# A Generalization of Calculus for Use with Continuous or Discrete Variables

**Jay Kaminsky**

**February 2013**

ABSTRACT: This document introduces a generalization of calculus that treats both continuous and discrete variables on an equal footing. This generalization of calculus was developed independently of the "Calculus on Time Scales'' literature but may be seen to have interesting overlap with it as well as with the "h-Calculus'' of the book *Quantum Calculus* by V. Kac and P. Cheung. As in the time scales literature, we first derive discrete analogues of all the common continuous calculus functions with an eye to maintaining as much similarity as possible between these discrete analogues and their continuous forebears. For example, in order to maintain the crucial property that the derivative of an exponential is a constant times itself, we replace the continuous exponential, $e^{ax}$, with a discrete function, $e_{\Delta x}(a,x) = [1+a\Delta x]^{\frac{x}{\Delta x}}$. Next, we develop a unified method of evaluating integrals of discrete variables. We discover that summations such as the Riemann Zeta Function, the Hurwitz Zeta Function, and the Digamma Function frequently appear in evaluating such integrals. Thus, we subsume these functions into a generalization of the natural logarithm, which we name "lnd(n,Δx,x)'', and evaluate many types of discrete variable integrals in terms of it. We provide a computer program, LNDX, to evaluate the lnd(n,Δx,x) function. LNDX is written in the UBASIC programming language, and it is downloadable from https://www.box.com/s/c59uvksuwomnw8rsyj92. Then, we develop a theory of control system analysis based on what we name a "$K_{\Delta x}$ Transform,'' which is related to the well-known Z Transform but has advantages beyond it. In closing, we highlight the fact that this document is structured somewhat like a textbook with many sample problems and solutions in the hope that it will be readily understood and found useful by those with even an undergraduate understanding of calculus and control systems.

## Author's Note

As mentioned in the Abstract, the research contained in this document was initially conceived and then largely developed independently of the mathematics and engineering research communities. However, it does have interesting overlap with two strands of contemporary research. For elaboration on this point as well as annotation of some key references in the literature, please read the Bibliographic Notes that preface the Bibliography at the end of this document.

# Table of Contents

# Introduction

Some years ago the thoughts came to mind, what would standard differential and integral Calculus (hereafter referred to simply as "Calculus") look like if $\Delta x$ did not go to zero and, if so, of what good would this non-infinitesimal $\Delta x$ calculus be? After some number of years, this calculus was developed and named Interval Calculus for the fact that its under function area is composed of $\Delta x$ intervals not of infinitesimal size. (The present Interval Calculus should not to be confused with the mathematics of the same name used to evaluate accumulated uncertainties in calculations.) Mathematically, this calculus worked. It found the area under curves that control engineers might call sample and hold shaped curves, the hold being for a $\Delta x$ interval. With care given to the selection of the Interval Calculus notation, its look and feel (mathematical manipulation) was found to be similar to that of Calculus. This led to its eventual application to the solution of difference equations and to difference equations in differential form. Later, from this work, the so-named $K_{\Delta x}$ Transform was conceptualized and developed.

This paper introduces both Interval Calculus and the $K_{\Delta x}$ Transform which are used to solve difference equations, difference equations in differential form, and problems involving sampled-data systems. Due to its similarity to Calculus, Interval Calculus needs little explanation. Its use becomes evident in the derivation and application of $K_{\Delta x}$ Transforms and in the many examples provided. The Interval Calculus functions which most often appear in the solution of non-infinitesimal $\Delta x$ equations and in the derivation of $K_{\Delta x}$ Transforms are presented in the tables in the Appendix.

A very important Interval Calculus function, the discrete lnd(n,$\Delta$x,x) function, has been derived and programmed in this paper. Amongst other things, this discrete function takes the place of the Calculus natural logarithm which has only limited application in discrete calculus.

It has been found that the discrete mathematics of Interval Calculus has applications to several advanced mathematics functions such as the Hurwitz Zeta Function, the Riemann Zeta Function, the Digamma Function, and the Polygamma Functions. The relationship of the Interval Calulus lnd(n,$\Delta$x,x) function to these functions is derived and presented in the Appendix at the end of this paper.

Below are some of the features and capabilities of Interval Calculus and the $K_{\Delta x}$ Transform:

1) Interval Calculus is a generalization of Calculus. Δx, in general, is not infinitesimal. If $\Delta x \to 0$, Interval Calculus is Calculus.

2) Interval Calculus is based on the function, $(1+a\Delta x)^{\frac{x}{\Delta x}}$ , where x = nΔx, n = integers, Δx = real or complex value.

   Comment - Calculus is based on the function, $e^{ax}$, where x = real or complex numbers. Difference mathematics is based on the function, $A^x$, where x = integers. A,a = real or complex constants.

3) Interval Calculus functions contain the x interval constant, Δx, where Δx may be any value, real or complex. For real values of x, the functions are sample and hold shaped waveforms. An example of an Interval Calculus function is shown below.

A Sample and Hold Shaped Interval Calculus function

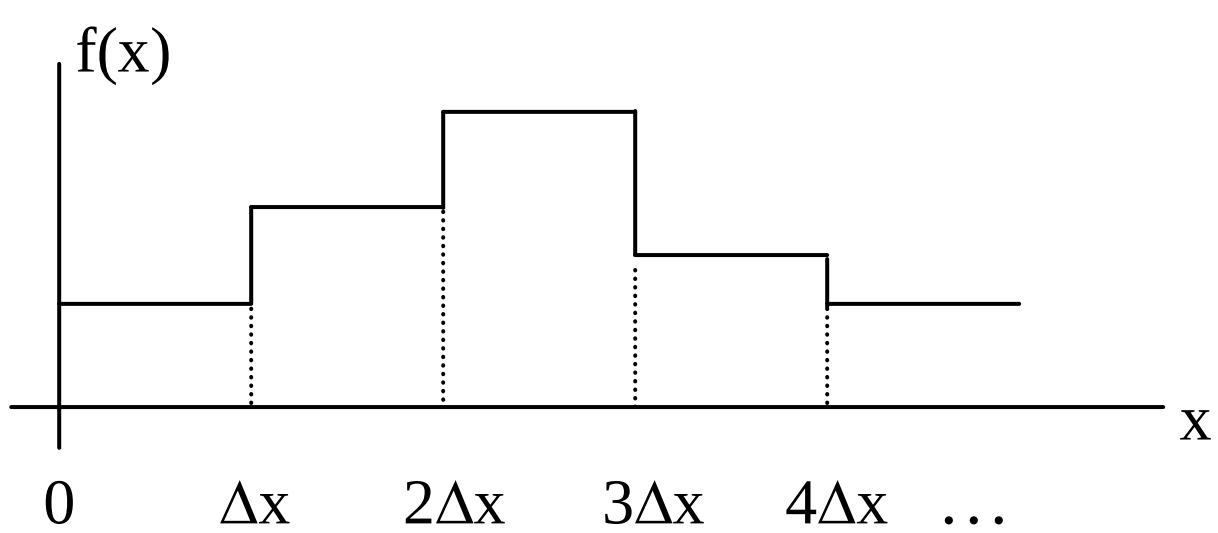

Note that the independent variable of an Interval Calculus function is discrete, x = nΔx , n = integers. Note, also, that as $\Delta x \to 0$ the Interval Calculus function becomes a continuous Calculus function where x = all real values.

4) Interval Calculus functions and equations are similar and related to Calculus functions. The similarities are in both notation and in mathematical operations. For example:

   a. $e_{\Delta x}(a,x) = (1+a\Delta x)^{\frac{x}{\Delta x}}$     Comment $\lim_{\Delta x \to 0}(1+a\Delta x)^{\frac{x}{\Delta x}} = e^{ax}$

   b. $\sin_{\Delta x}(a,x) = \dfrac{(1+ja\Delta x)^{\frac{x}{\Delta x}} - (1-ja\Delta x)^{\frac{x}{\Delta x}}}{2j}$

   c. $\cos_{\Delta x}(a,x) = \dfrac{(1+ja\Delta x)^{\frac{x}{\Delta x}} + (1-ja\Delta x)^{\frac{x}{\Delta x}}}{2}$

   d. $D_{\Delta x}\sin_{\Delta x}(a,x) = a\cos_{\Delta x}(a,x)$

   e. $D_{\Delta x}\cos_{\Delta x}(a,x) = -a\sin_{\Delta x}(a,x)$

   f. $D_{\Delta x}[u(x)v(x)] = v(x)D_{\Delta x}u(x) + D_{\Delta x}v(x)u(x+\Delta x)$

5) Interval Calculus performs discrete differentiation which is defined as follows:

$$D_{\Delta x}f(x) = \frac{f(x+\Delta x) - f(x)}{\Delta x}$$

6) Interval Calculus discrete differentiation is similar to Calculus differentiation
For example:

a. $D_{\Delta x}[x(x-\Delta x)] = 2x$
b. $D_{\Delta x}[x(x-\Delta x)(x-2\Delta x)] = 3x(x-\Delta x)$
c. $D_{\Delta x}e_{\Delta x}(a,x) = ae_{\Delta x}(a,x)$

7) Interval Calculus performs discrete integration which is defined as follows:

$$_{\Delta x}\int_{x_1}^{x_2} f(x)\Delta x = \Delta x \; _{\Delta x}\sum_{x=x_1}^{x_2-\Delta x} f(x)$$

8) Interval Calculus discrete integration calculates the area under sample and hold shaped waveforms.

9) Interval Calculus discrete integration is similar to Calculus integration
For example:

1. $_{\Delta x}\int_{x_1}^{x_2} x\,\Delta x = \frac{x(x-\Delta x)}{2}\Big|_{x_1}^{x_2}$
2. $_{\Delta x}\int_{x_1}^{x_2} e_{\Delta x}(a,x)\,\Delta x = \frac{e_{\Delta x}(a,x)}{a}\Big|_{x_1}^{x_2}$

10) Interval Calculus discrete integration can integrate by parts. The discrete integration by parts equation is as follows:

$$_{\Delta x}\int_{x_1}^{x_2} v(x)D_{\Delta x}u(x)\Delta x = u(x)v(x)\Big|_{x_1}^{x_2} - \; _{\Delta x}\int_{x_1}^{x_2} D_{\Delta x}v(x)u(x+\Delta x)\Delta x$$

11) Interval Calculus can be applied to differential difference equations, difference equations, area and summation evaluations, and calculations involving sampled-data sytems.

12) The Method of Undetermined Coefficients, the Method of Variation of Parameters, the Transform Method, and the Method of Related Functions are used by Interval Calculus to solve discrete differential difference equations.

13) Two series have been derived to calculate summations involving the function, $\frac{1}{x^n}$ . They are as follows:

1. $$\mathrm{lnd}(n,\Delta x,x) \approx \gamma + \ln\left(\frac{x}{\Delta x} - \frac{1}{2}\right) + \sum_{m=1}^{\infty} \frac{(2m-1)!\, C_m}{(2m+1)!\, 2^{2m} \left(\frac{x}{\Delta x} - \frac{1}{2}\right)^{2m}} \quad , \quad n = 1$$

2. $$\mathrm{lnd}(n,\Delta x,x) \approx - \sum_{m=0}^{\infty} \frac{\Gamma(n+2m-1)\left(\frac{\Delta x}{2}\right)^{2m} C_m}{\Gamma(n)(2m+1)! \left(x - \frac{\Delta x}{2}\right)^{n+2m-1}} + K \quad , \quad n \neq 1$$

The accuracy of both series increases rapidly for increasing $\left|\frac{x}{\Delta x}\right|$

<u>Comments</u> - These series are not derived from the Euler-Maclauren Series.

$$\sum_{\Delta x \; x=x_i}^{\pm\infty} \frac{1}{x^n} = \frac{1}{\Delta x} \mathrm{lnd}(n,\Delta x,x_i) \quad , \quad \mathrm{Re}(n)>1$$

14) A computer program, LNDX, has been written which uses the two equations presented in 13) to evaluate the lnd(n,Δx,x) function for all real or complex values of n, Δx, and x including the special case where n=1.

15) Summations such as the following can be evaluated using the lnd(n,Δx,x) function.

1. $$\sum_{\Delta x \; x=x_1}^{x_2} \frac{1}{x^n} = \pm \left[ -\frac{1}{\Delta x} \mathrm{lnd}(n,\Delta x,x) \Big|_{x_1}^{x_2+\Delta x} \right] , \quad + \text{ for } n\neq 1 \, , \; - \text{ for } n=1$$

2. $$\sum_{\Delta x \; x=x_i}^{\pm\infty} \frac{1}{x^n} = \frac{1}{\Delta x} \mathrm{lnd}(n,\Delta x,x_i) \quad , \quad \mathrm{Re}(n)>1$$

3. $$\sum_{\Delta x \; x=x_i}^{\pm\infty} (-1)^{\frac{x-x_i}{\Delta x}} \frac{1}{x^n} = \pm \left[ -\mathrm{lnd}(n, 2\Delta x, x_i+\Delta x) + \mathrm{lnd}(n,2\Delta x,x_i) \right] , \quad + \text{ for } n\neq 1 \, , \; - \text{ for } n=1 \, , \; \mathrm{Re}(n)>0$$

4. $$\sum_{\Delta x\; x=\mp\infty}^{\pm\infty} (-1)^{\frac{x-x_i}{\Delta x}} \frac{1}{x^n} = \pm \frac{1}{2\Delta x}\left[-\text{lnd}(n,2\Delta x,x_i+\Delta x)+\text{lnd}(n,2\Delta x,x_i)-\text{lnd}(n,-2\Delta x,x_i-2\Delta x)+\text{lnd}(n,-2\Delta x,x_i-\Delta x)\right]$$

$x = x_i + m\Delta x$, m = integers , + for $n\neq 1$ , – for n=1 , Re(n)>0

16) Series have been derived to evaluate the following summations:

1. $$\sum_{\Delta x\; x=x_1}^{x_2} f(x) = \frac{1}{\Delta x}\int_{x_1-\frac{\Delta x}{2}}^{x_2+\frac{\Delta x}{2}} f(x)dx + \sum_{m=1}^{\infty} B_m \Delta x^{2m-1} \frac{d^{2m-1}}{dx^{2m-1}} f(x)\Big|_{x_1-\frac{\Delta x}{2}}^{x_2+\frac{\Delta x}{2}}$$

2. $$\sum_{\Delta x\; x=x_1}^{x_2} (-1)^{\frac{x-x_1}{\Delta x}} f(x) = -\frac{1}{2} f(x)\Big|_{x_1}^{x_2+\Delta x} + \sum_{m=1}^{\infty} H_m (2\Delta x)^{2m-1} \frac{d^{2m-1}}{dx^{2m-1}} f(x)\Big|_{x_1}^{x_2+\Delta x}$$

Comment - These series differ from the Euler-Maclauren Series

17) Interval Calculus has a unique function, lnd(n,Δx,x). The lnd(n,Δx,x) function is to Interval Calculus what the natural log function is to Calculus. The lnd(n,Δx,x) function is defined as follows:

$$D_{\Delta x}\text{lnd}(n,\Delta x,x) \equiv D_{\Delta x}\ln_{\Delta x} x = +\frac{1}{x}\ ,\quad n=1$$

It has additional importance.

$$D_{\Delta x}\text{lnd}(n,\Delta x,x) = -\frac{1}{x^n}\ ,\qquad n\neq 1$$

18) The summation of the function, $\frac{1}{(x-a)^n}$, can be calculated using the lnd(n,Δx,x) function.

The calculation equations are as follows:

$$\sum_{\Delta x\; x=x_1}^{x_2-\Delta x} \frac{1}{(x-a)^n} = \pm\frac{1}{\Delta x}\,\text{lnd}(n,\Delta x,x-a)\Big|_{x_1}^{x_2}\ ,\ +\text{ for } n=1\ ,\ -\text{ for } n\neq 1$$

any $\frac{1}{0}$ term excluded

and

$$\sum_{\Delta x\; x=x_i}^{\pm\infty} \frac{1}{(x-a)^n} = \frac{1}{\Delta x}\,\text{lnd}(n,\Delta x,x_i-a)\ ,\ \text{Re}(n)>1$$

any $\frac{1}{0}$ term excluded

19) At least 60 functions and summations commonly used by engineers and scientists can be evaluated using the $\mathrm{lnd}(n,\Delta x,x)$ function.

20) The Riemann Zeta Function can be evaluated by the $\mathrm{lnd}(n,\Delta x,x)$ function. The equations are as follows:

$$\zeta(n) = \sum_{x=1}^{\infty}{}_{1}\frac{1}{x^n} = \mathrm{lnd}(n,1,1)\ , \quad \mathrm{Re}(n) > 1$$

$$\zeta(n) = \mathrm{lnd}(n,1,1)\ , \quad n \neq 1$$

21) The Hurwitz Zeta Function can be evaluated by the $\mathrm{lnd}(n,\Delta x,x)$ function. The equations are as follows:

$$\zeta(n,x_i) = \sum_{x=x_i}^{\infty}{}_{1}\frac{1}{x^n} = \mathrm{lnd}(n,1,x_i)\ , \quad \mathrm{Re}(n) > 1\ , \quad \text{any } \frac{1}{0} \text{ term is excluded}$$

$$\zeta(n,x_i) = \mathrm{lnd}(n,1,x_i)\ , \quad n \neq 1$$

22) Interval Operational Calculus is based on the $K_{\Delta x}$ Transform, a generalization of the Laplace Transform. This transform is related to the Laplace Transform and the Z Transform and can be used to solve differential difference equations or difference equations. The equation for the $K_{\Delta x}$ Transform is as follows:

$$K_{\Delta x}[f(x)] = {}_{\Delta x}\int_{0}^{\infty} (1+s\Delta x)^{-\left(\frac{x+\Delta x}{\Delta x}\right)} f(x)\Delta x = {}_{\Delta x}\int_{0}^{\infty} e_{\Delta x}(s,-x-\Delta x)f(x)\Delta x$$

<u>Note</u> – For $\Delta x \to 0$, the $K_{\Delta x}$ Transform becomes the Laplace Transform.

23) The Inverse $K_{\Delta x}$ Transform is as follows:

$$K_{\Delta x}^{-1}[F(s)] = F(x) = \frac{1}{2\pi j}\oint_{c} [1+s\Delta x]^{\frac{x}{\Delta x}-1} F(s)\, ds$$

$$s = \frac{e^{(\gamma+jw)\Delta x} - 1}{\Delta x}\ , \quad -\frac{\pi}{\Delta x} < w < +\frac{\pi}{\Delta x}$$

C is a circular contour in the complex plane defined by the above equation for s.
Residue theory can be used to find the Inverse $K_{\Delta x}$ Transform.

24) $K_{\Delta x}$ Transforms look similar to Laplace Transforms. Some additional $\Delta x$ terms may be added. If $\Delta x \to 0$, $K_{\Delta x}$ Transforms become Laplace Transforms.

25) The $K_{\Delta x}$ Transform is related to the Z Transform by the following substitutions:

1) $K_{\Delta x}$ Transform to Z Transform Conversion

$$Z[f(x)] = F(z) = \frac{z}{T} F(s)\Big|_{s = \frac{z-1}{T}}$$

$T = \Delta x$ sampling period

$x = nT$

$n = 0,1,2,3,\ldots$

$Z[f(x)] = F(z)$ Z Transform

$K_{\Delta x}[f(x)] = F(s)$ $K_{\Delta x}$ Transform

For conversion from the Inverse $K_{\Delta x}$ Transform to the Inverse Z Transform
The complex plane integration contour changes from

$$s = \frac{e^{(\gamma+jw)\Delta x} - 1}{\Delta x} \quad \text{to } z = e^{(\gamma+jw)T}$$

where

$\gamma$ = positive real constant

$$-\frac{\pi}{\Delta x} \le w < +\frac{\pi}{\Delta x}$$

2) Z Transform to $K_{\Delta t}$ Transform Conversion

$$K_{\Delta x}[f(x)] = F(s) = Tz^{-1} F(z)\Big|_{z = 1+s\Delta x}$$

$\Delta x = T$ sampling period

$x = n\Delta x$

$n = 0, 1, 2, 3, \ldots$

$K_{\Delta x}[f(x)] = F(s)$ $K_{\Delta x}$ Transform

$Z[f(x)] = F(z)$ Z Transform

For conversion from the Inverse Z Transform to the Inverse $K_{\Delta x}$ Transform
The complex plane integration contour changes from

$$z = e^{(\gamma+jw)T} \quad \text{to } s = \frac{e^{(\gamma+jw)\Delta x} - 1}{\Delta x}$$

where

$\gamma$ = positive real constant

$$-\frac{\pi}{\Delta x} \le w < +\frac{\pi}{\Delta x}$$

26) $K_{\Delta x}$ Transforms can be used to analyze discrete control system stability using Nyquist and Bode methods.

27) For $K_{\Delta x}$ Transform discrete control system stability analysis, control system stability occurs if all 1+A(s) zeros lie within a Critical Circle within the left half of the s plane. The Critical Circle has its center at $s = -\frac{1}{\Delta x}$ and has a radius of $\frac{1}{\Delta x}$.

<u>Comment</u> – For $\Delta x \to 0$, the Critical Circle becomes the left half of the s plane.

28) Interval Calculus can be used to calculate the area enclosed within a discrete closed contour in the complex plane.

29) Using $K_{\Delta x}$ Transforms to solve differential difference equations is very similar to using Laplace Transforms to solve differential equations. In fact, if a Calculus homogenous differential equation and an Interval Calculus homogeneous differential difference equation have the same value coefficients and initial conditions, their transform equations will look the same. Also, the solution functions become the same for $\Delta x \rightarrow 0$. See the table on the following page.

30) Interval Calculus may be considered to be an infinite set of calculi, a calculus for every $\Delta x$ interval value. Calculus is one element of this set, the calculus where $\Delta x \rightarrow 0$. Each calculus has its own unique functions. The functions, $e^{ax}$, sin(ax), cos(ax), ln(x), etc. belong to the Calculus element of the set. The functions of the finite $\Delta x$ non-Calculus elements of the set are $e_{\Delta x}(a,x)$, $\sin_{\Delta x}(b,x)$,$\cos_{\Delta x}(b,t)$, lnd(n,$\Delta x$,x) etc. A Calculus function may belong to a non-Calculus element of the set only if it is an identity of a function belonging to that non-Calculus element (For example:

$e^{bx} = e_{\Delta x}(a,x)$ where $b = \ln(1+a\Delta x)^{\frac{1}{\Delta x}}$, $x = 0, \Delta x, 2\Delta x, 3\Delta x, \ldots$).

**The Solution of a Homogeneous Differential Equation using the Laplace Transform**

For the homogeneous differential equation
$D^n y(x) + a_{n-1}D^{n-1}y(x) + a_{n-2}D^{n-2}y(x) + \ldots + a_1Dy(x) + a_0y(x) = 0$

With initial conditions
$y(0) = c_0,\ Dy(0) = c_1,\ \ldots,\ D^{n-2}y(0) = c_{n-2},\ D^{n-1}y(0) = c_{n-1}$

where

$$Dy(x) = \lim_{\Delta x\to 0} \frac{y(x+\Delta x)-y(x)}{\Delta x}$$

x = continuous variable

The Laplace Transform is:
$(s^n + a_{n-1}s^{n-1}+ a_{n-2}s^{n-2}+ \ldots + a_1s + a_0)y(x) = b_{n-1}s^{n-1} + b_{n-2}s^{n-2} + \ldots + b_1s + b_0$
the initial conditions are in this term

The solution functions are:

**The Solution of a Homogeneous Differential Difference Equation using the $K_{\Delta x}$ Transform (The Generalized Laplace Transform)**

For the homogeneous differential difference equation
$D_{\Delta x}{}^n y(x) + a_{n-1}D_{\Delta x}{}^{n-1}y(x) + a_{n-2}D_{\Delta x}{}^{n-2}y(x) + \ldots + a_1D_{\Delta x}y(x) + a_0y(x) = 0$

With initial conditions
$y(0) = c_0,\ D_{\Delta x}y(0) = c_1,\ \ldots,\ D_{\Delta x}{}^{n-2}y(0) = c_{n-2},\ D_{\Delta x}{}^{n-1}y(0) = c_{n-1}$

where

$$D_{\Delta x}y(x) = \frac{y(x+\Delta x)-y(x)}{\Delta x}$$

$x = 0, \Delta x, 2\Delta x, 3\Delta x, \ldots$

The $K_{\Delta x}$ Transform (The Generalized Laplace Transform) is:
$(s^n + a_{n-1}s^{n-1}+ a_{n-2}s^{n-2}+ \ldots + a_1s + a_0)y(x) = b_{n-1}s^{n-1} + b_{n-2}s^{n-2} + \ldots + b_1s + b_0$
the initial conditions are in this term

The solution functions are:

| Solution Function | Laplace Transform | Solution Function | $K_{\Delta x}$ Transform |
|---|---|---|---|
| e(ax) | $\frac{1}{s-a}$ | $e_{\Delta x}(a,x)$ | $\frac{1}{s-a}$ |
| sin(bx) | $\frac{b}{s^2+b^2}$ | $\sin_{\Delta x}(b,x)$ | $\frac{b}{s^2+b^2}$ |
| cos(bx) | $\frac{s}{s^2+b^2}$ | $\cos_{\Delta x}(b,x)$ | $\frac{s}{s^2+b^2}$ |
| e(ax)sinbx | $\frac{b}{(s-a)^2+b^2}$ | $e_{\Delta x}(a,x)\sin_{\Delta x}(\frac{b}{1+a\Delta x}, x)$<br>for $1+a\Delta x \neq 0$<br>$[b\Delta x]^{\frac{x}{\Delta x}} \sin \frac{\pi x}{2\Delta x}$<br>for $1+a\Delta x = 0$ | $\frac{b}{(s-a)^2+b^2}$ |
| e(ax)cosbx | $\frac{s-a}{(s-a)^2+b^2}$ | $e_{\Delta x}(a,x)\cos_{\Delta x}(\frac{b}{1+a\Delta x}, x)$<br>for $1+a\Delta x \neq 0$<br>$[b\Delta x]^{\frac{x}{\Delta x}} \cos \frac{\pi x}{2\Delta x}$<br>for $1+a\Delta x = 0$ | $\frac{s-a}{(s-a)^2+b^2}$ |

<u>Note</u> – For the functions specified, the transforms are the same. Tables of transforms are presented in the Appendix.

31) The $K_{\Delta x}$ Transform has some interesting characteristics which can be advantageous in solving differential difference equations.

1. $K_{\Delta x}$ transforms can be applied just as easily to difference equations as to differential difference equations.
2. Major Interval Calculus functions such as $e_{\Delta x}(a,x)$, $\sin_{\Delta x}(b,x)$, etc. are, in reality, sets of related functions of x which differ depending on the value of the constant, $\Delta x$. Interestingly, the $K_{\Delta x}$ Transform of these functions remains the same for all values of $\Delta x$.
3. To each Interval Calculus function, there is a related Calculus function (for example: e(ax) for $e_{\Delta x}(a,x)$, sin(bx) for $\sin_{\Delta x}(b,x)$ etc.) which has the same transform. This Calculus function is the Interval Calculus function where $\Delta x$ goes to zero.
4. Because discrete Interval Calculus functions have related Calculus functions and their transforms are equal, Calculus differential equations can be used to solve related differential difference equations.
5. The $K_{\Delta x}$ Transform is a generalized Laplace Transform which is applicable to discrete mathematics.

Consider a transform solution to a simple differential difference equation. See the table below.

| Find the solution to the differential difference equation, $D_{\Delta x}y(x) + ay(x) = 0$, where $\Delta x$ can equal any value including the case where $\Delta x \to 0$. Use both the $K_{\Delta x}$ Transform and the Laplace Transform for comparison purposes. | | | | |
|---|---|---|---|---|
| **Type of Transform** | **Equation Transform** | **Solution Transform** | **Solution** | **Comments** |
| $K_{\Delta x}$ Transform (Generalized Laplace Transform) | sy(s)-ay(s)=0 | $y(s) = y(0)\frac{1}{s+a}$ | $y(x) = y(0)e_{\Delta x}(-a,x)$ | $\Delta x \neq 0$<br>$e_{\Delta x}(-a,x) = (1-a\Delta x)^{\frac{x}{\Delta x}}$ |
| Laplace Transform | sy(s)-ay(s)=0 | $y(s) = y(0)\frac{1}{s+a}$ | $y(x) = y(0)e(-ax)$ | $\Delta x \to 0$<br>$e(-ax) = \lim_{\Delta x \to 0}(1-a\Delta x)^{\frac{x}{\Delta x}}$ |

32) Using Interval Calculus, differential difference equation stability, as a function of $\Delta x$, is easily found graphically when the equation coefficients are constant and not a function of $\Delta x$. Consider the simple differential difference equation in 30), $D_{\Delta x}y(x) + ay(x) = 0$.

$$\frac{y(x+\Delta x)-y(x)}{\Delta x} - ay(x) = 0 \quad , \quad x = 0, \Delta x, 2\Delta x, 3\Delta x, \ldots$$

Graphically, find the range of $\Delta x$ for which the above differential difference equation is stable.

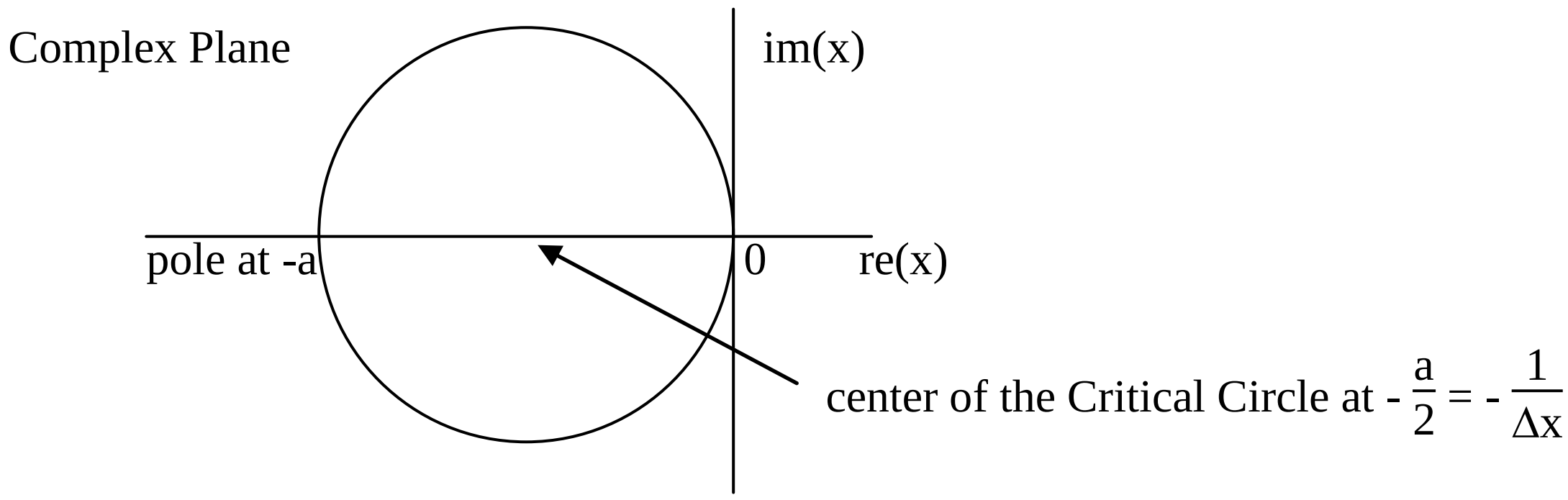

The radius of the Critical Circle is $\frac{1}{\Delta x} = \frac{a}{2}$. Then $\Delta x$ can not be greater than $\frac{2}{a}$ or the pole at -a will lie outside the circle of stability.

Thus the range of $\Delta x$ for stability is $0 \le \Delta x < \frac{2}{a}$

As a check, y(x), where $\Delta x = \frac{2}{a}$, should be oscillatory

$y(x) = e_{\Delta x}(-a,x) = (1-a\Delta x)^{\frac{x}{\Delta x}} = (1-a[\frac{2}{a}])^{\frac{x}{\Delta x}} = (-1)^{\frac{x}{\Delta x}}$ which is oscillatory good check

33) The $K_{\Delta t}$ Transform, the Z Transform, and the Laplace Transform have been shown to be very closely related.

1. $C^*(s) = F^*(s)\,G^*(s)$ $K_{\Delta t}$ Transform equation
2. $C^*(z) = F^*(z)\,G^*(z)$ Z Transform equation
3. $C^*(\mathit{s}) = F^*(\mathit{s})\,G^*(\mathit{s})$ Laplace Transform equation
4. $L[f^*(t)] = \frac{1}{\mathit{s}}\, sK_{\Delta t}[f^*(t)]\,\Big|_{s\,=\,\frac{e^{\mathit{s}\Delta t}-1}{\Delta t}}$ , $K_{\Delta t}$ Transform to Laplace Transform conversion
5. $L[f^*(t)] = \frac{1}{\mathit{s}}\,\frac{z-1}{z}\,Z[f^*(t)]\Big|_{z\,=\,e^{\mathit{s}T}}$ , Z Transform to Laplace Transform conversion

34) The $K_{\Delta x}$ Transform, the Z Transform, and the Laplace Transform have related Convolution Equations.

1. $$K_{\Delta x}[\underset{\Delta\lambda}{\int_0^x} f(x-\lambda-\Delta\lambda)g(\lambda)\Delta\lambda] = K_{\Delta x}[f(x)]K_{\Delta x}[g(x)]$$ , $K_{\Delta x}$ Transform Convolution Equations

   or

   $$K_{\Delta x}[\underset{\Delta\lambda}{\int_0^x} f(\lambda)g(x-\lambda-\Delta\lambda)\Delta\lambda] = K_{\Delta x}[f(x)]K_{\Delta x}[g(x)]$$

2. $$Z[\frac{1}{T}\underset{T}{\int_0^{x+T}} f(x-\lambda)g(\lambda)\Delta\lambda] = Z[f(x)]Z[g(x)]$$ , Z Transform Convolution Equations

   or

   $$Z[\frac{1}{T}\underset{T}{\int_0^{x+T}} f(\lambda)g(x-\lambda)\Delta\lambda] = Z[f(x)]Z[g(x)]$$

3. $$L[\underset{\Delta\lambda}{\int_0^x} f(x-\lambda)g(\lambda)\Delta\lambda] = L[f(x)]L[g(x)]$$ , Laplace Transform Convolution Equations

   or

   $$L[\underset{\Delta\lambda}{\int_0^x} f(\lambda)g(x-\lambda)\Delta\lambda] = L[f(x)]L[g(x)]$$

Comment - Discrete variable Duhamel formulas have also been derived.

35) Interval Calculus is a generalization of Calculus where the non-infinitesimal $\Delta x$ x interval has been reinserted.

Comment - Some applications of this generalization of Calculus to discrete mathematics are shown in Chapter 7.

36) There is a relationship between Interval Calculus and the Mathematics of Finite Differences. The commonly used Mathematics of Finite Differences operators, $H_{\Delta x}$ and $\Delta_{\Delta x}$, can be expressed in terms of the Interval Calculus discrete differentiation operator, $D_{\Delta x}$. The Interval Calculus discrete differentiation operator and the Mathematics of Finite Differences differential operator, $D_{\Delta x}$, are the same.

1) $H_{\Delta x}f(x) = f(x+\Delta x) = (\Delta x D_{\Delta x} + 1)f(x)$, Function variable forward increment
2) $\Delta_{\Delta x}f(x) = f(x+\Delta x) - f(x) = (H_{\Delta x} - 1)f(x) = \Delta x D_{\Delta x}f(x)$, Function difference
3) $D_{\Delta x}f(x) = \frac{f(x+\Delta x) - f(x)}{\Delta x} = \frac{\Delta_{\Delta x}f(x)}{\Delta x}$ , Discrete differentiation of a function

# CHAPTER 1

## The Solution of Difference Equations and Difference Equations in Differential Form Using Interval Calculus and the $K_{\Delta x}$ Transform

# Interval Calculus and the $K_{\Delta x}$ Transform

## Chapter 1: Commentary

Chapter 1 introduces a number of important and useful features of discrete Interval Calculus. These features are presented with some applications and explanation. Though some derivations are presented in this chapter, most of the important derivations are left to succeeding chapters where they can be presented in greater detail.

## Section 1.1: Notation and Basics

In the following discussion, when comparing Interval Calculus with Calculus, $f_{\Delta x}(x)$ is used to denote an Interval Calculus function (which contains the x interval, $\Delta x$) and f(x) is used to denote a Calculus function (where any derivation $\Delta x$ has gone to zero). Interval Calculus, a variation on the difference equation calculus that is commonly used, has one major difference. Instead of being structured on the functions $f(x) = A^x$ and $f(x) = A^{jx}$ , A being a constant and $j = \sqrt{-1}$ , it is structured on the respective functions $f_{\Delta x}(x) = (1+a\Delta x)^{\frac{x}{\Delta x}}$ and $f_{\Delta x}(x) = (1+ja\Delta x)^{\frac{x}{\Delta x}}$. There is a one to one relationship between each Interval Calculus function, $f_{\Delta x}(x)$, and a Calculus function, f(x). For example:

$f(x) = \lim_{\Delta x\to 0}(1+a\Delta x)^{\frac{x}{\Delta x}} = e^{ax}$ and $f(x) = \lim_{\Delta x\to 0}(1+ja\Delta x)^{\frac{x}{\Delta x}} = e^{jax}$. These functions, $f_{\Delta x}(x)$ and f(x) are related and become one in the same as $\Delta x\to 0$. Great care has been given to the naming and definition of the functions, $f_{\Delta x}(x)$, which arise from the solution of difference equations and difference equations of the differential form. For example, $f_{\Delta x}(x) = e_{\Delta x}(a,x) = (1+a\Delta x)^{\frac{x}{\Delta x}}$ which has as its $\Delta x\to 0$ related function, $f(x) = e^{ax}$ and $f_{\Delta x}(x) = \sin_{\Delta x}(a,x) = \dfrac{(1+ja\Delta x)^{\frac{x}{\Delta x}} - (1-ja\Delta x)^{\frac{x}{\Delta x}}}{2j}$ which has as its $\Delta x \to 0$ related function, $f(x) = \sin ax$ .

Differentiation in Interval Calculus is performed using the following relationship:

$$D_{\Delta x}f_{\Delta x}(x) = \frac{f_{\Delta x}(x+\Delta x) - f_{\Delta x}(x)}{\Delta x} \qquad (1.1\text{-}1)$$

The symbol, $D_{\Delta x,}$ is used as the differential operator. Difference equation derivative operations, when applied to the $f_{\Delta x}(x)$ functions, yield instantly recognizable results. For example, $D_{\Delta x}\sin_{\Delta x}(a,x) = a\cos_{\Delta x}(a,x)$ and $D_{\Delta x}\cos_{\Delta x}(a,x) = -a\sin_{\Delta x}(a,x)$. The solutions to discrete differential equations are recognizable too. For example, the discrete differential equation, $D^2_{\Delta x}y(x) – y(x) = 0$ has as a solution $k_1\sin_{\Delta x}(a,x) + k_2\cos_{\Delta x}(a,x)$. Even the $K_{\Delta x}$ Transform for this case would be recognizable as $y(s) = \frac{k_1 a}{s^2+a^2} + \frac{k_2 s}{s^2+a^2}$.

To maintain the strong similarities of Interval Calculus with Calculus, the following notation was adopted:

$$ {}_{\Delta x}\sum_{x=x_1}^{x_2} \equiv \sum_{x=\{X\}} \quad \text{where } \{X\} = x_1, x_1+\Delta x, x_1+2\Delta x, \ldots, x_2-\Delta x, x_2 \tag{1.1-2}$$

$$ {}_{\Delta x}\int_{x_1}^{x_2+\Delta x} \equiv (\Delta x)\left(\sum_{x=\{X\}}\right) \quad \text{where } \{X\} = x_1, x_1+\Delta x, x_1+2\Delta x, \ldots, x_2-\Delta x, x_2 \tag{1.1-3}$$

The leading subscript Δx indicates the interval used in the summation. Although other notation could certainly be chosen, it was found that the proposed notation leads to mathematical results with an appearance similar to standard Calculus. The leading subscript Δx, in the preferred form, may remain as a general value, Δx, or as a specific numerical value for Δx. These preferred notations are used within this paper. Note that integration is the inverse of differentiation.

## Section 1.2: The $K_{\Delta t}$ Transform

The $K_{\Delta x}$ Transform is derived in Chapter 5. Control system analysis is possible using the Interval Calculus $K_{\Delta x}$ Transform. Though the Z Transform provides a convenient way to evaluate and analyze discrete calculus problems, there is still reason for the use of the $K_{\Delta x}$ Transform instead. In particular, the form of the $K_{\Delta x}$ Transform lends itself well to discrete Interval Calculus integration and the $K_{\Delta x}$ Transform is very similar to the Laplace Transform. In fact, the $K_{\Delta x}$ Transform may be considered to be an extension of the Laplace Transform to discrete mathematics. The similarity between the $K_{\Delta x}$ transform and the Laplace Transform can be helpful in the analysis and solution of mathematical problems of a discrete variable. In this section and in the following, the $K_{\Delta x}$ Transform is presented with various applications.

The $K_{\Delta x}$ Transform, is now presented:

$$K_{\Delta x}[f(x)] = {}_{\Delta x}\int_{0}^{\infty} (1+s\Delta x)^{-(\frac{x+\Delta x}{\Delta x})} f(x)\Delta x \qquad (1.2\text{-}1)$$

**or**

$$K_{\Delta x}[f(x)] = {}_{1}\sum_{n=0}^{\infty} (1+s\Delta x)^{-n-1} f(n\Delta x)\Delta x \qquad (1.2\text{-}2)$$

Comment - Eq 1.2-1 and Eq 1.2-2 represent the sum of the $K_{\Delta x}$ Transforms of consecutive unit amplitude pulses of width $\Delta x$ multiplied by an f(x) weighting function. See the equation below from which both equations are derived.

$$K_{\Delta x}[f(x)] = {}_{1}\sum_{n=0}^{\infty} f(n\Delta x)K_{\Delta x}[U(x-n\Delta x)-U(x-n\Delta x-\Delta x)]$$

Described a different way, Eq 1.2-1 and Eq 1.2-2 represent the sum of the $K_{\Delta x}$ Transforms of consecutive unit area pulses of width $\Delta x$ multiplied by an f(x) weighting function. See the equation below from which both equations are derived.

$$K_{\Delta x}[f(x)] = \Delta x\,{}_{1}\sum_{n=0}^{\infty} f(n\Delta x)K_{\Delta x}[\frac{1}{\Delta x}\{U(x-n\Delta x)-U(x-n\Delta x-\Delta x)\}]$$

Note that the preferred notation described in the previous section is here used for integration. The transform may be used to conveniently solve difference equations and difference equations in differential form. It is similar to, and related to, the Laplace Transform. In particular, it has the following interesting characteristics:

1. The $K_{\Delta x}$ Transforms are functions not only of a variable s (similar to the Laplace Transform) but also of the x increment, $\Delta x$. Similarly, the functions of x which these transforms most often represent are functions not only of the variable, x, but of the x increment, $\Delta x$, also.

2. As $\Delta x \rightarrow 0$ the $K_{\Delta x}$ Transform becomes the Laplace Transform:

$$L[f(x)] = \int_0^{\infty} e^{-sx} f(x)dx \qquad (1.2\text{-}3)$$

Note the similarity between the two transforms.

3. Mathematical operations performed on the $K_{\Delta x}$ Transforms are primarily algebraic as they are with Laplace Transforms. The Interval Calculus derivative operation, $D_{\Delta x,}$ when applied to the $f_{\Delta x}(x)$ functions, yields instantly recognizable results. For example, $D_{\Delta x}\sin_{\Delta x}(a,x) = a\cos_{\Delta x}(a,x)$ and $D_{\Delta x}\cos_{\Delta x}(a,x) = -a\sin_{\Delta x}(a,x)$. Even the solutions to difference equations in differential form are instantly recognizable. For example, the difference equation in differential form, $D_{\Delta x}y(x) + ay(x) = 0$ has as a solution $ke_{\Delta x}(a,x.)$. For this case, the $K_{\Delta x}$ Transform would be $y(s) = \frac{k}{s+a}$.

4. Z Transforms have often been the method of preference for solving difference equations and difference equations in the differential form. The $K_{\Delta x}$ Transform is another option. Its close relationship to the Laplace Transform can be helpful for controls engineers. Instead of looking for poles in the complex right half plane, one would look for poles outside the circle with its center at $-\frac{1}{\Delta x}$ with a radius of $\frac{1}{\Delta x}$.

$K_{\Delta x}$ Transforms, Interval Calculus functions, $f_{\Delta x}(x)$, and the difference and differential operations on these functions are listed in tables in the appendix. A few proofs for the formulae provided in the appendix are shown in this paper. The proofs are not difficult and often proceed from two important sources. One source is Interval Calculus differentiation,

$$D_{\Delta x}u(x) = \frac{u(x+\Delta x) - u(x)}{\Delta x} \qquad (1.2\text{-}4)$$

with its inverse being used for Interval Calculus integration. The second source, used for Interval Calculus integration, is the integration by parts formula,

$$_{\Delta x}\int_{x_1}^{x_2} v(x)D_{\Delta x}u(x)\Delta x = u(x)v(x)\Big|_{x_1}^{x_2} - {}_{\Delta x}\int_{x_1}^{x_2} D_{\Delta x}v(x)u(x+\Delta x)\Delta x \qquad (1.2\text{-}5)$$

which is often very helpful. The two equations above, with the notation selected, should appear quite familiar. This is common in Interval Calculus. Note that if $\Delta x \to 0$ the first equation becomes the Calculus derivative and the second equation becomes the Calculus integration by parts equation. In general, if $f_{\Delta x}(x)$ is an Interval Calculus equation and f(x) is its related equation, then $\lim_{\Delta x \to 0} f_{\Delta x}(x) = f(x)$. In effect, Calculus is Interval Calculus where the limit of Δx goes to zero not at the beginning but at the end of any derivation process.

## Section 1.3: Examples

For demonstration purposes the derivative, integral, and $K_{\Delta x}$ Transform of the function $f_{\Delta x}(x) = e_{\Delta x}(a,x)$ will be derived. Then, to show the importance of the integration by parts equation in derivations, the $K_{\Delta x}$ Transform of the derivative of a function will be shown (i.e. $K_{\Delta x}[D_{\Delta x} f_{\Delta x}(x)]$).

Example 1.3-1 Find the derivative of $e_{\Delta x}(a,x)$

$$D_{\Delta x} e_{\Delta x}(a,x) = D_{\Delta x}(1+a\Delta x)^{\frac{x}{\Delta x}} = \frac{(1+a\Delta x)^{\frac{x+\Delta x}{\Delta x}} - (1+a\Delta x)^{\frac{x}{\Delta x}}}{\Delta x} \qquad (1.3\text{-}1)$$

$$D_{\Delta x} e_{\Delta x}(a,x) = \frac{(1+a\Delta x-1)(1+a\Delta x)^{\frac{x}{\Delta x}}}{\Delta x} \qquad (1.3\text{-}2)$$

$$D_{\Delta x} e_{\Delta x}(a,x) = a(1+a\Delta x)^{\frac{x}{\Delta x}} \qquad (1.3\text{-}3)$$

$$\mathbf{D_{\Delta x} e_{\Delta x}(a,x) = a e_{\Delta x}(a,x)} \qquad (1.3\text{-}4)$$

Note – Eq 1.3-4 is Equation #49 in Table 6 in the Appendix

Example 1.3-2  Find the integral of  $e_{\Delta x}(a,x)$

From Eq 1.3-3

$$_{\Delta x}\int D_{\Delta x}e_{\Delta x}(a,x)\Delta x \;=\; a\; _{\Delta x}\int e_{\Delta x}(a,x)\Delta x \;+\; c \tag{1.3-5}$$

$$\mathbf{_{\Delta x}\int e_{\Delta x}(a,x)\Delta x \;=\; \frac{1}{a}\, e_{\Delta x}(a,x) + k} \tag{1.3-6}$$

Note – Eq 1.3-6 is Equation #94 in Table 6 in the Appendix

or

$$\mathbf{_{\Delta x}\int (1+a\Delta x)^{\frac{x}{\Delta x}}\Delta x \;=\; \frac{1}{a}\,(1+a\Delta x)^{\frac{x}{\Delta x}} + k} \tag{1.3-7}$$

Example 1.3-3  Find the $K_{\Delta x}$ Transform of $e_{\Delta x}(a,x)$

$$e_{\Delta x}(a,x) = (1+a\Delta x)^{\frac{x}{\Delta x}} \tag{1.3-8}$$

$$K_{\Delta x}[f(x)] = \;_{\Delta x}\int_{0}^{\infty} (1+s\Delta x)^{-(\frac{x+\Delta x}{\Delta x})} f(x)\Delta x \tag{1.3-9}$$

$$K_{\Delta x}[e_{\Delta x}(a,x)] = \frac{1}{1+s\Delta x}\;_{\Delta x}\int_{0}^{\infty}\left(\frac{1+s\Delta x}{1+a\Delta x}\right)^{-(\frac{x}{\Delta x})}\Delta x \tag{1.3-10}$$

$$K_{\Delta x}[e_{\Delta x}(a,x)] = \frac{1}{1+s\Delta x}\;_{\Delta x}\int_{0}^{\infty}\left(1+\frac{(s-a)\Delta x}{1+a\Delta x}\right)^{-(\frac{x}{\Delta x})}\Delta x \tag{1.3-11}$$

Substituting b for the quantity, $\frac{s-a}{1+a\Delta x}$ in the following equation. This equation was obtained from Equation #99 in Table 6 in the Appendix.

$$_{\Delta x}\int (1+b\Delta x)^{-\frac{x}{\Delta x}}\Delta x = -\;\frac{1+b\Delta x}{b}(1+b\Delta x)^{-\frac{x}{\Delta x}} + k \tag{1.3-12}$$

$$K_{\Delta x}[e_{\Delta x}(a,x)] = \frac{1}{1+s\Delta x}\ {}_{\Delta x}\int_0^{\infty}\left(1+\frac{(s-a)\Delta x}{1+a\Delta x}\right)^{-\frac{x}{\Delta x}}\Delta x = -\frac{1}{1+s\Delta x}\,\frac{1+s\Delta x}{s-a}\left(1+\frac{(s-a)\Delta x}{1+a\Delta x}\right)^{-\frac{x}{\Delta x}}\Big|_0^{\infty} \qquad (1.3\text{-}13)$$

$$\mathbf{K_{\Delta x}[e_{\Delta x}(a,x)] = \frac{1}{s-a}} \qquad (1.3\text{-}14)$$

Note – Eq 1.3-14 is Equation #5 in Table 3 in the Appendix

<u>Example 1.3-4</u> Find the $K_{\Delta x}$ Transform of the derivative of a function ($K_{\Delta x}[D_{\Delta x}f(x)]$)

$$K_{\Delta x}[f(x)] = {}_{\Delta x}\int_0^{\infty}(1+s\Delta x)^{-(\frac{x+\Delta x}{\Delta x})}f(x)\Delta x \qquad (1.3\text{-}15)$$

Substituting

$$K_{\Delta x}[D_{\Delta x}f(x)] = {}_{\Delta x}\int_0^{\infty}(1+s\Delta x)^{-(\frac{x+\Delta x}{\Delta x})}D_{\Delta x}f(x)\Delta x \qquad (1.3\text{-}16)$$

Using the integration by parts equation

$${}_{\Delta x}\int_{x_1}^{x_2}v(x)D_{\Delta x}u(x)\Delta x = u(x)v(x)\Big|_{x1}^{x2} - {}_{\Delta x}\int_{x_1}^{x_2}D_{\Delta x}v(x)u(x+\Delta x)\Delta x \qquad (1.3\text{-}17)$$

Let $D_{\Delta x}u(x) = D_{\Delta x}f(x)$ $\qquad v(x) = (1+s\Delta x)^{-\frac{x+\Delta x}{\Delta x}}$

$u(x) = f(x)$ $\qquad D_{\Delta x}v(x) = \dfrac{(1+s\Delta x)^{-\frac{x+2\Delta x}{\Delta x}} - (1+s\Delta x)^{-\frac{x+\Delta x}{\Delta x}}}{\Delta x}$

$u(x+\Delta x) = f(x+\Delta x)$ $\qquad D_{\Delta x}v(x) = (1+s\Delta x)^{-\frac{x+\Delta x}{\Delta x}}\left(\frac{1}{1+s\Delta x}-1\right)\frac{1}{\Delta x}$

$$D_{\Delta x}v(x) = -s(1+s\Delta x)^{-\frac{x+2\Delta x}{\Delta x}}$$

From Eq 1.3-16 and Eq 1.3-17 and introducing the above substitutions

$$K_{\Delta x}[D_{\Delta x}f(x)] = f(x)(1+s\Delta x)^{-\frac{x+\Delta x}{\Delta x}}\Big|_0^\infty - {}_{\Delta x}\int_0^\infty -s(1+s\Delta x)^{-\frac{x+2\Delta x}{\Delta x}} f(x+\Delta x)\Delta x \qquad (1.3\text{-}18)$$

$$f(x+\Delta x) = f(x) + D_{\Delta x}f(x)\Delta x \qquad (1.3\text{-}19)$$

Substituting Eq 1.3-19 into Eq 1.3-18 and simplifying

$$K_{\Delta x}[D_{\Delta x}f(x)] = -\frac{f(0)}{1+s\Delta x} + \frac{s}{1+s\Delta x}\,{}_{\Delta x}\int_0^\infty (1+s\Delta x)^{-\frac{x+\Delta x}{\Delta x}} f(x)\Delta x$$

$$+\frac{s\Delta x}{1+s\Delta x}\,{}_{\Delta x}\int_0^\infty (1+s\Delta x)^{-\frac{x+\Delta x}{\Delta x}} D_{\Delta x} f(x)\Delta x \qquad (1.3\text{-}20)$$

Changing the notation

$$K_{\Delta x}[D_{\Delta x}f(x)] = -\frac{f(0)}{1+s\Delta x} + \frac{s}{1+s\Delta x}K_{\Delta x}[f(x)] + \frac{s\Delta x}{1+s\Delta x}K_{\Delta x}[D_{\Delta x}f(x)] \qquad (1.3\text{-}21)$$

$$\left(1 - \frac{s\Delta x}{1+s\Delta x}\right)K_{\Delta x}[D_{\Delta x}f(x)] = \frac{s}{1+s\Delta x}K_{\Delta x}[f(x)] - \frac{f(0)}{1+s\Delta x} \qquad (1.3\text{-}22)$$

Then

$$\mathbf{K_{\Delta x}[D_{\Delta x}f(x)] = sK_{\Delta x}[f(x)] - f(0)} \qquad (1.3\text{-}23)$$

Notes

1. Eq 1.3-23 above is Equation #7 in Table 2 in the Appendix.
2. For $\Delta x \to 0$, Eq 1.3-23 is recognized as the Calculus Laplace Transform for the derivative of a function.
3. The equation, $K_{\Delta x}[f(x+\Delta x)] = (1+s\Delta x)K_{\Delta x}[f(x)] - \Delta x f(0)$, used to directly solve difference equations, is derived in the same manner as Eq 1.3-23. This equation is Equation #11 in Table 2.

Later in this paper, using examples, difference equations and differential difference equations are solved by the $K_{\Delta x}$ Transform. The tables provided in the Appendix should provide sufficient support to derive and manipulate $K_{\Delta x}$ Transforms.

**Section 1.4: The lnd(n,Δx,x) Function and the Zeta Function**

Those who may wish to use Interval Calculus to find the area under sample and hold shaped waveforms will notice an equation, $\int_{\Delta x} \frac{1}{x^n}\,\Delta x = \pm\, \text{lnd}(n,\Delta x,x) + k$ , + for n = 1, – for n ≠ 1 (1.4-1)

This equation is Equation #89 in Table 6. This is a very important Interval Calculus equation. For Re(n)>1 and Δx→0, the function, $\int_{\Delta x\, x}^{\infty} \frac{1}{x^n}\,\Delta x = \text{lnd}(n,\Delta x,x)$ , has an equality in Calculus, $\int_{x}^{\infty} \frac{1}{x^n}\,dx = \frac{1}{n-1}\frac{1}{x^{n-1}}$, and for n=1 and Δx→0, the function, $\int_{\Delta x\, 1}^{x} \frac{1}{x}\,\Delta x = \text{lnd}(1,\Delta x,x) - \text{lnd}(1,\Delta x,1)$, has an equality in Calculus, $\int_{1}^{x} \frac{1}{x}\,dx = \ln x$. (See Table 9 Equation #49 in the Appendix). The function, lnd(n,Δx,x) can be used to evaluate a very important function, the General Zeta Function.

The lnd(n,Δx,x) Function

The lnd(n,Δx,x) function is defined as follows:

$$D_{\Delta x}\text{lnd}(n,\Delta x,x) \equiv D_{\Delta x}\ln_{\Delta x}x = +\frac{1}{x} , \quad n = 1 \qquad (1.4\text{-}2)$$

and

$$D_{\Delta x}\text{lnd}(n,\Delta x,x) = -\frac{1}{x^n} , \qquad n \neq 1 \qquad (1.4\text{-}3)$$

In Interval Calculus, the lnd(n,Δx,x) function holds the same position of importance as the natural logarithm in Calculus. Note the similarity.

The lnd(n,Δx,x) function is calculated using the program, LNDX. The LNDX program code is presented in the Calculation Programs section at the end of the Appendix.

The General Zeta Function

$$\zeta(n,\Delta x,x) = \frac{1}{\Delta x}\,\text{lnd}(n,\Delta x,x) , \quad \text{for all } n \qquad (1.4\text{-}4)$$

$$\zeta(n,\Delta x,x_i) = \frac{1}{\Delta x}\,\text{lnd}(n,\Delta x,x_i) = \sum_{\Delta x\; x=x_i}^{\pm\infty} \frac{1}{x^n} = \frac{1}{\Delta x}\int_{\Delta x\; x_i}^{\pm\infty} \frac{1}{x^n}\,\Delta x , \quad \text{Re}(n)>1 \qquad (1.4\text{-}5)$$

+∞ for Re(Δx)>0 or {Re(Δx)=0 and Im(Δx)>0}
–∞ for Re(Δx)<0 or {Re(Δx)=0 and Im(Δx)<0}

and

$$\zeta(n,\Delta x,x)\Big|_{x_1}^{x_2} = \frac{1}{\Delta x}\,\mathrm{lnd}(n,\Delta x,x)\Big|_{x_1}^{x_2} = \pm \sum_{\Delta x \; x=x_1}^{x_2-\Delta x} \frac{1}{x^n} = \pm \frac{1}{\Delta x} \int_{\Delta x \; x_1}^{x_2} \frac{1}{x^n}\,\Delta x\ , \quad + \text{ for } n = 1,\ - \text{ for } n \neq 1 \qquad (1.4\text{-}6)$$

where
- $\Delta x$ = x interval
- x = real or complex variable
- $x_i, x_1, x_2, \Delta x, n$ = real or complex constants
- Any summation term where x = 0 is excluded

Note – The relationship, $\lim_{x\to\infty} \mathrm{lnd}(n,\Delta x,x) \to 0$ where Re(n) > 1, is used in the derivation of Eq 1.4-5 from Eq 1.4-6

Special Cases of the General Zeta Function

The Hurwitz Zeta Function

For $\Delta x = 1$, $x_1 = x_i$, and $x_2 \to \infty$, the Hurwitz Zeta Function can be evaluated.

$$\zeta(n,x_i) = \mathrm{lnd}(n,1,x_i)\ , \quad n\neq 1 \qquad (1.4\text{-}7)$$

$$\zeta(n,x_i) = \mathrm{lnd}(n,1,x_i) = \sum_{1 \; x=x_i}^{\infty} \frac{1}{x^n} = \int_{1 \; x_i}^{\infty} \frac{1}{x^n}\,\Delta x\ , \quad \mathrm{Re}(n)>1 \qquad (1.4\text{-}8)$$

where
- x = real or complex variable
- $x_i$,n = real or complex constants
- Any summation term where x = 0 is excluded

The Riemann Zeta Function

For $\Delta x = 1$, $x_1 = 1$, and $x_2 \to \infty$, the Riemann Zeta Function can be evaluated.

$$\zeta(n) = \mathrm{lnd}(n,1,1)\ , \quad n\neq 1 \qquad (1.4\text{-}9)$$

$$\zeta(n) = \mathrm{lnd}(n,1,1) = \sum_{1 \; x=1}^{\infty} \frac{1}{x^n} = \int_{1 \; 1}^{\infty} \frac{1}{x^n}\,\Delta x\ , \quad \mathrm{Re}(n)>1 \qquad (1.4\text{-}10)$$

where
- x = real or complex variable
- n = real or complex constant

The N=1 Zeta Function

For $x_1 = \Delta x$, $x_2 = x_f$, n=1, the N=1 Zeta Function can be evaluated.

$$\zeta(1,\Delta x,x) = \frac{1}{\Delta x}\,\text{lnd}(1,\Delta x,x) \equiv \frac{1}{\Delta x}\,\ln_{\Delta x}x\,, \quad n=1 \tag{1.4-11}$$

$$\zeta(1,\Delta x,x_f) = \frac{1}{\Delta x}\,\text{lnd}(1,\Delta x,x_f) \equiv \frac{1}{\Delta x}\,\ln_{\Delta x}x_f = {}_{\Delta x}\!\!\sum_{x=\Delta x}^{x_f-\Delta x}\frac{1}{x} = \frac{1}{\Delta x}\,{}_{\Delta x}\!\!\int_{\Delta x}^{x_f}\frac{1}{x}\,\Delta x\,, \quad n=1 \tag{1.4-12}$$

and

$$\zeta(1,\Delta x,x)\Big|_{x_1}^{x_2} = \frac{1}{\Delta x}\,\text{lnd}(1,\Delta x,x)\Big|_{x_1}^{x_2} \equiv \frac{1}{\Delta x}\,\ln_{\Delta x}x\,\Big|_{x_1}^{x_2} = {}_{\Delta x}\!\!\sum_{x=x_1}^{x_2-\Delta x}\frac{1}{x} = \frac{1}{\Delta x}\,{}_{\Delta x}\!\!\int_{x_1}^{x_2}\frac{1}{x}\,\Delta x\,, \quad n=1 \tag{1.4-13}$$

where

$\Delta x$ = x interval
x = real or complex variable
$x_1, x_2, x_f, \Delta x$ = real or complex constants
Any summation term where x = 0 is excluded

$\text{lnd}(1,\Delta x,\Delta x) = 0$
$\text{lnd}(1,\Delta x,0) = 0$

$\ln_{\Delta x}x$ is an optional form of the $\text{lnd}(1,\Delta x,x)$ function used to emphasize its similarity to the natural logarithm. The above Zeta Functions are derived in Chapter 8.

The $\text{lnd}(n,\Delta x,x)$ Function can be calculated using two series which are functions of $n,\Delta x,x$ where $n,\Delta x,x$ may be real or complex values. The basic form of these two series is:

For n=1

$$\text{lnd}(n,\Delta x,x) \approx \gamma + \ln\left(\frac{x}{\Delta x} - \frac{1}{2}\right) + \sum_{m=1}^{\infty}\frac{(2m-1)!\,C_m}{(2m+1)!\,2^{2m}\left(\frac{x}{\Delta x} - \frac{1}{2}\right)^{2m}} \tag{1.4-12}$$

For n≠ 1

$$\text{lnd}(n,\Delta x,x) \approx -\sum_{m=0}^{\infty}\frac{\Gamma(n+2m-1)\left(\frac{\Delta x}{2}\right)^{2m} C_m}{\Gamma(n)(2m+1)!\left(x - \frac{\Delta x}{2}\right)^{n+2m-1}} + K \tag{1.4-14}$$

where

$$\Delta x \sum_{x=x_1}^{x_2} {}_{\Delta x} \frac{1}{x^n} = \pm \left. \mathrm{lnd}(n,\Delta x,x) \right|_{x_1}^{x_2+\Delta x} , \quad + \text{ for } n = 1, \ - \text{ for } n \neq 1 \tag{1.4-15}$$

any $\frac{1}{0}$ term is excluded

The accuracy of both series increases rapidly for increasing $\left|\frac{x}{\Delta x}\right|$

$x = x_i + r\Delta x$ ,  r = integers
$\Delta x$ = x increment
K = constant
$n,\Delta x,x,x_i$ = real or complex values
$\gamma$ = Euler's Constant = .577215664…
$C_0 = -1$
$C_1 = +1$
$C_2 = -\frac{7}{3}$
$C_3 = +\frac{31}{3}$
$C_4 = -\frac{381}{5}$
$C_5 = +\frac{2555}{3}$
$C_6 = -\frac{1414477}{105}$
$C_7 = +286685$
…

These constants can be calculated from Bernoulli constants. A computer program, LNDX, is available which calculates lnd(n,Δx,x). See Chapter 2 for the derivation and description of the lnd(n,Δx,x) function and its calculation program. The lnd(n,Δx,x) function is used to evaluate important mathematical functions, series, and discrete integrals. The use of the lnd(n,Δx,x) function is demonstrated in the Chapter 2 Solved Problems section at the end of Chapter 2.

## Section 1.5: The relationship between the $K_{\Delta x}$ Transform and the Z Transform

In Chapter 5, It is shown that the $K_{\Delta x}$ Transform is closely related to the Z Transform. It would be useful to have a convenient method for converting between these two transforms. Such a method does exist. Some of the calculation and conversion equations derived in Chapter 5 are presented on the following page.

Some useful $K_{\Delta x}$ Transform and Z Transform equations derived in Chapter 5 are:

1) $K_{\Delta x}$ Transform to Z Transform Conversion

$$Z[f(x)] = F(z) = \frac{z}{T} F(s)\Big|_{s = \frac{z-1}{T}} \qquad Z[f(x)] = F(z) \qquad \text{Z Transform} \quad (1.5\text{-}1)$$

$T = \Delta x$ sampling period

$x = nT$ $\qquad K_{\Delta x}[f(x)] = F(s) \qquad K_{\Delta x}$ Transform

$n = 0,1,2,3,\ldots$

For conversion from the Inverse $K_{\Delta x}$ Transform to the Inverse Z Transform
The complex plane integration contour changes from

$$s = \frac{e^{(\gamma+jw)\Delta x} - 1}{\Delta x} \text{ to } z = e^{(\gamma+jw)T} \qquad (1.5\text{-}2)$$

where

$\gamma$ = positive real constant

$-\frac{\pi}{\Delta x} \le w < \frac{\pi}{\Delta x}$

2) Z Transform to $K_{\Delta x}$ Transform Conversion

$$K_{\Delta x}[f(x)] = F(s) = \frac{\Delta x}{1+s\Delta x} F(z)\Big|_{z = 1+s\Delta x} \qquad K_{\Delta x}[f(x)] = F(s) \quad K_{\Delta x} \text{ Transform } (1.5\text{-}3)$$

$\Delta x = T$ sampling period

$x = n\Delta x$ $\qquad Z[f(x) = F(z) \qquad$ Z Transform

$n = 0, 1, 2, 3, \ldots$

For conversion from the Inverse Z Transform to the Inverse $K_{\Delta x}$ Transform
The complex plane integration contour changes from

$$z = e^{(\gamma+jw)T} \text{ to } s = \frac{e^{(\gamma+jw)\Delta x} - 1}{\Delta x} \qquad (1.5\text{-}4)$$

where

$\gamma$ = positive real constant

$-\frac{\pi}{\Delta x} \le w < \frac{\pi}{\Delta x}$

3) Z Transform Calculation Equation
(Interval Calculus form)

$$Z[f(x)] = F(z) = \frac{1}{T} {}_{T}\!\!\int_{0}^{\infty} z^{-\frac{x}{\Delta x}} f(x)\Delta x \qquad (5.4\text{-}5)$$

where

$x = n\Delta x$, $\quad n = 0,1,2,3,\ldots$

$T = \Delta x$

$\Delta x$ = sampling period, $x$ interval

Comments on Notation

1) In Z Transform literature, T is most often used to represent the sampling period. Either x or t is chosen as the variable.

2) The $K_{\Delta x}$ Transform designates its sampling period to be $\Delta x$. If the variable t is chosen instead of x, the $K_{\Delta t}$ Transform designates its sampling period to be $\Delta t$.

3) The f(x) functions for the Z Transform are Calculus functions such as $e^{ax}$, sinax, etc. The f(x) functions for the $K_{\Delta x}$ Transform are Interval Calculus discrete functions such as $e_{\Delta x}(a,x)$, $\sin_{\Delta x}(a,x)$ etc. The Calculus and Interval Calculus discrete functions are related.

For demonstation purposes, below are two conversion examples.

Example 1.5-1 Find the $K_{\Delta x}$ Transform for the function, $e^{ax}$. Calculate the $K_{\Delta x}$ Transform from the Z Transform of $e^{ax}$. Check the validity of the resulting $K_{\Delta x}$ Transform two ways.

$$f(x) = e^{ax} \tag{1.5-6}$$

$$Z[f(x)] = F(z) = \frac{z}{z - e^{aT}} \tag{1.5-7}$$

where

T = the sampling period

Rewriting Conversion Equation, Eq 1.5-3

$$K_{\Delta x}[f(x)] = f(s) = \frac{\Delta x}{1+s\Delta x} F(z)\Big|_{\substack{z = 1+s\Delta x \\ T = \Delta x}} \tag{1.5-8}$$

Substituting Eq 1.5-7 into Eq 1.5-8

$$f(s) = \frac{\Delta x}{1+s\Delta x}\left[\frac{z}{z - e^{aT}}\right]\Big|_{\substack{z = 1+s\Delta x \\ T = \Delta x}} \tag{1.5-9}$$

$$f(s) = \frac{\Delta x}{1+s\Delta x}\left[\frac{1+s\Delta x}{1+s\Delta x - e^{a\Delta x}}\right] = \frac{1}{s - \left(\frac{e^{a\Delta x}-1}{\Delta x}\right)}$$

$$f(s) = \frac{1}{s - \left(\frac{e^{a\Delta x}-1}{\Delta x}\right)} \tag{1.5-10}$$

Checking the validity of Eq 1.5-10

Check #1

From the $K_{\Delta x}$ Transform table, Table 3, in the Appendix

$$K_{\Delta x}[e_{\Delta x}(b,x)] = \frac{1}{s-b} \tag{1.5-11}$$

$$e_{\Delta x}(b,x) = (1+b\Delta x)^{\frac{x}{\Delta x}} \tag{1.5-12}$$

Let

$$b = \left(\frac{e^{a\Delta x}-1}{\Delta x}\right) \tag{1.5-13}$$

Substituting Eq 1.5-12 and Eq 1.5-13 into Eq 1.5-11

$$K_{\Delta x}\left[e_{\Delta x}\left(\frac{e^{a\Delta x}-1}{\Delta x},x\right)\right] = \frac{1}{s-\left(\frac{e^{a\Delta x}-1}{\Delta x}\right)} \tag{1.5-14}$$

From Eq 1.5-14 and Eq 1.5-12

$$K_{\Delta x}\left[e_{\Delta x}\left(\frac{e^{a\Delta x}-1}{\Delta x},x\right)\right] = K_{\Delta x}\left[\left(1+\left(\frac{e^{a\Delta x}-1}{\Delta x}\right)\Delta x\right)^{\frac{x}{\Delta x}}\right] = K_{\Delta x}[e^{ax}] = \frac{1}{s-\left(\frac{e^{a\Delta x}-1}{\Delta x}\right)} \tag{1.5-15}$$

$$K_{\Delta x}[e^{ax}] = \frac{1}{s-\left(\frac{e^{a\Delta x}-1}{\Delta x}\right)} \tag{1.5-16}$$

Then $\frac{1}{s-\left(\frac{e^{a\Delta x}-1}{\Delta x}\right)}$ is the $K_{\Delta x}$ Transform of $e^{ax}$ with a sample and hold period of $\Delta x$.

Checking the validity of Eq 1.5-8

Check #2

Rewriting Eq 1.2-2

$$K_{\Delta x}[f(x)] = {}_1\sum_{n=0}^{\infty}(1+s\Delta x)^{-n-1}f(n\Delta x)\Delta x \tag{1.5-17}$$

Substituting Eq 1.5-10 into Eq 1.5-17 to find f(x)

$$\frac{1}{s - (\frac{e^{a\Delta x} - 1}{\Delta x})} = {}_1\sum_{n=0}^{\infty} (1+s\Delta x)^{-n-1} f(n\Delta x)\Delta x \qquad (1.5\text{-}18)$$

Simplifying

$$\frac{\Delta x}{(1+s\Delta x) - e^{a\Delta x}} = \Delta x\, {}_1\sum_{n=0}^{\infty} (1+s\Delta x)^{-n-1} f(n\Delta x)$$

$$\frac{1}{(1+s\Delta x) - e^{a\Delta x}} = {}_1\sum_{n=0}^{\infty} (1+s\Delta x)^{-n-1} f(n\Delta x) \qquad (1.5\text{-}19)$$

Expanding the left side of Eq 1.5-19

$${}_1\sum_{n=0}^{\infty} f(n\Delta x)(1+s\Delta x)^{-n-1} = e^{a0\Delta x}(1+s\Delta x)^{-1} + e^{a1\Delta x}(1+s\Delta x)^{-2} + e^{a2\Delta x}(1+s\Delta x)^{-3} + e^{a3\Delta x}(1+s\Delta x)^{-4} + \ldots$$

(1.5-20)

Compare both sides of Eq 1.5-20

$f(0) = e^{a0\Delta x} = 1$
$f(\Delta x) = e^{a\Delta x}$
$f(2\Delta x) = e^{a2\Delta x}$
$f(3\Delta x) = e^{a3\Delta x}$
…

In General
$f(n\Delta x) = e^{an\Delta x}$ , $n = 0,1,2,3,\ldots$ (1.5-21)

Then $\frac{1}{s - (\frac{e^{a\Delta x} - 1}{\Delta x})}$ is the $K_{\Delta x}$ Transform of $e^{ax}$ with a sampling period of $\Delta x$.

.

Example 1.5-2 Derive the Inverse Z Transform from the Inverse $K_{\Delta x}$ Transform.

The Inverse $K_{\Delta x}$ Transform is:

$$f(x) = K_{\Delta x}^{-1}[f(s)] = \frac{1}{2\pi j} \oint_{c_1} [1+s\Delta x]^{\frac{x}{\Delta x}} f(s)\, ds \tag{1.5-22}$$

where

$$s = \frac{e^{(\gamma+jw)\Delta x} - 1}{\Delta x}\,, \quad -\frac{\pi}{\Delta x} \le w < \frac{\pi}{\Delta x}\,, \quad \text{This is the } c_1 \text{ contour} \tag{1.5-23}$$

$x = m\Delta x$ , $m$ = integers
$\Delta x$ = x increment
$f(x)$ = discrete function of x
$f(x) = 0$ for $x < 0$
$f(s)$ = function of s
$\gamma$, $w$ = real value constants
$\gamma > 0$

Rewriting Eq 1.5-1 and Eq 1.5-2

$K_{\Delta x}$ Transform to Z Transform Conversion

$$Z[f(x)] = F(z) = \frac{z}{T} f(s)\Big|_{s = \frac{z-1}{T}} \tag{1.5-24}$$

$T = \Delta x$ sampling period
$x = n\Delta x$
$n = 0,1,2,3,\ldots$

For inverse transform conversions, the complex plane integration contour changes from

$$s = \frac{e^{(\gamma+jw)\Delta x} - 1}{\Delta x} = \text{contour } c_1 \text{ to } z = e^{(\gamma+jw)T} = \text{contour } c_2 \tag{1.5-25}$$

where
$\gamma$ = positive real constant
$-\frac{\pi}{\Delta x} \le w < \frac{\pi}{\Delta x}$

Using Eq 1.5-24 and Eq 1.5-25 to derive the Inverse Z Transform from the Inverse $K_{\Delta x}$ Transform equations, Eq 1.5-22 and Eq 1.5-23.

$$f(x) = \frac{1}{2\pi j} \oint_{c_1} [1+s\Delta x]^{\frac{x}{\Delta x}} f(s)\, ds \tag{1.5-26}$$

$$x = n\Delta x \tag{1.5-27}$$

$n = 0,1,2,3,\ldots$ (1.5-28)

$T = \Delta x$ (1.5-29)

$s = \frac{z-1}{T}$ (1.5-30)

Substituting Eq 1.5-27 thru Eq 1.5-30 into Eq 1.5-26

$$f(nT) = \frac{1}{2\pi j} \oint_{c_2} [1+[(\frac{z-1}{T})T]^n \frac{T}{z} f(z) \frac{dz}{T} \quad (1.5\text{-}31)$$

Simplifying Eq 1.5-31

The Inverse Z Transform is:

$$f(nT) = \frac{1}{2\pi j} \oint_{c_2} f(z)\, z^{n-1} dz \quad (1.5\text{-}32)$$

where
$n = 0, 1, 2, 3, \ldots$
contour $c_2$ is $z = e^{(\eta+jw)T}$ in the z plane
$\gamma$ = positive real constant
$-\frac{\pi}{\Delta x} \le w < \frac{\pi}{\Delta x}$

Checking – Eq 1.5-32 is recognized as being the Inverse Z Transform.

## Section 1.6: The Inverse $K_{\Delta x}$ Transform

The $K_{\Delta x}$ Transform and the Inverse $K_{\Delta x}$ Transform are derived in Chapter 5. The relationships are as follows:

The $K_{\Delta x}$ Transform

$$f(s) = \int_{0}^{\infty} {}_{\Delta x} [1+s\Delta x]^{-\frac{x+\Delta x}{\Delta x}} f(x)\, \Delta x \quad (1.6\text{-}1)$$

The Inverse $K_{\Delta x}$ Transform

$$f(x) = K_{\Delta x}^{-1}[f(s)] = \frac{1}{2\pi j}\oint_c [1+s\Delta x]^{\frac{x}{\Delta x}} f(s)\, ds \qquad (1.6\text{-}2)$$

where

$s = \dfrac{e^{(\gamma+jw)\Delta x} - 1}{\Delta x}$ , $\quad -\dfrac{\pi}{\Delta x} < w \leq \dfrac{\pi}{\Delta x}$ , $\quad$ c contour

$x = m\Delta x$ , m = integers

$\Delta x$ = x increment

f(x) = discrete function of x

f(x) = 0 for x < 0

f(s) = function of s

$\gamma$, w = real value constants

$\gamma > 0$

$K_{\Delta x}^{-1}[f(s)] = f(x)$ = sum of the residues of $(1+s\Delta x)^{\frac{x}{\Delta x}} f(s)$ at all the poles of f(s) inside the complex plane contour, c, which is a circle in the left half of the complex plane with center at $-\dfrac{1}{\Delta x}$ and a radius of $\dfrac{1}{\Delta x}$ .

A demonstration of a $K_{\Delta x}$ Inverse Transform calculation using residue theory is shown in Section 3 of Chapter 4.

Of course, the inverse $K_{\Delta x}$ Transform can also be obtained from the $K_{\Delta x}$ Transform tables. In fact, this is usually done. In addition, there are two related methods of using a series to evaluate f(x) from its corresponding $K_{\Delta x}$ Transform, F(s). These methods are worthy of note and are presented below.

Series Method #1 for evaluating f(x) from $K_{\Delta x}[f(x)]$

Consider the following derived equation:

$$F(s) = K_{\Delta x}[f(x)] = \Delta x \sum_{n=0}^{\infty} f(n\Delta x)(1+s\Delta x)^{-n-1}) \equiv \Delta x \int_{\Delta x\, x=0}^{\infty} f(x)(1+s\Delta x)^{-(\frac{x+\Delta x}{\Delta x})} \quad , \; x = n\Delta x \; , \; n = 0,1,2,\ldots$$

(1.6-3)

Expanding the series of Eq 1.6-3

$$F(s) = \Delta x[f(0)(1+s\Delta x)^{-1} + f(1\Delta x)(1+s\Delta x)^{-2} + f(2\Delta x)(1+s\Delta x)^{-3} + f(3\Delta x)(1+s\Delta x)^{-4} + \ldots] \qquad (1.6\text{-}4)$$

where

f(x) = a function of x

$F(s) = K_{\Delta x}[f(x)]$ , the $K_{\Delta x}$ Transform of the function f(x)

$\Delta x$ = the x increment

Putting Eq 1.5-4 in a different form

$$\frac{F(s)}{\Delta x} = \frac{f(0)}{(1+s\Delta x)^1} + \frac{f(1\Delta x)}{(1+s\Delta x)^2} + \frac{f(2\Delta x)}{(1+s\Delta x)^3} + \frac{f(3\Delta x)}{(1+s\Delta x)^4} + \dots \qquad (1.6\text{-}5)$$

where

$f(x)$ = a function of x

$F(s) = K_{\Delta x}[f(x)]$ , the $K_{\Delta x}$ Transform of the function f(x)

$\Delta x$ = the x increment

It is observed that the coefficients of the Eq 1.6-5 series represent the values of f(x) for $x = n\Delta x$ , $n = 0,1,2,3,\dots$ Thus, by expanding a $K_{\Delta x}$ Transform, F(s), into a series of the form of Eq 1.6-5, $f(x) = K_{\Delta x}^{-1}[F(s)]$ can be evaluated. The $f(n\Delta x)$ n=0,1,2,3,… coefficients of Eq 1.6-5 are evaluated by the following method of taking the value of s to infinity.

Multiply both sides of Eq 1.6-5 by $1+s\Delta x$

$$(1+s\Delta x)\frac{F(s)}{\Delta x} = f(0) + \frac{f(1\Delta x)}{(1+s\Delta x)^1} + \frac{f(2\Delta x)}{(1+s\Delta x)^2} + \frac{f(3\Delta x)}{(1+s\Delta x)^3} + \dots \qquad (1.6\text{-}6)$$

Letting $s\rightarrow\infty$

$$(1+s\Delta x)\frac{F(s)}{\Delta x}\Big|_{s\rightarrow\infty} = f(0) + \frac{f(1\Delta x)}{(1+s\Delta x)^1} + \frac{f(2\Delta x)}{(1+s\Delta x)^2} + \frac{f(3\Delta x)}{(1+s\Delta x)^3} + \dots\Big|_{s\rightarrow\infty} = f(0) \qquad (1.6\text{-}7)$$

$$(1+s\Delta x)\frac{F(s)}{\Delta x}\Big|_{s\rightarrow\infty} = f(0) \qquad (1.6\text{-}8)$$

From Eq 1.6-6

$$(1+s\Delta x)\frac{F(s)}{\Delta x}\Big|_{s\rightarrow\infty} - f(0) = \frac{f(1\Delta x)}{(1+s\Delta x)^1} + \frac{f(2\Delta x)}{(1+s\Delta x)^2} + \frac{f(3\Delta x)}{(1+s\Delta x)^3} + \dots \qquad (1.6\text{-}9)$$

Applying the same process to Eq 1.6-9

Multiply both sides of Eq 1.6-9 by $1+s\Delta x$

$$(1+s\Delta x)\left[(1+s\Delta x)\frac{F(s)}{\Delta x}\Big|_{s\rightarrow\infty} - f(0)\right] = f(1\Delta x) + \frac{f(2\Delta x)}{(1+s\Delta x)^1} + \frac{f(3\Delta x)}{(1+s\Delta x)^2} + \dots \qquad (1.6\text{-}10)$$

Letting $s\rightarrow\infty$

$$(1+s\Delta x)\left[(1+s\Delta x)\frac{F(s)}{\Delta x}\Big|_{s\rightarrow\infty} - f(0)\right]\Big|_{s\rightarrow\infty} = f(1\Delta x) + \frac{f(2\Delta x)}{(1+s\Delta x)^1} + \frac{f(3\Delta x)}{(1+s\Delta x)^2} + \dots\Big|_{s\rightarrow\infty} = f(1\Delta x) \qquad (1.6\text{-}11)$$

$$(1+s\Delta x)\left[(1+s\Delta x)\frac{F(s)}{\Delta x}\Big|_{s\to\infty} - f(0)\right]\Big|_{s\to\infty} = f(1\Delta x) \tag{1.6-12}$$

From 1.6-10

$$(1+s\Delta x)\left[(1+s\Delta x)\frac{F(s)}{\Delta x}\Big|_{s\to\infty} - f(0)\right] - f(1\Delta x) = \frac{f(2\Delta x)}{(1+s\Delta x)^1} + \frac{f(3\Delta x)}{(1+s\Delta x)^2} + \ldots \tag{1.6-13}$$

Applying the same process to Eq 1.6-13

Multiply both sides of Eq 1.6-13 by $1+s\Delta x$

$$(1+s\Delta x)\left[(1+s\Delta x)\left[(1+s\Delta x)\frac{F(s)}{\Delta x}\Big|_{s\to\infty} - f(0)\right] - f(1\Delta x)\right] = f(2\Delta x) + \frac{f(3\Delta x)}{(1+s\Delta x)^1} + \ldots \tag{1.6-14}$$

Letting $s\to\infty$

$$(1+s\Delta x)\left[(1+s\Delta x)\left[(1+s\Delta x)\frac{F(s)}{\Delta x}\Big|_{s\to\infty} - f(0)\right] - f(1\Delta x)\right]\Big|_{s\to\infty} = f(2\Delta x) + \frac{f(3\Delta x)}{(1+s\Delta x)^1} + \ldots\Big|_{s\to\infty} = f(2\Delta x) \tag{1.6-15}$$

$$(1+s\Delta x)\left[(1+s\Delta x)\left[(1+s\Delta x)\frac{F(s)}{\Delta x}\Big|_{s\to\infty} - f(0)\right] - f(1\Delta x)\right]\Big|_{s\to\infty} = f(2\Delta x) \tag{1.6-16}$$

This process is repeated to find $f(3\Delta x), f(4\Delta x), f(5\Delta x), \ldots$

Thus,

Rewriting some of the previously derived equations

$$f(0) = (1+s\Delta x)\frac{F(s)}{\Delta x}\Big|_{s\to\infty} \tag{1.6-17}$$

$$f(1\Delta x) = (1+s\Delta x)\left[(1+s\Delta x)\frac{F(s)}{\Delta x}\Big|_{s\to\infty} - f(0)\right]\Big|_{s\to\infty} \tag{1.6-18}$$

$$f(2\Delta x) = (1+s\Delta x)\left[(1+s\Delta x)\left[(1+s\Delta x)\frac{F(s)}{\Delta x}\Big|_{s\to\infty} - f(0)\right] - f(1\Delta x)\right]\Big|_{s\to\infty} \tag{1.6-19}$$

…

Observing Eq 1.6-17 thru 1.6-19

In general

$$f(n\Delta x) = G_n(s)\Big|_{s\to\infty} \tag{1.6-20}$$

where

$$G_0(s) = (1+s\Delta x)\frac{F(s)}{\Delta x} \quad (1.6\text{-}21)$$

$$G_{n+1}(s) = (1+s\Delta x)[G_n(s) - G_n(s)|_{s\to\infty}\ ] \quad (1.6\text{-}22)$$

**<u>Series Method #1 for Evaluating f(x) from its corresponding $K_{\Delta x}$ Transform, $K_{\Delta x}[f(x)]$</u>**

**Evaluating $f(x) = K_{\Delta x}^{-1}[F(s)]$ given the $K_{\Delta x}$ Transform, $F(s) = K_{\Delta x}[f(x)]$**

$$\mathbf{f(n\Delta x) = G_n(s)|_{s\to\infty}\ ,\ \ x = n\Delta x,\ \ \ n = 0,1,2,3,\ldots} \quad (1.6\text{-}23)$$

**where**

$$\mathbf{G_0(s) = (1+s\Delta x)\frac{F(s)}{\Delta x}} \quad (1.6\text{-}24)$$

$$\mathbf{G_{n+1}(s) = (1+s\Delta x)[G_n(s) - G_n(s)|_{s\to\infty}\ ]} \quad (1.6\text{-}25)$$

Demonstate Series Method #1 for calculating the inverse $K_{\Delta x}$ Transform.

<u>Example 1.6-1</u> Evaluate $f(x) = K_{\Delta x}^{-1}[\frac{1}{s(s+.5)}]$ for x = 0, .2, .4, .6

Use Eq 1.6-23 thru Eq 1.6-25 to Evaluate f(x)

$$G_0(s) = \frac{1+s\Delta x}{\Delta x(s)(s+.5)} \quad (1.6\text{-}26)$$

$$G_0(s)\ |_{s\to\infty} = \frac{1+s\Delta x}{\Delta x(s)(s+.5)}\ |_{s\to\infty} = 0 \quad (1.6\text{-}27)$$

$$f(0) = G_0(s)\ |_{s\to\infty} = 0 \quad (1.6\text{-}28)$$

$$G_1(s) = (1+s\Delta x)[G_0(s) - G_0(s)\ |_{s\to\infty}\ ] \quad (1.6\text{-}29)$$

$$G_1(s) = (1+s\Delta x)\left[\frac{1+s\Delta x}{\Delta x(s)(s+.5)} - 0\right] = \frac{(1+s\Delta x)^2}{\Delta x(s)(s+.5)} \quad (1.6\text{-}30)$$

$\Delta x = .2$

$$G_1(s)\ |_{s\to\infty} = \left[\frac{(1+.2s)^2}{.2(s)(s+.5)}\right]|_{s\to\infty} = .2 \quad (1.6\text{-}31)$$

$$f(.2) = G_1(s)|_{s\to\infty} = .2 \quad (1.6\text{-}32)$$

$$G_2(s) = (1+s\Delta x)[G_1(s) - G_1(s)\ |_{s\to\infty}\ ] \quad (1.6\text{-}33)$$

$$G_2(s) = (1+s\Delta x)\left[\frac{(1+s\Delta x)^2}{\Delta x(s)(s+.5)} - .2\right] \quad (1.6\text{-}34)$$

$\Delta x = .2$

$$G_2(s) = (1+.2s)\left[\frac{(1+.2s)^2}{.2(s)(s+.5)} - .2\right] = (1+.2s)\left[\frac{1+.4s+.04s^2-.04s^2-.02s}{.2s(s+.5)}\right] = \left[\frac{1+.58s+.076s^2}{.2s(s+.5)}\right] \quad (1.6\text{-}35)$$

$$G_2(s)\Big|_{s\to\infty} = \left[\frac{1+.58s+.076s^2}{.2s(s+.5)}\right]\Big|_{s\to\infty} = .38 \quad (1.6\text{-}36)$$

$$f(.4) = G_2(s)\Big|_{s\to\infty} = .38 \quad (1.6\text{-}37)$$

$$G_3(s) = (1+s\Delta x)[G_2(s) - G_2(s)\,|_{s\to\infty}\,] \quad (1.6\text{-}38)$$

$$G_3(s) = (1+s\Delta x)\left[\frac{1+.58s+.076s^2}{.2s(s+.5)} - .38\right] \quad (1.6\text{-}39)$$

$\Delta x = .2$

$$G_3(s) = (1+.2s)\left[\frac{1+.58s+.076s^2}{.2(s)(s+.5)} - .38\right] = (1+.2s)\left[\frac{1+.58s+.076s^2-.076s^2-.038s}{.2s(s+.5)}\right] = \left[\frac{1+.742s+.1084s^2}{.2s(s+.5)}\right] \quad (1.6\text{-}40)$$

$$G_3(s)\Big|_{s\to\infty} = \left[\frac{1+.742s+.1084s^2}{.2s(s+.5)}\right]\Big|_{s\to\infty} = .542 \quad (1.5\text{-}41)$$

$$f(.6) = G_3(s)\Big|_{s\to\infty} = .542 \quad (1.6\text{-}42)$$

Then

$f(0) = 0$

$f(.2) = .2$

$f(.4) = .38$

$f(.6) = .542$

Checking

$$F(s) = [\frac{1}{s(s+.5)}] = \frac{2}{s} - \frac{2}{s+.5} \quad (1.6\text{-}43)$$

Taking the inverse $K_{\Delta x}$ Transform of F(s)

$$f(x) = K_{\Delta x}^{-1}[F(s)] = 2 - 2[1-.5(.2)]^{\frac{x}{.2}} = 2[1-(.9)^{5x}]$$

$$f(x) = 2[1-(.9)^{5x}] \quad (1.6\text{-}44)$$

$f(0) = 2[1-(.9)^0] = 0$

$f(.2) = 2[1-(.9)^1] = .2$

$f(.4) = 2[1-(.9)^2] = .38$

$f(.6) = 2[1-(.9)^3] = .542$ Good check

Series Method #2 for evaluating f(x) from $K_{\Delta x}[f(x)]$

Consider the following derived equation:

$$F(s) = K_{\Delta x}[f(x)] = \Delta x \sum_{n=0}^{\infty} f(n\Delta x)(1+s\Delta x)^{-n-1}) \equiv \int_{0}^{\infty}{}_{\Delta x} f(x)(1+s\Delta x)^{-(\frac{x+\Delta x}{\Delta x})} \Delta x \quad , \; x = n\Delta x \; , \; n = 0,1,2,\ldots \tag{1.6-45}$$

Expanding the series of Eq 1.6-45

$$F(s) = \Delta x[f(0)(1+s\Delta x)^{-1} + f(1\Delta x)(1+s\Delta x)^{-2} + f(2\Delta x)(1+s\Delta x)^{-3} + f(3\Delta x)(1+s\Delta x)^{-4} + \ldots] \tag{1.6-46}$$

where

f(x) = a function of x

$F(s) = K_{\Delta x}[f(x)]$ , the $K_{\Delta x}$ Transform of the function f(x)

$\Delta x$ = the x increment

It is observed that the coefficients of the series represent the values of f(x) for $x = n\Delta x$ , n = 0,1,2,3,… Thus, by expanding a $K_{\Delta x}$ Transform, F(s), into a series of the form of Eq 1.5-46, $f(x) = K_{\Delta x}^{-1}[F(s)]$ can be evaluated.

Simplifying the application of the method described above

$$s = \frac{(1+s\Delta x) - 1}{\Delta x} = \frac{p - 1}{\Delta x} \tag{1.6-47}$$

where

$$p = 1+s\Delta x \tag{1.6-48}$$

Substitute Eq 1.6-47 into Eq 1.6-46 and multiply both sides of Eq 1.6-46 by $\frac{p}{\Delta x}$ (1.6-49)

$$\frac{p}{\Delta x} F(s)\Big|_{s=\frac{p-1}{\Delta x}} = [f(0)p^{0} + f(1\Delta x)p^{-1} + f(2\Delta x)p^{-2} + f(3\Delta x)(p^{-3} + \ldots] \tag{1.6-50}$$

From Eq 1.6-50

Evaluating $f(x) = K_{\Delta x}^{-1}[F(s)]$ given the $K_{\Delta x}$ Transform, $F(s) = K_{\Delta x}[f(x)]$

$$\frac{p}{\Delta x} F(s)\Big|_{s=\frac{p-1}{\Delta x}} = \sum_{n=0}^{\infty} f(n\Delta x)p^{-n} \tag{1.6-51}$$

where

f(x) = a function of x

$F(s) = K_{\Delta x}[f(x)]$ , the $K_{\Delta x}$ Transform of the function f(x)

$\Delta x$ = the x increment

$x = n\Delta x$ , n=0,1,2,3,…

$p = 1+s\Delta x$

Referring to Eq 1.5-1 in the previous section, Section 1.5, an interesting relationship is found. Eq 1.5-1 is rewritten below.

$K_{\Delta x}$ Transform to Z Transform Conversion

$$Z[f(x)] = F(z) = \frac{z}{T} F(s)\Big|_{\substack{s = \frac{z-1}{T} \\ T = \Delta x}} \quad , \qquad Z[f(x)] = F(z) \quad \text{Z Transform} \qquad (1.6\text{-}52)$$

$$K_{\Delta x}[f(x)] = F(s) \quad K_{\Delta x} \text{ Transform}$$

From Eq 1.6-51 and Eq 1.6-52
Renaming p to z

$$\frac{z}{\Delta x} F(s)\Big|_{s = \frac{z-1}{\Delta x}} = Z[f(x)] = F(z) = \sum_{n=0}^{\infty} f(n\Delta x) z^{-n} \qquad (1.6\text{-}53)$$

where

- $f(x)$ = a function of x
- $F(s) = K_{\Delta x}[f(x)]$ , the $K_{\Delta x}$ Transform of the function, f(x)
- $F(z) = Z[f(x)]$ , the Z Transform of the function, f(x)
- $\Delta x$ = the x increment
- $x = n\Delta x$ , $n=0,1,2,3,\ldots$
- $z = 1+s\Delta x$

Note - The Z Transform uses the variable, T, to represent $\Delta x$

Thus

**Series Method #2 for Evaluating f(x) from its corresponding $K_{\Delta x}$ Transform**

**Evaluating $f(x) = K_{\Delta x}^{-1}[F(s)]$ given the $K_{\Delta x}$ Transform, $F(s) = K_{\Delta x}[f(x)]$**

$$\frac{z}{\Delta x} F(s)\Big|_{s = \frac{z-1}{\Delta x}} = Z[f(x)] = F(z) = \sum_{n=0}^{\infty} f(n\Delta x) z^{-n} \qquad (1.6\text{-}54)$$

**where**

- **$f(x)$ = a function of x**
- **$F(s) = K_{\Delta x}[f(x)]$ , the $K_{\Delta x}$ Transform of the function, f(x)**
- **$F(z) = Z[f(x)]$ , the Z Transform of the function, f(x)**
- **$\Delta x$ = the x increment**
- **$x = n\Delta x$ , $n = 0,1,2,3,\ldots$**
- **$z = 1+s\Delta x$**

**Note - The Z Transform uses the variable, T, to represent $\Delta x$**

Demonstate Series Method #2 for calculating the inverse $K_{\Delta x}$ Transform. Use the same problem solved by Series Method #1 in Example 1.6-1.

Example 1.6-2 Evaluate $f(x) = K_{\Delta x}^{-1}[\frac{1}{s(s+.5)}]$ for x = 0, .2, .4, .6

Use Eq 1.6-54 to Evaluate f(x)

Rewriting Eq 1.6-54

$$\frac{z}{\Delta x} F(s)\Big|_{s = \frac{z-1}{\Delta x}} = Z[f(x)] = F(z) = \sum_{n=0}^{\infty} f(n\Delta x) z^{-n}$$

$$F(s) = [\frac{1}{s(s+.5)}] \tag{1.6-55}$$

Substitute 1.6-55 into Eq 1.6-54

$$\frac{z}{\Delta x}[\frac{1}{s(s+.5)}]\Big|_{s = \frac{z-1}{\Delta x}} = F(z) = \sum_{n=0}^{\infty} f(n\Delta x) z^{-n} \tag{1.6-56}$$

Find F(z)

$$F(z) = \frac{z}{\Delta x}[\frac{1}{s(s+.5)}]\Big|_{s = \frac{z-1}{\Delta x}} = \frac{z}{\Delta x}[\frac{1}{(\frac{z-1}{\Delta x})(\frac{z-1}{\Delta x} + .5)}] = \frac{z\Delta x}{(z-1)(z-1+.5\Delta x)} \tag{1.6-57}$$

$\Delta x = .2$

$$F(z) = \frac{.2z}{(z-1)(z-.9)} = \frac{.2z}{z^2 - 1.9z + .9} \tag{1.6-58}$$

Expanding F(z) into a series in the form of Eq 1.6-56

$$F(z) = \frac{.2z}{z^2 - 1.9z + .9} = 0 + .2z^{-1} + .38z^{-2} + .542z^{-3} + \ldots \tag{1.6-59}$$

Substituting Eq 1.6-59 into 1.6-56

$$\sum_{n=0}^{\infty} f(n\Delta x) z^{-n} = 0 + .2z^{-1} + .38z^{-2} + .542z^{-3} + \ldots \tag{1.6-60}$$

Then comparing both sides of Eq 1.6-60

f(0) = 0
f(.2) = .2
f(.4) = .38
f(.6) = .542

Checking

$$F(s) = [\frac{1}{s(s+.5)}] = \frac{2}{s} - \frac{2}{s+.5} \qquad (1.6\text{-}61)$$

Taking the inverse $K_{\Delta x}$ Transform of F(s)

$$f(x) = K_{\Delta x}^{-1}[F(s)] = 2 - 2[1 - .5(.2)]^{\frac{x}{.2}} = 2[1-(.9)^{5x}]$$

$$f(x) = 2[1-(.9)^{5x}] \qquad (1.6\text{-}62)$$

Then

$f(0) = 2[1-(.9)^{0}] = 0$
$f(.2) = 2[1-(.9)^{1}] = .2$
$f(.4) = 2[1-(.9)^{2}] = .38$
$f(.6) = 2[1-(.9)^{3}] = .542$
Good check

## Section 1.7: Interval Calculus Area Calculation

Interval Calculus area calculation, of course, has similarities to area calculation in Calculus. If no poles are encountered within the range of integration, the following Inteval Calculus area calculation equation applies:

Area Calculation Equation

$$A = {}_{\Delta x}\int_{x_1}^{x_2} f(x)\Delta x \qquad (1.7\text{-}1)$$

where

$x = x_1, x_1+\Delta x, x_1+2\Delta x, x_1+3\Delta x, \ldots, x_2-\Delta x, x_2$

This equation is observed to be similar to the Calculus area calculation equation:

$$A = \int_{x_1}^{x_2} f(x)\Delta x \qquad (1.7\text{-}2)$$

where

$x_1 \le x \le x_2$

However, there are also dissimilarities. The Interval Calculus area calculation equation may also be represented as follows:

$$A = {}_{\Delta x}\int_{x_1}^{x_2} f(x)\Delta x = \Delta x \; {}_{\Delta x}\sum_{x=x_1}^{x_2-\Delta x} f(x) \qquad (1.7\text{-}3)$$

where

$x = x_1, x_1+\Delta x, x_1+2\Delta x, x_1+3\Delta x, \ldots, x_2-\Delta x, x_2$

Here, the Interval Calculus integral is shown to be equal to a summation. In Calculus, integrals are not normally shown to be equivalent to summations with an infinitesimal $\Delta x$.

A table of Interval Calculus discrete integrals is provided in Table 6 of the Appendix.

An example of a discrete integration area calculation is shown below in Example 1.7-1.

<u>Example 1.7-1</u>

Find the area, A, under the discrete sample and hold curve, $f(x) = \frac{1}{x(x-.5)}$, where $A = \int_{.5\,1}^{5} \frac{1}{x(x-.5)}\Delta x$ .

Using the integration table, Table 6, in the Appendix

$$\int_{\Delta x} \frac{1}{\prod_{n=a}^{b}(x-n\Delta x)}\Delta x = -\frac{1}{(b-a)\prod_{n=a+1}^{b}(x-n\Delta x)} + k \qquad (1.7\text{-}4)$$

where

$n = a, a+1, a+2,\ldots, b-1, b$

$b > a$

$x = x_0 + m\,\Delta x$, $m$ = integers

$b-a$ = one less than the order of the denominator polynomial of the function being integrated

$\Delta x$ = x interval

$k$ = constant of integration

For a = 0, b = 1

$$\int_{\Delta x} \frac{1}{x(x-\Delta x)}\Delta x = -\frac{1}{x-\Delta x} + k \qquad (1.7\text{-}5)$$

For $\Delta x = .5$

$x_1 = 1$ and $x_2 = 5$ , the integration limits

$$\int_{.5\,1}^{5} \frac{1}{x(x-.5)}\Delta x = -\frac{1}{x-.5}\Big|_{1}^{5} \qquad (1.7\text{-}6)$$

$$A = \int_{.5\,1}^{5} \frac{1}{x(x-.5)}\Delta x \qquad (1.7\text{-}7)$$

Solving the integral in Eq 1.7-7 using Eq 1.7-6

$$A = \int_{1}^{5} {}_{.5}\frac{1}{x(x-.5)}\Delta x = -\frac{1}{x-.5}\Big|_{1}^{5} = -\frac{1}{4.5} + \frac{1}{.5} = 1.7777777777 \qquad (1.7\text{-}8)$$

$$A = 1.7777777777 \qquad (1.7\text{-}9)$$

Checking the result of Eq 1.7-9

From Eq 1.7-3

$$A = \int_{1}^{5} {}_{.5}\frac{1}{x(x-.5)}\Delta x = .5 \sum_{x=1}^{4.5} {}_{.5}\frac{1}{x(x-.5)} \qquad (1.7\text{-}10)$$

where

x = 1, 1.5, 2, 2.5, 3, 3.5, 4, 4.5, 5

$$A = .5 \sum_{x=1}^{4.5} {}_{.5}\frac{1}{x(x-.5)} = .5\left[\frac{1}{1(.5)} + \frac{1}{1.5(1)} + \frac{1}{2(1.5)} + \frac{1}{2.5(2)} + \right. \qquad (1.7\text{-}11)$$

$$\left.\frac{1}{3(2.5)} + \frac{1}{3.5(3)} + \frac{1}{4(3.5)} + \frac{1}{4.5(4)}\right] = 1.7777777777$$

$$A = 1.7777777777 \qquad (1.7\text{-}12)$$

Good check

Dealing with poles that are encountered during integration or summation

In Calculus, it is possible that a value of x within the range of integration or summation may occur at a pole. Should this occur, most mathematical problem solving stops. The integral or sum has an infinite value. When it is acceptable to find the sum of just the finite terms of a summation or the area under the finite portions of a curve, multiple integrations or summations are performed on the finite values between the pole(s) and the lower and upper x limits. In Interval Calculus, the previous statement also applies but, interestingly, there are some equations which deal with the presence of poles. These equations are derived and presented below. If applicable, these equations exclude the poles in an integral or summation so that integration or summation can be performed on the remaining finite values.

Mean value mathematics will be used in the derivation of several of the equations. Here, the value of a function, f(x), at x is not calculated directly. It is calculated indirectly using the values, $(x+\varepsilon)$ and $(x-\varepsilon)$, on either side of x . For x not at a pole, it is rather evident that $f(x) = \lim_{\varepsilon\to 0}\frac{f(x+\varepsilon) + f(x-\varepsilon)}{2}$ . However, for x at a pole, $f(x) \neq \lim_{\varepsilon\to 0}\frac{f(x+\varepsilon) + f(x-\varepsilon)}{2}$ . The quantity $\lim_{\varepsilon\to 0}\frac{f(x+\varepsilon) + f(x-\varepsilon)}{2}$ is, nevertheless, useful and is defined as the mean value at a point, $MV[f(x)] = \lim_{\varepsilon\to 0}\frac{f(x+\varepsilon) + f(x-\varepsilon)}{2}$ . For each value of f(x), a unique mean value, MV[f(x)] can be calculated.
The following equation is valid for all x not occurring at a pole:

$$A = {}_{\Delta x}\int_{x_1}^{x_2} f(x)\Delta x = \Delta x \; {}_{\Delta x}\sum_{x=x_1}^{x_2-\Delta x} f(x) \tag{1.7-13}$$

$$x = x_1,\ x_1+\Delta x,\ x_1+2\Delta x,\ x_1+3\Delta x,\ \ldots,\ x_2-\Delta x,\ x_2$$

Consider the possibility of one or more values of x occurring at a pole.

Apply the mean value concept to Eq 1.7-13.

$$\frac{1}{2}\lim_{\varepsilon\to 0}\left[{}_{\Delta x}\int_{x_1+\varepsilon}^{x_2+\varepsilon} f(x)\Delta x + {}_{\Delta x}\int_{x_1-\varepsilon}^{x_2-\varepsilon} f(x)\Delta x\right] = \frac{1}{2}\lim_{\varepsilon\to 0}\left[\Delta x \; {}_{\Delta x}\sum_{x=x_1+\varepsilon}^{x_2+\varepsilon-\Delta x} f(x) + \Delta x \; {}_{\Delta x}\sum_{x=x_1-\varepsilon}^{x_2-\varepsilon-\Delta x} f(x)\right] \tag{1.7-14}$$

$$\frac{1}{2}\lim_{\varepsilon\to 0}\left[{}_{\Delta x}\int_{x_1+\varepsilon}^{x_2+\varepsilon} f(x)\Delta x + {}_{\Delta x}\int_{x_1-\varepsilon}^{x_2-\varepsilon} f(x)\Delta x\right] = \frac{1}{2}\lim_{\varepsilon\to 0}\left[\Delta x \; {}_{\Delta x}\sum_{x=x_1}^{x_2-\Delta x} f(x+\varepsilon) + \Delta x \; {}_{\Delta x}\sum_{x=x_1}^{x_2-\Delta x} f(x-\varepsilon)\right] \tag{1.7-15}$$

$$\frac{1}{2}\lim_{\varepsilon\to 0}\left[{}_{\Delta x}\int_{x_1+\varepsilon}^{x_2+\varepsilon} f(x)\Delta x + {}_{\Delta x}\int_{x_1-\varepsilon}^{x_2-\varepsilon} f(x)\Delta x\right] = \lim_{\varepsilon\to 0}\left[\Delta x \; {}_{\Delta x}\sum_{x=x_1}^{x_2-\Delta x} \frac{f(x+\varepsilon) + f(x-\varepsilon)}{2}\right] \tag{1.7-16}$$

Consider some x to occur at the poles where $x = \mathit{x}_1, \mathit{x}_2, \mathit{x}_3, \ldots, \mathit{x}_m$ , m = the number of poles within the range $x_1 < x < x_2$.

Separating the pole function values at $x$ from the non-pole function values

$$\frac{1}{2}\lim_{\varepsilon\to 0}\left[\ {}_{\Delta x}\!\int_{x_1+\varepsilon}^{x_2+\varepsilon} f(x)\Delta x + {}_{\Delta x}\!\int_{x_1-\varepsilon}^{x_2-\varepsilon} f(x)\Delta x\right] = \lim_{\varepsilon\to 0}\left[\Delta x\ {}_{\Delta x}\!\sum_{x=x_1}^{x_2-\Delta x} \frac{f(x+\varepsilon)+f(x-\varepsilon)}{2}\quad x \neq \text{any } x\right]$$

$$+ \lim_{\varepsilon\to 0}\left[\Delta x\ {}_{\Delta x}\!\sum_{r=1}^{m} \frac{f(x+\varepsilon)+f(x-\varepsilon)}{2}\quad x = \text{all } x\right] \qquad (1.7\text{-}17)$$

Let

$$MV[f(x)] = \lim_{\varepsilon\to 0}\left[\frac{f(x+\varepsilon)+f(x-\varepsilon)}{2}\right], \text{ the mean value of the function, f(x), at the point x} \qquad (1.7\text{-}18)$$

From Eq 1.7-17 and Eq 1.7-18

$$MV\left[{}_{\Delta x}\!\int_{x_1}^{x_2} f(x)\Delta x\right] = \underset{\text{any } \frac{1}{0} \text{ term excluded}}{\Delta x\ {}_{\Delta x}\!\sum_{x=x_1}^{x_2-\Delta x} f(x)} + \{\Delta x\textstyle\sum MV[f(x)]\ ,\ x = x_1, x_2, x_3, \ldots, x_m\} \qquad (1.7\text{-}19)$$

The sum, $\underset{\text{any } \frac{1}{0} \text{ term excluded}}{\Delta x\ {}_{\Delta x}\!\sum_{x=x_1}^{x_2-\Delta x} f(x)}$, in Eq 1.7-19 contains the f(x) pole values for x = $x_1, x_2, x_3, \ldots, x_m$ which have been excluded (i.e. set equal to 0).

Then, the sample and hold shaped curve area between the limits $x_1$ and $x_2$ where all f(x) values at a pole are excluded is:

$$A = \underset{\text{any } \frac{1}{0} \text{ term excluded}}{\Delta x\ {}_{\Delta x}\!\sum_{x=x_1}^{x_2-\Delta x} f(x)} \qquad (1.7\text{-}20)$$

Substituting Eq 1.7-20 into Eq 1.7-19

$$MV\left[{}_{\Delta x}\!\int_{x_1}^{x_2} f(x)\Delta x\right] = A + \{\Delta x\textstyle\sum MV[f(x)]\ ,\ x = x_1, x_2, x_3, \ldots, x_m\} \qquad (1.7\text{-}21)$$

Solving for A

$$A = MV[{}_{\Delta x}\int_{x_1}^{x_2} f(x)\Delta x] - \{\Delta x \sum MV[f(x)]\ ,\ x = x_1, x_2, x_3, \ldots, x_m\} \qquad (1.7\text{-}22)$$

On Eq 1.7-22, several conditions are imposed.

The values of $MV[{}_{\Delta x}\int_{x_1}^{x_2} f(x)\Delta x]$, $x_1, x_2,$ and $\{\sum MV[f(x)]\ ,\ x = x_1, x_2, x_3, \ldots, x_m\}$ must be finite.

For these conditions, A will be a finite value.

Then, from the derived equations above, Eq 1.7-23 has been obtained to calculate the area under a sample and hold curve where any existing poles are excluded.

$$A = \underset{\text{any } \frac{1}{0} \text{ term excluded}}{\Delta x_{\Delta x} \sum_{x=x_1}^{x_2-\Delta x} f(x)} = MV[{}_{\Delta x}\int_{x_1}^{x_2} f(x)\Delta x] - \{\Delta x \sum MV[f(x)]\ ,\ x = x_1, x_2, x_3, \ldots, x_m\} \qquad (1.7\text{-}23)$$

where

$$MV[f(x)] = \lim_{\varepsilon\to 0}\left[\frac{f(x+\varepsilon) + f(x-\varepsilon)}{2}\right] \qquad (1.7\text{-}24)$$

$x = x_1, x_1+\Delta x, x_1+2\Delta x, x_1+3\Delta x, \ldots, x_2-\Delta x, x_2$

$x_r$ = those values of x, $x_1<x<x_2$, for which f(x) poles occur

$r = 1,2,3,\ldots,m$

m = the number of poles between the limits $x_1$,$x_2$

$f(x_r) = \frac{k}{0}$ is excluded and given a 0 value, k = constant , $r = 1,2,3,\ldots,m$

For all x = real values

A = area under the f(x) curve between the limits $x_1$ and $x_2$ as

defined by the summation, $\Delta x_{\Delta x} \sum_{x=x_1}^{x_2-\Delta x} f(x)$, where if a summation term

value is $\frac{k}{0}$ at a pole at $x_r$, it is excluded

Necessary conditions for the calculation of A

The values of $MV[{}_{\Delta x}\int_{x_1}^{x_2} f(x)\Delta x]$, $x_1, x_2,$ and $\sum MV[f(x)]\ ,\ x = x_1, x_2, x_3, \ldots, x_m$, must be finite.

For these conditions, A will be a finite value.

$MV[f(x_r)]$ = the Mean Value of the function, f(x), calculated at $x=x_r$

$\Delta x$,x = real or complex values

If there are no poles between $x_1$ and $x_2$, the $\sum MV[f(x)]$ mean value term of the area equation, Eq 1.7-23, is left out.

Note 1 – If f(x) has one or more poles within the range $x_1$ to $x_2$, the integral, ${}_{\Delta x}\int_{x_1}^{x_2} f(x)\Delta x$ , will include these same non-finite values and be infinite. For some f(x) functions, the integral of the mean values of f(x), $MV[{}_{\Delta x}\int_{x_1}^{x_2} f(x)\Delta x]$, will be finite. In this case the integral of the mean values of f(x) will have the same values of f(x) for x not at a pole and will have a value, $MV[f(x_r)]$, for each value of $x = x_r$ at a pole. If the values of $MV[{}_{\Delta x}\int_{x_1}^{x_2} f(x)\Delta x]$ and all $MV(x_r)$ are finite, the area generated from all of the finite f(x) can be calculated. See Eq 1.7-23 above.

Note 2 – The necessary conditions that the values of $MV[{}_{\Delta x}\int_{x_1}^{x_2} f(x)\Delta x]$, $x_1$, $x_2$, and $\sum MV[f(x)]$ , $x = x_1, x_2, x_3, \ldots, x_m$, must be finite significantly limit the usefulness of Eq 1.7-23.

In addition to the above equations there are several other useful mean value equations which are derived and presented below.

Find the Mean Value of the function, $f(x) = \frac{1}{x-a}$ at $x = a$.

$$MV[f(x)] = \lim_{\varepsilon\to 0}\left[\frac{f(x+\varepsilon) + f(x-\varepsilon)}{2}\right]$$

$$MV[\frac{1}{x-a}] = \lim_{\varepsilon\to 0}\left[\frac{\frac{1}{x+\varepsilon-a} + \frac{1}{x-\varepsilon-a}}{2}\right]$$

$$MV[\frac{1}{x-a}]\,|_{x=a} = \lim_{\varepsilon\to 0}\frac{1}{2}\left[\frac{1}{\varepsilon} + \frac{1}{-\varepsilon}\right] = \lim_{\varepsilon\to 0}\frac{1}{2}\left[\frac{1}{\varepsilon} - \frac{1}{\varepsilon}\right] = 0$$

$$MV[\frac{1}{x-a}]\,|_{x=a} = 0 \qquad (1.7\text{-}25)$$

Find an expression for the Mean Value of the discrete integral, ${}_{\Delta x}\int_{x_1}^{x_2} \frac{1}{x-a}\Delta x$.

From Eq 1.7-23

$$MV\left[{}_{\Delta x}\int_{x_1}^{x_2} f(x)\Delta x\right] = \underset{\text{any }\frac{1}{0}\text{ term excluded}}{\Delta x_{\Delta x}\sum_{x=x_1}^{x_2-\Delta x} f(x)} + \{\Delta x \sum MV[f(x)]\ ,\ x = x_1, x_2, x_3, \ldots, x_m\} \tag{1.7-26}$$

Substituting $f(x) = \frac{1}{x-a}$ into Eq 1.7-26

$$MV\left[{}_{\Delta x}\int_{x_1}^{x_2} \frac{1}{x-a}\Delta x\right] = \underset{\text{any }\frac{1}{0}\text{ term excluded}}{\Delta x_{\Delta x}\sum_{x=x_1}^{x_2-\Delta x} \frac{1}{x-a}} + \Delta x\, MV\left[\frac{1}{x-a}\right]\Big|_{x=a} \tag{1.7-27}$$

Substituting Eq 1.7-25 into Eq 1.7-27

$$MV\left[{}_{\Delta x}\int_{x_1}^{x_2} \frac{1}{x-a}\Delta x\right] = \underset{\text{any }\frac{1}{0}\text{ term excluded}}{\Delta x_{\Delta x}\sum_{x=x_1}^{x_2-\Delta x} \frac{1}{x-a}} \tag{1.7-28}$$

In Chapter 2, a function, lnd(n,Δx,x), is derived and programmed. It is as follows:

$$\underset{\text{any }\frac{1}{0}\text{ term excluded}}{\Delta x_{\Delta x}\sum_{x=x_1}^{x_2-\Delta x} \frac{1}{(x-a)^n}} = \pm\, \mathrm{lnd}(n,\Delta x,x-a)\Big|_{x_1}^{x_2}\quad ,\quad + \text{ for } n = 1\quad ,\quad - \text{ for } n \neq 1 \tag{1.7-29}$$

where

$x = x_1, x_1+\Delta x, x_2+2\Delta x, x_1+3\Delta x, \ldots, x_2-\Delta x, x_2$

$x,\Delta x$ = real or complex variable

$\Delta x$ = x increment

$n,a$ = real or complex constant

Comment – The lnd(n,Δx,x) function of Eq 1.7-29 was derived such that any division by zero term is set equal to zero. This was done to exclude any pole in a summation. Note the Hurwitz Zeta function summation.

Substituting Eq 1.7-29 with n = 1 into Eq 1.7-28

$$MV\left[{}_{\Delta x}\int_{x_1}^{x_2} \frac{1}{x-a}\Delta x\right] = \underset{\text{any } \frac{1}{0} \text{ term excluded}}{\Delta x_{\Delta x}\sum_{x=x_1}^{x_2-\Delta x} \frac{1}{x-a}} = \text{lnd}(1,\Delta x,x-a)\Big|_{x_1}^{x_2} \qquad (1.7\text{-}30)$$

Note – Where $a$ is a value of x, the pole at x = $a$ has been excluded.

Thus, the equations which may be helpful in dealing with poles encountered in an integration or summation are as follows:

$$\mathbf{A} = \underset{\text{any } \frac{1}{0} \text{ term excluded}}{\Delta x_{\Delta x}\sum_{x=x_1}^{x_2-\Delta x} f(x)} = MV\left[{}_{\Delta x}\int_{x_1}^{x_2} f(x)\Delta x\right] - \{\Delta x\textstyle\sum MV[f(x)]\ ,\ x = x_1, x_2, x_3, \ldots, x_m\} \qquad (1.7\text{-}31)$$

$$MV[f(x)] = \lim_{\varepsilon\to 0}\left[\frac{f(x+\varepsilon) + f(x-\varepsilon)}{2}\right] \qquad (1.7\text{-}32)$$

**Conditions for the calculation of A**

**The values of** $MV\left[{}_{\Delta x}\int_{x_1}^{x_2} f(x)\Delta x\right]$, $x_1$, $x_2$, **and** $\sum MV[f(x)]\ ,\ x = x_1, x_2, x_3, \ldots, x_m$, **must be finite.**

**If there are no poles between $x_1$ and $x_2$, the $\sum MV[f(x)]$ mean value term of the area equation, Eq 1.7-31, is left out.**

**For** $f(x) = \frac{1}{(x-a)^n}$

$$\mathbf{A} = \underset{\text{any } \frac{1}{0} \text{ term excluded}}{\Delta x_{\Delta x}\sum_{x=x_1}^{x_2-\Delta x} \frac{1}{(x-a)^n}} = \pm\,\text{lnd}(n,\Delta x,x-a)\Big|_{x_1}^{x_2}\ ,\quad +\text{ for } n = 1\ ,\quad -\text{ for } n \neq 1 \qquad (1.7\text{-}33)$$

**For** $f(x) = \dfrac{1}{x-a}$

$$MV\left[{}_{\Delta x}\int_{x_1}^{x_2} \frac{1}{x-a}\,\Delta x\right] = \Delta x_{\Delta x}\sum_{x=x_1}^{x_2-\Delta x}\frac{1}{x-a} = \mathrm{lnd}(1,\Delta x,x-a)\Big|_{x_1}^{x_2} \qquad (1.7\text{-}34)$$

**any** $\frac{1}{0}$ **term excluded**

$$MV\left[\frac{1}{x-a}\right]\Big|_{x=a} = 0 \qquad (1.7\text{-}35)$$

**where**

$x = x_1, x_1+\Delta x, x_1+2\Delta x, x_1+3\Delta x, \ldots, x_2-\Delta x, x_2$

$x_r$ **= those values of x,** $x_1<x<x_2$, **for which f(x) poles occur**

$r = 1,2,3,\ldots,m$

**m = the number of poles between the limits** $x_1,x_2$

$f(x_r) = \frac{k}{0}$ **is excluded (i.e. given a 0 value), k = constant ,** $r = 1,2,3,\ldots,m$

**For all x = real values**

**A = area under the f(x) curve between the limits** $x_1$ **and** $x_2$ **as defined by the summation,** $\Delta x_{\Delta x}\sum_{x=x_1}^{x_2-\Delta x} f(x)$, **where if a summation term value is** $\frac{k}{0}$ **at a pole at** $x_r$, **it is excluded**

$MV[f(x_r)]$ **= the Mean Value of the function, f(x), calculated at** $x=x_r$

$x,\Delta x$ **= real or complex values**

**n,a = real or complex constants**

The previously derived equations, if applicable, provide a means to integrate or sum just the finite values within an integral or summation. These equations are used alone or in combination.

Below are some intregrals and summations to which the above specified equations are applicable.

1)
$$\int_{\Delta x} \frac{1}{\prod_{n=a}^{b}(x-n\Delta x)}\,\Delta x = -\frac{1}{(b-a)\prod_{n=a+1}^{b}(x-n\Delta x)} + k \qquad (1.7\text{-}36)$$

$n = a, a+1, a+2,\ldots, b-1, b$

$b > a$

$x = x_0 + m\,\Delta x$, $m$ = integers

b-a = one less than the order of the denominator polynomial of the function being integrated

2) $$\sum_{x=x_1}^{x_2}{}_{\Delta x}\ \frac{x^m + a_{m-1}x^{m-1} + a_{m-2}x^{m-2} + \dots + a_1x + a_0}{x^n + b_{n-1}x^{n-1} + b_{n-2}x^{n-2} + \dots + b_1x + b_0} \qquad (1.7\text{-}37)$$

where

$m,n$ = positive integers

$a_m,b_n$ = contants

$x = x_1, x_1{+}\Delta x, x_2{+}2\Delta x, x_1{+}3\Delta x, \dots , x_2{-}\Delta x, x_2$

$\Delta x$ = x increment

Note that the polynomial function of Eq 1.7-37 can be represented by a partial fraction expansion.

3) $${}_{\Delta x}\!\int_{x_1}^{x_2} \frac{1}{(x-a)^n}\,\Delta x \qquad (1.7\text{-}38)$$

where

$x = x_1, x_1{+}\Delta x, x_2{+}2\Delta x, x_1{+}3\Delta x, \dots , x_2{-}\Delta x, x_2$

$\Delta x$ = x increment

$n,a$ = real or complex constants

Comment – Discrete integrals and summations are related.

$${}_{\Delta x}\!\int_{x_1}^{x_2} f(x)\Delta x = \Delta x \sum_{x=x_1}^{x_2-\Delta x}{}_{\Delta x}\ f(x) \qquad (1.7\text{-}39)$$

where

$x = x_1, x_1{+}\Delta x, x_1{+}2\Delta x, x_1{+}3\Delta x, \dots , x_2{-}\Delta x, x_2$

The following example, Example 1.7-2, demonstates an application of the area calculation equations derived above.

Example 1.7-2

Find the integral, ${}_1\!\int_{-2}^{4} \frac{1}{x(x-1)}\Delta x$ , excluding the poles at x = 0,1 so that the area, A, attributable to the finite values between x = -2 and x = 4 can be obtained.

$$A = \int_{-2}^{4} {}_1 \frac{1}{x(x-1)}\Delta x = \sum_{x=-2}^{3} {}_1 \frac{1}{x(x-1)} \tag{1.7-40}$$

any $\frac{1}{0}$ term excluded any $\frac{1}{0}$ term excluded

Find the value of $MV[\int_{-2}^{4} {}_1 \frac{1}{x(x-1)}\Delta x]$

From the integration table, Table 6, in the Appendix

$$\int_{\Delta x} \frac{1}{x(x-\Delta x)}\Delta x = -\frac{1}{x-\Delta x} + k \tag{1.7-41}$$

$\Delta x = 1$

$$MV[\int_{-2}^{4} {}_1 \frac{1}{x(x-1)}\Delta x] = MV[-\frac{1}{x-1}\Big|_{-2}^{4}] = \frac{1}{2}\lim_{\varepsilon\to 0}[-\frac{1}{3+\varepsilon} + \frac{1}{-3+\varepsilon} - \frac{1}{3-\varepsilon} + \frac{1}{-3-\varepsilon}] = -\frac{2}{3} \tag{1.7-42}$$

Find the mean value of the function, $f(x) = \frac{1}{x(x-1)}$, at x = 0

$$MV[f(x)] = \lim_{\varepsilon\to 0}\left[\frac{f(x+\varepsilon) + f(x-\varepsilon)}{2}\right] \tag{1.7-43}$$

$$MV[f(0)] = \lim_{\varepsilon\to 0}\frac{1}{2}\left[\frac{1}{\varepsilon(\varepsilon-1)} + \frac{1}{-\varepsilon(-\varepsilon-1)}\right] = \lim_{\varepsilon\to 0}\frac{1}{2\varepsilon}\left[\frac{1}{(\varepsilon-1)} + \frac{1}{(\varepsilon+1)}\right] = \lim_{\varepsilon\to 0}\frac{1}{2\varepsilon}\left[\frac{2\varepsilon}{(\varepsilon^2-1)}\right] = -1 \tag{1.7-44}$$

Find the mean value of the function, $f(x) = \frac{1}{x(x-1)}$, at x = 1

$$MV[f(1)] = \lim_{\varepsilon\to 0}\frac{1}{2}\left[\frac{1}{(\varepsilon+1)\varepsilon} + \frac{1}{(1-\varepsilon)(-\varepsilon)}\right] = \lim_{\varepsilon\to 0}\frac{1}{2\varepsilon}\left[\frac{1}{(\varepsilon+1)} + \frac{1}{(\varepsilon-1)}\right] = \lim_{\varepsilon\to 0}\frac{1}{2\varepsilon}\left[\frac{2\varepsilon}{(\varepsilon^2-1)}\right] = -1 \tag{1.7-45}$$

Substitute the finite values obtained in Eq 1.7-42, Eq 1.7-44, and Eq 1,7-45 into the mean value area calculation equation, Eq 1.7-31, rewritten on the following page.

$$A = \underset{\text{any } \frac{1}{0} \text{ term excluded}}{\Delta x \sum_{\Delta x \; x=x_1}^{x_2-\Delta x} f(x)} = MV[{}_{\Delta x}\int_{x_1}^{x_2} f(x)\Delta x] - \{ \Delta x \sum MV[f(x)] \text{ , } x = x_1, x_2, x_3, \ldots, x_m \} \qquad (1.7\text{-}46)$$

$$A = \underset{\text{any } \frac{1}{0} \text{ term excluded}}{{}_1\sum_{x=-2}^{3} \frac{1}{x(x-1)}} = MV[{}_1\int_{-2}^{4} \frac{1}{x(x-1)}\Delta x] - \{\sum MV[f(x)] \text{ , } x = 0, 1]\} \text{ , } \Delta x = 1 \qquad (1.7\text{-}47)$$

$$A = \underset{\text{any } \frac{1}{0} \text{ term excluded}}{{}_1\int_{-2}^{4} \frac{1}{x(x-1)}\Delta x} = \underset{\text{any } \frac{1}{0} \text{ term excluded}}{{}_1\sum_{x=-2}^{3} \frac{1}{x(x-1)}} = -\frac{2}{3} - (-1-1) = 1\,\frac{1}{3} \qquad (1.7\text{-}48)$$

Checking the result of Eq 1.7-48

$$A = \underset{\text{any } \frac{1}{0} \text{ term excluded}}{{}_1\sum_{x=-2}^{3} \frac{1}{x(x-1)}} = \frac{1}{-2(-3)} + \frac{1}{-1(-2)} + 0 + 0 + \frac{1}{2(1)} + \frac{1}{3(2)} = \frac{1}{6} + \frac{1}{2} + \frac{1}{2} + \frac{1}{6} = 1\frac{1}{3} \qquad (1.7\text{-}48)$$

Good check

Two other applications of the area equations derived above, Example 1.7 and Example 1.8, are shown at the end of this chapter.

## Section 1.8: Application of the $K_{\Delta x}$ Transform to Control System Design

Control system design has depended heavily on the use of Operational Calculus to analyze and design stable control systems. Laplace Transform Methods for continuous time systems and Z Transform Methods for discrete time sampling systems have been used very effectively.
$K_{\Delta x}$ Transform system analysis methods can be added to the above mentioned methods.
The $K_{\Delta x}$ Transform, though it applies directly to discrete time systems, is remarkably similar to the Laplace Transform. $K_{\Delta x}$ Transforms become Laplace Transforms when the $\Delta x$ interval ($\Delta t$ interval for time variable systems) approaches zero. An advantage of $K_{\Delta x}$ Transform system design methods is that system stability can be clearly and easily observed for varying time intervals. System stability or instability is determined by observing any encirclements of the –1 point in the complex plane in a similar manner to that of Laplace Transform system analysis methods. In fact, the $K_{\Delta x}$ Transform system stability analysis methodology is just an extension of the Nyquist Criterion from continuous system analysis to discrete system analysis. The $K_{\Delta x}$ Transform is a major component of Operational Interval Calculus. Since Interval (Discrete) Calculus is a generalization of Calculus, it is strongly related to Calculus and can be learned and understood quite readily by those familiar with Calculus.

System stability analysis is primarily concerned with the location of the poles of a closed loop transfer function. Generally this transfer function is formed from an open loop transfer function which has its output fed back to its input for closed loop control. For stability, the poles of a transfer function must lie within a specified region of the s plane or z plane, whichever is appropriate. For Laplace Transform continuous time system analysis, all pole roots (of s) being within the left half of the s plane define a stable system. For Z Transform discrete (sampled) time system analysis, all pole roots (of z) being within the unit circle about the origin of the z plane define a stable system. For $K_{\Delta x}$ Transform discrete (sampled) time system analysis, all pole roots (of s) being within a $\frac{1}{\Delta x}$ radius circle centered at $-\frac{1}{\Delta x}$ in the left half of the s plane define a stable system.

Below, the two very important concepts of $K_{\Delta x}$ Transform discrete time system analysis are proved. The first concept relates to the fact that for closed loop system stability all closed loop transfer function poles must lie within a $\frac{1}{\Delta x}$ radius circle centered at $-\frac{1}{\Delta x}$ in the left half of the s plane. The second concept relates to the fact that a modified Nyquist Criterion identifies system stability and degree of stability.

**Proof of $K_{\Delta x}$ Transform discrete time system analysis Concept #1**
**Transfer function pole location for system stability**

1) <u>Consider a $K_{\Delta x}$ Transform function, F(s), with the denominator, $s^2+As+B$, having complex roots.</u>

$$F(s) = \frac{Cs+D}{s^2+As+B} \tag{1.8-1}$$

where

A,B,C,D = real constants

To find the inverse of the $K_{\Delta x}$ Transform function, F(s), change the form of Eq 1.8-1 to:

$$F(s) = \frac{K_1(s-a) + K_2w(1+a\Delta x)}{(s-a)^2 + (w[1+a\Delta x])^2} = \frac{Cs+D}{s^2+As+B} \tag{1.8-2}$$

$$m = 1+a\Delta x \tag{1.8-3}$$

Expanding Eq 1.8-2

$$F(s) = \frac{K_1s - K_1a + K_2wm}{s^2 - 2as + a^2 + w^2m^2} \tag{1.8-4}$$

Equating Eq 1.8-1 to Eq 1.8-4

$$s^2 + As + B = s^2 - 2as + a^2 + w^2m^2 \quad (1.8\text{-}5)$$

$$A = -2a$$

$$a = -\frac{A}{2} \quad (1.8\text{-}6)$$

$$B = a^2 + w^2m^2 \quad (1.8\text{-}7)$$

$$(wm)^2 = B - a^2 = B - (\frac{A}{2})^2$$

$$wm = \sqrt{B - (\frac{A}{2})^2} \quad (1.8\text{-}8)$$

$$Cs + D = K_1s + (-K_1a + K_2wm) \quad (1.8\text{-}9)$$

$$K_1 = C \quad (1.8\text{-}10)$$

$$D = -K_1a + K_2wm \quad (1.8\text{-}11)$$

Substituting Eq 8.1-6 and Eq 1.8-10 into Eq 1.8-11

$$D = \frac{CA}{2} + K_2wm$$

$$K_2wm = D - \frac{CA}{2} \quad (1.8\text{-}12)$$

Substituting Eq 1.8-8 into Eq 1.8-12

$$K_2 = \frac{D - \frac{CA}{2}}{\sqrt{B - (\frac{A}{2})^2}} \quad (1.8\text{-}13)$$

Find the roots of the denominator of F(s), $s^2 + As + B$

$$s = -\frac{A}{2} \pm \sqrt{(\frac{A}{2})^2 - B}$$

for complex roots, $B > (\frac{A}{2})^2$

$$s = -\frac{A}{2} \pm j\sqrt{B - (\frac{A}{2})^2} = a \pm bj \qquad (1.8\text{-}14)$$

As in Eq 18.6

$$a = -\frac{A}{2} \qquad (1.8\text{-}15)$$

$$b = \sqrt{B - (\frac{A}{2})^2} \qquad (1.8\text{-}16)$$

From Eq 1.8-15, Eq 1.8-16, Eq 1.8-3 and Eq 1.8-7

$$b = \sqrt{B - (\frac{A}{2})^2} = \sqrt{a^2 + w^2m^2 - (\frac{A}{2})^2} = \sqrt{(\frac{A}{2})^2 + w^2m^2 - (\frac{A}{2})^2} = wm = w(1+a\Delta x)$$

$$b = w(1+a\Delta x) \qquad (1.8\text{-}17)$$

$$w = \frac{b}{(1+a\Delta x)} \qquad (1.8\text{-}18)$$

Then

$$\mathbf{F(s)} = \frac{K_1(s-a) + K_2w(1+a\Delta x)}{(s-a)^2 + (w[1+a\Delta x])^2} = \frac{K_1(s-a) + K_2b}{(s-a)^2 + b^2} = \frac{Cs+D}{s^2+As+B} \qquad (1.8\text{-}19)$$

**where**

**$s = a \pm bj$, the roots of $s^2 + As + B$**

**$a,b,A,B,C,D,K_1,K_2$ = real constants**

$$a = -\frac{A}{2}$$

$$b = \sqrt{B - (\frac{A}{2})^2} = w(1+a\Delta x)$$

$$w = \frac{b}{(1+a\Delta x)}$$

$$K_1 = C$$

$$K_2 = \frac{D - \frac{CA}{2}}{\sqrt{B - (\frac{A}{2})^2}} = \frac{D + Ca}{b}$$

Find the inverse $K_{\Delta x}$ Transform of Eq 1.8-19

$$F(x) = e_{\Delta x}(a,x)[K_1\cos_{\Delta x}(\frac{b}{1+a\Delta x},x) + K_2\sin_{\Delta x}(\frac{b}{1+a\Delta x},x)] \;, \quad 1+a\Delta x \neq 0 \tag{1.8-20}$$

$$F(x) = [b\Delta x]^{\frac{x}{\Delta x}}[K_1\cos\frac{\pi x}{2\Delta x} + K_2\sin\frac{\pi x}{2\Delta x}] \;, \quad 1+a\Delta x = 0$$

where
$x = 0, \Delta x, 2\Delta x, 3\Delta x, 4\Delta x, \dots$

Eq 1.8-20 was obtained from the equations on lines 8 and 9 of the $K_{\Delta x}$ Transform Table, Table 3.

The functions, $e_{\Delta x}(a,x)\sin_{\Delta x}(\frac{b}{1+a\Delta x},x)$ and $e_{\Delta x}\cos_{\Delta x}(\frac{b}{1+a\Delta x},x)$, can be converted to equivalent functions containing sines and cosines using the equations on lines 15 and 18 of Table 4.

From Table 4

$$e_{\Delta x}(a,x)\sin_{\Delta x}(\frac{b}{1+a\Delta x},x) = [(1+a\Delta x)^2+(b\Delta x)^2]^{\frac{x}{2\Delta x}}\sin\beta\frac{x}{\Delta x} \tag{1.8-21}$$

$$e_{\Delta x}(a,x)\cos_{\Delta x}(\frac{b}{1+a\Delta x},x) = [(1+a\Delta x)^2+(b\Delta x)^2]^{\frac{x}{2\Delta x}}\cos\beta\frac{x}{\Delta x} \tag{1.8-22}$$

where

$$\beta = \begin{cases} \tan^{-1}\left|\frac{b\Delta t}{1+a\Delta t}\right| & \text{for } 1+a\Delta t > 0 \quad b\Delta t \geq 0 \\ -\tan^{-1}\left|\frac{b\Delta t}{1+a\Delta t}\right| & \text{for } 1+a\Delta t > 0 \quad b\Delta t < 0 \\ \pi - \tan^{-1}\left|\frac{b\Delta t}{1+a\Delta t}\right| & \text{for } 1+a\Delta t < 0 \quad b\Delta t \geq 0 \\ -\pi + \tan^{-1}\left|\frac{b\Delta t}{1+a\Delta t}\right| & \text{for } 1+a\Delta t < 0 \quad b\Delta t < 0 \end{cases} \tag{1.8-23}$$

$1+a\Delta t \neq 0$

$0 \leq \tan^{-1}\left|\frac{b\Delta t}{1+a\Delta t}\right| < \frac{\pi}{2}$

$x = 0, \Delta x, 2\Delta x, 3\Delta x, 4\Delta x, \dots$

Substituting Eq 1.8-21 and Eq 1.8-22 into Eq 1.8-20

$$F(x) = [(1+a\Delta x)^2+(b\Delta x)^2]^{\frac{x}{2\Delta x}}[K_1\cos\beta\frac{x}{\Delta x} + K_2\sin\beta\frac{x}{\Delta x}] \tag{1.8-24}$$

or

$$F(x) = [\sqrt{(1+a\Delta x)^2+(b\Delta x)^2}]^{\frac{x}{\Delta x}}[K_1\cos\beta\frac{x}{\Delta x} + K_2\sin\beta\frac{x}{\Delta x}] \qquad (1.8\text{-}25)$$

where

$$\beta = \begin{cases} \tan^{-1}\left|\frac{b\Delta t}{1+a\Delta t}\right| & \text{for } 1+a\Delta t > 0 \quad b\Delta t \geq 0 \\ -\tan^{-1}\left|\frac{b\Delta t}{1+a\Delta t}\right| & \text{for } 1+a\Delta t > 0 \quad b\Delta t < 0 \\ \pi - \tan^{-1}\left|\frac{b\Delta t}{1+a\Delta t}\right| & \text{for } 1+a\Delta t < 0 \quad b\Delta t \geq 0 \\ -\pi + \tan^{-1}\left|\frac{b\Delta t}{1+a\Delta t}\right| & \text{for } 1+a\Delta t < 0 \quad b\Delta t < 0 \end{cases} \qquad (1.8\text{-}26)$$

$1+a\Delta t \neq 0$

$0 \leq \tan^{-1}\left|\frac{b\Delta t}{1+a\Delta t}\right| < \frac{\pi}{2}$

$x = 0, \Delta x, 2\Delta x, 3\Delta x, 4\Delta x, \ldots$

Thus, collecting and presenting, together, the previously derived equations relatating to the inverse $K_{\Delta x}$ Transform of the $K_{\Delta x}$ Transform function, F(s)

$$\mathbf{F(s) = \frac{K_1(s\text{-}a) + K_2w(1+a\Delta x)}{(s\text{-}a)^2 + (w[1+a\Delta x])^2} = \frac{K_1(s\text{-}a) + K_2b}{(s\text{-}a)^2 + b^2} = \frac{Cs+D}{s^2+As+B}} \qquad (1.8\text{-}27)$$

**The roots of the denominator of Eq 1.8-27 are complex**

**$F(x) = F^{-1}(s)$ , The inverse $K_{\Delta x}$ Transform** (1.8-28)

$$\mathbf{F(x) = e_{\Delta x}(a,x)[K_1\cos_{\Delta x}(\frac{b}{1+a\Delta x},x) + K_2\sin_{\Delta x}(\frac{b}{1+a\Delta x},x)] , \quad 1+ax \neq 0} \qquad (1.8\text{-}29)$$

$$\mathbf{F(x) = [b\Delta x]^{\frac{x}{\Delta x}} [K_1\cos\frac{\pi x}{2\Delta x} + K_2\sin\frac{\pi x}{2\Delta x}] , \quad 1+a\Delta x = 0}$$

**or**

$$\mathbf{F(x) = [(1+a\Delta x)^2+b^2\Delta x^2]^{\frac{x}{2\Delta x}} [K_1\cos\beta\frac{x}{\Delta x} + K_2\sin\beta\frac{x}{\Delta x}]} \qquad (1.8\text{-}30)$$

**where**

**$x = 0, \Delta x, 2\Delta x, 3\Delta x, 4\Delta x, \ldots$**
**$\Delta x$ = x increment**
**$s = a \pm bj$ , the roots of the denominator of F(s), $s^2 + As + B$**
**$a,b,A,B,C,D,K_1,K_2$ = real constants**
**$a = -\frac{A}{2}$**

$$\mathbf{b = \sqrt{B - (\frac{A}{2})^2} = w(1+a\Delta x)}$$

$$\mathbf{w = \frac{b}{(1+a\Delta x)}}$$

$$\mathbf{K_1 = C}$$

$$\mathbf{K_2 = \frac{D - \frac{CA}{2}}{\sqrt{B - (\frac{A}{2})^2}} = \frac{D + Ca}{b}}$$

$$\boldsymbol{\beta = \begin{cases} \tan^{-1}\left|\frac{b\Delta t}{1+a\Delta t}\right| & \textbf{for } 1+a\Delta t > 0 \quad b\Delta t \ge 0 \\ -\tan^{-1}\left|\frac{b\Delta t}{1+a\Delta t}\right| & \textbf{for } 1+a\Delta t > 0 \quad b\Delta t < 0 \\ \pi - \tan^{-1}\left|\frac{b\Delta t}{1+a\Delta t}\right| & \textbf{for } 1+a\Delta t < 0 \quad b\Delta t \ge 0 \\ -\pi + \tan^{-1}\left|\frac{b\Delta t}{1+a\Delta t}\right| & \textbf{for } 1+a\Delta t < 0 \quad b\Delta t < 0 \end{cases}}$$

$$\mathbf{1+a\Delta t \ne 0}$$

$$\mathbf{0 \le \tan^{-1}\left|\frac{b\Delta t}{1+a\Delta t}\right| < \frac{\pi}{2}}$$

Consider the stability of the function, $F(x) = F^{-1}(s)$

In the case where the denominator of $F(s)$ (i.e. $s^2 + As + B$) has complex conjugate roots

From Eq 1.8-30

$$M = [(1+a\Delta x)^2+b^2\Delta x^2)]^{\frac{1}{2}} \tag{1.8-31}$$

For stability, $0 \le M < 1$ since M is raised to the power of $\frac{x}{\Delta x}$ (i.e. 0,1,2,3,…) and it is desired that $M^x \to 0$ as $x \to \infty$.

$$0 \le [(1+a\Delta x)^2+b^2\Delta x^2)]^{\frac{1}{2}} < 1 \tag{1.8-32}$$

From Eq 1.8-32

$$0 \le [(1+a\Delta x)^2 + (b\Delta x)^2] < 1$$

**Then for stability when the denominator roots of the $K_{\Delta x}$ Transform function, F(s), are complex conjugates**

$$\mathbf{0 \le [(1+a\Delta x)^2 + (b\Delta x)^2] < 1} \tag{1.8-33}$$

**where**
**$a \pm jb$ = denominator roots of the $K_{\Delta x}$ Transform function, F(s)**
**$\Delta x$ = x interval**

2) Consider a $K_{\Delta x}$ Transform function, F(s), with the denominator, $s^2+As+B$, having real roots.

$$F(s) = \frac{Cs+D}{s^2+As+B} \tag{1.8-34}$$

where
$A,B,C,D$ = real constants

To find the inverse of the $K_{\Delta x}$ Transform function, F(s), change the form of Eq 1.8-34 to:

$$F(s) = \frac{K_1}{s-a} + \frac{K_2}{s-b} \tag{1.8-35}$$

a,b are the real roots of the denominator of F(s), $s^2 + As + B$
$K_1,K_2$ = real constants

$$s = -\frac{A}{2} - \sqrt{(\frac{A}{2})^2 - B}\ ,\ -\frac{A}{2} + \sqrt{(\frac{A}{2})^2 - B} = a\ ,\ b \tag{1.8-36}$$

$$a = -\frac{A}{2} - \sqrt{(\frac{A}{2})^2 - B} \tag{1.8-37}$$

$$b = -\frac{A}{2} + \sqrt{(\frac{A}{2})^2 - B} \tag{1.8-38}$$

Expanding Eq 1.8-35

$$F(s) = \frac{K_1s - K_1b + K_2s - K_2a}{(s-a)(s-b)} = \frac{(K_1+K_2)s - (K_1b+K_2a)}{(s-a)(s-b)} \tag{1.8-39}$$

Equating Eq 1.8-34 to Eq 1.8-39

$$Cs + D = (K_1+K_2)s - (K_1b+K_2a) \tag{1.8-40}$$

$$K_1 + K_2 = C \tag{1.8-41}$$

$$-bK_1 - aK_2 = D \tag{1.8-42}$$

Solving the simultaneous equations, Eq 1.8-41 and Eq 1.8-42 for the constants $K_1$ and $K_2$

$$K_1 = \frac{D + aC}{a - b} \tag{1.8-43}$$

$$K_2 = -\frac{D + bC}{a - b} \qquad (1.8\text{-}44)$$

Find the inverse $K_{\Delta x}$ Transform of Eq 1.8-35

$$F(x) = K_1 e_{\Delta x}(a,x) + K_2 e_{\Delta x}(b,x) = K_1(1+a\Delta x)^{\frac{x}{\Delta x}} + K_1(1+b\Delta x)^{\frac{x}{\Delta x}}] \qquad (1.8\text{-}45)$$

where

$x = 0, \Delta x, 2\Delta x, 3\Delta x, 4\Delta x, \ldots$

Eq 1.8-45 was obtained from the equations on line 4 of the $K_{\Delta x}$ Transform Table, Table 3, and on line 2 of Table 4 of the Appendix.

Thus, collecting and presenting, together, the previously derived equations relatating to the inverse $K_{\Delta x}$ Transform of the $K_{\Delta x}$ Transform function, F(s)

$$\mathbf{F(s) = \frac{K_1}{s - a} + \frac{K_2}{s - b} = \frac{Cs+D}{s^2+As+B}} \qquad (1.8\text{-}45)$$

**The roots of the denominator of Eq 1.8-45 are real**

$$\mathbf{F(x) = F^{-1}(s)}\text{ , The inverse } \mathbf{K_{\Delta x}} \text{ Transform} \qquad (1.8\text{-}46)$$

$$\mathbf{F(x) = K_1 e_{\Delta x}(a,x) + K_2 e_{\Delta x}(b,x) = K_1(1+a\Delta x)^{\frac{x}{\Delta x}} + K_1(1+b\Delta x)^{\frac{x}{\Delta x}}]} \qquad (1.8\text{-}47)$$

**where**

**$x = 0, \Delta x, 2\Delta x, 3\Delta x, 4\Delta x, \ldots$**

**$\Delta x$ = x increment**

**s = a, b , the roots of $s^2 + As + B$**

**$a,b,A,B,C,D,K_1,K_2$ = real constants**

$$\mathbf{a = -\frac{A}{2} - \sqrt{\left(\frac{A}{2}\right)^2 - B}}$$

$$\mathbf{b = -\frac{A}{2} + \sqrt{\left(\frac{A}{2}\right)^2 - B}}$$

$$\mathbf{K_1 = \frac{D + aC}{a - b}}$$

$$\mathbf{K_2 = -\frac{D+bC}{a - b}}$$

Consider the stability of the function, $F(x) = F^{-1}(s)$

In the case where the denominator of F(s) (i.e. $s^2 + As + B$) has real roots

From Eq 1.8-47

$$M = 1+p\Delta x \tag{1.8-48}$$

where

$p$ = a root, a or b

For stability, $-1<M<1$ since M is raised to the power of $\frac{x}{\Delta x}$ (i.e. 0,1,2,3,…) and it is desired that $M^x \rightarrow 0$ as $x\rightarrow\infty$.

$$-1 < 1 + p\Delta x < 1 \tag{1.8-49}$$

$$-\frac{2}{\Delta x} < p < 0\ ,\ \ p = a \text{ or } b \tag{1.8-50}$$

**Then for stability when the denominator roots of F(s) are real**

$$\boldsymbol{-\frac{2}{\Delta x} < a < 0} \tag{1.8-51}$$

$$\boldsymbol{-\frac{2}{\Delta x} < b < 0} \tag{1.8-52}$$

**where**

**a,b = real roots of the denominator of the $K_{\Delta x}$ Transform function, F(s), Eq 1.8-45**

**$\Delta x$ = x increment**

In general, in control system design, a control system sustained oscillation is not sought. However, the condition for a control system sustained oscillation is known from Eq 1.8-30 which is rewritten below.

Rewriting Eq 1.8-30

$$F(x) = [(1+a\Delta x)^2+b^2\Delta x^2]^{\frac{x}{2\Delta x}}[K_1\cos\beta\frac{x}{\Delta x} + K_2\sin\beta\frac{x}{\Delta x}]$$

where

$K_1,K_2$ = real constants

$$\beta = \begin{cases} \tan^{-1}\left|\frac{b\Delta t}{1+a\Delta t}\right| & \text{for } 1+a\Delta t > 0 \quad b\Delta t \geq 0 \\ -\tan^{-1}\left|\frac{b\Delta t}{1+a\Delta t}\right| & \text{for } 1+a\Delta t > 0 \quad b\Delta t < 0 \\ \pi - \tan^{-1}\left|\frac{b\Delta t}{1+a\Delta t}\right| & \text{for } 1+a\Delta t < 0 \quad b\Delta t \geq 0 \\ -\pi + \tan^{-1}\left|\frac{b\Delta t}{1+a\Delta t}\right| & \text{for } 1+a\Delta t < 0 \quad b\Delta t < 0 \end{cases}$$

$1+a\Delta t \neq 0$

$0 \leq \tan^{-1}\left|\frac{b\Delta t}{1+a\Delta t}\right| < \frac{\pi}{2}$

$x = 0, \Delta x, 2\Delta x, 3\Delta x, 4\Delta x, \ldots$

$\Delta x$ = x increment

**An oscillatory result occurs for F(x), Eq 1.8-30, if:**

$$[(1+a\Delta x)^2+b^2\Delta x^2) = 1 \qquad (1.8\text{-}53)$$

**where**

**$\Delta x$ = x increment**

**$s = a \pm bj$ , the roots of the denominator of the $K_{\Delta x}$ Transform of F(x), $s^2+As+B$**

The system stability equations derived above can be presented graphically. This is shown below

<u>For complex roots, $s = a \pm jb$</u>

Rewriting Eq 1.8-33

For system stability

$0 \leq [(1+a\Delta x)^2 + (b\Delta x)^2] < 1$

Changing the form of Eq 1.8-33

Dividing Eq 1.8-33 by $\Delta x^2$

$$0 \leq [(\frac{1}{\Delta x} + a)^2 + (b)^2] < \frac{1}{\Delta x^2}$$

$$0 \leq [(a + \frac{1}{\Delta x})^2 + (b + 0)^2] < (\frac{1}{\Delta x})^2 \quad , \quad \Delta x \neq 0 \qquad (1.8\text{-}54)$$

Eq 1.8-54, which represents the condition for system stability when the denominator roots of F(s) are complex (a ± jb), is recognized as the equation of concentric circles of radius 0 to $\frac{1}{\Delta x}$ centered at the point, $-\frac{1}{\Delta x}$, in the left half of the complex plane.

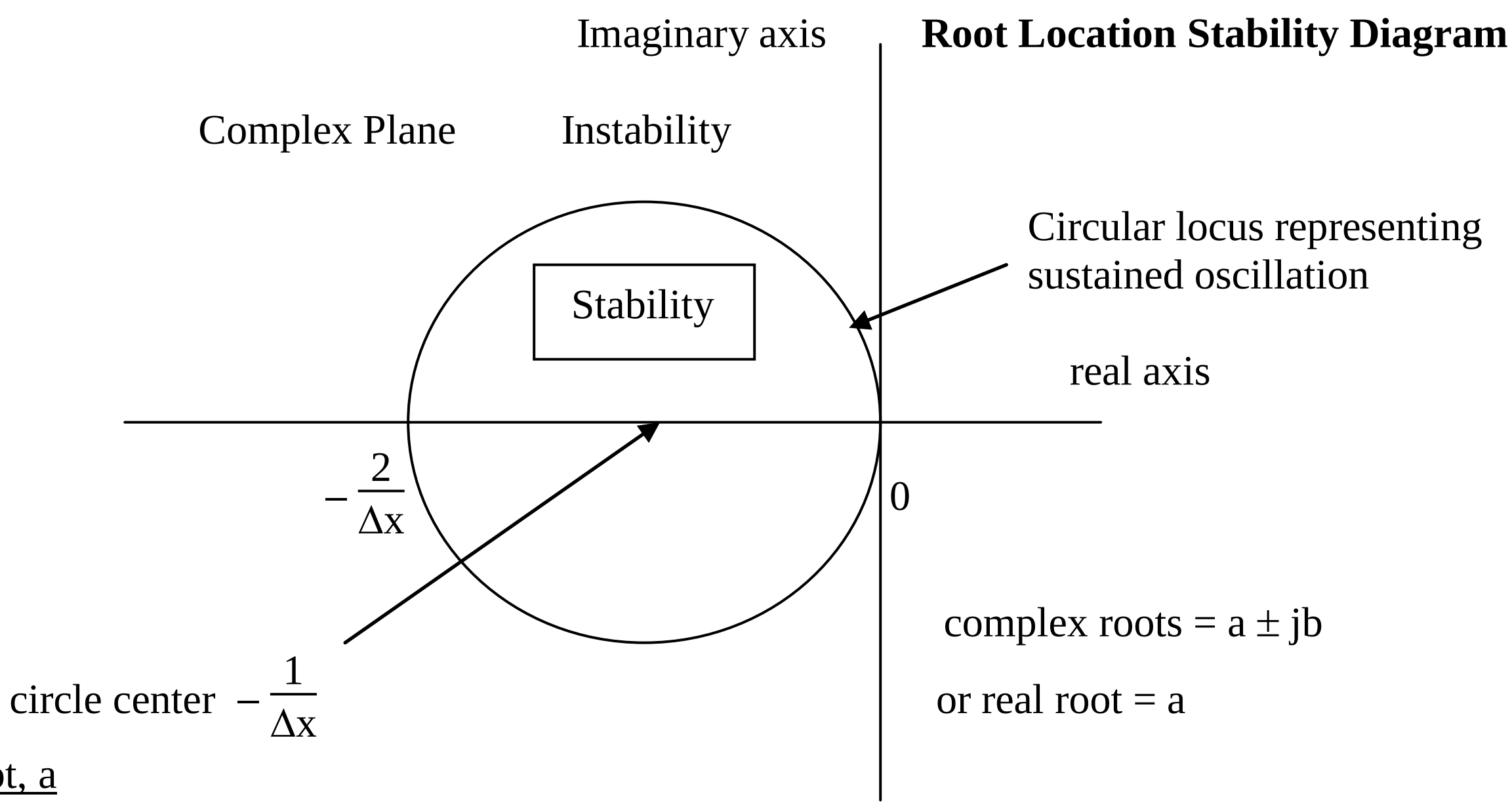

For a real root, a

Rewriting Eq 1.8-51
For system stability

$$-\frac{2}{\Delta x} < a < 0$$

It is observed that the system stability condition for real roots, Eq 1.8-51, is to lie on the negative real axis segment within the specified circular locus between 0 and $-\frac{2}{\Delta x}$.

Thus, for control system stability, a real or complex conjugate root pair in the $K_{\Delta x}$ Transform denominator must lie within the complex plane circle of radius, $\frac{1}{\Delta x}$, centered at the point, $-\frac{1}{\Delta x}$, in the left half of the complex plane. A complex conjugate root pair lying on the specified circular locus or a root at the point ($-\frac{2}{\Delta x}$,0) represents a control system in a purely oscillatory state. A root at the point (0,0) represents a control system in a steady state.

**Proof of $K_{\Delta x}$ Transform discrete variable system analysis Concept #2**
**The modification and use of the Nyquist Criterion to determine control system stability**

Modification of the Nyquist Criterion

To apply the Nyquist Criterion to Operational Interval Calculus control system design, a modification of this Criterion is required. Later, the need for the modification will become evident, but for clarity, only the required modification will be presented first. The modified Nyquist Criterion will be applicable to both continuous and discrete variable control system design applications.

## 1) Review of the Nyquist Criterion

Control System Stability is examined by seeking the presence of any zeros of the function, 1+A(s), in the right half of the s plane.

Consider an open loop Laplace Transform transfer function, A(s).

$$A(s) = \frac{K(s+a)(s+b)(s+c)\ldots}{(s+A)(s+B)(s+C)\ldots} = K\,\frac{A(s)\text{ zeros}}{A(s)\text{ poles}} \qquad (1.8\text{-}55)$$

where
K = real constant

Open loop

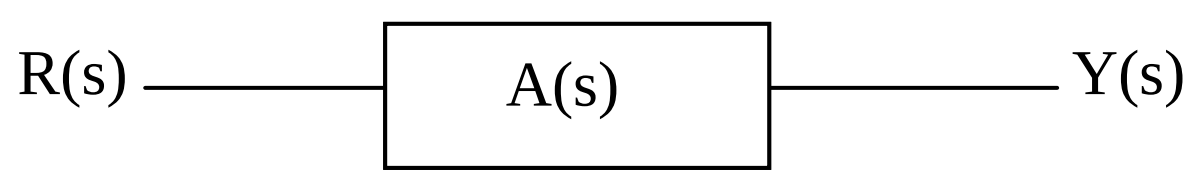

$$\frac{Y(s)}{R(s)} = A(s) \qquad (1.8\text{-}56)$$

where
R(s) = input transform
Y(s) = output transform

Feedback the output to the input

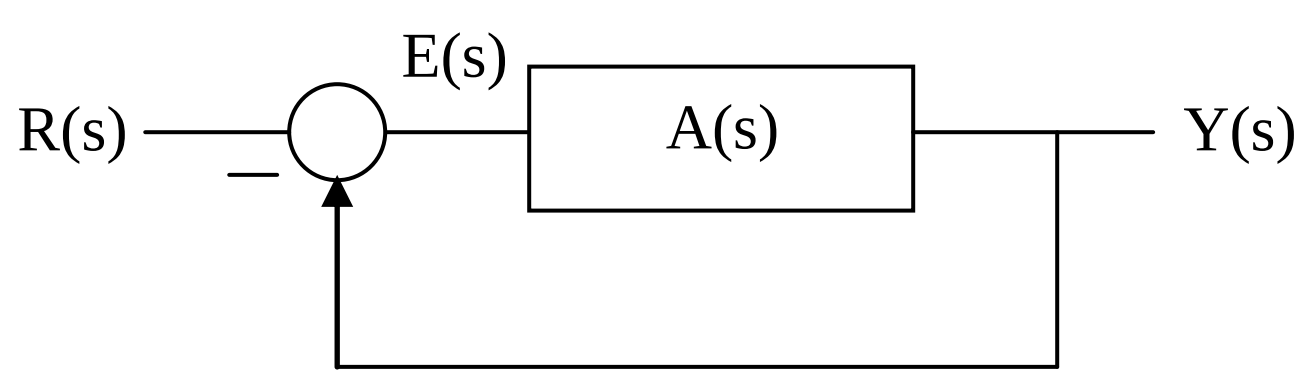

$$E(s) = R(s) - Y(s) \qquad (1.8\text{-}57)$$

where
E(s) = output error transform

$$Y(s) = E(s)A(s) \qquad (1.8\text{-}58)$$

From Eq 1.8-57 and Eq 1.8-58

$$\frac{E(s)}{R(s)} = \frac{1}{1 + A(s)} \qquad (1.8\text{-}59)$$

$$E(s) = \frac{R(s)}{1 + A(s)} \tag{1.8-60}$$

From Eq 1.8-55 and Eq 1.8-60

$$E(s) = \frac{R(s)}{1 + A(s)} = \frac{R(s)}{\dfrac{(s+A)(s+B)(s+C)\ldots+K(s+a)(s+b)(s+c)\ldots}{(s+A)(s+B)(s+C)\ldots}} \tag{1.8-61}$$

$$E(s) = \frac{R(s)}{1 + A(s)} = \frac{R(s)(s+A)(s+B)(s+C)\ldots}{(s+A)(s+B)(s+C)\ldots+ K(s+a)(s+b)(s+c)\ldots} \tag{1.8-62}$$

Often, the variable, x, represents time (t = x). It is desired that, in the time domain, the error approach a constant value, usually a very small value or zero. For this to occur in a continuous time control system, all the poles of the denominator of the transfer function of E(s), Eq 1.8-62, (which are the zeros of 1+ A(s)) must lie within the left half of the complex s plane.

Review the Nyquist Criterion, though perhaps, from a different perspective.

From Eq 1.8-61

$$1+A(s) = \frac{(s+A)(s+B)(s+C)\ldots+K(s+a)(s+b)(s+c)\ldots}{(s+A)(s+B)(s+C)\ldots} \tag{1.8-63}$$

Rewriting Eq 1.8-55

$$A(s) = \frac{K(s+a)(s+b)(s+c)\ldots}{(s+A)(s+B)(s+C)\ldots}$$

<u>Important facts concerning the functions A(s) and 1+A(s)</u>

1. The transfer functions, A(s) and 1+A(s), can have no more zeros than they have poles.
2. The functions 1+A(s) and A(s) have the same poles.
   (See Eq 1.8-55 and Eq 1.8-63.)
3. For control system stability, the zeros of 1+A(s) must be located in the left half of the s plane.

   These zeros are the poles of the E(s) transfer function, $\frac{1}{1+A(s)}$ .

   (See Eq 1.8-62 and Eq 1.8-63.)
4. The function 1+A(s) has the same number of zeros as it has poles. The numerator polynomial order is the same as the denominator polynomial order.

   (Refer to Eq 1.8-63 and fact 1.) above.)
5. If there are no zeros of 1+A(s) within the right half of the s plane, then the zeros must be located in the left half of the s plane or on its imaginary axis. The control system would then be stable or possibly oscillatory if a zero of 1+A(s) is located on the s plane imaginary axis.

Evaluate, in a clockwise direction, 1+A(s) for s = all points on the perimeter of the right half of the s plane. If a pole of 1+A(s) lies on the perimeter, semicircle around it very closely from within the right half of the s plane. A pole of 1+A(s) on the s plane imaginary axis is considered to be within the left half of the s plane.

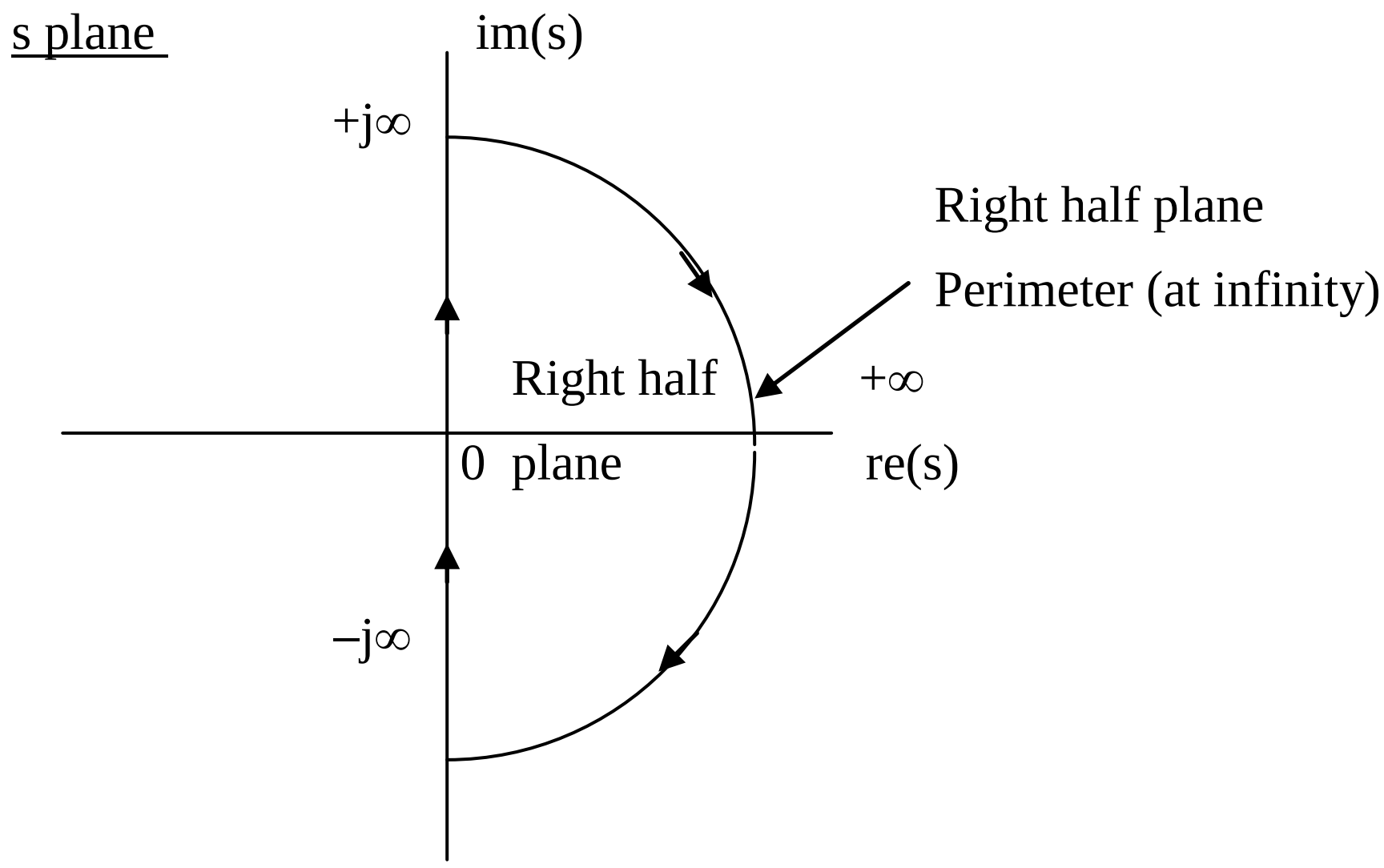

Plot 1+A(s) in the 1+A(s) plane

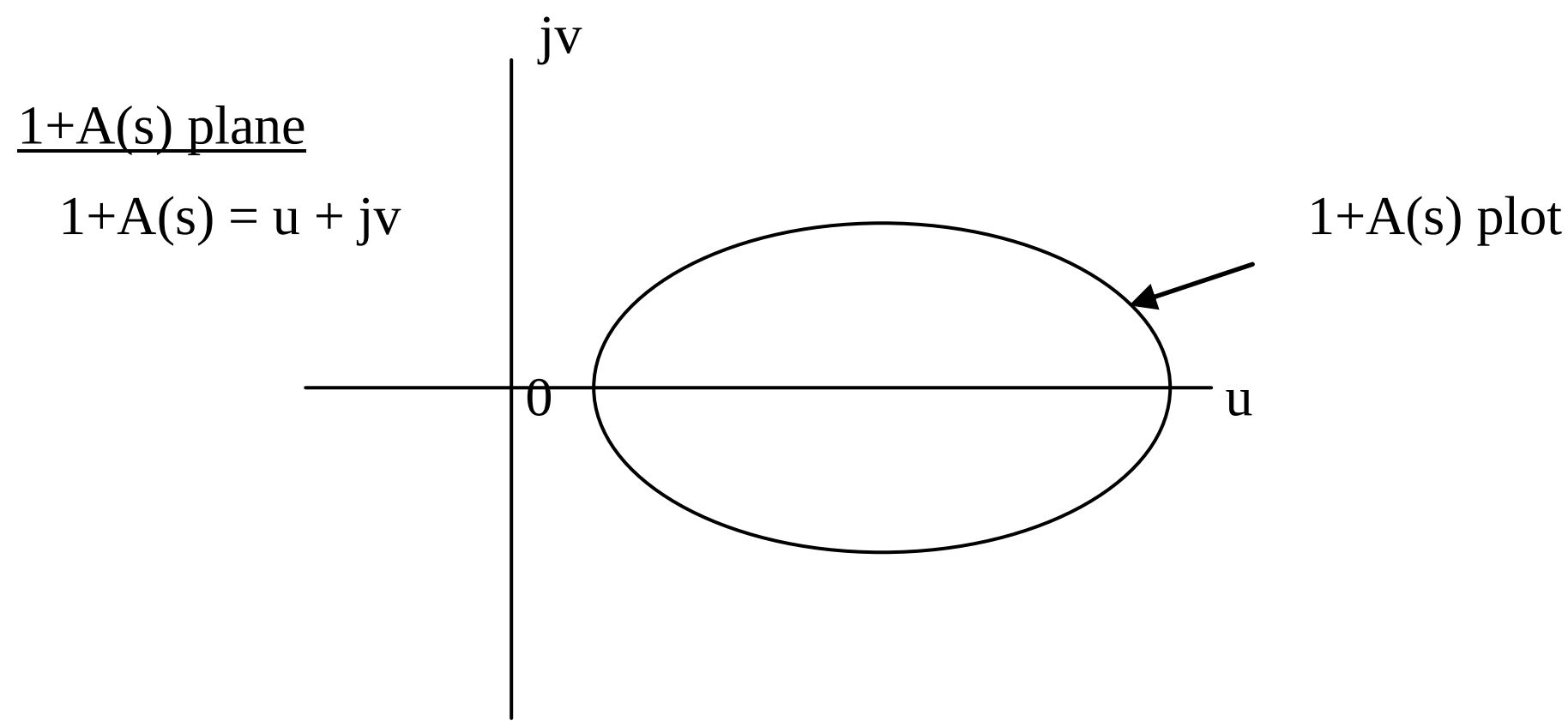

Observing the 1+A(s) plane plot:

If a pole of 1+A(s) is located within the right half of the s plane, a counterclockwise encirclement of the 0 point (origin) will occur in the 1+A(s) plane.

If a zero of 1+A(s) is located within the right half of the s plane, a clockwise encirclement of the 0 point (origin) will occur in the 1+A(s) plane.

If a zero of 1+A(s) is located on the imaginary axis of the s plane, a transition through the 0 point (origin) will occur in the 1+A(s) plane.

The above occurrences are additive so that the magnitude of the number of poles minus the number of zeros of 1+A(s) within the right half of the s plane equals the number of encirclements of the 1+A(s) plane 0 point. A plus sign indicates counterclockwise encirclements and a minus sign indicates clockwise encirclements. The number of transitions through the 0 point equals the number of zeros on the s plane imaginary axis.

The following table is developed to present the conditions for stability.

Let N= the number of poles and zeros of the function, 1+A(s), and

P = the number of poles of 1+A(s) and A(s) in the right half plane. P=0,1,2,3,…,N.

It is desired that the N zeros of the function, 1+A(s), all be located in the left half of the s plane.

A pole of the function, 1+A(s), on the imaginary axis of the s plane is considered to be within the left half of the s plane.

Control system stability will be determined from the number of of encirclements of the 0 point in the 1+A(s) plane.

The following table presents the requirements for all N zeros of 1+A(s) to lie within the left half of the s plane.

**Nyquist Criterion Stability Table**

| Number of poles of 1+A(s) in the left half of the s plane | Number of poles of 1+A(s) in the right half of the s plane | Counterclockwise encirclements of the 1+A(s) plane 0 point (origin) |
|---|---|---|
| **N** | **0** | **0** |
| **N-1** | **1** | **1** |
| **N-2** | **2** | **2** |
| **N-3** | **3** | **3** |
| **N-P** | **P** | **P** |

**Note - Each 1+A(s) plot transition through the 0 point of the 1+A(s) plane indicates a zero of 1+A(s) lying on the s plane imaginary axis.**

<u>There will be no zeros of 1+A(s) in the right half of the s plane if the following conditions are met.</u>

If there are no poles of 1+A(s) in the right half of the s plane, there must be no encirclements of the 1+A(s) plane 0 point.

If there are P poles of 1+A(s) in the right half of the s plane, there must be P counterclockwise encirclements of the 1+A(s) plane 0 point.

<u>There will be a zero of 1+A(s) on the s plane imaginary axis, for each transition through the 0 point in the 1+A(s) plane.</u>

It is useful and convenient not to plot 1+A(s) in the 1+A(s) plane but, rather, to plot A(s) in the A(s) plane. The plots are the same with all values off shifted by one real unit (i.e. 1). The critical 0 point of the 1+A(s) plane becomes the –1 point of the A(s) plane.

Plotting A(s) for the same s values used to plot 1+A(s)

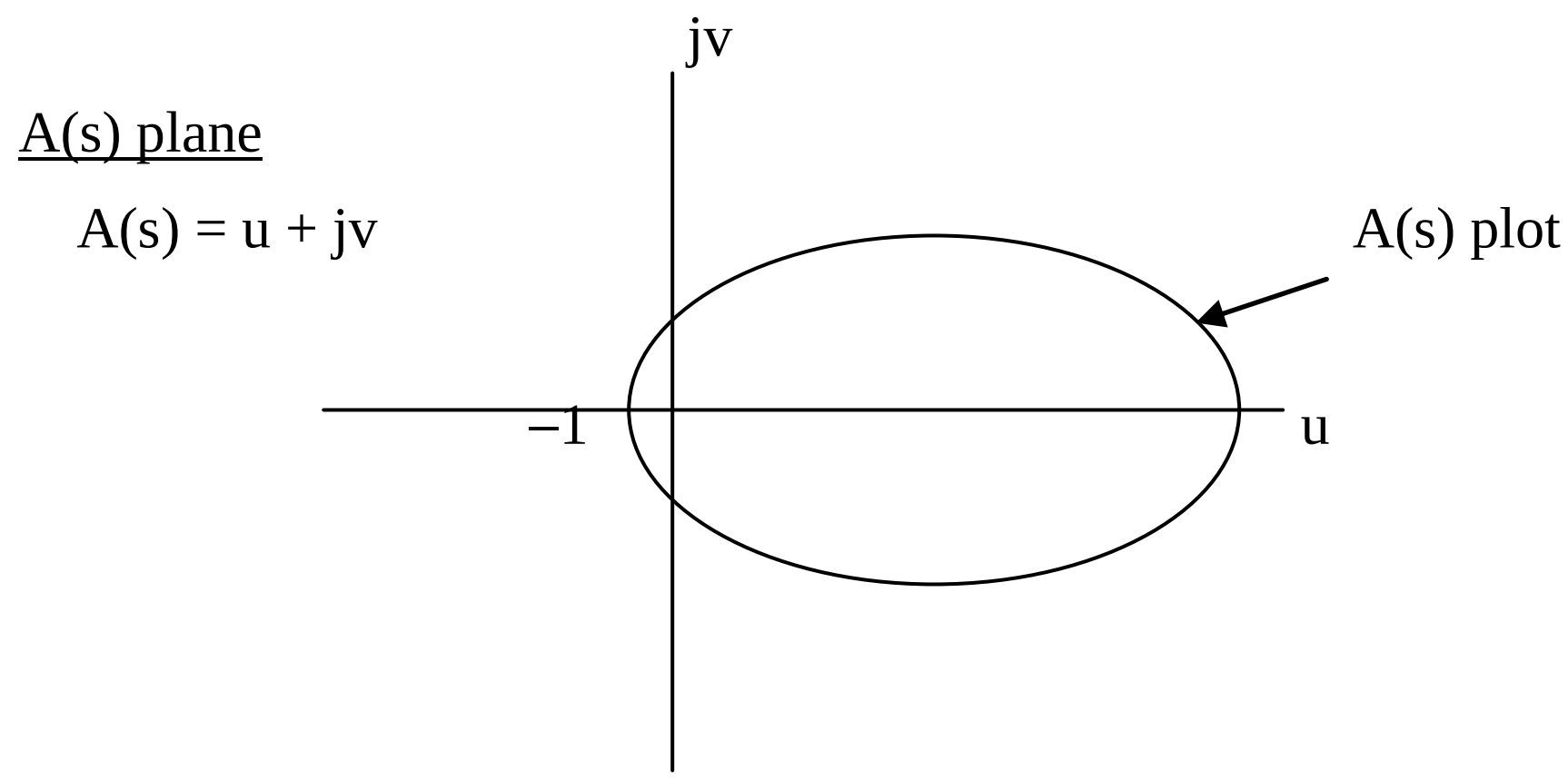

This A(s) plot is seen to be the 1+A(s) plot shifted left one real unit with the 1+A(s) plane 0 point becoming the A(s) plane –1 point. As with the 1+A(s) plot previously discussed, the A(s) plot will not have a zero in the right half of the s plane provided the following conditions are met:

<u>**There will be no zeros of 1+A(s) in the right half of the s plane if the following conditions are met.**</u>

> **If there are no poles of A(s) in the right half of the s plane then there must be no encirclements of the A(s) plane –1 point.**
>
> **If there are P (i.e. P=1,2,3,…,N) poles of A(s) in the right half of the s plane then there must be P counterclockwise encirclements of the A(s) plane -1 point.**

<u>**There will a zero of 1+A(s) on the s plane imaginary axis, for each transition through the –1 point in the 1+A(s) plane.**</u>

**<u>Note</u> – The poles of A(s) are the same as those of 1+A(s).**

<u>**For Stability, all zeros of 1+A(s) must be in the left half of the s plane**</u>

**The above conditions form the basis of the Nyquist Criterion.**

## <u>2) Description of the Modified Nyquist Criterion</u>

Control System Stability is examined by seeking the presence of any zeros of the function, 1+A(s), in the left half of the s plane.

The following derivation uses the same arguments used in obtaining the Nyquist Criterion with one major difference. A(s) is evaluated in a counterclockwise direction for all s on the perimeter of the left half of the s plane.

Consider an open loop $K_{\Delta x}$ Transform transfer function, A(s).

$$A(s) = \frac{K(s+a)(s+b)(s+c)\ldots}{(s+A)(s+B)(s+C)\ldots} = K\frac{A(s)\ zeros}{A(s)\ poles} \quad (1.8\text{-}64)$$

Note that $K_{\Delta x}$ Transform transfer functions look very much like Laplace Transform transfer functions. They are closely related.

Open loop

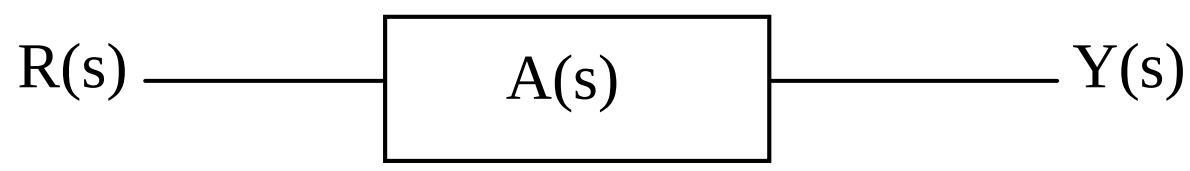

$$\frac{Y(s)}{R(s)} = A(s) \quad (1.8\text{-}65)$$

where
R(s) = input transform
Y(s) = output transform

Feedback the output to the input

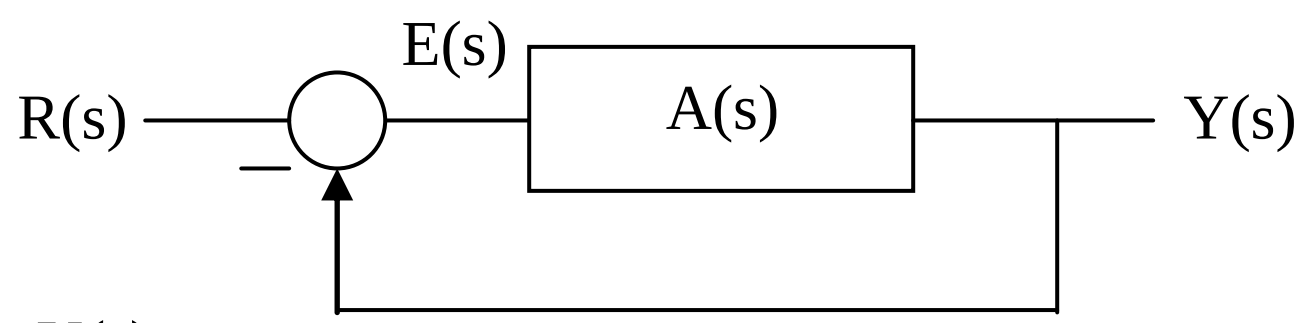

$$E(s) = R(s) - Y(s) \quad (1.8\text{-}66)$$

where
E(s) = output error transform

$$Y(s) = E(s)A(s) \quad (1.8\text{-}67)$$

From Eq 1.8-66 and Eq 1.8-67

$$\frac{E(s)}{R(s)} = \frac{1}{1 + A(s)} \quad (1.8\text{-}68)$$

$$E(s) = \frac{R(s)}{1 + A(s)} \quad (1.8\text{-}69)$$

From Eq 1.8-64 and Eq 1.8-69

$$E(s) = \frac{R(s)}{1 + A(s)} = \frac{R(s)}{\frac{(s+A)(s+B)(s+C)\ldots + K(s+a)(s+b)(s+c)\ldots}{(s+A)(s+B)(s+C)\ldots}} \quad (1.8\text{-}70)$$

$$E(s) = \frac{R(s)}{1 + A(s)} = \frac{R(s)(s+A)(s+B)(s+C)\ldots}{(s+A)(s+B)(s+C)\ldots + K(s+a)(s+b)(s+c)\ldots} \quad (1.8\text{-}71)$$

Note that up to this point, the $K_{\Delta x}$ Transform and the Laplace Transform equations and mathematical manipulations look exactly the same.

Most often, the variable, x, represents time (t = x). It is desired that, in the time domain, the error approach a constant value, usually a very small value or zero. For this to occur in a continuous time control system, all the poles of the denominator of the transfer function of E(s), Eq 1.8-71, (which are the zeros of 1+ A(s)) must lie within the left half of the complex s plane. To find if this is so, here, the Modified Nyquist Criterion is used.

Presenting the Modified Nyquist Criterion.

From Eq 1.8-70

$$1+A(s) = \frac{(s+A)(s+B)(s+C)\ldots+K(s+a)(s+b)(s+c)\ldots}{(s+A)(s+B)(s+C)\ldots} \qquad (1.8\text{-}72)$$

Rewriting Eq 1.8-64

$$A(s) = \frac{K(s+a)(s+b)(s+c)\ldots}{(s+A)(s+B)(s+C)\ldots}$$

Important facts concerning the functions A(s) and 1+A(s)

1. The transfer functions, 1+A(s) and A(s), can have no more zeros than they have poles.
2. The functions 1+A(s) and A(s) have the same poles.
   (See Eq 1.8-64 and Eq 1.8-72.)
3. For control system stability, the zeros of 1+A(s) must be located in the left half of the s plane.
   These zeros are the poles of the E(s) transfer function, $\frac{1}{1+A(s)}$ .
   (See Eq 1.8-71 and Eq 1.8-72.)
4. The function 1+A(s) has the same number of zeros as it has poles. The numerator polynomial order is the same as the denominator polynomial order.
   (Refer to Eq 1.8-72 and fact 1. above.)
5. Stability occurs when all the zeros of 1+A(s) are located in the left half of the s plane.

Evaluate, in a counterclockwise direction, 1+A(s) for s = all points on the perimeter of the left half of the s plane. If a pole of 1+A(s) lies on the perimeter, semicircle it very closely from outside the left half of the s plane. A pole of 1+A(s) on the s plane imaginary axis is considered to be within the left half of the s plane.

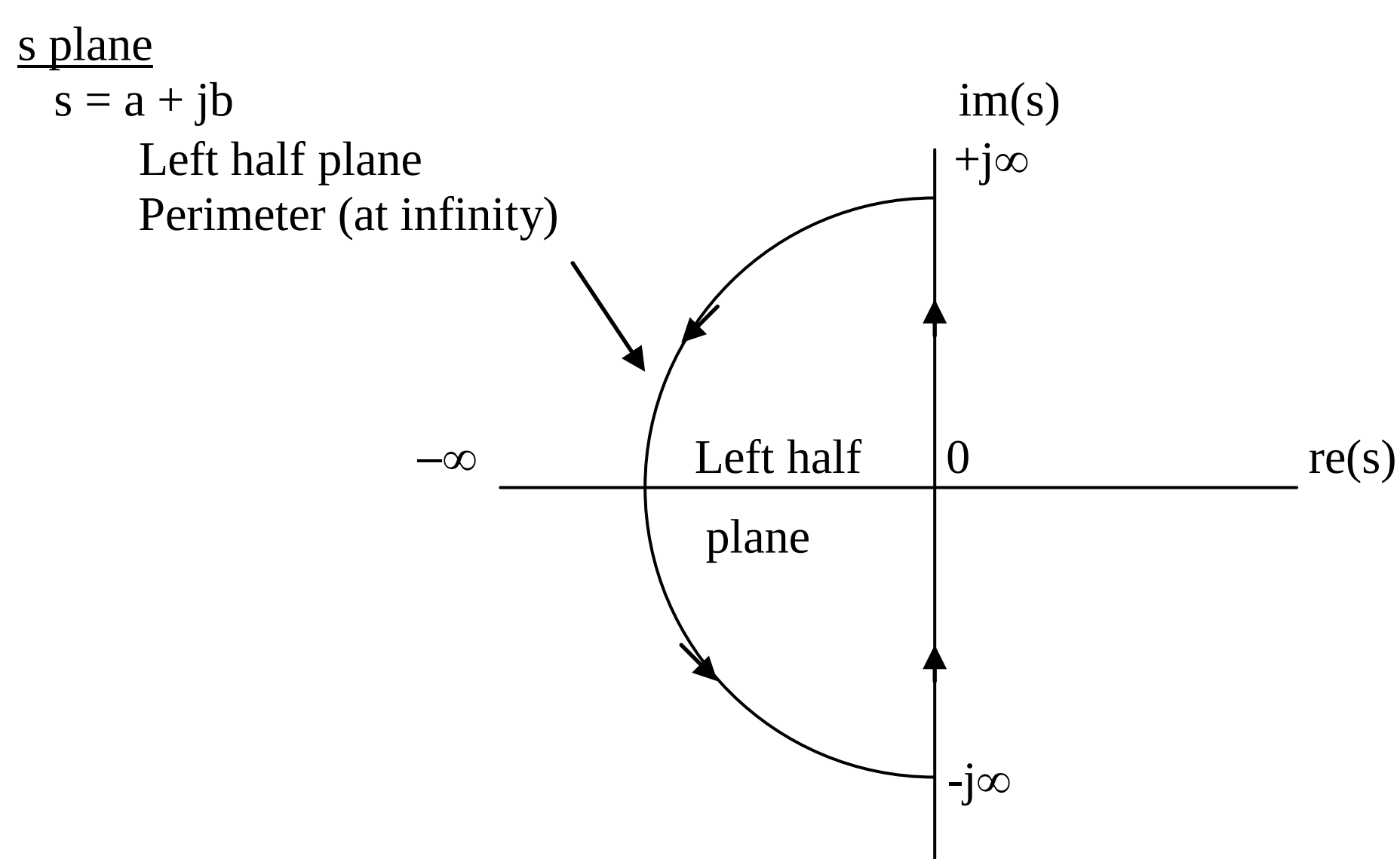

Plot 1+A(s) in the 1+A(s) plane

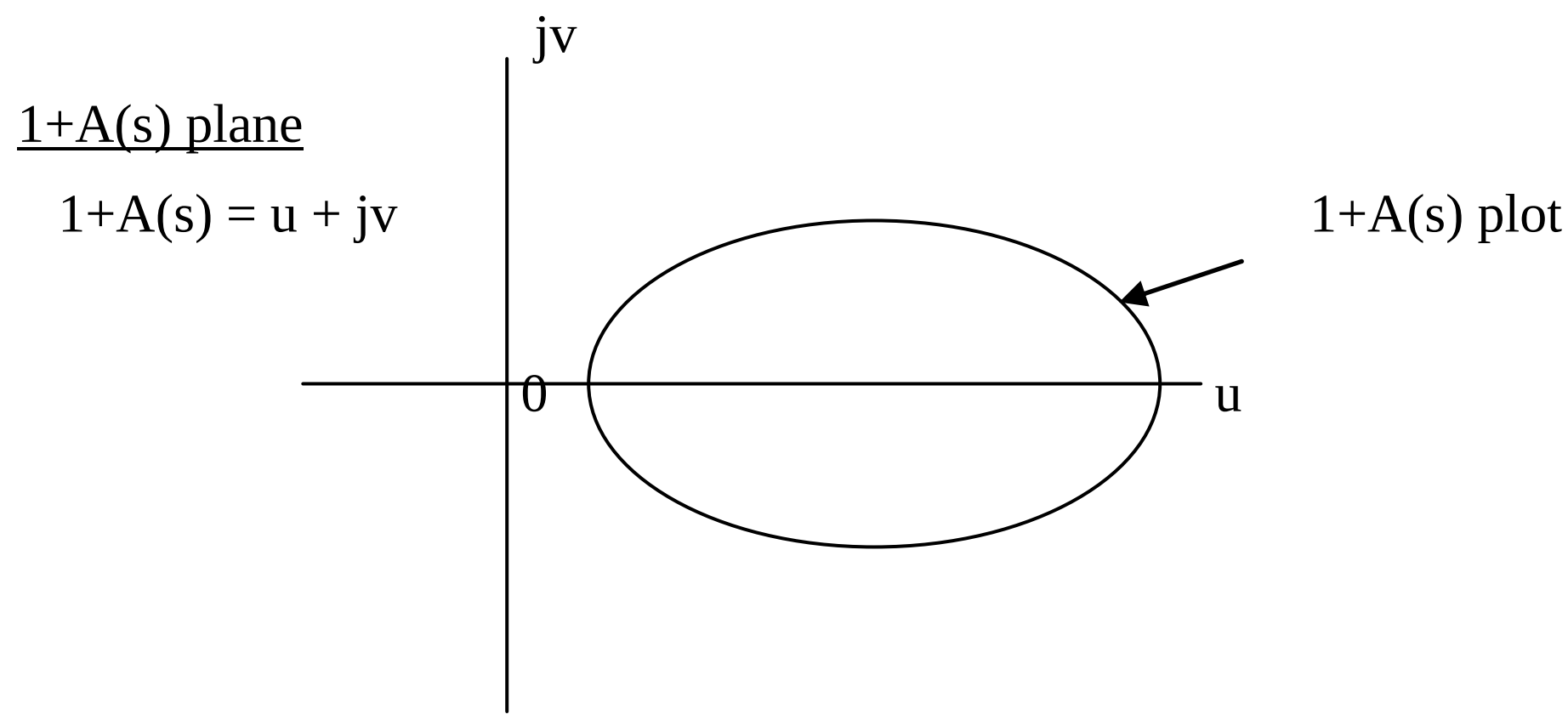

Observing the 1+A(s) plane plot which results from the evaluation of 1+A(s), counterclockwise, for all s on the perimeter of the left half plane:

For each pole of 1+A(s) located within the left half of the s plane, a clockwise encirclement of the 0 point (origin) will occur in the 1+A(s) plane.

For each zero of 1+A(s) located within the left half of the s plane, a counterclockwise encirclement of the 0 point (origin) will occur in the 1+A(s) plane.

For each zero of 1+A(s) located on the imaginary axis of the s plane, a transition through the 0 point (origin) will occur in the 1+A(s) plane.

The above occurrences are are additive so that the magnitude of the number of zeros minus the number of poles of 1+A(s) within the left half of the s plane equals the number of encirclements of the 1+A(s) plane 0 point. A plus sign indicates counterclockwise encirclements and a minus sign indicates clockwise encirclements. The number of transitions through the 0 point equals the number of zeros on the s plane imaginary axis.

As mentioned in fact 4 and fact 5 on the previous page, the function 1+A(s) has the same number of zeros as it has poles and for control system stability, the zeros of 1+A(s) must lie within the left half plane. Also, as mentioned in fact 2, The poles of 1+A(s) are the same as those of A(s). With these facts in mind, the following table is developed to present the conditions for stability.

Let N= the number of poles and zeros in the function,1+A(s), and the number of poles in the function, A(s).

P = the number of poles of 1+A(s) and A(s) in the right half of the s plane. P = 0,1,2,3,…,N.

It is desired that the N zeros of the function, 1+A(s), all be located in the left half of the s plane.

A pole of the function, 1+A(s), on the imaginary axis of the s plane is considered to be within the left half of the s plane.

Control system stability will be determined from the number of of encirclements of the 0 point in the 1+A(s) plane.

The following table presents the requirements for all N zeros of 1+A(s) to lie within the left half of the s plane.

**Modified Nyquist Criterion Stability Table**

| Number of poles of 1+A(s) in the left half of the s plane | Number of poles of 1+A(s) in the right half of the s plane | Counterclockwise encirclements of the 1+A(s) plane 0 point (origin) |
|---|---|---|
| N | 0 | 0 |
| N-1 | 1 | 1 |
| N-2 | 2 | 2 |
| N-3 | 3 | 3 |
| N-P | P | P |

**<u>Note</u> - Each 1+A(s) plot transition through the 0 point of the 1+A(s) plane indicates a zero of 1+A(s) lying on the s plane imaginary axis.**

Thus for control system stability, observing the 1+A(s) plane plot:

where

N = the number of poles of 1+A(s) (or A(s)) and the number of zeros of 1+A(s).

P = the number of poles of 1+A(s) (or A(s)) outside the left half of the s plane.

A pole on the imaginary axis of the s plane is considered to be within the left half of the s plane.

If N poles of 1+A(s) are located within the left half of the s plane, there must be no encirclement of the 1+A(s) plane 0 point (origin).

If N-P poles of 1+A(s) are located within the left half of the s plane, there must be P counterclockwise encirclements of the 1+A(s) plane 0 point (origin).

For each zero of 1+A(s) located on the imaginary axis of the s plane, there will be a transition through the 0 point (origin) of the 1+A(s) plane.

It is useful and convenient not to plot 1+A(s) in the 1+A(s) plane but, rather, to plot A(s) in the A(s) plane. The plots are the same with all values off shifted by one real unit (i.e. 1). In the A(s) plane, the critical 0 point of the 1+A(s) plane becomes the –1 point of the A(s) plane.

Plotting A(s) for the same s values used to plot 1+A(s)

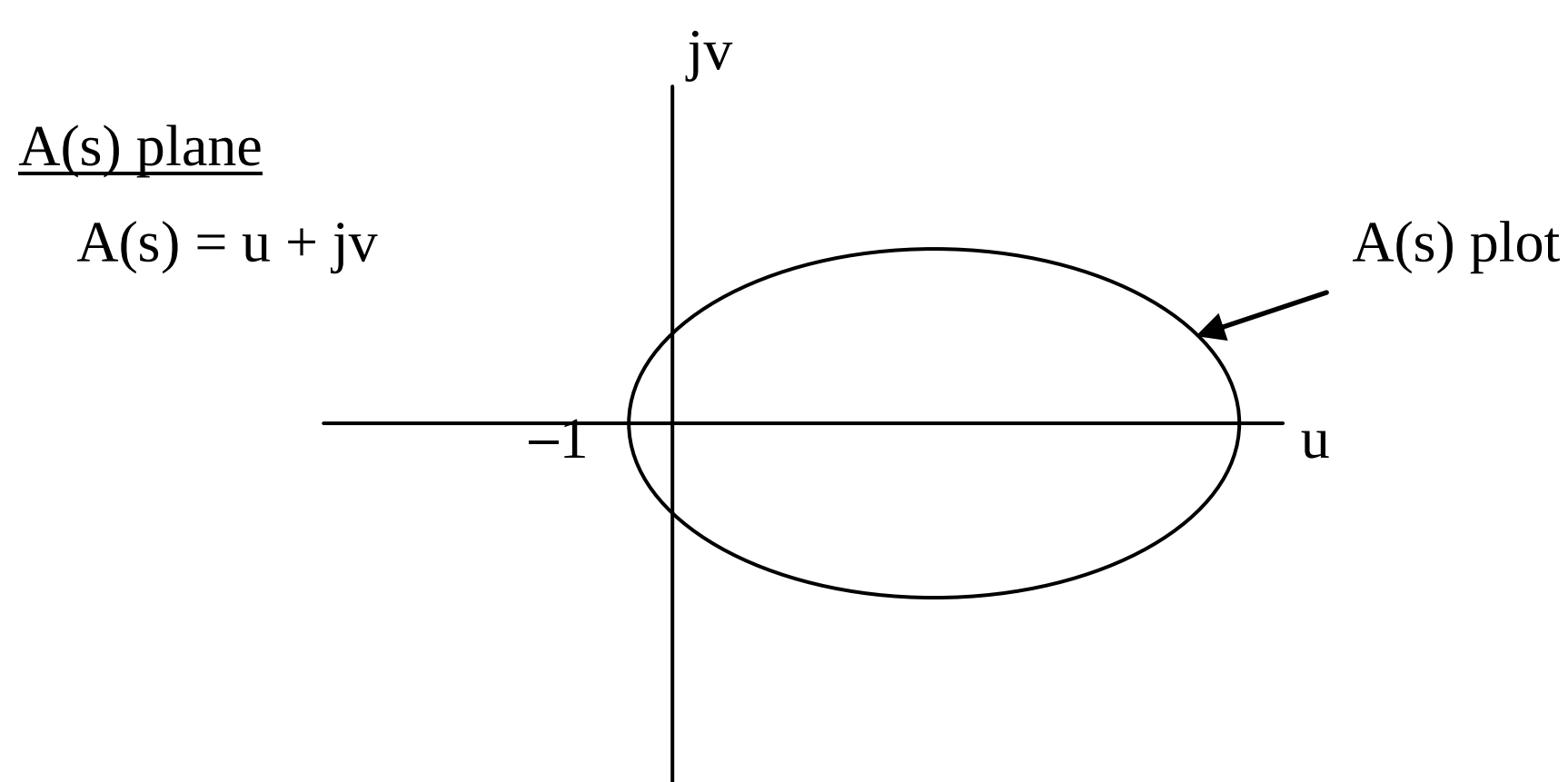

This A(s) plot is seen to be the 1+A(s) plot shifted left one real unit with the 1+A(s) plane 0 point becoming the A(s) plane –1 point. As with the 1+A(s) plot previously discussed, the A(s) plot will not have a zero in the right half of the s plane provided the following conditions are met:

**Evaluate, in a counterclockwise direction, A(s) for s = all points on the perimeter of the left half of the s plane. If a pole of A(s) lies on the perimeter, semicircle it very closely from outside the left half of the s plane. A pole of A(s) on the s plane imaginary axis is considered to be within the left half of the s plane.**

**<u>For control system stability</u>**

**where**

**N = the number of poles of 1+A(s) (or A(s)) and the number of zeros of 1+A(s).**

**P = the number of poles of 1+A(s) (or A(s)) outside the left half of the s plane.**

**If N poles of A(s) are located within the left half of the s plane, there must be no encirclement of the A(s) plane –1 point.**

**If N-P poles of A(s) are located within the left half of the s plane, there must be P counterclockwise encirclements of the A(s) plane –1 point.**

**For each zero of 1+A(s) located on the imaginary axis of the s plane, there will be a transition through the –1 point of the A(s) plane.**

**Note – The poles of A(s) are the same as those of 1+A(s).**

**The above conditions form the basis of the Modified Nyquist Criterion.**

Comparing the A(s) plots of both the Nyquist and Modified Nyquist Criterion

Consider the following diagram of the s plane

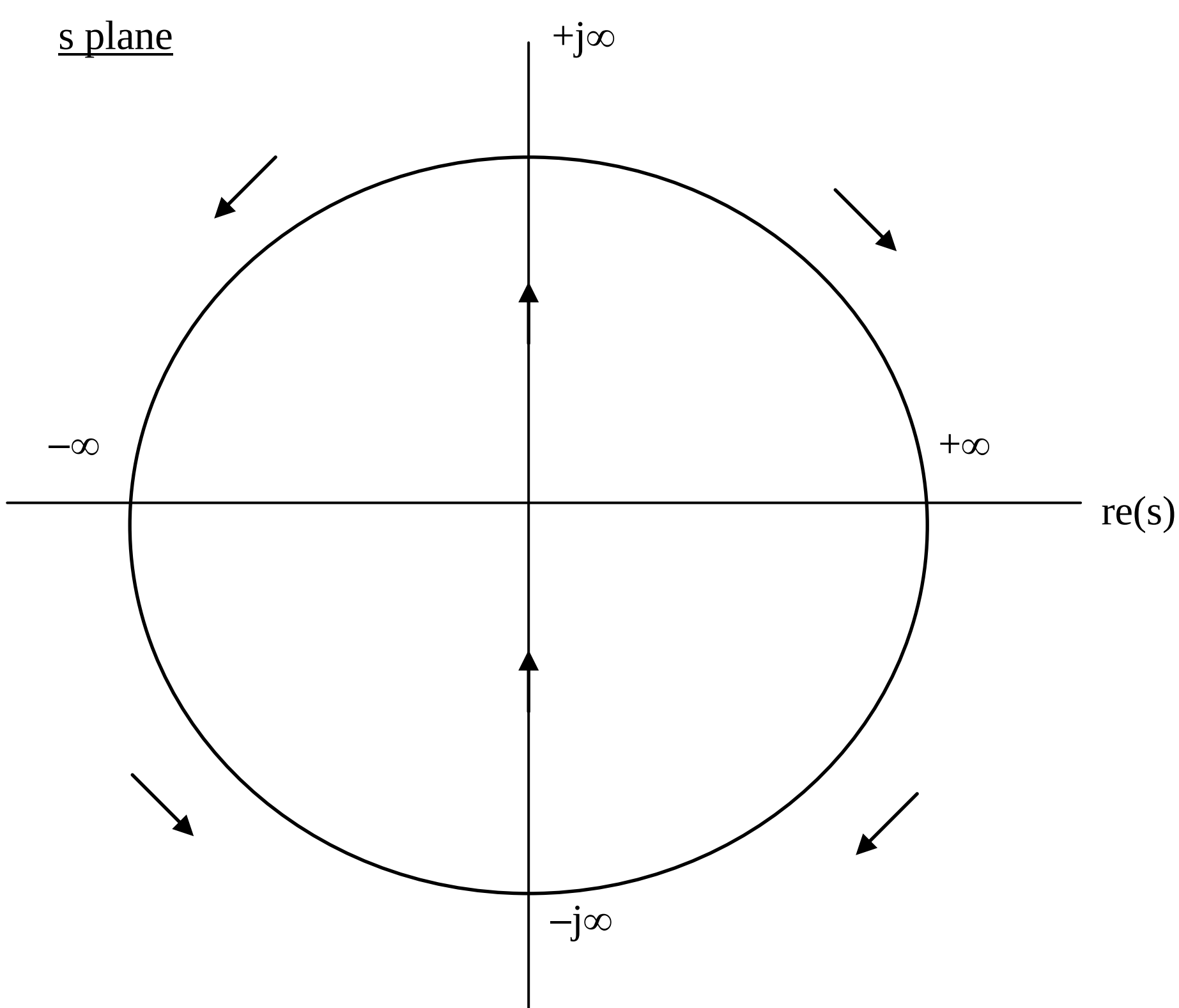

Consider the function, A(s), to have less zeros than poles. Evaluating A(s) for the Nyquist Criterion where s are the values along the periphery of the right half of the s plane, it is seen that the plot of A(s) is essentially the value of A(s) for all values of s on the s plane imaginary axis. The perimeter path other than the imaginary axis is at infinity and thus provides only multiple points of 0 value to the plot of A(s). Evaluating A(s) for the Modified Nyquist Criterion where s are the values along the periphery of the left half of the s plane, it is seen that the plot of A(s) is, again, essentially the value of A(s) for all values of s on the s plane imaginary axis. The perimeter path other than the imaginary axis is at infinity and thus provides only multiple points of 0 value to the plot of A(s). Thus, for the case where the function, A(s), has less zeros than poles, The A(s) plots of the Nyquist Criterion and the Modified Nyquist Criterion are exactly the same.

Consider the function, A(s), to have as many zeros as it has poles. Evaluating A(s) for the Nyquist Criterion where s are the values along the periphery of the right half of the s plane, it is seen that the plot of A(s) for finite values of s is the value of A(s) for all values of s on the s plane imaginary axis. The perimeter path other than the imaginary axis is at infinity so that along this path |A(s)| = 1. Evaluating A(s) for the Modified Nyquist Criterion where s are the values along the periphery of the left half of the s plane, it is seen that the plot of A(s) for finite values of s is the value of A(s) for all values of s on the s plane imaginary axis. The perimeter path other than the imaginary axis is at infinity so that along this path |A(s)| = 1. Thus, for infinite values of s, the Nyquist Criterion and Modified Nyquist Criterion evaluations of A(s) have magnitudes which are exactly the same, |A(s)| = 1. Where s is infinite, the phase change of A(s) for both Criterions must now be determined.

Let N = the number of poles of A(s)

N = the number of zeros of A(s)

P = the number of poles of A(s) in the right half of the s plane

Z = the number of zeros of A(s) in the right half of the s plane

Any pole or zero of A(s) on the imaginary axis is considered to be within the left half of the s plane.

N–P = the number of poles of A(s) in the left half plane

N–Z = the number of zeros of A(s) in the left half plane

P = 0,1,2,3,…,N

Z = 0,1,2,3,…,N

A clockwise evaluation of A(s) for the Nyquist Criterion, where s are the infinite values along the periphery of the right half of the s plane and P poles and Z zeros lie within the right half of the s plane, results in a counterclockwise phase change (at a magnitude of 1) of 180*(P-Z) degrees. A counterclockwise evaluation of A(s) for the Modified Nyquist Criterion, where s are the infinite values along the periphery of the left half of the s plane and N–P poles and N-Z zeros lie within the left half of the s plane, results in a counterclockwise phase change (at a magnitude of 1) of 180*(N-Z–[N–P]) degrees =180*(P-Z) degrees. Then, for infinite values of s, Nyquist Criterion and Modified Nyquist Criterion evaluations of A(s) have equal phase changes. Thus for the case where the function, A(s), has an equal number of zeros and poles, the A(s) plots of the Nyquist Criterion and the Modified Nyquist Criterion look exactly the same.

In summation, the Nyquist Criterion and the Modified Nyquist Criterion when applied to the same continuous variable transform function, A(s), yield identical results with respect to both the transform A(s) plane plot and the stability analysis. It should be noted that the stability criteria for both the Nyquist and Modified Criterion are equivalent even though the the wording of each may at first appear different. Compare the previously presented Nyquist Criterion Stability Table to the Modified Nyquist Criterion Stability Table.They are the same.

Thus, the Modified Nyquist Criterion for continuous variable control systems is as follows:

**The Modified Nyquist Criterion as it applies to continuous variable control systems ($\Delta x \rightarrow 0$)**

**Evaluate, in a counterclockwise direction, A(s) for s = all points on the perimeter of the left half of the s plane. If a pole of A(s) lies on the perimeter, semicircle it very closely from outside the left half of the s plane. A pole of A(s) on the s plane imaginary axis is considered to be within the left half of the s plane.**

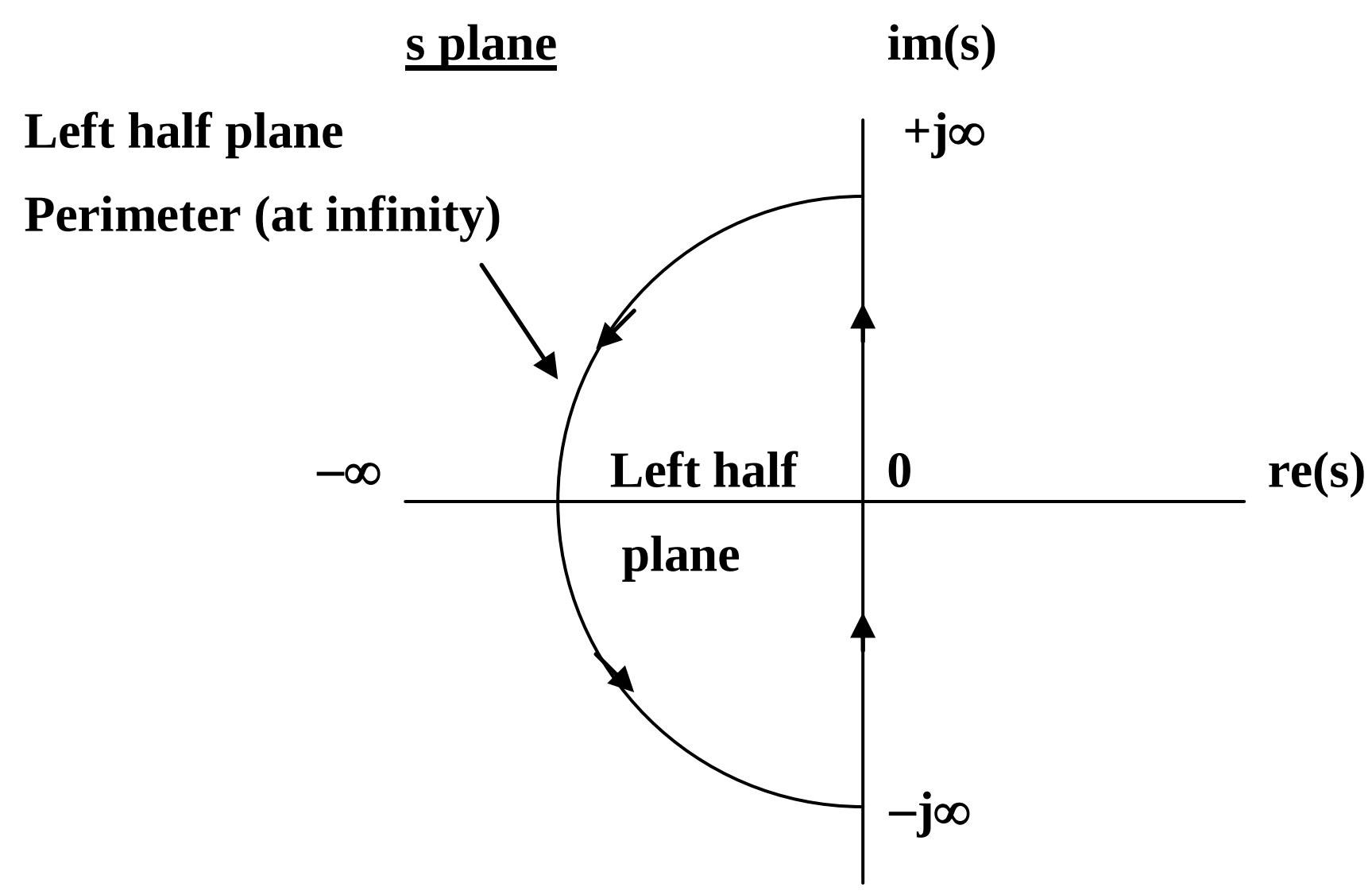

**For control system stability**

**where**

**N = the number of poles of 1+A(s) (or A(s)) and the number of zeros of 1+A(s).**

**P = the number of poles of 1+A(s) (or A(s)) outside the left half of the s plane.**

**If N poles of A(s) are located within the left half of the s plane, there must be no encirclement of the A(s) plane –1 point.**

**If N-P poles of A(s) are located within the left half of the s plane, there must be P counterclockwise encirclements of the A(s) plane –1 point.**

**For each zero of 1+A(s) located on the imaginary axis of the s plane, there will be a transition through the –1 point of the A(s) plane.**

**<u>Notes</u> – The poles of A(s) are the same as those of 1+A(s).**

**For the functions, A(s) and 1+A(s), there can be no more zeros than poles.**

**Comment – The Modified Nyquist Criterion as applied to continuous variable control systems yields exactly the same A(s) plot and stability analysis data as the Nyquist Criterion.**

**<u>3) Application of the Modified Nyquist Criterion to discrete variable control systems</u>**

The Modified Nyquist Criterion, as previously described, is readily adapted for use with discrete (sampled) independent variable control system analysis. Instead of evaluating the transfer function, A(s), along the perimeter of the left half of the s plane, A(s) is evaluated along the perimeter of a circle lying within the left half of the s plane. This circle, henceforth, will be named the stability Critical Circle.

Most often, the variable, x, represents time (t = x). It is desired that, in the time domain, feedback error approach a constant value, usually a very small value or zero. As previously proven, for this to occur in a discrete (sampled) time control system, all the poles of the denominator of the transfer function of E(s), Eq 1.8-71, (which are the zeros of 1+ A(s)) must lie in the left half of the complex s plane within the Critical Circle of radius $\frac{1}{\Delta x}$ centered at $-\frac{1}{\Delta x}$. With the exception of the change from the perimeter of the left half of the s plane for continuous variable system analysis to the Critical Circle within the s plane for discrete variable system analysis, all of the Modified Nyquist Criterion arguments apply equally well.

Thus, the Modified Nyquist Criterion for discrete variables is as follows:

**The Modified Nyquist Criterion as it applies to discrete (sampled) indepenent variable control systems**
**($\Delta x \neq 0$)**

**Evaluate, in a counterclockwise direction, A(s) for s = all points on the perimeter of the Critical Circle of radius $\frac{1}{\Delta x}$ centered at $-\frac{1}{\Delta x}$. If a pole of A(s) lies on the Critical Circle perimeter, semicircle it very closely from outside the Critical Circle. A pole of A(s) on the Critical Circle perimeter is considered to be within the Critical Circle.**

**Note - $\Delta x$ , the x increment, may be designated $\Delta t$, a time increment, for time related variables.**

**Diagram A**

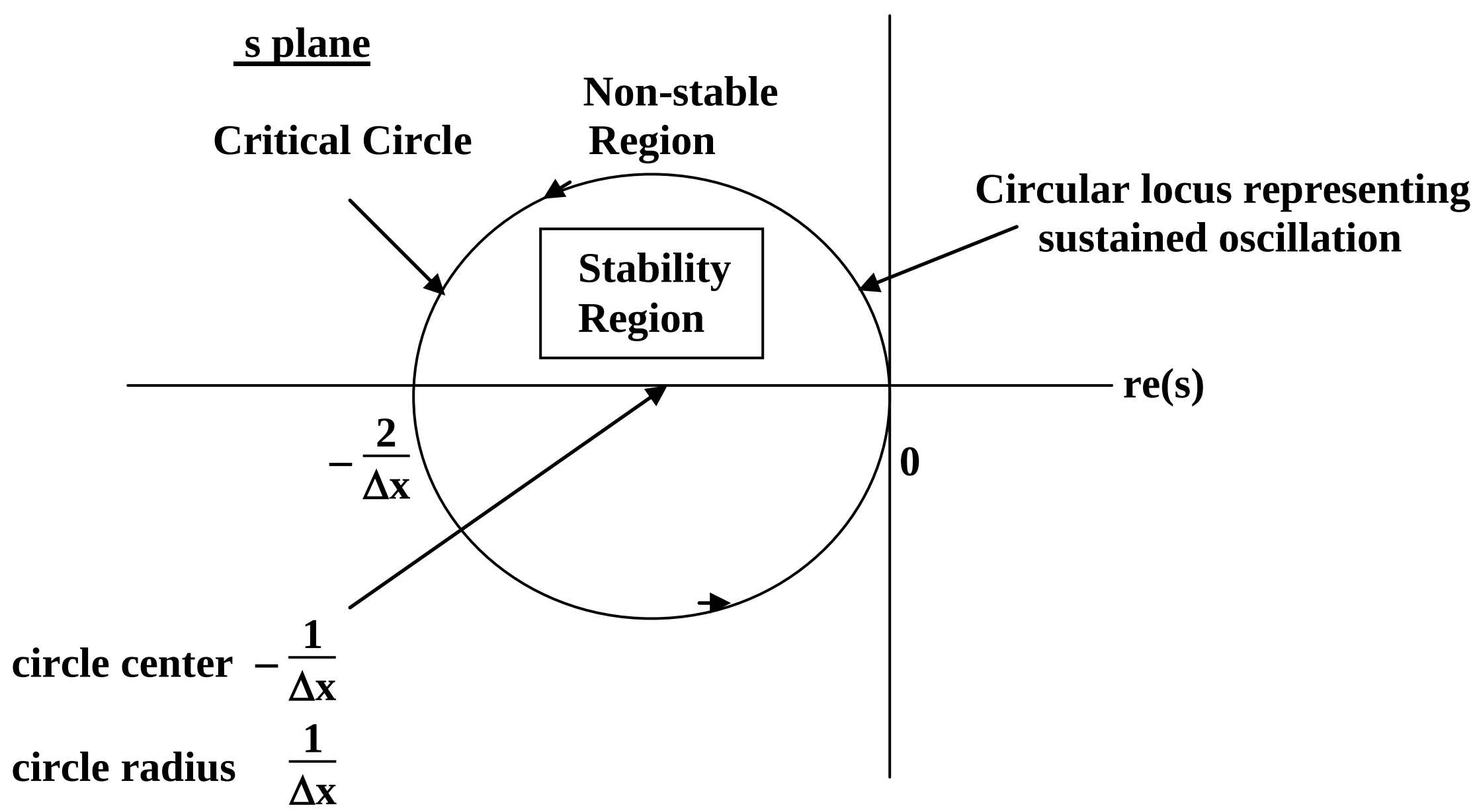

**For control system stability**

**where**

**N = the number of poles of 1+A(s) (or A(s)) and the number of zeros of 1+A(s).**

**P = the number of poles of 1+A(s) (or A(s)) located outside the Critical Circle which lies within the left half of the s plane.**

**If N poles of A(s) are located within the Critical Circle of the s plane, there must be no encirclement of the A(s) plane –1 point.**

**If N-P poles of A(s) are located within the Critical Circle of the s plane, there must be P counterclockwise encirclements of the A(s) plane –1 point.**

**For each zero of 1+A(s) located on the Critical Circle perimeter, there will be a transition through the –1 point of the A(s) plane.**

**Notes – The poles of A(s) are the same as those of 1+A(s).**
**For the functions, A(s) and 1+A(s), there can be no more zeros than poles.**

**4) Application of the Bilinear Transformation to the Modified Nyquist Criterion**

The Bilinear Transform, of interest, is as follows:

$$z = \frac{1 + w}{1 - w} \qquad (1.8\text{-}73)$$

The Bilinear Transformation does the following:

1. Maps the perimeter of the z plane unit circle into the imaginary axis of the w plane.
2. Maps the z plane area within the unit circle centered at 0 into the left half of the w plane.
3. Maps the z plane area outside the unit circle centered at 0 into the right half of the w plane.

Diagram B

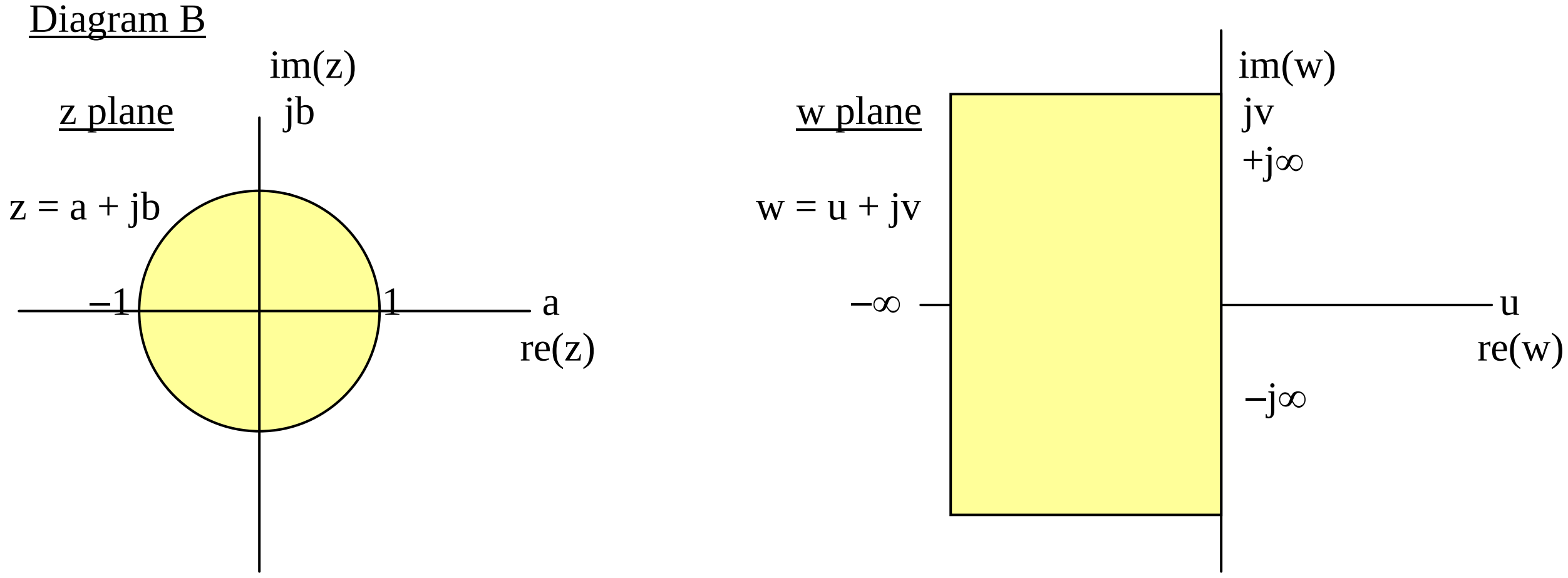

Referring to Diagram B and to Diagram A on the previous page, map the s plane Critical Circle to the z plane 0 centered unit circle. Then, map the resulting z plane to the w plane using the Bilinear Transformation

In particular, center the Critical Circle at 0 instead of at $-\frac{1}{\Delta x}$. This is done by moving all of the values of s to the right by $\frac{1}{\Delta x}$ by adding $\frac{1}{\Delta x}$ to s. Then, reduce the radius of the Critical Circle from $\frac{1}{\Delta x}$ to 1 by multiplying $(s + \frac{1}{\Delta x})$ by Δx. The result is the unit circle in the z plane. Finally, map this interim z plane into the w plane using the Bilinear Transformation.

$$\Delta x\left( s + \frac{1}{\Delta x} \right) = z \qquad (1.8\text{-}74)$$

Substituting Eq 1.8-73 into Eq 1.8-74

$$\Delta x\left( s + \frac{1}{\Delta x} \right) = \frac{1 + w}{1 - w} \qquad (1.8\text{-}75)$$

Solve for s

$$s\Delta x + 1 = \frac{1 + w}{1 - w}$$

$$s\Delta x = \frac{1 + w}{1 - w} - 1 = \frac{2w}{1-w} \quad s = \frac{2w}{\Delta x(1-w)} \qquad (1.8\text{-}76)$$

Eq 1.8-76 maps the interior of the s plane Critical Circle into the left half of the w plane, the exterior of the s plane Critical Circle into the right half of the w plane, and the circular perimeter of the s plane Critical Circle into the imaginary axis of the w plane. Substituting Eq 1.8-76 into a function of s creates a function of w to which the Bode methods of stability compensation can be applied.

In the application of the Modified Nyquist Criterion to discrete variable control systems, points on the perimeter of the Critical Circle were at first programmed as follows:

$$s = \frac{-1 + e^{j\theta}}{\Delta x} \ , \quad 0 \le \theta \le 2\pi \qquad (1.8\text{-}77)$$

The above equation, of course, works well. However, points on the perimeter of the Critical Circle may be calculated in a different and very interesting way. Eq 1.8-76 can be used.

$$s = \frac{2w}{\Delta x(1-w)} \ , \quad -j\infty < w < +j\infty \qquad (1.8\text{-}78)$$

where

$$\text{for} \quad s = 0\angle 0 \rightarrow -\frac{2}{\Delta x}\angle \pi^{-} \quad , \quad w = 0 \rightarrow j\infty$$

$$s = -\frac{2}{\Delta x}\angle \pi^{+} \rightarrow 0\angle 2\pi \quad , \quad w = -j\infty \rightarrow 0$$

The values of s calculated by Eq 1.8-77 and Eq 1.8-78 are the same.
The previously discussed methods for calculating Nyquist Criterion plots are listed below:

| # | Equations | Input Variable Commentary | Stability Criterion Used and Periphery |
|---|---|---|---|
| **1** | $s = \frac{-1 + e^{j\theta}}{\Delta x} \ , \quad 0 \le \theta \le 2\pi$ <br> Function = A(s) | The equation s values are on the perimeter of the Critical Circle in the left half of the s plane | Modified Nyquist Periphery – Critical Circle in the left half of the s plane |
| **2** | $s = \frac{2w}{\Delta x(1-w)} \ , \ -j\infty < w < +j\infty$ <br> Function = A(s) | The equation s values are on the perimeter of the Critical Circle in the left half of the s plane | Modified Nyquist Periphery – Critical Circle in the left half of the s plane |
| **3** | $-j\infty < w < +j\infty$ <br> Function = $A(s)\big\|_{s = \frac{2w}{\Delta x(1-w)}}$ | The specified w values are on the imaginary axis of the w plane | Modified Nyquist Periphery – Left half of the w plane |

Notes – 1. Methods #1 and #2 above differ only on how the Critical Circle perimeter values of s are calculated. In either case, the values of s are the same.
2. The equations of Methods #2 and #3 are the same even though their planes and peripheries of interest differ.
3. Method #3 is useful for the application of Bode's methods for stability compensation
4. The Nyquist plots of all three methods are the same.

In the following example, analyze the stability of a discrete variable control system using the Modified Nyquist Criterion.

Example 1.8-1 Investigate the stability of the specified discrete variable control system. Determine the system stability phase margin for the gain constant, K = 1. Find the value of the gain constant, K, at which the system is oscillatory. The x increment, $\Delta x$, is .1 and the initial values of y(x) are y(0)=y(.1)=y(.2)=0.

Note – Often the variable, x, represents time and, if so, is renamed, t.

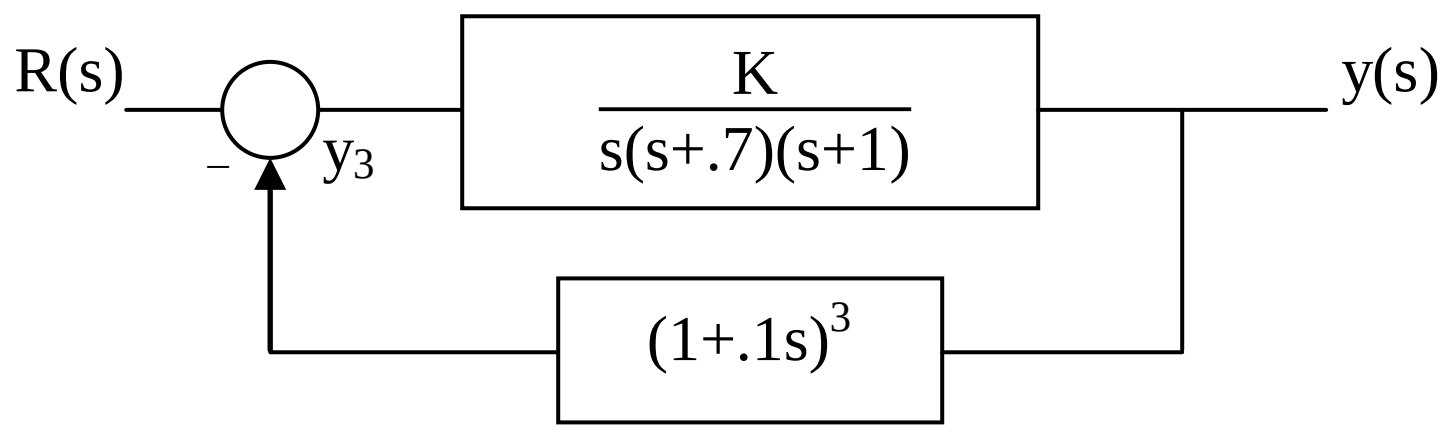

From the specified discrete feedback system diagram, write the transfer function for $\frac{y(s)}{R(s)}$.

$$\frac{y(s)}{R(s)} = \frac{\frac{K}{s(s+.7)(s+1)}}{1 + \frac{K(1+.1s)^3}{s(s+.7)(s+1)}} \tag{1.8-79}$$

$$A(s) = \frac{K(1+.1s)^3}{s(s+.7)(s+1)} \tag{1.8-80}$$

From the specified discrete feedback system diagram and the following $K_{\Delta x}$ Transforms, write the equation for y(x). Let R be a step input where R = 5.

$$K_{\Delta x}[D^n{}_{\Delta x}f(x)] = s^n K_{\Delta x}[f(x)] - s^{n-1}D^0{}_{\Delta x}f(0) - s^{n-2}D^1{}_{\Delta x}f(0) - s^{n-3}D^2{}_{\Delta x}f(0) - \ldots - s^0 D^{n-1}{}_{\Delta x}f(0)\ , \quad n = 1,2,3,\ldots \tag{1.8-81}$$

$$K_{\Delta x}[f(x+n\Delta x)] = (1+s\Delta x)^n K_{\Delta x}[f(x)] - \Delta x \sum_{m=1}^{n}(1+s\Delta x)^{n-m}f([n-m]\Delta x)\ , \quad n = 1,2,3,\ldots \tag{1.8-82}$$

$$[s(s+.7)(s+1)]y(s) = K[\frac{R}{s} - (1+.1s)^3y(s)] \tag{1.8-83}$$

$$[s^3 + 1.7s^2+.7s]y(s) = K[\frac{5}{s} - (1+.1s)^3y(s)] \tag{1.8-84}$$

The initial values are y(0)=y(.1)=y(.2)=0

y =y(x)

$x = n\Delta x$

$y_n = y(n\Delta x)$

$$D_{\Delta x}y(x) = \frac{y(x+\Delta x) - y(x)}{\Delta x} = \frac{y([n+1]\Delta x) - y(n\Delta x)}{\Delta x} = \frac{y_{n+1} - y_n}{\Delta x}$$ , The discrete derivative

From Eq 1.8-81, 1.8-82, and Eq 1.8-84

$\Delta x$ =.1

$$D_{\Delta x}{}^3 y(x) + 1.7D_{\Delta x}{}^2 y(x) + .7D_{\Delta x}{}^1 y(x) = K[5 - y(x+3\Delta x)] \tag{1.8-85}$$

$$\frac{y_{n+3}-3y_{n+2}+3y_{n+1}-y_n}{.1^3} + 1.7\left(\frac{y_{n+2}-2y_{n+1}+y_n}{.1^2}\right) + .7\left(\frac{y_{n+1}-y_n}{.1}\right) = 5K-Ky_{n+3} \tag{1.8-86}$$

$$y_{n+3}-3y_{n+2}+3y_{n+1}-y_n + .17(y_{n+2}-2y_{n+1}+y_n) + .007(y_{n+1}-y_n) = .005K-.001Ky_{n+3} \tag{1.8-87}$$

Combining the terms of Eq 1.8-87 and solving for $y_3$

$$y_{n+3} = \frac{2.83y_{n+2} -2.667y_{n+1}+.837y_n + .005K}{1+.001K} \tag{1.8-88}$$

A program to plot the Modified Nyquist Diagram of A(s) and to calculate the phase margin of A(s) was written. Also, a program to plot y(x) from a y(x) difference equation was written. These two programs appear in the Calculation Program section of the Appendix as programs 6 and 7.

From the A(s) equation, Eq 1.8-80, and the y(x) difference equation, Eq 1.8-88, The following plots were obtained.

The left half s plane Critical Circle has a radius of $\frac{1}{\Delta x} = \frac{1}{.1} = 10$ with a center at $s = -\frac{1}{\Delta x} = -\frac{1}{.1} = -10$ .

1) <u>For K=1, Δx =.1</u>

From Eq 1.8-80

$$A(s) = \frac{1(1+.1s)^3}{s(s+.7)(s+1)}$$

Evaluating A(s) around the left half s plane Critical Circle of radius 10 and center at s = –10

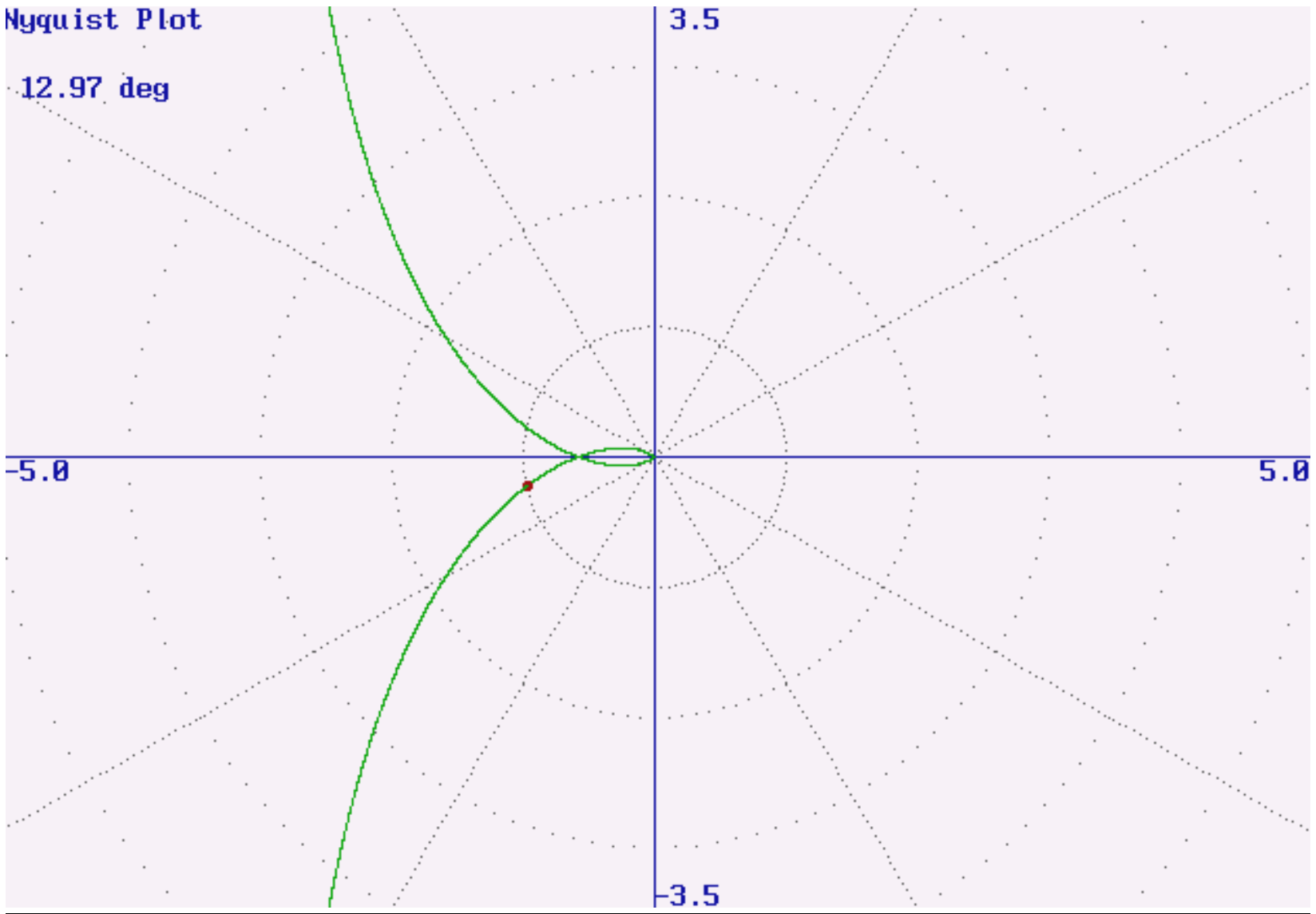

For K=1 and Δx = .1 a phase margin of 12.97 degrees was calculated.
Plotting the corresponding y(x) from the difference equation, Eq 1.8-88

$$y_{n+3} = \frac{2.83y_{n+2} - 2.667y_{n+1} + .837y_n + .005}{1.001}, \quad K=1 \text{ and } \Delta x = .1$$

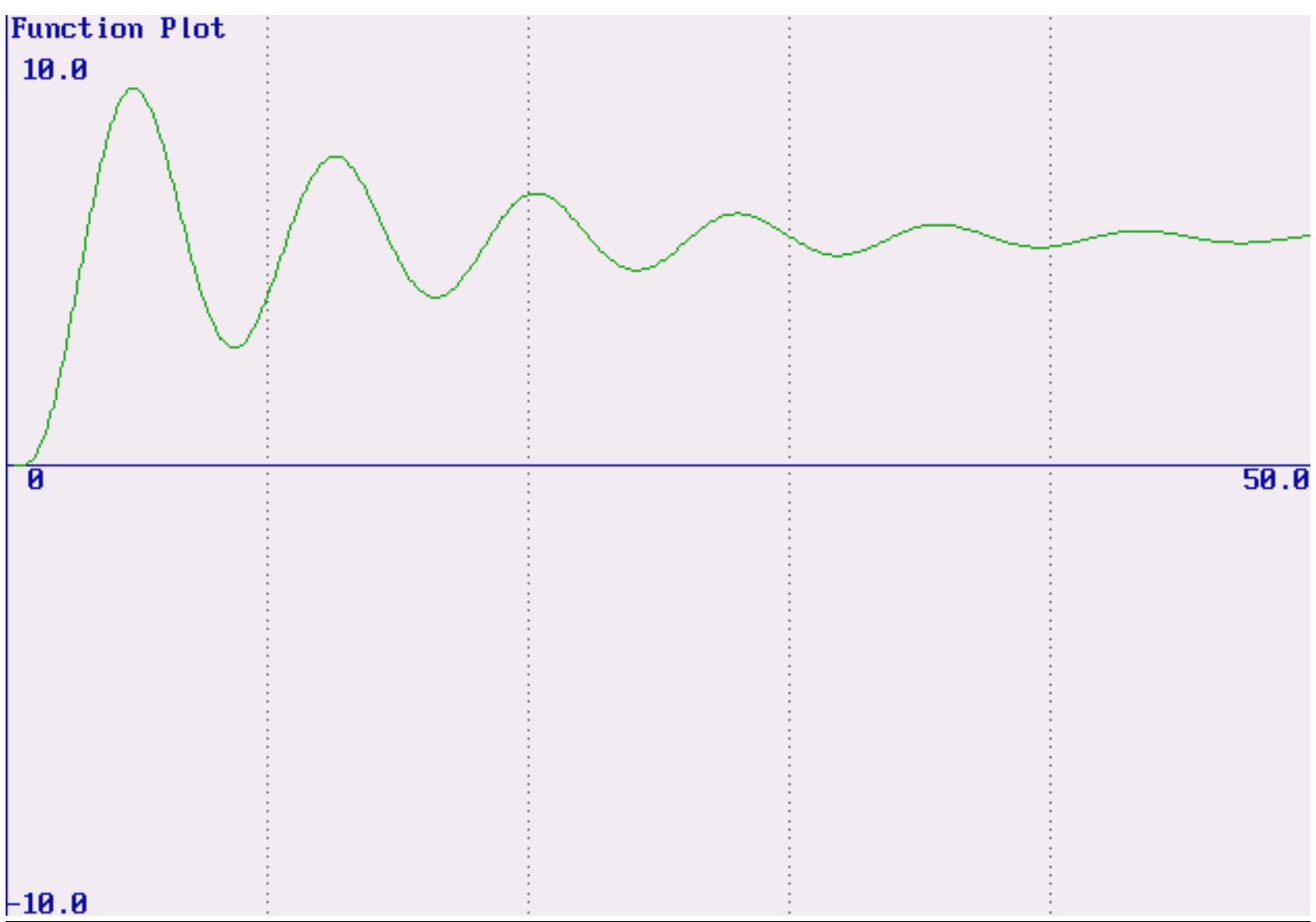

This y(x) plot displays stability but with a significant damped oscillation resulting from a step input of 5. This is what one would expect from a calculated phase margin of 12.97 degrees.

2) <u>For K=1.7, Δx =.1</u>

From Eq 1.8-80

$$A(s) = \frac{1.7(1+.1s)^3}{s(s+.7)(s+1)}$$

Evaluating A(s) around the left half s plane Critical Circle of radius 10 and center at s = –10

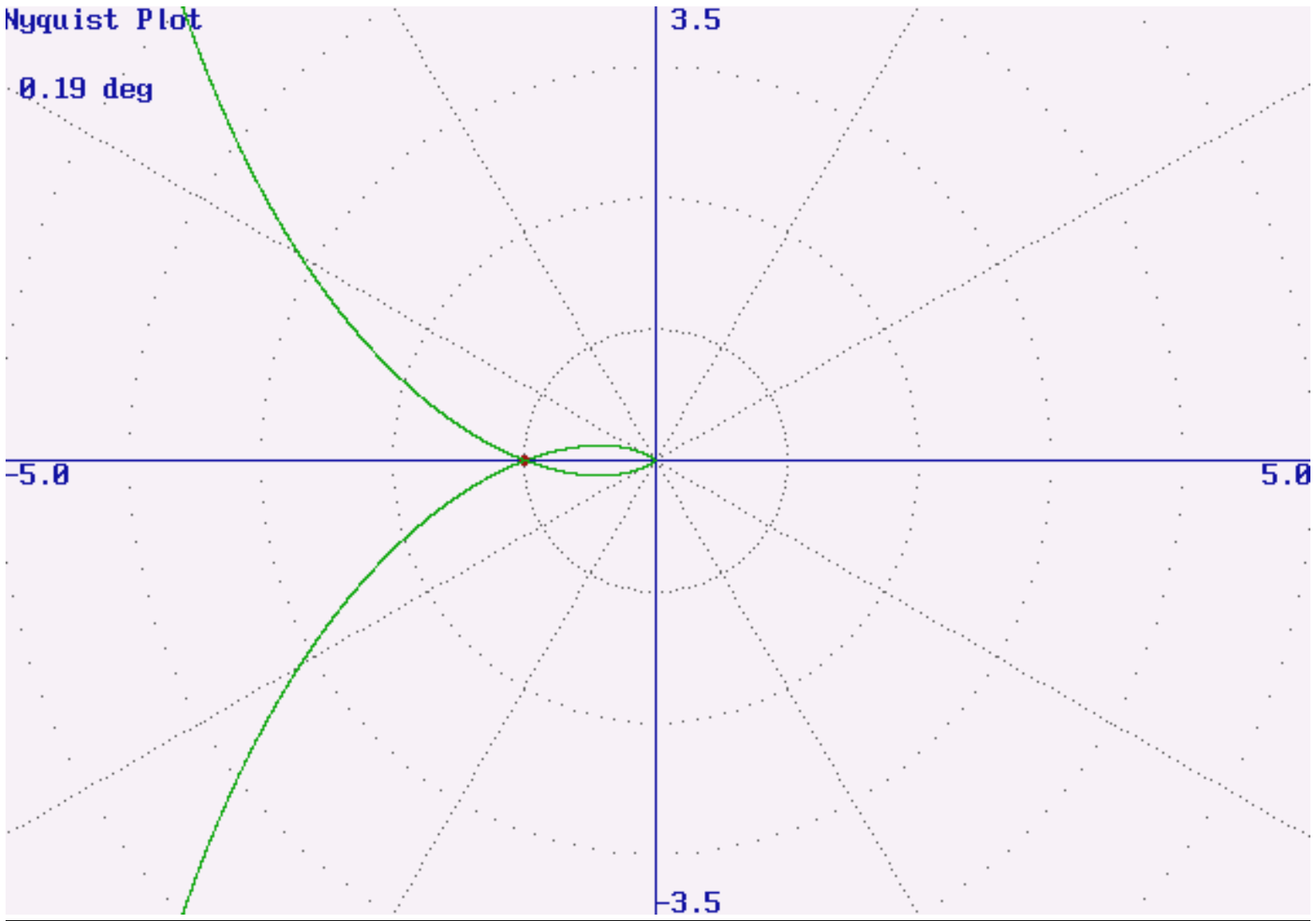

For K=1.7 and Δx = .1, the phase margin is calculated to be .19 degrees (virtually 0 degrees). This would indicate that the control system is oscillatory with a gain constant of K=1.7 and is on the verge of being unstable.

Plotting the corresponding y(x) from the difference equation, Eq 1.8-88

$$y_{n+3} = \frac{2.83y_{n+2} - 2.667y_{n+1} + .837y_n + .0085}{1.0017}, \qquad K = 1.7 \text{ and } \Delta x = .1$$

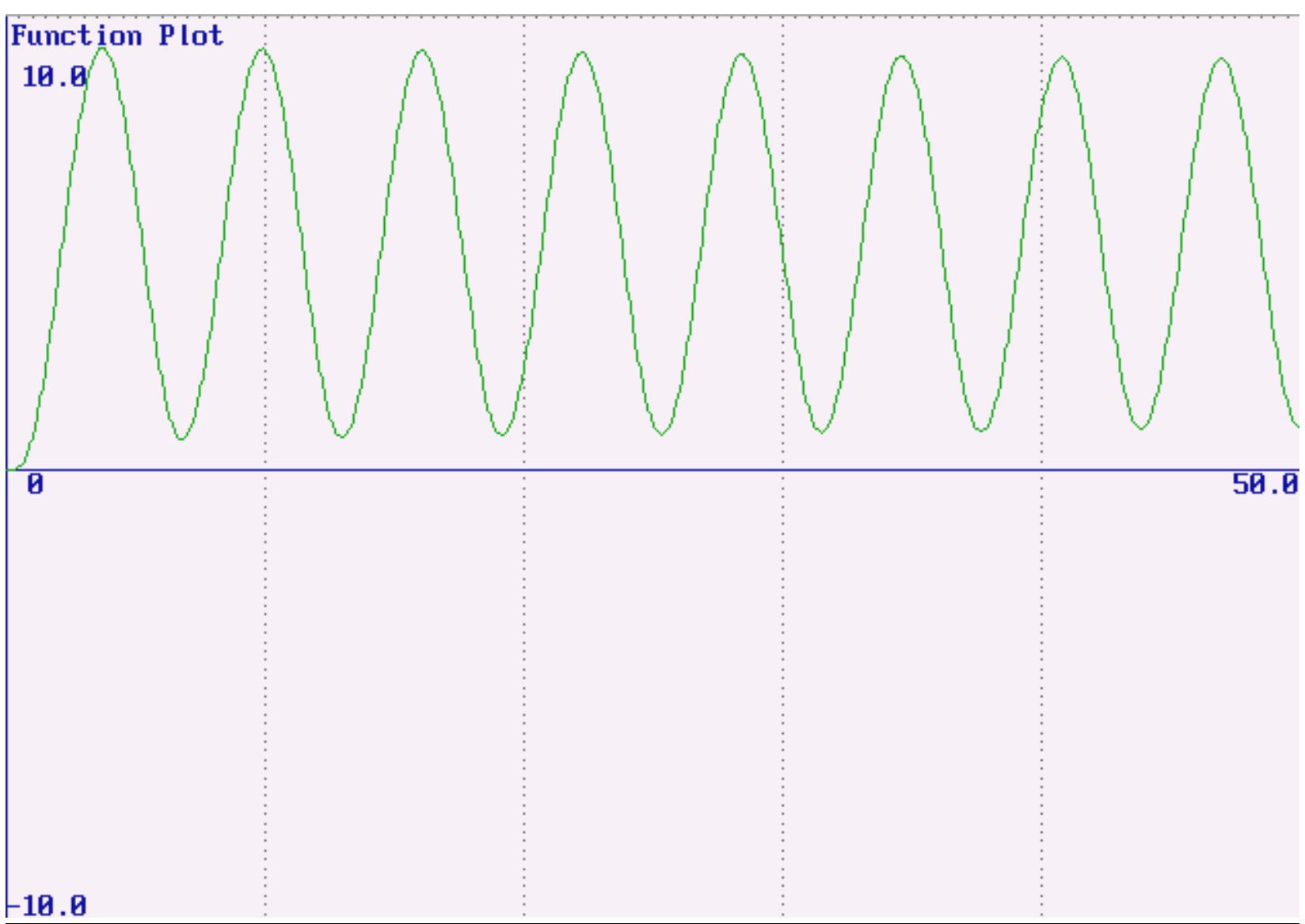

This y(x) plot displays a virtually sustained oscillation for a gain constant of K=1.7 resulting from a step input of 5. This is what one would expect from a calculated phase margin of .19 degrees.

# CHAPTER 1 SOLVED PROBLEMS

## Example 1.1

Find y(x) where $y(x+2\Delta x) - y(x) = 0$ , $\Delta x = 1$, $y(0) = 1$, $y(1) = -1$, $x = m\Delta x$ , m = integers. Solve for y(x) in two ways using the $K_{\Delta x}$ Transform.
From the table of $K_{\Delta x}$ Transform general equations (TABLE 2)

$$K_{\Delta x}[f(x+2\Delta x)] = (1+s\Delta x)^2 K_{\Delta x}[f(x)] - (1+s\Delta x)f(0)\Delta x - f(\Delta x)\Delta x \qquad 1)$$

$$y(x+2\Delta x) - y(x) = \quad 0 \qquad 2)$$

Taking the $K_{\Delta x}$ Transform of eq 2

$$K_{\Delta x}[y(x)+2] - K_{\Delta x}[y(x)] = 0 \qquad 3)$$

$$y(s) = K_{\Delta x}[y(x)] \qquad 4)$$

Substituting eq 1 and and eq 4 into eq 3 and introducing the y(x) initial conditions

$$(1+s)^2 y(s) - (1+s)(1)(1) - (-1)(1) – y(s) = 0 \qquad 5)$$

$$[(1+s)^2 – 1]y(s) – 1 – s + 1 = 0$$

$$(s^2+2s)y(s) = s$$

$$y(s) = \frac{1}{(s+2)} \qquad 6)$$

Solution 1

Finding the inverse $K_{\Delta x}$ Transform of eq 6

From the table of $K_{\Delta x}$ Transforms (TABLE 3)

$$K_{\Delta x}[\cos\frac{ax}{\Delta x}] = \frac{s + \frac{1-\cos a}{\Delta x}}{(s + \frac{1-\cos a}{\Delta x})^2 + (\frac{\sin a}{\Delta x})^2} \qquad 7)$$

For $\Delta x = 1$, $a = \pi$ , $x = 0, \Delta x, 2\Delta x, 3\Delta x, \dots$

Substituting into eq 7

$$K_{\Delta x}[\cos\pi x] = \frac{(s+2)}{(s+2)^2+0^2} = \frac{1}{(s+2)} \qquad 8)$$

$$y(s) = \frac{1}{s+2} = K_{\Delta x}[\cos\pi x] \qquad 9)$$

$$\mathbf{y(x) = \cos\pi x \ , \quad x = 0, 1, 2, 3, \dots} \qquad 10)$$

Checking eq 10 using the given difference equation and y(x) initial conditions

$y(0) = 1, \; y(1) = -1, \; \Delta x = 1$

$y(x+2) - y(x) = 0$

$\cos\pi(x+2) - \cos\pi x = \cos(\pi x+2\pi) - \cos\pi x = 0 \, , \; x = 0,1,2,3,\ldots$ Good check

$y(0) = \cos\pi(0) = 1 \, , \; y(1) = \cos\pi(1) = -1$ Good check

Solution 2

$$y(s) = \frac{1}{s+2}$$

From the tables of $K_{\Delta x}$ Transforms (TABLE 3) and Interval Calculus Functions (TABLE 4)

$$K_{\Delta x}[e_{\Delta x}(a,x)] = \frac{1}{s - a} \qquad 11)$$

$$e_{\Delta x}(a,x) = (1+a\Delta x)^{\frac{x}{\Delta x}} \qquad 12)$$

$a = -2, \; \Delta x = 1$

Substituting into Eq 11 and Eq 12

$$K_1[e_1(-2,x)] = \frac{1}{s + 2}$$

$$y(x) = e_1(-2,x) = (1-2(1))^{\frac{x}{1}} = (-1)^x$$

$$\mathbf{y(x) = (-1)^x \, , \; x = 0,1,2,3,\ldots} \qquad 13)$$

Checking Eq 13 using the given difference equation and y(x) initial conditions

$y(0) = 1, \; y(1) = -1, \; \Delta x = 1$

$y(x+2) - y(x) = 0$

$(-1)^{x+2} - (-1)^x = 0$

$(-1)^x(-1)^2 - (-1)^x = 0 \, , \; x = 0,1,2,3,\ldots$ Good check

$y(0) = (-1)^0 = 1 \, , \; y(1) = (-1)^1 = -1$ Good check

## Example 1.2

Find y(x) and z(x) from the following simultaneous difference equations and initial conditions using the $K_{\Delta x}$ Transform.

$y(x+\Delta x) + z(x) = 1$ , $y(0) = 1$ , $z(0) = 2$ , $\Delta x = .5$

$y(x) + z(x+\Delta x) = 0$ , $x = m\Delta x$ , $m = 0,1,2,3,\ldots$

From the table of $K_{\Delta x}$ Transform general equations (TABLE 2)

$$K_{\Delta x}[f(x+\Delta x)] = (1+s\Delta x)K_{\Delta x}[f(x)] - f(0)\Delta x \qquad 1)$$

$$y(x+\Delta x) + z(x) = 1 \qquad 2)$$

Taking the $K_{\Delta x}$ Transform of eq 2 with $\Delta x = .5$

$$K_{\Delta x}[y(x)+.5] + K_{\Delta x}[z(x)] = K_{\Delta x}[1] \qquad 3)$$

$$y(s) = K_{\Delta x}[y(x)]\ ,\ \ z(s) = K_{\Delta x}[z(x)] \qquad 4)$$

Substituting Eq 1 and Eq 4 into eq 3 with $\Delta x = .5$ and $y(0) = 1$

$$(1+.5s)y(s) – (1)(.5) +z(s) = K_{\Delta x}[1] \qquad 5)$$

From the table of $K_{\Delta x}$ Transforms (TABLE 3)

$$K_{\Delta x}[1] = \frac{1}{s} \qquad 6)$$

Substituting Eq 6 into Eq 5

$$(1+.5s)y(s) + z(s) = \frac{1}{s} + .5 \qquad 7)$$

Derived similarly from the equation $y(x) + z(x+\Delta x) = 0$ with $\Delta x = .5$ , $z(0) = 2$

$$y(s) + (1+.5s)z(s) = 1 \qquad 8)$$

Solving similataneous equations,Eq 7 and Eq 8

$$\begin{vmatrix} (1+.5s) & 1 \\ 1 & (1+.5s) \end{vmatrix} \begin{vmatrix} y(s) \\ z(s) \end{vmatrix} = \begin{vmatrix} \dfrac{1+.5s}{s} \\ 1 \end{vmatrix} \qquad 9)$$

$$[(1+.5s)^2 – 1]z(s) = 1+.5s - \frac{1+.5s}{s}$$

$$z(s) = \frac{(1+.5s)(s-1)}{s(.25s^2+s)}$$

$$z(s) = \frac{2(s+2)(s-1)}{s^2(s+4)} \qquad 10)$$

Finding the partial fraction expansion of Eq 10

$$z(s) = \frac{2(s+2)(s-1)}{s^2(s+4)} = \frac{A}{s^2} + \frac{B}{s} + \frac{C}{s+4} \qquad 11)$$

$$A = \frac{2(s+2)(s-1)}{s+4}\Big|_{s=0} = -1$$

$$C = \frac{2(s+2)(s-1)}{s^2}\Big|_{s=-4} = \frac{5}{4}$$

$$B = 2\frac{d}{ds}\frac{(s+2)(s-1)}{(s+4)}\Big|_{s=0} = \frac{3}{4}$$

Substituting into Eq 11

$$z(s) = -\frac{1}{s^2} + \frac{.75}{s} + \frac{1.25}{s+4} \qquad 12)$$

From the tables of $K_{\Delta x}$ Transforms (TABLE 3) and Discrete Functions (TABLE 4)

$$K_{\Delta x}[x] = \frac{1}{s^2},\ K_{\Delta x}[1] = \frac{1}{s},\ K_{\Delta x}[e_{\Delta x}(a,x)] = \frac{1}{s-a} \quad \{a = -4, \Delta x = .5\}$$

$$e_{\Delta x}(a,x) = (1+a\Delta x)^{\frac{x}{\Delta x}}$$

Finding the inverse $K_{\Delta x}$ Transform of z(s) in Eq 12

$$z(x) = -x +.75 + 1.25(1+(-4)(.5))^{\frac{x}{.5}}$$

$$\mathbf{z(x) = -x +.75 + 1.25(-1)^{2x}\ ,\quad x = .5m\ ,\ m = 0,1,2,3,\ldots} \qquad 13)$$

<u>Check 1</u>

$$z(0) = -0 + .75 + 1.25(-1)^0 = 2 \quad \text{Good check}$$

Finding y(s)

From Eq 8

$$y(s) = 1 - (1+.5s)z(s) \qquad 14)$$

Substituting Eq 10 into Eq 14

$$y(s) = 1 - .5(s+2)\left[\frac{2(s+2)(s-1)}{s^2(s^2+4)}\right] \qquad 15)$$

$$y(s) = \frac{s^2(s+4) - (s-1)(s+2)^2}{s^2(s+4)}$$

$$y(s) = \frac{s^3+4s^2 -s^3 -3s^2+4}{s^2(s+4)}$$

$$y(s) = \frac{s^2+4}{s^2(s+4)} \qquad 16)$$

Finding the partial fraction expansion of Eq 16

$$y(s) = \frac{s^2+4}{s^2(s+4)} = \frac{E}{s^2} + \frac{F}{s} + \frac{G}{s+4} \qquad 17)$$

$$E = \frac{s^2+4}{s+4}\Big|_{s=0} = 1$$

$$G = \frac{s^2+4}{s^2}\Big|_{s=-4} = \frac{5}{4}$$

$$F = \frac{d}{ds}\frac{s^2+4}{s+4}\Big|_{s=0} = -\frac{1}{4}$$

Substituting into Eq 17

$$y(s) = \frac{1}{s^2} - \frac{.25}{s} + \frac{1.25}{s+4} \qquad 18)$$

From the tables of $K_{\Delta x}$ Transforms (TABLE 3) and Discrete Functions (TABLE 4)

$$K_{\Delta x}[x] = \frac{1}{s^2}, \quad K_{\Delta x}[1] = \frac{1}{s}, \quad K_{\Delta x}[e_{\Delta x}(a,x)] = \frac{1}{s-a}) \quad \{a = -4, \Delta x = .5\}$$

$$e_{\Delta x}(a,x) = (1+a\Delta x)^{\frac{x}{\Delta x}}$$

Finding the inverse $K_{\Delta x}$ Transform of $y(s)$ in Eq 18

$$y(x) = x -.25 + 1.25(1+(-4)(.5))^{\frac{x}{.5}}$$

$$\mathbf{y(x) = x - .25 + 1.25(-1)^{2x}\ , \quad x = .5m\ , m = 0,1,2,3,\ldots} \qquad 19)$$

Check 2

$y(0) = 0 - .25 + 1.25(-1)^0 = 1$ Good check

Check 3

$z(x+.5) + y(x) = 0$ , $x = .5m$ , $m = 0,1,2,3,\ldots$

Substituting the values of z(x+.5) and y(x)

$-x - .5 +.75 + 1.25(-1)^{2x+1} + x - .25 + 1.25(-1)^{2x} = 0$ , $x = .5m$ , m = integers

$-.5 -.25 +.75 - 1.25 + 1.25 - x + x = 0$

$0 = 0$ Good check

Check 4

$y(x+.5) + z(x) = 1$ , $x = .5m$ , $m = 0,1,2,3,\ldots$

Substituting the values of y(x+.5) and z(x)

$x +.5 -.25 + 1.25(-1)^{2x+1} - x +.75 + 1.25(-1)^{2x} = 1$ , $x = .5m$ , $m = 0,1,2,3,\ldots$

$.5 -.25 +.75 - 1.25 + 1.25 + x - x = 1$

$1 = 1$ Good check

**Example 1.3** Find the Z Transform of f(x) = x using the interval calculus relationship,

$$Z[f(x)] = \frac{1}{T}\ {}_T\!\int_0^{\infty} z^{-(\frac{x}{\Delta x})} f(x)\,\Delta x\ .$$

$$Z[x] = \frac{1}{T}\ {}_T\!\int_0^{\infty} z^{-(\frac{x}{\Delta x})} x\,\Delta x \qquad 1)$$

Using integration by parts

$$_{\Delta x}\!\int_{x_1}^{x_2} v(x)D_{\Delta x}u(x)\Delta x = u(x)v(x)\Big|_{x1}^{x2} - {}_{\Delta x}\!\int_{x_1}^{x_2} D_{\Delta x}v(x)u(x+\Delta x)\Delta x \qquad 2)$$

From integration Table 6 in the Appendix

$$_{\Delta x}\!\int z^{-\frac{x}{\Delta x}}\Delta x = -\left(\frac{z\Delta x}{z-1}\right) z^{-\frac{x}{\Delta x}} + k \qquad 3)$$

Let $D_{\Delta x}u(x) = z^{-\frac{x}{\Delta x}}$ $\qquad v(x) = x$

$u(x) = -\left(\frac{z\Delta x}{z-1}\right) z^{-\frac{x}{\Delta x}}$ $\qquad D_{\Delta x}v(x) = 1$

$u(x+\Delta x) = -\left(\frac{\Delta x}{z-1}\right) z^{-\frac{x}{\Delta x}}$

$$Z[x] = \frac{1}{T}\ {}_T\!\int_0^{\infty} z^{-(\frac{x}{\Delta x})} x\,\Delta x = -\frac{1}{T}\, x\left(\frac{z\Delta x}{z-1}\right) z^{-\frac{x}{\Delta x}}\Bigg|_{0}^{\infty}{}_{T} + \frac{1}{T}\left(\frac{\Delta x}{z-1}\right)\ {}_T\!\int_0^{\infty} z^{-(\frac{x}{\Delta x})}\,\Delta x$$

$$Z[x] = 0 - \frac{1}{T}\left(\frac{\Delta x}{z-1}\right)\left(\frac{z\Delta x}{z-1}\right) z^{-\frac{x}{\Delta x}}\Bigg|_{0}^{\infty}{}_{T} = \frac{1}{T}\frac{T^2 z}{(z-1)^2}$$

$$\mathbf{Z[x] = \frac{Tz}{(z-1)^2}} \qquad 4)$$

**Example 1.4** Find the Z Transform, $Z[\sin\frac{ax}{\Delta x}]$ , given the $K_{\Delta x}$ Transform, $K_{\Delta x}[\sin\frac{ax}{\Delta x}]$ .

From transform Table 3 in the Appendix

$$K_{\Delta x}[\sin\frac{ax}{\Delta x}] = f(s) = \frac{\frac{\sin a}{\Delta x}}{(s + \frac{1-\cos a}{\Delta x})^2 + (\frac{\sin a}{\Delta x})^2} \qquad 1)$$

$$Z[f(x)] = F(z) = \frac{z}{T} f(s)\Big|_{\substack{s = \frac{z-1}{T} \\ T = \Delta x}} \qquad 2)$$

$$Z[\sin\frac{ax}{\Delta x}] = \frac{z}{T}\frac{\frac{\sin a}{T}}{(\frac{z-1}{T} + \frac{1-\cos a}{T})^2 + (\frac{\sin a}{T})^2}$$

$$Z[\sin\frac{ax}{\Delta x}] = \frac{z \sin a}{z^2-2z\cos a + \cos^2 a + \sin^2 a}$$

$$\mathbf{Z[\sin\frac{ax}{\Delta x}] = \frac{z \sin a}{z^2-2z\cos a +1}} \qquad 3)$$

**Example 1.5** Find the $K_{\Delta x}$ Transform, $K_{\Delta x}[x^3]$, given the Z Transform, $Z[x^3]$ .

$$Z[x^3] = F(z) = \frac{T^3 z(z^2+4z+1)}{(z-1)^4} \tag{1}$$

$$K_{\Delta x}[f(x)] = f(s) = \frac{\Delta x}{1+s\Delta x} F(z)\Big|_{\substack{z = 1+s\Delta x \\ T = \Delta x}} \tag{2}$$

$$K_{\Delta x}[x^3] = \frac{\Delta x}{1+s\Delta x} \; \frac{T^3 z(z^2+4z+1)}{(z-1)^4}\Big|_{\substack{z = 1+s\Delta x \\ T = \Delta x}}$$

$$K_{\Delta x}[x^3] = \frac{\Delta x}{1+s\Delta x} \; \frac{\Delta x^3(1+s\Delta x)(1+2s\Delta x+s^2\Delta x^2 +4+4s\Delta x+1)}{(1+s\Delta x-1)^4}$$

$$K_{\Delta x}[x^3] = \frac{6+6s\Delta x+s^2\Delta x^2}{s^4}$$

$$\mathbf{K_{\Delta x}[x^3] = \frac{6}{s^4} + \frac{6\Delta x}{s^3} + \frac{\Delta x^2}{s^2}} \tag{3}$$

**Example 1.6** For the $K_{\Delta x}$ Transform, $K_{\Delta x}[f(x)] = \frac{1}{s^2}$, find the values of f(x) for x = 0,.5,1,1.5.

$$K_{\Delta x}[f(x)] = \Delta x \sum_{n=0}^{\infty} f(n\Delta x)(1+s\Delta x)^{-n-1} \qquad 1)$$

$$x = n\Delta x \qquad 2)$$

$$\Delta x = .5 \qquad 3)$$

$$n = 0,1,2,3 \qquad 4)$$

$$\frac{1}{s^2} = \frac{\Delta x^2}{(s\Delta x)^2} = \frac{\Delta x^2}{(1+s\Delta x)^2-2s\Delta x-1} = \frac{\Delta x^2}{(1+s\Delta x)^2-2(1+s\Delta x)+1}$$

$$\frac{1}{s^2} = \Delta x\left[\frac{\Delta x}{(1+s\Delta x)^2-2(1+s\Delta x)+1}\right] \qquad 5)$$

Let

$$p = 1+s\Delta x \qquad 6)$$

Substituting Eq 6 into Eq 5 with $\Delta x = .5$

$$\frac{1}{s^2} = (.5)\left[\frac{.5}{p^2-2p+1}\right] \qquad 7)$$

Expanding Eq 7

$$\frac{1}{s^2} = (.5)(0p^{-1}+.5p^{-2}+1p^{-3}+1.5p^{-4}+\ldots) \qquad 8)$$

$$\frac{1}{s^2} = (.5)(0(1+s(.5))^{-1}+.5(1+s(.5))^{-2}+1(1+s(.5))^{-3}+1.5(1+s(.5))^{-4}+\ldots) \qquad 9)$$

Comparing Eq 9 to Eq 1

Then

f(0) = 0
f(.5) = .5
f(1) = 1
f(1.5) = 1.5

<u>Note</u> - $K_{\Delta x}[x] = \frac{1}{s^2}$

## <u>Example 1.7</u> A summation calculation when at least one value of x within the summation limits is at a pole

Evaluate the summation, $\sum_{x=-1.5}^{2}{}_{.5}\ \frac{x+.7}{(x+.5)(x-1)}$ . Exclude the poles at x = -.5,1 so that the sum attributable only to the finite values of $\frac{x+.7}{(x+.5)(x-1)}$ from x = -1.5 to x = 2 can be obtained.

Solve this example two ways.

<u>Solution 1</u>

Find the partial fraction expansion of the function, $f(x) = \frac{x+.7}{(x+.5)(x-1)}$

$$\frac{x+.7}{(x+.5)(x-1)} = \frac{A}{x+.5} + \frac{B}{x-1} \qquad 1)$$

Find A,B

$$Ax - A + Bx + .5B = x + .7 \qquad 2)$$

From Eq 2

$$A + B = 1 \qquad 3)$$

$$-A + .5B = .7 \qquad 4)$$

Add Eq 3 and Eq 4

$$B = \frac{1.7}{1.5} \qquad 5)$$

From Eq 3 and Eq 5

$$A = 1 - \frac{1.7}{1.5} = -\frac{.2}{1.5}$$

$$A = -\frac{.2}{1.5} \qquad 6)$$

From Eq 1, Eq 5, and Eq 6

$$\frac{x+.7}{(x+.5)(x-1)} = -\frac{.2}{1.5}\left[\frac{1}{x+.5}\right] + \frac{1.7}{1.5}\left[\frac{1}{x-1}\right] \qquad 7)$$

Finding the summation, $\displaystyle {}_{.5}\sum_{\substack{x=-1.5 \\ \text{any } \frac{1}{0} \text{ term excluded}}}^{2} \frac{x+.7}{(x+.5)(x-1)}$

Using the following equation to find the summation

$$\Delta x \; {}_{\Delta x}\sum_{\substack{x=x_1 \\ \text{any } \frac{1}{0} \text{ term excluded}}}^{x_2-\Delta x} f(x) = MV[{}_{\Delta x}\int_{x_1}^{x_2} f(x)\Delta x] - \{\, \Delta x \sum MV[f(x)] \;,\; x = x_1, x_2, x_3, \ldots, x_m \,\} \qquad 8)$$

$f(x) = \dfrac{x+.7}{(x+.5)(x-1)}$

$\Delta x = .5$

$x_1 = -.1.5$

$x_2 - \Delta x = 2$

Substituting into Eq 8

$$ {}_{.5}\sum_{\substack{x=-1.5 \\ \text{any } \frac{1}{0} \text{ term excluded}}}^{2} \frac{x+.7}{(x+.5)(x-1)} = \frac{1}{.5} MV[{}_{.5}\int_{-1.5}^{2.5} \frac{x+.7}{(x+.5)(x-1)} \Delta x] - \{ \sum MV[\frac{x+.7}{(x+.5)(x-1)}] \;,\; x = -.5, 1 \,\} \qquad 9)$$

Taking the integral of the partial fraction expansion of Eq 7

$$ {}_{.5}\int_{-1.5}^{2.5} \frac{x+.7}{(x+.5)(x-1)} \Delta x = -\frac{.2}{1.5}\, {}_{.5}\int_{-1.5}^{2.5} \frac{1}{x+.5} \Delta x \;+\; \frac{1.7}{1.5}\, {}_{.5}\int_{-1.5}^{2.5} \frac{1}{x-1} \Delta x \qquad 10)$$

Substituting Eq 10 into Eq 9

$$\sum_{x=-1.5}^{2}{}_{.5} \frac{x+.7}{(x+.5)(x-1)} = -\frac{.4}{1.5} MV[{}_{.5}\int_{-1.5}^{2.5} \frac{1}{x+.5}\Delta x] + \frac{3.4}{1.5} MV[{}_{.5}\int_{-1.5}^{2.5} \frac{1}{x-1}\Delta x]$$

any $\frac{1}{0}$ term excluded

$$-\{\sum MV[\frac{x+.7}{(x+.5)(x-1)}] , x = -.5, 1 \} \qquad 11)$$

Evaluate Eq 11

Using the equation:

$$MV[{}_{\Delta x}\int_{x_1}^{x_2} \frac{1}{x-a}\Delta x] = \text{lnd}(1,\Delta x,x-a) \Big|_{x_1}^{x_2} \qquad 12)$$

For $a = -.5$

$\Delta x = .5$

$x_1 = -.1.5$

$x_2 = 2.5$

$$MV[{}_{.5}\int_{-1.5}^{2.5} \frac{1}{x+.5}\Delta x] = \text{lnd}(1,.5,x+.5) \Big|_{-1.5}^{2.5} \qquad 13)$$

Calculating Eq 13 using the lnd(n,Δx,x) calculation program, LNDX

$$MV[{}_{.5}\int_{-1.5}^{2.5} \frac{1}{x+.5}\Delta x] = .7833333333 \qquad 14)$$

For $a = 1$

$\Delta x = .5$

$x_1 = -.1.5$

$x_2 = 2.5$

$$MV[_{.5}\int_{-1.5}^{2.5} \frac{1}{x-1}\,\Delta x] = \text{lnd}(1,.5,x-1)\Big|_{-1.5}^{2.5} \qquad 15)$$

Calculating Eq 15 using the lnd(n,Δx,x) calculation program, LNDX

$$MV[_{.5}\int_{-1.5}^{2.5} \frac{1}{x-1}\,\Delta x] = -.7833333333 \qquad 16)$$

Calculate $MV[\frac{x+.7}{(x+.5)(x-1)}]\,|_{x\,=\,-.5}$

$$MV[f(x)] = \lim_{\varepsilon\to 0}\left[\frac{f(x+\varepsilon)+f(x-\varepsilon)}{2}\right] \qquad 17)$$

$$MV[\frac{x+.7}{(x+.5)(x-1)}]\,|_{x\,=\,-.5} = \lim_{\varepsilon\to 0}\frac{1}{2}\left[\frac{.2+\varepsilon}{\varepsilon(\varepsilon-1.5)}+\frac{.2-\varepsilon}{-\varepsilon(-\varepsilon-1.5)}\right] = \lim_{\varepsilon\to 0}\frac{1}{2\varepsilon}\left[\frac{.2+\varepsilon}{(\varepsilon-1.5)}+\frac{.2-\varepsilon}{(\varepsilon+1.5)}\right]$$

$$MV[\frac{x+.7}{(x+.5)(x-1)}]\,|_{x\,=\,-.5} = \lim_{\varepsilon\to 0}\frac{1}{2\varepsilon}\left[\frac{(.2+\varepsilon)(\varepsilon+1.5)+(.2-\varepsilon)(\varepsilon-1.5)}{\varepsilon^2-(1.5)^2}\right] = \lim_{\varepsilon\to 0}\frac{1}{2\varepsilon}\left[\frac{3.4\varepsilon}{\varepsilon^2-(1.5)^2}\right]$$

$$MV[\frac{x+.7}{(x+.5)(x-1)}]\,|_{x\,=\,-.5} = -\frac{1.7}{(1.5)^2} \qquad 18)$$

Calculate $MV[\frac{x+.7}{(x+.5)(x-1)}]\,|_{x\,=\,1}$

$$MV[\frac{x+.7}{(x+.5)(x-1)}]\,|_{x\,=\,1} = \lim_{\varepsilon\to 0}\frac{1}{2}\left[\frac{1.7+\varepsilon}{\varepsilon(\varepsilon+1.5)}+\frac{1.7-\varepsilon}{-\varepsilon(-\varepsilon+1.5)}\right] = \lim_{\varepsilon\to 0}\frac{1}{2\varepsilon}\left[\frac{1.7+\varepsilon}{(\varepsilon+1.5)}+\frac{1.7-\varepsilon}{(\varepsilon-1.5)}\right]$$

$$MV[\frac{x+.7}{(x+.5)(x-1)}]\,|_{x\,=\,1} = \lim_{\varepsilon\to 0}\frac{1}{2\varepsilon}\left[\frac{(1.7+\varepsilon)(\varepsilon-1.5)+(1.7-\varepsilon)(\varepsilon+1.5)}{\varepsilon^2-(1.5)^2}\right] = \lim_{\varepsilon\to 0}\frac{1}{2\varepsilon}\left[\frac{.4\varepsilon}{\varepsilon^2-(1.5)^2}\right]$$

$$MV[\frac{x+.7}{(x+.5)(x-1)}]\,|_{x\,=\,1} = -\frac{.2}{(1.5)^2} \qquad 19)$$

From Eq 9, Eq 18 and Eq 19

$$\{\sum MV[\frac{x+.7}{(x+.5)(x-1)}]\,,\ x = -.5, 1\ \} = -\frac{1.7}{(1.5)^2}-\frac{.2}{(1.5)^2} = -\frac{1.9}{(1.5)^2}$$

$$\{ \sum MV[\frac{x+.7}{(x+.5)(x-1)}] , x = -.5, 1 \} = -\frac{1.9}{(1.5)^2} \qquad 20)$$

Substituting Eq 14, Eq 16, Eq 20 into Eq 11

$$\sum_{\substack{x=-1.5 \\ \text{any } \frac{1}{0} \text{ term} = 0}}^{2} {}_{.5} \frac{x+.7}{(x+.5)(x-1)} = -\frac{.4}{1.5}(.7833333333) + \frac{3.4}{1.5}(-.7833333333) - [-\frac{1.9}{(1.5)^2}]$$

$$\sum_{\substack{x=-1.5 \\ \text{any } \frac{1}{0} \text{ term excluded}}}^{2} {}_{.5} \frac{x+.7}{(x+.5)(x-1)} = -1.14 \qquad 21)$$

Checking the above result

$$\sum_{\substack{x=-1.5 \\ \text{any } \frac{1}{0} \text{ term excluded}}}^{2} {}_{.5} \frac{x+.7}{(x+.5)(x-1)} = \frac{-1.5+.7}{(-1.5+.5)(-1.5-1)} + \frac{-1+.7}{(-1+.5)(-1-1)} + 0 + \frac{0+.7}{(0+.5)(0-1)}$$

$$+ \frac{.5+.7}{(.5+.5)(.5-1)} + 0 + \frac{1.5+.7}{(1.5+.5)(1.5-1)} + \frac{2+.7}{(2+.5)(2-1)}$$

$$\sum_{\substack{x=-1.5 \\ \text{any } \frac{1}{0} \text{ term excluded}}}^{2} {}_{.5} \frac{x+.7}{(x+.5)(x-1)} = -.32 - .3 - 1.4 - 2.4 + 2.2 + 1.08 = -1.14 \qquad 22)$$

Good check

Solution 2

From Eq 1 thru Eq 7 in Solution 1 above, the partial fraction expansion of $\frac{x+.7}{(x+.5)(x-1)}$ is:

$$\frac{x+.7}{(x+.5)(x-1)} = -\frac{.2}{1.5}\left[\frac{1}{x+.5}\right] + \frac{1.7}{1.5}\left[\frac{1}{x-1}\right] \qquad 23)$$

Finding the summation,

$$\underset{\substack{x=-1.5 \\ \text{any } \frac{1}{0} \text{ term excluded}}}{\overset{2}{{}_{.5}\sum}} \frac{x+.7}{(x+.5)(x-1)} \qquad 24)$$

In the above summation, Eq 24, the terms, $\frac{x+.7}{(-.5+.5)(x-1)}$ and $\frac{x+.7}{(x+.5)(1-1)}$ are to be excluded.

From Eq 23 and Eq 24, the terms to be excluded are:

$$\frac{-.5+.7}{(-.5+.5)(-.5-1)} = \frac{.2}{(0)(-1.5)} = -\frac{.2}{1.5}\left[\frac{1}{-.5+.5}\right] + \frac{1.7}{1.5}\left[\frac{1}{-.5-1}\right] = -\frac{.2}{1.5}\left[\frac{1}{0}\right] - \frac{1.7}{(1.5)^2}$$

$$\frac{1+.7}{(1+.5)(1-1)} = \frac{1.7}{(1.5)(0)} = -\frac{.2}{1.5}\left[\frac{1}{1+.5}\right] + \frac{1.7}{1.5}\left[\frac{1}{1-1}\right] = -\frac{.2}{(1.5)^2} + \frac{1.7}{1.5}\left[\frac{1}{0}\right]$$

Note that when the partial fraction expansion is used, besides excluding the division by 0 term of the partial fraction expansion, another non division by zero term must also be excluded.

Using the partial fraction expansion

Substitute the partial fraction expansion of $\frac{x+.7}{(x+.5)(x-1)}$ into the summation of Eq 24

$$\underset{\substack{x=-1.5 \\ \text{any } \frac{1}{0} \text{ term excluded}}}{\overset{2}{{}_{.5}\sum}} \frac{x+.7}{(x+.5)(x-1)} = \underset{\substack{x=-1.5 \\ \text{any } \frac{1}{0} \text{ term excluded}}}{\overset{2}{{}_{.5}\sum}} \left(-\frac{.2}{1.5}\left[\frac{1}{x+.5}\right] + \frac{1.7}{1.5}\left[\frac{1}{x-1}\right]\right) - \underset{\text{non } \frac{1}{0} \text{ terms excluded}}{\left[\ -\frac{1.7}{(1.5)^2} - \frac{.2}{(1.5)^2}\ \right]}$$

$$\underset{\substack{x=-1.5 \\ \text{any } \frac{1}{0} \text{ term excluded}}}{\overset{2}{{}_{.5}\sum}} \frac{x+.7}{(x+.5)(x-1)} = -\frac{.2}{1.5} \underset{\substack{x=-1.5 \\ \text{any } \frac{1}{0} \text{ term excluded}}}{\overset{2}{{}_{.5}\sum}} \frac{1}{x+.5} + \frac{1.7}{1.5} \underset{\substack{x=-1.5 \\ \text{any } \frac{1}{0} \text{ term excluded}}}{\overset{2}{{}_{.5}\sum}} \frac{1}{x-1} + \frac{1.9}{(1.5)^2} \qquad 24)$$

Using the derived equation:

$$\sum_{\Delta x \; x=x_1}^{x_2} \frac{1}{(x-a)^n} \quad = \; \pm \frac{1}{\Delta x}\, \mathrm{lnd}(n,\Delta x,x\text{-}a)\, \Big|_{x_1}^{x_2+\Delta x} \quad , \quad + \text{ for } n = 1 \quad , \quad - \text{ for } n \neq 1 \qquad 25)$$

any $\frac{1}{0}$ term excluded

$n = 1$
$\Delta x = .5$
$x_1 = -1.5$
$x_2 - \Delta x = 2$

From Eq 24 and Eq 25

$$\sum_{.5 \; x=-1.5}^{2} \frac{x+.7}{(x+.5)(x-1)} \; = -\frac{.2}{1.5}\left(\frac{1}{.5}\right) \mathrm{lnd}(1,.5,x+.5)\Big|_{-1.5}^{2+.5} + \frac{1.7}{1.5}\left(\frac{1}{.5}\right) \mathrm{lnd}(1,.5,x-1)\Big|_{-1.5}^{2+.5} + \frac{1.9}{(1.5)^2} \qquad 26)$$

any $\frac{1}{0}$ term excluded

Note that the poles at x = -.5,1, the division by 0 terms, are excluded by Eq 25. The finite partial fraction terms at x = -.5,1 are excluded by the last term.

Calculating Eq 26 using the lnd(n,Δx,x) calculation program, LNDX

$$\sum_{.5 \; x=-1.5}^{2} \frac{x+.7}{(x+.5)(x-1)} \; = -\, .2666666666(.7833333333) + 2.266666666(-.7833333333) + .8444444444$$

any $\frac{1}{0}$ term excluded

$$\sum_{.5 \; x=-1.5}^{2} \frac{x+.7}{(x+.5)(x-1)} \; = \; -1.14 \qquad 27)$$

any $\frac{1}{0}$ term excluded

This is the same result obtained in Solution 1
Good check

## Example 1.8 Evaluation of the Hurwitz Zeta Function where a pole is encountered

Evaluate the Hurwitz Zeta function, $\zeta(1.2,-5)$

The Hurwitz Zeta Function is a generalization of the Riemann Zeta Function and is defined as follows:

$$\zeta(s,a) = \sum_{k=0}^{\infty} \frac{1}{(k+a)^s} \,, \quad \text{The definition of the Hurwitz Zeta Function} \qquad 1)$$

where

Any term with $k+a = 0$ is excluded

Changing the Hurwitz Zeta Function of Eq 1 to a more convenient equivalent form

$$\zeta(s,a) = {}_1\sum_{x=a}^{\infty} \frac{1}{x^s} \qquad 2)$$

where

Any term with $x = 0$ is excluded

Use Eq 3 and Eq 4 to evaluate $\zeta(n,x)$

$$\zeta(n,x) = \zeta(n,1,x) \qquad 3)$$

$$\zeta(n,\Delta x,x_i) = \frac{1}{\Delta x}\,\mathrm{lnd}(n,\Delta x,x_i) = {}_{\Delta x}\sum_{x=x_i}^{\pm\infty} \frac{1}{x^n} \,, \quad Re(n) > 1 \,, \quad \text{A form of the General Zeta Function} \qquad 4)$$

where

$+\infty$ for $Re(\Delta x)>0$ or $\{Re(\Delta x)=0$ and $Im(\Delta x)>0\}$

$-\infty$ for $Re(\Delta x)<0$ or $\{Re(\Delta x)=0$ and $Im(\Delta x)<0\}$

$x = x_i,\ x_i+\Delta x,\ x_i+2\Delta x,\ x_i+3\Delta x,\ \ldots$

$n,x,x_i,\Delta x$ = real or complex values

$\Delta x$ = x increment

Any term with $x = 0$ is excluded

Comment – The $\mathrm{lnd}(n,\Delta x,x)$ function excludes any summation term where $x = 0$.

For the Hurwitz Zeta Function

Let

$\Delta x = 1$

$n = s$

$x_i = a$

Substituting into Eq 3 and Eq 4

$$\zeta(s,a) = \zeta(s,1,a) = \sum_{x=a}^{\infty}{}_1 \frac{1}{x^s} = \text{lnd}(s,1,a) \qquad 5)$$

where

x = a, a+1, a+2, a+3, …

Re(s) > 1

Any term with x = 0 is excluded

Use Eq 5 to evaluate the Hurwitz Zeta function, ζ(1.2,-5)

$$\zeta(1.2,-5) = \sum_{x=-5}^{\infty}{}_1 \frac{1}{x^{1.2}} = \text{lnd}(1.2,1,-5) \qquad 6)$$

where

x = -5, -4, -3, -2, -1, 0, 1, 2, 3, 4, 5, 6 …

Any term with x = 0 is excluded

Using the lnd(n,Δx,x) calculation program to calculate Eq 6

$$\zeta(1.2,-5) = 3.9433912875609658 + 1.1974809673848110i \qquad 7)$$

Checking the above result

$$\zeta(1.2,-5) = \sum_{x=-5}^{\infty}{}_1 \frac{1}{x^{1.2}} = \sum_{x=-5}^{0}{}_1 \frac{1}{x^{1.2}} + \sum_{x=1}^{\infty}{}_1 \frac{1}{x^{1.2}} = \sum_{x=-5}^{0}{}_1 \frac{1}{x^{1.2}} + \zeta(1.2) \qquad 8)$$

where

ζ(1.2) = Riemann Zeta Function of 1.2

Find the Riemann Zeta Function of 1.2

Again using the lnd(n,Δx,x) calculation program, LNDX

$$\zeta(1.2) = \sum_{x=1}^{\infty}{}_1 \frac{1}{x^{1.2}} = \text{lnd}(1.2,1,1) = 5.5915824411777507 \qquad 9)$$

From Eq 8 and Eq 9

$$\zeta(1.2,-5) = \left[\frac{1}{(-5)^{1.2}} + \frac{1}{(-4)^{1.2}} + \frac{1}{(-3)^{1.2}} + \frac{1}{(-2)^{1.2}} + \frac{1}{(-1)^{1.2}} + 0\right] + 5.5915824411777507 \qquad 9)$$

Note that the pole at x = 0 is excluded

Using a computer to calculate the above summation

$$\zeta(1.2,-5) = -1.6481911536167848+1.1974809673848110i + 5.5915824411777507$$

$$\zeta(1.2,-5) = 3.9433912875609658 + 1.1974809673848110i \qquad 10)$$

Good check

# CHAPTER 2

## The lnd(n,Δx,x) Function

# Chapter 2

# The lnd(n,Δx,x) function

## Section 2.1: Description of the lnd(n,Δx,x) function

As previously mentioned, Interval Calculus is closely related to Calculus. Interval Calculus is Calculus where the x increment, Δx, is not infinitesimal. The Interval Calculus functions, in general, include the value of Δx as a variable in addition to the x variable. In fact, if Δx approaches zero as a limit (i.e. $\Delta x \to 0$), Interval Calculus and Calculus become one in the same.

Calculus did run into a particularly difficult problem early in its development, the integral of the function, $f(x) = \frac{1}{x}$ . The integral of this function turned out to be a rather special and complex function, the natural logarithm, ln(x). Interval Calculus, as it turns out, is not immune to this same Calculus difficulty. The discrete integral of this same function, $f(x) = \frac{1}{x} = \frac{1}{x^1}$ , is also a rather special and complex function, +lnd(1,Δx,x) (or alternately designated $\ln_{\Delta x} x$). Fortunately for Calculus, the integral of the function, $f(x) = \frac{1}{x^n}$ , did not present any difficulty provided the value of n was not equal to 1. A simple relationship existed for this integral. $\int \frac{1}{x^n} dx = \frac{x^{n+1}}{n+1} + k$ where k is the constant of integration. No function of such simplicity exists for the discrete integral of this function. The discrete integral of $f(x) = \frac{1}{x^n}$ where n is not equal to1 is a rather special function, -lnd(n,Δx,x). For clarity, writing in mathematical notation the two functions:

$$\int_{\Delta x} \frac{1}{x} \Delta x = +\text{lnd}(1,\Delta x,x) + k \quad \text{or} \quad \int_{\Delta x\, x_1}^{x_2} \frac{1}{x} \Delta x = +\text{lnd}(1,\Delta x,x) \Big|_{\Delta x\, x_1}^{x_2} , \quad n = 1 \qquad (2.1\text{-}1)$$

where k = constant of integration

and

$$\int_{\Delta x} \frac{1}{x^n} \Delta x = -\text{lnd}(n,\Delta x,x) + k \quad \text{or} \quad \int_{\Delta x\, x_1}^{x_2} \frac{1}{x^n} \Delta x = -\text{lnd}(n,\Delta x,x) \Big|_{\Delta x\, x_1}^{x_2} , \quad n \neq 1 \qquad (2.1\text{-}2)$$

where k = constant of integration

Merging the two above equations, 2.1-1 and 2.1-2

$$\int_{\Delta x} \frac{1}{x^n} \Delta x = \pm\, \text{lnd}(n,\Delta x,x) + k \quad \text{or} \quad \int_{\Delta x\, x_1}^{x_2} \frac{1}{x^n} \Delta x = \pm\, \text{lnd}(n,\Delta x,x) \Big|_{\Delta x\, x_1}^{x_2} , \quad + \text{ for } n = 1 , \quad - \text{ for } n \neq 1 \qquad (2.1\text{-}3)$$

where k = constant of integration

Equation 2.1-3 shows the equality of the Interval Calculus integral of $\frac{1}{x^n}$ to the function, lnd(n,Δx,x). Interestingly, the lnd(n,Δx,x) function has some additional significance in higher mathematics.

$$\Delta x \sum_{\Delta x \; x=x_1}^{x_2-\Delta x} \frac{1}{x^n} = {}_{\Delta x}\!\int_{x_1}^{x_2} \frac{1}{x^n}\,\Delta x \qquad (2.1\text{-}4)$$

From Eq 2.1-3 and Eq 2.1-4

$$\Delta x \sum_{\Delta x \; x=x_1}^{x_2-\Delta x} \frac{1}{x^n} = {}_{\Delta x}\!\int_{x_1}^{x_2} \frac{1}{x^n}\,\Delta x = \pm\, \mathrm{lnd}(n,\Delta x,x) \,\Big|_{\Delta x \; x_1}^{x_2} \;,\; + \text{ for } n = 1 \;,\; - \text{ for } n \neq 1$$

$$\sum_{\Delta x \; x=x_1}^{x_2-\Delta x} \frac{1}{x^n} = \pm \frac{1}{\Delta x}\, \mathrm{lnd}(n,\Delta x,x) \,\Big|_{\Delta x \; x_1}^{x_2} \;,\; + \text{ for } n = 1 \;,\; - \text{ for } n \neq 1 \qquad (2.1\text{-}5)$$

The summation of Eq 2.1-5 is recognized as being the Riemann Zeta Function for $\Delta x = 1$, $x_1 = 1$, $x_2 \rightarrow \infty$ and Re(n)>1

$$\sum_{1 \; x=1}^{\infty} \frac{1}{x^n} = -\mathrm{lnd}(n,\Delta x,x) \,\Big|\Big|_{1}^{\infty} = \mathrm{lnd}(n,1,1) \;,\quad n \neq 1 \qquad (2.1\text{-}6)$$

where
Re(n)>1

<u>Note</u> – For Re(n)>1, lnd(n,Δx,∞) = 0

$$\sum_{1 \; x=1}^{\infty} \frac{1}{x^n} = \mathrm{lnd}(n,1,1) = \zeta(n) \;,\quad n \neq 1 \qquad \text{The Riemann Zeta Function} \qquad (2.1\text{-}7)$$

where
Re(n)>1

The summation of eq 2.1-5 is recognized as being the Hurwitz Zeta Function for $\Delta x = 1$, $x_1 = x_i$, $x_2 \to \infty$ and Re(n)>1

$$\sum_{x=x_i}^{\infty}{}_{1} \frac{1}{x^n} = -\mathrm{lnd}(n,\Delta x,x)\Big|_{x_i}^{\infty} = \mathrm{lnd}(n,1,x_i)\ , \quad n \neq 1 \tag{2.1-8}$$

where

Re(n)>1

Note – For Re(n)>1, $\mathrm{lnd}(n,\Delta x,\infty) = 0$

$$\sum_{x=x_i}^{\infty}{}_{1} \frac{1}{x^n} = \mathrm{lnd}(n,1,x_i) = \zeta(n,x_i)\ , \quad n \neq 1 \qquad \text{The Hurwitz Zeta Function} \tag{2.1-9}$$

where

Re(n)>1

Note - In the Hurwitz Zeta Function summation, a division by zero term is given a zero value. This condition is built into the lnd(n,Δx,x) function.

Note - For $0 < \mathrm{Re}(n) \leq 1$, the summation, $\displaystyle\sum_{x=x_i}^{\infty}{}_{\Delta x} \frac{1}{x^n}$, may at first glance be considered to be useless since it is infinite. However, the summation, Eq 2.1-10, shown below is finite and does have mathematical value.

$$\sum_{x=x_i}^{\infty}{}_{2\Delta x} \frac{1}{x^n} - \sum_{x=x_i+\Delta x}^{\infty}{}_{2\Delta x} \frac{1}{x^n} = \sum_{x=x_i}^{\infty}{}_{\Delta x} (-1)^{\frac{x-x_i}{\Delta x}} \frac{1}{x^n} = \pm \frac{1}{2\Delta x}\left[-\mathrm{lnd}(n,2\Delta x,x_i) + \mathrm{lnd}(n,2\Delta x,x_i+\Delta x)\right] \tag{2.1-10}$$

where

+ is for n = 1

– is for n ≠ 1,

Eq 2.1-10 is derived from Eq 2.1-5.

Note – $\lim_{N\to\infty}[-\mathrm{lnd}(n,2\Delta x,N) + \mathrm{lnd}(n,\Delta x,N+\Delta x)] = \lim_{N\to\infty} \displaystyle\sum_{x=N}^{\infty}{}_{2\Delta x} \left[\frac{1}{x^n} - \frac{1}{(x+\Delta x)^n}\right] = 0$

Considering the fact that the following summation is a generalization of both the Riemann and Hurwitz Zeta Functions, in this paper it will be referred to as the General Zeta Function.

From Eq 2.1-5

$$\sum_{\Delta x \; x=x_1}^{x_2-\Delta x} \frac{1}{x^n} = \pm \frac{1}{\Delta x}\,\text{lnd}(n,\Delta x,x)\Big|_{x_1}^{x_2} = \pm\,\zeta(n,\Delta x,x)\Big|_{x_1}^{x_2}, \quad \text{The General Zeta Function} \qquad (2.1\text{-}11)$$

$$\zeta(n,\Delta x,x) = \frac{1}{\Delta x}\,\text{lnd}(n,\Delta x,x) \qquad (2.1\text{-}12)$$

where

+ for n=1 , – for $n\neq 1$

$x = x_1, x_1+\Delta x, x_1+2\Delta x, x_1+3\Delta x, \ldots , x_2-\Delta x, x_2$

$\Delta x$ = x increment

$n,\Delta x,x,x_1,x_2$ = real or complex values

$x_1,x_2$ may be infinite values

<u>Note</u> - n may be any real or complex value including 1, a special condition.

From Eq 2.1-11, the following very useful form of the General Zeta Function is derived:

For Re(n)>1, $\text{lnd}(n,\Delta x,\pm\infty) = 0$

Let

$x_1 = x_i$

$x_2 \to +\infty$ or $-\infty$

Re(n)>1

$$\sum_{\Delta x \; x=x_i}^{\pm\infty} \frac{1}{x^n} = -\frac{1}{\Delta x}\,\text{lnd}(n,\Delta x,x)\Big|_{x_i}^{\pm\infty} = -\,\zeta(n,\Delta x,x)\Big|_{x_i}^{\pm\infty} \qquad (2.1\text{-}13)$$

Then

A useful form of the General Zeta Function where Re(n)>1 is:

$$\zeta(n,\Delta x,x_i) = \sum_{\Delta x \; x=x_i}^{\pm\infty} \frac{1}{x^n} = \frac{1}{\Delta x}\,\text{lnd}(n,\Delta x,x_i)\,, \quad \text{Re}(n)>1 \qquad (2.1\text{-}14)$$

where

$x = x_i, x_i+\Delta x, x_i+2\Delta x, x_i+3\Delta x, \ldots$

$\Delta x$ = x increment

$n,\Delta x,x,x_i$ = real or complex values

$+\infty$ for Re(Δx)>0 or {Re(Δx)=0 and Im(Δx)>0}
$-\infty$ for Re(Δx)<0 or {Re(Δx)=0 and Im(Δx)<0}

Note – For more concerning the derivation of the General Zeta Function see Chapter 8.

## Section 2.2: Commentary and description of the two series used to calculate the lnd(n,Δx,x) function

Finding a function equivalent, lnd(n,Δx,x), for the discrete integral, $\int_{\Delta x} \frac{1}{x^n}\Delta x$, was a much harder task than initially anticipated. As it turned out, two series were required to calculate the integral. The first series, the series to calculate lnd(1,Δx,x) (also designated $\ln_{\Delta x}x$) where n=1 has two variables, Δx and x. The second series, the series to calculate lnd(n,Δx,x) n≠1 has three variables, n,Δx,x. Both series evaluate their respective lnd(n,Δx,x) function for either real or complex values of their variables. These two series were not derived from the Euler-Maclauren Sum Formula. The development of these two series will be shown in Section 2.10 and Section 2.11. For n=1 or n≠1 the appropriate series is selected to calculate the function, lnd(n,Δx,x). Both series do not calculate exact values for the lnd(n,Δx,x) function. However, both series rapidly converge to the exact value as the absolute value, $\left|\frac{x}{\Delta x}\right|$, increases in magnitude.

Note - The quantity, $\frac{x}{\Delta x}$, in the magnitude just referred to, is required to have a significant real component.

For n = 1

$$\int_{\Delta x} \frac{1}{x}\Delta x = \text{lnd}(1,\Delta x,x) + k \equiv \ln_{\Delta x}x + k \tag{2.2-1}$$

where
$x = x_i + m\Delta x$ , m = integers
$x_i$ = a value of x
Δx = x increment
x, Δx = real or complex values
k = real or complex constant of integration

Comment – It was early recognized that the function evaluation of the discrete integral, $\int_{\Delta x} \frac{1}{x}\Delta x$, had similarities to the natural logarithm, lnx. For this reason, the function was named, $\ln_{\Delta x}x$. Considerably later, this function was found to be a special case of the function, lnd(n,Δx,x) n≠1. It was then alternately designated as lnd(1,Δx,x).

For n ≠1

$$\int_{\Delta x} \frac{1}{x^n}\,\Delta x = -\mathrm{lnd}(n,\Delta x,x) + k\ , \quad k = \text{contant of integration} \tag{2.2-2}$$

where

$x = x_i + m\Delta x$ , m = integers

$x_i$ = a value of x

$\Delta x$ = x increment

x, $\Delta x$ = real or complex values

n = real or complex constant

k = real or complex constant of integration

The two series

Finding the anti-derivative of $D_{\Delta x}\mathrm{lnd}(1,\Delta x,x) = +\frac{1}{x}$ , n = 1

$$\int_{\Delta x} \frac{1}{x}\,\Delta x = \mathrm{lnd}(1,\Delta x,x) \approx \ln\!\left(\frac{x}{\Delta x} - \frac{1}{2}\right) + \sum_{m=1}^{\infty} \frac{(2m-1)!\, C_m}{(2m+1)!\, 2^{2m} \left(\frac{x}{\Delta x} - \frac{1}{2}\right)^{2m}} + \gamma \tag{2.2-3}$$

The accuracy of Eq 2.2-3 above increases rapidly for increasing $|\frac{x}{\Delta x}|$ .

Note that the constant of integration, γ, has only one value, that of Euler’s Constant (.577215664…).

The derivation of this series is shown in Section 2.10 of this chapter. Its programming was not difficult.

Finding the anti-derivative of $D_{\Delta x}\mathrm{lnd}(n,\Delta x,x) = -\frac{1}{x^n}$ , n ≠ 1

$$\int_{\Delta x} \frac{1}{x^n}\,\Delta x = -\mathrm{lnd}(n,\Delta x,x) + K \approx -\sum_{m=0}^{\infty} \frac{\Gamma(n+2m-1)\left(\frac{\Delta x}{2}\right)^{2m} C_m}{\Gamma(n)(2m+1)!\left(x - \frac{\Delta x}{2}\right)^{n+2m-1}} + K\ , \quad n\neq 1 \tag{2.2-4}$$

where

The accuracy of Eq 2.2-4 above increases rapidly for increasing $|\frac{x}{\Delta x}|$ .

The series function of Eq 2.3-4 was inexplicably found to change to an incorrect value in certain integration applications. However, this value change could accurately be compensated for by an appropriate change in the constant of integration, *K*. There is considerable discussion concerning this series characteristic in the following sections of this chapter.

The derivation of this series is shown in Section 2.11 of this chapter. Its programming was considerably complicated by the uncommon and perplexing series characteristic specified above. The programming of this series was difficult.

To provide some insight into the uncommon characteristics and the programming of the $\text{lnd}(n,\Delta x,x)$ $n \neq 1$ function of Eq 2.2-4, the following presentation is made.

<u>Description and definition of the lnd(n,Δx,x) n ≠ 1 function</u>

$$\int_{\Delta x\; x_1}^{x_2} \frac{1}{x^n}\,\Delta x = \Delta x \sum_{\Delta x\; x=x_1}^{x_2-\Delta x} \frac{1}{x^n} \tag{2.2-5}$$

From Eq 2.2-2

$$\int_{\Delta x\; x_1}^{x_2} \frac{1}{x^n}\,\Delta x = -\text{lnd}(n,\Delta x,x)\Big|_{x_1}^{x_2}, \quad n \neq 1 \tag{2.2-6}$$

From Eq 2.2-5 and Eq 2.2-6

$$-\text{lnd}(n,\Delta x,x)\Big|_{x_1}^{x_2} = \Delta x \sum_{\Delta x\; x=x_1}^{x_2-\Delta x} \frac{1}{x^n}, \quad n \neq 1 \tag{2.2-7}$$

From Eq 2.2-7

$$\boldsymbol{-\text{lnd}(n,\Delta x,x_2) + \text{lnd}(n,\Delta x,x_1) = \Delta x \sum_{\Delta x\; x=x_1}^{x_2-\Delta x} \frac{1}{x^n}} \tag{2.2-8}$$

**where**

**$x = x_1, x_1+\Delta x, x_1+2\Delta x, x_1+3\Delta x, \ldots, x_2-\Delta x, x_2$**

**$n \neq 1$**

**$\Delta x$ = x increment**

Substituting into Eq 2.2-7

$x_1 = x_i$
$x_2 \rightarrow \pm\infty$
+ for Re(Δx)>0 or {Re(Δx)=0 and Im(Δx)>0}
– for Re(Δx)<0 or {Re(Δx)=0 and Im(Δx)<0}
and
Let Re(n) > 1

$$-\mathrm{lnd}(n,\Delta x,x)\Big|_{x_i}^{\pm\infty} = \Delta x \sum_{\Delta x\; x=x_i}^{\pm\infty} \frac{1}{x^n}, \quad n \neq 1 \tag{2.2-9}$$

From Eq 2.2-9

$$-\mathrm{lnd}(n,\Delta x,\pm\infty) + \mathrm{lnd}(n,\Delta x,x_i) = \Delta x \sum_{\Delta x\; x=x_i}^{\pm\infty} \frac{1}{x^n}, \quad n \neq 1 \tag{2.2-10}$$

Let x→±∞ in Eq 2.2-4

$$\mathrm{lnd}(n,\Delta x,\pm\infty) = 0 \tag{2.2-11}$$

From Eq 2.2-10, Eq 2.2-11 and Eq 2.2-5

$$\mathbf{lnd(n,\Delta x,x_i) = \Delta x \sum_{\Delta x\; x=x_i}^{\pm\infty} \frac{1}{x^n} = {}_{\Delta x}\!\int_{x_i}^{\pm\infty} \frac{1}{x^n}\, \Delta x} \tag{2.2-12}$$

**where**

**$x = x_i, x_i+\Delta x, x_i+2\Delta x, x_i+3\Delta x + x_i+4\Delta x, x_i+5\Delta x, \ldots$**

**Re(n) > 1**

**Δx = x increment**

## The x locus and the x locus line

A summation is performed between two limits

$$\Delta x \sum_{\substack{x=x_1 \\ \Delta x}}^{x_2} \frac{1}{x^n} = \pm\, \mathrm{lnd}(n,\Delta x,x)\,\Big|_{x_1}^{x_2+\Delta x} = \pm[\mathrm{lnd}(n,\Delta x,x_2+\Delta x) - \mathrm{lnd}(n,\Delta x,x_1)] \qquad (2.2\text{-}13)$$

**where**

**+ for n = 1**

**– for n ≠ 1**

**$x = x_1, x_1+\Delta x, x_1+2\Delta x, x_1+3\Delta x + x_1+4\Delta x, x_1+5\Delta x, \ldots, x_2-\Delta x, x_2$**

**<u>Note</u> – For Eq 2.2-13, Re(n) can be any real value and n can equal 1**

There is a commonality between all of the summation x values. When plotted in the complex plane, the x values form a line of equally spaced points. See the following diagram, Diagram 2.2-1.

<u>Diagram 2.2-1: An x locus and x locus line in the complex plane</u>

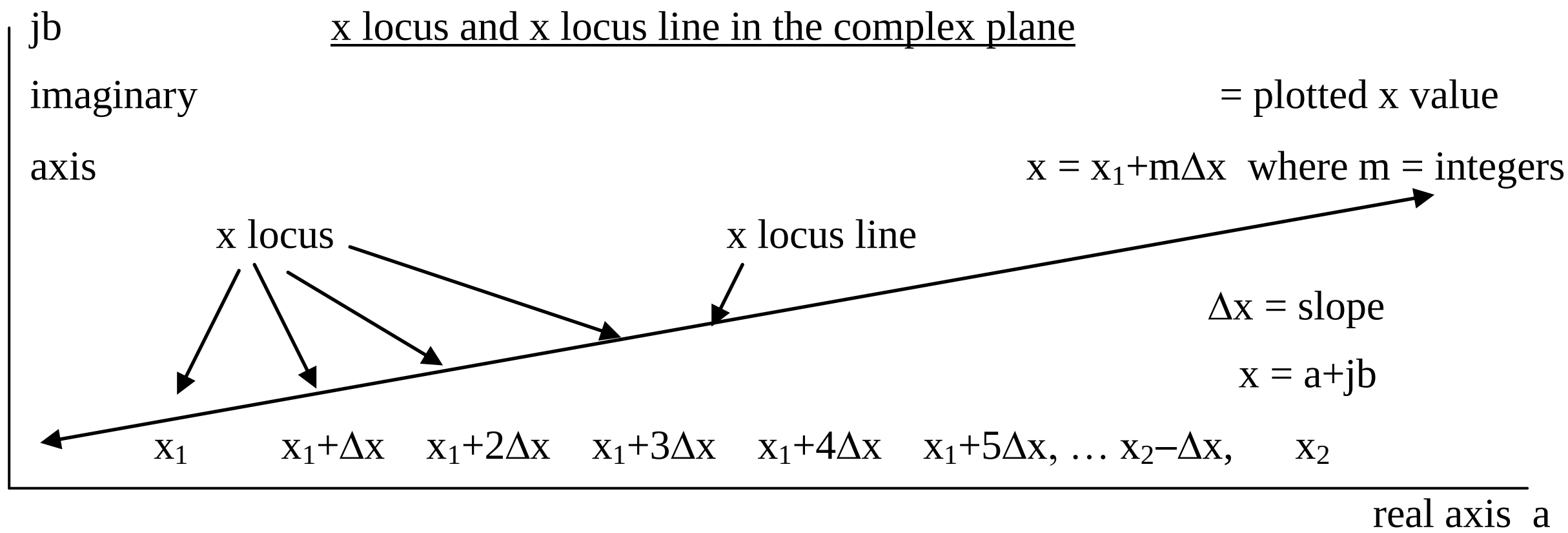

The x locus shown in Diagram 2.2-1 is formed by plotting all the values x in the complex plane ($x = x_1 + m\Delta x$ where m = integers). The x locus line is obtained by drawing a line through any two summation values of x. If an extended summation or a partial summation is desired, the x values of these summations will be in the same x locus and on the same x locus line as the original summation. The x locus line, just described, not only is indicative of the x values of a summation but also provides some insight into a perplexing characteristic involving the derived lnd(n,Δx,x) n ≠ 1 equation, Eq 2.2-4. This odd and uncommon characteristic greatly complicated the programming of Eq 2.2-4.

An example of an lnd(n,Δx,x) n≠1 function x locus line in the complex plane is shown in Diagram 2.2-2 below.

Rewriting Eq 2.2-12

$$\text{lnd}(n,\Delta x,x_i) = \Delta x \sum_{\substack{\Delta x \\ x=x_i}}^{\infty} \frac{1}{x^n} = \int_{\substack{\Delta x \\ x_i}}^{\infty} \frac{1}{x^n}\,\Delta x \quad \text{, the lnd}(n,\Delta x,x) \text{ function where Re}(n)>0$$

Diagram 2.2-2: An lnd(n,Δx,x) n≠1 function x locus line in the complex plane (Refer to Eq 2.2-12)

An lnd(n,Δx,x) n≠1 function x locus line in the complex plane where Re(Δx)>0

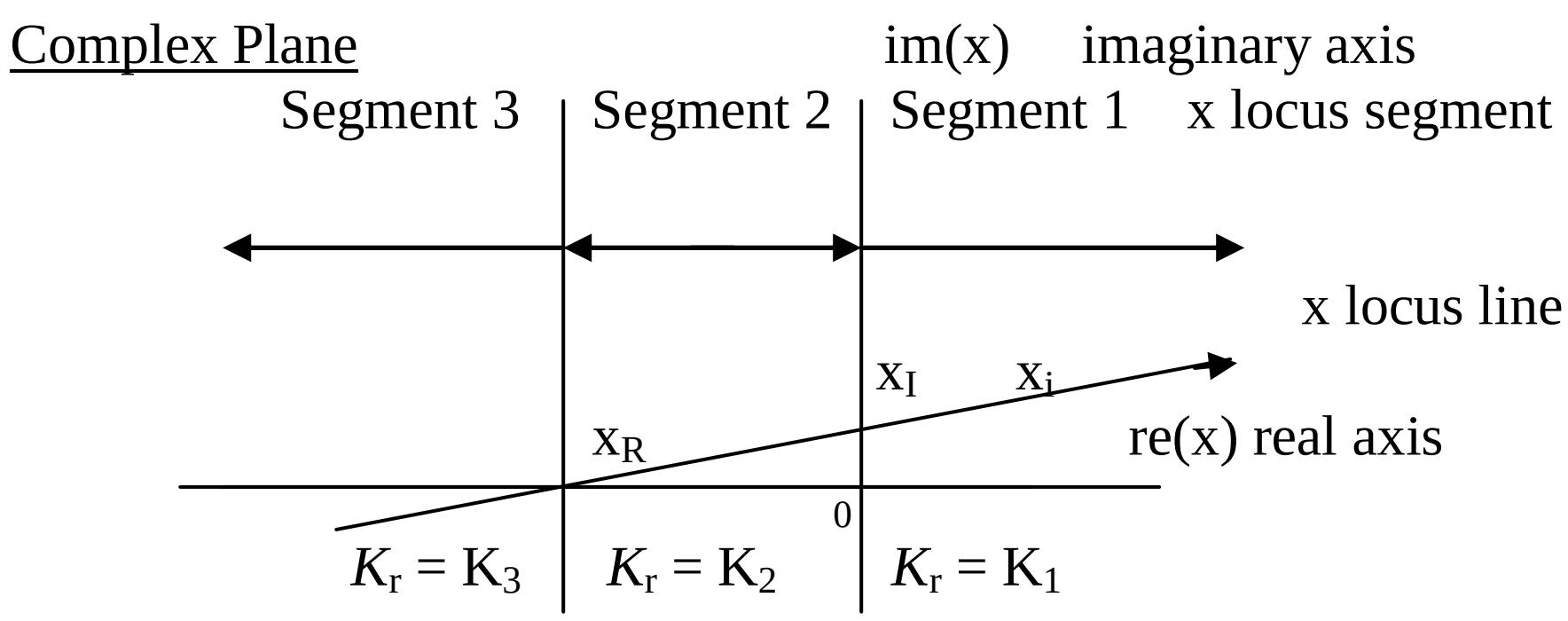

$x_R$ = real axis x locus crossover point
$x_I$ = imaginary axis x locus crossover point
$x = x_i + m\Delta x$ , m=1,2,3,… This is the x locus
$x_i$ = a value of x, the summation initial value x on the x locus line
(See Eq 2.2-12 above)
Re(Δx)>0

Note1 – The arrowhead is placed on the x locus line to indicate the direction of summation.

Note2 – This x locus line diagram is one of many. For Re(Δx)<0, Re(Δx) = 0, Im(Δx) = 0, and different values of $x_i$ and Δx there are other x locus line diagrams. The relationship of this x locus line diagram and the others to the characteristics of their respective lnd(n,Δx,x) n≠1 functions will be explained in the following sections.

Diagram 2.2-2 presents an example of an x locus line in the complex plane. It shows the complex plane divided into three segments. The vertical line through the point, $x_I$, where the x locus line intersects the imaginary axis forms one dividing line. This vertical line, of course, is the imaginary axis. The vertical line through the point, $x_R$, where the x locus line intersects the real axis forms another dividing line. Complex plane Segment 1 is the rightmost segment. Segment 2 is the next complex plane segment to the left and Segment 3 is the leftmost complex plane segment. This diagram, its complex plane segmentation, and x locus line are important to the resolution of the problem previously mentioned. Diagram 2.2-2 and its significance will be described in the following sections.

**Section 2.3: The discovery of a major lnd(n,Δx,x) n≠1 function difficulty**

Once both lnd(n,Δx,x) series were derived, Eq 2.2-3 and Eq 2.2-4, the task of programming them was begun. The lnd(1,Δx,x) Series programming was not difficult. This series was programmed and subsequently tested for proper functionality. In a relatively short time, it was functioning well. In a similar manner the programming of the lnd(n,Δx,x) n≠1 Series was begun. The fact that the lnd(n,Δx,x) n≠1 Series would at all times require a large $|\frac{x}{\Delta x}|$ value to obtain high series accuracy did not cause any undue programming difficulty. To overcome this problem, should a low value of $|\frac{x}{\Delta x}|$ result from the values of x and Δx entered for evaluation, code was written for the computer to term by term calculate the initial part of the summation until a sufficiently large $|\frac{x}{\Delta x}|$ value was obtained. Next, the series was used to evaluate the remainder of the summation. Finally, the two resulting quantities were summed together to obtain the desired summation. This approach worked well. To test the lnd(n,Δx,x) n≠1 program, the x locus line described above and the relationship of Eq 2.2-12 were used. First, a set of values for $x_i$, Δx, and n were chosen. Using these selected values, the program calculated the lnd(n,Δx,x) n≠1 Series value and then its corresponding summation value. If the two values obtained agreed to a large number of decimal places, the test was considered to be a success. This test was repeated again and again by changing the value of $x_i$ to that of another x value on the same x locus line. When the series of tests were all completed satisfactorily, a new lnd(n,Δx,x) n≠1 function was selected and the process was repeated. Initially, it appeared that the lnd(n,Δx,x) n≠1 calculation program was working well. But in time, it became evident that something was not right. Amongst all of the many good tests, a bad test would result. This was very perplexing. An in depth examination of the derivation of the lnd(n,Δx,x) n≠1 Series showed nothing to be wrong. A subsequent review of the programming showed nothing to be wrong but, obviously, something was wrong. Only extensive testing, as described above, indicated the problem. For many values of n, Δx, and $x_i$, the lnd(n,Δx,x) n≠1 calculation program provided correct values to excellent accuracy. However, it was found that for some combinations of n, Δx, and $x_i$ testing along their respective x locus lines produced good and bad values. In the case where Re(Δx)>0 (Refer to Diagram 2.2-2) the lnd(n,Δx,x) n≠1 function values calculated for $x_i$ placed in complex plane Segment 1 were always correct. However, in Segment 2, in Segment 3, or in both of these segments, the values sometimes would not be correct. This characteristic of the lnd(n,Δx,x) n≠1 Series was baffling. The programming was checked so many times and found to be correct that the problem had to involve the series. How could it be that the series, itself, in most cases provided accurate results and then in some small number of cases did not? After much more investigation, it was found that the intersection of the x locus line with the complex plane real and imaginary axes provided an indication as to where calculation difficulties might occur. This is the reason for the three segments shown in Diagram 2.2-2. Additional tedious investigation showed something rather remarkable. Moving left from Segment 1 to Segment 2 to Segment 3 once an incorrect value was calculated close to a dividing line and x locus line intersection, the following calculated values remained incorrect. The value transition between the calculationed good and bad values was surprising abrupt. Only in Segment 1 were the values always calculated correctly. More tedious investigation showed something else quite remarkable. If within a segment where the calculated values were incorrect, an appropriate constant was added to the lnd(n,Δx,x) n≠1 Series, the values calculated in the segment would be correct. From the observations just described, the following conclusion was reached.

Conclusion

The lnd(n,$\Delta$x,x) n $\neq$1 Series function as the anti-derivative of the function, $-\frac{1}{x^n}$ , (where ($D_{\Delta x}$lnd(n,$\Delta$x,x) = $-\frac{1}{x^n}$ , n $\neq$1) is more complicated than it at first appears. True to its derivation, the lnd(n,$\Delta$x,x) n $\neq$1 Series function does have a discrete derivative equal to $-\frac{1}{x^n}$ . However when a series function of specified n and $\Delta$x is plotted for the values of x on its $x=x_i+m\Delta x$, m=integers determined complex plane x locus line, the resulting curve will have either of two forms, a continuous form or a non-continuous form. See Diagram 2.3-1 on the following page.

Diagram 2.3-1: An example of the two lnd(n,Δx,x) n ≠1 Series function forms

Form #1 Continuous Form (The desired anti-derivative form)

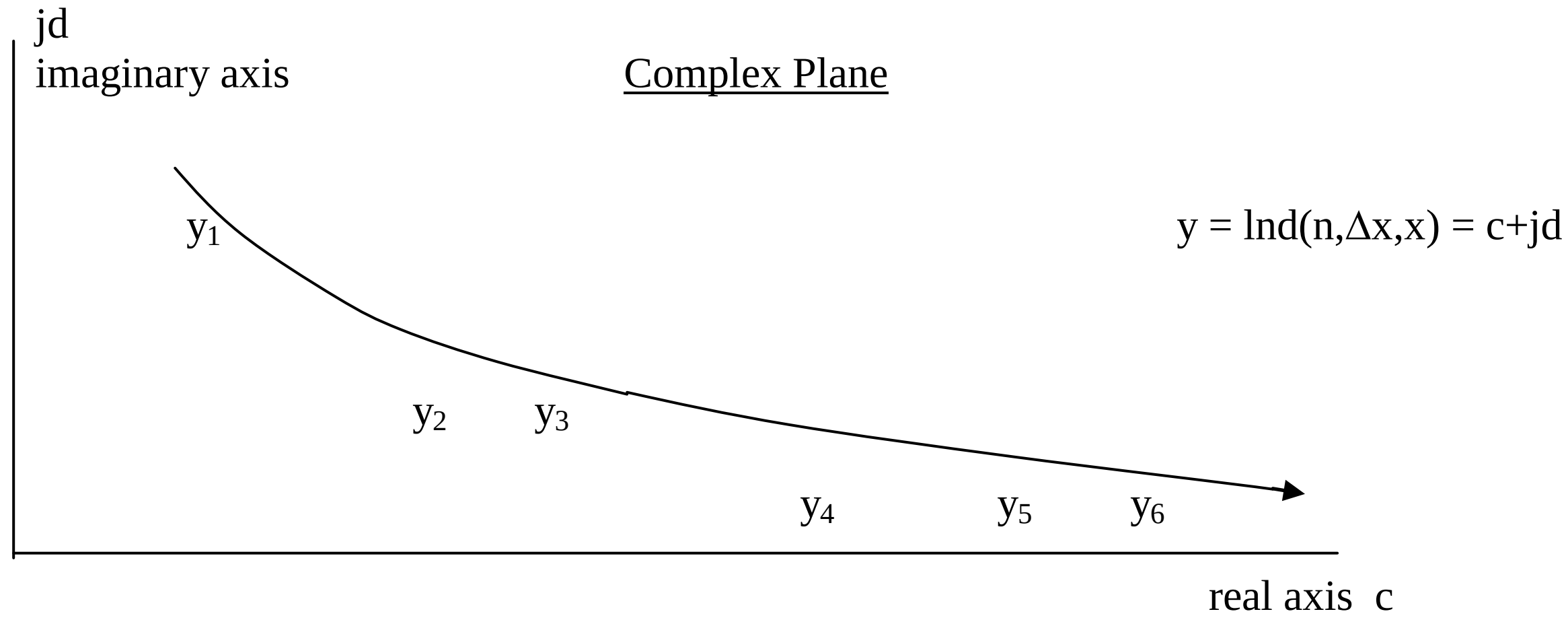

Form #2 Non-continuous Form (The problematic anti-derivative form)

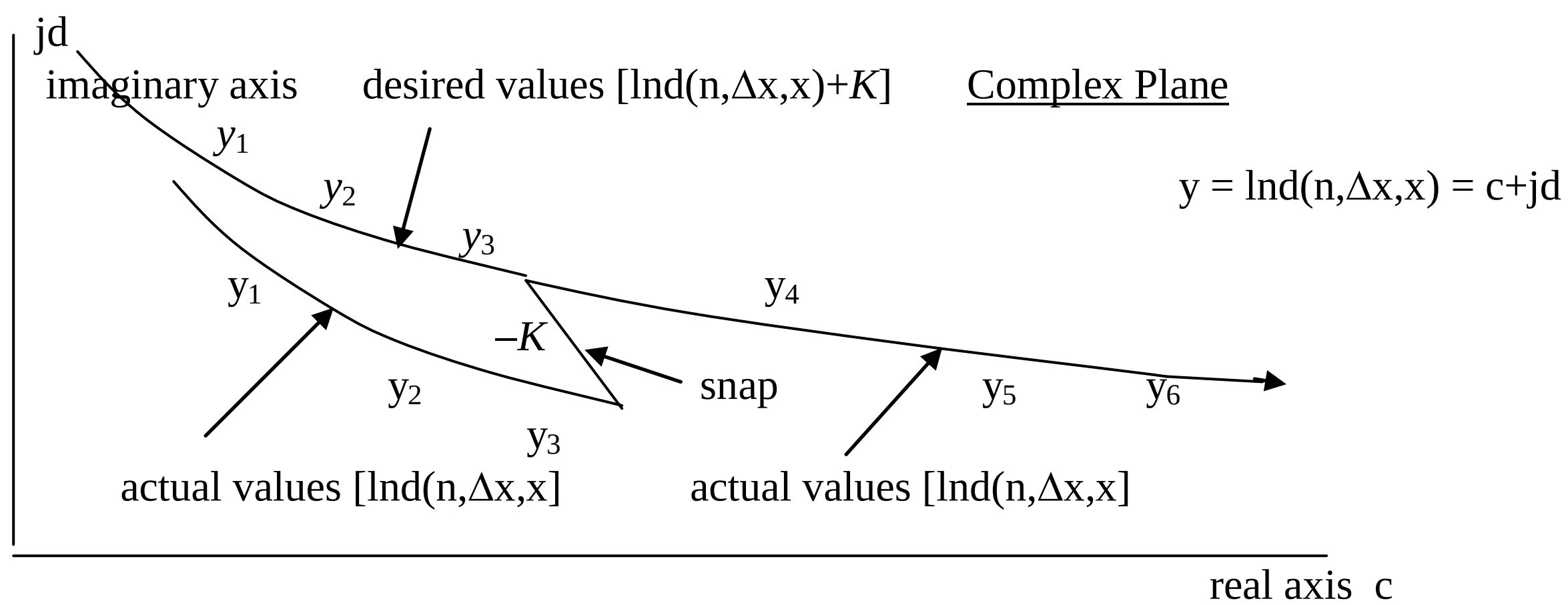

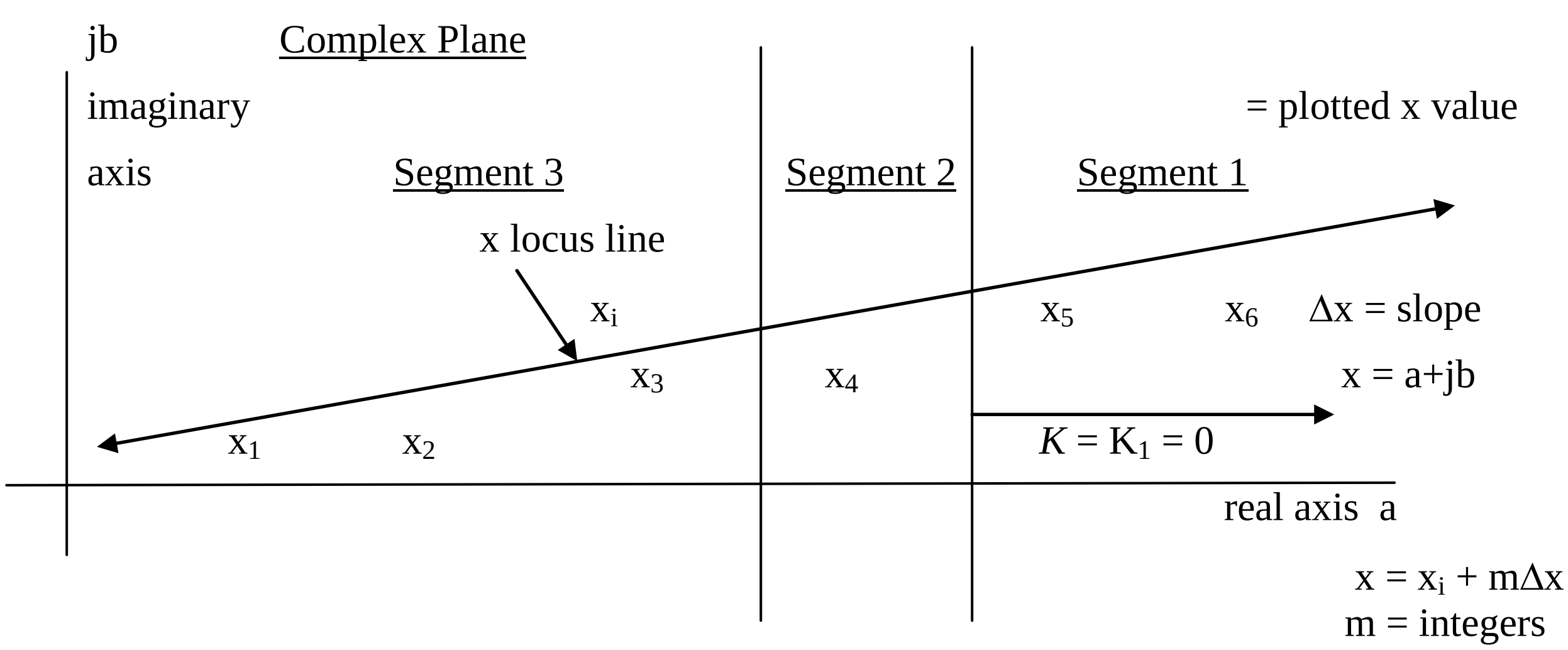

Observing Diagram 2.3-1

Form #1 shows an anti-derivative (integral) of the function, $-\frac{1}{x^n}$ where n≠1, which is a continuous x locus as indicated by the x locus line through the plotted x values. Such a result would provide the correct evaluation of the summation of Eq 2.1-5, the General Zeta Function. Form #2 shows an anti-derivative which is not a continuous x locus as indicated by the x locus line through the plotted x values. Such a result, depending on the location of x on the x locus line, might not provide the correct evaluation of the summation of Eq 2.1-5, the General Zeta Function. Comparing Form #2 to Form #1 it is observed that their discrete derivatives (rates of change) are the same. This is properly the result of the derivation of Eq 2.2-4, the lnd(n,Δx,x) n≠1 Series function. However, for some values of x, n, and Δx, this series will not correctly evaluate the summation of
Eq 2.1-5. Should this anti-derivative form result, a considerable summation evaluation problem is encountered. The results of the summation calculation become unreliable. When and how this result occurs is unknown. No doubt, a full mathematical understanding of this characteristic of the lnd(n,Δx,x) n≠1 Series could be very challenging. A review and comparison of Form #1 and Form #2 provides a means to obtain a "work around" resolution to this problem. It is observed that Form #1 and Form #2 would be the same if the effect labeled "snap" in the Form #2 anti-derivative plot could be removed. It can. The "snap" effect is seen to be the introduction of a constant value at some value of x. This "snap" effect is analogous to the presence of a Heavyside Step Function in a Calculus equation. "Snap" in this paper is defined as a very rapid transition between two real or complex finite values. (The word snap was selected from the terminology "snap action switch". It was felt that the commonly used word, "step" would be somewhat confusing when referring to the transition between two complex values.) By adding a constant value, *K*, in the location where a function snap of −*K* occurs, a Form #1 result will occur, the same result sought to properly evaluate the summation of Eq 2.1-5. It is this value of *K* which has been placed in Eq 2.2-4. However, contrary to what is considered normal, this value of *K* may change over the range of x depending upon the presence of a lnd(n,Δx,x) n≠1 Series function snap condition.

The results of an investigation into the lnd(n,Δx,x) n≠1 Series function snap characteristic

It was discovered that if the lnd(n,Δx,x) n≠1 Series unexpectedly changed value, it would be in very close proximity to the intersection of its linear x locus line and the real and imaginary axes of the complex plane. The change in the value, should it occur, is exceedingly rapid. The series value changes abruptly as would occur in a snap action switch. This unique series characteristic made programming of the lnd(n,Δx,x) function far more complicated than would otherwise be expected. To clarify what has already been stated, a more mathematical explanation is presented below:

The lnd(1,Δx,x) ≡ $\ln_{\Delta x}x$ Series is defined as:

$$D_{\Delta x}\ln_{\Delta x}x = +\frac{1}{x} \quad , \qquad n=1 \qquad (2.3\text{-}1)$$

Integrating (finding the anti-derivative of) Eq 2.3-1

$$\ln_{\Delta x}x \approx \ln\left(\frac{x}{\Delta x}-\frac{1}{2}\right)+\sum_{m=1}^{\infty}\frac{(2m-1)!\,C_m}{(2m+1)!\,2^{2m}\left(\frac{x}{\Delta x}-\frac{1}{2}\right)^{2m}}+\gamma \qquad (2.3\text{-}2)$$

The accuracy of Eq 2.3-2 above increases rapidly for increasing $\left|\frac{x}{\Delta x}\right|$.

Note that the constant of integration, $\gamma$, Euler's Constant (.577215664…), has only one value for all x, that of Euler's Constant

The $lnd(n,\Delta x,x)$ $n\neq1$ Series is defined as:

$$D_{\Delta x}lnd(n,\Delta x,x) = -\frac{1}{x^n} \quad , \quad n\neq1 \qquad (2.3\text{-}3)$$

Integrating (finding the anti-derivative of) Eq 2.3-3

$$lnd(n,\Delta x,x) \approx -\sum_{m=0}^{\infty}\frac{\Gamma(n+2m-1)\left(\frac{\Delta x}{2}\right)^{2m}C_m}{\Gamma(n)(2m+1)!\left(x-\frac{\Delta x}{2}\right)^{n+2m-1}} + K_r \quad , \quad n\neq1 \qquad (2.3\text{-}4)$$

where

r = 1, 2 or 3, the x locus segment designation number

The accuracy of Eq 2.3-4 above increases rapidly for increasing $\left|\frac{x}{\Delta x}\right|$.

Note - The constant of integration, $K_r$, is known to be a function of n as well as r but this is not shown in the constant notation.

Consider the $lnd(n,\Delta x,x)$ $n\neq1$ function x locus line Diagrams, Diagram 2.3-2 and Diagram 2.3-3, shown below.

Diagram 2.3-2: An $lnd(n,\Delta x,x)$ $n\neq1$ function x locus line in the complex plane where $Re(\Delta x)>0$

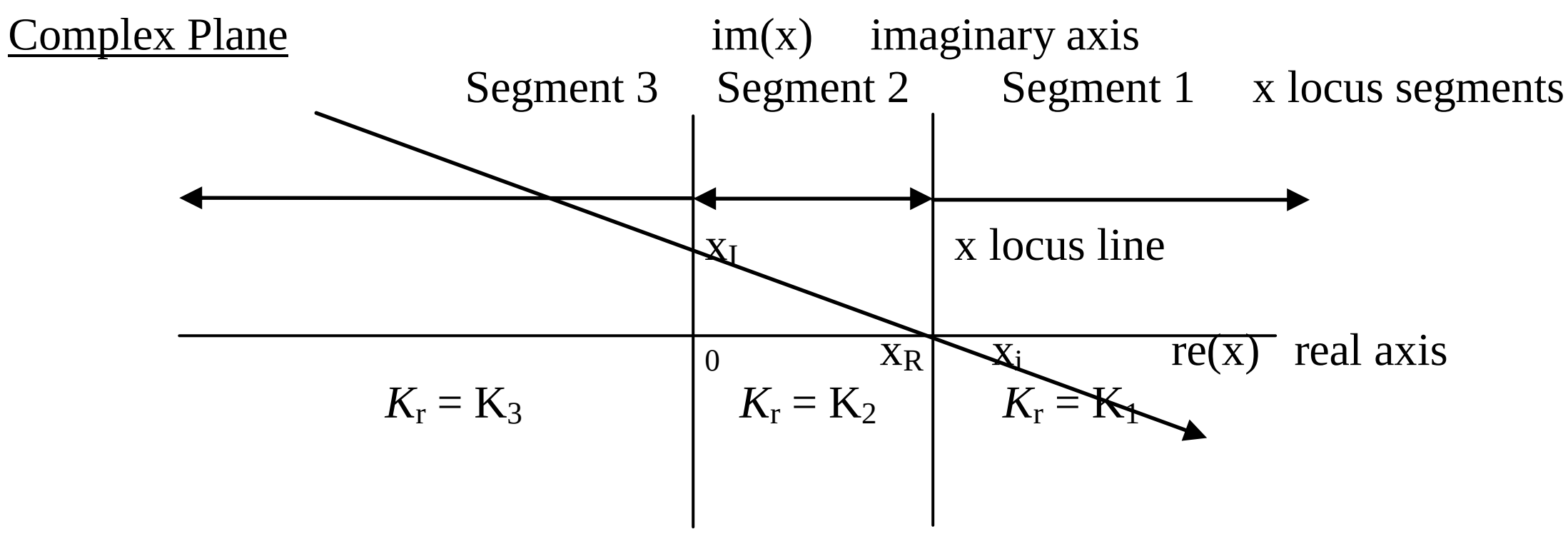

$x_R$ = real axis x locus crossover point
$x_I$ = imaginary axis x locus crossover point
$x = x_i + m\Delta x$ , m=0,1,2,3,… This is the x locus
$x_i$ = a value of x, the summation initial value of x on the x locus line
(See Eq 2.2-12)
$Re(\Delta x)>1$

Diagram 2.3-3: An lnd(n,Δx,x) n≠1 function x locus line in the complex plane where $Re(\Delta x)=0$ and $Im(\Delta x)>0$

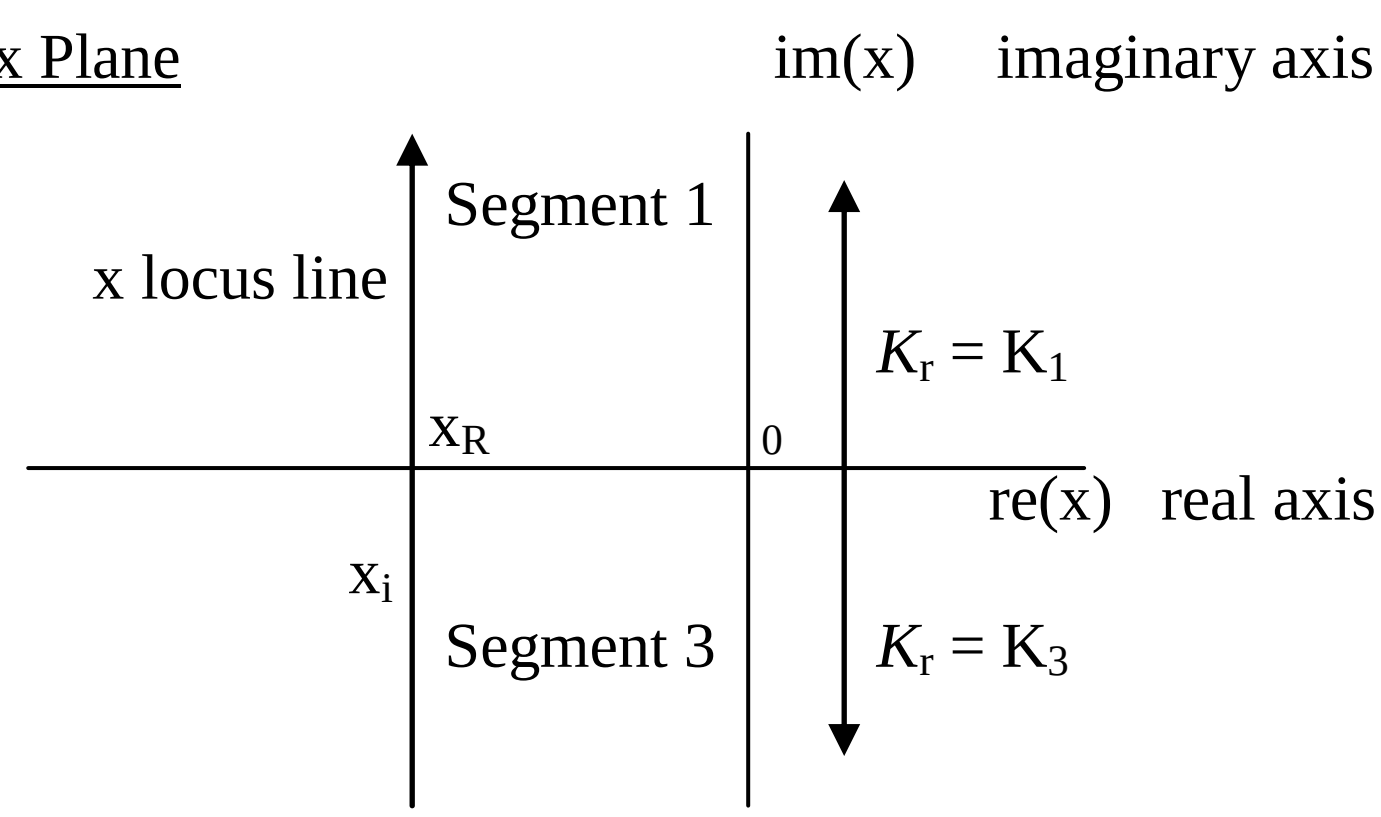

$x_R$ = real axis x locus crossover point
$x = x_i + m\Delta x$ , m=0,1,2,3,… This is the x locus
$x_i$ = a value of x, the summation initial value of x on the x locus line
(See Eq 2.2-12)
$Re(\Delta x)=0$ and $Im(\Delta x)>0$

Note that in this case there is no Segment 2 in the lnd(n,Δx,x) n≠1 function x locus line diagram.

In Diagram 2.3-2, where $Re(\Delta x)>1$, the represented lnd(n,Δx,x) n≠1 Series (See Eq 2.3-4 ) may have as many as three values over all x for its constant of integration, $K_r$, (i.e. $K_r = K_1, K_2, K_3$) or at minimum it may have just one value over all x for its contant of integration, $K_r$, (i.e. $K_r = K_1=K_2=K_3$). Lnd(n,Δx,x) n≠1 Series snap (i.e. snap action transitions) may occur close to $x_R$ and $x_I$, where there are x locus and complex plane axis crossovers. In Diagram 2.3-3, where $Re(\Delta x)=0$ and $Im(\Delta x)>0$, the represented lnd(n,Δx,x) n≠1 Series (See Eq 2.3-4 ) may have as many as two values over all x for its constant of integration, $K_r$, (i.e. $K_r = K_1, K_3$) or at minimum it may have just one value over all x for its contant of integration, $K_r$, (i.e. $K_r = K_1 = K_3$). Lnd(n,Δx,x) n≠1 Series snap (i.e. snap action transitions) may occur close to $x_R$ and $x_I$, (at least one of these crossovers will exist) where there is an x locus and complex plane axis crossover. Why and when snap occurs is not known, but it does occur. To establish if series snap occurs and to correct for it, a considerable amount of programming code had to be written. In particular, the x locus line axis crossovers had to first be identified. Then an evaluation of the lnd(n,Δx,x) function at a point on each side of the axis crossover had to be performed. If the values were found to be inconsistent, this was an indication that series snap had occurred at the crossover boundry. These same calculated values were then used to obtain an integration constant correction. More details concerning the lnd(n,Δx,x) n≠1 Series snap characteritic will be presented in the following sections.

As previously mentioned, the development of the program to calculate the lnd(n,Δx,x) function was slow and tedious, but it was accomplished. The code for calculating the lnd(n,Δx,x) function is written in the UBASIC Programming Language (an advanced form of the Basic Programming Language designed for engineers and mathematicians). UBASIC is run in the Microsoft MSDOS computer environment. The UBASIC Programming Language has all of the features necessary for the calculation of the lnd(n,Δx,x) function. UBASIC performs real and complex evaluations automatically for all algebraic operations and all function computations. Also, the accuracy of the UBASIC computations is at the discretion of the programmer. Typically, all program development calculations were performed with at least 60 places of accuracy.

The UBASIC program for the calculation of the lnd(n,Δx,x) function, LNDX, and a program to calculate the lnd(n,Δx,x) n≠1 Series $C_n$ constants, CNCALC, appear in the PROGRAMS TO CALCULATE LND(N,ΔX,X) AND FORMULA CONSTANTS section at the end of the Appendix. The lnd(n,Δx,x) calculation program retains the coding used for its development. This code may be helpful if further mathematical research is desired to better understand the two series' unique characteristics.

## Section 2.4: The lnd(n,Δx,x) n≠1 Series "snap" characteristic

The series snap characteristic observed during the use of the lnd(n,Δx,x) n≠1 Series was perplexing. The results of an lnd(n,Δx,x) n≠1 Series calculation were often correct but sometimes not. When it became apparent that something about the series was not understood, coding for debugging purposes was added to the lnd(n,Δx,x) n≠1 Series calculation program to observe the program calculations in detail.

Through observation of many lnd(n,Δx,x) n≠1 function computer calculations and the investigation of lnd(n,Δx,x) n≠1 mathematical relationships, the following facts were obtained.

1) The lnd(1,Δx,x) Series (also called the $\ln_{\Delta x}x$ Series), behaves in a straight forward manner. Its constant of integration is one value for all x, Euler's Constant.

2) The lnd(n,Δx,x) n≠1 Series is functional but problematic. For all x, its constant of integration may be unique or it may have as many as six values, three for Re(Δx)>0 and three for Re(Δx)<0.

3) When the the lnd(n,Δx,x) n≠1 Series is used to calculate the lnd(n,Δx,x) function where n≠1, the lnd(n,Δx,x) Series' value can change abruptly at some points yielding an erroneous calculation result. This rather odd characteristic is called " snap" due to the value's surprisingly rapid transition.

4) Plotting the values of x in the complex plane is helpful for finding the x values at which lnd(n,Δx,x) n≠1 Series snap may take place. Series snap has only been observed to occur at an x value very close to either the x locus line real axis intersection or the x locus line imaginary axis intersection. Function snap has not been observed to occur anywhere other than at an x value very close to an x locus line real or imaginary axis crossover.

5) There is no lnd(n,Δx,x) n≠1 Series snap for all x values located on the portion of the x locus line beyond all axis intersections in the direction of summation. For these x values, the constant of integration, *K* , equals zero. (See Eq 2.2-4 and Eq 2.2-12)

6) If series snap is detected, it can be corrected so that the calculated value of lnd(n,Δx,x) n≠1 is valid for all x.

7) Lnd(n,Δx,x) n≠1 Series snap often does not occur but it may.

8) If series snap does not occur, the lnd(n,Δx,x) n≠1 Series constant of integration, *K*, is unique and is equal to 0.

9) Lnd(n,Δx,x) series snap apparently does not occur when n is a negative integer.

10) There is not, at this time, a definitive mathematical method to predict whether series snap will occur given particular values of n, Δx, and x.

11) The developed lnd(n,Δx,x) function computation program, LNDX, has provided correct high accuracy results for an exceeding large number of computations.

**The lnd(n,Δx,x) n≠1 Series Snap Hypothesis**

**Lnd(n,Δx,x) n≠1 Series snap will occur, if it occurs at all, only in close proximity to an x locus line transition across a complex plane axis and, except where series snap occurs, the lnd(n,Δx,x) n≠1 Series constant of integration will not change.**

Comment

The above hypothesis statement has been observed to be true in a great many trials. However, a proof of this statement has not yet been found.

The determination of the value of the lnd(n,Δx,x) n≠1 Series constant of integration is very important. The methodology by which to calculate the constant of integration value is presented in the next section, Section 2.5.

## Section 2.5: The lnd(n,Δx,x) n≠1 Series constants of integration and their calculation

There are six possible values for the lnd(n,Δx,x) n≠1 Series constant of integration. There are as many as three constant values for Re(Δx) > 0 or two constant values for the case where Re(Δx)=0 and Im(Δx)>0. Also, there are as many as three constant values for Re(Δx) < 0 or two constant values for the case where Re(Δx)=0 and Im(Δx)<0. The sign of Re(Δx), or if Re(Δx)=0 the sign of Im(Δx), determines the direction of summation from the initial value of x, $x_i$, along the lnd(n,Δx,x) function complex plane x locus line. The x locus diagrams in Diagram 2.5-1 and Diagram 2.5-2 below show the lnd(n,Δx,x) n≠1 function summation direction and also the axis crossover segments of several x loci. These segments are designated as 1, 2, and 3. As shown, it is possible that there is no segment 2. For each x locus segment there are two lnd(n,Δx,x) n≠1 Series constant of integration values. A constant value may be equal to the constant value of another segment or it may not. The integration constant value for each segment is designated as $K_r$ or $k_r$ where r is the segment designation number, 1, 2, or 3. $K_r$ designates the constants for Re(Δx)>0 or Re(Δx)=0 and Im(Δx)>0. $k_r$ designates the constants for Re(Δx)<0 or Re(Δx)=0 and Im(Δx)<0.

Applicable Equations for the following diagrams, Diagram 2.5-1 and Diagram 2.5-2

$$\mathrm{lnd}(n,\Delta x,x) \approx -\sum_{m=0}^{\infty} \frac{\Gamma(n+2m-1)\left(\frac{\Delta x}{2}\right)^{2m} C_m}{\Gamma(n)(2m+1)!\left(x-\frac{\Delta x}{2}\right)^{n+2m-1}} + \begin{cases} K_r \\ k_r \end{cases}, \quad n\neq 1 \qquad (2.5\text{-}1)$$

$K_r$ = constants of integration for Re(Δx) > 0 or Re(Δx) = 0 and Im(Δx) > 0
$k_r$ = constants of integration for Re(Δx) < 0 or Re(Δx) = 0 and Im(Δx) < 0
$C_m$ = series constants
r = 1, 2, or 3 , The x locus segment designations

Accuracy increases rapidly as $\left|\frac{x}{\Delta x}\right|$ increases in value.

**$K_1 = 0$** **$k_3 = 0$** <u>Comment</u> – The proof of $K_1$ and $k_3$ being equal to 0 depends on the validity of the lnd(n,Δx,x) n≠1 Series Snap Hypothesis. The proof is in the following section, Section 2.6.

$$\Delta x \sum_{\Delta x\; x=x_1}^{x_2} \frac{1}{x^n} = \int_{\Delta x\; x_1}^{x_2+\Delta x} \frac{1}{x^n}\,\Delta x = -\mathrm{lnd}(n,\Delta x,x)\,\Big|_{x_1}^{x_2+\Delta x}, \quad n\neq 1 \qquad (2.5\text{-}2)$$

where

$x = x_1, x_1+\Delta x, x_1+2\Delta x, x_1+3\Delta x, \ldots, x_2-\Delta x, x_2$

or

$x = x_1+m\Delta x$ , $m=0,1,2,3,\ldots,M$

$M = \frac{x_2 - x_1}{\Delta x}$

$x_1,x_2$ may be infinite

Δx = x increment

n,Δx,x = real or complex values

Diagram 2.5-1: Some complex plane x loci line diagrams for Eq 2.5-1 where $Re(\Delta x) > 0$ or $Re(\Delta x) = 0$ and $Im(\Delta x) > 0$

The x locus is $x = x_i + m\Delta x$ , m = integers , $Re(\Delta x) > 0$ or $Re(\Delta x) = 0$ and $Im(\Delta x) > 0$, which includes the lnd(n,Δx,x) n≠1 summation x values plotted in the complex plane (See Eq 2.2-12). The x locus line determined from the x locus is an infinitely long line through all the x values plotted in the complex plane. The locus of the summation x values from $x_1$ to $x_2$ is on a segment of the x locus line. Summation progresses from $x_1$ to $x_2$, here, left to right or bottom to top in the complex plane.

.

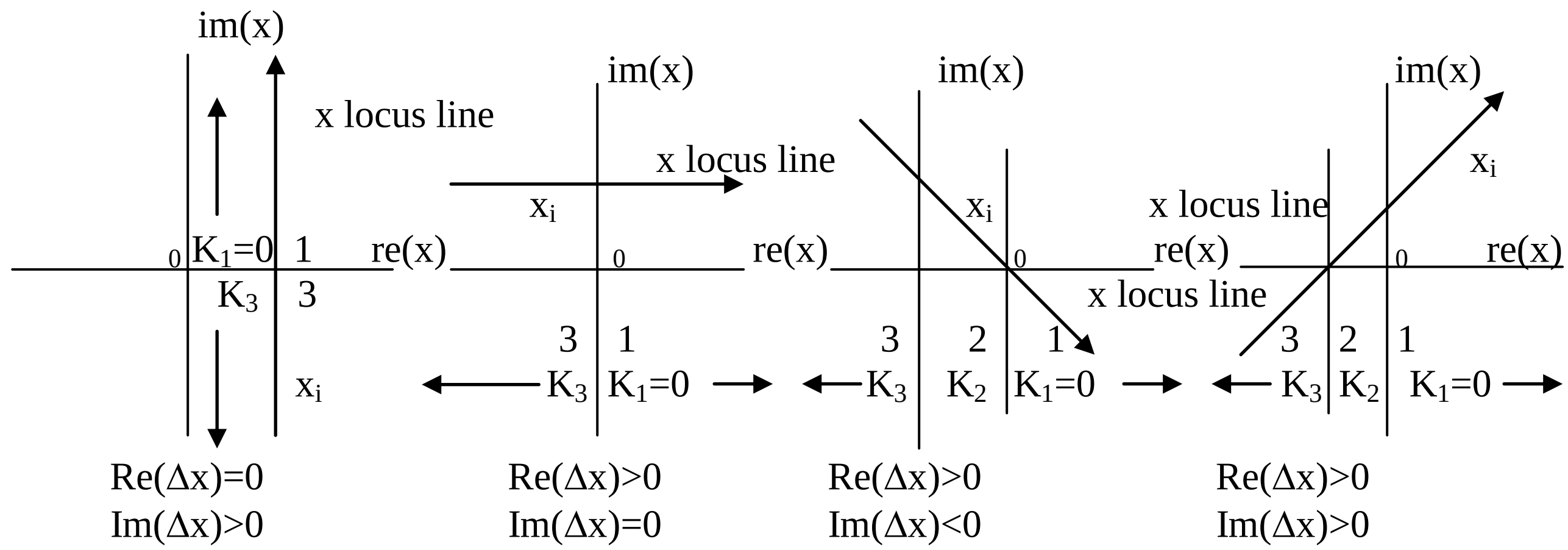

$Re(\Delta x) > 0$ or $Re(\Delta x) = 0$ and $Im(\Delta x) > 0$

$x = x_i + m\Delta x$ , m = integers The value of x goes from left to right in the complex plane

$x_i$ = a value of x, the summation initial value of x (See Eq 2.2-12)

1,2,3 = x locus segment designations (Note that for some x loci there may not be a segment 2)

1 is the rightmost segment of the x locus and 3 is the leftmost segment of the x locus.

The complex plane axes divide the x locus into two or three segments.

If lnd(n,Δx,x) snap occurs, it will occur only in close proximity to an x locus line transition across a complex plane axis. (See the lnd(n,Δx,x) n≠1 Series Snap Hypothesis previously stated.)

Diagram 2.5-2: Some complex plane x loci line diagrams for Eq 2.5-1 where $Re(\Delta x) < 0$ or $Re(\Delta x) = 0$ and $Im(\Delta x) < 0$

The x locus is $x = x_i + m\Delta x$ , m = integers , $Re(\Delta x) < 0$ or $Re(\Delta x) = 0$ and $Im(\Delta x) > 0$, which includes the lnd(n,Δx,x) n≠1 summation x values plotted in the complex plane (See Eq 2.2-12). The x locus line determined from the x locus is an infinitely long line through all of the x values plotted in the complex plane. The locus of the summation x values from $x_1$ to $x_2$ is on a segment of the x locus line. Summation progresses from $x_1$ to $x_2$, here, right to left or top to bottom in the complex plane.

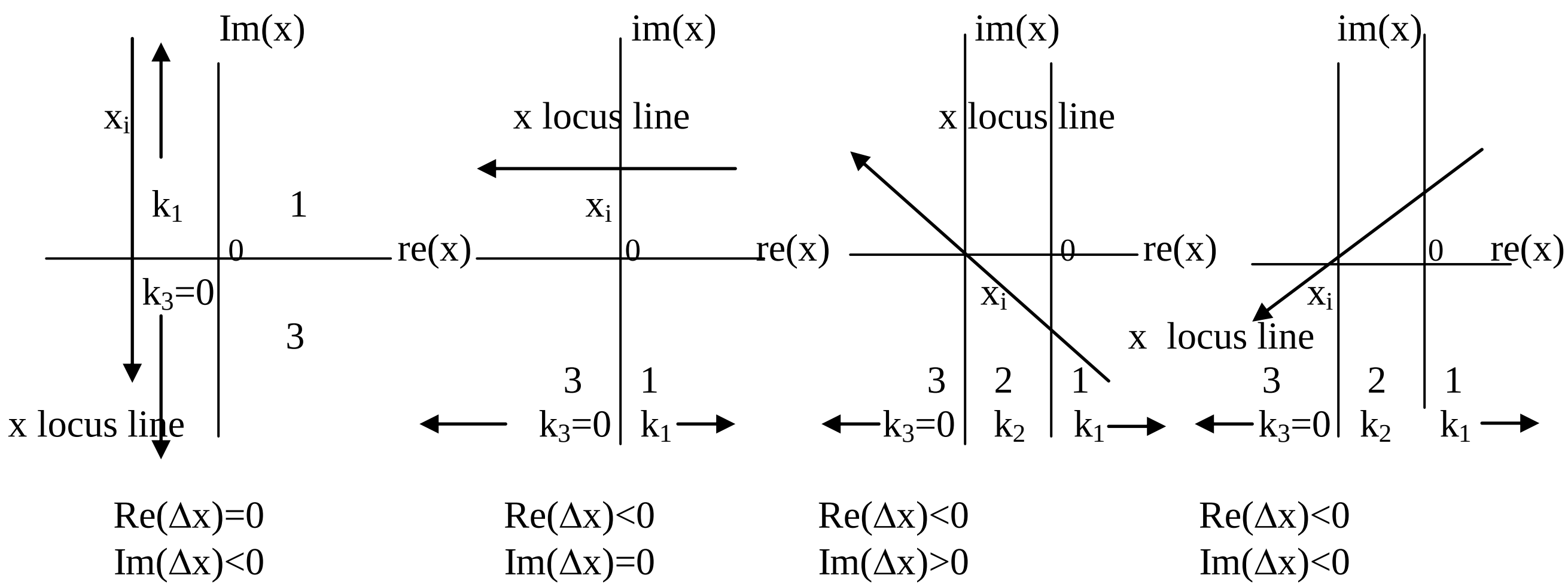

$Re(\Delta x) < 0$ or $Re(\Delta x) = 0$ and $Im(\Delta x) < 0$

$x = x_i + m\Delta x$ , m = integers The value of x goes from right to left in the complex plane

$x_i$ = a value of x, the summation initial value of x (See Eq 2.2-12)

1,2,3 = x locus segment designations (Note that for some x loci there may not be a segment 2)

For $Re(\Delta x)\neq 0$, 1 is the rightmost segment of the x locus and 3 is the leftmost segment of the x locus.

For $Re(\Delta x)=0$ and $Im(\Delta x)\neq 0$, 1 is the uppermost segment of the x locus and 3 is the lower segment of the x locus.

The complex plane axes divide the x locus into two or three segments.

If lnd(n,Δx,x) snap occurs, it will occur only in close proximity to an x locus transition across a complex plane axis. (See the lnd(n,Δx,x) n≠1 Series Snap Hypothesis previously stated.)

Diagram 2.5-1 and Diagram 2.5-2 identify the x regions where constants of integration of different values (as many as three) may exist. It would be very useful to know the values of these constants given the values of n, Δx, and $x_i$. Some equations will now be derived to calculate the value of the Eq 2.5-1 constants of integration designated $K_1$, $K_2$, $K_3$, $k_1$, $k_2$, and $k_3$.

Derivation of some equations to calculate the value of the Eq 2.5-1 constants of integration for a specified x locus

Derive some very useful equations relating to the constants of integration of the function, lnd(n,Δx,x) n≠1 for a specified summation x locus in the complex plane.The specified x locus line is the straight line through all of the summation x points plotted in the complex plane. See the following diagrams, Diagram 2.5-3 and Diagram 2.5-4.

Diagram 2.5-3: Representative x locus and x locus line in the complex plane for Re(Δx) ≠ 0

Re(Δx) ≠ 0

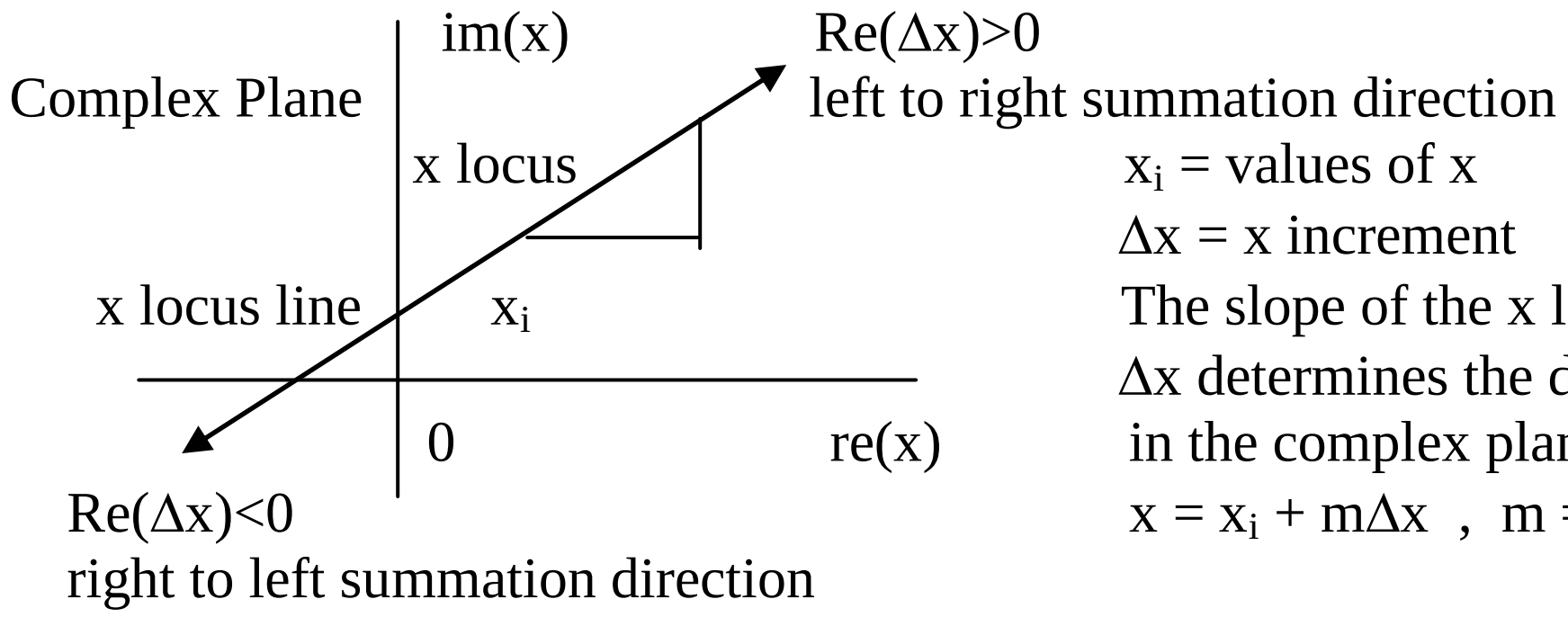

$x_i$ = values of x
Δx = x increment
The slope of the x locus line = Δx
Δx determines the direction of summation in the complex plane
$x = x_i + m\Delta x$ , m = integers

Left to right summation direction where Re(Δx)>0
Right to left summation direction where Re(Δx)<0

$$\Delta x \sum_{\substack{\Delta x \\ x=x_i}}^{\pm\infty} \frac{1}{x^n} = \text{lnd}(n,\Delta x,x_i) \tag{2.5-3}$$

where
$x = x_i + m\Delta x$
m = 0,1,2,3,… for Re(Δx)>0
m = 0,-1,-2,-3,… for Re(Δx)<0
$x_i$ to +∞ for Re(Δx)>0
$x_i$ to -∞ for Re(Δx)<0
Re(n)>1

$$\Delta x \sum_{\substack{\Delta x \\ x=\mp\infty}}^{\pm\infty} \frac{1}{x^n} = \text{lnd}(n,\Delta x,x_i) - \text{lnd}(n,-\Delta x,x_i-\Delta x) \tag{2.5-4}$$

where
$x = x_i + m\Delta x$
m = integers
-∞ to +∞ for Re(Δx)>0
+∞ to -∞ for Re(Δx)<0
Re(n)>1

Note – Eq 2.5-4 is derived from Eq 2.5-3

Diagram 2.5-4: Representative x locus and x locus line in the complex plane for $Re(\Delta x) = 0$

$Re(\Delta x) = 0$

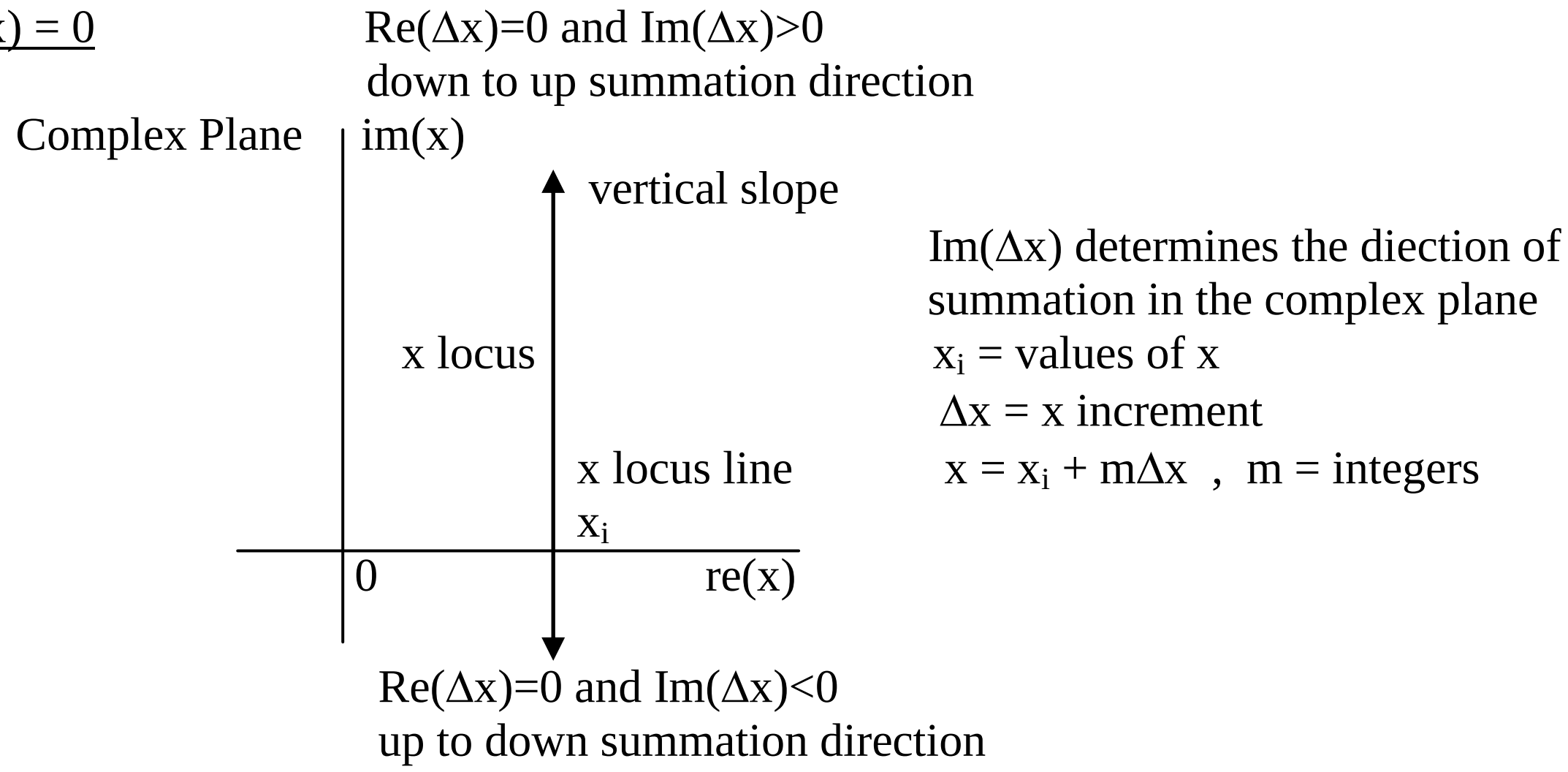

Down to up summation direction where $Re(\Delta x)=0$ and $Im(\Delta x)>0$
Up to down summation direction where $Re(\Delta x)=0$ and $Im(\Delta x)<0$

$$\Delta x \sum_{\substack{\Delta x \\ x=x_i}}^{\pm\infty} \frac{1}{x^n} = \mathrm{lnd}(n,\Delta x,x_i) \qquad (2.5\text{-}5)$$

where
$x = x_i + m\Delta x$
$m = 0,1,2,3,\ldots$ for $Re(\Delta x)=0$ and $Im(\Delta x)>0$
$m = 0,-1,-2,-3,\ldots$ for $Re(\Delta x)=0$ and $Im(\Delta x)<0$
$x_i$ to $+\infty$ for $Re(\Delta x)=0$ and $Im(\Delta x)>0$
$x_i$ to $-\infty$ for for $Re(\Delta x)=0$ and $Im(\Delta x)<0$
$Re(n)>1$

$$\Delta x \sum_{\substack{\Delta x \\ x=\mp\infty}}^{\pm\infty} \frac{1}{x^n} = \mathrm{lnd}(n,\Delta x,x_i) - \mathrm{lnd}(n,-\Delta x,x_i-\Delta x) \qquad (2.5\text{-}6)$$

where
$x = x_i + m\Delta x$
$m$ = integers
$-\infty$ to $+\infty$ for $Re(\Delta x)=0$ and $Im(\Delta x)>0$
$+\infty$ to $-\infty$ for $Re(\Delta x)=0$ and $Im(\Delta x)<0$
$Re(n)>1$

Note – Eq 2.5-4 is derived from Eq 2.5-5

Definitions and comments concerning the x increment, $\Delta x$

(Refer to Diagram 2.5-3 and Diagram 2.5-4)

1. $\Delta x$ = summation x increment
2. $\Delta x = \boldsymbol{\Delta x}$ or $-\boldsymbol{\Delta x}$
3. $\boldsymbol{\Delta x}$ = all $\Delta x$ values where $Re(\Delta x)>0$ or $\{Re(\Delta x)=0$ and $Im(\Delta x)>0\}$
4. $\boldsymbol{\Delta x}$ is $\Delta x$, the summation x increment, where $Re(\Delta x)>0$ or $\{Re(\Delta x)=0$ and $Im(\Delta x)>0\}$
5. $-\boldsymbol{\Delta x}$ = all $\Delta x$ values where $Re(\Delta x)<0$ or $\{Re(\Delta x)=0$ and $Im(\Delta x)<0\}$
6. $-\boldsymbol{\Delta x}$ is $\Delta x$, the summation x increment, where $Re(\Delta x)<0$ or $\{Re(\Delta x)=0$ and $Im(\Delta x)<0\}$
7. The x locus line is the straight line in the complex plane through all of the plotted summation x points
8. The slope of a summation x locus line in the complex plane is determined by $\Delta x$.
9. $\{x_i$ and $\boldsymbol{\Delta x}\}$ or $\{x_i$ and $-\boldsymbol{\Delta x}\}$ determine the same x locus line in the complex plane.
10. When $Re(\Delta x)>0$, the summation progresses along its complex plane x locus line from left to right or when $Re(\Delta x)=0$ and $Im(\Delta x)>0$ from down to up.
11. When $Re(\Delta x)<0$, the summation progresses along its complex plane x locus line from right to left or when $Re(\Delta x)=0$ and $Im(\Delta x)<0$ from up to down.

$$\Delta x = \boldsymbol{\Delta x} \quad (Re(\Delta x)>0 \text{ or } \{Re(\Delta x)=0 \text{ and } Im(\Delta x)>0\}) \tag{2.5-7}$$

From Eq 2.5-1 and Eq 2.5-7

$$\mathrm{lnd}(n, \boldsymbol{\Delta x}, x) \approx -\sum_{m=0}^{\infty} \frac{\Gamma(n+2m-1)\left(\frac{\boldsymbol{\Delta x}}{2}\right)^{2m} C_m}{\Gamma(n)(2m+1)!\left(x-\frac{\boldsymbol{\Delta x}}{2}\right)^{n+2m-1}} + K_r\ , \quad n\neq 1 \tag{2.5-8}$$

Accuracy increases rapidly as $|\frac{x}{\Delta x}|$ increases in value.

r = 1, 2, or 3 , The x locus segment designations

$Re(\Delta x)>0$ or $\{Re(\Delta x)=0$ and $Im(\Delta x)>0\}$ so the constant of integration is $K_r$

$$\Delta x = -\boldsymbol{\Delta x} \quad (Re(\Delta x)<0 \text{ or } \{Re(\Delta x)=0 \text{ and } Im(\Delta x)<0\}) \tag{2.5-9}$$

From Eq 2.5-1 and Eq 2.5-9

$$\mathrm{lnd}(n, -\boldsymbol{\Delta x}, x) \approx -\sum_{m=0}^{\infty} \frac{\Gamma(n+2m-1)\left(\frac{\boldsymbol{\Delta x}}{2}\right)^{2m} C_m}{\Gamma(n)(2m+1)!\left(x+\frac{\boldsymbol{\Delta x}}{2}\right)^{n+2m-1}} + k_r\ , \quad n\neq 1 \tag{2.5-10}$$

Accuracy increases rapidly as $|\frac{x}{\Delta x}|$ increases in value.

$r = 1, 2, \text{ or } 3$
$Re(-\Delta x)<0$ or $\{Re(\Delta x)=0$ and $Im(-\Delta x)<0\}$ so the constant of integration is $k_r$

Change the value of x in Eq 2.5-10 to $x-\Delta\mathbf{x}$

$$lnd(n,-\Delta x,x-\Delta\mathbf{x}) \approx -\sum_{m=0}^{\infty} \frac{\Gamma(n+2m-1)\left(\frac{\Delta\mathbf{x}}{2}\right)^{2m} C_m}{\Gamma(n)(2m+1)!\left(x-\frac{\Delta\mathbf{x}}{2}\right)^{n+2m-1}} + k_r \;,\quad n\neq 1 \tag{2.5-11}$$

Subtracting Eq 2.5-11 from Eq 2.5-8 the series vanishes

$$lnd(n,\Delta\mathbf{x},x) - lnd(n,-\Delta\mathbf{x},x-\Delta\mathbf{x}) \approx K_r - k_r \;,\quad n\neq 1 \;,\; r = 1,2,3 \tag{2.5-12}$$

where
$x = x_i + m\Delta x$ , m = integers , This is an x locus in the complex plane
$x_i$ = a value of x
$\Delta x$ = x increment
$Re(\Delta x)>0$ or $\{Re(\Delta x)=0$ and $Im(\Delta x)>0\}$
x = real or complex values
$n,\Delta x,x_i,K_r,k_r$ = real or complex constants
The x locus line is the straight line in the complex plane through all of the plotted x points
Accuracy increases rapidly as $\left|\frac{x}{\Delta x}\right|$ increases in value.

$$\Delta x = -\Delta\mathbf{x} \quad (Re(\Delta x)<0 \text{ or } \{Re(\Delta x)=0 \text{ and } Im(\Delta x)<0\}) \tag{2.5-13}$$

From Eq 2.5-1 and Eq 2.5-13

$$lnd(n,-\Delta\mathbf{x},x) \approx -\sum_{m=0}^{\infty} \frac{\Gamma(n+2m-1)\left(\frac{\Delta\mathbf{x}}{2}\right)^{2m} C_m}{\Gamma(n)(2m+1)!\left(x+\frac{\Delta\mathbf{x}}{2}\right)^{n+2m-1}} + k_r \;,\quad n\neq 1 \tag{2.5-14}$$

Accuracy increases rapidly as $\left|\frac{x}{\Delta x}\right|$ increases in value.

$r = 1, 2, \text{ or } 3$
$Re(\Delta x)<0$ or $\{Re(\Delta x)=0$ and $Im(\Delta x)<0\}$ so the constant of integration is $k_r$

$$\Delta x = \Delta\mathbf{x} \quad (Re(\Delta x)>0 \text{ or } \{Re(\Delta x)=0 \text{ and } Im(\Delta x)>0\}) \tag{2.5-15}$$

From Eq 2.5-1 and Eq 2.5-15

$$\text{lnd}(n, \Delta x, x) \approx -\sum_{m=0}^{\infty} \frac{\Gamma(n+2m-1)\left(\frac{\Delta x}{2}\right)^{2m} C_m}{\Gamma(n)(2m+1)!\left(x - \frac{\Delta x}{2}\right)^{n+2m-1}} + K_r \;,\; n \neq 1 \tag{2.5-16}$$

Accuracy increases rapidly as $|\frac{x}{\Delta x}|$ increases in value.

r = 1, 2, or 3

Re(Δx)>0 or {Re(Δx)=0 and Im(Δx)>0} so the constant of integration is $K_r$

Change the value of x in Eq 2.5-16 to x+Δx

$$\text{lnd}(n, \Delta x, x+\Delta x) \approx -\sum_{m=0}^{\infty} \frac{\Gamma(n+2m-1)\left(\frac{\Delta x}{2}\right)^{2m} C_m}{\Gamma(n)(2m+1)!\left(x + \frac{\Delta x}{2}\right)^{n+2m-1}} + K_r \;,\; n \neq 1 \tag{2.5-17}$$

Subtracting Eq 2.5-17 from Eq 2.5-15 the series vanishes

$$\text{lnd}(n,-\Delta x,x) - \text{lnd}(n, \Delta x,x+\Delta x) \approx -(K_r - r_r), \quad n \neq 1 \;,\; r = 1,2,3 \tag{2.5-18}$$

where

$x = x_i + m\Delta x$ ,  m = integers ,  This is an x locus in the complex plane

$x_i$ = a value of x

Δx = x increment

Re(Δx)<0 or {Re(Δx)=0 and Im(Δx)<0}

x = real or complex values

$n,\Delta x,x_i,K_r,k_r$ = real or complex constants

The x locus line is the straight line in the complex plane through all of the plotted x points

Accuracy increases rapidly as $|\frac{x}{\Delta x}|$ increases in value.

Rewriting Eq 2.5-12 and Eq 2.5-18

$$\text{lnd}(n, [\Delta x],x) - \text{lnd}(n,-[\Delta x],x-[\Delta x]) \approx +(K_r - k_r) \;,\; n \neq 1 \;,\; r = 1,2,3 \tag{2.5-19}$$

$$\text{lnd}(n,[-\Delta x],x) - \text{lnd}(n,-[-\Delta x],x-[-\Delta x]) \approx -(K_r - r_r), \quad n \neq 1 \;,\; r = 1,2,3 \tag{2.5-20}$$

It is possible to merge Eq 2.5-19 and Eq 2.5-20 into one equation. This equation is shown in Eq 2.5-21 below.

$$\mathbf{lnd(n,\Delta x,x) - lnd(n,-\Delta x,x-\Delta x) \approx \pm(K_r - k_r),\quad n\neq 1\ ,\ r = 1,2,3} \qquad (2.5\text{-}21)$$

**where**

**$x = x_i + m\Delta x$ , m = integers , This is an x locus in the complex plane**

**$x_i$ = a value of x**

**$\Delta x$ = x increment**

**+ for $Re(\Delta x)>0$ or $\{Re(\Delta x)=0$ and $Im(\Delta x)>0\}$**

**– for $Re(\Delta x)<0$ or $\{Re(\Delta x)=0$ and $Im(\Delta x)<0\}$**

**x = real or complex values**

**$n,\Delta x,x_i,K_r,k_r$ = real or complex constants**

**The x locus line is the straight line in the complex plane through all of the plotted x points**

**Accuracy increases rapidly as $|\frac{x}{\Delta x}|$ increases in value.**

For a large absolute value, $|\frac{x}{\Delta x}|$ , it is seen that Eq 2.5-14 is accurate. The question arises as to what happens when the absolute value, $|\frac{x}{\Delta x}|$ , is not large. Eq 2.5-14 is problematic for this reason. It would be very helpful if the accuracy restriction on Eq 2.5-14 could be removed. This can and will be done below.

Find another relationship involving the quantity, lnd(n,Δx,x) – lnd(n,-Δx,x-Δx).
Consider summing both ways along an x locus line in the complex plane. See the following diagram, Diagram 2.5-5.

Diagram 2.5-5: An x locus line in the complex plane

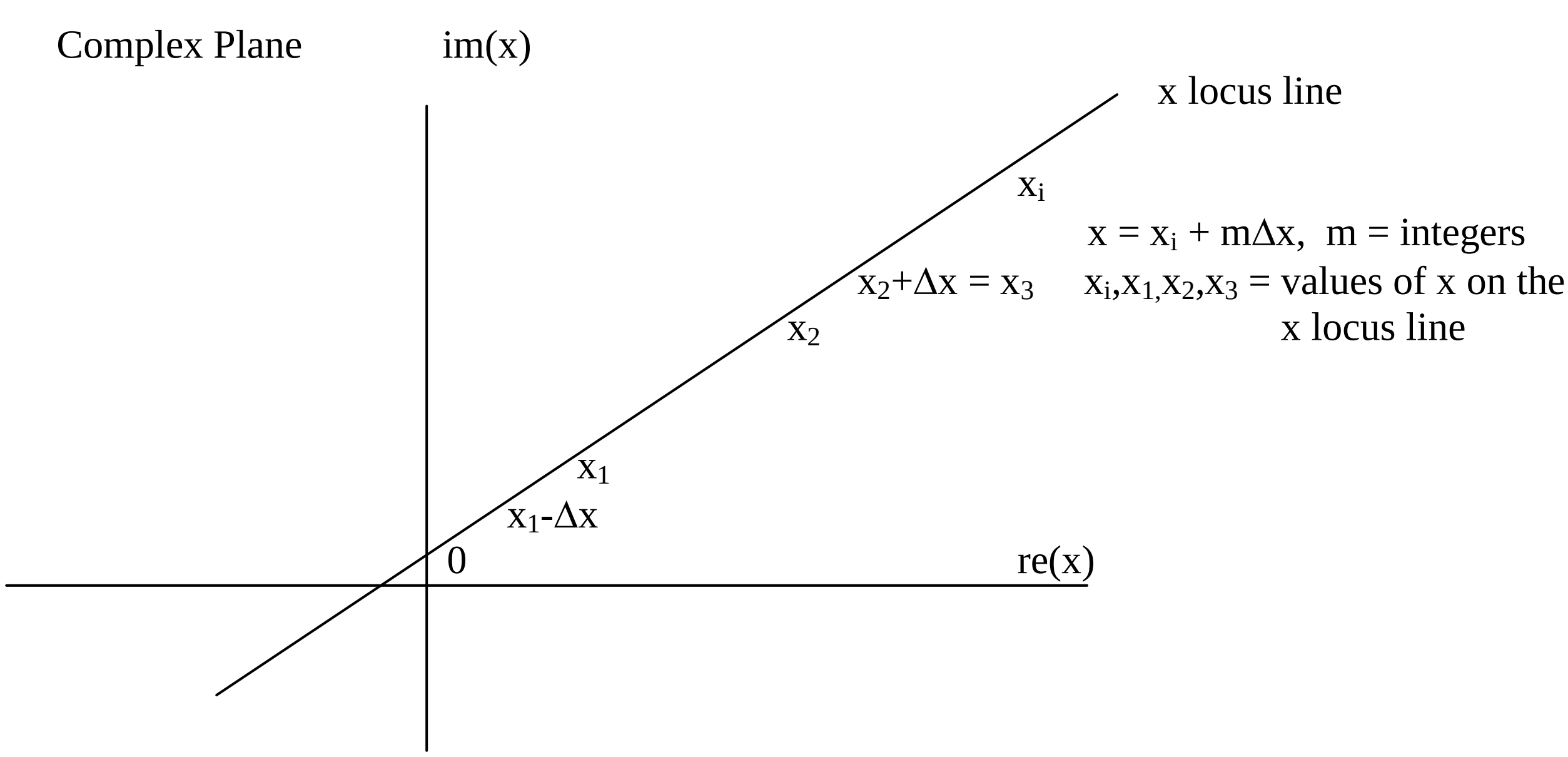

$$\Delta x \sum_{x=x_1}^{x_2}{}_{\Delta x}\frac{1}{x^n} = {}_{\Delta x}\int_{x_1}^{x_2+\Delta x}\frac{1}{x^n}\,\Delta x = [\mathrm{lnd}(n,\Delta x,x_2+\Delta x) - \mathrm{lnd}(n,\Delta x,x_1)]\ ,\ \ n\neq 1 \tag{2.5-22}$$

$$-\Delta x \sum_{x=x_2}^{x_1}{}_{-\Delta x}\frac{1}{x^n} = {}_{-\Delta x}\int_{x_2}^{x_1-\Delta x}\frac{1}{x^n}\,\Delta x = [\mathrm{lnd}(n,-\Delta x,x_1-\Delta x) - \mathrm{lnd}(n,-\Delta x,x_2)]\ ,\ \ n\neq 1 \tag{2.5-23}$$

But

$$\sum_{x=x_1}^{x_2}{}_{\Delta x}\frac{1}{x^n} = \sum_{x=x_2}^{x_1}{}_{-\Delta x}\frac{1}{x^n} \tag{2.5-24}$$

Adding Eq 2.5-22 and Eq 2.5-23

$$0 = [\mathrm{lnd}(n,\Delta x,x_2+\Delta x) - \mathrm{lnd}(n,\Delta x,x_1)] + [\mathrm{lnd}(n,-\Delta x,x_1-\Delta x) - \mathrm{lnd}(n,-\Delta x,x_2)] \tag{2.5-25}$$

Rearranging terms

$$[\mathrm{lnd}(n,\Delta x,x_1) - \mathrm{lnd}(n,-\Delta x,x_1-\Delta x)] = [\mathrm{lnd}(n,\Delta x,x_2+\Delta x) - \mathrm{lnd}(n,-\Delta x,x_2)] \tag{2.5-26}$$

$$x_3 = x_2 + \Delta x \tag{2.5-27}$$

$$x_2 = x_3 - \Delta x \tag{2.5-28}$$

Substituting Eq 2.5-27 and Eq 2.5-28 into Eq 2.5-26

$$\mathrm{lnd}(n,\Delta x,x_1) - \mathrm{lnd}(n,-\Delta x,x_1-\Delta x) = \mathrm{lnd}(n,\Delta x,x_3) - \mathrm{lnd}(n,-\Delta x,x_3-\Delta x) \tag{2.5-29}$$

Since $x_1$ and $x_3$ are any points on the x locus and the left and right sides of Eq 2.5-29 are equal,

Then

$$\mathrm{lnd}(n,\Delta x,x) - \mathrm{lnd}(n,-\Delta x,x-\Delta x) = \text{constant} \tag{2.5-30}$$

where

$x = x_i + m\Delta x$ ,  m = integers ,  This is the x locus

$x_i$ = a value of x

$\Delta x$ = x increment

$n,\Delta x,x.x_i$ = real or complex values

Note

1) The quantity, lnd(n,Δx,x) – lnd(n,-Δx,x-Δx), equals the same value for any value of x on an x locus line where n is a contant value
2) The constant value remains the same for a specific n and x locus line
3) The relationship of Eq 2.5-30 is valid for all values of n including the special case where n=1.

Let

$\Delta x = \Delta\mathbf{x} \quad (Re(\Delta x)>0 \text{ or } \{Re(\Delta x)=0 \text{ and } Im(\Delta x)>0\})$ (2.5-31)

From Eq 2.5-30 and Eq 2.5-31

$lnd(n,\Delta\mathbf{x},x) - lnd(n,-\Delta\mathbf{x},x-\Delta\mathbf{x}) = K$ (2.5-32)

where

K = real or complex value constant

Compare Eq 2.5-32 to Eq 2.5-12

$lnd(n,\Delta\mathbf{x},x) - lnd(n,-\Delta\mathbf{x},x-\Delta\mathbf{x}) = K$ (2.5-33)

where

K = constant

$lnd(n,\Delta\mathbf{x},x) - lnd(n,-\Delta\mathbf{x},x-\Delta\mathbf{x}) \approx K_{r-}k_r \;,\quad n\neq 1 \;,\quad r = 1,2,3$ (2.5-34)

where

Accuracy increases rapidly as $|\frac{x}{\Delta x}|$ increases in value.

For $|\frac{x}{\Delta x}|$ large, $lnd(n, \Delta\mathbf{x},x) - lnd(n,- \Delta\mathbf{x},x-\Delta\mathbf{x}) = K_{r-}k_r$ then $K = K_{r-}k_r$

$K = K_{r-}k_r$ (2.5-35)

Substituting Eq 2.5-35 into Eq 2.5-32

$lnd(n,\Delta\mathbf{x},x) - lnd(n,-\Delta\mathbf{x},x-\Delta\mathbf{x}) = K_{r-}k_r$ (2.5-36)

Let

$\Delta x = -\Delta\mathbf{x} \quad (Re(\Delta x)<0 \text{ or } \{Re(\Delta x)=0 \text{ and } Im(\Delta x)<0\})$ (2.5-37)

From Eq 2.5-30 and Eq 2.5-37

$lnd(n,-\Delta\mathbf{x},x) - lnd(n,\Delta\mathbf{x},x+\Delta\mathbf{x}) = c$ (2.5-38)

where

c = real or complex value constant

Compare Eq 2.5-38 to Eq 2.5-18

$$\mathrm{lnd}(n,-\Delta x,x) - \mathrm{lnd}(n,\Delta x,x+\Delta x) = c \qquad (2.5\text{-}39)$$

where

c = constant

and

$$\mathrm{lnd}(n,-\Delta x,x) - \mathrm{lnd}(n,\Delta x,x+\Delta x) \approx -(K_{r} - r_{r}),\quad n\neq 1\ ,\ r = 1,2,3 \qquad (2.5\text{-}40)$$

Accuracy increases rapidly as $|\frac{x}{\Delta x}|$ increases in value.

For $|\frac{x}{\Delta x}|$ large, $\mathrm{lnd}(n,-\Delta x,x) - \mathrm{lnd}(n,\Delta x,x+\Delta x) = -(K_{r} - k_{r})$ then $c = -(K_{r} - k_{r})$

$$c = -(K_{r} - k_{r}) \qquad (2.5\text{-}41)$$

Substituting Eq 2.5-41 into Eq 2.5-38

$$\mathrm{lnd}(n,-\Delta x,x) - \mathrm{lnd}(n,\Delta x,x+\Delta x) = -(K_{r} - k_{r}) \qquad (2.5\text{-}42)$$

Rewriting Eq 2.5-36 and Eq 2.5-42

$$\mathrm{lnd}(n,[\Delta x],x) - \mathrm{lnd}(n,-[\Delta x],x-[\Delta x]) = +(K_{r} - k_{r})\quad ,\quad n\neq 1\ ,\ r = 1,2,3 \qquad (2.5\text{-}43)$$

$$\mathrm{lnd}(n,[-\Delta x],x) - \mathrm{lnd}(n,-[-\Delta x],x-[-\Delta x]) = -(K_{r} - r_{r}),\quad n\neq 1\ ,\ r = 1,2,3 \qquad (2.5\text{-}44)$$

It is possible to merge Eq 2.5-43 and Eq 2.5-44 into one equation. The set of values of $\Delta x$ is equal to the set of values of $\Delta x$ and $-\Delta x$. This equation is shown in Eq 2.5-45 below.

$$\mathbf{lnd(n,\Delta x,x) - lnd(n,-\Delta x,x-\Delta x) = \pm(K_{r} - k_{r})} \qquad (2.5\text{-}45)$$

**where**

**$x = x_i + m\Delta x$ , m = integers , This is an x locus in the complex plane**

**$x_i$ = a value of x**

**$\Delta x$ = x increment**

**+ for $Re(\Delta x)>0$ or $\{Re(\Delta x)=0$ and $Im(\Delta x)>0\}$**

**– for $Re(\Delta x)<0$ or $\{Re(\Delta x)=0$ and $Im(\Delta x)<0\}$**

**x = real or complex values**

**$n\neq 1$**

**r = 1,2,3**

**$n,\Delta x,x_i,K_r,k_r$ = real or complex constants**

**The x locus line is the straight line in the complex plane through all of the plotted x points**

From Eq 2.5-43 and Eq 2.5-44

$$\mathrm{lnd}(n,[\Delta x],x) - \mathrm{lnd}(n,-[\Delta x],x-[\Delta x]) = -[\mathrm{lnd}(n,[-\Delta x],x) - \mathrm{lnd}(n,-[-\Delta x],x-[-\Delta x])] \qquad (2.5\text{-}46)$$

Rearranging terms

$[\text{lnd}(n,[-\Delta x],x) - \text{lnd}(n,-[-\Delta x],x-[-\Delta x])] = -[\text{lnd}(n,[\Delta x],x) - \text{lnd}(n,-[\Delta x],x-[\Delta x])]$ (2.5-47)

It is possible to merge Eq 2.5-46 and Eq 2.5-47 into one equation. The set of values of $\Delta x$ is equal to the set of values of $\Delta x$ and $-\Delta x$. This equation is shown in Eq 2.5-48 below.

**$\text{lnd}(n,\Delta x,x) - \text{lnd}(n,-\Delta x,x-\Delta x) = -[\text{lnd}(n,-\Delta x,x) - \text{lnd}(n,\Delta x,x+\Delta x)]$** (2.5-48)

**where**

**$x = x_i + m\Delta x$ , $m$ = integers , This is an x locus in the complex plane**

**$x_i$ = a value of x**

**$\Delta x$ = x increment**

**x = real or complex values**

**$n,\Delta x,x_i,K_r,k_r$ = real or complex constants**

**The x locus line is the straight line in the complex plane through all of the plotted x points**

From the previously derived equations other useful relationships can be obtained. Some of these relationships are as follows:

1) $K_1$ is known to equal zero.

$K_1 = 0$ (2.5-49)

Then from Eq 2.5-45 and Eq 2.5-49

**$k_1 = \pm[\text{lnd}(n,\Delta x,x) - \text{lnd}(n,-\Delta x,x-\Delta x)], \ n\neq 1$** (2.5-50)

**where**

**$x = x_i + m\Delta x$ , $m$ = integers , This is an x locus in the complex plane**

**$x_i$ = a value of x**

**$\Delta x$ = x increment**

**– for $Re(\Delta x)>0$ or $\{Re(\Delta x)=0$ and $Im(\Delta x)>0\}$**

**+ for $Re(\Delta x)<0$ or $\{Re(\Delta x)=0$ and $Im(\Delta x)<0\}$**

**x = real or complex values**

**$n,\Delta x,x_i,K_r,k_r$ = real or complex constants**

**The x locus line is the straight line in the complex plane through all of the plotted x points**

2) $k_3$ is known to equal zero.

$k_3 = 0$ (2.5-51)

Then from Eq 2.5-45 and Eq 2.5-51

**$K_3 = \pm[\text{lnd}(n,\Delta x,x) - \text{lnd}(n,-\Delta x,x-\Delta x)], \ n\neq 1$** (2.5-52)

**where**

**$x = x_i + m\Delta x$ ,    m = integers  ,   This is an x locus in the complex plane**

**$x_i$ = a value of x**

**$\Delta x$ = x increment**

**+ for $Re(\Delta x)>0$ or $\{Re(\Delta x)=0$ and $Im(\Delta x)>0\}$**

**– for $Re(\Delta x)<0$ or $\{Re(\Delta x)=0$ and $Im(\Delta x)<0\}$**

**x = real or complex values**

**$n,\Delta x,x_i,K_r,k_r$ = real or complex constants**

**The x locus line is the straight line in the complex plane through all of the plotted x points.**

3) From Eq 2.5-50 and Eq 2.5-52 a relationship between $K_3$ and $k_1$ is found

$$\mathbf{K_3 = -k_1} \tag{2.5-53}$$

4) Finding another relationship involving the function, $\text{lnd}(n,\Delta x,x) - \text{lnd}(n,-\Delta x,x-\Delta x)$

Find a relationship for $\displaystyle {}_{\Delta x}\!\sum_{x=-\infty}^{+\infty} \frac{1}{x^n}$

Where

$x = x_i + m\Delta x$ ,  m = integers

$x_i$ = a value of x

$\Delta x$ = x increment

$Re(n)>1$

$\Delta x$ is $\Delta x$, the summation x increment, where $Re(\Delta x)>0$ or $\{Re(\Delta x)=0$ and $Im(\Delta x)>0\}$.

The summation proceeds along its complex plane x locus line from left to right or down to up.

$$\Delta x \, {}_{\Delta x}\!\sum_{x=x_i}^{+\infty} \frac{1}{x^n} = \text{lnd}(n,\Delta x,x_i) = \Delta x \left[ \frac{1}{x_i^n} + \frac{1}{(x_i+\Delta x)^n} + \frac{1}{(x_i+2\Delta x)^n} + \frac{1}{(x_i+3\Delta x)^n} + \ldots \right] \tag{2.5-54}$$

Eq 2.5-54 sums the function, $\frac{1}{x^n}$ , where $x = x_i,\ x_i+\Delta x,\ x_i+2\Delta x,\ x_i+3\Delta x,\ x_i+4\Delta x,\ \ldots$

$$-\Delta x \, {}_{-\Delta x}\!\sum_{x=x_i-\Delta x}^{-\infty} \frac{1}{x^n} = \text{lnd}(n,-\Delta x,x_i-\Delta x) = -\Delta x \left[ \ldots + \frac{1}{(x_i-3\Delta x)^n} + \frac{1}{(x_i-2\Delta x)^n} + \frac{1}{(x_i-\Delta x)^n} \right] \tag{2.5-55}$$

Eq 2.5-55 sums the function, $\frac{1}{x^n}$ , where $x = x_i-\Delta x,\ x_i-2\Delta x,\ x_i-3\Delta x,\ x_i-4\Delta x,\ \ldots$

Subtracting Eq 2.5-56 from Eq 2.5-55

$$\Delta x \sum_{\substack{\Delta x \\ x=-\infty}}^{+\infty} \frac{1}{x^n} = \text{lnd}(n,\Delta x,x_i) - \text{lnd}(n,-\Delta x,x_i-\Delta x) \qquad (2.5\text{-}56)$$

where
$x = x_i + m\Delta x$ , m = integers
$x_i$ = a value of x
$\Delta x$ = x increment
Re(n) > 1
$x,x_i,\Delta x,n$ = real or complex values

Note that the summation of Eq 2.5-56 proceeds along its complex plane x locus line from left to right or down to up.

Find a relationship for $\displaystyle \sum_{\substack{-\Delta x \\ x=+\infty}}^{-\infty} \frac{1}{x^n}$

where
$x = x_i + m\Delta x$ , m = integers
$x_i$ = a value of x
$\Delta x$ = x increment
Re(n)>1
$-\Delta x$ is $\Delta x$, the summation x increment, where $Re(\Delta x)<0$ or{$Re(\Delta x)=0$ and $Im(\Delta x)<0$}.
The summation proceeds along its complex plane x locus line from right to left or up to down.

$$-\Delta x \sum_{\substack{-\Delta x \\ x=x_i}}^{-\infty} \frac{1}{x^n} = \text{lnd}(n,-\Delta x,x_i) = -\Delta x \left[ \frac{1}{x_i^n} + \frac{1}{(x_i-\Delta x)^n} + \frac{1}{(x_i-2\Delta x)^n} + \frac{1}{(x_i-3\Delta x)^n} + \ \ldots\right] \qquad (2.5\text{-}57)$$

Eq 2.5-57 sums the function, $\frac{1}{x^n}$ , where $x = x_i, x_i-\Delta x, x_i-2\Delta x, x_i-3\Delta x, x_i-4\Delta x, \ldots$

$$\Delta x \sum_{\substack{\Delta x \\ x=x_i+\Delta x}}^{+\infty} \frac{1}{x^n} = \text{lnd}(n,\Delta x,x_i+\Delta x) = \Delta x \left[\ldots + \frac{1}{(x_i+\Delta x)^n} + \frac{1}{(x_i+2\Delta x)^n} + \frac{1}{(x_i+3\Delta x)^n}\right] \qquad (2.5\text{-}58)$$

Eq 2.5-58 sums the function, $\frac{1}{x^n}$ , where $x = x_i+\Delta x, x_i+2\Delta x, x_i+3\Delta x, x_i+4\Delta x, \ldots$

Subtracting Eq 2.5-58 from Eq 2.5-57

$$-\Delta x \sum_{\substack{-\Delta x \\ x=+\infty}}^{-\infty} \frac{1}{x^n} = \mathrm{lnd}(n,-\Delta x,x_i) - \mathrm{lnd}(n,\Delta x,x_i+\Delta x) \qquad (2.5\text{-}59)$$

where
$x = x_i + m\Delta x$ ,   m = integers
$x_i$ = a value of x
$\Delta x$ = x increment
Re(n) > 1
$x,x_i,\Delta x,n$ = real or complex values

Note that the summation of Eq 2.5-59 proceeds along its complex plane x locus line from right to left or up to down.

Rewriting Eq 2.5-56 and Eq 2.5-59

$$[\Delta x] \sum_{\substack{\Delta x \\ x=-\infty}}^{+\infty} \frac{1}{x^n} = \mathrm{lnd}(n,[\Delta x],x_i) - \mathrm{lnd}(n,-[\Delta x],x_i-[\Delta x]) \qquad (2.5\text{-}60)$$

$$[-\Delta x] \sum_{\substack{-\Delta x \\ x=+\infty}}^{-\infty} \frac{1}{x^n} = \mathrm{lnd}(n,[-\Delta x],x_i) - \mathrm{lnd}(n,\Delta x,x_i-[-\Delta x]) \qquad (2.5\text{-}61)$$

Also,

$$\sum_{\substack{\Delta x \\ x=-\infty}}^{+\infty} \frac{1}{x^n} = \sum_{\substack{-\Delta x \\ x=+\infty}}^{-\infty} \frac{1}{x^n} \quad , \quad x = x_i + m\Delta x \ , \quad m = \text{integers} \qquad (2.5\text{-}62)$$

It is possible using Eq 2.5-62 to merge Eq 2.5-60 and Eq 2.5-61 into one equation. The set of values of $\Delta x$ is equal to the set of values of $\Delta x$ and $-\Delta x$. This equation is shown in Eq 2.5-63 below.

$$\mathbf{\Delta x \sum_{\substack{\Delta x \\ x=\mp\infty}}^{\pm\infty} \frac{1}{x^n} = lnd(n,\Delta x,x_i) - lnd(n,-\Delta x,x_i-\Delta x)} \qquad (2.5\text{-}63)$$

**or**

$$\sum_{\substack{\Delta x \\ x=\mp\infty}}^{\pm\infty} \frac{1}{x^n} = \frac{1}{\Delta x}\left[\text{lnd}(n,\Delta x,x_i) - \text{lnd}(n,-\Delta x,x_i-\Delta x)\right] \tag{2.5-64}$$

**where**

**$x = x_i + m\Delta x$ , m = integers , This is an x locus in the complex plane**

**$x_i$ = a value of x**

**$\Delta x$ = x increment**

**$-\infty$ to $+\infty$ for $Re(\Delta x)>0$ or $\{Re(\Delta x)=0$ and $Im(\Delta x)>0\}$**

**$+\infty$ to $-\infty$ for $Re(\Delta x)<0$ or $\{Re(\Delta x)=0$ and $Im(\Delta x)<0\}$**

**x = real or complex values**

**$n,\Delta x,x_i,K_r,k_r$ = real or complex constants**

**$Re(n)>1$**

**The x locus line is the straight line in the complex plane through all of the plotted x points.**

From Eq 2.5-45 with $x = x_i$

$$\text{lnd}(n,\Delta x,x_i) - \text{lnd}(n,-\Delta x,x_i-\Delta x) = \pm(K_r - k_r),\quad n\neq 1\ ,\ r = 1,2,3 \tag{2.5-65}$$

Substituting Eq 2.5-65 into Eq 2.5-63

$$\pm(K_r - k_r) = \text{lnd}(n,\Delta x,x_i) - \text{lnd}(n,-\Delta x,x_i-\Delta x) = \Delta x \sum_{\substack{\Delta x \\ x=\mp\infty}}^{\pm\infty} \frac{1}{x^n} = \text{a constant irrespective of } r \tag{2.5-66}$$

**where**

**$x = x_i + m\Delta x$ , m = integers , This is an x locus in the complex plane**

**$x_i$ = a value of x**

**$\Delta x$ = x increment**

**+ and $-\infty$ to $+\infty$ for $Re(\Delta x)>0$ or $\{Re(\Delta x)=0$ and $Im(\Delta x)>0\}$**

**– and $+\infty$ to $-\infty$ for $Re(\Delta x)<0$ or $\{Re(\Delta x)=0$ and $Im(\Delta x)<0\}$**

**x = real or complex values**

**$n,\Delta x,x,x_i,K_r,k_r$ = real or complex constants**

**$Re(n)>1$**

**r = 1,2,3**

**The x locus line is the straight line in the complex plane through all of the plotted x points.**

Notes - x can be represented by a linear x locus and an x locus line in the complex plane.

$\sum_{\substack{\Delta x \\ x=\mp\infty}}^{\pm\infty} \frac{1}{x^n}$ is the same value for any value of $x_i$, and $\Delta x$ on an x locus line.

$$\Delta x \sum_{\Delta x \; x=\mp\infty}^{\pm\infty} \frac{1}{x^n} = \Delta x \left[\ldots + \frac{1}{(x_i-2\Delta x)^n} + \frac{1}{(x_i-\Delta x)^n} + \frac{1}{x_i^n} + \frac{1}{(x_i+\Delta x)^n} + \frac{1}{(x_i+2\Delta x)^n} + \ldots\right] \qquad (2.5\text{-}67)$$

where

$x = x_i + m\Delta x$ ,  m = integers
$x_i$ = a value of x
$\Delta x$ = x increment
Re(n) > 1
$x, x_i, \Delta x, n$ = real or complex values

**Though the lnd(n,Δx,x) n≠1 constants of integration, $K_1, K_2, K_3, k_1, k_2, k_3$, may change along an x locus line in the complex plane, the value of $K_r - k_r = \pm[\text{lnd}(n,\Delta x,x) - \text{lnd}(n,-\Delta x,x-\Delta x)]$ will not change. This non-varying value will also be $K_r - k_r = \pm \left[\Delta x \sum_{\Delta x \; x=\mp\infty}^{\pm\infty} \frac{1}{x^n}\right]$.**

**where**

**$x = x_i + m\Delta x$ ,  m = integers ,  This is an x locus in the complex plane**
**$x_i$ = a value of x**
**Δx = x increment**
**+ and -∞ to +∞ for Re(Δx)>0 or {Re(Δx)=0 and Im(Δx)>0}**
**– and +∞ to -∞ for Re(Δx)<0 or {Re(Δx)=0 and Im(Δx)<0}**
**x = real or complex values**
**r = 1,2,3 , The x locus segment designations**
**Re(n)>1**

<u>Note</u> - $\sum_{\Delta x \; x=-\infty}^{+\infty} \frac{1}{x^n}$ and $\sum_{-\Delta x \; x=+\infty}^{-\infty} \frac{1}{x^n}$ have the same value.

From Eq 2.5-66

$$\pm(K_1 - k_1) = \pm(K_2 - k_2) \pm = \pm(K_3 - k_3) = \Delta x \sum_{\Delta x \; x=\mp\infty}^{\pm\infty} \frac{1}{x^n} = \text{lnd}(n,\Delta x,x_i) - \text{lnd}(n,-\Delta x,x_i-\Delta x) \qquad (2.5\text{-}68)$$

where

$x = x_i + m\Delta x$ ,  m = integers ,  This is an x locus in the complex plane
$x_i$ = a value of x
Δx = x increment
+ and -∞ to +∞ for Re(Δx)>0 or {Re(Δx)=0 and Im(Δx)>0}

– and +∞ to -∞ for Re(Δx)<0 or {Re(Δx)=0 and Im(Δx)<0}
x = real or complex values
$n,\Delta x,x_i,K_r,k_r$ = real or complex constants
Re(n)>1
r = 1,2,3

The x locus line is the straight line in the complex plane through all of the plotted x points.

From Eq 2.5-68 and Eq 2.5-49

Substituting Eq 2.5-49, $K_1 = 0$, into Eq 2.5-68

$$\pm(0 - k_1) = \Delta x \sum_{\substack{x=\mp\infty \\ \Delta x}}^{\pm\infty} \frac{1}{x^n}$$

$$\mathbf{k_1 = \pm \left[-\Delta x \sum_{\substack{x=\mp\infty \\ \Delta x}}^{\pm\infty} \frac{1}{x^n}\right] \ , \quad Re(n)>1} \qquad (2.5\text{-}69)$$

**where**

**$x = x_i + m\Delta x$ , m = integers , This is an x locus in the complex plane**
**$x_i$ = a value of x**
**Δx = x increment**
**+ and -∞ to +∞ for Re(Δx)>0 or {Re(Δx)=0 and Im(Δx)>0}**
**– and +∞ to -∞ for Re(Δx)<0 or {Re(Δx)=0 and Im(Δx)<0}**
**x = real or complex values**
**$n,\Delta x,x_i,K_r,k_r$ = real or complex constants**
**r = 1,2,3 , The x locus segment designations**
**Re(n)>1**

**The x locus line is the straight line in the complex plane through all of the plotted x points.**

From Eq 2.5-68 and Eq 2.5-51

Substituting Eq 2.5-51, $k_3 = 0$, into Eq 2.5-68

$$\pm(K_3 - 0) = \Delta x \sum_{\substack{x=\mp\infty \\ \Delta x}}^{\pm\infty} \frac{1}{x^n}$$

$$\mathbf{K_3 = \pm \left[\Delta x \sum_{\substack{x=\mp\infty \\ \Delta x}}^{\pm\infty} \frac{1}{x^n}\right] \ , \quad Re(n)>1}$$

**where** (2.5-70)

**$x = x_i + m\Delta x$ ,  m = integers ,  This is an x locus in the complex plane**
**$x_i$ = a value of x**
**$\Delta x$ = x increment**
**+ and $-\infty$ to $+\infty$ for $Re(\Delta x)>0$ or $\{Re(\Delta x)=0$ and $Im(\Delta x)>0\}$**
**– and $+\infty$ to $-\infty$ for $Re(\Delta x)<0$ or $\{Re(\Delta x)=0$ and $Im(\Delta x)<0\}$**
**x = real or complex values**

**$n,\Delta x,x_i,K_r,k_r$ = real or complex constants**
**r = 1,2,3 ,  The x locus segment designations**

From Eq 2.5-68

$$K_2 - k_2 = \pm\left[\Delta x \sum_{\Delta x \; x=\mp\infty}^{\pm\infty} \frac{1}{x^n}\right] \;,\quad Re(n)>1 \tag{2.5-71}$$

**where**

**$x = x_i + m\Delta x$ ,  m = integers ,  This is an x locus in the complex plane**
**$x_i$ = a value of x**
**$\Delta x$ = x increment**
**+ and $-\infty$ to $+\infty$ for $Re(\Delta x)>0$ or $\{Re(\Delta x)=0$ and $Im(\Delta x)>0\}$**
**– and $+\infty$ to $-\infty$ for $Re(\Delta x)<0$ or $\{Re(\Delta x)=0$ and $Im(\Delta x)<0\}$**
**x = real or complex values**
**$n,\Delta x,x_i,K_r,k_r$ = real or complex constants**
**r = 1,2,3 ,  The x locus segment designations**

Diagram 2.5-6 below is an example of an Eq 2.5-68 x locus line in the complex plane.

<u>Diagram 2.5-6: An example of an x locus line in the complex plane</u>

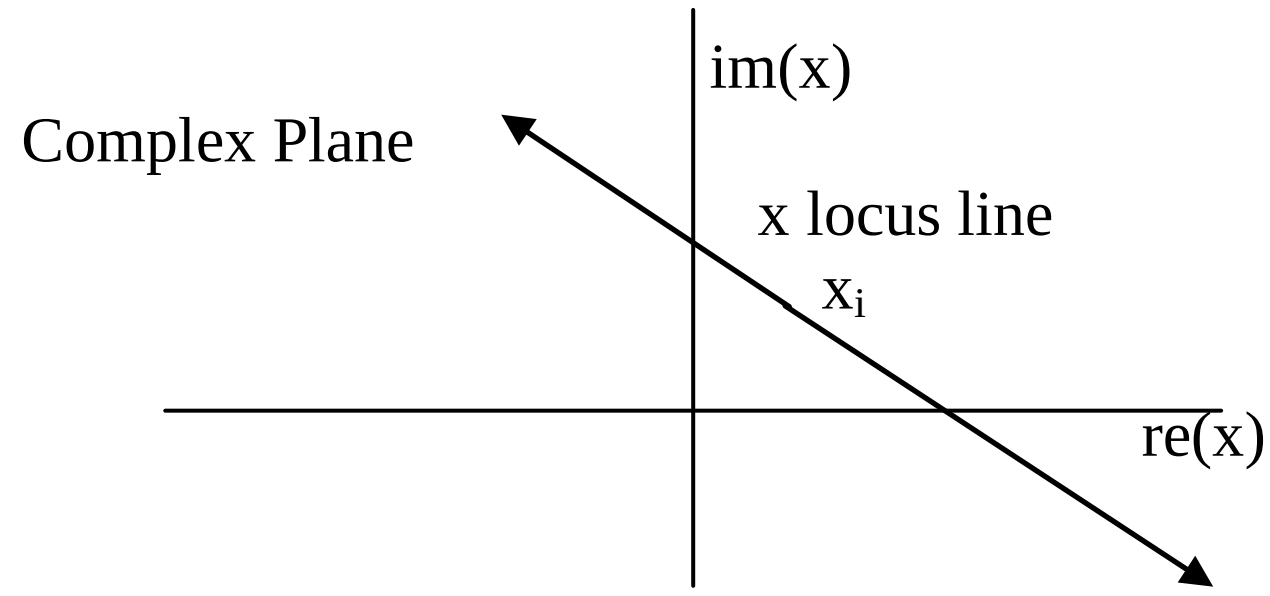

5) The $lnd(n,\Delta x,x)$ $n\neq1$ constant of integration can be obtained from Eq 2.5-1

Rewriting Eq 2.5-1

$$\ln d(n,\Delta x,x) \approx -\sum_{m=0}^{\infty} \frac{\Gamma(n+2m-1)\left(\frac{\Delta x}{2}\right)^{2m} C_m}{\Gamma(n)(2m+1)!\left(x - \frac{\Delta x}{2}\right)^{n+2m-1}} + \begin{cases} K_r \\ k_r \end{cases} , \quad n \neq 1 \qquad (2.5\text{-}72)$$

$K_r$ = constants of integration for $Re(\Delta x) > 0$ or $\{Re(\Delta x) = 0$ and $Im(\Delta x) > 0\}$
$k_r$ = constants of integration for $Re(\Delta x) < 0$ or $\{Re(\Delta x) = 0$ and $Im(\Delta x) < 0\}$
r = 1, 2, or 3 , The x locus segment designations
Accuracy increases rapidly as $|\frac{x}{\Delta x}|$ increases in value

Let $\ln d_f(n,\Delta x,x)$ represent the above summation only

$$\ln d_f(n,\Delta x,x) = -\sum_{m=0}^{\infty} \frac{\Gamma(n+2m-1)\left(\frac{\Delta x}{2}\right)^{2m} C_m}{\Gamma(n)(2m+1)!\left(x - \frac{\Delta x}{2}\right)^{n+2m-1}} , \quad n \neq 1 \qquad (2.5\text{-}73)$$

Then from Eq 2.5-72 and Eq 2.5-73

$$\left.\begin{matrix} \mathbf{K_r} \\ \mathbf{k_r} \end{matrix}\right\} \approx \mathbf{\ln d(n,\Delta x,x) - \ln d_f(n,\Delta x,x)} \qquad (2.5\text{-}74)$$

**x is in the x locus segment, r**

**where**

**$x = x_i + m\Delta x$ , m = integer , This is an x locus in the complex plane**
**$x_i$ = value of x**
**x = real or complex values**
**$K_r$ = constants of integration for $Re(\Delta x) > 0$ or $Re(\Delta x) = 0$ and $Im(\Delta x) > 0$**
**$k_r$ = constants of integration for $Re(\Delta x) < 0$ or $Re(\Delta x) = 0$ and $Im(\Delta x) < 0$**
**r = 1, 2, or 3 , The x locus segment designations**
**$n \neq 1$**

**The absolute value, $|\frac{x}{\Delta x}|$ , must be large for good accuracy**

This equation is only useful if $\ln d(n,\Delta x,x)$ is known.

Merging Eq 2.5-3 and Eq 2.5-5

$$\Delta x \sum_{\substack{\Delta x \\ x=x_i}}^{\pm\infty} \frac{1}{x^n} = \text{lnd}(n,\Delta x,x_i) \tag{2.5-75}$$

where

$x = x_i + m\Delta x$

$x_i$ = initial value of x

$\Delta x$ = x increment

$m = 0,1,2,3,\ldots$ and $x_i$ to $+\infty$ for $Re(\Delta x)>0$ or $\{Re(\Delta x)=0$ and $Im(\Delta x)>0\}$

$m = 0,-1,-2,-3,\ldots$ and $x_i$ to $-\infty$ for $Re(\Delta x)<0$ or $\{Re(\Delta x)=0$ and $Im(\Delta x)<0\}$

$Re(n)>1$

Substituting Eq 2.5-75 into Eq 2.5-74

$$\left.\begin{matrix} \mathbf{K_r} \\ \mathbf{k_r} \end{matrix}\right\} = \Delta x \sum_{\substack{\Delta x \\ x=x_i}}^{\pm\infty} \frac{1}{x^n} - \text{lnd}_f(n,\Delta x,x_i) \ , \quad Re(n)>1$$

**$x_i$ is in the x locus segment, r**

**where**

**$x = x_i + m\Delta x$, This is an x locus in the complex plane**

**$x_i$ = initial value of x**

**$\Delta x$ = x increment**

**$m = 0,1,2,3,\ldots$ and $x_i$ to $+\infty$ for $Re(\Delta x)>0$ or $\{Re(\Delta x)=0$ and $Im(\Delta x)>0\}$**

**$m = 0,-1,-2,-3,\ldots$ and $x_i$ to $-\infty$ for $Re(\Delta x)<0$ or $\{Re(\Delta x)=0$ and $Im(\Delta x)<0\}$**

**$K_r$ = constant of integration for $Re(\Delta x)>0$ or $\{Re(\Delta x) = 0$ and $Im(\Delta x) > 0\}$**

**$k_r$ = constant of integration for $Re(\Delta x)<0$ or $\{Re(\Delta x) = 0$ and $Im(\Delta x) < 0\}$**

**r = 1, 2, or 3 , x locus segment designation**

**$n \neq 1$**

**Accuracy increases rapidly as $\left|\frac{x_i}{\Delta x}\right|$ increases in value.**

Another method for finding the constants of integration will be presented in Section 2.7.

6) Eq 2.5-45 provides a necessary conditon for no lnd(n,Δx,x) n≠1 function snap to occur for all x.

Rewriting Eq 2.5-45

$$\text{lnd}(n,\Delta x,x) - \text{lnd}(n,-\Delta x,x-\Delta x) = \pm(K_r - k_r)$$

where

$x = x_i + m\Delta x$ , m = integers , This is an x locus in the complex plane

$x_i$ = a value of x
$\Delta x$ = x increment
+ for $Re(\Delta x)>0$ or $\{Re(\Delta x)=0$ and $Im(\Delta x)>0\}$
– for $Re(\Delta x)<0$ or $\{Re(\Delta x)=0$ and $Im(\Delta x)<0\}$
x = real or complex values
$n\neq 1$ ,
r = 1,2,3
$n,\Delta x,x_i,K_r,k_r$ = real or complex constants

The x locus line is the straight line in the complex plane through all of the plotted x points.

If the value of $lnd(n,\Delta x,x) - lnd(n,-\Delta x,x-\Delta x) = 0$ for an x value on an x locus line, and knowing $K_1=k_3=0$, then $K_1=k_1=K_3=k_3=0$ and $K_2=k_2$. When $K_1,K_2,K_3,k_1,k_2,k_3$ all equal 0 (no function snap) this condition is satisfied.

Then

**A necessary condition that $lnd(n,\Delta x,x)$ $n\neq 1$ Series snap not occur for all x on an x locus line is that $lnd(n,\Delta x,x) - lnd(n,-\Delta x,x-\Delta x) = 0$ for an x on that x locus line.**

7) Though there is no formal proof, it appears that for n = -1,-2,-3,… $lnd(n,\Delta x,x) - lnd(n,\Delta x,x-\Delta x) = 0$.

Rewriting the $lnd(n,\Delta x,x)$ $n\neq 1$ Series constant of integration evaluation equations previously derived.

**The following equations can be used to find the $lnd(n,\Delta x,x)$ $n\neq 1$ function constants of integration, $K_1,K_2,K_3,k_1,k_2,k_3$.**

**<u>For all values of R(n)</u>**

$$K_1 = k_3 = 0 \tag{2.5-76}$$
$$K_3 = -k_1 = \pm[lnd(n,\Delta x,x_i) - lnd(n,-\Delta x,x_i-\Delta x)] \tag{2.5-77}$$
$$K_2 - k_2 = \pm[lnd(n,\Delta x,x_i) - lnd(n,-\Delta x,x_i-\Delta x)] \tag{2.5-78}$$
$$K_r - k_r = \pm[lnd(n,\Delta x,x_i) - lnd(n,-\Delta x,x_i-\Delta x)] \;,\quad r = 1,2,\text{ or }3 \tag{2.5-79}$$

$$\left.\begin{matrix}K_r\\ k_r\end{matrix}\right\}\approx lnd(n,\Delta x,x) - lnd_f(n,\Delta x,x) \tag{2.5-80}$$

**x is in the x locus segment, r**

**Accuracy increases rapidly as $|\frac{x_i}{\Delta x}|$ increases in value.**

**where**
**$x = x_i + m\Delta x$ , m = integer , This is an x locus in the complex plane**
**$x_i$ = value of x**
**$\Delta x$ = x increment**

**$x$ = real or complex values**
**+ for Re(Δx)>0 or {Re(Δx)=0 and Im(Δx)>0}**
**– for Re(Δx)<0 or {Re(Δx)=0 and Im(Δx)<0}**
**$K_r$ = constants of integration for Re(Δx) > 0 or Re(Δx) = 0 and Im(Δx) > 0**
**$k_r$ = constants of integration for Re(Δx) < 0 or Re(Δx) = 0 and Im(Δx) < 0**
**r = 1, 2, or 3 , The x locus segment designations**
**$n \neq 1$**
**$n,\Delta x,x_i,K_r,k_r$ = real or complex constants**

**The x locus line is the straight line in the complex plane through all of the plotted x points.**

**A necessary condition that lnd(n,Δx,x) n≠1 Series snap not occur for all x on an x locus line is that lnd(n,Δx,x) – lnd(n,-Δx,x-Δx) = 0 for an x on that x locus line.**

Comment

To find the values of the constants $K_2$,$K_3$,$k_1$, and $k_2$, using these equations, the lnd(n,Δx,x) n≠1 function must first be calculated. However, in order to calculate the lnd(n,Δx,x) n≠1 function these same constant values must be known. As a result of this difficulty, Eq 2.5-77, Eq 2.5-78, and Eq 2.5-80 have limited use. The values of these four constants can be found in another way. This different calculation method will be presented in Section 2.8.

**For Re(n)>1**

$$K_1 = k_3 = 0 \tag{2.5-81}$$

$$k_1 = \pm\left[-\Delta x \sum_{\Delta x\ x=\mp\infty}^{\pm\infty} \frac{1}{x^n}\right] \tag{2.5-82}$$

$$K_3 = \pm\left[\Delta x \sum_{\Delta x\ x=\mp\infty}^{\pm\infty} \frac{1}{x^n}\right] \tag{2.5-83}$$

$$K_2 - k_2 = \pm\left[\Delta x \sum_{\Delta x\ x=\mp\infty}^{\pm\infty} \frac{1}{x^n}\right] \tag{2.5-84}$$

$$K_r - k_r = \pm\left[\Delta x \sum_{\Delta x\ x=\mp\infty}^{\pm\infty} \frac{1}{x^n}\right]\ , \quad r = 1,2,\text{ or }3 \tag{2.5-85}$$

**where**
**$x = x_i + m\Delta x$ , m = integers , This is an x locus in the complex plane**
**$x_i$ = a value of x**

**$\Delta x$ = x increment**
**+ and -∞ to +∞ for $Re(\Delta x)>0$ or $\{Re(\Delta x)=0$ and $Im(\Delta x)>0\}$**
**– and +∞ to -∞ for $Re(\Delta x)<0$ or $\{Re(\Delta x)=0$ and $Im(\Delta x)<0\}$**
**x = real or complex values**
**$n,\Delta x,x_i,K_r,k_r$ = real or complex constants**
**r = 1,2,3 , The x locus segment designations**
**$Re(n)>1$**

**The x locus line is the straight line in the complex plane through all of the plotted x points.**

$$\left.\begin{matrix} K_r \\ k_r \end{matrix}\right\} = \Delta x \sum_{\Delta x \; x=x_i}^{\pm\infty} \frac{1}{x^n} - \mathrm{lnd}_f(n,\Delta x,x_i) \qquad (2.5\text{-}86)$$

**x is in the x locus segment, r**

**where**

**$x = x_i + m\Delta x$, This is an x locus in the complex plane**
**$x_i$ = initial value of x**
**$\Delta x$ = x increment**
**$m = 0,1,2,3,\ldots$ and $x_i$ to +∞ for $Re(\Delta x)>0$ or $\{Re(\Delta x)=0$ and $Im(\Delta x)>0\}$**
**$m = 0,-1,-2,-3,\ldots$ and $x_i$ to -∞ for $Re(\Delta x)<0$ or $\{Re(\Delta x)=0$ and $Im(\Delta x)<0\}$**
**$K_r$ = constant of integration for $Re(\Delta x)>0$ or $\{Re(\Delta x) = 0$ and $Im(\Delta x) > 0\}$**
**$k_r$ = constant of integration for $Re(\Delta x)<0$ or $\{Re(\Delta x) = 0$ and $Im(\Delta x) < 0\}$**
**r = 1, 2, or 3 , x locus segment designation**
**$n \neq 1$**
**$Re(n)>1$**

**Accuracy increases rapidly as $\left|\frac{x_i}{\Delta x}\right|$ increases in value.**

## Section 2.6: Evaluation of the constants of integration, $K_1$ and $k_3$

Constant of integration, $K_1$

The relationship for the $\mathrm{lnd}(n,\Delta x,x)$ $n\neq 1$ function, where $x_i$ is in segment 1 of the x locus line, is:

$$\int_{\Delta x \; x_i}^{\infty} \frac{1}{x^n}\,\Delta x = \Delta x \sum_{\Delta x \; x=x_i}^{\infty} \frac{1}{x^n} = \mathrm{lnd}(n,\Delta x,x_i) \approx \mathrm{lnd}_f(n,\Delta x,x_i) + K_1 \qquad (2.6\text{-}1)$$

where

$K_1$ = constant of integration for all x on segment 1 of the complex plane x locus line
where $Re(\Delta x)>0$ or $\{Re(\Delta x)=0$ and $Im(\Delta x)>0\}$

$x_i$ = the summation initial value of x, a value of x on segment 1 of the x locus line in the complex plane

$x = x_i + m\Delta x$ , $m=1,2,3,\ldots$ , The x locus

$x = x_i, x_i+\Delta x, x_i+2\Delta x, x_i+3\Delta x, \ldots$

$Re(\Delta x)>0$ or $\{Re(\Delta x)=0$ and $Im(\Delta x)>0\}$

$\Delta x$ = x increment

Accuracy increases rapidly as $\left|\frac{x_i}{\Delta x}\right|$ increases in value

Rewriting Eq 2.5-73

$$lnd_f(n,\Delta x,x) = -\sum_{m=0}^{\infty} \frac{\Gamma(n+2m-1)\left(\frac{\Delta x}{2}\right)^{2m} C_m}{\Gamma(n)(2m+1)!\left(x - \frac{\Delta x}{2}\right)^{n+2m-1}} \text{ , } n\neq 1 \qquad (2.6\text{-}2)$$

For the summation, $\sum_{\Delta x, x=x_i}^{\infty} \frac{1}{x^n}$ , some examples of x locus lines, which will be used to describe the constant, $K_1$, are shown in Diagram 2.6-1 below.

Diagram 2.6-1: The location in the complex plane of x locus line segment 1

The darkened arrows represent segment 1 of the x locus line. For $Re(\Delta x)>0$, x goes from left to right along the x locus line in the summation process from x to $\infty$. For $Re(\Delta x)=0$ and $Im(\Delta x)<0$, x goes from bottom to top along the x locus line in the summation process from x to $\infty$. The $lnd(n,\Delta x,x)$ $n\neq 1$ Series constant of integration value, $K_1$, of x locus line segment 1, will be shown to be equal to 0. $K_1$ retains its zero value for all $n,\Delta x$ where $x_i$ is located on segment 1 of the x locus line.

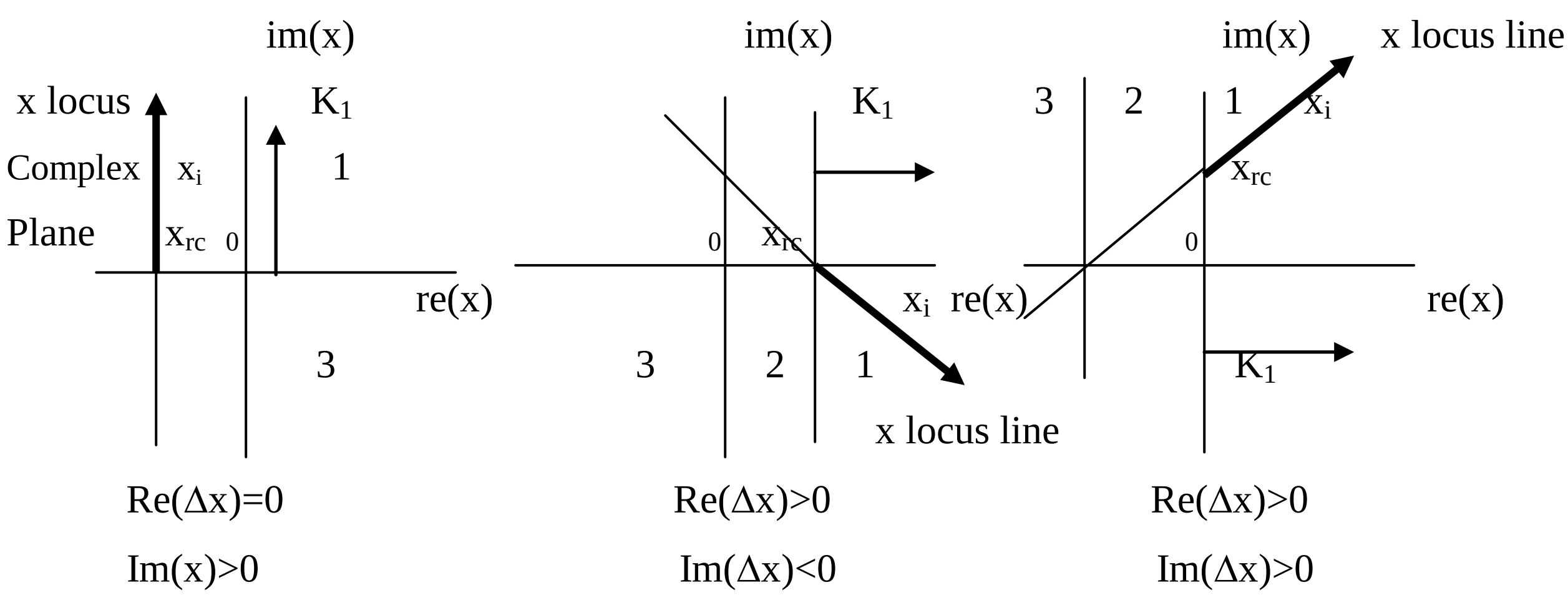

$x_{rc}$ = rightmost or only x locus line and complex plane axis crossover point

$x_{rcr}$ = real part of $x_{rc}$

$Re(\Delta x)>0$ or $\{Re(\Delta x)=0$ and $im(\Delta x)>0\}$

$x = x_i + m\Delta x$ ,   m = 1,2,3,…   This is the x locus

$x = x_i, x_i+\Delta x, x_i+2\Delta x, x_i+3\Delta x, \ldots$

1,2,3 represent x locus line segment designations

x locus line segment 1 is where: $Re(\Delta x)>0$, $Re(x)>x_{rcr}$

or

$Re(\Delta x)=0$, $Im(\Delta x)>0$, $Im(x)>0$

Show that $K_1 = 0$.

1) Consider the case where $Re(\Delta x)>0$ or $\{Re(\Delta x)=0$ and $Im(\Delta x)>0\}$, $1<Re(n)$ , and $n\neq 1$

$$lnd(n,\Delta x,x_i) = \Delta x \sum_{\Delta x\ x=x_i}^{\infty} \frac{1}{x^n} = \int_{\Delta x\ x_i}^{\infty} \frac{1}{x^n}\,\Delta x \qquad (2.6\text{-}3)$$

where

$x = x_i + m\Delta x$,   m = 0,1,2,3, …   This is the x locus

$\Delta x$ = x increment

$1<Re(n)$

Eq 2.6-3 above has been shown to be valid

$lnd_f(n,\Delta x,x)$ = the $lnd(n,\Delta x,x)$ $n\neq 1$ Series formula only (See Eq 2.6-2)

For $x_i$ on segment 1 of the complex plane x locus line, the constant of integration is $K_1$ where $Re(\Delta x)>0$ or $\{Re(\Delta x)=0$ and $im(\Delta x)>0\}$

$$\int_{\Delta x\ x_i}^{\infty} \frac{1}{x^n}\,\Delta x = lnd(n,\Delta x,x_i) \approx lnd_f(n,\Delta x,x_i) + K_1 \ ,\ Re(\Delta x)>0 \qquad (2.6\text{-}4)$$

Accuracy increases rapidly as $\left|\frac{x_i}{\Delta x}\right|$ increases in value

Find $K_1$

$x_i$ is a value of x on x locus line segment 1

$$\lim_{x_i\to\infty} \int_{\Delta x\ x_i}^{\infty} \frac{1}{x^n}\,\Delta x = \lim_{x_i\to\infty} \left[\Delta x \sum_{\Delta x\ x=x_i}^{\infty} \frac{1}{x^n}\right] = \lim_{x_i\to\infty} lnd(n,\Delta x,x_i) \to 0 \qquad (2.6\text{-}5)$$

From Eq 2.6-2

$$\lim_{x_i\to\infty} \text{lnd}_f(n,\Delta x,x_i)\to 0 \qquad (2.6\text{-}6)$$

From Eq 2.6-4, Eq 2.6-5 and Eq 2.6-6 where $x\to\infty$

$$0 = 0 + K_1$$
$$K_1 = 0 \qquad (2.6\text{-}7)$$

The function $\text{lnd}(n,\Delta x,x_i)$ has a constant of integration, $K_1 = 0$ for $x_i\to\infty$. Since $x_i$ is a value of x on segment 1 of the x locus line, the $\text{lnd}(n,\Delta x,x)$ $n\neq 1$ Series Snap Hypothesis specifies that for all other x on segment 1 the constant of integration must be the same, $K_1 = 0$.

Thus from the above equations:

**For $x_i$ on segment 1 of the x locus line**

$$\mathbf{lnd(n,\Delta x,x_i) = \Delta x \sum_{\Delta x\ x=x_i}^{\infty} \frac{1}{x^n} = \int_{\Delta x\ x_i}^{\infty} \frac{1}{x^n}\,\Delta x}$$

**where**

**$x = x_i + m\Delta x$, $m = 0,1,2,3,\ldots$ This is the x locus**

**$\Delta x$ = x increment**

**$1<Re(n)$**

**$n \neq 1$**

**$Re(\Delta x)>0$ or $\{Re(\Delta x)=0$ and $Im(\Delta x)>0\}$**

$$\mathbf{\Delta x \sum_{\Delta x\ x=x_i}^{\infty} \frac{1}{x^n} = \int_{\Delta x\ x_i}^{\infty} \frac{1}{x^n}\,\Delta x = lnd(n,\Delta x,x_i) \approx lnd_f(n,\Delta x,x_i) + K_1}$$

**Accuracy increases rapidly as $\left|\frac{x_i}{\Delta x}\right|$ increases in value**

**$K_1 = 0$**

2) Consider the case where $Re(\Delta x)>0$ or $\{Re(\Delta x)=0$ and $Im(\Delta x)>0\}$, $0<Re(n)$, and $n\neq 1$

$$\int_{\Delta x\ x_i}^{\infty} (-1)^{\frac{x-x_i}{\Delta x}} \frac{1}{x^n}\,\Delta x = \Delta x \sum_{\Delta x\ x=x_i}^{\infty} (-1)^{\frac{x-x_i}{\Delta x}} \frac{1}{x^n} = \text{lnd}(n,2\Delta x,x_i) - \text{lnd}(n,2\Delta x,x_i+\Delta x) \qquad (2.6\text{-}8)$$

where

$x = x_i + m\Delta x$, $m = 0,1,2,3,\ldots$ This is the x locus

$\Delta x$ = x increment

$0 < Re(n)$

$n \neq 1$

Eq 2.6-8 above has been shown to be valid

In segment 1 of the complex plane x locus line, the constant of integration is $K_1$ where $Re(\Delta x)>0$ or $Re(\Delta x)=0$ and $im(\Delta x)>0$}

$$lnd(n,\Delta x,x) \approx lnd_f(n,\Delta x,x) + K_1 \qquad (2.6\text{-}9)$$

Accuracy increases rapidly as $|\frac{x}{\Delta x}|$ increases in value

Find $K_1$

$x_i$ is a value of x on x locus line segment 1

Substitute Eq 2.6-9 into Eq 2.6-8

$$\underset{\Delta x}{\int_{x_i}^{\infty}} (-1)^{\frac{x-x_i}{\Delta x}} \frac{1}{x^n}\,\Delta x \approx lnd_f(n,2\Delta x,x_i) + K_1 - lnd(n,2\Delta x,x_i+\Delta x) - K_1 \qquad (2.6\text{-}10)$$

Accuracy increases rapidly as $|\frac{x_i}{\Delta x}|$ increases in value

The constant of integration, $K_1$, cancels out of Eq 2.6-10. Therefore, the value of the integral in Eq 2.6-10 is not a function of $K_1$. Then, the value of $K_1$, in this case, may be chosen arbitrarily. $K_1$ is chosen to have a value of 0.

Then

$$K_1 = 0 \qquad (2.6\text{-}11)$$

Thus

**For $x_i$ on segment 1 of the x locus line**

$$\underset{\Delta x}{\int_{x_i}^{\infty}} (-1)^{\frac{x-x_i}{\Delta x}} \frac{1}{x^n}\,\Delta x = \Delta x \underset{\Delta x}{\sum_{x=x_i}^{\infty}} (-1)^{\frac{x-x_i}{\Delta x}} \frac{1}{x^n} = lnd(n,2\Delta x,x_i) - lnd(n,2\Delta x,x_i+\Delta x)$$

**where**

**$x = x_i + m\Delta x$, $m = 0,1,2,3,\ldots$ This is the x locus**

**$\Delta x$ = x increment**

**$Re(\Delta x)>0$ or $\{Re(\Delta x)=0$ and $Im(\Delta x)>0\}$**

**$0 < Re(n)$**

**$n \neq 1$**

**$lnd(n,\Delta x,x_i) \approx lnd_f(n,\Delta x,x_i) + K_1$**

**Accuracy increases rapidly as $|\frac{x_i}{\Delta x}|$ increases in value**

**$K_1 = 0$**

3) Consider the case where $Re(\Delta x)>0$ or $\{Re(\Delta x)=0$ and $Im(\Delta x)>0\}$, $n \neq 1$

$$\int_{\Delta x\ x_1}^{x_2} \frac{1}{x^n}\,\Delta x = \Delta x \sum_{\Delta x\ x=x_1}^{x_2} \frac{1}{x^n} = -lnd(n,\Delta x,x)\Big|_{x_1}^{x_2} = -lnd(n,\Delta x,x_2) + lnd(n,\Delta x,x_1) \qquad (2.6\text{-}12)$$

$x = x_1, x_1+\Delta x, x_1+2\Delta x, x_1+3\Delta x, x_2-\Delta x, x_2$

$x = x_1 + m\Delta x, \quad m=1,2,3,\ldots,\frac{x_2 - x_1}{\Delta x}$ , This is the x locus

$\Delta x$ = x increment

$n\neq 1$

Eq 2.6-12 above has been shown to be valid

For x on segment 1 of the complex plane x locus line, the constant of integration is $K_1$ where $Re(\Delta x)>0$ or $\{Re(\Delta x)=0$ and $im(\Delta x)>0\}$

All points, $x_1$ thru $x_2$, are on x locus segment 1

Rewriting Eq 2.6-9

$lnd(n,\Delta x,x) \approx lnd_f(n,\Delta x,x) + K_1$

Accuracy increases rapidly as $|\frac{x}{\Delta x}|$ increases in value

Find $K_1$

Substitute Eq 2.6-9 into Eq 2.6-12

$$\int_{\Delta x\ x_1}^{x_2} \frac{1}{x^n}\,\Delta x \approx -lnd_f(n,\Delta x,x_2) - K_1 + lnd_f(n,\Delta x,x_1) + K_1 \qquad (2.6\text{-}13)$$

Accuracy increases rapidly as $|\frac{x_1}{\Delta x}|$ and $|\frac{x_2}{\Delta x}|$ increase in value

The constant of integration, $K_1$, cancels out of Eq 2.6-13. Therefore, the value of the integral in Eq 2.6-13 is not a function of $K_1$. Then, the value of $K_1$, in this case, may be chosen arbitrarily. $K_1$ is chosen to have a value of 0.

Then

$$K_1 = 0 \tag{2.6-14}$$

Thus

**For x = $x_1$ thru $x_2$ on segment 1 of the x locus line**

$$ {}_{\Delta x}\int_{x_1}^{x_2} \frac{1}{x^n}\,\Delta x = {}_{\Delta x}\sum_{x=x_1}^{x_2} \frac{1}{x^n} = -\text{lnd}(n,\Delta x,x)\Big|_{x_1}^{x_2} = -\text{lnd}(n,\Delta x,x_2) + \text{lnd}(n,\Delta x,x_1) $$

**where**

$x = x_1,\ x_1+\Delta x,\ x_1+2\Delta x,\ x_1+3\Delta x,\ x_2-\Delta x,\ x_2$

$x = x_1 + m\Delta x,\quad m=1,2,3,\ldots,\frac{x_2 - x_1}{\Delta x}$ , **This is the x locus**

$\Delta x$ **= x increment**

$Re(\Delta x)>0$ **or** $\{Re(\Delta x)=0$ **and** $Im(\Delta x)>0\}$

$n\neq 1$

$\text{lnd}(n,\Delta x,x) \approx \text{lnd}_f(n,\Delta x,x) + K_1$

**Accuracy increases rapidly as** $|\frac{x_1}{\Delta x}|$ **and** $|\frac{x_2}{\Delta x}|$ **increase in value**

$K_1 = 0$

From Eq 2.6-7, Eq 2.6-11, and Eq 2.6-14

**For x values on x locus segment 1 where** $Re(\Delta x)>0$ **or** $\{Re(\Delta x)=0$ **and** $Im(\Delta x)>0\}$ , $\frac{1}{x^n}$ **sums or integrals may be calculated using the relationship,**

$$\text{lnd}(n,\Delta x,x) \approx \text{lnd}_f(n,\Delta x,x) + K_1 ,\quad n\neq 1$$

$$K_1 = 0 \tag{2.6-15}$$

**Accuracy increases rapidly as** $|\frac{x}{\Delta x}|$ **increases in value.**

Constant of integration, $k_3$

The relationship for the lnd(n,Δx,x) n≠1 function, where $x_i$ is in segment 3 of the x locus line, is:

$$\int_{\Delta x\, x_i}^{-\infty} \frac{1}{x^n}\Delta x = \Delta x \sum_{\Delta x\, x=x_i}^{-\infty} \frac{1}{x^n} = \text{lnd}(n,\Delta x,x_i) \approx \text{lnd}_f(n,\Delta x,x_i) + k_3 \qquad (2.6\text{-}16)$$

where

$k_3$ = constant of integration for all x on segment 3 of the complex plane x locus where Re(Δx)<0 or {Re(Δx)=0 and Im(Δx)<0

$x_i$ = the summation initial value of x, a value of x on segment 3 of the x locus line in the complex plane

$x = x_i + m\Delta x$ , m=1,2,3,… , This is the x locus

$x = x_i, x_i+\Delta x, x_i+2\Delta x, x_i+3\Delta x, \ldots$

Re(Δx)<0 or {Re(Δx)=0 and Im(Δx)<0

Δx = x increment

Accuracy increases rapidly as $\left|\frac{x_i}{\Delta x}\right|$ increases in value

Rewriting Eq 2.6-2

$$\text{lnd}_f(n,\Delta x,x) = -\sum_{m=0}^{\infty} \frac{\Gamma(n+2m-1)\left(\frac{\Delta x}{2}\right)^{2m} C_m}{\Gamma(n)(2m+1)!\left(x - \frac{\Delta x}{2}\right)^{n+2m-1}} , \ n\neq 1$$

For the summation, $\sum_{\Delta x\, x=x_i}^{-\infty} \frac{1}{x^n}$ , some examples of x locus lines which will be used to describe the constant, $k_3$, are shown in the following diagram, Diagram 2.6-2.

Diagram 2.6-2: The location in the complex plane of x locus line segment 3

The darkened arrows represent segment 3 of the x locus. For Re(Δx)<0, x goes from right to left along the x locus line in the summation process from x to -∞. For Re(Δx)=0 and Im(Δx)<0, x goes from top to bottom along the x locus line in the summation process from x to -∞. The lnd(n,Δx,x) n≠1 Series constant of integration value, $k_3$, of x locus line segment 3, will be shown to be equal to 0. $k_3$ retains its zero value for all n,Δx where $x_i$ is located on segment 3 of the x locus line.

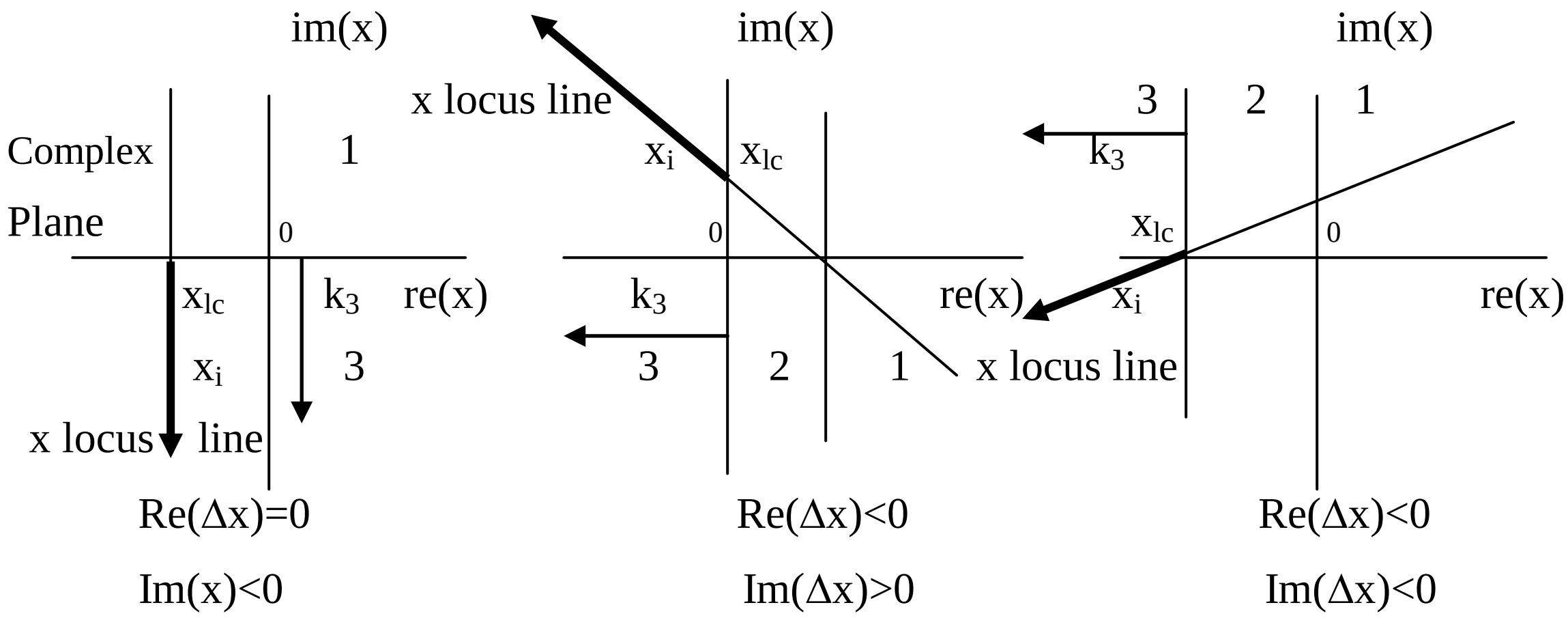

$x_{lc}$ = leftmost or only x locus and complex plane axis crossover point

$x_{lcr}$ = the real part of $x_{lc}$

Re(Δx)<0 or {Re(Δx)=0 and im(Δx)<0}

$x = x_i + m\Delta x$ ,   m = 1,2,3,…     This is the x locus

$x = x_i, x_i+\Delta x, x_i+2\Delta x, x_i+3\Delta x, \ldots$

1,2,3 represent x locus segment designations

x locus line segment 3 is where:   Re(Δx)<0,  Re(x)<$x_{lcr}$
or
Re(Δx)=0, Im(Δx)<0,  Im(x)<0

Show that $k_3 = 0$

1) Consider the case where Re(Δx)<0 or {Re(Δx)=0 and Im(Δx)<0}, 1<Re(n) , and n≠1

$$\mathrm{lnd}(n,\Delta x,x_i) = \Delta x \sum_{\Delta x \; x=x_i}^{-\infty} \frac{1}{x^n} = \int_{\Delta x \; x_i}^{-\infty} \frac{1}{x^n}\,\Delta x \qquad (2.6\text{-}17)$$

where

$x = x_i + m\Delta x$,   m = 0,1,2,3, …     This is the x locus

Δx = x increment

1<Re(n)

Eq 2.6-17 above has been shown to be valid

$\mathrm{lnd}_f(n,\Delta x,x)$ = the lnd(n,Δx,x) n≠1 Series formula only (See Eq 2.6-2)

For x on segment 3 of the complex plane x locus line, the constant of integration is $k_3$ where $Re(\Delta x)<0$ or $\{Re(\Delta x)=0$ and $im(\Delta x)<0\}$

$$\int_{\Delta x\ x_i}^{-\infty} \frac{1}{x^n}\Delta x = lnd(n,\Delta x,x_i) \approx lnd_f(n,\Delta x,x_i) + k_3\ ,\quad Re(\Delta x)<0 \tag{2.6-18}$$

Accuracy increases rapidly as $|\frac{x_i}{\Delta x}|$ increases in value

Find $k_3$

$x_i$ is a value of x on x locus line segment 3

$$\lim_{x_i\to-\infty}\ \int_{\Delta x\ x_i}^{-\infty} \frac{1}{x^n}\Delta x = \lim_{x_i\to-\infty}\ [\Delta x \sum_{\Delta x\ x=x_i}^{-\infty} \frac{1}{x^n}] = \lim_{x_i\to-\infty} lnd(n,\Delta x,x_i) \to 0 \tag{2.6-19}$$

From Eq 2.6-2

$$\lim_{x_i\to-\infty} lnd_f(n,\Delta x,x_i) \to 0 \tag{2.6-20}$$

From Eq 2.6-18, Eq 2.6-19 and Eq 2.6-20 where $x\to-\infty$

$$0 = 0 + k_3$$

$$k_3 = 0 \tag{2.6-21}$$

The function $lnd(n,\Delta x,x_i)$ has a constant of integration, $k_3 = 0$ for $x_i\to-\infty$. Since $x_i$ is a value of x on segment 3 of the x locus line, the $lnd(n,\Delta x,x)$ $n\neq1$ Series Snap Hypothesis specifies that for all other x on segment 3 the constant of integration must be the same, $k_3 = 0$.

Thus from the above equations

**For x on segment 3 of the x locus line**

$$\mathbf{lnd(n,\Delta x,x_i) = \Delta x \sum_{\Delta x\ x=x_i}^{-\infty} \frac{1}{x^n} = \int_{\Delta x\ x_i}^{-\infty} \frac{1}{x^n}\Delta x}$$

**where**

**$x = x_i + m\Delta x$,   $m = 0,1,2,3, \ldots$     This is the x locus**

**$\Delta x$ = x increment**

**$1<Re(n)$**

**$n \neq 1$**
**$Re(\Delta x)<0$ or $\{Re(\Delta x)=0$ and $Im(\Delta x)<0\}$**

$$\Delta x \sum_{\Delta x \; x=x_i}^{-\infty} \frac{1}{x^n} = \int_{\Delta x \; x_i}^{-\infty} \frac{1}{x^n}\,\Delta x = \mathrm{lnd}(n,\Delta x,x_i) \approx \mathrm{lnd}_f(n,\Delta x,x_i) + k_3$$

**$k_3 = 0$**

**Accuracy increases rapidly as $\left|\frac{x_i}{\Delta x}\right|$ increases in value**

2) Consider the case where $Re(\Delta x)<0$ or $\{Re(\Delta x)=0$ and $Im(\Delta x)<0\}$, $0<Re(n)$, and $n\neq 1$

$$\int_{\Delta x \; x_i}^{-\infty} (-1)^{\frac{x-x_i}{\Delta x}} \frac{1}{x^n}\,\Delta x = \Delta x \sum_{\Delta x \; x=x_i}^{-\infty} (-1)^{\frac{x-x_i}{\Delta x}} \frac{1}{x^n} = \mathrm{lnd}(n,2\Delta x,x_i) - \mathrm{lnd}(n,2\Delta x,x_i+\Delta x) \tag{2.6-22}$$

where

$x = x_i + m\Delta x$, $m = 0,1,2,3, \ldots$ This is the x locus

$\Delta x$ = x increment

$0 < Re(n)$

$n \neq 1$

Eq 2.6-22 above has been shown to be valid

In segment 3 of the complex plane x locus line, the constant of integration is $k_3$ where $Re(\Delta x)<0$ or $Re(\Delta x)=0$ and $im(\Delta x)<0\}$

$$\mathrm{lnd}(n,\Delta x,x) \approx \mathrm{lnd}_f(n,\Delta x,x) + k_3 \tag{2.6-23}$$

Accuracy increases rapidly as $\left|\frac{x}{\Delta x}\right|$ increases in value

Find $k_3$

$x_i$ is a value of x on the x locus line segment 3

Substitute Eq 2.6-23 into Eq 2.6-22

$$\int_{\Delta x \; x_i}^{-\infty} (-1)^{\frac{x-x_i}{\Delta x}} \frac{1}{x^n}\,\Delta x \approx \mathrm{lnd}_f(n,2\Delta x,x_i) + k_3 - \mathrm{lnd}(n,2\Delta x,x_i+\Delta x) - k_3 \tag{2.6-24}$$

Accuracy increases rapidly as $\left|\frac{x_i}{\Delta x}\right|$ increases in value

The constant of integration, $k_3$, cancels out of Eq 2.6-23. Therefore, the value of the integral in Eq 2.6-24 is not a function of $k_3$. Then, the value of $k_3$, in this case, may be chosen arbitrarily. $k_3$ is chosen to have a value of 0.

Then

$$k_3 = 0 \tag{2.6-25}$$

Thus

**For x on segment 3 of the x locus line**

$$\mathbf{{}_{\Delta x}\int_{x_i}^{-\infty} (-1)^{\frac{x-x_i}{\Delta x}} \frac{1}{x^n} \Delta x = \Delta x \sum_{\Delta x \; x=x_i}^{-\infty} (-1)^{\frac{x-x_i}{\Delta x}} \frac{1}{x^n} = lnd(n,2\Delta x,x_i) - lnd(n,2\Delta x,x_i+\Delta x)}$$

**where**

**$x = x_i + m\Delta x$, $m = 0,1,2,3,\dots$ This is the x locus**

**$\Delta x$ = x increment**

**$Re(\Delta x)<0$ or $\{Re(\Delta x)=0$ and $Im(\Delta x)<0\}$**

**$0 < Re(n)$**

**$n \neq 1$**

**$lnd(n,\Delta x,x_i) \approx lnd_f(n,\Delta x,x_i) + k_3$**

**$k_3 = 0$**

**Accuracy increases rapidly as $\left|\frac{x_i}{\Delta x}\right|$ increases in value**

3) Consider the case where $Re(\Delta x)<0$ or $\{Re(\Delta x)=0$ and $Im(\Delta x)<0\}$, $n \neq 1$

$${}_{\Delta x}\int_{x_1}^{x_2} \frac{1}{x^n} \Delta x = \Delta x \sum_{\Delta x \; x=x_1}^{x_2} \frac{1}{x^n} = -lnd(n,\Delta x,x)\Big|_{x_1}^{x_2} = -lnd(n,\Delta x,x_2) + lnd(n,\Delta x,x_1) \tag{2.6-26}$$

$x = x_1, x_1+\Delta x, x_1+2\Delta x, x_1+3\Delta x, x_2-\Delta x, x_2$

$x = x_1 + m\Delta x$, $m=1,2,3,\dots, \frac{x_2-x_1}{\Delta x}$ , This is the x locus

$\Delta x$ = x increment

$n \neq 1$

Eq 2.6-26 above has been shown to be valid

For x on segment 3 of the complex plane x locus line, the constant of integration is $k_3$ where Re(Δx)<0 or {Re(Δx)=0 and im(Δx)<0}

All points are on x locus segment 3

Rewriting Eq 2.6-23

$$\mathrm{lnd}(n,\Delta x,x) \approx \mathrm{lnd_f}(n,\Delta x,x) + k_3$$

Accuracy increases rapidly as $|\frac{x}{\Delta x}|$ increases in value

Find $k_3$

Substitute Eq 2.6-23 into Eq 2.6-26

$$\int_{\Delta x\, x_1}^{x_2} \frac{1}{x^n}\,\Delta x \approx -\mathrm{lnd_f}(n,\Delta x,x_2) - k_3 + \mathrm{lnd_f}(n,\Delta x,x_1) + k_3 \qquad (2.6\text{-}27)$$

Accuracy increases rapidly as $|\frac{x_1}{\Delta x}|$ and $|\frac{x_2}{\Delta x}|$ increase in value

The constant of integration, $k_3$, cancels out of Eq 2.6-27. Therefore, the value of the integral in Eq 2.6-26 is not a function of $k_3$. Then, the value of $k_3$, in this case, may be chosen arbitrarily. $k_3$ is chosen to have a value of 0.

Then

$$k_3 = 0 \qquad (2.6\text{-}28)$$

Thus

**For x = $x_1$ thru $x_2$ on segment 3 of the x locus line**

$$\int_{\Delta x\, x_1}^{x_2} \frac{1}{x^n}\,\Delta x = \Delta x \sum_{\Delta x\, x=x_1}^{x_2} \frac{1}{x^n} = -\mathrm{lnd}(n,\Delta x,x)\Big|_{x_1}^{x_2+\Delta x} = -\mathrm{lnd}(n,\Delta x,x_2+\Delta x) + \mathrm{lnd}(n,\Delta x,x_1)$$

**$x = x_1, x_1+\Delta x, x_1+2\Delta x, x_1+3\Delta x, x_2-\Delta x, x_2$**

**$x = x_1 + m\Delta x$, $m=1,2,3,\ldots,\frac{x_2-x_1}{\Delta x}$ ,** **This is the x locus**

**Δx = x increment**

**Re(Δx)<0 or {Re(Δx)=0 and Im(Δx)<0}**

**$n \neq 1$**

**$lnd(n,\Delta x,x) \approx lnd_f(n,\Delta x,x) + k_3$**

**$k_3 = 0$**

**Accuracy increases rapidly as $|\frac{x}{\Delta x}|$ increases in value**

From Eq 2.6-21, Eq 2.6-25, and Eq 2.6-28

**For x values on x locus segment 3 where Re(Δx)<0 or {Re(Δx)=0 and Im(Δx)<0} , $\frac{1}{x^n}$ sums or integrals may be calculated using the relationship,**

**$lnd(n,\Delta x,x) \approx lnd_f(n,\Delta x,x) + k_3$ , n≠1**

**$k_3 = 0$** (2.6-29)

**Accuracy increases rapidly as $|\frac{x}{\Delta x}|$ increases in value.**

<u>A Characteristic of each complex plane x locus</u>

From Eq 2.5-45

$$lnd(n,\Delta x,x) - lnd(n,-\Delta x,x-\Delta x) = \pm(K_r - k_r) \qquad (2.6\text{-}30)$$

where

$x = x_i + m\Delta x$ , m = integers , This is an x locus in the complex plane

$x_i$ = a value of x

Δx = x increment

+ for Re(Δx)>0 or {Re(Δx)=0 and Im(Δx)>0}

– for Re(Δx)<0 or {Re(Δx)=0 and Im(Δx)<0}

x = real or complex values

r = 1,2,3

n≠1

$n,x_i,K_r,k_r$ = real or complex constants

Eq 2.6-30 can be rewritten in a different form

$$K_r - k_r = \begin{cases} lnd(n,\Delta x,x) - lnd(n,-\Delta x,x-\Delta x) & \text{for } Re(\Delta x)>0 \text{ or } \{Re(\Delta x)=0 \text{ and } Im(\Delta x)>0\} \\ -[lnd(n,-\Delta x,x) - lnd(n,\Delta x,x+\Delta x)] & \text{for } Re(\Delta x)<0 \text{ or } \{Re(\Delta x)=0 \text{ and } Im(\Delta x)<0\} \end{cases} \qquad (2.6\text{-}31)$$

Note that $x_i$ together with Δx and $x_i$ together with -Δx determines the same x locus and x locus line in the complex plane.

Rewriting Eq 2.5-48

$$lnd(n,\Delta x,x) - lnd(n,-\Delta x,x-\Delta x) = -[lnd(n,-\Delta x,x) - lnd(n,\Delta x,x+\Delta x)] \qquad (2.6\text{-}32)$$

From Eq 2.6-31 and Eq 2.6-32 it is seen that the constant, $K_r - k_r$, is the same value along an x locus and x locus line with a specified value of n irrespective of complex plane segment or summation direction.

$\mathbf{K_r - k_r = constant}$ (2.6-33)

Note - The condition, Re(Δx)>0 or {Re(Δx)=0 and Im(Δx)>0}, is associated with a summation to +∞ and motion along an x locus line in the complex plane from left to right or from down to up. The condition, Re(Δx)<0 or {Re(Δx)=0 and Im(Δx)<0}, is associated with a summation to -∞ and motion along an x locus line in the complex plane from right to left or from up to down.

## Conclusion

Diagram 2.6-3: Locations where $\text{lnd}(n,\Delta x,x) \approx \text{lnd}_f(n,\Delta x,x)$

**For $x_i$, a value of x, on the following highlighted x locus segments, $\text{lnd}(n,\Delta x,x_i) \approx \text{lnd}_f(n,\Delta x,x_i)$ ($K_1=k_3=0$).**

**Accuracy increases rapidly as $|\frac{x_i}{\Delta x}|$ increases in value.**

For Re(Δx)>0 and {Re(Δx)=0 and Im(Δx)>0

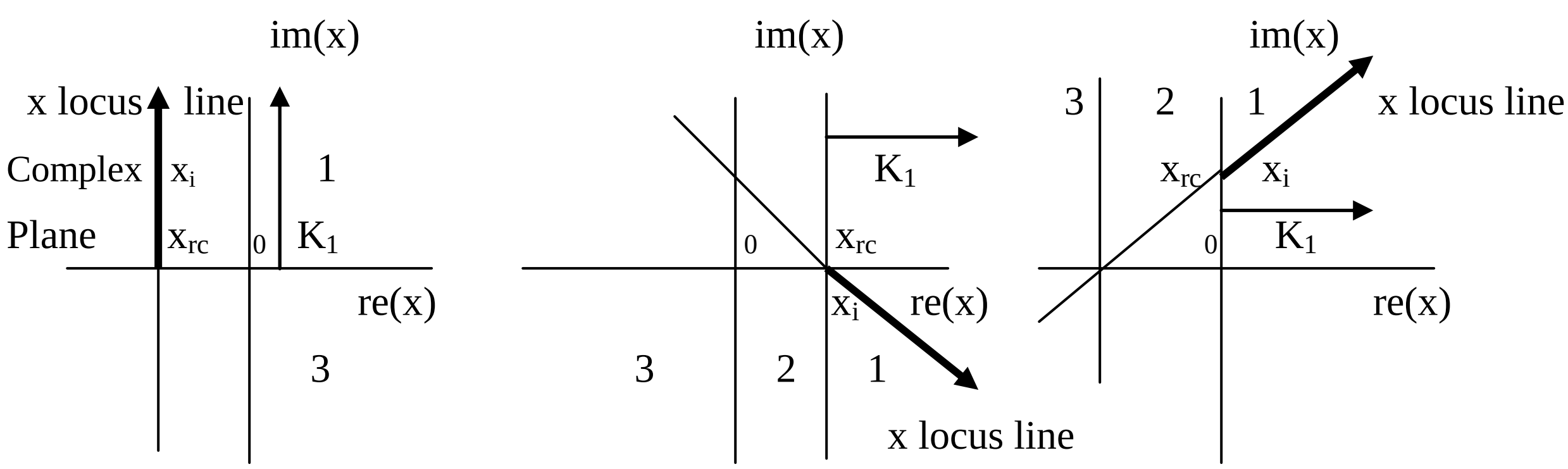

Re(Δx)=0 Re(Δx)>0 Re(Δx)>0

Im(x)>0 Im(Δx)<0 Im(Δx)>0

For Re(Δx)<0 and {Re(Δx)=0 and Im(Δx)}<0

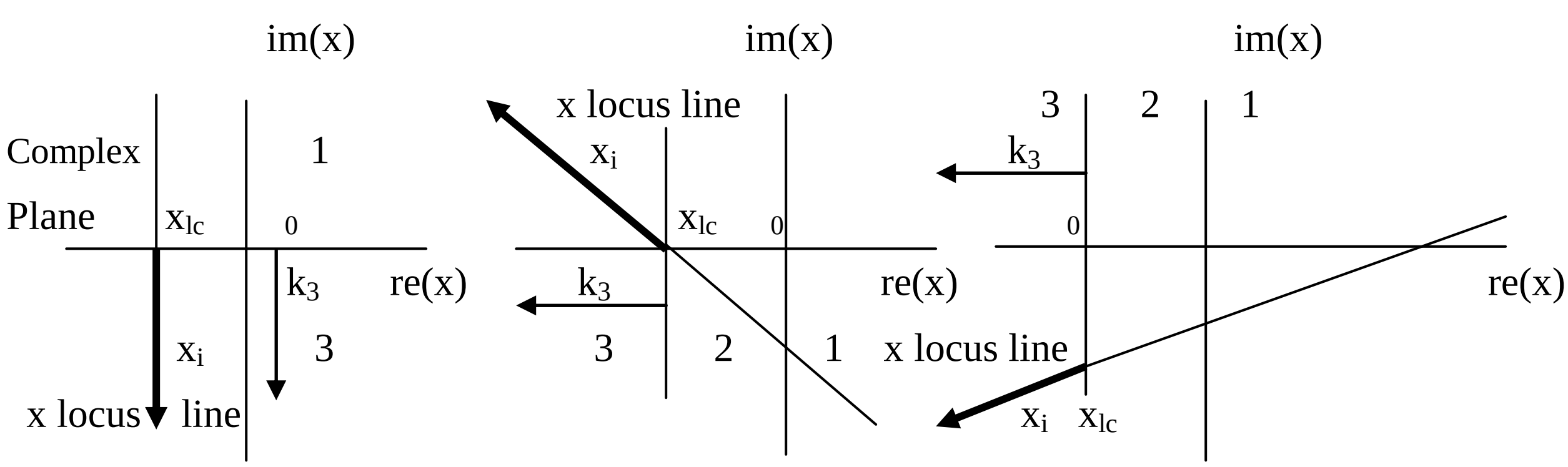

Re(Δx)=0 Re(Δx)<0 Re(Δx)<0

Im(x)<0 Im(Δx)>0 Im(Δx)<0

$x_{rc}$ = rightmost or only x locus line and complex plane axis crossover point

$x_{lc}$ =leftmost or only x locus ine and complex plane axis crossover point

$x_{rcr}$ = real part of $x_{rc}$

$x_{lcr}$ = real part of $x_{lc}$

$x = x_i + m\Delta x$ , m = 1,2,3,… This is the x locus

$x = x_i, x_i+\Delta x, x_i+2\Delta x, x_i+3\Delta x, \ldots$

1,2,3 represent x locus segment designations

x locus segment 1 is where: $\mathrm{Re}(\Delta x)>0$, $\mathrm{Re}(x)>x_{rcr}$
or
$\mathrm{Re}(\Delta x)=0$, $\mathrm{Im}(\Delta x)>0$, $\mathrm{Im}(x)>0$

x locus segment 3 is where: $\mathrm{Re}(\Delta x)<0$, $\mathrm{Re}(x)<x_{lcr}$
or
$\mathrm{Re}(\Delta x)=0$, $\mathrm{Im}(\Delta x)<0$, $\mathrm{Im}(x)<0$

The above conclusion is very important. With the relationship, $\mathrm{lnd}(n,\Delta x,x) \approx \mathrm{lnd}_f(n,\Delta x,x)$, valid in the regions and under the conditions shown, the function $\mathrm{lnd}(n,\Delta x,x)$ $n\neq1$ can be calculated for all x. Two methods are available to do this. The method which does not require any knowledge of the constants previously mentioned, $K_2$, $K_3$, $k_1$, and $k_2$, is presented in the following section, Section 2.7. The second method which requires the calculated values of $K_2$, $K_3$, $k_1$, and $k_2$ is presented in Section 2.8. The second method requires more computer coding but is generally faster.

Example 2.6-1

Put to use the following equations and concepts derived and described in the previous sections of Chapter 2. Demonstrate the validity of these equations and concepts.

For this demonstration, arbitrarily let:
$n = 3.1-1.2i$
$\Delta x = 1+i$
$x = -1000 + m(1+i)$ , m = integers This is a linear x locus in the complex plane

**Equations**

**E1)** $$\Delta x \sum_{\Delta x \; x=x_1}^{x_2} \frac{1}{x^n} = \int_{\Delta x \; x_1}^{x_2} \frac{1}{x^n}\Delta x = -\mathrm{lnd}(n,\Delta x,x)\Big|_{x_1}^{x_2+\Delta x} = -\mathrm{lnd}(n,\Delta x,x_2+\Delta x) + \mathrm{lnd}(n,\Delta x,x_1)]$$

$x = x_1, x_1+\Delta x, x_1+2\Delta x, x_1+3\Delta x, \ldots, x_2-\Delta x, x_2$ (from Eq 2.2-13)

$x_1,x_2$ can be infinite values

**E2)** $$\Delta x \sum_{\Delta x \; x=x_i}^{\pm\infty} \frac{1}{x^n} = \int_{\Delta x \; x_i}^{\pm\infty} \frac{1}{x^n}\Delta x = \text{lnd}(n,\Delta x,x_i) \approx \text{lnd}_f(n,\Delta x,x_i) + \begin{cases} K_r \\ k_r \end{cases}$$

$x = x_i, x_i+\Delta x, x_i+2\Delta x, x_i+3\Delta x + x_i+4\Delta x, x_i+5\Delta x, \ldots$ (from Eq 2.2-12, Eq 2.5-74, Eq 2.5-75)

$Re(n) > 1$

Accuracy increases rapidly as $|\frac{x}{\Delta x}|$ increases in value.

$K_r, k_r$ = constants of integration

r = 1, 2, or 3

**E3)** $$\Delta x \sum_{\Delta x \; x=\mp\infty}^{\pm\infty} \frac{1}{x^n} = \text{lnd}(n,\Delta x,x_i) - \text{lnd}(n,-\Delta x,x_i-\Delta x) = \pm (K_r - k_r)$$ (from Eq 2.5-66)

$Re(n)>1$

**E4)** $$\text{lnd}_f(n,\Delta x,x) = -\sum_{m=0}^{\infty} \frac{\Gamma(n+2m-1)\left(\frac{\Delta x}{2}\right)^{2m} C_m}{\Gamma(n)(2m+1)!\left(x - \frac{\Delta x}{2}\right)^{n+2m-1}}$$ (from Eq 2.5-73)

**E5)** $$\text{lnd}(n,\Delta x,x) \approx \text{lnd}_f(n,\Delta x,x) + \begin{cases} K_r \\ k_r \end{cases},$$ x is within the x locus segment, r

The absolute value, $|\frac{x}{\Delta x}|$, must be large for good accuracy (from Eq 2.5-73 and Eq 2.5-74)

**E6)** $$\left.\begin{matrix} K_r \\ k_r \end{matrix}\right\} \approx \text{lnd}(n,\Delta x,x) - \text{lnd}_f(n,\Delta x,x),$$ x is within the x locus segment, r (from Eq 2.5-74)

The absolute value, $|\frac{x}{\Delta x}|$, must be large for good accuracy

**E7)** $$\text{lnd}(n,\Delta x,x_i) = \Delta x \sum_{\Delta x \; x=x_i}^{x_p-\Delta x} \frac{1}{x^n} + \text{lnd}_f(n,\Delta x,x_p) + \begin{cases} K_r \\ k_r \end{cases}$$ (from Eq 2.7-5)

**E8)** $K_1 = 0$ (from Eq 2.6-15)

**E9)** $k_3 = 0$ (from Eq 2.6-29)

**E10)** $K_3 = -k_1$ (from Eq 2.5-77)

**E11)** $\ln d(n,\Delta x,x) - \ln d(n,-\Delta x,x-\Delta x) = \pm(K_r - k_r)$ (from Eq 2.5-45)

**E12)** $\ln d(n,\Delta x,x) - \ln d(n,-\Delta x,x-\Delta x) = -[\ln d(n,-\Delta x,x) - \ln d(n,\Delta x,x+\Delta x)]$ (from Eq 2.5-48)

**E13)** $$k_1 = \pm\left[-\Delta x \sum_{\Delta x\; x=\mp\infty}^{\pm\infty} \frac{1}{x^n}\right], \quad Re(n)>1$$ (from Eq 2.5-82)

**E14)** $$K_3 = \pm\left[\Delta x \sum_{\Delta x\; x=\mp\infty}^{\pm\infty} \frac{1}{x^n}\right], \quad Re(n)>1$$ (from Eq 2.5-83)

**E15)** $$K_2 - k_2 = \pm\left[\Delta x \sum_{\Delta x\; x=\mp\infty}^{\pm\infty} \frac{1}{x^n}\right], \quad Re(n)>1$$ (from Eq 2.5-84)

**E16)** $$K_r - k_r = \pm\left[\Delta x \sum_{\Delta x\; x=\mp\infty}^{\pm\infty} \frac{1}{x^n}\right], \quad Re(n)>1, \quad r = 1,2, \text{ or } 3$$ (from Eq 2.5-85)

**E17)** $$\left.\begin{matrix} K_r \\ k_r \end{matrix}\right\} = \Delta x \sum_{\Delta x\; x=x_i}^{\pm\infty} \frac{1}{x^n} - \ln d_f(n,\Delta x,x_i), \quad Re(n)>1, \quad r = 1,2, \text{ or } 3$$ (from Eq 2.5-86)

$x_i$ is within the x locus segment, r

Accuracy increases rapidly as $|\frac{x}{\Delta x}|$ increases in value.

$x_i$ to $+\infty$ for $Re(\Delta x)>0$ or $\{Re(\Delta x)=0$ and $Im(\Delta x)>0\}$
$x_i$ to $-\infty$ for $Re(\Delta x)<0$ or $\{Re(\Delta x)=0$ and $Im(\Delta x)<0\}$

where

$x = x_i + m\Delta x$, m = integers, This is an x locus in the complex plane
$x_i,x_p$ = values of x
$\Delta x$ = x increment
+ for $Re(\Delta x)>0$ or $\{Re(\Delta x)=0$ and $Im(\Delta x)>0\}$
– for $Re(\Delta x)<0$ or $\{Re(\Delta x)=0$ and $Im(\Delta x)<0\}$
$-\infty$ to $+\infty$ for $Re(\Delta x)>0$ or $\{Re(\Delta x)=0$ and $Im(\Delta x)>0\}$
$+\infty$ to $-\infty$ for $Re(\Delta x)<0$ or $\{Re(\Delta x)=0$ and $Im(\Delta x)<0\}$
x = real or complex values
$n,\Delta x,x_i,x_p Kr,k_r$ = real or complex constants
r = 1,2,3 , The x locus segment designations
$K_r$ = constant of integration for $Re(\Delta x)>0$ or $\{Re(\Delta x) = 0$ and $Im(\Delta x) > 0\}$
$k_r$ = constant of integration for $Re(\Delta x)<0$ or $\{Re(\Delta x) = 0$ and $Im(\Delta x) < 0\}$
$n \neq 1$

The x locus line is the straight line in the complex plane through all of the plotted x points.

**Concepts**

**C1) The x locus and the x locus line in the complex plane**

**C2) Function snap – An exceeding abrupt transition between two successive values, often from a correct to an incorrect value**

**C3) Snap Hypothesis – lnd(n,Δx,x) n≠1 Series snap will occur, if it occurs at all, only at an x locus transition across a complex plane axis and, except where series snap occurs, the lnd(n,Δx,x) n≠1 Series constant of integration will not change.**

**C4) Over all x, the series function, $\mathbf{lnd_f(n,\Delta x,x)}$, can have as many as six constants of integration, $\mathbf{K_1,K_2,K_3,k_1,k_2,k_3}$.**

**Diagram 2.6-4 x locus and x locus line in the complex plane**

$$\Delta x \sum_{\substack{\Delta x \\ x=x_i}}^{\pm\infty} \frac{1}{x^n} = \int_{\substack{\Delta x \\ x_i}}^{\pm\infty} \frac{1}{x^n}\,\Delta x = \mathrm{lnd}(n,\Delta x,x_i) \approx \mathrm{lnd_f}(n,\Delta x,x_i) + \begin{cases} K_r \\ k_r \end{cases},\quad x_i = \text{an x locus value},\ \mathrm{Re}(n)>1$$

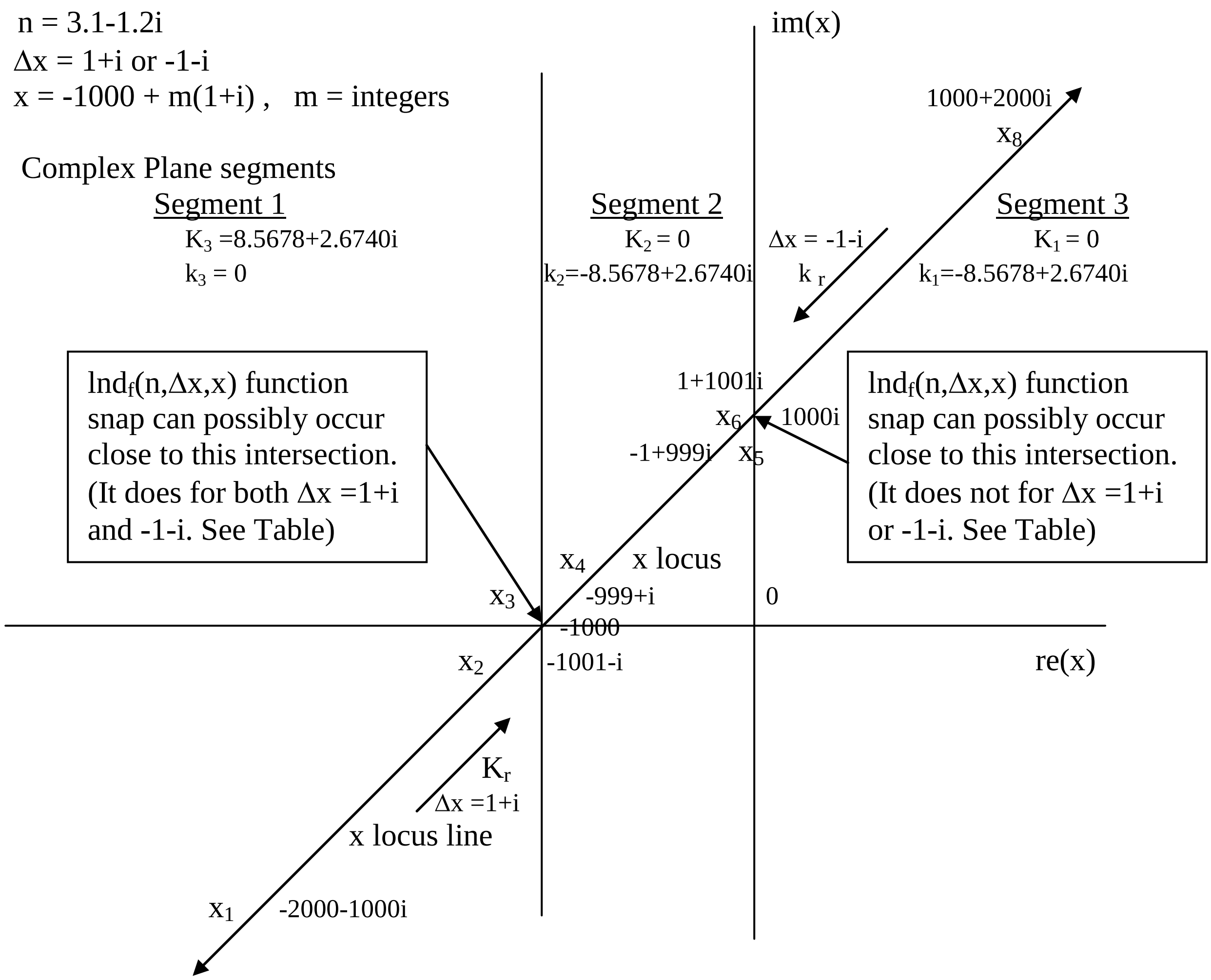

$K_1, K_2, K_3$ are the complex plane segment $lnd_f(n,\Delta x,x)$ constants of integration for those summations moving upward along the x locus line where $\Delta x = 1+i$.
$k_1, k_2, k_3$ are the complex plane segment $lnd_f(n,\Delta x,x)$ constants of integration for those summations moving downward along the x locus line where $\Delta x = -1-i$.

Note – The $lnd_f(n,\Delta x,x)$ constants of integration specified in the above diagram were calculated using the $lnd(n,\Delta x,x)$ calculation program, LNDX. The programming code for this program is shown in the Calculation Programs section at the end of the Appendix.

To demonstrate the validity of Equations E1 thru E12 and Concepts C1 thru C4 calculations will be made and verified using other mathematical means.

**Demonstration #1**

Demonstate the validity of Equations E1 thru E9 and Concepts C1 thru C4.

Equations E4 thru E9 and Concepts C1 Thru C4 have been incorporated into the UBASIC $lnd(n,\Delta x,x)$ calculation program, LNDX, which is presented in the Calculation Programs section at the end of the Appendix.

The program, LNDX, will be used to calculate equations, E1 thru E3 using n = 3.1-1.2i, $\Delta x = 1+i$ or -1-i, the x locus and x locus line shown in the previous diagram, and various values for $x_i$. Also, a computer calculation check of the results of the LNDX program will be preformed to verify accuracy. All values of x will lie on the x locus line shown in the previous diagram.

Comment - The computer calculaton check of the LNDX program sums the terms of the summation being calculated one by one. To sum to infinity in this way and obtain a high precision result is not practical. An excessive amount of computer time would be necessary. The checks performed below verify the results of the LNDX program from four to seven digits of accuracy. This is sufficient to verify the functionality of the LNDX computer program. It is known that the LNDX program provides high precision results.

**1)** Evaluate the summation, $\displaystyle\sum_{\substack{1+i \\ x=\mp\infty}}^{\pm\infty} \frac{1}{x^{3.1-1.2i}}$ , using Eq E3

where

n = 3.1-1.2i
$\Delta x = 1+i$
$x_i = -989+11i$
$x = x_i + m\Delta x$ ,   m = integers

From Eq E3

$$\sum_{\substack{\Delta x \\ x=\mp\infty}}^{\pm\infty} \frac{1}{x^{n}} = \frac{1}{\Delta x}\left[\mathrm{lnd}(n,\Delta x,x_i) - \mathrm{lnd}(n,-\Delta x,x_i-\Delta x)\right]$$

Substituting values

$$\sum_{\substack{1+i \\ x=\mp\infty}}^{\pm\infty} \frac{1}{x^{3.1-1.2i}} = \frac{1}{1+i}\left[\mathrm{lnd}(3.1-1.2i,1+i,-989+11i) - \mathrm{lnd}(3.1-1.2i,-1-i,-990+10i)\right]$$

Using the LNDX program to calculate the lnd(n,Δx,x) functions

$$\sum_{\substack{1+i \\ x=\mp\infty}}^{\pm\infty} \frac{1}{x^{3.1-1.2i}} = \frac{1}{1+i}\left[8.562801571620273 - 2.670413187174567i\right]x10^{-6}$$

$$\sum_{\substack{1+i \\ x=\mp\infty}}^{\pm\infty} \frac{1}{x^{3.1-1.2i}} = \left[2.946194192222852 - 5.616607379397420i\right]x10^{-6}$$

checking

For $x = -989+11i + m(1+i)$ , $m$ = the integers from –100,000 to +100,000

$$\sum_{\substack{1+i \\ x=-100989-99989i}}^{99011+100011i} \frac{1}{x^{3.1-1.2i}} = \left[2.94626 - 5.61660i\right]x10^{-6}$$ Good check

**2)** Evaluate the summation, $\sum_{\substack{1+i \\ x=1102-102i}}^{+\infty} \frac{1}{x^{3.1-1.2i}}$ , using Eq E2

where

$n = 3.1-1.2i$

$\Delta x = 1+i$

$x_i = -1102-102i$

$x = x_i + m\Delta x$ ,   $m = 0,1,2,3,\ldots$

From Eq E2

$$\Delta x \sum_{\Delta x \; x=x_i}^{\pm\infty} \frac{1}{x^n} = \frac{1}{\Delta x} \text{lnd}(n,\Delta x,x_{i)}$$

Substituting values

$$\sum_{1+i \;\; x=-1102-102i}^{+\infty} \frac{1}{x^{3.1-1.2i}} = \frac{1}{1+i} [\text{lnd}(3.1\text{-}1.2i,1+i,-1102\text{-}102i)]$$

Using the LNDX program to calculate the lnd(n,Δx,x) function

$$\sum_{1+i \;\; x=-1102-102i}^{+\infty} \frac{1}{x^{3.1-1.2i}} = \frac{1}{1+i} [2.520092396540544-2.947396943161950i]x10^{-6}$$

$$\sum_{1+i \;\; x=-1102-102i}^{+\infty} \frac{1}{x^{3.1-1.2i}} = [1.112676351112174-1.407416045428370i]x10^{-6}$$

checking

For x = -1102-102i + m(1+i) , m = the integers from 0 to 200,000

$$\sum_{1+i \;\; x=-1102-102i}^{198898+199898i} \frac{1}{x^{3.1-1.2i}} = [1.112676 - 1.407416i]x10^{-6}$$ Good check

**3)** Evaluate the summation, $\sum_{-1-i \;\; x=1000i}^{-\infty} \frac{1}{x^{3.1-1.2i}}$ , using Eq E2

where

$n = 3.1-1.2i$

$\Delta x = -1-i$

$x_i = 1000i$

$x = x_i + m\Delta x$ , m = 0,1,2,3,…

From Eq E2

$$\Delta x \; {}_{\Delta x}\sum_{x=x_i}^{\pm\infty} \frac{1}{x^n} = \frac{1}{\Delta x}\, \mathrm{lnd}(n,\Delta x,x_{i)}$$

Substituting values

$${}_{-1-i}\sum_{x=1000i}^{-\infty} \frac{1}{x^{3.1-1.2i}} = \frac{1}{1+i}\,[\mathrm{lnd}(3.1\text{-}1.2i,\text{-}1\text{-}i,1000i)]$$

Using the LNDX program to calculate the lnd(n,Δx,x) function

$${}_{-1-i}\sum_{x=1000i}^{-\infty} \frac{1}{x^{3.1-1.2i}} = \frac{1}{1+i}\,[-8.540262768319782+2.648481133711169i]\mathrm{x}10^{-6}$$

$${}_{-1-i}\sum_{x=1000i}^{-\infty} \frac{1}{x^{3.1-1.2i}} = [2.945890817304306+5.594371951015475i]\mathrm{x}10^{-6}$$

checking

For $x = 1000i + m(-1-i)$ , m = the integers from 0 to 200,000

$${}_{-1-i}\sum_{x=1000i}^{-200000-199000i} \frac{1}{x^{3.1-1.2i}} = [2.9459 - 5.5943i]\mathrm{x}10^{-6} \qquad \text{Good check}$$

**4)** Evaluate the summation, ${}_{1+i}\sum_{x=-1001-i}^{1+1001i} \frac{1}{x^{3.1-1.2i}}$ , using Eq E1

where

$n = 3.1-1.2i$

$\Delta x = 1+i$

$x_1 = -1001-i$

$x_2 = 1+1001i$

$x = x_1 + m\Delta x$ , m = the integers from 0 to 1002 (i.e. $\frac{x_2-x_1}{\Delta x}$ )

From Eq E1

$$\Delta x \sum_{\Delta x\, x=x_1}^{x_2} \frac{1}{x^n} = -\mathrm{lnd}(n,\Delta x,x)\Big|_{x_1}^{x_2+\Delta x} = -\mathrm{lnd}(n,\Delta x,x_2+\Delta x) + \mathrm{lnd}(n,\Delta x,x_1)]$$

Substituting values

$$\sum_{1+i\; x=-1001-i}^{1+1001i} \frac{1}{x^{3.1-1.2i}} = \frac{1}{1+i}\,[-\mathrm{lnd}(3.1-1.2i,1+i,2+1002i) + \mathrm{lnd}(3.1-1.2i,1+i,-1001-i]$$

Using the LNDX program to calculate the lnd(n,Δx,x) functions

$$\sum_{1+i\; x=-1001-i}^{1+1001i} \frac{1}{x^{3.1-1.2i}} = \frac{1}{1+i}\,[.5152442710026848+2.450653840064823]x10^{-8}$$

$$\sum_{1+i\; x=-1001-i}^{1+1001i} \frac{1}{x^{3.1-1.2i}} = [1.482949055533754+.9677047845310694]x10^{-8}$$

checking

$$\sum_{1+i\; x=-1001-i}^{1+1001i} \frac{1}{x^{3.1-1.2i}} = [1.482949055533754+.9677047845310694]x10^{-8} \qquad \text{Good check}$$

The LNDX program has calculated the values of the lnd(n,Δx,x) functions correctly to good accuracy. The validity of Equations E1 thru E9, Concepts C1 thru C4 and the LNDX program itself has been demonstrated.

**Demonstration #2**

Demonstate the validity of Equation E10

$K_3 = -k_1$.

Note columns G and H in the following table of calculated values relating to the diagram x locus. The value of $K_3$ = 8.5628-2.6704i and the value of $k_1$ = -8.5628-2.6704i , $K_3 = -\,k_1$. The validity of Equation E10 has been demonstrated.

**Demonstration #3**

Demonstate the validity of Equation E11

$lnd(n,\Delta x,x) - lnd(n,-\Delta x,x-\Delta x) = \pm(K_r - k_r)$

where

$r = 1,2,3$

$+$ for $Re(\Delta x)>0$ or $\{Re(\Delta x)=0$ and $Im(\Delta x)>0\}$

$-$ is for $Re(\Delta x)<0$ or $\{Re(\Delta x)=0$ and $Im(\Delta x)<0\}$

$n,\Delta x,K_r,k_r$ = real or complex constants

Note columns C, E and F in the following table of calculated values relating to the diagram x locus. The first half of columns C, E, and F represents the condition where $\Delta x = 1+i$ and the second half of columns C,E, and F represents the condition where $\Delta x = -(1+i)$. Where $Re(1+i)>0$, $+(K_r - k_r)$, is seen to equal $lnd(n,\Delta x,x) - lnd(n,-\Delta x,x-\Delta x)$ and where $Re(-1-i)<0$, $-(K_r - k_r)$ is seen to equal

$lnd(n,\Delta x,x) - lnd(n,-\Delta x,x-\Delta x)$. The validity of Equation E11 has been demonstrated.

Comment - It is interesting to note that for all points on an x locus line with a specified value of n, the value of $K_r - k_r$ is the same. This is in agreement with Eq 2.6-33 which states that $K_r - k_r$ = constant value all along an x locus line with a specified value of n.

**Demonstration #4**

Demonstate the validity of Equation E3.

$$\Delta x \sum_{\substack{x=\mp\infty \\ \Delta x}}^{\pm\infty} \frac{1}{x^n} = lnd(n,\Delta x,x_i) - lnd(n,-\Delta x,x_i-\Delta x) = \pm(K_r - k_r)$$

The validity of $lnd(n,\Delta x,x_i) - lnd(n,-\Delta x,x_i-\Delta x) = \pm(K_r - k_r)$ has been shown in Demonstration #3.

Demonstrate the validity of

$$\Delta x \sum_{\substack{x=\mp\infty \\ \Delta x}}^{\pm\infty} \frac{1}{x^n} = lnd(n,\Delta x,x_i) - lnd(n,-\Delta x,x_i-\Delta x)$$

where

$x = x_i + m\Delta x$ , $m$ = integers , This is an x locus in the complex plane

$x_i$ = any value of x in the x locus

$\Delta x$ = x increment

$-\infty$ to $+\infty$ for $Re(\Delta x)>0$ or $\{Re(\Delta x)=0$ and $Im(\Delta x)>0\}$

$+\infty$ to $-\infty$ for $Re(\Delta x)<0$ or $\{Re(\Delta x)=0$ and $Im(\Delta x)<0\}$

x = real or complex values

$n,\Delta x,x_i,K_r,k_r$ = real or complex constants

$Re(n)>1$

$r = 1,2,3$

Note columns B,C, D, F and J in the following table of calculated values relating to the diagram x locus. For the same values of n, Δx, x the calculated values in columns F and J are seen to be equal. The validity of Equation E3 has been demonstrated.

**Demonstration #5**

Demonstate the validity of Equation E12.

$$\mathrm{lnd}(n,\Delta x,x) - \mathrm{lnd}(n,-\Delta x,x-\Delta x) = -[\mathrm{lnd}(n,-\Delta x,x) - \mathrm{lnd}(n,\Delta x,x+\Delta x)]$$

where

$x = x_i + m\Delta x$ ,  m = integers , This is an x locus in the complex plane
$x_i$ = a value of x
Δx = x increment
$n \neq 1$
n,Δx = real or complex constants

Note columns B,C,D and F in the following table of calculated values relating to the diagram x locus. For a specified value of n and x it is observed that the calculated value of lnd(n,Δx,x) – lnd(n,-Δx,x-Δx) reverses sign if the value of Δx reverses sign. The validity of Equation E3 has been demonstrated.

**Demonstration #6**

Demonstate the validity of Equation E13.

$$k_1 = \pm\left[-\Delta x \sum_{\Delta x \; x=\mp\infty}^{\pm\infty} \frac{1}{x^n}\right], \quad \mathrm{Re}(n)>1$$

where

$x = x_i + m\Delta x$ ,  m = integers , This is an x locus in the complex plane
$x_i$ = a value of x
Δx = x increment
+ and -∞ to +∞ for Re(Δx)>0 or {Re(Δx)=0 and Im(Δx)>0}
– and +∞ to -∞ for Re(Δx)<0 or {Re(Δx)=0 and Im(Δx)<0}
$n \neq 1$
$n,\Delta x,k_1$ = real or complex constants

Re(n)>1
Δx = -1-i
Re(Δx)<0

$$k_1 = \Delta x \sum_{\Delta x \; +\infty}^{-\infty} \frac{1}{x^n}$$

In the table, the calculated value in J14 is $\Delta x \sum_{\Delta x \; +\infty}^{-\infty} \frac{1}{x^n} = [-8.5628+2.6704i]\mathrm{x}10^{-6}$.

Substituting this value into the above equation, $k_1 = [-8.5628+2.6704i]x10^{-6}$.
This value agrees with the $\text{lnd}(n,\Delta x,x) - \text{lnd}_f(n,\Delta x,x)$ calculated value for $k_1$ in H14.
The validity of Equation E13 has been demonstrated.

**Demonstration #7**

Demonstate the validity of Equation E14.

$$K_3 = \pm\left[\Delta x \sum_{\Delta x \; x=\mp\infty}^{\pm\infty} \frac{1}{x^n}\right] , \quad Re(n)>1$$

where

$x = x_i + m\Delta x$ , $m$ = integers , This is an x locus in the complex plane
$x_i$ = a value of x
$\Delta x$ = x increment
+ and $-\infty$ to $+\infty$ for $Re(\Delta x)>0$ or $\{Re(\Delta x)=0$ and $Im(\Delta x)>0\}$
– and $+\infty$ to $-\infty$ for $Re(\Delta x)<0$ or $\{Re(\Delta x)=0$ and $Im(\Delta x)<0\}$
$n \neq 1$
$n,\Delta x,K_3$ = real or complex constants

$Re(n)>1$
$\Delta x = 1+i$
$Re(\Delta x)>0$

$$K_3 = \Delta x \sum_{\Delta x \; -\infty}^{+\infty} \frac{1}{x^n}$$

In the table, the calculated value in J1 is $\Delta x \sum_{\Delta x \; -\infty}^{+\infty} \frac{1}{x^n} = [8.5628+2.6704i]x10^{-6}$.

Substituting this value into the above equation, $K_3 = [8.5628+2.6704i]x10^{-6}$.
This value agrees with the $\text{lnd}(n,\Delta x,x) - \text{lnd}_f(n,\Delta x,x)$ calculated value for $K_3$ in H1.
The validity of Equation E14 has been demonstrated.

**Demonstration #8**

Demonstate the validity of Equation E15.

$$K_2 - k_2 = \pm\left[\Delta x \sum_{\Delta x \; x=\mp\infty}^{\pm\infty} \frac{1}{x^n}\right] , \quad Re(n)>1$$

where

$x = x_i + m\Delta x$ , $m$ = integers , This is an x locus in the complex plane
$x_i$ = a value of x
$\Delta x$ = x increment

+ and $-\infty$ to $+\infty$ for $Re(\Delta x)>0$ or $\{Re(\Delta x)=0$ and $Im(\Delta x)>0\}$
– and $+\infty$ to $-\infty$ for $Re(\Delta x)<0$ or $\{Re(\Delta x)=0$ and $Im(\Delta x)<0\}$
$n \neq 1$
$n,\Delta x,K_2,k_2$ = real or complex constants

$Re(n)>1$
$\Delta x = 1+i$
$Re(\Delta x)>0$

$$K_2 - k_2 = \Delta x \sum_{\Delta x\, -\infty}^{+\infty} \frac{1}{x^n}$$

In the table, the calculated value in J5 is $\Delta x \sum_{\Delta x\, -\infty}^{+\infty} \frac{1}{x^n} = [8.5628+2.6704i]x10^{-6}$.

Substituting this value into the above equation, $K_2 - k_2 = [8.5628+2.6704i]x10^{-6}$.
This value agrees with the $lnd(n,\Delta x,x) - lnd_f(n,\Delta x,x)$ calculated value for $K_2 - k_2$ in F5.
The validity of Equation E15 has been demonstrated.

**Demonstration #9**
Demonstate the validity of Equation E16.

$$K_r - k_r = \pm [\Delta x \sum_{\Delta x\, x=\mp\infty}^{\pm\infty} \frac{1}{x^n}] , \quad Re(n)>1$$

where
$x = x_i + m\Delta x$ , $m$ = integers , This is an x locus in the complex plane
$x_i$ = a value of x
$\Delta x$ = x increment
+ and $-\infty$ to $+\infty$ for $Re(\Delta x)>0$ or $\{Re(\Delta x)=0$ and $Im(\Delta x)>0\}$
– and $+\infty$ to $-\infty$ for $Re(\Delta x)<0$ or $\{Re(\Delta x)=0$ and $Im(\Delta x)<0\}$
$n \neq 1$
$n,\Delta x$ = real or complex constants
$r = 1,2,3$

$Re(n)<1$
$\Delta x = -1-i$
$Re(\Delta x)<0$

$$K_r - k_r = \Delta x \sum_{\Delta x\, +\infty}^{-\infty} \frac{1}{x^n}$$

In the table, the calculated values in J9 thru J16 are $\Delta x \sum_{\Delta x}^{+\infty}_{-\infty} \frac{1}{x^n} = [-8.5628+2.6704i]x10^{-6}$.

Substituting these values into the above equation, $K_r - k_r = [-8.5628+2.6704i]x10^{-6}$.
These values agree with the $lnd(n,\Delta x,x) - lnd_f(n,\Delta x,x)$ calculated values for $K_r - k_r$ in F9 thru F16.
The validity of Equation E16 has been demonstrated.

**Demonstration #10**

Demonstate the validity of Equation E17.

$$\left.\begin{matrix}K_r\\k_r\end{matrix}\right\} = \Delta x \sum_{\Delta x \; x=x_i}^{\pm\infty} \frac{1}{x^n} - lnd_f(n,\Delta x,x_i)\ , \quad Re(n)>1$$

where

$x = x_i + m\Delta x$ , $m$ = integers , This is an x locus in the complex plane
$x_i$ = values of x
$x_i$ is within the x locus segment, r
$x_i$ to $+\infty$ for $Re(\Delta x)>0$ or $\{Re(\Delta x)=0$ and $Im(\Delta x)>0\}$
$x_i$ to $-\infty$ for $Re(\Delta x)<0$ or $\{Re(\Delta x)=0$ and $Im(\Delta x)<0\}$
$\Delta x$ = x increment
x = real or complex values
$n,\Delta x,x_i,K_r,k_r$ = real or complex constants
r = 1,2,3 , The x locus segment designations
$K_r$ = constant of integration for $Re(\Delta x)>0$ or $\{Re(\Delta x) = 0$ and $Im(\Delta x) > 0\}$
$k_r$ = constant of integration for $Re(\Delta x)<0$ or $\{Re(\Delta x) = 0$ and $Im(\Delta x) < 0\}$
$n \neq 1$
Accuracy increases rapidly as $\left|\frac{x}{\Delta x}\right|$ increases in value.

In the table, the calculated values in column K are $\left.\begin{matrix}K_r\\k_r\end{matrix}\right\} = \Delta x \sum_{\Delta x \; x=x_i}^{\pm\infty} \frac{1}{x^n} - lnd_f(n,\Delta x,x_i)$ where $x_i$ are the x values in column D. In the table, the calculated values in column H are $\left.\begin{matrix}K_r\\k_r\end{matrix}\right\} = lnd(n,\Delta x,x) - lnd_f(n,\Delta x,x)$ . Note that the calculated values in column K and in column H are the same. The validity of Equation E17 has been demonstrated.

## Table - Calculations relating to the x locus diagram

| Row No. | Pnt | n | Δx | x | $\pm(K_r-k_r)$ $lnd(n,\Delta x,x) - lnd(n,-\Delta x,x-\Delta x)$ $x10^{-6}$ | | $K_r$ or $k_r$ $lnd(n,\Delta x,x) - lnd_f(n,\Delta x,x)$ $x10^{-6}$ | | $\pm(K_r-k_r)$ $\Delta x \sum_{x=\mp\infty}^{\pm\infty} \frac{1}{x^n}$ $x=x_i+m\Delta x$ , m = integers; $x_i$ = x in column D $x10^{-6}$ | | $K_r$ or $k_r$ (See column G) $\Delta x \sum_{x=x_i}^{\pm\infty} \frac{1}{x^n} - lnd_f(n,\Delta x,x_i)$ Re(n)>1; $x_i$ = x in column D $x10^{-6}$ |
|---|---|---|---|---|---|---|---|---|---|---|---|
| | A | B | C | D | E | F | G | H | I | J | K |
| **1** | $x_1$ | 3.1-1.2i | 1+i | -2000-1000i | $K_3$-$k_3$ | 8.5628-2.6704i | $K_3$ | 8.5628-2.6704i | $K_3$-$k_3$ | 8.5628-2.6704i | 8.5628-2.6704i |
| **2** | $x_2$ | 3.1-1.2i | 1+i | -1001-i | $K_3$-$k_3$ | 8.5628-2.6704i | $K_3$ | 8.5628-2.6704i | $K_3$-$k_3$ | 8.5628-2.6704i | 8.5628-2.6704i |
| **3** | $x_3$ | 3.1-1.2i | 1+i | -1000 | $K_3$-$k_3$ | 8.5628-2.6704i | $K_3$ | 8.5628-2.6704i | $K_3$-$k_3$ | 8.5628-2.6704i | 8.5628-2.6704i |
| **4** | $x_4$ | 3.1-1.2i | 1+i | -999+i | $K_2$-$k_2$ | 8.5628-2.6704i | $K_2$ | 0 snap | $K_2$-$k_2$ | 8.5628-2.6704i | 0 snap |
| **5** | $x_5$ | 3.1-1.2i | 1+i | -1+999i | $K_2$-$k_2$ | 8.5628-2.6704i | $K_2$ | 0 | $K_2$-$k_2$ | 8.5628-2.6704i | 0 |
| **6** | $x_6$ | 3.1-1.2i | 1+i | 1000i | $K_2$-$k_2$ | 8.5628-2.6704i | $K_2$ | 0 | $K_2$-$k_2$ | 8.5628-2.6704i | 0 |
| **7** | $x_7$ | 3.1-1.2i | 1+i | 1+1001i | $K_1$-$k_1$ | 8.5628-2.6704i | $K_1$ | 0 | $K_1$-$k_1$ | 8.5628-2.6704i | 0 |
| **8** | $x_8$ | 3.1-1.2i | 1+i | 1000+2000i | $K_1$-$k_1$ | 8.5628-2.6704i | $K_1$ | 0 | $K_1$-$k_1$ | 8.5628-2.6704i | 0 |
| **9** | $x_1$ | 3.1-1.2i | -1-i | -2000-1000i | $k_3$-$K_3$ | -8.5628+2.6704i | $k_3$ | 0 | $k_3$-$K_3$ | -8.5628+2.6704i | 0 |
| **10** | $x_2$ | 3.1-1.2i | -1-i | -1001-i | $k_3$-$K_3$ | -8.5628+2.6704i | $k_3$ | 0 snap | $k_3$-$K_3$ | -8.5628+2.6704i | 0 snap |
| **11** | $x_3$ | 3.1-1.2i | -1-i | -1000 | $k_2$-$K_2$ | -8.5628+2.6704i | $k_2$ | -8.5628+2.6704i | $k_2$-$K_2$ | -8.5628+2.6704i | -8.5628+2.6704i |
| **12** | $x_4$ | 3.1-1.2i | -1-i | -999+i | $k_2$-$K_2$ | -8.5628+2.6704i | $k_2$ | -8.5628+2.6704i | $k_2$-$K_2$ | -8.5628+2.6704i | -8.5628+2.6704i |
| **13** | $x_5$ | 3.1-1.2i | -1-i | -1+999i | $k_2$-$K_2$ | -8.5628+2.6704i | $k_2$ | -8.5628+2.6704i | $k_2$-$K_2$ | -8.5628+2.6704i | -8.5628+2.6704i |
| **14** | $x_6$ | 3.1-1.2i | -1-i | 1000i | $k_1$-$K_1$ | -8.5628+2.6704i | $k_1$ | -8.5628+2.6704i | $k_1$-$K_1$ | -8.5628+2.6704i | -8.5628+2.6704i |
| **15** | $x_7$ | 3.1-1.2i | -1-i | 1+10001i | $k_1$-$K_1$ | -8.5628+2.6704i | $k_1$ | -8.5628+2.6704i | $k_1$-$K_1$ | -8.5628+2.6704i | -8.5628+2.6704i |
| **16** | $x_8$ | 3.1-1.2i | -1-i | 1000+2000i | $k_1$-$K_1$ | -8.5628+2.6704i | $k_1$ | -8.5628+2.6704i | $k_1$-$K_1$ | -8.5628+2.6704i | -8.5628+2.6704i |

+ and $K_r$ are used when Re(Δx)>0 or {Re(Δx)=0 and Im(x)>0}

– and $k_r$ are used when Re(Δx)<0 or {Re(Δx)=0 and Im(x)<0}

r = 1,2,3

Thus, from the calculations of this example, the validity of Equations E1 thru E17, the Concepts C1 thru C4 and the functionality of the LNDX lnd(n,Δx,x) function calculation program have been demonstrated.

Comment - The LNDX program calculates but does not print out the $\text{lnd}_f(n,\Delta x,x)$ function and the lnd(n,Δx,x) function constants of integration, $K_1,K_2,K_3,k_1,k_2,k_3$. To obtain these values, some additional code was added to the LNDX program to print them out. This modified LNDX program is called FINDK.

**Section 2.7: Method 1 for the the calculation of the function, lnd(n,Δx,x) n≠1**

In the previous section, Section 2.6, it was shown that the function, lnd(n,Δx,x) n≠1 is very nearly equal to the series, $\text{lnd}_f(n,\Delta x,x)$, under certain conditions (where $K_1=k_3=0$). Pertinent function relationships and the equality conditions are restated below:

$$\text{lnd}(n,\Delta x,x_i) = \Delta x \sum_{\substack{\Delta x \\ x=x_i}}^{\pm\infty} \frac{1}{x^n} \approx \text{lnd}_f(n,\Delta x,x_i) \tag{2.7-1}$$

where

$x = x_i + m\Delta x$, m = integers , This is the x locus

$x_i$ = summation initial value of x

Δx = x increment

n≠1

+ for Re(Δx)>0 or {Re(Δx)=0 and Im(Δx)>0}

– for Re(Δx)>0 or {Re(Δx)=0 and Im(Δx)>0}

The value of $\left|\frac{x_i}{\Delta x}\right|$ is large

The $\text{lnd}_f(n,\Delta x,x)$ Series

$$\text{lnd}_f(n,\Delta x,x) \approx -\sum_{m=0}^{\infty} \frac{\Gamma(n+2m-1)\left(\frac{\Delta x}{2}\right)^{2m} C_m}{\Gamma(n)(2m+1)!\left(x-\frac{\Delta x}{2}\right)^{n+2m-1}}, \quad n\neq 1 \tag{2.7-2}$$

where

$C_m$ = Series constants

The accuracy of Eq 2.7-2 above increases rapidly for increasing $|\frac{x}{\Delta x}|$ .

$$\text{lnd}(n,\Delta x,x_i) = \Delta x \sum_{\substack{\Delta x \\ x=x_i}}^{\pm\infty} \frac{1}{x^n} \approx \text{lnd}_f(n,\Delta x,x_i) \tag{2.7-3}$$

The accuracy of Eq 2.7-3 increases rapidly as $|\frac{x_i}{\Delta x}|$ increases in value**.**

where

For segment 1 on the $\text{lnd}(n,\Delta x,x)$ summation x locus line

$Re(\Delta x)>0$, $Re(x)>x_{rcr}$

or

$Re(\Delta x)=0$, $Im(\Delta x)>0$, $Im(x)>0$

$x_i$ = value of x on segment 1 of the $\text{lnd}(n,\Delta x,x)$ summation x locus line

$x_{rc}$ = rightmost or only x locus line and complex plane axis crossover point

$x_{rcr}$ = real part of $x_{rc}$

For segment 3 on the $\text{lnd}(n,\Delta x,x)$ summation x locus line

$Re(\Delta x)<0$, $Re(x)<x_{lcr}$

or

$Re(\Delta x)=0$, $Im(\Delta x)<0$, $Im(x)<0$

$x_i$ = value of x on segment 3 of the $\text{lnd}(n,\Delta x,x)$ summation x locus line

$x_{lc}$ = leftmost or only x locus line and complex plane axis crossover point

$x_{lcr}$ = real part of $x_{lc}$

See Diagram 2.6-3

To obtain a very near equality relationship in Eq 2.7-3 for all $x_i$, the accuracy requirement involving the quantity, $|\frac{x_i}{\Delta x}|$, must be dealt with. Much thought and many trials, yielded an algorithm which provides the very near equality desired. This algorithm assures that the value of $x_i$ used in the $\text{lnd}_f(n,\Delta x,x)$ n≠1 Series has real and imaginary values greater than a minimum considered necessary for high accuracy. These minimum real and imaginary values are variable and a function of Δx. To implement the algorithm, the x locus segments where $K_1$ and $k_3$ equal zero are divided into two parts, 1a and 1b and 3a and 3b (See Diagram 2.7-1). Only $x_i$ values on the x locus segments 1b and 3b are suitable for use with the lnd(n,Δx,x) n≠1 Series. All $x_i$ values on x locus line segments 1b and 3b will be free of potential $\text{lnd}_f(n,\Delta x,x_i)$ n≠1 Series snap irregularities and will yield an $|\frac{x_i}{\Delta x}|$ value which is adequately large for high accuracy. The algorithm selects a G value for two conditions. They are as follows:

For Re(Δx)≠0

$G = N\text{Re}(\Delta x)$

For Re(Δx)=0

$G = N\text{Im}(\Delta x)$

$N \geq 100$, This value for N has been found to obtain accurate results

Diagram 2.7-1 below shows how these G values are used in selecting x locus $x_i$ values that will yield accurate lnd(n,Δx,x) n≠1 function values from $\text{lnd}_f(n,\Delta x,x)$ Series calculations. Values of $x_i$ on the highlighted x locus line segments may be used directly in the lnd(n,Δx,x) n≠1 Series to obtain highly accurate values of the function, lnd(n,Δ,x) n≠1. In this case:

$$\text{lnd}(n,\Delta x,x_i) = \text{lnd}_f(n,\Delta x,x_i) = \Delta x \sum_{\Delta x \; x=x_i}^{\pm\infty} \frac{1}{x^n}, \quad n\neq 1 \qquad (2.7\text{-}4)$$

However, the desired value of $x_i$ may be on the x locus line but not on the segment which is highlighted. This condition is handled using the following equation:

$$\text{lnd}(n,\Delta x,x_i) = \Delta x \sum_{\Delta x \; x=x_i}^{x_p-\Delta x} \frac{1}{x^n} + \text{lnd}_f(n,\Delta x,x_p) = \Delta x \sum_{\Delta x \; x=x_i}^{\pm\infty} \frac{1}{x^n}, \quad n\neq 1 \qquad (2.7\text{-}5)$$

where

$x = x_i + m\Delta x$, m = integers , This is the x locus

$x_p$ = an x value on the x locus segment 1b for Re(Δx)>0 or {Re(Δx)=0 and Im(Δx)>0} or an x value on the x locus segment 3b for Re(Δx)<0 or {Re(Δx)=0 and Im(Δx)<0}

$\Delta x$ = x increment
+ for $Re(\Delta x)>0$ or $\{Re(\Delta x)=0$ and $Im(\Delta x)>0\}$
– for $Re(\Delta x)<0$ or $\{Re(\Delta x)=0$ and $Im(\Delta x)<0\}$

The summation term, in effect, jumps over regions of the x locus line which may be subject to lnd(n,Δx,x) n≠1 Series snap irregularities or low $|\frac{x}{\Delta x}|$ which degrades lnd(n,Δx,x) n≠1 Series accuracy.

In Diagram 2.7-1 below, there are six figures, Fig 1 thru Fig 6, representing various examples of x loci. On each x locus line is marked an $x_i$ value. The equation to calculate the function, lnd(n,Δx,x) n≠1 is shown at the bottom of each figure.

In Fig 1, Fig 2, Fig 5, and Fig 6 lnd(n,Δx,x) is calculated using Eq 2.7-5. In Fig 3 and Fig 4 lnd(n,Δx,x) is able to be calculated directly using Eq 2.7-4.

Diagram 2.7-1: Various examples of x locii

For $Re(\Delta x)>0$ or $\{Re(\Delta x)=0$ and $Im(\Delta x)>0$

Fig 1

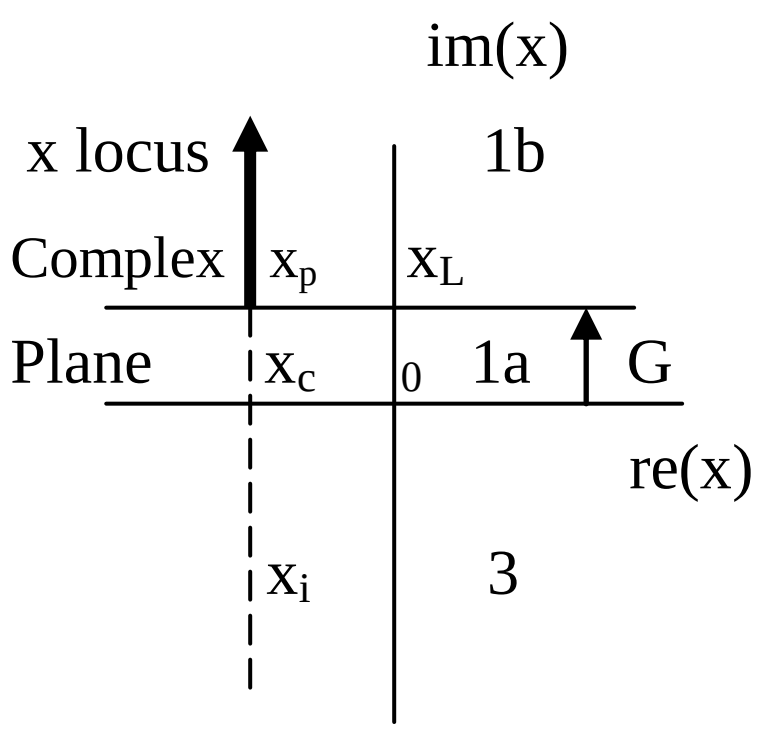

$Re(\Delta x)=0$ and $Im(x)>0$

$G = NIm(\Delta x)$

$x_L = 0+jG$

Fig 2

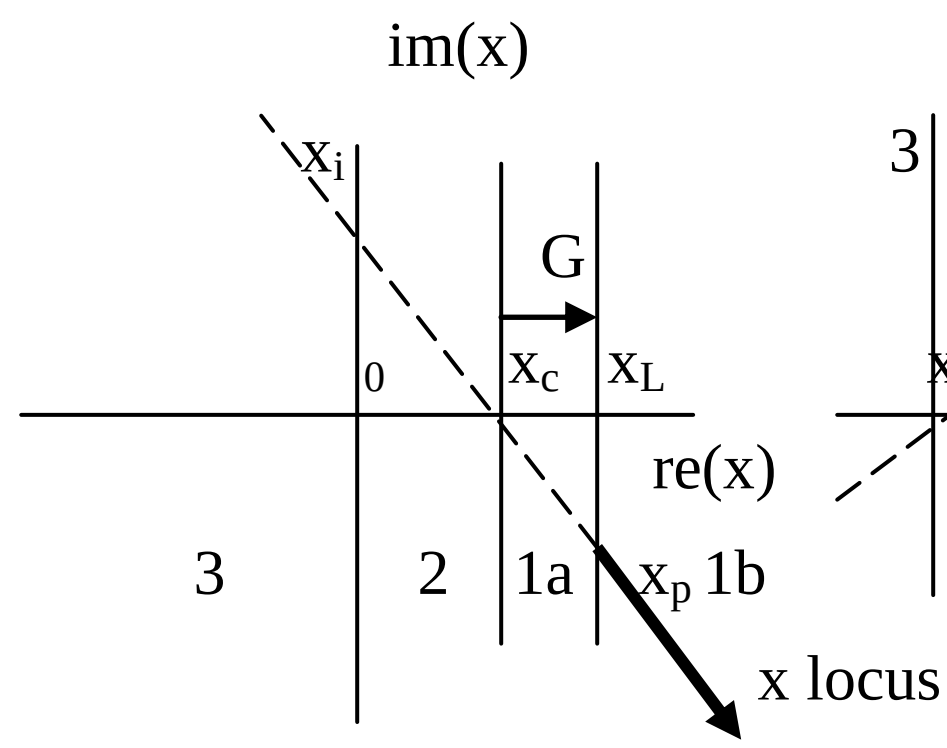

$Re(\Delta x)>0$

$G = NRe(\Delta x)$

$x_L = x_c+G$

Fig 3

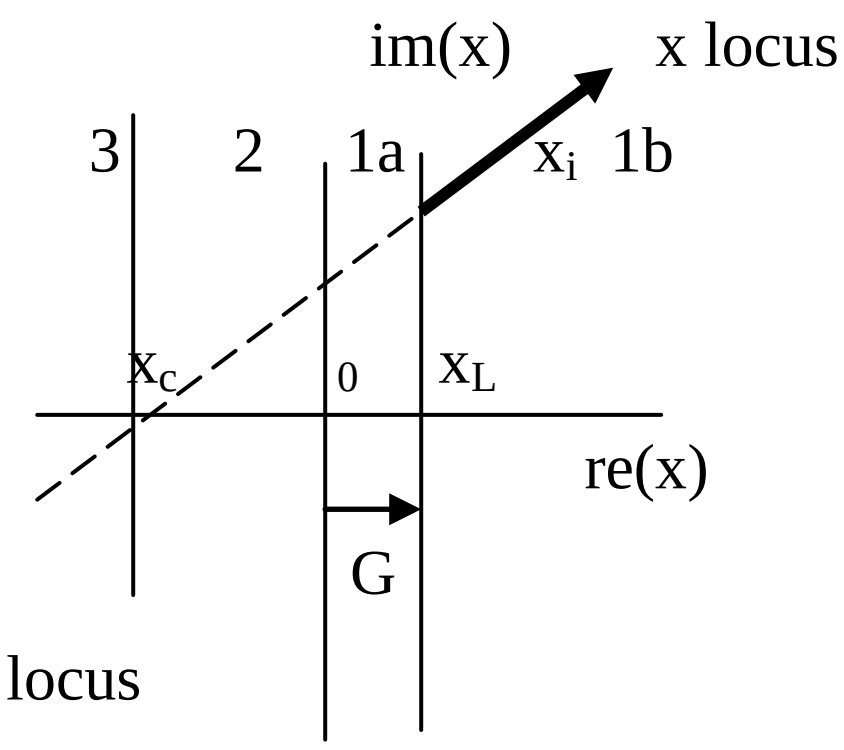

$Re(\Delta x)>0$

$G = NRe(\Delta x)$

$x_L = 0+G$

$$lnd(n,\Delta x,x_i) = \Delta x \sum_{x=x_i}^{x_p-\Delta x}{}_{\Delta x} \frac{1}{x^n} + lnd_f(n,\Delta x,x_p) \quad lnd(n,\Delta x,x_i) = \Delta x \sum_{x=x_i}^{x_p-\Delta x}{}_{\Delta x} \frac{1}{x^n} + lnd_f(n,\Delta x,x_p) \quad lnd(n,\Delta x,x_i) = lnd_f(n,\Delta x,x_i)$$

For $Re(\Delta x)<0$ or $\{Re(\Delta x)=0$ and $Im(\Delta x)<0$

Fig 4

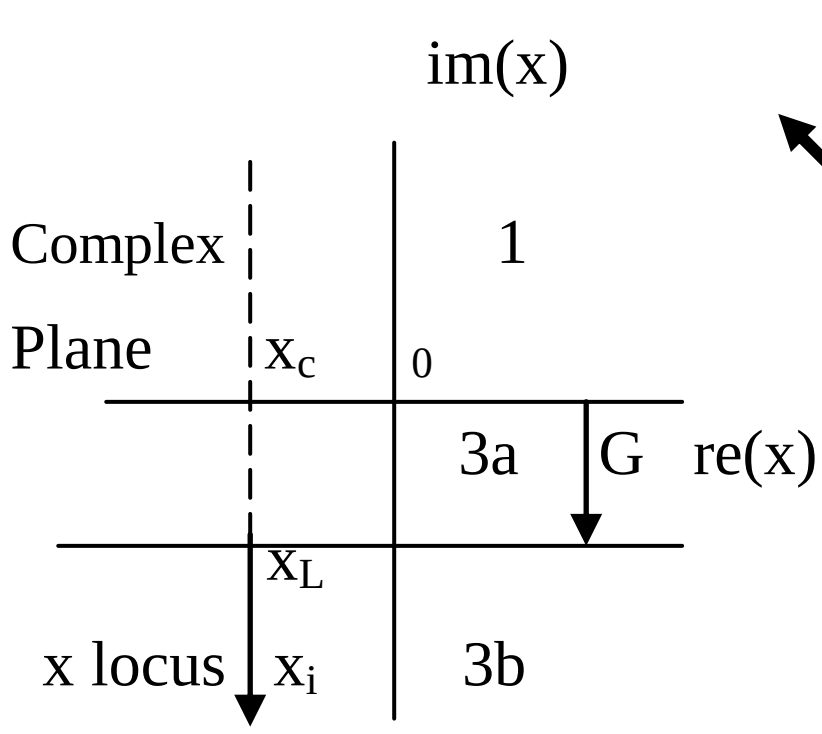

$Re(\Delta x)=0$ and $Im(x)<0$

$G = NIm(\Delta x)$

$x_L = 0+jG$

Fig 5

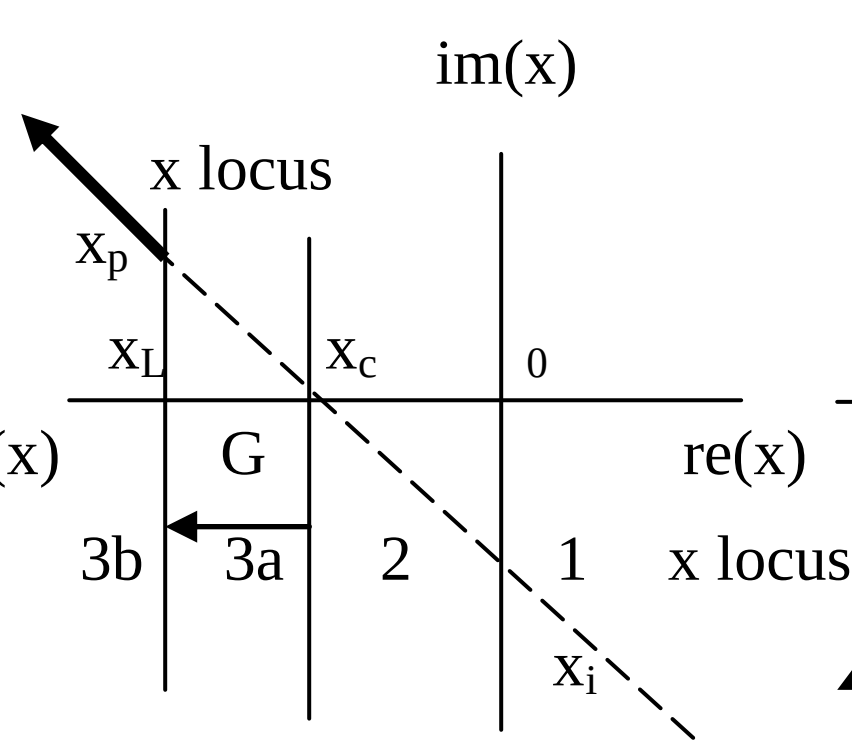

$Re(\Delta x)<0$

$G = NRe(\Delta x)$

$x_L = \ x_c+G$

Fig 6

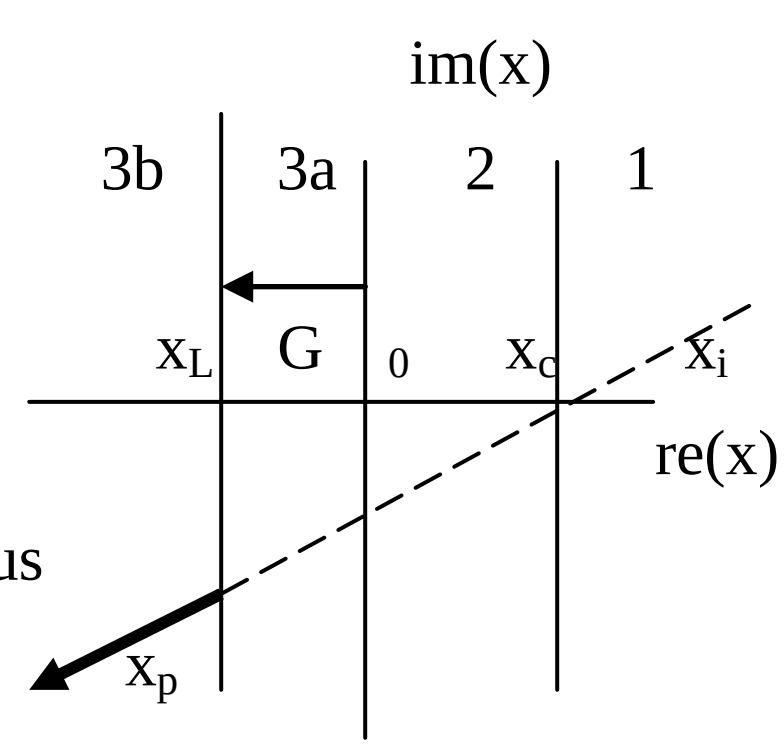

$Re(\Delta x),<0$

$G = NRe(\Delta x)$

$x_L = 0+G$

$$lnd(n,\Delta x,x_i) = lnd_f(n,\Delta x,x_i) \quad lnd(n,\Delta x,x_i) = \Delta x \sum_{x=x_i}^{x_p-\Delta x}{}_{\Delta x} \frac{1}{x^n} + lnd_f(n,\Delta x,x_p) \quad lnd(n,\Delta x,x_i) = \Delta x \sum_{x=x_i}^{x_p-\Delta x}{}_{\Delta x} \frac{1}{x^n} + lnd_f(n,\Delta x,x_p)$$

$N \geq 100$,  This value for N has been found to obtain accurate results

$x = x_i + m\Delta x$ ,  $m = 1,2,3,\ldots$  This is the x locus

$x = x_i, x_i+\Delta x, x_i+2\Delta x, x_i+3\Delta x, \ldots$

1a,1b,2,3a,3b represent x locus segment designations

$x_p$ = an x value on the x locus segment 1b for $Re(\Delta x)>0$ or $\{Re(\Delta x)=0$ and $Im(\Delta x)>0\}$ or

an x value on the x locus segment 3b for $Re(\Delta x)<0$ or $\{Re(\Delta x)=0$ and $Im(\Delta x)<0\}$

x locus line segment 1b is where: $Re(\Delta x)>0,\ Re(x)>x_L$
or
$Re(\Delta x)=0,\ Im(\Delta x)>0,\ Im(x)>x_L$

x locus line segment 3b is where: $Re(\Delta x)<0,\ Re(x)<x_L$
or
$Re(\Delta x)=0,\ Im(\Delta x)<0,\ Im(x)<x_L$

The lnd(n,Δx,x) n≠1 calculation equation, Eq 2.7-5, while being very useful, has a significant drawback. For $\left|\frac{Re(x) - Re(x_p)}{Re(\Delta x)}\right|$ being a very large number, the summation part of Eq 2.7-5 will require many calculations. Many computer calculations can take an excessive amount of time. To eliminate this shortcoming, a Method 2 has been devised to calculate the function, lnd(n,Δx,x) n≠1, which does use Eq 2.7-5 but only with short summations. This Method 2 is more complex than Method 1. It requires that the snap phenomenon involving the lnd(n,Δx,x) n≠1 Series be confronted even though it is not presently fully understood. The method of calculation is based on the lnd(n,Δx,x) Function Snap Hypothesis. This hypothesis states that “lnd(n,Δx,x) n≠1 Series snap will occur, if it occurs at all, only at an x locus transition across a complex plane axis and, except where series snap occurs, the lnd(n,Δx,x) n≠1 Series constant of integration will not change.” In this method, the constants of integration, $K_2$,$K_3$,$k_1$, and $k_2$, previously defined, have to be evaluated. This Method 2 will be described in the next section, Section 2.8.

Method 1, for the calculation of lnd(n,Δx,x) n≠1, uses Eq 2.7-4 and Eq 2.7-5 to calculate lnd(n,Δx,x) n≠1. The conditions and logic just described determine which of the two equations should be used given a particular value of x and Δx. The UBASIC language program, LNDX, for calculating the function lnd(n,Δx,x), which is placed at the end of the Appendix, implements Method 1. However, Method 1 is not the default method. Method 2 is the default method. To use Method 1 to calculate the function, lnd(n,Δx,x) n≠1, the program F2 flag must be changed from a value of 0 to a value of 1. This is done on program line 90.

## Section 2.8: Method 2 for the the calculation of the function, lnd(n,Δx,x) n≠1

This method for the calculation of the function, lnd(n,Δx,x) where n≠1, is relatively complex. One would initially think that the implementation of the lnd(n,Δx,x) n≠1 Series would be straight forward. It is not. Though the series appears to be rather common in appearance, it has characteristics which are most certainly not common. In fact, to obtain accurate computational results from the lnd(n,Δx,x) n≠1 Series, one must accept a hypothesis, the previously mentioned "lnd(n,Δx,x) n≠1 Series Snap Hypothesis", that is yet unproven. Also, the lnd(n,Δx,x) n≠1 Series places requirements on the magnitude of the quantity, $|\frac{x}{\Delta x}|$ . The quantity, $|\frac{x}{\Delta x}|$ , must be large and the quantity, $\frac{x}{\Delta x}$ , where Re(Δx)≠0, should have a significantly large real component. The larger the quantities, the more accurate the calculated value of the function, lnd(n,Δx,x).

Through considerable observation, it was found that the calculation of the function, $D_{\Delta x}$lnd(n,Δx,x) held no surprises. However, the discrete integral of this function, $_{\Delta x}\int D_{\Delta x}\text{lnd}(n,\Delta x,x)$ , did.

$_{\Delta x}\int D_{\Delta x}\text{lnd}(n,\Delta x,x) = \text{lnd}(n,\Delta x,x) + K$ where K is the constant of integration. Knowing Calculus, one would expect that the constant, K, would have a single value for all x. This turns out not to be so. For all x, the constant of integration may have as many as six values to keep the value of lnd(n,Δx,x) consistently accurate over all x given a value of n, Δx, and an arbitrary value of x, $x_i$.
See Diagram 2.8-1 below. As a result of extensive observation it was found that the locations where the constant of integration might change, could be determined. However, whether the constant of integration would change at these locations could not be determined from the values of n,Δx, and x. What is known about the behavior of the lnd(n,Δx,x) n≠1 Series forms the basis of the lnd(n,Δx,x) n≠1 Series Snap Hypothesis. The term "snap" was selected to describe the abrupt change in the series' value should a change occur. The lnd(n,Δx,x) n≠1 Series Snap Hypothesis is rewritten below. The x locus referred to is the straight line locus of all x values plotted in the complex plane.

<u>The lnd(n,Δx,x) n≠1 Series Snap Hypothesis</u>

Lnd(n,Δx,x) n≠1 Series snap will occur, if it occurs at all, only at an x locus transition across a complex plane axis and, except where series snap occurs, the lnd(n,Δx,x) n≠1 Series constant of integration will not change.

Diagram 2.8-1: Problematic locations on an x locus line

An example of an x locus line

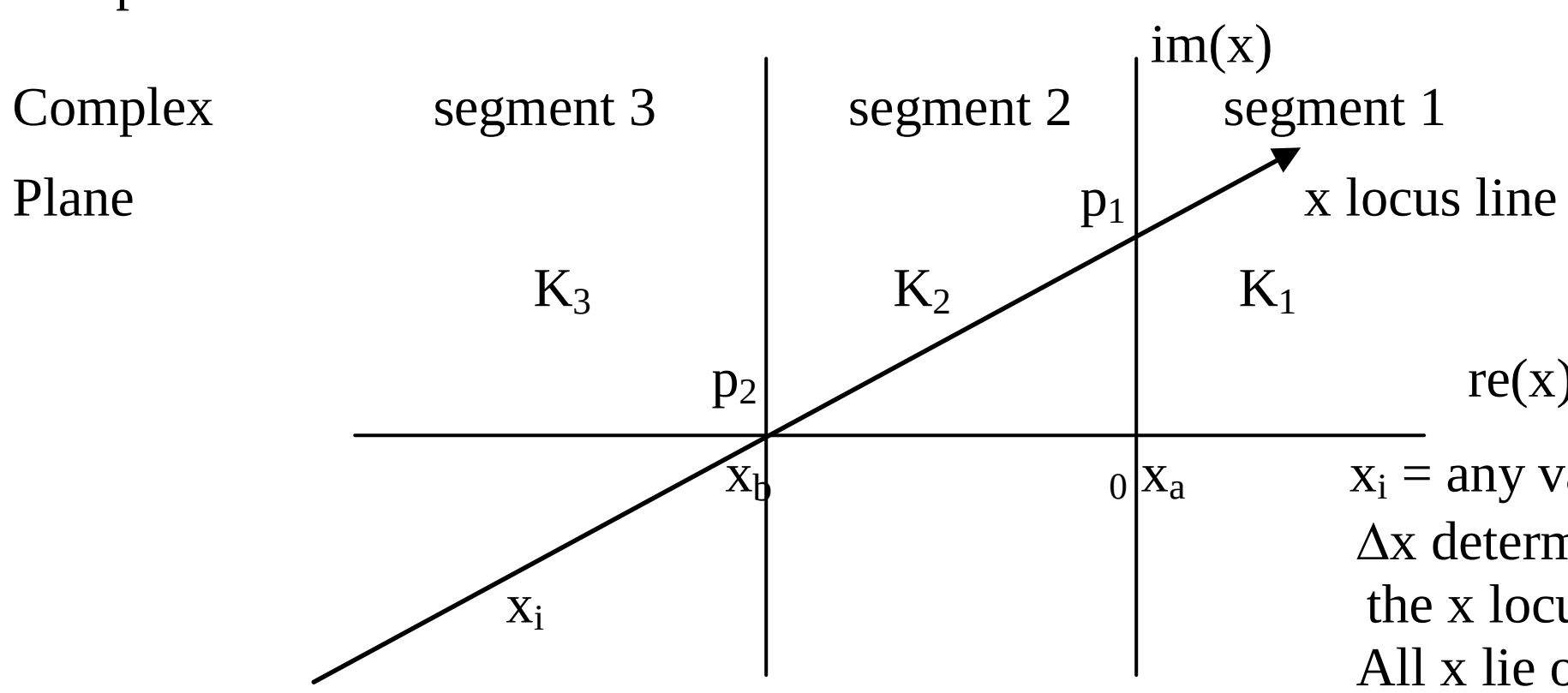

$x_i$ = any value of x
$\Delta x$ determines the slope of the x locus line
All x lie on the x locus line
$x_i$ and $\Delta x$ or $x_i$ and $-\Delta x$ determine the same x locus line
$p_1,p_2$ = complex plane axis crossover points
$Re(p_1) = x_a$
$Re(p_2) = x_b$
$K_1,K_2,K_3$ = x locus line segment constants of integration

Note – The points, $p_1$ and $p_2$, are x locus line points but are not necessarily x locus points.

The vertical lines intersecting the re(x) axis at $x_a$ and $x_b$ intersect the x locus line at a transition point, $p_1$ or $p_2$, where the x values on either side have a real or complex component of opposite sign. One of these lines is collinear with the im(x) axis. The other line passes through the x locus line re(x) axis crossover point.

A discussion concerning the values of the lnd(n,Δx,x) n≠1 Series constant of integration has been given in Sections 2.5, 2.6, and 2.7.

Potential lnd(n,Δx,x) n≠1 Series calculation problems

The lnd(n,Δx,x) n≠1 Series should not be evaluated for an x value at or near the transition points, $p_1$ or $p_2$. Calculations using x values where Re(x) are nearly equal to $x_a$ or $x_b$ are not reliably accurate.

1) There could be an lnd(n,Δx,x) n≠1 Series snap value change at the x locus line points $x=p_1$ or $x=p_2$ where $Re(x) = x_a$ or $x_b$. This value change would most likely be to an incorrect value.

2) In some cases the quantities $\left|\frac{x}{\Delta x}\right|$ and the real component of $\frac{x}{\Delta x}$ are not sufficiently large.

   In this region, the real or imaginary part of the quantity, $\frac{x}{\Delta x}$, is small.

To properly use the lnd(n,Δx,x) n≠1 Series to calculate the lnd(n,Δx,x) function where n≠1, it is necessary to address all of the important concerns mentioned above. This can and has been done in the lnd(n,Δx,x) calculation program, LNDX, presented at the end of the Appendix. However, it has been done at the expense of a considerable amount of additional coding. Though the Method 2 program code for calculating the lnd(n,Δx,x) n≠1 Series function is more extensive then that of Method 1 for calculating this same function (described in Section 2.7), its computation time can be much shorter.

Method 2 for calculating the lnd(n,Δx,x) function where n≠1, in addition to using the lnd(n,Δx,x) n≠1 Series, uses a considerable amount of support software. The support software is used primarily for two reasons. The first reason is to evaluate the proper constant of integration value for the given n, Δx, and x values. The second reason is to avoid the use of x values in the series which may degrade its accuracy or cause it to fail. A description of Method 2 follows.

Method 2 Description

Method 2 for the calculation of lnd(n,Δx,x), where n≠1, requires that the function,

$$\mathrm{lnd}(n,\Delta x,x_i) = \Delta x \sum_{\substack{\Delta x \\ x=x_i}}^{\pm\infty} \frac{1}{x^n}$$

, be plotted in the complex plane and labeled in a specific way. The various steps to properly represent this function in the complex plane are listed below. From the resulting plot diagram, a computer program to calculate the function, lnd(n,Δx,x) n≠1 for all Δx, x, and n with the exception of n=1 can be written.

lnd(n,Δx,x) n≠1 function complex plane plot description for Re(Δx) ≠0

1) Plot in the complex plane all x values of the function, $x = x_i + m\Delta x$ where m=integers. This is an x locus. These points will all lie on a straight line, the x locus line.

2) Draw a line through all of the x points to represent the x locus straight line.

3) For the desired value of Δx, indicate the direction of the summation along the x locus with an arrowhead placed on the x locus line pointing in the direction of the summation. For Re(Δx)>0, the arrowhead points from left to right. For Re(Δx)<0, the arrowhead points from right to left.

4) Identify the one or two x locus line complex plane axis crossover points. There will be one crossover point if the x locus line passes through the complex plane axis origin.

5) Draw a vertical line through each crossover point dividing the x locus line into two or three segments. One of these vertical lines may be collinear with the imaginary axis. Label the vertical line intersections with the re(x) axis from right to left $x_a$ then $x_b$ for Re(Δx)>0 or $x_b$ then $x_a$ for Re(x)<0. If there is only one re(x) axis intersection, it is labeled, $x_a$.

7) The rightmost x locus line segment is segment 1 and the leftmost line segment is segment 3. If there are three x locus line segments, the middle segment is segment 2.

8) If the arrowhead points to the right, the x locus line segment constants of integration are designated $K_1$, $K_2$, and $K_3$. The subscript number indicates on which x locus line segment the constant of integration applies. In the case where there are only two segments, there will not be a constant of integration, $K_2$. $K_1$ has been proven to always have a zero value.

9) If the arrowhead points to the left, the x locus line segment constants of integration are designated $k_1$,$k_2$, and $k_3$. The subscript number indicates on which x locus line segment the constant of integration applies. In the case where there are only two segments, there will not be a constant of integration, $k_2$. $k_1$ has been proven to always have a zero value.

10) Draw two parallel lines equidistant on either side of the vertical lines $x_a$ and $x_b$ (if it exists) described in 5). Label the four vertical line intersections with the re(x) axis from right to left $x_1$,$x_2$,$x_3$,$x_4$ consecutively for Re(Δx)<0 or $x_4$,$x_3$,$x_2$,$x_1$ consecutively for Re(Δx)<0. Where there is no $x_b$ there is no $x_3$ and $x_4$.

See Diagram 2.8-2 for some examples of x locus lines with proper labeling. The evaluation of $x_a$,$x_b$ and $x_1$ thru $x_4$ will be described on the following pages.

<u>lnd(n,Δx,x) n≠1 function complex plane plot description for Re(Δx) = 0</u>

1) Plot in the complex plane all x values of the function, $x = x_i + m\Delta x$ where m=integers. This is an x locus. These points will all lie on a vertical straight line, the x locus line.

2) Draw a line through all of the x points to represent the x locus vertical straight line.

3) For the desired value of Δx, indicate the direction of the summation along the x locus with an arrowhead placed on the x locus line pointing in the direction of the summation. For Im(Δx)>0, the arrowhead points from down to up. For Im(Δx)<0, the arrowhead points from up to down.

4) The x locus vertical straight line is perpendicular to the re(x) axis which divides it into two segments. The uppermost x locus line segment is segment 1 and the bottommost line segment is segment 3.

5) If the arrowhead points upward, the x locus line segment constants of integration are designated $K_1$ and $K_3$ (there is no $K_2$). The subscript number indicates on which x locus line segment the constant of integration applies. $K_1$ has been proven to always have a zero value.

6) If the arrowhead points downward, the x locus line segment constants of integration are designated $k_1$ and $k_3$ (there is no $k_2$). The subscript number indicates on which x locus line segment the constant of integration applies. $k_3$ has been proven to always have a zero value.

7) Equidistant from the re(x) axis draw two lines parallel to the re(x) axis. These two lines intersect the im(x) axis.

8) Label the complex plane axis origin, $x_a$. Label the parallel line intersections with the im(x) axis from up to down, $x_1$ then $x_2$ for Im(Δx)>0 or $x_2$ then $x_1$ for Im(Δx)<0.

See Diagram 2.8-2 for some examples of x locus lines with proper labeling. The evaluation of $x_1$ and $x_2$ will be described on the following pages. $x_a = 0$.

Diagram 2.8-2: Some examples of x locus lines with proper labeling

The x locus problem areas needing special attention are highlighted below. The mathematical methods to calculate all specified values of x and the constants of integration, $K_2$,$K_3$,$k_1$, and $k_2$, ($K_1$=$k_3$=0) will later be shown. The real values of the axis x locus line crossover points are $x_a$ and $x_b$.

Fig 1

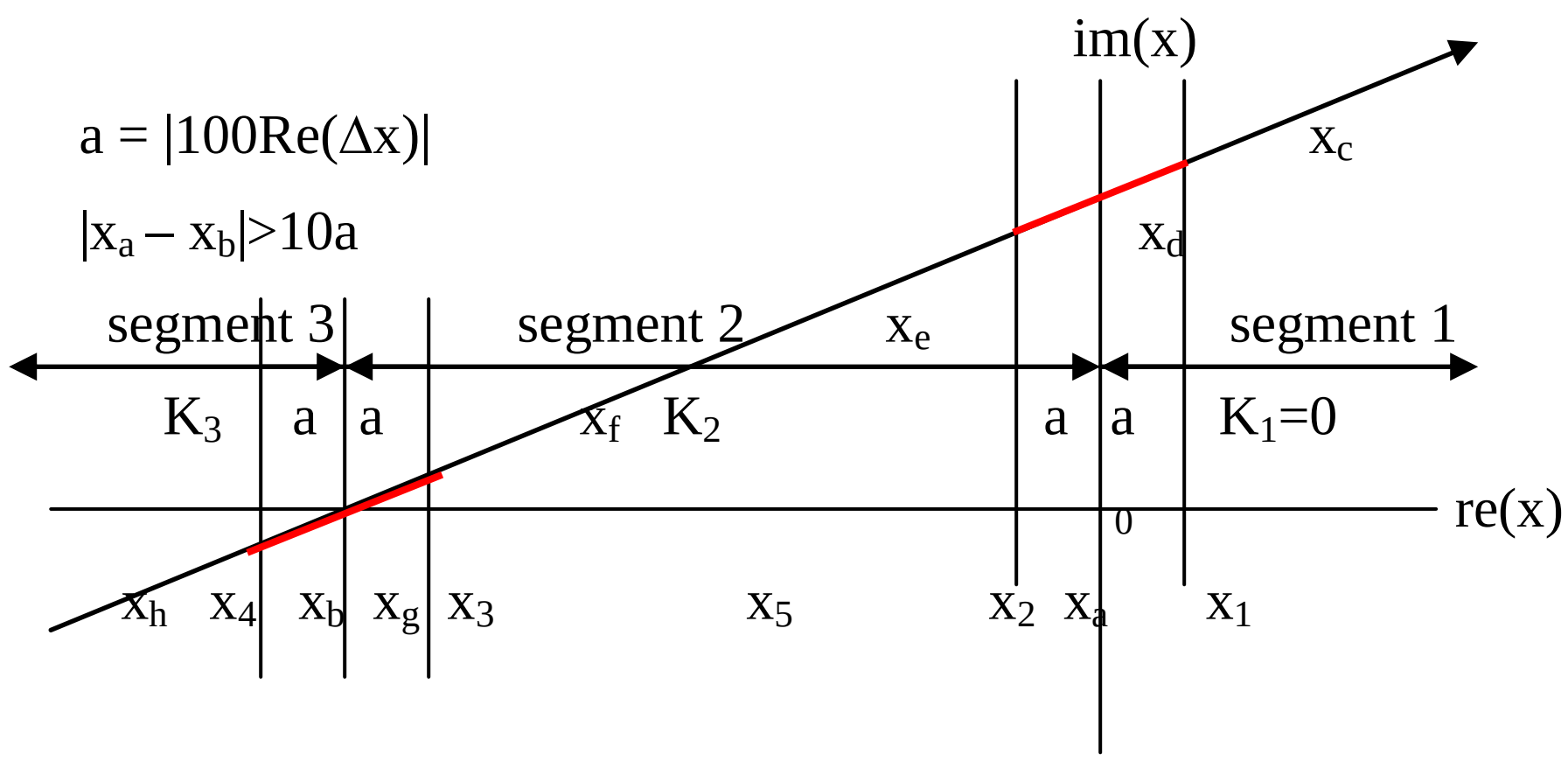

$x_a$ and $x_b$ are not close

$x_5 = x_a - 14a$

$Re(\Delta x) > 0$

Fig 2

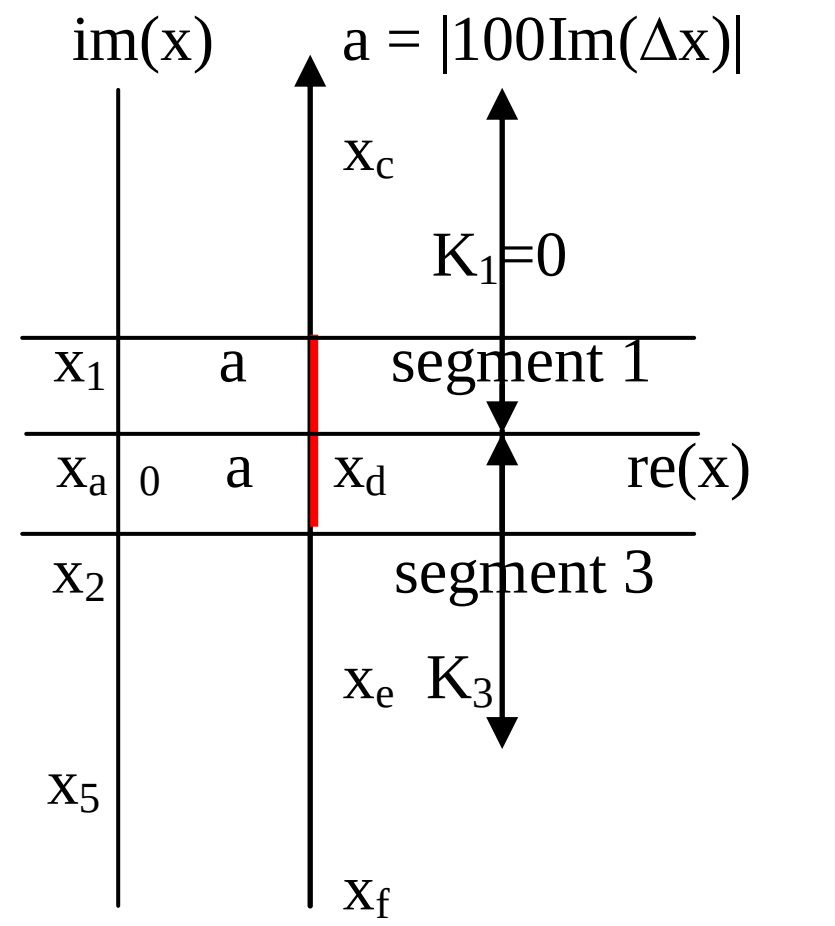

$x_5 = x_a - 14aj$ , $x_a = 0$

$Re(\Delta x) = 0$ and $Im(\Delta x)>0$

= points in various positions on the x locus line or on a complex plane axis

Fig 3

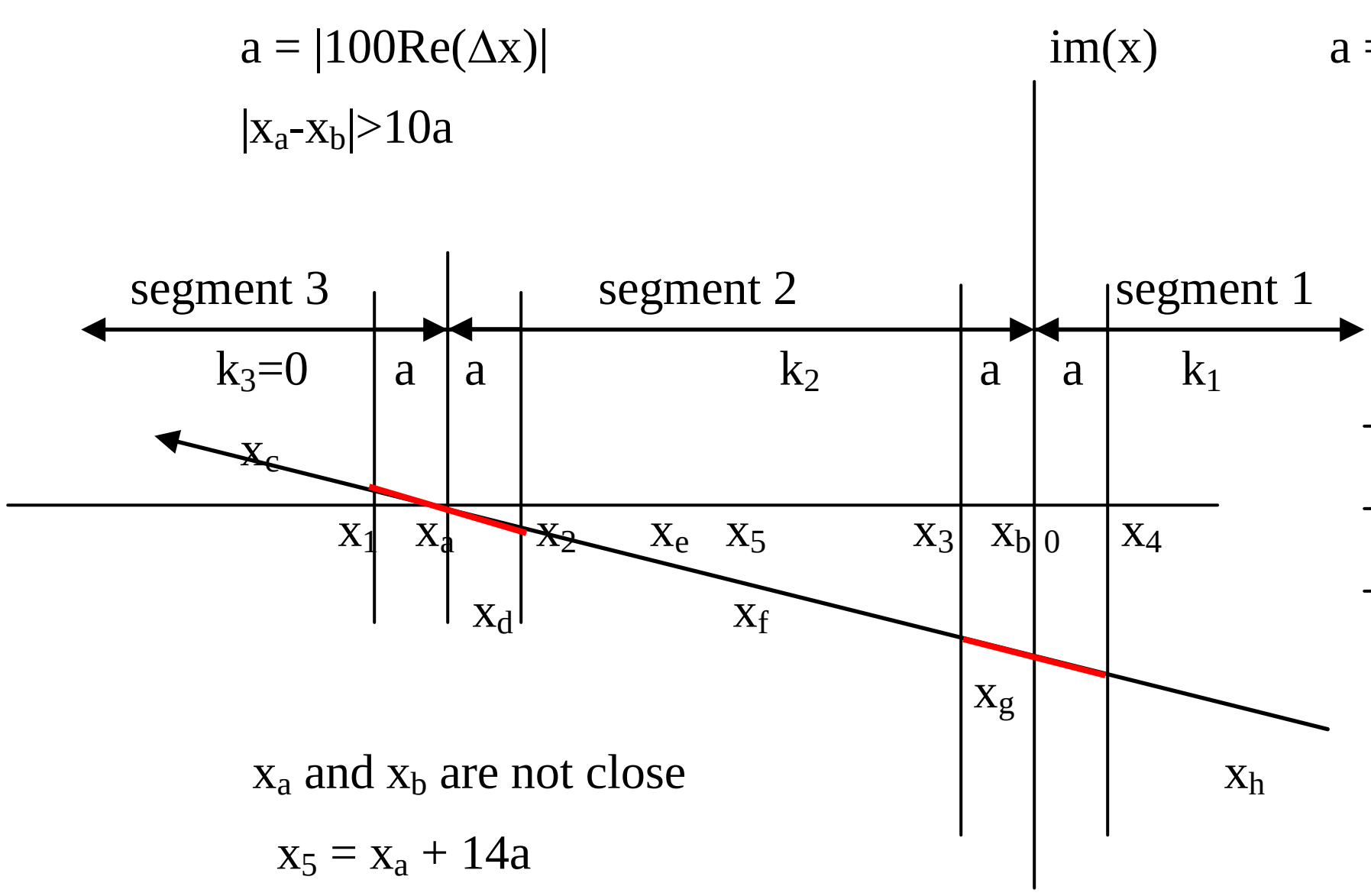

$x_a$ and $x_b$ are not close

$x_5 = x_a + 14a$

$Re(\Delta x) < 0$

Fig 4

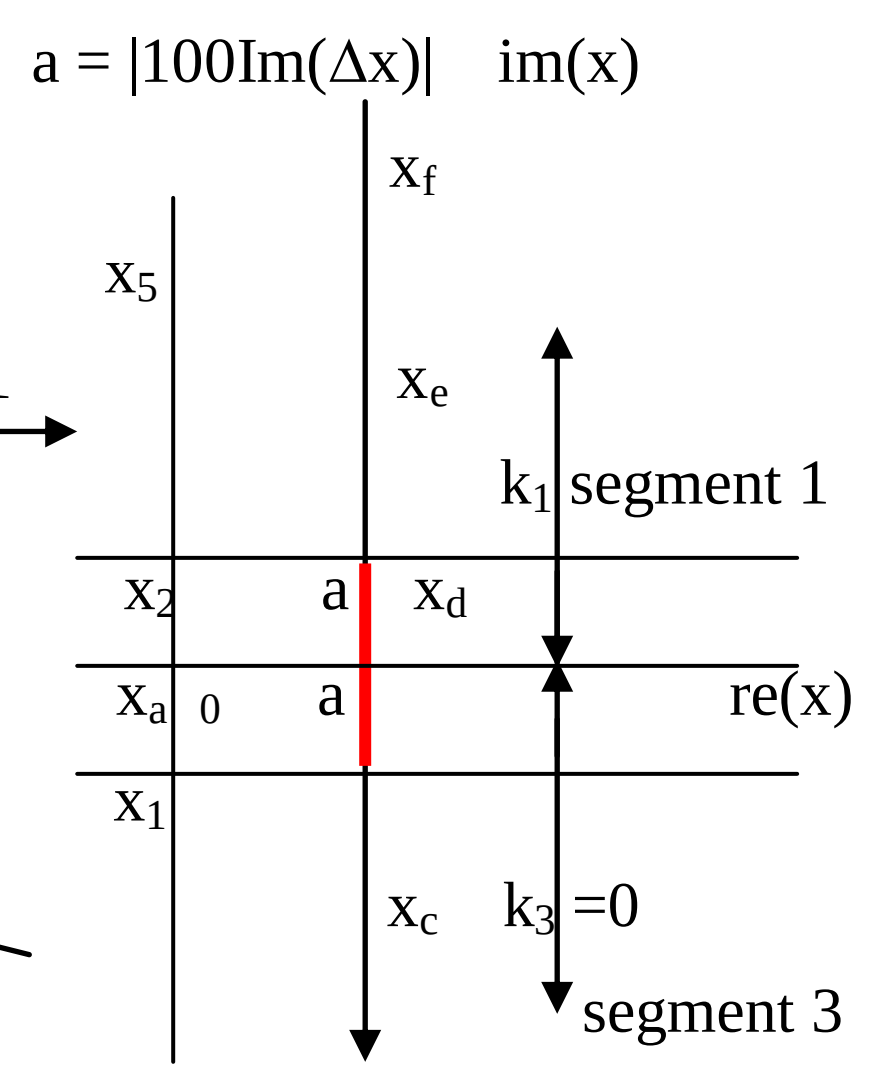

$x_5 = x_a + 14aj$ , $x_a = 0$

$Re(\Delta x) = 0$ and $Im(\Delta x)<0$

Some Special Cases

a = |100Re(Δx)| Fig 5

$x_5 = x_a + 14a$

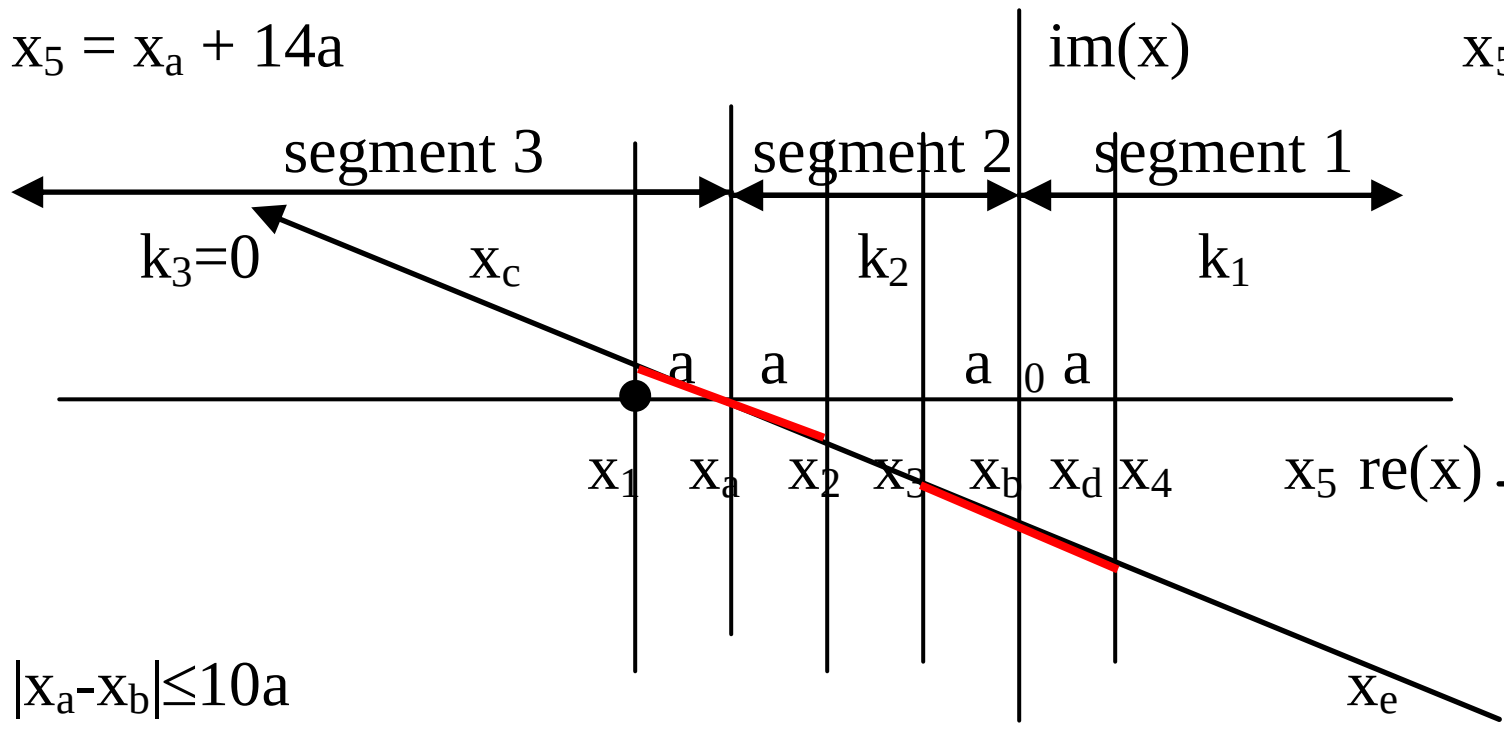

$|x_a - x_b| \le 10a$

$x_a$ and $x_b$ are close

Re(Δx) < 0

a = |100Re(Δx)| Fig 6

$x_5 = x_a - 14a$

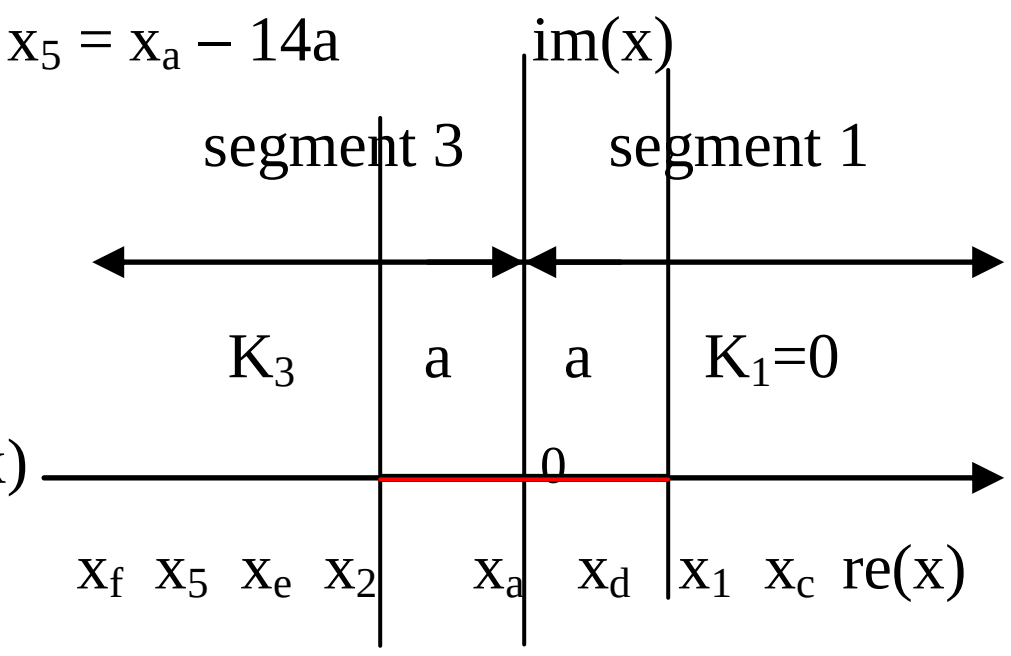

no segment 2

Re(Δx) > 0

Im(Δx) = 0

Based on the examples of x locus diagrams above, these are the methods, equations, and quantities necessary to calculate the function, lnd(n,Δx,x) n≠1, for all values of n,Δx, and x using the lnd(n,Δx,x) n≠1 Series.

1) Finding the x locus line complex plane x axis crossover point, $x_c$

xc = the calculated crossover point xc is on the x locus line but may or may not be a value of the

set of all x (i.e. $x = x_i + m\Delta x$ where $x_i$ = a value of x and m = integers)

For Re(Δx) ≠ 0

$$xc = \mathrm{Re}(x) - \mathrm{Im}(x) * \left[\frac{\mathrm{Re}(\Delta x)}{\mathrm{Im}(\Delta x)}\right] \qquad (2.8\text{-}1)$$

where

x = any x locus value

Δx = x interval

For Re(Δx) = 0

$$xc = \mathrm{Re}(x) + j0 \qquad (2.8\text{-}2)$$

where

x = any x locus value

Δx = x interval

2) Defining the points, $x_a,x_b$

For {Re(Δx)>0 and xc>0}or {Re(Δx)<0 and xc<0}

$x_a$ = xc = x locus line re(x) axis crossover point (2.8-3)

$x_b$ = 0 = real part of the x locus line im(x) axis crossover point (2.8-4)

For {Re(Δx)>0 and xc<0} or {Re(Δx)<0 and xc>0}or

{Re(Δx)=0 and xc>0}or {Re(Δx)=0 and xc<0} (2.8-5)

$x_a$ = 0 = real part of the x locus line im(x) axis crossover point

$x_b$ = xc = x locus line re(x) axis crossover point (2.8-6)

For xc=0

$x_a = x_b$ = xc = 0 = complex plane axis origin (2.8-7)

3) Defining the computation safety margin value, a

a = |100Re(Δx)| for use when Re(Δx)≠0 (2.8-8)

a = |100Im(Δx)| for use when Re(Δx)=0 (2.8-9)

Note - This length value has been found to work well in avoiding lnd(n,Δx,x) n≠1 Series low accuracy regions.

4) Find the values, $x_1,x_2,x_3, x_4$

For Re(Δx)>0

$x_1 = x_a + a$ (2.8-10)

$x_2 = x_a - a$ (2.8-11)

$x_3 = x_b + a$ (2.8-12)

$x_4 = x_b - a$ (2.8-13)

For Re(Δx)<0

$x_1 = x_a - a$ (2.8-14)

$x_2 = x_a + a$ (2.8-15)

$x_3 = x_b - a$ (2.8-16)

$x_4 = x_b + a$ (2.8-17)

Note - Some x loci do not have an $x_3$ to $x_4$ region

For Re(Δx)=0 and Im(Δx)>0

$x_1 = x_a + ja\,,\quad x_a = 0$ (2.8-18)

$x_2 = x_a - ja\,,\quad x_a = 0$ (2.8-19)

For Re(Δx)=0 and Im(Δx)<0

$x_1 = x_a - ja\,,\quad x_a = 0$ (2.8-20)

$x_2 = x_a + ja\,,\quad x_a = 0$ (2.8-21)

5) Calculating the lnd(n,Δx,$x_i$) n≠1 function when $x_i$ is located in the x locus line region between $x_1$ and $x_2$ or $x_3$ and $x_4$

$$\mathrm{lnd}(n,\Delta x,x_i) = \Delta x \sum_{x=x_i}^{\mathit{x}_m-\Delta x} \frac{1}{x^n} + \mathrm{lnd}(n,\Delta x,\mathit{x}_m) \qquad (2.8\text{-}22)$$

where

$x_i$ = an x locus value on the x locus line

m = 1 for $x_i$ between $x_1$ and $x_2$

m = 3 for $x_i$ between $x_3$ and $x_4$

$\mathit{x}_m$ = the x locus value where Re($\mathit{x}_m$) = $x_m$

Δx = x interval

The x locus line regions between $x_1$ and $x_2$ and between $x_3$ and $x_4$ are potentially problematic for obtaining an accuarate calculation using the lnd(n,Δx,x) n≠ Series.

The summation part of the equation skips over the potentially problematic x locus line region. The values $x_1$ or $x_3$, at the end of this region are satifactory for calculating a highly accurate value of lnd(n,Δx,x) n≠1 using the series.

Note - It is necessary to find an x locus value, $\mathit{x}_m$, where Re($\mathit{x}_m$) ≈ $x_m$. The mathematical method to do this is shown below.

6) Find the x locus value, $\mathit{x}_m$, where Re($\mathit{x}_m$) ≈ $x_m$.

For Re(Δx)≠0

$$\mathit{x}_m = x_i + \mathrm{int}\left[\frac{x_m - \mathrm{re}(x_i)}{\mathrm{re}(\Delta x)}\right]\Delta x\,,\qquad m = 1,2,3\ldots \qquad (2.8\text{-}23)$$

where

$x_i$ = an x locus value on the x locus line

int(p) = integer value of p, a real number

For $Re(\Delta x)=0$

$$\mathit{x}_m = x_i + \text{int}\left[\frac{x_m - \text{im}(x_i)}{\text{im}(\Delta x)}\right]\Delta x\,, \qquad m = 1,2,3\ldots \tag{2.8-24}$$

where

$x_i$ = an x locus value on the x locus line

int(p) = integer value of p, a real number

7) <u>Calculating the constants of integration, $K_2$, $K_3$, $k_1$, and $k_2$</u>

$$c_1 = \Delta x \sum_{x=\mathit{x}_2}^{\mathit{x}_1-\Delta x} \frac{1}{x^n} + \text{lnd}_f(n,\Delta x,\mathit{x}_1) - \text{lnd}_f(n,\Delta x,\mathit{x}_2) \tag{2.8-25}$$

$$c_2 = \Delta x \sum_{x=\mathit{x}_4}^{\mathit{x}_3-\Delta x} \frac{1}{x^n} + \text{lnd}_f(n,\Delta x,\mathit{x}_3) - \text{lnd}_f(n,\Delta x,\mathit{x}_4) \tag{2.8-26}$$

where

$\mathit{x}_m$ = the x locus value where $Re(\mathit{x}_m) \approx x_m$

$\Delta x$ = x interval

$$\text{lnd}_f(n,\Delta x,x) = -\sum_{m=0}^{\infty} \frac{\Gamma(n+2m-1)\left(\frac{\Delta x}{2}\right)^{2m} C_m}{\Gamma(n)(2m+1)!\left(x - \frac{\Delta x}{2}\right)^{n+2m-1}}\,, \quad n\neq 1 \tag{2.8-27}$$

For good accuracy the magnitude, $\left|\frac{x}{\Delta x}\right|$, must be large.

For finding the constants, $c_1$ and $c_2$, using Eq 2.8-25 and Eq 2.8-26, the magnitude, $\left|\frac{x}{\Delta x}\right|$, is large.

<u>Note</u> - The relative position of the $x_1$, $x_2$, $x_3$, and $x_4$ points change as a function of $\Delta x$ characteristics. See the previously presented lnd(n,Δx,x) n≠1 function complex plane plot descriptions.

Table 2.8-1 Calculation of the constants of integration

| No. | Δx Characteristics | Constants of Integration | Sample x locus in Diagram 2.8-2 |
|---|---|---|---|
| 1 | All Δx | $K_1 = 0$<br>$k_3 = 0$ | See Fig 1 thru Fig 6 |
| 2 | Re(Δx)>0 | $K_2 = c_1$<br>$K_3 = c_1 + c_2$ | Fig 1 |
| 3 | Re(Δx)<0 | $k_2 = c_1$<br>$k_1 = c_1 + c_2$ | Fig 3 |
| 4 | Re(Δx)=0 and Im(Δx)>0 | $K_3 = c_1$ | Fig 2 |
| 5 | Re(Δx)=0 and Im(Δx)<0 | $k_1 = c_1$ | Fig 4 |
| 6 | Re(Δx)>0 and Im(Δx)=0 | $K_3 = c_1$ | Fig 6 |
| 7 | Re(Δx)<0 and Im(Δx)=0 | $k_1 = c_1$ | Related to Fig 6<br>Opposite summation direction |
| 8 | Re(Δx)<0 and<br>$\lvert Re(x_a)–Re(x_b)\rvert \le 10a$ | $k_1 = \text{lnd}(n,\Delta x,x_5) – \text{lnd}_f(n,\Delta x,x_5)$ | Fig 5 |
| 9 | Re(Δx)>0 and<br>$\lvert Re(x_a)–Re(x_b)\rvert \le 10a$ | $K_3 = \text{lnd}(n,\Delta x,x_5) – \text{lnd}_f(n,\Delta x,x_5)$ | Related to Fig 5<br>Opposite summation direction |

The derivation of some of the constants of integration shown in the table above is shown below. The equations presented are based on the lnd(n,Δx,x) n≠1 Series Snap Hypothesis.

For Re(Δx)>0 (row 2)

Refer to Diagram 2.8-2 Fig 1

$$K_2 + \text{lnd}_f(n,\Delta x,x_2) = \Delta x \sum_{x=x_2}^{x_1-\Delta x} \frac{1}{x^n} + K_1 + \text{lnd}_f(n,\Delta x,x_1) \tag{2.8-28}$$

$$K_1 = 0 \tag{2.8-29}$$

$$K_2 = \Delta x \sum_{x=x_2}^{x_1-\Delta x} \frac{1}{x^n} + \text{lnd}_f(n,\Delta x,x_1) - \text{lnd}_f(n,\Delta x,x_2) \tag{2.8-30}$$

$$K_3 + \ln d_f(n,\Delta x,x_4) = \Delta x \sum_{x=x_4}^{x_3-\Delta x} \frac{1}{x^n} + K_2 + \ln d_f(n,\Delta x,x_3) \quad (2.8\text{-}31)$$

$$K_3 = \Delta x \sum_{x=x_4}^{x_3-\Delta x} \frac{1}{x^n} + \ln d_f(n,\Delta x,x_3) - \ln d_f(n,\Delta x,x_4) + K_2 \quad (2.8\text{-}32)$$

Substituting

$$K_3 = [\Delta x \sum_{x=x_4}^{x_3-\Delta x} \frac{1}{x^n} + \ln d_f(n,\Delta x,x_3) - \ln d_f(n,\Delta x,x_4)] + [\Delta x \sum_{x=x_2}^{x_1-\Delta x} \frac{1}{x^n} + \ln d_f(n,\Delta x,x_1) - \ln d_f(n,\Delta x,x_2)] = c_1 + c_2$$

(2.8-33)

Then

$$c_1 = \Delta x \sum_{x=x_2}^{x_1-\Delta x} \frac{1}{x^n} + \ln d_f(n,\Delta x,x_1) - \ln d_f(n,\Delta x,x_2) \quad (2.8\text{-}34)$$

$$c_2 = \Delta x \sum_{x=x_4}^{x_3-\Delta x} \frac{1}{x^n} + \ln d_f(n,\Delta x,x_3) - \ln d_f(n,\Delta x,x_4)] \quad (2.8\text{-}35)$$

where

$x_m$ = the x locus value where $Re(x_m) \approx x_m$

$\Delta x$ = x interval

$$K_2 = c_1 \quad (2.8\text{-}36)$$

$$K_3 = c_1 + c_2 \quad (2.8\text{-}37)$$

For re(Δx)<0 (row 3)

Refer to Diagram 2.8-2 Fig 3

$$k_2 + \ln d_f(n,\Delta x,x_2) = \Delta x \sum_{x=x_2}^{x_1-\Delta x} \frac{1}{x^n} + k_3 + \ln d_f(n,\Delta x,x_1) \quad (2.8\text{-}38)$$

$$k_3 = 0 \quad (2.8\text{-}39)$$

$$k_2 = \Delta x \sum_{x=x_2}^{x_1-\Delta x} \frac{1}{x^n} + \ln d_f(n,\Delta x,x_1) - \ln d_f(n,\Delta x,x_2) \qquad (2.8\text{-}40)$$

$$k_1 + \ln d_f(n,\Delta x,x_4) = \Delta x \sum_{x=x_4}^{x_3-\Delta x} \frac{1}{x^n} + k_2 + \ln d_f(n,\Delta x,x_3) \qquad (2.8\text{-}41)$$

$$k_1 = \Delta x \sum_{x=x_4}^{x_3-\Delta x} \frac{1}{x^n} + \ln d_f(n,\Delta x,x_3) - \ln d_f(n,\Delta x,x_4) + k_2 \qquad (2.8\text{-}42)$$

Substituting

$$k_1 = \left[\Delta x \sum_{x=x_4}^{x_3-\Delta x} \frac{1}{x^n} + \ln d_f(n,\Delta x,x_3) - \ln d_f(n,\Delta x,x_4)\right] + \left[\Delta x \sum_{x=x_2}^{x_1-\Delta x} \frac{1}{x^n} + \ln d_f(n,\Delta x,x_1) - \ln d_f(n,\Delta x,x_2)\right] = c_2 + c_1$$

Then (2.8-43)

$$c_1 = \Delta x \sum_{x=x_2}^{x_1-\Delta x} \frac{1}{x^n} + \ln d_f(n,\Delta x,x_1) - \ln d_f(n,\Delta x,x_2) \qquad (2.8\text{-}44)$$

$$c_2 = \Delta x \sum_{x=x_4}^{x_3-\Delta x} \frac{1}{x^n} + \ln d_f(n,\Delta x,x_3) - \ln d_f(n,\Delta x,x_4) \qquad (2.8\text{-}45)$$

where

$x_m$ = the x locus value where $\mathrm{Re}(x_m) \approx x_m$

$\Delta x$ = x interval

$$k_2 = c_1 \qquad (2.8\text{-}46)$$

$$k_1 = c_1 + c_2 \qquad (2.8\text{-}47)$$

For re(Δx)=0 and Im(Δx)>0   (row 4)

Refer to Diagram 2.8-2 Fig 2

$$K_3 + \ln d_f(n,\Delta x,x_2) = \Delta x \sum_{x=x_2}^{x_1-\Delta x} \frac{1}{x^n} + K_1 + \ln d_f(n,\Delta x,x_1) \tag{2.8-48}$$

$$K_1 = 0 \tag{2.8-49}$$

$$K_3 = \Delta x \sum_{x=x_2}^{x_1-\Delta x} \frac{1}{x^n} + \ln d_f(n,\Delta x,x_1) - \ln d_f(n,\Delta x,x_2) = c_1 \tag{2.8-50}$$

Then

$$c_1 = \Delta x \sum_{x=x_2}^{x_1-\Delta x} \frac{1}{x^n} + \ln d_f(n,\Delta x,x_1) - \ln d_f(n,\Delta x,x_2) \tag{2.8-51}$$

where

$x_m$ = the x locus value where $\mathrm{Re}(x_m) \approx \mathrm{x}_m$

Δx = x interval

$$K_3 = c_1 \tag{2.8-52}$$

For re(Δx)=0 and Im(Δx)<0   (row 5)

Refer to Diagram 2.8-2 Fig 4

$$k_1 + \ln d_f(n,\Delta x,x_2) = \Delta x \sum_{x=x_2}^{x_1-\Delta x} \frac{1}{x^n} + \ln d_f(n,\Delta x,x_1) + k_3 \tag{2.8-53}$$

$$k_3 = 0 \tag{2.8-54}$$

$$k_1 = \Delta x \sum_{x=x_2}^{x_1-\Delta x} \frac{1}{x^n} + \ln d_f(n,\Delta x,x_1) - \ln d_f(n,\Delta x,x_2) = c_1 \tag{2.8-55}$$

Then

$$c_1 = \Delta x \sum_{x=x_2}^{x_1-\Delta x} \frac{1}{x^n} + \ln d_f(n,\Delta x,x_1) - \ln d_f(n,\Delta x,x_2) \qquad (2.8\text{-}56)$$

where

$x_m$ = the x locus value where $\mathrm{Re}(x_m) \approx x_m$

$\Delta x$ = x interval

$$k_1 = c_1 \qquad (2.8\text{-}57)$$

For re(Δx)>0 and Im(Δx)=0 (row 6)

Refer to Diagram 2.8-2 Fig 6

$$K_3 + \ln d_f(n,\Delta x,x_2) = \Delta x \sum_{x=x_2}^{x_1-\Delta x} \frac{1}{x^n} + \ln d_f(n,\Delta x,x_1) + K_1 \qquad (2.8\text{-}58)$$

$$K_1 = 0 \qquad (2.8\text{-}59)$$

$$K_3 = \Delta x \sum_{x=x_2}^{x_1-\Delta x} \frac{1}{x^n} + \ln d_f(n,\Delta x,x_1) - \ln d_f(n,\Delta x,x_2) = c_1 \qquad (2.8\text{-}60)$$

Then

$$c_1 = \Delta x \sum_{x=x_2}^{x_1-\Delta x} \frac{1}{x^n} + \ln d_f(n,\Delta x,x_1) - \ln d_f(n,\Delta x,x_2) \qquad (2.8\text{-}61)$$

where

$x_m$ = the x locus value where $\mathrm{Re}(x_m) \approx x_m$

$\Delta x$ = x interval

$$K_3 = c_1 \qquad (2.8\text{-}62)$$

For re(Δx)<0 and |Re($x_a$)–Re($x_b$)|≤10a   (row 8)

Refer to Diagram 2.8-5 Fig 5

$$\text{lnd}(n,\Delta x,x_5) = \text{lnd}_f(n,\Delta x,x_5) + k_1 \qquad (2.8\text{-}63)$$

Then

$$k_1 = \text{lnd}(n,\Delta x,x_5) - \text{lnd}_f(n,\Delta x,x_5) \qquad (2.8\text{-}64)$$

After the necessary constants of integration have been calculated, the lnd(n,Δx,x) n≠1 function can be calculated for the specified values of n, Δx, and x. To clarify the the way in which this is done, the function, lnd(n,Δx,x) where n≠1, will be calculated for various points located on the x locus lines plotted in Fig1 thru Fig 6 of  Diagram 2.8-2. Note the table below. The Program Condition Number referred to in the table identifies which of the five calculation procedures (which are based on the presence of certain conditions) is used in the lnd(n,Δx,x) calculation program, LNDX, presented at the end of the Appendix. The program will specify which procedure is used in the lnd(n,Δx,x) n≠1 calculation by designating it as Cn = 1, 2, 3, 4, or 5. The procedures differ depending on the Δx condition and where the x value of the lnd(n,Δx,x) function is positioned on the complex plane x locus line.

Table 2.8-2  lnd(n,Δx,x) n≠1 function calculations from the x locus diagrams of Diagram 2.8-2

| **Fig no.** | **Δx Condition** | **$x_i$ value** | $\int_{\Delta x\, x_i}^{\pm\infty} \frac{1}{x^n}\Delta x = \Delta x \sum_{\Delta x\, x=x_i}^{\pm\infty} \frac{1}{x^n}$ = lnd(n,Δx,$x_i$) | **Program Condition no.** | **Comments** |
|---|---|---|---|---|---|
| 1 | re(Δx)>0 | $x_c$ | $\text{lnd}_f(n,\Delta x,x_c)$ | 1 | $x_c$ is on segment 1<br>Re($x_c$) ≥ $x_1$ |
| 1 | re(Δx)>0 | $x_d$ | $\Delta x \sum_{\Delta x\, x=x_d}^{x_1-\Delta x} \frac{1}{x^n} + \text{lnd}_f(n,\Delta x,x_1)$ | 1 | Must sum out of a problem region |
| 1 | re(Δx)>0 | $x_e$ | $\Delta x \sum_{\Delta x\, x=x_e}^{x_1-\Delta x} \frac{1}{x^n} + \text{lnd}_f(n,\Delta x,x_1)$ | 1 | Sum over a problem region<br>Re($x_e$) ≥ $x_5$ |
| 1 | re(Δx)>0 | $x_f$ | $\text{lnd}_f(n,\Delta x,x_f) + K_2$ | 3 | $x_f$ is on segment 2<br>Re($x_f$) < $x_5$ |

| Fig no. | Δx Condition | $x_i$ value | $\int_{x_i}^{\pm\infty}{}_{\Delta x}\frac{1}{x^n}\,\Delta x = \Delta x\sum_{\Delta x,\,x=x_i}^{\pm\infty}\frac{1}{x^n} = \mathrm{lnd}(n,\Delta x,x_i)$ | Program Condition no. | Comments |
|---|---|---|---|---|---|
| 1 | re(Δx)>0 | $x_g$ | $\Delta x\sum_{\Delta x,\,x=x_g}^{x_3-\Delta x}\frac{1}{x^n} + \mathrm{lnd}_f(n,\Delta x,x_3) + K_2$ | 4 | Must sum out of a Problem region |
| 1 | re(Δx)>0 | $x_h$ | $\mathrm{lnd}_f(n,\Delta x,x_h) + K_3$ | 5 | $x_h$ is on segment 3<br>$\mathrm{Re}(x_h) \le x_4$ |
| 2 | re(Δx)=0 and im(Δx)>0 | $x_c$ | $\mathrm{lnd}_f(n,\Delta x,x_c)$ | 1 | $x_c$ is on segment 1<br>$\mathrm{Im}(x_c) \ge \mathrm{Im}(x_1)$ |
| 2 | re(Δx)=0 and im(Δx)>0 | $x_d$ | $\Delta x\sum_{\Delta x,\,x=x_d}^{x_1-\Delta x}\frac{1}{x^n} + \mathrm{lnd}_f(n,\Delta x,x_1)$ | 1 | Must sum out of a problem region |
| 2 | re(Δx)=0 and im(Δx)>0 | $x_e$ | $\Delta x\sum_{\Delta x,\,x=x_e}^{x_1-\Delta x}\frac{1}{x^n} + \mathrm{lnd}_f(n,\Delta x,x_1)$ | 1 | Sum over a problem region<br>The summation is short<br>$\mathrm{Im}(x_e) \ge \mathrm{Im}(x_5)$ |
| 2 | re(Δx)=0 and im(Δx)>0 | $x_f$ | $\mathrm{Lnd}_f(n,\Delta x,x_f) + K_3$ | 2 | $x_f$ is on segment 3<br>$\mathrm{Im}(x_f) < \mathrm{Im}(x_5)$ |
| 3 | re(Δx)<0 | $x_c$ | $\mathrm{lnd}_f(n,\Delta x,x_c)$ | 1 | $x_c$ is on segment 3<br>$\mathrm{Re}(x_c) \le x_1$ |
| 3 | re(Δx)<0 | $x_d$ | $\Delta x\sum_{\Delta x,\,x=x_d}^{x_1-\Delta x}\frac{1}{x^n} + \mathrm{lnd}_f(n,\Delta x,x_1)$ | 1 | Must sum out of a problem region |
| 3 | re(Δx)<0 | $x_e$ | $\Delta x\sum_{\Delta x,\,x=x_e}^{x_1-\Delta x}\frac{1}{x^n} + \mathrm{lnd}_f(n,\Delta x,x_1)$ | 1 | Sum over a problem region<br>The summation is short<br>$\mathrm{Re}(x_e) \le x_5$ |

| **Fig no.** | **Δx Condition** | **$x_i$ value** | $\displaystyle {}_{\Delta x}\!\int_{x_i}^{\pm\infty} \frac{1}{x^n}\,\Delta x = \Delta x \, {}_{\Delta x}\!\sum_{x=x_i}^{\pm\infty} \frac{1}{x^n}$ $= \mathrm{lnd}(n,\Delta x,x_i)$ | **Program Condition no.** | **Comments** |
|---|---|---|---|---|---|
| 3 | re(Δx)<0 | $x_f$ | $\mathrm{lnd}_f(n,\Delta x,x_f) + k_2$ | 3 | $x_f$ is on segment 2<br>$\mathrm{Re}(x_f) > x_5$ |
| 3 | re(Δx)<0 | $x_g$ | $\displaystyle \Delta x \, {}_{\Delta x}\!\sum_{x=x_g}^{x_3-\Delta x} \frac{1}{x^n} + \mathrm{lnd}_f(n,\Delta x,x_3) + k_2$ | 4 | Must sum out of a problem region |
| 3 | re(Δx)<0 | $x_h$ | $\mathrm{lnd}_f(n,\Delta x,x_h) + k_1$ | 5 | $x_h$ is on segment 1<br>$\mathrm{Re}(x_h) \geq x_4$ |
| 4 | re(Δx)=0 and im(Δx)<0 | $x_c$ | $\mathrm{lnd}_f(n,\Delta x,x_c)$ | 1 | $x_c$ is on segment 3<br>$\mathrm{Im}(x_c) \leq \mathrm{Im}(x_1)$ |
| 4 | re(Δx)=0 and im(Δx)<0 | $x_d$ | $\displaystyle \Delta x \, {}_{\Delta x}\!\sum_{x=x_d}^{x_1-\Delta x} \frac{1}{x^n} + \mathrm{lnd}_f(n,\Delta x,x_1)$ | 1 | Must sum out of a problem region |
| 4 | re(Δx)=0 and im(Δx)<0 | $x_e$ | $\displaystyle \Delta x \, {}_{\Delta x}\!\sum_{x=x_e}^{x_1-\Delta x} \frac{1}{x^n} + \mathrm{lnd}_f(n,\Delta x,x_1)$ | 1 | Sum over a problem region<br>The summation is short<br>$\mathrm{Im}(x_e) \leq \mathrm{Im}(x_5)$ |
| 4 | re(Δx)=0 and im(Δx)<0 | $x_f$ | $\mathrm{lnd}_f(n,\Delta x,x_f) + k_1$ | 2 | $x_f$ is on segment 1<br>$\mathrm{Im}(x_f) > \mathrm{Im}(x_5)$ |
| 5 | re(Δx)<0 | $x_c$ | $\mathrm{lnd}_f(n,\Delta x,x_c)$ | 1 | $x_c$ is on segment 3<br>$\mathrm{Re}(x_c) \leq x_1$ |
| 5 | re(Δx)<0 | $x_d$ | $\displaystyle \Delta x \, {}_{\Delta x}\!\sum_{x=x_d}^{x_1-\Delta x} \frac{1}{x^n} + \mathrm{lnd}_f(n,\Delta x,x_1)$ | 1 | Must sum out of a problem region or the region extending to $x_5$<br>The summation is short<br>$\mathrm{Re}(x_d) \leq x_5$ |

| **Fig no.** | **$\Delta$x Condition** | **$x_i$ value** | $\underset{\Delta x}{\int_{x_i}^{\pm\infty}} \frac{1}{x^n}\,\Delta x = \Delta x \underset{\Delta x}{\sum_{x=x_i}^{\pm\infty}} \frac{1}{x^n}$ $= \text{lnd}(n,\Delta x,x_i)$ | **Program Condition no.** | **Comments** |
|---|---|---|---|---|---|
| 5 | re($\Delta$x)<0 | $x_f$ | $\text{lnd}_f(n,\Delta x,x_e) + k_1$ | 2 | $x_e$ is on segment 1<br>$\text{Re}(x_e) > x_5$ |
| 6 | re($\Delta$x)>0 | $x_c$ | $\text{lnd}_f(n,\Delta x,x_c)$ | 1 | $x_c$ is on segment 1<br>$x_c \geq x_1$ |
| 6 | re($\Delta$x)>0 | $x_d$ | $\Delta x \underset{\Delta x}{\sum_{x=x_d}^{x_1-\Delta x}} \frac{1}{x^n} + \text{lnd}_f(n,\Delta x,x_1)$ | 1 | Must sum out of a problem region |
| 6 | re($\Delta$x)>0 | $x_e$ | $\Delta x \underset{\Delta x}{\sum_{x=x_e}^{x_1-\Delta x}} \frac{1}{x^n} + \text{lnd}_f(n,\Delta x,x_1)$ | 1 | Sum over a problem region<br>The summation is short<br>$x_e \geq x_5$ |
| 6 | re($\Delta$x)>0 | $x_f$ | $\text{lnd}_f(n,\Delta x,x_f) + K_3$ | 2 | $x_f$ is on segment3<br>$x_f < x_5$ |

$x_m$ = the x locus value where $\text{Re}(x_m) \approx x_m$ for $\text{Re}(\Delta x)\neq 0$ or $\text{Im}(x_m) \approx \text{Im}(x_m)$ for $\text{Re}(\Delta x)=0$

The Method 2 calculation of the lnd(n,$\Delta$x,x) n$\neq$1 function is more complex than one would first think. Firstly, it is based on a hypothesis, the lnd(n,$\Delta$x,x) Function Snap Hypothesis, which is not yet proven. Secondly, the lnd(n,$\Delta$x,x) n$\neq$1 Series has problematic characteristics which must be managed. The position of x on the complex plane x locus is important and the magnitude, $\left|\frac{x}{\Delta x}\right|$, is important when determining how to calculate the lnd(n,$\Delta$x,x) n$\neq$1 function. The explanation above endeavors to clarify the methodology behind the Method 2 lnd(n,$\Delta$x,x) n$\neq$1 function calculation. As is seen, different conditions require a different computational approach. The lnd(n,$\Delta$x,x) calculation program, LNDX, in the Appendix incorporates all of these various computational approaches into one program to calculate the value of the function, lnd(n,$\Delta$x,x) n$\neq$1, given a specific real or complex value for n, $\Delta$x, and x.The program, of course, calculates the function, lnd(n,$\Delta$x,x), where n=1. However, the methodlogy to do this is different from that of Method 2 and is not nearly as complex. A short discription of this methodology is found in the following section, Section 2.9.

To see an application of the LNDX program for a case where n=1 and a case where n$\neq$1, see Example 2.1 and Example 2.2 in the Chapter 2 Solved Problems section at the end of this Chapter.

## Section 2.9: The Method for the calculation of the function, $\mathrm{lnd}(1,\Delta x,x) \equiv \mathrm{lnd}_{\Delta x}x$

The $\ln_{\Delta x}x$ Series is functionally different from the lnd(n,Δx,x) n≠1 Series. As previously pointed out, the lnd(n,Δx,x) n≠1 Series is highly irregular and considerable care must be given to the way in which it is used. Though the $\ln_{\Delta x}x$ Series even shares the same series constants with the lnd(n,Δx,x) n≠1 Series, its use is, nevertheless, straightforward. Using the $\ln_{\Delta x}x$ Series to calculate the function, $\ln_{\Delta x}x \equiv \mathrm{lnd}(1,\Delta x,x)$ is relatively easy. Note the $\ln_{\Delta x}x$ Series presented below:

The $\ln_{\Delta x}x$ Series

$$\ln_{\Delta x}x \equiv \mathrm{lnd}(1,\Delta x,x) \approx \gamma + \ln\left(\frac{x}{\Delta x} - \frac{1}{2}\right) + \sum_{m=1}^{\infty} \frac{(2m-1)!\, C_m}{(2m+1)!\, 2^{2m}\left(\frac{x}{\Delta x} - \frac{1}{2}\right)^{2m}} \tag{2.9-1}$$

where

accuracy improves rapidly for increasing $\left|\frac{x}{\Delta x}\right|$

n,Δx,x,*x* are real or complex values

$\gamma$ = Euler's Constant, .577215664

$x = x_i + n\Delta x$ ,  n = integers

$x_i$ = a value of x

$\Delta x$ = x increment

$C_m$ = series coefficients

The only complication in calculating the $\ln_{\Delta x}x$ Series is the handling of the accuaracy requirement that the value of the quantity, $\left|\frac{x}{\Delta x}\right|$ , must be large. For some selected combinations of x and Δx, the quantity, $\left|\frac{x}{\Delta x}\right|$ , may not be large. To handle instances such as this and to assure high calculation accuracy for all x and Δx, a method for calculating the function, $\ln_{\Delta x}x$, will be derived below.

Derive a method to accuarately calculate the function, $\ln_{\Delta x}x$.

$$\ln_{\Delta x}x = \ln_1\left(\frac{x}{\Delta x}\right) \tag{2.9-2}$$

$$R = \frac{x}{\Delta x} \tag{2.9-3}$$

$$\int_{R}^{R+N} {}_1\frac{1}{r}\,\Delta r = \ln_1 r \,\Big|_{R}^{R+N} = \ln_1(R+N) - \ln_1 R = (1)\ {}_1\!\sum_{r=R}^{R+N-1} \frac{1}{r} \qquad (2.9\text{-}4)$$

where

$N = 1,2,3,\ldots$

From Eq 2.9-4

$$\ln_1(R+N) - \ln_1 R = {}_1\!\sum_{r=R}^{R+N-1} \frac{1}{r}$$

Rearranging terms

$$\ln_1 R = \ln_1(R+N) - {}_1\!\sum_{r=R}^{R+N-1} \frac{1}{r} \qquad (2.9\text{-}5)$$

or

$$\ln_1 R = \ln_1(R+N) - {}_{-1}\!\sum_{r=R+N-1}^{R} \frac{1}{r} \qquad (2.9\text{-}6)$$

From Eq 2.9-6

$$\ln_1 R = \ln_1(R+N) - {}_1\!\sum_{n=1}^{N} \frac{1}{R+N-n} \qquad (2.9\text{-}7)$$

From Eq 2.9-2, Eq 2.9-3, and Eq 2.9-7

$$\ln_{\Delta x} x = \ln_1\left(\frac{x}{\Delta x}\right) = \ln_1\left(\frac{x}{\Delta x} + N\right) - {}_1\!\sum_{n=1}^{N} \frac{1}{\frac{x}{\Delta x} + N - n} \qquad (2.9\text{-}8)$$

where

$N$ = large positive integer to assure good accuracy when the $\ln_{\Delta x} x$ Series, Eq 2.9-1, is used to calculate the function, $\ln_{\Delta x} x$

From Eq 2.9-8 an algorithm can be devised to calculate $\ln_{\Delta x}x$ to good accuracy for all x and Δx. The algorithm used in the lnd(n,Δx,x) calculation program, LNDX, at the end of the Appendix is as follows:

For $\mathrm{Re}(\frac{x}{\Delta x}) > 0$

If the value, $|\frac{x}{\Delta x}|$, is large, Eq 2.9-1 is used directly to calculate $\ln_{\Delta x}x$

If the value, $|\frac{x}{\Delta x}|$ , is not large, Eq 2.9-1 is used to calculate $\ln_{\Delta x}x$ indirectly using Eq 2.9-8 where

$$N = 20000 - \mathrm{int}(\mathrm{Re}(\frac{x}{\Delta x})) \tag{2.9-9}$$

where

int(p) = integer value of p

Note that using this value for N, the value, $|\frac{x}{\Delta x} + N|$ will be large and thus the calculated value of $\ln_1(\frac{x}{\Delta x} + N)$ will be accurate.

For $\mathrm{Re}(\frac{x}{\Delta x}) < 0$

Using a derived equation involving the function, $\ln_{\Delta x}x$

$$\ln_{\Delta x}x - \ln_{\Delta x}(\Delta x - x) = \begin{cases} 0\,, & \frac{x}{\Delta x} = \text{integer} \\ -\pi\cot\mathrm{an}\frac{\pi x}{\Delta x}\,, & \frac{x}{\Delta x} \neq \text{integer} \end{cases} \tag{2.9-10}$$

From Eq 2.9-2 and Eq 2.9-10

$$\ln_{\Delta x}x = \ln_1(\frac{x}{\Delta x}) = \ln_1(1-\frac{x}{\Delta x}) \qquad \text{for } \frac{x}{\Delta x} = \text{integer} \tag{2.9-11}$$

$$\ln_{\Delta x}x = \ln_1(\frac{x}{\Delta x}) = \ln_1(1-\frac{x}{\Delta x}) - \pi\cot\mathrm{an}\pi\frac{x}{\Delta x} \qquad \text{for } \frac{x}{\Delta x} \neq \text{integer} \tag{2.9-12}$$

Note that $\mathrm{Re}(1-\frac{x}{\Delta x}) > 0$

Since, $\mathrm{Re}(1-\frac{x}{\Delta x}) > 0$, $\ln_1(1-\frac{x}{\Delta x})$ can be calculated using Eq 2.9-1, or Eq 2.9-8. Then, $\ln_{\Delta x}x$ can be calculated using Eq 2.9-11 or Eq 2.9-12.

## Section 2.10: Derivation of the $\ln_{\Delta x}$ Series

The $\ln_{\Delta x}x$ Series is derived using the following procedure. Firstly, selected functions are converted to convergent series using the Maclauren asymptotic expansion. The series expansions of these functions are then added to one another to obtain a series expression for the discrete derivative of the function, $\ln_{\Delta x}x$, which is $\frac{1}{x}$. From this expression, the $\ln_{\Delta x}x$ Series is obtained. The derivation of the $\ln_{\Delta x}x$ Series is as follows:

$$z = \frac{x}{\Delta x} \tag{2.10-1}$$

$$\ln_{\Delta x}x = \ln_1(\frac{x}{\Delta x}) = \ln_1 z \tag{2.10-2}$$

The following function, $F_1(z)$, is selected since it has been previously found to be a good representation of the function, $\ln_1 z$, for large z. However, a more accurate (especially for smaller z) expression for $\ln_1 z$ is now being sought.

$$F_1(z) = \ln(z - \frac{1}{2}) + \gamma \tag{2.10-3}$$

where

$\gamma$ = Euler's Constant, .577215664…

Taking the discrete derivative of Eq 2.10-3

$$D_1F_1(z) = D_1[\ln(z - \frac{1}{2}) + \gamma] = \frac{\ln(z+1-\frac{1}{2}) + \gamma - \ln(z-\frac{1}{2}) - \gamma}{1} = \ln(z+\frac{1}{2}) - \ln(z-\frac{1}{2}) = \ln(\frac{2z+1}{2z-1}) \tag{2.10-4}$$

The function Maclauren asymptotic expansions which follow can be very time consuming to compute. The symbolic logic computation program, MAXIMA, was used to greatly reduce time and effort.

$$D_1[\ln(z - \frac{1}{2}) + \gamma] = \ln(\frac{2z+1}{2z-1}) = 2[\frac{1}{2z} + \frac{1}{3}(\frac{1}{2z})^3 + \frac{1}{5}(\frac{1}{2z})^5 + \frac{1}{7}(\frac{1}{2z})^7 + \ldots], \quad |z| > \frac{1}{2} \tag{2.10-5}$$

This series and all of the other series to follow are convergent for $|z| >> 1$

From Eq 2.10-5

$$D_1[\ln(z - \tfrac{1}{2}) + \gamma] = \frac{1}{z} + \frac{1}{12z^3} + \frac{1}{80z^5} + \frac{1}{448z^7} + \frac{1}{2304z^9} + \frac{1}{11264z^{11}} + \frac{1}{53248z^{13}} + \frac{1}{245760z^{15}} + \dots \quad (2.10\text{-}6)$$

$$D_1\ln_1 z = \frac{1}{z} \quad (2.10\text{-}7)$$

Substituting Eq 2.10-7 into Eq 2.10-6

$$D_1[\ln(z - \tfrac{1}{2}) + \gamma] = D_1\ln_1 z + \frac{1}{12z^3} + \frac{1}{80z^5} + \frac{1}{448z^7} + \frac{1}{2304z^9} + \frac{1}{11264z^{11}} + \frac{1}{53248z^{13}} + \frac{1}{245760z^{15}} + \dots \quad (2.10\text{-}8)$$

Rearranging the terms of Eq 2.10-8

$$D_1\ln_1 z = D_1[\ln(z - \tfrac{1}{2}) + \gamma] - \frac{1}{12z^3} - \frac{1}{80z^5} - \frac{1}{448z^7} - \frac{1}{2304z^9} - \frac{1}{11264z^{11}} - \frac{1}{53248z^{13}} - \frac{1}{245760z^{15}} - \dots \quad (2.10\text{-}9)$$

Note that the most significant error term of $D_1\ln_1 z$ decreases as $z^{-3}$

For z large and integrating both sides of Eq 2.10-9

$$\ln_1 z \approx \ln(z-\tfrac{1}{2}) + \gamma \ , \quad \gamma = \text{constant of integration, Euler's Constant} \quad (2.10\text{-}10)$$

The accuracy of Eq 2.10-10 increases as |z| increases in value

Selecting another function

$$F_2(z) = \frac{1}{(z - \frac{1}{2})^2} \quad (2.10\text{-}11)$$

$$D_1F_2(z) = D_1[\frac{1}{(z - \frac{1}{2})^2}] = \frac{1}{(z + \frac{1}{2})^2} - \frac{1}{(z - \frac{1}{2})^2} \quad (2.10\text{-}12)$$

Expanding Eq 2.10-12 into a Maclauren asymptotic series

$$D_1[\frac{1}{(z - \frac{1}{2})^2}] = -\frac{2}{z^3} - \frac{1}{5z^5} - \frac{3}{8z^7} - \frac{1}{8z^9} - \frac{5}{128z^{11}} - \frac{3}{256z^{13}} - \frac{7}{2048z^{15}} - \dots \quad (2.10\text{-}13)$$

Multiplying Eq 2.10-13 by $+\frac{\frac{1}{12}}{2} = +\frac{1}{24}$ and adding it to Eq 2.10-6

$$D_1[\ln(z - \tfrac{1}{2}) + \gamma + \frac{1}{24(z - \frac{1}{2})^2}] = \frac{1}{z} - \frac{7}{240z^5} - \frac{3}{224z^7} - \frac{11}{2304z^9} - \frac{13}{8448z^{11}} - \frac{25}{53248z^{13}} - \frac{17}{122880z^{15}} - \dots \quad (2.10\text{-}14)$$

Substituting Eq 2.10-7 into Eq 2.10-14

$$D_1[\ln(z-\frac{1}{2}) + \gamma + \frac{1}{24(z-\frac{1}{2})^2}] = D_1\ln_1 z \;-\frac{7}{240z^5} - \frac{3}{224z^7} - \frac{11}{2304z^9} - \frac{13}{8448z^{11}} - \frac{25}{53248z^{13}} - \frac{17}{122880z^{15}} - \ldots \quad (2.10\text{-}15)$$

Rearranging the terms of Eq 2.10-15

$$D_1\ln_1 z \; = D_1[\ln(z-\frac{1}{2}) + \gamma + \frac{1}{24(z-\frac{1}{2})^2}] = \; + \frac{7}{240z^5} + \frac{3}{224z^7} + \frac{11}{2304z^9} + \frac{13}{8448z^{11}} + \frac{25}{53248z^{13}} + \frac{17}{122880z^{15}} + \ldots$$

(2.10-16)

Note that the most significant error term of $D_1\ln_1 z$ decreases as $z^{-5}$

For large z and integrating both sides of Eq 2.10-16

$$\ln_1 z \approx \ln(z-\frac{1}{2}) + \gamma + \frac{1}{24(z-\frac{1}{2})^2} \quad , \quad \gamma = \text{constant of integration, Euler's Constant} \qquad (2.10\text{-}17)$$

The accuracy of Eq 2.10-17 increases as |z| increases in value

Compare Eq 2.10-16 to Eq 2.10-9 which is rewritten below

$$D_1\ln_1 z = D_1[\ln(z - \frac{1}{2}) + \gamma] - \frac{1}{12z^3} - \frac{1}{80z^5} - \frac{1}{448z^7} - \frac{1}{2304z^9} - \frac{1}{11264z^{11}} - \frac{1}{53248z^{13}} - \frac{1}{245760z^{15}} - \ldots$$

Observing the error terms, it is noted that the Eq 2.10-16 representation of the discrete derivative of $\ln_1 z$ is a more accurate representation than that of Eq 2.10-9 for large z. Therefore the expression for $\ln_1 z$ of Eq 2.10-17 will be a more accurate representation of $\ln_1 z$ than Eq 2.10-10. The process shown is continued using the following functions:

$$F_p(x) = \frac{1}{(z - \frac{1}{2})^{2(p-1)}} \quad , \qquad p = 2,3,4,5, \ldots \qquad (2.10\text{-}18)$$

After each process iteration, a more accurate representation of $D_1\ln_1 z$ and $\ln_1 z$ is found. The final result is the $\ln_{\Delta x}x$ Series (where $z = \frac{x}{\Delta x}$ ).

On the following page is a table which shows in tabular form the $\ln_{\Delta x}x$ Series derivation for nine terms of the series. The series constants, $K_1$ thru $K_7$ are evaluated. The other series constants can be evaluated using the same methodology presented or may be evaluated from the $C_m$ constants of the lnd(n,Δx,x) n≠1 Series. The derivation of the lnd(n,Δx,x) n≠1 Series is shown in the following section, Section 2.11

| # | **TABLE 2.10-1 Partial Series** | $K_m$ | $K_m$ | $z^{-1}$ | $z^{-3}$ | $z^{-5}$ | $z^{-7}$ | $z^{-9}$ | $z^{-11}$ | $z^{-13}$ | $z^{-15}$ |
|---|---|---|---|---|---|---|---|---|---|---|---|
| **1** | $D_1[\ln(z - \frac{1}{2}) + \gamma]$ | | | $+\frac{1}{z}$ | $+\frac{1}{12z^3}$ | $+\frac{1}{80z^5}$ | $+\frac{1}{448z^7}$ | $+\frac{1}{2304z^9}$ | $+\frac{1}{11264z^{11}}$ | $+\frac{1}{53248z^{13}}$ | $+\frac{1}{245760z^{15}}$ |
| **2** | $D_1[\frac{1}{(z - \frac{1}{2})^2}]$ | | | | $-\frac{2}{z^3}$ | $-\frac{1}{5z^5}$ | $-\frac{3}{8z^7}$ | $-\frac{1}{8z^9}$ | $-\frac{5}{128z^{11}}$ | $-\frac{3}{256z^{13}}$ | $-\frac{7}{2048z^{15}}$ |
| **3** | $D_1[\frac{1}{24(z - \frac{1}{2})^2}]$ | $\frac{\frac{+1}{12}}{2}$ | $K_1 = +\frac{1}{24}$ | | $-\frac{12}{z^3}$ | $-\frac{1}{24z^5}$ | $-\frac{1}{64z^7}$ | $-\frac{1}{192z^9}$ | $-\frac{5}{3072z^{11}}$ | $-\frac{1}{2048z^{13}}$ | $-\frac{7}{49152z^{15}}$ |
| **4** | $D_1[\ln(z - \frac{1}{2}) + \gamma + \frac{1}{24(z - \frac{1}{2})^2}]$ | | | $+\frac{1}{z}$ | | $-\frac{7}{240z^5}$ | $-\frac{3}{224z^7}$ | $-\frac{11}{2304z^9}$ | $-\frac{13}{8448z^{11}}$ | $-\frac{25}{53248z^{13}}$ | $-\frac{17}{122880z^{15}}$ |
| **5** | $D_1[\frac{1}{(z - \frac{1}{2})^4}]$ | | | | | $-\frac{4}{5z^5}$ | $-\frac{5}{z^7}$ | $-\frac{7}{2z^9}$ | $-\frac{15}{8z^{11}}$ | $-\frac{55}{64z^{13}}$ | $-\frac{91}{256z^{15}}$ |
| **6** | $D_1[\frac{-7}{960(z - \frac{1}{2})^4}]$ | $\frac{\frac{-7}{240}}{4}$ | $K_2 = -\frac{7}{960}$ | | | $+\frac{7}{240z^5}$ | $+\frac{7}{192z^7}$ | $+\frac{49}{1920z^9}$ | $+\frac{7}{512z^{11}}$ | $+\frac{77}{12288z^{13}}$ | $+\frac{637}{245760x^{15}}$ |

| # | **TABLE 2.10-1 Partial Series** | $K_m$ | $K_m$ | $z^{-1}$ | $z^{-3}$ | $z^{-5}$ | $z^{-7}$ | $z^{-9}$ | $z^{-11}$ | $z^{-13}$ | $z^{-15}$ |
|---|---|---|---|---|---|---|---|---|---|---|---|
| **7** | $D_1[\ln(z-\frac{1}{2})+\gamma+\sum_{m=1}^{2}\frac{K_m}{(z-\frac{1}{2})^{2m}}]$ | | | $+\frac{1}{z}$ | | | $+\frac{31}{1344z^7}$ | $+\frac{239}{11520z^9}$ | $+\frac{205}{16896z^{11}}$ | $+\frac{463}{79872z^{13}}$ | $+\frac{201}{81920x^{15}}$ |
| **8** | $D_1[\frac{1}{(z-\frac{1}{2})^6}]$ | | | | | | $-\frac{6}{z^7}$ | $-\frac{14}{z^9}$ | $-\frac{63}{4z^{11}}$ | $-\frac{99}{8z^{13}}$ | $-\frac{1001}{128z^{15}}$ |
| **9** | $D_1[\frac{+31}{8064(z-\frac{1}{2})^6}]$ | $\frac{\frac{+31}{1344}}{6}$ | $K_3=+\frac{31}{8064}$ | | | | $-\frac{31}{1344z^7}$ | $-\frac{31}{576z^9}$ | $-\frac{31}{512z^{11}}$ | $-\frac{341}{7168z^{13}}$ | $-\frac{4433}{147456z^{15}}$ |
| **10** | $D_1[\ln(z-\frac{1}{2})+\gamma+\sum_{m=1}^{3}\frac{K_m}{(z-\frac{1}{2})^{2m}}]$ | | | $+\frac{1}{z}$ | | | | $-\frac{127}{3840z^9}$ | $-\frac{409}{8448z^{11}}$ | $-\frac{23357}{559104z^{13}}$ | $-\frac{5089}{184320z^{15}}$ |
| **11** | $D_1[\frac{1}{(z-\frac{1}{2})^8}]$ | | | | | | | $-\frac{8}{z^9}$ | $-\frac{30}{z^{11}}$ | $-\frac{99}{2z^{13}}$ | $-\frac{429}{8z^{15}}$ |
| **12** | $D_1[\frac{-127}{30720(z-\frac{1}{2})^8}]$ | $\frac{\frac{-127}{3840}}{8}$ | $K_4=-\frac{127}{30720}$ | | | | | $+\frac{127}{3840z^9}$ | $+\frac{127}{1024z^{11}}$ | $+\frac{4191}{20480z^{13}}$ | $+\frac{18161}{81920x^{15}}$ |

| # | **TABLE 2.10-1 Partial Series** | $K_m$ | $\mathbf{K_m}$ | $z^{-1}$ | $z^{-3}$ | $z^{-5}$ | $z^{-7}$ | $z^{-9}$ | $z^{-11}$ | $z^{-13}$ | $z^{-15}$ |
|---|---|---|---|---|---|---|---|---|---|---|---|
| **13** | $D_1[\ln(z-\frac{1}{2})+\gamma+\sum_{m=1}^{4}\frac{K_m}{(z-\frac{1}{2})^{2m}}]$ | | | $+\frac{1}{z}$ | | | | | $+\frac{2555}{33792z^{11}}$ | $+\frac{910573}{5591040z^{13}}$ | $+\frac{143093}{737280x^{15}}$ |
| **14** | $D_1[\frac{1}{(z-\frac{1}{2})^{10}}]$ | | | | | | | | $-\frac{10}{z^{11}}$ | $-\frac{55}{z^{13}}$ | $-\frac{1001}{8z^{15}}$ |
| **15** | $D_1[\frac{+511}{67584(z-\frac{1}{2})^{10}}]$ | $\frac{+\frac{2555}{33792}}{10}$ | $K_5=+\frac{511}{67584}$ | | | | | | $-\frac{2555}{33792z^{11}}$ | $-\frac{2555}{6144z^{13}}$ | $-\frac{46501}{49152z^{15}}$ |
| **16** | $D_1[\ln(z-\frac{1}{2})+\gamma+\sum_{m=1}^{5}\frac{K_m}{(z-\frac{1}{2})^{2m}}]$ | | | $+\frac{1}{z}$ | | | | | | $-\frac{1414477}{5591040z^{13}}$ | $-\frac{277211}{368640z^{15}}$ |
| **17** | $D_1[\frac{1}{(z-\frac{1}{2})^{12}}]$ | | | | | | | | | $-\frac{12}{z^{13}}$ | $-\frac{91}{z^{15}}$ |

| # | **<u>TABLE 2.10-1</u> Partial Series** | ***$K_m$*** | **$K_m$** | **$z^{-1}$** | **$z^{-3}$** | **$z^{-5}$** | **$z^{-7}$** | **$z^{-9}$** | **$z^{-11}$** | **$z^{-13}$** | **$z^{-15}$** |
|---|---|---|---|---|---|---|---|---|---|---|---|
| **18** | $D_1[\frac{-1414477}{67092480(z-\frac{1}{2})^{12}}]$ | $\frac{\frac{-1414477}{5591040}}{12}$ | $K_6 = -\frac{1414477}{67092480}$ | | | | | | | $+\frac{1414477}{5591040z^{13}}$ | $+\frac{1414477}{737280z^{15}}$ |
| **19** | $D_1[\ln(z-\frac{1}{2}) + \gamma + \sum_{m=1}^{6}\frac{K_m}{(z-\frac{1}{2})^{2m}}]$ | | | $+\frac{1}{z}$ | | | | | | | $+\frac{57337}{49152z^{15}}$ |
| **20** | $D_1[\frac{1}{(z-\frac{1}{2})^{14}}]$ | | | | | | | | | | $-\frac{14}{z^{15}}$ |
| **21** | $D_1[\frac{+8191}{98304(z-\frac{1}{2})^{14}}]$ | $\frac{\frac{+57337}{49152}}{14}$ | $K_7 = +\frac{8191}{98304}$ | | | | | | | | $-\frac{57337}{49152z^{15}}$ |
| **22** | $D_1[\ln(z-\frac{1}{2}) + \gamma + \sum_{m=1}^{7}\frac{K_m}{(z-\frac{1}{2})^{2m}}]$ | | | $+\frac{1}{z}$ | | | | | | | |
| **23** | . . . | | | | | | | | | | |

From Table 2.10-1 line 22 for large |z|, The following equation is obtained:

$$D_1[\ln(z - \frac{1}{2}) + \gamma + \sum_{m=1}^{7} \frac{K_m}{(z - \frac{1}{2})^{2m}}] \approx \frac{1}{z} \qquad (2.10\text{-}19)$$

Rewriting Eq 2.10-7

$$D_1 \ln_1 z = \frac{1}{z}$$

Substituting Eq 2.10-7 into Eq 2.10-19

$$D_1[\ln(z - \frac{1}{2}) + \gamma + \sum_{m=1}^{7} \frac{K_m}{(z - \frac{1}{2})^{2m}}] \approx D_1 \ln_1 z \qquad (2.10\text{-}20)$$

Integrating both sides of Eq 2.10-20 and rearranging terms

$$\ln_1 z \approx \ln(z - \frac{1}{2}) + \gamma + \sum_{m=1}^{7} \frac{K_m}{(z - \frac{1}{2})^{2m}} + \quad , \quad \gamma = \text{constant of integration, Euler's Constant} \qquad (2.10\text{-}21)$$

Expanding Eq 2.10-21

$$\ln_1 z \approx \ln(z - \frac{1}{2}) + \gamma + \frac{1}{24(z - \frac{1}{2})^2} - \frac{7}{960(z - \frac{1}{2})^4} + \frac{31}{8064(z - \frac{1}{2})^6} - \frac{127}{30720(z - \frac{1}{2})^6} + \frac{511}{67584(z - \frac{1}{2})^{10}}$$
$$- \frac{1414477}{67092480(z - \frac{1}{2})^{12}} + \frac{8191}{98304(z - \frac{1}{2})^{14}} - \dots \qquad (2.10\text{-}22)$$

The accuracy of Eq 2.10-22 increases rapidly as |z| increases in value

And in General

$$\ln_1 z \approx \ln(z - \frac{1}{2}) + \gamma + \sum_{m=1}^{\infty} \frac{K_m}{(z - \frac{1}{2})^{2m}} \qquad (2.10\text{-}23)$$

The accuracy of Eq 2.10-23 increases rapidly as $|z|$ increases in value

$$z = \frac{x}{\Delta x} \qquad (2.10\text{-}24)$$

$$\ln_1 z = \ln_1(\frac{x}{\Delta x}) = \ln_{\Delta x} x \qquad (2.10\text{-}25)$$

The $\ln_{\Delta x} x$ Series is

$$\ln_{\Delta x} x \approx \ln(\frac{x}{\Delta x} - \frac{1}{2}) + \gamma + \sum_{m=1}^{\infty} \frac{K_m}{(\frac{x}{\Delta x} - \frac{1}{2})^{2m}} \qquad (2.10\text{-}26)$$

The accuracy of Eq 2.10-26 increases rapidly as $|\frac{x}{\Delta x}|$ increases in value

where

$\gamma$ = constant of integration, Euler's Constant, .577215664

$K_m = \ln_{\Delta x}$ Series constants

$m = 1,2,3,\ldots$

$K_1 = +\frac{1}{24}$ $\qquad K_5 = +\frac{511}{67584}$

$K_2 = -\frac{7}{960}$ $\qquad K_6 = -\frac{1414477}{67092480}$

$K_3 = +\frac{31}{8064}$ $\qquad K_7 = +\frac{8191}{98304}$

$K_4 = -\frac{127}{30720}$

The $\ln_{\Delta x}$ Series $K_m$ constants can be calculated from the $C_m$ constants. The relationship is as follows:

$$K_m = \frac{(2m-1)!\,C_m}{2^m(2m+1)!} = \frac{C_m}{2^m(2m+1)(2m)}\ ,\ \ m = 1,2,3,\ldots \tag{2.10-27}$$

There is another equation representation for the function, $\ln_{\Delta x}x$. It is derived from Eq 2.10-26 and Eq 2.10-27. It is as follows:

The $\ln_{\Delta x}x$ Series is

$$\ln_{\Delta x}x \approx \gamma + \ln\left(\frac{x}{\Delta x} - \frac{1}{2}\right) + \sum_{m=1}^{\infty} \frac{(2m-1)!\,C_m}{(2m+1)!\,2^{2m}\left(\frac{x}{\Delta x} - \frac{1}{2}\right)^{2m}} \tag{2.10-28}$$

accuracy improves rapidly for increasing $|\frac{x}{\Delta x}|$

$\gamma$ = constant of integration, Euler's Constant, .577215664

$C_m$ = series coefficients

$C_1 = +1$ $\qquad C_5 = +\frac{2555}{3}$

$C_2 = -\frac{7}{3}$ $\qquad C_6 = -\frac{1414477}{105}$

$C_3 = +\frac{31}{3}$ $\qquad C_7 = +286685$

$C_4 = -\frac{381}{5}$ $\qquad \ldots$

A method to calculate the $C_m$ constants from the Bernoulli Constants is given in Table 10 in the Appendix.

Other forms of the $\ln_{\Delta x}x$ Series are shown in Table 7 in the Appendix.

## Section 2.11: Derivation of the lnd(n,Δx,x) n≠1 Series

The lnd(n,Δx,x) n≠1 Series is derived using a procedure similar to that used by the $\ln_{\Delta x}x$ Series in Section 2.10. Firstly, selected functions are converted to convergent series using the Maclauren asymptotic expansion. The series expansions of these functions are then added to one another to obtain a series expression for the discrete derivative of the function, lnd(n,Δx,x) n≠1, which is $\frac{1}{x^n}$ where n≠1. From this expression, the lnd(n,Δx,x) n≠1 Series is obtained. The derivation of the lnd(n,Δx,x) n≠1 Series is as follows:

The following function, $F_0(n,\Delta x,x)$, is selected since it has been previously found to be a good representation of the function, lnd(n,Δx,x) n≠1, for large $|\frac{x}{\Delta x}|$. However, a more accurate (especially for smaller $|\frac{x}{\Delta x}|$) expression for the function, lnd(n,Δx,x) n≠1, is now being sought.

$$F_0(n,\Delta x,x) = -\frac{1}{(n-1)(x-\frac{\Delta x}{2})^{n-1}} \quad , \ n\neq 1 \tag{2.11-1}$$

Taking the discrete derivative of Eq 2.11-1 with respect to Δx

$$D_{\Delta x}F_0(n,\Delta x,x) = D_{\Delta x}[-\frac{1}{(n-1)(x-\frac{\Delta x}{2})^{n-1}}] \tag{2.11-2}$$

The function Maclauren asymptotic expansions which follow can be very time consuming to compute. The symbolic logic computation program, MAXIMA, was used to greatly reduce time and effort.

$$D_{\Delta x}F_0(n,\Delta x,x) = D_{\Delta x}[-\frac{1}{(n-1)(x-\frac{\Delta x}{2})^{n-1}}] = -\frac{1}{\Delta x(n-1)}\left[\frac{1}{(x+\frac{\Delta x}{2})^{n-1}} - \frac{1}{(x-\frac{\Delta x}{2})^{n-1}}\right] \tag{2.11-3}$$

From Eq 2.11-3

$$D_{\Delta x}F_0(n,\Delta x,x) = \frac{1}{x^n} + \frac{\Delta x^2 \prod_{r=0}^{1}(n+r)}{24\,x^{n+2}} + \frac{\Delta x^4 \prod_{r=0}^{3}(n+r)}{1920\,x^{n+4}} + \frac{\Delta x^6 \prod_{r=0}^{5}(n+r)}{322560\,x^{n+6}} + \frac{\Delta x^8 \prod_{r=0}^{7}(n+r)}{92897280\,x^{n+8}} + \frac{\Delta x^{10} \prod_{r=0}^{9}(n+r)}{40874803200\,x^{n+10}}$$

$$+ \frac{\Delta x^{12} \prod_{r=0}^{11}(n+r)}{25505877196800\,x^{12}} + \frac{\Delta x^{14} \prod_{r=0}^{13}(n+r)}{21424936845312000\,x^{14}} + \dots \qquad (2.11\text{-}4)$$

This series and all of the other series to follow are convergent for $\left|\frac{x}{\Delta x}\right| >> 1$

By definition

$$D_{\Delta x}[\mathrm{lnd}(n,\Delta x,x)] = -\frac{1}{x^n} \quad , \quad n \neq 1 \qquad (2.11\text{-}5)$$

Substituting Eq 2.11-1 and Eq 2.11-5 into Eq 2.11-4

$$D_{\Delta x}\left[-\frac{1}{(n-1)(x-\frac{\Delta x}{2})^{n-1}}\right] = -D_{\Delta x}[\mathrm{lnd}(n,\Delta x,x)] + \frac{\Delta x^2 \prod_{r=0}^{1}(n+r)}{24\,x^{n+2}} + \frac{\Delta x^4 \prod_{r=0}^{3}(n+r)}{1920\,x^{n+4}} + \frac{\Delta x^6 \prod_{r=0}^{5}(n+r)}{322560\,x^{n+6}}$$

$$+ \frac{\Delta x^8 \prod_{r=0}^{7}(n+r)}{92897280\,x^{n+8}} + \frac{\Delta x^{10} \prod_{r=0}^{9}(n+r)}{40874803200\,x^{n+10}} + \frac{\Delta x^{12} \prod_{r=0}^{11}(n+r)}{25505877196800\,x^{12}}$$

$$+ \frac{\Delta x^{14} \prod_{r=0}^{13}(n+r)}{21424936845312000\,x^{14}} + \dots \quad , \quad n \neq 1 \qquad (2.11\text{-}6)$$

Rearranging the terms of Eq 2.11-6

$$D_{\Delta x}[\mathrm{lnd}(n,\Delta x,x)] = D_{\Delta x}\left[\frac{1}{(n-1)(x-\frac{\Delta x}{2})^{n-1}}\right] + \frac{\Delta x^2 \prod_{r=0}^{1}(n+r)}{24\, x^{n+2}} + \frac{\Delta x^4 \prod_{r=0}^{3}(n+r)}{1920\, x^{n+4}} + \frac{\Delta x^6 \prod_{r=0}^{5}(n+r)}{322560\, x^{n+6}}$$

$$+ \frac{\Delta x^8 \prod_{r=0}^{7}(n+r)}{92897280\, x^{n+8}} + \frac{\Delta x^{10} \prod_{r=0}^{9}(n+r)}{40874803200\, x^{n+10}} + \frac{\Delta x^{12} \prod_{r=0}^{11}(n+r)}{25505877196800\, x^{12}}$$

$$+ \frac{\Delta x^{14} \prod_{r=0}^{13}(n+r)}{21424936845312000\, x^{14}} + \ldots \quad , \quad n \neq 1 \qquad (2.11\text{-}7)$$

Note that the most significant error term of $D_{\Delta x}[\mathrm{lnd}(n,\Delta x,x)]$, where $n\neq 1$, decreases as $x^{n+2}$

For $\left|\frac{x}{\Delta x}\right|$ large and integrating both sides of Eq 2.11-7

$$\mathrm{lnd}(n,\Delta x,x) \approx \frac{1}{(n-1)(x-\frac{\Delta x}{2})^{n-1}} + k \ , \qquad n \neq 1 \ , \quad k = \text{constant of integration} \qquad (2.11\text{-}8)$$

The accuracy of Eq 2.11-8 increases as $\left|\frac{x}{\Delta x}\right|$ increases in value

Selecting another function

$$F_1(n,\Delta x,x) = -\frac{1}{(n+1)(x-\frac{\Delta x}{2})^{n+1}} \qquad (2.11\text{-}9)$$

Taking the discrete derivative of Eq 2.11-9 with respect to $\Delta x$

$$D_{\Delta x}F_1(n,\Delta x,x) = D_{\Delta x}\left[-\frac{1}{(n+1)(x-\frac{\Delta x}{2})^{n+1}}\right] = -\frac{1}{\Delta x(n+1)}\left[\frac{1}{(x+\frac{\Delta x}{x})^{n+1}} - \frac{1}{(x-\frac{\Delta x}{x})^{n+1}}\right] \qquad (2.11\text{-}10)$$

Expanding Eq 2.11-10 into a Maclauren asymptotic series

$$D_{\Delta x}\left[-\frac{1}{(n+1)(x-\frac{\Delta x}{2})^{n+1}}\right] = \frac{1}{x^{n+2}} + \frac{\Delta x^2 \prod_{r=2}^{3}(n+r)}{24\, x^{n+4}} + \frac{\Delta x^4 \prod_{r=2}^{5}(n+r)}{1920\, x^{n+6}} + \frac{\Delta x^6 \prod_{r=2}^{7}(n+r)}{322560\, x^{n+8}} + \frac{\Delta x^8 \prod_{r=2}^{9}(n+r)}{92897280\, x^{n+10}}$$

$$+ \frac{\Delta x^{10} \prod_{r=2}^{11}(n+r)}{40874803200\, x^{n+12}} + \frac{\Delta x^{12} \prod_{r=2}^{13}(n+r)}{25505877196800\, x^{n+14}} \qquad (2.11\text{-}11)$$

Multiplying Eq 2.11-11 by $\dfrac{-\Delta x^2 \prod_{r=0}^{1}(n+r)}{24}$ and adding it to Eq 2.11-7

$$D_{\Delta x}[\text{lnd}(n,\Delta x,x)] + D_{\Delta x}\left[\frac{\Delta x^2 n}{24(x-\frac{\Delta x}{2})^{n+1}}\right] = D_{\Delta x}\left[\frac{1}{(n-1)(x-\frac{\Delta x}{2})^{n-1}}\right] - \frac{7\Delta x^4 \prod_{r=0}^{3}(n+r)}{5760\, x^{n+4}} - \frac{\Delta x^6 \prod_{r=0}^{5}(n+r)}{53760\, x^{n+6}}$$

$$- \frac{11\Delta x^8 \prod_{r=0}^{7}(n+r)}{92897280\, x^{n+8}} - \frac{13\Delta x^{10} \prod_{r=0}^{9}(n+r)}{30656102400\, x^{n+10}} - \frac{\Delta x^{12} \prod_{r=0}^{11}(n+r)}{1020235087872\, x^{n+12}} - \frac{-17\Delta x^{14} \prod_{r=0}^{13}(n+r)}{10712468422656000\, x^{n+14}}$$

where $n \neq 1$ (2.11-12)

Rearranging the terms of Eq 2.11-12

$$D_{\Delta x}[\text{lnd}(n,\Delta x,x)] = D_{\Delta x}\left[\frac{1}{(n-1)(x-\frac{\Delta x}{2})^{n-1}} - \frac{\Delta x^2 n}{24(x-\frac{\Delta x}{2})^{n+1}}\right] - \frac{7\Delta x^4 \prod_{r=0}^{3}(n+r)}{5760\, x^{n+4}} - \frac{\Delta x^6 \prod_{r=0}^{5}(n+r)}{53760\, x^{n+6}}$$

$$- \frac{11\Delta x^8 \prod_{r=0}^{7}(n+r)}{92897280\, x^{n+8}} - \frac{13\Delta x^{10} \prod_{r=0}^{9}(n+r)}{30656102400}\, x^{n+10} - \frac{\Delta x^{12} \prod_{r=0}^{11}(n+r)}{1020235087872\, x^{n+12}} - \frac{-17\Delta x^{14} \prod_{r=0}^{13}(n+r)}{10712468422656000\, x^{n+14}}$$

where $n \neq 1$ (2.11-13)

Note that the most significant error term of $D_{\Delta}x[\text{lnd}(n,\Delta x,x)]$, where $n\neq1$, decreases as $x^{n+4}$

For large $\left|\frac{x}{\Delta x}\right|$ and integrating both sides of Eq 2.11-13

$$\text{lnd}(n,\Delta x,x) \approx \frac{1}{(n-1)(x-\frac{\Delta x}{2})^{n-1}} - \frac{\Delta x^2 n}{24(x-\frac{\Delta x}{2})^{n+1}} + k \ , \quad n\neq 1 \ , \quad k = \text{constant of integration} \qquad (2.11\text{-}14)$$

The accuracy of Eq 2.11-14 increases as $\left|\frac{x}{\Delta x}\right|$ increases in value

Compare Eq 2.11-13 to Eq 2.11-7 which is rewritten below

$$D_{\Delta x}[\text{lnd}(n,\Delta x,x)] = D_{\Delta x}\left[\frac{1}{(n-1)(x-\frac{\Delta x}{2})^{n-1}}\right] + \frac{\Delta x^2 \prod_{r=0}^{1}(n+r)}{24\,x^{n+2}} + \frac{\Delta x^4 \prod_{r=0}^{3}(n+r)}{1920\,x^{n+4}} + \frac{\Delta x^6 \prod_{r=0}^{5}(n+r)}{322560\,x^{n+6}}$$

$$+ \frac{\Delta x^8 \prod_{r=0}^{7}(n+r)}{92897280\,x^{n+8}} + \frac{\Delta x^{10} \prod_{r=0}^{9}(n+r)}{40874803200\,x^{n+10}} + \frac{\Delta x^{12} \prod_{r=0}^{11}(n+r)}{25505877196800\,x^{12}}$$

$$+ \frac{\Delta x^{14} \prod_{r=0}^{13}(n+r)}{21424936845312000\,x^{14}} + \dots \quad , \quad n\neq 1$$

Observing the error terms, it is noted that the Eq 2.11-13 representation of the discrete derivative of lnd(n,Δx,x) is a more accurate representation than that of Eq 2.11-7 for large $\left|\frac{x}{\Delta x}\right|$. Therefore the expression for lnd(n,Δx,x) of Eq 2.11-14 will be a more accurate representation of lnd(n,Δx,x) than Eq 2.11-8. The process shown is continued using the following functions:

$$F_p(n,\Delta x,x) = -\frac{1}{(n+2p-1)(x-\frac{\Delta x}{2})^{n+2p-1}} \ , \qquad p = 2,3,4,5,\dots \qquad (2.11\text{-}15)$$

After each process iteration, a more accurate representation of $D_{\Delta x}$[lnd(n,Δx,x)] and lnd(n,Δx,x) is found. The final result is the lnd(n,Δx,x) n≠1 Series.

On the following page is a table which shows in tabular form the lnd(n,Δx,x) n≠1 Series derivation for eight terms of the series. The series constants, $K_0$ thru $K_7$ are evaluated. The other series constants can be evaluated using the same methodology presented or may be evaluated from the $C_m$ constants. A method to calculate the $C_m$ constants from the Bernoulli Constants is given in Table 10 in the Appendix.

| | **TABLE 2.11-1** **Partial Series / $K_m$** | $\frac{1}{x^n}$ | $\frac{1}{x^{n+2}}$ | $\frac{1}{x^{n+4}}$ | $\frac{1}{x^{n+6}}$ | $\frac{1}{x^{n+8}}$ | $\frac{1}{x^{n+10}}$ | $\frac{1}{x^{n+12}}$ | $\frac{1}{x^{n+14}}$ |
|---|---|---|---|---|---|---|---|---|---|
| **1** | $D_{\Delta x}[-\frac{1}{(n-1)(x-\frac{\Delta x}{2})^{n-1}}]$ | 1 | $\frac{\Delta x^2 \prod_{r=0}^{1}(n+r)}{24}$ | $\frac{\Delta x^4 \prod_{r=0}^{3}(n+r)}{1920}$ | $\frac{\Delta x^6 \prod_{r=0}^{5}(n+r)}{322560}$ | $\frac{\Delta x^8 \prod_{r=0}^{7}(n+r)}{92897280}$ | $\frac{\Delta x^{10} \prod_{r=0}^{9}(n+r)}{40874803200}$ | $\frac{\Delta x^{12} \prod_{r=0}^{11}(n+r)}{25505877196800}$ | $\frac{\Delta x^{14} \prod_{r=0}^{13}(n+r)}{21424936845312000}$ |
| **2** | $D_{\Delta x}[-\frac{1}{(n+1)(x-\frac{\Delta x}{2})^{n+1}}]$ Then multiply by $\frac{-\Delta x^2 \prod_{r=0}^{1}(n+r)}{24}$ | | 1 | $\frac{\Delta x^2 \prod_{r=2}^{3}(n+r)}{24}$ | $\frac{\Delta x^4 \prod_{r=2}^{5}(n+r)}{1920}$ | $\frac{\Delta x^6 \prod_{r=2}^{7}(n+r)}{322560}$ | $\frac{\Delta x^8 \prod_{r=2}^{9}(n+r)}{92897280}$ | $\frac{\Delta x^{10} \prod_{r=2}^{11}(n+r)}{40874803200}$ | $\frac{\Delta x^{12} \prod_{r=2}^{13}(n+r)}{25505877196800}$ |
| **3** | $D_{\Delta x}[\frac{\Delta x^2 \prod_{r=0}^{0}(n+r)}{24(x-\frac{\Delta x}{2})^{n+1}}]$ $K_1 = +\frac{1}{24}$ | | $\frac{-\Delta x^2 \prod_{r=0}^{1}(n+r)}{24}$ | $\frac{-\Delta x^4 \prod_{r=0}^{3}(n+r)}{576}$ | $\frac{-\Delta x^6 \prod_{r=0}^{5}(n+r)}{46080}$ | $\frac{-\Delta x^8 \prod_{r=0}^{7}(n+r)}{7741440}$ | $\frac{-\Delta x^{10} \prod_{r=0}^{9}(n+r)}{2229534720}$ | $\frac{-\Delta x^{12} \prod_{r=0}^{11}(n+r)}{980995276800}$ | $\frac{-\Delta x^{14} \prod_{r=0}^{13}(n+r)}{612141052723200}$ |

| # | **TABLE 2.11-1** **Partial Series / $K_m$** | $\frac{1}{x^n}$ | $\frac{1}{x^{n+2}}$ | $\frac{1}{x^{n+4}}$ | $\frac{1}{x^{n+6}}$ | $\frac{1}{x^{n+8}}$ | $\frac{1}{x^{n+10}}$ | $\frac{1}{x^{n+12}}$ | $\frac{1}{x^{n+14}}$ |
|---|---|---|---|---|---|---|---|---|---|
| **4** | $D_{\Delta x}\left[ -\frac{1}{(n-1)(x-\frac{\Delta x}{2})^{n-1}} + \sum_{m=1}^{1} \frac{K_m \Delta x^{2m} \prod_{r=0}^{2(m-1)}(n+r)}{(x-\frac{\Delta x}{2})^{n+2m-1}} \right]$ | 1 | | $\frac{-7\Delta x^4 \prod_{r=0}^{3}(n+r)}{5760}$ | $\frac{-\Delta x^6 \prod_{r=0}^{5}(n+r)}{53760}$ | $\frac{-11\Delta x^8 \prod_{r=0}^{7}(n+r)}{92897280}$ | $\frac{-13\Delta x^{10} \prod_{r=0}^{9}(n+r)}{30656102400}$ | $\frac{-\Delta x^{12} \prod_{r=0}^{11}(n+r)}{1020235087872}$ | $\frac{-17\Delta x^{14} \prod_{r=0}^{13}(n+r)}{10712468422656000}$ |
| **5** | $D_{\Delta x}\left[ -\frac{1}{(n+3)(x-\frac{\Delta x}{2})^{n+3}} \right]$ Then multiply by $\frac{7\Delta x^4 \prod_{r=0}^{3}(n+r)}{5760}$ | | | 1 | $\frac{\Delta x^2 \prod_{r=4}^{5}(n+r)}{24}$ | $\frac{\Delta x^4 \prod_{r=4}^{7}(n+r)}{1920}$ | $\frac{\Delta x^6 \prod_{r=4}^{9}(n+r)}{322560}$ | $\frac{\Delta x^8 \prod_{r=4}^{11}(n+r)}{92897280}$ | $\frac{\Delta x^{10} \prod_{r=4}^{13}(n+r)}{40874803200}$ |
| **6** | $D_{\Delta x}\left[ \frac{-7\Delta x^4 \prod_{r=0}^{2}(n+r)}{5760(x-\frac{\Delta x}{2})^{n+3}} \right]$ $K_2 = -\frac{7}{5760}$ | | | $\frac{7\Delta x^4 \prod_{r=0}^{3}(n+r)}{5760}$ | $\frac{7\Delta x^6 \prod_{r=0}^{5}(n+r)}{138240}$ | $\frac{7\Delta x^8 \prod_{r=0}^{7}(n+r)}{11059200}$ | $\frac{\Delta x^{10} \prod_{r=0}^{9}(n+r)}{265420800}$ | $\frac{\Delta x^{12} \prod_{r=0}^{11}(n+r)}{76441190400}$ | $\frac{\Delta x^{14} \prod_{r=0}^{13}(n+r)}{33634123776000}$ |

| # | **TABLE 2.11-1** **Partial Series / $K_m$** | $\frac{1}{x^n}$ | $\frac{1}{x^{n+2}}$ | $\frac{1}{x^{n+4}}$ | $\frac{1}{x^{n+6}}$ | $\frac{1}{x^{n+8}}$ | $\frac{1}{x^{n+10}}$ | $\frac{1}{x^{n+12}}$ | $\frac{1}{x^{n+14}}$ |
|---|---|---|---|---|---|---|---|---|---|
| **7** | $D_{\Delta x}[-\frac{1}{(n-1)(x-\frac{\Delta x}{2})^{n-1}}$ $+\sum_{m=1}^{2}\frac{K_m\Delta x^{2m}\prod_{r=0}^{2(m-1)}(n+r)}{(x-\frac{\Delta x}{2})^{n+2m-1}}]$ | 1 | | | $\frac{31\Delta x^6\prod_{r=0}^{5}(n+r)}{967680}$ | $\frac{239\Delta x^8\prod_{r=0}^{7}(n+r)}{464486400}$ | $\frac{41\Delta x^{10}\prod_{r=0}^{9}(n+r)}{12262440960}$ | $\frac{463\Delta x^{12}\prod_{r=0}^{11}(n+r)}{38258815795200}$ | $\frac{67\Delta x^{14}\prod_{r=0}^{13}(n+r)}{2380548538368000}$ |
| **8** | $D_{\Delta x}[-\frac{1}{(n+5)(x-\frac{\Delta x}{2})^{n+5}}]$ Then multiply by $\frac{-31\Delta x^6\prod_{r=0}^{5}(n+r)}{967680}$ | | | | 1 | $\frac{\Delta x^2\prod_{r=6}^{7}(n+r)}{24}$ | $\frac{\Delta x^4\prod_{r=6}^{9}(n+r)}{1920}$ | $\frac{\Delta x^6\prod_{r=6}^{11}(n+r)}{322560}$ | $\frac{\Delta x^8\prod_{r=6}^{13}(n+r)}{92897280}$ |

| # | **TABLE 2.11-1** **Partial Series / $K_m$** | $\frac{1}{x^n}$ | $\frac{1}{x^{n+2}}$ | $\frac{1}{x^{n+4}}$ | $\frac{1}{x^{n+6}}$ | $\frac{1}{x^{n+8}}$ | $\frac{1}{x^{n+10}}$ | $\frac{1}{x^{n+12}}$ | $\frac{1}{x^{n+14}}$ |
|---|---|---|---|---|---|---|---|---|---|
| **9** | $D_{\Delta x}\left[ \frac{31\Delta x^6 \prod_{r=0}^{4}(n+r)}{967680(x - \frac{\Delta x}{2})^{n+5}} \right]$ <br><br> $K_3 = + \frac{31}{967680}$ | | | | $\frac{-31\Delta x^6 \prod_{r=0}^{5}(n+r)}{967680}$ | $\frac{-31\Delta x^8 \prod_{r=0}^{7}(n+r)}{23224320}$ | $\frac{-31\Delta x^{10} \prod_{r=0}^{9}(n+r)}{1857945600}$ | $\frac{-31\Delta x^{12} \prod_{r=0}^{11}(n+r)}{312134860800}$ | $\frac{-31\Delta x^{14} \prod_{r=0}^{13}(n+r)}{89894839910400}$ |
| **10** | $D_{\Delta x}\left[ -\frac{1}{(n-1)(x - \frac{\Delta x}{2})^{n-1}} + \sum_{m=1}^{3} \frac{K_m \Delta x^{2m} \prod_{r=0}^{2(m-1)}(n+r)}{(x - \frac{\Delta x}{2})^{n+2m-1}} \right]$ | 1 | | | | $\frac{-127\Delta x^8 \prod_{r=0}^{7}(n+r)}{154828800}$ | $\frac{-409\Delta x^{10} \prod_{r=0}^{9}(n+r)}{30656102400}$ | $\frac{-23357\Delta x^{12} \prod_{r=0}^{11}(n+r)}{267811710566400}$ | $\frac{-727\Delta x^{14} \prod_{r=0}^{13}(n+r)}{2295528947712000}$ |
| **11** | $D_{\Delta x}\left[ -\frac{1}{(n+7)(x - \frac{\Delta x}{2})^{n+7}} \right]$ <br><br> Then multiply by <br><br> $\frac{+127\Delta x^8 \prod_{r=0}^{7}(n+r)}{154828800}$ | | | | | 1 | $\frac{\Delta x^2 \prod_{r=8}^{9}(n+r)}{24}$ | $\frac{\Delta x^4 \prod_{r=8}^{11}(n+r)}{1920}$ | $\frac{\Delta x^6 \prod_{r=8}^{13}(n+r)}{322560}$ |

| # | **TABLE 2.11-1** **Partial Series / $K_m$** | $\frac{1}{x^n}$ | $\frac{1}{x^{n+2}}$ | $\frac{1}{x^{n+4}}$ | $\frac{1}{x^{n+6}}$ | $\frac{1}{x^{n+8}}$ | $\frac{1}{x^{n+10}}$ | $\frac{1}{x^{n+12}}$ | $\frac{1}{x^{n+14}}$ |
|---|---|---|---|---|---|---|---|---|---|
| **12** | $D_{\Delta x}\left[\frac{-127\Delta x^8 \prod_{r=0}^{6}(n+r)}{154828800(x-\frac{\Delta x}{2})^{n+7}}\right]$ $K_4 = -\frac{127}{154828800}$ | | | | | $\frac{127\Delta x^8 \prod_{r=0}^{7}(n+r)}{154828800}$ | $\frac{127\Delta x^{10} \prod_{r=0}^{9}(n+r)}{3715891200}$ | $\frac{127\Delta x^{12} \prod_{r=0}^{11}(n+r)}{29271296000}$ | $\frac{127\Delta x^{14} \prod_{r=0}^{13}(n+r)}{49941577728000}$ |
| **13** | $D_{\Delta x}\left[-\frac{1}{(n-1)(x-\frac{\Delta x}{2})^{n-1}} + \sum_{m=1}^{4}\frac{K_m\Delta x^{2m}\prod_{r=0}^{2(m-1)}(n+r)}{(x-\frac{\Delta x}{2})^{n+2m-1}}\right]$ | 1 | | | | | $\frac{73\Delta x^{10} \prod_{r=0}^{9}(n+r)}{3503554560}$ | $\frac{910573\Delta x^{12} \prod_{r=0}^{11}(n+r)}{2678117105664000}$ | $\frac{143093\Delta x^{14} \prod_{r=0}^{13}(n+r)}{64274810535936000}$ |
| **14** | $D_{\Delta x}\left[-\frac{1}{(n+9)(x-\frac{\Delta x}{2})^{n+9}}\right]$ Then multiply by $-\frac{73\Delta x^{10}\prod_{r=0}^{9}(n+r)}{3503554560}$ | | | | | | 1 | $\frac{\Delta x^2 \prod_{r=10}^{11}(n+r)}{24}$ | $\frac{\Delta x^4 \prod_{r=10}^{13}(n+r)}{1920}$ |

| # | **TABLE 2.11-1** **Partial Series / $K_m$** | $\frac{1}{x^n}$ | $\frac{1}{x^{n+2}}$ | $\frac{1}{x^{n+4}}$ | $\frac{1}{x^{n+6}}$ | $\frac{1}{x^{n+8}}$ | $\frac{1}{x^{n+10}}$ | $\frac{1}{x^{n+12}}$ | $\frac{1}{x^{n+14}}$ |
|---|---|---|---|---|---|---|---|---|---|
| **15** | $D_{\Delta x}\left[\frac{+73\Delta x^{10}\prod_{r=0}^{8}(n+r)}{3503554560(x-\frac{\Delta x}{2})^{n+9}}\right]$ <br> $K_5 = +\frac{73}{3503554560}$ | | | | | | $\frac{-73\Delta x^{10}\prod_{r=0}^{9}(n+r)}{3503554560}$ | $\frac{-73\Delta x^{12}\prod_{r=0}^{11}(n+r)}{84085309440}$ | $\frac{-73\Delta x^{14}\prod_{r=0}^{13}(n+r)}{6726824755200}$ |
| **16** | $D_{\Delta x}\left[-\frac{1}{(n-1)(x-\frac{\Delta x}{2})^{n-1}} + \sum_{m=1}^{5}\frac{K_m\Delta x^{2m}\prod_{r=0}^{2(m-1)}(n+r)}{(x-\frac{\Delta x}{2})^{n+2m-1}}\right]$ | 1 | | | | | | $\frac{-1414477\Delta x^{12}\prod_{r=0}^{11}(n+r)}{2678117105664000}$ | $\frac{-25201\Delta x^{14}\prod_{r=0}^{13}(n+r)}{2921582297088000}$ |
| **17** | $D_{\Delta x}\left[-\frac{1}{(n+11)(x-\frac{\Delta x}{2})^{n+11}}\right]$ <br> Then multiply by <br> $\frac{+1414477\Delta x^{12}\prod_{r=0}^{11}(n+r)}{2678117105664000}$ | | | | | | | 1 | $\frac{\Delta x^{2}\prod_{r=12}^{13}(n+r)}{24}$ |

| # | **TABLE 2.11-1** Partial Series / $K_m$ | $\frac{1}{x^n}$ | $\frac{1}{x^{n+2}}$ | $\frac{1}{x^{n+4}}$ | $\frac{1}{x^{n+6}}$ | $\frac{1}{x^{n+8}}$ | $\frac{1}{x^{n+10}}$ | $\frac{1}{x^{n+12}}$ | $\frac{1}{x^{n+14}}$ |
|---|---|---|---|---|---|---|---|---|---|
| **18** | $D_{\Delta x}[\frac{-1414477\Delta x^{12}\prod_{r=0}^{10}(n+r)}{2678117105664000(x-\frac{\Delta x}{2})^{n+11}}]$ <br> $K_6 = -\frac{1414477}{2678117105664000}$ | | | | | | | $\frac{1414477\Delta x^{12}\prod_{r=0}^{11}(n+r)}{2678117105664000}$ | $\frac{1414477\Delta x^{14}\prod_{r=0}^{13}(n+r)}{64274810535936000}$ |
| **19** | $D_{\Delta x}[-\frac{1}{(n-1)(x-\frac{\Delta x}{2})^{n-1}} + \sum_{m=1}^{6}\frac{K_m\Delta x^{2m}\prod_{r=0}^{2(m-1)}(n+r)}{(x-\frac{\Delta x}{2})^{n+2m-1}}]$ | 1 | | | | | | | $\frac{8191\Delta x^{14}\prod_{r=0}^{13}(n+r)}{612141052723200}$ |
| **20** | $D_{\Delta x}[-\frac{1}{(n+13)(x-\frac{\Delta x}{2})^{n+13}}]$ <br> Then multiply by <br> $\frac{-8191\Delta x^{14}\prod_{r=0}^{13}(n+r)}{612141052723200}$ | | | | | | | | 1 |

| # | **TABLE 2.11-1** **Partial Series / $K_m$** | $\frac{1}{x^n}$ | $\frac{1}{x^{n+2}}$ | $\frac{1}{x^{n+4}}$ | $\frac{1}{x^{n+6}}$ | $\frac{1}{x^{n+8}}$ | $\frac{1}{x^{n+10}}$ | $\frac{1}{x^{n+12}}$ | $\frac{1}{x^{n+14}}$ |
|---|---|---|---|---|---|---|---|---|---|
| **21** | $D_{\Delta x}\left[\frac{+8191\Delta x^{14}\prod_{r=0}^{12}(n+r)}{61241052723200(x-\frac{\Delta x}{2})^{n+13}}\right]$ <br> $K_7 = +\frac{8191}{612141052723200}$ | | | | | | | | $\frac{-8191\Delta x^{14}\prod_{r=0}^{13}(n+r)}{612141052723200}$ |
| **22** | $D_{\Delta x}\left[-\frac{1}{(n-1)(x-\frac{\Delta x}{2})^{n-1}}\right.$ <br> $\left.+\sum_{m=1}^{7}\frac{K_m\Delta x^{2m}\prod_{r=0}^{2(m-1)}(n+r)}{(x-\frac{\Delta x}{2})^{n+2m-1}}\right]$ | 1 | | | | | | | |
| **23** | . <br> . <br> . | | | | | | | | |

From Table 2.11-1 line 22 for large $\left|\frac{x}{\Delta x}\right|$, The following equation is obtained:

$$D_{\Delta x}\left[-\frac{1}{(n-1)(x-\frac{\Delta x}{2})^{n-1}}+\sum_{m=1}^{7}\frac{K_m\Delta x^{2m}\prod_{r=0}^{2(m-1)}(n+r)}{(x-\frac{\Delta x}{2})^{n+2m-1}}\right]\approx\frac{1}{x^{n}} \qquad (2.11\text{-}16)$$

By definition

$$D_{\Delta x}\mathrm{lnd}(n,\Delta x,x)=-\frac{1}{x^{n}}\quad,\quad n\neq 1 \qquad (2.11\text{-}17)$$

Substituting Eq 2.11-17 into Eq 2.11-16

$$D_{\Delta x}\left[-\frac{1}{(n-1)(x-\frac{\Delta x}{2})^{n-1}}+\sum_{m=1}^{7}\frac{K_m\Delta x^{2m}\prod_{r=0}^{2(m-1)}(n+r)}{(x-\frac{\Delta x}{2})^{n+2m-1}}\right]\approx-D_{\Delta x}\mathrm{lnd}(n,\Delta x,x) \qquad (2.11\text{-}18)$$

Integrating both sides of EQ 2.11-18 and rearranging terms

$$\mathrm{lnd}(n,\Delta x,x)\approx\frac{1}{(n-1)(x-\frac{\Delta x}{2})^{n-1}}-\sum_{m=1}^{7}\frac{K_m\Delta x^{2m}\prod_{r=0}^{2(m-1)}(n+r)}{(x-\frac{\Delta x}{2})^{n+2m-1}}+k \qquad (2.11\text{-}19)$$

where

$n\neq 1$

k = constant of integration

From Eq 2.11-19

The lnd(n,Δx,x) n≠1 Series is

$$\mathrm{lnd}(n,\Delta x,x) \approx -\frac{K_0}{(n-1)(x-\frac{\Delta x}{2})^{n-1}} - \sum_{m=1}^{\infty} \frac{K_m \Delta x^{2m} \prod_{r=0}^{2(m-1)}(n+r)}{(x-\frac{\Delta x}{2})^{n+2m-1}} + k \qquad (2.11\text{-}20)$$

where

$n \neq 1$

k = constant of integration

$K_m$ = Series constants

m = 1,2,3,…

$K_0 = -1$ $\qquad K_4 = -\frac{127}{154828800}$

$K_1 = +\frac{1}{24}$ $\qquad K_5 = +\frac{73}{3503554560}$

$K_2 = -\frac{7}{5760}$ $\qquad K_6 = -\frac{1414477}{2678117105664000}$

$K_3 = +\frac{31}{967680}$ $\qquad K_7 = +\frac{8191}{612141052723200}$

…

The lnd(n,Δx,x) n≠1 Series $K_m$ constants can be calculated from the $C_m$ constants.
The relationship is as follows:

$$K_m = \frac{C_m}{(2m+1)!2^{2m}} \quad , \quad m = 0,1,2,3,\ldots \qquad (2.11\text{-}21)$$

The accuracy of Eq 2.11-20 increases rapidly as $|\frac{x}{\Delta x}|$ increases in value

There is another equation representation for the function, lnd(n,Δx,x) n≠1. It is derived from Eq 2.11-20 and Eq 2.11-21. Its form is more concise than that of Eq 2.11-20. It is as follows:

The lnd(n,Δx,x) n≠1 Series is:

$$\mathrm{lnd}(n,\Delta x,x) \approx -\sum_{m=0}^{\infty} \frac{\Gamma(n+2m-1)\left(\frac{\Delta x}{2}\right)^{2m} C_m}{\Gamma(n)(2m+1)!\left(x-\frac{\Delta x}{2}\right)^{n+2m-1}} + k \qquad (2.11\text{-}22)$$

accuracy improves rapidly for increasing $|\frac{x}{\Delta x}|$

where

$n \neq 1$

$k$ = constant of integration

$C_m$ = Series constants

$m = 1,2,3,\ldots$

$C_0 = -1$ $\qquad C_4 = -\frac{381}{5}$

$C_1 = +1$ $\qquad C_5 = +\frac{2555}{3}$

$C_2 = -\frac{7}{3}$ $\qquad C_6 = -\frac{1414477}{105}$

$C_3 = +\frac{31}{3}$ $\qquad C_7 = +286685$

. . .

A method to calculate the $C_m$ constants from the Bernoulli Constants is given in Table 10 in the Appendix.

Other forms of the lnd(n,Δx,x) n≠1 Series are shown in Table 7 in the Appendix.

## Section 2.12: Merging the lnd(n,Δx,x) n≠1 Series and the $\ln_{\Delta x}$ Series

In Section 2.10 the $\ln_{\Delta x}x$ Series was derived and in Section 2.11 the lnd(n,Δx,x) n≠1 Series was derived. From these individual derivations it would seem that there are two unique series needed to represent the lnd(n,Δx,x) function for all n,Δx, and x. The lnd(n,Δx,x) n≠1 series would be the general series while the the $\ln_{\Delta x}x$ Series ( $\ln_{\Delta x}x \equiv \text{lnd}(1,\Delta x,x)$ ) would be the special case series where n=1. This is not entirely true since the $\ln_{\Delta x}x$ Series can be derived from the lnd(n,Δx,x) n≠1 Series where n→1 as a limit. This will be shown below.

Derive the $\ln_{\Delta x}x$ Series from the lnd(n,Δx,x) n≠1 Series where n→1 as a limit

The lnd(n,Δx,x) n≠1 Series is:

$$\text{lnd}(n,\Delta x,x) \approx -\sum_{m=0}^{\infty} \frac{\Gamma(n+2m-1)\left(\frac{\Delta x}{2}\right)^{2m} C_m}{\Gamma(n)(2m+1)!\left(x-\frac{\Delta x}{2}\right)^{n+2m-1}} + k \qquad (2.12\text{-}1)$$

accuracy improves rapidly for increasing $\left|\frac{x}{\Delta x}\right|$

where

$n \neq 1$

k = constant of integration

Δx = x increment

$C_m$ = Series constants, m = 1,2,3,…

$C_0 = -1$ $\qquad$ $C_4 = -\frac{381}{5}$

$C_1 = +1$ $\qquad$ $C_5 = +\frac{2555}{3}$

$C_2 = -\frac{7}{3}$ $\qquad$ $C_6 = -\frac{1414477}{105}$

$C_3 = +\frac{31}{3}$ $\qquad$ $C_7 = +286685$ …

Differentiating and changing the form of Eq 2.12-1

$$D_{\Delta x}\text{lnd}(n,\Delta x,x) \approx D_{\Delta x}\left[ \frac{1}{(n-1)(x-\frac{\Delta x}{2})^{n-1}} - \sum_{m=1}^{\infty} \frac{\Gamma(n+2m-1)\left(\frac{\Delta x}{2}\right)^{2m} C_m}{\Gamma(n)(2m+1)!\left(x-\frac{\Delta x}{2}\right)^{n+2m-1}} \right] \tag{2.12-2}$$

The accuracy of Eq 2.12-2 increases rapidly as $\left|\frac{x}{\Delta x}\right|$ increases in value

$$D_{\Delta x}\text{lnd}(n,\Delta x,x) \approx D_{\Delta x}\left[ \frac{1}{(n-1)(x-\frac{\Delta x}{2})^{n-1}} \right] - \sum_{m=1}^{\infty} \frac{\Gamma(n+2m-1)\Delta x^{2m} Cm}{\Gamma(n)(2m+1)!2^{2m}} D_{\Delta x}\left[ \frac{1}{(x-\frac{\Delta x}{2})^{n+2m-1}} \right] \tag{2.12-3}$$

By definition

$$D_{\Delta x}\text{lnd}(n,\Delta x,x) = -\frac{1}{x^{n}} \tag{2.12-4}$$

Taking the discrete derivative and the limit as $n\to 1$ of Eq 2.12-3

$$\lim_{n\to 1}\left[-\frac{1}{x^{n}}\right] \approx \lim_{n\to 1}\left[\frac{1}{\Delta x(n-1)}\right]\left[\frac{1}{(x+\frac{\Delta x}{2})^{n-1}} - \frac{1}{(x-\frac{\Delta x}{2})^{n-1}}\right]$$

$$- \sum_{m=1}^{\infty} \frac{\Gamma(n+2m-1)\Delta x^{2m} C_m}{\Gamma(n)(2m+1)!2^{2m}} \lim_{n\to 1} \frac{1}{\Delta x}\left[\frac{1}{(x+\frac{\Delta x}{2})^{n+2m-1}} - \frac{1}{(x-\frac{\Delta x}{2})^{n+2m-1}}\right] \tag{2.12-5}$$

From Eq 2.12-5

$$-\frac{1}{x} \approx \frac{1}{\Delta x}\lim_{n\to 1}\left[\frac{1}{(n-1)}\right]\left[\frac{1}{(x+\frac{\Delta x}{2})^{n-1}} - \frac{1}{(x-\frac{\Delta x}{2})^{n-1}}\right]$$

$$- \sum_{m=1}^{\infty} \frac{\Gamma(n+2m-1)\Delta x^{2m} C_m}{\Gamma(n)(2m+1)!2^{2m}} \frac{1}{\Delta x}\left[\frac{1}{(x+\frac{\Delta x}{2})^{2m}} - \frac{1}{(x-\frac{\Delta x}{2})^{2m}}\right] \tag{2.12-6}$$

For n = 1

$$\frac{\Gamma(n+2m-1)}{\Gamma(n)} = \frac{\Gamma(2m)}{\Gamma(1)} = (2m-1)! \tag{2.12-7}$$

Using the Maxima symbolic logic program to find the limit term in Eq 2.12-6

(L'hospital's Rule is used by the Maxima symbolic logic program.)

$$\lim_{n\to 1}\left[\frac{1}{(n-1)}\right]\left[\frac{1}{(x+\frac{\Delta x}{2})^{n-1}} - \frac{1}{(x-\frac{\Delta x}{2})^{n-1}}\right] = -\left[\ln(2x+\Delta x) - \ln(2x-\Delta x)\right] \tag{2.12-8}$$

$$-\left[\ln(2x+\Delta x) - \ln(2x-\Delta x)\right] = -\left[\ln(2\Delta x)(\frac{x}{\Delta x}+\frac{1}{2}) - \ln(2\Delta x)(\frac{x}{\Delta x}-\frac{1}{2})\right.$$

$$-\left[\ln(2x+\Delta x) - \ln(2x-\Delta x)\right] = -\ln\left(\frac{\frac{x}{\Delta x}+\frac{1}{2}}{\frac{x}{\Delta x}-\frac{1}{2}}\right) = -\Delta x D_{\Delta x}\ln(\frac{x}{\Delta x}-\frac{1}{2}) \tag{2.12-9}$$

$$\frac{1}{\Delta x}\lim_{n\to 1}\left[\frac{1}{(n-1)}\right]\left[\frac{1}{(x+\frac{\Delta x}{2})^{n-1}} - \frac{1}{(x-\frac{\Delta x}{2})^{n-1}}\right] = -\frac{1}{\Delta x}\ln\left(\frac{\frac{x}{\Delta x}+\frac{1}{2}}{\frac{x}{\Delta x}-\frac{\Delta x}{2}}\right) = -D_{\Delta x}\ln(\frac{x}{\Delta x}-\frac{1}{2})$$

(2.12-10)

$$\lim_{n\to 1} D_{\Delta x}\left[\frac{1}{(n-1)(x-\frac{\Delta x}{2})^{n-1}}\right] = -D_{\Delta x}\ln(\frac{x}{\Delta x}-\frac{1}{2}) \tag{2.12-11}$$

Eq 2.12-11 could also have been obtained by observing the first series of Table 2.10-1 in Section 2.10 and the first series of Table 2.11-1 in Section 2.11.

For $n \to 1$

$$z = \frac{x}{\Delta x} \tag{2.12-12}$$

$$\frac{1}{\Delta x}D_1[\ln(z-\frac{1}{2})+\gamma] = D_{\Delta x}\left[-\frac{1}{(n-1)(x-\frac{\Delta x}{2})^{n-1}}\right] = +\frac{1}{x} + \frac{\Delta x^2}{12x^3} + \frac{\Delta x^4}{80x^5} + \frac{\Delta x^6}{448x^7} + \frac{\Delta x^8}{2304x^9} + \ldots$$

(2.12-13)

Then

$$-\,D_{\Delta x}\ln\left(\frac{x}{\Delta x}-\frac{1}{2}\right)=\lim_{n\to 1}D_{\Delta x}\left[\frac{1}{(n-1)(x-\frac{\Delta x}{2})^{n-1}}\right] \qquad (2.12\text{-}14)$$

$$\frac{1}{\Delta x}\left[\frac{1}{(x+\frac{\Delta x}{2})^{2m}}-\frac{1}{(x-\frac{\Delta x}{2})^{2m}}\right]=D_{\Delta x}\left[\frac{1}{(x-\frac{\Delta x}{2})^{2m}}\right] \qquad (2.12\text{-}15)$$

By definition

$$D_{\Delta x}\ln_{\Delta x}x\equiv \text{lnd}(1,\Delta x,x)=+\frac{1}{x} \qquad (2.12\text{-}16)$$

Substituting Eq 2.12-7, Eq 2.12-10 Eq 2.12-15, and Eq 2.12-16 into Eq 2.12-6

$$-\,D_{\Delta x}\ln_{\Delta x}x\approx -\,D_{\Delta x}\ln\left(\frac{x}{\Delta x}-\frac{1}{2}\right)-\sum_{m=1}^{\infty}\frac{(2m-1)!\Delta x^{2m}C_m}{(2m+1)!2^{2m}}D_{\Delta x}\left[\frac{1}{(x-\frac{\Delta x}{2})^{2m}}\right] \qquad (2.12\text{-}17)$$

Integrating both sides of Eq 2.12-17

$$\ln_{\Delta x}x\approx\ln\left(\frac{x}{\Delta x}-\frac{1}{2}\right)+\gamma+\sum_{m=1}^{\infty}\frac{(2m-1)!\Delta x^{2m}C_m}{(2m+1)!2^{2m}(x-\frac{\Delta x}{2})^{2m}} \qquad (2.12\text{-}18)$$

Simplifying Eq 2.12-18

**This is the $\ln_{\Delta x}x$ Series**

$$\boldsymbol{\ln_{\Delta x}x\approx\ln\left(\frac{x}{\Delta x}-\frac{1}{2}\right)+\gamma+\sum_{m=1}^{\infty}\frac{(2m-1)!C_m}{(2m+1)!2^{2m}(\frac{x}{\Delta x}-\frac{1}{2})^{2m}}} \qquad (2.12\text{-}19)$$

where

$\gamma$ = constant of integration, Euler's Constant .577215664…

$\Delta x$ = x increment

$C_m$ = $C_m$ constants

The accuracy of Eq 2.12-19 increases rapidly as $\left|\frac{x}{\Delta x}\right|$ increases in value

Note – The constant of integration was chosen to have the value of Euler's Constant so that

$$\ln_{\Delta x}x \;=\; {}_{\Delta x}\!\int_{\Delta x}^{x}\frac{1}{x}\,\Delta x \;=\; \Delta x \;{}_{\Delta x}\!\sum_{x=\Delta x}^{x-\Delta x}\frac{1}{x} \;=\; {}_{1}\!\sum_{x=1}^{M-1}\frac{1}{x}\quad,\quad x = M\Delta x\;,\quad M = 2,3,4,5,\ldots \qquad (2.12\text{-}20)$$

The derivations of the lnd(n,Δx,x) n≠1 Series and the $\ln_{\Delta x}x$ Series, until this section, were obtained independently of one another. The relationship of these series to each other was not entirely clear. Certainly, the $\ln_{\Delta x}x \equiv$ lnd(1,Δx,x) function was known to be a special case of the lnd(n,Δx,x) function and the lnd(n,Δx,x) n≠1 function was known to be the general part of the lnd(n,Δx,x) function. It was also known that together, both of these functions, the lnd(n,Δx,x) n≠1 function and the $\ln_{\Delta x}x$ function could be used to represent the lnd(n,Δx,x) function. However, in this section, some additional useful knowledge has been obtained. The $\ln_{\Delta x}x$ function series can be derived from the lnd(n,Δx,x) n≠1 function series. The two functions are seen to be more closely related than was previously realized. In fact, from this discovery, it now appears that one series can be used to represent the function, lnd(n,Δx,x). In effect, both the lnd(n,Δx,x) n≠1 Series and the $\ln_{\Delta x}x$ Series can be logically merged into an lnd(n,Δx,x) Series. The series equation representation may not be common, but it is logically sound and functional. It must represent the function, lnd(n,Δx,x), for all n including the special case (obtained from a limit) where n=1. This lnd(n,Δx,x) function series representation is as follows:

Rewriting the pertinent equations required for the lnd(n,Δx,x) Series derivation

$$D_{\Delta x}\mathrm{lnd}(n,\Delta x,x) \;=\; -\,\frac{1}{x^{n}}\;,\qquad n\neq 1$$

$$D_{\Delta x}\ln_{\Delta x}x \equiv \mathrm{lnd}(1,\Delta x,x) \;=\; +\,\frac{1}{x}\;,\quad n=1$$

$$D_{\Delta x}\left[\frac{1}{(r-1)(x-\frac{\Delta x}{2})^{r-1}}\right]\Bigg|_{r\to 1} \;=\; -\;D_{\Delta x}\ln\!\left(\frac{x}{\Delta x}-\frac{1}{2}\right)$$

From Eq 3

Letting r→n and integrating

$$-\alpha(n)\mathrm{lnd}(n,\Delta x,x) \approx {}_{\Delta x}\!\!\int\left\{D_{\Delta x}\left[\frac{1}{(r-1)(x-\frac{\Delta x}{2})^{r-1}}\right]\Bigg|_{r\to n} - \sum_{m=1}^{\infty}\frac{\Gamma(r+2m-1)\Delta x^{2m}C_m}{\Gamma(r)(2m+1)!2^{2m}}\,D_{\Delta x}\left[\frac{1}{(x-\frac{\Delta x}{2})^{r+2m-1}}\right]\Bigg|_{r\to n}\right\}\Delta x$$

(2.12-21)

where

$$\alpha(n) = \begin{cases} 1 & n = 1 \\ -1 & n \neq 1 \end{cases}$$

$$-\alpha(n)\text{lnd}(n,\Delta x,x) \approx \int_{\Delta x} \{D_{\Delta x}\left[\frac{1}{(r-1)(x-\frac{\Delta x}{2})^{r-1}}\right]\Big|_{r\to n} - \sum_{m=1}^{\infty} \frac{\Gamma(n+2m-1)(\frac{\Delta x}{2})^{2m} Cm}{\Gamma(n)(2m+1)!} D_{\Delta x}\left[\frac{1}{(x-\frac{\Delta x}{2})^{n+2m-1}}\right]\}\Delta x$$

(2.12-22)

**The lnd(n,Δx,x) Series**

$$\mathbf{lnd(n,\Delta x,x) \approx \left[\frac{1+\alpha(n)}{2}\right]\left[\ln\left(\frac{x}{\Delta x}-\frac{1}{2}\right)+\gamma\right]+\left[\frac{1-\alpha(n)}{2}\right]\left[\frac{1}{(n-1)(x-\frac{\Delta x}{2})^{n-1}}+K\right]}$$

$$\mathbf{+\ \alpha(n)\sum_{m=1}^{\infty}\frac{\Gamma(n+2m-1)\left(\frac{\Delta x}{2}\right)^{2m} C_m}{\Gamma(n)(2m+1)!\left(x-\frac{\Delta x}{2}\right)^{n+2m-1}}} \qquad (2.12\text{-}23)$$

**accuracy improves rapidly for increasing** $\mathbf{\left|\frac{x}{\Delta x}\right|}$

where $\alpha(n) = \begin{cases} 1 & n = 1 \\ -1 & n \neq 1 \end{cases}$

n,Δx,x = real or complex values
K = constant of integration for n≠1
γ = constant of integration for n=1, Euler's Constant .577215664…
Δx = x increment
$C_m$ = Series constants, m = 1,2,3,…

$C_1 = +1$ $\qquad$ $C_5 = +\frac{2555}{3}$

$C_2 = -\frac{7}{3}$ $\qquad$ $C_6 = -\frac{1414477}{105}$

$C_3 = +\frac{31}{3}$ $\qquad$ $C_7 = +286685$

$C_4 = -\frac{381}{5}$ $\qquad$ …

# CHAPTER 2 SOLVED PROBLEMS

## Example 2.1 - Evaluation of a polynomial summation using partial fractions

Evaluate the sum, $\displaystyle {}_{.5}\sum_{x=1}^{3}\frac{(x+1)^2}{x^2+1}$

From a partial fraction expansion

$$\frac{(x+1)^2}{x^2+1} = 1 + \frac{1}{x-j} + \frac{1}{x+j} \qquad 1)$$

$$ {}_{.5}\sum_{x=1}^{3}\frac{(x+1)^2}{x^2+1} = {}_{.5}\sum_{x=1}^{3}\left(1+\frac{1}{x-j}+\frac{1}{x+j}\right) = \frac{1}{.5}\int_{1}^{3+.5}\Delta x + \frac{1}{.5}\int_{1}^{3+.5}\frac{\Delta x}{x-j} + \frac{1}{.5}\int_{1}^{3+.5}\frac{\Delta x}{x+j} \qquad 2)$$

$$ {}_{.5}\sum_{x=1}^{3}\frac{(x+1)^2}{x^2+1} = \frac{1}{.5}\left[ -\,\mathrm{lnd}(0,.5,x)\Big|_{1}^{3.5} + \mathrm{lnd}(1,.5,x-j)\Big|_{1}^{3.5} + \mathrm{lnd}(1,.5,x+j)\Big|_{1}^{3.5} \right] \qquad 3)$$

Note – The $\mathrm{lnd}(n,\Delta x,x)$ function is preceded by a – for $n \neq 1$ and by a + if $n = 1$

$$ {}_{.5}\sum_{x=1}^{3}\frac{(x+1)^2}{x^2+1} = 2\,[\, -\,\mathrm{lnd}(0,.5,3.5) + \mathrm{lnd}(0,.5,1) + \mathrm{lnd}(1,.5,3.5-j) - \mathrm{lnd}(1,.5,1-j)\,] +$$

$$2\,[\, \mathrm{lnd}(1,.5,3.5+j) - \mathrm{lnd}(1,.5,1+j)\,]$$

$$ {}_{.5}\sum_{x=1}^{3}\frac{(x+1)^2}{x^2+1} = 5 + (\, 2.0063660477453580 + 1.2456233421750663i\,) +$$

$$(\, 2.0063660477453580 - 1.2456233421750663i\,)$$

$$\boldsymbol{ {}_{.5}\sum_{x=1}^{3}\frac{(x+1)^2}{x^2+1} = 9.0127320954907162} \qquad 4)$$

Checking $\displaystyle {}_{.5}\sum_{x=1}^{3}\frac{(x+1)^2}{x^2+1} = \frac{2^2}{2} + \frac{2.5^2}{3.25} + \frac{3^2}{5} + \frac{3.5^2}{7.25} + \frac{4^2}{10} = 9.0127320954907162$ Good check

**<u>Example 2.2</u> - Evaluation of an alternating sign summation**

Evaluate the sum, $\sum_{x=0}^{\infty}(-1)^x\left(\frac{1}{(2x+1)^7}\right) = 1-\frac{1}{3^7}+\frac{1}{5^7}-\frac{1}{7^7}+\frac{1}{9^7}-\frac{1}{11^7}+\dots$

$$\sum_{x=0}^{\infty}(-1)^x\left(\frac{1}{(2x+1)^7}\right) = \sum_{x=0}^{\infty}(-1)^x\left(\frac{1}{2^7}\right)\left(\frac{1}{(x+.5)^7}\right) \qquad 1)$$

$$\sum_{x=0}^{\infty}(-1)^x\left(\frac{1}{(2x+1)^7}\right) = \frac{1}{128}\left[\sum_{x=0}^{\infty}{}_2\frac{1}{(x+.5)^7} - \sum_{x=1}^{\infty}{}_2\frac{1}{(x+.5)^7}\right] \qquad 2)$$

$$\sum_{x=0}^{\infty}(-1)^x\left(\frac{1}{(2x+1)^7}\right) = \frac{1}{128}\left[\frac{1}{2}\,\mathrm{lnd}(7,2,x+.5)\Big|_{x=0} - \frac{1}{2}\,\mathrm{lnd}(7,2,x+.5)\Big|_{x=1}\right]$$

$$\sum_{x=0}^{\infty}(-1)^x\left(\frac{1}{(2x+1)^7}\right) = \frac{1}{256}\left[\,\mathrm{lnd}(7,2,.5) - \mathrm{lnd}(7,2,1.5)\,\right]$$

$$\sum_{x=0}^{\infty}(-1)^x\left(\frac{1}{(2x+1)^7}\right) = \frac{255.88595401997821683}{256}$$

$$\boldsymbol{\sum_{x=0}^{\infty}(-1)^x\left(\frac{1}{(2x+1)^7}\right) = .9995545078905399} \qquad 3)$$

Checking

Using a computer program to calculate the summation

$\sum_{x=0}^{\infty}(-1)^x\left(\frac{1}{(2x+1)^7}\right) = .9995545078905399$ Good check

**<u>Example 2.3</u> - Evaluation of tan(ax)**

<u>1) Find the value of tan(2.73)</u>

$$\tan(ax) = \frac{\pi}{\text{lnd}(1,1,1-\frac{ax}{\pi}) - \text{lnd}(1,1,\frac{ax}{\pi})} \quad 1)$$

$$\tan(2.73) = \frac{\pi}{\text{lnd}(1,1,1-\frac{2.73}{\pi}) - \text{lnd}(1,1,\frac{2.73}{\pi})} \quad 2)$$

$$\tan(2.73) = \frac{\pi}{-7.43573175290067708225 + .23892690272816824579}$$

**tan(2.73) = - .43652603050845385099** 3)

Checking

tan(2.73) = - .43652603050845385099 Good check
(exact)

<u>2) Find the value of tan(-2.1 – 3.7i )</u>

$$\tan(ax) = \frac{\pi}{\text{lnd}(1,1,1-\frac{ax}{\pi}) - \text{lnd}(1,1,\frac{ax}{\pi})} \quad 4)$$

$$\tan(-2.1.3.7i) = \frac{\pi}{\text{lnd}(1,1,1-\frac{-2.1-3.7i}{\pi}) - \text{lnd}(1,1,\frac{-2.1-3.7i}{\pi})} \quad 5)$$

$$\tan(-2.1-3.7i) = \frac{\pi}{[1.08424001429150699008 + .77436573948229227571i \; -1.08089463414945316032 + 2.36534279371204579220i]}$$

$$\tan(-2.1-3.7i) = \frac{\pi}{.00334538014205382975 + 3.13970853319433806792i}$$

**tan(-2.1-3.7i ) = .00106614478671118932 – 1.00059895806122610947i** 6)

Checking

tan(-2.1-3.7i ) = .00106614478671118932 – 1.00059895806122610947i Good check
(exact)

**<u>Example 2.4</u> - Evaluation of $\tan_{\Delta x}(a,x)$**

Find the value of $\tan_{.5}(2,1+i)$

$$\tan_{\Delta x}(a,x) = \frac{\pi}{\text{lnd}(1,1,1-\frac{bx}{\pi}) - \text{lnd}(1,1,\frac{bx}{\pi})} \qquad 1)$$

$$b = \frac{\tan^{-1}a\Delta x}{\Delta x} \qquad 2)$$

$$b = \frac{\tan^{-1}1}{.5} = \frac{\pi}{4(.5)} = \frac{\pi}{2} \qquad 3)$$

$$\tan_{.5}(2,1+i) = \frac{\pi}{\text{lnd}(1,1,1-\frac{\pi[1+i]}{2\pi}) - \text{lnd}(1,1,\frac{\pi[1+i]}{2\pi})} \qquad 4)$$

$$\tan_{.5}(2,1+i) = \frac{\pi}{\text{lnd}(1,1,1-\frac{1+i}{2}) - \text{lnd}(1,1,\frac{1+i}{2})}$$

$$\tan_{.5}(2,1+i) = \frac{\pi}{0 - 2.88131903995502918531i}$$

$$\mathbf{\tan_{.5}(2,1+i) = 1.09033141072736823003i} \qquad 5)$$

## <u>Example 2.5</u> - Evaluation of the Gamma Function

<u>1) Find the value of $\Gamma(7)$</u>

$$\Gamma(x) = \sqrt{2\pi}\, e^{\lim_{n\to 0}\left[\frac{lnd(n,1,x) - lnd(0,1,x)}{n}\right]} \qquad 1)$$

$x = 7$

let $n = 10^{-17}$

$$\Gamma(7) = \sqrt{2\pi}\, e^{\frac{lnd(10^{-17},1,7) - lnd(0,1,7)}{10^{-17}}} \qquad 2)$$

$$\Gamma(7) = \sqrt{2\pi}\, e^{\frac{5.6603126788805428196\ 10^{-17}}{10^{-17}}}$$

**$\Gamma(7) = 719.9999999999999$** 3)

Checking

$\Gamma(7) = 6! = 720$ Good check

<u>2) Find the value of $\Gamma(1+i)$</u>

$$\Gamma(x) = \sqrt{2\pi}\, e^{\lim_{n\to 0}\left[\frac{lnd(n,1,x) - lnd(0,1,x)}{n}\right]} \qquad 4)$$

$x = 1+i$

let $n = 10^{-17}$

$$\Gamma(1+i) = \sqrt{2\pi}\, e^{\frac{lnd(10^{-17},1,1+i) - lnd(0,1,1+i)}{10^{-17}}} \qquad 5)$$

$$\Gamma(x) = \sqrt{2\pi}\, e^{\frac{[-1.569861732506529091 - .301640320467533195i]\ 10^{-17}}{10^{-17}}}$$

**$\Gamma(1+i) = .4980156681183560 - .1549498283018106i$** 6)

Checking using Gamma Function Tables

$\ln\Gamma(1+i) = -.650923199302 - .301640320468i$

$\Gamma(1+i) = e^{-.650923199302 - .301640320468i}$

$\Gamma(1+i) = .498015668118 - .154949828302i$ Good Check

**Example 2.6 - Evaluation of the summation,** $\displaystyle {}_{3}\sum_{x=1}^{\infty} \frac{1}{x} - {}_{3}\sum_{x=2}^{\infty} \frac{1}{x}$

Evaluate the summation, $\displaystyle {}_{3}\sum_{x=1}^{\infty} \frac{1}{x} - {}_{3}\sum_{x=2}^{\infty} \frac{1}{x} = 1 - \frac{1}{2} + \frac{1}{4} - \frac{1}{5} + \frac{1}{7} - \frac{1}{8} + \frac{1}{10} - \frac{1}{11} + \dots$ , using the lnd(n,Δx,x) function.

The lnd(n,Δx,x) function evaluates the following summation as shown

$$ {}_{\Delta x}\sum_{x=x_1}^{x_2} \frac{1}{x+a} = \frac{1}{\Delta x} \operatorname{lnd}(1,\Delta x,x+a)\Big|_{x_1}^{x_2+\Delta x} $$

$N \to \infty$
$\Delta x = 3$

$$ {}_{3}\sum_{x=1}^{\infty} \frac{1}{x} - {}_{3}\sum_{x=2}^{\infty} \frac{1}{x} = \lim_{N\to\infty} \left[ {}_{3}\sum_{x=1}^{N} \frac{1}{x} - {}_{3}\sum_{x=2}^{N+1} \frac{1}{x} \right] \qquad 1) $$

$$ {}_{3}\sum_{x=1}^{\infty} \frac{1}{x} - {}_{3}\sum_{x=2}^{\infty} \frac{1}{x} = \lim_{N\to\infty} \frac{1}{3} \left[ \operatorname{lnd}(1,3,x)\,)\Big|_{1}^{N+3} - \operatorname{lnd}(1,3,x)\,)\Big|_{2}^{N+4} \right] \qquad 2) $$

$$ {}_{3}\sum_{x=1}^{\infty} \frac{1}{x} - {}_{3}\sum_{x=2}^{\infty} \frac{1}{x} = \frac{1}{3} [- \operatorname{lnd}(1,3,1) + \operatorname{lnd}(1,3,2) ] + \frac{1}{3} \lim_{N\to\infty} [\operatorname{lnd}(1,3,N+3)\,) - \operatorname{lnd}(1,3,N+4) ] $$

But

$\lim_{N\to\infty} [\operatorname{lnd}(1,3,N+3)\,) - \operatorname{lnd}(1,3,N+4)] \to 0$

Then

$$ {}_{3}\sum_{x=1}^{\infty} \frac{1}{x} - {}_{3}\sum_{x=2}^{\infty} \frac{1}{x} = \frac{1}{3} [- \operatorname{lnd}(1,3,1) + \operatorname{lnd}(1,3,2) ] \qquad 3) $$

Evaluating using the lnd(n,Δx,x) computer computation program, LNDX

$$ {}_3\sum_{x=1}^{\infty}\frac{1}{x} - {}_3\sum_{x=2}^{\infty}\frac{1}{x} = \frac{1}{3}[2.55481811511927346238990698 - .74101875088505561179582872\ ] $$

$$ {}_3\sum_{x=1}^{\infty}\frac{1}{x} - {}_3\sum_{x=2}^{\infty}\frac{1}{x} = .60459978807807261686469275\ldots = \frac{\pi}{3\sqrt{3}}\ldots $$

$$ {}_3\sum_{x=1}^{\infty}\frac{1}{x} - {}_3\sum_{x=2}^{\infty}\frac{1}{x} = \mathbf{1 - \frac{1}{2} + \frac{1}{4} - \frac{1}{5} + \frac{1}{7} - \frac{1}{8} + \frac{1}{10} - \frac{1}{11} + \ldots = \frac{\pi}{3\sqrt{3}}\ldots} \qquad 4) $$

Checking using a computer solution to the summation, ${}_3\sum_{x=1}^{10000000}(\frac{1}{x} - \frac{1}{x+1})$

$$ {}_3\sum_{x=1}^{10000000}(\frac{1}{x} - \frac{1}{x+1}) = .6045997547447 $$

good check

This summation is slow to converge.

## <u>Example 2.7</u> - Evaluation of the Riemann and Hurwitz Zeta Functions

### Riemann Zeta Function

<u>1) Evaluate ζ(i)</u>

$\zeta(n) = \text{lnd}(n,1,1)$ 1)

$\zeta(i) = \text{lnd}(i,1,1) = .0033002236 - .4181554491i$

$\boldsymbol{\zeta(i) = .0033002236 - .4181554491i}$ 2)

Comment: This answer and the ones to follow will, arbitrarily, be to 10 place accuracy.

<u>2) Evaluate ζ(1+i)</u>

$\zeta(n) = \text{lnd}(n,1,1)$ 3)

$\zeta(1+i) = \text{lnd}(1+i,1,1) = .5821580597 - .9268485643i$

$\boldsymbol{\zeta(1+i) = .5821580597 - .9268485643i}$ 4)

<u>3) Evaluate ζ(-11.234)</u>

$\zeta(n) = \text{lnd}(n,1,1)$ 5)

$\zeta(-11.234) = \text{lnd}(-11.234,1,1) = .02272911368$

$\boldsymbol{\zeta(-11.234) = .02272911368}$ 6)

### Hurwitz Zeta Function

<u>1) Evaluate ζ(2,.25) + ζ(2,.75)</u>

$\zeta(n,x) = \text{lnd}(n,1,x)$ 7)

$\zeta(2,.25) + \zeta(2,.75) = \text{lnd}(2,1,.25) + \text{lnd}(2,1,.75)$ 8)

$\zeta(2,.25) + \zeta(2,.75) = 17.19732915450711073927 + 2.54187964767160649839$

$\zeta(2,.25) + \zeta(2,.75) = 19.73920880217871723766$

$\boldsymbol{\zeta(2,.25) + \zeta(2,.75) = 2\pi^2}$ 9)

## **Example 2.8** - Evaluation of the Digamma and Polygamma Functions

### Digamma Function

Find the value of $\Psi(1.2+3.1i)$

$\Psi(x) = \text{lnd}(1,1,x) - \gamma$ 1)

where

$\gamma = .5772157\ldots$ , Euler's Constant

$\Psi(1.2+3.1i) = \text{lnd}(1,1,1.2+3.1i) - \gamma$

$\Psi(1.2+3.1i)=1.72971244071551957500+1.34688217783110624669i-.57721566490153286060$

$\mathbf{\Psi(1.2+3.1i)= 1.15249677581398671439 + 1.34688217783110624669i}$ 2)

### Polygamma Function

Find $\psi^{(2)}(1.49)$

$\psi^{(n)}(x) = (-1)^{n+1}n!\zeta(n+1,x)$ , $n = 1, 2, 3, \ldots$ 3)

$\psi^{(2)}(1.49) = (-1)^3 2!\zeta(3,1.49) = -2\zeta(3,1.49) = -2\text{lnd}(3,1,1.49)$

$\psi^{(2)}(1.49) = -2\text{lnd}(3,1,1.49) = -2(.42153156879978985652)$

$\mathbf{\psi^{(2)}(1.49) = - .84306313759957971305}$ 4)

Find $\psi^{(3)}(1.49)$

$\psi^{(n)}(x) = (-1)^{n+1}n!\zeta(n+1,x)$ , $n = 1, 2, 3, \ldots$ 5)

$\psi^{(3)}(1.49) = (-1)^4 3!\zeta(4,1.49) = 6\zeta(4,1.49) = 6\text{lnd}(4,1,1.49)$

$\psi^{(3)}(1.49) = 6\text{lnd}(4,1,1.49) = 6(.24073275616715942218)$

$\mathbf{\psi^{(3)}(1.49) = 1.44439653700295653312}$ 6)

## **Example 2.9** - Evaluation of the summation, $\sum_{x=4.6}^{8.6} \ln(1+x)$

Find the value of the sum, $\sum_{x=4.6}^{8.6} \ln(1+x)$

Use the following equation

$$\sum_{x=x_1}^{x_2} \ln(1+x) = \left[ x(\ln x-1) + \frac{\mathrm{lnd}(1,1,x+.5)}{2} + \sum_{n=2}^{\infty} \frac{\mathrm{lnd}(2n-2,1,x+.5)}{(2n-1)(2)^{2n-2}} - \frac{\mathrm{lnd}(2n-1,1,x+.5)}{(2n-1)(2)^{2n-1}} \right] \Bigg|_{x_1}^{x_2+1} \qquad 1)$$

$\mathrm{Re}(\frac{1}{x}) \geq -1$ , $\frac{1}{x} \neq -1$

Accuracy increases rapidly for increasing |x|

Let $x_1 = 4.6$
$x_2 = 8.6$

$$\sum_{x=4.6}^{8.6} \ln(1+x) = \left[ x(\ln x-1) + \frac{\mathrm{lnd}(1,1,x+.5)}{2} + \sum_{n=2}^{\infty} \frac{\mathrm{lnd}(2n-2,1,x+.5)}{(2n-1)(2)^{2n-2}} - \frac{\mathrm{lnd}(2n-1,1,x+.5)}{(2n-1)(2)^{2n-1}} \right] \Bigg|_{4.6}^{9.6} \qquad 2)$$

$$\sum_{x=4.6}^{8.6} \ln(1+x) = [9.6(\ln(9.6-1) - 4.6(\ln(4.6)-1] + \frac{1}{2}\mathrm{lnd}(1,1,x+.5) \Big|_{4.6}^{9.6}$$

$$+ \frac{1}{12}\mathrm{lnd}(2,1,x+.5) \Big|_{4.6}^{9.6} - \frac{1}{24}\mathrm{lnd}(3,1,x+.5) \Big|_{4.6}^{9.6}$$

$$+ \frac{1}{80}\mathrm{lnd}(4,1,x+.5) \Big|_{4.6}^{9.6} - \frac{1}{160}\mathrm{lnd}(5,1,x+.5) \Big|_{4.6}^{9.6} + \ldots$$

Writing values to 12 decimal places

| | Term Values | Partial Sum |
|---|---|---|
| $\sum_{x=4.6}^{8.6} \ln(1+x)=$ | +12.112925745348 | 12.112925745348 |
| | – 2.419858996077 | 9.693066749271 |
| | + .367102414019 | 10.060169163290 |
| | – .009372999707 | 10.050796163583 |
| | + .000747788160 | 10.051543951743 |
| | – .000037150532 | 10.051506801211 |
| | + .000003177271 | 10.051509978482 |

$$\sum_{x=4.6}^{8.6} \ln(1+x) = \underline{10.051509}978$$ (8 place accuracy for 6 series terms) 3)

Checking

$$\sum_{x=4.6}^{8.6} \ln(1+x) = \ln(5.6) + \ln(6.6) + \ln(7.6) + \ln(8.6) + \ln(9.6)$$

$$\sum_{x=4.6}^{8.6} \ln(1+x) = 10.051509795799$$ Good check

(exact value)

## **Example 2.10** - **Derivation of an interesting Zeta Function relationship**

Find an interesting Zeta Function relationship

Consider the following series:

$$\frac{1}{x-1} = \frac{1}{x} + \frac{1}{x^2} + \frac{1}{x^3} + \frac{1}{x^4} + \frac{1}{x^5} + \ldots \quad , \qquad |x| > 1 \tag{1}$$

$$D_{\Delta x}\mathrm{lnd}(1,\Delta x,x-a) = +\frac{1}{x-a} \quad , \quad n = 1 \tag{2}$$

$$D_{\Delta x}\mathrm{lnd}(n,\Delta x,x-a) = -\frac{1}{(x-a)^n} \quad , \quad n \neq 1 \tag{3}$$

where

a = constant

$$\frac{1}{x-1} = +\left(\frac{1}{x}\right) - \left(-\frac{1}{x^2}\right) - \left(-\frac{1}{x^3}\right) - \left(-\frac{1}{x^4}\right) - \left(-\frac{1}{x^5}\right) - \ldots \tag{4}$$

From Eq 2 thru Eq 4

$$\frac{1}{x-1} = D_{\Delta x}\mathrm{lnd}(1,\Delta x,x) - D_{\Delta x}\mathrm{lnd}(2,\Delta x,x) - D_{\Delta x}\mathrm{lnd}(3,\Delta x,x) - D_{\Delta x}\mathrm{lnd}(4,\Delta x,x) - D_{\Delta x}\mathrm{lnd}(5,\Delta x,x) - \ldots$$

$$D_{\Delta x}\mathrm{lnd}(1,\Delta x,x-1) = D_{\Delta x}[\mathrm{lnd}(1,\Delta x,x) - \mathrm{lnd}(2,\Delta x,x) - \mathrm{lnd}(3,\Delta x,x) - \mathrm{lnd}(4,\Delta x,x) - \mathrm{lnd}(5,\Delta x,x) - \ldots] \tag{5}$$

Integrate both sides of the above equation with a constant of integration equal to zero

$$\mathrm{lnd}(1,\Delta x,x-1) = \mathrm{lnd}(1,\Delta x,x) - \mathrm{lnd}(2,\Delta x,x) - \mathrm{lnd}(3,\Delta x,x) - \mathrm{lnd}(4,\Delta x,x) - \mathrm{lnd}(5,\Delta x,x) - \ldots \tag{6}$$

Rearranging terms

$$\mathrm{lnd}(1,\Delta x,x) - \mathrm{lnd}(1,\Delta x,x-1) = \mathrm{lnd}(2,\Delta x,x) + \mathrm{lnd}(3,\Delta x,x) + \mathrm{lnd}(4,\Delta x,x) + \mathrm{lnd}(5,\Delta x,x) + \ldots$$

Then

$$\mathrm{lnd}(1,\Delta x,x) - \mathrm{lnd}(1,\Delta x,x-1) = \sum_{n=2}^{\infty} \mathrm{lnd}(n,\Delta x,x) \tag{7}$$

Use the above equation to find an interesting Zeta Function relationship

Let $\Delta x = 1$

$x = 2$

$$\text{lnd}(1,1,2) - \text{lnd}(1,1,1) = \sum_{n=2}^{\infty} \text{lnd}(n,1,2) \qquad 8)$$

$$\sum_{\Delta x \; x=x_i}^{\infty} \frac{1}{x^n} = \frac{1}{\Delta x} \text{lnd}(n,\Delta x,x_i) \;, \quad \text{for } \text{Re}(n) > 1 \qquad 9)$$

For n = 2,3,4,…
$\Delta x = 1$

$$\sum_{1\; x=x_i}^{\infty} \frac{1}{x^n} = \text{lnd}(n,1,x_i) \qquad 10)$$

$$1 + \sum_{1\; x=2}^{\infty} \frac{1}{x^n} = \sum_{1\; x=1}^{\infty} \frac{1}{x^n} \qquad 11)$$

From Eq 10 and Eq 11

$$1 + \text{lnd}(n,1,2) = \text{lnd}(n,1,1) = \zeta(n) \qquad 12)$$

where

$$\zeta(n) = \sum_{1\; x=1}^{\infty} \frac{1}{x^n} \;, \quad \text{The Riemann Zeta Function}$$

From Eq 12

$$\text{lnd}(n,\Delta x,2) = \zeta(n) - 1 \qquad 13)$$

Substituting Eq 13 into Eq 8

$$\text{lnd}(1,1,2) - \text{lnd}(1,1,1) = \sum_{n=2}^{\infty} \text{lnd}(n,1,2) = \sum_{n=2}^{\infty} [\zeta(n) - 1]$$

$$\text{lnd}(1,1,2) - \text{lnd}(1,1,1) = \sum_{n=2}^{\infty} [\zeta(n) - 1] \qquad 14)$$

Evaluating using the lnd(n,Δx,x) computer computation program, LNDX

$\text{lnd}(1,1,2) - \text{lnd}(1,1,1) = 1 - 0 = 1$ 15)

Thus

$$\sum_{n=2}^{\infty} [\zeta(n) - 1] = 1 \qquad 16)$$

Checking the above equation

| | | |
|---|---|---|
| 1 | $\zeta(2) - 1 =$ | .64493 |
| 2 | $\zeta(3) - 1 =$ | .20205 |
| 3 | $\zeta(4) - 1 =$ | .08232 |
| 4 | $\zeta(5) - 1 =$ | .03692 |
| 5 | $\zeta(6) - 1 =$ | .01734 |
| 6 | $\zeta(7) - 1 =$ | .00834 |
| 7 | $\zeta(8) - 1 =$ | .00407 |
| 8 | $\zeta(9) - 1 =$ | .00200 |
| 9 | $\zeta(10) - 1 =$ | .00099 |
| 10 | $\zeta(11) - 1 =$ | .00049 |
| | | .99945 (for 10 terms, 1 = exact value) |

**Example 2.11** - **A calculation using the Discrete Calculus Summation Equation**

Evaluate the summation, $\sum_{\frac{\pi}{4}\, x=0}^{\frac{\pi}{2}} \cos x$

Evaluating using the Discrete Calculus Summation Equation

$$\sum_{\Delta x\, x=x_1}^{x_2} f(x) = \frac{1}{\Delta x} \int_{\Delta x\, x_1-\frac{\Delta x}{2}}^{x_2+\frac{\Delta x}{2}} f(x)\Delta x + \sum_{n=1}^{\infty} b_m\, \Delta x^{2m-1} \frac{d^{2m-1}}{dx^{2m-1}} f(x)\Big|_{x_1-\frac{\Delta x}{2}}^{x_2+\frac{\Delta x}{2}} \qquad 1)$$

$$b_m = \frac{-C_m}{(2m+1)!2^{2m}} \qquad 2)$$

Let $f(x) = \cos x$

$\Delta x = \frac{\pi}{4}$

$x_1 = 0$

$x_2 = \frac{\pi}{2}$

Substituting

$$\sum_{\frac{\pi}{4}\, x=0}^{\frac{\pi}{2}} \cos x = \frac{1}{\frac{\pi}{4}} \int_{\frac{\pi}{4}\, 0-\frac{\pi}{8}}^{\frac{\pi}{2}+\frac{\pi}{8}} f(x)\Delta x + \sum_{n=1}^{\infty} b_m\, \Delta x^{2m-1} \frac{d^{2m-1}}{dx^{2m-1}} f(x)\Big|_{0-\frac{8}{2}}^{\frac{\pi}{2}+\frac{\pi}{8}} \qquad 3)$$

$$\sum_{\frac{\pi}{4}\, x=0}^{\frac{\pi}{2}} \cos x = \frac{4}{\pi} \sin x\Big|_{-\frac{\pi}{8}}^{\frac{5\pi}{8}} + \frac{1}{24}\left(\frac{\pi}{4}\right)^1 \sin x\Big|_{-\frac{\pi}{8}}^{\frac{5\pi}{8}} + \left(\frac{7}{5760}\right)\left(\frac{\pi}{4}\right)^3 \sin x\Big|_{-\frac{\pi}{8}}^{\frac{5\pi}{8}} +$$

$$\left(\frac{31}{967680}\right)\left(\frac{\pi}{4}\right)^5 \sin x\Big|_{-\frac{\pi}{8}}^{\frac{5\pi}{8}} + \left(\frac{127}{154828800}\right)\left(\frac{\pi}{4}\right)^7 \sin x\Big|_{-\frac{\pi}{8}}^{\frac{5\pi}{8}} +$$

$$\left(\frac{73}{3503554560}\right)\left(\frac{\pi}{4}\right)^9 \sin x\Big|_{-\frac{\pi}{8}}^{\frac{5\pi}{8}} + \left(\frac{1414477}{2678117105664000}\right)\left(\frac{\pi}{4}\right)^{11} \sin x\Big|_{-\frac{\pi}{8}}^{\frac{5\pi}{8}} +$$

$$\left(\frac{8191}{612141052723200}\right)\left(\frac{\pi}{4}\right)^{13} \sin x\Big|_{-\frac{\pi}{8}}^{\frac{5\pi}{8}} + \ldots \qquad 4)$$

$${}_{\frac{\pi}{4}}\sum_{x=0}^{\frac{\pi}{2}} \cos x =$$

| | Term Values | Partial Sum | Place Accuracy |
|---|---|---|---|
| 1 | 1.66356703456770220055 | 1.66356703456770220055 | 1 |
| 2 | .04275717304070962036 | 1.70632480760773162595 | 3 |
| 3 | .00076926424024611682 | 1.70709407184797774277 | 4 |
| 4 | .00001250862806468674 | 1.70710658047604242952 | 7 |
| 5 | .00000019756567220325 | 1.70710677804171463278 | 8 |
| 6 | .00000000309565996843 | 1.70710678113737460121 | 11 |
| 7 | .00000000004840445842 | 1.70710678118577905964 | 12 |
| 8 | .00000000000075645695 | 1.70710678118653551659 | 14 |

$${}_{\frac{\pi}{4}}\sum_{x=0}^{\frac{\pi}{2}} \mathbf{\cos x} = \mathbf{1.70710678118653551659} \qquad 5)$$

(14 place accuracy for 8 series terms)

Checking

$${}_{\frac{\pi}{4}}\sum_{x=0}^{\frac{\pi}{2}} \cos x = 1 + \frac{\sqrt{2}}{2} + 0 = 1.70710678118654752440$$

(exact value)

## <u>Example 2.12</u> - A calculation using the Alternating Sign Discrete Calculus Summation Equation

Evaluate the summation, $\displaystyle {}_{\frac{1}{3}}\sum_{x=0}^{5} (-1)^{3x} \sin 2x$

Evaluating using the Alternating Sign Discrete Calculus Summation Equation

$$\sum_{\substack{\Delta x \\ x=x_1}}^{x_2} (-1)^{\frac{x-x_1}{\Delta x}} f(x) = -\frac{1}{2} f(x) \Big|_{x_1}^{x_2+\Delta x} + \sum_{n=1}^{\infty} h_m (2\Delta x)^{2m-1} \frac{d^{2m-1}}{dx^{2m-1}} f(x) \Big|_{x_1}^{x_2+\Delta x} \qquad 1)$$

$$h_m = \left[ \frac{B_{2m} + \dfrac{C_m}{(2m+1)!2^{2m}}}{(2m)!} \right], \quad m = 1,2,3,\ldots \qquad 2)$$

$x_2 = x_1 + (2p-1)$ ,   $p = 1,2,3,\ldots$

Let $f(x) = \sin 2x$

$\Delta x = \frac{1}{3}$

$x_1 = 0$

$x_2 = 5$

$x_2 + \Delta x = 5 + \frac{1}{3} = \frac{16}{3}$

Substituting

$$ {}_{\frac{1}{3}}\sum_{x=0}^{5} (-1)^{3x} \sin 2x = -\frac{1}{2} \sin 2x \Big|_{0}^{\frac{16}{3}} + \sum_{n=1}^{\infty} h_m \left(\frac{2}{3}\right)^{2m-1} \frac{d^{2m-1}}{dx^{2m-1}} f(x) \Big|_{0}^{\frac{16}{3}}$$

$$ {}_{\frac{1}{3}}\sum_{x=0}^{5} (-1)^{3x} \sin 2x = -\frac{1}{2} \sin 2x \Big|_{0}^{\frac{16}{3}} + \frac{1}{8}\left(\frac{2}{3}\right)^1 (2)^1 \cos 2x \Big|_{0}^{\frac{16}{3}} + \left(\frac{1}{384}\right)\left(\frac{2}{3}\right)^3 (2)^3 \cos 2x \Big|_{0}^{\frac{16}{3}}) +$$

$$\left(\frac{1}{15360}\right)\left(\frac{2}{3}\right)^5 (2)^5 \cos 2x \Big|_{0}^{\frac{16}{3}} + \left(\frac{17}{10321920}\right)\left(\frac{2}{3}\right)^7 (2)^7 \cos 2x \Big|_{0}^{\frac{16}{3}} +$$

$$\left(\frac{31}{743178240}\right)\left(\frac{2}{3}\right)^{9}(2)^{9}\cos 2x\,\Big|_{0}^{\frac{16}{3}} + \left(\frac{691}{653996851200}\right)\left(\frac{2}{3}\right)^{11}(2)^{11}\cos 2x\,\Big|_{0}^{\frac{16}{3}} +$$

$$\left(\frac{5461}{204047017574400}\right)\left(\frac{2}{3}\right)^{13}(2)^{13}\cos 2x\,\Big|_{0}^{\frac{16}{3}} + \ldots \qquad 3)$$

$${}_{\frac{1}{3}}\sum_{x=0}^{5}(-1)^{3x}\sin 2x =$$

| | Term Values | Partial Sum | Place Accuracy |
|---|---|---|---|
| 1 | +.4731978784190540 | .4731978784905400 | 0 |
| 2 | - .2205015663958965 | $\underline{.2}526963120231575$ | 1 |
| 3 | - .0081667246813295 | $\underline{.244}5295873418280$ | 3 |
| 4 | - .0003629655413924 | $\underline{.2441}666218004355$ | 4 |
| 5 | - .0000163238470996 | $\underline{.2441}502979533359$ | 4 |
| 6 | - .0000007349880320 | $\underline{.2441495}629653038$ | 7 |
| 7 | - .0000000330972127 | $\underline{.24414952}98687910$ | 8 |
| 8 | - .0000000014904191 | $\underline{.2441495283}776718$ | 10 |

$$\mathbf{{}_{\frac{1}{3}}\sum_{x=0}^{5}(-1)^{3x}\sin 2x = \underline{.2441495283}776718} \qquad 4)$$

(10 place accuracy for 8 series terms)

Checking

Using a computer program to exactly calculate, ${}_{\frac{1}{3}}\sum_{x=0}^{5}(-1)^{3x}\sin 2x$

$${}_{\frac{1}{3}}\sum_{x=0}^{5}(-1)^{3x}\sin 2x = .2441495283073910169 \qquad \text{Good check}$$

(exact value)

# CHAPTER 3

# Area Calculation Using Discrete Closed Contour Summation and Integration in the Complex Plane

## Section 3.1: Description of Area Calculation Using Discrete Closed Contour Summation/Integration in the Complex Plane

Discrete mathematics in the complex plane has some unique features which are now specified. The complex plane grid is composed of squares. The allowed values, $z = x+jy$, occur only at a corner of a complex plane grid square. As a result, any plotted complex plane curves, which are also called contours, have a two dimentional “sample and hold” shape. Contours in the complex plane pass only through grid corner points (i.e. z points) and progress only along the horizontal or vertical sides of grid squares. All closed contours encircle an integer number of grid squares.

Each grid square is $\delta$ units in length. The area of each square is $\delta^2$ square units in area. See Diagrams 3.1-1 thru 3.1-2 below.

Diagram 3.1-1 The complex plane with a square grid

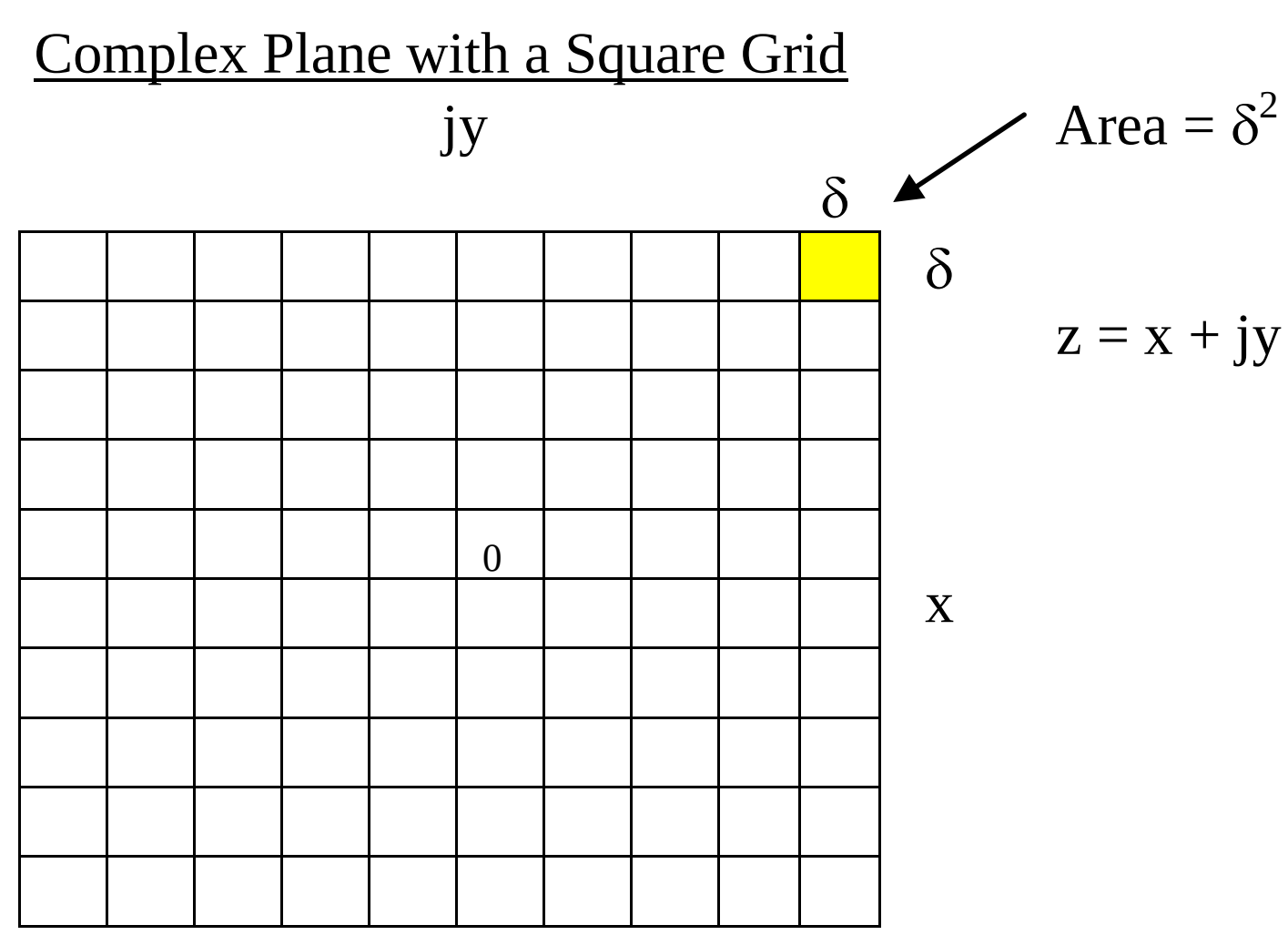

In the complex plane grid, each square side is considered to be a vector. The four vector sides are positioned head to tail in either a clockwise or counterclockwise direction forming a closed loop (i.e. the square). See Diagram 3.1-2 below.

Diagram 3.1-2 The complex plane grid square

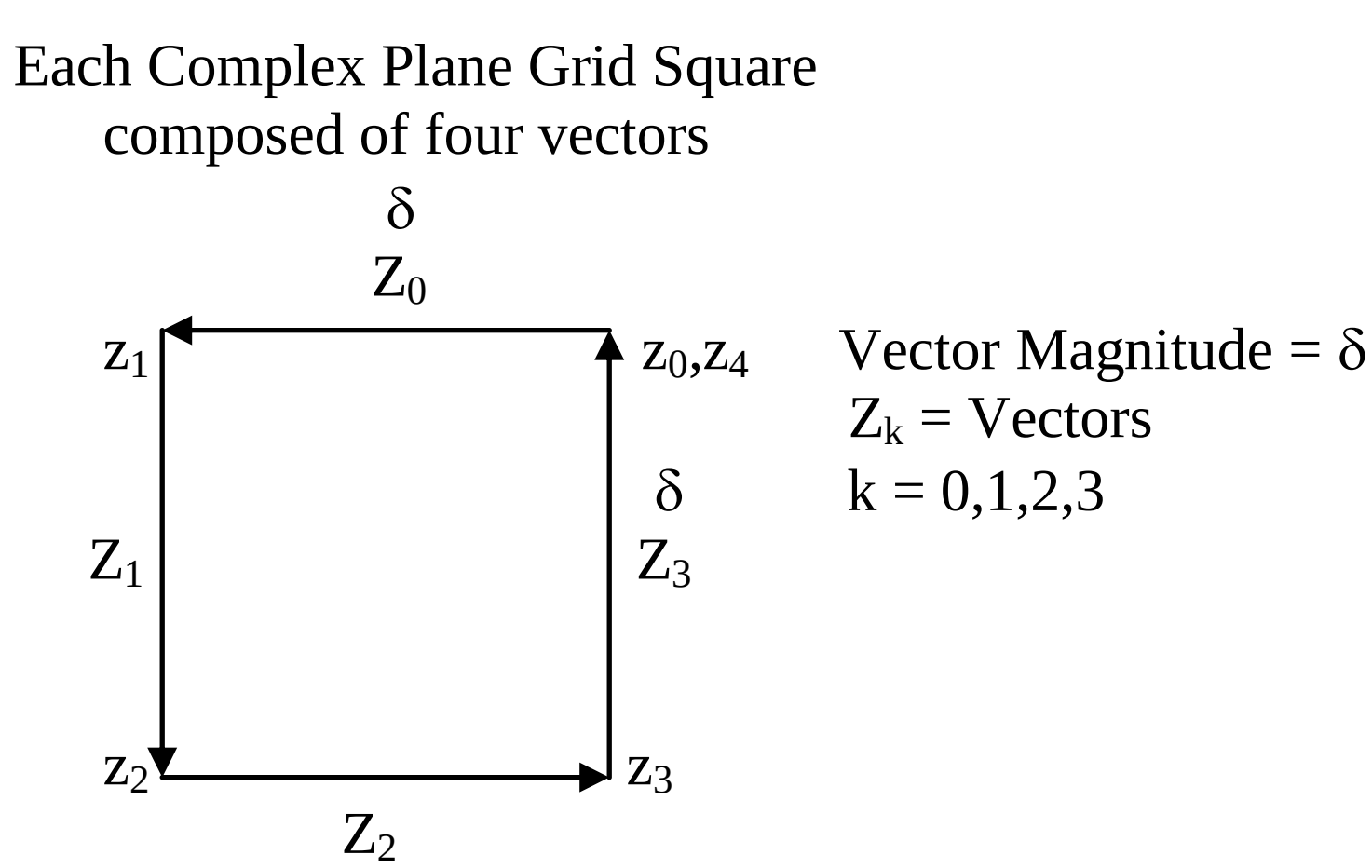

Note – The grid square shown in diagram 3.2 is composed of four vectors positioned head to tail in a counter-clockwise direction. The positioning of the vectors may also be head to tail in a clockwise direction.

Placing the vector squares of Diagram 3.1-2 together into a grouping of squares, internal vectors (i.e. those vectors not on the grouping outside perimeter contour) are seen to cancel. See Diagram 3.1-3 below.

Diagram 3.1-3 An example of a closed contour in the complex plane

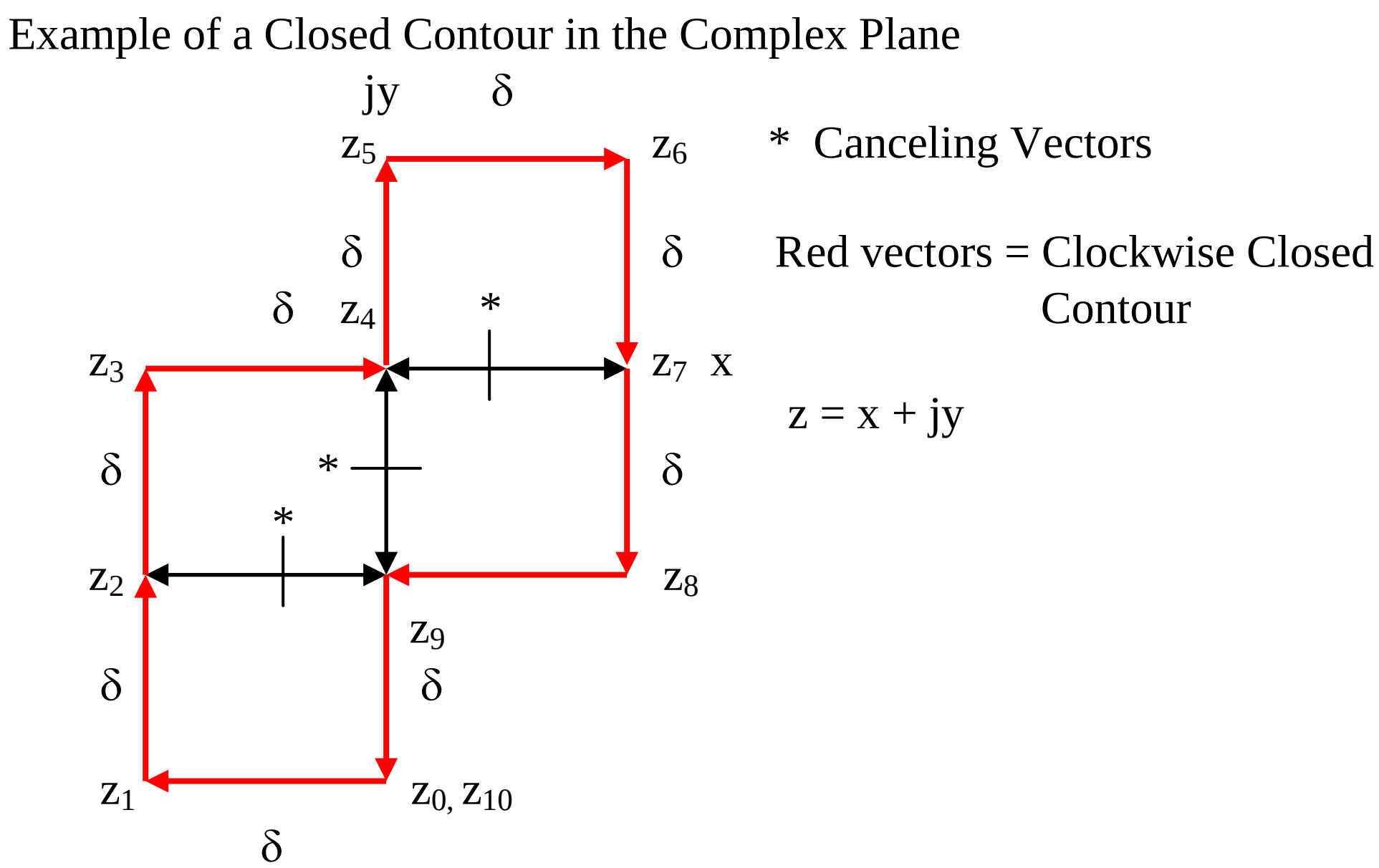

Note – The grid squares shown are composed of four vectors positioned head to tail in a clockwise direction. The vector head to tail positioning may also be in a counterclockwise direction.

It has been found that each summation, presented below, evaluated along the four vector sides of a complex plane grid square, either $\sum_{m=0}^{3} \left(\frac{z_{n+1} + z_n}{2}\right)^3 (z_{n+1} - z_n)$ for clockwise vectors or

$\sum_{m=0}^{3} \left(\frac{z_{n+1} + z_n}{2}\right)^3 (z_{n+1} - z_n)$ for counterclockwise vectors, has the same value no matter where the grid square is located in the complex plane. Changing the direction of vector rotation from clockwise to counterclockwise or vice versa changes the sign of the value. The summations shown above are of equal value but of opposite sign. Each summation can be related to the area of a grid square. By evaluating and totaling all of the summation terms along the vectors of the grid squares within a complex plane contour, the area within the contour can be obtained. Interestingly, only those summation terms associated with the grid square vectors forming a closed complex plane perimeter contour need to be evaluated and summed to obtain the area within the contour. This mathematical process is called contour integration or contour summation. Those grid square summation terms associated with vectors located within a closed complex plane contour cancel to zero. See Diagram 3.3. The * mark identifies the canceling vectors and their associated canceling summation terms. The derivation of the above two summations is shown in the following section, Section 3.2.

## Section 3.2: Derivation of Area Calculation Equations which use Discrete Closed Contour Summation/Integration in the Complex Plane

Consider the complex plane to be composed of δxδ, xy grid squares. $\Delta x = \delta$ and $\Delta y = \delta$. Also, randomly select a complex plane grid square. At each of the four corners of this grid square specify a grid value, $z_0$, $z_1$, $z_2$, and $z_3$. In the derivation which follows, $z_0$ is sometimes referred to as $z_4$ ($z_4 = z_0$). As described in Section 3.1, this grid square is composed of four vectors, each forming a side of the square. The four equal in magnitude grid square vectors are positioned head to tail either in a clockwise or counterclockwise direction. See Diagram 3.2-1 and Diagram 3.2-2 presented below. In the derivation to follow, the vector direction will be considered to be counterclockwise.

Diagram 3.2-1 A complex plane grid square

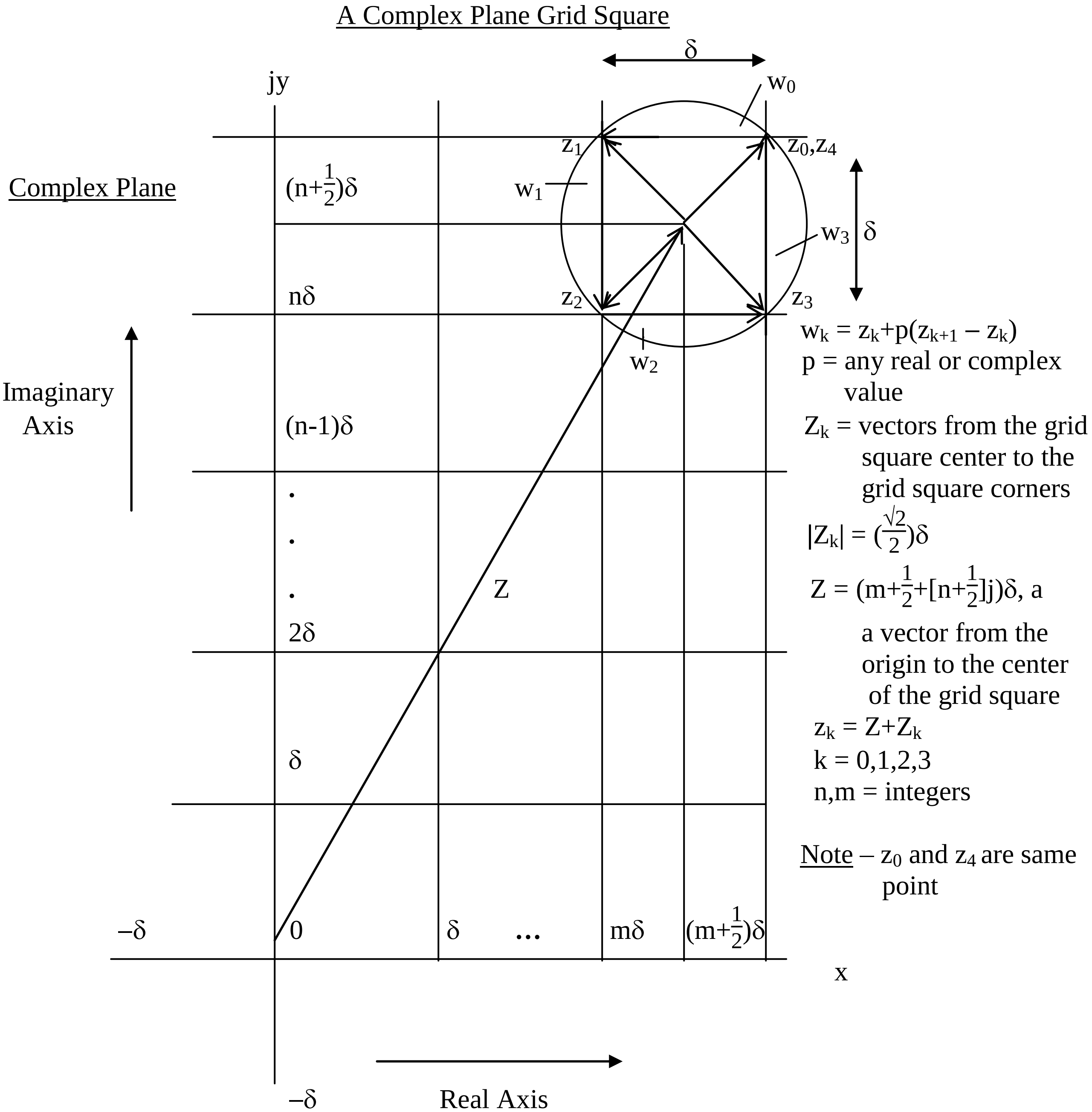

Diagram 3.2-2 The complex plane grid square in detail

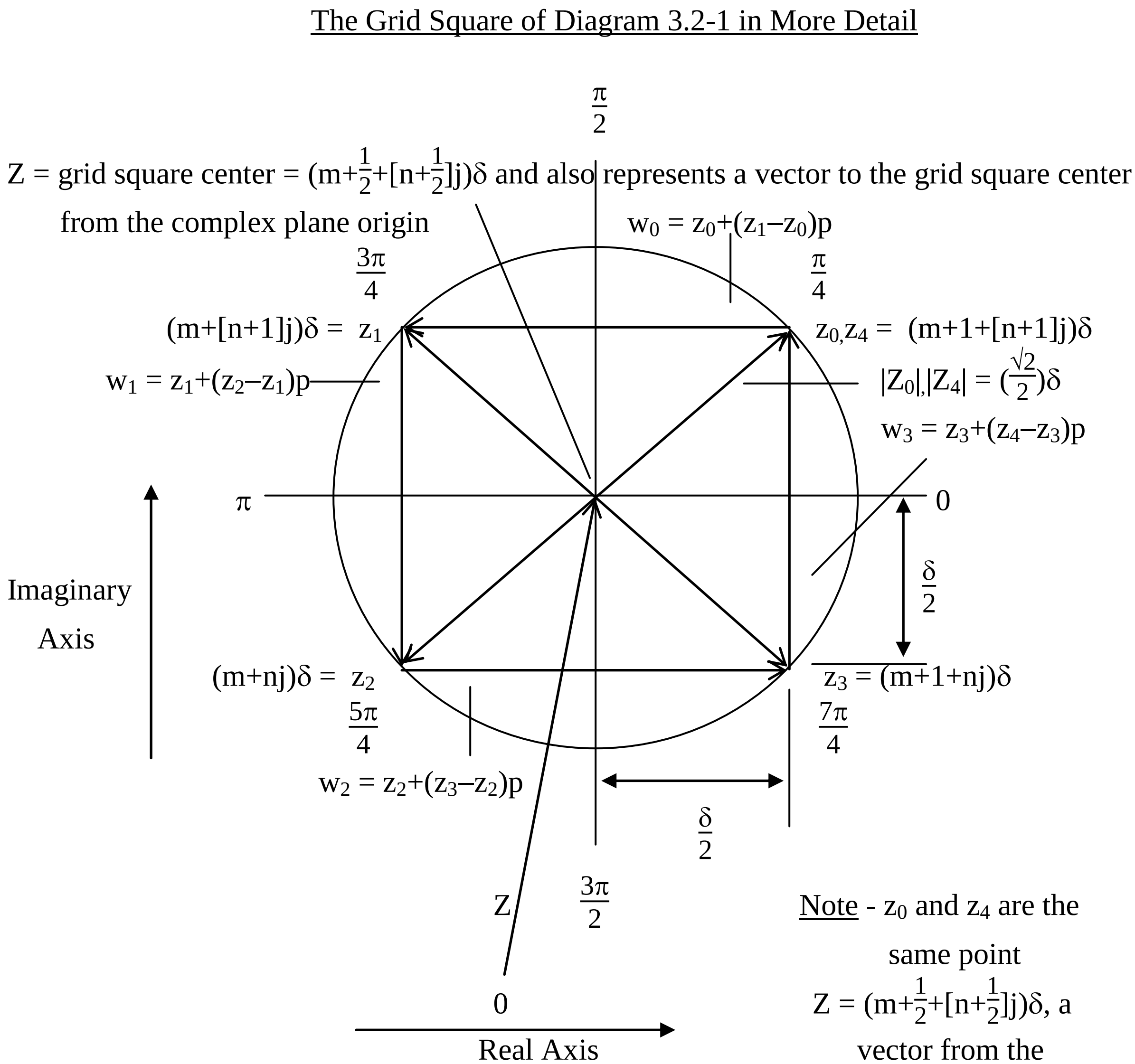

Note - $z_0$ and $z_4$ are the same point

$Z = (m+\frac{1}{2}+[n+\frac{1}{2}]j)\delta$, a vector from the origin to the center of the grid square

$Z_k$ = vectors from the grid square center to the grid square corners

$Z_k = (\frac{\sqrt{2}}{2})\delta e^{j(\frac{\pi}{4} + \frac{k\pi}{2})}$

$z_k = Z+Z_k$

$w_k = z_k+p(z_{k+1} - z_k)$

p = any real or complex value

$0 \leq p \leq 1$ for a point on a grid square vector

k = 0,1,2,3

n,m = integers

Using Diagram 3.2-1 and Diagram 3.2-2 show that the circular complex plane grid square vector function sum, $\sum_{k=0}^{3} ( z_k + p[z_{k+1} - z_k] )^3 (z_{k+1} - z_k)$, has the same value no matter where in the complex plane the grid square lies.

Derivation

$$z_k = (m+\frac{1}{2} + [n+\frac{1}{2}]j)\delta + (\frac{\sqrt{2}}{2})\delta e^{j(\frac{\pi}{4} + \frac{k\pi}{2})} \qquad (3.2\text{-}1)$$

where

m,n = integers

$-\infty \le m \le +\infty$

$-\infty \le n \le +\infty$

k = 0,1,2,3

$$w_k = z_k + p(z_{k+1} - z_k) \qquad (3.2\text{-}2)$$

where

p = any real or complex value

$0 \le p \le 1$ for a point on a grid square vector

k = 0,1,2,3

$$Z = (m+\frac{1}{2} + [n+\frac{1}{2}]j)\delta \qquad (3.2\text{-}3)$$

From Eq 3.2-1 and Eq 3.2-3

$$z_k = Z + (\frac{\sqrt{2}}{2})\delta e^{j(\frac{\pi}{4} + \frac{k\pi}{2})} \qquad (3.2\text{-}4)$$

Substituting Eq 3.2-4 into Eq 3.2-2

$$w_k = Z + (\frac{\sqrt{2}}{2})\delta e^{j(\frac{\pi}{4} + \frac{k\pi}{2})} + p\{Z + (\frac{\sqrt{2}}{2})\delta e^{j(\frac{\pi}{4} + \frac{[k+1]\pi}{2})} - Z - (\frac{\sqrt{2}}{2})\delta e^{j(\frac{\pi}{4} + \frac{k\pi}{2})}\} \qquad (3.2\text{-}5)$$

Simplifying Eq 3.2-5

$$w_k = Z + (\frac{\sqrt{2}}{2})\delta e^{j(\frac{\pi}{4} + \frac{k\pi}{2})} + (\frac{\sqrt{2}}{2})\delta p[e^{j(\frac{\pi}{4} + \frac{k\pi}{2})} e^{\frac{j\pi}{2}} - e^{j(\frac{\pi}{4} + \frac{k\pi}{2})}] \qquad (3.2\text{-}6)$$

$$w_k = Z + (\frac{\sqrt{2}}{2})\delta e^{j(\frac{\pi}{4} + \frac{k\pi}{2})}[ 1+p(e^{\frac{j\pi}{2}} - 1) ] \qquad (3.2\text{-}7)$$

Diagram 3.2-3  Vector diagram

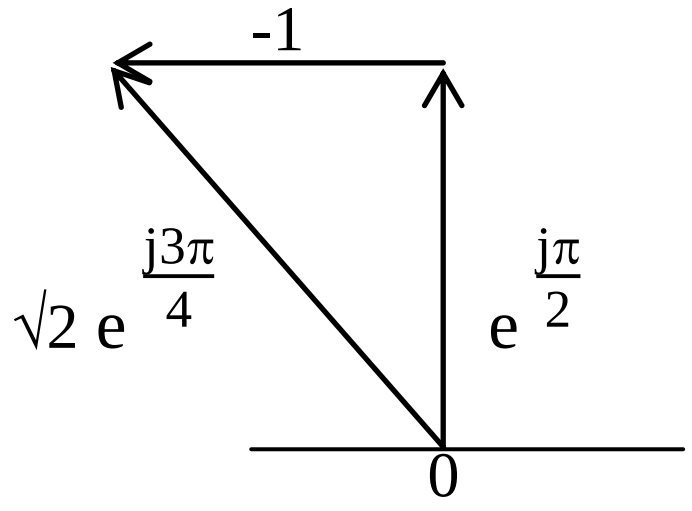

$$e^{\frac{j\pi}{2}} - 1 = \sqrt{2}\, e^{\frac{j3\pi}{4}} = j - 1 \qquad (3.2\text{-}8)$$

Substituting Eq 3.2-8 into Eq 3.2-7

$$w_k = Z + \left(\frac{\sqrt{2}}{2}\right)\delta e^{j\left(\frac{\pi}{4} + \frac{k\pi}{2}\right)}\left[\, 1 + p\sqrt{2}e^{\frac{j3\pi}{4}}\, \right]$$

Simplifying

$$w_k = Z + \delta e^{j\left(\frac{\pi}{4} + \frac{k\pi}{2}\right)}\left[\, \frac{\sqrt{2}}{2} + pe^{\frac{j3\pi}{4}}\, \right] \qquad (3.2\text{-}9)$$

Define the following complex constants

From Eq 3.2-9

$$C = \delta\left( \frac{\sqrt{2}}{2} + pe^{\frac{j3\pi}{4}} \right) e^{\frac{j\pi}{4}} \qquad (3.2\text{-}10)$$

From Eq 3.2-4

$$B = \left(\frac{\sqrt{2}}{2}\right)\delta e^{\frac{j\pi}{4}} \qquad (3.2\text{-}11)$$

Substituting Eq 3.2-11 into Eq 3.2-4

$$z_k = Z + Be^{\frac{jk\pi}{2}} \qquad (3.2\text{-}12)$$

Substituting Eq 3.2-10 into Eq 3.2-9

$$w_k = Z + Ce^{\frac{jk\pi}{2}} \qquad (3.2\text{-}13)$$

Consider the discrete circular integral, $\alpha_r$

$$\alpha_r = \oint_{0}^{4} w_k^{\,r}(z_{k+1} - z_k)\Delta k \qquad (3.2\text{-}14)$$

where

$w_k^{\,r} = (z_k + p[z_{k+1} - z_k]\,)^r$ , a function of the vector, $z_{k+1} - z_k$, $z_{k+1}$ at the vector head, $z_k$ at the tail

$\Delta k = 1$

$k = 0,1,2,3$

$r = 0,1,2,3,\ldots$

The summation below is equivalent to the above integral

$$\alpha_r = \sum_{k=0}^{3} w_k^{\,r}(z_{k+1} - z_k) \qquad (3.2\text{-}15)$$

where

$k = 0,1,2,3$

Find $\alpha_r$ by substituting Eq 3.2-12 and Eq 3.2-13 into Eq 3.2-14 and integrating counterclockwise along the vector sides of a grid square. See Diagram 3.7 below.

<u>Diagram 3.2-4 A complex plane grid square</u>

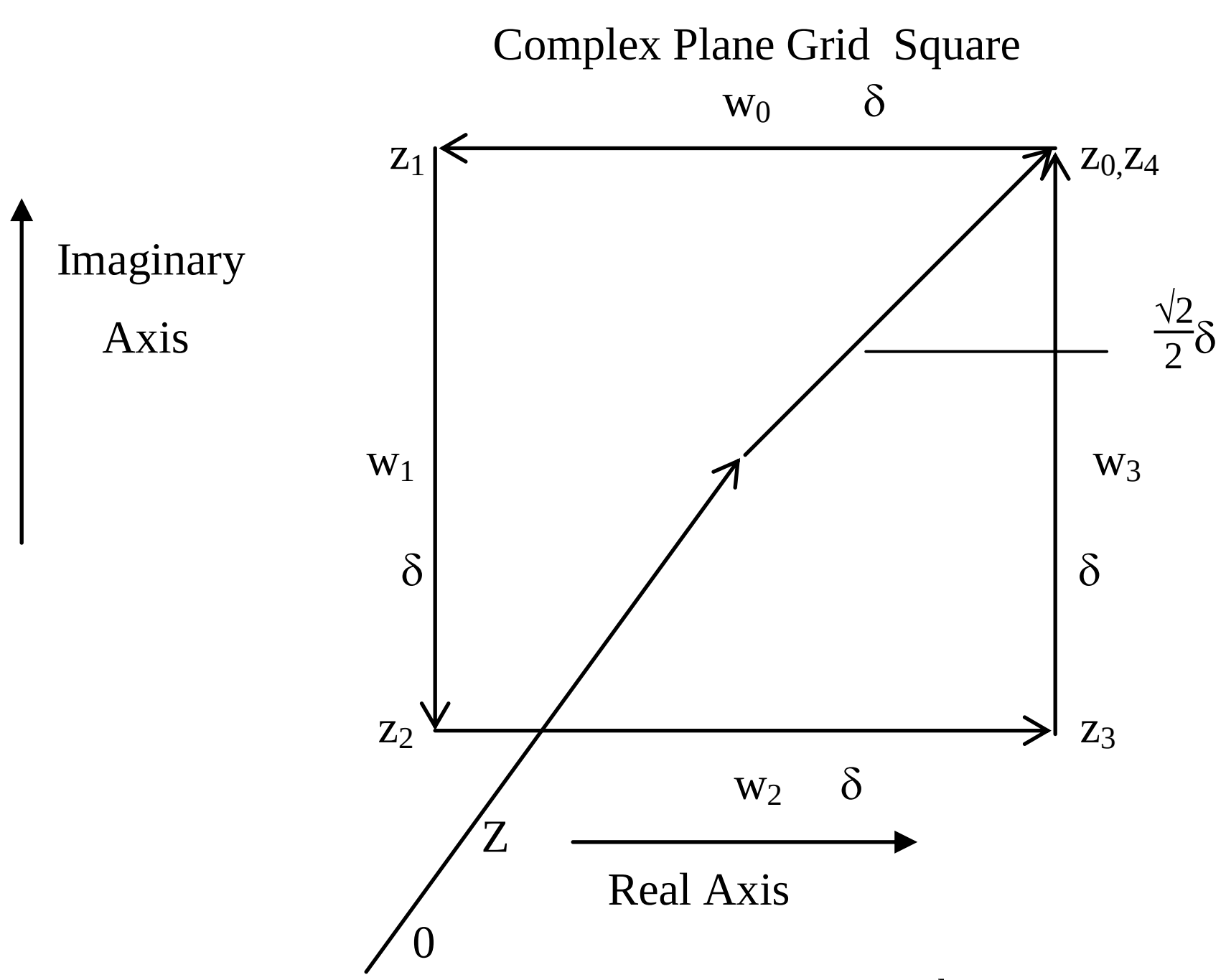

$$z_k = Z + \left(\frac{\sqrt{2}}{2}\right)\delta e^{\,j\left(\frac{\pi}{4} + \frac{k\pi}{2}\right)}$$

$w_k = z_k + (z_{k+1} - z_k)p$

$p$ = any real or complex value

$0 \le p \le 1$ for a point on a grid square vector
$k = 0,1,2,3$

$$\alpha_r = {}_1\oint_0^4 (Z + Ce^{\frac{jk\pi}{2}})^r ( Z + Be^{\frac{j[k+1]\pi}{2}} - Z - Be^{\frac{jk\pi}{2}} )\Delta k \qquad (3.2\text{-}16)$$

Simplifying Eq 3.2-16

$$\alpha_r = B_1\oint_0^4 (Z + Ce^{\frac{jk\pi}{2}})^r e^{\frac{jk\pi}{2}} ( e^{\frac{j\pi}{2}} - 1) \Delta k$$

$$\alpha_r = B( e^{\frac{j\pi}{2}} - 1) \, {}_1\oint_0^4 (Z + Ce^{\frac{jk\pi}{2}})^r e^{\frac{jk\pi}{2}} \Delta k \qquad (3.2\text{-}17)$$

Evaluate the quantity, $B( e^{\frac{j\pi}{2}} - 1)$

From Eq 3.2-8 and Eq 3.2-11

$$B( e^{\frac{j\pi}{2}} - 1) = (\tfrac{\sqrt{2}}{2})\delta e^{\frac{j\pi}{4}} ( e^{\frac{j\pi}{2}} - 1) = (\tfrac{\sqrt{2}}{2})\delta e^{\frac{j\pi}{4}} \sqrt{2}e^{\frac{j3\pi}{4}} = \delta e^{j\pi} = -\delta$$

$$B( e^{\frac{j\pi}{2}} - 1) = -\delta \qquad (3.2\text{-}18)$$

Substituting Eq 3.2-18 into Eq 3.2-17

$$\alpha_r = -\delta \, {}_1\oint_0^4 (Z + Ce^{\frac{jk\pi}{2}})^r e^{\frac{jk\pi}{2}} \Delta k \qquad (3.2\text{-}19)$$

Changing the form of Eq 3.2-19

$$\alpha_r = -\delta \, {}_1\oint_0^4 (1 + \tfrac{C}{Z}e^{\frac{jk\pi}{2}})^r Z^r e^{\frac{jk\pi}{2}} \Delta k \qquad (3.2\text{-}20)$$

Expand the term, $(1 + \frac{C}{Z}e^{\frac{jk\pi}{2}})^r$, in Eq 3.2-20 using the Binomial Series Expansion

The Binomial Expansion

$$(1+q)^r = 1 + \frac{r}{1!}q + \frac{r(r-1)}{2!}q^2 + \frac{r(r-1)(r-2)}{3!}q^3 + \ldots + \frac{r(r-1)(r-2)\ldots(r-n+1)}{n!}q^n + \ldots \qquad (3.2\text{-}21)$$

For r = a positive integer, 0,1,2,3,… , the Binomial Series Expansion is exact and has r +1 terms.

$$(1+\frac{C}{Z}e^{\frac{jk\pi}{2}})^r = 1 + \frac{r}{1!}\frac{C}{Z}e^{\frac{jk\pi}{2}} + \frac{r(r-1)}{2!}\frac{C^2}{Z^2}e^{jk\pi} + \frac{r(r-1)(r-2)}{3!}\frac{C^3}{Z^3}e^{\frac{j3k\pi}{2}} + \ldots + \frac{C^{r-1}}{1!Z^{r-1}}e^{\frac{j(r-1)k\pi}{2}} + \frac{C^r}{Z^r}e^{\frac{jrk\pi}{2}} \qquad (3.2\text{-}22)$$

Substituting Eq 3.2-22 into Eq 3.2-20

$$\alpha_r = -\delta\, {}_1\oint_0^4 [1 + \frac{r}{1!}\frac{C}{Z}e^{\frac{jk\pi}{2}} + \frac{r(r-1)}{2!}\frac{C^2}{Z^2}e^{jk\pi} + \frac{r(r-1)(r-2)}{3!}\frac{C^3}{Z^3}e^{\frac{j3k\pi}{2}} + \ldots$$

$$+ \frac{C^{r-1}}{1!Z^{r-1}}e^{\frac{j(r-1)k\pi}{2}} + \frac{C^r}{Z^r}e^{\frac{jrk\pi}{2}}]\, Z^r e^{\frac{jk\pi}{2}}\,\Delta k \qquad (3.2\text{-}23)$$

Simplifying Eq 3.2-23

$$\alpha_r = -\delta\, {}_1\oint_0^4 [Z^r e^{\frac{jk\pi}{2}} + Z^{r-1}\frac{r}{1!}C e^{jk\pi} + Z^{r-2}\frac{r(r-1)}{2!}C^2 e^{\frac{j3k\pi}{2}} + Z^{r-3}\frac{r(r-1)(r-2)}{3!}C^3 e^{j2k\pi} + \ldots$$

$$+ Z\frac{r}{1!}C^{r-1}e^{\frac{jrk\pi}{2}} + C^r e^{\frac{j(r+1)k\pi}{2}}]\Delta k \qquad (3.2\text{-}24)$$

Integrating using the following two relationships

$$ {}_1\int_L^H e^{ak}\Delta k = {}_1\sum_{k=L}^{H-1} e^{ak}\Delta k \qquad (3.2\text{-}25)$$

or

$$ {}_1\int_L^H e^{ak}\Delta k = \frac{1}{e^a-1}e^{ak}\Big|_{1\,L}^{H} = \frac{1}{e^a-1}[e^{aH} - e^{aL}], \qquad e^a-1 \neq 0 \qquad (3.2\text{-}26)$$

$$\alpha_r = -\delta[\frac{Z^r}{e^{\frac{j\pi}{2}}-1}e^{\frac{jk\pi}{2}}\Big|_0^4 + \frac{r}{1!}C\frac{Z^{r-1}}{e^{j\pi}-1}e^{jk\pi}\Big|_0^4 + \frac{r(r-1)}{2!}C^2\frac{Z^{r-2}}{e^{\frac{j3\pi}{2}}-1}e^{\frac{j3k\pi}{2}}\Big|_0^4 +$$

$$\frac{r(r-1)(r-2)}{3!} C^3 Z^{r-3} {}_1\sum_{k=0}^{3} e^{j2\pi k} + \frac{r(r-1)(r-2)(r-3)}{4!} C^4 \frac{Z^{r-4}}{e^{\frac{j5\pi}{2}}-1} e^{\frac{j5k\pi}{2}} \Big|_0^4 +$$

$$\frac{r(r-1)(r-2)(r-3)(r-4)}{5!} C^5 \frac{Z^{r-5}}{e^{j3\pi}-1} e^{j3k\pi} \Big|_0^4 + \frac{r(r-1)(r-2)(r-3)(r-4)(r-5)}{6!} C^6 \frac{Z^{r-6}}{e^{\frac{j7\pi}{2}}-1} e^{\frac{j7k\pi}{2}} \Big|_0^4 +$$

$$\frac{r(r-1)(r-2)(r-3)(r-4)(r-5)(r-6)}{7!} C^7 Z^{r-7} {}_1\sum_{k=0}^{3} e^{j4\pi k} + \ldots ] \qquad (3.2\text{-}27)$$

Then

From Eq 3.2-14

$$\alpha_r = {}_1\oint_0^4 w_k^r (z_{k+1} - z_k)\Delta k$$

$\Delta k = 1$
$k = 0,1,2,3$
$r = 0,1,2,3,\ldots$

Note the following diagram and Eq 3.2-27

Diagram 3.2-5 A complex plane grid square

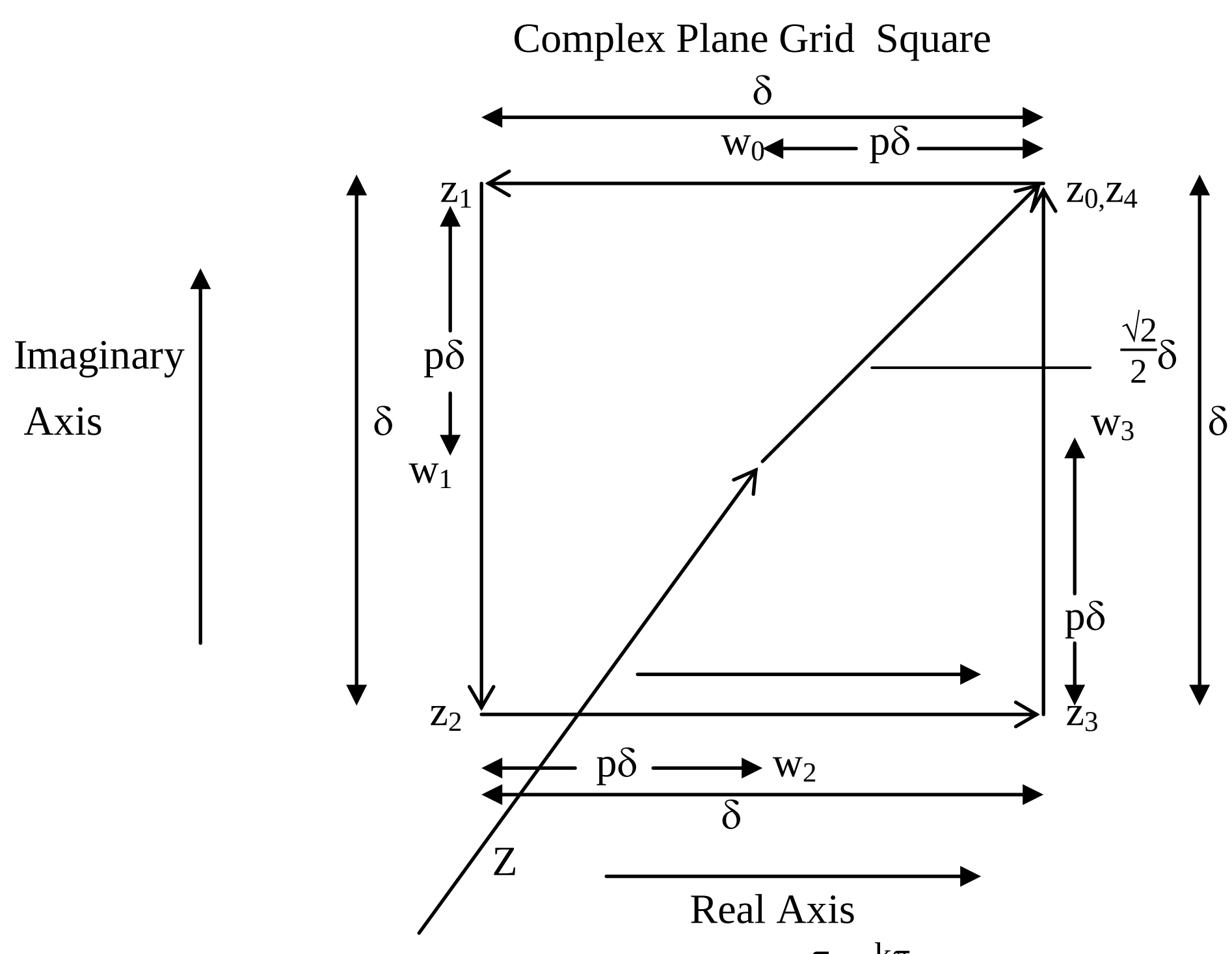

$$z_k = Z + \left(\frac{\sqrt{2}}{2}\right)\delta e^{j\left(\frac{\pi}{4} + \frac{k\pi}{2}\right)}$$

$w_k = z_k + p(z_{k+1} - z_k)$
p = any real or complex value
$0 \le p \le 1$ for a point on a grid square vector
Z = vector from the origin to the center of a grid square

$$ {}_rC_m = \frac{\prod_{s=0}^{m-1}(r-s)}{m!} \qquad (3.2\text{-}28) $$

From Eq 3.2-27

For r = 0,1,2

$$ \alpha_r = 0 \qquad (3.2\text{-}29) $$

For r = 3,4,5,6

$$ \alpha_r = -4\,\delta\,{}_rC_3\,C^3\,Z^{r-3} \qquad (3.2\text{-}30) $$

For r = 7,8,9,10

$$ \alpha_r = -4\,\delta\,(\,{}_rC_3\,C^3\,Z^{r-3} + {}_rC_7\,C^7\,Z^{r-7}) \qquad (3.2\text{-}31) $$

For r = 11,12,13,14

$$ \alpha_r = -4\,\delta\,(\,{}_rC_3\,C^3\,Z^{r-3} + {}_rC_7\,C^7\,Z^{r-7} + {}_rC_{11}\,C^{11}\,Z^{r-11}) \qquad (3.2\text{-}32) $$

.
.
.

In General

For r = 0,1,2

$$ \alpha_r = 0 \qquad (3.2\text{-}33) $$

For n = 1,2,3,…

r = 4n-1, 4n, 4n+1, 4n+2 ,

$$\alpha_r = -4\,\delta \sum_{p=1}^{n} {}_rC_{4p-1} C^{4p-1} Z^{r-(4p-1)} \qquad (3.2\text{-}34)$$

Referring to Equation 3.2-34 something very important is noticed. For n = 1 and r = 3, $\alpha_3 = -4\,\delta\, {}_3C_3\, C^3 Z^{3-3} = -4\delta C^3$. $\alpha_3$ is seen to have a non-zero value where Z is not a factor in that value. This occurs only for n = 1 and r = 3. This unique occurrence has considerable mathematical importance. For this condition, the integral now specified around any grid square in the complex plane (See Diagram 3.2-1 and Diagram 3.2-5) has the same value irrespective of its location in the plane.

$$\alpha_3 = {}_1\oint_0^4 w_k^{\,3}(z_{k+1} - z_k)\Delta k = {}_1\oint_0^4 (z_k + p[z_{k+1} - z_k])^3 (z_{k+1} - z_k)\Delta k = -4\delta C^3 \qquad (3.2\text{-}35)$$

Evaluating the complex constant C

From Eq 3.2-10

$$C = \delta\left(\frac{\sqrt{2}}{2} + pe^{\frac{j3\pi}{4}}\right) e^{\frac{j\pi}{4}}$$

$$e^{\frac{j\pi}{4}} = \cos\frac{\pi}{4} + j\sin\frac{\pi}{4} = \frac{1}{\sqrt{2}}(1+j) \qquad (3.2\text{-}36)$$

From Eq 3.2-8, Eq 3.2-10 and Eq 3.2-36

$$C = \delta\left(\frac{\sqrt{2}}{2} + p\frac{j-1}{\sqrt{2}}\right)\left(\frac{1}{\sqrt{2}}\right)(j+1)$$

Simplying

$$C = \delta\left(\frac{1}{2}j + \frac{1}{2} - p\right) = \frac{\delta}{2}(1 - 2p + j) \qquad (3.2\text{-}37)$$

Then

Substituting Eq 3.2-37 into Eq 3.2-35

$${}_1\oint_0^4 (z_k + p[z_{k+1} - z_k])^3 (z_{k+1} - z_k)\Delta k = -\frac{\delta^4}{2}(1-2p+j)^3 \;, \quad \Delta k = 1 \qquad (3.2\text{-}38)$$

p = any real or complex value

Using a derivation similar to that used to obtain Eq 3.2-38, the clockwise integral around a complex plane grid square can be obtained. It is as follows:

$$_1\oint_0^4 (z_k + p[z_{k+1} - z_k])^3 (z_{k+1} - z_k)\Delta k = +\frac{\delta^4}{2}(1-2p+j)^3 \;, \quad \Delta k=1 \qquad (3.2\text{-}39)$$

p = any real or complex value

Changing the direction of integration changes the sign of the result.

In discrete Interval Calculus, integrations can be written equivalently as summations as shown below. The contour summation and contour integral yield equivalent results. However, the summation and the integral calculation operations differ. Contour summation adds each evaluation of the summation function to yield the final result. Contour Integration, as in Calculus, mathematically manipulates the integrated function and then evaluates the resulting function to yield the final result.

From Eq 3.2-38
The counterclockwise integral/summation around any grid square is:

$$_1\oint_0^4 ( z_k + p[z_{k+1} - z_k] )^3 (z_{k+1} - z_k)\Delta k = {}_1\oint\!\!\sum_{k=0}^{3} ( z_k + p[z_{k+1} - z_k] )^3 (z_{k+1} - z_k) = -\frac{\delta^4}{2}(1-2p+j)^3 \qquad (3.2\text{-}40)$$

p = any real or complex value

From Eq 3.2-39
The clockwise integral/summation around any grid square is:

$$_1\oint_0^4 ( z_k + p[z_{k+1} - z_k] )^3 (z_{k+1} - z_k)\Delta k = {}_1\oint\!\!\sum_{k=0}^{3} ( z_k + p[z_{k+1} - z_k] )^3 (z_{k+1} - z_k) = +\frac{\delta^4}{2}(1-2p+j)^3 \qquad (3.2\text{-}41)$$

p = any real or complex value

An appropriate selection for the value of p in Eq 3.2-40 and Eq 3.2-41 simplifies and makes these summations exceedingly useful. The value chosen is $p = \frac{1}{2}$.

For $p = \frac{1}{2}$, the following summations result:

From Eq 3.2-40
The counterclockwise integration/summation around any grid square is:

$$ {}_1\oint_0^4 \left(\frac{z_{k+1} + z_k}{2}\right)^3 (z_{k+1} - z_k)\Delta k = {}_1\sum_{k=0}^{3} \left(\frac{z_{k+1} + z_k}{2}\right)^3 (z_{k+1} - z_k) = +\frac{j\delta^4}{2} \qquad (3.2\text{-}42)$$

From Eq 3.2-41
The clockwise summation around any grid square is:

$$ {}_1\oint_0^4 \left(\frac{z_{k+1} + z_k}{2}\right)^3 (z_{k+1} - z_k)\Delta k = {}_1\sum_{k=0}^{3} \left(\frac{z_{k+1} + z_k}{2}\right)^3 (z_{k+1} - z_k) = -\frac{j\delta^4}{2} \qquad (3.2\text{-}43)$$

Note Diagram 3.2-6 below. Grid square sides not on the circular perimeter contour (specified in black) are seen to be formed by two equal and opposite vectors from the adjacent grid squares. The * mark shown in the diagram indicates these canceling vectors. Notice, also, from Eq 3.2-42 and Eq 3.2-43 that the evaluations associated with these canceling vectors are the negative of each other and also cancel. Both this property and the property whereby the summation has the same circular summation value around any grid square makes equations Eq 3.2-42 and Eq 3.2-43 exceedingly useful. In particular, these equations can be used to find the area within a closed contour in the complex plane.

<u>Diagram 3.2-6 A discrete closed contour in the complex plane</u>

An example of a discrete closed contour in the complex plane
(Counterclockwise vector direction)

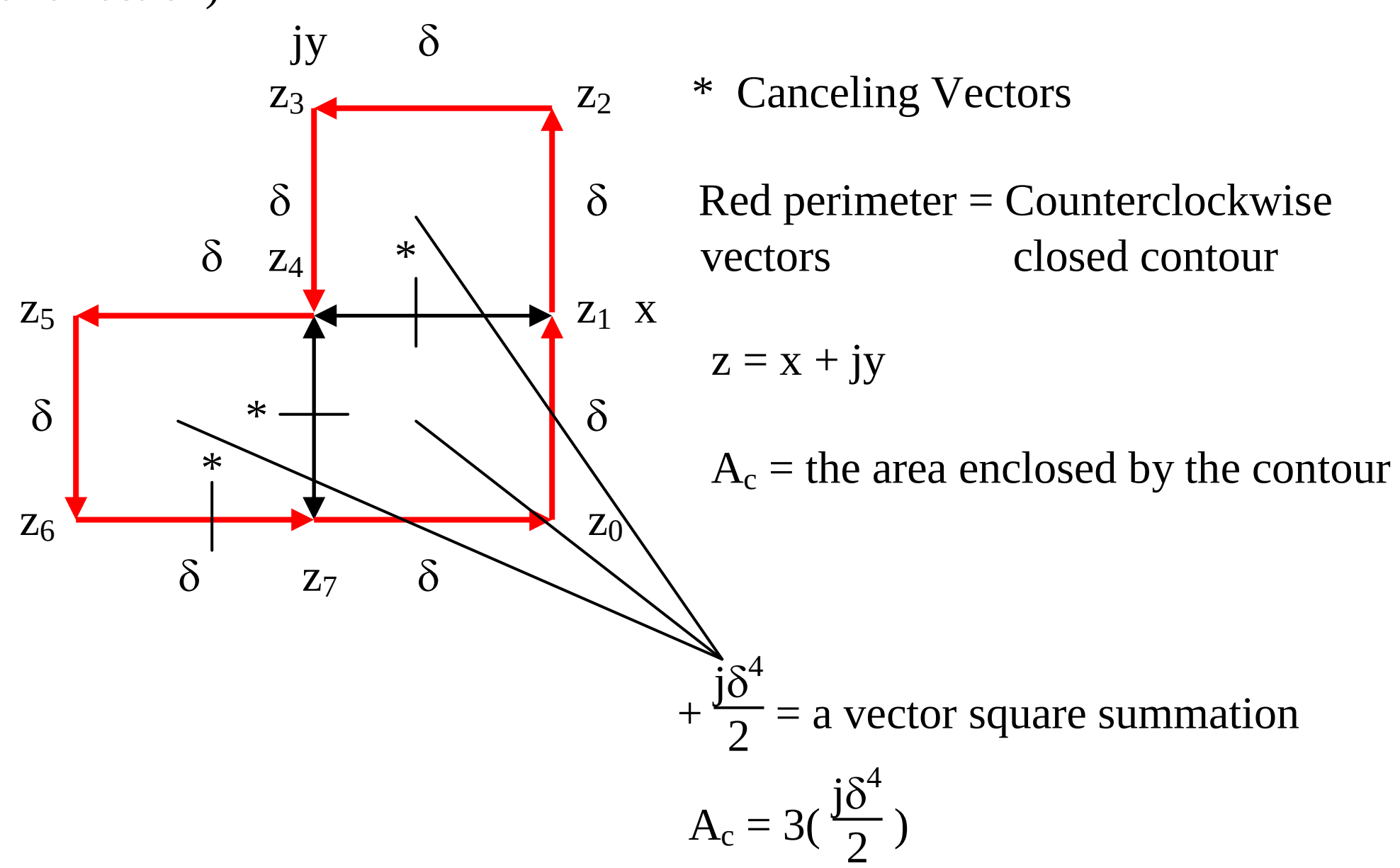

By totaling all grid square summations (Eq 3.2-42) of the grid squares contained within a counterclockwise complex plane contour, the value $+N\frac{j\delta^4}{2}$ will be obtained. N equals the number of grid squares enclosed within the complex plane contour. By totaling all grid square summations (Eq 3.2-43) of the grid squares contained within a clockwise complex plane contour, the value $-N\frac{j\delta^4}{2}$ will be obtained. Interestingly, only those evaluations associated with the grid square vectors forming a closed complex plane perimeter contour need to be performed and summed (closed contour integration/summation) to obtain the area within the contour. Those vectors located within a closed complex plane contour cancel to zero. Refer to Diagram 3.2-6.

Since the value obtained by totaling all the summations of the grid squares within a closed complex plane contour is $+N\frac{j\delta^4}{2}$ for counterclockwise summations and $-N\frac{j\delta^4}{2}$ for clockwise summations, by the argument above, these same results can be obtained by summing just the summation terms associated with the closed contour perimeter vectors (note Diagram 3.2-6).
The notation adopted to represent just the perimeter contour vector summation evaluations is as follows:

For counterclockwise complex plane closed contour summation and integration.

$$\text{Contour Summation} = {}_1\circlearrowleft\!\!\!\!\!\!\sum_{k=0}^{v-1} \left(\frac{z_{k+1}+z_k}{2}\right)^3 (z_{k+1}-z_k) \tag{3.2-44}$$

$$\text{Contour Integration} = {}_1\oint_0^{v}\!\!\circlearrowleft \left(\frac{z_{k+1}+z_k}{2}\right)^3 (z_{k+1}-z_k)\Delta k \tag{3.2-45}$$

where
v = the number of vectors composing the contour

For clockwise complex plane closed contour summation and integration. They are equivalent.

$$\text{Contour Summation} = {}_1\circlearrowright\!\!\!\!\!\!\sum_{k=0}^{v-1} \left(\frac{z_{k+1}+z_k}{2}\right)^3 (z_{k+1}-z_k) \tag{3.2-46}$$

$$\text{Contour Integration} = {}_1\oint_0^{v}\!\!\circlearrowright \left(\frac{z_{k+1}+z_k}{2}\right)^3 (z_{k+1}-z_k)\Delta k \tag{3.2-47}$$

where
v = the number of vectors composing the contour

From the equations, Eq 3.42 thru Eq 3.47, and the explanation of Diagram 3.9

The counterclockwise integration/summation around a complex plane discrete closed contour is:

$$\oint_{0}^{v}\left(\frac{z_{k+1}+z_k}{2}\right)^3(z_{k+1}-z_k)\Delta k = \sum_{k=0}^{v-1}\left(\frac{z_{k+1}+z_k}{2}\right)^3(z_{k+1}-z_k) = +N\left(\frac{j\delta^4}{2}\right) \quad (3.2\text{-}48)$$

The clockwise summation around around a complex plane discrete closed contour is:

$$\oint_{0}^{v}\left(\frac{z_{k+1}+z_k}{2}\right)^3(z_{k+1}-z_k)\Delta k = \sum_{k=0}^{v-1}\left(\frac{z_{k+1}+z_k}{2}\right)^3(z_{k+1}-z_k) = -N\left(\frac{j\delta^4}{2}\right) \quad (3.2\text{-}49)$$

Calculate the area of a complex plane grid square

$$A_g = \delta^2 \quad (3.2\text{-}50)$$

Calculate the area encircled by a discrete closed contour in the complex plane

$$A_c = N\delta^2 = NA_g \quad (3.2\text{-}51)$$

Find the number of grid squares encircled by a complex plane discrete closed contour
From Eq 3.2-48 and Eq 3.2-50

$$\mathbf{N = -\frac{2j}{A_g^2}\oint_{0}^{v}\left(\frac{z_{k+1}+z_k}{2}\right)^3(z_{k+1}-z_k)\Delta k = -\frac{2j}{A_g^2}\sum_{k=0}^{v-1}\left(\frac{z_{k+1}+z_k}{2}\right)^3(z_{k+1}-z_k)} \quad (3.2\text{-}52)$$

From Eq 3.2-49 and Eq 3.2-50

$$\mathbf{N = \frac{2j}{A_g^2}\oint_{0}^{v}\left(\frac{z_{k+1}+z_k}{2}\right)^3(z_{k+1}-z_k)\Delta k = \frac{2j}{A_g^2}\sum_{k=0}^{v-1}\left(\frac{z_{k+1}+z_k}{2}\right)^3(z_{k+1}-z_k)} \quad (3.2\text{-}53)$$

**<u>Note</u> -** $\mathbf{\frac{z_{k+1}+z_k}{2} = z_k + \frac{z_{k+1}-z_k}{2}}$ **, an alternate form which may be used for calculation**

**where**

**N = the number of grid squares encircled by a discrete closed contour in the complex plane**
**$A_g$ = the area of a complex plane grid square**
**v = the number of counterclockwise or clockwise pointing vectors comprising the complex plane discrete closed contour being integrated or summed**
**$z_k$ = the coordinate value of the kth point on the closed contour**
**$\Delta k = 1$**
**$z_{k+1} - z_k$ = a contour grid square vector value, $z_{k+1}$ at the vector head, $z_k$ at the tail**

Find the area enclosed in a closed contour in the complex plane

From Eq 3.2-48 and Eq 3.2-51

$$A_c = -\frac{2j}{A_g}\,\frac{1}{\,}\oint_0^v \left(\frac{z_{k+1}+z_k}{2}\right)^3 (z_{k+1}-z_k)\Delta k = -\frac{2j}{A_g}\,\frac{1}{\,}\sum_{k=0}^{v-1} \left(\frac{z_{k+1}+z_k}{2}\right)^3 (z_{k+1}-z_k) \qquad (3.2\text{-}54)$$

From Eq 3.2-49 and Eq 3.2-51

$$A_c = \frac{2j}{A_g}\,\frac{1}{\,}\oint_0^v \left(\frac{z_{k+1}+z_k}{2}\right)^3 (z_{k+1}-z_k)\Delta k = \frac{2j}{A_g}\,\frac{1}{\,}\sum_{k=0}^{v-1} \left(\frac{z_{k+1}+z_k}{2}\right)^3 (z_{k+1}-z_k) \qquad (3.2\text{-}55)$$

**Note** - $\frac{z_{k+1}+z_k}{2} = z_k + \frac{z_{k+1}-z_k}{2}$ , **an alternate form which may be used for calculation**

From Eq 3.2-51

$$N = \frac{A_c}{A_g} \qquad (3.2\text{-}56)$$

**where**

**N = the number of grid squares enclosed in a discrete closed contour in the complex plane**

**$A_c$ = the area enclosed in a discrete closed contour in the complex plane**

**$A_g$ = the area of a complex plane grid square**

**v = the number of counterclockwise or clockwise pointing vectors comprising the complex plane discrete closed contour being integrated or summed**

**$z_k$ = the coordinate value of the kth point on the closed contour**

**$\Delta k = 1$**

**$z_{k+1} - z_k$ = a contour grid square vector value, $z_{k+1}$ at the vector head, $z_k$ at the tail**

**Section 3.3: Area Calculation Using Discrete Closed Contour Summation in the Complex Plane**

In Calculus, mathematical manipulation of functions is the means employed to exactly evaluate areas. This is the process known as integration. To exactly evaluate areas, the summation process is not used. The summation of an infinite number of function calculations multiplied by infinitesimal intervals is not a practical process for finding area. Thus, in Calculus, the summation process is used only to obtain area approximations. However, in discrete Interval Calculus, exact area evaluations may be obtained in both ways, by integration and by summation. Complex plane closed contour summation is discussed below.

In Interval Calculus, the summation area equations of Eq 3.2-54 and Eq 3.2-55 derived in Section 3.2 can be directly applied to the calculation of the area enclosed within a discrete closed contour in the complex plane. Discrete closed contours composed of counterclockwise vectors have their enclosed areas calculated by the summation equation of Eq 3.2-54. Discrete closed contours composed of clockwise vectors have their enclosed areas calculated by the summation equation of Eq 3.2-55. The use of these equations is best demonstated by an example. See Example 3.3-1 below.

Example 3.3-1

Find the area, $A_c$, enclosed within the complex plane discrete closed contour presented in Diagram 3.3-1 below. Use the following complex plane clockwise vector closed contour summation equation to calculate $A_c$.

$$A_c = \frac{2j}{A_g} \, {}_1\!\!\sum_{k=0}^{v-1} \left(\frac{z_{k+1} + z_k}{2}\right)^3 (z_{k+1} - z_k) \qquad (3.3\text{-}1)$$

where

$A_c$ = the area enclosed in a discrete closed contour in the complex plane

$A_g$ = the area of a complex plane grid square

$v$ = the number of clockwise pointing vectors comprising the complex plane discrete closed contour

$k = 0, 1, 2, 3,\ldots, v-1$

$z_k$ = the coordinate value of the kth point on the closed contour

$\Delta k = 1$

$z_{k+1} - z_k$ = a contour grid square vector value, $z_{k+1}$ at the vector head, $z_k$ at the tail

Diagram 3.3-1  An example of a complex plane discrete closed contour

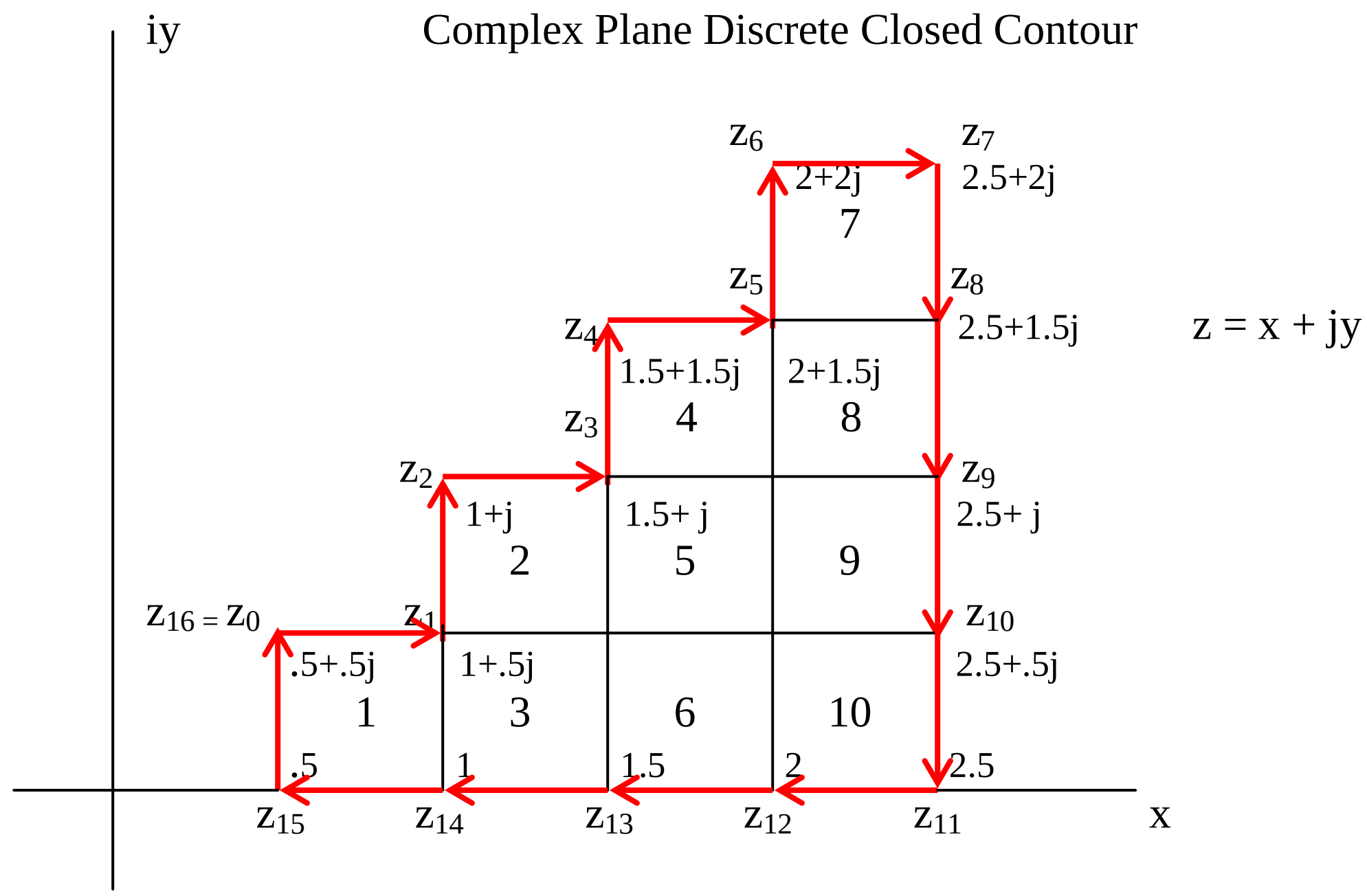

To calculate the contour enclosed area, $A_{c,}$ the following UBASIC program was written.

```
10   'This program performs a closed contour clockwise summation
20   'to find the area within a discrete complex plane closed contour
30   AG=0.25 'The area of a grid square
40   M=15 'One less than the number of vectors forming the contour
50   dim Z(M+1)
60   'A list of the points defining a closed complex plane contour
70   Z(0)=0.5+0.5#i
80   Z(1)=1+0.5#i
90   Z(2)=1+#i
100  Z(3)=1.5+#i
110  Z(4)=1.5+1.5#i
120  Z(5)=2+1.5#i
130  Z(6)=2+2#i
140  Z(7)=2.5+2#i
150  Z(8)=2.5+1.5#i
160  Z(9)=2.5+#i
170  Z(10)=2.5+0.5#i
180  Z(11)=2.5
190  Z(12)=2
200  Z(13)=1.5
210  Z(14)=1
220  Z(15)=0.5
230  Z(16)=Z(0)
240  S=0
```

```
250   for K=0 to M
260   S=S+(((Z(K+1)+Z(K))/2)^3)*(Z(K+1)-Z(K))
270   next K
280   AC=(2*#i/AG)*S
290   print "The area within the closed contour is ";AC
```

The evaluation obtained by running the above program is as follows:

**The area within the closed contour is 2.5**

Checking the above result

Referring to the contour diagram, Diagram 3.3-1

$$A_c = 10 \text{ grid squares} * \frac{.25}{\text{grid square}} = 2.5$$

Good check

As can be seen, there is a calculation required for each vector of the contour. The resulting calculated values are summed to obtain the contour enclosed area, $A_{c.}$ This area is exact. This same area can also be obtained using discrete closed contour integration. The complex plane discrete closed contour integration process for finding area will be discussed in the following section, Section 3.4.

## Section 3.4: Area Calculation Using Discrete Closed Contour Integration in the Complex Plane

In Interval Calculus one has a choice of processes for calculating the area enclosed within a closed discrete contour in the complex plane. The summation process described in Section 3.3 above is direct and straight forward but its computational requirements can be excessive. On the other hand, the discrete integration process is often less direct and straight forward but its computational requirements can be substantially less. The nature of the area to be calculated determines the most expedient calculation methodology. As previously mentioned, integration requires the mathematical manipulation of functions to evaluate areas. The integration process thus requires that the area enclosed within a complex plane closed contour be defined by an integrable mathematical function that describes the closed contour. A single function which fully describes a closed complex plane contour and is conveniently integrable often does not exist. As a result, it is usual that a closed contour be described by several different functions. Below are derived two discrete closed contour integral equations. These equations can be used together with different functions to integrate along a complex plane discrete closed contour to evaluate the area within.

From Eq 3.2-54 and Eq 3.2-55

$$A_c = \frac{2j}{A_g} \; {}_1\oint_0^v \left(\frac{z_{k+1} + z_k}{2}\right)^3 (z_{k+1} - z_k)\Delta k \; = -\frac{2j}{A_g} \; {}_1\oint_0^v \left(\frac{z_{k+1} + z_k}{2}\right)^3 (z_{k+1} - z_k)\Delta k \qquad (3.4\text{-}1)$$

Putting Eq 3.4-1 into a form more convenient for integration

$$\frac{z_{k+1} + z_k}{2} = z_k + \frac{z_{k+1} - z_k}{2} \qquad (3.4\text{-}2)$$

$$\Delta z = z_{k+1} - z_k \qquad (3.4\text{-}3)$$

Substituting Eq 3.4-2 and Eq 3.4-3 into Eq 3.4-1

The complex plane contour area calculation equations are:

$$\mathbf{A_c = \frac{2j}{A_g} \; {}_1\oint_0^v \left[\left( z_k + \frac{\Delta z}{2}\right)^3 \Delta z\right] \Delta k} \qquad (3.4\text{-}4)$$

$$\mathbf{A_c = -\frac{2j}{A_g} \; {}_1\oint_0^v \left[\left( z_k + \frac{\Delta z}{2}\right)^3 \Delta z\right] \Delta k} \qquad (3.4\text{-}5)$$

**where**

**$A_c$ = the area enclosed by the vectors forming a discrete closed contour in the complex plane**

**v = the number of horizontal and vertical vectors forming the discrete closed contour**

**$z_k$ = f(k) , function(s) of k**

**$z_k$ = the coordinate value of the contour kth vector tail point**

**$\Delta z$ = real or imaginary vector values**

$$\Delta z = \begin{cases} +\delta \text{ for a horizontal left to right pointing vector} \\ -\delta \text{ for a horizontal right to left pointing vector} \\ +j\delta \text{ for a vertical down to up pointing vector} \\ -j\delta \text{ for a vertical up to down pointing vector} \end{cases}$$

**$z_k + \Delta z$ = the coordinates of the contour kth vector head point**

**$A_g$ = $\delta x \delta$ , the area of a complex plane grid**

**$\delta$ = the length of the side of a complex plane grid square**

**$\Delta k = 1$**

Some additional explanation concerning functions and contours is as follows:

1) The function variable, k, identifies the tail points of all vectors forming the contour.
   a) In discrete Interval Calculus, points can only occur at a complex plane grid square corner.
   b) All vectors have the same length, the length of a side of a grid square, δ.
   c) Vectors are placed head to tail to form a complex plane contour.

2) The function variable, k, is an integer.
3) $k = 0,1,2,3,\ldots,v$ , v is the number of vectors forming a complex plane contour.
4) There is a unique function variable integer value to identify each vector.
5) A function specifies the complex plane coordinates of the point designated by its respective variable value.
7) Discrete contours in the complex plane are composed of equal length vectors placed only in horizontal or vertical positions collinear with a side of a grid square.
8) One or more functions are used to describe all the points on a complex plane discrete closed contour.
9) Functions are selected not only to identify the tail coordinates of contour vectors but also to be easily integrable. Linear functions are often used.

Closed contours composed of clockwise vectors have their enclosed areas calculated by the integral equation of Eq 3.4-4. Closed contours composed of counterclockwise vectors have their enclosed areas calculated by the integral equation of Eq 3.4-5. The use of these equations is best demonstated by an example. See Example 3.4-1 below.

Example 3.4-1

Find the area, $A_{c,}$ enclosed within the complex plane closed contour presented below. See Diagram 3.4-1 below. Use the following complex plane closed contour integration equation to calculate $A_c$.

$$A_c = \frac{2j}{A_g} \oint_{0}^{v} {}_1 \left[\left( z_k + \frac{\Delta z}{2}\right)^3 \Delta z\right] \Delta k \qquad (3.4\text{-}6)$$

where

$A_c$ = the area enclosed by the vectors forming a discrete closed contour in the complex plane
$v$ = the number of horizontal and vertical vectors forming the discrete closed contour
$z_k$ = $f(k)$ , function(s) of k
$z_k$ = the coordinates of the contour kth vector tail point
$\Delta z$ = real or imaginary value

$$\Delta z = \begin{cases} +\delta \text{ for a horizontal left to right pointing vector} \\ -\delta \text{ for a horizontal right to left pointing vector} \\ +j\delta \text{ for a vertical down to up pointing vector} \\ -j\delta \text{ for a vertical up to down pointing vector} \end{cases}$$

$z_k + \Delta z$ = the coordinates of the contour kth vector head point
$A_g$ = $\delta x \delta$ , the area of a complex plane grid square
$\delta$ = the length of the side of a complex plane grid square
$\Delta k = 1$

Note that the complex plane discrete closed contour of this example, Diagram 3.4-1, is the same as the complex plane closed contour of Example 3.3-1, Diagram 3.3-1. Note, also, that the notation relating to the closed contour vectors has been changed. This change has been made to facilitate the integration process. The contour vectors have been divided into four groups of four vectors, A, B, C, and D. A different linear equation has been selected to represent the vector tail points of the vectors in each group.

Diagram 3.4-1 A Complex plane discrete closed contour

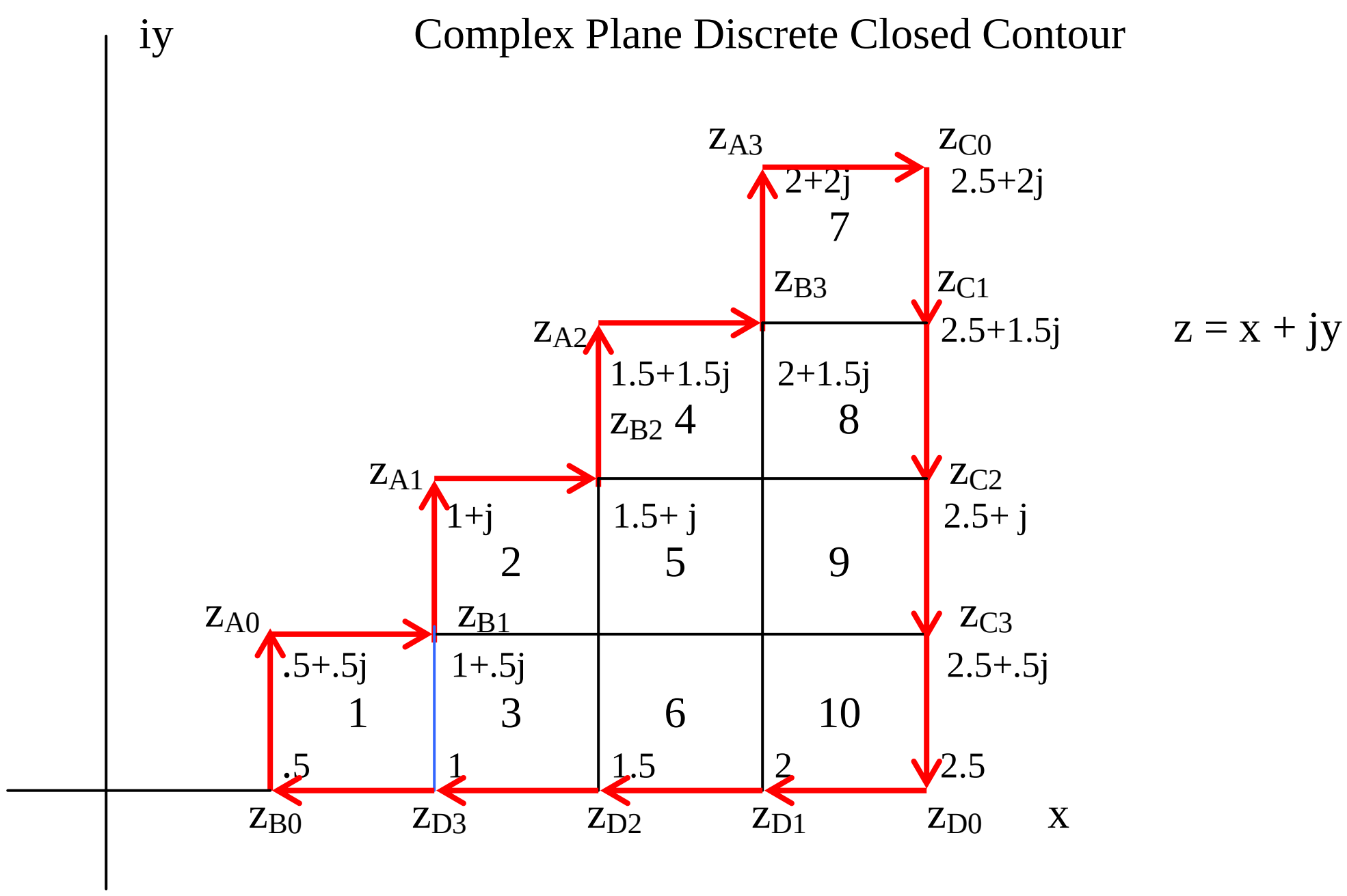

The closed contour vectors are grouped into four groups of four vectors

Referring to Diagram 3.4-1 above, calculate the area within the above closed contour using Eq 3.4-6 and four functions collectively representing all closed contour vectors.

The area, $A_c$, is found by adding the four vector group integrals

$A_c = G_A + G_B + G_C + G_D$

Vector Group A

$z_k = .5+.5j +k(.5+.5j)$ , $k = 0,1,2,3$ (3.4-7)

$\Delta z = .5$ (3.4-8)

$u = 4$ , the number of vectors in Group A (3.4-9)

$A_g = .5x.5 = .25$ , the area of a grid square (3.4-10)

From Eq 3.4-6

$$G_A = \frac{2j}{A_g}\,\frac{1}{\,}\int_0^u [( z_k + \frac{\Delta z}{2})^3 \Delta z]\, \Delta k$$ , the Group A contour segment integral equation (3.4-11)

Substituting Eq 3.4-7 thru Eq 3.4-10 into Eq 3.4-11

$$G_A = \frac{2j}{.25}\ {}_1\int_0^4 [( .5+.5j +k(.5+.5j) + \frac{.5}{2} )^3 .5]\ \Delta k = 4j\ {}_1\int_0^4 [.75+.5j +k(.5+.5j)]^3\ \Delta k$$

$$G_A = \frac{j(1+j)^3}{2}\ {}_1\int_0^4 [\frac{.75+.5j}{.5+.5j} + k]^3 \Delta k$$

$$G_A = (-1-j)\ {}_1\int_0^4 [k + 1.25-.25j]^3 \Delta k \qquad (3.4\text{-}12)$$

Changing variables

Let $p = k + (1.25-.25j)$

$\Delta p = \Delta k = 1$

for $k = 0$ $p = 1.25-.25j$

for $k = 4$ $p = 5.25-.25j$

Substituting into Eq 3.4-12

$$G_A = (-1-j)\ {}_1\int_{1.25-.25j}^{5.25-.25j} p^3 \Delta p \qquad (3.4\text{-}13)$$

From the Integral equations in Table 6 in the Appendix

$${}_{\Delta x}\int x^3 \Delta x = \frac{x(x-\Delta x)(x-2\Delta x)(x-3\Delta x)}{4} + x(x-\Delta x)(x-2\Delta x)\Delta x + \frac{x(x-\Delta x)}{2}\Delta x^2 + k \qquad (3.4\text{-}14)$$

Using Eq 3.4-14 to integrate Eq 3.4-13

$x = p$

$\Delta x = \Delta p = 1$

$$G_A = (-1-j)\ {}_1\int_{1.25-.25j}^{5.25-.25j} p^3 \Delta p = (-1-j)[\frac{p(p-1)(p-2)(p-3)}{4}\Big|_{1.25-.25j}^{5.25-.25j} + p(p-1)(p-2)\Big|_{1.25-.25j}^{5.25-.25j} + \frac{p(p-1)}{2}\Big|_{1.25-.25j}^{5.25-.25j}] \qquad (3.4\text{-}15)$$

$$G_A = -148.75-96j \qquad (3.4\text{-}16)$$

Vector Group B

$z_k = .5 + k(.5+.5j)$ ,  $k = 0,1,2,3$ (3.4-17)

$\Delta z = .5j$ (3.4-18)

$u = 4$ ,  the number of vectors in Group B (3.4-19)

$A_g = .5x.5 = .25$ ,  the area of a grid square (3.4-20)

From Eq 3.4-6

$$G_B = \frac{2j}{A_g}\ {}_1\!\int_0^u [(z_k + \frac{\Delta z}{2})^3 \Delta z]\ \Delta k$$ ,  the Group B contour segment integral equation (3.4-21)

Substituting Eq 3.4-17 thru Eq 3.4-20 into Eq 3.4-21

$$G_B = \frac{2j}{.25}\ {}_1\!\int_0^4 [( .5+k(.5+.5j) + \frac{.5j}{2})^3 .5j]\ \Delta k = -4\ {}_1\!\int_0^4 [.5+.25j + k(.5+.5j)]^3\ \Delta k$$

$$G_B = \frac{-(1+j)^3}{2}\ {}_1\!\int_0^4 [\frac{.5+.25j}{.5+.5j} + k]^3 \Delta k$$

$$G_B = (1-j)\ {}_1\!\int_0^4 [k + .75-.25j]^3 \Delta k$$ (3.4-22)

Changing variables

Let  $p = k + (.75-.25j)$

$\Delta p = \Delta k = 1$

for $k = 0$  $p = .75-.25j$

for $k = 4$  $p = 4.75-.25j$

Substituting into Eq 3.4-22

$$G_B = (-1-j)\ {}_1\!\int_{.75-.25j}^{4.75-.25j} p^3 \Delta p$$ (3.4-23)

From the Integral equations in Table 6 in the Appendix

$${}_{\Delta x}\!\int x^3 \Delta x = \frac{x(x-\Delta x)(x-2\Delta x)(x-3\Delta x)}{4} + x(x-\Delta x)(x-2\Delta x)\Delta x + \frac{x(x-\Delta x)}{2}\Delta x^2 + k$$ (3.4-24)

Using Eq 3.4-24 to integrate Eq 3.4-23

$x = p$

$\Delta x = \Delta p = 1$

$$G_B = (1-j)\int_{1\ \ .75-.25j}^{4.75-.25j} p^3 \Delta p = (1-j)\left[\frac{p(p-1)(p-2)(p-3)}{4}\Big|_{.75-.25j}^{4.75-.25j} + p(p-1)(p-2)\Big|_{.75-.25j}^{4.75-.25j} + \frac{p(p-1)}{2}\Big|_{.75-.25j}^{4.75-.25j}\right] \tag{3.4-25}$$

$$G_B = 58.75-96.5j \tag{3.4-26}$$

<u>Vector Group C</u>

$$z_k = 2.5+2j+k(-.5j)\ , \quad k = 0,1,2,3 \tag{3.4-27}$$

$$\Delta z = -.5j \tag{3.4-28}$$

$$u = 4\ , \quad \text{the number of vectors in Group C} \tag{3.4-29}$$

$$A_g = .5x.5 = .25\ , \quad \text{the area of a grid square} \tag{3.4-30}$$

From Eq 3.4-6

$$G_C = \frac{2j}{A_g}\ {}_1\!\int_0^u \left[\left(z_k + \frac{\Delta z}{2}\right)^3 \Delta z\right] \Delta k\ , \quad \text{the Group C contour segment integral equation} \tag{3.4-31}$$

Substituting Eq 3.4-27 thru Eq 3.4-30 into Eq 3.4-31

$$G_C = \frac{2j}{.25}\ {}_1\!\int_0^4 \left[\left(2.5+2j+k(-.5j) + \frac{-.5j}{2}\right)^3(-.5j)\right] \Delta k = 4\ {}_1\!\int_0^4 [2.5+1.75j+k(-.5j)]^3\ \Delta k$$

$$G_C = -\frac{j^3}{2}\ {}_1\!\int_0^4 \left[\frac{2.5+1.75j}{-.5j} + k\right]^3 \Delta k$$

$$G_C = \frac{j}{2}\ {}_1\!\int_0^4 [k + (-3.5+5j)]^3 \Delta k \tag{3.4-32}$$

Changing variables

Let $p = k - 3.5+5j$

$\Delta p = \Delta k = 1$

for $k = 0$ $p = -3.5+5j$

for $k = 4$ $p = .5+5j$

Substituting into Eq 3.4-32

$$G_C = \frac{j}{2} \; {}_1\!\int_{-3.5+5j}^{.5+5j} p^3 \Delta p \qquad (3.4\text{-}33)$$

From the Integral equations in Table 6 in the Appendix

$$_{\Delta x}\!\int x^3 \Delta x = \frac{x(x-\Delta x)(x-2\Delta x)(x-3\Delta x)}{4} + x(x-\Delta x)(x-2\Delta x)\Delta x + \frac{x(x-\Delta x)}{2}\Delta x^2 + k \qquad (3.4\text{-}34)$$

Using Eq 3.4-34 to integrate Eq 3.4-33

$x = p$

$\Delta x = \Delta p = 1$

$$G_C = \frac{j}{2} \; {}_1\!\int_{-3.5+5j}^{.5+5j} p^3 \Delta p = -\frac{j}{2}\left[\frac{p(p-1)(p-2)(p-3)}{4}\Big|_{-3.5+5j}^{.5+5j} + p(p-1)(p-2)\Big|_{-3.5+5j}^{.5+5j} + \frac{p(p-1)}{2}\Big|_{-3.5+5j}^{.5+5j}\right]$$

(3.4-35)

$$G_C = 92.5+269j \qquad (3.4\text{-}36)$$

<u>Vector Group D</u>

$$z_k = 2.5+k(-.5)\;, \quad k = 0,1,2,3 \qquad (3.4\text{-}37)$$

$$\Delta z = -.5 \qquad (3.4\text{-}38)$$

$u = 4$ , the number of vectors in Group D (3.4-39)

$A_g = .5x.5 = .25$ , the area of a grid square (3.4-40)

From Eq 3.4-6

$$G_D = \frac{2j}{A_g} \; {}_1\!\int_0^u \left[\left(z_k + \frac{\Delta z}{2}\right)^3 \Delta z\right] \Delta k$$ , the Group D contour segment integral equation (3.4-41)

Substituting Eq 3.4-37 thru Eq 3.4-40 into Eq 3.4-41

$$G_D = \frac{2j}{.25}\ {}_1\!\int_0^4 \left[\left(2.5+k(-.5) + \frac{-.5}{2}\right)^3(-.5)\right]\Delta k = -4j\ {}_1\!\int_0^4 [2.25+k(-.5)]^3\,\Delta k$$

$$G_D = \frac{j}{2}\ {}_1\!\int_0^4 \left[\frac{2.25}{-.5} + k\right]^3\Delta k$$

$$G_D = \frac{j}{2}\ {}_1\!\int_0^4 [k - 4.5]^3\Delta k \qquad (3.4\text{-}42)$$

Changing variables

Let $p = k - 4.5$

$\Delta p = \Delta k = 1$

for $k = 0$ $p = -4.5$

for $k = 4$ $p = -.5$

Substituting into Eq 3.4-32

$$G_D = \frac{j}{2}\ {}_1\!\int_{-4.5}^{-.5} p^3\Delta p \qquad (3.4\text{-}43)$$

From the Integral equations in Table 6 in the Appendix

$$_{\Delta x}\!\int x^3\Delta x = \frac{x(x-\Delta x)(x-2\Delta x)(x-3\Delta x)}{4} + x(x-\Delta x)(x-2\Delta x)\Delta x + \frac{x(x-\Delta x)}{2}\Delta x^2 + k \qquad (3.4\text{-}44)$$

Using Eq 3.4-44 to integrate Eq 3.4-43

$x = p$

$\Delta x = \Delta p = 1$

$$G_D = \frac{j}{2}\ {}_1\!\int_{-4.5}^{-.5} p^3\Delta p = \frac{j}{2}\left[\frac{p(p-1)(p-2)(p-3)}{4}\Big|_{-4.5}^{-.5} + p(p-1)(p-2)\Big|_{-4.5}^{-.5} + \frac{p(p-1)}{2}\Big|_{-4.5}^{-.5}\right]$$

(3.4-45)

$G_D = -76.5j$

(3.4-46)

Then

$A_c$ is found by adding the four group integrals

$$A_c = G_A + G_B + G_C + G_D = (-148.75-96j) + (58.75-96.5j) + (92.5+269j) + (-76.5j) \qquad (3.4\text{-}47)$$

$$A_c = 2.5 \qquad (3.4\text{-}48)$$

**The area within the discrete closed contour is 2.5**

Checking the above result

Referring to the contour diagram, Diagram 3.4-1

$$A_c = 10 \text{ grid squares} * \frac{.25}{\text{grid square}} = 2.5$$

Good check

It is obvious that the integration process can save a considerable amount of computational time if a complex plane contour is composed of many vectors  There is another method for finding the area within a complex plane closed contour. This method , which uses integration in its derivation, is derived in the next section, Section 3.5.

## Section 3.5:  Area Calculation Using Discrete Complex Plane Closed Contour Corner Points

Often, in a typical discrete closed contour in the complex plane, many vectors are placed head to tail in straight lines. When calculating contour enclosed area, it would be useful to have an equation to apply to such multi-vector straight lines. This equation is derived below.

For clockwise contour integration
From Eq 3.4-4 and Eq 3.4-5

$$A_c = \frac{2j}{A_g}\ {}_1\oint_0^v [( z_k + \frac{\Delta z}{2})^3 \Delta z]\ \Delta k = -\frac{2j}{A_g}\ {}_1\oint_0^v [( z_k + \frac{\Delta z}{2})^3 \Delta z]\ \Delta k \qquad (3.5\text{-}1)$$

The closed contour area equations of Eq 3.5-1 above are functions of k where k = 0,1,2,3,…,v and z is a function of k. $\Delta k = 1$.

Change the form of Eq (3.5-1)

$$\mathbf{A_c = \frac{2j}{A_g}\ {}_{\Delta z}\oint ( z + \frac{\Delta z}{2})^3 \Delta z\ =\ -\frac{2j}{A_g}\ {}_{\Delta z}\oint ( z + \frac{\Delta z}{2})^3 \Delta z} \qquad (3.5\text{-}2)$$

**where**

**$A_c$ = the area enclosed by the vectors forming a discrete closed contour in the complex plane**
**z = the complex plane coordinates of the tail points of the discrete contour vectors**
**Δz = real or imaginary vector value**

$$\Delta z = \begin{cases} +\delta \text{ for a horizontal left to right pointing vector} \\ -\delta \text{ for a horizontal right to left pointing vector} \\ +j\delta \text{ for a vertical down to up pointing vector} \\ -j\delta \text{ for a vertical up to down pointing vector} \end{cases}$$

**z + Δz = the coordinates of the head points of the discrete contour vectors**
**$A_g$ = δxδ , the area of a complex plane grid square**
**δ = the length of the side of a complex plane grid square**

**Should a numerical value be substituted for the integral subscript symbol Δz, it would be δ = |Δz|.**

The closed contour area equations of Eq 3.5-2 above are functions of z

Note that the Eq 3.5-2 discrete closed contour integral is shown above without specifically designating the contour beginning and ending vectors as in Eq 3.5-1. With the integral functions now in terms of z, it no longer adds clarity to the closed contour integration to specify the beginning and ending contour vectors. The initial and final points on a discrete closed contour are one in the same.

Diagram 3.5-1

Complex plane contour straight line segment

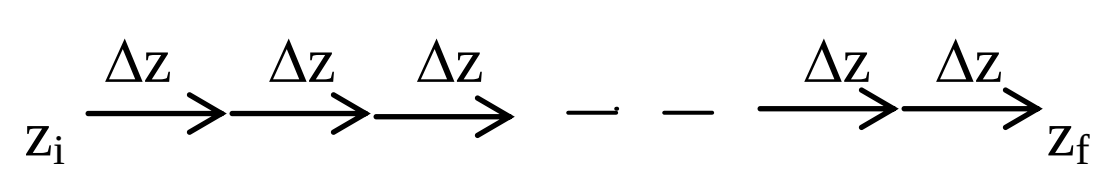

Discrete contours are formed only from horizontal or vertical vectors. The above contour in-line vectors are shown pointing to the right. They may also point to the left, point upward, or point downward.

Using Eq 3.5-2, integrate along a contour straight line segment composed of head to tail vectors.
See Diagram 3.12 above.

Integrate to calculate G

$$G = \underset{\Delta z}{\int_{z_i}^{z_f}} \left( z + \frac{\Delta z}{2} \right)^3 \Delta z \qquad (3.5\text{-}3)$$

where

G = a straight line vector group integral
$\Delta z$ = real or imaginary constant value

$$\Delta z = \begin{cases} +\delta \text{ for horizontal left to right pointing vectors} \\ -\delta \text{ for horizontal right to left pointing vectors} \\ +j\delta \text{ for vertical down to up pointing vectors} \\ -j\delta \text{ for vertical up to down pointing vectors} \end{cases}$$

$\delta$ = the length of the side of a complex plane grid square
$z$ = the complex plane coordinates of the tail points of the in-line contour vectors
$z + \Delta z$ = the coordinates of the head points of the contour vectors
$z_i$ = the initial value of z
$z_f$ = the final value of z

Evaluate the integral of Eq 3.5-3 using Eq 3.4-44 which is rewritten below

$$_{\Delta x}\int x^3 \Delta x = \frac{x(x-\Delta x)(x-2\Delta x)(x-3\Delta x)}{4} + x(x-\Delta x)(x-2\Delta x)\Delta x + \frac{x(x-\Delta x)}{2}\Delta x^2 + k \tag{3.5-4}$$

Let $$x = z + \frac{\Delta z}{2} \tag{3.5-5}$$

$$\Delta x = \Delta z \tag{3.5-6}$$

From Eq 3.5-3 thru 3.5-6

$$_{\Delta z}\int_{z_i}^{z_f} \left( z + \frac{\Delta z}{2} \right)^3 \Delta z = \frac{(z + \frac{\Delta z}{2})(z + \frac{\Delta z}{2} - \Delta z)(z + \frac{\Delta z}{2} - 2\Delta z)(z + \frac{\Delta z}{2} - 3\Delta z)}{4}$$

$$+ (z + \tfrac{\Delta z}{2})(z + \tfrac{\Delta z}{2} - \Delta z)(z + \tfrac{\Delta z}{2} - 2\Delta z)\Delta z + \frac{(z + \frac{\Delta z}{2})(z + \frac{\Delta z}{2} - \Delta z)}{2}\Delta z^2 \Big|_{z_i}^{z_f} \tag{3.5-7}$$

simplifying Eq 3.5-7

$$_{\Delta z}\int_{z_i}^{z_f} \left( z + \frac{\Delta z}{2} \right)^3 \Delta z = \frac{(z + \frac{\Delta z}{2})(z - \frac{\Delta z}{2})(z - \frac{3\Delta z}{2})(z - \frac{5\Delta z}{2})}{4}$$

$$+ (z + \tfrac{\Delta z}{2})(z - \tfrac{\Delta z}{2})(z - \tfrac{3\Delta z}{2})\Delta z + \frac{(z + \frac{\Delta z}{2})(z - \frac{\Delta z}{2})}{2}\Delta z^2 \Big|_{z_i}^{z_f} \tag{3.5-8}$$

Expanding the three right side terms of Eq 3.5-8

$$\Delta z \int_{z_i}^{z_f} ( z + \frac{\Delta z}{2} )^3 \Delta z = ( \frac{1}{4} z^4 - z^3\Delta z + \frac{7}{8} z^2\Delta z + \frac{1}{4} z\Delta z^3 - \frac{15}{64} \Delta z^4 ) + ( z^3\Delta z - \frac{3}{2} z^2\Delta z^2 - \frac{1}{4} z\Delta z^3 + \frac{3}{8}\Delta z^4 )$$

$$+ ( \frac{1}{2} z^2\Delta z^2 - \frac{1}{8} \Delta z^4 ) \Big|_{z_i}^{z_f} \qquad (3.5\text{-}9)$$

Combining like terms of Eq 3.5-9

$$\Delta z \int_{z_i}^{z_f} ( z + \frac{\Delta z}{2} )^3 \Delta z = \frac{1}{4} z^4 - \frac{1}{8} z^2\Delta z^2 + \frac{1}{64} \Delta z^4 \Big|_{z_i}^{z_f} \qquad (3.5\text{-}10)$$

Then

For complex plane integratation along a contour straight line segment composed of head to tail vectors See Diagram 3.5-1.

$$\Delta z \int_{z_i}^{z_f} ( z + \frac{\Delta z}{2} )^3 \Delta z = \frac{1}{4} z^4 - \frac{1}{8} z^2\Delta z^2 \Big|_{z_i}^{z_f} = \frac{z^2}{8} ( 2z^2 - \Delta z^2 ) \Big|_{z_i}^{z_f} \qquad (3.5\text{-}11)$$

**where**

**z = the complex plane coordinates of the tail points of the contour vectors**

**Δz = real or imaginary constant value**

$$\Delta z = \begin{cases} +\delta \text{ for horizontal left to right pointing vectors} \\ -\delta \text{ for horizontal right to left pointing vectors} \\ +j\delta \text{ for vertical down to up pointing vectors} \\ -j\delta \text{ for vertical up to down pointing vectors} \end{cases}$$

**δ = the length of the side of a complex plane grid square**

**z + Δz = the coordinates of the head points of the contour vectors**

**$z_i$ = the initial value of z**

**$z_f$ = the final value of z**

To better understand Eq 3.4-11 and its use in finding the area within a complex plane discrete closed contour, the following four examples, Example 3.5-1 thru Example 3.5-4 are provided. These four examples will calculate the area enclosed within the same discrete closed contour. However, for each example, the starting point on the contour and the direction of integration (clockwise or counterclockwise) will be different.

Example 3.5-1

Find the area enclosed within the following discrete closed contour.
Integrate in a clockwise direction starting a discrete contour corner point from which integration proceeds in a vertical direction.

Diagram 3.5-2  A complex plane discrete closed contour

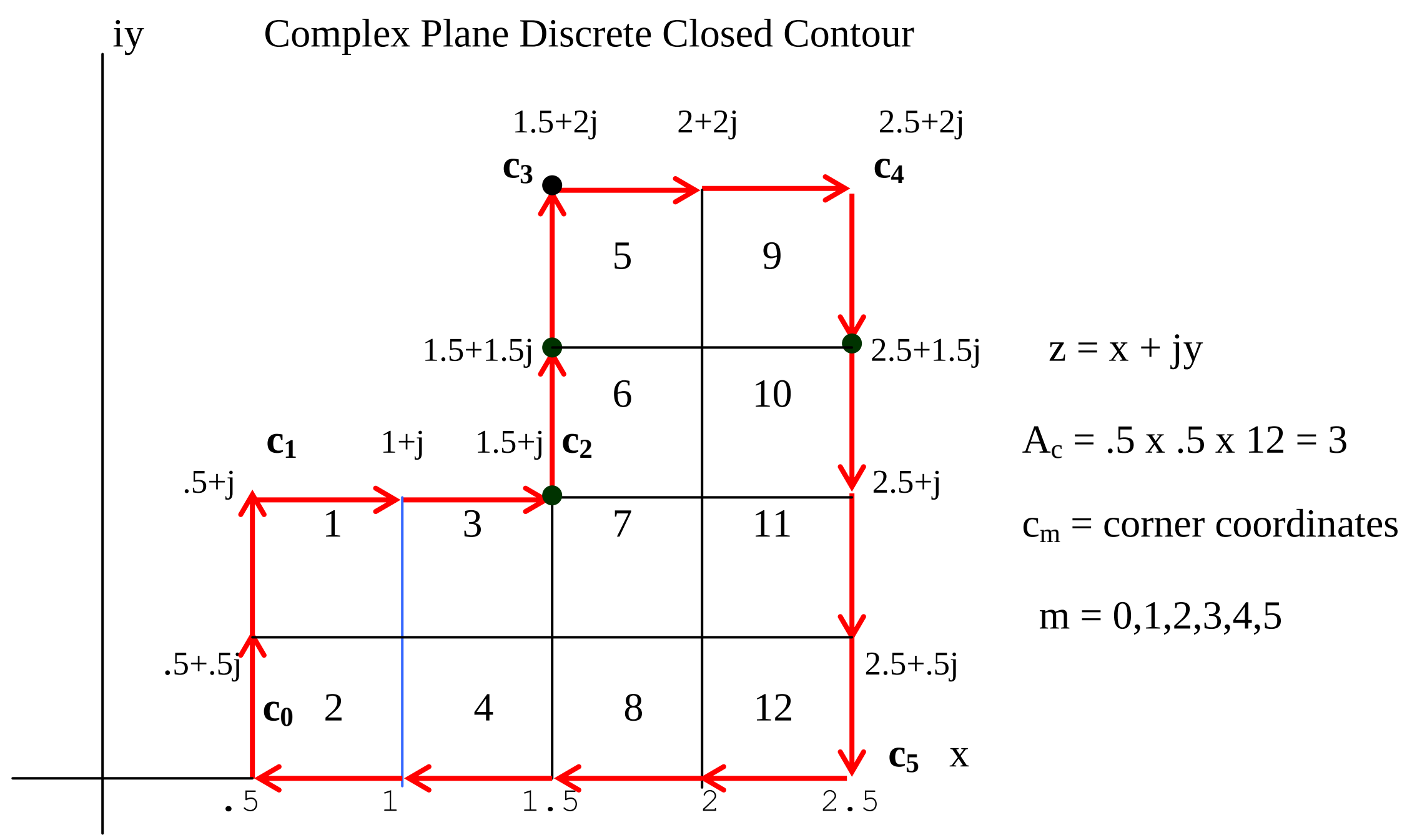

Apply Eq 3.5-2 and Eq 3.5-11 to the calculation of the area enclosed within the complex plane discrete closed contour of Diagram 3.5-2.

Rewriting Eq 3.5-2 (clockwise direction)

$$A_c = \frac{2j}{A_g} \oint_{\Delta z} ( z + \frac{\Delta z}{2} )^3 \Delta z \tag{3.5-12}$$

The discrete closed clockwise contour of Diagram 3.5-2 is composed of six in-line vector segments.

From Eq 3.5-12

$$A_c = \frac{2j}{A_g} \oint_{\Delta z} ( z + \frac{\Delta z}{2} )^3 \Delta z = \frac{2j}{A_g} [ \int_{\Delta z\, c_0}^{c_1} ( z + \frac{\Delta z}{2} )^3 \Delta z + \int_{\Delta z\, c_1}^{c_2} ( z + \frac{\Delta z}{2} )^3 \Delta z + \int_{\Delta z\, c_2}^{c_3} ( z + \frac{\Delta z}{2} )^3 \Delta z$$

$$+ \int_{\Delta z\, c_3}^{c_4} ( z + \frac{\Delta z}{2} )^3 \Delta z + \int_{\Delta z\, c_4}^{c_5} ( z + \frac{\Delta z}{2} )^3 \Delta z + \int_{\Delta z\, c_5}^{c_0} ( z + \frac{\Delta z}{2} )^3 \Delta z ]$$

(3.5-13)

Substituting Eq 3.5-11 into Eq 3.5-13

$$A_c = \frac{2j}{A_g} \oint_{\Delta z} \left( z + \frac{\Delta z}{2} \right)^3 \Delta z = \frac{2j}{A_g} \Big[ \left( \frac{1}{4} z^4 - \frac{1}{8} z^2 \{+j\delta\}^2 \right) \Big|_{c_0}^{c_1} + \left( \frac{1}{4} z^4 - \frac{1}{8} z^2 \{+\delta\}^2 \right) \Big|_{c_1}^{c_2}$$

$$+ \left( \frac{1}{4} z^4 - \frac{1}{8} z^2 \{+j\delta\}^2 \right) \Big|_{c_2}^{c_3} + \left( \frac{1}{4} z^4 - \frac{1}{8} z^2 \{+\delta\}^2 \right) \Big|_{c_3}^{c_4}$$

$$+ \left( \frac{1}{4} z^4 - \frac{1}{8} z^2 \{-j\delta\}^2 \right) \Big|_{c_4}^{c_5} + \left( \frac{1}{4} z^4 - \frac{1}{8} z^2 \{-\delta\}^2 \right) \Big|_{c_5}^{c_0} \Big] \qquad (3.5\text{-}14)$$

From Eq 3.5-14, combine the z terms and substitute the limits

$$A_c = \frac{2j}{A_g} \oint_{\Delta z} \left( z + \frac{\Delta z}{2} \right)^3 \Delta z = \frac{j}{2A_g} [ (c_1^4 - c_0^4 + c_2^4 - c_1^4 + c_3^4 - c_2^4 + c_4^4 - c_3^4 + c_5^4 - c_4^4 + c_0^4 - c_5^4 ]$$

$$- \frac{j}{4A_g} [ (+j\delta)^2 c_1^2 - (+j\delta)^2 c_0^2 + (+\delta)^2 c_2^2 - (+\delta)^2 c_1^2 + (+j\delta)^2 c_3^2 - (+j\delta)^2 c_2^2 ]$$

$$- \frac{j}{4A_g} [ (+\delta)^2 c_4^2 - (+\delta)^2 c_3^2 + (-j\delta)^2 c_5^2 - (-j\delta)^2 c_4^2 + (-\delta)^2 c_0^2 - (-\delta)^2 c_5^2 ]$$

(3.5-15)

Simplifying Eq 3.5-15

The first term on the right side of Eq 3.5-15 is seen to equal zero

$A_g = \delta x \delta$

$$A_c = \frac{2j}{A_g} \oint_{\Delta z} \left( z + \frac{\Delta z}{2} \right)^3 \Delta z = - \frac{jA_g}{4A_g} [ -c_1^2 + c_0^2 + c_2^2 - c_1^2 - c_3^2 + c_2^2 + c_4^2 - c_3^2 - c_5^2 + c_4^2 + c_0^2 - c_5^2 ]$$

(3.5-16)

Simplifying Eq 3.5-16

$$A_c = \frac{2j}{A_g} \oint_{\Delta z} \left( z + \frac{\Delta z}{2} \right)^3 \Delta z = - \frac{jA_g}{4A_g} [ -2c_1^2 + 2c_0^2 + 2c_2^2 - 2c_3^2 + 2c_4^2 - 2c_5^2 ]$$

$$A_c = \frac{2j}{A_g} \oint_{\Delta z} \left( z + \frac{\Delta z}{2} \right)^3 \Delta z = - \frac{j}{2} [ +c_0^2 - c_1^2 + c_2^2 - c_3^2 + c_4^2 - c_5^2 ] \qquad (3.5\text{-}17)$$

where

1. The contour vector which has as its tail point the initial point, $c_0$, is vertical.
2. The contour is traversed in a clockwise direction.

Evaluate the area, $A_c$, enclosed within the complex plane discrete closed contour of Diagram 3.5-2

From Diagram 3.5-2

$c_0 = .5$
$c_1 = .5+j$
$c_2 = 1.5+j$
$c_3 = 1.5+2j$
$c_4 = 2.5+2j$
$c_5 = 2.5$

Substituting the above values for $c_0$ thru $c_5$ into Eq 3.5-17

$$A_c = \frac{2j}{A_g} \oint_{\Delta z} \left( z + \frac{\Delta z}{2} \right)^3 \Delta z = -\frac{j}{2} \left[ +(.5)^2 - (.5+j)^2 + (1.5+j)^2 - (1.5+2j)^2 + (2.5+2j)^2 - (2.5)^2 \right] \qquad (3.5\text{-}18)$$

$$A_c = \frac{2j}{A_g} \oint_{\Delta z} \left( z + \frac{\Delta z}{2} \right)^3 \Delta z = 3 \qquad (3.5\text{-}19)$$

Checking the above value for $A_c$

Referring to Diagram 3.5-2

$A_c = .5 \text{ x } .5 \text{ x } 12 = 3$
Good check

<u>Comment</u> - The integration initial point is arbitrarily designated as $c_0$. The following points, progressing along the contour in the clockwise direction, are successively designated $c_1, c_2, c_3, \ldots$

The validity and the usefulness of Eq 3.5-2 and Eq 3.5-11 has been shown.

Example 3.5-2

Find the area enclosed within the following discrete closed contour.
Integrate in a clockwise direction starting a discrete contour corner point from which integration proceeds in a horizontal direction.

Diagram 3.5-3 A complex plane discrete closed contour

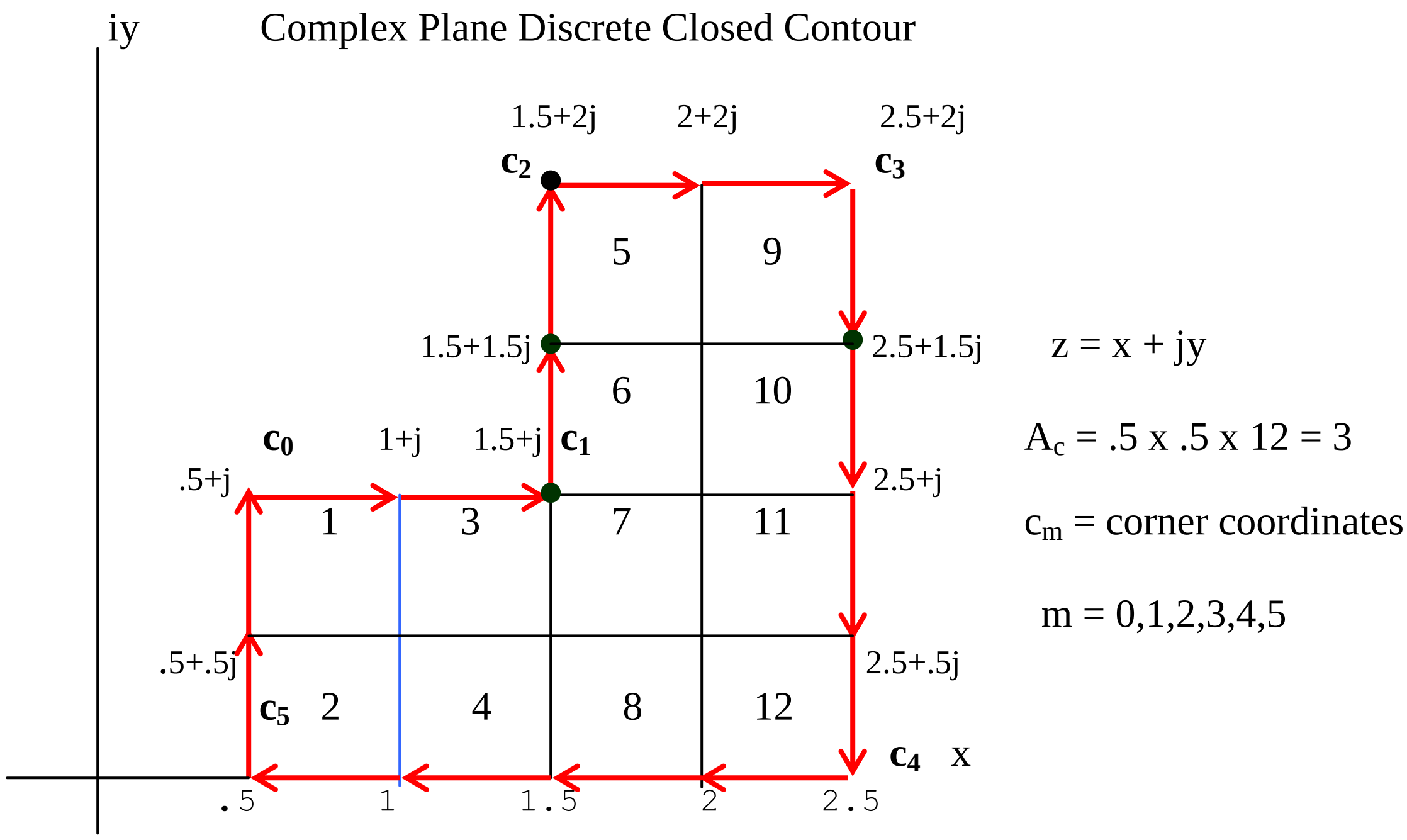

Apply Eq 3.5-2 and Eq 3.5-11 to the calculation of the area enclosed within the complex plane discrete closed contour of Diagram 3.5-3.

Rewriting Eq 3.5-2 (clockwise direction)

$$A_c = \frac{2j}{A_g} \oint_{\Delta z} ( z + \frac{\Delta z}{2} )^3 \Delta z \qquad (3.5\text{-}20)$$

The discrete closed clockwise contour of Diagram 3.5-3 is composed of six in-line vector segments.

From Eq 3.5-20

$$A_c = \frac{2j}{A_g} \oint_{\Delta z} ( z + \frac{\Delta z}{2} )^3 \Delta z = \frac{2j}{A_g} [ \int_{\Delta z \, c_0}^{c_1} ( z + \frac{\Delta z}{2} )^3 \Delta z + \int_{\Delta z \, c_1}^{c_2} ( z + \frac{\Delta z}{2} )^3 \Delta z + \int_{\Delta z \, c_2}^{c_3} ( z + \frac{\Delta z}{2} )^3 \Delta z$$

$$+ \int_{\Delta z \, c_3}^{c_4} ( z + \frac{\Delta z}{2} )^3 \Delta z + \int_{\Delta z \, c_4}^{c_5} ( z + \frac{\Delta z}{2} )^3 \Delta z + \int_{\Delta z \, c_5}^{c_0} ( z + \frac{\Delta z}{2} )^3 \Delta z \, ]$$

(3.5-21)

Substituting Eq 3.5-11 into Eq 3.5-21

$$A_c = \frac{2j}{A_g} \oint_{\Delta z} ( z + \frac{\Delta z}{2} )^3 \Delta z = \frac{2j}{A_g} [ ( \frac{1}{4} z^4 - \frac{1}{8} z^2\{+\delta\}^2 ) \Big|_{c_0}^{c_1} + ( \frac{1}{4} z^4 - \frac{1}{8} z^2\{+j\delta\}^2 ) \Big|_{c_1}^{c_2}$$

$$+ ( \frac{1}{4} z^4 - \frac{1}{8} z^2\{+\delta\}^2 ) \Big|_{c_2}^{c_3} + ( \frac{1}{4} z^4 - \frac{1}{8} z^2\{-j\delta\}^2 ) \Big|_{c_3}^{c_4}$$

$$+ ( \frac{1}{4} z^4 - \frac{1}{8} z^2\{-\delta\}^2 ) \Big|_{c_4}^{c_5} + ( \frac{1}{4} z^4 - \frac{1}{8} z^2\{+j\delta\}^2 ) \Big|_{c_5}^{c_0} ] \qquad (3.5\text{-}22)$$

From Eq 3.5-22, combine the z terms and substitute the limits

$$A_c = \frac{2j}{A_g} \oint_{\Delta z} ( z + \frac{\Delta z}{2} )^3 \Delta z = \frac{j}{2A_g} [ (c_1^4 - c_0^4 + c_2^4 - c_1^4 + c_3^4 - c_2^4 + c_4^4 - c_3^4 + c_5^4 - c_4^4 + c_0^4 - c_5^4 ]$$

$$- \frac{j}{4A_g} [ (+\delta)^2 c_1^2 - (+\delta)^2 c_0^2 + (+j\delta)^2 c_2^2 - (+j\delta)^2 c_1^2 + (+\delta)^2 c_3^2 - (+\delta)^2 c_2^2 ]$$

$$- \frac{j}{4A_g} [ (-j\delta)^2 c_4^2 - (-j\delta)^2 c_3^2 + (-\delta)^2 c_5^2 - (-\delta)^2 c_4^2 + (+j\delta)^2 c_0^2 - (+j\delta)^2 c_5^2 ]$$

(3.5-23)

Simplifying Eq 3.5-23

The first term on the right side of Eq 3.5-23 is seen to equal zero

$A_g = \delta x \delta$

$$A_c = \frac{2j}{A_g} \oint_{\Delta z} ( z + \frac{\Delta z}{2} )^3 \Delta z = - \frac{jA_g}{4A_g} [ c_1^2 - c_0^2 - c_2^2 + c_1^2 + c_3^2 - c_2^2 - c_4^2 + c_3^2 + c_5^2 - c_4^2 - c_0^2 + c_5^2 ]$$

(3.5-24)

Simplifying Eq 3.5-24

$$A_c = \frac{2j}{A_g} \oint_{\Delta z} ( z + \frac{\Delta z}{2} )^3 \Delta z = - \frac{jA_g}{4A_g} [ +2c_1^2 - 2c_0^2 - 2c_2^2 + 2c_3^2 - 2c_4^2 + 2c_5^2 ]$$

$$A_c = \frac{2j}{A_g} \oint_{\Delta z} ( z + \frac{\Delta z}{2} )^3 \Delta z = + \frac{j}{2} [ +c_0^2 + c_1^2 - c_2^2 + c_3^2 - c_4^2 + c_5^2 ] \qquad (3.5\text{-}25)$$

where

1. The contour vector which has as its tail point the initial point, $c_0$, is horizontal.
2. The contour is traversed in a clockwise direction.

Evaluate the area, $A_c$, enclosed within the complex plane discrete closed contour of Diagram 3.5-3

From Diagram 3.5-3

$c_0 = .5+j$
$c_1 = 1.5+j$
$c_2 = 1.5+2j$
$c_3 = 2.5+2j$
$c_4 = 2.5$
$c_5 = .5$

Substituting the above values for $c_0$ thru $c_5$ into Eq 3.5-25

$$A_c = \frac{2j}{A_g} \oint_{\Delta z} \left( z + \frac{\Delta z}{2} \right)^3 \Delta z = \frac{j}{2} [\, +(.5+j)^2 - (1.5+j)^2 + (1.5+2j)^2 - (2.5+2j)^2 + (2.5)^2 - (.5)^2 \,] \qquad (3.5\text{-}26)$$

$$A_c = \frac{2j}{A_g} \oint_{\Delta z} \left( z + \frac{\Delta z}{2} \right)^3 \Delta z = 3 \qquad (3.5\text{-}27)$$

Checking the above value for $A_c$

Referring to Diagram 3.5-3

$A_c = .5 \text{ x } .5 \text{ x } 12 = 3$
Good check

<u>Comment</u> - The integration initial point is arbitrarily designated as $c_0$. The following points, progressing along the contour in the clockwise direction, are successively designated $c_1, c_2, c_3, \ldots$

The validity and the usefulness of Eq 3.5-2 and Eq 3.5-11 has been shown.

Example 3.5-3

Find the area enclosed within the following discrete closed contour.
Integrate in a counterclockwise direction starting a discrete contour corner point from which integration proceeds in a horizontal direction.

Diagram 3.5-4  A complex plane discrete closed contour

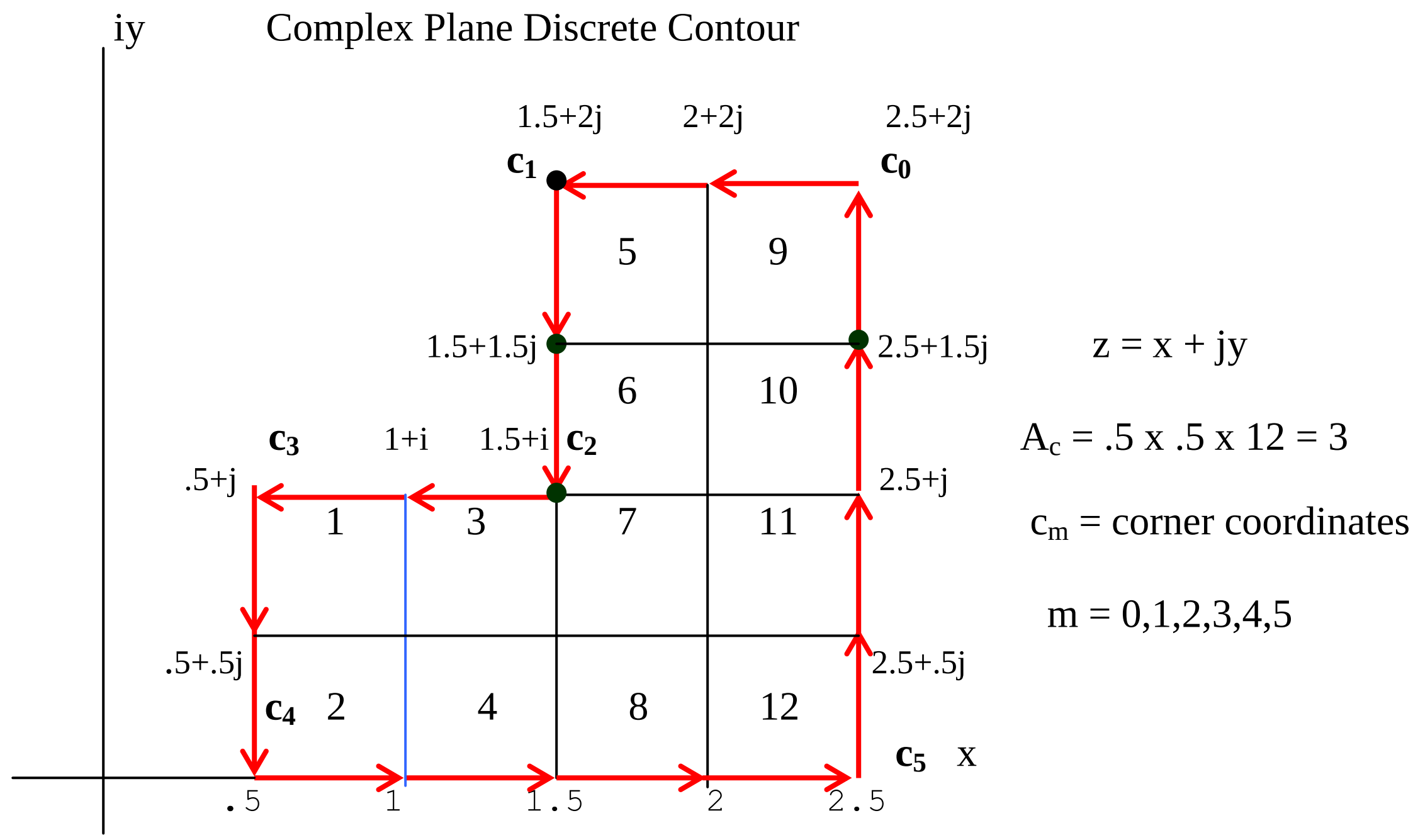

Apply Eq 3.5-2 and Eq 3.5-11 to the calculation of the area enclosed within the complex plane discrete closed contour of Diagram 3.5-4.

Rewriting Eq 3.5-2 (counterclockwise direction)

$$A_c = -\frac{2j}{A_g} \underset{\Delta z}{\oint} (z + \frac{\Delta z}{2})^3 \Delta z \qquad (3.5\text{-}28)$$

The discrete closed clockwise contour of Diagram 3.5-4 is composed of six in-line vector segments.

From Eq 3.5-28

$$A_c = -\frac{2j}{A_g} \underset{\Delta z}{\oint} (z + \frac{\Delta z}{2})^3 \Delta z = -\frac{2j}{A_g} [ \underset{\Delta z}{\int_{c_0}^{c_1}} (z + \frac{\Delta z}{2})^3 \Delta z + \underset{\Delta z}{\int_{c_1}^{c_2}} (z + \frac{\Delta z}{2})^3 \Delta z + \underset{\Delta z}{\int_{c_2}^{c_3}} (z + \frac{\Delta z}{2})^3 \Delta z$$

$$+ \underset{\Delta z}{\int_{c_3}^{c_4}} (z + \frac{\Delta z}{2})^3 \Delta z + \underset{\Delta z}{\int_{c_4}^{c_5}} (z + \frac{\Delta z}{2})^3 \Delta z + \underset{\Delta z}{\int_{c_5}^{c_0}} (z + \frac{\Delta z}{2})^3 \Delta z ] \qquad (3.5\text{-}29)$$

Substituting Eq 3.5-11 into Eq 3.5-29

$$A_c = -\frac{2j}{A_g}\oint_{\Delta z}\left(z+\frac{\Delta z}{2}\right)^3\Delta z = -\frac{2j}{A_g}\Big[\left(\frac{1}{4}z^4-\frac{1}{8}z^2\{-\delta\}^2\right)\Big|_{c_0}^{c_1} + \left(\frac{1}{4}z^4-\frac{1}{8}z^2\{-j\delta\}^2\right)\Big|_{c_1}^{c_2}$$
$$+\left(\frac{1}{4}z^4-\frac{1}{8}z^2\{-\delta\}^2\right)\Big|_{c_2}^{c_3} + \left(\frac{1}{4}z^4-\frac{1}{8}z^2\{-j\delta\}^2\right)\Big|_{c_3}^{c_4}$$
$$+\left(\frac{1}{4}z^4-\frac{1}{8}z^2\{+\delta\}^2\right)\Big|_{c_4}^{c_5} + \left(\frac{1}{4}z^4-\frac{1}{8}z^2\{+j\delta\}^2\right)\Big|_{c_5}^{c_0}\Big] \quad (3.5\text{-}30)$$

From Eq 3.5-30, combine the z terms and substitute the limits

$$A_c = -\frac{2j}{A_g}\oint_{\Delta z}\left(z+\frac{\Delta z}{2}\right)^3\Delta z = -\frac{j}{2A_g}\left[(c_1^4-c_0^4+c_2^4-c_1^4+c_3^4-c_2^4+c_4^4-c_3^4+c_5^4-c_4^4+c_0^4-c_5^4\right]$$
$$+\frac{j}{4A_g}\left[(-\delta)^2c_1^2-(-\delta)^2c_0^2+(-j\delta)^2c_2^2-(-j\delta)^2c_1^2+(-\delta)^2c_3^2-(-\delta)^2c_2^2\right]$$
$$+\frac{j}{4A_g}\left[(-j\delta)^2c_4^2-(-j\delta)^2c_3^2+(+\delta)^2c_5^2-(+\delta)^2c_4^2+(+j\delta)^2c_0^2-(+j\delta)^2c_5^2\right] \quad (3.5\text{-}31)$$

Simplifying Eq 3.5-31

The first term on the right side of Eq 3.5-31 is seen to equal zero

$A_g = \delta x \delta$

$$A_c = -\frac{2j}{A_g}\oint_{\Delta z}\left(z+\frac{\Delta z}{2}\right)^3\Delta z = +\frac{jA_g}{4A_g}\left[+c_1^2-c_0^2-c_2^2+c_1^2+c_3^2-c_2^2-c_4^2+c_3^2+c_5^2-c_4^2-c_0^2+c_5^2\right] \quad (3.5\text{-}32)$$

Simplifying Eq 3.5-32

$$A_c = -\frac{2j}{A_g}\oint_{\Delta z}\left(z+\frac{\Delta z}{2}\right)^3\Delta z = +\frac{jA_g}{4A_g}\left[+2c_1^2-2c_0^2-2c_2^2+2c_3^2-2c_4^2+2c_5^2\right]$$

$$A_c = -\frac{2j}{A_g}\oint_{\Delta z}\left(z+\frac{\Delta z}{2}\right)^3\Delta z = -\frac{j}{2}\left[+c_0^2-c_1^2+c_2^2-c_3^2+c_4^2-c_5^2\right] \quad (3.5\text{-}33)$$

where
1. The contour vector which has as its tail point the initial point, $c_0$, is horizontal.
2. The contour is traversed in a counterclockwise direction.

Evaluate the area, $A_c$, enclosed within the complex plane discrete closed contour of Diagram 3.5-4

From Diagram 3.5-4

$c_0 = 2.5+2j$
$c_1 = 1.5+2j$
$c_2 = 1.5+j$
$c_3 = .5+j$
$c_4 = .5$
$c_5 = 2.5$

Substituting the above values for $c_0$ thru $c_5$ into Eq 3.5-33

$$A_c = -\frac{2j}{A_g} \oint_{\Delta z} \left( z + \frac{\Delta z}{2} \right)^3 \Delta z = -\frac{j}{2} \left[ +(2.5+2j)^2-(1.5+2j)^2+(1.5+j)^2-(.5+j)^2+(.5)^2-(2.5)^2 \right] \quad (3.5\text{-}34)$$

$$A_c = -\frac{2j}{A_g} \oint_{\Delta z} \left( z + \frac{\Delta z}{2} \right)^3 \Delta z = 3 \quad (3.5\text{-}35)$$

Checking the above value for $A_c$

Referring to Diagram 3.5-4

$A_c = .5 \text{ x } .5 \text{ x } 12 = 3$

Good check

<u>Comment</u> - The integration initial point is arbitrarily designated as $c_0$. The following points, progressing along the contour in the counterclockwise direction, are successively designated $c_1, c_2, c_3, \ldots$

The validity and the usefulness of Eq 3.5-2 and Eq 3.5-11 has again been shown.

Example 3.5-4

Find the area enclosed within the following discrete closed contour.
Integrate in a counterclockwise direction starting a discrete contour corner point from which integration proceeds in a vertical direction.

Diagram 3.5-5 A complex plane discrete closed contour

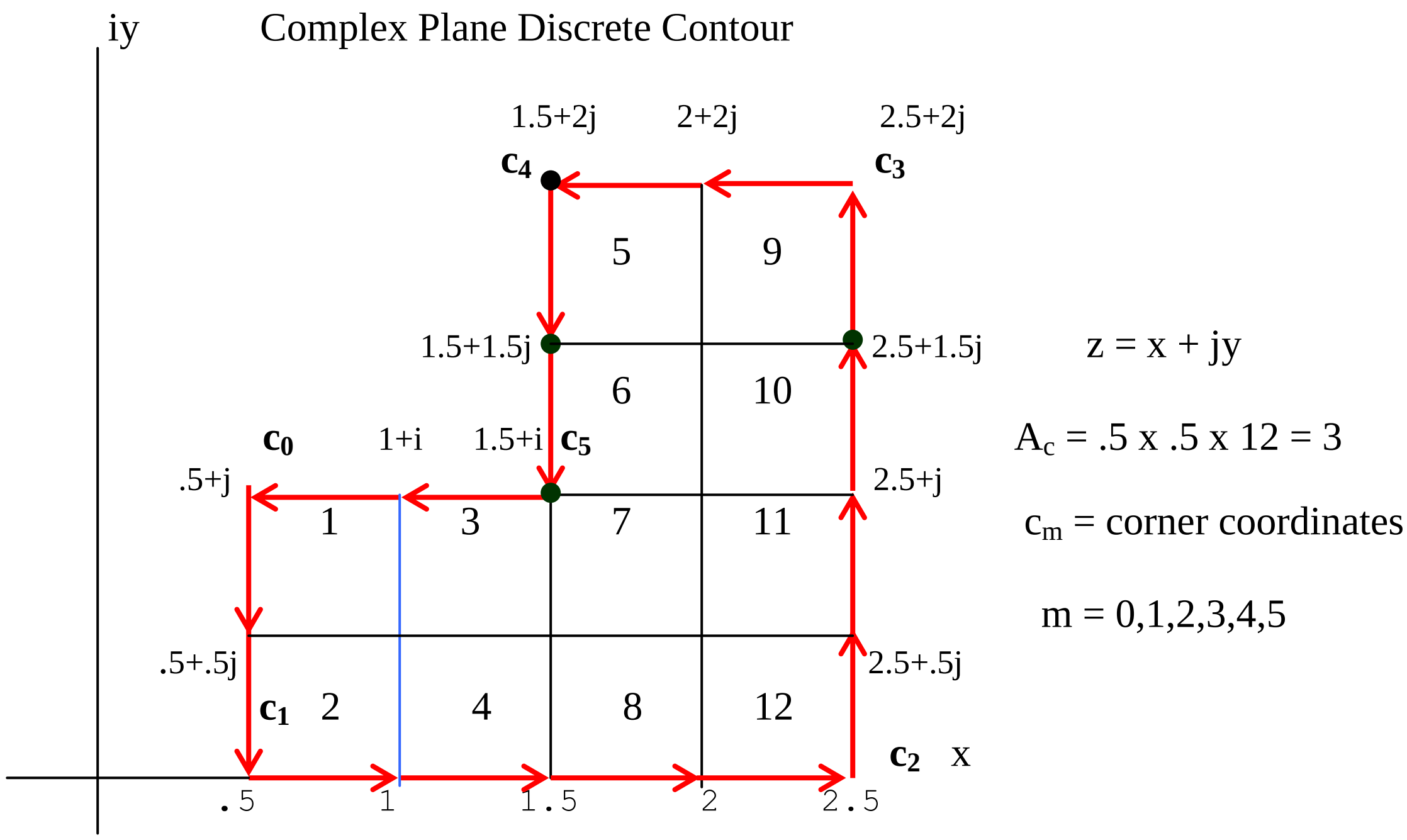

Apply Eq 3.5-2 and Eq 3.5-11 to the calculation of the area enclosed within the complex plane discrete closed contour of Diagram 3.5-5.

Rewriting Eq 3.5-2 (counterclockwise direction)

$$A_c = -\frac{2j}{A_g} \underset{\Delta z}{\oint} (z + \frac{\Delta z}{2})^3 \Delta z \tag{3.5-36}$$

The discrete closed clockwise contour of Diagram 3.5-5 is composed of six in-line vector segments.

From Eq 3.5-36

$$A_c = -\frac{2j}{A_g} \underset{\Delta z}{\oint} (z + \frac{\Delta z}{2})^3 \Delta z = -\frac{2j}{A_g} [ \underset{\Delta z}{\int_{c_0}^{c_1}} (z + \frac{\Delta z}{2})^3 \Delta z + \underset{\Delta z}{\int_{c_1}^{c_2}} (z + \frac{\Delta z}{2})^3 \Delta z + \underset{\Delta z}{\int_{c_2}^{c_3}} (z + \frac{\Delta z}{2})^3 \Delta z$$

$$+ \underset{\Delta z}{\int_{c_3}^{c_4}} (z + \frac{\Delta z}{2})^3 \Delta z + \underset{\Delta z}{\int_{c_4}^{c_5}} (z + \frac{\Delta z}{2})^3 \Delta z + \underset{\Delta z}{\int_{c_5}^{c_0}} (z + \frac{\Delta z}{2})^3 \Delta z ] \tag{3.5-37}$$

Substituting Eq 3.5-11 into Eq 3.5-37

$$A_c = -\frac{2j}{A_g}\oint_{\Delta z}\left(z+\frac{\Delta z}{2}\right)^3\Delta z = -\frac{2j}{A_g}\left[\left(\frac{1}{4}z^4-\frac{1}{8}z^2\{-j\delta\}^2\right)\Big|_{c_0}^{c_1}+\left(\frac{1}{4}z^4-\frac{1}{8}z^2\{+\delta\}^2\right)\Big|_{c_1}^{c_2}\right.$$
$$+\left(\frac{1}{4}z^4-\frac{1}{8}z^2\{+j\delta\}^2\right)\Big|_{c_2}^{c_3}+\left(\frac{1}{4}z^4-\frac{1}{8}z^2\{-\delta\}^2\right)\Big|_{c_3}^{c_4}$$
$$\left.+\left(\frac{1}{4}z^4-\frac{1}{8}z^2\{-j\delta\}^2\right)\Big|_{c_4}^{c_5}+\left(\frac{1}{4}z^4-\frac{1}{8}z^2\{-\delta\}^2\right)\Big|_{c_5}^{c_0}\right] \tag{3.5-38}$$

From Eq 3.5-38, combine the z terms and substitute the limits

$$A_c = -\frac{2j}{A_g}\oint_{\Delta z}\left(z+\frac{\Delta z}{2}\right)^3\Delta z = -\frac{j}{2A_g}\left[\left(c_1^4-c_0^4+c_2^4-c_1^4+c_3^4-c_2^4+c_4^4-c_3^4+c_5^4-c_4^4+c_0^4-c_5^4\right]\right.$$
$$+\frac{j}{4A_g}\left[(-j\delta)^2c_1^2-(-j\delta)^2c_0^2+(+\delta)^2c_2^2-(+\delta)^2c_1^2+(+j\delta)^2c_3^2-(+j\delta)^2c_2^2\right]$$
$$+\frac{j}{4A_g}\left[(-\delta)^2c_4^2-(-\delta)^2c_3^2+(+j\delta)^2c_5^2-(+j\delta)^2c_4^2+(-\delta)^2c_0^2-(-j\delta)^2c_5^2\right] \tag{3.5-39}$$

Simplifying Eq 3.5-39

The first term on the right side of Eq 3.5-39 is seen to equal zero

$A_g = \delta x \delta$

$$A_c = -\frac{2j}{A_g}\oint_{\Delta z}\left(z+\frac{\Delta z}{2}\right)^3\Delta z = +\frac{jA_g}{4A_g}\left[-c_1^2+c_0^2+c_2^2-c_1^2-c_3^2+c_2^2+c_4^2-c_3^2-c_5^2+c_4^2+c_0^2-c_5^2\right] \tag{3.5-40}$$

Simplifying Eq 3.5-40

$$A_c = -\frac{2j}{A_g}\oint_{\Delta z}\left(z+\frac{\Delta z}{2}\right)^3\Delta z = +\frac{jA_g}{4A_g}\left[-2c_1^2+2c_0^2+2c_2^2-2c_3^2+2c_4^2-2c_5^2\right]$$

$$A_c = -\frac{2j}{A_g}\oint_{\Delta z}\left(z+\frac{\Delta z}{2}\right)^3\Delta z = +\frac{j}{2}\left[+c_0^2-c_1^2+c_2^2-c_3^2+c_4^2-c_5^2\right] \tag{3.5-41}$$

where

1. The contour vector which has as its tail point the initial point, $c_0$, is vertical.
2. The contour is traversed in a counterclockwise direction.

Evaluate the area, $A_c$, enclosed within the complex plane discrete closed contour of Diagram 3.5-5

From Diagram 3.5-5

$c_0 = .5+j$
$c_1 = .5$
$c_2 = 2.5$
$c_3 = 2.5+2j$
$c_4 = 1.5+2j$
$c_5 = 1.5+j$

Substituting the above values for $c_0$ thru $c_5$ into Eq 3.5-41

$$A_c = -\frac{2j}{A_g} \underset{\Delta z}{\oint} \left( z + \frac{\Delta z}{2} \right)^3 \Delta z = \frac{j}{2} \left[ +(.5+j)^2 - (.5)^2 + (2.5)^2 - (2.5+2j)^2 + (1.5+2j)^2 - (1.5+j)^2 \right] \qquad (3.5\text{-}42)$$

$$A_c = -\frac{2j}{A_g} \underset{\Delta z}{\oint} \left( z + \frac{\Delta z}{2} \right)^3 \Delta z = 3 \qquad (3.5\text{-}43)$$

Checking the above value for $A_c$

Referring to Diagram 3.5-5

$A_c = .5 \text{ x } .5 \text{ x } 12 = 3$

Good check

<u>Comment</u> - The integration initial point is arbitrarily designated as $c_0$. The following points, progressing along the contour in the counterclockwise direction, are successively designated $c_1, c_2, c_3, \ldots$

The validity and the usefulness of Eq 3.5-2 and Eq 3.5-11 has again been shown.

Note that by picking any corner point on a discrete contour as the initial point, $c_0$, and selecting either summation direction, clockwise or counterclockwise, the resulting calculated area within the discrete contour is the same positive value.

From Eq 3.5-17 and Eq 3.5-25, Eq 3.5-33, and Eq 3.5-41 a general equation for finding the area enclosed within a discrete complex plane contour can be derived.

From Eq 3.5-17

$$A_c = +\frac{2j}{A_g} \oint_{\Delta z} \left( z + \frac{\Delta z}{2} \right)^3 \Delta z = (-1)\left[ \frac{j}{2} (c_0^2 - c_1^2 + c_2^2 - c_3^2 + c_4^2 - c_5^2) \right] \tag{3.5-44}$$

where

1. The contour vector which has as its tail point the initial point, $c_0$, is vertical.
2. The contour is traversed in a clockwise direction.

From Eq 3.5-25

$$A_c = +\frac{2j}{A_g} \oint_{\Delta z} \left( z + \frac{\Delta z}{2} \right)^3 \Delta z = (+1)\left[ \frac{j}{2} (c_0^2 - c_1^2 + c_2^2 - c_3^2 + c_4^2 - c_5^2) \right] \tag{3.5-45}$$

where

1. The contour vector which has as its tail point the initial point, $c_0$, is horizontal.
2. The contour is traversed in a clockwise direction.

Rewriting Eq 3.5-33

$$A_c = -\frac{2j}{A_g} \oint_{\Delta z} \left( z + \frac{\Delta z}{2} \right)^3 \Delta z = (-1)\left[ \frac{j}{2} (c_0^2 - c_1^2 + c_2^2 - c_3^2 + c_4^2 - c_5^2) \right] \tag{3.5-46}$$

where

1. The contour vector which has as its tail point the initial point, $c_0$, is horizontal.
2. The contour is traversed in a counterclockwise direction.

From Eq 3.5-41

$$A_c = -\frac{2j}{A_g} \oint_{\Delta z} \left( z + \frac{\Delta z}{2} \right)^3 \Delta z = (+1)\left[ \frac{j}{2} (c_0^2 - c_1^2 + c_2^2 - c_3^2 + c_4^2 - c_5^2) \right] \tag{3.5-47}$$

where

1. The contour vector which has as its tail point the initial point, $c_0$, is vertical.
2. The contour is traversed in a counterclockwise direction.

Reviewing the derivations of Eq 3.5-44 thru Eq 3.5-47, the sign of the function, $\pm\left[\frac{j}{2}(c_0^2-c_1^2+c_2^2-c_3^2+c_4^2-c_5^2)\right]$, is seen to be a function of the product of two signs. One sign is related to the direction of the complex plane discrete closed contour vector summation. This sign is specified as follows:

$$\text{closed contour vector summation sign} = \begin{cases} + & \text{for clockwise closed contour summation} \\ - & \text{for counterclockwise closed contour summation} \end{cases} \qquad (3.5\text{-}48)$$

The other sign is related to the direction of the initial vector of the summation through the function, $\frac{\Delta z^2}{\delta^2}$. As previously specified:

$$\Delta z = \begin{cases} +\delta \text{ for horizontal left to right pointing vectors} \\ -\delta \text{ for horizontal right to left pointing vectors} \\ +j\delta \text{ for vertical down to up pointing vectors} \\ -j\delta \text{ for vertical up to down pointing vectors} \end{cases} \qquad (3.5\text{-}49)$$

Thus the signs associated with the function, $\frac{\Delta z^2}{\delta^2}$, are:

$$\frac{\Delta z^2}{\delta^2}\text{ sign} = \begin{cases} + \text{ for horizontal left to right pointing vectors} \\ + \text{ for horizontal right to left pointing vectors} \\ - \text{ for vertical down to up pointing vectors} \\ - \text{ for vertical up to down pointing vectors} \end{cases} \qquad (3.5\text{-}50)$$

or simplifying

$$\frac{\Delta z^2}{\delta^2}\text{ sign} = \begin{cases} + \text{ for horizontal pointing vectors} \\ - \text{ for vertical pointing vectors} \end{cases} \qquad (3.5\text{-}51)$$

Finally, from Eq 3.5-44 thru Eq 3.5-47, their mathematical derivations, and their application in Examples 3.5-1 thru Example 3.5-4, the following observations have been made.

1. The area, $A_c$, is real and positive.
2. The quantity, $\frac{j}{2}(c_0^2-c_1^2+c_2^2-c_3^2+c_4^2-c_5^2)$, changes sign as a function of the direction that the discrete contour is traversed (i.e. clockwise or counterclockwise) and the position of the vector which has as its tail point the initial point, $c_0$ (i.e. horizontal or vertical). Note the signs in Eq 3.5-44, Eq 3.5-45, Eq 3.5-46, and Eq 3.5-47.
3. The area enclosed within a discrete closed contour in the complex plane can be obtained solely from the coordinates of the corner points of the contour.

4. The equation for finding the enclosed contour area is not a function of the complex plane grid dimensions. They cancel out.
5. The minimum number of corner points in a discrete complex plane closed contour is 4, the number of corner points in a rectangle. The number of corner points in any discrete complex plane contour is an even number.
6. The number of horizontal sides and the number of vertical sides of any discrete complex plane discrete closed contour is the same. This can be easily proven.
7. Though the integral equation, Eq 3.5-11, has two terms, a $z^4$ term and a $z^2$ term, in a discrete closed contour the $z^4$ terms cancel out leaving only the $z^2$ terms.
8. The $\frac{\Delta z^2}{\delta^2}$ term in the derivation of Eq 3.5-43 thru Eq 3.5-47 indicates contour vector direction by a sign, $(\pm 1)^2 = +1$ for a horizontal direction and $(\pm j)^2 = -1$ for a vertical direction.

From Eq 3.5-44, Eq 3.5-45, Eq 3.5-46, Eq 3.5-47, and the above observations a general equation to calculate the area enclosed within any discrete complex plane contour can be derived. This equation is presented below.

**Discrete complex plane contour area equation**

$$\mathbf{A_c = +\frac{2j}{A_g}\oint_{\Delta z}\!\!\circlearrowright (z+\frac{\Delta z}{2})^3 \Delta z = -\frac{2j}{A_g}\oint_{\Delta z}\!\!\circlearrowleft (z+\frac{\Delta z}{2})^3 \Delta z = \alpha\,\frac{j}{2}\sum_{p=0}^{N}(-1)^p c_p^{\,2}} \qquad (3.5\text{-}52)$$

**where**

$$\alpha = \begin{cases} +1 & \text{for clockwise contour vector summation with initial vector horizontal} \\ +1 & \text{for counterclockwise contour vector summation with initial vector vertical} \\ -1 & \text{for clockwise contour vector summation with initial vector vertical} \\ -1 & \text{for counterclockwise contour vector summation with initial vector horizontal} \end{cases}$$

**$A_c$ = The area enclosed within a complex plane discrete closed contour**
**N = 3,5,7,9,11,…**
**N = the number of the discrete complex plane contour corner points minus one**
**$c_p$ = the coordinates of the corner points of the discrete complex plane contour**
**p = 0,1,2,3,4,…,N**

<u>**Comments**</u> **- 1. The value of the summation, $\frac{j}{2}\sum_{p=0}^{N}(-1)^p c_p^{\,2}$ , is always real.**

**2. The value of α is such that the area, $A_c$, is always positive.**

Since the area, $A_c$, is always real and positive, Eq 3.5-52 can be simplified. This general equation to calculate the area enclosed within any discrete complex plane contour is presented below.

**Discrete complex plane contour area equation**

$$\mathbf{A_c = +\frac{2j}{A_g}\oint_{\Delta z}\!\!\circlearrowright (z+\frac{\Delta z}{2})^3 \Delta z = -\frac{2j}{A_g}\oint_{\Delta z}\!\!\circlearrowleft (z+\frac{\Delta z}{2})^3 \Delta z = \left|\frac{j}{2}\sum_{p=0}^{N}(-1)^p c_p^{\,2}\right|} \qquad (3.5\text{-}53)$$

**where**

**$A_c$ = The area enclosed within a complex plane discrete closed contour**
**N = 3,5,7,9,11,…**
**N = the number of the discrete complex plane contour corner points minus one**
**$c_p$ = the coordinates of the corner points of the discrete complex plane contour**
**p = 0,1,2,3,4,…,N**
**|r| = absolute value of r**

**The contour integration can begin at any of the contour corner points. The integration initial point is designated as $c_0$. The following points, progressing along the contour in the direction of integration, are successively designated $c_1, c_2, c_3, \ldots$**

Interestingly, the complex plane discrete closed contour area equation stated above shows the contour enclosed area to be solely a function of the complex plane coordinates of the contour corner points.

Comments - 1. The value of the summation, $\frac{j}{2}\sum_{p=0}^{N}(-1)^p c_p^{\,2}$ , is real.

2. The above summation is negative if the contour between $c_0$ and $c_1$ is horizontal and p increases counterclockwise around the contour or if the contour between $c_0$ and $c_1$ is vertical and p increases clockwise around the contour else the summation is positive.

Example 3.5-5

Find the area enclosed within the complex plane discrete closed contour shown in Diagram 3.5-6 below. Use Eq 3.5-53 to calculate the area, $A_c$.

Diagram 3.5-6 A complex plane closed contour

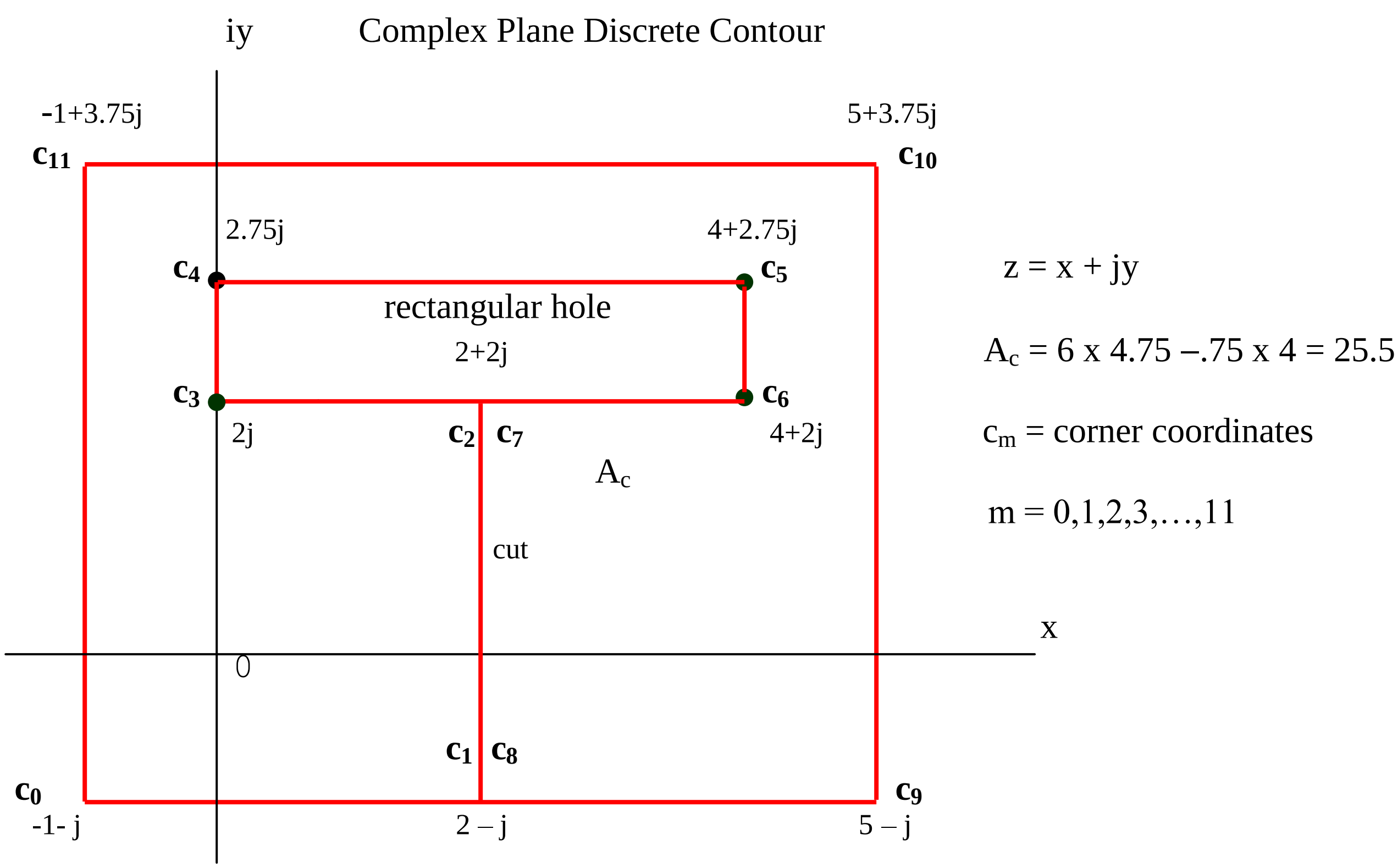

$A_c$, the area enclosed within the complex plane discrete closed contour of Diagram 3.5-6, is calculated from Eq 3.5-53 rewritten below.

$$A_c = \left| \frac{j}{2} \sum_{p=0}^{N} (-1)^p c_p^{\,2} \right|$$

where

$A_c$ = The area enclosed within a discrete complex plane closed contour

N = 3,5,7,9,11,…

N = the number of the discrete complex plane contour corner points minus one

$c_p$ = the coordinates of the corner points of the discrete complex plane contour

p = 0,1,2,3,4,…,N

|r| = absolute value of r

Note - The integration initial point is designated as $c_0$. The following points, progressing along the contour in the direction of integration, are successively designated $c_1, c_2, c_3, \ldots$

The following UBASIC program uses the above equation to calculate the area enclosed within the complex plane discrete closed contour shown in Diagram 3.5-6.

```
 10   dim C(20)
 20   C(0)=-1-#i
 30   C(1)=2-#i
 40   C(2)=2+2#i
 50   C(3)=2#i
 60   C(4)=2.75#i
 70   C(5)=4+2.75#i
 80   C(6)=4+2#i
 90   C(7)=2+2#i
100   C(8)=2-#i
110   C(9)=5-#i
120   C(10)=5+3.75#i
130   C(11)=-1+3.75#i
140   S=0
150   for K=0 to 11
160   S=S+((-1)^K)*C(K)^2
170   next K
180   W=abs((#i/2)*S)
190   print "The closed contour enclosed area is ";W
```

The evaluation obtained by running the above program is as follows:

**The closed contour enclosed area is 25.5**

Checking the above result

Referring to the contour diagram, Diagram 3.5-6

$A_c$ = 6 x 4.75 –.75 x 4 = 25.5 Good check

## Section 3.6: The Derivation of Additional Equations to Calculate the Area of a Closed Contour in the Complex Plane

There are a number of useful complex plane closed contour area calculation equations. These equations include and are derived from Eq 3.5-52. All equations are presented in Table 3.6-1 on the following page. Equation variables, constants, and some additional equations are defined at the end of the table. The area calculation equations listed will be derived in this section, Section 3.6.

The equations in Table 3.6-1 apply to three distict complex plane closed contour shapes. All of these three complex plane closed contour shapes consist of straight line vectors connected head to tail forming a closed loop. Corner points are designated as $z_n$, complex coordinates where $n = 1, 2, 3, \ldots, N$ with N being the total number of corner points that define the closed contour.The first shape, Shape S1, to which Eq 3.5-52 directly applies, consists of interconnected alternating horizontal and vertical vectors of finite length. For this case, N is always an even number. The second shape, Shape S2, consists of interconnected vectors of finite length with any slope. For this case, the number of corner points, N, need not be an even number. Shape S1 is a special case of Shape S2. The third shape, Shape S3, consists of interconnected vectors of infinitessimal length with any slope. For this case, the number, N, is infinite and the closed loop is formed by a continuous closed curve. Shape S3 is a special case of Shape S2 where N becomes infinitely large.

**Table 3.6-1** **Complex Plane Closed Contour Area Calculation Equations**

| # | Closed Contour Area Calculation Equation | Initial Direction of Summation/ Integration | Comments | Equation Applicable to Complex Plane Closed Contours Composed of | | |
|---|---|---|---|---|---|---|
| | | | | Vertical/Horizontal Straight Line Vectors Only S1 | Straight Line Vectors with any Slope S2 | Continuous Curves S3 |
| 1 | $A = \frac{ja}{2}\sum_{n=1}^{N}(-1)^{n-1}z_n^2$ | Horizontal or vertical then clockwise or counterclockwise | Use when all sides of a complex plane closed contour are alternating head to tail horizontal and vertical vectors<br><br>From this basic equation all of the following equations can be derived. | Yes | No | No |
| 2 | $A = -\frac{ja}{2}\sum_{n=1}^{\frac{N}{2}}(z_{2n}^2 - z_{2n-1}^2)$ | Horizontal or vertical then clockwise or counterclockwise | Use when all sides of a complex plane closed contour are alternating head to tail horizontal and vertical vectors | Yes | No | No |

| # | Closed Contour Area Calculation Equation | Initial Direction of Summation/ Integration | Comments | Equation Applicable to Complex Plane Closed Contours Composed of | | |
|---|---|---|---|---|---|---|
| | | | | **Vertical/Horizontal Straight Line Vectors Only S1** | **Straight Line Vectors with any Slope S2** | **Continuous Curves S3** |
| **3** | $A = -a \sum_{n=1}^{\frac{N}{2}} (x_{2n-1}y_{2n-1} - x_{2n}y_{2n})$ <br> or <br> $A = -a \sum_{n=1}^{\frac{N}{2}} \begin{vmatrix} x_{2n-1} & y_{2n} \\ x_{2n} & y_{2n-1} \end{vmatrix}$ | Horizontal or vertical then clockwise or counterclockwise | Use when all sides of a complex plane closed contour are alternating head to tail horizontal and vertical vectors | Yes | No | No |
| **4** | $A = a \sum_{n=1}^{\frac{N}{2}} y_{2n-1}\Delta x_{2n-1}$ | Horizontal only then clockwise or counterclockwise | Use when all sides of a complex plane closed contour are alternating head to tail horizontal and vertical vectors | Yes | No | No |
| **5** | $A = a \sum_{n=1}^{\frac{N}{2}} x_{2n-1}\Delta y_{2n-1}$ | Vertical only then clockwise or counterclockwise | Use when all sides of a complex plane closed contour are alternating head to tail horizontal and vertical vectors | Yes | No | No |

| # | Closed Contour Area Calculation Equation | Initial Direction of Summation/ Integration | Comments | Equation Applicable to Complex Plane Closed Contours Composed of | | |
|---|---|---|---|---|---|---|
| | | | | **Vertical/Horizontal Straight Line Vectors Only S1** | **Straight Line Vectors with any Slope S2** | **Continuous Curves S3** |
| **6** | $A = a \sum_{n=1}^{\frac{N}{2}} y_{2n} \Delta x_{2n-1}$ | Horizontal only then clockwise or counterclockwise | Use when all sides of a complex plane closed contour are alternating head to tail horizontal and vertical vectors | Yes | No | No |
| **7** | $A = a \sum_{n=1}^{\frac{N}{2}} x_{2n} \Delta y_{2n-1}$ | Vertical only then clockwise or counterclockwise | Use when all sides of a complex plane closed contour are alternating head to tail horizontal and vertical vectors | Yes | No | No |
| **8** | $A = -ja \sum_{n=1}^{\frac{N}{2}} \bar{z}_{2n-1} \Delta x_{2n-1}$ | Horizontal only then clockwise or counterclockwise | Use when all sides of a complex plane closed contour are alternating head to tail horizontal and vertical vectors | Yes | No | No |
| **9** | $A = a \sum_{n=1}^{\frac{N}{2}} \bar{z}_{2n-1} \Delta y_{2n-1}$ | Vertical only then clockwise or counterclockwise | Use when all sides of a complex plane closed contour are alternating head to tail horizontal and vertical vectors | Yes | No | No |

| # | Closed Contour Area Calculation Equation | Initial Direction of Summation/ Integration | Comments | Equation Applicable to Complex Plane Closed Contours Composed of | | |
|---|---|---|---|---|---|---|
| | | | | **Vertical/Horizontal Straight Line Vectors Only S1** | **Straight Line Vectors with any Slope S2** | **Continuous Curves S3** |
| **10** | $A = \frac{a}{2}\sum_{n=1}^{\frac{N}{2}}(y_{2n-1}\Delta x_{2n-1} - x_{2n}\Delta y_{2n})$ | Horizontal only then clockwise or counterclockwise | Use when all sides of a complex plane closed contour are alternating head to tail horizontal and vertical vectors | Yes | No | No |
| **11** | $A = \frac{a}{2}\sum_{n=1}^{\frac{N}{2}}(x_{2n-1}\Delta y_{2n-1} - y_{2n}\Delta x_{2n})$ | Vertical only then clockwise or counterclockwise | Use when all sides of a complex plane closed contour are alternating head to tail horizontal and vertical vectors | Yes | No | No |
| **12** | $A = -\frac{ja}{2}\sum_{n=1}^{\frac{N}{2}}[(x_{2n+1}+jy_{2n-1})^2 - (x_{2n-1}+jy_{2n-1})^2]$ | Horizontal only then clockwise or counterclockwise | Use when all sides of a complex plane closed contour are alternating head to tail horizontal and vertical vectors | Yes | No | No |
| **13** | $A = -\frac{ja}{2}\sum_{n=1}^{\frac{N}{2}}[(x_{2n-1}+jy_{2n+1})^2 - (x_{2n-1}+jy_{2n-1})^2]$ | Vertical only then clockwise or counterclockwise | Use when all sides of a complex plane closed contour are alternating head to tail horizontal and vertical vectors | Yes | No | No |

| # | Closed Contour Area Calculation Equation | Initial Direction of Summation/ Integration | Comments | Equation Applicable to Complex Plane Closed Contours Composed of | | |
|---|---|---|---|---|---|---|
| | | | | **Vertical/Horizontal Straight Line Vectors Only S1** | **Straight Line Vectors with any Slope S2** | **Continuous Curves S3** |
| **14** | $A = \frac{jb}{4}\sum_{n=1}^{N}[(x_{n+1}+jy_n)^2 - (x_n+jy_{n+1})^2]$ | Clockwise or counterclockwise | Use when all sides of a complex plane closed contour are head to tail vectors of any slope | Yes | Yes | Can use for approximation by applying to selected points on a continuous closed contour |
| **15** | $A = -b\sum_{n=1}^{N}\bar{y}_n\Delta x_n$ | Clockwise or counterclockwise | Use when all sides of a complex plane closed contour are head to tail vectors of any slope | Yes | Yes | Can use for approximation by applying to selected points on a continuous closed contour |
| **16** | $A = b\sum_{n=1}^{N}\bar{x}_n\Delta y_n$ | Clockwise or counterclockwise | Use when all sides of a complex plane closed contour are head to tail vectors of any slope | Yes | Yes | Can use for approximation by applying to selected points on a continuous closed contour |
| **17** | $A = jb\sum_{n=1}^{N}\bar{z}_n\Delta x_n$ | Clockwise or counterclockwise | Use when all sides of a complex plane closed contour are head to tail vectors of any slope | Yes | Yes | Can use for approximation by applying to selected points on a continuous closed contour |
| **18** | $A = b\sum_{n=1}^{N}\bar{z}_n\Delta y_n$ | Clockwise or counterclockwise | Use when all sides of a complex plane closed contour are head to tail vectors of any slope | Yes | Yes | Can use for approximation by applying to selected points on a continuous closed contour |

| # | Closed Contour Area Calculation Equation | Initial Direction of Summation/ Integration | Comments | Equation Applicable to Complex Plane Closed Contours Composed of | | |
|---|---|---|---|---|---|---|
| | | | | **Vertical/Horizontal Straight Line Vectors Only S1** | **Straight Line Vectors with any Slope S2** | **Continuous Curves S3** |
| **19** | $A = \frac{jb}{2}\sum_{n=1}^{N} \bar{z}_n \Delta z_n^*$ | Clockwise or counterclockwise | Use when all sides of a complex plane closed contour are head to tail vectors of any slope | Yes | Yes | Can use for approximation by applying to selected points on a continous closed contour |
| **20** | $A = -\frac{jb}{2}\sum_{n=1}^{N} \bar{z}_n^* \Delta z_n$ | Clockwise or counterclockwise | Use when all sides of a complex plane closed contour are head to tail vectors of any slope | Yes | Yes | Can use for approximation by applying to selected points on a continuous closed contour |
| **21** | $A = \frac{b}{2}\sum_{n=1}^{N}(x_n y_{n+1} - x_{n+1} y_n)$ or $A = \frac{b}{2}\sum_{n=1}^{N}\begin{vmatrix} x_n & y_n \\ x_{n+1} & y_{n+1}\end{vmatrix}$ | Clockwise or counterclockwise | Use when all sides of a complex plane closed contour are head to tail vectors of any slope | Yes | Yes | Can use for approximation by applying to selected points on a continuous closed contour |
| **22** | $A = \frac{b}{2}\sum_{n=1}^{N}(\bar{x}_n \Delta y_n - \bar{y}_n \Delta x_n)$ | Clockwise or counterclockwise | Use when all sides of a complex plane closed contour are head to tail vectors of any slope | Yes | Yes | Can use for approximation by applying to selected points on a continous closed contour |

| # | Closed Contour Area Calculation Equation | Initial Direction of Summation/ Integration | Comments | Equation Applicable to Complex Plane Closed Contours Composed of | | |
|---|---|---|---|---|---|---|
| | | | | **Vertical/Horizontal Straight Line Vectors Only S1** | **Straight Line Vectors with any Slope S2** | **Continuous Curves S3** |
| **23** | $A = \frac{b}{2}\sum_{n=1}^{N}(x_{n+1}\Delta y_n - y_{n+1}\Delta x_n)$ | Clockwise or counterclockwise | Use when all sides of a complex plane closed contour are head to tail vectors of any slope | Yes | Yes | Can use for approximation by applying to selected points on a continous closed contour |
| **24** | $A = \frac{b}{2}\sum_{n=1}^{N}(x_n\Delta y_n - y_n\Delta x_n)$ | Clockwise or counterclockwise | Use when all sides of a complex plane closed contour are head to tail vectors of any slope<br><br>This is the discrete form of Green's Theorem for the calculation of area | Yes | Yes | Can use for approximation by applying to selected points on a continous closed contour |
| **25** | $A = \frac{b}{2}\oint_{c}[xdy - ydx]$ | Clockwise or counterclockwise | Use when a complex plane closed contour is a continuous closed curve<br><br>This is Green's Theorem for the calculation of area | No | No | Yes |

| # | Closed Contour Area Calculation Equation | Initial Direction of Summation/ Integration | Comments | Equation Applicable to Complex Plane Closed Contours Composed of | | |
|---|---|---|---|---|---|---|
| | | | | **Vertical/Horizontal Straight Line Vectors Only S1** | **Straight Line Vectors with any Slope S2** | **Continuous Curves S3** |
| **26** | $A = -b\oint_c ydx$ | Clockwise or counterclockwise | Use when a complex plane closed contour is a continuous closed curve | No | No | Yes |
| **27** | $A = b\oint_c xdy$ | Clockwise or counterclockwise | Use when a complex plane closed contour is a continuous closed curve | No | No | Yes |
| **28** | $A = jb\oint_c zdx$ | Clockwise or counterclockwise | Use when a complex plane closed contour is a continuous closed curve | No | No | Yes |
| **29** | $A = b\oint_c zdy$ | Clockwise or counterclockwise | Use when a complex plane closed contour is a continuous closed curve | No | No | Yes |
| **30** | $A = \frac{jb}{2}\oint_c zdz^*$ | Clockwise or counterclockwise | Use when a complex plane closed contour is a continuous closed curve | No | No | Yes |

| # | Closed Contour Area Calculation Equation | Initial Direction of Summation/ Integration | Comments | Equation Applicable to Complex Plane Closed Contours Composed of | | |
|---|---|---|---|---|---|---|
| | | | | **Vertical/Horizontal Straight Line Vectors Only S1** | **Straight Line Vectors with any Slope S2** | **Continuous Curves S3** |
| **31** | $A = -\frac{jb}{2}\oint_{c} z^{*}dz$ | Clockwise or counterclockwise | Use when a complex plane closed contour is a continuous closed curve | No | No | Yes |

Variable and Constant Definitions

Summation and integration direction constants a and b

$$a = \begin{cases} +1 & \text{For initial direction horizontal then clockwise} \\ -1 & \text{For initial direction horizontal then counterclockwise} \\ -1 & \text{For initial direction vertical then clockwise} \\ +1 & \text{For initial direction vertical then counterclockwise} \end{cases}$$

$$b = \begin{cases} -1 & \text{For clockwise direction} \\ +1 & \text{For counterclockwise direction} \end{cases}$$

Categorization of Complex Plane Closed Contour Shapes

1. S1 – Discrete complex plane closed contours composed entirely of alternating horizontal and vertical vectors connected head to tail forming a closed loop
2. S2 – Discrete complex plane closed contours composed of vectors of any slope connected head to tail forming a closed loop
3. S3 – Continuous complex plane closed contours composed of infinitesimal vectors of any slope connected head to tail forming a closed loop

Examples of S1, S2, and S3 Complex Plane Closed Contours

S1 S2 S3

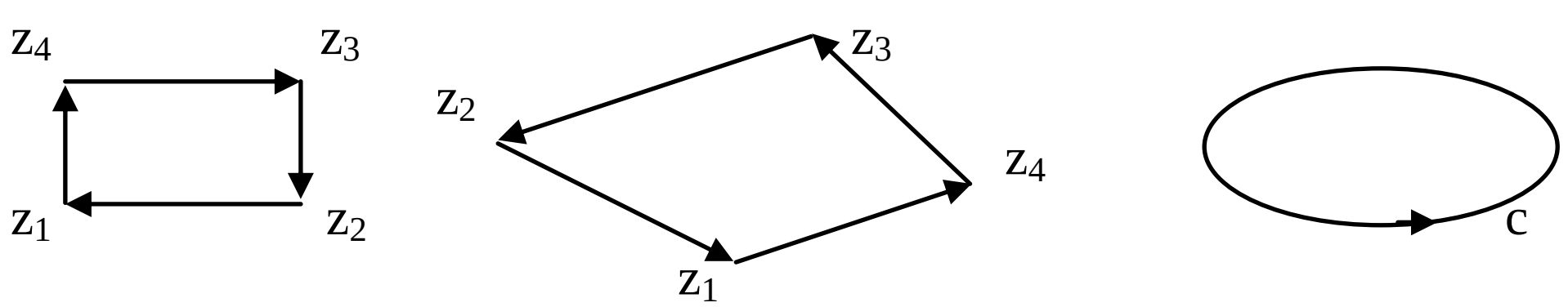

$z_1,z_2,z_3,z_4$ = Closed contour corner points defining the S1 and S2 discrete closed contours
c = The S3 continous closed contour

Clockwise Discrete Closed Contour | Counterclockwise Discrete Closed Contour | Counterclockwise Continuous Closed Contour

Comment 1 - N = the number of corner points that define a discrete closed contour in the complex plane
For S1 closed contours N = 4, 6, 8, 10, …
For S2 closed contours N = 3, 4, 5, 6, …
For S3 closed contours N is infinite

Comment 2 - An S1 and S3 closed contour is a special case of the S2 closed contour.

$$A = \frac{1}{2}(A_H + A_V)$$

where

$A$ = the area enclosed within a complex plane S2 closed contour

$A_H, A_V$ = the areas of two complex plane S1 closed contours constructed from the S2 closed contour corner points

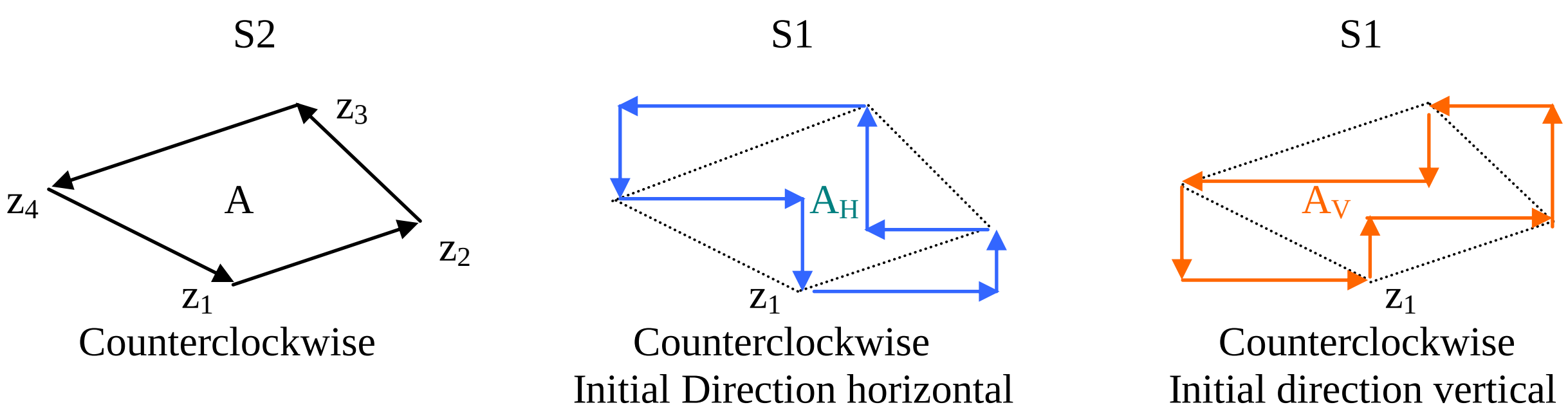

It has been shown that the area enclosed within an S2 closed contour can be calculated from half the sum of the areas enclosed within two S1 closed contours constructed from the same corner points of the S2 closed contour. Note the above complex plane S2 and two S1 closed contours.

Notation and Variable Definitions

$m$ = defined positive integers

$x_m$ = real values

$y_m$ = real values

$z_m = x_m + jy_m$ , the set of complex plane corner points defining a discrete closed contour

Note 1 – $z_m$ can be represented by a discrete function of m. For example, $z_m = e^{\frac{j2\pi(m-1)}{M}}$ where m = 1, 2, 3,…, M

$\bar{x}_m = \frac{x_{m+1} + x_m}{2}$ , $\bar{y}_m = \frac{y_{m+1} + y_m}{2}$ , $\bar{z}_m = \frac{z_{m+1} + z_m}{2}$

$\Delta x_m = x_{m+1} - x_m$ , $\Delta y_m = y_{m+1} - y_m$ , $\Delta z_m = z_{m+1} - z_m$

$x_{N+1} = x_1$ , $y_{N+1} = y_1$ , $z_{N+1} = z_1$

$x(\theta)$ = A continuous real value function of $\theta$

$y(\theta)$ = A continuous real value function of $\theta$

$z(\theta) = x(\theta) + jy(\theta)$ , A continuous complex value function of $\theta$ defining a closed contour in the complex plane (ex: $z(\theta) = e^{j\theta}$)

A = The area enclosed within an S1, S2, or S3 closed contour in the complex plane. The area, A, is always a positive value

$A_H$ = The area enclosed within the S1 closed contour constructed from an S2 closed contour by starting at an initial S2 corner point, $z_1$, and then drawing alternately horizontal then vertical vectors connecting the S2 corner points. The two S1 closed contours contructed from the same S2 corner points must have the same sense, clockwise or counterclockwise.

$A_V$ = The area enclosed within the S1 closed contour constructed from an S2 closed contour by starting at an initial S2 corner point, $z_1$, and then drawing alternately vertical then horizontal vectors connecting the S2 corner points. The two S1 closed contours contructed from the same S2 corner points must have the same sense, clockwise or counterclockwise.

N = The number of corner points on the complex plane closed contour (i.e. The number of corner points defining the complex plane closed contour).

* = complex conjugate designation (ex: $[a+jb]^* = a-jb$ or $[e^{j\theta}]^* = e^{-j\theta}$ )

$\oint_c$ = The integration sign representing complex plane integration along the continuous closed contour, c

<u>Note 2</u> - Consecutive contour points, designated $z_1$, $z_2$, $z_3$, …, are connected by straight line vectors. These vectors are connected head to tail to form a closed loop.

Comment 3 - There are several interesting closed contour summations and integrals that equal zero. They are as follows:

For S1 complex plane closed contours

1. $\sum_{n=1}^{\frac{N}{2}} [\,(x_{2n}{}^2 - x_{2n-1}{}^2) = 0$ , $x_N = x_1$ , $x_{2n+1} = x_{2n}$ , $y_{2n} = y_{2n-1}$ , Initial direction of summation = horizontal

2. $\sum_{n=1}^{\frac{N}{2}} [\,(y_{2n}{}^2 - y_{2n-1}{}^2) = 0$ , $y_N = y_1$ , $y_{2n+1} = y_{2n}$ , $x_{2n} = x_{2n-1}$ , Initial direction of summation = vertical

For S2 complex plane closed contours

3. $\sum_{n=1}^{N} (x_{n+1}{}^2 - x_n{}^2) = 0$ , $x_{N+1} = x_1$

4. $\sum_{n=1}^{N} (y_{n+1}{}^2 - y_n{}^2) = 0$ , $y_{N+1} = y_1$

5. $\sum_{n=1}^{N} \overline{x}_n \Delta x_n = 0$

6. $\sum_{n=1}^{N} \overline{y}_n \Delta y_n = 0$

7. $\sum_{n=1}^{N} \overline{z}_n \Delta z_n = 0$

8. $\sum_{n=1}^{N} \bar{z}_n^{\,*} \Delta z_n^{\,*} = 0$

For S3 complex plane closed contours

9. $\oint_c x(\theta)dx(\theta) = 0$

10. $\oint_c y(\theta)dy(\theta) = 0$

11. $\oint_c z(\theta)dz(\theta) = 0$

12. $\oint_c z(\theta)^* dz(\theta)^* = 0$

<u>Comment 4</u> - The closed contour summations previously listed have another form, a discrete integral form. The discrete integral form is particularly useful when a function is available to describe the closed contour for which the enclosed area is to be calculated.

For example: $A = \frac{jb}{2}\sum_{n=1}^{N} \bar{z}_n \Delta z_n^{\,*} = \frac{jb}{2}\,{}_1\!\int_{1}^{N+1} \bar{z}_n \Delta z_n^{\,*}$ .

Below is a listing of the thirty-one complex plane closed contour area calculation equations which are presented in Table 3.6-1. These equations are derived in the following pages in the order shown.

**Listing of the Derived Complex Plane Closed Contour Area Calculation Equations**

Equations to calculate the area of S1 complex plane closed contours
(closed contours composed only of horizontal and vertical vectors)

1. $\mathbf{A = \frac{ja}{2} \sum_{n=1}^{N} (-1)^{n-1} z_n^2}$ , This basic equation is derived in Sections 3.1 thru 3.5.

2. $\mathbf{A = -\frac{ja}{2} \sum_{n=1}^{\frac{N}{2}} ( z_{2n}^2 - z_{2n-1}^2)}$

3. $\mathbf{A = a \sum_{n=1}^{\frac{N}{2}} y_{2n-1} \Delta x_{2n-1}}$ , where the inital direction of summation is horizontal.

4. $\mathbf{A = a \sum_{n=1}^{\frac{N}{2}} y_{2n} \Delta x_{2n-1}}$ , where the inital direction of summation is horizontal.

5. $\mathbf{A = a \sum_{n=1}^{\frac{N}{2}} x_{2n-1} \Delta y_{2n-1}}$ , where the initial direction of summation is vertical.

6. $\mathbf{A = a \sum_{n=1}^{\frac{N}{2}} x_{2n} \Delta y_{2n-1}}$ , where the initial direction of summation is vertical.

7. $\mathbf{A = -a \sum_{n=1}^{\frac{N}{2}} (x_{2n-1} y_{2n-1} - x_{2n} y_{2n})}$, where the initial direction of summation is horizontal or vertical.

8. $\mathbf{A = -ja \sum_{n=1}^{\frac{N}{2}} \bar{z}_{2n-1} \Delta x_{2n-1}}$ , where the initial direction of summation is horizontal.

9. $\mathbf{A = a \sum_{n=1}^{\frac{N}{2}} \overline{z}_{2n-1} \Delta y_{2n-1}}$ , where the initial direction of summation is vertical.

10. $\mathbf{A = a \sum_{n=1}^{\frac{N}{2}} (y_{2n-1} \Delta x_{2n-1} - x_{2n} \Delta y_{2n})}$, where the initial direction of summation is horizontal.

11. $\mathbf{A = a \sum_{n=1}^{\frac{N}{2}} (x_{2n-1} \Delta y_{2n-1} - y_{2n} \Delta x_{2n})}$, where the initial direction of summation is vertical.

12. $\mathbf{A = - \frac{ja}{2} \sum_{n=1}^{\frac{N}{2}} [(x_{2n+1} + jy_{2n-1})^2 - (x_{2n-1} + jy_{2n-1})^2]}$, where the initial direction of summation is horizontal.

13. $\mathbf{A = - \frac{ja}{2} \sum_{n=1}^{\frac{N}{2}} [(x_{2n-1} + jy_{2n+1})^2 - (x_{2n-1} + jy_{2n-1})^2]}$, where the initial direction of summation is vertical.

<u>Equations to calculate the area of S2 complex plane closed contours</u>
<u>(closed contours composed of vectors of any slope)</u>

<u>Note</u> – S1 complex plane closed contours are a special case of S2 complex plane closed contours

14. $\mathbf{A = \frac{b}{2} \sum_{n=1}^{N} (x_n \Delta y_n - y_n \Delta x_n)}$

15. $\mathbf{A = \frac{b}{2} \sum_{n=1}^{N} (x_n y_{n+1} - x_{n+1} y_n)}$

16. $\mathbf{A = \frac{jb}{4} \sum_{n=1}^{N} [(x_{n+1} + jy_n)^2 - (x_n + jy_{n+1})^2]}$

17. $\mathbf{A = \frac{jb}{2} \sum_{n=1}^{N} \overline{z}_n \Delta z_n^{*}}$

18. $\mathbf{A = -\frac{jb}{2}\sum_{n=1}^{N} \bar{z}_n^{*}\Delta z_n}$

19. $\mathbf{A = jb\sum_{n=1}^{N} \bar{z}_n \Delta x_n}$

20. $\mathbf{A = -b\sum_{n=1}^{N} \bar{y}_n \Delta x_n}$

21. $\mathbf{A = b\sum_{n=1}^{N} \bar{z}_n \Delta y_n}$

22. $\mathbf{A = b\sum_{n=1}^{N} \bar{x}_n \Delta y_n}$

23. $\mathbf{A = \frac{b}{2}\sum_{n=1}^{N} (\bar{x}_n \Delta y_n - \bar{y}_n \Delta x_n)}$

24. $\mathbf{A = \frac{b}{2}\sum_{n=1}^{N} (x_{n+1}\Delta y_n - y_{n+1}\Delta x_n)}$

<u>Equations to calculate the area of S3 complex plane closed contours</u>
<u>(closed contours composed of infinitesimal vectors of any slope)</u>

<u>Note</u> – S3 complex plane closed contours are continuous curve closed contours

25. $\mathbf{A = \frac{b}{2}\oint_{c} [xdy - ydx]}$, Green's area calculation equation

26. $\mathbf{A = -b\oint_{c} ydx}$

27. $\mathbf{A = b\oint_{c} xdy}$

28. $\mathbf{A = jb \oint_c z dx}$

29. $\mathbf{A = b \oint_c z dy}$

30. $\mathbf{A = \frac{jb}{2} \oint_c z dz^*}$

31. $\mathbf{A = - \frac{jb}{2} \oint_c z^* dz}$

The derivations of the above thirty-one complex plane closed contour area calculation equations presented in Table 3.6-1 and listed above are shown on the following pages.

## Derivation of thirty-one equations to calculate the area of complex plane closed contours

For the derivations which follow, Eq 3.5-52 is rewritten below with several notational modifications for improved clarity and applicability.

**Fundamental discrete complex plane contour area equation**

$$\mathbf{1.}\quad A = \frac{ja}{2}\sum_{n=1}^{N}(-1)^{n-1}z_n^{\,2} \qquad (3.6\text{-}1)$$

**where**

**The initial direction of summation is horizontal or vertical then clockwise or counterclockwise**

**A = The area enclosed within a complex plane closed contour of shape S1**
**The closed contour sides are composed of alternating horizontal and vertical head to tail connected vectors**

$$a = \begin{cases} +1 & \text{for clockwise contour vector summation with initial vector horizontal} \\ -1 & \text{for counterclockwise contour vector summation with initial vector horizontal} \\ -1 & \text{for clockwise contour vector summation with initial vector vertical} \\ +1 & \text{for counterclockwise contour vector summation with initial vector vertical} \end{cases}$$

**N = 4,6,8,10,…**
**N = the number of the discrete complex plane contour corner points, an even number**
**$z_n$ = the coordinates of the corner points of the discrete complex plane contour**
**n = 1,2,3,4,…,N**

<u>**Comments**</u> **- 1. The value of the summation, $\frac{j}{2}\sum_{n=1}^{N}(-1)^{n-1}z_n^{\,2}$ , is always real.**

**2. The value of the constant, a (+1 or –1), is such that the area, A, is always positive.**

**3. Eq 3.6-1 is directly applicable to those complex plane closed contours previously described as being of the shape, S1.**

<u>Derive equations to calculate the area within S1 closed contours in the complex plane.</u>

Consider the following diagram of an S1 closed contour in the complex plane
In this case, the initial vector direction from the point, $z_1$, is horizontal

<u>Diagram 3.6-1 Several points on an S1 Closed Contour where the initial direction of summation is horizontal</u>

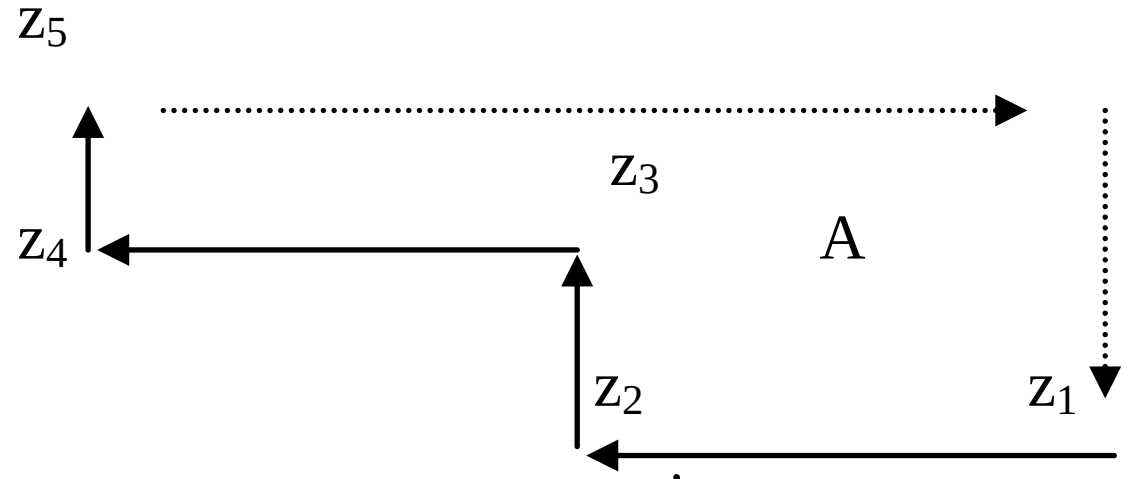

Complex Plane

$z_n = x_n + jy_n$ for n = 1 ,2, 3, …, N
N = number of corner points in the S1 closed contour
A = area within the closed contour

$z_{2n} = x_{2n}+jy_{2n-1}$ $\quad z_{2n-1} = x_{2n-1}+jy_{2n-1}$ $\quad n = 1, 2, 3,\ldots, \frac{N}{2}$

Change the form of Eq 3.6-1

**2.** $$\mathbf{A = -\frac{ja}{2}\sum_{n=1}^{\frac{N}{2}}(z_{2n}^2 - z_{2n-1}^2)} \qquad (3.6\text{-}2)$$

**where**

**The initial direction of summation is horizontal or vertical then clockwise or counterclockwise**

**A = The area enclosed within a complex plane closed contour of shape S1**
**The closed contour sides are composed of alternating horizontal and vertical head to tail connected vectors**

$$\mathbf{a} = \begin{cases} +1 & \text{for clockwise contour vector summation with initial vector horizontal} \\ -1 & \text{for counterclockwise contour vector summation with initial vector horizontal} \\ -1 & \text{for clockwise contour vector summation with initial vector vertical} \\ +1 & \text{for counterclockwise contour vector summation with initial vector vertical} \end{cases}$$

**N = 4,6,8,10,…**

**N = the number of the discrete complex plane contour corner points, an even number**

**$z_n$ = the coordinates of the corner points of the discrete complex plane contour**

**n = 1,2,3,4,…,N**

**3. Derive the equation,** $\mathbf{A = a\sum_{n=1}^{\frac{N}{2}} y_{2n-1}\Delta x_{2n-1}}$ , where the inital direction of summation is horizontal.

Noting Diagram 3.6-1

$$z_{2n} = x_{2n}+jy_{2n-1} \qquad (3.6\text{-}3)$$

$$z_{2n-1} = x_{2n-1}+jy_{2n-1} \qquad (3.6\text{-}4)$$

Substituting Eq 3.6-3 and Eq 3.6-4 into Eq 3.6-2

$$A = -\frac{ja}{2}\sum_{n=1}^{\frac{N}{2}}[\,(x_{2n}+jy_{2n-1})^2 - (x_{2n-1}+jy_{2n-1})^2\,] \qquad (3.6\text{-}5)$$

Simplifying Eq 3.6-5

$$A = -\frac{ja}{2}\sum_{n=1}^{\frac{N}{2}}[\,x_{2n}^2 + 2jx_{2n}y_{2n-1} - y_{2n-1}^2 - x_{2n-1}^2 - 2jx_{2n-1}y_{2n-1} + y_{2n-1}^2\,] \qquad (3.6\text{-}6)$$

$$A = -\frac{ja}{2}\sum_{n=1}^{\frac{N}{2}}[(x_{2n}^2 - x_{2n-1}^2) + a\sum_{n=1}^{\frac{N}{2}}(x_{2n} - x_{2n-1})\, y_{2n-1} \qquad (3.6\text{-}7)$$

$$\sum_{n=1}^{\frac{N}{2}}[(x_{2n}^2 - x_{2n-1}^2) = x_2^2 - x_1^2 + x_4^2 - x_3^2 + x_6^2 - x_5^2 + \ldots + x_N^2 - x_{N-1}^2 \qquad (3.6\text{-}8)$$

Noting Diagram 3.6-1

$$x_{2n+1} = x_{2n} \qquad (3.6\text{-}9)$$

$$x_N = x_1 \qquad (3.6\text{-}10)$$

$$y_{2n} = y_{2n-1} \qquad (3.6\text{-}11)$$

From Eq 3.6-8 thru Eq 3.6-10

$$\sum_{n=1}^{\frac{N}{2}}[(x_{2n}^2 - x_{2n-1}^2) = x_3^2 - x_N^2 + x_5^2 - x_3^2 + x_{N-1}^2 - x_5^2 + \ldots + x_N^2 - x_{N-1}^2 = 0 \qquad (3.6\text{-}12)$$

$$\sum_{n=1}^{\frac{N}{2}}[(x_{2n}^2 - x_{2n-1}^2) = 0 \qquad (3.6\text{-}13)$$

Substituting Eq 3.6-13 into Eq 3.6-7

$$A = a\sum_{n=1}^{\frac{N}{2}} y_{2n-1}(x_{2n} - x_{2n-1}) = a\sum_{n=1}^{\frac{N}{2}} y_{2n-1}\,\Delta x_{2n-1} \qquad (3.6\text{-}14)$$

$$A = a\sum_{n=1}^{\frac{N}{2}} y_{2n-1}\,\Delta x_{2n-1} \qquad (3.6\text{-}15)$$

Then from Eq 3.6-15

$$\mathbf{A = a \sum_{n=1}^{\frac{N}{2}} y_{2n-1} \Delta x_{2n-1}} \qquad (3.6\text{-}16)$$

**where**

**The initial direction of summation is horizontal then clockwise or counterclockwise**

**A = The area enclosed within a complex plane closed contour of shape S1**
**The closed contour sides are composed of alternating horizontal and vertical head to tail connected vectors**

$$\mathbf{a} = \begin{cases} +1 & \textbf{for clockwise contour vector summation with initial vector horizontal} \\ -1 & \textbf{for counterclockwise contour vector summation with initial vector horizontal} \\ -1 & \textbf{for clockwise contour vector summation with initial vector vertical} \\ +1 & \textbf{for counterclockwise contour vector summation with initial vector vertical} \end{cases}$$

**N = 4,6,8,10,…**

**N = the number of the discrete complex plane contour corner points, an even number**

**$z_n = x_n + jy_n$ , the coordinates of the corner points of the discrete complex plane contour**

**$\Delta x_{2n-1} = x_{2n} - x_{2n-1}$**

**4. Derive the equation, $A = a \sum_{n=1}^{\frac{N}{2}} y_{2n} \Delta x_{2n-1}$** , where the inital direction of summation is horizontal.

From Eq 3.6-16 and Eq 3.6-11

$$A = a \sum_{n=1}^{\frac{N}{2}} y_{2n} \Delta x_{2n-1} \qquad (3.6\text{-}17)$$

Then

$$\mathbf{A = a \sum_{n=1}^{\frac{N}{2}} y_{2n} \Delta x_{2n-1}} \qquad (3.6\text{-}18)$$

**where**

**The initial direction of summation is horizontal then clockwise or counterclockwise**

**A = The area enclosed within a complex plane closed contour of shape S1**
**The closed contour sides are composed of alternating horizontal and vertical head to tail connected vectors**

$$\mathbf{a} = \begin{cases} +1 & \text{for clockwise contour vector summation with initial vector horizontal} \\ -1 & \text{for counterclockwise contour vector summation with initial vector horizontal} \\ -1 & \text{for clockwise contour vector summation with initial vector vertical} \\ +1 & \text{for counterclockwise contour vector summation with initial vector vertical} \end{cases}$$

**N = 4,6,8,10,…**
**N = the number of the discrete complex plane contour corner points, an even number**
**$z_n = x_n + jy_n$ , the coordinates of the corner points of the discrete complex plane contour**
**$\Delta x_{2n-1} = x_{2n} - x_{2n-1}$**

Diagram 3.6-2 Several points on an S1 Closed Contour where the initial direction of summation is vertical

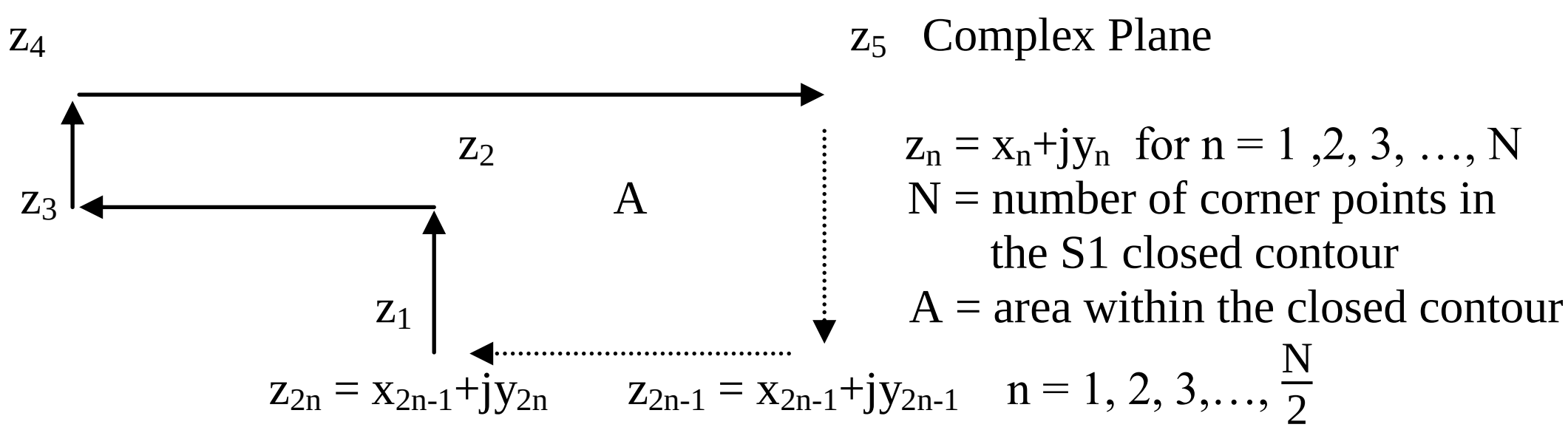

**5. Derive the equation, $A = a \sum_{n=1}^{\frac{N}{2}} x_{2n-1} \Delta y_{2n-1}$** , where the initial direction of summation is vertical.

Noting Diagram 3.6-2

$$z_{2n} = x_{2n-1} + jy_{2n} \tag{3.6-19}$$

$$z_{2n-1} = x_{2n-1} + jy_{2n-1} \tag{3.6-20}$$

Substituting Eq 3.6-19 and Eq 3.6-20 into Eq 3.6-2

$$A = -\frac{ja}{2} \sum_{n=1}^{\frac{N}{2}} [\,(x_{2n-1} + jy_{2n})^2 - (x_{2n-1} + jy_{2n-1})^2\,] \tag{3.6-21}$$

Simplifying Eq 3.6-21

$$A = -\frac{ja}{2} \sum_{n=1}^{\frac{N}{2}} [\,x_{2n-1}^2 + 2jx_{2n-1}y_{2n} - y_{2n}^2 - x_{2n-1}^2 - 2jx_{2n-1}y_{2n-1} + y_{2n-1}^2\,] \tag{3.6-22}$$

$$A = \frac{ja}{2} \sum_{n=1}^{\frac{N}{2}} [\,(y_{2n}{}^2 - y_{2n-1}{}^2) + a \sum_{n=1}^{\frac{N}{2}} (y_{2n} - y_{2n-1})\, x_{2n-1} \qquad (3.6\text{-}23)$$

$$\sum_{n=1}^{\frac{N}{2}} [\,(y_{2n}{}^2 - y_{2n-1}{}^2) \;= y_2{}^2 - y_1{}^2 + y_4{}^2 - y_3{}^2 + \; y_6{}^2 - y_5{}^2 + \ldots + \; y_N{}^2 - y_{N-1}{}^2 \qquad (3.6\text{-}24)$$

Noting Diagram 3.6-2

$$y_{2n+1} = y_{2n} \qquad (3.6\text{-}25)$$

$$y_N = y_1 \qquad (3.6\text{-}26)$$

$$x_{2n} = x_{2n-1} \qquad (3.6\text{-}27)$$

From Eq 3.6-24 thru Eq 3.6-26

$$\sum_{n=1}^{\frac{N}{2}} [\,(y_{2n}{}^2 - y_{2n-1}{}^2) \;= y_3{}^2 - y_N{}^2 + y_5{}^2 - y_3{}^2 + \; y_{N-1}{}^2 - y_5{}^2 + \ldots + \; y_N{}^2 - y_{N-1}{}^2 = 0 \qquad (3.6\text{-}28)$$

$$\sum_{n=1}^{\frac{N}{2}} [\,(y_{2n}{}^2 - y_{2n-1}{}^2) = 0 \qquad (3.6\text{-}29)$$

Substituting Eq 3.6-29 into Eq 3.6-23

$$A = a \sum_{n=1}^{\frac{N}{2}} x_{2n-1}\,(y_{2n} - y_{2n-1}) \;= a \sum_{n=1}^{\frac{N}{2}} x_{2n-1}\,\Delta y_{2n-1} \qquad (3.6\text{-}30)$$

$$A = a \sum_{n=1}^{\frac{N}{2}} x_{2n-1}\,\Delta y_{2n-1} \qquad (3.6\text{-}31)$$

Then from Eq 3.6-31

$$\mathbf{A = a \sum_{n=1}^{\frac{N}{2}} x_{2n-1} \Delta y_{2n-1}} \qquad (3.6\text{-}32)$$

**where**

**The initial direction of summation is vertical then clockwise or counterclockwise**

**A = The area enclosed within a complex plane closed contour of shape S1**
**The closed contour sides are composed of alternating horizontal and vertical head to tail connected vectors**

$$\mathbf{a} = \begin{cases} +1 & \textbf{for clockwise contour vector summation with initial vector horizontal} \\ -1 & \textbf{for counterclockwise contour vector summation with initial vector horizontal} \\ -1 & \textbf{for clockwise contour vector summation with initial vector vertical} \\ +1 & \textbf{for counterclockwise contour vector summation with initial vector vertical} \end{cases}$$

**N = 4,6,8,10,…**

**N = the number of the discrete complex plane contour corner points, an even number**

$\mathbf{z_n = x_n + jy_n}$ **, the coordinates of the corner points of the discrete complex plane contour**

$\mathbf{\Delta y_{2n-1} = y_{2n} - y_{2n-1}}$

**6. Derive the equation,** $\mathbf{A = a \sum_{n=1}^{\frac{N}{2}} x_{2n} \Delta y_{2n-1}}$ , where the initial direction of summation is vertical.

From Eq 3.6-27 and Eq 3.6-32

$$A = a \sum_{n=1}^{\frac{N}{2}} x_{2n} \Delta y_{2n-1} \qquad (3.6\text{-}33)$$

Then

$$\mathbf{A = a \sum_{n=1}^{\frac{N}{2}} x_{2n} \Delta y_{2n-1}} \qquad (3.6\text{-}34)$$

**where**

**The initial direction of summation is vertical then clockwise or counterclockwise**

**A = The area enclosed within a complex plane closed contour of shape S1**
**The closed contour sides are composed of alternating horizontal and vertical head to tail connected vectors**

$$\mathbf{a} = \begin{cases} +1 & \textbf{for clockwise contour vector summation with initial vector horizontal} \\ -1 & \textbf{for counterclockwise contour vector summation with initial vector horizontal} \\ -1 & \textbf{for clockwise contour vector summation with initial vector vertical} \\ +1 & \textbf{for counterclockwise contour vector summation with initial vector vertical} \end{cases}$$

**N = 4,6,8,10,…**

**N = the number of the discrete complex plane contour corner points, an even number**
**$z_n = x_n + jy_n$ , the coordinates of the corner points of the discrete complex plane contour**
**$\Delta y_{2n-1} = y_{2n} - y_{2n-1}$**

**7. Derive the equation, $A = -a \sum_{n=1}^{\frac{N}{2}} (x_{2n-1}y_{2n-1} - x_{2n}y_{2n})$,** where the initial direction of summation is horizontal or vertical.

$$z_{2n} = x_{2n} + jy_{2n} \tag{3.6-35}$$

$$z_{2n-1} = x_{2n-1} + jy_{2n-1} \tag{3.6-36}$$

Substituting Eq 3.6-35 and Eq 3.6-36 into Eq 3.6-2

$$A = -\frac{ja}{2} \sum_{n=1}^{\frac{N}{2}} [\, x_{2n}^2 + 2jx_{2n}y_{2n} - y_{2n}^2 - x_{2n-1}^2 - 2jx_{2n-1}y_{2n-1} + y_{2n-1}^2 \,] \tag{3.6-37}$$

Simplifying Eq 3.6-37

$$A = -\frac{ja}{2} \sum_{n=1}^{\frac{N}{2}} [\, (x_{2n}^2 - x_{2n-1}^2) - (y_{2n}^2 - y_{2n-1}^2) + 2j(x_{2n}y_{2n} - x_{2n-1}y_{2n-1}) \,] \tag{3.6-38}$$

Then from Eq 3.6-38, Eq 3.6-13 and Eq 3.6-29

$$\mathbf{A = -a \sum_{n=1}^{\frac{N}{2}} [\, (x_{2n-1}y_{2n-1} - x_{2n}y_{2n}) \,]} \tag{3.6-39}$$

**or**

$$\mathbf{A = -a \sum_{n=1}^{\frac{N}{2}} \begin{vmatrix} x_{2n-1} & y_{2n} \\ x_{2n} & y_{2n-1} \end{vmatrix}} \tag{3.6-40}$$

**where**

**The initial direction of summation is horizontal or vertical then clockwise or counterclockwise**
**A = The area enclosed within a complex plane closed contour of shape S1**
**The closed contour sides are composed of alternating horizontal and vertical head to tail connected vectors**

**a =**
- **+1 for clockwise contour vector summation with initial vector horizontal**
- **–1 for counterclockwise contour vector summation with initial vector horizontal**
- **–1 for clockwise contour vector summation with initial vector vertical**
- **+1 for counterclockwise contour vector summation with initial vector vertical**

**N = 4,6,8,10,…**

**N = the number of the discrete complex plane contour corner points, an even number**
**$z_n = x_n+jy_n$ , the coordinates of the corner points of the discrete complex plane contour**
**n = 1,2,3,4,…,N**

**8. Derive the equation, $A = -ja\sum_{n=1}^{\frac{N}{2}} \bar{z}_{2n-1}\Delta x_{2n-1}$** , where the initial direction of summation is horizontal.

Rewriting Eq 3.6-5

$$A = -\frac{ja}{2}\sum_{n=1}^{\frac{N}{2}}[(x_{2n}+jy_{2n-1})^2 - (x_{2n-1}+jy_{2n-1})^2] \quad (3.6\text{-}41)$$

Simplifying Eq 3.6-41

$$A = -\frac{ja}{2}\sum_{n=1}^{\frac{N}{2}}[x_{2n}^2 + 2jx_{2n}y_{2n-1} - y_{2n-1}^2 - x_{2n-1}^2 - 2jx_{2n-1}y_{2n-1} + y_{2n-1}^2] \quad (3.6\text{-}42)$$

$$A = -\frac{ja}{2}\sum_{n=1}^{\frac{N}{2}}[(x_{2n}^2 - x_{2n-1}^2) + 2j(x_{2n} - x_{2n-1})\, y_{2n-1} \quad (3.6\text{-}43)$$

$$A = -\frac{ja}{2}\sum_{n=1}^{\frac{N}{2}}[(x_{2n} - x_{2n-1})[\, x_{2n} + jy_{2n-1} + x_{2n-1} + jy_{2n-1}] \quad (3.6\text{-}44)$$

From Eq 3.6-11

$$y_{2n} = y_{2n-1} \quad (3.6\text{-}45)$$

Substituting Eq 3.6-45 into Eq 3.6-44

$$A = -\frac{ja}{2}\sum_{n=1}^{\frac{N}{2}}[(x_{2n} - x_{2n-1})[\,(x_{2n} + jy_{2n}) + (x_{2n-1} + jy_{2n-1})] = -ja\sum_{n=1}^{\frac{N}{2}}(x_{2n} - x_{2n-1})\left[\frac{z_{2n} + z_{2n-1}}{2}\right] \quad (3.6\text{-}46)$$

$$A = -ja\sum_{n=1}^{\frac{N}{2}} \bar{z}_{2n-1}\Delta x_{2n-1} \quad (3.6\text{-}47)$$

Then

$$\mathbf{A = - ja \sum_{n=1}^{\frac{N}{2}} \bar{z}_{2n-1}\Delta x_{2n-1}} \qquad (3.6\text{-}48)$$

**where**

**The initial direction of summation is horizontal then clockwise or counterclockwise**

**A = The area enclosed within a complex plane closed contour of shape S1**
**The closed contour sides are composed of alternating horizontal and vertical head to tail connected vectors**

$$\mathbf{a} = \begin{cases} +1 & \textbf{for clockwise contour vector summation with initial vector horizontal} \\ -1 & \textbf{for counterclockwise contour vector summation with initial vector horizontal} \\ -1 & \textbf{for clockwise contour vector summation with initial vector vertical} \\ +1 & \textbf{for counterclockwise contour vector summation with initial vector vertical} \end{cases}$$

**N = 4,6,8,10,…**

$\mathbf{z_{n-1} = x_{2n-1}+jy_{2n-1}}$ **, the coordinates of the corner points of the discrete complex plane contour**

$$\mathbf{\bar{z}_{2n-1} = \frac{z_{2n} + z_{2n-1}}{2}}$$

$$\mathbf{n = 1,2,3,4,\ldots,\frac{N}{2}}$$

**9. Derive the equation,** $\mathbf{A = a \sum_{n=1}^{\frac{N}{2}} \bar{z}_{2n-1}\Delta y_{2n-1}}$ , where the initial direction of summation is vertical.

Rewriting Eq 3.6-21

$$A = -\frac{ja}{2}\sum_{n=1}^{\frac{N}{2}}[(x_{2n-1}+jy_{2n})^2 - (x_{2n-1}+jy_{2n-1})^2] \qquad (3.6\text{-}49)$$

Simplifying Eq 3.6-49

$$A = -\frac{ja}{2}\sum_{n=1}^{\frac{N}{2}}[x_{2n-1}^2 + 2jx_{2n-1}y_{2n} - y_{2n}^2 - x_{2n-1}^2 - 2jx_{2n-1}y_{2n-1} + y_{2n-1}^2] \qquad (3.6\text{-}50)$$

$$A = -\frac{ja}{2}\sum_{n=1}^{\frac{N}{2}}[-(y_{2n}^2 - y_{2n-1}^2) + 2j(y_{2n} - y_{2n-1})\,x_{2n-1} \qquad (3.6\text{-}51)$$

$$A = \frac{ja}{2}\sum_{n=1}^{\frac{N}{2}}[(y_{2n} - y_{2n-1})[\,y_{2n} + y_{2n-1} - jx_{2n-1} - jx_{2n-1}\,] \tag{3.6-52}$$

$$A = \frac{a}{2}\sum_{n=1}^{\frac{N}{2}}[(y_{2n} - y_{2n-1})[\,jy_{2n} + jy_{2n-1} + x_{2n-1} + x_{2n-1}\,] \tag{3.6-53}$$

From Eq 3.6-27

$$x_{2n} = x_{2n-1} \tag{3.6-54}$$

From Eq 3.6-53 and Eq 3.6-54

$$A = \frac{a}{2}\sum_{n=1}^{\frac{N}{2}}(y_{2n} - y_{2n-1})[\,(x_{2n} + jy_{2n}) + (x_{2n-1} + jy_{2n-1})\,] = a\sum_{n=1}^{\frac{N}{2}}(y_{2n} - y_{2n-1})[\,\frac{z_{2n} + z_{2n-1}}{2}\,] \tag{3.6-55}$$

$$A = a\sum_{n=1}^{\frac{N}{2}}\bar{z}_{2n-1}\Delta y_{2n-1} \tag{3.6-56}$$

Then

$$\mathbf{A = a\sum_{n=1}^{\frac{N}{2}}\bar{z}_{2n-1}\Delta y_{2n-1}} \tag{3.6-57}$$

**where**

**The initial direction of summation is vertical then clockwise or counterclockwise**

**A = The area enclosed within a complex plane closed contour of shape S1**
**The closed contour sides are composed of alternating horizontal and vertical head to tail connected vectors**

**a =**
- **+1 for clockwise contour vector summation with initial vector horizontal**
- **–1 for counterclockwise contour vector summation with initial vector horizontal**
- **–1 for clockwise contour vector summation with initial vector vertical**
- **+1 for counterclockwise contour vector summation with initial vector vertical**

**N = 4,6,8,10,…**

**$z_{n-1} = x_{2n-1}+jy_{2n-1}$, the coordinates of the corner points of the discrete complex plane contour**

**$\bar{z}_{2n-1} = \frac{z_{2n} + z_{2n-1}}{2}$**

**$n = 1,2,3,4,\ldots,\frac{N}{2}$**

**10. Derive the equation, $A = a \sum_{n=1}^{\frac{N}{2}} (y_{2n-1}\Delta x_{2n-1} - x_{2n}\Delta y_{2n})$,** where the initial direction of summation is horizontal.

Refer to Diagram 3.6-1

From Eq 3.6-16

$$A = a \sum_{n=1}^{\frac{N}{2}} y_{2n-1}\Delta x_{2n-1} \quad \text{where the initial direction of summation is horizontal} \qquad (3.6\text{-}58)$$

From Eq 3.6-31

$$A = a \sum_{n=1}^{\frac{N}{2}} x_{2n-1}\Delta y_{2n-1} \quad \text{where the initial direction of summation is vertical} \qquad (3.6\text{-}59)$$

where

$$a = \begin{cases} +1 & \text{for clockwise contour vector summation with initial vector horizontal} \\ -1 & \text{for counterclockwise contour vector summation with initial vector horizontal} \\ -1 & \text{for clockwise contour vector summation with initial vector vertical} \\ +1 & \text{for counterclockwise contour vector summation with initial vector vertical} \end{cases}$$

Applying Eq 3.6-58 and Eq 3.6-59 to Diagram 3.1

$$A = \frac{1}{2}\left[\, a \sum_{n=1}^{\frac{N}{2}} y_{2n-1}\Delta x_{2n-1} + (-a) \sum_{n=1}^{\frac{N}{2}} x_{2n}\Delta y_{2n} \,\right] \qquad (3.6\text{-}60)$$

where the initial direction of summation is horizontal

The value of the constant, a, in Eq 3.6-57 and Eq 3.6-58 depends on the listing of the four summation conditions specified on the previous page. The value of the constant, a, will be either +1 or –1. In Eq 3.6-59, for the same values of n, the first summation has an initial horizontal direction while the second summation has an initial vertical direction. Then, since both summations of Eq 3.6-59 must have the same circular direction, clockwise or counterclockwise, the constant of the second summation becomes –a.

Simplifying Eq 3.6-60

$$A = \frac{a}{2}\sum_{n=1}^{\frac{N}{2}} (y_{2n-1}\Delta x_{2n-1} - x_{2n}\Delta y_{2n}) \qquad (3.6\text{-}61)$$

Then

$$\mathbf{A = \frac{a}{2}\sum_{n=1}^{\frac{N}{2}} (y_{2n-1}\Delta x_{2n-1} - x_{2n}\Delta y_{2n})} \qquad (3.6\text{-}62)$$

**where**

**The initial direction of summation is horizontal then clockwise or counterclockwise**

**A = The area enclosed within a complex plane closed contour of shape S1**
**The closed contour sides are composed of alternating horizontal and vertical head to tail connected vectors**

**a =**
- **+1 for clockwise contour vector summation with initial vector horizontal**
- **–1 for counterclockwise contour vector summation with initial vector horizontal**
- **–1 for clockwise contour vector summation with initial vector vertical**
- **+1 for counterclockwise contour vector summation with initial vector vertical**

**N = 4,6,8,10,…**

**$z_n = x_n + jy_n$, the coordinates of the corner points of the discrete complex plane contour**

**n = 1,2,3,4,…,N**

**11. Derive the equation, $A = a\sum_{n=1}^{\frac{N}{2}} (x_{2n-1}\Delta y_{2n-1} - y_{2n}\Delta x_{2n})$**, where the initial direction of summation is vertical.

Refer to Diagram 3.2

Applying Eq 3.6-58 and Eq 3.6-59 to Diagram 3.2

$$A = \frac{1}{2}\left[ a\sum_{n=1}^{\frac{N}{2}} x_{2n-1}\Delta y_{2n-1} + (-a)\sum_{n=1}^{\frac{N}{2}} y_{2n}\Delta x_{2n} \right] \qquad (3.6\text{-}63)$$

where the initial direction of summation is vertical

The value of the constant, a, in Eq 3.6-58 and Eq 3.6-57 depends on the listing of the four summation conditions specified on the previous page. The value of the constant, a, will be either +1 or –1. In Eq 3.6-61, for the same values of n, the first summation has an initial vertical direction while the second summation has an initial horizontal direction. Then, since both summations of Eq 3.6-62 must have the same circular direction, clockwise or counterclockwise, the constant of the second summation becomes –a.

Simplifying Eq 3.6-63

$$A = \frac{a}{2}\sum_{n=1}^{\frac{N}{2}} (x_{2n-1}\Delta y_{2n-1} - y_{2n}\Delta x_{2n})$$

(3.6-64)

Then

$$\mathbf{A = \frac{a}{2}\sum_{n=1}^{\frac{N}{2}} (x_{2n-1}\Delta y_{2n-1} - y_{2n}\Delta x_{2n})} \qquad (3.6\text{-}65)$$

**where**

**The initial direction of summation is vertical then clockwise or counterclockwise**

**A = The area enclosed within a complex plane closed contour of shape S1**
**The closed contour sides are composed of alternating horizontal and vertical head to tail connected vectors**

$$\mathbf{a} = \begin{cases} +1 & \text{for clockwise contour vector summation with initial vector horizontal} \\ -1 & \text{for counterclockwise contour vector summation with initial vector horizontal} \\ -1 & \text{for clockwise contour vector summation with initial vector vertical} \\ +1 & \text{for counterclockwise contour vector summation with initial vector vertical} \end{cases}$$

**N = 4,6,8,10,…**

**$z_n = x_n + jy_n$, the coordinates of the corner points of the discrete complex plane contour**

**n = 1,2,3,4,…,N**

**12. Derive the equation, $A = -\frac{ja}{2}\sum_{n=1}^{\frac{N}{2}}[(x_{2n+1}+jy_{2n-1})^2 - (x_{2n-1}+jy_{2n-1})^2]$,** where the initial direction of summation is horizontal.

Refer to Diagram 3.6-1

From Eq 3.6-9

$$x_{2n+1} = x_{2n} \qquad (3.6\text{-}66)$$

From Eq 3.6-5

$$A = -\frac{ja}{2}\sum_{n=1}^{\frac{N}{2}}[(x_{2n}+jy_{2n-1})^2 - (x_{2n-1}+jy_{2n-1})^2] \qquad (3.6\text{-}67)$$

Substituting Eq 3.6-66 into Eq 3.6-67

$$A = -\frac{ja}{2}\sum_{n=1}^{\frac{N}{2}}[(x_{2n+1}+jy_{2n-1})^2 - (x_{2n-1}+jy_{2n-1})^2] \qquad (3.6\text{-}68)$$

Then

$$\mathbf{A = -\frac{ja}{2}\sum_{n=1}^{\frac{N}{2}}[(x_{2n+1}+jy_{2n-1})^2 - (x_{2n-1}+jy_{2n-1})^2]} \qquad (3.6\text{-}69)$$

**where**

**The initial direction of summation is horizontal then clockwise or counterclockwise**

**A = The area enclosed within a complex plane closed contour of shape S1**
**The closed contour sides are composed of alternating horizontal and vertical head to tail connected vectors**

**a =**
- **+1 for clockwise contour vector summation with initial vector horizontal**
- **−1 for counterclockwise contour vector summation with initial vector horizontal**
- **−1 for clockwise contour vector summation with initial vector vertical**
- **+1 for counterclockwise contour vector summation with initial vector vertical**

**N = 4,6,8,10,…**

**$z_n = x_n+jy_n$ , the coordinates of the corner points of the discrete complex plane contour**

**n = 1,2,3,4,…,N**

**13. Derive the equation, $A = -\frac{ja}{2}\sum_{n=1}^{\frac{N}{2}}[(x_{2n-1}+jy_{2n+1})^2 - (x_{2n-1}+jy_{2n-1})^2]$,** where the initial direction of summation is vertical.

Refer to Diagram 3.6-2

From Eq 3.6-25

$$y_{2n+1} = y_{2n} \qquad (3.6\text{-}70)$$

From Eq 3.6-21

$$A = -\frac{ja}{2}\sum_{n=1}^{\frac{N}{2}}[(x_{2n-1}+jy_{2n})^2 - (x_{2n-1}+jy_{2n-1})^2] \qquad (3.6\text{-}71)$$

Substituting Eq 3.6-70 into Eq 3.6-71

$$A = -\frac{ja}{2} \sum_{n=1}^{\frac{N}{2}} [(x_{2n-1}+jy_{2n+1})^2 - (x_{2n-1}+jy_{2n-1})^2] \qquad (3.6\text{-}72)$$

Then

$$\mathbf{A = -\frac{ja}{2} \sum_{n=1}^{\frac{N}{2}} [(x_{2n-1}+jy_{2n+1})^2 - (x_{2n-1}+jy_{2n-1})^2]} \qquad (3.6\text{-}73)$$

**where**

**The initial direction of summation is vertical then clockwise or counterclockwise**

**A = The area enclosed within a complex plane closed contour of shape S1**
**The closed contour sides are composed of alternating horizontal and vertical head to tail connected vectors**

**a =**
- **+1 for clockwise contour vector summation with initial vector horizontal**
- **–1 for counterclockwise contour vector summation with initial vector horizontal**
- **–1 for clockwise contour vector summation with initial vector vertical**
- **+1 for counterclockwise contour vector summation with initial vector vertical**

**N = 4,6,8,10,…**

**$z_n = x_n+jy_n$, the coordinates of the corner points of the discrete complex plane contour**

**n = 1,2,3,4,…,N**

The S1 complex plane closed contour area calculation equations derived above, at first glance, seem to have a limited use. These equations are applicable to complex plane closed contours with only horizontal and vertical vector sides. Such closed contours are a small part of those most commonly encountered. The closed contours most often encountered are of the S2 (vector sides of any slope) and S3 (continuous curve) shape. However, it will be shown that these same S1 area calculation equations can be used to derive equations applicable to the calculation of the areas of both the S2 and S3 closed contours. Note the following diagram, Diagram 3.6-3.

Diagram 3.6-3 Example of an S2 closed contour with its two related S1 closed contours

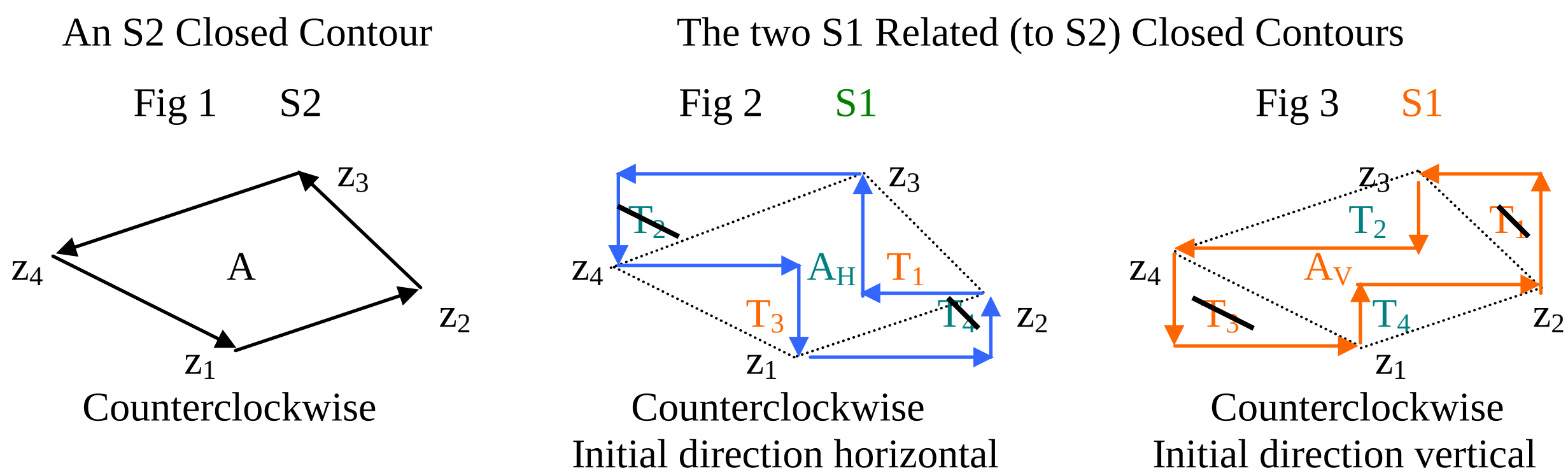

where
A = the area of a complex plane S2 closed contour
$A_H$,$A_V$ = the areas of the two complex plane S1 closed contours constructed from the S2 closed contour corner points

$z_n$ = the S2 closed contour corner points
n = 1, 2, 3, …, N
N = the number of S2 complex plane closed contour corner points
$z_1$ = the initial corner point of the S2 closed contour. The circular direction, counterclockwise or clockwise, of the head to tail closed contour vectors is usually indicated by the circular direction of the increasing $z_n$ subscripts.

Note – The circular direction of all vectors in Diagram 3.6-3 is shown as counterclockwise. This vector circular direction may also be clockwise.

One will readily notice from the above diagrams that the area of an S2 closed contour can be obtained from its related S1 closed contours with the addition or subtraction of right triangles. Since the previously derived S1 equations will calculate the area of an S1 complex plane closed contour and the area equation for a right triangle is well known, a number of equations for the area calculation of an S2 and S3 complex plane closed contour can be derived.

In Diagram 3.6-3, the two S1 complex plane closed contours are said to be related to the S2 complex plane closed contour because they have been constructed from it. Note that the vectors of the two related S1 closed contours of Fig 1 and Fig 2 have the same circular direction (here counterclockwise) as the S2 closed contour of Fig 1. This is required. The vector circular direction can be counterclockwise, as is shown, or clockwise. The two related S1 closed contours are constructed by connecting all consecutive corner points of the S2 closed contour by two vectors, one horizontal and the other vertical. Note Fig 2 and Fig 3. If the initial vector from the starting corner point, $z_1$, is horizontal, the related S1 closed contour of Fig 2 is obtained. If the initial vector from the starting corner point, $z_1$, is vertical, the related S1 closed contour of Fig 3 is obtained. In this way, the two S1 related closed contours, Fig 2 and Fig 3, of the S2 closed contour, Fig 1, are constructed. This construction procedure can be applied to any S2 closed contour to obtain its two related S1 closed contours as demonstrated in Diagram 3.6-3.

An important formula to calculate the area of an S2 complex plane closed contour from the area of its two related S1 complex plane closed contours can be derived using the three closed contours presented in Diagram 3.6-3. The derivation is fairly obvious and can be performed visually. Note that each of the two S1 complex plane closed contours, when compared to the S2 closed contour from which they were constructed, form right triangles. There will be as many right triangles as there are S2 closed contour corner points. Comparing the two related S1 closed contours to the S2 closed contour, it is noticed that an excessive right triangle in one S1 closed contour is matched by a deficiency of this same (congruent) right triangle in the other. Rearranging these right triangles amongst the two S1 closed contours results in the formation of two S2 closed contours. The area of the two related S1 complex plane closed contours is observed to be twice the area of the S2 complex plane closed contour. Thus, the formula to evaluate the area of an S2 complex plane closed contour from the area of its two related S1 complex plane closed contours is easily written as follows:

The formula to calculate the area of an S2 complex plane closed contour from the area of its two related S1 complex plane closed contours is:

$$\mathbf{A = \frac{1}{2}(A_H + A_V)} \tag{3.6-74}$$

**where**

**A = the area of an S2 complex plane closed contour with consecutive corner points designated $z_1$, $z_2$, $z_3$ …**

**$A_H$ = the area of the S1 complex plane closed contour that is related to the S2 complex plane closed contour where the initial vector direction from the $z_1$ corner point is horizontal**

**$A_V$ = the area of the S1 complex plane closed contour that is related to the S2 complex plane closed contour where the initial vector direction from the $z_1$ corner point is vertical**

The two S1 Closed Contours related to the S2 Closed Contour

S2 Closed Contour | S1 Closed Contour | S1 Closed Contour

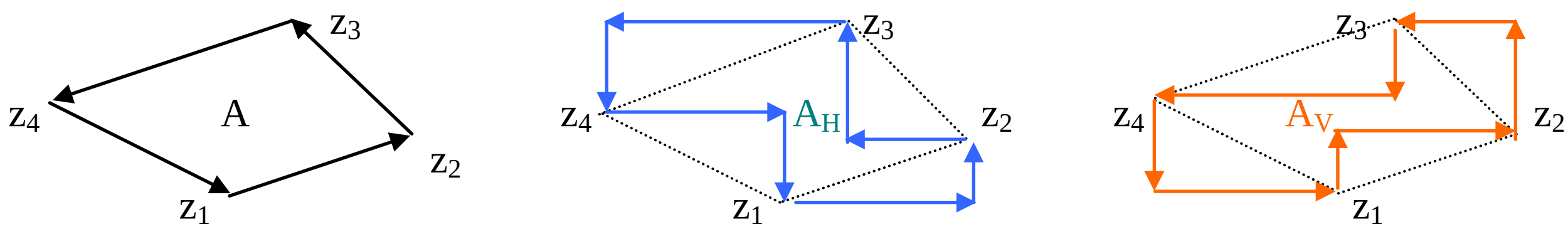

Initial vector direction horizontal | Initial vector direction vertical

Note – The vector circular directions of all three complex plane closed contours must be the same. Here, the shown circular direction is counterclockwise. The vector circular direction may also be clockwise.

Comment – Eq 3.6-74 is useful for the derivation of area calculation equations for S2 and S3 complex plane closed contours from the area calculation equations derived for S1 complex plane closed contours.

In the following derivations, several S2 and S3 complex plane closed contour area calculation equations will be derived from equations derived to calculate the area of S1 complex plane closed contours.
Eq 3.6-74 will be used.

**14. Derive the equation, $A = \frac{b}{2}\sum_{n=1}^{N}(x_n\Delta y_n - y_n\Delta x_n)$.**

From Eq 3.6-16

$$A_H = a\sum_{n=1}^{\frac{N}{2}} y_{2n-1}\,\Delta x_{2n-1}\ ,\quad \text{where the inital direction of summation is horizontal.} \qquad (3.6\text{-}75)$$

Change the notation of Eq 3.6-75 to match that used by S2 complex plane closed contour area calculation equations. Since S2 closed contours have half the number of corner points as their two S1 related closed contours, redefine the value of N to the value, $\frac{N}{2}$.

From Eq 3.6-75

$$A_H = a \sum_{n=1}^{N} y_n \Delta x_n \,, \quad \text{where the inital direction of summation is horizontal.} \qquad (3.6\text{-}76)$$

where

$$a = \begin{cases} +1 & \text{For initial direction horizontal then clockwise} \\ -1 & \text{For initial direction horizontal then counterclockwise} \\ -1 & \text{For initial direction vertical then clockwise} \\ +1 & \text{For initial direction vertical then counterclockwise} \end{cases}$$

$$b = \begin{cases} -1 & \text{For clockwise direction} \\ +1 & \text{For counterclockwise direction} \end{cases}$$

For the initial direction of summation being horizontal

$$a = \begin{cases} +1 & \text{For clockwise direction} \\ -1 & \text{For counterclockwise direction} \end{cases}$$

Then

$$a = -b \qquad (3.6\text{-}77)$$

Substituting Eq 3.6-77 into Eq 3.6-76

$$A_H = -b \sum_{n=1}^{N} y_n \Delta x_n \,, \quad \text{where the inital direction of summation is horizontal.} \qquad (3.6\text{-}78)$$

From Eq 3.6-32

$$A_v = a \sum_{n=1}^{\frac{N}{2}} x_{2n-1} \Delta y_{2n-1} \,, \text{ where the initial direction of summation is vertical.} \qquad (3.6\text{-}79)$$

Change the notation of Eq 3.6-79 to match that used by S2 complex plane closed contour area calculation equations. Since S2 closed contours have half the number of corner points as their two S1 related closed contours, redefine the value of N to the value, $\frac{N}{2}$.

From Eq 3.6-79

$$A_v = a \sum_{n=1}^{N} x_n \Delta y_n \,, \quad \text{where the inital direction of summation is horizontal.} \qquad (3.6\text{-}80)$$

For the initial direction of summation being vertical

$$a = \begin{cases} -1 & \text{For clockwise direction} \\ +1 & \text{For counterclockwise direction} \end{cases}$$

$$b = \begin{cases} -1 & \text{For clockwise direction} \\ +1 & \text{For counterclockwise direction} \end{cases}$$

Then

$$a = b \qquad (3.6\text{-}81)$$

Substituting Eq 3.6-81 into Eq 3.6-80

$$A_v = b \sum_{n=1}^{N} x_n \Delta y_n \,, \quad \text{where the initial direction of summation is vertical.} \qquad (3.6\text{-}82)$$

From Eq 3.6-74, Eq 3.6-78, and Eq 3.6-82

$$A = \frac{1}{2}(A_H + A_V) = \frac{1}{2}\left[-b \sum_{n=1}^{N} y_n \Delta x_n + b \sum_{n=1}^{N} x_n \Delta y_n\right] = \frac{b}{2} \sum_{n=1}^{N} (x_n \Delta y_n - y_n \Delta x_n) \qquad (3.6\text{-}83)$$

$$A = \frac{b}{2} \sum_{n=1}^{N} (x_n \Delta y_n - y_n \Delta x_n) \qquad (3.6\text{-}84)$$

Note that the above equation is the discrete form of Green’s Theorm for the calculation of area.

Then from Eq 3.6-84

$$\mathbf{A = \frac{b}{2} \sum_{n=1}^{N} (x_n \Delta y_n - y_n \Delta x_n)} \qquad (3.6\text{-}85)$$

**where**

**A = The area enclosed within a complex plane closed contour of shape S2**
**The closed contour sides are composed of head to tail connected vectors of any slope**

$$\mathbf{b = \begin{cases} -1 & \text{For clockwise direction} \\ +1 & \text{For counterclockwise direction} \end{cases}}$$

**$z_n = x_n + jy_n$ , the coordinates of the corner points of the discrete complex plane contour**
**n = 1,2,3,4,…,N**

**15. Derive the equation,** $\mathbf{A = \frac{b}{2}\sum_{n=1}^{N}(x_n y_{n+1} - x_{n+1} y_n)}$.

From Eq 3.6-85

$$A = \frac{b}{2}\sum_{n=1}^{N}(x_n \Delta y_n - y_n \Delta x_n) \qquad (3.6\text{-}86)$$

Expanding Eq 3.6-86

$$A = \frac{b}{2}\sum_{n=1}^{N}[x_n(y_{n+1} - y_n) - y_n(x_{n+1} - x_n)] = \frac{b}{2}\sum_{n=1}^{N}(x_n y_{n+1} - x_{n+1} y_n) \qquad (3.6\text{-}87)$$

$$A = \frac{b}{2}\sum_{n=1}^{N}(x_n y_{n+1} - x_{n+1} y_n) \qquad (3.6\text{-}88)$$

Then from Eq 3.6-88

$$\mathbf{A = \frac{b}{2}\sum_{n=1}^{N}(x_n y_{n+1} - x_{n+1} y_n)} \qquad (3.6\text{-}89)$$

**or**

$$\mathbf{A = \frac{b}{2}\sum_{n=1}^{N}\begin{vmatrix} x_n & y_n \\ x_{n+1} & y_{n+1} \end{vmatrix}} \qquad (3.6\text{-}90)$$

**where**

**A = The area enclosed within a complex plane closed contour of shape S2**
**The closed contour sides are composed of head to tail connected vectors of any slope**

$\mathbf{b = \begin{cases} -1 & \text{For clockwise direction} \\ +1 & \text{For counterclockwise direction} \end{cases}}$

$\mathbf{z_n = x_n + jy_n}$ **, the coordinates of the corner points of the discrete complex plane contour**

**n = 1,2,3,4,…,N**

**16. Derive the equation,** $\mathbf{A = \frac{jb}{4}\sum_{n=1}^{N}[(x_{n+1}+jy_n)^2 - (x_n+jy_{n+1})^2]}$.

From Eq 3.6-68

$$A_H = -\frac{ja}{2}\sum_{n=1}^{\frac{N}{2}}[(x_{2n+1}+jy_{2n-1})^2 - (x_{2n-1}+jy_{2n-1})^2],\ \text{where the inital direction of summation is horizontal.} \quad (3.6\text{-}91)$$

Change the notation of Eq 3.6-90 to match that used by S2 complex plane closed contour area calculation equations. Since S2 closed contours have half the number of corner points as their two S1 related closed contours, redefine the value of N to the value, $\frac{N}{2}$.

From Eq 3.6-91

$$A_H = -\frac{ja}{2}\sum_{n=1}^{N}[(x_{n+1}+jy_n)^2 - (x_n+jy_n)^2],\ \text{where the inital direction of summation is horizontal.} \quad (3.6\text{-}92)$$

where

$$a = \begin{cases} +1 & \text{For initial direction horizontal then clockwise} \\ -1 & \text{For initial direction horizontal then counterclockwise} \\ -1 & \text{For initial direction vertical then clockwise} \\ +1 & \text{For initial direction vertical then counterclockwise} \end{cases}$$

$$b = \begin{cases} -1 & \text{For clockwise direction} \\ +1 & \text{For counterclockwise direction} \end{cases}$$

For the initial direction of summation being horizontal

$$a = \begin{cases} +1 & \text{For clockwise direction} \\ -1 & \text{For counterclockwise direction} \end{cases}$$

Then

$$a = -b \quad (3.6\text{-}93)$$

Substituting Eq 3.6-93 into Eq 3.6-92

$$A_H = \frac{jb}{2}\sum_{n=1}^{N}[(x_{n+1}+jy_n)^2 - (x_n+jy_n)^2],\ \text{where the inital direction of summation is horizontal.} \quad (3.6\text{-}94)$$

From Eq 3.6-72

$$A_v = -\frac{ja}{2}\sum_{n=1}^{\frac{N}{2}}[(x_{2n-1}+jy_{2n+1})^2 - (x_{2n-1}+jy_{2n-1})^2], \text{ where the initial direction of summation is vertical.} \quad (3.6\text{-}95)$$

Change the notation of Eq 3.6-95 to match that used by S2 complex plane closed contour area calculation equations. Since S2 closed contours have half the number of corner points as their two S1 related closed contours, redefine the value of N to the value, $\frac{N}{2}$.

From Eq 3.6-95

$$A_v = -\frac{ja}{2}\sum_{n=1}^{\frac{N}{2}}[(x_{2n-1}+jy_{2n+1})^2 - (x_{2n-1}+jy_{2n-1})^2], \text{ where the initial direction of summation is vertical.} \quad (3.6\text{-}96)$$

For the initial direction of summation being vertical

$$a = \begin{cases} -1 & \text{For clockwise direction} \\ +1 & \text{For counterclockwise direction} \end{cases}$$

$$b = \begin{cases} -1 & \text{For clockwise direction} \\ +1 & \text{For counterclockwise direction} \end{cases}$$

Then

$$a = b \quad (3.6\text{-}97)$$

Substituting Eq 3.6-97 into Eq 3.6-96

$$A_v = -\frac{jb}{2}\sum_{n=1}^{N}[(x_n+jy_{n+1})^2 - (x_n+jy_n)^2], \text{ where the initial direction of summation is vertical.} \quad (3.6\text{-}98)$$

From Eq 3.6-74, Eq 3.6-94, and Eq 3.6-98

$$A = \frac{1}{2}(A_H + A_V) = \frac{1}{2}\{\frac{jb}{2}\sum_{n=1}^{N}[(x_{n+1}+jy_n)^2 - (x_n+jy_n)^2] - \frac{jb}{2}\sum_{n=1}^{N}[(x_n+jy_{n+1})^2 - (x_n+jy_n)^2]\} \quad (3.6\text{-}99)$$

$$A = \frac{jb}{4}\sum_{n=1}^{N}\{[(x_{n+1}+jy_n)^2 - (x_n+jy_n)^2] - [(x_n+jy_{n+1})^2 - (x_n+jy_n)^2]\} \quad (3.6\text{-}100)$$

$$A = \frac{jb}{4}\sum_{n=1}^{N}[(x_{n+1}+jy_n)^2 - (x_n+jy_{n+1})^2] \qquad (3.6\text{-}101)$$

Then from Eq 3.6-101

$$\mathbf{A = \frac{jb}{4}\sum_{n=1}^{N}[(x_{n+1}+jy_n)^2 - (x_n+jy_{n+1})^2]} \qquad (3.6\text{-}102)$$

**where**

**A = The area enclosed within a complex plane closed contour of shape S2**
**The closed contour sides are composed of head to tail connected vectors of any slope**

$$\mathbf{b = \begin{cases} -1 & \text{For clockwise direction} \\ +1 & \text{For counterclockwise direction} \end{cases}}$$

**$z_n = x_n+jy_n$, the coordinates of the corner points of the discrete complex plane contour**
**n = 1,2,3,4,…,N**

**7. Derive the equation, $A = \frac{jb}{2}\sum_{n=1}^{N}\bar{z}_n\Delta z_n^*$ .**

From Eq 3.6-48

$$A_H = -ja\sum_{n=1}^{\frac{N}{2}}\bar{z}_{2n-1}\Delta x_{2n-1} \quad , \quad \text{where the inital direction of summation is horizontal.} \qquad (3.6\text{-}103)$$

Change the notation of Eq 3.6-103 to match that used by S2 complex plane closed contour area calculation equations. Since S2 closed contours have half the number of corner points as their two S1 related closed contours, redefine the value of N to the value, $\frac{N}{2}$.

From Eq 3.6-103

$$A_H = -ja\sum_{n=1}^{N}\bar{z}_n\Delta x_n \quad , \quad \text{where the inital direction of summation is horizontal.} \qquad (3.6\text{-}104)$$

where

$$a = \begin{cases} +1 & \text{For initial direction horizontal then clockwise} \\ -1 & \text{For initial direction horizontal then counterclockwise} \\ -1 & \text{For initial direction vertical then clockwise} \\ +1 & \text{For initial direction vertical then counterclockwise} \end{cases}$$

$$b = \begin{cases} -1 & \text{For clockwise direction} \\ +1 & \text{For counterclockwise direction} \end{cases}$$

For the initial direction of summation being horizontal

$$a = \begin{cases} +1 & \text{For clockwise direction} \\ -1 & \text{For counterclockwise direction} \end{cases}$$

Then

$$a = -b \qquad (3.6\text{-}105)$$

Substituting Eq 3.6-105 into Eq 3.6-104

$$A_H = jb\sum_{n=1}^{N} \bar{z}_n \Delta x_n \ , \quad \text{where the inital direction of summation is horizontal.} \qquad (3.6\text{-}106)$$

From Eq 3.6-57

$$A_v = a\sum_{n=1}^{\frac{N}{2}} \bar{z}_{2n-1} \Delta y_{2n-1} \ , \text{ where the initial direction of summation is vertical.} \qquad (3.6\text{-}107)$$

Change the notation of Eq 3.6-107 to match that used by S2 complex plane closed contour area calculation equations. Since S2 closed contours have half the number of corner points as their two S1 related closed contours, redefine the value of N to the value, $\frac{N}{2}$.

From Eq 3.6-107

$$A_v = a\sum_{n=1}^{N} \bar{z}_n \Delta y_n \ , \quad \text{where the inital direction of summation is horizontal.} \qquad (3.6\text{-}108)$$

For the initial direction of summation being vertical

$$a = \begin{cases} -1 & \text{For clockwise direction} \\ +1 & \text{For counterclockwise direction} \end{cases}$$

$$b = \begin{cases} -1 & \text{For clockwise direction} \\ +1 & \text{For counterclockwise direction} \end{cases}$$

Then

$$a = b \qquad (3.6\text{-}109)$$

Substituting Eq 3.6-109 into Eq 3.6-108

$$A_v = b\sum_{n=1}^{N} \bar{z}_n \Delta y_n \quad , \quad \text{where the initial direction of summation is vertical.} \qquad (3.6\text{-}110)$$

From Eq 3.6-74, Eq 3.6-106, and Eq 3.6-110

$$A = \frac{1}{2}(A_H + A_V) = \frac{1}{2}\left[jb\sum_{n=1}^{N} \bar{z}_n \Delta x_n + b\sum_{n=1}^{N} \bar{z}_n \Delta y_n \right] = \frac{jb}{2}\sum_{n=1}^{N} \bar{z}_n(\Delta x_n - j\Delta y_n) \qquad (3.6\text{-}111)$$

$$A = \frac{jb}{2}\sum_{n=1}^{N} \bar{z}_n(\Delta x_n - j\Delta y_n) \qquad (3.6\text{-}112)$$

$$\Delta z_n^{*} = \Delta x_n - j\Delta y_n \qquad (3.6\text{-}113)$$

Substituting Eq 3.6-113 into Eq 3.6-112

$$A = \frac{jb}{2}\sum_{n=1}^{N} \bar{z}_n \Delta z_n^{*} \qquad (3.6\text{-}114)$$

Then

$$\mathbf{A = \frac{jb}{2}\sum_{n=1}^{N} \bar{z}_n \Delta z_n^{*}} \qquad (3.6\text{-}115)$$

**where**

**$A$ = The area enclosed within a complex plane closed contour of shape S2**
**The closed contour sides are composed of head to tail connected vectors of any slope**

**$b = \begin{cases} -1 & \text{For clockwise direction} \\ +1 & \text{For counterclockwise direction} \end{cases}$**

**$z_n = x_n + jy_n$ , the coordinates of the corner points of the discrete complex plane contour**

**$z_n^{*} = x_n - jy_n$**

**$\bar{z}_n = \dfrac{z_{n+1} - z_n}{2}$**

**$n = 1,2,3,4,\ldots,N$**

**18. Derive the equation,** $\mathbf{A = -\frac{jb}{2}\sum_{n=1}^{N}\overline{z}_n^{*}\Delta z_n}$.

From Eq 3.6-115

$$A = \frac{jb}{2}\sum_{n=1}^{N}\overline{z}_n\Delta z_n^{*} \qquad (3.6\text{-}116)$$

The equation for the complex conjugate of a product is as follows:

$$(uv)^{*} = u^{*}v^{*} \qquad (3.6\text{-}117)$$

Using Eq 3.6-117 and taking the complex conjugate of both sides of Eq 3.6-116
Since A is real, $A^{*} = A$

$$(\overline{z}_n\Delta z_n^{*})^{*} = \overline{z}_n^{*}\Delta z_n \ , \quad (\frac{jb}{2})^{*} = -\frac{jb}{2} \qquad (3.6\text{-}118)$$

Substituting Eq 3.6-118 into Eq 3.6-116

$$A = -\frac{jb}{2}\sum_{n=1}^{N}\overline{z}_n^{*}\Delta z_n \qquad (3.6\text{-}119)$$

Then

$$\mathbf{A = -\frac{jb}{2}\sum_{n=1}^{N}\overline{z}_n^{*}\Delta z_n} \qquad (3.6\text{-}120)$$

**where**

**A = The area enclosed within a complex plane closed contour of shape S2**
**The closed contour sides are composed of head to tail connected vectors of any slope**

$$\mathbf{b = \begin{cases} -1 & \text{For clockwise direction} \\ +1 & \text{For counterclockwise direction} \end{cases}}$$

$\mathbf{z_n = x_n + jy_n}$ **, the coordinates of the corner points of the discrete complex plane contour**

$$\mathbf{\overline{z}_n^{*} = (\frac{z_{n+1} - z_n}{2})^{*} = \frac{z_{n+1}^{*} - z_n^{*}}{2}}$$

$\mathbf{n = 1,2,3,4,\ldots,N}$

There are several more complex plane closed loop area calculation equations to be derived. As before, the derivations of these equations involve the area of an S1 complex plane closed contour and the addition or subtraction of the areas of right triangles in order to calculate the area of an S2 complex plane closed contour. Note the two diagrams which follow, Diagram 3.6-4 and Diagram 3.6-5. However, the method of derivation of the following equations differs from that previously employed. Eq 3.6-74 is not used.

Rewriting the S1 complex plane closed contour area equation, Eq 3.6-2. The variable, z, has been replaced by c for convenience and clarity in the derivations that follow.

$$A_B = -\frac{ja}{2}\sum_{n=1}^{\frac{N}{2}}\left( c_{2n}^{2} - c_{2n-1}^{2}\right) \qquad (3.6\text{-}121)$$

where

The initial direction of summation is horizontal or vertical then clockwise or counterclockwise

$A_B$ = The area enclosed within a complex plane closed contour of shape S1
The closed contour sides are composed of alternating horizontal and vertical head to tail connected vectors

$$a = \begin{cases} +1 & \text{for clockwise contour vector summation with initial vector horizontal} \\ -1 & \text{for counterclockwise contour vector summation with initial vector horizontal} \\ -1 & \text{for clockwise contour vector summation with initial vector vertical} \\ +1 & \text{for counterclockwise contour vector summation with initial vector vertical} \end{cases}$$

N = 4,6,8,10,…

N = the number of the discrete complex plane contour corner points, an even number

$c_n$ = the coordinates of the corner points of the discrete complex plane contour

n = 1,2,3,4,…,N

Consider Diagram 3.6-4 below. The black line complex plane closed contour is a closed contour of the shape, S1. The red line complex plane closed contour is a complex plane closed contour of the shape, S2. The following derivation is to find the area equation for closed contour of the shape, S2. Closed contours of the shape S2 have sides consisting of straight line vectors of any slope.

<u>Diagram 3.6-4  An example of Complex Plane Closed Contours of the Shape S1 and S2</u>

Complex Plane Closed Contours of the form S1 and S2

The Initial vector direction of the S1 closed contour is horizontal then counterclockwise

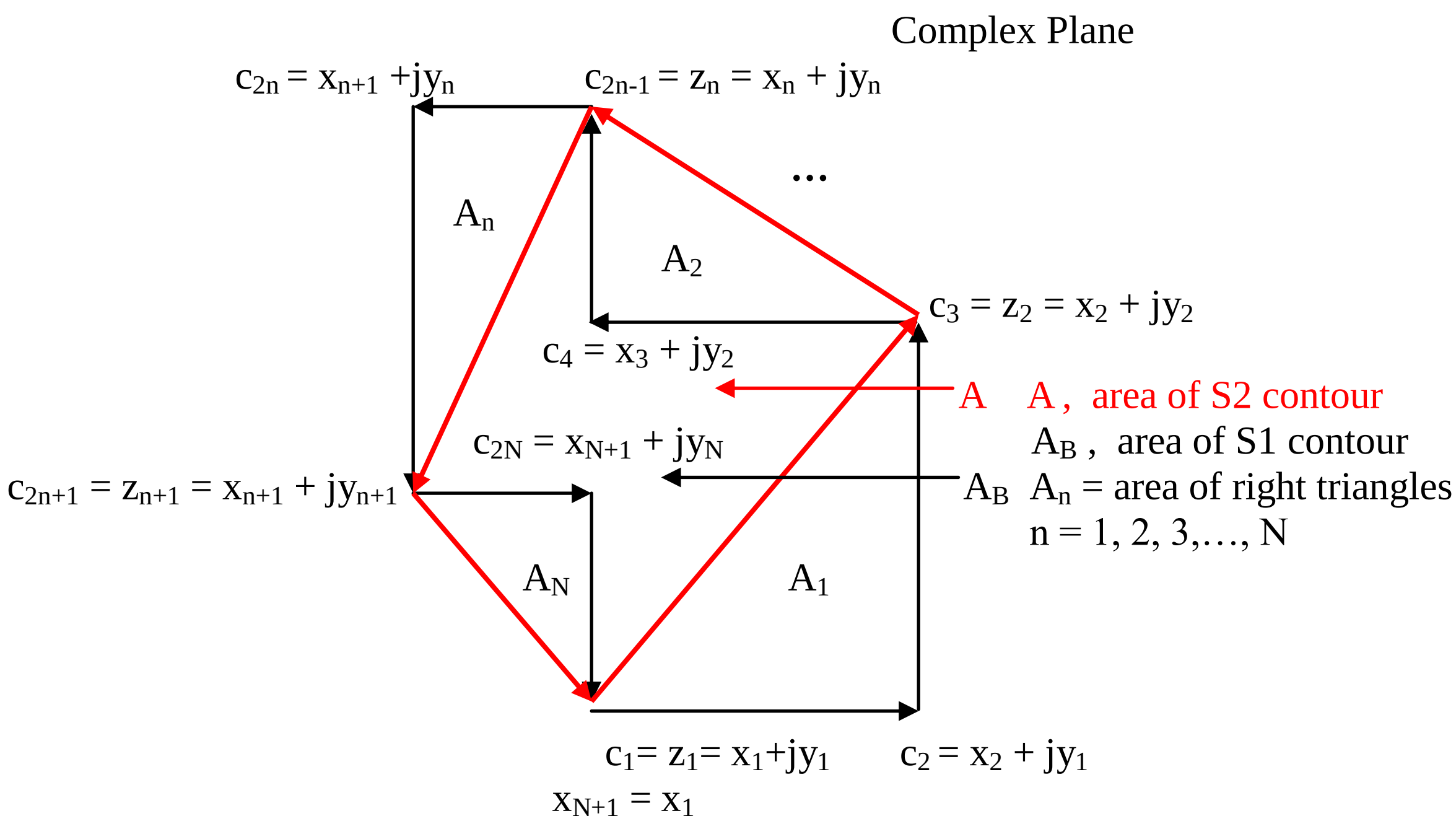

where

$A_B$ = the area within the black lined S1 closed contour

A = the area within the S2 redlined closed contour
N = the number of corner points defining the S2 redlined closed contour
2N = the number of corner points defining the black lined S1 closed contour
n = 1, 2, 3,…, N

<u>Comment</u> – In Diagram 3.6-4 the vector sides of the complex plane closed contours, S1 and S2, are shown in the counterclockwise summation direction.

Note that the area, A, can be calculated by properly adding or subtracting the area of the right triangles, $A_n$, where n = 1, 2 , 3,…, N to or from the the area, $A_B$. Since the black lined closed contour is of the shape, S1, its area can be calculated using Eq  3.6-121.

In this case, the area, $A_B$, is calculated with summation initially horizontal then counterclockwise or clockwise.

**19. Derive the equation, $A = jb \sum_{n=1}^{N} \bar{z}_n \Delta x_n$**

$$A = A_B + \sum_{n=1}^{N} A_n \qquad (3.6\text{-}122)$$

Substituting Eq 3.6-121into Eq 3.6-122

$$A = -\frac{ja}{2} \sum_{n=1}^{\frac{N}{2}} ( c_{2n}^{2} - c_{2n-1}^{2}) + \sum_{n=1}^{N} A_n \qquad (3.6\text{-}123)$$

From Diagram 3.6-4

$$c_{2n-1} = x_n + jy_n \qquad (3.6\text{-}124)$$

$$c_{2n} = x_{n+1} + jy_n \qquad (3.6\text{-}125)$$

$$A_n = a\sum_{n=1}^{N}\left[\frac{-j(x_{n+1} - x_n)(jy_{n+1} - jy_n)}{2}\right] \qquad (3.6\text{-}126)$$

For horizontal initial summation as in this case, a = -1 for counterclockwise summation and a = +1 for clockwise summation. See the definition of the constant, a, for Eq 3.6-121. The constant, a, determines, as a function of a counterclockwise or clockwise summation direction, whether the area of the right triangle, $A_n$, is added or subtracted from the area, $A_B$.

Substituting Eq 3.6-124 thru Eq 3.6-126 into Eq 3.6-123

$$A = -\frac{ja}{2} \sum_{n=1}^{N} [(x_{n+1}+jy_n)^{2} - (x_n+jy_n)^2] + a\sum_{n=1}^{N}\left[\frac{-j(x_{n+1} - x_n)(jy_{n+1} - jy_n)}{2}\right] \qquad (3.6\text{-}127)$$

Simplifying Eq 3.6-127

$$A = -\frac{ja}{2} \sum_{n=1}^{N} [(x_{n+1}+jy_n)^{2} - (x_n+jy_n)^2 + (x_{n+1} - x_n)(jy_{n+1} - jy_n)] \qquad (3.6\text{-}128)$$

$$A = -\frac{ja}{2}\sum_{n=1}^{N}[x_{n+1}^2 + 2jx_{n+1}y_n - y_n^2 - x_n^2 - 2jx_ny_n + y_n^2 + (x_{n+1} - x_n)(jy_{n+1} - jy_n)] \quad (3.6\text{-}129)$$

$$A = -\frac{ja}{2}\sum_{n=1}^{N}[(x_{n+1}^2 - x_n^2) + 2jy_n(x_{n+1} - x_n) + (x_{n+1} - x_n)(jy_{n+1} - jy_n)] \quad (3.6\text{-}130)$$

Factoring Eq 3.6-130

$$A = -\frac{ja}{2}\sum_{n=1}^{N}(x_{n+1} - x_n)[x_{n+1} + x_n + 2jy_n + jy_{n+1} - jy_n] \quad (3.6\text{-}131)$$

$$A = -\frac{ja}{2}\sum_{n=1}^{N}(x_{n+1} - x_n)[(x_{n+1} + jy_{n+1}) + (x_n + jy_n)] \quad (3.6\text{-}132)$$

$$z_n = x_n + jy_n \quad (3.6\text{-}133)$$

$$z_{n+1} = x_{n+1} + jy_{n+1} \quad (3.6\text{-}134)$$

$$\Delta x_n = x_{n+1} - x_n \quad (3.6\text{-}135)$$

$$\bar{z}_n = \frac{z_{n+1} + z_n}{2} \quad (3.6\text{-}136)$$

Substituting Eq 3.6-133 thru 3.6-136 into Eq 3.6-132

$$A = -\frac{ja}{2}\sum_{n=1}^{N}(x_{n+1} - x_n)(z_{n+1} + z_n) = -ja\sum_{n=1}^{N}(x_{n+1} - x_n)\,\bar{z}_n = -ja\sum_{n=1}^{N}\bar{z}_n\Delta x_n \quad (3.6\text{-}137)$$

From Eq 3.6-137

$$A = -ja\sum_{n=1}^{N}\bar{z}_n\Delta x_n \quad (3.6\text{-}138)$$

where

$$a = \begin{cases} +1 & \text{For initial direction horizontal then clockwise} \\ -1 & \text{For initial direction horizontal then counterclockwise} \\ -1 & \text{For initial direction vertical then clockwise} \\ +1 & \text{For initial direction vertical then counterclockwise} \end{cases}$$

$$b = \begin{cases} -1 & \text{For clockwise direction} \\ +1 & \text{For counterclockwise direction} \end{cases}$$

For the initial direction of summation being horizontal

$$a = \begin{cases} +1 & \text{For clockwise direction} \\ -1 & \text{For counterclockwise direction} \end{cases}$$

Then

$$a = -b \qquad (3.6\text{-}139)$$

Substituting Eq 3.6-139 into Eq 3.6-138

$$A_R = jb \sum_{n=1}^{N} \bar{z}_n \Delta x_n \qquad (3.6\text{-}140)$$

Then

$$\mathbf{A = jb \sum_{n=1}^{N} \bar{z}_n \, \Delta x_n} \qquad (3.6\text{-}141)$$

**where**

**A = Area enclosed within a complex plane closed contour of shape S2**

**N = the number of corner points defining the S2 complex plane closed contour**

$$\mathbf{b} = \begin{cases} \mathbf{-1} & \textbf{For clockwise direction} \\ \mathbf{+1} & \textbf{For counterclockwise direction} \end{cases}$$

$\mathbf{z_n = x_n + jy_n}$

$\mathbf{\bar{z}_n = \dfrac{z_{n+1} + z_n}{2}}$

$\mathbf{\Delta x_n = x_{n+1} - x_n}$

**n = 1, 2, 3,…, N**

**20. Derive the equation,** $\mathbf{A = -b \sum_{n=1}^{N} \bar{y}_n \, \Delta x_n}$

Rewriting Eq 3.6-131 with a = -b from Eq 3.6-139

$$A = \frac{jb}{2} \sum_{n=1}^{N} (x_{n+1} - x_n)[x_{n+1} + x_n + 2jy_n + jy_{n+1} - jy_n] \qquad (3.6\text{-}142)$$

Rearranging the terms of Eq 3.6-142

$$A = \frac{jb}{2}\sum_{n=1}^{N}[(x_{n+1} - x_n)(x_{n+1} + x_n) + j(x_{n+1} - x_n)(y_{n+1} + y_n)] \qquad (3.6\text{-}143)$$

$$A = \frac{jb}{2}\sum_{n=1}^{N}(x_{n+1}^2 - x_n^2) - b\sum_{n=1}^{N}(x_{n+1} - x_n)\left(\frac{y_{n+1} + y_n}{2}\right) \qquad (3.6\text{-}144)$$

Since $A_B$ and all $A_n$ are real, from Eq 3.6-122 A is real.

Then from Eq 3.6-144

$$\sum_{n=1}^{N}(x_{n+1}^2 - x_n^2) = 0 \;, \quad x_{N+1} = x_1 \qquad (3.6\text{-}145)$$

and

$$A = -b\sum_{n=1}^{N}(x_{n+1} - x_n)\left(\frac{y_{n+1} + y_n}{2}\right) = -b\sum_{n=1}^{N}\bar{y}_n(x_{n+1} - x_n) = -b\sum_{n=1}^{N}\bar{y}_n\Delta x_n \qquad (3.6\text{-}146)$$

Then

$$\mathbf{A = -b\sum_{n=1}^{N}\bar{y}_n\,\Delta x_n} \qquad (3.6\text{-}147)$$

**where**

**A = Area enclosed within a complex plane closed contour of shape S2**

**N = the number of corner points defining the S2 complex plane closed contour**

$$\mathbf{b} = \begin{cases} -1 & \textbf{For clockwise direction} \\ +1 & \textbf{For counterclockwise direction} \end{cases}$$

$$\bar{y}_n = \frac{y_{n+1} + y_n}{2}$$

$$\Delta x_n = x_{n+1} - x_n$$

**n = 1, 2, 3,…, N**

Diagram 3.6-5  An example of Complex Plane Closed Contours of the Shape S1 and S2

Complex Plane Closed Contours of the form S1 and S2
The Initial vector direction of the S1 closed contour is vertical then counterclockwise

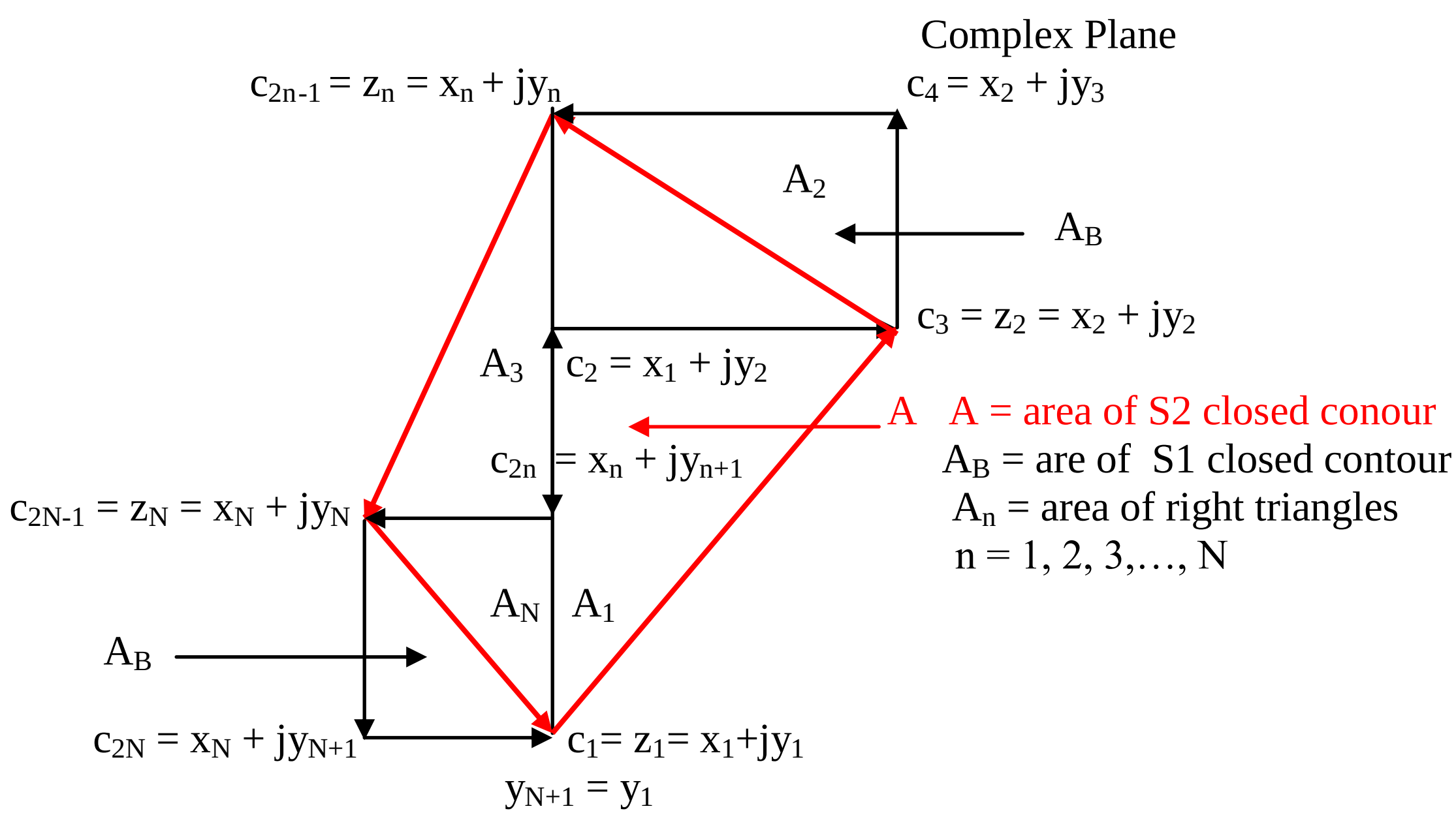

where
$A_B$ = the area within the black lined S1 closed contour
A = the area within the S2 redlined closed contour
N = the number of corner points defining the S2 redlined closed contour
2N = the number of corner points defining the black lined S1 closed contour
n = 1, 2, 3,…, N

Comment – In Diagram 3.6-5 the vector sides of the complex plane closed contours, S1 and S2, are shown in the counterclockwise summation direction.

Note that the area, A, can be calculated by properly adding or subtracting the area of the right triangles, $A_n$, where n = 1, 2 , 3,…, N to or from the area, $A_B$. Since the black lined closed contour is of the shape, S1, its area can be calculated using Eq  3.6-121.

In this case, the area, $A_B$, is calculated with summation initially vertical then counterclockwise or clockwise.

**21.  Derive the equation,** $\mathbf{A = b\sum_{n=1}^{N} \bar{z}_n \Delta y_n}$

$$A = A_B + \sum_{n=1}^{N} A_n \qquad (3.6\text{-}148)$$

Substituting Eq 3.6-121into Eq 3.6-148

$$A = -\frac{ja}{2}\sum_{n=1}^{\frac{N}{2}}( c_{2n}^2 - c_{2n-1}^2) + \sum_{n=1}^{N} A_n \qquad (3.6\text{-}149)$$

From Diagram 3.6-5

$$c_{2n-1} = x_n + jy_n \qquad (3.6\text{-}150)$$

$$c_{2n} = x_n + jy_{n+1} \qquad (3.6\text{-}151)$$

$$A_n = a\sum_{n=1}^{N}[\frac{-j(x_{n+1} - x_n)(jy_{n+1} - jy_n)}{2}] \qquad (3.6\text{-}152)$$

For vertical initial summation as in this case, a = +1 for counterclockwise summation and a = –1 for clockwise summation. See the definition of the constant, a, for Eq 3.6-121. The constant, a, determines, as a function of a counterclockwise or clockwise summation direction, whether the area of the right triangle, $A_n$, is added or subtracted from the area, $A_B$.

Substituting Eq 3.6-150 thru Eq 3.6-152 into Eq 3.6-149

$$A = -\frac{ja}{2}\sum_{n=1}^{N}[(x_n+jy_{n+1})^2 - (x_n+jy_n)^2] + a\sum_{n=1}^{N}[\frac{-j(x_{n+1} - x_n)(jy_{n+1} - jy_n)}{2}] \qquad (3.6\text{-}153)$$

Simplifying Eq 3.6-153

$$A = -\frac{ja}{2}\sum_{n=1}^{N}[(x_n+jy_{n+1})^2 - (x_n+jy_n)^2 + (x_{n+1} - x_n)(jy_{n+1} - jy_n)] \qquad (3.6\text{-}154)$$

$$A = -\frac{ja}{2}\sum_{n=1}^{N}[x_n^2 + 2jx_ny_{n+1} - y_{n+1}^2 - x_n^2 - 2jx_ny_n + y_n^2 + (x_{n+1} - x_n)(jy_{n+1} - jy_n)] \qquad (3.6\text{-}155)$$

$$A_R = -\frac{ja}{2}\sum_{n=1}^{N}[-(y_{n+1}^2 - y_n^2) + 2jx_n(y_{n+1} - y_n) + (x_{n+1} - x_n)(jy_{n+1} - jy_n)] \qquad (3.6\text{-}156)$$

Factoring Eq 3.6-156

$$A = -\frac{ja}{2}\sum_{n=1}^{N}(y_{n+1} - y_n)[-y_{n+1} - y_n + 2jx_n + jx_{n+1} - jx_n] \qquad (3.6\text{-}157)$$

$$A = -\frac{ja}{2}\sum_{n=1}^{N}(y_{n+1} - y_n)[-y_{n+1} + jx_{n+1} - y_n + jx_n] \qquad (3.6\text{-}158)$$

$$A = \frac{a}{2}\sum_{n=1}^{N}(y_{n+1} - y_n)[x_{n+1} + jy_{n+1} + x_n + jy_n] \qquad (3.6\text{-}159)$$

$$z_n = x_n + jy_n \qquad (3.6\text{-}160)$$

$$z_{n+1} = x_{n+1} + jy_{n+1} \qquad (3.6\text{-}161)$$

$$\Delta y = y_{n+1} - y_n \qquad (3.6\text{-}162)$$

$$\bar{z}_n = \frac{z_{n+1} + z_n}{2} \qquad (3.6\text{-}163)$$

Substituting Eq 3.6-160 thru 3.6-163 into Eq 3.6-159

$$A = \frac{a}{2}\sum_{n=1}^{N}(y_{n+1} - y_n)(z_{n+1} + z_n) = a\sum_{n=1}^{N}(y_{n+1} - y_n)\,\bar{z}_n = a\sum_{n=1}^{N}\bar{z}_n\Delta y_n \qquad (3.6\text{-}164)$$

$$A = a\sum_{n=1}^{N}\bar{z}_n\Delta y_n \qquad (3.6\text{-}165)$$

where

$$a = \begin{cases} +1 & \text{For initial direction horizontal then clockwise} \\ -1 & \text{For initial direction horizontal then counterclockwise} \\ -1 & \text{For initial direction vertical then clockwise} \\ +1 & \text{For initial direction vertical then counterclockwise} \end{cases}$$

$$b = \begin{cases} -1 & \text{For clockwise direction} \\ +1 & \text{For counterclockwise direction} \end{cases}$$

For the initial direction of summation being vertical

$$a = \begin{cases} -1 & \text{For clockwise direction} \\ +1 & \text{For counterclockwise direction} \end{cases}$$

Then

$$a = b \qquad (3.6\text{-}166)$$

Substituting Eq 3.6-166 into Eq 3.6-165

$$A = b\sum_{n=1}^{N} \bar{z}_n \Delta y_n \qquad (3.6\text{-}167)$$

Then

$$\mathbf{A = b\sum_{n=1}^{N} \bar{z}_n \, \Delta y_n} \qquad (3.6\text{-}168)$$

**where**

**A = Area enclosed within a complex plane closed contour of shape S2**

**N = the number of corner points defining the S2 complex plane closed contour**

$$\mathbf{b = \begin{cases} -1 & \text{For clockwise direction} \\ +1 & \text{For counterclockwise direction} \end{cases}}$$

$\mathbf{z_n = x_n + jy_n}$

$\mathbf{\bar{z}_n = \dfrac{z_{n+1} + z_n}{2}}$

$\mathbf{\Delta y_n = y_{n+1} - y_n}$

**n = 1, 2, 3,…, N**

**22. Derive the equation,** $\mathbf{A = b\sum_{n=1}^{N} \bar{x}_n \Delta y_n}$

Rewriting Eq 3.6-159 with a = b from Eq 3.6-166

$$A = \frac{b}{2}\sum_{n=1}^{N}(y_{n+1} - y_n)[x_{n+1} + jy_{n+1} + x_n + jy_n] \qquad (3.6\text{-}169)$$

Rearranging the terms of Eq 3.6-169

$$A = \frac{b}{2}\sum_{n=1}^{N}(y_{n+1} - y_n)[j(y_{n+1} + y_n) + (x_{n+1} + x_n)] \qquad (3.6\text{-}170)$$

$$A = \frac{jb}{2}\sum_{n=1}^{N}(y_{n+1}^2 - y_n^2) + b\sum_{n=1}^{N}(y_{n+1} - y_n)\left(\frac{x_{n+1} + x_n}{2}\right) \qquad (3.6\text{-}171)$$

Since $A_B$ and all $A_n$ are real, from Eq 3.6-148 A is real.

Then from Eq 3.6-171

$$\sum_{n=1}^{N} (y_{n+1}^2 - y_n^2) = 0 \quad , \quad y_{N+1} = y_1 \tag{3.6-172}$$

and

$$A = b\sum_{n=1}^{N} (y_{n+1} - y_n)\left(\frac{x_{n+1} + x_n}{2}\right) = b\sum_{n=1}^{N} \bar{x}_n(y_{n+1} - y_n) = b\sum_{n=1}^{N} \bar{x}_n\Delta y_n \tag{3.6-173}$$

Then

$$\mathbf{A = b\sum_{n=1}^{N} \bar{x}_n \, \Delta y_n} \tag{3.6-174}$$

**where**

**A = Area enclosed within a complex plane closed contour of shape S2**
**N = the number of corner points defining the S2 complex plane closed contour**

$$\mathbf{b} = \begin{cases} \mathbf{-1} & \textbf{For clockwise direction} \\ \mathbf{+1} & \textbf{For counterclockwise direction} \end{cases}$$

$$\mathbf{\bar{x}_n = \frac{x_{n+1} + x_n}{2}}$$

$\mathbf{\Delta y_n = y_{n+1} - y_n}$
**n = 1, 2, 3,…, N**

**23. Derive the equation,** $\mathbf{A = \frac{b}{2}\sum_{n=1}^{N} (\bar{x}_n \Delta y_n - \bar{y}_n \Delta x_n)}$

From Eq 3.6-174

$$A = b\sum_{n=1}^{N} \bar{x}_n \, \Delta y_n \tag{3.6-175}$$

From Eq 3.6-147

$$A = -b\sum_{n=1}^{N} \bar{y}_n \, \Delta x_n \tag{3.6-176}$$

Adding Eq 3.6-175 and Eq 3.6-176

$$A = \frac{b}{2}\sum_{n=1}^{N}(\bar{x}_n\,\Delta y_n - \bar{y}_n\,\Delta x_n) \qquad (3.6\text{-}177)$$

Then

$$\mathbf{A = \frac{b}{2}\sum_{n=1}^{N}(\bar{x}_n\,\Delta y_n - \bar{y}_n\,\Delta x_n)} \qquad (3.6\text{-}178)$$

**where**

**A = Area enclosed within a complex plane closed contour of shape S2**

**N = the number of corner points defining the S2 complex plane closed contour**

$$\mathbf{b} = \begin{cases} \mathbf{-1} & \textbf{For clockwise direction} \\ \mathbf{+1} & \textbf{For counterclockwise direction} \end{cases}$$

$$\bar{x}_n = \frac{x_{n+1} + x_n}{2}$$

$$\bar{y}_n = \frac{y_{n+1} + y_n}{2}$$

$$\Delta x_n = x_{n+1} - x_n$$

$$\Delta y_n = y_{n+1} - y_n$$

**n = 1, 2, 3,…, N**

**24. Derive the equation,** $\mathbf{A = \frac{b}{2}\sum_{n=1}^{N}(x_{n+1}\Delta y_n - y_{n+1}\Delta x_n)}$

From Eq 3.6-178

$$A = \frac{b}{2}\sum_{n=1}^{N}(\bar{x}_n\,\Delta y_n - \bar{y}_n\,\Delta x_n) \qquad (3.6\text{-}179)$$

Expanding Eq 3.6-179

$$A = \frac{b}{4}\sum_{n=1}^{N}[(x_{n+1} + x_n)\,\Delta y_n - (y_{n+1} + y_n)\,\Delta x_n] \qquad (3.6\text{-}180)$$

Rearranging terms

$$2A = \frac{b}{2}\sum_{n=1}^{N}(x_n\Delta y_n - y_n\Delta x_n) + \frac{b}{2}\sum_{n=1}^{N}(x_{n+1}\Delta y_n - y_{n+1}\Delta x_n) \qquad (3.6\text{-}181)$$

But from Eq 3.6-85

$$A = \frac{b}{2}\sum_{n=1}^{N}(x_n\Delta y_n - y_n\Delta x_n) \qquad (3.6\text{-}182)$$

Substituting Eq 3.6-182 into Eq 3.6-181

$$2A = A + \frac{b}{2}\sum_{n=1}^{N}(x_{n+1}\Delta y_n - y_{n+1}\Delta x_n) \qquad (3.6\text{-}183)$$

$$A = \frac{b}{2}\sum_{n=1}^{N}(x_{n+1}\Delta y_n - y_{n+1}\Delta x_n) \qquad (3.6\text{-}184)$$

Then

$$\mathbf{A = \frac{b}{2}\sum_{n=1}^{N}(x_{n+1}\Delta y_n - y_{n+1}\Delta x_n)} \qquad (3.6\text{-}185)$$

**where**

**A = Area enclosed within a complex plane closed contour of shape S2**

**N = the number of corner points defining the S2 complex plane closed contour**

$$\mathbf{b = \begin{cases} -1 & \text{For clockwise direction} \\ +1 & \text{For counterclockwise direction} \end{cases}}$$

$\mathbf{\Delta x_n = x_{n+1} - x_n}$

$\mathbf{\Delta y_n = y_{n+1} - y_n}$

**n = 1, 2, 3,…, N**

<u>Derivation of equations which calculate the area of S3 continuous complex plane closed contours</u>

The derivations which follow are of equations that calculate the area of continuous complex plane closed contours. These continuous complex plane closed contours have been categorized as S3 complex plane closed contours. As previously mentioned, S3 complex plane closed contours are a special case of S2 complex plane closed contours. S2 closed contours consist of head to tail connected complex plane vectors which form a closed loop. The enclosed area is said to be the area of the closed loop. The S2 vectors are of any slope and of finite number. S3 complex plane closed contours are a special case of S2 complex plane closed contours where the number of vectors forming the closed contour and the number of corner points is infinite. Thus, S3 closed contours have the appearance of a continuous closed curve. As will be seen below, The area calculation equations of S3 complex plane closed contours can be easily derived from their related S2 complex plane closed contours. In particular, increments of length approach zero, sets of values and discrete functions become continuous, and summations become integrals.

**25. Derive the equation, $\mathbf{A = \frac{b}{2}\oint_c [xdy - ydx]}$**, Green’s area calculation equation

Rewriting Eq 3.6-85

$$A = \frac{b}{2}\sum_{n=1}^{N}(x_n\Delta y_n - y_n\Delta x_n) \qquad (3.6\text{-}186)$$

where

$x_n = x(n)$ , a discrete function of n

$y_n = y(n)$ , a discrete function of n

$n = 1, 2, 3,\ldots, N$

N = the number of corner points in the complex plane closed loop

$\Delta x_n = x_{n+1} - x_n$

$\Delta y_n = y_{n+1} - y_n$

$$b = \begin{cases} -1 & \text{For clockwise direction} \\ +1 & \text{For counterclockwise direction} \end{cases}$$

For $N \to \infty$

The summation becomes a closed contour integral

$x_n = x(n)$, a discrete function of n becomes $x = x(\theta)$, a continuous function of $\theta$

$y_n = y(n)$, a discrete function of n becomes $y = y(\theta)$, a continuous function of $\theta$

$0 \le \theta < 2\pi$

$x_n \to x$

$y_n \to y$

$\Delta x_n \to dx$

$\Delta y_n \to dy$

Substituting into Eq 3.6-186

$$A = \frac{b}{2}\oint_c [xdy - ydx] \qquad (3.6\text{-}187)$$

Then

$$\mathbf{A = \frac{b}{2}\oint_c [xdy - ydx]} \qquad (3.6\text{-}188)$$

**where**

**A = Area enclosed within a continuous curve (S3) complex plane closed contour**

**c = closed S3 contour in the complex plane**

$$\mathbf{b = \begin{cases} -1 & \text{For clockwise integration direction} \\ +1 & \text{For counterclockwise integration direction} \end{cases}}$$

See Diagram 3.6-6 below

<u>Diagram 3.6-6  Formation of a complex plane S3 continuous curve closed contour from an S2 complex plane closed contour where the number of its corner points becomes infinite, $N \to \infty$</u>

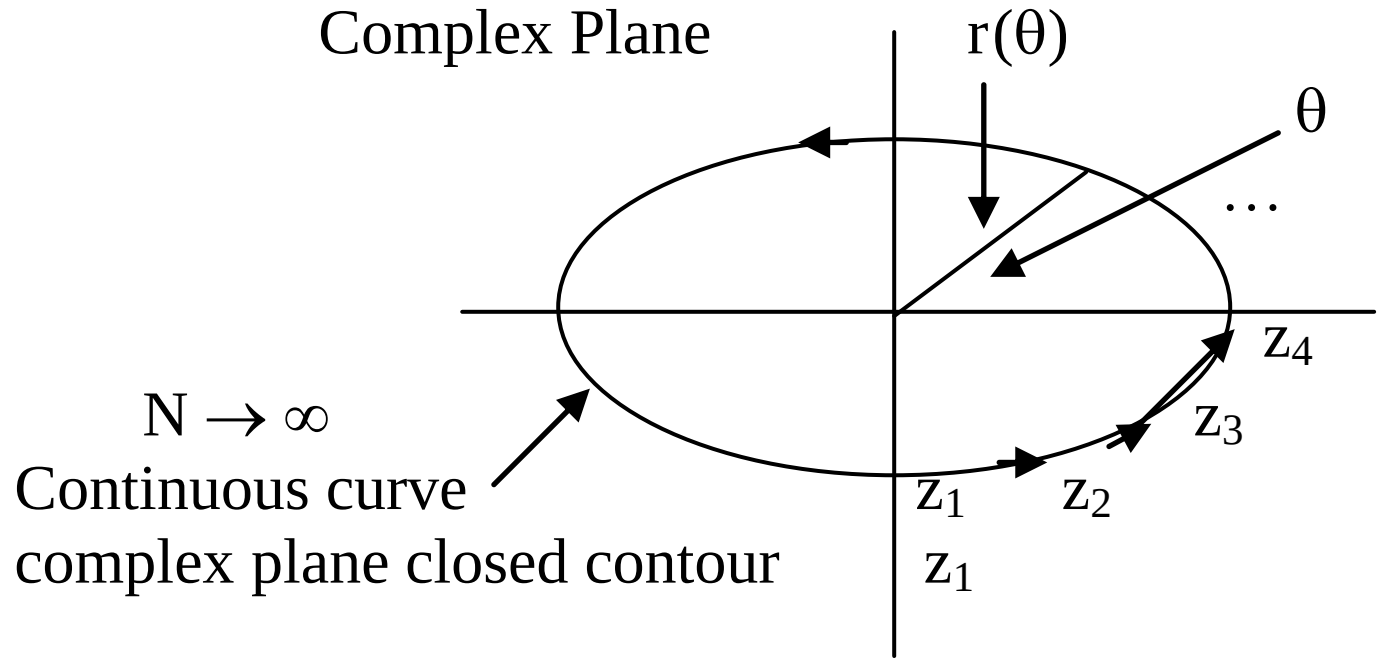

$z_n = z(n)$, discrete function of n
$z = z(\theta)$, continuous function of θ
$z = x(\theta) + jy(\theta) = x + jy$
$z_n = x_n + jy_n$
$z_n$ = corner points of a complex plane closed contour
$n = 1, 2, 3, \ldots, N$
N = number of corner points defining the complex plane closed contour
$0 \le \theta < 2\pi$

For $N \to \infty$, the number of corner points defining the complex plane closed contour becomes infinite thereby forming a continuous curve complex plane closed contour

**26. Derive the equation, $A = -b \oint_c y\,dx$**

Rewriting Eq 3.6-147

$$A = -b\sum_{n=1}^{N} \bar{y}_n \Delta x_n \tag{3.6-189}$$

where

$x_n = x(n)$ , a discrete function of n
$y_n = y(n)$ , a discrete function of n
$n = 1, 2, 3, \ldots, N$
N = the number of corner points in the complex plane closed loop

$$\bar{y}_n = \frac{y_{n+1} - y_n}{2}$$

$\Delta x_n = x_{n+1} - x_n$

$$b = \begin{cases} -1 & \text{For clockwise direction} \\ +1 & \text{For counterclockwise direction} \end{cases}$$

For $N \to \infty$

The summation becomes a closed contour integral
$x_n = x(n)$, a discrete function of n becomes $x = x(\theta)$, a continuous function of θ
$y_n = y(n)$, a discrete function of n becomes $y = y(\theta)$, a continuous function of θ
$0 \le \theta < 2\pi$
$x_n \to x$

$y_n \to y$
$\Delta x_n \to dx$
$\bar{y}_n \to y$

Substituting into Eq 3.6-189

$$A = -b \oint_c y\,dx \tag{3.6-190}$$

Then

$$\mathbf{A = -b \oint_c y\,dx} \tag{3.6-191}$$

**where**

**A = Area enclosed within a continuous curve (S3) complex plane closed contour**
**c = closed S3 contour in the complex plane**

$$\mathbf{b = \begin{cases} -1 & \text{For clockwise integration direction} \\ +1 & \text{For counterclockwise integration direction} \end{cases}}$$

**27. Derive the equation, $A = b \oint_c x\,dy$**

Rewriting Eq 3.6-174

$$A = b \sum_{n=1}^{N} \bar{x}_n \Delta y_n \tag{3.6-192}$$

where

$x_n = x(n)$ , a discrete function of n
$y_n = y(n)$ , a discrete function of n
$n = 1, 2, 3, \ldots, N$
N = the number of corner points in the complex plane closed loop

$$\bar{x}_n = \frac{x_{n+1} - x_n}{2}$$

$\Delta y_n = y_{n+1} - y_n$

$$b = \begin{cases} -1 & \text{For clockwise direction} \\ +1 & \text{For counterclockwise direction} \end{cases}$$

For $N \to \infty$

The summation becomes a closed contour integral
$x_n = x(n)$, a discrete function of n becomes $x = x(\theta)$, a continuous function of $\theta$
$y_n = y(n)$, a discrete function of n becomes $y = y(\theta)$, a continuous function of $\theta$
$0 \le \theta < 2\pi$
$x_n \to x$

$y_n \rightarrow y$
$\Delta y_n \rightarrow dy$
$\bar{x}_n \rightarrow x$

Substituting into Eq 3.6-192

$$A = b \oint_c x\,dy \qquad (3.6\text{-}193)$$

Then

$$\mathbf{A = b \oint_c x\,dy} \qquad (3.6\text{-}194)$$

**where**

**A = Area enclosed within a continuous curve (S3) complex plane closed contour**
**c = closed S3 contour in the complex plane**

$$\mathbf{b} = \begin{cases} -1 & \text{For clockwise integration direction} \\ +1 & \text{For counterclockwise integration direction} \end{cases}$$

**28. Derive the equation, $A = jb \oint_c z\,dx$**

Rewriting Eq 3.6-141

$$A = jb \sum_{n=1}^{N} \bar{z}_n \Delta x_n \qquad (3.6\text{-}195)$$

where

$x_n = x(n)$ , a discrete function of n
$y_n = y(n)$ , a discrete function of n
$n = 1, 2, 3,\ldots, N$
N = the number of corner points in the complex plane closed loop
$z_n = x_n + jy_n$
$\bar{z}_n = \dfrac{z_{n+1} - z_n}{2}$
$\Delta x_n = x_{n+1} - x_n$

$$b = \begin{cases} -1 & \text{For clockwise direction} \\ +1 & \text{For counterclockwise direction} \end{cases}$$

For $N \rightarrow \infty$

The summation becomes a closed contour integral
$x_n = x(n)$, a discrete function of n becomes $x = x(\theta)$, a continuous function of $\theta$
$y_n = y(n)$, a discrete function of n becomes $y = y(\theta)$, a continuous function of $\theta$
$0 \le \theta < 2\pi$

$x_n \rightarrow x$
$y_n \rightarrow y$
$z_n \rightarrow z$
$\Delta x_n \rightarrow dx$
$\bar{z}_n \rightarrow z$

Substituting into Eq 3.6-195

$$A = jb \oint_c z\,dx \qquad (3.6\text{-}196)$$

Then

$$\mathbf{A = jb \oint_c z\,dx} \qquad (3.6\text{-}197)$$

**where**
**A = Area enclosed within a continuous curve (S3) complex plane closed contour**
**c = closed S3 contour in the complex plane**
$$\mathbf{b = \begin{cases} -1 & \text{For clockwise integration direction} \\ +1 & \text{For counterclockwise integration direction} \end{cases}}$$

**29. Derive the equation,** $\mathbf{A = b \oint_c z\,dy}$

Rewriting Eq 3.6-168

$$A = b \sum_{n=1}^{N} \bar{z}_n \Delta y_n \qquad (3.6\text{-}198)$$

where
$x_n = x(n)$ , a discrete function of n
$y_n = y(n)$ , a discrete function of n
$n = 1, 2, 3,\ldots, N$
N = the number of corner points in the complex plane closed loop
$z_n = x_n + jy_n$
$\bar{z}_n = \dfrac{z_{n+1} - z_n}{2}$
$\Delta y_n = y_{n+1} - y_n$
$$b = \begin{cases} -1 & \text{For clockwise direction} \\ +1 & \text{For counterclockwise direction} \end{cases}$$

For $N \rightarrow \infty$
The summation becomes a closed contour integral
$x_n = x(n)$, a discrete function of n becomes $x = x(\theta)$, a continuous function of θ
$y_n = y(n)$, a discrete function of n becomes $y = y(\theta)$, a continuous function of θ

$0 \le \theta < 2\pi$
$x_n \to x$
$y_n \to y$
$z_n \to z$
$\Delta y_n \to dy$
$\bar{z}_n \to z$

Substituting into Eq 3.6-198

$$A = b \oint_c z\,dy \tag{3.6-199}$$

Then

$$\mathbf{A = b \oint_c z\,dy} \tag{3.6-200}$$

**where**

**A = Area enclosed within a continuous curve (S3) complex plane closed contour**
**c = closed S3 contour in the complex plane**

$$\mathbf{b} = \begin{cases} -1 & \textbf{For clockwise integration direction} \\ +1 & \textbf{For counterclockwise integration direction} \end{cases}$$

**30. Derive the equation,** $\mathbf{A = \frac{jb}{2} \oint_c z\,dz^*}$

Rewriting Eq 3.6-115

$$A = \frac{jb}{2} \sum_{n=1}^{N} \bar{z}_n \Delta z_n^* \tag{3.6-201}$$

where

$x_n = x(n)$ , a discrete function of n
$y_n = y(n)$ , a discrete function of n
$n = 1, 2, 3,\ldots, N$
N = the number of corner points in the complex plane closed loop
$z_n = x_n + jy_n$
$\bar{z}_n = \dfrac{z_{n+1} - z_n}{2}$
$z_n^*$ = complex conjugate of $z_n$
$\Delta z_n^* = z_{n+1}^* - z_n^*$

$$b = \begin{cases} -1 & \text{For clockwise direction} \\ +1 & \text{For counterclockwise direction} \end{cases}$$

For $N \to \infty$

The summation becomes a closed contour integral
$x_n = x(n)$, a discrete function of n becomes $x = x(\theta)$, a continuous function of $\theta$
$y_n = y(n)$, a discrete function of n becomes $y = y(\theta)$, a continuous function of $\theta$

$0 \le \theta < 2\pi$
$x_n \to x$
$y_n \to y$
$z_n \to z$
$\Delta z_n^* \to dz^*$
$\overline{z}_n \to z$

Substituting into Eq 3.6-201

$$A = \frac{jb}{2} \oint_c z\,dz^* \qquad (3.6\text{-}202)$$

Then

$$\mathbf{A = \frac{jb}{2} \oint_c z\,dz^*} \qquad (3.6\text{-}203)$$

**where**

**A = Area enclosed within a continuous curve (S3) complex plane closed contour**
**c = closed S3 contour in the complex plane**

$$\mathbf{b = \begin{cases} -1 & \text{For clockwise integration direction} \\ +1 & \text{For counterclockwise integration direction} \end{cases}}$$

**31. Derive the equation, $A = -\frac{jb}{2} \oint_c z^* dz$**

Rewriting Eq 3.6-120

$$A = \frac{jb}{2} \sum_{n=1}^{N} \overline{z}_n^* \, \Delta z_n \qquad (3.6\text{-}204)$$

where

$x_n = x(n)$ , a discrete function of n
$y_n = y(n)$ , a discrete function of n
$n = 1, 2, 3,\ldots, N$
N = the number of corner points in the complex plane closed loop
$z_n = x_n + jy_n$
$\overline{z}_n^* = \dfrac{z_{n+1}^* - z_n^*}{2}$
$z_n^*$ = complex conjugate of $z_n$
$\Delta z_n = z_{n+1} - z_n$

$$b = \begin{cases} -1 & \text{For clockwise direction} \\ +1 & \text{For counterclockwise direction} \end{cases}$$

For $N \to \infty$

The summation becomes a closed contour integral
$x_n = x(n)$, a discrete function of n becomes $x = x(\theta)$, a continuous function of $\theta$

$y_n = y(n)$, a discrete function of n becomes $y = y(\theta)$, a continuous function of $\theta$
$0 \le \theta < 2\pi$
$x_n \rightarrow x$
$y_n \rightarrow y$
$z_n \rightarrow z$
$\Delta z_n \rightarrow dz$
$\overline{z}_n^{\,*} \rightarrow z^*$

Substituting into Eq 3.6-204

$$A = -\frac{jb}{2} \oint_c z^* dz \qquad (3.6\text{-}205)$$

Then

$$\mathbf{A = -\frac{jb}{2} \oint_c z^* dz} \qquad (3.6\text{-}206)$$

**where**

**A = Area enclosed within a continuous curve (S3) complex plane closed contour**

**c = closed S3 contour in the complex plane**

$$\mathbf{b = \begin{cases} -1 & \text{For clockwise integration direction} \\ +1 & \text{For counterclockwise integration direction} \end{cases}}$$

The following two diagrams, Diagram 3.6-7 and Diagram 3.6-8, are used to demonstrate the application of some of the complex plane closed contour area calculation equations previously derived.

Diagram 3.6-7  An Example of Complex Plane S1 (black) and S2 (red) Closed Contours

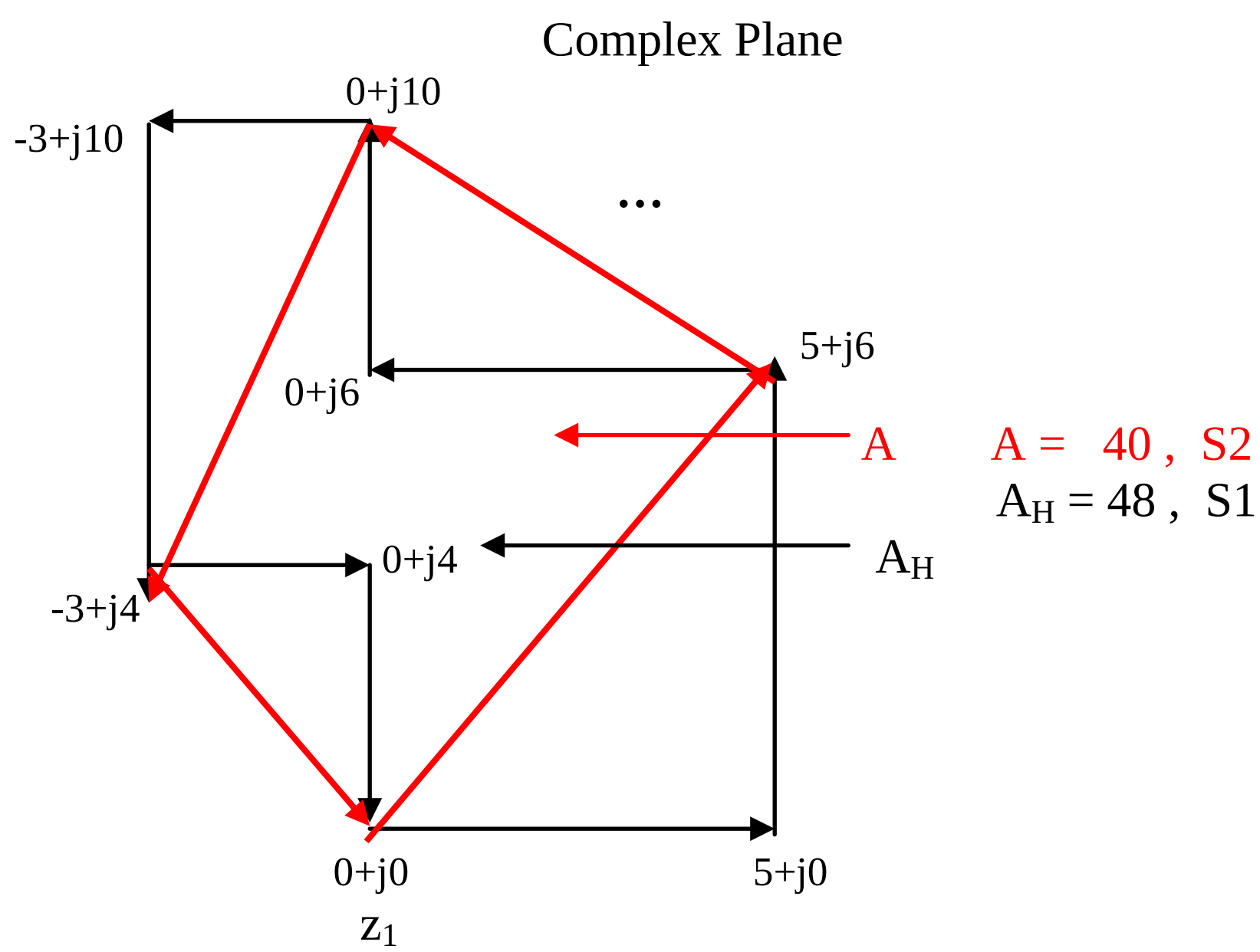

Note that the initial vector direction from the starting point point, $z_1$, is shown to be horizontal then counterclockwise. If the vector directions are reversed, the initial vector direction from the starting point point, $z_1$, is vertical then clockwise.

Diagram 3.6-8  Another Example of Complex Plane S1 (black) and S2 (red) Closed Contours

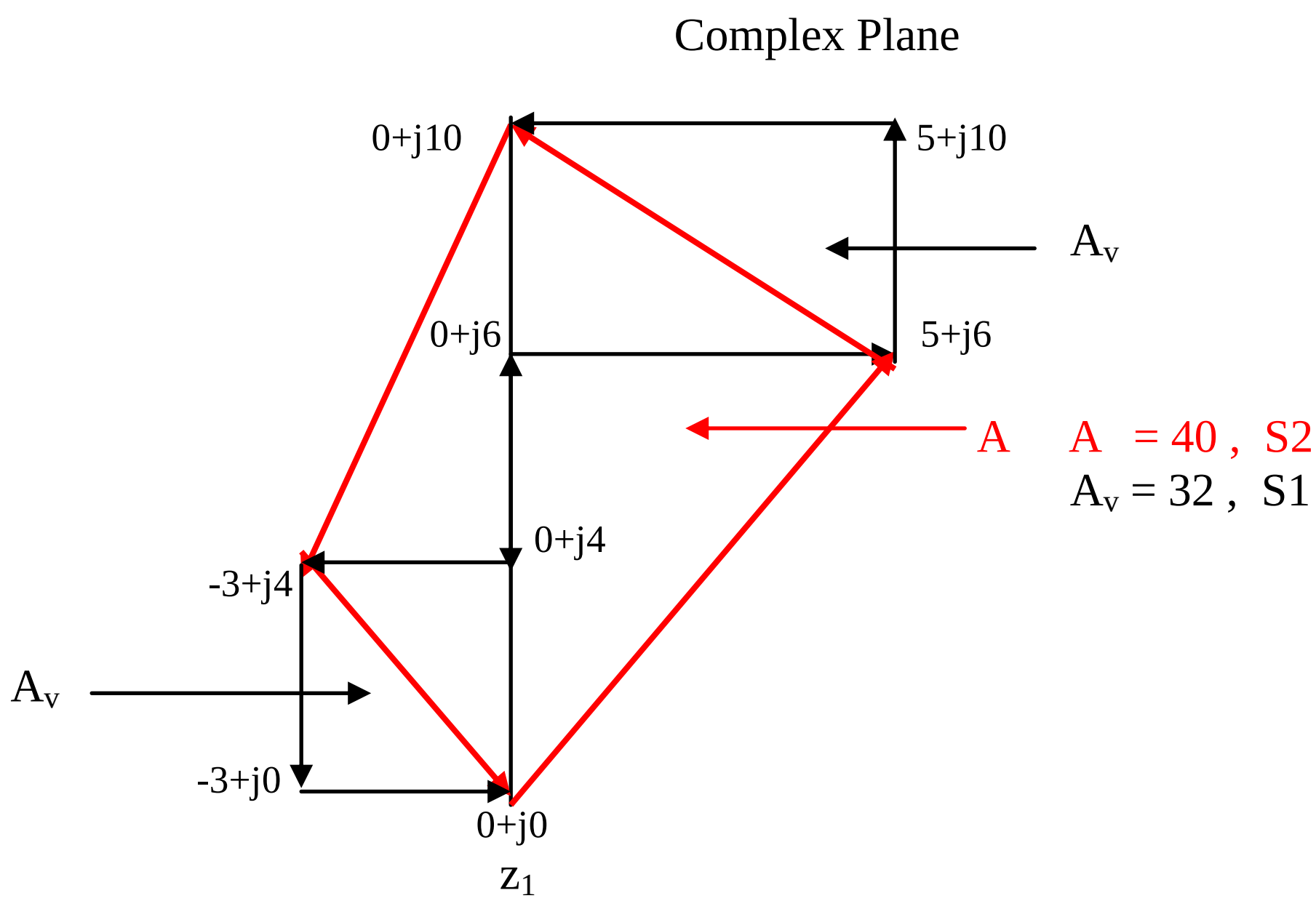

Note that the initial vector direction from the starting point point, $z_1$, is shown to be vertical then counterclockwise. If the vector directions are reversed, the initial vector direction from the starting point point, $z_1$, is horizontal then clockwise.

Example 3.1

Find the area within the black lined S1 closed contour in Diagram 3.6-7 using the equation,

$A = \frac{ja}{2}\sum_{p=1}^{P}(-1)^{p-1}c_p^{\,2}$ . Let the initial direction of summation be horizontal then counterclockwise.

$$A_H = \frac{ja}{2}\sum_{p=1}^{P}(-1)^{p-1}c_p^{\,2} \qquad (3.6\text{-}207)$$

where

$$a = \begin{cases} +1 & \text{for clockwise contour vector summation with initial vector horizontal} \\ -1 & \text{for counterclockwise contour vector summation with initial vector horizontal} \\ -1 & \text{for clockwise contour vector summation with initial vector vertical} \\ +1 & \text{for counterclockwise contour vector summation with initial vector vertical} \end{cases}$$

$c_p$ = the coordinates of the corner points of an S1 complex plane closed contour

P = the number of corner points in an S1 complex plane closed contour

p = 1,2,3,4,…,P

$$A_H = \frac{j(-1)}{2}\sum_{p=1}^{8}(-1)^{p-1}c_p^{\,2} = -\frac{j}{2}[(0+j0)^2 - (5+j0)^2 + (5+j6)^2 - (0+j6)^2 + (0+j10)^2 - (-3+j10)^2$$
$$+ (-3+j4)^2 - (0+j4)^2] \qquad (3.6\text{-}208)$$

$$A_H = -\frac{j}{2}[0 - 25 - 11 + j60 + 36 - 100 + 91 + j60 - 7 - j24 + 16] = -\frac{j}{2}[0 + j96] = 48 \qquad (3.6\text{-}209)$$

$A_H = 48$ Good check ($A_H$ exact = 48) (3.6-210)

<u>Example 3.2</u>

Find the area within the black lined S1 closed contour in Diagram 3.6-8 using the equation,

$A = \frac{ja}{2}\sum_{p=1}^{P}(-1)^{p-1}c_p^{\,2}$ . Let the initial direction of summation be vertical then counterclockwise.

$$A_v = \frac{ja}{2}\sum_{p=1}^{P}(-1)^{p-1}c_p^{\,2} \qquad (3.6\text{-}211)$$

where

$$a = \begin{cases} +1 & \text{for clockwise contour vector summation with initial vector horizontal} \\ -1 & \text{for counterclockwise contour vector summation with initial vector horizontal} \\ -1 & \text{for clockwise contour vector summation with initial vector vertical} \\ +1 & \text{for counterclockwise contour vector summation with initial vector vertical} \end{cases}$$

$c_p$ = the coordinates of the corner points of an S1 complex plane closed contour
P = the number of corner points in an S1 complex plane closed contour
p = 1,2,3,4,…,P

$$A_v = \frac{j(+1)}{2}\sum_{p=1}^{8}(-1)^{p-1}c_p^{\,2} = -\frac{j}{2}[(0+j0)^2 - (0+j6)^2 + (5+j6)^2 - (5+j10)^2 + (0+j10)^2 - (0+j4)^2$$
$$+ (-3+j4)^2 - (-3+j0)^2] \qquad (3.6\text{-}212)$$

$$A_v = -\frac{j}{2}[0 + 36 - 11 + j60 + 75 - j100 - 100 + 16 - 7 - j24 - 9] = \frac{j}{2}[0 - j64] = 32 \qquad (3.6\text{-}213)$$

$A_v = 32$ Good check ($A_v$ exact = 32) (3.6-214)

Example 3.3

Find the area within the redlined S2 closed contour in Diagram 3.6-7 and Diagram 3.6-8.

Using Eq 3.6-74, $A = \frac{1}{2}(A_H + A_V)$, and the results of Example 3.6-1 and Example 3.6-2, the area within the redlined S2 closed contour in Diagram 3.6-7 and Diagram 3.6-8 can be calculated.

From Example 3.6-1

$$A_H = 48 \qquad (3.6\text{-}215)$$

From Example 3.6-2

$$A_v = 32 \qquad (3.6\text{-}216)$$

Then

$$A = \frac{1}{2}(A_H + A_V) = \frac{1}{2}(48 + 32) = 40 \qquad (3.6\text{-}217)$$

$A = 40$ Good check ($A$ exact $= 40$) (3.6-218)

Example 3.4

Find the area within the black lined S1 closed contour in Diagram 3.6-7 using the equation, $A = a\sum_{n=1}^{\frac{N}{2}} \bar{z}_{2n-1}\Delta y_{2n-1}$ . Let the initial direction of summation be vertical then clockwise.

$$A = a\sum_{n=1}^{\frac{N}{2}} \bar{z}_{2n-1}\Delta y_{2n-1} \qquad (3.6\text{-}219)$$

where

$$a = \begin{cases} +1 & \text{for clockwise contour vector summation with initial vector horizontal} \\ -1 & \text{for counterclockwise contour vector summation with initial vector horizontal} \\ -1 & \text{for clockwise contour vector summation with initial vector vertical} \\ +1 & \text{for counterclockwise contour vector summation with initial vector vertical} \end{cases}$$

$$A = (-1)\sum_{n=1}^{4} \overline{z}_{2n-1}\Delta y_{2n-1} = -\frac{1}{2}[(0+j4)(4)+(-6+j14)(6)+(0+j16)(-4)+(10+j6)(-6)] = -\frac{1}{2}(-96) = 48 \quad (3.6\text{-}220)$$

$$A = -\frac{j}{2}[0 + j16 - 36 + j84 + 0 - j64 - 60 - j36] = -\frac{j}{2}[-96 + j0] = 48 \quad (3.6\text{-}221)$$

A = 48 Good check (A exact = 48) (3.6-222)

Example 3.5

Find the area within the redlined S2 closed contour in Diagram 3.6-8 using the equation,

$A = -b\sum_{n=1}^{N} \overline{y}_n\Delta x_n$ . Let the direction of summation be counterclockwise.

$$A = -b\sum_{n=1}^{N} \overline{y}_n\Delta x_n \quad (3.6\text{-}223)$$

where

$$b = \begin{cases} -1 & \text{For clockwise direction} \\ +1 & \text{For counterclockwise direction} \end{cases}$$

$$A = -(+1)\sum_{n=1}^{4} \overline{y}_n\Delta x_n = -\frac{1}{2}[(0+6)(5)+(6+10)(-5)+(10+4)(-3)+(4+0)(3)] \quad (3.6\text{-}224)$$

$$A = -\frac{1}{2}(30 - 80 - 42 + 12) = -\frac{1}{2}(-80) = 40 \quad (3.6\text{-}225)$$

A = 40 Good check (A exact = 40) (3.6-226)

Example 3.6

Find the area within the redlined S2 closed contour in Diagram 3.6-8 using $A = jb\sum_{n=1}^{N} \overline{z}_n\Delta x_n$ .

Let the direction of summation be clockwise.

$$A = jb\sum_{n=1}^{N} \bar{z}_n \Delta x_n \tag{3.6-227}$$

where

$$b = \begin{bmatrix} -1 & \text{For clockwise direction} \\ +1 & \text{For counterclockwise direction} \end{bmatrix}$$

$$A = j(-1)\sum_{n=1}^{4} \bar{z}_n \Delta x_n = -\frac{j}{2}[(-3+j4)(-3)+(-3+j14)(3)+(5+j16)(5)+(5+j6)(-5)] \tag{3.6-228}$$

$$A = -\frac{j}{2}[9 - j12 - 9 + j42 + 25 + j80 - 25 - j30] = -\frac{j}{2}(j80) = 40 \tag{3.6-229}$$

$$A = 40 \quad \text{Good check} \quad (A \text{ exact} = 40) \tag{3.6-230}$$

Example 3.7

Find the area of the S3 closed contour, $z = 2e^{j\theta}$ where $0 \le \theta < 2\pi$, using the equation, $A = \frac{jb}{2}\oint_c z\, dz^*$. Let the direction of summation be counterclockwise.

$$A = \frac{jb}{2}\oint_c z\, dz^*. \tag{3.6-231}$$

where

$$b = \begin{bmatrix} -1 & \text{For clockwise direction} \\ +1 & \text{For counterclockwise direction} \end{bmatrix}$$

$$z = 2e^{j\theta}, \quad 0 \le \theta < 2\pi \tag{3.6-232}$$

$$z^* = 2e^{-j\theta} \tag{3.6-233}$$

$$dz^* = -j2e^{-j\theta}d\theta \tag{3.6-234}$$

Substitute Eq 3.6-232 and Eq 3.6-234 into Eq 3.6-231

$$A = \frac{jb}{2}\oint_c z\, dz^* = \frac{j(+1)}{2}\int_0^{2\pi} 2e^{j\theta}(-j2e^{-j\theta})d\theta = 2\int_0^{2\pi} d\theta = 4\pi \tag{3.6-235}$$

$$A = 4\pi \quad \text{Good Check} \quad (A \text{ exact} = \pi[2]^2 = 4\pi) \tag{3.6-236}$$

Example 3.8

Find the area of the S3 closed contour, $z = 2e^{j\theta}$ where $0 \geq \theta > -2\pi$, using the equation $A = b \oint_c x\,dy$.

Let the direction of summation be clockwise.

$$A = b \oint_c x\,dy. \qquad (3.6\text{-}237)$$

where

$$b = \begin{cases} -1 & \text{For clockwise direction} \\ +1 & \text{For counterclockwise direction} \end{cases}$$

where $0 \geq \theta > -2\pi$

$$z = x + jy = 2e^{j\theta} = 2\cos\theta + j2\sin\theta \qquad (3.6\text{-}236)$$

$$x = 2\cos\theta \qquad (3.6\text{-}237)$$

$$y = 2\sin\theta \qquad (3.6\text{-}238)$$

$$dy = 2\cos\theta d\theta \qquad (3.6\text{-}239)$$

Substituting Eq 3.6-237 and Eq 3.6-239 into Eq 3.6-237

$$A = (-1) \int_0^{-2\pi} 2\cos\theta(2\cos\theta d\theta) = -4 \int_0^{-2\pi} \cos^2\theta d\theta = -4 \int_0^{-2\pi} \frac{1+\cos 2\theta}{2} d\theta = -2 \int_0^{-2\pi} d\theta - 2 \int_0^{-2\pi} \cos 2\theta d\theta \qquad (3.6\text{-}240)$$

$$A = -2(-2\pi) + \sin 2\theta \Big|_0^{-2\pi} = 4\pi \qquad (3.6\text{-}241)$$

$$A = 4\pi \qquad \text{Good Check} \qquad (A \text{ exact} = \pi[2]^2 = 4\pi) \qquad (3.6\text{-}242)$$

# CHAPTER 4

# The Solution of Discrete Calculus Differential Difference Equations

## Section 4.1: Discrete Calculus Differential Difference Equation Solution Overview

In the previous chapters, the mathematical methodology called Interval Calculus has been developed. It is a discrete calculus with many similarities to Calculus. Its notation is similar, its operations are similar, and its application to the solution of mathematical problems is similar. In fact, this Interval Calculus becomes Calculus if the independent variable increment, Δx, approaches zero. Where Calculus is applied to the solution of differential equations and areas involving functions of a continuous variable, x, Interval Calculus is applied to the solution of differential difference equations and areas involving functions of a discrete independent variable, x.

In the discussion which follows, the terms Interval Calculus and Discrete Calculus may be used interchangeably. Interval Calculus is Discrete Calculus. However, it is a discrete calculus which uses the particular notation and functions previous presented. Refer to the Appendix for tables, definitions, and listings.

In the previous chapters many problems have been solved for demonstative purposes. However, an in-depth description of the use of Interval Calculus in the solution of differential difference equations has not yet been presented. This chapter, Chapter 4, will describe four methods used for the solution of Discrete Calculus differential difference equations. The first is the solution of differential difference equations using functions with undetermined coeffients. The second method is the $K_{\Delta x}$ Transform method. The third method is a method which uses related functions and the forth method is the method of variation of parameters. These differential difference equation solution methods will be described in the sections which follow.

## Section 4.2: The solution of differential difference equations using the Method of Undetermined Coefficients

### Definition and Description of the Interval Calculus Differential Difference Equation

The differential difference equation is similar to the Calculus differential equation. However, the derivatives are discrete derivatives and the functions are Interval Calculus discrete functions.

The differential difference equation is:

$$D_{\Delta x}{}^{n}f(x) + a_{n-1}D_{\Delta x}{}^{n-1}f(x) + a_{n-2}D_{\Delta x}{}^{n-2}f(x) + \ldots + a_1D_{\Delta x}f(x) + a_0f(x) = Q(x) \qquad (4.2\text{-}1)$$

where

$n = 0,1,2,3,\ldots$

The discrete derivatives are based on the discrete calculus differential difference operation,

$$D_{\Delta x}f(x) = \frac{f(x+\Delta x) - f(x)}{\Delta x} \qquad (4.2\text{-}2)$$

where

$f(x)$ = function of x

$\Delta x$ = x increment

$x = x_0 + m\Delta x$ , $m$ = integers

$x_0$ = initial value of x

$D_{\Delta x}$ = Interval Calculus differential operator

Successive discrete differentiation of Eq 4.2-2 yields the following relationships:

$$D_{\Delta x}^{2}\,[f(x)] = D_{\Delta x}[D_{\Delta x}f(x)] = \frac{\dfrac{f(x+2\Delta x) - f(x+\Delta x)}{\Delta x} - \dfrac{f(x+\Delta x) - f(x)}{\Delta x}}{\Delta x}$$

$$D_{\Delta x}^{2}\,[f(x)] = \frac{(x+2\Delta x) - 2f(x+\Delta x) + f(x)}{\Delta x^{2}} \qquad (4.2\text{-}3)$$

$$D_{\Delta x}^{3}\,[f(x)] = D_{\Delta x}\,[D_{\Delta x}^{2}\,[f(x)]] = \frac{\dfrac{f(x+3\Delta x) - 2f(x+2\Delta x) + f(x+\Delta x)}{\Delta x^{2}} - \dfrac{f(x+2\Delta x) - 2f(x+\Delta x) + f(x)}{\Delta x^{2}}}{\Delta x}$$

$$D_{\Delta x}^{3}\,[f(x)] = \frac{(x+3\Delta x) - 3f(x+2\Delta x) + 3f(x+\Delta x) - f(x)}{\Delta x^{3}} \qquad (4.2\text{-}4)$$

The coefficients of Eq 4.2-2 thru Eq 4.2-4 are recognized as the Binomial Expansion Coefficients. These coefficients are generated by the expansion of the function, $(1-x)^n$ where n = 0,1,2,3…, or the following formula.

$$_{n}C_{m} = \frac{n!}{m!(n-m)!}\ , \quad m = 0,1,2,3,\ \dots\ , n \qquad (4.2\text{-}5)$$

where
n = 0,1,2,3,…
$m \le n$

From Eq 4.2-2 thru 5.2-5 a general expression for the nth discrete derivative of a function, f(x), can be derived.

$$D_{\Delta x}^{n}\,[f(x)] = \frac{1}{\Delta x^{n}} \sum_{m=0}^{n} (-1)^{m}\,{}_{n}C_{m}\, f(x+[n-m]\Delta x) \qquad (4.2\text{-}6)$$

where

$_{n}C_{m} = \dfrac{n!}{m!(n-m)!}$ , the combinations of n things take m at a time

f(x) = function of x

$D_{\Delta x}^{n}$ = discrete nth derivative operator with an x increment of $\Delta x$

$\Delta x$ = x increment

$x = x_0 + m\Delta x$ , m = integers

$x_0$ = initial value of x

n = the order of the discrete derivative

Rewriting Eq 4.2-1

$D_{\Delta x}{}^{n}f(x) + a_{n-1}D_{\Delta x}{}^{n-1}f(x) + a_{n-2}D_{\Delta x}{}^{n-2}f(x) + \ldots + a_1D_{\Delta x}f(x) + a_0f(x) = Q(x)$

where

$n = 0,1,2,3,\ldots$

$Q(x),f(x)$ = discrete functions of x

$x = x_0 + m\Delta x$ , m = integers

$x_0$ = initial value of x

$a_p$ = real value constant coefficients

$p = 0,1,2,3,\ldots, n-1$

$D_{\Delta x}{}^{n}f(x)$ = nth discrete derivative of f(x)

Given the constant coefficients and the function, Q(x), a method for solving the differential difference equation, Eq 4.2-1, for the function, f(x), will now be sought.

For demonstration, consider Eq 4.2-1 where the order, n, of the equation is two (i.e. n = 2).

$$D_{\Delta x}{}^{2}f(x) + a_1D_{\Delta x}f(x) + a_0 = Q(x) \qquad (4.2\text{-}7)$$

Substituting Eq 4.2-2 and Eq 4.2-3 into Eq 4.2-7

$$\frac{(x+2\Delta x) - 2f(x+\Delta x) + f(x)}{\Delta x^2} + a_1[\,\frac{f(x+\Delta x) - f(x)}{\Delta x}\,] + a_0f(x) = Q(x) \qquad (4.2\text{-}8)$$

Solving Eq 4.2-8 for $f(x+2\Delta x)$ and simplifying

$f(x+2\Delta x) - 2f(x+\Delta x) + f(x) + a_1\Delta xf(x+\Delta x) - a_1\Delta xf(x) + a_0\Delta x^2f(x) = \Delta x^2Q(x)$

$f(x+2\Delta x) + [-2+a_1\Delta x]f(x+\Delta x) + [1-a_1\Delta x+a_0\Delta x^2]f(x) = \Delta x^2Q(x)$

$f(x+2\Delta x) - A_1f(x+\Delta x) - A_0f(x) = BQ(x)$

where

$A_1 = 2 - a_1\Delta x$

$A_0 = -1 + a_1\Delta x - a_0\Delta x^2$

$B = \Delta x^2$

$$f(x+2\Delta x) = A_1f(x+\Delta x) + A_0f(x) + BQ(x) \qquad (4.2\text{-}9)$$

Note that Eq 4.2-9 is a difference equation equivalent to Eq 4.2-7, the differential equation from which it was derived. Note, also, that Eq 4.2-9 is a recursive equation which is the solution to Eq 4.2-7. To use Eq 4.2-9 to find the value of f(x) at all values of x, $x = x_0 + m\Delta x$ , m=integers, two initial conditions, $f(x_0)$ and $f(x_0+\Delta x)$ or equivalently $f(x_0)$ and $D_{\Delta x}(x_0)$ where $D_{\Delta x}(x_0) = \dfrac{f(x_0+\Delta x) - f(x_0)}{\Delta x}$ , must first be known.

One can extrapolate the results of this second order (i.e. n = 2) demonstration example to a differential difference equation of any order in the form of Eq 4.2-10 below. The more general conclusions are as follows:

1) For any differential difference equation of the form:

$D_{\Delta x}{}^{n}f(x) + a_{n-1}D_{\Delta x}{}^{n-1}f(x) + a_{n-2}D_{\Delta x}{}^{n-2}f(x) + \ldots + a_1D_{\Delta x}f(x) + a_0f(x) = Q(x)$

(From Eq 4.2-1)

there is an equivalent equation obtained by substituting Eq 4.2-6 into Eq 4.2-1

$$f(x+n\Delta x) = A_{n-1}f(x+[n-1]\Delta x + A_{n-2}f(x+[n-2]\Delta x + A_{n-3}f(x+[n-3]\Delta x + \ldots + A_1f(x+\Delta x) + A_0f(x) + BQ(x) \quad (4.2\text{-}10)$$

where

$a_p$, $A_p$, B = real value constants

p = 0,1,2,3,…, n-1

n = the order of the differential difference equation

$\Delta x$ = x increment

2) Eq 4.2-10 is not only an equivalent equation for Eq 4.2-1 but also a recursive solution for Eq 4.2-1.

3) As can be seen from Eq 4.2-10, any solution of Eq 4.2-1, f(x), for all $x = x_0 + m\Delta x$, m=integers must allow for the specification of n independent initial conditions. These initial conditions may be in the form:

$f(x_0)$, $f(x_0+\Delta x)$, $f(x_0+2\Delta x)$, …

or

$f(x_0)$, $D_{\Delta x}f(x_0)$, $D_{\Delta x}{}^{2}f(x_0)$, …

or

A combination of the above so long as the initial conditions are independent of each other

For example, $f(x_0)$, $f(x_0+\Delta x)$, and $D_{\Delta x}f(x_0)$ are not independent initial conditions.

One initial condition is expressible in terms of two of the other initial conditions, $D_{\Delta x}f(x_0) = \dfrac{f(x_0+\Delta x)-f(x_0)}{\Delta x}$ .

As can be seen from Eq 4.2-10, a recursion method is an option to solve the differential difference equation, Eq 4.2-1, for f(x). However, if $\Delta x$ is small or the values of x are large, the recursion method to find f(x) may be impractical. In this case, another method for finding a function solution for f(x) in equation Eq 4.2-1 would be required. Such a method will now be discussed.

Note – Observation 3, above, is a very important requirement for obtaining a general solution for Eq 4.2-1.

**The discrete function, $(1+a\Delta x)^{\frac{x}{\Delta x}}$**

Consider the discrete function,

$$g(x) = (1+a\Delta x)^{\frac{x}{\Delta x}} \qquad (4.2\text{-}11)$$

where

$\Delta x$ = x increment

a = index of the function, g(x) (See note 1 below)

Notes - 1) The constant, a, within the function, g(x), is of considerable importance as will soon be seen. For this reason, it has been given a name. This name is "index " for "internal exponent".

2) $\lim_{\Delta x\to 0} g(x) = \lim_{\Delta x\to 0} (1+a\Delta x)^{\frac{x}{\Delta x}} = e^{ax} = (e^a)^x$

For $\lim_{\Delta x\to 0} g(x)$, the index of g(x) becomes the exponent of e, the base of the natural log.

3) Referring to the function, $g(x) = \lim_{\Delta x\to 0} (1+a\Delta x)^{\frac{x}{\Delta x}}$, the value, $(1+a\Delta x)^{\frac{1}{\Delta x}}$, is the base of the function, g(x). This base is raised to the power of x.

The function, g(x), has a very important property which is exceedingly useful in the solution of differential difference equations.

Consider the discrete derivative of $g(x) = (1+a\Delta x)^{\frac{x}{\Delta x}}$

$$D_{\Delta x}g(x) = \frac{g(x+\Delta x)-g(x)}{\Delta x} = \frac{(1+a\Delta x)^{\frac{x+\Delta x}{\Delta x}} - (1+a\Delta x)^{\frac{x}{\Delta x}}}{\Delta x} = \frac{[(1+a\Delta x)-1](1+a\Delta x)^{\frac{x}{\Delta x}}}{\Delta x}$$

$$D_{\Delta x}g(x) = a(1+a\Delta x)^{\frac{x}{\Delta x}} \qquad (4.2\text{-}12)$$

Consider the second derivative of g(x)

$$D_{\Delta x}[D_{\Delta x}g(x)] = D_{\Delta x}{}^2g(x) = \frac{a(1+a\Delta x)^{\frac{x+\Delta x}{\Delta x}} - a(1+a\Delta x)^{\frac{x}{\Delta x}}}{\Delta x} = \frac{a[(1+a\Delta x)-1](1+a\Delta x)^{\frac{x}{\Delta x}}}{\Delta x}$$

$$D_{\Delta x}{}^2g(x) = a^2(1+a\Delta x)^{\frac{x}{\Delta x}} \qquad (4.2\text{-}13)$$

.

.

In general

$$D_{\Delta x}{}^n g(x) = D_{\Delta x}{}^n (1+a\Delta x)^{\frac{x}{\Delta x}} = a^n(1+a\Delta x)^{\frac{x}{\Delta x}} \qquad (4.2\text{-}14)$$

Consider the function, f(x), of the differential difference equation, Eq 4.2-1, to be composed of two parts. Let

$$f(x) = f_C + f_P \tag{4.2-15}$$

where

f(x) = the general solution to the differential difference equation

$f_C(x)$ = the complementary solution to the differential difference equation and the general solution to its related homogeneous equation. (See Eq 4.2-17 below.)

$f_P(x)$ = the particular solution to the differential difference equation

Differential Difference Equation Solutions

1) The Differential Difference Equation General Solution

$f(x) = f_C(x) + f_P(x)$

$$D_{\Delta x}{}^{n}f(x) + a_{n-1}D_{\Delta x}{}^{n-1}f(x) + a_{n-2}D_{\Delta x}{}^{n-2}f(x) + \ldots + a_1D_{\Delta x}f(x) + a_0f(x) = Q(x) \tag{4.2-16}$$

2) The Differential Difference Equation Complementary Solution and also the General Solution to the Differential Difference Equation Related Homogeneous Equation where Q(x) is given a zero value

$f_C(x)$

$$D_{\Delta x}{}^{n}f_C(x) + a_{n-1}D_{\Delta x}{}^{n-1}f_C(x) + a_{n-2}D_{\Delta x}{}^{n-2}f_C(x) + \ldots + a_1D_{\Delta x}f_C(x) + a_0f_C(x) = 0 \tag{4.2-17}$$

3) The Differential Difference Equation Particular Solution

$f_P(x)$

$$D_{\Delta x}{}^{n}f_P(x) + a_{n-1}D_{\Delta x}{}^{n-1}f_P(x) + a_{n-2}D_{\Delta x}{}^{n-2}f_P(x) + \ldots + a_1D_{\Delta x}f_P(x) + a_0f_P(x) = Q(x) \tag{4.2-18}$$

**Derivation of the Complementary Solution of a Differential Difference Equation with Unique Roots**

Find the complementary solution, $f_C(x)$, of the differential difference equation, Eq 4.2-16 when the characteristic polynomial, h(r), has no multiple roots.

Investigate the function, $k(1+r\Delta x)^{\frac{x}{\Delta x}}$, to see if it satisfies the related homogeneous differential difference equation, Eq 4.2-17. Let $f_C = k(1+r\Delta x)^{\frac{x}{\Delta x}}$. The values k and r are real or complex constants.

Substituting $f_C = k(1+r\Delta x)^{\frac{x}{\Delta x}}$ into Eq 4.2-17

$$kr^n(1+r\Delta x)^{\frac{x}{\Delta x}} + ka_{n-1}r^{n-1}(1+r\Delta x)^{\frac{x}{\Delta x}} + ka_{n-2}\,r^{n-2}(1+r\Delta x)^{\frac{x}{\Delta x}} + \ldots + ka_1r(1+r\Delta x)^{\frac{x}{\Delta x}} + ka_0(1+r\Delta x)^{\frac{x}{\Delta x}} = 0 \tag{4.2-19}$$

$$k(1+r\Delta x)^{\frac{x}{\Delta x}} [\, r^n + a_{n-1}r^{n-1} + a_{n-2}\, r^{n-2} + \ldots + a_1 r + a_0 \,] = 0 \qquad (4.2\text{-}20)$$

From Eq 4.2-20

$$h(r) = r^n + a_{n-1}r^{n-1} + a_{n-2}\, r^{n-2} + \ldots + a_1 r + a_0 = 0 \qquad 4.2\text{-}21)$$

The above polynomial, h(r), is called the characteristic polynomial associated with the differential difference equations, Eq 4.2-16 thru Eq 4.2-18. From this polynomial, a function which satisfies the related homogenous equation, Eq 4.17 can be found. Referring to Eq 4.2-21, it is observed that there are n values of r which will satisfy the related homogeneous equation, Eq 4.17.

Eq 4.2-21 has n roots

Factoring Eq 4.2-21

$$h(r) = (r-r_n)\,(r-r_{n-1})\,(r-r_{n-2})\ldots\,(r-r_2)\,(r-r_1) = 0 \qquad (4.2\text{-}22)$$

The n roots of h(r) are $r_{1,}\, r_2,\, r_3,\, \ldots,\, r_{n-1,}\, r_n$ . These roots may be real or complex values.

There are then n $k_p(1+r_p\Delta x)^{\frac{x}{\Delta x}}$ functions where p = 1,2,3,…,n that satisfy the related homogeneous equation, 4.2-17. Since each of these n functions has a constant available for the specification of an initial condition, the linear combination of all n of these functions provides a general solution for the homogeneous equation, 4.2-17. This is in accordance with Observation 3 of the general conclusions previously specified that requires a general solution to a differential difference equation to have the ability to have specified n initial conditions. Of course, these initial conditions must be independent of one another. This will be discussed in the following paragraphs .

Then

**The complementary solution of the differential difference equation, Eq 4.2-16, where all roots of the characteristic polynomial, h(r), are unique is:**

$$f_C(x) = k_1(1+r_1\Delta x)^{\frac{x}{\Delta x}} + k_2(1+r_2\Delta x)^{\frac{x}{\Delta x}} + k_3(1+r_3\Delta x)^{\frac{x}{\Delta x}} + \ldots + k_{n-1}(1+r_{n-1}\Delta x)^{\frac{x}{\Delta x}} + k_n(1+r_n\Delta x)^{\frac{x}{\Delta x}}$$

**or**

$$f_C(x) = k_1 e_{\Delta x}(r_1,x) + k_2 e_{\Delta x}(r_2,x) + k_3 e_{\Delta x}(r_3,x) + \ldots + k_{n-1} e_{\Delta x}(r_{n-1},x) + k_n e_{\Delta x}(r_n,x) \qquad (4.2\text{-}23)$$

**where**

**$r_p$, p=1,2,3,…,n are the roots of the characteristic polynomial, h(r)**

$$h(r) = r^n + a_{n-1}r^{n-1} + a_{n-2}\, r^{n-2} + \ldots + a_1 r + a_0 = (r-r_n)\,(r-r_{n-1})\,(r-r_{n-2})\ldots\,(r-r_2)\,(r-r_1) = 0$$

**None of the roots have the same value.**

**The roots, $r_1$ thru $r_n$ may be real or complex conjugate pair values.**

Individually, each term of the function, $f_C(x)$, satisfies the related homogeneous equation, Eq 4.2-17, but only their sum fulfills the requirements of being a complementary solution. Firstly, the function, $f_C(x)$, satisfies the related homogeneous equation. Secondly, it satisfies the general solution requirement that there be allowed n independent initial conditions to be specified. The n constants, $k_1$ thru $k_n$, satisfy this condition. Thus, the function, $f_C(x)$, is the general solution to the related homogeneous equation and is designated as the complementary solution. If there are roots of multiplicity two or more, terms could be combined to form a single term. No longer would there be n independent k constants to represent n initial conditions. In this case,
Eq 4.2-23, would not be the sought after complementary solution. General solutions do exist for homogeneous differential difference equations with roots of multiplicity two or more. The derivation of these general solutions follow.

## **Derivation of the Complementary Solution of a Differential Difference Equation with Multiple Roots**

Find the complementary solution, $f_C(x)$, of the differential difference equation, Eq 4.2-16 when the characteristic polynomial, h(r), has a root of multiplicity two.

From Eq 4.2-17, the related homogeneous equation of the differential difference equation is:

$$D_{\Delta x}{}^{n} f_C(x) + a_{n-1}D_{\Delta x}{}^{n-1} f_C(x) + a_{n-2}D_{\Delta x}{}^{n-2} f_C(x) + \ldots + a_1 D_{\Delta x} f_C(x) + a_0 f_C(x) = 0 \qquad (4.2\text{-}24)$$

where

$$h(r) = r^n + a_{n-1}r^{n-1} + a_{n-2}\, r^{n-2} + \ldots + a_2 r^2 + a_1 r + a_0 = (r-r_n)(r-r_{n-1})(r-r_{n-2})\ldots(r-r_3)(r-a)(r-a) = 0 \qquad (4.2\text{-}25)$$

Two of the roots have the same value, $r_1 = r_2 = a$
$r_p$, $p=1,2,3,\ldots,n$ are the roots of the characteristic polynomial, h(r)
The roots, $r_1$ thru $r_n$ may be real or complex conjugate pair values.

Change the form of Eq 4.2-24 to represent the presence a double root with a value, r = a.

$$(D_{\Delta x}-r_n)(D_{\Delta x}-r_{n-1})(D_{\Delta x}-r_{n-2})\ldots(D_{\Delta x}-a)(D_{\Delta x}-a)\, f_C(x) = 0 \qquad (4.2\text{-}26)$$

Investigate to see if the function, $v(x)e_{\Delta x}(a,x)$, satisfies Eq 4.2-26 where v(x) is some to be determined function of x.

Substitute the function, $v(x)e_{\Delta x}(a,x)$, into Eq 4.2-26 for $f_C(x)$

$$(D_{\Delta x}-r_n)(D_{\Delta x}-r_{n-1})(D_{\Delta x}-r_{n-2})\ldots(D_{\Delta x}-r_3)[(D_{\Delta x}{}^2 - 2aD_{\Delta x} + a^2)v(x)e_{\Delta x}(a,x)] = 0 \qquad (4.2\text{-}27)$$

If

$$(D_{\Delta x}{}^2 - 2aD_{\Delta x} + a^2)v(x)e_{\Delta x}(a,x) = (D_{\Delta x}-a)^2 v(x)e_{\Delta x}(a,x)] = 0 \qquad (4.2\text{-}28)$$

the function $f_C(x) = v(x)e_{\Delta x}(a,x)$ would satisfy the homogeneous equation, Eq 4.2-27, which is a form of Eq 4.2-24 where two roots of h(r) have the same value.

Evaluate Eq 4.2-28.

These Interval Calculus functions and equations will be used in the evaluations and derivations which follow.

$v(x)$ = function of x

$$e_{\Delta x}(r,x) = (1+r\Delta x)^{\frac{x}{\Delta x}}$$

$$e_{\Delta x}(r,x+n\Delta x) = (1+r\Delta x)^n\, e_{\Delta x}(r,x)$$

$$[x]_{\Delta x}^{0} = 1$$

$$[x]_{\Delta x}^{h} = x(x-\Delta x)(x-2\Delta x)(x-3\Delta x)\ldots(x-[h-1]\Delta x)\ , \quad h = 1,2,3,\ldots$$

$$D_{\Delta x}{}^{n}\, e_{\Delta x}(r,x) = r^n e_{\Delta x}(r,x) = r^n(1+r\Delta x)^{\frac{x}{\Delta x}}$$

$$D_{\Delta x}[x]_{\Delta x}^{n} = \begin{cases} 0 & \text{for } n = 0 \\ n[x]_{\Delta x}^{n-1} & \text{for } n = 1,2,3,\ldots \end{cases}$$

$$D_{\Delta x}[u(x)v(x)] = [D_{\Delta x}u(x)]v(x) + [D_{\Delta x}v(x)]u(x+\Delta x)$$

Find $D_{\Delta x}{}^{0}[v(x)e_{\Delta x}(a,x)]$

$$D_{\Delta x}{}^{0}[v(x)e_{\Delta x}(a,x)] = v(x)e_{\Delta x}(a,x) \qquad (4.2\text{-}29)$$

Find $D_{\Delta x}{}^{1}[v(x)e_{\Delta x}(a,x)]$

Using the function product derivative equation

$$D_{\Delta x}{}^{1}[v(x)e_{\Delta x}(a,x)] = ae_{\Delta x}(a,x)v(x) + D_{\Delta x}v(x)e_{\Delta x}(a,x+\Delta x) \qquad (4.2\text{-}30)$$

Find $D_{\Delta x}{}^{2}[v(x)e_{\Delta x}(a,x)]$

$$D_{\Delta x}{}^{2}[v(x)e_{\Delta x}(a,x)] = D_{\Delta x}[D_{\Delta x}{}^{1}[v(x)e_{\Delta x}(a,x)]]$$

Using the function product derivative equation again

$$D_{\Delta x}{}^{2}[v(x)e_{\Delta x}(a,x)] = a^2e_{\Delta x}(a,x)v(x) + aD_{\Delta x}v(x)e_{\Delta x}(a,x+\Delta x) + ae_{\Delta x}(a,x+\Delta x)D_{\Delta x}v(x) + D_{\Delta x}{}^{2}v(x)e_{\Delta x}(a,x+2\Delta x) \qquad (4.2\text{-}31)$$

Substituting Eq 4.2-29 thru 4.2-31 into Eq 4.2-28

$(D_{\Delta x}–a)^2v(x)e_{\Delta x}(a,x) = [v(x)e_{\Delta x}(a,x)] –2a[ae_{\Delta x}(a,x)v(x) + D_{\Delta x}v(x)e_{\Delta x}(a,x+\Delta x)] + [a^2e_{\Delta x}(a,x)v(x) + aD_{\Delta x}v(x)e_{\Delta x}(a,x+\Delta x) + ae_{\Delta x}(a,x+\Delta x)D_{\Delta x}v(x) + D_{\Delta x}{}^2v(x)e_{\Delta x}(a,x+2\Delta x)] = 0$

Simplifying using the Interval Calculus equations listed above.

$$(D_{\Delta x}–a)^2v(x)e_{\Delta x}(a,x) = [a^2v(x)–2a^2v(x)–2a(1+a\Delta x)D_{\Delta x}v(x)+a^2v(x)+a(1+a\Delta x)D_{\Delta x}v(x)+a(1+a\Delta x)D_{\Delta x}v(x) +(1+a\Delta x)^2D_{\Delta x}{}^2v(x)]e_{\Delta x}(a,x) = 0 \qquad (4.2\text{-}32)$$

After canceling terms in Eq 4.2-32

$$(D_{\Delta x}–a)^2v(x)e_{\Delta x}(a,x) = (D_{\Delta x}{}^2 – 2aD_{\Delta x} + a^2)v(x)e_{\Delta x}(a,x) = (1+a\Delta x)^2e_{\Delta x}(a,x)D_{\Delta x}{}^2v(x) = 0 \qquad (4.2\text{-}33)$$

From Eq 4.2-27, Eq 4.2-28 and Eq 4.2-33

$$(D_{\Delta x}-r_n)(D_{\Delta x}-r_{n-1})(D_{\Delta x}-r_{n-2})\ldots(D_{\Delta x}-r_3)[(D_{\Delta x}{}^2 – 2aD_{\Delta x} + a^2)v(x)e_{\Delta x}(a,x)] =$$
$$(D_{\Delta x}-r_n)(D_{\Delta x}-r_{n-1})(D_{\Delta x}-r_{n-2})\ldots(D_{\Delta x}-r_3)[(D_{\Delta x}–a)^2v(x)e_{\Delta x}(a,x)] =$$
$$(D_{\Delta x}-r_n)(D_{\Delta x}-r_{n-1})(D_{\Delta x}-r_{n-2})\ldots(D_{\Delta x}-r_3)[(1+a\Delta x)^2e_{\Delta x}(a,x)D_{\Delta x}{}^2v(x)] = 0 \qquad (4.2\text{-}34)$$

To find the general solution for Eq 4.2-34, $D_{\Delta x}{}^2v(x)$ must equal zero and v(x) must introduce two constants, one for each of the two r = a roots, in order to obtain two initial conditions.

Since $D_{\Delta x}{}^2(k_2+k_1[x]_{\Delta x}^{1}) = 0$,

$$v(x) = k_2+k_1[x]_{\Delta x}^{1} \qquad (4.2\text{-}35)$$

where

$k_1,k_2$ = constants

$[x]_{\Delta x}^{1} = x$

From Eq 4.2-34 and Eq 4.2-35

$$(D_{\Delta x}-r_n)(D_{\Delta x}-r_{n-1})(D_{\Delta x}-r_{n-2})\ldots (D_{\Delta x}-r_3)[(D_{\Delta x}{}^2 – 2aD_{\Delta x} + a^2) [v(x)e_{\Delta x}(a,x)] =$$
$$(D_{\Delta x}-r_n)(D_{\Delta x}-r_{n-1})(D_{\Delta x}-r_{n-2})\ldots (D_{\Delta x}-r_3)(D_{\Delta x}–a)(D_{\Delta x}–a)[( k_2+k_1[x]_{\Delta x}^{1})e_{\Delta x}(a,x)] = 0 \qquad (4.2\text{-}36)$$

Then the function,

$$( k_2+k_1[x]_{\Delta x}^{1})e_{\Delta x}(a,x) \qquad (4.2\text{-}37)$$

satisfies the related homogeneous equation, Eq 4.2-17, when two of its characteristic polynomial ,h(r), roots have the same value, a.

Thus

**The complementary solution of a differential difference equation with a characteristic polynomial, h(r), with two roots of equal value, a, is as follows:**

**From the related homogeneous equation:**

$$\mathbf{D_{\Delta x}^{n} f_C(x) + a_{n-1}D_{\Delta x}^{n-1} f_C(x) + a_{n-2}D_{\Delta x}^{n-2} f_C(x) + \ldots + a_1D_{\Delta x} f_C(x) + a_0 f_C(x) = 0}$$

**or**

$$\mathbf{(D_{\Delta x}-r_n)(D_{\Delta x}-r_{n-1})(D_{\Delta x}-r_{n-2})\ldots(D_{\Delta x}-r_3)(D_{\Delta x}-a)(D_{\Delta x}-a) f_C(x) = 0}$$

**The differential difference equation complementary solution, $f_C(x)$, is:**

$$\mathbf{f_C(x) = k_ne_{\Delta x}(r_n,x) + k_{n-1}e_{\Delta x}(r_{n-1},x) + \ldots + k_3e_{\Delta x}(r_3,x) + k_2e_{\Delta x}(a,x) + k_1[x]_{\Delta x}^{1}e_{\Delta x}(a,x)}$$

**or**

$$\mathbf{f_C(x) = k_n(1+r_n\Delta x)^{\frac{x}{\Delta x}} + k_{n-1}(1+r_{n-1}\Delta x)^{\frac{x}{\Delta x}} + \ldots + k_3(1+r_3\Delta x)^{\frac{x}{\Delta x}} + k_2(1+a\Delta x)^{\frac{x}{\Delta x}} + k_1[x]_{\Delta x}^{1}(1+a\Delta x)^{\frac{x}{\Delta x}}}$$

**where**

**n = the order of the homogeneous equation**

**$h(r) = r^n + a_{n-1}r^{n-1} + a_{n-2} r^{n-2} + \ldots + a_2r^2 + a_1r + a_0 = (r-r_n)(r-r_{n-1})(r-r_{n-2})\ldots(r-r_3)(r-a)(r-a) = 0$**

**$r_p$, p=1,2,3,…,n are the roots of the characteristic polynomial, h(r)**

**Two of the roots have the same value, $r_1 = r_2 = a$**

**The roots, $r_1$ thru $r_n$ may be real or complex conjugate pair values.**

**$[x]_{\Delta x}^{0} = 1$**

**$[x]_{\Delta x}^{1} = x$**

Find the complementary solution, $f_C(x)$, of the differential difference equation, Eq 4.2-16 when the characteristic polynomial, h(r), has a root of multiplicity three.

From Eq 4.2-17, the related homogeneous equation of the differential difference equation is:

$$D_{\Delta x}^{n}f_C(x) + a_{n-1}D_{\Delta x}^{n-1}f_C(x) + a_{n-2}D_{\Delta x}^{n-2}f_C(x) + \ldots + a_1D_{\Delta x}f_C(x) + a_0 f_C(x) = 0 \qquad (4.2\text{-}38)$$

where

$$h(r) = r^n + a_{n-1}r^{n-1} + a_{n-2} r^{n-2} + \ldots + a_2r^2 + a_1r + a_0 = (r-r_n)(r-r_{n-1})(r-r_{n-2})\ldots(r-r_4)(r-a)(r-a)(r-a) = 0 \qquad (4.2\text{-}39)$$

Three of the roots have the same value, $r_1 = r_2 = r_3 = a$

$r_p$, p=1,2,3,…,n are the roots of the characteristic polynomial, h(r)

The roots, $r_1$ thru $r_n$ may be real or complex conjugate pair values.

Change the form of Eq 4.2-38 to represent the presence of a triple root with a value, r = a.

$$(D_{\Delta x}-r_n)(D_{\Delta x}-r_{n-1})(D_{\Delta x}-r_{n-2})\ldots(D_{\Delta x}-a)(D_{\Delta x}-a)(D_{\Delta x}-a) f_C(x) = 0 \qquad (4.2\text{-}40)$$

From Eq 4.2-33

$(D_{\Delta x}–a)^2v(x)e_{\Delta x}(a,x) = (D_{\Delta x}{}^2 – 2aD_{\Delta x} + a^2)v(x)e_{\Delta x}(a,x) = (1+a\Delta x)^2e_{\Delta x}(a,x)\ D_{\Delta x}{}^2v(x) = 0$

Introducing another root, r = a

$$(D_{\Delta x}–a)^3v(x)e_{\Delta x}(a,x) = (D_{\Delta x}–a)(D_{\Delta x}{}^2 – 2aD_{\Delta x} + a^2)v(x)e_{\Delta x}(a,x) = (D_{\Delta x}–a)(1+a\Delta x)^2e_{\Delta x}(a,x)D_{\Delta x}{}^2v(x) = 0 \qquad (4.2\text{-}41)$$

From Eq 4.2-41 and Eq 4.2-34

$(D_{\Delta x}-r_n)(D_{\Delta x}-r_{n-1})(D_{\Delta x}-r_{n-2})\ldots(D_{\Delta x}-r_4)(D_{\Delta x}-a)[(D_{\Delta x}{}^2 – 2aD_{\Delta x} + a^2)v(x)e_{\Delta x}(a,x)] =$
$(D_{\Delta x}-r_n)(D_{\Delta x}-r_{n-1})(D_{\Delta x}-r_{n-2})\ldots(D_{\Delta x}-r_4)[(D_{\Delta x}–a)^3v(x)e_{\Delta x}(a,x)] =$
$(D_{\Delta x}-r_n)(D_{\Delta x}-r_{n-1})(D_{\Delta x}-r_{n-2})\ldots(D_{\Delta x}-r_4)[(D_{\Delta x}–a)(1+a\Delta x)^2e_{\Delta x}(a,x)D_{\Delta x}{}^2v(x)] = 0$ (4.2-42)

Eq 4.2-42 is satisfied for:

$$(D_{\Delta x}–a)(1+a\Delta x)^2e_{\Delta x}(a,x)D_{\Delta x}{}^2v(x) = 0 \qquad (4.2\text{-}43)$$

Find the function, v(x), which satisfies Eq 4.2-42
$(D_{\Delta x}–a)^3v(x)e_{\Delta x}(a,x) =$
$(D_{\Delta x}–a)(1+a\Delta x)^2e_{\Delta x}(a,x)D_{\Delta x}{}^2v(x) = (1+a\Delta x)^2[(D_{\Delta x}–a)e_{\Delta x}(a,x)\ D_{\Delta x}{}^2v(x)] =$
$(1+a\Delta x)^2[ae_{\Delta x}(a,x)D_{\Delta x}{}^2v(x) + D_{\Delta x}{}^3v(x)e_{\Delta x}(a,x+\Delta x) – ae_{\Delta x}(a,x)D_{\Delta x}{}^2v(x)] =$
$(1+a\Delta x)^3e_{\Delta x}(a,x)D_{\Delta x}{}^3v(x) = 0$
(4.2-44)

To find the general solution for Eq 4.2-42, $D_{\Delta x}{}^3v(x)$ must equal zero and v(x) must introduce three constants, one for each of the three r = a roots, in order to obtain three initial conditions.

Since $D_{\Delta x}{}^3(k_3+k_2[x]^1_{\Delta x} + k_1[x]^2_{\Delta x}) = 0$,

$$v(x) = k_3+k_2[x]^1_{\Delta x} + k_1[x]^2_{\Delta x} \qquad (4.2\text{-}45)$$

where

$k_1,k_2,k_3$ = constants

$[x]^1_{\Delta x} = x$

$[x]^2_{\Delta x} = x(x-\Delta x)$

Then the function,

$$(k_3+k_2[x]^1_{\Delta x} + k_1[x]^2_{\Delta x})e_{\Delta x}(a,x) \qquad (4.2\text{-}46)$$

satisfies the related homogeneous equation, Eq 4.2-38, when three of its characteristic polynomial ,h(r), roots have the same value, a.

Thus

**The complementary solution of a differential difference equation with a characteristic polynomial, h(r), with three roots of equal value, a, is as follows:**

**From the related homogeneous equation:**

$$D_{\Delta x}^{n} f_C(x) + a_{n-1}D_{\Delta x}^{n-1} f_C(x) + a_{n-2}D_{\Delta x}^{n-2} f_C(x) + \ldots + a_1D_{\Delta x} f_C(x) + a_0 f_C(x) = 0$$

**or**

$$(D_{\Delta x}-r_n)(D_{\Delta x}-r_{n-1})(D_{\Delta x}-r_{n-2})\ldots(D_{\Delta x}-r_4)(D_{\Delta x}-a)(D_{\Delta x}-a)(D_{\Delta x}-a) f_C(x) = 0$$

**The differential difference equation complementary solution, $f_C(x)$, is:**

$$f_C(x) = k_ne_{\Delta x}(r_n,x) + k_{n-1}e_{\Delta x}(r_{n-1},x) + \ldots + k_4e_{\Delta x}(r_4,x) + k_3e_{\Delta x}(a,x) + k_2[x]_{\Delta x}^{1}e_{\Delta x}(a,x) + k_1[x]_{\Delta x}^{2}e_{\Delta x}(a,x)$$

**or**

$$f_C(x) = k_n(1+r_n\Delta x)^{\frac{x}{\Delta x}} + k_{n-1}(1+r_{n-1}\Delta x)^{\frac{x}{\Delta x}} + \ldots + k_4(1+r_4\Delta x)^{\frac{x}{\Delta x}} + k_3(1+a\Delta x)^{\frac{x}{\Delta x}} + k_2[x]_{\Delta x}^{1}(1+a\Delta x)^{\frac{x}{\Delta x}} +$$

$$k_1[x]_{\Delta x}^{2}(1+a\Delta x)^{\frac{x}{\Delta x}}$$

**where**

**n = the order of the homogeneous equation**
**$h(r) = r^n + a_{n-1}r^{n-1} + a_{n-2}r^{n-2} + \ldots + a_2r^2 + a_1r + a_0 = (r-r_n)(r-r_{n-1})(r-r_{n-2})\ldots(r-r_4)(r-a)(r-a)(r-a) = 0$**
**$r_p$, p=1,2,3,…,n are the roots of the characteristic polynomial, h(r)**
**Three of the roots have the same value, $r_1 = r_2 = r_3 = a$**
**The roots, $r_1$ thru $r_n$ may be real or complex conjugate pair values.**

$$[x]_{\Delta x}^{1} = x$$

$$[x]_{\Delta x}^{2}e(a,x) = x(x-\Delta x)$$

From Eq 4.2-33, Eq 4.2-44 and a continuation of the calculation process

In general

For m = 2,3,4,…

$$(D_{\Delta x}-a)^{m}v(x)e_{\Delta x}(a,x) = (1+a\Delta x)^{m}e_{\Delta x}(a,x)D_{\Delta x}^{m}v(x) = 0 \qquad (4.2\text{-}47)$$

Then

$$v(x) = \sum_{p=1}^{m} k_{m-p+1}[x]_{\Delta x}^{p-1} \qquad (4.2\text{-}48)$$

where

$k_1,k_2,k_3\ldots k_m$ = real or complex value constants
m = the multiplicity of the repeated roots of the homogeneous equation characteristic polynomial, h(r)

$m = 2,3,4\ldots$

$[x]_{\Delta x}^{0} = 1$

$[x]_{\Delta x}^{h} = x(x-\Delta x)(x-2\Delta x)(x-3\Delta x)\ldots(x-[h-1]\Delta x)\ , \quad h = 1,2,3,\ldots$

Substituting Eq 4.2-48 into Eq 4.2-47

$$(D_{\Delta x}-a)^m\{e_{\Delta x}(a,x)\sum_{p=1}^{m} k_{m-p+1}[x]_{\Delta x}^{p-1}\} = 0 \qquad (4.2\text{-}49)$$

Note that Eq 4.2-49 is valid for m = 1 where all roots are unique

where

$k_1,k_2,k_3\ldots k_m$ = real or complex value constants

m = the multiplicity of the repeated roots of the homogeneous equation characteristic polynomial, h(r)

$m = 1,2,3,4\ldots$

r = a is the root of multiplicity m

$[x]_{\Delta x}^{0} = 1$

$[x]_{\Delta x}^{h} = x(x-\Delta x)(x-2\Delta x)(x-3\Delta x)\ldots(x-[h-1]\Delta x)\ , \quad h = 1,2,3,\ldots$

Note that the function, $(\sum_{p=1}^{m} k_{m-p+1}[x]_{\Delta x}^{p-1})e_{\Delta x}(a,x)$, which satisfies a homogeneous differential difference equation with roots of multiplicity m, has m independent constants. These constants are necessary to provided the required number of initial conditions to obtain a homogeneous equation general solution. The required number of general solution initial conditions is equal to the order of the homogeneous equation, n.

Thus, in general

**The complementary solution of a differential difference equation with a characteristic polynomial, h(x), with m roots of equal value, a, is as follows:**

**From the related homogeneous equation:**

$$\mathbf{D_{\Delta x}^{n} f_C(x) + a_{n-1}D_{\Delta x}^{n-1} f_C(x) + a_{n-2}D_{\Delta x}^{n-2} f_C(x) + \ldots + a_1D_{\Delta x} f_C(x) + a_0 f_C(x) = 0}$$

**or**

$$\mathbf{(D_{\Delta x}-r_n)(D_{\Delta x}-r_{n-1})(D_{\Delta x}-r_{n-2}) \ldots (D_{\Delta x}-r_{m+1})(D_{\Delta x}-a)^m f_C(x) = 0}$$

**The differential difference equation complementary solution, $f_C(x)$, is:**

$$f_C(x) = k_n e(r_n,x) + k_{n-1}e(r_{n-1},x) + k_{n-2}e(r_{n-2},x) + \ldots + k_{m+1}e(r_{m+1},x) + \left(\sum_{p=1}^{m} k_{m+1-p}[x]_{\Delta x}^{p-1}\right)e_{\Delta x}(a,x)$$

**or**

$$f_C(x) = k_n(1+r_n\Delta x)^{\frac{x}{\Delta x}} + k_{n-1}(1+r_{n-1}\Delta x)^{\frac{x}{\Delta x}} + k_{n-2}(1+r_{n-2}\Delta x)^{\frac{x}{\Delta x}} + \ldots + k_{m+1}(1+r_{m+1}\Delta x)^{\frac{x}{\Delta x}} +$$

$$\left(\sum_{p=1}^{m} k_{m+1-p}[x]_{\Delta x}^{p-1}(1+a\Delta x)^{\frac{x}{\Delta x}}\right.$$

**where**

**n = the order of the homogeneous equation**

$$h(r) = r^n + a_{n-1}r^{n-1} + a_{n-2}r^{n-2} + \ldots + a_2r^2 + a_1r + a_0 = (r-r_n)(r-r_{n-1})(r-r_{n-2})\ldots(r-r_4)(r-a)(r-a)(r-a) = 0$$

**$r_p$, p=1,2,3,…,n are the roots of the characteristic polynomial, h(r)**

**m of the roots have the same value, $r_1$ thru $r_m$ have the value, a**

**The roots, $r_1$ thru $r_n$ may be real or complex conjugate pair values.**

$$[x]_{\Delta x}^{0} = 1$$

$$[x]_{\Delta x}^{h} = x(x-\Delta x)(x-2\Delta x)(x-3\Delta x)\ldots(x-[h-1]\Delta x)\ , \quad h = 1,2,3,\ldots$$

The results of the previous derivations are presented in a more useful form in Table 4.2-1 below.

**Table 4.2-1 The General Solution to Homogeneous Differential Difference Equations**

For the homogeneous differential difference equation:

$$D_{\Delta x}^{\ n}f_C(x) + a_{n-1}D_{\Delta x}^{\ n-1}f_C(x) + a_{n-2}D_{\Delta x}^{\ n-2}f_C(x) + \ldots + a_1D_{\Delta x}f_C(x) + a_0f_C(x) = 0$$

with the characteristic polynomial:

$$h(r) = r^n + a_{n-1}r^{n-1} + a_{n-2}r^{n-2} + \ldots + a_2r^2 + a_1r + a_0 = (r-r_n)(r-r_{n-1})(r-r_{n-2})\ldots(r-r_3)(r-r_2)(r-r_1) = 0$$

$r = r_1, r_2, r_3, \ldots, r_{n-1}, r_n$ , the n roots of h(r)

$$f_C(x) = \left[\sum_{m=1}^{n} k_m w_m(x)\right]e_{\Delta x}(r_m,x)$$

or

$$f_C(x) = \left[\sum_{m=1}^{n} k_m w_m(x)\right](1+r_m\Delta x)^{\frac{x}{\Delta x}}$$

$$w_m(x) = \begin{cases} [x]_{\Delta x}^{0} = 1 & \text{For non-multiple roots or the 1}^{\text{st}}\text{ of a multiple root} \\ [x]_{\Delta x}^{1} & \text{For the 2}^{\text{nd}}\text{ of a multiple root} \\ [x]_{\Delta x}^{2} & \text{For the 3}^{\text{rd}}\text{ of a multiple root} \\ [x]_{\Delta x}^{3} & \text{For the 4}^{\text{th}}\text{ of a multiple root} \\ \cdot & \\ \cdot & \\ \cdot & \end{cases}$$

where

$x = x_0 + p\Delta x$ , $p$ = integers

$x_0$ = initial value of x

$\Delta x$ = x increment

$a_{n-1},a_{n-1},\ldots,a_0$ = real constants

n = order of the homogeneous equation

$r_1,r_2,\ldots,r_n$ = real or complex roots of the characteristic polynomial, h(r)

$k_m$ = the n constants of $f_C(x)$ which allow the specification of n initial conditions

$m = 1,2,3,\ldots,n$

$f_C(x)$ = general solution to the homogeneous equation

$D_{\Delta x}q(x) = \dfrac{q(x+\Delta x) - q(x)}{\Delta x}$ , discrete derivative

$[x]_{\Delta x}^{0} = 1$

$[x]_{\Delta x}^{h} = \prod_{u=1}^{h}(x-[u-1]\Delta x)$ , $h = 1,2,3,\ldots$

$e_{\Delta x}(r,x) = (1+r\Delta x)^{\frac{x}{\Delta x}}$

**Derivation of the homogeneous differential difference equation real value functions associated with complex conjugate roots**

Two of the roots of the homogeneous differential difference equation characteristic polynomial may be complex conjugate pairs, $a \pm jb$. Then the functions $[1+(a+jb)\Delta x]^{\frac{x}{\Delta x}}$ and $[1+(a-jb)\Delta x]^{\frac{x}{\Delta x}}$ each satisfy the homogeneous differential difference equation. These two functions obviously yield complex results.

However, these functions can be linearly combined and simplified to yield real value functions which also satisfy the homogeneous equation.

1) Consider the case where two roots of the characteristic polynomial are the complex conjugate pairs, 0+jb, 0-jb.

For the roots 0±jb, the functions, $[1+jb\Delta x]^{\frac{x}{\Delta x}}$ and $[1-jb\Delta x]^{\frac{x}{\Delta x}}$ each satisfy the differential difference homogeneous equation. These two functions yield complex results. However, they can be linearly combined to form two other functions which are real value functions.

Combining the above two functions to form two other functions

$$\sin_{\Delta x}(b,x) = \frac{[1+jb\Delta x]^{\frac{x}{\Delta x}} - [1-jb\Delta x]^{\frac{x}{\Delta x}}}{2j} \quad , \qquad \text{roots } jb, -jb \tag{4.2-51}$$

Note that since $\frac{x}{\Delta x}$ = positive integer where x = 0, Δx, 2Δx, … expanding both terms yields a result with all real expansion terms canceling. Dividing by j produces a real result.

and

$$\cos_{\Delta x}(b,x) = \frac{[1+jb\Delta x]^{\frac{x}{\Delta x}} + [1-jb\Delta x]^{\frac{x}{\Delta x}}}{2} \quad , \qquad \text{roots } jb, -jb \tag{4.2-52}$$

Note that since $\frac{x}{\Delta x}$ = positive integer where x = 0, Δx, 2Δx, … expanding both terms yields a result with all imaginary expansion terms canceling leaving a real result.

Since the functions, $[1+jb\Delta x]^{\frac{x}{\Delta x}}$ and $[1-jb\Delta x]^{\frac{x}{\Delta x}}$, both satisfy the homogeneous equation, their linear combinations do also. Then the real value functions, $\sin_{\Delta x}(b,x)$ and $\cos_{\Delta x}(b,x)$ can be used instead of the functions, $[1+jb\Delta x]^{\frac{x}{\Delta x}}$ and $[1-jb\Delta x]^{\frac{x}{\Delta x}}$.

Note that the functions, $[1+jb\Delta x]^{\frac{x}{\Delta x}}$ and $[1-jb\Delta x]^{\frac{x}{\Delta x}}$, can be expressed as $e_{\Delta x}(jb,x)$ and $e_{\Delta x}(-jb,x)$ respectively.

2) Consider the case where two roots of the characteristic polynomial are the complex conjugate pairs, a+jb, a-jb where 1+aΔx = 0.

For the roots a ± jb where 1+aΔx = 0, the functions, $[+jb\Delta x]^{\frac{x}{\Delta x}}$ and $[-jb\Delta x]^{\frac{x}{\Delta x}}$ each satisfy the differential difference homogeneous equation. These two functions yield complex results. However, they can be linearly combined to form two other functions which are real value funcions. Since each function satisfies the homogeneous equation, the linear combinations of these functions satisfy the homogeneous equation also.

Finding another form for the function, $[+jb\Delta x]^{\frac{x}{\Delta x}}$

$$[+jb\Delta x]^{\frac{x}{\Delta x}} = [b\Delta x]^{\frac{x}{\Delta x}}[+j]^{\frac{x}{\Delta x}} = [b\Delta x]^{\frac{x}{\Delta x}}[e^{\ln(j)}]^{\frac{x}{\Delta x}} = [b\Delta x]^{\frac{x}{\Delta x}}[e^{\frac{\pi}{2}j}]^{\frac{x}{\Delta x}} = [b\Delta x]^{\frac{x}{\Delta x}}\, e^{j\frac{\pi x}{2\Delta x}}$$

$$[+jb\Delta x]^{\frac{x}{\Delta x}} = [b\Delta x]^{\frac{x}{\Delta x}}\, e^{j\frac{\pi x}{2\Delta x}} \qquad (4.2\text{-}53)$$

Finding another form for the function, $[-jb\Delta x]^{\frac{x}{\Delta x}}$

$$[-jb\Delta x]^{\frac{x}{\Delta x}} = [b\Delta x]^{\frac{x}{\Delta x}}[-j]^{\frac{x}{\Delta x}} = [b\Delta x]^{\frac{x}{\Delta x}}[e^{\ln(-j)}]^{\frac{x}{\Delta x}} = [b\Delta x]^{\frac{x}{\Delta x}}[e^{-\frac{\pi}{2}j}]^{\frac{x}{\Delta x}} = [b\Delta x]^{\frac{x}{\Delta x}}\, e^{-j\frac{\pi x}{2\Delta x}}$$

$$[-jb\Delta x]^{\frac{x}{\Delta x}} = [b\Delta x]^{\frac{x}{\Delta x}}\, e^{-j\frac{\pi x}{2\Delta x}} \qquad (4.2\text{-}54)$$

Subtracting Eq 4.2-54 from Eq 4.2-53

$$\frac{[+jb\Delta x]^{\frac{x}{\Delta x}} - [-jb\Delta x]^{\frac{x}{\Delta x}}}{2j} = [b\Delta x]^{\frac{x}{\Delta x}}\left[\frac{e^{j\frac{\pi x}{2\Delta x}} - e^{-j\frac{\pi x}{2\Delta x}}}{2j}\right] = [b\Delta x]^{\frac{x}{\Delta x}} \sin\frac{\pi x}{2\Delta x}$$

$$\frac{[+jb\Delta x]^{\frac{x}{\Delta x}} - [-jb\Delta x]^{\frac{x}{\Delta x}}}{2j} = [b\Delta x]^{\frac{x}{\Delta x}} \sin\frac{\pi x}{2\Delta x}, \quad \text{roots } a\pm jb \quad \text{where } 1+a\Delta x = 0 \qquad (4.2\text{-}55)$$

Adding Eq 4.2-54 to Eq 4.2-53

$$\frac{[+jb\Delta x]^{\frac{x}{\Delta x}} + [-jb\Delta x]^{\frac{x}{\Delta x}}}{2} = [b\Delta x]^{\frac{x}{\Delta x}}\left[\frac{e^{j\frac{\pi x}{2\Delta x}} + e^{-j\frac{\pi x}{2\Delta x}}}{2}\right] = [b\Delta x]^{\frac{x}{\Delta x}} \cos\frac{\pi x}{2\Delta x}$$

$$\frac{[+jb\Delta x]^{\frac{x}{\Delta x}} + [-jb\Delta x]^{\frac{x}{\Delta x}}}{2} = [b\Delta x]^{\frac{x}{\Delta x}} \cos\frac{\pi x}{2\Delta x}, \quad \text{roots } a\pm jb \quad \text{where } 1+a\Delta x = 0 \qquad (4.2\text{-}56)$$

Since the functions, $[+jb\Delta x]^{\frac{x}{\Delta x}}$ and $[-jb\Delta x]^{\frac{x}{\Delta x}}$, both satisfy the homogeneous equation, their linear combinations do also. Then the real value functions, $[b\Delta x]^{\frac{x}{\Delta x}} \sin\frac{\pi x}{2\Delta x}$ and $[b\Delta x]^{\frac{x}{\Delta x}} \cos\frac{\pi x}{2\Delta x}$ can be used instead of the functions, $[+jb\Delta x]^{\frac{x}{\Delta x}}$ and $[-jb\Delta x]^{\frac{x}{\Delta x}}$.

Note that the functions, $[+jb\Delta x]^{\frac{x}{\Delta x}}$ and $[-jb\Delta x]^{\frac{x}{\Delta x}}$, can be expressed as $[b\Delta x]^{\frac{x}{\Delta x}}\, e^{j\frac{\pi x}{2\Delta x}}$ and $[b\Delta x]^{\frac{x}{\Delta x}}\, e^{-j\frac{\pi x}{2\Delta x}}$ respectively.

3) Consider the case where two roots of the characteristic polynomial are the complex conjugate pairs, a+jb, a-jb where $1+a\Delta x \neq 0$.

For the roots $a \pm jb$ where $1+a\Delta x \neq 0$, the functions, $[1+(a+jb)\Delta x]^{\frac{x}{\Delta x}}$ and $[1+(a-jb)\Delta x]^{\frac{x}{\Delta x}}$ each satisfy the differential difference homogeneous equation. These two functions yield complex results. However, they can be linearly combined to form two other functions which are real value functions. Since both functions satisfy the homogeneous equation, their linear combinations satisfy the homogeneous equation also.

Combining the above two functions to form two other functions

$$\frac{[1+(a+jb)\Delta x]^{\frac{x}{\Delta x}} - [1+(a-jb)\Delta x]^{\frac{x}{\Delta x}}}{2j} = [1+a\Delta x]^{\frac{x}{\Delta x}} \left\{ \frac{\left[\frac{1+(a+jb)\Delta x}{1+a\Delta x}\right]^{\frac{x}{\Delta x}} - \left[\frac{1+(a-jb)\Delta x}{1+a\Delta x}\right]^{\frac{x}{\Delta x}}}{2j} \right\} \qquad (4.2\text{-}57)$$

$$\frac{1+(a+jb)\Delta x}{1+a\Delta x} = 1 + \frac{jb\Delta x}{1+a\Delta x} \qquad (4.2\text{-}58)$$

$$\frac{1+(a-jb)\Delta x}{1+a\Delta x} = 1 - \frac{jb\Delta x}{1+a\Delta x} \qquad (4.2\text{-}59)$$

Substituting Eq 4.2-58 and Eq 4.2-59 into Eq 4.2-57

$$\frac{[1+(a+jb)\Delta x]^{\frac{x}{\Delta x}} - [1+(a-jb)\Delta x]^{\frac{x}{\Delta x}}}{2j} = [1+a\Delta x]^{\frac{x}{\Delta x}} \left\{ \frac{\left[1 + \frac{jb\Delta x}{1+a\Delta x}\right]^{\frac{x}{\Delta x}} - \left[1 - \frac{jb\Delta x}{1+a\Delta x}\right]^{\frac{x}{\Delta x}}}{2j} \right\} \qquad (4.2\text{-}60)$$

$$\frac{[1+(a+jb)\Delta x]^{\frac{x}{\Delta x}} - [1+(a-jb)\Delta x]^{\frac{x}{\Delta x}}}{2j} = [1+a\Delta x]^{\frac{x}{\Delta x}} \sin_{\Delta x}\left(\frac{b}{1+a\Delta x}, x\right), \quad \text{roots } a \pm jb \text{ where } 1+a\Delta x \neq 0$$

(4.2-61)

Note that the function, $[1+a\Delta x]^{\frac{x}{\Delta x}}$, and the function, $\sin_{\Delta x}\left(\frac{b}{1+a\Delta x}, x\right)$, are both real value functions so that their product, $[1+a\Delta x]^{\frac{x}{\Delta x}} \sin_{\Delta x}\left(\frac{b}{1+a\Delta x}, x\right)$, is a real value function.

Also,

$$\frac{[1+(a+jb)\Delta x]^{\frac{x}{\Delta x}} + [1+(a-jb)\Delta x]^{\frac{x}{\Delta x}}}{2} = [1+a\Delta x]^{\frac{x}{\Delta x}} \left\{ \frac{\left[\frac{1+(a+jb)\Delta x}{1+a\Delta x}\right]^{\frac{x}{\Delta x}} + \left[\frac{1+(a-jb)\Delta x}{1+a\Delta x}\right]^{\frac{x}{\Delta x}}}{2} \right\} \qquad (4.2\text{-}62)$$

Substituting Eq 4.2-58 and Eq 4.2-59 into Eq 4.2-62

$$\frac{[1+(a+jb)\Delta x]^{\frac{x}{\Delta x}} + [1+(a-jb)\Delta x]^{\frac{x}{\Delta x}}}{2} = [1+a\Delta x]^{\frac{x}{\Delta x}} \left\{ \frac{\left[1 + \frac{jb\Delta x}{1+a\Delta x}\right]^{\frac{x}{\Delta x}} + \left[1 - \frac{jb\Delta x}{1+a\Delta x}\right]^{\frac{x}{\Delta x}}}{2} \right\} \qquad (4.2\text{-}63)$$

$$\frac{[1+(a+jb)\Delta x]^{\frac{x}{\Delta x}} + [1+(a-jb)\Delta x]^{\frac{x}{\Delta x}}}{2} = [1+a\Delta x]^{\frac{x}{\Delta x}} \cos_{\Delta x}\left( \frac{b}{1+a\Delta x} , x \right), \quad \text{roots } a \pm jb \text{ where } 1+a\Delta x \neq 0 \tag{4.2-64}$$

Note that the function, $[1+a\Delta x]^{\frac{x}{\Delta x}}$, and the function, $\cos_{\Delta x}( \frac{b}{1+a\Delta x} , x )$, are both real value functions so that their product, $[1+a\Delta x]^{\frac{x}{\Delta x}} \cos_{\Delta x}( \frac{b}{1+a\Delta x} , x )$, is a real value function.

The results of the previous function derivations are presented in Table 4.2-2 below.

**Table 4.2-2 Homogeneous Differential Difference Equation Real Value Functions Associated with Complex Conjugate Roots**

Real value functions which satisfy the differential difference homogeneous equation,

$D_{\Delta x}{}^{n}f_C(x) + a_{n-1}D_{\Delta x}{}^{n-1}f_C(x) + a_{n-2}D_{\Delta x}{}^{n-2}f_C(x) + \ldots + a_1D_{\Delta x}f_C(x) + a_0f_C(x) = 0,$

where the characteristic polynomial,

$h(r) = r^n + a_{n-1}r^{n-1} + a_{n-2}\, r^{n-2} + \ldots + a_2r^2 + a_1r + a_{0},$

has a complex conjugate root pair of the form, $a \pm jb$

<u>For a complex conjugate root pair, $0 \pm jb$</u>

$$\sin_{\Delta x}(b,x) = \frac{[1+jb\Delta x]^{\frac{x}{\Delta x}} - [1-jb\Delta x]^{\frac{x}{\Delta x}}}{2j}$$

$$\cos_{\Delta x}(b,x) = \frac{[1+jb\Delta x]^{\frac{x}{\Delta x}} + [1-jb\Delta x]^{\frac{x}{\Delta x}}}{2}$$

<u>For a complex conjugate root pair,</u> $-\frac{1}{\Delta x} \pm jb$

$$[b\Delta x]^{\frac{x}{\Delta x}} \sin\frac{\pi x}{2\Delta x} = \frac{[+jb\Delta x]^{\frac{x}{\Delta x}} - [-jb\Delta x]^{\frac{x}{\Delta x}}}{2j}$$

$$[b\Delta x]^{\frac{x}{\Delta x}} \cos\frac{\pi x}{2\Delta x} = \frac{[+jb\Delta x]^{\frac{x}{\Delta x}} + [-jb\Delta x]^{\frac{x}{\Delta x}}}{2}$$

For a complex conjugate root pair, $a \pm jb$, where $a \neq -\frac{1}{\Delta x}$

$$[1+a\Delta x]^{\frac{x}{\Delta x}} \sin_{\Delta x}\left(\frac{b}{1+a\Delta x}, x\right) = \frac{[1+(a+jb)\Delta x]^{\frac{x}{\Delta x}} - [1+(a-jb)\Delta x]^{\frac{x}{\Delta x}}}{2j}$$

$$[1+a\Delta x]^{\frac{x}{\Delta x}} \cos_{\Delta x}\left(\frac{b}{1+a\Delta x}, x\right) = \frac{[1+(a+jb)\Delta x]^{\frac{x}{\Delta x}} + [1+(a-jb)\Delta x]^{\frac{x}{\Delta x}}}{2}$$

Introducing the Functions of Table 4.2-2 into the Homogeneous Equation General Solution presented in Table 4.2-1

Some of the functions representing the general solution of a homogeneous differential difference equation in Table 4.2-1, while valid, are not necessarily the most convenient to use. In the case of single real roots of the characteristic polynomial, there is no problem. The functions yield real values. However, in the case of complex conjugate root pairs, $a\pm jb$, the function representing each of the two roots yields complex values. It is true that when these functions are combined, the imaginary components of both cancel. However, mathematically manipulating these function's real and imaginary value components can be time consuming. The real value functions derived and then presented in Table 4.2-2 can be very helpful in resolving this difficulty.

Consider the function, $F(x)$, which is the sum of the two general solution functions associated with one of the m $a\pm jb$ complex conjugate root pairs.

$$F(x) = C_p[x]_{\Delta x}^{p-1}[1+(a+jb)\Delta x]^{\frac{x}{\Delta x}} + D_p[x]_{\Delta x}^{p-1}[1+(a-jb)\Delta x]^{\frac{x}{\Delta x}} \qquad (4.2\text{-}65)$$

where

This $a\pm jb$ complex conjugate root pair is the pth root pair of a root pair of multiplicity m
$m = 1,2,3,\ldots$

$p = 1$ for a single complex conjugate root pair
$p = 1$ for the 1st conjugate root pair of a root pair of multiplicity m
$p = 2$ for the 2nd conjugate root pair of a root pair of multiplicity m
$p = 3$ for the 3rd conjugate root pair of a root pair of multiplicity m
…
$a,b$ = real constants
$C_p, D_p$ = complex constants
$\Delta x$ = x increment
$x = x_0 + p\Delta x$ , $p$ = integers
$x_0$ = initial value of x

From Eq 4.2-65

$$F(x) = [x]_{\Delta x}^{p-1}\{C_p[1+(a+jb)\Delta x]^{\frac{x}{\Delta x}} + D_p[1+(a-jb)\Delta x]^{\frac{x}{\Delta x}}\} \tag{4.2-66}$$

Let

$$C_p = \frac{A_p}{2} - j\frac{B_p}{2} \tag{4.2-67}$$

$$D_p = \frac{A_p}{2} + j\frac{B_p}{2} \tag{4.2-68}$$

where

$A_p, B_p$ = real constants

Substituting Eq 4.2-67 and Eq 4.2-68 into Eq 4.2-66

$$F(x) = [x]_{\Delta x}^{p-1}\{(\frac{A_p}{2} - j\frac{B_p}{2})[1+(a+jb)\Delta x]^{\frac{x}{\Delta x}} + (\frac{A_p}{2} + j\frac{B_p}{2})[1+(a-jb)\Delta x]^{\frac{x}{\Delta x}}\} \tag{4.2-69}$$

Rearranging Eq 4.2-66

$$F(x) = [x]_{\Delta x}^{p-1}\{A_p(\frac{[1+(a+jb)\Delta x]^{\frac{x}{\Delta x}} + [1+(a-jb)\Delta x]^{\frac{x}{\Delta x}}}{2}) + B_p(\frac{[1+(a+jb)\Delta x]^{\frac{x}{\Delta x}} - [1+(a-jb)\Delta x]^{\frac{x}{\Delta x}}}{2j})\} \tag{4.2-70}$$

There are three possibilities for the complex conjugate root pair, $a \pm jb$, $a = 0$, $a = -\frac{1}{\Delta x}$ and $a \neq -\frac{1}{\Delta x}$. The function, F(x), for each of the above three possibilities is specified below.

For $a \neq -\frac{1}{\Delta x}$.

---

From Table 4.2-2

$$[1+a\Delta x]^{\frac{x}{\Delta x}} \sin(\frac{b}{1+a\Delta x}, x) = \frac{[1+(a+jb)\Delta x]^{\frac{x}{\Delta x}} - [1+(a-jb)\Delta x]^{\frac{x}{\Delta x}}}{2j}, \quad \text{for } a \neq -\frac{1}{\Delta x}. \tag{4.2-71}$$

$$[1+a\Delta x]^{\frac{x}{\Delta x}} \cos(\frac{b}{1+a\Delta x}, x) = \frac{[1+(a+jb)\Delta x]^{\frac{x}{\Delta x}} + [1+(a-jb)\Delta x]^{\frac{x}{\Delta x}}}{2}, \quad \text{for } a \neq -\frac{1}{\Delta x} \tag{4.2-72}$$

Substituting Eq 4.2-71 and Eq 4.2-72 into Eq 4.2-70

For a complex conjugate root pair, $a \pm jb$, where $a \neq -\frac{1}{\Delta x}$

$$F_1(x) = [x]_{\Delta x}^{p-1}\{A_p[1+a\Delta x]^{\frac{x}{\Delta x}} \cos_{\Delta x}(\frac{b}{1+a\Delta x}, x) + B_p[1+a\Delta x]^{\frac{x}{\Delta x}} \sin_{\Delta x}(\frac{b}{1+a\Delta x}, x)\}$$

or since $e_{\Delta x}(a,x) = [1+a\Delta x]^{\frac{x}{\Delta x}}$

$$F_1(x) = [x]_{\Delta x}^{p-1}\{A_p e_{\Delta x}(a,x) \cos_{\Delta x}(\frac{b}{1+a\Delta x}, x) + B_p e_{\Delta x}(a,x) \sin_{\Delta x}(\frac{b}{1+a\Delta x}, x)\} \qquad (4.2\text{-}73)$$

where

$A_p = C_p + D_p$

$B_p = j(C_p - D_p)$

For a = 0.

From Eq 4.2-73 where a = 0

For a complex conjugate root pair, $0 \pm jb$

$$F_2(x) = [x]_{\Delta x}^{p-1}\{A_p \cos_{\Delta x}(b,x) + B_p \sin_{\Delta x}(b,x)\} \quad , \qquad (4.2\text{-}74)$$

where

$A_p = C_p + D_p$

$B_p = j(C_p - D_p)$

For $a = -\frac{1}{\Delta x}$.

From Eq 4.2-70 where $a = -\frac{1}{\Delta x}$.

$$F(x) = [x]_{\Delta x}^{p-1}\{A_p(\frac{[+jb\Delta x]^{\frac{x}{\Delta x}} + [-jb\Delta x]^{\frac{x}{\Delta x}}}{2}) + B_p(\frac{[+jb\Delta x]^{\frac{x}{\Delta x}} - [-jb\Delta x]^{\frac{x}{\Delta x}}}{2j})\} \qquad (4.2\text{-}75)$$

From Table 4.2-2

$$[b\Delta x]^{\frac{x}{\Delta x}} \sin\frac{\pi x}{2\Delta x} = \frac{[+jb\Delta x]^{\frac{x}{\Delta x}} - [-jb\Delta x]^{\frac{x}{\Delta x}}}{2j} \quad , \qquad \text{for } a = -\frac{1}{\Delta x} \qquad (4.2\text{-}76)$$

$$[b\Delta x]^{\frac{x}{\Delta x}} \cos\frac{\pi x}{2\Delta x} = \frac{[+jb\Delta x]^{\frac{x}{\Delta x}} + [-jb\Delta x]^{\frac{x}{\Delta x}}}{2} \quad , \qquad \text{for } a = -\frac{1}{\Delta x} \tag{4.2-77}$$

Substituting Eq 4.2-76 and Eq 4.2-77 into Eq 4.2-75

For a complex conjugate root pair, $-\frac{1}{\Delta x} \pm jb$

$$F_3(x) = [x]^{p-1}_{\Delta x}\{A_p[b\Delta x]^{\frac{x}{\Delta x}} \cos\frac{\pi x}{2\Delta x} + B_p[b\Delta x]^{\frac{x}{\Delta x}} \sin\frac{\pi x}{2\Delta x}\} \tag{4.2-78}$$

Note - The functions, $F_1(x)$, $F_2(x)$, and $F_3(x)$ yield real values for real x and real $A_p$ and $B_p$. For $A_p$ and $B_p$ to be real, $C_p$ and $D_p$ must be complex conjugate root pairs.

## The Real Function General Solutions to Homogeneous Differential Difference Equations

From the previous derivations, the real value function general solutions to homogeneous differential difference equations have been obtained. These general solutions to homogeneous differential difference equations are presented in Table 4.2-3 on the following page.

**Table 4.2-3 The real value function general solution of homogeneous differential difference equations**

For the differential difference homogeneous equation:

$$D_{\Delta x}{}^{n}f_C(x) + a_{n-1}D_{\Delta x}{}^{n-1}f_C(x) + a_{n-2}D_{\Delta x}{}^{n-2}f_C(x) + \ldots + a_1D_{\Delta x}f_C(x) + a_0f_C(x) = 0$$

where

$n$ = order of the homogeneous equation
$a_{n-1}, a_{n-1}, \ldots, a_0$ = real constants
$f_C(x)$ = general solution to the homogeneous differential difference equation
$D_{\Delta x}{}^{n}f(x)$ = nth discrete derivative of the function, $f(x)$
$x = x_0 + p\Delta x$ , $p$ = integers
$x_0$ = initial value of x
$\Delta x$ = x increment

with the characteristic polynomial:

$$h(r) = r^n + a_{n-1}r^{n-1} + a_{n-2}\,r^{n-2} + \ldots + a_2r^2 + a_1r + a_0 = (r-r_n)(r-r_{n-1})(r-r_{n-2})\ldots(r-r_3)(r-r_2)(r-r_1) = 0$$

with n roots

$r = r_1, r_2, r_3, \ldots, r_{n-1}, r_n$

group the roots and catagorize them in the following way:

v = Each unique root or complex conjugate root pair value
m = multiplicity of each unique root or complex conjugate root pair value, m = 1,2,3,…
g = number of groups of unique root/root pair values ($g \le n$)

Note – A root of multiplicity 1 is a single root, a root of multiplicity 2 is a double root, etc.

W(v,m,x) = real value general solution function associated with the root(s), v, of multiplicity, m

$f_C(x) = \sum_{s=1}^{g} W_s(v_s,m_s,x)$ , the general solution to the homogeneous differential difference equation

In the following table

$a,b,A_pB_p$ = real constants

$[x]_{\Delta x}^{0} = 1$

$[x]_{\Delta x}^{q} = \sum_{u=1}^{q}(x-[u-1]\Delta x)$ , $q = 1,2,3,\ldots$

$x = x_0 + c\Delta x$ , c = integers

$x_0$ = initial value of x

$\Delta x$ = x increment

| # | v Root/Root Pair Value | W(v,m,x) Function | W(v,m,x) Calculation Function | Comments |
|---|---|---|---|---|
| **1** | a | $[\sum_{p=1}^{m} A_p[x]_{\Delta x}^{p-1}]\, e_{\Delta x}(a,x)$ | $e_{\Delta x}(a,x) = [1+a\Delta x]^{\frac{x}{\Delta x}}$ | a is a single real root |
| **2** | $0 \pm jb$ | $[\sum_{p=1}^{m} A_p[x]_{\Delta x}^{p-1}]\, \sin_{\Delta x}(b,x) + [\sum_{p=1}^{m} B_p[x]_{\Delta x}^{p-1}]\, \cos_{\Delta x}(b,x)$ | $\sin_{\Delta x}(b,x) = \frac{[1+jb\Delta x]^{\frac{x}{\Delta x}} - [1-jb\Delta x]^{\frac{x}{\Delta x}}}{2j}$ <br> $\cos_{\Delta x}(b,x) = \frac{[1+jb\Delta x]^{\frac{x}{\Delta x}} + [1-jb\Delta x]^{\frac{x}{\Delta x}}}{2}$ | $a \pm jb$ is a complex conjugate root pair <br> $a = 0$ |
| **3** | $-\frac{1}{\Delta x} \pm jb$ | $[\sum_{p=1}^{m} A_p[x]_{\Delta x}^{p-1}][b\Delta x]^{\frac{x}{\Delta x}}\sin\frac{\pi x}{2\Delta x} + [\sum_{p=1}^{m} B_p[x]_{\Delta x}^{p-1}][b\Delta x]^{\frac{x}{\Delta x}}\cos\frac{\pi x}{2\Delta x}$ | $[b\Delta x]^{\frac{x}{\Delta x}}\sin\frac{\pi x}{2\Delta x} = \frac{[+jb\Delta x]^{\frac{x}{\Delta x}} - [-jb\Delta x]^{\frac{x}{\Delta x}}}{2j}$ <br> $[b\Delta x]^{\frac{x}{\Delta x}}\cos\frac{\pi x}{2\Delta x} = \frac{[+jb\Delta x]^{\frac{x}{\Delta x}} + [-jb\Delta x]^{\frac{x}{\Delta x}}}{2}$ | $a \pm jb$ is a complex conjugate root pair <br> $a = -\frac{1}{\Delta x}$ |
| **4** | $a \pm jb$<br>$a \neq -\frac{1}{\Delta x}$ | $[\sum_{p=1}^{m} A_p[x]_{\Delta x}^{p-1}]\, e_{\Delta x}(a,x)\, \sin_{\Delta x}(\frac{b}{1+a\Delta x}, x) +$<br>$[\sum_{p=1}^{m} B_p[x]_{\Delta x}^{p-1}]\, e_{\Delta x}(a,x)\, \cos_{\Delta x}(\frac{b}{1+a\Delta x}, x)$ | $e_{\Delta x}(a,x)\, \sin_{\Delta x}(\frac{b}{1+a\Delta x}, x) = \frac{[1+(a+jb)\Delta x]^{\frac{x}{\Delta x}} - [1+(a-jb)\Delta x]^{\frac{x}{\Delta x}}}{2j}$<br>$e_{\Delta x}(a,x)\, \cos_{\Delta x}(\frac{b}{1+a\Delta x}, x) = \frac{[1+(a+jb)\Delta x]^{\frac{x}{\Delta x}} + [1+(a-jb)\Delta x]^{\frac{x}{\Delta x}}}{2}$ | $a \pm jb$ is a complex conjugate root pair<br>$a \neq -\frac{1}{\Delta x}$ |

The general solution, f(x) to a differential difference equation is composed of two components, the complementary solution, $f_C(x)$, and the particular solution, $f_P(x)$. Table 4.2-3, which is presented on the previous pages, provides the methodology to obtain the complementary solution, $f_C(x)$. It is now necessary that a methodology be derived to obtain the particular solution, $f_P(x)$. The derivation of this methodology follows.

**The Differential Difference Equation Particular Solution**

For convenient reference, the differential difference equations associated with a differential equation general solution are rewritten below. (From Eq 4.2-16, Eq 4.2-17, and Eq 4.2-18)

1) The Differential Difference Equation General Solution
$f(x) = f_C(x) + f_P(x)$

$$D_{\Delta x}{}^{n}f(x) + a_{n-1}D_{\Delta x}{}^{n-1}f(x) + a_{n-2}D_{\Delta x}{}^{n-2}f(x) + \ldots + a_1D_{\Delta x}f(x) + a_0f(x) = Q(x) \qquad (4.2\text{-}79)$$

2) The Related Homogeneous Equation Complementary Solution
$f_C(x)$

$$D_{\Delta x}{}^{n}f_C(x) + a_{n-1}D_{\Delta x}{}^{n-1}f_C(x) + a_{n-2}D_{\Delta x}{}^{n-2}f_C(x) + \ldots + a_1D_{\Delta x}f_C(x) + a_0f_C(x) = 0 \qquad (4.2\text{-}80)$$

3) The Differential Difference Equation Particular Solution
$f_P(x)$

$$D_{\Delta x}{}^{n}f_P(x) + a_{n-1}D_{\Delta x}{}^{n-1}f_P(x) + a_{n-2}D_{\Delta x}{}^{n-2}f_P(x) + \ldots + a_1D_{\Delta x}f_P(x) + a_0f_P(x) = Q(x) \qquad (4.2\text{-}81)$$

Note – f(x), $f_C(x)$, and $f_P(x)$ are the general solutions of their respective equations.

$f(x) = f_C(x) + f_P(x)$ , the general solution to differential difference equation, Eq 4.2-79 (4.2-82)

Substituting Eq 4.2-82 into Eq 4.2-79

$$D_{\Delta x}{}^{n}[f_C(x)+f_P(x)] + a_{n-1}D_{\Delta x}{}^{n-1}[f_C(x)+f_P(x)] + a_{n-2}D_{\Delta x}{}^{n-2}[f_C(x)+f_P(x)] + \ldots + \qquad (4.2\text{-}83)$$
$$a_1D_{\Delta x}[f_C(x)+f_P(x)] + a_0[f_C(x)+f_P(x)] = Q(x)$$

Expanding Eq 4.2-83

$$[D_{\Delta x}{}^{n} f_C(x) + a_{n-1}D_{\Delta x}{}^{n-1}f_C(x) + a_{n-2}D_{\Delta x}{}^{n-2}f_C(x) + \ldots + a_1D_{\Delta x}f_C(x) + a_0f_C(x)] + \qquad (4.2\text{-}84)$$
$$[D_{\Delta x}{}^{n} f_P(x) + a_{n-1}D_{\Delta x}{}^{n-1}f_P(x) + a_{n-2}D_{\Delta x}{}^{n-2}f_P(x) + \ldots + a_1D_{\Delta x}f_P(x) + a_0f_P(x)] = Q(x)$$

But the first term of Eq 4.2-84 is the related homogeneous equation and $f_C(x)$ is its solution. Then,

$$[D_{\Delta x}{}^{n} f_C(x) + a_{n-1}D_{\Delta x}{}^{n-1}f_C(x) + a_{n-2}D_{\Delta x}{}^{n-2}f_C(x) + \ldots + a_1D_{\Delta x}f_C(x) + a_0f_C(x)] = 0 \qquad (4.2\text{-}85)$$

Note that Eq 4.2-85 is Eq 4.2-80

Substituting Eq 4.2-85 into Eq 4.2-84

$D_{\Delta x}{}^{n} f_P(x) + a_{n-1}D_{\Delta x}{}^{n-1}f_P(x) + a_{n-2}D_{\Delta x}{}^{n-2}f_P(x) + \ldots + a_1D_{\Delta x}f_P(x) + a_0f_P(x) = Q(x)$ (4.2-86)

Note that Eq 4.2-86 is Eq 4.2-81

The objective now is to find the function, $f_P(x)$. $f_P(x)$, the particular solution of Eq 4.2-79, is that function, the sum of Inteval Calculus functions, which when operated on in accordance with the left side of a specified differential difference equation of the form of Eq 4.2-79, yields the function, Q(x).

Finding a function, $f_P(x)$ for any function Q(x) could be difficult. However, a convenient method to find $f_P(x)$ for the most commonly used functions represented by Q(x) can be found with minimum difficulty.

In the following table, Table 4.2-4, those functions which are most commonly represented by Q(x) are listed. Some of these functions are not themselves easily differentiated but their identities are. For these functions, their identities are provided.

**Table 4.2-4 Commonly used Q(x) functions for which a Particular Solution is required**

| # | Commonly used functions represented by Q(x) | Interval Calculus Function Identities |
|---|---|---|
| **1** | c a constant | |
| **2** | $[x]^{n}_{\Delta x}$ equivalent to $\prod_{m=1}^{n}(x-[m-1]\Delta x)$<br>n = 1,2,3,… | |
| **3** | x | $[x]^{1}_{\Delta x}$ |
| **4** | $x^2$ | $[x]^{2}_{\Delta x} + \Delta x[x]^{1}_{\Delta x}$ |
| **5** | $x^3$ | $[x]^{3}_{\Delta x} + (3\Delta x)[x]^{2}_{\Delta x} + (\Delta x^2)[x]^{1}_{\Delta x}$ |
| **6** | $e_{\Delta x}(a,x)$ | |
| **7** | $\sin_{\Delta x}(b,x)$ | |
| **8** | $\cos_{\Delta x}(b,x)$ | |
| **9** | $e_{\Delta x}(a,x)\ \sin_{\Delta x}(b,x)$ | |
| **10** | $e_{\Delta x}(a,x)\ \cos_{\Delta x}(b,x)$ | |
| **11** | $\sinh_{\Delta x}(b,x)$ | |
| **12** | $\cosh_{\Delta x}(b,x)$ | |
| **13** | $[x]_{\Delta x}^{n}\, e_{\Delta x}(a,x)$ | |
| **14** | $[x]_{\Delta x}^{n}\, \sin_{\Delta x}(b,x)$ | |
| **15** | $[x]_{\Delta x}^{n} \cos_{\Delta x}(b,x)$ | |
| **16** | $[x]_{\Delta x}^{n} e_{\Delta x}(a,x)\, \sin_{\Delta x}(b,x)$ | |
| **17** | $[x]_{\Delta x}^{n} e_{\Delta x}(a,x)\, \cos_{\Delta x}(b,x)$ | |
| **18** | $A^x$ | $e_{\Delta x}( \frac{A^{\Delta x}-1}{\Delta x}, x )$ |
| **19** | $e^{ax}$ | $e_{\Delta x} ( \frac{e^{a\Delta x}-1}{\Delta x}, x )$ |
| **20** | sinbx | $e_{\Delta x}( \frac{\cos b\Delta x -1}{\Delta x}, x) \sin_{\Delta x}( \frac{\tan b\Delta x}{\Delta x}, x )$ <br> where <br> $\frac{x}{\Delta x}$ = integer <br> $\cos b\Delta x \neq 0$ |
| **21** | cosbx | $e_{\Delta x}( \frac{\cos b\Delta x -1}{\Delta x}, x) \cos_{\Delta x}( \frac{\tan b\Delta x}{\Delta x}, x )$ <br> where <br> $\frac{x}{\Delta x}$ = integer <br> $\cos b\Delta x \neq 0$ |
| **22** | $e^{ax}\sin bx$ | $e_{\Delta x}(\frac{e^{a\Delta x}\cos b\Delta x-1}{\Delta x}, x) \sin_{\Delta x}( \frac{\tan b\Delta x}{\Delta x}, x )$ <br> where <br> $\frac{x}{\Delta x}$ = integer <br> $\cos b\Delta x \neq 0$ |
| **23** | $e^{ax}\cos bx$ | $e_{\Delta x}(\frac{e^{a\Delta x}\cos b\Delta x-1}{\Delta x}, x)\ \cos_{\Delta x}(\frac{\tan b\Delta x}{\Delta x}, x)$ <br> where <br> $\frac{x}{\Delta x}$ = integer <br> $\cos b\Delta x \neq 0$ |
| **24** | $\sinh bx$ | $e_{\Delta x}(\frac{\cosh b\Delta x -1}{\Delta x}, x)\ \sinh_{\Delta x}(\frac{\tanh b\Delta x}{\Delta x}, x)$ |
| **25** | $\cosh bx$ | $e_{\Delta x}(\frac{\cosh b\Delta x -1}{\Delta x}, x)\ \cosh_{\Delta x}(\frac{\tanh b\Delta x}{\Delta x}, x)$ |
| **26** | $[x]_{\Delta x}^{n} e^{ax}$ | $[x]_{\Delta x}^{n}\ e_{\Delta x}(\frac{e^{a\Delta x}-1}{\Delta x}, x)$ |
| **27** | $\sin\frac{\pi x}{2\Delta x}$ | |
| **28** | $\cos\frac{\pi x}{2\Delta x}$ | |

Note – For generality, all of the functions in Table 4.2-4 above should be multiplied by a constant, K.

The Interval Calculus identities for $e^{ax}$, sinbx, cosbx, etc. presented in Table 4.2-4 can be used to expedite the discrete differentiation of these functions. Discrete differentiation of these functions can be awkward and time consuming. By converting these functions into their identities, discrete differentiation of the identities becomes simple and direct thereby reducing differentiation mathematical manipulation effort. To better understand these Interval Calculus identities, the identities for the functions, $e^{ax}$ and cosbx will be derived below.

Derivation of the $e^{ax}$ Interval Calculus Identity

$$e^{ax} = [1 + e^{a\Delta x} - 1\ ]^{\frac{x}{\Delta x}} = [1 + (\frac{e^{a\Delta x}-1}{\Delta x})\Delta x]^{\frac{x}{\Delta x}} = e_{\Delta x}(\frac{e^{a\Delta x}-1}{\Delta x}, x) \qquad (4.2\text{-}87)$$

$$e^{ax} = e_{\Delta x}(\frac{e^{a\Delta x}-1}{\Delta x}, x) \qquad (4.2\text{-}88)$$

Derivation of the cosbx Interval Calculus Identity

$$\cos bx = \frac{e^{jbx} + e^{-jbx}}{2} = \frac{[1 + e^{jb\Delta x} - 1]^{\frac{x}{\Delta x}} + [1 + e^{-jb\Delta x} - 1]^{\frac{x}{\Delta x}}}{2} \tag{4.2-89}$$

$$e^{jbx} = [e^{jb\Delta x}]^{\frac{x}{\Delta x}} \quad \text{for} \quad \frac{x}{\Delta x} = \text{integer} \tag{4.2-90}$$

In use, this Interval Calculus identity will have an $\frac{x}{\Delta x}$ = integer since $x = 0, \Delta x, 2\Delta x, 3\Delta x, \ldots$

and $\frac{x}{\Delta x} = 0, 1, 2, 3, \ldots$

$$\cos bx = \frac{[1 + (\frac{e^{jb\Delta x} - 1}{\Delta x})\Delta x]^{\frac{x}{\Delta x}} + [1 + (\frac{e^{-jb\Delta x} - 1}{\Delta x})\Delta x]^{\frac{x}{\Delta x}}}{2} \tag{4.2-91}$$

$$\cos bx = \frac{[1 + (\frac{\cos bx - 1 + j\sin bx}{\Delta x})\Delta x]^{\frac{x}{\Delta x}} + [1 + (\frac{\cos bx - 1 - j\sin bx}{\Delta x})\Delta x]^{\frac{x}{\Delta x}}}{2} \tag{4.2-92}$$

$$\cos bx = \frac{[1 + (\frac{\cos bx - 1}{\Delta x} + \frac{j\sin bx}{\Delta x})\Delta x]^{\frac{x}{\Delta x}} + [1 + (\frac{\cos bx - 1}{\Delta x} - \frac{j\sin bx}{\Delta x})\Delta x]^{\frac{x}{\Delta x}}}{2} \tag{4.2-93}$$

Let $u = \frac{\cos bx - 1}{\Delta x}$

$v = \frac{\sin bx}{\Delta x}$

Substituting into Eq 4.2-93

$$\cos bx = \frac{[1 + (u + jv)\Delta x]^{\frac{x}{\Delta x}} + [1 + (u - jv)\Delta x]^{\frac{x}{\Delta x}}}{2} \tag{4.2-94}$$

$$\cos bx = [1+u\Delta x]^{\frac{x}{\Delta x}} \left( \frac{[\frac{1 + u\Delta x + jv\Delta x}{1+u\Delta x}]^{\frac{x}{\Delta x}} + [\frac{1 + u\Delta x - jv\Delta x}{1+u\Delta x}]^{\frac{x}{\Delta x}}}{2} \right) \tag{4.2-95}$$

$$\cos bx = [1+u\Delta x]^{\frac{x}{\Delta x}} \left( \frac{[1 + j\frac{v\Delta x}{1+u\Delta x}]^{\frac{x}{\Delta x}} + [1 - j\frac{v\Delta x}{1+u\Delta x}]^{\frac{x}{\Delta x}}}{2} \right) \tag{4.2-96}$$

$$\cos bx = [1+u\Delta x]^{\frac{x}{\Delta x}} \cos_{\Delta x}\left( \frac{v}{1+u\Delta x}, x \right) \tag{4.2-97}$$

Substituting $u = \frac{\cos bx - 1}{\Delta x}$ and $v = \frac{\sin bx}{\Delta x}$ into Eq 4.2-97

$$\cos bx = [1 + (\frac{\cos bx - 1}{\Delta x})\Delta x]^{\frac{x}{\Delta x}} \cos_{\Delta x}(\frac{\sin b\Delta x}{\Delta x \cos b\Delta x}, x) \tag{4.2-98}$$

$$\cos bx = [1 + (\frac{\cos bx - 1}{\Delta x})\Delta x]^{\frac{x}{\Delta x}} \cos_{\Delta x}(\frac{\tan b\Delta x}{\Delta x}, x)$$

$$\cos bx = e_{\Delta x}(\frac{\cos bx - 1}{\Delta x}, x)\cos_{\Delta x}(\frac{\tan b\Delta x}{\Delta x}, x) \tag{4.2-99}$$

where

$\frac{x}{\Delta x}$ = integer

$\cos b\Delta x \neq 0$

The Interval Calculus identities for the functions, $A^x$, sinbx, $e^{ax}$sinbx, $e^{ax}$cosbx, sinhbx, and coshbx can be derived in a similar manner.

Note – The identities of the functions, sinbx, cosbx, $e^{ax}$sinbx, and $e^{ax}$cosbx have a restriction that $\frac{x}{\Delta x}$ = integer. The identities for the functions, $A^x$, $e^{ax}$, sinhbx, and coshbx do not have this restriction.

Interestingly, all of the functions in Table 4.2-4 can be represented by either of two general functions. These functions are:

$$q_1(x) = [\sum_{p=1}^{m} A_p[x]_{\Delta x}^{p-1}]\, e_{\Delta x}(u,x)\, \sin_{\Delta x}(w,x) \;+\; [\sum_{p=1}^{m} B_p[x]_{\Delta x}^{p-1}]\, e_{\Delta x}(u,x)\, \cos_{\Delta x}(w,x) \qquad (4.2\text{-}100)$$

where

$A_p, B_p, u, w$ = real value constants

or

$$q_2(x) = [\sum_{p=1}^{m} A_p[x]_{\Delta x}^{p-1}][w\Delta x]^{\frac{x}{\Delta x}}\sin\frac{\pi x}{2\Delta x} + [\sum_{p=1}^{m} B_p[x]_{\Delta x}^{p-1}][w\Delta x]^{\frac{x}{\Delta x}}\cos\frac{\pi x}{2\Delta x} \qquad (4.2\text{-}101)$$

where

$A_p, B_p, w$ = real value constants

For demonstation purposes, show the relationship of various functions of Table 4.2-4 to Eq 4.2-100 or Eq 4.2-101.

| Function | Relationship to Eq 4.2-100 | $q_1(x)$ |
|---|---|---|
| $e_{\Delta x}(a,x)$ | $m = 1$, $u = a$, $w = 0$, $B_1 = K$ | $Ke_{\Delta x}(a,x)$ |
| $\sin_{\Delta x}(b,x)$ | $m = 1$, $u = 0$, $A_1 = K$, $B_1 = 0$ <br> $w = b$ | $K\sin_{\Delta x}(b,x)$ |
| $\cos_{\Delta x}(b,x)$ | $m = 1$, $u = 0$, $A_1 = 0$, $B_1 = K$ <br> $w = b$ | $K\cos_{\Delta x}(b,x)$ |
| $[x]_{\Delta x}^{1} e_{\Delta x}(a,x)\sin_{\Delta x}(b,x)$ | $m = 2$, $u = a$, $w = b$, $A_1 = 0$ <br> $A_2 = K$, $B_1 = 0$, $B_2 = 0$ | $K[x]_{\Delta x}^{1} e_{\Delta x}(a,x)\sin_{\Delta x}(b,x)$ |
| $e^{ax}\cos bx =$ <br> $e_{\Delta x}(\frac{e^{a\Delta x}\cos b\Delta x - 1}{\Delta x}, x)\cos_{\Delta x}(\frac{\tan b\Delta x}{\Delta x}, x)$ | $m = 1$, $u = \frac{e^{a\Delta x}\cos b\Delta x - 1}{\Delta x}$ <br> $A_1 = 0$, $B_1 = K$, $w = \frac{\tan b\Delta x}{\Delta x}$ | $Ke_{\Delta x}(\frac{e^{a\Delta x}\cos b\Delta x - 1}{\Delta x}, x)\cos_{\Delta x}(\frac{\tan b\Delta x}{\Delta x}, x)$ |
| $x^2 =$ <br> $[x]_{\Delta x}^{2} + \Delta x[x]_{\Delta x}^{1}$ | $m = 3$, $u = 0$, $w = 0$, $B_1 = 0$ <br> $B_2 = K\Delta x$, $B_3 = K$ | $K([x]_{\Delta x}^{2} + \Delta x[x]_{\Delta x}^{1})$ |

| Function | Relationship to Eq 4.2-101 | $q_2(x)$ |
|---|---|---|
| $\sin\frac{\pi x}{2\Delta x}$ | $m = 1$ , $A_1 = K$ , $B_1 = 0$<br>$w = \frac{1}{\Delta x}$ | $K\sin\frac{\pi x}{2\Delta x}$ |

Rewriting the general form of a differential difference equation, Eq 4.2-79

$$D_{\Delta x}{}^{n}f(x) + a_{n-1}D_{\Delta x}{}^{n-1}f(x) + a_{n-2}D_{\Delta x}{}^{n-2}f(x) + \ldots + a_1D_{\Delta x}f(x) + a_0f(x) = Q(x)$$

where

$n$ = order of the homogeneous equation

$a_{n-1},a_{n-1},\ldots,a_0$ = real constants

$f(x) = f_C(x) + f_P(x)$ , the differential difference equation general solution

$f_C(x)$ = differential difference equation complementary solution

$f_P(x)$ = differential difference equation particular solution

Since Eq 4.2-100 and Eq 4.2-101 can represent so many useful functions, especially those in Table 4.2-4, let the function, Q(x), in the differential difference equation, Eq 4.2-79, be equal to Eq 4.2-100 or Eq 4.2-101. Then, finding the particular solution for the resulting two differential difference equations will yield the particular solutions for all the differential difference equations that both Eq 4.2-100 and Eq 4.2-101 can represent.

Let

$$\underline{Q(x) = q_1(x)} \tag{4.2-102}$$

Changing the form of Eq 4.2-79

$$(D_{\Delta x}-r_n)(D_{\Delta x}-r_{n-1})(D_{\Delta x}-r_{n-2})\ldots(D_{\Delta x}-r_3)(D_{\Delta x}-r_2)(D_{\Delta x}-r_1)f(x) = Q(x) \tag{4.2-103}$$

where

$h(r) = r^n + a_{n-1}r^{n-1} + a_{n-2}r^{n-2} + \ldots + a_2r^2 + a_1r + a_0 = (r-r_n)(r-r_{n-1})(r-r_{n-2})\ldots(r-r_3)(r-r_2)(r-r_1) = 0$

$r = r_1, r_2, r_3, \ldots, r_{n-1}, r_n$ , the n roots of h(r)

Substituting Eq 4.2-102 into Eq 4.2-103

$$(D_{\Delta x}-r_n)(D_{\Delta x}-r_{n-1})(D_{\Delta x}-r_{n-2})\ldots(D_{\Delta x}-r_3)(D_{\Delta x}-r_2)(D_{\Delta x}-r_1)f(x) = q_1(x) \tag{4.2-104}$$

Find the particular solution to the above differential difference equation, Eq 4.2-104.

Consider the following homogeneous differential difference equation:

$$[(D_{\Delta x}-\{a+jb\})(D_{\Delta x}-\{a-jb\})]^m F(x) = 0 \tag{4.2-105}$$

where

$m$ = multiplicity of the complex conjugate root pair, a+jb and a–jb

a,b = real constants

$$a \neq -\frac{1}{\Delta x}$$

From Table 4.2-3, the general solution to Eq 4.2-105 is easily found.

For $a \pm jb$ where $a \neq -\frac{1}{\Delta x}$

$$F(x) = [\sum_{p=1}^{m} A_p[x]_{\Delta x}^{p-1}]\, e_{\Delta x}(a,x)\, \sin_{\Delta x}(\frac{b}{1+a\Delta x}, x) + [\sum_{p=1}^{m} B_p[x]_{\Delta x}^{p-1}]\, e_{\Delta x}(a,x)\, \cos_{\Delta x}(\frac{b}{1+a\Delta x}, x) \qquad (4.2\text{-}106)$$

Rewriting the equation for $q_1(x)$, Eq 4.2-100

$$q_1(x) = [\sum_{p=1}^{m} A_p[x]_{\Delta x}^{p-1}]\, e_{\Delta x}(u,x)\, \sin_{\Delta x}(w,x) + [\sum_{p=1}^{m} B_p[x]_{\Delta x}^{p-1}]\, e_{\Delta x}(u,x)\, \cos_{\Delta x}(w,x)$$

$F(x)$ is seen to be very similar to $q_1(x)$

Let

$$w = \frac{b}{1+a\Delta x}$$

$$u = a$$

Solving for a and b in terms of u and w

$$a = u \qquad (4.2\text{-}107)$$

$$b = w(1+u\Delta x) \qquad (4.2\text{-}108)$$

Substituting Eq 4.2-107 and Eq 4.2-108 into Eq 4.2-105 and Eq 4.2-106

$$[(D_{\Delta x}-\{u+jw(1+u\Delta x)\})(D_{\Delta x}-\{u-jw(1-u\Delta x)\})]^m\, F(x) = 0 \qquad (4.2\text{-}109)$$

$$F(x) = [\sum_{p=1}^{m} A_p[x]_{\Delta x}^{p-1}]\, e_{\Delta x}(u,x)\, \sin_{\Delta x}(w,x) + [\sum_{p=1}^{m} B_p[x]_{\Delta x}^{p-1}]\, e_{\Delta x}(u,x)\, \cos_{\Delta x}(w,x) = q_1(x) \qquad (4.2\text{-}110)$$

Note that $F(x)$ is now equal to $q_1(x)$

From Eq 4.2-109 and Eq 4.2-110 it is seen that

$$q_1(x) = [\sum_{p=1}^{m} A_p[x]_{\Delta x}^{p-1}]\ e_{\Delta x}(u,x)\ \sin_{\Delta x}(w,x) + [\sum_{p=1}^{m} B_p[x]_{\Delta x}^{p-1}]\ e_{\Delta x}(u,x)\ \cos_{\Delta x}(w,x) \qquad (4.2\text{-}111)$$

is the general solution to the homogeneous differential difference equation

$$[(D_{\Delta x}-\{u+jw(1+u\Delta x)\})(D_{\Delta x}-\{u-jw(1-u\Delta x)\})]^{m}\ q_1(x) = 0 \qquad (4.2\text{-}112)$$

where

$A_p,B_p,u,w$ = real constants

$m$ = multiplicity of the complex conjugate root pair, $u+jw(1+u\Delta x)$ and $u-jw(1+u\Delta x)$

$\Delta x$ = x increment

$u \neq -\frac{1}{\Delta x}$

The roots $u+jw(1+u\Delta x)$ and $u-jw(1+u\Delta x)$ are henceforth referred to as the related roots of $q_1(x)$.

Note - m may also be determined from the order of the summation polynomials in Eq 4.2-111.

$m = N + 1$

$N$ = the order of the summation polynomials in $q_1(x)$

Rewriting Eq 4.2-104

$(D_{\Delta x}-r_n)(D_{\Delta x}-r_{n-1})(D_{\Delta x}-r_{n-2})\ldots(D_{\Delta x}-r_3)(D_{\Delta x}-r_2)(D_{\Delta x}-r_1)f(x) = q_1(x)$

In the following derivation, the related roots of $q_1(x)$ , $u+jw(1+u\Delta x)$ and $u-jw(1+u\Delta x)$, are not equal to any of the roots in Eq 4.2-104.

Differentiating each side of Eq 4.2-104 using the derivative operator,

$[(D_{\Delta x}-\{u+jw(1+u\Delta x)\})(D_{\Delta x}-\{u-jw(1-u\Delta x)\})]^{m}$

$$[(D_{\Delta x}-\{u+jw(1+u\Delta x)\})(D_{\Delta x}-\{u-jw(1-u\Delta x)\})]^{m}\ (D_{\Delta x}-r_n)(D_{\Delta x}-r_{n-1})(D_{\Delta x}-r_{n-2})\ldots(D_{\Delta x}-r_3)(D_{\Delta x}-r_2)(D_{\Delta x}-r_1)f(x)$$
$$= [(D_{\Delta x}-\{u+jw(1+u\Delta x)\})(D_{\Delta x}-\{u-jw(1-u\Delta x)\})]^{m}\ q_1(x) \qquad (4.2\text{-}113)$$

From Eq 4.2-112 and 4.2-113

$$\{[(D_{\Delta x}-\{u+jw(1+u\Delta x)\})(D_{\Delta x}-\{u-jw(1-u\Delta x)\})]^{m}\}\{(D_{\Delta x}-r_n)(D_{\Delta x}-r_{n-1})(D_{\Delta x}-r_{n-2})\ldots(D_{\Delta x}-r_3)(D_{\Delta x}-r_2)(D_{\Delta x}-r_1)\}f(x) = 0 \qquad (4.2\text{-}114)$$

Eq 4.2-104, a non-homgeneous equation, has now been converted into a homogeneous equation, Eq 4.2-114. The function, f(x), the general solution to both equations, remains the same.

As previously stated, the function f(x), the general solution of Eq 4.2-79, has two components, a complementary solution, $f_C(x)$, and a particular solution, $f_P(x)$.

$$f(x) = f_C(x) + f_P(x) \qquad (4.2\text{-}115)$$

Comparing Eq 4.2-114 to Eq 4.2-104, the Eq 4.2-114 term, $(D_{\Delta x}\text{-}r_n)(D_{\Delta x}\text{-}r_{n-1})(D_{\Delta x}\text{-}r_{n-2})\ldots(D_{\Delta x}\text{-}r_3)(D_{\Delta x}\text{-}r_2)(D_{\Delta x}\text{-}r_1)$ is immediately recognized as the portion of Eq 4.2-114, that contains roots which will generate in f(x) the complementary solution, $f_C(x)$ (from Eq 4.2-104, the related homogeneous equation is $(D_{\Delta x}\text{-}r_n)(D_{\Delta x}\text{-}r_{n-1})(D_{\Delta x}\text{-}r_{n-2})\ldots(D_{\Delta x}\text{-}r_3)(D_{\Delta x}\text{-}r_2)(D_{\Delta x}\text{-}r_1)f_C(x) = 0$). The roots of the remaining Eq 4.2-114 term, $[(D_{\Delta x}-\{u+jw(1+u\Delta x)\})(D_{\Delta x}-\{u-jw(1-u\Delta x)\})]^m$, must then be those that will yield the particular solution, $f_P(x)$ (this is seen from Eq 4.2-112 and Eq 4.2-111). Then, the particular solution to Eq 4.2-104 has been found. It is $q_1(x)$ with the appropriate $A_p$ and $B_p$ constant values provided. Q(x), as previously shown, is also formed from $q_1(x)$ but, in general, with other values for the constants, $A_p$ and $B_p$.

Thus

**For the differential difference equation, 4.2-79, rewritten below**

$$D_{\Delta x}{}^n f(x) + a_{n-1}D_{\Delta x}{}^{n-1}f(x) + a_{n-2}D_{\Delta x}{}^{n-2}f(x) + \ldots + a_1D_{\Delta x}f(x) + a_0f(x) = Q(x)$$

**where**

$$h(r) = r^n + a_{n-1}r^{n-1} + a_{n-2}r^{n-2} + \ldots + a_2r^2 + a_1r + a_0 = (r-r_n)(r-r_{n-1})(r-r_{n-2})\ldots(r-r_3)(r-r_2)(r-r_1)$$

**if** $Q(x) = [\sum_{p=1}^{m} A_p[x]_{\Delta x}^{p-1}]\, e_{\Delta x}(u,x)\, \sin_{\Delta x}(w,x) + [\sum_{p=1}^{m} B_p[x]_{\Delta x}^{p-1}]\, e_{\Delta x}(u,x)\, \cos_{\Delta x}(w,x)$ **, the particular solution to Eq 4.2-79 is:** (4.2-116)

$$f_P(x) = [\sum_{p=1}^{m} \mathcal{A}_p[x]_{\Delta x}^{p-1}]\, e_{\Delta x}(u,x)\, \sin_{\Delta x}(w,x) + [\sum_{p=1}^{m} \mathcal{B}_p[x]_{\Delta x}^{p-1}]\, e_{\Delta x}(u,x)\, \cos_{\Delta x}(w,x) \qquad (4.2\text{-}117)$$

**where**

**n = order of the differential difference equation**

$A_p,B_p,\mathcal{A}_p,\mathcal{B}_p,a_{n-1}\ldots a_0,u,w$ **= real value constants**

$f_P(x)$ **= differential difference equation particular solution**

**m = N + 1**

**N = the order of the summation polynomials of Q(x)**

**Δx = x interval**

$u \neq -\frac{1}{\Delta x}$

**The related roots of the Q(x) function, u+jw(1+uΔx) and u–jw(1+uΔx), are not equal to any of the roots of the characteristic polynomial, h(r).**

As previously shown, the Q(x) function, Eq 4.2-116, with the proper selection of its constants, $A_p$,$B_p$,u,w, and m, can represent many useful functions.
The function, Q(x), of the differential difference equation, Eq 4.2-79, may also equal the function,

$$q_2(x) = [\sum_{p=1}^{m} A_p[x]_{\Delta x}^{p-1}][w\Delta x]^{\frac{x}{\Delta x}}\sin\frac{\pi x}{2\Delta x} + [\sum_{p=1}^{m} B_p[x]_{\Delta x}^{p-1}][w\Delta x]^{\frac{x}{\Delta x}}\cos\frac{\pi x}{2\Delta x} \qquad (4.2\text{-}118)$$

where

$A_p$,$B_p$,w = real value constants

$\Delta x$ = x interval

This case will now be considered. The argument used above for $Q(x) = q_1(x)$ will now be repeated for

<u>$Q(x) = q_2(x)$</u>.

$$Q(x) = q_2(x) \qquad (4.2\text{-}119)$$

Changing the form of Eq 4.2-79

$$(D_{\Delta x}-r_n)(D_{\Delta x}-r_{n-1})(D_{\Delta x}-r_{n-2})\ldots(D_{\Delta x}-r_3)(D_{\Delta x}-r_2)(D_{\Delta x}-r_1)f(x) = Q(x) \qquad (4.2\text{-}120)$$

where

$h(r) = r^n + a_{n-1}r^{n-1} + a_{n-2}r^{n-2} + \ldots + a_2r^2 + a_1r + a_0 = (r-r_n)(r-r_{n-1})(r-r_{n-2})\ldots(r-r_3)(r-r_2)(r-r_1) = 0$

$r = r_1, r_2, r_3, \ldots, r_{n-1}, r_n$ , the n roots of h(r)

Substituting Eq 4.2-119 into Eq 4.2-120

$$(D_{\Delta x}-r_n)(D_{\Delta x}-r_{n-1})(D_{\Delta x}-r_{n-2})\ldots(D_{\Delta x}-r_3)(D_{\Delta x}-r_2)(D_{\Delta x}-r_1)f(x) = q_2(x) \qquad (4.2\text{-}121)$$

Find the particular solution to the above differential difference equation, Eq 4.2-121.

Consider the following homogeneous differential difference equation:

$$[(D_{\Delta x}-\{a+jb\})(D_{\Delta x}-\{a-jb\})]^m F(x) = 0 \qquad (4.2\text{-}122)$$

where

m = multiplicity of the complex conjugate root pair, a+jb, a–jb

a,b = real constants

$\Delta x$ = x increment

$a = -\frac{1}{\Delta x}$

From Table 4.2-3, the general solution to Eq 4.2-122 is easily found.

For $a \pm jb$ where $a = -\frac{1}{\Delta x}$

$$F(x) = [\sum_{p=1}^{m} A_p[x]_{\Delta x}^{p-1}][b\Delta x]^{\frac{x}{\Delta x}} \sin\frac{\pi x}{2\Delta x} + [\sum_{p=1}^{m} B_p[x]_{\Delta x}^{p-1}][b\Delta x]^{\frac{x}{\Delta x}} \cos\frac{\pi x}{2\Delta x} \qquad (4.2\text{-}123)$$

Rewriting the equation for $q_2(x)$, Eq 4.2-118

$$q_2(x) = [\sum_{p=1}^{m} A_p[x]_{\Delta x}^{p-1}][w\Delta x]^{\frac{x}{\Delta x}} \sin\frac{\pi x}{2\Delta x} + [\sum_{p=1}^{m} B_p[x]_{\Delta x}^{p-1}][w\Delta x]^{\frac{x}{\Delta x}} \cos\frac{\pi x}{2\Delta x}$$

$F(x)$ is seen to be very similar to $q_2(x)$

Let

$$b = w \qquad (4.2\text{-}124)$$

Substituting Eq 4.2-124 and $a = -\frac{1}{\Delta x}$ into Eq 4.2-122 and Eq 4.2-123

$$[(D_{\Delta x} - \{-\frac{1}{\Delta x} + jw\})(D_{\Delta x} - \{-\frac{1}{\Delta x} - jw\})]^m F(x) = 0 \qquad (4.2\text{-}125)$$

$$F(x) = [\sum_{p=1}^{m} A_p[x]_{\Delta x}^{p-1}][w\Delta x]^{\frac{x}{\Delta x}} \sin\frac{\pi x}{2\Delta x} + [\sum_{p=1}^{m} B_p[x]_{\Delta x}^{p-1}][w\Delta x]^{\frac{x}{\Delta x}} \cos\frac{\pi x}{2\Delta x} = q_2(x) \qquad (4.2\text{-}126)$$

Note that $F(x)$ is now equal to $q_2(x)$

Then from Eq 4.2-125 and Eq 4.2-126

$$q_2(x) = [\sum_{p=1}^{m} A_p[x]_{\Delta x}^{p-1}][w\Delta x]^{\frac{x}{\Delta x}} \sin\frac{\pi x}{2\Delta x} + [\sum_{p=1}^{m} B_p[x]_{\Delta x}^{p-1}][w\Delta x]^{\frac{x}{\Delta x}} \cos\frac{\pi x}{2\Delta x} \qquad (4.2\text{-}127)$$

is the general solution to the homogeneous differential difference equation

$$[(D_{\Delta x} - \{-\frac{1}{\Delta x} + jw\})(D_{\Delta x} - \{-\frac{1}{\Delta x} - jw\})]^m q_2(x) = 0 \qquad (4.2\text{-}128)$$

where

$A_p, B_p, w$ = real constants

m = multiplicity of the complex conjugate root pair, $-\frac{1}{\Delta x} + jw$ and $-\frac{1}{\Delta x} - jw$

$\Delta x$ = x increment

$u = -\frac{1}{\Delta x}$

The roots $-\frac{1}{\Delta x} + jw$ and $-\frac{1}{\Delta x} - jw$ are henceforth referred to as the related roots of $q_2(x)$.

Note - m may also be determined from the order of the summation polynomials in Eq 4.2-127.

m = N + 1

N = the order of the summation polynomials in $q_2(x)$

Rewriting Eq 4.2-121

$(D_{\Delta x}-r_n)(D_{\Delta x}-r_{n-1})(D_{\Delta x}-r_{n-2})\ldots(D_{\Delta x}-r_3)(D_{\Delta x}-r_2)(D_{\Delta x}-r_1)f(x) = q_2(x)$

In the following derivation, the related roots of $q_1(x)$ , $-\frac{1}{\Delta x} + jw$ and $-\frac{1}{\Delta x} - jw$ , are not equal to any of the roots in Eq 4.2-121.

Differentiating each side of Eq 4.2-121 using the derivative operator,

$[(D_{\Delta x} - \{-\frac{1}{\Delta x} + jw\})(D_{\Delta x} - \{-\frac{1}{\Delta x} - jw\})]^m$

$$[(D_{\Delta x} - \{-\frac{1}{\Delta x} + jw\})(D_{\Delta x} - \{-\frac{1}{\Delta x} - jw\})]^m (D_{\Delta x}-r_n)(D_{\Delta x}-r_{n-1})(D_{\Delta x}-r_{n-2})\ldots(D_{\Delta x}-r_3)(D_{\Delta x}-r_2)(D_{\Delta x}-r_1)f(x) =$$
$$[(D_{\Delta x} - \{-\frac{1}{\Delta x} + jw\})(D_{\Delta x} - \{-\frac{1}{\Delta x} - jw\})]^m F(x)\, q_2(x) \qquad (4.2\text{-}129)$$

From Eq 4.2-128 and 4.2-129

$$[(D_{\Delta x} - \{-\frac{1}{\Delta x} + jw\})(D_{\Delta x} - \{-\frac{1}{\Delta x} - jw\})]^m (D_{\Delta x}-r_n)(D_{\Delta x}-r_{n-1})(D_{\Delta x}-r_{n-2})\ldots(D_{\Delta x}-r_3)(D_{\Delta x}-r_2)(D_{\Delta x}-r_1)f(x) = 0$$
$$(4.2\text{-}130)$$

Eq 4.2-121, a non-homgeneous equation, has now been converted into a homogeneous equation, Eq 4.2-130. The function, f(x), the general solution to both equations, remains the same.

As previously stated, the function f(x), the general solution of Eq 4.2-108, has two components, a complementary solution, $f_C(x)$, and a particular solution, $f_P(x)$.

$$f(x) = f_C(x) + f_P(x) \qquad (4.2\text{-}131)$$

Comparing Eq 4.2-130 to Eq 4.2-121, the Eq 4.2-130 term, $(D_{\Delta x}\text{-}r_n)(D_{\Delta x}\text{-}r_{n-1})(D_{\Delta x}\text{-}r_{n-2})\ldots(D_{\Delta x}\text{-}r_3)(D_{\Delta x}\text{-}r_2)(D_{\Delta x}\text{-}r_1)$ is immediately recognized as the portion of Eq 4.2-130, that contains roots which will generate in f(x) the complementary solution, $f_C(x)$ (from Eq 4.2-121, the related homogeneous equation is $(D_{\Delta x}\text{-}r_n)(D_{\Delta x}\text{-}r_{n-1})(D_{\Delta x}\text{-}r_{n-2})\ldots(D_{\Delta x}\text{-}r_3)(D_{\Delta x}\text{-}r_2)(D_{\Delta x}\text{-}r_1)f_C(x) = 0$). The roots of the remaining Eq 4.2-130 term, $[(D_{\Delta x}-\{u+jw(1+u\Delta x)\})(D_{\Delta x}-\{u-jw(1-u\Delta x)\})]^m$, must then be those that will yield the particular solution, $f_P(x)$ (this is seen from Eq 4.2-128 and Eq 4.2-127). Then, the particular solution to Eq 4.2-121 has been found. It is $q_2(x)$ with the appropriate $A_p$ and $B_p$ constant values provided. Q(x), as previously shown, is also formed from $q_2(x)$ but, in general, with other values for the constants, $A_p$ and $B_p$.

Thus

**For the differential difference equation, 4.2-79, rewritten below**

$$D_{\Delta x}^{n}f(x) + a_{n-1}D_{\Delta x}^{n-1}f(x) + a_{n-2}D_{\Delta x}^{n-2}f(x) + \ldots + a_1D_{\Delta x}f(x) + a_0f(x) = Q(x)$$

**where**

$$h(r) = r^n + a_{n-1}r^{n-1} + a_{n-2}r^{n-2} + \ldots + a_2r^2 + a_1r + a_0 = (r-r_n)(r-r_{n-1})(r-r_{n-2})\ldots(r-r_3)(r-r_2)(r-r_1)$$

**if** $Q(x) = [\sum_{p=1}^{m} A_p[x]_{\Delta x}^{p-1}][w\Delta x]^{\frac{x}{\Delta x}}\sin\frac{\pi x}{2\Delta x} + [\sum_{p=1}^{m} B_p[x]_{\Delta x}^{p-1}][w\Delta x]^{\frac{x}{\Delta x}}\cos\frac{\pi x}{2\Delta x}$**, the particular solution to Eq 4.2-79 is:** (4.2-132)

$$f_P(x) = [\sum_{p=1}^{m} \mathcal{A}_p[x]_{\Delta x}^{p-1}][w\Delta x]^{\frac{x}{\Delta x}}\sin\frac{\pi x}{2\Delta x} + [\sum_{p=1}^{m} \mathcal{B}_p[x]_{\Delta x}^{p-1}][w\Delta x]^{\frac{x}{\Delta x}}\cos\frac{\pi x}{2\Delta x} \quad (4.2\text{-}133)$$

**where**

**n = order of the differential difference equation**

$A_p, B_p, \mathcal{A}_p, \mathcal{B}_p, a_{n-1}\ldots a_0, w$ **= real value constants**

$f_P(x)$ **= differential difference equation particular solution**

**m = N + 1**

**N = the order of the summation polynomials of Q(x)**

**$\Delta x$ = x increment**

**In the above derivation, the related roots of the Q(x) function,** $-\frac{1}{\Delta x} + jw$ **and** $-\frac{1}{\Delta x} - jw$**, are not equal to any of the roots of the characteristic polynomial, h(r).**

As previously shown, the Q(x) function, Eq 4.2-132, with the proper selection of its constants, $A_p$, $B_p$, w, and m, can represent many useful functions.

In the following derivation, the related root of multiplicity, m, of the $Q(x) = q_1(x)$ function is equal to a root of the characteristic polynomial, h(r) of multiplicity, $m_h$.

$Q(x) = q_1(x)$ (4.2-134)

Changing the form of Eq 4.2-79

$$(D_{\Delta x}-r_n)(D_{\Delta x}-r_{n-1})(D_{\Delta x}-r_{n-2})\ldots(D_{\Delta x}-r_3)(D_{\Delta x}-r_2)(D_{\Delta x}-r_1)f(x) = Q(x) \quad (4.2\text{-}135)$$

where

$h(r) = r^n + a_{n-1}r^{n-1} + a_{n-2}r^{n-2} + \ldots + a_2r^2 + a_1r + a_0 = (r-r_n)(r-r_{n-1})(r-r_{n-2})\ldots(r-r_3)(r-r_2)(r-r_1) = 0$

$r = r_1, r_2, r_3, \ldots, r_{n-1}, r_n$, the n roots of h(r)

Substituting Eq 4.2-134 into Eq 4.2-135

$$(D_{\Delta x}-r_n)(D_{\Delta x}-r_{n-1})(D_{\Delta x}-r_{n-2})\ldots(D_{\Delta x}-r_3)(D_{\Delta x}-r_2)(D_{\Delta x}-r_1)f(x) = q_1(x) \quad (4.2\text{-}136)$$

From Eq 4.2-111 and Eq 4.2-112

$$q_1(x) = \left[\sum_{p=1}^{m} A_p[x]_{\Delta x}^{p-1}\right] e_{\Delta x}(u,x)\sin_{\Delta x}(w,x) + \left[\sum_{p=1}^{m} B_p[x]_{\Delta x}^{p-1}\right] e_{\Delta x}(u,x)\cos_{\Delta x}(w,x) \quad (4.2\text{-}137)$$

is the general solution to the homogeneous differential difference equation

$$[(D_{\Delta x}-\{u+jw(1+u\Delta x)\})(D_{\Delta x}-\{u-jw(1-u\Delta x)\})]^m q_1(x) = 0 \quad (4.2\text{-}138)$$

where

$A_p, B_p, u, w$ = real constants

m = multiplicity of the complex conjugate root pair, $u+jw(1+u\Delta x)$ and $u-jw(1+u\Delta x)$

$\Delta x$ = x increment

$u \neq -\frac{1}{\Delta x}$

Let several of the h(r) roots, $r_1$ thru $r_n$, be equal to some multiple, $m_h$, of the same complex conjugate root pair, $u+jw(1+u\Delta x)$ and $u-jw(1+u\Delta x)$ as in Eq 4.2-138.

Several of the h(r) roots are then

$$[(D_{\Delta x}-\{u+jw(1+u\Delta x)\})(D_{\Delta x}-\{u-jw(1-u\Delta x)\})]^{m_h} \quad (4.2\text{-}139)$$

Substituting Eq 4.2-139 into Eq 4.2-136

$[(D_{\Delta x}-\{u+jw(1+u\Delta x)\})(D_{\Delta x}-\{u-jw(1-u\Delta x)\})]^{m_h}(D_{\Delta x}-r_v)(D_{\Delta x}-r_{v-1})(D_{\Delta x}-r_{v-2})\ldots(D_{\Delta x}-r_3)(D_{\Delta x}-r_2)$

$(D_{\Delta x}\text{-}r_1)f(x) = q_1(x)$ (4.2-140)

where

n = order of the differential difference equation, Eq 4.2-140, and also the number of roots in its related homogeneous differential difference equation

$v = n - 2m_h$

The related roots of $q_1$ are a multiple, m, of the complex conjugate root pair, $u+jw(1+u\Delta x)$ and $u-jw(1+u\Delta x)$.

Take the derivative of both sides of Eq 4.2-140 using the derivative operator,

$[(D_{\Delta x}-\{u+jw(1+u\Delta x)\})(D_{\Delta x}-\{u-jw(1-u\Delta x)\})]^m$

$$[(D_{\Delta x}-\{u+jw(1+u\Delta x)\})(D_{\Delta x}-\{u-jw(1-u\Delta x)\})]^m[(D_{\Delta x}-\{u+jw(1+u\Delta x)\})(D_{\Delta x}-\{u-jw(1-u\Delta x)\})]^{m_h}(D_{\Delta x}\text{-}r_v)$$
$$(D_{\Delta x}\text{-}r_{v-1})(D_{\Delta x}\text{-}r_{v-2})\ldots(D_{\Delta x}\text{-}r_3)(D_{\Delta x}\text{-}r_2)(D_{\Delta x}\text{-}r_1)f(x) = [(D_{\Delta x}-\{u+jw(1+u\Delta x)\})(D_{\Delta x}-\{u-jw(1-u\Delta x)\})]^m q_1(x)$$
(4.2-141)

From Eq 4.2-138 and Eq 4.2-141

$$[(D_{\Delta x}-\{u+jw(1+u\Delta x)\})(D_{\Delta x}-\{u-jw(1-u\Delta x)\})]^{m+m_h}(D_{\Delta x}\text{-}r_v)(D_{\Delta x}\text{-}r_{v-1})(D_{\Delta x}\text{-}r_{v-2})\ldots(D_{\Delta x}\text{-}r_3)(D_{\Delta x}\text{-}r_2)(D_{\Delta x}\text{-}r_1)f(x) = 0$$
(4.2-142)

Eq 4.2-140, a non-homgeneous equation, has now been converted into a homogeneous equation, Eq 4.2-142. The function, f(x), the general solution to both equations, remains the same.

The undetermined coefficient solution to Eq 4.2-142 yields f(x) as:

$$f(x) = [\sum_{p=1}^{m+m_h} A_p[x]_{\Delta x}^{p-1}]\, e_{\Delta x}(u,x)\sin_{\Delta x}(w,x) + [\sum_{p=1}^{m+m_h} B_p[x]_{\Delta x}^{p-1}]\, e_{\Delta x}(u,x)\cos_{\Delta x}(w,x) + \qquad (4.2\text{-}143)$$
$$K_v e^{r_v x} + K_{v-1}e^{r_{v-1}x} + K_{v-2}e^{r_{v-2}x} + \ldots + K_3 e^{r_3 x} + K_2 e^{r_2 x} + K_1 e^{r_1 x}$$

or

$$f(x) = \{\, [\sum_{p=m_h+1}^{m+m_h} A_p[x]_{\Delta x}^{p-1}]\, e_{\Delta x}(u,x)\sin_{\Delta x}(w,x) + [\sum_{p=m_h+1}^{m+m_h} B_p[x]_{\Delta x}^{p-1}]\, e_{\Delta x}(u,x)\cos_{\Delta x}(w,x)\,\} + \qquad (4.2\text{-}144)$$
$$\{\, [\sum_{p=1}^{m_h} A_p[x]_{\Delta x}^{p-1}]\, e_{\Delta x}(u,x)\sin_{\Delta x}(w,x) + [\sum_{p=1}^{m_h} B_p[x]_{\Delta x}^{p-1}]\, e_{\Delta x}(u,x)\cos_{\Delta x}(w,x) +$$
$$K_v e^{r_v x} + K_{v-1}e^{r_{v-1}x} + K_{v-2}e^{r_{v-2}x} + \ldots + K_3 e^{r_3 x} + K_2 e^{r_2 x} + K_1 e^{r_1 x}\,\}$$

where

$A_p, B_p$,u,w = real constants

$K_1$ thru $K_v$ = real or complex constants

$m_h$ = the multiplicity of the h(r) characteristic polynomial complex conjugate roots, u+jw(1+uΔx) and u–jw(1+uΔx) which are equal to the related roots of $q_1(x)$

m = the multiplicity of the $q_1(x)$ related complex conjugate roots, u+jw(1+uΔx) and u–jw(1+uΔx)

Δx = x increment

$u \neq -\dfrac{1}{\Delta x}$

Identify the particular solution in the above equation, Eq 4.2-144.

The function f(x), the general solution of Eq 4.2-140, has two components, a complementary solution, $f_C(x)$, and a particular solution, $f_P(x)$.

$$f(x) = f_C(x) + f_P(x) \qquad (4.2\text{-}145)$$

The second term of Eq 4.2-144 is readily recognized as the complementary solution, $f_C(x)$. Then, the first term must be the particular solution, $f_P(x)$.

Thus

**For the differential difference equation, 4.2-79, rewritten below**

$$\mathbf{D_{\Delta x}^{n}f(x) + a_{n-1}D_{\Delta x}^{n-1}f(x) + a_{n-2}D_{\Delta x}^{n-2}f(x) + \ldots + a_1D_{\Delta x}f(x) + a_0f(x) = Q(x)} \qquad (4.2\text{-}146)$$

**where**

$$\mathbf{h(r) = r^n + a_{n-1}r^{n-1} + a_{n-2}\,r^{n-2} + \ldots + a_2r^2 + a_1r + a_0 = (r-r_n)(r-r_{n-1})(r-r_{n-2})\ldots(r-r_3)(r-r_2)(r-r_1)}$$

**if** $\mathbf{Q(x) = [\sum_{p=1}^{m} A_p[x]_{\Delta x}^{p-1}]\, e_{\Delta x}(u,x)\, \sin_{\Delta x}(w,x) + [\sum_{p=1}^{m} B_p[x]_{\Delta x}^{p-1}]\, e_{\Delta x}(u,x)\, \cos_{\Delta x}(w,x)}$**, the particular solution to Eq 4.2-146 is:** (4.2-147)

$$\mathbf{f_P(x) = [\sum_{p=m_h+1}^{m+m_h} \mathit{A}_p[x]_{\Delta x}^{p-1}]\, e_{\Delta x}(u,x)\, \sin_{\Delta x}(w,x) + [\sum_{p=m_h+1}^{m+m_h} \mathit{B}_p[x]_{\Delta x}^{p-1}]\, e_{\Delta x}(u,x)\, \cos_{\Delta x}(w,x)} \qquad (4.2\text{-}148)$$

**where**

**n = order of the differential difference equation**

$\mathbf{A_p, B_p, \mathit{A}_p, \mathit{B}_p, a_{n-1}\ldots a_0}$**,u,w = real value constants**

$\mathbf{f_P(x)}$ **= differential difference equation particular solution**

**$m_h$ = the multiplicity of the h(r) characteristic polynomial complex conjugate roots, u+jw(1+uΔx) and u−jw(1+uΔx) which are equal to the related roots of Q(x)**

**m = the multiplicity of the Q(x) related complex conjugate roots, u+jw(1+uΔx) and u−jw(1+uΔx)**

**Δx = x increment**

**$u \neq -\frac{1}{\Delta x}$**

**Note – The particular solution equation, Eq 4.2-128, is more general than Eq 4.2-117.**

As previously shown, the Q(x) function presented above, with the proper selection of its constants, $A_p$,$B_p$,u,w,m and $m_h$, can represent many useful functions.

In the following derivation, a root of multiplicity, $m_p$, of the Q(x) function is equal to a root of the characteristic polynomial, h(r) of multiplicity, $m_h$.

$Q(x) = q_2(x)$ (4.2-149)

Changing the form of Eq 4.2-79

$(D_{\Delta x}-r_n)(D_{\Delta x}-r_{n-1})(D_{\Delta x}-r_{n-2})\ldots(D_{\Delta x}-r_3)(D_{\Delta x}-r_2)(D_{\Delta x}-r_1)f(x) = Q(x)$ (4.2-150)

where

$h(r) = r^n + a_{n-1}r^{n-1} + a_{n-2}r^{n-2} + \ldots + a_2r^2 + a_1r + a_0 = (r-r_n)(r-r_{n-1})(r-r_{n-2})\ldots(r-r_3)(r-r_2)(r-r_1) = 0$

$r = r_1, r_2, r_3, \ldots, r_{n-1}, r_n$ , the n roots of h(r)

Substituting Eq 4.2-149 into Eq 4.2-150

$(D_{\Delta x}-r_n)(D_{\Delta x}-r_{n-1})(D_{\Delta x}-r_{n-2})\ldots(D_{\Delta x}-r_3)(D_{\Delta x}-r_2)(D_{\Delta x}-r_1)f(x) = q_2(x)$ (4.2-151)

From Eq 4.2-127 and Eq 4.2-128

$$q_2(x) = [\sum_{p=1}^{m} A_p[x]_{\Delta x}^{p-1}][w\Delta x]^{\frac{x}{\Delta x}}\sin\frac{\pi x}{2\Delta x} + [\sum_{p=1}^{m} B_p[x]_{\Delta x}^{p-1}][w\Delta x]^{\frac{x}{\Delta x}}\cos\frac{\pi x}{2\Delta x} \qquad (4.2\text{-}152)$$

is the general solution to the homogeneous differential difference equation

$$[(D_{\Delta x} - \{-\frac{1}{\Delta x} + jw\})(D_{\Delta x} - \{-\frac{1}{\Delta x} - jw\})]^m q_2(x) = 0 \qquad (4.2\text{-}153)$$

where

$A_p$,$B_p$,u,w = real constants

$m =$ the multiplicity of the $q_2(x)$ complex conjugate root pair, $-\frac{1}{\Delta x} + jw$ and $-\frac{1}{\Delta x} - jw$

$\Delta x$ = x increment

Let several of the h(r) roots, $r_1$ thru $r_n$ be equal to some multiple, $m_h$, of the same complex conjugate root pair, $-\frac{1}{\Delta x} + jw$ and $-\frac{1}{\Delta x} - jw$

Several of the h(r) roots are then

$$[(D_{\Delta x} - \{-\frac{1}{\Delta x} + jw\})(D_{\Delta x} - \{-\frac{1}{\Delta x} - jw\})]^{m_h} \qquad (4.2\text{-}154)$$

Substituting Eq 4.2-154 into Eq 4.2-151

$$[(D_{\Delta x} - \{-\frac{1}{\Delta x} + jw\})(D_{\Delta x} - \{-\frac{1}{\Delta x} - jw\})]^{m_h}(D_{\Delta x}-r_v)(D_{\Delta x}-r_{v-1})(D_{\Delta x}-r_{v-2})\ldots(D_{\Delta x}-r_3)(D_{\Delta x}-r_2) = q_2(x)$$

where (4.2-155)

$n$ = order of the differential difference equation, Eq 4.2-150, and also the number of roots in its related homogeneous differential difference equation

$v = n - 2m_h$

The related roots of $q_2$ are a multiple, m, of the complex conjugate root pair, $-\frac{1}{\Delta x} + jw$ and $-\frac{1}{\Delta x} - jw$.

Take the derivative of both sides of Eq 4.2-155 using the derivative operator,

$$[(D_{\Delta x} - \{-\frac{1}{\Delta x} + jw\})(D_{\Delta x} - \{-\frac{1}{\Delta x} - jw\})]^{m} \qquad (4.2\text{-}156)$$

$$[(D_{\Delta x} - \{-\frac{1}{\Delta x} + jw\})(D_{\Delta x} - \{-\frac{1}{\Delta x} - jw\})]^{m}[(D_{\Delta x} - \{-\frac{1}{\Delta x} + jw\})(D_{\Delta x} - \{-\frac{1}{\Delta x} - jw\})]^{m_h}(D_{\Delta x}-r_v)$$

$$(D_{\Delta x}-r_{v-1})(D_{\Delta x}-r_{v-2})\ldots(D_{\Delta x}-r_3)(D_{\Delta x}-r_2)(D_{\Delta x}-r_1)f(x) = [(D_{\Delta x} - \{-\frac{1}{\Delta x} + jw\})(D_{\Delta x} - \{-\frac{1}{\Delta x} - jw\})]^{m}q_2(x)$$

(4.2-157)

From Eq 4.2-153 and Eq 4.2-157

$$[(D_{\Delta x} - \{-\frac{1}{\Delta x} + jw\})(D_{\Delta x} - \{-\frac{1}{\Delta x} - jw\})]^{m+m_h}(D_{\Delta x}-r_v)(D_{\Delta x}-r_{v-1})(D_{\Delta x}-r_{v-2})\ldots(D_{\Delta x}-r_3)(D_{\Delta x}-r_2)(D_{\Delta x}-r_1)f(x) = 0$$

(4.2-158)

Eq 4.2-155, a non-homgeneous equation, has now been converted into a homogeneous equation, Eq 4.2-158. The function, f(x), the general solution to both equations, remains the same.

The undetermined coefficient solution to Eq 4.2-158 yields f(x) as:

$$f(x) = [\sum_{p=1}^{m+m_h} A_p[x]_{\Delta x}^{p-1}][w\Delta x]^{\frac{x}{\Delta x}}\sin\frac{\pi x}{2\Delta x} + [\sum_{p=1}^{m+m_h} B_p[x]_{\Delta x}^{p-1}][w\Delta x]^{\frac{x}{\Delta x}}\cos\frac{\pi x}{2\Delta x} + \qquad (4.2\text{-}159)$$

$$K_v e^{r_v x} + K_{v-1}e^{r_{v-1}x} + K_{v-2}e^{r_{v-2}x} + \ldots + K_3e^{r_3 x} + K_2e^{r_2 x} + K_1e^{r_1 x}$$

or

$$f(x) = \{ [\sum_{p=m_h+1}^{m+m_h} A_p[x]_{\Delta x}^{p-1}][w\Delta x]^{\frac{x}{\Delta x}}\sin\frac{\pi x}{2\Delta x} + [\sum_{p=m_h+1}^{m+m_h} B_p[x]_{\Delta x}^{p-1}][w\Delta x]^{\frac{x}{\Delta x}}\cos\frac{\pi x}{2\Delta x} \} + \qquad (4.2\text{-}160)$$

$$\{ [\sum_{p=1}^{m_h} A_p[x]_{\Delta x}^{p-1}][w\Delta x]^{\frac{x}{\Delta x}}\sin\frac{\pi x}{2\Delta x} + [\sum_{p=1}^{m_h} B_p[x]_{\Delta x}^{p-1}][w\Delta x]^{\frac{x}{\Delta x}}\cos\frac{\pi x}{2\Delta x} +$$

$$K_v e^{r_v x} + K_{v-1}e^{r_{v-1}x} + K_{v-2}e^{r_{v-2}x} + \ldots + K_3e^{r_3 x} + K_2e^{r_2 x} + K_1e^{r_1 x} \}$$

where

$A_p, B_p$, u,w = real constants

$K_1$ thru $K_v$ = real or complex constants

$m_h$ = the multiplicity of the h(r) characteristic polynomial complex conjugate roots, $-\frac{1}{\Delta x} + jw$ and $-\frac{1}{\Delta x} - jw$, which are equal to the related roots of $q_2(x)$

m = the multiplicity of the $q_2(x)$ complex conjugate roots, $-\frac{1}{\Delta x} + jw$ and $-\frac{1}{\Delta x} - jw$ .

$\Delta x$ = x increment

Identify the particular solution in equation, Eq 4.2-160.

The function f(x), the general solution of Eq 4.2-160, has two components, a complementary solution, $f_C(x)$, and a particular solution, $f_P(x)$.

$$f(x) = f_C(x) + f_P(x) \qquad (4.2\text{-}161)$$

The second term of Eq 4.2-147 is readily recognized as the complementary solution, $f_C(x)$. Then, the first term must be the particular solution, $f_P(x)$.

Thus

**For the differential difference equation, 4.2-79, rewritten below**

$$\mathbf{D_{\Delta x}^{n}f(x) + a_{n-1}D_{\Delta x}^{n-1}f(x) + a_{n-2}D_{\Delta x}^{n-2}f(x) + \ldots + a_1D_{\Delta x}f(x) + a_0f(x) = Q(x),} \qquad (4.2\text{-}162)$$

**where**

$$h(r) = r^n + a_{n-1}r^{n-1} + a_{n-2}\, r^{n-2} + \ldots + a_2r^2 + a_1r + a_0 = (r-r_n)(r-r_{n-1})(r-r_{n-2})\ldots(r-r_3)(r-r_2)(r-r_1)$$

**if** $Q(x) = [\sum_{p=1}^{m} A_p[x]_{\Delta x}^{p-1}][w\Delta x]^{\frac{x}{\Delta x}} \sin\frac{\pi x}{2\Delta x} + [\sum_{p=1}^{m} B_p[x]_{\Delta x}^{p-1}][w\Delta x]^{\frac{x}{\Delta x}} \cos\frac{\pi x}{2\Delta x}$**, the particular solution to Eq 4.2-162 is:** (4.2-163)

$$f_P(x) = [\sum_{p=m_h+1}^{m+m_h} \mathcal{A}_p[x]_{\Delta x}^{p-1}][w\Delta x]^{\frac{x}{\Delta x}} \sin\frac{\pi x}{2\Delta x} + [\sum_{p=m_h+1}^{m+m_h} \mathcal{B}_p[x]_{\Delta x}^{p-1}][w\Delta x]^{\frac{x}{\Delta x}} \cos\frac{\pi x}{2\Delta x} \qquad (4.2\text{-}164)$$

**where**

**n = order of the differential difference equation**

$A_p, B_p, \mathcal{A}_p, \mathcal{B}_p, a_{n-1}\ldots a_0, w$ **= real value constants**

$f_P(x)$ **= differential difference equation particular solution**

$m_h$ **= the multiplicity of the h(r) characteristic polynomial complex conjugate roots,** $-\frac{1}{\Delta x} + jw$ **and** $-\frac{1}{\Delta x} - jw$**, which are equal to the related roots of Q(x)**

**m = the multiplicity of the Q(x) complex conjugate related roots,** $-\frac{1}{\Delta x} + jw$ **and** $-\frac{1}{\Delta x} - jw$

$\Delta x$ **= x increment**

**Note – The particular solution equation, Eq 4.2-164, is more general than Eq 4.2-133.**

The previous equations for obtaining a particular solution to a differential difference equation can be made more understandable and convenient to use with a change in notation.

Let

N = order of the summation polynomials in Q(x)

m = N + 1

$N_h$ = the number of times the complex conjugate roots, $u+jw(1+u\Delta x)$ and $u-jw(1+u\Delta x)$, appear in the characteristic polynomial, h(r)

$N_h = m_h$

Using the above notation and relationships, modify and simplify Eq 4.2-147 and Eq 4.2-148

**For the differential difference equation**

$$D_{\Delta x}{}^n f(x) + a_{n-1}D_{\Delta x}{}^{n-1}f(x) + a_{n-2}D_{\Delta x}{}^{n-2}f(x) + \ldots + a_1D_{\Delta x}f(x) + a_0f(x) = Q(x), \qquad (4.2\text{-}165)$$

**where**

$$h(r) = r^n + a_{n-1}r^{n-1} + a_{n-2}\, r^{n-2} + \ldots + a_2r^2 + a_1r + a_0 = (r-r_n)(r-r_{n-1})(r-r_{n-2})\ldots(r-r_3)(r-r_2)(r-r_1)$$

**if** $Q(x) = [\sum_{p=0}^{N} A_p[x]_{\Delta x}^{p}]\, e_{\Delta x}(u,x) \sin_{\Delta x}(w,x) + [\sum_{p=0}^{N} B_p[x]_{\Delta x}^{p}]\, e_{\Delta x}(u,x) \cos_{\Delta x}(w,x)$**, the particular solution to Eq 4.2-165 is:** (4.2-166)

$$f_P(x) = [\sum_{p=N_h}^{N+N_h} \mathcal{A}_p[x]_{\Delta x}^{p}]\, e_{\Delta x}(u,x) \sin_{\Delta x}(w,x) + [\sum_{p=N_h}^{N+N_h} \mathcal{B}_p[x]_{\Delta x}^{p}]\, e_{\Delta x}(u,x) \cos_{\Delta x}(w,x) \quad (4.2\text{-}167)$$

**where**

**$n$ = order of the differential difference equation**

**$N$ = order of the summation polynomials in $Q(x)$**

**$N_h$ = the number of times the $Q(x)$ complex conjugate related roots, $u+jw(1+u\Delta x)$ and $u-jw(1+u\Delta x)$, appear in the characteristic polynomial, $h(r)$**

**$A_p, B_p, \mathcal{A}_p, \mathcal{B}_p, a_{n-1} \ldots a_0, u, w$ = real value constants**

**$f_P(x)$ = differential difference equation particular solution**

**$\Delta x$ = x increment**

$u \neq -\frac{1}{\Delta x}$

**If $Q(x) = Q_1(x) + Q_2(x) + Q_3(x) + \ldots + Q_g(x)$, then $f_P(x)$ is the sum of the particular solutions of each $Q_r(x)$ where $r = 1,2,3,\ldots,g$.**

**<u>Comment</u> – The form of these equations is, in general, the most convenient to use.**

Let

$N$ = the order of the summation polynomials in $Q(x)$

$m = N + 1$

$N_h$ = the number of times the complex conjugate roots, $-\frac{1}{\Delta x} + jw$ and $-\frac{1}{\Delta x} - jw$ appear in the characteristic polynomial, $h(r)$

$N_h = m_h$

Using the above notation and relationships, modify and simplify Eq 4.2-163 and Eq 4.2-164

**For the differential difference equation**

$$D_{\Delta x}^{n}f(x) + a_{n-1}D_{\Delta x}^{n-1}f(x) + a_{n-2}D_{\Delta x}^{n-2}f(x) + \ldots + a_1 D_{\Delta x}f(x) + a_0 f(x) = Q(x), \quad (4.2\text{-}168)$$

**where**

$$h(r) = r^n + a_{n-1}r^{n-1} + a_{n-2}r^{n-2} + \ldots + a_2 r^2 + a_1 r + a_0 = (r-r_n)(r-r_{n-1})(r-r_{n-2})\ldots(r-r_3)(r-r_2)(r-r_1)$$

**if** $Q(x) = [\sum_{p=0}^{N} A_p[x]_{\Delta x}^{p}][w\Delta x]^{\frac{x}{\Delta x}} \sin\frac{\pi x}{2\Delta x} + [\sum_{p=0}^{N} B_p[x]_{\Delta x}^{p}][w\Delta x]^{\frac{x}{\Delta x}} \cos\frac{\pi x}{2\Delta x}$**, the particular solution to Eq 4.2-168 is:** (4.2-169)

$$f_P(x) = [\sum_{p=N_h}^{N+N_h} \mathcal{A}_p[x]_{\Delta x}^{p}][w\Delta x]^{\frac{x}{\Delta x}} \sin\frac{\pi x}{2\Delta x} + [\sum_{p=N_h}^{N+N_h} \mathcal{B}_p[x]_{\Delta x}^{p}][w\Delta x]^{\frac{x}{\Delta x}} \cos\frac{\pi x}{2\Delta x} \quad (4.2\text{-}170)$$

**where**

**n = order of the differential difference equation**

**N = the order of the summation polynomials in Q(x)**

$N_h$ **= the number of times the Q(x) complex conjugate related roots,** $-\frac{1}{\Delta x} + jw$ **and** $-\frac{1}{\Delta x} - jw$**, appear in the characteristic polynomial, h(r)**

$A_p, B_p, \mathcal{A}_p, \mathcal{B}_p, a_{n-1} \ldots a_0, w$ **= real value constants**

$f_P(x)$ **= differential difference equation particular solution**

$\Delta x$ **= x increment**

**If** $Q(x) = Q_1(x) + Q_2(x) + Q_3(x) + \ldots + Q_g(x)$**, then** $f_P(x)$ **is the sum of the particular solutions of each** $Q_r(x)$ **where r = 1,2,3,…,g.**

**Comment – The form of these equations is, in general, the most convenient to use.**

The differential difference equation relationships previously derived in this section are presented in the following table, Table 4.2-5. This table provides the means to solve a general differential difference equation of the form specified. The solution method shown here is called The Method of Undetermined Coefficients. This method will provide the same results as the $K_{\Delta x}$ Transform Method previously described. The choice of method used is at the descretion of the user. Other methods for solving differential difference equations are available. They will be described later.

**Table 4.2-5 Solving Differential Difference Equations using a Method of Undetermined Coefficients**

| |
|---|
| **f(x) - The Solution of a Differential Difference Equation** |
| For the differential difference equation |
| $D_{\Delta x}^{n}f(x) + a_{n-1}D_{\Delta x}^{n-1}f(x) + a_{n-2}D_{\Delta x}^{n-2}f(x) + \ldots + a_1D_{\Delta x}f(x) + a_0f(x) = Q(x)$ |
| $f(x) = f_C(x) + f_P(x)$, the differential difference equation general solution |
| with the characteristic polynomial, h(r), |

$h(r) = r^n + a_{n-1}r^{n-1} + a_{n-2}\,r^{n-2} + \ldots + a_2r^2 + a_1r + a_0 = (r-r_n)(r-r_{n-1})(r-r_{n-2})\ldots(r-r_3)(r-r_2)(r-r_1)$

**fc(x) - The Complementary Solution of the Differential Difference Equation**

For the related homogeneous differential difference equation

$$D_{\Delta x}^{\ n}f_C(x) + a_{n-1}D_{\Delta x}^{\ n-1}f_C(x) + a_{n-2}D_{\Delta x}^{\ n-2}f_C(x) + \ldots + a_1D_{\Delta x}f_C(x) + a_0f_C(x) = 0$$

where
- $n$ = order of the homogeneous equation
- $a_{n-1,}a_{n-1},\ldots,a_0$ = real constants
- $f_C(x)$ = complementary solution, the general solution to the related homogeneous differential difference equation
- $D_{\Delta x}^{\ n}f(x)$ = nth discrete derivative of the function, f(x)
- $x = x_i + p\Delta x$ ,   p = integers
- $x_i$ = initial value of x
- $\Delta x$ = x increment

with the characteristic polynomial:

$$h(r) = r^n + a_{n-1}r^{n-1} + a_{n-2}\,r^{n-2} + \ldots + a_2r^2 + a_1r + a_0 = (r-r_n)(r-r_{n-1})(r-r_{n-2})\ldots(r-r_3)(r-r_2)(r-r_1) = 0$$

with n roots

$r = r_1, r_2, r_3, \ldots, r_{n-1,} r_n$

group the roots and catagorize them in the following way:

v = one of the different unique root or complex conjugate root pair values
m = multiplicity of each different unique root or complex comjugate root pair value,  m = 1,2,3,…
g = number of groups of different unique root and root pair values  ($g \le n$)

Note – A root of multiplicity 1 is a single root, a root of multiplicity 2 is a double root, etc.

W(v,m,x) = real value general solution function associated with the root(s), v, of multiplicity, m

The W(v,m,x) functions are shown in the table below

$$f_C(x) = \sum_{s=1}^{g} W_s(v_s,m_s,x)$$ , the general solution to the homogeneous differential difference equation

In the following table
- $a,b,A_pB_p$ = real constants
- $[x]_{\Delta x}^{0} = 1$

$$[x]_{\Delta x}^{q} = \sum_{u=1}^{q}(x-[u-1]\Delta x)\ , \ \ q = 1,2,3,\ldots$$

$x = x_i + p\Delta x$ , $p$ =integers

$x_i$ = initial value of x

$\Delta x$ = x increment

The $A_p$, $B_p$ constants in the differential difference equation general solution, f(x), are evaluated from the problem specified initial conditions:

$f(0), f(\Delta x), f(2\Delta x), f(3\Delta x), \ldots, f(n\Delta x)$

or

$f(0), D_{\Delta x}f(0), D_{\Delta x}^{2}f(0), D_{\Delta x}^{3}f(0), \ldots, D_{\Delta x}^{n}f(0)$

The Possible h(r) = 0 Root Outcomes

| v Root/Root Pair Value | W(v,m,x) Function | W(v,m,x) Calculation Function |
|---|---|---|
| a | $[\sum_{p=1}^{m} A_p[x]_{\Delta x}^{p-1}]\ e_{\Delta x}(a,x)$ | $e_{\Delta x}(a,x) = [1+a\Delta x]^{\frac{x}{\Delta x}}$ |
| $0 \pm jb$ | $[\sum_{p=1}^{m} A_p[x]_{\Delta x}^{p-1}]\ \sin_{\Delta x}(b,x)\ +$ $[\sum_{p=1}^{m} B_p[x]_{\Delta x}^{p-1}]\ \cos_{\Delta x}(b,x)$ | $\sin_{\Delta x}(b,x) = \frac{[1+jb\Delta x]^{\frac{x}{\Delta x}} - [1-jb\Delta x]^{\frac{x}{\Delta x}}}{2j}$ $\cos_{\Delta x}(b,x) = \frac{[1+jb\Delta x]^{\frac{x}{\Delta x}} + [1-jb\Delta x]^{\frac{x}{\Delta x}}}{2}$ |
| $-\frac{1}{\Delta x} \pm jb$ | $[\sum_{p=1}^{m} A_p[x]_{\Delta x}^{p-1}][b\Delta x]^{\frac{x}{\Delta x}}\sin\frac{\pi x}{2\Delta x}\ +$ $[\sum_{p=1}^{m} B_p[x]_{\Delta x}^{p-1}][b\Delta x]^{\frac{x}{\Delta x}}\cos\frac{\pi x}{2\Delta x}$ | $[b\Delta x]^{\frac{x}{\Delta x}}\sin\frac{\pi x}{2\Delta x} = \frac{[+jb\Delta x]^{\frac{x}{\Delta x}} - [-jb\Delta x]^{\frac{x}{\Delta x}}}{2j}$ $[b\Delta x]^{\frac{x}{\Delta x}}\cos\frac{\pi x}{2\Delta x} = \frac{[+jb\Delta x]^{\frac{x}{\Delta x}} + [-jb\Delta x]^{\frac{x}{\Delta x}}}{2}$ |
| $a \pm jb$ <br> $a \neq -\frac{1}{\Delta x}$ | $[\sum_{p=1}^{m} A_p[x]_{\Delta x}^{p-1}]\, e_{\Delta x}(a,x)\, \sin_{\Delta x}(\frac{b}{1+a\Delta x}, x)$ + <br> $[\sum_{p=1}^{m} B_p[x]_{\Delta x}^{p-1}]\, e_{\Delta x}(a,x)\, \cos_{\Delta x}(\frac{b}{1+a\Delta x}, x)$ | $e_{\Delta x}(a,x)\, \sin_{\Delta x}(\frac{b}{1+a\Delta x}, x) = \frac{[1+(a+jb)\Delta x]^{\frac{x}{\Delta x}} - [1+(a-jb)\Delta x]^{\frac{x}{\Delta x}}}{2j}$ <br> $e_{\Delta x}(a,x)\, \cos_{\Delta x}(\frac{b}{1+a\Delta x}, x) = \frac{[1+(a+jb)\Delta x]^{\frac{x}{\Delta x}} + [1+(a-jb)\Delta x]^{\frac{x}{\Delta x}}}{2}$ |

**$f_P(x)$ - The Particular Solution of the Differential Difference Equation**

if $Q(x) = [\sum_{p=0}^{N} A_p[x]_{\Delta x}^{p}]\, e_{\Delta x}(u,x)\, \sin_{\Delta x}(w,x) + [\sum_{p=0}^{N} B_p[x]_{\Delta x}^{p}]\, e_{\Delta x}(u,x)\, \cos_{\Delta x}(w,x)$ , the particular solution to this equation is:

$$f_P(x) = [\sum_{p=N_h}^{N+N_h} \mathcal{A}_p[x]_{\Delta x}^{p}]\, e_{\Delta x}(u,x)\, \sin_{\Delta x}(w,x) + [\sum_{p=N_h}^{N+N_h} \mathcal{B}_p[x]_{\Delta x}^{p}]\, e_{\Delta x}(u,x)\, \cos_{\Delta x}(w,x)$$

where

$u \neq -\frac{1}{\Delta x}$

N = the polynomial order of the summations in Q(x)

$N_h$ = the number of times the Q(x) complex conjugate related roots, $u+jw(1+u\Delta x)$ and $u-jw(1+u\Delta x)$, appear in the characteristic polynomial, h(r)

$A_p, B_p, \mathcal{A}_p, \mathcal{B}_p, u, w$ = real value constants

n = order of the differential difference equation

$f_P(x)$ = particular solution to a differential difference equation

$\Delta x$ = x increment

or

if $Q(x) = [\sum_{p=0}^{N} A_p[x]_{\Delta x}^{p}][w\Delta x]^{\frac{x}{\Delta x}} \sin\frac{\pi x}{2\Delta x} + [\sum_{p=0}^{N} B_p[x]_{\Delta x}^{p}][w\Delta x]^{\frac{x}{\Delta x}} \cos\frac{\pi x}{2\Delta x}$ , the particular solution to this equation is:

$$f_P(x) = \left[\sum_{p=N_h}^{N+N_h} \mathcal{A}_p[x]_{\Delta x}^{p}\right][w\Delta x]^{\frac{x}{\Delta x}} \sin\frac{\pi x}{2\Delta x} + \left[\sum_{p=N_h}^{N+N_h} \mathcal{B}_p[x]_{\Delta x}^{p}\right][w\Delta x]^{\frac{x}{\Delta x}} \cos\frac{\pi x}{2\Delta x}$$

where

N = the polynomial order of the summations in Q(x)

$N_h$ = the number of times the Q(x) complex conjugate related roots, $-\frac{1}{\Delta x} + jw$ and $-\frac{1}{\Delta x} - jw$ appear in the characteristic polynomial, h(r)

$A_p, B_p, \mathcal{A}_p, \mathcal{B}_p$, w = real value constants

n = order of the differential difference equation

$f_P(x)$ = particular solution to a differential difference equation

$\Delta x$ = x increment

If $Q(x) = Q_1(x) + Q_2(x) + Q_3(x) + \ldots + Q_g(x)$, then $f_P(x)$ is the sum of the particular solutions of each $Q_r(x)$ where r = 1,2,3,…,g.

The undetermined coefficients of the $f_P(x)$ function are evaluated by introducing the obtained $f_P(x)$ function into the differential difference equation and solving for the coefficient values which will yield an identity.

From Eq 4.2-166, Eq 4.2-167, Eq 4.2-169, Eq 4.2-170, Table 4.2-6 has been developed and is presented below. In the first column are listed functions commonly represented by Q(x) in the general differential difference equation, $D_{\Delta x}{}^{n}f(x) + a_{n-1}D_{\Delta x}{}^{n-1}f(x) + a_{n-2}D_{\Delta x}{}^{n-2}f(x) + \ldots + a_1D_{\Delta x}f(x) + a_0f(x) = Q(x)$. In the second column is listed the undetermined coefficient particular solution function used for each function in Q(x). In the third column is the related root(s) for each function in Q(x).

**Table 4.2-6 Commonly used Q(x) functions and their corresponding Undetermined Coefficient Particular Solution functions**

| # | For these Interval Calculus Functions in Q(x)<br><br>The function specified should be replaced in Q(x) by its identity if its identity is listed.<br><br>$A_p, B_p, k_p$, A, a, b, k, u, w, $\Delta x$ are constants | Put these Interval Calculus functions with undetermined coefficients in the particular solution, $f_P(x)$.<br><br>$N_h$ = the number of times the Q(x) related root(s) appear in the characteristic polynomial, h(r) | Related Root(s) |
|---|---|---|---|
| **1** | $\sum_{p=0}^{N} k_p[x]_{\Delta x}^{p}$ | $\sum_{p=N_h}^{N+N_h} A_p[x]_{\Delta x}^{p}$ | r = 0 |
| **1a** | k | $\sum_{p=N_h}^{N_h} A_p[x]_{\Delta x}^{p}$ | r = 0 |
| **1b** | kx<br><br>$= k[x]_{\Delta x}^{1}$ | $\sum_{p=N_h}^{1+N_h} A_p[x]_{\Delta x}^{p}$ | r = 0 |
| **1c** | $kx^2$<br><br>$= k\Delta x[x]_{\Delta x}^{1} + k[x]_{\Delta x}^{2}$ | $\sum_{p=N_h}^{2+N_h} A_p[x]_{\Delta x}^{p}$ | r = 0 |
| **1d** | $kx^3$<br><br>$= k(\Delta x^2)\,[x]_{\Delta x}^{1} + k(3\Delta x)\,[x]_{\Delta x}^{2} + k[x]_{\Delta x}^{3}$ | $\sum_{p=N_h}^{3+N_h} A_p[x]_{\Delta x}^{p}$ | r = 0 |
| **1e** | $k[x]_{\Delta x}^{n}$<br><br>n = 1,2,3,… | $\sum_{p=N_h}^{n+N_h} A_p[x]_{\Delta x}^{p}$ | r = 0 |
| **2** | $[\sum_{p=0}^{N} k_p[x]_{\Delta x}^{p}]\ e_{\Delta x}(u,x)$ | $[\sum_{p=N_h}^{N+N_h} A_p[x]_{\Delta x}^{p}]\ e_{\Delta x}(u,x)$ | r = u<br><br>$u \neq -\frac{1}{\Delta x}$ |

| # | For these Interval Calculus Functions in Q(x)<br><br>The function specified should be replaced in Q(x) by its identity if its identity is listed.<br><br>$A_p, B_p, k_p$, A, a, b, k, u, w, Δx are constants | Put these Interval Calculus functions with undetermined coefficients in the particular solution, $f_P(x)$.<br><br>$N_h$ = the number of times the Q(x) related root(s) appear in the characteristic polynomial, h(r) | Related Root(s) |
|---|---|---|---|
| **2a** | $ke_{\Delta x}(u,x)$ | $[\sum_{p=N_h}^{N_h} A_p[x]_{\Delta x}^{p}]\ e_{\Delta x}(u,x)$ | r = u<br>$u \neq -\frac{1}{\Delta x}$ |
| **2b** | $kxe_{\Delta x}(u,x)$ | $[\sum_{p=N_h}^{1+N_h} A_p[x]_{\Delta x}^{p}]\ e_{\Delta x}(u,x)$ | r = u<br>$u \neq -\frac{1}{\Delta x}$ |
| **2c** | $k[x]_{\Delta x}^{n} e_{\Delta x}(u,x)$ | $[\sum_{p=N_h}^{n+N_h} A_p[x]_{\Delta x}^{p}]\ e_{\Delta x}(u,x)$ | r = u<br>$u \neq -\frac{1}{\Delta x}$ |
| **3** | $[\sum_{p=0}^{N} k_p[x]_{\Delta x}^{p}]\ \sin_{\Delta x}(w,x)$ | $[\sum_{p=N_h}^{N+N_h} A_p[x]_{\Delta x}^{p}]\ \sin_{\Delta x}(w,x)\ +$<br>$[\sum_{p=N_h}^{N+N_h} B_p[x]_{\Delta x}^{p}]\ \cos_{\Delta x}(w,x)$ | r = ±jw |
| **3a** | $k\sin_{\Delta x}(w,x)$ | $[\sum_{p=N_h}^{N_h} A_p[x]_{\Delta x}^{p}]\ \sin_{\Delta x}(w,x)\ +$<br>$[\sum_{p=N_h}^{N_h} B_p[x]_{\Delta x}^{p}]\ \cos_{\Delta x}(w,x)$ | r = ±jw |
| **3b** | $kx\sin_{\Delta x}(w,x)$ | $[\sum_{p=N_h}^{1+N_h} A_p[x]_{\Delta x}^{p}]\ \sin_{\Delta x}(w,x)\ +$ | r = ±jw |

| # | For these Interval Calculus Functions in Q(x)<br><br>The function specified should be replaced in Q(x) by its identity if its identity is listed.<br><br>$A_p, B_p, k_p$, A, a, b, k, u, w, Δx are constants | Put these Interval Calculus functions with undetermined coefficients in the particular solution, $f_P(x)$.<br><br>$N_h$ = the number of times the Q(x) related root(s) appear in the characteristic polynomial, h(r) | Related Root(s) |
|---|---|---|---|
| | | $[\sum_{p=N_h}^{1+N_h} B_p [x]_{\Delta x}^{p}] \cos_{\Delta x}(w,x)$ | |
| **3c** | $k[x]_{\Delta x}^{n} \sin_{\Delta x}(w,x)$ | $[\sum_{p=N_h}^{n+N_h} A_p [x]_{\Delta x}^{p}] \sin_{\Delta x}(w,x)$ +<br>$[\sum_{p=N_h}^{n+N_h} B_p [x]_{\Delta x}^{p}] \cos_{\Delta x}(w,x)$ | $r = \pm jw$ |
| **4** | $[\sum_{p=0}^{N} k_p [x]_{\Delta x}^{p}] \cos_{\Delta x}(w,x)$ | $[\sum_{p=N_h}^{N+N_h} A_p [x]_{\Delta x}^{p}] \sin_{\Delta x}(w,x)$ +<br>$[\sum_{p=N_h}^{N+N_h} B_p [x]_{\Delta x}^{p}] \cos_{\Delta x}(w,x)$ | $r = \pm jw$ |
| **4a** | $k\cos_{\Delta x}(w,x)$ | $[\sum_{p=N_h}^{N_h} A_p [x]_{\Delta x}^{p}] \sin_{\Delta x}(w,x)$ +<br>$[\sum_{p=N_h}^{N_h} B_p [x]_{\Delta x}^{p}] \cos_{\Delta x}(w,x)$ | $r = \pm jw$ |
| **4b** | $kx\cos_{\Delta x}(w,x)$ | $[\sum_{p=N_h}^{1+N_h} A_p [x]_{\Delta x}^{p}] \sin_{\Delta x}(w,x)$ + | $r = \pm jw$ |

| # | **For these Interval Calculus Functions in Q(x)**<br><br>**The function specified should be replaced in Q(x) by its identity if its identity is listed.**<br>**$A_p$, $B_p$, $k_p$, A, a, b, k, u, w, Δx are constants** | **Put these Interval Calculus functions with undetermined coefficients in the particular solution, $f_P(x)$.**<br><br>**$N_h$ = the number of times the Q(x) related root(s) appear in the characteristic polynomial, h(r)** | **Related Root(s)** |
|---|---|---|---|
| | | $[\sum_{p=N_h}^{1+N_h} B_p[x]_{\Delta x}^{p}]\ \cos_{\Delta x}(w,x)$ | |
| **4c** | $k[x]_{\Delta x}^{n}\ \cos_{\Delta x}(w,x)$ | $[\sum_{p=N_h}^{n+N_h} A_p[x]_{\Delta x}^{p}]\ \sin_{\Delta x}(w,x)\ +$<br>$[\sum_{p=N_h}^{n+N_h} B_p[x]_{\Delta x}^{p}]\ \cos_{\Delta x}(w,x)$ | $r = \pm jw$ |
| **5** | $[\sum_{p=0}^{N} k_p[x]_{\Delta x}^{p}]\ e_{\Delta x}(u,x)\ \sin_{\Delta x}(w,x)$ | $[\sum_{p=N_h}^{N+N_h} A_p[x]_{\Delta x}^{p}]\ e_{\Delta x}(u,x)\ \sin_{\Delta x}(w,x)\ +$<br>$[\sum_{p=N_h}^{N+N_h} B_p[x]_{\Delta x}^{p}]\ e_{\Delta x}(u,x)\ \cos_{\Delta x}(w,x)$ | $r = u \pm jw(1+u\Delta x)$<br>$u \neq -\frac{1}{\Delta x}$ |
| **5a** | $ke_{\Delta x}(u,x)\ \sin_{\Delta x}(w,x)$ | $[\sum_{p=N_h}^{N_h} A_p[x]_{\Delta x}^{p}]\ e_{\Delta x}(u,x)\ \sin_{\Delta x}(w,x)\ +$<br>$[\sum_{p=N_h}^{N_h} B_p[x]_{\Delta x}^{p}]\ e_{\Delta x}(u,x)\ \cos_{\Delta x}(w,x)$ | $r = u \pm jw(1+u\Delta x)$<br>$u \neq -\frac{1}{\Delta x}$ |
| **5b** | $kxe_{\Delta x}(u,x)\ \sin_{\Delta x}(w,x)$ | $[\sum_{p=N_h}^{1+N_h} A_p[x]_{\Delta x}^{p}]\ e_{\Delta x}(u,x)\ \sin_{\Delta x}(w,x)\ +$ | $r = u \pm jw(1+u\Delta x)$<br>$u \neq -\frac{1}{\Delta x}$ |

| # | For these Interval Calculus Functions in Q(x)<br><br>The function specified should be replaced in Q(x) by its identity if its identity is listed.<br><br>$A_p$, $B_p$, $k_p$, A, a, b, k, u, w, Δx are constants | Put these Interval Calculus functions with undetermined coefficients in the particular solution, $f_P(x)$.<br><br>$N_h$ = the number of times the Q(x) related root(s) appear in the characteristic polynomial, h(r) | Related Root(s) |
|---|---|---|---|
| | | $[\sum_{p=N_h}^{1+N_h} B_p[x]_{\Delta x}^{p}]\ e_{\Delta x}(u,x)\ \cos_{\Delta x}(w,x)$ | |
| **5c** | $k[x]_{\Delta x}^{n}\ e_{\Delta x}(u,x)\ \sin_{\Delta x}(w,x)$ | $[\sum_{p=N_h}^{n+N_h} A_p[x]_{\Delta x}^{p}]\ e_{\Delta x}(u,x)\ \sin_{\Delta x}(w,x)\ +$<br>$[\sum_{p=N_h}^{n+N_h} B_p[x]_{\Delta x}^{p}]\ e_{\Delta x}(u,x)\ \cos_{\Delta x}(w,x)$ | $r = u \pm jw(1+u\Delta x)$<br>$u \neq -\frac{1}{\Delta x}$ |
| **6** | $[\sum_{p=0}^{N} k_p[x]_{\Delta x}^{p}]\ e_{\Delta x}(u,x)\ \cos_{\Delta x}(w,x)$ | $[\sum_{p=N_h}^{N+N_h} A_p[x]_{\Delta x}^{p}]\ e_{\Delta x}(u,x)\ \sin_{\Delta x}(w,x)\ +$<br>$[\sum_{p=N_h}^{N+N_h} B_p[x]_{\Delta x}^{p}]\ e_{\Delta x}(u,x)\ \cos_{\Delta x}(w,x)$ | $r = u \pm jw(1+u\Delta x)$<br>$u \neq -\frac{1}{\Delta x}$ |
| **6a** | $k e_{\Delta x}(u,x)\ \cos_{\Delta x}(w,x)$ | $[\sum_{p=N_h}^{N_h} A_p[x]_{\Delta x}^{p}]\ e_{\Delta x}(u,x)\ \sin_{\Delta x}(w,x)\ +$<br>$[\sum_{p=N_h}^{N_h} B_p[x]_{\Delta x}^{p}]\ e_{\Delta x}(u,x)\ \cos_{\Delta x}(w,x)$ | $r = u \pm jw(1+u\Delta x)$<br>$u \neq -\frac{1}{\Delta x}$ |
| **6b** | $kx e_{\Delta x}(u,x)\ \cos_{\Delta x}(w,x)$ | $[\sum_{p=N_h}^{1+N_h} A_p[x]_{\Delta x}^{p}]\ e_{\Delta x}(u,x)\ \sin_{\Delta x}(w,x)\ +$ | $r = u \pm jw(1+u\Delta x)$<br>$u \neq -\frac{1}{\Delta x}$ |

| # | For these Interval Calculus Functions in Q(x)<br><br>The function specified should be replaced in Q(x) by its identity if its identity is listed.<br><br>$A_p, B_p, k_p$, A, a, b, k, u, w, Δx are constants | Put these Interval Calculus functions with undetermined coefficients in the particular solution, $f_P(x)$.<br><br>$N_h$ = the number of times the Q(x) related root(s) appear in the characteristic polynomial, h(r) | Related Root(s) |
|---|---|---|---|
| | | $[\sum_{p=N_h}^{1+N_h} B_p[x]^{p}_{\Delta x}]\ e_{\Delta x}(u,x)\ \cos_{\Delta x}(w,x)$ | |
| **6c** | $k[x]^{n}_{\Delta x}\ e_{\Delta x}(u,x)\ \cos_{\Delta x}(w,x)$ | $[\sum_{p=N_h}^{n+N_h} A_p[x]^{p}_{\Delta x}]\ e_{\Delta x}(u,x)\ \sin_{\Delta x}(w,x)\ +$<br>$[\sum_{p=N_h}^{n+N_h} B_p[x]^{p}_{\Delta x}]\ e_{\Delta x}(u,x)\ \cos_{\Delta x}(w,x)$ | $r = u \pm jw(1+u\Delta x)$<br>$u \neq -\frac{1}{\Delta x}$ |
| **7** | $[\sum_{p=0}^{N} k_p[x]^{p}_{\Delta x}][w\Delta x]^{\frac{x}{\Delta x}} \sin\frac{\pi x}{2\Delta x}$ | $[\sum_{p=N_h}^{N+N_h} A_p[x]^{p}_{\Delta x}][w\Delta x]^{\frac{x}{\Delta x}} \sin\frac{\pi x}{2\Delta x}\ +$<br>$[\sum_{p=N_h}^{N+N_h} B_p[x]^{p}_{\Delta x}][w\Delta x]^{\frac{x}{\Delta x}} \cos\frac{\pi x}{2\Delta x}$ | $r = -\frac{1}{\Delta x} \pm jw$ |
| **7a** | $k[w\Delta x]^{\frac{x}{\Delta x}} \sin\frac{\pi x}{2\Delta x}$ | $[\sum_{p=N_h}^{N_h} A_p[x]^{p}_{\Delta x}][w\Delta x]^{\frac{x}{\Delta x}} \sin\frac{\pi x}{2\Delta x}\ +$<br>$[\sum_{p=N_h}^{N_h} B_p[x]^{p}_{\Delta x}][w\Delta x]^{\frac{x}{\Delta x}} \cos\frac{\pi x}{2\Delta x}$ | $r = -\frac{1}{\Delta x} \pm jw$ |
| **7b** | $kx[w\Delta x]^{\frac{x}{\Delta x}} \sin\frac{\pi x}{2\Delta x}$ | $[\sum_{p=N_h}^{1+N_h} A_p[x]^{p}_{\Delta x}][w\Delta x]^{\frac{x}{\Delta x}} \sin\frac{\pi x}{2\Delta x}\ +$ | $r = -\frac{1}{\Delta x} \pm jw$ |

| # | For these Interval Calculus Functions in Q(x)<br><br>The function specified should be replaced in Q(x) by its identity if its identity is listed.<br><br>$A_p, B_p, k_p$, A, a, b, k, u, w, $\Delta x$ are constants | Put these Interval Calculus functions with undetermined coefficients in the particular solution, $f_P(x)$.<br><br>$N_h$ = the number of times the Q(x) related root(s) appear in the characteristic polynomial, h(r) | Related Root(s) |
|---|---|---|---|
| | | $[\sum_{p=N_h}^{1+N_h} B_p[x]_{\Delta x}^{p}][w\Delta x]^{\frac{x}{\Delta x}}\cos\frac{\pi x}{2\Delta x}$ | |
| **7c** | $k[x]_{\Delta x}^{n}[w\Delta x]^{\frac{x}{\Delta x}}\sin\frac{\pi x}{2\Delta x}$ | $[\sum_{p=N_h}^{n+N_h} A_p[x]_{\Delta x}^{p}][w\Delta x]^{\frac{x}{\Delta x}}\sin\frac{\pi x}{2\Delta x}$ +<br>$[\sum_{p=N_h}^{n+N_h} B_p[x]_{\Delta x}^{p}][w\Delta x]^{\frac{x}{\Delta x}}\cos\frac{\pi x}{2\Delta x}$ | $r = -\frac{1}{\Delta x} \pm jw$ |
| **8** | $[\sum_{p=0}^{N} k_p[x]_{\Delta x}^{p}][w\Delta x]^{\frac{x}{\Delta x}}\cos\frac{\pi x}{2\Delta x}$ | $[\sum_{p=N_h}^{N+N_h} A_p[x]_{\Delta x}^{p}][w\Delta x]^{\frac{x}{\Delta x}}\sin\frac{\pi x}{2\Delta x}$ +<br>$[\sum_{p=N_h}^{N+N_h} B_p[x]_{\Delta x}^{p}][w\Delta x]^{\frac{x}{\Delta x}}\cos\frac{\pi x}{2\Delta x}$ | $r = -\frac{1}{\Delta x} \pm jw$ |
| **8a** | $k[w\Delta x]^{\frac{x}{\Delta x}}\cos\frac{\pi x}{2\Delta x}$ | $[\sum_{p=N_h}^{N_h} A_p[x]_{\Delta x}^{p}][w\Delta x]^{\frac{x}{\Delta x}}\sin\frac{\pi x}{2\Delta x}$ +<br>$[\sum_{p=N_h}^{N_h} B_p[x]_{\Delta x}^{p}][w\Delta x]^{\frac{x}{\Delta x}}\cos\frac{\pi x}{2\Delta x}$ | $r = -\frac{1}{\Delta x} \pm jw$ |
| **8b** | $kx[w\Delta x]^{\frac{x}{\Delta x}}\cos\frac{\pi x}{2\Delta x}$ | $[\sum_{p=N_h}^{1+N_h} A_p[x]_{\Delta x}^{p}][w\Delta x]^{\frac{x}{\Delta x}}\sin\frac{\pi x}{2\Delta x}$ + | $r = -\frac{1}{\Delta x} \pm jw$ |

| # | For these Interval Calculus Functions in Q(x)<br><br>The function specified should be replaced in Q(x) by its identity if its identity is listed.<br><br>$A_p, B_p, k_p$, A, a, b, k, u, w, Δx are constants | Put these Interval Calculus functions with undetermined coefficients in the particular solution, $f_P(x)$.<br><br>$N_h$ = the number of times the Q(x) related root(s) appear in the characteristic polynomial, h(r) | Related Root(s) |
|---|---|---|---|
| | | $[\sum_{p=N_h}^{1+N_h} B_p[x]_{\Delta x}^{p}][w\Delta x]^{\frac{x}{\Delta x}} \cos\frac{\pi x}{2\Delta x}$ | |
| **8c** | $k[x]_{\Delta x}^{n}[w\Delta x]^{\frac{x}{\Delta x}}\cos\frac{\pi x}{2\Delta x}$ | $[\sum_{p=N_h}^{n+N_h} A_p[x]_{\Delta x}^{p}][w\Delta x]^{\frac{x}{\Delta x}} \sin\frac{\pi x}{2\Delta x}$ +<br>$[\sum_{p=N_h}^{n+N_h} B_p[x]_{\Delta x}^{p}][w\Delta x]^{\frac{x}{\Delta x}} \cos\frac{\pi x}{2\Delta x}$ | $r = -\frac{1}{\Delta x} \pm jw$ |
| **9** | $[\sum_{p=0}^{N} k_p[x]_{\Delta x}^{p}] \sinh_{\Delta x}(w,x)$ | $[\sum_{p=N_h}^{N+N_h} A_p[x]_{\Delta x}^{p}] \sinh_{\Delta x}(w,x)$ +<br>$[\sum_{p=N_h}^{N+N_h} B_p[x]_{\Delta x}^{p}] \cosh_{\Delta x}(w,x)$ | $r = \pm w$ |
| **9a** | $k\sinh_{\Delta x}(w,x)$ | $[\sum_{p=N_h}^{N_h} A_p[x]_{\Delta x}^{p}] \sinh_{\Delta x}(w,x)$ +<br>$[\sum_{p=N_h}^{N_h} B_p[x]_{\Delta x}^{p}] \cosh_{\Delta x}(w,x)$ | $r = \pm w$ |
| **9b** | $kx\sinh_{\Delta x}(w,x)$ | $[\sum_{p=N_h}^{1+N_h} A_p[x]_{\Delta x}^{p}] \sinh_{\Delta x}(w,x)$ + | $r = \pm w$ |

| # | For these Interval Calculus Functions in Q(x)<br><br>The function specified should be replaced in Q(x) by its identity if its identity is listed.<br><br>$A_p, B_p, k_p$, A, a, b, k, u, w, $\Delta x$ are constants | Put these Interval Calculus functions with undetermined coefficients in the particular solution, $f_P(x)$.<br><br>$N_h$ = the number of times the Q(x) related root(s) appear in the characteristic polynomial, h(r) | Related Root(s) |
|---|---|---|---|
| | | $[\sum_{p=N_h}^{n+N_h} B_p[x]_{\Delta x}^{p}]\cos_{\Delta x}(w,x)$ | |
| **9c** | $k[x]_{\Delta x}^{n}\sinh_{\Delta x}(w,x)$ | $[\sum_{p=N_h}^{n+N_h} A_p[x]_{\Delta x}^{p}]\sinh_{\Delta x}(w,x)$ +<br>$[\sum_{p=N_h}^{n+N_h} B_p[x]_{\Delta x}^{p}]\cosh_{\Delta x}(w,x)$ | $r = \pm w$ |
| **10** | $[\sum_{p=0}^{N} k_p[x]_{\Delta x}^{p}]\cosh_{\Delta x}(w,x)$ | $[\sum_{p=N_h}^{N+N_h} A_p[x]_{\Delta x}^{p}]\sinh_{\Delta x}(w,x)$ +<br>$[\sum_{p=N_h}^{N+N_h} B_p[x]_{\Delta x}^{p}]\cosh_{\Delta x}(w,x)$ | $r = \pm w$ |
| **10a** | $k\cosh_{\Delta x}(w,x)$ | $[\sum_{p=N_h}^{N_h} A_p[x]_{\Delta x}^{p}]\sinh_{\Delta x}(w,x)$ +<br>$[\sum_{p=N_h}^{N_h} B_p[x]_{\Delta x}^{p}]\cosh_{\Delta x}(w,x)$ | $r = \pm w$ |
| **10b** | $kx\text{coshh}_{\Delta x}(w,x)$ | $[\sum_{p=N_h}^{1+N_h} A_p[x]_{\Delta x}^{p}]\sinh_{\Delta x}(w,x)$ + | $r = \pm w$ |

| # | For these Interval Calculus Functions in Q(x)<br><br>The function specified should be replaced in Q(x) by its identity if its identity is listed.<br><br>$A_p, B_p, k_p$, A, a, b, k, u, w, $\Delta x$ are constants | Put these Interval Calculus functions with undetermined coefficients in the particular solution, $f_P(x)$.<br><br>$N_h$ = the number of times the Q(x) related root(s) appear in the characteristic polynomial, h(r) | Related Root(s) |
|---|---|---|---|
| | | $[\sum_{p=N_h}^{n+N_h} B_p[x]_{\Delta x}^{p}] \cosh_{\Delta x}(w,x)$ | |
| **10c** | $k[x]_{\Delta x}^{n} \cosh_{\Delta x}(w,x)$ | $[\sum_{p=N_h}^{n+N_h} A_p[x]_{\Delta x}^{p}] \sinh_{\Delta x}(w,x) +$<br>$[\sum_{p=N_h}^{n+N_h} B_p[x]_{\Delta x}^{p}] \cosh_{\Delta x}(w,x)$ | $r = \pm w$ |
| **11** | $[\sum_{p=0}^{N} k_p[x]_{\Delta x}^{p}] A^x$<br><br>$= [\sum_{p=0}^{N} k_p[x]_{\Delta x}^{p}] e_{\Delta x}(\frac{A^{\Delta x}-1}{\Delta x}, x)$ | $[\sum_{p=N_h}^{N+N_h} A_p[x]_{\Delta x}^{p}] e_{\Delta x}(\frac{A^{\Delta x}-1}{\Delta x}, x)$ | $r = \frac{A^{\Delta x}-1}{\Delta x}$ |
| **11a** | $kA^x$<br><br>$= ke_{\Delta x}(\frac{A^{\Delta x}-1}{\Delta x}, x)$ | $[\sum_{p=N_h}^{N_h} A_p[x]_{\Delta x}^{p}] e_{\Delta x}(\frac{A^{\Delta x}-1}{\Delta x}, x)$ | $r = \frac{A^{\Delta x}-1}{\Delta x}$ |
| **11b** | $kxA^x$<br><br>$= kxe_{\Delta x}(\frac{A^{\Delta x}-1}{\Delta x}, x)$ | $[\sum_{p=N_h}^{1+N_h} A_p[x]_{\Delta x}^{p}] e_{\Delta x}(\frac{A^{\Delta x}-1}{\Delta x}, x)$ | $r = \frac{A^{\Delta x}-1}{\Delta x}$ |
| **12** | $[\sum_{p=0}^{N} k_p[x]_{\Delta x}^{p}] e^{ax}$ | $[\sum_{p=N_h}^{N+N_h} A_p[x]_{\Delta x}^{p}] e_{\Delta x}(\frac{e^{a\Delta x}-1}{\Delta x}, x)$ | $r = \frac{e^{a\Delta x}-1}{\Delta x}$ |

| # | For these Interval Calculus Functions in Q(x)<br><br>The function specified should be replaced in Q(x) by its identity if its identity is listed.<br><br>$A_p, B_p, k_p$, A, a, b, k, u, w, Δx are constants | Put these Interval Calculus functions with undetermined coefficients in the particular solution, $f_P(x)$.<br><br>$N_h$ = the number of times the Q(x) related root(s) appear in the characteristic polynomial, h(r) | Related Root(s) |
|---|---|---|---|
| | $= [\sum_{p=0}^{N} k_p[x]_{\Delta x}^{p}]\, e_{\Delta x}(\frac{e^{a\Delta x}-1}{\Delta x}, x)$ | | |
| **12a** | $ke^{ax}$<br><br>$= ke_{\Delta x}(\frac{e^{a\Delta x}-1}{\Delta x}, x)$ | $[\sum_{p=N_h}^{N_h} A_p[x]_{\Delta x}^{p}]\, e_{\Delta x}(\frac{e^{a\Delta x}-1}{\Delta x}, x)$ | $r = \frac{e^{a\Delta x}-1}{\Delta x}$ |
| **12b** | $kxe^{ax}$<br><br>$= kxe_{\Delta x}(\frac{e^{a\Delta x}-1}{\Delta x}, x)$ | $[\sum_{p=N_h}^{1+N_h} A_p[x]_{\Delta x}^{p}]\, e_{\Delta x}(\frac{e^{a\Delta x}-1}{\Delta x}, x)$ | $r = \frac{e^{a\Delta x}-1}{\Delta x}$ |

| # | For these Interval Calculus Functions in Q(x)<br><br>The function specified should be replaced in Q(x) by its identity if its identity is listed.<br><br>k = constant | Put these Interval Calculus functions with undetermined coefficients in the particular solution, $f_P(x)$.<br><br>$N_h$ = the number of times the Q(x) related root(s) appear in the characteristic polynomial, h(r) | Related Root(s) |
|---|---|---|---|
| **13** | $[\sum_{p=0}^{N} k_p[x]_{\Delta x}^{p}]\ \sin bx$<br><br>$= [\sum_{p=0}^{N} k_p[x]_{\Delta x}^{p}]\ e_{\Delta x}(\frac{\cos b\Delta x - 1}{\Delta x}, x)\ \sin_{\Delta x}(\frac{\tan b\Delta x}{\Delta x}, x)$<br><br>where<br><br>$\frac{x}{\Delta x}$ = integer<br><br>$\cos b\Delta x \neq 0$ | $[\sum_{p=N_h}^{N+N_h} A_p[x]_{\Delta x}^{p}]\ e_{\Delta x}(\frac{\cos b\Delta x-1}{\Delta x}, x)\ \sin_{\Delta x}(\frac{\tan b\Delta x}{\Delta x}, x) +$<br><br>$[\sum_{p=N_h}^{N+N_h} B_p[x]_{\Delta x}^{p}]\ e_{\Delta x}(\frac{\cos b\Delta x-1}{\Delta x}, x)\ \cos_{\Delta x}(\frac{\tan b\Delta x}{\Delta x}, x)$ | $r = \frac{\cos b\Delta x - 1}{\Delta x} \pm j\,\frac{\sin b\Delta x}{\Delta x}$<br><br>$\cos b\Delta x \neq 0$ |
| **13a** | ksinbx<br><br>$= k e_{\Delta x}(\frac{\cos b\Delta x - 1}{\Delta x}, x)\ \sin_{\Delta x}(\frac{\tan b\Delta x}{\Delta x}, x)$<br><br>where<br><br>$\frac{x}{\Delta x}$ = integer<br><br>$\cos b\Delta x \neq 0$ | $[\sum_{p=N_h}^{N_h} A_p[x]_{\Delta x}^{p}]\ e_{\Delta x}(\frac{\cos b\Delta x-1}{\Delta x}, x)\ \sin_{\Delta x}(\frac{\tan b\Delta x}{\Delta x}, x) +$<br><br>$[\sum_{p=N_h}^{N_h} B_p[x]_{\Delta x}^{p}]\ e_{\Delta x}(\frac{\cos b\Delta x-1}{\Delta x}, x)\ \cos_{\Delta x}(\frac{\tan b\Delta x}{\Delta x}, x)$ | $r = \frac{\cos b\Delta x - 1}{\Delta x} \pm j\,\frac{\sin b\Delta x}{\Delta x}$<br><br>$\cos b\Delta x \neq 0$ |

| # | For these Interval Calculus Functions in Q(x)<br><br>The function specified should be replaced in Q(x) by its identity if its identity is listed.<br><br>k = constant | Put these Interval Calculus functions with undetermined coefficients in the particular solution, $f_P(x)$.<br><br>$N_h$ = the number of times the Q(x) related root(s) appear in the characteristic polynomial, h(r) | Related Root(s) |
|---|---|---|---|
| **13b** | kxsinbx<br><br>$= kxe_{\Delta x}(\frac{\cos b\Delta x - 1}{\Delta x}, x)\ \sin_{\Delta x}(\frac{\tan b\Delta x}{\Delta x}, x)$<br><br>where<br><br>$\frac{x}{\Delta x}$ = integer<br><br>$\cos b\Delta x \neq 0$ | $[\sum_{p=N_h}^{1+N_h} A_p[x]_{\Delta x}^{p}]\ e_{\Delta x}(\frac{\cos b\Delta x-1}{\Delta x}, x)\ \sin_{\Delta x}(\frac{\tan b\Delta x}{\Delta x}, x) +$<br><br>$[\sum_{p=N_h}^{1+N_h} B_p[x]_{\Delta x}^{p}]\ e_{\Delta x}(\frac{\cos b\Delta x-1}{\Delta x}, x)\ \cos_{\Delta x}(\frac{\tan b\Delta x}{\Delta x}, x)$ | $r = \frac{\cos b\Delta x - 1}{\Delta x} \pm j\ \frac{\sin b\Delta x}{\Delta x}$<br><br>$\cos b\Delta x \neq 0$ |
| **14** | $[\sum_{p=0}^{N} k_p[x]_{\Delta x}^{p}]\ \cos bx$<br><br>$= [\sum_{p=0}^{N} k_p[x]_{\Delta x}^{p}]\ e_{\Delta x}(\frac{\cos b\Delta x - 1}{\Delta x}, x)\ \cos_{\Delta x}(\frac{\tan b\Delta x}{\Delta x}, x)$<br><br>where<br><br>$\frac{x}{\Delta x}$ = integer<br><br>$\cos b\Delta x \neq 0$ | $[\sum_{p=N_h}^{N+N_h} A_p[x]_{\Delta x}^{p}]\ e_{\Delta x}(\frac{\cos b\Delta x-1}{\Delta x}, x)\ \sin_{\Delta x}(\frac{\tan b\Delta x}{\Delta x}, x) +$<br><br>$[\sum_{p=N_h}^{N+N_h} B_p[x]_{\Delta x}^{p}]\ e_{\Delta x}(\frac{\cos b\Delta x-1}{\Delta x}, x)\ \cos_{\Delta x}(\frac{\tan b\Delta x}{\Delta x}, x)$ | $r = \frac{\cos b\Delta x - 1}{\Delta x} \pm j\ \frac{\sin b\Delta x}{\Delta x}$<br><br>$\cos b\Delta x \neq 0$ |

| # | For these Interval Calculus Functions in Q(x)<br><br>The function specified should be replaced in Q(x) by its identity if its identity is listed.<br>k = constant | Put these Interval Calculus functions with undetermined coefficients in the particular solution, $f_P(x)$.<br>$N_h$ = the number of times the Q(x) related root(s) appear in the characteristic polynomial, h(r) | Related Root(s) |
|---|---|---|---|
| **14a** | kcosbx<br>$= ke_{\Delta x}(\frac{\cos b\Delta x - 1}{\Delta x}, x)\cos_{\Delta x}(\frac{\tan b\Delta x}{\Delta x}, x)$<br>where<br>$\frac{x}{\Delta x}$ = integer<br>$\cos b\Delta x \neq 0$ | $[\sum_{p=N_h}^{N_h} A_p[x]_{\Delta x}^{p}]\, e_{\Delta x}(\frac{\cos b\Delta x-1}{\Delta x}, x)\sin_{\Delta x}(\frac{\tan b\Delta x}{\Delta x}, x) +$<br>$[\sum_{p=N_h}^{N_h} B_p[x]_{\Delta x}^{p}]\, e_{\Delta x}(\frac{\cos b\Delta x-1}{\Delta x}, x)\cos_{\Delta x}(\frac{\tan b\Delta x}{\Delta x}, x)$ | $r = \frac{\cos b\Delta x - 1}{\Delta x} \pm j\frac{\sin b\Delta x}{\Delta x}$<br>$\cos b\Delta x \neq 0$ |
| **14b** | kxcosbx<br>$= kxe_{\Delta x}(\frac{\cos b\Delta x - 1}{\Delta x}, x)\cos_{\Delta x}(\frac{\tan b\Delta x}{\Delta x}, x)$<br>where<br>$\frac{x}{\Delta x}$ = integer<br>$\cos b\Delta x \neq 0$ | $[\sum_{p=N_h}^{1+N_h} A_p[x]_{\Delta x}^{p}]\, e_{\Delta x}(\frac{\cos b\Delta x-1}{\Delta x}, x)\sin_{\Delta x}(\frac{\tan b\Delta x}{\Delta x}, x) +$<br>$[\sum_{p=N_h}^{1+N_h} B_p[x]_{\Delta x}^{p}]\, e_{\Delta x}(\frac{\cos b\Delta x-1}{\Delta x}, x)\cos_{\Delta x}(\frac{\tan b\Delta x}{\Delta x}, x)$ | $r = \frac{\cos b\Delta x - 1}{\Delta x} \pm j\frac{\sin b\Delta x}{\Delta x}$<br>$\cos b\Delta x \neq 0$ |

| # | For these Interval Calculus Functions in Q(x)<br><br>The function specified should be replaced in Q(x) by its identity if its identity is listed.<br><br>k = constant | Put these Interval Calculus functions with undetermined coefficients in the particular solution, $f_P(x)$.<br><br>$N_h$ = the number of times the Q(x) related root(s) appear in the characteristic polynomial, h(r) | Related Root(s) |
|---|---|---|---|
| **15** | $[\sum_{p=0}^{N} k_p[x]_{\Delta x}^{p}]\ e^{ax}\sin bx$<br><br>$= [\sum_{p=0}^{N} k_p[x]_{\Delta x}^{p}]\ e_{\Delta x}(\frac{e^{a\Delta x}\cos b\Delta x - 1}{\Delta x}, x)\ \sin_{\Delta x}(\frac{\tan b\Delta x}{\Delta x}, x)$<br><br>where<br><br>$\frac{x}{\Delta x}$ = integer<br><br>$\cos b\Delta x \neq 0$ | $[\sum_{p=N_h}^{N+N_h} A_p[x]_{\Delta x}^{p}]\ e_{\Delta x}(\frac{e^{a\Delta x}\cos b\Delta x-1}{\Delta x}, x)\ \sin_{\Delta x}(\frac{\tan b\Delta x}{\Delta x}, x) +$<br><br>$[\sum_{p=N_h}^{N+N_h} B_p[x]_{\Delta x}^{p}]\ e_{\Delta x}(\frac{e^{a\Delta x}\cos b\Delta x-1}{\Delta x}, x)\ \cos_{\Delta x}(\frac{\tan b\Delta x}{\Delta x}, x)$ | $r = \frac{e^{a\Delta x}\cos b\Delta x - 1}{\Delta x} \pm j\ \frac{e^{a\Delta x}\sin b\Delta x}{\Delta x}$<br><br>$\cos b\Delta x \neq 0$ |
| **15a** | $ke^{ax}\sin bx$<br><br>$= ke_{\Delta x}(\frac{e^{a\Delta x}\cos b\Delta x-1}{\Delta x}, x)\ \sin_{\Delta x}(\frac{\tan b\Delta x}{\Delta x}, x)$<br><br>where<br><br>$\frac{x}{\Delta x}$ = integer<br><br>$\cos b\Delta x \neq 0$ | $[\sum_{p=N_h}^{N_h} A_p[x]_{\Delta x}^{p}]\ e_{\Delta x}(\frac{e^{a\Delta x}\cos b\Delta x-1}{\Delta x}, x)\ \sin_{\Delta x}(\frac{\tan b\Delta x}{\Delta x}, x) +$<br><br>$[\sum_{p=N_h}^{N_h} B_p[x]_{\Delta x}^{p}]\ e_{\Delta x}(\frac{e^{a\Delta x}\cos b\Delta x-1}{\Delta x}, x)\ \cos_{\Delta x}(\frac{\tan b\Delta x}{\Delta x}, x)$ | $r = \frac{e^{a\Delta x}\cos b\Delta x - 1}{\Delta x} \pm j\ \frac{e^{a\Delta x}\sin b\Delta x}{\Delta x}$<br><br>$\cos b\Delta x \neq 0$ |

| # | For these Interval Calculus Functions in Q(x)<br><br>**The function specified should be replaced in Q(x) by its identity if its identity is listed.**<br><br>**k = constant** | Put these Interval Calculus functions with undetermined coefficients in the particular solution, $f_P(x)$.<br><br>**$N_h$ = the number of times the Q(x) related root(s) appear in the characteristic polynomial, h(r)** | Related Root(s) |
|---|---|---|---|
| **15b** | $kxe^{ax}\sin bx$<br><br>$= kxe_{\Delta x}(\frac{e^{a\Delta x}\cos b\Delta x - 1}{\Delta x}, x)\ \sin_{\Delta x}(\frac{\tan b\Delta x}{\Delta x}, x)$<br><br>where<br><br>$\frac{x}{\Delta x}$ = integer<br><br>$\cos b\Delta x \neq 0$ | $[\sum_{p=N_h}^{1+N_h} A_p[x]_{\Delta x}^{p}]\ e_{\Delta x}(\frac{e^{a\Delta x}\cos b\Delta x - 1}{\Delta x}, x)\ \sin_{\Delta x}(\frac{\tan b\Delta x}{\Delta x}, x) +$<br><br>$[\sum_{p=N_h}^{1+N_h} B_p[x]_{\Delta x}^{p}]\ e_{\Delta x}(\frac{e^{a\Delta x}\cos b\Delta x - 1}{\Delta x}, x)\ \cos_{\Delta x}(\frac{\tan b\Delta x}{\Delta x}, x)$ | $r = \frac{e^{a\Delta x}\cos b\Delta x - 1}{\Delta x} \pm j\ \frac{e^{a\Delta x}\sin b\Delta x}{\Delta x}$<br><br>$\cos b\Delta x \neq 0$ |
| **16** | $[\sum_{p=0}^{N} k_p[x]_{\Delta x}^{p}]\ e^{ax}\cos bx$<br><br>=<br><br>$[\sum_{p=0}^{N} k_p[x]_{\Delta x}^{p}]\ e_{\Delta x}(\frac{e^{a\Delta x}\cos b\Delta x - 1}{\Delta x}, x)\ \cos_{\Delta x}(\frac{\tan b\Delta x}{\Delta x}, x)$<br><br>where<br><br>$\frac{x}{\Delta x}$ = integer<br><br>$\cos b\Delta x \neq 0$ | $[\sum_{p=N_h}^{N+N_h} A_p[x]_{\Delta x}^{p}]\ e_{\Delta x}(\frac{e^{a\Delta x}\cos b\Delta x - 1}{\Delta x}, x)\ \sin_{\Delta x}(\frac{\tan b\Delta x}{\Delta x}, x) +$<br><br>$[\sum_{p=N_h}^{N+N_h} B_p[x]_{\Delta x}^{p}]\ e_{\Delta x}(\frac{e^{a\Delta x}\cos b\Delta x - 1}{\Delta x}, x)\ \cos_{\Delta x}(\frac{\tan b\Delta x}{\Delta x}, x)$ | $r = \frac{e^{a\Delta x}\cos b\Delta x - 1}{\Delta x} \pm j\ \frac{e^{a\Delta x}\sin b\Delta x}{\Delta x}$<br><br>$\cos b\Delta x \neq 0$ |

| # | For these Interval Calculus Functions in Q(x)<br><br>The function specified should be replaced in Q(x) by its identity if its identity is listed.<br>k = constant | Put these Interval Calculus functions with undetermined coefficients in the particular solution, $f_P(x)$.<br>$N_h$ = the number of times the Q(x) related root(s) appear in the characteristic polynomial, h(r) | Related Root(s) |
|---|---|---|---|
| **16a** | $ke^{ax}\cos bx$<br>$= ke_{\Delta x}(\frac{e^{a\Delta x}\cos b\Delta x-1}{\Delta x}, x)\cos_{\Delta x}(\frac{\tan b\Delta x}{\Delta x}, x)$<br>where<br>$\frac{x}{\Delta x}$ = integer<br>$\cos b\Delta x \neq 0$ | $[\sum_{p=N_h}^{N_h} A_p[x]_{\Delta x}^{p}]\, e_{\Delta x}(\frac{e^{a\Delta x}\cos b\Delta x-1}{\Delta x}, x)\sin_{\Delta x}(\frac{\tan b\Delta x}{\Delta x}, x) +$<br>$[\sum_{p=N_h}^{N_h} B_p[x]_{\Delta x}^{p}]\, e_{\Delta x}(\frac{e^{a\Delta x}\cos b\Delta x-1}{\Delta x}, x)\cos_{\Delta x}(\frac{\tan b\Delta x}{\Delta x}, x)$ | $r = \frac{e^{a\Delta x}\cos b\Delta x - 1}{\Delta x} \pm j\,\frac{e^{a\Delta x}\sin b\Delta x}{\Delta x}$<br>$\cos b\Delta x \neq 0$ |
| **16b** | $kxe^{ax}\cos bx$<br>$= kxe_{\Delta x}(\frac{e^{a\Delta x}\cos b\Delta x-1}{\Delta x}, x)\cos_{\Delta x}(\frac{\tan b\Delta x}{\Delta x}, x)$<br>where<br>$\frac{x}{\Delta x}$ = integer<br>$\cos b\Delta x \neq 0$ | $[\sum_{p=N_h}^{1+N_h} A_p[x]_{\Delta x}^{p}]\, e_{\Delta x}(\frac{e^{a\Delta x}\cos b\Delta x-1}{\Delta x}, x)\sin_{\Delta x}(\frac{\tan b\Delta x}{\Delta x}, x) +$<br>$[\sum_{p=N_h}^{1+N_h} B_p[x]_{\Delta x}^{p}]\, e_{\Delta x}(\frac{e^{a\Delta x}\cos b\Delta x-1}{\Delta x}, x)\cos_{\Delta x}(\frac{\tan b\Delta x}{\Delta x}, x)$ | $r = \frac{e^{a\Delta x}\cos b\Delta x - 1}{\Delta x} \pm j\,\frac{e^{a\Delta x}\sin b\Delta x}{\Delta x}$<br>$\cos b\Delta x \neq 0$ |

| # | For these Interval Calculus Functions in Q(x)<br><br>The function specified should be replaced in Q(x) by its identity if its identity is listed.<br><br>k = constant | Put these Interval Calculus functions with undetermined coefficients in the particular solution, $f_P(x)$.<br><br>$N_h$ = the number of times the Q(x) related root(s) appear in the characteristic polynomial, h(r) | Related Root(s) |
|---|---|---|---|
| **17** | $[\sum_{p=0}^{N} k_p[x]_{\Delta x}^{p}] \sinh bx$<br><br>$= [\sum_{p=0}^{N} k_p[x]_{\Delta x}^{p}] e_{\Delta x}(\frac{\cosh b\Delta x - 1}{\Delta x}, x) \sinh_{\Delta x}(\frac{\tanh b\Delta x}{\Delta x}, x)$ | $[\sum_{p=N_h}^{N+N_h} A_p[x]_{\Delta x}^{p}] e_{\Delta x}(\frac{\cosh b\Delta x-1}{\Delta x}, x) \sinh_{\Delta x}(\frac{\tanh b\Delta x}{\Delta x}, x) +$<br><br>$[\sum_{p=N_h}^{N+N_h} B_p[x]_{\Delta x}^{p}] e_{\Delta x}(\frac{\cosh b\Delta x-1}{\Delta x}, x) \cosh_{\Delta x}(\frac{\tanh b\Delta x}{\Delta x}, x)$ | $r = \frac{\cosh b\Delta x - 1}{\Delta x} \pm \frac{\sinh b\Delta x}{\Delta x}$ |
| **17a** | ksinhbx<br><br>$= k e_{\Delta x}(\frac{\cosh b\Delta x - 1}{\Delta x}, x) \sinh_{\Delta x}(\frac{\tanh b\Delta x}{\Delta x}, x)$ | $[\sum_{p=N_h}^{N_h} A_p[x]_{\Delta x}^{p}] e_{\Delta x}(\frac{\cosh b\Delta x-1}{\Delta x}, x) \sinh_{\Delta x}(\frac{\tanh b\Delta x}{\Delta x}, x) +$<br><br>$[\sum_{p=N_h}^{N_h} B_p[x]_{\Delta x}^{p}] e_{\Delta x}(\frac{\cosh b\Delta x-1}{\Delta x}, x) \cosh_{\Delta x}(\frac{\tanh b\Delta x}{\Delta x}, x)$ | $r = \frac{\cosh b\Delta x - 1}{\Delta x} \pm \frac{\sinh b\Delta x}{\Delta x}$ |

| # | For these Interval Calculus Functions in Q(x)<br><br>The function specified should be replaced in Q(x) by its identity if its identity is listed.<br><br>k = constant | Put these Interval Calculus functions with undetermined coefficients in the particular solution, $f_P(x)$.<br><br>$N_h$ = the number of times the Q(x) related root(s) appear in the characteristic polynomial, h(r) | Related Root(s) |
|---|---|---|---|
| **17b** | kxsinhbx<br><br>$= kxe_{\Delta x}(\frac{\cosh b\Delta x - 1}{\Delta x}, x)\sinh_{\Delta x}(\frac{\tanh b\Delta x}{\Delta x}, x)$ | $[\sum_{p=N_h}^{1+N_h} A_p[x]_{\Delta x}^{p}]\, e_{\Delta x}(\frac{\cosh b\Delta x - 1}{\Delta x}, x)\sinh_{\Delta x}(\frac{\tanh b\Delta x}{\Delta x}, x) +$<br><br>$[\sum_{p=N_h}^{1+N_h} B_p[x]_{\Delta x}^{p}]\, e_{\Delta x}(\frac{\cosh b\Delta x - 1}{\Delta x}, x)\cosh_{\Delta x}(\frac{\tanh b\Delta x}{\Delta x}, x)$ | $r = \frac{\cosh b\Delta x - 1}{\Delta x} \pm \frac{\sinh b\Delta x}{\Delta x}$ |
| **18** | $[\sum_{p=0}^{N} k_p[x]_{\Delta x}^{p}]\cosh bx$<br><br>=<br><br>$[\sum_{p=0}^{N} k_p[x]_{\Delta x}^{p}]\, e_{\Delta x}(\frac{\cosh b\Delta x - 1}{\Delta x}, x)\cosh_{\Delta x}(\frac{\tanh b\Delta x}{\Delta x}, x)$ | $[\sum_{p=N_h}^{N+N_h} A_p[x]_{\Delta x}^{p}]\, e_{\Delta x}(\frac{\cosh b\Delta x - 1}{\Delta x}, x)\sinh_{\Delta x}(\frac{\tanh b\Delta x}{\Delta x}, x) +$<br><br>$[\sum_{p=N_h}^{N+N_h} B_p[x]_{\Delta x}^{p}]\, e_{\Delta x}(\frac{\cosh b\Delta x - 1}{\Delta x}, x)\cosh_{\Delta x}(\frac{\tanh b\Delta x}{\Delta x}, x)$ | $r = \frac{\cosh b\Delta x - 1}{\Delta x} \pm \frac{\sinh b\Delta x}{\Delta x}$ |

| # | For these Interval Calculus Functions in Q(x)<br><br>The function specified should be replaced in Q(x) by its identity if its identity is listed.<br><br>k = constant | Put these Interval Calculus functions with undetermined coefficients in the particular solution, $f_P(x)$.<br><br>$N_h$ = the number of times the Q(x) related root(s) appear in the characteristic polynomial, h(r) | Related Root(s) |
|---|---|---|---|
| **18a** | kcoshbx<br><br>$= ke_{\Delta x}(\frac{\cosh b\Delta x - 1}{\Delta x}, x) \cosh_{\Delta x}(\frac{\tanh b\Delta x}{\Delta x}, x)$ | $[\sum_{p=N_h}^{N_h} A_p[x]_{\Delta x}^{p}]\, e_{\Delta x}(\frac{\cosh b\Delta x-1}{\Delta x}, x) \sinh_{\Delta x}(\frac{\tanh b\Delta x}{\Delta x}, x) +$<br><br>$[\sum_{p=N_h}^{N_h} B_p[x]_{\Delta x}^{p}]\, e_{\Delta x}(\frac{\cosh b\Delta x-1}{\Delta x}, x) \cosh_{\Delta x}(\frac{\tanh b\Delta x}{\Delta x}, x)$ | $r = \frac{\cosh b\Delta x - 1}{\Delta x} \pm \frac{\sinh b\Delta x}{\Delta x}$ |
| **18b** | kxcoshbx<br><br>$= kxe_{\Delta x}(\frac{\cosh b\Delta x - 1}{\Delta x}, x) \cosh_{\Delta x}(\frac{\tanh b\Delta x}{\Delta x}, x)$ | $[\sum_{p=N_h}^{1+N_h} A_p[x]_{\Delta x}^{p}]\, e_{\Delta x}(\frac{\cosh b\Delta x-1}{\Delta x}, x) \sinh_{\Delta x}(\frac{\tanh b\Delta x}{\Delta x}, x) +$<br><br>$[\sum_{p=N_h}^{1+N_h} B_p[x]_{\Delta x}^{p}]\, e_{\Delta x}(\frac{\cosh b\Delta x-1}{\Delta x}, x) \cosh_{\Delta x}(\frac{\tanh b\Delta x}{\Delta x}, x)$ | $r = \frac{\cosh b\Delta x - 1}{\Delta x} \pm \frac{\sinh b\Delta x}{\Delta x}$ |

The following two examples are presented to show the application of the just derived Method of Undetermined Coefficients to the solution of differential difference equations.

<u>Example 4.2-1</u>

Solve the differential difference equation, $D_{\Delta x}y(x) + y(x) = xe^{-x}$, where $\Delta x = 2$ and $y(0) = 3$. Use the Method of Undetermined Coefficients.

$$D_{\Delta x}y(x) + y(x) = xe^{-x}, \qquad \Delta x = 2, \ y(0) = 3. \tag{1}$$

Find the related homogeneous equation complementary solution

$$D_{\Delta x}y(x) + y(x) = 0 \ , \quad \text{the related homogeneous differential difference equation} \tag{2}$$

Let

$$y(x) = Ke_{\Delta x}(a,x) = K(1+a\Delta x)^{\frac{x}{\Delta x}} \tag{3}$$

Substituting this y(x) function into the related homogenous differential difference equation, Eq 2

$$Kae_2(a,x) + Ke_2(a,x) = (1 + a)Ke_2(a,x) = 0 \ , \quad \Delta x = 2$$

$$a = -1$$

Substituting a = -1 into Eq 3

The complementary solution to the differential difference equation, Eq 1 is then

$$y_c(x) = Ke_2(-1,x) \tag{4}$$

Find the particular solution to the differential difference equation, Eq 1.

The function, $xe^{-x}$, should first be converted into its identity. This will facilitate the necessary discrete differentiation.

$$e^{cx} = e_{\Delta x}\left( \frac{e^{c\Delta x}-1}{\Delta x}, x \right)$$

$$c = -1, \quad \Delta x = 2$$

$$xe^{-x} = xe_2\left( \frac{e^{-2}-1}{2}, x \right) = xe_2(-.43233235,x) \tag{5}$$

Substituting into Eq 1

$$D_{\Delta x}y(x) + y(x) = xe_2(-.43233235,x) \tag{6}$$

Note that all of the functions in Eq 6 are discrete Interval Calculus functions

From the undetermined coefficient particular solution tables

$$y_p(x) = [\sum_{p=N_h}^{N+N_h} A_p[x]_{\Delta x}^{p}]\, e_{\Delta x}(\frac{e^{a\Delta x}-1}{\Delta x}, x) \qquad 7)$$

for the function

$$[\sum_{p=0}^{N} k_p[x]_{\Delta x}^{p}]\, e^{ax} = [\sum_{p=0}^{N} k_p[x]_{\Delta x}^{p}]\, e_{\Delta x}(\frac{e^{a\Delta x}-1}{\Delta x}, x) \qquad 8)$$

and the related root is

$$r = \frac{e^{a\Delta x}-1}{\Delta x}$$

For a = -1 and $\Delta x = 2$

The related root is r = -.432333235
The related homogeneous equation root is r = -1

Since the roots are not the same
$N_h = 0$

N = 1, the polynomial power, from Eq 6 and Eq 8 ( i.e. $x^1 e_2(-.43233235,x)$ )

From Eq 7

$$y_p(x) = [\sum_{p=0}^{1} A_p[x]_{\Delta x}^{p}]\, e_{\Delta x}(-.43233235,x)$$

$$y_p(x) = (A_0 + A_1 x)\, e_{\Delta x}(-.43233235,x) \qquad 9)$$

Substituting Eq 9, into Eq 6 to solve for the constants, $A_0$ , $A_1$

Using the derivative of a function product formula

$D_{\Delta x}[u(x)v(x) = v(x)D_{\Delta x}u(x) + D_{\Delta x}v(x)u(x+\Delta x)$ where $u(x) = (A_0 + A_1x)$ , $v(x) = e_{\Delta x}(-.43233235,x)$

$D_2y_p(x) = A_1e_2(-.43233235,x) + (-.43233235)e_2(-.43233235,x)(A_0 + A_1[x + 2])$

$D_2 y_p(x) = (A_1 - .43233235A_0 \; - .43233235A_1 x \; - .86466470A_1)e_2(-.43233235,x)$ 10)

Substituting Eq 9 and Eq 10 into Eq 6

$(A_1 - .43233235A_0 \; - .43233235A_1 x \; - .86466470A_1)e_2(-.43233235,x) + (A_0 + A_1 x)\, e_{\Delta x}( -.43233235,x) = xe_2(-.43233235,x)$

$(A_1 - .43233235A_0 - .86466470A_1 + A_0)\, e_2(-.43233235,x) + (A_1 - .43233235A_1)xe_2(-.43233235,x) = xe_2(-.43233235,x)$

Simplifying and equating constants

$.56766765A_0 + .13533530A_1 = 0$

$.56766765A_1 = 1$

Solving for $A_0$, $A_1$

$A_1 = 1.76159413$

$A_0 = \frac{-.13533530(1.76159413)}{.56766765} = -.41997438$

$A_0 = \; = -.41997438$

Substituting the values of $A_0$, $A_1$ into Eq 9

$y_p(x) = (-.41997438 + 1.76159413x)\, e_2(-.43233235,x)$ 11)

$y(x) = y_c(x) + y_p(x)$ 12)

Substituting Eq 4 and Eq 11 into Eq 12

$y(x) = Ke_2(-1,x) + (-.41997438 + 1.76159413x)\, e_2(-.43233235,x)$ 13)

Find K from the initial condition y(0) = 3

$3 = Ke_2(-1,0) - 41997438e_2(-.43233235,0)$

$e_{\Delta x}(c,0) = (1+c\Delta x)^{\frac{0}{\Delta x}} = 1$

$3 = K - .41997438$

$K = 3.41997438$

Substituting the value of K into Eq 13

Then

$$\mathbf{y(x) = 3.14199743e_2(-1,x) + (-.41997438 + 1.76159413x)\, e_2(-.43233235,x)} \qquad 14)$$

is the solution to the differential difference equation

$D_{\Delta x}y(x) + y(x) = xe^{-x} = xe_2(-.43233235, x)$

where

$x = 0, 2, 4, 6, \ldots$

If desired, Eq 14 can be changed to another form

$$e_{\Delta x}(c,x) = (1+c\Delta x)^{\frac{x}{\Delta x}}$$

$$e_2(-1,x) = [1+(-1)(2)]^{\frac{x}{2}} = (-1)^{\frac{x}{2}}$$

From Eq 5

$e_2(-.43233235,x) = e^{-x}$

Substituting these equalities into Eq 14

$$\mathbf{y(x) = 3.14199743(-1)^{\frac{x}{2}} + (-.41997438 + 1.76159413x)\, e^{-x}} \qquad 15)$$

Check Eq 14 and Eq 15

From Eq 1

$D_{\Delta x}y(x) + y(x) - xe^{-x} = 0$

$$\frac{y(x+\Delta x)-y(x)}{\Delta x} + y(x) - xe^{-x} = 0$$

For

$x = 2$

$\Delta x = 2$

$$\frac{y(4)-y(2)}{2} + y(2) - xe^{-2} = 0$$

Finding y(2) and y(4) from Eq 15

y(2) = -2.7220231
y(4) =  3.2633642

Substituting

$$\frac{3.2633642 - (-2.7220231)}{2} + (-2.7220231) - .27067056 = 0$$

0 = 0        Good check

Example 4.2-2

Solve the differential difference equation, $D_{\Delta x}^2\, y(x) + D_{\Delta x}\, y(x) = 3$, where $\Delta x = 2$, $y(0) = 1$ and $y(2) = 4$. Use the Method of Undetermined Coefficients.

$$D_{\Delta x}^2\, y(x) + D_{\Delta x}\, y(x) = 3 \qquad 1)$$

where

$\Delta x = 2$

$y(0) = 1$

$y(2) = 4$

Find the related homogenous equation complementary solution to Eq 1

$$D_{\Delta x}^2\, y(x) + D_{\Delta x}\, y(x) = 0 \qquad 2)$$

Let

$$y(x) = Ke_{\Delta x}(a,x) = K(1+a\Delta x)^{\frac{x}{\Delta x}} \qquad 3)$$

Substitute Eq 3 into Eq 2

$$Ka^2e_{\Delta x}(a,x) + Kae_{\Delta x}(a,x) = 0$$

$$(a^2 + a)\, Ke_{\Delta x}(a,x) = 0$$

$$a = 0, -1$$

The complementary solution to Eq 1 is:

$$y_c(x) = K_1e_{\Delta x}(-1,x) + K_2e_{\Delta x}(0,x)$$

$$y_c(x) = K_1e_{\Delta x}(-1,x) + K_2 \qquad 4)$$

Find the particular solution to Eq 1

From the particular solution tables where Q(x) is a constant (here it is 3)

$$\text{For } Q(x) = \sum_{p=0}^{N} k_p[x]_{\Delta x}^{p} \qquad 5)$$

The undetermined coefficient $y_p(x)$ is:

$$y_p(x) = \sum_{p=N_h}^{N+N_h} A_p[x]_{\Delta x}^{p} \qquad 6)$$

$N = 0$,    from Eq 1 for $Q(x) = 3$
$\Delta x = 2$

The related root is given as $r = 0$ but this root is equal to one of the complementary solution roots, $r = 0$.

Then
$N_h = 1$

Substituting the values for $N$, $N_h$, $\Delta x$ into Eq 6

$$y_p(x) = \sum_{p=1}^{1} A_p[x]_2^p = A_1[x]_2^1 = A_1 x$$

$$y_p(x) = A_1 x \qquad 7)$$

Substitute Eq 7 into Eq 1

$$D_{\Delta x}^{\ 2}(A_1 x) + D_{\Delta x}(A_1 x) = 3$$

$$0 + A_1 = 3$$

$$A_1 = 3$$

Substituting the value of $A_1$ into Eq 7

$$y_p(x) = 3x \qquad 8)$$

Find $y(x)$

$$y(x) = y_c(x) + y_p(x) \qquad 9)$$

From Eq 4, Eq8, and Eq 9

$$y(x) = K_1 e_{\Delta x}(-1,x) + K_2 + 3x \qquad 10)$$

From the initial conditions find $K_1$,$K_2$

$$y(0) = 1 = K_1 e_{\Delta x}(-1,0) + K_2$$

$$K_1 + K_2 = 1$$

$$y(2) = 4 = K_1 e_2(-1,2) + K_2 + 6\ , \quad \Delta x = 2$$

$$e_2(-1,2) = [1+(-1)(2)]^{\frac{2}{2}} = -1$$

$4 = -K_1 + K_2 + 6$

$-K_1 + K_2 = -2$
$K_1 + K_2 = 1$

$2K_2 = -1$

$K_2 = -\frac{1}{2}$

$K_1 = \frac{3}{2}$

Substituting_ the values of $K_1$,$K_2$ into Eq 10

$$\mathbf{y(x) = \frac{3}{2}\, e_2(-1,x) - \frac{1}{2} + 3x} \qquad 11)$$

Since

$e_2(-1,x) = [1+(-1)(2)]^{\frac{x}{2}} = (-1)^{\frac{x}{2}}$

$$\mathbf{y(x) = \frac{3}{2}\,(-1)^{\frac{x}{2}} - \frac{1}{2} + 3x} \qquad 12)$$

Both Eq 11 and Eq 12 are equivalent forms of the solution to the differential difference equation

$D_{\Delta x}{}^{2}\, y(x) + D_{\Delta x}\, y(x) = 3$, where $\Delta x = 2$, $y(0) = 1$ and $y(2) = 4$.

Checking Eq 11 and Eq 12

$y(0) = \frac{3}{2}(-1)^{0} - \frac{1}{2} + 0 = \frac{3}{2} - \frac{1}{2} = 1$ Good check

$y(2) = \frac{3}{2}(-1)^{1} - \frac{1}{2} + (3)(2) = -\frac{3}{2} - \frac{1}{2} + 6 = 4$ Good check

Substituting Eq 11 into Eq 1

$(-1)^2\frac{3}{2}\, e_2(-1,x) + (-1)\,\frac{3}{2}\, e_2(-1,x) + 3 = 3$

$3 = 3$ Good check

This Method of Undetermined Coefficients for solving differential difference equations is one of several methods. The other methods for solving differential difference equations are described in the sections which follow.

## Section 4.3: The solution of differential difference equations using the $K_{\Delta x}$ Transform Method

The $K_{\Delta x}$ Transform Method for solving differential difference equations has been described in Chapter 1. The $K_{\Delta x}$ Transform was previously introduced to show that discrete Interval Calculus also includes, as does Calculus, an operational method for solving differential equations. Originally, the $K_{\Delta x}$ Transform Method was found through observation, superficial logic, common sense, and trial and error. Of course, mathematics, in general, can not be structured on such a foundation. However, such an analysis can provide invaluable insight. Perhaps in some measure by luck, the type of analysis just described was very effective in deriving the $K_{\Delta x}$ Transform. The $K_{\Delta x}$ Transform, as derived, was found to work well in all applications in which it was used and tested. It was not until much later that the $K_{\Delta x}$ Transform was rigorously derived by adapting the Fourier Series and the Fourier Transform for discrete mathematics use. The derivation of the $K_{\Delta x}$ Transform Transform (and its closely related $J_{\Delta x}$ Transform) is shown in Chapter 5. After studying the derivation of the $K_{\Delta x}$ Transform, mathematicians may prefer to call it an extension of the Laplace Transform for discrete mathematics use. In fact, the Laplace Transform is the $K_{\Delta x}$ Transform where $\Delta x$ approaches zero. Depending on the perspective, the $K_{\Delta x}$ Transform can be considered to be the Laplace Transform generalized to accommodate functions of a non-continuous variable (i.e. a discrete variable of the form $x = m\Delta x$ where m is an integer). In any case, this paper will continue to use the originally proposed $K_{\Delta x}$ Transform designation.

In Chapter 1, the $K_{\Delta x}$ Transform was shown to have the ability to solve difference equations, differential difference equations, and engineering control system related problems. In this chapter, Chapter 4, some additional description of the $K_{\Delta x}$ transform will be provided.

From its derivation in Chapter 5, the $K_{\Delta x}$ Transform and its inverse are expressed as follows:

Table 4.3-1

The $K_{\Delta x}$ Transform

$$F(s) = \underset{0}{\overset{\infty}{\int}}{}_{\Delta x} [1+s\Delta x]^{-\frac{x+\Delta x}{\Delta x}} F(t)\, \Delta x$$

The Inverse $K_{\Delta x}$ Transform

$$F(x) = \frac{1}{2\pi j} \oint_c [1+s\Delta x]^{\frac{x}{\Delta x}} F(s)\, ds$$

$$s = \frac{e^{(\gamma+jw)\Delta x} - 1}{\Delta x}, \quad -\frac{\pi}{\Delta t} \le w < \frac{\pi}{\Delta t}$$

where

$x = m\Delta x$, m = integers

$\Delta x$ = x increment

$F(x)$ = discrete function of x

$F(x) = 0$ for $x < 0$

Complex Plane

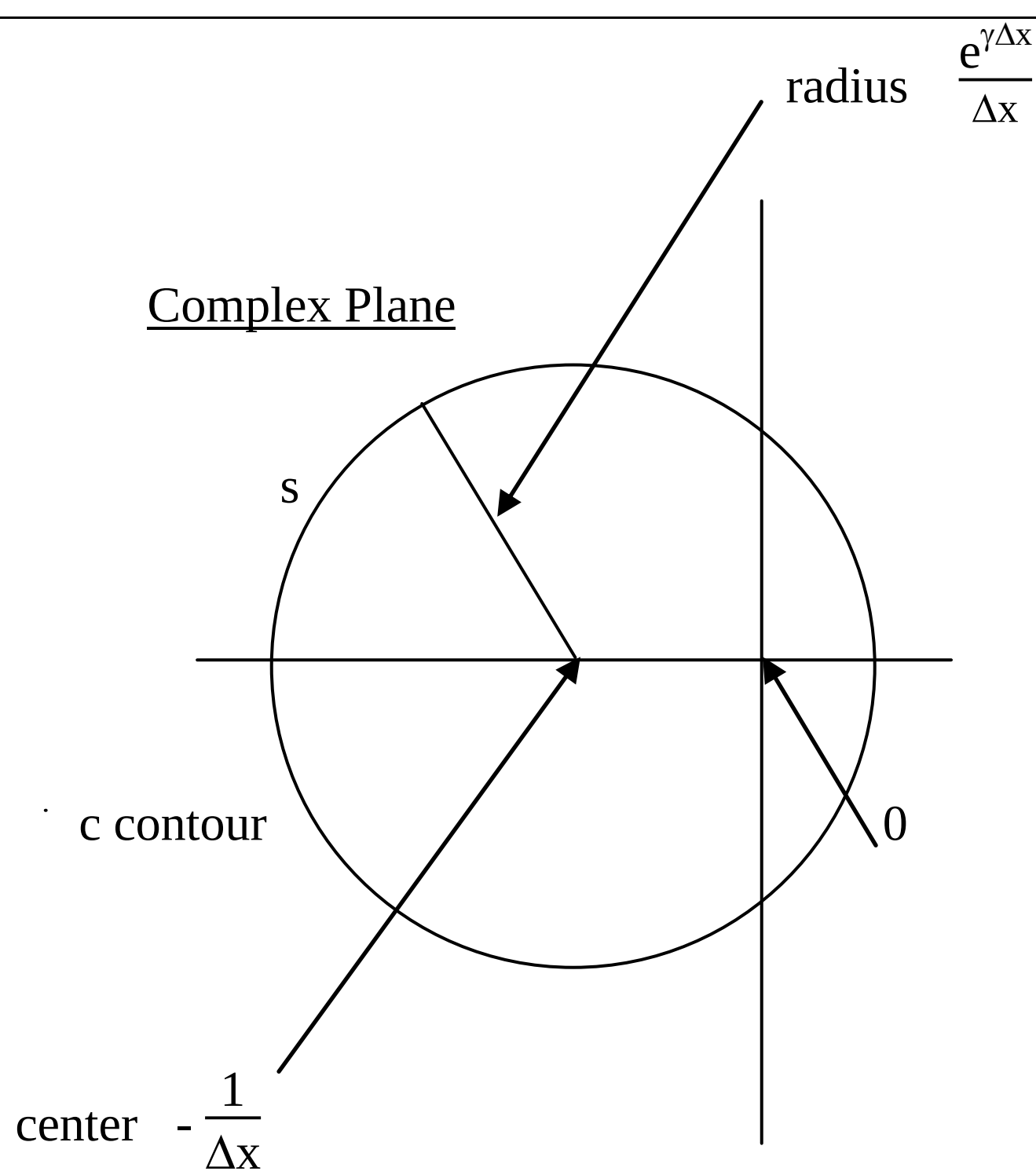

F(s) = function of s

$\gamma$, w = real value constants

$\gamma > 0$

C, the complex plane contour of integration, is a circle of radius, $\frac{e^{\gamma\Delta x}}{\Delta x}$, with center at $-\frac{1}{\Delta x}$.

$\gamma$ is chosen so that the contour encloses all poles of F(s).

Note - The $K_{\Delta x}$ Transform becomes the Laplace Transform for $\Delta x \to 0$.

To demonstrate the use of both the $K_{\Delta x}$ Transform and the Inverse $K_{\Delta x}$ Transform several examples will now be presented.

Example 4.3-1 Find the $K_{\Delta x}$ Transform of the function $e_{\Delta x}(a,x)$.

$$e_{\Delta x}(a,x) = (1+a\Delta x)^{\frac{x}{\Delta x}} \tag{4.3-1}$$

$$K_{\Delta x}[f(x)] = {}_{\Delta x}\int_0^\infty (1+s\Delta x)^{-(\frac{x+\Delta x}{\Delta x})} f(x)\Delta x \tag{4.3-2}$$

Let

$$f(x) = e_{\Delta x}(a,x) = (1+a\Delta x)^{\frac{x}{\Delta x}} \tag{4.3-3}$$

Substituting Eq 4.3-3 into the $K_{\Delta x}$ Transform equation, Eq 4.3-2

$$K_{\Delta x}[e_{\Delta x}(a,x)] = \frac{1}{1+s\Delta x}\ {}_{\Delta x}\int_0^\infty \left(\frac{1+s\Delta x}{1+a\Delta x}\right)^{-(\frac{x}{\Delta x})} \Delta x \tag{4.3-4}$$

$$\frac{1+s\Delta x}{1+a\Delta x} = 1 + \frac{(s-a)\Delta x}{1+a\Delta x} \tag{4.3-5}$$

Substituting Eq 4.3-5 into Eq 4.3-4

$$K_{\Delta x}[e_{\Delta x}(a,x)] = \frac{1}{1+s\Delta x}\ {}_{\Delta x}\int_0^\infty \left(1+\frac{(s-a)\Delta x}{1+a\Delta x}\right)^{-(\frac{x}{\Delta x})} \Delta x \tag{4.3-6}$$

Substitute the quantity, $\frac{s-a}{1+a\Delta x}$, for b in the following equation. This integral equation was obtained from the integral equation table, Table 6, in the Appendix.

$$_{\Delta x}\int (1+b\Delta x)^{-\frac{x}{\Delta x}}\,\Delta x = -\,\frac{1+b\Delta x}{b}(1+b\Delta x)^{-\frac{x}{\Delta x}} + k \qquad (4.3\text{-}7)$$

$$K_{\Delta x}[e_{\Delta x}(a,x)] = \frac{1}{1+s\Delta x}\;{}_{\Delta x}\!\int_0^{\infty}\left(1+\frac{(s-a)\Delta x}{1+a\Delta x}\right)^{-\frac{x}{\Delta x}}\Delta x = -\,\frac{1}{1+s\Delta x}\,\frac{1+s\Delta x}{s-a}\left(1+\frac{(s-a)\Delta x}{1+a\Delta x}\right)^{-\frac{x}{\Delta x}}\Big|_0^{\infty} \qquad (4.3\text{-}8)$$

$$\mathbf{K_{\Delta x}[e_{\Delta x}(a,x)] = \frac{1}{s-a}} \qquad (4.3\text{-}9)$$

Eq 4.3-10 specifies the $K_{\Delta x}$ Transform of $e_{\Delta x}(a,x)$

Example 4.3-2 Find the Inverse $K_{\Delta x}$ Transform of the function, $f(s) = \frac{1}{s-a}$, using the Inverse $K_{\Delta x}$ Transform Integral.

$$K_{\Delta x}^{-1}[f(s)] = F(x) = \frac{1}{2\pi j}\oint_c [1+s\Delta x]^{\frac{x}{\Delta x}} f(s)\,ds \text{ , The Inverse } K_{\Delta x} \text{ Transform} \qquad (4.3\text{-}10)$$

Notes – The x variable of F(x) is discrete where x = 0, Δx, 2Δx, 3Δx, …

The s variable of the complex plane closed contour integral of Eq 4.3-10 is continuous along the closed contour, c.

$\frac{x}{\Delta x}$ = 0, 1, 2, 3, …

The closed contour, c, in the complex plane is shown in Table 4.3-1. The positive constant, γ, is made large enough to have the closed contour, c, encircle all poles of $(1+s\Delta x)^{\frac{x}{\Delta x}}$ f(s).

To calculate the integral of Eq 4.3-10, the theory of residues will be used.

Provided that the function of s, $[1+s\Delta x]^{\frac{x}{\Delta x}}$ f(s), can be expanded into a convergent Laurent Series for each pole of f(s)

$$K_{\Delta x}^{-1}[f(s)] = \frac{1}{2\pi j}\oint_c [1+s\Delta x]^{\frac{x}{\Delta x}} f(s)\,ds = \sum_{p=1}^{P} R_p \tag{4.3-11}$$

$$R = \lim_{s\to a}\frac{1}{(n-1)!}\frac{d^{n-1}}{ds^{n-1}}[(s-r)^n (1+s\Delta x)^{\frac{x}{\Delta x}} f(s)]$$ , the residue calculation formula for a pole at s = r

where (4.3-12)

$R$ = the residue of a pole of $(1+s\Delta x)^{\frac{x}{\Delta x}} f(s)$

$\sum_{p=1}^{P} R_p$ = the sum of the P residues of the P poles of $(1+s\Delta x)^{\frac{x}{\Delta x}} f(s)$

Let $f(s) = \frac{1}{s-a}$

$r = a$

$P = 1$ , $(1+s\Delta x)^{\frac{x}{\Delta x}} \frac{1}{s-a}$ has one first order pole at s = a

Substituting the specified f(s), r, and P into Eq 4.3-11

$$K_{\Delta x}^{-1}[\frac{1}{s-a}] = \frac{1}{2\pi j}\oint_c [1+s\Delta x]^{\frac{x}{\Delta x}} \frac{1}{s-a}\,ds = \sum_{p=1}^{1} R_p = R_1 \tag{4.3-13}$$

$n = 1$

<u>Note</u> - $(1+s\Delta x)^{\frac{x}{\Delta x}} \frac{1}{s-a}$ has a convergent Laurent Series which is:

$$\frac{(1+s\Delta x)^{\frac{x}{\Delta x}}}{s-a} = e_{\Delta x}(a,x)\sum_{n=0}^{\infty}\frac{1}{n!}\frac{[x]_{\Delta x}^{n}}{(1+a\Delta x)^n}\frac{1}{(s-a)^{1-n}}$$

Since there is only one pole to be encircled by the complex plane closed contour, c, there is only one residue required, $R_1$, the residue of the first order pole at s = a.

Substituting the specified f(s) and r into Eq 4.3-12

$$R_1 = \lim_{s\to a}[(s-a)(1+s\Delta x)^{\frac{x}{\Delta x}}\frac{1}{s-a}] = (1+a\Delta x)^{\frac{x}{\Delta x}} = e_{\Delta x}(a,x) \tag{4.3-14}$$

Substituting $R_1$ into Eq 4.3-13

$$\mathbf{K_{\Delta x}^{-1}[\frac{1}{s-a}] = \frac{1}{2\pi j} \oint_c [1+s\Delta x]^{\frac{x}{\Delta x}} \frac{1}{s-a} ds = e_{\Delta x}(a,x)} \tag{4.3-15}$$

Eq 4.3-15 specifies the Inverse $K_{\Delta x}$ Transform of $\frac{1}{s-a}$ .

Example 4.3-2 above calculated the Inverse $K_{\Delta x}$ Transform of the function, $f(s) = \frac{1}{s-a}$ , using the Inverse $K_{\Delta x}$ Transform Integral. Residue theory was used. In particular, the necessary residue (at s = a) of the function, $[1+s\Delta x]^{\frac{x}{\Delta x}} \frac{1}{s-a}$ , was calculated using the residue calculation formula. Residues of a function can also be calculated directly from the function's convergent Laurent Series expansion. This will be demonstrated below in Example 4.3-3.

<u>Example 4.3-3</u> Find the Inverse $K_{\Delta x}$ Transform of the function, $f(s) = \frac{1}{(s-a)^m}$ , m = 1,2,3,… , using the Inverse $K_{\Delta x}$ Transform Integral. Evaluate the integral directly from the convergent Laurent Series expansion of the function, $[1+s\Delta x]^{\frac{x}{\Delta x}} \frac{1}{(s-a)^m}$ .

$$K_{\Delta x}^{-1}[f(s)] = F(x) = \frac{1}{2\pi j} \oint_c [1+s\Delta x]^{\frac{x}{\Delta x}} f(s)\, ds \text{ , The Inverse } K_{\Delta x} \text{ Transform} \tag{4.3-16}$$

<u>Notes</u> – The x variable of F(x) is discrete where x = 0, Δx, 2Δx, 3Δx, …

The s variable of the complex plane closed contour integral of Eq 4.3-10 is continuous along the closed contour, c.

$\frac{x}{\Delta x} = 0, 1, 2, 3, \ldots$

The closed contour, c, in the complex plane is shown in Table 4.3-1. The positive constant, γ, is made large enough to have the closed contour, c, encircle all poles of $(1+s\Delta x)^{\frac{x}{\Delta x}} f(s)$.

$$f(s) = \frac{1}{(s-a)^m} \tag{4.3-17}$$

Substituting Eq 4.3-17 into Eq 4.3-16

$$K_{\Delta x}^{-1}[f(s)] = F(x) = \frac{1}{2\pi j} \oint_c [1+s\Delta x]^{\frac{x}{\Delta x}} \frac{1}{(s-a)^m} \, ds \tag{4.3-18}$$

To calculate the integral of Eq 4.3-18, a convergent Laurent Series of the function,

$[1+s\Delta x]^{\frac{x}{\Delta x}} \frac{1}{(s-a)^m}$, will now be found.

$$H(s) = \frac{[1+s\Delta x]^{\frac{x}{\Delta x}}}{(s-a)^m} \quad , \quad m = 1,2,3,\ldots \tag{4.3-19}$$

Let

$$u = s-a \tag{4.3-20}$$

$$s = a+u \tag{4.3-21}$$

$$e_{\Delta x}(s,x) = [1+s\Delta x]^{\frac{x}{\Delta x}} \tag{4.3-22}$$

Substituting Eq 4.3-20 thru Eq 4.3-22 into Eq 4.3-19

$$H(s) = \frac{[1+(a+u)\Delta x]^{\frac{x}{\Delta x}}}{u^m} = \frac{1}{u^m} e_{\Delta x}(a+u,x) \tag{4.3-23}$$

$$e_{\Delta x}(a+u,x) = e_{\Delta x}(a,x)\, e_{\Delta x}\left(\frac{u}{1+a\Delta x},x\right) \tag{4.3-24}$$

Substituting Eq 4.3-24 int Eq 4.3-23

$$H(s) = \frac{1}{u^m} e_{\Delta x}(a,x)\, e_{\Delta x}\left(\frac{u}{1+a\Delta x},x\right) \tag{4.3-25}$$

$$e_{\Delta x}(r,x) = \sum_{n=0}^{\infty} \frac{r^n}{n!} [x]_{\Delta x}^{n} \, , \quad \text{convergent series} \tag{4.3-26}$$

Substituting Eq 4.3-26 into Eq 4.3-25 with $r = \frac{u}{1+a\Delta x}$ and u = s-a

$$H(s) = \frac{1}{(s-a)^m} e_{\Delta x}(a,x) \sum_{n=0}^{\infty} \frac{\left(\frac{s-a}{1+a\Delta x}\right)^n}{n!} [x]_{\Delta x}^{n} \tag{4.3-27}$$

Simplifying Eq 4.3- 27

$$H(s) = e_{\Delta x}(a,x) \sum_{n=0}^{\infty} \frac{[x]_{\Delta x}^{n}}{n!\,(1+a\Delta x)^n} \frac{1}{(s-a)^{m-n}} \tag{4.3-28}$$

Substituting Eq 4.3-19 into 4.3-28

$$\frac{[1+s\Delta x]^{\frac{x}{\Delta x}}}{(s-a)^m} = e_{\Delta x}(a,x) \sum_{n=0}^{\infty} \frac{[x]_{\Delta x}^{n}}{n!\,(1+a\Delta x)^n} \frac{1}{(s-a)^{m-n}} \,, \qquad \text{The Laurent Series} \tag{4.3-29}$$

where

m = 1,2,3,…

m = the order of the pole at s = a

Expanding Eq 4.3-29

$$\frac{[1+s\Delta x]^{\frac{x}{\Delta x}}}{(s-a)^m} = e_{\Delta x}(a,x) \left[ \frac{1}{(s-a)^m} + \frac{[x]_{\Delta x}^{1}}{1!\,(1+a\Delta x)^1} \frac{1}{(s-a)^{m-1}} + \frac{[x]_{\Delta x}^{2}}{2!\,(1+a\Delta x)^2} \frac{1}{(s-a)^{m-2}} + \frac{[x]_{\Delta x}^{3}}{3!\,(1+a\Delta x)^3} \frac{1}{(s-a)^{m-3}} + \dots \right. \tag{4.3-30}$$

where

m = 1,2,3,… ,

m = the order of the pole at s = a

Using residue theory

Since there is only one single or multiple order pole (at s = a) in the function, $\frac{[1+s\Delta x]^{\frac{x}{\Delta x}}}{(s-a)^m}$, there is only one residue, $R_m$.

This residue of the function, $\frac{[1+s\Delta x]^{\frac{x}{\Delta x}}}{(s-a)^m}$, can be found from its convergent Laurent Series, Eq 4.3-30. The residue, $R_m$, is the coefficient of the $\frac{1}{(s-a)^1}$ series term. This coefficient changes as a function of m.

$$\frac{1}{2\pi j}\oint_c [1+s\Delta x]^{\frac{x}{\Delta x}}\frac{1}{(s-a)^m}\,ds = R_m \quad \text{(The sum of the residues, in this case one residue, is } R_m.) \qquad (4.3\text{-}31)$$

$$\frac{1}{2\pi j}\oint_c [1+s\Delta x]^{\frac{x}{\Delta x}}\frac{1}{(s-a)^1}\,ds = R_1 = e_{\Delta x}(a,x)\,, \qquad \text{for } m = 1 \qquad (4.3\text{-}32)$$

$$\frac{1}{2\pi j}\oint_c [1+s\Delta x]^{\frac{x}{\Delta x}}\frac{1}{(s-a)^2}\,ds = R_2 = \frac{[x]_{\Delta x}^{1}\, e_{\Delta x}(a,x)}{1!\,(1+a\Delta x)^1}\,, \qquad \text{for } m = 2 \qquad (4.3\text{-}33)$$

$$\frac{1}{2\pi j}\oint_c [1+s\Delta x]^{\frac{x}{\Delta x}}\frac{1}{(s-a)^3}\,ds = R_3 = \frac{[x]_{\Delta x}^{2}\, e_{\Delta x}(a,x)}{2!\,(1+a\Delta x)^2}\,, \qquad \text{for } m = 3 \qquad (4.3\text{-}34)$$

$$\frac{1}{2\pi j}\oint_c [1+s\Delta x]^{\frac{x}{\Delta x}}\frac{1}{(s-a)^4}\,ds = R_4 = \frac{[x]_{\Delta x}^{3}\, e_{\Delta x}(a,x)}{3!\,(1+a\Delta x)^3}\,, \qquad \text{for } m = 4 \qquad (4.3\text{-}35)$$

In general, from Eq 4.3-16, Eq 4.3-17 and Eq 4.3-32 thru Eq 4.3-35

$$\mathbf{K_{\Delta x}^{-1}\left[\frac{1}{(s-a)^m}\right] = F(x) = \frac{1}{2\pi j}\oint_c [1+s\Delta x]^{\frac{x}{\Delta x}}\frac{1}{(s-a)^m}\,ds = \frac{[x]_{\Delta x}^{m-1}\, e_{\Delta x}(a,x)}{(m-1)!\,(1+a\Delta x)^{m-1}}} \qquad (4.3\text{-}36)$$

**where**

**m = 1,2,3,…**

**$\Delta x$ = x increment**

**x = 0, $\Delta x$, $2\Delta x$, $3\Delta x$, …**

Comment – This result agrees with the equation in the $K_{\Delta x}$ Transform table, Table 3, in the Appendix. Eq 4.3-36 represents the Inverse $K_{\Delta x}$ Transform of the function, $f(s) = \frac{1}{(s-a)^m}$ , for m = 1,2,3,…

Checking

Eq 4.3-36 can be used to find the Inverse $K_{\Delta x}$ Transform of the function, $f(s) = \frac{1}{s-a}$, using the Inverse $K_{\Delta x}$ Transform Integral.

For m = 1

$$K_{\Delta x}^{-1}[\frac{1}{(s-a)}] = F(x) = \frac{1}{2\pi j}\oint_c [1+s\Delta x]^{\frac{x}{\Delta x}}\frac{1}{(s-a)}\,ds = e_{\Delta x}(a,x) \qquad (4.3\text{-}37)$$

Note that Eq 4.3-37 and Eq 4.3-15 are the same.

Good check

The method of finding the Inverse $K_{\Delta x}$ Transform using its complex plane closed contour integral definition, while available, is not often used. Complex mathematics calculations can be tedious. Calculations using the $K_{\Delta x}$ Transform integral are generally easier to perform. By producing a table of commonly used functions of x and their corresponding transform functions of s, an inverse transform for a specific function of s can be easily obtained.

Transform tables are quite efficient for finding both $K_{\Delta x}$ Transforms and Inverse $K_{\Delta x}$ Transforms. For this reason, transform tables are commonly used when problem solving employs the use of $K_{\Delta x}$ Transform methods. Table 2 and Table 3 in the Appendix are $K_{\Delta x}$ Transform Tables. Some of the $K_{\Delta x}$ transforms from each table are presented below.

# TABLE 2

## $K_{\Delta x}$ Transform General Equations

| | |
|---|---|
| **1** | $K_{\Delta x}[f(x)] = {}_{\Delta x}\int_0^{\infty} (1+s\Delta x)^{-(\frac{x+\Delta x}{\Delta x})} f(x)\Delta x$ , The $K_{\Delta x}$ Transform |
| **1a** | $K_{\Delta x}[f(x)] = \int_0^{\infty} e_{\Delta x}(s,-x-\Delta x)f(x)\Delta x$ |
| **1b** | $K_{\Delta x}[f(x)] = \Delta x\sum_{n=0}^{\infty} f(n\Delta x)K_{\Delta x}[\frac{1}{\Delta x}\{U(x-n\Delta x)-U(x-n\Delta x-\Delta x)\}]$ Sum of Unit Area Pulse $K_{\Delta x}$ Transforms |
| **1c** | $K_{\Delta x}[f(x)] = \sum_{n=0}^{\infty} f(n\Delta x)[\frac{(1+s\Delta x)^{-n} - (1+s\Delta x)^{-n-1}}{s}]$ |
| **1d** | $K_{\Delta x}[f(x)] = \Delta x \sum_{n=0}^{\infty} f(n\Delta x)(1+s\Delta x)^{-n-1}) \equiv \Delta x \sum_{\Delta x \; x=0}^{\infty} f(x)(1+s\Delta x)^{-(\frac{x+\Delta x}{\Delta x})}$ , $x = n\Delta x$ |
| **2** | $K_{\Delta x}[cf(x)] = cK_{\Delta x}[f(x)]$ |
| **3** | $K_{\Delta x}[f(x) + g(x)] = K_{\Delta x}[f(x)] + K_{\Delta x}[g(x)]$ |
| **4** | $K_{\Delta x}[D_{\Delta x}f(x)] = sK_{\Delta x}[f(x)] - f(0)$ |
| ***5** | $K_{\Delta x}[D^n{}_{\Delta x}f(x)] = s^n K_{\Delta x}[f(x)] - s^{n-1}D^0{}_{\Delta x}f(0) - s^{n-2}D^1{}_{\Delta x}f(0) - s^{n-3}D^2{}_{\Delta x}f(0) - \ldots - s^0 D^{n-1}{}_{\Delta x}f(0)$ $n = 1,2,3,\ldots$ |
| **6** | $K_{\Delta x}[f(x+\Delta x)] = (1+s\Delta x)K_{\Delta x}[f(x)] - \Delta x f(0)$ |
| **7** | $K_{\Delta x}[f(x+n\Delta x)] = (1+s\Delta x)^n K_{\Delta x}[f(x)] - \Delta x \sum_{m=1}^{n}(1+s\Delta x)^{n-m} f([m-1]\Delta x)$ , $n = 1,2,3,\ldots$ |

* A $K_{\Delta x}$ Transform used in the solution of the differential difference equation of Example 4.3-4

# TABLE 3

## $K_{\Delta x}$ Transforms

## (Compared to the Lapace Transform)

$$K_{\Delta x}[f(x)] = {}_{\Delta x}\int_{0}^{\infty} (1+s\Delta x)^{-(\frac{x+\Delta x}{\Delta x})} f(x)\Delta x$$

| # | $f_{\Delta x}(x)$ | $K_{\Delta x}[f_{\Delta x}(x)]$ | f(x) | L[f(x)] |
|---|---|---|---|---|
| **1** | 1 | $\frac{1}{s}$ | 1 | $\frac{1}{s}$ |
| **2** | c | $\frac{c}{s}$ | c | $\frac{c}{s}$ |
| **3** | x | $\frac{1}{s^2}$ | x | $\frac{1}{s^2}$ |

| # | $f_{\Delta x}(x)$ | $K_{\Delta x}[f_{\Delta x}(x)]$ | f(x) | L[f(x)] |
|---|---|---|---|---|
| *4 | $e_{\Delta x}(a,x)$ | $\frac{1}{s-a}$ root $s = a$ | $e^{ax}$ | $\frac{1}{s-a}$ root $s = a$ |
| *5 | $xe_{\Delta x}(a,x)$ | $\frac{1+a\Delta x}{(s-a)^2}$ | $xe^{ax}$ | $\frac{1}{(s-a)^2}$ |
| 6 | $\sin_{\Delta x}(b,x)$ | $\frac{b}{s^2+b^2}$ roots $s = jb, -jb$ | sinbx | $\frac{b}{s^2+b^2}$ roots $s = jb, -jb$ |
| 7 | $\cos_{\Delta x}(b,x)$ | $\frac{s}{s^2+b^2}$ roots $s = jb, -jb$ | cosbx | $\frac{s}{s^2+b^2}$ roots $s = jb, -jb$ |
| 8 | $(1+a\Delta x)^{\frac{x}{\Delta x}}\sin_{\Delta x}(\frac{b}{1+a\Delta x}, x)$ | $\frac{b}{(s-a)^2+b^2}$ roots $s = a+jb, a-jb$ | $e^{ax}\sin bx$ | $\frac{b}{(s-a)^2+b^2}$ roots $s = a+jb, a-jb$ |
| 9 | $(1+a\Delta x)^{\frac{x}{\Delta x}}\cos_{\Delta x}(\frac{b}{1+a\Delta x}, x)$ | $\frac{s-a}{(s-a)^2+b^2}$ roots $s = a+jb, a-jb$ | $e^{ax}\cos bx$ | $\frac{s-a}{(s-a)^2+b^2}$ roots $s = a+jb, a-jb$ |
| *10 | $[x]_{\Delta x}^{n}\, e^{ax}$ or $[\prod_{m=1}^{n}(x-[m-1]\Delta x)]e^{ax}$ $n = 1,2,3, \ldots$ | $\frac{e^{na\Delta x}\, n!}{(s-\frac{e^{a\Delta x}-1}{\Delta x})^{n+1}}$ | $x^n e^{ax}$ | $\frac{n!}{(s-a)^{n+1}}$ |

* $K_{\Delta x}$ Transforms used in the solution of the differential difference equation of Example 4.3-4

Example 4.3-4 below is solved using $K_{\Delta x}$ Transform methodology. All necessary transforms for the solution of the given first order differential difference equation are obtained from $K_{\Delta x}$ Transform Table 2 and Table 3 which are found in the Appendix but, in part, are presented above. The convenience of using these tables is obvious. Note that after finding the roots of the polynomial multiplier of y(s), the differential difference equation complementary solution and particular solution are obtained from simple algebraic manipulations and transform conversions.

Example 4.3-4 Use the $K_{\Delta x}$ Transform Method to solve the following differential difference equation.

$$D_{\Delta x}y(x) + y(x) = xe^{-x}, \qquad \Delta x = 2,\ y(0) = 3. \qquad 1)$$

Apply $K_{\Delta x}$ Transforms from the $K_{\Delta x}$ Transform Tables in the Appendix to Eq 1

$$K_{\Delta x}[D_{\Delta x}y(x)] + K_{\Delta x}[y(x)] = K_{\Delta x}[xe^{-x}]$$

$$sy(s) - y(0) + y(s) = \frac{e^{a\Delta x}}{\left(s - \frac{e^{a\Delta x}-1}{\Delta x}\right)^2} \qquad 2)$$

where

$a = -1$

$\Delta x = 2$

$y(0) = 3$

Substituting these values into Eq 2

$$sy(s) - 3 + y(s) = \frac{e^{-2}}{\left(s - \frac{e^{-2}-1}{2}\right)^2} \qquad 3)$$

$e^{-2} = .135335283$

$\frac{e^{-2}-1}{2} = -.432442358$

Substituting these two values into Eq 3

$$(s + 1)y(s) = 3 + \frac{.135335283}{(s + .432332358)^2}$$

$$y(s) = \frac{3}{(s + 1)} + \frac{.135335283}{(s + 1)(s + .432332358)^2} \qquad 4)$$

Finding a partial fraction expansion of y(s)

$$\frac{.135335283}{(s + 1)(s + .432332358)^2} = \frac{A}{(s + 1)} + \frac{B}{(s + .432332358)^2} + \frac{C}{(s + .432332358)} \qquad 5)$$

$$A = \frac{.135335283}{(s + .432332358)^2}\bigg|_{s = -1} = \frac{.135335283}{(-1 + .432332358)^2} = .41997434$$

$$B = \frac{.135335283}{(s + 1)}\bigg|_{s = -.432332358} = \frac{.135335283}{(-.432332358 + 1)} = .238405843$$

$$C = \frac{d}{ds}\left[\frac{.135335283}{(s+1)}\right]\Big|_{s = -.432332358} = \left[-\frac{.135335283}{(s+1)^2}\right]\Big|_{s = -.432332358}$$

$$C = -\frac{.135335283}{(-.432332358 + 1)^2} = -.41997434$$

Substituting the values of A,B,C into Eq 5 and then substituting Eq 5 into Eq 4

$$y(s) = \frac{3}{s+1} + \frac{.41997434}{s+1} - \frac{.41997434}{s + .432332358} + \frac{.238405843}{(s + .432332358)^2} \qquad 6)$$

Find the inverse $K_{\Delta x}$ Transform of Eq 6

From the $K_{\Delta x}$ Transform Tables

$$K_{\Delta x}[e_{\Delta x}(a,x)] = \frac{1}{s-a}$$

$$K_{\Delta x}[xe_{\Delta x}(a,x)] = \frac{1 + a\Delta x}{(s-a)^2}$$

Putting Eq 6 into a form which will facilitate taking the inverse $K_{\Delta x}$ Transform

$a = -.43233235$
$\Delta x = 2$

$1+a\Delta x = 1 - .432332358(2) = .1353335284$

$$y(s) = \frac{3.41997434}{s+1} - \frac{.41997434}{s + .432332358} + 1.761594138\left[\frac{.135335284}{(s + .432332358)^2}\right] \qquad 7)$$

Taking the inverse $K_{\Delta x}$ Transform of Eq 7

$y(x) = 3.41997434e_2(-1,x) - .41997434e_2(-.432332358,x) + 1.761594138xe_2(-.432332358,x)$

**$y(x) = 3.4199743e_2(-1,x) + ( - .41997434 + 1.76159413x)e_2(-.43233235,x)$** 8)

Eq 8 is the solution to the differential difference equation:

$D_{\Delta x}y(x) + y(x) = xe^{-x}$
where
$\Delta x = 2$
$y(0) = 3.$
$x = 0,2,4,6,\ldots$

If desired, Eq 8 can be changed to another form

$$e_{\Delta x}(c,x) = (1+c\Delta x)^{\frac{x}{\Delta x}} \qquad 9)$$

$$e_2(-1,x) = [1+(-1)(2)]^{\frac{x}{2}} = [-1]^{\frac{x}{2}} \qquad 10)$$

From Eq 9

$$e_{\Delta x}(c,x) = \left[e^{\ln(1+c\Delta x)}\right]^{\frac{x}{\Delta x}}$$

Then

$$e_2(-.432332358,x) = \left[e^{\ln(1+\{-.432332358\}\{2\})}\right]^{\frac{x}{2}} = [\,e^{-2}\,]^{\frac{x}{2}} = e^{-x}$$

$$e_2(-.432332358,x) = e^{-x} \qquad 11)$$

Substituting Eq 10 and Eq 11 into Eq 8

$$\mathbf{y(x) = 3.4199743[-1]^{\frac{x}{2}} + ( - .41997434 + 1.76159413x)e^{-x}} \qquad 12)$$

Checking

Eq 8 and Eq 12 are the same as Eq 14 and Eq 15 respectively calculated in Example 4.2-1. Example 4.2-1 solved this same differential difference equation using the Method of Undetermined Coefficients.

Good check

There are additional methods for solving differential difference equations. The Method of Variation of Parameters is presented in the following section, Section 4.4.

## Section 4.4: The solution of differential difference equations using the Method of Variation of Parameters

The Method of Variation of Parameters used in Calculus to obtain the particular solution of a differential equation is also applicable to finding the particular solution of a differential difference equation. However, all discrete differentiation and integration must be performed in accordance with the previously developed operations of Interval Calculus discrete mathematics. The Method of Variation of Parameters for finding the particular solution of a differential difference equation is presented below.

Given an nth order differential difference equation to be solved, the Method of Variation of Parameters firstly requires that its complementary solution, $y_c(x)$, be obtained. Each initial condition constant of the complementary solution is then replaced by an undetermined function of the differential equation independent variable. The equation, $y_p(x)$, resulting from this modification to the complementary solution is now considered to be the particular solution of the differential difference equation. Of course, the undetermined functions, which are part of this equation, must be determined. The procedure for doing is as follows:

1. Differentiate $y_p(x)$ to obtain $D_{\Delta x}y_p(x)$. Within $D_{\Delta x}y_p(x)$ there will be terms which contain, $c_i'(x)$ (i = 1,2,3,…,n), the first derivatives of the n undetermined functions. Set this combination of terms to zero.
2. Differentiate $D_{\Delta x}y_p(x)$ to obtain $D_{\Delta x}^2(x)$. Again there will be terms which contain, $c_i'(x)$. Set this combination of terms to zero as before.
3. Continue the differentiation of $y_p(x)$ process up to and including $D_{\Delta x}^{n-1}y_p(x)$.
4. Differentiate $D_{\Delta x}^{n-1}y_p(x)$ to obtain $D_{\Delta x}^{n}y_p(x)$.
5. Substitute all $D_{\Delta x}^{m}y_p(x)$, m = 0,1,2,3,…n into the differential difference equation. The resulting equation will not contain any $c_i(x)$ terms. They will cancel out due to the manner in which $y_p(x)$ was derived from $y_c(x)$, the complementary solution.
6. The equation obtained from 5 and the equations equated to zero that were obtained from 1 to 3 yield a system of n linear equations with n unknowns, $c_i'(x)$, (i = 1,2,3,…,n). These equations may be solved to determine $c_i'(x)$, (i = 1,2,3,…,n).
7. Using discrete integration, obtain $c_i(x)$ from $c_i'(x)$ (i = 1,2,3,…,n) determined in 6.
8. Substitute the values of $c_i(x)$, (i = 1,2,3,…,n) obtained in 7 into $y_p(x)$.

The particular solution of the differential difference equation has now been obtained using the Method of Variation of Parameters.

The Method of Variation of Parameters gets progressively more difficult as the order, n, of the differential difference equation increases since the method requires the solution of a system n linear equations in n unknows.

The two examples, Example 4.4-1 and Example 4.4-2, which follow demonstrate the application of the Method of Variation of Parameters process described above.

Example 4.4-1 Use the Method of Variation of Parameters to solve the following differential difference equation:

$D_{\Delta x}y(x) + y(x) = xe^{-x}$, $\Delta x = 2$, $y(0) = 3$, $x = 0,2,4,6,\ldots$ 1)

Find the related homogeneous equation complementary solution, $y_c(x)$.

$D_{\Delta x}y_c(x) + y_c(x) = 0$ , the related homogeneous differential difference equation 2)

Let

$y_c(x) = ce_{\Delta x}(a,x) = c(1+a\Delta x)^{\frac{x}{\Delta x}}$ 3)

Substituting this $y_c(x)$ function into the related homogenous differential difference equation, Eq 2

$cae_{\Delta x}(a,x) + ce_{\Delta x}(a,x) = (a + 1)ce_{\Delta x}(a,x) = 0$

$a = -1$

Substituting $a = -1$ into Eq 3

The complementary solution to the differential difference equation, Eq 1 is then

$$y_c(x) = ce_2(-1,x) \qquad 4)$$

Find the particular solution to the differential difference equation, Eq 1.

The function, $e^{-x}$, should first be converted into its identity. This will facilitate the necessary discrete differentiation.

Using the relationship:

$$e^{cx} = e_{\Delta x}\left(\frac{e^{c\Delta x}-1}{\Delta x},x\right)$$

$c = -1, \quad \Delta x = 2$

$$e^{-x} = e_2\left(\frac{e^{-2}-1}{2},x\right) = e_2(-.43233235,x) \qquad 5)$$

Substituting into Eq 1

$$D_2y(x) + y(x) = xe_2(-.43233235,x) \qquad 6)$$

Note that all of the functions in Eq 6 are discrete Interval Calculus functions.
Consider the constant, c, in Eq 4 to be a function of x.

$$c = c(x) \qquad 7)$$

From Eq 4 and Eq 7, let the particular solution, $y_p(x)$, be as follows:

$$y_p(x) = c(x)e_2(-1,x) \qquad 8)$$

Take the discrete derivative of Eq 8.
Use the following derivative of the product of two functions formula.

$D_{\Delta x}[u(x)v(x)] = v(x)D_{\Delta x}u(x) + D_{\Delta x}v(x)u(x+\Delta x)$ where $u(x) = e_2(-1,x)$ , $v(x) = c(x)$

The following notation will be used to represent the discrete derivative.

$f^{'}(x) = D_{\Delta x}f(x)$

Differentiating Eq 8

$y_p'(x) = -c(x)e_2(-1,x) + c'(x)e_2(-1,x+2)$ 9)

Using the relationship:
$e_{\Delta x}(b,x+\Delta x) = (1+b\Delta x)e_{\Delta x}(b,x)$ 10)

$e_2(-1,x+2) = [1+(-1)(2)]e_2(-1,x) = -e_2(-1,x)$

Substituting into Eq 9

$y_p'(x) = -c(x)e_2(-1,x) - c'(x)e_2(-1,x)$ 11)

Substituting Eq 8 and Eq 11 into Eq 6

$-c(x)e_2(-1,x) - c'(x)e_2(-1,x) + c(x)e_2(-1,x) = xe_2(-.43233235,x)$

$c'(x)e_2(-1,x) = -xe_2(-.43233235,x)$

$$c'(x) = \frac{-xe_2(-.43233235,x)}{e_2(-1,x)}$$

$c'(x) = -x[e_2(-.43233235,x)e_2(-1,-x)]$ 12)

Using the relationship:
$e_{\Delta x}(a,x)e_{\Delta x}(b,-x) = e_{\Delta x}(\frac{a-b}{1+b\Delta x},x)$ 13)

$$e_2(-.43233235,x)e_2(-1,-x) = e_2(\frac{-.432332358-(-1)}{1+(-1)(2)},x)$$

$e_2(-.43233235,x)e_2(-1,-x) = e_2(-.567667642,x)$ 14)

Substituting Eq 14 into Eq 12

$c'(x) = -xe_2(-.567667642,x)$ 15)

Integrate Eq 15 to find c(x).
Use the following integration formula for the product of two functions.

${}_{\Delta x}\int v(x)D_{\Delta x}u(x)\Delta x = u(x)v(x) - {}_{\Delta x}\int D_{\Delta x}v(x)u(x+\Delta x)\Delta x + K$ 16)

Let

$v(x) = x$ $\quad D_{\Delta x}v(x) = 1$

$D_{\Delta x}u(x) = e_2(-.567667642,x)$ $\quad u(x) = \frac{e_2(-.567667642,x)}{-.567667642} = -1.761594155\, e_2(-.567667642,x)$

From Eq 15 and Eq 16

$$c(x) = -_{\Delta x}\int xe_2(-.567667642,x)\Delta x = 1.761594155xe_2(-.567667642,x) -$$
$$1.761594155\ _{\Delta x}\int e_2(-.567667642,x+2)\Delta x$$

Using the relationship:

$$e_{\Delta x}(a,x+\Delta x) = (1+a\Delta x)e_{\Delta x}(a,x) \qquad 17)$$

$$c(x) = 1.761594155xe_2(-.567667642,x) - (1.761594155)[1+(-.567667642)(2)]\ _{\Delta x}\int e_2(-.567667642,x)\Delta x$$

$$c(x) = 1.761594155xe_2(-.567667642,x) - \frac{(1.761594155)(-.135335284)}{-.567667642}\,e_2(-.567667642,x)$$

$$c(x) = 1.761594155xe_2(-.567667642,x) - .419974343\ e_2(-.567667642,x) \qquad 18)$$

Substituting Eq 18 into Eq 8

$$y_p(x) = (-.419974343 + 1.761594155x)[e_2(-.567667642,x)\ e_2(-1,x)] \qquad 19)$$

Using the relationship:

$$e_{\Delta x}(a,x)e_{\Delta x}(b,x) = e_{\Delta x}(a+b+ab\Delta x\ ,x) \qquad 20)$$

$$e_2(-.567667642,x)e_{\Delta x}(-1,x) = e_{\Delta x}(-.567667642 - 1 + [-.567667642][-1][2],\ x)$$

$$e_2(-.567667642,x)e_{\Delta x}(-1,x) = e_2(-.432332358\ ,x)$$

Substituting into Eq 19

$$y_p(x) = (-.419974343 + 1.761594155x)e_2(-.432332358\ ,x) \qquad 21)$$

$$y(x) = y_c(x) + y_p(x) \qquad 22)$$

Substituting Eq 4 and Eq 21 into Eq 22

$$y(x) = ce_2(-1,x) + (-.419974343 + 1.761594155x)e_2(-.432332358\ ,x) \qquad 23)$$

Find c from the initial condition, y(0) = 3

$$e_{\Delta x}(a,0) = 1$$

$$3 = c - .419974343$$

$$c = 3.419974343$$

Substituting this value for c into Eq 23

$$\mathbf{y(x) = 3.419974343e_2(-1,x) + (-.419974343 + 1.761594155x)e_2(-.43233235\,,x)} \qquad 24)$$

Eq 24 can be changed to another form

Using the relationship:

$$e_{\Delta x}(a,x) = [1+a\Delta x]^{\frac{x}{\Delta x}} \qquad 25)$$

$$e_2(-1,x) = [1+(-1)(2)]^{\frac{x}{2}}$$

$$e_2(-1,x) = [-1]^{\frac{x}{2}}$$

From Eq 5

$$e^{-x} = e_2(-.43233235,x)$$

Substituting into Eq 24

$$\mathbf{y(x) = 3.419974343[-1]^{\frac{x}{2}} + (-.419974343 + 1.761594155x)e^{-x}} \qquad 26)$$

Checking Eq 26

y(0) = 3.419974343-.419974343 = 3 Good check

Substitute Eq 26 into Eq 1

$$\frac{3.419974343[-1]^{\frac{x+2}{2}} + (-.419974343 + 1.761594155\{x+2\})e^{-(x+2)}}{2}$$

$$-\frac{3.419974343[-1]^{\frac{x}{2}} + (-.419974343 + 1.761594155x)e^{-x}}{2}$$

$$+\ 3.419974343[-1]^{\frac{x}{2}} + (-.419974343 + 1.761594155x)e^{-x} = xe^{-x}$$

$$+\frac{-3.419974343[-1]^{\frac{x}{2}}}{2} + \frac{-.419974343e^{-2}e^{-x}}{2} + \frac{1.761594155(2)e^{-2}e^{-x}}{2} + \frac{1.761594155e^{-2}xe^{-x}}{2}$$

$$-\frac{3.419974343[-1]^{\frac{x}{2}}}{2}+\frac{.419974343e^{-x}}{2}-\frac{1.761594155xe^{-x}}{2}$$

$$+\,3.419974343[-1]^{\frac{x}{2}}\ -.419974343e^{-x}+1.761594155xe^{-x}=xe^{-x}$$

$$[\frac{-.419974343e^{-2}}{2}+\frac{1.761594155(2)e^{-2}}{2}+\frac{.419974343}{2}-.419974343\ ]e^{-x}+$$

$$[\frac{1.761594155e^{-2}}{2}-\frac{1.761594155}{2}+1.761594155]xe^{-x}=xe^{-x}$$

$$0e^{-x}+1xe^{-x}=xe^{-x}$$

$xe^{-x} = xe^{-x}$ Good check

Example 4.4-2 Use the Method of Variation of Parameters to solve the following differential difference equation:

$$D_{\Delta x}{}^2\, y(x)+D_{\Delta x}y(x)=3,\qquad \Delta x=2,\ \ y(0)=1,\ \ y(2)=4,\ \ x=0,2,4,6,\ldots \tag{1}$$

Find the related homogeneous equation complementary solution, $y_c(x)$.

$$D_{\Delta x}{}^2 y_c(x)+D_{\Delta x}y_c(x)=0\ , \quad \text{the related homogeneous differential difference equation} \tag{2}$$

Let

$$y_c(x)=ce_{\Delta x}(a,x)=c(1+a\Delta x)^{\frac{x}{\Delta x}} \tag{3}$$

Substituting this $y_c(x)$ function into the related homogenous differential difference equation, Eq 2
$ca^2e_{\Delta x}(a,x) + cae_{\Delta x}(a,x) = (a^2 + a)ce_{\Delta x}(a,x) = a(a + 1)ce_{\Delta x}(a,x) = 0$

$a = 0,-1$

Substituting $a = 0,-1$ into Eq 3
The complementary solution to the differential difference equation, Eq 1 is then

$$y_c(x)=c_1e_2(-1,x)+c_2e_2(0,x)=c_1e_2(-1,x)+c_2$$

$$y_c(x)=c_1e_2(-1,x)+c_2 \tag{4}$$

Find the particular solution to the differential difference equation, Eq 1.

Consider the constants, $c_1,c_2$, in Eq 4 to be a function of x

$c_1 = c_1(x)$ 5)
$c_2 = c_2(x)$ 6)

From Eq 4, Eq 5, and Eq 6, let the particular solution, $y_p(x)$, be as follows:

$y_p(x) = c_1(x)e_2(-1,x) + c_2(x)$ 7)

Take the discrete derivative of Eq 7.
Use the following formula for the derivative of the product of two functions.

$D_{\Delta x}[u(x)v(x)] = v(x)D_{\Delta x}u(x) + D_{\Delta x}v(x)u(x+\Delta x)$ where $u(x) = e_2(-1,x)$ , $v(x) = c_1(x)$ 8)

The following notation will be used to represent the discrete derivative.

$f'(x) = D_{\Delta x}f(x)$

Differentiating Eq 7

$y_p'(x) = -c_1(x)e_2(-1,x) + c_1'(x)e_2(-1,x+2) + c_2'(x)$ 9)

Using the relationship:
$e_{\Delta x}(b,x+\Delta x) = (1+b\Delta x)e_{\Delta x}(b,x)$ 10)

$e_2(-1,x+2) = [1+(-1)(2)]e_2(-1,x) = -e_2(-1,x)$ 11)

Substituting into Eq 9

$y_p'(x) = -c_1(x)e_2(-1,x) - c_1'(x)e_2(-1,x) + c_2'(x)$ 12)

Let
$-c_1'(x)e_2(-1,x) + c_2'(x) = 0$ 13)
$y_p'(x) = -c_1(x)e_2(-1,x)$ 14)

Taking the derivative of Eq 13 using Eq 8

$y_p''(x) = c_1(x)e_2(-1,x) - c_1'(x)e_2(-1,x+2)$ 15)

Substituting Eq 11 into Eq 15

$y_p''(x) = c_1(x)e_2(-1,x) + c_1'(x)e_2(-1,x)$ 16)

Substituting Eq 14 and Eq 16 into Eq 1

$c_1(x)e_2(-1,x) + c_1^{'}(x)e_2(-1,x) - c_1(x)e_2(-1,x) = 3$

$c_1^{'}(x)e_2(-1,x) = 3$ 17)

$c_1^{'}(x) = \frac{3}{e_2(-1,x)} = 3e_2(-1,-x)$

$c_1^{'}(x) = 3e_2(-1,-x)$ 18)

Integrating Eq 18

$c_1(x) = 3 {}_{2}\int e_2(-1,-x)\Delta x$ 19)

Using the relationship:

${}_{\Delta x}\int e_{\Delta x}(a,-x)\Delta x = -\frac{1+a\Delta x}{a}e_{\Delta x}(a,-x)$ 20)

$c_1(x) = 3[-\frac{1+(-1)(2)}{-1}]\, e_2(-1,-x) = -3e_2(-1,-x)$

$c_1(x) = -3e_2(-1,-x)$ 21)

From Eq 13

$c_2^{'}(x) = c_1^{'}(x)e_2(-1,x)$ 22)

Substituting Eq 18 into Eq 22

$c_2^{'}(x) = 3e_2(-1,-x)e_2(-1,x) = 3$

$c_2^{'}(x) = 3$

Integrating

$c_2(x) = 3x$ 23)

Substituting Eq 21 and Eq 23 into Eq 7

$y_p(x) = -3e_2(-1,-x)e_2(-1,x) + 3x$

$y_p(x) = 3x - 3$ 24)

$y(x) = y_c(x) + y_p(x)$ 25)

From Eq 4, Eq 24, and Eq 25)

$$y(x) = c_1 e_2(-1,x) + c_2 + 3x - 3 \tag{26}$$

Let
$K_1 = c_1$
$K_2 = c_2 - 3$

Substituting into Eq 26

$$y(x) = K_1 e_2(-1,x) + K_2 + 3x \tag{27}$$

Eq 26 can be put into another form

$$e_2(-1,x) = [1+(-1)(2)]^{\frac{x}{2}} = [-1]^{\frac{x}{2}}$$

Substituting into Eq 27

$$y(x) = K_1[-1]^{\frac{x}{2}} + K_2 + 3x \tag{28}$$

Evaluate $K_1$ and $K_2$ from the initial conditions, $y(0) = 1$ and $y(2) = 4$

$$y(0) = 1 = K_1[-1]^0 + K_2 + 0$$

$$K_1 + K_2 = 1 \tag{29}$$

$$y(2) = 4 = K_1[-1]^1 + K_2 + 3(2)$$

$$-K_1 + K_2 = -2 \tag{30}$$

Adding Eq 29 and Eq 30

$$2K_2 = -1$$

$$K_2 = -\frac{1}{2}$$

$$K_1 = 1 - \left(-\frac{1}{2}\right) = \frac{3}{2}$$

$$K_1 = \frac{3}{2}$$

Substituting the values of $K_1$ and $K_2$ into Eq 27 and Eq 28

$$\mathbf{y(x) = \frac{3}{2}\, e_2(-1,x) - \frac{1}{2} + 3x} \qquad 31)$$

**or**

$$\mathbf{y(x) = \frac{3}{2}\, [-1]^{\frac{x}{2}} - \frac{1}{2} + 3x} \qquad 32)$$

Checking Eq 31 and Eq 32

$y(0) = \frac{3}{2} - \frac{1}{2} = 1$ good check

$y(2) = \frac{3}{2}\, [-1]^{\frac{2}{2}} - \frac{1}{2} + 3(2) = -\frac{3}{2} - \frac{1}{2} + 6 = 4$

$y(2) = 4$ Good check

Substituting Eq 31 into Eq 1

$\frac{3}{2}\,(-1)^2\, e_2(-1,x) + (-1)\,\frac{3}{2}\, e_2(-1,x) + 3 = 3$

$3 = 3$ Good check

There is one more method for solving differential difference equations to be described, the Method of Related Functions. It appears in the following section, Section 4.5.

## Section 4.5: The solution of differential difference equations using the Method of Related Functions

If one reviews the $K_{\Delta x}$ Transforms and Laplace Transforms of Table 2 and Table 3 in the Appendix, some things which are very interesting (and perhaps perplexing) are observed. The $K_{\Delta x}$ Transforms of some discrete functions , functions of $\Delta x$ and x, are found to be the same as the Laplace Transforms of other functions, functions only of x. This, of course, could occur by coincidence if $K_{\Delta x}$ Transforms and Laplace Transforms were totally unrelated. However, there is a relationship. The two functions that share the common transform become identical if $\Delta x$ goes to zero. A logical conclusion is that Laplace Transforms are actually $K_{\Delta x}$ Transforms where $\Delta x$ approaches zero. This, in fact, is the case since the $K_{\Delta x}$ Transform integral becomes the Laplace Transform integral as $\Delta x$ approaches zero.

$$\lim_{\Delta x\to 0}\ {}_{\Delta x}\int_0^{\infty} (1+s\Delta x)^{-(\frac{x+\Delta x}{\Delta x})} f(x)\Delta x = \int_0^{\infty} e^{-sx} f(x)dx \,, \quad \text{Transform equations} \qquad (4.5\text{-}1)$$

where

$$\lim_{\Delta x\to 0}(1+s\Delta x)^{-\frac{x}{\Delta x}} = e^{-sx}$$

$\Delta x$ = x increment

Having established that the $K_{\Delta x}$ Transform and the Laplace Transform are related, another perplexing observation is made. Many of the common discrete functions, functions of $\Delta x$ and x, not only have the same transform as a Calculus function, but also have the same transform irrespective of the value of $\Delta x$. These functions are called the Primary Interval Calculus Functions. See the following function comparison table, Table 4.5-1. The functions, $e_{\Delta x}(a,x)$ and $e^{ax}$ with a transform of $\frac{1}{s-a}$ are an example of a Primary Interval Calculus Function pair. The Primary Interval Calculus Function pairs that have been identified are shown at the beginning of Table 4.5-1. In addition to these functions there are others that have similar but not identical transforms, for example, the discrete $e^{ax}$ function with its $K_{\Delta x}$ Transform, $\frac{1}{s-\frac{e^{a\Delta x}-1}{\Delta x}}$ , and the Calculus $e^{ax}$ function with its Laplace Transform, $\frac{1}{s-a}$ . However, if we compare the transforms of the discrete $e^{ax}$ function to the Calculus function, $e^{\frac{e^{a\Delta x}-1}{\Delta x}x}$ , they are found to be identical, $\frac{1}{s-\frac{e^{a\Delta x}-1}{\Delta x}}$ . Functions such as this have been identified and added below the Primary Interval Calculus Functions in Table 4.5-1.

Table 4.5-1 is a table of function pairs, one function a discrete calculus function and the other function a Calculus function. Both functions have the exact same transform. Besides having the same transform, it will be noted that the paired functions that share the common transform become identical when $\Delta x$ goes to zero. The function pairs listed in Table 4.5-1 are referred to as Related Functions. Pertinent facts concerning these listed functions and their transforms have been included in Table 4.5-1.

It is interesting that these function relationships exist, but can they be put to use? They can. Their potential use is somewhat surprising and perhaps contrary to logic. From these relationships it appears that, at least for some functions, discrete differential difference equations can be solved using Calculus (which typically is not considered to be discrete mathematics). The methodology to use Calculus to solve discrete differential difference equations is the topic of this section, Section 4.5. In this section, the solution method, specified as the Method of Related Equations, will be presented.

In Table 4.5-1, each identified discrete Interval Calculus function and Calculus function pair is listed together with its common transform. This table was developed as part of an initial investigation to identify any similarities between Inteval Calculus and Calculus. After much observation and study, a result not initially anticipated was obtained. Calculus was found to be a special case of discrete Interval Calculus. This result implies that Calculus is not a continuous mathematics as previously thought but, in actually, a discrete mathematics. True, with the independent variable interval so small ($\Delta x \to 0$), visually it is continuous, but not mathematically. The independent variable inteval is some value, it can not be zero ($\Delta x \neq 0$). In the development of Calculus, the independent variable interval, $\Delta x$, is so small relative to all other terms of a mathematical expression that it need not be shown. In most cases, this is not unreasonable. However, in some few cases, it must be remembered that although $\Delta x$ is not shown, it is, nevertheless, there. This mathematical analysis has been one of those few cases. Table 4.5-1 is rewritten in a more useful form in Table 4.5-2. In Table 4.5-2 Calculus functions are presented as discrete Interval Calculus functions where their $\Delta x$ value is infinitesimal (i.e. $\Delta x \to 0$). Table 4.5-2 directly follows Table 4.5-1.

# Discrete Calculus/Calculus Related Functions Table 4.5-1

<u>Comments</u> – 1) $L[f(x)] = \lim_{\Delta x \to 0} K_{\Delta x}[f(x)]$ , $L[f(x)]$ is the Laplace Transform, $K_{\Delta x}[f(x)]$ is the $K_{\Delta x}$ Transform

2) $f_{\Delta x}(x)$ are discrete calculus functions, $x = x_0 + m\Delta x$ , m = integers

3) f(x) are Calculus functions, $x = x_0 + m\Delta x$ , m = integers, $\Delta x \to 0$

4) $\lim_{\Delta x \to 0} f_{\Delta x}(x) = \lim_{\Delta x \to 0} f(x)$. For the Primary Interval Calculus Functions, $\lim_{\Delta x \to 0} f_{\Delta x}(x) = f(x)$.

5) If $L[f(x)]$ is the same as $K_{\Delta x}[f_{\Delta x}(x)]$, the functions, f(x) and $f_{\Delta x}(x)$ are said to be related.

6) The derivative operator $\frac{d^n}{dx^n} \equiv \lim_{\Delta x \to 0} D^n{}_{\Delta x}$ where n = 1,2,3,…

7) The Calculus differential equation, $\frac{d^n}{dx^n} y(x) + A_{n-1}\frac{d^{n-1}}{dx^{n-1}} y(x) + A_{n-2}\frac{d^{n-2}}{dx^{n-2}} y(x) + \ldots + A_1\frac{d}{dx} y(x) + A_0 y(x)$, is related to the discrete calculus equation, $D_{\Delta x}{}^n y(x) + A_{n-1}D_{\Delta x}{}^{n-1}y(x) + A_{n-2}D_{\Delta x}{}^{n-2}y(x) + \ldots + A_1 D_{\Delta x}y(x) + A_0 y(x)$, because their Laplace and $K_{\Delta x}$ Transforms are the same.

8) The $K_{\Delta x}$ Transform can be considered to be a Generalized Laplace Transform or the Laplace Transform can be considered to be the $K_{\Delta x}$ Transform where $\Delta x \to 0$.

**Table 4.5-1 A Listing of Discrete Interval Calculus and Calculus Functions having a common $K_{\Delta x}$/Laplace Transform**

$A_m, B_m, a, b, c, x_0$ are constants

| # | f(x) **Calculus Functions** $\lim_{\Delta x\to 0} f(x) = \lim_{\Delta x\to 0} f_{\Delta x}(x)$, $\Delta x \neq 0$; $x = x_0 + m\Delta x$, m = integers, $\Delta x \to 0$ | $f_{\Delta x}(x)$ **Discrete Calculus (Interval Calculus) Functions** | $f_{\Delta x}(x)$ **Calculation Equations** $x = x_0 + m\Delta x$, m = integers, $\Delta x \neq 0$ | $K_{\Delta x}[f_{\Delta x}(x)]$ L[f(x)] |
|---|---|---|---|---|
| **1** | $\frac{d^n}{dx^n} y(x) + A_{n-1}\frac{d^{n-1}}{dx^{n-1}} y(x) + A_{n-2}\frac{d^{n-2}}{dx^{n-2}} y(x) + \ldots + A_1\frac{d}{dx} y(x) + A_0\, y(x)$ | $D_{\Delta x}{}^n y(x) + A_{n-1}D_{\Delta x}{}^{n-1}y(x) + A_{n-2}D_{\Delta x}{}^{n-2}y(x) + \ldots + A_1 D_{\Delta x}y(x) + A_0 y(x)$ | --- | $s^n y(s) + A_{n-1}s^{n-1}y(s) + A_{n-2}s^{n-2}y(s) + \ldots + A_1 s y(s) + A_0 y(s) - B_{n-1}s^{n-1} - B_{n-2}s^{n-2} - B_{n-3}s^{n-3} - \ldots - B_1 s - B_0$<br>$B_m$, m=0,1,2,…,n-1 are intial condition constants |
| **2** | $e^{ax}$ | $e_{\Delta x}(a,x)$ | $(1+a\Delta x)^{\frac{x}{\Delta x}}$<br>or<br>$e^{[\frac{\ln(1+a\Delta x)}{\Delta x}]x}$ | $\frac{1}{s-a}$<br>root s = a |
| **3** | sinbx | $\sin_{\Delta x}(b,x)$ | $\frac{(1+jb\Delta x)^{\frac{x}{\Delta x}} - (1-jb\Delta x)^{\frac{x}{\Delta x}}}{2j}$<br>or<br>$(1+[b\Delta x]^2)^{\frac{x}{2\Delta x}} \sin\beta\frac{x}{\Delta x}$<br>or | $\frac{b}{s^2+b^2}$<br>roots s = jb, -jb |

| # | f(x)<br>**Calculus Functions**<br>$\lim_{\Delta x \to 0} f(x) = \lim_{\Delta x \to 0} f_{\Delta x}(x)$, $\Delta x \neq 0$<br>$x = x_0 + m\Delta x$, m = integers, $\Delta x \to 0$ | $f_{\Delta x}(x)$<br>**Discrete Calculus (Interval Calculus) Functions** | $f_{\Delta x}(x)$<br>**Calculation Equations**<br>$x = x_0 + m\Delta x$, m = integers, $\Delta x \neq 0$ | $K_{\Delta x}[f_{\Delta x}(x)]$<br>$L[f(x)]$ |
|---|---|---|---|---|
| | | | $(\sec\beta)^{\frac{x}{\Delta x}} \sin\beta\frac{x}{\Delta x}$<br><br>$\beta = \tan^{-1} b\Delta x$ | |
| **4** | cosbx | $\cos_{\Delta x}(b,x)$ | $\dfrac{(1+jb\Delta x)^{\frac{x}{\Delta x}} + (1-jb\Delta x)^{\frac{x}{\Delta x}}}{2}$<br>or<br>$(1+[b\Delta x]^2)^{\frac{x}{2\Delta x}} \cos\beta\frac{x}{\Delta x}$<br>or<br>$(\sec\beta)^{\frac{x}{\Delta x}} \cos\beta\frac{x}{\Delta x}$<br><br>$\beta = \tan^{-1} b\Delta x$ | $\dfrac{s}{s^2+b^2}$<br>roots $s = jb, -jb$ |

| # | f(x)<br>Calculus Functions<br>$\lim_{\Delta x \to 0} f(x) = \lim_{\Delta x \to 0} f_{\Delta x}(x)$, $\Delta x \neq 0$<br>$x = x_0 + m\Delta x$ , m = integers , $\Delta x \to 0$ | $f_{\Delta x}(x)$<br>Discrete Calculus (Interval Calculus) Functions | $f_{\Delta x}(x)$<br>Calculation Equations<br>$x = x_0 + m\Delta x$ , m = integers , $\Delta x \neq 0$ | $K_{\Delta x}[f_{\Delta x}(x)]$<br>L[f(x)] |
|---|---|---|---|---|
| **5** | $e^{ax}\sin bx$ | <u>$1+a\Delta x \neq 0$</u> $(a \neq -\frac{1}{\Delta x})$<br><br>$e_{\Delta x}(a,x)\sin_{\Delta x}(\frac{b}{1+a\Delta x}, x)$<br><br>or<br><br>$(1+a\Delta x)^{\frac{x}{\Delta x}}\sin_{\Delta x}(\frac{b}{1+a\Delta x}, x)$<br><br>or<br><br><u>$1+a\Delta x = 0$</u> $(a = -\frac{1}{\Delta x})$<br><br>$[b\Delta x]^{\frac{x}{\Delta x}} \sin\frac{\pi x}{2\Delta x}$ | <u>$1+a\Delta x \neq 0$</u> $(a \neq -\frac{1}{\Delta x})$<br><br>$[(1+a\Delta x)^{\frac{x}{\Delta x}}]*$<br>$[\frac{(1+j\frac{b\Delta x}{1+a\Delta x})^{\frac{x}{\Delta x}} - (1-j\frac{b\Delta x}{1+a\Delta x})^{\frac{x}{\Delta x}}}{2j}]$<br>or<br>$[\sqrt{(1+a\Delta x)^2 + (b\Delta x)^2}\ ]^{\frac{x}{\Delta x}} \sin\frac{\beta x}{\Delta x}$<br>or<br>$(1+a\Delta x)^{\frac{x}{\Delta x}}(1+[\frac{b\Delta x}{1+a\Delta x}]^2)^{\frac{x}{2\Delta x}} \sin\frac{\beta x}{\Delta x}$<br><br>$\beta = \begin{cases} \tan^{-1}\frac{b\Delta x}{1+a\Delta x} & \text{for } 1+a\Delta x > 0 \\ \pi+\tan^{-1}\frac{b\Delta x}{1+a\Delta x} & \text{for } 1+a\Delta x < 0 \end{cases}$<br><br><u>$1+a\Delta x = 0$</u> $(a = -\frac{1}{\Delta x})$<br><br>$[b\Delta x]^{\frac{x}{\Delta x}} \sin\frac{\pi x}{2\Delta x}$ | $\frac{b}{(s-a)^2+b^2}$<br>roots $s = a+jb, a-jb$ |

| # | $f(x)$ **Calculus Functions** $\lim_{\Delta x\to 0} f(x)=\lim_{\Delta x\to 0} f_{\Delta x}(x)$, $\Delta x \neq 0$ <br> $x = x_0 + m\Delta x$, m = integers, $\Delta x \to 0$ | $f_{\Delta x}(x)$ **Discrete Calculus (Interval Calculus) Functions** | $f_{\Delta x}(x)$ **Calculation Equations** <br> $x = x_0 + m\Delta x$, m = integers, $\Delta x \neq 0$ | $K_{\Delta x}[f_{\Delta x}(x)]$ <br> $L[f(x)]$ |
|---|---|---|---|---|
| **6** | $e^{ax}\cos bx$ | $\underline{1+a\Delta x \neq 0}$ $\quad (a \neq -\frac{1}{\Delta x})$ <br><br> $e_{\Delta x}(a,x)\cos_{\Delta x}(\frac{b}{1+a\Delta x}, x)$ <br><br> or <br><br> $(1+a\Delta x)^{\frac{x}{\Delta x}}\cos_{\Delta x}(\frac{b}{1+a\Delta x}, x)$ <br><br> or <br><br> $\underline{1+a\Delta x = 0}$ $\quad (a = -\frac{1}{\Delta x})$ <br><br> $[b\Delta x]^{\frac{x}{\Delta x}} \cos\frac{\pi x}{2\Delta x}$ | $\underline{1+a\Delta x \neq 0}$ $\quad (a \neq -\frac{1}{\Delta x})$ <br><br> $[(1+a\Delta x)^{\frac{x}{\Delta x}}]*$ <br> $[\frac{(1+j\frac{b\Delta x}{1+a\Delta x})^{\frac{x}{\Delta x}} + (1-j\frac{b\Delta x}{1+a\Delta x})^{\frac{x}{\Delta x}}}{2}]$ <br><br> or <br><br> $[\sqrt{(1+a\Delta x)^2 + (b\Delta x)^2}\ ]^{\frac{x}{\Delta x}} \cos\frac{\beta x}{\Delta x}$ <br><br> or <br><br> $(1+a\Delta x)^{\frac{x}{\Delta x}} (1+[\frac{b\Delta x}{1+a\Delta x}]^2)^{\frac{x}{2\Delta x}} \cos\frac{\beta x}{\Delta x}$ <br><br> $\beta = \begin{cases} \tan^{-1}\frac{b\Delta x}{1+a\Delta x} & \text{for } 1+a\Delta x > 0 \\ \pi+\tan^{-1}\frac{b\Delta x}{1+a\Delta x} & \text{for } 1+a\Delta x < 0 \end{cases}$ <br><br> $\underline{1+a\Delta x = 0}$ $\quad (a = -\frac{1}{\Delta x})$ <br><br> $[b\Delta x]^{\frac{x}{\Delta x}} \cos\frac{\pi x}{2\Delta x}$ | $\frac{s-a}{(s-a)^2+b^2}$ <br> roots $s = a+jb,\ a-jb$ |

| # | f(x) Calculus Functions $\lim_{\Delta x \to 0} f(x) = \lim_{\Delta x \to 0} f_{\Delta x}(x)$, $\Delta x \neq 0$ $x = x_0 + m\Delta x$, m = integers, $\Delta x \to 0$ | $f_{\Delta x}(x)$ Discrete Calculus (Interval Calculus) Functions | $f_{\Delta x}(x)$ Calculation Equations $x = x_0 + m\Delta x$, m = integers, $\Delta x \neq 0$ | $K_{\Delta x}[f_{\Delta x}(x)]$ L[f(x)] |
|---|---|---|---|---|
| **7** | sinhbx | $\sinh_{\Delta x}(b,x)$ | $\frac{(1+b\Delta x)^{\frac{x}{\Delta x}} - (1-b\Delta x)^{\frac{x}{\Delta x}}}{2}$ or $(1-[b\Delta x]^2)^{\frac{x}{2\Delta x}} \sinh\frac{\beta x}{\Delta x}$ $\beta = \tanh^{-1} b\Delta x$ | $\frac{b}{s^2-b^2}$ roots s = b, -b |
| **8** | coshbx | $\cosh_{\Delta x}(b,x)$ | $\frac{(1+b\Delta x)^{\frac{x}{\Delta x}} + (1-b\Delta x)^{\frac{x}{\Delta x}}}{2}$ or $(1-[b\Delta x]^2)^{\frac{x}{2\Delta x}} \cosh\frac{\beta x}{\Delta x}$ $\beta = \tanh^{-1} b\Delta x$ | $\frac{s}{s^2-b^2}$ roots s = b, -b |

| # | f(x) Calculus Functions $\lim_{\Delta x \to 0} f(x) = \lim_{\Delta x \to 0} f_{\Delta x}(x)$, $\Delta x \neq 0$ $x = x_0 + m\Delta x$, m = integers, $\Delta x \to 0$ | $f_{\Delta x}(x)$ Discrete Calculus (Interval Calculus) Functions | $f_{\Delta x}(x)$ Calculation Equations $x = x_0 + m\Delta x$, m = integers, $\Delta x \neq 0$ | $K_{\Delta x}[f_{\Delta x}(x)]$ L[f(x)] |
|---|---|---|---|---|
| **9** | $x^n e^{ax}$ <br> n = 1,2,3, … | $[x]_{\Delta x}^{n} e_{\Delta x}(a,x-n\Delta x)$ <br> or <br> $[x]_{\Delta x}^{n} e_{\Delta x}(a,-n\Delta x)e_{\Delta x}(a,x)$ <br> or <br> $[x]_{\Delta x}^{n}(1+a\Delta x)^{-n}e_{\Delta x}(a,x)$ | $[\prod_{m=1}^{n}(x-[m-1]\Delta x)][(1+a\Delta x)^{\frac{x-n\Delta x}{\Delta x}}]$ | $\frac{n!}{(s-a)^{n+1}}$ |
| **10** | xsinbx | $x\sin_{\Delta x}(b,x-\Delta x)$ <br> or <br> $\frac{-[b\Delta x]x\cos_{\Delta x}(b,x)+x\sin_{\Delta x}(b,x)}{1+[b\Delta x]^2}$ | To calculate $\sin_{\Delta x}(b,x)$ and $\cos_{\Delta x}(b,x)$ see rows 3 and 4 | $\frac{2bs}{(s^2+b^2)^2}$ |
| **11** | xcosbx | $x\cos_{\Delta x}(b,x-\Delta x)$ <br> or <br> $\frac{x\cos_{\Delta x}(b,x)+[b\Delta x]x\sin_{\Delta x}(b,x)}{1+[b\Delta x]^2}$ | To calculate $\sin_{\Delta x}(b,x)$ and $\cos_{\Delta x}(b,x)$ see rows 3 and 4 | $\frac{s^2-b^2}{(s^2+b^2)^2}$ |
| **12** | xsinbx+[bΔx]xcosbx | $x\sin_{\Delta x}(b,x)$ | To calculate $\sin_{\Delta x}(b,x)$ see row 3 | $\frac{2bs+[b\Delta x](s^2-b^2)}{(s^2+b^2)^2}$ |
| **13** | xcosbx–[bΔx]xsinbx | $x\cos_{\Delta x}(b,x)$ | To calculate $\cos_{\Delta x}(b,x)$ see row 4 | $\frac{(s^2-b^2)-[b\Delta x]2bs}{(s^2+b^2)^2}$ |

| # | f(x) Calculus Functions $\lim_{\Delta x\to 0} f(x) = \lim_{\Delta x\to 0} f_{\Delta x}(x)$, $\Delta x \neq 0$ $x = x_0 + m\Delta x$, m = integers, $\Delta x \to 0$ | $f_{\Delta x}(x)$ Discrete Calculus (Interval Calculus) Functions | $f_{\Delta x}(x)$ Calculation Equations $x = x_0 + m\Delta x$, m = integers, $\Delta x \neq 0$ | $K_{\Delta x}[f_{\Delta x}(x)]$ L[f(x)] |
|---|---|---|---|---|
| **14** | $x^n \sin bx$ | $[x]_{\Delta x}^{n}\ \sin_{\Delta x}(b,x-n\Delta x)$ or $[x]_{\Delta x}^{n}\cos_{\Delta x}(b,-n\Delta x)\sin_{\Delta x}(b,x)+$ $[x]_{\Delta x}^{n}\sin_{\Delta x}(b,-n\Delta x)\cos_{\Delta x}(b,x)$ | $[\prod_{m=1}^{n}(x-[m-1]\Delta x)] *$ $[\frac{(1+jb\Delta x)^{\frac{x-n\Delta x}{\Delta x}} - (1-jb\Delta x)^{\frac{x-n\Delta x}{\Delta x}}}{2j}]$ | $\frac{n!}{2j}[\frac{(s+jb)^{n+1} - (s-jb)^{n+1}}{(s^2+b^2)^{n+1}}]$ |
| **15** | $x^n \cos bx$ | $[x]_{\Delta x}^{n}\ \cos_{\Delta x}(b,x-n\Delta x)$ or $[x]_{\Delta x}^{n}\cos_{\Delta x}(b,-n\Delta x)\cos_{\Delta x}(b,x)-$ $[x]_{\Delta x}^{n}\sin_{\Delta x}(b,-n\Delta x)\sin_{\Delta x}(b,x)$ | $[\prod_{m=1}^{n}(x-[m-1]\Delta x)]$ $[\frac{(1+jb\Delta x)^{\frac{x-n\Delta x}{\Delta x}} + (1-jb\Delta x)^{\frac{x-n\Delta x}{\Delta x}}}{2}]$ | $\frac{n!}{2}[\frac{(s+jb)^{n+1} + (s-jb)^{n+1}}{(s^2+b^2)^{n+1}}]$ |
| **16** | $x\sinh bx$ | $x\sinh_{\Delta x}(b,x-\Delta x)$ or $\frac{-[b\Delta x]x\cosh_{\Delta x}(b,x)+x\sinh_{\Delta x}(b,x)}{1-[b\Delta x]^2}$ | To calculate $\sinh_{\Delta x}(b,x)$ and $\cosh_{\Delta x}(b,x)$ see rows 7 and 8 | $\frac{2bs}{(s^2-b^2)^2}$ |
| **17** | $x\cosh bx$ | $x\cosh_{\Delta x}(b,x-\Delta x)$ or $\frac{x\cosh_{\Delta x}(b,x)-[b\Delta x]x\sinh_{\Delta x}(b,x)}{1-[b\Delta x]^2}$ | To calculate $\sinh_{\Delta x}(b,x)$ and $\cosh_{\Delta x}(b,x)$ see rows 7 and 8 | $\frac{s^2+b^2}{(s^2-b^2)^2}$ |

| # | $f(x)$<br>**Calculus Functions**<br>$\lim_{\Delta x\to 0} f(x)=\lim_{\Delta x\to 0} f_{\Delta x}(x)$, $\Delta x \neq 0$<br>$x = x_0 + m\Delta x$ , m = integers , $\Delta x \to 0$ | $f_{\Delta x}(x)$<br>**Discrete Calculus (Interval Calculus) Functions** | $f_{\Delta x}(x)$<br>**Calculation Equations**<br>$x = x_0 + m\Delta x$ , m = integers , $\Delta x \neq 0$ | $K_{\Delta x}[f_{\Delta x}(x)]$<br>$L[f(x)]$ |
|---|---|---|---|---|
| **18** | $x\sinh bx+[b\Delta x]x\cosh bx$ | $x\sinh_{\Delta x}(b,x)$ | To calculate $\sinh_{\Delta x}(b,x)$ see row 7 | $\frac{2bs+[b\Delta x](s^2+b^2)}{(s^2-b^2)^2}$ |
| **19** | $x\cosh bx+[b\Delta x]x\sinh bx$ | $x\cosh_{\Delta x}(b,x)$ | To calculate $\cosh_{\Delta x}(b,x)$ see row 8 | $\frac{(s^2+b^2)+[b\Delta x]2bs}{(s^2-b^2)^2}$ |
| **20** | $x^n$<br>n = 1,2,3, … | $[x]_{\Delta x}^{n}$ | $\prod_{m=1}^{n}(x-[m-1]\Delta x)$ | $\frac{n!}{s^{n+1}}$ |
| **21** | $x^n e^{\frac{e^{a\Delta x}-1}{\Delta x}x}$<br>n = 1,2,3, … | $[x]_{\Delta x}^{n}\, e^{a(x-n\Delta x)}$ | $[\prod_{m=1}^{n}(x-[m-1]\Delta x)]e^{a(x-n\Delta x)}$ | $\frac{n!}{(s-\frac{e^{a\Delta x}-1}{\Delta x})^{n+1}}$ |
| **22** | c | c | c | $\frac{c}{s}$ |
| **23** | x | $[x]_{\Delta x}^{1}$ | x | $\frac{1}{s^2}$ |
| **24** | $e^{\frac{e^{a\Delta x}-1}{\Delta x}x}$ | $e^{ax}$ | $e^{ax}$<br>or<br>$e_{\Delta x}(\frac{e^{a\Delta x}-1}{\Delta x}, x)$<br>an equivalent function for $e^{ax}$ | $\frac{1}{s-\frac{e^{a\Delta x}-1}{\Delta x}}$ |

| # | f(x) **Calculus Functions** $\lim_{\Delta x \to 0} f(x) = \lim_{\Delta x \to 0} f_{\Delta x}(x)$, $\Delta x \neq 0$ $x = x_0 + m\Delta x$, m = integers, $\Delta x \to 0$ | $f_{\Delta x}(x)$ **Discrete Calculus (Interval Calculus) Functions** | $f_{\Delta x}(x)$ **Calculation Equations** $x = x_0 + m\Delta x$, m = integers, $\Delta x \neq 0$ | $K_{\Delta x}[f_{\Delta x}(x)]$ $L[f(x)]$ |
|---|---|---|---|---|
| **25** | $e^{\frac{A^{\Delta x}-1}{\Delta x}x}$ | $A^x$ | $A^x$ <br> or <br> $e_{\Delta x}(\frac{A^{\Delta x}-1}{\Delta x}, x)$ <br> an equivalent function for $A^x$ | $\frac{1}{s - \frac{A^{\Delta x}-1}{\Delta x}}$ |
| **26** | $e^{\frac{\cos b\Delta x - 1}{\Delta x}x}\sin(\frac{\sin b\Delta x}{\Delta x}x)$ | sinbx | sinbx <br> or <br> $e_{\Delta x}(\frac{\cos b\Delta x\ \text{-}1}{\Delta x}, x)\ \sin_{\Delta x}(\frac{\tan b\Delta x}{\Delta x}, x)$ <br> where <br> $\frac{x}{\Delta x}$ = integer <br> $\cos(b\Delta x) \neq 0$ <br> an equivalent function for sinbx except where $\cos b\Delta x = 0$ | $\frac{\frac{\sin b\Delta x}{\Delta x}}{(s - \frac{\cos b\Delta x - 1}{\Delta x})^2 + (\frac{\sin b\Delta x}{\Delta x})^2}$ |

| # | f(x) Calculus Functions $\lim_{\Delta x\to 0} f(x) = \lim_{\Delta x\to 0} f_{\Delta x}(x)$, $\Delta x \neq 0$ $x = x_0 + m\Delta x$, m = integers, $\Delta x \to 0$ | $f_{\Delta x}(x)$ Discrete Calculus (Interval Calculus) Functions | $f_{\Delta x}(x)$ Calculation Equations $x = x_0 + m\Delta x$, m = integers, $\Delta x \neq 0$ | $K_{\Delta x}[f_{\Delta x}(x)]$ L[f(x)] |
|---|---|---|---|---|
| **27** | $e^{\frac{\cos b\Delta x-1}{\Delta x}x}\cos(\frac{\sin b\Delta x}{\Delta x}x)$ | cosbx | cosbx<br>or<br>$e_{\Delta x}(\frac{\cos b\Delta x - 1}{\Delta x}, x)\cos_{\Delta x}(\frac{\tan b\Delta x}{\Delta x}, x)$<br>where<br>$\frac{x}{\Delta x}$ = integer<br>$\cos b\Delta x \neq 0$<br>an equivalent function for cosbx except where $\cos b\Delta x = 0$ | $\frac{s - \frac{\cos b\Delta x-1}{\Delta x}}{(s - \frac{\cos b\Delta x-1}{\Delta x})^2 + (\frac{\sin b\Delta x}{\Delta x})^2}$ |
| **28** | $e^{\frac{e^{a\Delta x}\cos b\Delta x-1}{\Delta x}x}\sin(\frac{e^{a\Delta x}\sin b\Delta x}{\Delta x}x)$ | $e^{ax}\sin bx$ | $e^{ax}\sin bx$<br>or<br>$e_{\Delta x}(\frac{e^{a\Delta x}\cos b\Delta x-1}{\Delta x}, x)\sin_{\Delta x}(\frac{\tan b\Delta x}{\Delta x}, x)$<br>where<br>$\frac{x}{\Delta x}$ = integer<br>$\cos b\Delta x \neq 0$<br>an equivalent function for $e^{ax}\sin bx$ except where $\cos b\Delta x = 0$ | $\frac{\frac{e^{a\Delta x}\sin b\Delta x}{\Delta x}}{(s - \frac{e^{a\Delta x}\cos b\Delta x-1}{\Delta x})^2 + (\frac{e^{a\Delta x}\sin b\Delta x}{\Delta x})^2}$ |

| # | f(x) Calculus Functions $\lim_{\Delta x \to 0} f(x) = \lim_{\Delta x \to 0} f_{\Delta x}(x)$, $\Delta x \neq 0$ $x = x_0 + m\Delta x$, m = integers, $\Delta x \to 0$ | $f_{\Delta x}(x)$ Discrete Calculus (Interval Calculus) Functions | $f_{\Delta x}(x)$ Calculation Equations $x = x_0 + m\Delta x$, m = integers, $\Delta x \neq 0$ | $K_{\Delta x}[f_{\Delta x}(x)]$ L[f(x)] |
|---|---|---|---|---|
| **29** | $e^{\frac{e^{a\Delta x}\cos b\Delta x - 1}{\Delta x}x}\cos(\frac{e^{a\Delta x}\sin b\Delta x}{\Delta x}x)$ | $e^{ax}\cos bx$ | $e^{ax}\cos bx$<br>or<br>$e_{\Delta x}(\frac{e^{a\Delta x}\cos b\Delta x - 1}{\Delta x}, x)\cos_{\Delta x}(\frac{\tan b\Delta x}{\Delta x}, x)$<br>where<br>$\frac{x}{\Delta x}$ = integer<br>$\cos(b\Delta x) \neq 0$<br>an equivalent function for $e^{ax}\cos bx$ except where $\cos b\Delta x = 0$ | $\dfrac{s - \frac{e^{a\Delta x}\cos b\Delta x - 1}{\Delta x}}{(s - \frac{e^{a\Delta x}\cos b\Delta x - 1}{\Delta x})^2 + (\frac{e^{a\Delta x}\sin b\Delta x}{\Delta x})^2}$ |
| **30** | $e^{\frac{\cosh b\Delta x - 1}{\Delta x}x}\sinh(\frac{\sinh b\Delta x}{\Delta x}x)$ | sinhbx | sinhbx<br>or<br>$e_{\Delta x}(\frac{\cosh b\Delta x - 1}{\Delta x}, x)\sinh_{\Delta x}(\frac{\tanh b\Delta x}{\Delta x}, x)$<br>an equivalent function for sinhbx | $\dfrac{\frac{\sinh b\Delta x}{\Delta x}}{(s - \frac{\cosh b\Delta x - 1}{\Delta x})^2 - (\frac{\sinh b\Delta x}{\Delta x})^2}$ |
| **31** | $e^{\frac{\cosh b\Delta x - 1}{\Delta x}x}\cosh(\frac{\sinh b\Delta x}{\Delta x}, x)$ | coshbx | coshbx<br>or<br>$e_{\Delta x}(\frac{\cosh b\Delta x - 1}{\Delta x}, x)\cosh_{\Delta x}(\frac{\tanh b\Delta x}{\Delta x}, x)$<br>an equivalent function for coshbx | $\dfrac{s - \frac{\cosh b\Delta x - 1}{\Delta x}}{(s - \frac{\cosh b\Delta x - 1}{\Delta x})^2 - (\frac{\sinh b\Delta x}{\Delta x})^2}$ |

Comment - Where $u = \frac{\cos b\Delta x - 1}{\Delta x}$ , $v = \frac{\sin b\Delta x}{\Delta x}$ and $v\Delta x = \sqrt{1 - (1 + u\Delta x)^2}$ , $b = \frac{1}{\Delta x} \tan^{-1} \left(\frac{v\Delta x}{1+u\Delta x}\right)$ for $1+u\Delta x > 0$

or $b = \frac{1}{\Delta x} [\pi + \tan^{-1} (\frac{v\Delta x}{1+u\Delta x})]$ for $1+u\Delta x < 0$

**Table 4.5-2 Related Functions, functions with the same Laplace and $K_{\Delta x}$ Transform**

$A_m, B_m, a, b, c$ are constants

| | Discrete Functions | Laplace/$K_{\Delta x}$ Transform |
|---|---|---|
| **1** | Discrete calculus<br>$y(x)$<br>Calculus<br>$y(x)$<br>See Note 9 | $Y(s)$ |
| **2** | Discrete calculus<br>$D_{\Delta x}{}^{n}y(x) + A_{n-1}D_{\Delta x}{}^{n-1}y(x) + A_{n-2}D_{\Delta x}{}^{n-2}y(x) + \ldots + A_1D_{\Delta x}y(x) + A_0y(x)$<br>Calculus<br>$\frac{d^n}{dx^n}y(x) + A_{n-1}\frac{d^{n-1}}{dx^{n-1}}y(x) + A_{n-2}\frac{d^{n-2}}{dx^{n-2}}y(x) + \ldots + A_1\frac{d}{dx}y(x) + A_0 y(x)$<br>See Note 10 | $s^n y(s) + A_{n-1}s^{n-1}y(s) + A_{n-2}s^{n-2}y(s) + \ldots + A_1sy(s) + A_0y(s) - B_{n-1}s^{n-1} - B_{n-2}s^{n-2} - B_{n-3}s^{n-3} - \ldots - B_1s - B_0$<br>$B_m$, m=0,1,2,…,n-1 are intial condition constants |
| **3** | Discrete calculus<br>$c$<br>Calculus<br>$c$ | $\frac{c}{s}$ |
| **4** | Discrete calculus<br>$x = [x]_{\Delta x}^{1}$<br>Calculus<br>$x$ | $\frac{1}{s^2}$ |
| **5** | Discrete calculus<br>$[x]_{\Delta x}^{n}$<br>Calculus<br>$x^n$<br>n = 0,1,2,3, … | $\frac{n!}{s^{n+1}}$ |
| **6** | Discrete calculus<br>$e_{\Delta x}(a,x)$<br>Calculus<br>$e^{ax}$ | $\frac{1}{s - a}$<br>root s = a |
| **7** | Discrete calculus<br>$\sin_{\Delta x}(b,x)$<br>Calculus<br>$\sin bx$ | $\frac{b}{s^2+b^2}$<br>roots s = jb, -jb |

| | Discrete Functions | Laplace/$K_{\Delta x}$ Transform |
|---|---|---|
| **8** | Discrete calculus<br>$\cos_{\Delta x}(b,x)$<br>Calculus<br>$\cos bx$ | $\frac{s}{s^2+b^2}$<br>roots $s = jb, -jb$ |
| **9** | Discrete calculus<br>$a \neq -\frac{1}{\Delta x}$<br>$e_{\Delta x}(a,x)\ \sin_{\Delta x}(\frac{b}{1+a\Delta x}, x)$<br>$a = -\frac{1}{\Delta x}$<br>$[b\Delta x]^{\frac{x}{\Delta x}} \sin\frac{\pi x}{2\Delta x}$<br>Calculus<br>$e^{ax}\sin bx$ | $\frac{b}{(s-a)^2+b^2}$<br>roots $s = a+jb, a-jb$ |
| **10** | Discrete calculus<br>$a \neq -\frac{1}{\Delta x}$<br>$e_{\Delta x}(a,x)\ \cos_{\Delta x}(\frac{b}{1+a\Delta x}, x)$<br>$a = -\frac{1}{\Delta x}$<br>$[b\Delta x]^{\frac{x}{\Delta x}} \cos\frac{\pi x}{2\Delta x}$<br>Calculus<br>$e^{ax}\cos bx$ | $\frac{s-a}{(s-a)^2+b^2}$<br>roots $s = a+jb, a-jb$ |
| **11** | Discrete calculus<br>$\sinh_{\Delta x}(b,x)$<br>Calculus<br>$\sinh bx$ | $\frac{b}{s^2-b^2}$<br>roots $s = b, -b$ |
| **12** | Discrete calculus<br>$\cosh_{\Delta x}(b,x)$<br>Calculus<br>$\cosh bx$ | $\frac{s}{s^2-b^2}$<br>roots $s = b, -b$ |
| **13** | Discrete calculus<br>$x\sin_{\Delta x}(b,x)$<br>Calculus<br>$x\sin bx+[b\Delta x]x\cos bx$ | $\frac{2bs+[b\Delta x](s^2-b^2)}{(s^2+b^2)^2}$ |

| | Discrete Functions | Laplace/$K_{\Delta x}$ Transform |
|---|---|---|
| **14** | Discrete calculus<br>$x\cos_{\Delta x}(b,x)$<br>Calculus<br>$x\cos bx - [b\Delta x]x\sin bx$ | $\frac{(s^2-b^2)-[b\Delta x]2bs}{(s^2+b^2)^2}$ |
| **15** | Discrete calculus<br>$x\sin_{\Delta x}(b,x-\Delta x)$<br>or<br>$\frac{-[b\Delta x]x\cos_{\Delta x}(b,x)+x\sin_{\Delta x}(b,x)}{1+[b\Delta x]^2}$<br>Calculus<br>$x\sin bx$ | $\frac{2bs}{(s^2+b^2)^2}$ |
| **16** | Discrete calculus<br>$x\cos_{\Delta x}(b,x-\Delta x)$<br>or<br>$\frac{x\cos_{\Delta x}(b,x)+[b\Delta x]x\sin_{\Delta x}(b,x)}{1+[b\Delta x]^2}$<br>Calculus<br>$x\cos bx$ | $\frac{s^2-b^2}{(s^2+b^2)^2}$ |
| **17** | Discrete calculus<br>$[x]^n_{\Delta x}\, e_{\Delta x}(a,x-n\Delta x)$<br>or<br>$\frac{1}{(1+a\Delta x)^n}[x]^n_{\Delta x}\, e_{\Delta x}(a,x)$<br>Calculus<br>$x^n e^{ax}$<br>$n = 0,1,2,3, \ldots$ | $\frac{n!}{(s-a)^{n+1}}$ |
| **18** | Discrete calculus<br>$[x]^n_{\Delta x}\ \sin_{\Delta x}(b,x-n\Delta x)$<br>or | $\frac{n!}{2j}\left[\frac{(s+jb)^{n+1}-(s-jb)^{n+1}}{(s^2+b^2)^{n+1}}\right]$ |

| | Discrete Functions | Laplace/$K_{\Delta x}$ Transform |
|---|---|---|
| | $[x]^{n}_{\Delta x}\cos_{\Delta x}(b,-n\Delta x)\sin_{\Delta x}(b,x) +$<br>$[x]^{n}_{\Delta x}\sin_{\Delta x}(b,-n\Delta x)\cos_{\Delta x}(b,x)$<br><u>Calculus</u><br>$x^{n}\sin bx$<br>$n = 0,1,2,3, \ldots$ | |
| **19** | <u>Discrete calculus</u><br>$[x]^{n}_{\Delta x}\cos_{\Delta x}(b, x-n\Delta x)$<br>or<br>$[x]^{n}_{\Delta x}\cos_{\Delta x}(b,-n\Delta x)\cos_{\Delta x}(b,x) -$<br>$[x]^{n}_{\Delta x}\sin_{\Delta x}(b,-n\Delta x)\sin_{\Delta x}(b,x)$<br><u>Calculus</u><br>$x^{n}\cos bx$<br>$n = 0,1,2,3, \ldots$ | $\frac{n!}{2}\left[\frac{(s+jb)^{n+1}+(s-jb)^{n+1}}{(s^2+b^2)^{n+1}}\right]$ |
| **20** | <u>Discrete calculus</u><br>$x\sinh_{\Delta x}(b,x-\Delta x)$<br>or<br>$\frac{-[b\Delta x]x\cosh_{\Delta x}(b,x)+x\sinh_{\Delta x}(b,x)}{1-[b\Delta x]^2}$<br><u>Calculus</u><br>$x\sinh bx$ | $\frac{2bs}{(s^2-b^2)^2}$ |
| **21** | <u>Discrete calculus</u><br>$x\cosh_{\Delta x}(b,x-\Delta x)$<br>or<br>$\frac{x\cosh_{\Delta x}(b,x)-[b\Delta x]x\sinh_{\Delta x}(b,x)}{1-[b\Delta x]^2}$<br><u>Calculus</u><br>$x\cosh bx$ | $\frac{s^2+b^2}{(s^2-b^2)^2}$ |

| | Discrete Functions | Laplace/$K_{\Delta x}$ Transform |
|---|---|---|
| **22** | Discrete calculus<br>$x\sinh_{\Delta x}(b,x)$<br>Calculus<br>$x\sinh bx+[b\Delta x]x\cosh bx$ | $\frac{2bs+[b\Delta x](s^2+b^2)}{(s^2-b^2)^2}$ |
| **23** | Discrete calculus<br>$x\cos_{\Delta x}(b,x)$<br>Calculus<br>$x\cosh bx+[b\Delta x]x\sinh bx$ | $\frac{(s^2+b^2)+[b\Delta x]2bs}{(s^2-b^2)^2}$ |
| **24** | Discrete calculus<br>$e^{ax}$<br>Calculus<br>$e^{\frac{e^{a\Delta x}-1}{\Delta x}x}$<br>See Note 1 at the end of this table | $\frac{1}{s-\frac{e^{a\Delta x}-1}{\Delta x}}$ |
| **25** | Discrete calculus<br>$A^x$<br>Calculus<br>$e^{\frac{A^{\Delta x}-1}{\Delta x}x}$<br>See Note 2 at the end of this table | $\frac{1}{s-\frac{A^{\Delta x}-1}{\Delta x}}$ |
| **26** | Discrete calculus<br>$[x]_{\Delta x}^{n}\, e^{a(x-n\Delta x)}$<br>or<br>$e^{-an\Delta x}[x]_{\Delta x}^{n}\, e^{ax}$<br>Calculus<br>$x^n e^{\frac{e^{a\Delta x}-1}{\Delta x}x}$<br>$n = 0,1,2,3, \ldots$ | $\frac{n!}{(s-\frac{e^{a\Delta x}-1}{\Delta x})^{n+1}}$ |
| **27** | Discrete calculus<br>$\sin bx$<br>Calculus<br>$e^{\frac{\cos b\Delta x-1}{\Delta x}x}\sin(\frac{\sin b\Delta x}{\Delta x}x)$<br>See Note 3 at the end of this table | $\frac{\frac{\sin b\Delta x}{\Delta x}}{(s-\frac{\cos b\Delta x-1}{\Delta x})^2+(\frac{\sin b\Delta x}{\Delta x})^2}$ |

| | Discrete Functions | Laplace/$K_{\Delta x}$ Transform |
|---|---|---|
| **28** | Discrete calculus<br>$\cos bx$<br>Calculus<br>$e^{\frac{\cos b\Delta x-1}{\Delta x}x}\cos(\frac{\sin b\Delta x}{\Delta x}x)$<br>See Note 4 at the end of this table | $\dfrac{s-\frac{\cos b\Delta x-1}{\Delta x}}{(s-\frac{\cos b\Delta x-1}{\Delta x})^2+(\frac{\sin b\Delta x}{\Delta x})^2}$ |
| **29** | Discrete calculus<br>$e^{ax}\sin bx$<br>Calculus<br>$e^{\frac{e^{a\Delta x}\cos b\Delta x-1}{\Delta x}x}\sin(\frac{e^{a\Delta x}\sin b\Delta x}{\Delta x}x)$<br>See Note 5 at the end of this table | $\dfrac{\frac{e^{a\Delta x}\sin b\Delta x}{\Delta x}}{(s-\frac{e^{a\Delta x}\cos b\Delta x-1}{\Delta x})^2+(\frac{e^{a\Delta x}\sin b\Delta x}{\Delta x})^2}$ |
| **30** | Discrete calculus<br>$e^{ax}\cos bx$<br>Calculus<br>$e^{\frac{e^{a\Delta x}\cos(b\Delta x)-1}{\Delta x}x}\cos(\frac{e^{a\Delta x}\sin(b\Delta x)}{\Delta x}x)$<br>See Note 6 at the end of this table | $\dfrac{s-\frac{e^{a\Delta x}\cos b\Delta x-1}{\Delta x}}{(s-\frac{e^{a\Delta x}\cos b\Delta x-1}{\Delta x})^2+(\frac{e^{a\Delta x}\sin b\Delta x}{\Delta x})^2}$ |
| **31** | Discrete calculus<br>$\sinh bx$<br>Calculus<br>$e^{\frac{\cosh b\Delta x-1}{\Delta x}x}\sinh(\frac{\sinh b\Delta x}{\Delta x}x)$<br>See Note 7 at the end of this table | $\dfrac{\frac{\sinh b\Delta x}{\Delta x}}{(s-\frac{\cosh b\Delta x-1}{\Delta x})^2-(\frac{\sinh b\Delta x}{\Delta x})^2}$ |
| **32** | Discrete calculus<br>$\cosh bx$<br>Calculus<br>$e^{\frac{\cosh b\Delta x-1}{\Delta x}x}\cosh(\frac{\sinh b\Delta x}{\Delta x}x)$<br>See Note 8 at the end of this table | $\dfrac{s-\frac{\cosh b\Delta x-1}{\Delta x}}{(s-\frac{\cosh b\Delta x-1}{\Delta x})^2-(\frac{\sinh b\Delta x}{\Delta x})^2}$ |

Notes

Notes 1 thru 8 are identities

1. $e^{ax} = e_{\Delta x}(\frac{e^{a\Delta x}-1}{\Delta x}, x)$

2. $A^x = e_{\Delta x}(\frac{A^{\Delta x}-1}{\Delta x}, x)$

3. $\sin bx = e_{\Delta x}(\frac{\cos b\Delta x - 1}{\Delta x}, x) \sin_{\Delta x}(\frac{\tan b\Delta x}{\Delta x}, x)$ , $\frac{x}{\Delta x}$ = integer , $\cos(b\Delta x) \neq 0$

4. $\cos bx = e_{\Delta x}(\frac{\cos b\Delta x - 1}{\Delta x}, x) \cos_{\Delta x}(\frac{\tan b\Delta x}{\Delta x}, x)$ , $\frac{x}{\Delta x}$ = integer , $\cos(b\Delta x) \neq 0$

5. $e^{ax} \sin bx = e_{\Delta x}(\frac{e^{a\Delta x}\cos b\Delta x - 1}{\Delta x}, x)\sin_{\Delta x}(\frac{\tan b\Delta x}{\Delta x}, x)$ , $\frac{x}{\Delta x}$ = integer , $\cos(b\Delta x) \neq 0$

6. $e^{ax} \cos bx = e_{\Delta x}(\frac{e^{a\Delta x}\cos b\Delta x - 1}{\Delta x}, x)\sin_{\Delta x}(\frac{\tan b\Delta x}{\Delta x}, x)$ , $\frac{x}{\Delta x}$ = integer , $\cos(b\Delta x) \neq 0$

7. $\sinh bx = e_{\Delta x}(\frac{\cosh b\Delta x - 1}{\Delta x}, x) \sinh_{\Delta x}(\frac{\tanh b\Delta x}{\Delta x}, x)$

8. $\cosh bx = e_{\Delta x}(\frac{\cosh b\Delta x - 1}{\Delta x}, x) \cosh_{\Delta x}(\frac{\tanh b\Delta x}{\Delta x}, x)$

9. The table upper discrete calculus function, y(x), is associated with discrete differential equations and the lower Calculus function, *y*(x), is associated with Calculus differential equations. The relationship between the two functions is $\lim_{\Delta x \to 0} f_1(x) = \lim_{\Delta x \to 0} f_2(x)$.
10. For the Related Functions discrete calculus and Calculus differential equations, the value of the respective initial conditions is the same (i.e. y(0) = *y*(0), $D_{\Delta x}y(0) = \frac{d}{dx}y(0)$, $D_{\Delta x}^2 y(0) = \frac{d^2}{dx^2}y(0)$, etc.).

Observing Table 4.5-2, there are two important things which should be noted. The first is that the discrete calculus function and Calculus function which are paired share a common Laplace/$K_{\Delta x}$ Transform. The second thing is that both functions become identical if Δx goes to zero. As a result of these relationships, these paired functions are designated as Related Functions. It is these related functions which will be used in the Method of Related Functions. This method will be used in the solution of discrete differential difference equations.

Description of the Method of Related Functions

Using the Method of Related Functions, a discrete calculus differential equation can be solved by first converting it into its related Calculus differential equation using the functions in Table 4.5-2. This related Calculus differential equation is then solved by means of any of the usual Calculus methods. When the solution function to the Calculus differential equation is obtained, it is converted into its related discrete calculus function again using the related functions of Table 4.5-2. The resulting function is the solution to the original discrete calculus differential equation to be solved.

The Method of Related Functions is shown in detail in the block diagram which follows. In this block diagram the function relationships to the Laplace/$K_{\Delta x}$ Transform are shown in dotted lines. However, in the application of this method, the transform relationships are not directly involved.

<u>Block Diagram 4.5-1</u> The Method of Related Functions used to solve differential difference equations

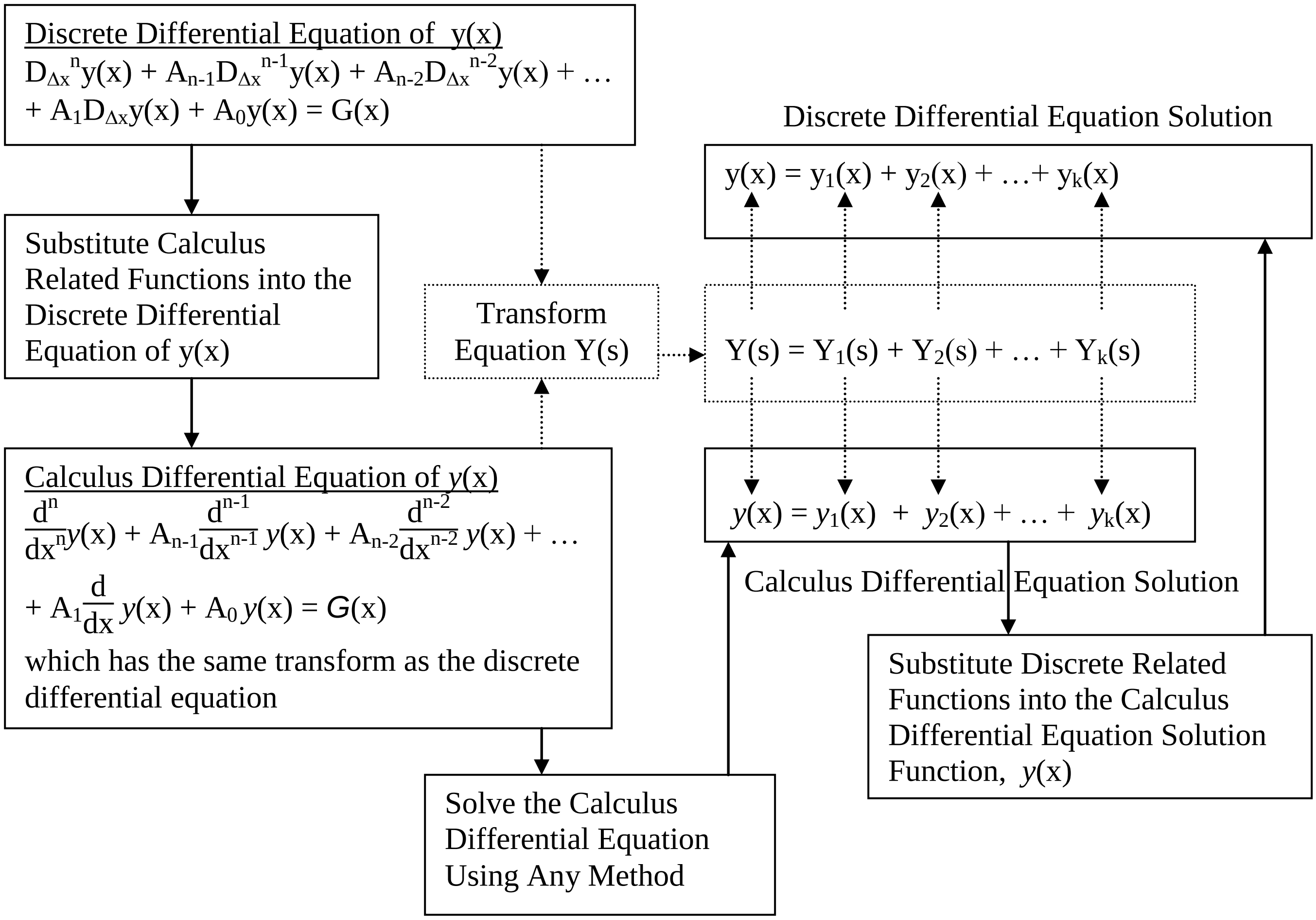

<u>Notes</u> - 1) $[\mathrm{D_{\Delta x}}^{n}\mathrm{y(x)} + \mathrm{A_{n-1}D_{\Delta x}}^{n-1}\mathrm{y(x)} + \mathrm{A_{n-2}D_{\Delta x}}^{n-2}\mathrm{y(x)} + \ldots + \mathrm{A_1D_{\Delta x}y(x)} + \mathrm{A_0y(x)},$
$\frac{d^n}{dx^n}y(x) + A_{n-1}\frac{d^{n-1}}{dx^{n-1}}y(x) + A_{n-2}\frac{d^{n-2}}{dx^{n-2}}y(x) + \ldots + A_1\frac{d}{dx}y(x) + A_0\,y(x)]$ are Related Functions.

2) $[\mathrm{G(x)}, \mathit{G}(\mathrm{x})]$, $[\mathrm{y(x)}, y(\mathrm{x})]$, $[\mathrm{y_m(x)}, y_m(\mathrm{x})\ \ m = 1,2,3,\ldots,k]$ are Related Functions. They have the same Laplace/$\mathrm{K_{\Delta x}}$ Transform but are not equal functions.

3) Related Function pairs are listed in Table 4.5-2.

4) The Related Functions Method can only solve discrete differential equations (also called differential difference equations) which contain functions that can be found in the Related Functions listing, Table 4.5-2.

5) The block diagram solid lines represent the Method of Related Functions process. The block diagram dotted lines show the relationship of the Related Functions method to the transform solution method. The Related Functions solution process itself does not involve the use of transforms unless the Laplace Transform, a special case of the $\mathrm{K_{\Delta x}}$ Transform, is used to find the solution to the related Calculus differential equation.

6) The Related Function substitution specified in the method process is simply the replacement of one function of a related function pair with the other.

8) The symbol, $\overrightarrow{\leftarrow}$, has been selected to indicate a Related Functions relationship. For example, $e^{ax} \overrightarrow{\leftarrow} e^{\frac{e^{a\Delta x}-1}{\Delta x}x}$ is a Related Functions relationship. The order in which the functions are written is discrete calculus function on the left, Calculus function on the right.

On the following page there are two examples, Example 4.5-1 and Example 4.5-2, that show the application of the Method of Related Functions to the solution of a discrete calculus differential equation. The solution process is in accordance with Block Diagram 4.5-1 presented on the previous page.

Example 4.5-1 Use of the Method of Related Functions to solve the differential difference equation

$$D_{\Delta x}y(x) + y(x) = xe^{-x}, \qquad \Delta x = 2,\ y(0) = 3,\ x = 0,2,4,6,\ldots \tag{1}$$

Using the Method of Related Functions convert Eq 1 into a related Calculus equation.

From the table of related functions (Discrete Function $\rightleftarrows$ Calculus Function)

The following related functions are reversible.

$$y(x) \rightleftarrows Y(x) \tag{2}$$

$$D_{\Delta x} \rightleftarrows \frac{d}{dx} \tag{3}$$

$$y(0) = Y(0) \tag{4}$$

$$e_{\Delta x}(a,x) \rightleftarrows e^{ax} \tag{5}$$

$$xe^{a(x-\Delta x)} \rightleftarrows xe^{\frac{e^{a\Delta x}-1}{\Delta x}x} \tag{6}$$

Evaluating Eq 6

For $a = -1$
$\Delta x = 2$

$$e^{a\Delta x} = e^{-2} = .135335283$$

$$\frac{e^{a\Delta x}-1}{\Delta x} = \frac{e^{-2}-1}{2} = -.432332358$$

$$xe^{ax} \rightleftarrows e^{a\Delta x}\, xe^{\frac{e^{a\Delta x}-1}{\Delta x}x}$$

$$xe^{-x} \rightleftarrows .135335283\, xe^{-.432332358x} \tag{7}$$

From the related functions of Eq 2 and Eq 7, convert the discrete differential difference equation, Eq 1, into its related Calculus differential equation.

$$\frac{d}{dx}Y(x) + Y(x) = .135335283\, xe^{-.432332358x}, \quad Y(0) = y(0) = 3 \tag{8}$$

Using Calculus solve Eq 8

The homogeneous differential equation for Eq 8 is:

$$\frac{d}{dx} Y_c(x) + Y_c(x) = 0 \qquad 9)$$

Solve Eq 9 for the related Calculus differential equation complementary solution, $Y_c(x)$.

Let

$$Y_c(x) = ce^{bx} \qquad 10)$$

Substituting this $Y_c(x)$ function into the homogenous differential equation, Eq 9

$$cbe^{bx} + ce^{bx} = (b + 1)ce^{bx} = 0$$

$$b = -1$$

Substituting b = -1 into Eq 10
The complementary solution to the differential equation, Eq 8, is:

$$Y_c(x) = ce^{-x} \qquad 11)$$

Find the particular solution to the differential equation, Eq 8, using the Method of Undetermined Coefficients

This method specifies a particular solution:

$$Y_p(x) = (K_1x + K_2)e^{-.432332358x} \qquad 12)$$

Substitute Eq 12 into Eq 8

$$K_1e^{-.432332358x} - .432332358(K_1x + K_2)e^{-.432332358x} + (K_1x + K_2)e^{-.432332358x} = .135335283xe^{-.432332358x}$$

Combining terms

$$e^{-.432332358x}[K_1 - .432332358K_2 + K_2] + xe^{-.432332358x}[-.432332358K_1 + K_1] = .135335283xe^{-.432332358x}$$

$$(1-.432332358)K_1 = .135335283$$

$$K_1 = .238405843$$

$$K_1 + (1-.432332358)K_2 = 0$$

$$K_2 = \frac{-.238405843}{1-.432332358}$$

$$K_2 = -.419974339$$

Substituting the constants $K_1$ and $K_2$ into Eq 12

$Y_p(x) = (.238405843x - .419974339)e^{-.432332358x}$ 13)

$Y(x) = Y_c(x) + Y_p(x)$ 14)

Substituting Eq 11 and Eq 13 into Eq 14

$Y(x) = ce^{-x} + (-\ .419974339 + .238405843x)e^{-.432332358x}$ 15)

Find the constant, c, using the initial condition, Y(0) = 3

$3 = c - .419974339$

$c = 3.419974339$

Substituting the value of c into Eq 15

$Y(x) = 3.419974339e^{-x} + (-\ .419974339 + .238405843x)e^{-.432332358x}$ 16)

This is the related Calculus function to y(x).

Using the related functions of Eq 2 thru Eq 7 and the following related functions and equations, the Calculus function, Y(x), can be converted into its related discrete Interval Calculus function, y(x).

From the related functions table

$e_{\Delta x}(a,x) = (1+a\Delta x)^{\frac{x}{\Delta x}} = e^{\frac{x}{\Delta x}\ln(1+a\Delta x)} \rightleftarrows e^{ax}$ 17)

For a = -.432332358
$\Delta x = 2$

Substiituting into Eq 17

$e_2(-.432332358,x) = (1+[-.432332358][2])^{\frac{x}{2}} = e^{\frac{\ln[1+(-.432332358)(2)]}{2}x} = e^{-x} \rightleftarrows e^{-.432332358x}$

$e^{-x} \rightleftarrows e^{-.432332358x}$ 18)

From Eq 7

$\frac{1}{.135335283} xe^{-x} \rightleftarrows xe^{-.432332358x}$ 19)

From Eq 5

$$e_{\Delta x}(a,x) = (1+a\Delta x)^{\frac{x}{\Delta x}} \rightleftarrows e^{ax}$$

a = -1
Δx = 2

$$e_2(-1,x) = [1+(-1)(2)]^{\frac{x}{2}} = [-1]^{\frac{x}{2}} \rightleftarrows e^{-x} \qquad 20)$$

$$[-1]^{\frac{x}{2}} \rightleftarrows e^{-x} \qquad 21)$$

From Eq 2

$$y(x) \rightleftarrows Y(x) \qquad 22)$$

Substituting the related functions of Eq 18. Eq 19, Eq 21 and Eq 22 into Eq 16

$$y(x) = 3.419974339[-1]^{\frac{x}{2}} + (-\ .419974339 + \frac{.238405843}{.135335283}x)e^{-x}$$

$$\mathbf{y(x) = 3.419974339[-1]^{\frac{x}{2}} + (-\ .419974339 + 1.761594151x)e^{-x}} \qquad 23)$$

Eq 23 can be put into another form

From Eq 20

$$[-1]^{\frac{x}{2}} = e_2(-1,x) \qquad 24)$$

$$e^{ax} = (1 + \frac{e^{a\Delta x}-1}{\Delta x}\Delta x)^{\frac{x}{\Delta x}} = e_{\Delta x}(\frac{e^{a\Delta x}-1}{\Delta x},x) \qquad 25)$$

For a = -1
Δx = 2

$$e^{-x} = e_2(\frac{e^{-2}-1}{2},x)$$

$$e^{-x} = e_2(-.432332358,x) \qquad 26)$$

Substituting Eq 24 and Eq 26 into Eq 23

$$\mathbf{y(x) = 3.419974339e_2(-1,x) + (-\ .419974339 + 1.761594151x)e_2(-.432332358,x)} \qquad 27)$$

Checking Eq 23 and Eq 27

y(0) = 3.4199743393-.419974339 = 3     Good check

Substitute Eq 23 into Eq 1

$$\frac{3.419974339[-1]^{\frac{x+2}{2}} + (-.419974339 + 1.761594151\{x+2\})e^{-(x+2)}}{2}$$

$$-\frac{3.419974339[-1]^{\frac{x}{2}} + (-.419974339 + 1.761594151x)e^{-x}}{2}$$

$$+ 3.419974339[-1]^{\frac{x}{2}} + (-.419974339 + 1.761594151x)e^{-x} = xe^{-x}$$

$$+\frac{-3.419974339[-1]^{\frac{x}{2}}}{2} + \frac{-.419974339e^{-2}e^{-x}}{2} + \frac{1.761594151(2)e^{-2}e^{-x}}{2} + \frac{1.761594151e^{-2}xe^{-x}}{2}$$

$$-\frac{3.419974339[-1]^{\frac{x}{2}}}{2} + \frac{.419974339e^{-x}}{2} - \frac{1.761594151xe^{-x}}{2}$$

$$+ 3.419974339[-1]^{\frac{x}{2}} - .419974339e^{-x} + 1.761594151xe^{-x} = xe^{-x}$$

$$\left[\frac{-.419974339e^{-2}}{2} + \frac{1.761594151(2)e^{-2}}{2} + \frac{.419974339}{2} - .419974339\right]e^{-x} +$$

$$\left[\frac{1.761594151e^{-2}}{2} - \frac{1.761594151}{2} + 1.761594151\right]xe^{-x} = xe^{-x}$$

$0e^{-x} + 1xe^{-x} = xe^{-x}$

$xe^{-x} = xe^{-x}$     Good check

Example 4.5-2 Use of the Method of Related Functions to solve the differential difference equation:

$$D_{\Delta x}^{2}y(x) + D_{\Delta x}y(x) = x(x-\Delta x) = [x]_{\Delta x}^{2}\ , \quad \Delta x = 2,\ y(0) = 1,\ \ y(2) = 5,\ x = 0,2,4,6,\ldots \tag{1}$$

Using the Method of Related Functions convert Eq 1 into a related Calculus equation

From the table of related functions (Discrete Function $\rightleftarrows$ Calculus Function)

The following function relationships are reversible.

$$y(x) \rightleftarrows Y(x) \tag{2}$$

$$D_{\Delta x} \rightleftarrows \frac{d}{dx} \tag{3}$$

$$y(0) = Y(0) \tag{4}$$

$$y^{'}(0) = Y^{'}(0) \quad \text{i.e.} \quad D_{\Delta x}y(x)|_{x=0} = \frac{d}{dx}Y(x)|_{x=0} \tag{5}$$

$$e_{\Delta x}(a,x) \rightleftarrows e^{ax} \tag{6}$$

$$[x]_{\Delta x}^{n} \rightleftarrows x^{n} \quad n = 0,1,2,3,\ldots \tag{7}$$

Converting Eq 1 to its related Calculus equation using the relationships specified above

$$\frac{d^2}{dx^2}Y(x) + \frac{d}{dx}Y(x) = x^2\ , \qquad Y(0) = 1,\ \ Y^{'}(0) = y^{'}(0) = \frac{5-1}{2} = 2 \tag{8}$$

The related homogeneous equation for Eq 8 is:

$$\frac{d^2}{dx^2}Y_c(x) + \frac{d}{dx}Y_c(x) = 0 \tag{9}$$

Find the complementary solution, $Y_c(x)$, for Eq 8

Let

$$Y_c(x) = ce^{bx} \tag{10}$$

Substituting into Eq 9

$$b^2ce^{bx} + bce^{bx} = 0$$

$b(b+1)ce^{bx} = 0$

$b = 0,-1$

Substituting into Eq 10

$Y_c(x) = c_1e^{-1x} + c_2e^{0x} = c_1e^{-x} + c_2$

$Y_c(x) = c_1e^{-x} + c_2$ 11)

Find the particular solution for Eq 8 using the Method of Undetermined Coefficients

Note that a root of the homogeneous equation is 0 and so is that of the function, $x^2$ (i.e. $x^2e^{0x}$). The Method of Undetermined Coeffients would normally specify a particular solution of $Y_p(x) = K_1x^2+K_2x+K_3$. However, because of the common root, 0, $K_1x^2+K_2x+K_3$ must be multiplied by x.

Then

$Y_p(x) = K_1x^3+K_2x^2+K_3x$ 12)

Substituting Eq 12 into Eq 8

$6K_1x +2K_2 + 3K_1x^2 + 2K_2x + K_3 = x^2$

Combining terms

$(3K_1 – 1)x^2 +(6K_1 + 2K_2)x + (2K_2 + K_3) = 0$

$K_1 = \frac{1}{3}$

$6(\frac{1}{3}) + 2K_2 = 0$

$K_2 = -1$

$2(-1) + K_3 = 0$

$K_3 = 2$

Substituting the values $K_1$,$K_2$,$K_3$ into Eq 12

$Y_p(x) = \frac{1}{3}x^3 – x^2 + 2x$ 13)

$Y(x) = Y_c(x) + Y_p(x)$ 14)

Substituting Eq 11 and Eq 13 into Eq 14

$$Y(x) = c_1e^{-x} + c_2 + \frac{1}{3}x^3 - x^2 + 2x \qquad 15)$$

Find the constants $c_1$ and $c_2$ from the initial conditions $Y(0) = 1$ and $Y^{'}(0) = 2$

$Y(0) = 1 = c_1 + c_2$

$c_1 + c_2 = 1$

$Y^{'}(0) = 2 = - c_1 + 2$

$c_1 = 0$
$c_2 = 1$

Substituting the values of $c_1$ and $c_2$ into Eq 15

$$Y(x) = \frac{1}{3}x^3 - x^2 + 2x + 1 \qquad 15)$$

This function is the Calculus related function of the discrete Interval Calculus function Eq 1

Checking Eq 15

$Y(0) = 1$ Good check

$Y^{'}(0) = x^2 - 2x + 2 \,|_{x=0} = 2$ Good check

Substituting Eq 15 into Eq 8

$2x - 2 + x^2 - 2x + 2 = x^2$

$x^2 = x^2$ Good check

Convert Y(x) to y(x)

From Eq 7

$[x]^{n}_{\Delta x} \rightleftarrows x^n \quad n = 0,1,2,3,\ldots$

$$[x]^{n}_{\Delta x} = \prod_{m=1}^{n}(x-[m-1]\Delta x)$$

$[x]^{0}_{\Delta x} = 1$

Substituting these previous relationships into Eq 15

$$y(x) = \frac{1}{3}[x]_2^3 - [x]_2^2 + 2[x]_2^1 + 1 \;, \qquad \Delta x = 2 \qquad\qquad 16)$$

or

$$y(x) = \frac{1}{3}\,x(x-2)(x-4) - x(x-2) + 2x + 1 \qquad\qquad 17)$$

Checking Eq 16 and Eq 17

$y(0) = 1$ Good check

$y(2) = 4 + 1 = 5$ Good check

Substitute Eq 16 into Eq 1 and take the necessary derivatives using the derivative formula,
$D_{\Delta x}[x]_{\Delta x}^{n} = n[x]_{\Delta x}^{n-1}$

$$\frac{3(2)}{3}[x]_2^1 - 2(1)[x]_2^0 + 0 + 0 + \frac{3}{3}[x]_2^2 - 2[x]_2^1 + 2 = [x]_2^2$$

$2x - 2 + x(x-2) - 2x + 2 = x(x-2)$

$x(x-2) = x(x-2)$ Good check

<u>Comment</u> – The Method of Related Functions can also be used to perform discrete integration. To obtain a discrete integral, $y(x) = {}_{\Delta x}\!\int f_{\Delta x}(x)\Delta x$, the following procedure is used:

1. If the function to be integrated is not a discrete function, it is replaced by its discrete function identity.
2. Using the Related Functions Table, Table 4.5-2, the Calculus function, $f(x)$, is obtained by replacing the discrete functions comprising $f_{\Delta x}(x)$ with their related Calculus functions. $f_{\Delta x}(x) \rightarrow f(x)$.
3. The function, $F(x)$, is obtained using Calculus to integrate $f(x)$. $F(x) = \int f(x)dx$.
4. Using the Related Functions Table, Table 4.5-2, the discrete calculus function, $y(x)$, is obtained by replacing the Calculus functions comprising $F(x)$ with their related discrete calculus functions. $F(x) \rightarrow y(x)$.

# CHAPTER 5

# Discrete Variable Operational Calculus

## Derivation of the Discrete Calculus $J_{\Delta x}$ and $K_{\Delta x}$ Transforms, the Conversion relationships between the $K_{\Delta x}$ and Z Transforms, and some useful Interval Calculus Functions and Equalities

## Section 5.1: Derivation of the $J_{\Delta x}$ and $K_{\Delta x}$ Transforms

Initial Comment – Because transforms such as the Laplace Transform and the Z Transform are most often applied to time domain problems, the variable selected for the following derivations and explanations will be t instead of x.

Interval Calculus is a discrete calculus developed for the solution of mathematical problems involving discrete variables. It would be very useful to have the ability to conveniently apply the Z Transform as a part of its discrete mathematical calculations. However, a significant difficulty is encountered. Interval Calculus discrete differentiation and integration are most easily performed on Interval Calculus functions such as $e_{\Delta t}(a,t)$, $\sin_{\Delta x}(a,t)$, $\cos_{\Delta x}(a,t)$, etc.but the Z Transform, as derived, uses Calculus functions such as $e^{at}$, sin(at), cos(at), etc. Fortunately, this mismatch can be overcome. The Z Transform can be modified. The resulting two modifications, the so called $J_{\Delta x}$ and $K_{\Delta x}$ Transforms, are compatible with all Interval Calculus functions and consistent with all Interval Calculus notational conventions.

The derivation of the $J_{\Delta x}$ Transform from the Z Transform is performed below.

The Z Transform

$$F(z) = \sum_{n=0}^{\infty} f(nT) z^{-n} \qquad (5.1\text{-}1)$$

The Inverse Z Transform

$$f(nT) = \frac{1}{2\pi j} \oint_{c_1} F(z)\, z^{n-1} dz \qquad (5.1\text{-}2)$$

The $c_1$ contour is:

$$z = e^{(\gamma + jw)T} \text{, the complex plane circular contour, } c_1 \qquad (5.1\text{-}3)$$

where

- $n = 0,1,2,3, \ldots$
- The contour $c_1$ is a circular contour in the Z Plane
- $\gamma$ = positive real constant
- $-\frac{\pi}{T} \le w < \frac{\pi}{T}$
- f(t) = function of t
- F(z) = Z[f(t)] , the Z Transform of f(t)

Figure 5.1-1

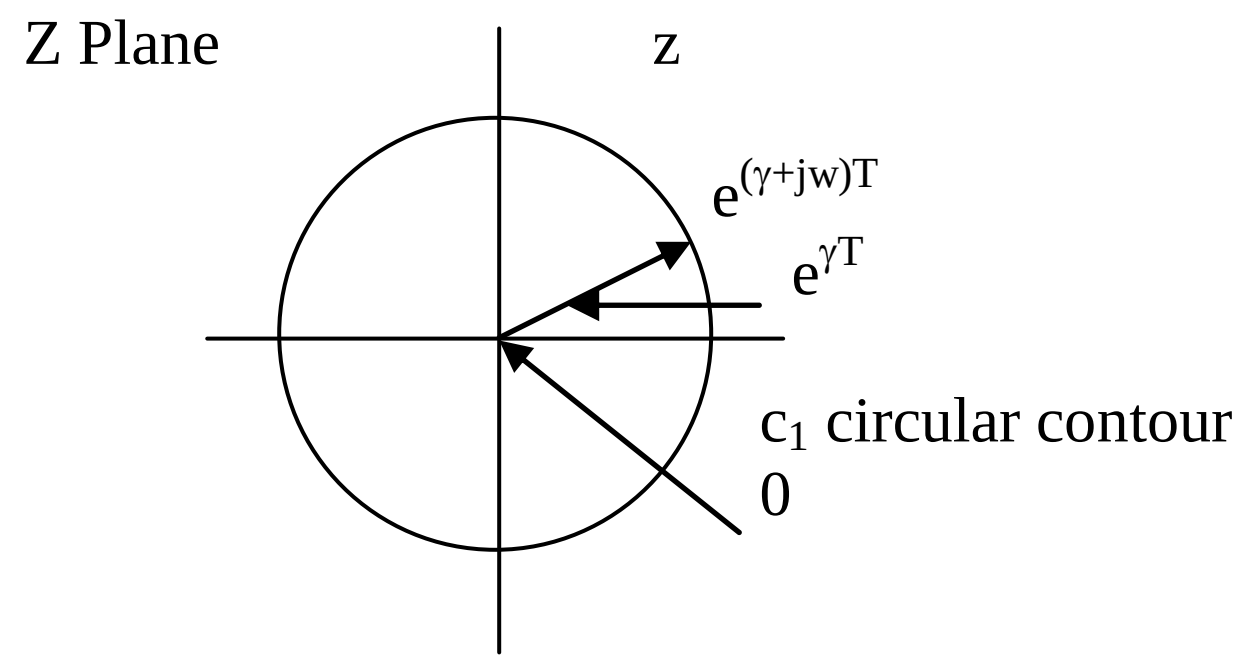

In the following derivation of the $J_{\Delta t}$ and $K_{\Delta t}$ Transforms, Interval Calculus notation will be used.

$\Delta t = T$ (5.1-4)

$t = n\Delta t \;,\quad n = 0,1,2,3,\ldots$ (5.1-5)

or

$t = 0,\ \Delta t,\ 2\Delta t,\ 3\Delta t,\ \ldots$ (5.1-6)

Map the Z Plane of the Z Transform into the S Plane using the linear transformation,

$$s = \frac{1}{\Delta t} z - \frac{1}{\Delta t} = \frac{z-1}{\Delta t}$$

$s = \dfrac{z-1}{\Delta t}$ , The linear conformal mapping transform (5.1-7)

Substituting Eq 5.1-3 into Eq .1-7

$s = \dfrac{e^{(\gamma+jw)\Delta t}-1}{\Delta t}$ , the complex plane circular contour, $c_2$ (5.1-8)

where

The contour $c_2$ is a circular contour in the S Plane

$\gamma$ = positive real constant

$-\dfrac{\pi}{\Delta t} \le w < \dfrac{\pi}{\Delta t}$

The linear conformal mapping transform, Eq 5.1-7, maps the $c_1$ circular Z Plane locus of Figure 1.5-1 into the $c_2$ circular S Plane locus of Figure 1.5-2 below.

Figure 5.1-2

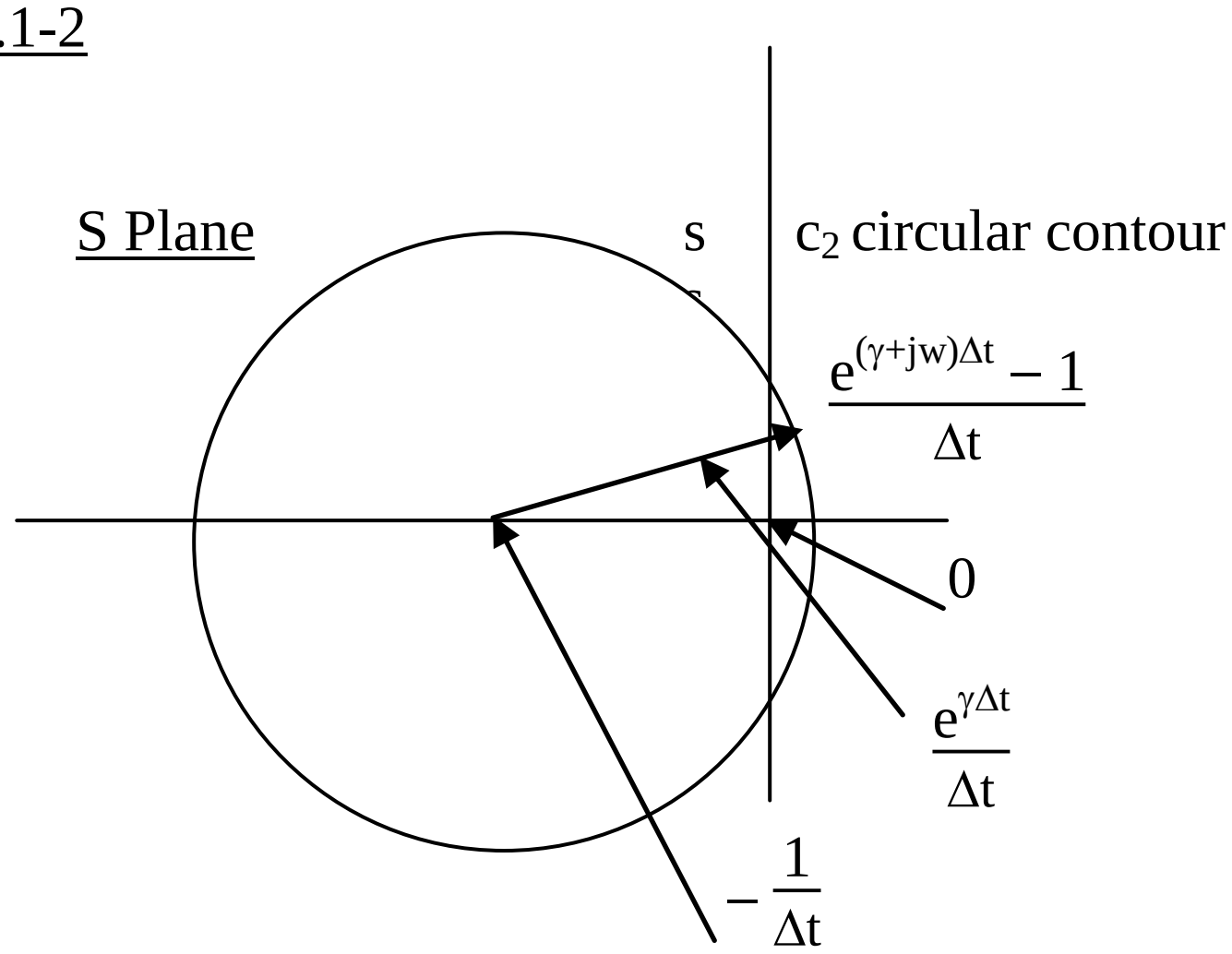

From Eq 5.1-7

$$z = 1+s\Delta t \tag{5.1-9}$$

From Eq 5.1-1, Eq 5.1-5, and 5.1-9

$$F(z)\Big|_{z = 1+s\Delta t} = \sum_{\Delta t \, t=0}^{\infty} f(t)(1+s\Delta t)^{-\frac{t}{\Delta t}} \tag{5.1-10}$$

where
- $t = n\Delta t$ , $n = 0,1,2,3,\ldots$
- $t = 0, \Delta t, 2\Delta t, 3\Delta t, \ldots$
- $f(t)$ = function of t
- $F(z)$ = Z Transform of $f(t)$
- $\Delta t$ = sampling period, the t increment

Let

$$F(s) = \Delta t F(z)\Big|_{z = 1+s\Delta t} = \Delta t \sum_{\Delta t \, t=0}^{\infty} f(t)(1+s\Delta t)^{-\frac{t}{\Delta t}} = \sum_{\Delta t \, t=0}^{\infty} f(t)(1+s\Delta t)^{-\frac{t}{\Delta t}} \Delta t \tag{5.1-11}$$

$$F(s) = \sum_{\Delta t \, t=0}^{\infty} f(t)(1+s\Delta t)^{-\frac{t}{\Delta t}} \Delta t \tag{5.1-12}$$

Changing the Interval Calculus summation notation of Eq 5.1-12 to its equivalent Interval Calculus integration notation.

$$F(s) = \int_{0}^{\infty} {}_{\Delta t}\, f(t)(1+s\Delta t)^{-\frac{t}{\Delta t}}\,\Delta t \qquad (5.1\text{-}13)$$

Note from Eq 5.1-10 thru Eq 5.1-13 that F(s), due to the introduction of Δt, represents an area calculation.

Comment – Eq 5.1-13 is compatible with the Interval Calculus discrete integration and differentiation of the very important Interval Calculus function, $(1+c\Delta t)^{\frac{t}{\Delta t}}$. The Interval Calculus discrete functions, $e_{\Delta t}(a,t)$, $\sin_{\Delta t}(b,t)$, $\cos_{\Delta t}(b,t)$, $\sinh_{\Delta t}(b,t)$, $\cosh_{\Delta t}(b,t)$, and several others are based on this function.

Comment – The function, $(1+s\Delta t)^{-\frac{t}{\Delta t}}$, is equal to $e^{-(\gamma+jw)t}$

$$(1+s\Delta t)^{-\frac{t}{\Delta t}} = \left(1+\left[\frac{e^{(\gamma+jw)\Delta t}-1}{\Delta t}\right]\Delta t\right)^{-\frac{t}{\Delta t}} = e^{-(\gamma+jw)t}$$

Comment – $\lim_{\Delta t\to 0} (1+s\Delta t)^{-\frac{t}{\Delta t}} = e^{-st}$

$$\lim_{\Delta t\to 0} (1+s\Delta t)^{-\frac{t}{\Delta t}} = \lim_{\Delta t\to 0} e^{-t\frac{\ln(1+s\Delta t)}{\Delta t}} = \lim_{\Delta t\to 0} e^{-t\frac{s}{1+s\Delta t}} = e^{-st}$$

Comment – For $\Delta t \to 0$ Eq 5.1-13 is the Laplace Transform

$$\lim_{\Delta t\to 0} \int_{0}^{\infty} {}_{\Delta t}\, f(t)(1+s\Delta t)^{-\frac{t}{\Delta t}}\,\Delta t = \int_{0}^{\infty} f(t)e^{-st}dt$$

Eq 5.1-13 has arbitrarily been named the $J_{\Delta t}$ Transform. The subscript, Δt, represents the value of the sampling (sample and hold) interval. Note this transform's similarity to the Laplace Transform.

Map the Z Plane Inverse Z Transform, Eq 5.1-2, into S Plane

Rewriting the Inverse Z Transform equation, Eq 5.1-2

$$f(nT) = \frac{1}{2\pi j} \oint_{c_1} F(z)\, z^{n-1} dz$$

From Eq 5.1-7

$$\Delta t ds = dz \qquad (5.1\text{-}14)$$

From Eq 5.1-2, Eq 5.1-5, Eq 5.1-8, Eq 5.1-9 and Eq 5.1-14

$$f(t) = \frac{1}{2\pi j} \oint_{c_2} F(z)|_{z\,=\,1+s\Delta t}\, (1+s\Delta t)^{\frac{t}{\Delta t}-1} \Delta t ds \qquad (5.1\text{-}15)$$

From Eq Eq 5.1-15, Eq 5.1-11, and Eq 5.1-13

$$f(t) = \frac{1}{2\pi j} \oint_{c_2} \frac{1}{\Delta t} \left[ {}_{\Delta t}\!\int_0^{\infty} f(t)(1+s\Delta t)^{-\frac{t}{\Delta t}} \Delta t \right] (1+s\Delta t)^{\frac{t}{\Delta t}-1} \Delta t ds \qquad (5.1\text{-}15)$$

$$f(t) = \frac{1}{2\pi j} \oint_{c_2} \left[ {}_{\Delta t}\!\int_0^{\infty} f(t)(1+s\Delta t)^{-\frac{t}{\Delta t}} \Delta t \right] (1+s\Delta t)^{\frac{t}{\Delta t}-1} ds \qquad (5.1\text{-}16)$$

Rewriting Eq 5.1-13

$$\mathsf{F}(s) = {}_{\Delta t}\!\int_0^{\infty} f(t)(1+s\Delta t)^{-\frac{t}{\Delta t}} \Delta t \qquad (5.1\text{-}17)$$

, the $J_{\Delta t}$ Transform

Substituting Eq 5.1-17 into Eq 5.1-16

The Inverse $J_{\Delta t}$ Transform is:

$$f(t) = \frac{1}{2\pi j} \oint_{c_2} \mathsf{F}(s)(1+s\Delta t)^{\frac{t}{\Delta t}-1} ds \qquad (5.1\text{-}18)$$

The $J_{\Delta t}$ Transform and the Inverse $J_{\Delta t}$ Transform have been derived.

From Eq 5.1-17 and Eq 5.1-18

**The $J_{\Delta t}$ Transform is:**

$$\mathbf{F(s) = J_{\Delta t}[f(t)] = \underset{\Delta t}{\int_0^\infty} f(t)(1+s\Delta t)^{-\frac{t}{\Delta t}}\Delta t} \qquad (5.1\text{-}19)$$

**The Inverse $J_{\Delta t}$ Transform is:**

$$\mathbf{f(t) = J_{\Delta t}^{-1}[F(s)] = \frac{1}{2\pi j}\oint_{c_2} F(s)(1+s\Delta t)^{\frac{t}{\Delta t}-1} ds} \qquad (5.1\text{-}20)$$

**where**

**$s = \dfrac{e^{(\gamma+jw)\Delta t} - 1}{\Delta t}$ , the complex plane circular contour, $c_2$**

**The contour $c_2$ is a circular contour in the S Plane**

**$\gamma$ = positive real constant**

**$-\dfrac{\pi}{\Delta t} \le w < \dfrac{\pi}{\Delta t}$**

**f(t) = function of t**

**F(s) = $J_{\Delta t}$[f(t)] , a function of s**

**$\Delta t$ = sampling interval, t increment**

**$t = n\Delta t$ , n = 0,1,2,3,…**

**$t = 0, \Delta t, 2\Delta t, 3\Delta t, \ldots$**

The derivation of the $K_{\Delta x}$ Transform from the Z Transform is performed below.

Rewriting Eq 5.1-16

$$f(t) = \frac{1}{2\pi j}\oint_{c_2} [\underset{\Delta t}{\int_0^\infty} f(t)(1+s\Delta t)^{-\frac{t}{\Delta t}}\Delta t\,](1+s\Delta t)^{\frac{t}{\Delta t}-1} ds$$

The term, $(1+s\Delta t)^{-1}$, can be moved into the discrete integral since it is not a function of t.

$$f(t) = \frac{1}{2\pi j}\oint_{c_2} [\underset{\Delta t}{\int_0^\infty} f(t)(1+s\Delta t)^{-\frac{t}{\Delta t}-1}\Delta t\,](1+s\Delta t)^{\frac{t}{\Delta t}} ds \qquad (1.5\text{-}21)$$

$$f(t) = \frac{1}{2\pi j} \oint_{c_2} \left[ \int_0^\infty f(t)(1+s\Delta t)^{-\frac{t+\Delta t}{\Delta t}} \Delta t \right] (1+s\Delta t)^{\frac{t}{\Delta t}} ds \qquad (5.1\text{-}22)$$

Let

$$F(s) = \int_0^\infty f(t)(1+s\Delta t)^{-\frac{t+\Delta t}{\Delta t}} \Delta t \qquad (5.1\text{-}23)$$

Comment – For $\Delta t \rightarrow 0$ Eq 5.1-23 is the Laplace Transform

$$\lim_{\Delta t \rightarrow 0} \int_0^\infty f(t)(1+s\Delta t)^{-\frac{t+\Delta t}{\Delta t}} \Delta t = \int_0^\infty f(t) e^{-st} dt$$

Eq 5.1-23 has arbitrarily been named the $K_{\Delta t}$ Transform. The subscript, $\Delta t$, represents the value of the sampling (sample and hold) interval. This transform is closely related to the $J_{\Delta t}$ Transform, $F(s) = (1+s\Delta t)^{-1}\mathsf{F}(s)$.

Substituting Eq 5.1-23 into Eq 5.1-22

The Inverse $K_{\Delta t}$ Transform is:

$$f(t) = \frac{1}{2\pi j} \oint_{c_2} F(s)(1+s\Delta t)^{\frac{t}{\Delta t}} ds \qquad (5.1\text{-}24)$$

The $K_{\Delta x}$ Transform and Inverse $K_{\Delta x}$ Transform have been derived.

From Eq 5.1-23 and Eq 5.1-24

**The $K_{\Delta t}$ Transform is:**

$$\mathbf{F(s) = K_{\Delta t}[f(t)] = \int_0^\infty f(t)(1+s\Delta t)^{-\frac{t+\Delta t}{\Delta t}} \Delta t} \qquad (5.1\text{-}25)$$

**The Inverse $K_{\Delta t}$ Transform is:**

$$\mathbf{f(t) = K_{\Delta t}^{-1}[F(s)] = \frac{1}{2\pi j} \oint_{c_2} F(s)(1+s\Delta t)^{\frac{t}{\Delta t}} ds} \qquad (5.1\text{-}26)$$

**where**

$$s = \frac{e^{(\gamma+jw)\Delta t} - 1}{\Delta t}$$ **, the complex plane circular contour, $c_2$**

**The contour $c_2$ is a circular contour in the S Plane**

**$\gamma$ = positive real constant**

$$-\frac{\pi}{\Delta t} \le w < \frac{\pi}{\Delta t}$$

**f(t) = function of t**

**$F(s) = K_{\Delta t}[f(t)]$ , a function of s**

**$\Delta t$ = sampling interval, t increment**

**$t = n\Delta t$ , n = 0,1,2,3,…**

**$t = 0, \Delta t, 2\Delta t, 3\Delta t, \ldots$**

Comparing Eq 5.1-19 with Eq 5.1-25

$$J_{\Delta t}[f(t)] = (1+s\Delta t)K_{\Delta t}[f(t)] \qquad (5.1\text{-}27)$$

The previously derived Transforms are presented in Table 5.1-1 on the following page.

**Table 5.1-1**

<table>
<tr>
<td>**1**</td>
<td>

The $K_{\Delta t}$ Transform

-

$$F(s) = K_{\Delta t}[f(t)] = \int_{0}^{\infty}{}_{\Delta t} [1+s\Delta t]^{-\frac{t+\Delta t}{\Delta t}} f(t)\, \Delta t$$

Complex Plane

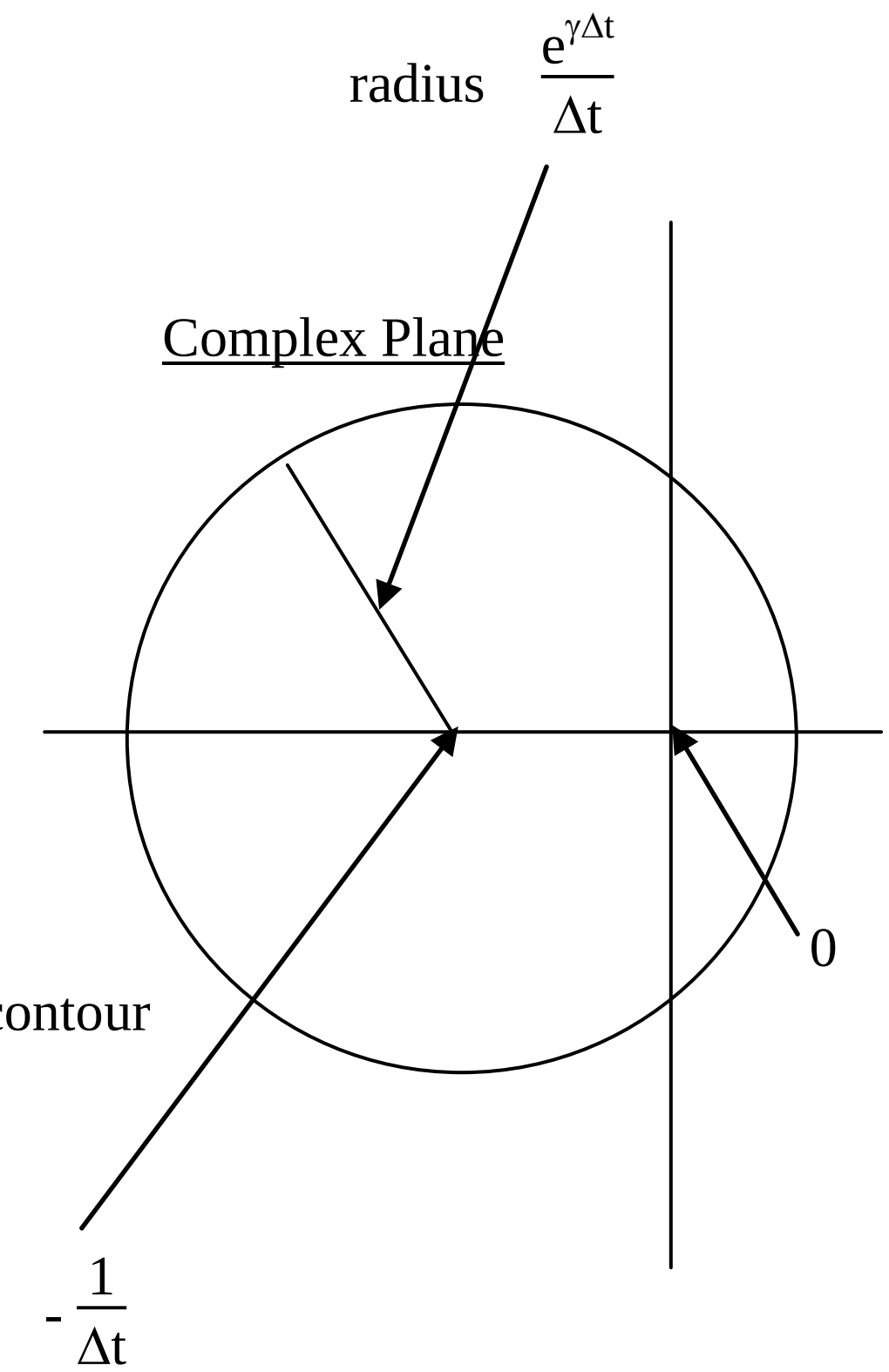

The Inverse $K_{\Delta x}$ Transform

$$f(t) = K_{\Delta t}^{-1}[F(s)] = \frac{1}{2\pi j} \oint_{c} [1+s\Delta t]^{\frac{t}{\Delta t}} F(s)\, ds$$

$$s = \frac{e^{(\gamma+jw)\Delta t} - 1}{\Delta t}\,, \quad -\frac{\pi}{\Delta t} \le w < \frac{\pi}{\Delta t}$$

where

$t = m\Delta t$ , $m = 0,1,2,3,\ldots$

$\Delta t$ = sampling interval, t increment

$f(t)$ = function of t

$f(t) = 0$ for $t < 0$

$F(s) = K_{\Delta t}[f(t)]$, a function of s

$\gamma$, w = real value constants

$\gamma > 0$

C, the complex plane contour of integration, is a circle of radius, $\frac{e^{\gamma\Delta t}}{\Delta t}$ , with center at $-\frac{1}{\Delta t}$

$\gamma$ is chosen so that the contour encloses all poles of F(s). If so, the f(t) complex plane closed contour integral can be evaluated using residue theorem methodology. The function of s, $[1+s\Delta t]^{\frac{t}{\Delta t}} F(s)$, for many F(s), can be expressed as a convergent Laurent Series for each pole as required by residue theory.

Notes - The $K_{\Delta t}$ Transform becomes the Laplace Transform for $\Delta t \to 0$. The c contour of the $K_{\Delta t}$ and $J_{\Delta t}$ Transforms is the same.

The $K_{\Delta t}$ Transform is related to the $J_{\Delta t}$ Transform, $K_{\Delta t}[f(t)] = \frac{1}{1+s\Delta t} J_{\Delta t}[f(t)]$.

</td>
<td>

The $K_{\Delta t}$ Transform is closely related to the $J_{\Delta x}$ Transform,

$K_{\Delta x}[f(t)] = \frac{J_{\Delta x}[f(t)]}{(1+s\Delta t)^{-}}$

</td>
</tr>
</table>

<table>
<tr>
<td>2</td>
<td>
<u>The $J_{\Delta t}$ Transform</u><br><br>
$$F(s) = J_{\Delta t}[f(t)] = \int_{0}^{\infty} [1+s\Delta t]^{-\frac{t}{\Delta t}} f(t)\,\Delta t$$
<br><br>
<u>The Inverse $J_{\Delta t}$ Transform</u><br>
$$f(t) = J_{\Delta t}^{-1}[F(s)] = \frac{1}{2\pi j}\oint_{c} [1+s\Delta t]^{\frac{t}{\Delta t}-1} F(s)\,ds$$
$$s = \frac{e^{(\gamma+jw)\Delta t} - 1}{\Delta t}, \quad -\frac{\pi}{\Delta t} \le w < \frac{\pi}{\Delta t}$$
where<br>
$t = m\Delta t$ , $m = 0,1,2,3,\ldots$<br>
$\Delta t$ = sampling period, t increment<br>
$f(t)$ = function of t<br>
$F(s) = J_{\Delta t}[f(t)]$ , a function of s<br>
$f(t) = 0$ for $t < 0$<br>
$\gamma$, w = real value constants<br>
$\gamma > 0$<br>
C, the complex plane contour of integration, is a circle of radius, $\frac{e^{\gamma\Delta t}}{\Delta t}$ , with center at $-\frac{1}{\Delta t}$<br><br>
$\gamma$ is chosen so that the contour encloses all poles of F(s). If so, the f(t) complex plane closed contour integral can be evaluated using residue theorem methodology. The function of s, $[1+s\Delta t]^{\frac{t}{\Delta t}-1} F(s)$, for many F(s), can be expressed as a convergent Laurent Series for each pole as required by residue theory.<br><br>
<u>Notes</u> - The $J_{\Delta t}$ Transform becomes the Laplace Transform for $\Delta t \to 0$.<br>
The c contour of the $J_{\Delta t}$ and $K_{\Delta t}$ Transforms is the same. See the c contour diagram in Section 1 above.<br><br>
The $J_{\Delta t}$ Transform is related to the $K_{\Delta t}$ Transform, $J_{\Delta t}[f(t)] = (1+s\Delta t)\,K_{\Delta t}[f(t)]$.
</td>
<td>The $J_{\Delta t}$ Transform is closely related to the $K_{\Delta x}$ Transform,<br>$J_{\Delta x}[f(t)] = (1+s\Delta x)K_{\Delta x}[f(t)]$</td>
</tr>
</table>

<u>Comment</u> - There is another way to derive the $J_{\Delta t}$ and $K_{\Delta t}$ Transform equations. Rather than obtaining the $J_{\Delta t}$ and $K_{\Delta t}$ Transform equations from the Z Transform, the Laplace Transform is used. This derivation method is presented in the following example, Example 5.1-1.

Example 5.1 Derive the $J_{\Delta t}$ Transform and $K_{\Delta t}$ Transform equations from the Laplace Transform using Interval Calculus discrete integration.

Derive the $J_{\Delta t}$ Transform from the Laplace Transform using Interval Calculus discrete integration.

$$G(p) = \int_{0}^{\infty} g(t)e^{-pt}\,dt \quad , \quad \text{The Laplace Transform} \qquad 1)$$

where

$g(t)$ and $e^{-pt}$ are continuous functions of t

Note that the Laplace Transform continuous integral of Eq 1 is the area under the continuous function, $g(t)e^{-pt}$, where $0 \le t < \infty$.

Consider f(t) to be a discrete sample and hold function where any change of value occurs only at $t = 0, \Delta t, 2\Delta t, 3\Delta t, \ldots, \infty$. $\Delta t$ = t interval.

For example

f(t) A Discrete Sample and Hold Waveform

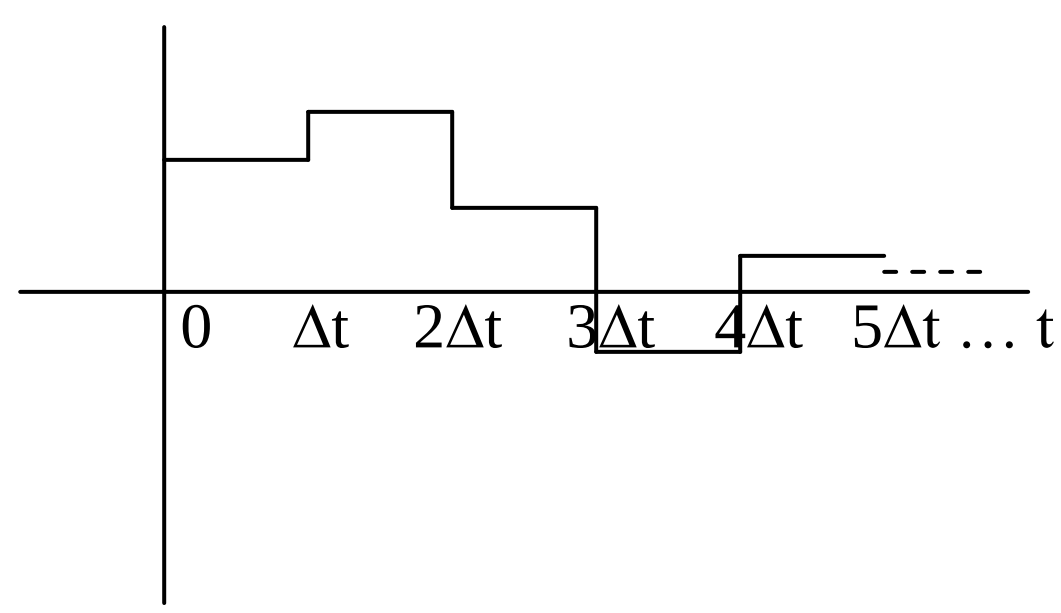

From Eq 1 using Interval Calculus discrete integration

$$F(p) = {}_{\Delta t}\int_{0}^{\infty} f(t)e^{-pt}\Delta t = \Delta t \sum_{n=0}^{\infty} f(n\Delta t)e^{-pn\Delta t} = T \sum_{n=0}^{\infty} f(n\Delta t)z^{-n} \qquad 2)$$

where

$\Delta t$ = t interval

$t = 0, \Delta t, 2\Delta t, 3\Delta t, \ldots, \infty$

$T = \Delta t$

f(t) = discrete sample and hold waveform function

$z = e^{p\Delta t}$

F(p) = Laplace Transform of the discrete sample and hold function, f(t).

$$F(z) = \sum_{n=0}^{\infty} f(n\Delta t)z^{-n} \quad , \quad \text{The Z Transform}$$

Note that the Interval Calculus discrete integral of Eq 2, the Laplace Transform of a discrete sample and hold waveform function, f(t), is the area under the discrete function, $f(t)e^{-pt}$, where $t = 0, \Delta t, 2\Delta t, 3\Delta t, \ldots, \infty$.

Note, also, that the Interval Calculus discrete integral of Eq 2 is closely related to the Z Transform.

For the convenience of discrete integration, use the following discrete function identity:

$$e^{-pt} = [1+(\frac{e^{p\Delta t}-1}{\Delta t})\Delta t]^{-\frac{t}{\Delta t}} \tag{3}$$

Let

$$s = \frac{e^{p\Delta t}-1}{\Delta t} \tag{4}$$

Substituting Eq 4 into Eq 3

$$e^{-pt} = [1+s\Delta t]^{-\frac{t}{\Delta t}} \tag{5}$$

From Eq 2 and Eq 5

$$J(s) = {}_{\Delta t}\int_{0}^{\infty} f(t)\ [1+s\Delta t]^{-\frac{t}{\Delta t}}\ \Delta t \tag{6}$$

The $J_{\Delta t}$ Transform, J(s), is related to the Laplace Transform, F(p), through the relationship, $s = \frac{e^{p\Delta t}-1}{\Delta t}$.

Then

**The $J_{\Delta t}$ Transform is:**

$$\mathbf{J_{\Delta t}[f(t)] = J(s) = {}_{\Delta t}\int_{0}^{\infty} f(t)\,[1+s\Delta t]^{-\frac{t}{\Delta t}}\,\Delta t}\ ,\quad \textbf{The } \mathbf{J_{\Delta t}} \textbf{ Transform} \tag{7}$$

**where**

**f(t) = discrete sample and hold waveform function**
**Δt = t interval**
**t = nΔt**
**n = 0, 1, 2, 3, …**

Find the function, f(t), from the $J_{\Delta t}$ Transform, J(s).

From Eq 7

$$J_{\Delta t}[f(t)] = J(s) = {}_{\Delta t}\int_{0}^{\infty} f(t)\ [1+s\Delta t]^{-\frac{t}{\Delta t}}\ \Delta t \equiv {}_{\Delta t}\sum_{t=0}^{\infty} f(t)\ [1+s\Delta t]^{-\frac{t}{\Delta t}}\ \Delta t\ = \sum_{n=0}^{\infty} f(nt)\Delta t\ [1+s\Delta t]^{-\frac{n\Delta t}{\Delta t}} \tag{8}$$

$$J_{\Delta t}[f(t)] = J(s) = {}_{\Delta t}\int_{0}^{\infty} f(t)\,[1+s\Delta t]^{-\frac{t}{\Delta t}}\,\Delta t \equiv \Delta t\,{}_{\Delta t}\sum_{t=0}^{\infty} f(t)\,[1+s\Delta t]^{-\frac{t}{\Delta t}} = \sum_{n=0}^{\infty} f(nt)\Delta t\,[1+s\Delta t]^{-n} \qquad 9)$$

where

$f(n\Delta t) = f(t)$ = discrete sample and hold waveform function

$\Delta t$ = t interval

$t = n\Delta t$

$n = 0, 1, 2, 3, \ldots$

Expanding Eq 9

$$J_{\Delta t}[f(t] = J(s) = f(0)\Delta t + f(\Delta t)\Delta t\,[1+s\Delta t]^{-1} + f(2\Delta t)\Delta t\,[1+s\Delta t]^{-2} + f(3\Delta t)\Delta t\,[1+s\Delta t]^{-3} + \ldots \qquad 10)$$

$$J_{\Delta t}[f(t] = J(s) = c_0 + c_1\,[1+s\Delta t]^{-1} + c_2\,[1+s\Delta t]^{-2} + c_3\,[1+s\Delta t]^{-3} + \ldots \qquad 11)$$

$$f(n\Delta t) = \frac{c_n}{\Delta t}\ , \quad n = 0,1,2,3, \ldots \qquad 12)$$

Eq 11 and Eq 12 above provide a means to evaluate from a $J_{\Delta t}$ Transform, J(s), the function which generated it, f(t).

Thus

**Evaluating f(t) from J(s) = $J_{\Delta t}$[f(t)] using a J(s) expansion**

$$\mathbf{J(s) = c_0 + c_1\,[1+s\Delta t]^{-1} + c_2\,[1+s\Delta t]^{-2} + c_3\,[1+s\Delta t]^{-3} + \ldots} \qquad 13)$$

$$\mathbf{f(t) = f(n\Delta t) = \frac{c_n}{\Delta t}\ , \quad n = 0,1,2,3, \ldots} \qquad 14)$$

**where**

**$c_n$ = series constants**

**$f(t) = f(n\Delta t)$ = discrete sample and hold waveform function represented by J(s), the $J_{\Delta t}$ Transform**

**$\Delta t$ = t increment**

**$t = n\Delta t$**

**$n = 0,1,2,3, \ldots$**

Find the Inverse $J_{\Delta t}$ Transform

Derive the Inverse $J_{\Delta t}$ Transform from Eq 10

Rewriting Eq 10

$$J_{\Delta t}[f(t)] = J(s) = f(0)\Delta t + f(\Delta t)\Delta t\,[1+s\Delta t]^{-1} + f(2\Delta t)\Delta t\,[1+s\Delta t]^{-2} + f(3\Delta t)\Delta t\,[1+s\Delta t]^{-3} + \dots$$

Multiply both sides of Eq 10 by $[1+s\Delta t]^{n-1}$ where $n = 0,1,2,3, \dots$

$$[1+s\Delta t]^{n-1} J(s) = f(0)\Delta t[1+s\Delta t]^{n-1} + f(\Delta t)\Delta t\,[1+s\Delta t]^{n-2} + f(2\Delta t)\Delta t\,[1+s\Delta t]^{n-3}$$

$$+ f(3\Delta t)\Delta t\,[1+s\Delta t]^{n-4} + \dots \qquad 15)$$

Change the form of Eq 15

$$[1+s\Delta t]^{n-1} J(s) = \frac{f(0)}{\Delta t^{-n}} \frac{1}{[s+\frac{1}{\Delta t}]^{1-n}} + \frac{f(\Delta t)}{\Delta t^{1-n}} \frac{1}{[s+\frac{1}{\Delta t}]^{2-n}} + \frac{f(2\Delta t)}{\Delta t^{2-n}} \frac{1}{[s+\frac{1}{\Delta t}]^{3-n}} + \frac{f(3\Delta t)}{\Delta t^{3-n}} \frac{1}{[s+\frac{1}{\Delta t}]^{4-n}} + \dots \qquad 16)$$

Perform contour integration along the complex plane closed contour, C, which circles the pole at $-\frac{1}{\Delta x}$. The complex plane closed contour, C, is defined below.

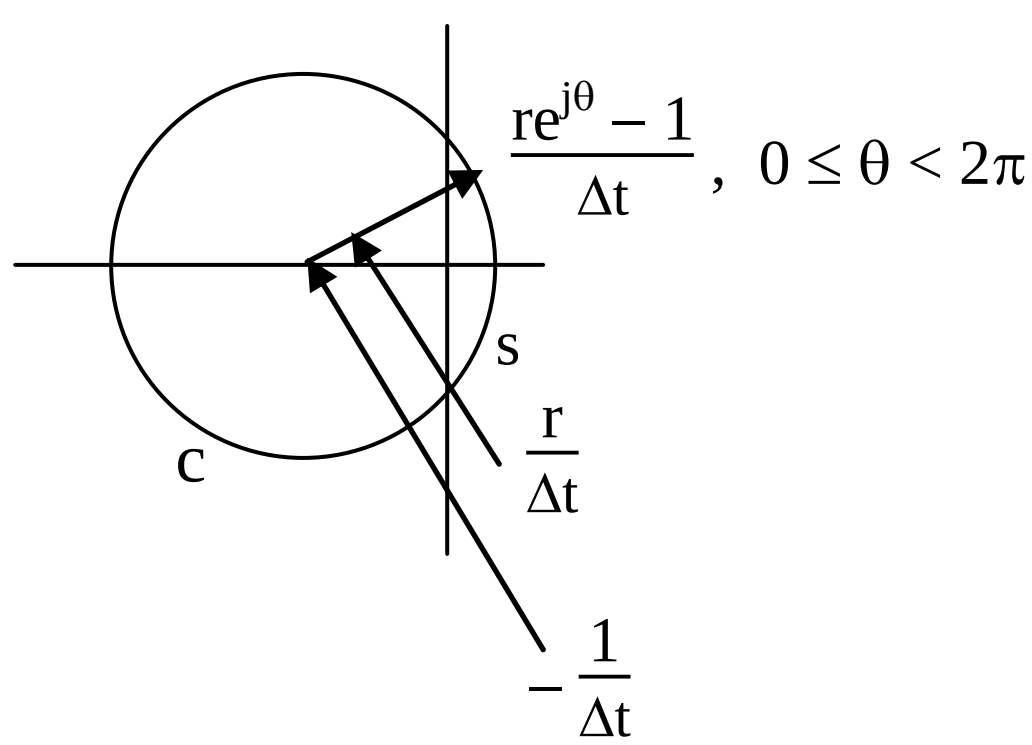

The constant, r, is chosen so that all poles of J(s) lie within the circular complex plane closed contour, c.

$$\frac{1}{2\pi j} \oint_c [1+s\Delta t]^{n-1} J(s)ds = \frac{f(0)}{\Delta t^{-n}} \frac{1}{2\pi j} \oint_c \frac{1}{[s+\frac{1}{\Delta t}]^{1-n}} ds + \frac{f(\Delta t)}{\Delta t^{1-n}} \frac{1}{2\pi j} \oint_c \frac{1}{[s+\frac{1}{\Delta t}]^{2-n}} ds$$

$$+ \frac{f(2\Delta t)}{\Delta t^{2-n}} \frac{1}{2\pi j} \oint_c \frac{1}{[s+\frac{1}{\Delta t}]^{3-n}} ds + \frac{f(3\Delta t)}{\Delta t^{3-n}} \frac{1}{2\pi j} \oint_c \frac{1}{[s+\frac{1}{\Delta t}]^{4-n}} ds + \dots \qquad 17)$$

Using Cauchy's Integral Formulas, evaluate Eq 17 for n = 0, 1, 2, 3, …

Cauchy's Integral Formulas are:

$$\oint_c P(s)ds = 0 \qquad 18)$$

$$\frac{1}{2\pi j}\oint_c \frac{P(s)}{(s-a)^{r+1}}ds = \frac{1}{r!}\frac{d^r}{ds^r}P(s)\Big|_{s=a} \qquad 19)$$

where
r = 0, 1, 2, 3, …
P(s) has no poles within the closed contour, c
The pole at s = a is within the closed contour, c

Substituting n = 0 into Eq 17

$$\frac{1}{2\pi j}\oint_c [1+s\Delta t]^{0-1}J(s)ds = \frac{f(0)}{1}\frac{1}{2\pi j}\oint_c \frac{1}{[s+\frac{1}{\Delta t}]}ds + \frac{f(\Delta t)}{\Delta t}\frac{1}{2\pi j}\oint_c \frac{1}{[s+\frac{1}{\Delta t}]^2}ds + \frac{f(2\Delta t)}{\Delta t^2}\frac{1}{2\pi j}\oint_c \frac{1}{[s+\frac{1}{\Delta t}]^3}ds$$

$$+\frac{f(3\Delta t)}{\Delta t^3}\frac{1}{2\pi j}\oint_c \frac{1}{[s+\frac{1}{\Delta t}]^4}ds + \ldots = f(0) + 0 + 0 + 0 + \ldots \qquad 20)$$

$$\frac{1}{2\pi j}\oint_c [1+s\Delta t]^{0-1}J(s)ds = f(0) \qquad 21)$$

Substituting n = 1 into Eq 17

$$\frac{1}{2\pi j}\oint_c [1+s\Delta t]^{1-1}J(s)ds = \frac{f(0)}{\Delta t^{-1}}\frac{1}{2\pi j}\oint_c ds + f(\Delta t)\frac{1}{2\pi j}\oint_c \frac{1}{[s+\frac{1}{\Delta t}]}ds + \frac{f(2\Delta t)}{\Delta t}\frac{1}{2\pi j}\oint_c \frac{1}{[s+\frac{1}{\Delta t}]^2}ds$$

$$+\frac{f(3\Delta t)}{\Delta t^2}\frac{1}{2\pi j}\oint_c \frac{1}{[s+\frac{1}{\Delta t}]^3}ds + \ldots = 0 + f(\Delta t) + 0 + 0 + \ldots \qquad 22)$$

$$\frac{1}{2\pi j}\oint_c [1+s\Delta t]^{1-1}J(s)ds = f(\Delta t) \qquad 23)$$

Substituting n = 2 into Eq 17

$$\frac{1}{2\pi j}\oint_c [1+s\Delta t]^{2-1} J(s)ds = \frac{f(0)}{\Delta t^{-2}}\frac{1}{2\pi j}\oint_c [s+\frac{1}{\Delta t}]\, ds + \frac{f(\Delta t)}{\Delta t^{-1}}\frac{1}{2\pi j}\oint_c ds + f(2\Delta t)\frac{1}{2\pi j}\oint_c \frac{1}{[s+\frac{1}{\Delta t}]} ds$$

$$+\frac{f(3\Delta t)}{\Delta t}\frac{1}{2\pi j}\oint_c \frac{1}{[s+\frac{1}{\Delta t}]^2} ds + \ldots = 0 + 0 + f(2\Delta t) + 0 + 0 + \ldots \qquad 24)$$

$$\frac{1}{2\pi j}\oint_c [1+s\Delta t]^{2-1} J(s)ds = f(2\Delta t) \qquad 25)$$

.

.

.

Observing the previous evaluations of the contour integral, $\frac{1}{2\pi j}\oint_c [1+s\Delta t]^{n-1} J(s)ds$, for n = 0, 1, 2, 3, ...

$$\frac{1}{2\pi j}\oint_c [1+s\Delta t]^{n-1} J(s)ds = f(n\Delta t) \qquad 26)$$

Then

**The Inverse $J_{\Delta t}$ Transform is:**

$$\mathbf{f(n\Delta t) = \frac{1}{2\pi j}\oint_c J(s)[1+s\Delta t]^{n-1} ds} \qquad 27)$$

**or**

$$\mathbf{f(t) = \frac{1}{2\pi j}\oint_c J(s)[1+s\Delta t]^{\frac{t}{\Delta t}-1} ds} \qquad 28)$$

**where**

**$J(s) = J_{\Delta t}[f(n\Delta t)]$ , The $J_{\Delta t}$ Transform**
**$f(n\Delta t) = f(t)$ = discrete sample and hold waveform function**
**$\Delta t$ = t increment**
**$t = n\Delta t$**
**n = 0, 1, 2, 3, …**

$$s = \frac{re^{j\theta} - 1}{\Delta t}, \; 0 \le \theta < 2\pi$$

**r = a constant such that all of the poles of J(s) lie within the circular complex plane contour, c**

**c = the circular complex plane closed contour which is shown below.**

**The complex plane closed contour, C, is defined below.**

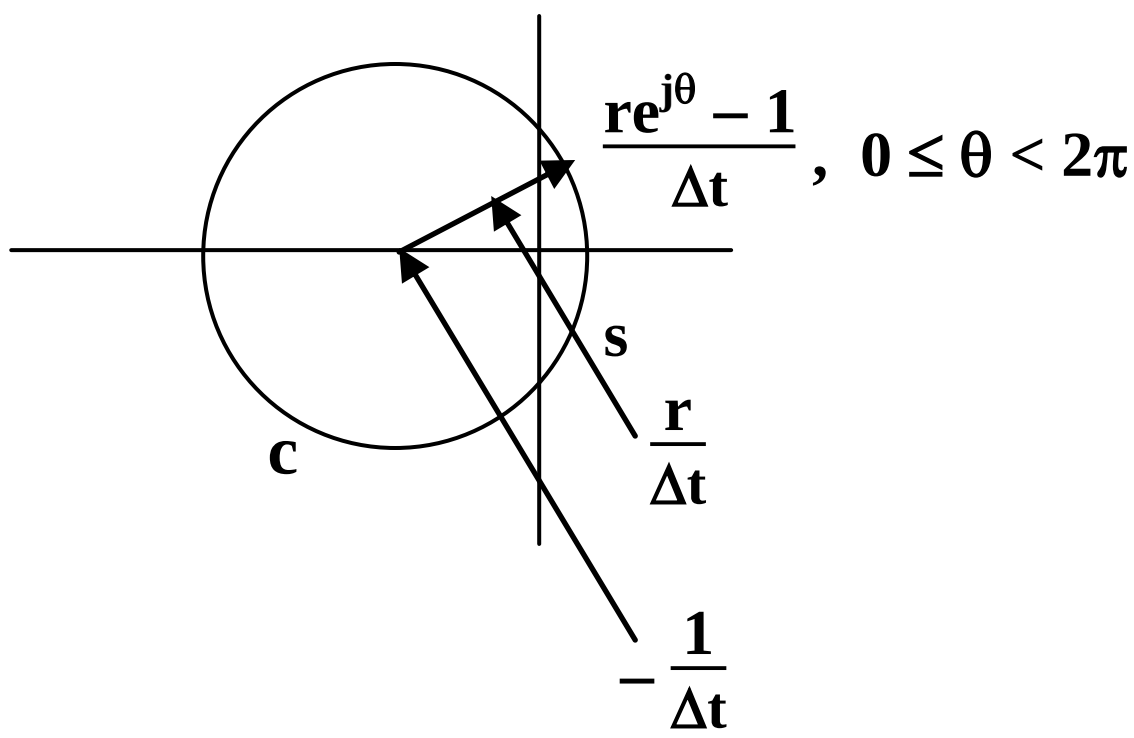

To calculate the closed contour integral of Eq 27 or the closed contour integral of Eq 28 where $n = \frac{t}{\Delta t}$,

the theory of residues can be used provided that the function of s, $J(s)[1+s\Delta t]^{n-1}$, can be expanded into a convergent Laurent Series for each pole of J(s)

$$f(n\Delta t) = J_{\Delta t}^{-1}J_{\Delta t}[f(n\Delta t)] = J_{\Delta t}^{-1}J(s) = \frac{1}{2\pi j}\oint_c J(s)(1+s\Delta t)^{n-1}ds = \sum_{p=1}^{P} R_p \qquad 29)$$

$f(n\Delta t) = J_{\Delta t}^{-1}J_{\Delta t}[f(n\Delta t)]$ = sum of the residues of $J(s)[1+s\Delta t]^{n-1}$ at all of the poles of J(s), all of which lie within the circular complex plane closed contour, c, with center at $-\frac{1}{\Delta x}$ and radius, $\frac{r}{\Delta t}$.

$$R = \lim_{s\to r}\frac{1}{(m-1)!}\frac{d^{m-1}}{ds^{m-1}}[(s-r)^m J(s)(1+s\Delta t)^{n-1}]$$ , the residue calculation formula for a pole at s = r 30)

where

$R$ = the residue of a pole of $J(s)(1+s\Delta t)^{n-1}$

$\sum_{p=1}^{P} R_p$ = the sum of the P residues of the P poles of $J(s)(1+s\Delta t)^{n-1}$

Derive the $K_{\Delta t}$ Transform from the $J_{\Delta t}$ Transform.

By definition

$$\mathbf{K_{\Delta t}[f(t)] = K(s) = (1+s\Delta t)^{-1}J_{\Delta t}[f(t)] = (1+s\Delta t)^{-1}J(s)} \qquad 31)$$

Multiplying both sides of Eq 6 by $(1+s\Delta t)^{-1}$

$$(1+s\Delta t)^{-1}J_{\Delta t}[f(t)] = (1+s\Delta t)^{-1}J(s) = {}_{\Delta t}\int_0^{\infty} f(t)(1+s\Delta t)^{-1}[1+s\Delta t]^{-\frac{t}{\Delta t}}\Delta t \qquad 32)$$

From Eq 31 and Eq 32

$$K_{\Delta t}[f(t)] = K(s) = {}_{\Delta t}\int_0^{\infty} f(t)\,[1+s\Delta t]^{-\frac{t+\Delta t}{\Delta t}}\,\Delta t \qquad 33)$$

Then

**The $K_{\Delta t}$ Transform is:**

$$\mathbf{K_{\Delta t}[f(t)] = K(s) = {}_{\Delta t}\int_0^{\infty} f(t)\,[1+s\Delta t]^{-\frac{t+\Delta t}{\Delta t}}\,\Delta t} \ , \quad \textbf{The } \mathbf{K_{\Delta t}} \textbf{ Transform} \qquad 34)$$

**where**

**$f(t)$ = discrete sample and hold waveform function**

**$\Delta t$ = t interval**

**$t = 0, \Delta t, 2\Delta t, 3\Delta t, \ldots$**

Find the Inverse $K_{\Delta t}$ Transform

Derive the Inverse $K_{\Delta t}$ Transform from the Inverse $J_{\Delta t}$ Transform

Rewriting Eq 26

$$f(n\Delta t) = \frac{1}{2\pi j}\oint_c [1+s\Delta t]^{n-1}J(s)ds \qquad 35)$$

From Eq 31

$$J(s) = (1+s\Delta t)K(s) \qquad 36)$$

Substitute Eq 36 into Eq 35

$$f(n\Delta t) = \frac{1}{2\pi j} \oint_c [1+s\Delta t]^{n-1}(1+s\Delta t)K(s)ds \qquad 37)$$

$$f(n\Delta t) = \frac{1}{2\pi j} \oint_c [1+s\Delta t]^{n} K(s)ds = \frac{1}{2\pi j} \oint_c [1+s\Delta t]^{\frac{t}{\Delta t}} K(s)ds \qquad 38)$$

where
$t = n\Delta t$
$n = 0, 1, 2, 3, \ldots$

Then

**The Inverse $K_{\Delta t}$ Transform is:**

$$f(n\Delta t) = K_{\Delta t}^{-1}[K(s)] = \frac{1}{2\pi j} \oint_c K(s)[1+s\Delta t]^{n} ds \qquad 39)$$

**or**

$$f(t) = K_{\Delta t}^{-1}[K(s)] = \frac{1}{2\pi j} \oint_c K(s)(1+s\Delta t)^{\frac{t}{\Delta t}} ds \qquad 40)$$

**where**

**$K(s) = K_{\Delta t}[f(t)]$ , The $K_{\Delta t}$ Transform**
**$f(n\Delta t) = f(t)$ = discrete sample and hold waveform function**
**$\Delta t$ = t increment**
**$n = 0, 1, 2, 3, \ldots$**
**$s = \frac{re^{j\theta} - 1}{\Delta t}$ , $0 \le \theta < 2\pi$**
**r = a constant such that all of the poles of J(s) lie within the circular complex plane closed contour, c**
**c = The circular closed contour in the complex plane which is shown below.**

**S Plane** **Circular Contour**

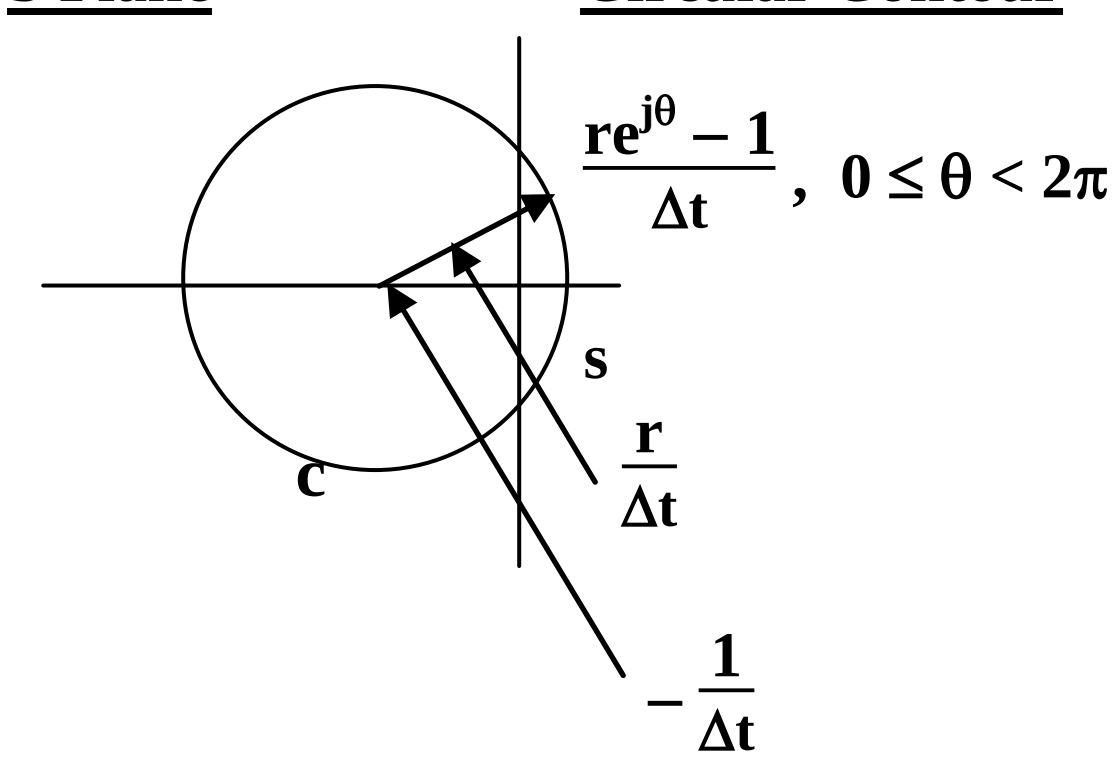

To calculate the closed contour integral of Eq 39 where $n = \frac{t}{\Delta t}$ or the closed contour integral of Eq 40,

the theory of residues can be used provided that the function of s, $K(s)[1+s\Delta t]^{\frac{t}{\Delta t}}$, can be expanded into a convergent Laurent Series for each pole of K(s)

$$f(n\Delta t) = K_{\Delta t}^{-1}K_{\Delta t}[f(n\Delta t)] = K_{\Delta t}^{-1}K(s) = \frac{1}{2\pi j}\oint_c K(s)(1+s\Delta t)^{\frac{t}{\Delta t}}ds = \sum_{p=1}^{P} R_p \qquad 41)$$

$f(n\Delta t) = K_{\Delta t}^{-1}K_{\Delta t}[f(n\Delta t)]$ = sum of the residues of $K(s)[1+s\Delta t]^{\frac{t}{\Delta t}}$at all of the poles of K(s), all of which lie within the circular complex plane contour, c, with center at $-\frac{1}{\Delta x}$ and radius, $\frac{r}{\Delta t}$ .

$$R = \lim_{s\to r}\frac{1}{(m-1)!}\frac{d^{m-1}}{ds^{m-1}}[(s-r)^m K(s)(1+s\Delta t)^{\frac{t}{\Delta t}}]$$ , the residue calculation formula for a pole at s = r 42)

where

$R$ = the residue of a pole of $K(s)(1+s\Delta t)^{\frac{t}{\Delta t}}$

$\sum_{p=1}^{P} R_p$ = the sum of the P residues of the P poles of $K(s)(1+s\Delta t)^{\frac{t}{\Delta t}}$

$t = 0, \Delta t, 2\Delta t, 3\Delta t, \ldots$

## Section 5.2: Investigating the characteristics of the $J_{\Delta x}$ and $K_{\Delta x}$ Transforms

Similarities to the Laplace Transform

As can be seen from Eq 5.1-27, there is only a small difference between the $J_{\Delta t}$ Transform and the $K_{\Delta t}$ Transform, a factor of $1+s\Delta t$. Either transform could be selected for use. However, the $K_{\Delta t}$ Transform is a better choice. The $K_{\Delta t}$ Transform has a greater similarity to the Laplace Transform. For some applications, this greater similarity to the Laplace Transform can be advantageous as in the case of Related Functions. These related functions, a Calculus function and an Interval Calculus function which have the same transform, can be used to solve differential difference equations. See Section 4.5 in Chapter 4. To demonstate, consider two related functions, the Calculus function, $e^{at}$, and the discrete Interval Calculus function, $e_{\Delta t}(a,t)$.

The Laplace Transform of $e^{at}$ is $\dfrac{1}{s-a}$ .

The $K_{\Delta t}$ Transform of $e_{\Delta t}(a,t)$ is $\dfrac{1}{s-a}$ .

The $J_{\Delta t}$ Transform of $e_{\Delta t}(a,t)$ is $\dfrac{1+s\Delta t}{s-a}$ .

Note that the $K_{\Delta t}$ Transform of $e_{\Delta t}(a,t)$ has a greater similarity to the Laplace Transform than the $J_{\Delta t}$ Transform. It must be kept in mind that the Laplace Transform s and the Interval Calculus s, though they appear to be the same, with the exception of the unique case where $\Delta t \to 0$, they are not.

The Laplace Transform s is equal to $\gamma+jw$ , $-\infty < w < +\infty$ .

The Interval Calculus s is equal to $s = \dfrac{e^{(\gamma+jw)\Delta t} - 1}{\Delta t}$ , $-\dfrac{\pi}{\Delta t} \le w < \dfrac{\pi}{\Delta t}$ .

.

For $\Delta t \to 0$ both s are equal to $\gamma+jw$ where $-\infty < w < +\infty$ .

Further investigation shows another significant similarity between the $K_{\Delta t}$ and $J_{\Delta t}$ Transform and the Laplace Transform, in particular where $\Delta t \to 0$.

$$K_{\Delta t}[f(t)] = \int_{\Delta t\,0}^{\infty} f(t)(1+s\Delta t)^{-(\frac{t+\Delta t}{\Delta t})}\Delta t = \frac{1}{1+s\Delta t}\int_{\Delta t\,0}^{\infty} f(t)(1+s\Delta t)^{-\frac{t}{\Delta t}}\Delta t = \frac{1}{1+s\Delta t}J_{\Delta t}[f(t)] \qquad (5.2\text{-}1)$$

where

$t = n\Delta t$ , $n = 0,1,2,\ldots$

$f(t)$ = function of t

Taking the limit of Eq 5.2-1 as $\Delta t \rightarrow 0$

$$\lim_{\Delta t\rightarrow 0} K_{\Delta t}[f(t)] = \lim_{\Delta t\rightarrow 0} \int_{0}^{\infty} f(t)(1+s\Delta t)^{-(\frac{t+\Delta t}{\Delta t})}\Delta t \;=\; \lim_{\Delta t\rightarrow 0} \frac{1}{1+s\Delta t}\int_{0}^{\infty} f(t)(1+s\Delta t)^{-\frac{t}{\Delta t}}\Delta t =$$

$$\lim_{\Delta t\rightarrow 0} \frac{1}{1+s\Delta t} J_{\Delta t}[f(t)] \tag{5.2-2}$$

Find $\lim_{\Delta t\rightarrow 0} (1+s\Delta t)^{-\frac{t}{\Delta t}}$

$$\lim_{\Delta t\rightarrow 0} (1+s\Delta t)^{-\frac{t}{\Delta t}} = \lim_{\Delta t\rightarrow 0} e^{-t\frac{\ln(1+s\Delta t)}{\Delta t}} = \lim_{\Delta t\rightarrow 0} e^{-t\frac{s}{1+s\Delta t}} = e^{-st}$$

$$\lim_{\Delta t\rightarrow 0} (1+s\Delta t)^{-\frac{t}{\Delta t}} = e^{-st} \tag{5.2-3}$$

From Eq 5.2-2 and Eq 5.2-3

$$\lim_{\Delta t\rightarrow 0} K_{\Delta t}[f(t)] = \lim_{\Delta t\rightarrow 0} \frac{1}{1+s\Delta t} J_{\Delta t}[f(t)] = \int_{0}^{\infty} f(t)e^{-st}dt \tag{5.2-4}$$

Then

$$\lim_{\Delta t\rightarrow 0} K_{\Delta t}[f(t)] = \lim_{\Delta t\rightarrow 0} J_{\Delta t}[f(t)] = \int_{0}^{\infty} f(t)e^{-st}dt \;\;,\;\; \text{The Laplace Transform} \tag{5.2-5}$$

The limit as $\Delta t \rightarrow 0$ of both the $K_{\Delta t}$ Transform and the $J_{\Delta t}$ Transform is the Laplace Transform.

When using Interval Calculus function notation, even the appearance of the $K_{\Delta t}$ Transform and the $J_{\Delta t}$ Transform is similar.

Rewriting the $K_{\Delta t}$ and the $J_{\Delta t}$ Transforms

$$K_{\Delta t}[f(t)] = \int_{0}^{\infty} f(t)(1+s\Delta t)^{-(\frac{t+\Delta t}{\Delta t})}\Delta t \;,\quad \text{The } K_{\Delta t} \text{ Transform} \tag{5.2-6}$$

$$J_{\Delta t}[f(t)] = \int_{0}^{\infty} f(t)(1+s\Delta t)^{-\frac{t}{\Delta t}}\Delta t \;,\qquad \text{The } J_{\Delta t} \text{ Transform} \tag{5.2-7}$$

The Interval Calculus function, $e_{\Delta t}(a,-t)$, is defined as follows:

$$e_{\Delta t}(a,-t) = (1+s\Delta x)^{-\frac{t}{\Delta t}} \qquad (5.2\text{-}8)$$

From Eq 5.2-6, Eq 5.2-7, and Eq 5.2-8

$$K_{\Delta t}[f(t)] = {}_{\Delta t}\int_0^{\infty} f(t)e_{\Delta t}(s,-t-\Delta t)\Delta t \ , \quad \text{The } K_{\Delta t} \text{ Transform} \qquad (5.2\text{-}9)$$

$$J_{\Delta t}[f(t)] = {}_{\Delta t}\int_0^{\infty} f(t)e_{\Delta t}(s,-t)\Delta t \ , \quad \text{The } J_{\Delta t} \text{ Transform} \qquad (5.2\text{-}10)$$

$$L[f(t)] = \int_0^{\infty} f(t)e^{-st}dt \ , \quad \text{The Laplace Transform} \qquad (5.2\text{-}11)$$

Note the similarity of the above three transform equations

As previously stated, either the $K_{\Delta t}$ Transform or the $J_{\Delta t}$ Transform can be selected for use. Both transforms are derived from either the Z Transform or the Laplace Transform and are very similar. However, for the reasons previously mentioned, the selection of the $K_{\Delta t}$ Transform is a better choice. For this reason, in the discussions and derivations which follow, only the $K_{\Delta t}$ Transform will be considered.

## Section 5.3: Analyzing the $K_{\Delta t}$ Transform

Another form of the $K_{\Delta t}$ Transform can be written

Rewriting the derived $K_{\Delta t}$ Transform

$$F(s) = K_{\Delta t}[f(t)] = {}_{\Delta t}\int_0^{\infty} [1+s\Delta t]^{-\frac{t+\Delta t}{\Delta t}} f(t)\ \Delta t \qquad (5.3\text{-}1)$$

where
- $t = n\Delta t$ ,   $n = 0,1,2,3,\ldots$
- $\Delta t$ = sampling interval, t increment
- $f(t)$ = function of t
- $f(t) = 0$   for $t < 0$

From Eq 5.3-1

$$F(s) = K_{\Delta t}[f(n\Delta t)] = {}_{1}\sum_{n=0}^{\infty} [1+s\Delta t]^{-n-1}[f(n\Delta t)\Delta t] \qquad (5.3\text{-}2)$$

In Eq 5.3-2, the term [f(nΔt)Δt] is readily recognized as the nth area segment of a sample and hold waveform, f(t), obtained from sampling the function, f(t), where t = nΔt ,   n = 0,1,2,3,…
Note Diagram 5.3-1 below.

Diagram 5.3-1:  Functions of t

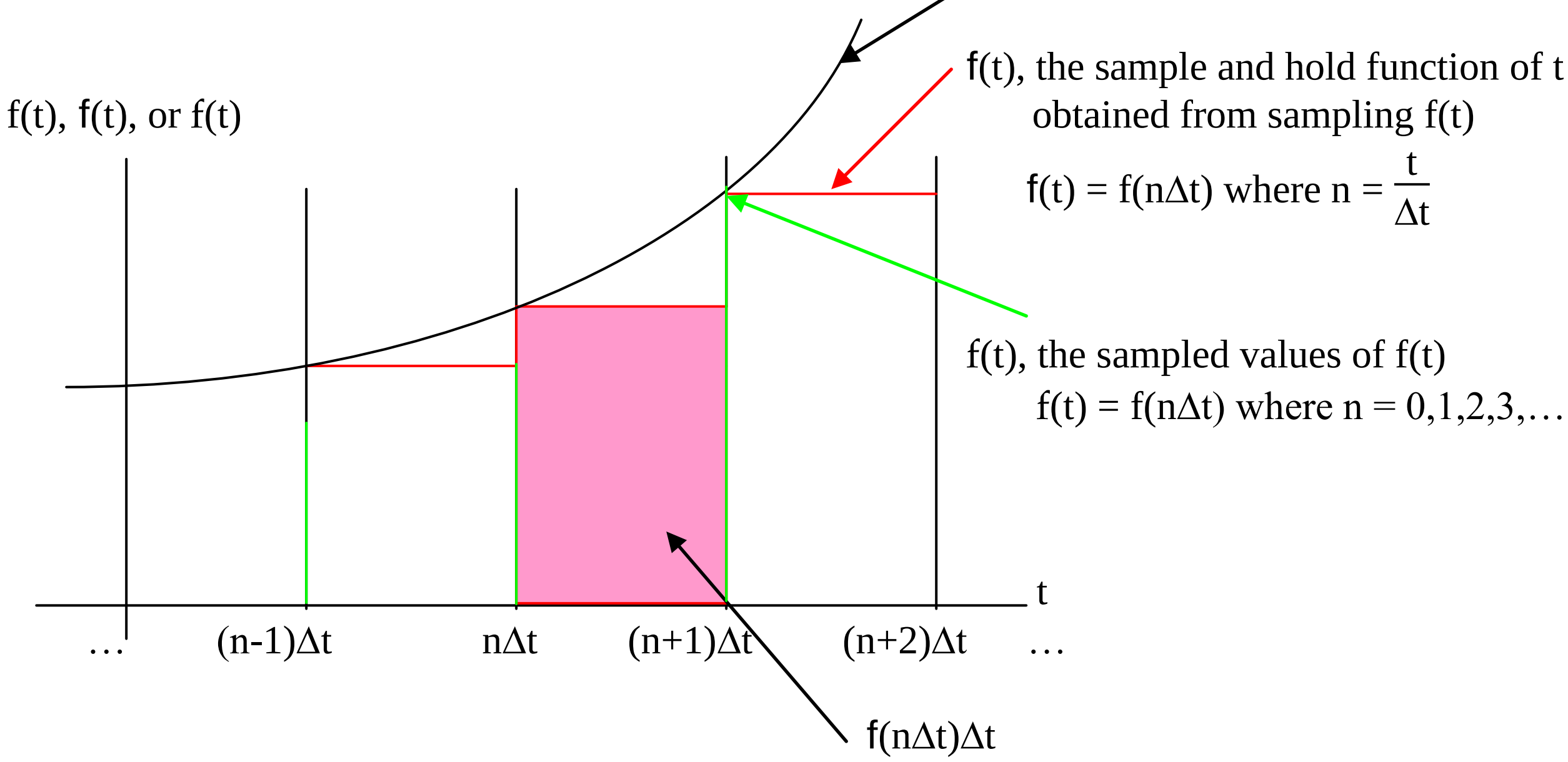

Comments – In Diagram 5.3-1 there are specified three different functions of t, f(t), f(t), and f(t). All of these functions have the same value for t = nΔt where n = 0,1,2,3,… When used by the $K_{\Delta t}$ Transform equations, these functions produce the same result since their values are the same at the sampling times, t = nΔt where n = 0,1,2,3…

The question arises, "What does the term, $[1+s\Delta t]^{-n-1}$, in Eq 5.3-2 represent?"

Consider a Unit Area Pulse and find its $K_{\Delta t}$ Transform. See Diagram 5.3-2 below.

Diagram 5.3-2:  The Unit Area Pulse

$$\text{Unit Area Pulse} = \frac{1}{\Delta t}\,\Delta t = 1$$

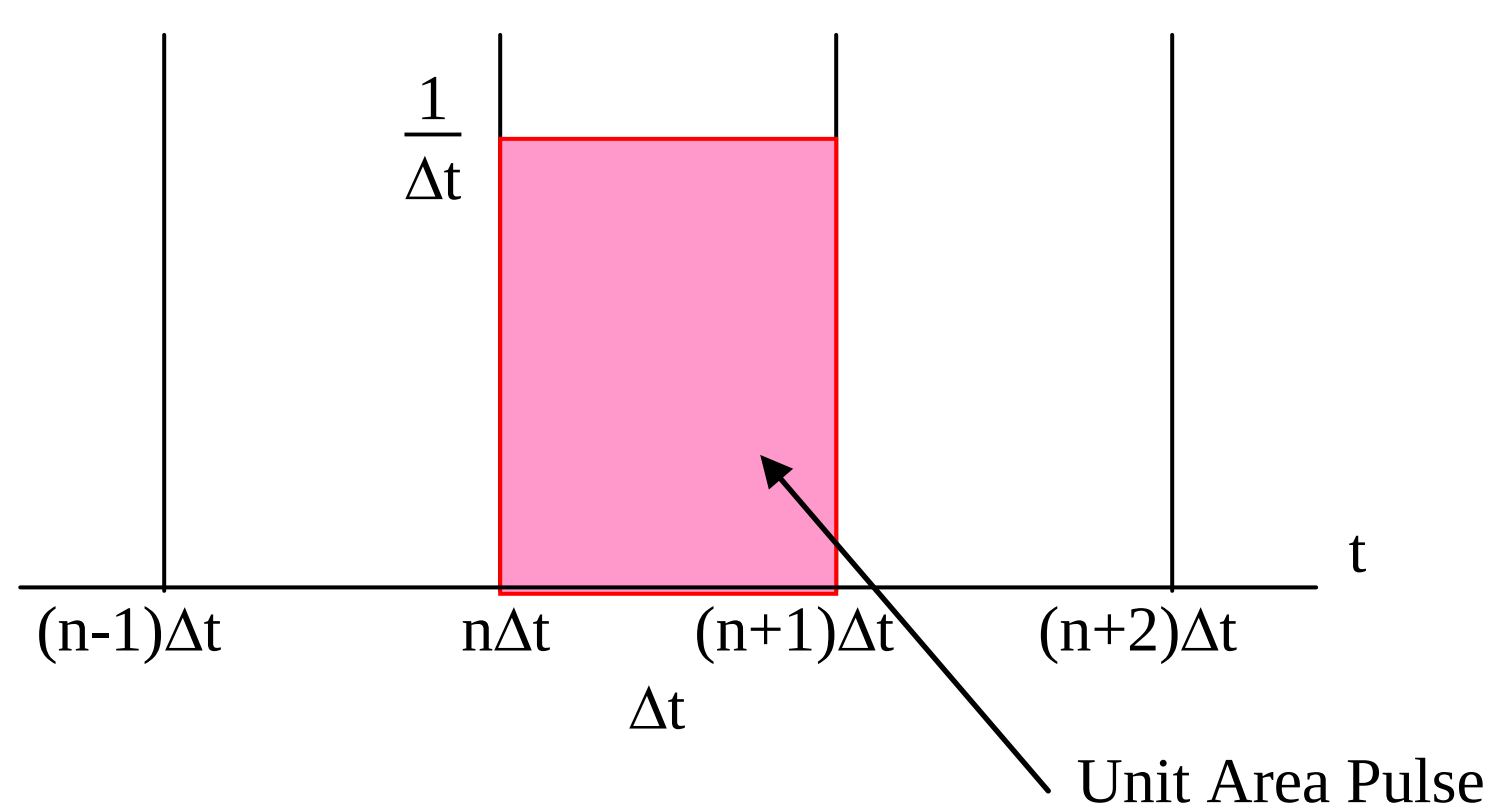

Let $f(t) = \frac{1}{\Delta t}$ , a constant value

Find the $K_{\Delta t}$ Transform of the function, $f(t) = \frac{1}{\Delta t}$ , from $t = n\Delta t$ to $(n+1)\Delta t$ , n = 0, 1, 2, 3, … , a unit area pulse.

$$K_{\Delta t}[\text{Unit Area Pulse}] = {}_{\Delta t}\int_{n\Delta t}^{(n+1)\Delta t} \frac{1}{\Delta t}\,[1+s\Delta t]^{-\frac{t+\Delta t}{\Delta t}}\,\Delta t \tag{5.3-3}$$

$$K_{\Delta t}[\text{Unit Area Pulse}] = \frac{1}{1+s\Delta t}\;{}_{\Delta t}\int_{n\Delta t}^{(n+1)\Delta t} \frac{1}{\Delta t}\,[1+s\Delta t]^{-\frac{t}{\Delta t}}\,\Delta t = \frac{1}{\Delta t(1+s\Delta t)}\;{}_{\Delta t}\int_{n\Delta t}^{(n+1)\Delta t} e_{\Delta t}(s,-t)\,\Delta t \tag{5.3-4}$$

From Table 6, the table of discrete integrals in the Appendix

$$_{\Delta x}\int e_{\Delta x}(a,-x)\Delta x = -\,\frac{1+a\Delta x}{a}e_{\Delta x}(a,-x) + k \tag{5.3-5}$$

Using Eq 5.3-5 to evaluate Eq 5.3-4

$a = s$
$x = t$

$$K_{\Delta t}[\text{Unit Area Pulse}] = \frac{1}{\Delta t(1+s\Delta t)}\,(-1)\,\frac{1+s\Delta t}{s}e_{\Delta t}(s,-t)\,\Big|_{n\Delta t}^{(n+1)\Delta t} \tag{5.3-6}$$

$$e_{\Delta t}(s,-t) = (1+s\Delta t)^{-\frac{t}{\Delta t}} \tag{5.3-7}$$

From Eq 5.3-6 and Eq 5.3-7

$$K_{\Delta t}[\text{Unit Area Pulse}] = -\,\frac{1}{s\Delta t}\,(1+s\Delta t)^{-\frac{t}{\Delta t}}\Big|_{n\Delta t}^{(n+1)\Delta t} = -\,\frac{1}{s\Delta t}\,[(1+s\Delta t)^{-(n+1)} - (1+s\Delta t)^{-n}\,] \tag{5.3-8}$$

$$K_{\Delta t}[\text{Unit Area Pulse}] = -\,\frac{1}{s\Delta t}\,(1+s\Delta t)^{-n}\,[\frac{1}{1+s\Delta t} - 1)] = -\,\frac{1}{s\Delta t}\,(1+s\Delta t)^{-n}\,[\frac{-s\Delta t}{1+s\Delta t}] \tag{5.3-9}$$

$K_{\Delta t}$[Unit Area Pulse] = $(1+s\Delta t)^{-n-1}$ (5.3-10)

where

$n = \dfrac{t}{\Delta t}$

n = 0,1,2,3,…

t = 0, Δt, 2Δt, 3Δt,…

Δt = sampling interval, t increment

From Eq 5.3-10, the term, $(1+s\Delta t)^{-n-1}$ or $(1+s\Delta t)^{-\frac{t+\Delta t}{\Delta t}}$, has been identified as the $K_{\Delta t}$ Transform of a unit area pulse positioned between t = nΔt and t = (n+1)Δt.

Thus, the $K_{\Delta x}$ Transform of Eq 5.3-1 or Eq 5.3-2 is seen to be a summation (or equivalently, a discrete integral) of the product of two terms over each discrete time interval, t = 0,Δt,2Δt,3Δt,…, ∞. One term is the $K_{\Delta t}$ Transform of a unit area pulse in the time interval. The other term is the area under the discrete sample and hold function, f(t), in the time interval. By expanding the $K_{\Delta t}$ Transform of Eq 5.3-2, one can evaluate the values of f(t) for all t where t = 0,Δt,2Δt,3Δt,…

<u>Comment</u> - In later explanations, the term $(1+s\Delta t)^{-n-1}\Delta t$ or $(1+s\Delta t)^{-\frac{t+\Delta t}{\Delta t}}\Delta t$ may be used. These terms are equivalent and are the $K_{\Delta t}$ Transform of a unit amplitude pulse of width Δt initiated at t = nΔt and ending at t = (n+1)Δt. n = 0,1,2,3,… , t = nΔt

Finding the values of the function, f(t), from its $K_{\Delta t}$ Transform

From Eq 5.3-2

$$F(s) = K_{\Delta t}[f(t)] = {}_1\sum_{n=0}^{\infty} [1+s\Delta t]^{-n-1} f(n\Delta t)\Delta t \qquad (5.3\text{-}11)$$

where

t = nΔt, n = 0,1,2,3,…

Δt = t interval, sampling interval

f(t) = a function of t

$F(s) = K_{\Delta t}[f(t)]$ , the $K_{\Delta t}$ Transform of the function f(t)

Expanding the summation of Eq 5.3-2 into a series

$$F(s) = \Delta t\,[f(0)(1+s\Delta t)^{-1} + f(1\Delta t)(1+s\Delta t)^{-2} + f(2\Delta t)(1+s\Delta t)^{-3} + f(3\Delta t)(1+s\Delta t)^{-4} + \ldots] \qquad (5.3\text{-}12)$$

$$F(s) = f(0)[\frac{\Delta t}{(1+s\Delta t)^1}] + f(\Delta t)[\frac{\Delta t}{(1+s\Delta t)^2}] + f(2\Delta t)[\frac{\Delta t}{(1+s\Delta t)^3}] + f(3\Delta t)[\frac{\Delta t}{(1+s\Delta t)^4}] + \ldots] \qquad (5.3\text{-}13)$$

where

$t = n\Delta t$, $n = 0,1,2,3,\ldots$

$\Delta t$ = t interval, sampling interval

$f(t)$ = function of t

$F(s) = K_{\Delta t}[f(t)]$ , the $K_{\Delta t}$ Transform of the function f(t)

From Eq 5.3-13, the series expansion of Eq 5.3-2, the values of $f(n\Delta t)$ where $n = 0,1,2,3,\ldots$ are found to be the coefficients of the terms, $[\frac{\Delta t}{(1+s\Delta t)^{n+1}}]$.

Control system analysis is possible using the Interval Calculus $K_{\Delta x}$ Transform. Though the Z Transform provides a convenient way to evaluate and analyze discrete calculus problems, there is still reason for the use of the $K_{\Delta x}$ Transform instead. In particular, the form of the $K_{\Delta x}$ Transform lends itself well to discrete Interval Calculus integration and the $K_{\Delta x}$ Transform is very similar to the Laplace Transform. In fact, the $K_{\Delta x}$ Transform may be considered to be an extension of the Laplace Transform to discrete mathematics. The similarity between the $K_{\Delta x}$ transform and the Laplace Transform can be helpful in the analysis and solution of mathematical problems of a discrete variable.

## Section 5.4: The relationship between the $K_{\Delta x}$ Transform and the Z Transform

In Section 1 of this chapter, the $K_{\Delta t}$ Transform was derived from the Z Transform. Therefore, the two transforms are directly related. It would be useful to have a convenient method for converting between these two transforms. Such a method has been found and is presented below.

The derived transform conversion equations are as follows:

1) K$\Delta$x Transform to Z Transform Conversion

$$Z[f(t)] = F(z) = \frac{z}{T} F(s)\Big|_{\substack{s = \frac{z-1}{T} \\ T = \Delta t \\ t = 0,\, T,\, 2T,\, 3T,\, \ldots}} \qquad \begin{matrix} Z[f(t)] = F(z) & \text{Z Transform} \\ K_{\Delta x}[f(t)] = F(s) & K_{\Delta t}\text{ Transform} \end{matrix} \qquad (5.4\text{-}1)$$

For inverse transform conversions, the complex plane integration contour changes from

$$s = \frac{e^{(\gamma+jw)\Delta t} - 1}{\Delta t} \ \text{ to } z = e^{(\gamma+jw)T} \qquad (5.4\text{-}2)$$

where

$\gamma$ = positive real constant

$-\frac{\pi}{\Delta t} \leq w < \frac{\pi}{\Delta t}$

2) Z Transform to $K_{\Delta x}$ Transform Conversion

$$K_{\Delta t}[f(t)] = F(s) = \frac{\Delta t}{1+s\Delta t} F(z)\Big|_{\substack{z = 1+s\Delta t \\ \Delta t = T \\ t = 0, \Delta t, 2\Delta t, 3\Delta t, \dots}} \qquad K_{\Delta t}[f(t)] = F(s) \quad K_{\Delta t} \text{ Transform} \quad (5.4\text{-}3)$$

$$Z[f(t)] = F(z) \quad \text{Z Transform}$$

For inverse transform conversions, the complex plane integration contour changes from

$$z = e^{(\gamma+jw)T} \text{ to } s = \frac{e^{(\gamma+jw)\Delta t} - 1}{\Delta t} \qquad (5.4\text{-}4)$$

where

$\gamma$ = positive real constant

$$-\frac{\pi}{\Delta t} \le w < \frac{\pi}{\Delta t}$$

The derivation of the above two conversion relationships follows.

The definitions of the Z and $K_{\Delta t}$ Transforms are rewritten below. Comparing both transforms, the following observations are made. The transforms are series which represent the sampling of a continuous function, f(t), at a sampling interval of $\Delta t = T$. Both transforms involve unit functions multiplied by a weighting function of $f(n\Delta t)$ where n = 0,1,2,3,… The unit function for the Z Transform is the unit impulse function which has a Laplace Transform of $e^{-nTs}$. The unit function for the the $K_{\Delta t}$ Transform is the unit amplitude pulse function which has a $K_{\Delta t}$ Transform of $(1+s\Delta t)^{-n-1}\Delta t$. From this similarity, the specified conversion relationships, Eq 5.4-1 and Eq 5.4-3, are derived:

The definition of the $K_{\Delta t}$ Transform is as follows:

From Eq 5.3-2

$$K_{\Delta t}[f(t)] = {}_{1}\sum_{n=0}^{\infty} f(n\Delta t)[1+s\Delta t]^{-n-1}\Delta t = {}_{\Delta t}\sum_{t=0}^{\infty} f(t)[1+s\Delta t]^{-(\frac{t+\Delta t}{\Delta t})}\Delta t = {}_{\Delta t}\int_{0}^{\infty} (1+s\Delta t)^{-(\frac{t+\Delta t}{\Delta t})} f(t)\Delta t \qquad (5.4\text{-}5)$$

where $t = n\Delta t$ , n = 0,1,2,3,…

As shown in the first summation, the $K_{\Delta t}$ Transform is the sum of a weighting function, $f(n\Delta t)$, times its corresponding unit amplitude pulse $K_{\Delta t}$ Transform, $(1+s\Delta t)^{-n-1}\Delta t$. The $K_{\Delta t}$ Transform of a unit area pulse is $(1+s\Delta t)^{-n-1}$.

The definition of the Z Transform is as follows:

$$Z[f(t)] = \sum_{n=0}^{\infty} f(nT)e^{-nTs} = \sum_{n=0}^{\infty} f(nT)z^{-n} \ , \quad z = e^{Ts} \qquad (5.4\text{-}6)$$

As shown in the first summation, the Z transform is the sum of a weighting function, f(nT), times the Laplace Transform of its corresponding unit impulse, $e^{-nTs}$.

The literature describing the Z Transform commonly uses the letter, T, to represent the sampling interval. The $K_{\Delta t}$ Transform uses the designation, $\Delta t$, to represent this same sampling interval.
Thus, $\Delta t = T$.

$$\Delta t = T \tag{5.4-7}$$

Eq 5.4-5 and Eq 5.4-6 are transforms of the same function, f(t)

From Eq 5.4-5

$$K_{\Delta t}[f(t)] = \Delta t \sum_{n=0}^{\infty} f(n\Delta t)(1+s\Delta t)^{-n-1} \tag{5.4-8}$$

From Eq 5.4-6 and Eq 5.4-7

$$Z[f(t)] = \sum_{n=0}^{\infty} f(n\Delta t)z^{-n} \tag{5.4-9}$$

Note that there is a similarity between Eq 5.4-8 and Eq 5.4-9.

$$K_{\Delta t}[f(t)] = \Delta t \sum_{n=0}^{\infty} f(n\Delta t)(1+s\Delta t)^{-n-1} = f(0)(1+s\Delta t)^{-1}\Delta t + f(\Delta t)(1+s\Delta t)^{-1}\Delta t + f(2\Delta t)(1+s\Delta t)^{-2}\Delta t + \ldots$$

$$Z[f(t)] = \sum_{n=0}^{\infty} f(n\Delta t)z^{-n} = f(0)z^{0} + f(\Delta t)z^{-1} + f(2\Delta t)z^{-2} + \ldots$$

From Eq 5.4-8 and Eq 5.4-9, find the relationship to convert the $K_{\Delta t}$ Transform of a function to the Z Transform of the same function.

From Eq 5.4-8

$$K_{\Delta t}[f(t)] = \frac{\Delta t}{1+s\Delta t} \sum_{n=0}^{\infty} f(n\Delta t)(1+s\Delta t)^{-n} \tag{5.4-10}$$

$$\sum_{n=0}^{\infty} f(n\Delta t)(1+s\Delta t)^{-n} = \frac{1+s\Delta t}{\Delta t} K_{\Delta t}[f(t)] \tag{5.4-11}$$

$$Z[f(t)] = F(z) , \quad \text{The Z Transform} \tag{5.4-12}$$

$$K_{\Delta t}[f(t)] = F(s) , \quad \text{The } K_{\Delta t} \text{ Transform} \tag{5.4-13}$$

Let

$$s = \frac{z-1}{\Delta t} \tag{5.4-14}$$

From Eq 5.4-7 and Eq 5.4-11 thru Eq 1.6-14

$$\sum_{n=0}^{\infty} f(n\Delta t)(1+s\Delta t)^{-n} \Big|_{s=\frac{z-1}{\Delta t}} = \sum_{n=0}^{\infty} f(nT)z^{-n} = F(z) = \frac{(1+s\Delta t)}{\Delta t} K_{\Delta t}[f(t)] \Big|_{s=\frac{z-1}{\Delta t}} \;,\; T = \Delta t \qquad (5.4\text{-}15)$$

Note that the relationship, $s = \frac{z-1}{\Delta t}$, has converted the $K_{\Delta t}$ Transform to the Z Transform. This was to be expected. This relationship is the linear conformal mapping transform that was used to derive the $K_{\Delta t}$ Transform from the Z Transform. See Eq 5.1-7.

$$\sum_{n=0}^{\infty} f(nT)z^{-n} = Z[f(t)] = F(z) = \frac{z}{T} F(s) \Big|_{s=\frac{z-1}{T}} \;,\; T = \Delta t \qquad (5.4\text{-}16)$$

Then, the relationship to convert the $K_{\Delta t}$ Transform of f(t) to the corresponding Z Transform of f(t) is as follows:

$$\mathbf{Z[f(t)] = F(z) = \frac{z}{T} F(s)\Big|_{\substack{s=\frac{z-1}{T} \\ T=\Delta t}}} \qquad (5.4\text{-}17)$$

**where**

**Z[f(t)] = F(z)** **The Z Transform**

**$K_{\Delta x}$[f(t)] = F(s)** **The $K_{\Delta x}$ Transform**

From Eq 5.4-7, Eq 5.4-10, and Eq 5.4-12 thru Eq 5.4-14, find the relationship to convert the Z Transform of a function to the $K_{\Delta t}$ Transform of the same function.

$$K_{\Delta t}[f(t)] = \frac{\Delta t}{1+s\Delta t} \sum_{n=0}^{\infty} f(n\Delta t)(1+s\Delta t)^{-n} \qquad (5.4\text{-}18)$$

$$z = 1+s\Delta t \qquad (5.4\text{-}19)$$

$$K_{\Delta t}[f(t)]\Big|_{z=1+s\Delta t} = \frac{\Delta t}{1+s\Delta t} \sum_{n=0}^{\infty} f(n\Delta t)(1+s\Delta t)^{-n} \Big|_{z=1+s\Delta t} = \frac{\Delta t}{1+s\Delta t} \sum_{n=0}^{\infty} f(n\Delta t)z^{-n} \Big|_{z=1+s\Delta t}$$

$$= \frac{\Delta t}{1+s\Delta t} F(z) \Big|_{z=1+s\Delta t} \qquad (5.4\text{-}20)$$

$$K_{\Delta t}[f(t)] = F(s) = \frac{\Delta t}{1+s\Delta t} F(z) \Big|_{z=1+s\Delta t} \;,\; T = \Delta t \qquad (5.4\text{-}21)$$

Then, the relationship to convert the Z Transform of f(t) to the corresponding $K_{\Delta x}$ Transform of f(t) is as follows:

$$\mathbf{K_{\Delta t}[f(t)] = F(s) = \frac{\Delta t}{1+s\Delta t} F(z)\Big|_{\substack{z\, =\, 1+s\Delta t \\ \Delta t\, =\, T}}} \qquad (5.4\text{-}22)$$

**where**

**Z[f(t)] = F(z)** **The Z Transform**

$\mathbf{K_{\Delta t}[f(t)] = F(s)}$ **The $K_{\Delta x}$ Transform**

Though the previous conversions provide the Z and $K_{\Delta t}$ Transforms for the same function, f(t), it should be remembered that the Transforms themselves do not represent the same thing. The $K_{\Delta t}$ Transform represents a series of unit amplitude pulses and the Z Transform represents a series of unit area impulses with the same weighting, f(0), f(Δt), f(2Δt), … from the same function, f(t). In control system analysis, the $K_{\Delta t}$ Transform would be a good match where sample and hold samplers are used and the Z Transform would be a good match where impulse samplers are used.

For converting between the Inverse $K_{\Delta t}$ Transform and the Z Transform, refer to Section 5.1. In this section, the $K_{\Delta t}$ Transform and Z Transform contours are specified. See Figure 5.1-1 and Figure 5.1-2. Comparing these two transform complex plane contours, the following contour changes are determined.

**For conversion from the Inverse $K_{\Delta t}$ Transform to the Inverse Z Transform**

**The complex plane integration contour changes from**

$$\mathbf{s = \frac{e^{(\gamma+jw)\Delta t} - 1}{\Delta t} \text{ to } z = e^{(\gamma+jw)T}} \qquad (5.4\text{-}23)$$

**where**

**γ = positive real constant**

$$\mathbf{-\frac{\pi}{\Delta t} \le w < \frac{\pi}{\Delta t}}$$

**For conversion from the Inverse Z Transform to the Inverse $K_{\Delta t}$ Transform**

**The complex plane integration contour changes from**

$$\mathbf{z = e^{(\gamma+jw)T} \text{ to } s = \frac{e^{(\gamma+jw)\Delta t} - 1}{\Delta t}} \qquad (5.4\text{-}24)$$

**where**

**γ = positive real constant**

$$\mathbf{-\frac{\pi}{\Delta t} \le w < \frac{\pi}{\Delta t}}$$

**Thus, from the previously derived equations, the $K_{\Delta t}$ and Z Transform conversion equations are as follows:**

**1) $K_{\Delta t}$ Transform to Z Transform Conversion**

$$Z[f(t)] = F(z) = \frac{z}{T} F(s)\Big|_{s = \frac{z-1}{T}} \qquad Z[f(t)] = F(z) \quad \text{Z Transform} \qquad (5.4\text{-}25)$$

$T = \Delta t$ sampling period

$t = nT$ $\qquad K_{\Delta x}[f(t)] = F(s) \quad K_{\Delta t}$ Transform

$n = 0,1,2,3,\ldots$

**For conversion from the Inverse $K_{\Delta t}$ Transform to the Inverse Z Transform**
**The complex plane integration contour changes from**

$$s = \frac{e^{(\gamma+jw)\Delta t} - 1}{\Delta t} \quad \text{to } z = e^{(\gamma+jw)T} \qquad (5.4\text{-}26)$$

**where**

$\gamma$ = positive real constant

$$-\frac{\pi}{\Delta t} \le w < \frac{\pi}{\Delta t}$$

**2) <u>Z Transform to $K_{\Delta t}$ Transform Conversion</u>**

$$K_{\Delta t}[f(t)] = F(s) = \frac{\Delta t}{1+s\Delta t} F(z)\Big|_{z = 1+s\Delta t} \qquad K_{\Delta x}[f(t)] = F(s) \quad K_{\Delta x} \text{ Transform} \qquad (5.4\text{-}27)$$

$\Delta t = T$ sampling period

$t = n\Delta t$ $\qquad Z[f(t) = F(z) \quad$ Z Transform

$n = 0, 1, 2, 3, \ldots$

**For conversion from the Inverse Z Transform to the Inverse $K_{\Delta t}$ Transform**
**The complex plane integration contour changes from**

$$z = e^{(\gamma+jw)T} \quad \text{to } s = \frac{e^{(\gamma+jw)\Delta t} - 1}{\Delta t} \qquad (5.4\text{-}28)$$

**where**

$\gamma$ = positive real constant

$$-\frac{\pi}{\Delta t} \le w < \frac{\pi}{\Delta t}$$

The previously derived Transform conversion equations, Eq 5.4-24 thru Eq 5.4-28, can be used to derive other important transform equations.

Derive three equations to obtain the Z transform of a function, f(t).

From Eq 5.4-5

$$K_{\Delta t}[f(t)] = \int_{0}^{\infty}{}_{\Delta t} (1+s\Delta t)^{-(\frac{t+\Delta t}{\Delta t})} f(t)\Delta t \qquad (5.4\text{-}29)$$

From Eq 5.4-25

$$Z[f(t)] = F(z) = \frac{z}{T}\, F(s)\Big|_{\substack{s = \frac{z-1}{T} \\ \Delta x = T}} = \frac{z}{T}\, K_{\Delta t}[f(t)]\Big|_{\substack{s = \frac{z-1}{T} \\ \Delta x = T}} \qquad (5.4\text{-}30)$$

Use Eq 5.4-30 to convert the $K_{\Delta t}$ Transform of f(t), Eq 5.4-29, to the Z Transform of f(t)

$$Z[f(t)] = \frac{z}{T}\, K_{\Delta t}[f(t)]\Big|_{\substack{s = \frac{z-1}{T} \\ \Delta x = T}} = \frac{z}{T} \int_{0}^{\infty}{}_{T} (1 + [\frac{z-1}{\Delta t}]\Delta t)^{-(\frac{t+\Delta t}{\Delta t})} f(t)\Delta t \qquad (5.4\text{-}31)$$

$$Z[f(t)] = \frac{z}{T} \int_{0}^{\infty}{}_{T} z^{-\frac{t}{\Delta t}-1} f(t)\Delta t \quad , \ \Delta t = T \qquad (5.4\text{-}32)$$

$$Z[f(t)] = \frac{1}{T} \int_{0}^{\infty}{}_{T} z^{-\frac{t}{\Delta t}} f(t)\Delta t \quad , \ \Delta t = T \qquad (5.4\text{-}33)$$

Then, the first equation derived to calculate the Z Transform is:

Equation #1

$$\mathbf{Z[f(t)] = F(z) = \frac{1}{T} \int_{0}^{\infty}{}_{T} z^{-\frac{t}{\Delta t}} f(t)\Delta t} \qquad (5.4\text{-}34)$$

**where**

**$\Delta t = T$**

**$\Delta t$ = sampling period, t interval**

$$t = n\Delta t \ , \quad n = 0,1,2,3,\ldots \qquad (5.4\text{-}35)$$

Substituting Eq 5.4-35 into Eq 5.4-34 and changing the integration variable from t to n

$$Z[f(t)] = \frac{\Delta t}{T} \int_{0}^{\infty}{}_{1} z^{-n} f(n\Delta t)\Delta n \qquad (5.4\text{-}36)$$

$\Delta t = T$

Then, the second equation derived to calculate the Z Transform is:

Equation #2

$$\mathbf{Z[f(t)] = F(z) = {}_1\!\int_0^\infty z^{-n}\, f(n\Delta t)\Delta n} \qquad (5.4\text{-}37)$$

**where**

**$\Delta t = T$**

**$t = n\Delta t$ ,   $n = 0,1,2,3,\ldots$**

**$\Delta t$ = sampling period,  t interval**

From Eq 5.4-34

$$Z[f(t)] = \frac{1}{T}\,{}_T\!\int_0^\infty z^{-\frac{t}{\Delta t}}\, f(t)\Delta t = \frac{1}{T}\,{}_T\!\int_0^\infty (1 + [\tfrac{z-1}{T}]\Delta t)^{-\frac{t}{\Delta t}}\, f(t)\Delta t = \frac{1}{T}\,{}_T\!\int_0^\infty (1+a\Delta t)^{-\frac{t}{\Delta t}}\, f(t)\Delta t \qquad (5.4\text{-}38)$$

where

$a = \frac{z-1}{T}$

$\Delta t = T$

$\Delta t$ = sampling period,  t interval

Then, the third equation derived to calculate the Z Transform is:

Equation #3

$$\mathbf{Z[f(t)] = F(z) = \frac{1}{T}\,{}_T\!\int_0^\infty (1+a\Delta t)^{-\frac{t}{\Delta t}}\, f(t)\Delta t = \frac{1}{T}\,{}_T\!\int_0^\infty e_{\Delta t}(a,-t)\, f(t)\Delta t} \qquad (5.4\text{-}39)$$

**where**

**$a = \frac{z-1}{T}$**

**$\Delta t = T$**

**$\Delta t$ = sampling period,  t interval**

**$e_{\Delta t}(a,-t) = (1+a\Delta t)^{-\frac{t}{\Delta t}}$**

Comment – To help solve Eq 5.4-39, the following important Interval Calculus equation is available:

$${}_{\Delta t}\!\int e_{\Delta t}(a,-t)\Delta x = -\frac{1+a\Delta t}{a} e_{\Delta x}(a,-t) + k$$

Eq 5.4-34, Eq 5.4-37, and Eq 5.4-39 are all valid equations that can be used to calculate the Z Transform of a function, f(t). Depending on the application, one or another of these equations may be more convenient to use.

To demonstrate the use of some of the previously derived Interval Calculus equations, the following two examples are provided.

Example 5.4-1  Find the Z Transform of $f(t) = e^{at}$ using the derived Interval Calculus equation, Eq 5.4-34.

Writing Eq 5.4-34

$$Z[f(t)] = F(z) = \frac{1}{T}\int_{0}^{\infty} z^{-\frac{t}{\Delta t}} f(t)\Delta t \quad (5.4\text{-}40)$$

where

$\Delta t = T$

$\Delta t$ = sampling period, t interval

$$f(t) = e^{at} \quad (5.4\text{-}41)$$

Substituting Eq 5.4-41 into Eq 5.4-40

$$Z[e^{at}] = F(z) = \frac{1}{T}\int_{0}^{\infty} z^{-\frac{t}{\Delta t}} e^{at}\Delta t = \frac{1}{T}\int_{0}^{\infty} z^{-\frac{t}{\Delta t}} [e^{a\Delta t}]^{\frac{t}{\Delta t}}\Delta t = \frac{1}{T}\int_{0}^{\infty} [\frac{z}{e^{a\Delta t}}]^{-\frac{t}{\Delta t}}\Delta t$$

$$Z[e^{at}] = F(z) = \frac{1}{T}\int_{0}^{\infty} [\frac{z}{e^{a\Delta t}}]^{-\frac{t}{\Delta t}}\Delta t \quad (5.4\text{-}42)$$

Let

$$1+b\Delta t = \frac{z}{e^{a\Delta t}} \quad (5.4\text{-}43)$$

Solving for b

$$b = \frac{ze^{-a\Delta t} - 1}{\Delta t} \quad (5.4\text{-}44)$$

Substituting Eq 5.4-43 into Eq 5.4-42

$$Z[e^{at}] = F(z) = \frac{1}{T}\int_{0}^{\infty} [1+b\Delta t]^{-\frac{t}{\Delta t}}\Delta t \quad (5.4\text{-}45)$$

From the discrete integration table, Table 6, in the Appendix

$$\int_{\Delta t} (1+b\Delta t)^{-\frac{t}{\Delta t}}\,\Delta t = -\frac{1+b\Delta t}{b}(1+b\Delta t)^{-\frac{t}{\Delta t}} + k \tag{5.4-46}$$

From Eq 5.4-44 thru Eq 5.4-46

$$Z[e^{at}] = \frac{1}{T}\int_{T\ 0}^{\infty}[1+(\frac{ze^{-a\Delta t}-1}{\Delta t})\Delta t]^{-\frac{t}{\Delta t}}\,\Delta t = -\frac{1}{T}\frac{1+\frac{ze^{-a\Delta t}-1}{\Delta t}\Delta t}{\frac{ze^{-a\Delta t}-1}{\Delta t}}(1+\frac{ze^{-a\Delta t}-1}{\Delta t}\Delta t)^{-\frac{t}{\Delta t}}\Big|_0^{\infty} \tag{5.4-47}$$

Simplifying

$\Delta t = T$

$$Z[e^{at}] = -\frac{ze^{-a\Delta t}}{ze^{-a\Delta t}-1}(ze^{-a\Delta t})^{-\frac{t}{\Delta t}}\Big|_0^{\infty} = \frac{z}{z-e^{a\Delta t}} = \frac{z}{z-e^{aT}} \tag{5.4-48}$$

Then

$$\mathbf{Z[e^{at}] = \frac{z}{z-e^{aT}}} \tag{5.4-49}$$

<u>Example 5.4-2</u> Find the Z Transform of $f(t) = e^{at}$ from the $K_{\Delta t}$ Transform of $f(t) = e^{at}$ using the derived Interval Calculus equation, Eq 5.4-25.

Writing Eq 5.4-25, the $K_{\Delta t}$ Transform to Z Transform Conversion equation

$$Z[f(t)] = F(z) = \frac{z}{T}f(s)\Big|_{s=\frac{z-1}{T}} \tag{5.4-50}$$

$\Delta t = T$ sampling period

$t = n\Delta t$

$n = 0,1,2,3,\ldots$

$Z[f(t)] = F(z)$ Z Transform

$K_{\Delta x}[f(t)] = F(s)$ $K_{\Delta t}$ Transform

$$f(t) = e^{at} \tag{5.4-51}$$

From Table 3, a $K_{\Delta t}$ Transform table, in the Appendix

$$K_{\Delta t}[e^{at}] = \frac{1}{s-(\frac{e^{a\Delta t}-1}{\Delta t})} \tag{5.4-52}$$

From Eq 5.4-50 thru Eq 5.4- 52

$$Z[e^{at}] = \frac{z}{T}K_{\Delta t}[e^{at}]\Big|_{s=\frac{z-1}{T}} \tag{5.4-53}$$

$\Delta t = T$ sampling period

$$Z[e^{at}] = \frac{z}{T}\frac{1}{s-\left(\frac{e^{a\Delta t}-1}{\Delta t}\right)}\Big|_{s=\frac{z-1}{T}} \quad \Delta t = T \text{ sampling period} \tag{5.4-54}$$

$$Z[e^{at}] = \frac{z}{T}\frac{1}{\frac{z-1}{T}-\left(\frac{e^{aT}-1}{T}\right)} = \frac{z}{z-e^{aT}} \tag{5.4-55}$$

$$\mathbf{Z[e^{at}] = \frac{z}{z-e^{aT}}} \tag{5.4-56}$$

Three additional examples of the use of $K_{\Delta t}$ and Z Transform calculation and conversion equations are shown in Example 1.3, Example 1.4, and Example 1.5 of the Chapter 1 Solved Problems Section located at the end of Chapter 1.

Some final comments concerning the Z Transform, and the $K_{\Delta t}$ Transform

1. The $K_{\Delta t}$ Transform is derived from the Z Transform and there exist simple relationships which convert the results of one transform to the other.

2. The $K_{\Delta t}$ Transform is applied to discrete Interval Calculus functions such as $e_{\Delta t}(a,x)$, $\sin_{\Delta t}(a,x)$, etc.

3. Interval Calculus discrete functions such as $e_{\Delta t}(a,t) = (1+a\Delta t)^{\frac{t}{\Delta t}}$ are actually a related set of functions which differ as a function of $\Delta t$. The $K_{\Delta t}$ Transform of each function of the set is the same, $K_{\Delta t}[e_{\Delta t}(a,t)] = \frac{1}{s-a}$. The set of functions, $e_{\Delta t}(a,t)$, is the solution to the set of related differential difference equations, $D_{\Delta t}y_{\Delta t}(t) + ay_{\Delta t}(t) = \frac{y_{\Delta t}(t+\Delta t) - y_{\Delta t}(t)}{\Delta t} + ay_{\Delta t}(t) = 0$.

4. The $K_{\Delta t}$ Transform is primarily applied to differential difference equations but, like the Z Transform, may also be applied to difference equations.

5. The Inverse Z Transform and the Inverse $K_{\Delta t}$ Transform are discrete functions where the independent variable, t, is defined only at specific equally spaced intervals, $t = 0, \Delta t, 2\Delta t, 3\Delta t, \ldots$

6. Both the Z Transform and the $K_{\Delta t}$ Transform can be expanded into series, the coefficients of which are the values of f(t) at $t = 0, \Delta t, 2\Delta t, 3\Delta t, \ldots$

7. Though the $K_{\Delta t}$ Transform is derived from the Z Transform, it is more similar to the Laplace Transform. In fact, the $K_{\Delta t}$ Transform becomes the Laplace Transform when the value of $\Delta t$ is infinitesimal.

## Section 5.5: Application of the $K_{\Delta t}$ to Z and Z to $K_{\Delta t}$ Transform relationships to discrete variable control system analysis

In the previous sections of this chapter the following relationships were derived:

1) $K_{\Delta x}$ Transform to Z Transform Conversion

$$Z[f(t)] = F(z) = \frac{z}{T} F(s)\Big|_{s = \frac{z-1}{T}} \qquad Z[f(t)] = F(z) \quad \text{Z Transform} \tag{5.5-1}$$

$T = \Delta t$

$t = 0, T, 2T, 3T, \ldots \qquad K_{\Delta x}[f(t)] = F(s) \quad K_{\Delta t}$ Transform

For inverse transform conversions, the complex plane integration contour changes from

$$s = \frac{e^{(\gamma+jw)\Delta t} - 1}{\Delta t} \quad \text{to } z = e^{(\gamma+jw)T} \tag{5.5-2}$$

where

$\gamma$ = positive real constant

$-\frac{\pi}{\Delta t} \leq w < \frac{\pi}{\Delta t}$

2) <u>Z Transform to $K_{\Delta x}$ Transform Conversion</u>

$$K_{\Delta t}[f(t)] = F(s) = \frac{\Delta t}{1+s\Delta t} F(z)\Big|_{z = 1+s\Delta t} \qquad K_{\Delta t}[f(t)] = F(s) \quad K_{\Delta t} \text{ Transform} \tag{5.5-3}$$

$\Delta t = T$

$t = 0, \Delta t, 2\Delta t, 3\Delta t, \ldots \qquad Z[f(t)] = F(z) \quad$ Z Transform

For inverse transform conversions, the complex plane integration contour changes from

$$z = e^{(\gamma+jw)T} \text{ to } s = \frac{e^{(\gamma+jw)\Delta t} - 1}{\Delta t} \tag{5.5-4}$$

where

$\gamma$ = positive real constant

$-\frac{\pi}{\Delta t} \leq w < \frac{\pi}{\Delta t}$

These above relationships between the $K_{\Delta t}$ Transform and the Z Transform can be used to redefine, using $K_{\Delta t}$ Transforms, a control system defined by Z Transforms and vice versa.

<u>Z Transform to $K_{\Delta t}$ Transform Conversion</u>

Consider the following general feedback control system block diagram, Block Diagram 5.5-1.

<u>Block Diagram 5.5-1</u>

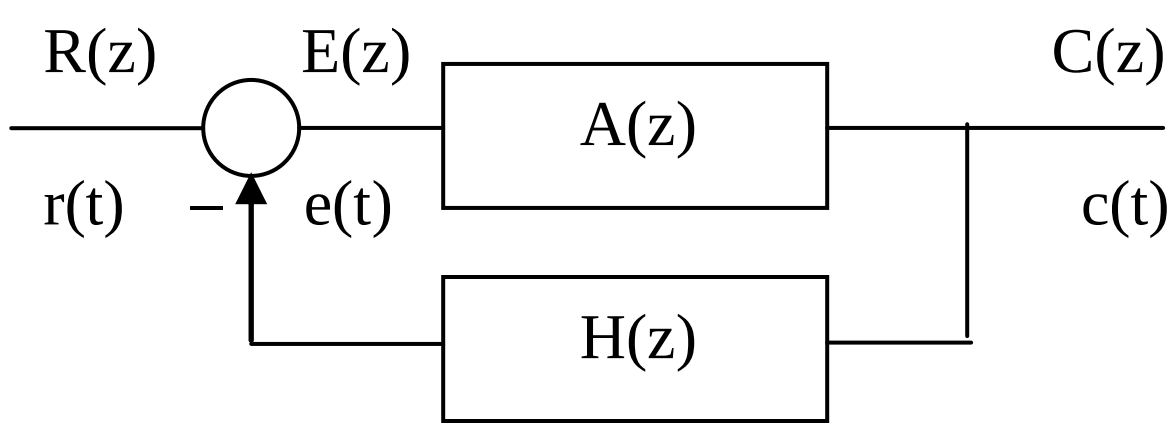

$C(z) = Z[c(t)]$, The Z Transform of $c(t)$
$R(z) = Z[r(t)]$, The Z Transform of $r(t)$
$E(z) = Z[e(t)]$, The Z Transform of $e(t)$
$t = 0, T, 2T, 3T, \ldots$

Combining transfer functions A(z) and H(z) into a single transfer function, G(z)

$$G(z) = \frac{A(z)}{1+A(z)H(z)}$$

With respect to the input, R(z), and the output, C(z), Block Diagram 5.5-2 below is the equivalent of Block Diagram 5.5-1.

Block Diagram 5.5-2

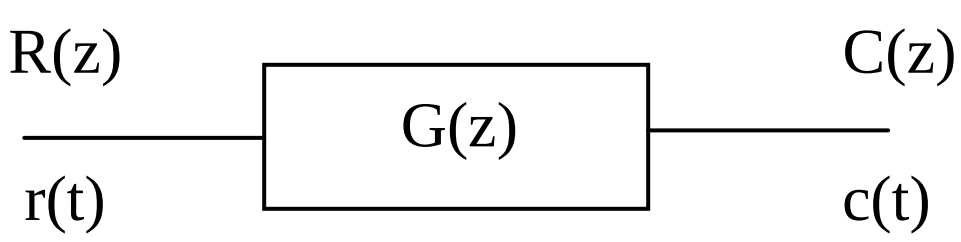

C(z) = Z[c(t)] , The Z Transform of c(t)
R(z) = Z[r(t)] , The Z Transform of r(t)
G(z) = Z Transform transfer function
t = 0, T, 2T, 3T, …

$$C(z) = G(z)R(z)$$ , Equation representing Block Diagram 5.5-2 (5.5-5)

Use Eq 5.5-3 , the Z Transform to $K_{\Delta t}$ Transform Convervion Relationship, to convert Eq 5.5-5 to an equivalent $K_{\Delta t}$ Transform equation.

Rewriting Eq 5.5-3

$$K_{\Delta t}[f(t)] = F(s) = \frac{\Delta t}{1+s\Delta t} F(z)\Big|_{\substack{z = 1+s\Delta t \\ \Delta t = T}}$$

$K_{\Delta t}[f(t)] = F(s)$ $K_{\Delta t}$ Transform

t = 0, Δt, 2Δt, 3Δt, … $Z[f(t)] = F(z)$ Z Transform

From Eq 5.5-3 and Eq 5.5-5

$$\frac{\Delta t}{1+s\Delta t} C(z)\Big|_{\substack{z = 1+s\Delta t \\ \Delta t = T}} = \frac{\Delta t}{1+s\Delta t} [G(z)R(z)]\Big|_{\substack{z = 1+s\Delta t \\ \Delta t = T}} \tag{5.5-6}$$

$$\frac{\Delta t}{1+s\Delta t} C(z)\Big|_{\substack{z = 1+s\Delta t \\ \Delta t = T}} = [G(z)\Big|_{\substack{z = 1+s\Delta t \\ \Delta t = T}}][\frac{\Delta t}{1+s\Delta t} [R(z)]\Big|_{\substack{z = 1+s\Delta t \\ \Delta t = T}}] \tag{5.5-7}$$

$$C(s) = G(z)\Big|_{\substack{z = 1+s\Delta t \\ \Delta t = T}} R(s) \tag{5.5-8}$$

$$C(s) = \mathit{G}(s)R(s) \tag{5.5-9}$$

where

$\mathit{G}(s) = G(z)\Big|_{\substack{z = 1+s\Delta t \\ \Delta t = T}}$ , $K_{\Delta t}$ Transform transfer function

$C(s) = K_{\Delta t}[c(t)]$ , $K_{\Delta t}$ Transform of c(t)

$R(s) = K_{\Delta t}[r(t)]$ , $K_{\Delta t}$ Transform of r(t)

Block Diagram 5.5-3 below represents Eq 5.5-9 and is the $K_{\Delta t}$ Transform equivalent of the Z Transform Block Diagram 5.5-2.

Block Diagram 5.5-3

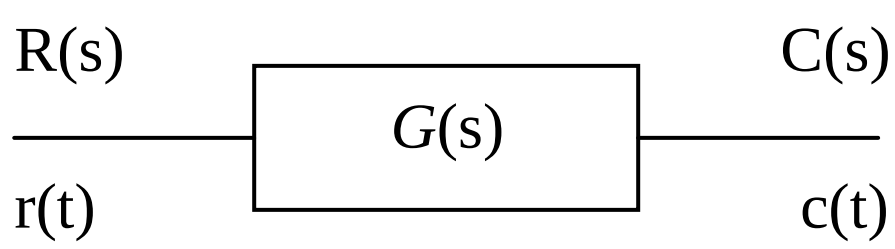

$C(s) = K_{\Delta t}[c(t)]$ , The $K_{\Delta t}$ Transform of c(t)
$R(s) = K_{\Delta t}[r(t)]$ , The $K_{\Delta t}$ Transform of r(t)
$G(s) = G(z)\Big|_{\substack{z = 1+s\Delta t \\ \Delta t = T \\ t = 0,\ \Delta t,\ 2\Delta t,\ 3\Delta t,\ \ldots}}$ , $K_{\Delta t}$ Transform transfer function

c(t) can be calculated using either Z Transforms or $K_{\Delta t}$ Transforms using Eq 5.5-5 or Eq 5.5-9

Then

**The Z Transform to $K_{\Delta t}$ Transform Transfer Function Conversion is:**

$$\frac{C(z)}{R(z)} = G(z) \text{ to } \frac{C(s)}{R(s)} = G(z)\Big|_{\substack{z = 1+s\Delta t \\ \Delta t = T \\ t = 0,\ \Delta x,\ 2\Delta x,\ 3\Delta x,\ \ldots}} \quad \textbf{Equivalent Z and } K_{\Delta t} \textbf{ Transform transfer functions} \quad (5.5\text{-}10)$$

**where**

**$C(z) = Z[c(t)]$ , The Z Transform of c(t)**
**$R(z) = Z[r(t)]$ , The Z Transform of r(t)**
**G(z) = Z Transform transfer function**
**$C(s) = K_{\Delta t}[c(t)]$ , The $K_{\Delta t}$ Transform of c(t)**
**$R(s) = K_{\Delta t}[r(t)]$ , The $K_{\Delta t}$ Transform of r(t)**
**$G(z)\Big|_{\substack{z = 1+s\Delta t \\ \Delta t = T \\ t = 0,\ \Delta t,\ 2\Delta t,\ 3\Delta t,\ \ldots}} = G(s)$, $K_{\Delta t}$ Transform transfer function**

$K_{\Delta t}$ Transform to Z Transform Conversion

Consider the following general feedback control system block diagram, Block Diagram 5.5-4.

Block Diagram 5.5-4

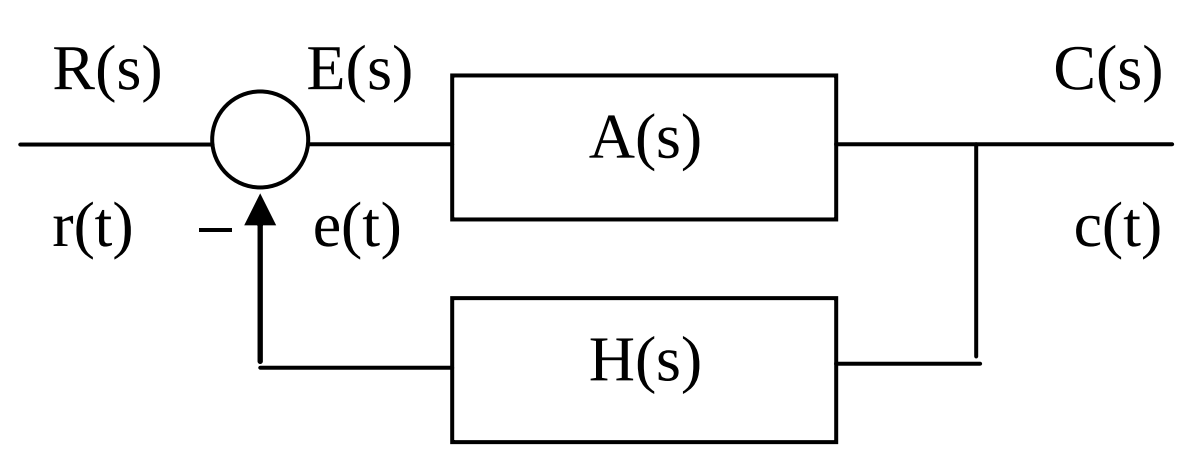

$C(s) = K_{\Delta t}[c(t)]$ , The $K_{\Delta t}$ Transform of c(t)
$R(s) = K_{\Delta t}[r(t)]$ , The $K_{\Delta t}$ Transform of r(t)
$E(s) = K_{\Delta t}[e(t)]$ , The $K_{\Delta t}$ Transform of e(t)
$t = 0, \Delta t, 2\Delta t, 3\Delta t, \ldots$

Combining transfer functions A(s) and H(s) into a single transfer function, G(s)

$$G(s) = \frac{A(s)}{1+A(s)H(s)}$$

With respect to the input, R(s), and the output, C(s), Block Diagram 5.5-5 below is the equivalent of Block Diagram 5.5-4.

<u>Block Diagram 5.5-5</u>

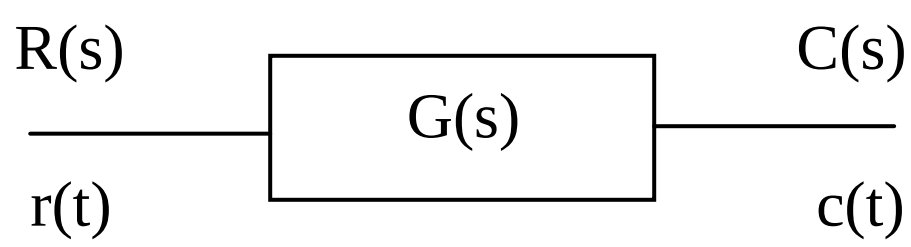

$C(s) = K_{\Delta t}[c(t)]$ , The $K_{\Delta t}$ Transform of c(t)
$R(s) = K_{\Delta t}[r(t)]$ , The $K_{\Delta t}$ Transform of r(t)
$G(s) = K_{\Delta t}$ Transform transfer function
$t = 0, \Delta t, 2\Delta t, 3\Delta t, \ldots$

$C(s) = G(s)R(s)$ , Equation representing Block Diagram 5.5-5 (5.5-11)

Use Eq 5.5-1, the $K_{\Delta t}$ Transform to Z Transform Convervion Relationship, to convert Eq 5.5-11 to an equivalent Z Transform equation.

Rewriting Eq 5.5-1

$$Z[f(t)] = F(z) = \frac{z}{T} F(s)\Big|_{\substack{s = \frac{z-1}{T} \\ T = \Delta t \\ t = 0, T, 2T, 3T, \ldots}}$$

$Z[f(t)] = F(z)$ Z Transform

$K_{\Delta x}[f(t)] = F(s)$ $K_{\Delta t}$ Transform

From Eq 5.5-11 and Eq 5.5-1

$$\frac{z}{T} C(s)\Big|_{\substack{s = \frac{z-1}{T} \\ T = \Delta t}} = \frac{z}{T}[G(s)R(s)]\Big|_{\substack{s = \frac{z-1}{T} \\ T = \Delta t}} \tag{5.5-12}$$

$$\frac{z}{T} C(s)\Big|_{\substack{s = \frac{z-1}{T} \\ T = \Delta t}} = [G(s)\Big|_{\substack{s = \frac{z-1}{T} \\ T = \Delta t}}][\frac{z}{T}[R(s)]\Big|_{\substack{s = \frac{z-1}{T} \\ T = \Delta t}}] \tag{5.5-13}$$

$$C(z) = G(s)\Big|_{\substack{s = \frac{z-1}{T} \\ T = \Delta t}} R(z) \tag{5.5-14}$$

$$C(z) = \mathcal{G}(z)R(z) \tag{5.5-15}$$

where

$G(z) = G(s)\Big|_{s = \frac{z-1}{T}}$ ,   Z Transform transfer function

$T = \Delta t$

$t = 0, T, 2T, 3T, \ldots$

$G(s) =$ $K_{\Delta t}$ Transform transfer function

$C(z) = Z[c(t)]$ ,   Z Transform of $c(t)$

$R(z) = Z[r(t)]$ ,   Z Transform of $r(t)$

$C(s) = K_{\Delta t}[c(t)]$ , $K_{\Delta t}$ Transform of $c(t)$

$R(s) = K_{\Delta t}[r(t)]$ , $K_{\Delta t}$ Transform of $r(t)$

Block Diagram 5.5-6 below represents Eq 5.5-15 and is the $K_{\Delta t}$ Transform equivalent of the Z Transform Block Diagram 5.5-5.

Block Diagram 5.5-6

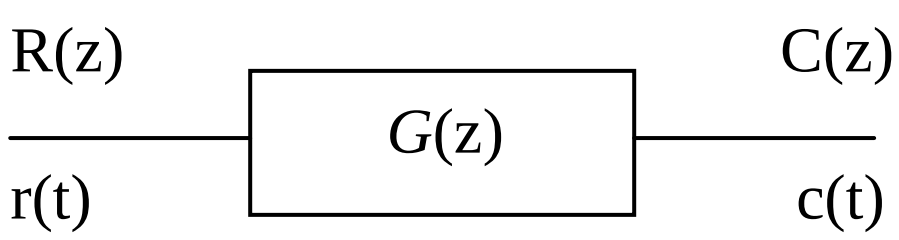

$C(z) = Z[c(t)]$ ,   The Z Transform of $c(t)$

$R(z) = Z[r(t)]$ ,   The Z Transform of $r(t)$

$G(z) = G(s)\Big|_{s = \frac{z-1}{T}}$ ,   Z Transform transfer function

$T = \Delta t$

$t = 0, T, 2T, 3T, \ldots$

c(t) can be calculated using either $K_{\Delta t}$ Transforms or Z Transforms using Eq 5.5-11 or Eq 5.5-15

Then

**The $K_{\Delta t}$ Transform to Z Transform Transfer Function Conversion is:**

$$\frac{C(s)}{R(s)} = G(s) \;\text{ to }\; \frac{C(z)}{R(z)} = G(s)\Big|_{s = \frac{z-1}{T}} \quad (5.5\text{-}16)$$

**, Equivalent Z and $K_{\Delta t}$ Transform transfer functions** (5.5-16)

**$T = \Delta t$**

**$t = 0, T, 2T, 3T, \ldots$**

**where**

**$C(s) = K_{\Delta t}[c(t)]$ , The $K_{\Delta t}$ Transform of $c(t)$**

**$R(s) = K_{\Delta t}[r(t)]$ , The $K_{\Delta t}$ Transform of $r(t)$**

**$G(s) = K_{\Delta t}$ Transform transfer function**

**$C(z) = Z[c(t)]$ , The Z Transform of $c(t)$**

**$R(z) = Z[r(t)]$ , The Z Transform of $r(t)$**

**$G(s)\Big|_{s = \frac{z-1}{T}} = G(z)$ , Z Transform transfer function**

**$T = \Delta t$**

**$t = 0, T, 2T, 3T, \ldots$**

**Discrete variable control systems can be analyzed using either Z Transforms or $K_{\Delta t}$ Transforms. If desired, one transform can be converted to the other using the following conversion equations:**

**The Z Transform to $K_{\Delta t}$ Transform Transfer Function Conversion is:**

$$\frac{C(z)}{R(z)} = G(z) \text{ to } \frac{C(s)}{R(s)} = G(z)\Big|_{\substack{z = 1+s\Delta t \\ \Delta t = T \\ t = 0, \Delta x, 2\Delta x, 3\Delta x, \ldots}} \quad \text{Equivalent Z and } K_{\Delta t} \text{ Transform transfer functions} \quad (5.5\text{-}17)$$

**The $K_{\Delta t}$ Transform to Z Transform Transfer Function Conversion is:**

$$\frac{C(s)}{R(s)} = G(s) \text{ to } \frac{C(z)}{R(z)} = G(s)\Big|_{\substack{s = \frac{z-1}{T} \\ T = \Delta t \\ t = 0, T, 2T, 3T, \ldots}} , \quad \text{Equivalent Z and } K_{\Delta t} \text{ Transform transfer functions} \quad (5.5\text{-}18)$$

**where**

**$C(s) = K_{\Delta t}[c(t)]$ , The $K_{\Delta t}$ Transform of c(t)**

**$R(s) = K_{\Delta t}[r(t)]$ , The $K_{\Delta t}$ Transform of r(t)**

**$G(z)\Big|_{\substack{z = 1+s\Delta t \\ \Delta t = T}} = G(s)$, $K_{\Delta t}$ Transform transfer function**

**$C(z) = Z[c(t)]$ , The Z Transform of c(t)**

**$R(z) = Z[r(t)]$ , The Z Transform of r(t)**

**$G(s)\Big|_{\substack{s = \frac{z-1}{T} \\ T = \Delta t}} = G(z)$ , Z Transform transfer function**

On the following page, Example 5.5-1 is presented which demonstrates a feedback control system analyzed by both the Z and $K_{\Delta t}$ Transforms.

Example 5.5-1 Solve the following Z Transform defined closed loop system for c(t) using $K_{\Delta t}$ Transforms. As a check, solve again for c(t) using Z Transforms.

Closed Loop Control System Z Transform Diagram

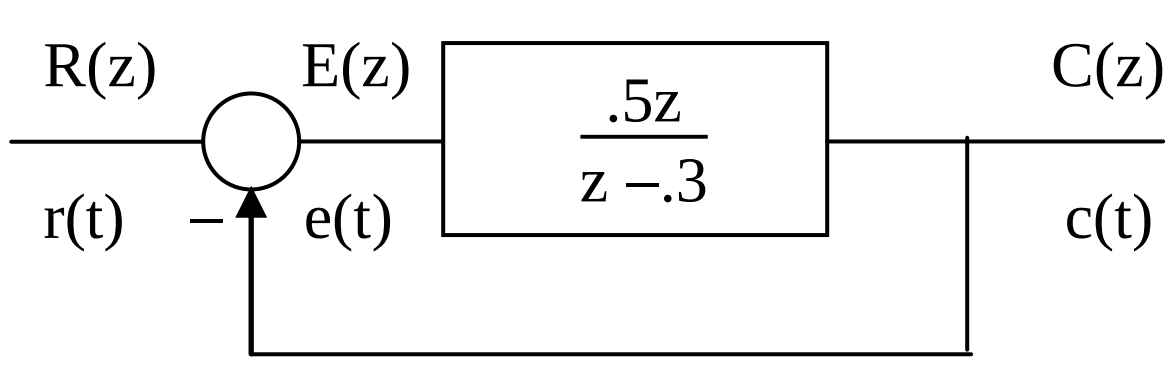

t = 0, T, 2T, 3T, …

r(t) = U(t) , unit step at t = 0

$R(z) = \frac{z}{z-1}$ , Z Transform of U(t)

T = .2

Initial conditions = 0

1) Find c(t) using $K_{\Delta t}$ Transforms

Convert the above closed loop control system Z Transform diagram to an equivalent $K_{\Delta t}$ Transform control system diagram.

$$\frac{C(z)}{E(z)} = \frac{.5z}{z - .3} \qquad (5.5\text{-}19)$$

Use the following Z Transform to $K_{\Delta t}$ Transform Transfer Function Conversion equations to convert $\frac{C(z)}{E(z)}$ to its equivalent $\frac{C(s)}{E(s)}$

$$\frac{O(z)}{I(z)} = F(z)\,, \quad \text{Z Transform equation} \qquad (5.5\text{-}20)$$

$$\frac{O(s)}{I(s)} = F(z)\Big|_{\substack{z = 1+s\Delta t \\ \Delta t = T \\ t = 0,\ \Delta x,\ 2\Delta x,\ 3\Delta x,\ \ldots}}\,, \quad \text{Equivalent } K_{\Delta t} \text{ Transform equation} \qquad (5.5\text{-}21)$$

From Eq 5.5-19 thru Eq 5.5-21
Calculate the equivalent $K_{\Delta t}$ Transform equation for Eq 5.5-19

$$\frac{C(s)}{E(s)} = \frac{.5z}{z - .3}\Big|_{\substack{z = 1+s\Delta t \\ \Delta t = T \\ t = 0,\ \Delta x,\ 2\Delta x,\ 3\Delta x,\ \ldots}}\,, \qquad (5.5\text{-}22)$$

$$\Delta t = T = .2 \qquad (5.5\text{-}23)$$

Simplifying Eq 5.5-22

$$\frac{C(s)}{E(s)} = \frac{.5(1+.2s)}{(1+.2s) - .3} = \frac{(.1s+.5)}{.2s + .7} = \frac{.1(s+5)}{.2(s+3.5)}$$

$$\frac{C(s)}{E(s)} = \frac{.5(s+5)}{(s+3.5)}\,, \quad \text{The equivalent } K_{\Delta t} \text{ Transform equation for Eq 5.5-19} \qquad (5.5\text{-}24)$$

Then the equivalent $K_{\Delta t}$ Transform system diagram is as follows

Closed Loop Control System $K_{\Delta t}$ Transform Diagram

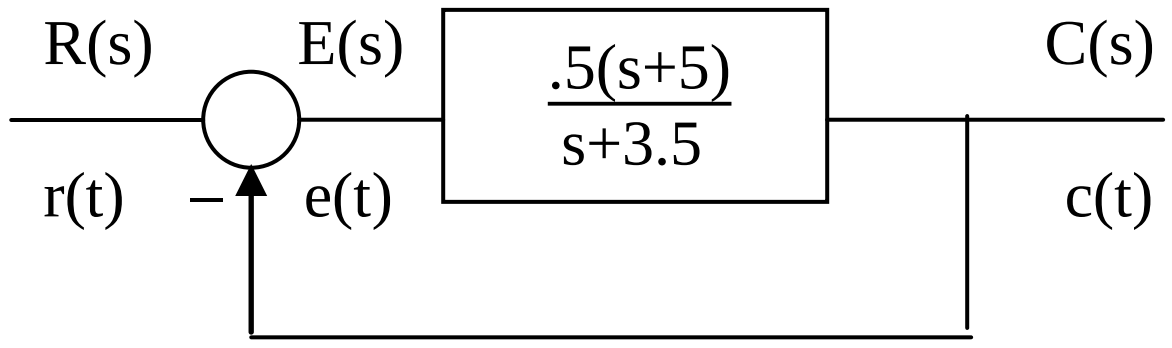

$t = 0, \Delta t, 2\Delta t, 3\Delta t, \ldots$

$r(t) = U(t)$ , unit step at $t = 0$

$R(s) = \frac{1}{s}$ , $K_{\Delta t}$ Transform of $U(t)$

$\Delta t = .2$

Find C(s)

$$C(s) = \frac{\frac{.5(s+5)}{s+3.5}}{1 + \frac{.5(s+5)}{s+3.5}} R(s) = \frac{\frac{.5(s+5)}{s+3.5}}{1 + \frac{.5(s+5)}{s+3.5}} \frac{1}{s} \qquad (5.5\text{-}25)$$

$$C(s) = \frac{.5(s+5)}{1.5s+6} \frac{1}{s} = \frac{1}{3} \frac{(s+5)}{s+4} \frac{1}{s}$$

$$C(s) = \frac{1}{3} \frac{(s+5)}{s(s+4)} \qquad (5.5\text{-}26)$$

Expand Eq 5.5-26 in a partial fraction expansion

$$C(s) = \frac{1}{3} \frac{(s+5)}{s(s+4)} = \frac{A}{s} + \frac{B}{s+4} \qquad (5.5\text{-}27)$$

$$A = \frac{1}{3} \frac{(s+5)}{(s+4)} \Big|_{s=0} = \frac{5}{3(4)} = \frac{5}{12}$$

$$A = \frac{5}{12} \qquad (5.5\text{-}28)$$

$$B = \frac{1}{3} \frac{(s+5)}{s} \Big|_{s=-4} = \frac{-4+5}{3(-4)} = -\frac{1}{12}$$

$$B = -\frac{1}{12} \qquad (5.5\text{-}29)$$

Substituting Eq 5.5-28 and Eq 5.5-29 into Eq 5.5-27

$$C(s) = \frac{5}{12} \frac{1}{s} - \frac{1}{12} \frac{1}{s+4} \qquad (5.5\text{-}30)$$

Taking the Inverse $K_{\Delta t}$ Transform of Eq 5.5-30

$$K_{\Delta t}[1] = \frac{1}{s} \qquad (5.5\text{-}31)$$

$$K_{\Delta t}[e_{\Delta t}(a,t)] = \frac{1}{s-a} \qquad (5.5\text{-}32)$$

$$\Delta t = .2 \qquad (5.5\text{-}33)$$

From Eq 5.5-30 thru Eq 5.5-33

$$c(t) = \frac{5}{12} - \frac{1}{12} e_{.2}(-4,t) \qquad (5.5\text{-}34)$$

$$e_{\Delta t}(a,t) = (1+a\Delta t)^{\frac{t}{\Delta t}} \qquad (5.5\text{-}35)$$

From Eq 5.5-33 thru 5.5-35

$$c(t) = \frac{5}{12} - \frac{1}{12} [1-4(.2)]^{\frac{t}{.2}} = \frac{5}{12} - \frac{1}{12} [.2]^{5t} \qquad (5.5\text{-}36)$$

$$c(t) = \frac{5}{12} - \frac{1}{12} [\frac{1}{5}]^{5t} \qquad (5.5\text{-}37)$$

Check the above result using a Z Transform solution of the control system defined by the Closed Loop Control System Z Transform Diagram.

2) Find c(t) using Z Transforms

Find C(z)

$$C(z) = \frac{\frac{.5z}{z-.3}}{1+\frac{.5z}{z-.3}} R(z) = \frac{\frac{.5z}{z-.3}}{1+\frac{.5z}{z-.3}} \frac{z}{z-1} \qquad (5.5\text{-}38)$$

$$C(z) = \frac{.5z}{z-.3+.5z} \frac{z}{z-1} = \frac{.5z}{1.5z-.3} \frac{z}{z-1}$$

$$C(z) = \frac{1}{3} \frac{z}{z-\frac{1}{5}} \frac{z}{z-1} = \frac{z}{3} \frac{z}{(z-\frac{1}{5})(z-1)}$$

$$C(z) = \frac{z}{3} \frac{z}{(z-1)(z-\frac{1}{5})} \qquad (5.5\text{-}39)$$

Expand Eq 5.5-39 in a partial fraction expansion

$$C(z) = \frac{z}{3} \frac{z}{(z - \frac{1}{5})(z\text{-}1)} = \frac{z}{3} [\frac{A}{z - 1} + \frac{B}{z - \frac{1}{5}} ] \qquad (5.5\text{-}40)$$

$$A = \frac{z}{z - \frac{1}{5}} |_{z=1} = \frac{1}{\frac{4}{5}} = \frac{5}{4}$$

$$A = \frac{5}{4} \qquad (5.5\text{-}41)$$

$$B = \frac{z}{z - 1} |_{s = \frac{1}{5}} = \frac{\frac{1}{5}}{-\frac{4}{5}} = -\frac{1}{4}$$

$$B = -\frac{1}{12} \qquad (5.5\text{-}42)$$

Substituting Eq 5.5-41 and Eq 5.5-42 into Eq 5.5-40

$$C(z) = \frac{5}{4(3)} \frac{z}{z - 1} - \frac{1}{4(3)} \frac{z}{z - \frac{1}{5}} \qquad (5.5\text{-}43)$$

$$C(z) = \frac{5}{12} \frac{z}{z - 1} - \frac{1}{12} \frac{z}{z - \frac{1}{5}} \qquad (5.5\text{-}44)$$

Taking the Inverse Z Transform of Eq 5.5-44

$$Z[1] = \frac{z}{z - 1} \qquad (5.5\text{-}45)$$

$$Z[a^{wt}] = \frac{z}{z - a^{wT}} \qquad (5.5\text{-}46)$$

$$T = .2 \qquad (5.5\text{-}47)$$

$$a^{wt} = [a^{wT}]^{\frac{t}{T}} = [\frac{1}{5}]^{\frac{t}{.2}} = [\frac{1}{5}]^{5t} \qquad (5.5\text{-}48)$$

From 5.5-44 thru 5.5-48

$$c(t) = \frac{5}{12} - \frac{1}{12}\left[\frac{1}{5}\right]^{5t} \qquad (5.5\text{-}46)$$

This is the same result as Eq 5.5-37.
Good check

## Section 5.6: Several Additional Methods to Obtain and Evaluate an Inverse $K_{\Delta t}$ Transform

As derived, the $K_{\Delta x}$ Transform equation is one of two closely related equations. The other equation is the Inverse $K_{\Delta x}$ Transform equation. Both equations are derived in Section 5.1. The two equations are shown below.

1. <u>The $K_{\Delta t}$ Transform Equation</u>

$$F(s) = K_{\Delta x}[f(t)] = {}_{\Delta x}\int_{0}^{\infty} (1+s\Delta t)^{-\left(\frac{t+\Delta t}{\Delta t}\right)} f(t)\Delta t$$

2. <u>The Inverse $K_{\Delta t}$ Transform Equation</u>

$$f(t) = K_{\Delta t}^{-1}[F(s)] = \frac{1}{2\pi j}\oint_{c} [1+s\Delta t]^{\frac{t}{\Delta t}} F(s)\, ds$$

The $K_{\Delta x}$ Transform Equation converts a function of t, f(t), into a function of s, F(s). The Inverse $K_{\Delta t}$ Transform Equation reconverts the function of s, F(s), into the original function of t, f(t). The discrete integral of the $K_{\Delta t}$ Transform equation is usually not difficult to solve. The function, f(t), is a real value function and a table of discrete integrals is available in the Appendix. On the other hand, the Inverse $K_{\Delta t}$ Transform Equation can present some obstacles to the user. The function of s, F(s), is a complex value function and integration in the complex plane is more difficult. The use of residue theory is required. Integration is possible, but for many, fairly difficult to perform. Where the form of a $K_{\Delta t}$ Transform is that of a polynomial divided by another polynomial, Heavyside's method for finding the necessary residues for integration simplifies the integration process considerably. A polynomial fraction $K_{\Delta t}$ Transform, for which an inverse is required, is often expanded into a partial fraction expansion after its denominator roots are determined. A table of $K_{\Delta t}$ Transforms and their corresponding f(t) functions is then used to find the required $K_{\Delta t}$ Transform inversions. This table is also found in the Appendix. It must be noted that not all $K_{\Delta t}$ Transforms are of a polynomial fraction form. When this is the case, the use of residue theory can be difficult. The following methods for evaluating f(t) from a $K_{\Delta t}$ Transform, $F(s) = K_{\Delta t}[f(t)]$ can be applied to any $K_{\Delta t}$ Transform. Should polynomials be involved, no polynomial roots need be calculated.

This, in some cases, may be an advantage. As will be seen in the following discussions and derivations, asymptotic expansions of $K_{\Delta t}$ Transforms or modified $K_{\Delta t}$ Transforms can be very useful.

The $K_{\Delta t}$ Transform Unit Pulse Series

The $K_{\Delta t}$ Transform Unit Pulse Series is derived directly from the definition of the $K_{\Delta t}$ Transform rewritten below.

Definition of the $K_{\Delta t}$ Transform

$$F(s) = K_{\Delta t}[f(t)] = \int_{0}^{\infty} f(t)(1+s\Delta t)^{-(\frac{t+\Delta t}{\Delta t})} \Delta t \equiv \sum_{n=0}^{\infty} f(n\Delta t)\Delta t(1+s\Delta t)^{-n-1} , \quad t = n\Delta t , \; n = 0,1,2,3, \ldots \quad (5.6\text{-}1)$$

where

$f(t) =$ a function of t

$F(s) = K_{\Delta t}[f(t)]$ , the $K_{\Delta t}$ Transform of the function f(t)

$\Delta t =$ the t increment

$(1+s\Delta t)^{-n-1} =$ the $K_{\Delta t}$ Transform of a unit area pulse where the pulse rises at $t = n\Delta t$ and falls at $t = (n+1)\Delta t$

$n = 0,1,2,3,\ldots$

The $K_{\Delta t}$ Transform Unit Pulse Series is derived by expanding the $K_{\Delta t}$ Transform summation shown in Eq 5.6-1.

The unit area pulse, $\frac{1}{\Delta t}[\, U(t-n\Delta t) - U(t-[n+1]\Delta t \,]$, which has a $K_{\Delta t}$ Transform of $(1+s\Delta t)^{-1-n}$, can be converted into a unit amplitude pulse, $U(t-n\Delta t) - U(t-[n+1]\Delta t$, by multiplying it by $\Delta x$. The unit amplitude pulse has a $K_{\Delta t}$ Transform of $\Delta t(1+s\Delta t)^{-1-n}$.

**The $K_{\Delta t}$ Transform Unit Pulse Series**

$$\mathbf{F(s) = \sum_{n=0}^{\infty} f(n\Delta t)[\Delta t(1+s\Delta t)^{-n-1}] = f(0)[\Delta t(1+s\Delta t)^{-1}] + f(\Delta t)[\Delta t(1+s\Delta t)^{-2}] + f(2\Delta t)[\Delta t(1+s\Delta t)^{-3}] + f(3\Delta t)[\Delta t(1+s\Delta t)^{-4}] + \ldots} \quad (5.6\text{-}2)$$

**where**

**$f(t) =$ a function of t**

**$F(s) = K_{\Delta t}[f(t)]$ , the $K_{\Delta t}$ Transform of the function f(t)**

**$\Delta t =$ the t increment**

**$\Delta t(1+s\Delta t)^{-n-1} =$ the $K_{\Delta t}$ Transform of a unit area pulse where the pulse rises at $t = n\Delta t$ and falls at $t = (n+1)\Delta t$**

**$n = 0,1,2,3,\ldots$**

The coefficients of the $K_{\Delta t}$ Transform Unit Pulse Series, Eq 5.6-2, represent the values of f(t) for $t = n\Delta t$ , $n = 0,1,2,3,\ldots$

Thus, by expanding a $K_{\Delta t}$ Transform, $F(s) = K_{\Delta t}[f(t)]$, as shown in Eq 5.6-2, $f(t) = K_{\Delta x}^{-1}[F(s)]$ can be evaluated.

## The Modified $K_{\Delta t}$ Transform Unit Pulse Asymptotic Series

The $K_{\Delta t}$ Transform Unit Pulse Series, as defined in Eq 5.6-2, has a difficulty in application. The $K_{\Delta t}$ Transform, $F(s) = K_{\Delta x}[f(t)]$, is a function of s and not of $(1+s\Delta t)$. A means is needed to expand F(s) into a series of $(1+s\Delta t)^{-n-1}$ terms where n = 0,1,2,3,…

A modification to the $K_{\Delta t}$ Transform, F(s), can be made. This modification is shown below.

Let

$$s = \frac{(1+s\Delta t) - 1}{\Delta t} = \frac{p - 1}{\Delta t} \qquad (5.6\text{-}3)$$

where

$$p = 1+s\Delta t \qquad (5.6\text{-}4)$$

Substitute Eq 5.6-3 and Eq 5.6-4 into the $K_{\Delta t}$ Transform Unit Pulse Series of Eq 5.6-2 and divide both sides of Eq 5.6-2 by $\Delta t$.

The modified $K_{\Delta t}$ Transform Unit Pulse Series is:

$$\frac{1}{\Delta t} F(s)\Big|_{s = \frac{p-1}{\Delta t}} = M(p) = \sum_{n=0}^{\infty} f(n\Delta t)p^{-n-1} = f(0)p^{-1} + f(\Delta t)p^{-2} + f(2\Delta t)p^{-3} + f(3\Delta t)p^{-4} + f(4\Delta t)p^{-5} + \ldots \qquad (5.6\text{-}5)$$

where

$F(s) = K_{\Delta t}[f(t)]$ , the $K_{\Delta t}$ Transform of the function f(t)
$M(p)$ = Modified $K_{\Delta t}$ Transform
$f(t)$ = a function of t
$\Delta t$ = the t increment
$t = n\Delta t$ , $n = 0,1,2,3,\ldots$
$p = 1+s\Delta t$
$\Delta t\, p^{-n-1}$ = the $K_{\Delta t}$ Transform of a unit pulse where the pulse rises at $t = n\Delta t$ and falls at $t = (n+1)\Delta t$

If the the form of F(s) is that of a polynomial fraction as shown below in Eq 5.6-6, the form of the function M(p) will also be that of a polynomial fraction. For both functions, F(s) and M(p), the maximum order of the numerator polynomial is one less than the order of the denominator polynomial.

A polynomial fraction $K_{\Delta t}$ Transform

$$F(s) = K_{\Delta t}[f(t)] = \frac{a_{m-1}s^{m-1} + a_{m-2}s^{m-2} + a_{m-3}s^{m-3} + \ldots + a_1 s^1 + a_0 s^0}{s^m + b_{m-1}s^{m-1} + b_{m-2}s^{m-2} + b_{m-3}s^{m-3} + \ldots + b_1 s^1 + b_0 s^0} \qquad (5.6\text{-}6)$$

where

$a_c, b_c$ = constants , c = 1,2,3, …

m = order of the denominator polynomial of the $K_{\Delta t}$ Transform, F(s)
The maximum order of the numerator polynomial is one less than the order of the denominator polynomial.

In this case, by dividing the denominator polynomial of M(p) into the numerator polynomial of M(p), an asymptotic series is generated as follows:

$$M(p) = \sum_{n=0}^{\infty} f(n\Delta t)p^{-n-1} = a_1p^{-1} + a_2p^{-2} + a_3p^{-3} + a_4p^{-4} + a_5p^{-5} + \dots \qquad (5.6\text{-}7)$$

where
$a_n$ = constants
n = 1,2,3,…

Comparing Eq 5.6-7 to Eq 5.6-5 the coefficients of the asymptotic series are seen to be the values of f(t) at t = 0,Δt,2Δt,3Δt, …

$$a_{n+1} = f(n\Delta t) \ , \quad n = 0,1,2,3,\dots \qquad (5.6\text{-}8)$$

Some $K_{\Delta t}$ Transforms are not of a polynomial form, for example $F(s) = \frac{1}{\sqrt{s^2+1}}$. If the form of a $K_{\Delta t}$ Transform is not that of a polynomial fraction, neither is the form of the function, M(p).

In this case, an asymptotic series similar to Eq 5.6-7 can be generated using a computer program that expands functions into asymptotic series. The computer generated M(p) asymptotic expansion generates the following asymptotic series:

$$M(p) = \sum_{n=0}^{\infty} f(n\Delta t)p^{-n-1} = a_1p^{-1} + a_2p^{-2} + a_3p^{-3} + a_4p^{-4} + a_5p^{-5} + \dots \qquad (5.6\text{-}9)$$

where
$a_n$ = constants
n = 1,2,3,…

Comparing Eq 5.6-9 to Eq 5.6-5 the coefficients of the asymptotic series are seen to be the values of f(t) at t = 0,Δt,2Δt,3Δt, …

$$a_{n+1} = f(n\Delta t) \ , \quad n = 0,1,2,3,\dots \qquad (5,6\text{-}10)$$

Various computer programs are available to expand functions into asymptotic series. A symbolic logic program such as Maxima can be used or an internet site such as Wolframalpha.com can also be used. Without a computer, the required mathematical computations can be rather tedious.

Then, for any $K_{\Delta t}$ Transform, the modified $K_{\Delta t}$ Transform Unit Pulse Series of Eq 5.6-5 will provide evaluations of f(t) at t = 0,Δt,2Δt,3Δt, … for a $K_{\Delta t}$ Transform, $F(s) = K_{\Delta t}[f(t)]$.

The series of Eq 5.6-5 has been named The Modified $K_{\Delta t}$ Transform Unit Pulse Asymptotic Series. Using this series, derived from the $K_{\Delta t}$ Transform, $F(s) = K_{\Delta t}[f(t)]$, the value of $f(t) = K_{\Delta t}^{-1}[F(s)]$ at t = 0,Δt,2Δt,3Δt,… can be found.

Then
From Eq 5.6-5

The Modified $K_{\Delta t}$ Transform Unit Pulse Asyptotic Series is as follows:

**The Modified $K_{\Delta t}$ Transform Unit Pulse Asymptotic Series**

$$\mathbf{M(p) = \sum_{n=0}^{\infty} f(n\Delta t)p^{-n-1} = f(0)p^{-1} + f(1\Delta t)p^{-2} + f(2\Delta t)p^{-3} + f(3\Delta t)p^{-4} + f(4\Delta t)p^{-5} + \ldots}$$

**where**

$$\mathbf{M(p) = \frac{1}{\Delta t} F(s)\Big|_{s = \frac{p-1}{\Delta t}}} \qquad (5.6\text{-}11)$$

**$F(s) = K_{\Delta t}[f(t)]$ , the $K_{\Delta t}$ Transform of the function f(t)**
**M(p) = Modified $K_{\Delta t}$ Transform**
**f(t) = a function of t**
**$\Delta t$ = the t increment**
**$t = n\Delta t$ , $n = 0,1,2,3,\ldots$**
**$p = 1+s\Delta t$**
**$\Delta t\, p^{-n-1}$ = the $K_{\Delta t}$ Transform of a unit pulse where the pulse rises at $t = n\Delta t$ and falls at $t = (n+1)\Delta t$**

**The Modified $K_{\Delta t}$ Transform Unit Pulse Asymptotic Series, which is derived from a $K_{\Delta t}$ Transform, F(s), finds the values of f(t) where $f(t) = K_{\Delta t}^{-1}[F(s)]$. The coefficients of the Modified $K_{\Delta t}$ Transform Unit Pulse Asymptotic Series represent the values of f(t) at $t = n\Delta t$ , $n = 0,1,2,3,\ldots$**

**Where M(p) is a polynomial fraction, the Modified $K_{\Delta t}$ Transform Unit Pulse Asymptotic Series can be generated by dividing the denominator of M(p) into the numerator of M(p). In general, a computer program that expands functions into asymptotic series can always be used.**

The following two examples, Example 5.6-1 and Example 5.6-2, provide demonstrations of the use of the Modified $K_{\Delta t}$ Transform Unit Pulse Asymptotic Series to find the values of $f(t) = K_{\Delta t}^{-1}[F(s)]$ from its $K_{\Delta t}$ Transform, $F(s) = K_{\Delta t}[f(t)]$.

<u>Example 5.6-1</u> Find the values of f(t) at t = 0, .5, 1, and 1.5 given the $K_{\Delta t}$ Transform, $F(s) = K_{\Delta t}[f(t)] = \frac{2(s-1)}{s(s-2)}$ . Use the Modified $K_{\Delta t}$ Transform Unit Pulse Asymptotic Series.

$$F(s) = K_{\Delta t}[f(t)] = \frac{2(s-1)}{s(s-2)} \qquad (5.6\text{-}12)$$

$$\Delta t = .5 \qquad (5.6\text{-}13)$$

$$M(p) = \frac{1}{\Delta t} F(s)\Big|_{s = \frac{p-1}{\Delta t}} = \sum_{n=0}^{\infty} f(n\Delta t)p^{-n-1} \qquad (5.6\text{-}14)$$

Find the function, M(p)

$$M(p) = \frac{1}{\Delta t} F(s)\Big|_{s = \frac{p-1}{\Delta t}} = \frac{1}{.5}\frac{2(s-1)}{s(s-2)}\Big|_{s = \frac{p-1}{.5}} = \frac{4(s-1)}{s(s-2)}\Big|_{s = 2p-2} \qquad (5.6\text{-}15)$$

$$M(p) = \frac{4(2p-2-1)}{(2p-2)(2p-4)} = \frac{2p-3}{(p-1)(p-2)} = \frac{2p-3}{p^2-3p+2} \qquad (5.6\text{-}16)$$

$$M(p) = \frac{2p-3}{p^2-3p+2} \qquad (5.6\text{-}17)$$

Expand M(p) by dividing its denominator into its numerator

$$M(p) = 2p^{-1} + 3p^{-2} + 5p^{-3} + 9p^{-4} + 17p^{-5} + 33p^{-6} + \ldots \qquad (5.6\text{-}18)$$

Comparing the coefficients of the expansion of Eq 5.6-18 to that of Eq 5.6-14

The values of f(t) at t = 0,1,2,3, … are found to be:

f(0) = 2
f(.5) = 3
f(1) = 5
f(1.5) = 9

Check the above values

Expanding the $K_{\Delta t}$ Transform into a partial fraction expansion

$$K_{\Delta t}[f(t)] = \frac{2(s-1)}{s(s-2)} = \frac{A}{s} + \frac{B}{s-2} \qquad (5.6\text{-}19)$$

$$A = \frac{2s(s-1)}{s(s-2)}\Big|_{s = 0} = 1 \qquad (5.6\text{-}20)$$

$$B = \frac{2(s-2)(s-1)}{s(s-2)}\Big|_{s = 2} = 1 \qquad (5.6\text{-}21)$$

Substiituting Eq 5.6-20 and Eq 5.6-21 into Eq 5.6-19

$$K_{\Delta t}[f(t)] = \frac{1}{s} + \frac{1}{s-2} \qquad (5.6\text{-}22)$$

Taking the Inverse $K_{\Delta t}$ Transform

$$f(t) = 1 + e_{.5}(2,t) = 1 + (1+2[.5])^{\frac{t}{.5}} = 1 + 2^{2t} \tag{5.6-23}$$

$$f(t) = 1 + 2^{2t} \tag{5.6-24}$$

The values of f(t) at t = 0,.5,1, and 1.5 are as follows:

f(0) = 2
f(.5) = 3
f(1) = 5
f(1.5) = 9

These are the same results previously obtained.

Good check

Example 5.6-2 Find the values of f(t) at t = 0, .2, .4, and .6 given the $K_{\Delta t}$ Transform, $F(s) = K_{\Delta t}[f(t)] = \frac{1}{\sqrt{s^2+1}}$ . Use the Modified $K_{\Delta t}$ Transform Unit Pulse Asymptotic Series.

$$F(s) = K_{\Delta t}[f(t)] = \frac{1}{\sqrt{s^2+1}} \tag{5.6-25}$$

$$\Delta t = .2 \tag{5.6-26}$$

$$M(p) = \frac{1}{\Delta t}F(s)\Big|_{s=\frac{p-1}{\Delta t}} = \sum_{n=0}^{\infty} f(n\Delta t)p^{-n-1} \tag{5.6-27}$$

Find the function, M(p)

$$M(p) = \frac{1}{\Delta t}F(s)\Big|_{s=\frac{p-1}{\Delta t}} = \frac{1}{.2}\frac{1}{\sqrt{s^2+1}}\Big|_{s=\frac{p-1}{.2}} = \frac{5}{\sqrt{s^2+1}}\Big|_{s=5p-5} \tag{5.6-28}$$

$$M(p) = \frac{5}{\sqrt{(5p-5)^2+1}} = \frac{5}{\sqrt{25p^2-50p+25+1}} = \frac{5}{\sqrt{25p^2-50p+26}} = \frac{1}{\sqrt{p^2-2p+1.04}} = \sqrt{\frac{1}{p^2-2p+1.04}} \tag{5.6-29}$$

$$M(p) = \sqrt{\frac{1}{p^2-2p+1.04}} \tag{5.6-30}$$

Expanding M(p) into an asymptotic series using the symbolic logic computer program, Maxima

$$M(p) = 1p^{-1} + 1p^{-2} + \frac{49}{50}p^{-3} + \frac{47}{50}p^{-4} + \frac{4403}{5000}p^{-5} + \frac{803}{1000}p^{-6} + \ldots \qquad (5.6\text{-}31)$$

Comparing Eq 5.6-31 to Eq 5.6-27

$f(0) = 1$

$f(.2) = 1$

$f(.4) = \frac{49}{50} = .98$

$f(.6) = \frac{47}{50} = .94$

This same problem is calculated by a different method in Example 5.6-6. The values calculated in Example 5.6-6 are as follows:

$f(0) = 1$

$f(.2) = 1$

$f(.4) = .98$

$f(.6) = .94$

Good check

The Z Transform Unit Impulse Asymptotic Series

As shown previously in Section 1.6, The Z Transform of a function of t can be obtained from the $K_{\Delta t}$ Transform of that same function of t. Instead of using the The Modified $K_{\Delta t}$ Transform Unit Pulse Series, an expansion of the Z Transform into a $z^{-1}$ power series, the Z Transform Unit Impulse Series, can be used. The $K_{\Delta t}$ Transform, $F(s) = K_{\Delta t}[f(t)]$, for which the values of the function, f(t), are to be found, can first be converted into its equivalent Z Transform, Z[f(t)]. Then, this Z Transform can be expanded into an symptotic series to obtain the values of the function, f(t). The equation for the $K_{\Delta t}$ Transform to Z transform conversion and series expansion is as follows:

**The Z Transform Unit Impulse Asymptotic Series**

$$\mathbf{F(z) = \frac{z}{\Delta t} F(s)\Big|_{s = \frac{z-1}{\Delta t}} = \sum_{n=0}^{\infty} f(n\Delta t)z^{-n} = f(0) + f(1\Delta t)z^{-1} + f(2\Delta t)z^{-2} + f(3\Delta t)(z^{-3} + f(4\Delta t)z^{-3} + \ldots} \qquad (5.6\text{-}32)$$

**where**

**$F(s) = K_{\Delta t}[f(t)]$ , the $K_{\Delta t}$ Transform of the function, f(t)**

**$F(z) = Z[f(t)]$ , the Z Transform of the function, f(t)**

**f(t) = a function of t**
**$\Delta t$ = the t increment**
**$t = n\Delta t$ , $n = 0,1,2,3,\ldots$**
**$z = 1+s\Delta t$**

**For the $K_{\Delta t}$ Transform, F(s), the maximum order of the numerator polynomial is one less than the order of the denominator polynomial. For the Z Transform, F(z), the maximum order of the numerator polynomial is the same as the order of the denominator polynomial.**

**The Z Transform Unit Impulse Series, which is derived from a $K_{\Delta t}$ Transform, F(s), finds the values of f(t) where $f(t) = Z^{-1}[F(z)]$. The coefficients of the Z Transform Unit Impulse Series represent the values of f(t) at $t = n\Delta t$ , $n = 0,1,2,3,\ldots$**

Note the similarity of the Z Transform Unit Impulse Asymptotic Series to the the Modified $K_{\Delta t}$ Transform Unit Pulse Series.

The following example, Example 5.6-3, provides a demonstration of the use of the Z Transform Unit Impulse Asymptotic Series to find the values of f(t) from its $K_{\Delta t}$ Transform, $F(s) = K_{\Delta t}[f(t)]$.

Example 5.6-3 Find the values of f(t) at t = 0, .5, 1, and 1.5 given the $K_{\Delta t}$ Transform, $F(s) = K_{\Delta t}[f(t)] = \frac{2(s-1)}{s(s-2)}$. Use the Z Transform Unit Impulse Asymptotic Series.

$$F(s) = K_{\Delta t}[f(t)] = \frac{2(s-1)}{s(s-2)} \qquad (5.6\text{-}33)$$

$$\Delta t = .5 \qquad (5.6\text{-}34)$$

$$F(z) = \frac{z}{\Delta t}F(s)\Big|_{s=\frac{z-1}{\Delta t}} = \sum_{n=0}^{\infty} f(n\Delta t)z^{-n} = f(0) + f(\Delta t)z^{-1} + f(2\Delta t)z^{-2} + f(3\Delta t)(z^{-3} + f(4\Delta t)z^{-3} + \ldots \qquad (5.6\text{-}35)$$

Find the the Z Transform, $F(z) = Z[f(t)]$, from the $K_{\Delta t}$ Transform, $F(s) = K_{\Delta t}[f(t)]$

$$F(z) = \frac{z}{\Delta t}F(s)\Big|_{s=\frac{z-1}{\Delta t}} = \frac{z}{.5}\frac{2(s-1)}{s(s-2)}\Big|_{s=\frac{z-1}{.5}} = \frac{4z(s-1)}{s(s-2)}\Big|_{s=2z-2} = \frac{4z(s-1)}{s(s-2)}\Big|_{s=2z-2} \qquad (5.6\text{-}36)$$

$$F(z) = \frac{4z(2z-3-1)}{(2z-2)(2z-4)} = \frac{4z(2z-2-1)}{(2z-2)(2z-4)} = \frac{z(2z-3)}{(z-1)(z-2)} = \frac{z(2z-3)}{z^2-3z+2} \qquad (5.6\text{-}37)$$

$$F(z) = \frac{2z^2-3z}{z^2-3z+2} \qquad (5.6\text{-}38)$$

Expand F(z) into an asymptotic series by dividing its denominator into its numerator or by using a symbolic logic computer program that expands functions into asymptotic series.

The asymptotic expansion of F(z) is as follows:

$$F(z) = 2 + 3z^{-1} + 5z^{-2} + 9z^{-3} + 17z^{-4} + 33z^{-5} + \ldots\ldots \qquad (5.6\text{-}39)$$

Comparing the coefficients of the expansion of Eq 5.6-39 to that of Eq 5.6-35

The values of f(t) for t = 0,.5,1, and 1.5 are as follows:

f(0) = 2
f(.5) = 3
f(1) = 5
f(1.5) = 9

These are the same values calculated in Example 5.6-1.

f(0) = 2
f(.5) = 3
f(1) = 5
f(1.5) = 9

$$f(t) = 1 + e_{.5}(2,t) = 1 + (1+2[.5])^{\frac{t}{.5}} = 1 + 2^{2t} \qquad (5.6\text{-}40)$$

Good check

<u>The $K_{\Delta t}$ Transform Asymptotic Series</u>

The $K_{\Delta t}$ Transform Asymptotic Series is the series that results from an aymptotic expansion of the $K_{\Delta t}$ Transform, F(s). The asymptotic expansion of the $K_{\Delta t}$ Transform is performed by division or a computer program as previously explained.

The resulting asymptotic series is as follows:

$$F(s) = K_{\Delta t}[f(t)] = \sum_{n=1}^{N} c_n s^{-n} + R_{N+1}(s) = c_1 s^{-1} + c_2 s^{-2} + c_3 s^{-3} + c_4 s^{-4} + c_5 s^{-5} + \ldots + R_{N+1}(s) \qquad (5.6\text{-}41)$$

where
$R_{N+1}(s)$ = remainder term
$N \to \infty$

$R_{N+1}(s) \to ks^{-N-1}$ where $s \to \infty$ , k = the value of the largest term of the remainder numerator divided by the largest term of the remainder denominator. (5.6-42)

The values of the coefficients of this $K_{\Delta t}$ Transform asynptotic series can be very useful. They evaluate all of the discrete derivatives of f(t) at t = 0.

Find the values of $c_n$ , n = 1,2,3,… in Eq 5.6-41.

Find the value of $c_1$

$$f(0) = \lim_{s\to\infty} sF(s) \tag{5.6-43}$$

From Eq 5.6-41 thru Eq 5.6-43

$$f(0) = \lim_{s\to\infty} sF(s) = \lim_{s\to\infty}[\ c_1 + c_2s^{-1} + c_3s^{-2} + c_4s^{-3} + c_5s^{-4} + \ldots + ks^{-N}] \tag{5.6-44}$$

where

$f(0) = f(t)|_{t=0}$

$N \to \infty$

$$c_1 = f(0) \tag{5.6-45}$$

Find the value of $c_2$

$$sF(s) - f(0) = K_{\Delta t}[D_{\Delta t}{}^1 f(t)] \tag{5.6-46}$$

where

$D_{\Delta t}{}^1 f(t) = \dfrac{f(t + \Delta t) - f(t)}{\Delta t}$ , the first discrete derivative of $f(t)$

From Eq 5.6-41, Eq 42, Eq 5.6-45 and Eq 5.6-46

$$sF(s) - f(0) = K_{\Delta t}[D_{\Delta t}{}^1 f(t)] = [f(0) + c_2s^{-1} + c_3s^{-2} + c_4s^{-3} + c_5s^{-4} + \ldots + ks^{-N}] - f(0) \tag{5.6-47}$$

$$K_{\Delta t}[D_{\Delta t}{}^1 f(t)] = [c_2s^{-1} + c_3s^{-2} + c_4s^{-3} + c_5s^{-4} + \ldots + ks^{-N}] \tag{5.6-48}$$

$$D_{\Delta t}{}^1 f(0) = \lim_{s\to\infty} sK_{\Delta t}[D_{\Delta t}{}^1 f(t)] \tag{5.6-49}$$

From Eq 5.6-48 and 5.6-49

$$D_{\Delta t}{}^1 f(0) = \lim_{s\to\infty} sK_{\Delta t}[D_{\Delta t}{}^1 f(t)] = \lim_{s\to\infty} [c_2 + c_3s^{-1} + c_4s^{-2} + c_5s^{-3} + \ldots + ks^{-N+1}] \tag{5.6-50}$$

$$c_2 = D_{\Delta t}{}^1 f(0) \tag{5.6-51}$$

Find the value of $c_3$

$$s^2F(s) - f(0)s - D_{\Delta t}{}^1 f(0) = K_{\Delta t}[D_{\Delta t}{}^2 f(t)] \tag{5.6-52}$$

where

$D_{\Delta t}{}^2 f(t)$ = the second discrete derivative of $f(t)$

From Eq 5.6-41, Eq 42, Eq 5.6-45, Eq 5.6-51 and Eq 5.6-52

$$s^2F(s) - f(0)s - D_{\Delta t}{}^1 f(0) = K_{\Delta t}[D_{\Delta t}{}^2 f(t)] = [f(0)s + D_{\Delta t}{}^1 f(0) + c_3s^{-1} + c_4s^{-2} + c_5s^{-3} + \ldots + ks^{-N+1}] - f(0)s - D_{\Delta t}{}^1 f(0) \tag{5.6-53}$$

$$K_{\Delta t}[D_{\Delta t}{}^2 f(t)] = [c_3s^{-1} + c_4s^{-2} + c_5s^{-3} + \ldots + ks^{-N+1}] \tag{5.6-54}$$

$D_{\Delta t}{}^2 f(0) = \lim_{s\to\infty} sK_{\Delta t}[D_{\Delta t}{}^2 f(t)]$ (5.6-55)

From Eq 5.6-54 and 5.6-55

$D_{\Delta t}{}^2 f(0) = \lim_{s\to\infty} sK_{\Delta t}[D_{\Delta t}{}^2 f(t)] = \lim_{s\to\infty} [c_3 + c_4 s^{-1} + c_5 s^{-2} + \ldots + ks^{-N+2}]$ (5.6-56)

$c_3 = D_{\Delta t}{}^2 f(0)$ (5.6-57)

.
.
.

From Eq 5.6-41, Eq 5.6-45, Eq 5.6-51, and Eq 5.6-57

$$F(s) = K_{\Delta t}[f(t)] = \sum_{n=1}^{N} D_{\Delta t}{}^{n-1} f(0)s^{-n} + R_{N+1}(s) = f(0)s^{-1} + D_{\Delta t}{}^1 f(0)s^{-2} + D_{\Delta t}{}^2 f(0)s^{-3} + D_{\Delta t}{}^3 f(0)s^{-4} + \ldots + R_{N+1}(s)$$ (5.6-58)

where
- $N \to \infty$
- $D_{\Delta t}{}^m f(0) = D_{\Delta t}{}^m f(t)|_{t=0}$, $m = 0,1,2,3, \ldots$
- $F(s) = K_{\Delta t}[f(t)]$ = the $K_{\Delta t}$ Transform of f(t)
- $R_{N+1}(s)$ = remainder term
- $R_{N+1}(s) \to 0$ for $N \to \infty$

Note – The values of of the discrete derivatives remain the same for all values of $\Delta t$.
This can be seen from the following example:

$$D_{\Delta t} e_{\Delta t}(a,x) = D_{\Delta t}(1+a\Delta t)^{\frac{t}{\Delta t}} = \frac{(1+a\Delta t)^{\frac{t+\Delta t}{\Delta t}} - (1+a\Delta t)^{\frac{t}{\Delta t}}}{\Delta t} = \frac{(1+a\Delta t)^{\frac{t}{\Delta t}}(1+a\Delta t - 1)}{\Delta t} = a(1+a\Delta t)^{\frac{t}{\Delta t}} = ae_{\Delta t}(a,x)$$

Then
The $K_{\Delta t}$ Transform Asymptotic Series is defined as follows:

**The $K_{\Delta t}$ Transform Asymptotic Series**

$$\mathbf{F(s) = K_{\Delta t}[f(t)] = \sum_{n=1}^{\infty} D_{\Delta t}{}^{n-1} f(0)s^{-n} = f(0)s^{-1} + D_{\Delta t}{}^1 f(0)s^{-2} + D_{\Delta t}{}^2 f(0)s^{-3} + D_{\Delta t}{}^3 f(0)s^{-4} + \ldots}$$

**where** (5.6-59)
- **$D_{\Delta t}{}^m f(0) = D_{\Delta t}{}^m f(t)|_{t=0}$ , $m = 0,1,2,3,\ldots$ , Discrete derivatives of f(t) evaluated at t = 0**
- **f(t) = function of t**
- **$F(s) = K_{\Delta t}[f(t)] = K_{\Delta t}$ Transform of the function, f(t)**

**The asymptotic expansion of the $K_{\Delta t}$ Transform is obtained by division or by the use of a computer program such as Maxima. The internet site, Wolframalpha.com can also be used.**

**Note – The $K_{\Delta t}$ Transform Asymptotic Series is a series expansion of a $K_{\Delta t}$ Transform at x = ∞.**

Note the following derivation and the relationship of an asymptotic expansion of a $K_{\Delta t}$ Transform, $F(s) = K_{\Delta t}[f(t)]$, to the discrete Maclaurin Series of the function, f(t).

F(s) can be expanded into an asymptotic series of the following form.

$$F(s) = K_{\Delta t}[f(t)] \;=\; \frac{a_1}{s} + \frac{a_2}{s^2} + \frac{a_3}{s^3} + \frac{a_4}{s^4} + \frac{a_5}{s^5} + \dots \qquad (5.6\text{-}60)$$

where

$a_n$ = constants , $n = 1,2,3, \dots$

Changing the form of Eq 5.6-60

$$F(s) = K_{\Delta t}[f(t)] \;=\; \frac{a_1}{0!}\left(\frac{0!}{s}\right) + \frac{a_2}{1!}\left(\frac{1!}{s^2}\right) + \frac{a_3}{2!}\left(\frac{2!}{s^3}\right) + \frac{a_4}{3!}\left(\frac{3!}{s^4}\right) + \frac{a_5}{4!}\left(\frac{4!}{s^5}\right) + \dots \qquad (5.6\text{-}61)$$

$$K_{\Delta t}^{-1}[F(s)] = f(t) \qquad (5.6\text{-}62)$$

From Eq 5.6-61 and Eq 5.6-62

$$f(t) = K_{\Delta t}^{-1}[F(s)] \;=\; \frac{a_1}{0!}K_{\Delta t}^{-1}[\frac{0!}{s}] + \frac{a_2}{1!}K_{\Delta t}^{-1}[\frac{1!}{s^2}] + \frac{a_3}{2!}K_{\Delta t}^{-1}[\frac{2!}{s^3}] + \frac{a_4}{3!}K_{\Delta t}^{-1}[\frac{3!}{s^4}] + \frac{a_5}{4!}K_{\Delta t}^{-1}[\frac{4!}{s^5}] + \dots \qquad (5.6\text{-}63)$$

$$K_{\Delta t}^{-1}[\frac{n!}{s^{n+1}}] = [t]_{\Delta t}^{n} \qquad (5.6\text{-}64)$$

$$[t]_{\Delta t}^{n} = 1 \quad , \quad n = 0 \qquad (5.6\text{-}65)$$

$$[t]_{\Delta t}^{n} = \prod_{m=1}^{n}(t-[m-1]\Delta t) \quad , \quad n = 1,2,3, \dots \qquad (5.6\text{-}66)$$

From Eq 5.6-63 thru 5.6-66

$$f(t) = K_{\Delta t}^{-1}[F(s)] \;=\; \frac{a_1}{0!}[t]_{\Delta t}^{0} + \frac{a_2}{1!}[t]_{\Delta t}^{1} + \frac{a_3}{2!}[t]_{\Delta t}^{2} + \frac{a_4}{3!}[t]_{\Delta t}^{3} + \frac{a_5}{4!}[t]_{\Delta t}^{4} + \dots \qquad (5.6\text{-}67)$$

Expanding Eq 5.6-67 using Eq 5.6-66

$$f(t) = \frac{a_1}{0!}\,1 + \frac{a_2}{1!}\,t + \frac{a_3}{2!}\,t(t-\Delta t) + \frac{a_4}{3!}\,t(t-\Delta t)(t-2\Delta t) + \frac{a_5}{4!}\,t(t-\Delta t)(t-2\Delta t)(t-3\Delta t) + \dots \qquad (5.6\text{-}68)$$

Eq 5.6-68 above is recognized as the discrete Maclaurin Series of f(t).

The Discrete Maclaurin Series is:

$$f(t) = \sum_{m=0}^{\infty} b_m [t]_{\Delta t}^{m} \qquad (5.6\text{-}69)$$

$$b_m = \frac{1}{m!} D_{\Delta t}{}^{m} f(t)\Big|_{t=0}\ , \text{ initial condition discrete derivatives} \qquad (5.6\text{-}70)$$

where

$$D_{\Delta t} f(t) = \frac{f(t+\Delta t) - f(t)}{\Delta t}\ , \text{ discrete derivative of } f(t)$$

$$[t]_{\Delta t}^{m} = 1\ ,\ m = 0$$

$$[t]_{\Delta t}^{m} = \prod_{r=1}^{m} (t-[r-1]\Delta t)\ ,\ m = 1,2,3,\dots$$

$t = 0, \Delta t, 2\Delta t, 3\Delta t, \dots$

$D_{\Delta t}{}^{m} f(t)$ = mth discrete derivative of f(t)

f(t) = function of t

Then

From the expansion of the $K_{\Delta t}$ Transform, $F(s) = K_{\Delta t}[f(t)]$, into an asymptotic series where

$$F(s) = K_{\Delta t}[f(t)] = \frac{a_1}{s} + \frac{a_2}{s^2} + \frac{a_3}{s^3} + \frac{a_4}{s^4} + \frac{a_5}{s^5} + \dots\ ,\quad a_n = \text{constants}\ ,\quad n = 1,2,3,\dots \qquad (5.6\text{-}71)$$

the discrete derivatives of the function, f(t), can be derived.

$$D_{\Delta t}{}^{n-1} f(t)|_{t=0} = a_n\ ,\quad n = 1,2,3,\dots\ ,\ \text{Discrete derivatives of } f(t) \text{ at } t = 0 \qquad (5.6\text{-}72)$$

where

$$D_{\Delta t}{}^{0} f(t)|_{t=0} = f(0)$$

The coefficients of a $K_{\Delta t}$ Transform, $F(s) = K_{\Delta t}[f(t)]$, asymptotic series evaluate the discrete derivatives of f(t) at t = 0 . From these discrete derivatives, the values of f(t) at $t = 0, \Delta t, 2\Delta t, 3\Delta t, \dots$ can be obtained from the definitions of the discrete derivatives shown below.

$$D_{\Delta t}{}^{0} f(0) = f(0) \qquad (5.6\text{-}73)$$

$$D_{\Delta t}{}^{1} f(0) = \frac{f(\Delta t) - f(0)}{\Delta t} \qquad (5.6\text{-}74)$$

$$D_{\Delta t}{}^{2} f(0) = \frac{f(2\Delta t) - 2f(\Delta t) + f(0)}{\Delta t^2} \qquad (5.6\text{-}75)$$

$$D_{\Delta t}^{3}f(0) = \frac{f(3\Delta t) - 3f(2\Delta t) + 3f(\Delta t) - f(0)}{\Delta t^{3}} \qquad (5.6\text{-}76)$$

$$D_{\Delta t}^{4}f(0) = \frac{f(4\Delta t) - 4f(3\Delta t) + 6f(2\Delta t) - 4f(\Delta t) + f(0)}{\Delta t^{4}} \qquad (5.6\text{-}77)$$

.
.
.

The following example, Example 5.6-4, provides a demonstration of the use of the $K_{\Delta t}$ Transform, $F(s) = K_{\Delta t}[f(t)]$, Asymptotic Series to find the values of f(t).

Example 5.6-4  Find the values of f(t) for t = 0, .5, 1, and 1.5 given the $K_{\Delta t}$ Transform, $F(s) = K_{\Delta t}[f(t)] = \frac{2(s-1)}{s(s-2)}$. Use the $K_{\Delta t}$ Transform Asymptotic Series

$$F(s) = K_{\Delta t}[f(t)] = \frac{2(s-1)}{s(s-2)} = \frac{2s-2}{s^2-2s} \qquad (5.6\text{-}78)$$

Dividing the denominator of F(s) into the numerator of F(s)

$$F(s) = K_{\Delta t}[f(t)] = \frac{2s-2}{s^2-2s} = 2s^{-1} + 2s^{-2} + 4s^{-3} + 8s^{-4} + 16s^{-5} + \ldots \qquad (5.6\text{-}79)$$

The $K_{\Delta t}$ Transform Asymptotic Series is:

$$F(s) = K_{\Delta t}[f(t)] = f(0)s^{-1} + D_{\Delta t}^{1}f(0)s^{-2} + D_{\Delta t}^{2}f(0)s^{-3} + D_{\Delta t}^{3}f(0)s^{-4} + D_{\Delta t}^{4}f(0)s^{-5} + \ldots \qquad (5.6\text{-}80)$$

Comparing Eq 5.6-79 to Eq 5.6-80

$f(0) = 2$ (5.6-81)
$D_{\Delta t}^{1}f(0) = 2$ (5.6-82)
$D_{\Delta t}^{2}f(0) = 4$ (5.6-83)
$D_{\Delta t}^{3}f(0) = 8$ (5.6-84)

$\Delta t = .5$ (5.6-85)

Find f(.5)

$$D_{\Delta t}^{1}f(0) = \frac{f(\Delta t) - f(0)}{\Delta t} \qquad (5.6\text{-}86)$$

From Eq 81, Eq 82, and Eq 85

$$2 = \frac{f(.5) - 2}{.5} \tag{5.6-87}$$

$$f(.5) - 2 = 1 \tag{5.6-88}$$

$$f(.5) = 3 \tag{5.6-89}$$

Find f(1)

$$D_{\Delta t}^{2} f(0) = \frac{f(2\Delta t) - 2f(\Delta t) + f(0)}{\Delta t^2} \tag{5.6-90}$$

From Eq 81, Eq 83, Eq 85, Eq 89

$$4 = \frac{f(1) - 2(3) + 2}{.5^2} \tag{5.6-91}$$

$$f(1) - 6 + 2 = 1 \tag{5.6-92}$$

$$f(1) = 5 \tag{5.6-93}$$

Find f(1.5)

$$D_{\Delta t}^{3} f(0) = \frac{f(3\Delta t) - 3f(2\Delta t) + 3f(\Delta t) - f(0)}{\Delta t^3} \tag{5.6-94}$$

From Eq 81, Eq 84, Eq 85, Eq 89, and Eq 93

$$8 = \frac{f(1.5) - 3(5) + 3(3) - 2}{.5^3} \tag{5.6-95}$$

$$f(1.5) - 15 + 9 - 2 = 1 \tag{5.6-96}$$

$$f(1.5) = 9 \tag{5.6-97}$$

From Eq 81, Eq 89, Eq 93, and Eq 97

$f(0) = 2$
$f(.5) = 3$
$f(1) = 5$
$f(1.5) = 9$

These are the same values calculated in Example 5.6-1 and Example 5.6-3

f(0) = 2
f(.5) = 3
f(1) = 5
f(1.5) = 9

$$f(t) = 1 + e_{.5}(2,t) = 1 + (1+2[.5])^{\frac{t}{.5}} = 1 + 2^{2t} \quad (5.6\text{-}98)$$

Good check

After the discrete derivatives of a f(t) are known, there is another way to find the values of f(t) at t = 0,Δt,2Δt,3Δt,… The discrete Maclaurin Series can be used.

Evaluating $f(t) = K_{\Delta t}^{-1}[F(s)]$ using the $K_{\Delta t}$ Transform Asymptotic Series and the Discrete Maclaurin Series

If the $K_{\Delta t}$ Transform Asymptotic Series is used to obtain the initial condition discrete derivatives of f(t) from the $K_{\Delta t}$ Transform, $F(s) = K_{\Delta t}[f(t)]$, f(t) can be evaluated at t = 0,Δt,2Δt,3Δt,… using the discrete Maclaurin Series.

The discrete Maclaurin Series is as follows:

The Discrete Maclaurin Series

$$f(t) = \sum_{m=0}^{\infty} b_m [t]_{\Delta t}^{m} \quad (5.6\text{-}99)$$

$$b_m = \frac{1}{m!} D_{\Delta t}^{m} f(t)\Big|_{t=0} \text{ , initial condition discrete derivatives} \quad (5.6\text{-}100)$$

where

$$D_{\Delta t}f(t) = \frac{f(t+\Delta t) - f(t)}{\Delta t} \text{ , discrete derivative of f(t)}$$

$$[t]_{\Delta t}^{m} = 1 \text{ , } m = 0$$

$$[t]_{\Delta t}^{m} = \prod_{r=1}^{m}(t-[r-1]\Delta t) \text{ , } m = 1,2,3,\ldots$$

t = 0, Δx, 2Δx, 3Δx, …

$D_{\Delta t}^{m}f(t)$ = mth discrete derivative of f(t)

f(t) = function of t

The initial condition discrete derivatives of f(t) are first found from the $K_{\Delta t}$ Transform, $F(s) = K_{\Delta t}[f(t)]$, using a $K_{\Delta t}$ Transform expansion into an asyptotic series.
The resulting $K_{\Delta t}$ Transform Asyptotic Series is:

$$F(s) = K_{\Delta t}[f(t)] = \sum_{n=1}^{N} c_n s^{-n} = c_1 s^{-1} + c_2 s^{-2} + c_3 s^{-3} + c_4 s^{-4} + c_5 s^{-5} + \ldots \qquad (5.6\text{-}101)$$

$$D_{\Delta t}^{\,n-1} f(t)|_{t=0} = c_n, \quad n = 1,2,3,\ldots \quad , \text{ Discrete derivatives of } f(t) \text{ at } t = 0 \qquad (5.6\text{-}102)$$

where

$$D_{\Delta t}^{\,0} f(t)|_{t=0} = f(0)$$

Finding f(t)

The evaluated values of the initial condition discrete derivatives and the value of $\Delta t$ are entered into a discrete Maclaurin Series.

$$f(t) = \sum_{m=0}^{\infty} \frac{1}{m!} D_{\Delta t}^{\,m} f(t)\Big|_{t=0} [t]_{\Delta t}^{m} \qquad (5.6\text{-}103)$$

Expanding Eq 5.6-103

$$f(t) = f(0) + \frac{D_{\Delta t}^{\,1} f(0)}{1!} t + \frac{D_{\Delta t}^{\,2} f(0)}{2!} t(t-\Delta t) + \frac{D_{\Delta t}^{\,3} f(0)}{3!} t(t-\Delta t)(t-2\Delta t) + \frac{D_{\Delta t}^{\,4} f(0)}{4!} t(t-\Delta t)(t-2\Delta t)(t-3\Delta t) + \ldots \qquad (5.6\text{-}104)$$

The above discrete Maclaurin Series is used to evaluate f(t)

If the discrete Maclaurin Series is truncated at the rth term, the series will calculate exact values for f(t) at $t = 0, \Delta t, 2\Delta t, 3\Delta t, \ldots, (r-1)\Delta t$. In general, for larger values of t there will be some error and series convergence must be investigated.

The following examples, Example 5.6-5 and Example 5.6-6 , provide a demonstration of the use of the $K_{\Delta t}$ Transform Asymptotic Series and the discrete Maclaurin Series to find the values of f(t) at $t = 0, \Delta t, 2\Delta t, 3\Delta t, \ldots$ from the $K_{\Delta t}$ Transform, $F(s) = K_{\Delta t}[f(t)]$.

<u>Example 5.6-5</u> Find the values of f(t) for t = 0, .5, 1, and 1.5 given the $K_{\Delta t}$ Transform, $F(s) = K_{\Delta t}[f(t)] = \frac{2(s-1)}{s(s-2)}$. Use the $K_{\Delta t}$ Transform Asymptotic Series and the discrete Maclaurin Series.

$$F(s) = K_{\Delta t}[f(t)] = \frac{2(s-1)}{s(s-2)} = \frac{2s-2}{s^2-2s} \qquad (5.6\text{-}105)$$

Divide the denominator of F(s) into the numerator of F(s)

$$F(s) = K_{\Delta t}[f(t)] = \frac{2s-2}{s^2-2s} = 2s^{-1} + 2s^{-2} + 4s^{-3} + 8s^{-4} + 16s^{-5} + \dots \tag{5.6-106}$$

The $K_{\Delta t}$ Transform Asymptotic Series is:

$$F(s) = K_{\Delta t}[f(t)] = f(0)s^{-1} + D_{\Delta t}^{1}f(0)s^{-2} + D_{\Delta t}^{2}f(0)s^{-3} + D_{\Delta t}^{3}f(0)s^{-4} + D_{\Delta t}^{4}f(0)s^{-5} + \dots \tag{5.6-107}$$

Comparing Eq 5.6-107 to Eq 5.6-106

$$f(0) = 2 \tag{5.6-108}$$
$$D_{\Delta t}^{1}f(0) = 2 \tag{5.6-109}$$
$$D_{\Delta t}^{2}f(0) = 4 \tag{5.6-110}$$
$$D_{\Delta t}^{3}f(0) = 8 \tag{5.6-111}$$
$$D_{\Delta t}^{4}f(0) = 16 \tag{5.6-112}$$

$$\Delta t = .5 \tag{5.6-113}$$

Use a discrete Maclaurin Series to evaluate the function f(t)

$$f(t) = f(0) + \frac{D_{\Delta t}^{1}f(0)}{1!}t + \frac{D_{\Delta t}^{2}f(0)}{2!}t(t-\Delta t) + \frac{D_{\Delta t}^{3}f(0)}{3!}t(t-\Delta t)(t-2\Delta t) + \frac{D_{\Delta t}^{4}f(0)}{4!}t(t-\Delta t)(t-2\Delta t)(t-3\Delta t) + \dots \tag{5.6-114}$$

From Eq 5.6-108 thru Eq 5.6-114

$$f(t) = 2 + 2t + 2t(t-.5) + \frac{4}{3}t(t-.5)(t-1) + \frac{2}{3}t(t-.5)(t-1)(t-1.5) + \dots \tag{5.6-115}$$

$$f(0) = 2 \tag{5.6-116}$$

$$f(.5) = 2+1 = 3 \tag{5.6-117}$$

$$f(1) = 2 + 2(1) + 2(1)(.5) = 5 \tag{5.6-118}$$

$$f(1.5) = 2 + 2(1.5) + 2(1.5)(1) + \frac{4}{3}(1.5)(1)(.5) = 2 + 3 + 3 + 1 = 9 \tag{5.6-119}$$

From Eq 116 thru Eq 119

f(0) = 2
f(.5) = 3
f(1) = 5
f(1.5) = 9

These are the same values calculated in Example 5.6-1, Example 5.6-3, and Example 5.6-4

$f(0) = 2$
$f(.5) = 3$
$f(1) = 5$
$f(1.5) = 9$

$$f(t) = 1 + e_{.5}(2,t) = 1 + (1+2[.5])^{\frac{t}{.5}} = 1 + 2^{2t} \qquad (5.6\text{-}120)$$

Good check

Example 5.6-6 Find the values of f(t) for t = 0, .2, .4, and .6 given the $K_{\Delta t}$ Transform, $F(s) = K_{\Delta t}[f(t)] = \frac{1}{\sqrt{s^2+1}}$. Use the $K_{\Delta t}$ Transform Asymptotic Series and the discrete Maclaurin Series.

$$F(s) = K_{\Delta t}[f(t)] = \frac{1}{\sqrt{s^2+1}} = \sqrt{\frac{1}{s^2+1}} \qquad (5.6\text{-}121)$$

Expand F(s) into a $K_{\Delta t}$ Transform Asymptotic Series

$$F(s) = K_{\Delta t}[f(t)] = \sqrt{\frac{1}{s^2+1}} = 1s^{-1} - \frac{1}{2}s^{-3} + \frac{3}{8}s^{-5} - \frac{5}{16}s^{-7} + \ldots \qquad (5.6\text{-}122)$$

The above asymptotic expansion of F(s), Eq 5.6-122, was obtained using the internet site, wolframalpha.com.
The $K_{\Delta t}$ Transform Asymptotic Series is:

$$F(s) = K_{\Delta t}[f(t)] = f(0)s^{-1} + D_{\Delta t}{}^{1}f(0)s^{-2} + D_{\Delta t}{}^{2}f(0)s^{-3} + D_{\Delta t}{}^{3}f(0)s^{-4} + D_{\Delta t}{}^{4}f(0)s^{-5} + \ldots \qquad (5.6\text{-}123)$$

Comparing Eq 5.6-123 to Eq 5.6-122

$$f(0) = 1 \qquad (5.6\text{-}124)$$

$$D_{\Delta t}{}^{2}f(0) = -\frac{1}{2} \qquad (5.6\text{-}125)$$

$$D_{\Delta t}{}^{4}f(0) = \frac{3}{8} \qquad (5.6\text{-}126)$$

$$D_{\Delta t}{}^{6}f(0) = -\frac{5}{16} \qquad (5.6\text{-}127)$$

$$\Delta t = .2 \qquad (5.6\text{-}128)$$

Use a discrete Maclaurin Series to evaluate the function f(t)

$$f(t) = f(0) + \frac{D_{\Delta t}{}^{1}f(0)}{1!}t + \frac{D_{\Delta t}{}^{2}f(0)}{2!}t(t-\Delta t) + \frac{D_{\Delta t}{}^{3}f(0)}{3!}t(t-\Delta t)(t-2\Delta t) + \frac{D_{\Delta t}{}^{4}f(0)}{4!}t(t-\Delta t)(t-2\Delta t)(t-3\Delta t) + \ldots \qquad (5.6\text{-}129)$$

From Eq 5.6-124 thru Eq 5.6-126, Eq 128, and Eq 129

$$f(t) = 1 + \frac{1}{2(2!)} t(t-.2) + \frac{3}{8(4!)} t(t-.2)(t-.4)(t-.6) + \ldots \qquad (5.6\text{-}130)$$

$$f(0) = 1 \qquad (5.6\text{-}131)$$

$$f(.2) = 1 \qquad (5.6\text{-}132)$$

$$f(.4) = 1 - \frac{1}{4}(.4)(.2) = 1 - .02 = .98 \qquad (5.6\text{-}133)$$

$$f(.6) = 2 - \frac{1}{4}(.6)(.4) = 1 - .06 = .94 \qquad (5.6\text{-}134)$$

From Eq 5.6-131 thru Eq 5.6-134

$f(0) = 1$
$f(.2) = 1$
$f(.4) = .98$
$f(.6) = .94$

Checking

In Example 5.6-2 this same problem was calculated in a different way.
The above values are the same values calculated in Example 5.6-2.

$f(0) = 1$
$f(.2) = 1$
$f(.4) = .98$
$f(.6) = .94$

Good check

Comment - This calculation method derives a Discrete Maclaurin Series representation of $f(t) = K_{\Delta t}^{-1}[K_{\Delta t}[f(t)]]$. In this form, a Discrete Maclaurin Series representation of the discrete derivatives and integral of f(t) can be easily obtained using discrete differentiation and discrete integration.

There is another method for finding the Inverse of a $K_{\Delta t}$ Transform using the $K_{\Delta t}$ Transform Asymptotic Series. This method is derived from the Discrete Maclaurin Series.

Evaluating $f(t) = K_{\Delta t}^{-1}[F(s)]$ using the $K_{\Delta t}$ Transform Asymptotic Series and a formula derived from the Discrete Maclaurin Series

If the $K_{\Delta t}$ Transform Asymptotic Series is used to obtain the initial condition discrete derivatives of f(t) from the $K_{\Delta t}$ Transform, $F(s) = K_{\Delta t}[f(t)]$, f(t) can be evaluated at $t = 0,\Delta t,2\Delta t,3\Delta t,\ldots$ using the formula derived below.

$$K_{\Delta t}[f(t)] = F(s) \qquad (5.6\text{-}135)$$

Changing the form of F(s) to that of an Asymptotic Maclaurin Series

$$K_{\Delta t}[f(t)] = F(s) = a_1 s^{-1} + a_2 s^{-2} + a_3 s^{-3} + a_4 s^{-4} + \dots \tag{5.6-136}$$

The $K_{\Delta t}$ Transform Asymptotic Series for $K_{\Delta t}[\, f(t)]$ is:

$$K_{\Delta t}[f(t)] = f(0)s^{-1} + D_{\Delta t}{}^{1}f(0)s^{-2} + D_{\Delta t}{}^{2}f(0)s^{-3} + D_{\Delta t}{}^{3}f(0)s^{-4} + D_{\Delta t}{}^{4}f(0)s^{-5} + \dots \tag{5.6-137}$$

Comparing Eq 5.6-137 to Eq 5.6-136

The initial conditions of f(t) are:

$$f(0) = a_1 \tag{5.6-138}$$
$$D_{\Delta t}{}^{1}f(0) = a_2 \tag{5.6-139}$$
$$D_{\Delta t}{}^{2}f(0) = a_3 \tag{5.6-140}$$
$$D_{\Delta t}{}^{3}f(0) = a_4 \tag{5.6-141}$$
$$D_{\Delta t}{}^{4}f(0) = a_5 \tag{5.6-142}$$
…

The Discrete Maclaurin Series is:

$$f(t) = \sum_{p=0}^{\infty} \frac{1}{p!} D_{\Delta t}{}^{p} f(0) \prod_{r=1}^{p} (t-[\,r-1]\Delta t) \tag{5.6-143}$$

where

$$D_{\Delta t}f(t) = \frac{f(t+\Delta t) - f(t)}{\Delta t}\,, \text{ discrete derivative of } f(t)$$

$$[t]_{\Delta t}^{0} = \prod_{r=1}^{0}(t-[\,r-1]\Delta t) = 1$$

$$[t]_{\Delta t}^{m} = \prod_{r=1}^{m}(t-[r-1]\Delta t) \;, \; m = 1,2,3,\dots$$

$t = 0, \Delta t, 2\Delta t, 3\Delta t, \dots$

$D_{\Delta t}{}^{p}f(t)$ = pth discrete derivative of f(t)

f(t) = function of t

Expanding Eq 5.6-143

$$f(t) = f(0) + \frac{D_{\Delta t}\, f(0)}{1!}\, t + \frac{D_{\Delta t}{}^{2}f(0)}{2!}\, t(t-\Delta t) + \frac{D_{\Delta t}{}^{3}f(0)}{3!}\, t(t-\Delta t)(t-2\Delta t) + \dots \tag{5.6-144}$$

where

$t = 0, \Delta t, 2\Delta t, 3\Delta t, \dots$

Substituting into Eq 5.6-144 the initial conditions of f(t)

$$f(t) = a_1 + a_2 t + \frac{a_3}{2!} t(t-\Delta t) + \frac{a_4}{3!} t(t-\Delta t)(t-2\Delta t) + \frac{a_5}{4!} t(t-\Delta t)(t-2\Delta t)(t-3\Delta t) + \ldots \tag{5.6-145}$$

where

$t = 0, \Delta t, 2\Delta t, 3\Delta t, \ldots$

From Eq 5.6-145

For $t = 0$

$$f(0) = a_1 \tag{5.6-146}$$

For $t = \Delta t$

$$f(\Delta t) = a_1 + a_2 \Delta t \tag{5.6-147}$$

For $t = 2\Delta t$

$$f(2\Delta t) = a_1 + a_2 2\Delta t + a_3 \Delta t^2 \tag{5.6-148}$$

For $t = 3\Delta t$

$$f(3\Delta t) = a_1 + a_2 3\Delta t + a_3 3\Delta t^2 + a_4 \Delta t^3 \tag{5.6-149}$$

For $t = 4\Delta t$

$$f(4\Delta t) = a_1 + a_2 4\Delta t + a_3 6\Delta t^2 + a_4 4\Delta t^3 + a_5 \Delta t^4 \tag{5.6-150}$$

.
.
.

In general

$$f(n\Delta t) = \sum_{p=0}^{n} [{}_nC_p \Delta t^p]\, a_{p+1} \tag{5.6-151}$$

where

$n = 0, 1, 2, 3, \ldots$

$p = 0, 1, 2, 3, \ldots, n$

$nC_p = \frac{n!}{r!(n-r)!}$

$0! = 1$

Then

**The Inverse $K_{\Delta t}$ Transform Formula for a function, f(t), is:**

$$\mathbf{K_{\Delta t}^{-1}[K_{\Delta t}[f(t)]] = f(n\Delta t) = \sum_{p=0}^{n}[{}_nC_p\Delta t^p]\, a_{p+1}} \qquad (5.6\text{-}152)$$

**where**

**$K_{\Delta t}[f(t)] = F(s) = a_1s^{-1} + a_2s^{-2} + a_3s^{-3} + a_4s^{-4} + \dots$ , Asymptotic Series form of F(s)**
**$F(s) = K_{\Delta t}[f(t)] = K_{\Delta t}$ Transform of the function, f(t)**
**$K_{\Delta t}^{-1}[F(s)] = f(t) =$ Inverse $K_{\Delta t}$ Transform of the function, F(s)**
**$n = 0, 1, 2, 3, 4, \dots$**
**$p = 0, 1, 2, 3, \dots , n$**
**$nC_p = \dfrac{n!}{r!(n-r)!}$**
**$0! = 1$**
**f(t) = function of t**
**$\Delta t$ = interval between successive values of t**
**$t = 0, \Delta t, 2\Delta t, 3\Delta t, \dots$**
**$D_{\Delta t}^{m}f(0) = a_{m+1}$ , $m = 0,1,2,3,\dots$ , Initial conditions of f(t)**

Example 5.6-7 below uses the Inverse $K_{\Delta t}$ Transform Formula , Eq 5.6-152, to solve Example 5.6-6 in a different way.

Example 5.6-7 Find the values of f(t) for t = 0, .2, .4, and .6 given the $K_{\Delta t}$ Transform,
$F(s) = K_{\Delta t}[f(t)] = \dfrac{1}{\sqrt{s^2+1}}$ . Use the Inverse $K_{\Delta t}$ Transform Formula.

$$F(s) = K_{\Delta t}[f(t)] = \frac{1}{\sqrt{s^2+1}} = \sqrt{\frac{1}{s^2+1}} \qquad (5.6\text{-}153)$$

Expand F(s) into a $K_{\Delta t}$ Transform Asymptotic Series

$$F(s) = K_{\Delta t}[f(t)] = \sqrt{\frac{1}{s^2+1}} = 1s^{-1} - \frac{1}{2}s^{-3} + \frac{3}{8}s^{-5} - \frac{5}{16}s^{-7} + \dots \qquad (5.6\text{-}154)$$

The above asymptotic expansion of F(s), Eq 5.6-154, was obtained using the Internet site, wolframalpha.com.

The Inverse $K_{\Delta t}$ Transform Formula for a function, f(t), is:

$$K_{\Delta t}^{-1}[K_{\Delta t}[f(t)]] = f(n\Delta t) = \sum_{p=0}^{n}[{}_nC_p\Delta t^p]\, a_{p+1} \qquad (5.6\text{-}155)$$

$$K_{\Delta t}[f(t)] = F(s) = a_1s^{-1} + a_2s^{-2} + a_3s^{-3} + a_4s^{-4} + \dots \text{ , Asymptotic Series form of F(s)} \qquad (5.6\text{-}156)$$

where

$F(s) = K_{\Delta t}[f(t)] = K_{\Delta t}$ Transform of the function, f(t)

$K_{\Delta t}^{-1}[F(s)] = f(t) =$ Inverse $K_{\Delta t}$ Transform of the function, F(s)

n = 0, 1, 2, 3, 4, …

p = 0, 1, 2, 3, … , n

$nC_p = \dfrac{n!}{r!(n-r)!}$

0! = 1

f(t) = function of t

$\Delta t$ = interval between successive values of t

$t = 0, \Delta t, 2\Delta t, 3\Delta t, \ldots$

$D_{\Delta t}{}^{m} f(0) = a_{m+1}$ , m = 0,1,2,3,… , Initial conditions of f(t)

Comparing Eq 5.6-156 to Eq 5.6-154

$$a_1 = 1 \tag{5.6-157}$$

$$a_2 = 0 \tag{5.6-158}$$

$$a_3 = -\frac{1}{2} \tag{5.6-159}$$

$$a_4 = 0 \tag{5.6-160}$$

$$a_5 = \frac{3}{8} \tag{5.6-161}$$

$$\Delta t = .2 \tag{5.6-162}$$

From the Inverse $K_{\Delta t}$ Transform Formula, Eq 5.6-155

$$f(0) = a_1 \tag{5.6-163}$$

$$f(0) = 1 \tag{5.6-164}$$

$$f(\Delta t) = a_1 + a_2\Delta t \tag{5.6-165}$$

$$f(.2) = 1 + (0)(.2) = 1 \tag{5.6-166}$$

$$f(2\Delta t) = a_1 + a_2 2\Delta t + a_3\Delta t^2 \tag{5.6-167}$$

$$f(.4) = 1 + 0 + (-.5)(.04) = .98 \tag{5.6-168}$$

$$f(3\Delta t) = a_1 + a_2 3\Delta t + a_3 3\Delta t^2 + a_4\Delta t^3 \tag{5.6-169}$$

$$f(.6) = 1 + 0 + (-.5)3(.04) + 0 = .94 \tag{5.6-170}$$

Then

From Eq 5.6-164, Eq 5.6-166, Eq 5.6-168, and Eq 5.6-170

f(0) = 1
f(.2) = 1
f(.4) = .98
f(.6) = .94

Checking

In Example 5.6-6 this same problem was solved in a different way.
The above values are the same values calculated in Example 5.6-6.

f(0) = 1
f(.2) = 1
f(.4) = .98
f(.6) = .94

Good check

<u>Evaluating $f(t) = L^{-1}[F(s)]$ using the Laplace Transform Asymptotic Series and the Maclaurin Series</u>

As previously stated, the Laplace Transform is a $K_{\Delta t}$ Transform where $\Delta t$ is an infinitesimal value, $\Delta t \to 0$. Thus, the inverse transform methods previously derived and discussed can also be used for Laplace Transforms.
From Eq 5.6-59

$\Delta t \to 0$

**The Laplace Transform Asymptotic Series is:**

$$\mathbf{F(s) = L[f(t)] = \sum_{n=1}^{\infty} \frac{d^{n-1}}{dt^{n-1}} f(0)s^{-n} = f(0)s^{-1} + \frac{d}{dt} f(0)s^{-2} + \frac{d^2}{dt^2} f(0)s^{-3} + \frac{d^3}{dt^3} f(0)s^{-4} + \dots} \qquad (5.6\text{-}171)$$

**where**

$$\mathbf{\frac{d^m}{dt^m} f(0) = \frac{d^m}{dt^m} f(t)\Big|_{t=0}\ ,\quad m = 0,1,2,3,\dots\ ,\quad \text{Derivatives of } f(t) \text{ evaluated at } t = 0}$$

**f(t) = function of t**

**F(s) = L[f(t)] = Laplace Transform of the function, f(t)**

First, the Laplace Transform, F(s) =L[f(t)], is expanded into a Laplace Transform Asymptotic Series.

The resulting Laplace Transform asymptotic series is:

$$F(s) = L[f(t)] = \sum_{n=1}^{\infty} c_n s^{-n} = c_1 s^{-1} + c_2 s^{-2} + c_3 s^{-3} + c_4 s^{-4} + c_5 s^{-5} + \dots \quad (5.6\text{-}172)$$

$$\frac{d^{n-1}}{dt^{n-1}} f(t)|_{t=0} = c_n, \quad n = 1,2,3,\dots \quad (5.6\text{-}173)$$

where

$$\frac{d^0}{dt^0} f(t)|_{t=0} = f(0)$$

Finding f(t)

The evaluated values of the initial condition derivatives are entered into a Maclaurin Series.

$$f(t) = \sum_{m=0}^{\infty} \frac{1}{m!} \frac{d^m}{dt^m} f(t)|_{t=0} t^m = \sum_{m=0}^{\infty} \frac{1}{m!} f^{(m)}(0)\, t^m \quad (5.6\text{-}174)$$

Expanding Eq 5.6-174

$$f(t) = f(t) = f(0) + \frac{1}{1!} f^{(1)}(0)\, t + \frac{1}{2!} f^{(2)}(0)\, t^2 + \frac{1}{3!} f^{(3)}(0)\, t^3 + \frac{1}{4!} f^{(4)}(0)\, t^4 + \dots \quad (5.6\text{-}175)$$

The above Maclaurin Series is used to evaluate f(t).
The following example, Example 5.6-8, provides a demonstration of the use of a Laplace Transform asymptotic series expansion and a Maclaurin Series to find the inverse of a Laplace Transform.

<u>Example 5.6-8</u> Find the values of f(t) at t = 1,2,3,4,and 5 given the Laplace Transform,
$F(s) = L[f(t)] = \frac{1}{\sqrt{s^2+1}}$, the Laplace Transform of the $J_0(t)$ Bessel Function. Show the values of f(t) to eight decimal places. Use a Laplace Transform Asymptotic Series and a Maclaurin Series.

$$F(s) = L[f(t)] = \frac{1}{\sqrt{s^2+1}} = \sqrt{\frac{1}{s^2+1}} \quad (5.6\text{-}176)$$

Expanding F(s) into an asymptotic series

$$F(s) = L[f(t)] = \sqrt{\frac{1}{s^2+1}} = 1s^{-1} - \frac{1}{2} s^{-3} + \frac{3}{8} s^{-5} - \frac{5}{16} s^{-7} + \frac{35}{128} s^{-9} - \frac{63}{256} s^{-11} + \frac{231}{1024} s^{-13} - \frac{429}{2048} s^{-15} + \dots$$
$$+ \frac{6435}{32768} s^{-17} - \frac{12155}{65536} s^{-19} + \frac{46189}{262144} s^{-21} - \frac{88179}{524288} s^{-23} + \dots \quad (5.6\text{-}177)$$

The above expansion of F(s) into an asymptotic series was calculated using the Maxima symbolic logic computer program.

Writing the Laplace Transform Asymptotic Series

$$F(s) = L[f(t)] = \sum_{n=1}^{\infty} \frac{d^{n-1}}{dt^{n-1}} f(0)s^{-n} = f(0)s^{-1} + \frac{d}{dt} f(0)s^{-2} + \frac{d^2}{dt^2} f(0)s^{-3} + \frac{d^3}{dt^3} f(0)s^{-4} + \dots \qquad (5.6\text{-}178)$$

where

$\frac{d^m}{dt^m} f(0) = \frac{d^m}{dt^m} f(t)|_{t=0}$ , m = 0,1,2,3,… , Derivatives of f(t) evaluated at t = 0

f(t) = function of t

F(s) = L[f(t)] = Laplace Transform of the function, f(t)

Comparing the series for F(s), Eq 5.6-177, to the Laplace Transform Asymptotic Series, Eq 5.6-178, the derivatives of f(t) at t = 0 are evaluated:

$f(0) = 1 \qquad \frac{d^2}{dt^2} f(0) = -\frac{1}{2} \qquad \frac{d^4}{dt^4} f(0) = +\frac{3}{8} \qquad \frac{d^6}{dt^6} f(0) = -\frac{5}{16} \qquad \frac{d^8}{dt^8} f(0) = +\frac{35}{128}$

$\frac{d^{10}}{dt^{10}} f(0) = -\frac{63}{256} \qquad \frac{d^{12}}{dt^{12}} f(0) = +\frac{231}{1024} \qquad \frac{d^{14}}{dt^{14}} f(0) = -\frac{429}{2048} \qquad \frac{d^{16}}{dt^{16}} f(0) = +\frac{6435}{32768}$

$\frac{d^{18}}{dt^{18}} f(0) = -\frac{12155}{65536} \qquad \frac{d^{20}}{dt^{20}} f(0) = +\frac{46189}{262144} \qquad \frac{d^{22}}{dt^{22}} f(0) = -\frac{88179}{524288} \qquad \dots$

Placing the above initial condition derivatives of f(t) into a Maclaurin Series

$$f(t) = 1 - \frac{1}{2(2!)} t^2 + \frac{3}{8(4!)} t^4 - \frac{5}{16(6!)} t^6 + \frac{35}{128(8!)} t^8 - \frac{63}{256(10!)} t^{10} + \frac{231}{1024(12!)} t^{12} - \frac{429}{2048(14!)} t^{14}$$
$$+ \frac{6435}{32768(16!)} t^{16} - \frac{12155}{65536(18!)} t^{18} + \frac{46189}{262144(20!)} t^{20} - \frac{88179}{524288(22!)} t^{22} + \dots \qquad (5.6\text{-}179)$$

From Eq 5.6-179

f(1) = .76519768 , $f(1)_{exact}$ = .76519768 Good check

f(2) = .22389077 , $f(2)_{exact}$ = .22389077 Good check

f(3) = -.26005195 , $f(3)_{exact}$ = -.26005195 Good check
f(4) = -.39714980 , $f(4)_{exact}$ = -.39714980 Good check

f(5) = -.17759677 , $f(5)_{exact}$ = -.17759678 Good check

Evaluating $K_{\Delta t}$ Transforms using the $K_{\Delta t}$ Transform Asymptotic Series and Laplace Transforms using the Laplace Transform Asymptotic Series

The $K_{\Delta t}$ Transform Asymptotic Series is:

$$F(s) = K_{\Delta t}[f(t)] = \sum_{n=1}^{\infty} D_{\Delta t}^{n-1}f(0)s^{-n} = f(0)s^{-1} + D_{\Delta t}^{1}f(0)s^{-2} + D_{\Delta t}^{2}f(0)s^{-3} + D_{\Delta t}^{3}f(0)s^{-4} + \ldots \qquad (5.6\text{-}180)$$

where

$D_{\Delta t}^{m}f(0) = D_{\Delta t}^{m}f(t)|_{t=0}$ , m = 0,1,2,3,… , Discrete derivatives of f(t) evaluated at t = 0

f(t) = function of t

$F(s) = K_{\Delta t}[f(t)] = K_{\Delta t}$ Transform of the function, f(t)

The Laplace Transform Asymptotic Series is:

$$F(s) = K_{\Delta t}[f(t)] = \sum_{n=1}^{\infty} \frac{d^{n-1}}{dt^{n-1}} f(0)s^{-n} = f(0)s^{-1} + \frac{d}{dt} f(0)s^{-2} + \frac{d^2}{dt^2} f(0)s^{-3} + \frac{d^3}{dt^3} f(0)s^{-4} + \ldots \qquad (5.6\text{-}181)$$

where

$\frac{d^m}{dt^m} f(0) = \frac{d^m}{dt^m} f(t)|_{t=0}$ , m = 0,1,2,3,… , Derivatives of f(t) evaluated at t = 0

f(t) = function of t

F(s) = L[f(t)] = Laplace Transform of the function, f(t)

Note that the Laplace Transform Asyptotic Series is the $K_{\Delta t}$ Transform Asymptotic Series where $\Delta t$ is an infinitessimal increment of the t variable.

From a given function, f(t), its t = 0 discrete derivatives or derivatives (where $\Delta t$ is infinitesimal) can be determined. The derivative values obtained when inserted into the appropriate transform asymptotic series will determine the transform of f(t), $F(s) = K_{\Delta t}[f(t)]$ where $\Delta t$ is not infinitesimal and F(s) = L[f(t)] where $\Delta t$ is infinitesimal. The F(s) transform determined will be in an open ended series form. In many cases, the series obtained can be represented in a closed function form.

The following example, Example 5.6-9, provides a demonstration of the use of the $K_{\Delta t}$ Transform or Laplace Transform Asymptotic Series to obtain the transform of a given function of t, f(t).

Example 5.6-9 Given several functions of t, f(t), find the requested transforms. Use the appropriate $K_{\Delta t}$ Transform or Laplace Transform Asymptotic Series.

1) Find the $K_{\Delta t}$ Transform of the discrete function, $f(t) = t(t-\Delta t)(t-2\Delta t)$

Finding the discrete derivatives of f(t) evaluated at t = 0

$$D_{\Delta t}^{0}f(0) = f(0) = t(t-\Delta t)(t-2\Delta t)|_{t=0} = 0 \qquad (5.6\text{-}182)$$

$D_{\Delta t}{}^{1}f(0) = 3t(t\text{-}\Delta t)|_{t=0} = 0$ (5.6-183)

$D_{\Delta t}{}^{2}f(0) = 3(2)(t)|_{t=0} = 0$ (5.6-184)

$D_{\Delta t}{}^{3}f(0) = 3(2)(1) = 3! = 6$ (5.6-185)

Substituting the values, Eq 5.6-182 thru Eq 5.6-185, into The $K_{\Delta t}$ Transform Asymptotic Series, Eq 5.6-180

$$F(s) = K_{\Delta t}[f(t)] = 0s^{-1} + 0s^{-2} + 0s^{-3} + 3!s^{-4} \tag{5.6-186}$$

Then

$$F(s) = K_{\Delta t}[t(t\text{-}\Delta t)(t\text{-}2\Delta t)] = \frac{3!}{s^4}\ , \quad \text{This is the correct } K_{\Delta t} \text{ Transform} \tag{5.6-187}$$

2) Find the $K_{\Delta t}$ Transform of the discrete function, $f(t) = e_{\Delta t}(a,t)$

Finding the discrete derivatives of f(t) evaluated at t = 0

$D_{\Delta t}{}^{0}f(0) = f(0) = e_{\Delta t}(a,t)|_{t=0} = 1$ (5.6-188)

$D_{\Delta t}{}^{1}f(0) = ae_{\Delta t}(a,t)|_{t=0} = a$ (5.6-189)

$D_{\Delta t}{}^{2}f(0) = a^2e_{\Delta t}(a,t)|_{t=0} = a^2$ (5.6-190)

$D_{\Delta t}{}^{3}f(0) = a^3e_{\Delta t}(a,t)|_{t=0} = a^3$ (5.6-191)

.
.
.

$D_{\Delta t}{}^{n}f(0) = a^ne_{\Delta t}(a,t)|_{t=0} = a^n\ , \quad n = 0,1,2,3, \ldots$ (5.6-192)

Rewriting $K_{\Delta t}$ Transform Asymptotic Series

$$F(s) = K_{\Delta t}[f(t)] = \sum_{n=1}^{\infty} D_{\Delta t}{}^{n-1}f(0)s^{-n} = f(0)s^{-1} + D_{\Delta t}{}^{1}f(0)s^{-2} + D_{\Delta t}{}^{2}f(0)s^{-3} + D_{\Delta t}{}^{3}f(0)s^{-4} + \ldots \tag{5.6-193}$$

where

$D_{\Delta t}{}^{m}f(0) = D_{\Delta t}{}^{m}f(t)|_{t=0}$ , m = 0,1,2,3,… , Discrete derivatives of f(t) evaluated at t = 0
f(t) = function of t
$F(s) = K_{\Delta t}[f(t)] = K_{\Delta t}$ Transform of the function, f(t)

Substituting the values of Eq 5.6-192 into The $K_{\Delta t}$ Transform Asymptotic Series, Eq 5.6-193

$$F(s) = K_{\Delta t}[f(t)] = 1s^{-1} + as^{-2} + a^2s^{-3} + a^3s^{-4} + a^4s^{-5} + \ldots \tag{5.6-194}$$

$$\frac{s^{-1}}{1 - as^{-1}} = \frac{\frac{1}{s}}{1 - \frac{a}{s}} = \frac{1}{s-a} = 1s^{-1} + as^{-2} + a^2s^{-3} + a^3s^{-4} + a^4s^{-5} + \ldots \tag{5.6-195}$$

Substituting Eq 5.6-195 into Eq 5.6-194

$$F(s) = K_{\Delta t}[f(t)] = \frac{1}{s-a} \tag{5.6-196}$$

Then

$$F(s) = K_{\Delta t}[e_{\Delta t}(a,t)] = \frac{1}{s-a} \ , \quad \text{This is the correct } K_{\Delta t} \text{ Transform} \tag{5.6-197}$$

Comment - Since $K_{\Delta t}$ Transforms become Laplace Transforms and the discrete functions become their related Calculus functions when $\Delta t$ is infinitessimal, $L[e^{at}] = \frac{1}{s-a}$ .

3) Find the LaplaceTransform of the function, f(t) = sinat

Finding the derivatives of f(t) evaluated at t = 0

$$\frac{d^0}{dt^0} f(0) = f(0) = \sin at \,|_{t=0} = 0 \tag{5.6-198}$$

$$\frac{d^1}{dt^1} f(0) = a\cos at|_{t=0} = a \tag{5.6-199}$$

$$\frac{d^2}{dt^2} f(0) = -a^2 \sin at|_{t=0} = 0 \tag{5.6-200}$$

$$\frac{d^3}{dt^3} f(0) = -a^3 \cos at|_{t=0} = -a^3 \tag{5.6-201}$$

.
.
.

$$\frac{d^n}{dt^n} f(0) = 0 \ , \qquad n = 0,2,4,6, \dots \tag{5.6-202}$$

$$\frac{d^n}{dt^n} f(0) = (-1)^{\frac{n-1}{2}} a^n \ , \quad n = 1,3,5,7, \dots \tag{5.6-203}$$

Rewriting the Laplace Transform Asymptotic Series

$$F(s) = K_{\Delta t}[f(t)] = \sum_{n=1}^{\infty} \frac{d^{n-1}}{dt^{n-1}} f(0)s^{-n} = f(0)s^{-1} + \frac{d}{dt} f(0)s^{-2} + \frac{d^2}{dt^2} f(0)s^{-3} + \frac{d^3}{dt^3} f(0)s^{-4} + \dots \tag{5.6-204}$$

where

$\frac{d^m}{dt^m} f(0) = \frac{d^m}{dt^m} f(t)|_{t=0}$ , m = 0,1,2,3,… , Derivatives of f(t) evaluated at t = 0

f(t) = function of t

$F(s) = L[f(t)]$ = Laplace Transform of the function, $f(t)$

Substituting the values of Eq 5.6-202 and Eq 5.6-203 into the Laplace Transform Asymptotic Series, Eq 5.6-204

$$F(s) = L[f(t)] = 0s^{-1} + as^{-2} + 0s^{-3} - a^3s^{-4} + 0s^{-5} + a^5s^{-6} + 0s^{-7} - a^7s^{-8} + \dots \tag{5.6-205}$$

$$F(s) = L[f(t)] = as^{-2} - a^3s^{-4} + a^5s^{-6} - a^7s^{-8} + \dots \tag{5.6-206}$$

$$\frac{as^{-2}}{1 + a^2s^{-2}} = \frac{\frac{a}{s^2}}{1 + \frac{a^2}{s^2}} = \frac{a}{s^2 + a^2} = as^{-2} - a^3s^{-4} + a^5s^{-6} - a^7s^{-8} + \dots \tag{5.6-207}$$

Substituting Eq 5.6-207 into Eq 5.6-206

$$F(s) = L[f(t)] = \frac{a}{s^2 + a^2} \tag{5.6-208}$$

Then

$$F(s) = L[\sin at] = \frac{a}{s^2 + a^2}\ , \quad \text{This is the correct Laplace Transform} \tag{5.6-209}$$

**Section 5.7: The use of Calculus functions in Interval Calculus**

Previously, Interval Calculus discrete calculations have used only discrete Interval Calculus functions such as $e_{\Delta t}(a,t)$, $\sin_{\Delta t}(b,t)$, and $\cos_{\Delta t}(b,t)$. However, Calculus funcions such as $e^{at}$, sinbt, and cosbt can also be used. In Calculus, the Calculus functions are defined as functions of a continuous variable, t. In Interval Calculus, any Calculus functions used are defined only at equally spaced discrete values of the variable, t, for example at $t = 0, \Delta t, 2\Delta t, 3\Delta t, \dots$

Calculus functions and Interval Calculus functions are related. Consider the following derivation.

Show that the Interval Calculus function, $e_{\Delta t}(a,t)$, and the Calculus function, $e^{\alpha t}$, are related. $\alpha$ and a are constants.

$$e_{\Delta t}(\alpha,t) = (1+\alpha\Delta t)^{\frac{t}{\Delta t}} = e^{\frac{\ln(1+\alpha\Delta t)}{\Delta t}t} = e^{at} \tag{5.7-1}$$

Then

$$\mathbf{e_{\Delta t}(\alpha,t) = e^{at}} \tag{5.7-2}$$

**where**

**$\alpha$, a = constants**

$$\mathbf{a = \frac{\ln(1+\alpha\Delta t)}{\Delta t}}$$

**$\Delta t$ = sampling interval, t increment**

**$t = 0, \Delta t, 2\Delta t, 3\Delta t, \dots$**

Also

From Eq 5.7-2

$$\alpha = \frac{e^{a\Delta t} - 1}{\Delta t} \qquad (5.7\text{-}3)$$

From Eq 5.7-2 and Eq 5.7-3

$$e^{at} = e_{\Delta t}(\alpha,t) = [1+\alpha\Delta t]^{\frac{t}{\Delta t}} = [1 + (\frac{e^{a\Delta t}-1}{\Delta t})\Delta t]^{\frac{t}{\Delta t}} = e_{\Delta t}(\frac{e^{a\Delta t}-1}{\Delta t}, t) \qquad (5.7\text{-}4)$$

Then

$$\mathbf{e^{at} = e_{\Delta t}(\alpha, t)} \qquad (5.7\text{-}5)$$

**where**

**$\alpha$, a = constants**

$$\boldsymbol{\alpha = \frac{e^{a\Delta t} - 1}{\Delta t}}$$

**$\Delta t$ = sampling interval, t increment**

**t = 0, $\Delta t$, $2\Delta t$, $3\Delta t$, …**

Both Eq 5.7-2 and Eq 5.7-5 show two relationships between the Interval Calculus function, $e_{\Delta t}(a,t)$, and the Calculus function, $e^{\alpha t}$. In fact, these two relationships are equalities that in Interval Calculus are referred to as identities. Identities are defined as pairs of functions, one an Interval Calculus function and the other a Calculus function, that yield equal results for all values of their common discrete independent variable. Since each function of a pair can be expressed in terms of an expression containing the other, for each pair of functions there are two identities. A Table of several Interval Calculus/Calculus function identities is provided in Table 5.7-1 which appears below. Two examples involving identities, Example 5.7-1 and Example 5.7-2, are presented following Table 5.7-1.

Note the derived identities, Eq 5.7-2 and Eq 5.7-5, in rows 1 and 8 of Table 5.7-1.

Comment – There is another relationship between Interval Calculus functions and Calculus functions. If the $\Delta t$ value of an Interval Calculus function becomes infinitessimal (i.e. $\Delta t \to 0$), a related Calculus function results. For example, $\lim_{\Delta t\to 0}e_{\Delta t}(a,t) = e^{at}$, $\lim_{\Delta t\to 0}\sin_{\Delta t}(b,t) = \sin bt$, $\lim_{\Delta t\to 0}\cos_{\Delta t}(b,t) = \cos bt$, etc.

## Table 5.7-1 Interval Calculus Function/Calculus Function Identities

| # | Calculus Function<br>$t = 0, \Delta t, 2\Delta t, 3\Delta t, \ldots$ | Interval Calculus Function<br>$t = 0, \Delta t, 2\Delta t, 3\Delta t, \ldots$ |
|---|---|---|
| **1** | $e^{at}$ | $e_{\Delta t}(\frac{e^{a\Delta t}-1}{\Delta t}, t)$<br>or<br>$[1 + (\frac{e^{a\Delta t}-1}{\Delta t})\Delta t]^{\frac{t}{\Delta t}}$ |
| **2** | $\sin bt$ | $[\cos b\Delta t]^{\frac{t}{\Delta t}}\sin_{\Delta t}(\frac{\tan b\Delta t}{\Delta t}, t)$<br>or<br>$e_{\Delta t}(\frac{\cos b\Delta t - 1}{\Delta t}, t)\sin_{\Delta t}(\frac{\tan b\Delta t}{\Delta t}, t)$<br>$\cos b\Delta t \neq 0$ |
| **3** | $\cos bt$ | $[\cos b\Delta t]^{\frac{t}{\Delta t}}\cos_{\Delta t}(\frac{\tan b\Delta t}{\Delta t}, t)$<br>or<br>$e_{\Delta t}(\frac{\cos b\Delta t - 1}{\Delta t}, t)\cos_{\Delta t}(\frac{\tan b\Delta t}{\Delta t}, t)$<br>$\cos b\Delta t \neq 0$ |
| **4** | $e^{at}\sin bt$ | $[e^{a\Delta t}\cos b\Delta t]^{\frac{t}{\Delta t}}\sin_{\Delta t}(\frac{\tan b\Delta t}{\Delta t}, t)$<br>or<br>$e_{\Delta t}(\frac{e^{a\Delta t}\cos b\Delta t - 1}{\Delta t}, t)\sin_{\Delta t}(\frac{\tan b\Delta t}{\Delta t}, t)$<br>$\cos b\Delta t \neq 0$ |
| **5** | $e^{at}\cos bt$ | $[e^{a\Delta t}\cos b\Delta t]^{\frac{t}{\Delta t}}\cos_{\Delta t}(\frac{\tan b\Delta t}{\Delta t}, t)$<br>or<br>$e_{\Delta t}(\frac{e^{a\Delta t}\cos b\Delta t - 1}{\Delta t}, t)\cos_{\Delta t}(\frac{\tan b\Delta t}{\Delta t}, t)$<br>$\cos b\Delta t \neq 0$ |
| **6** | $\sinh bt$ | $[\cosh b\Delta t]^{\frac{t}{\Delta t}}\sinh_{\Delta t}(\frac{\tanh b\Delta t}{\Delta t}, t)$<br>or<br>$e_{\Delta t}(\frac{\cosh b\Delta t - 1}{\Delta t}, t)\sinh_{\Delta x}(\frac{\tanh b\Delta t}{\Delta t}, t)$ |

| # | Calculus Function<br>t = 0, Δt, 2Δt, 3Δt, … | Interval Calculus Function<br>t = 0, Δt, 2Δt, 3Δt, … |
|---|---|---|
| | | |
| 7 | $\cosh bt$ | $[\cosh(b\Delta t)]^{\frac{t}{\Delta t}}\cosh_{\Delta t}(\frac{\tanh b\Delta t}{\Delta t}, t)$<br>or<br>$e_{\Delta t}(\frac{\cosh b\Delta t - 1}{\Delta t}, t)\cosh_{\Delta t}(\frac{\tanh b\Delta t}{\Delta t}, t)$ |
| 8 | $e^{at}$<br>$a = \frac{\ln(1+\alpha\Delta t)}{\Delta t}$ | $e_{\Delta t}(\alpha,t)$ |
| 9 | $(\sec\beta)^{\frac{t}{\Delta t}}\sin\beta\frac{t}{\Delta t}$<br>$\beta = \tan^{-1}b\Delta t$ | $\sin_{\Delta t}(b,t)$ |
| 10 | $(\sec\beta)^{\frac{t}{\Delta t}}\cos\beta\frac{t}{\Delta t}$<br>$\beta = \tan^{-1}b\Delta t$ | $\cos_{\Delta t}(b,t)$ |
| 11 | $[(1+a\Delta t)^2+(b\Delta t)^2]^{\frac{t}{2\Delta t}}\sin\beta\frac{t}{\Delta t}$<br>or<br>$e_{\Delta t}(a,t)(\sec\beta)^{\frac{t}{\Delta t}}\sin\beta\frac{t}{\Delta t}$<br>$\beta = \begin{cases} \tan^{-1}\left\lvert\frac{b\Delta t}{1+a\Delta t}\right\rvert & \text{for } 1+a\Delta t > 0 \quad b\Delta t \geq 0 \\ -\tan^{-1}\left\lvert\frac{b\Delta t}{1+a\Delta t}\right\rvert & \text{for } 1+a\Delta t > 0 \quad b\Delta t < 0 \\ \pi - \tan^{-1}\left\lvert\frac{b\Delta t}{1+a\Delta t}\right\rvert & \text{for } 1+a\Delta t < 0 \quad b\Delta t \geq 0 \\ -\pi + \tan^{-1}\left\lvert\frac{b\Delta t}{1+a\Delta t}\right\rvert & \text{for } 1+a\Delta t < 0 \quad b\Delta t < 0 \end{cases}$<br>$1+a\Delta t \neq 0$<br>$0 \leq \tan^{-1}\left\lvert\frac{b\Delta t}{1+a\Delta t}\right\rvert < \frac{\pi}{2}$ | $e_{\Delta t}(a,t)\sin_{\Delta t}(\frac{b}{1+a\Delta t}, t)$ |
| 12 | $[(1+a\Delta t)^2+(b\Delta t)^2]^{\frac{t}{2\Delta t}}\cos\beta\frac{t}{\Delta t}$<br>or | $e_{\Delta t}(a,t)\cos_{\Delta t}(\frac{b}{1+a\Delta t}, t)$ |

| # | Calculus Function<br>t = 0, Δt, 2Δt, 3Δt, … | Interval Calculus Function<br>t = 0, Δt, 2Δt, 3Δt, … |
|---|---|---|
| | $e_{\Delta t}(a,t)(\sec\beta)^{\frac{t}{\Delta t}}\cos\beta\frac{t}{\Delta t}$<br>$\beta = \begin{cases} \tan^{-1}\left\lvert\frac{b\Delta t}{1+a\Delta t}\right\rvert & \text{for } 1+a\Delta t > 0 \quad b\Delta t \geq 0 \\ -\tan^{-1}\left\lvert\frac{b\Delta t}{1+a\Delta t}\right\rvert & \text{for } 1+a\Delta t > 0 \quad b\Delta t < 0 \\ \pi - \tan^{-1}\left\lvert\frac{b\Delta t}{1+a\Delta t}\right\rvert & \text{for } 1+a\Delta t < 0 \quad b\Delta t \geq 0 \\ -\pi + \tan^{-1}\left\lvert\frac{b\Delta t}{1+a\Delta t}\right\rvert & \text{for } 1+a\Delta t < 0 \quad b\Delta t < 0 \end{cases}$<br>$1+a\Delta t \neq 0$<br>$0 \leq \tan^{-1}\left\lvert\frac{b\Delta t}{1+a\Delta t}\right\rvert < \frac{\pi}{2}$ | |
| **13** | $(1-[b\Delta t]^2)^{\frac{x}{2\Delta t}}\sinh\beta\frac{t}{\Delta t}$<br>$\beta = \tanh^{-1}b\Delta t$ | $\sinh_{\Delta t}(b,t)$ |
| **14** | $(1-[b\Delta t]^2)^{\frac{t}{\Delta t}}\cosh\beta\frac{t}{\Delta t}$<br>$\beta = \tanh^{-1}b\Delta t$ | $\cosh_{\Delta t}(b,t)$ |

Identities such as those shown in Table 5.1 above are used in a number of equation derivations.

Example 5.7-1 presents the derivation of an identity involving the discrete Interval Calculus function, $e_{\Delta t}(a,t)\sin_{\Delta t}(\frac{b}{1+a\Delta t}, t)$.

<u>Example 5.7-1</u> Find an identity involving the Interval Calculus function, $e_{\Delta t}(a,t)\sin_{\Delta t}(\frac{b}{1+a\Delta t}, t)$ .

$$e_{\Delta t}(a,t)\sin_{\Delta t}(\frac{b}{1+a\Delta t}, t) = (1+a\Delta t)^{\frac{t}{\Delta t}}\sin_{\Delta t}(\frac{b\Delta t}{1+a\Delta t}, t) = \frac{(1+a\Delta t + jb\Delta t)^{\frac{t}{\Delta t}} - (1+a\Delta t - jb\Delta t)^{\frac{t}{\Delta t}}}{2j} \qquad 1)$$

$$e_{\Delta t}(a,t)\sin_{\Delta t}(\frac{b}{1+a\Delta t}, t) = (1+a\Delta t)^{\frac{t}{\Delta t}}\left[\frac{(1 + \frac{jb\Delta t}{1+a\Delta t})^{\frac{t}{\Delta t}} - (1 - \frac{jb\Delta t}{1+a\Delta t})^{\frac{t}{\Delta t}}}{2j}\right] \qquad 2)$$

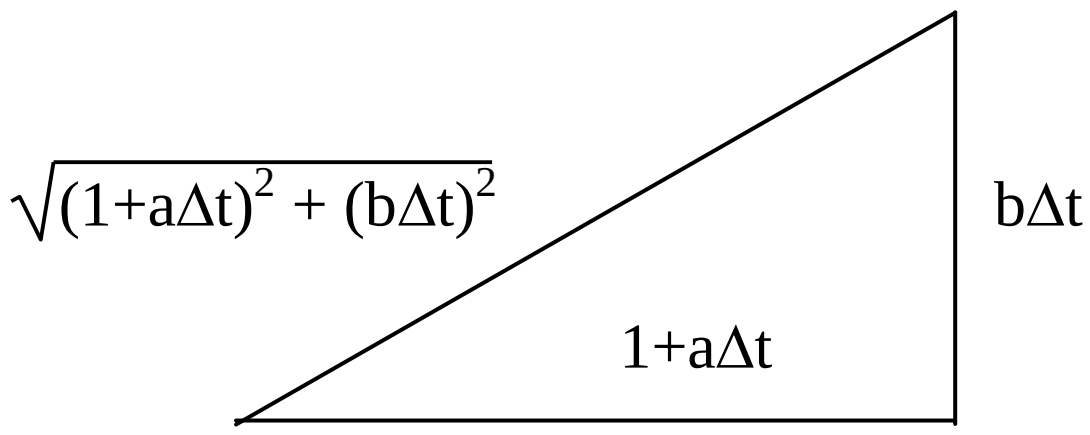

From the above diagram

$$e^{j\beta} = \cos\beta + j\sin\beta = \frac{1+a\Delta t}{\sqrt{(1+a\Delta t)^2 + (b\Delta t)^2}} + j\,\frac{b\Delta t}{\sqrt{(1+a\Delta t)^2 + (b\Delta t)^2}} \qquad 3)$$

$$e^{-j\beta} = \cos\beta + j\sin\beta = \frac{1+a\Delta t}{\sqrt{(1+a\Delta t)^2 + (b\Delta t)^2}} - j\,\frac{b\Delta t}{\sqrt{(1+a\Delta t)^2 + (b\Delta t)^2}} \qquad 4)$$

$$\tan\beta = \frac{b\Delta t}{1+a\Delta t} \qquad 5)$$

From Eq 1 thru Eq 4

$$e_{\Delta t}(a,t)\sin_{\Delta t}\left(\frac{b}{1+a\Delta t}, t\right) = \frac{1}{2j}\left[\sqrt{(1+a\Delta t)^2 + (b\Delta t)^2}\right]^{\frac{t}{\Delta t}}\left[\frac{1+a\Delta t}{\sqrt{(1+a\Delta t)^2 + (b\Delta t)^2}} + j\,\frac{b\Delta t}{\sqrt{(1+a\Delta t)^2 + (b\Delta t)^2}}\right]^{\frac{t}{\Delta t}} -$$

$$\frac{1}{2j}\left[\sqrt{(1+a\Delta t)^2 + (b\Delta t)^2}\right]^{\frac{t}{\Delta t}}\left[\frac{1+a\Delta t}{\sqrt{(1+a\Delta t)^2 + (b\Delta t)^2}} - j\,\frac{b\Delta t}{\sqrt{(1+a\Delta t)^2 + (b\Delta t)^2}}\right]^{\frac{t}{\Delta t}} \qquad 6)$$

$$e_{\Delta t}(a,t)\sin_{\Delta t}\left(\frac{b}{1+a\Delta t}, t\right) = \left[\sqrt{(1+a\Delta t)^2 + (b\Delta t)^2}\right]^{\frac{t}{\Delta t}}\left[\frac{e^{\frac{j\beta t}{\Delta t}} - e^{\frac{-j\beta t}{\Delta t}}}{2j}\right] \qquad 7)$$

$$e_{\Delta t}(a,t)\sin_{\Delta t}\left(\frac{b}{1+a\Delta t}, t\right) = \left[\sqrt{(1+a\Delta t)^2 + (b\Delta t)^2}\right]^{\frac{t}{\Delta t}} \sin\frac{\beta t}{\Delta t} \qquad 8)$$

where

$$\beta = \tan^{-1}\frac{b\Delta t}{1+a\Delta t}$$

The value of β will depend on the sign of (1+aΔt) and bΔt.

Four Quadrants

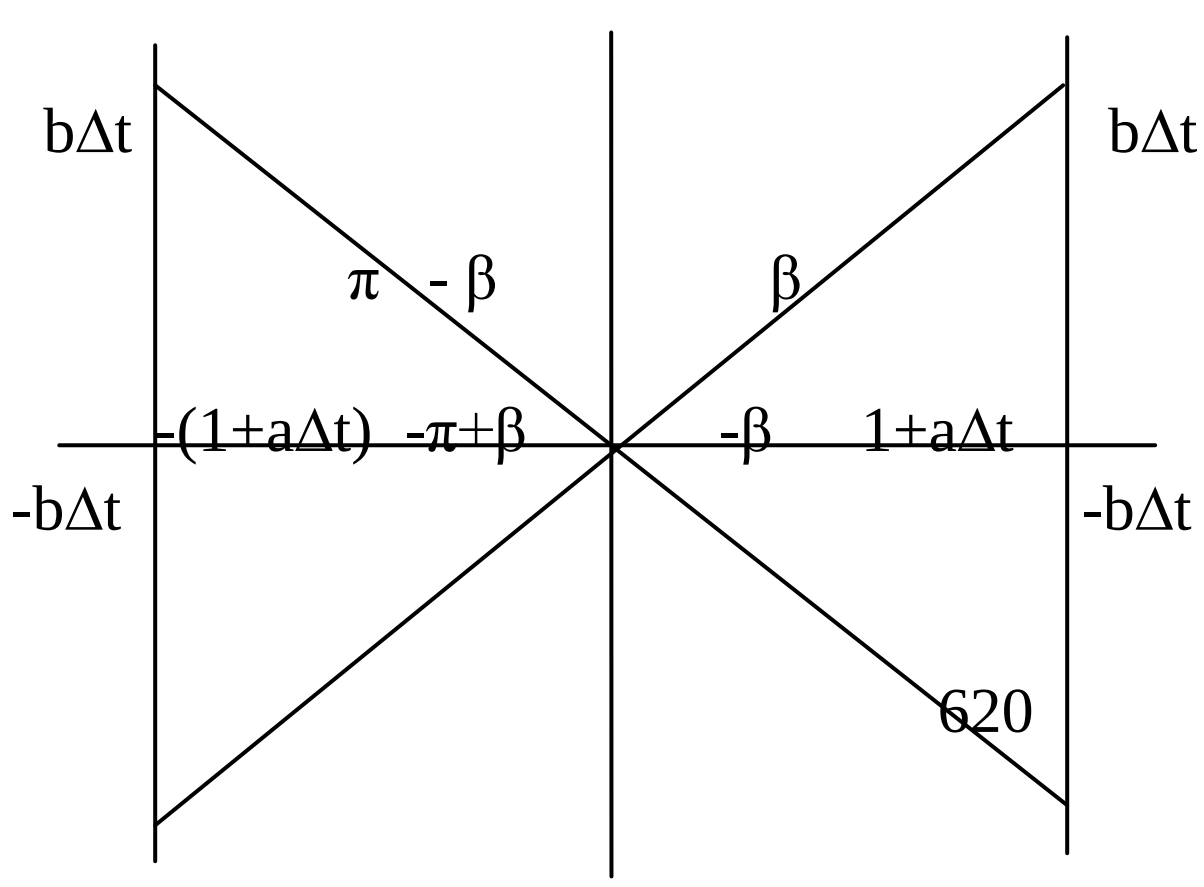

or

$$\beta = \begin{cases} \tan^{-1}\left|\frac{b\Delta t}{1+a\Delta t}\right| & \text{for } 1+a\Delta t > 0 \quad b\Delta t \geq 0 \\ -\tan^{-1}\left|\frac{b\Delta t}{1+a\Delta t}\right| & \text{for } 1+a\Delta t > 0 \quad b\Delta t < 0 \\ \pi - \tan^{-1}\left|\frac{b\Delta t}{1+a\Delta t}\right| & \text{for } 1+a\Delta t < 0 \quad b\Delta t \geq 0 \\ -\pi + \tan^{-1}\left|\frac{b\Delta t}{1+a\Delta t}\right| & \text{for } 1+a\Delta t < 0 \quad b\Delta t < 0 \end{cases}$$

$1+a\Delta t \neq 0$

$0 \leq \tan^{-1}\left|\frac{b\Delta t}{1+a\Delta t}\right| < \frac{\pi}{2}$

Then

**An identity for the discrete Interval Calculus function, $e_{\Delta t}(a,t)\sin_{\Delta t}(\frac{b}{1+a\Delta t}, t)$, is:**

$$\mathbf{e_{\Delta t}(a,t)\sin_{\Delta t}(\frac{b}{1+a\Delta t}, t) = \left[\sqrt{(1+a\Delta t)^2 + (b\Delta t)^2}\,\right]^{\frac{t}{\Delta t}} \sin\frac{\beta t}{\Delta t}} \qquad 9)$$

**where**

$$\boldsymbol{\beta} = \begin{cases} \tan^{-1}\left|\frac{b\Delta t}{1+a\Delta t}\right| & \text{for } 1+a\Delta t > 0 \quad b\Delta t \geq 0 \\ -\tan^{-1}\left|\frac{b\Delta t}{1+a\Delta t}\right| & \text{for } 1+a\Delta t > 0 \quad b\Delta t < 0 \\ \pi - \tan^{-1}\left|\frac{b\Delta t}{1+a\Delta t}\right| & \text{for } 1+a\Delta t < 0 \quad b\Delta t \geq 0 \\ -\pi + \tan^{-1}\left|\frac{b\Delta t}{1+a\Delta t}\right| & \text{for } 1+a\Delta t < 0 \quad b\Delta t < 0 \end{cases}$$

**$1+a\Delta t \neq 0$**

**$0 \leq \tan^{-1}\left|\frac{b\Delta t}{1+a\Delta t}\right| < \frac{\pi}{2}$**

**$\Delta t$ = sampling interval, t increment**

**$t = 0, \Delta t, 2\Delta t, 3\Delta t, \ldots$**

Example 5.7-2 is an example of the use of an Interval Calculus identity.

<u>Example 5.7-2</u> Find the inverse discrete derivative (i.e. the discrete integral) of $e^{2t}$ where $\Delta t = .5$.

Using an identity for $e^{2t}$ where $\Delta t = .5$

$$e^{2t} = (1 + [\frac{e^{2\Delta t}-1}{\Delta t}]\Delta t)^{\frac{t}{\Delta t}} = (1 + [\frac{e^{1}-1}{.5}]\Delta t)^{2t} = (1 + [3.436563657][.5])^{2t} \qquad 1)$$

$$_{\Delta t}\int (1 + a\Delta t)^{\frac{t}{\Delta t}}\Delta t = \frac{(1 + a\Delta t)^{\frac{t}{\Delta t}}}{a} + k \qquad 2)$$

where
k = constant of integration

From Eq 1 and Eq 2

$a = 3.43656365$
$\Delta t = .5$

$$_{.5}\int (2.718281828)^{2t}\Delta t = \frac{(2.718281828)^{2t}}{3.43656365} + k = .29098835e^{2t} + k \qquad 3)$$

Then

The inverse discrete derivative, the discrete integral, of $e^{2t}$ where $\Delta t = .5$ is:

$$_{.5}\int e^{2t}\Delta t = .29098835e^{2t} + k \qquad 4)$$

Checking

Assume the solution for the inverse discrete derivative of $e^{2t}$ is $be^{2t} + k$
where $\Delta t = .5$ and b,k = constants.

Taking the discrete derivative of $be^{2t}$ where $\Delta t = .5$

$$D_{.5}[be^{2t} + k] = \frac{be^{2(t+.5)} + k - be^{2t} - k}{.5} = \frac{(be - b)e^{2t}}{.5} = b\frac{(e-1)}{.5}e^{2t} = b(3.436563657)e^{2t}$$

Note that k can be any constant value

$b(3.436563657) = 1$

$b = .29098835$

$be^{2t} = .29098835e^{2t}$

Then the inverse discrete derivative of $e^{2t}$ is $.29098835e^{2t} + k$.

Good check

Previously, the use of $K_{\Delta t}$ Transforms was associated primarily with Interval Calculus functions such as $e_{\Delta t}(a,x)$, $\sin_{\Delta t}(b,t)$, and $\cos_{\Delta t}(b,t)$. The Calculus functions, such as $e^{at}$, sinbt, and cosbt, were associated with the Laplace Transform or, equivalently, the $K_{\Delta t}$ Transform where $\Delta t \rightarrow 0$. It will now be shown that Calculus functions can also be used with $K_{\Delta t}$ Transforms where $\Delta t$ is a finite interval and $t = 0, \Delta t, 2\Delta t, 3\Delta t, \ldots$ This would be similar to Z Transforms where $\Delta t$ (often referred to as T) is a finite interval and the associated sampled time domain functions are Calculus functions.

Consider the following derivation.

Derive the $K_{\Delta t}$ Transform for $e^{at}$ where $t = 0, \Delta t, 2\Delta t, 3\Delta t, \ldots$

The discrete Interval Calculus function, $e_{\Delta t}(\alpha,t)$, has a $K_{\Delta t}$ Transform which is stated in Eq 5.7-6 below.

$$K_{\Delta t}[e_{\Delta t}(\alpha,t)] = \frac{1}{s+\alpha} \qquad (5.7\text{-}6)$$

From Eq 5.7-5 and Eq 5.7-6

$$K_{\Delta t}[e^{at}] = K_{\Delta t}[e_{\Delta t}(\alpha,t)] = \frac{1}{s-\alpha} \quad \text{where } \alpha = \frac{e^{a\Delta t} - 1}{\Delta t} \qquad (5.7\text{-}7)$$

$$K_{\Delta t}[e^{at}] = \frac{1}{s - \frac{e^{a\Delta t} - 1}{\Delta t}} \qquad (5.7\text{-}8)$$

Then

**The $K_{\Delta t}$ Transform of $e^{at}$ is:**

$$\mathbf{K_{\Delta t}[e^{at}] = \frac{1}{s - \frac{e^{a\Delta t} - 1}{\Delta t}}} \qquad (5.7\text{-}9)$$

**where**

**a = constant**

**$\Delta t$ = sampling interval, t increment**

**$t = 0, \Delta t, 2\Delta t, 3\Delta t, \ldots$**

The above derived equation relates the Calculus exponential function, $e^{at}$, to its $K_{\Delta t}$ Transform. Note that this equation, Eq 5.7-9, is found in the following table, Table 5.7-2, in row 2.

Derive the inverse $K_{\Delta t}$ Transform of $\frac{1}{s-a}$ in terms of the exponential function $e^{\alpha t}$

From Eq 5.7-7 redefining the variable, $\alpha$, as a and the variable, a, as $\alpha$

$$K_{\Delta t}^{-1}\left[\frac{1}{s-a}\right] = e^{\alpha t} \qquad (5.7\text{-}10)$$

$$a = \frac{e^{\alpha\Delta t} - 1}{\Delta t} \tag{5.7-11}$$

Solving for $\alpha$

$$\alpha = \frac{\ln(1+a\Delta t)}{\Delta t}$$

Then

**The Inverse $K_{\Delta t}$ Transform of $\frac{1}{s-a}$ is:**

$$\mathbf{K_{\Delta t}^{-1}\left[\frac{1}{s-a}\right] = e^{\alpha t}} \tag{5.7-12}$$

**where**

**$\alpha$, a = constants**

$$\boldsymbol{\alpha = \frac{\ln(1+a\Delta t)}{\Delta t}}$$

**$\Delta t$ = sampling interval, t increment**

**t = 0, $\Delta t$, $2\Delta t$, $3\Delta t$, …**

The above derived equation relates the $K_{\Delta t}$ Transform, $\frac{1}{s-a}$, to its Calculus exponential function, $e^{\alpha t}$.

Note that this equation, Eq 5.7-12, is found in the following table, Table 5.7-3, in row 2.

In a similar manner, equations relating the Calculus functions, sinbx, cosbx, etc. to their $K_{\Delta t}$ Transforms can also be found. These relationships are also tabulated in Table 5.7-2 and Table 5.7-3.

The user of $K_{\Delta t}$ Transforms has the option of using either Interval Calculus functions or Calculus functions in his/her mathematical analysis. Either set of functions will provide the correct results.

See Example 5.7-3.

**Table 5.7-2 Conversions of Calculus Function Laplace Transforms to Equivalent Function Discrete $K_{\Delta t}$ Transforms**

| # | $K_{\Delta t}$ Transform $\Delta t \to 0$ Laplace Transform F(t) where $0 \le t < \infty$ | $K_{\Delta t}$ Transform Generalized Laplace Transform F(t) where $t = 0, \Delta t, 2\Delta t, 3\Delta t, \ldots$ | F(t) Calculus Function |
|---|---|---|---|
| **1** | $\frac{1}{s}$ | $\frac{1}{s}$ | 1 |
| **2** | $\frac{1}{s-a}$ | $\frac{1}{s - \frac{e^{a\Delta t}-1}{\Delta t}}$ | $e^{at}$ |
| **3** | $\frac{b}{s^2+b^2}$ | $\frac{\frac{\sin b\Delta t}{\Delta t}}{(s + \frac{1-\cos b\Delta t}{\Delta t})^2 + (\frac{\sin b\Delta t}{\Delta t})^2}$ | sinbt |
| **4** | $\frac{b}{s^2-b^2}$ | $\frac{\frac{\sinh b\Delta t}{\Delta t}}{(s + \frac{1-\cosh b\Delta t}{\Delta t})^2 - (\frac{\sinh b\Delta t}{\Delta t})^2}$ | sinhbt |
| **5** | $\frac{s}{s^2+b^2}$ | $\frac{s + \frac{1-\cos b\Delta t}{\Delta t}}{(s + \frac{1-\cos b\Delta t}{\Delta t})^2 + (\frac{\sin b\Delta t}{\Delta t})^2}$ | cosbt |
| **6** | $\frac{s}{s^2-b^2}$ | $\frac{s + \frac{1-\cosh b\Delta t}{\Delta t}}{(s + \frac{1-\cosh b\Delta t}{\Delta t})^2 - (\frac{\sinh b\Delta t}{\Delta t})^2}$ | coshbt |

| # | $K_{\Delta t}$ Transform $\Delta t \to 0$ Laplace Transform F(t) where $0 \le t < \infty$ | $K_{\Delta t}$ Transform Generalized Laplace Transform F(t) where $t = 0, \Delta t, 2\Delta t, 3\Delta t, \ldots$ | F(t) Calculus Function |
|---|---|---|---|
| **7** | $\frac{b}{(s-a)^2 + b^2}$ | $\frac{\frac{e^{a\Delta t}\sin b\Delta t}{\Delta t}}{(s - \frac{e^{a\Delta t}\cos b\Delta t - 1}{\Delta t})^2 + (\frac{e^{a\Delta t}\sin b\Delta t}{\Delta t})^2}$ | $e^{at}\sin bt$ |
| **8** | $\frac{b}{(s-a)^2 - b^2}$ | $\frac{\frac{e^{a\Delta t}\sinh b\Delta t}{\Delta t}}{(s - \frac{e^{a\Delta t}\cosh b\Delta t - 1}{\Delta t})^2 - (\frac{e^{a\Delta t}\sinh b\Delta t}{\Delta t})^2}$ | $e^{at}\sinh bt$ |
| **9** | $\frac{s-a}{(s-a)^2 + b^2}$ | $\frac{s - \frac{e^{a\Delta t}\cos b\Delta t - 1}{\Delta t}}{(s - \frac{e^{a\Delta t}\cos b\Delta t - 1}{\Delta t})^2 + (\frac{e^{a\Delta t}\sin b\Delta t}{\Delta t})^2}$ | $e^{at}\cos bt$ |
| **10** | $\frac{s-a}{(s-a)^2 - b^2}$ | $\frac{s - \frac{e^{a\Delta t}\cosh b\Delta t - 1}{\Delta t}}{(s - \frac{e^{a\Delta t}\cosh b\Delta t - 1}{\Delta t})^2 - (\frac{e^{a\Delta t}\sinh b\Delta t}{\Delta t})^2}$ | $e^{at}\cosh bt$ |
| **11** | $e^{-sn\Delta t}f(s)$ | $(1+s\Delta t)^{-n}f(s)$ | $U(t-n\Delta t)F(t-n\Delta t)$ $n = 0,1,2,3,\ldots$ |
| **12** | $\frac{1}{s^2}$ | $\frac{1}{s^2}$ | $t$ |
| **13** | $\frac{2}{s^3}$ | $\frac{2}{s^3} + \frac{\Delta t}{s^2}$ | $t^2$ |
| **14** | $\frac{6}{s^4}$ | $\frac{6}{s^4} + \frac{6\Delta t}{s^3} + \frac{\Delta t^2}{s^2}$ | $t^3$ |
| **15** | $\frac{1}{(s-a)^2}$ | $\frac{e^{a\Delta t}}{(s - \frac{e^{a\Delta t} - 1}{\Delta t})^2}$ | $te^{at}$ |
| **16** | $\frac{2}{(s-a)^3}$ | $\frac{e^{a\Delta t}(s\Delta t + e^{a\Delta t} + 1)}{(s - \frac{e^{a\Delta t} - 1}{\Delta t})^3}$ | $t^2e^{at}$ |

| # | $K_{\Delta t}$ Transform $\Delta t \rightarrow 0$<br>Laplace Transform<br>F(t) where $0 \le t < \infty$ | $K_{\Delta t}$ Transform<br>Generalized Laplace Transform<br>F(t) where $t = 0, \Delta t, 2\Delta t, 3\Delta t, \ldots$ | F(t)<br>Calculus<br>Function |
|---|---|---|---|
| **17** | $\frac{2bs}{[\ s^2 + b^2\ ]^2}$ | $\frac{\frac{\sin b\Delta t}{\Delta t}(2s+\Delta t s^2)}{[(s + \frac{1-\cos b\Delta t}{\Delta t})^2 + (\frac{\sin b\Delta t}{\Delta t})^2]^2}$ | tsinbt |
| **18** | $\frac{2bs}{[\ s^2 - b^2\ ]^2}$ | $\frac{\frac{\sinh b\Delta t}{\Delta t}(2s+\Delta t s^2)}{[(s + \frac{1-\cosh b\Delta t}{\Delta t})^2 - (\frac{\sinh b\Delta t}{\Delta t})^2]^2}$ | tsinhbt |
| **19** | $\frac{s^2 - b^2}{[\ s^2 + b^2\ ]^2}$ | $\frac{[(s + \frac{1-\cos b\Delta t}{\Delta t})^2 - (\frac{\sin b\Delta t}{\Delta t})^2] - (\frac{1-\cos b\Delta t}{\Delta t})[\Delta t s^2 + 4s + 2(\frac{1-\cos b\Delta t}{\Delta t})]}{[(s + \frac{1-\cos b\Delta t}{\Delta t})^2 + (\frac{\sin b\Delta t}{\Delta t})^2]^2}$ | tcosbt |
| **20** | $\frac{s^2 + b^2}{[\ s^2 - b^2\ ]^2}$ | $\frac{[(s + \frac{1-\cosh b\Delta t}{\Delta t})^2 + (\frac{\sinh b\Delta t}{\Delta t})^2] - (\frac{1-\cosh b\Delta t}{\Delta t})[\Delta t s^2 + 4s + 2(\frac{1-\cosh b\Delta t}{\Delta t})]}{[(s + \frac{1-\cosh b\Delta t}{\Delta t})^2 - (\frac{\sinh b\Delta t}{\Delta t})^2]^2}$ | tcoshbt |

## Table 5.7-3 Conversions of $K_{\Delta t}$ Transforms to Calculus Functions

| # | $K_{\Delta t}$ Transform Generalized Laplace Transform of F(t) | F(t) Calculus Functions<br><br>F(t) where t = 0, Δt, 2Δt, … | Coefficient Definitions |
|---|---|---|---|
| 1 | $\frac{1}{s}$ | 1 | |
| 2 | $\frac{1}{s-a}$ | $e^{\alpha t}$<br><br>Comment – The Equivalent Interval Calculus Function is:<br><br>$e_{\Delta t}(\frac{e^{\alpha\Delta t}-1}{\Delta t}, t)$ | Denominator polynomial a real root<br><br>$\alpha = \frac{\ln(1+a\Delta t)}{\Delta t}$ |
| 3 | $\frac{b}{(s-a)^2+b^2}$ | $e^{\alpha t}\sin\beta t$<br><br>Comment – The Equivalent non-Calculus Function is:<br><br>$(1+a\Delta t)^{\frac{t}{\Delta t}}\sin_{\Delta t}(\frac{b\Delta t}{1+a\Delta t}, t)$ | Denominator polynomial complex conjugate roots = a ± jb<br><br>$\alpha = \frac{1}{2\Delta t}\ln[(1+a\Delta t)^2 + (b\Delta t)^2]$<br><br>$\beta = \begin{cases} \frac{1}{\Delta t}\tan^{-1}\left\lvert\frac{b\Delta t}{1+a\Delta t}\right\rvert & \text{for } 1+a\Delta t > 0 \quad b\Delta t \ge 0 \\ -\frac{1}{\Delta t}\tan^{-1}\left\lvert\frac{b\Delta t}{1+a\Delta t}\right\rvert & \text{for } 1+a\Delta t > 0 \quad b\Delta t < 0 \\ \frac{1}{\Delta t}[\pi - \tan^{-1}\left\lvert\frac{b\Delta t}{1+a\Delta t}\right\rvert] & \text{for } 1+a\Delta t < 0 \quad b\Delta t \ge 0 \\ \frac{1}{\Delta t}[-\pi + \tan^{-1}\left\lvert\frac{b\Delta t}{1+a\Delta t}\right\rvert] & \text{for } 1+a\Delta t < 0 \quad b\Delta t < 0 \end{cases}$<br><br>$1+a\Delta t \ne 0$<br><br>$0 \le \tan^{-1}\left\lvert\frac{b\Delta t}{1+a\Delta t}\right\rvert < \frac{\pi}{2}$<br><br>Note – For sinβt, $(1+a\Delta t)^2 + (b\Delta t)^2 = 1$ |

| # | $K_{\Delta t}$ Transform Generalized Laplace Transform of F(t) | F(t) Calculus Functions<br><br>F(t) where t = 0, Δt, 2Δt, … | Coefficient Definitions |
|---|---|---|---|
| **4** | $\frac{s-a}{(s-a)^2+b^2}$ | $e^{\alpha t}\cos\beta t$<br><br>Comment – The Equivalent non-Calculus Function is:<br>$(1+a\Delta t)^{\frac{t}{\Delta t}}\cos_{\Delta t}(\frac{b\Delta t}{1+a\Delta t}, t)$ | Denominator polynomial complex conjugate roots = a ± jb<br>$\alpha = \frac{1}{2\Delta t}\ln[(1+a\Delta t)^2 + (b\Delta t)^2]$<br><br>$\beta = \left[\begin{array}{ll} \frac{1}{\Delta t}\tan^{-1}\left\lvert\frac{b\Delta t}{1+a\Delta t}\right\rvert & \text{for } 1+a\Delta t > 0 \quad b\Delta t \ge 0 \\ -\frac{1}{\Delta t}\tan^{-1}\left\lvert\frac{b\Delta t}{1+a\Delta t}\right\rvert & \text{for } 1+a\Delta t > 0 \quad b\Delta t < 0 \\ \frac{1}{\Delta t}[\pi - \tan^{-1}\left\lvert\frac{b\Delta t}{1+a\Delta t}\right\rvert] & \text{for } 1+a\Delta t < 0 \quad b\Delta t \ge 0 \\ \frac{1}{\Delta t}[-\pi + \tan^{-1}\left\lvert\frac{b\Delta t}{1+a\Delta t}\right\rvert] & \text{for } 1+a\Delta t < 0 \quad b\Delta t < 0 \end{array}\right.$<br><br>$1+a\Delta t \ne 0$<br>$0 \le \tan^{-1}\left\lvert\frac{b\Delta x}{1+a\Delta x}\right\rvert < \frac{\pi}{2}$<br><br>Note – For cosβt, $(1+a\Delta t)^2 + (b\Delta t)^2 = 1$ |
| **5** | $(1+s\Delta t)^{-n}F(s)$ | $f(t-\Delta t)U(t-n\Delta t)$ | The Unit Step Function, U(t-nΔt)<br>$U(t-n\Delta t) = \left\{\begin{array}{ll} 1 & t \ge n\Delta t \\ 0 & t < n\Delta t \end{array}\right\}$<br><br>n = 0, 1, 2, 3, …<br>F(s) = $K_{\Delta t}$ Transform of f(t)<br><br>Comment<br>To find c(t) where $c(t) = K_{\Delta t}^{-1}[(1+s\Delta t)^{-n}F(s)]$<br>First find f(t) where $f(t) = K_{\Delta t}^{-1}[F(s)]$<br>Then $c(t) = f(t-n\Delta t)U(t-n\Delta t)$ |
| **6** | $\frac{n!}{s^{n+1}}$ | $\prod_{m=1}^{n}(t-[m-1]\Delta t)$ | n = 1,2,3,… |

| # | $K_{\Delta t}$ Transform Generalized Laplace Transform of F(t) | F(t) Calculus Functions<br><br>F(t) where t = 0, Δt, 2Δt, … | Coefficient Definitions |
|---|---|---|---|
| **7** | $\frac{(1+a\Delta t)^n\, n!}{(s-a)^{n+1}}$ | $e^{\alpha t}\prod_{m=1}^{n}(t-[m-1]\Delta t)$ | Denominator polynomial a multiple real root<br>$\alpha = \frac{\ln(1+a\Delta t)}{\Delta t}$<br>n = 1,2,3,… |

Example 5.7-3 is an example of the solution of a differential difference equation with the solution obtained in two forms. One solution is in the form of Interval Calculus functions and the other solution is in the form of equivalent Calculus functions. The two solutions are identities.

Example 5.7-3 Solve the differential difference equation, $D_{\Delta t}^2\, y + D_{\Delta t}y + 2y = 7$. Obtain the solution in two forms, one solution in terms of Interval Calculus functions and the other equivalent solution in terms of Calculus functions. Use $K_{\Delta t}$ Transforms to obtain the solutions. y(0) = 0, $y'(0) = 0$, and Δt = .1. Find y(1).

$$D_{\Delta t}^2\, y + D_{\Delta t}y + 2y = 7 \qquad 1)$$

Take the $K_{\Delta t}$ Transform of Eq 1

$$s^2y(s) - sy(s) - 0 + sy(s) - 0 + 2y(s) = \frac{7}{s} \qquad 2)$$

$$s^2 + sy(s) + 2y(s) = (s^2 + s + 2)y(s) = \frac{7}{s} \qquad 3)$$

$$y(s) = \frac{7}{s(s^2 + s + 2)} = \frac{A}{s} + \frac{Bs + C}{s^2 + s + 2} \qquad 4)$$

Find A

Multiplying each term of Eq 4 by s and setting s = 0

$$A = \frac{7}{s^2 + s + 2}\Big|_{s=0} = \frac{7}{2} = 3.5$$

$$A = 3.5 \qquad 5)$$

From Eq 4 and Eq 5

$$y(s) = \frac{7}{s(s^2 + s + 2)} = \frac{3.5}{s} + \frac{Bs + C}{s^2 + s + 2} = \frac{3.5s^2 + 3.5s + 7 + Bs^2 + Cs}{s(s^2 + s + 2)} \qquad 6)$$

Find B

3.5 + B = 0

$$B = -3.5 \qquad 7)$$

Find C

3.5 + C = 0

$$C = -3.5 \qquad 8)$$

Substituting Eq 7 and Eq 8 into Eq 6

$$y(s) = \frac{7}{s(s^2 + s + 2)} = \frac{3.5}{s} - 3.5\,\frac{s + 1}{s^2 + s + 2} = \frac{3.5}{s} - 3.5\,\frac{s + 1}{(s + \frac{1}{2})^2 + (\frac{\sqrt{7}}{2})^2} \qquad 9)$$

$$y(s) = \frac{3.5}{s} - 3.5\,\frac{s + 1}{(s + \frac{1}{2})^2 + (\frac{\sqrt{7}}{2})^2} \qquad 10)$$

Changing the form of Eq 10 in order to find the inverse $K_{\Delta t}$ transform of y(s)

$$y(s) = 3.5\left[\frac{1}{s} - \frac{s - (-\frac{1}{2})}{(s - \{-\frac{1}{2}\})^2 + (\frac{\sqrt{7}}{2})^2} - \frac{\frac{1}{2}}{(s - \{-\frac{1}{2}\})^2 + (\frac{\sqrt{7}}{2})^2}\right] \qquad 11)$$

$$y(s) = 3.5\left[\frac{1}{s} - \frac{s - (-\frac{1}{2})}{(s - \{-\frac{1}{2}\})^2 + (\frac{\sqrt{7}}{2})^2} - \frac{1}{\sqrt{7}}\,\frac{\frac{\sqrt{7}}{2}}{(s - \{-\frac{1}{2}\})^2 + (\frac{\sqrt{7}}{2})^2}\right] \qquad 12)$$

Finding the inverse $K_{\Delta t}$ Transform of y(s), $K_{\Delta t}^{-1}[y(s)]$

<u>Find an Interval Calculus function solution to Eq 12</u>

From Table 3 in the Appendix

$$K_{\Delta t}^{-1}\left[\frac{b}{(s-a)^2 + b^2}\right] = (1+a\Delta t)^{\frac{t}{\Delta t}} \sin_{\Delta t}\left(\frac{b}{1+a\Delta t}, t\right) \qquad 13)$$

$$K_{\Delta t}^{-1}[\frac{s-a}{(s-a)^2 + b^2}] = (1+a\Delta t)^{\frac{t}{\Delta t}} \cos_{\Delta t}(\frac{b}{1+a\Delta t}, t) \qquad 14)$$

$$K_{\Delta t}^{-1}[\frac{1}{s}] = 1 \qquad 15)$$

Substituting Eq 13 thru Eq 15 into Eq 12

$$y(t) = 3.5\left[1 - (1+a\Delta t)^{\frac{t}{\Delta t}} \cos_{\Delta t}(\frac{b}{1+a\Delta t}, t) - \frac{1}{\sqrt{7}}(1+a\Delta t)^{\frac{t}{\Delta t}} \sin_{\Delta t}(\frac{\frac{\sqrt{7}}{2}}{1+a\Delta t}, t)\right] \qquad 16)$$

Let $a = -\frac{1}{2}$

$b = \frac{\sqrt{7}}{2}$

$\Delta t = .1$

$$y(t) = 3.5\left[1 - (1 - \tfrac{1}{2}\{.1\})^{\frac{t}{.1}} \cos_{.1}(\frac{\frac{\sqrt{7}}{2}}{1 + \frac{1}{2}\{.1\}}, t) - \frac{1}{\sqrt{7}}(1 - \tfrac{1}{2}\{.1\})^{\frac{t}{.1}} \sin_{.1}(\frac{\frac{\sqrt{7}}{2}}{1 + \frac{1}{2}\{.1\}}, t)\right] \qquad 17)$$

$$\mathbf{y(t) = 3.5\,[1 - (.95)^{10t}\cos_{.1}(1.39250069,t) - .377964473\,(.95)^{10t}\sin_{.1}(1.39250069,t)]} \qquad 18)$$

Find y(1)

From Eq 18

$$y(t) = 3.5\,[1 - (.95)^{10}\cos_{.1}(1.39250069,1) - .377964473\,(.95)^{10}\sin_{.1}(1.39250069,1)] \qquad 19)$$

Using the computer program, ICFNCALC, to calculate the Interval Calculus functions, $\cos_{\Delta t}(b,t)$ and $\sin_{\Delta t}(b,t)$

$$y(s) = 3.5\,[1 - (.95)^{10}(.20485800) - .377964473\,(.95)^{10}(1.0815563)] \qquad 20)$$

$$y(1) = 2.2140522 \qquad 21)$$

<u>Find a Calculus function solution to Eq 12</u>

From Table 5.7-3 or Table 3b in the Appendix

$$K_{\Delta t}^{-1}[\frac{b}{(s-a)^2 + b^2}] = e^{\alpha t}\sin\beta t \qquad 22)$$

$$K_{\Delta t}^{-1}\left[\frac{s-a}{(s-a)^2+b^2}\right] = e^{\alpha t}\cos\beta t \qquad 23)$$

$$\alpha = \frac{1}{2\Delta t}\ln[(1+a\Delta t)^2 + (b\Delta t)^2] \qquad 24)$$

$$\beta = \begin{cases} \frac{1}{\Delta t}\tan^{-1}\left|\frac{b\Delta t}{1+a\Delta t}\right| & \text{for } 1+a\Delta t > 0 \quad b\Delta t \ge 0 \\ -\frac{1}{\Delta t}\tan^{-1}\left|\frac{b\Delta t}{1+a\Delta t}\right| & \text{for } 1+a\Delta t > 0 \quad b\Delta t < 0 \\ \frac{1}{\Delta t}\left[\pi - \tan^{-1}\left|\frac{b\Delta t}{1+a\Delta t}\right|\right] & \text{for } 1+a\Delta t < 0 \quad b\Delta t \ge 0 \\ \frac{1}{\Delta t}\left[-\pi + \tan^{-1}\left|\frac{b\Delta t}{1+a\Delta t}\right|\right] & \text{for } 1+a\Delta t < 0 \quad b\Delta t < 0 \end{cases} \qquad 25)$$

$1+a\Delta t \neq 0$

$0 \le \tan^{-1}\left|\frac{b\Delta x}{1+a\Delta x}\right| < \frac{\pi}{2}$

$$K_{\Delta t}^{-1}\left[\frac{1}{s}\right] = 1 \qquad 26)$$

Find the inverse $K_{\Delta t}$ Transform of y(s)

Substituting Eq 22 thru Eq 26 into Eq 12

$$y(t) = 3.5\left[1 - e^{\alpha t}\cos\beta t - \frac{1}{\sqrt{7}}e^{\alpha t}\sin\beta t\right] \qquad 27)$$

where

$$\alpha = \frac{1}{2\Delta t}\ln[(1+a\Delta t)^2 + (b\Delta t)^2]$$

$$\beta = \begin{cases} \frac{1}{\Delta t}\tan^{-1}\left|\frac{b\Delta t}{1+a\Delta t}\right| & \text{for } 1+a\Delta t > 0 \quad b\Delta t \ge 0 \\ -\frac{1}{\Delta t}\tan^{-1}\left|\frac{b\Delta t}{1+a\Delta t}\right| & \text{for } 1+a\Delta t > 0 \quad b\Delta t < 0 \\ \frac{1}{\Delta t}\left[\pi - \tan^{-1}\left|\frac{b\Delta t}{1+a\Delta t}\right|\right] & \text{for } 1+a\Delta t < 0 \quad b\Delta t \ge 0 \\ \frac{1}{\Delta t}\left[-\pi + \tan^{-1}\left|\frac{b\Delta t}{1+a\Delta t}\right|\right] & \text{for } 1+a\Delta t < 0 \quad b\Delta t < 0 \end{cases}$$

$1+a\Delta t \neq 0$

$0 \le \tan^{-1}\left|\frac{b\Delta x}{1+a\Delta x}\right| < \frac{\pi}{2}$

Let $a = -\frac{1}{2}$

$b = \frac{\sqrt{7}}{2}$

$\Delta t = .1$

$$\alpha = \frac{1}{2(.1)} \ln[(1 - \frac{1}{2}\{.1\})^2 + (\frac{\sqrt{7}}{2}\{.1\})^2] = 5 \ln[(.95)^2 + (.132287565)^2] \quad 28)$$

$$\alpha = -.416908045 \quad 29)$$

$$1+a\Delta t = 1 - \frac{1}{2}\{.1\}) = .95 > 0 \quad 30)$$

$$b\Delta t = \frac{\sqrt{7}}{2}\{.1\}) = .132287565 > 0 \quad 31)$$

$$\beta = \frac{1}{\Delta t} \tan^{-1}\left|\frac{b\Delta x}{1+a\Delta x}\right| = \frac{1}{.1} \tan^{-1}\left|\frac{.132287565}{.95}\right| = 1.3836035 \quad 32)$$

$$\beta = 1.3836035 \quad 33)$$

Substituting Eq 29 and Eq 33 into Eq 27

$$\mathbf{y(t) = 3.5 \left[1 - e^{-.416908045t}\cos(1.3836035t) - .377964473\, e^{-.416908045t}\sin(1.3836035t)\right]} \quad 34)$$

find y(1)

$$y(t) = 3.5 [1 - e^{-.416908045}\cos(1.3836035) - .377964473\, e^{-.416908045}\sin(1.3836035)] \quad 35)$$

$$y(1) = 2.2140522 \quad 36)$$

Note that the values calculated for y(1) by the Interval Calculus function solution, Eq 18, and the Calculus function solution, Eq 34, are the same.The two solutions are identities.

## Section 5.8: A demonstration of problem solving using both Interval Calculus and Calculus functions

For demonstration, consider the following Calculus and Interval Calculus functions

Diagram 5.8-1

Calculus Function

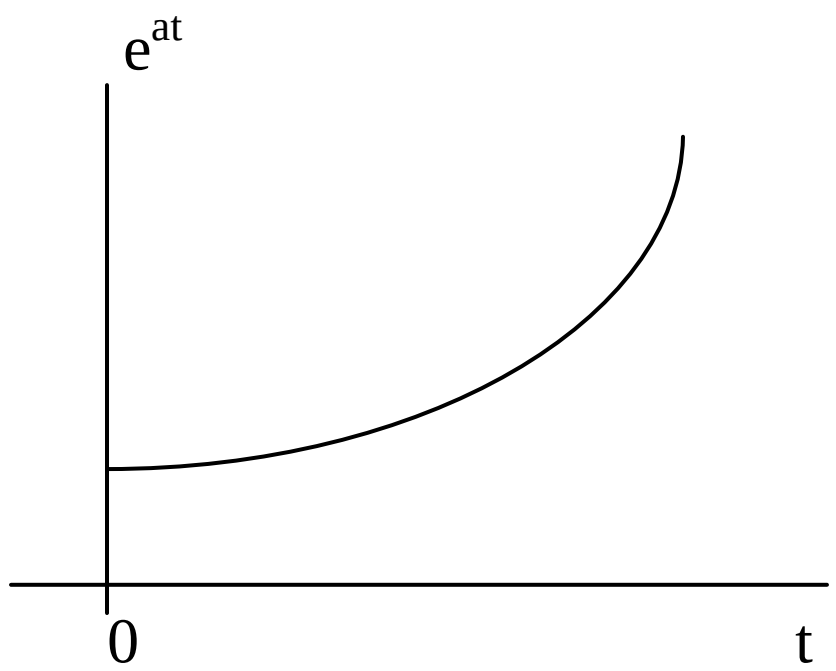

Smooth continuous curve

$e^{at}$ solves differential equations

Diagram 5.8-2

Interval Calculus Discrete Function

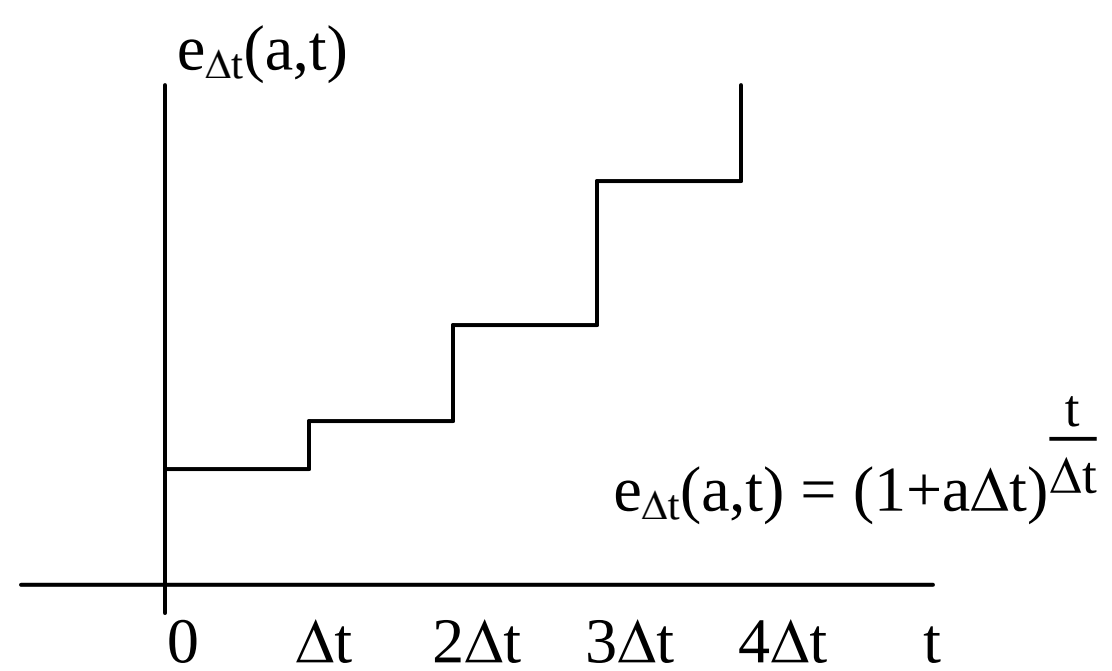

Sample and hold shaped continuous curve

$e_{\Delta t}(a,t)$ solves differential difference equations

The four integrations below calculate the area under the two functions specified above.

Integration of the Calculus function, $e^{at}$

Calculation of the area under the continuous curve of the Calculus function, $e^{at}$

$$A = \int_{t_1}^{t_2} \mathbf{e^{at}\,dt} = \frac{e^{at}}{a}\Big|_{t_1}^{t_2} = \frac{1}{a}\left[e^{at_2} - e^{at_1}\right] \qquad (5.8\text{-}1)$$

where

$$e^{at} = \lim_{\Delta t \to 0}(1+a\Delta t)^{\frac{t}{\Delta t}}$$

$t_1 \le t \le t_2$

$a = \text{constant}$

Discrete integration of the discrete Interval Calculus function, $e_{\Delta t}(a,t)$

Calculation of the area under the sample and hold shaped continuous curve of the discrete Interval Calculus function, $e_{\Delta t}(a,t)$

$$A = {}_{\Delta t}\int_{t_1}^{t_2} \mathbf{e_{\Delta t}(a,t)\,\Delta t} = \Delta t \; {}_{\Delta t}\sum_{t=t_1}^{t_2-\Delta t} e_{\Delta t}(a,t) = \frac{e_{\Delta t}(a,t)}{a}\Big|_{t_1}^{t_2} = \frac{1}{a}\left[e_{\Delta t}(a,t_2) - e_{\Delta t}(a,t_1)\right] \qquad (5.8\text{-}2)$$

where

$$e_{\Delta t}(a,t) = (1+a\Delta t)^{\frac{t}{\Delta t}}$$

$t,\ t_1+3\Delta t,\ \ldots,\ t_2-\Delta t,\ t_2$

a = constant
$\Delta t$ = interval between successive values of t

Integration of the discrete Interval Calculus function, $e_{\Delta t}(a,t)$

Calculation of the area under the sample and hold continuous curve of the discrete Interval Calculus function, $e_{\Delta t}(a,t)$

$$A = \int_{t_1}^{t_2} \mathbf{e_{\Delta t}(a,t)\,dt} = {}_{\Delta t}\!\int_{t_1}^{t_2} e_{\Delta t}(a,t)\,\Delta t = \Delta t \sum_{\Delta t\; t=t_1}^{t_2-\Delta t} e_{\Delta t}(a,t) = \frac{e_{\Delta t}(a,t)}{a}\Big|_{t_1}^{t_2} = \frac{1}{a}[e_{\Delta t}(a,t_2) - e_{\Delta t}(a,t_1)] \qquad (5.8\text{-}3)$$

where

$e_{\Delta t}(a,t) = (1+a\Delta t)^{\frac{t}{\Delta t}}$
$t_1 \le t \le t_2$
a = constant
$\Delta t$ = increment of t

Discrete integration of the Calculus function, $e^{at}$, evaluated at equally spaced discrete values of t

Calculation of the area under the sample and hold shaped continuous curve of the Calculus function, $e^{at}$, evaluated at equally spaced discrete values of t. The Interval Calculus identity for $e^{at}$ is used.

$$\mathbf{A = {}_{\Delta t}\!\int_{t_1}^{t_2} e^{at}\,\Delta t} = \Delta t \sum_{\Delta t\; t=t_1}^{t_2-\Delta t} e^{at} = {}_{\Delta t}\!\int_{t_1}^{t_2} e_{\Delta t}(\alpha, t)\,\Delta t = \frac{e_{\Delta t}(\alpha, t)}{\alpha}\Big|_{t_1}^{t_2} = \frac{1}{\alpha}[e_{\Delta t}(\alpha, t_2) - e_{\Delta t}(\alpha, t_1)] \qquad (5.8\text{-}4)$$

where

$e^{at} = e_{\Delta t}(\alpha, t)$ , Interval Calculus identity
$e_{\Delta t}(\alpha,t) = (1+\alpha\Delta t)^{\frac{t}{\Delta t}}$
$\alpha = \dfrac{e^{a\Delta t} - 1}{\Delta t}$
$t = t_1, t_1+\Delta t, t_1+2\Delta t, t_1+3\Delta t, \ldots, t_2-\Delta t, t_2$
a, $\alpha$ = constants
$\Delta t$ = interval between successive values of t

Eq 5.8-1 thru Eq 5.8-4 are integrations of the functions, $e^{at}$ or $e_{\Delta t}(a,t)$, in order to calculate the areas under the specified curves. Eq 5.8-1 is an example of a Calculus integration of a continuous Calculus function. Eq 5.8-2 is an example of a discrete integration of a continuous Interval Calculus discrete function. These two integrations are straight forward. The type of integration matches the type of function to be integrated (i.e. the subscripts match). The integration formulas can be obtained from Calculus or Interval Calculus integration tables. Interval Calculus integration formula tables appear in the Appendix of this paper.

Eq 5.8-3 is an example of a Calculus integration of a discrete Interval Calculus function. It is an integration of a discrete function with a sample and hold continuous curve. As shown, this integration is equal to the discrete integral of the same discrete function. Such an integration is shown in Eq 5.8-2. Thus, Eq 5.8-3 is an integration that can also be performed in a straight forward manner. Eq 5.8-4 is different. It represents the discrete integration of a continuous Calculus function. This type of integration is necessary to sum values of a continuous Calculus function at equally spaced discrete values of t. The mismatch between the discrete integration and the continuous Calculus function must be resolved. For this integration, there are no integration tables available. There is, however, a straight forward methodology for performing this type of integration. It involves the use of an Interval Calculus identity. An Interval Calculus identity can be substituted for the Calculus function to be integrated and then an Interval Calculus discrete integration similar to Eq 5.8-2 can be performed. Note that in the integration procedure of Eq 5.8-4, an Interval Calculus identity for $e^{at}$ is used. Interval Calculus identities appear in several tables in the Appendix of this paper. In particular, note Table 5 and Table 5a. Some examples of the use of an $e^{at}$ identity in the solution of problems involving sample and hold shaped functions, summations, and discrete derivatives are shown below.

Example 5.8 -1 Find the area under the sampled exponential function, $e^{2t}$, from t = 0 thru 4 where $\Delta t = .5$.

$$A = .5 \;{}_{.5}\!\sum_{t=0}^{3.5} e^{2t} = {}_{.5}\!\int_{0}^{4} e^{2t}\,\Delta t \qquad 1)$$

$e^{2t}$ is sampled at t = 0, .5, 1, 1.5, …, 3.5 and the sampled value is held for a period of .5 (i.e. $\Delta t = .5$)

Using an Interval Calculus identity for $e^{at}$ to facilitate the integration operation

$$e^{at} = e_{\Delta t}(\alpha, t)\ , \quad \text{Interval Calculus identity for } e^{at} \qquad 2)$$

where

$\alpha = \dfrac{e^{a\Delta t} - 1}{\Delta t}$

$t_i$ = initial value of t

$t_f$ = final value of t

$t = t_i, t_i+\Delta t, t_i+2\Delta t, t_i+3\Delta t, \ldots, t_f-\Delta t, t_f$

a, $\alpha$ = constants

$\Delta t$ = interval between successive values of t

$$a = 2 \qquad 3)$$

$$\Delta t = .5 \qquad 4)$$

$$t_i = 0 \qquad 5)$$

$$t_f = 4 \qquad 6)$$

$$\alpha = \frac{e^{a\Delta t} - 1}{\Delta t} = \frac{e - 1}{.5} = 3.436563657 \qquad 7)$$

$$e^{2t} = e_{.5}(3.436563657, t) \qquad 8)$$

Substituting Eq 8 into Eq 1

$$A = {}_{.5}\int_0^4 e^{2t}\Delta t = {}_{.5}\int_0^4 e_{.5}(3.436563657, t)\,\Delta t = \frac{e_{.5}(3.436563657, t)}{3.436563657}\Big|_0^4 \qquad 9)$$

$$A = \frac{1}{3.436563657}\,[e_{.5}(3.436563657, 4) - e_{.5}(3.436563657, 0)] \qquad 10)$$

$$e_{\Delta t}(a,t) = (1+a\Delta t)^{\frac{t}{\Delta t}} \qquad 11)$$

$$e_{.5}(3.436563657, 0) = (1+3.436563657[.5])^{\frac{0}{.5}} = (2.718281829)^0 = 1 \qquad 12)$$

$$e_{.5}(3.436563657, 4) = (1+3.436563657[.5])^{\frac{4}{.5}} = (2.718281829)^8 = 2980.957992 \qquad 13)$$

Substituting Eq 12 and Eq 13 into Eq 10

$$A = \frac{1}{3.436563657}\,[2980.957992 - 1] = 867.1330693 \qquad 14)$$

$$A = {}_{.5}\int_0^4 e^{2t}\Delta t = 867.1330693 \qquad 15)$$

Finding the value of the summation, ${}_{.5}\sum_{t=0}^{3.5} e^{2t}$

Using the value of A from Eq 15

.

$${}_{.5}\sum_{t=0}^{3.5} e^{2t} = \frac{1}{.5}\,{}_{.5}\int_0^4 e^{2t}\Delta t = \frac{A}{.5} = 2(867.1330693) = 1734.266136 \qquad 16)$$

$$\boldsymbol{{}_{.5}\sum_{t=0}^{3.5} e^{2t} = 1734.26613613} \qquad 17)$$

Checking

$${}_{.5}\sum_{t=0}^{3.5} e^{2t} = e^0 + e^1 + e^2 + e^3 + e^4 + e^5 + e^6 + e^7 = 1734.266136 \qquad 18)$$

$$\sum_{t=0}^{3.5}{}_{.5}\, e^{2t} = 1734.266136 \qquad 19)$$

Good check

Interval Calculus identities are also useful in discrete differentiation. A demonstration of this is shown in Example 5.8-2 and Example 5.8-3 below.

Example 5.8-2 Find the discrete derivative $D_{\Delta t}{}^3e^{at}$ using the Interval Calculus identity for the Calculus function,$e^{at}$

The discrete derivative of the Interval Calculus function, $e_{\Delta t}(b,t)$, is:

$$D_{\Delta t}e_{\Delta t}(b,t) = be_{\Delta t}(b,t) \qquad 1)$$

The interval Calculus identity for the Calculus function, $e^{at}$, is:

$$e^{at} = e_{\Delta t}(\frac{e^{a\Delta t}-1}{\Delta t}, t) \qquad 2)$$

The values of t are discrete values separated by equal intervals of $\Delta t$

From Eq 1 and Eq 2

$$D_{\Delta t}{}^3e^{at} = D_{\Delta t}{}^3e_{\Delta t}(\frac{e^{a\Delta t}-1}{\Delta t}, t) = \frac{(e^{a\Delta t}-1)}{\Delta t}D_{\Delta t}{}^2e_{\Delta t}(\frac{e^{a\Delta t}-1}{\Delta t}, t) = \left(\frac{(e^{a\Delta t}-1)}{\Delta t}\right)^2 D_{\Delta t}e_{\Delta t}(\frac{e^{a\Delta t}-1}{\Delta t}, t) \qquad 3)$$

$$D_{\Delta t}{}^3e^{at} = \left(\frac{(e^{a\Delta t}-1)}{\Delta t}\right)^3 e_{\Delta t}(\frac{e^{a\Delta t}-1}{\Delta t}, t) = \left(\frac{(e^{a\Delta t}-1)}{\Delta t}\right)^3 e^{at} \qquad 4)$$

$$\mathbf{D_{\Delta t}{}^3e^{at} = \left(\frac{(e^{a\Delta t}-1)}{\Delta t}\right)^3 e^{at}} \qquad 5)$$

Compare the above calculation of $D_{\Delta t}{}^3e^{at}$ to the calculation obtained using the definition of the discrete derivative

$$D_{\Delta t}y(t) = \frac{y(t+\Delta t) - y(t)}{\Delta t}, \quad \text{Definition of the discrete derivative} \qquad 6)$$

$$D_{\Delta t}{}^3e^{at} = D_{\Delta t}{}^2[D_{\Delta t}e^{at}] = D_{\Delta t}{}^2[\frac{e^{a(t+\Delta t)}-e^{at}}{\Delta t}] = D_{\Delta t}{}^2[\frac{e^{at}(e^{a\Delta t}-1)}{\Delta t}] = \frac{(e^{a\Delta t}-1)}{\Delta t}D_{\Delta t}{}^2[e^{at}] = \frac{(e^{a\Delta t}-1)}{\Delta t}D_{\Delta t}[D_{\Delta t}e^{at}] \qquad 7)$$

$$D_{\Delta t}{}^3e^{at} = \frac{(e^{a\Delta t}-1)}{\Delta t}D_{\Delta t}[D_{\Delta t}e^{at}] = \frac{(e^{a\Delta t}-1)}{\Delta t}D_{\Delta t}[\frac{e^{a(t+\Delta t)}-e^{at}}{\Delta t}] = \frac{(e^{a\Delta t}-1)}{\Delta t}D_{\Delta t}[\frac{e^{at}(e^{a\Delta t}-1)}{\Delta t}] = \left(\frac{(e^{a\Delta t}-1)}{\Delta t}\right)^2 D_{\Delta t}[e^{at}] \qquad 8)$$

$$D_{\Delta t}{}^3e^{at} = \left(\frac{(e^{a\Delta t}-1)}{\Delta t}\right)^2 D_{\Delta t}[e^{at}] = \left(\frac{(e^{a\Delta t}-1)}{\Delta t}\right)^2[\frac{e^{a(t+\Delta t)}-e^{at}}{\Delta t}] = \left(\frac{(e^{a\Delta t}-1)}{\Delta t}\right)^2[\frac{e^{at}(e^{a\Delta t}-1)}{\Delta t}] = \left(\frac{(e^{a\Delta t}-1)}{\Delta t}\right)^3 e^{at} \qquad 9)$$

$$\mathbf{D_{\Delta t}{}^3 e^{at} = \left(\frac{(e^{a\Delta t}-1)}{\Delta t}\right)^3 e^{at}} \qquad 10)$$

Both methods obtain the same result. However, the use of the $e^{at}$ Interval Calculus identity is seen to yield a solution with less mathematical manipulation.

Example 5.8-3 Find the solution to the differential equation, $(DD_{.9} + 1)y(t) = t$, using the Interval Calculus identity for the Calculus function, $e^{at}$

$(DD_{.9} + 1)y(t) = t$ 1)

where

$D = \frac{d}{dt}$

Find the homogeneous equation solution, $y_h(t)$, for Eq 1

$(DD_{.9} + 1)y_h(t) = 0$ , Homogeneous equation of Eq 1 2)

$D_{.9}Dy(t) = DD_{.9}y(t)$ 3)

Note - $D_{.9}D = DD_{.9}$ is a discrete derivative operation equality

Substitute Eq 3 into Eq 2

$(D_{.9}D + 1)y_h(t) = 0$ 4)

Assume

$y_h(t) = Ke^{at}$ 5)

where

$K$ = constant

From Eq 4 and Eq 5

$aD_{.9}e^{at} + e^{at} = 0$ 6)

Using the Interval Calculus identity for the Calculus function, $e^{at}$

$e^{at} = e_{\Delta t}(\frac{e^{a\Delta t} - 1}{\Delta t}, \Delta t)$ , The Interval Calculus $e^{at}$ identity 7)

$\Delta t = .9$ 8)

From Eq 6 thru Eq 8

$aD_{.9}e_{.9}(\frac{e^{.9a} - 1}{.9}, x) + e^{at} = 0$ 9)

$$a(\frac{e^{.9a}-1}{.9})e_{.9}(\frac{e^{.9a}-1}{.9}, x) + e^{at} = 0 \qquad 10)$$

$$a(\frac{e^{.9a}-1}{.9})e^{at} + e^{at} = 0 \qquad 11)$$

$$a(\frac{e^{.9a}-1}{.9}) + 1 = 0 \qquad 12)$$

From Eq 12 find the constant, a

Using the internet site, WolframAlpha to find the constant, a

$$a = .210461 \pm .944384j \qquad 13)$$

The solution of Eq 2, the homogeneous equation solution of Eq 1, is:

$$y_h(t) = c_1 e^{(.210461+.944384j)t} + c_2 e^{(.210461-.944384j)t} \qquad 14)$$

or

$$y_h(t) = K_1 e^{.210461t}\sin .944384t + K_2 e^{.210461t}\cos .944384t \qquad 15)$$

Find the particular solution, $y_p$, of Eq 1

$$y_p(t) = t \qquad 16)$$

Verifying Eq 16

Substituting Eq 16 into Eq 1

$$DD_{.9}[\,t] + t = t \qquad 17)$$

$$D[1] + t = t \quad , \quad D = \frac{d}{dt} \qquad 18)$$

$$t = t \qquad \text{good check} \qquad 19)$$

$$y(t) = y_h(t) + y_p(t) \qquad 20)$$

Then the solution to the discrete differential equation, $(DD_{.9} + 1)y(t) = t$, Eq 1, is:

$$\mathbf{y(t) = K_1 e^{.210461t}\sin .944384t + K_2 e^{.210461t}\cos .944384t + t} \qquad 21)$$

In this section and in Section 5.7, Interval Calculus, a generalization of Calculus, is shown to use both Calculus functions and its own Interval Calculus functions to solve problems involving discrete variables. In Section 5.7 a problem involving the solution of a differential difference equation is solved in two different ways. See Example 5.7-3. One solution uses Calculus functions and the other solution uses

Interval Calculus functions to obtain the same result. Though the forms of the two solutions differ, they can be shown to be equivalent Interval Calculus identities. In this section, Interval Calculus discrete integration of a Calculus function is demonstrated.

In the following section, Section 5.9, the equations to calculate the Fourier Series coefficients of a sample and hold shaped waveform are derived. The derivation involves the use of both Calculus and discrete Interval Calculus. The Fourier Series equations obtained include both Calculus and discrete Interval Calculus functions.

**Section 5.9: Derivation of Fourier Series equations generalized for use in the Fourier Series expansion of discrete sample and hold shaped waveforms**

The Fourier Series coefficient evaluation equations are valid for continuous functions. Sample and hold shaped waveforms can be integrated over a single period but the summation of many interval integrations may be tedious. Also, for sample and hold shaped waveforms, the Fourier Series coefficient evaluation equations do not provide a convenient closed form result. Interval Calculus provides a way to integrate discrete sample and hold shaped waveform functions using a single Interval Calculus discrete integral. The derivation which follows generalizes the Fourier Series coefficient evaluation equations for use in the integration of sample and hold shaped waveform functions. All Interval Calculus functions represent sample and hold shaped waveforms.

Derivation of the Fourier Series equations used in the expansion of sample and hold shaped waveforms

Diagram 5.9-1 - An Interval Calculus sample and hold shaped waveform, f(t)

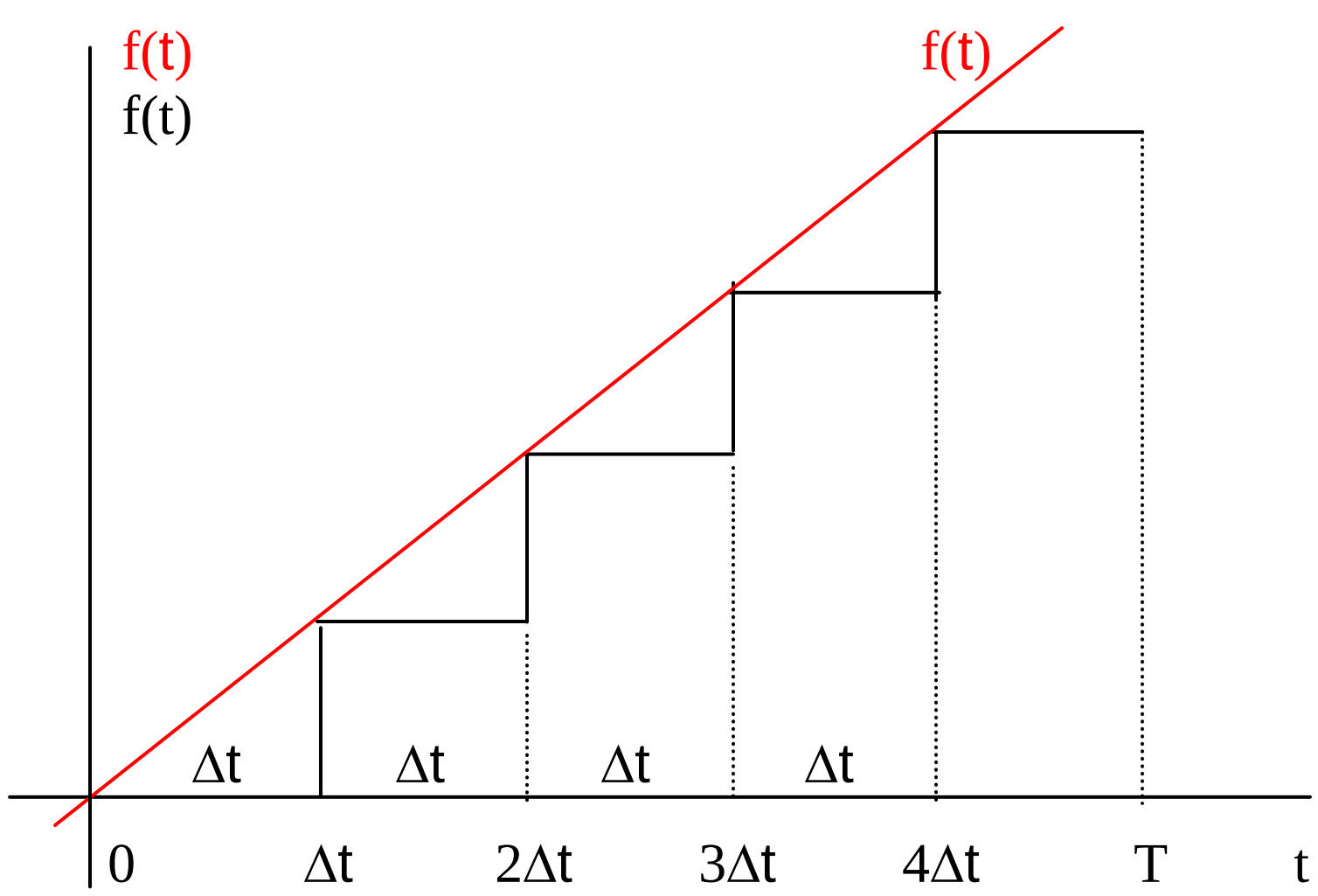

$\Delta t$ = interval between successive discrete values of t

$$\Delta t = \frac{T}{m}$$

T = waveform period
m = number of intervals within the period, 0 to T
m = positive integer
f(t) = continuous function of t where $0 \le t \le T$

f(t) = discrete sample and hold shaped waveform function of t where t = 0, Δt, 2Δt, 3Δt, …, T-Δt, T
f(t) = f(t) where t = t

Comments – t are discrete values of t separated by a common interval , Δt.

The values of t, t = 0, Δt, 2Δt, 3Δt, …, T-Δt, T, are a subset of t.

f(t), a discrete function of t, is equal to f(t), a continuous function of t, where t = 0, Δt, 2Δt, 3Δt, …, T-Δt, T .

Fourier Series Equations where t is continuous, $0 \le t \le T$

$$F(t) = \frac{a_0}{2} + \sum_{n=1}^{\infty} \left( a_n \cos\frac{2\pi n}{T} t + b_n \sin\frac{2\pi n}{T} t \right), \quad n = 1, 2, 3, \ldots \qquad (5.9\text{-}1)$$

$$a_n = \frac{2}{T}\int_0^T f(t)\cos\frac{2\pi n}{T} t\, dt \qquad (5.9\text{-}2)$$

$$b_n = \frac{2}{T}\int_0^T f(t)\sin\frac{2\pi n}{T} t\, dt \qquad (5.9\text{-}3)$$

$$a_0 = \frac{2}{T}\int_0^T f(t)\, dt \qquad (5.9\text{-}4)$$

Find $a_0$

$$a_0 = \frac{2}{T}\int_0^T f(t)\, dt = \frac{2}{T}{}_{\Delta t}\sum_{t=0}^{T-\Delta t} f(t)\Delta t = \frac{2}{T}{}_{\Delta t}\int_0^T f(t)\Delta t \qquad (5.9\text{-}5)$$

$$\mathbf{a_0 = \frac{2}{T}{}_{\Delta t}\int_0^T f(t)\Delta t} \qquad (5.9\text{-}6)$$

**where**

**f(t) = discrete sample and hold shaped waveform function of t**
**t = 0, Δt, 2Δt, 3Δt, …, T-Δt, T**
**Δt = interval between successive values of t**
**T = waveform period**

Find $a_n$ , n = 0, 1, 2, 3, …

$$a_n = \frac{2}{T}\int_0^T f(t)\cos\frac{2\pi n}{T}t\,dt = \frac{2}{T}\sum_{t=0}^{T-\Delta t}{}_{\Delta t}\; f(t)\int_t^{t+\Delta t}\cos\frac{2\pi n}{T}t\,dt \qquad (5.9\text{-}7)$$

$$f(t)\int_t^{t+\Delta t}\cos\frac{2\pi n}{T}t\,dt = \frac{f(t)T}{2\pi n}\sin\frac{2\pi n}{T}t\,\Big|_t^{t+\Delta t} \qquad (5.9\text{-}8)$$

$$f(t)\int_t^{t+\Delta t}\cos\frac{2\pi n}{T}t\,dt = \frac{f(t)T}{2\pi n}[\sin\frac{2\pi n}{T}(t+\Delta t) - \sin\frac{2\pi n}{T}t) \qquad (5.9\text{-}9)$$

Substituting Eq 5.9-9 into Eq 5.9-7

$$a_n = \frac{2}{T}\sum_{t=0}^{T-\Delta t}{}_{\Delta t}\;\frac{f(t)T}{2\pi n}[\sin\frac{2\pi n}{T}(t+\Delta t) - \sin\frac{2\pi n}{T}t] \qquad (5.9\text{-}10)$$

Simplifying

$$a_n = \frac{1}{\pi n}\sum_{t=0}^{T-\Delta t}{}_{\Delta t}\; f(t)[\sin\frac{2\pi n}{T}(t+\Delta t) - \sin\frac{2\pi n}{T}t] \qquad (5.9\text{-}11)$$

$$a_n = \frac{1}{\pi n}\sum_{t=0}^{T-\Delta t}{}_{\Delta t}\; f(t)\Delta\sin\frac{2\pi n}{T}t \qquad (5.9\text{-}12)$$

where

$\Delta$ = Difference operator

Changing the form of Eq 5.9-12 to its discrete integral form

$$\mathbf{a_n = \frac{1}{\pi n}{}_{\Delta t}\int_0^T f(t)\Delta\sin\frac{2\pi n}{T}t = \frac{1}{\pi n}\sum_{t=0}^{T-\Delta t}{}_{\Delta t}\; f(t)\Delta\sin\frac{2\pi n}{T}t} \qquad (5.9\text{-}13)$$

**where**

**$n = 1, 2, 3, \ldots$**

**$f(t)$ = discrete function of t**

**$t = 0, \Delta t, 2\Delta t, 3\Delta t, \ldots, T-\Delta t, T$**

**$\Delta t$ = interval between successive values of t**

**$\Delta$ = Difference operator**

**$\Delta g(t) = g(t+\Delta t) - g(t)$**

**T = waveform period**

$$\sin\frac{2\pi nt}{T} = \frac{e^{j\frac{2\pi nt}{T}} - e^{-j\frac{2\pi nt}{T}}}{2j} \qquad (5.9\text{-}14)$$

Substituting Eq 5.9-14 into Eq 5.9-13

The complex form of Eq 5.9-13 is as follows:

$$\mathbf{a_n = \frac{1}{2\pi nj}\left[\int_{0\;\Delta t}^{T} f(t)\Delta e^{j\frac{2\pi nt}{T}} - \int_{0\;\Delta t}^{T} f(t)\Delta e^{-j\frac{2\pi nt}{T}}\right]} \qquad (5.9\text{-}15)$$

**where**

**$n = 1, 2, 3, \ldots$**

**$f(t)$ = discrete function of t**

**$t = 0, \Delta t, 2\Delta t, 3\Delta t, \ldots, T-\Delta t, T$**

**$\Delta t$ = interval between successive values of t**

**$\Delta$ = Difference operator**

**T = waveform period**

Using Interval Calculus identities to facilitate the integrations in Eq 5.9-13 and Eq 5.9-15

Writing the Interval Calculus identity for $\sin\frac{2\pi nt}{T}$

$$\sin\frac{2\pi nt}{T} = e_{\Delta t}\left(\frac{\cos\frac{2\pi n}{T}\Delta t - 1}{\Delta t}, t\right)\sin_{\Delta t}\left(\frac{\tan\frac{2\pi n}{T}\Delta t}{\Delta t}, t\right) \qquad (5.9\text{-}16)$$

Substituting Eq 5.9-16 into Eq 5.9-13

$$\mathbf{a_n = \frac{1}{\pi n}\int_{0\;\Delta t}^{T} f(t)\,\Delta e_{\Delta t}\left(\frac{\cos\frac{2\pi n}{T}\Delta t - 1}{\Delta t}, t\right)\sin_{\Delta t}\left(\frac{\tan\frac{2\pi n}{T}\Delta t}{\Delta t}, t\right)} \qquad (5.9\text{-}17)$$

**where**

**n = 1, 2, 3, …**

**f(t) = discrete function of t**

**t = 0, Δt, 2Δt, 3Δt, …, T-Δt, T**

**Δt = interval between successive values of t**

**Δ = Difference operator**

**T = waveform period**

Writing the Interval Calculus identities for $e^{j\frac{2\pi nt}{T}}$ and $e^{-j\frac{2\pi nt}{T}}$

$$e^{j\frac{2\pi nt}{T}} = e_{\Delta t}\left(\frac{e^{j\frac{2\pi n}{T}\Delta t} - 1}{\Delta t}, t\right) \tag{5.9-18}$$

$$e^{-j\frac{2\pi nt}{T}} = e_{\Delta t}\left(\frac{e^{-j\frac{2\pi n}{T}\Delta t} - 1}{\Delta t}, t\right) \tag{5.9-19}$$

Substituting Eq 5.9-18 and Eq 5.9-19 into Eq 5.9-15

$$\mathbf{a_n = \frac{1}{2\pi nj}\left[\int_{0}^{T} f(t)\Delta e_{\Delta t}\left(\frac{e^{j\frac{2\pi n}{T}\Delta t} - 1}{\Delta t}, t\right)\Delta t - \int_{0}^{T} f(t)\Delta e_{\Delta t}\left(\frac{e^{-j\frac{2\pi n}{T}\Delta t} - 1}{\Delta t}, t\right)\Delta t\right]} \tag{5.9-20}$$

**where**

**n = 1, 2, 3, …**

**f(t) = discrete function of t**

**t = 0, Δt, 2Δt, 3Δt, …, T-Δt, T**

**Δt = interval between successive values of t**

**Δ = Difference operator**

$\left(\frac{e^{j\frac{2\pi n}{T}\Delta t}-1}{\Delta t}\right)$ **= constant**

$\left(\frac{e^{-j\frac{2\pi n}{T}\Delta t}-1}{\Delta t}\right)$ **= constant**

**T = waveform period**

Thus

**The equations to calculate $a_n$ are as follows:**

$$1)\ \ a_n = \frac{1}{\pi n}\int_{0}^{T}{}_{\Delta t} f(t)\Delta \sin\frac{2\pi n t}{T} = \frac{1}{\pi n}\sum_{t=0}^{T-\Delta t}{}_{\Delta t} f(t)\,\Delta \sin\frac{2\pi n}{T}t \qquad (5.9\text{-}21)$$

$$2)\ \ a_n = \frac{1}{2\pi n j}\left[\int_{0}^{T}{}_{\Delta t} f(t)\Delta e^{j\frac{2\pi n t}{T}} - \int_{0}^{T}{}_{\Delta t} f(t)\Delta e^{-j\frac{2\pi n t}{T}}\right] \qquad (5.9\text{-}22)$$

$$3)\ \ a_n = \frac{1}{\pi n}\int_{0}^{T}{}_{\Delta t} f(t)\Delta\left[e_{\Delta t}\left(\frac{\cos\frac{2\pi n}{T}\Delta t - 1}{\Delta t}, t\right)\sin_{\Delta t}\left(\frac{\tan\frac{2\pi n}{T}\Delta t}{\Delta t}, t\right)\right] \qquad (5.9\text{-}23)$$

$$4)\ \ a_n = \frac{1}{2\pi n j}\left[\int_{0}^{T}{}_{\Delta t} f(t)\Delta e_{\Delta t}\left(\frac{e^{j\frac{2\pi n}{T}\Delta t} - 1}{\Delta t}, t\right) - \int_{0}^{T}{}_{\Delta t} f(t)\Delta e_{\Delta t}\left(\frac{e^{-j\frac{2\pi n}{T}\Delta t} - 1}{\Delta t}, t\right)\right] \qquad (5.9\text{-}24)$$

**where**

**$n = 1, 2, 3, \ldots$**

**$\Delta t = \frac{T}{m}$ = interval between successive values of t**

**T = waveform period**

**m = number of intervals within the range, 0 to T**

**m = integer**

**f(t) = discrete function of t**

**$t = 0, \Delta t, 2\Delta t, 3\Delta t, \ldots, T-\Delta t, T$**

**$\left(\frac{e^{j\frac{2\pi n}{T}\Delta t} - 1}{\Delta t}\right)$ = constant**

**$\left(\frac{e^{-j\frac{2\pi n}{T}\Delta t} - 1}{\Delta t}\right)$ = constant**

**$\Delta$ = Difference operator**

**$D_{\Delta t} = \frac{\Delta}{\Delta t}$ = discrete derivative operator**

**$\Delta g(t) = g(t+\Delta t) - g(t)$**

**$D_{\Delta t}\, g(t) = \frac{g(t+\Delta t) - g(t)}{\Delta t}$**

**<u>Comment</u> – The previous equations containing the difference operator, Δ, can be rewritten using the discrete derivative operator, $D_{\Delta t}$. For example, the equation,**

$$\mathbf{a_n = \frac{1}{\pi n}\int_0^T{}_{\Delta t} f(t)\Delta\sin\frac{2\pi nt}{T}}, \text{ can be rewritten, } \mathbf{a_n = \frac{1}{\pi n}\int_0^T{}_{\Delta t} f(t)D_{\Delta t}\sin\frac{2\pi nt}{T}\Delta t}\ .$$

Check Eq 5.9-23

If $\Delta t \to 0$, Eq 5.9-23 should become the Fourier Series equation, $a_n = \frac{2}{T}\int_0^T f(t)\cos\frac{2\pi n}{T}t\,dt$.

From Eq 5.9-23 letting $\Delta t \to 0$

$$a_n = \frac{1}{\pi n}\int_0^T f(t)d[e^{0t}\sin\frac{2\pi n}{T}t] = \frac{1}{\pi n}\int_0^T f(t)d[\sin\frac{2\pi n}{T}t]$$

$$a_n = \frac{2\pi n}{T\pi n}\int_0^T f(t)\frac{2\pi n}{T}\cos\frac{2\pi n}{T}t\,dt = \frac{2}{T}\int_0^T f(t)\cos\frac{2\pi n}{T}t\,dt$$

$$a_n = \frac{2}{T}\int_0^T f(t)\cos\frac{2\pi n}{T}t\,dt$$

The desired result has been obtained

Good check

<u>Find $b_n$ , n = 0, 1, 2, 3, …</u>

$$b_n = \frac{2}{T}\int_0^T f(t)\sin\frac{2\pi n}{T}t\,dt = \frac{2}{T}\sum_{t=0}^{T-\Delta t}{}_{\Delta t}\ f(t)\int_t^{t+\Delta t}\sin\frac{2\pi n}{T}t\,dt \tag{5.9-25}$$

$$f(t)\int_t^{t+\Delta t}\sin\frac{2\pi n}{T}t\,dt = -\frac{f(t)T}{2\pi n}\cos\frac{2\pi n}{T}t\,\Big|_t^{t+\Delta t} \tag{5.9-26}$$

$$f(t)\int_{t}^{t+\Delta t}\sin\frac{2\pi n}{T}t\,dt = -\frac{f(t)T}{2\pi n}\left[\cos\frac{2\pi n}{T}(t+\Delta t) - \cos\frac{2\pi n}{T}t\right) \qquad (5.9\text{-}27)$$

Substituting Eq 5.9-27 into Eq 5.9-25

$$b_n = \frac{2}{T}\sum_{\Delta t}^{T-\Delta t}{}_{t=0} -\frac{f(t)T}{2\pi n}\left[\cos\frac{2\pi n}{T}(t+\Delta t) - \cos\frac{2\pi n}{T}t\right] \qquad (5.9\text{-}28)$$

Simplifying

$$b_n = -\frac{1}{\pi n}\sum_{\Delta t}^{T-\Delta t}{}_{t=0} f(t)\left[\cos\frac{2\pi n}{T}(t+\Delta t) - \cos\frac{2\pi n}{T}t\right] \qquad (5.9\text{-}29)$$

$$b_n = -\frac{1}{\pi n}\sum_{\Delta t}^{T-\Delta t}{}_{t=0} f(t)\,\Delta\cos\frac{2\pi n}{T}t \qquad (5.9\text{-}30)$$

Changing the form of Eq 5.9-30 to its discrete integral form

$$\mathbf{b_n = -\frac{1}{\pi n}\int_{\Delta t\,0}^{T} f(t)\,\Delta\cos\frac{2\pi n}{T}t \;=\; -\frac{1}{\pi n}\sum_{\Delta t}^{T-\Delta t}{}_{t=0} f(t)\,\Delta\cos\frac{2\pi n}{T}t} \qquad (5.9\text{-}31)$$

**where**

**n = 1, 2, 3, …**
**f(t) = discrete function of t**
**t = 0, Δt, 2Δt, 3Δt, …, T-Δt, T**
**Δt = interval between successive values of t**
**Δ = Difference operator**
**T = waveform period**

$$\cos\frac{2\pi nt}{T} = \frac{e^{j\frac{2\pi nt}{T}} + e^{-j\frac{2\pi nt}{T}}}{2} \qquad (5.9\text{-}32)$$

Substituting Eq 5.9-32 into Eq 5.9-31

The complex form of Eq 5.9-31 is as follows:

$$\mathbf{b_n = -\frac{1}{2\pi n}\left[\int_{\Delta t\,0}^{T} f(t)\Delta e^{j\frac{2\pi nt}{T}} + \int_{\Delta t\,0}^{T} f(t)\Delta e^{-j\frac{2\pi nt}{T}}\right]} \tag{5.9-33}$$

**where**

**n = 1, 2, 3, …**

**f(t) = discrete function of t**

**t = 0, Δt, 2Δt, 3Δt, …, T-Δt, T**

**Δt = interval between successive values of t**

**Δ = Difference operator**

**T = waveform period**

Using Interval Calculus identities to facilitate the integrations in Eq 5.9-31 and Eq 5.9-33

Writing the Interval Calculus identity for $\cos\frac{2\pi nt}{T}$

$$\cos\frac{2\pi nt}{T} = e_{\Delta t}\left(\frac{\cos\frac{2\pi n}{T}\Delta t - 1}{\Delta t}, t\right)\cos_{\Delta t}\left(\frac{\tan\frac{2\pi n}{T}\Delta t}{\Delta t}, t\right) \tag{5.9-34}$$

Substituting Eq 5.9-34 into Eq 5.9-31

$$\mathbf{b_n = -\frac{1}{\pi n}\int_{\Delta t\,0}^{T} f(t)\,\Delta e_{\Delta t}\left(\frac{\cos\frac{2\pi n}{T}\Delta t - 1}{\Delta t}, t\right)\cos_{\Delta t}\left(\frac{\tan\frac{2\pi n}{T}\Delta t}{\Delta t}, t\right)} \tag{5.9-35}$$

**where**

**n = 1, 2, 3, …**

**f(t) = discrete function of t**

**t = 0, Δt, 2Δt, 3Δt, …, T-Δt, T**

**Δt = interval between successive values of t**

**Δ = Difference operator**

**T = waveform period**

Writing the Interval Calculus identities for $e^{j\frac{2\pi nt}{T}}$ and $e^{-j\frac{2\pi nt}{T}}$

$$e^{j\frac{2\pi nt}{T}} = e_{\Delta t}\left(\frac{e^{j\frac{2\pi n}{T}\Delta t} - 1}{\Delta t}, t\right) \tag{5.9-36}$$

$$e^{-j\frac{2\pi nt}{T}} = e_{\Delta t}\left(\frac{e^{-j\frac{2\pi n}{T}\Delta t} - 1}{\Delta t}, t\right) \tag{5.9-37}$$

Substituting Eq 5.9-36 and Eq 5.9-37 into Eq 5.9-33

$$\mathbf{b_n = -\frac{1}{2\pi nj}\left[\int_{0\ \Delta t}^{T} f(t)\Delta e_{\Delta t}\left(\frac{e^{j\frac{2\pi n}{T}\Delta t} - 1}{\Delta t}, t\right)\Delta t + \int_{0\ \Delta t}^{T} f(t)\Delta e_{\Delta t}\left(\frac{e^{-j\frac{2\pi n}{T}\Delta t} - 1}{\Delta t}, t\right)\Delta t\right]} \tag{5.9-38}$$

**where**

**n = 1, 2, 3, …**

**f(t) = discrete function of t**

**t = 0, Δt, 2Δt, 3Δt, …, T-Δt, T**

**Δt = interval between successive values of t**

**Δ = Difference operator**

$$\left(\frac{e^{j\frac{2\pi n}{T}\Delta t} - 1}{\Delta t}\right) = \text{constant}$$

$$\left(\frac{e^{-j\frac{2\pi n}{T}\Delta t} - 1}{\Delta t}\right) = \text{constant}$$

**T = waveform period**

Thus

**The equations to calculate $b_n$ are as follows:**

$$\mathbf{1)\ \ b_n = -\frac{1}{\pi n}\int_{0\ \Delta t}^{T} f(t)\Delta\cos\frac{2\pi nt}{T} = -\frac{1}{\pi n}\sum_{t=0\ \Delta t}^{T} f(t)\Delta\cos\frac{2\pi nt}{T}} \tag{5.9-39}$$

$$\mathbf{2)\ \ b_n = -\frac{1}{2\pi n}\left[\int_{0\ \Delta t}^{T} f(t)\Delta e^{j\frac{2\pi nt}{T}} + \int_{0\ \Delta t}^{T} f(t)\Delta e^{-j\frac{2\pi nt}{T}}\right]} \tag{5.9-40}$$

$$\mathbf{3)\ \ b_n = -\frac{1}{\pi n}\int_{0\ \Delta t}^{T} f(t)\Delta\left[e_{\Delta t}\left(\frac{\cos\frac{2\pi n}{T}\Delta t - 1}{\Delta t}, t\right)\cos_{\Delta t}\left(\frac{\tan\frac{2\pi n}{T}\Delta t}{\Delta t}, t\right)\right]} \tag{5.9-41}$$

$$4)\ \ \mathbf{b_n = -\frac{1}{2\pi n}\left[\ \int_{\Delta t\ 0}^{T} f(t)\Delta e_{\Delta t}\left(\frac{e^{j\frac{2\pi n}{T}\Delta t}-1}{\Delta t}, t\right) + \int_{\Delta t\ 0}^{T} f(t)\Delta e_{\Delta t}\left(\frac{e^{-j\frac{2\pi n}{T}\Delta t}-1}{\Delta t}, t\right)\right]} \qquad (5.9\text{-}42)$$

**where**

**n = 1, 2, 3, …**

$\mathbf{\Delta t = \frac{T}{m}}$ **= interval between successive values of t**

**T = waveform period**

**m = number of intervals within the range, 0 to T**

**m = integer**

**f(t) = discrete function of t**

**t = 0, Δt, 2Δt, 3Δt, …, T-Δt, T**

$\mathbf{\left(\frac{e^{j\frac{2\pi n}{T}\Delta t}-1}{\Delta t}\right)}$ **= constant**

$\mathbf{\left(\frac{e^{-j\frac{2\pi n}{T}\Delta t}-1}{\Delta t}\right)}$ **= constant**

**Δ = Difference operator**

$\mathbf{D_{\Delta t} = \frac{\Delta}{\Delta t}}$ **= discrete derivative**

**Δg(t) = g(t+Δt) – g(t)**

$\mathbf{D_{\Delta t}\, g(t) = \frac{g(t+\Delta t) - g(t)}{\Delta t}}$

**<u>Comment</u> – The previous equations containing the difference operator, Δ, can be rewritten using the discrete derivative operator, $\mathbf{D_{\Delta t}}$. For example, the equation,**

$$\mathbf{b_n = -\frac{1}{\pi n}\int_{\Delta t\ 0}^{T} f(t)\Delta\cos\frac{2\pi nt}{T}}, \textbf{ can be rewritten, } \mathbf{b_n = -\frac{1}{\pi n}\int_{\Delta t\ 0}^{T} f(t)D_{\Delta t}\cos\frac{2\pi nt}{T}\Delta t}\ .$$

Check Eq 5.9-41

If $\Delta t \to 0$, Eq 5.9-41 should become the Fourier Series equation, $b_n = \frac{2}{T}\int_0^T f(t)\sin\frac{2\pi n}{T}t\ dt$.

From Eq 5.9-41 letting $\Delta t \to 0$

$$b_n = -\frac{1}{\pi n}\int_0^T f(t)d[e^{0t}\cos\frac{2\pi n}{T}t] = -\frac{1}{\pi n}\int_0^T f(t)d[\cos\frac{2\pi n}{T}t]$$

$$b_n = \frac{2\pi n}{T\pi n}\int_0^T f(t)\frac{2\pi n}{T}\sin\frac{2\pi n}{T}t\,dt = \frac{2}{T}\int_0^T f(t)\sin\frac{2\pi n}{T}t\,dt$$

$$b_n = \frac{2}{T}\int_0^T f(t)\sin\frac{2\pi n}{T}t\,dt$$

The desired result has been obtained

Good check

The discrete integration process involved in the use of the above derived equations is facilitated by the tables of discrete integrals, discrete derivatives, and identities in the Appendix.

Thus

**The Fourier Series for sample and hold shaped functions, f(t), is:**

**Diagram 5.9-2 - An example of a discrete Interval Calculus sample and hold shaped waveform**

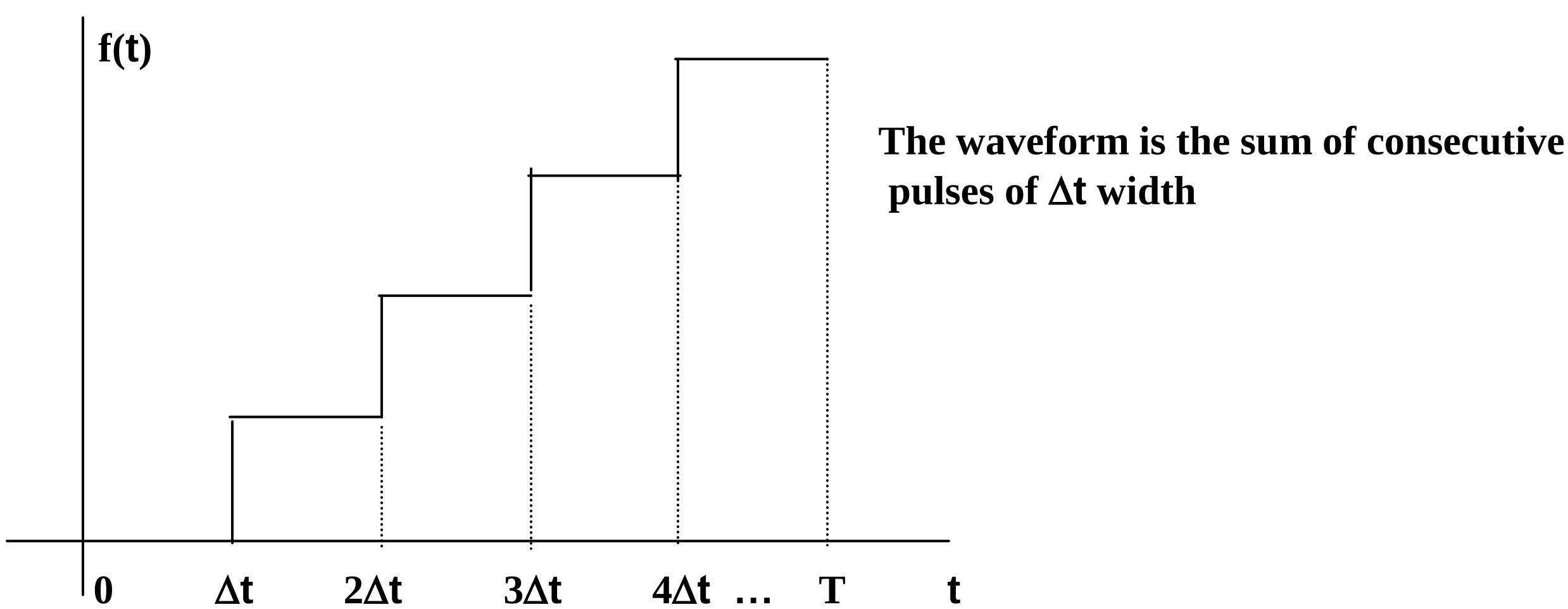

**Fourier Series of Sample and Hold Shaped Waveforms**

**1)** $$\mathbf{F(t) = \frac{a_0}{2} + \sum_{n=1}^{\infty} \left( a_n\cos\frac{2\pi n}{T}t + b_n\sin\frac{2\pi n}{T}t \right)} \qquad (5.9\text{-}43)$$

**Fourier Series Coefficient Calculation Equations Generalized for Sample and Hold Shaped Waveforms**

2) $$a_0 = \frac{2}{T}\,{}_{\Delta t}\!\int_0^T f(t)\Delta t \qquad (5.9\text{-}44)$$

3) $$a_n = \frac{1}{\pi n}\,{}_{\Delta t}\!\int_0^T f(t)\Delta\sin\frac{2\pi nt}{T} = \frac{1}{\pi n}\,{}_{\Delta t}\!\sum_{t=0}^{T-\Delta t} f(t)\,\Delta\sin\frac{2\pi n}{T}t \qquad (5.9\text{-}45)$$

4) $$a_n = \frac{1}{2\pi nj}\left[\,{}_{\Delta t}\!\int_0^T f(t)\Delta e^{j\frac{2\pi nt}{T}} - {}_{\Delta t}\!\int_0^T f(t)\Delta e^{-j\frac{2\pi nt}{T}}\right] \qquad (5.9\text{-}46)$$

5) $$a_n = \frac{1}{\pi n}\,{}_{\Delta t}\!\int_0^T f(t)\Delta\left[e_{\Delta t}\left(\frac{\cos\frac{2\pi n}{T}\Delta t - 1}{\Delta t}, t\right)\sin_{\Delta t}\left(\frac{\tan\frac{2\pi n}{T}\Delta t}{\Delta t}, t\right)\right] \qquad (5.9\text{-}47)$$

6) $$a_n = \frac{1}{2\pi nj}\left[\,{}_{\Delta t}\!\int_0^T f(t)\Delta e_{\Delta t}\left(\frac{e^{j\frac{2\pi n}{T}\Delta t} - 1}{\Delta t}, t\right) - {}_{\Delta t}\!\int_0^T f(t)\Delta e_{\Delta t}\left(\frac{e^{-j\frac{2\pi n}{T}\Delta t} - 1}{\Delta t}, t\right)\right] \qquad (5.9\text{-}48)$$

7) $$b_n = -\frac{1}{\pi n}\,{}_{\Delta t}\!\int_0^T f(t)\Delta\cos\frac{2\pi nt}{T} = -\frac{1}{\pi n}\,{}_{\Delta t}\!\sum_{t=0}^{T-\Delta t} f(t)\Delta\cos\frac{2\pi nt}{T} \qquad (5.9\text{-}49)$$

8) $$b_n = -\frac{1}{2\pi n}\left[\,{}_{\Delta t}\!\int_0^T f(t)\Delta e^{j\frac{2\pi nt}{T}} + {}_{\Delta t}\!\int_0^T f(t)\Delta e^{-j\frac{2\pi nt}{T}}\right] \qquad (5.9\text{-}50)$$

9) $$b_n = -\frac{1}{\pi n}\,{}_{\Delta t}\!\int_0^T f(t)\Delta\left[e_{\Delta t}\left(\frac{\cos\frac{2\pi n}{T}\Delta t - 1}{\Delta t}, t\right)\cos_{\Delta t}\left(\frac{\tan\frac{2\pi n}{T}\Delta t}{\Delta t}, t\right)\right] \qquad (5.9\text{-}51)$$

10) $$b_n = -\frac{1}{2\pi n}\left[\,{}_{\Delta t}\!\int_0^T f(t)\Delta e_{\Delta t}\left(\frac{e^{j\frac{2\pi n}{T}\Delta t} - 1}{\Delta t}, t\right) + {}_{\Delta t}\!\int_0^T f(t)\Delta e_{\Delta t}\left(\frac{e^{-j\frac{2\pi n}{T}\Delta t} - 1}{\Delta t}, t\right)\right] \qquad (5.9\text{-}52)$$

**where**

$t = 0, \Delta t, 2\Delta t, 3\Delta t, \ldots, T-\Delta t, T$

$\Delta t = \frac{T}{m}$ = **interval between successive values of t**

$n = 1, 2, 3, \ldots$

**T = waveform period**

**m = number of intervals within the range, 0 to T**

**m = positive integer**

**f(t) = discrete sample and hold shaped waveform function of t**

$$\left(\frac{e^{j\frac{2\pi n}{T}\Delta t} - 1}{\Delta t}\right) = \text{constant}$$

$$\left(\frac{e^{-j\frac{2\pi n}{T}\Delta t} - 1}{\Delta t}\right) = \text{constant}$$

**Δ = Difference operator**

$D_{\Delta t} = \frac{\Delta}{\Delta t}$ = **discrete derivative**

$\Delta g(t) = g(t+\Delta t) - g(t)$

$$D_{\Delta t}\, g(t) = \frac{g(t+\Delta t) - g(t)}{\Delta t}$$

**<u>Comment</u> – The above equations containing the difference operator, Δ, can be rewritten using the discrete derivative operator, $D_{\Delta t}$. For example, the equation,**

$$b_n = -\frac{1}{\pi n}\ {}_{\Delta t}\!\int_0^T f(t)\Delta\cos\frac{2\pi nt}{T}\text{, can be rewritten, } b_n = -\frac{1}{\pi n}\ {}_{\Delta t}\!\int_0^T f(t)D_{\Delta t}\cos\frac{2\pi nt}{T}\Delta t\ .$$

For a demonstration of the use of the previously derived Fourier Series equations that have been generalized for use with sample and hold shaped waveforms. See Example 5.9-1 on the following page.

Example 5.9-1 – The expansion of a sample and hold shaped waveform into a Fourier Series

Find the Fourier Series for the following sample and hold shaped waveform. Use the Fourier Series coefficient calculation equations generalized for use with sample and hold shaped waveforms.

Diagram 5.9-3 - A sample and hold shaped waveform where f(t) = t and t = 0, 1, 2

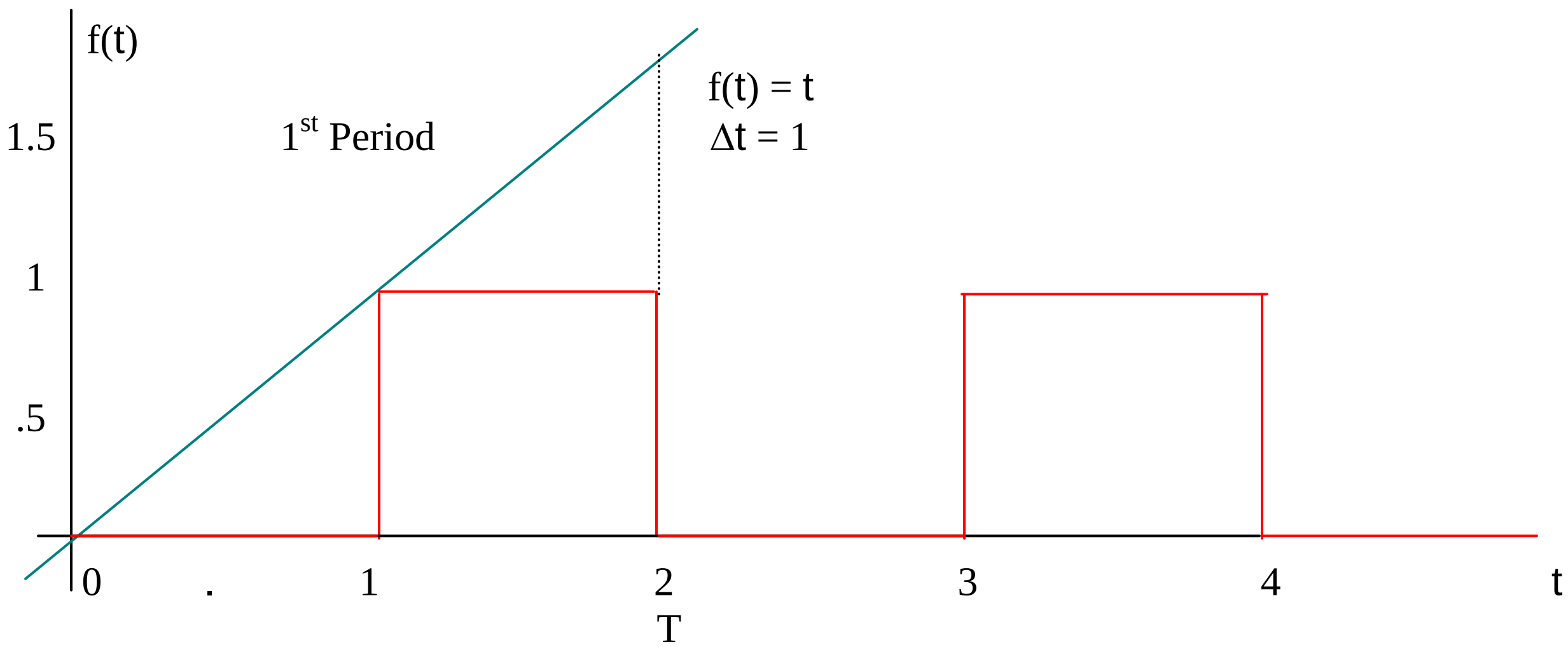

Periodic Sample and Hold Shaped Waveform

$f(t) = t$ (5.9-53)

$t = 0, 1, 2, \ldots$ (5.9-54)

$\Delta t = 1$ , the sampling interval (5.9-55)

$T = 2$ (5.9-56)

$$F(t) = \frac{a_0}{2} + \sum_{n=1}^{\infty} \left( a_n \cos\frac{2\pi n}{T}\, t + b_n \sin\frac{2\pi n}{T}\, t \right) \quad (5.9\text{-}57)$$

$$a_0 = \frac{2}{T}\ {}_{\Delta t}\!\int_0^T f(t)\Delta t \quad (5.9\text{-}58)$$

$$a_n = \frac{1}{\pi n}\ {}_{\Delta t}\!\int_0^T f(t)\Delta\left[e_{\Delta t}\left(\frac{\cos\frac{2\pi n}{T}\Delta t - 1}{\Delta t}, t\right)\sin_{\Delta t}\left(\frac{\tan\frac{2\pi n}{T}\Delta t}{\Delta t}, t\right)\right] \quad (5.9\text{-}59)$$

$$b_n = -\frac{1}{\pi n}\ {}_{\Delta t}\!\int_0^T f(t)\Delta\left[e_{\Delta t}\left(\frac{\cos\frac{2\pi n}{T}\Delta t - 1}{\Delta t}, t\right)\cos_{\Delta t}\left(\frac{\tan\frac{2\pi n}{T}\Delta t}{\Delta t}, t\right)\right] \quad (5.9\text{-}60)$$

Find $a_0$ for the sample and hold shaped waveform of Diagram 5.9-3

From Eq 5.9-53 thru 5.9-56 and Eq 5.9-58

$$a_0 = \frac{2}{T}\, {}_{\Delta t}\!\int_0^T f(t)\Delta t = \frac{2}{2}\, {}_1\!\int_0^2 t\ \Delta t = \frac{t(t\text{-}\Delta t)}{2}\Big|_0^2 1 = \frac{t(t\text{-}1)}{2}\Big|_0^2 = \frac{2(1)}{2} - 0 = 1 \qquad (5.9\text{-}61)$$

$$a_0 = 1 \qquad (5.9\text{-}62)$$

Checking

Referring to Diagram 5.9-3

$a_0 = \frac{2}{2}(1)(1) = 1$ Good check

Find $b_0$ for the sample and hold shaped waveform of Diagram 5.9-3

From Eq 5.9-53 thru 5.9-56 and Eq 5.9-60

$$b_n = -\frac{1}{\pi n}\, {}_{\Delta t}\!\int_0^T f(t)\Delta[e_{\Delta t}(\frac{\cos\frac{2\pi n}{T}\Delta t - 1}{\Delta t}, t)\cos_{\Delta t}(\frac{\tan\frac{2\pi n}{T}\Delta t}{\Delta t}, t)] \qquad (5.9\text{-}63)$$

$n = 1, 2, 3, \ldots$
$f(t) = t$
$T = 2$
$\Delta t = 1$

$$b_n = -\frac{1}{\pi n}\, {}_{\Delta t}\!\int_0^2 t\,\Delta[e_1(\cos\pi n - 1, t)\cos_1(\tan\pi n, t)] \qquad (5.9\text{-}64)$$

$$b_n = -\frac{1}{\pi n}\, {}_{\Delta t}\!\int_0^2 t\, D_1[e_1(\cos\pi n - 1, t)\cos_1(\tan\pi n, t)]\,\Delta t \qquad (5.9\text{-}65)$$

For $n = 2, 4, 6, \ldots$

$$b_n = -\frac{1}{\pi n}\, {}_{\Delta t}\!\int_0^2 t\, D_1[e_1(0, t)\cos_1(0, t)\}\,\Delta t \qquad (5.9\text{-}66)$$

$$e_1(0,t) = [1+0(1)]^{\frac{t}{1}} = 1 \qquad (5.9\text{-}67)$$

$$\cos_1(0,t) = \frac{e_1(j0,t) + e_1(-j0,t)}{2} = \frac{1+1}{2} = 1 \qquad (5.9\text{-}68)$$

Substituting Eq 5.9-67 and Eq 5.9-68 into Eq 5.9-66

$$b_n = -\frac{1}{\pi n}\ {}_{\Delta t}\!\int_0^2 t\, D_1[1]\,\Delta t = -\frac{1}{\pi n}\ {}_{\Delta t}\!\int_0^2 t(0)\Delta t = 0 \qquad (5.9\text{-}69)$$

$$\mathbf{b_n = 0\ ,\quad for\ n = 2, 4, 6, \ldots} \qquad (5.9\text{-}70)$$

From Eq 5.9-65

For n = 1, 3, 5, …

$$b_n = -\frac{1}{\pi n}\ {}_{\Delta t}\!\int_0^2 t\, D_1[e_1(-2, t)\cos_1(0, t)]\,\Delta t \qquad (5.9\text{-}71)$$

Substituting Eq 5.9-68 into Eq 5.5-71

$$b_n = -\frac{1}{\pi n}\ {}_{\Delta t}\!\int_0^2 t\, D_1[e_1(-2, t)]\,\Delta t \qquad (5.9\text{-}72)$$

$$b_n = -\frac{1}{\pi n}\ {}_{\Delta t}\!\int_0^2 t\,(-2)[e_1(-2, t)]\,\Delta t \qquad (5.9\text{-}73)$$

$$b_n = \frac{2}{\pi n}\ {}_{\Delta t}\!\int_0^2 t\, e_1(-2, t)\,\Delta t \qquad (5.9\text{-}74)$$

Integrating Eq 5.9-74 using integration by parts

$${}_{\Delta t}\!\int_{t_1}^{t_2} v(t)D_{\Delta t}u(t)\Delta t = u(t)v(t)\Big|_{t_1}^{t_2} - {}_{\Delta t}\!\int_{t_1}^{t_2} D_{\Delta t}v(t)u(t+\Delta t)\Delta t \qquad (5.9\text{-}75)$$

$$\Delta t = 1 \qquad (5.9\text{-}76)$$

let

$v(t) = t \qquad D_1 v(t) = 1$

$D_1 u(t) = e_1(-2,t) \qquad u(t) = \frac{e_1(-2,t)}{-2}$

$$u(t+\Delta t) = \frac{e_1(-2,t+1)}{-2}$$

From Eq 5.9-74 thru Eq 5.9-76 and the above equations

$$b_n = \frac{2}{\pi n} \int_0^2 t\, e_1(-2,t)\, \Delta t = \frac{2}{\pi n} \frac{t e_1(-2,t)}{-2} \Big|_0^2 - \frac{2}{\pi n} \int_0^2 (1) \frac{e_1(-2,t+1)}{-2} \Delta t \qquad (5.9\text{-}77)$$

Changing variable and simplifying

$$\tau = t + 1 \qquad (5.9\text{-}78)$$

$$b_n = \frac{2}{\pi n} \int_0^2 t\, e_1(-2,t)\, \Delta t = \frac{2}{\pi n} \frac{2e_1(-2,2)}{-2} + \frac{1}{\pi n} \int_1^3 e_1(-2,\tau)\, \Delta\tau \qquad (5.9\text{-}79)$$

$$b_n = -\frac{2}{\pi n} [1 - 2(1)]^{\frac{2}{1}} + \frac{1}{\pi n} \frac{e_1(-2,\tau)}{-2} \Big|_1^3 \qquad (5.9\text{-}80)$$

$$b_n = -\frac{2}{\pi n} [1 - 2(1)]^{\frac{2}{1}} + \frac{1}{\pi n} \frac{[1-2(1)]^{\frac{3}{1}} - [1-2(1)]^{\frac{1}{1}}}{-2} \qquad (5.9\text{-}81)$$

$$b_n = -\frac{2}{\pi n} (-1)^2 - \frac{1}{2\pi n} [(-1)^3 - (-1)^1] = -\frac{2}{\pi n} - \frac{1}{2\pi n} [-1 + 1] = -\frac{2}{\pi n} \qquad (5.9\text{-}82)$$

$$\mathbf{b_n = -\frac{2}{\pi n}\ ,\quad for\ n = 1, 3, 5, \ldots} \qquad (5.9\text{-}83)$$

Then

$$\mathbf{b_n} = \begin{cases} \mathbf{0} & \textbf{for n = 2 4 6 …} \\ -\mathbf{\frac{2}{\pi n}} & \textbf{for n = 1 3 5 …} \end{cases} \qquad (5.9\text{-}84)$$

Checking

$$b_n = \frac{2}{T}\int_0^T f(t)\sin\frac{2\pi n}{T}t\,dt$$

n = 1, 2, 3, …
T = 2
f(t) = 0 for $0 \le t < 1$
f(t) = 1 for $1 \le t < 2$

Substituting

$$b_n = \frac{2}{2}\int_0^1 (0)\sin\frac{2\pi n}{2}t\,dt + \frac{2}{2}\int_1^2 (1)\sin\frac{2\pi n}{2}t\,dt$$

$$b_n = -\frac{1}{\pi n}\cos\pi nt\,\Big|_1^2 = -\frac{1}{\pi n}[\cos 2\pi n - \cos\pi n] = -\frac{1}{\pi n}[1 - \cos\pi n]$$

$$b_n = -\frac{1}{\pi n}[1 - \cos\pi n]$$

Then

For n = 1, 2, 3, 4, …

$$b_n = \begin{cases} 0 & \text{for } n = 2\ \ 4\ \ 6\ \ldots \\ -\dfrac{2}{\pi n} & \text{for } n = 1\ \ 3\ \ 5\ \ldots \end{cases} \qquad \text{Good check}$$

<u>Find $a_n$ for the sample and hold shaped waveform of Diagram 5.9-3</u>

$$a_n = \frac{1}{\pi n}\ {}_{\Delta t}\!\!\int_0^T f(t)\Delta[e_{\Delta t}(\frac{\cos\frac{2\pi n}{T}\Delta t - 1}{\Delta t}, t)\sin_{\Delta t}(\frac{\tan\frac{2\pi n}{T}\Delta t}{\Delta t}, t)] \qquad (5.9\text{-}85)$$

n = 1, 2, 3, …
f(t) = t
T = 2

$\Delta t = 1$

$$a_n = \frac{1}{\pi n}\,{}_1\!\int_0^2 t\,\Delta[e_1(\cos\pi n - 1, t)\sin_1(\tan\pi n, t)] \qquad (5.9\text{-}86)$$

$$a_n = \frac{1}{\pi n}\,{}_1\!\int_0^2 t\,D_1[e_1(\cos\pi n - 1, t)\sin_1(\tan\pi n, t)]\,\Delta t \qquad (5.9\text{-}87)$$

where

$D_1 = \frac{\Delta}{\Delta t}$, the discrete derivative operator

For n = 1, 2, 3, 4, …

$$\tan\pi n = 0 \qquad (5.9\text{-}88)$$

$$\sin_{\Delta t}(0,t) = 0 \qquad (5.9\text{-}89)$$

From Eq 5.9-87 thru Eq 5.9-89)

$$a_n = \frac{1}{\pi n}\,{}_1\!\int_0^2 t\,D_1[0]]\,\Delta t = 0 \qquad (5.9\text{-}90)$$

$$a_n = 0 \qquad (5.9\text{-}91)$$

Then

$$\mathbf{a_n = 0\ ,\ \ for\ n = 1, 2, 3, 4, \ldots} \qquad (5.9\text{-}92)$$

Checking

$$a_n = \frac{2}{T}\int_0^T f(t)\cos\frac{2\pi n}{T}t\,dt$$

$n = 1, 2, 3, 4, \ldots$
$T = 2$
$f(t) = 0 \quad \text{for } 0 \le t < 1$
$f(t) = 1 \quad \text{for } 1 \le t < 2$

Substituting

$$a_n = \frac{2}{2}\int_0^1 (0)\cos\frac{2\pi n}{2}t\,dt + \frac{2}{2}\int_1^2 (1)\cos\frac{2\pi n}{2}t\,dt$$

$$a_n = \frac{1}{\pi n}\sin\pi nt\,\Big|_1^2 = \frac{1}{\pi n}[\sin 2\pi n - \sin\pi n]$$

Then

For n = 1, 2, 3, 4, …

$a_n = 0$

Good check

Thus

**The Fourier Series Expansion of the sample and hold shaped waveform of Diagram 5.9-3 is:**

$$\mathbf{F(t) = \frac{a_0}{2} + \sum_{n=1}^{\infty}\left( a_n\cos\frac{2\pi n}{T}\,t + b_n\sin\frac{2\pi n}{T}\,t \right)} \qquad (5.9\text{-}93)$$

**where**

**F(t) = the sample and hold shaped waveform of Diagram 5.9-3**

$\mathbf{a_0 = 1}$

$\mathbf{a_n = 0}$ **for n = 1, 2, 3, …**

$$\mathbf{b_n = \begin{cases} 0 & \text{for n = 2 4 6 …} \\ -\dfrac{2}{\pi n} & \text{for n = 1 3 5 …} \end{cases}}$$

This section, Section 5.9, has shown that discrete sample and hold shaped waveforms can be expanded into Fourier Series using generalized Fourier Series coefficient calculation equations. Interval Calculus discrete integration is successfully employed to calculate the necessary series coefficients. In the following section, Interval Calculus discrete integration will again be used but this time to derive and evaluate convolution integration equations for the $K_{\Delta t}$ Transform and the Z Transform. As with the Laplace Transform Convolution Equation, the $K_{\Delta t}$ Transform Convolution Equation and the Z Transform Convolution Equation can be put to good mathematical use.

## Section 5.10: The derivation of the $K_{\Delta t}$ Transform Convolution Equation and the Z Transform Convolution Equation

Of considerable importance in Operational Calculus is the Laplace Transform Convolution Equation which is shown below.

The Laplace Transform Convolution Equation is:

$$L\left[\int_0^t f(t)g(t-\lambda)\Delta\lambda\right] = L\left[\int_0^t f(t-\lambda)g(\lambda)\Delta\lambda\right] = L[f(t)]\,L[g(t)] \quad , \quad \text{Laplace Convolution Equation} \qquad (5.10\text{-}1)$$

where

$f(t),g(t)$ = functions of t

$f(t-\lambda)$ = function of t and $\lambda$

$L[h(t)] = \int_0^\infty e^{-st}\,h(t)dt$ , The Laplace Transform of a function of t

$L[f(t)]$ , $L[g(t)]$ = Laplace Transforms, functions of s

The Laplace Convolution Equation is used in mathematics for many important and useful calculations. This would include obtaining inverse Laplace Transforms, finding solutions to differential equations, and finding the output of Laplace Transform transfer functions. The considerable mathematical capabilities of the Laplace Transform Convolution Equations are available to discrete mathematics $K_{\Delta t}$ Transforms and Z Transforms also. $K_{\Delta t}$ Transform and Z Transform Convolution Equations, generalizations of the Laplace Transform Equation, do exist.They are presented below.

The $K_{\Delta t}$ Transform Convolution Equation is:

$$K_{\Delta t}\left[{}_{\Delta\lambda}\int_0^t f(t)g(t-\lambda-\Delta\lambda)\Delta\lambda\right] = K_{\Delta t}\left[{}_{\Delta\lambda}\int_0^t f(t-\lambda-\Delta\lambda)g(\lambda)\Delta\lambda\right] = K_{\Delta t}[f(t)]\,K_{\Delta t}[g(t)] \qquad (5.10\text{-}2)$$

where

$f(t),g(t)$ = discrete functions of t

$f(t-\lambda-\Delta\lambda)$ = discrete function of t, $\lambda$, and $\Delta\lambda$

$\Delta t = \Delta\lambda$ , Interval between consecutive values of the independent discrete variables, t and $\lambda$

$t = 0, \Delta t, 2\Delta t, 3\Delta t, \ldots$

$\lambda = 0, \Delta\lambda, 2\Delta\lambda, 3\Delta\lambda, \ldots, t-\Delta\lambda, t$

$K_{\Delta t}[h(t)] = {}_{\Delta t}\int_0^\infty (1+s\Delta t)^{-\left(\frac{t+\Delta t}{\Delta t}\right)}h(t)\Delta t$ , The $K_{\Delta t}$ Transform of a function of t

$K_{\Delta t}[f(t)]$ , $K_{\Delta t}[g(t)]$ = $K_{\Delta t}$ Transforms, functions of s

Note that the $K_{\Delta t}$ Transform Convolution Equation becomes the Laplace Transform Convolution Equation for $\Delta\lambda \to 0$. This is to be expected since Interval Calculus becomes Calculus when the independent variable interval between values (in this case $\Delta\lambda$) becomes infinitessimal.

The Z Transform Convolution Equation is:

$$Z[\frac{1}{T}\, {}_T\!\int_0^{t+T} f(t)g(t-\lambda)\Delta\lambda] = Z[\frac{1}{T}\, {}_T\!\int_0^{t+T} f(t-\lambda)g(\lambda)\Delta\lambda] = Z[f(t)]\, Z[g(t)] \qquad (5.10\text{-}3)$$

where

$f(t),g(t)$ = discrete functions of t

$f(t-\lambda)$ = discrete function of t and $\lambda$

$\Delta t = \Delta\lambda = T$ ,   Interval between consecutive values of the independent discrete variables, t and $\lambda$

$t = 0, \Delta t, 2\Delta t, 3\Delta t, \ldots$

$\lambda = 0, \Delta\lambda, 2\Delta\lambda, 3\Delta\lambda, \ldots, t-\Delta\lambda, t, t+T$

$Z[h(t)] = {}_{\Delta t}\!\int_0^{\infty} z^{-(\frac{t}{\Delta t})} h(t)\Delta t$ ,     Z Transform of a function of t

$Z[f(t)]$ , $Z[g(t)]$ = Z Transforms, functions of z

The $K_{\Delta t}$ Transform and Z Transform Convolution Equations presented above provide additional ways to solve problems involving discrete variables (i.e. where $t = 0, \Delta t, 2\Delta t, 3\Delta t, \ldots$). As always, the problem solution methodology selected is at the discretion of the engineer or mathematician. The use of these equations requires the application of  Interval Calculus. The integral used in these equations is the Interval Calculus discrete integral. Only when the value of the discrete variable interval, $\Delta\lambda$, $\Delta t$ or T, becomes an infinitessimal value does this Interval Calculus discrete integral become the commonly used Calculus integral. The Interval Calculus discrete integral is defined as follows:

The Definition of the Interval Calculus Discrete Integral is:

$${}_{\Delta t}\!\int_{t_1}^{t_2} f(t)\Delta t \;\equiv\; \Delta t \, {}_{\Delta t}\!\sum_{t=t_1}^{t_2-\Delta t} f(t) \qquad (5.10\text{-}4)$$

where

$f(t)$ = discrete function of t

$t = t_1, t_1+\Delta t, t_1+2\Delta t, \ldots, t_2-\Delta t, t_2$

$\Delta t$ = Interval between consecutive values of the discrete independent variable, t

To facilitate the calculation of  Interval Calculus discrete integrals, there are evaluation formulas provided in Table 6 in the Appendix.

The derivation of the $K_{\Delta t}$ Transform Convolution Equation, Eq 5.10-2, follows:

<u>Derivation of the $K_{\Delta t}$ Transform Convolution Equation</u>

f(t), g(t) = sample and hold shaped waveform functions

Consider the following discrete function.

$$\int_{0}^{t}{}_{\Delta\lambda} f(t-\lambda-\Delta\lambda)g(\lambda)\Delta\lambda] \;=\; \text{a discrete function of } \lambda \text{ and } t \qquad (5.10\text{-}5)$$

where

$f(t-\lambda-\Delta\lambda)$ = discrete function of t, $\lambda$, and $\Delta\lambda$

$g(\lambda)$ = discrete function of $\lambda$

$\Delta t = \Delta\lambda$

$t = 0, \Delta t, 2\Delta t, 3\Delta t, \ldots$

$\lambda = 0, \Delta\lambda, 2\Delta\lambda, 3\Delta\lambda, \ldots, t-\Delta\lambda, t$

$$K_{\Delta t}[h(t)]\;) = \int_{0}^{\infty}{}_{\Delta t} [1+s\Delta t]^{-\left(\frac{t+\Delta t}{\Delta t}\right)} h(t)\Delta t \;\;,\quad K_{\Delta t} \text{ Transform of the function, } h(t) \qquad (5.10\text{-}6)$$

Take the $K_{\Delta t}$ Transform of Eq 5.10-5

$$K_{\Delta t}\left[\int_{0}^{t}{}_{\Delta\lambda} f(t-\lambda-\Delta\lambda)g(\lambda)\Delta\lambda\right] = \int_{0}^{\infty}{}_{\Delta t} \left[\int_{0}^{t}{}_{\Delta\lambda} f(t-\lambda-\Delta\lambda)g(\lambda)\Delta\lambda\right][1+s\Delta t]^{-\left(\frac{t+\Delta t}{\Delta t}\right)}\Delta t \qquad (5.10\text{-}7)$$

<u>Comment</u> – For discrete integrals, the independent discrete variable range is:
lower limit ≤ variable < upper limit.
In the direction of increasing discrete variable value, the specified discrete integration lower limit is the first value in range and the specified discrete integration upper limit is the first value out of range.

$$\Delta\lambda = \Delta t \qquad (5.10\text{-}8)$$

$$\lambda = 0, \Delta\lambda, 2\Delta\lambda, 3\Delta\lambda, \ldots, t-\Delta\lambda, t \qquad (5.10\text{-}9)$$

$$U(t-\lambda-\Delta\lambda) = \begin{cases} 1 & \text{for } \lambda < t \\ 0 & \text{for } \lambda \ge t \end{cases} \qquad (5.10\text{-}10)$$

Then

$$f(t-\lambda-\Delta\lambda)g(\lambda)U(t-\lambda-\Delta\lambda) = \begin{cases} f(t-\lambda-\Delta\lambda)g(\lambda) & \text{for } \lambda < t \\ 0 & \text{for } \lambda \geq t \end{cases} \tag{5.10-11}$$

Since the product of Eq 5.10-11 is 0 for all $\lambda \geq t$, the inner integration in Eq 5.10-7 can be extended to $\infty$ if the factor $U(t-\lambda-\Delta\lambda)$ is inserted into the integrand.

$$K_{\Delta t}[{}_{\Delta\lambda}\int_0^t f(t-\lambda-\Delta\lambda)g(\lambda)\Delta\lambda] = {}_{\Delta t}\int_0^\infty [{}_{\Delta\lambda}\int_0^\infty f(t-\lambda-\Delta\lambda)g(\lambda)U(t-\lambda-\Delta\lambda)\Delta\lambda][1+s\Delta t]^{-(\frac{t+\Delta t}{\Delta t})}\Delta t \tag{5.10-12}$$

Our usual assumptions about the functions we transform are sufficient to permit the order of integration in Eq 5.10-12 to be interchanged.

$$K_{\Delta t}[{}_{\Delta\lambda}\int_0^t f(t-\lambda-\Delta\lambda)g(\lambda)\Delta\lambda] = {}_{\Delta\lambda}\int_0^\infty [{}_{\Delta t}\int_0^\infty f(t-\lambda-\Delta\lambda)g(\lambda)U(t-\lambda-\Delta\lambda)[1+s\Delta t]^{-(\frac{t+\Delta t}{\Delta t})}\Delta t]\Delta\lambda \tag{5.10-13}$$

$$K_{\Delta t}[{}_{\Delta\lambda}\int_0^t f(t-\lambda-\Delta\lambda)g(\lambda)\Delta\lambda] = {}_{\Delta\lambda}\int_0^\infty g(\lambda)[{}_{\Delta t}\int_0^\infty f(t-\lambda-\Delta\lambda)\, U(t-\lambda-\Delta\lambda)[1+s\Delta t]^{-(\frac{t+\Delta t}{\Delta t})}\Delta t]\Delta\lambda \tag{5.10-14}$$

Due to the term, $U(t-\lambda-\Delta\lambda)$, the integrand of the inner integral of Eq 5.10-14 is 0 for $t < \lambda+\Delta\lambda$

$$K_{\Delta t}[{}_{\Delta\lambda}\int_0^t f(t-\lambda-\Delta\lambda)g(\lambda)\Delta\lambda] = {}_{\Delta\lambda}\int_0^\infty g(\lambda)[{}_{\Delta t}\int_{\lambda+\Delta\lambda}^\infty f(t-\lambda-\Delta\lambda)\, [1+s\Delta t]^{-(\frac{t+\Delta t}{\Delta t})}\Delta t]\Delta\lambda \tag{5.10-15}$$

In the inner integral on the right side of Eq 5.10-15
Let

$$\tau = t-\lambda-\Delta\lambda \tag{5.10-16}$$

From Eq 5.10-16

$$\Delta t = \Delta\tau \tag{5.10-17}$$

$$t = \tau + \lambda + \Delta\lambda \tag{5.10-18}$$

From Eq 5.10-16, finding the limits of τ from the limits of t

For $t = \lambda+\Delta\lambda$ , $\tau = 0$ (5.10-19)

For $t = \infty$ , $\tau = \infty$ (5.10-20)

Substituting Eq 5.10-17 thru Eq 5.10-20 into Eq 5.10-15

$$K_{\Delta t}[{}_{\Delta\lambda}\int_0^t f(t-\lambda-\Delta\lambda)g(\lambda)\Delta\lambda] = {}_{\Delta\lambda}\int_0^\infty g(\lambda)[{}_{\Delta\tau}\int_0^\infty f(\tau)\,[1+s\Delta\tau]^{-\left(\frac{\tau+\lambda+\Delta\lambda+\Delta\tau}{\Delta\tau}\right)}\Delta\tau]\Delta\lambda \tag{5.10-21}$$

Rearranging Eq 5.10-21

$$K_{\Delta t}[{}_{\Delta\lambda}\int_0^t f(t-\lambda-\Delta\lambda)g(\lambda)\Delta\lambda] = {}_{\Delta\lambda}\int_0^\infty g(\lambda)[1+s\Delta\tau]^{-\left(\frac{\lambda+\Delta\lambda}{\Delta\tau}\right)}\,[{}_{\Delta\tau}\int_0^\infty f(\tau)\,[1+s\Delta\tau]^{-\left(\frac{\tau+\Delta\tau}{\Delta\tau}\right)}\Delta\tau]\Delta\lambda \tag{5.10-22}$$

Note From Eq 5.10-8 and Eq 5.10-17

$$\Delta\tau = \Delta\lambda = \Delta t \tag{5.10-23}$$

From Eq 5.10-22 and Eq 5.10-23

$$K_{\Delta t}[{}_{\Delta\lambda}\int_0^t f(t-\lambda-\Delta\lambda)g(\lambda)\Delta\lambda] = {}_{\Delta\lambda}\int_0^\infty g(\lambda)[1+s\Delta\lambda]^{-\left(\frac{\lambda+\Delta\lambda}{\Delta\lambda}\right)}\,[{}_{\Delta\tau}\int_0^\infty f(\tau)\,[1+s\Delta\tau]^{-\left(\frac{\tau+\Delta\tau}{\Delta\tau}\right)}\Delta\tau]\Delta\lambda \tag{5.10-24}$$

Rearranging Eq 5.10-24

$$K_{\Delta t}[{}_{\Delta\lambda}\int_0^t f(t-\lambda-\Delta\lambda)g(\lambda)\Delta\lambda] = [{}_{\Delta\lambda}\int_0^\infty g(\lambda)[1+s\Delta\lambda]^{-\left(\frac{\lambda+\Delta\lambda}{\Delta\lambda}\right)}\Delta\lambda][{}_{\Delta\tau}\int_0^\infty f(\tau)\,[1+s\Delta\tau]^{-\left(\frac{\tau+\Delta\tau}{\Delta\tau}\right)}\Delta\tau] \tag{5.10-25}$$

Note that the right side of Eq 5.10-25 is the product of two $K_{\Delta t}$ Transforms.

$$K_{\Delta t}[\ {}_{\Delta\lambda}\int_0^t f(t-\lambda-\Delta\lambda)g(\lambda)\Delta\lambda] = K_{\Delta\lambda}[g(\lambda)]\,K_{\Delta\tau}[f(\tau)] \qquad (5.10\text{-}26)$$

Since $\Delta\tau = \Delta\lambda = \Delta t$, for clarity, the $\lambda$ and $\tau$ variables on the right side of Eq 5.10-26 can be changed to t.

Then

$$K_{\Delta t}[\ {}_{\Delta\lambda}\int_0^t f(t-\lambda-\Delta\lambda)g(\lambda)\Delta\lambda] = K_{\Delta t}[g(t)]\,K_{\Delta t}[f(t)] \qquad (5.10\text{-}27)$$

or

$$K_{\Delta t}[\ {}_{\Delta\lambda}\int_0^t f(t-\lambda-\Delta\lambda)g(\lambda)\Delta\lambda] = K_{\Delta t}[f(t)]\,K_{\Delta t}[g(t)] \qquad (5.10\text{-}28)$$

or interchanging the function designations

$$K_{\Delta t}[\ {}_{\Delta\lambda}\int_0^t f(\lambda)g(t-\lambda-\Delta\lambda)\Delta\lambda] = K_{\Delta t}[f(t)]\,K_{\Delta t}[g(t)] \qquad (5.10\text{-}29)$$

Thus

**The $K_{\Delta t}$ Transform Convolution Equation is:**

$$\mathbf{K_{\Delta t}[\ {}_{\Delta\lambda}\int_0^t f(t-\lambda-\Delta\lambda)g(\lambda)\Delta\lambda] = K_{\Delta t}[f(t)]\,K_{\Delta t}[g(t)]} \qquad (5.10\text{-}30)$$

**or**

$$\mathbf{K_{\Delta t}[\ {}_{\Delta\lambda}\int_0^t f(\lambda)g(t-\lambda-\Delta\lambda)\Delta\lambda] = K_{\Delta t}[f(t)]\,K_{\Delta t}[g(t)]} \qquad (5.10\text{-}31)$$

**where**

**$f(t),g(t)$ = discrete functions of t**

**$f(t-\lambda-\Delta\lambda)$ = discrete function of t, $\lambda$, and $\Delta\lambda$**

**$\Delta t = \Delta\lambda$**

**$t = 0, \Delta t, 2\Delta t, 3\Delta t, \ldots$**

**$\lambda = 0, \Delta\lambda, 2\Delta\lambda, 3\Delta\lambda, \ldots, t-\Delta\lambda, t$**

$$\mathbf{K_{\Delta t}[h(t)] = {}_{\Delta t}\int_0^\infty (1+s\Delta t)^{-\left(\frac{t+\Delta t}{\Delta t}\right)} h(t)\Delta t}$$ **, The $K_{\Delta t}$ Transform of a function of t**

**$K_{\Delta t}[f(t)]$ , $K_{\Delta t}[g(t)]$ = $K_{\Delta t}$ Transforms, functions of s**

**Comment - By definition, the discrete integral, $_{\Delta\lambda}\int_0^t h(\lambda)\Delta\lambda$ , is the sum, $_{\Delta\lambda}\sum_{\lambda=0}^{t-\Delta\lambda} h(\lambda)\Delta\lambda$ . Though the integration is from $\lambda = 0$ thru t, the summation that defines this integration is from $\lambda = 0$ thru $t-\Delta\lambda$.**

**Eq 5.10-30 can be rewritten in other useful forms**

$$_{\Delta\lambda}\int_0^t f(t-\lambda-\Delta\lambda)g(\lambda)\Delta\lambda = K_{\Delta t}^{-1}[\, K_{\Delta t}[f(t)]\, K_{\Delta t}[g(t)]\, ] \qquad (5.10\text{-}32)$$

**or**

$$_{\Delta\lambda}\int_0^t f(t-\lambda-\Delta\lambda)g(\lambda)\Delta\lambda = K_{\Delta t}^{-1}[F(s)G(s)] \qquad (5.10\text{-}33)$$

**where**

**$F(s) = K_{\Delta t}[f(t)]$**
**$G(s) = K[g(t)]$**
**$K_{\Delta t}^{-1}[H(s)]$ = Inverse $K_{\Delta t}$ Transform of a discrete function, H(s)**

**And Eq 5.10-31 can also be rewritten in other useful forms**

$$_{\Delta\lambda}\int_0^t f(\lambda)g(t-\lambda-\Delta\lambda)\,\Delta\lambda = K_{\Delta t}^{-1}[\, K_{\Delta t}[f(t)]\, K_{\Delta t}[g(t)]\, ] \qquad (5.10\text{-}34)$$

**or**

$$_{\Delta\lambda}\int_0^t f(\lambda)g(t-\lambda-\Delta\lambda)\Delta\lambda = K_{\Delta t}^{-1}[F(s)G(s)] \qquad (5.10\text{-}35)$$

**where**

**$F(s) = K_{\Delta t}[f(t)]$**
**$G(s) = K[g(t)]$**
**$K_{\Delta t}^{-1}[H(s)]$ = Inverse $K_{\Delta t}$ Transform of a discrete function, H(s)**

Below is an analysis of the $K_{\Delta t}$ Transform Convolution Discrete Integral, $_{\Delta\lambda}\int_0^t f(\lambda)g(t-\lambda-\Delta\lambda)\,\Delta\lambda$, that will make it more understandable.

Consider the following transfer function block diagram of a $K\Delta t$ Transform Transfer Function response to an input Unit Amplitude Pulse of $\Delta t$ width. The Unit Amplitude Pulse input starts at $t = n\Delta t$ and ends at $t = [n+1]\Delta t$ where $n = 0, 1, 2, 3, ...$ The $K_{\Delta t}$ Transform Transfer Function response to the Unit Amplitude Pulse input is $c(t,n) = K_{\Delta t}^{-1}[R(s)G(s)]$.

$g(t) = K_{\Delta t}^{-1}[G(s)]$

$r(t,n) = U(t-n\Delta t) - U(t-[n+1]\Delta t)$ $\qquad$ $c(t,n) = g(t-[n+1]\Delta t)U(t-[n+1]\Delta t)\Delta t$

$K_{\Delta t}$ Transform
Transfer Function

$R(s) = (1+s\Delta t)^{-(n+1)}\Delta t$ $\qquad$ $G(s)$ $\qquad$ $C(s) = R(s)G(s) = (1+s\Delta t)^{-(n+1)}G(s)\Delta t$

Input $\qquad$ Transfer Function $\qquad$ Output

$$c(t,n) = K_{\Delta t}^{-1}[(1+s\Delta t)^{-(n+1)}G(s)\Delta t] \tag{5.10-36}$$

$$c(t,n) = g(t-[n+1]\Delta t)U(t-[n+1]\Delta t)\Delta t \tag{5.10-37}$$
$K_{\Delta t}$ Transform Transfer Function Response to an input Unit Amplitude Pulse of $\Delta t$ width

where

$r(t,n) = U(t-n\Delta t) - U(t-[n+1]\Delta t)$ , Unit Amplitude Pulse Input

$R(s) = \dfrac{(1+s\Delta t)^{-n} - (1+s\Delta t)^{-(n+1)}}{s} = (1+s\Delta t)^{-(n+1)}\Delta t$ , Unit Amplitude Pulse Input $K_{\Delta t}$ Transform

$G(s) = K_{\Delta t}$ Transform Transfer Function

$g(t) = K_{\Delta t}^{-1}[G(s)]$ , Inverse $K_{\Delta t}$ Transform of $G(s)$

$C(s) = K_{\Delta t}[c(t)] = R(s)G(s) = (1+s\Delta t)^{-(n+1)}G(s)\Delta t$ , Output $K_{\Delta t}$ Transform

$c(t,n) = g(t-[n+1]\Delta t)U(t-[n+1]\Delta t)\Delta t$ , Output response to an input Unit Amplitude Pulse

$\Delta t$ = Interval between successive discrete values of t

$t = 0, \Delta t, 2\Delta t, 3\Delta t, \ldots$ , discrete values of t

$n = 0, 1, 2, 3, \ldots$

From the above analysis, the $K_{\Delta t}$ Transform Transfer Function, $G(s)$, response to an input Unit Amplitude Pulse, $r(t,n) = U(t-n\Delta t) - U(t-[n+1]\Delta t)$, is $c(t,n) = g(t-[n+1]\Delta t)U(t-[n+1]\Delta t)\Delta t$ where n denotes the interval, $n\Delta t \le t \le t = [n+1]\Delta t$, in which the Unit Amplitude Pulse occurs.

More generally, consider the input to the $K_{\Delta t}$ Transform function, G(s), to be a discrete sample and hold shaped waveform, f(t). All discrete functions used in Interval Calculus are sample and hold shaped waveform functions. Note Diagram 5.10-1 below.

Diagram 5.10-1  A diagram of a sample and hold shaped waveform, f(t)

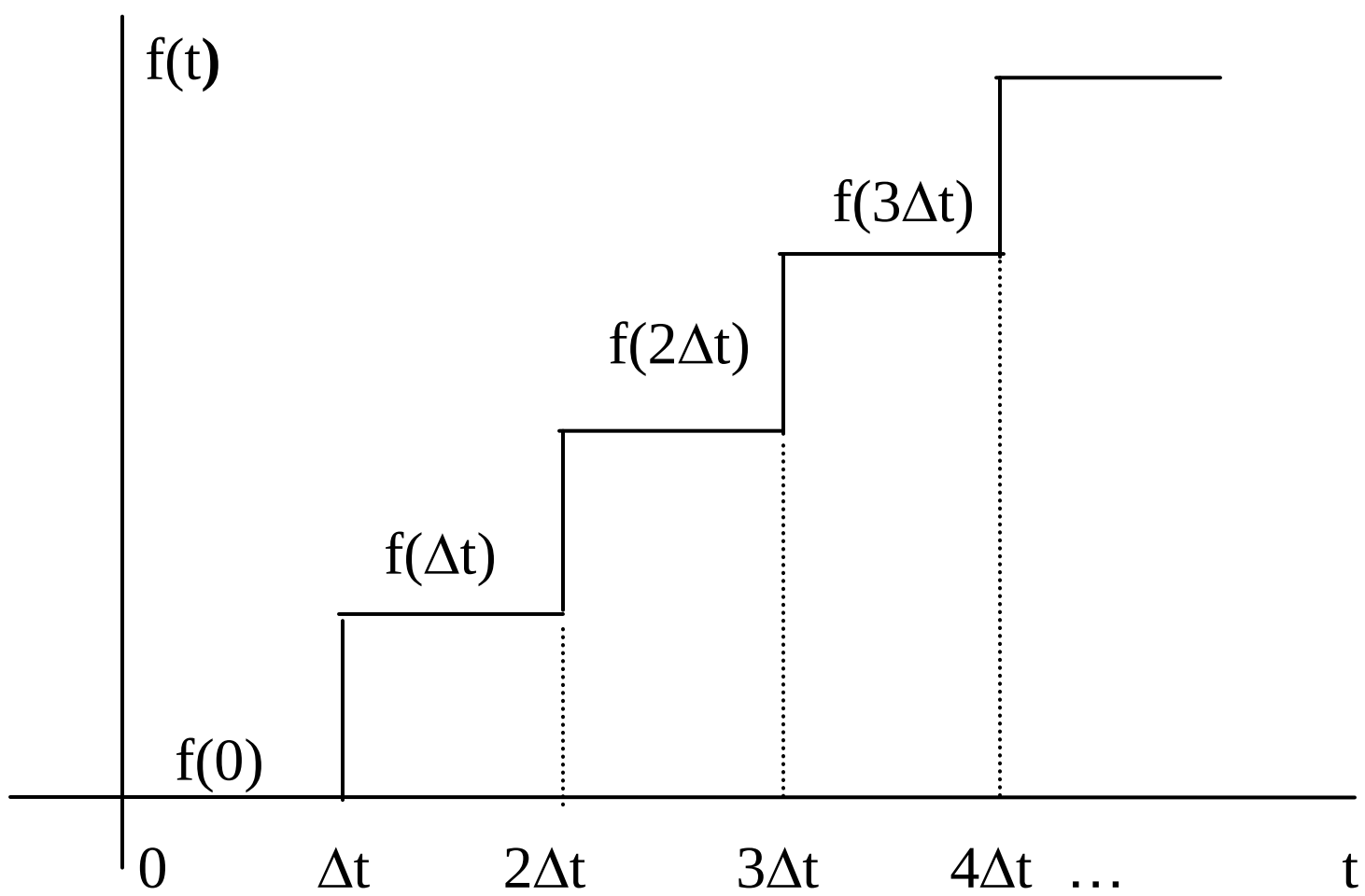

The sample and hold shaped waveform, f(t), can be considered to be a series of consecutive pulses of Δt width where the amplitude of the nth consecutive pulse is f(nΔt) , n = 0, 1, 2, 3, … This can be expressed mathematically as follows:

f(t) = Infinite series of pulses
f(t) = Summation of Unit Amplitude Pulses, U(t-nΔt) – U(t-[n+1]Δt), each weighted to have the value, f(nΔt), in the nth interval, $n\Delta t \le t \le [n+1]\Delta t$

$$f(t) = {}_1\sum_{n=0}^{\infty} f(n\Delta t) r(t,n) = {}_1\sum_{n=0}^{\infty} f(n\Delta t)[\, U(t-n\Delta t) - U(t-[n+1]\Delta t)\, ] \qquad (5.10\text{-}38)$$

where
n = 0, 1, 2, 3, …  ,  f(t) interval designation
r(t,n) = U(t-nΔt) – U(t-[n+1]Δt) ,  nth interval Unit Amplitude Pulse
Δt = width of all Unit Amplitude Pulses, r(t,n)
f(nΔt) = the amplitude of the f(t) consecutive pulse in the nth interval, $n\Delta t \le t \le [n+1\}\Delta t$

Consider the following transfer function block diagram of a $K_{\Delta t}$ Transform Transfer Function response to the input function, f(t).

$f(t) = {}_1\sum_{n=0}^{\infty} f(n\Delta t)[\, U(t-n\Delta t)-U(t-[n+1]\Delta t)\,]$ $\quad g(t) = K_{\Delta t}^{-1}[G(s)]$ $\quad O(t) = K_{\Delta t}^{-1}[\, K_{\Delta t}[f(t)]\, K_{\Delta t}[g(t)]\,] =$

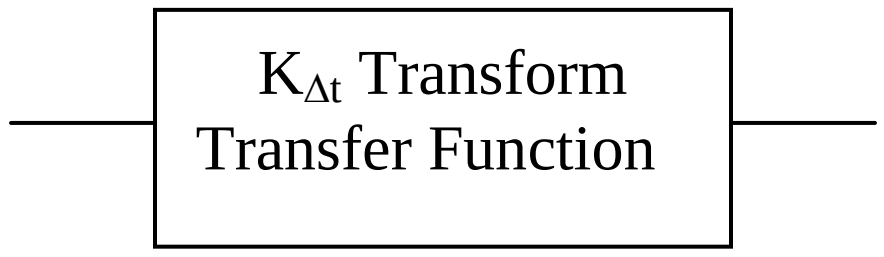

${}_1\sum_{n=0}^{\infty} f(n\Delta t)g(t-[n+1]\Delta t)U(t-[n+1]\Delta t)\Delta t$

$K_{\Delta t}[f(t)] = {}_1\sum_{n=0}^{\infty} f(n\Delta t)(1+s\Delta t)^{-(n+1)}\Delta t$ $\quad G(s)$ $\quad O(s) = K_{\Delta t}[f(t)]\, K_{\Delta t}[g(t)]$

Input $\qquad$ Transfer Function $\qquad$ Output

O(t) is the $K_{\Delta t}$ Transform Transfer Function Response to a series of weighted Unit Amplitude Pulse inputs which are equivalent to the input function, f(t).

$$O(t) = K_{\Delta t}^{-1}[\, K_{\Delta t}[f(t)]\, K_{\Delta t}[g(t)]\,] = {}_1\sum_{n=0}^{\infty} f(n\Delta t)g(t-[n+1]\Delta t)U(t-[n+1]\Delta t)\Delta t \qquad (5.10\text{-}39)$$

where
- f(t) = Input
- F(s) = Input $K_{\Delta t}$ Transform
- G(s) = $K_{\Delta t}$ Transform Transfer Function
- $g(t) = K_{\Delta t}^{-1}[G(s)]$, Inverse $K_{\Delta t}$ Transform of G(s)
- $O(s) = K_{\Delta t}[O(t)] = R(s)G(s)$, Output $K_{\Delta t}$ Transform
- O(t) = Output response to the Input, f(t)
- $\Delta t$ = Interval between successive discrete values of t
- $t = 0, \Delta t, 2\Delta t, 3\Delta t, \ldots$ = discrete independent variable values
- $n = 0, 1, 2, 3, \ldots$

Then

The $K_{\Delta t}$ Transform Function output response, O(t) to the input function f(t) is:

$$O(t) = {}_1\sum_{n=0}^{\infty} f(n\Delta t)g(t-[n+1]\Delta t)U(t-[n+1]\Delta t)\Delta t = K_{\Delta t}^{-1}[\, K_{\Delta t}[f(t)]\, K_{\Delta t}[g(t)]\,] \qquad (5.10\text{-}40)$$

Let

$$\lambda = n\Delta t \qquad (5.10\text{-}41)$$

$$\Delta\lambda = \Delta t \qquad (5.10\text{-}42)$$

$$n = 0, 1, 2, 3, \ldots \qquad (5.10\text{-}43)$$

$$\lambda = 0, \Delta t, 2\Delta t, 3\Delta t, \ldots \qquad (5.10\text{-}44)$$

Substituting Eq 5.10-41 thru Eq 5.10-44 into Eq 5.10-40

$$O(t) = {}_{\Delta\lambda}\sum_{\lambda=0}^{\infty} f(\lambda)g(t-\lambda-\Delta\lambda)U(t-\lambda-\Delta\lambda)\Delta\lambda = K_{\Delta t}^{-1}[\, K_{\Delta t}[f(t)]\, K_{\Delta t}[g(t)]\,] \qquad (5.10\text{-}45)$$

Since

$$U(t-\lambda-\Delta\lambda) = \begin{cases} 1 & \text{for } \lambda < t \\ 0 & \text{for } \lambda \geq t \end{cases} \qquad (5.10\text{-}46)$$

Eq 5.9-45 can be simplified

$$O(t) = {}_{\Delta\lambda}\sum_{\lambda=0}^{t-\Delta t} f(\lambda)g(t-\lambda-\Delta\lambda)\Delta\lambda = K_{\Delta t}^{-1}[\, K_{\Delta t}[f(t)]\, K_{\Delta t}[g(t)]\,] \qquad (5.10\text{-}47)$$

Changing the form of Eq 5.9-47 to its discrete integral form

$$O(t) = {}_{\Delta\lambda}\int_{0}^{t} f(\lambda)g(t-\lambda-\Delta\lambda)\Delta\lambda = K_{\Delta t}^{-1}[\, K_{\Delta t}[f(t)]\, K_{\Delta t}[g(t)]\,] \qquad (5.10\text{-}48)$$

From Eq 5.10-34 the $K_{\Delta t}$ Transform Convolution Equation is:

$$ {}_{\Delta\lambda}\int_{0}^{t} f(\lambda)g(t-\lambda-\Delta\lambda)\Delta\lambda = K_{\Delta t}^{-1}[\, K_{\Delta t}[f(t)]\, K_{\Delta t}[g(t)]\,] = O(t) \qquad (5.10\text{-}49)$$

Comparing Eq 5.9-48 to the $K_{\Delta t}$ Transform Convolution Equation, Eq 5.9-49 the two equations are seen to be the same. Thus, the $K_{\Delta t}$ Transform Convolution integral can now be described as follows:

1. The $K_{\Delta t}$ Transform Convolution integral, $O(t) = {}_{\Delta\lambda}\int_{0}^{t} f(\lambda)g(t-\lambda-\Delta\lambda)\Delta\lambda$ , is the sum of all G(s) Transfer Function responses, ${}_{1}\sum_{n=0}^{\infty} f(n\Delta t)g(t-[n+1]\Delta t)U(t-[n+1]\Delta t)\Delta t$, to all $f(n\Delta t)$ weighted Unit Amplitude Pulses, ${}_{1}\sum_{n=0}^{\infty} f(n\Delta t)[\,U(t-n\Delta t)-U(t-[n+1]\Delta t)\,]$, that together comprise the input function, f(t). $\lambda = n\Delta t$, $\Delta\lambda = \Delta t$, and n = 0, 1, 2, 3, ...

2. The $K_{\Delta t}$ Transform Convolution integral term, $g(t-\lambda-\Delta\lambda)\Delta\lambda$, if more completely written, $g(t-\lambda-\Delta\lambda)U(t-\lambda-\Delta\lambda)\Delta\lambda$ (see Eq 5.10-45), represents the output response of a $K_{\Delta t}$ Transform Transfer Function, G(s), to an input Unit Amplitude Pulse of $\Delta\lambda$ width initiated at $\lambda$ and ending at $\lambda+\Delta\lambda$ where $\lambda = 0, \Delta\lambda, 2\Delta\lambda, 3\Delta\lambda, \ldots, t-\Delta\lambda$. $\lambda = n\Delta t$, $\Delta\lambda = \Delta t$, and n = 0, 1, 2, 3, ... , $\frac{t}{\Delta\lambda} -1$

3. The $K_{\Delta t}$ Transform Convolution integral term, $f(\lambda) = f(n\Delta t)$, is the value of the sample and hold shaped waveform, f(t), in its nth interval. $f(\lambda)$ is the weighting placed on each consecutive pulse of the f(t) infinite series of Unit Amplitude Pulses of $\Delta\lambda$ width where $\lambda = 0, \Delta\lambda, 2\Delta\lambda, 3\Delta\lambda, \ldots, t-\Delta\lambda$ and $n = 0, 1, 2, 3, \ldots, \frac{t}{\Delta\lambda} -1$ . $\lambda = n\Delta t$, $\Delta\lambda = \Delta t$

To demonstrate the use of the $K_{\Delta t}$ Transform Convolution Equation, several examples are provided below.

<u>Example 5.10-1</u> Find the Inverse $K_{\Delta t}$ Transform, $K_{\Delta t}^{-1}[\frac{1}{s(s-a)}]$, using the $K_{\Delta t}$ Transform Convolution Equation. The constant, a, is any real value.

The $K_{\Delta t}$ Transform Convolution Equations is:

$$K_{\Delta t}[\underset{\Delta\lambda}{\int_0^t} f(t-\lambda-\Delta\lambda)g(\lambda)\Delta\lambda] = K_{\Delta t}[f(t)]\,K_{\Delta t}[g(t)] \qquad 1)$$

where

$\Delta t = \Delta\lambda$

$\lambda = 0, \Delta\lambda, 2\Delta\lambda, 3\Delta\lambda, \ldots, t-\Delta\lambda, t$

$f(t),g(t)$ = discrete functions of t

$f(t-\lambda-\Delta\lambda)$ = discrete function of t, $\lambda$, and $\Delta\lambda$

$t = 0, \Delta t, 2\Delta t, 3\Delta t, \ldots$

$K_{\Delta t}[f(t)] = F(s) = K_{\Delta t}$ Transform of f(t)

$$K_{\Delta t}[f(t)] = \underset{\Delta t}{\int_0^\infty} [1+s\Delta t]^{-\frac{t+\Delta t}{\Delta t}} f(t)\Delta t \quad = K_{\Delta t} \text{ Transform of the function, } f(t)$$

Find $K_{\Delta t}^{-1}[\frac{1}{s(s-a)}]$

Let

$$K_{\Delta t}[f(t)] = \frac{1}{s-a} \qquad 2)$$

$$K_{\Delta t}[g(t)] = \frac{1}{s} \qquad 3)$$

$$K_{\Delta t}[f(t)]K_{\Delta t}[g(t)] = [\frac{1}{s(s-a)}] \qquad 4)$$

Writing the inverse $K_{\Delta t}$ Transforms of Eq 2 and Eq 3

$$f(t) = e_{\Delta t}(a,t) = [1+a\Delta t]^{\frac{t}{\Delta t}} \qquad 5)$$
$$g(t) = U(t) = 1 \qquad 6)$$

Substituting Eq 2 thru Eq 6 into Eq 1 and taking the Inverse $K_{\Delta t}$ Transform

$$K_{\Delta t}^{-1}[\frac{1}{s(s-a)}] = {}_{\Delta\lambda}\int_0^t [1+a\Delta\lambda]^{\frac{t-\lambda-\Delta\lambda}{\Delta\lambda}} (1)\, \Delta\lambda \tag{7}$$

$$K_{\Delta t}^{-1}[\frac{1}{s(s-a)}] = [1+a\Delta\lambda]^{\frac{t}{\Delta\lambda}-1}[{}_{\Delta\lambda}\int_0^t [1+a\Delta\lambda]^{-\frac{\lambda}{\Delta\lambda}} \Delta\lambda\, ] \tag{8}$$

From the integration table, Table 6, in the Appendix

$$_{\Delta x}\int (1+a\Delta x)^{-\frac{x}{\Delta x}} \Delta x = -\frac{1+a\Delta x}{a} (1+a\Delta x)^{-\frac{x}{\Delta x}} + k \tag{9}$$

Integrating the right side of Eq 8 using Eq 9

$$K_{\Delta t}^{-1}[\frac{1}{s(s-a)}] = [1+a\Delta\lambda]^{\frac{t}{\Delta\lambda}-1}(-\frac{1+a\Delta\lambda}{a})\, [1+a\Delta\lambda]^{-\frac{\lambda}{\Delta\lambda}} \Big|_0^t \tag{10}$$

$$K_{\Delta t}^{-1}[\frac{1}{s(s-a)}] = [1+a\Delta\lambda]^{\frac{t}{\Delta\lambda}-1}(-\frac{1+a\Delta\lambda}{a})(\, [1+a\Delta\lambda]^{-\frac{t}{\Delta\lambda}} -1) \tag{11}$$

$$\Delta t = \Delta\lambda \tag{12}$$

Substituting Eq 12 into Eq 11

$$K_{\Delta t}^{-1}[\frac{1}{s(s-a)}] = [1+a\Delta t]^{\frac{t}{\Delta t}-1}(-\frac{1+a\Delta t}{a})(\, [1+a\Delta t]^{-\frac{t}{\Delta t}} -1) \tag{13}$$

$$K_{\Delta t}^{-1}[\frac{1}{s(s-a)}] = -\frac{1}{a}\, [1+a\Delta t]^{\frac{t}{\Delta t}}(\, [1+a\Delta t]^{-\frac{t}{\Delta t}} -1) \tag{14}$$

$$K_{\Delta t}^{-1}[\frac{1}{s(s-a)}] = -\frac{1}{a}\, (\, 1 - [1+a\Delta t]^{\frac{t}{\Delta t}}\, ) \tag{15}$$

$$K_{\Delta t}^{-1}[\frac{1}{s(s-a)}] = \frac{[1+a\Delta t]^{\frac{t}{\Delta t}} - 1}{a} \tag{16}$$

$$e_{\Delta t}(a,t) = [1+a\Delta t]^{\frac{t}{\Delta t}} \tag{17}$$

Then

Substituting Eq 17 into Eq 16

$$\mathbf{K_{\Delta t}^{-1}[\,\frac{1}{s(s-a)}\,] = \frac{e_{\Delta t}(a,t) - 1}{a}} \qquad 18)$$

Checking Eq 18

Taking the $K_{\Delta t}$ Transform of Eq 18

$$\frac{1}{s(s-a)} = \frac{1}{a}\left(\frac{1}{s-a} - \frac{1}{s}\right) = \frac{s-s+a}{as(s-a)} = \frac{1}{s(s-a)}$$

Good check

<u>Example 5.10-2</u> Find the output of a unit amplitude pulse into an integrator.Use the $K_{\Delta t}$ Transform Convolution Equation.

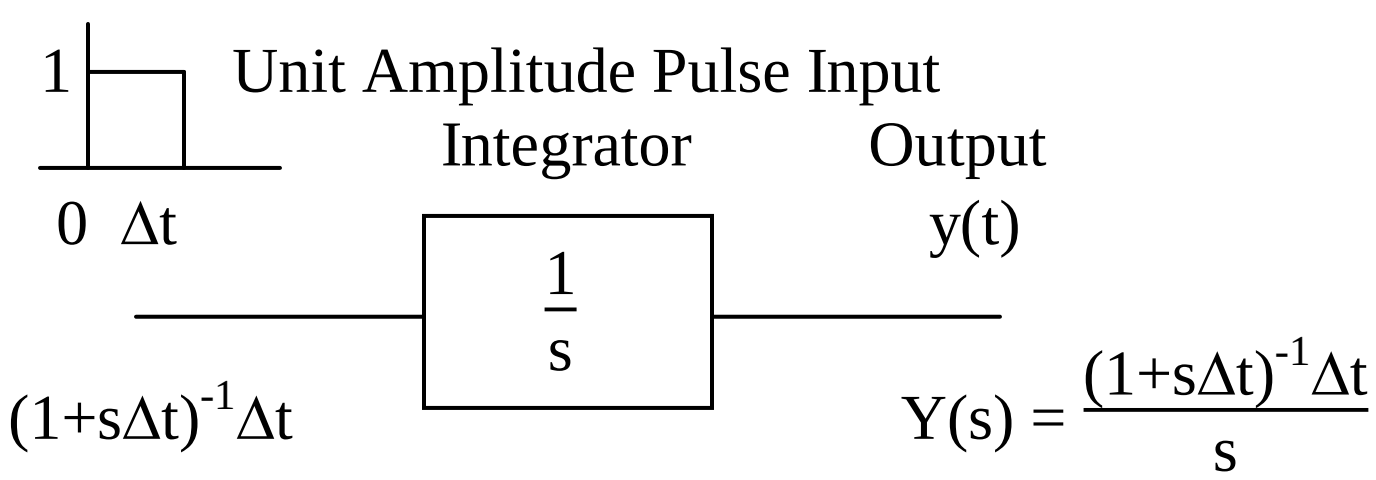

The $K_{\Delta t}$ Transform Convolution Equations is:

$$K_{\Delta t}[\int_{0,\,\Delta\lambda}^{t} f(t-\lambda-\Delta\lambda)g(\lambda)\Delta\lambda] = K_{\Delta t}[f(t)]K_{\Delta t}[g(t)] \qquad 1)$$

where

$\Delta t = \Delta\lambda$

$\lambda = 0, \Delta\lambda, 2\Delta\lambda, 3\Delta\lambda, \ldots, t-\Delta\lambda, t$

f(t),g(t) = discrete functions of t

$f(t-\lambda-\Delta\lambda)$ = discrete function of t, $\lambda$, and $\Delta\lambda$

$t = 0, \Delta t, 2\Delta t, 3\Delta t, \ldots$

$K_{\Delta t}[f(t)] = F(s) = K_{\Delta t}$ Transform of f(t)

$$K_{\Delta t}[f(t)] = \int_{0,\,\Delta t}^{\infty} [1+s\Delta t]^{-\frac{t+\Delta t}{\Delta t}} f(t)\Delta t \;=\; K_{\Delta t} \text{ Transform of the function, f(t)}$$

Find $y(t) = K_{\Delta t}^{-1}[\frac{(1+s\Delta t)^{-1}\Delta t}{s}]$ using the $K_{\Delta t}$ Transform Convolution Equation, Eq 1.

Let

$K_{\Delta t}[f(t)] = \frac{1}{s}$ , The $K_{\Delta t}$ Transform of the integrator 2)

$K_{\Delta t}[g(t)] = (1+s\Delta t)^{-1}\Delta t$ , The $K_{\Delta t}$ Transform of the Unit amplitude pulse input 3)

$K_{\Delta t}[f(t)]\, K_{\Delta t}[g(t)] = Y(s) = \frac{(1+s\Delta t)^{-1}\Delta t}{s}$ , The $K_{\Delta t}$ Transform of the unit amplitude pulse output 4)

Finding the Inverse $K_{\Delta t}$ Transforms of Eq 2 and Eq 3

$$f(t) = K_{\Delta t}^{-1}[\frac{1}{s}] = U(t) \quad 5)$$

$$g(t) = K_{\Delta t}^{-1}[(1+s\Delta t)^{-1}\Delta t] = [U(t) - U(t-\Delta t)] \quad 6)$$

From Eq 1 and Eq 4 thru Eq 6

$$y(t) = K_{\Delta t}^{-1}[\frac{1}{s}(1+s\Delta t)^{-1}\Delta t] = {}_{\Delta\lambda}\int_0^t U(t-\lambda-\Delta\lambda)[U(\lambda) - U(\lambda-\Delta\lambda)]\Delta\lambda \quad 7)$$

$$y(t) = {}_{\Delta\lambda}\int_0^t U(t-\lambda-\Delta\lambda)U(\lambda)\Delta\lambda - {}_{\Delta\lambda}\int_0^t U(t-\lambda-\Delta\lambda)U(\lambda-\Delta\lambda)]\Delta\lambda \quad 8)$$

$$y(t) = {}_{\Delta\lambda}\int_0^t \Delta\lambda - {}_{\Delta\lambda}\int_{\Delta\lambda}^t \Delta\lambda \quad 9)$$

$$y(t) = \lambda \Big|_0^t - \lambda \Big|_{\Delta\lambda}^t \quad 10)$$

$$y(t) = t - t + \Delta\lambda = \Delta\lambda \quad 11)$$

$$\Delta\lambda = \Delta t \quad 12)$$

Substituting Eq 12 into Eq 11

$$\mathbf{y(t) = \Delta t} \quad 13)$$

Checking Eq 13

$$y(t) = K_{\Delta t}^{-1}[\frac{1}{s}(1+s\Delta t)^{-1}\Delta t\ ] = K_{\Delta t}^{-1}[\frac{1}{s(s+\frac{1}{\Delta t})}] = K_{\Delta t}^{-1}[\frac{\Delta t}{s} - \frac{\Delta t}{s+\frac{1}{\Delta t}}]$$

Taking the Inverse $K_{\Delta t}$ Transform

$$y(t) = \Delta t - \Delta t[1 - \frac{1}{\Delta t}\Delta t]^{\frac{t}{\Delta t}} = \Delta t$$

$y(t) = \Delta t$

Good check

Derivation of the Z Transform Convolution Equation

f(t), g(t) = sample and hold shaped waveform functions

Consider the following discrete function.

$$\underset{\Delta\lambda}{\int_{0}^{t+T}} f(t-\lambda)g(\lambda)\Delta\lambda] \; = \text{a discrete function of } \lambda \text{ and } t \qquad (5.10\text{-}50)$$

where

- $f(t-\lambda)$ = discrete function of t and $\lambda$
- $g(\lambda)$ = discrete function of $\lambda$
- $\Delta t = \Delta\lambda = T$
- $t = 0, \Delta t, 2\Delta t, 3\Delta t, \ldots$
- $\lambda = 0, \Delta\lambda, 2\Delta\lambda, 3\Delta\lambda, \ldots, t\text{-}\Delta\lambda, t, t\text{+}T$

$$Z[h(t)]\;) = \frac{1}{T}\underset{T}{\int_{0}^{\infty}} z^{-(\frac{t}{\Delta t})} h(t)\Delta t \;\;, \quad \text{Z Transform} \qquad (5.10\text{-}51)$$

where

- $t = 0, \Delta t, 2\Delta t, 3\Delta t, \ldots, \infty$
- $\Delta t = T$
- $\Delta t$ = the interval between successive values of t
- $h(t)$ = discrete function of t

Take the Z Transform of Eq 5.10-50

$$Z[\underset{\Delta\lambda}{\int_{0}^{t+T}} f(t-\lambda)g(\lambda)\Delta\lambda] \; = \; \frac{1}{T}\underset{T}{\int_{0}^{\infty}} [\underset{\Delta\lambda}{\int_{0}^{t+T}} f(t-\lambda)g(\lambda)\Delta\lambda] z^{-(\frac{t}{\Delta t})}\Delta t \qquad (5.10\text{-}52)$$

<u>Comment</u> – For discrete integrals, the discrete variable range is: lower limit ≤ variable < upper limit.

$$\Delta\lambda = \Delta t = T \qquad (5.10\text{-}53)$$

$$\lambda = 0, \Delta\lambda, 2\Delta\lambda, 3\Delta\lambda, \ldots, t\text{-}\Delta\lambda, t, t\text{+}T \qquad (5.10\text{-}54)$$

$$U(t\text{-}\lambda) = \begin{cases} 1 & \text{for } \lambda < t\text{+}T \\ 0 & \text{for } \lambda \geq t\text{+}T \end{cases} \qquad (5.10\text{-}55)$$

Then

$$f(t-\lambda)g(\lambda)U(t\text{-}\lambda) = \begin{cases} f(t-\lambda)g(\lambda) & \text{for } \lambda < t\text{+}T \\ 0 & \text{for } \lambda \geq t\text{+}T \end{cases} \qquad (5.10\text{-}56)$$

Since the product of Eq 5.10-56 is 0 for all $\lambda \geq t+T$, the inner integration in Eq 5.10-52 can be extended to $\infty$ if the factor $U(t-\lambda)$ is inserted into the integrand.

$$Z_T\Big[\int_{0}^{t+T}{}_{\Delta\lambda} f(t-\lambda)g(\lambda)\Delta\lambda\Big] = \frac{1}{T}\int_{0}^{\infty}{}_{T}\Big[\int_{0}^{\infty}{}_{\Delta\lambda} f(t-\lambda)g(\lambda)U(t-\lambda)\Delta\lambda\Big]z^{-\left(\frac{t}{\Delta t}\right)}\Delta t \tag{5.10-57}$$

Our usual assumptions about the functions we transform are sufficient to permit the order of integration in Eq 5.10-57 to be interchanged.

$$Z_T\Big[\int_{0}^{t+T}{}_{\Delta\lambda} f(t-\lambda)g(\lambda)\Delta\lambda\Big] = \frac{1}{T}\int_{0}^{\infty}{}_{\Delta\lambda}\Big[\int_{0}^{\infty}{}_{T} f(t-\lambda)g(\lambda)U(t-\lambda)z^{-\left(\frac{t}{\Delta t}\right)}\Delta t\Big]\Delta\lambda \tag{5.10-58}$$

$$Z_T\Big[\int_{0}^{t+T}{}_{\Delta\lambda} f(t-\lambda)g(\lambda)\Delta\lambda\Big] = \frac{1}{T}\int_{0}^{\infty}{}_{\Delta\lambda} g(\lambda)\Big[\int_{0}^{\infty}{}_{T} f(t-\lambda)U(t-\lambda)z^{-\left(\frac{t}{\Delta t}\right)}\Delta t\Big]\Delta\lambda \tag{5.10-59}$$

Due to the term, $U(t-\lambda)$, the integrand of the inner integral of Eq 5.10-59 is 0 for $t < \lambda$

$$Z_T\Big[\int_{0}^{t+T}{}_{\Delta\lambda} f(t-\lambda)g(\lambda)\Delta\lambda\Big] = \frac{1}{T}\int_{0}^{\infty}{}_{\Delta\lambda} g(\lambda)\Big[\int_{\lambda}^{\infty}{}_{\Delta t} f(t-\lambda)z^{-\left(\frac{t}{\Delta t}\right)}\Delta t\Big]\Delta\lambda \tag{5.10-60}$$

In the inner integral on the right side of Eq 5.10-60

Let

$$\tau = t-\lambda \tag{5.10-61}$$

From Eq 5.10-61

$$\Delta t = \Delta\tau \tag{5.10-62}$$

$$t = \tau + \lambda \tag{5.10-63}$$

From Eq 5.10-63, finding the limits of $\tau$ from the limits of t

$$\text{For } t = \lambda\,,\ \ \tau = 0 \tag{5.10-64}$$

$$\text{For } t = \infty\,,\ \ \tau = \infty \tag{5.10-65}$$

Substituting Eq 5.10-62 thru Eq 5.10-65 into Eq 5.10-60

$$Z[{}_{\Delta\lambda}\int_0^{t+T} f(t-\lambda)g(\lambda)\Delta\lambda] = \frac{1}{T}{}_{\Delta\lambda}\int_0^{\infty} g(\lambda)[{}_{\Delta\tau}\int_0^{\infty} f(\tau)z^{-(\frac{\tau+\lambda}{\Delta\tau})}\Delta\tau]\Delta\lambda \tag{5.10-66}$$

Rearranging Eq 5.10-66

$$Z[{}_{\Delta\lambda}\int_0^{t+T} f(t-\lambda)g(\lambda)\Delta\lambda] = \frac{1}{T}{}_{\Delta\lambda}\int_0^{\infty} g(\lambda)z^{-(\frac{\lambda}{\Delta\tau})}[{}_{\Delta\tau}\int_0^{\infty} f(\tau)z^{-(\frac{\tau}{\Delta\tau})}\Delta\tau]\Delta\lambda \tag{5.10-67}$$

From Eq 5.10-53 and Eq 5.10-62

$$\Delta\tau = \Delta\lambda = \Delta t = T \tag{5.10-68}$$

From Eq 5.10-67 and Eq 5.10-68

$$Z[{}_{T}\int_0^{t+T} f(t-\lambda)g(\lambda)\Delta\lambda] = \frac{1}{T}{}_{\Delta\lambda}\int_0^{\infty} g(\lambda)z^{-(\frac{\lambda}{\Delta\lambda})}[{}_{\Delta\tau}\int_0^{\infty} f(\tau)z^{-(\frac{\tau}{\Delta\tau})}\Delta\tau]\Delta\lambda \tag{5.10-69}$$

Rearranging Eq 5.10-69 and multiplying both sides by $\frac{1}{T}$

$$\frac{1}{T}Z[{}_{T}\int_0^{t} f(t-\lambda)g(\lambda)\Delta\lambda] = [\frac{1}{T}{}_{\Delta\lambda}\int_0^{\infty} g(\lambda)z^{-(\frac{\lambda}{\Delta\lambda})}\Delta\lambda][\frac{1}{T}{}_{\Delta\tau}\int_0^{\infty} f(\tau)z^{-(\frac{\tau}{\Delta\tau})}\Delta\tau] \tag{5.10-70}$$

Note that the right side of Eq 5.10-70 is the product of two Z Transforms.

$$Z[\frac{1}{T}{}_{T}\int_0^{t+T} f(t-\lambda)g(\lambda)\Delta\lambda] = Z[g(\lambda)]\,Z[f(\tau)] \tag{5.10-71}$$

Since $\Delta\tau = \Delta\lambda = \Delta t$, for clarity, the $\lambda$ and $\tau$ variables on the right side of Eq 5.10-71 can be changed to t.

Then

$$Z[\frac{1}{T}{}_T\int_0^{t+T} f(t-\lambda)g(\lambda)\Delta\lambda] = Z[g(t)]\,Z[f(t)] \qquad (5.10\text{-}72)$$

or

$$Z[\frac{1}{T}{}_T\int_0^{t+T} f(t-\lambda)g(\lambda)\Delta\lambda] = Z[f(t)]\,Z[g(t)] \qquad (5.10\text{-}73)$$

or interchanging the function designations

$$Z[\frac{1}{T}{}_T\int_0^{t+T} f(\lambda)g(t-\lambda)\Delta\lambda] = Z[f(t)]\,Z[g(t)] \qquad (5.10\text{-}74)$$

Thus

**The Z Transform Convolution Equation is:**

$$\mathbf{Z[\frac{1}{T}{}_T\int_0^{t+T} f(t-\lambda)g(\lambda)\Delta\lambda] = Z[f(t)]\,Z[g(t)]} \qquad (5.10\text{-}75)$$

**or**

$$\mathbf{Z[\frac{1}{T}{}_T\int_0^{t+T} f(\lambda)g(t-\lambda)\Delta\lambda] = Z[f(t)]\,Z[g(t)]} \qquad (5.10\text{-}76)$$

**where**

**f(t),g(t) = discrete functions of t**

**$f(t-\lambda)$ = discrete function of t and $\lambda$**

**$\Delta t = \Delta\lambda = T$**

**$t = 0, \Delta t, 2\Delta t, 3\Delta t, \ldots$**

**$\lambda = 0, \Delta\lambda, 2\Delta\lambda, 3\Delta\lambda, \ldots, t-\Delta\lambda, t, t+T$**

**$Z[h(t)] = {}_{\Delta t}\int_0^{\infty} z^{-(\frac{t}{\Delta t})} h(t)\Delta t$ ,     Z Transform of a function of t**

**Z[f(t)] , Z[g(t)] = Z Transforms, functions of z**

**Eq 5.10-75 can be rewritten in other useful forms**

$$\mathbf{\frac{1}{T}{}_T\int_0^{t+T} f(t-\lambda)g(\lambda)\Delta\lambda = Z^{-1}[\,Z[f(t)]\,Z[g(t)]\,]} \qquad (5.10\text{-}77)$$

**or**

$$\frac{1}{T}\,{}_T\!\int_0^{t+T} f(\lambda)g(t-\lambda)\Delta\lambda \;=\; Z^{-1}[\, Z[f(t)]\, Z[g(t)]\, ] \qquad (5.10\text{-}78)$$

**or**

$$\frac{1}{T}\,{}_T\!\int_0^{t+T} f(t-\lambda)g(\lambda)\Delta\lambda \;=\; Z^{-1}[F(z)G(z)] \qquad (5.10\text{-}79)$$

**or**

$$\frac{1}{T}\,{}_T\!\int_0^{t+T} f(\lambda)g(t-\lambda)\Delta\lambda \;=\; Z^{-1}[F(z)G(z)] \qquad (5.10\text{-}80)$$

**where**

**$f(t) = Z^{-1}[F(z)]$**
**$g(t) = Z^{-1}[G(z)]$**
**$Z^{-1}[H(z)]$ = Inverse Z Transform of a function, H(z)**

__Comment__ - **By definition, the discrete integral, ${}_{\Delta\lambda}\!\int_0^t h(\lambda)\Delta\lambda$ , is the sum, ${}_{\Delta\lambda}\sum_{\lambda=0}^{t-\Delta\lambda} h(\lambda)\Delta\lambda$ . Though the integration is from $\lambda = 0$ thru t, the summation that defines this integration is from $\lambda = 0$ thru $t-\Delta\lambda$.**

There is another form of the Z Transform Convolution Equation that can be derived directly from the $K_{\Delta t}$ Transform Convolution Equation, Eq 5.10-30 and Eq 5.10-31, with the use of the $K_{\Delta t}$ Transform to Z Transform Conversion Equation, Eq 5.5-1. It is derived as follows:

Rewriting The $K_{\Delta x}$ Transform to Z Transform Conversion Equation, Eq 5.5-1

$$Z[f(t)] = F(z) = \frac{z}{T} F(s)\Big|_{s = \frac{z-1}{T}} \qquad Z[f(t)] = F(z) \quad \text{Z Transform} \qquad (5.10\text{-}81)$$

$T = \Delta t$
$t = 0, T, 2T, 3T, \ldots$ $\qquad K_{\Delta x}[f(t)] = F(s) \quad K_{\Delta t}$ Transform

Rewriting the $K_{\Delta t}$ Transform Convolution Equation from Eq 5.10-30 and Eq 5.10-31

$$K_{\Delta t}[{}_{\Delta\lambda}\!\int_0^t f(t-\lambda-\Delta\lambda)g(\lambda)\Delta\lambda] = K_{\Delta t}[{}_{\Delta\lambda}\!\int_0^t f(\lambda)g(t-\lambda-\Delta\lambda)\Delta\lambda] = K_{\Delta t}[f(t)]\, K_{\Delta t}[g(t)] \qquad (5.10\text{-}82)$$

Applying Eq 5.10-81 to Eq 5.10-82

$$\frac{z}{T}\frac{z}{T}K_{\Delta t}\Big[\sum_{\Delta\lambda}\int_0^t f(t-\lambda-\Delta\lambda)g(\lambda)\Delta\lambda\Big]\Big|_{\substack{s=\frac{z-1}{T}\\ T=\Delta t\\ t=0,T,2T,3T,\dots}} = \frac{z}{T}\frac{z}{T}K_{\Delta t}\Big[\sum_{\Delta\lambda}\int_0^t f(t-\lambda-\Delta\lambda)g(\lambda)\Delta\lambda\Big]\Big|_{\substack{s=\frac{z-1}{T}\\ T=\Delta t\\ t=0,T,2T,3T,\dots}}$$

$$= \frac{z}{T}K_{\Delta t}[f(t)]\frac{z}{T}K_{\Delta t}[g(t)]\ \Big|_{\substack{s=\frac{z-1}{T}\\ T=\Delta t\\ t=0,T,2T,3T,\dots}} \qquad (5.10\text{-}83)$$

$$Z\Big[\frac{1}{T}\sum_{\Delta\lambda}\int_0^t f(t-\lambda-\Delta\lambda)g(\lambda)\Delta\lambda\Big] = Z\Big[\frac{1}{T}\sum_{\Delta\lambda}\int_0^t f(\lambda)g(t-\lambda-\Delta\lambda)\Delta\lambda\Big] = \frac{Z[f(t)]\,Z[g(t)]}{z} \qquad (5.10\text{-}84)$$

Then

**Another form of the Z Transform Convolution Equation is:**

$$\mathbf{Z\Big[\frac{1}{T}\sum_{T}\int_0^t f(t-\lambda-\Delta\lambda)g(\lambda)\Delta\lambda\Big] = \frac{Z[f(t)]Z[g(t)]}{z} = \frac{F(z)G(z)}{z}} \qquad (5.10\text{-}85)$$

**and**

$$\mathbf{Z\Big[\frac{1}{T}\sum_{T}\int_0^t f(\lambda)g(t-\lambda-\Delta\lambda)\Delta\lambda\Big] = \frac{Z[f(t)]Z[g(t)]}{z} = \frac{F(z)G(z)}{z}} \qquad (5.10\text{-}86)$$

**where**

**$f(t),g(t)$ = discrete functions of t**

**$f(t-\lambda-\Delta\lambda)$ = discrete function of t, $\lambda$, and $\Delta\lambda$**

**$\Delta\lambda = \Delta t = \Delta\lambda = T$**

**$t = n\Delta t$ ,    $n = 0, 1, 2, 3, \dots, \frac{t}{\Delta t}-1, \frac{t}{\Delta t}$**

**$\lambda = 0, \Delta\lambda, 2\Delta\lambda, 3\Delta\lambda, \dots, t-\Delta\lambda, t$**

**$Z[h(t)]\ ) = \frac{1}{T}\sum_{T}\int_0^\infty z^{-(\frac{t}{\Delta t})}h(t)\Delta t$ ,    The Z Transform of a function of t**

**$F(z) = Z[f(t)]$ = Z Transform, a function of z**

**$G(z) = Z[g(t)]$ = Z Transform, a function of z**

Below is an analysis of the Z Transform Convolution Equation, $\frac{1}{T}{}_T\int_{0}^{t+T} f(\lambda)g(t-\lambda)\Delta\lambda$, that will make it more understandable.

Consider the following transfer function block diagram of a Z Transform Transfer Function response to an input function, f(t).

f(t) = Function of t evaluated only at equally spaced discrete values of t, t = 0, T, 2T, 3T, …
δ(t) = Unit impulse function.

$f(t) = {}_1\sum_{n=0}^{\infty} f(n\Delta t)\delta(t-n\Delta t)$ $\quad g(t) = Z^{-1}[G(z)]$ $\quad O(t) = Z^{-1}[\,Z[f(t)]\,Z[g(t)]\,] =$

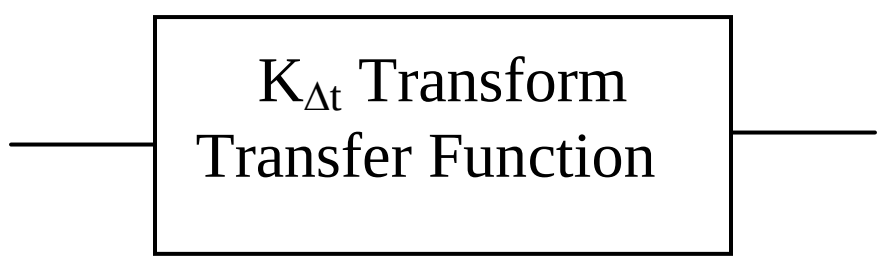

$${}_1\sum_{n=0}^{\infty} f(n\Delta t)g(t-n\Delta t)U(t-n\Delta t)$$

$Z[f(t)] = {}_1\sum_{n=0}^{\infty} f(n\Delta t)z^{-n}$ $\quad G(z) = Z[g(t)]$ $\quad O(z) = Z[f(t)]\,G(z) = {}_1\sum_{n=0}^{\infty} f(n\Delta t)z^{-n}G(z)$

Input Transfer Function Output

O(t) is the Z Transform Transfer Function response to a weighted series of unit impulses which is equivalent to the function, f(t), where t = 0, T, 2T, 3T, …

$$O(t) = Z^{-1}[\,Z[f(t)]\,Z[g(t)]\,] = {}_1\sum_{n=0}^{\infty} f(n\Delta t)g(t-n\Delta t)U(t-n\Delta t)\Delta t \qquad (5.10\text{-}87)$$

where
f(t) = Input function

$Z[f(t)] = {}_1\sum_{n=0}^{\infty} f(n\Delta t)z^{-n}$ , Input function Z Transform
G(z) = Z Transform Transfer Function
$g(t) = Z^{-1}[G(z)]$ , Inverse Z Transform of G(z)
O(z) = Z[O(t)] = Z[f(t)] G(z) , Output function $K_{\Delta t}$ Transform
O(t) = , Output response to the input function, f(t)
Δt = T = Interval between successive discrete values of t
t = 0, Δt, 2Δt, 3Δt, … = discrete independent variable values
n = 0, 1, 2, 3, …

Then, the Z Transform function output response, O(t), to the input function f(t) is:

$$O(t) = {}_1\sum_{n=0}^{\infty} f(n\Delta t)g(t-n\Delta t)U(t-n\Delta t) = Z^{-1}[\,Z[f(t)]\,Z[g(t)]\,] \qquad (5.10\text{-}88)$$

Let

$$\lambda = n\Delta t \qquad (5.10\text{-}89)$$

$$\Delta\lambda = \Delta t = T \qquad (5.10\text{-}90)$$

$$n = 0, 1, 2, 3, \ldots \qquad (5.10\text{-}91)$$

$$\lambda = 0, \Delta t, 2\Delta t, 3\Delta t, \ldots \qquad (5.10\text{-}92)$$

Substituting Eq 5.10-89 thru Eq 5.10-92 into Eq 5.10-88

$$O(t) = {}_{\Delta\lambda}\sum_{\lambda=0}^{\infty} f(\lambda)g(t-\lambda)U(t-\lambda) = Z^{-1}[\ Z[f(t)]\, Z[g(t)]\ ] \qquad (5.10\text{-}93)$$

Since

$$U(t-\lambda) = \begin{cases} 1 & \text{for } \lambda \le t \\ 0 & \text{for } \lambda > t \end{cases} \qquad (5.10\text{-}94)$$

Eq 5.9-93 can be simplified

$$O(t) = {}_{\Delta\lambda}\sum_{\lambda=0}^{t} f(\lambda)g(t-\lambda) = Z^{-1}[\ Z[f(t)]\, Z[g(t)]\ ] \qquad (5.10\text{-}95)$$

Changing the form of Eq 5.10-95 to its discrete integral form with $\Delta\lambda = T$

$$O(t) = \frac{1}{T}{}_{T}\int_{0}^{t+T} f(\lambda)g(t-\lambda)\Delta\lambda = Z^{-1}[\ Z[f(t)]\, Z[g(t)]\ ] \qquad (5.10\text{-}96)$$

From Eq 5.10-78 the $K_{\Delta t}$ Transform Convolution Equation is:

$$\frac{1}{T}{}_{T}\int_{0}^{t+T} f(\lambda)g(t-\lambda)\Delta\lambda = Z^{-1}[\ Z[f(t)]\, Z[g(t)]\ ] = O(t) \qquad (5.10\text{-}97)$$

Comparing Eq 5.9-96 to the Z Transform Convolution Equation, Eq 5.10-97, the two equations are seen to be the same. Thus, the Z Transform Convolution integral can now be described as follows:

1. The $K_{\Delta t}$ Transform Convolution integral, $O(t) = \frac{1}{T}{}_T\int_0^{t+T} f(\lambda)g(t-\lambda)\Delta\lambda$ , is the sum of all G(z) Transfer

   Function responses, ${}_1\sum_{n=0}^{\infty} f(n\Delta t)g(t-n\Delta t)U(t-n\Delta t)$, to all $f(n\Delta t)$ weighted unit impulses,

   ${}_1\sum_{n=0}^{\infty} f(n\Delta t)\delta(t-n\Delta t)$, that together comprise the input function, f(t). $\lambda = n\Delta t$, $\Delta\lambda = \Delta t$, and n = 0, 1, 2, 3, ...

2. The Z Transform Convolution integral term, $g(t-\lambda)$, if more completely written, $g(t-\lambda)U(t-\lambda)$ (see Eq 5.10-89), represents the output response of a Z Transform Transfer Function, G(z), to an input unit area impulse at $t = \lambda$.

3. The Z Transform Convolution integral term, $f(\lambda) = f(n\Delta t)$, is the value of the impulse series of f(t) in its nth interval. $f(\lambda)$ is the weighting placed on each consecutive unit area impulse of the f(t) infinite series where $\lambda = 0, \Delta\lambda, 2\Delta\lambda, 3\Delta\lambda, \ldots, t$ and $n = 0, 1, 2, 3, \ldots, \frac{t}{\Delta\lambda}$ . $\lambda = n\Delta t$, $\Delta\lambda = \Delta t$

To demonstrate the use of the Z Transform Convolution Equation, several examples are provided below.

<u>Example 5.10-3</u> Find the Inverse Z Transform, $Z^{-1}[\frac{z^2}{(z-1)^4}]$, using the Z Transform Convolution Equation.

The Z Transform Convolution Equation is:

$$Z[\frac{1}{T}{}_T\int_0^{t+T} f(t-\lambda)g(\lambda)\Delta\lambda] = Z[f(t)]Z[g(t)] \qquad 1)$$

where
- $\Delta\lambda = \Delta t = T$
- $\lambda = 0, \Delta\lambda, 2\Delta\lambda, 3\Delta\lambda, \ldots, t-\Delta\lambda, t, t+T$
- f(t),g(t) = discrete functions of t
- $f(t-\lambda)$ = discrete function of t and $\lambda$
- $t = 0, \Delta t, 2\Delta t, 3\Delta t, \ldots$
- Z[f(t)] = F(z) = Z Transform of f(t)
- $Z[f(t)] = \frac{1}{T}{}_T\int_0^{\infty} z^{-(\frac{t}{\Delta t})} f(t)\Delta t$ = the Z Transform of the function, f(t)

Find $Z^{-1}[\frac{z^2}{(z-1)^4}]$

Let

$$Z[f(t)] = \frac{z}{(z-1)^2} = \frac{1}{T}\frac{Tz}{(z-1)^2} \qquad 2)$$

$$Z[g(t)] = \frac{z}{(z-1)^2} = \frac{1}{T}\frac{Tz}{(z-1)^2} \qquad 3)$$

$$Z[f(t)]Z[g(t)]) = [\frac{z}{(z-1)^2}]^2 = \frac{z^2}{(z-1)^4} \qquad 4)$$

Writing the inverse Z Transforms of Eq 2) and Eq 3)

$$f(t) = \frac{t}{T} \qquad 5)$$

$$g(t) = \frac{t}{T} \qquad 6)$$

Substituting Eq 2 thru Eq 6 into Eq 1 and taking the Inverse Z Transform

$$Z^{-1}[\frac{z^2}{(z-1)^4}] = \frac{1}{T}\,{}_T\!\int_0^{t+T}\frac{t-\lambda}{T}\frac{\lambda}{T}\,\Delta\lambda = \frac{1}{T^3}\,{}_T\!\int_0^{t+T} t\lambda\,\Delta\lambda - \frac{1}{T^3}\,{}_T\!\int_0^{t+T}\lambda^2\,\Delta\lambda \qquad 7)$$

Evaluating Eq 7 using discrete integration

$$Z^{-1}[\frac{z^2}{(z-1)^4}] = \frac{1}{T^3}\,t\,{}_T\!\int_0^{t+T}\lambda\,\Delta\lambda - \frac{1}{T^3}\,{}_T\!\int_0^{t+T}[\lambda(\lambda-\Delta\lambda) + \lambda\Delta\lambda]\,\Delta\lambda \qquad 8)$$

$$Z^{-1}[\frac{z^2}{(z-1)^4}] = \frac{1}{T^3}\,t\,{}_T\!\int_0^{t+T}\lambda\,\Delta\lambda - \frac{1}{T^3}\,{}_T\!\int_0^{t+T}\lambda(\lambda-\Delta\lambda)\,\Delta\lambda - \frac{1}{T^3}\Delta\lambda\,{}_T\!\int_0^{t+T}\lambda\,\Delta\lambda \qquad 9)$$

$$Z^{-1}[\frac{z^2}{(z-1)^4}] = \frac{1}{T^3}[\frac{t\lambda(\lambda-\Delta\lambda)}{2}\Big|_0^{t+T} - \frac{\lambda(\lambda-\Delta\lambda)(\lambda-2\Delta\lambda)}{3}\Big|_0^{t+T} - \frac{\Delta\lambda\lambda(\lambda-\Delta\lambda)}{2}\Big|_0^{t+T}] \qquad 10)$$

$$T = \Delta\lambda \qquad 11)$$

Substituting Eq 11 into Eq 10

$$Z^{-1}[\frac{z^2}{(z-1)^4}] = \frac{1}{T^3}[\frac{t\lambda(\lambda-T)}{2}\Big|_0^{t+T} - \frac{\lambda(\lambda-T)(\lambda-2T)}{3}\Big|_0^{t+T} - \frac{T\lambda(\lambda-T)}{2}\Big|_0^{t+T}] \qquad 12)$$

$$Z^{-1}[\frac{z^2}{(z-1)^4}] = \frac{1}{T^3}[\frac{t^2(t+T)}{2} - \frac{(t+T)t(t-T)}{3} - \frac{T(t+T)t}{2}] \qquad 13)$$

$$Z^{-1}[\frac{z^2}{(z-1)^4}] = \frac{1}{T^3}[\frac{t^3}{2} + \frac{Tt^2}{2} - \frac{t^3}{3} + \frac{T^2t}{3} - \frac{Tt^2}{2} - \frac{T^2t}{2}] \qquad 14)$$

$$Z^{-1}[\frac{z^2}{(z-1)^4}] = \frac{1}{T^3}[\frac{t^3}{6} - \frac{T^2t}{6}] = \frac{1}{6T^3}[t^3 - T^2t)] \qquad 15)$$

Then

$$\mathbf{Z^{-1}[\frac{z^2}{(z-1)^4}] = \frac{1}{6}[(\frac{t}{T})^3 - \frac{t}{T}]} \qquad 16)$$

Checking Eq 16

Taking the Z Transform of both sides of Eq 16

$$\frac{z^2}{(z-1)^4} = Z[\frac{1}{6}\{(\frac{t}{T})^3 - \frac{t}{T}\}]$$

From the Z Transform Tables

$$Z^{-1}[t] = \frac{Tz}{(z-1)^2}$$

$$Z^{-1}[t^3] = \frac{T^3z(z^2+4z+1)}{(z-1)^4}$$

Substituting

$$\frac{z^2}{(z-1)^4} = \frac{1}{6}[\frac{T^3z(z^2+4z+1)}{T^3(z-1)^4} - \frac{Tz}{T(z-1)^2}]$$

$$\frac{z^2}{(z-1)^4} = \frac{1}{6}[\frac{z(z^2+4z+1)}{(z-1)^4} - \frac{z}{(z-1)^2}]$$

$$\frac{z^2}{(z-1)^4} = \frac{1}{6}[\frac{z^3+4z^2+z - z(z^2-2z+1)}{(z-1)^4}]$$

$$\frac{z^2}{(z-1)^4} = \frac{1}{6}\left[\frac{z^3+4z^2+z - z^3+2z^2-z}{(z-1)^4}\right]$$

$$\frac{z^2}{(z-1)^4} = \frac{z^2}{(z-1)^4}$$

Good check

Example 5.10-4 Find the particular solution to the following difference equation. Solve using Z Transforms and the Z Transform Convolution Equation. $y(0) = 0$ and $\Delta t = T$

$y(t+\Delta t) + 3y(t) = 1 \quad , \quad t = 0, \Delta t, 2\Delta t, 3\Delta t, \ldots \quad ,$ Difference equation 1)

Taking the Z Transform of Eq 1

$Z[y(t+\Delta t)] + 3Z[y(t)] = Z[1]$ 2)

$zY(z) - zy(0) + 3Y(z) = \frac{z}{z-1}$ 3)

$(z+3)Y(z) - 0z = \frac{z}{z-1}$ 4)

$Y(z) = Z[y(t)] = [\frac{1}{z+3}]\,[\frac{z}{z-1}]$ 5)

Taking the Inverse Z Transform of Eq 5

$y(t) = Z^{-1}[\,[\frac{1}{z+3}]\,[\frac{z}{z-1}]\,]$ 6)

The Z Transform Convolution Equation is:

$Z[\frac{1}{T}\,{}_T\!\int_0^{t+T} f(t-\lambda)g(\lambda)\Delta\lambda] = Z[f(t)]Z[g(t)]$ 7)

where

$\Delta\lambda = \Delta t = T$
$\lambda = 0, \Delta\lambda, 2\Delta\lambda, 3\Delta\lambda, \ldots, t-\Delta\lambda, t, t+T$
$f(t),g(t)$ = discrete functions of t
$f(t-\lambda)$ = discrete function of t and $\lambda$
$t = 0, \Delta t, 2\Delta t, 3\Delta t, \ldots$
$Z[f(t)] = F(z)$ = Z Transform of f(t)

$$Z[f(t)] = \frac{1}{T}\ {}_{T}\int_{0}^{\infty} z^{-\left(\frac{t}{\Delta t}\right)} f(t)\Delta t = \text{the Z Transform of the function, f(t)}$$

Find y(t) using the Z Transform Convolution Equation, Eq 7.

From Eq 6 and Eq 7

Let

$$Z[f(t)] = \frac{1}{z+3} = z^{-1}\frac{z}{z+3} \qquad 8)$$

Finding the Inverse Z Transform of Eq 8 using the Z Transform tables where $n = \frac{t}{T}$

$$f(t) = (-3)^{\frac{t}{T}-1}U(\tfrac{t}{T}-1) = (-3)^{\frac{t-T}{T}}U(\tfrac{t-T}{T}) \qquad 9)$$

Let

$$Z[g(t)] = \frac{z}{z-1} \qquad 10)$$

Finding the Inverse Z Transform of Eq 10 using the Z Transform tables

$$g(t) = 1 \qquad 11)$$

Taking the Inverse Z Transform of Eq 7

$$\frac{1}{T}\ {}_{T}\int_{0}^{t+T} f(t-\lambda)g(\lambda)\Delta\lambda = Z^{-1}[\ Z[f(t)]Z[g(t)]\ ] \qquad 12)$$

Substituting Eq 8 thru Eq 11 into Eq 12

$$\frac{1}{T}\ {}_{T}\int_{0}^{t+T} (-3)^{\frac{t-\lambda-T}{T}}U(\tfrac{t-\lambda-T}{T})(1)\Delta\lambda = Z^{-1}[\ [\frac{1}{z+3}][\frac{z}{z-1}]\ ] \qquad 13)$$

where

$\lambda = 0, \Delta\lambda, 2\Delta\lambda, 3\Delta\lambda, \ldots, t-\lambda, t, t+T$

$\Delta\lambda = T$

From Eq 6 and Eq 13

$$y(t) = \frac{1}{T} {}_T\int_0^{t+T} (-3)^{\frac{t-\lambda-T}{T}} U\left(\frac{t-\lambda-T}{T}\right)\Delta\lambda \qquad 14)$$

From the function, $U(\frac{t-\lambda-T}{T})$, it is seen that the integration includes the summation of terms from $\lambda = 0$ thru $\lambda = t-T$. Thus, the actual upper limit of the integration is $\lambda = t-T+T = t$.

Then, simplifying Eq 14

$$y(t) = \frac{1}{T} {}_T\int_0^{t} (-3)^{\frac{t-\lambda-T}{T}} \Delta\lambda \qquad 15)$$

$$y(t) = \frac{1}{T} (-3)^{\frac{t}{T}-1} {}_T\int_0^{t} (-3)^{\frac{-\lambda}{T}} \Delta\lambda \qquad 16)$$

Integrating Eq 16

Use the following integration formula found in Table 6 of the Appendix

$$_{\Delta x}\int c^{-\frac{x}{\Delta x}} \Delta x = \frac{\Delta x}{c^{-1}-1} c^{-\frac{x}{\Delta x}} + k \qquad 17)$$

where
c = constant
$c^{-1}-1 \neq 0$

From Eq 16 and Eq 17

$$y(t) = \frac{1}{T} (-3)^{\frac{t}{T}-1} {}_T\int_0^{t} (-3)^{-\frac{\lambda}{T}} \Delta\lambda = \frac{1}{T} (-3)^{\frac{t}{T}-1} \frac{T}{(-3)^{-1}-1} (-3)^{-\frac{\lambda}{T}} \Big|_0^t \qquad 18)$$

Simplifying Eq 18

$$y(t) = \frac{1}{T} (-3)^{\frac{t}{T}} \left(\frac{1}{-3}\right) \frac{(-3)T}{4} \left[(-3)^{-\frac{t}{T}} - 1\right] = \frac{1}{4}\left[1 - (-3)^{\frac{t}{T}}\right] \qquad 19)$$

Then

$$\mathbf{y(t) = \frac{1}{4}\,[1 - (-3)^{\frac{t}{T}}]} \qquad 20)$$

Checking Eq 20

$\Delta t = T$

$y(t+T) + 3y(t) = 1$

Substituting Eq 20 into the above difference equation

$$\frac{1}{4}\,[1 - (-3)^{\frac{t}{T}+1}] + (3)\,\frac{1}{4}\,[1 - (-3)^{\frac{t}{T}}] = 1$$

$$\frac{1}{4} + \frac{3}{4}\,(-3)^{\frac{t}{T}} + \frac{3}{4} - \frac{3}{4}\,(-3)^{\frac{t}{T}} = 1$$

$1 = 1$

Good check

<u>Example 5.10-5</u> Find the Inverse Z Transform, $Z^{-1}[\frac{1}{z-1}\,C(z)]$

$$\text{Find } Z^{-1}[\frac{1}{z-1}\,C(z)] \qquad (1)$$

Writing the following Z Transform Convolution Equation

$$\frac{1}{T}\,{}_{T}\!\int_{0}^{t+T} f(t-\lambda)g(\lambda)\Delta\lambda = Z^{-1}\,[F(z)G(z)] \qquad (2)$$

where

$\Delta\lambda = \Delta t = T$

From Eq 1 and Eq 2 where $F(z) = \frac{1}{z-1}$ , $G(z) = C(z)$

$$f(t) = U(t-T) \qquad (3)$$

$$g(t) = c(t) \qquad (4)$$

$$Z^{-1}\,[\frac{1}{z-1}\,C(z)] = \frac{1}{T}\,{}_{T}\!\int_{0}^{t+T} U(t-\lambda-\Delta\lambda)c(\lambda)\Delta\lambda = \frac{1}{T}\,{}_{T}\!\int_{0}^{t} U(t-\lambda-\Delta\lambda)c(\lambda)\Delta\lambda = \frac{1}{T}\,{}_{T}\!\int_{0}^{t} c(\lambda)\Delta\lambda \qquad (5)$$

where

$\Delta\lambda = \Delta t = T$

Comment - The last value of the above discrete integral summation is for $\lambda = t-T$

Checking

$$\mathbf{Z^{-1}[\frac{1}{z-1}C(z)] = \frac{1}{T}\,{}_T\!\int_0^t c(\lambda)\Delta\lambda} \qquad (6)$$

**where**

$\mathbf{\Delta\lambda = \Delta t = T}$

$\mathbf{c(\lambda) = Z^{-1}[C(z)]}$

Let $C(z) = \frac{Tz}{z-1}$

$c(t) = T$

$$Z^{-1}[\frac{Tz}{(z-1)^2}] = \frac{1}{T}\,{}_T\!\int_0^t c(\lambda)\Delta\lambda = \frac{1}{T}\,{}_T\!\int_0^t T\Delta\lambda = {}_T\!\int_0^t \Delta\lambda = \lambda\Big|_0^t = t$$

$$Z^{-1}[\frac{Tz}{(z-1)^2}] = t$$

Good check

In this section $K_{\Delta t}$ Transform and Z Transform Convolution Equations were derived and applied to several demonstration problems. Like the Laplace Transform Convolution Equations, these equations were found to be very useful in the solution of challenging mathematical problems. In addition to these convolution equations, there are several other related equations that are of considerable use. They are the Laplace Transform and Z Transform Duhamel Equations. In the following section, Section 5.11, Duhamel Equations will be derived for both the $K_{\Delta t}$ Transform and the Z Transform.

## Section 5.11: The Interval Calculus Duhamel Formulas

Of significant importance in Operational Calculus are the Duhamel Equations that are shown below.

The Laplace Transform Duhamel Equations

System Diagram

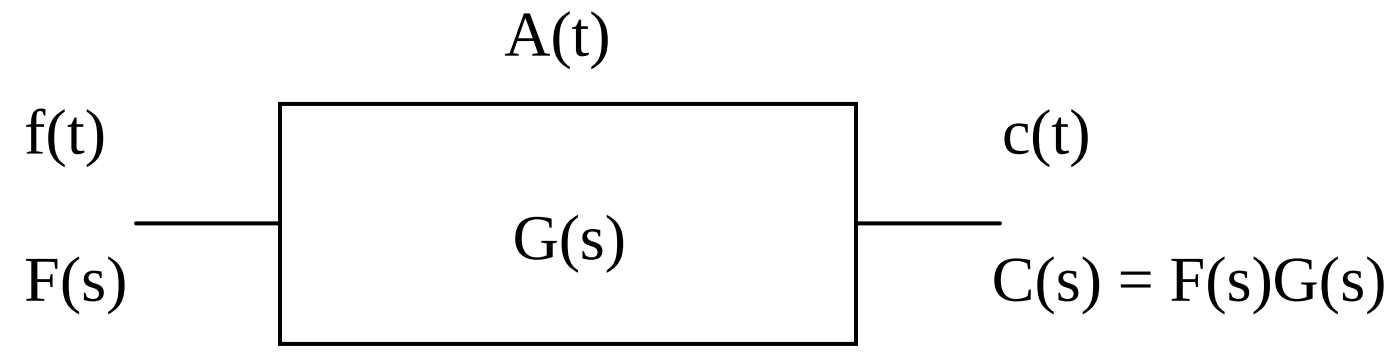

$$A(s) = \frac{1}{s}\,G(s)$$

Input Transfer Function Output

1. $c(t) = A(0)f(t) + \int_0^t f(\lambda)A'(t-\lambda)d\lambda$ , $A(0) = 0$ since the system is, by definition, initially passive. (5.11-1)

See Eq 5.11-3.

2. $c(t) = A(0)f(t) + \int_0^t f(t-\lambda)A'(\lambda)d\lambda$ , $A(0) = 0$ since the system is, by definition, initially passive. (5.11-2)

See Eq 5.11-4.

3. $c(t) = \int_0^t f(\lambda)A'(t-\lambda)d\lambda$ (5.11-3)

4. $c(t) = \int_0^t f(t-\lambda)A'(\lambda)d\lambda$ (5.11-4)

5. $c(t) = f(0)A(t) + \int_0^t f'(\lambda)A(t-\lambda)d\lambda$ (5.11-5)

6. $c(t) = f(0)A(t) + \int_0^t f'(t-\lambda)A(\lambda)d\lambda$ (5.11-6)

7. $c(t) = \frac{d}{dt}\int_0^t f(\lambda)A(t-\lambda)d\lambda$ (5.11-7)

8. $c(t) = \frac{d}{dt}\int_0^t f(t-\lambda)A(\lambda)d\lambda$ (5.11-8)

where

$f(t)$ = System input excitation

$c(t)$ = System response to the input excitation, $f(t)$

$A(t)$ = System response to a unit step function

$A(s)$ = Laplace Transform of the system response to a unit step function initiated at $t = 0$.

$G(s)$ = System Laplace Transform transfer function

$\frac{1}{s}$ = Laplace Transform of a unit step function initiated at $t = 0$

$A(s) = \frac{1}{s} G(s)$

In this digital age, the sampling of system variables is common. The results of this type of variable sampling can be represented as sample and hold shaped waveforms. Interval Calculus efficiently manipulates and analyzes these type of waveforms. Thus, in the derivations that follow, the Duhamel Equations specified above will be generalized for use in Interval Calculus.

The following equations are the $K_{\Delta t}$ Transform Duhamel Equations that will later be derived. Note that they appear to be very similar to the previously listed Laplace Transform Durhamel Equations. However, the mathematical operations that will be necessary for equation solution will be those of discrete variable Interval Calculus.

Derivation of the Interval Calculus Duhamel Equations using $K_{\Delta t}$ Transforms

For the analysis of discrete variable systems, it has been shown in the previous section, Section 5.10, that the derived $K_{\Delta t}$ Transform Convolution Equation can be very useful. However, there  are other equations derived from the $K_{\Delta t}$ Transform Convolution Equation that can also be very useful. These equations are derived in this section and named the Interval Calculus Duhamel Equations since they are generalizations of the Duhamel Equations developed by the French mathematician, J.M.C. Duhamel. Durhamel’s Equations are derived from the Laplace Transform and are used in the analysis of continuous variable systems. However, the Interval Calculus Durhamel Equations are derived either by the use of $K_{\Delta t}$ Transforms or the related Z Transforms and are used in the analysis of discrete variable systems. An advantage of using Durhamel Equations for analysis is the use of the transfer function, $A(s)$ or $A(z)$. The function, $A(t)$, used in the equations can be easily obtained. For a given system with an initially passive  transfer function, $G(s)$ or $G(z)$, $A(t)$ can be generated and recorded by simply exciting the system with a unit step function input initiated at $t = 0$. For many types of inputs, this could be performed with a simple switch and a constant voltage or current source.

The Interval Calculus Duhamel Equations listed below are the same whether derived by $K_{\Delta t}$ Transforms or Z Transforms. In this section each of the eight Interval Calculus Duhamel Equations is derived twice, initially using $K_{\Delta t}$ Transforms and later using Z Transforms.

Interval Calculus Duhamel Equations

$K_{\Delta t}$ Transform System Diagram

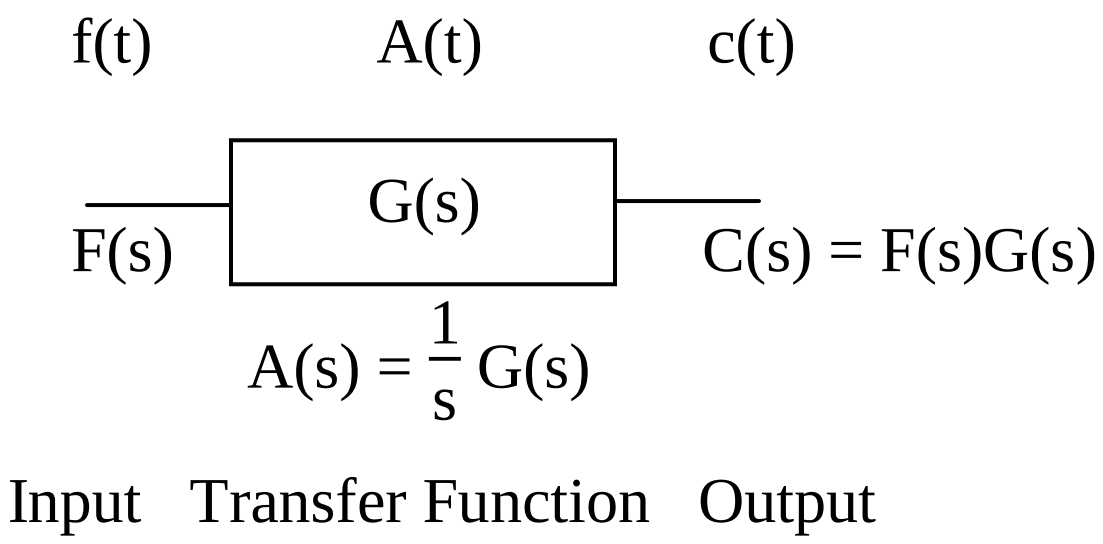

Z Transform System Diagram

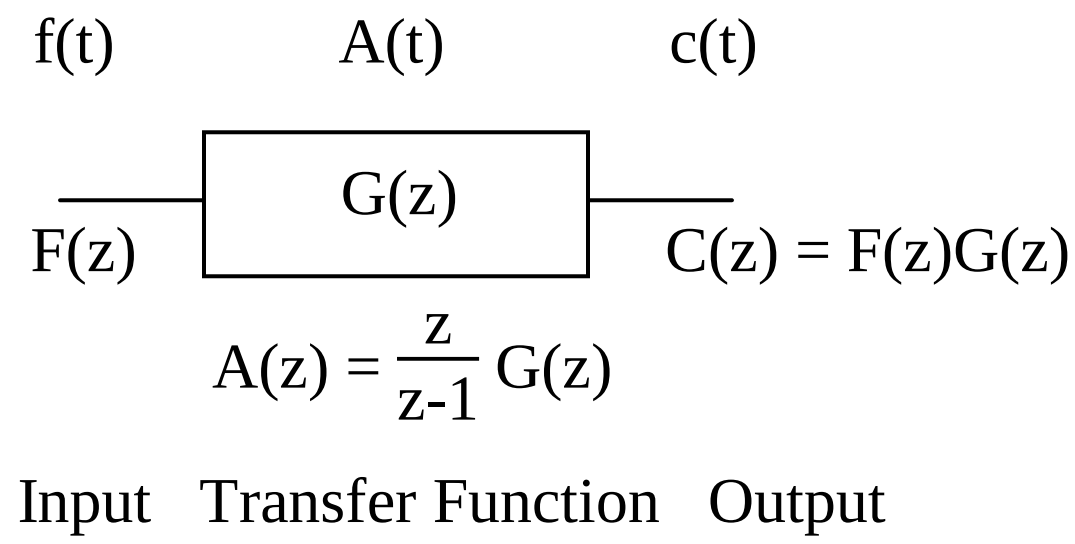

The following $K_{\Delta t}$ Transform and Z Transform Duhamel Equations make it possible to obtain the system response to a general input function, f(t), if the unit step function system response, A(t), is known.

1. $c(t) = A(0)f(t) + {}_{\Delta\lambda}\int_0^t f(\lambda)A^{'}(t-\lambda-\Delta\lambda)\Delta\lambda$ , $\quad C(s) = F(s)[\frac{sG(s)}{s}] = F(s)[sA(s)]$ (5.11-9)

$$C(z) = F(z)[\frac{z-1}{z}][\frac{z}{z-1}G(z)] = F(z)[\frac{z-1}{z}A(z)]$$

Since A(t) is by definition the response of a system that is initially passive, A(0) = 0.
See Eq 5.11-11.

2. $c(t) = A(0)f(t) + {}_{\Delta\lambda}\int_0^t f(t-\lambda-\Delta\lambda)A^{'}(\lambda)\Delta\lambda$ , $\quad C(s) = F(s)\ [\frac{sG(s)}{s}] = F(s)[sA(s)]$ (5.11-10)

$$C(z) = F(z)[\frac{z-1}{z}][\frac{z}{z-1}G(z)] = F(z)[\frac{z-1}{z}A(z)]$$

Since A(t) is by definition the response of a system that is initially passive, A(0) = 0.

See Eq 5.11-12.

3. $c(t) = {}_{\Delta\lambda}\int_0^t f(\lambda)A^{'}(t-\lambda-\Delta\lambda)\Delta\lambda$ $\quad C(s) = F(s)\ [\frac{sG(s)}{s}] = F(s)[sA(s)]$ (5.11-11)

$$C(z) = F(z)[\frac{z-1}{z}][\frac{z}{z-1}G(z)] = F(z)[\frac{z-1}{z}A(z)]$$

4. $c(t) = {}_{\Delta\lambda}\int_0^t f(t-\lambda-\Delta\lambda)A^{'}(\lambda)\Delta\lambda$ $\quad C(s) = F(s)\ [\frac{sG(s)}{s}] = F(s)[sA(s)]$ (5.11-12)

$$C(z) = F(z)[\frac{z-1}{z}][\frac{z}{z-1}G(z)] = F(z)[\frac{z-1}{z}A(z)]$$

5. $c(t) = f(0)A(t) + \underset{\Delta\lambda}{\int_0^t} f'(\lambda)A(t-\lambda-\Delta\lambda)\Delta\lambda$ $\qquad C(s) = [sF(s)]\frac{G(s)}{s} = [sF(s)]A(s)$ (5.11-13)

$$C(z) = [\frac{z-1}{z}F(z)][\frac{z}{z-1}G(z)] = [\frac{z-1}{z}F(z)]A(z)$$

6. $c(t) = f(0)A(t) + \underset{\Delta\lambda}{\int_0^t} f'(t-\lambda-\Delta\lambda)A(\lambda)\Delta\lambda$ $\qquad C(s) = [sF(s)]\frac{G(s)}{s} = [sF(s)]A(s)$ (5.11-14)

$$C(z) = [\frac{z-1}{z}F(z)][\frac{z}{z-1}G(z)] = [\frac{z-1}{z}F(z)]A(z)$$

7. $c(t) = D_{\Delta t}\underset{\Delta\lambda}{\int_0^t} f(\lambda)A(t-\lambda-\Delta\lambda)\Delta\lambda$ $\qquad C(s) = s[F(s)\frac{G(s)}{s}] = s[F(s)A(s)]$ (5.11-15)

$$C(z) = \frac{z-1}{z}[F(z)\frac{z}{z-1}G(z)] = \frac{z-1}{z}[F(z)A(z)]$$

8. $c(t) = D_{\Delta t}\underset{\Delta\lambda}{\int_0^t} f(t-\lambda-\Delta\lambda)A(\lambda)\Delta\lambda$ $\qquad C(s) = s[F(s)\frac{G(s)}{s}] = s[F(s)A(s)]$ (5.11-16)

$$C(z) = \frac{z-1}{z}[F(z)\frac{z}{z-1}G(z)] = \frac{z-1}{z}[F(z)A(z)]$$

where

$F(s) = K_{\Delta t}$ Transform of the system input excitation function, f(t)
f(t) = System input excitation function
$G(s)$ = System $K_{\Delta t}$ Transform transfer function
$A(s) = \frac{1}{s}G(s) = K_{\Delta t}$ Transform of the system response to a unit step function
A(t) = System response to a unit step function
A(0) = 0, The system, by definition, is initally passive
$A'(t) = D_{\Delta t}A(t)$ = Discrete derivative of A(t)
$C(s) = K_{\Delta t}$ Transform of the system output function, c(t)
c(t) = System output function
F(z) = Z Transform of the system input excitation function, f(t)
G(z) = System Z Transform transfer function
$A(z) = \frac{z}{z-1}G(z)$ = Z Transform of the system response to a unit step function
C(z) = Z Transform of the system output function, c(t)
$\Delta\lambda = \Delta t$
$\lambda = n\Delta t, \quad n = 0, 1, 2, 3, \ldots, \frac{t}{\Delta t}-1, t$
$\lambda = 0, \Delta\lambda, 2\Delta\lambda, 3\Delta\lambda, \ldots, t-\Delta\lambda, t$

$f'(\lambda) = D_{\Delta\lambda}f(\lambda) = \frac{f(\lambda+\Delta\lambda)-f(\lambda)}{\Delta\lambda}$ = discrete derivative of $f(\lambda)$

$\Delta f(\lambda) = f(\lambda+\Delta\lambda)-f(\lambda) = f'(\lambda)\Delta\lambda$

Comment - By definition, the discrete integral, $\Delta\lambda\int_0^t h(\lambda)\Delta\lambda$, is the sum, $\Delta\lambda\sum_{\lambda=0}^{t-\Delta\lambda} h(\lambda)\Delta\lambda$. Though the integration is from $\lambda = 0$ thru t, the summation that defines this integration is from $\lambda = 0$ thru $t-\Delta\lambda$.

The functions previously specified, f(t), $f'(t)$, c(t), are discrete variable sample and hold shaped waveforms. They can be represented by the following two series.

Two series representations of the sample and hold shaped waveform function, r(t)

Using unit amplitude $\Delta t$ width pulses

$$r(t) = {}_1\sum_{n=0}^{\frac{t}{\Delta t}-1} r(n\Delta t)[U(t-n\Delta t)-U(t-[n+1]\Delta t)],$$ A series of unit amplitude pulses of $\Delta t$ width weighted by $r(n\Delta t)$

(5.11-17)

or in a different form where $\lambda = n\Delta t$ and $\Delta\lambda = \Delta t$

$$r(t) = {}_{\Delta\lambda}\sum_{\lambda=0}^{t-\Delta\lambda} r(\lambda)[U(t-\lambda) - U(t-\lambda-\Delta\lambda)],$$ A series of unit amplitude pulses of $\Delta\lambda$ width weighted by $r(\lambda)$

(5.11-18)

Using unit step functions

$$r(t) = f(0)U(t) + {}_1\sum_{n=0}^{\frac{t}{\Delta t}-1} r'(n\Delta t)U(t-[n+1]\Delta t)\Delta t,$$ A series of unit step functions weighted by $r'(n\Delta t)$ (5.11-19)

or in a different form where $\lambda = n\Delta t$ and $\Delta\lambda = \Delta t$

$$r(t) = f(0)U(t) + {}_{\Delta\lambda}\sum_{\lambda=0}^{t-\Delta\lambda} r'(\lambda)U(t-\lambda-\Delta\lambda)\Delta\lambda,$$ A series of unit step functions weighted by $r'(\lambda)$ (5.11-20)

System response to a unit amplitude $\Delta t$ width pulse input excitation

$c(t) = K_{\Delta t}^{-1}[(1+s\Delta t)^{-1-n}\Delta t\, A(s)] = A(t-n\Delta t-\Delta t)\Delta t = A(t-\lambda-\Delta\lambda)\Delta\lambda$ = System response to a unit amplitude $\Delta t$ width pulse input excitation initiated at $t = n\Delta t = \lambda$ (5.11-21)

where
- $n = 0, 1, 2, 3, \ldots$
- $\lambda = 0, \Delta\lambda, 2\Delta\lambda, 3\Delta\lambda, \ldots$
- $\lambda = n\Delta t$
- $\Delta\lambda = \Delta t$

System response to a unit step function input excitation

$c(t) = K_{\Delta t}^{-1}[\frac{(1+s\Delta t)^{-n}}{s} G(s)] = K_{\Delta t}^{-1}[(1+s\Delta t)^{-n}A(s)] = A(t-n\Delta t) = A(t-\lambda)$ = System response to a unit step input excitation at $t = n\Delta t = \lambda$ (5.11-22)

where
- $n = 0, 1, 2, 3, \ldots$
- $\lambda = 0, \Delta\lambda, 2\Delta\lambda, 3\Delta\lambda, \ldots$
- $\lambda = n\Delta t$
- $\Delta\lambda = \Delta t$

Derivation of the Interval Calculus Duhamel Equations using $K_{\Delta t}$ Transforms

For the $K_{\Delta t}$ Transform equation, $C(s) = F(s)[sA(s)]$

System Diagram

f(t) A(t) c(t)

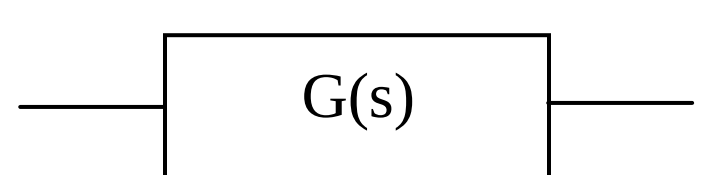

F(s) $C(s) = F(s)G(s) = F(s)s\frac{G(s)}{s} = F(s)[sA(s)]$

$A(s) = \frac{1}{s} G(s)$

Input Transfer Function Output

$$C(s) = F(s)[sA(s)] = F(s)[sA(s) - A(0) + A(0)] = A(0)F(s) + F(s)[sA(s) - A(0)] \quad (5.11\text{-}23)$$

$$C(s) = A(0)F(s) + K_{\Delta t}[f(t)]\, K_{\Delta t}[\, A'(t)] \quad (5.11\text{-}24)$$

where

$$A'(t) = D_{\Delta t}A(t) = \frac{A(t-\Delta t) - A(t)}{\Delta t} = \text{Discrete derivative of } A(t)$$

$$F(s) = K_{\Delta t}[f(t)]$$

Taking the Inverse $K_{\Delta t}$ Transform of Eq 5.11-24 and using the $K_{\Delta t}$ Transform Convolution Equations, Eq 5.10-32 and Eq 5.10-34, to find c(t)

$$c(t) = A(0)f(t) + K_{\Delta t}^{-1}[\, K_{\Delta t}[f(t)]\, K_{\Delta t}[\, A'(t)]\, ] = A(0)f(t) + {}_{\Delta\lambda}\int_0^t f(t-\lambda-\Delta\lambda)A'(\lambda)\Delta\lambda \qquad (5.11\text{-}25)$$

and

$$c(t) = A(0)f(t) + {}_{\Delta\lambda}\int_0^t f(\lambda)A'(t-\lambda-\Delta\lambda)\Delta\lambda \qquad (5.11\text{-}26)$$

Then

$$\mathbf{c(t) = A(0)f(t) + {}_{\Delta\lambda}\int_0^t f(t-\lambda-\Delta\lambda)A'(\lambda)\Delta\lambda} \qquad (5.11\text{-}27)$$

**where**

**f(t) = System input excitation function**
**c(t) = System output function**
**A(t) = System response to a unit step function input**
**A(0) = 0, The system, by definition, is initally passive**
**$A'(t) = D_{\Delta t}A(t)$ = discrete derivative of A(t)**
**$A'(s) = sA(s) - A(0) = sA(s) = s\,\frac{G(s)}{s} = G(s)$**
**$A(s) = \frac{G(s)}{s}$**
**$\Delta\lambda = \Delta t$**
**$\lambda = n\Delta t, \quad n = 0, 1, 2, 3, \ldots, \frac{t}{\Delta t}-1, \frac{t}{\Delta t}$**
**$\lambda = 0, \Delta\lambda, 2\Delta\lambda, 3\Delta\lambda, \ldots, t-\Delta\lambda, t$**
**c(t), f(t) = Sample and hold shaped waveforms, Interval Calculus discrete variable functions**

Since A(t) is by definition the response of a system that is initially passive, A(0) = 0.

$$\mathbf{c(t) = {}_{\Delta\lambda}\int_0^t f(t-\lambda-\Delta\lambda)A'(\lambda)\Delta\lambda} \qquad (5.11\text{-}28)$$

**where**

**f(t) = System input excitation function**
**c(t) = System output function**
**A(t) = System response to a unit step function input**
**A(0) = 0, The system, by definition, is initally passive**

**$A'(t) = D_{\Delta t}A(t)$ = discrete derivative of A(t)**

**$A'(s) = sA(s) - A(0) = sA(s) = s\,\dfrac{G(s)}{s} = G(s)$**

**$A(s) = \dfrac{G(s)}{s}$**

**$\Delta\lambda = \Delta t$**

**$\lambda = n\Delta t, \quad n = 0, 1, 2, 3, \ldots, \dfrac{t}{\Delta t}-1, \dfrac{t}{\Delta t}$**

**$\lambda = 0, \Delta\lambda, 2\Delta\lambda, 3\Delta\lambda, \ldots, t-\Delta\lambda, t$**

**c(t), f(t) = Sample and hold shaped waveforms, Interval Calculus discrete variable functions**

**and**

$$c(t) = A(0)f(t) + {}_{\Delta\lambda}\int_{0}^{t} f(\lambda)A'(t-\lambda-\Delta\lambda)\,\Delta\lambda \qquad (5.11\text{-}29)$$

**where**

**f(t) = System input excitation function ( See Eq 5.11-18 where f(t) = r(t) )**
**c(t) = System output function**
**A(t) = System response to a unit step function input excitation**
**A(0) = 0, The system, by definition, is initally passive**
**$A'(t) = D_{\Delta t}A(t)$ = discrete derivative of A(t)**

**$A'(s) = sA(s) - A(0) = sA(s) = s\,\dfrac{G(s)}{s} = G(s)$**

**$A(s) = \dfrac{G(s)}{s}$**

**$\Delta\lambda = \Delta t$**

**$\lambda = n\Delta t, \quad n = 0, 1, 2, 3, \ldots, \dfrac{t}{\Delta t}-1, \dfrac{t}{\Delta t}$**

**$\lambda = 0, \Delta\lambda, 2\Delta\lambda, 3\Delta\lambda, \ldots, t-\Delta\lambda, t$**

**$A'(t-\lambda-\Delta\lambda)\Delta\lambda$ = System response to a unit amplitude $\Delta\lambda$ width pulse input initiated at $t = \lambda$ ( See Eq 5.11-21 where $A'(t)$ is substituted for A(t))**

**c(t), f(t) = Sample and hold shaped waveforms, Interval Calculus discrete variable functions**

Since A(t) is by definition the response of a system that is initially passive, A(0) = 0.

$$c(t) = {}_{\Delta\lambda}\int_{0}^{t} f(\lambda)A'(t-\lambda-\Delta\lambda)\,\Delta\lambda \qquad (5.11\text{-}30)$$

**where**

**f(t) = System input function ( See Eq 5.11-18 where f(t) = r(t) )**
**c(t) = System output function**
**A(t) = System response to a unit step function input**
**A(0) = 0, The system, by definition, is initally passive**
**$A'(t) = D_{\Delta t}A(t)$ = discrete derivative of A(t)**

$$\mathbf{A'(s) = sA(s) - A(0) = sA(s) = s\frac{G(s)}{s} = G(s)}$$

$$\mathbf{A(s) = \frac{G(s)}{s}}$$

$$\mathbf{\Delta\lambda = \Delta t}$$

$$\mathbf{\lambda = n\Delta t, \quad n = 0, 1, 2, 3, \ldots, \frac{t}{\Delta t}-1, \frac{t}{\Delta t}}$$

$$\mathbf{\lambda = 0, \Delta\lambda, 2\Delta\lambda, 3\Delta\lambda, \ldots, t-\Delta\lambda, t}$$

**$A'(t-\lambda-\Delta\lambda)\Delta\lambda$ = System response to a unit amplitude $\Delta\lambda$ width pulse input initiated at t = $\lambda$**
**( See Eq 5.11-21 where $A'(t)$ is substituted for A(t))**

**c(t), f(t) = Sample and hold shaped waveforms, Interval Calculus discrete variable functions**

For the $K_{\Delta t}$ Transform equation, C(s) = [sF(s)]A(s)

System Diagram

f(t) A(t) c(t)

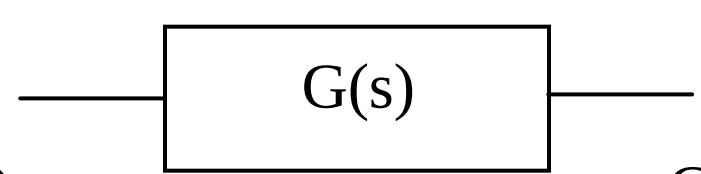

F(s) $C(s) = F(s)G(s) = sF(s)\frac{G(s)}{s} = [sF(s)]A(s)$

$A(s) = \frac{1}{s}G(s)$

Input Transfer Function Output

$$C(s) = F(s)G(s) = [sF(s)]A(s) = [sF(s) - f(0)]A(s) + f(0)A(s) = K_{\Delta t}[f'(t)]A(s) + f(0)A(s) \qquad (5.11\text{-}31)$$

$$C(s) = f(0)A(s) + K_{\Delta t}[f'(t)]\,K_{\Delta t}[A(t)] \qquad (5.11\text{-}32)$$

where

$$f'(t) = D_{\Delta t}f(t) = \frac{f(t-\Delta t) - f(t)}{\Delta t} = \text{Discrete derivative of } f(t)$$

$$A(s) = K_{\Delta t}[A(t)]$$

Taking the Inverse $K_{\Delta t}$ Transform of Eq 5.11-32 and using the $K_{\Delta t}$ Transform Convolution Equations, Eq 5.10-32 and Eq 5.10-34, to find c(t)

$$c(t) = f(0)A(t) + K_{\Delta t}^{-1}[\,K_{\Delta t}[f'(t)]K_{\Delta t}[A(t)]\,] = f(0)A(t) + {}_{\Delta\lambda}\int_0^t f'(t-\lambda-\Delta\lambda)A(\lambda)\Delta\lambda \qquad (5.11\text{-}33)$$

and

$$c(t) = f(0)A(t) + {}_{\Delta\lambda}\int_0^t f'(t-\lambda-\Delta\lambda)A(\lambda)\Delta\lambda \qquad (5.11\text{-}34)$$

Then

$$c(t) = f(0)A(t) + {}_{\Delta\lambda}\int_{0}^{t} f'(t-\lambda-\Delta\lambda)A(\lambda)\Delta\lambda \qquad (5.11\text{-}35)$$

**where**

**$f(t)$ = System input excitation function**
**$c(t)$ = System output function**
**$A(t)$ = System response to a unit step function input**
**$f'(t) = D_{\Delta t}f(t)$ = discrete derivative of $f(t)$**

$A(s) = \frac{G(s)}{s}$

$\Delta\lambda = \Delta t$

$\lambda = n\Delta t, \quad n = 0, 1, 2, 3, \ldots, \frac{t}{\Delta t}-1, \frac{t}{\Delta t}$

$\lambda = 0, \Delta\lambda, 2\Delta\lambda, 3\Delta\lambda, \ldots, t-\Delta\lambda, t$

**$c(t), f(t)$ = Sample and hold shaped waveforms, Interval Calculus discrete variable functions**

**and**

$$c(t) = f(0)A(t) + {}_{\Delta\lambda}\int_{0}^{t} f'(\lambda)A(t-\lambda-\Delta\lambda)\Delta\lambda \qquad (5.11\text{-}36)$$

**where**

**$f(t)$ = System input excitation function ( See Eq 5.11-20 where $f(t) = r(t)$ )**
**$c(t)$ = System output function**
**$A(t)$ = System response to a unit step function input**
**$f'(t) = D_{\Delta t}f(t)$ = discrete derivative of $f(t)$**

$A(s) = \frac{G(s)}{s}$

$\Delta\lambda = \Delta t$

$\lambda = n\Delta t, \quad n = 0, 1, 2, 3, \ldots, \frac{t}{\Delta t}-1, \frac{t}{\Delta t}$

$\lambda = 0, \Delta\lambda, 2\Delta\lambda, 3\Delta\lambda, \ldots, t-\Delta\lambda, t$

**$A(t-n\Delta t-\Delta t)$ = System response to a unit step function initiated at $t = [n+1]\Delta t$,**

**( See Eq 5.11-22), $\quad n = 0, 1, 2, 3, \ldots, \frac{t}{\Delta t}-1$**

$f'(n\Delta t)\Delta t = \Delta f(n\Delta t) = [f(n\Delta t+\Delta t)-f(n\Delta t)]$

**$c(t), f(t)$ = Sample and hold shaped waveforms, Interval Calculus discrete variable functions**

For the $K_{\Delta t}$ Transform equation, $C(s) = s[F(s)A(s)]$

System Diagram

f(t) A(t) c(t)

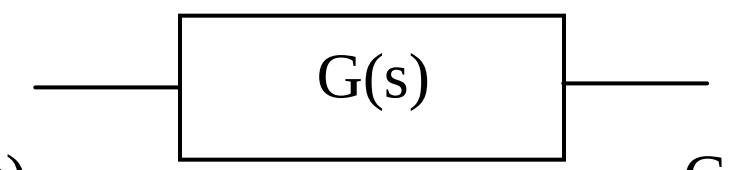

F(s) $C(s) = F(s)G(s) = s[F(s)\frac{G(s)}{s}] = s[F(s)A(s)]$

$A(s) = \frac{1}{s}G(s)$

Input Transfer Function Output

Find c(t)

$$C(s) = s[F(s)A(s)] \tag{5.11-37}$$

Dividing both sides of Eq 5.11-37 by s

$$\frac{C(s)}{s} = F(s)A(s) \tag{5.11-38}$$

Taking the $K_{\Delta t}$ Transform of Eq 5.11-38

$$_{\Delta t}\int_0^t c(t)\Delta t = K_{\Delta t}^{-1}[F(s)A(s)] \tag{5.11-39}$$

Taking the discrete derivative of each side of Eq 5.11-39

$$D_{\Delta t}[\,_{\Delta t}\int_0^t c(t)\Delta t] = D_{\Delta t}[K_{\Delta t}^{-1}[F(s)A(s)]] \tag{5.11-40}$$

$$c(t) = D_{\Delta t}[K_{\Delta t}^{-1}[K_{\Delta t}[f(t)]K_{\Delta t}[A(t)]]] \tag{5.11-41}$$

where

$K_{\Delta t}[f(t)] = F(s)$

$K_{\Delta t}[A(t)] = A(s)$

Rewriting Eq 5.10-32 and Eq 5.10-34 using the functions f(t) and A(t)

$$_{\Delta\lambda}\int_0^t f(t-\lambda-\Delta\lambda)A(\lambda)\Delta\lambda = [K_{\Delta t}^{-1}[K_{\Delta t}[f(t)]K_{\Delta t}[A(t)]]] \tag{5.11-42}$$

$$_{\Delta\lambda}\int_0^t f(\lambda)A(t-\lambda-\Delta\lambda)\Delta\lambda = [K_{\Delta t}^{-1}[K_{\Delta t}[f(t)]K_{\Delta t}[A(t)]]] \tag{5.11-43}$$

Substituting the $K_{\Delta t}$ Transform Convolution Equations of Eq 5.10-42 and Eq 5.10-43 into Eq 5.11-41

Then

$$\mathbf{c(t) = D_{\Delta t}[\,{}_{\Delta\lambda}\!\int_0^t f(t-\lambda-\Delta\lambda)A(\lambda)\Delta\lambda\,]} \tag{5.11-44}$$

**where**

$\mathbf{f(t)}$ **= System input excitation function**

$\mathbf{c(t)}$ **= System output function**

$\mathbf{A(t)}$ **= System response to a unit step function input**

$$\mathbf{D_{\Delta t}\,h(t) = \frac{h(t+\Delta t) - h(t)}{\Delta t}} \textbf{ = discrete derivative of h(t)}$$

$$\mathbf{A(s) = \frac{G(s)}{s}}$$

$$\mathbf{\Delta\lambda = \Delta t}$$

$$\mathbf{\lambda = n\Delta t, \quad n = 0, 1, 2, 3, \ldots, \frac{t}{\Delta t}-1, \frac{t}{\Delta t}}$$

$$\mathbf{\lambda = 0, \Delta\lambda, 2\Delta\lambda, 3\Delta\lambda, \ldots, t-\Delta\lambda, t}$$

$\mathbf{c(t), f(t)}$ **= Sample and hold shaped waveforms, Interval Calculus discrete variable functions**

$$\mathbf{c(t) = D_{\Delta t}[\,{}_{\Delta\lambda}\!\int_0^t f(\lambda)A(t-\lambda-\Delta\lambda)\Delta\lambda\,]} \tag{5.11-45}$$

**where**

$\mathbf{f(t)}$ **= System input excitation function ( See Eq 5.11-20 where f(t) = r(t) )**

$\mathbf{c(t)}$ **= System output function**

$\mathbf{A(t)}$ **= System response to a unit step function input**

$$\mathbf{D_{\Delta t}\,h(t) = \frac{h(t+\Delta t) - h(t)}{\Delta t}} \textbf{ = discrete derivative of h(t)}$$

$$\mathbf{A(s) = \frac{G(s)}{s}}$$

$$\mathbf{\Delta\lambda = \Delta t}$$

$$\mathbf{\lambda = n\Delta t, \quad n = 0, 1, 2, 3, \ldots, \frac{t}{\Delta t}-1, \frac{t}{\Delta t}}$$

$$\mathbf{\lambda = 0, \Delta\lambda, 2\Delta\lambda, 3\Delta\lambda, \ldots, t-\Delta\lambda, t}$$

$\mathbf{A(t-n\Delta t-\Delta t)}$ **= System response to a unit step function initiated at t = [n+1]Δt,**

**( See Eq 5.11-22),** $\mathbf{n = 0, 1, 2, 3, \ldots, \frac{t}{\Delta t}-1}$

$\mathbf{c(t), f(t)}$ **= Sample and hold shaped waveforms, Interval Calculus discrete variable functions**

Thus all of the eight $K_{\Delta t}$ Transform Duhamel Equations previously listed are now derived.

Derivation of the Interval Calculus Duhamel Equations using Z Transforms

The eight previously listed Interval Calculus Duhamel Equations can be derived using Z Transforms instead of $K_{\Delta t}$ Transforms. For those who wish to see how this would be done, the previous derivations are repeated below using Z Transforms.

For the Z Transform equation, $C(z) = F(z)[\frac{z-1}{z} A(z)]$

System Diagram

f(t) A(t) c(t)

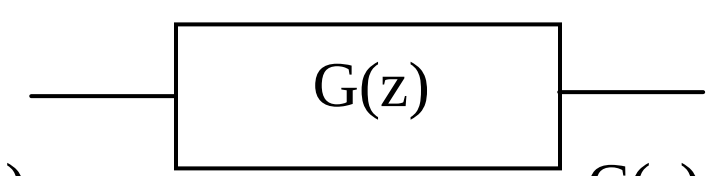

F(z) $C(z) = F(z)G(z) = F(z)\frac{z-1}{z}[\frac{z}{z-1} G(z)] = F(z)[\frac{z-1}{z} A(z)]$

$A(z) = \frac{z}{z-1} G(z)$

Input Transfer Function Output

Find c(t)

$$Z[A(t+T)] = zA(z) - zA(0) \quad (5.11\text{-}46)$$

$$A^{'}(z) = Z[A^{'}(t)] = Z[\Delta A(t)] = Z[A(t+T) - A(t)] = zA(z) - zA(0) - A(z) = (z-1)A(z) - zA(0) \quad (5.11\text{-}47)$$

$$C(z) = F(z)G(z) = F(z)\frac{z-1}{z}[\frac{z}{z-1} G(z)] = F(z)[\frac{z-1}{z} A(z)] \quad (5.11\text{-}48)$$

$$C(z) = \frac{F(z)}{z}[\{(z-1)A(z) - zA(0)\} + zA(0)] \quad (5.11\text{-}49)$$

Substituting Eq 5.11-47 into Eq 5.11-49

$$C(z) = \frac{F(z)}{z}[A^{'}(z) + zA(0)] \quad (5.11\text{-}50)$$

Changing the form of Eq 5.11-50

$$C(z) = A(0)F(z) + \frac{F(z)A^{'}(z)}{z} \quad (5.11\text{-}51)$$

Taking the Inverse Z Transform of Eq 5.11-51

$$c(t) = Z^{-1}[C(z)] = A(0)f(t) + Z^{-1}[\frac{F(z)A^{'}(z)}{z}] \quad (5.11\text{-}52)$$

Rewriting the Z Transform Convolution Equations, Eq 5.10-85 and Eq 5.10-86

$$Z[\frac{1}{T}\,{}_{T}\!\int_{0}^{t} f(t-\lambda-\Delta\lambda)g(\lambda)\Delta\lambda] = \frac{Z[f(t)]Z[g(t)]}{z} = \frac{F(z)G(z)}{z} \quad (5.11\text{-}53)$$

$$Z[\frac{1}{T}\,{}_{T}\!\int_{0}^{t} f(\lambda)g(t-\lambda-\Delta\lambda)\Delta\lambda] = \frac{Z[f(t)]Z[g(t)]}{z} = \frac{F(z)G(z)}{z} \quad (5.11\text{-}54)$$

From Eq 5.11-52 and Eq 5.11-53

$$c(t) = A(0)f(t) + \frac{1}{T}\,{}_{T}\!\int_{0}^{t} f(t-\lambda-\Delta\lambda)A^{'}(\lambda)\Delta\lambda] \quad (5.11\text{-}55)$$

where
$A^{'}(\lambda) = \Delta A(\lambda) = A(\lambda+\Delta\lambda) - A(\lambda)$
$\Delta\lambda = \Delta t = T$

$$\Delta\lambda = \Delta t = T \quad (5.11\text{-}56)$$

$$A^{'}(\lambda) = \Delta A(\lambda) = A(\lambda+\Delta\lambda) - A(\lambda) \quad (5.11\text{-}57)$$

$$A'(\lambda) = D_{\Delta\lambda}A(\lambda) = \frac{\Delta A(\lambda)}{\Delta\lambda} = \frac{A^{'}(\lambda)}{T} \quad (5.11\text{-}58)$$

From Eq 5.11-55, Eq 5.11-56 and Eq 5.11-58

$$c(t) = A(0)f(t) + {}_{\Delta\lambda}\!\int_{0}^{t} f(t-\lambda-\Delta\lambda)A'(\lambda)\Delta\lambda \quad (5.11\text{-}59)$$

Then

$$\mathbf{c(t) = A(0)f(t) + {}_{\Delta\lambda}\!\int_{0}^{t} f(t-\lambda-\Delta\lambda)A'(\lambda)\Delta\lambda} \quad (5.11\text{-}60)$$

**where**
**f(t) = System input excitation function**
**c(t) = System output function**

**A(t) = System response to a unit step function input**
**A(0) = 0, The system, by definition, is initally passive**
$$A^{'}(t) = D_{\Delta t}A(t) = \frac{\Delta A(t)}{\Delta t} = \frac{A(t+\Delta t) - A(t)}{\Delta t} = \text{discrete derivative of } A(t)$$
$$A(z) = \frac{z}{z-1}G(z)$$
$$\Delta\lambda = \Delta t = T$$
$$\lambda = n\Delta t, \quad n = 0, 1, 2, 3, \ldots, \frac{t}{\Delta t}-1, \frac{t}{\Delta t}$$
$$\lambda = 0, \Delta\lambda, 2\Delta\lambda, 3\Delta\lambda, \ldots, t-\Delta\lambda, t$$
**c(t), f(t) = Sample and hold shaped waveforms, Interval Calculus discrete variable functions**

Since A(t) is by definition the response of a system that is initially passive, A(0) = 0.

$$c(t) = {}_{\Delta\lambda}\int_{0}^{t} f(t-\lambda-\Delta\lambda)A^{'}(\lambda)\Delta\lambda \qquad (5.11\text{-}61)$$

**where**
**f(t) = System input excitation function**
**c(t) = System output function**
**A(t) = System response to a unit step function input**
**A(0) = 0, The system, by definition, is initally passive**
$$A^{'}(t) = D_{\Delta t}A(t) = \frac{\Delta A(t)}{\Delta t} = A^{'}(t) = \frac{A(t+\Delta t) - A(t)}{\Delta t} = \text{discrete derivative of } A(t)$$
$$A(z) = \frac{z}{z-1}G(z)$$
$$\Delta\lambda = \Delta t = T$$
$$\lambda = n\Delta t, \quad n = 0, 1, 2, 3, \ldots, \frac{t}{\Delta t}-1, \frac{t}{\Delta t}$$
$$\lambda = 0, \Delta\lambda, 2\Delta\lambda, 3\Delta\lambda, \ldots, t-\Delta\lambda, t$$
**c(t), f(t) = Sample and hold shaped waveforms, Interval Calculus discrete variable functions**

From Eq 5.11-52 and Eq 5.11-54

$$c(t) = A(0)F(t) + \frac{1}{T}\,{}_{T}\int_{0}^{t} f(\lambda)A^{'}(t-\lambda-\Delta\lambda)\Delta\lambda] \qquad (5.11\text{-}62)$$

where
$A^{'}(\lambda) = \Delta A(\lambda) = A(\lambda+\Delta\lambda) - A(\lambda)$
$\Delta\lambda = \Delta t = T$

$$\Delta\lambda = \Delta t = T \qquad (5.11\text{-}63)$$

$$A^{'}(\lambda) = \Delta A(\lambda) = A(\lambda+\Delta\lambda) - A(\lambda) \qquad (5.11\text{-}64)$$

$$A^{'}(\lambda) = D_{\Delta\lambda}A(\lambda) = \frac{\Delta A(\lambda)}{\Delta\lambda} = \frac{A^{'}(\lambda)}{T} \qquad (5.11\text{-}65)$$

From Eq 5.11-62, Eq 5.11-63 and Eq 5.11-65

$$c(t) = A(0)F(t) + {}_{\Delta\lambda}\int_{0}^{t} f(t-\lambda-\Delta\lambda)A'(\lambda)\Delta\lambda \qquad (5.11\text{-}66)$$

Then

$$\mathbf{c(t) = A(0)f(t) + {}_{\Delta\lambda}\int_{0}^{t} f(\lambda)A^{'}(t-\lambda-\Delta\lambda)\Delta\lambda} \qquad (5.11\text{-}67)$$

**where**

**f(t) = System input excitation function ( See Eq 5.11-18 where f(t) = r(t) )**
**c(t) = System output function**
**A(t) = System response to a unit step function input**
**A(0) = 0, The system, by definition, is initally passive**

$$\mathbf{A^{'}(t) = D_{\Delta t}A(t) = \frac{\Delta A(t)}{\Delta t} = \frac{A(t+\Delta t) - A(t)}{\Delta t}} \text{ = discrete derivative of A(t)}$$

$$\mathbf{A(z) = \frac{z}{z-1}G(z)}$$

$$\mathbf{\Delta\lambda = \Delta t = T}$$

$$\mathbf{\lambda = n\Delta t, \quad n = 0, 1, 2, 3, \ldots, \frac{t}{\Delta t}-1, \frac{t}{\Delta t}}$$

$$\mathbf{\lambda = 0, \Delta\lambda, 2\Delta\lambda, 3\Delta\lambda, \ldots, t-\Delta\lambda, t}$$

**$A^{'}(t-\lambda-\Delta\lambda)\Delta\lambda$ = System response to a unit amplitude $\Delta\lambda$ width pulse input initiated at t = $\lambda$ ( See Eq 5.11-21 where $A^{'}(t)$ is substituted for A(t))**

**c(t), f(t) = Sample and hold shaped waveforms, Interval Calculus discrete variable functions**

Since A(t) is by definition the response of a system that is initially passive, A(0) = 0.

$$\mathbf{c(t) = {}_{\Delta\lambda}\int_{0}^{t} f(\lambda)A^{'}(t-\lambda-\Delta\lambda)\Delta\lambda} \qquad (5.11\text{-}68)$$

**where**

**f(t) = System input excitation function ( See Eq 5.11-18 where f(t) = r(t) )**
**c(t) = System output function**
**A(t) = System response to a unit step function input**
**A(0) = 0, The system, by definition, is initally passive**

$$\mathbf{A^{'}(t) = D_{\Delta t}A(t) = \frac{\Delta A(t)}{\Delta t} = \frac{A(t+\Delta t) - A(t)}{\Delta t}} \text{ = discrete derivative of A(t)}$$

$$\mathbf{A(z) = \frac{z}{z-1}G(z)}$$

**$\Delta\lambda = \Delta t = T$**

**$\lambda = n\Delta t, \quad n = 0, 1, 2, 3, \ldots, \frac{t}{\Delta t} - 1, \frac{t}{\Delta t}$**

**$\lambda = 0, \Delta\lambda, 2\Delta\lambda, 3\Delta\lambda, \ldots, t-\Delta\lambda, t$**

**$A'(t-\lambda-\Delta\lambda)\Delta\lambda$ = System response to a unit amplitude $\Delta\lambda$ width pulse input initiated at $t = \lambda$**
**( See Eq 5.11-21 where $A'(t)$ is substituted for $A(t)$)**

**$c(t), f(t)$ = Sample and hold shaped waveforms, Interval Calculus discrete variable functions**

For the Z Transform equation, $C(z) = [\frac{z-1}{z} F(z)]A(z)$

System Diagram

f(t) A(t) c(t)

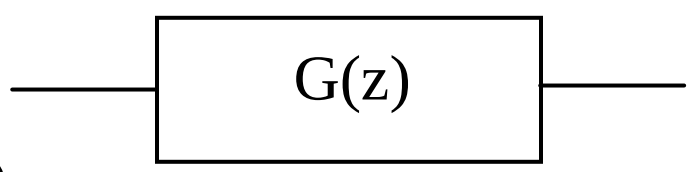

F(z) $C(z) = [\frac{z-1}{z} F(z)][\frac{z}{z-1} G(z)] = [\frac{z-1}{z} F(z)]A(z)$

$A(z) = \frac{z}{z-1} G(z)$

Input Transfer Function Output

$$C(z) = F(z)G(z) = [\frac{z-1}{z} F(z)]A(z) = [\frac{z-1}{z} F(z) - \frac{z}{z} f(0)]A(z) + f(0)A(z) \quad (5.11\text{-}69)$$

$$C(z) = [(z-1)F(z) - z\, f(0)]\frac{A(z)}{z} + f(0)A(z) \quad (5.11\text{-}70)$$

From Eq 5.11-70 and Eq 5.11-47

$$C(z) = f(0)A(z) + \frac{F'(z)A(z)}{z} \quad (5.11\text{-}71)$$

where

$F'(z) = Z[F'(t)] = Z[\Delta f(t)] = Z[f(t+T) - f(t)]$

Taking the Inverse Z Transform of Eq 5.11-71

$$c(t) = Z^{-1}[C(z)] = f(0)A(t) + Z^{-1}[\frac{F'(z)A(z)}{z}] \quad (5.11\text{-}72)$$

Rewriting the Z Transform Convolution Equations, Eq 5.11-53 and Eq 5.11-54

$$Z[\frac{1}{T}\,{}_{T}\int_{0}^{t} f(t-\lambda-\Delta\lambda)g(\lambda)\Delta\lambda] = \frac{Z[f(t)]Z[g(t)]}{z} = \frac{F(z)G(z)}{z} \qquad (5.11\text{-}73)$$

$$Z[\frac{1}{T}\,{}_{T}\int_{0}^{t} f(\lambda)g(t-\lambda-\Delta\lambda)\Delta\lambda] = \frac{Z[f(t)]Z[g(t)]}{z} = \frac{F(z)G(z)}{z} \qquad (5.11\text{-}74)$$

From Eq 5.11-72 and Eq 5.11-73

$$c(t) = f(0)A(t) + \frac{1}{T}\,{}_{T}\int_{0}^{t} f^{'}(t-\lambda-\Delta\lambda)A(\lambda)\Delta\lambda] \qquad (5.11\text{-}75)$$

where
- $f^{'}(t) = \Delta f(t) = f(t+T) - f(t)$
- $\Delta\lambda = \Delta t = T$
- f(t) = System input excitation function
- c(t) = System output function
- A(t) = System response to a unit step function input
- A(0) = 0, The system, by definition, is initally passive

$$f^{'}(\lambda) = D_{\Delta\lambda}f(\lambda) = \frac{f(\lambda+\Delta\lambda) - f(\lambda)}{\Delta\lambda} \qquad (5.11\text{-}76)$$

$$f^{'}(\lambda) = \Delta f(\lambda) = f(\lambda+\Delta\lambda) - f(\lambda) \qquad (5.11\text{-}77)$$

$$f^{'}(\lambda) = D_{\Delta\lambda}f(\lambda) = \frac{\Delta f(\lambda)}{\Delta\lambda} = \frac{f^{'}(\lambda)}{T} \qquad (5.11\text{-}78)$$

From Eq 5.11-75, thru Eq 5.11-78

$$c(t) = f(0)A(t) + {}_{\Delta\lambda}\int_{0}^{t} f'(t-\lambda-\Delta\lambda)A(\lambda)\Delta\lambda \qquad (5.11\text{-}79)$$

Then

$$\mathbf{c(t) = f(0)A(t) + {}_{\Delta\lambda}\int_{0}^{t} f'(t-\lambda-\Delta\lambda)A(\lambda)\Delta\lambda} \qquad (5.11\text{-}80)$$

**where**

$$\mathbf{f'(t) = \frac{\Delta f(t)}{\Delta t} = \frac{f(t+\Delta t) - f(t)}{\Delta t}}$$ **= discrete derivative of f(t)**

$$\mathbf{A(z) = \frac{z}{z-1}G(z)}$$

$$\mathbf{\Delta\lambda = \Delta t = T}$$

$$\mathbf{\lambda = n\Delta t, \quad n = 0, 1, 2, 3, \ldots, \frac{t}{\Delta t}-1, \frac{t}{\Delta t}}$$

$$\mathbf{\lambda = 0, \Delta\lambda, 2\Delta\lambda, 3\Delta\lambda, \ldots, t-\Delta\lambda, t}$$

**c(t), f(t) = Sample and hold shaped waveforms, Interval Calculus discrete variable functions**

From Eq 5.11-72 and Eq 5.11-74

$$c(t) = f(0)A(t) + \frac{1}{T}\,{}_T\!\int_0^t f^{'}(\lambda)A(t-\lambda-\Delta\lambda)\Delta\lambda] \qquad (5.11\text{-}81)$$

where

$f^{'}(t) = \Delta f(t) = f(t+T) - f(t)$

$\Delta\lambda = \Delta t = T$

f(t) = System input excitation function

c(t) = System output function

A(t) = System response to a unit step function input

A(0) = 0, The system, by definition, is initally passive

$$f^{'}(\lambda) = D_{\Delta\lambda}f(\lambda) = \frac{f(\lambda+\Delta\lambda) - f(\lambda)}{\Delta\lambda} \qquad (5.11\text{-}82)$$

$$f^{'}(\lambda) = \Delta f(\lambda) = f(\lambda+\Delta\lambda) - f(\lambda) \qquad (5.11\text{-}83)$$

$$f^{'}(\lambda) = D_{\Delta\lambda}f(\lambda) = \frac{\Delta f(\lambda)}{\Delta\lambda} = \frac{f^{'}(\lambda)}{T} \qquad (5.11\text{-}84)$$

From Eq 5.11-81, thru Eq 5.11-84

$$c(t) = f(0)A(t) + {}_{\Delta\lambda}\!\int_0^t f'(\lambda)A(t-\lambda-\Delta\lambda)\Delta\lambda \qquad (5.11\text{-}85)$$

Then

$$\mathbf{c(t) = f(0)A(t) + {}_{\Delta\lambda}\!\int_0^t f'(\lambda)A(t-\lambda-\Delta\lambda)\Delta\lambda} \qquad (5.11\text{-}86)$$

**where**

**f(t) = System input excitation function See Eq 5.11-20 where f(t) = r(t))**

**c(t) = System output function**

**A(t) = System response to a unit step function input**

**$f'(t) = D_{\Delta t}f(t) = \frac{\Delta f(t)}{\Delta t} = \frac{f(t+\Delta t) - f(t)}{\Delta t}$ = discrete derivative of f(t)**

**$A(z) = \frac{z}{z-1}G(z)$**

**$\Delta\lambda = \Delta t = T$**

**$\lambda = n\Delta t, \quad n = 0, 1, 2, 3, \ldots, \frac{t}{\Delta t}-1, \frac{t}{\Delta t}$**

**$\lambda = 0, \Delta\lambda, 2\Delta\lambda, 3\Delta\lambda, \ldots, t-\Delta\lambda, t$**

**$A(t-n\Delta t-\Delta t)$ = System response to a unit step function initiated at $t = [n+1]\Delta t$,**

**( See Eq 5.11-22), $\quad n = 0, 1, 2, 3, \ldots, \frac{t}{\Delta t}-1$**

**$f'(n\Delta t)\Delta t = \Delta f(n\Delta t) = [f(n\Delta t+\Delta t)-f(n\Delta t)]$**

**c(t), f(t) = Sample and hold shaped waveforms, Interval Calculus discrete variable functions**

<u>For the Z Transform equation, $C(z) = \frac{z-1}{z}[F(z)A(z)]$</u>

System Diagram

f(t) A(t) c(t)

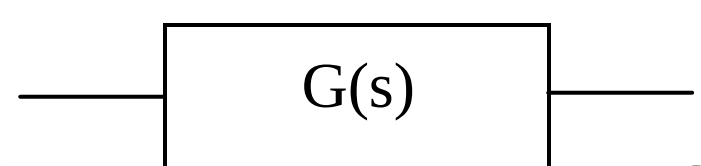

F(z) $C(z) = \frac{z-1}{z}[F(z)\frac{z}{z-1}G(z)] = \frac{z-1}{z}[F(z)A(z)]$

$A(z) = \frac{z}{z-1}G(z)$

Input Transfer Function Output

Find c(t)

$$C(z) = \frac{z-1}{z}[F(z)\frac{z}{z-1}G(z)] = \frac{z-1}{z}[F(z)A(z)] \tag{5.11-87}$$

$$C(z) = (z-1)\left[\frac{F(z)A(z)}{z}\right] \tag{5.11-88}$$

After changing variables substitute Eq 5.11-73 into Eq 5.11-88

$$C(z) = (z-1)Z\left[\frac{1}{T}\,{}_T\!\int_0^t f(t-\lambda-\Delta\lambda)A(\lambda)\Delta\lambda\right] \tag{5.11-89}$$

Rewriting Eq 5.11-89 in a different form

$$C(z) = (z-1)P(z) \tag{5.11-90}$$

$$P(z) = Z[\frac{1}{T}\,{}_{T}\int_{0}^{t} f(t-\lambda-\Delta\lambda)A(\lambda)\Delta\lambda] \tag{5.11-91}$$

Taking the Inverse Z Transform of Eq 5.11-91

$$P(t) = \frac{1}{T}\,{}_{T}\int_{0}^{t} f(t-\lambda-\Delta\lambda)A(\lambda)\Delta\lambda \tag{5.11-92}$$

For t = 0

$$P(0) = 0 \tag{5.11-93}$$

From Eq 5.11-88

$$C(z) = [(z-1)P(z) - zP(0)] + zP(0) = (z-1)P(z) - zP(0) \tag{5.11-94}$$

From Eq 5.11-47 and Eq 5.11-94

$$Z[\Delta P(t)] = (z-1)P(z) - zP(0) \tag{5.11-95}$$

Substituting Eq 5.11-95 into Eq 5.11-94

$$C(z) = Z[\Delta P(t)] \tag{5.11-96}$$

Taking the Inverse Z Transform of Eq 5.11-96

$$c(t) = \Delta P(t) \tag{5.11-97}$$

Substituting Eq 5.11-92 into Eq 5.11-97

$$c(t) = \frac{\Delta}{T}\,{}_{T}\int_{0}^{t} f(t-\lambda-\Delta\lambda)A(\lambda)\Delta\lambda \tag{5.11-98}$$

$$D_{\Delta t} = \frac{\Delta}{T} \tag{5.11-99}$$

Substituting Eq 5.11-99 into Eq 5.11-98

$$c(t) = D_{\Delta t}\,{}_{T}\int_{0}^{t} f(t-\lambda-\Delta\lambda)A(\lambda)\Delta\lambda \tag{5.11-100}$$

From Eq 5.11-73 and Eq 5.11-74

$${}_{T}\int_{0}^{t} f(\lambda)A(t-\lambda-\Delta\lambda)\Delta\lambda = {}_{T}\int_{0}^{t} f(t-\lambda-\Delta\lambda)A(\lambda)\Delta\lambda \tag{5.11-101}$$

Substituting Eq 5.11-101 into Eq 5.11-100

$$c(t) = D_{\Delta t}\ {}_{T}\!\int_{0}^{t} f(\lambda)A(t-\lambda-\Delta\lambda)\Delta\lambda \qquad (5.11\text{-}102)$$

Then

From Eq 5.11-100 and Eq 5.11-102

$$\mathbf{c(t) = D_{\Delta t}\ {}_{T}\!\int_{0}^{t} f(t-\lambda-\Delta\lambda)A(\lambda)\Delta\lambda} \qquad (5.11\text{-}103)$$

**where**

**f(t) = System input excitation function**

**c(t) = System output function**

**A(t) = System response to a unit step function input**

$D_{\Delta t}\,h(t) = \dfrac{h(t+\Delta t)-h(t)}{\Delta t}$ **= discrete derivative of h(t)**

$A(s) = \dfrac{G(s)}{s}$

$\Delta\lambda = \Delta t$

$\lambda = n\Delta t, \quad n = 0, 1, 2, 3, \ldots, \frac{t}{\Delta t}-1, \frac{t}{\Delta t}$

$\lambda = 0, \Delta\lambda, 2\Delta\lambda, 3\Delta\lambda, \ldots, t-\Delta\lambda, t$

**c(t), f(t) = Sample and hold shaped waveforms, Interval Calculus discrete variable functions**

**and**

$$\mathbf{c(t) = D_{\Delta t}[\ {}_{\Delta\lambda}\!\int_{0}^{t} f(\lambda)A(t-\lambda-\Delta\lambda)\Delta\lambda]} \qquad (5.11\text{-}101)$$

**where**

**f(t) = System input excitation function ( See Eq 5.11-20 where f(t) = r(t) )**

**c(t) = System output function**

**A(t) = System response to a unit step function input**

$D_{\Delta t}\,h(t) = \dfrac{h(t+\Delta t)-h(t)}{\Delta t}$ **= discrete derivative of h(t)**

$A(s) = \dfrac{G(s)}{s}$

$\Delta\lambda = \Delta t$

$\lambda = n\Delta t, \quad n = 0, 1, 2, 3, \ldots, \frac{t}{\Delta t}-1, \frac{t}{\Delta t}$

$\lambda = 0, \Delta\lambda, 2\Delta\lambda, 3\Delta\lambda, \ldots, t-\Delta\lambda, t$

$A(t-\lambda-\Delta\lambda)\Delta\lambda$ **= System response to a unit amplitude** $\Delta\lambda$ **width pulse input initiated at** $t = \lambda$ **( See Eq 5.11-21)**

**c(t), f(t) = Sample and hold shaped waveforms, Interval Calculus discrete variable functions**

Thus all of the eight Z Transform Duhamel Equations previously listed are now derived. Note that they are the same as the $K_{\Delta t}$ Transform Duhamel Equations.

Example 5.11-1 Use several of the previously derived Interval Calculus Duhamel equations to find the output of a system, c(t), where the input f(t) = U(t), and A(t) = t(t-Δt).

Using $c(t) = {}_{\Delta\lambda}\int_0^t f(\lambda)A'(t-\lambda-\Delta\lambda)\,\Delta\lambda$ 1)

$A(t) = t(t-\Delta t)$ 2)

$A' = 2t$ 3)

$f(t) = 1$ 4)

$c(t) = 2\,{}_{\Delta\lambda}\int_0^t (1)A'(t-\lambda-\Delta\lambda)\,\Delta\lambda$ 5)

$c(t) = 2[t\,{}_{\Delta\lambda}\int_0^t \Delta\lambda - {}_{\Delta\lambda}\int_0^t \lambda\Delta\lambda - \Delta\lambda\,{}_{\Delta\lambda}\int_0^t \Delta\lambda\,] = 2[t\lambda\Big|_0^t - \frac{\lambda(\lambda-\Delta\lambda)}{2}\Big|_0^t - \Delta\lambda\lambda\Big|_0^t\,]$ 6)

$c(t) = 2[t^2 - \frac{t(t-\Delta t)}{2} - t\Delta t] = 2[t^2 - \frac{t^2}{2} + \frac{t\Delta t}{2} - t\Delta t]\ ,\quad \Delta t = \Delta\lambda$ 7)

$c(t) = 2[\,\frac{t^2}{2} - \frac{t\Delta t}{2}\,]$ 8)

$\mathbf{c(t) = t(t-\Delta t)}$ 9)

Using $c(t) = D_{\Delta t}[\,{}_{\Delta\lambda}\int_0^t f(t-\lambda-\Delta\lambda)A(\lambda)\,\Delta\lambda\,]$ 10)

$A(t) = t(t-\Delta t)$ 11)

$f(t) = 1$ 12)

$c(t) = D_{\Delta t}[\,{}_{\Delta\lambda}\int_0^t (1)\lambda(\lambda-\Delta\lambda)\,\Delta\lambda\,] = D_{\Delta t}\,[\frac{\lambda(\lambda-\Delta\lambda)(\lambda-2\Delta\lambda)}{3}]\Big|_0^t$ 13)

$c(t) = D_{\Delta t}[\,{}_{\Delta\lambda}\int_0^t (1)\lambda(\lambda-\Delta\lambda)\,\Delta\lambda\,] = D_{\Delta t}\,[\frac{t(t-\Delta t)(t-2\Delta t)}{3}]\ ,\quad \Delta t = \Delta\lambda$ 14)

$$\mathbf{c(t) = t(t\text{-}\Delta t)} \qquad 15)$$

Good check

$$\text{Using } c(t) = D_{\Delta t}[\,{}_{\Delta\lambda}\int_0^t f(\lambda)A(t-\lambda-\Delta\lambda)\,\Delta\lambda\,] \qquad 16)$$

$$A(t) = t(t-\Delta t) \qquad 17)$$

$$f(t) = 1 \qquad 18)$$

$$c(t) = D_{\Delta t}[\,{}_{\Delta\lambda}\int_0^t (1)(t-\lambda-\Delta\lambda)(t-\lambda-2\Delta\lambda)\,\Delta\lambda\,] \qquad 19)$$

$$c(t) = D_{\Delta t}[\,{}_{\Delta\lambda}\int_0^t (\,t^2 - t\lambda - 2t\Delta\lambda - \lambda t + \lambda^2 + 2\lambda\Delta\lambda - \Delta\lambda t + \Delta\lambda\lambda + 2\Delta\lambda^2)\Delta\lambda\,] \qquad 20)$$

$$c(t) = D_{\Delta t}[\,{}_{\Delta\lambda}\int_0^t (\,t^2 - 2t\lambda - 3t\Delta\lambda + 3\lambda\Delta\lambda + \lambda^2 + 2\Delta\lambda^2)\Delta\lambda\,] \qquad 21)$$

$$c(t) = D_{\Delta t}[t^2\lambda\,\Big|_0^t - \frac{2t\lambda(\lambda-\Delta\lambda)}{2}\Big|_0^t - 3t\Delta\lambda\lambda\,\Big|_0^t + \frac{3\Delta\lambda\lambda(\lambda-\Delta\lambda)}{2}\Big|_0^t + 2\Delta\lambda^2\lambda\,\Big|_0^t + {}_{\Delta\lambda}\int_0^t \{\lambda(\lambda-\Delta\lambda) + \lambda\Delta\lambda\}\Delta\lambda] \qquad 22)$$

$$c(t) = D_{\Delta t}[t^2\lambda\,\Big|_0^t - \frac{2t\lambda(\lambda-\Delta\lambda)}{2}\Big|_0^t - 3t\Delta\lambda\lambda\,\Big|_0^t + \frac{3\Delta\lambda\lambda(\lambda-\Delta\lambda)}{2}\Big|_0^t + \frac{\lambda(\lambda-\Delta\lambda)(\lambda-2\Delta\lambda)}{3}\Big|_0^t + \frac{\Delta\lambda\lambda(\lambda-\Delta\lambda)}{2}\Big|_0^t + 2\Delta\lambda^2\lambda\,\Big|_0^t\,] \qquad 23)$$

$$c(t) = D_{\Delta t}[\,t^3 - t^2(t-\Delta t) - 3t^2\Delta t + \frac{3tt(t-\Delta t)}{2} + \frac{t(t-\Delta t)(t-2\Delta\lambda)}{3} + \frac{\Delta tt(t-\Delta t)}{2} + 2\Delta t^2 t\,]\,,\ \ \Delta t = \Delta\lambda \qquad 24)$$

$$c(t) = D_{\Delta t}[-2t^2\Delta t] + 3\Delta tt + t(t-\Delta t) + \Delta tt + 2\Delta t^2 \qquad 25)$$

$$c(t) = D_{\Delta t}[-2\Delta t\{t(t-\Delta t) + t\Delta t\}] + 3\Delta tt + t(t-\Delta t) + \Delta tt + 2\Delta t^2 \qquad 26)$$

$$c(t) = -4\Delta tt - 2\Delta t^2 + 3\Delta tt + t^2 - t\Delta t + \Delta tt + 2\Delta t^2 \qquad 27)$$

$$c(t) = t^2 - t\Delta t \qquad 28)$$

$$\mathbf{c(t) = t(t\text{-}\Delta t)} \qquad 29)$$

Good check

Using $c(t) = f(0)A(t) + \int_{0}^{t} {}_{\Delta\lambda} f'(\lambda)A(t-\lambda-\Delta\lambda)\,\Delta\lambda$ 30)

$A(t) = t(t-\Delta t)$ 31)

$f(t) = 1$ 32)

$\Delta\lambda = \Delta t$ 33)

$c(t) = (1)[t(t-\Delta t)] + \int_{0}^{t} {}_{\Delta\lambda} (0)(t-\lambda-\Delta\lambda)(t-\lambda-2\Delta t)\Delta\lambda$ 34)

**$c(t) = t(t-\Delta t)$** 35)

Good check

Using $c(t) = \int_{0}^{t} {}_{\Delta\lambda} f(t-\lambda-\Delta\lambda)A'(\lambda)\,\Delta\lambda$ 36)

$A(t) = t(t-\Delta t)$ 37)

$A'(t) = 2t$ 38)

$f(t) = 1$ 39)

$\Delta\lambda = \Delta t$ 40)

$c(t) = \int_{0}^{t} {}_{\Delta\lambda} (1)(2\lambda)\Delta\lambda = \frac{2\lambda(\lambda-\Delta\lambda)}{2}\Big|_{0}^{t} = t(t-\Delta t)$ 41)

**$c(t) = t(t-\Delta t)$** 42)

Good check

## Section 5.12: Mathematical analysis of sampled-data systems using Interval Calculus

Presently, a very common approach to analyzing sampled-data control systems is to use the Z Transform to describe actions measured in discrete time. Over many years, Z Transform analysis has been found to be very effective in solving problems involving discrete time. With the advent of the use of digital systems, and thus variable sampling, Z Transform control system analysis has become more and more in demand. Control system designers need to be able to design feedback systems which are stable within given specifications. Since the use of Z Transforms and Laplace Transforms (from which Z Transfroms are derived) has been used so effectively, one might wonder why anyone would seek another discrete time analytical methodology. The answer is two fold. Firstly, If such a mathematical methodology exists, it should be made known. Mathematics research is performed just for this reason. Secondly, there is always the possibility that a new mathematical methodology might provide features that increase insight for more effective problem solving. With this in mind, the Interval Calculus previously derived and described is put to use in analyzing discrete time variable (sampled-data) control systems.

Calculus is a special case of Interval Calculus where the discrete time interval, $\Delta t$, is infinitessimal. Thus, one might conclude that Interval Calculus, itself, might have the potential to provide the same analysis commonly performed by the Z Transform. This has been found to be the case. $K_{\Delta t}$ Transform analysis can be used instead of Z Transform analysis. However, as part of the mathematical research performed to establish this fact, something very useful was found. Because of the close relationship between the $K_{\Delta t}$ Transform and the Z Transform, Interval Calculus was found to be useful in simplifying some Z Transform mathematical manipulations. Note Section 5.10 and Section 5.11 where the Z Transform Convolution Equation is derived and the Interval Calculus Duhamel Equations are derived using Z Transforms. In fact, in Section 5.4 and in Section 5.5, a methodology to convert between the Interval Calculus $K_{\Delta t}$ Transform and the Z Transform is derived.

There are a number of Interval Calculus characteristics that can be advantageous to sampled-data analysis. Note the listing below.

1. Interval Calculus deals with functions that have a sample and hold shaped waveform. Often in the analysis of sampled data systems it is functions of this type that are being analyzed.

2. Z Transforms change a sampled continuous time function into a sequence of unit area impulses weighted with the value of the continuous time function at the sampling instants. Mathematically, this is fine but, to some, it may be visually confusing. Interval Calculus $K_{\Delta t}$ Transforms change a sampled continuous time function into a sequence of unit amplitude pulses of $\Delta t$ width weighted by the value of the continuous time function at the sampling instants. This sequence has the appearance of a sample and hold shaped waveform of the original continuous time function.

3. Interval Calculus notation was made to be similar to that of Calculus. For those familiar with Calculus, it is not difficult to understand Interval Calculus operations.

4. Interval Calculus $K_{\Delta t}$ Transforms are very similar in appearance to Laplace Transforms. Easier visual checking of transforms could result in less careless errors being made.

5. Interval Calculus problem solutions are often composed of Interval Calculus discrete functions. These functions can be manipulated directly using discrete integration and differentiation.

6. Calculus functions such as $e^{at}$ and sinat can be used in Interval Calculus. However, to simplify mathematical manipulations on these functions, such as discrete variable differentiation and discrete variable integration, their Interval Calculus identify functions are often used.

7. There exist conversion relationships that allow an Interval Calculus analysis to be changed into an equivalent Z Transform analysis and vice versa.

8. Nyquist and Bode stability analysis can be applied to Interval Calculus $K_{\Delta t}$ Transforms. The left half plane area of stability is closely related to the left half plane of Laplace Transforms. Stability exists within a critical circle in the left half plane with a center at the point $-\frac{1}{\Delta t}$ and with a radius of $\frac{1}{\Delta t}$.

9. Due to the relationship of Calculus (where $\Delta t$ is an infinitesimal value) to Interval Calculus (where $\Delta t$ need not be an infinitesimal value), it is fairly easy in the complex plane to visualize the effect of $\Delta t$ on system stability.

10. The sample and hold function resulting from the sampling of a Calculus function is the Interval Calculus identity of that function.

11. Interval Calculus integrates the area under functions with sample and hold shaped waveforms in a manor similar to the way Calculus integrates the area under continuous time functions. This discrete variable integration capability provides many mathematical manipulation simplifications.

Below, various methods for the use of Interval Calculus in the analysis of sampled-data systems will be derived and presented.

<u>The use of Interval Calculus to describe and analyze sampled-data systems</u>

In previous sections where derivations have been performed, the notation, f(t), was used for both continuous and discrete time functions. Which type of function was being referenced was readily apparent from the context of the derivation. In the following derivations and explanations it is necessary to distinguish between functions of a continuous time variable and functions of a discrete time variable. For continuous time the function designation, f(t), is used and for discrete time the function designation, $f^*(t)$, is used. In the following derivations where f(t) and $f^*(t)$ are both used, $f^*(t)$ represents the waveform resulting from the sample and hold sampling of f(t).

Consider a sample and hold switch as shown in Diagram 5.12-1 below.

Diagram 5.12-1  A Sample and Hold Switch

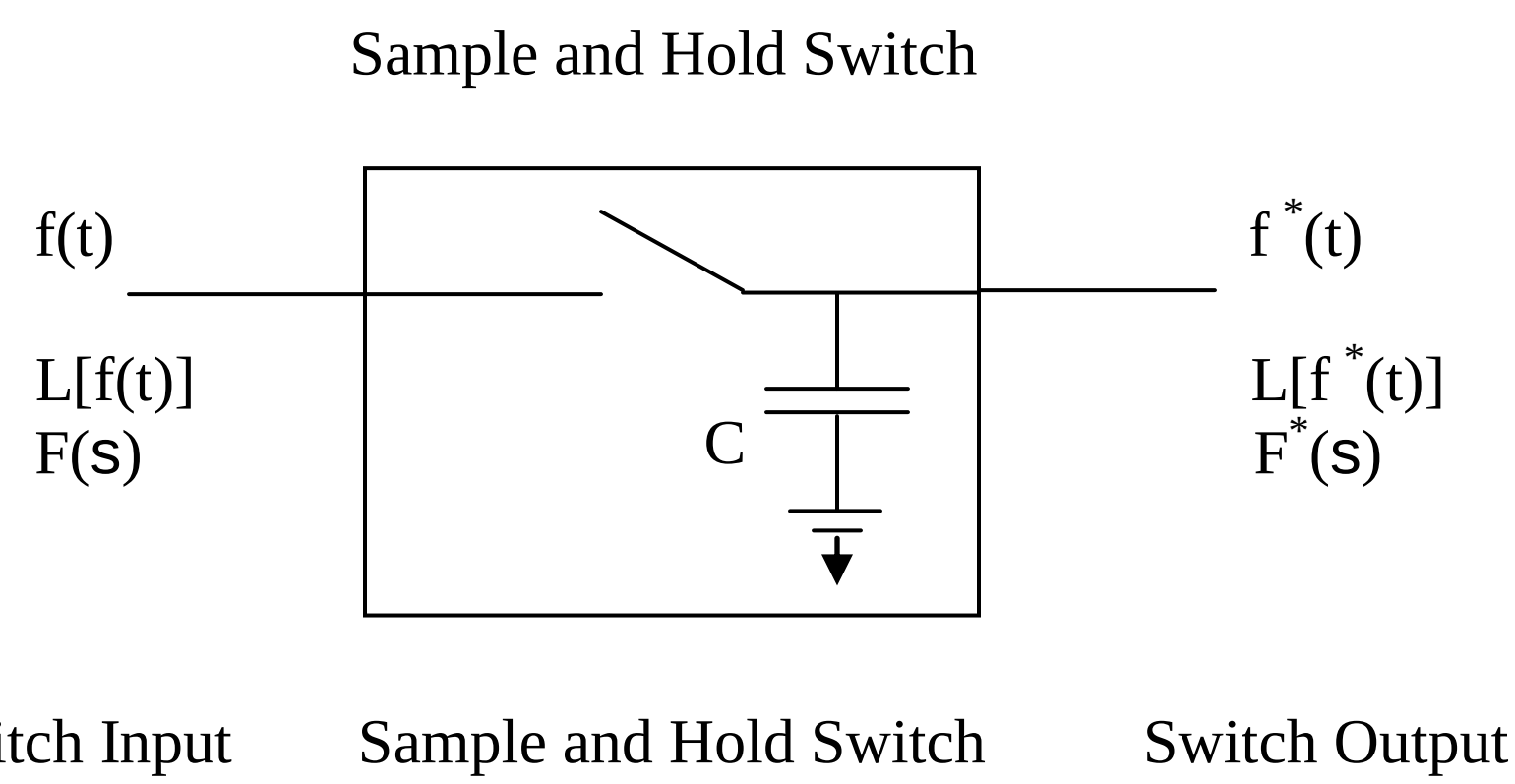

where

- $f(t)$ = Switch input function, a continuous function of time, t
- $f^*(t)$ = Switch output function where output values change only at t = 0, Δt, 2Δt, 3Δt, …
- $L[f(t)] = F(s)$ , The Laplace Transform of the switch input function, f (t)
- $L[f^*(t)] = F^*(s)$ , The Laplace Transform of the switch output function, $f^*(t)$
- C is the sample and hold mechanism depicted in the diagram as a capacitor, C, connected to ground. It could be a digital to analog converter, etc.
- s = The Laplace Transform function variable
- s = The $K_{\Delta t}$ Transform function variable

Derivation of the Laplace Transform of the output of a sample and hold switch using the $K_{\Delta t}$ Transform

Note Diagram 5.12-1 above.

$$f^*(t) = {}_1\sum_{n=0}^{\infty} f(n\Delta t)[U(t-n\Delta t) - U(t-[n+1]\Delta t)] \qquad (5.12\text{-}1)$$

Taking the Laplace Transform of Eq 5.12-1

$$L[f^*(t)] = {}_1\sum_{n=0}^{\infty} f(n\Delta t)\left[\frac{e^{-n\Delta t s}}{s} - \frac{e^{-[n+1]\Delta t s}}{s}\right] \qquad (5.12\text{-}2)$$

Simplifying

$$L[f^*(t)] = \frac{1}{s}\,{}_1\sum_{n=0}^{\infty} f(n\Delta t)e^{-n\Delta t s}(1 - e^{-\Delta t s}) \qquad (5.12\text{-}3)$$

$$L[f^*(t)] = \frac{1 - e^{-\Delta t s}}{s} \; {}_1\sum_{n=0}^{\infty} f(n\Delta t)e^{-n\Delta t s} \qquad (5.12\text{-}4)$$

Changing the form of Eq 5.12-4 to that of discrete integration

$$L[f^*(t)] = \frac{1 - e^{-\Delta t s}}{s\Delta t} \; {}_1\int_0^{\infty} f(n\Delta t)e^{-n\Delta t s}\,\Delta t \qquad (5.12\text{-}5)$$

Let

$$t = n\Delta t \qquad (5.12\text{-}6)$$

From Eq 5.12-5 and Eq 5.12-6

$$\mathbf{L[f^*(t)] = \frac{1 - e^{-s\Delta t}}{s\Delta t} \; {}_{\Delta t}\int_0^{\infty} f(t)e^{-st}\,\Delta t} \qquad (5.12\text{-}7)$$

**where**

**$f(t)$ = Switch input function, a continuous function of time, t**

**$f^*(t)$ = Switch output discrete function of time**

**$L[f^*(t)]$ = Laplace Transform of the switch output function, $f^*(t)$**

**$\Delta t$ = Sampling interval**

**$t = 0, \Delta t, 2\Delta t, 3\Delta t, \dots$**

To facilitate the discrete integration of Eq 5.12-7, the Interval Calculus identity of $e^{-st}$ is used.

$$e^{-st} = \left[1 + \frac{e^{s\Delta t} - 1}{\Delta t}\Delta t\right]^{-\frac{t}{\Delta t}}, \quad t = 0, \Delta t, 2\Delta t, 3\Delta t, \dots \qquad (5.12\text{-}8)$$

Substituting Eq 5.12-8 into Eq 5.12-7

$$L[f^*(t)] = \frac{1 - e^{-s\Delta t}}{s\Delta t} \; {}_{\Delta t}\int_0^{\infty} f(t)\left[1 + \frac{e^{s\Delta t} - 1}{\Delta t}\Delta t\right]^{-\frac{t}{\Delta t}} \Delta t \qquad (5.12\text{-}9)$$

Changing the form of Eq 5.12-9

$$L[f^*(t)] = \frac{1 - e^{-s\Delta t}}{s\Delta t}\left[1 + \frac{e^{s\Delta t} - 1}{\Delta t}\Delta t\right] \; {}_{\Delta t}\int_0^{\infty} f(t)\left[1 + \frac{e^{s\Delta t} - 1}{\Delta t}\Delta t\right]^{-\frac{t}{\Delta t}-1} \Delta t \qquad (5.12\text{-}10)$$

Simplifying

$$L[f^*(t)] = \frac{1 - e^{-\mathit{s}\Delta t}}{\mathit{s}\Delta t}\, e^{\mathit{s}\Delta t} \int_0^\infty {}_{\Delta t} f(t)\, [1 + \frac{e^{\mathit{s}\Delta t} - 1}{\Delta t}\Delta t]^{-\frac{t + \Delta t}{\Delta t}} \Delta t \qquad (5.12\text{-}11)$$

$$L[f^*(t)] = \frac{e^{\mathit{s}\Delta t} - 1}{\mathit{s}\Delta t} \int_0^\infty {}_{\Delta t} f(t)\, [1 + \frac{e^{\mathit{s}\Delta t} - 1}{\Delta t}\Delta t]^{-\frac{t + \Delta t}{\Delta t}} \Delta t \qquad (5.12\text{-}12)$$

Let

$$s = \frac{e^{\mathit{s}\Delta t} - 1}{\Delta t} \qquad (5.12\text{-}13)$$

Note – $\mathit{s}$ and s are two different variables

Substituting Eq 5.12-13 into Eq 5.12-12

$$L[f^*(t)] = \frac{1}{\mathit{s}}\, s \int_0^\infty {}_{\Delta t} f(t)\, [1 + s\Delta t]^{-\frac{t + \Delta t}{\Delta t}} \Delta t \qquad (5.12\text{-}14)$$

The integral in Eq 5.12-14 is recognized as being the $K_{\Delta t}$ Transform of f(t)

$$L[f^*(t)] = \frac{1}{\mathit{s}}\, sK_{\Delta t}[f(t)]\, \Big|_{s = \frac{e^{\mathit{s}\Delta t} - 1}{\Delta t}} \qquad (5.12\text{-}15)$$

, The $K_{\Delta t}$ Transform, $K_{\Delta t}[f(t)]$, is a function of s

Changing the form of Eq 5.12-15

$$L[f^*(t)] = \frac{1}{\mathit{s}}\, [\, sK_{\Delta t}[f(t)] - f(0) + f(0) \Big|_{s = \frac{e^{\mathit{s}\Delta t} - 1}{\Delta t}} ] \qquad (5.12\text{-}16)$$

$$K_{\Delta t}[D_{\Delta t}f(t)] = sK_{\Delta t}[f(t)] - f(0) \qquad (5.12\text{-}17)$$

Substituting Eq 5.12-17 into Eq 5.12-16

$$L[f^*(t)] = \frac{1}{\mathit{s}}\, [\, K_{\Delta t}[D_{\Delta t}f(t)]\, \Big|_{s = \frac{e^{\mathit{s}\Delta t} - 1}{\Delta t}} + f(0)\, ] \qquad (5.12\text{-}18)$$

$$L[f^*(t)] = \frac{f(0)}{s} + \frac{1}{s}[\, K_{\Delta t}[D_{\Delta t}F(s)] \,|_{s = \frac{e^{s\Delta t} - 1}{\Delta t}} ] \qquad (5.12\text{-}19)$$

$$f^*(t) = L^{-1}[F^*(s)] \qquad (5.12\text{-}20)$$

where

$L^{-1}[F(s)]$ = Inverse Laplace Transform of $F(s)$

Then

From Eq 5.12-15, Eq 5.12-19, and Eq 5.12-20

$$\mathbf{L[f^*(t)] = \frac{1}{s} sK_{\Delta t}[f(t)] \,|_{s = \frac{e^{s\Delta t} - 1}{\Delta t}}} \text{ , } \mathbf{K_{\Delta t} \text{ Transforms are a function of s}} \qquad (5.12\text{-}21)$$

**or**

$$\mathbf{F^*(s) = \frac{1}{s} sK_{\Delta t}[L^{-1}[F(s)]] \,|_{s = \frac{e^{s\Delta t} - 1}{\Delta t}}} \text{ , } \mathbf{K_{\Delta t} \text{ Transforms are a function of s}} \qquad (5.12\text{-}22)$$

**or**

$$\mathbf{L[f^*(t)] = \frac{f(0)}{s} + \frac{1}{s}[\, K_{\Delta t}[D_{\Delta t}f(t)] \,|_{s = \frac{e^{s\Delta t} - 1}{\Delta t}} ]} \text{ , } \mathbf{K_{\Delta t} \text{ Transforms are a function of s}} \qquad (5.12\text{-}23)$$

**where**

**f(t) = Switch input function, a continuous function of time, t**

**$L[f(t)] = F(s)$ , the Laplace Transform of the switch input function, f(t), a function of $s$**

**$K_{\Delta t}[f(t)] = K_{\Delta t}$ Transform of the input function f(t) , a function of s**

**$f^*(t)$ = Switch output discrete function of time where t = 0, Δt, 2Δt, 3Δt, …, a sample and hold shaped waveform**

**$L[f^*(t)] = F^*(s)$ , the Laplace Transform of the switch output function, $f^*(t)$, a function of $s$**

**t = 0, Δt, 2Δt, 3Δt, …**

**Δt = Sampling interval**

**$s$ = The Laplace Transform function variable**

**s = The $K_{\Delta t}$ Transform function variable**

**<u>Note</u> – The Laplace Transform of the output of a sample and hold switch, $L[f^*(t)]$, is a function of the $K_{\Delta t}$ Transform of the input function, f(t), to the sample and hold switch.**

From Eq 5.12-2 and Eq 5.12-21

$$L[f^*(t)] = \frac{1}{s} sK_{\Delta t}[f(t)]\,\Big|_{s = \frac{e^{s\Delta t} - 1}{\Delta t}} = \sum_{n=0}^{\infty} f(n\Delta t)\left[\frac{e^{-ns\Delta t}}{s} - \frac{e^{-(n+1)s\Delta t}}{s}\right] \qquad (5.12\text{-}24)$$

The Inverse Laplace Transform of Eq 5.12-24 is shown in the following diagram, Diagram 5.12-2.

Diagram 5.12-2 The output of a sample and hold switch (Laplace Transform system analysis)

Sampled f(t)

$$f^*(t) = L^{-1}[\,L[f^*(t)]] = L^{-1}\left[\frac{1}{s} sK_{\Delta t}[f(t)]\Big|_{s = \frac{e^{s\Delta t} - 1}{\Delta t}}\right]$$

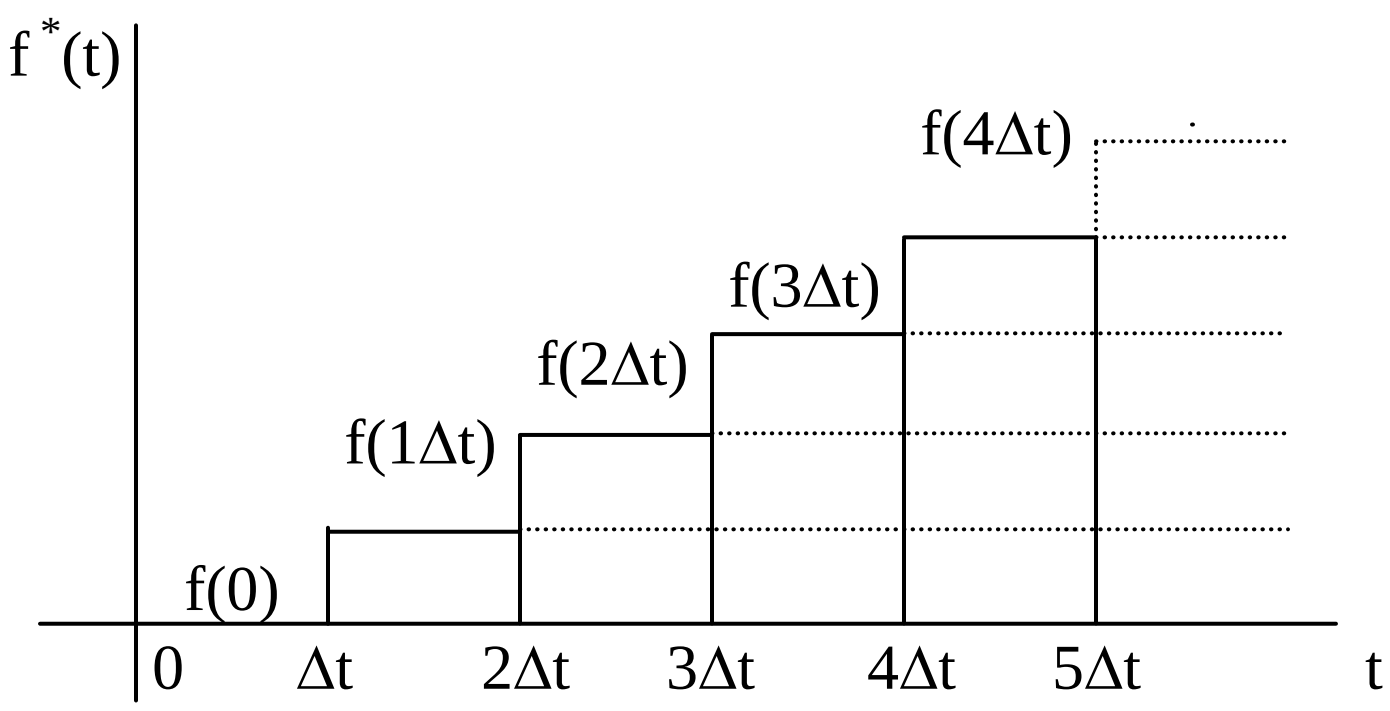

The following three Examples, Example 5.12-1 thru Example 5.12-3, demonstrate the use of the Laplace Transform Equations, Eq 5.12-21 thru Eq 5.12-23, in the solution of problems involving sample and hold sampled-data systems.

Introducing sample and hold sampling into a continuous time system using a Laplace Transform Methodology

Sample and hold sampling can be introduced into continuous time systems using the previously described sample and hold switch. The sample and hold switch is described by any one of the equations, Eq 5.12-21 thru Eq 5.12-23. Demonstrations of sample and hold sampling using a Laplace Transform Methodology are shown in Example 5.12-1 thru Example 5.12-3.

Example 5.12-1

Find the output of a sample and hold switch, $f^*(t)$, for a switch input of f(t) = t using Eq 5.12-23.

Sample and hold switch

$f(t) = t$ 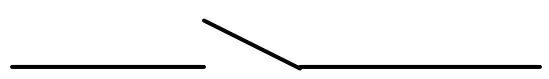 $f^*(t)$

$$L[f^*(t)] = \frac{f(0)}{s} + \frac{1}{s}\left[ K_{\Delta t}[D_{\Delta t}f(t)]\,\Big|_{s = \frac{e^{s\Delta t} - 1}{\Delta t}} \right] \qquad 1)$$

$$f(t) = t \qquad 2)$$

$$f(0) = 0 \qquad 3)$$

$$D_{\Delta t}[t] = D_{\Delta t}\,[t]^{1}_{\Delta x} = 1 \qquad 4)$$

$$K_{\Delta t}[1] = \frac{1}{s} \qquad 5)$$

Substituting into Eq 1)

$$L[f^*(t)] = \frac{0}{s} + \frac{1}{s}\left[\frac{1}{s}\right]\Big|_{s = \frac{e^{s\Delta t} - 1}{\Delta t}} = \frac{1}{s}\left[\frac{1}{\frac{e^{s\Delta t} - 1}{\Delta t}}\right] = \frac{\Delta t}{s(\,e^{s\Delta t} - 1\,)} \qquad 6)$$

Then

$$L[f^*(t)] = \frac{\Delta t}{s(\,e^{s\Delta t} - 1\,)} \qquad 7)$$

Expanding Eq 7)

$$L[f^*(t)] = e^{-s\Delta t}\frac{\Delta t}{s} + e^{-2s\Delta t}\frac{\Delta t}{s} + e^{-3s\Delta t}\frac{\Delta t}{s} + e^{-4s\Delta t}\frac{\Delta t}{s} + e^{-5s\Delta t}\frac{\Delta t}{s} + \dots \qquad 8)$$

Taking the Inverse Laplace Transform of Eq 8

$$f^*(t) = L^{-1}[F^*(s)] = 0U(t) + \Delta tU(t-\Delta t) + \Delta tU(t-2\Delta t) + \Delta tU(t-3\Delta t) + \Delta tU(t-4\Delta t) + \Delta tU(t-5\Delta t) + \dots \qquad 9)$$

$$\mathbf{f^*(t) = {}_{1}\sum_{n=0}^{\infty} \Delta t\, U(t-[n+1]\Delta t)} \qquad 10)$$

**or**

**$f^*(t) = t$ , $t = 0, \Delta t, 2\Delta t, 3\Delta t, \dots$ , t is a sample and hold shaped waveform** 11)

The above result has been checked using $K_{\Delta t}$ Transforms in Example 5.12-6 part 1.
Good check

Diagram 5.12-3 Diagram of the output of sample and hold switch, $f^*(t)$, with a switch input of $f(t) = t$

Sample and hold switch output, $f^*(t)$

$$f^*(t) = {}_1\sum_{n=0}^{\infty} \Delta t\, U(t-[n+1]\Delta t)$$

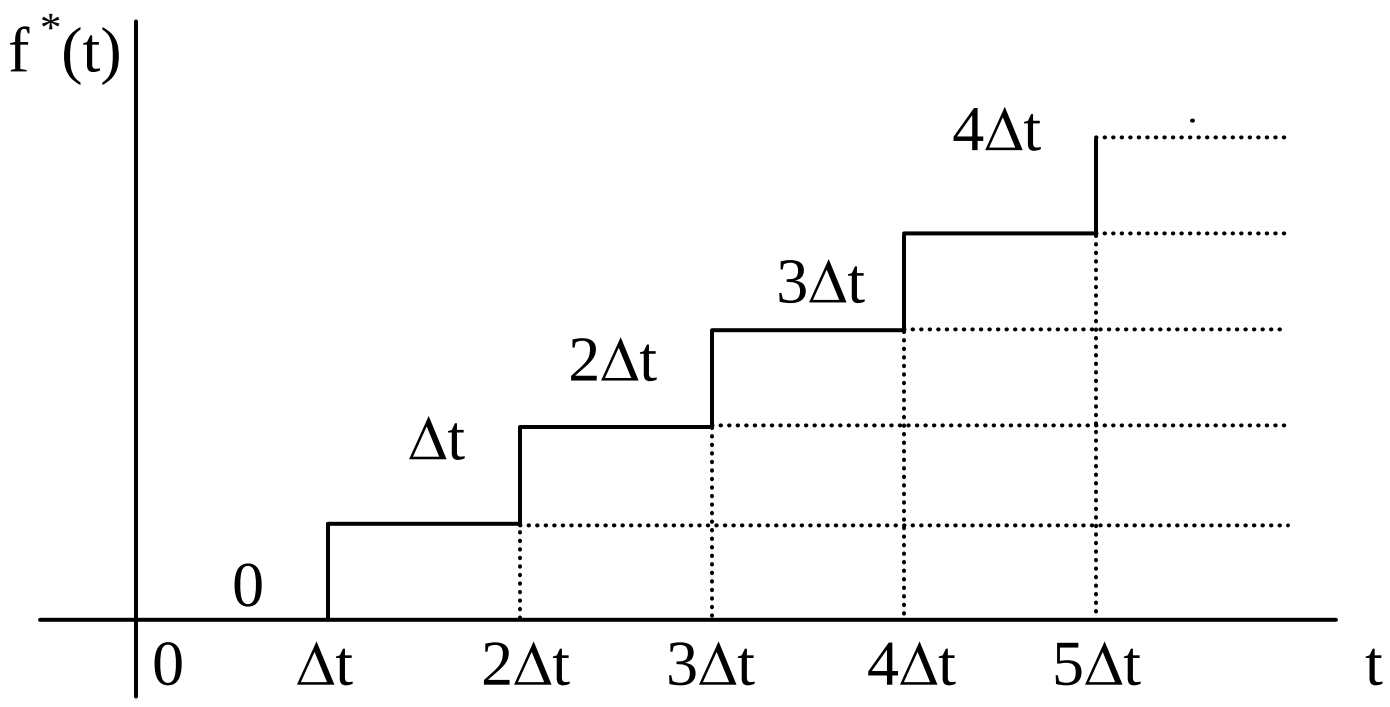

$f^*(t) = t$ , $t = 0, \Delta t, 2\Delta t, 3\Delta t, \ldots$ , a sample and hold shaped waveform

Example 5.12-2

Find the output of an integrator, c(t), where the input of the integrator is the output, $f^*(t)$, of the sample and hold switch of Example 5.1-1. The integrator is initially passive.

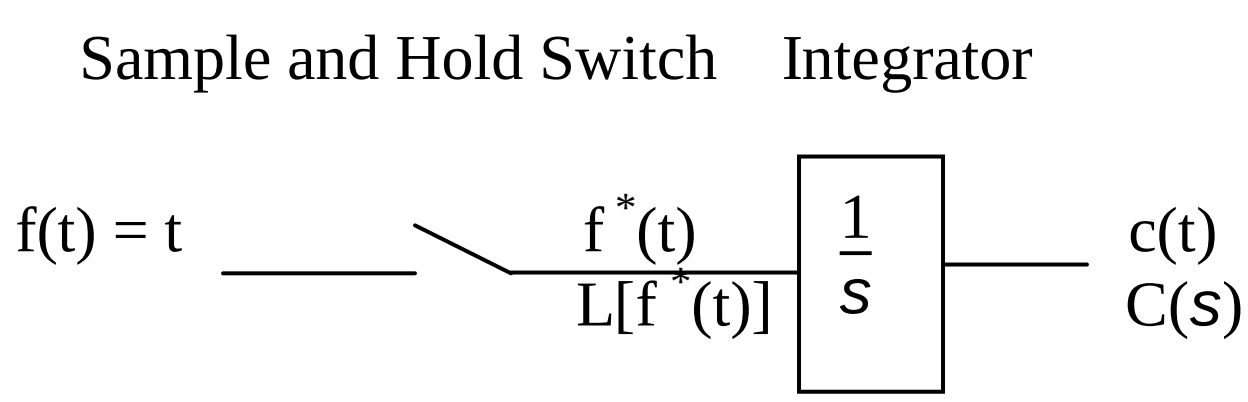

From Eq 5.12-23

$$L[f^*(t)] = \frac{f(0)}{s} + \frac{1}{s}\left[\, K_{\Delta t}[D_{\Delta t} f(t)]\, \Big|_{s = \frac{e^{s\Delta t} - 1}{\Delta t}} \right] \tag{1}$$

$$f(t) = t \tag{2}$$

$$f(0) = 0 \tag{3}$$

$$D_{\Delta t}[t] = D_{\Delta t}\,[t]_{\Delta x}^{1} = 1 \tag{4}$$

$$K_{\Delta t}[1] = \frac{1}{s} \tag{5}$$

Substituting into Eq 1)

$$L[f^*(t)] = \frac{0}{s} + \frac{1}{s}\left[\frac{1}{s}\Big|_{s = \frac{e^{s\Delta t} - 1}{\Delta t}}\right] = \frac{1}{s}\left[\frac{1}{\frac{e^{s\Delta t} - 1}{\Delta t}}\right] = \frac{\Delta t}{s(\, e^{s\Delta t} - 1\,)} \tag{6}$$

$$L[f^*(t)] = \frac{\Delta t}{s(e^{s\Delta t} - 1)} \qquad 7)$$

$$C(s) = \frac{1}{s} L[f^*(t)] \qquad 8)$$

$$C(s) = \frac{\Delta t}{s^2(e^{s\Delta t} - 1)} \qquad 9)$$

Expanding Eq 9)

$$C(s) = e^{-s\Delta t}\frac{\Delta t}{s^2} + e^{-2s\Delta t}\frac{\Delta t}{s^2} + e^{-3s\Delta t}\frac{\Delta t}{s^2} + e^{-4s\Delta t}\frac{\Delta t}{s^2} + e^{-5s\Delta t}\frac{\Delta t}{s^2} + \ldots \qquad 10)$$

Taking the Inverse Laplace Transform of Eq 10

$$c(t) = 0U(t) + \Delta t(t-\Delta t)U(t-\Delta t) + \Delta t(t-2\Delta t)U(t-2\Delta t) + \Delta t(t-3\Delta t)U(t-3\Delta t) \qquad 11)$$
$$+ \Delta t(t-4\Delta t)U(t-4\Delta t) + (t-5\Delta t)U(t-5\Delta t) + \ldots$$

$$\mathbf{c(t) = {}_1\sum_{n=0}^{\infty} \Delta t\,(t-[n+1]\Delta t)U(t-[n+1]\Delta t)} \qquad 12)$$

**or**

$$\mathbf{c(t) = \frac{t(t-\Delta t)}{2}, \quad t = 0, \Delta t, 2\Delta t, 3\Delta t, \ldots} \qquad 13)$$

The above result has been checked using $K_{\Delta t}$ Transforms in Example 5.12-6 part 2.

Good check

Diagram 5.12-4  Diagram of the output, c(t), with a sample and hold switch input, f(t) = t

Sampled system output, c(t)

$$c(t) = {}_1\sum_{n=0}^{\infty} \Delta t\,(t-[n+1]\Delta t)U(t-[n+1]\Delta t)$$

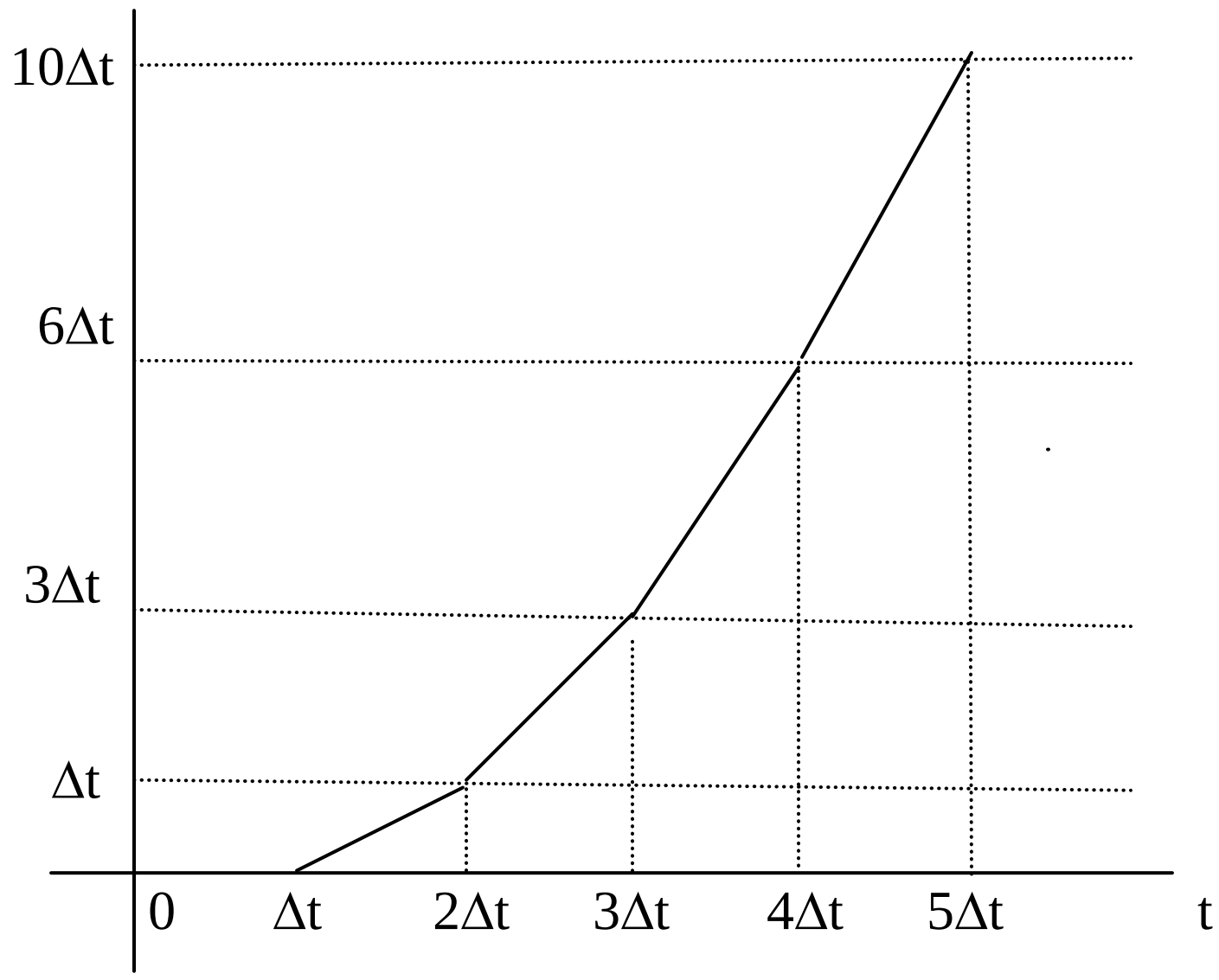

$$c(t) = \frac{t(t-\Delta t)}{2}, \quad t = 0, \Delta t, 3\Delta t, 3\Delta t, \ldots$$

Example 5.12-3

Find the output function, c(t), for the following system where the input is U(t). Use Laplace Transforms. The sample and hold switches are synchronized and the integrators are initially passive.

Three sample and hold switches and two integrators

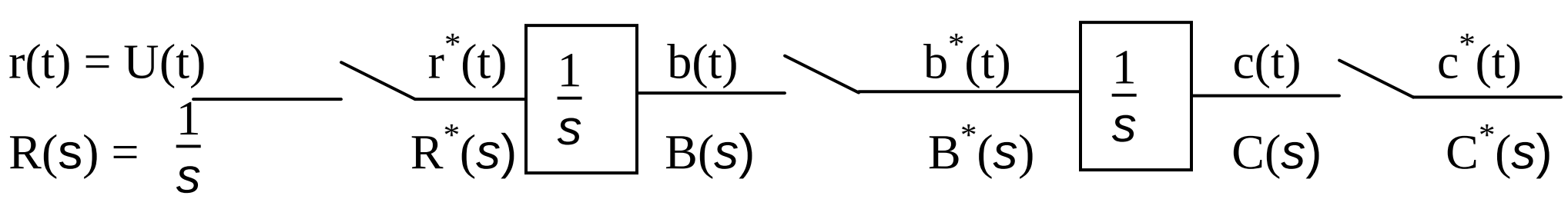

From Eq 5.12-22

$$R^*(s) = \frac{1}{s}\, sK_{\Delta t}[L^{-1}[R(s)]]\,\Big|_{s\,=\,\frac{e^{s\Delta t}-1}{\Delta t}} = \frac{1}{s}\, sK_{\Delta t}[U(t)]\,\Big|_{s\,=\,\frac{e^{s\Delta t}-1}{\Delta t}} = \frac{1}{s}\, s\, \frac{1}{s}\,\Big|_{s\,=\,\frac{e^{s\Delta t}-1}{\Delta t}} = \frac{1}{s} \qquad 1)$$

$$R^*(s) = \frac{1}{s} \qquad 2)$$

$$B(s) = \frac{1}{s} R^*(s) = \frac{1}{s}\frac{1}{s} = \frac{1}{s^2} \qquad 3)$$

$$B^*(s) = \frac{1}{s} sK_{\Delta t}[L^{-1}[\frac{1}{s^2}] \Big|_{s = \frac{e^{s\Delta t} - 1}{\Delta t}} = \frac{1}{s} sK_{\Delta t}[t] \Big|_{s = \frac{e^{s\Delta t} - 1}{\Delta t}} = \frac{1}{s} s \frac{1}{s^2} \Big|_{s = \frac{e^{s\Delta t} - 1}{\Delta t}} = \frac{\Delta t}{s(e^{s\Delta t} -1)} \qquad 4)$$

$$B^*(s) = \frac{\Delta t}{s(e^{s\Delta t} -1)} \qquad 5)$$

$$C(s) = \frac{1}{s} B*(s) = [\frac{1}{s}] \frac{\Delta t}{s(e^{s\Delta t} -1)} = \frac{\Delta t}{s^2(e^{s\Delta t} -1)} \qquad 6)$$

$$C(s) = \frac{\Delta t}{s^2(e^{s\Delta t} -1)} \qquad 7)$$

Expanding Eq 7

$$C(s) = \Delta t[ \frac{e^{-s\Delta t}}{s^2} + \frac{e^{-2s\Delta t}}{s^2} + \frac{e^{-3s\Delta t}}{s^2} + \frac{e^{-4s\Delta t}}{s^2} + \ldots ] \qquad 8)$$

Taking the Inverse Laplace Transform of Eq 8)

$$c(t) = L^{-1}[C(s)] = \Delta t[ (t-\Delta t)U(t-\Delta T) + (t-2\Delta t)U(t-2\Delta T) + (t-3\Delta t)U(t-3\Delta T) + (t-4\Delta t)U(t-4\Delta T) +\ldots] \qquad 9)$$

$$C^*(s) = \frac{1}{s} sK_{\Delta t}[c(t)] \Big|_{s = \frac{e^{s\Delta t} - 1}{\Delta t}} \qquad 10)$$

$$C^*(s) = \frac{\Delta t}{s} sK_{\Delta t}[t-\Delta t)U(t-\Delta T) + (t-2\Delta t)U(t-2\Delta T) + (t-3\Delta t)U(t-3\Delta T) + (t-4\Delta t)U(t-4\Delta T) +\ldots] \Big|_{s = \frac{e^{s\Delta t} - 1}{\Delta t}} \qquad 11)$$

$$C^*(s) = \frac{\Delta t}{s} s [\frac{(1+s\Delta t)^{-1}}{s^2} + \frac{(1+s\Delta t)^{-2}}{s^2} + \frac{(1+s\Delta t)^{-3}}{s^2} \; \frac{(1+s\Delta t)^{-4}}{s^2} + \ldots] \Big|_{s = \frac{e^{s\Delta t} - 1}{\Delta t}} \qquad 12)$$

$$C^*(s) = \frac{\Delta t^2}{s} ( \frac{1}{e^{s\Delta t} -1} ) [ e^{-s\Delta t} + e^{-2s\Delta t} + e^{-3s\Delta t} + e^{-4s\Delta t} +\ldots] \qquad 13)$$

$$C^*(s) = \frac{\Delta t^2}{s} ( \frac{1}{e^{s\Delta t}-1} )^2 \qquad 14)$$

Expanding $( \frac{1}{e^{s\Delta t}-1} )^2$

$$( \frac{1}{e^{s\Delta t}-1} )^2 = e^{-2s\Delta t} + 2e^{-3s\Delta t} + 3e^{-4s\Delta t} + 4e^{-5s\Delta t} + 5e^{-6s\Delta t} + \ldots \qquad 15)$$

Substituting Eq 15 into Eq 14

$$C^*(s) = \Delta t^2\left[\frac{e^{-2s\Delta t}}{s}{}^{t} + 2\frac{e^{-3s\Delta t}}{s} + 3\frac{e^{-4s\Delta t}}{s} + 4\frac{e^{-5s\Delta t}}{s} + 5\frac{e^{-6s\Delta t}}{s} + \ldots\right] \qquad 16)$$

Taking the Laplace Transform of Eq 16

$$c^*(t) = \Delta t^2[\ 1U(t-2\Delta t) + 2U(t-3\Delta t) + 3U(t-4\Delta t) + 4U(t-5\Delta t) + \ldots\ ] \qquad 17)$$

$$c^*(t) = \Delta t^2 \sum_{n=0}^{\infty}(n+1)U(t-[n+2]\Delta t) \qquad 18)$$

Calculate c(t) from Eq 18 for t = 0, Δt, 2Δt, 3Δt, 4Δt, 5Δt, 6Δt

| t | c(t) |
|---|---|
| 0 | 0 |
| $\Delta t$ | 0 |
| $2\Delta t$ | $\Delta t^2$ |
| $3\Delta t$ | $3\Delta t^2$ |
| $4\Delta t$ | $6\Delta t^2$ |
| $5\Delta t$ | $10\Delta t^2$ |
| $6\Delta t$ | $15\Delta t^2$ |

$$\mathbf{c^*(t) = \Delta t^2 \sum_{n=0}^{\infty}(n+1)U(t-[n+2]\Delta t)} \qquad 19)$$

**or**

$$\mathbf{c^*(t) = \frac{t(t-\Delta t)}{2}\ ,\quad t = 0, \Delta t, 2\Delta t, 3\Delta t, 4\Delta t, 5\Delta t, 6\Delta t} \qquad 20)$$

A check of the above result has been obtained using Example 5.12-6 part 2.

Good check

The output of a sample and hold switch can also be described using the Z Transform. Note the derivation below.

Derivation of the Laplace Transform of the output of a sample and hold switch using the Z Transform

Sample and hold switch

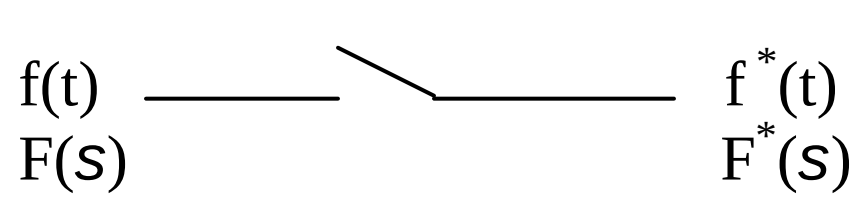

$$f^{*}(t) = {}_{1}\sum_{n=0}^{\infty} f(n\Delta t)[U(t-n\Delta t) - U(t-[n+1]\Delta t)] \qquad (5.12\text{-}25)$$

Taking the Laplace Transform of Eq 5.12-25

$$L[f^{*}(t)] = {}_{1}\sum_{n=0}^{\infty} f(n\Delta t)[\frac{e^{-n\Delta t s}}{s} - \frac{e^{-(n+1)\Delta t s}}{s}] \qquad (5.12\text{-}26)$$

Simplifying

$$L[f^{*}(t)] = \frac{1}{s}\,{}_{1}\sum_{n=0}^{\infty} f(n\Delta t)e^{-n\Delta t s}\,(1 - e^{-\Delta t s}) \qquad (5.12\text{-}27)$$

$$L[f^{*}(t)] = \frac{1 - e^{-\Delta t s}}{s}\,{}_{1}\sum_{n=0}^{\infty} f(n\Delta t)e^{-n\Delta t s} \qquad (5.12\text{-}28)$$

Changing the form of Eq 5.12-28 to that of discrete integration

$$L[f^{*}(t)] = \frac{1 - e^{-\Delta t s}}{s\Delta t}\,{}_{1}\int_{0}^{\infty} f(n\Delta t)e^{-n\Delta t s}\,\Delta t \qquad (5.12\text{-}29)$$

Let

$$t = n\Delta t \qquad (5.12\text{-}30)$$

From Eq 5.12-29 and Eq 5.12-30

$$\mathbf{L[f^{*}(t)] = \frac{1 - e^{-s\Delta t}}{s\Delta t}\,{}_{\Delta t}\int_{0}^{\infty} f(t)\,e^{-st}\,\Delta t} \qquad (5.12\text{-}31)$$

**where**

**f(t) = Switch input function, a function of time, t**

**$f^{*}(t)$ = Switch sample and hold output function of time, t**

**$L[f^{*}(t)]$ , the Laplace Transform of the switch output function, $f^{*}(t)$**

**t = 0, Δt, 2Δt, 3Δt, …**

**Δt = Sampling interval**

To facilitate the discrete integration of Eq 5.12-31 use the Interval Calculus identity of $e^{-st}$

$$e^{-st} = [1 + \frac{e^{s\Delta t} - 1}{\Delta t}\Delta t]^{-\frac{t}{\Delta t}}, \quad t = 0, \Delta t, 2\Delta t, 3\Delta t, \ldots \tag{5.12-32}$$

Substituting Eq 5.12-32 into Eq 5.12-31

$$L[f^*(t)] = \frac{1 - e^{-s\Delta t}}{s\Delta t} \, {}_{\Delta t}\!\int_0^{\infty} f(t)\,[1 + \frac{e^{s\Delta t} - 1}{\Delta t}\Delta t]^{-\frac{t}{\Delta t}}\,\Delta t \tag{5.12-33}$$

Let

$$s = \frac{e^{s\Delta t} - 1}{\Delta t} \tag{5.12-34}$$

<u>Note</u> – $s$ and s are two different variables

Substituting Eq 5.12-34 into Eq 5.12-33

$$L[f^*(t)] = \frac{1 - e^{-s\Delta t}}{s\Delta t} \, {}_{\Delta t}\!\int_0^{\infty} f(t)\,[1+s\Delta t]^{-\frac{t}{\Delta t}}\,\Delta t \tag{5.12-35}$$

The integral in Eq 5.12-35 is recognized as being the $J_{\Delta t}$ Transform of f(t)

$$L[f^*(t)] = \frac{1 - e^{-s\Delta t}}{s\Delta t} J_{\Delta t}[f(t)]\,\Big|_{s = \frac{e^{s\Delta t} - 1}{\Delta t}}, \quad K_{\Delta t}[f(t)] \text{ is a function of s} \tag{5.12-36}$$

Then

$$\mathbf{L[f^*(t)] = \frac{1 - e^{-s\Delta t}}{s\Delta t} J_{\Delta t}[f(t)]\,\Big|_{s = \frac{e^{s\Delta t} - 1}{\Delta t}}, \quad J_{\Delta t}[f(t)] \text{ is a function of s}} \tag{5.12-37}$$

**where**

**f(t) = Switch input function, a continuous function of time, t**
**L[f(t)] = F($s$) , the Laplace Transform of the switch input function, f (t), a function of $s$**
**$J_{\Delta t}$[f(t)] = $J_{\Delta t}$ Transform of the input function, f(t), a function of s**
**$f^*$(t) = Switch sample and hold output function of time, t**
**L[$f^*$(t)] = $F^*$($s$) , the Laplace Transform of the switch output function, $f^*$(t), a function of $s$**
**Δt = Sampling interval**
**t = 0, Δt, 2Δt, 3Δt, …**
**$s$ = The Laplace Transform function variable**
**s = The $K_{\Delta t}$ Transform function variable**

**<u>Note</u> – The Laplace Transform of the output of a sample and hold switch, L[$f^*$(t)], is a function of the $J_{\Delta t}$ Transform of the input, f(t), to the sample and hold switch.**

From Eq 5.12-36 and Eq 5.12-26

$$L[f^*(t)] = \frac{1 - e^{-s\Delta t}}{s\Delta t} J_{\Delta t}[f(t)] \Big|_{s = \frac{e^{s\Delta t} - 1}{\Delta t}} = \sum_{n=0}^{\infty} f(n\Delta t)\left[\frac{e^{-ns\Delta t}}{s} - \frac{e^{-(n+1)s\Delta t}}{s}\right] \quad (5.12\text{-}38)$$

Let

$$z = e^{s\Delta t} \quad (5.12\text{-}39)$$

$$T = \Delta t \quad (5.12\text{-}40)$$

Substituting Eq 5.12-39 and Eq 5.12-40 into Eq 5.12-38

$$L[f^*(t)] = \frac{1 - z^{-1}}{sT} J_{\Delta t}[f(t)] \Big|_{s = \frac{z - 1}{T}} \quad (5.12\text{-}41)$$

$$Z[f(t)] = \frac{1}{T} J_{\Delta t}[f(t)] \Big|_{s = \frac{z - 1}{T}} \quad (5.12\text{-}42)$$

where
$T = \Delta t$
$t = nT$
$n = 0, 1, 2, 3, \ldots$

Substituting 5.12-42 into Eq 5.12-41

$$L[f^*(t)] = \frac{1 - z^{-1}}{s} Z[f(t)] = \frac{1}{s}\,\frac{z - 1}{z} Z[f(t)] \quad (5.12\text{-}43)$$

$$L[f^*(t)] = \frac{1}{s}\,\frac{z - 1}{z} Z[f(t)] \quad (5.12\text{-}44)$$

Then

$$\mathbf{L[f^*(t)] = \frac{1}{s}\,\frac{z - 1}{z} Z[f(t)]\Big|_{z = e^{sT}}} \quad (5.12\text{-}45)$$

**or**

$$\mathbf{F^*(s) = \frac{1}{s}\,\frac{z - 1}{z} Z[L^{-1}[F(s)]\Big|_{z = e^{sT}}} \quad (5.12\text{-}46)$$

**where**
**f(t) = Sample and hold switch input function, a continuous function of time, t**
**Z[f(t)] = F(z) , The Z Transform of the switch input function, f (t), a function of z**
**$f^*(t)$ = Sample and hold switch output function, a sample and hold function of time, t**
**$L[f^*(t)] = F^*(s)$ , the Laplace Transform of the switch output function, $f^*(t)$, a function of s**
**$T = \Delta t$ = Sampling interval**

**$t = 0, T, 2T, 3T, \ldots$**
**$F(s)$ = Laplace Transform of the switch input function, f(t)**
**$F^*(s)$ = Laplace Transform of the switch sample and hold output function, $f^*(t)$**
**z = The Z Transform function variable**
**$z = e^{sT}$**
**$s$ = Laplace Transform variable**

<u>**Note**</u> **– The Laplace Transform of the output of a sample and hold switch, $L[f^*(t)]$, is a function of the Z Transform of the input, f(t), to the sample and hold switch.**

There is another useful equation to calculate Laplace Transform of the sample and hold waveform, $f^*(t)$, from the Z Transform of f(t).

$$\Delta f(t) = f(t+\Delta t) - f(t) \qquad (5.12\text{-}47)$$

$$Z[\Delta f(t)] = (z-1)Z[f(t)] - zf(0) \qquad (5.12\text{-}48)$$

Dividing both sides of Eq 5.12-48 by z

$$\frac{1}{z} Z[\Delta f(t)] = \frac{z-1}{z} Z[f(t)] - f(0) \qquad (5.12\text{-}49)$$

Rearranging the terms of Eq 5.12-49

$$\frac{z-1}{z} Z[f(t)] = f(0) + \frac{1}{z} Z[\Delta f(t)] \qquad (5.12\text{-}50)$$

$$\frac{z-1}{z} Z[f(t)] = f(0) + Tz^{-1} Z[\frac{\Delta f(t)}{T}] \qquad (5.12\text{-}51)$$

Multiplying both sides of Eq 5.12-51 by $\frac{1}{s}$ and defining z as $z = e^{sT}$

$$L[f^*(t)] = \frac{1}{s}\,\frac{z-1}{z}\,Z[f(t)]\Big|_{z = e^{sT}} = \frac{1}{s}\,[\,f(0) + Tz^{-1} Z[\frac{\Delta f(t)}{T}]\,\Big|_{z = e^{sT}}\,] \qquad (5.12\text{-}52)$$

$$L[f^*(t)] = \frac{1}{s}[\,f(0) + Tz^{-1} Z[\frac{\Delta f(t)}{T}]\,\Big|_{z = e^{sT}}\,] \qquad (5.12\text{-}53)$$

$$D_T[f(t)] = \frac{\Delta f(t)}{T} \qquad (5.12\text{-}54)$$

where
$D_T[f(t)]$ = Discrete derivative of f(t) with a sampling interval of T

Substituting Eq 5.12-54 into Eq 5.12-53

$$L[f^{*}(t)] = \frac{1}{s}\,[\, f(0) + Tz^{-1}\, Z[D_T[f(t)]\,]\,|_{z = e^{sT}}\,] \qquad (5.12\text{-}55)$$

Then

$$\mathbf{L[f^{*}(t)] = \frac{1}{s}\,[\, f(0) + Tz^{-1}\, Z[D_T[f(t)]\,]\,|_{z = e^{sT}}\,]} \qquad (5.12\text{-}56)$$

**where**

**$f(t)$ = Sample and hold switch input function, a continuous function of time, t**
**$Z[f(t)] = F(z)$, The Z Transform of the switch input function, f (t), a function of z**
**$f^{*}(t)$ = Sample and hold switch output function, a sample and hold function of time, t**
**$L[f^{*}(t)] = F^{*}(s)$, the Laplace Transform of the switch output function, $f^{*}(t)$, a function of $s$**
**$T = \Delta t$ = Sampling interval**
**$t = 0, T, 2T, 3T, \ldots$**
**$F(s)$ = Laplace Transform of the switch input function, f(t)**
**$F^{*}(s)$ = Laplace Transform of the switch sample and hold output function, $f^{*}(t)$**
**z = The Z Transform function variable**
**$z = e^{sT}$**
**$D_T[f(t)] = \dfrac{f(t+\Delta t) - f(t)}{T}$, Discrete derivative of f(t) with a sampling interval of T**
**$s$ = Laplace Transform variable**

**<u>Note</u> – The Laplace Transform of the output of a sample and hold switch, $L[f^{*}(t)]$, is a function of the Z Transform of the input, f(t), to the sample and hold switch.**

Some characteristics of $K_{\Delta t}$ Transform sample and hold switches

$K_{\Delta t}$ Transform Sample and Hold Switch

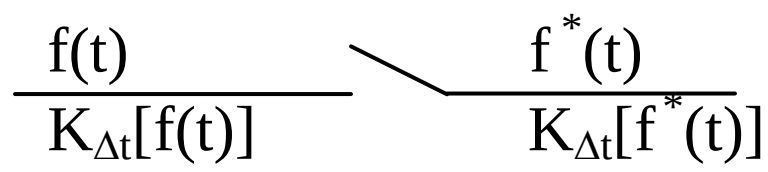

$f(t)$ = Switch input function
$f^*(t)$ = Switch output function, a sample and hold shaped waveform
$0 \le t < \infty$
$\Delta t$ = interval between samples
$f^*(t) = f(n\Delta t)$ for $n\Delta t \le t < [n+1]\Delta t$ , $n = 0, 1, 2, 3, \ldots$
$K_{\Delta t}[f^*(t)] = K_{\Delta t}[f(t)]$

Sample and Hold Shaped Waveform

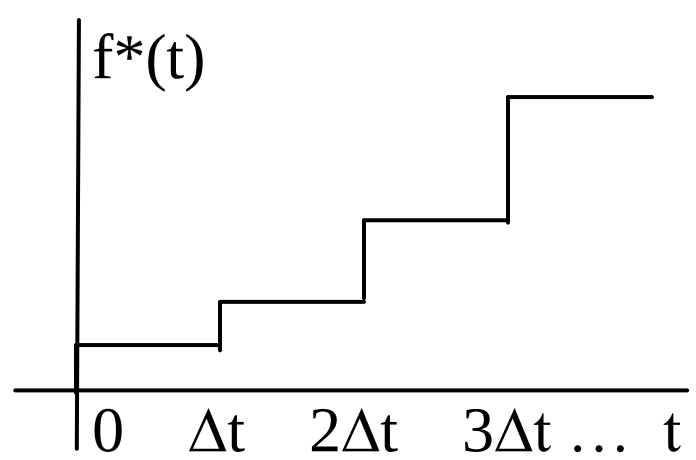

1) $f^*(t) = f(t)$ at $t = 0, \Delta t, 2\Delta t, 3\Delta t, \ldots$
2) $K_{\Delta t}[f^*(t)] = K_{\Delta t}[f(t)]$. If two functions have the same values at $t = 0, \Delta t, 2\Delta t, 3\Delta t, \ldots$, their $K_{\Delta t}$ Transforms are the same.
3) If the input of a $K_{\Delta t}$ Transform sample and hold switch is designated as $f(t) = F(t)$, the switch's output is designatd as $f^*(t) = F(t)$ also. However, it is understood that the two waveforms generally are not the same. The switch output is a sample and hold shaped waveform with values equal to the switch input values at the sampling instants $t = 0, \Delta t, 2\Delta t, 3\Delta t, \ldots$
4) If two or more synchronously actived sample and hold switches are connected in series, the combination of switches acts as one switch.

Example 5.12-4

Find the output of a sample and hold switch, $f^*(t)$, for a switch input of $f(t) = t$ using Eq 5.12-45.

Sample and Hold Switch

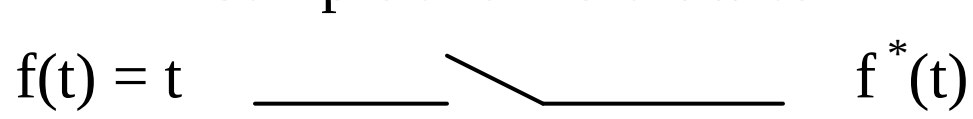

$$L[f^*(t)] = \frac{1}{s}\,\frac{z-1}{z}\,Z[f(t)]\Big|_{z = e^{sT}} \qquad 1)$$

$$f(t) = t \qquad 2)$$

$$Z[t] = \frac{Tz}{(z-1)^2} \qquad 3)$$

Substituting Eq 3 into Eq 1

$$[f^*(t)] = \frac{1}{s}\,\frac{z-1}{z}\,\frac{Tz}{(z-1)^2}\Big|_{z = e^{sT}} = \frac{1}{s}\,\frac{T}{(z-1)}\Big|_{z = e^{sT}} \qquad 4)$$

$$L[f^*(t)] = \frac{T}{s(\,e^{sT} - 1\,)} \qquad 5)$$

Expanding Eq 5

$$L[f^*(t)] = e^{-sT}\frac{T}{s} + e^{-2sT}\frac{T}{s} + e^{-3sT}\frac{T}{s} + e^{-4sT}\frac{T}{s} + e^{-5sT}\frac{T}{s} + \ldots \qquad 6)$$

Taking the Inverse Laplace Transform of Eq 6

$$f^*(t) = L^{-1}[L[f^*(t)]] = TU(t-T) + TU(t-2T) + TU(t-3T) + TU(t-4T) + TU(t-5T) + \ldots \qquad 7)$$

$$\mathbf{f^*(t) = {}_1\sum_{n=0}^{\infty} T\,U(t-[n+1]T)} \qquad 8)$$

<u>Diagram 5.12-5  Diagram of the output of a sample and hold switch, $f^*(t)$, with a switch input of $f(t) = t$</u>

<u>Sample and hold switch output, $f^*(t)$</u>

$$f^*(t) = L^{-1}[L[f^*(t)]] = {}_1\sum_{n=0}^{\infty} T\,U(t-[n+1]T)$$

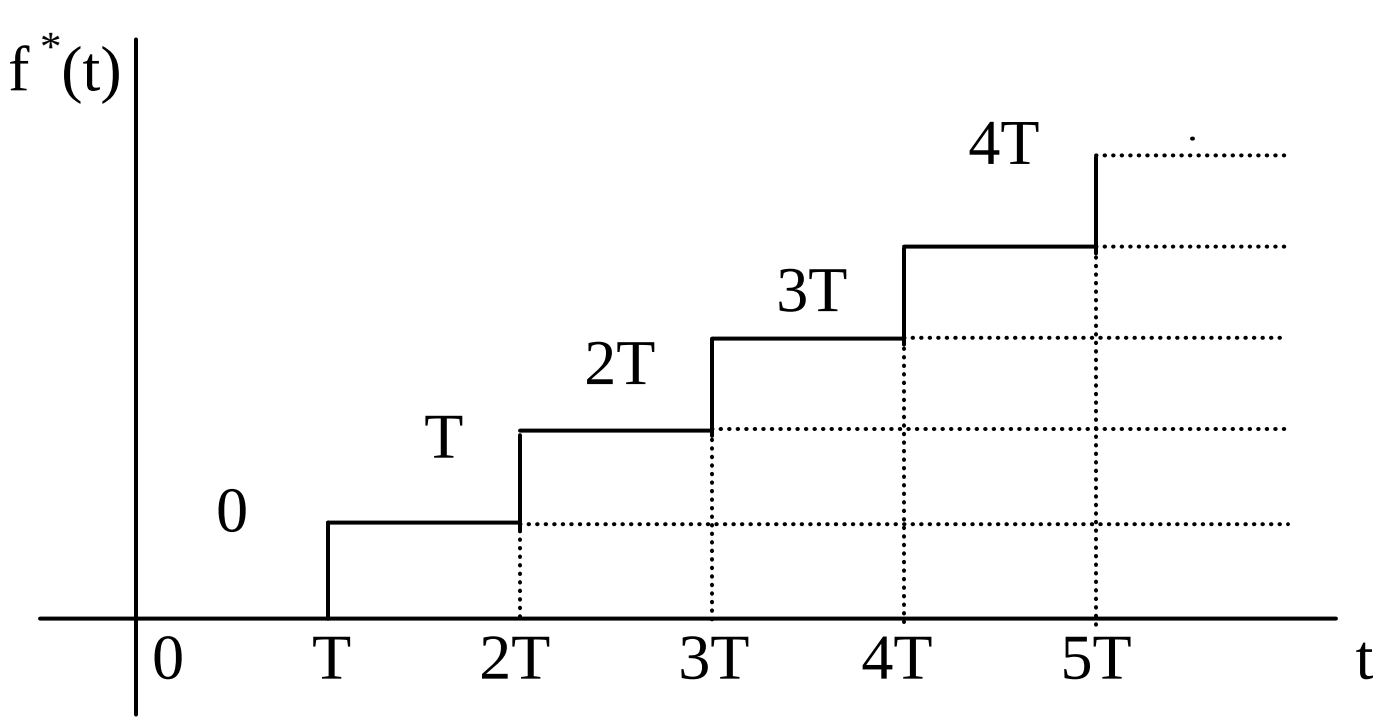

Note that Example 5.12-4 and Example 5.12-1 are the same though solved by different means. Note, also, the results, as expected, are the same. Compare Diagram 5.12-5 to Diagram 5.12-3. $T = \Delta t$.

Example 5.12-5

Find the output of a sample and hold switch, $f^*(t)$, for a switch input of $f(t) = e^{at}$ using Eq 5.12-45.

Sample and Hold Switch

$f(t) = t$ 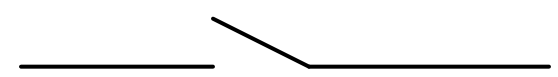 $f^*(t)$

$$L[f^*(t)] = \frac{1}{s}\,\frac{z-1}{z}\,Z[f(t)]\Big|_{z = e^{sT}} \qquad 1)$$

$$f(t) = e^{at} \qquad 2)$$

$$Z[t] = \frac{z}{z - e^{at}} \qquad 3)$$

Substituting Eq 3 into Eq 1

$$L[f^*(t)] = \frac{1}{s}\,\frac{z-1}{z}\,\frac{z}{z-e^{at}}\Big|_{z = e^{sT}} = \frac{1}{s}\,\frac{z-1}{z-e^{at}}\Big|_{z = e^{sT}} \qquad 4)$$

$$L[f^*(t)] = \frac{1}{s}\,\frac{z-1}{z-e^{at}}\Big|_{z = e^{sT}} \qquad 5)$$

Expanding Eq 5

$$L[f^*(t)] = \frac{1}{s} + e^{aT}\frac{e^{-sT}}{s} + e^{2aT}\frac{e^{-2sT}}{s} + e^{3aT}\frac{e^{-3sT}}{s} + e^{4aT}\frac{e^{-4sT}}{s} + \ldots \qquad 6)$$

Taking the Inverse Laplace Transform of Eq 6

$$f^*(t) = L^{-1}[L[f^*(t)]] = U(t) + e^{at}U(t-T) + e^{2at}U(t-2T) + e^{3at}U(t-3T) + e^{4at}U(t-4T) + \ldots \qquad 7)$$

Then

$$\mathbf{f^*(t) = {}_1\sum_{n=0}^{\infty} e^{nat}\,U(t-nT)} \qquad 8)$$

Diagram 5.12-6 Diagram of the output of a sample and hold switch, $f^*(t)$, with a switch input of $f(t) = t$

Sampled $f(t) = e^{at}$

$$f^*(t) = L^{-1}[L[f^*(t)]] = {}_1\sum_{n=0}^{\infty} e^{nat}\,U(t-nT)$$

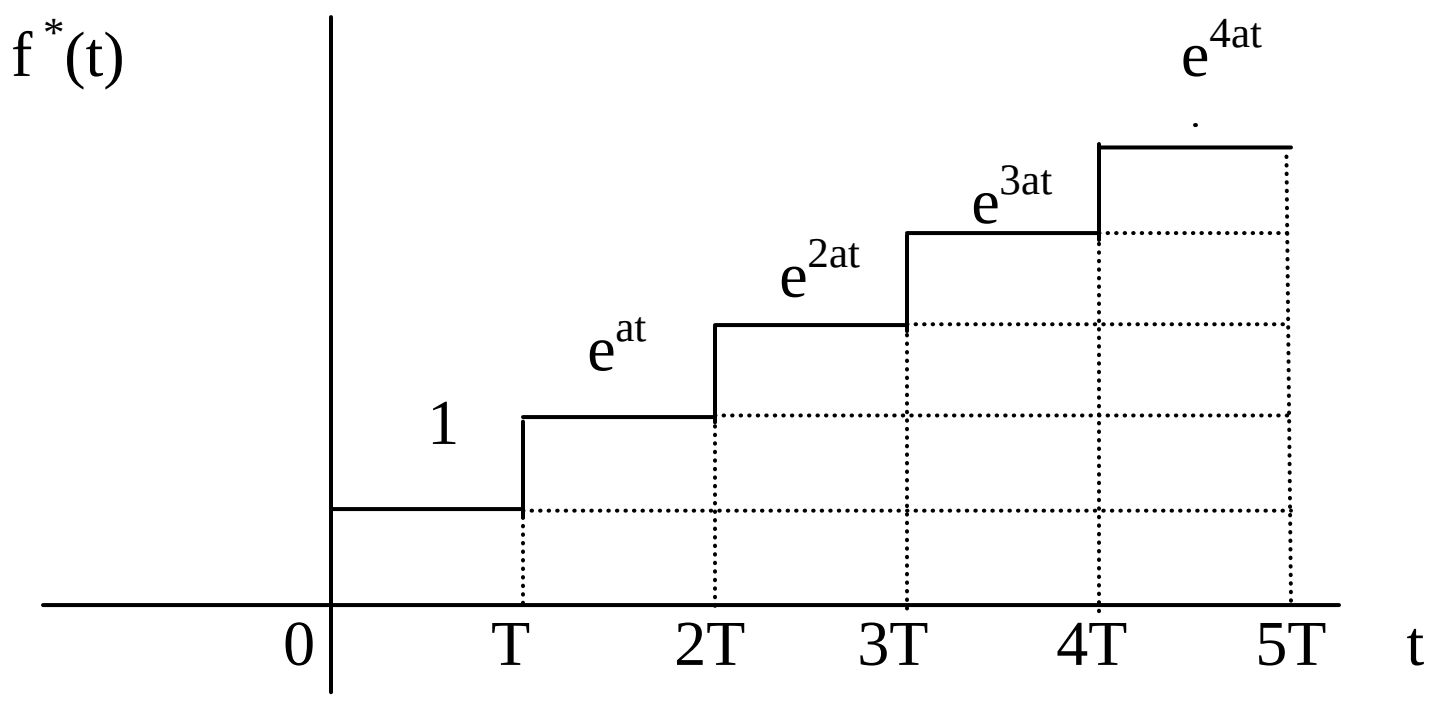

Eq 5.12-21 and Eq 5.12-45 convert the $K_{\Delta t}$ Transform and Z Transform of a function, f(t), into the Laplace Transform of its sampled sample and hold waveform, $L[f^*(t)]$.

Eq 5.12-21 and Eq 5.12-45 are rewritten below.

$$L[f^*(t)] = \frac{1}{s}\, sK_{\Delta t}[f(t)]\,\Big|_{s = \frac{e^{s\Delta t} - 1}{\Delta t}} \;,\quad K_{\Delta t} \text{ Transforms are a function of s} \tag{5.12-57}$$

$$L[f^*(t)] = \frac{1}{s}\,\frac{z - 1}{z}\, Z[f(t)]\Big|_{z = e^{sT}} \;,\quad \text{Z Transforms are a function of z} \tag{5.12-58}$$

For example

Let

$$f(t) = U(t) \tag{5.12-59}$$

$$L[f^*(t)] = \frac{1}{s}\, sK_{\Delta t}[f(t)]\,\Big|_{s = \frac{e^{s\Delta t} - 1}{\Delta t}} = \frac{1}{s}\, s(\frac{1}{s})\Big|_{s = \frac{e^{s\Delta t} - 1}{\Delta t}} = \frac{1}{s} \tag{5.12-60}$$

$$f^*(t) = U(t) \qquad \text{Good check} \tag{5.12-61}$$

$$L[f^*(t)] = \frac{1}{s}\,\frac{z - 1}{z}\, Z[f(t)]\Big|_{z = e^{sT}} = \frac{1}{s}\,\frac{z - 1}{z}\,\frac{z}{z - 1}\Big|_{z = e^{sT}} = \frac{1}{s} \tag{5.12-62}$$

$$f^*(t) = U(t) \qquad \text{Good check} \tag{5.12-63}$$

Note that Eq 5.12-61 and Eq 5.12-63 are the same.

For another example see Example 5.12-1 and Example 5.12-4 where f(t) = t and T = Δt

Let

$$f(t) = t \qquad (5.12\text{-}64)$$

$$L[f^*(t)] = \frac{1}{s} sK_{\Delta t}[f(t)]\Big|_{s = \frac{e^{s\Delta t}-1}{\Delta t}} = \frac{1}{s} s \frac{1}{s^2}\Big|_{s = \frac{e^{s\Delta t}-1}{\Delta t}} = \frac{1}{s}\frac{1}{s}\Big|_{s = \frac{e^{s\Delta t}-1}{\Delta t}} = \frac{\Delta t}{s(e^{s\Delta t}-1)} \qquad (5.12\text{-}65)$$

Expanding Eq 5.12-65

$$L[f^*(t)] = e^{-s\Delta t}\frac{\Delta t}{s} + e^{-2s\Delta t}\frac{\Delta t}{s} + e^{-3s\Delta t}\frac{\Delta t}{s} + e^{-4s\Delta t}\frac{\Delta t}{s} + e^{-5s\Delta t}\frac{\Delta t}{s} + \ldots \qquad (5.12\text{-}66)$$

$$f^*(t) = {}_1\sum_{n=0}^{\infty} \Delta t\, U(t-[n+1]\Delta t) \qquad (5.12\text{-}67)$$

$$L[f^*(t)] = \frac{1}{s}\frac{z-1}{z} Z[f(t)]\Big|_{z = e^{sT}} = \frac{1}{s}\frac{z-1}{z}\frac{Tz}{(z-1)^2}\Big|_{z = e^{sT}} = \frac{1}{s}\frac{Tz}{z-1}\Big|_{z = e^{sT}} = \frac{T}{s(e^{sT}-1)} \qquad (5.12\text{-}68)$$

Expanding Eq 5.12-68

$$L[f^*(t)] = e^{-sT}\frac{T}{s} + e^{-2sT}\frac{T}{s} + e^{-3sT}\frac{T}{s} + e^{-4sT}\frac{T}{s} + e^{-5sT}\frac{T}{s} + \ldots \qquad (5.12\text{-}69)$$

$$f^*(t) = {}_1\sum_{n=0}^{\infty} T\, U(t-[n+1]T) \qquad (5.12\text{-}70)$$

Note that Eq 5.12-67 and Eq 5.12-70 are the same. T = Δt.  Good check

<u>Introducing sample and hold sampling into a continuous time system using a $K_{\Delta t}$ Transform Methodology</u>

It may be remembered that the $K_{\Delta t}$ Transform was originally developed to solve differential difference equations, however, it has been found that the $K_{\Delta t}$ Transform is also useful in the analysis of sampled-data systems where sample and hold sampling is used. A derivation and demonstration of this capability of the $K_{\Delta t}$ Transform follows.

Consider the following sample and hold switch.

Sample and Hold Switch

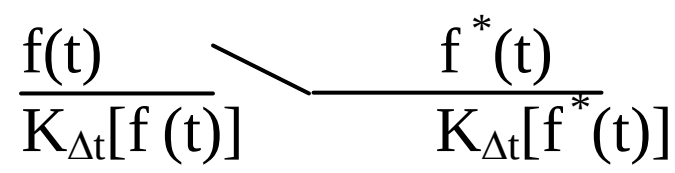

$K_{\Delta x}[f^*(t)] = K_{\Delta x}[f(t)]$ (5.12-71)

where

$f(t)$ = Switch input function

$f^*(t)$ = Switch output function, a sample and hold shaped waveform

$0 \le t < \infty$

$\Delta t$ = interval between samples

$f^*(t) = f(n\Delta t)$ for $n\Delta t \le t < [n+1]\Delta t$, $n = 0, 1, 2, 3, \ldots$

Sample and Hold Shaped Waveform

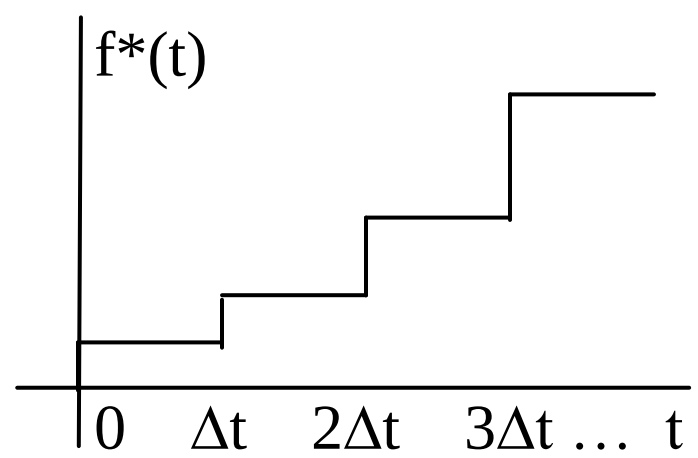

Derive, using $K_{\Delta t}$ Transforms, the sample and hold wavefrom, $f^*(t)$, that is obtained by the sample and hold sampling of $f(t)$.

$$f^*(t) = {}_1\sum_{n=0}^{\infty} f(n\Delta t)[U(t-n\Delta t) - U(t-[n+1]\Delta t)] \qquad (5.12\text{-}72)$$

where

$n = 0, 1, 2, 3, \ldots$

$f^*(t)$ = Series of unit amplitude pulses of width $\Delta t$ with an amplitude weighting of $f(n\Delta t)$

Taking the $K_{\Delta t}$ Transform of Eq 5.12-72

$$K_{\Delta t}[f^*(t)] = {}_1\sum_{n=0}^{\infty} f(n\Delta t)\left[\frac{(1+s\Delta t)^{-n}}{s} - \frac{(1+s\Delta t)^{-n-1}}{s}\right] = {}_1\sum_{n=0}^{\infty} f(n\Delta t)\,\frac{(1+s\Delta t)^{-n}}{s}\left[1 - \frac{1}{(1+s\Delta t)}\right] \qquad (5.12\text{-}73)$$

$$K_{\Delta t}[f^*(t)] = {}_1\sum_{n=0}^{\infty} f(n\Delta t)\,\frac{(1+s\Delta t)^{-n}}{s}\left[1 - \frac{1}{(1+s\Delta t)}\right] = {}_1\sum_{n=0}^{\infty} f(n\Delta t)\,\frac{(1+s\Delta t)^{-n}}{s}\left[\frac{s\Delta t}{(1+s\Delta t)}\right] \qquad (5.12\text{-}74)$$

$$K_{\Delta t}[f^*(t)] = \sum_{n=0}^{\infty} f(n\Delta t)\,(1+s\Delta t)^{-n-1}\Delta t \;\equiv\; {}_{\Delta t}\int_0^{\infty} f(t)(1+s\Delta t)^{-\left(\frac{t+\Delta t}{\Delta t}\right)}\Delta t = K_{\Delta t}[f(t)] \qquad (5.12\text{-}75)$$

where

$t = n\Delta t$

$n = 0, 1, 2, 3, \ldots$

$$K_{\Delta t}[f^*(t)] = K_{\Delta t}[f(t)] \qquad (5.12\text{-}76)$$

It should be noted that even though $f^*(t)$ is generally not equal to f(t) for all t > 0, $f^*(t) = f(t)$ at the sampling times, $t = 0, \Delta t, 2\Delta t, 3\Delta t, \ldots$ Since a $K_{\Delta t}$ Transform sample and hold switch calculation is based only on the sampling times, the input and output functions of a sample and hold switch have the same $K_{\Delta t}$ Transform.

Comment - $Z[f^*(t)] = Z[f(t)]$. Since the Z Transform switch calculation is based only on the sampling times, the input and output functions of a Z Transform sampler have the same Z Transform.

From Eq 5.12-75, The output of a $K_{\Delta t}$ Transform sample and hold switch is shown in Diagram 5.12-7 below.

Diagram 5.12-7 The output of a sample and hold switch ($K_{\Delta t}$ Transform system analysis)

Sampled f(t)

$$f^*(t) = K_{\Delta t}^{-1}\left[\sum_{n=0}^{\infty} f(n\Delta t)\,(1+s\Delta t)^{-n-1}\Delta t\right]$$

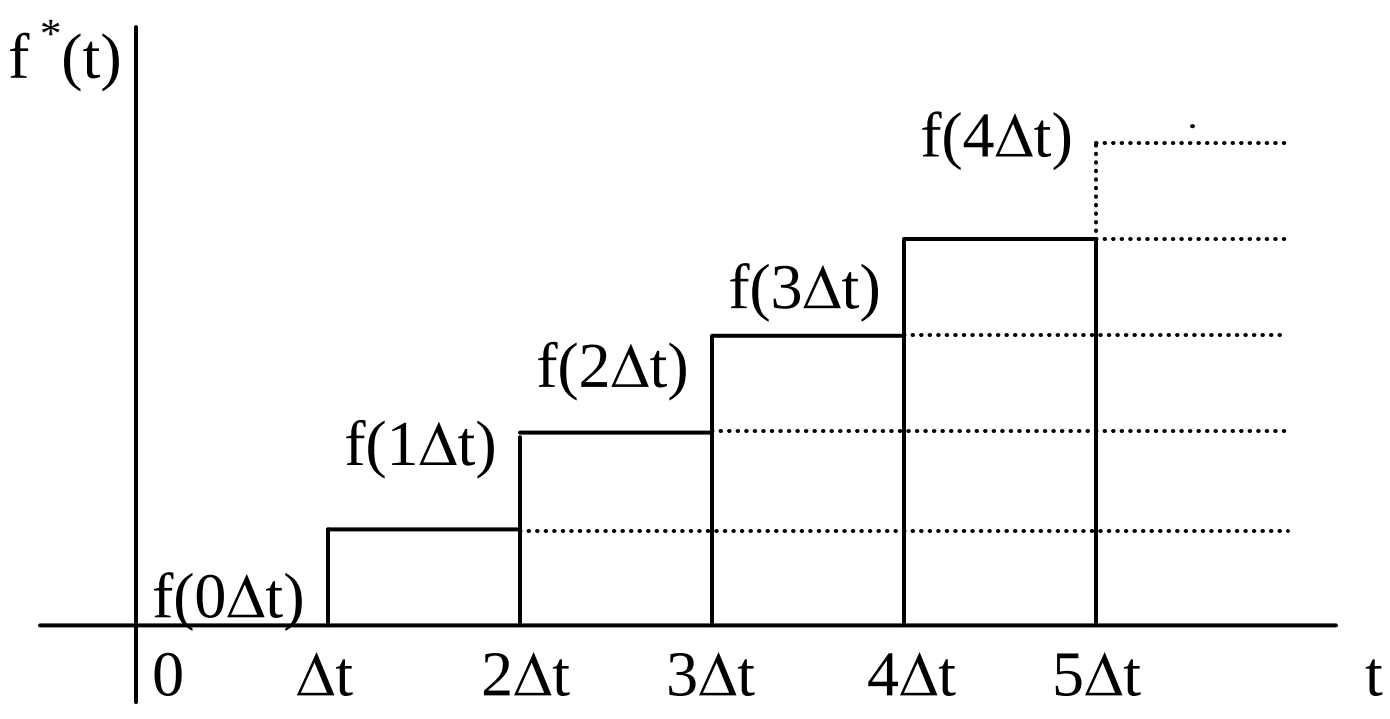

The term, $(1+s\Delta t)^{-n-1}\Delta t$ where $n = 0, 1, 2, 3, \ldots$ is the $K_{\Delta t}$ Transform of a unit amplitude pulse of width $\Delta t$ extending from $t = n\Delta t$ to $t = [n+1]\Delta t$. The $K_{\Delta t}$ Transform of the function, $f^*(t)$, is a series of consecutive unit amplitude pulses of $\Delta t$ width with their amplitudes weighted by the value, $f(n\Delta t)$. Thus, the resulting waveform represented is a sample and hold shaped waveform formed by the sampled values obtained at $t = 0, \Delta t, 2\Delta t, 3\Delta t, \ldots$

Note that the sample and hold switch diagrams of Diagram 5.12-7 and Diagram 5.12-2 are the same.

<u>Sampled f(t)</u>

$$f^*(t) = K_{\Delta t}^{-1}\left[\sum_{n=0}^{\infty} f(n\Delta t)\,(1+s\Delta t)^{-n-1}\,\Delta t\right] \qquad (5.12\text{-}77)$$

and

$$f^*(t) = L^{-1}[\,L[f^*(t)]] = L^{-1}\left[\frac{1}{\mathsf{s}}\,sK_{\Delta t}[f(t)]\Big|_{s\,=\,\frac{e^{\mathsf{s}\Delta t}-1}{\Delta t}}\right] \qquad (5.12\text{-}78)$$

<u>Derivation of several K$\Delta$t Transform, Z Transform, and Laplace Transform conversion equations</u>

Derive an equation to convert the $K_{\Delta t}$ Transform of a sample and hold shaped waveform function into its Laplace Transform equivalent.

Rewriting Eq 5.12-76

$$K_{\Delta t}[f^*(t)] = K_{\Delta t}[f(t)] \qquad (5.12\text{-}79)$$

Rewriting Eq 5.12-21

$$L[f^*(t)] = \frac{1}{\mathsf{s}}\,sK_{\Delta t}[f(t)]\,\Big|_{s\,=\,\frac{e^{\mathsf{s}\Delta t}-1}{\Delta t}}\ ,\quad K_{\Delta t}\text{ Transforms are a function of s} \qquad (5.12\text{-}80)$$

Substituting Eq 5.12-79 into Eq 5.12-80

$$L[f^*(t)] = \frac{1}{\mathsf{s}}\,sK_{\Delta t}[f^*(t)]\,\Big|_{s\,=\,\frac{e^{\mathsf{s}\Delta t}-1}{\Delta t}}\ ,\quad K_{\Delta t}\text{ Transforms are a function of s} \qquad (5.12\text{-}81)$$

Then

**The equation to convert the $K_{\Delta t}$ Transform of a sample and hold shaped waveform function into its Laplace Transform equivalent is:**

$$\mathbf{L[f^*(t)] = \frac{1}{\mathsf{s}}\,sK_{\Delta t}[f^*(t)]\,\Big|_{s\,=\,\frac{e^{\mathsf{s}\Delta t}-1}{\Delta t}}\ ,\quad K_{\Delta t}\text{ Transforms are a function of s}} \qquad (5.12\text{-}82)$$

**where**

**$\Delta t$ = sampling interval**
**t = 0, $\Delta$t, 2$\Delta$t, 3$\Delta$t, …**
**$f^*(t)$ = Sample and hold shaped waveform function**
**s = $K_{\Delta t}$ Transform variable**
**$\mathsf{s}$ = Laplace Transform variable**

Derive an equation to convert the Laplace Transform of a sample and hold shaped waveform function into its $K_{\Delta t}$ Transform equivalent

From Eq 5.12-82

$$L[f^*(t)] = \frac{1}{s} sK_{\Delta t}[f^*(t)] \tag{5.12-83}$$

$$e^{s\Delta t} = 1+s\Delta t \tag{5.12-84}$$

From Eq 5.12-83 and Eq 5.12-84

$$K_{\Delta t}[f^*(t)] = \frac{1}{s}\,[sL[f^*(t)]]\,|_{e^{s\Delta t} = 1+s\Delta t} \tag{5.12-85}$$

Then

**The equation to convert the Laplace Transform of a sample and hold shaped waveform function into its $K_{\Delta t}$ Transform equivalent is:**

$$\mathbf{K_{\Delta t}[f^*(t)] = \frac{1}{s}\,[sL[f^*(t)]]\,|_{e^{s\Delta t} = 1+s\Delta t}} \tag{5.12-86}$$

**where**

**$\Delta t$ = sampling interval**

**$t = 0, \Delta t, 2\Delta t, 3\Delta t, \ldots$**

**$f^*(t)$ = Sample and hold shaped waveform function**

**$s = K_{\Delta t}$ Transform variable**

**$s = \frac{1}{\Delta t}\ln(1+s\Delta t)$ , Laplace Transform variable**

Derive an equation to convert the Z Transform of a sample and hold shaped waveform function into its Laplace Transform equivalent.

Rewriting Eq 5.12-76

$$K_{\Delta t}[f^*(t)] = K_{\Delta t}[f(t)] \tag{5.12-87 }$$

Using the $K_{\Delta t}$ Transform to Z Transform conversion equation, Eq 1.6-52, on Eq 5.12-87

$$Z[f^*(t)] = Z[f(t)] \tag{5.12-88}$$

Rewriting Eq 5.12-45

$$L[f^*(t)] = \frac{1}{s}\,\frac{z-1}{z}\,Z[f(t)]|_{z = e^{sT}} \tag{5.12-89}$$

Substituting Eq 5.12-88 into Eq 5.12-89

$$L[f^*(t)] = \frac{1}{s}\,\frac{z-1}{z}\,Z[f^*(t)]\Big|_{z = e^{sT}} \tag{5.12-90}$$

Then

**The equation to convert the Z Transform of a sample and hold shaped waveform function into its Laplace Transform equivalent is:**

$$\mathbf{L[f^*(t)] = \frac{1}{s}\,\frac{z-1}{z}\,Z[f^*(t)]\Big|_{z = e^{sT}}} \tag{5.12-91}$$

**where**

**T = Δt = sampling interval**

**t = 0, T, 2T, 3T, …**

**$f^*(t)$ = Sample and hold shaped waveform function**

**$z = e^{sT}$ = Z Transform variable**

**$s$ = Laplace Transform variable**

Derive an equation to convert the Laplace Transform of a sample and hold shaped waveform function into its Z Transform equivalent.

From Eq 5.12-91

$$L[f^*(t)] = \frac{1}{s}\,\frac{z-1}{z}\,Z[f^*(t)]\Big|_{z = e^{sT}} \tag{5.12-92}$$

From Eq 5.12-92

$$Z[f^*(t)] = \frac{z}{z-1}\,[sL[f^*(t)]]\,)]\Big|_{e^{sT} = z} \tag{5.12-93}$$

Then

**The equation to convert the Laplace Transform of a sample and hold shaped waveform function into its Z Transform equivalent is:**

$$\mathbf{Z[f^*(t)] = \frac{z}{z-1}\,[sL[f^*(t)]]\,)]\Big|_{e^{sT} = z}} \tag{5.12-94}$$

**where**

**T = Δt = sampling interval**

**t = 0, T, 2T, 3T, …**

**$f^*(t)$ = Sample and hold shaped waveform function**

**$z = e^{sT}$ = Z Transform variable**

**$s$ = Laplace Transform variable**

Derivation of the very important KΔt Transform equation, $K_{\Delta t}[c(t)] = K_{\Delta t}[g(t)]\, K_{\Delta t}[f(t)]$

Below is a listing of the variables used in the derivations that follow.

Listing of derivation variables

- f(t) = system input function
- $F(s) = K_{\Delta t}[f(t)]$ , The $K_{\Delta t}$ Transform of the system input function, f(t)
- G(s) = System Laplace Transform transfer function
- $g(t) = L^{-1}[G(s)]$ , Inverse Laplace Transform of the system Laplace transfer function, G(s)
- c(t) = System output function
- C(s) = L[c(t)] , Laplace Transform of the output function, c(t)
- $f^*(t)$ = Sampled system input function
- $f^*(t)$ = f(t) at t = 0, Δt, 2Δt, 3Δt, … , a sample and hold waveform
- $F^*(s) = K_{\Delta t}[f^*(t)]$ , The $K_{\Delta t}$ Transform of $f^*(t)$
- $G^*(s)$ = System $K_{\Delta t}$ Transform transfer function
- $g^*(t) = K_{\Delta t}^{-1}[G^*(s)]$ , Inverse Laplace Transform of the system $K_{\Delta t}$ Transform transfer function, $G^*(s)$
- $c^*(t)$ = Sampled system output function
- $c^*(t)$ = c(t) at t = 0, Δt, 2Δt, 3Δt, … , a sample and hold waveform
- $C^*(s) = K_{\Delta t}[c^*(t)]$ , $K_{\Delta t}$ Transform of the sampled output function, $c^*(t)$

Consider a basic continuous time system that is specified in Diagram 5.12-8 below.

Diagram 5.12-8  Initially passive system transfer function diagram

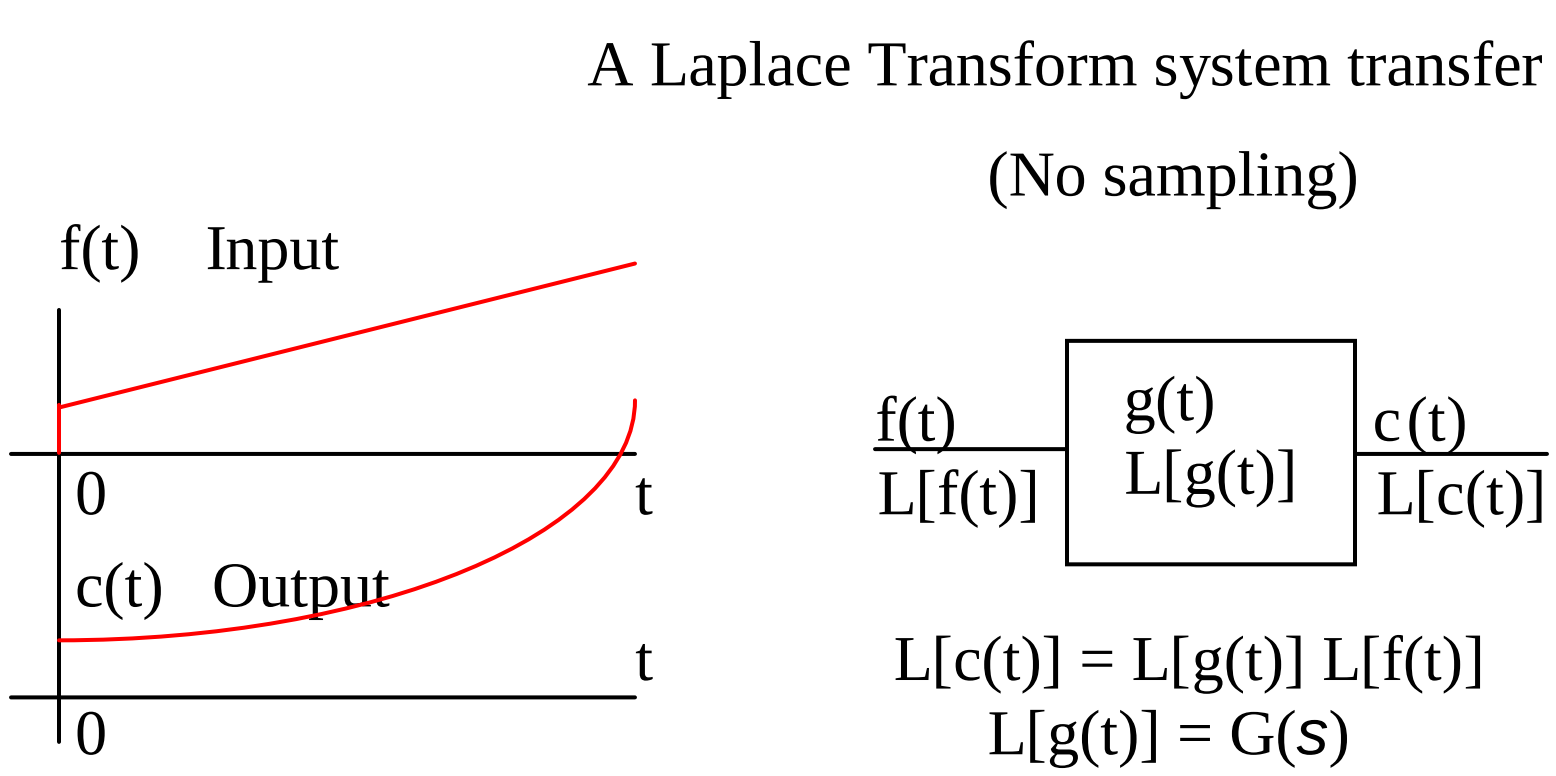

The input to the system shown in Diagram 5.12-8 is f(t) and its output is c(t). The system is represented by a Laplace Transform transfer function, L[g(t)] = G(s). It is desired that this same system respond to a sample and hold shaped input waveform, $f^*(t)$, obtained from the sample and hold sampling of the input function, f(t), where $f^*(t)$ = f(t) at the sampling instants. It is also desired that the sample and hold sampling of the system output, c(t), be a sample and hold shaped waveform, $c^*(t)$, where $c^*(t)$ = c(t) at the sampling instants. The sampling switches are synchronized and the sampling interval is Δt. Refer to the following sampled-data system diagram, Diagram 5.12-9.

Diagram 5.12-9  The system of Diagram 5.12-8 with the sample and hold sampling of the system input and output introduced

Sample and Hold Laplace Transform sampled system transfer function

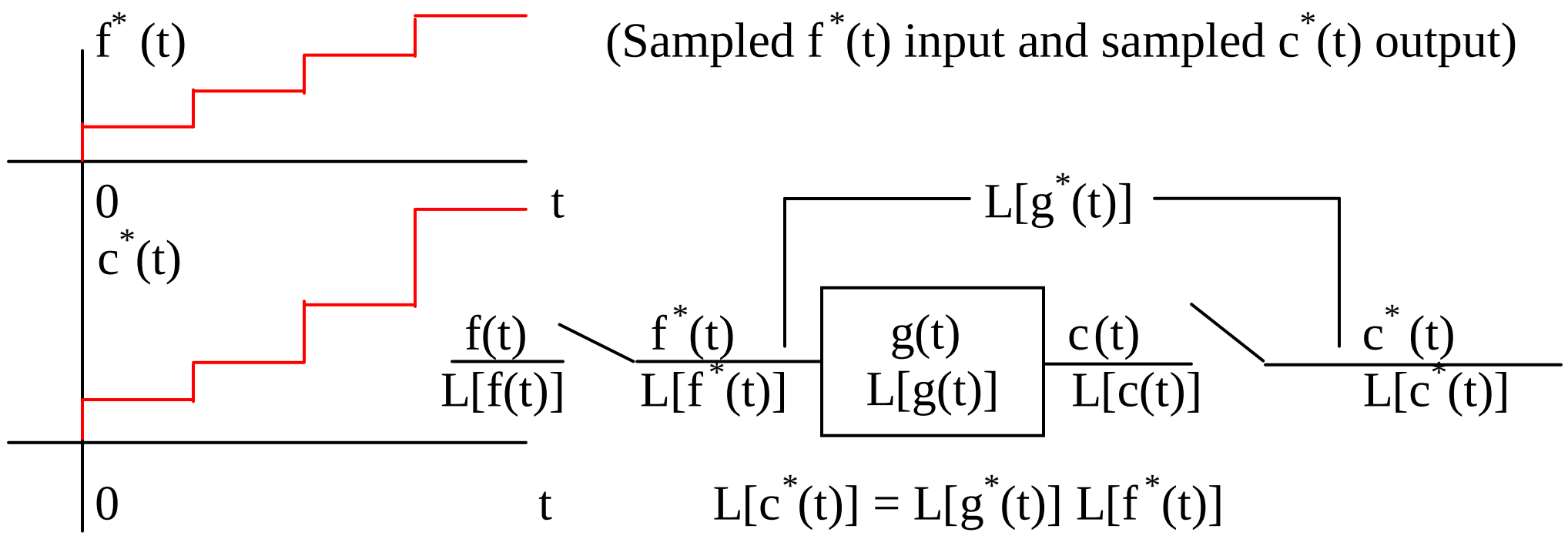

$$L[c^*(t)] = L[g^*(t)]\, L[f^*(t)]$$
$$C(s) = G(s)\, F(s)$$

The sample and hold switches are synchronized

Though the analysis of systems with sampled inputs and outputs can be performed using Laplace Transforms, perhaps using one of the equations derived above, the use of Z Transforms is most often used. Z Transforms have some important advantages. There are many technical books describing in depth the derivation of the Z Transform and its use. An excellent textbook for learning about Z Transforms is “Digital & Sampled-data Control Systems” written by Julius T. Tou and published by McGraw-Hill. The Z Transform sampled system transfer function for the system shown in Diagram 5.12-8 and Eq 5.12-9 is shown below in Diagram 5.12-10

Diagram 5.12-10  Z Transform system diagram of the system shown in Diagram 5.12-8 and Diagram 5.12-8

Z Transform sampled system transfer function

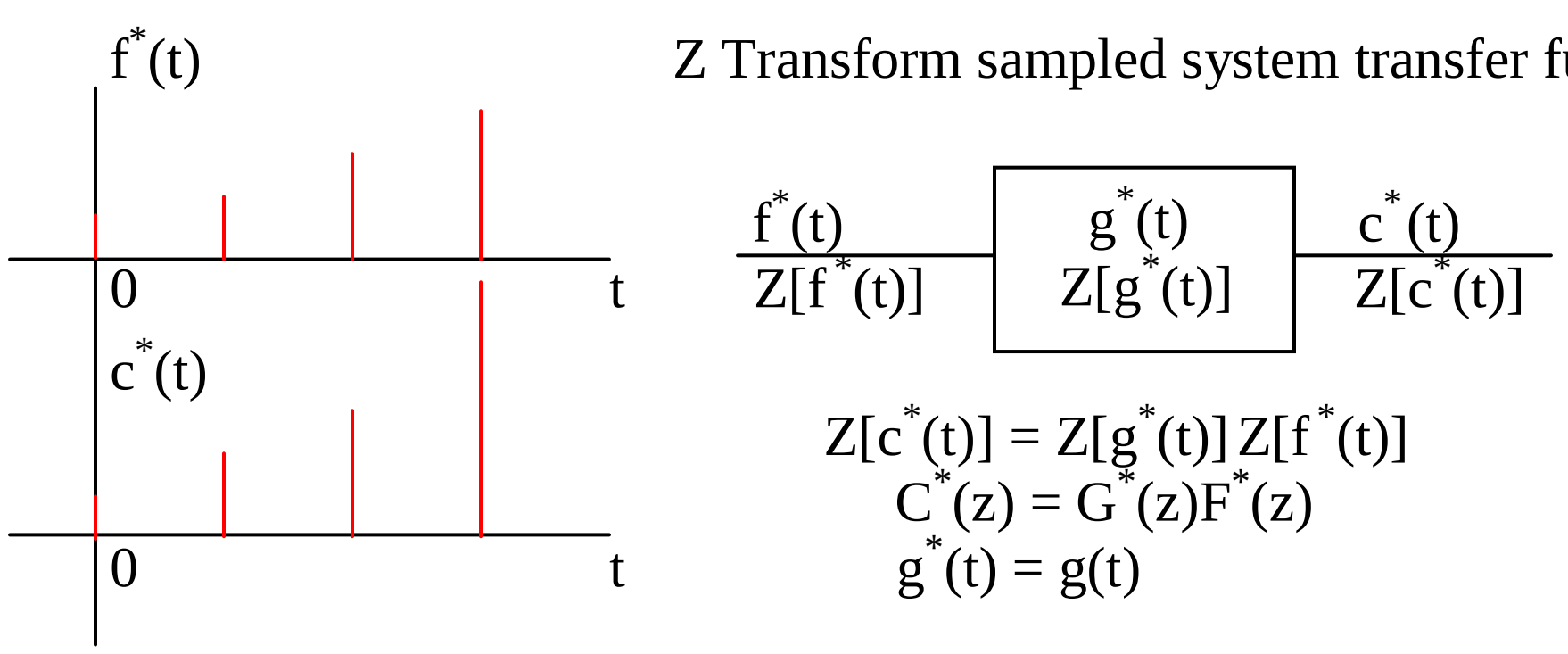

$$Z[c^*(t)] = Z[g^*(t)]\, Z[f^*(t)]$$
$$C^*(z) = G^*(z) F^*(z)$$
$$g^*(t) = g(t)$$

It is necessary to derive the $K_{\Delta t}$ Transform sample and hold sampled system transfer function equation (similar to the Z Transform equation, $Z[c^*(t)] = Z[g^*(t)]\, Z[f^*(t)]$, for the system shown in Diagram 5.12-8 and Diagram 5.12-9. See Diagram 5.12-11 below.

Diagram 5.12-11  The sampling of the input and output of the system described in Diagram 5.12-8

Sample and Hold $K_{\Delta t}$ Transform sampled system transfer function

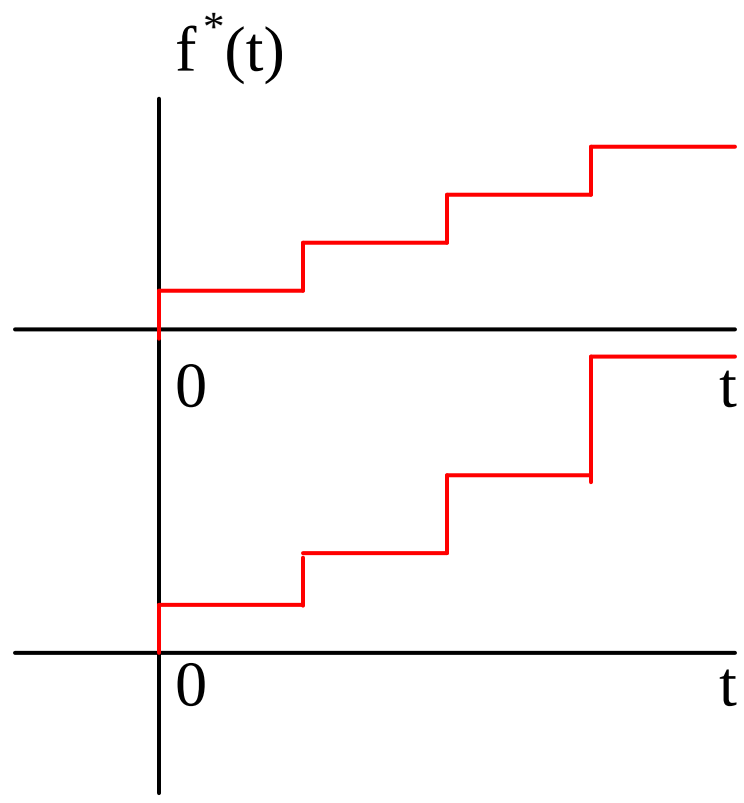

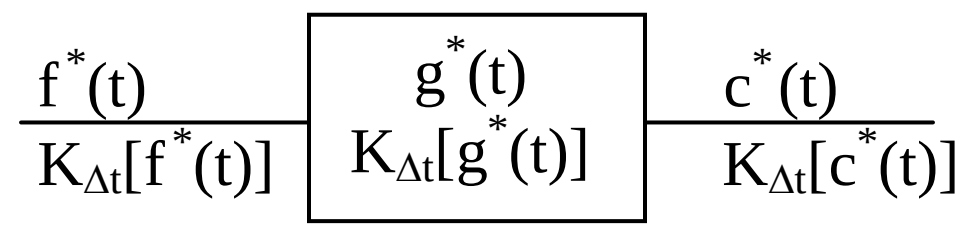

$$K_{\Delta t}[c^*(t)] = K_{\Delta t}[g^*(t)]\,K_{\Delta t}[f^*(t)]$$
$$C^*(s) = G^*(s)F^*(s)$$
$$g^*(t) = ?$$

It is necessary to determine the $K_{\Delta t}$ Transform sample and hold sampled-data system transfer function, $g^*(t)$.

From Eq 5.10-86

Z Transform Convolution Integral

$$Z[\frac{1}{T}\,{}_T\!\int_0^t f^*(\lambda)g^*(t-\lambda-\Delta\lambda)\Delta\lambda] = \frac{Z[f^*(t)]Z[g^*(t)]}{z} = \frac{F^*(z)G^*(z)}{z} \qquad (5.12\text{-}95)$$

where
$f^*(t), g^*(t)$ = discrete functions of t
$f^*(t-\lambda-\Delta\lambda)$ = discrete function of t, $\lambda$, and $\Delta\lambda$
$\Delta\lambda = T = \Delta t$
$t = n\Delta t$ ,    $n = 0, 1, 2, 3, \ldots, \frac{t}{\Delta t}-1, \frac{t}{\Delta t}$
$\lambda = 0, \Delta\lambda, 2\Delta\lambda, 3\Delta\lambda, \ldots, t-\Delta\lambda, t$

$Z[h(t)]\,) = \frac{1}{T}\,{}_T\!\int_0^\infty z^{-(\frac{t}{\Delta t})}h(t)\Delta t$ ,    The Z Transform of a function of t
$F^*(z) = Z[f^*(t)]$ = Z Transform of the function $f^*(t)$
$G^*(z) = Z[g^*(t)]$ = Z Transform of the function $g^*(t)$
$Z = e^{sT}$ ,   Z Transform variable
$s$ = Laplace Transform variable

Multiplying both sides of Eq 5.12-95 by T

$$Z[{}_T\int_0^t f^*(\lambda)g^*(t-\lambda-\Delta\lambda)\Delta\lambda] = Z[f^*(t)]Tz^{-1}Z[g^*(t)] = F^*(z)[Tz^{-1}G^*(z)] = Z[c^*(t)] \quad (5.12\text{-}96)$$

$$c^*(t) = {}_T\int_0^t f^*(\lambda)g^*(t-\lambda-\Delta\lambda)\Delta\lambda \quad (5.12\text{-}97)$$

$$Z[c^*(t)] = Z[f^*(t)]Tz^{-1}Z[g^*(t)] = F^*(z)[Tz^{-1}G^*(z)] = C^*(z) \quad (5.12\text{-}98)$$

Changing the form of Eq 5.12-98 by multiplying both sides by $Tz^{-1}$ where $z = e^{sT}$, $z = 1+s\Delta t$, $s = \frac{e^{s\Delta t}-1}{\Delta t}$, and $T = \Delta t$

$$Tz^{-1}Z[c^*(t)]\,|_{z = 1+s\Delta t} = [Tz^{-1}Z[f^*(t)]][Tz^{-1}Z[g^*(t)]]\,|_{z = e^{sT}} \quad (5.12\text{-}99)$$

$$Tz^{-1}Z[c^*(t)]\,|_{z = 1+s\Delta t} = [Tz^{-1}Z[f^*(t)]\,|_{z = 1+s\Delta t}][Tz^{-1}Z[g^*(t)]\,|_{z = e^{sT}}] \quad (5.12\text{-}100)$$

For the Z Transform, the function, $g^*(t)$, has been shown to be equal to g(t).

$$g^*(t) = g(t) \quad (5.12\text{-}101)$$

Substituting Eq 5.12-101 into Eq 5.12-100

$$Tz^{-1}Z[c^*(t)]\,|_{z = e^{sT}} = [Tz^{-1}Z[f^*(t)]\,|_{z = e^{sT}}][Tz^{-1}Z[g(t)]\,|_{z = e^{sT}}] \quad (5.12\text{-}102)$$

Writing the Z Transform to $K_{\Delta t}$ Transform conversion equation

$$K_{\Delta t}[p^*(t)] = Tz^{-1}Z[p^*(t)]|_{z = 1+s\Delta t} \quad (5.12\text{-}103)$$

$\Delta t = T$ sampling period

$t = n\Delta t$

$n = 0, 1, 2, 3, \dots$

$K_{\Delta x}[p^*(t)] = K_{\Delta x}$ Transform of $p^*(t)$

$Z[p^*(t)] =$ Z Transform of $p^*(t)$

Any function with the same values at $t = 0, \Delta t, 2\Delta t, 3\Delta t, \dots$ will have the same Z Transform. For example, $p^*(t)$ and p(t) would have the same Z Transform.

Applying the conversion relationship of Eq 5.12-8103 to Eq 5.12-102

$$K_{\Delta t}[c^*(t)] = K_{\Delta t}[f^*(t)]\,K_{\Delta t}[g(t)] \quad (5.12\text{-}104)$$

or in a different form

$$C^*(s) = F^*(s)\,G^*(s) \quad (5.12\text{-}105)$$

$T = \Delta t$ (5.12-106)

From Eq 5.12-102 thru Eq 5.12-106

$$G^*(s) = K_{\Delta t}[g(t)] = Tz^{-1}Z[g(t)]\,\Big|_{z\,=\,1+sT} = Tz^{-1}G^*(z)\Big|_{z\,=\,1+sT} \quad (5.12\text{-}107)$$

$$G^*(z) = Z[g(t)] = Z[g^*(t)] \quad (5.12\text{-}108)$$

From Diagram 5.12-8

$$g(t) = L^{-1}[G(s)] \quad (5.12\text{-}109)$$

Substituting Eq 5.12-109 into Eq 5.12-107

$$G^*(s) = K_{\Delta t}[L^{-1}[G(s)]] \quad (5.12\text{-}110)$$

To find $G^*(s)$ from $G(s)$ there is a table in the Appendix that may help. See TABLE 3a, The Conversion of Calculus Function Laplace Transforms to Equivalent Function $K_{\Delta t}$ Transforms.

With $G^*(s)$ defined as shown in Eq 5.12-110, Eq 5.12-105 defines the system of Diagram 5.12-8 where the input and output are sampled using synchronized sample and hold switches. See Diagram 5.12-9.

The $K_{\Delta t}$ Transform equations, Eq 5.12-105 and Eq 5.12-110, are to sample and hold sampled systems as the Z Transform equations, $C^*(z) = G^*(z)F^*(z)$ and $G^*(z) = Z[L^{-1}[G(s)]]$ are to impulse sampled systems. The definition of g(t) comes from the system transfer function shown in Diagram 5.12-8 where $g(t) = L^{-1}[G(s)]$.

Since $K_{\Delta t}$ Transforms become Laplace Transforms for $\Delta t \rightarrow 0$, in the limit, Eq 5.12-104 becomes $L[c(t)] = L[f(t)]\,L[g(t)]$ , the equation for the unsampled system.

From Eq 5.12-98, Eq 5.12-101, Eq 5.12-104 and Eq 5.12-105

**Equivalent $K_{\Delta t}$ Transform and Z Transform Equations**

$$\mathbf{Z[c^*(t)] = Z[f(t)]Tz^{-1}Z[g(t)]} \quad (5.12\text{-}111)$$

$$\mathbf{K_{\Delta t}[c^*(t)] = K_{\Delta t}[f^*(t)]\,K_{\Delta t}[g(t)]} \quad (5.12\text{-}112)$$

**or**

$$\mathbf{C^*(z) = F^*(z)Tz^{-1}\,G^*(z)} \quad (5.12\text{-}113)$$

$$\mathbf{C^*(s) = F^*(s)\,G^*(s)} \quad (5.12\text{-}114)$$

**Note the equivalency of the above two transform equations, Eq 5.12-111 and Eq 5.12-112. Taking the appropriate inverse transforms, for a given $f^*(t)$ and g(t), the same $c^*(t)$ is obtained.**

Then

From Eq 5.12-111 thru Eq 5.12-114 the following sample and hold transfer functions are equivalent.

**Equivalent Sample and Hold $K_{\Delta t}$ Transform and Z Transform Transfer Functions #1**

**$K_{\Delta t}$ Transform Transfer Function Diagram** **Z Transform Transfer Function Diagram**

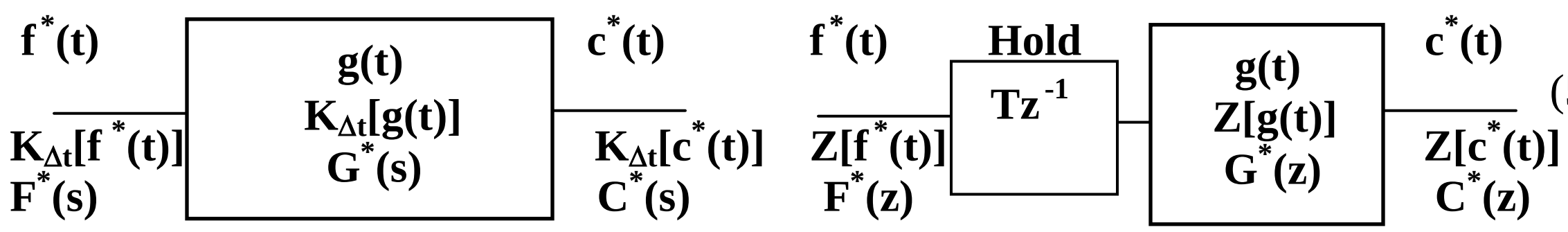

(5.12-115)

$$G^*(s) = K_{\Delta t}[L^{-1}[G(s)]] \qquad G^*(z) = Z[L^{-1}[G(s)]]$$

**For the same input to both transfer functions the same output will be obtained.**

To better understand the Z Transform $Tz^{-1}$ hold term, consider the following comparison between the equivalent $K_{\Delta t}$ Transform Transfer Function and Z Transform Transfer Function diagrams shown above. For the $K_{\Delta t}$ Transform system let $f^*(t)$ be a unit amplitutde pulse of width $T = \Delta t$ and the transfer function be an integrator. For the Z Transform system let $f^*(t)$ be a unit area impulse and its transfer function also be an integrator. Compare the outputs of both systems.

$K_{\Delta t}$ Transform Transfer Function Diagram Z Transform Transfer Function Diagram

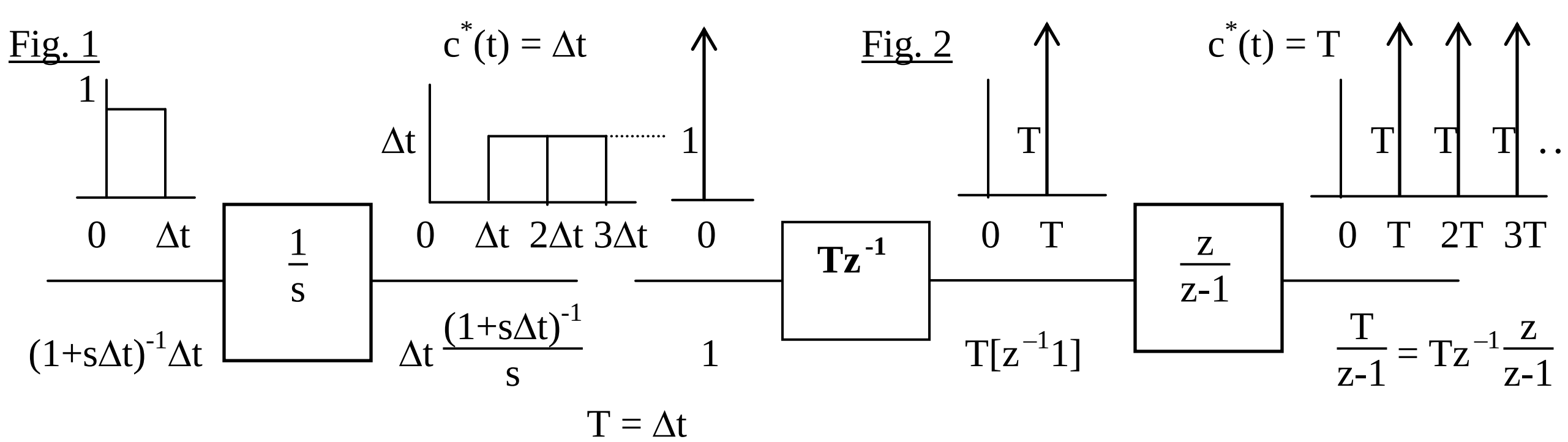

The first thing that should be observed from the above diagrams is that there is a significant difference between the $K_{\Delta t}$ Transform and the Z Transform. The $K_{\Delta t}$ Transform is based on weighted unit amplitude pulses of width $\Delta t$ and the Z Transform is based on weighted unit area impulses. The second thing that should be observed from the above diagrams is that the outputs of both systems represent the same function, $c^*(t) = T = \Delta t$, even though the input of the $K_{\Delta t}$ Transform system is a unit amplitude pulse of $\Delta t$ width and the input of the Z Transform system is a unit area impulse. It is evident that the Z Transform term, $Tz^{-1}$, introduces the sample and hold operation into the Z Transform system since the output results are the same.

Referring to Fig. 2 above, it is seen that the Z Transform transfer function, $Tz^{-1}$, integrates an input impulse over one time interval, T, to yield an output impulse value of T with a delay of T. The Z Transform representing a unit amplitude pulse of width T is $[Tz^{-1}]1$. Note the check of this conclusion below.

Show that the Z Transform of a unit amplitude pulse of duration T is $[Tz^{-1}]1$.

From the definition of the $K_{\Delta t}$ Transform

$$s = \frac{e^{\mathit{s}\Delta t}-1}{\Delta t} \tag{5.12-116}$$

where

$\Delta t$ = interval between consecutive values of t

s = $K_{\Delta t}$ Transform variable

$\mathit{s}$ = Laplace Transform variable

$$z = e^{\mathit{s}\Delta t} \tag{5.12-117}$$

From Eq 5.12-116 and Eq 5.12-117

$$z = 1+s\Delta t \tag{5.12-118}$$

$$T = \Delta t \tag{5.12-119}$$

Substituting Eq 5.118 and Eq 5.12-119 into $[Tz^{-1}]1$

$$[Tz^{-1}]1 = (1+s\Delta t)^{-1}\Delta t \tag{5.12-120}$$

$(1+s\Delta t)^{-1}\Delta t$ is recognized as being the $K_{\Delta t}$ Transform of a unit amplitude pulse of width T initiated at t = 0.

Thus

**$Tz^{-1}$ = Z Transform of a unit amplitude pulse of width T initiated at t = 0.** (5.12-121)

From the above derivations, sample and hold sampling of the input and output of the system in Diagram 5.12-8 is shown in the following diagram, Diagram 5.12-12.

## Diagram 5.12-12 K∆t Transform sample and hold sampling of the input and output to a continuous time system

$K_{\Delta t}[c^*(t)] = K_{\Delta t}[g(t)]\,K_{\Delta t}[f^*(t)]$ or $C^*(s) = G^*(s)F^*(s)$

Fig. 1 represents the sample and hold sampled system of Fig. 2. The inputs of both Fig. 1 and Fig. 2 are the same for t = 0, ∆t, 2∆t, 3∆t, …

$K_{\Delta t}$ Transform system transfer function

**Fig. 1**

Sample and hold sampling

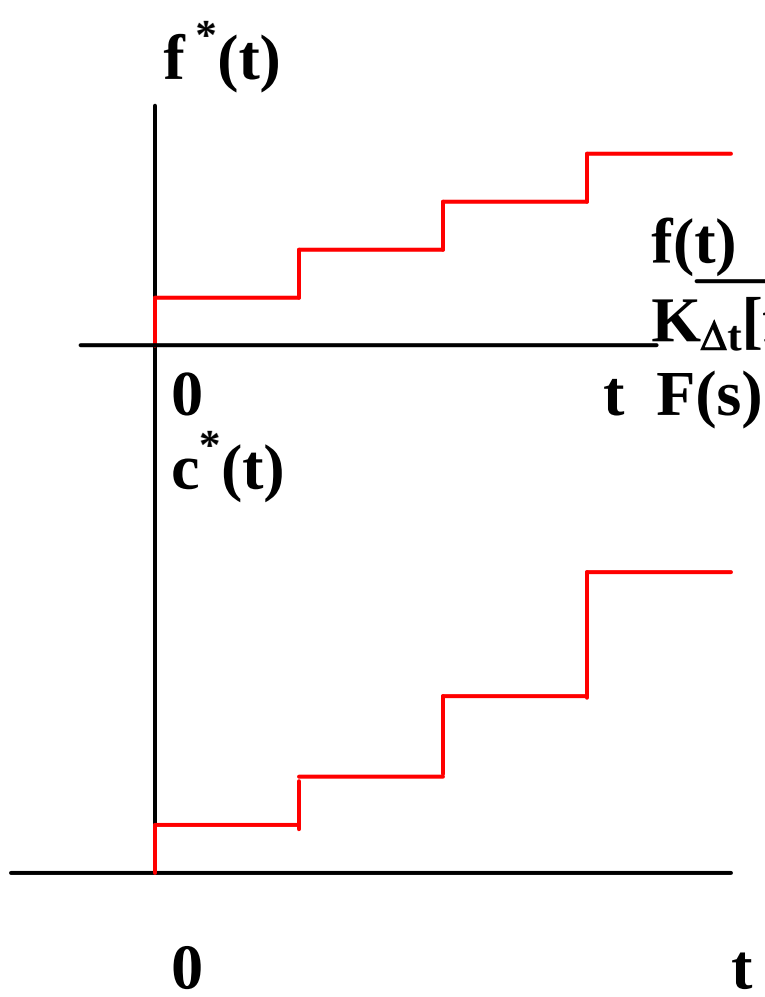

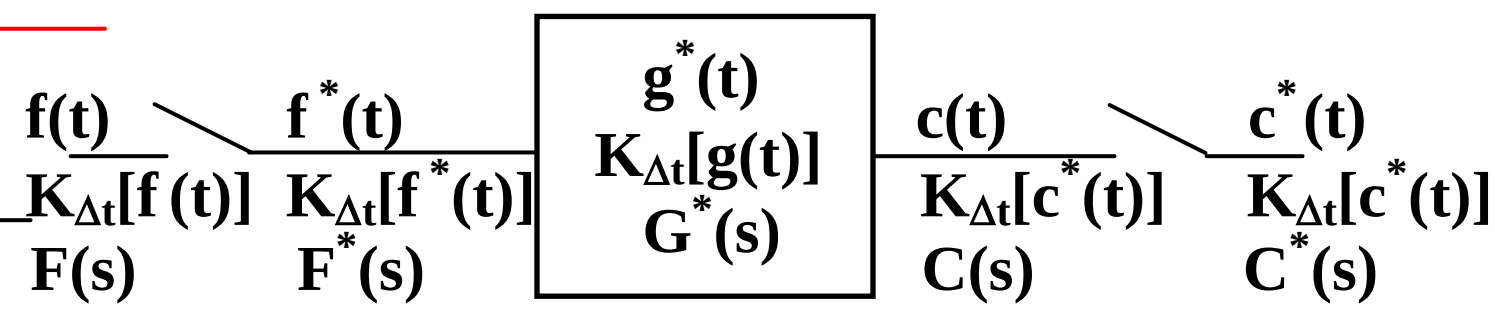

$K_{\Delta t}[c^*(t)] = K_{\Delta t}[g^*(t)]\,K_{\Delta t}[f^*(t)]$

$C^*(s) = G^*(s)F^*(s)$

$G^*(s) = K_{\Delta t}[L^{-1}[G(s)]]$

$g^*(t) = g(t) = L^{-1}[G(s)]$

$G^*(s) = K_{\Delta t}[g(t)] = Tz^{-1}Z[g(t)]\,\big|_{z\,=\,1+s\Delta t}$

$G^*(s) = Tz^{-1}G^*(z)\,\big|_{z\,=\,1+s\Delta t}$

$K_{\Delta t}[f^*(t)] = K_{\Delta t}[f(t)]$

The two switches shown in Fig. 1 above are synchronous sample and hold switches.

Laplace Transform system transfer function

**Fig. 2**

No sampling

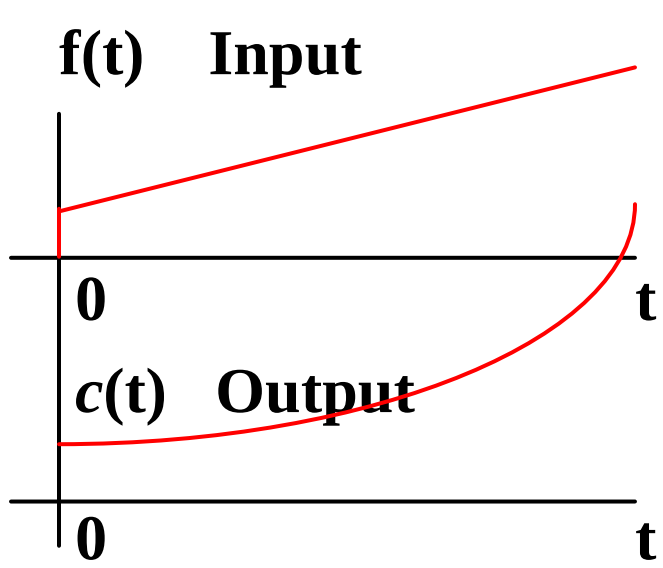

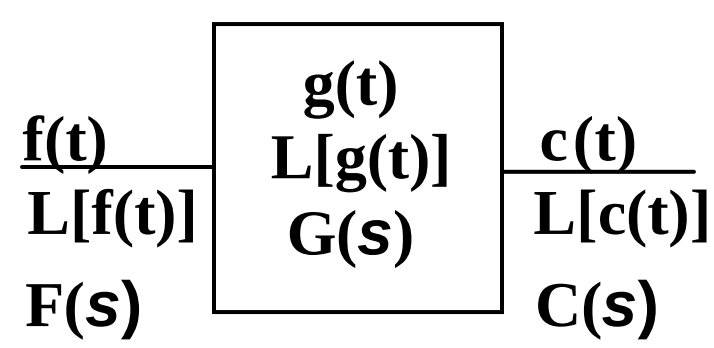

$L[c(t)] = L[g(t)]\,L[f(t)]$

$C(s) = G(s)F(s)$

$G(s) = L[g(t)]$

$g(t) = L^{-1}[G(s)]$

where

- $f(t)$ = system unsampled input function
- $f^*(t)$ = system sample and hold sampled input function
- $F^*(s) = K_{\Delta t}[f^*(t)]$, The $K_{\Delta t}$ Transform of the sampled system input function, $f^*(t)$
- $g(t) = L^{-1}[G(s)]$, Inverse Laplace Transform of the system Laplace transfer function, $G(s)$
- $c(t)$ = System unsampled output function
- $c^*(t)$ = System sample and hold sampled output function
- $G^*(s) = K_{\Delta t}[g(t)]$ = System $K_{\Delta t}$ Transform transfer function

**$g^*(t) = K_{\Delta t}^{-1}[G^*(s)]$ , Inverse Laplace Transform of the system $K_{\Delta t}$ Transform transfer function, $G^*(s)$**
**$C^*(s) = K_{\Delta t}[c^*(t)]$ , $K_{\Delta t}$ Transform of the sampled output function, $c^*(t)$**
**$g^*(t) = g(t) = L^{-1}[G(\mathsf{s})]$**
**$G^*(s) = K_{\Delta t}[L^{-1}[G(\mathsf{s})]]$**
**$F(\mathsf{s}) = L[f(t)]$ , The Laplace Transform of the system input function, $f(t)$**
**$C(\mathsf{s}) = L[c(t)]$ , The Laplace Transform of the system output function, $c(t)$**
**$G(\mathsf{s}) = L[g(t)]$ = System Laplace Transform transfer function**
**$G^*(z) = Z[g(t)]$ = Z Transform of the function $g(t)$**
**$z = e^{\mathsf{s}T} = 1+s\Delta t$**
**$s = K_{\Delta t}$ Transform variable**
**$\mathsf{s}$ = Laplace Transform variable**

**To find $G^*(s)$ from $G(\mathsf{s})$ there is a table in the Appendix that may help. See TABLE 3a, The Conversion of Calculus Function Laplace Transforms to Equivalent Function $K_{\Delta t}$ Transforms.**

For the sake of completeness, there is another set of equivalent sample and hold $K_{\Delta t}$ Transform and Z Transform transfer functions which should be shown. This second set of transfer functions is presented below.

**Equivalent Sample and Hold $K_{\Delta t}$ Transform and Z Transform Transfer Functions #2**

**$K_{\Delta t}$ Transform Transfer Function Diagram     Z Transform Transfer Function Diagram**

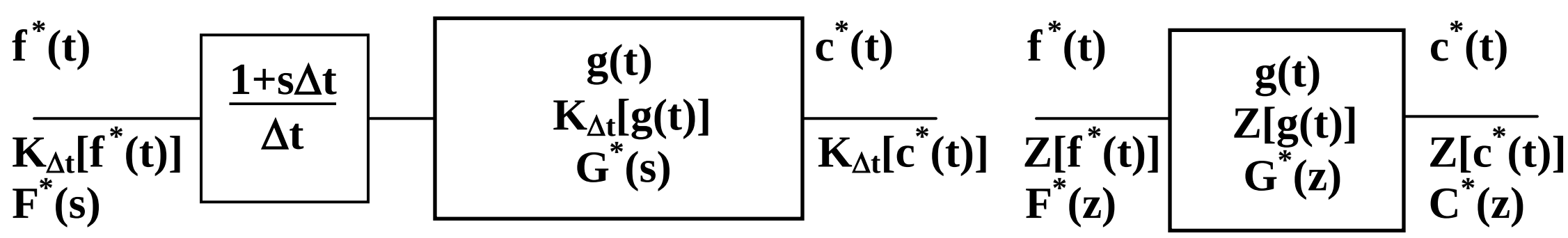

(5.12-122)

**$G^*(s) = K_{\Delta t}[L^{-1}[G(\mathsf{s})]]$          $G^*(z) = Z[L^{-1}[G(\mathsf{s})]]$**

**$\dfrac{1+s\Delta t}{\Delta t} = \dfrac{z}{T}\bigg|_{\substack{z = 1+s\Delta t \\ T = \Delta t}}$ , This term cancels the $K_{\Delta t}$ Transform sample and hold.**

**For the same input to both transfer functions the same output will be obtained.**

There is one other important thing that should be pointed out. Because a $K_{\Delta t}$ Transform system analysis deals with weighted unit amplitude pulses of $\Delta t$ width and a Z Transform system analysis deals with weighted unit area impulses, there is a difference. This difference will be seen in the results of an analysis of the same system where both the $K_{\Delta t}$ Transform methodology and Z Transform methodology are used without $Tz^{-1}$ or $\dfrac{1+s\Delta t}{\Delta t}$ compensation. Consider the following simple system which is analyzed in both ways.

$K_{\Delta t}$ Transform System Analysis (Sample and Hold Sampling)

Z Transform System Analysis (Impulse Sampling)

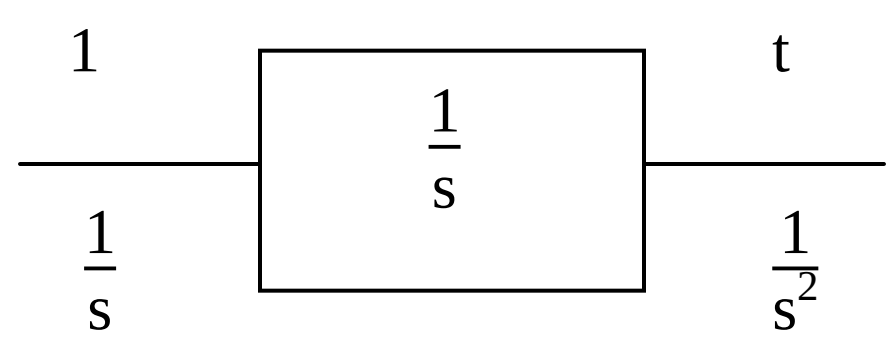

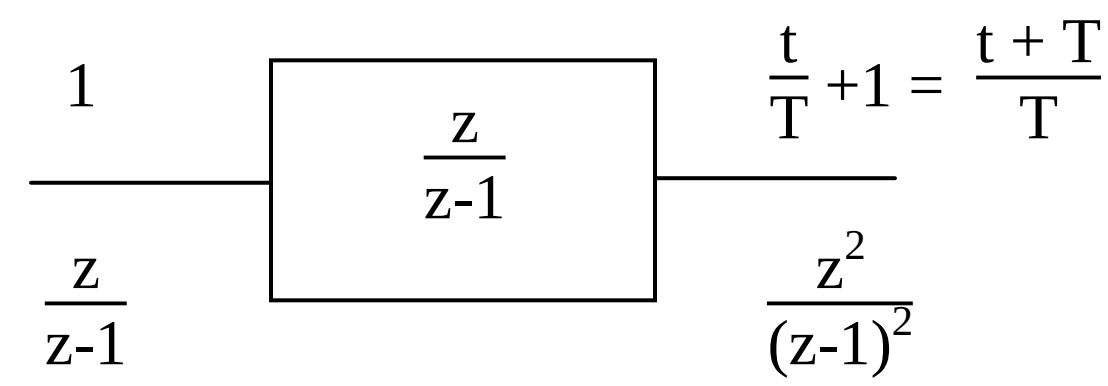

An integrator

$$\left(\frac{1}{s}\right)\left(\frac{1}{s}\right) = \frac{1}{s^2}$$

An integrator

$$\left(\frac{z}{z-1}\right)\left(\frac{z}{z-1}\right) = \left(\frac{z}{T}\right)\left(\frac{Tz}{(z-1)^2}\right)$$

$T = \Delta t$
$t = n\Delta t = nT$
$n = 0, 1, 2, 3, \ldots$

Note that for the same inputs the outputs differ.

Fig.1 System output

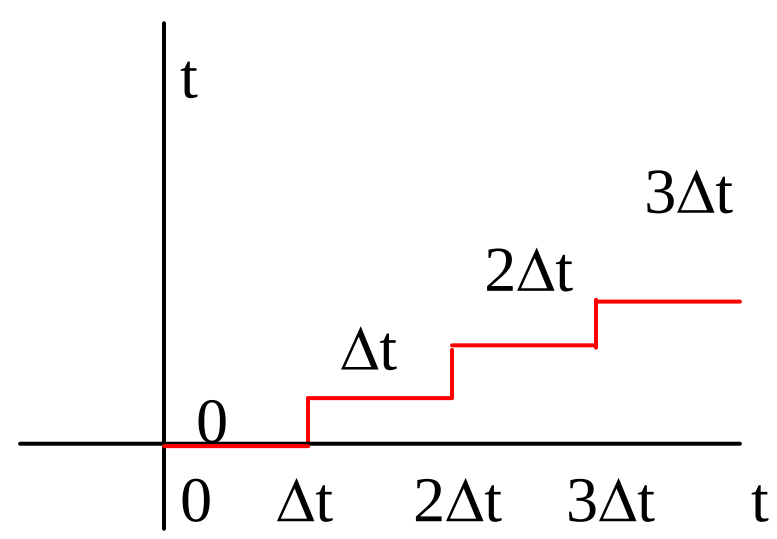

Series of weighted unit amplitude pulses of $\Delta t$ width

Fig. 2 System output

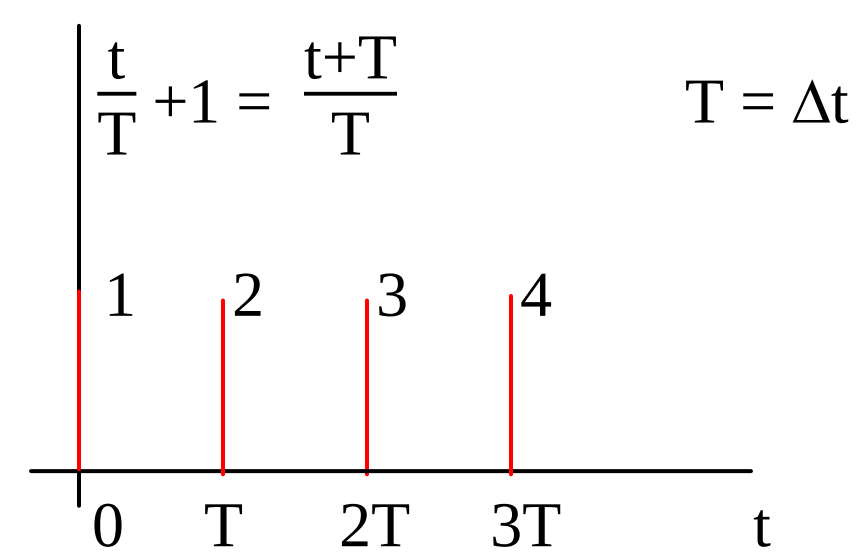

Series of weighted unit area impulses with an interval between impulses of T

From a comparison of the $K_{\Delta t}$ Transform and Z Transform calculation methodologies, some observed differences are listed below. Note Fig. 1 and Fig. 2.

1) Relative to the $K_{\Delta t}$ Transform system output, the Z Transform system output is shifted T to the left and divided by T

2) The $K_{\Delta t}$ Transform system output is a sample and hold shaped waveform whereas the Z Transform system output is a series of weighted unit area impulses.

3) By multiplying the Z Transform output by $Tz^{-1}$ it can be changed to equal the $K_{\Delta t}$ Transform output. Shift right by T then multiply by T.

4) By multiplying the $K_{\Delta t}$ Transform output by $\frac{(1+s\Delta t)^1}{\Delta t}$ it can be changed to equal the Z Transform output. Shift left by $\Delta t$ then divide by $\Delta t$.

5) The $K_{\Delta t}$ Transform system output shows an output function where its value changes only at $t = 0, \Delta t, 2\Delta t, 3\Delta t, \ldots$ and retains that value until a value change occurs. The resulting waveform is what would be expected from sample and hold sampling. The Z Transform system output shows an output function where its value is defined only at $t = 0, T, 2T, 3T, \ldots\ T = \Delta t$.

6) In appearance $K_{\Delta t}$ Transforms are more similar to Laplace Transforms than Z Transforms. In fact, $K_{\Delta t}$ Transforms are generalizations of Laplace Transforms for use where $\Delta t$ is not infinitesimal. If $\Delta t$ becomes infinitesimal, $K_{\Delta t}$ Transforms become Laplace Transforms.

7) The output function of the $K_{\Delta t}$ Transform analysis matches the Laplace Transform analysis of the same system in Example 5.12-3 where $L^{-1}[B^*(s)] = L^{-1}[\frac{\Delta t}{s(e^{s\Delta t}-1)}] = t$ , $t = 0, \Delta t, 2\Delta t, 3\Delta t, \ldots.$

<u>Some comments concerning Interval Calculus and the K$\Delta$t Transform</u>

Interval Calculus $K_{\Delta t}$ Transforms do not represent a sequence of unit inpulses as do Z Transforms. They represent a series of weighted unit amplitude pulses of $\Delta t$ width time increments. The sum total of these weighted unit amplitude pulses can be used to form a sample and hold representation of a continuous function of time with sampling at uniform intervals, at $t = 0, \Delta t, 2\Delta t, 3\Delta t, \ldots$ For an example of this see Diagram 5.12-12 above. Both the Z transform and the $K_{\Delta t}$ Transform can be used to analyze sampled-data systems. The solutions obtained will numerically be the same though the form of the solutions may not be the same. Typically, the solutions obtained using Z Transforms will be in terms of Calculus functions and the solutions obtained using $K_{\Delta t}$ Transforms will be in terms of Interval Calculus functions. Both solutions will be related, they will be Interval Calculus/Calculus identities. If desired, $K_{\Delta t}$ Transforms can be expressed in terms of Calculus functions that are defined at the discrete discrete values of t, at $t = 0, \Delta t, 2\Delta t, 3\Delta t, \ldots$ Tables 2, 3, 3a, and 3b in the Appendix identify $K_{\Delta t}$ Transforms and their corresponding Interval Calculus or Calculus functions

In this paper the term, Interval Calculus, has generally been used to refer to a discrete calculus where $\Delta t$ is a finite value. The term, Calculus, is used to refer to the commonly used Calculus of Newton and Leibniz where $\Delta t \to 0$. Initially in the development of Interval Calculus, these two terms provided a clear and convenient way to delineate between the type of calculus being considered, continuous or discrete. Later in the development of Interval Calculus something rather surprising was realized. Calculus was found to be a discrete mathematics where $\Delta t$ is an infinitesimal value (but not 0) and thus a subset of Interval Calculus where $\Delta t$ can be an infinitesimal value. Therefore, Interval Calculus includes Calculus. As a result of this, one can correctly refer to Calculus as being a part of Interval Calculus or Interval Calculus as beeing a generalization of Calculus. For clarity, it has continued to be advantageous to distinguish between continuous and discrete calculus using the terms, Calculus and Interval Calculus. However, the relationship just described between the two calculi should always be kept in mind. Since Calculus is a subset of Interval Calculus, another surprising fact has been realized. Interval Calculus is composed of an infinite number of calculus subsets, one for every value of $\Delta t$.

Below are several examples of the use of the $K_{\Delta t}$ Transform in the analysis of sampled-data systems.

<u>Example 5.12-6</u>

Find the outputs of the following two sampled-data systems. Use the $K_{\Delta t}$ Transform methodology.

1) Find the output, $f^*(t)$, for the sample and hold switch shown below.

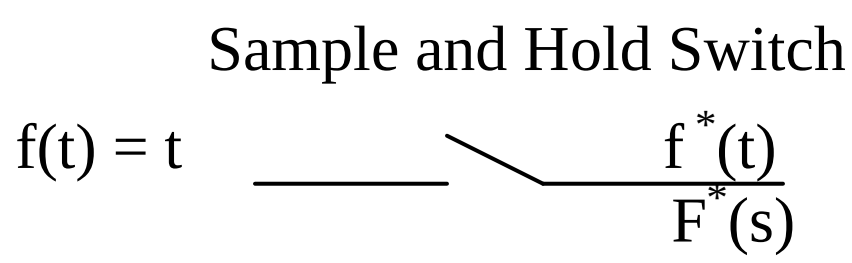

$$F^*(s) = K_{\Delta t}[f(t)] = K_{\Delta t}[t] = \frac{1}{s^2} \qquad 1)$$

Taking the Inverse $K_{\Delta t}$ Transform of Eq 1

**$f^*(t) = t$ ,   $t = 0, \Delta t, 2\Delta t, 3\Delta t, \ldots$ ,    t is a sample and hold shaped waveform**

The above result has been checked using Laplace Transforms in Example 5.12-1.

2) Find the output of the following integrator system with synchronously sampled input and output. The two switches are synchronous sample and hold switches. The integrator is initially passive.

<u>Fig. 1</u>

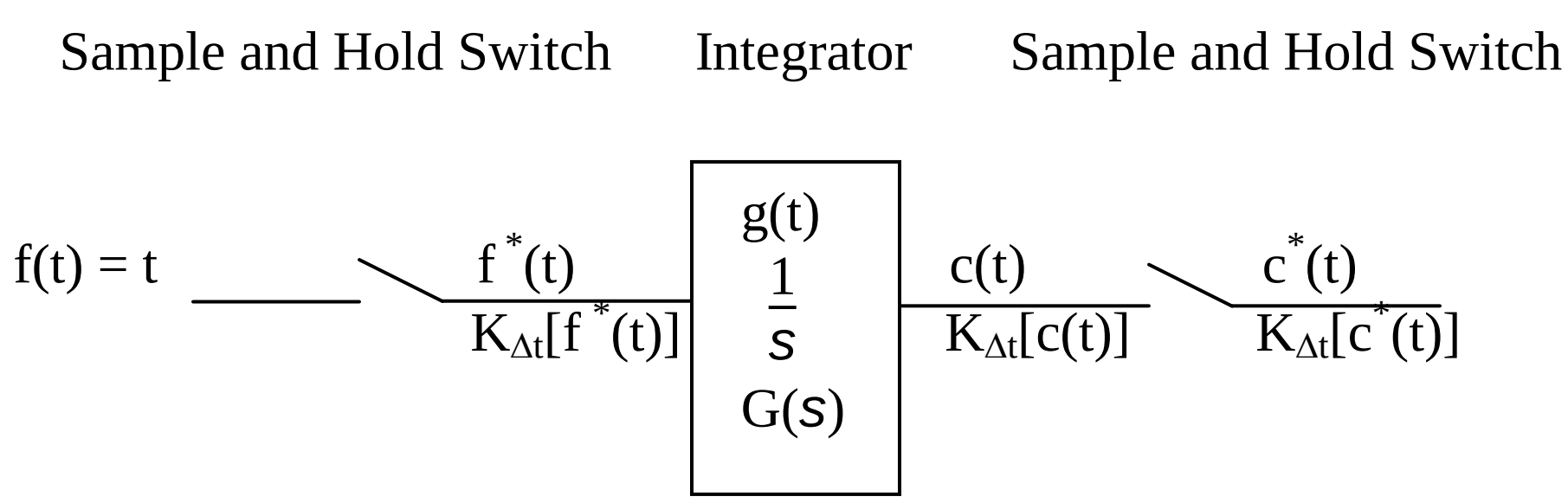

Redrawing the sampled-data system of Fig. 1 using only $K_{\Delta t}$ Transform notation.

$$g(t) = L^{-1}[G(s)] = L^{-1}[\frac{1}{s}] = U(t) \qquad 2)$$

$$G^*(s) = K_{\Delta t}[g(t)] = K_{\Delta t}[U(t)] = \frac{1}{s} \qquad 3)$$

<u>Fig. 2</u>

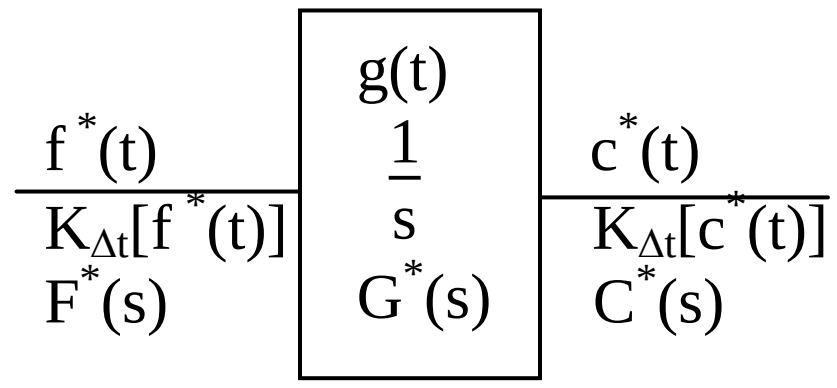

$f(t) = t$ 4)

$F^*(s) = K_{\Delta t}[f^*(t)] = K_{\Delta t}[f(t)] = K_{\Delta t}[t] = \frac{1}{s^2}$ 5)

$C^*(s) = G^*(s)F^*(s)$ 6)

Substituting Eq 3 and Eq 5 into Eq 6

$K_{\Delta t}[c^*(t)] = \frac{1}{s}\frac{1}{s^2} = \frac{1}{s^3}$ 7)

From Table 3 in the Appendix

$K_{\Delta t}^{-1}[\frac{n!}{s^{n+1}}] = \prod_{m=1}^{n}(t-[m-1]\Delta t) = [t]_{\Delta t}^{n}$ , n = 0, 1, 2, 3, … 8)

Taking the Inverse $K_{\Delta t}$ Transform of Eq 7 using Eq 8

$c^*(t) = \frac{1}{2}K_{\Delta t}^{-1}[\frac{2}{s^3}] = \frac{[t]_{\Delta t}^{2}}{2} = \frac{t(t-\Delta t)}{2}$ 9)

$\mathbf{c^*(t) = \frac{t(t-\Delta t)}{2}}$**,  t = 0, Δt, 2Δt, 3Δt, 4Δt, 5Δt, 6Δt** 10)

Referring to Example 5.12-2, the above result is seen to be correct. Good check

Example 5.12-7

For the following sampled-data feedback system find the maximum sampling rate for which the system is stable. The sample and hold switches are synchronized and the system is initially passive. Solve this problem using $K_{\Delta t}$ Transforms.

A Sampled-data Feedback System

$G(s) = \frac{3}{s+2}$

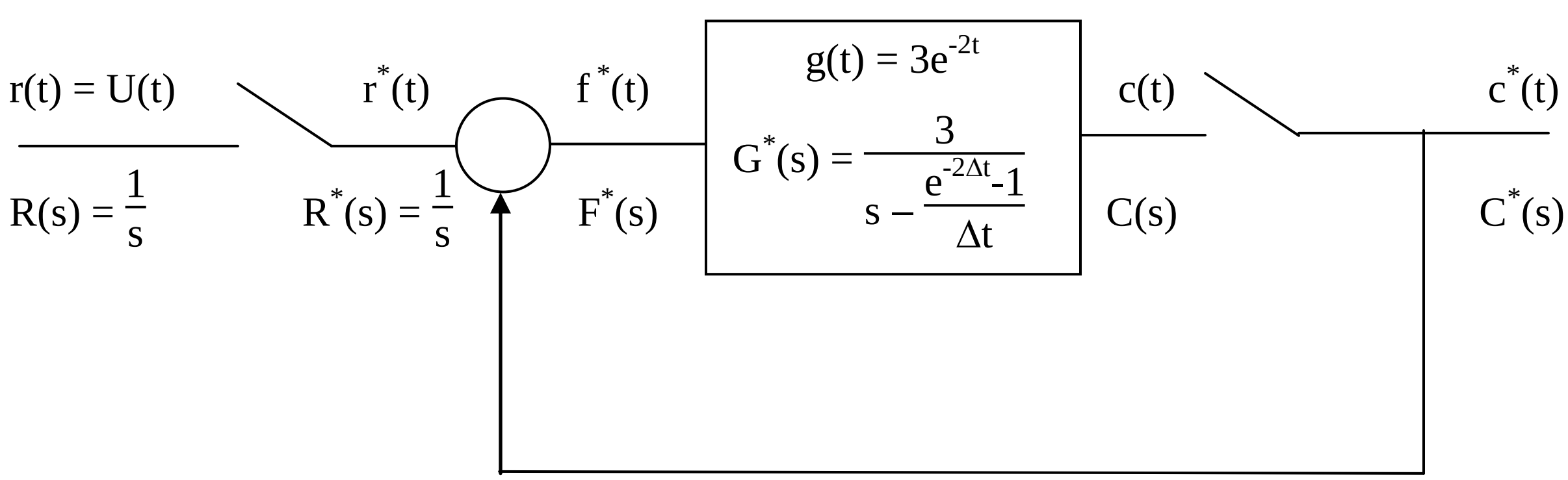

$$g(t) = L^{-1}[G(s)] = L^{-1}[\frac{3}{s+2}] = 3e^{-2t} \quad 1)$$

$$G^*(s) = K_{\Delta t}[g(t)] = K_{\Delta t}[3e^{-2t}] = \frac{3}{s - \frac{e^{-2\Delta t} - 1}{\Delta t}} \quad 2)$$

$$F^*(s) = R^*(s) - C^*(s) \quad 3)$$

$$C^*(s) = G^*(s)\,F^*(s) \quad 4)$$

Substituting Eq 4 into Eq 3

$$\frac{C^*(s)}{G^*(s)} = R^*(s) - C^*(s) \quad 5)$$

Simplifying Eq 5

$$\frac{C^*(s)}{G^*(s)} + C^*(s) = R^*(s) \quad 6)$$

$$C^*(s)[\frac{1}{G^*(s)} + 1] = R^*(s) \quad 7)$$

$$\frac{C^*(s)}{R^*(s)} = \frac{G^*(s)}{1+G^*(s)} \quad 8)$$

Substituting Eq 2 into Eq 8

$$\frac{C^*(s)}{R^*(s)} = \frac{\frac{3}{s - \frac{e^{-2\Delta t}-1}{\Delta t}}}{1 + \frac{3}{s - \frac{e^{-2\Delta t}-1}{\Delta t}}} \quad 9)$$

Simplifying Eq 9

$$\frac{C^*(s)}{R^*(s)} = \frac{3}{s - [\frac{e^{-2\Delta t}-1}{\Delta t} - 3]} \quad 10)$$

The denominator root, $s = \frac{e^{-2\Delta t}-1}{\Delta t} - 3$, must lie within the critical circle in the left half of the complex plane. The critical circle has its center at $-\frac{1}{\Delta t}$ and has a radius of $\frac{1}{\Delta t}$.

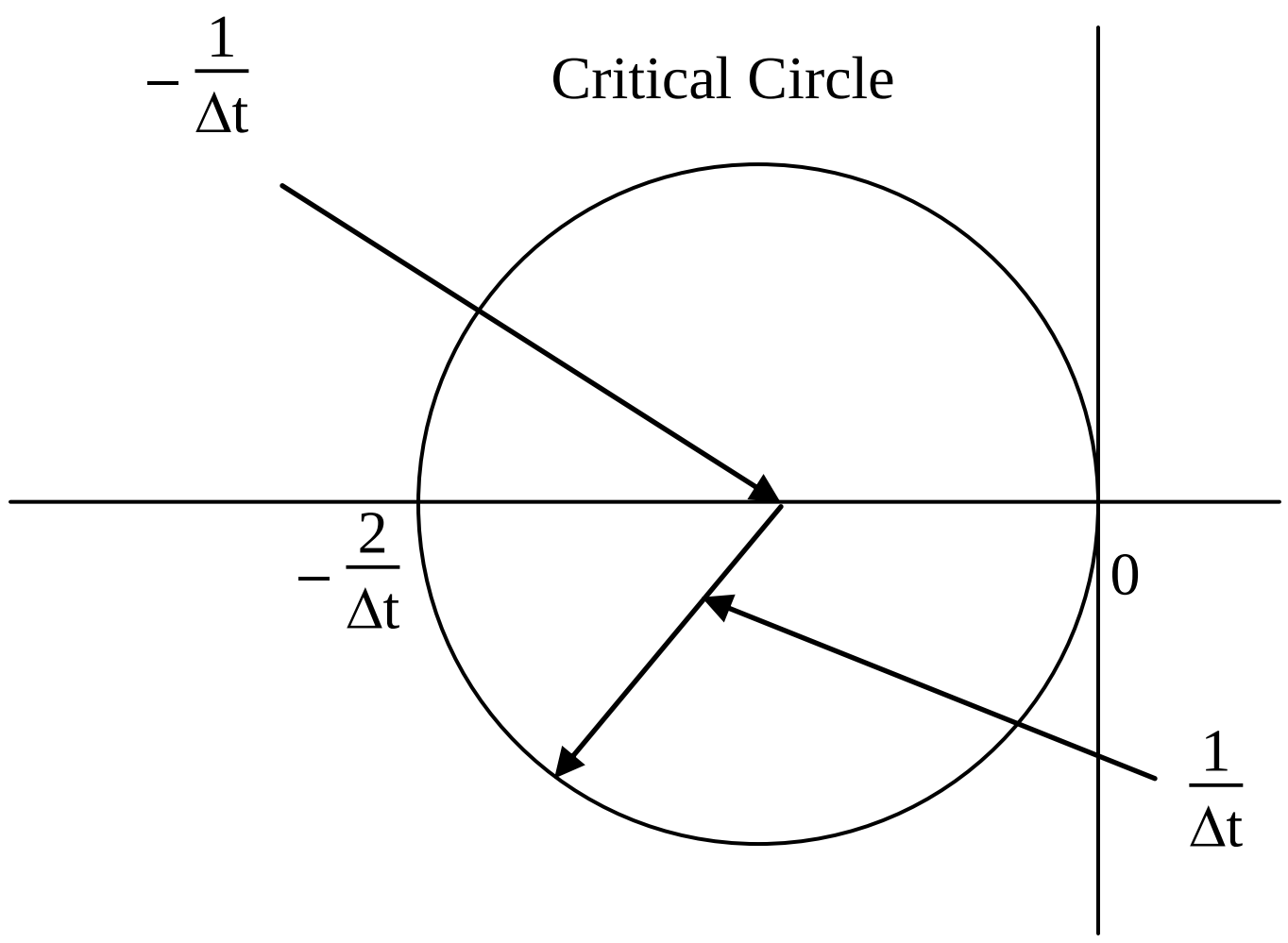

To find the maximum value of $\Delta t$ for system stability

$$-\frac{2}{\Delta t} < \frac{e^{-2\Delta t}-1}{\Delta t} - 3 < 0 \qquad 11)$$

Simplifying

$$-1 < e^{-2\Delta t} - 3\Delta t < 1 \qquad 12)$$

From Eq 12

Finding the limiting stability condition (oscillation)

$$-1 = e^{-2\Delta t} - 3\Delta t \qquad 13)$$

Solving Eq 13 using Newton's Method

$$\Delta t = .46488222 \qquad 14)$$

$$\text{Sampling rate} > \frac{1}{\Delta t} = \frac{1}{.46488222} = 2.151082483 \text{ samples/sec} \qquad 15)$$

Sampling rate > 2.151082483 samples/sec 16)

Then

**For system stability**

**$\Delta t$ < .46488222 sec/sample** 17)

**Sampling rate > 2.151082483 samples/sec** 18)

This example has been reanalyzed in Example 5.12-8 using Z Transforms. The results are the same. Good check

Using another check

Find the system oscillatory output at the minimum sampling rate of 2.151082483 samples/sec (.46488222 sec/sample).

The $K_{\Delta t}$ Transform of the input and output of a sample and hold switch are the same.

$$R^*(s) = R(s) = \frac{1}{s} \qquad 19)$$

$$R^*(s) = \frac{1}{s} \qquad 20)$$

From Eq 10

$$C^*(s) = R^*(s)\,\frac{3}{s - [\frac{e^{-2\Delta t}-1}{\Delta t} - 3]} = \frac{1}{s}\,[\,\frac{3}{s - [\frac{e^{-2\Delta t}-1}{\Delta t} - 3]}\,] = \frac{3}{s\,(s + 4.30216496)} \qquad 21)$$

$$C^*(s) = \frac{3}{s\,(s + 4.30216496)} = \frac{A}{s} + \frac{B}{s + 4.30216496} \qquad 22)$$

$$A = \frac{3}{(s + 4.30216496)}\Big|_{s = 0} = \frac{3}{4.30216496} = .69732333 \qquad 23)$$

$$B = \frac{3}{s}\Big|_{s = -4.30216496} = -\frac{3}{4.30216496} = -.69732333 \qquad 24)$$

Substituting Eq 23 and Eq 24 into Eq 22

$$C^*(s) = .69732333\,[\,\frac{1}{s} - \frac{1}{s + 4.30216496}\,] \qquad 25)$$

$$K_{\Delta t}[\,(1 - b\Delta t)^{\frac{t}{\Delta t}}\,] = \frac{1}{s + b}\,, \quad b = \text{constant} \qquad 26)$$

$$K_{\Delta t}[1] = \frac{1}{s} \qquad 27)$$

Taking the Inverse $K_{\Delta t}$ Transform of Eq 25 using Eq 26 and Eq 27

$$c^*(t) = .69732333\,[1 - [\,1 - (4.30216496)(\,.46488222)\,]^{\frac{t}{\Delta t}}\,] \qquad 28)$$

$$c^*(t) = .69732333\,[1 - (\,-1\,)^{\frac{t}{\Delta t}}\,] \qquad 29)$$

Then

$$\mathbf{c^*(t) = .69732333\,[1 - (-1)^{\frac{t}{\Delta t}}]} \qquad 30)$$

**where**

**t = 0, Δt, 2Δt, 3Δt, …**

**Δt = .46488222 sec/sample**

**$c^*(t)$ = system sample and hold shaped waveform output**

At the minimum sampling rate of 2.151082483 samples/sec (.46488222 sec/sample ) the system output waveform is oscillatory as expected.

Good check

Previously the $K_{\Delta t}$ Transform equation, $C^*(s) = G^*(s)F^*(s)$ was derived and demonstrated. A similar equation can be derived for the Z Transform, $C^*(z) = \mathcal{G}^*(z)F^*(z)$. The derivation of this equation follows.

Derivation of the Z Transform sample and hold system equation, $C^*(z) = \mathcal{G}^*(z)F^*(z)$

From Eq 5.12-98 and Eq 5.12-101

$$Z[c^*(t)] = Z[f^*(t)]Tz^{-1}Z[g(t)] = F^*(z)[Tz^{-1}G^*(z)] \qquad (5.12\text{-}123)$$

$$g(t) = L^{-1}[G(s)] \qquad (5.12\text{-}124)$$

$$G^*(z) = Z[g(t)] = Z[L^{-1}[G(s)]] \qquad (5.12\text{-}125)$$

$$F^*(z) = Z[f^*(t)] \qquad (5.12\text{-}126)$$

$$C^*(z) = Z[c^*(t)] \qquad (5.12\text{-}127)$$

Let

$$\mathcal{G}^*(z) = Tz^{-1}G^*(z) \qquad (5.12\text{-}128)$$

Substituting Eq 5.12-125 thru Eq 5.12-128 into Eq 5.12-123

$$C^*(z) = F^*(z)\,\mathcal{G}^*(z) \qquad (5.12\text{-}129)$$

Then

$$\mathbf{Z[c^*(t)] = Tz^{-1}Z[g(t)]\,Z[f^*(t)]} \qquad (5.12\text{-}130)$$

**or**

$$\mathbf{C^*(z) = F^*(z)\,\mathcal{G}^*(z) = Tz^{-1}G^*(z)F^*(z)} \qquad (5.12\text{-}131)$$

**where**

**$g(t) = L^{-1}[G(s)]$**

**$G^*(z) = Z[g(t)]$**

Fig. 1 represents the sample and hold sampled system of Fig. 2.

**Z Transform system transfer function**

**Fig. 1**

**Sample and hold sampling**

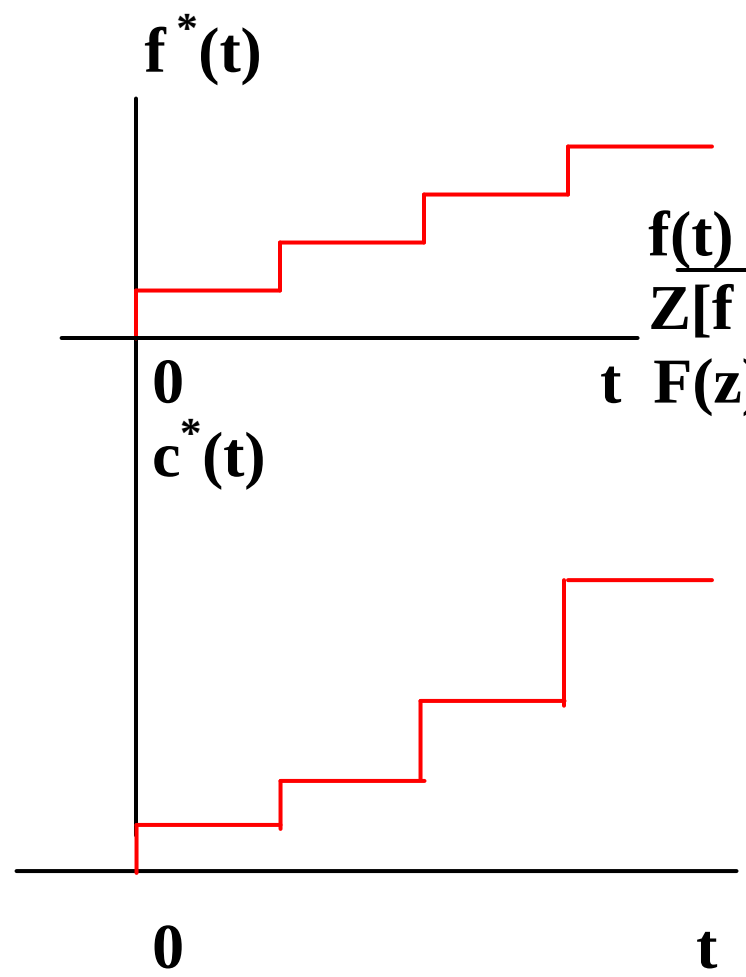

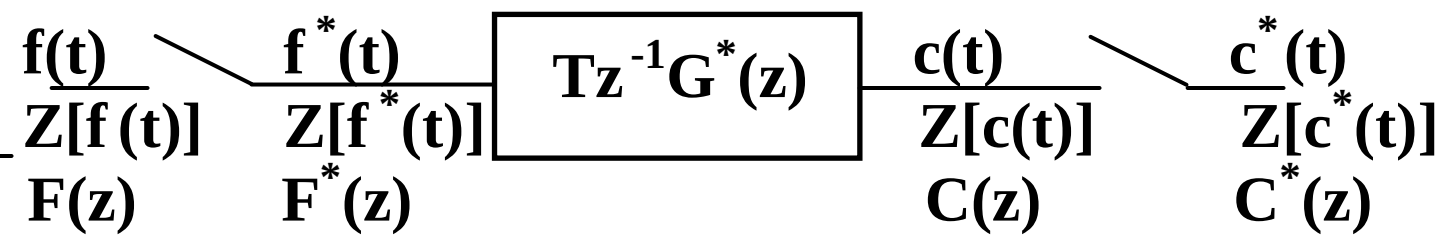

$Z[c^*(t)] = Tz^{-1}Z[g(t)]\,Z[f^*(t)]$

$C^*(z) = G^*(z)F^*(z) = Tz^{-1}G^*(z)F^*(z)$

$F^*(z) = Z[f^*(t)]$

$G^*(z) = Z[L^{-1}[G(s)]]$

$g(t) = L^{-1}[G(s)]$

The two switches shown in Fig. 1 above are synchronous sample and hold switches.

**Laplace Transform system transfer function**

**Fig. 2**

**No sampling**

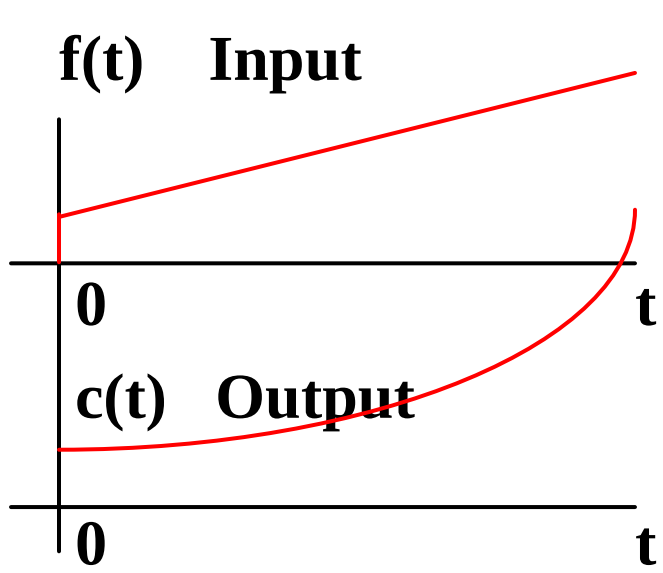

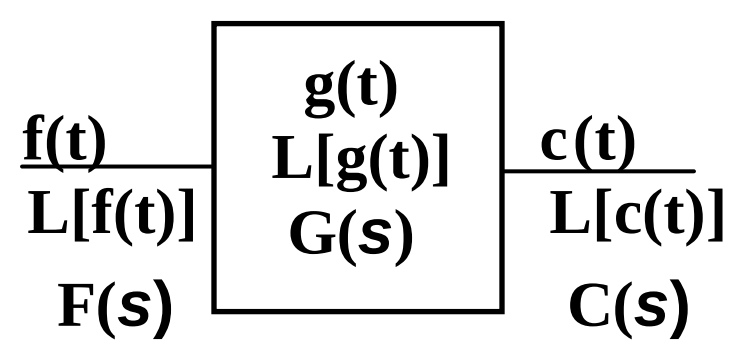

$L[c(t)] = L[g(t)]\,L[f(t)]$

$C(s) = G(s)F(s)$

$G(s) = L[g(t)]$

$g(t) = L^{-1}[G(s)]$

where

$f(t)$ = system unsampled input function

$F(z) = Z[f(t)]$, The Z Transform of the system input function, $f(t)$

$g(t) = L^{-1}[G(s)]$, Inverse Laplace Transform of the system Laplace transfer function, $G(s)$

$G(s) = L[g(t)]$ = System Laplace Transform transfer function

$c(t)$ = System unsampled output function

$C(z) = Z[c(t)]$, Z Transform of the output function, $c(t)$

$f^*(t)$ = system sample and hold sampled input function

$F^*(z) = Z[f^*(t)]$, The Z Transform of the sampled input function, $f^*(t)$

$c^*(t)$ = System sample and hold sampled output function

$C^*(z) = Z[c^*(t)]$, Z Transform of the sampled output function, $c^*(t)$

$G^*(z) = Z[g(t)]$

$\Delta t = T$ = Sampling interval

$z = e^{sT}$, Z Transform variable

$s$ = Laplace Transform variable

It has previously been shown that the term, $Tz^{-1}$, is associated with sample and hold sampling in a Z Transform analyzed system. This term can be shown to be derived directly from the definition of the $K_{\Delta t}$ Transform.

Consider the definition of the $K_{\Delta t}$Transform, the summation of weighted unit amplitude pulses of width, Δt.

$$K_{\Delta t}[f(t)] = {}_{\Delta x}\int_0^{\infty} (1+s\Delta t)^{-(\frac{t+\Delta t}{\Delta t})} f(t)\Delta t \ , \qquad \text{The } K_{\Delta t} \text{ Transform} \tag{5.12-132}$$

The value of s is defined as: $s = \dfrac{e^{s\Delta t} - 1}{\Delta t}$ , $-\dfrac{\pi}{\Delta t} \le w < \dfrac{\pi}{\Delta t}$ , $s = \gamma + jw$

where γ is any positive real value which makes s an indefinitely large value.

Substitute the equation, $s = \dfrac{e^{s\Delta t} - 1}{\Delta t}$ , into Eq 5.12-132 and T = Δt

$$K_{\Delta t}[f(t)] = {}_{\Delta x}\int_0^{\infty} e^{-s\Delta t(\frac{t+\Delta t}{\Delta t})} f(t)\Delta t = e^{-sT}\,{}_{\Delta x}\int_0^{\infty} e^{-sT\frac{t}{\Delta t}} f(t)\Delta t = Te^{-sT}\frac{1}{T}\,{}_{T}\int_0^{\infty} e^{-sT\frac{t}{\Delta t}} f(t)\Delta t \tag{5.12-133}$$

$$K_{\Delta t}[f(t)] = Te^{-sT}\frac{1}{T}\,{}_{T}\int_0^{\infty} e^{-sT\frac{t}{\Delta t}} f(t)\Delta t \tag{5.12-134}$$

$$z = e^{sT} \tag{5.12-135}$$

Substituting Eq 5.12-135 into Eq 5.12-134

$$K_{\Delta t}[f(t)] = Tz^{-1}\frac{1}{T}\,{}_{T}\int_0^{\infty} z^{-\frac{t}{\Delta t}} f(t)\Delta t \tag{5.12-136}$$

$$Z[f(t)] = \frac{1}{T}\,{}_{T}\int_0^{\infty} z^{-\frac{t}{\Delta t}} f(t)\Delta t \ , \quad \text{Definition of the Z Transform} \tag{5.12-137}$$

Substituting Eq 5.12-137 into Eq 5.12-1136

**$K_{\Delta t}[f(t)] = Tz^{-1}Z[f(t)]$** (5.12-138)

**where**

**T = Δt = sampling period**

**$z = e^{sT} = 1+s\Delta t$**

The $Tz^{-1}$ term converts a Z Transform unit area impulse into a unit amplitude pulse of width T

From the above investigation, the impulse transfer function to convert an impulse to a unit amplitude pulse of width T is:

Impulse to unit amplitude pulse transfer function

Unit area impulse input    Unit amplitude pulse of width T output

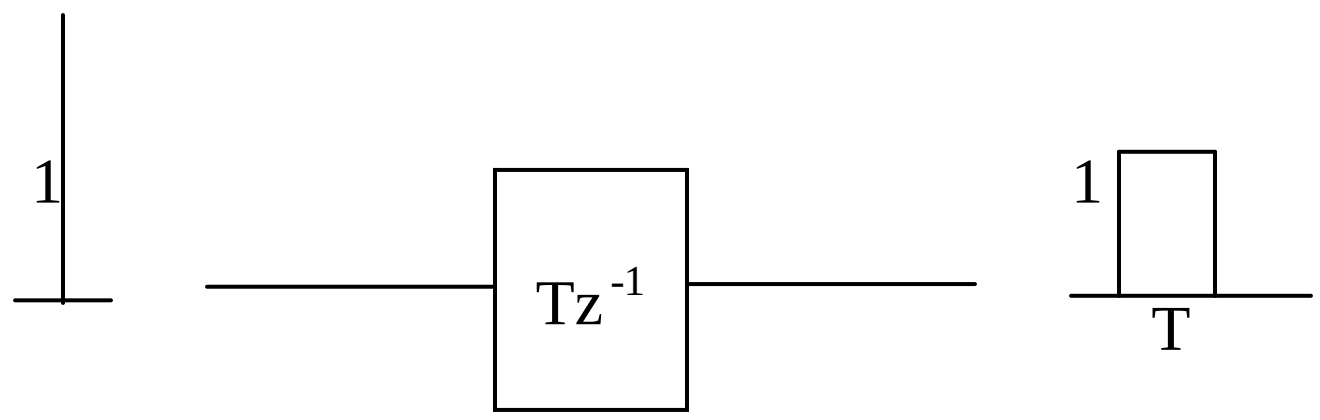

From Eq 5.12-138

Fig. 1 and Fig. 2 are equivalent

Z Transform transfer function

Fig. 1

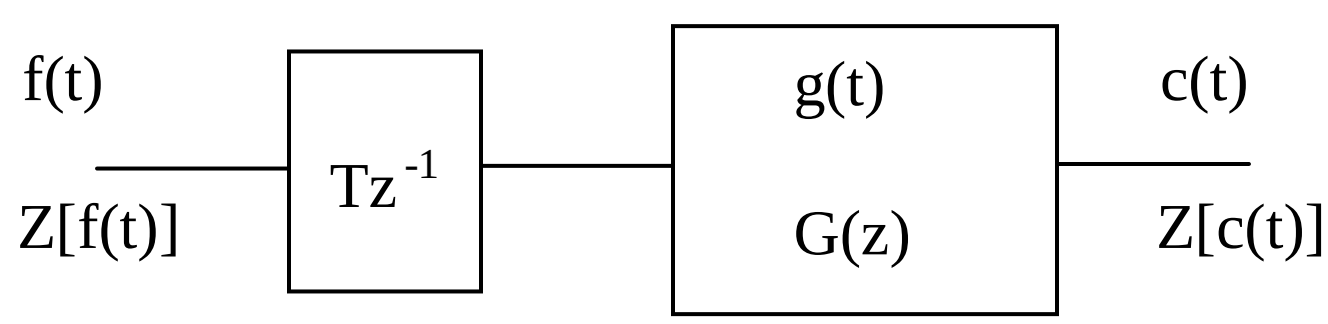

$G(z) = Z[g(t)]$
$g(t) = L^{-1}[G(s)]$

$K_{\Delta t}$ Transform Transfer function

Fig. 2

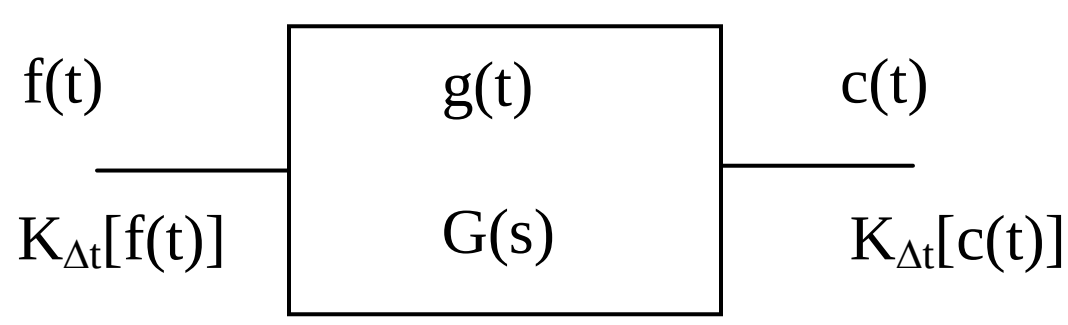

$G(s) = K_{\Delta t}[g(t)]$
$g(t) = L^{-1}[G(s)]$

The following example, Example 5.12-8, is a demonstration of the use of a sample and hold transfer function in a Z Transform system analysis. Example 5.12-8 is a reanalysis of the system of Example 5.12-7 using Z Transforms instead of $K_{\Delta t}$ Transforms.

Example 5.12-8

For the following sampled-data feedback system find the maximum sampling rate for which the system is stable. The sample and hold switches are synchronized and the system is initially passive. Solve this problem using Z Transforms.

A Sampled-data Feedback System

$$G(s) = \frac{3}{s+2}$$

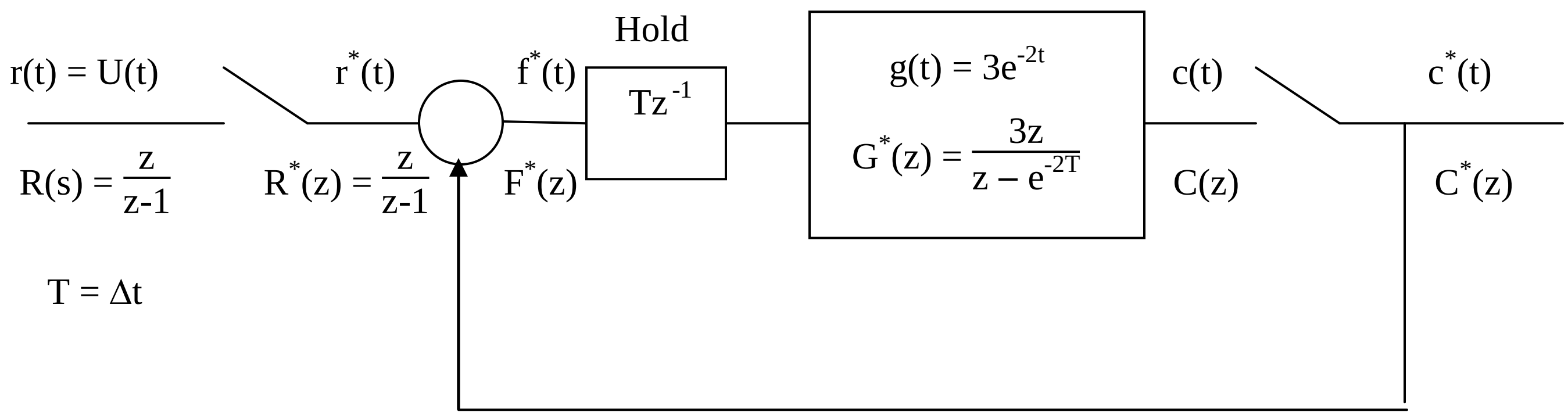

$$g(t) = L^{-1}[G(s)] = L^{-1}[\frac{3}{s+2}] = 3e^{-2t} \qquad 1)$$

$$G^*(z) = Z[g(t)] = Z[3e^{-2t}] = \frac{3z}{z - e^{-2T}} \qquad 2)$$

$$F^*(z) = R^*(z) - C^*(z) \qquad 3)$$

$$C^*(z) = Tz^{-1}G^*(z)\,F^*(z) \qquad 4)$$

Substituting Eq 4 into Eq 3

$$\frac{C^*(z)}{Tz^{-1}G^*(z)} = R^*(z) - C^*(z) \qquad 5)$$

Simplifying

$$C^*(z)[\,1 + \frac{1}{Tz^{-1}G^*(z)}\,] = R^*(z) \qquad 6)$$

$$\frac{C^*(z)}{R^*(z)} = \frac{Tz^{-1}G^*(z)}{1 + Tz^{-1}G^*(z)} \qquad 7)$$

Substituting Eq 2) into Eq 7

$$\frac{C^*(s)}{R^*(s)} = \frac{\dfrac{3T}{z - e^{-2T}}}{1 + \dfrac{3T}{z - e^{-2T}}} \qquad 8)$$

Simplifying Eq 7

$$\frac{C^*(s)}{R^*(s)} = \frac{3T}{z - [e^{-2T} - 3T]} \qquad 9)$$

The denominator root, $z = e^{-2T} - 3T$, must lie within the critical circle of the complex plane. The critical circle has its center at 0 and has a radius of 1.

Critical Circle

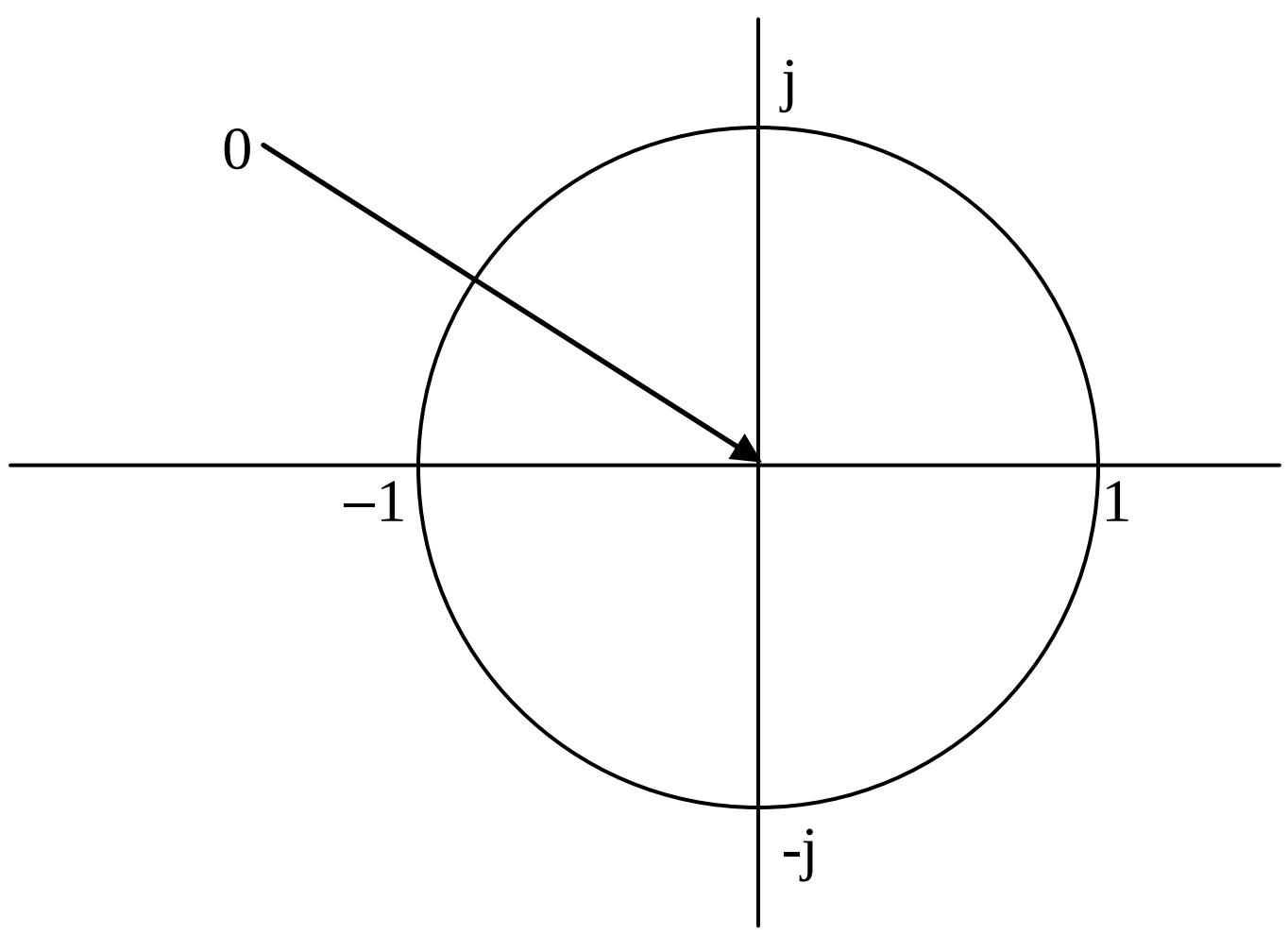

To find the maximum value of T for system stability

$$-1 = e^{-2T} - 3T < 1 \qquad 10)$$

From Eq 10

$$-1 = e^{-2T} - 3T \qquad 11)$$

Solving Eq 11 using Newton's Method

$$T = .46488222 \qquad 12)$$

$$\text{Sampling rate} > \frac{1}{T} = \frac{1}{.46488222} = 2.151082483 \text{ samples/sec} \qquad 13)$$

$$\text{Sampling rate} > 2.151082483 \text{ samples/sec} \qquad 14)$$

Then

**For system stability**

**T < .46488222 sec/sample** 15)

**Sampling rate > 2.151082483 samples/sec** 16)

The above result agrees with the $K_{\Delta t}$ Transform calculation of Example 5.12-7

Good check

The derivations above show a very close relationship between the $K_{\Delta t}$ Transform and the Z Transform.

<u>Derivation of the Laplace Transform sample and hold system equation, $C^*(s) = G^*(s)F^*(s)$</u>

Laplace Transform of a non-sampled system

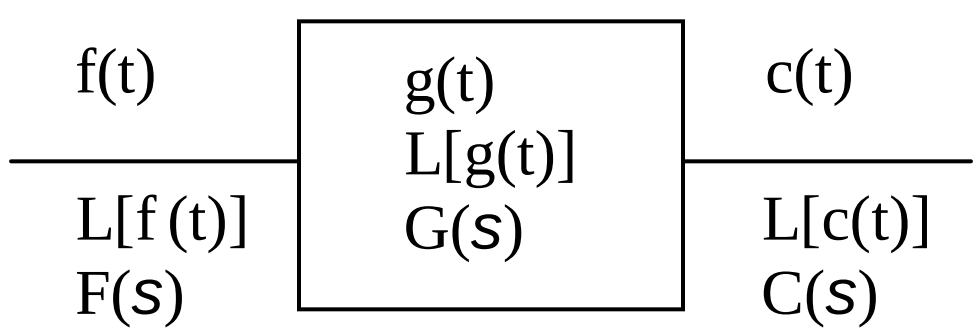

$g(t) = L^{-1}[G(s)]$

$K_{\Delta t}$ Transform of the sample and hold sampled system

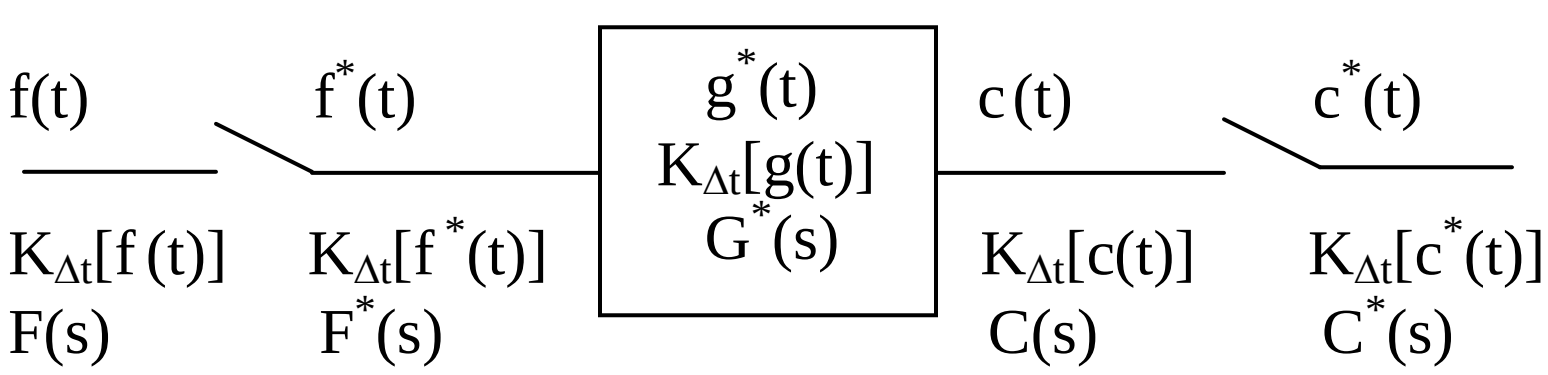

$G^*(s) = K_{\Delta t}[L^{-1}[G(s)]]$

$g^*(t) = g(t) = L^{-1}[G(s)]$

$K_{\Delta t}[c^*(t)] = K_{\Delta t}[g^*(t)]\, K_{\Delta t}[f^*(t)]$

$C^*(s) = G^*(s)F^*(s)$

$G^*(s) = Tz^{-1}G^*(z)\big|_{z\,=\,1+sT}$

$G^*(s) = K_{\Delta t}[g(t)] = Tz^{-1}Z[g(t)]\,\big|_{z\,=\,1+sT}$

Laplace Transform of the sample and hold sampled system

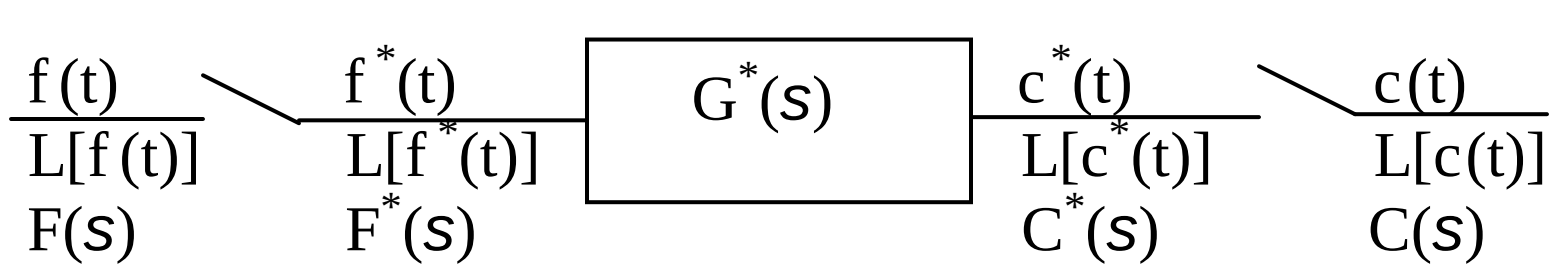

The switches shown above are synchronous sample and hold switches.

From Eq 5.12-104

$$K_{\Delta t}[c^*(t)] = K_{\Delta t}[g(t)]\,K_{\Delta t}[f^*(t)] \qquad (5.12\text{-}139)$$

where

$g(t) = L^{-1}[G(s)]$

Rewriting the sample and hold switch equations relating switch output to switch input

From Eq 5.12-21

$$L[f^*(t)] = \frac{1}{s}\,sK_{\Delta t}[f(t)]\,\Big|_{s = \frac{e^{s\Delta t}-1}{\Delta t}} \quad , \quad K_{\Delta t}\text{ Transforms are a function of s} \qquad (5.12\text{-}140)$$

From Eq 5.12-79

The following equation relates the $K_{\Delta t}$ Transform of the input of a sample and hold switch to the $K_{\Delta t}$ Transform of its output

$$K_{\Delta t}[f^*(t)] = K_{\Delta t}[f(t)] \qquad (5.12\text{-}141)$$

Comment – The functions, $f(t)$ and $f^*(t)$, have the same values at $t = 0, \Delta t, 2\Delta t, 3\Delta t, \ldots$ Thus, these functions, which have different waveforms, have the same $K_{\Delta t}$ Transform. The $K_{\Delta t}$ Transform is calculated only from the function values at $t = 0, \Delta t, 2\Delta t, 3\Delta t, \ldots$

Substituting Eq 5.12-141 into Eq 5.12-140

**Equation to convert the $K_{\Delta t}$ Transform of a sample and hold shaped waveform function into its equivalent Laplace Transform**

$$\mathbf{L[f^*(t)] = \frac{1}{s}\,sK_{\Delta t}[f^*(t)]\,\Big|_{s = \frac{e^{s\Delta t}-1}{\Delta t}}} \qquad (5.12\text{-}142)$$

Multiply both sides of Eq 5.12-139 by $\frac{1}{s}$ s then let $s = \frac{e^{s\Delta t}-1}{\Delta t}$

$$\frac{1}{s}\,sK_{\Delta t}[c^*(t)]\,\Big|_{s = \frac{e^{s\Delta t}-1}{\Delta t}} = \left(K_{\Delta t}[g(t)]\Big|_{s = \frac{e^{s\Delta t}-1}{\Delta t}}\right)\left(\frac{1}{s}\,sK_{\Delta t}[f^*(t)]\Big|_{s = \frac{e^{s\Delta t}-1}{\Delta t}}\right) \qquad (5.12\text{-}143)$$

From Eq 5.12-142 and Eq 5.12-143

$$L[c^*(t)] = \left[K_{\Delta t}[g(t)]\,\Big|_{s = \frac{e^{s\Delta t}-1}{\Delta t}}\right] L[f^*(t)] \qquad (5.12\text{-}144)$$

or

Writing Eq 5.12-144 in a different way

$$L[c^*(t)] = K_{\Delta t}[g(t)]\,L[f^*(t)] \qquad (5.12\text{-}145)$$

where

$$s = \frac{e^{s\Delta t} - 1}{\Delta t}$$

$$C^*(s) = L[c^*(t)] \tag{5.12-146}$$

$$F^*(s) = L[f^*(t)] \tag{5.12-147}$$

$$G^*(s) = K_{\Delta t}[g(t)]\Big|_{s = \frac{e^{s\Delta t} - 1}{\Delta t}} = G^*(s)\Big|_{s = \frac{e^{s\Delta t} - 1}{\Delta t}} \tag{5.12-148}$$

$$g(t) = L^{-1}[G(s)] \tag{5.12-149}$$

Substituting Eq 5.12-149 into Eq 148

$$G^*(s) = K_{\Delta t}[L^{-1}[G(s)]]\Big|_{s = \frac{e^{s\Delta t} - 1}{\Delta t}} \tag{5.12-150}$$

From Eq 5.12-144 and Eq 5.12-146 thru Eq 5.12-150

$$C^*(s) = G^*(s)\,F^*(s) \tag{5.12-151}$$

Thus

$$\mathbf{L[c^*(t)] = \left[K_{\Delta t}[g(t)]\Big|_{s = \frac{e^{s\Delta t} - 1}{\Delta t}}\right] L[f^*(t)]\,,\quad K_{\Delta t}[g(t)]\ \text{is a function of s}} \tag{5.12-152}$$

**or**

$$\mathbf{C^*(s) = G^*(s)\,F^*(s)} \tag{5.12-153}$$

**where**

**Fig. 1 represents the sample and hold sampled system of Fig. 2.**

**$K_{\Delta t}$ Transform system transfer function**

**Fig. 1** **Sample and hold sampling**

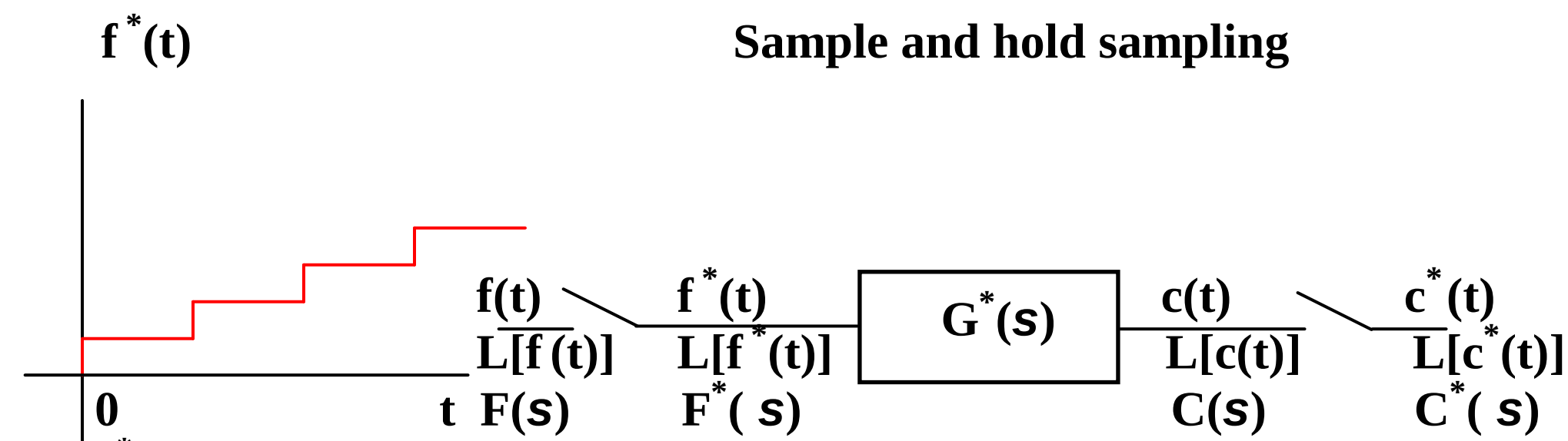

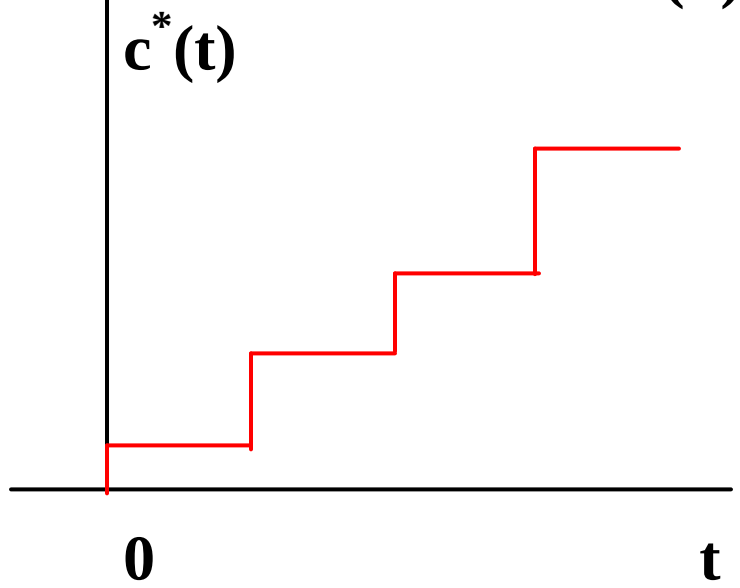

$$\mathbf{C^*(s) = G^*(s)F^*(s)}$$

$$\mathbf{F^*(s) = L[f^*(t)] = \frac{1}{s}\, sK_{\Delta t}[L^{-1}[F(s)]]\Big|_{s = \frac{e^{s\Delta t} - 1}{\Delta t}}}$$

$$\mathbf{G^*(s) = K_{\Delta t}[L^{-1}[G(s)]]\Big|_{s = \frac{e^{s\Delta t} - 1}{\Delta t}}}$$

$$\mathbf{G^*(s) = Tz^{-1}Z[L^{-1}[G(s)]]\Big|_{z = e^{s\Delta t}}}$$

$$G^*(s) = K_{\Delta t}[g(t)]\Big|_{s = \frac{e^{s\Delta t} - 1}{\Delta t}} = G^*(s)\Big|_{s = \frac{e^{s\Delta t} - 1}{\Delta t}}$$

**The two switches shown in the previous Fig. 1 are synchronous sample and hold switches.**

**Laplace Transform system transfer function**

**Fig. 2**

**No sampling**

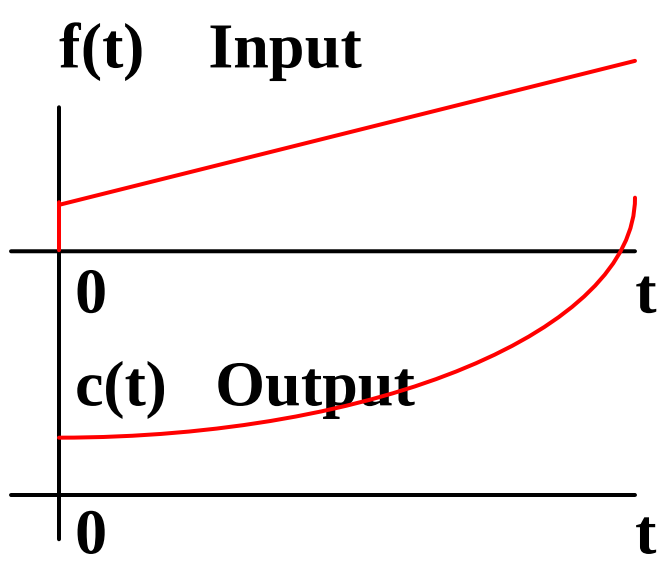

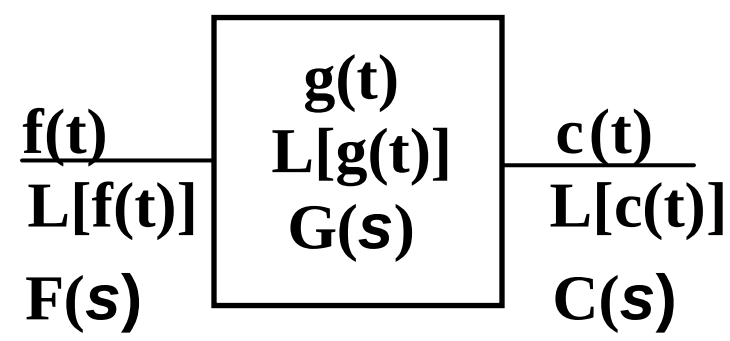

$C(s) = G(s)F(s)$
$F(s) = L[f(t)]$
$G(s) = L[g(t)]$
$g(t) = L^{-1}[G(s)]$

**f(t) = system unsampled input function**
**$F(s) = L[f(t)]$, The Laplace Transform of the system input function, f(t)**
**$g(t) = L^{-1}[G(s)]$, Inverse Laplace Transform of the system Laplace transfer function, G(s)**
**$G(s) = L[g(t)]$ = System Laplace Transform transfer function**
**c(t) = System unsampled output function**
**$C(s) = L[c(t)]$, Laplace Transform of the output function, c(t)**
**$f^*(t)$ = system sample and hold sampled input function**
**$F^*(s) = L[f^*(t)] = \frac{1}{s} sK_{\Delta t}[L^{-1}[F(s)]]$, The Laplace Transform of the sampled input function, $f^*(t)$**
**$c^*(t)$ = System sample and hold sampled output function**
**$C^*(s) = L[c^*(t)]$, Laplace Transform of the sampled output function, $c^*(t)$**

$$G^*(s) = K_{\Delta t}[L^{-1}[G(s)]]\Big|_{s = \frac{e^{s\Delta t} - 1}{\Delta t}}$$

**or**

$$G^*(s) = Tz^{-1}Z[L^{-1}[G(s)]]\Big|_{z = e^{s\Delta t}}$$

$$s = \frac{e^{s\Delta t} - 1}{\Delta t}$$

**$\Delta t = T$ = Sampling interval**
**s = $K_{\Delta t}$ Transform variable**
**s = Laplace Transform variable**

**To find the $K_{\Delta t}$ transform of a function there are Tables 2, 3, 3a, and 3b in the Appendix that may help.**

If a $K_{\Delta t}$ Transform analysis has been performed on a sample and hold sampled system and the system output, $K_{\Delta t}[c^*(t)]$, is known, its equivalent Laplace Transform, $L[c^*(t)]$, is easily found. Note the derivation below.

From Eq 5.12-153

$$C^*(s) = G^*(s)\,F^*(s) \qquad (5.12\text{-}154)$$

where

$$g(t) = L^{-1}[G(s)]$$
$$C^*(s) = L[c^*(t)]$$
$$G^*(s) = K_{\Delta t}[g(t)]\Big|_{s = \frac{e^{s\Delta t}-1}{\Delta t}} = G^*(s)\Big|_{s = \frac{e^{s\Delta t}-1}{\Delta t}}$$
$$F^*(s) = L[f^*(t)] = \frac{1}{s}\,sK_{\Delta t}[f^*(t)]$$
$$s = \frac{e^{s\Delta t}-1}{\Delta t}$$

From Eq 5.12-142

$$F^*(s) = L[f^*(t)] = \frac{1}{s}\,sK_{\Delta t}[f^*(t)] \qquad (5.1\text{-}155)$$

where

$$s = \frac{e^{s\Delta t}-1}{\Delta t}$$

Substituting Eq 5.12-155 into Eq 5.12-154

$$C^*(s) = [\frac{1}{s}\,s]\,G^*(s)\,K_{\Delta t}[f^*(t)] = [\frac{1}{s}\,s]\,G^*(s)\,F^*(s) = [\frac{1}{s}\,s]\,C^*(s) \qquad (5.12\text{-}156)$$

Then

**The conversion of the $K_{\Delta t}$ Transform of the output of a system to its Laplace Transform equivalent is:**

$$\mathbf{C^*(s) = [\frac{s}{s}]\,C^*(s)\,\Big|_{s = \frac{e^{s\Delta t}-1}{\Delta t}}} \qquad (5.12\text{-}157)$$

**where**

**$C^*(s) = L[c^*(t)]$**

**$C^*(s) = K_{\Delta t}[c^*(t)]$**

**$s = K_{\Delta t}$ Transform variable**

**$s$ = Laplace Transform variable**

Example 5.12-9 below demonstates the use of Eq 5.12-154 in the analysis of a sample and hold sampled system using Laplace Transforms. The system shown has been previously analyzed using other means.

<u>Example 5.12-9</u>

Find the output, $c^*(t)$, for the sample and hold sampled-data system shown below. The sample and hold switches are synchronized and the integrators are initially passive. Use the previously derived equation, $C^*(s) = G^*(s)\,F^*(s)$, to find $c^*(t)$.

Sample and hold system with two integrators

<u>Fig. 1</u>

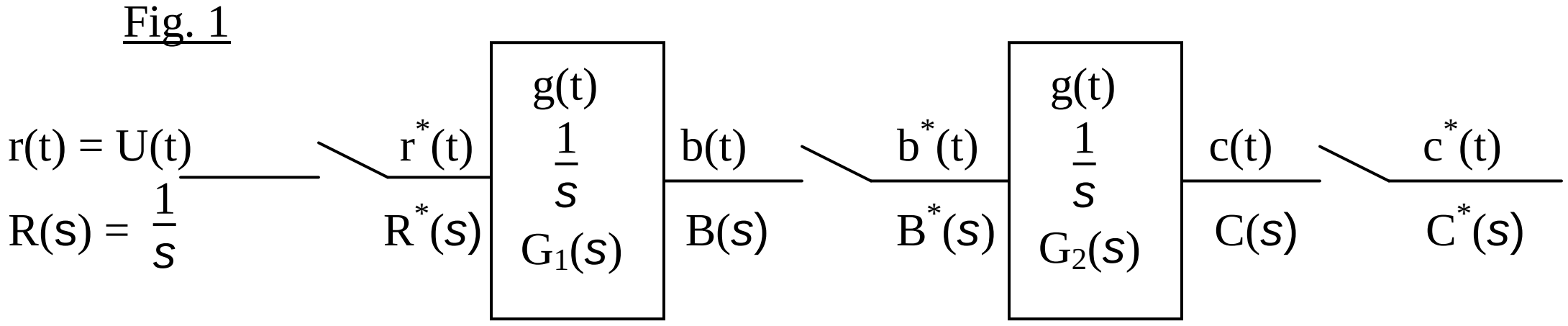

Changing the form of Fig. 1

<u>Fig. 2</u>

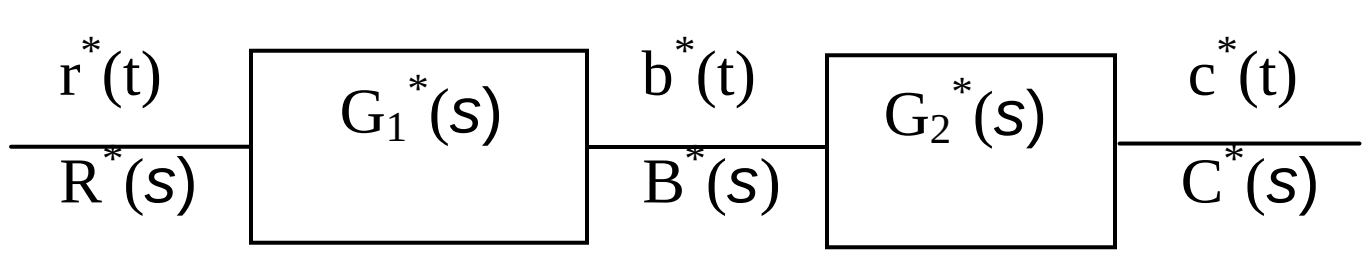

$$s = \frac{e^{s\Delta t} - 1}{\Delta t}$$

$$g(t) = L^{-1}[G_1(s)] = L^{-1}[G_2(s)] = L^{-1}[\,\frac{1}{s}\,] = U(t) \qquad 1)$$

$$G_1^*(s) = G_1^*(s) = K_{\Delta t}[g(t)] = K_{\Delta t}[U(t)] = \frac{1}{s} \qquad 2)$$

where

$$s = \frac{e^{s\Delta t} - 1}{\Delta t}$$

Using the sample and hold switch equation to find $R^*(s)$

$$R^*(s) = \frac{1}{s}\, sK_{\Delta t}[L^{-1}[R(s)]\,\Big|_{s\,=\,\frac{e^{s\Delta t} - 1}{\Delta t}} \qquad 3)$$

$$R(s) = \frac{1}{s} \qquad 4)$$

Substituting Eq 4 into Eq 3

$$R^*(s) = \frac{1}{s} sK_{\Delta t}[L^{-1}[\frac{1}{s}] \,|_{s = \frac{e^{s\Delta t} - 1}{\Delta t}} = R^*(s) = \frac{1}{s} sK_{\Delta t}[U(t)] \,|_{s = \frac{e^{s\Delta t} - 1}{\Delta t}} = \frac{1}{s}\frac{s}{s} \,|_{s = \frac{e^{s\Delta t} - 1}{\Delta t}} = \frac{1}{s} \qquad 5)$$

$$R^*(s) = \frac{1}{s} \qquad 6)$$

$$\frac{B^*(s)}{R^*(s)} = G_1{}^*(s) = \frac{1}{s} \qquad 7)$$

$$\frac{C^*(s)}{B^*(s)} = G_2{}^*(s) = \frac{1}{s} \qquad 8)$$

Multiplying Eq 7 and Eq 8

$$\frac{B^*(s)}{R^*(s)}\frac{C^*(s)}{B^*(s)} = \frac{C^*(s)}{R^*(s)} = G_1{}^*(s)\,G_2{}^*(s) = \frac{1}{s}\frac{1}{s} = \frac{1}{s^2} \qquad 9)$$

$$\frac{C^*(s)}{R^*(s)} = G_1{}^*(s)\,G_2{}^*(s) = \frac{1}{s^2} \qquad 10)$$

where

$$s = \frac{e^{s\Delta t} - 1}{\Delta t}$$

From Eq 10 and Eq 6

$$C^*(s) = \frac{1}{s}\left[\frac{1}{s^2}\,|_{s = \frac{e^{s\Delta t} - 1}{\Delta t}}\right] = \frac{\Delta t^2}{s\,(e^{s\Delta t} - 1)^2} \qquad 11)$$

$$C^*(s) = \frac{\Delta t^2}{s}\left(\frac{1}{e^{s\Delta t} - 1}\right)^2 \qquad 12)$$

Expanding $\left(\frac{1}{e^{s\Delta t}-1}\right)^2$

$$\left(\frac{1}{e^{s\Delta t}-1}\right)^2 = e^{-2s\Delta t} + 2e^{-3s\Delta t} + 3e^{-4s\Delta t} + 4e^{-5s\Delta t} + 5e^{-6s\Delta t} + \ldots \qquad 13)$$

Substituting Eq 14 into Eq 13

$$C^*(s) = \Delta t^2\left[\frac{e^{-2s\Delta t}}{s} + 2\frac{e^{-3s\Delta t}}{s} + 3\frac{e^{-4s\Delta t}}{s} + 4\frac{e^{-5s\Delta t}}{s} + 5\frac{e^{-6s\Delta t}}{s} + \ldots\right] \qquad 14)$$

Taking the Laplace Transform of Eq 14

$$c^*(t) = \Delta t^2[\ 1U(t-2\Delta t) + 2U(t-3\Delta t) + 3U(t-4\Delta t) + 4U(t-5\Delta t) + \dots\ ] \qquad 15)$$

$$c^*(t) = \Delta t^2 \sum_{n=0}^{\infty}(n+1)U(t-[n+2]\Delta t) \qquad 16)$$

Calculate $c^*(t)$ from Eq 17 for t = 0, Δt, 2Δt, 3Δt, 4Δt, 5Δt, 6Δt

| $t$ | $c^*(t)$ |
|---|---|
| 0 | 0 |
| $\Delta t$ | 0 |
| $2\Delta t$ | $\Delta t^2$ |
| $3\Delta t$ | $3\Delta t^2$ |
| $4\Delta t$ | $6\Delta t^2$ |
| $5\Delta t$ | $10\Delta t^2$ |
| $6\Delta t$ | $15\Delta t^2$ |

$$\mathbf{c^*(t) = \Delta t^2 \sum_{n=0}^{\infty}(n+1)U(t-[n+2]\Delta t)} \qquad 17)$$

**or**

$$\mathbf{c^*(t) = \frac{t(t-\Delta t)}{2}} \qquad 18)$$

The above result is checked by the results of Example 5.12-3 and can also be checked using the $K_{\Delta t}$ Transform equation, $C^*(s) = G^*(s)R^*(s)$. This check follows.

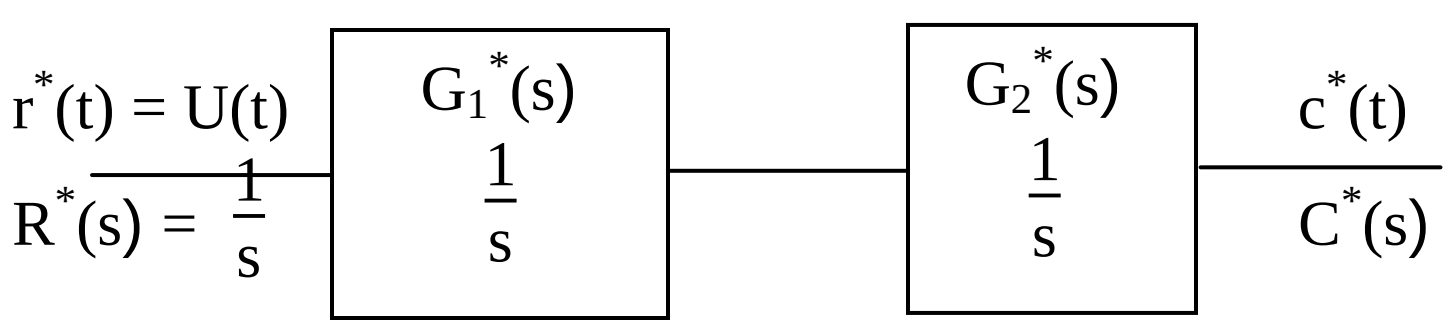

$$R^*(s) = \frac{1}{s} \qquad 19)$$

$$G^*(s) = \frac{1}{s}\frac{1}{s} = \frac{1}{s^2} \qquad 20)$$

$$C^*(s) = G^*(s)R^*(s) \qquad 21)$$

$$C^*(s) = \left(\frac{1}{s^2}\right)\left(\frac{1}{s}\right) = \frac{1}{s^3} \qquad 22)$$

$$C^*(s) = \frac{1}{s^3} \qquad 23)$$

From Table 3 in the Appendix

$$K_{\Delta t}^{-1}\left[\frac{n!}{s^{n+1}}\right] = \prod_{m=1}^{n}(t-[m-1]\Delta t) = [t]_{\Delta t}^{n}\ , \quad n = 0, 1, 2, 3, \ldots \qquad 24)$$

Taking the Inverse $K_{\Delta t}$ Transform of Eq 23 using Eq 24

$$c^*(t) = \frac{1}{2}K_{\Delta t}^{-1}\left[\frac{2}{s^3}\right] = \frac{[t]_{\Delta t}^{2}}{2} = \frac{t(t-\Delta t)}{2} \qquad 25)$$

$$\mathbf{c^*(t) = \frac{t(t-\Delta t)}{2}} \qquad 26)$$

Good check

Previously, the Laplace Transform system equation, $C^*(s) = G^*(s)\,F^*(s)$, was derived with $G^*(s)$ defined in terms of a $K_{\Delta t}$ Transform where $G^*(s) = G^*(s)\big|_{s = \frac{e^{s\Delta t}-1}{\Delta t}}$. It is also possible to define $G^*(s)$ in terms of a Z Transform where $G^*(s) = [Tz^{-1}G^*(z)]\big|_{z = e^{s\Delta t}}$. The derivation of this equation follows.

<u>Derivation of the Laplace Transform system equation, $C^*(s) = G^*(s)F^*(s)$, where $G^*(s) = Tz^{-1}G^*(z)$</u>

Laplace Transform of a non-sampled system

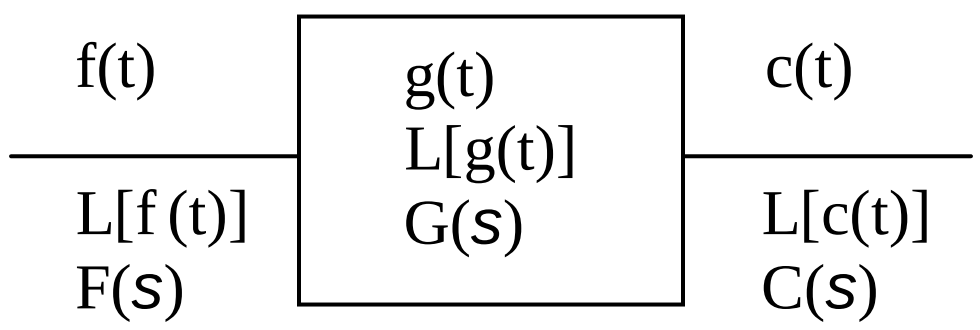

Z Transform of an impulse sampled system

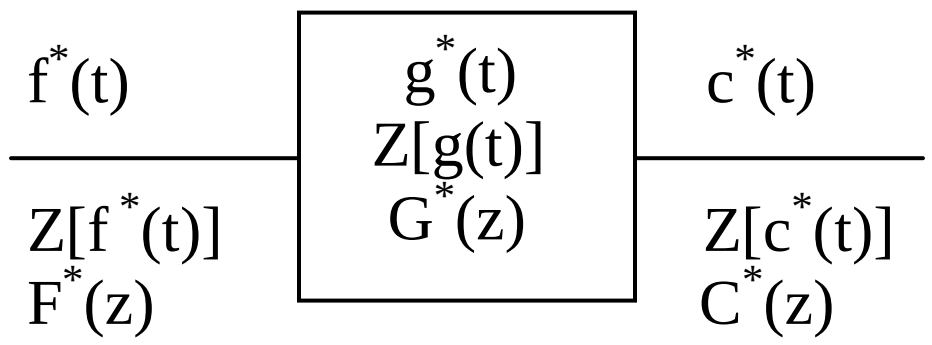

$G^*(z) = Z[L^{-1}[G(s)]]$

$g^*(t) = g(t) = L^{-1}[G(s)]$

$Z[c^*(t)] = Z[g^*(t)]\,Z[f^*(t)]$

$C^*(z) = G^*(z)F^*(z)$

$G^*(z) = Z[g(t)]$

$z = e^{sT}$

Laplace Transform of a sample and hold sampled system

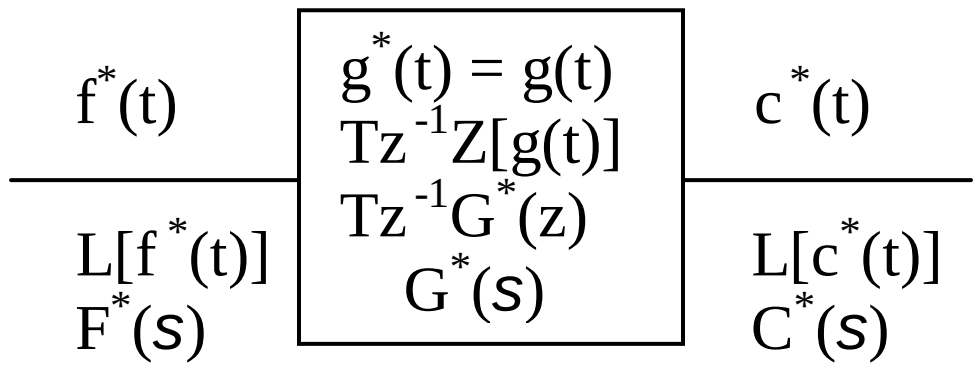

$G^*(s) = Tz^{-1}G^*(z)]$
$z = e^{sT}$
$T = \Delta t$
$G^*(z) = Z[L^{-1}[G(s)]]$
$G^*(z) = Z[g(t)]$
$g^*(t) = g(t) = L^{-1}[G(s)]$
$C^*(s) = [Tz^{-1}G^*(z)]F^*(s)$

From Eq 5.12-152

$$L[c^*(t)] = \left[K_{\Delta t}[g(t)]\,\Big|_{s = \frac{e^{s\Delta t} - 1}{\Delta t}}\right] L[f^*(t)] \qquad (5.12\text{-}158)$$

$$C^*(s) = G^*(s)\,F^*(s) \qquad (5.12\text{-}159)$$

$$G^*(s) = G^*(s) = K_{\Delta t}[g(t)] \qquad (5.12\text{-}160)$$

where

$$s = \frac{e^{s\Delta t} - 1}{\Delta t}$$

$$s = \frac{z - 1}{\Delta t}$$

$$z = e^{s\Delta t}$$

$$z = 1+s\Delta t$$

From Eq 5.12-107

$$K_{\Delta t}[g(t)] = Tz^{-1}Z[g(t)]\,\Big|_{z = 1+s\Delta t} = Tz^{-1}G^*(z)\,\Big|_{z = 1+s\Delta t} \qquad (5.12\text{-}161)$$

Substituting Eq 5.12-161 into Eq 5.12-160

$$G^*(s) = K_{\Delta t}[g(t)] = Tz^{-1}G^*(z) \qquad (5.12\text{-}162)$$

$$G^*(s) = Tz^{-1}G^*(z) \qquad (5.12\text{-}163)$$

where

$$z = e^{sT}$$

Substituting Eq 5.12-163 into Eq 5.12-159

$$C^*(s) = [Tz^{-1}G^*(z)]\,F^*(s) \qquad (5.12\text{-}164)$$

where

$$z = e^{sT}$$

Thus

$$\mathbf{C^*(s) = G^*(s)\,F^*(s)} \qquad (5.12\text{-}165)$$

**where**

$G^*(s) = Tz^{-1}G^*(z)$

$G^*(z) = Z[g(t)] = Z[L^{-1}[G(s)]]$

$z = e^{sT}$

$g^*(t) = g(t) = L^{-1}[G(s)]$

$C^*(s) = L[c^*(t)]$ , **Laplace Transform of the sampled output function,** $c^*(t)$

$F^*(s) = L[f^*(t)]$ , **The Laplace Transform of the sampled input function,** $f^*(t)$

$G(s) = L[g(t)] =$ **System Laplace Transform transfer function**

$G^*(z) = Z[g(t)] =$ **Z Transform of the function** $g(t)$

$T = \Delta t =$ **sampling interval**

$s =$ **Laplace Transform variable**

$F(s) = L[f(t)]$ , **The Laplace Transform of the input function,** $f(t)$

$F^*(s) = L[f^*(t)] = \dfrac{1}{s}\,\dfrac{z-1}{z}\,Z[L^{-1}[F(s)]]$ , **The Laplace Transform of the sampled input function,** $f(t)$

An example of the analysis of a system using the Laplace Transform with the above Z Transform methodology is shown in Example 5.12-8 below.

Example 5.12-10

Find the output, $c^*(t)$ of the integrator system with synchronously sampled input and output that is shown below. The two switches are synchronous sample and hold switches. The integrator is initially passive. Use a Z Transform methodology to find $c^*(t)$.

Fig. 1 Sample and Hold Switch Integrator Sample and Hold Switch

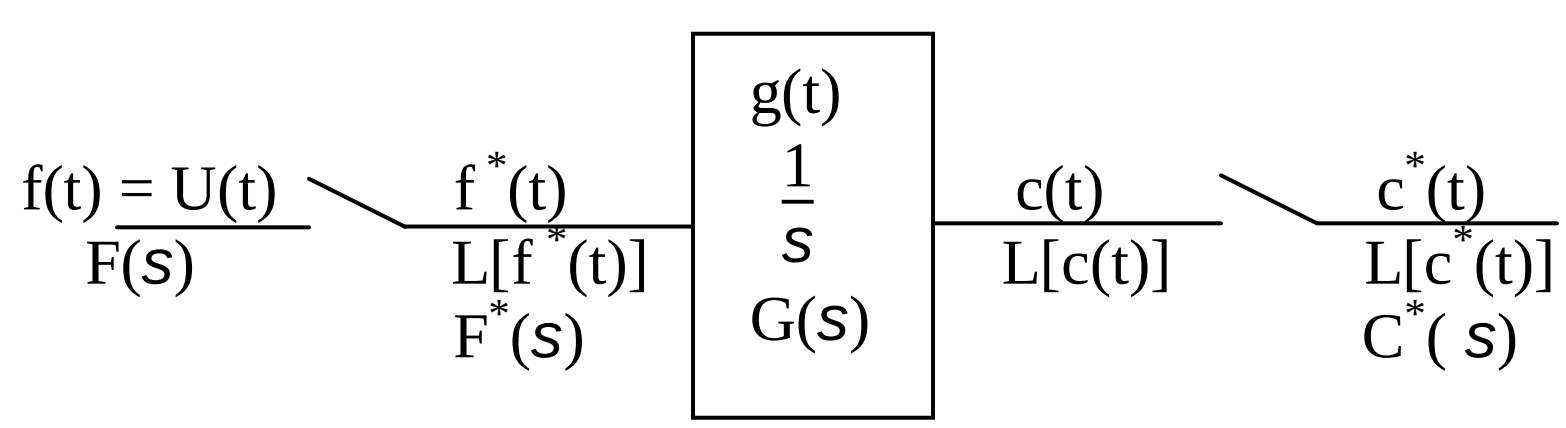

Redrawing the sampled-data system of Fig. 1 using only Laplace Transform and Z Transform notation.

Fig. 2

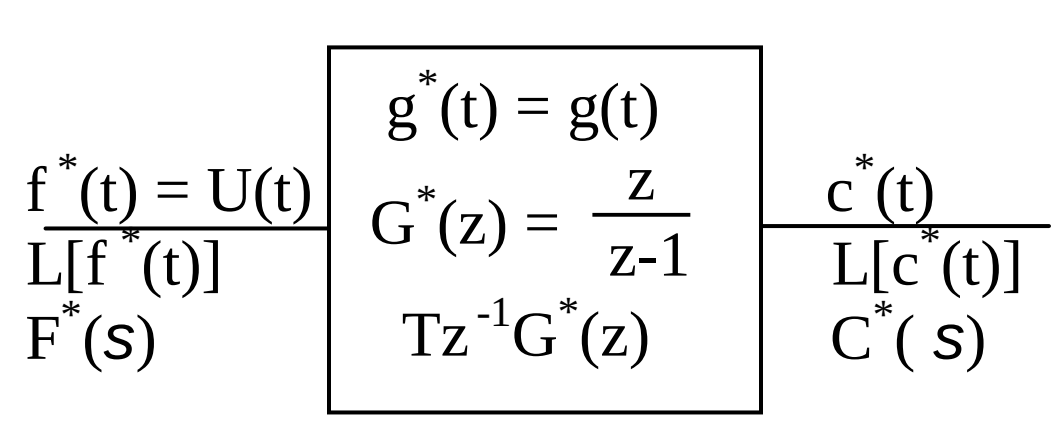

$$C^*(s) = [Tz^{-1}G^*(z)]\,F^*(s) \qquad 1)$$

$$g(t) = L^{-1}[G(s)] = L^{-1}[\frac{1}{s}] = U(t) \qquad 2)$$

$$G^*(z) = Z[g(t)] = Z[U(t)] = \frac{z}{z-1} \qquad 3)$$

$$Tz^{-1}G^*(z) = Tz^{-1}\frac{z}{z-1} = \frac{T}{z-1} \qquad 4)$$

$$Tz^{-1}G^*(z) = \frac{T}{z-1} \qquad 5)$$

$$z = e^{sT} \qquad 6)$$

Substituting Eq 6 into Eq 5)

$$Tz^{-1}G^*(z) = \frac{T}{e^{sT}-1} \qquad 7)$$

$$F(s) = L[f(t)] = L[U(t)] = \frac{1}{s} \qquad 8)$$

$$F^*(s) = L[f^*(t)] = \frac{1}{s}\,\frac{z-1}{z}\,Z[L^{-1}[F(s)] \qquad 9)$$

Substituting Eq 8 into Eq 9 to find $F^*(s)$

$$F^*(s) = \frac{1}{s}\,\frac{z-1}{z}\,Z[L^{-1}[\frac{1}{s}] = \frac{1}{s}\,\frac{z-1}{z}\,Z[U(t)] = \frac{1}{s}\,\frac{z-1}{z}\,\frac{z}{z-1} = \frac{1}{s} \qquad 10)$$

$$F^*(s) = \frac{1}{s} \qquad 11)$$

Substituting Eq 11 and Eq 7 into Eq 1

$$C^*(s) = [Tz^{-1}G^*(z)]\,F^*(s) \qquad 12)$$

$$C^*(s) = \frac{T}{e^{sT}-1}\,[\frac{1}{s}] \qquad 13)$$

$$C^*(s) = \frac{T}{s(e^{sT}-1)} \qquad 14)$$

$$C^*(s) = L[f^*(t)] = \frac{T}{s(e^{sT} - 1)} \qquad 15)$$

Expanding Eq 15

$$L[f^*(t)] = e^{-sT}\frac{T}{s} + e^{-2sT}\frac{T}{s} + e^{-3sT}\frac{T}{s} + e^{-4sT}\frac{T}{s} + e^{-5sT}\frac{T}{s} + \ldots \qquad 16)$$

Taking the Inverse Laplace Transform of Eq 16

$$f^*(t) = L^{-1}[L[f^*(t)]] = 0U(t) + TU(t-T) + TU(t-2T) + TU(t-3T) + TU(t-4T) + TU(t-5T) + \ldots \qquad 17)$$

Then

$$\mathbf{f^*(t) = {}_1\sum_{n=0}^{\infty} T\,U(t-[n+1]T)} \qquad 18)$$

**or**

$$\mathbf{f^*(t) = nT \quad \text{where } nT \le t < [n+1]T\,, \quad n = 0, 1, 2, 3, \ldots} \qquad 19)$$

**where**

**$0 \le t < \infty$**

**T = interval between samples**

**or**

$$\mathbf{f^*(t) = t} \qquad 20)$$

**where**

**$f^*(t)$ = Sample and hold shaped waveform**

**t = 0, T, 2T, 3T, …**

Checking the previous result using $K_{\Delta t}$ Transforms

Redrawing the sampled-data system of Fig. 1 using $K_{\Delta t}$ Transform notation notation.

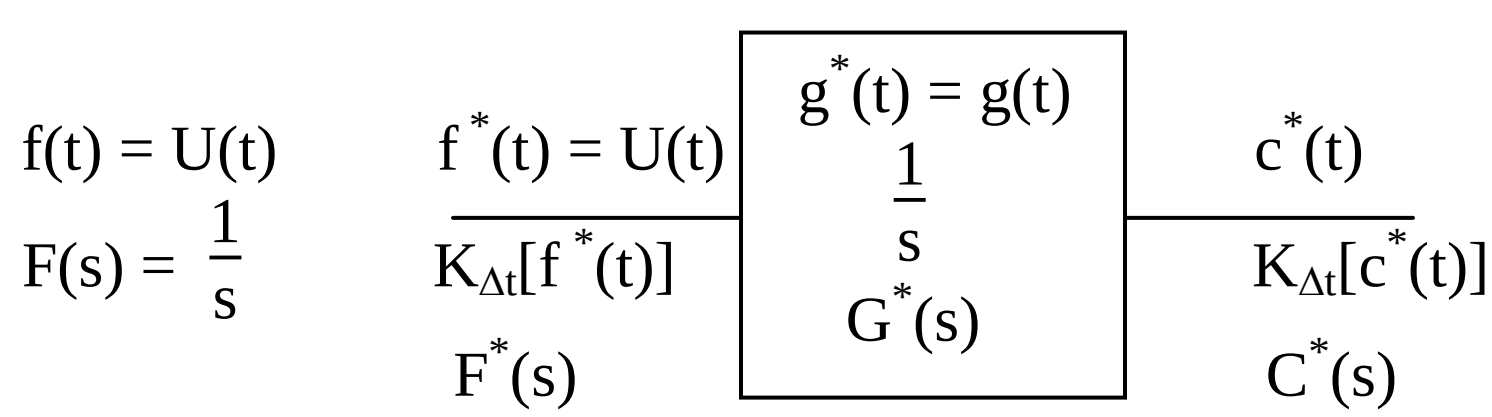

The input and output $K_{\Delta t}$ Transforms of a $K_{\Delta t}$ Transform sample and hold switch are the same.

$$C^*(s) = G^*(s)\, F^*(s) \tag{21}$$

$$F^*(s) = F(s) = \frac{1}{s} \tag{22}$$

$$g(t) = L^{-1}[G(s)] = L^{-1}[\frac{1}{s}] = U(t) \tag{23}$$

$$G^*(s) = K_{\Delta t}[g(t)] = K_{\Delta t}[U(t)] = \frac{1}{s} \tag{24}$$

$$G^*(s) = \frac{1}{s} \tag{25}$$

Substituting Eq 25 and Eq 24 into Eq 21

$$C^*(s) = \frac{1}{s}\,\frac{1}{s} = \frac{1}{s^2} \tag{26}$$

$$C^*(s) = \frac{1}{s^2} \tag{27}$$

From the $K_{\Delta t}$ Transform Table 3 in the Appendix

$$K_{\Delta t}[t] = \frac{1}{s^2} \tag{28}$$

From Eq 27 and Eq 26

$$C^*(s) = K_{\Delta t}[t] \tag{29}$$

Taking the inverse $K_{\Delta t}$ Transform of Eq 29

**$c^*(t) = t$ , $t = 0, \Delta t, 2\Delta t, 3\Delta t \ldots$** 30)

**where**

**$c^*(t)$ is a sample and hold shaped waveform**

Good check

Find the Laplace Transform of $c^*(t)$ from the $K_{\Delta t}$ Transform of $c^*(t)$

Rewriting Eq 5.12-157

$$C^*(s) = \frac{1}{s}\, sC^*(s)\,\Big|_{s = \frac{e^{s\Delta t} - 1}{\Delta t}} \tag{31}$$

From Eq 31 and Eq 27

$$C^*(s) = \frac{1}{s} s \frac{1}{s^2} \Big|_{s = \frac{e^{s\Delta t} - 1}{\Delta t}} = \frac{1}{s}\frac{1}{s} \Big|_{s = \frac{e^{s\Delta t} - 1}{\Delta t}} = \frac{1}{s}\frac{\Delta t}{e^{s\Delta t} - 1} \qquad 32)$$

$$\mathbf{C^*(s) = \frac{\Delta t}{s(e^{s\Delta t} - 1)}} \qquad 33)$$

Eq 33 is the same as Eq 14 where $\Delta t = T$

Good check

Comment

From the previous derivations and demonstrations it is evident that the $K_{\Delta t}$ Transform, the Z Transform, and the Laplace Transform are very closely related.

## Section 5.13: The $e_{\Delta t}(a,t)$ variable subscript identity

**The $e_{\Delta t}(a,t)$ variable subscript identity is defined as follows:**

$$\mathbf{e_{\Delta t}(a,t) = e_{m\Delta t}(\frac{(1+a\Delta t)^m - 1}{m\Delta t},t)\ ,\ m = 1,2,3,\ldots\ ,\ t = 0,\Delta t,2\Delta t,3\Delta t,\ldots} \qquad (5.13\text{-}1)$$

To show that Eq 5.13-1 represents an identity, expand the function, $e_{m\Delta t}(\frac{(1+a\Delta t)^m - 1}{m\Delta t},t)$ .

The definition of the function, $e_{\Delta t}(b,t)$, is:

$$e_{\Delta t}(b,t) = (1+b\Delta t)^{\frac{t}{\Delta t}}\ ,\quad t = 0,\Delta t,2\Delta t,3\Delta t,\ldots \qquad (5.13\text{-}2)$$

$$\Delta t = m\Delta t,\ \ m = 1,2,3,\ldots \qquad (5.13\text{-}3)$$

$$b = \frac{(1+a\Delta t)^m - 1}{m\Delta t} \qquad (5.13\text{-}4)$$

Substituting Eq 5.13-3 and Eq 5.13-4 into Eq 5.13-2

$$e_{m\Delta t}(\frac{(1+a\Delta t)^m - 1}{m\Delta t},t) = [1 + \{\frac{(1+a\Delta t)^m - 1}{m\Delta t}\}m\Delta t]^{\frac{t}{m\Delta t}} = (1+a\Delta t)^{\frac{t}{\Delta t}} = e_{\Delta t}(a,t) \qquad (5.13\text{-}5)$$

An identity is shown.

As an identity, the function, $e_{m\Delta t}(\frac{(1+a\Delta t)^m-1}{m\Delta t},t)$ where m = 1,2,3,… , can replace the function, $e_{\Delta t}(a,t)$, in any equation. This ability, to have an equivalent function for $e_{\Delta t}(a,t)$ but with a different subscript value, can be put to good use in matching some Interval Calculus operations to Interval Calculus functions. For example, $D_2e_1(3,t)$ can not be evaluated using the discrete derivative table, Table 6, in the Appendix. However, its equivalent, $D_2e_2(7.5,t)$, can be evaluated using this discrete derivative table. The formulas in the discrete derivative table in the Appendix are only for discrete derivatives of discrete functions with the same subscripts. Three examples, Example 5.13-1 thru Example 5.13-3, are presented below to demonstate the usefulness of the $e_{\Delta t}(a,t)$ variable subscript identity in matching the subscripts of an Interval Calculus operation/function pair.

Example 5.13-1 Calculation of $D_2e_1$ using the $e_{\Delta t}(a,t)$ variable subscript identity

Find $F(t) = D_2e_1(3,t)$ using the $e_{\Delta t}(a,t)$ variable subscript identity

$F(t) = D_2e_1(3,t)$ 1)

Using the $e_{\Delta t}(a,t)$ variable subscript identity,

$e_{\Delta t}(a,t) = e_{m\Delta t}(\frac{(1+a\Delta t)^m-1}{m\Delta t},t)$ , m = 1,2,3,… , t = 0,Δt,2Δt,3Δt,… 2)

$m = 2$ 3)

$\Delta t = 1$ 4)

$a = 3$ 5)

Substituting Eq 3 thru Eq 5 into Eq 2

$e_1(3,t) = e_2(\frac{(1+3)^2-1}{2},t) = e_2(7.5,t)$ 6)

$e_1(3,t) = e_2(7.5,t)$ 7)

Substituting Eq 7 into Eq 1

$F(t) = D_2e_1(3,t) = D_2e_2(7.5,t)$ 8)

From the discrete derivative table, Table 6, in the Appendix

$D_{\Delta x}e_{\Delta x}(a,x) = ae_{\Delta x}(a,x)$ 9)

From Eq 8 and Eq 9

$F(t) = D_2e_1(3,t) = D_2e_2(7.5,t) = 7.5e_2(7.5,t)$ 10)

From Eq 10 and Eq 7

$$F(t) = D_2e_1(3,t) = D_2e_2(7.5,t) = 7.5e_2(7.5,t) = 7.5e_1(3,t) \qquad 11)$$

Then

$$\mathbf{F(t) = D_2e_1(3,t) = 7.5e_1(3,t)\ ,\ \ t = 0,\Delta t,2\Delta t,3\Delta t,\ldots} \qquad 12)$$

Checking the result of Eq 12 using the definition of the discrete derivative and the function, $e_{\Delta t}(a,t)$

$$D_{2\Delta t}f(t) = \frac{f(t+2\Delta t) - f(t)}{2\Delta t} \qquad 13)$$

$$f(t) = e_{\Delta t}(a,t) = [1+a\Delta t]^{\frac{t}{\Delta t}} \qquad 14)$$

Rewriting the values, m = 2, $\Delta t = 1$, a = 3

From Eq 13 and Eq 14 and substituting the above values

$$F(t) = D_2e_1(3,t) = \frac{e_1(3,t+2) - e_1(3,t)}{2} = \frac{[1+3(1)]^{\frac{t+2}{1}} - [1+3(1)]^{\frac{t}{1}}}{2} \qquad 15)$$

$$F(t) = [\frac{4^2-1}{2}][1+3(1)]^{\frac{t}{1}} = 7.5e_1(3,t) \qquad 16)$$

$$\mathbf{F(t) = D_2e_1(3,t) = 7.5e_1(3,t)\ ,\ \ t = 0,\Delta t,2\Delta t,3\Delta t,\ldots} \qquad 17)$$

Good check

Example 5.13-2 Calculation of $D_{2\Delta t}sin_{\Delta t}(b,t)$ using the $e_{\Delta t}(a,t)$ variable subscript identity

Calculate $F(t) = D_{2\Delta t}sin_{\Delta t}(b,t)$ using the $e_{\Delta t}(a,t)$ variable subscript identity.

$$e_{\Delta t}(a,t) = e_{m\Delta t}(\frac{(1+a\Delta t)^m-1}{m\Delta t},t) \text{ where } m = 1,2,3,\ldots\ ,\ \ t = 0,\Delta t,2\Delta t,3\Delta t,\ldots \qquad 1)$$

$$F(t) = D_{2\Delta t}sin_{\Delta t}(b,t) \qquad 2)$$

$$sin_{\Delta t}(b,t) = \frac{e_{\Delta t}(jb,t) - e_{\Delta t}(-jb,t)}{2j} \qquad 3)$$

Substituting Eq 3 into Eq 2

$$F(t) = D_{2\Delta t}\, sin_{\Delta t}(b,t) = D_{2\Delta t}[\ \frac{e_{\Delta t}(jb,t) - e_{\Delta t}(-jb,t)}{2j}\ ] \qquad 4)$$

Using the $e_{\Delta t}(a,t)$ variable subscript identity, Eq 1

$m = 2$ 5)

$$e_{\Delta t}(jb,t) = e_{2\Delta t}\left(\frac{(1+jb\Delta t)^2-1}{2\Delta t},t\right) \tag{6}$$

Substituting Eq 6 into Eq 4

$$F(t) = D_{2\Delta t}\sin_{\Delta t}(b,t) = D_{2\Delta t}\left[\frac{e_{2\Delta t}\left(\frac{(1+jb\Delta t)^2-1}{2\Delta t},t\right) - e_{2\Delta t}\left(\frac{(1-jb\Delta t)^2-1}{2\Delta t},t\right)}{2j}\right] \tag{7}$$

$$F(t) = D_{2\Delta t}\sin_{\Delta t}(b,t) = \frac{\frac{(1+jb\Delta t)^2-1}{2\Delta t}\,e_{2\Delta t}\left(\frac{(1+jb\Delta t)^2-1}{2\Delta t},t\right) - \frac{(1-jb\Delta t)^2-1}{2\Delta t}\,e_{2\Delta t}\left(\frac{(1-jb\Delta t)^2-1}{2\Delta t},t\right)}{2j} \tag{8}$$

Simplifying Eq 8

$$F(t) = D_{2\Delta t}\sin_{\Delta t}(b,t) = \frac{1}{2j}\left[\frac{j2b\Delta t-b^2\Delta t^2}{2\Delta t}\,e_{2\Delta t}\left(\frac{(1+jb\Delta t)^2-1}{2\Delta t},t\right) - \frac{-j2b\Delta t-b^2\Delta t^2}{2\Delta t}\,e_{2\Delta t}\left(\frac{(1-jb\Delta t)^2-1}{2\Delta t},t\right)\right] \tag{9}$$

$$F(t) = D_{2\Delta t}\sin_{\Delta t}(b,t) = \frac{1}{4j}\left[(j2b - b^2\Delta t)\,e_{2\Delta t}\left(\frac{(1+jb\Delta t)^2-1}{2\Delta t},t\right) + (j2b + b^2\Delta t)\,e_{2\Delta t}\left(\frac{(1-jb\Delta t)^2-1}{2\Delta t},t\right)\right] \tag{10}$$

$$\begin{aligned} F(t) = D_{2\Delta t}\sin_{\Delta t}(b,t) = {} & \frac{1}{2}b\left[e_{2\Delta t}\left(\frac{(1+jb\Delta t)^2-1}{2\Delta t},t\right) + e_{2\Delta t}\left(\frac{(1-jb\Delta t)^2-1}{2\Delta t},t\right)\right] \\ & - \frac{b^2\Delta t}{4j}\left[e_{2\Delta t}\left(\frac{(1+jb\Delta t)^2-1}{2\Delta t},t\right) - e_{2\Delta t}\left(\frac{(1-jb\Delta t)^2-1}{2\Delta t},t\right)\right] \end{aligned} \tag{11}$$

From Eq 11 and Eq 6

$$F(t) = D_{2\Delta t}\sin_{\Delta t}(b,t) = b\left[\frac{e_{\Delta t}(jb,t) + e_{\Delta t}(jb,t)}{2}\right] - \frac{b^2\Delta t}{2}\left[\frac{e_{\Delta t}(jb,t) - e_{\Delta t}(jb,t)}{2j}\right] \tag{12}$$

Then

$$\mathbf{F(t) = D_{2\Delta t}\sin_{\Delta t}(b,t) = b\cos_{\Delta t}(b,t) - \frac{b^2\Delta t}{2}\sin_{\Delta t}(b,t)} \tag{13}$$

Eq 13 is checked by Example 5.14-2 in Section 5.14

Good check

Example 5.13-3 Calculation of ${}_{2\Delta t}\int_{t_1}^{t_2} \sin_{\Delta t}(b,t)\Delta t$ using the $e_{\Delta t}(a,t)$ variable subscript identity

Calculate ${}_{2\Delta t}\int_{t_1}^{t_2} \sin_{\Delta t}(b,t)\Delta t$ using the $e_{\Delta t}(a,t)$ variable subscript identity.

The $e_{\Delta t}(a,t)$ variable subscript identity is:

$$e_{\Delta t}(a,t) = e_{m\Delta t}\left(\frac{(1+a\Delta t)^m-1}{m\Delta t},t\right) \text{ where } m = 1,2,3,\ldots ,\quad t = 0,\Delta t,2\Delta t,3\Delta t,\ldots \qquad 1)$$

Finding a $2\Delta t$ subscript identity for $\sin_{\Delta t}(b,t)$

$$\sin_{\Delta t}(b,t) = \frac{e_{\Delta t}(jb,t) - e_{\Delta t}(-jb,t)}{2j} \qquad 2)$$

Substituting Eq 1 into Eq 2 where m = 2

$$\sin_{\Delta t}(b,t) = \left[\frac{e_{2\Delta t}\left(\frac{(1+jb\Delta t)^2-1}{2\Delta t},t\right) - e_{2\Delta t}\left(\frac{(1-jb\Delta t)^2-1}{2\Delta t},t\right)}{2j}\right] \qquad 3)$$

Integrating Eq 3

$${}_{2\Delta t}\int_{t_1}^{t_2} \sin_{\Delta t}(b,t)\Delta t = \frac{1}{2}\,{}_{2\Delta t}\int_{t_1}^{t_2} \sin_{\Delta t}(b,t)2\Delta t = \frac{1}{4j}\left[{}_{2\Delta t}\int_{t_1}^{t_2} e_{2\Delta t}\left(\frac{(1+jb\Delta t)^2-1}{2\Delta t},t\right)(2\Delta t) - {}_{2\Delta t}\int_{t_1}^{t_2} e_{2\Delta t}\left(\frac{(1-jb\Delta t)^2-1}{2\Delta t},t\right)(2\Delta t)\right] \qquad 4)$$

$${}_{2\Delta t}\int_{t_1}^{t_2} \sin_{\Delta t}(b,t)\Delta t = \frac{1}{4j}\left[\frac{1}{\frac{(1+jb\Delta t)^2-1}{2\Delta t}} e_{2\Delta t}\left(\frac{(1+jb\Delta t)^2-1}{2\Delta t},t\right) - \frac{1}{\frac{(1-jb\Delta t)^2-1}{2\Delta t}} e_{2\Delta t}\left(\frac{(1-jb\Delta t)^2-1}{2\Delta t},t\right)\right]\Big|_{t_1}^{t_2} \qquad 5)$$

Substituting Eq 1 into Eq 5 and simplifying

$${}_{2\Delta t}\int_{t_1}^{t_2} \sin_{\Delta t}(b,t)\Delta t = \frac{1}{4j}\left[\frac{2\Delta t}{(1+jb\Delta t)^2-1} e_{2\Delta t}(jb,t) - \frac{2\Delta t}{(1-jb\Delta t)^2-1} e_{2\Delta t}(-jb,t)\right]\Big|_{t_1}^{t_2} \qquad 6)$$

$${}_{2\Delta t}\int_{t_1}^{t_2} \sin_{\Delta t}(b,t)\Delta t = \frac{1}{4j}\left[\frac{2\Delta t}{2jb\Delta t-b^2\Delta t^2} e_{2\Delta t}(jb,t) - \frac{2\Delta t}{-2jb\Delta t-b^2\Delta t^2} e_{2\Delta t}(-jb,t)\right]\Big|_{t_1}^{t_2} \qquad 7)$$

Canceling terms in Eq 7 and simplifying

$$\underset{2\Delta t}{\int_{t_1}^{t_2}} \sin_{\Delta t}(b,t)\Delta t = \frac{1}{2}\left[\frac{1}{-2b-jb^2\Delta t} e_{2\Delta t}(jb,t) - \frac{1}{2b-jb^2\Delta t} e_{2\Delta t}(-jb,t)\right]\Big|_{t_1}^{t_2} \qquad 8)$$

$$\underset{2\Delta t}{\int_{t_1}^{t_2}} \sin_{\Delta t}(b,t)\Delta t = \frac{1}{2}\left[\frac{-2b+jb^2\Delta t}{4b^2+b^4\Delta t^2} e_{2\Delta t}(jb,t) - \frac{2b+jb^2\Delta t}{4b^2+b^4\Delta t^2} e_{2\Delta t}(-jb,t)\right]\Big|_{t_1}^{t_2} \qquad 9)$$

$$\underset{2\Delta t}{\int_{t_1}^{t_2}} \sin_{\Delta t}(b,t)\Delta t = \frac{1}{2}\left[\frac{-2+jb\Delta t}{4b+b^3\Delta t^2} e_{2\Delta t}(jb,t) + \frac{-2-jb\Delta t}{4b+b^3\Delta t^2} e_{2\Delta t}(-jb,t)\right]\Big|_{t_1}^{t_2} \qquad 10)$$

$$\underset{2\Delta t}{\int_{t_1}^{t_2}} \sin_{\Delta t}(b,t)\Delta t = \frac{-2}{4b+b^3\Delta t^2}\,\frac{e_{2\Delta t}(jb,t) + e_{\Delta t}(-jb,t)}{2}\Big|_{t_1}^{t_2} + \frac{-b\Delta t}{4b+b^3\Delta t^2}\,\frac{e_{2\Delta t}(jb,t) - e_{\Delta t}(-jb,t)}{2j}\Big|_{t_1}^{t_2} \qquad 11)$$

Then

$$\boldsymbol{\underset{2\Delta t}{\int_{t_1}^{t_2}} \sin_{\Delta t}(b,t)\Delta t = -\frac{2}{4b+b^3\Delta t^2}\cos_{\Delta t}(b,t)\Big|_{t_1}^{t_2} - \frac{b\Delta t}{4b+b^3\Delta t^2}\sin_{\Delta t}(b,t)\Big|_{t_1}^{t_2}, \quad t = 0, 2\Delta t, 4\Delta t, 6\Delta t, \ldots} \qquad 12)$$

<u>Note</u> – The integration operation establishes the discrete time as being $t = 0, 2\Delta t, 4\Delta t, 6\Delta t, \ldots$

Checking Eq 12

Let
$\Delta t = .5$
$b = .8$
$t_1 = 0$
$t_2 = 3$

Substituting the above values into Eq 12

$$\underset{1}{\int_{0}^{3}} \sin_{.5}(.8,t)\Delta t = -\frac{2}{4(.8)+(.8)^3(.5)^2}\cos_{.5}(.8,t) - \frac{.8(.5)}{4(.8)+(.8)^3(.5)^2}\sin_{.5}(.8,t)\Big|_{0}^{3} \qquad 13)$$

$$\underset{1}{\int_{0}^{3}} \sin_{.5}(.8,t)\Delta t = -.600961538\cos_{.5}(.8,t) - .120192307\sin_{.5}(.8,t)\Big|_{0}^{3} \qquad 14)$$

$$\int_{0}^{3} {}_1\sin_{.5}(.8,t)\Delta t = -.600961538\,[\cos_{.5}(.8,3) - \cos_{.5}(.8,0)] - .120192307\,[\sin_{.5}(.8,3) - \sin_{.5}(.8,0)] \quad 15)$$

Using the DXFUN1 program to calculate the $\sin_{\Delta t}(b,t)$ and $\cos_{\Delta t}(b,t)$ functions

$$\int_{0}^{3} {}_1\sin_{.5}(.8,t)\Delta t = -.600961538\,[-1.020096 - 1] - .120192307\,[1.18145 - 0] \quad 16)$$

$$\int_{0}^{3} {}_1\sin_{.5}(.8,t)\Delta t = 1.072 \quad 17)$$

Verifying Eq 17

Finding $\int_{0}^{3} {}_1\sin_{.5}(.8,t)\Delta t$ using a series where $\Delta t = .5$ and $2\Delta t = 1$

$$\int_{0}^{3} {}_1\sin_{.5}(.8,t)\Delta t = {}_1\sum_{t=0}^{2} \sin_{.5}(.8,t)(.5) = .5\,{}_1\sum_{t=0}^{2} \sin_{.5}(.8,t) = .5[\sin_{.5}(.8,0) + \sin_{.5}(.8,1) + \sin_{.5}(.8,2)] \quad 18)$$

Using the DXFUN1 program to calculate the $\sin_{\Delta t}(b,t)$ function

$$\int_{0}^{3} {}_1\sin_{.5}(.8,t)\Delta t = .5[0 + .8 + 1.344] = 1.072 \quad 19)$$

$$\int_{0}^{3} {}_1\sin_{.5}(.8,t)\Delta t = .5[0 + .8 + 1.344] = 1.072 \quad 20)$$

Good check

<u>Comment</u> – The following two integrals are derived using the same $e_{\Delta t}$ variable subscript identitiy methodology that has been used above.

1. $${}_{2\Delta x}\int e_{\Delta x}(a,x)\Delta x = \frac{1}{2a+a^2\Delta t}\, e_{\Delta x}(a,x) + k \quad , \quad x = 0, 2\Delta x, 4\Delta x, 6\Delta x, \ldots$$

2. $${}_{2\Delta x}\int \cos_{\Delta x}(b,x)\Delta x = \frac{2}{4b+b^3\Delta x^2}\sin_{\Delta x}(b,x) - \frac{b\Delta x}{4b+b^3\Delta x^2}\cos_{\Delta x}(b,x) + k \; , \; x = 0, 2\Delta x, 4\Delta x, 6\Delta x, \ldots$$

## Section 5.14: Discrete derivative equalities

Besides the use of the $e_{\Delta t}(a,t)$ variable subscript identity described in Section 5.13, there is another way to evaluate discrete derivatives of functions where both subscripts are not the same. This method uses a discrete derivative equality. Discrete derivative equalities make it possible express a $D^n_{m\Delta t}$ discrete derivative in terms of a series of $D^p_{\Delta t}$ discrete derivatives (p=1,2,3,…) where the derivative subscripts are the same as the function subscript. Then the discrete derivative table in the Appendix, Table 6, can be used for evaluation. The formulas in the discrete derivative table are only for discrete derivatives of discrete functions where both have the same subscripts.

Below are listed some derived discrete derivative equalities that can be used to calculate derivatives of functions where the discrete time values of the derivative are a subset of those of the function.

1. $D_{2\Delta t} = D_{\Delta t} + \frac{1}{2}\Delta t\, D_{\Delta t}^{\,2}$

2. $D_{3\Delta t} = D_{\Delta t} + \Delta t D_{\Delta t}^{\,2} + \frac{1}{3}\Delta t^2\, D_{\Delta t}^{\,3}$

3. $D_{4\Delta t} = D_{\Delta t} + \frac{3}{2}\Delta t D_{\Delta t}^{\,2} + \Delta t^2 D_{\Delta t}^{\,3} + \frac{1}{4}\Delta t^3 D_{\Delta t}^{\,4}$

4. $$D_{m\Delta t} = D_{\Delta t} + \sum_{n=2}^{\infty} \frac{\Delta t^{n-1}}{n!}\prod_{p=1}^{n-1}(m-p) D_{\Delta t}^{\,n} = D_{\Delta t} + \frac{(m-1)}{2!}\Delta t D_{\Delta t}^{\,2} + \frac{(m-1)(m-2)}{3!}\Delta t^2 D_{\Delta t}^{\,3} + \frac{(m-1)(m-2)(m-3)}{4!}\Delta t^3 D_{\Delta t}^{\,4} + \ldots \quad , \quad m = 2,3,4,\ldots$$

5. $D_{2\Delta t}^{\,2} = D_{\Delta t}^{\,2} + \Delta t D_{\Delta x}^{\,3} + \frac{1}{4}\Delta t^2\, D_{\Delta t}^{\,4}$

where

$$D_{\Delta t}f(t) = \frac{f(t+\Delta t) - f(t)}{\Delta t} \quad , \quad t = 0, \Delta t, 2\Delta t, 3\Delta t, 4\Delta t, \ldots$$

$$D_{2\Delta t}f(t) = \frac{f(t+2\Delta t) - f(t)}{2\Delta t} \quad , \quad t = 0, 2\Delta t, 4\Delta t, 6\Delta t, 8\Delta t, \ldots$$

$$D_{3\Delta t}f(t) = \frac{f(t+3\Delta t) - f(t)}{3\Delta t} \quad , \quad t = 0, 3\Delta t, 6\Delta t, 9\Delta t, 12\Delta t, \ldots$$

$$D_{4\Delta t}f(t) = \frac{f(t+4\Delta t) - f(t)}{4\Delta t} \quad , \quad t = 0, 4\Delta t, 8\Delta t, 12\Delta t, 16\Delta t, \ldots$$

$$D_{2\Delta t}^{\,2}f(t) = \frac{f(t+4\Delta t) - 2f(t+2\Delta t) + f(t)}{(2\Delta t)^2} \quad , \quad t = 0, 2\Delta t, 4\Delta t, 6\Delta t, 8\Delta t, \ldots$$

Comment - $D_{m\Delta t}^{\,n}$ is alternately represented as $s_{m\Delta t}^{\,n}$

Two derivations are shown below. In Example 5.14-1 a derivation of the discrete derivative equality, $D_{2\Delta x} = D_{\Delta x} + \frac{\Delta x}{2} D_{\Delta x}^2$, is shown and in Example 5.14-2 the derivation of the discrete derivative equality,

$$D_{m\Delta t} = D_{\Delta t} + \sum_{n=2}^{\infty} \frac{\Delta t^{n-1}}{n!} \prod_{p=1}^{n-1} (m-p) D_{\Delta t}^{n}$$ where m = 2,3,4, … , is shown.

<u>Example 5.14-1</u> Derivation of the discrete derivative equality, $D_{2\Delta x} = D_{\Delta x} + \frac{\Delta x}{2} D_{\Delta x}^2$

Derive the discrete derivative equality, $D_{2\Delta t} = D_{\Delta t} + \frac{\Delta t}{2} D_{\Delta t}^2$

$$D_{2\Delta t}y = \frac{y_{2\Delta t} - y_0}{2\Delta t} \tag{1}$$

$$D_{2\Delta t}y = \frac{\Delta t}{2}\left[\frac{y_{2\Delta t} - y_0}{\Delta t^2}\right] \tag{2}$$

$$D^2_{\Delta t}y = \frac{y_{2\Delta t} - 2y_{\Delta t} + y_0}{\Delta t^2} \tag{3}$$

$$D_{\Delta t}y = \frac{y_{\Delta t} - y_0}{\Delta t} \tag{4}$$

From Eq 2 thru Eq 4

$$D_{2\Delta t}y = \frac{\Delta t}{2}\left[\frac{y_{2\Delta t} - 2y_{\Delta t} + 2y_0 - y_0}{\Delta t^2}\right] + \frac{\Delta t}{2}\left[\frac{2y_{\Delta t} - 2y_0}{\Delta t^2}\right] = \frac{\Delta t}{2} D_{\Delta x}^2 y + D_{\Delta t}y \tag{5}$$

$$D_{2\Delta t}y = D_{\Delta t}\, y + \frac{\Delta t}{2} D_{\Delta x}^2 y \tag{6}$$

Then

**The $D_{2\Delta t}$ discrete derivative equality is:**

$$\mathbf{D_{2\Delta t} = D_{\Delta t} + \frac{\Delta t}{2} D_{\Delta x}^2} \tag{7}$$

Example 5.14-2 Derivation of the discrete derivative equality, $D_{m\Delta t} = D_{\Delta t} + \sum_{n=2}^{\infty} \frac{\Delta t^{n-1}}{n!} \prod_{p=1}^{n-1} (m-p) D_{\Delta t}^{n}$

where $m = 2,3,4, \ldots$

Derive the discrete derivative equality, $D_{m\Delta t} = D_{\Delta t} + \sum_{n=2}^{\infty} \frac{\Delta t^{n-1}}{n!} \prod_{p=1}^{n-1} (m-p) D_{\Delta t}^{n}$ where $m = 2,3,4, \ldots$

Using the Discrete Variable Maclaurin Series

$$f(t) = \frac{1}{0!} f(0) + \frac{1}{1!} D_{\Delta t} f(0)t + \frac{1}{2!} D_{\Delta t}^{2} f(0)t(t-\Delta t) + \frac{1}{3!} D_{\Delta t}^{3} f(0)\, t(t-\Delta t)(t-2\Delta t) + \ldots \qquad 1)$$

Let

$$t = m\Delta t \ , \quad m = 2,3,4, \ldots \qquad 2)$$

$$\frac{f(mt) - f(0)}{m\Delta t} = D_{\Delta t} f(0) + \frac{1}{2!} D_{\Delta t}^{2} f(0)(m-1)\Delta t + \frac{1}{3!} D_{\Delta t}^{3} f(0)(m-1)(m-2)\Delta t^{2} + \ldots \qquad 3)$$

$$D_{m\Delta t} f(0) = D_{\Delta t} f(0) + \frac{(m-1)}{2!} \Delta t\, D_{\Delta t}^{2} f(0) + \frac{(m-1)(m-2)}{3!} \Delta t^{2} D_{\Delta t}^{3} f(0) + \ldots \qquad 4)$$

From Eq 4 write the discrete derivative equality

$$D_{m\Delta t} = D_{\Delta t} + \frac{(m-1)}{2!} \Delta t\, D_{\Delta t}^{2} + \frac{(m-1)(m-2)}{3!} \Delta t^{2} D_{\Delta t}^{3} + \ldots \qquad 5)$$

Then

Expanding Eq 5

**The $D_{m\Delta t}$ discrete derivative equality is:**

$$\mathbf{D_{m\Delta t} = D_{\Delta t} + \sum_{n=2}^{\infty} \frac{\Delta t^{n-1}}{n!} \prod_{p=1}^{n-1} (m-p) D_{\Delta t}^{n} \ , \quad m = 2,3,4, \ldots} \qquad 6)$$

Discrete derivative equalities such as those derived in Example 5.14-1 and Example 5.14-2 above are used very efficiently to calculate derivatives of functions where the discrete time values of the derivative are a subset of those of the function. An example of the use of a discrete derivative equality in the solution of a problem is shown in the following example, Example 5.14-3.

Example 5.14-3 Calculation of $F(t) = D_{2\Delta t}\sin_{\Delta t}(b,t)$ using a discrete derivative equality

Calculate $F(t) = D_{2\Delta t}\sin_{\Delta t}(b,t)$ using the discrete derivative equality, $D_{2\Delta t} = D_{\Delta t} + \frac{\Delta t}{2} D_{\Delta t}^2$

$$F(t) = D_{2\Delta t}\sin_{\Delta t}(b,t) \qquad 1)$$

$$D_{2\Delta t} = D_{\Delta t} + \frac{\Delta t}{2} D_{\Delta t}^2 \qquad 2)$$

Substituting Eq 2 into Eq 1

$$F(t) = (D_{\Delta t} + \frac{\Delta t}{2} D_{\Delta t}^2)\sin_{\Delta t}(b,t) \qquad 3)$$

Differentiating

$$F(t) = b\cos_{\Delta t}(b,t) - \frac{\Delta t}{2} b^2\sin_{\Delta t}(b,t) \qquad 4)$$

Then

$$\mathbf{D_{2\Delta t}\sin_{\Delta t}(b,t) = b\cos_{\Delta t}(b,t) - \frac{\Delta t}{2} b^2\sin_{\Delta t}(b,t)} \qquad 5)$$

Checking Eq 5 using the definition of the discrete derivative

$$F(t) = D_{2\Delta t}\sin_{\Delta t}(b,t) \qquad 6)$$

$$D_{2\Delta t}f(t) = \frac{f(t+2\Delta t) - f(t)}{2\Delta t} \qquad 7)$$

From Eq 6 and Eq 7

$$F(t) = D_{2\Delta t}\sin_{\Delta t}(b,t) = \frac{\sin_{\Delta t}(b,t+2\Delta t) - \sin_{\Delta t}(b,t)}{2\Delta t} \qquad 8)$$

Using an identity from Table 5 in the Appendix

$$\sin_{\Delta x}(a,x+y) = \sin_{\Delta x}(a,x)\cos_{\Delta x}(a,y) + \sin_{\Delta x}(a,y)\cos_{\Delta x}(a,x) \qquad 9)$$

$$a = b \qquad 10)$$

$$x = t \qquad 11)$$

$$y = 2\Delta t \qquad 12)$$

Substituting Eq 10 thru Eq 12 into Eq 9

$$\sin_{\Delta t}(b,t+2\Delta t) = \sin_{\Delta t}(b,t)\cos_{\Delta t}(b,2\Delta t) + \sin_{\Delta t}(b,2\Delta t)\cos_{\Delta t}(b,t) \qquad 13)$$

$$\sin_{\Delta t}(b,2\Delta t) = \frac{(1+jb\Delta t)^{\frac{2\Delta t}{\Delta t}} - (1-jb\Delta t)^{\frac{2\Delta t}{\Delta t}}}{2j} = \frac{(1+jb\Delta t)^2 - (1-jb\Delta t)^2}{2j} = \frac{1+2jb\Delta t-b^2\Delta t^2-1+2jb\Delta t+b^2\Delta t^2}{2j} \qquad 14)$$

$$\sin_{\Delta t}(b,2\Delta t) = 2b\Delta t \qquad 15)$$

$$\cos_{\Delta t}(b,2\Delta t) = \frac{(1+jb\Delta t)^{\frac{2\Delta t}{\Delta t}} + (1-jb\Delta t)^{\frac{2\Delta t}{\Delta t}}}{2} = \frac{(1+jb\Delta t)^2 + (1-jb\Delta t)^2}{2} = \frac{1+2jb\Delta t-b^2\Delta t^2+1-2jb\Delta t-b^2\Delta t^2}{2} \qquad 16)$$

$$\cos_{\Delta t}(b,2\Delta t) = 1-b^2\Delta t^2 \qquad 17)$$

Substituting Eq 15 and Eq 17 into Eq 13

$$\sin_{\Delta t}(b,t+2\Delta t) = (1-b^2\Delta t^2)\sin_{\Delta t}(b,t) + (2b\Delta t)\cos_{\Delta t}(b,t) \qquad 18)$$

Substituting Eq 18 into Eq 8

$$F(t) = D_{2\Delta t}\sin_{\Delta t}(b,t) = \frac{(1-b^2\Delta t^2)\sin_{\Delta t}(b,t) + (2b\Delta t)\cos_{\Delta t}(b,t) - \sin_{\Delta t}(b,t)}{2\Delta t} \qquad 19)$$

$$F(t) = b\cos_{\Delta t}(b,t) - \frac{\Delta t}{2}b^2\sin_{\Delta t}(b,t) \qquad 20)$$

Then

$$D_{2\Delta t}\sin_{\Delta t}(b,t) = b\cos_{\Delta t}(b,t) - \frac{\Delta t}{2}b^2\sin_{\Delta t}(b,t) \qquad 21)$$

Good check

## Section 5.15: Demonstration of the use of discrete derivative equalities in the analysis of differential difference equations and sampled-data systems

The use of discrete derivative equations in the analysis of differential difference equations

Discrete derivative equalities can be used in the solution of differential difference equations where derivative operators such as $D_{2\Delta t}$ and $D_{3\Delta t}^{2}$ appear. A demonstration of the solution of such a differential difference equation is shown in the following example, Example 5.15-1.

Example 5.15-1 A demonstration of the use of discrete derivative equalities in the solution of differential difference equations

Solve the following related differential difference equations, $D_{\Delta t}y(t) + 2y(t) = 0$ where $\Delta t = 2$, $t = 0, 2, 4, 6, \dots$ and $D_{2\Delta t}y(t) + 2y(t) = 0$ where $\Delta t = 1$, $t = 0, 1, 2, 3, \dots$ using $K_{\Delta t}$ Transforms. $y(0) = 2$.

1. Solve $D_{\Delta t}y(x) + 2y(t) = 0$ where $\Delta t = 2$

$$D_2y(t) + 2y(t) = 0 \tag{1}$$

Taking the $K_{\Delta t}$ Transform of Eq 1

$$sY(s) - y(0) + 2Y(s) = 0 \tag{2}$$

$$y(0) = 2 \tag{3}$$

From Eq 2 and Eq 3

$$Y(s) = \frac{2}{s+2} \tag{4}$$

Finding the Inverse $K_{\Delta t}$ Transform of Eq 4

$$K_{\Delta t}^{-1}\left[\frac{1}{s+a}\right] = [1 - a\Delta t]^{\frac{t}{\Delta t}} \tag{5}$$

$$\Delta t = 2 \tag{6}$$

From Eq 4 thru Eq 6

$$y(t) = 2[1 - 2(2)]^{\frac{t}{2}} \tag{7}$$

Then

$$\mathbf{y(t) = 2[-3]^{\frac{t}{2}}, \quad t = 0, 2, 4, 6, \dots} \tag{8}$$

From Eq 8

Table #1

| t | y(t) |
|---|---|
| 0 | 2 |
| 2 | -6 |
| 4 | 18 |
| 6 | -54 |
| 8 | 162 |
| … | |

2. Solve $D_{2\Delta t}y(x) + 2y(t) = 0$ where $\Delta t = 1$

$D_{2\Delta t}y(t) + 2y(t) = 0 \quad , \quad t = 0, \Delta t, 2\Delta t, 3\Delta t, \ldots$ 9)

$D_{2\Delta t}\, y(t) = \frac{\Delta t}{2} D_{\Delta t}{}^2 y(t) + D_{\Delta t}\, y(t)$ , A discrete derivative identity 10)

Substituting Eq 10 into Eq 9

$\frac{\Delta t}{2} D_{\Delta t}{}^2 y(t) + D_{\Delta t}\, y(t) + 2y(t) = 0$ 11)

Taking the $K_{\Delta t}$ Transform of Eq 11

$\frac{\Delta t}{2}[s^2Y(s) - sy(0) - D_{\Delta t}y(0)] + sY(s) - y(0) + 2Y(s) = 0$ 11)

$\Delta t = 1$ 12)

$y(0) = 2$ 13)

$D_{\Delta t}y(0) = 0$ 14)

Substituting Eq 12 thru Eq 14 into Eq 11

$\frac{1}{2} s^2Y(s) - s + sY(s) - 2 + 2Y(s) = 0$ 15)

Simplifying

$(\frac{1}{2} s^2 + s + 2)Y(s) = s + 2$ 16)

$Y(s) = \frac{2s + 4}{s^2 + 2s + 4} = \frac{2s+4}{(s + 1+j\sqrt{3})(s + 1-j\sqrt{3})} = \frac{A}{s + 1+j\sqrt{3}} + \frac{B}{s + 1-j\sqrt{3}}$ 17)

Finding the Inverse $K_{\Delta t}$ Transform of Eq 17

$$A = \frac{2s+4}{s + 1-j\sqrt{3}}\Big|_{s = -1-j\sqrt{3}} = \frac{-2-2j\sqrt{3} + 4}{-2j\sqrt{3}} = \frac{1-j\sqrt{3}}{-j\sqrt{3}} = 1 + \frac{j}{\sqrt{3}} = 1 + \frac{j\sqrt{3}}{3} \quad 18)$$

$$B = \frac{2s+4}{s + 1+j\sqrt{3}}\Big|_{s = -1+j\sqrt{3}} = \frac{-2+2j\sqrt{3} + 4}{2j\sqrt{3}} = \frac{1+j\sqrt{3}}{j\sqrt{3}} = 1 - \frac{j}{\sqrt{3}} = 1 - \frac{j\sqrt{3}}{3} \quad 19)$$

Substituting Eq 18 and Eq 19 into Eq 17

$$Y(s) = \frac{1 + \frac{j\sqrt{3}}{3}}{s + 1+j\sqrt{3}} + \frac{1 - \frac{j\sqrt{3}}{3}}{s + 1-j\sqrt{3}} \quad 20)$$

Finding the Inverse $K_{\Delta t}$ Transform of Eq 20 using Eq 5

$$y(t) = [1 + \frac{j\sqrt{3}}{3}][1 - (1+j\sqrt{3})(1)]^{\frac{t}{1}} + [1 - \frac{j\sqrt{3}}{3}][1 - (1-j\sqrt{3})(1)]^{\frac{t}{1}} \quad 21)$$

Then

$$\mathbf{y(t) = [1 + \frac{j\sqrt{3}}{3}][-j\sqrt{3})]^{t} + [1 - \frac{j\sqrt{3}}{3}][\, j\sqrt{3}]^{t}, \quad t = 0, 2, 3, \ldots} \quad 22)$$

or

Simplifying Eq 22

$$y(t) = [(-1)^t + \frac{j\sqrt{3}}{3}(-1)^t + 1 - \frac{j\sqrt{3}}{3}](\, j\sqrt{3})^t \quad 23)$$

$$\mathbf{y(t) = \{ [(-1)^t + 1] + \frac{j\sqrt{3}}{3}[(-1)^t - 1] \}(\, j\sqrt{3})^t, \quad t = 0, 1, 2, 3, \ldots} \quad 24)$$

Evaluating Eq 24 for t = 0, 1, 2, 3, 4, 5, 6, 7, 8

$$y(0) = 2(1) = 2 \quad 25)$$

$$y(1) = -2\frac{j\sqrt{3}}{3}(\, j\sqrt{3}) = 2 \quad 26)$$

$$y(2) = 2(\, j\sqrt{3})^2 = -6 \quad 27)$$

$$y(3) = -2\frac{j\sqrt{3}}{3}(\, j\sqrt{3})^3 = -6 \quad 28)$$

$y(4) = ( 2( j\sqrt{3})^4 = 18$ 29)

$y(5) = -2 \frac{j\sqrt{3}}{3} ( j\sqrt{3})^5 = 18$ 30)

$y(6) = 2( j\sqrt{3})^6 = -54$ 31)

$y(7) = ( -2 \frac{j\sqrt{3}}{3} ( j\sqrt{3})^7 = -54$ 32)

$y(8) = 2( j\sqrt{3})^8 = 162$ 33)

Then

<u>Table #2</u>

| t | y(t) |
|---|---|
| 0 | 2 |
| 1 | 2 |
| 2 | -6 |
| 2 | -6 |
| 4 | 18 |
| 5 | 18 |
| 6 | -54 |
| 7 | -54 |
| 8 | 162 |
| … | |

<u>Comment</u> - As expected, the values of y(t) in Solutions 1 and 2 are the same for t = 0, 2, 4, 6, 8.

Checking Solution 2 using Z Transforms

From Eq 17

$$Y(s) = \frac{2s + 4}{s^2 + 2s + 4} \qquad 34)$$

Converting Y(s) into its equivalent Z Transform using the $K_{\Delta t}$ Transform to Z Transform Conversion Equation

$$Y(z) = \frac{z}{T} Y(s)\Big|_{\substack{s = \frac{z-1}{T} \\ T = \Delta t}} , \qquad Z[y(t)] = Y(z) \quad \text{Z Transform} \qquad 35)$$

$$K_{\Delta t}[f(t)] = F(s) \quad K_{\Delta t} \text{ Transform}$$

$$Y(z) = \frac{z}{T}\frac{2s + 4}{s^2 + 2s + 4}\Big|_{\substack{s = \frac{z-1}{T} \\ T = \Delta t}} \qquad 36)$$

$$\Delta t = T = 1 \qquad 37)$$

$$Y(z) = \frac{z[2(z-1) + 4]}{(z-1)^2 + 2(z-1) + 4} = \frac{z[2z + 2]}{z^2 - 2z + 1 + 2z - 2 + 4} = \frac{2z^2 + 2z}{z^2 + 3} \qquad 38)$$

$$Y(z) = \frac{2z^2 + 2z}{z^2 + 3} \qquad 39)$$

Dividing the denominator of Y(z) into the numerator of Y(z)

$$Y(z) = 2z^0 + 2z^{-1} - 6z^{-2} - 6z^{-3} + 18z^{-4} + 18z^{-5} - 54z^{-6} - 54z^{-7} + 162z^{-8} + 162z^{-9} - \ldots \qquad 40)$$

From the above Z Transform Series Table #3 is obtained

Table #3

| t | y(t) |
|---|---|
| 0 | 2 |
| 1 | 2 |
| 2 | -6 |
| 2 | -6 |
| 4 | 18 |
| 5 | 18 |
| 6 | -54 |
| 7 | -54 |
| 8 | 162 |
| … | |

Note that Table #3 is the same as Table #2 for t = 0, 2, 4, 6, 8, … However, unlike Table #2, Table #3 defines the values of y(t) for t = 1, 3, 5, 7, …

Good check

The use of discrete derivative equalities in the analysis of sampled-data systems

Discrete derivative equalities can be used in the analysis of sampled-data systems where there is an operator mismatch. For example, consider the system $K_{\Delta t}$ Transform, $C(s) = \frac{1}{s s_{2\Delta t}}$. This transform has an operator mismatch. The s operator is defined as $sc(t) = \frac{c(t+\Delta t) - c(t)}{\Delta t}$, a difference of the function, c(t), over a Δt interval. However, the $s_{2\Delta t}$ operator is defined differently, it is defined as $s_{2\Delta t}c(t) = \frac{c(t+2\Delta t) - c(t)}{2\Delta t}$, a difference of the function, c(t), over a 2Δt interval. Taking the inverse of a $K_{\Delta t}$ Transform that includes dissimilar operators is facilitated by the use of discrete derivative equalities. Using the discrete derivative equality, $s_{2\Delta t} = s + \frac{\Delta t}{2} s^2$, C(s) becomes $C(s) = \frac{1}{\frac{\Delta t}{2} s^3 + s^2}$. Taking the inverse of this $K_{\Delta t}$ Transform is not difficult.

Discrete derivative equalities may be helpful in analyzing single clock sampled-data systems that have samplers with different sampling periods all derived from the same single clock. In systems such as this, an operator such as $s_{2\Delta t}$, $s_{3\Delta t}$, $s_{\Delta t}^{\,2}$, etc. could appear in a system transform together with s operators. An s difference operator would be associated with a sampler with a sample and hold period of Δt, an $s_{2\Delta t}$ difference operator would be associated with a sampler with a sample and hold period of 2Δt, etc. A discrete derivative equality can be used to convert the system transform into an equivalent transform that is a function of s only. Once this is done, mathematical analysis proceeds with little difficulty. To demonstrate the use discrete derivative equalities in sampled-data system analysis, the following three examples, Example 5.15-2 thru Example 5.15-4, are provided.

In each of the following three demonstration examples, the sampled-data system has samplers with different sampling periods.

Example 5.14-2 A demonstration of the analysis of a single clock sampled-data system with two sample and hold samplers where the period of one sampler is twice that of the other

Demonstrate the analysis of a sampled-data system with two samplers with different sampling periods. One sampling period is twice that of the other. Use the discrete derivative identity, $s_{2\Delta t} = s + \frac{\Delta t}{2} s^2$. The system is initially passive. Find c(t) at t = 0, Δt, 2Δt, 3Δt, 4Δt, 5Δt, 6Δt, 7Δt and 8Δt.

The above operators are defined as follows:

$$sy(t) = \frac{y(t+\Delta t) - y(t)}{\Delta t}$$

$$s^2y(t) = \frac{y(t+2\Delta t) - 2y(t+\Delta t) + y(t)}{\Delta t^2}$$

$$s_{2\Delta t}\, y(t) = sy(t) + \frac{\Delta t}{2} s^2 y(t) = \frac{y(t+2\Delta t) - y(t)}{2\Delta t}$$

where

y(t) = a function of t

t = 0, Δt, 2Δt, 3Δt, …

Consider the following sampled-data system, an integrator with a unit step input.

<u>Fig #1</u>

Integrator

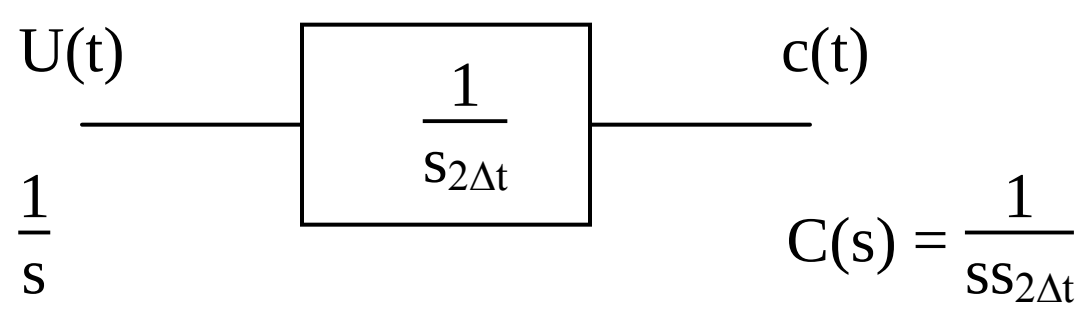

The system is initially passive

$$s_{2\Delta t} = s + \frac{\Delta t}{2} s^2 \tag{1}$$

$$C(s) = \frac{1}{s s_{2\Delta t}} \tag{2}$$

Substituting Eq 1 into Eq 2

$$C(s) = \frac{1}{s s_{2\Delta t}} = \frac{1}{s(\frac{\Delta t}{2} s^2 + s)} = \frac{\frac{2}{\Delta t}}{s^2(s + \frac{2}{\Delta t})} \tag{3}$$

From Fig #1, Eq 1, and Eq 2

<u>Fig #2</u>

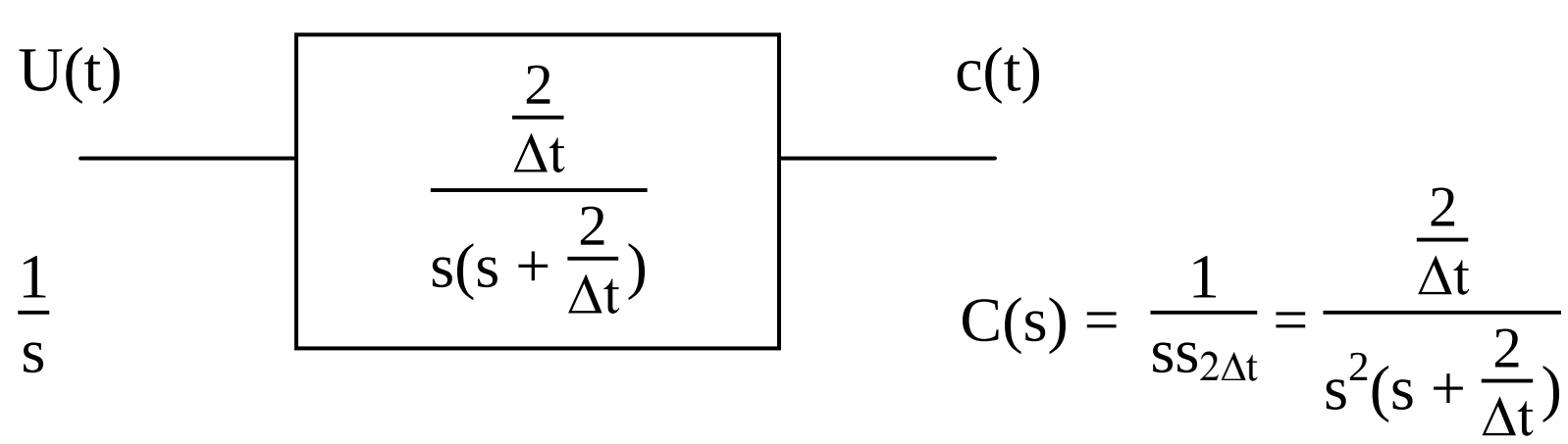

$$C(s) = \frac{\frac{2}{\Delta t}}{s^2(s + \frac{2}{\Delta t})} = \frac{A}{s} + \frac{B}{s^2} + \frac{C}{s + \frac{2}{\Delta t}} \tag{4}$$

Find the constants A, B, and C

$$B = \frac{\frac{2}{\Delta t}}{(s + \frac{2}{\Delta t})}\Big|_{s \to 0} = 1 \tag{5}$$

$$C = \frac{\frac{2}{\Delta t}}{s^2}\Big|_{s\rightarrow -\frac{2}{\Delta t}} = \frac{\Delta t}{2} \qquad 6)$$

$$A = \frac{d}{ds}\left[\frac{\frac{2}{\Delta t}}{(s + \frac{2}{\Delta t})}\right]\Big|_{s\rightarrow 0} = \left[\frac{-\frac{2}{\Delta t}}{(s + \frac{2}{\Delta t})^2}\right]\Big|_{s\rightarrow 0} = -\frac{\Delta t}{2} \qquad 7)$$

Substituting Eq 5 thru Eq 7 into Eq 4

$$C(s) = \frac{\frac{2}{\Delta t}}{s^2(s + \frac{2}{\Delta t})} = -\frac{\Delta t}{2}\frac{1}{s} + \frac{1}{s^2} + \frac{\Delta t}{2}\frac{1}{s + \frac{2}{\Delta t}} \qquad 8)$$

Finding the Inverse $K_{\Delta t}$ Transform of Eq 8

$$K_{\Delta t}^{-1}\left[\frac{1}{s+a}\right] = (1-a\Delta t)^{\frac{t}{\Delta t}} \qquad 9)$$

From Eq 8 and Eq 9

$$c(t) = -\frac{\Delta t}{2} + t + \frac{\Delta t}{2}\left[1 - \frac{2}{\Delta t}\Delta t\right]^{\frac{t}{\Delta t}} \qquad 10)$$

$$c(t) = -\frac{\Delta t}{2} + t + \frac{\Delta t}{2}(-1)^{\frac{t}{\Delta t}} \qquad 11)$$

Then

$$\mathbf{c(t) = t - \frac{\Delta t}{2}\left[1 - (-1)^{\frac{t}{\Delta t}}\right]} \qquad 12)$$

From Eq 12 find c(t) for t = 0, Δt, 2Δt, 3Δt, 4Δt, 5Δt, 6Δt, 7Δt, and 8Δt

Table #1

System output, c(t) vs t

| t | c(t) |
|---|---|
| 0 | 0 |
| Δt | 0 |
| 2Δt | 2Δt |
| 3Δt | 2Δt |
| 4Δt | 4Δt |
| 5Δt | 4Δt |
| 6Δt | 6Δt |
| 7Δt | 6Δt |
| 8Δt | 8Δt |

Plotting c(t) from Table #1

Fig # 3

c(t) vs t plot

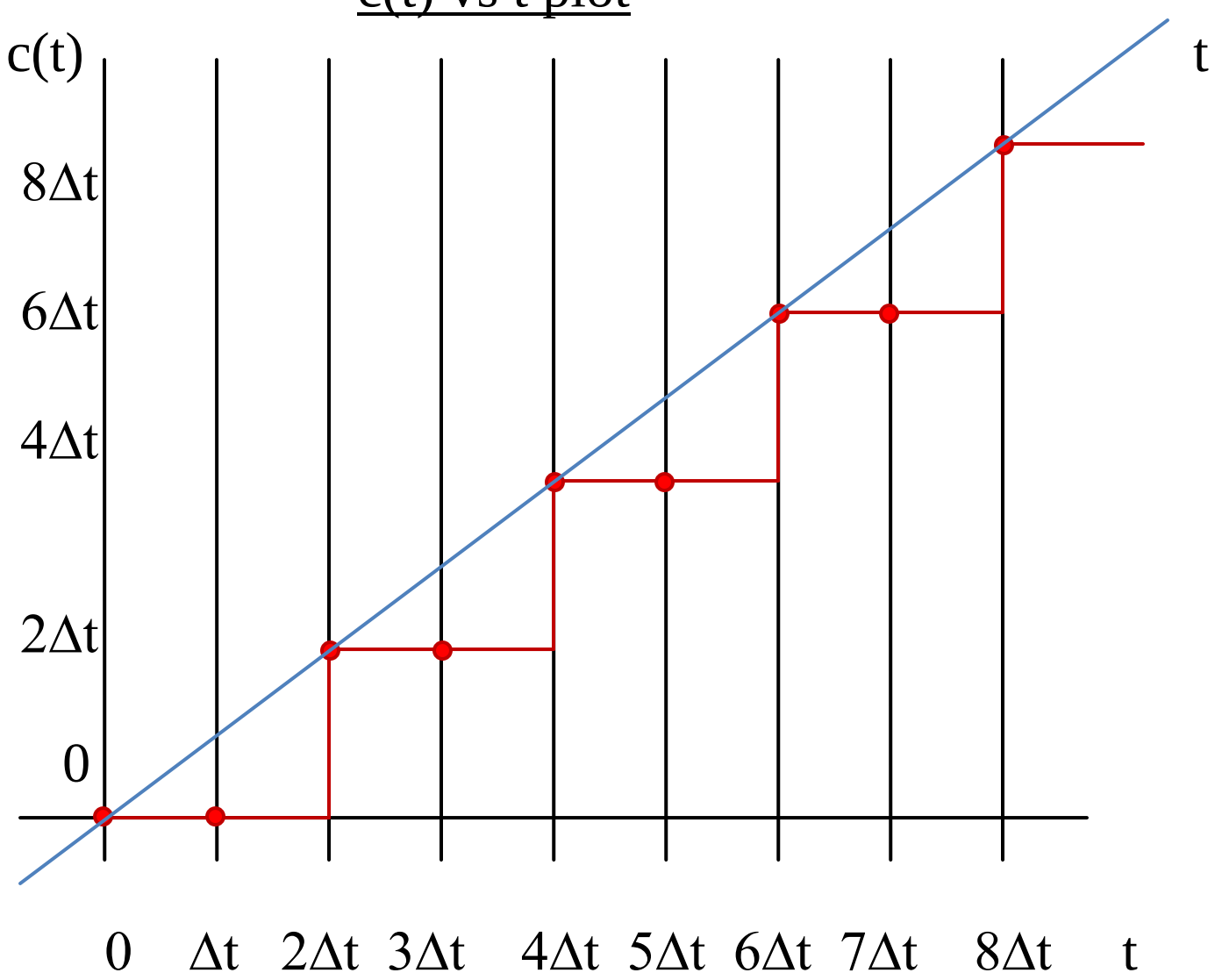

$$c(t) = t - \frac{\Delta t}{2}\left[1 - (-1)^{\frac{t}{\Delta t}}\right]$$

$t = 0, \Delta t, 2\Delta t, 3\Delta t, \ldots$

Note from Fig #3 that system output sampling occurs at the beginning of alternate system clock cycles (i.e. at t = 0, 2Δt, 4Δt, 6Δt, … ) and that the sampling period is 2Δt.

The following system diagram clarifies system diagram Fig #2 by showing the system input and output sample and hold switches.

Fig #4

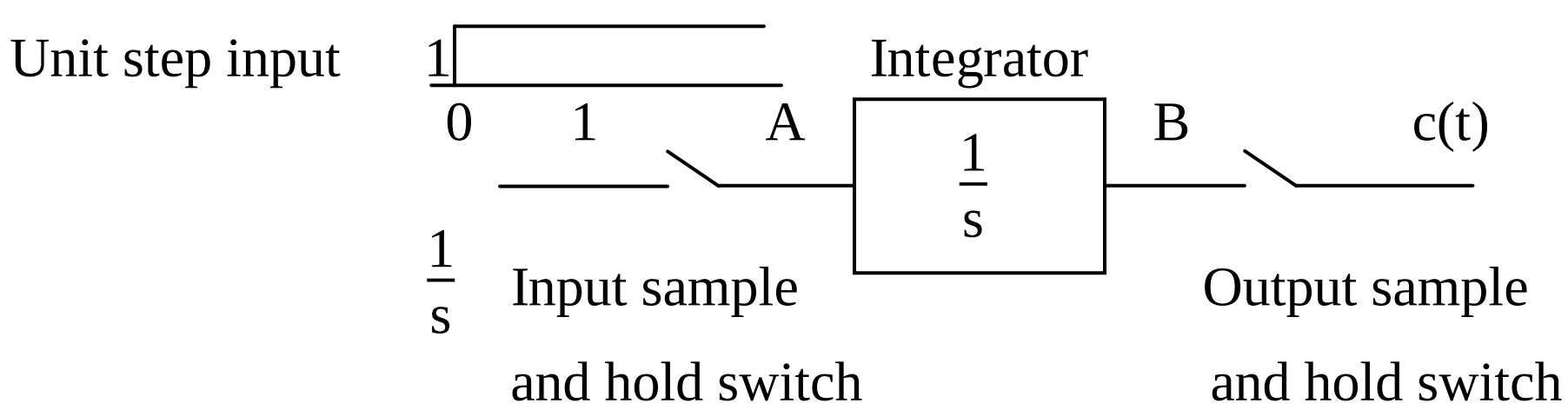

$$c(t) = t - \frac{\Delta t}{2}\left[1 - (-1)^{\frac{t}{\Delta t}}\right]$$

$t = 0, \Delta t, 2\Delta t, 3\Delta t, \ldots$

The input sampler samples at every system clock cycle (i.e. sampling period = Δt). The input to the integrator at A is 1.

The output sampler samples at every other system clock cycle (i.e. sampling period = 2Δt).

The system clock cycle period is Δt.

The combination of the integrator, $\frac{1}{s}$, together with the input and output sample and hold switches, represents the integrator, $\frac{1}{s_{2\Delta t}}$.

Showing the operation of the system of Fig #4 in response to a unit step input

Table #2

| t | A | B | c(t) | Comments | |
|---|---|---|---|---|---|
| 0 | 1 | 0Δt | 0Δt | A = 1 | Input sample of 1 at t = 0 |
| Δt | 1 | 1Δt | 0Δt | B = Δt, c(Δt) = 0Δt, | No output sample at t = Δt |
| 2Δt | 1 | 2Δt | 2Δt | B = Δt+Δt, c(2Δt) = 2Δt, | Output sample at t = 2Δt |
| 3Δt | 1 | 3Δt | 2Δt | B = 2Δt+Δt, c(3Δt) = 2Δt | No output sample at t = 3Δt |
| 4Δt | 1 | 4Δt | 4Δt | B = 3Δt+Δt, c(4Δt) = 4Δt | Output sample at t = 4Δt |
| 5Δt | 1 | 5Δt | 4Δt | B = 4Δt+Δt, c(5Δt) = 4Δt | No output sample at t = 5Δt |
| 6Δt | 1 | 6Δt | 6Δt | B = 5Δt+Δt, c(6Δt) = 6Δt | Output sample at t = 6Δt |
| 7Δt | 1 | 7Δt | 6Δt | B = 6Δt+Δt, c(7Δt) = 6Δt | No output sample at t = 7Δt |
| 8Δt | 1 | 8Δt | 8Δt | B = 7Δt+Δt, c(8Δt) = 8Δt | Output sample at t = 8t |

Note that c(t) in Table #2 is the same as c(t) in Table #1

Good check

Comment − If the integrator of Fig #1 was $\frac{1}{s_{3\Delta t}}$ , using the discrete derivative identity,

$s_{3\Delta t} = s + \Delta t s^2 + \frac{\Delta t^2}{3} s^3$, the output sample and hold switch of Fig #4 would have a sample

and hold period of 3Δt where $c(t) = t - \Delta t + \Delta t[\frac{2}{3-j\sqrt{3}}(\frac{-1-j\sqrt{3}}{2})^{\frac{t}{\Delta t}} + \frac{2}{3+j\sqrt{3}}(\frac{-1+j\sqrt{3}}{2})^{\frac{t}{\Delta t}}]$ .

Example 5.14-3 The analysis of a sampled-data double integrator single clock system with the sampling period of one of the samplers being twice the sampling period of the others

Find the output response of the following sampled-data system to a unit step input. The system is initially passive.

Fig #1

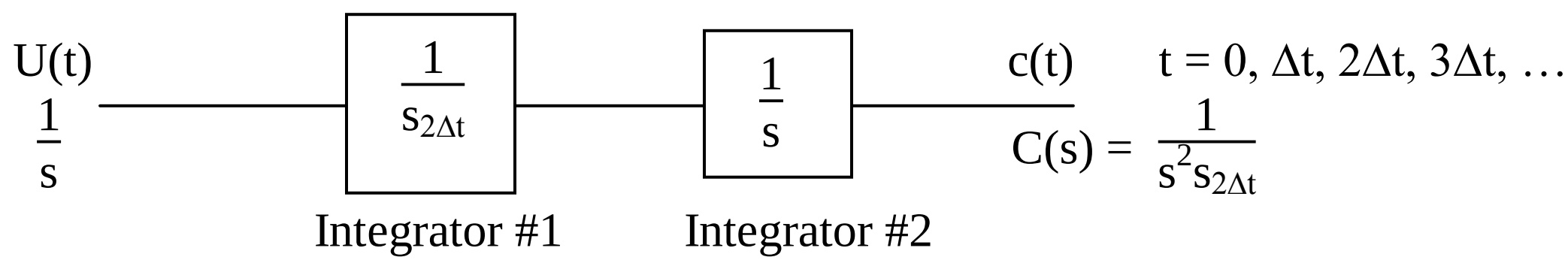

$\frac{1}{s_{2\Delta t}}$ = $K_{\Delta t}$ Transform of Integrator #1

$\frac{1}{s}$ = $K_{\Delta t}$ Transform of Integrator #2

$t = 0, \Delta t, 2\Delta t, 3\Delta t, \ldots$

Definitions

The above operators are defined as follows:

$$sy(t) = \frac{y(t+\Delta t) - y(t)}{\Delta t} , \qquad t = 0, \Delta t, 2\Delta t, 3\Delta t, \ldots$$

$$s_{2\Delta t}y(t) = \frac{y(t+2\Delta t) - y(t)}{2\Delta t} , \qquad t = 0, \Delta t, 2\Delta t, 3\Delta t, \ldots$$

$$s^2y(t) = \frac{y(t+2\Delta t) - 2y(t+\Delta t) + y(t)}{\Delta t^2} , \qquad t = 0, \Delta t, 2\Delta t, 3\Delta t, \ldots$$

$s_{2\Delta t} = s + \frac{\Delta t}{2}s^2$ , A discrete derivative equality 1)

Substituting $s_{2\Delta t} = s + \frac{\Delta t}{2}s^2$ into Fig #1 and simplifying

Fig #2

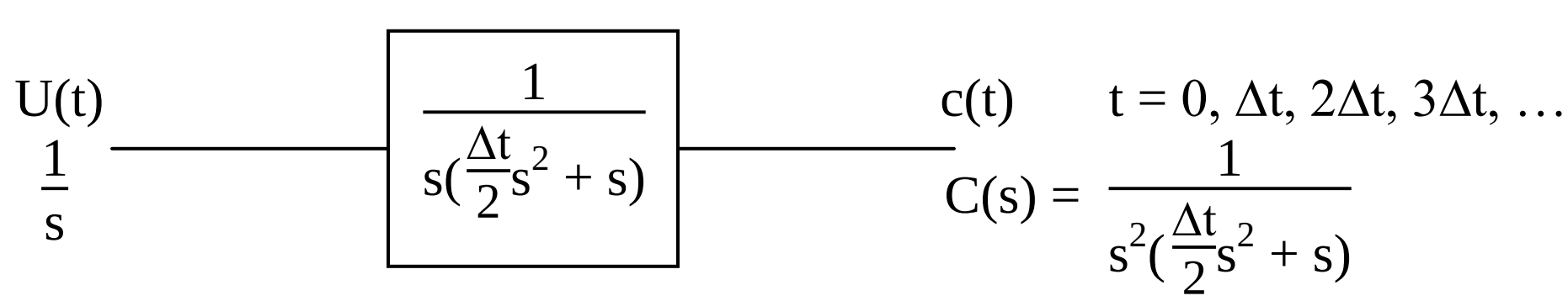

From Fig #2

$$C(s) = \frac{1}{s^2(\frac{\Delta t}{2}s^2 + s)} = \frac{1}{s^3(\frac{\Delta t}{2}s + 1)} = \frac{\frac{2}{\Delta t}}{s^3(s + \frac{2}{\Delta t})} = \frac{A}{s} + \frac{B}{s^2} + \frac{C}{s^3} + \frac{D}{s + \frac{2}{\Delta t}} \quad 2)$$

Finding the constants A, B, C, and D

$$C = \frac{\frac{2}{\Delta t}}{(s + \frac{2}{\Delta t})}\Big|_{s=0} = 1 \quad 3)$$

$$B = \frac{d}{ds}\left[\frac{\frac{2}{\Delta t}}{(s + \frac{2}{\Delta t})}\right]\Big|_{s=0} = \frac{-\frac{2}{\Delta t}}{(s + \frac{2}{\Delta t})^2}\Big|_{s=0} = -\frac{\Delta t}{2} \quad 4)$$

$$A = \frac{1}{2!}\frac{d^2}{ds^2}\left[\frac{\frac{2}{\Delta t}}{(s + \frac{2}{\Delta t})}\right]\Big|_{s=0} = \frac{1}{2}\frac{d}{ds}\frac{d}{ds}\left[\frac{\frac{2}{\Delta t}}{(s + \frac{2}{\Delta t})}\right]\Big|_{s=0} = \frac{1}{2}\frac{d}{ds}\frac{-\frac{2}{\Delta t}}{(s + \frac{2}{\Delta t})^2}\Big|_{s=0} = \frac{1}{2}\left(\frac{-2}{\Delta t}\right)\frac{-2(s + \frac{2}{\Delta t})}{(s + \frac{2}{\Delta t})^4}\Big|_{s=0} \quad 5)$$

$$A = \frac{\frac{2}{\Delta t}}{(s + \frac{2}{\Delta t})^3}\Big|_{s=0} = \frac{\Delta t^2}{4} \quad 6)$$

$$D = \frac{\frac{2}{\Delta t}}{s^3}\Big|_{s=-\frac{2}{\Delta t}} = -\frac{\Delta t^2}{4} \quad 7)$$

Substituting Eq 3, Eq 4, Eq 6, and Eq 7 into Eq 2

$$C(s) = \frac{\Delta t^2}{4}\frac{1}{s} - \frac{\Delta t}{2}\frac{1}{s^2} + \frac{1}{s^3} - \frac{\Delta t^2}{4}\frac{1}{s + \frac{2}{\Delta t}} \quad 8)$$

$$K_{\Delta t}^{-1}\left[\frac{1}{s}\right] = 1 \quad 9)$$

$$K_{\Delta t}^{-1}\left[\frac{1}{s^2}\right] = t \quad 10)$$

$$K_{\Delta t}^{-1}\left[\frac{1}{s^3}\right] = \frac{1}{2}t(t-\Delta t) \quad 11)$$

$$K_{\Delta t}^{-1}[\frac{1}{s+a}] = (1-a\Delta t)^{\frac{t}{\Delta t}} \qquad 12)$$

Finding the Inverse $K_{\Delta t}$ Transform of Eq 8 using Eq 9 thru Eq 12

$$c(t) = \frac{\Delta t^2}{4} - \frac{\Delta t}{2} t + \frac{1}{2} t(t-\Delta t) - \frac{\Delta t^2}{4} [-1]^{\frac{t}{\Delta t}} \quad , \quad t = 0, \Delta t, 2\Delta t, 3\Delta t, \ldots \qquad 13)$$

Simplifying Eq 13

$$\mathbf{c(t) =} \frac{\Delta t^2}{4} [\, 1 - [-1]^{\frac{t}{\Delta t}}] + \frac{t^2}{2} - \frac{\Delta t}{2} t - \frac{\Delta t}{2} t \qquad 14)$$

Then

$$\mathbf{c(t) = \frac{\Delta t^2}{4} [\, 1 - [-1]^{\frac{t}{\Delta t}}] + \frac{1}{2} t(t - 2\Delta t) \quad , \quad t = 0, \Delta t, 2\Delta t, 3\Delta t, \ldots} \qquad 15)$$

Evaluating Eq 15 for $t = 0, \Delta t, 2\Delta t, 3\Delta t, 4\Delta t, 5\Delta t, 6\Delta t, 7\Delta t$, and $8\Delta t$

$$c(0) = 0 \qquad 16)$$

$$c(\Delta t) = \frac{\Delta t^2}{2} - \frac{\Delta t^2}{2} = 0 \qquad 17)$$

$$c(2\Delta t) = 0 + 0 = 0 \qquad 18)$$

$$c(3\Delta t) = \frac{\Delta t^2}{2} + \frac{3\Delta t^2}{2} = 2\Delta t^2 \qquad 19)$$

$$c(4\Delta t) = 0 + \frac{1}{2} 4\Delta t(2\Delta t) = 4\Delta t^2 \qquad 20)$$

$$c(5\Delta t) = \frac{\Delta t^2}{2} + \frac{1}{2} 5\Delta t(3\Delta t) = 8\Delta t^2 \qquad 21)$$

$$c(6\Delta t) = 0 + \frac{1}{2} 6\Delta t(4\Delta t) = 12\Delta t^2 \qquad 22)$$

$$c(7\Delta t) = \frac{\Delta t^2}{2} + \frac{1}{2} 7\Delta t(5\Delta t) = 18\Delta t^2 \qquad 23)$$

$$c(8\Delta t) = 0 + \frac{1}{2} 8\Delta t(6\Delta t) = 24\Delta t^2 \qquad 24)$$

The previously calculated values are placed in the following table, Table #1.

Table #1

c(t) vs t

| t | c(t) |
|---|---|
| 0 | 0 |
| $\Delta t$ | 0 |
| $2\Delta t$ | 0 |
| $3\Delta t$ | $2\Delta t^2$ |
| $4\Delta t$ | $4\Delta t^2$ |
| $5\Delta t$ | $8\Delta t^2$ |
| $6\Delta t$ | $12\Delta t^2$ |
| $7\Delta t$ | $18\Delta t^2$ |
| $8\Delta t$ | $24\Delta t^2$ |

Rewriting Fig #1 in more detail showing the system sample and hold sampling switches and checking Eq 15

Fig #3

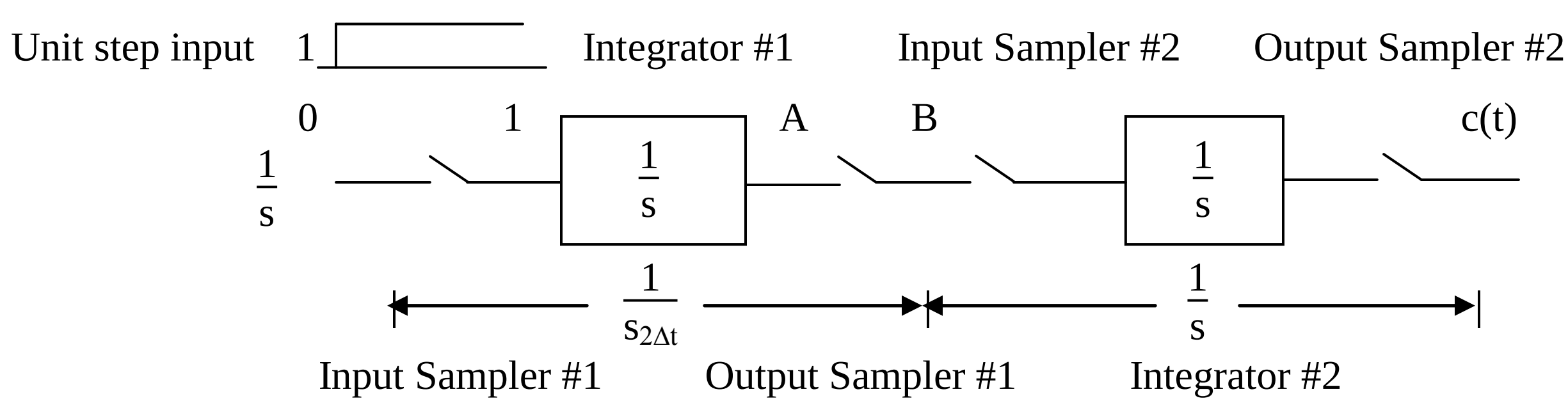

Samplers #1 and #2 are sample and hold switches.
Input Sampler #1 period = $\Delta t$ (i.e. sampling at every system clock cycle)
Output Sampler #1 period = $2\Delta t$ (i.e. sampling at every other system clock cycle)
Input Sampler #2 period = $\Delta t$ (i.e. sampling at every system clock cycle)
Output Sampler #2 period = $\Delta t$ (i.e. sampling at every system clock cycle)
The output waveform of Integrator #1 changes value every $2\Delta t$
$t = 0, \Delta t, 2\Delta t, 3\Delta t, \ldots$
The system clock cycle period is $\Delta t$.

Visually consider the operation of the Fig #3 system in response to a unit step input. Note the following table, Table #2, which lists the values at various points in the system as a function of t.

Table #2

| t | A | B | c(t) | Comments |
|---|---|---|---|---|
| 0 | $0\Delta t$ | $0\Delta t$ | $0\Delta t^2$ | Input Sampler #1 puts 1 at the input of Integrator #1 |
| $\Delta t$ | $1\Delta t$ | $0\Delta t$ | $0\Delta t^2$ | $A = \Delta t$ |
| $2\Delta t$ | $2\Delta t$ | $2\Delta t$ | $0\Delta t^2$ | $A = 2\Delta t$, $B = A = 2\Delta t$ (Output sampler #1 period is $2\Delta t$) |
| $3\Delta t$ | $3\Delta t$ | $2\Delta t$ | $2\Delta t^2$ | $A = 3\Delta t$, $B = 2\Delta t$, $C = c(3\Delta t) = 2\Delta t(\Delta t) = 2\Delta t^2$ |
| $4\Delta t$ | $4\Delta t$ | $4\Delta t$ | $4\Delta t^2$ | $A = 4\Delta t$, $B = 4\Delta t$, $C = c(4\Delta t) = 2\Delta t^2 + 2\Delta t(\Delta t) = 4\Delta t^2$ |
| $5\Delta t$ | $5\Delta t$ | $4\Delta t$ | $8\Delta t^2$ | $A = 5\Delta t$, $B = 4\Delta t$, $C = c(5\Delta t) = 4\Delta t^2 + 4\Delta t(\Delta t) = 8\Delta t^2$ |
| $6\Delta t$ | $6\Delta t$ | $6\Delta t$ | $12\Delta t^2$ | $A = 6\Delta t$, $B = 6\Delta t$, $C = c(6\Delta t) = 8\Delta t^2 + 4\Delta t(\Delta t) = 12\Delta t^2$ |
| $7\Delta t$ | $7\Delta t$ | $6\Delta t$ | $18\Delta t^2$ | $A = 7\Delta t$, $B = 6\Delta t$, $C = c(7\Delta t) = 12\Delta t^2 + 6\Delta t(\Delta t) = 18\Delta t^2$ |
| $8\Delta t$ | $8\Delta t$ | $8\Delta t$ | $24\Delta t^2$ | $A = 8\Delta t$, $B = 8\Delta t$, $C = c(8\Delta t) = 18\Delta t^2 + 6\Delta t(\Delta t) = 24\Delta t^2$ |

Note that c(t) in Table #2 is the same as c(t) in Table #1.

Good check

Example 5.14-4 The analysis of a sampled-data closed loop single clock system where the sampling period of one of the two samplers is twice the sampling period of the other

Demonstrate the analysis of a sampled-data closed loop system with two sample and hold samplers. The sampling period of one sampler is twice the sampling period of the other. Use the discrete derivative equality, $s_{2\Delta t} = s + \frac{\Delta t}{2} s^2$. The system is initially passive. Find c(t) at t = 0, Δt, 2Δt, 3Δt, 4Δt, 5Δt, 6Δt, and 7Δt and check the calculated values of c(t) using various other methods. Also, find the maximum system clock period, Δt, for which the system is stable.

The above specified operators are defined as follows:

$$sy(t) = \frac{y(t+\Delta t) - y(t)}{\Delta t}$$

$$s^2y(t) = \frac{y(t+2\Delta t) - 2y(t+\Delta t) + y(t)}{\Delta t^2}$$

$$s_{2\Delta t}y(t) = sy(t) + \frac{\Delta t}{2}s^2y(t) = \frac{y(t+2\Delta t) - y(t)}{2\Delta t}$$

where
y(t) = a function of t
t = 0, Δt, 2Δt, 3Δt, …

Consider the following closed loop integrator sampled-data system with a unit step input.

Fig #1

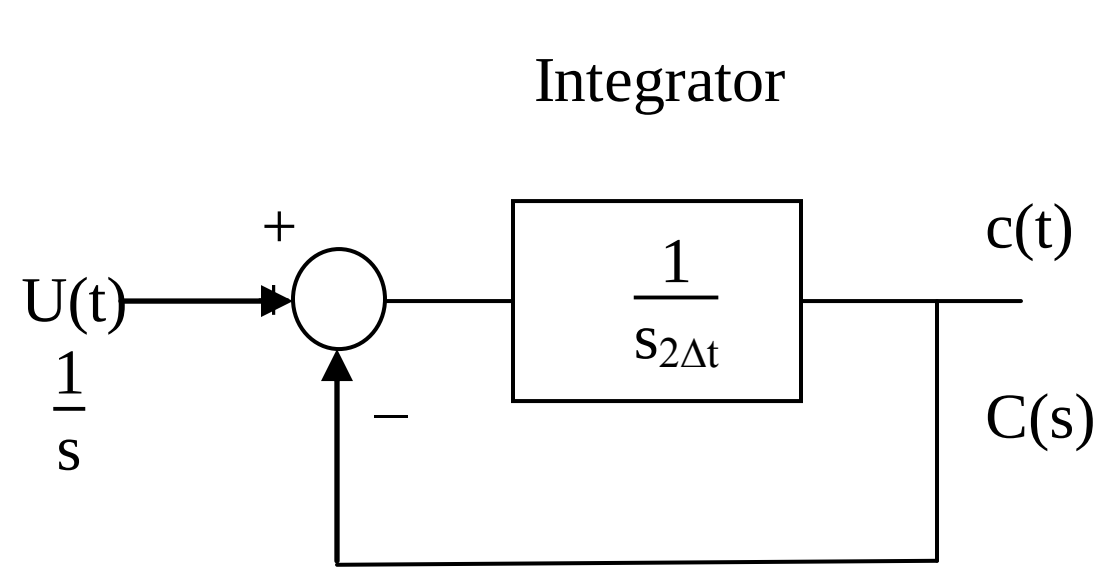

System is initially passive

$$C(s) = \frac{1}{s}\left(\frac{\frac{1}{s_{\Delta t}}}{1 + \frac{1}{s_{\Delta t}}}\right) \qquad 1)$$

$$C(s) = \frac{1}{s}\left(\frac{1}{s_{\Delta t} + 1}\right) \qquad 2)$$

$$s_{2\Delta t} = s + \frac{\Delta t}{2}s^2 \quad \text{, A discrete derivative equality} \qquad 3)$$

Substituting Eq 3 into Eq 2

$$C(s) = \frac{1}{s} \frac{1}{\frac{\Delta t}{2} s^2 + s + 1} \qquad 4)$$

Changing the form of Eq 4

$$C(s) = \frac{1}{s} \frac{\frac{2}{\Delta t}}{s^2 + \frac{2}{\Delta t}s + \frac{2}{\Delta t}} \qquad 5)$$

The roots of the denominator of Eq 5 are found using the quadratic equation.

$$C(s) = \frac{\frac{2}{\Delta t}}{s(\frac{1}{\Delta t} + \frac{1}{\Delta t}\sqrt{1\text{-}2\Delta t})(\frac{1}{\Delta t} - \frac{1}{\Delta t}\sqrt{1\text{-}2\Delta t})} \qquad 6)$$

Expanding Eq 6 into a partial fraction expansion

$$C(s) = \frac{\frac{2}{\Delta t}}{s(s + \frac{1}{\Delta t} + \frac{1}{\Delta t}\sqrt{1\text{-}2\Delta t})(s + \frac{1}{\Delta t} - \frac{1}{\Delta t}\sqrt{1\text{-}2\Delta t})} = \frac{A}{s} + \frac{B}{(s + \frac{1}{\Delta t} + \frac{1}{\Delta t}\sqrt{1\text{-}2\Delta t})} + \frac{C}{(s + \frac{1}{\Delta t} - \frac{1}{\Delta t}\sqrt{1\text{-}2\Delta t})} \qquad 7)$$

From Eq 5 and Eq 7

$$A = \frac{\frac{2}{\Delta t}}{s^2 + \frac{2}{\Delta t}s + \frac{2}{\Delta t}} \Big|_{s = 0} = 1 \qquad 8)$$

From Eq 7

$$B = \frac{\frac{2}{\Delta t}}{s(s + \frac{1}{\Delta t} - \frac{1}{\Delta t}\sqrt{1\text{-}2\Delta t})} \Big|_{s = -\frac{1}{\Delta t} - \frac{1}{\Delta t}\sqrt{1\text{-}2\Delta t}} = \frac{\frac{2}{\Delta t}}{(-\frac{1}{\Delta t} - \frac{1}{\Delta t}\sqrt{1\text{-}2\Delta t})(-\frac{2}{\Delta t}\sqrt{1\text{-}2\Delta t})} \qquad 9)$$

$$B = \frac{\Delta t}{(1 + \sqrt{1\text{-}2\Delta t})(\sqrt{1\text{-}2\Delta t})} = \frac{\Delta t}{\sqrt{1\text{-}2\Delta t} + 1\text{-}2\Delta t} \qquad 10)$$

$$B = \frac{\Delta t}{1\text{-}2\Delta t + \sqrt{1\text{-}2\Delta t}} \qquad 11)$$

$$C = \frac{\frac{2}{\Delta t}}{s(s + \frac{1}{\Delta t} + \frac{1}{\Delta t}\sqrt{1-2\Delta t})}\Big|_{s = -\frac{1}{\Delta t} + \frac{1}{\Delta t}\sqrt{1-2\Delta t}} = \frac{\frac{2}{\Delta t}}{(-\frac{1}{\Delta t} + \frac{1}{\Delta t}\sqrt{1-2\Delta t})(\frac{2}{\Delta t}\sqrt{1-2\Delta t})} \quad 12)$$

$$C = \frac{\Delta t}{(-1 + \sqrt{1-2\Delta t})(\sqrt{1-2\Delta t})} = \frac{\Delta t}{-\sqrt{1-2\Delta t} + 1-2\Delta t} \quad 13)$$

$$C = \frac{\Delta t}{1-2\Delta t - \sqrt{1-2\Delta t}} \quad 14)$$

Substituting Eq 8, Eq 11, and Eq 14 into Eq 7

$$C(s) = \frac{1}{s} + \frac{\Delta t}{1-2\Delta t + \sqrt{1-2\Delta t}}\,\frac{1}{s + \frac{1}{\Delta t} + \frac{1}{\Delta t}\sqrt{1-2\Delta t}} + \frac{\Delta t}{1-2\Delta t - \sqrt{1-2\Delta t}}\,\frac{1}{s + \frac{1}{\Delta t} - \frac{1}{\Delta t}\sqrt{1-2\Delta t}} \quad 15)$$

Find the Inverse $K_{\Delta t}$ Transform of Eq 15

$$K_{\Delta t}^{-1}[\frac{1}{s+b}] = (1-b\Delta t)^{\frac{t}{\Delta t}} \quad 16)$$

From Eq 15 and Eq 16

$$c(t) = 1 + \frac{\Delta t}{1-2\Delta t + \sqrt{1-2\Delta t}}(1 - \frac{1 + \sqrt{1-2\Delta t}}{\Delta t}\Delta t)^{\frac{t}{\Delta t}} + \frac{\Delta t}{1-2\Delta t - \sqrt{1-2\Delta t}}(1 - \frac{1 - \sqrt{1-2\Delta t}}{\Delta t}\Delta t)^{\frac{t}{\Delta t}} \quad 17)$$

$$c(t) = 1 + \frac{\Delta t}{1-2\Delta t + \sqrt{1-2\Delta t}}(-\sqrt{1-2\Delta t})^{\frac{t}{\Delta t}} + \frac{\Delta t}{1-2\Delta t - \sqrt{1-2\Delta t}}(\sqrt{1-2\Delta t})^{\frac{t}{\Delta t}} \quad 18)$$

Simplifying Eq 18

Let

$$a = 1-2\Delta t \quad 19)$$

Substituting Eq 19 into Eq 18

$$c(t) = 1 + \frac{\Delta t}{a + \sqrt{a}}(-\sqrt{a})^{\frac{t}{\Delta t}} + \frac{\Delta t}{a - \sqrt{a}}(\sqrt{a})^{\frac{t}{\Delta t}} \quad 20)$$

Simplifying Eq 20

$$c(t) = 1 + \frac{\Delta t(a - \sqrt{a})}{a^2 - a}(-\sqrt{a})^{\frac{t}{\Delta t}} + \frac{\Delta t(a + \sqrt{a})}{a^2 - a}(\sqrt{a})^{\frac{t}{\Delta t}} = 1 + \frac{\Delta t(a - \sqrt{a})}{a(a-1)}(-\sqrt{a})^{\frac{t}{\Delta t}} + \frac{\Delta t(a + \sqrt{a})}{a(a-1)}(\sqrt{a})^{\frac{t}{\Delta t}} \quad 21)$$

$a(a-1) = a(1-2\Delta t-1) = a(-2\Delta t)$ 22)

Substituting Eq 22 into Eq 21

$$c(t) = 1 + \frac{\Delta t(a - \sqrt{a})}{a(-2\Delta t)}(-\sqrt{a})^{\frac{t}{\Delta t}} + \frac{\Delta t(a + \sqrt{a})}{a(-2\Delta t)}(\sqrt{a})^{\frac{t}{\Delta t}} \quad 23)$$

Then

$$\mathbf{c(t) = 1 - \frac{1}{2a}\left[(a - \sqrt{a})(-\sqrt{a})^{\frac{t}{\Delta t}} + (a + \sqrt{a})(\sqrt{a})^{\frac{t}{\Delta t}}\right]} \quad 24)$$

**where**

**$a = 1-2\Delta t$**

Calculating c(t) for t = 0, Δt, 2Δt, 3Δt, …

$$c(0) = 1 - \frac{1}{2a}[a - \sqrt{a} + a + \sqrt{a}] = 1 - 1 = 0 \quad 25)$$

$$c(\Delta t) = 1 - \frac{1}{2a}[(a - \sqrt{a})(-\sqrt{a}) + (a + \sqrt{a})(\sqrt{a})] = 1 - \frac{1}{2a}[(-a\sqrt{a} + a + a\sqrt{a} + a] = 1 - 1 = 0 \quad 26)$$

$$c(2\Delta t) = 1 - \frac{1}{2a}[(a - \sqrt{a})(-\sqrt{a})^2 + (a + \sqrt{a})(\sqrt{a})^2] = 1 - \frac{1}{2a}[(a^2 - a\sqrt{a} + a^2 + a\sqrt{a}] = 1 - a \quad 27)$$

$$c(3\Delta t) = 1 - \frac{1}{2a}[(a - \sqrt{a})(-\sqrt{a})^3 + (a + \sqrt{a})(\sqrt{a})^3] = 1 - \frac{1}{2a}[(-a^2\sqrt{a} + a^2 + a\sqrt{a} + a^2] = 1 - a \quad 28)$$

$$c(4\Delta t) = 1 - \frac{1}{2a}[(a - \sqrt{a})(-\sqrt{a})^4 + (a + \sqrt{a})(\sqrt{a})^4] = 1 - \frac{1}{2a}[(a^3 - a^2\sqrt{a} + a^3 + a^2\sqrt{a}] = 1 - a^2 \quad 29)$$

$$c(5\Delta t) = 1 - \frac{1}{2a}[(a - \sqrt{a})(-\sqrt{a})^5 + (a + \sqrt{a})(\sqrt{a})^5] = 1 - \frac{1}{2a}[(-a^3\sqrt{a} + a^3 + a^3\sqrt{a} + a^3] = 1 - a^2 \quad 30)$$

$$c(6\Delta t) = 1 - \frac{1}{2a}[(a - \sqrt{a})(-\sqrt{a})^6 + (a + \sqrt{a})(\sqrt{a})^6] = 1 - \frac{1}{2a}[(a^4 - a^3\sqrt{a} + a^4 + a^3\sqrt{a}] = 1 - a^3 \quad 31)$$

$$c(7\Delta t) = 1 - \frac{1}{2a}[(a - \sqrt{a})(-\sqrt{a})^7 + (a + \sqrt{a})(\sqrt{a})^7] = 1 - \frac{1}{2a}[(-a^4\sqrt{a} + a^4 + a^4\sqrt{a} + a^4] = 1 - a^3 \quad 32)$$

Putting the previously calculated values of c(t) into a table

Table #1

c(t) vs t

| t | c(t) |
|---|---|
| 0 | 0 |
| $\Delta t$ | 0 |
| $2\Delta t$ | $1 - (1-2\Delta t)^1$ |
| $3\Delta t$ | $1 - (1-2\Delta t)^1$ |
| $4\Delta t$ | $1 - (1-2\Delta t)^2$ |
| $5\Delta t$ | $1 - (1-2\Delta t)^2$ |
| $6\Delta t$ | $1 - (1-2\Delta t)^3$ |
| $7\Delta t$ | $1 - (1-2\Delta t)^3$ |

In general

$$C(n\Delta t) = C([n+1]\Delta t) = 1 - a^{\frac{n}{2}}, \quad n = 0, 2, 4, 6, \ldots \qquad 33)$$

$$\mathbf{C(n\Delta t) = C([n+1]\Delta t) = 1 - (1-2\Delta t)^{\frac{n}{2}}, \quad n = 0, 2, 4, 6, \ldots} \qquad 34)$$

Finding the maximum system clock period, $\Delta t$, for which the system is stable

Noting Eq 34 it is seen that the system is stable for $|1-2\Delta t| < 1$, $\Delta t \neq 0$

Then

**For system stability $0 < \Delta t < 1$** 35)

The system diagram of Fig # 1 can be depicted in a different way showing the system samplers. See Fig #2 below.

Fig #2

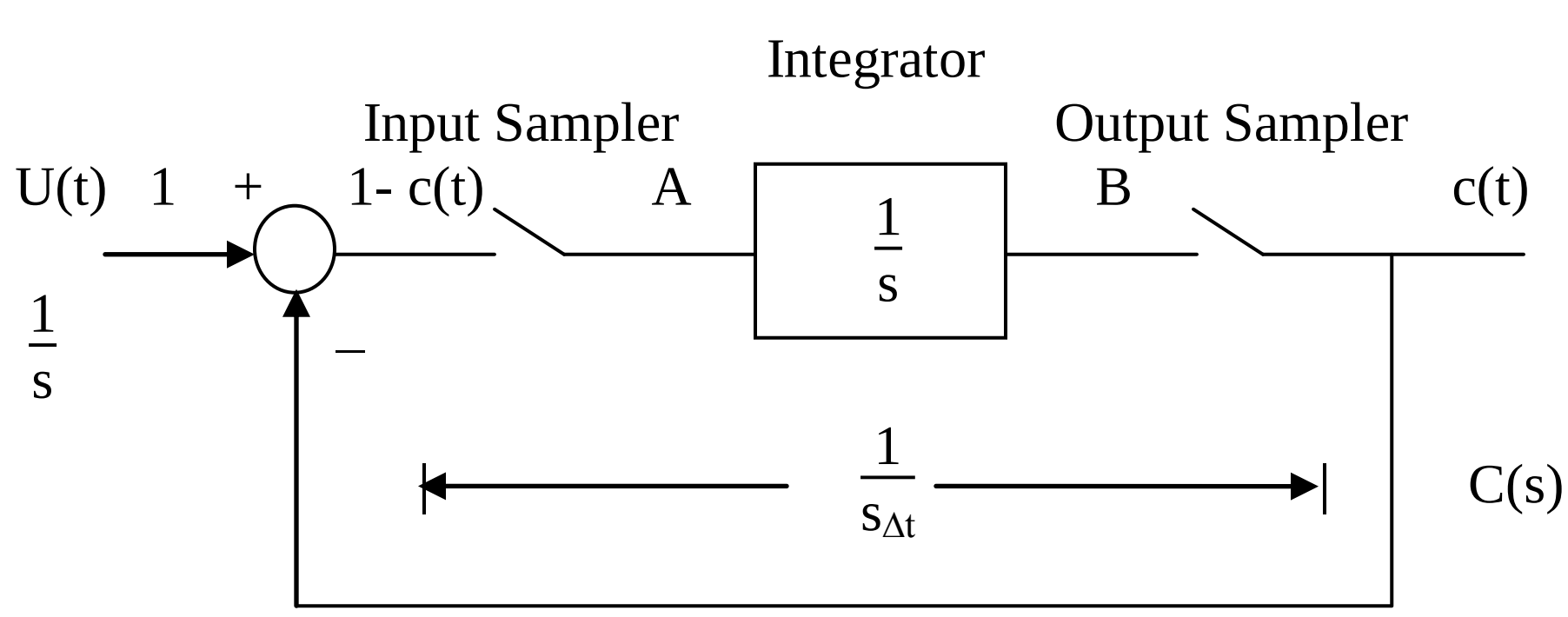

The system is initially passive

The Input Sampler period is $\Delta t$ (i.e. sampling at every system clock cycle)
The Output Sampler period is $2\Delta t$ (i.e. sampling at every other system clock cycle)
$t = 0, \Delta t, 2\Delta t, 3\Delta t, \ldots$
The system clock cycle period is $\Delta t$.

Checking Eq 34

Check #1

Check the calculated values of c(t) by observing the operation of the system of Fig #2 as t increases. Note Table #2 below.

Table #2

| t | A | B | c(t) |
|---|---|---|---|
| 0 | 1 | 0 | 0 |
| $\Delta t$ | 1 | $1(\Delta t) = \Delta t$<br>$\Delta t$ | 0 |
| $2\Delta t$ | $1-2\Delta t$ | $\Delta t + 1(\Delta t)$<br>$2\Delta t$ | $2\Delta t$<br>$1-(1-2\Delta t)^1$ |
| $3\Delta t$ | $1-2\Delta t$ | $2\Delta t + \Delta t(1-2\Delta t)$<br>$3\Delta t-2\Delta t^2$ | $2\Delta t$<br>$1-(1-2\Delta t)^1$ |
| $4\Delta t$ | $1-(4\Delta t-4\Delta t^2)$<br>$1-4\Delta t+4\Delta t^2$ | $3\Delta t-2\Delta t^2+\Delta t(1-2\Delta t)$<br>$4\Delta t-4\Delta t^2$ | $4\Delta t-4\Delta t^2$<br>$1-(1-2\Delta t)^2$ |
| $5\Delta t$ | $1-(4\Delta t-4\Delta t^2)$<br>$1-4\Delta t+4\Delta t^2$ | $4\Delta t-4\Delta t^2+\Delta t(1-4\Delta t+4\Delta t^2)$<br>$5\Delta t-8\Delta t^2+4\Delta t^3$ | $4\Delta t-4\Delta t^2$<br>$1-(1-2\Delta t)^2$ |
| $6\Delta t$ | $1-(6\Delta t-12\Delta t^2+8\Delta t^3)$<br>$1-6\Delta t+12\Delta t^2-8\Delta t^3$ | $5\Delta t-8\Delta t^2+4\Delta t^3+\Delta t(1-4\Delta t+4\Delta t^2)$<br>$6\Delta t-12\Delta t^2+8\Delta t^3$ | $6\Delta t-12\Delta t^2+8\Delta t^3$<br>$1-(1-2\Delta t)^3$ |
| $7\Delta t$ | | | $6\Delta t-12\Delta t^2+8\Delta t^3$<br>$1-(1-2\Delta t)^3$ |

The c(t) values in Table #2 are the same as those in Table #1.

Good check

Check #2

Check the calculated values of c(t) using the $K_{\Delta t}$ Transform Asymptotic Series and the Discrete Maclaurin Series.

From the system $K_{\Delta t}$ Transform, Eq 4

$$C(s) = \frac{1}{s}\,\frac{1}{\frac{\Delta t}{2}s^2 + s + 1} \qquad 36)$$

$$C(s) = \frac{1}{\frac{\Delta t}{2}s^3 + s^2 + s} \qquad 37)$$

Dividing the denominator of Eq 37 into the numerator of Eq 37

$$\begin{array}{r|l} & 0s^{-1} + 0s^{-2} + \frac{2}{\Delta t}s^{-3} - \frac{4}{\Delta t^2}s^{-4} + \frac{8-4\Delta t}{\Delta t^3}s^{-5} + \dots \\ \hline \frac{\Delta t}{2}s^3 + s^2 + s & 1 \end{array} \qquad 38)$$

From the Asymptotic Series of Eq 38

$$c(0) = 0 \qquad 39)$$

$$D_{\Delta t}\, c(0) = 0 \qquad 40)$$

$$D_{\Delta t}^{2}c(0) = \frac{2}{\Delta t} \qquad 41)$$

$$D_{\Delta t}^{3}c(0) = -\frac{4}{\Delta t^2} \qquad 42)$$

$$D_{\Delta t}^{4}c(0) = \frac{8-4\Delta t}{\Delta t^3} \qquad 43)$$

…

Using a Discrete Maclaurin Series to represent c(t)

$$c(t) = c(0) + \frac{D_{\Delta t}^{1}c(0)}{1!}t + \frac{D_{\Delta t}^{2}c(0)}{2!}t(t-\Delta t) + \frac{D_{\Delta t}^{3}c(0)}{3!}t(t-\Delta t)(t-2\Delta t) + \frac{D_{\Delta t}^{4}c(0)}{4!}t(t-\Delta t)(t-2\Delta t)(t-3\Delta t) + \dots \qquad 44)$$

Substituting Eq 39 thru Eq 43 into Eq 44 and simplifying

$$c(t) = \frac{1}{\Delta t} t(t-\Delta t) - \frac{2}{3\Delta t^2} t(t-\Delta t)(t-2\Delta t) + \frac{8-4\Delta t}{24\Delta t^3} t(t-\Delta t)(t-2\Delta t)(t-3\Delta t) + \ldots \qquad 45)$$

From Eq 45 calculate c(t) for t = 0, Δt, 2Δt, 3Δt, and 4Δt

$$c(0) = 0 \qquad 46)$$

$$c(\Delta t) = 0 \qquad 47)$$

$$c(2\Delta t) = 2\Delta t = 1 - (1-2\Delta t) \qquad 48)$$

$$c(3\Delta t) = 6\Delta t - \frac{2}{3\Delta t^2} 3!\Delta t^3 = 6\Delta t - 4\Delta t = 2\Delta t = 1 - (1-2\Delta t) \qquad 49)$$

$$c(4\Delta t) = 4(3)\Delta t - \frac{2}{3\Delta t^2} 4(3)(2)\Delta t^3 + \frac{8-4\Delta t}{24\Delta t^3} 4(3)(2)(1)\Delta t^4 = 12\Delta t - 16\Delta t + 8\Delta t - 4\Delta t^2 = 4\Delta t - 4\Delta t^2 \qquad 50)$$

$$c(4\Delta t) = 4\Delta t - 4\Delta t^2 = 1 - (1-2\Delta t)^2 \qquad 51)$$

The above values of c(t) calculated from the $K_{\Delta t}$ Transform Series and Discrete Maclaurin Series are shown below in Table #3.

Table #3

c(t) vs t

| t | c(t) |
|---|---|
| 0 | 0 |
| Δt | 0 |
| 2Δt | $1 - (1-2\Delta t)^1$ |
| 3Δt | $1 - (1-2\Delta t)^1$ |
| 4Δt | $1 - (1-2\Delta t)^2$ |

The values of c(t) shown in Table #3 are the same as those shown in Table #1.

Good check

Section 5.16: The use of Z Transforms with discrete Interval Calculus functions

In the solution of difference equations and problems involving discrete variables, the Z Transform has been successfully and efficiently used for many years. This paper has introduced the $K_{\Delta t}$ Transform which is a suitable substitute for the Z Transform. One might think that a mathematician or engineer would have only these two transforms to choose from. However, this is not the case. There is a third choice, a hybrid combination of both the Z and $K_{\Delta t}$ transforms. The Z Transform can be used with the discrete functions of the $K_{\Delta t}$Transform instead of those of Calculus. In this section, this third choice transform methodology will be presented together with an explanation of its advantages.

During the development of Interval Calculus, a problem was recognized and eventually resolved. While the solution of a differential difference equation, $D_{\Delta t}y(t) + y(t) = e_{\Delta t}(a,t)$, was easy to solve, the differential difference equation, $D_{\Delta t}y(t) + y(t) = e^{at}$, was not. The problem encountered in the latter equation is the presence of a function mismatch. The function, $e^{at}$, is an Interval Calculus function but for an infinitesimal value of $\Delta t$ ($e^{at} = \lim_{\Delta t\to 0}(1+a\Delta t)^{\frac{t}{\Delta t}}$). The discrete derivative interval value of $\Delta t$ is not infinitesimal. In spite of the $\Delta t$ interval mismatch, the equation, $D_{\Delta t}y(t) + y(t) = e^{at}$, is not unreasonable and should have a solution for y(t). The difficulty was resolved through the conceptualization and derivation of Interval Calculus function identities.

Consider the derivation of the $e^{at}$ Interval Calculus identity which is shown below.

$$e_{\Delta t}\left(\frac{e^{a\Delta t}-1}{\Delta t},t\right) = \left(1 + \frac{e^{a\Delta t}-1}{\Delta t}\Delta t\right)^{\frac{t}{\Delta t}} = e^{a\Delta t\frac{t}{\Delta t}} = e^{at} \qquad (5.16\text{-}1)$$

$$e^{at} = e_{\Delta t}\left(\frac{e^{a\Delta t}-1}{\Delta t},t\right) \qquad (5.16\text{-}2)$$

$e_{\Delta t}\left(\frac{e^{a\Delta t}-1}{\Delta t},t\right)$ is the Interval Calculus identity for $e^{at}$.

At the discrete values of $t = 0,\Delta t,2\Delta t,3\Delta t,\ldots$, $e^{at}$ and $e_{\Delta t}\left(\frac{e^{a\Delta t}-1}{\Delta t},t\right)$ have the same values.

Then

$D_{\Delta t}y(t) + y(t) = e_{\Delta t}(a,t)$ and $D_{\Delta t}y(t) + y(t) = e_{\Delta t}\left(\frac{e^{a\Delta t}-1}{\Delta t},t\right)$ have the same solution for y(t) where the latter equation has no Interval Calculus discrete function mismatch. This equation can be solved without difficulty.

It is thus seen that Interval Calculus has a way of representing Calculus functions in terms of Interval Calculus discrete functions. On the following page is Table 5.16-1, a listing of some Interval Calculus function identities.

Table 5.16-1 Some Interval Calculus Identity Functions

| No. | Calculus Function $t = 0, \Delta t, 2\Delta t, 3\Delta t, \ldots$ | Equivalent Interval Calculus Identity Function $t = 0, \Delta t, 2\Delta t, 3\Delta t, \ldots$ |
|---|---|---|
| 1 | $e^{at}$ | $e_{\Delta t}(\frac{e^{a\Delta t}-1}{\Delta t},t)$ |
| 2 | sinbt | $e_{\Delta t}(\frac{\cos b\Delta t-1}{\Delta t}, t)\sin_{\Delta t}(\frac{\tan b\Delta t}{\Delta t}, t)$ |
| 3 | cosbt | $e_{\Delta t}(\frac{\cos b\Delta t-1}{\Delta t}, t)\cos_{\Delta t}(\frac{\tan b\Delta t}{\Delta t}, t)$ |
| 4 | $e^{at}\sin bt$ | $e_{\Delta t}(\frac{e^{a\Delta t}\cos b\Delta t-1}{\Delta t}, t)\sin_{\Delta t}(\frac{\tan b\Delta t}{\Delta t}, t)$ |
| 5 | $e^{at}\cos bt$ | $e_{\Delta t}(\frac{e^{a\Delta t}\cos b\Delta t-1}{\Delta t}, t)\cos_{\Delta t}(\frac{\tan b\Delta t}{\Delta t}, t)$ |
| 6 | sinhbt | $e_{\Delta t}(\frac{\cosh b\Delta t-1}{\Delta t}, t)\sinh_{\Delta t}(\frac{\tanh b\Delta t}{\Delta t}, t)$ |
| 7 | coshbt | $e_{\Delta t}(\frac{\cosh b\Delta t-1}{\Delta t}, t)\cosh_{\Delta t}(\frac{\tanh b\Delta t}{\Delta t}, t)$ |

A more complete listing of Interval Calculus identities is found in the Appendix in Table 5, Interval Calculus Equations and Identities.

Calculus functions, when used with discrete variables and Z Transforms, can not be efficiently manipulated using Calculus derivatives and integrals. However, Interval Calculus identities when substituted for their equivalent Calculus functions can be efficiently manipulated using Interval Calculus discrete derivatives and discrete integrals. As a result, additional useful information can be conveniently obtained. Some functions that can be derived from an Interval Calculus identity are listed below.

Where h(t) is an Interval Calculus discrete identity function substituted for its equivalent Calculus function:

$t = 0, \Delta t, 2\Delta t, 3\Delta t, \ldots$ (5.16-3)

$\Delta h(t) = \Delta t D_{\Delta t} h(t) = h(t+\Delta t) - h(t)$ , The change of h(t) between successive values of t h(t) forward difference (5.16-4)

$\Delta^2 h(t) = \Delta t^2 D_{\Delta t}{}^2 h(t) = h(t+2\Delta t) - 2h(t+\Delta t) + h(t)$ , Second order h(t) forward difference (5.16-5)

$D_{\Delta t} h(t) = \frac{h(t+\Delta t) - h(t)}{\Delta t}$ , The rate of change of h(t) between successive values of t Interval Calculus discrete derivative of h(t) (5.16-6)

$D_{\Delta t}{}^2 h(t) = \frac{h(t+2\Delta t) - 2h(t+\Delta t) + h(t)}{\Delta t^2}$ , Interval Calculus second order discrete derivative of h(t) (5.16-7)

$$\Delta t \int_{t=t_1}^{t_2} h(t)\,\Delta t = \Delta t \sum_{t=t_1}^{t_2-\Delta t} h(t)\,\Delta t \quad , \quad \text{The area under the sample and hold shaped waveform of } h(t) \qquad (5.16\text{-}8)$$

between $t = t_1$ and $t = t_2$ obtained from the discrete integral of $h(t)$

Considering the above discussion, one might question as to whether Interval Calculus identity functions could be used instead of Calculus functions with Z Transforms. Due to their equivalence, such a substitution should be possible. It is. On the following page is Table 5.16-2, a listing of the Z Transforms of many Interval Calculus functions. The derivation of the Z Transforms of Interval Calculus functions is not difficult. They are easily derived from previously derived $K_{\Delta t}$ Transforms using the $K_{\Delta t}$ Transform to Z Transform Conversion equation rewritten below:

$$Z[f(t)] = \frac{z}{T} K_{\Delta t}[f(t)]\Big|_{\substack{s = \frac{z-1}{T} \\ T = \Delta t \\ t = nT,\ n = 0,1,2,3,\ldots}} \qquad Z[f(t)] = F(z) \quad \text{Z Transform} \qquad (5.16\text{-}9)$$

$$K_{\Delta t}[f(t)] = f(s) \quad K_{\Delta t} \text{ Transform}$$

**Table 5.16-2** **Z Transforms of Some Interval Calculus Functions**

| # | F(t) Interval Calculus Functions<br>$t = 0, \Delta t, 2\Delta t \ldots$ | $K_{\Delta t}$ Transform of F(t)<br>$t = 0, \Delta t, 2\Delta t, \ldots$ | Z Transform of F(t)<br>$T = \Delta t,\ \ t = 0, T, 2T, \ldots$ |
|---|---|---|---|
| **1** | $[t]_{\Delta t}^{n}$<br>or<br>$\prod_{m=1}^{n}(t-[m-1]\Delta t)$<br>$n = 1,2,3,\ldots$ | $\frac{n!}{s^{n+1}}$ | $\frac{n!T^n z}{(z-1)^{n+1}}$ |
| **2** | $e_{\Delta t}(a,t)$ | $\frac{1}{s-a}$ | $\frac{z}{z-(1+aT)}$ |
| **3** | $te_{\Delta t}(a,t)$ | $\frac{1+a\Delta t}{(s-a)^2}$ | $\frac{(1+aT)Tz}{z^2-2(1+aT)z+(1+aT)^2}$ |
| **4** | $[t]_{\Delta t}^{n}\ e_{\Delta t}(a,t)$<br>or<br>$e_{\Delta x}(a,t)\ \prod_{m=1}^{n}(t-[m-1]\Delta t)$<br>$n=1,2,3,\ldots$ | $\frac{(1+a\Delta t)^n\, n!}{(s-a)^{n+1}}$ | $\frac{n!(1+aT)^n T^n z}{[z-(1+aT)]^{n+1}}$ |
| **5** | $\sin_{\Delta t}(b,t)$ | $\frac{b}{s^2+b^2}$ | $\frac{bTz}{z^2-2z+(1+b^2T^2)}$ |

| # | F(t)<br>Interval Calculus Functions<br><br>t = 0, Δt, 2Δt… | $K_{\Delta t}$ Transform<br>of F(t)<br><br>t = 0, Δt, 2Δt, … | Z Transform<br>of F(t)<br><br>T = Δt,   t = 0, T, 2T, … |
|---|---|---|---|
| **6** | $t\sin_{\Delta t}(b,t)$ | $\frac{2bs + b\Delta t(s^2-b^2)}{(s^2+b^2)^2}$ | $\frac{T^2bz[z^2 - (1+b^2T^2)]}{[z^2-2z+(1+b^2T^2)]^2}$ |
| **7** | $\cos_{\Delta t}(b,t)$ | $\frac{s}{s^2+b^2}$ | $\frac{z(z-1)}{z^2-2z+(1+b^2T^2)}$ |
| **8** | $t\cos_{\Delta t}(b,t)$ | $\frac{(s^2-b^2)-2\Delta tb^2s}{(s^2+b^2)^2}$ | $\frac{Tz[z^2-2(1+b^2T^2)z+(1+b^2T^2)]}{[z^2-2z+(1+b^2T^2)]^2}$ |
| **9** | $e_{\Delta t}(a,t)\sin_{\Delta t}(\frac{b}{1+a\Delta t},t)$<br>$a \neq -\frac{1}{\Delta t}$ | $\frac{b}{(s-a)^2+b^2}$ | $\frac{bTz}{z^2-2[1+aT]z+[(1+aT)^2+b^2T^2]}$ |
| **10** | $e_{\Delta t}(a,t)\cos_{\Delta t}(\frac{b}{1+a\Delta t},t)$<br>$a \neq -\frac{1}{\Delta t}$ | $\frac{s-a}{(s-a)^2+b^2}$ | $\frac{z^2- (1+aT)z}{z^2-2[1+aT]z+[(1+aT)^2+b^2T^2]}$ |
| **11** | $\sinh_{\Delta t}(b,t)$ | $\frac{b}{s^2-b^2}$ | $\frac{bTz}{z^2-2z+(1-b^2T^2)}$ |
| **12** | $t\sinh_{\Delta t}(b,t)$ | $\frac{2bs + b\Delta t(s^2+b^2)}{(s^2-b^2)^2}$ | $\frac{T^2bz[z^2 - (1-b^2T^2)]}{[z^2-2z+(1-b^2T^2)]^2}$ |
| **13** | $\cosh_{\Delta t}(b,t)$ | $\frac{s}{s^2-b^2}$ | $\frac{z(z-1)}{z^2-2z+(1-b^2T^2)}$ |
| **14** | $t\cosh_{\Delta t}(b,t)$ | $\frac{(s^2+b^2)+2\Delta tb^2s}{(s^2-b^2)^2}$ | $\frac{Tz[z^2-2(1-b^2T^2)z+(1-b^2T^2)]}{[z^2-2z+(1-b^2T^2)]^2}$ |

| # | F(t)<br>Interval Calculus Functions<br><br>**t = 0, Δt, 2Δt…** | $K_{\Delta t}$ Transform of F(t)<br><br>**t = 0, Δt, 2Δt, …** | Z Transform of F(t)<br><br>**T = Δt, t = 0, T, 2T, …** |
|---|---|---|---|
| **15** | $K_{\Delta t}[f(t-n\Delta t)U(t-n\Delta t)]$<br><br>$U(t-n\Delta t) = \begin{Bmatrix} 1 & t \geq n\Delta t \\ 0 & t < n\Delta t \end{Bmatrix}$<br><br>n = 0,1,2,3,…<br><br>Unit Step Function | $(1+s\Delta t)^{-n}K_{\Delta t}[f(t)]$ | $z^{-n}Z[f(t)]$<br><br>n = 0,1,2,3,… |
| **16** | $U(t-n\Delta t)-U(t-[n+1]\Delta t))$<br>n = 0,1,2,3,…<br><br>Unit Amplitude Pulse<br><br>Pulse Interval = Δt<br>Pulse Ampitude = 1<br>Pulse initiation at t = nΔt<br>Pulse ends at t = (n+1)Δt | $(1+s\Delta t)^{-n-1}\Delta t$<br>or<br>$(1+s\Delta t)^{-(\frac{t+\Delta t}{\Delta t})}\Delta t$<br>or<br>$\frac{1}{s}[(1+s\Delta t)^{-n}-(1+s\Delta t)^{-n-1}]$<br>t = nΔt | $Tz^{-n-1}$<br><br>n = 0,1,2,3,… |
| **17** | f(t) | $\int_{0}^{\infty}{}_{\Delta t}\,(1+s\Delta t)^{-(\frac{t+\Delta t}{\Delta t})}f(t)\Delta t$ | $\frac{1}{T}\int_{0}^{\infty}{}_{T}\,z^{-\frac{t}{\Delta t}}f(t)\Delta t$ |
| **18** | tf(t) | $-(1+s\Delta t)\frac{d}{ds}K_{\Delta t}[f(t)] - \Delta t K_{\Delta t}[f(t)]$ | $-Tz\frac{d}{dz}Z[f(t)]$<br><br>or<br><br>$-Tz^2\frac{d}{dz}[\frac{Z[f(t)]}{z}] - Tz[f(t)]]$ |

| # | F(t)<br>Interval Calculus Functions<br><br>$t = 0, \Delta t, 2\Delta t\ldots$ | $K_{\Delta t}$ Transform<br>of F(t)<br><br>$t = 0, \Delta t, 2\Delta t, \ldots$ | Z Transform<br>of F(t)<br><br>$T = \Delta t, \quad t = 0, T, 2T, \ldots$ |
|---|---|---|---|
| **19** | $D_{\Delta t}f(t)$ | $sK_{\Delta t}[f(t)] - f(0)$ | $\frac{z-1}{T} Z[f(t)] - \frac{z}{T} f(0)$<br><br>$D_T f(t) = \frac{f(t+T) - f(t)}{T}$ |
| **20** | $D^n{}_{\Delta t}f(t)$ | $s^n K_{\Delta t}[f(t)] - s^{n-1}D^0{}_{\Delta t}f(0) - s^{n-2}D^1{}_{\Delta t}f(0) - s^{n-3}D^2{}_{\Delta t}f(0) - \ldots - s^0 D^{n-1}{}_{\Delta t}f(0)$<br><br>or<br><br>$s^n K_{\Delta t}[f(t)] - \sum_{m=1}^{n} s^{n-m}D_{\Delta t}{}^{m-1}f(0)$<br><br>$n = 1,2,3,\ldots$ | $(\frac{z-1}{T})^n Z[f(t)] - \sum_{m=1}^{n} \frac{z}{T} (\frac{z-1}{T})^{n-m}D_{\Delta t}{}^{m-1}f(0)$<br><br>$n = 1,2,3,\ldots$ |
| **21** | ${}_{\Delta t}\int_0^t f(t)\Delta t$ | $\frac{1}{s} K_{\Delta t}[f(t)]$ | $\frac{T}{z-1} Z[f(t)]$ |

$K_{\Delta t}$ Transform to Z Transform Conversion

$$Z[f(t)] = \frac{z}{T} K_{\Delta t}[f(t)]\Big|_{s = \frac{z-1}{T}}$$

$T = \Delta t$

$t = nT, \quad n = 0,1,2,3,\ldots$

$Z[f(t)] = F(z)$ Z Transform

$K_{\Delta t}[f(t)] = f(s)$ $K_{\Delta t}$ Transform

Interval Calculus discrete functions are used to represent sample and hold shaped waveforms even when used with Z Transforms. Note Fig #1 and Fig #2 below.

Representation of the Inverse Z Transform using Calculus functions and Interval Calculus functions

Fig #1

Inverse Z Transform
Usual Calculus function representation

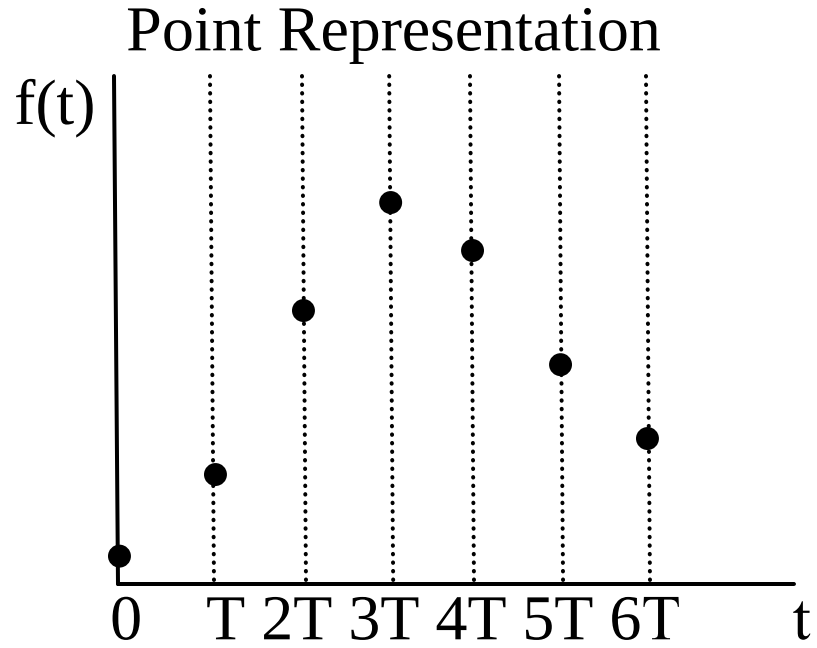

$T = \Delta t$

Fig #2

Inverse Z Transform
Interval Calculus identity function representation

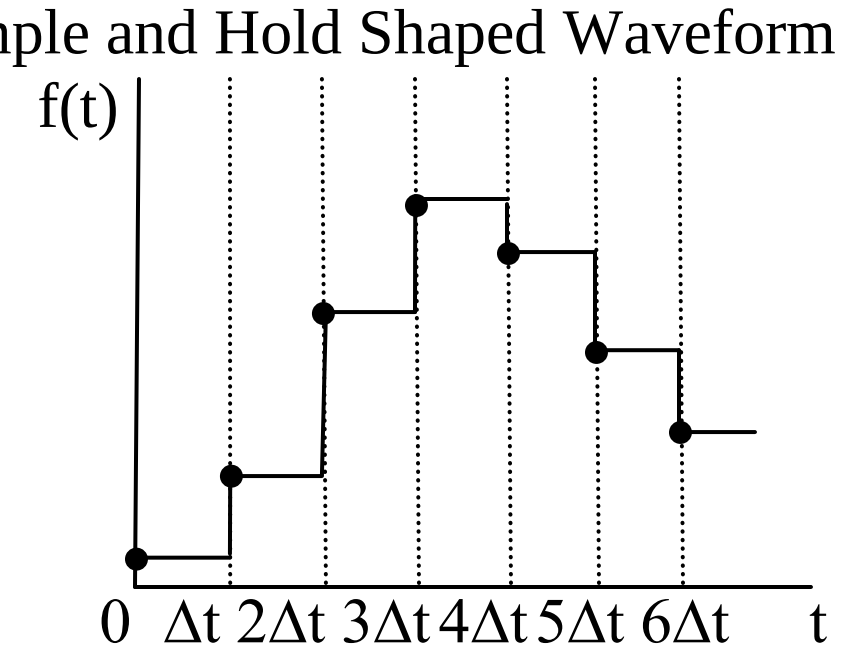

Note – Finding the area under sample and hold shaped waveforms is possible using Interval Calculus discrete (variable) integration

To demonstrate the use of Interval Calculus functions with Z Transforms, two examples are presented below, Example 5.16-1 and Example 5.16-2

Example 5.16-1 Analysis of a sampled data feedback system

For the following sampled data feedback system find the following:

1. The output system response, c(t), at t = .5
2. The change in value of c(t) from t = .1 to t = .2 , $\Delta c(t)$ at t = .1
3. The rate of change of c(t) from t = .1 to t = .2, $D_{\Delta t}c(t)$, at t = .1
4. The area under c(t) from t = 0 to t = 1

The system is initially passive and $T = \Delta t = .1$, $t = 0, .1, .2, .3, \ldots$

Find c(t) in terms of Interval Calculus discrete functions using Z Transforms.

Fig #3 Sampled Data Feedback System Diagram

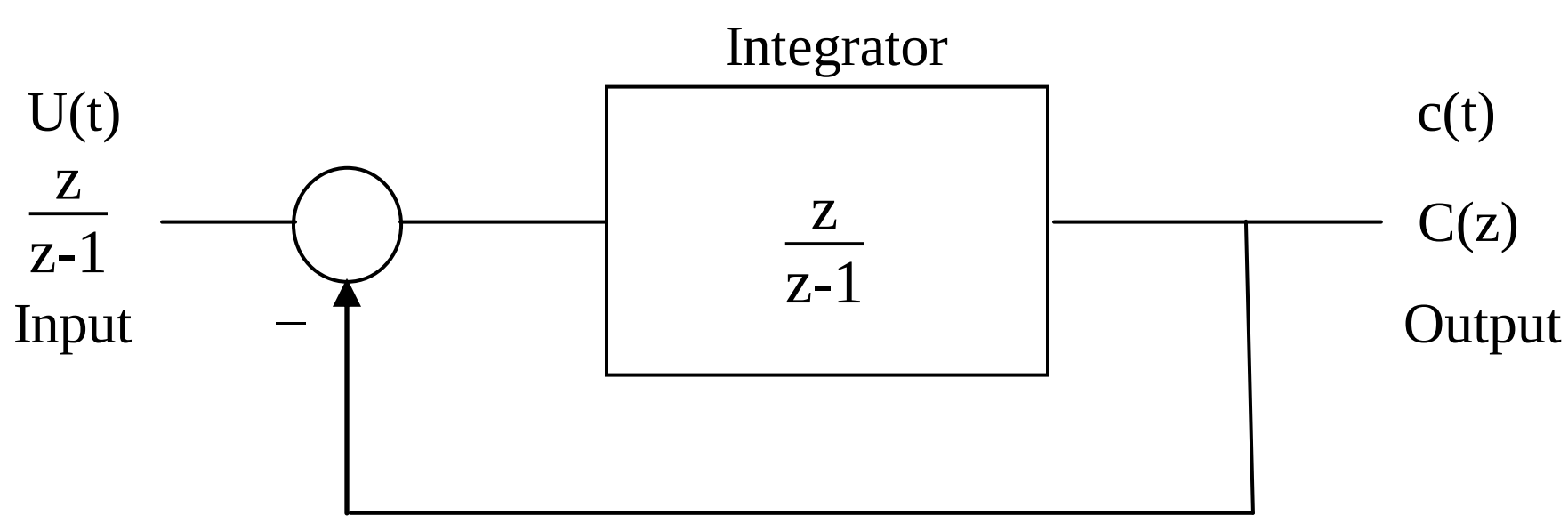

From Fig #3

1. Find c(t) at t = .5

$$C(z) = \frac{z}{z-1}\,\frac{\frac{z}{z-1}}{1+\frac{z}{z-1}} = \frac{z}{(z-1)(z-1+z)} = \frac{z}{2(z-1)(z-.5)} = \frac{.5z}{(z-1)(z-.5)} \qquad 1)$$

$$C(z) = \frac{.5z}{(z-1)(z-.5)} \qquad 2)$$

$$C(z) = \frac{.5z}{(z-1)(z-.5)} = \frac{A}{z-1} + \frac{B}{z-.5} = \frac{Az-.5A+Bz-B}{(z-1)(z-.5)} = \frac{(A+B)z-(.5A+B)}{(z-1)(z-.5)} \qquad 3)$$

From Eq 3

$A + B = .5$ 4)

$.5A + B = 0$ 5)

From Eq 4 and Eq 5

$A - .5A = .5$ 6)

$A = 1$ 7)

$B = -.5$ 8)

From Eq 3, Eq 7, and Eq 8

$$C(z) = \frac{1}{z-1} - \frac{.5}{z-.5} \qquad 9)$$

$$C(z) = z^{-1}\left[\frac{z}{z-1}\right] - z^{-1}\left[\frac{.5z}{z-.5}\right] \qquad 10)$$

Finding the Inverse Z Transform of Eq 10 using Table 5.16-1

$$Z[e_{\Delta t}(a,t)] = \frac{z}{z-(1+aT)} \qquad 11)$$

$T = \Delta t = .1$ 12)

$$Z[e_{.1}(a,t)] = \frac{z}{z-(1+.1a)} = \frac{.5z}{z-.5} \qquad 13)$$

$1+.1a = .5$ 14)

$a = -5$ 15)

Finding c(t) from Eq 10, Eq 13, and Eq 15

$c(t) = U(t-.1) - .5e_{.1}(-5,t-.1)U(t-.1)\ , \quad t = .1, .2, .3, \ldots$ 16)

Find the value of c(t) at t = .5

Using the computer program DXFUN11 to calculate the function, $e_{\Delta t}(a,t)$

$c(.5) = 1 - .5e_{.1}(-5,.4) = 1 - .0625 = .9375$ 17)

**c(.5) = .96875** 18)

Checking Eq 18

From Eq 2

Dividing

$$C(z) = \frac{.5z}{(z-1)(z-.5)} = \frac{.5z}{z^2-1.5z+.5} = 0z^{-0}+.5z^{-1}+.75z^{-2}+.875z^{-3}+.9375z^{-4}+\underline{.96875}z^{-5}+\ldots \qquad 19)$$

$t = .5$ 20)

c(t) for t = .5 is the coefficient of the sixth series term of Eq 19

$c(.5) = .96875$ 21)

Good check

<u>2. Find the change in value of c(t) from t = .1 to t = .2 , Δc(t) at t = .1</u>

Rewriting Eq 16

$c(t) = U(t-.1) - .5e_{.1}(-5,t-.1)U(t-.1)\ , \quad t = .1, .2, .3, \ldots$ 22)

The change of c(t) between successive values of t is the forward difference, Δc(t).

$\Delta c(t) = \Delta t D_{\Delta t}\, c(t)$ 23)

For the difference of c(t) between t = .1 and t = .2 where Δt = T = .1

$\Delta c(t) = .1D_{.1}\,[1 - .5e_{.1}(-5,t-.1)]$ 24)

$\Delta c(.1) = .1D_{.1}\ [1 - .5e_{.1}(-5,t-.1)] = .1(-5)(-.5)e_{.1}(-5,.1-.1) = .25e_{.1}(-5,0) = .25(1) = .25$ 25)

$\mathbf{\Delta c(.1) = c(.2) - c(.1) = .25}$ 26)

Checking Eq 26

$\Delta c(.1) = c(.2) - c(.1) = [U(.2-.1) - .5e_{.1}(-5,.2-.1)U(.2-.1)] - [U(.1-.1) - .5e_{.1}(-5,.1-.1)U(.1-.1)]$ 27)

$c(t) = c(.2) - c(.1) = [1 - .5e_{.1}(-5,.1)] - [1 - .5e_{.1}(-5,0)]$ 28)

Using the computer program, DXFUN11 to calculate the function, $e_{\Delta t}(a,t)$

$c(t) = c(.2) - c(.1) = 1 - .5(.5) - 1 + .5(1) = .5 - .25 = .25$ 29)

$c(t) = c(.2) - c(.1) = .25$ 30)

Good check

3. Find the rate of change of c(t) from t = .1 to t = .2 , $D_{\Delta t}\,c(t)$ at t = .1

Rewriting Eq 16

$c(t) = U(t-.1) - .5e_{.1}(-5,t-.1)U(t-.1)\ , \quad t = .1, .2, .3, \ldots$ 31)

$$D_{\Delta t}\,c(t) = \frac{c(t+\Delta t) - c(t)}{\Delta t}$$ 32)

The rate of change of c(t) between successive values of t is $D_{\Delta t}c(t)$.

For the rate of change of c(t) between t = .1 and t = .2 where $\Delta t = T = .1$

$D_{.1}c(t) = D_{.1}\ [1 - .5e_{.1}(-5,.t-.1)]$ 33)

Using the computer program, DXFUN11 to calculate the function, $e_{\Delta t}(a,t)$

$D_{.1}c(.1) = (-5)(-.5)e_{.1}(-5,.1-.1) = 2.5e_{.1}(-5,0) = 2.5(1) = 2.5$ 34)

$\mathbf{D_{.1}c(.1) = 2.5}$ 35)

Checking Eq 35

From Eq 31 and Eq 32

For the rate of change of c(t) between t = .1 and t = .2 where $\Delta t = T = .1$

$$D_{.1}c\,(.1) = \frac{c(.2) - c(.1)}{.1} = \frac{[U(.2\text{-}.1) - .5e_{.1}(\text{-}5,.2\text{-}.1)U(.2\text{-}.1)] - [U(.1\text{-}.1) - .5e_{.1}(\text{-}5,.1\text{-}.1)U(.1\text{-}.1)]}{.1} \qquad 36)$$

$$D_{.1}c\,(.1) = \frac{c(.2) - c(.1)}{.1} = \frac{[1 - .5e_{.1}(\text{-}5,.1)] - [1 - .5e_{.1}(\text{-}5,0)]}{.1} \qquad 37)$$

Using the computer program, DXFUN11 to calculate the function, $e_{\Delta t}(a,t)$

$$D_{.1}c(.1) = \frac{c(.2) - c(.1)}{.1} = \frac{1 - .5(.5) - 1 + .5(1)}{.1} = \frac{.25}{.1} = 2.5 \qquad 38)$$

$$D_{.1}c(.1) = \frac{c(.2) - c(.1)}{.1} = 2.5 \qquad 39)$$

Good check

4. Find the area under the sample and hold shaped waveform of c(t) from t = 0 to t = 1

Rewriting Eq 16

$$c(t) = U(t\text{-}.1) - .5e_{.1}(\text{-}5,t\text{-}.1)U(t\text{-}.1)\,, \quad t = .1, .2, .3, \ldots \qquad 40)$$

$$_{.1}\int_{t=0}^{1} c(t)\,\Delta t = {}_{.1}\int_{t=0}^{1}[\,U(t\text{-}.1) - .5e_{.1}(\text{-}5,t\text{-}.1)U(t\text{-}.1)\,]\,\Delta t = 0 + {}_{.1}\int_{t=.1}^{1}[\,U(t\text{-}.1) - .5e_{.1}(\text{-}5,t\text{-}.1)U(t\text{-}.1)\,]\,\Delta t \qquad 41)$$

$$_{.1}\int_{t=0}^{1} c(t)\,\Delta t = {}_{.1}\int_{t=.1}^{1}\Delta t - .5\,{}_{.1}\int_{t=.1}^{1} e_{.1}(\text{-}5,t\text{-}.1)\,\Delta t = t\,\Big|_{.1}^{1} - .5\,\frac{e_{.1}(\text{-}5,t\text{-}.1)}{\text{-}5}\,\Big|_{.1}^{1} \qquad 42)$$

$$_{.1}\int_{t=0}^{1} c(t)\,\Delta t = +\,[1 - .1] + .1[e_{.1}(\text{-}5,.9) - e_{.1}(\text{-}5,0)] \qquad 43)$$

Using the computer program, DXFUN11 to calculate the function, $e_{\Delta t}(a,t)$

$$_{.1}\int_{t=0}^{1} c(t)\,\Delta t = [1 - .1] + .1\,[\,.001953125 - 1] = .8001953125 \qquad 44)$$

$$.1\int_{t=0}^{1} c(t)\,\Delta t = .8001953125 \qquad 45)$$

Checking Eq 45

$$.1\int_{t=0}^{1} c(t)\,\Delta t = {}_{.1}\sum_{t=0}^{.9} c(t)\,\Delta t = .1\ {}_{.1}\sum_{t=0}^{.9} [U(t-.1) - .5e_{.1}(-5,t-.1)U(t-.1)] = .1\ {}_{.1}\sum_{t=.1}^{.9} [1 - .5e_{.1}(-5,t-.1)] \qquad 46)$$

Using the computer program, DXFUN11 to calculate the function, $e_{\Delta t}(a,t)$

$$.1\int_{t=0}^{1} c(t)\,\Delta t = .1(0 + .5 + .75 + .875 + .9375 + .96875 + .984375 +.9921875 + .99609375 + .998046875) \qquad 47)$$

$$.1\int_{t=0}^{1} c(t)\,\Delta t = .8001953125 \qquad 48)$$

Good check

The previous example, Example 5.16-1, was an analysis of a sampled data feedback system using Z Transforms with Interval Calculus discrete functions. The following example, Example 5.16-2, is the calculation of a second order difference equation using Z Transforms with Interval Calculus discrete functions.

Example 5.16-2  The solution of a second order difference equation

For the following second order difference equation, y(t+.4) +.1y(t+.2) + .2y(t) = 5 find the following:

1. y(t) at t = 1
2. The change in value of y(t) from t = 1 to t = 1.2 , Δy(t) at t = 1
3. The rate of change of y(t) from t = 1 to t = 1.2, $D_{.2}y(t)$, at t = 1
4. The area under y(t) from t = 0 to t = 1.2, $_{.2}\int_{0}^{1.2} y(t)\,\Delta t$

The initial conditions are y(0) = 0, y(.2) = .5 and T = Δt = .2,   t = 0, .2, .4, .6, …

Find y(t) in terms of Interval Calculus discrete functions using Z Transforms.

$$y(t+.4) +.1y(t+.2) + .2y(t) = 5\ ,\quad y(0) = 0,\ y(.2) =.5 \qquad 1)$$

Find the Z Transform of Eq 1

$$Y(z) = Z[y(t)] \qquad 2)$$

$$Z[y(t+T)] = zY(z) - zy(0) \qquad 3)$$

$$Z(y(t+2T)] = z^2Y(z) - z^2y(0) - zy(T) \qquad 4)$$

$$Z(1) = \frac{z}{z-1} \qquad 5)$$

$$T = .2 \qquad 6)$$

Taking the Z Transform of Eq 1 using Eq 2 thru Eq 6

$$z^2Y(z) - z^2y(0) - zy(.2) + .1[zY(z) - zy(0)] + .2Y(z) = \frac{5z}{z-1} \qquad 7)$$

Simplifying

$$z^2Y(z) - .5z + .1zY(z) + .2Y(z) = \frac{5z}{z-1} \qquad 8)$$

$$(z^2+.1z+.2)Y(z) = \frac{5z}{z-1} + .5z = \frac{.5z^2 + 4.5z}{z-1} = \frac{.5z(z + 9)}{z-1} \qquad 9)$$

$$Y(z) = \frac{.5z(z + 9)}{(z-1)(z^2+.1z+.2)} \qquad 10)$$

Finding the Inverse Z Transform of Eq 10

$$Y(z) = \frac{.5z(z + 9)}{(z-1)(z^2+.1z+.2)} = \frac{A}{z-1} + \frac{Bz+C}{z^2+.1z+.2} \qquad 11)$$

From Eq 11

$$A = \frac{.5z(z + 9)}{(z^2+.1z+.2)}\Big|_{z=1} \qquad 12)$$

$$A = \frac{.5(1)(10)}{1+.1+.2} = \frac{.5(1)(10)}{1+.1+.2} = \frac{5}{1.3} \qquad 13)$$

$$A+B = .5 \qquad 14)$$

$$B = .5 - \frac{5}{1.3} = -\frac{4.35}{1.3} \qquad 15)$$

$$.2A - C = 0 \qquad 16)$$

$$C = .2A = .2(\frac{5}{1.3}) = \frac{1}{1.3} \qquad 17)$$

Substituting Eq 12, Eq 15, and Eq 17 into Eq 11

$$Y(z) = \frac{5}{1.3}\frac{1}{z-1} - \frac{1}{1.3}\frac{4.35z-1}{z^2+.1z+.2} \qquad 18)$$

Changing the form of Eq 18

$$Y(z) = \frac{5}{1.3} z^{-1}[\frac{z}{z-1}] - \frac{4.35}{1.3}\frac{z}{z^2+.1z+.2} + \frac{1}{1.3} z^{-1} [\frac{z}{z^2+.1z+.2}] \qquad 19)$$

From Table 5.16-2

$$Z[e_{\Delta t}(a,t)\sin_{\Delta t}(\frac{b}{1+a\Delta t},t)] = \frac{bTz}{z^2-2[1+aT]z+[(1+aT)^2+b^2T^2]} \qquad 20)$$

From Eq 19, Eq 20, and Eq 6

$$-2(1+aT) = -2(1+.2a) = .1 \qquad 21)$$

$$1+.2a = -.05 \qquad 22)$$

$$a = -\frac{1.05}{.2} = -5.25 \qquad 23)$$

$$(1+aT)^2+b^2T^2 = [1+(-5.25)(.2)]^2 + .04b^2 = .2 \qquad 24)$$

$$b^2 = \frac{.2 - (-.05)^2}{.04} = 4.9375 \qquad 25)$$

$$b = 2.2220486 \qquad 26)$$

$$bT = 2.2220486*.2 = .44440972 \qquad 27)$$

$$\frac{b}{1+aT} = \frac{2.2220486}{1+(-5.25).2} = -44.440972 \qquad 28)$$

From Eq 20, Eq 23, Eq 26, and Eq 27

$$Z[e_{.2}(-5.25,t)\sin_{.2}(-44.440972,t)] = \frac{.44440972z}{z^2+.1z+.2} \qquad 29)$$

From Eq 19, Eq 27

$$Y(z) = \frac{5}{1.3} z^{-1}[\frac{z}{z-1}] - \frac{4.35}{1.3(.44440972)} \frac{.44440972z}{z^2+.1z+.2} + \frac{1}{1.3(.44440972)} z^{-1} \frac{.44440972z}{z^2+.1z+.2} \qquad 30)$$

$$Y(z) = \frac{5}{1.3} z^{-1}[\frac{z}{z-1}] - 7.52943442 \frac{.44440972z}{z^2+.1z+.2} + 1.73090446\, z^{-1} \frac{.44440972z}{z^2+.1z+.2} \qquad 31)$$

Finding the Inverse Z Transform of Eq 31 using Eq 29

$$\mathbf{y(t) = \frac{5}{1.3} U(t-.2) - 7.52943442\, e_{.2}(-5.25,t)\sin_{.2}(-44.440972,t) +}$$
$$\mathbf{1.73090446\, e_{.2}(-5.25,t-.2)\sin_{.2}(-44.440972,t-.2)U(t-.2)} \qquad 32)$$

<u>1. Find y(t) at t = 1</u>

Using Eq 32 find y(t) where t = 1

$$y(1) = \frac{5}{1.3} - 7.52943442\, e_{.2}(-5.25,1)\sin_{.2}(-44.440972,1) +$$
$$1.73090446\, e_{.2}(-5.25,.8)\sin_{.2}(-44.440972,.8) \qquad 33)$$

Using the computer program, DXFUN11, to calculate the functions, $e_{\Delta t}(a,t)$ and $\sin_{\Delta t}(b,t)$

$$y(1) = \frac{5}{1.3} - 7.52943442\,(\,.01515437\,) + 1.73090446\,(.01733197) \qquad 34)$$

**y(1) = 3.76205** 35)

Checking Eq 35

Rewriting Eq 1 as a recursion equation

$$y(t+.4) = 5 - .1y(t+.2) - .2y(t) \quad , \quad y(0) = 0,\ y(.2) = .5 \qquad 36)$$

Using a computer to find y(.4), y(.6), y(.8), …

**<u>Table 5.16-3</u>**

| t | y(t) | |
|---|---|---|
| 0 | 0 | |
| .2 | .5 | |
| .4 | 4.95 | |
| .6 | 4.405 | |
| .8 | 3.5695 | |
| **1.0** | **3.76205** | |
| 1.2 | 3.909895 | Good check |

2. <u>Find the change in value of y(t) from t = 1 to t = 1.2 , Δy(t) at t = 1</u>

$\Delta y(t) = \Delta t D_{\Delta t} y(t)$ 37)

$\Delta y(t) = \Delta t D_{\Delta t}[\frac{5}{1.3} U(t-.2) - 7.52943442\, e_{.2}(-5.25,t)\sin_{.2}(-44.440972,t) +$
$1.73090446\, e_{.2}(-5.25,t-.2)\sin_{.2}(-44.440972,t-.2)U(t-.2)]$ 38)

for $t \geq .2$

$\Delta y(t) = \Delta t D_{\Delta t}[\frac{5}{1.3} - 7.52943442\, e_{.2}(-5.25,t)\sin_{.2}(-44.440972,t) +$
$1.73090446\, e_{.2}(-5.25,t-.2)\sin_{.2}(-44.440972,t-.2)]$ 39)

$D_{\Delta\tau}[v(\tau)u(\tau)] = v(\tau)D_{\Delta\tau}u(\tau) + D_{\Delta\tau}v(\tau)u(\tau+\Delta\tau)$ 40)

Let

$u(\tau) = \sin_{.2}(-44.440972,\tau)$ 41)

$v(\tau) = e_{.2}(-5.25, \tau)$ 42)

From Eq 40 thru Eq 42

$D_{\Delta\tau}[e_{.2}(-5.25, \tau)\sin_{.2}(-44.440972,\tau)] = e_{.2}(-5.25, \tau)D_{\Delta\tau}\sin_{.2}(-44.440972,\tau) +$
$D_{\Delta\tau}e_{.2}(-5.25, \tau)\sin_{.2}(-44.440972,\tau+\Delta\tau)$ 43)

$D_{\Delta\tau}[e_{.2}(-5.25, \tau)\sin_{.2}(-44.440972,\tau)] = -44.440972 e_{.2}(-5.25, \tau)\cos_{.2}(-44.440972,\tau)$
$- 5.25 e_{.2}(-5.25, \tau)\sin_{.2}(-44.440972,\tau+\Delta\tau)$ 44)

From Eq 39

$\Delta y(t) = - 7.52943442 \Delta t D_{\Delta t}[e_{.2}(-5.25,t)\sin_{.2}(-44.440972,t)]$
$+ 1.73090446 \Delta t D_{\Delta t}[e_{.2}(-5.25,t-.2)\sin_{.2}(-44.440972,t-.2)]$ 45)

From Eq 44 and Eq 45

$\Delta t = .2$ 46)

For the first term $\tau = t$ and $\Delta\tau = \Delta t = .2$
For the second term $\tau = t-.2$ and $\Delta\tau = \Delta(t-.2) = \Delta t = .2$

$\Delta y(t) = - 7.52943442(.2)[-44.440972 e_{.2}(-5.25, t)\cos_{.2}(-44.440972,t)$
$- 5.25 e_{.2}(-5.25, t)\sin_{.2}(-44.440972,t+.2)]$
$+ 1.73090446(.2)[-44.440972 e_{.2}(-5.25, t-.2)\cos_{.2}(-44.440972,t-.2)$
$- 5.25 e_{.2}(-5.25, t-.2)\sin_{.2}(-44.440972,t-.2+.2)]$ 47)

Simplifying Eq 47

$$\mathbf{\Delta y(t) = +66.92307684e_{.2}(-5.25,\ t)\cos_{.2}(-44.440972,t) +7.90590614e_{.2}(-5.25,\ t)\sin_{.2}(-44.440972,t+.2) - 15.38461532e_{.2}(-5.25,\ t-.2)\cos_{.2}(-44.440972,t-.2) - 1.81744968e_{.2}(-5.25,\ t-.2)\sin_{.2}(-44.440972,t)} \quad 48)$$

Find $\Delta y(t)$ for $t = 1$

$$\Delta y(1) = +66.92307684e_{.2}(-5.25,\ 1)\cos_{.2}(-44.440972,1) +7.90590614e_{.2}(-5.25,\ 1)\sin_{.2}(-44.440972,1.2) - 15.38461532e_{.2}(-5.25,\ .8)\cos_{.2}(-44.440972,.8) - 1.81744968e_{.2}(-5.25,\ .8)\sin_{.2}(-44.440972,1) \quad 49)$$

Calculating $e_{\Delta t}(a,t)$, $\sin_{\Delta t}(b,t)$, $\cos_{\Delta t}(b,t)$, $e_{\Delta t}(a,t)\sin_{\Delta t}(b,t)$ and $e_{\Delta t}(a,t)\cos_{\Delta t}(b,t)$ using the program, DXFUN11

$$\Delta y(1) = +66.92307684(-.00950500) +7.90590614(-.0000003125)(-318835.6645) - 15.38461532(.03605000) - 1.81744968(.00000625)(-48493.9882) \quad 50)$$

$$\Delta y(1) = -.63610384+.78771401-.55461538+.55084614 = .14784 \quad 51)$$

$$\mathbf{\Delta y(1) = .14784} \quad 52)$$

Checking Eq 52

From Table 5.16-3

$$y(1) = 3.76205 \quad 53)$$

$$y(1.2) = 3.909895 \quad 54)$$

$$\Delta y(1) = 3.909895-3.76205 = .147845 \quad 55)$$

Good check

3. Find the rate of change of y(t) from t = 1 to t = 1.2, $D_{.2}y(t)$, at t = 1

$$D_{\Delta t}y(t) = \frac{\Delta y(t)}{\Delta t} \quad 56)$$

From Eq 48 and Eq 56

$$D_{.2}y(t) = \frac{\Delta y(t)}{.2} = 5[+66.92307684e_{.2}(-5.25, t)\cos_{.2}(-44.440972,t)$$
$$+7.90590614e_{.2}(-5.25, t)\sin_{.2}(-44.440972,t+.2)$$
$$-15.38461532e_{.2}(-5.25, t-.2)\cos_{.2}(-44.440972,t-.2)$$
$$-1.81744968e_{.2}(-5.25, t-.2)\sin_{.2}(-44.440972,t)] \qquad 57)$$

$$\mathbf{D_{.2}y(t) = +334.6153842e_{.2}(-5.25, t)\cos_{.2}(-44.440972,t)}$$
$$\mathbf{+39.5295307e_{.2}(-5.25, t)\sin_{.2}(-44.440972,t+.2)}$$
$$\mathbf{-76.9230766e_{.2}(-5.25, t-.2)\cos_{.2}(-44.440972,t-.2)}$$
$$\mathbf{-9.0872484e_{.2}(-5.25, t-.2)\sin_{.2}(-44.440972,t)]} \qquad 58)$$

Find $D_{.2}y(1)$

From Eq 52 and Eq 57

$$D_{.2}y(1) = \frac{\Delta y(1)}{.2} = \frac{.14784}{.2} = .7392 \qquad 59)$$

$$\mathbf{D_{.2}y(1) = .7392} \qquad 60)$$

Checking Eq 60

From Eq 55 and Eq 57

$$D_{.2}y(1) = \frac{\Delta y(1)}{.2} = \frac{.147845}{.2} = .7392 \qquad 61)$$

$$D_{.2}y(1) = .739225 \qquad 62)$$

Good check

4. <u>Find the area under y(t) from t = 0 to t = 1.2,</u> $\,_{.2}\int_{0}^{1.2} y(t)\Delta t$

Rewriting Eq 32

$$y(t) = \frac{5}{1.3}U(t-.2) - 7.52943442\, e_{.2}(-5.25,t)\sin_{.2}(-44.440972,t) +$$
$$1.73090446\, e_{.2}(-5.25,t-.2)\sin_{.2}(-44.440972,t-.2)U(t-.2) \qquad 63)$$

Use discrete integration to find $\,_{.2}\int_{0}^{1.2} y(t)\Delta t$

Integrating Eq 63

$$ {}_{.2}\int_{0}^{1.2} y(t)\Delta t = 0 + \frac{5}{1.3}\,{}_{.2}\int_{.2}^{1.2} 1\,\Delta t - 7.52943442\,{}_{.2}\int_{.2}^{1.2} e_{.2}(-5.25,t)\sin_{.2}(-44.440972,t)\Delta t\ + $$

$$ 1.73090446\,{}_{.2}\int_{.2}^{1.2} e_{.2}(-5.25,t\text{-}.2)\sin_{.2}(-44.440972,t\text{-}.2)\Delta t \qquad 64) $$

From Table 6 in the Appendix

$$ {}_{\Delta t}\int e_{\Delta t}(a,t)\sin_{\Delta t}(b,t)\Delta t = \frac{e_{\Delta t}(a,t)}{a^2+b^2(1+a\Delta t)^2}\,[a\sin_{\Delta t}(b,t) - b(1+a\Delta t)\cos_{\Delta t}(b,t)] + k \qquad 65) $$

Using Eq 65 to evaluate Eq 64

$$ {}_{.2}\int_{0}^{1.2} y(t)\Delta t = \frac{5}{1.3}\,t\,\Big|_{.2}^{1.2} - 7.52943442\,\frac{e_{\Delta t}(a,t)}{a^2+b^2(1+a\Delta t)^2}\,[a\sin_{\Delta t}(b,t) - b(1+a\Delta t)\cos_{\Delta t}(b,t)]\,\Big|_{.2}^{1.2} + $$

$$ 1.73090446\,\frac{e_{\Delta t}(a,t)}{a^2+b^2(1+a\Delta t)^2}\,[a\sin_{\Delta t}(b,t\text{-}.2) - b(1+a\Delta t)\cos_{\Delta t}(b,t\text{-}.2)]\,\Big|_{.2}^{1.2} \qquad 66) $$

Substituting into Eq 66 the values from Eq 64

$$ a = -5.25 \qquad 67) $$

$$ b = -44.440972 \qquad 68) $$

$$ \Delta t = .2 \qquad 69) $$

$$ a^2+b^2(1+a\Delta t)^2 = (-5.25)^2+(-44.440972)^2[1-5.25(.2)]^2 = 32.5 \qquad 70) $$

$$ b(1+a\Delta t) = -44.440972[1-5.25(.2)] = 2.2220486 \qquad 71) $$

$$ {}_{.2}\int_{0}^{1.2} y(t)\Delta t = \frac{5}{1.3}\,t\,\Big|_{.2}^{1.2} - \frac{7.52943442(-5.25)}{32.5}\,e_{.2}(-5.25,t)\sin_{.2}(-44.040972,t)\,\Big|_{.2}^{1.2} $$

$$ + \frac{7.52943442(2.2220486)}{32.5}\,e_{.2}(-5.25,t)\cos_{.2}(-44.440972,t)\,\Big|_{.2}^{1.2} $$

$$ + \frac{1.73090446(-5.25)}{32.5}\,e_{.2}(-5.25,t)\sin_{.2}(-44.040972,t\text{-}.2)\,\Big|_{.2}^{1.2} $$

$$ - \frac{1.73090446(2.2220486)}{32.5}\,e_{.2}(-5.25,t)\cos_{.2}(-44.440972,t\text{-}.2)\,\Big|_{.2}^{1.2} \qquad 72) $$

$$\int_{0}^{1.2}{}_{.2}\, y(t)\Delta t = \frac{5}{1.3} t \Big|_{.2}^{1.2} + 1.21629325\, e_{.2}(-5.25,t)\sin_{.2}(-44.440972,t) \Big|_{.2}^{1.2}$$

$$+ .51479289\, e_{.2}(-5.25,t)\cos_{.2}(-44.440972,t) \Big|_{.2}^{1.2}$$

$$- .27960764\, e_{.2}(-5.25,t-.2)\sin_{.2}(-44.040972,t-.2) \Big|_{.2}^{1.2}$$

$$- .11834319\, e_{.2}(-5.25,t-.2)\cos_{.2}(-44.440972,t-.2) \Big|_{.2}^{1.2} \qquad 73)$$

Calculating $e_{\Delta t}(a,t)$, $\sin_{\Delta t}(b,t)$, $\cos_{\Delta t}(b,t)$, $e_{\Delta t}(a,t)\sin_{\Delta t}(b,t)$ and $e_{\Delta t}(a,t)\cos_{\Delta t}(b,t)$ using the program, DXFUN11

$$\int_{0}^{1.2}{}_{.2}\, y(t)\Delta t = \frac{5}{1.3}[1.2-.2] + 1.21629325[-.00498183-.44440972] + .51479289\,[-.00625950+.05]$$

$$- .27960764[.01515437-0] - .11834319[-.00950450-1] = 3.43731 \qquad 74)$$

$$\boldsymbol{\int_{0}^{1.2}{}_{.2}\, y(t)\Delta t = 3.43731} \qquad 75)$$

Checking Eq 75

Using a computer to find A(1.2)

Rewriting Table 3

| t | y(t) |
|---|---|
| 0 | 0 |
| .2 | .5 |
| .4 | 4.95 |
| .6 | 4.405 |
| .8 | 3.5695 |
| 1.0 | 3.76205 |
| 1.2 | 3.909895 |

A(1.2) = .2(0+.5+4.95+4.405+3.5695+3.76205) = 3.43731 76)

A(1.2) = 3.43731 77)

Good check

The previous two examples, Example 5.16-1 and Example 5.16-2, have demonstrated the use of Z Transforms with Interval Calculus discrete variable functions, functions with sample and hold shaped waveforms. Having the the Z Transform solutions in terms of Interval Calculus discrete functions facilitates the the calculation of differences, rates of change, and under curve areas. This same capability is available using $K_{\Delta t}$ Transforms where solutions are also in terms of Interval Calculus discrete functions. In addition to these two specified methods of analysis, there is a third. A typical Z Transform analysis can be performed yielding a solution in terms of Calculus functions that are then replaced by their Interval Calculus identities. Such an analysis is shown below in Example 5.16-3.

Example 5.16-3 The solution of a system using Z Transforms where the solution Calculus functions are replaced by their equivalent Interval Calculus identities

For the system specified below, obtain Interval Calculus discrete functions to represent:

1. y(t)
2. The change in the value of y(t)
3. The rate of change in the value of y(t)
4. The area under the sample and hold shaped curve of y(t) , $\underset{\Delta t}{\int_{t_1}^{t_2}} y(t)\Delta t$

The system shown below is initially passive.

Find y(t) in terms of Interval Calculus discrete functions using a typical Z Transform analysis. Replace the Calculus function solution with its equivalent Interval Calculus identities.

System Diagram

$r(t) = t^2$ Sample and Hold Integrator $y(t) = \frac{t^2}{2} - \frac{t}{2}$

$$Z[r(t)] = \frac{Tz}{(z-1)^2} \longrightarrow \boxed{\frac{T}{z-1}} \longrightarrow Z[y(t)] = \frac{T^2z}{(z-1)^3}$$

Input Output

1. Find y(t)

Note - $\frac{T}{z-1} = [Tz^{-1}][\frac{z}{z-1}]$ = Sample and Hold Integrator

$$Z[y(t)] = [\frac{Tz}{(z-1)^2}][\frac{T}{z-1}] = \frac{T^2z}{(z-1)^3} \qquad 1)$$

Assume the Inverse Z Transform of y(t) has a $t^2$ and t term

$$Z[y(t)] = \frac{T^2z}{(z-1)^3} = \frac{AT^2z(z+1)}{(z-1)^3} + \frac{BTz}{(z-1)^2} \qquad 2)$$

$$Z[y(t)] = \frac{T^2z}{(z-1)^3} = \frac{AT^2z(z+1)}{(z-1)^3} + \frac{BTz(z-1)}{(z-1)^3} = \frac{AT^2z(z+1) + BTz(z-1)}{(z-1)^3} \qquad 3)$$

Finding the constants A,B

$$T^2z = AT^2z^2 + AT^2z + BTz^2 - BTz \quad 4)$$

$$0 = AT^2 + BT \quad 5)$$

$$B = -AT \quad 6)$$

$$T^2 = AT^2 - BT \quad 7)$$

Substituting Eq 6 into Eq 7

$$T^2 = AT^2 + AT^2 \quad 8)$$

$$A = \frac{1}{2} \quad 9)$$

Substituting Eq 9 into Eq 6

$$B = -\frac{T}{2} \quad 10)$$

Substituting Eq 9 and Eq 10 into Eq 2

$$Z[y(t)] = \frac{1}{2}\frac{T^2z(z+1)}{(z-1)^3} - \frac{T}{2}\frac{Tz}{(z-1)^2} \quad 11)$$

Find the Inverse Z Transform of Eq 11

$$\mathbf{y(t) = \frac{t^2}{2} - \frac{Tt}{2}} \quad 12)$$

To find $\Delta y(t)$, $D_{\Delta t}y(t)$, and $\int_{t_1}^{t_2}{}_{\Delta t}\, y(t)\Delta t$ substitute Interval Calculus identities into Eq 12

$$T = \Delta t \quad 13)$$

Interval Calculus identities

$$t^2 = [t]_{\Delta t}^{2} + \Delta t[t]_{\Delta t}^{1} = t(t-\Delta t) + \Delta t t \quad 14)$$

$$t = [t]_{\Delta t}^{1} \quad 15)$$

Substituting Eq 14 and Eq 15 into Eq 12

$$y(t) = \frac{1}{2}[t]_{\Delta t}^{2} + \frac{\Delta t}{2}[t]_{\Delta t}^{1} - \frac{\Delta t}{2}[t]_{\Delta t}^{1} = \frac{1}{2}[t]_{\Delta t}^{2} \qquad 16)$$

$$\mathbf{y(t) = \frac{1}{2}[t]_{\Delta t}^{2} = \frac{t(t-\Delta t)}{2}} \qquad 17)$$

Checking Eq 17

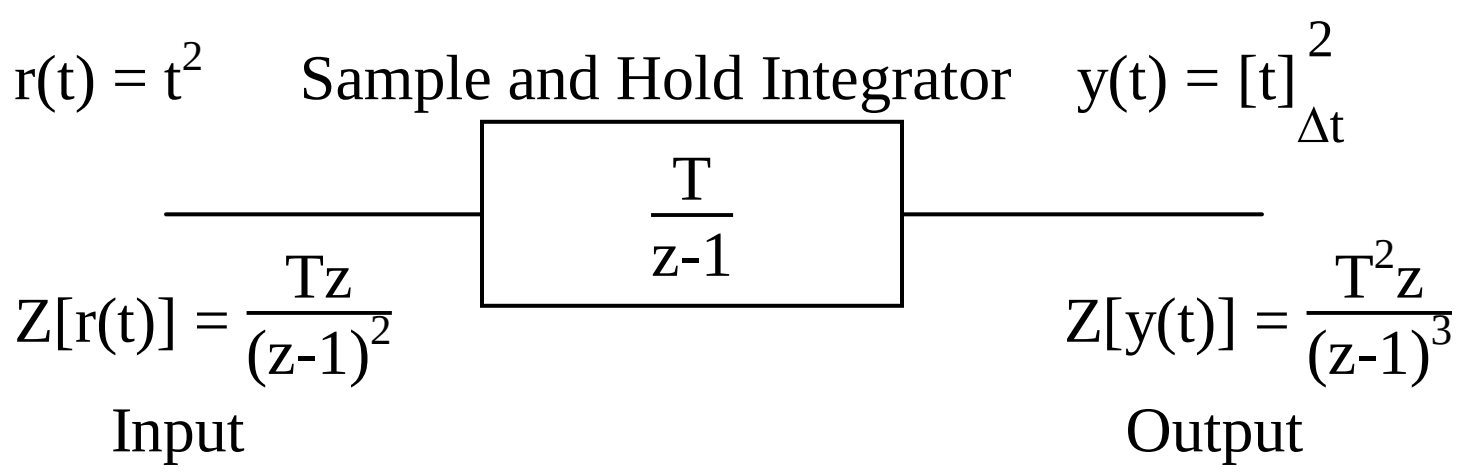

$$Z[y(t)] = \frac{T^2 z}{(z-1)^3} \qquad 18)$$

From Table 5.16-2

$$[t]_{\Delta t}^{n} = \prod_{m=1}^{n}(t-[m-1]\Delta t) \qquad 19)$$

$n = 1,2,3,\ldots$

$$Z[\,[t]_{\Delta t}^{n}\,] = \frac{n!T^n z}{(z-1)^{n+1}} \qquad 20)$$

From Eq 18

$$Z[y(t)] = \frac{1}{2}\frac{2T^2 z}{(z-1)^3} \qquad 21)$$

From Eq 18 thru Eq 20

$$y(t) = \frac{1}{2}[t]_{\Delta t}^{2} = \frac{t(t-\Delta t)}{2} \qquad 22)$$

Good check

2. Find the change in value of y(t) between the discrete values of t

$$\Delta y(t) = \Delta t D_{\Delta t} y(t) = \Delta t D_{\Delta t}[\frac{1}{2}[t]_{\Delta t}^{2}] = \Delta t [t]_{\Delta t}^{1} = \Delta t t \qquad 23)$$

$$\mathbf{\Delta y(t) = \Delta t t} \qquad 24)$$

3. Find the rate of change of y(t) between the discrete values of t

$$D_{\Delta t}y(t) = D_{\Delta t}[t]^{2}_{\Delta t} = D_{\Delta t}[\frac{1}{2}[t]^{2}_{\Delta t}] = [t]^{1}_{\Delta t} = t \qquad 25)$$

$$\mathbf{D_{\Delta t}y(t) = t} \qquad 26)$$

4. Find the area under the sample and hold shaped curve of y(t) , ${}_{\Delta t}\int_{t_1}^{t_2} y(t)\Delta t$

$${}_{\Delta t}\int_{t_1}^{t_2} y(t)\Delta t = {}_{\Delta t}\int_{t_1}^{t_2} \frac{1}{2}[t]^{2}_{\Delta t}\Delta t = \frac{1}{6}[t]^{3}_{\Delta t}\Big|_{t_1}^{t_2} = \frac{t(t-\Delta t)(t-2\Delta t)}{6}\Big|_{t_1}^{t_2} \qquad 27)$$

$$\mathbf{{}_{\Delta t}\int_{t_1}^{t_2} y(t)\Delta t = \frac{t(t-\Delta t)(t-2\Delta t)}{6}\Big|_{t_1}^{t_2}} \qquad 28)$$

The previous three examples, Example 5.16-1 thru Example 5.16-3, have shown the use of Z Transforms with Interval Calculus discrete functions in the solution of discrete variable problems. It is observed that the use of Interval Calculus discrete functions to express y(t) greatly facilitates the calculation of differences, rates of change and under curve areas relating to y(t).

## Section 5.17: The Modified $K_{\Delta t}$ Transform

When analyzing a system using Z Transforms, system responses to a series of equally spaced input impulses are evaluated. Usually the responses are evaluated at the instants of time , t = 0, T, 2T, 3T, … For some systems knowing only the system outputs at these times is adequate for a suitable understanding of the system and its response characteristics. However, in some systems with slower sampling rates, knowing the system output responses only at t = 0, T, 2T, 3T, … may be insufficient and even misleading. For example, consider the block diagram of the system shown below.

Fig #1 Demonstation System

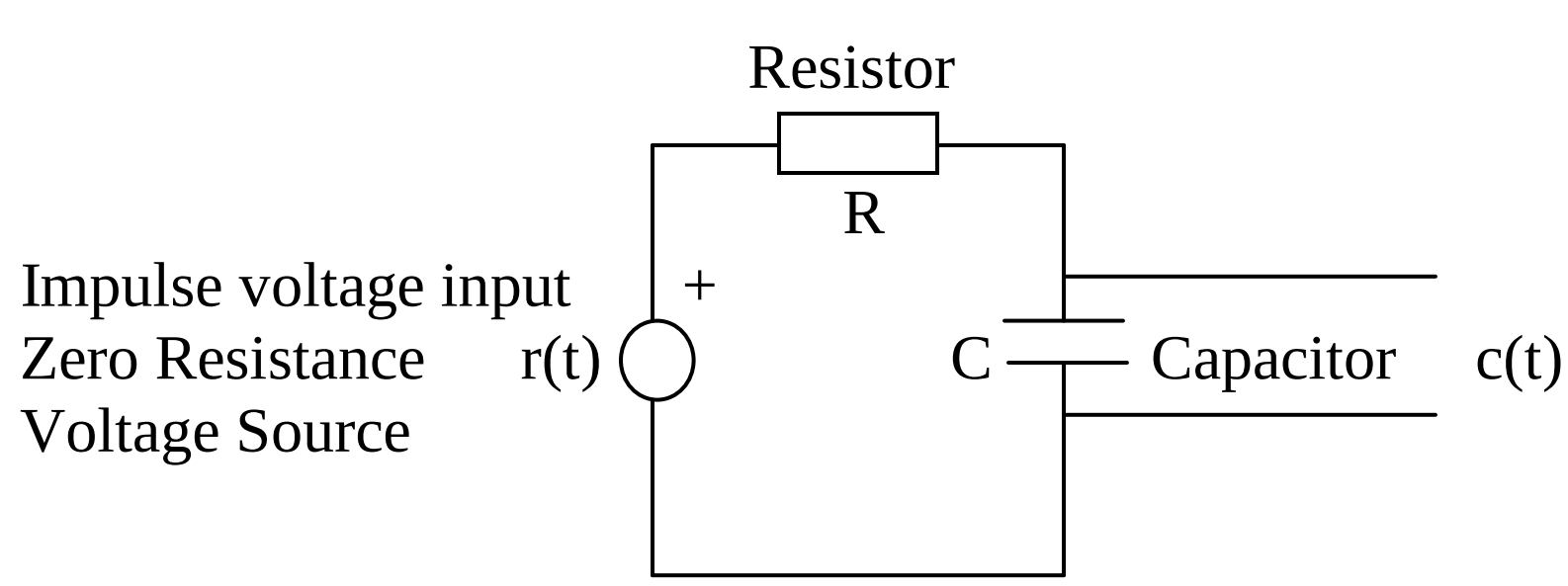

r(t) = | | | | … , a series of input unit area impulses
0 T 2T 3T

c(t) = output voltage across the capacitor

T = interval between successive input voltage impulses
c(0) = 0 , initial condition of c(t)

<u>Comment</u> - In the system diagram of Fig#1, c(t), the voltage across the capacitor, is defined for all values of t even though the system input, r(t), is a discrete input, equally spaced weighted impulses for a Z Transform analysis or equally spaced weighted unit amplitude pulses for a $K_{\Delta t}$ Transform analysis. When using a Z Transform or $K_{\Delta t}$ Transform analysis, c(t) is calculated only at the instants of time defining the system discrete input, r(t). In the calculations and discussions that follow, the system output, c(t), will be defined with all allowable values of t specified.

$$\frac{C(s)}{R(s)} = \frac{\frac{1}{Cs}}{\frac{1}{Cs} + R} = \frac{1}{RCs+1} = \frac{1}{RC}\frac{1}{s + \frac{1}{RC}} = \frac{a}{s+a} \tag{5.17-1}$$

$$\frac{C(s)}{R(s)} = \frac{a}{s+a} \tag{5.17-2}$$

where

$$a = \frac{1}{RC}$$

<u>Fig #2</u> Z Transform of the Demonstation System

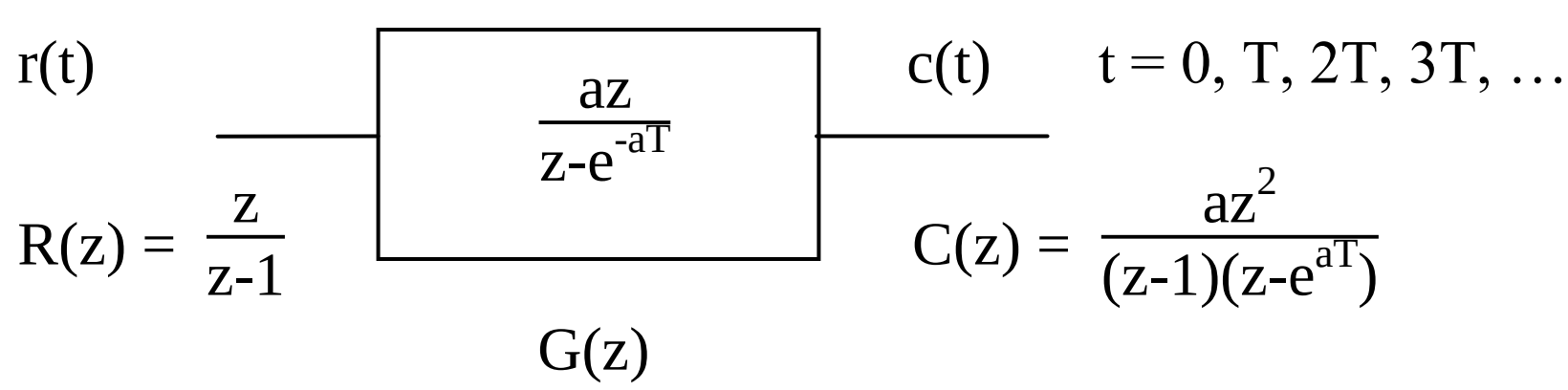

r(t) = | | | | … , a series of input voltage unit area impulses
0 T 2T 3T

From Fig #2

$$C(z) = \frac{a}{1-e^{-aT}}\left(\frac{z}{z-1} - \frac{e^{-aT}z}{z-e^{-aT}}\right) \tag{5.17-3}$$

The Inverse Z Transform of Eq 5.17-3 is:

$$c(t) = \frac{a}{1- e^{-aT}}\left(1 - e^{-aT}e^{-at}\right), \quad t = 0, T, 2T, 3T, \dots \tag{5.17-4}$$

$$aT = \frac{T}{RC} > 1 \tag{5.17-5}$$

From Eq 5.17-4 one would expect the system output, c(t), to look like Curve 1 shown in Fig #3 presented below. However, Curve 1 is deceptive. Looking at the system diagram of Fig #1 it is realized that the capacitor output voltage, c(t), must rise sharply when the energy from each individual voltage impulse is transferred to the capacitor. Following this, capacitor exponential voltage decay would be expected as capacitor current drains back through the system resistor. Observe curve 2 in Fig #3. It should be noted that Eq 5.17-4 does not show any of the capacitor voltage drop due to the capacitor's current discharge through the system resistor. Eq 5.17-4 represents the capacitor voltage only at $t = 0, T, 2T, 3T, \ldots$ just following impulse capacitor charging. Note Fig #3 shown on the following page.

Fig #3 Perceived and Actual System Output

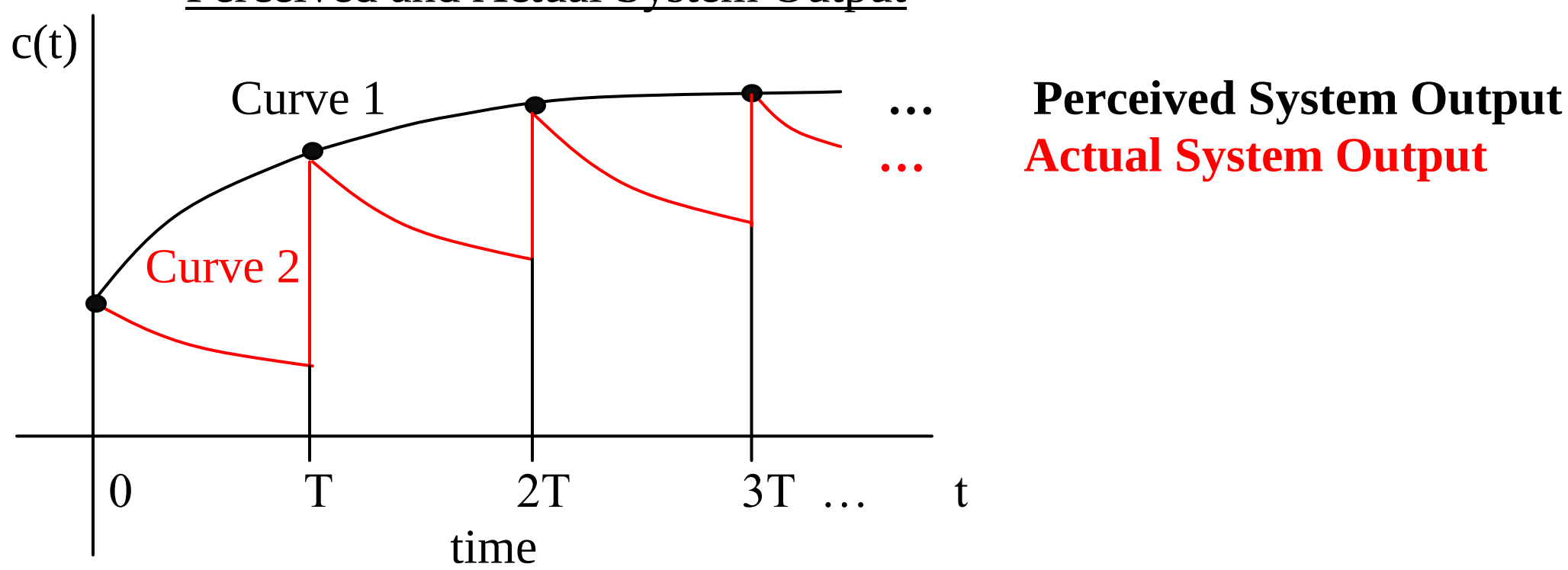

For systems such as this demonstration system, the Modified Z Transform was developed to evaluate the output values occurring between $t = 0, T, 2T, 3T, \ldots$ The Modified Z Transform is well known in the field of sampled-data control system analysis. Its derivation can be found in a number of text books. A recommended text book for finding the derivation of the Modified Z Transform is listed in the Reference Book Listing at the end of this paper. It is the text book listed as "Digital and Sampled-Data Control Systems", by Julius Tou. A Modified $K_{\Delta t}$ Transform has also been developed to evaluate output values occurring between $t = 0, T, 2T, 3T, \ldots$ After a discription and comparison of the Z and $K_{\Delta t}$ Transforms, a derivation of the Modified $K_{\Delta t}$ Transform is given.

Note the following series definitions of the Z Transform and the $K_{\Delta t}$ Transform.

Z Transform

$$Z[f(t)] = f(0)z^0 + f(T)z^{-1} + f(2T)z^{-2} + f(3T)z^{-3} + f(4T)z^{-4} + \ldots \qquad (5.17\text{-}6)$$

where

$f(t) =$ a function of t

$t = 0, T, 2T, 3T, 4T, \ldots$

$T =$ the interval between successive discrete values of t

$Z[f(t)] =$ the Z Transform of the function, f(t)

$e^{snT} =$ unit area impulses initiated at $t = nT$, $n = 0, 1, 2, 3, 4, \ldots$

$K_{\Delta t}$ Transform

$$K_{\Delta t}[f(t)] = f(0)(1+s\Delta t)^{-1}\Delta t + f(\Delta t)(1+s\Delta t)^{-2}\Delta t + f(2\Delta t)(1+s\Delta t)^{-3}\Delta t + f(3\Delta t)(1+s\Delta t)^{-4}\Delta t + \ldots \qquad (5.17\text{-}7)$$

where

$f(t) =$ a function of t

$t = 0, \Delta t, 2\Delta t, 3\Delta t, 4\Delta t, \ldots$

$\Delta t =$ the interval between successive discrete values of t

$K_{\Delta t}[f(t)] =$ the $K_{\Delta t}$ Transform of the function f(t)

$(1+s\Delta t)^{-n}\Delta t =$ unit amplitude pulse of $\Delta t$ duration initiated at $t = n\Delta t$, $n = 1,2,3,4,\ldots$

Clearly from Eq 5.17-6 and Eq 5.17-7 the Z Transform and the $K_{\Delta t}$ Transform are different. However, there is an important relationship. They can both represent the same discrete function, f(t), and as shown in the Z Transform to $K_{\Delta t}$ Transform Conversion Equation previously derived, the Z Transform of f(t) can be easily converted into the $K_{\Delta t}$ Transform of f(t). The Z Transform to $K_{\Delta t}$ Transform Conversion Equation is written below.

The Z Transform to $K_{\Delta t}$ Transform Conversion Equation is:

$$K_{\Delta t}[f(t)] = F(s) = Tz^{-1}\, F(z)\Big|_{z\,=\,1+s\Delta t} \qquad K_{\Delta t}[f(t)] = F(s) \quad K_{\Delta t} \text{ Transform} \qquad (5.17\text{-}8)$$

$\Delta t = T$ sampling period

$t = n\Delta t$ $\qquad Z[f(t)] = F(z)$ $\quad$ Z Transform

$n = 0, 1, 2, 3, \ldots$

<u>Comment</u> – It should be noted that there is a $K_{\Delta t}$ Transform to Z Transform Conversion Equation also.

The $K_{\Delta t}$ Transform to Z Transform Conversion Equation is:

$$Z[f(t)] = F(z) = \frac{1+s\Delta t}{\Delta t}\, F(s)\Big|_{s\,=\,\frac{z-1}{T}} \qquad Z[f(t)] = F(z) \quad \text{Z Transform} \qquad (5.17\text{-}9)$$

$T = \Delta t$ sampling period

$t = nT$ $\qquad K_{\Delta t}[f(t)] = F(s)$ $\quad$ $K_{\Delta t}$ Transform

$n = 0,1,2,3,\ldots$

Note that if Eq 5.17-8 is applied to Eq 5.17-6, Eq 5.17-7 results.

Consider the following example, Example 5.17-1.

<u>Example 5.17-1 The conversion of the Z Transform, Z[t], to the $K_{\Delta t}$ Transform, $K_{\Delta t}[t]$</u>

Convert the Z Transform of t, $Z[t] = \frac{Tz}{(z-1)^2}$, into the $K_{\Delta t}$ Transform of t, $K_{\Delta t}[t]$. Use The Z Transform to $K_{\Delta t}$ Transform Conversion Equation.

$$Z[t] = \frac{Tz}{(z-1)^2}\,, \quad \text{The Z Transform of t} \qquad 1)$$

The $K_{\Delta t}$ Transform to Z Transform Conversion Equation is:

$$K_{\Delta t}[f(t)] = F(s) = Tz^{-1}\, F(z)\Big|_{z\,=\,1+s\Delta t} \qquad K_{\Delta t}[f(t)] = F(s) \quad K_{\Delta t} \text{ Transform} \qquad 2)$$

$\Delta t = T$ sampling period

$t = n\Delta t$ $\qquad Z[f(t)] = F(z)$ $\quad$ Z Transform

$n = 0, 1, 2, 3, \ldots$

Applying Eq 2 to Eq 1

$$K_{\Delta t}[t] = Tz^{-1}Z[t]\Big|_{\substack{z\,=\,1+s\Delta t \\ \Delta t\,=\,T \text{ sampling period} \\ t\,=\,n\Delta t \\ n\,=\,0,\,1,\,2,\,3,\,\ldots}} = Tz^{-1}\frac{Tz}{(z-1)^2}\Big|_{\substack{z\,=\,1+s\Delta t \\ \Delta t\,=\,T \text{ sampling period} \\ t\,=\,n\Delta t \\ n\,=\,0,1,2,3,\ldots}} \qquad 3)$$

$$K_{\Delta t}[t] = \frac{\Delta t^2}{(1+s\Delta t-1)^2} = \frac{\Delta t^2}{s^2\Delta t^2} = \frac{1}{s^2} \qquad 4)$$

Then

$$\mathbf{K_{\Delta t}[t] = \frac{1}{s^2}}\ ,\ \textbf{The } \mathbf{K_{\Delta t}} \textbf{ Transform of t} \qquad 5)$$

Eq 5 is known to be the $K_{\Delta t}$ Transform of t

It is seen in Example 5.17-1 that the inverse transform, $K_{\Delta t}^{-1}[\frac{1}{s^2}]$ and $Z^{-1}[\frac{Tz}{(z-1)^2}]$, both represent the same function, t. However, the manner in which the function, t, assigns a value to each discrete value of time differs for each transform.

For the Z Transform

A series of weighted unit area impulses

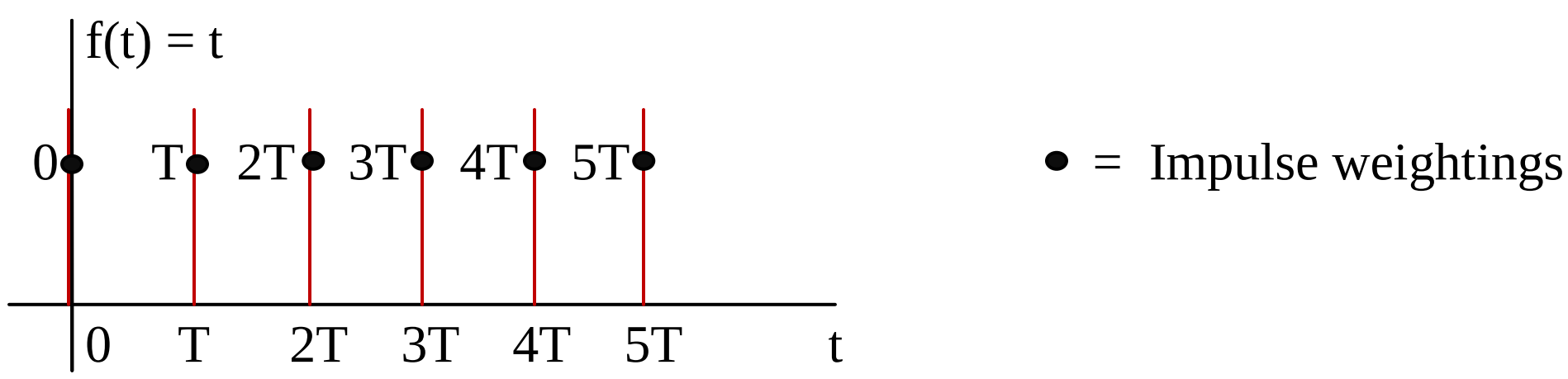

For the $K_{\Delta t}$ Transform

A series of weighted unit amplitude pulses of $\Delta t$ duration

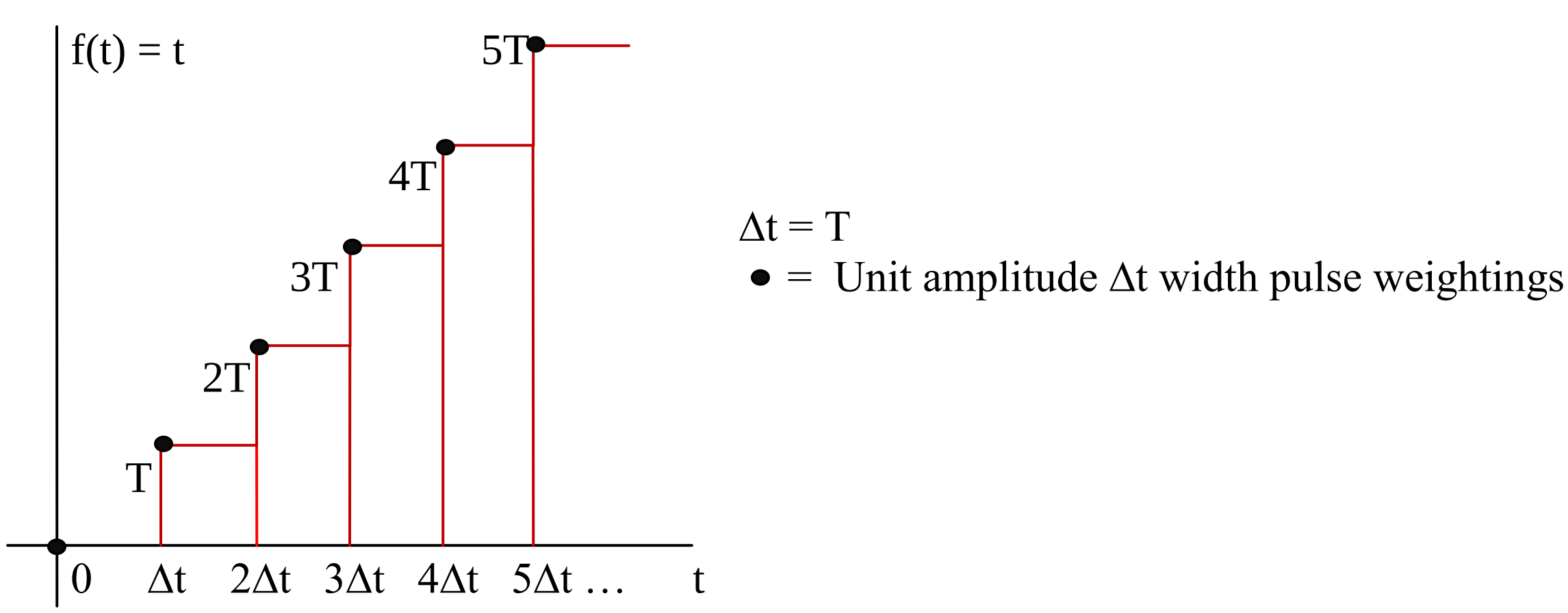

In a manner similar to that shown in Example 5.17-1 above, the $K_{\Delta t}$ Transform to Z Transform Conversion Equation can be applied to a Z Transform system block diagram composed of an input function, a transfer function, and an output function. This type of application is shown below.

Consider the following Z Transform system block diagram.

Fig #4 Z Transform System Block Diagram

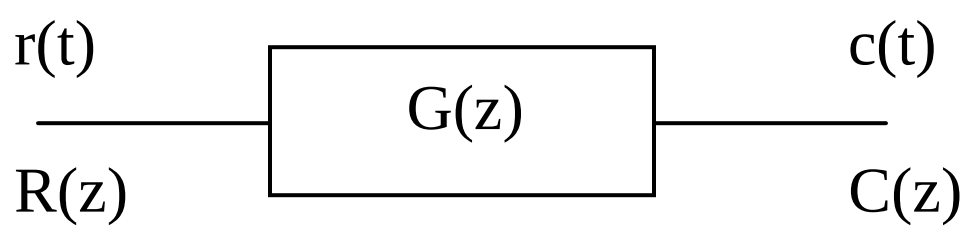

$t = 0, T, 2T, 3T \ldots$

$C(z) = G(z)R(z)$

Input function Transfer function Output function

From Fig #4

$$C(z) = G(z)R(z) \qquad (5.17\text{-}10)$$

where

R(z) = The Z Transform of the system input function, r(t)

C(z) = The Z Transform of the system output function, c(t)

G(z) = Z Transform system transfer function

Applying the Z Transform to $K_{\Delta t}$ Transform Conversion Equation, Eq 5.17-8, to Eq 5.17-10

$$C(s) = Tz^{-1}\,C(z)\Big|_{\substack{z\,=\,1+s\Delta t \\ \Delta t\,=\,T\ \ \text{sampling period} \\ t\,=\,n\Delta t \\ n\,=\,0,\,1,\,2,\,3,\,\ldots}} = Tz^{-1}\,G(z)R(z)\Big|_{\substack{z\,=\,1+s\Delta t \\ \Delta t\,=\,T\ \ \text{sampling period} \\ t\,=\,n\Delta t \\ n\,=\,0,\,1,\,2,\,3,\,\ldots}} \qquad (5.17\text{-}11)$$

Rearranging and simplifying Eq 5.17-11

$$C(s) = Tz^{-1}\,C(z)\Big|_{\substack{z\,=\,1+s\Delta t \\ \Delta t\,=\,T\ \ \text{sampling period}}} = G(z)\,\Big|_{\substack{z\,=\,1+s\Delta t \\ \Delta t\,=\,T\ \ \text{sampling period}}} \; Tz^{-1}R(z)\Big|_{\substack{z\,=\,1+s\Delta t \\ \Delta t\,=\,T\ \ \text{sampling period}}} \qquad (5.17\text{-}12)$$

where

$t = n\Delta t$

$n = 0, 1, 2, 3, \ldots$

$$C(s) = G(z)\,\Big|_{\substack{z\,=\,1+s\Delta t \\ \Delta t\,=\,T\ \ \text{sampling period}}} \; R(s) \qquad (5.17\text{-}13)$$

Let

$$G(s) = G(z)\,\Big|_{\substack{z\,=\,1+s\Delta t \\ \Delta t\,=\,T\ \ \text{sampling period}}} \qquad (5.17\text{-}14)$$

Substituting Eq 5.17-14 into Eq 5.17-13

$$C(s) = G(s)R(s) \qquad (5.17\text{-}15)$$

where

$R(s) = Tz^{-1}R(z)\Big|_{\substack{z = 1+s\Delta t \\ \Delta t = T \text{ sampling period}}}$ , $K_{\Delta t}$ Transform of the system input, r(t)

$C(s) = Tz^{-1}\, C(z)\Big|_{\substack{z = 1+s\Delta t \\ \Delta t = T \text{ sampling period}}}$ , $K_{\Delta t}$ Transform of the system output, c(t)

$G(s) = G(z)\,\Big|_{\substack{z = 1+s\Delta t \\ \Delta t = T \text{ sampling period}}}$ , $K_{\Delta t}$ Transform of the system transfer function

R(z) = The Z Transform of the system input, r(t)
C(z) = The Z Transform of the system output, c(t)
G(z) = Z Transform of the system transfer function
$t = n\Delta t$
$n = 0, 1, 2, 3, \ldots$

Eq 5.17-15 describes the following $K_{\Delta t}$ Transform system

Fig #5 $K_{\Delta t}$ Transform System Block Diagram

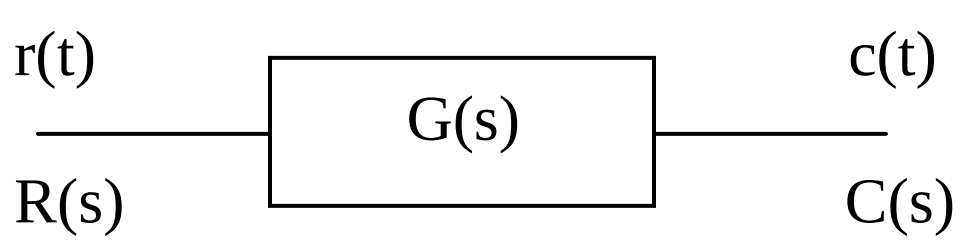

$t = 0, \Delta t, 2\Delta t, 3\Delta t, \ldots$

$C(s) = G(s)R(s)$

Input function Transfer function Output function

Comparing the system of Fig #5 (derived from the system of Fig #4) to the system of Fig#4, something interesting is observed. The two systems are the same. In both systems, an input of r(t) results in an output of c(t). The apparent difference is in the transform methodology selected to analyze the system. In Fig #4 the Z transform methodology uses weighted impulses to represent the discrete values of r(t) and c(t) at the discrete values of t, t = 0, T, 2T, 3T, … Whereas, in Fig #5 the $K_{\Delta t}$ Transform methodology uses weighted unit amplitude pulses of $\Delta t = T$ duration to represent the discrete values of r(t) and c(t) at the discrete values of t, $t = 0, \Delta t, 2\Delta t, 3\Delta t, \ldots$ An analysis using either of these two methodologies provides the same results.

Comment – If a $K_{\Delta t}$ Transform system is originally given, an equivalent Z Transform system could be derived from it using the $K_{\Delta t}$ Transform to Z Transform Conversion Equation in a derivation similar to the one shown above.

$$C(z) = G(z)R(z) \qquad (5.17\text{-}16)$$

where

$R(z) = \dfrac{1+s\Delta t}{\Delta t}\, R(s)\Big|_{\substack{s = \frac{z-1}{T} \\ T = \Delta t \text{ sampling period}}}$ , Z Transform system input

$$C(z) = = \frac{1+s\Delta t}{\Delta t}\,C(s)\Big|_{s = \frac{z-1}{T}} \quad \text{, Z Transform system output}$$

$T = \Delta t$ sampling period

$$G(z) = G(s)\Big|_{s = \frac{z-1}{T}} \quad \text{, Z Transform system transfer function}$$

$T = \Delta t$ sampling period

R(s) , The $K_{\Delta t}$ Transform of the system input
C(s) , The $K_{\Delta t}$ Transform of the system output
G(s) , $K_{\Delta t}$ Transform of the system transfer function
R(z) , The Z Transform of the system input
C(z) , The Z Transform of the system output
G(z) , Z Transform of the system transfer function
$t = n\Delta t$
$n = 0, 1, 2, 3, \ldots$

Considering that the Z Transform has a Modified Z Transform that can be used to calculate system output values between t = 0, T, 2T, 3T, …, one would think that a Modified $K_{\Delta t}$ Transform should exist for the same purpose. It does. The derivation of the Modified $K_{\Delta t}$ Transform is presented below.

The Modified $K_{\Delta t}$ Transform is the $K_{\Delta t}$ Transform equivalent of the Modified Z Transform. It is derived from the Modified Z Transform and is used in a similar manner. Depending on the discrete mathematics problem being analyzed, its Z Transform or $K_{\Delta t}$ Transform solution may be adequate or it may not. Both transform solutions, y(t), are accurate at the times, 0, Δt, 2Δt, 3Δt, … but nothing is known of the values of y(t) between the stated time increments. In many cases, where Δt is a sufficiently small value, knowing the value of y(t) at t = 0, Δt, 2Δt, 3Δt, … is adequate. Sufficiently accurate in-between values can be found using interpolation. If, however, Δt is not small or the function, y(t), is particularly variable, the values obtained for y(t) from a Z Transform or $K_{\Delta t}$ Transform analysis may be inadequate and even misleading (Note the system output waveform shown in Fig#3). The Modified Z Transform and Modified $K_{\Delta t}$ Transform have been developed to resolve this problem.

The Modified Z Transform is shown below.

<u>The Modified Z Transform of y(t) and related equations are:</u>

$$Z_M[y(t)] = Y(z,m) = z^{-1}\sum_{k=0}^{\infty} y(kT+mT)z^{-k} = \sum_{k=0}^{\infty} y(kT+[m-1]T)z^{-k} \text{ , Modified Z Transform of y(t)} \quad (5.17\text{-}17)$$

$$y(kT+[m-1]T) = Z^{-1}[\, Z_M[y(t)]\,] \quad (5.17\text{-}18)$$

where
$0 \le m < 1$
$Z_M[y(t)] = Y(z,m)$ , Modified Z Transform of the function, y(t), a function of z and m
Modified Z Transform discrete values of t = mT, T+mT, 2T+mT, 3T+mT, …

<u>Comment</u> - The Modified Z Transform Equation, Eq 5.17-17, is described and derived in the textbook, Digital and Sampled-Data Control Systems, by Julius T. Tou Copyrighted 1959 McGraw-Hill Book Company Inc., New York, Toronto, London

$$Z[y(t)] = Y(z) = \sum_{k=0}^{\infty} y(kT)z^{-1} \text{ , Z Transform of the function, y(t)} \quad (5.17\text{-}19)$$

From Eq 5.17-17 and Eq 5.17-19

**The definition of the Modified Z Transform is:**

$$\mathbf{Z_M[y(t)] = z^{-1}Z[y(kT+mT)] = z^{-1} \sum_{k=0}^{\infty} y(kT+mT)z^{-1} = Z[y(kT+[m-1]T)] = \sum_{k=0}^{\infty} y(kT+[m-1]T)z^{-k}} \quad (5.17\text{-}20)$$

**where**

**$0 \leq m < 1$**

**$Z_M[y(t)] = Y(z,m)$ = Modified Z Transform of the function, y(t), a function of z and m**

**$Z[y(t)]$ = Z Transform of y(t)**

**Modified Z Transform discrete values of t = mT, T+mT, 2T+mT, 3T+mT, …**

**T = Interval between successive discrete values of t**

**<u>Comment</u> - If the Modified Z Transform is expanded into a series, $\sum_{n=0}^{\infty} a_k z^{-k}$ ,**

**for $a_k z^{-k}$, $a_k = y(kT+[m-1]T)$ where k = 1, 2, 3, …**

Deriving the Modified $K_{\Delta t}$ Transform from the Z Transform

Using the Z Transform to $K_{\Delta x}$ Transform Conversion Equation, derived in Section 5.4, to relate the Modified Z Transform to a $K_{\Delta t}$ Transform

The Z Transform to $K_{\Delta t}$ Transform Conversion Equation is:

$$K_{\Delta t}[f(t)] = F(s) = Tz^{-1}F(z)\Big|_{\substack{z = 1+s\Delta t \\ \Delta t = T \\ t = 0, \Delta t, 2\Delta t, 3\Delta t, \dots}} \qquad K_{\Delta t}[f(t)] = F(s) \quad K_{\Delta t} \text{ Transform} \qquad Z[f(t)] = F(z) \quad \text{Z Transform} \quad (5.17\text{-}21)$$

From Eq 5.17-20 and Eq 5.17-21

$$TZ_M[y(t)]\Big|_{\substack{z = 1+s\Delta t \\ \Delta t = T}} = Tz^{-1}Z[y(t+mT)]\Big|_{\substack{z = 1+s\Delta t \\ \Delta t = T}} = K_{\Delta t}[y(t+m\Delta t)] = K_{\Delta t M}[y(t)] \quad (5.17\text{-}22)$$

$$K_{\Delta t M}[y(t)] = K_{\Delta t}[y(t+m\Delta t)] = TZ_M[y(t)]\Big|_{\substack{z = 1+s\Delta t \\ \Delta t = T}} \quad (5.17\text{-}23)$$

$$K_{\Delta t}[y(t)] = \sum_{k=0}^{\infty} y(k\Delta t)(1+s\Delta t)^{-k-1}\Delta t \text{ , Definition of the } K_{\Delta t} \text{ Transform of y(t)} \quad (5.17\text{-}24)$$

From Eq 5.17-24

$$K_{\Delta t}[y(t+m\Delta t)] = \sum_{k=0}^{\infty} y(k\Delta t+m\Delta t)(1+s\Delta t)^{-k-1}\Delta t$$ , Definition of the $K_{\Delta t}$ Transform of y(t+mΔt) (5.17-25)

Then

From Eq 5.17-23 and Eq 5.17-25

**The definition of the Modified $K_{\Delta t}$ Transform is:**

$$K_{\Delta t_M}[y(t)] = K_{\Delta t}[y(k\Delta t+m\Delta t)] = \sum_{k=0}^{\infty} y(k\Delta t+m\Delta t)(1+s\Delta t)^{-k-1}\Delta t = TZ_M[y(t)]\Big|_{\substack{z = 1+s\Delta t \\ \Delta t = T}} \qquad (5.17\text{-}26)$$

**where**

- **$0 \leq m < 1$**
- **$K_{\Delta t_M}[y(t)] = Y(s,m)$ , Modified $K_{\Delta t}$ Transform of the function, y(t), a function of s and m**
- **$K_{\Delta t}[y(t)]$ = $K_{\Delta t}$ Transform of y(t)**
- **$Z_M[y(t)] = Y(z,m)$ , Modified Z Transform of the function, y(t), a function of z and m**
- **Z[y(t)] = Z Transform of y(t)**
- **Modified $K_{\Delta t}$ Transform discrete values of t = mΔt, Δt+mΔt, 2Δt+mΔt, 3Δt+mΔt, …**
- **Δt = T = Interval between successive values of t**

**<u>Comment</u> - If the Modified $K_{\Delta t}$ Transform is expanded into a series, $\sum_{k=0}^{\infty} a_k(1+s\Delta t)^{-k-1}\Delta t$ ,**

**for $a_k(1+s\Delta t)^{-k-1}\Delta t$, $a_k = y(k\Delta t+m\Delta t)$ where k = 0, 1, 2, 3, …**

From Eq 5.17-26

**The Conversion from the Modified Z Transform to the Modified $K_{\Delta t}$ Transform is:**

$$K_{\Delta t_M}[y(t)] = TZ_M[y(t)]\Big|_{\substack{z = 1+s\Delta t \\ \Delta t = T}} \qquad (5.17\text{-}27)$$

**where**

- **$0 \leq m < 1$**
- **$K_{\Delta t_M}[y(t)] = Y(s,m)$ , Modified $K_{\Delta t}$ Transform of the function, y(t), a function of s and m**
- **$K_{\Delta t}[y(t)]$ = $K_{\Delta t}$ Transform of y(t)**
- **$Z_M[y(t)] = Y(z,m)$ , Modified Z Transform of the function, y(t), a function of z and m**
- **Z[y(t)] = Z Transform of y(t)**
- **Modified $K_{\Delta t}$ Transform discrete values of t = mΔt, Δt+mΔt, 2Δt+mΔt, 3Δt+mΔt, …**
- **Δt = T = Interval between successive values of t**

From Eq 5.17-27

**The Conversion from the Modified $K_{\Delta t}$ Transform to the Modified Z Transform is:**

$$Z_M[y(t)] = z^{-1}Z[y(t+mT)] = \frac{1}{T}K_{\Delta t_M}[y(t)]\Big|_{\substack{s = \frac{z-1}{T} \\ T = \Delta t}} \qquad (5.17\text{-}28)$$

**where**

**$0 \leq m < 1$**
**$Z_M[y(t)] = Y(z,m)$ , Modified Z Transform of the function, y(t), a function of z and m**
**$Z[y(t)]$ = Z Transform of y(t)**
**$K_{\Delta t_M}[y(t)] = Y(s,m)$ , Modified $K_{\Delta t}$ Transform of the function, y(t), a function of s and m**
**$K_{\Delta t}[y(t)] = K_{\Delta t}$ Transform of y(t)**
**Modified Z Transform discrete values of t = T, T+mT, T+mT, 3T+mT, …**
**$T = \Delta t$ = Interval between successive values of t**

Using The Conversion from the Modified Z Transform to the Modified $K_{\Delta t}$ Transform Equation, Eq 5.17-27, a table of Modified $K_{\Delta t}$ Transforms can be obtained. A table of Modified $K_{\Delta t}$ Transforms is provided on the following page.

# Modified $K_{\Delta t}$ Transforms

| # | Calculus Function Laplace Transforms | y(t) Calculus Functions Interval Calculus Functions $t = 0, \Delta t, 2\Delta t \ldots$ | $K_{\Delta t}$ Transform of y(t) Y(s) $t = 0, \Delta t, 2\Delta t, \ldots$ | Modified $K_{\Delta t}$ Transform of y(t) Y(s,m) $0 \le m < 1$, $t = m\Delta t, \Delta t + m\Delta t, 2\Delta t + m\Delta t, \ldots$ |
|---|---|---|---|---|
| **1** | $\frac{1}{s}$ | $1$ | $\frac{1}{s}$ | $\frac{1}{s}$ |
| **2** | $\frac{1}{s^2}$ | $t$ | $\frac{1}{s^2}$ | $\frac{1}{s^2} + \frac{m\Delta t}{s}$ |
| **3** | | $t(t-\Delta t)$ | $\frac{2}{s^3}$ | $\frac{2}{s^3} + \frac{2m\Delta t}{s^2} + \frac{m(m-1)\Delta t^2}{s}$ |
| **4** | $\frac{2}{s^3}$ | $t^2$ | $\frac{2}{s^3} + \frac{\Delta x}{s^2}$ | $\frac{2}{s^3} + \frac{(2m+1)\Delta t}{s^2} + \frac{m^2\Delta t^2}{s}$ |
| **5** | | $e_{\Delta t}(a,t)$ | $\frac{1}{s-a}$ | $\frac{(1+a\Delta t)^m}{s-a}$ |
| **6** | $\frac{1}{s-a}$ | $e^{at}$<br>or<br>$e_{\Delta t}(\frac{e^{a\Delta t}-1}{\Delta t},t)$ | $\frac{1}{s-\frac{e^{a\Delta t}-1}{\Delta t}}$ | $\frac{e^{am\Delta t}}{s-\frac{e^{a\Delta t}-1}{\Delta t}}$ |
| **7** | | $\sin_{\Delta t}(b,t)$ | $\frac{b}{s^2+b^2}$ | $\frac{\sin_{\Delta t}(b,m\Delta t)s + b\cos_{\Delta t}(b,m\Delta t)}{s^2+b^2}$ |
| **8** | $\frac{b}{s^2+b^2}$ | $\sin bt$<br>or<br>$e_{\Delta t}(\frac{\cos b\Delta t - 1}{\Delta t}, t)\sin_{\Delta t}(\frac{\tan b\Delta t}{\Delta t}, t)$ | $\frac{\frac{\sin b\Delta t}{\Delta t}}{(s + \frac{1-\cos b\Delta t}{\Delta t})^2 + (\frac{\sin b\Delta t}{\Delta t})^2}$ | $\frac{\cos bm\Delta t(\frac{\sin b\Delta t}{\Delta t}) + \sin bm\Delta t(s + \frac{1-\cos b\Delta t}{\Delta t})}{s^2 + 2(\frac{1-\cos b\Delta t}{\Delta t})s + 2(\frac{1-\cos b\Delta t}{\Delta t^2})}$ |
| **9** | | $\cos_{\Delta t}(b,t)$ | $\frac{s}{s^2+b^2}$ | $\frac{\cos_{\Delta t}(b,m\Delta t)s - b\sin_{\Delta t}(b,m\Delta t)}{s^2+b^2}$ |

| # | Calculus Function Laplace Transforms | y(t)<br>Calculus Functions<br>Interval Calculus Functions<br>t = 0, Δt, 2Δt… | $K_{\Delta t}$ Transform of y(t)<br>Y(s)<br>t = 0, Δt, 2Δt, … | Modified $K_{\Delta t}$ Transform of y(t)<br>Y(s,m)<br>$0 \le m < 1$, t = mΔt, Δt+mΔt, 2Δt+mΔt, … |
|---|---|---|---|---|
| **10** | $\frac{s}{s^2+b^2}$ | $\cos bt$<br>or<br>$e_{\Delta t}(\frac{\cos b\Delta t -1}{\Delta t}, x)\ \cos_{\Delta t}(\frac{\tan b\Delta t}{\Delta t}, t)$ | $\frac{s + \frac{1-\cos b\Delta t}{\Delta t}}{(s + \frac{1-\cos b\Delta t}{\Delta t})^2 + (\frac{\sin b\Delta t}{\Delta t})^2}$ | $\frac{\cos bm\Delta t(s + \frac{1-\cos b\Delta t}{\Delta t}) - \sin bm\Delta t(\frac{\sin b\Delta t}{\Delta t})}{s^2 + 2(\frac{1-\cos b\Delta t}{\Delta t})s + 2(\frac{1-\cos b\Delta t}{\Delta t^2})}$ |
| **11** | | $e_{\Delta t}(a,t)\sin_{\Delta t}(\frac{b}{1+a\Delta t},t)$<br>$a \neq -\frac{1}{\Delta t}$ | $\frac{b}{(s-a)^2+b^2}$ | $\frac{(1+a\Delta t)^m\,[b\cos_{\Delta t}(\frac{b}{1+a\Delta t},m\Delta t) + \sin_{\Delta t}(\frac{b}{1+a\Delta t},m\Delta t)(s-a)]}{(s-a)^2 + b^2}$ |
| **12** | $\frac{b}{(s-a)^2 + b^2}$ | $e^{at}\sin bt$<br>or<br>$e_{\Delta t}(\frac{e^{a\Delta t}\cos b\Delta t-1}{\Delta t}, t)\sin_{\Delta t}(\frac{\tan b\Delta t}{\Delta t}, t)$ | $\frac{\frac{e^{a\Delta t}\sin b\Delta t}{\Delta t}}{(s - \frac{e^{a\Delta t}\cos b\Delta t-1}{\Delta t})^2 + (\frac{e^{a\Delta t}\sin b\Delta t}{\Delta t})^2}$ | $\frac{e^{am\Delta t}[\cos m\Delta t(\frac{e^{a\Delta t}\sin b\Delta t}{\Delta t}) + \sin bm\Delta t(s + \frac{1-e^{a\Delta t}\cos b\Delta t}{\Delta t})]}{s^2 + 2(\frac{1-e^{a\Delta t}\cos b\Delta t}{\Delta t})s + (\frac{1-2e^{a\Delta t}\cos b\Delta t+e^{2a\Delta t}}{\Delta t^2})}$ |
| **13** | | $e_{\Delta t}(a,t)\cos_{\Delta t}(\frac{b}{1+a\Delta t},t)$<br>$a \neq -\frac{1}{\Delta t}$ | $\frac{s-a}{(s-a)^2+b^2}$ | $\frac{(1+a\Delta t)^m\,[\cos_{\Delta t}(\frac{b}{1+a\Delta t},m\Delta t)(s-a) - b\sin_{\Delta t}(\frac{b}{1+a\Delta t},m\Delta t)]}{(s-a)^2 + b^2}$ |
| **14** | $\frac{s-a}{(s-a)^2 + b^2}$ | $e^{at}\cos bt$<br>or<br>$e_{\Delta t}(\frac{e^{a\Delta t}\cos b\Delta t-1}{\Delta t}, t)\cos_{\Delta t}(\frac{\tan b\Delta t}{\Delta t}, t)$ | $\frac{s - \frac{e^{a\Delta t}\cos b\Delta t-1}{\Delta t}}{(s - \frac{e^{a\Delta t}\cos b\Delta t-1}{\Delta t})^2 + (\frac{e^{a\Delta t}\sin b\Delta t}{\Delta t})^2}$ | $\frac{e^{am\Delta t}[\cos bm\Delta t(s + \frac{1-e^{a\Delta t}\cos b\Delta t}{\Delta t}) - \sin bm\Delta t(\frac{e^{a\Delta t}\sin b\Delta t}{\Delta t})]}{s^2 + 2(\frac{1-e^{a\Delta t}\cos b\Delta t}{\Delta t})s + (\frac{1-2e^{a\Delta t}\cos b\Delta t+e^{2a\Delta t}}{\Delta t^2})}$ |

1. $K\Delta t_M[y(t)] = K_{\Delta t}[y(t+m\Delta t)]$ — Definition of the Modified $K_{\Delta t}$ Transform

2. $K\Delta t_M[y(t)] = Tz^{-1}Z[y(t+mT)] \big|_{z = 1+s\Delta t,\ \Delta t = T \text{ sampling period}}$ — Modified $K_{\Delta t}$ Transform from the Z Transform

   or

   $K\Delta t_M[y(t)] = TZ[y(t+[m-1]T)] \big|_{z = 1+s\Delta t,\ \Delta t = T \text{ sampling period}}$ — Modified $K_{\Delta t}$ Transform from the Z Transform

3. $K\Delta t_M[y(t)] = TZ_M[y(t)] \big|_{z = 1+s\Delta t,\ \Delta t = T \text{ sampling period}}$ — Modified $K_{\Delta t}$ Transform from the Modified Z Transform

4. $Z_M[y(t)] = z^{-1}Z[y(t+mT)]$ — Definition of the Modified Z Transform

   or

   $Z_M[y(t)] = Z[y(t+[m-1]T)]$ — Definition of the Modified Z Transform

5. $Z_M[y(t)] = \frac{1}{T} K\Delta t_M[y(t)] \big|_{s = \frac{z-1}{T},\ T = \Delta t \text{ sampling period}}$ — Modified Z Transform from the Modified $K_{\Delta t}$ Transform

where

$0 \leq m < 1$

$\Delta t = T$

Below are three examples, Example 5.17-1 thru Example 5.17-3, involving the use of the Modified $K_{\Delta t}$ Transform.

Example 5.17-1 Derivation of the Modified $K_{\Delta t}$ Transform of t, $K_{\Delta t_M}[t]$, from the Modified Z Transform of t, $Z_M[t]$.

Derive the Modified $K_{\Delta t}$ Transform of t, $K_{\Delta t_M}[t]$, from theModified Z Transform of t, $Z_M[t]$.

$$K_{\Delta t_M}[y(t)] = TZ_M[y(t)]\,\Big|_{\substack{z = 1+s\Delta t \\ \Delta t = T \text{ sampling period}}} \qquad \text{Modified } K_{\Delta t} \text{ Transform from Modified Z Transform} \quad 1)$$

$$y(t) = t \quad 2)$$

From the Modified Z Transform table in a CRC Handbook of Mathematical Tables

$$Z_M[t] = \frac{mT}{z-1} + \frac{T}{(z-1)^2}\,, \quad 0 \le m < 1 \quad 3)$$

From Eq 1 thru Eq 3

$$K_{\Delta t_M}[t] = TZ_M[t]\,\Big|_{\substack{z = 1+s\Delta t \\ \Delta t = T \text{ sampling period}}} = T\left[\frac{mT}{z-1} + \frac{T}{(z-1)^2}\right]\Big|_{\substack{z = 1+s\Delta t \\ \Delta t = T \text{ sampling period}}} \quad 4)$$

$$K_{\Delta t_M}[t] = \Delta t\left[\frac{m\Delta t}{s\Delta t} + \frac{\Delta t}{\Delta t^2 s^2}\right] = \frac{m\Delta t}{s} + \frac{1}{s^2} \quad 5)$$

Then

$$\mathbf{K_{\Delta t_M}[t] = \frac{1}{s^2} + \frac{m\Delta t}{s}\,, \quad 0 \le m < 1} \quad 6)$$

Checking Eq 6

$$K_{\Delta t_M}[y(t)] = K_{\Delta t}[y(t+m\Delta t]\ \ , \ \text{Definition of the Modified } K_{\Delta t} \text{ Transform} \quad 7)$$

$$K_{\Delta t_M}[t] = K_{\Delta t}[t+m\Delta t] = \frac{1}{s^2} + \frac{m\Delta t}{s}\,, \quad 0 \le m < 1 \quad 8)$$

Good check

Example 5.17-2 From the Modified $K\Delta t$ Transform, $K\Delta t_M[t] = \frac{1}{s^2} + \frac{m\Delta t}{s}$ where m= .5 and $\Delta t$ = .2, find the values of the function, y(t) = t, at t = .1, .3, .4, and .7

$K\Delta t_M[y(t)] = K_{\Delta t}[y(t+m\Delta t)] = Y(s,m)$ , Definition of the Modified $K_{\Delta t}$ Transform 1)

where

$t = 0, \Delta t, 2\Delta t, 3\Delta t, \ldots$

$Y(s,m)$ = Modified $K_{\Delta t}$ Transform of y(t), a function of s and m

Comment – The $K_{\Delta t}$ Transform of $y(t+m\Delta t)$, $K_{\Delta t}[y(t+m\Delta t)]$, is the Modified $K_{\Delta t}$ Transform of y(t), $K\Delta t_M[y(t)]$. Its inverse $K_{\Delta t}$ Transform yields $y(t+m\Delta t)$ at $t = 0, \Delta t, 2\Delta t, 3\Delta t, \ldots$

$y(t) = t$ 2)

$\Delta t = .2$ 3)

$m = .5$ 4)

$K\Delta t_M[t] = \frac{1}{s^2} + \frac{m\Delta t}{s}$ 5)

Substituting Eq 3 and Eq 4 into Eq 5

$K_{.2_M}[t] = \frac{1}{s^2} + \frac{.1}{s}$ 6)

Changing the form of $K\Delta t_M[t]$, Eq 5, to that of an Asymptotic Maclaurin Series

$K_{.2_M}[t] = \frac{1}{s^2} + \frac{.1}{s} = .1s^{-1} + 1s^{-2} + 0s^{-3} + 0s^{-4} + 0s^{-5} + \ldots$ 7)

$K\Delta t_M[y(t)] = K_{\Delta t}[y(t+m\Delta t)]$ , $t = 0, \Delta t, 2\Delta t, 3\Delta t, \ldots$ 8)

Finding the values of y(t+.1) where t = 0, .2, .4, .6, … using the Inverse $K_{\Delta t}$ Transform Formula

Comment – The Inverse $K_{\Delta t}$ Transform Formula was derived in Section 5.6

The Inverse $K_{\Delta t}$ Transform Formula is:

$K_{\Delta t}^{-1}[K_{\Delta t}[f(t)]] = f(n\Delta t) = \sum_{p=0}^{n} [{}_nC_p \Delta t^p]\, a_{p+1}$ 9)

where

$K_{\Delta t}[f(t)] = F(s) = a_1 s^{-1} + a_2 s^{-2} + a_3 s^{-3} + a_4 s^{-4} + \ldots$ , Asymptotic Series form of F(s)

$F(s) = K_{\Delta t}[f(t)]$ = $K_{\Delta t}$ Transform of the function, f(t)

$K_{\Delta t}^{-1}[F(s)] = f(t)$ = Inverse $K_{\Delta t}$ Transform of the function, F(s)

$n = 0, 1, 2, 3, 4, \ldots$

$p = 0, 1, 2, 3, \ldots, n$

$$ {}_nC_p = \frac{n!}{r!(n-r)!} $$

$0! = 1$

$f(t)$ = function of t

$\Delta t$ = interval between successive values of t

$t = 0, \Delta t, 2\Delta t, 3\Delta t, \ldots$

$D_{\Delta t}{}^{m}f(0) = a_{m+1}$ , $m = 0,1,2,3,\ldots$ , Initial conditions of f(t)

$$ f(t) = y(t+m\Delta t) \qquad 10) $$

where

$t = n\Delta t$

$n = 0, 1, 2, 3, \ldots$

From Eq 3, Eq 4, Eq 9, and Eq 10,

$$ y(.2n+.1) = \sum_{p=0}^{n} [{}_nC_p(.2)^p]\, a_{p+1} \ , \quad n = 0, 1, 2, 3, \ldots \qquad 11) $$

From Eq 7 and Eq 9

$$ a_1 = .1 \qquad 12) $$
$$ a_2 = 1 \qquad 13) $$
$$ a_n = 0 \ , \ n = 3, 4, 5, \ldots \qquad 14) $$

From Eq 11 thru Eq 14

For n = 0

$$ y(.1) = (1)a_1 = (1)(.1) = .1 \qquad 15) $$

For n = 1

$$ y(.3) = (1)a_1 + (1)(.2)a_2 = (1)(.1) + (1)(.2)(1) = .3 \qquad 16) $$

For n = 2

$$ y(.5) = (1)a_1 + (2)(.2)a_2 + (1)(.2)^2a_3 = (1)(.1) + (2)(.2)(1) + (1)(.2)^2(0) = .5 \qquad 17) $$

For n = 3

$$ y(.7) = (1)a_1 + (3)(.2)a_2 + (3)(.2)^2a_3 + (1)(.2)^3a_4 = (1)(.1) + (3)(.2)(1) + (3)(.2)^2(0) + (1)(.2)^3(0) = .7 \qquad 18) $$

.

.

.

$$nC_1 = \frac{n!}{1!(n-1)!} = n \quad \text{for } n = 1,2,3,\ldots \qquad 19)$$

In general

$$y(.1) = .1 \qquad 20)$$

$$y(.2n+.1) = (1)a_1 + n(.2)a_2 + 0 = (1)(.1) + n(.2)(1) = .1 + .2n\ , \quad \text{for } n = 1,2,3,\ldots \qquad 21)$$

From Eq 20 and Eq 21

$$y(.2n+.1) = .2n + .1\ , \quad \text{for } n = 0,1,2,3,\ldots \qquad 21)$$

or

**$y(t) = t\ ,\ \ t = .1, .3, .5, .7, \ldots$**

Example 5.17-3 Find c(t) , the output of the following Z Transform system block diagram, using $K_{\Delta t}$ Transforms. Find the values of c(t) at t = 0, .05, .1, .15 and .2 in response to a unit step input. The initial condition for c(t) is c(0) = 0, a = 5 and T = .1.

Fig #1 Z Transform System Block Diagram

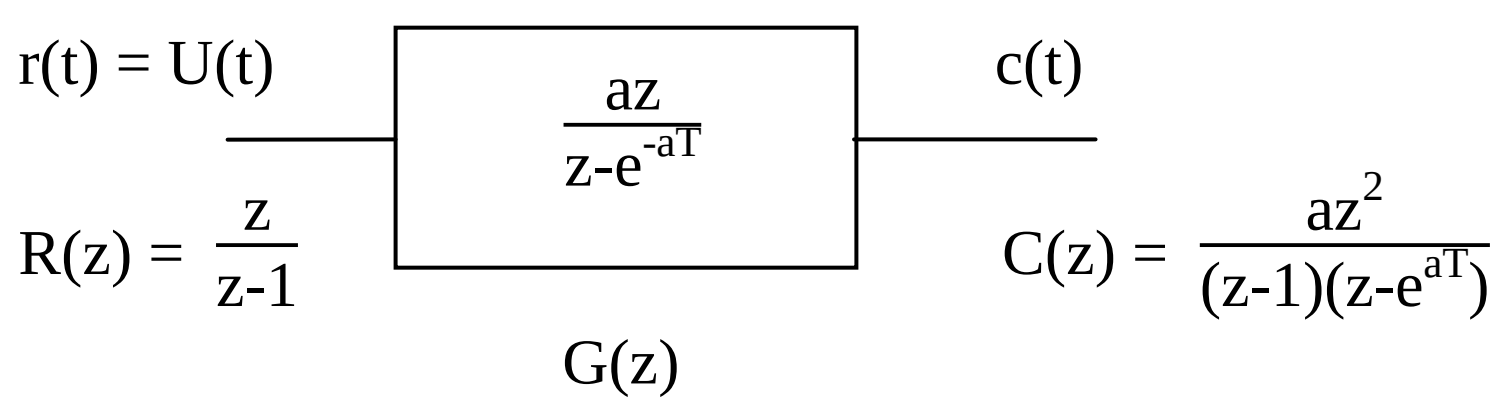

r(t) = | | | | … , a series of input unit value impulses
0 T 2T 3T

Convert the above Z Transform system block diagram into its equivalent $K_{\Delta t}$ Transform system block diagram.

From the Z Transform system block diagram

$$C(z) = \frac{z}{z-1}\,\frac{az}{z-e^{-aT}} = \frac{az^2}{(z-1)(z-e^{-aT})} \qquad 1)$$

Writing the Z Transform to $K_{\Delta t}$ Transform Conversion Equation

$$K_{\Delta t}[f(t)] = F(s) = Tz^{-1}F(z)\Big|_{\substack{z = 1+s\Delta t \\ \Delta t = T \\ t = 0,\ \Delta t,\ 2\Delta t,\ 3\Delta t,\ \ldots}} \qquad 2)$$

$K_{\Delta t}[f(t)] = F(s)$ $K_{\Delta t}$ Transform

$Z[f(t)] = F(z)$ Z Transform

Applying the Z Transform to $K_{\Delta t}$ Transform Conversion Equation to Eq 1

$$Tz^{-1}C(z)\Big|_{\substack{z = 1+s\Delta t \\ \Delta t = T}} = Tz^{-1}\frac{az^2}{(z-1)(z-e^{-aT})}\Big|_{\substack{z = 1+s\Delta t \\ \Delta t = T}} = [Tz^{-1}\frac{z}{z-1}]\Big|_{\substack{z = 1+s\Delta t \\ \Delta t = T}} [\frac{az}{(z-e^{-aT})}]\Big|_{\substack{z = 1+s\Delta t \\ \Delta t = T}} \qquad 3)$$

Simplifying

$$C(s) = \frac{1}{s}[\frac{a(1+s\Delta t)}{1+s\Delta t-e^{-at}}] = \frac{a(1+s\Delta t)}{\Delta ts(s-\frac{e^{-a\Delta t}-1}{\Delta t})} \qquad 4)$$

$a = 5$ 5)

$\Delta t = T = .1$ 6)

Substituting Eq 5 and Eq 6 into Eq 4

$$C(s) = \frac{5(1+.1s)}{.1s(s-\frac{e^{-.5}-1}{.1})} = \frac{50(1+.1s)}{s(s+.393469340)} \qquad 7)$$

$$C(s) = \frac{50(1+.1s)}{s(s+.393469340)} \quad , \quad K_{\Delta t} \text{ Transform of } c(t) \qquad 8)$$

Fig #2 Equivalent $K_{\Delta t}$ Transform System Block Diagram

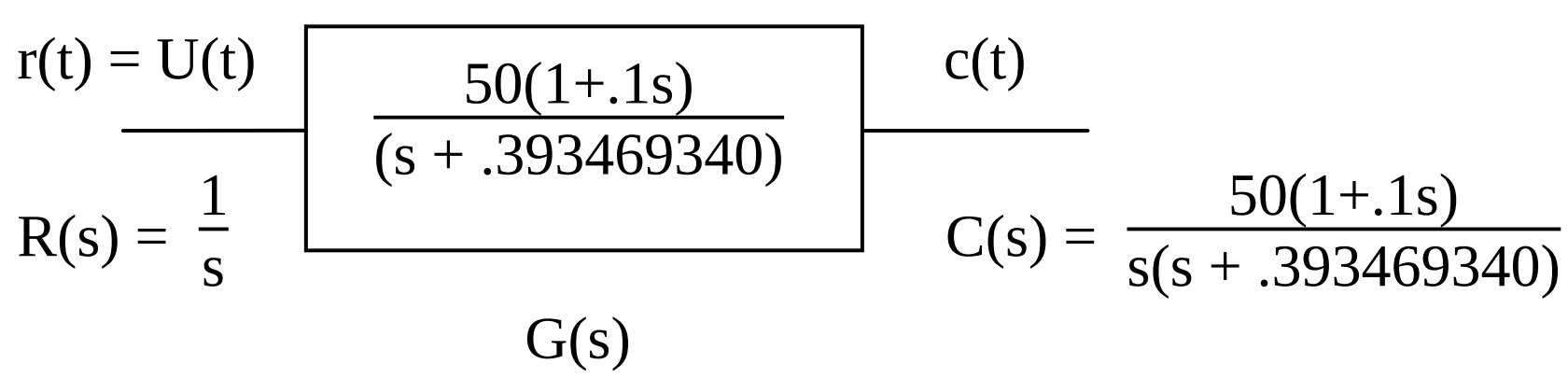

r(t) = ⎸‾‾⎸‾‾⎸‾‾⎸ … , a series of input unit amplitude pulses of $\Delta t$ duration that form a unit step initiated at t = 0

0 $\Delta t$ $2\Delta t$ $3\Delta t$

Convert Fig #2 into an Equivalent System Diagram using Modified $K_{\Delta t}$ Transforms.

From the Modified $K_{\Delta t}$ Transform Table shown in this section or in the Appendix

The $K_{\Delta t}$ Transform, $\dfrac{1}{s-\frac{e^{-a\Delta t}-1}{\Delta t}}$, becomes the Modified $K_{\Delta t}$ Transform, $\dfrac{e^{-am\Delta t}}{s-\frac{e^{-a\Delta t}-1}{\Delta t}}$,

where $0 \le m < 1$. 9)

$a\Delta t = 5(.1) = .5$

From Fig #2 and Eq 9

Fig #3    Equivalent Modified $K_{\Delta t}$ Transform System Block Diagram

$r(t) = U(t)$ → $G(s) = \dfrac{50(1+.1s)e^{-.5m}}{(s + .393469340)}$ → $c(t+.1m)$

$R(s) = \dfrac{1}{s}$    $C(s,m) = \dfrac{50(1+.1s)e^{-.5m}}{s(s + .393469340)}$

$r(t) =$ unit step initiated at $t = 0$    $0 \le m < 1$

$$\frac{1}{s(s + .393469340)} = \frac{A}{s} + \frac{B}{s+3.93469340} \qquad 10)$$

$$A = \frac{1}{s+3.93469340}\Big|_{s = 0} = \frac{1}{3.93469340} \qquad 11)$$

$$B = \frac{1}{s}\Big|_{s = -3.93469340} = -\frac{1}{3.93469340} \qquad 12)$$

Substituting Eq 11 and Eq 12 into Eq 10

$$\frac{1}{s(s + .393469340)} = \frac{1}{3.93469340}\left[\frac{1}{s} - \frac{1}{s+3.93469340}\right] \qquad 13)$$

$$C(s,m) = \frac{50(1+.1s)e^{-.5m}}{s(s + .393469340)} = \frac{50(1+.1s)e^{-.5m}}{3.93469340}\left[\frac{1}{s} - \frac{1}{s+3.93469340}\right] \qquad 14)$$

$$C(s,m) = \frac{50(1+.1s)e^{-.5m}}{3.93469340}\left[\frac{1}{s} - \frac{1}{s+3.93469340}\right] \qquad 15)$$

Finding the Inverse $K_{\Delta t}$ Transform of Eq 15

$$c(t+.1m) = 12.7074704e^{-.5m}\left[1 - \{1 - 3.93469340(.1)\}^{\frac{t+.1}{.1}}\right] \qquad 16)$$

Then

$$\mathbf{c(t+.1m) = 12.7074704e^{-.5m}\left[1 - (.60653066)^{10t+1}\right], \quad t = 0, .1, .2, .3, \ldots, \quad 0 \le m < 1} \qquad 17)$$

Calculating c(t) for t = 0, .05, .1 using Eq 17

For t = 0 and m = 0

$c(0) = 12.7074704(1 - .60653066) = 5$ 18)

**c(0) = 5** 19)

For t = 0 and m = .5

$c(.05) = 12.7074704e^{-.25}(1 - .60653066) = 3.89400391$ 20)

**c(.05) = 3.89400391** 21)

For t = 0 and m →1

$c(.1) = 12.7074704e^{-.5}(1 - .60653066) = 3.03265329$ 22)

**c(.1) = 3.03265329** 23)

Calculating c(t) for t = .1, .15, .2 using Eq 17

For t = .1 and m = 0

$c(.1) = 12.7074704(1 - .60653066^{2}) = 8.03265329$ 24)

**c(.1) = 8.03265329** 25)

For t = .1 and m = .5

$c(.15) = 12.7074704e^{-.25}(1 - .60653066^{2}) = 6.255836679$ 26)

**c(.15) = 6.25583667** 27)

For t = .1 and m → 1

$c(.2) = 12.7074704e^{-.5}(1 - .60653066^{2}) = 4.87205050$ 28)

**c(.2) = 4.87205050** 29)

Fig #4 System Output

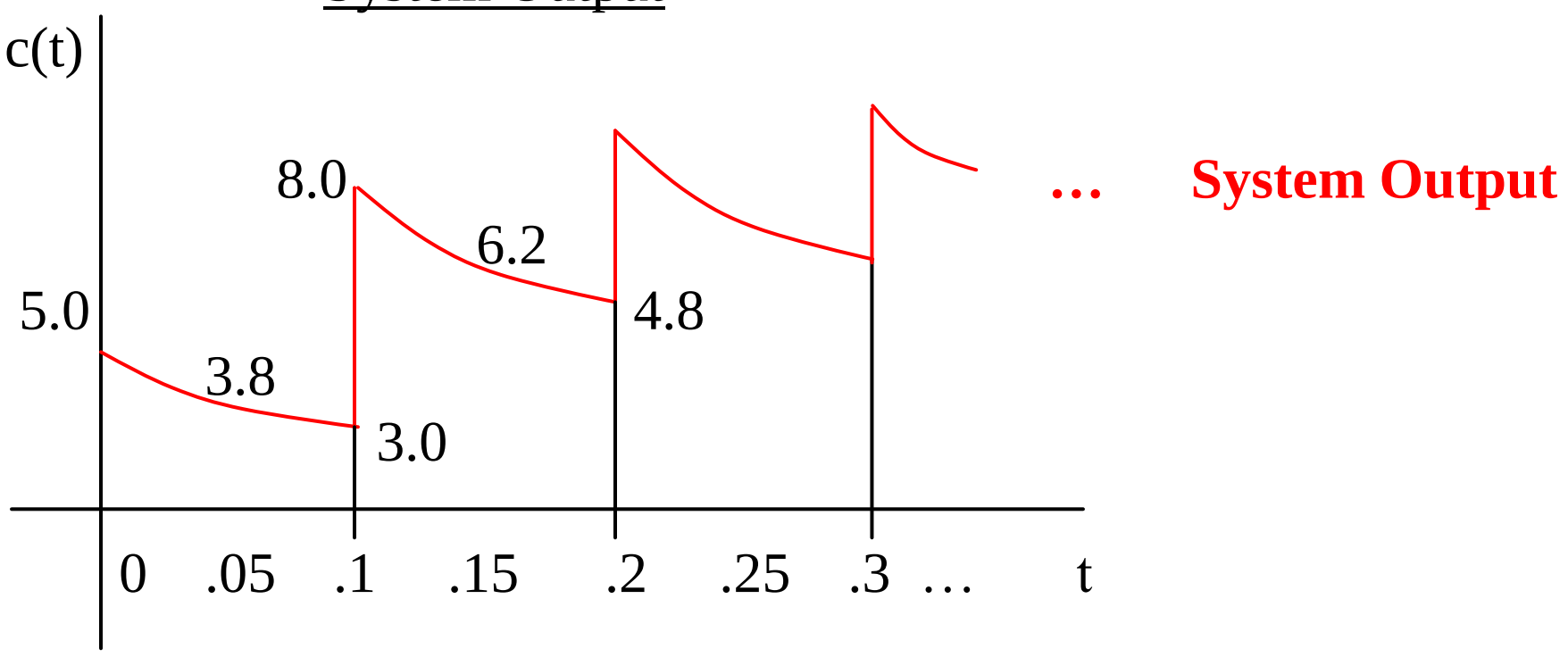

Check Eq 17 by finding c(t) for t = 0, .05, .1, .15, and .2 using a Modified Z Transform analysis of the system block diagram of Fig #1. The Modified $K_{\Delta t}$ Transform system analysis just completed and the following Modified Z Transform system analysis should provide the same results. T = .1 and a = 5

Rewriting the Z Transform System Block Diagram of Fig #1

Z Transform System Block Diagram of Fig#1

r(t) = U(t) — [ $G(z) = \frac{az}{z-e^{-aT}}$ ] — c(t)

$R(z) = \frac{z}{z-1}$  $C(z) = \frac{az^2}{(z-1)(z-e^{-aT})}$

G(z)

r(t) = | | | | … , a series of input unit value impulses
0 T 2T 3T

From the Modified Z Transform Table shown in the textbook, "Digital and Sampled-Data Control Systems" by Julius T. Tou (See the Reference Book Listing at the end of the Appendix)

The Z Transform, $\frac{z}{z-e^{-aT}}$, becomes the Modified Z Transform, $\frac{e^{-amT}}{z-e^{-aT}}$, where $0 \le m < 1$. 30)

Converting the Z Transform System Block Diagram of Fig#1 to a Modified Z Transform System Block Diagram shown below in Fig# 5

Fig #5 Modified Z Transform System Block Diagram

r(t) = U(t) — [ $G(z) = \frac{ae^{-amT}}{z-e^{-aT}}$ ] — c(t)

$R(z) = \frac{z}{z-1}$  $C(z,m) = \frac{aze^{-amT}}{(z-1)(z-e^{-aT})}$

G(z)  $0 \le m < 1$

r(t) = | | | | … , a series of input unit value impulses
0 T 2T 3T

$$C(z,m) = \frac{aze^{-amT}}{(z-1)(z-e^{-aT})} \quad 31)$$

$a = 5$ 32)

$T = \Delta t = .1$ 33)

$aT = .5$ 34)

Substituting Eq 32 thru Eq 34 into Eq 31 and expanding

$$C(z,m) = \frac{5ze^{-.5m}}{(z-1)(z-e^{-.5})} = 5e^{-.5m}\left[\frac{z}{(z-1)(z-e^{-.5})}\right] = 5e^{-.5m}\left[\frac{A}{z-1} + \frac{B}{z-e^{-.5}}\right] \quad 35)$$

$$A = \frac{z}{z-e^{-aT}}\Big|_{z=1} = \frac{1}{1-e^{-aT}} \quad 36)$$

$$B = \frac{z}{z-1}\Big|_{z=e^{-.5}} = \frac{e^{-.5}}{e^{-.5}-1} = -\frac{e^{-.5}}{1-e^{-.5}} \quad 37)$$

Substituting Eq 36 and Eq 37 into Eq 3

$$C(z,m) = \frac{5e^{-.5m}}{1-e^{-.5}}\left[\frac{1}{z-1} - \frac{e^{-.5}}{z-e^{-.5}}\right] \quad 36)$$

or

$$C(z,m) = \frac{5e^{-.5m}z^{-1}}{1-e^{-.5}}\left[\frac{z}{z-1} - \frac{e^{-.5}z}{z-e^{-.5}}\right] \;,\quad T = .1 \;,\quad \text{The Modified Z Transform of C(z)} \quad 37)$$

Finding the Inverse Z Transform of Eq 37

$$c(t,m) = c(t+[m-1].1) = \frac{5e^{-.5m}}{1-e^{-.5}}\left[1 - e^{-.5}e^{-5(t-.1)}\right] \;,\quad t = 0, .1, .2, .3, \ldots \quad 38)$$

$$c(t+[m-1].1) = \frac{5e^{-.5m}}{1-e^{-.5}}\left[1 - e^{-5t}\right] \;,\quad t = 0, .1, .2, .3, \ldots \quad 39)$$

Then

$$\mathbf{c(t+[m-1].1) = 12.7074704e^{-.5m}\left[1 - e^{-5t}\right] \;,\quad t = 0\ .1, .2, .3, \ldots \;,\quad 0 \le m < 1} \quad 40)$$

Calculating c(t) for t = 0, .05, .1 using Eq 38

For t = .1 and m = 0

$$c(0) = 12.7074704(1 - .60653066) = 5 \quad 41)$$

$$\mathbf{c(0) = 5} \quad 42)$$

For t = .1 and m = .5

$$c(.05) = 12.7074704e^{-.25}(1 - .60653066) = 3.89400391 \quad 43)$$

$$\mathbf{c(.05) = 3.89400391} \quad 44)$$

For t = .1 and m →1

$c(.1) = 12.7074704e^{-.5}(1 - .60653066) = 3.03265329$ 45)

**c(.1) = 3.03265329** 46)

Calculating c(t) for t = .1, .15, .2 using Eq 17

For t = .2 and m = 0

$c(.1) = 12.7074704(1 - .367879441) = 8.03265329$ 47)

**c(.1) = 8.03265329** 48)

For t = .2 and m = .5

$c(.15) = 12.7074704e^{-.25}(1 - .367879441) = 6.25583667$ 49)

**c(.15) = 6.25583667** 50)

For t = .2 and m → 1

$c(.2) = 12.7074704e^{-.5}(1 - .367879441) = 4.87205050$ 51)

**c(.2) = 4.87205050**

Note that Eq 17 and Eq 40 are equivalent though in a different form. 52)

Good check

# CHAPTER 6

## Miscellaneous Derivations

**Section 6.1: Derivation of the Interval Calculus discrete integration by parts formula**

Consider the following two rectangles of equal size.

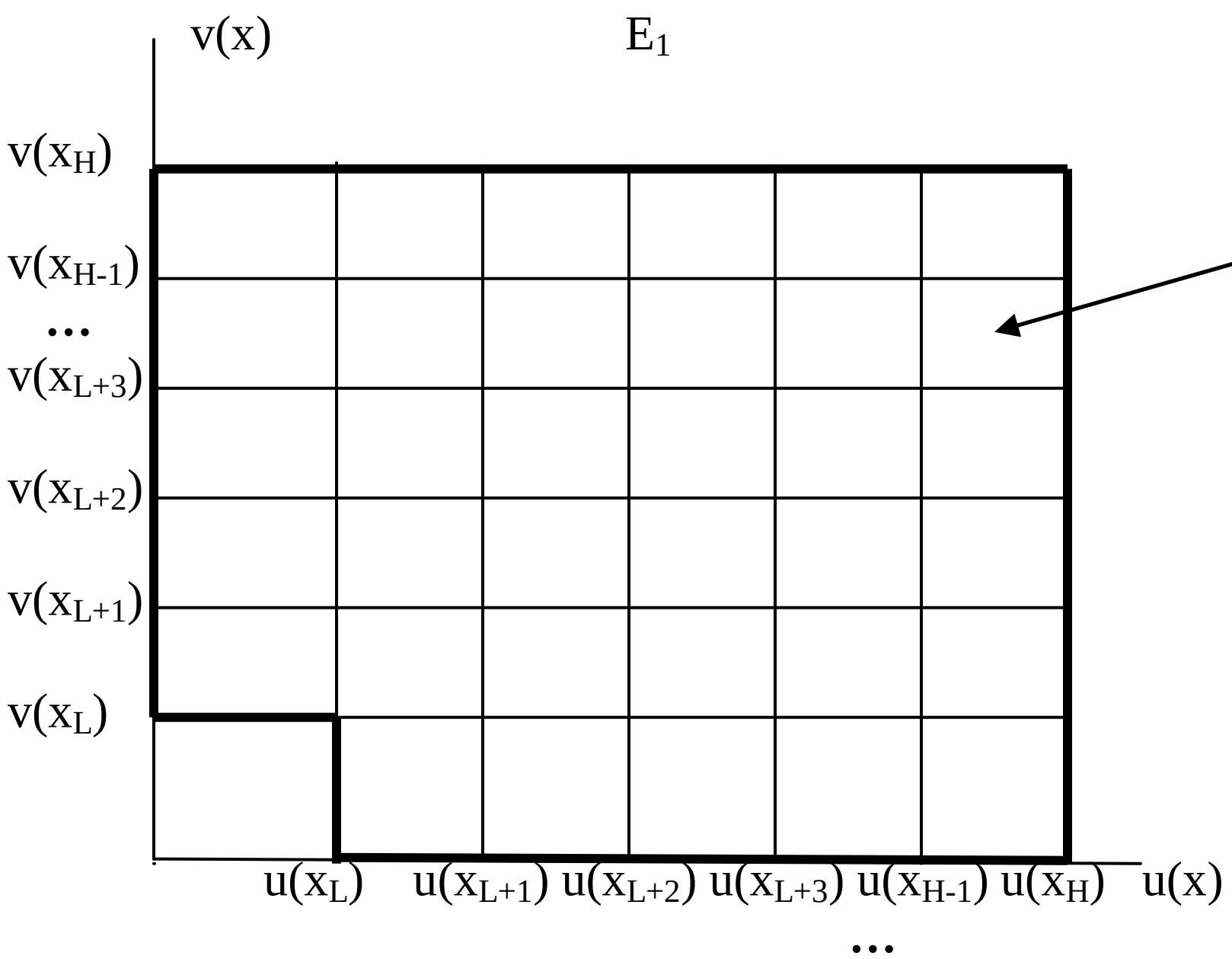

$$E_1 = {}_1\sum_{n=L}^{H}[u(x_{n+1})v(x_{n+1}) - u(x_n)v(x_n)]$$

$E_1$ = Area within the dark black lines

$$E_2 = {}_1\sum_{n=L}^{H-1}[v(x_{n+1}) - v(x_n)]\, u(x_{n+1}) + {}_1\sum_{n=L}^{H-1}[u(x_{n+1}) - u(x_n)]\, v(x_n)$$

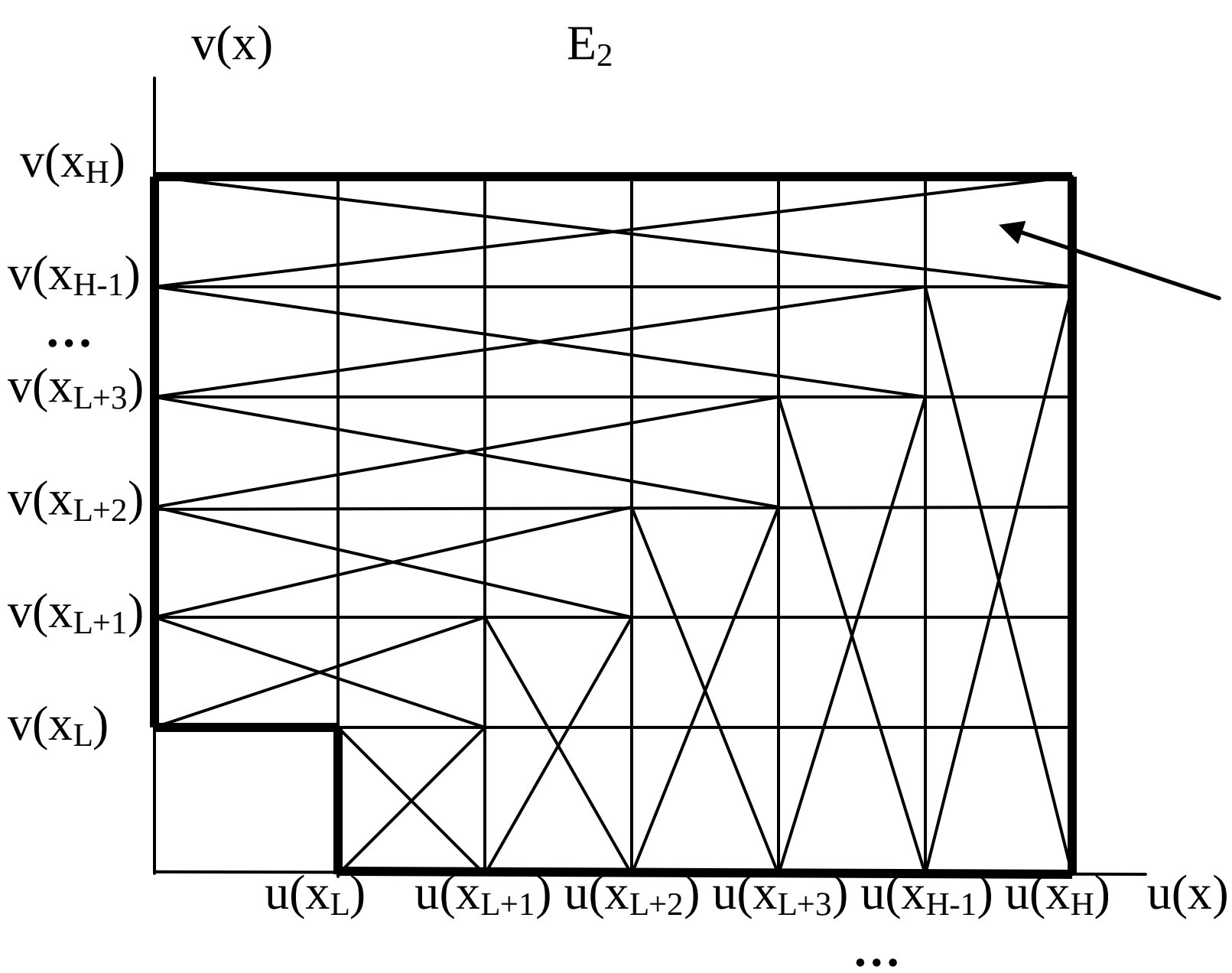

$E_2$ = Area within the dark black lines

$$E_1 = {}_1\sum_{n=L}^{H}[u(x_{n+1})v(x_{n+1}) - u(x_n)v(x_n)] \qquad (6.1\text{-}1)$$

Canceling terms

$$E_1 = {}_1\sum_{n=L}^{H}[u(x_{n+1})v(x_{n+1}) - u(x_n)v(x_n)] = u(x_n)v(x_n)\Big|_{L}^{H} = u(x_H)v(x_H) - u(x_L)v(x_L)$$

$$E_1 = u(x_n)v(x_n)\Big|_{L}^{H} \qquad (6.1\text{-}2)$$

$$E_2 = [v(x_{L+1})-v(x_L)]u(x_{L+1}) + [u(x_{L+1})-u(x_L)]v(x_L) + [v(x_{L+2})-v(x_{L+1})]u(x_{L+2}) + [u(x_{L+2})-u(x_{L+1})]v(x_{L+1}) + [v(x_{L+3})-v(x_{L+2})]u(x_{L+3}) + [u(x_{L+3})-u(x_{L+2})]v(x_{L+2}) + \ldots + [v(x_{H-1})-v(x_{H-2})]u(x_{H-1}) + [u(x_{H-1})-u(x_{H-2})]v(x_{H-2}) + [v(x_H)-v(x_{H-1})]u(x_H) + [u(x_H)-u(x_{H-1})]v(x_{H-1})$$

$$E_2 = {}_1\sum_{n=L}^{H-1}[v(x_{n+1}) - v(x_n)]\, u(x_{n+1}) + {}_1\sum_{n=L}^{H-1}[u(x_{n+1}) - u(x_n)]\, v(x_n) \qquad (6.1\text{-}3)$$

Let

$$x_{n+1} - x_n = \Delta x \ , \text{ an x increment} \qquad (6.1\text{-}4)$$

$$x_{n+1} = x_n + \Delta x \qquad (6.1\text{-}5)$$

From Eq 6.1-3 and Eq 6.1-5

$$E_2 = {}_1\sum_{n=L}^{H-1}\left[\frac{v(x_n+\Delta x) - v(x_n)}{\Delta x}\right] u(x_n+\Delta x)\, \Delta x + {}_1\sum_{n=L}^{H-1}\left[\frac{u(x_n+\Delta x) - u(x_n)}{\Delta x}\right] v(x_n)\, \Delta x \qquad (6.1\text{-}6)$$

$$D_{\Delta x}f(x_n) = \frac{f(x_n+\Delta x)-f(x_n)}{\Delta x}\ , \quad f(x_n) = \text{a function of x evaluated at } x_n\ ,\ \ n = 1,2,3,\ldots \qquad (6.1\text{-}7)$$

From Eq 6.1-6 and Eq 6.1-7

$$E_2 = {}_1\sum_{n=L}^{H-1}D_{\Delta x}v(x_n)\, u(x_n+\Delta x)\, \Delta x + {}_1\sum_{n=L}^{H-1}D_{\Delta x}u(x_n)\, v(x_n)\, \Delta x \qquad (6.1\text{-}8)$$

Eq 6.1-8 can be written in the discrete integral form where $x = x_L, x_{L+1}, x_{L+2}, x_{L+3}, \dots, x_{H-1}, x_H$

$$E_2 = {}_{\Delta x}\int_{x_L}^{x_H} D_{\Delta x}v(x)\, u(x+\Delta x)\, \Delta x + {}_{\Delta x}\int_{x_L}^{x_H} D_{\Delta x}u(x)\, v(x)\, \Delta x \qquad (6.1\text{-}9)$$

Note from the diagrams of $E_1$ and $E_2$ above that $E_1 = E_2$.

$$E_1 = E_2 \qquad (6.1\text{-}10)$$

From Eq 6.1-2 where $x = x_n$

$$E_1 = u(x_n)v(x_n)\Big|_L^H = u(x)v(x)\Big|_{x_L}^{x_H} \qquad (6.1\text{-}11)$$

Equating Eq 6.1-11 and Eq 6.1-9

$$u(x)v(x)\Big|_{x_L}^{x_H} = {}_{\Delta x}\int_{x_L}^{x_H} D_{\Delta x}v(x)\, u(x+\Delta x)\, \Delta x + {}_{\Delta x}\int_{x_L}^{x_H} D_{\Delta x}u(x)\, v(x)\, \Delta x \qquad (6.1\text{-}12)$$

Thus
Rearranging the terms of Eq 6.1-12

**The Interval Calculus discrete integration by parts formula is:**

$$\mathbf{{}_{\Delta x}\int_{x_L}^{x_H} D_{\Delta x}u(x)\, v(x)\, \Delta x = u(x)v(x)\Big|_{x_L}^{x_H} - {}_{\Delta x}\int_{x_L}^{x_H} D_{\Delta x}v(x)\, u(x+\Delta x)\, \Delta x} \qquad (6.1\text{-}13)$$

## Section 6.2: Derivation of the discrete derivative of the product of two functions

From the Interval Calculus discrete integration by parts formula derived in Section 6.1

$$_{\Delta x}\int_{x_L}^{x_H} D_{\Delta x}u(x)\, v(x)\, \Delta x = u(x)v(x)\Big|_{x_L}^{x_H} - {}_{\Delta x}\int_{x_L}^{x_H} D_{\Delta x}v(x)\, u(x+\Delta x)\, \Delta x \qquad (6.2\text{-}1)$$

Changing the form Eq 6.2-1 to its indefinite form

$$_{\Delta x}\int D_{\Delta x}u(x)\, v(x)\, \Delta x = u(x)v(x) - {}_{\Delta x}\int D_{\Delta x}v(x)\, u(x+\Delta x)\, \Delta x + k \qquad (6.2\text{-}2)$$

where

k = constant of integration

Rearranging Eq 6.2-2 and differentiating

$$D_{\Delta x}[u(x)v(x)] = D_{\Delta x}\, {}_{\Delta x}\int D_{\Delta x}u(x)\, v(x)\, \Delta x + D_{\Delta x}\, {}_{\Delta x}\int D_{\Delta x}v(x)\, u(x+\Delta x)\, \Delta x + D_{\Delta x}k \qquad (6.2\text{-}3)$$

Thus
From Eq 6.2-3

**The discrete derivative of a product of two functions is:**

$$\mathbf{D_{\Delta x}[u(x)v(x)] = D_{\Delta x}u(x)\, v(x) + D_{\Delta x}v(x)\, u(x+\Delta x)} \qquad (6.2\text{-}4)$$

## Section 6.3: Derivation of the discrete derivative of the division of two functions

$$D_{\Delta x}f(x) = \frac{f(x+\Delta x) - f(x)}{\Delta x} \quad , \quad \text{discrete derivative of } f(x) \qquad (6.3\text{-}1)$$

where

$f(x)$ = function of x

$\Delta x$ = x increment

$x = m\Delta x$ , m = integers

Take the derivative of the division of two functions, $\frac{u(x)}{v(x)}$ .

$$D_{\Delta x}[\, \frac{u(x)}{v(x)} \,] = \frac{1}{\Delta x}\, [\, \frac{u(x+\Delta x)}{v(x+\Delta x)} - \frac{u(x)}{v(x)} \,] \qquad (6.3\text{-}2)$$

$$D_{\Delta x}[\, \frac{u(x)}{v(x)} \,] = \frac{1}{\Delta x}\, [\, \frac{u(x+\Delta x)v(x) - u(x)v(x+\Delta x)}{v(x)v(x+\Delta x)} \,] \qquad (6.3\text{-}3)$$

$$D_{\Delta x}[\, \frac{u(x)}{v(x)} \,] = \frac{1}{\Delta x}\, [\, \frac{u(x+\Delta x)v(x) - u(x)v(x+\Delta x) + u(x)v(x) - u(x)v(x)}{v(x)v(x+\Delta x)} \,] \qquad (6.3\text{-}4)$$

$$D_{\Delta x}[\, \frac{u(x)}{v(x)} \,] = \frac{v(x)[\frac{u(x+\Delta x) - u(x)}{\Delta x}] - u(x)[\frac{v(x+\Delta x) - v(x)}{\Delta x}]}{v(x)v(x+\Delta x)} \qquad (6.3\text{-}5)$$

Thus

From Eq 6.3-5 and Eq 6.3-1

**The discrete derivative of a division of two functions is:**

$$\mathbf{D_{\Delta x}[\, \frac{u(x)}{v(x)} \,] = \frac{v(x)D_{\Delta x}u(x) - u(x)D_{\Delta x}v(x)}{v(x)v(x+\Delta x)}} \qquad (6.3\text{-}6)$$

## Section 6.4: Derivation of the Discrete Function Chain Rule

$F(v)$ = function of $v$
$g_{\Delta x}(x)$ = discrete Interval Calculus function of $x,\Delta x$
$v = g_{\Delta x}(x)$
$\Delta v = \Delta x D_{\Delta x} g_{\Delta x}(x)$
$\Delta f(x) = f(x+\Delta x) - f(x)$

$D_{\Delta x} = \dfrac{\Delta}{\Delta x}$

$D_{\Delta x} f(x) = \dfrac{\Delta f(x)}{\Delta x} = \dfrac{f(x+\Delta x) - f(x)}{\Delta x}$ = discrete derivative of $f(x)$

$$\left[\frac{\Delta F(v)}{\Delta v}\right] = \frac{\frac{\Delta F(v)}{\Delta x}}{\frac{\Delta v}{\Delta x}} = \left[\frac{\Delta F(v)}{\Delta x}\right]\left[\frac{1}{\frac{\Delta v}{\Delta x}}\right] \tag{6.4-1}$$

$$\left[\frac{\Delta F(v)}{\Delta v}\right]\Big|_{\substack{v=g_{\Delta x}(x) \\ \Delta v = \Delta x D_{\Delta x} g_{\Delta x}(x)}} = \frac{\Delta}{\Delta x}\left[F(v)\,\Big|_{v=g_{\Delta x}(x)}\right]\left[\frac{1}{\frac{\Delta v}{\Delta x}}\right]\Big|_{v=g_{\Delta x}(x)} \tag{6.4-2}$$

$$\left[\frac{\Delta F(v)}{\Delta v}\right]\Big|_{\substack{v=g_{\Delta x}(x) \\ \Delta v = \Delta x D_{\Delta x} g_{\Delta x}(x)}} = \frac{\Delta}{\Delta x}\left[F(v)\,\Big|_{v=g_{\Delta x}(x)}\right]\left[\frac{1}{\frac{\Delta g_{\Delta x}(x)}{\Delta x}}\right] \tag{6.4-3}$$

$$\left[D_{\Delta v} F(v)\right]\Big|_{\substack{v=g_{\Delta x}(x) \\ \Delta v = \Delta x D_{\Delta x} g_{\Delta x}(x)}} = D_{\Delta x}\left[F(v)\,\Big|_{v=g_{\Delta x}(x)}\right]\left[\frac{1}{D_{\Delta x} g_{\Delta x}(x)}\right] \tag{6.4-4}$$

Rearranging the terms of Eq 6.4-4

$$D_{\Delta x}\left[F(v)\,\Big|_{v=g_{\Delta x}(x)}\right] = \left[D_{\Delta v} F(v)\right]\Big|_{\substack{v=g_{\Delta x}(x) \\ \Delta v = \Delta x D_{\Delta x} g_{\Delta x}(x)}} \left[D_{\Delta x} g_{\Delta x}(x)\right] \tag{6.4-5}$$

The discrete derivative is found either from its definition, $D_{\Delta x}f(x) = \dfrac{f(x+\Delta x) - f(x)}{\Delta x}$, or from the discrete differentiation of a discrete function identity. For example:

$$D_{\Delta x}v^2 = \frac{(v+\Delta v)^2 - v^2}{\Delta v} = \frac{v^2+2v\Delta v+\Delta v^2-v^2}{\Delta v} = 2v+\Delta x$$

or using the $v^2$ discrete function identity

$$v^2 = v(v-\Delta v)+v\Delta v$$

$$D_{\Delta v}v^2 = D_{\Delta v}[v(v-\Delta v)+v\Delta v] = 2v+\Delta v$$

Then

**The Discrete Function Chain Rule is:**

$$\mathbf{D_{\Delta x}\left[F(v)\,\Big|_{v=g_{\Delta x}(x)}\right] = \left[D_{\Delta v}F(v)\right]\Big|_{\substack{v=g_{\Delta x}(x)\\ \Delta v = \Delta x D_{\Delta x}g_{\Delta x}(x)}}\left[D_{\Delta x}g_{\Delta x}(x)\right]} \qquad (6.4\text{-}6)$$

**where**

**$F(v)$ = function of v**

**$g_{\Delta x}(x)$ = discrete Interval Calculus function of x**

**$v = g_{\Delta x}(x)$**

**$\Delta v = \Delta x D_{\Delta x}g_{\Delta x}(x)$**

**$\Delta x,\Delta v$ = interval increments**

**$D_{\Delta x}$ = the discrete derivative with respect to $\Delta x$**

## Section 6.5: The derivation of the $e_{\Delta x}(a,x)$ , $e^{ax}$ identity relationship

$$e^{ax} = e^{ax} \qquad (6.5\text{-}1)$$

$$e^{ax} = [1 + (\frac{e^{a\Delta x} - 1}{\Delta x})\Delta x]^{\frac{x}{\Delta x}} \qquad (6.5\text{-}2)$$

where

$x = m\Delta x$ , m = integers

$\Delta x$ = x increment

$$e_{\Delta x}(c,x) = [1 + c\Delta x]^{\frac{x}{\Delta x}} \qquad (6.5\text{-}3)$$

From Eq 6.5-2 and Eq 6.5-3

**The $e_{\Delta x}(a,x)$ , $e^{ax}$ identity relationship is:**

$$\mathbf{e^{ax} = e_{\Delta x}\left(\frac{e^{a\Delta x}-1}{\Delta x},x\right)} \tag{6.5-4}$$

Comment - The function, $e_{\Delta x}(\frac{e^{a\Delta x}-1}{\Delta x},x)$ , is commonly used in discrete calculus to represent the function, $e^{ax}$. Both functions have the same value at $x = m\Delta x$ , m = integers.

**Section 6.6: Derivation of the function,** $\displaystyle \mathbf{\sum_{\Delta x\; x=\mp\infty}^{\pm\infty}\frac{1}{x^n} = \frac{-\pi}{(n-1)!\Delta x}\frac{d^{n-1}}{dx^{n-1}}\cot\left(\frac{\pi x}{\Delta x}\right)\Big|_{x=-x_i}}$

**$x = x_i + m\Delta x$ , m = integers**

Consider the summation

$$S = \sum_{\Delta x\; x=\mp\infty}^{\pm\infty}\frac{1}{x^n} \tag{6.6-1}$$

where $x = x_i + m\Delta x$ , m = integers

$x_i$ = value of x

n = 1,2,3,…

The expanded series of Eq 6.6-1 is:

$$S = \ldots + \frac{1}{(x_i-3\Delta x)^n} + \frac{1}{(x_i-2\Delta x)^n} + \frac{1}{(x_i-\Delta x)^n} + \frac{1}{(x_i)^n} + \frac{1}{(x_i+\Delta x)^n} + \frac{1}{(x_i+2\Delta x)^n} + \frac{1}{(x_i+3\Delta x)^n} + \ldots \tag{6.6-2}$$

Changing the form of the summation of Eq 6.6-1

$$S = \sum_{1\; r=\mp\infty}^{\pm\infty}\frac{1}{(x_i+r\Delta x)^n}\ ,\quad \text{where } r = \text{integers},\quad n = 1,2,3,\ldots$$

$$S = \frac{1}{(\Delta x)^n}\sum_{r=\mp\infty}^{\pm\infty}\frac{1}{\left(r+\frac{x_i}{\Delta x}\right)^n}\ ,\quad \text{where } r = \text{integers},\quad n = 1,2,3,\ldots \tag{6.6-3}$$

Note that for Eq 6.6-3, $r = 0, \pm 1, \pm 2, \pm 3, \ldots$ and the summation function has a single pole at $-\frac{x}{\Delta x}$ of order n.

This summation can be evaluated using complex variable mathematics residue theory.

<u>Comment</u> – Summing from $-\infty$ to $+\infty$ or from $+\infty$ to $-\infty$ yields the same result.

In general, if a summation series of the form, $\sum_{r=-\infty}^{+\infty} f(r)$ , converges, then:

$$\sum_{r=-\infty}^{+\infty} f(r) = -\{\text{sum of the residues of } \pi\cot\pi z\ f(z) \text{ at all the poles of } f(z)\} \qquad (6.6\text{-}4)$$

> Eq 6.6-4, was derived in a book entitled, Schaum's Outline Series Theory and Problems of Complex Variables, written by Murray R. Spiegel, printed by McGraw-Hill Book Company, and copyright in 1964.

The series of Eq 6.6-3 is seen to be convergent for all n, n = 1,2,3,… provided $-\frac{x_i}{\Delta x} \neq$ integer.

Thus, the relationship of Eq 6.6-4 may be utilized to evaluate the summation of Eq 6.6-3.

Find the residue of $\dfrac{\pi\cot\pi r}{(r + \frac{x_i}{\Delta x})^n}$ at the pole $r = -\frac{x}{\Delta x}$

To find the residue, R, for all n, the residue calculation formula is used

$$R = \lim_{r \to -\frac{x_i}{\Delta x}} \frac{1}{(n-1)!} \frac{d^{n-1}}{dr^{n-1}} \frac{(r + \frac{x_i}{\Delta x})^n \pi \cot\pi r}{(r + \frac{x_i}{\Delta x})^n} = \frac{\pi}{(n-1)!} \frac{d^{n-1}}{dr^{n-1}} \frac{(r + \frac{x_i}{\Delta x})^n \cot\pi r}{(r + \frac{x_i}{\Delta x})^n} \Bigg|_{r \to -\frac{x_i}{\Delta x}} \qquad (6.6\text{-}5)$$

Canceling

$$R = \frac{\pi}{(n-1)!} \frac{d^{n-1}}{dr^{n-1}} \cot\pi r \,\Big|_{r \to -\frac{x_i}{\Delta x}} \qquad (6.6\text{-}6)$$

From 6.6-3, 6.6-4 and 6.6-6

$$S = \frac{1}{(\Delta x)^n} \sum_{r=\mp\infty}^{\pm\infty} \frac{1}{( r + \frac{x}{\Delta x} )^n} = -\frac{1}{(\Delta x)^n} \frac{\pi}{(n-1)!} \frac{d^{n-1}}{dr^{n-1}} \cot \pi r \Big|_{r \to -\frac{x_i}{\Delta x}} \qquad (6.6\text{-}7)$$

Substituting Eq 6.6-1 into Eq 6.6-7

$$\sum_{\substack{x=\mp\infty \\ \Delta x}}^{\pm\infty} \frac{1}{x^n} = -\frac{1}{(\Delta x)^n} \frac{\pi}{(n-1)!} \frac{d^{n-1}}{dr^{n-1}} \cot \pi r \Big|_{r \to -\frac{x_i}{\Delta x}} \qquad (6.6\text{-}8)$$

where

$x = x_i + m\Delta x$ , $m$ = integers

$x_i$ = value of x

$n = 1,2,3,\ldots$

$-\frac{x_i}{\Delta x} \neq$ integer.

Let

$$r = \frac{x}{\Delta x} \qquad (6.6\text{-}9)$$

Using the Chain rule

$$\frac{d}{dr} = \frac{1}{\frac{dr}{dx}} \frac{d}{dx} = \Delta x \frac{d}{dx}$$

$$\frac{d}{dr} = \Delta x \frac{d}{dx}$$

$$\frac{d}{dr}\frac{d}{dr} = \Delta x \frac{d}{dx} \Delta x \frac{d}{dx} = \Delta x^2 \frac{d}{dx}$$

$$\frac{d^2}{dr^2} = \Delta x^2 \frac{d^2}{dx^2}$$

In general

$$\frac{d^{n-1}}{dr^{n-1}} = \Delta x^{n-1} \frac{d^{n-1}}{dx^{n-1}} , \quad n = 1,2,3,\ldots \qquad (6.6\text{-}10)$$

$$r = \frac{x}{\Delta x} \rightarrow - \frac{x_i}{\Delta x} \quad \text{then} \quad x \rightarrow -x_i \tag{6.6-11}$$

Substituting Eq 6.6-10 and Eq 6.6-11 into Eq 6.6-8

$$\sum_{\Delta x \; x=\mp\infty}^{\pm\infty} \frac{1}{x^n} = \frac{-\pi}{(n-1)!\Delta x} \frac{d^{n-1}}{dx^{n-1}} \cot\frac{\pi x}{\Delta x}\Big|_{x\rightarrow -x_i} \tag{6.6-12}$$

**where**

**$x = x_i + m\Delta x$ , $m$ = integers**

**$x_i$ = value of x**

**$n = 1,2,3,\ldots$**

**$-\frac{x_i}{\Delta x} \neq$ integer.**

# CHAPTER 7

# Demonstration of Interval Calculus' Unique and Useful Mathematical Methods

## Section 7.1: Introduction to Chapter 7

Interval Calculus is a discrete calculus with a specific unique notation. It has been developed as a mathematical tool to solve applied mathematical problems efficiently, perhaps even more efficiently than other mathematical methods commonly used. Of course, the mathematical method of preference depends on the user. This chapter is written to point out some of the unique and useful features of Interval Calculus in the solution of various types of mathematical problems.

This chapter is composed of problem examples solved by what is considered to be interesting Interval Calculus mathematical methods. The derivation of these methods can be found in the previous chapters. Only the application of these methods is presented here. The heading of each example will describe the problem to be solved. In the previous chapters many problems have been presented and solved for demonstration purposes. Some of these problem solutions reappear here as examples. The emphasis here will be to point out the features of Interval Calculus which may provide a useful alternative or improved method for solving problems presently solved by other mathematical methodologies. The order in which the problem examples are presented is of no particular significance.

## Section 7.2: Problem Examples

Below is a listing of problem examples solved in this chapter using Interval Calculus methods.

Example 7.1 Evaluation of the summation, $\displaystyle {}_{\Delta x}\sum_{x=x_1}^{x_2}\frac{1}{(x+a)^n}$ , using the derived lnd(n,Δx,x) function

Example 7.2 Evaluation of the summation, $\displaystyle {}_{3}\sum_{x=1}^{\infty}\frac{1}{x} - {}_{3}\sum_{x=2}^{\infty}\frac{1}{x}$ , using the lnd(n,Δx,x) function

Example 7.3 Evaluation of a polynomial summation

Example 7.4 Evaluation of useful mathematical functions using the lnd(n,Δx,x) function

Example 7.5 Evaluation of the Riemann and Hurwitz Zeta Functions using the derived lnd(n,Δx,x) Function

Example 7.6 Derivation of an interesting Zeta Function relationship

Example 7.7 Use of the lnd(n,Δx,x) function to calculate the Gamma Function

Example 7.8 Calculation of an alternating sign summation using the derived Alternating Sign Discrete Calculus Summation Equation

Example 7.9 A calculation using the derived Discrete Calculus Summation Equation

Example 7.10 Finding the Z Transform of a function using Interval Calculus

Example 7.11 Using the complex plane area calculation equation, $A_c = \left| \frac{j}{2} \sum_{p=0}^{N} (-1)^p c_p^{\,2} \right|$, to calculate the area within a discrete closed contour in the complex plane

Example 7.12 Using the $K_{\Delta x}$ Transform to solve difference equations

Example 7.13 Solving servomechanism control problems using the $K_{\Delta x}$ Transform and Nyquist Methods

Example 7.14 Evaluating a summation which has no division by zero terms

Example 7.15 Evaluating a summation which has a division by zero term

Example 7.16 A demonstration of discrete differentiation

Example 7.17 The solution of a differential difference equation using four different methods

Solution 1 Use of the Method of Undetermined Coefficients to solve a differential difference equation

Solution 2 Use of the $K_{\Delta x}$ Transform Method to solve a differential difference equation

Solution 3 Use of the Method of Variation of Parameters to solve a differential difference equation

Solution 4 Use of the Method of Related Functions to solve a differential difference equation

Example 7.18 Finding the $K_{\Delta x}$ Transform, $K_{\Delta x}[e_{\Delta x}(a,x)] = \frac{1}{s-a}$ two different ways, and the Inverse $K_{\Delta x}$ Transform, $K_{\Delta x}^{-1}[\frac{1}{s-a}] = e_{\Delta x}(a,x)$ using the Inverse $K_{\Delta x}$ Transform

Example 7.19 Determining for what range of $\Delta x$ a control system with a differential difference equation, $D_{\Delta x}^{2} y_{\Delta x}(x) + 20 D_{\Delta x} y_{\Delta x}(x) + 125 y_{\Delta x}(x) = 0$, is stable

Example 7.20 Evaluating a Z Transform defined closed loop system using $K_{\Delta t}$ Transforms

Example 7.21 Evaluating in detail several Zeta Function summations, $\underset{\Delta x}{\sum_{x=x_1}^{x_2}} \frac{1}{x^n}$ ( $x_1$,$x_2$ may be infinite), using the lnd(n,Δx,x) function calculation program, LNDX, and some pertinent equations, and demonstrating the validity of all of the equations and concepts presented.

Example 7.22 Finding the derivative, $D_{\Delta x}[\ln_{\Delta x}{}^2 x]$, using the Discrete Function Chain Rule and the Discrete Derivative of the Product of two Functions Equation

**Example 7.1** Evaluation of the summation, $\displaystyle {}_{\Delta x}\sum_{x=x_1}^{x_2} \frac{1}{(x+a)^n}$, using the derived lnd(n,Δx,x) function

$$ {}_{\Delta x}\sum_{x=x_1}^{x_2} \frac{1}{(x+a)^n} = \frac{1}{\Delta x}\ {}_{\Delta x}\int_{x_1}^{x_2+\Delta x} \frac{1}{(x+a)^n}\,\Delta x = -\frac{1}{\Delta x}\,\mathrm{lnd}(n,\Delta x,x+a)\Big|_{x_1}^{x_2+\Delta x} \quad , \quad n \neq 1 $$

Problem Description

Find the summation, $\displaystyle {}_{1+i}\sum_{x=-1-3i}^{2} \frac{1}{(x-2+i)^{2.3+3.7i}}$, using the lnd(n,Δx,x) function.

Let , $\Delta x = 1+i$, $a = -2+i$, $x_1 = -1-3i$, $x_2 = 2$

$$ {}_{1+i}\sum_{x=-1-3i}^{2} \frac{1}{(x-2+i)^{2.3+3.7i}} = \frac{1}{1+i}\ {}_{1+i}\int_{-1-3i}^{3+i} \frac{1}{(x-2+i)^{2.3+3.7i}}\,\Delta x = -\frac{1}{1+i}\,\mathrm{lnd}(2.3+3.7i,1+i,x-2+i)\Big|_{-1-3i}^{3+i} $$

$$ {}_{1+i}\sum_{x=-1-3i}^{2} \frac{1}{(x-2+i)^{2.3+3.7i}} = \frac{1}{1+i}\,[-\,\mathrm{lnd}(2.3+3.7i,1+i,1+2i) + \mathrm{lnd}(2.3+3.7i,1+i,-3-2i)] $$

Evaluating using the lnd(n,Δx,x) computer computation program, LNDX

$$ {}_{1+i}\sum_{x=-1-3i}^{2} \frac{1}{(x-2+i)^{2.3+3.7i}} = \frac{1}{1+i}\,[(-3.0442574254041-12.2176183323067i)+ $$
$$ (+155626.32956516287-24853.3785053376592i)] $$

$$ {}_{1+i}\sum_{x=-1-3i}^{2} \frac{1}{(x-2+i)^{2.3+3.7i}} = 65378.8445920337525-90244.4407157037186i $$

Checking

$$\sum_{x=-1-3i}^{2}{}_{1+i}\frac{1}{(x-2+i)^{2.3+3.7i}} = \frac{1}{(-3-2i)^{2.3+3.7i}} + \frac{1}{(-2-i)^{2.3+3.7i}} + \frac{1}{(-1)^{2.3+3.7i}} + \frac{1}{(i)^{2.3+3.7i}}$$

$$\sum_{x=-1-3i}^{2}{}_{1+i}\frac{1}{(x-2+i)^{2.3+3.7i}} = 65378.8445920337525-90244.4407157037186i$$

Good check

**Example 7.2** Evaluation of the summation, $\sum_{3\,x=1}^{\infty} \frac{1}{x} - \sum_{3\,x=2}^{\infty} \frac{1}{x}$, using the lnd(n,Δx,x) function

Problem Description

Evaluate the summation, $\sum_{3\,x=1}^{\infty} \frac{1}{x} - \sum_{3\,x=2}^{\infty} \frac{1}{x} = 1 - \frac{1}{2} + \frac{1}{4} - \frac{1}{5} + \frac{1}{7} - \frac{1}{8} + \frac{1}{10} - \frac{1}{11} + \ldots$, using the lnd(n,Δx,x) function.

The lnd(n,Δx,x) function evaluates the following summation as shown

$$\sum_{\Delta x\; x=x_1}^{x_2} \frac{1}{x+a} = \frac{1}{\Delta x} \text{lnd}(1,\Delta x,x+a)\Big|_{x_1}^{x_2+\Delta x}$$

$N \rightarrow \infty$
$\Delta x = 3$

$$\sum_{3\,x=1}^{\infty} \frac{1}{x} - \sum_{3\,x=2}^{\infty} \frac{1}{x} = \lim_{N\rightarrow\infty} \left[ \sum_{3\,x=1}^{N} \frac{1}{x} - \sum_{3\,x=2}^{N+1} \frac{1}{x} \right]$$

$$\sum_{3\,x=1}^{\infty} \frac{1}{x} - \sum_{3\,x=2}^{\infty} \frac{1}{x} = \lim_{N\rightarrow\infty} \frac{1}{3} \left[ \text{lnd}(1,3,x)\,\Big|_{1}^{N+3} - \text{lnd}(1,3,x)\,\Big|_{2}^{N+4} \right]$$

$$\sum_{3\,x=1}^{\infty} \frac{1}{x} - \sum_{3\,x=2}^{\infty} \frac{1}{x} = \frac{1}{3} [-\text{lnd}(1,3,1) + \text{lnd}(1,3,2)] + \frac{1}{3} \lim_{N\rightarrow\infty} [\text{lnd}(1,3,N+3)) - \text{lnd}(1,3,N+4)]$$

But

$\lim_{N\rightarrow\infty} [\text{lnd}(1,3,N+3)) - \text{lnd}(1,3,N+4)] \rightarrow 0$

Then

$$\sum_{3\,x=1}^{\infty} \frac{1}{x} - \sum_{3\,x=2}^{\infty} \frac{1}{x} = \frac{1}{3} [-\text{lnd}(1,3,1) + \text{lnd}(1,3,2)]$$

Evaluating using the lnd(n,Δx,x) computer computation program, LNDX

$$\sum_{x=1}^{\infty}{}_3\frac{1}{x} - \sum_{x=2}^{\infty}{}_3\frac{1}{x} = \frac{1}{3}[2.5548181151192734623899069 8 - .7410187508850556117958287 2 ]$$

$$\sum_{x=1}^{\infty}{}_3\frac{1}{x} - \sum_{x=2}^{\infty}{}_3\frac{1}{x} = .6045997880780726168646927 5\ldots = \frac{\pi}{3\sqrt{3}}\ldots$$

$$\sum_{x=1}^{\infty}{}_3\frac{1}{x} - \sum_{x=2}^{\infty}{}_3\frac{1}{x} = 1 - \frac{1}{2} + \frac{1}{4} - \frac{1}{5} + \frac{1}{7} - \frac{1}{8} + \frac{1}{10} - \frac{1}{11} + \ldots = \frac{\pi}{3\sqrt{3}}\ldots$$

Checking using a computer solution to the summation, $\sum_{x=1}^{10000000}{}_3(\frac{1}{x} - \frac{1}{x+1})$

$$\sum_{x=1}^{10000000}{}_3(\frac{1}{x} - \frac{1}{x+1}) = .6045997547447$$

good check

This summation is slow to converge.

**Example 7.3** Evaluation of a polynomial summation

Problem Description

Evaluate the sum, $\sum_{.5\ x=1}^{3}\frac{(x+1)^2}{x^2+1}$, using the derived lnd(n,Δx,x) function.

Expanding using a partial fraction expansion

$$\frac{(x+1)^2}{x^2+1} = 1 + \frac{1}{x-j} + \frac{1}{x+j}$$

$$\sum_{.5\ x=1}^{3}\frac{(x+1)^2}{x^2+1} = \sum_{.5\ x=1}^{3}\left(1 + \frac{1}{x-j} + \frac{1}{x+j}\right) = \frac{1}{.5}\ {}_{.5}\!\int_{1}^{3+.5}\Delta x + \frac{1}{.5}\ {}_{.5}\!\int_{1}^{3+.5}\frac{\Delta x}{x-j} + \frac{1}{.5}\ {}_{.5}\!\int_{1}^{3+.5}\frac{\Delta x}{x+j}$$

The lnd(n,Δx,x) function evaluates the following summation and integral as shown

$$\sum_{\Delta x\ x=x_1}^{x_2}\frac{1}{(x+a)^n} = \frac{1}{\Delta x}\ {}_{\Delta x}\!\int_{x_1}^{x_2+\Delta x}\frac{1}{(x+a)^n}\,\Delta x = \pm\frac{1}{\Delta x}\,\mathrm{lnd}(n,\Delta x,x+a)\Big|_{x_1}^{x_2+\Delta x}\ , \quad -\text{ for } n \neq 1,\ +\text{ for } n = 1$$

$$\sum_{.5\ x=1}^{3}\frac{(x+1)^2}{x^2+1} = \frac{1}{.5}\left[\ -\mathrm{lnd}(0,.5,x)\Big|_{1}^{3.5} + \mathrm{lnd}(1,.5,x-j)\Big|_{1}^{3.5} + \mathrm{lnd}(1,.5,x+j)\Big|_{1}^{3.5}\ \right]$$

Note – The lnd(n,Δx,x) function is preceded by a – for n ≠ 1 and by a + for n = 1

$$\sum_{.5\ x=1}^{3}\frac{(x+1)^2}{x^2+1} = 2\,[\ -\mathrm{lnd}(0,.5,3.5) + \mathrm{lnd}(0,.5,1) + \mathrm{lnd}(1,.5,3.5-j) - \mathrm{lnd}(1,.5,1-j)\ ] +$$

$$2\,[\ \mathrm{lnd}(1,.5,3.5+j) - \mathrm{lnd}(1,.5,1+j)\ ]$$

Evaluating using the lnd(n,Δx,x) computer computation program, LNDX

$$\sum_{.5\ x=1}^{3}\frac{(x+1)^2}{x^2+1} = 5 + (\ 2.0063660477453580 + 1.2456233421750663i\ ) +$$

$$(\ 2.0063660477453580 - 1.2456233421750663i\ )$$

$$\sum_{\substack{x=1 \\ .5}}^{3} \frac{(x+1)^2}{x^2+1} = 9.0127320954907162$$

Checking

$$\sum_{\substack{x=1 \\ .5}}^{3} \frac{(x+1)^2}{x^2+1} = \frac{2^2}{2} + \frac{2.5^2}{3.25} + \frac{3^2}{5} + \frac{3.5^2}{7.25} + \frac{4^2}{10} = 9.0127320954907162$$

Good check

**Example 7.4** Evaluation of useful mathematical functions using the lnd(n,Δx,x) function

Problem Description

List some of the functions evaluated by the lnd(n,Δx,x) function.

## List of Functions Evaluated by the lnd(n, Δx, x) Function

**(where n, Δx, x are real or complex values)**

| # | FUNCTION | LND(n,Δx,x) EQUATION | DESCRIPTION |
|---|---|---|---|
| **1** | $\sum_{\substack{x=x_1 \\ \Delta x}}^{x_2} \frac{1}{x}$ | $\frac{1}{\Delta x} \text{lnd}(1,\Delta x,x) \Big\vert_{x_1}^{x_2+\Delta x}$<br>Any summation term with x = 0 is excluded | Summation |
| **2** | $\sum_{\substack{x=x_1 \\ \Delta x}}^{x_2} \frac{1}{x^n}$ | $-\frac{1}{\Delta x} \text{lnd}(n,\Delta x,x) \Big\vert_{x_1}^{x_2+\Delta x}, \quad n \neq 1$<br>Any summation term with x = 0 is excluded | Summation |
| **3** | $\sum_{\substack{x=x_1 \\ \Delta x}}^{x_2} (-1)^{\frac{x-x_1}{\Delta x}} \frac{1}{x^n}$ | $\frac{\alpha(n)}{2\Delta x} \left[ -\text{lnd}(n,2\Delta x,x) \Big\vert_{x_1}^{x_2+\Delta x} + \text{lnd}(n,2\Delta x,x) \Big\vert_{x_1+\Delta x}^{x_2+2\Delta x} \right]$<br>$x = x_1, x_1+\Delta x, x_1+2\Delta x, x_1+3\Delta x, \ldots, x_2-\Delta x, x_2$<br>$x_2 = x_1 + (2m-1)\Delta x, \quad m = 1,2,3,\ldots$<br>Δx = x increment<br>n = any value<br>$\alpha(n) = \begin{cases} 1 & n \neq 1 \\ -1 & n = 1 \end{cases}$<br>Any summation term with x = 0 is excluded | Summation |
| **4** | $\sum_{\Delta x,\ x=x_i}^{\pm\infty} \frac{1}{x^n}$ | $\frac{1}{\Delta x}\,\mathrm{lnd}(n,\Delta x,x_i)$<br>$+\infty$ for $Re(\Delta x)>0$ or [ $Re(\Delta x)=0$ and $Im(\Delta x)>0$ ]<br>$-\infty$ for $Re(\Delta x)<0$ or [ $Re(\Delta x)=0$ and $Im(\Delta x)<0$ ]<br>$x = x_i + m\Delta x$, $m = 0,1,2,3,\ldots$<br>$Re(n)>1$<br><br>Any summation term with $x = 0$ is excluded | Summation |
| **5** | $\sum_{\Delta x,\ x=x_i}^{\pm\infty} (-1)^{\frac{x-x_i}{\Delta x}} \frac{1}{x^n}$ | $\frac{\alpha(n)}{2\Delta x}$ [ $-\mathrm{lnd}(n,2\Delta x,x_i+\Delta x) + \mathrm{lnd}(n,2\Delta x,x_i)$ ]<br>n = all values<br>$\alpha(n) = -1$ for $n = 1$<br>$\alpha(n) = +1$ for $n \neq 1$<br>$+\infty$ for $Re(\Delta x)>0$ or [ $Re(\Delta x)=0$ and $Im(\Delta x)>0$ ]<br>$-\infty$ for $Re(\Delta x)<0$ or [ $Re(\Delta x)=0$ and $Im(\Delta x)<0$ ]<br>$x = x_i + m\Delta x$, $m = 0,1,2,3,\ldots$<br>$Re(n)>0$<br><br>Any summation term with $x = 0$ is excluded | Summation |
| **6** | $\sum_{\Delta x,\ x=\mp\infty}^{\pm\infty} \frac{1}{x^n}$ | $\frac{\alpha(n)}{\Delta x}$ [$\mathrm{lnd}(n,\Delta x,x_i) - \mathrm{lnd}(n,-\Delta x,x_i-\Delta x)$]<br>$\alpha(n) = -1$ for $n = 1$<br>$\alpha(n) = +1$ for $n \neq 1$<br>$-\infty$ to $+\infty$ for $Re(\Delta x)>0$ or [$Re(\Delta x)=0$ and $Im(\Delta x)>0$]<br>$+\infty$ to $-\infty$ for $Re(\Delta x)<0$ or [$Re(\Delta x)=0$ and $Im(\Delta x)<0$]<br>$x = x_i+m\Delta x$, $m = 0,1,2,3,\ldots$<br>$x_i$ = value of x<br>$n,x,x_i,\Delta x$ = real or complex values<br>$Re(n)\geq 1$<br><br>Any summation term with $x = 0$ is excluded<br><br><u>Comment</u> – If the value of $x_i$ is changed to $x_i + r\Delta x$, r = integers, the summation value remains the same. | Summation |
| **7** | $\sum_{\Delta x \; x=\mp\infty}^{\pm\infty} (-1)^{\frac{x-x_i}{\Delta x}} \frac{1}{x^n}$ | $\frac{\alpha(n)}{2\Delta x}$ [-lnd(n,2Δx,$x_i$+Δx)+lnd(n,2Δx,$x_i$)-lnd(n,-2Δx,$x_i$-2Δx)+lnd(n,-2Δx,$x_i$-Δx)]<br>α (n) = -1 for n = 1<br>α (n) = +1 for n ≠ 1<br>-∞ to +∞ for Re(Δx)>0 or [Re(Δx)=0 and Im(Δx)>0]<br>+∞ to -∞ for Re(Δx)<0 or [Re(Δx)=0 and Im(Δx)<0]<br>$x = x_i + m\Delta x$, m = 0,1,2,3,…<br>$x_i$ = value of x<br>Re(n)>0<br><br>Any summation term with x = 0 is excluded<br><br>Comment – If the value of $x_i$ is changed to $x_i$ + rΔx, r = integers, the summation value remains the same. | Summation |
| **8** | $\sum_{x=x_1}^{x_2} \ln(1+x)$ | $\left[ x(\ln x-1) + \frac{\text{lnd}(1,1,x+\frac{1}{2})}{2} + \sum_{n=2}^{\infty} \left\{ \frac{\text{lnd}(2n-2,1,x+\frac{1}{2})}{(2n-1)(2)^{2n-2}} - \frac{\text{lnd}(2n-1,1,x+\frac{1}{2})}{(2n-1)(2)^{2n-1}} \right\} \right] \Big|_{x_1}^{x_2+1}$<br>Re(1/x) ≥ –1, 1/x ≠ -1 | Natural Log Summation |
| **9** | ln(x) | 1) $\lim_{m\to\infty}$ [ lnd(1,$\frac{x-1}{m}$,x) – lnd(1,$\frac{x-1}{m}$,1) ]<br>$\Delta x = \frac{x-1}{m}$, m = positive integers, m→∞<br><br>For x = positive real or complex values<br><br>2) $\lim_{m\to\infty}$ [ lnd(1,$\frac{x+1}{m}$,x) – lnd(1,$\frac{x+1}{m}$, -1) ] + πi<br>$\Delta x = \frac{x+1}{m}$, m = positive integers, m→∞<br><br>For x = negative real values<br>x ≠ 0 | Natural Log |
| **10** | Γ(x) | $\sqrt{2\pi}\, e^{\lim_{n\to 0} \frac{\text{lnd}(n,1,x)-\text{lnd}(0,1,x)}{n}}$<br>n = real number | Gamma Function |
| **11** | LnΓ(x) | $\lim_{n\to 0} \frac{\text{lnd}(n,1,x) - \text{lnd}(0,1,x)}{n} + \frac{\ln 2\pi}{2}$<br>n = real number, n→0 | Natural Log of the Gamma Function |
| **12** | ψ(x) | lnd(1,1,x) - γ<br>where γ = Euler's constant | Bigamma Function |
| **13** | $\psi^{(n)}(x)$ | $(-1)^{n+1}$ n! lnd(n+1,1,x)<br>where n = 1, 2, 3, 4, … | Polygamma Functions |
| **14** | $\zeta(n) = \sum_{x=1}^{\infty} \frac{1}{x^n}$ | $\text{lnd}(n,1,1)$<br>where $\text{Re}(n)>1$ | Riemann Zeta Function |
| **15** | $\zeta(n,x) = \sum_{x=x_i}^{\infty} \frac{1}{x^n}$ | $\text{lnd}(n,1,x)$<br>where $\text{Re}(n)>1$<br>Any summation term with x = 0 is excluded | Hurwitz Zeta Function |
| **16** | $\beta(w,z)$ | $\sqrt{2\pi}\, e^{\lim_{n\to 0} \frac{\text{lnd}(n,1,w)-\text{lnd}(0,1,w)+\text{lnd}(n,1,z)-\text{lnd}(0,1,z)-\text{lnd}(n,1,w+z)+\text{lnd}(0,1,w+z)}{n}}$ | Beta Function |
| **17** | $B_n$ | $-n\,[\text{lnd}(1-n,1,1)]$<br>where n = 2, 4, 6, 8, … | Bernoulli Constants |
| **18** | $\tan(ax)$ | $\dfrac{\pi}{\text{lnd}(1,1,1-\frac{ax}{\pi}) - \text{lnd}(1,1,\frac{ax}{\pi})}$<br>for $ax = n\pi$□ n=integer, $\tan(ax)=0$ | Tangent Function |
| **19** | $\tanh(ax)$ | $\dfrac{-j\pi}{\text{lnd}(1,1,1-\frac{jax}{\pi}) - \text{lnd}(1,1,\frac{jax}{\pi})}$<br>for $ax = 0$□ $\tanh(ax)=0$ | Hyperbolic Tangent Function |
| **20** | $X^N$ | $\alpha(N)[\,\text{lnd}(-N,1,x+1) - \text{lnd}(-N,1,x)\,]$<br>$\alpha(N) = +1$ for $N = -1$<br>$\alpha(N) = -1$ for $N \neq -1$ | Variable to a constant power |
| **21** | $\pi$ | $-\text{lnd}(1,4,4m+5) + \text{lnd}(1,-4,4m+1)$<br>m = integer | Pi |
| **22** | $\gamma$ | $\lim_{n\to 0}\,[\,\text{lnd}(1+n,1,1) - \frac{1}{n}\,]$<br>n = real number | Euler's Constant |
| **23** | $\tan_{\Delta x}(a,x)$ | $\dfrac{\pi}{\text{lnd}(1,1,1-\frac{bx}{\pi}) - \text{lnd}(1,1,\frac{bx}{\pi})}$<br>$b = \dfrac{\tan^{-1}(a\Delta x)}{\Delta x}$<br>for $bx = n\pi$ n=integer, $\tan_{\Delta x}(a,x)=0$ | Discrete Tangent Function |
| **24** | $\tanh_{\Delta x}(a,x)$ | $\dfrac{-j\pi}{\text{lnd}(1,1,1-\frac{bx}{\pi}) - \text{lnd}(1,1,\frac{bx}{\pi})}$<br>$b = \dfrac{\tan^{-1}(ja\Delta x)}{\Delta x}$<br>for $bx = 0$□ $\tanh_{\Delta x}(a,x)=0$ | Discrete Hyperbolic Tangent Function |

Note -A computer program, LNDX, has been written to calculate the function, lnd(n,Δx,x).

**Example 7.5** - Evaluation of the Riemann and Hurwitz Zeta Functions using the derived lnd(n,Δx,x) function

Riemann Zeta Function

Problem Description

1) Evaluate ζ(i)

ζ(n) = lnd(n,1,1)

Evaluating using the lnd(n,Δx,x) computer computation program

n = i

ζ(i) = lnd(i,1,1) = .0033002236 - .4181554491i

ζ(i) = .0033002236 - .4181554491i

Comment - This answer and the ones to follow will, arbitrarily, be to 10 place accuracy.

Checking using a Zeta Function calculation program, LNDX

ζ(i) = .0033002236 - .4181554491i
good check

2) Evaluate ζ(1+i)

ζ(n) = lnd(n,1,1)

Evaluating using the lnd(n,Δx,x) computer computation program, LNDX

n = 1+i

ζ(1+i) = lnd(1+i,1,1) = .5821580597 - .9268485643i

ζ(1+i) = .5821580597 - .9268485643i

Checking using a Zeta Function calculation program

ζ(1+i) = .5821580597 - .9268485643i
good check

3) Evaluate $\zeta(-11.234)$

$\zeta(n) = \text{lnd}(n,1,1)$

Evaluating using the lnd(n,Δx,x) computer computation program, LNDX

$n = -11.234$

$\zeta(-11.234) = \text{lnd}(-11.234,1,1) = .02272911368$

$\zeta(-11.234) = .02272911368$

Checking using a Zeta Function calculation program

$\zeta(-11.234) = .02272911368$

good check

Hurwitz Zeta Function

Problem Description

Evaluate $\zeta(2,.25) + \zeta(2,.75)$

$\zeta(n,x) = \text{lnd}(n,1,x)$

Evaluating using the lnd(n,Δx,x) computer computation program, LNDX

$\zeta(2,.25) + \zeta(2,.75) = \text{lnd}(2,1,.25) + \text{lnd}(2,1,.75)$

$\zeta(2,.25) + \zeta(2,.75) = 17.19732915450711073927 + 2.54187964767160649839$

$\zeta(2,.25) + \zeta(2,.75) = 19.73920880217871723766$

$\zeta(2,.25) + \zeta(2,.75) = 2\pi^2$

Checking the above result using calculations obtained from the internet site, functions.wolfram.com

$\zeta(2,.25) + \zeta(2,.75) = 2\pi^2$

good check

**Example 7.6** Derivation of an interesting Zeta Function relationship

Problem Description

Find an interesting Zeta Function relationship

Consider the following series:

$$\frac{1}{x-1} = \frac{1}{x} + \frac{1}{x^2} + \frac{1}{x^3} + \frac{1}{x^4} + \frac{1}{x^5} + \ldots \quad , \qquad |x| > 1 \tag{1}$$

$$D_{\Delta x}\mathrm{lnd}(1,\Delta x,x\text{-}a) = +\frac{1}{x-a} \quad , \quad n = 1 \tag{2}$$

$$D_{\Delta x}\mathrm{lnd}(n,\Delta x,x\text{-}a) = -\frac{1}{(x-a)^n} \quad , \quad n \neq 1 \tag{3}$$

where

a = constant

$$\frac{1}{x-1} = +\left(\frac{1}{x}\right) - \left(-\frac{1}{x^2}\right) - \left(-\frac{1}{x^3}\right) - \left(-\frac{1}{x^4}\right) - \left(-\frac{1}{x^5}\right) - \ldots \tag{4}$$

From Eq 2 thru Eq 4

$$\frac{1}{x-1} = D_{\Delta x}\mathrm{lnd}(1,\Delta x,x) - D_{\Delta x}\mathrm{lnd}(2,\Delta x,x) - D_{\Delta x}\mathrm{lnd}(3,\Delta x,x) - D_{\Delta x}\mathrm{lnd}(4,\Delta x,x) - D_{\Delta x}\mathrm{lnd}(5,\Delta x,x) - \ldots$$

$$D_{\Delta x}\mathrm{lnd}(1,\Delta x,x\text{-}1) = D_{\Delta x}[\mathrm{lnd}(1,\Delta x,x) - \mathrm{lnd}(2,\Delta x,x) - \mathrm{lnd}(3,\Delta x,x) - \mathrm{lnd}(4,\Delta x,x) - \mathrm{lnd}(5,\Delta x,x) - \ldots] \tag{5}$$

Integrate both sides of the above equation with a constant of integration equal to zero

$$\mathrm{lnd}(1,\Delta x,x\text{-}1) = \mathrm{lnd}(1,\Delta x,x) - \mathrm{lnd}(2,\Delta x,x) - \mathrm{lnd}(3,\Delta x,x) - \mathrm{lnd}(4,\Delta x,x) - \mathrm{lnd}(5,\Delta x,x) - \ldots \tag{6}$$

Rearranging terms

$$\mathrm{lnd}(1,\Delta x,x) - \mathrm{lnd}(1,\Delta x,x\text{-}1) = \mathrm{lnd}(2,\Delta x,x) + \mathrm{lnd}(3,\Delta x,x) + \mathrm{lnd}(4,\Delta x,x) + \mathrm{lnd}(5,\Delta x,x) + \ldots$$

Then

$$\mathrm{lnd}(1,\Delta x,x) - \mathrm{lnd}(1,\Delta x,x\text{-}1) = \sum_{n=2}^{\infty} \mathrm{lnd}(n,\Delta x,x) \tag{7}$$

Use the above equation to find an interesting Zeta Function relationship

Let $\Delta x = 1$

$x = 2$

$$\text{lnd}(1,1,2) - \text{lnd}(1,1,1) = \sum_{n=2}^{\infty} \text{lnd}(n,1,2) \qquad 8)$$

$$\sum_{\Delta x \; x=x_i}^{\infty} \frac{1}{x^n} = \frac{1}{\Delta x}\,\text{lnd}(n,\Delta x,x_i) \;, \quad \text{for } \operatorname{Re}(n) > 1 \qquad 9)$$

For $n = 2,3,4,\ldots$

$\Delta x = 1$

$$\sum_{1 \; x=x_i}^{\infty} \frac{1}{x^n} = \text{lnd}(n,1,x_i) \qquad 10)$$

$$1 + \sum_{1 \; x=2}^{\infty} \frac{1}{x^n} = \sum_{1 \; x=1}^{\infty} \frac{1}{x^n} \qquad 11)$$

From Eq 10 and Eq 11

$$1 + \text{lnd}(n,1,2) = \text{lnd}(n,1,1) = \zeta(n) \qquad 12)$$

where

$$\zeta(n) = \sum_{1 \; x=1}^{\infty} \frac{1}{x^n}, \quad \text{The Riemann Zeta Function}$$

From Eq 12

$$\text{lnd}(n,\Delta x,2) = \zeta(n) - 1 \qquad 13)$$

Substituting Eq 13 into Eq 8

$$\text{lnd}(1,1,2) - \text{lnd}(1,1,1) = \sum_{n=2}^{\infty} \text{lnd}(n,1,2) = \sum_{n=2}^{\infty} [\zeta(n) - 1]$$

$$\text{lnd}(1,1,2) - \text{lnd}(1,1,1) = \sum_{n=2}^{\infty} [\zeta(n) - 1] \qquad 14)$$

Evaluating using the lnd(n,Δx,x) computer computation program, LNDX

$$\text{lnd}(1,1,2) - \text{lnd}(1,1,1) = 1 - 0 = 1 \qquad 15)$$

Thus

$$\sum_{n=2}^{\infty} [\zeta(n) - 1] = 1 \qquad 16)$$

Checking the above equation

| | | |
|---|---|---|
| 1 | $\zeta(2) - 1 =$ | .64493 |
| 2 | $\zeta(3) - 1 =$ | .20205 |
| 3 | $\zeta(4) - 1 =$ | .08232 |
| 4 | $\zeta(5) - 1 =$ | .03692 |
| 5 | $\zeta(6) - 1 =$ | .01734 |
| 6 | $\zeta(7) - 1 =$ | .00834 |
| 8 | $\zeta(8) - 1 =$ | .00407 |
| 8 | $\zeta(9) - 1 =$ | .00200 |
| 9 | $\zeta(10) - 1 =$ | .00099 |
| 10 | $\zeta(11) - 1 =$ | .00049 |
| | | .99945 (for 10 terms, 1 = exact value) |

**Example 7.7** Use of the lnd(n,Δx,x) function to calculate the Gamma Function

Problem Description

Find the value of $\Gamma(1+i)$

$$\Gamma(x) = \sqrt{2\pi}\, e^{\lim_{n\to 0}\left[ \frac{\text{lnd}(n,1,x) - \text{lnd}(0,1,x)}{n} \right]}$$

Evaluating using the lnd(n,Δx,x) computer computation program, LNDX

$x = 1+i$

let $n = 10^{-17}$

$$\Gamma(1+i) = \sqrt{2\pi}\, e^{\frac{\text{lnd}(10^{-17},1,1+i) - \text{lnd}(0,1,1+i)}{10^{-17}}}$$

$$\Gamma(x) = \sqrt{2\pi}\, e^{\frac{[-1.569861732506529091 - .301640320467533195i]\, 10^{-17}}{10^{-17}}}$$

$\Gamma(1+i) = .4980156681183560 - .1549498283018106i$

Checking using Gamma Function Tables

$\ln\Gamma(1+i) = -.650923199302 - .301640320468i$

$\Gamma(1+i) = e^{-.650923199302 - .301640320468i}$

$\Gamma(1+i) = .498015668118 - .154949828302i$

good check

**Example 7.8** Calculation of an alternating sign summation using the derived Alternating Sign Discrete Calculus Summation Equation

Problem Description

Evaluate the summation, $\sum_{\frac{1}{3}\, x=0}^{5} (-1)^{3x} \sin 2x$ , using the Alternating Sign Discrete Calculus Summation Equation.

$$\sum_{\Delta x\, x=x_1}^{x_2} (-1)^{\frac{x-x_1}{\Delta x}} f(x) = -\frac{1}{2} f(x) \Big|_{x_1}^{x_2+\Delta x} + \sum_{m=1}^{\infty} h_m (2\Delta x)^{2m-1} \frac{d^{2m-1}}{dx^{2m-1}} f(x)\Big|_{x_1}^{x_2+\Delta x}$$

$$h_m = \left[ \frac{B_{2m} + \frac{C_m}{(2m+1)!2^{2m}}}{(2m)!} \right], \quad m = 1,2,3,\ldots$$

where

$B_{2m}$ = Bernoulli Constants

$C_m$ = $C_m$ Constants

$f(x)$ = function of x

$\Delta x$ = x increment

$$x_2 = x_1 + (2p-1)\Delta x\ , \quad p = 1,2,3,\ldots$$

Let $f(x) = \sin 2x$

$\Delta x = \frac{1}{3}$

$x_1 = 0$

$x_2 = 5$

$x_2 + \Delta x = 5 + \frac{1}{3} = \frac{16}{3}$

Substituting

$$\sum_{\frac{1}{3}\, x=0}^{5} (-1)^{3x} \sin 2x = -\frac{1}{2} \sin 2x\Big|_{0}^{\frac{16}{3}} + \sum_{m=1}^{\infty} h_m \left(\frac{2}{3}\right)^{2m-1} \frac{d^{2m-1}}{dx^{2m-1}} f(x)\Big|_{0}^{\frac{16}{3}}$$

$$\sum_{\frac{1}{3}\, x=0}^{5} (-1)^{3x} \sin 2x = -\frac{1}{2} \sin 2x\Big|_{0}^{\frac{16}{3}} + \frac{1}{8}\left(\frac{2}{3}\right)^1 (2)^1 \cos 2x\Big|_{0}^{\frac{16}{3}} + \left(-\frac{1}{384}\right)\left(\frac{2}{3}\right)^3 (-2)^3 \cos 2x\Big|_{0}^{\frac{16}{3}}\Big) +$$

$$(\frac{1}{15360})(\frac{2}{3})^5(2)^5\cos 2x\Big|_0^{\frac{16}{3}} + (-\frac{17}{10321920})(\frac{2}{3})^7(-2)^7\cos 2x\Big|_0^{\frac{16}{3}} +$$

$$(\frac{31}{743178240})(\frac{2}{3})^9(2)^9\cos 2x\Big|_0^{\frac{16}{3}} + (-\frac{691}{653996851200})(\frac{2}{3})^{11}(-2)^{11}\cos 2x\Big|_0^{\frac{16}{3}} +$$

$$(\frac{5461}{204047017574400})(\frac{2}{3})^{13}(2)^{13}\cos 2x\Big|_0^{\frac{16}{3}} + \ldots$$

$$\sum_{\frac{1}{3}\ x=0}^{5}(-1)^{3x}\sin 2x =$$

| | Term Values | Partial Sum | Place Accuracy |
|---|---|---|---|
| 1 | +.4731978784190540 | .4731978784905400 | 0 |
| 2 | - .2205015663958965 | .2526963120231575 | 1 |
| 3 | - .0081667246813295 | .2445295873418280 | 3 |
| 4 | - .0003629655413924 | .2441666218004355 | 4 |
| 5 | - .0000163238470996 | .2441502979533359 | 4 |
| 6 | - .0000007349880320 | .2441495629653038 | 7 |
| 7 | - .0000000330972127 | .2441495298687910 | 8 |
| 8 | - .0000000014904191 | .2441495283776718 | 10 |

$$\sum_{\frac{1}{3}\ x=0}^{5}(-1)^{3x}\sin 2x = \underline{.2441495283}776718$$

(10 place accuracy for 8 series terms)

Checking

Using a computer program to exactly calculate, $\sum_{\frac{1}{3}\ x=0}^{5}(-1)^{3x}\sin 2x$

$$\sum_{\frac{1}{3}\ x=0}^{5}(-1)^{3x}\sin 2x = .2441495283073910169$$ good check

(exact value)

**Example 7.9** A calculation using the derived Discrete Calculus Summation Equation

Note – The Discrete Calculus Summation Equation is not the Euler-Maclaurin Sum Formula. The Discrete Calculus Summation Equation and the Euler-Maclaurin Sum Formula are both used to derive the Alternating Sign Discrete Calculus Summation Equation presented in the previous example.

Problem Description

Evaluate the summation, $\displaystyle\sum_{\frac{\pi}{4}\; x=0}^{\frac{\pi}{2}} \cos x$ , using the Discrete Calculus Summation Equation.

$$\sum_{\Delta x\; x=x_1}^{x_2} f(x) = \frac{1}{\Delta x}\int_{\Delta x\; x_1-\frac{\Delta x}{2}}^{x_2+\frac{\Delta x}{2}} f(x)\Delta x + \sum_{m=1}^{\infty} b_m\, \Delta x^{2m-1} \frac{d^{2m-1}}{dx^{2m-1}} f(x)\Big|_{x_1-\frac{\Delta x}{2}}^{x_2+\frac{\Delta x}{2}}$$

$$b_m = \frac{-C_m}{(2m+1)!2^{2m}}$$

where
$C_m = C_m$ Constants
$f(x)$ = function of x
$\Delta x$ = x increment

Let $f(x) = \cos x$

$\Delta x = \frac{\pi}{4}$

$x_1 = 0$

$x_2 = \frac{\pi}{2}$

Substituting

$$\sum_{\frac{\pi}{4}\; x=0}^{\frac{\pi}{2}} \cos x = \frac{1}{\frac{\pi}{4}}\int_{\frac{\pi}{4}\; 0-\frac{\pi}{8}}^{\frac{\pi}{2}+\frac{\pi}{8}} f(x)\Delta x + \sum_{m=1}^{\infty} b_m\, \Delta x^{2m-1} \frac{d^{2m-1}}{dx^{2m-1}} f(x)\Big|_{0-\frac{8}{2}}^{\frac{\pi}{2}+\frac{\pi}{8}}$$

$$\sum_{\frac{\pi}{4}\; x=0}^{\frac{\pi}{2}} \cos x = \frac{4}{\pi}\sin x\Big|_{-\frac{\pi}{8}}^{\frac{5\pi}{8}} + (-\frac{1}{24})(\frac{\pi}{4})^1(-\sin x)\Big|_{-\frac{\pi}{8}}^{\frac{5\pi}{8}} + (\frac{7}{5760})(\frac{\pi}{4})^3\sin x\Big|_{-\frac{\pi}{8}}^{\frac{5\pi}{8}} +$$

$$\left(-\frac{31}{967680}\right)\left(\frac{\pi}{4}\right)^5(-\sin x)\Big|_{-\frac{\pi}{8}}^{\frac{5\pi}{8}} + \left(\frac{127}{154828800}\right)\left(\frac{\pi}{4}\right)^7\sin x\Big|_{-\frac{\pi}{8}}^{\frac{5\pi}{8}} +$$

$$\left(-\frac{73}{3503554560}\right)\left(\frac{\pi}{4}\right)^9(-\sin x)\Big|_{-\frac{\pi}{8}}^{\frac{5\pi}{8}} + \left(\frac{1414477}{2678117105664000}\right)\left(\frac{\pi}{4}\right)^{11}\sin x\Big|_{-\frac{\pi}{8}}^{\frac{5\pi}{8}} +$$

$$\left(-\frac{8191}{612141052723200}\right)\left(\frac{\pi}{4}\right)^{13}(-\sin x)\Big|_{-\frac{\pi}{8}}^{\frac{5\pi}{8}} + \ldots$$

$${}_{\frac{\pi}{4}}\sum_{x=0}^{\frac{\pi}{2}} \cos x =$$

| | Term Values | Partial Sum | Place Accuracy |
|---|---|---|---|
| 1 | 1.66356703456770220055 | 1.66356703456770220055 | 1 |
| 2 | .04275717304070962036 | 1.70632480760773162595 | 3 |
| 3 | .00076926424024611682 | 1.70709407184797774277 | 4 |
| 4 | .00001250862806468674 | 1.70710658047604242952 | 7 |
| 5 | .00000019756567220325 | 1.70710677804171463278 | 8 |
| 6 | .00000000309565996843 | 1.70710678113737460121 | 11 |
| 7 | .00000000004840445842 | 1.70710678118577905964 | 12 |
| 8 | .00000000000075645695 | 1.70710678118653551659 | 14 |

$${}_{\frac{\pi}{4}}\sum_{x=0}^{\frac{\pi}{2}} \cos x = \underline{1.70710678118653}551659$$

(14 place accuracy for 8 series terms)

Checking

$${}_{\frac{\pi}{4}}\sum_{x=0}^{\frac{\pi}{2}} \cos x = 1 + \frac{\sqrt{2}}{2} + 0 = 1.70710678118654752440 \quad \text{good check}$$

(exact value)

**Example 7.10** Finding the Z Transform of a function using Interval Calculus

Problem Description .

Find the Z Transform of f(x) = x using the derived interval calculus relationship,

$$Z[f(x)] = \frac{1}{T}\ {}_T\!\int_0^{\infty} z^{-(\frac{x}{\Delta x})} f(x)\,\Delta x\ .$$

$$Z[x] = \frac{1}{T}\ {}_T\!\int_0^{\infty} z^{-(\frac{x}{\Delta x})}\, x\,\Delta x$$

Using integration by parts and the integration table

$${}_{\Delta x}\!\int_{x_1}^{x_2} v(x) D_{\Delta x} u(x)\Delta x = u(x)v(x)\Big|_{x1}^{x2} - {}_{\Delta x}\!\int_{x_1}^{x_2} D_{\Delta x} v(x) u(x+\Delta x)\Delta x$$

$${}_{\Delta x}\!\int z^{-\frac{x}{\Delta x}}\,\Delta x = -\left(\frac{z\Delta x}{z-1}\right) z^{-\frac{x}{\Delta x}} + k$$

Let $D_{\Delta x}u(x) = z^{-\frac{x}{\Delta x}}$ $\qquad v(x) = x$

$$u(x) = -\left(\frac{z\Delta x}{z-1}\right) z^{-\frac{x}{\Delta x}} \qquad D_{\Delta x}v(x) = 1$$

$$u(x+\Delta x) = -\left(\frac{\Delta x}{z-1}\right) z^{-\frac{x}{\Delta x}}$$

$$Z[x] = \frac{1}{T}\ {}_T\!\int_0^{\infty} z^{-(\frac{x}{\Delta x})}\, x\,\Delta x = -\frac{1}{T}\, x\left(\frac{z\Delta x}{z-1}\right) z^{-\frac{x}{\Delta x}}\Bigg|_{0}^{\infty\ T} + \frac{1}{T}\left(\frac{\Delta x}{z-1}\right)\ {}_T\!\int_0^{\infty} z^{-(\frac{x}{\Delta x})}\,\Delta x$$

$$Z[x] = 0 - \frac{1}{T}\left(\frac{\Delta x}{z-1}\right)\left(\frac{z\Delta x}{z-1}\right) z^{-\frac{x}{\Delta x}}\Bigg|_{0}^{\infty\ T} = \frac{1}{T}\,\frac{T^2 z}{(z-1)^2}\ , \quad \Delta x = T$$

$$Z[x] = \frac{Tz}{(z-1)^2}$$

This is the Z Transform of the function, x

Good check

**Example 7.11** Use the complex plane area calculation equation, $A_c = \left| \frac{j}{2} \sum_{p=0}^{N} (-1)^p c_p^{\,2} \right|$, to calculate the area within a discrete closed contour in the complex plane

Problem Description

Find the area enclosed within the complex plane discrete closed contour shown in Diagram 6.1 below.

Use the complex plane area calculation equation, $A_c = \left| \frac{j}{2} \sum_{p=0}^{N} (-1)^p c_p^{\,2} \right|$, to calculate the area.

Diagram 7.1

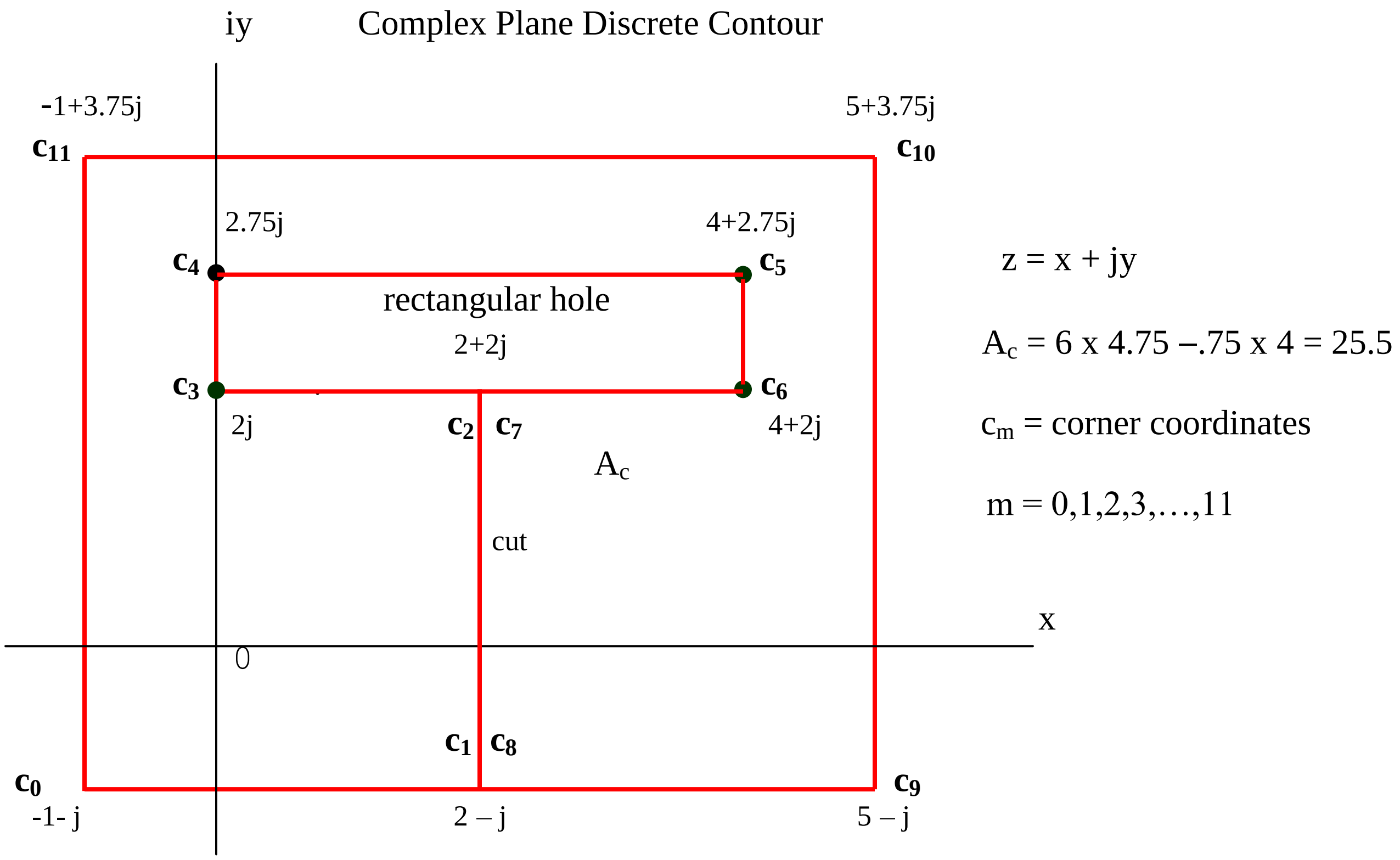

$A_c$, the area enclosed within the complex plane discrete closed contour of Diagram 7.1, is calculated using the following equation.

$$A_c = \left| \frac{j}{2} \sum_{p=0}^{N} (-1)^p c_p^{\,2} \right|$$

where

$A_c$ = The area enclosed within a discrete complex plane closed contour

$N = 3,5,7,9,11,\ldots$

N = the number of the discrete complex plane closed contour corner points minus one
$c_p$ = the coordinates of the corner points of the discrete complex plane closed contour
p = 0,1, 2, 3, 4,…, N

Note - The integration initial point is designated as $c_0$. The following points, progressing along the contour in the direction of integration, are successively designated $c_1, c_2, c_3, \ldots, c_N$

The following UBASIC program uses the above equation to calculate the area enclosed within the complex plane discrete closed contour shown in Diagram 7.1.

```
10   dim C(20)
20   C(0)=-1-#i
30   C(1)=2-#i
40   C(2)=2+2#i
50   C(3)=2#i
60   C(4)=2.75#i
70   C(5)=4+2.75#i
80   C(6)=4+2#i
90   C(7)=2+2#i
100  C(8)=2-#i
110  C(9)=5-#i
120  C(10)=5+3.75#i
130  C(11)=-1+3.75#i
140  S=0
150  for K=0 to 11
160  S=S+((-1)^K)*C(K)^2
170  next K
180  W=abs((#i/2)*S)
190  print "The closed contour enclosed area is ";W
```

The evaluation obtained by running the above program is as follows:

**The closed contour enclosed area is 25.5**

Checking the above result

Referring to the contour diagram, Diagram 7.1

$A_c = 6 \text{ x } 4.75 – .75 \text{ x } 4 = 25.5$      good check

**Example 7.12** Use of the $K_{\Delta x}$ Transform to solve difference equations

Problem Description

Find y(x) where $y(x+2\Delta x) - y(x) = 0$ , $\Delta x = 1$, $y(0) = 1$, $y(1) = -1$, $x = m\Delta x$ , $m = 0,1,2,3,\ldots$
Solve for y(x) using the $K_{\Delta x}$ Transform.

From the table of $K_{\Delta x}$ Transform general equations (TABLE 2)

$$K_{\Delta x}[y(x+2\Delta x)] = (1+s\Delta x)^2 K_{\Delta x}[y(x)] - (1+s\Delta x)y(0)\Delta x - y(\Delta x)\Delta x \qquad 1)$$

$$y(x+2\Delta x) - y(x) = \quad 0 \qquad 2)$$

Taking the $K_{\Delta x}$ Transform of Eq 2

$$K_{\Delta x}[y(x)+2\Delta x] - K_{\Delta x}[y(x)] = 0 \qquad 3)$$

$$y(s) = K_{\Delta x}[y(x)] \qquad 4)$$

Substituting Eq 1 and Eq 4 into Eq 3 and introducing the y(x) initial conditions

$$(1+s)^2 y(s) - (1+s)(1)(1) - (-1)(1) - y(s) = 0$$

$$[(1+s)^2 - 1]y(s) - 1 - s + 1 = 0$$

$$(s^2+2s)y(s) = s$$

$$y(s) = \frac{1}{(s+2)}$$

From the tables of $K_{\Delta x}$ Transforms (TABLE 3) and Interval Calculus Functions (TABLE 4)

$$K_{\Delta x}[e_{\Delta x}(a,x)] = K_{\Delta x}[(1+a\Delta x)^{\frac{x}{\Delta x}}] = \frac{1}{s - a}$$

$$a = -2, \ \Delta x = 1$$

Substituting

$$K_1[e_1(-2,x)] = K_1[(1-2(1))^{\frac{x}{1}}] = \frac{1}{s + 2} = y(s)$$

$$y(s) = K_1[e_1(-2,x)] = K_1[(1-2(1))^{\frac{x}{1}}]$$

Taking the inverse $K_{\Delta x}$ Transform

$$y(x) = e_1(-2,x) = (1-2(1))^{\frac{x}{1}} = (-1)^x$$

$$y(x) = (-1)^x , \ x = 0,1,2,3,\ldots \qquad 5)$$

Checking Eq 5 using the given difference equation and the y(x) initial conditions

$y(0) = 1, \; y(1) = -1, \; \Delta x = 1$

$y(x+2) - y(x) = 0$

$y(x) = (-1)^x$

$(-1)^{x+2} - (-1)^x = 0$

$(-1)^x(-1)^2 - (-1)^x = 0 , \quad x = 0,1,2,3,\ldots$

good check

$y(0) = (-1)^0 = 1 , \; y(1) = (-1)^1 = -1$

good check

<u>Comment</u> - The commonly used transform for solving difference equations is the Z Transform. As can be seen in this example, the $K_{\Delta x}$ Transform can be used also. The solution, of course, will be the same using either transform. However, an advantage of the $K_{\Delta x}$ Transform over the Z Transform may be its similarity to the Laplace Transform. It has the look and feel of the Laplace Transform which is familiar to most engineers. In fact, the $K_{\Delta x}$ Transform becomes the Laplace Transform when $\Delta x \rightarrow 0$.

**Example 7.13** Solving servomechanism control problems using the $K_{\Delta x}$ Transform and Nyquist methods

Problem Description

Investigate the stability of the specified discrete variable control system using the $K_{\Delta x}$ Transform and Nyquist Criteria methods. Determine the system stability phase margin for the gain constant, K = 1. Find the value of the gain constant, K, at which the system is oscillatory. The x increment, Δx, is .1 and the initial values of y(x) are y(0)=y(.1)=y(.2)=0.

Note – Often the variable, x, represents time and, if so, a variable, t, may be used instead.

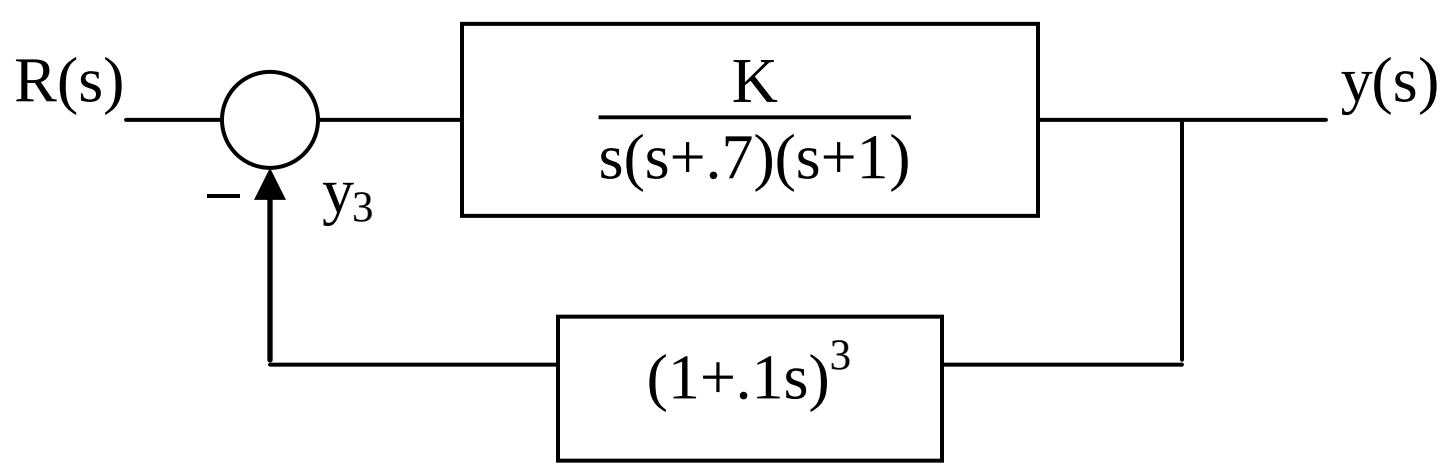

From the specified discrete feedback system diagram, write the transfer function for $\frac{y(s)}{R(s)}$.

$$\frac{y(s)}{R(s)} = \frac{\frac{K}{s(s+.7)(s+1)}}{1 + \frac{K(1+.1s)^3}{s(s+.7)(s+1)}}$$

$$A(s) = \frac{K(1+.1s)^3}{s(s+.7)(s+1)} \qquad 1)$$

From the specified discrete feedback system diagram and the following $K_{\Delta x}$ Transforms, write the equation for y(x). Let R be a step input where R = 5.

$$K_{\Delta x}[D^n{}_{\Delta x}f(x)] = s^n K_{\Delta x}[f(x)] - s^{n-1}D^0{}_{\Delta x}f(0) - s^{n-2}D^1{}_{\Delta x}f(0) - s^{n-3}D^2{}_{\Delta x}f(0) - \ldots - s^0 D^{n-1}{}_{\Delta x}f(0), \quad n = 1,2,3,\ldots \qquad 2)$$

$$K_{\Delta x}[f(x+n\Delta x)] = (1+s\Delta x)^n K_{\Delta x}[f(x)] - \Delta x \sum_{m=1}^{n}(1+s\Delta x)^{n-m}f([n-m]\Delta x), \quad n = 1,2,3,\ldots \qquad 3)$$

$$[s(s+.7)(s+1)]y(s) = K[\frac{R}{s} - (1+.1s)^3 y(s)]$$

$$[s^3 + 1.7s^2 + .7s]y(s) = K[\frac{5}{s} - (1+.1s)^3 y(s)] \qquad 4)$$

The initial values are y(0)=y(.1)=y(.2)=0

y =y(x)

x = nΔx

$y_n = y(n\Delta x)$

$D_{\Delta x}y(x) = \dfrac{y(x+\Delta x) - y(x)}{\Delta x} = \dfrac{y([n+1]\Delta x) - y(n\Delta x)}{\Delta x} = \dfrac{y_{n+1} - y_n}{\Delta x}$, The discrete derivative

From Eq 2, 3, and Eq 4

$\Delta x = .1$

$D_{\Delta x}^3 y(x) + 1.7D_{\Delta x}^2 y(x) + .7D_{\Delta x}^1 y(x) = K[5 - y(x+3\Delta x)]$

$$\frac{y_{n+3}-3y_{n+2}+3y_{n+1}-y_n}{.1^3} + 1.7\left(\frac{y_{n+2}-2y_{n+1}+y_n}{.1^2}\right) + .7\left(\frac{y_{n+1}-y_n}{.1}\right) = 5K-Ky_{n+3}$$

$$y_{n+3}-3y_{n+2}+3y_{n+1}-y_n + .17(y_{n+2}-2y_{n+1}+y_n) + .007(y_{n+1}-y_n) = .005K-.001Ky_{n+3} \qquad 5)$$

Combining the terms of Eq 5 and solving for $y_3$

$$y_{n+3} = \frac{2.83y_{n+2} - 2.667y_{n+1} + .837y_n + .005K}{1+.001K} \qquad 6)$$

A program to plot the Modified Nyquist Diagram of A(s) and to calculate the phase margin of A(s) was written. Also, a program to plot y(x) from a y(x) difference equation was written. These two programs appear in the Calculation Program section of the Appendix as programs 6 and 7.

From the A(s) equation, Eq 1, and the y(x) difference equation, Eq 6, the following plots were obtained.

The left half s plane Critical Circle has a radius of $\dfrac{1}{\Delta x} = \dfrac{1}{.1} = 10$ with a center at $s = -\dfrac{1}{\Delta x} = -\dfrac{1}{.1} = -10$ .

1) <u>For K=1, Δx =.1</u>

From Eq 1

$$A(s) = \frac{1(1+.1s)^3}{s(s+.7)(s+1)}$$

Evaluating A(s) around the left half s plane Critical Circle of radius 10 and center at s = –10

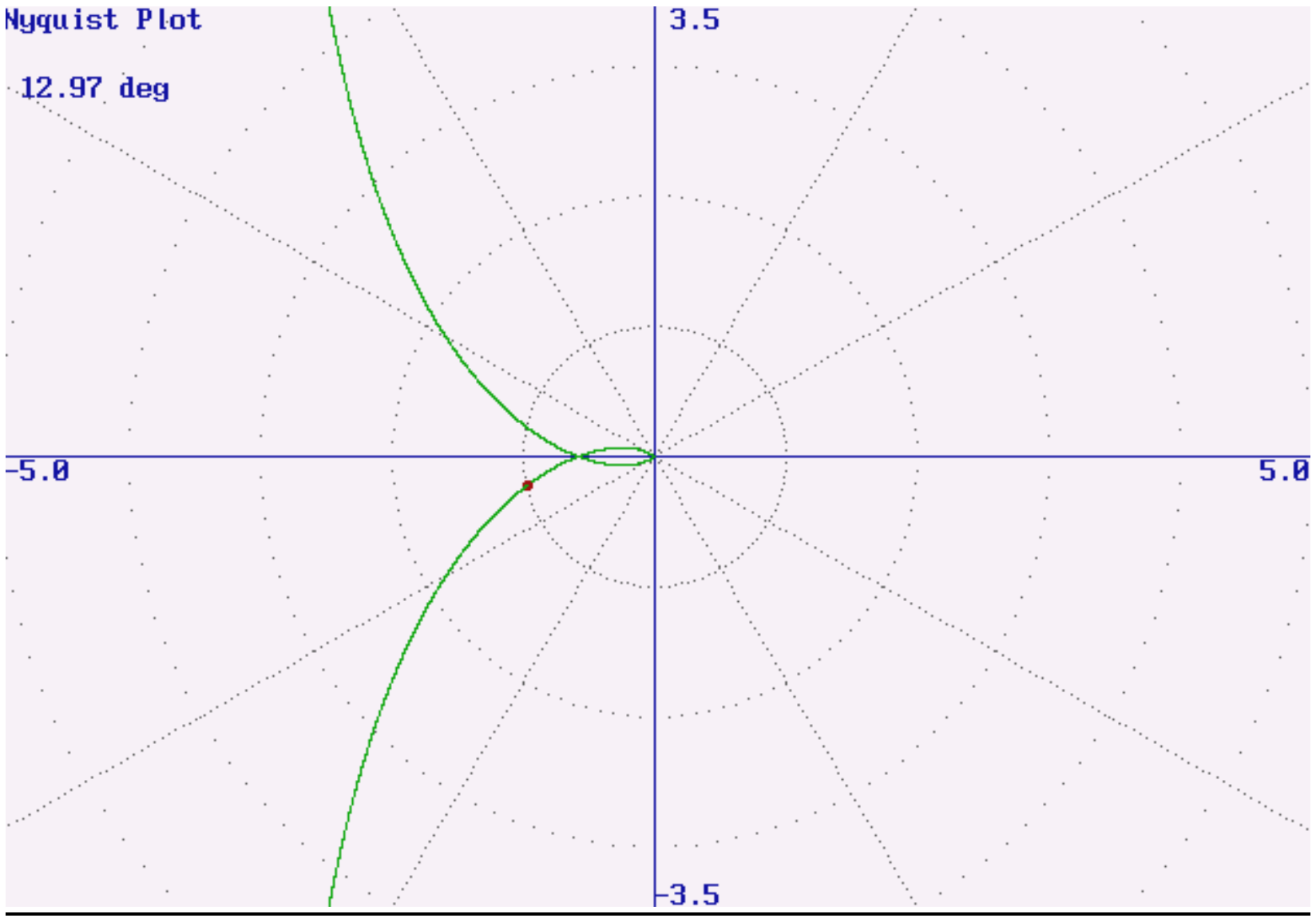

For K=1 and Δx = .1 a phase margin of 12.97 degrees was calculated.

Plotting the corresponding y(x) from the difference equation, Eq 6

$$y_{n+3} = \frac{2.83y_{n+2} - 2.667y_{n+1} + .837y_n + .005}{1.001}, \quad K=1 \text{ and } \Delta x = .1$$

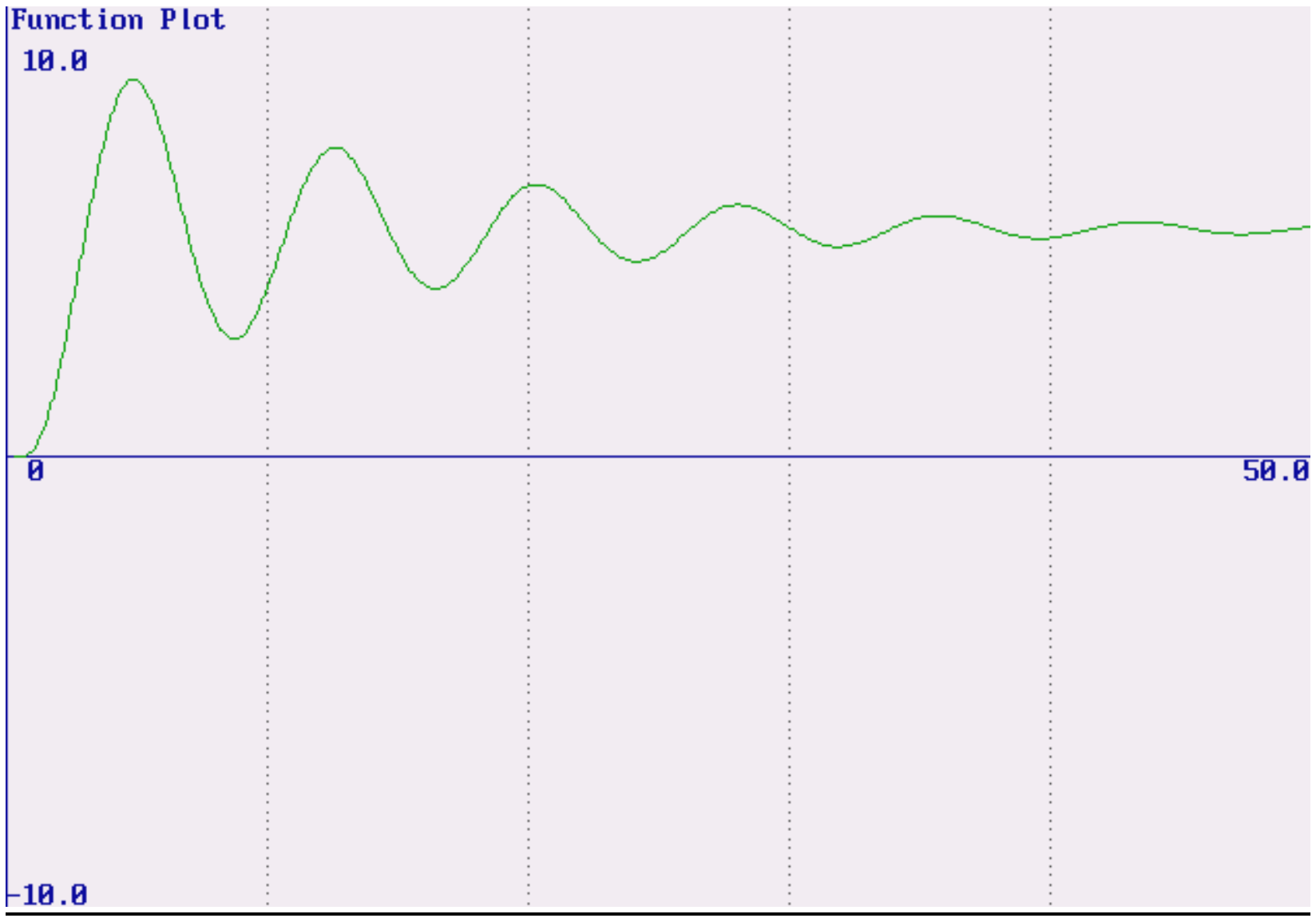

This y(x) plot displays stability but with a significant damped oscillation resulting from a step input of 5. This is what one would expect from a calculated phase margin of 12.97 degrees.

2) <u>For K=1.7, Δx =.1</u>

From Eq 1

$$A(s) = \frac{1.7(1+.1s)^3}{s(s+.7)(s+1)}$$

Evaluating A(s) around the left half s plane Critical Circle of radius 10 and center at s = –10

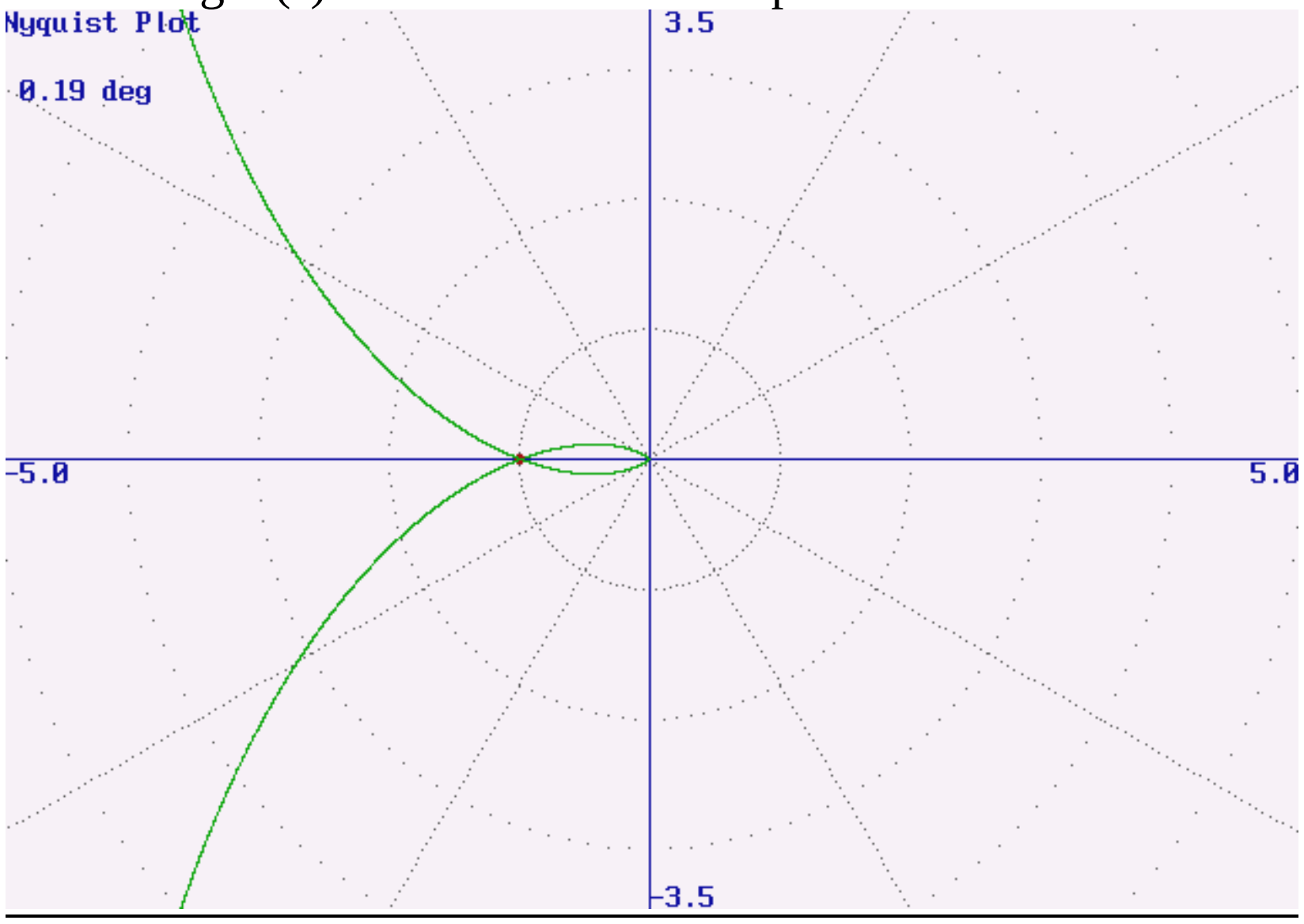

For K=1.7 and Δx = .1, the phase margin is calculated to be .19 degrees (virtually 0 degrees). This would indicate that the control system is oscillatory with a gain constant of K=1.7 and is on the verge of being unstable.

Plotting the corresponding y(x) from the difference equation, Eq 6

$$y_{n+3} = \frac{2.83y_{n+2} - 2.667y_{n+1} + .837y_n + .0085}{1.0017}, \qquad K = 1.7 \text{ and } \Delta x = .1$$

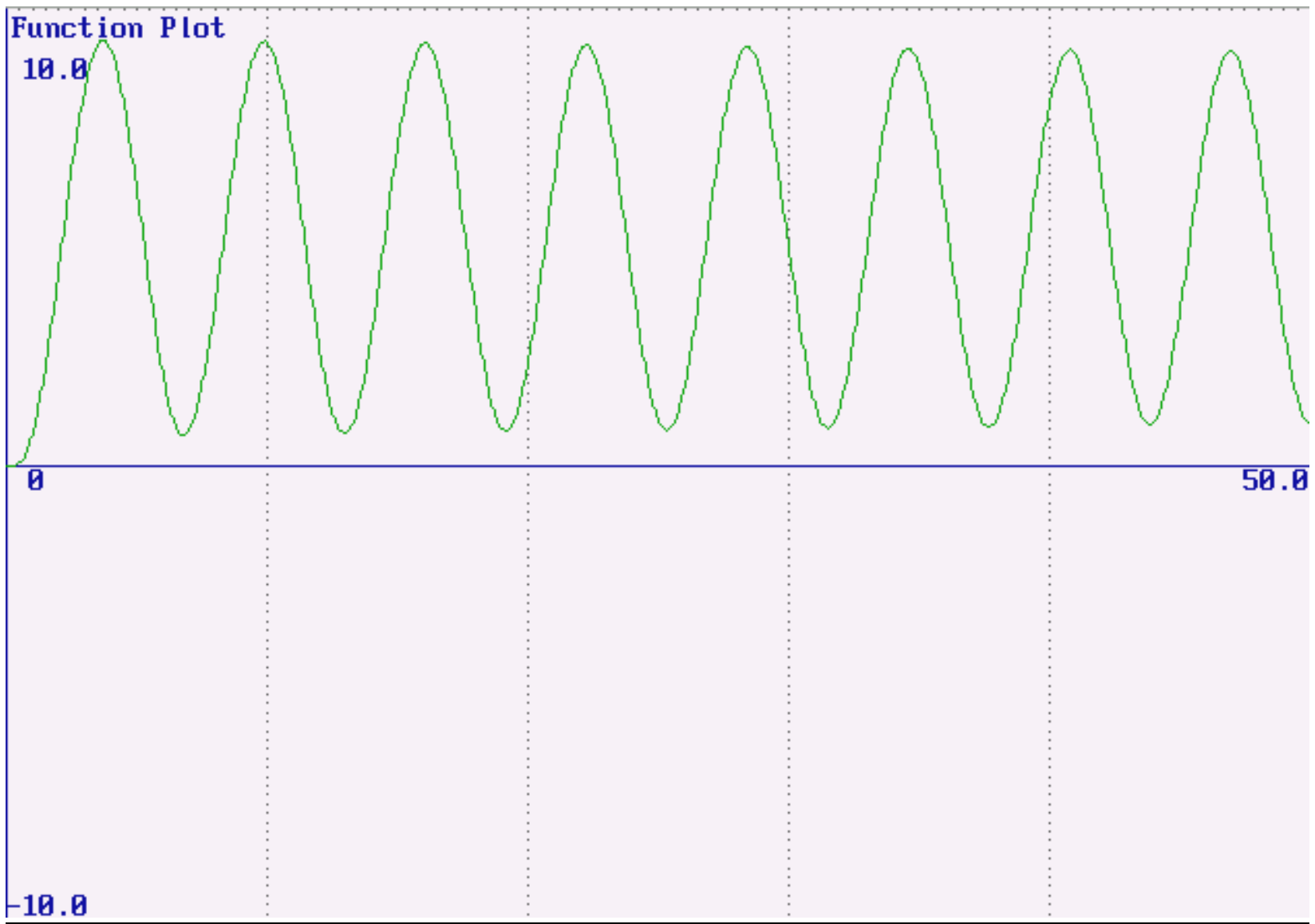

This y(x) plot displays a virtually sustained oscillation for a gain constant of K=1.7 resulting from a step input of 5. This is what one would expect from a calculated phase margin of .19 degrees.

**Example 7.14** Evaluation of a summation which has no division by zero term

Problem Description

Evaluate the summation, $\displaystyle {}_{.5}\sum_{x=-.9}^{1.6}\frac{5x-2}{x-.5}$ where $\Delta x = .5$, $x = -.9,-.4,.1,.6,\ldots,1.6$

The summation equations to be used, when no division by zero term is involved, are as follows:

$$ {}_{\Delta x}\sum_{x=x_1}^{x_2} f(x) = \frac{1}{\Delta x}\ {}_{\Delta x}\int_{x_1}^{x_2+\Delta x} f(x)\Delta x \qquad 1)$$

and

$$ {}_{\Delta x}\int_{x_1}^{x_2}\frac{1}{(x-a)^n}\Delta x = \pm\, \mathrm{lnd}(n,\Delta x,x-a)\Big|_{x_1}^{x_2}, \qquad + \text{ for } n = 1,\ - \text{ for } n \neq 1 \qquad 2)$$

where

$x = x_1, x_1+\Delta x, x_1+2\Delta x, x_1+3\Delta x, \ldots, x_2$

$$ {}_{.5}\sum_{x=-.9}^{1.6}\frac{5x-2}{x-.5} = \frac{1}{.5}\ {}_{.5}\int_{-.9}^{2.1}\frac{5x-2}{x-.5}\Delta x, \qquad \Delta x = .5$$

$$ {}_{.5}\sum_{x=-.9}^{1.6}\frac{5x-2}{x-.5} = 2\ {}_{.5}\int_{-.9}^{2.1}\frac{5x-2}{x-.5}\Delta x$$

$$\frac{5x+2}{x-.5} = \frac{.5}{x-.5} + 5$$

Substituting

$$ {}_{.5}\sum_{x=-.9}^{1.6}\frac{5x-2}{x-.5} = 2\ {}_{.5}\int_{-.9}^{2.1}\frac{5x-2}{x-.5}\Delta x = 2\ {}_{.5}\int_{-.9}^{2.1}\frac{.5}{x-.5}\Delta x + 2\ {}_{.5}\int_{-.9}^{2.1}5\Delta x \qquad 3)$$

$$\sum_{x=-.9}^{1.6} {}_{.5}\frac{5x-2}{x-.5} = {}_{.5}\int_{-.9}^{2.1} \frac{1}{x-.5}\,\Delta x + 10\,{}_{.5}\int_{-.9}^{2.1} \Delta x \qquad 4)$$

Integrating Eq 4

Eq 2 is used to integrate the integral, ${}_{.5}\int_{-.9}^{2.1} \frac{1}{x-.5}\,\Delta x$ ,   n = 1

$$\sum_{x=-.9}^{1.6} {}_{.5}\frac{5x-2}{x-.5} = \text{lnd}(1,.5,x-.5)\Big|_{-.9}^{2.1} + 10x\Big|_{-.9}^{2.1}$$

$$\sum_{x=-.9}^{1.6} {}_{.5}\frac{5x-2}{x-.5} = [\text{lnd}(1,.5,1.6) - \text{lnd}(1,.5,-1.4)] + 30 \qquad 5)$$

Using the lnd(n,Δx,x) calculation program, LNDX, to calculate Eq 5

$$\sum_{x=-.9}^{1.6} {}_{.5}\frac{5x-2}{x-.5} = [1.576054556188 + 2.549125818992] + 30 = 4.125180375180 + 30$$

$$\sum_{x=-.9}^{1.6} {}_{.5}\frac{5x-2}{x-.5} = 34.125180375180 \qquad 6)$$

Checking

$$\sum_{x=-.9}^{1.6} {}_{.5}\frac{5x-2}{x-.5} = \frac{-4.5-2}{-.9-.5} + \frac{-2-2}{-.4-.5} + \frac{.5-2}{.1-.5} + \frac{3-2}{.6-.5} + \frac{5.5-2}{1.1-.5} + \frac{8-2}{1.6-.5}$$

$$\sum_{x=-.9}^{1.6}{}_{.5}\ \frac{5x-2}{x-.5} = \frac{-6.5}{-1.4} + \frac{-4}{-.9} + \frac{-1.5}{-.4} + \frac{1}{.1} + \frac{3.5}{.6} + \frac{6}{1.1}$$

$$\sum_{x=-.9}^{1.6}{}_{.5}\ \frac{5x-2}{x-.5} = 34.125180375180$$

Good check

Note - If it is desired to find the area under the "sample and hold" shaped curve of the discrete function, $f(x) = \frac{5x-2}{x-.5}$, where $-.9 \le x \le 2.1$ and $\Delta x = .5$, this area can be easily calculated.

$$A = \Delta x \sum_{x=x_1}^{x_2}{}_{\Delta x} f(x) = (.5) \sum_{x=-.9}^{1.6}{}_{.5}\ \frac{5x-2}{x-.5} = (.5)(\ 34.125180375180) = 17.062590187590$$

**Example 7.15** Evaluation of a summation which has a division by zero term

Problem Description

Evaluate the Hurwitz Zeta function, $\zeta(1.2,-5)$

The Hurwitz Zeta Function is a generalization of the Riemann Zeta Function and is defined as follows:

$$\zeta(s,a) = \sum_{k=0}^{\infty} \frac{1}{(k+a)^s} \text{ , } \quad \text{The definition of the Hurwitz Zeta Function} \qquad 1)$$

where

Any term with k+a = 0 is excluded

Changing the Hurwitz Zeta Function of Eq 1 to a more convenient equivalent form

$$\zeta(s,a) = {}_1\!\sum_{x=a}^{\infty} \frac{1}{x^s} \qquad 2)$$

where

Any term with x = 0 is excluded

Use Eq 3 and Eq 4 to evaluate $\zeta(n,x)$

$$\zeta(n,x) = \zeta(n,1,x) \qquad 3)$$

$$\zeta(n,\Delta x,x_i) = \frac{1}{\Delta x}\,\mathrm{lnd}(n,\Delta x,x_i) = {}_{\Delta x}\!\sum_{x=x_i}^{\pm\infty} \frac{1}{x^n} \text{ , } \quad \mathrm{Re}(n) > 1 \text{ , } \quad \text{A form of the General Zeta Function} \qquad 4)$$

where

$+\infty$ for $\mathrm{Re}(\Delta x)>0$ or $\{\mathrm{Re}(\Delta x)=0 \text{ and } \mathrm{Im}(\Delta x)>0\}$

$-\infty$ for $\mathrm{Re}(\Delta x)<0$ or $\{\mathrm{Re}(\Delta x)=0 \text{ and } \mathrm{Im}(\Delta x)<0\}$

$x = x_i, x_i+\Delta x, x_i+2\Delta x, x_i+3\Delta x, \ldots$

$n,x,x_i,\Delta x$ = real or complex values

$\Delta x$ = x increment

Any term with x = 0 is excluded

Comment – The $\mathrm{lnd}(n,\Delta x,x)$ function excludes any summation term where x = 0.

For the Hurwitz Zeta Function

Let

$\Delta x = 1$

$n = s$

$x_i = a$

Substituting into Eq 3 and Eq 4

$$\zeta(s,a) = \zeta(s,1,a) = \sum_{x=a}^{\infty} {}_1 \frac{1}{x^s} = \text{lnd}(s,1,a) \qquad 5)$$

where

$x = a, a+1, a+2, a+3, \ldots$

$Re(s) > 1$

Any term with $x = 0$ is excluded

Use Eq 5 to evaluate the Hurwitz Zeta function, $\zeta(1.2,-5)$

$$\zeta(1.2,-5) = \sum_{x=-5}^{\infty} {}_1 \frac{1}{x^{1.2}} = \text{lnd}(1.2,1,-5) \qquad 6)$$

where

$x = -5, -4, -3, -2, -1, 0, 1, 2, 3, 4, 5, 6 \ldots$

Any term with $x = 0$ is excluded

Using the lnd$(n,\Delta x,x)$ calculation program to calculate Eq 6

$$\zeta(1.2,-5) = 3.9433912875609658 + 1.1974809673848110i \qquad 7)$$

Checking the above result

$$\zeta(1.2,-5) = \sum_{x=-5}^{\infty} {}_1 \frac{1}{x^{1.2}} = \sum_{x=-5}^{0} {}_1 \frac{1}{x^{1.2}} + \sum_{x=1}^{\infty} {}_1 \frac{1}{x^{1.2}} = \sum_{x=-5}^{0} {}_1 \frac{1}{x^{1.2}} + \zeta(1.2) \qquad 8)$$

where

$\zeta(1.2)$ = Riemann Zeta Function of 1.2

Find the Riemann Zeta Function of 1.2

Again using the lnd(n,Δx,x) calculation program, LNDX

$$\zeta(1.2) = \sum_{x=1}^{\infty} \frac{1}{x^{1.2}} = \text{lnd}(1.2,1,1) = 5.5915824411777507 \qquad 9)$$

From Eq 8 and Eq 9

$$\zeta(1.2,-5) = \left[\frac{1}{(-5)^{1.2}} + \frac{1}{(-4)^{1.2}} + \frac{1}{(-3)^{1.2}} + \frac{1}{(-2)^{1.2}} + \frac{1}{(-1)^{1.2}}\right] + 5.5915824411777507 \qquad 10)$$

Note that the pole at x = 0 is excluded

Using a computer to calculate the above summation

$\zeta(1.2,-5) = -1.6481911536167848+1.1974809673848110i + 5.5915824411777507$

$$\zeta(1.2,-5) = 3.9433912875609658 + 1.1974809673848110i \qquad 11)$$

Good check

**Example 7.16** A demonstration of discrete differentiation

Problem Description

Find the discrete derivative of the following function, $A(x) = \dfrac{(x+7)(x+3)(x-1)(x-5)}{(x-1)(x-5)}$. Find the discrete derivative four different ways.

Some equations which will be used:

1. $[x]_{\Delta x}^{n} = \prod_{m=1}^{n}(x-(m-1)\Delta x)$ , n = 1,2,3… , polynomial discrete function
2. $[x]_{\Delta x}^{0} = 1$
3. $x = x + a$ , a = constant
4. $D_{\Delta x}[x]_{\Delta x}^{n} = n[x]_{\Delta x}^{n-1}$ for n = 1,2,3… , discrete derivative of the discrete function, $[x]_{\Delta x}^{n}$
5. $D_{\Delta x}[u(x)v(x)] = v(x)D_{\Delta x}u(x) + D_{\Delta x}v(x)u(x+\Delta x)$ , discrete derivative of the product of two functions
6. $D_{\Delta x}\left[\dfrac{u(x)}{v(x)}\right] = \dfrac{v(x)D_{\Delta x}u(x) - u(x)D_{\Delta x}v(x)}{v(x)v(x+\Delta x)}$ , discrete derivative of the division of two functions

**1)** Find $D_{\Delta x}A(x)$ by first canceling and then using the definition of a discrete derivative,

$$D_{\Delta x}f(x) = \frac{f(x+\Delta x) - f(x)}{\Delta x} .$$

$A(x) = (x+7)(x+3)$ , $\Delta x = 4$

$$D_4A(x) = \frac{(x+11)(x+7) - (x+7)(x+3)}{4}$$

$$D_4A(x) = \frac{(x+7)(x+11-x-3)}{4} = 2(x+7)$$

$\underline{D_4A(x) = 2(x+7)}$

**2)** Find $D_{\Delta x}A(x)$ by first canceling and then using the formula for taking the discrete derivative, $D_{\Delta x}[x]^{n}_{\Delta x}$.

$A(x) = (x+7)(x+3) = [x+7]^{2}_{4}$ , $\Delta x = 4$, product terms = 2

$D_4A(x) = D_4[x+7]^{2}_{4} = 2\,[x+7]^{2-1}_{4} = 2\,[x+7]^{1}_{4} = 2(x+7)$

$\underline{D_4A(x) = 2(x+7)}$

**3)** Find $D_{\Delta x}A(x)$ by first canceling and then using the formula for taking the discrete derivative of a product of two functions, $D_{\Delta x}[u(x)v(x)] = v(x)D_{\Delta x}u(x) + D_{\Delta x}v(x)u(x+\Delta x)$.

$A(x) = (x+7)(x+3)$ , $\Delta x = 4$

$u(x) = (x+7) = [x+7]^{1}_{4}$ , product terms = 1

$v(x) = (x+3) = [x+3]^{1}_{4}$ , product terms = 1

Substituting into the derivative of the product of two functions equation

$D_4A(x) = D_4[\,[x+7]^{1}_{4}\,[x+3]^{1}_{4}\,] = [x+3]^{1}_{4}\,D_4[x+7]^{1}_{4} + D_4[x+3]^{1}_{4}\,[x+7+4]^{1}_{4}$

$D_4A(x) = (x+3)(1) + (1)(x+11) = 2x+14 = 2(x+7)$

$\underline{D_4A(x) = 2(x+7)}$

**4)** Find $D_{\Delta x}A(x)$ using the equation for taking the discrete derivative of a division of two functions, $D_{\Delta x}\left[\frac{u(x)}{v(x)}\right] = \frac{v(x)D_{\Delta x}u(x) - u(x)D_{\Delta x}v(x)}{v(x)v(x+\Delta x)}$.

$A(x) = \frac{(x+7)(x+3)(x-1)(x-5)}{(x-1)(x-5)} = \frac{[x+7]^{4}_{4}}{[x-1]^{2}_{4}} = \frac{u(x)}{v(x)}$ , $\Delta x = 4$ , product terms = 4,2

$$D_4A(x) = D_{\Delta x}\left[ \frac{[x+7]_4^4}{[x-1]_4^2} \right] = \frac{[x-1]_4^2\, D_4[x+7]_4^4 - [x+7]_4^4\, D_4[x-1]_4^2}{[x-1]_4^2\, [x-1+4]_4^2}$$

$$D_4A(x) = \frac{[x-1]_4^2\, 4[x+7]_4^3 - [x+7]_4^4\, 2[x-1]_4^1}{[x-1]_4^2\, [x+3]_4^2}$$

$$D_4A(x) = \frac{(x-1)(x-5)(4)(x+7)(x+3)(x-1) - (x+7)(x+3)(x-1)(x-5)(2)(x-1)}{(x-1)(x-5)(x+3)(x-1)}$$

Canceling terms

$$D_4A(x) = 4(x+7) - 2(x+7) = 2(x+7)$$

$\underline{D_4A(x) = 2(x+7)}$

Good checks

**Example 7.17** The solution of a differential difference equation using four different methods

Problem Description

Solve the differential difference equation, $D_{\Delta x}y(x) + y(x) = xe^{-x}$, where $\Delta x = 2$ and $y(0) = 3$ using four different methods.

**Solution 1 Use of the Method of Undetermined Coefficients to solve a differential difference equation**

$$D_{\Delta x}y(x) + y(x) = xe^{-x}, \qquad \Delta x = 2, \; y(0) = 3. \qquad 1)$$

Find the related homogeneous equation complementary solution

$$D_{\Delta x}y(x) + y(x) = 0 \;, \quad \text{the related homogeneous differential difference equation} \qquad 2)$$

Let

$$y(x) = Ke_{\Delta x}(a,x) = K(1+a\Delta x)^{\frac{x}{\Delta x}} \qquad 3)$$

Substituting this y(x) function into the related homogenous differential difference equation, Eq 2

$$Kae_2(a,x) + Ke_2(a,x) = (1 + a)Ke_2(a,x) = 0 \;, \quad \Delta x = 2$$

$$a = -1$$

Substituting a = -1 into Eq 3

The complementary solution to the differential difference equation, Eq 1 is then

$$y_c(x) = Ke_2(-1,x) \qquad 4)$$

Find the particular solution to the differential difference equation, Eq 1.

The function, $xe^{-x}$, should first be converted into its identity. This will facilitate the necessary discrete differentiation.

$$e^{cx} = e_{\Delta x}\left( \frac{e^{c\Delta x}-1}{\Delta x}, x \right)$$

$$c = -1, \quad \Delta x = 2$$

$$xe^{-x} = xe_2\left( \frac{e^{-2}-1}{2}, x \right) = xe_2(-.43233235,x) \qquad 5)$$

Substituting into Eq 1

$$D_{\Delta x}y(x) + y(x) = xe_2(-.43233235,x) \qquad 6)$$

Note that all of the functions in Eq 6 are discrete Interval Calculus functions

From the undetermined coefficient particular solution tables

$$y_p(x) = [\sum_{p=N_h}^{N+N_h} A_p[x]_{\Delta x}^{p}]\, e_{\Delta x}(\frac{e^{a\Delta x}-1}{\Delta x}, x) \qquad 7)$$

for the function

$$[\sum_{p=0}^{N} k_p[x]_{\Delta x}^{p}]\, e^{ax} = [\sum_{p=0}^{N} k_p[x]_{\Delta x}^{p}]\, e_{\Delta x}(\frac{e^{a\Delta x}-1}{\Delta x}, x) \qquad 8)$$

and the related root is

$$r = \frac{e^{a\Delta x}-1}{\Delta x}$$

For a = -1 and $\Delta x = 2$

The related root is r = -.432333235
The related homogeneous equation root is r = -1

Since the roots are not the same
$N_h = 0$

N = 1, the polynomial power, from Eq 6 and Eq 8 ( i.e. $x^1 e_2(-.43233235,x)$ )

From Eq 7

$$y_p(x) = [\sum_{p=0}^{1} A_p[x]_{\Delta x}^{p}]\, e_{\Delta x}( -.43233235,x )$$

$$y_p(x) = (A_0 + A_1 x)\, e_{\Delta x}( -.43233235,x) \qquad 9)$$

Substituting Eq 9, into Eq 6 to solve for the constants, $A_0$ , $A_1$

Using the derivative of a function product formula

$D_{\Delta x}[u(x)v(x) = v(x)D_{\Delta x}u(x) + D_{\Delta x}v(x)u(x+\Delta x)$ where $u(x) = (A_0 + A_1x)$ , $v(x) = e_{\Delta x}( -.43233235,x)$

$D_2y_p(x) = A_1e_2(-.43233235,x) + (-.43233235)e_2(-.43233235,x)(A_0 + A_1[x + 2])$

$D_2y_p(x) = (A_1 - .43233235A_0 - .43233235A_1x - .86466470A_1)e_2(-.43233235,x)$ 10)

Substituting Eq 9 and Eq 10 into Eq 6

$(A_1 - .43233235A_0 - .43233235A_1x - .86466470A_1)e_2(-.43233235,x) + (A_0 + A_1x)\, e_{\Delta x}( -43233235,x)$
$= xe_2(-.43233235,x)$

$(A_1 - .43233235A_0 - .86466470A_1 + A_0)\, e_2(-.43233235,x) + (A_1 - .43233235A_1)xe_2(-.43233235,x) =$
$xe_2(-.43233235,x)$

Simplifying and equating constants

$.56766765A_0 + .13533530A_1 = 0$

$.56766765A_1 = 1$

Solving for $A_0$, $A_1$

$A_1 = 1.76159413$

$$A_0 = \frac{-.13533530(1.76159413)}{.56766765} = -.41997438$$

$A_0 = \ = -.41997438$

Substituting the values of $A_0$, $A_1$ into Eq 9

$y_p(x) = (-.41997438 + 1.76159413x)\, e_2(-.43233235,x)$ 11)

$y(x) = y_c(x) + y_p(x)$ 12)

Substituting Eq 4 and Eq 11 into Eq 12

$y(x) = Ke_2(-1,x) + (-.41997438 + 1.76159413x)\, e_2(-.43233235,x)$ 13)

Find K from the initial condition y(0) = 3

$3 = Ke_2(-1,0) - 41997438e_2(-.43233235,0)$

$e_{\Delta x}(c,0) = (1+c\Delta x)^{\frac{0}{\Delta x}} = 1$

$3 = K - .41997438$

$K = 3.41997438$

Substituting the value of K into Eq 13

Then

$$\mathbf{y(x) = 3.14199743e_2(-1,x) + (-.41997438 + 1.76159413x)\, e_2(-.43233235,x)} \qquad 14)$$

is the solution to the differential difference equation

$D_{\Delta x}y(x) + y(x) = xe^{-x} = xe_2(-.43233235, x)$

where

$x = 0, 2, 4, 6, \ldots$

If desired, Eq 14 can be changed to another form

$$e_{\Delta x}(c,x) = (1+c\Delta x)^{\frac{x}{\Delta x}}$$

$$e_2(-1,x) = [1+(-1)(2)]^{\frac{x}{2}} = (-1)^{\frac{x}{2}}$$

From Eq 5

$e_2(-.43233235,x) = e^{-x}$

Substituting these equalities into Eq 14

$$\mathbf{y(x) = 3.14199743(-1)^{\frac{x}{2}} + (-.41997438 + 1.76159413x)\, e^{-x}} \qquad 15)$$

Check Eq 14 and Eq 15

From Eq 1

$D_{\Delta x}y(x) + y(x) - xe^{-x} = 0$

$$\frac{y(x+\Delta x)-y(x)}{\Delta x} + y(x) - xe^{-x} = 0$$

For

$x = 2$

$\Delta x = 2$

$$\frac{y(4)-y(2)}{2} + y(2) - xe^{-2} = 0$$

Finding y(2) and y(4) from Eq 15

y(2) = -2.7220231
y(4) =  3.2633642

Substituting

$$\frac{3.2633642 - (-2.7220231)}{2} + (-2.7220231) - .27067056 = 0$$

0 = 0        Good check

**Solution 2  Use of the $K_{\Delta x}$ Transform Method to solve a differential difference equation**

$$D_{\Delta x}y(x) + y(x) = xe^{-x}, \qquad \Delta x = 2, \ y(0) = 3. \qquad 1)$$

Apply $K_{\Delta x}$ Transforms from the $K_{\Delta x}$ Transform Tables in the Appendix to Eq 1

$$K_{\Delta x}[D_{\Delta x}y(x)] + K_{\Delta x}[y(x)] = K_{\Delta x}[xe^{-x}]$$

$$sy(s) - y(0) + y(s) = \frac{e^{a\Delta x}}{\left(s - \frac{e^{a\Delta x}-1}{\Delta x}\right)^2} \qquad 2)$$

where
a = -1
$\Delta x = 2$
y(0) = 3

Substituting these values into Eq 2

$$sy(s) - 3 + y(s) = \frac{e^{-2}}{\left(s - \frac{e^{-2}-1}{2}\right)^2} \qquad 3)$$

$e^{-2} = .135335283$

$$\frac{e^{-2}-1}{2} = -.432442358$$

Substituting these two values into Eq 3

$$(s + 1)y(s) = 3 + \frac{.135335283}{(s + .432332358)^2}$$

$$y(s) = \frac{3}{(s + 1)} + \frac{.135335283}{(s + 1)(s + .432332358)^2} \qquad 4)$$

$$\frac{.135335283}{(s + 1)(s + .432332358)^2} = \frac{A}{(s + 1)} + \frac{B}{(s + .432332358)^2} + \frac{C}{(s + .432332358)} \qquad 5)$$

$$A = \frac{.135335283}{(s + .432332358)^2}\Big|_{s = -1} = \frac{.135335283}{(-1 + .432332358)^2} = .41997434$$

$$B = \frac{.135335283}{(s + 1)}\Big|_{s = -.432332358} = \frac{.135335283}{(-.432332358 + 1)} = .238405843$$

$$C = \frac{d}{ds}\Big[\frac{.135335283}{(s + 1)}\Big]\Big|_{s = -.432332358} = \Big[-\frac{.135335283}{(s + 1)^2}\Big]\Big|_{s = -.432332358}$$

$$C = -\frac{.135335283}{( -.432332358 + 1)^2} = -.41997434$$

Substituting the values of A,B,C into Eq 5 and then substituting Eq 5 into Eq 4

$$y(s) = \frac{3}{s + 1} + \frac{.41997434}{s + 1} - \frac{.41997434}{s + .432332358} + \frac{.238405843}{(s + .432332358)^2} \qquad 6)$$

Find the inverse $K_{\Delta x}$ Transform of Eq 6

From the $K_{\Delta x}$ Transform Tables

$$K_{\Delta x}[e_{\Delta x}(a,x)] = \frac{1}{s - a}$$

$$K_{\Delta x}[xe_{\Delta x}(a,x)] = \frac{1 + a\Delta x}{(s - a)^2}$$

Putting Eq 6 into a form which will facilitate taking the inverse $K_{\Delta x}$ Transform

$a = -.43233235$
$\Delta x = 2$

$1+a\Delta x = 1 - .432332358(2) = .1353335284$

$$y(s) = \frac{3.41997434}{s + 1} - \frac{.41997434}{s + .432332358} + 1.761594138\Big[\frac{.135335284}{(s + .432332358)^2}\Big] \qquad 7)$$

Taking the inverse $K_{\Delta x}$ Transform of Eq 7

$y(x) = 3.41997434e_2(-1,x) - .41997434e_2(-.432332358,x) + 1.761594138x\ e_2(-.432332358,x)$

$$\mathbf{y(x) = 3.41997434e_2(-1,x) + ( - .41997434 + 1.761594138x)e_2(-.432332358,x)} \qquad 8)$$

Eq 8 is the solution to the differential difference equation:

$D_{\Delta x}y(x) + y(x) = xe^{-x}$

where
$\Delta x = 2$
$y(0) = 3.$
$x = 0,2,4,6,\ldots$

If desired, Eq 8 can be changed to another form

$$e_{\Delta x}(c,x) = (1+c\Delta x)^{\frac{x}{\Delta x}} \qquad 9)$$

$$e_2(-1,x) = [1+(-1)(2)]^{\frac{x}{2}} = [-1]^{\frac{x}{2}} \qquad 10)$$

From Eq 9

$$e_{\Delta x}(c,x) = \left[e^{\ln(1+c\Delta x)}\right]^{\frac{x}{\Delta x}}$$

Then

$$e_2(-.432332358,x) = \left[e^{\ln(1+\{-.432332358\}\{2\})}\right]^{\frac{x}{2}} = [\ e^{-2}\ ]^{\frac{x}{2}} = e^{-x}$$

$$e_2(-.432332358,x) = e^{-x} \qquad 11)$$

Substituting Eq 10 and Eq 11 into Eq 8

$$\mathbf{y(x) = 3.41997434[-1]^{\frac{x}{2}} + ( - .41997434 + 1.761594138x)e^{-x}} \qquad 12)$$

Note - The $K_{\Delta x}$ Transform solution results are the same as those obtained in Solution 1 using the Undetermined Coefficient Method.

**Solution 3 Use of the Method of Variation of Parameters to solve a differential difference equation**

$D_{\Delta x}y(x) + y(x) = xe^{-x}$, $\Delta x = 2$, $y(0) = 3$, $x = 0,2,4,6,\ldots$ 1)

Find the related homogeneous equation complementary solution, $y_c(x)$.

$D_{\Delta x}y_c(x) + y_c(x) = 0$ , the related homogeneous differential difference equation 2)

Let

$y_c(x) = ce_{\Delta x}(a,x) = c(1+a\Delta x)^{\frac{x}{\Delta x}}$ 3)

Substituting this $y_c(x)$ function into the related homogenous differential difference equation, Eq 2

$cae_{\Delta x}(a,x) + ce_{\Delta x}(a,x) = (a + 1)ce_{\Delta x}(a,x) = 0$

$a = -1$

Substituting $a = -1$ into Eq 3
The complementary solution to the differential difference equation, Eq 1 is then

$y_c(x) = ce_2(-1,x)$ 4)
Find the particular solution to the differential difference equation, Eq 1.

The function, $e^{-x}$, should first be converted into its identity. This will facilitate the necessary discrete differentiation.

Using the relationship:

$$e^{cx} = e_{\Delta x}\left(\frac{e^{c\Delta x}-1}{\Delta x},x\right)$$

$c = -1$, $\Delta x = 2$

$e^{-x} = e_2\left(\frac{e^{-2}-1}{2},x\right) = e_2(-.43233235,x)$ 5)

Substituting into Eq 1

$D_2y(x) + y(x) = xe_2(-.43233235,x)$ 6)

Note that all of the functions in Eq 6 are discrete Interval Calculus functions

Consider the constant, c, in Eq 4 to be a function of x

$c = c(x)$ 7)

From Eq 4 and Eq 7, let the particular solution, $y_p(x)$, be as follows:

$$y_p(x) = c(x)e_2(-1,x) \qquad 8)$$

Take the derivative of Eq 8.
Use the following derivative of the product of two functions formula.

$D_{\Delta x}[u(x)v(x)] = v(x)D_{\Delta x}u(x) + D_{\Delta x}v(x)u(x+\Delta x)$ where $u(x) = e_2(-1,x)$ , $v(x) = c(x)$

The following notation will be used to represent the discrete derivative.

$$f'(x) = D_{\Delta x}f(x)$$

Differentiating Eq 8

$$y_p'(x) = -c(x)e_2(-1,x) + c'(x)e_2(-1,x+2) \qquad 9)$$

Using the relationship:

$$e_{\Delta x}(b,x+\Delta x) = (1+b\Delta x)e_{\Delta x}(b,x) \qquad 10)$$

$$e_2(-1,x+2) = [1+(-1)(2)]e_2(-1,x) = -e_2(-1,x)$$

Substituting into Eq 9

$$y_p'(x) = -c(x)e_2(-1,x) - c'(x)e_2(-1,x) \qquad 11)$$

Substituting Eq 8 and Eq 11 into Eq 6

$$-c(x)e_2(-1,x) - c'(x)e_2(-1,x) + c(x)e_2(-1,x) = xe_2(-.43233235,x)$$

$$c'(x)e_2(-1,x) = -xe_2(-.43233235,x)$$

$$c'(x) = \frac{-xe_2(-.43233235,x)}{e_2(-1,x)}$$

$$c'(x) = -x[e_2(-.43233235,x)e_2(-1,-x)] \qquad 12)$$

Using the relationship:

$$e_{\Delta x}(a,x)e_{\Delta x}(b,-x) = e_{\Delta x}\left(\frac{a-b}{1+b\Delta x}, x\right) \qquad 13)$$

$$e_2(-.43233235,x)e_2(-1,-x) = e_2\left(\frac{-.432332358-(-1)}{1+(-1)(2)}, x\right)$$

$$e_2(-.43233235,x)e_2(-1,-x) = e_2(-.567667642,x) \qquad 14)$$

Substituting Eq 14 into Eq 12

$c'(x) = -xe_2(-.567667642,x)$ 15)

Integrate Eq 15 to find c(x).
Use the following integration formula for the product of two functions.

$_{\Delta x}\int v(x)D_{\Delta x}u(x)\Delta x = u(x)v(x) - {}_{\Delta x}\int D_{\Delta x}v(x)u(x+\Delta x)\Delta x + K$ 16)

Let

$v(x) = x \qquad D_{\Delta x}v(x) = 1$

$D_{\Delta x}u(x) = e_2(-.567667642,x) \qquad u(x) = \dfrac{e_2(-.567667642,x)}{-.567667642} = -1.761594155\, e_2(-.567667642,x)$

From Eq 15 and Eq 16

$$c(x) = -{}_{\Delta x}\int xe_2(-.567667642,x)\Delta x = 1.761594155xe_2(-.567667642,x) - 1.761594155\ {}_{\Delta x}\int e_2(-.567667642,x+2)\Delta x$$

Using the relationship:
$e_{\Delta x}(a,x+\Delta x) = (1+a\Delta x)e_{\Delta x}(a,x)$ 17)

$$c(x) = 1.761594155xe_2(-.567667642,x) - (1.761594155)[1+(-.567667642)(2)]\ {}_{\Delta x}\int e_2(-.567667642,x)\Delta x$$

$$c(x) = 1.761594155xe_2(-.567667642,x) - \frac{(1.761594155)(-.135335284)}{-.567667642}\, e_2(-.567667642,x)$$

$c(x) = 1.761594155xe_2(-.567667642,x) - .419974343\, e_2(-.567667642,x)$ 18)

Substituting Eq 18 into Eq 8

$y_p(x) = (-.419974343 + 1.761594155x)[e_2(-.567667642,x)\, e_2(-1,x)]$ 19)

Using the relationship:
$e_{\Delta x}(a,x)e_{\Delta x}(b,x) = e_{\Delta x}(a+b+ab\Delta x\,,x)$ 20)

$e_2(-.567667642,x)e_{\Delta x}(-1,x) = e_2(-1.567667642+1.135335285\,,x)$

$e_2(-.567667642,x)e_{\Delta x}(-1,x) = e_2(-.432332358\,,x)$

Substituting into Eq 19

$y_p(x) = (-.419974343 + 1.761594155x)e_2(-.432332358\,,x)$ 21)

$y(x) = y_c(x) + y_p(x)$ 22)

Substituting Eq 4 and Eq 21 into Eq 22

$y(x) = ce_2(-1,x) + (-.419974343 + 1.761594155x)e_2(-.432332358\,,x)$ 23)

Find c from the initial condition, y(0) = 3

$e_{\Delta x}(a,0) = 1$

$3 = c - .419974343$

$c = 3.419974343$

Substituting this value for c into Eq 23

$\mathbf{y(x) = 3.419974343e_2(-1,x) + (-.419974343 + 1.761594155x)e_2(-.43233235\,,x)}$ 24)

Eq 24 can be changed to another form

Using the relationship:

$e_{\Delta x}(a,x) = [1+a\Delta x]^{\frac{x}{\Delta x}}$ 25)

$e_2(-1,x) = [1+(-1)(2)]^{\frac{x}{2}}$

$e_2(-1,x) = [-1]^{\frac{x}{2}}$

From Eq 5

$e^{-x} = e_2(-.43233235,x)$

Substituting into Eq 24

$\mathbf{y(x) = 3.419974343[-1]^{\frac{x}{2}} + (-.419974343 + 1.761594155x)e^{-x}}$ 26)

Checking Eq 26

y(0) = 3.419974343-.419974343 = 3 good check

Substitute Eq 26 into Eq 1

$$\frac{3.419974343[-1]^{\frac{x+2}{2}} + (-.419974343 + 1.761594155\{x+2\})e^{-(x+2)}}{2}$$

$$-\frac{3.419974343[-1]^{\frac{x}{2}} + (-.419974343 + 1.761594155x)e^{-x}}{2}$$

$$+ 3.419974343[-1]^{\frac{x}{2}} + (-.419974343 + 1.761594155x)e^{-x} = e^{-x}$$

$$+\frac{-3.419974343[-1]^{\frac{x}{2}}}{2} + \frac{-.419974343e^{-2}e^{-x}}{2} + \frac{1.761594155(2)e^{-2}e^{-x}}{2} + \frac{1.761594155e^{-2}xe^{-x}}{2}$$

$$-\frac{3.419974343[-1]^{\frac{x}{2}}}{2} + \frac{.419974343e^{-x}}{2} - \frac{1.761594155xe^{-x}}{2}$$

$$+ 3.419974343[-1]^{\frac{x}{2}} - .419974343e^{-x} + 1.761594155xe^{-x} = e^{-x}$$

$$[\frac{-.419974343e^{-2}}{2} + \frac{1.761594155(2)e^{-2}}{2} + \frac{.419974343}{2} - .419974343\ ]e^{-x} +$$
$$[\frac{1.761594155e^{-2}}{2} - \frac{1.761594155}{2} + 1.761594155]xe^{-x} = xe^{-x}$$

$0e^{-x} + 1xe^{-x} = xe^{-x}$

$xe^{-x} = xe^{-x}$

Good check

**Solution 4  Use of the Method of Related Functions to solve a differential difference equation**

$D_{\Delta x}y(x) + y(x) = xe^{-x}, \qquad \Delta x = 2,\ y(0) = 3,\ x = 0,2,4,6,\dots$ 1)

Using the Method of Related Functions convert Eq 1 into a related Calculus equation.

From the table of related functions (Discrete Function $\rightleftarrows$ Calculus Function)

The following related functions are reversible.

$y(x) \rightleftarrows Y(x)$ 2)

$D_{\Delta x} \rightleftarrows \frac{d}{dx}$ 3)

$y(0) = Y(0)$ 4)

$e_{\Delta x}(a,x) \rightleftarrows e^{ax}$ 5)

$xe^{a(x-\Delta x)} \rightleftarrows xe^{\frac{e^{a\Delta x}-1}{\Delta x}x}$ 6)

Evaluating Eq 6

For $a = -1$
$\Delta x = 2$

$e^{a\Delta x} = e^{-2} = .135335283$

$\frac{e^{a\Delta x}-1}{\Delta x} = \frac{e^{-2}-1}{2} = -.432332358$

$xe^{ax} \rightleftarrows e^{a\Delta x}xe^{\frac{e^{a\Delta x}-1}{\Delta x}x}$

$xe^{-x} \rightleftarrows .135335283\,xe^{-.432332358x}$ 7)

From the related functions of Eq 2 and Eq 7, convert the discrete differential difference equation, Eq 1, into its related Calculus differential equation.

$\frac{d}{dx}Y(x) + Y(x) = .135335283\,xe^{-.432332358x}, \quad Y(0) = y(0) = 3$ 8)

Using Calculus solve Eq 8

The homogeneous differential equation for Eq 8 is:

$$\frac{d}{dx} Y_c(x) + Y_c(x) = 0 \qquad 9)$$

Solve Eq 9 for the related Calculus differential equation complementary solution, $Y_c(x)$.

Let

$$Y_c(x) = ce^{bx} \qquad 10)$$

Substituting this $Y_c(x)$ function into the homogenous differential equation, Eq 9

$$cbe^{bx} + ce^{bx} = (b + 1)ce^{bx} = 0$$

$$b = -1$$

Substituting b = -1 into Eq 10
The complementary solution to the differential equation, Eq 8, is:

$$Y_c(x) = ce^{-x} \qquad 11)$$

Find the particular solution to the differential equation, Eq 8, using the Method of Undetermined Coefficients

This method specifies a particular solution:

$$Y_p(x) = (K_1x + K_2)e^{-.432332358x} \qquad 12)$$

Substitute Eq 12 into Eq 8

$$K_1e^{-.432332358x} - .432332358(K_1x + K_2)e^{-.432332358x} + (K_1x + K_2)e^{-.432332358x} = .135335283xe^{-.432332358x}$$

Combining terms

$$e^{-.432332358x}[K_1 - .432332358K_2 + K_2] + xe^{-.432332358x}[-.432332358K_1 + K_1] = .135335283xe^{-.432332358x}$$

$$(1-.432332358)K_1 = .135335283$$

$$K_1 = .238405843$$

$$K_1 + (1-.432332358)K_2 = 0$$

$$K_2 = \frac{-.238405843}{1-.432332358}$$

$$K_2 = -.419974339$$

Substituting the constants $K_1$ and $K_2$ into Eq 12

$Y_p(x) = (.238405843x - .419974339)e^{-.432332358x}$ 13)

$Y(x) = Y_c(x) + Y_p(x)$ 14)

Substituting Eq 11 and Eq 13 into Eq 14

$Y(x) = ce^{-x} + (- .419974339 + .238405843x)e^{-.432332358x}$ 15)

Find the constant, c, using the initial condition, Y(0) = 3

$3 = c - .419974339$

$c = 3.419974339$

Substituting the value of c into Eq 15

$Y(x) = 3.419974339e^{-x} + (- .419974339 + .238405843x)e^{-.432332358x}$ 16)

This is the related Calculus function to y(x).

Using the related functions of Eq 2 thru Eq 7 and the following related functions and equations, the Calculus function, Y(x), can be converted into its related discrete Interval Calculus function, y(x).

From the related functions table

$e_{\Delta x}(a,x) = (1+a\Delta x)^{\frac{x}{\Delta x}} = e^{\frac{x}{\Delta x}\ln(1+a\Delta x)} \rightleftarrows e^{ax}$ 17)

For $a = -.432332358$
$\Delta x = 2$

Substiituting into Eq 17

$e_2(-.432332358,x) = (1+ [-.432332358][2])^{\frac{x}{2}} = e^{\frac{\ln[1+(-.432332358)(2)]}{2}x} = e^{-x} \rightleftarrows e^{-.432332358x}$

$e^{-x} \rightleftarrows e^{-.432332358x}$ 18)

From Eq 7

$\frac{1}{.135335283}xe^{-x} \rightleftarrows xe^{-.432332358x}$ 19)

From Eq 5

$$e_{\Delta x}(a,x) = (1+a\Delta x)^{\frac{x}{\Delta x}} \rightleftarrows e^{ax}$$

a = -1
Δx = 2

$$e_2(-1,x) = [1+(-1)(2)]^{\frac{x}{2}} = [-1]^{\frac{x}{2}} \rightleftarrows e^{-x} \qquad 20)$$

$$[-1]^{\frac{x}{2}} \rightleftarrows e^{-x} \qquad 21)$$

From Eq 2

$$y(x) \rightleftarrows Y(x) \qquad 22)$$

Substituting the related functions of Eq 18. Eq 19, Eq 21 and Eq 22 into Eq 16

$$y(x) = 3.419974339[-1]^{\frac{x}{2}} + (-\ .419974339 + \frac{.238405843}{.135335283}x)e^{-x}$$

$$\mathbf{y(x) = 3.419974339[-1]^{\frac{x}{2}} + (-\ .419974339 + 1.761594151x)e^{-x}} \qquad 23)$$

Eq 23 can be put into another form

From Eq 20

$$[-1]^{\frac{x}{2}} = e_2(-1,x) \qquad 24)$$

$$e^{ax} = (1 + \frac{e^{a\Delta x}-1}{\Delta x}\Delta x)^{\frac{x}{\Delta x}} = e_{\Delta x}(\frac{e^{a\Delta x}-1}{\Delta x},x) \qquad 25)$$

For a = -1
Δx = 2

$$e^{-x} = e_2(\frac{e^{-2}-1}{2},x)$$

$$e^{-x} = e_2(-.432332358,x) \qquad 26)$$

Substituting Eq 24 and Eq 26 into Eq 23

$$\mathbf{y(x) = 3.419974339e_2(-1,x) + (-\ .419974339 + 1.761594151x)e_2(-.432332358,x)} \qquad 27)$$

Checking Eq 23 and Eq 27

y(0) = 3.4199743393-.419974339 = 3     Good check

Substitute Eq 23 into Eq 1

$$\frac{3.419974339[-1]^{\frac{x+2}{2}} + (-.419974339 + 1.761594151\{x+2\})e^{-(x+2)}}{2}$$

$$-\frac{3.419974339[-1]^{\frac{x}{2}} + (-.419974339 + 1.761594151x)e^{-x}}{2}$$

$$+ 3.419974339[-1]^{\frac{x}{2}} + (-.419974339 + 1.761594151x)e^{-x} = xe^{-x}$$

$$+\frac{-3.419974339[-1]^{\frac{x}{2}}}{2} + \frac{-.419974339e^{-2}e^{-x}}{2} + \frac{1.761594151(2)e^{-2}e^{-x}}{2} + \frac{1.761594151e^{-2}xe^{-x}}{2}$$

$$-\frac{3.419974339[-1]^{\frac{x}{2}}}{2} + \frac{.419974339e^{-x}}{2} - \frac{1.761594151xe^{-x}}{2}$$

$$+ 3.419974339[-1]^{\frac{x}{2}} - .419974339e^{-x} + 1.761594151xe^{-x} = xe^{-x}$$

$$[\frac{-.419974339e^{-2}}{2} + \frac{1.761594151(2)e^{-2}}{2} + \frac{.419974339}{2} - .419974339\ ]e^{-x} +$$

$$[\frac{1.761594151e^{-2}}{2} - \frac{1.761594151}{2} + 1.761594151]xe^{-x} = xe^{-x}$$

$0e^{-x} + 1xe^{-x} = xe^{-x}$

$xe^{-x} = xe^{-x}$

Good check

**Example 7.18** Find the $K_{\Delta x}$ Transform, $K_{\Delta x}[e_{\Delta x}(a,x)] = \frac{1}{s-a}$ two different ways, and the Inverse $K_{\Delta x}$ Transform, $K_{\Delta x}^{-1}[\frac{1}{s-a}] = e_{\Delta x}(a,x)$ using the Inverse $K_{\Delta x}$ Transform.

Solution 1 Find the $K_{\Delta x}$ Transform of the function $e_{\Delta x}(a,x)$ using the formula,

$$K_{\Delta x}[f(x)] = {}_{\Delta x}\int_0^\infty (1+s\Delta x)^{-(\frac{x+\Delta x}{\Delta x})} f(x)\Delta x.$$

$$e_{\Delta x}(a,x) = (1+a\Delta x)^{\frac{x}{\Delta x}} \qquad 1)$$

$$K_{\Delta x}[f(x)] = {}_{\Delta x}\int_0^\infty (1+s\Delta x)^{-(\frac{x+\Delta x}{\Delta x})} f(x)\Delta x \qquad 2)$$

Let

$$f(x) = e_{\Delta x}(a,x) = (1+a\Delta x)^{\frac{x}{\Delta x}} \qquad 3)$$

Substituting Eq 3 into the $K_{\Delta x}$ Transform equation, Eq 2

$$K_{\Delta x}[e_{\Delta x}(a,x)] = \frac{1}{1+s\Delta x}\ {}_{\Delta x}\int_0^\infty \left(\frac{1+s\Delta x}{1+a\Delta x}\right)^{-(\frac{x}{\Delta x})} \Delta x \qquad 4)$$

$$\frac{1+s\Delta x}{1+a\Delta x} = 1 + \frac{(s-a)\Delta x}{1+a\Delta x} \qquad 5)$$

Substituting Eq 5 into Eq 4

$$K_{\Delta x}[e_{\Delta x}(a,x)] = \frac{1}{1+s\Delta x}\ {}_{\Delta x}\int_0^\infty \left(1+\frac{(s-a)\Delta x}{1+a\Delta x}\right)^{-(\frac{x}{\Delta x})} \Delta x \qquad 6)$$

Substitute the quantity, $\frac{s-a}{1+a\Delta x}$, for b in the following equation. This integral equation was obtained from the integral equation table, Table 6, in the Appendix.

$$ {}_{\Delta x}\int (1+b\Delta x)^{-\frac{x}{\Delta x}}\Delta x = -\frac{1+b\Delta x}{b}(1+b\Delta x)^{-\frac{x}{\Delta x}} + k \qquad 7) $$

$$ K_{\Delta x}[e_{\Delta x}(a,x)] = \frac{1}{1+s\Delta x}\,{}_{\Delta x}\int_{0}^{\infty}\left(1+\frac{(s-a)\Delta x}{1+a\Delta x}\right)^{-\frac{x}{\Delta x}}\Delta x = -\frac{1}{1+s\Delta x}\,\frac{1+s\Delta x}{s-a}\left(1+\frac{(s-a)\Delta x}{1+a\Delta x}\right)^{-\frac{x}{\Delta x}}\Big|_{0}^{\infty} \qquad 8) $$

$$ \mathbf{K_{\Delta x}[e_{\Delta x}(a,x)] = \frac{1}{s-a}} \qquad 9) $$

Eq 9 specifies the $K_{\Delta x}$ Transform of $e_{\Delta x}(a,x)$

Solution 2 Find the $K_{\Delta x}$ Transform of the function $e_{\Delta x}(a,x)$ using the formula,

$$ K_{\Delta x}[f(x)] = \int_{0}^{\infty} e_{\Delta x}(s,-x-\Delta x)f(x)\Delta x. $$

$$ K_{\Delta x}[f(x)] = \int_{0}^{\infty} e_{\Delta x}(s,-x-\Delta x)f(x)\Delta x \qquad 10) $$

From Table 5 in the Appendix

$$ e_{\Delta x}(s,-x-\Delta x) = \frac{1}{1+s\Delta x}\, e_{\Delta x}(s,-x) \qquad 11) $$

Substituting Eq 11 into Eq 10

$$ K_{\Delta x}[f(x)] = \frac{1}{1+s\Delta x}\int_{0}^{\infty} e_{\Delta x}(s,-x)f(x)\Delta x \qquad 12) $$

Let

$$ f(x) = e_{\Delta x}(a,x) \qquad 13) $$

Substitute Eq 13 into Eq 12)

$$ K_{\Delta x}[e_{\Delta x}(a,x)] = \frac{1}{1+s\Delta x}\int_{0}^{\infty} e_{\Delta x}(a,x)e_{\Delta x}(s,-x)\Delta x \qquad 14) $$

From Table 5 in the Appendix

$$e_{\Delta x}(a,x)e_{\Delta x}(s,-x) = e_{\Delta x}(\frac{a-s}{1+s\Delta x},x) \qquad 15)$$

Substituting Eq 15 into Eq 14

$$K_{\Delta x}[e_{\Delta x}(a,x)] = \frac{1}{1+s\Delta x}\int_{0}^{\infty} e_{\Delta x}(\frac{a-s}{1+s\Delta x},x)\Delta x \qquad 16)$$

From Table 6 in the Appendix

$$_{\Delta x}\int e_{\Delta x}(c,x)\Delta x = \frac{1}{c}e_{\Delta x}(c,x) + k \qquad 17)$$

From Eq 16 and Eq 17

$$K_{\Delta x}[e_{\Delta x}(a,x)] = (\frac{1}{1+s\Delta x})\,(\frac{1+s\Delta x}{a-s})\,e_{\Delta x}(\frac{a-s}{1+s\Delta x},x)\,\Big|_{0}^{\infty} \qquad 18)$$

Comment - In the derivation of the $K_{\Delta x}$ Transform and its inverse, the value of s is defined as:

$$s = \frac{e^{(\gamma+jw)\Delta x} - 1}{\Delta x}\,, \quad -\frac{\pi}{\Delta x} \leq w < \frac{\pi}{\Delta x}$$

where $\gamma$ is any positive real value which makes s an indefinitely large value.

$$K_{\Delta x}[e_{\Delta x}(a,x)] = (\frac{1}{1+s\Delta x})\,(\frac{1+s\Delta x}{a-s})\,e_{\Delta x}(\frac{a-s}{1+s\Delta x},x)\,\Big|_{0}^{\infty} = \frac{1}{a-s}\,(1+\frac{a-s}{1+s\Delta x}\Delta x)^{\frac{x}{\Delta x}}\,\Big|_{0}^{\infty} = \frac{1}{a-s}\,(\frac{1+a\Delta x}{1+s\Delta x})^{\frac{x}{\Delta x}}\,\Big|_{0}^{\infty}$$

$\gamma$ is chosen so that $\left|\frac{1+a\Delta x}{1+s\Delta x}\right| < 1$

$$K_{\Delta x}[e_{\Delta x}(a,x)] = \frac{1}{a-s}\,(0-1) = \frac{1}{s-a} \qquad 19)$$

$$\mathbf{K_{\Delta x}[e_{\Delta x}(a,x)] = \frac{1}{s-a}} \qquad 20)$$

Eq 20 specifies the $K_{\Delta x}$ Transform of $e_{\Delta x}(a,x)$

Find the Inverse $K_{\Delta x}$ Transform of the function, $f(s) = \frac{1}{s-a}$, using the Inverse $K_{\Delta x}$ Transform Integral.

$$K_{\Delta x}^{-1}[f(s)] = F(x) = \frac{1}{2\pi j} \oint_c [1+s\Delta x]^{\frac{x}{\Delta x}} f(s)\, ds \text{ , The Inverse } K_{\Delta x} \text{ Transform} \qquad 21)$$

Notes – The x variable of F(x) is discrete where $x = 0, \Delta x, 2\Delta x, 3\Delta x, \ldots$

The s variable of the complex plane closed contour integral of Eq 10 is continuous along the closed contour, c.

$\frac{x}{\Delta x} = 0, 1, 2, 3, \ldots$

The closed contour, c, in the complex plane is a circle of radius, $\frac{e^{\gamma\Delta t}}{\Delta t}$, with center at $-\frac{1}{\Delta t}$.

The positive constant, $\gamma$, is made large enough to have the closed contour, c, encircle all poles of $(1+s\Delta x)^{\frac{x}{\Delta x}} f(s)$.

To calculate the integral of Eq 20, the theory of residues will be used.

Provided that the function of s, $[1+s\Delta x]^{\frac{x}{\Delta x}} f(s)$, can be expanded into a convergent Laurent Series for each pole of f(s)

$$K_{\Delta x}^{-1}[f(s)] = \frac{1}{2\pi j} \oint_c [1+s\Delta x]^{\frac{x}{\Delta x}} f(s)\, ds = \sum_{p=1}^{P} R_p \qquad 22)$$

$$R = \lim_{s\to r} \frac{1}{(n-1)!} \frac{d^{n-1}}{ds^{n-1}} [(s-r)^n (1+s\Delta x)^{\frac{x}{\Delta x}} f(s)] \text{ , the residue calculation formula for a pole at } s = r$$

where 23)

$R$ = the residue of a pole of $(1+s\Delta x)^{\frac{x}{\Delta x}} f(s)$

$\sum_{p=1}^{P} R_p$ = the sum of the P residues of the P poles of $(1+s\Delta x)^{\frac{x}{\Delta x}} f(s)$

Let $f(s) = \frac{1}{s-a}$

$r = a$

$P = 1$ , $(1+s\Delta x)^{\frac{x}{\Delta x}} \frac{1}{s-a}$ has one pole which is at $s = a$

Substituting the specified f(s), r, and P into Eq 22 and Eq 23

$$K_{\Delta x}^{-1}[\frac{1}{s-a}] = \frac{1}{2\pi j} \oint_c [1+s\Delta x]^{\frac{x}{\Delta x}} \frac{1}{s-a}\, ds = \sum_{p=1}^{1} R_p = R_1 \qquad 24)$$

$n = 1$ , the pole at $s = a$ is a first order pole

<u>Note</u> - $(1+s\Delta x)^{\frac{x}{\Delta x}} \frac{1}{s-a}$ has a convergent Laurent Series which is:

$$\frac{(1+s\Delta x)^{\frac{x}{\Delta x}}}{s-a} = e_{\Delta x}(a,x) \sum_{n=0}^{\infty} \frac{1}{n!} \frac{[x]_{\Delta x}^{n}}{(1+a\Delta x)^{n}} \frac{1}{(s-a)^{1-n}}$$

Since there is only one pole to be encircled by the complex plane closed contour, c, there is only one residue required, $R_{1,}$ the residue at $s = a$.

From Eq 23

$$R_1 = \lim_{s\to a} [(s-a)\,(1+s\Delta x)^{\frac{x}{\Delta x}} \frac{1}{s-a}] = (1+a\Delta x)^{\frac{x}{\Delta x}} = e_{\Delta x}(a,x) \qquad 25)$$

Substituting $R_1$ into Eq 24

$$\mathbf{K_{\Delta x}^{-1}[\frac{1}{s-a}] = \frac{1}{2\pi j} \oint_c [1+s\Delta x]^{\frac{x}{\Delta x}} \frac{1}{s-a}\, ds = e_{\Delta x}(a,x)} \qquad 26)$$

Eq 26 specifies the Inverse $K_{\Delta x}$ Transform of $\frac{1}{s-a}$ .

**Example 7.19** Determine for what range of Δx a control system with a differential difference equation, $D_{\Delta x}^2 y_{\Delta x}(x) + 20D_{\Delta x}y_{\Delta x}(x) + 125y_{\Delta x}(x) = 0$, is stable.

Comment – The differential difference equation, $D_{\Delta x}^2 y_{\Delta x}(x) + 20D_{\Delta x}y_{\Delta x}(x) + 125y_{\Delta x}(x) = 0$ is actually a set of related differential difference equations, one for every value of Δx.

$$D_{\Delta x}^2 y_{\Delta x}(x) + 20D_{\Delta x}y_{\Delta x}(x) + 125y_{\Delta x}(x) = 0 \qquad 1)$$

Find the roots of the characteristic polynomial, $s^2 + 20s + 125 = 0$

$$s = \frac{-20 \pm \sqrt{(20)^2 - 4(125)}}{2(1)} = -10 \pm 5j \qquad 2)$$

The roots are in the left half of the complex plane so for $\Delta x \to 0$ the control system is stable.

Plot the roots, $-10 \pm 5j$, in the complex plane and graphically find the Δx interval at which the control system goes unstable. Use the Interval Calculus Modified Nyquist Criteria. For stability, the Interval Calculus Modified Nyquist Criteria requires all roots to be within the complex plane critical circle of radius, $\frac{1}{\Delta x}$, with center at the point $(-\frac{1}{\Delta x}, 0)$.

Draw a circle through the points (-10,5), (-10,-5), (0,0) with the center somewhere on the complex plane negative real axis. Find the centerpoint. The center point value is equal to $-\frac{1}{\Delta x}$.

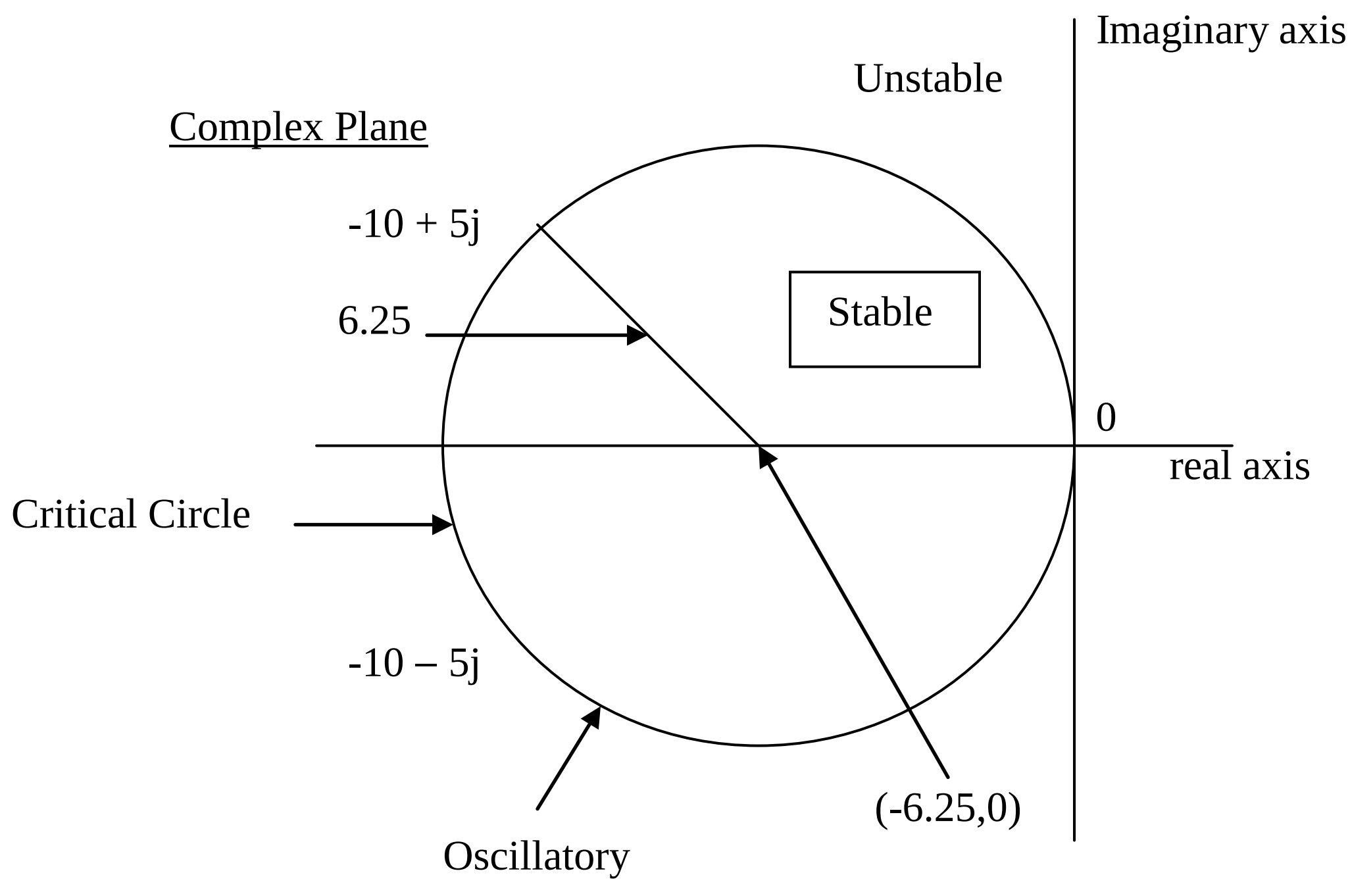

$$-\frac{1}{\Delta x} = -6.25 = -\frac{25}{4}$$ 3)

$$\Delta x = \frac{4}{25} = .16$$ 4)

The control system is oscillatory for $\Delta x = .16$ 5)
The maximum $\Delta x$ interval for which the control system is stable is $\Delta x < .16$ 6)

Comment - This method works well since the poles of the $K_{\Delta x}$ Transform of the specified differential difference equation are not functions of $\Delta x$.

Thus, the control system is stable for $0 < \Delta x < .16$. 7)

Checking the previously calculated value of $\Delta x$, $\Delta x = .16$.

As previously derived, the solution to the differential difference equation,
$D_{\Delta x}^2 y_{\Delta x}(x) + AD_{\Delta x}y_{\Delta x}(x) + By_{\Delta x}(x) = 0$, is:

$$y_{\Delta x}(x) = \left[\sqrt{(1+a\Delta x)^2 + (b\Delta x)^2}\ \right]^{\frac{x}{\Delta x}} \left[ K_1\cos\frac{\beta x}{\Delta x} + K_2\sin\frac{\beta x}{\Delta x} \right]$$

where

$x = m\Delta x$ ,  $m = 0,1,2,3,\ldots$

$K_1,K_2,a,b$ = real constants

$$\beta = \begin{cases} \tan^{-1}\frac{b\Delta x}{1+a\Delta x} & \text{for } 1+a\Delta x \ge 0 \\ \pi+\tan^{-1}\frac{b\Delta x}{1+a\Delta x} & \text{for } 1+a\Delta x < 0 \end{cases}$$

$a+jb$, $a-jb$ = the roots of the characteristic polynomial $s^2 + As + B$

$\Delta x$ = x increment , $\Delta x > 0$

For an oscillatory condition (the stability limit), the quantity, $\sqrt{(1+a\Delta x)^2 + (b\Delta x)^2}$ , must equal 1.

$a = -10$
$b = 5$
$\Delta x = .16$ (the calculated value of $\Delta x$)

$$\sqrt{\{1+(-10)(.16)\}^2 + \{(5)(.16)\}^2} = \sqrt{(-.6)^2 + (.8)^2} = 1$$

Good check

**Example 7.20** Evaluating a Z Transform defined closed loop system using $K_{\Delta t}$ Transforms

Problem Description

Solve the following Z Transform defined closed loop system for c(t) using $K_{\Delta t}$ Transforms.
As a check, solve again for c(t) using Z Transforms.

Closed Loop Control System Z Transform Diagram

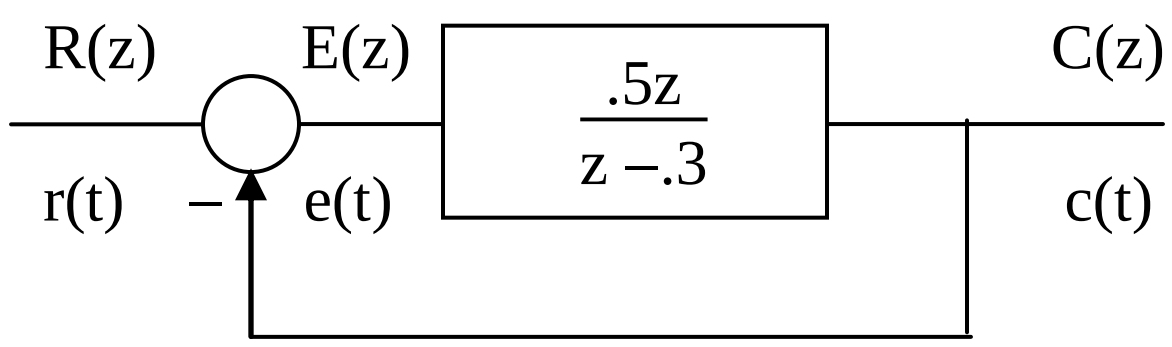

t = 0, T, 2T, 3T, …

r(t) = U(t) , unit step at t = 0

$R(z) = \dfrac{z}{z-1}$ , Z Transform of U(t)

T = .2

Initial conditions = 0

1) Find c(t) using $K_{\Delta t}$ Transforms

Convert the above closed loop control system Z Transform diagram to an equivalent $K_{\Delta t}$ Transform control system diagram.

$$\frac{C(z)}{E(z)} = \frac{.5z}{z-.3} \qquad 1)$$

Use the following Z Transform to $K_{\Delta t}$ Transform Transfer Function Conversion equations to convert $\dfrac{C(z)}{E(z)}$ to its equivalent $\dfrac{C(s)}{E(s)}$

$$\frac{O(z)}{I(z)} = F(z) \text{ , Z Transform equation} \qquad 2)$$

$$\frac{O(s)}{I(s)} = F(z)\Big|_{\substack{z = 1+s\Delta t \\ \Delta t = T}} \text{ , Equivalent } K_{\Delta t} \text{ Transform equation} \qquad 3)$$

t = 0, Δx, 2Δx, 3Δx, …

From Eq 1 thru Eq 3
Calculate the equivalent $K_{\Delta t}$ Transform equation for Eq 1

$$\frac{C(s)}{E(s)} = \frac{.5z}{z-.3}\Big|_{\substack{z = 1+s\Delta t \\ \Delta t = T \\ t = 0,\ \Delta x,\ 2\Delta x,\ 3\Delta x,\ \ldots}} \quad , \qquad 4)$$

$$\Delta t = T = .2 \qquad 5)$$

Simplifying Eq 4

$$\frac{C(s)}{E(s)} = \frac{.5(1+.2s)}{(1+.2s)-.3} = \frac{(.1s+.5)}{.2s+.7} = \frac{.1(s+5)}{.2(s+3.5)}$$

$$\frac{C(s)}{E(s)} = \frac{.5(s+5)}{(s+3.5)} \quad , \text{ The equivalent } K_{\Delta t} \text{ Transform equation for Eq 1} \qquad 6)$$

Then the equivalent $K_{\Delta t}$ Transform system diagram is as follows

Closed Loop Control System $K_{\Delta t}$ Transform Diagram

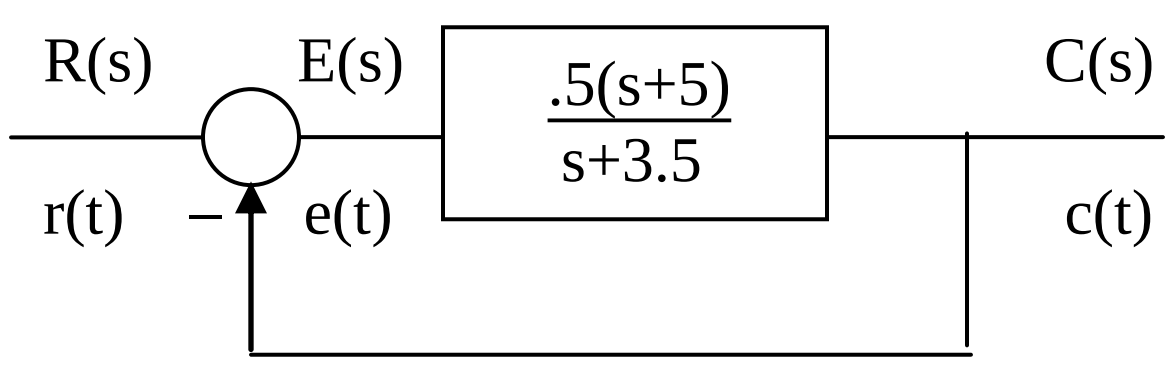

$t = 0,\ \Delta t,\ 2\Delta t,\ 3\Delta t,\ \ldots$

$r(t) = U(t)$ , unit step at $t = 0$

$R(s) = \frac{1}{s}$ , $K_{\Delta t}$ Transform of $U(t)$

$\Delta t = .2$

Find C(s)

$$C(s) = \frac{\frac{.5(s+5)}{s+3.5}}{1 + \frac{.5(s+5)}{s+3.5}} R(s) = \frac{\frac{.5(s+5)}{s+3.5}}{1 + \frac{.5(s+5)}{s+3.5}} \frac{1}{s} \qquad 7)$$

$$C(s) = \frac{.5(s+5)}{1.5s+6}\frac{1}{s} = \frac{1}{3}\frac{(s+5)}{s+4}\frac{1}{s}$$

$$C(s) = \frac{1}{3}\frac{(s+5)}{s(s+4)} \qquad 8)$$

Expand Eq 8 in a partial fraction expansion

$$C(s) = \frac{1}{3}\frac{(s+5)}{s(s+4)} = \frac{A}{s} + \frac{B}{s+4} \qquad 9)$$

$$A = \frac{1}{3}\frac{(s+5)}{(s+4)}\Big|_{s=0} = \frac{5}{3(4)} = \frac{5}{12}$$

$$A = \frac{5}{12} \qquad 10)$$

$$B = \frac{1}{3}\frac{(s+5)}{s}\Big|_{s=-4} = \frac{-4+5}{3(-4)} = -\frac{1}{12}$$

$$B = -\frac{1}{12} \qquad 11)$$

Substituting Eq10 and Eq 11 into Eq 9

$$C(s) = \frac{5}{12}\frac{1}{s} - \frac{1}{12}\frac{1}{s+4} \qquad 12)$$

Taking the Inverse $K_{\Delta t}$ Transform of Eq 12

$$K_{\Delta t}[1] = \frac{1}{s} \qquad 13)$$

$$K_{\Delta t}[e_{\Delta t}(a,t)] = \frac{1}{s-a} \qquad 14)$$

$$\Delta t = .2 \qquad 15)$$

From Eq 13 thru Eq 15

$$c(t) = \frac{5}{12} - \frac{1}{12} e_{.2}(-4,t) \qquad 16)$$

$$e_{\Delta t}(a,t) = (1+a\Delta t)^{\frac{t}{\Delta t}} \qquad 17)$$

From Eq 15 thru 17

$$c(t) = \frac{5}{12} - \frac{1}{12}[1-4(.2)]^{\frac{t}{.2}} = \frac{5}{12} - \frac{1}{12}[.2]^{5t} \qquad 18)$$

$$c(t) = \frac{5}{12} - \frac{1}{12}\left[\frac{1}{5}\right]^{5t} \qquad 19)$$

Check the above result using a Z Transform evaluation of the control system defined by the Closed Loop Control System Z Transform Diagram.

2) Find c(t) using Z Transforms

Find C(z)

$$C(z) = \frac{\frac{.5z}{z-.3}}{1 + \frac{.5z}{z-.3}} R(z) = \frac{\frac{.5z}{z-.3}}{1 + \frac{.5z}{z-.3}} \frac{z}{z-1} \qquad 20)$$

$$C(z) = \frac{.5z}{z-.3+.5z} \frac{z}{z-1} = \frac{.5z}{1.5z-.3} \frac{z}{z-1}$$

$$C(z) = \frac{1}{3} \frac{z}{z-\frac{1}{5}} \frac{z}{z-1} = \frac{z}{3} \frac{z}{(z-\frac{1}{5})(z-1)}$$

$$C(z) = \frac{z}{3} \frac{z}{(z-1)(z-\frac{1}{5})} \qquad 21)$$

Expand Eq 21 in a partial fraction expansion

$$C(z) = \frac{z}{3} \frac{z}{(z-\frac{1}{5})(z-1)} = \frac{z}{3}\left[\frac{A}{z-1} + \frac{B}{z-\frac{1}{5}}\right] \qquad 22)$$

$$A = \frac{z}{z-\frac{1}{5}}\Big|_{z=1} = \frac{1}{\frac{4}{5}} = \frac{5}{4}$$

$$A = \frac{5}{4} \qquad 23)$$

$$B = \frac{z}{z-1}\Big|_{s=\frac{1}{5}} = \frac{\frac{1}{5}}{-\frac{4}{5}} = -\frac{1}{4}$$

$$B = -\frac{1}{12} \qquad 24)$$

Substituting Eq 23 and Eq 24 into Eq 22

$$C(z) = \frac{5}{4(3)} \frac{z}{z-1} - \frac{1}{4(3)} \frac{z}{z - \frac{1}{5}} \qquad 25)$$

$$C(z) = \frac{5}{12} \frac{z}{z-1} - \frac{1}{12} \frac{z}{z - \frac{1}{5}} \qquad 26)$$

Taking the Inverse Z Transform of Eq 26

$$Z[1] = \frac{z}{z-1} \qquad 27)$$

$$Z[a^{wt}] = \frac{z}{z - a^{wT}} \qquad 28)$$

$$T = .2 \qquad 29)$$

$$a^{wt} = [a^{wT}]^{\frac{t}{T}} = [\frac{1}{5}]^{\frac{t}{.2}} = [\frac{1}{5}]^{5t} \qquad 30)$$

From 5.5-44 thru 5.5-48

$$c(t) = \frac{5}{12} - \frac{1}{12} [\frac{1}{5}]^{5t} \qquad 31)$$

This is the same result as Eq 19.

Good check

**Example 7.21** Evaluating in detail several Zeta Function summations, $\sum_{\Delta x,\, x=x_1}^{x_2} \frac{1}{x^n}$ ( $x_1,x_2$ may be infinite), using the lnd(n,Δx,x) function calculation program, LNDX, and some pertinent equations. Demonstrate the validity of all of the equations and concepts presented.

For this demonstration, arbitrarily let:
n = 3.1-1.2i
Δx = 1+i
x = -1000 + m(1+i) , m = integers This is a linear x locus in the complex plane

Equations

E1) $$\Delta x \sum_{\Delta x,\, x=x_1}^{x_2} \frac{1}{x^n} = \int_{\Delta x,\, x_1}^{x_2} \frac{1}{x^n}\,\Delta x = -\mathrm{lnd}(n,\Delta x,x)\Big|_{x_1}^{x_2+\Delta x} = -\mathrm{lnd}(n,\Delta x,x_2+\Delta x) + \mathrm{lnd}(n,\Delta x,x_1)]$$

$x = x_1, x_1+\Delta x, x_1+2\Delta x, x_1+3\Delta x, \ldots, x_2-\Delta x, x_2$

$x_1,x_2$ can be infinite values

E2) $$\Delta x \sum_{\Delta x,\, x=x_i}^{\pm\infty} \frac{1}{x^n} = \int_{\Delta x,\, x_i}^{\pm\infty} \frac{1}{x^n}\,\Delta x = \mathrm{lnd}(n,\Delta x,x_i) \approx \mathrm{lnd}_f(n,\Delta x,x_i) + \begin{cases} K_r \\ k_r \end{cases}$$

$x = x_i, x_i+\Delta x, x_i+2\Delta x, x_i+3\Delta x + x_i+4\Delta x, x_i+5\Delta x, \ldots$

Re(n) > 1

Accuracy increases rapidly as $|\frac{x}{\Delta x}|$ increases in value.

$K_r,k_r$ = constants of integration
r = 1, 2, or 3

E3) $$\Delta x \sum_{\Delta x,\, x=\mp\infty}^{\pm\infty} \frac{1}{x^n} = \mathrm{lnd}(n,\Delta x,x_i) - \mathrm{lnd}(n,-\Delta x,x_i-\Delta x) = \pm (K_r - k_r)$$

Re(n)>1

E4) $$\ln d_f(n,\Delta x,x) = -\sum_{m=0}^{\infty}\frac{\Gamma(n+2m-1)\left(\frac{\Delta x}{2}\right)^{2m} C_m}{\Gamma(n)(2m+1)!\left(x-\frac{\Delta x}{2}\right)^{n+2m-1}}$$

E5) $\ln d(n,\Delta x,x) \approx \ln d_f(n,\Delta x,x) + \begin{cases} K_r \\ k_r \end{cases}$ , x is within the x locus segment, r

The absolute value, $|\frac{x}{\Delta x}|$ , must be large for good accuracy

E6) $\left.\begin{matrix} K_r \\ k_r \end{matrix}\right\} \approx \ln d(n,\Delta x,x) - \ln d_f(n,\Delta x,x)$ , x is within the x locus segment, r

The absolute value, $|\frac{x}{\Delta x}|$ , must be large for good accuracy

E7) $$\ln d(n,\Delta x,x_i) = \Delta x \sum_{\Delta x \; x=x_i}^{x_p-\Delta x}\frac{1}{x^n} + \ln d_f(n,\Delta x,x_p) + \begin{cases} K_r \\ k_r \end{cases}$$

The absolute value, $|\frac{x}{\Delta x}|$ , must be large for good accuracy

E8) $K_1 = 0$

E9) $k_3 = 0$

E10) $K_3 = -k_1$

E11) $\ln d(n,\Delta x,x) - \ln d(n,-\Delta x,x-\Delta x) = \pm(K_r - k_r)$

E12) $\ln d(n,\Delta x,x) - \ln d(n,-\Delta x,x-\Delta x) = -[\ln d(n,-\Delta x,x) - \ln d(n,\Delta x,x+\Delta x)]$

E13) $$k_1 = \pm\left[-\Delta x \sum_{\Delta x \; x=\mp\infty}^{\pm\infty}\frac{1}{x^n}\right], \quad Re(n)>1$$

E14) $$K_3 = \pm\left[\Delta x \sum_{\Delta x \; x=\mp\infty}^{\pm\infty}\frac{1}{x^n}\right], \quad Re(n)>1$$

E15) $K_2 - k_2 = \pm [\Delta x \sum_{x=\mp\infty}^{\pm\infty} {}_{\Delta x} \frac{1}{x^n}]$ , Re(n)>1

E16) $K_r - k_r = \pm [\Delta x \sum_{x=\mp\infty}^{\pm\infty} {}_{\Delta x} \frac{1}{x^n}]$ , Re(n)>1 , r = 1,2, or 3

E17) $\left.\begin{matrix} K_r \\ k_r \end{matrix}\right\} = \Delta x \sum_{x=x_i}^{\pm\infty} {}_{\Delta x} \frac{1}{x^n} - \ln d_f(n,\Delta x,x_i)$ , Re(n)>1 , r = 1,2, or 3

$x_i$ is within the x locus segment, r

Accuracy increases rapidly as $|\frac{x}{\Delta x}|$ increases in value.

$x_i$ to $+\infty$ for Re(Δx)>0 or {Re(Δx)=0 and Im(Δx)>0}
$x_i$ to $-\infty$ for Re(Δx)<0 or {Re(Δx)=0 and Im(Δx)<0}

where

$x = x_i + m\Delta x$ , m = integers , This is an x locus in the complex plane
$x_i, x_p$ = values of x
Δx = x increment
+ for Re(Δx)>0 or {Re(Δx)=0 and Im(Δx)>0}
– for Re(Δx)<0 or {Re(Δx)=0 and Im(Δx)<0}
$-\infty$ to $+\infty$ for Re(Δx)>0 or {Re(Δx)=0 and Im(Δx)>0}
$+\infty$ to $-\infty$ for Re(Δx)<0 or {Re(Δx)=0 and Im(Δx)<0}
x = real or complex values
$n,\Delta x,x_i,x_p Kr,k_r$ = real or complex constants
r = 1,2,3 , The x locus segment designations
$K_r$ = constant of integration for Re(Δx)>0 or {Re(Δx) = 0 and Im(Δx) > 0}
$k_r$ = constant of integration for Re(Δx)<0 or {Re(Δx) = 0 and Im(Δx) < 0}
$n \neq 1$

The x locus line is the straight line in the complex plane through all of the plotted x points.

Concepts

C1) The x locus and the x locus line in the complex plane

C2) Function snap – An exceeding abrupt transition between two successive values, often from a correct to an incorrect value

C3) Snap Hypothesis – lnd(n,Δx,x) n≠1 Series snap will occur, if it occurs at all, only at an x locus transition across a complex plane axis and, except where series snap occurs, the lnd(n,Δx,x) n≠1 Series constant of integration will not change.

C4) Over all x, the series function, $lnd_f(n,\Delta x,x)$, can have as many as six constants of integration, $K_1,K_2,K_3,k_1,k_2,k_3$.

Diagram - x locus and x locus line in the complex plane

$$\Delta x \sum_{\substack{\Delta x \\ x=x_i}}^{\pm\infty} \frac{1}{x^n} = \int_{\Delta x\, x_i}^{\pm\infty} \frac{1}{x^n}\Delta x = lnd(n,\Delta x,x_i) \approx lnd_f(n,\Delta x,x_i) + \begin{cases} K_r \\ k_r \end{cases} ,\quad x_i = \text{an x locus value} ,\ Re(n)>1$$

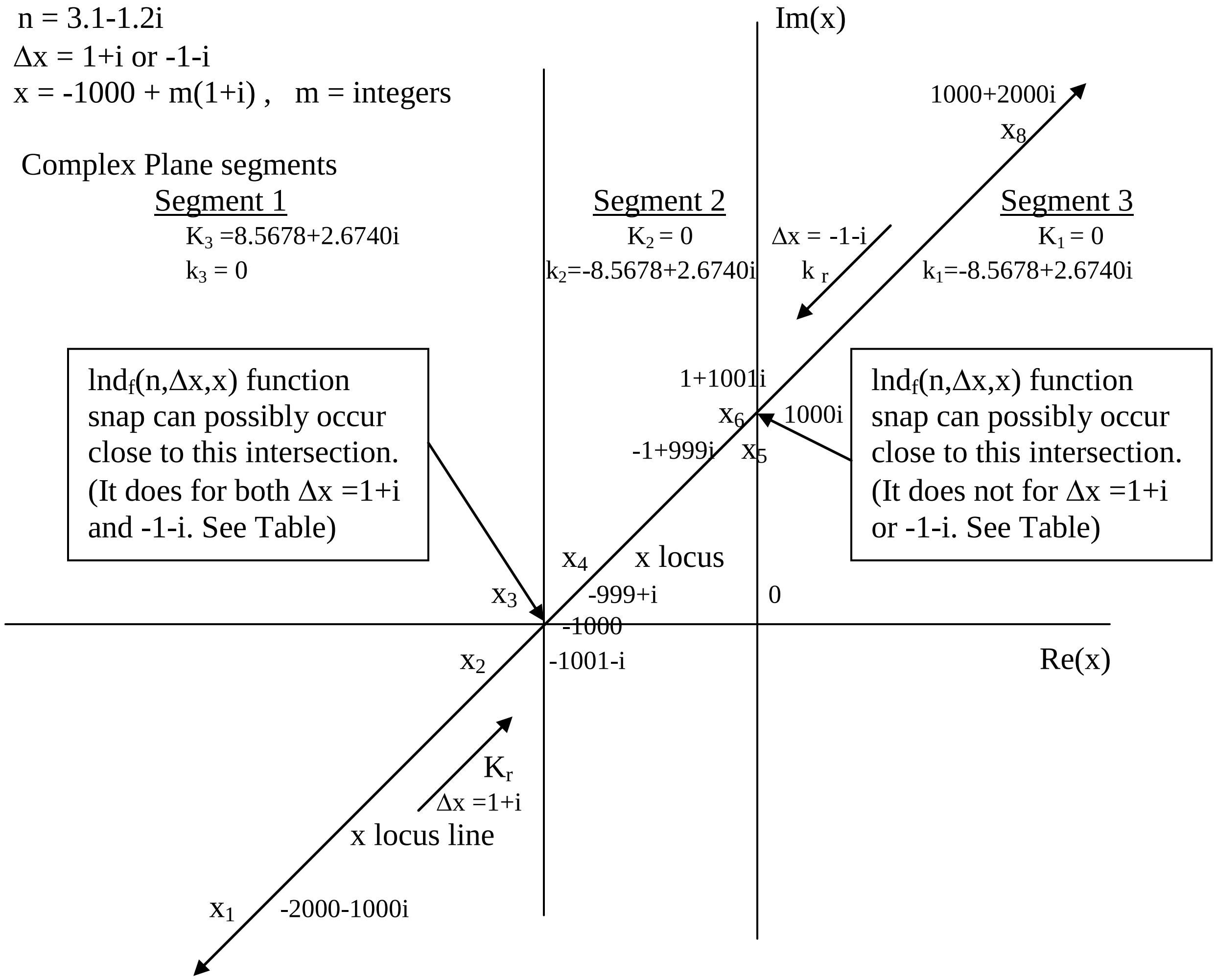

$K_1,K_2,K_3$ are the complex plane segment $lnd_f(n,\Delta x,x)$ constants of integration for those summations moving upward along the x locus line where $\Delta x = 1+i$.
$k_1,k_2,k_3$ are the complex plane segment $lnd_f(n,\Delta x,x)$ constants of integration for those summations moving downward along the x locus line where $\Delta x = -1-i$.

Note – The $lnd_f(n,\Delta x,x)$ constants of integration specified in the above diagram were calculated using the $lnd(n,\Delta x,x)$ calculation program, LNDX. The programming code for this program is shown in the Calculation Programs section at the end of the Appendix.

To demonstrate the validity of Equations E1 thru E12 and Concepts C1 thru C4 calculations will be made and verified using other mathematical means.

**Demonstration #1**

Demonstate the validity of Equations E1 thru E9 and Concepts C1 thru C4.

Equations E4 thru E9 and Concepts C1 Thru C4 have been incorporated into the UBASIC lnd(n,Δx,x) calculation program, LNDX. The program, LNDX, will be used to calculate equations, E1 thru E3 using n = 3.1-1.2i, Δx = 1+i or -1-i, the x locus and x locus line shown in the previous diagram, and various values for $x_i$. Also, a computer calculation check of the results of the LNDX program will be preformed to verify accuracy. All values of x will lie on the x locus line shown in the previous diagram.

Comment - The computer calculaton check of the LNDX program sums the terms of the summation being calculated one by one. To sum to infinity in this way and obtain a high precision result is not practical. An excessive amount of computer time would be necessary. The checks performed below verify the results of the LNDX program from four to seven digits of accuracy. This is sufficient to verify the functionality of the LNDX computer program. It is known that the LNDX program provides high precision results.

1) Evaluate the summation, $\displaystyle\sum_{x=\mp\infty}^{\pm\infty}{}_{1+i}\,\frac{1}{x^{3.1-1.2i}}$ , using Eq E3

where

$n = 3.1-1.2i$

$\Delta x = 1+i$

$x_i = -989+11i$

$x = x_i + m\Delta x$ , m = integers

From Eq E3

$$\sum_{x=\mp\infty}^{\pm\infty}{}_{\Delta x}\,\frac{1}{x^{n}} = \frac{1}{\Delta x}\,[\mathrm{lnd}(n,\Delta x,x_i) - \mathrm{lnd}(n,-\Delta x,x_i-\Delta x)]$$

Substituting values

$$\sum_{x=\mp\infty}^{\pm\infty}{}_{1+i}\,\frac{1}{x^{3.1-1.2i}} = \frac{1}{1+i}\,[\mathrm{lnd}(3.1-1.2i,1+i,-989+11i) - \mathrm{lnd}(3.1-1.2i,-1-i,-990+10i]$$

Using the LNDX program to calculate the lnd(n,Δx,x) functions

$$\sum_{1+i \; x=\mp\infty}^{\pm\infty} \frac{1}{x^{3.1-1.2i}} = \frac{1}{1+i}[8.562801571620273 - 2.670413187174567i]x10^{-6}$$

$$\sum_{1+i \; x=\mp\infty}^{\pm\infty} \frac{1}{x^{3.1-1.2i}} = [2.946194192222852 - 5.616607379397420i]x10^{-6}$$

checking

For $x = -989+11i + m(1+i)$ , m = the integers from –100,000 to +100,000

$$\sum_{1+i \; x=-100989-99989i}^{99011+100011i} \frac{1}{x^{3.1-1.2i}} = [2.94626 - 5.61660i]x10^{-6} \qquad \text{Good check}$$

2) Evaluate the summation, $\sum_{1+i \; x=1102-102i}^{+\infty} \frac{1}{x^{3.1-1.2i}}$ , using Eq E2

where

$n = 3.1-1.2i$

$\Delta x = 1+i$

$x_i = -1102-102i$

$x = x_i + m\Delta x$ , $m = 0,1,2,3,\ldots$

From Eq E2

$$\Delta x \sum_{\Delta x \; x=x_i}^{\pm\infty} \frac{1}{x^n} = \frac{1}{\Delta x}\,\mathrm{lnd}(n,\Delta x,x_i)$$

Substituting values

$$\sum_{1+i \; x=-1102-102i}^{+\infty} \frac{1}{x^{3.1-1.2i}} = \frac{1}{1+i}[\mathrm{lnd}(3.1-1.2i,1+i,-1102-102i)]$$

Using the LNDX program to calculate the lnd(n,Δx,x) function

$$\sum_{x=-1102-102i}^{+\infty}{}_{1+i}\ \frac{1}{x^{3.1-1.2i}} = \frac{1}{1+i}\,[2.520092396540544-2.947396943161950i]\mathrm{x}10^{-6}$$

$$\sum_{x=-1102-102i}^{+\infty}{}_{1+i}\ \frac{1}{x^{3.1-1.2i}} = [1.112676351112174-1.407416045428370i]\mathrm{x}10^{-6}$$

checking

For $x = -1102-102i + m(1+i)$ , $m$ = the integers from 0 to 200,000

$$\sum_{x=-1102-102i}^{198898+199898i}{}_{1+i}\ \frac{1}{x^{3.1-1.2i}} = [1.112676 - 1.407416i]\mathrm{x}10^{-6} \qquad \text{Good check}$$

3) Evaluate the summation, $\sum_{x=1000i}^{-\infty}{}_{-1-i}\ \frac{1}{x^{3.1-1.2i}}$ , using Eq E2

where
$n = 3.1-1.2i$
$\Delta x = -1-i$
$x_i = 1000i$
$x = x_i + m\Delta x$ , $m = 0,1,2,3,\ldots$

From Eq E2

$$\Delta x \sum_{x=x_i}^{\pm\infty}{}_{\Delta x}\ \frac{1}{x^n} = \frac{1}{\Delta x}\,\mathrm{lnd}(n,\Delta x,x_i)$$

Substituting values

$$\sum_{x=1000i}^{-\infty}{}_{-1-i}\ \frac{1}{x^{3.1-1.2i}} = \frac{1}{1+i}\,[\mathrm{lnd}(3.1-1.2i,-1-i,1000i)]$$

Using the LNDX program to calculate the lnd(n,Δx,x) function

$$ {}_{-1-i}\sum_{x=1000i}^{-\infty} \frac{1}{x^{3.1-1.2i}} = \frac{1}{1+i}\,[-8.540262768319782+2.648481133711169i]\mathrm{x}10^{-6} $$

$$ {}_{-1-i}\sum_{x=1000i}^{-\infty} \frac{1}{x^{3.1-1.2i}} = [2.945890817304306+5.594371951015475i]\mathrm{x}10^{-6} $$

checking

For $x = 1000i + m(-1-i)$ , m = the integers from 0 to 200,000

$$ {}_{-1-i}\sum_{x=1000i}^{-200000-199000i} \frac{1}{x^{3.1-1.2i}} = [2.9459 - 5.5943i]\mathrm{x}10^{-6} \qquad \text{Good check} $$

4) Evaluate the summation, $\displaystyle {}_{1+i}\sum_{x=-1001-i}^{1+1001i} \frac{1}{x^{3.1-1.2i}}$ , using Eq E1

where

$n = 3.1-1.2i$

$\Delta x = 1+i$

$x_1 = -1001-i$

$x_2 = 1+1001i$

$x = x_1 + m\Delta x$ , m = the integers from 0 to 1002 (i.e. $\frac{x_2-x_1}{\Delta x}$ )

From Eq E1

$$ \Delta x \, {}_{\Delta x}\sum_{x=x_1}^{x_2} \frac{1}{x^n} = -\mathrm{lnd}(n,\Delta x,x)\,\Big|_{x_1}^{x_2+\Delta x} = -\mathrm{lnd}(n,\Delta x,x_2+\Delta x) + \mathrm{lnd}(n,\Delta x,x_1)] $$

Substituting values

$$ {}_{1+i}\sum_{x=-1001-i}^{1+1001i} \frac{1}{x^{3.1-1.2i}} = \frac{1}{1+i}\,[-\mathrm{lnd}(3.1-1.2i,1+i,2+1002i) + \mathrm{lnd}(3.1-1.2i,1+i,-1001-i] $$

Using the LNDX program to calculate the lnd(n,Δx,x) functions

$$\sum_{1+i\,x=-1001-i}^{1+1001i} \frac{1}{x^{3.1-1.2i}} = \frac{1}{1+i}[.5152442710026848+2.450653840064823]x10^{-8}$$

$$\sum_{1+i\,x=-1001-i}^{1+1001i} \frac{1}{x^{3.1-1.2i}} = [1.482949055533754+.9677047845310694]x10^{-8}$$

checking

$$\sum_{1+i\,x=-1001-i}^{1+1001i} \frac{1}{x^{3.1-1.2i}} = [1.482949055533754+.9677047845310694]x10^{-8} \qquad \text{Good check}$$

The LNDX program has calculated the values of the lnd(n,Δx,x) functions correctly to good accuracy. The validity of Equations E1 thru E9, Concepts C1 thru C4 and the LNDX program itself has been demonstrated.

<u>Demonstration #2</u>
Demonstate the validity of Equation E10

$K_3 = -k_1$.

Note columns G and H in the following table of calculated values relating to the diagram x locus. The value of $K_3$ = 8.5628-2.6704i and the value of $k_1$ = -8.5628-2.6704i , $K_3 = - k_1$. The validity of Equation E10 has been demonstrated.

<u>Demonstration #3</u>
Demonstate the validity of Equation E11

lnd(n,Δx,x) – lnd(n,-Δx,x-Δx) = $\pm(K_r - k_r)$

where

r = 1,2,3

\+ for Re(Δx)>0 or {Re(Δx)=0 and Im(Δx)>0}

– is for Re(Δx)<0 or {Re(Δx)=0 and Im(Δx)<0}

n,Δx,$K_r$,$k_r$ = real or complex constants

Note columns C, E and F in the following table of calculated values relating to the diagram x locus. The first half of columns C, E, and F represents the condition where Δx = 1+i and the second half of columns C,E, and F represents the condition where Δx = -(1+i). Where Re(1+i)>0, $+(K_r - k_r)$, is seen to equal lnd(n,Δx,x) – lnd(n,-Δx,x-Δx) and where Re(-1-i)<0, $-(K_r - k_r)$ is seen to equal lnd(n,Δx,x) – lnd(n,-Δx,x-Δx). The validity of Equation E11 has been demonstrated.

Comment - From Eq E11 and Eq E12 it can be shown that for all points on an x locus line with a specified value of n, the value of $K_r - k_r$ is the same. This can be seen by observing columns E and F in the following table.

Demonstration #4

Demonstate the validity of Equation E3.

$$\Delta x \sum_{\substack{\Delta x \\ x=\mp\infty}}^{\pm\infty} \frac{1}{x^n} = \mathrm{lnd}(n,\Delta x,x_i) - \mathrm{lnd}(n,-\Delta x,x_i-\Delta x) = \pm\,(K_r - k_r)$$

The validity of $\mathrm{lnd}(n,\Delta x,x_i) - \mathrm{lnd}(n,-\Delta x,x_i-\Delta x) = \pm\,(K_r - k_r)$ has been shown in Demonstration #3.

Demonstrate the validity of

$$\Delta x \sum_{\substack{\Delta x \\ x=\mp\infty}}^{\pm\infty} \frac{1}{x^n} = \mathrm{lnd}(n,\Delta x,x_i) - \mathrm{lnd}(n,-\Delta x,x_i-\Delta x)$$

where

$x = x_i + m\Delta x$ , m = integers , This is an x locus in the complex plane

$x_i$ = any value of x in the x locus

$\Delta x$ = x increment

$-\infty$ to $+\infty$ for $\mathrm{Re}(\Delta x)>0$ or $\{\mathrm{Re}(\Delta x)=0$ and $\mathrm{Im}(\Delta x)>0\}$

$+\infty$ to $-\infty$ for $\mathrm{Re}(\Delta x)<0$ or $\{\mathrm{Re}(\Delta x)=0$ and $\mathrm{Im}(\Delta x)<0\}$

x = real or complex values

$n,\Delta x,x_i,K_r,k_r$ = real or complex constants

$\mathrm{Re}(n)>1$

r = 1,2,3

Note columns B,C, D, F and J in the following table of calculated values relating to the diagram x locus. For the same values of n, Δx, x the calculated values in columns F and J are seen to be equal. The validity of Equation E3 has been demonstrated.

Demonstration #5

Demonstate the validity of Equation E12.

$$\mathrm{lnd}(n,\Delta x,x) - \mathrm{lnd}(n,-\Delta x,x-\Delta x) = -[\mathrm{lnd}(n,-\Delta x,x) - \mathrm{lnd}(n,\Delta x,x+\Delta x)]$$

where

$x = x_i + m\Delta x$ , m = integers , This is an x locus in the complex plane

$x_i$ = a value of x

$\Delta x$ = x increment

$n \neq 1$

$n,\Delta x$ = real or complex constants

Note columns B,C,D and F in the following table of calculated values relating to the diagram x locus. For a specified value of n and x it is observed that the calculated value of $\ln d(n,\Delta x,x) - \ln d(n,-\Delta x,x-\Delta x)$ reverses sign if the value of $\Delta x$ reverses sign. The validity of Equation E3 has been demonstrated.

Demonstration #6

Demonstate the validity of Equation E13.

$$k_1 = \pm [-\Delta x \sum_{\Delta x \; x=\mp\infty}^{\pm\infty} \frac{1}{x^n} ] \; , \quad Re(n)>1$$

where

$x = x_i + m\Delta x$ , m = integers , This is an x locus in the complex plane

$x_i$ = a value of x

$\Delta x$ = x increment

+ and $-\infty$ to $+\infty$ for $Re(\Delta x)>0$ or $\{Re(\Delta x)=0$ and $Im(\Delta x)>0\}$

– and $+\infty$ to $-\infty$ for $Re(\Delta x)<0$ or $\{Re(\Delta x)=0$ and $Im(\Delta x)<0\}$

$n \neq 1$

$n,\Delta x,k_1$ = real or complex constants

$Re(n)>1$

$\Delta x = -1-i$

$Re(\Delta x)<0$

$$k_1 = \Delta x \sum_{\Delta x \; +\infty}^{-\infty} \frac{1}{x^n}$$

In the table, the calculated value in J14 is $\Delta x \sum_{\Delta x \; +\infty}^{-\infty} \frac{1}{x^n} = [-8.5628+2.6704i]x10^{-6}$.

Substituting this value into the above equation, $k_1 = [-8.5628+2.6704i]x10^{-6}$.

This value agrees with the $\ln d(n,\Delta x,x) - \ln d_f(n,\Delta x,x)$ calculated value for $k_1$ in H14.

The validity of Equation E13 has been demonstrated.

Demonstration #7

Demonstate the validity of Equation E14.

$$K_3 = \pm [\Delta x \sum_{\Delta x \; x=\mp\infty}^{\pm\infty} \frac{1}{x^n} ] \; , \quad Re(n)>1$$

where

$x = x_i + m\Delta x$ , m = integers , This is an x locus in the complex plane

$x_i$ = a value of x

$\Delta x$ = x increment
+ and $-\infty$ to $+\infty$ for $Re(\Delta x)>0$ or $\{Re(\Delta x)=0$ and $Im(\Delta x)>0\}$
– and $+\infty$ to $-\infty$ for $Re(\Delta x)<0$ or $\{Re(\Delta x)=0$ and $Im(\Delta x)<0\}$
$n \neq 1$
$n,\Delta x,K_3$ = real or complex constants

$Re(n)>1$
$\Delta x = 1+i$
$Re(\Delta x)>0$

$$K_3 = \Delta x \sum_{\Delta x \; -\infty}^{+\infty} \frac{1}{x^n}$$

In the table, the calculated value in J1 is $\Delta x \sum_{\Delta x \; -\infty}^{+\infty} \frac{1}{x^n} = [8.5628+2.6704i]x10^{-6}$.

Substituting this value into the above equation, $K_3 = [8.5628+2.6704i]x10^{-6}$.
This value agrees with the $lnd(n,\Delta x,x) - lnd_f(n,\Delta x,x)$ calculated value for $K_3$ in H1.
The validity of Equation E14 has been demonstrated.

<u>Demonstration #8</u>
Demonstate the validity of Equation E15.

$$K_2 - k_2 = \pm [\Delta x \sum_{\Delta x \; x=\mp\infty}^{\pm\infty} \frac{1}{x^n} ] , \quad Re(n)>1$$

where
$x = x_i + m\Delta x$ , m = integers , This is an x locus in the complex plane
$x_i$ = a value of x
$\Delta x$ = x increment
+ and $-\infty$ to $+\infty$ for $Re(\Delta x)>0$ or $\{Re(\Delta x)=0$ and $Im(\Delta x)>0\}$
– and $+\infty$ to $-\infty$ for $Re(\Delta x)<0$ or $\{Re(\Delta x)=0$ and $Im(\Delta x)<0\}$
$n \neq 1$
$n,\Delta x,K_2,k_2$ = real or complex constants

$Re(n)>1$
$\Delta x = 1+i$
$Re(\Delta x)>0$

$$K_2 - k_2 = \Delta x \sum_{\Delta x \; -\infty}^{+\infty} \frac{1}{x^n}$$

In the table, the calculated value in J5 is $\Delta x \sum_{\Delta x}^{+\infty}_{-\infty} \frac{1}{x^n} = [8.5628+2.6704i]x10^{-6}$.

Substituting this value into the above equation, $K_2 - k_2 = [8.5628+2.6704i]x10^{-6}$.
This value agrees with the $lnd(n,\Delta x,x) - lnd_f(n,\Delta x,x)$ calculated value for $K_2 - k_2$ in F5.
The validity of Equation E15 has been demonstrated.

Demonstration #9
Demonstate the validity of Equation E16.

$$K_r - k_r = \pm [\Delta x \sum_{\Delta x \; x=\mp\infty}^{\pm\infty} \frac{1}{x^n}] , \quad Re(n)>1$$

where
- $x = x_i + m\Delta x$ , m = integers , This is an x locus in the complex plane
- $x_i$ = a value of x
- $\Delta x$ = x increment
- + and $-\infty$ to $+\infty$ for $Re(\Delta x)>0$ or $\{Re(\Delta x)=0$ and $Im(\Delta x)>0\}$
- – and $+\infty$ to $-\infty$ for $Re(\Delta x)<0$ or $\{Re(\Delta x)=0$ and $Im(\Delta x)<0\}$
- $n \neq 1$
- $n,\Delta x$ = real or complex constants
- r = 1,2,3

$Re(n)<1$
$\Delta x = -1-i$
$Re(\Delta x)<0$

$$K_r - k_r = \Delta x \sum_{\Delta x \; +\infty}^{-\infty} \frac{1}{x^n}$$

In the table, the calculated values in J9 thru J16 are $\Delta x \sum_{\Delta x}^{+\infty}_{-\infty} \frac{1}{x^n} = [-8.5628+2.6704i]x10^{-6}$.

Substituting these values into the above equation, $K_r - k_r = [-8.5628+2.6704i]x10^{-6}$.
These values agree with the $lnd(n,\Delta x,x) - lnd_f(n,\Delta x,x)$ calculated values for $K_r - k_r$ in F9 thru F16.
The validity of Equation E16 has been demonstrated.

Demonstration #10

Demonstate the validity of Equation E17.

$$\left.\begin{matrix} K_r \\ k_r \end{matrix}\right\} = \Delta x \sum_{\Delta x \; x=x_i}^{\pm\infty} \frac{1}{x^n} - \mathrm{lnd}_f(n,\Delta x,x_i) \ , \quad \mathrm{Re}(n)>1$$

where

$x = x_i + m\Delta x$ ,    m = integers  ,   This is an x locus in the complex plane
$x_i$ = values of x
$x_i$ is within the x locus segment, r
$x_i$ to $+\infty$   for $\mathrm{Re}(\Delta x)>0$ or $\{\mathrm{Re}(\Delta x)=0$ and $\mathrm{Im}(\Delta x)>0\}$
$x_i$ to $-\infty$    for $\mathrm{Re}(\Delta x)<0$ or $\{\mathrm{Re}(\Delta x)=0$ and $\mathrm{Im}(\Delta x)<0\}$
$\Delta x$ = x increment
x = real or complex values
$n,\Delta x,x_i,K_r,k_r$ = real or complex constants
r = 1,2,3 ,   The x locus segment designations
$K_r$ = constant of integration for $\mathrm{Re}(\Delta x)>0$ or $\{\mathrm{Re}(\Delta x) = 0$ and $\mathrm{Im}(\Delta x) > 0\}$
$k_r$ = constant of integration for $\mathrm{Re}(\Delta x)<0$ or $\{\mathrm{Re}(\Delta x) = 0$ and $\mathrm{Im}(\Delta x) < 0\}$
$n \neq 1$
Accuracy increases rapidly as $\left|\frac{x}{\Delta x}\right|$ increases in value.

In the table, the calculated values in column K are $\left.\begin{matrix} K_r \\ k_r \end{matrix}\right\} = \Delta x \sum_{\Delta x \; x=x_i}^{\pm\infty} \frac{1}{x^n} - \mathrm{lnd}_f(n,\Delta x,x_i)$ where $x_i$ are the x values in column D. In the table, the calculated values in column H are $\left.\begin{matrix} K_r \\ k_r \end{matrix}\right\} = \mathrm{lnd}(n,\Delta x,x) - \mathrm{lnd}_f(n,\Delta x,x)$ .
Note that the calculated values in column K and in column H are the same. The validity of Equation E17 has been demonstrated.

**Table - Calculations relating to the x locus diagram**

| Row No. | Pnt | n | Δx | x | ±($K_r$-$k_r$) lnd(n,Δx,x) –lnd(n,-Δx,x-Δx) x10$^{-6}$ | | $K_r$ or $k_r$ lnd(n,Δx,x) –$lnd_f$(n,Δx,x) x10$^{-6}$ | | ±($K_r$-$k_r$) $\Delta x \sum_{\Delta x,\, x=\mp\infty}^{\pm\infty} \frac{1}{x^n}$ $x=x_i+m\Delta x$ , m = integers; $x_i$ = x in column D x10$^{-6}$ | | $K_r$ or $k_r$ (See column G) $\Delta x \sum_{\Delta x,\, x=x_i}^{\pm\infty} \frac{1}{x^n} - lnd_f(n,\Delta x,x_i)$ Re(n)>1; $x_i$ = x in column D x10$^{-6}$ |
|---|---|---|---|---|---|---|---|---|---|---|---|
| | A | B | C | D | E | F | G | H | I | J | K |
| 1 | $x_1$ | 3.1-1.2i | 1+i | -2000-1000i | $K_3$-$k_3$ | 8.5628-2.6704i | $K_3$ | 8.5628-2.6704i | $K_3$-$k_3$ | 8.5628-2.6704i | 8.5628-2.6704i |
| 2 | $x_2$ | 3.1-1.2i | 1+i | -1001-i | $K_3$-$k_3$ | 8.5628-2.6704i | $K_3$ | 8.5628-2.6704i | $K_3$-$k_3$ | 8.5628-2.6704i | 8.5628-2.6704i |
| 3 | $x_3$ | 3.1-1.2i | 1+i | -1000 | $K_3$-$k_3$ | 8.5628-2.6704i | $K_3$ | 8.5628-2.6704i | $K_3$-$k_3$ | 8.5628-2.6704i | 8.5628-2.6704i |
| 4 | $x_4$ | 3.1-1.2i | 1+i | -999+i | $K_2$-$k_2$ | 8.5628-2.6704i | $K_2$ | 0 snap | $K_2$-$k_2$ | 8.5628-2.6704i | 0 snap |
| 5 | $x_5$ | 3.1-1.2i | 1+i | -1+999i | $K_2$-$k_2$ | 8.5628-2.6704i | $K_2$ | 0 | $K_2$-$k_2$ | 8.5628-2.6704i | 0 |
| 6 | $x_6$ | 3.1-1.2i | 1+i | 1000i | $K_2$-$k_2$ | 8.5628-2.6704i | $K_2$ | 0 | $K_2$-$k_2$ | 8.5628-2.6704i | 0 |
| 7 | $x_7$ | 3.1-1.2i | 1+i | 1+1001i | $K_1$-$k_1$ | 8.5628-2.6704i | $K_1$ | 0 | $K_1$-$k_1$ | 8.5628-2.6704i | 0 |
| 8 | $x_8$ | 3.1-1.2i | 1+i | 1000+2000i | $K_1$-$k_1$ | 8.5628-2.6704i | $K_1$ | 0 | $K_1$-$k_1$ | 8.5628-2.6704i | 0 |
| 9 | $x_1$ | 3.1-1.2i | -1-i | -2000-1000i | $k_3$-$K_3$ | -8.5628+2.6704i | $k_3$ | 0 | $k_3$-$K_3$ | -8.5628+2.6704i | 0 |
| 10 | $x_2$ | 3.1-1.2i | -1-i | -1001-i | $k_3$-$K_3$ | -8.5628+2.6704i | $k_3$ | 0 snap | $k_3$-$K_3$ | -8.5628+2.6704i | 0 snap |
| 11 | $x_3$ | 3.1-1.2i | -1-i | -1000 | $k_2$-$K_2$ | -8.5628+2.6704i | $k_2$ | -8.5628+2.6704i | $k_2$-$K_2$ | -8.5628+2.6704i | -8.5628+2.6704i |
| 12 | $x_4$ | 3.1-1.2i | -1-i | -999+i | $k_2$-$K_2$ | -8.5628+2.6704i | $k_2$ | -8.5628+2.6704i | $k_2$-$K_2$ | -8.5628+2.6704i | -8.5628+2.6704i |
| 13 | $x_5$ | 3.1-1.2i | -1-i | -1+999i | $k_2$-$K_2$ | -8.5628+2.6704i | $k_2$ | -8.5628+2.6704i | $k_2$-$K_2$ | -8.5628+2.6704i | -8.5628+2.6704i |
| 14 | $x_6$ | 3.1-1.2i | -1-i | 1000i | $k_1$-$K_1$ | -8.5628+2.6704i | $k_1$ | -8.5628+2.6704i | $k_1$-$K_1$ | -8.5628+2.6704i | -8.5628+2.6704i |
| 15 | $x_7$ | 3.1-1.2i | -1-i | 1+10001i | $k_1$-$K_1$ | -8.5628+2.6704i | $k_1$ | -8.5628+2.6704i | $k_1$-$K_1$ | -8.5628+2.6704i | -8.5628+2.6704i |
| 16 | $x_8$ | 3.1-1.2i | -1-i | 1000+2000i | $k_1$-$K_1$ | -8.5628+2.6704i | $k_1$ | -8.5628+2.6704i | $k_1$-$K_1$ | -8.5628+2.6704i | -8.5628+2.6704i |

+ and $K_r$ are used when Re(Δx)>0 or {Re(Δx)=0 and Im(x)>0}

– and $k_r$ are used when Re(Δx)<0 or {Re(Δx)=0 and Im(x)<0}

r = 1,2,3

Thus, from the calculations of this example, the validity of Equations E1 thru E17, the Concepts C1 thru C4 and the functionality of the LNDX lnd(n,Δx,x) function calculation program have been demonstrated.

Comment - The LNDX program calculates but does not print out the $lnd_f$(n,Δx,x) function and the lnd(n,Δx,x) function constants of integration, $K_1$,$K_2$,$K_3$,$k_1$,$k_2$,$k_3$. To obtain these values, some additional code was added to the LNDX program to print them out. This modified LNDX program is called FINDK.

**Example 7.22** Finding the derivative, $D_{\Delta x}[\ln_{\Delta x}{}^2 x]$, using the Discrete Function Chain Rule and the Discrete Derivative of the Product of two Functions Equation

Solution #1

Find the discrete derivative of the discrete function, $\ln_{\Delta x}{}^2 x$, using the Discrete Function Chain Rule

The Discrete Function Chain Rule is as follows:

$$D_{\Delta x}\left[F(v)\,\Big|_{v=g_{\Delta x}(x)}\right] = \left[D_{\Delta v}F(v)\right]\Big|_{\substack{v=g_{\Delta x}(x)\\ \Delta v = \Delta x D_{\Delta x} g_{\Delta x}(x)}} \left[D_{\Delta x}g_{\Delta x}(x)\right] \qquad 1)$$

where
- $F(v)$ = function of v
- $g_{\Delta x}(x)$ = discrete Interval Calculus function of x
- $v = g_{\Delta x}(x)$
- $\Delta x, \Delta v$ = interval increments
- $\Delta v = \Delta x D_{\Delta x} g_{\Delta x}(x)$
- $D_{\Delta x}$ = the discrete derivative with respect to $\Delta x$

Let

$F(v) = v^2$

$D_{\Delta v}v^2 = \dfrac{(v+\Delta v)^2 - v^2}{\Delta v} = \dfrac{v^2+2v\Delta v+\Delta v^2-v^2}{\Delta v} = 2v+\Delta v$

$v = g_{\Delta x}(x) = \text{lnd}(1,\Delta x,x) \equiv \ln_{\Delta x}x$

$\Delta v = \Delta x D_{\Delta x} g_{\Delta x}(x) = \Delta x D_{\Delta x}\ln_{\Delta x}x = \Delta x(\frac{1}{x}) = \dfrac{\Delta x}{x}$

$D_{\Delta x}g_{\Delta x}(x) = D_{\Delta x}\ln_{\Delta x}x = \dfrac{1}{x}$

Substituting into Eq 1

$$D_{\Delta x}[\ln_{\Delta x}{}^2 x] = [2v+\Delta v]\Big|_{\substack{v=\ln_{\Delta x}x\\ \Delta v = \frac{\Delta x}{x}}} \left[\frac{1}{x}\right] \qquad 2)$$

$$D_{\Delta x}[\ln_{\Delta x}{}^2 x] = [2\ln_{\Delta x}x + \frac{\Delta x}{x}]\,\frac{1}{x} = [\ln_{\Delta x}x + (\ln_{\Delta x}x + \frac{\Delta x}{x})]\,\frac{1}{x} \qquad 3)$$

But

$$\ln_{\Delta x}x + \frac{\Delta x}{x} = \ln_{\Delta x}(x+\Delta x) \qquad 4)$$

Substituting Eq 4 into Eq 3

$$\mathbf{D_{\Delta x}[\ln_{\Delta x}{}^2 x] = [\ln_{\Delta x}x + \ln_{\Delta x}(x+\Delta x)]\,\frac{1}{x}} \qquad 5)$$

Note

$$_{\Delta x}\int_{\Delta x}^{x}\frac{1}{x}\Delta x = \Delta x\,_{\Delta x}\sum_{x=\Delta x}^{x-\Delta x}\frac{1}{x}\Delta x = \ln_{\Delta x}x\,\Big|_{\Delta x}^{x} = \ln_{\Delta x}x$$

$$_{\Delta x}\int_{\Delta x}^{x+\Delta x}\frac{1}{x}\Delta x = \Delta x\,_{\Delta x}\sum_{x=\Delta x}^{x}\frac{1}{x}\Delta x = \ln_{\Delta x}x\,\Big|_{\Delta x}^{x+\Delta x} = \ln_{\Delta x}(x+\Delta x)$$

$$\ln_{\Delta x}(x+\Delta x) - \ln_{\Delta x}x = \Delta x\,_{\Delta x}\sum_{x=\Delta x}^{x}\frac{1}{x}\Delta x - \Delta x\,_{\Delta x}\sum_{x=\Delta x}^{x-\Delta x}\frac{1}{x}\Delta x = \frac{\Delta x}{x}$$

$$\ln_{\Delta x}x + \frac{\Delta x}{x} = \ln_{\Delta x}(x+\Delta x)$$

This is the proof of Eq 4

Solution #2

Find the discrete derivative of the discrete function, $\ln_{\Delta x}{}^{2}x$, using the Discrete Derivative of the Product of two Functions Equation.

The Discrete Derivative of the Product of two Functions Equation is:

$$D_{\Delta x}[u(x)v(x)] = D_{\Delta x}u(x)v(x) + D_{\Delta x}v(x)u(x+\Delta x) \qquad 6)$$

Let

$u(x) = \ln_{\Delta x}x$

$v(x) = \ln_{\Delta x}x$

$D_{\Delta x}u(x) = \frac{1}{x}$

$D_{\Delta x}u(x) = \frac{1}{x}$

Substituting into Eq 6

$$D_{\Delta x}[\ln_{\Delta x}{}^{2}x] = \frac{1}{x}\ln_{\Delta x}x + \frac{1}{x}\ln_{\Delta x}(x+\Delta x) \qquad 6)$$

$$\mathbf{D_{\Delta x}[\ln_{\Delta x}{}^{2}x] = [\ln_{\Delta x}x + \ln_{\Delta x}(x+\Delta x)]\,\frac{1}{x}} \qquad 7)$$

Both methods for finding $D_{\Delta x}\ln_{\Delta x}{}^{2}x$ are seen to obtain the same results

# CHAPTER 8

# The General Zeta Function, the lnd(n,Δx,x) Function, and other Related Functions

# Chapter 8

# The General Zeta Function, the lnd(n,Δx,x) Function, and other Related Functions

## Section 8.1: Origin of the General Zeta Function

In Chapter 2, the so-called lnd(n,Δx,x) function was derived. The reason for its derivation was to resolve a difficulty encountered in the development of Interval Calculus. The difficulty involved the evaluation of the discrete integral, ${}_{\Delta x}\int_{x_1}^{x_2} \frac{1}{x^n}\,\Delta x$ , which is equal to the summation, $\Delta x \sum_{\Delta x\; x=x_1}^{x_2-\Delta x} \frac{1}{x^n}$ .

In Calculus where the x interval, Δx, is infinitessimal, the integral, $\int_{x_1}^{x_2} \frac{1}{x^n}\,dx$, evaluation functions are well known. For n=1 the integral evaluation function is the natural logarithm where $\int_{x_1}^{x_2} \frac{1}{x}\,dx = \ln x \Big|_{x_1}^{x_2}$.

For n≠1 the integral evaluation function is $\frac{x^{-n+1}}{-n+1}$ where $\int_{x_1}^{x_2} \frac{1}{x^n}\,dx = \frac{x^{-n+1}}{-n+1} \Big|_{x_1}^{x_2}$.

When Δx is not infinitessimal, neither of these two functions can be used.
What must be used, instead, is the discrete function, lnd(n,Δx,x), where

$$ {}_{\Delta x}\int_{x_1}^{x_2} \frac{1}{x^n}\,\Delta x = \Delta x \sum_{\Delta x\; x=x_1}^{x_2-\Delta x} \frac{1}{x^n}\,\Delta x = \pm\, \mathrm{lnd}(n,\Delta x,x) \Big|_{x_1}^{x_2}, \quad - \text{ for } n \neq 1,\ + \text{ for } n = 1, $$

$x = x_1, x_1+\Delta x, x_1+2\Delta x, \ldots , x_2-\Delta x, x_2$ .

To evaluate the discrete integral, ${}_{\Delta x}\int_{x_1}^{x_2} \frac{1}{x^n}\,\Delta x$, the lnd(n,Δx,x) function was derived from its following definition:

$$ D_{\Delta x}\mathrm{lnd}(n,\Delta x,x) = \pm \frac{1}{x^n}, \quad - \text{ for } n \neq 1, \quad + \text{ for } n=1 \qquad (8.1\text{-}1) $$

where
- Δx = x interval
- x = real or complex variable
- $x_1, x_2, \Delta x, n$ = real or complex constants
- Any summation term where x = 0 is excluded

Integrating Eq 8.1-1 using discrete Interval Calculus Integration

$$_{\Delta x}\int_{x_1}^{x_2} \frac{1}{x^n}\,\Delta x = \Delta x \,{}_{\Delta x}\sum_{x=x_1}^{x_2-\Delta x} \frac{1}{x^n} = \pm\, \mathrm{lnd}(n,\Delta x,x)\Big|_{x_1}^{x_2}, \quad -\text{ for } n \neq 1, \quad +\text{ for } n = 1 \tag{8.1-2}$$

where

$$_{\Delta x}\sum_{x=x_1}^{x_1-\Delta x} \frac{1}{x^n} = 0$$

$\Delta x$ = x interval

x = real or complex variable

$x_1, x_2, \Delta x, n$ = real or complex constants

Any summation term where x = 0 is excluded

As seen in the previous chapters the lnd(n,Δx,x) function properly evaluates the discrete integral,

$_{\Delta x}\int_{x_1}^{x_2} \frac{1}{x^n}\,\Delta x$ and the summation, $\Delta x \,{}_{\Delta x}\sum_{x=x_1}^{x_2} \frac{1}{x^n}$ .

From Eq 8.1-2

$$_{\Delta x}\int_{x_1}^{x_2} \frac{1}{x^n}\,\Delta x = \pm\, \mathrm{lnd}(n,\Delta x,x)\Big|_{x_1}^{x_2}, \quad -\text{ for } n \neq 1, \quad +\text{ for } n = 1 \tag{8.1-3}$$

and

$$_{\Delta x}\sum_{x=x_1}^{x_2} \frac{1}{x^n} = \pm \frac{1}{\Delta x}\, \mathrm{lnd}(n,\Delta x,x)\Big|_{x_1}^{x_2+\Delta x}, \quad -\text{ for } n \neq 1, \quad +\text{ for } n = 1 \tag{8.1-4}$$

In Chapter 2 the series for calculating the lnd(n,Δx,x) function is derived. The derived series is rewritten on the following page.

The lnd(n,Δx,x) Series is:

$$\mathrm{lnd}(n,\Delta x,x) \approx \left[\frac{1+\alpha(n)}{2}\right]\left[\ln\left(\frac{x}{\Delta x}-\frac{1}{2}\right)+\gamma\right]+\left[\frac{1-\alpha(n)}{2}\right]\left[\frac{1}{(n-1)(x-\frac{\Delta x}{2})^{n-1}}+K\right]$$

$$+\ \alpha(n)\sum_{m=1}^{\infty}\frac{\Gamma(n+2m-1)\left(\frac{\Delta x}{2}\right)^{2m}C_m}{\Gamma(n)(2m+1)!\left(x-\frac{\Delta x}{2}\right)^{n+2m-1}} \qquad (8.1\text{-}5)$$

Accuracy improves rapidly for increasing $|\frac{x}{\Delta x}|$

where

$\alpha(n) = \begin{cases} 1 & n = 1 \\ -1 & n \neq 1 \end{cases}$

n,Δx,x = real or complex values

K = constant of integration for n≠1

γ = constant of integration for n=1, Euler's Constant .577215664…

Δx = x increment

$C_m$ = Series constants, m = 1,2,3,…

$C_1 = +1$ $\qquad C_5 = +\frac{2555}{3}$

$C_2 = -\frac{7}{3}$ $\qquad C_6 = -\frac{1414477}{105}$

$C_3 = +\frac{31}{3}$ $\qquad C_7 = +286685$

$C_4 = -\frac{381}{5}$ $\qquad$ …

As mentioned in Chapter 2, the lnd(n,Δx,x) Series, though appearing relatively typical, is not. Interestingly, the convergence of the series is quite rapid, however, a series oddity named "snap" requires a significant amount of programming code to reliably maintain computational accuracy. The programming of the lnd(n,Δx,x) evaluation series of Eq 8.1-5 is described at length in Chapter 2.
The lnd(n,Δx,x) evaluation program, written in the UBASIC programming language and run on the Microsoft MS-DOS operating system, is named, LNDX. It is provided on a CD disk located in a pouch on the back cover of this paper.

After the development of Interval Calculus in Chapter 1 through Chapter 7 an interesting similarity was observed. It was realized that the summation in Eq 8.1-2 could be converted into a generalization of the Hurwitz Zeta Function. This generalization of the Hurwitz Zeta Function, called the General Zeta Function, is derived and described in the following section, Section 8.2. Section 8.2 also shows the relationship between the General Zeta Function, the Hurwitz Zeta Function, the Riemann Zeta Function, the Digamma Function, the Polygamma Functions, and the Interval Calculus lnd(n,Δx,x) Function.

## Section 8.2: Derivation of the General Zeta Function, the Hurwitz Zeta Function, the Riemann Zeta Function, the Digamma Function, and the Polygamma Functions in terms of the lnd(n,Δx,x) Function

The General Zeta Function

The General Zeta Function is defined in terms of the lnd(n,Δx,x) Function.

**The General Zeta Function is:**

$$\boldsymbol{\zeta(n,\Delta x,x) = \frac{1}{\Delta x}\,\mathrm{lnd}(n,\Delta x,x)} \qquad (8.2\text{-}1)$$

**This General Zeta Function evaluation equation is valid for all real or complex values of n, Δx, and x.**

From 8.1-2

$$\mathrm{lnd}(n,\Delta x,x)\Big|_{x_1}^{x_2} = \pm\,\Delta x \sum_{\Delta x\; x=x_1}^{x_2-\Delta x} \frac{1}{x^n} = \pm \int_{\Delta x\; x_1}^{x_2} \frac{1}{x^n}\,\Delta x\,, \quad -\text{ for } n \neq 1,\ +\text{ for } n = 1 \qquad (8.2\text{-}2)$$

where

$$\sum_{\Delta x\; x=x_1}^{x_1-\Delta x} \frac{1}{x^n} = 0$$

x = real or complex variable

$x_1,x_2,\Delta x,n$ = real or complex constants

Δx = x interval

Any summation term where x = 0 is excluded

The lnd(n,Δx,x) function represents an area as shown from the discrete integral above.

From Eq 8.2-1 and Eq 8.2-2

$$\zeta(n,\Delta x,x)\Big|_{x_1}^{x_2} = \pm \sum_{\Delta x\; x=x_1}^{x_2-\Delta x} \frac{1}{x^n} = \pm \frac{1}{\Delta x} \int_{\Delta x\; x_1}^{x_2} \frac{1}{x^n}\,\Delta x \;= \frac{1}{\Delta x}\,\text{lnd}(n,\Delta x,x)\Big|_{x_1}^{x_2}, \;-\text{ for } n \neq 1,\; +\text{ for } n = 1 \quad (8.2\text{-}3)$$

where

$$\sum_{\Delta x\; x=x_1}^{x_1-\Delta x} \frac{1}{x^n} = 0$$

x = real or complex variable

$x_1,x_2,\Delta x,n$ = real or complex constants

$\Delta x$ = x interval

Any summation term where x = 0 is excluded

The Genral Zeta Function, $\zeta(n,\Delta x,x)$, represents a summation as shown from the summation of Eq 8.2-3.

From Eq 8.2-3

**The General Zeta Function Summation Between Finite Limits Equation is:**

$$\sum_{\Delta x\; x=x_1}^{x_2-\Delta x} \frac{1}{x^n} = \frac{1}{\Delta x} \int_{\Delta x\; x_1}^{x_2} \frac{1}{x^n}\,\Delta x \;= \pm\,\zeta(n,\Delta x,x)\Big|_{x_1}^{x_2} = \pm \frac{1}{\Delta x}\,\text{lnd}(n,\Delta x,x)\Big|_{x_1}^{x_2}, \;-\text{ for } n \neq 1,\; +\text{ for } n = 1 \quad (8.2\text{-}4)$$

**where**

$$\sum_{\Delta x\; x=x_1}^{x_1-\Delta x} \frac{1}{x^n} = 0$$

**x = real or complex variable**

**$x_1,x_2,\Delta x,n$ = real or complex constants**

**$\Delta x$ = x interval**

**Any summation term where x = 0 is excluded**

**or**

$$\zeta(n,\Delta x,x)\Big|_{x_1}^{x_2} = \frac{1}{\Delta x}\,\text{lnd}(n,\Delta x,x)\Big|_{x_1}^{x_2} \pm \sum_{\Delta x\; x=x_1}^{x_2-\Delta x} \frac{1}{x^n} = \pm \frac{1}{\Delta x} \int_{\Delta x\; x_1}^{x_2} \frac{1}{x^n}\,\Delta x \;,\; -\text{ for } n \neq 1,\; +\text{ for } n = 1 \quad (8.2\text{-}5)$$

**where**

$$\underset{\Delta x}{\sum_{x=x_1}^{x_1-\Delta x}} \frac{1}{x^n} = 0$$

**x = real or complex variable**

**$x_1$,$x_2$,Δx,n = real or complex constants**

**Δx = x interval**

**Any summation term where x = 0 is excluded**

Derive the General Zeta Function Definite Discrete Integral Between Finite Limits Equation.

From Eq 8.2-4

**The General Zeta Function Definite Discrete Integral Between Finite Limits is:**

$$\underset{\Delta x}{\int_{x_1}^{x_2}} \frac{1}{x^n}\,\Delta x = \Delta x \underset{\Delta x}{\sum_{x=x_1}^{x_2-\Delta x}} \frac{1}{x^n} = \pm\,\Delta x\zeta(n,\Delta x,x)\,\Big|_{x_1}^{x_2} = \pm\,\mathrm{lnd}(n,\Delta x,x)\,\Big|_{x_1}^{x_2}\,,\quad -\text{ for } n\neq 1,\ +\text{ for } n=1 \qquad (8.2\text{-}6)$$

**where**

$$\underset{\Delta x}{\sum_{x=x_1}^{x_1-\Delta x}} \frac{1}{x} = 0$$

**x = real or complex variable**

**$x_1$,$x_2$,Δx,n = real or complex constants**

**Δx = x increment**

**Any summation term where x = 0 is excluded**

Derive the General Zeta Function Indefinite Discrete Integral

From Eq 8.2-6

**The General Zeta Function Indefinite Discrete Integral is:**

$$\underset{\Delta x}{\int} \frac{1}{x^n}\,\Delta x = \pm\,\Delta x\zeta(n,\Delta x,x) + k(n,\Delta x) = \pm\,\mathrm{lnd}(n,\Delta x,x) + k(n,\Delta x)\ ,\quad -\text{ for } n\neq 1,\ +\text{ for } n=1 \qquad (8.2\text{-}7)$$

**where**

**x = real or complex variable**

**Δx,n = real or complex constants**

**Δx = x increment**

**k = constant of integration, a function of n,Δx**

Derive The General Zeta Function Summation to Infinity Equation.

From Eq 8.2-5

$n\neq1$

Let $x_2\rightarrow\pm\infty$, Re(n)>1
$+\infty$ for Re(Δx)>0 or [ Re(Δx)=0 and Im(Δx)>0 ]
$-\infty$ for Re(Δx)<0 or [ Re(Δx)=0 and Im(Δx)<0 ]

$$\zeta(n,\Delta x,x)\Big|_{x_1}^{\pm\infty} = -\sum_{\Delta x\, x=x_1}^{\pm\infty}\frac{1}{x^n} = \frac{1}{\Delta x}\,\mathrm{lnd}(n,\Delta x,x)\Big|_{x_1}^{\pm\infty} \qquad (8.2\text{-}8)$$

It is observed from the implementation of Eq 8.1-5 in the program, LNDX, that $\mathrm{lnd}(n,\Delta x,x)\big|_{x\rightarrow\pm\infty}\rightarrow 0$ when Re(n)>1. Then, from Eq 8.2-1, $\zeta(n,\Delta x,x)\big|_{x\rightarrow\pm\infty}\rightarrow 0$ when Re(n)>1

Simplifying Eq 8.2-8

$$\zeta(n,\Delta x,\pm\infty) - \zeta(n,\Delta x,x_1) = -\sum_{\Delta x\, x=x_1}^{\pm\infty}\frac{1}{x^n} = \frac{1}{\Delta x}\,\mathrm{lnd}(n,\Delta x,\pm\infty) - \frac{1}{\Delta x}\,\mathrm{lnd}(n,\Delta x,x_1) \qquad (8.2\text{-}9)$$

$$-\zeta(n,\Delta x,x_1) = -\sum_{\Delta x\, x=x_1}^{\pm\infty}\frac{1}{x^n} = -\frac{1}{\Delta x}\,\mathrm{lnd}(n,\Delta x,x_1) \qquad (8.2\text{-}10)$$

From Eq 8.2-10

**The General Zeta Function Summation to Infinity Equation is:**

$$\boldsymbol{\zeta(n,\Delta x,x_i) = \sum_{\Delta x\, x=x_i}^{\pm\infty}\frac{1}{x^n} = \frac{1}{\Delta x}\,\mathrm{lnd}(n,\Delta x,x_i)\,,\quad \mathrm{Re}(n)>1} \qquad (8.2\text{-}11)$$

**where**
**x = real or complex variable**
**$x_i$,n,Δx = real or complex constants**
**Δx = x increment**
**$+\infty$ for Re(Δx)>0 or [ Re(Δx)=0 and Im(Δx)>0 ]**
**$-\infty$ for Re(Δx)<0 or [ Re(Δx)=0 and Im(Δx)<0 ]**
**Any summation term where x = 0 is excluded**

Derive the General Zeta Function Definite Discrete Integral to Infinity

$$\int_{\Delta x\, x_1}^{x_2} \frac{1}{x^n}\Delta x = \Delta x \sum_{\Delta x\, x=x_1}^{x_2-\Delta x} \frac{1}{x^n} \qquad (8.2\text{-}12)$$

for $x_2 \rightarrow \pm\infty$

$$\int_{\Delta x\, x_1}^{\pm\infty} \frac{1}{x^n}\Delta x = \Delta x \sum_{\Delta x\, x=x_1}^{\pm\infty} \frac{1}{x^n} \qquad (8.2\text{-}13)$$

From Eq 8.2-11 and Eq 8.2-13

$$\int_{\Delta x\, x_1}^{\pm\infty} \frac{1}{x^n}\Delta x = \Delta x \sum_{\Delta x\, x=x_1}^{\pm\infty} \frac{1}{x^n} = \Delta x\,\zeta(n,\Delta x,x_i) = \mathrm{lnd}(n,\Delta x,x_i)\ , \quad \mathrm{Re}(n)>1 \qquad (8.2\text{-}14)$$

From Eq 8.2-14

**The General Zeta Function Definite Discrete Integral to Infinity is:**

$$\int_{\Delta x\, x=x_i}^{\pm\infty} \frac{1}{x^n}\Delta x = \Delta x \sum_{\Delta x\, x=x_i}^{\pm\infty} \frac{1}{x^n} = \Delta x\,\zeta(n,\Delta x,x_i) = \mathrm{lnd}(n,\Delta x,x_i)\ , \quad \mathrm{Re}(n)>1 \qquad (8.2\text{-}15)$$

**where**

**x = real or complex variable**

**$x_i$,n,Δx = real or complex constants**

**Δx = x increment**

**+∞ for Re(Δx)>0 or [ Re(Δx)=0 and Im(Δx)>0 ]**

**-∞ for Re(Δx)<0 or [ Re(Δx)=0 and Im(Δx)<0 ]**

**Any summation term where x = 0 is excluded**

Comment - The following equation, which relates the lnd(n,Δx,x) function to the Hurwitz Zeta Function, may be used.

$$\mathrm{lnd}(n,\Delta x,x) = \frac{1}{(\Delta x)^{n-1}}\,\mathrm{lnd}(n,\frac{x}{\Delta x})\ , \qquad (8.2\text{-}16)$$

x,Δx = real values with the condition that

x ≠ positive real value when

Δx = negative real value

n = real or complex value

or

x = Δx = real or complex values

n = real or complex value

or

$x,\Delta x$ = real or complex values
$n$ = integer

Derive the General Zeta Function discrete derivative.

The discrete derivative of the lnd(n,Δx,x) function is:

$$D_{\Delta x}\text{lnd}(n,\Delta x,x) = \pm\frac{1}{x^n}\,,\ \ -\text{ for } n\neq 1,\ +\text{ for } n=1 \qquad (8.2\text{-}17)$$

where
x = real or complex variable
n,Δx = real or complex constant
Δx = x increment

Comment -It is from Eq 8.2-17 that the lnd(n,Δx,x) function is derived.

Multiplying each side of Eq 8.2-17 by $\frac{1}{\Delta x}$

$$D_{\Delta x}\frac{1}{\Delta x}\text{lnd}(n,\Delta x,x) = \pm\frac{1}{\Delta x}\frac{1}{x^n}\,,\quad -\text{ for } n\neq 1,\ +\text{ for } n=1 \qquad (8.2\text{-}18)$$

Rewriting Eq 8.2-1

$$\zeta(n,\Delta x,x) = \frac{1}{\Delta x}\text{lnd}(n,\Delta x,x)$$

From Eq 8.2-1 and Eq 8.2-18

$$D_{\Delta x}\zeta(n,\Delta x,x) = \pm\frac{1}{\Delta x}\frac{1}{x^n}\,,\quad -\text{ for } n\neq 1,\ +\text{ for } n=1 \qquad (8.2\text{-}19)$$

From Eq 8.2-19

**The General Zeta Function Discrete Derivative is:**

$$\mathbf{D_{\Delta x}\zeta(n,\Delta x,x) = \pm\frac{1}{\Delta x}\frac{1}{x^n}\,,\quad -\text{ for } n\neq 1,\ +\text{ for } n=1} \qquad (8.2\text{-}20)$$

**where**
**x = real or complex variable**
**n,Δx = real or complex constants**
**Δx = x increment**

Derive the General Zeta Function term by term relationship.

Expanding Eq 8.2-20

$$\frac{\zeta(n,\Delta x,x+\Delta x) - \zeta(n,\Delta x,x)}{\Delta x} = \pm \frac{1}{\Delta x}\frac{1}{x^n}\ , \quad - \text{for } n\neq 1,\ + \text{for } n=1 \tag{8.2-21}$$

From Eq 8.2-21

**The General Zeta Function term to term relationship is:**

$$\boldsymbol{\zeta(n,\Delta x,x+\Delta x) - \zeta(n,\Delta x,x) = \pm \frac{1}{x^n}\ , \quad - \text{for } n\neq 1,\ + \text{for } n=1} \tag{8.2-22}$$

**where**

**x = real or complex variable**

**Δx,n = real or complex constants**

**Δx = x increment**

Derive the General Zeta Function Recursion Equation.

Rearranging Eq 8.2-22

**The General Zeta Function Recursion Equation is:**

$$\boldsymbol{\zeta(n,\Delta x,x+\Delta x) = \zeta(n,\Delta x,x) \pm \frac{1}{x^{\,n}}\ , \ - \text{for } n\neq 1,\ + \text{for } n=1,} \tag{8.2-23}$$

**where**

**x = real or complex variable**

**Δx,n = real or complex constants**

**Δx = x increment**

Equations relating to the Hurwitz Zeta Function, the Riemann Zeta Function, the Digammma Function, and the Polygamma Functions are derived from the General Zeta Function equations specified above.

When $n \neq 1$, the General Zeta Function is found to be related to the Hurwitz Zeta Function, the Riemann Zeta Function, and the Polygamma Functions.

<u>The Hurwitz Zeta Function</u>

The Hurwitz Zeta Function is a special case of the General Zeta Function where $\Delta x=1$ and $n\neq 1$.

Derive the relationship of the Hurwitz Zeta Function to the lnd(n,Δx,x) Function.

Rewriting Eq 8.2-1

$$\zeta(n,\Delta x,x) = \frac{1}{\Delta x}\,\text{lnd}(n,\Delta x,x)$$

From Eq 8.2-1 with $\Delta x = 1$, $n \neq 1$

$$\zeta(n,x) = \zeta(n,1,x) = \text{lnd}(n,1,x)\ ,\ n \neq 1 \tag{8.2-24}$$

From Eq 8.2-24

**The relationship of the Hurwitz Zeta Function to the lnd(n,Δx,x) function**

$$\boldsymbol{\zeta(n,x) = \zeta(n,1,x) = \text{lnd}(n,1,x)\ ,\ \ n \neq 1} \tag{8.2-25}$$

**where**
**x = real or complex variable**
**n = real or complex constant**

Derive the Hurwitz Zeta Function Summation Between Finite Limits Equation.

From the General Zeta Function Summation Between Finite Limits Equation, Eq 8.2-4 and Eq 8.2-5

Let $\Delta x = 1$
$n \neq 1$

**The Hurwitz Zeta Function Summation Between Finite Limits is:**

$$\sum_{1\ x=x_1}^{x_2} \frac{1}{x^n} = -\,\zeta(n,x)\Big|_{x_1}^{x_2+1} = -\,\zeta(n,1,x)\Big|_{x_1}^{x_2+1} = -\,\text{lnd}(n,1,x)\Big|_{x_1}^{x_2+1}\ ,\ \ n \neq 1 \tag{8.2-26}$$

**where**
**x = real or complex variable**
**$x_1, x_2$, n = real or complex constants**
**Any summation term where x = 0 is excluded**

**or**

$$\zeta(n,x)\Big|_{x_1}^{x_2} = -\sum_{1\ x=x_1}^{x_2-1} \frac{1}{x^n} = \zeta(n,1,x)\Big|_{x_1}^{x_2} = \text{lnd}(n,1,x)\Big|_{x_1}^{x_2}\ ,\ \ n\neq 1 \tag{8.2-27}$$

**where**

$$\sum_{\Delta x\ x=x_1}^{x_1-1} \frac{1}{x^n} = 0$$

**x = real or complex variable**
**$x_1, x_2$, n = real or complex constants**
**Any summation term where x = 0 is excluded**

Derive the Hurwitz Zeta Function Definite Discrete Integral Between Finite Limits Equation.

From the General Zeta Function Definite Discrete Integral Between Finite Limits Equation, Eq 8.2-6

Let $\Delta x = 1$
$n \neq 1$

**The Hurwitz Zeta Function Definite Discrete Integral Between Finite Limits is:**

$$ {}_1\!\int_{x_1}^{x_2} \frac{1}{x^n}\,\Delta x = {}_1\!\sum_{x=x_1}^{x_2-1} \frac{1}{x^n} = -\,\zeta(n,x)\Big|_{x_1}^{x_2} = -\,\zeta(n,1,x)\Big|_{x_1}^{x_2} = -\,\mathrm{lnd}(n,1,x)\Big|_{x_1}^{x_2}\,, \quad n\neq 1 \qquad (8.2\text{-}28)$$

**where**

$$ {}_1\!\sum_{x=x_1}^{x_1-1} \frac{1}{x} = 0 $$

**x = real or complex variable**
**$x_1$,$x_2$,n = real or complex constants**
**$\Delta x = 1$, x increment**
**Any summation term where x = 0 is excluded**

Derive the Hurwitz Zeta Function Indefinite Discrete Integral.

From Eq 8.2-28

**The Hurwitz Zeta Function Indefinite Discrete Integral is:**

$$ {}_1\!\int \frac{1}{x^n}\,\Delta x = -\,\zeta(n,x) + k(n) = -\,\zeta(n,1,x) + k(n) = -\,\mathrm{lnd}(n,1,x) + k(n)\,, \quad n\neq 1 \qquad (8.2\text{-}29)$$

**where**

**x = real or complex variable**
**n = real or complex constant**
**$\Delta x = 1$, x increment**
**k = constant of integration, a function of n**

Derive the Hurwitz Zeta Function Summation to Infinity Equation.

From the General Zeta Function Summation to Infinity Equation, Eq 8.2-11

Let $\Delta x = 1$
$n \neq 1$

**The Hurwitz Zeta Function Summation to Infinity is:**

$$\sum_{x=x_i}^{\infty} {}_1\frac{1}{x^n} = \zeta(n,x_i) = \zeta(n,1,x_i) = \operatorname{lnd}(n,1,x_i)\,, \quad \operatorname{Re}(n)>1 \qquad (8.2\text{-}30)$$

**where**

**x = real or complex variable**

**$x_i$,n = real or complex constants**

**$\Delta x$ = 1, x increment**

**Any summation term where x = 0 is excluded**

Derive the Hurwitz Zeta Function Definite Discrete Integral to Infinity Equation.

From the General Zeta Function Definite Discrete Integral to Infinity Equation, Eq 8.2-15

Let $\Delta x = 1$
$n \neq 1$

**The Hurwitz Zeta Function Definite Discrete Integral to Infinity is:**

$$\int_{x=x_i}^{\infty} {}_1\frac{1}{x^n}\,\Delta x = \sum_{x=x_i}^{\infty} {}_1\frac{1}{x^n} = \zeta(n,x_i) = \zeta(n,1,x_i) = \operatorname{lnd}(n,1,x_i)\,, \quad \operatorname{Re}(n)>1 \qquad (8.2\text{-}31)$$

**where**

**x = real or complex variable**

**$x_i$,n= real or complex constants**

**$\Delta x$ = 1, x increment**

**Any summation term where x = 0 is excluded**

Derive the Hurwitz Zeta Function Discrete Derivative.

From the General Zeta Function Discrete Derivative Equation, Eq 8.2-20

Let $\Delta x = 1$
$n \neq 1$

**The Hurwitz Zeta Function Discrete Derivative is:** (8.2-32)

$$D_1\zeta(n,x) = -\frac{1}{x^n}\,, \quad n\neq 1$$

**where**

**x = real or complex variable**

**n= real or complex constant**

Derive the Hurwitz Zeta Function term to term relationship.

From the General Zeta Function term to term relationship, Eq 8.2-22

Let $\Delta x = 1$
$n \neq 1$

**The Hurwitz Zeta Function term to term relationship is:**

$$\boldsymbol{\zeta(n,x+1) - \zeta(n,x) = -\frac{1}{x^n}, \quad n\neq 1} \tag{8.2-33}$$

**where**
**x = real or complex variable**
**n = real or complex constant**

Derive the Hurwitz Zeta Function Recursion Equation.

From the General Zeta Function Recursion Equation, Eq 8.2-23

Let $\Delta x = 1$
$n \neq 1$

**The Hurwitz Zeta Function Recursion Equation is:**

$$\boldsymbol{\zeta(n,x+1) = \zeta(n,x) - \frac{1}{x^n}, \quad n\neq 1} \tag{8.2-34}$$

**where**
**x = real or complex variable**
**n = real or complex constant**

The Riemann Zeta Function

Derive the relationship of the Riemann Zeta Function to the lnd(n,Δx,x) function

The General Zeta Function, Eq 8.2-1, is:

$$\zeta(n,\Delta x,x) = \frac{1}{\Delta x}\,\mathrm{lnd}(n,\Delta x,x)$$

From Eq 8.2-1

$\Delta x = 1$
$x = 1$
$n \neq 1$

$$\zeta(n) = \zeta(n,1,1) = \mathrm{lnd}(n,1,1)\,, \quad n \neq 1 \tag{8.2-35}$$

From Eq 8.2-35

**The relationship of the Riemann Zeta Function to the lnd(n,Δx,x) Function is:**

$$\zeta(n) = \zeta(n,1,1) = \mathrm{lnd}(n,1,1)\ ,\quad n \neq 1 \tag{8.2-36}$$

**where**
**n = real or complex constant**

Derive the Riemann Zeta Function Summation to Infinity Equation.

The Riemann Zeta Function Summation to Infinity is a special case of the General Zeta Function Summation to Infinity where $\Delta x=1$, $x_i=1$ and $n\neq1$.

From the General Zeta Function Summation to Infinity Equation, Eq 8.2-11

Let $\Delta x = 1$
$x_i = 1$
$n \neq 1$

**The Riemann Zeta Function Summation to Infinity**

$$\zeta(n) = {}_1\sum_{x=1}^{\infty} \frac{1}{x^n} = \zeta(n,1,1) = \mathrm{lnd}(n,1,1)\ ,\quad \mathrm{Re}(n)>1 \tag{8.2-37}$$

**where**
**x = real or complex variable**
**n = real or complex constant**

Derive the Riemann Zeta Function Definite Discrete Integral to Infinity Equation.

From the General Zeta Function Definite Discrete Integral to Infinity Equation, Eq 8.2-15

Let $\Delta x = 1$
$x_i = 1$
$n \neq 1$

**The Riemann Zeta Function Definite Discrete Integral to Infinity is:**

$${}_1\int_{x=1}^{\infty} \frac{1}{x^n}\,\Delta x = {}_1\sum_{x=1}^{\infty} \frac{1}{x^n} = \zeta(n) = \zeta(n,1,1) = \mathrm{lnd}(n,1,1)\ ,\quad \mathrm{Re}(n)>1 \tag{8.2-38}$$

**where**
**x = real or complex variable**
**n = real or complex constant**
**Δx = 1, x increment**

The Digamma Function

When n = 1, the General Zeta Function is found to be related to the Digamma Function.

Show the relationship of the Digamma Function to the General Zeta Function and to the lnd(n,Δx,x) Function.

Derive the Digamma term to term relationship

$$\Gamma(x+1) = x\Gamma(x)\ ,\ \ x \neq 0,-1,-2,-3,\dots\ ,\quad \Gamma(x) \text{ is infinite for } x = 0,-1,-2,-3,\dots \tag{8.2-39}$$

Taking the natural logarithm of both sides of Eq 8.2-39

$$\ln\Gamma(x+1) = \ln x + \ln\Gamma(x) \tag{8.2-40}$$

Taking the derivative of both sides of Eq 8.2-40

$$\frac{d}{dx}\ln\Gamma(x+1) = \frac{1}{x} + \frac{d}{dx}\ln\Gamma(x) \tag{8.2-41}$$

$$\frac{d}{dx}\ln\Gamma(x) = \frac{\Gamma(x)'}{\Gamma(x)} \tag{8.2-42}$$

Rearranging Eq 8.2-41

$$\frac{d}{dx}\ln\Gamma(x+1) - \frac{d}{dx}\ln\Gamma(x) = \frac{1}{x}\ , \tag{8.2-43}$$

Let

$$\psi(x) = \frac{d}{dx}\ln\Gamma(x) = \frac{\Gamma(x)'}{\Gamma(x)}\ ,\ \ \text{Digamma Function definition} \tag{8.2-44}$$

From Eq 8.2-43 and Eq 8.2-44

$$\psi(x+1) - \psi(x) = \frac{1}{x}\ ,\ \ \text{The Digamma term to term relationship} \tag{8.2-45}$$

From Eq 8.2-45

**The Digamma Function term to term relationship is:**

$$\boldsymbol{\psi(x+1) - \psi(x) = \frac{1}{x}} \tag{8.2-46}$$

**where**

**x = real or complex variable**

Derive the Digamma Recursion Equation

Rearranging Eq 8.2-46

**The Digamma Recursion Equation is:**

$$\boldsymbol{\psi(x+1) = \psi(x) + \frac{1}{x}}\,, \quad \textbf{The Digamma Recursion Equation} \qquad (8.2\text{-}47)$$

**where**
**x = real or complex variable**

Derive the Digamma Function discrete derivative

From Eq 8.2-46

**The Digamma Function discrete derivative**

$$\boldsymbol{D_1\psi(x) = \frac{1}{x}} \qquad (8.2\text{-}48)$$

**where**
**x = real or complex variable**
**n = real or complex constant**

Find a relationship of the Digamma Function to the lnd(n,Δx,x) Function.

Consider the derivative of the lnd(n,Δx,x) function

$$D_{\Delta x}\text{lnd}(n,\Delta x,x) = \pm\frac{1}{x^n}\ ,\ \ +\text{ for }n=1,\ \ -\text{ for }n\neq 1 \qquad (8.2\text{-}49)$$

Let n = 1
Δx = 1

$$D_1\text{lnd}(1,1,x) = \frac{1}{x} \qquad (8.2\text{-}50)$$

$$\text{lnd}(1,1,x+1) - \text{lnd}(1,1,x) = \frac{1}{x} \qquad (8.2\text{-}51)$$

More generally

$$[\text{lnd}(1,1,x+1) + k] - [\text{lnd}(1,1,x) + k] = \frac{1}{x} \qquad (8.2\text{-}52)$$

where
k = constant

Comparing Eq 8.2-52 to Eq 8.2-46

$$\psi(x) = \mathrm{lnd}(1,1,x) + k \tag{8.2-53}$$

Evaluate the constant, k

From the definition of $\psi(x)$ for x = 1,2,3,…

$$\psi(x) = {}_1\sum_{r=1}^{x-1} \frac{1}{r} - \gamma \tag{8.2-54}$$

where

$${}_1\sum_{r=1}^{0} \frac{1}{r} = 0$$

r = 1 ,2 ,3 ,…, x-1

$\gamma$ = Euler's Constant, .5772157…

$$\psi(1) = -\gamma \tag{8.2-55}$$

$$\mathrm{lnd}(1,1,1) = 0 \tag{8.2-56}$$

From Eq 8.2-53, Eq 8.2-55, and Eq 8.2-56

$$-\gamma = 0 + k \tag{8.2-57}$$

$$k = -\gamma \tag{8.2-58}$$

Substituting Eq 8.2-58 into Eq 8.2-53

$$\psi(x) = \mathrm{lnd}(1,1,x) - \gamma \tag{8.2-59}$$

Then from Eq 8.2-59

**A relationship of the Digamma Function to the lnd(n,Δx,x) Function is:**

$$\boldsymbol{\psi(x) = \mathrm{lnd}(1,1,x) - \gamma\,,\quad x \neq 0,-1,-2,-3,\ldots} \tag{8.2-60}$$

**where**

**x = real or complex values**

**$\gamma$ = Euler's Constant, .5772157…**

**$\psi(x)$ is infinite for x = 0,-1,-2,-3,…**

Derive another relationship of the Digamma Function to the lnd(n,Δx,x) Function

From Eq 8.2-60

$$\psi(\frac{x}{\Delta x}) = \mathrm{lnd}(1,1,\frac{x}{\Delta x}) - \gamma\,, \quad x \neq 0,-\Delta x,-2\Delta x,-3\Delta x,\ldots \tag{8.2-61}$$

$$\mathrm{lnd}(1,\Delta x,x) = \mathrm{lnd}(1,1,\frac{x}{\Delta x}) \tag{8.2-62}$$

$$\mathrm{lnd}(1,\Delta x,x) = \Delta x\, \zeta(1,\Delta x,x) \tag{8.2-63}$$

From Eq 8.2-61 thru 8.2-63

$$\psi(\frac{x}{\Delta x}) = \mathrm{lnd}(1,\Delta x,x) - \gamma = \mathrm{lnd}(1,1,\frac{x}{\Delta x}) - \gamma = \Delta x\zeta(1,\Delta x,x) - \gamma\,, \quad x \neq 0,-\Delta x,-2\Delta x,-3\Delta x,\ldots \tag{8.2-64}$$

Then from Eq 8.2-64

**Another relationship of the Digamma Function to the lnd(n,Δx,x) Function is:**

$$\boldsymbol{\psi(\frac{x}{\Delta x}) = \mathrm{lnd}(1,\Delta x,x) - \gamma = \mathrm{lnd}(1,1,\frac{x}{\Delta x}) - \gamma = \Delta x\zeta(1,\Delta x,x) - \gamma\,, \quad x \neq 0,-\Delta x,-2\Delta x,-3\Delta x,\ldots} \tag{8.2-65}$$

**where**
**x = real or complex values**
**γ = Euler's Constant, .5772157…**
**ψ(x) is infinite for x = 0,-1,-2,-3,…**
**Δx = x increment**
**ζ(1,Δx,x) is the N=1 Zeta Function**

**<u>Comments</u> – The condition that x ≠ 0,-1,-2,-3,… has been added to Eq 5.2-60 since the Digamma Function has first order poles at these values of x. At these values of x, where $\psi(x) = \frac{1}{0}$, the value of ψ(x) is infinite.**

**The lnd(1,1,x) function minus Euler's number is equal to the Digamma Function for all values of x except for x = 0,-1,-2,-3,… See Eq 5.2-60. For x = 0,-1,-2,-3,…, the lnd(1,1,x) function differs from the ψ(x) function. The equation used to calculate lnd(1,1,x) when x = 0,-1,-2,-3,… is lnd(1,1,x) = lnd(1,1,1–x). In addition, any summation calculated using the lnd(1,1,x) function, excludes any summation term with a division by zero. These differences were introduced into the lnd(1,1,x) function to make possible the summation of integers along the real axis of the complex plane.**

**For x = 0,-1,-2,-3,…, if the program, LNDX, is selected to find ψ(x) using Eq 8.2-60, the resulting calculated value for ψ(x) is ψ(x) = ψ(1–x).**

Derive the relationship of the N=1 Zeta Function to the Digamma Function.
.
From Eq 8.2-65

$$\psi(\frac{x}{\Delta x}) = \text{lnd}(1,\Delta x,x) - \gamma\,, \quad x \neq 0,-\Delta x,-2\Delta x,-3\Delta x,\ldots \tag{8.2-66}$$

Rearranging terms

$$\text{lnd}(1,\Delta x,x) = \psi(\frac{x}{\Delta x}) + \gamma\,, \quad x \neq 0,-\Delta x,-2\Delta x,-3\Delta x,\ldots \tag{8.2-67}$$

$$\frac{1}{\Delta x}\,\text{lnd}(1,\Delta x,x) = \frac{1}{\Delta x}\,[\psi(\frac{x}{\Delta x}) + \gamma]\,, \quad x \neq 0,-\Delta x,-2\Delta x,-3\Delta x,\ldots \tag{8.2-68}$$

$$\zeta(1,\Delta x,x) = \frac{1}{\Delta x}\,\text{lnd}(1,\Delta x,x) \tag{8.2-69}$$

Substituting Eq 8.2-69 into Eq 8.2-68

**The relationship of the N=1 Zeta Function to the Digamma Function is:**

$$\boldsymbol{\zeta(1,\Delta x,x) = \frac{1}{\Delta x}\,[\psi(\frac{x}{\Delta x}) + \gamma]\,, \quad x \neq 0,-\Delta x,-2\Delta x,-3\Delta x,\ldots} \tag{8.2-70}$$

**where**

**x = real or complex variable**

**Δx = x interval**

**γ = Euler's Constant, .5772157…**

Derive the Digamma Function Summation Between Finite Limits equation.

Rewriting the General Zeta Function Summation Between Finite Limits Equation, Eq 8.2-4

$$\sum_{\Delta x\; x=x_1}^{x_2-\Delta x} \frac{1}{x^n} = \pm\frac{1}{\Delta x}\,\text{lnd}(n,\Delta x,x)\Big|_{x_1}^{x_2}\,, \quad - \text{ for } n \neq 1\,, \; + \text{ for } n = 1 \tag{8.2-71}$$

where

$$\sum_{\Delta x\; x=x_1}^{x_1-\Delta x} \frac{1}{x^n} = 0$$

x = real or complex variable

$x_1, x_2, \Delta x, n$ = real or complex constants

Any summation term where x = 0 is excluded

Substitute n=1 into Eq 8.2-71

$$\sum_{\Delta x \; x=x_1}^{x_2-\Delta x} \frac{1}{x} = \frac{1}{\Delta x} \,\mathrm{lnd}(1,\Delta x,x)\Big|_{x_1}^{x_2} \tag{8.2-72}$$

Substituting Eq 8.2-67 into Eq 8.2-72

$$\sum_{\Delta x \; x=x_1}^{x_2-\Delta x} \frac{1}{x} = \frac{1}{\Delta x} \,\mathrm{lnd}(1,\Delta x,x)\Big|_{x_1}^{x_2} = \frac{1}{\Delta x}\left[\psi\left(\frac{x}{\Delta x}\right) + \gamma\right]\Big|_{x_1}^{x_2} = \frac{1}{\Delta x}\,\psi\left(\frac{x}{\Delta x}\right)\Big|_{x_1}^{x_2}, \quad x_1,x_2 \neq 0,-\Delta x,-2\Delta x,-3\Delta x,\ldots$$

(8.2-73)

<u>Comment</u> - Eq 8.2-72 is more general than Eq 8.2-73 since it has no restrictions on the value of x.

From Eq 8.2-73

$$\sum_{\Delta x \; x=x_1}^{x_2} \frac{1}{x} = \frac{1}{\Delta x}\,\psi\left(\frac{x}{\Delta x}\right)\Big|_{x_1}^{x_2+\Delta x} = \frac{1}{\Delta x} \,\mathrm{lnd}(1,\Delta x,x)\Big|_{x_1}^{x_2+\Delta x}, \quad x_1,x_2 \neq 0,-\Delta x,-2\Delta x,-3\Delta x,\ldots \tag{8.2-74}$$

where

$x_1 \neq 0,-\Delta x,-2\Delta x,-3\Delta x,\ldots$

$x_2 \neq -\Delta x,-2\Delta x,-3\Delta x,\ldots$

<u>Note</u> - From Eq 8.2-73 where $x_1=1$, $x_2=n$, $\Delta x=1$ with $\psi(1) = -\gamma$

$$\psi(n) = \sum_{1 \; x=1}^{n-1} \frac{1}{x} - \gamma, \quad \text{Definition of the Digamma Function for } n = 1,2,3,\ldots \tag{8.2-75}$$

From Eq 8.2-62, Eq 8.2-63, and Eq 8.2-74

**The Digamma Function Summation Between Finite Limits is:**

$$\sum_{\Delta x\, x=x_1}^{x_2} \frac{1}{x} = \frac{1}{\Delta x}\psi\left(\frac{x}{\Delta x}\right)\Big|_{x_1}^{x_2+\Delta x} = \zeta\left(1,1,\frac{x}{\Delta x}\right)\Big|_{x_1}^{x_2+\Delta x} = \frac{1}{\Delta x}\,\mathrm{lnd}\left(1,1,\frac{x}{\Delta x}\right)\Big|_{x_1}^{x_2+\Delta x} = \frac{1}{\Delta x}\,\mathrm{lnd}(1,\Delta x,x)\Big|_{x_1}^{x_2+\Delta x}$$

**where** (8.2-76)

**$x$ = real or complex variable**

**$x_1, x_2, \Delta x$ = real or complex constants**

**$\Delta x$ = x increment**

**$x_1 \neq 0,-\Delta x,-2\Delta x,-3\Delta x,\ldots$ , These $x_1,x_2$ exclusions apply only if the Digamma Function is used.**

**$x_2 \neq -\Delta x,-2\Delta x,-3\Delta x,\ldots$ , If the Digamma Function is not used, any term of the summation with a division by zero will be excluded.**

**or**

$$\psi\left(\frac{x}{\Delta x}\right)\Big|_{x_1}^{x_2} = \Delta x \sum_{\Delta x\, x=x_1}^{x_2-\Delta x} \frac{1}{x} = \Delta x\,\zeta\left(1,1,\frac{x}{\Delta x}\right)\Big|_{x_1}^{x_2} = \mathrm{lnd}\left(1,1,\frac{x}{\Delta x}\right)\Big|_{x_1}^{x_2} = \mathrm{lnd}(1,\Delta x,x)\Big|_{x_1}^{x_2} \qquad (8.2\text{-}77)$$

**where**

$$\sum_{\Delta x\, x=x_1}^{x_1-\Delta x} \frac{1}{x} = 0$$

**$x$ = real or complex variable**

**$x_1, x_2, \Delta x, n$ = real or complex constants**

. **$\Delta x$ = x increment**

**$x_1, x_2 \neq 0,-\Delta x,-2\Delta x,-3\Delta x,\ldots$ , These $x_1,x_2$ exclusions apply only if the Digamma Function is used. If the Digamma Function is not used, any term of the summation with a division by zero will be excluded.**

Derive the Digamma Definite Discrete Integral Between Finite Limits Equation.

$$\int_{\Delta x\, x_1}^{x_2} \frac{1}{x}\,\Delta x = \Delta x \sum_{\Delta x\, x=x_1}^{x_2-\Delta x} \frac{1}{x} \qquad (8.2\text{-}78)$$

Substituting Eq 8.2-78 into Eq 8.2-77

**The Digamma Function Definite Discrete Integral Between Finite Limits is:**

$$\int_{\Delta x\, x_1}^{x_2} \frac{1}{x}\Delta x = \Delta x \sum_{\Delta x\, x=x_1}^{x_2-\Delta x} \frac{1}{x} = \psi(\frac{x}{\Delta x})\Big|_{x_1}^{x_2} = \Delta x\zeta(1,1,\frac{x}{\Delta x})\Big|_{x_1}^{x_2} = \text{lnd}(1,1,\frac{x}{\Delta x})\Big|_{x_1}^{x_2} = \text{lnd}(1,\Delta x,x)\Big|_{x_1}^{x_2} \qquad (8.2\text{-}79)$$

**where**

$$\sum_{\Delta x\, x=x_1}^{x_1-\Delta x} \frac{1}{x} = 0$$

**$x$ = real or complex variable**

**$x_1,x_2,\Delta x$ = real or complex constants**

**$\Delta x$ = x increment**

**$x_1,x_2 \neq 0,-\Delta x,-2\Delta x,-3\Delta x,\ldots$, These $x_1,x_2$ exclusions apply only if the Digamma Function is used. If the Digamma Function is not used, any term of the integral with a division by zero will be excluded.**

Derive the Digamma Function indefinite Discrete Integral Equation.

From Eq 8.2-79

**The Digamma Function Indefinite Discrete Integral Equation is:**

$$\int_{\Delta x} \frac{1}{x}\Delta x = \psi(\frac{x}{\Delta x}) + k(\Delta x) = \Delta x\zeta(1,1,\frac{x}{\Delta x}) + k(\Delta x) = \text{lnd}(1,1,\frac{x}{\Delta x}) + k(\Delta x) = \text{lnd}(1,\Delta x,x) + k(\Delta x) \qquad (8.2\text{-}80)$$

**where**

**$x$ = real or complex variable**

**$\Delta x$ = real or complex constant**

**$\Delta x$ = x increment**

**$x \neq 0,-\Delta x,-2\Delta x,-3\Delta x,\ldots$ , This x value exclusion applies only if the Digamma Function is used. If the Digamma Function is not used, any term of the integral with a division by zero will be excluded.**

**$k$ = constant of integration, a function of $\Delta x$**

The Polygamma Functions

Derive the relationship of the Polygamma Functions to the lnd(n,Δx,x) Function.

$$\psi^{(m)}(x) = \frac{d^m}{dx^m}\psi(x) = \frac{d^{m+1}}{dx^{m+1}}\ln\Gamma(x) \text{ , Polygamma Functions definition} \qquad (8.2\text{-}81)$$

where

$m = 1,2,3,\ldots$

$x$ = real or complex values

Note – Eq 8.2-81 shows the relationship of the Polygamma Functions to the Digamma Function.

Rewriting Eq 8.2-59

$\psi(x) = \text{lnd}(1,1,x) - \gamma\ ,\quad x \neq 0,-1,-2,-3,\ldots$ (8.2-82)

where
$x$ = real or complex values
$\gamma$ = Euler's Constant, .5772157…

Take the mth derivative of Eq 8.2-82 using the following equation, Eq 8.2-83

$\frac{d^m}{dx^m}\text{lnd}(1,\Delta x,x) = (-1)^{m+1}m!\,\text{lnd}(1+m,\Delta x,x)$ (8.2-83)

where
$m = 1,2,3,\ldots$
$\Delta x$ = x increment

<u>Note</u> – Eq 8.2-83 is taken from Table 8 row 7 in the Appendix

From Eq 8.2-82 and Eq 8.2-83 with $\Delta x = 1$ and $m=1,2,3,\ldots$

$\psi^{(m)}(x) = \frac{d^m}{dx^m}\psi(x) = \frac{d^m}{dx^m}[\ \text{lnd}(1,1,x) - \gamma\ ] = \frac{d^m}{dx^m}\text{lnd}(1,1,x)$ (8.2-84)

$\psi^{(m)}(x) = (-1)^{m+1}m!\,\text{lnd}(1+m,1,x)$, The Polygamma Functions in terms of the lnd(n,Δx,x) Function (8.2-85)

Rewriting Eq 8.2-25

$\zeta(n,x) = \zeta(n,1,x) = \text{lnd}(n,1,x)\ ,\ \ n \neq 1$ (8.2-86)

where
$x$ = real or complex variable
$x,n$ = real or complex constants

Substituting Eq 8.2-86 into Eq 8.2-85

$\psi^{(m)}(x) = (-1)^{m+1}m!\,\text{lnd}(1+m,1,x) = (-1)^{m+1}m!\,\zeta(1+m,x)$, The Polygamma Functions (8.2-87)

Then from Eq 8.2-87

**The relationship of the Polygamma Functions to the lnd(n,Δx,x) function is:**

$\boldsymbol{\psi^{(m)}(x) = (-1)^{m+1}m!\ \zeta(1+m,x) = (-1)^{m+1}m!\ \zeta(1+m,1,x) = (-1)^{m+1}m!\ \text{lnd}(1+m,1,x)}$ (8.2-88)

**where**
**$m = 1,2,3,\ldots$**
**$x$ = real or complex value**

Derive the Digamma Function Summation Between Finite Limits Equation.

Rewriting Eq 8.2-77, the Digamma Function Summation Between Finite Limits Equation.

$$\psi(\frac{x}{\Delta x})\Big|_{x_1}^{x_2} = \Delta x \sum_{\Delta x \; x=x_1}^{x_2-\Delta x} \frac{1}{x} = \Delta x\,\zeta(1,1,\frac{x}{\Delta x})\Big|_{x_1}^{x_2} = \operatorname{lnd}(1,1,\frac{x}{\Delta x})\Big|_{x_1}^{x_2} = \operatorname{lnd}(1,\Delta x,x)\Big|_{x_1}^{x_2} \qquad (8.2\text{-}89)$$

where

$$\sum_{\Delta x \; x=x_1}^{x_1-\Delta x} \frac{1}{x} = 0$$

$x$ = real or complex variable

$x_1,x_2,\Delta x,n$ = real or complex constants

. $\Delta x$ = x increment

$x_1,x_2 \neq 0,-\Delta x,-2\Delta x,-3\Delta x,\ldots$ , These $x_1,x_2$ exclusions apply only if the Digamma Function is used. If the Digamma Function is not used, any term of the summation with a division by zero will be excluded.

From Eq 8.2-89

$$\psi(\frac{x}{\Delta x})\Big|_{x_1}^{x_2} = \operatorname{lnd}(1,\Delta x,x)\Big|_{x_1}^{x_2} \qquad (8.2\text{-}90)$$

The derivative of the lnd(n,Δx,x) function where n=1is as follows:

$$\frac{d^m}{dx^m}\operatorname{lnd}(1,\Delta x,x) = (-1)^{m+1}m!\operatorname{lnd}(1+m,\Delta x,x) \qquad (8.2\text{-}91)$$

Applying the derivative equation, Eq 8.2-91, to Eq 8.2-90

$$\frac{d^m}{dx^m}\psi(\frac{x}{\Delta x})\Big|_{x_1}^{x_2} = \frac{1}{\Delta x^m}\psi^{(m)}(\frac{x}{\Delta x})\Big|_{x_1}^{x_2} = \frac{d^m}{dx^m}\operatorname{lnd}(1,\Delta x,x)\Big|_{x_1}^{x_2} = (-1)^{m+1}m!\operatorname{lnd}(1+m,\Delta x,x)\Big|_{x_1}^{x_2} \qquad (8.2\text{-}92)$$

Note - $\psi^{(m)}(\frac{x}{\Delta x})$ is the function obtained by the mth derivative of $\psi(\frac{x}{\Delta x})$.

$$\Delta x \sum_{\Delta x \; x=x_1}^{x_2-\Delta x} \frac{1}{x^n} = {}_{\Delta x}\!\int_{x=x_1}^{x_2} \frac{1}{x^n}\,\Delta x = -\operatorname{lnd}(n,\Delta x,x)\Big|_{x_1}^{x_2} = -\Delta x\,\zeta(n,\Delta x,x)\Big|_{x_1}^{x_2}, \quad n \neq 1 \qquad (8.2\text{-}93)$$

From Eq 8.2-93

Let n = 1+m where m = 1,2,3,…

$$\Delta x \sum_{\Delta x\; x=x_1}^{x_2-\Delta x} \frac{1}{x^{1+m}} = \int_{\Delta x\; x=x_1}^{x_2} \frac{1}{x^{1+m}}\,\Delta x = -\,\mathrm{lnd}(1+m,\Delta x,x)\Big|_{x_1}^{x_2} = -\,\Delta x\,\zeta(1+m,\Delta x,x)\Big|_{x_1}^{x_2} \qquad (8.2\text{-}94)$$

From Eq 8.2-92

$$\mathrm{lnd}(1+m,\Delta x,x)\Big|_{x_1}^{x_2} = \frac{(-1)^{m+1}}{m!\Delta x^m}\,\psi^{(m)}\!\left(\frac{x}{\Delta x}\right)\Big|_{x_1}^{x_2} \qquad (8.2\text{-}95)$$

Substituting Eq 8.2-95 into Eq 8.2-94 and dividing by Δx

$$\sum_{\Delta x\; x=x_1}^{x_2-\Delta x} \frac{1}{x^{1+m}} = \frac{1}{\Delta x}\int_{\Delta x\; x=x_1}^{x_2} \frac{1}{x^{1+m}}\,\Delta x = -\frac{1}{\Delta x}\,\mathrm{lnd}(1+m,\Delta x,x)\Big|_{x_1}^{x_2} = \frac{(-1)^m}{m!\Delta x^{m+1}}\,\psi^{(m)}\!\left(\frac{x}{\Delta x}\right)\Big|_{x_1}^{x_2} - \zeta(1+m,\Delta x,x)\Big|_{x_1}^{x_2}$$

(8.2-96)

where

m = 1,2,3,…

Then

From Eq 8.2-95

**The Polygamma Functions Summation Between Finite Limits is:**

$$\sum_{\Delta x\; x=x_1}^{x_2} \frac{1}{x^{1+m}} = \frac{(-1)^m}{m!\Delta x^{m+1}}\,\psi^{(m)}\!\left(\frac{x}{\Delta x}\right)\Big|_{x_1}^{x_2+\Delta x} = -\zeta(1+m,\Delta x,x)\Big|_{x_1}^{x_2+\Delta x} = -\frac{1}{\Delta x}\,\mathrm{lnd}(1+m,\Delta x,x)\Big|_{x_1}^{x_2+\Delta x} \qquad (8.2\text{-}97)$$

**where**

**x = real or complex variable**

**$x_1$,$x_2$,Δx = real or complex constants**

**Δx = x increment**

**m = 1,2,3,…**

**Any summation term where x = 0 is excluded**

**or**

$$\psi^{(m)}\!\left(\frac{x}{\Delta x}\right)\Big|_{x_1}^{x_2} = (-1)^{m+1} m!\Delta x^{m+1}\,\zeta(1+m,\Delta x,x)\Big|_{x_1}^{x_2} = (-1)^m m!\Delta x^{m+1} \sum_{\Delta x\; x=x_1}^{x_2-\Delta x} \frac{1}{x^{1+m}} = (-1)^{m+1} m!\Delta x^m\,\mathrm{lnd}(1+m,\Delta x,x)\Big|_{x_1}^{x_2}$$

**where** (8.2-98)

$$\Delta x\sum_{x=x_1}^{x_1-\Delta x}\frac{1}{x^{1+m}} = 0$$

**x = real or complex variable**

**$x_1$,$x_2$,Δx = real or complex constants**

**Δx = x increment**

**m = 1,2,3,…**

**Any summation term where x = 0 is excluded**

Derive the Polygamma Functions Definite Discrete Integral Between Finite Limits Equation.

From Eq 8.2-96

**The Polygamma Functions Definite Discrete Integral Between Finite Limits is:**

$$\int_{\Delta x\, x_1}^{x_2}\frac{1}{x^{1+m}}\,\Delta x = \Delta x\sum_{\Delta x\, x=x_1}^{x_2-\Delta x}\frac{1}{x^{1+m}} = \frac{(-1)^m}{m!\Delta x^m}\psi^{(m)}(\frac{x}{\Delta x})\Big|_{x_1}^{x_2} - \Delta x\zeta(1+m,\Delta x,x)\Big|_{x_1}^{x_2} = -\,\mathrm{lnd}(1+m,\Delta x,x)\Big|_{x_1}^{x_2}$$

**where** (8.2-99)

$$\Delta x\sum_{x=x_1}^{x_1-\Delta x}\frac{1}{x} = 0$$

**x = real or complex variable**

**$x_1$,$x_2$,Δx = real or complex constants**

**Δx = x increment**

**m = 1,2,3,…**

**Any summation term where x = 0 is excluded**

Derive the Polygamma Functions Indefinite Discrete Integral Between Finite Limits Equation.

From Eq 8.2-99

**The Polygamma Functions Indefinte Discrete Integral is:**

$$\int_{\Delta x}\frac{1}{x^{1+m}}\,\Delta x = \frac{(-1)^m}{m!\Delta x^m}\psi^{(m)}(\frac{x}{\Delta x}) + k(\Delta x) = -\,\Delta x\zeta(1+m,\Delta x,x) + k(\Delta x) = -\,\mathrm{lnd}(1+m,\Delta x,x) + k(\Delta x)$$

**where** (8.2-100)

**x = real or complex variable**

**Δx = real or complex constant**

**Δx = x increment**

**m = 1,2,3,…**

**k = constant of integration, a function of Δx**

Derive the Polygamma Functions Summation to Infinity Equation.

From Eq 8.2-97

$m = 1,2,3,\ldots$

Let $x_2 \to \pm\infty$

$+\infty$ for $Re(\Delta x)>0$ or [ $Re(\Delta x)=0$ and $Im(\Delta x)>0$ ]

$-\infty$ for $Re(\Delta x)<0$ or [ $Re(\Delta x)=0$ and $Im(\Delta x)<0$ ]

$$\sum_{\Delta x \; x=x_1}^{\pm\infty} \frac{1}{x^{1+m}} = \frac{(-1)^m}{m!\Delta x^{m+1}} \psi^{(m)}(\frac{x}{\Delta x})\Big|_{x_1}^{\pm\infty} = -\zeta(1+m,\Delta x,x)\Big|_{x_1}^{\pm\infty} = -\frac{1}{\Delta x} \text{lnd}(1+m,\Delta x,x)\Big|_{x_1}^{\pm\infty} \qquad (8.2\text{-}101)$$

$$\text{lnd}(n,\Delta x,x)\Big|_{x\to\pm\infty} \to 0 \text{ when } Re(n)>1 \qquad (8.2\text{-}102)$$

Rewriting Eq 8.2-88

The relationship of the Polygamma Functions to the lnd(n,Δx,x) is:

$$\psi^{(m)}(x) = (-1)^{m+1}m!\,\zeta(1+m,x) = (-1)^{m+1}m!\,\zeta(1+m,1,x) = (-1)^{m+1}m!\,\text{lnd}(1+m,1,x) \qquad (8.2\text{-}103)$$

where

$m = 1,2,3,\ldots$

$x$ = real or complex value

From Eq 8.2-101 and Eq 8.2-103

Let $n = 1+m$

$m = 1,2,3,\ldots$

$Re(1+m) > 1$

$$\lim\Big|_{x\to\pm\infty} \psi^{(m)}(x) = (-1)^{m+1}m![\lim_{x\to\pm\infty} \text{lnd}(1+m,1,x)] \to 0 \qquad (8.2\text{-}104)$$

$$\text{lnd}(1+m,\Delta x,x)\Big|_{x\to\pm\infty} \to 0 \qquad (8.2\text{-}105)$$

From Eq 101, Eq 104, and Eq 105

$$\sum_{\Delta x \; x=x_1}^{\pm\infty} \frac{1}{x^{1+m}} = \frac{(-1)^{m+1}}{m!\Delta x^{m+1}} \psi^{(m)}(\frac{x_1}{\Delta x}) = \zeta(1+m,\Delta x,x_1) = \frac{1}{\Delta x} \text{lnd}(1+m,\Delta x,x_1)) \qquad (8.2\text{-}106)$$

Then

From Eq 8.2-106

**The Polygamma Function Summation to Infinity is:**

$$\sum_{\Delta x \; x=x_i}^{\pm\infty} \frac{1}{x^{1+m}} = \frac{(-1)^{m+1}}{m!\Delta x^{m+1}}\psi^{(m)}\left(\frac{x_i}{\Delta x}\right) = \zeta(1+m,\Delta x,x_i) = \frac{1}{\Delta x}\,\mathrm{lnd}(1+m,\Delta x,x_i) \qquad (8.2\text{-}107)$$

**where**

**x = real or complex variable**

**$x_i$,m,Δx = real or complex constants**

**+∞ for Re(Δx)>0 or [ Re(Δx)=0 and Im(Δx)>0 ]**

**–∞ for Re(Δx)<0 or [ Re(Δx)=0 and Im(Δx)<0 ]**

**m = 1,2,3,…**

**Δx = x increment**

**Any summation term where x = 0 is excluded**

Derive the Polygamma Functions Definite Discrete Integral to Infinity Equation.

$$\int_{\Delta x \; x_i}^{\pm\infty} \frac{1}{x^{1+m}}\,\Delta x = \Delta x \sum_{\Delta x \; x=x_i}^{\pm\infty} \frac{1}{x^{1+m}} \qquad (8.2\text{-}108)$$

Multiplying Eq 8.2-107 by Δx

$$\Delta x \sum_{\Delta x \; x=x_i}^{\pm\infty} \frac{1}{x^{1+m}} = \frac{(-1)^{m+1}}{m!\Delta x^{m}}\psi^{(m)}\left(\frac{x_i}{\Delta x}\right) = \Delta x\zeta(1+m,\Delta x,x_i) = \mathrm{lnd}(1+m,\Delta x,x_i) \qquad (8.2\text{-}109)$$

Substituting Eq 8.2-108 into Eq 8.2-109

$$\int_{\Delta x \; x_i}^{\pm\infty} \frac{1}{x^{1+m}}\,\Delta x = \Delta x \sum_{\Delta x \; x=x_i}^{\pm\infty} \frac{1}{x^{1+m}} = \frac{(-1)^{m+1}}{m!\Delta x^{m}}\psi^{(m)}\left(\frac{x_i}{\Delta x}\right) = \Delta x\zeta(1+m,\Delta x,x_i) = \mathrm{lnd}(1+m,\Delta x,x_i) \qquad (8.2\text{-}110)$$

Then

From Eq 8.2-110

**The Polygamma Functions Definite Discrete Integral to Infinity is:**

$$\underset{\Delta x}{\int_{x_i}^{\pm\infty}} \frac{1}{x^{1+m}}\Delta x = \frac{(-1)^{m+1}}{m!\Delta x^m}\psi^{(m)}\left(\frac{x_i}{\Delta x}\right) = \Delta x \underset{\Delta x}{\sum_{x=x_i}^{\pm\infty}} \frac{1}{x^{1+m}} = \Delta x\zeta(1+m,\Delta x,x_i) = \mathrm{lnd}(1+m,\Delta x,x_i) \tag{8.2-111}$$

**where**

**x = real or complex variable**

**$x_i$,m,Δx = real or complex constants**

**+∞ for Re(Δx)>0 or [ Re(Δx)=0 and Im(Δx)>0 ]**

**–∞ for Re(Δx)<0 or [ Re(Δx)=0 and Im(Δx)<0 ]**

**m = 1,2,3,…**

**Δx = x increment**

**Any summation term where x = 0 is excluded**

Derive the Polygamma Functions term by term relationship.

From Eq 8.2-45, the Digamma Function term to term relationship

$$\psi(x+1) - \psi(x) = \frac{1}{x} \tag{8.2-112}$$

where

x = real or complex variable

Taking m successive derivatives of Eq 8.2-112

$$\frac{d^m}{dx^m}[\psi(x+1) - \psi(x)] = \frac{d^m}{dx^m}\left[\frac{1}{x}\right] \tag{8.2-113}$$

$$\psi^{(m)}(x+1) - \psi^{(m)}(x) = (-1)^m m!\frac{1}{x^{m+1}} \tag{8.2-114}$$

Then

From Eq 8.2-114

**The Polygamma Functions term to term relationship is:**

$$\boldsymbol{\psi^{(m)}(x+1) - \psi^{(m)}(x) = (-1)^m m!\frac{1}{x^{m+1}}} \tag{8.2-115}$$

**where**

**m = 1,2,3,…**

**x = real or complex value**

Derive the Polygamma Functions discrete derivative

$$D_1\psi^{(m)}(x) = \psi^{(m)}(x+1) - \psi^{(m)}(x) \tag{8.2-116}$$

Substituting Eq 8.2-116 into Eq 8.2-115

$$D_1\psi^{(m)}(x) = (-1)^m m!\frac{1}{x^{m+1}} \tag{8.2-117}$$

Then

From Eq 8.2-117

**The Polygamma Functions discrete derivative is:** (8.2-118)

$$\mathbf{D_1\psi^{(m)}(x) = (-1)^m m!\frac{1}{x^{m+1}}}$$

**where**

**m = 1,2,3,…**

**x = real or complex variable**

Derive the Polygamma Recursion Equation.

Rearranging the terms of Eq 8.2-115

**The Polygamma Functions Recursion Equation is:**

$$\mathbf{\psi^{(m)}(x+1) = \psi^{(m)}(x) + (-1)^m m!\frac{1}{x^{m+1}}} \tag{8.2-119}$$

**where**

**m = 1,2,3,…**

**x = real or complex value**

A listing of the above derived equations is provided in Table 8.2-1 on the following page.

**Table 8.2-1 A listing of Zeta Function, Digamma Function, Polygamma Function, and lnd(n,Δx,x) Function Equations**

| | **Zeta Function, Digamma Function, Polygamma Function, and lnd(n,Δx,x) Function Equations** |
|---|---|
| | **The General Zeta Function** |
| 1 | The relationship of the General Zeta Function to the lnd(n,Δx,x) function<br><br>$\zeta(n,\Delta x,x) = \frac{1}{\Delta x}\,\mathrm{lnd}(n,\Delta x,x)$<br><br>where<br>x = real or complex variable<br>n,Δx = real or complex constants<br>Δx = x increment |
| 1a | The General Zeta Function Summation Between Finite Limits<br><br>$\sum_{\Delta x\; x=x_1}^{x_2} \frac{1}{x^n} = \pm\, \zeta(n,\Delta x,x)\Big\|_{x_1}^{x_2+\Delta x} = \pm \frac{1}{\Delta x}\,\mathrm{lnd}(n,\Delta x,x)\Big\|_{x_1}^{x_2+\Delta x}$ , – for n≠1, + for n=1<br><br>where<br>x = real or complex variable<br>$x_1, x_2$,Δx,n = real or complex constants<br>Δx = x increment<br>Any summation term where x = 0 is excluded<br><br>or<br><br>$\zeta(n,\Delta x,x)\Big\|_{x_1}^{x_2} = \pm \sum_{\Delta x\; x=x_1}^{x_2-\Delta x} \frac{1}{x^n} = \frac{1}{\Delta x}\,\mathrm{lnd}(n,\Delta x,x)\Big\|_{x_1}^{x_2}$ , – for n≠1, + for n=1<br><br>where<br>$\sum_{\Delta x\; x=x_1}^{x_1-\Delta x} \frac{1}{x^n} = 0$<br>x = real or complex variable<br>$x_1, x_2$,Δx,n = real or complex constants<br>Δx = x increment<br>Any summation term where x = 0 is excluded |

| | **Zeta Function, Digamma Function, Polygamma Function, and lnd(n,Δx,x) Function Equations** |
|---|---|
| | The following equation, which relates the lnd(n,Δx,x) Function to the Hurwitz Zeta Function, may be used.<br><br>$\text{lnd}(n,\Delta x,x) = \frac{1}{\Delta x^{n-1}}\zeta(n,\frac{x}{\Delta x})$<br><br>where<br>x,Δx = real values with the condition that x ≠ positive real value when Δx = negative real value<br>n = complex values<br>or<br>x = Δx = complex values<br>n = complex value<br>or<br>x,Δx = complex values<br>n = integer<br><br>ζ(n,x) = Hurwitz Zeta Function , n ≠ 1<br><br>Comments – A complex value, x+jy, can be a real value, x, an imaginary value, jy, or an integer, N+j0.<br><br>For value combinations not specified above, the equality of the stated equation may not be valid. |
| 1b | The General Zeta Function Definite Discrete Integral Between Finite Limits<br><br>$_{\Delta x}\int_{x_1}^{x_2}\frac{1}{x^n}\Delta x = \Delta x \, {}_{\Delta x}\sum_{x=x_1}^{x_2-\Delta x}\frac{1}{x^n} = \pm\,\Delta x\zeta(n,\Delta x,x)\Big\|_{x_1}^{x_2} = \pm\text{lnd}(n,\Delta x,x)\Big\|_{x_1}^{x_2}$ , – for n≠1, + for n=1<br><br>where<br>$_{\Delta x}\sum_{x=x_1}^{x_1-\Delta x}\frac{1}{x} = 0$<br>x = real or complex variable<br>$x_1,x_2,\Delta x,n$ = real or complex constants<br>Δx = x increment<br>Any summation term where x = 0 is excluded |

| | **Zeta Function, Digamma Function, Polygamma Function, and lnd(n,Δx,x) Function Equations** |
|---|---|
| 1c | <u>The General Zeta Function Indefinite Discrete Integral</u><br>$\int_{\Delta x} \frac{1}{x^n}\,\Delta x = \pm\,\Delta x\,\zeta(n,\Delta x,x) + k(n,\Delta x) = \pm\,\text{lnd}(n,\Delta x,x) + k(n,\Delta x)$ , $-$ for $n\neq 1$, $+$ for $n=1$<br>where<br>x = real or complex variable<br>Δx,n = real or complex constants<br>Δx = x increment<br>k = constant of integration, a function of n,Δx |
| 1d | <u>The General Zeta Function Summation to Infinity</u><br>$\sum_{\Delta x,\ x=x_i}^{\pm\infty} \frac{1}{x^n} = \zeta(n,\Delta x,x_i) = \frac{1}{\Delta x}\,\text{lnd}(n,\Delta x,x_i)$ , $\text{Re}(n)>1$<br>where<br>x = real or complex variable<br>$x_i$,n,Δx = real or complex constants<br>Δx = x increment<br>$+\infty$ for $\text{Re}(\Delta x)>0$ or [ $\text{Re}(\Delta x)=0$ and $\text{Im}(\Delta x)>0$ ]<br>$-\infty$ for $\text{Re}(\Delta x)<0$ or [ $\text{Re}(\Delta x)=0$ and $\text{Im}(\Delta x)<0$ ]<br>Any summation term where x = 0 is excluded |
| 1e | <u>The General Zeta Function Definite Discrete Integral to Infinity</u><br>$\int_{\Delta x,\ x=x_i}^{\pm\infty} \frac{1}{x^n}\,\Delta x = \Delta x \sum_{\Delta x,\ x=x_i}^{\pm\infty} \frac{1}{x^n} = \Delta x\,\zeta(n,\Delta x,x_i) = \text{lnd}(n,\Delta x,x_i)$ , $\text{Re}(n)>1$<br>where<br>x = real or complex variable<br>$x_i$,n,Δx = real or complex constants<br>Δx = x increment<br>$+\infty$ for $\text{Re}(\Delta x)>0$ or [ $\text{Re}(\Delta x)=0$ and $\text{Im}(\Delta x)>0$ ]<br>$-\infty$ for $\text{Re}(\Delta x)<0$ or [ $\text{Re}(\Delta x)=0$ and $\text{Im}(\Delta x)<0$ ]<br>Any summation term where x = 0 is excluded |

| | Zeta Function, Digamma Function, Polygamma Function, and lnd(n,Δx,x) Function Equations |
|---|---|
| 1f | The General Zeta Function Discrete Derivative<br><br>$D_{\Delta x}\,\zeta(n,\Delta x,x) = \pm\frac{1}{\Delta x}\frac{1}{x^n}$,   – for n≠1,  + for n=1<br><br>where<br>x = real or complex variable<br>Δx,n = real or complex constants<br>Δx = x increment |
| 1g | The General Zeta Function term to term relationship<br><br>$\zeta(n,\Delta x,x+\Delta x) - \zeta(n,\Delta x,x) = \pm\frac{1}{x^n}$,   – for n≠1,  + for n=1<br><br>where<br>x = real or complex variable<br>Δx,n = real or complex constants<br>Δx = x increment |
| 1h | The General Zeta Function Recursion Equation<br><br>$\zeta(n,\Delta x,x+\Delta x) = \zeta(n,\Delta x,x) \pm \frac{1}{x^n}$,  – for n≠1,  + for n=1,<br><br>where<br>x = real or complex variable<br>Δx,n = real or complex constants<br>Δx = x increment |
| | **The Hurwitz Zeta Function** |
| 2 | The relationship of the Hurwitz Zeta Function to the lnd(n,Δx,x) function<br><br>$\zeta(n,x) = \zeta(n,1,x) = \mathrm{lnd}(n,1,x)$ ,  n ≠ 1<br><br>where<br>x = real or complex variable<br>n = real or complex constant |
| 2a | The Hurwitz Zeta Function Summation Between Finite Limits<br><br>$\sum_{1\ x=x_1}^{x_2}\frac{1}{x^n} = -\zeta(n,x)\Big\vert_{x_1}^{x_2+1} = -\zeta(n,1,x)\Big\vert_{x_1}^{x_2+1} = -\mathrm{lnd}(n,1,x)\Big\vert_{x_1}^{x_2+1}$ ,  n ≠ 1<br><br>where<br>x = real or complex variable<br>$x_1, x_2$,n = real or complex constants<br>Any summation term where x = 0 is excluded<br><br>or |

| | **Zeta Function, Digamma Function, Polygamma Function, and lnd(n,Δx,x) Function Equations** |
|---|---|
| | $$\zeta(n,x)\Big\|_{x_1}^{x_2} = -\,{}_{1}\!\sum_{x=x_1}^{x_2-1}\frac{1}{x^n} = \zeta(n,1,x)\Big\|_{x_1}^{x_2} = \text{lnd}(n,1,x)\Big\|_{x_1}^{x_2}\,, \quad n\neq 1$$ <br> where <br> $${}_{\Delta x}\!\sum_{x=x_1}^{x_1-1}\frac{1}{x^n} = 0$$ <br> x = real or complex variable <br> $x_1, x_2$, n = real or complex constants <br> Any summation term where x = 0 is excluded |
| 2b | <u>The Hurwitz Zeta Function Definite Discrete Integral Between Finite Limits</u> <br> $${}_{1}\!\int_{x_1}^{x_2}\frac{1}{x^n}\,\Delta x = {}_{1}\!\sum_{x=x_1}^{x_2-1}\frac{1}{x^n} = -\,\zeta(n,x)\Big\|_{x_1}^{x_2} = -\,\zeta(n,1,x)\Big\|_{x_1}^{x_2} = -\,\text{lnd}(n,1,x)\Big\|_{x_1}^{x_2}\,, \quad n\neq 1$$ <br> where <br> $${}_{1}\!\sum_{x=x_1}^{x_1-1}\frac{1}{x} = 0$$ <br> x = real or complex variable <br> $x_1, x_2$, n = real or complex constants <br> Δx = 1, x increment <br> Any summation term where x = 0 is excluded |
| 2c | <u>The Hurwitz Zeta Function Indefinite Discrete Integral</u> <br> $${}_{1}\!\int\frac{1}{x^n}\,\Delta x = -\,\zeta(n,x) + k(n) = -\,\zeta(n,1,x) + k(n) = -\,\text{lnd}(n,1,x) + k(n)\,, \quad n\neq 1$$ <br> where <br> x = real or complex variable <br> n = real or complex constant <br> Δx = 1, x increment <br> k = constant of integration, a function of n |

| | **Zeta Function, Digamma Function, Polygamma Function, and lnd(n,Δx,x) Function Equations** |
|---|---|
| 2d | The Hurwitz Zeta Function Summation to Infinity<br><br>$$\sum_{x=x_i}^{\infty}{}_1 \frac{1}{x^n} = \zeta(n,x_i) = \zeta(n,1,x_i) = \text{lnd}(n,1,x_i)\,, \quad \text{Re}(n)>1$$<br>where<br>x = real or complex variable<br>$x_i$,n = real or complex constants<br>Δx = 1, x increment<br>Any summation term where x = 0 is excluded |
| 2e | The Hurwitz Zeta Function Definite Discrete Integral to Infinity<br><br>$$\int_{x=x_i}^{\infty}{}_1 \frac{1}{x^n}\Delta x = \sum_{x=x_i}^{\infty}{}_1 \frac{1}{x^n} = \zeta(n,x_i) = \zeta(n,1,x_i) = \text{lnd}(n,1,x_i)\,, \quad \text{Re}(n)>1$$<br>where<br>x = real or complex variable<br>$x_i$,n= real or complex constants<br>Δx = 1, x increment<br>Any summation term where x = 0 is excluded |
| 2f | The Hurwitz Zeta Function Discrete Derivative<br><br>$$D_1\zeta(n,x) = -\frac{1}{x^n}\,, \quad n\neq 1$$<br>where<br>x = real or complex variable<br>n = real or complex constant |
| 2g | The Hurwitz Zeta Function term to term relationship<br><br>$$\zeta(n,x+1) - \zeta(n,x) = -\frac{1}{x^n}\,, \quad n\neq 1$$<br>where<br>x = real or complex variable<br>n = real or complex constant |

| | **Zeta Function, Digamma Function, Polygamma Function, and lnd(n,Δx,x) Function Equations** |
|---|---|
| 2h | <u>The Hurwitz Zeta Function Recursion Equation</u><br><br>$\zeta(n,x+1) = \zeta(n,x) - \frac{1}{x^n}$ , $n \neq 1$<br>where<br>x = real or complex variable<br>n = real or complex constant |
| | **The Riemann Zeta Function** |
| 3 | <u>The relationship of the Riemann Zeta Function to the lnd(n,Δx,x) function</u><br><br>$\zeta(n) = \zeta(n,1,1) = \text{lnd}(n,1,1)$ , $n \neq 1$<br>where<br>n = real or complex constant |
| 3a | <u>The Riemann Zeta Function Summation to Infinity</u><br><br>$\zeta(n) = \sum_{x=1}^{\infty} \frac{1}{x^n} = \zeta(n,1,1) = \text{lnd}(n,1,1)$ , $\text{Re}(n)>1$<br>where<br>x = real or complex variable<br>n = real or complex constant |
| 3b | <u>The Riemann Zeta Function Definite Discrete Integral to Infinity</u><br><br>$\int_{x=1}^{\infty} \frac{1}{x^n} \Delta x = \sum_{x=1}^{\infty} \frac{1}{x^n} = \zeta(n) = \zeta(n,1,1) = \text{lnd}(n,1,1)$ , $\text{Re}(n)>1$<br>where<br>x = real or complex variable<br>n = real or complex constant<br>Δx = 1, x increment |
| | **The Polygamma Functions** |
| 4 | <u>The relationship of the Polygamma Functions to the lnd(n,Δx,x) function</u><br><br>$\psi^{(m)}(x) = (-1)^{m+1} m! \zeta(1+m,x) = (-1)^{m+1} m! \zeta(1+m,1,x) = (-1)^{m+1} m! \text{lnd}(1+m,1,x)$<br>where<br>m = 1,2,3,…<br>x = real or complex value |

| | Zeta Function, Digamma Function, Polygamma Function, and lnd(n,Δx,x) Function Equations |
|---|---|
| 4a | The Polygamma Functions Summation Between Finite Limits<br>$$\sum_{\Delta x\, x=x_1}^{x_2} \frac{1}{x^{1+m}} = \frac{(-1)^m}{m!\Delta x^{m+1}}\psi^{(m)}(\frac{x}{\Delta x})\Big\|_{x_1}^{x_2+\Delta x} = -\zeta(1+m,\Delta x,x)\Big\|_{x_1}^{x_2+\Delta x} = -\frac{1}{\Delta x}\,\mathrm{lnd}(1+m,\Delta x,x)\Big\|_{x_1}^{x_2+\Delta x}$$<br>where<br>x = real or complex variable<br>$x_1, x_2, \Delta x$ = real or complex constants<br>Δx = x increment<br>m = 1,2,3,…<br>Any summation term where x = 0 is excluded<br><br>or<br>$$\psi^{(m)}(\frac{x}{\Delta x})\Big\|_{x_1}^{x_2} = (-1)^{m+1} m!\Delta x^{m+1}\zeta(1+m,\Delta x,x)\Big\|_{x_1}^{x_2} = (-1)^m m!\Delta x^{m+1}\sum_{\Delta x\, x=x_1}^{x_2-\Delta x}\frac{1}{x^{1+m}} = (-1)^{m+1} m!\Delta x^m\,\mathrm{lnd}(1+m,\Delta x,x)\Big\|_{x_1}^{x_2}$$<br>where<br>$$\sum_{\Delta x\, x=x_1}^{x_1-\Delta x}\frac{1}{x^{1+m}} = 0$$<br>x = real or complex variable<br>$x_1, x_2, \Delta x$ = real or complex constants<br>Δx = x increment<br>m = 1,2,3,…<br>Any summation term where x = 0 is excluded |
| 4b | The Polygamma Functions Definite Discrete Integral Between Finite Limits<br>$$\int_{\Delta x\, x_1}^{x_2}\frac{1}{x^{1+m}}\,\Delta x = \Delta x\sum_{\Delta x\, x=x_1}^{x_2-\Delta x}\frac{1}{x^{1+m}} = \frac{(-1)^m}{m!\Delta x^m}\psi^{(m)}(\frac{x}{\Delta x})\Big\|_{x_1}^{x_2} - \Delta x\,\zeta(1+m,\Delta x,x)\Big\|_{x_1}^{x_2} = -\,\mathrm{lnd}(1+m,\Delta x,x)\Big\|_{x_1}^{x_2}$$<br>where<br>$$\sum_{\Delta x\, x=x_1}^{x_1-\Delta x}\frac{1}{x} = 0$$<br>x = real or complex variable<br>$x_1, x_2, \Delta x$ = real or complex constants<br>Δx = x increment<br>m = 1,2,3,…<br>Any summation term where x = 0 is excluded |

| | **Zeta Function, Digamma Function, Polygamma Function, and lnd(n,Δx,x) Function Equations** |
|---|---|
| 4c | The Polygamma Functions Indefinte Discrete Integral<br><br>$\int_{\Delta x} \frac{1}{x^{1+m}}\,\Delta x = \frac{(-1)^m}{m!\Delta x^m}\,\psi^{(m)}(\frac{x}{\Delta x}) + k(\Delta x) = -\,\Delta x\,\zeta(1+m,\Delta x,x) + k(\Delta x) = -\,\mathrm{lnd}(1+m,\Delta x,x) + k(\Delta x)$<br>where<br>x = real or complex variable<br>Δx = real or complex constant<br>Δx = x increment<br>m = 1,2,3,…<br>k = constant of integration, a function of Δx |
| 4d | The Polygamma Function Summation to Infinity<br><br>$\sum_{\Delta x\; x=x_i}^{\pm\infty} \frac{1}{x^{1+m}} = \frac{(-1)^{m+1}}{m!\Delta x^{m+1}}\,\psi^{(m)}(\frac{x_i}{\Delta x}) = \zeta(1+m,\Delta x,x_i) = \frac{1}{\Delta x}\,\mathrm{lnd}(1+m,\Delta x,x_i)$<br>where<br>x = real or complex variable<br>$x_i$,m,Δx = real or complex constants<br>+∞ for Re(Δx)>0 or [ Re(Δx)=0 and Im(Δx)>0 ]<br>–∞ for Re(Δx)<0 or [ Re(Δx)=0 and Im(Δx)<0 ]<br>m = 1,2,3,…<br>Δx = x increment<br>Any summation term where x = 0 is excluded |
| 4e | The Polygamma Functions Definite Discrete Integral to Infinity<br><br>$\int_{\Delta x\; x_i}^{\pm\infty} \frac{1}{x^{1+m}}\,\Delta x = \frac{(-1)^{m+1}}{m!\Delta x^m}\,\psi^{(m)}(\frac{x_i}{\Delta x}) = \Delta x \sum_{\Delta x\; x=x_i}^{\pm\infty} \frac{1}{x^{1+m}} = \Delta x\,\zeta(1+m,\Delta x,x_i) = \mathrm{lnd}(1+m,\Delta x,x_i)$<br>where<br>x = real or complex variable<br>$x_i$,m,Δx = real or complex constants<br>+∞ for Re(Δx)>0 or [ Re(Δx)=0 and Im(Δx)>0 ]<br>–∞ for Re(Δx)<0 or [ Re(Δx)=0 and Im(Δx)<0 ]<br>m = 1,2,3,…<br>Δx = x increment<br>Any summation term where x = 0 is excluded |

| | **Zeta Function, Digamma Function, Polygamma Function, and lnd(n,Δx,x) Function Equations** |
|---|---|
| 4f | The Polygamma Functions Discrete Derivative<br><br>$D_1\psi^{(m)}(x) = (-1)^m m!\frac{1}{x^{m+1}}$<br>where<br>m = 1,2,3,…<br>x = real or complex variable |
| 4g | The Polygamma Functions term to term relationship<br>$\psi^{(m)}(x+1) - \psi^{(m)}(x) = (-1)^m m!\frac{1}{x^{m+1}}$<br>where<br>m = 1,2,3,…<br>x = real or complex value |
| 4h | The Polygamma Functions Recursion Equation<br><br>$\psi^{(m)}(x+1) = \psi^{(m)}(x) + (-1)^m m!\frac{1}{x^{m+1}}$<br>where<br>m = 1,2,3,…<br>x = real or complex value |
| | **The Digamma Function** |
| 5 | The relationship of the Digamma Function to the lnd(n,Δx,x) function<br><br>$\psi(x) = \mathrm{lnd}(1,1,x) - \gamma$ ,   $x \neq 0,-1,-2,-3,\ldots$<br>or<br>$\psi(\frac{x}{\Delta x}) = \mathrm{lnd}(1,\Delta x,x) - \gamma = \mathrm{lnd}(1,1,\frac{x}{\Delta x}) - \gamma = \Delta x\,\zeta(1,\Delta x,x) - \gamma$ ,   $x \neq 0,-\Delta x,-2\Delta x,-3\Delta x,\ldots$<br><br>where<br>x = real or complex values<br>Δx = x increment<br>γ = Euler’s Constant, .5772157…<br>ψ(x) is infinite for x = 0,-1,-2,-3,…<br>$\psi(\frac{x}{\Delta x})$ is infinite for $x = 0,-\Delta x,-2\Delta x,-3\Delta x,\ldots$<br>ζ(1,Δx,x) is the N=1 Zeta Function<br><br>Comments – The Digamma Function, ψ(x), has first order poles at x ≠ 0,-1,-2,-3,….<br>At these values of x, where $\psi(x) = \frac{1}{0}$ , the value of ψ(x) is infinite.<br><br>The lnd(1,1,x) function minus Euler’s number is equal to the Digamma Function for all values of x except for x = 0,-1,-2,-3,… For x = 0,-1,-2,-3,…, the lnd(1,1,x) function |

differs from the ψ(x) function. The equation used to calculate the lnd(1,1,x) function when
x = 0,-1,-2,-3,… is lnd(1,1,x) = lnd(1,1,1−x). In addition, any summation calculated using the lnd(1,1,x) function, excludes any summation term with a division by zero. These differences were introduced into the lnd(1,1,x) function to make possible the summation of integers along the real axis of the complex plane.

For x = 0,-1,-2,-3,…, if the program, LNDX, is selected to find ψ(x) using the equation, ψ(x) = lnd(1,1,x) − γ, the resulting calculated value for ψ(x) is ψ(x) = ψ(1−x).

5a The Digamma Function Summation Between Finite Limits

$$\sum_{\Delta x\; x=x_1}^{x_2} \frac{1}{x} = \frac{1}{\Delta x}\,\psi\left(\frac{x}{\Delta x}\right)\Big|_{x_1}^{x_2+\Delta x} = \zeta\left(1,1,\frac{x}{\Delta x}\right)\Big|_{x_1}^{x_2+\Delta x} = \frac{1}{\Delta x}\,\mathrm{lnd}\left(1,1,\frac{x}{\Delta x}\right)\Big|_{x_1}^{x_2+\Delta x} = \frac{1}{\Delta x}\,\mathrm{lnd}(1,\Delta x,x)\Big|_{x_1}^{x_2+\Delta x}$$

where
x = real or complex variable
$x_1,x_2,\Delta x$ = real or complex constants
Δx = x increment
$x_1 \neq 0,-\Delta x,-2\Delta x,-3\Delta x,\ldots$ , These $x_1,x_2$ exclusions apply only if the Digamma Function is used.
$x_2 \neq -\Delta x,-2\Delta x,-3\Delta x,\ldots$ , If the Digamma Function is not used, any term of the summation with a division by zero will be excluded.

or

$$\psi\left(\frac{x}{\Delta x}\right)\Big|_{x_1}^{x_2} = \Delta x \sum_{\Delta x\; x=x_1}^{x_2-\Delta x} \frac{1}{x} = \Delta x\,\zeta\left(1,1,\frac{x}{\Delta x}\right)\Big|_{x_1}^{x_2} = \mathrm{lnd}\left(1,1,\frac{x}{\Delta x}\right)\Big|_{x_1}^{x_2} = \mathrm{lnd}(1,\Delta x,x)\Big|_{x_1}^{x_2}$$

where

$$\sum_{\Delta x\; x=x_1}^{x_1-\Delta x} \frac{1}{x} = 0$$

x = real or complex variable
$x_1,x_2,\Delta x,n$ = real or complex constants
. Δx = x increment
$x_1,x_2 \neq 0,-\Delta x,-2\Delta x,-3\Delta x,\ldots$ , These $x_1,x_2$ exclusions apply only if the Digamma Function is used. If the Digamma Function is not used, any term of the summation with a division by zero will be excluded.

| | **Zeta Function, Digamma Function, Polygamma Function, and lnd(n,Δx,x) Function Equations** |
|---|---|
| 5b | The Digamma Function Definite Discrete Integral Between Finite Limits<br><br>$${}_{\Delta x}\int_{x_1}^{x_2} \frac{1}{x}\,\Delta x = \Delta x \,{}_{\Delta x}\sum_{x=x_1}^{x_2-\Delta x} \frac{1}{x} = \psi\left(\frac{x}{\Delta x}\right)\Big|_{x_1}^{x_2} = \Delta x\,\zeta\left(1,1,\frac{x}{\Delta x}\right)\Big|_{x_1}^{x_2} = \text{lnd}\left(1,1,\frac{x}{\Delta x}\right)\Big|_{x_1}^{x_2} = \text{lnd}(1,\Delta x,x)\Big|_{x_1}^{x_2}$$<br>where<br>$${}_{\Delta x}\sum_{x=x_1}^{x_1-\Delta x} \frac{1}{x} = 0$$<br>x = real or complex variable<br>$x_1,x_2,\Delta x$ = real or complex constants<br>Δx = x increment<br>$x_1,x_2 \neq 0,-\Delta x,-2\Delta x,-3\Delta x,\ldots$, These $x_1,x_2$ exclusions apply only if the Digamma Function is used. If the Digamma Function is not used, any term of the integral with a division by zero will be excluded. |
| 5c | The Digamma Function Indefinite Discrete Integral<br><br>$${}_{\Delta x}\int \frac{1}{x}\,\Delta x = \psi\left(\frac{x}{\Delta x}\right) + k(\Delta x) = \Delta x\,\zeta\left(1,1,\frac{x}{\Delta x}\right) + k(\Delta x) = \text{lnd}\left(1,1,\frac{x}{\Delta x}\right) + k(\Delta x) = \text{lnd}(1,\Delta x,x) + k(\Delta x)$$<br>where<br>x = real or complex variable<br>Δx = real or complex constant<br>Δx = x increment<br>$x \neq 0,-\Delta x,-2\Delta x,-3\Delta x,\ldots$ , This x value exclusion applies only if the Digamma Function is used. If the Digamma Function is not used, any term of the integral with a division by zero will be excluded.<br>k = constant of integration, a function of Δx |
| 5d | The Digamma Function Discrete Derivative<br><br>$$D_1\psi(x) = \frac{1}{x}$$<br>where<br>x = real or complex variable<br>n = real or complex constant |

| | **Zeta Function, Digamma Function, Polygamma Function, and lnd(n,Δx,x) Function Equations** |
|---|---|
| 5e | The Digamma Function term to term relationship<br><br>$\psi(x+1) - \psi(x) = \frac{1}{x}$<br>where<br>x = real or complex variable |
| 5f | The Digamma Function Recursion Equation<br>$\psi(x+1) = \psi(x) + \frac{1}{x}$<br>where<br>x = real or complex variable |

## Section 8.3: Demonstration of the use of the lnd(n,Δx,x) function to evaluate the Hurwitz Zeta Function, the Riemann Zeta Function, and the Digamma Function

To demonstrate the usefulness of the lnd(n,Δx,x) function in the evaluation of the Hurwitz Zeta Function, the Riemann Zeta Function, and the Digamma Function, a table of calculated values will be presented below. In this table, various evaluations of the Hurwitz Zeta Function, the Riemann Zeta Function, and the Digamma Function are obtained using the lnd(n,Δx,x) function and other calculation methods.

The use of the Section 8.2 derived equations together with the written lnd(n,Δx,x) calculation program, LNDX, is quite efficient in obtaining function evaluations. There are other computer programs available such as Mathematica which will calculate the Zeta Functions, Polygamma Functions, Digamma Function, etc. However, the use of the lnd(n,Δx,x) function together with Interval Calculus can provide a mathematical advantage. For several demonstrations, see Example 8.3-1 thru Example 8.3-4 at the end of this section. There are many relationships involving the lnd(n,Δx,x) function specified in the Appendix that can be used to simplify calculations involving summations and other discrete mathematical functions.

Consider the following three equations involving the lnd(n,Δx,x) function.

The Hurwitz Zeta Function is:

$\zeta(n,x) = \zeta(n,1,x) = \text{lnd}(n,1,x)\ ,\ \ n \neq 1$ (8.3-1)

where
x = real or complex variable
x,n = real or complex constants
Any summation term where x = 0 is excluded

The Riemann Zeta Function is:

$\zeta(n) = \zeta(n,1,1) = \text{lnd}(n,1,1)$ , $n \neq 1$ (8.3-2)

where
- $x$ = real or complex variable
- $n$ = real or complex constants
- Any summation term where $x = 0$ is excluded

The Digamma Function

$\psi(x) = \zeta(1,1,x) - \gamma = \text{lnd}(1,1,x) - \gamma$ (8.3-3)

where
- $x$ = real or complex values
- $\gamma$ = Euler's Constant, .5772157…

The above three functions will be evaluated in two different ways. The functions will be evaluated using the lnd(n,Δx,x) function, as specified above in the three function equations, and also by using tables or a computer program. The evaluations obtained are placed in Table 8.3-1.

The functional characteristics of the lnd(n,Δx,x) function change for different sets of n values. These sets are specified below and include all real and complex values of n.

<u>Sets of n values</u>

$\text{Re}(n) > 1$
$0 \leq \text{Re}(n) \leq 1$ and $n \neq 0$ and $n \neq 1$
$n = 1$
$n = 0$
$\text{Re}(n) < 0$

In Table 8.3-1 on the following page, three values are calculated using the lnd(n,Δx,x) function for each of the five sets of n listed above. The accuracy of these calculations is then checked using calculations obtained from tables or a computer program.

**Table 8.3-1 Checking the Hurwitz Zeta Function, Riemann Zeta Function, and Digamma lnd(n,Δx,x) function evaluation equations, Eq 8.3-1 thru Eq 8.3-3, Δx = 1**

| # | n | x | Function evaluations obtained using the lnd(n,Δx,x) function | Function evaluations obtained using tables or a computer program | N value set | Check | Function |
|---|---|---|---|---|---|---|---|
| **1** | 2 | 1 | 1.6449340668 | 1.6449340668 | Re(n)>1 | Good | Riemann Zeta Function |
| **2** | 2+i | 3 | .802461390929-.607367926173i | .802461390929-.607367926173i | Re(n)>1 | Good | Hurwitz Zeta Function |
| **3** | 1.1-2i | 1.5 | -.032339285896+.564141251383i | -.032339285896+.564141251383i | Re(n)>1 | Good | Hurwitz Zeta Function |
| **4** | 1+i | 2 | -.417841940247-.926848564330i | -.417841940247-.926848564330i | 0 ≤ Re(n) ≤ 1 and n ≠ 0 and n ≠ 1 | Good | Hurwitz Zeta Function |
| **5** | i | .5 | .266423576237+.098602725982i | .266423576237+.098602725982i | 0 ≤ Re(n) ≤ 1 and n ≠ 0 and n ≠ 1 | Good | Hurwitz Zeta Function |
| **6** | .5 | 1.5 | -2.019112205794 | -2.019112205794 | 0 ≤ Re(n) ≤ 1 and n ≠ 0 and n ≠ 1 | Good | Hurwitz Zeta Function |
| **7** | 1 | 1.5 | .0364899740 | .0364899740 | n = 1 | Good | Digamma Function |
| **8** | 1 | 96 | 4.5591308160 | 4.5591308160 | n = 1 | Good | Digamma Function |
| **9** | 1 | 1+4.4i | 1.48593+1.45716i | 1.48593+1.45716i | n = 1 | Good | Digamma Function |
| **11** | 0 | 1 | -.5 | -.5 | n = 0 | Good | Riemann Zeta Function |
| **12** | 0 | 2.5-i | -2+i | -2+i | n = 0 | Good | Hurwitz Zeta Function |
| **13** | 0 | 3.2+i | -2.7-i | -2.7-i | n = 0 | Good | Hurwitz Zeta Function |
| **14** | -1 | 1 | -.083333333333 | -.083333333333 | Re(n) < 0 | Good | Riemann Zeta Function |
| **15** | -4 | -1 | 1 | 1 | Re(n) < 0 | Good | Hurwitz Zeta Function |
| **16** | -1-2.1i | .5 | -.141276933256-.146828469024i | -.141276933256-.146828469024i | Re(n) < 0 | Good | Hurwitz Zeta Function |

Note that the lnd(n,Δx,x) function obtains valid results for all five sets of n.

Table 8.3-1 is a simple demonstation of the usefulness and accuracy of the discrete Interval Calculus function, lnd(n,Δx,x), when dealing with advanced mathematical functions such as the Hurwitz Zeta Function, the Riemann Zeta Function, and the Digamma Function. It should be pointed out that the lnd(n,Δx,x) function is defined and functional for all real or complex values of x, Δx, and n including the very special case where n = 1.

The lnd(n,Δx,x) function, being an Interval Calculus function, can be conveniently mathematically manipulated to obtain desired calculation results. This will be demonstated in the following four examples.

**Example 8.1** Express the Hurwitz Zeta Function, $\zeta(n,\frac{3}{2})$ in terms of the Riemann Zeta Function, $\zeta(n)$.

$$\sum_{x=1}^{\infty}{}_{1}\frac{1}{x^n} = \frac{1}{(1)^n} + \frac{1}{(2)^n} + \frac{1}{(3)^n} + \frac{1}{(4)^n} + \frac{1}{(5)^n} + \frac{1}{(6)^n} + \frac{1}{(7)^n} + \frac{1}{(8)^n} + \frac{1}{(9)^n} + \ldots \tag{8.3-4}$$

$$\sum_{x=1}^{\infty}{}_{1}\frac{1}{x^n} = 1 + \left[\frac{1}{(2)^n} + \frac{1}{(4)^n} + \frac{1}{(6)^n} + \frac{1}{(8)^n} + \ldots\right] + \left[\frac{1}{(3)^n} + \frac{1}{(5)^n} + \frac{1}{(7)^n} + \frac{1}{(9)^n} + \ldots\right] \tag{8.3-5}$$

$$\sum_{x=1}^{\infty}{}_{1}\frac{1}{x^n} = 1 + \sum_{x=2}^{\infty}{}_{2}\frac{1}{x^n} + \sum_{x=3}^{\infty}{}_{2}\frac{1}{x^n} \tag{8.3-6}$$

$$\sum_{x=1}^{\infty}{}_{\Delta x}\frac{1}{x^n} = \frac{1}{\Delta x}\,\mathrm{lnd}(n,\Delta x,x) \tag{8.3-7}$$

From Eq 8.3-6 and Eq 8.3-7

$$\frac{1}{1}\,\mathrm{lnd}(n,1,1) = 1 + \frac{1}{2}\,\mathrm{lnd}(n,2,2) + \frac{1}{2}\,\mathrm{lnd}(n,2,3) \tag{8.3-8}$$

$$\mathrm{lnd}(n,\Delta x,x) = \frac{1}{\Delta x^{n-1}}\,\mathrm{lnd}(n,1,\frac{x}{\Delta x}), \quad x,\Delta x = \text{real values with the condition that} \tag{8.3-9}$$

$x \neq$ positive real value when $\Delta x$ = negative real value

From Eq 8.3-8 and Eq 8.3-9

$$\mathrm{lnd}(n,1,1) = 1 + \frac{1}{2(2)^{n-1}}\,\mathrm{lnd}(n,1,1) + \frac{1}{2(2)^{n-1}}\,\mathrm{lnd}(n,1,\frac{3}{2}) \tag{8.3-10}$$

$$\text{lnd}(n,1,1) = 1 + \frac{1}{(2)^n}\text{lnd}(n,1,1) + \frac{1}{(2)^n}\text{lnd}(n,1,\tfrac{3}{2}) \qquad (8.3\text{-}11)$$

Simplifying

$$\text{lnd}(n,1,\tfrac{3}{2}) = (2)^n\text{lnd}(n,1,1) - \text{lnd}(n,1,1) - (2)^n \qquad (8.3\text{-}12)$$

$$\text{lnd}(n,1,\tfrac{3}{2}) = (2^n - 1)\text{lnd}(n,1,1) - (2)^n \qquad (8.3\text{-}13)$$

$\zeta(n,x) = \zeta(n,1,x) = \text{lnd}(n,1,x)$ , $n \neq 1$ , The Hurwitz Zeta Function evaluation equation (8.3-14)

$\zeta(n) = \zeta(n,1,1) = \text{lnd}(n,1,1)$ , The Riemann Zeta Function evaluation equation (8.3-15)

Substituting Eq 8.3-14 and Eq 8.3-15 into Eq 8.3-13

$$\boldsymbol{\zeta(n,\tfrac{3}{2}) = (2^n - 1)\zeta(n) - (2)^n} \qquad (8.3\text{-}16)$$

Checking

The above equation can be found on the internet at wolfram.com (http://functions.wolfram.com/ZetaFunctionsandPolylogarithms)

**<u>Example 8.2</u>** Evaluate the summation, $\sum_{k=0}^{\infty}\frac{(-1)^k}{zk+1}$, in terms of the Digamma Function where z is real and z > 0.

$$\sum_{k=0}^{\infty}\frac{(-1)^k}{zk+1} = \frac{1}{1} - \frac{1}{z+1} + \frac{1}{2z+1} - \frac{1}{3z+1} + \frac{1}{4z+1} - \frac{1}{5z+1} + \frac{1}{6z+1} - \frac{1}{7z+1} + \frac{1}{8z+1} - \ldots \qquad (8.3\text{-}17)$$

$$\sum_{k=0}^{\infty}\frac{(-1)^k}{zk+1} = \left[\, 1 + \frac{1}{2z+1} + \frac{1}{4z+1} + \frac{1}{6z+1} + \frac{1}{8z+1} + \ldots \right] - \left[\, \frac{1}{z+1} + \frac{1}{3z+1} + \frac{1}{5z+1} + \frac{1}{7z+1} + \ldots \right] \qquad (8.3\text{-}18)$$

From Eq 8.3-18

$$\sum_{k=0}^{\infty}\frac{(-1)^k}{zk+1} = \lim_{N\to\infty}\left[\ {}_{2z}\sum_{x=0}^{N}\frac{1}{x+1} - {}_{2z}\sum_{x=z}^{N}\frac{1}{x+1}\ \right] \qquad (8.3\text{-}19)$$

$$\sum_{\Delta x \; x=x_1}^{x_2} \frac{1}{x^n} = \pm \frac{1}{\Delta x} \mathrm{lnd}(n,\Delta x,x) \Big|_{x_1}^{x_2} \quad , \quad + \text{ for } n=1, \; - \text{ for } n\neq 1 \tag{8.3-20}$$

From Eq 8.3-19 and Eq 8.3-20

$$\sum_{k=0}^{\infty} \frac{(-1)^k}{zk+1} = \lim_{N\to\infty}\Big[ \frac{1}{2z} \mathrm{lnd}(1,2z,x+1) \Big|_0^N - \frac{1}{2z} \mathrm{lnd}(1,2z,x+1) \Big|_z^N \Big] \tag{8.3-21}$$

$$\sum_{k=0}^{\infty} \frac{(-1)^k}{zk+1} = \lim_{N\to\infty}\Big[ \frac{1}{2z} \mathrm{lnd}(1,2z,N+1) - \frac{1}{2z} \mathrm{lnd}(1,2z,1) - \frac{1}{2z} \mathrm{lnd}(1,2z,N+1) + \frac{1}{2z} \mathrm{lnd}(1,2z,z+1)\Big] \tag{8.3-22}$$

$$\boldsymbol{\sum_{k=0}^{\infty} \frac{(-1)^k}{zk+1} = -\frac{1}{2z} \mathrm{lnd}(1,2z,1) + \frac{1}{2z} \mathrm{lnd}(1,2z,z+1)} \tag{8.3-23}$$

$$\mathrm{lnd}(n,\Delta x,x) = \frac{1}{\Delta x^{\,n-1}} \mathrm{lnd}(n,1,\frac{x}{\Delta x}) \;, \quad x,\Delta x = \text{real values with the condition that} \tag{8.3-24}$$

$x \neq$ positive real value when $\Delta x$ = negative real value

From Eq 8.3-23 and Eq 8.3-24 where $\Delta x = 2z$, $n=1$

$$\sum_{k=0}^{\infty} \frac{(-1)^k}{zk+1} = -\frac{1}{2z} \Big[ \mathrm{lnd}(1,1,\frac{1}{2z}) - \mathrm{lnd}(1,1,\frac{z+1}{2z}) \Big] \tag{8.3-25}$$

$$\sum_{k=0}^{\infty} \frac{(-1)^k}{zk+1} = -\frac{1}{2z} \Big[ \mathrm{lnd}(1,1,\frac{1}{2z}) - \gamma - \mathrm{lnd}(1,1,\frac{z+1}{2z}) + \gamma \Big] \tag{8.3-26}$$

$$\psi(x) = \zeta(1,1,x) - \gamma = \mathrm{lnd}(1,1,x) - \gamma \;, \quad x \neq 0,-1,-2,-3,\dots , \qquad \text{The Digamma Function} \tag{8.3-27}$$

From Eq 8.3-26 and Eq 8.3-27

$$\boldsymbol{\sum_{k=0}^{\infty} \frac{(-1)^k}{zk+1} = -\frac{1}{2z} \Big[ \psi(\frac{1}{2z}) - \psi(\frac{z+1}{2z}) \Big]} \tag{8.3-28}$$

Checking

The above equation can be found on the internet at Wolfram MathWorld (http://mathworld.wolfram.com/DigammFunction.html)

**Example 8.3** Evaluate the summation, $\sum_{k=0}^{\infty}\frac{(-1)^k}{(zk+1)^n}$, where $n\neq 1$, z is real, z > 0 in terms of the Hurwitz Zeta Function.

$$\sum_{k=0}^{\infty}\frac{(-1)^k}{(zk+1)^n} = \frac{1}{(1)^n}-\frac{1}{(z+1)^n}+\frac{1}{(2z+1)^n}-\frac{1}{(3z+1)^n}+\frac{1}{(4z+1)^n}-\frac{1}{(5z+1)^n}+\frac{1}{(6z+1)^n}-\frac{1}{(7z+1)^n}+\frac{1}{(8z+1)^n}-\dots \tag{8.3-29}$$

$$\sum_{k=0}^{\infty}\frac{(-1)^k}{(zk+1)^n} = [\,1+\frac{1}{(2z+1)^n}+\frac{1}{(4z+1)^n}+\frac{1}{(6z+1)^n}+\frac{1}{(8z+1)^n}+\dots\,]-[\,\frac{1}{(z+1)^n}+\frac{1}{(3z+1)^n}+\frac{1}{(5z+1)^n}+\frac{1}{(7z+1)^n}+\dots\,] \tag{8.3-30}$$

From Eq 8.3-30

$$\sum_{k=0}^{\infty}\frac{(-1)^k}{(zk+1)^n} = [\,{}_{2z}\sum_{x=0}^{\infty}\frac{1}{(x+1)^n} - {}_{2z}\sum_{x=z}^{\infty}\frac{1}{(x+1)^n}\,] \tag{8.3-31}$$

$${}_{\Delta x}\sum_{x=x_i}^{\infty}\frac{1}{(x+a)^n} = \frac{1}{\Delta x}\,\mathrm{lnd}(n,\Delta x,x_i+a)\,,\quad n\neq 1 \tag{8.3-32}$$

From Eq 8.3-31 and Eq 8.3-32 where $\Delta x = 2z$, $a = 1$

$$\sum_{k=0}^{\infty}\frac{(-1)^k}{(zk+1)^n} = [\,\frac{1}{2z}\,\mathrm{lnd}(n,2z,x+1)\Big|_{x=0} - \frac{1}{2z}\,\mathrm{lnd}(n,2z,x+1)\Big|_{x=z}\,] \tag{8.3-33}$$

$$\sum_{k=0}^{\infty}\frac{(-1)^k}{(zk+1)^n} = [\,\frac{1}{2z}\,\mathrm{lnd}(n,2z,1) - \frac{1}{2z}\,\mathrm{lnd}(n,2z,z+1)] \tag{8.3-34}$$

$$\boldsymbol{\sum_{k=0}^{\infty}\frac{(-1)^k}{(zk+1)^n} = \frac{1}{2z}\,\mathrm{lnd}(n,2z,1) + \frac{1}{2z}\,\mathrm{lnd}(n,2z,z+1)} \tag{8.3-35}$$

$$\mathrm{lnd}(n,\Delta x,x) = \frac{1}{\Delta x^{n-1}}\,\mathrm{lnd}(n,1,\frac{x}{\Delta x})\,,$$ $x,\Delta x$ = real values with the condition that $x \neq$ positive real value when $\Delta x$ = negative real value (8.3-36)

From Eq 8.3-35 and Eq 8.3-36 where $\Delta x = 2z$

$$\sum_{k=0}^{\infty} \frac{(-1)^k}{(zk+1)^n} = \frac{1}{2z}\left(\frac{1}{2z}\right)^{n-1} \left[ \operatorname{lnd}(n,1,\frac{1}{2z}) - \operatorname{lnd}(n,1,\frac{z+1}{2z}) \right] \tag{8.3-37}$$

$$\sum_{k=0}^{\infty} \frac{(-1)^k}{(zk+1)^n} = \left(\frac{1}{2z}\right)^{n} \left[ \operatorname{lnd}(n,1,\frac{1}{2z}) - \operatorname{lnd}(n,1,\frac{z+1}{2z}) \right] \tag{8.3-38}$$

$\zeta(n,x) = \zeta(n,1,x) = \operatorname{lnd}(n,1,x)$ , The Hurwitz Zeta Function evaluation equation (8.3-39)

From Eq 8.3-38 and Eq 8.3-39

$$\boldsymbol{\sum_{k=0}^{\infty} \frac{(-1)^k}{(zk+1)^n} = \left(\frac{1}{2z}\right)^{n} \left[ \zeta(n,\frac{1}{2z}) - \zeta(n,\frac{z+1}{2z}) \right]} \tag{8.3-40}$$

Check Eq 8.3-40

Let z = 2, n = 3

$$\sum_{k=0}^{\infty} \frac{(-1)^k}{(2k+1)^3} = \left(\frac{1}{4}\right)^{3} \left[ \zeta(3,\frac{1}{4}) - \zeta(3,\frac{3}{4}) \right] = \frac{1}{64} \left[ \zeta(3,\frac{1}{4}) - \zeta(3,\frac{3}{4}) \right]$$

$$\sum_{k=0}^{\infty} \frac{(-1)^k}{(2k+1)^3} = \frac{1}{64} \left[ 64.6638699687684601 - 2.6513166081688198 \right]$$

$$\sum_{k=0}^{\infty} \frac{(-1)^k}{(2k+1)^3} = .9689461462593693$$

Checking the above value using a computer summation program

$$\sum_{k=0}^{1000000} \frac{(-1)^k}{(2k+1)^3} = .9689461462593693$$

Good check

**Note – The summation equations, Eq 8.3-28 and Eq 8.3-40, appear at first glance to be decidedly different. This is because these two equations involve two different functions, the Digamma Function and the Hurwitz Zeta Function. However, if the summations are expressed in terms of the lnd(n,Δx,x) function, these summations are seen to be closely related. From Eq 8.3-23 and Eq 8.3-35 the following single equation can be written.**

$$\sum_{k=0}^{\infty} \frac{(-1)^k}{(zk+1)^n} = \pm\left[ \frac{1}{2z}\,\mathrm{lnd}(n,2z,1) - \frac{1}{2z}\,\mathrm{lnd}(n,2z,z+1) \right] \text{, + for } n\neq 1\text{ , – for } n=1 \qquad (8.3\text{-}41)$$

**Example 8.4** Derive an equation to evaluate the Hurwitz Zeta Function, ζ(0,x), using Interval Calculus.

$$D_{\Delta x}\mathrm{lnd}(n,\Delta x,x) = \pm\frac{1}{x^n} \text{ , + for } n=1\text{ , – for } n\neq 1 \qquad (8.3\text{-}42)$$

$n = 0\ \ (n\neq 1)$
$\Delta x = 1$

$$D_1\mathrm{lnd}(0,1,x) = D_1\zeta(0,x) = -\frac{1}{x^0} \qquad (8.3\text{-}43)$$

$$D_1\zeta(0,x) = -1 \qquad (8.3\text{-}44)$$

Using discrete integration to integrate Eq 8.3-44

$$\zeta(0,x) = -{}_1\!\int 1\,\Delta x + k = -x + k \qquad (8.3\text{-}45)$$

$$\zeta(0,x) = -x + k \qquad (8.3\text{-}46)$$

Find the constant of integration, k

Let x = 2

Using the LNDX computer program to calculate ζ(0,2) = lnd(0,1,2)

$$\zeta(0,2) = \mathrm{lnd}(0,1,2) = -1.5 \qquad (8.3\text{-}47)$$

From Eq 8.3-46 and Eq 8.3-47

$$k = \zeta(0,2) - 2 = -1.5 + 2 = .5 \qquad (8.3\text{-}48)$$

$$k = \frac{1}{2} \qquad (8.3\text{-}49)$$

Substituting Eq 8.3-49 into Eq 8.3-46

$$\boldsymbol{\zeta(0,x) = \frac{1}{2} - x} \qquad (8.3\text{-}50)$$

<u>Comment</u> - The above equation can be found on the internet at Wolfram MathWorld (http://mathworld.wolfram.com/HurwitzZetaFunction.html)

**<u>Example 8.5</u>** Derive an equation to evaluate the Hurwitz Zeta Function, ζ(-2,x), using Interval Calculus.

$$D_{\Delta x}\text{lnd}(n,\Delta x,x) = \pm\frac{1}{x^n}\ ,\ \ + \text{ for } n=1\ ,\ \ - \text{ for } n\neq 1 \qquad (8.3\text{-}51)$$

$n = -2\ \ (n\neq 1)$
$\Delta x = 1$

$$D_1\text{lnd}(n,1,x) = D_1\zeta(-2,x) = -\frac{1}{x^{-2}} \qquad (8.3\text{-}52)$$

$$D_1\zeta(-2,x) = -x^2 \qquad (8.3\text{-}53)$$

Using discrete integration to integrate Eq 8.3-53

$$\zeta(-2,x) = -{}_1\!\int x^2\,\Delta x + k = -{}_1\!\int [x(x-1) + x]\,\Delta x + k \qquad (8.3\text{-}54)$$

$$\zeta(-2,x) = -{}_1\!\int x^2\,\Delta x + k = -\frac{x(x-1)(x-2)}{3} - \frac{x(x-1)}{2} + k \qquad (8.3\text{-}55)$$

$$\zeta(-2,x) = -\frac{x(x-1)(x-2)}{3} - \frac{x(x-1)}{2} + k \qquad (8.3\text{-}56)$$

Find the constant of integration, k

$$\text{Let } x = 0 \qquad (8.3\text{-}57)$$

Using the LNDX computer program to calculate ζ(-2,0) = lnd(-2,1,0)

$$\zeta(-2,0) = \text{lnd}(-2,1,0) = 0 \qquad (8.3\text{-}58)$$

Substituting Eq 8.3-58 and Eq 8.3-57 into Eq 8.3-56

$$k = \zeta(-2,0) + \frac{0(0-1)(0-2)}{3} + \frac{0(0-1)}{2} = 0 + 0 + 0 = 0 \qquad (8.3\text{-}59)$$

$$k = 0 \qquad (8.3\text{-}60)$$

Substituting Eq 8.3-60 into Eq 8.3-56

$$\boldsymbol{\zeta(-2,x) = -\frac{x(x-1)(x-2)}{3} - \frac{x(x-1)}{2}} \qquad (8.3\text{-}61)$$

Checking Eq 8.3-61

Let x = 1+ i

Substituting into Eq 8.3-61

$$\zeta(-2,1+i) = -\frac{(1+i)(1+i-1)(1+i-2)}{3} - \frac{(1+i)(1+i-1)}{2} = \frac{1}{2} + \frac{1}{6}i$$

$$\zeta(-2,1+i) = \frac{1}{2} + \frac{1}{6}i$$

Using the LNDX computer program to calculate ζ(-2,1+i) = lnd(-2,1,1+i)

$$\zeta(-2,1+i) = \frac{1}{2} + \frac{1}{6}i$$

Good Check

Comment – The derivation of Eq 8.3-61 could have been obtained more directly from the derived general equation, $\mathrm{lnd}(-2,\Delta x,x) = -\frac{x(x-\Delta x)(x-2\Delta x)}{3} - \frac{x(x-\Delta x)}{2}$, where Δx=1 and ζ(-2,x) = lnd(-2,1,x).

**Example 8.6** Evaluate four summations of the function, $\frac{1}{x^n}$, using Zeta Function and related equations.

Summation #1

Evaluate the summation, $\underset{1+i}{\sum_{x=-3-1.2i}^{4+5.8i}} \frac{1}{x^{.5}}$ , using the General Zeta Function summation between finite limits equation.

The General Zeta Function summation between finite limits equation is:

$$\underset{\Delta x}{\sum_{x=x_1}^{x_2}} \frac{1}{x^n} = \pm [\, \zeta(n,\Delta x,x_2+\Delta x) - \zeta(n,\Delta x,x_1) \,] = \pm \frac{1}{\Delta x} [\mathrm{lnd}(n,\Delta x,x_2+\Delta x) - \mathrm{lnd}(n,\Delta x,x_1)] \qquad (8.3\text{-}62)$$

where
- $+$ for $n = 1$, $-$ for $n \neq 1$
- $x$ = real or complex variable
- $x_1, x_2, \Delta x, n$ = real or complex constants
- Any summation term where $x = 0$ is excluded

Let $x_1 = -3-1.2i$
$x_2 = 4+5.8i$
$\Delta x = 1+i$
$n = .5$

$$\underset{1+i}{\sum_{x=-3-1.2i}^{4+5.8i}} \frac{1}{x^{.5}} = \frac{1}{1+i} [-\mathrm{lnd}(.5,1+i,5+6.8i) + \mathrm{lnd}(.5,1+i,-3-1.2i)] \qquad (8.3\text{-}63)$$

$$\underset{1+i}{\sum_{x=-3-1.2i}^{4+5.8i}} \frac{1}{x^{.5}} = \frac{1}{1+i} [\, 4.9500567561719026 + 2.5466870715060025i \,] + \frac{1}{1+i} [\, -1.3292804270644887 - 1.0709836690132466i \,] \qquad (8.3\text{-}64)$$

$$\boldsymbol{\underset{1+i}{\sum_{x=-3-1.2i}^{4+5.8i}} \frac{1}{x^{.5}} = 2.5482398658000848 - 1.0725364633073289i} \qquad (8.3\text{-}65)$$

checking

$$\sum_{x=-3-1.2i}^{4+5.8i}{}_{1+i}\frac{1}{x^{.5}} = \frac{1}{(-3-1.2i)^{.5}} + \frac{1}{(-2-.2i)^{.5}} + \frac{1}{(-1+.8i)^{.5}} + \frac{1}{(1.8i)^{.5}} + \frac{1}{(1+2.8i)^{.5}} +$$

$$\frac{1}{(2+3.8i)^{.5}} + \frac{1}{(3+4.8i)^{.5}} + \frac{1}{(4+5.8i)^{.5}} \qquad (8.3\text{-}66)$$

$$\sum_{x=-3-1.2i}^{4+5.8i}{}_{1+i}\frac{1}{x^{.5}} = 2.5482398658000848 - 1.0725364633073289i \qquad (8.3\text{-}67)$$

Good check

Summation #2

Evaluate the summation, $\sum_{x=-4}^{6}{}_{1}\frac{1}{x^{.5}}$ , using the Hurwitz Zeta Function summation between finite limits equation.

The Hurwitz Zeta Function summation between finite limits equation is:

$$\sum_{x=x_1}^{x_2}{}_{1}\frac{1}{x^{n}} = -\zeta(n,x_2+1) + \zeta(n,x_1) \qquad (8.3\text{-}68)$$

where
x = real or complex variable
$x_1, x_2, n$ = real or complex constants
Any summation term where x = 0 is excluded

Let $x_1 = -4$
$x_2 = 6$
$n = .5$

$$\sum_{x=-4}^{6}{}_{1}\frac{1}{x^{.5}} = -\zeta(.5,7) + \zeta(.5,-4) \qquad (8.3\text{-}69)$$

$$ {}_{1}\sum_{x=-4}^{6}\frac{1}{x^{.5}} = -(-5.1002734451495810) + (-1.4603545088095868-2.7844570503761732i) \qquad (8.3\text{-}70) $$

$$ \boldsymbol{{}_{1}\sum_{x=-4}^{6}\frac{1}{x^{.5}} = 3.6399189363399942-2.7844570503761732i} \qquad (8.3\text{-}71) $$

checking

$$ {}_{1}\sum_{x=-4}^{6}\frac{1}{x^{.5}} = \frac{1}{(-4)^{.5}} + \frac{1}{(-3)^{.5}} + \frac{1}{(-2)^{.5}} + \frac{1}{(-1)^{.5}} + \frac{1}{(1)^{.5}} + \frac{1}{(2)^{.5}} + \frac{1}{(3)^{.5}} + \frac{1}{(4)^{.5}} + \frac{1}{(5)^{.5}} + \frac{1}{(6)^{.5}} $$

<u>Note</u> - The x = 0 summation term has been excluded

$$ {}_{1}\sum_{x=-4}^{6}\frac{1}{x^{.5}} = 3.6399189363399942-2.7844570503761732i $$

Good check

<u>Summation #3</u>

Evaluate the summation, ${}_{1}\sum_{x=-5.09}^{.91}\frac{1}{x}$ , using the Digamma Function summation between finite limits equation.

The Digamma Function summation between finite limits equation is:

$$ {}_{\Delta x}\sum_{x=x_1}^{x_2}\frac{1}{x} = \frac{1}{\Delta x}\left[\psi\left(\frac{x_2}{\Delta x}+1\right) - \psi\left(\frac{x_1}{\Delta x}\right)\right] \qquad (8.3\text{-}72) $$

where
- $x$ = real or complex variable
- $x_1, x_2, \Delta x$ = real or complex constants
- $x_1 \neq 0, -\Delta x, -2\Delta x, -3\Delta x, \ldots$
- $x_2 \neq -\Delta x, -2\Delta x, -3\Delta x, \ldots$
- $x_2 = x_1 + m\Delta x$ ,   $m = 0,1,2,3,\ldots$
- $\Delta x$ = x increment

Let $x_1 = -5.09$
$x_2 = .91$
$\Delta x = 1$

$$ {}_{1}\sum_{x=-5.09}^{.91} \frac{1}{x} = \psi(1.91) - \psi(-5.09) \tag{8.3-73}$$

$$ {}_{1}\sum_{x=-5.09}^{.91} \frac{1}{x} = .3630410646488811 - 12.5357381962403155 \tag{8.3-74}$$

$$ \mathbf{ {}_{1}\sum_{x=-5.09}^{.91} \frac{1}{x} = -12.1726971315914304} \tag{8.3-75}$$

checking

$$ {}_{1}\sum_{x=-5.09}^{.91} \frac{1}{x} = \frac{1}{-5.09} + \frac{1}{-4.09} + \frac{1}{-3.09} + \frac{1}{-2.09} + \frac{1}{-1.09} + \frac{1}{-.09} + \frac{1}{.91}$$

$$ {}_{1}\sum_{x=-5.09}^{.91} \frac{1}{x} = -12.1726971315914304$$

Good check

Summation #4

Evaluate the summation, ${}_{1}\sum_{x=-3}^{4} \frac{1}{x}$ , using the General Zeta Function summation between finite limits equation. Exclude any summation term where x = 0.

The General Zeta Function summation between finite limits equation is:

$$\sum_{\Delta x \; x=x_1}^{x_2} \frac{1}{x^n} = \pm [\,\zeta(n,\Delta x,x_2+\Delta x) - \zeta(n,\Delta x,x_1)\,] = \pm \frac{1}{\Delta x} [\mathrm{lnd}(n,\Delta x,x_2+\Delta x) - \mathrm{lnd}(n,\Delta x,x_1)] \qquad (8.3\text{-}76)$$

where
- \+ for $n = 1$, $-$ for $n \neq 1$
- $x$ = real or complex variable
- $x_1, x_2, \Delta x, n$ = real or complex constants
- Any summation term where $x = 0$ is excluded

Let $x_1 = -3$
$x_2 = 4$
$\Delta x = 1$
$n = 1$

$$\sum_{1 \; x=-3}^{4} \frac{1}{x} = \mathrm{lnd}(1,1,5) - \mathrm{lnd}(1,1,-3) = 2\frac{1}{12} - 1\frac{5}{6} = \frac{25}{12} - \frac{22}{12} = \frac{3}{12} = \frac{1}{4} \qquad (8.3\text{-}77)$$

$$\boldsymbol{\sum_{1 \; x=-3}^{4} \frac{1}{x} = \frac{1}{4}} \qquad (8.3\text{-}78)$$

checking

$$\sum_{1 \; x=-3}^{4} \frac{1}{x} = \frac{1}{-3} + \frac{1}{-2} + \frac{1}{-1} + \frac{1}{1} + \frac{1}{2} + \frac{1}{3} + \frac{1}{4}$$

<u>Note</u> - The $x = 0$ summation term has been excluded

$$\sum_{1 \; x=-3}^{4} \frac{1}{x} = \frac{1}{4}$$

Good check

## A Final Comment

When the work on this paper was begun, the objective was to see what Calculus would be like if its independent variable was discrete. The resulting calculus would be applicable to discrete mathematics. In order not to confuse Calculus (where $\Delta x$ is infinitessimal) with this new calculus (where $\Delta x$ is finite), it was necessary to select for this new calculus a unique name. The name, Interval Calculus, was chosen. During the early years of development, it was believed that Calculus and its related Interval Calculus were similar but different. Calculus was a calculus with a continuous independent variable and Interval Calculus was a related calculus with a discrete independent variable. Both calculi were related by the fact that when the $\Delta x$ value of Interval Calculus was made to approach zero, Interval Calculus became Calculus. This characteristic of Interval Calculus was, of course, interesting. Early in its development, Interval Calculus was not considered to be a generalization of Calculus. It was thought to be a related calculus. However, later in the development of Interval Calculus, something surprising was observed. The Laplace Transform of the Calculus function, $e^{ax}$, and the $K_{\Delta x}$ Transform of the Interval Calculus functions, $e_{\Delta x}(a,x)$, was seen to have the same transform, $\frac{1}{s\text{-}a}$ . It was known that the $K_{\Delta x}$ Transform became the Laplace Transform if $\Delta x$ approached zero and that the $\lim_{\Delta x\to 0} e_{\Delta x}(a,x) = e^{ax}$. Could this mean that Calculus is a subset of Interval Calculus where $\Delta x$ is a very small value (i.e. $\Delta x\to 0$)? Comparing the Laplace Transforms of sin(ax) and cos(ax) to the $K_{\Delta x}$ Transforms of $\sin_{\Delta x}(a,x)$ and $\cos_{\Delta x}(a,x)$ respectively, their transforms were also found to be same. In fact, it has been found that all Interval Calculus functions, when $\Delta x\to 0$, become their related Calculus functions and have the same transform. As a result, it appears that Calculus is a subset of Interval Calculus based on the similarity of transforms. If this is true, it would mean that Calculus is a discrete mathematics as is Interval Calculus. Could this be so? Probably it is. To try to clarify this, consider the Calculus first order derivative,

$$\frac{dy(x)}{dx} = \lim_{\Delta x\to 0}\frac{\Delta y(x)}{\Delta x} = \lim_{\Delta x\to 0}\frac{y(x+\Delta x)-y(x)}{\Delta x}$$ which can be rewritten as

$$\frac{dy(x)}{dx} = \lim_{\Delta x\to 0}\frac{\Delta y(x)}{\Delta x} = \lim_{\Delta x\to \alpha}\frac{y(x+\Delta x)-y(x)}{\Delta x}$$ where $\alpha\to 0$ but $\alpha \neq 0$.

By this definition, where $\Delta x$ is a very small non-zero value, the Calculus derivative is a discrete derivative. It is strongly believed that Interval Calculus can correctly be called a generalization of Calculus.

# Appendix

## Interval Calculus Tables of Equations and Calculation Programs

# Table 1
## Interval Calculus Notation Definitions

| | Notation & Symbols | Description |
|---|---|---|
| **1** | x,s | Variables |
| **2** | Re(x) | real value of x |
| **3** | Im(x) | imaginary value of x |
| **4** | $x_1,x_2$ | values of x |
| **5** | $\Delta$ | difference operator |
| **6** | $\Delta x$ | x increment<br>$\Delta x = x_2 - x_1$ |
| **7** | $\Delta f(x)$ | $\Delta f(x) = f(x+\Delta x) - f(x)$ |
| **8** | f(x),g(x),u(x),v(x) | Continuous functions of x<br>x variable is continuous |
| **8a** | $f^*(x),g^*(x),u^*(x),v^*(x)$ | Discrete functions of x<br>Impulse function<br>or<br>Sample and hold shaped waveform function<br>x variable is discrete , often $x = 0,\Delta x,2\Delta x,3\Delta x, \ldots$<br><br>Comment – Discrete functions have values defined only at specific values of x or can change value only at specific values of x. |
| **8b** | F(s),G(s),F(z),G(z) | Transform functions of s,z representing continuous functions |
| **8c** | $F^*(s),G^*(s),F^*(z),G^*(z)$ | Transform functions of s,z representing discrete functions |
| **9** | a,b,c,k,$\beta$ | Constants |

| | Notation & Symbols | Description |
|---|---|---|
| **10** | $j = \sqrt{-1}$ | |
| **11** | $\sum$ , ${}_{\Delta x}\sum$ , ${}_{\Delta k}\oint\!\!\!\sum$ , ${}_{\Delta k}\circlearrowleft\!\!\!\sum$ | Interval Calculus summation symbols |
| **12** | $\sum_{x=x_1}^{x_2} f(x)$ | summation symbol with limits and $\Delta x = 1$<br>$\sum_{x=x_1}^{x_2} f(x) = f(x_1) + f(x_1+1) + f(x_1+2) + \ldots$<br>$+ f(x_2-1) + f(x_2)$<br><br>$x_2-x_1$ = integer |
| **13** | ${}_{\Delta x}\sum_{x=x_1}^{x_2} f(x)$ | summation symbol with limits and $\Delta x$<br><br>${}_{\Delta x}\sum_{x=x_1}^{x_2} f(x) = f(x_1) + f(x_1+\Delta x) + f(x_1+2\Delta x) + \ldots +$<br>$f(x_2-\Delta x) + f(x_2)$<br>${}_{\Delta x}\sum_{x=x_1}^{x_2} f(x) = \frac{1}{\Delta x}\ {}_{\Delta x}\int_{x_1}^{x_2} f(x)\, \Delta x$<br>$\frac{x_2-x_1}{\Delta x}$ = integer |
| **14** | ${}_{1}\oint\!\!\!\sum_{k=0}^{v-1} f(z_k,z_{k+1})$<br>where<br>$\Delta k = 1$<br>v = the number of counterclockwise pointing vectors comprising the complex plane closed contour being integrated<br>$z_k$ = points on the closed contour<br>$f(z_k,z_{k+1})$ = a function of consecutive points (vectors) on the closed contour<br>k = 0,1,2,3,…,v-1<br>$z_{k+1} - z_k$ = a contour grid square vector value, $z_{k+1}$ at the vector head, $z_k$ at the tail | Interval Calculus complex plane closed contour summation. The summation is performed in a counterclockwise direction. |

| | Notation & Symbols | Description |
|---|---|---|
| **15** | ${}_{1}\circlearrowright\!\!\!\!\sum_{k=0}^{v-1} f(z_k,z_{k+1})$<br>where<br>$\Delta k = 1$<br>v = the number of clockwise pointing vectors comprising the complex plane closed contour being integrated<br>$z_k$ = points on the closed contour<br>$f(z_k,z_{k+1})$ = a function of consecutive points (vectors) on the closed contour<br>$k = 0,1,2,3,\ldots,v-1$<br>$z_{k+1} - z_k$ = a contour grid square vector value, $z_{k+1}$ at the vector head, $z_k$ at the tail | Interval Calculus complex plane closed contour summation. The summation is performed in a clockwise direction. |
| **16** | $\prod$ | Product symbol |
| **17** | $\prod_{x=x_1}^{x_2} f(x)$ | Product symbol with limits and $\Delta x = 1$<br>$\prod_{x=x_1}^{x_2} f(x) = f(x_1)f(x_1+1)f(x_1+2)\ldots f(x_2-1)f(x_2)$ |
| **18** | $D_{\Delta x}$ | Interval Calculus discrete derivative operator<br>$D_{\Delta x} = \frac{\Delta}{\Delta x}$, $\Delta x$ = interval between x values<br>$\lim_{\Delta x\to 0} D_{\Delta x} = \frac{d}{dx}$ |
| **19** | $D_{\Delta x}f(x)$ | $D_{\Delta x}f(x) = \frac{f(x+\Delta x) - f(x)}{\Delta x} = \frac{\Delta f(x)}{\Delta x}$ |
| **20** | $D_{-\Delta x}$ | Interval Calculus discrete derivative operator ($-\Delta x$) |
| **21** | $D_{-\Delta x}f(x)$ | $D_{-\Delta x}f(x) = \frac{f(x-\Delta x) - f(x)}{\Delta x}$ |

| | Notation & Symbols | Description |
|---|---|---|
| **22** | $s = s_{\Delta x} = \dfrac{e^{s\Delta x} - 1}{\Delta x}$, $-\dfrac{\pi}{\Delta x} \le w < \dfrac{\pi}{\Delta x}$, $s = \gamma + jw$<br>where<br>$x = 0, \Delta x, 2\Delta x, 3\Delta x, \ldots$<br>$\gamma$ = any positive real value which makes s an indefinitely large value. | $K_{\Delta x}$ Transform variable, s |
| **23** | $\vartheta_{\Delta x}$ | Interval Calculus discrete partial derivative operator<br>$\lim_{\Delta x \to 0} \vartheta_{\Delta x} = \dfrac{\vartheta}{\vartheta x}$ |
| **24** | $\vartheta_{\Delta x} f(x,y)$ | $\vartheta_{\Delta x} f(x,y) = \dfrac{f(x+\Delta x,y) - f(x,y)}{\Delta x}$ |
| **25** | $H_{\Delta x}$ | Incremental $x+\Delta x$ operator |
| **26** | $H_{\Delta x} f(x)$ | $H_{\Delta x} f(x) = f(x+\Delta x)$ |
| **27** | $H_{-\Delta x}$ | Incremental $x-\Delta x$ operator |
| **28** | $H_{-\Delta x} f(x)$ | $H_{-\Delta x} f(x) = f(x-\Delta x)$ |
| **29** | $sf(x) = D_{\Delta x} f(x) = \dfrac{f(x+\Delta x)-f(x)}{\Delta x}$<br>For Interval Calculus<br>$sf(x) = \dfrac{df(x)}{dx} = \lim_{\Delta x \to 0} D_{\Delta x} f(x)$<br>$= \lim_{\Delta x \to 0} \dfrac{f(x+\Delta x)-f(x)}{\Delta x}$<br>For Calculus | Interval Calculus and Calculus discrete derivative transform operator (s) |
| **30** | $\rightarrow$<br>$a \rightarrow b$ | Variable change symbol<br>a approaches b |
| **31** | $\rightleftarrows$<br>$y(x) \rightleftarrows Y(x)$<br>(Discrete Function $\rightleftarrows$ Calculus Function) | Function exchange symbol for Related Functions<br>The function exchange is reversible<br>No equality is implied |

| | Notation & Symbols | Description |
|---|---|---|
| **32** | $_{\Delta x}\int$ , $_{\Delta k}\oint$ , $_{\Delta k}\oint$ | Interval Calculus integral symbols<br>Δx,Δk are the interval between successive discrete values of x,k<br><br><u>Comment</u> – The Interval Calculus integral is sometimes called a discrete integral |
| **33** | $_{\Delta x}\int f(x)\Delta x + k$ | Interval Calculus indefinite integral |
| **34** | $_{\Delta x}\int_{x_1}^{x_2} f(x)\Delta x \equiv {}_{\Delta x}\sum_{x=x_1}^{x=x_2-\Delta x} f(x)\Delta x$ | Interval Calculus definite integral |
| **35** | $[x]_{\Delta x}^{n}$ | Symbol for the Interval Calculus function,<br>$[x]_{\Delta x}^{n} = x(x-\Delta x)(x-2\Delta x)(x-3\Delta x)\ldots(x-[n-1]\Delta x)$<br>$n = 1,2,3,\ldots$ , the number of product terms<br>$[x]_{\Delta x}^{0} = 1$ |
| **36** | $MV[f(x)]\,\rvert_{x=x} = \lim_{\varepsilon\to 0}\left[\frac{f(x+\varepsilon) + f(x-\varepsilon)}{2}\right]$<br>where<br>$MV[f(x)]\,\rvert_{x=x}$ = the mean value of the function, f(x), calculated at *x* | Mean Value of a function<br><br>If f(*x*) is at a pole, $MV[f(x)]\,\rvert_{x=x}$ may be finite or infinite |

| | Notation & Symbols | Description |
|---|---|---|
| **37** | $_{1}\oint_{0}^{v} f(z_k,z_{k+1})\Delta k$<br>where<br>$\Delta k = 1$<br>v = the number of counterclockwise pointing vectors comprising the complex plane closed contour being integrated<br>$z_k$ = coordinate of the kth point on the closed contour<br>$f(z_k,z_{k+1})$ = a function of consecutive points (vectors) on the closed contour<br>k = 0,1,2,3,…,v , point designations<br>$z_{k+1} - z_k$ = a contour grid square vector value, $z_{k+1}$ at the vector head, $z_k$ at the tail<br>or | Interval Calculus complex plane closed contour integral. Integration is performed in a counterclockwise direction |
| | $_{\Delta z}\oint f(z)\Delta z$<br>where<br>z = the complex plane coordinates of the tail points of the closed contour counterclockwise pointing vectors<br>$\Delta z$ = real or imaginary value<br>$z + \Delta z$ = the coordinates of the head points of the closed contour counterclockwise pointing vectors<br>f(z) = a function of z<br><br>Should a numerical value be substituted for the integral symbol $\Delta z$, it would be the absolute value of $\Delta z$. | Interval Calculus complex plane closed contour integral. Integration is performed in a counter clockwise direction |

| | Notation & Symbols | Description |
|---|---|---|
| **38** | $\oint_{1\ 0}^{v} f(z_k,z_{k+1})\Delta k$<br>where<br>$\Delta k = 1$<br>v = the number of clockwise pointing vectors comprising the complex plane closed contour being integrated<br>$z_k$ = points on the closed contour<br>$f(z_k,z_{k+1})$ = a function of consecutive points (vectors) on the closed contour<br>k = 0,1,2,3,…,v<br>$z_{k+1} - z_k$ = a contour grid square vector value, $z_{k+1}$ at the vector head, $z_k$ at the tail<br>or<br><br>$\oint_{\Delta z} f(z)\Delta z$<br>where<br>z = the complex plane coordinates of the tail points of the closed contour clockwise pointing vectors<br>$\Delta z$ = real or imaginary value<br>$z + \Delta z$ = the coordinates of the head points of the closed contour clockwise pointing vectors<br>f(z) = a function of z<br><br>Should a numerical value be substituted for the integral symbol $\Delta z$, it would be the absolute value of $\Delta z$. | Interval Calculus complex plane closed contour integral. Integration is performed in a clockwise direction |

| | Notation & Symbols | Description |
|---|---|---|
| **39** | A = Area | Area = area under the f(x) curve between the limits $x_1$ and $x_2$ as defined by the sum $\Delta x \underset{\Delta x}{\sum_{x=x_1}^{x_2}} f(x)$ where if any f(x) term is $f(x) = \frac{1}{0}$ at a pole, this term is excluded<br>__________________________________________<br>Area with no poles within the integration interval<br>$A = {}_{\Delta x}\int_{x_1}^{x_2} f(x)\Delta x = \Delta x \underset{\Delta x}{\sum_{x=x_1}^{x_2-\Delta x}} f(x)$<br>$x = x_1, x_1+\Delta x, x_1+2\Delta x, x_1+3\Delta x, \dots, x_2-\Delta x, x_2$<br>Interval calculus integration of the area, A, under the function f(x) between the limits $x_1,x_2$<br>Note – for $Re(x_2-x_1) < 0$ the value of $Re(\Delta x) < 0$<br>for $Re(x_2-x_1) > 0$ the value of $Re(\Delta x) > 0$<br>This equation applies if there are no f(x) poles within $x_1 \le x \le x_2$<br>__________________________________________<br>Area with poles within the integration interval<br>All $\frac{1}{0}$ terms are excluded<br>$A = {}_{\Delta x}\int_{x_1}^{x_2} f(x)\Delta x - \{\Delta x \sum MV[f(x)]\,,\ x = x_1, x_2, x_3, \dots, x_m\}$<br>$MV[f(x)] = \lim_{\varepsilon \to 0}\left[\frac{f(x+\varepsilon) + f(x-\varepsilon)}{2}\right]$, Mean Value of f(x) at x<br>$A = \Delta x \underset{\Delta x}{\sum_{x=x_1}^{x_2-\Delta x}} f(x)$, where all $f(x_m) = \frac{1}{0}$ terms are excluded<br>$x = x_1, x_1+\Delta x, x_1+2\Delta x, x_1+3\Delta x, \dots, x_2-\Delta x, x_2$<br>$x_m$ = x values at the f(x) poles<br>Interval calculus integration of the area under the function f(x) between the limits $x_1,x_2$<br>Note – for $Re(x_2-x_1) < 0$ the value of $Re(\Delta x) < 0$<br>for $Re(x_2-x_1) > 0$ the value of $Re(\Delta x) > 0$<br>This equation applies if there are f(x) poles within $x_1<x<x_2$<br>There can be no poles at the limits $x_1, x_2$<br>__________________________________________<br>$A = \Delta x \underset{\Delta x}{\sum_{x=x_1}^{x_2-\Delta x}} \frac{1}{(x-a)^n} = \pm lnd(n,\Delta x,x-a)\Big\vert_{x_1}^{x_2}$ + for n=1, – for n≠1<br>$f(x) = \frac{1}{(x-a)^n}$<br>$x = x_1, x_1+\Delta x, x_1+2\Delta x, x_1+3\Delta x, \dots, x_2-\Delta x, x_2$<br>If x=a, the term, $f(a) = \frac{1}{0}$, is excluded<br>a = constant<br>There can be no poles at the limits, $x_1$ and $x_2$ |

| | Notation & Symbols | Description |
|---|---|---|
| **40** | $u(\Delta x,x)\Big\vert_{x_1}^{x_2}\Delta x_p = u(\Delta x_p,x_2) - u(\Delta x_p,x_1)$ | $u(\Delta x,x)$ = function of $\Delta x$ and x<br>$\Delta x_p$ = the value of $\Delta x$<br>$x_1$ = lower limit of x<br>$x_2$ = upper limit of x |
| **41** | $\frac{d}{dx}$ | Continuous derivative symbol |
| **42** | $\frac{d}{dx}f(x)$ | Continuous derivative of f(x)<br>$\frac{d}{dx}f(x) = \lim_{\Delta x\to 0} D_{\Delta x}f(x)$<br>or more precisely<br>$\frac{d}{dx}f(x) = \text{let}_{\Delta x=\varepsilon} D_{\Delta x}f(x)$<br>where<br>$\varepsilon\to 0$ , $\varepsilon\neq 0$ |
| **43** | $\int$ | Continuous integral symbol |
| **44** | $\int f(x)\Delta x + k$ | Continuous integral<br>$\int f(x)\Delta x = \lim_{\Delta x\to 0} {}_{\Delta x}\!\int f(x)\Delta x$ |
| **45** | $e_{\Delta x}x$ | Function of x, $\Delta x$ related to $e^x$<br>Function defined in Table 4 |
| **46** | $\sin_{\Delta x}x$ | Function of x, $\Delta x$ related to sinx<br>Function defined in Table 4 |
| **47** | $\cos_{\Delta x}x$ | Function of x, $\Delta x$ related to cosx<br>Function defined in Table 4 |
| **48** | $\sinh_{\Delta x}x$ | Function of x, $\Delta x$ related to sinhx<br>Function defined in Table 4 |
| **49** | $\cosh_{\Delta x}x$ | Function of x, $\Delta x$ related to coshx<br>Function defined in Table 4 |
| **50** | $e_{\Delta x}(a,x)$ | Function of a,x, $\Delta x$ related to $e^{ax}$<br>Function defined in Table 4 |
| **51** | $\sin_{\Delta x}(a,x)$ | Function of a,x, $\Delta x$ related to sinax<br>Function defined in Table 4 |
| **52** | $\cos_{\Delta x}(a,x)$ | Function of a,x, $\Delta x$ related to cosax<br>Function defined in Table 4 |
| **53** | $\sinh_{\Delta x}(a,x)$ | Function of a,x, $\Delta x$ related to sinhax<br>Function defined in Table 4 |

| | Notation & Symbols | Description |
|---|---|---|
| **54** | $\cosh_{\Delta x}(a,x)$ | Function of a,x, Δx related to coshax<br>Function defined in Table 4 |
| **55** | $F_{\Delta x}(x)$<br>$F(x)$<br>$F(x) = \lim_{\Delta x\to 0} F_{\Delta x}(x)$<br>or more precisely<br>$F(x) = \text{let}_{\Delta x=\varepsilon} F_{\Delta x}(x)$<br>where<br>$\varepsilon\to 0$ , $\varepsilon\neq 0$ | Function of x,Δx (includes the above functions)<br>Related function of x |
| **56** | $K_{\Delta x}[f(x)]$ | $K_{\Delta x}$ Transform of f(x)<br>$K_{\Delta x}[f(x)] = {}_{\Delta x}\int_0^{\infty} (1+s\Delta x)^{-(\frac{x+\Delta x}{\Delta x})} f(x)\Delta x$ |
| **57** | $L[f(x)]$ | Laplace Transform of f(x)<br>$L[f(x)] = \lim_{\Delta x\to 0} K_{\Delta x}[f(x)]$ |
| **58** | $U(x-a)$ | 1 for $x \geq a$ , 0 for $x < a$<br>Unit Step |
| | $\text{lnd}(n,\Delta x,x)$ | The lnd(n,Δx,x) Function<br>The lnd(n,Δx,x) function is defined by:<br>$D_{\Delta x}\text{lnd}(n,\Delta x,x) \equiv D_{\Delta x}\ln_{\Delta x}x = +\frac{1}{x}$ , n =1<br>and<br>$D_{\Delta x}\text{lnd}(n,\Delta x,x) = -\frac{1}{x^n}$ , $n \neq 1$<br><br>The lnd(n,Δx,x) Function is calculated using the program, LNDX. The LNDX program code is presented in the Calculation Programs section at the end of the Appendix. |

| | Notation & Symbols | Description |
|---|---|---|
| **59** | $\zeta(n,\Delta x,x)$ | General Zeta Function<br>$\zeta(n,\Delta x,x) = \frac{1}{\Delta x}\,\mathrm{lnd}(n,\Delta x,x)$<br>$\zeta(n,\Delta x,x_i) = \frac{1}{\Delta x}\,\mathrm{lnd}(n,\Delta x,x_i) = {}_{\Delta x}\sum_{x=x_i}^{\pm\infty}\frac{1}{x^n}$ , Re(n)>1<br>$\zeta(n,\Delta x,x)\Big\vert_{x_1}^{x_2} = \frac{1}{\Delta x}\,\mathrm{lnd}(n,\Delta x,x)\Big\vert_{x_1}^{x_2} = \pm\,{}_{\Delta x}\sum_{x=x_1}^{x_2-\Delta x}\frac{1}{x^n}$ ,<br>+ for n = 1, – for n ≠ 1, Re(n)>1<br><br>if x=0, the summation term, $\frac{1}{0}$ , is excluded |
| **60** | $\zeta(n,x)$ | Hurwitz Zeta Function<br>$\zeta(n,x_i) = \mathrm{lnd}(n,1,x_i)$ , $\Delta x=1$, $n\neq 1$<br>$\zeta(n,x_i) = \mathrm{lnd}(n,1,x_i) = {}_{1}\sum_{x=x_i}^{\infty}\frac{1}{x^n}$ , Re(n)>1<br>if x=0, the summation term, $\frac{1}{0}$ , is excluded<br><br>A special case of the General Zeta Function |
| **61** | $\zeta(n)$ | Riemann Zeta Function<br>$\zeta(n) = \mathrm{lnd}(n,1,1)$ , $\Delta x=1$, $n\neq 1$<br>$\zeta(n) = \mathrm{lnd}(n,1,1) = {}_{1}\sum_{x=1}^{\infty}\frac{1}{x^n}$ , Re(n)>1<br><br>A special case of the General Zeta Function |

| | Notation & Symbols | Description |
|---|---|---|
| **62** | $\zeta(1,\Delta x,x)$ | N=1 Zeta Function<br><br>$\zeta(1,\Delta x,x) = \frac{1}{\Delta x}\,\mathrm{lnd}(1,\Delta x,x) \equiv \frac{1}{\Delta x}\ln_{\Delta x}x$<br><br>$\zeta(1,\Delta x,x_f) = \frac{1}{\Delta x}\,\mathrm{lnd}(1,\Delta x,x_f) \equiv \frac{1}{\Delta x}\ln_{\Delta x}x_f = \sum_{\Delta x \; x=\Delta x}^{x_f-\Delta x}\frac{1}{x}$<br><br>$\zeta(1,\Delta x,x)\Big\|_{x_1}^{x_2} = \frac{1}{\Delta x}\,\mathrm{lnd}(1,\Delta x,x)\Big\|_{x_1}^{x_2} \equiv \frac{1}{\Delta x}\ln_{\Delta x}x\Big\|_{x_1}^{x_2} = \sum_{\Delta x \; x=x_1}^{x_2-\Delta x}\frac{1}{x}$<br><br>$\mathrm{lnd}(1,\Delta x,\Delta x) = 0$<br>$\mathrm{lnd}(1,\Delta x,0) = 0$<br><br>$\ln_{\Delta x}x$ is an optional form used to emphasize the similarity to the natural logarithm.<br><br>A special case of the General Zeta Function |
| **63** | $\psi(x)$ | Digamma Function<br><br>$\psi(x) = \mathrm{lnd}(1,1,x) - \gamma\,,\quad x \neq 0,-1,-2,-3,\ldots$<br><br>$\sum_{\Delta x \; x=x_1}^{x_2}\frac{1}{x} = \frac{1}{\Delta x}\,\psi(\frac{x}{\Delta x})\Big\|_{x_1}^{x_2+\Delta x}$<br><br>$D_1\psi(x) = \frac{1}{x}$ |
| **64** | $\psi(x)^{(m)}\,,\quad m = 1,2,3,\ldots$ | Polygamma Functions<br><br>$\psi^{(m)}(x) = (-1)^{m+1}m!\,\mathrm{lnd}(1+m,1,x)$<br><br>$\sum_{\Delta x \; x=x_1}^{x_2}\frac{1}{x^{1+m}} = \frac{(-1)^m}{m!\Delta x^{m+1}}\,\psi^{(m)}(\frac{x}{\Delta x})\Big\|_{x_1}^{x_2+\Delta x}$<br><br>$D_1\psi^{(m)}(x) = (-1)^m m!\frac{1}{x^{m+1}}$ |

# TABLE 2

## $K_{\Delta x}$ Transform General Equations

| | |
|---|---|
| **1** | The $K_{\Delta x}$ Transform<br>$K_{\Delta x}[f(x)] = {}_{\Delta x}\int_0^{\infty} (1+s\Delta x)^{-(\frac{x+\Delta x}{\Delta x})} f(x)\Delta x$<br>The value of s is defined as: $s = s_{\Delta x} = \frac{e^{s\Delta x} - 1}{\Delta x}$, $-\frac{\pi}{\Delta x} \le w < \frac{\pi}{\Delta x}$, $s = \gamma + jw$<br>where<br>$\Delta x$ = interval between consecutive discrete values of x<br>$x = 0, \Delta x, 2\Delta x, 3\Delta x, \ldots$<br>$\gamma$ is any positive real value which makes s an indefinitely large value. |
| **1a** | $K_{\Delta x}[f(x)] = {}_{\Delta x}\int_0^{\infty} e_{\Delta x}(s,-x-\Delta x)f(x)\Delta x$ |
| **1b** | $K_{\Delta x}[f(x)] = \frac{1}{1+s\Delta x}\int_0^{\infty} e_{\Delta x}(s,-x)f(x)\Delta x$ |
| **1c** | $K_{\Delta x}[f(x)] = \Delta x \sum_{n=0}^{\infty} f(n\Delta x)K_{\Delta x}[\frac{1}{\Delta x}\{U(x-n\Delta x)-U(x-n\Delta x-\Delta x)\}]$<br>Sum of the $K_{\Delta x}$ Transforms of consecutive Unit Amplitude Pulses of $\Delta t$ width |
| **1d** | $K_{\Delta x}[f(x)] = f(s) = \Delta x \sum_{n=0}^{\infty} f(n\Delta x)(1+s\Delta x)^{-n-1} = \Delta x \, {}_{\Delta x}\sum_{x=0}^{\infty} f(x)(1+s\Delta x)^{-(\frac{x+\Delta x}{\Delta x})}$ , $x = n\Delta x$<br><u>Note</u> – This equation can be used to find f(x) from a $K_{\Delta x}$ Transform, f(s) |
| **1e** | $K_{\Delta x}[f(x)] = \Delta x \, {}_{1}\int_0^{\infty} [1+s\Delta x]^{-n-1} f(n\Delta x)\,\Delta n$ , $x = n\Delta x$ |
| **2** | The Inverse $K_{\Delta x}$ Transform<br>$f(x) = K_{\Delta x}^{-1}[F(s)] = \frac{1}{2\pi j}\oint_c F(s)(1+s\Delta x)^{\frac{x}{\Delta x}} ds$<br>or<br>$f(x) = \sum_{p=1}^{P} R_p$ = the sum of the P residues of the P poles of $F(s)(1+s\Delta x)^{\frac{x}{\Delta x}}$ |

$$R = \lim_{s\to r} \frac{1}{(m-1)!} \frac{d^{m-1}}{ds^{m-1}} [(s-r)^m F(s)(1+s\Delta x)^{\frac{x}{\Delta x}}]$$ , Residue calculation formula for a pole at s = r

where

$F(s) = K_{\Delta x}[f(x)]$ , $K_{\Delta x}$ Transform of f(x)

R = the residue of a pole of $F(s)(1+s\Delta x)^{\frac{x}{\Delta x}}$

$\sum_{p=1}^{P} R_p$ = the sum of the P residues of the P poles of $F(s)(1+s\Delta x)^{\frac{x}{\Delta x}}$

$x = 0, \Delta x, 2\Delta x, 3\Delta x, \ldots$

c = The circular closed contour in the complex plane which is shown below.

S Plane Circular Contour

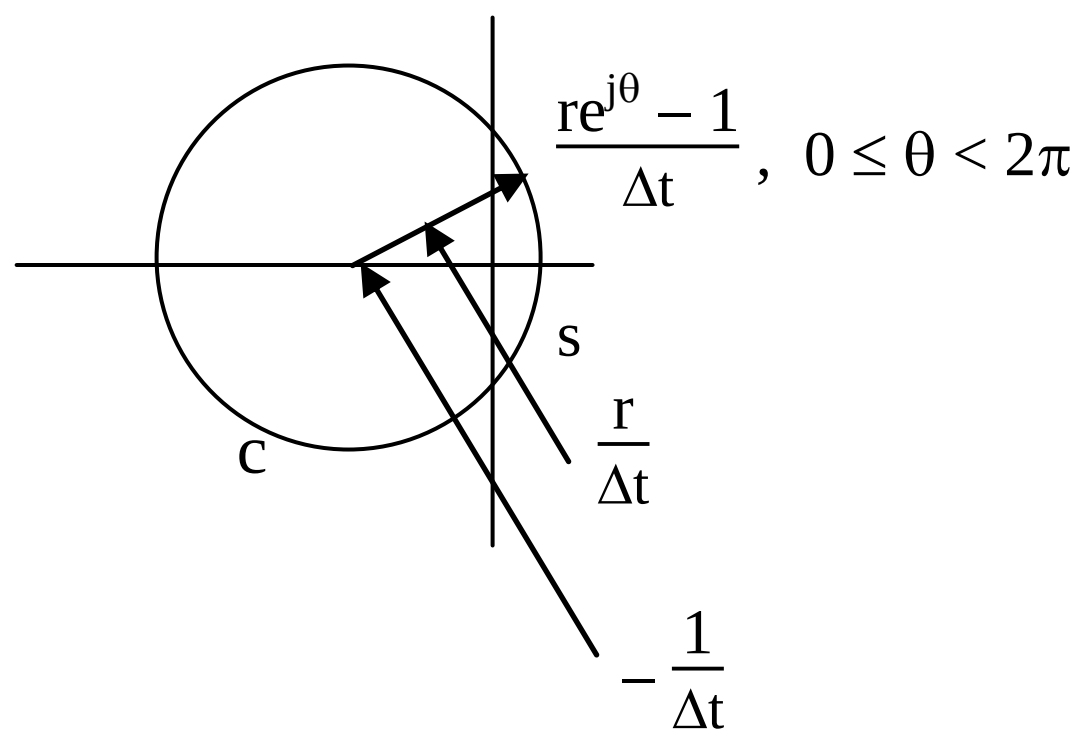

**3**

Relationship between the $K_{\Delta t}$ Transform and the Z Transform

$$K_{\Delta t}[f(t)] = {}_{\Delta t}\int_0^\infty (1+s\Delta t)^{-(\frac{t+\Delta t}{\Delta t})} f(t)\Delta t = {}_{\Delta t}\int_0^\infty e^{-s\Delta t(\frac{t+\Delta t}{\Delta t})} f(t)\Delta t = e^{-sT} {}_{\Delta t}\int_0^\infty e^{-sT\frac{t}{\Delta t}} f(t)\Delta t$$

$$K_{\Delta t}[f(t)] = Te^{-sT} \frac{1}{T} {}_{T}\int_0^\infty e^{-sT\frac{t}{\Delta t}} f(t)\Delta t = Tz^{-1}[\frac{1}{T} {}_{T}\int_0^\infty z^{-\frac{t}{\Delta t}} f(t)\Delta t] = Tz^{-1}Z[f(t)]$$

where

$z = e^{sT}$

$T = \Delta t$ = sampling period

$e^{sT} = 1+s\Delta t$

$Tz^{-1}$ = impulse transfer function of a T sampling period hold

Note - In engineering problems the chosen variable is often t instead of x

| | |
|---|---|
| **4** | Relationship between the $K_{\Delta t}$ Transform and the Laplace Transform<br>$\lim_{\Delta x\to 0} K_{\Delta x}[f(x)] = \lim_{\Delta x\to 0} {}_{\Delta x}\int_{0}^{\infty} (1+s\Delta x)^{-(\frac{x+\Delta x}{\Delta x})} f(x)\Delta x = \lim_{\Delta x\to 0} {}_{\Delta x}\int_{0}^{\infty} (e^{s\Delta x})^{-(\frac{x+\Delta x}{\Delta x})} f(x)\Delta x$<br>$\lim_{\Delta x\to 0} K_{\Delta x}[f(x)] = \lim_{\Delta x\to 0} {}_{\Delta x}\int_{0}^{\infty} (e^{s\Delta x})^{-(\frac{x+\Delta x}{\Delta x})} f(x)\Delta x = \int_{0}^{\infty} e^{-sx} f(x)dx) = L[f(x)]$<br>where<br>$\Delta x$ = interval between consecutive discrete values of x<br>$x = 0, \Delta x, 2\Delta x, 3\Delta x, \ldots$<br>$e^{s\Delta x} = 1+s\Delta x$ |
| **5** | $K_{\Delta x}[cf(x)] = cK_{\Delta x}[f(x)]$ |
| **6** | $K_{\Delta x}[f(x) + g(x)] = K_{\Delta x}[f(x)] + K_{\Delta x}[g(x)]$ |
| **7** | $K_{\Delta x}[D_{\Delta x}f(x)] = sK_{\Delta x}[f(x)] - f(0)$<br><br>The related Z Transform Equation is:<br>$Z[D_T f(x)] = \frac{z-1}{T} Z[f(x)] - \frac{z}{T} f(0)$<br>where<br>$D_{\Delta x}f(x) = \frac{f(x+\Delta x) - f(x)}{\Delta x}$<br>$T = \Delta x$ |
| **8** | $K_{\Delta x}[D^2_{\Delta x}f(x)] = s^2K_{\Delta x}[f(x)] - sf(0) - D_{\Delta x}f(0)$<br><br>The related Z Transform Equation is:<br>$Z[D^2_T f(x)] = (\frac{z-1}{T})^2 Z[f(x)] - \frac{z}{T}(\frac{z-1}{T})f(0) - \frac{z}{T} D_T f(0)$<br>where<br>$D_{\Delta x}f(x) = \frac{f(x+\Delta x) - f(x)}{\Delta x}$<br>$T = \Delta x$ |
| **9** | $K_{\Delta x}[D^n{}_{\Delta x}f(x)] = s^n K_{\Delta x}[f(x)] - s^{n-1}D^0{}_{\Delta x}f(0) - s^{n-2}D^1{}_{\Delta x}f(0) - s^{n-3}D^2{}_{\Delta x}f(0) - \ldots - s^0 D^{n-1}{}_{\Delta x}f(0)$ , n = 1,2,3,…<br>or<br>$K_{\Delta x}[D^n{}_{\Delta x}f(x)] = s^n K_{\Delta x}[f(x)] - \sum_{m=1}^{n} s^{n-m}D_{\Delta x}{}^{m-1}f(0)$ , n = 1,2,3,…<br>The related Z Transform Equation is:<br>$Z[D_T{}^n] = (\frac{z-1}{T})^n \, Z[f(x)] - \sum_{m=1}^{n} \frac{z}{T}(\frac{z-1}{T})^{n-m}D_T{}^{m-1}f(0)$ , n = 1,2,3,… , $T = \Delta x$ |
| **10** | $K_{\Delta x}[\Delta f(x)] = \Delta x s K_{\Delta x}[f(x)] - \Delta x f(0)$<br>The related Z Transform Equation is:<br>$Z[\Delta f(x)] = (z-1)Z[f(x)] - zf(0)$<br>where<br>$\Delta f(x) = f(x+\Delta x) - f(x)$<br>$\Delta = \Delta x D_{\Delta x}$ |
| **11** | $K_{\Delta x}[f(x+\Delta x] = (1+s\Delta x)K_{\Delta x}[f(x)] - \Delta x f(0)$ |
| **12** | $K_{\Delta x}[f(x+2\Delta x] = (1+s\Delta x)^2 K_{\Delta x}[f(x)] - (1+s\Delta x)\Delta x f(0) - \Delta x f(\Delta x)$ |
| **13** | $K_{\Delta x}[f(x+m\Delta x] = (1+s\Delta x)^m K_{\Delta x}[f(x)] - \Delta x \sum_{r=0}^{m-1}(1+s\Delta x)^{m-1-r}f(r\Delta x)$ , m =1,2,3,… |
| **14** | $\Delta f(x) = \Delta x D_{\Delta x} f(x)$<br>$\Delta f(x) = f(x+\Delta x) - f(x)$ |
| **15** | $\Delta^2 f(x) = \Delta x^2 D_{\Delta x}{}^2 f(x)$ |
| **16** | $\Delta^n f(x) = \Delta x^n D_{\Delta x}{}^n f(x)$ , n = 1,2,3, … |
| **17** | $K_{\Delta x}[b(x)] = s\varphi(s)$<br>$b(x) = D_{\Delta x}K_{\Delta x}{}^{-1}[\varphi(s)]$<br>where<br>$x = 0, \Delta x, 2\Delta x, 3\Delta x, \ldots, \infty$<br>b(x) is finite for all x |
| | $\sum_{0}^{\infty}{}_{\Delta x}\int (1+s\Delta x)^{-(\frac{x+\Delta x}{\Delta x})} b(x)\,\Delta x$ is convergent<br>$\lim_{s\to\infty} K_{\Delta x}[b(x)] = 0$ |
| **18** | $K_{\Delta x}[f(x+\Delta x)] = (1+s\Delta x)K_{\Delta x}[f(x)] - \Delta x f(0)$ |
| **19** | $K_{\Delta x}[f(x+2\Delta x)] = (1+s\Delta x)^2 K_{\Delta x}[f(x)] - \Delta x(1+s\Delta x)f(0) - \Delta x f(\Delta x)$ |
| **20** | $K_{\Delta x}[f(x+n\Delta x)] = (1+s\Delta x)^n K_{\Delta x}[f(x)] - \Delta x \sum_{m=1}^{n}(1+s\Delta x)^{n-m} f([m-1]\Delta x)$ , $n = 1,2,3$<br><br>The related Z Transform Equation is:<br>$Z[f(x+nT)] = z^n Z[f(x)] - \sum_{m=0}^{n-1} z^{n-m} f(mT)$ , $n = 1,2,3\ldots$ , $T = \Delta x$ |
| **21** | $K_{\Delta x}[{}_{\Delta x}\int_0^x f(x)\Delta x] = \frac{1}{s} K_{\Delta x}[f(x)]$<br><br>The related Z Transform Equation is:<br>$Z[\frac{1}{T}{}_T\int_0^x f(x)\Delta x] = \frac{1}{z-1} Z[f(x)]$ , $\Delta x = T$ |
| **22** | $K_{\Delta x}[f(x-n\Delta x)U(x-n\Delta x)] = (1+s\Delta x)^{-n}F(s)$<br>$U(x-n\Delta x) = \begin{Bmatrix} 1 & x \geq n\Delta x \\ 0 & x < n\Delta x \end{Bmatrix}$ , Unit Step Function<br><br>$n = 0, 1, 2, 3, \ldots$<br>$F(s) = K_{\Delta x}$ Transform of $f(x)$<br><br><u>Comment</u> – To find $c(x)$ where $c(x) = K_{\Delta x}^{-1}[(1+s\Delta x)^{-n}F(s)]$<br>First find $f(x)$ where $f(x) = K_{\Delta x}^{-1}[F(s)]$<br>Then $c(x) = f(x-n\Delta x)U(x-n\Delta x)$ |
| **23** | $K_{\Delta x}[(1+a\Delta x)^{\frac{x}{\Delta x}} f(x)] = \frac{1}{1+a\Delta x} K_{\Delta x}[f(x)]\ \|_{s\to\frac{s-a}{1+a\Delta x}}$ |
| **24** | $K_{\Delta x}[b^{-\frac{x}{\Delta x}}f(x)] = b\,K_{\Delta x}[f(x)]\big\|_{s\rightarrow\, bs + \frac{b-1}{\Delta x}}$ |
| **25** | Function with an alternating sign<br>$K_{\Delta x}[(-1)^{-\frac{x}{\Delta x}}f(x)] = -K_{\Delta x}[f(x)]\big\|_{s\rightarrow\, -s - \frac{2}{\Delta x}}$<br>The related Z Transform Equation is:<br>$Z[(-1)^{-\frac{x}{\Delta x}}f(x)] = Z[f(x)]\big\|_{z\rightarrow\, -z}$ |
| **26** | $K_{\Delta t}[\frac{1}{c}\, f(\frac{t}{c})] = K_{\Delta t}[f(t)]\big\|_{s\rightarrow cs}$ |
| **27** | $K_{\Delta x}[e^{\frac{ax}{\Delta x}}f(x)] = \frac{1}{e^{a}}\, K_{\Delta x}[f(x)]\,\big\|_{s\rightarrow \frac{s + \frac{1-e^{a}}{\Delta x}}{e^{a}}}$ |
| **28** | $K_{\Delta x}[xf(x)] = -(1+s\Delta x)\frac{d}{ds}\, K_{\Delta x}[f(x)] - \Delta x K_{\Delta x}[f(x)]$ |
| **29** | $K_{\Delta x}[x^2 f(x)] = (1+s\Delta x)^2\frac{d^2}{ds^2}\, K_{\Delta x}[f(x)] + 3\Delta x(1+s\Delta x)\frac{d}{ds}\, K_{\Delta x}[f(x)] + \Delta x^2\, K_{\Delta x}[f(x)]$ |
| **30** | Periodic function, f(x)<br>$K_{\Delta x}[f(x)] = \dfrac{\Delta x \int_{0}^{N\Delta x} (1+s\Delta x)^{-(\frac{x+\Delta x}{\Delta x})} f(x)\Delta x}{1 - (1+s\Delta x)^{-N}}$ , $f(x) = f(x+N\Delta x)$ , $N\Delta x$ = period<br>N = positive integer |
| **31** | $K_{\Delta x}[\{\prod_{m=1}^{n}(x+m\Delta x)\}f(x)] = (-1)^n(1+s\Delta x)^n\frac{d^n}{ds^n}K_{\Delta x}[f(x)]$<br>$n = 1,2,3,\ldots$ |
| **32** | Multiplication of x and Δx by a constant (b)<br>For $H(x,\Delta x)$ = an Interval Calculus function of x and Δx<br>$h(s) = K_{\Delta x}[H(x,\Delta x)]$<br>$H(bx,b\Delta x) = K_{\Delta x}^{-1}[\frac{1}{b}\, h(s)\big\|_{s\rightarrow \frac{s}{b},\ \Delta x\rightarrow b\Delta x}]$<br>$x = 0, \Delta x, 2\Delta x, 3\Delta x, \ldots$<br>$H(x,\Delta x)$ functions: $e_{\Delta x}(a,x)$, $\sin_{\Delta x}(a,x)$, $x(x-\Delta x)$<br>$H(bx,b\Delta x)$ functions: $e_{b\Delta x}(a,bx) = e_{\Delta x}(ba,x)$, $\sin_{b\Delta x}(a,bx) = \sin_{\Delta x}(ba,x)$, $bx(bx-b\Delta x) = b^2x(x-\Delta x)$ |
| **33** | Initial-value theorem<br>$\lim_{x\to 0} f(x) = \lim_{s\to\infty} sf(s)$ , if the indicated limits exist<br>$f(s) = K_{\Delta x}[f(x)]$ |
| **34** | Final-value theorem<br>$\lim_{x\to\infty} f(x) = \lim_{s\to 0} sf(s)$ , if the indicated limits exist<br>$f(s) = K_{\Delta x}[f(x)]$ |
| **35** | $\lim_{s\to\infty} f(s) = 0$<br>$f(s) = K_{\Delta x}[f(x)]$ |
| **36** | $f(n\Delta x) = \dfrac{1}{\Delta x(n+1)!} \lim_{p\to 0} \dfrac{d^{n+1}}{dp^{n+1}} F(p)$ , $x=n\Delta x$ , $n = 0,1,2,3,\ldots$<br>$f(s) = K_{\Delta x}[f(x)]$<br>$F(p) = f(s)\Big\vert_{s = \frac{1-p}{p\Delta x}}$ |
| **37** | $K_{\Delta x}$ Transform to Z Transform Conversion<br>$Z[f(x)] = \dfrac{1+s\Delta x}{\Delta x} K_{\Delta x}[f(x)]\Big\vert_{s = \frac{z-1}{T}}$ $Z[f(x)] = F(z)$ Z Transform<br>$T = \Delta x$<br>$x = nT$ , $n = 0,1,2,3,\ldots$ $K_{\Delta x}[f(x)] = f(s)$ $K_{\Delta x}$ Transform |
| **38** | Z Transform to $K_{\Delta x}$ Transform Conversion<br>$K_{\Delta x}[f(x)] = Tz^{-1} Z[f(x)]\Big\vert_{z = 1+s\Delta x}$ $K_{\Delta x}[f(x)] = f(s)$ $K_{\Delta x}$ Transform<br>$\Delta x = T$<br>$x = n\Delta x$ , $n = 0,1,2,3,\ldots$ $Z[f(x)] = F(z)$ Z Transform |
| **39** | Z Transform to $J_{\Delta x}$ Transform Conversion<br>$J_{\Delta x}[f(x)] = \Delta x\, Z[f(x)]\Big\vert_{z = 1+s\Delta x}$ $J_{\Delta x}[f(x)] = f(s)$ $J_{\Delta x}$ Transform<br>$\Delta x = T$<br>$x = n\Delta x$ , $n = 0,1,2,3,\ldots$ $Z[f(x)] = F(z)$ Z Transform |
| **40** | $J_{\Delta x}$ Transform to Z Transform Conversion<br>$Z[f(x)] = \dfrac{1}{T} J_{\Delta x}[f(x)]\Big\vert_{s = \frac{z-1}{T}}$ , $Z[f(x)] = F(z)$ Z Transform<br>$T = \Delta x$<br>$x = nT$ , $n = 0,1,2,3,\ldots$ $J_{\Delta x}[f(x)] = g(s)$ $J_{\Delta x}$ Transform |
| **41** | $K_{\Delta x}$ Transform to $J_{\Delta x}$ Transform Conversion<br>$J_{\Delta x}[f(x)] = (1+s\Delta x)\, K_{\Delta x}[f(x)]$ $J_{\Delta x}[f(x)] = g(s)$ $J_{\Delta x}$ Transform<br>$x = n\Delta x$ , $n = 0,1,2,3,\ldots$ $K_{\Delta x}[f(x)] = f(s)$ $K_{\Delta x}$ Transform |
| **42** | $J_{\Delta x}$ Transform to $K_{\Delta x}$ Transform Conversion<br>$K_{\Delta x}[f(x)] = (1+s\Delta x)^{-1} J_{\Delta x}[f(x)]$ $K_{\Delta x}[f(x)] = f(s)$ $K_{\Delta x}$ Transform<br>$x = n\Delta x$ , $n = 0,1,2,3,\ldots$ $J_{\Delta x}[f(x)] = g(s)$ $J_{\Delta x}$ Transform |

**43** $K_{\Delta x}$ Transform

$$K_{\Delta x}[f(x)] = {}_{\Delta x}\sum_{x=0}^{\infty} (1+s\Delta x)^{-\left(\frac{x+\Delta x}{\Delta x}\right)} f(x)\Delta x = {}_{\Delta x}\int_{0}^{\infty} (1+s\Delta x)^{-\left(\frac{x+\Delta x}{\Delta x}\right)} f(x)\Delta x$$

**44** $J_{\Delta x}$ Transform

$$J_{\Delta x}[f(x)] = {}_{\Delta x}\sum_{x=0}^{\infty} (1+s\Delta x)^{-\left(\frac{x}{\Delta x}\right)} f(x)\Delta x = {}_{\Delta x}\int_{0}^{\infty} (1+s\Delta x)^{-\left(\frac{x}{\Delta x}\right)} f(x)\Delta x$$

**45** Z Transform

$$Z[f(x)] = \frac{1}{T}\,{}_{T}\sum_{x=0}^{\infty} z^{-\left(\frac{x}{\Delta x}\right)} f(x)\Delta x = \frac{1}{T}\,{}_{T}\int_{0}^{\infty} z^{-\left(\frac{x}{\Delta x}\right)} f(x)\Delta x \qquad Z[f(x)] = F(z) \quad \text{Z Transform}$$

$x = 0, \Delta x, 2\Delta x, 3\Delta x, \dots, \infty$

$\Delta x = T$

**46** Inverse Z Transform

$Z^{-1}[F(z)] = f(x) =$ sum of the residues of $z^{\frac{x}{T}-1}F(z)$ $\qquad Z^{-1}[F(z)] = f(x)$ Inverse Z Transform

$x = 0, T, 2T, 3T, \dots$

**47** Z Transform to $K_{\Delta x}$ Transform Transfer Function Diagram Conversion

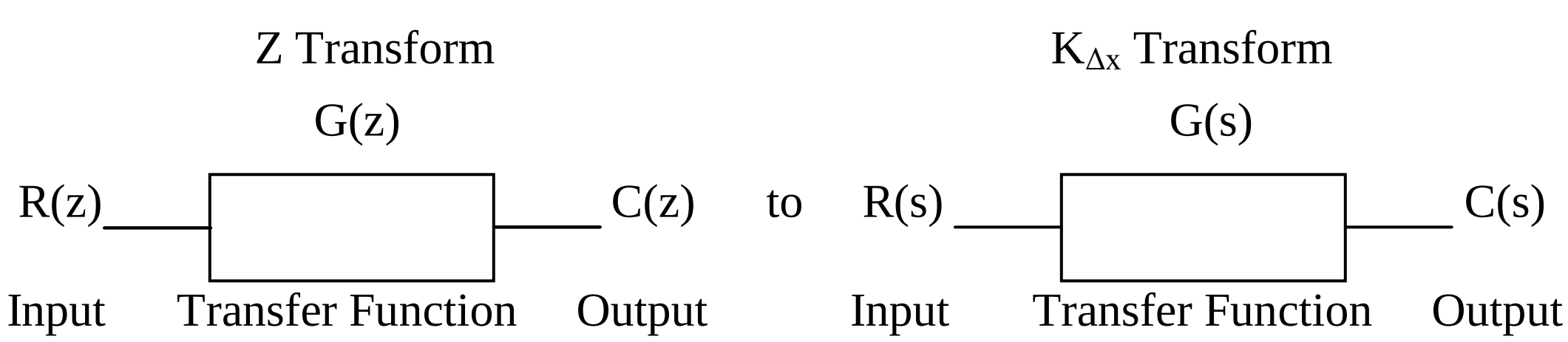

$R(s) = Tz^{-1}R(z)\Big|_{\substack{z = 1+s\Delta x \\ \Delta x = T \\ x = n\Delta x,\ n = 0,1,2,3,\dots}}$

$C(s) = Tz^{-1}C(z)\Big|_{\substack{z = 1+s\Delta x \\ \Delta x = T \\ x = n\Delta x,\ n = 0,1,2,3,\dots}}$

$G(s) = G(z)\Big|_{\substack{z = 1+s\Delta x \\ \Delta x = T \\ x = n\Delta x,\ n = 0,1,2,3,\dots}}$

$C(z) = R(z)G(z)$ to $C(s) = R(s)G(s)$ Equivalent Z and $K_{\Delta x}$ Transform functions

where

$C(z)$ = Z Transform of the output function, $c(x)$

R(z) = Z Transform of the input function, r(x)
G(z) = Z Transform transfer function
C(s) = $K_{\Delta x}$ Transform of the output function, c(x)
R(s) = $K_{\Delta x}$ Transform of the input function, r(x)
G(s) = $K_{\Delta x}$ Transform transfer function

This conversion allows one to use $K_{\Delta t}$ Transforms instead of Z Transforms to obtain the same calculation result.

**48** $K_{\Delta x}$ Transform to Z Transform Transfer Function Diagram Conversion

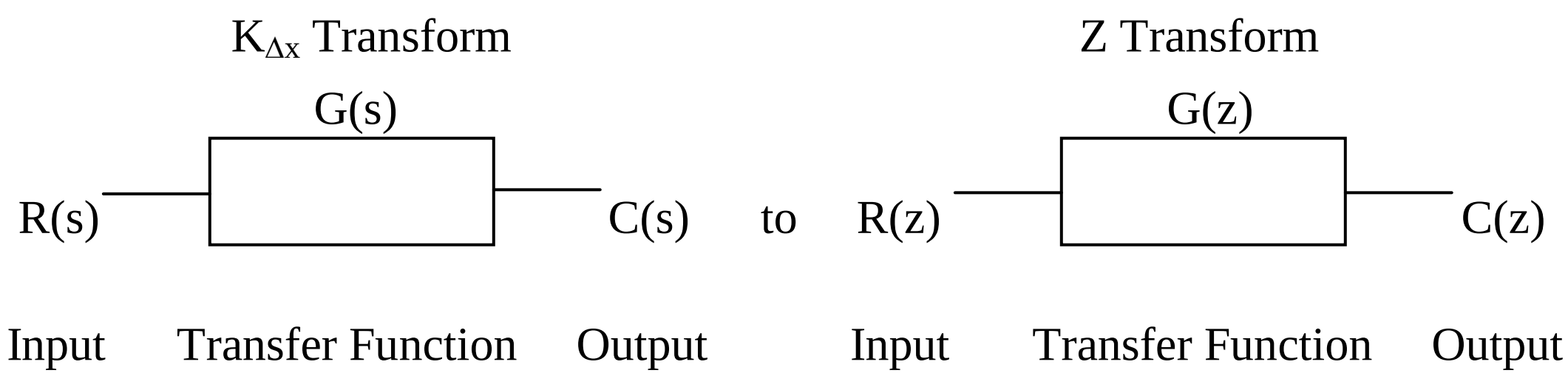

$$R(z) = \frac{1+s\Delta x}{\Delta x} R(s)\Big|_{s = \frac{z-1}{T}}$$
$T = \Delta x$
$x = nT$ ,  $n = 0,1,2,3,\ldots$

$$C(z) = \frac{1+s\Delta x}{\Delta x} C(s)\Big|_{s = \frac{z-1}{T}}$$
$T = \Delta x$
$x = nT$ ,  $n = 0,1,2,3,\ldots$

$$G(z) = G(s)\Big|_{s = \frac{z-1}{T}}$$
$T = \Delta x$
$x = nT$ ,  $n = 0,1,2,3,\ldots$

C(s) = R(s)G(s)   to   C(z) = R(z)G(z)     Equivalent Z and $K_{\Delta x}$ Transform functions

where
C(z) = Z Transform of the output function, c(x)
R(z) = Z Transform of the input function, r(x)
G(z) = Z Transform transfer function
C(s) = $K_{\Delta x}$ Transform of the output function, c(x)
R(s) = $K_{\Delta x}$ Transform of the input function, r(x)
G(s) = $K_{\Delta x}$ Transform transfer function

This conversion allows one to use Z Transforms instead of $K_{\Delta t}$ Transforms to obtain the same calculation result.

| | |
|---|---|
| **49** | $K_{\Delta x}$ Transform Convolution Equation<br>$K_{\Delta x}[\,{}_{\Delta\lambda}\int_{0}^{x} f(x-\lambda-\Delta\lambda)g(\lambda)\Delta\lambda] = K_{\Delta x}[f(x)]K_{\Delta x}[g(x)]$<br>or<br>$K_{\Delta x}[\,{}_{\Delta\lambda}\int_{0}^{x} f(\lambda)g(x-\lambda-\Delta\lambda)\Delta\lambda] = K_{\Delta x}[f(x)]K_{\Delta x}[g(x)]$<br>where<br>$\Delta\lambda = \Delta x$<br>$\lambda = 0, \Delta\lambda, 2\Delta\lambda, 3\Delta\lambda, \ldots, x-\Delta\lambda, x$<br>$f(x),g(x)$ = functions of x<br>$x = 0, \Delta x, 2\Delta x, 3\Delta x, \ldots$<br>$K_{\Delta x}[f(x)] = f(s) = K_{\Delta x}$ Transform of f(x)<br>$K_{\Delta x}[f(x)] = {}_{\Delta x}\int_{0}^{\infty} (1+s\Delta x)^{-(\frac{x+\Delta x}{\Delta x})} f(x)\Delta x$ |
| **50** | Z Transform Convolution Equations<br>$Z[\frac{1}{T}\,{}_{T}\int_{0}^{x+T} f(x-\lambda)g(\lambda)\Delta\lambda] = Z[f(x)]Z[g(x)]$ , $\lambda = 0, \Delta\lambda, 2\Delta\lambda, 3\Delta\lambda, \ldots, x-\Delta\lambda, x, x+T$<br>or<br>$Z[\frac{1}{T}\,{}_{T}\int_{0}^{x+T} f(\lambda)g(x-\lambda)\Delta\lambda] = Z[f(x)]Z[g(x)]$ , $\lambda = 0, \Delta\lambda, 2\Delta\lambda, 3\Delta\lambda, \ldots, x-\Delta\lambda, x, x+T$<br>or<br>$Z[\frac{1}{T}\,{}_{T}\int_{0}^{x} f(x-\lambda-\Delta\lambda)g(\lambda)\Delta\lambda] = \frac{Z[f(x)]Z[g(x)]}{z}$ , $\lambda = 0, \Delta\lambda, 2\Delta\lambda, 3\Delta\lambda, \ldots, x-\Delta\lambda, x$<br>or<br>$Z[\frac{1}{T}\,{}_{T}\int_{0}^{x} f(\lambda)g(x-\lambda-\Delta\lambda)\Delta\lambda] = \frac{Z[f(x)]Z[g(x)]}{z}$ , $\lambda = 0, \Delta\lambda, 2\Delta\lambda, 3\Delta\lambda, \ldots, x-\Delta\lambda, x$<br>where<br>$\Delta\lambda = \Delta x = T$<br>$f(x),g(x)$ = functions of x<br>$x = 0, \Delta x, 2\Delta x, 3\Delta x, \ldots$<br>$Z[f(x)] = F(z)$ = Z Transform of f(x) |

$$Z[f(x)] = \frac{1}{T}\ {}_{T}\!\int_{0}^{\infty} z^{-(\frac{x}{\Delta x})} f(x)\Delta x$$

**51** <u>Interval Calculus Duhamel Equations</u>

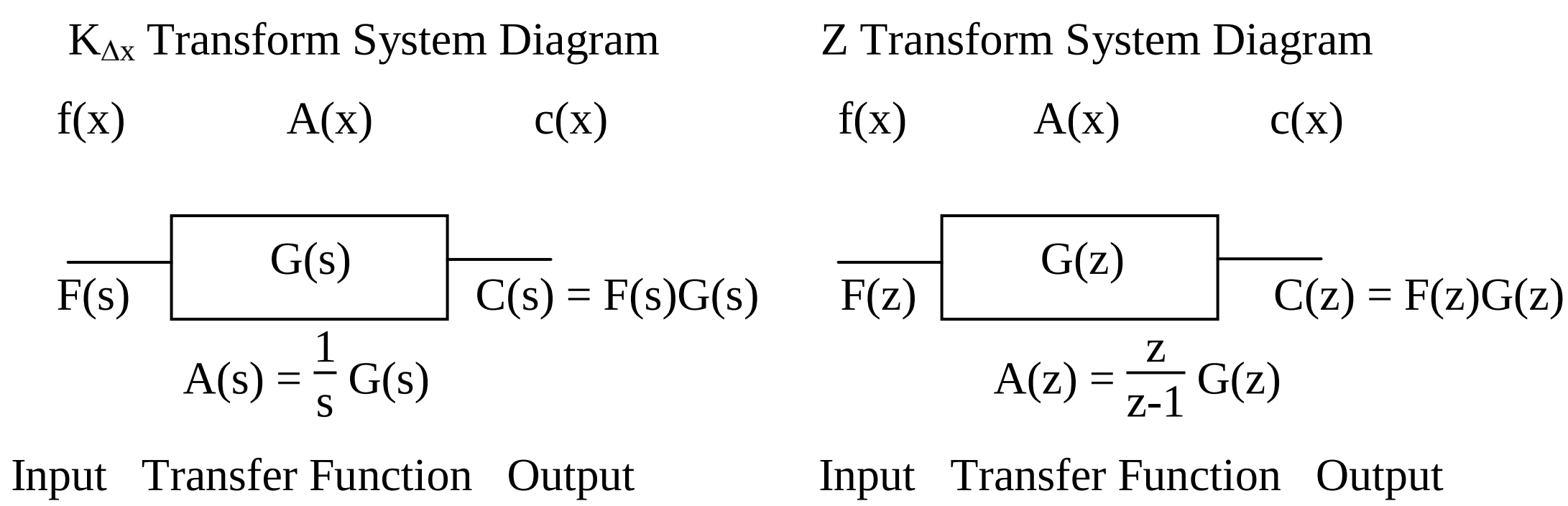

The following $K_{\Delta x}$ Transform and Z Transform Duhamel Equations make it possible to obtain the system response to a general input function, f(x), if the unit step function system response, A(x), is known.

1. $c(x) = A(0)f(x) + {}_{\Delta\lambda}\!\int_{0}^{x} f(\lambda)A'(x-\lambda-\Delta\lambda)\Delta\lambda$ , $\quad C(s) = F(s)[\frac{sG(s)}{s}] = F(s)[sA(s)]$

   $C(z) = F(z)[\frac{z-1}{z}][\frac{z}{z-1}G(z)] = F(z)[\ \frac{z-1}{z}\ A(z)]$

   Since A(x) is by definition the response of a system that is initially passive, A(0) = 0.

2. $c(x) = A(0)f(x) + {}_{\Delta\lambda}\!\int_{0}^{x} f(x-\lambda-\Delta\lambda)A'(\lambda)\Delta\lambda$ , $\quad C(s) = F(s)\ [\frac{sG(s)}{s}] = F(s)[sA(s)]$

   $C(z) = F(z)[\frac{z-1}{z}][\frac{z}{z-1}G(z)] = F(z)[\ \frac{z-1}{z}\ A(z)]$

   Since A(x) is by definition the response of a system that is initially passive, A(0) = 0.

3. $c(x) = {}_{\Delta\lambda}\!\int_{0}^{x} f(\lambda)A'(x-\lambda-\Delta\lambda)\Delta\lambda$ , $\quad C(s) = F(s)\ [\frac{sG(s)}{s}] = F(s)[sA(s)]$

   $C(z) = F(z)[\frac{z-1}{z}][\frac{z}{z-1}G(z)] = F(z)[\ \frac{z-1}{z}\ A(z)]$

4. $c(x) = {}_{\Delta\lambda}\!\int_{0}^{x} f(x-\lambda-\Delta\lambda)A'(\lambda)\Delta\lambda$ , $\quad C(s) = F(s)\ [\frac{sG(s)}{s}] = F(s)[sA(s)]$

$$C(z) = F(z)[\frac{z-1}{z}][\frac{z}{z-1}G(z)] = F(z)[\frac{z-1}{z}A(z)]$$

5. $c(x) = f(0)A(x) + {}_{\Delta\lambda}\int_0^x f'(\lambda)A(x-\lambda-\Delta\lambda)\Delta\lambda$ , $C(s) = [sF(s)]\frac{G(s)}{s} = [sF(s)]A(s)$

$$C(z) = [\frac{z-1}{z}F(z)][\frac{z}{z-1}G(z)] = [\frac{z-1}{z}F(z)]A(z)$$

6. $c(t) = f(0)A(x) + {}_{\Delta\lambda}\int_0^x f'(x-\lambda-\Delta\lambda)A(\lambda)\Delta\lambda$ , $C(s) = [sF(s)]\frac{G(s)}{s} = [sF(s)]A(s)$

$$C(z) = [\frac{z-1}{z}F(z)][\frac{z}{z-1}G(z)] = [\frac{z-1}{z}F(z)]A(z)$$

7. $c(x) = D_{\Delta x}\,{}_{\Delta\lambda}\int_0^x f(\lambda)A(x-\lambda-\Delta\lambda)\Delta\lambda$ $C(s) = s[F(s)\frac{G(s)}{s}] = s[F(s)A(s)]$

$$C(z) = \frac{z-1}{z}[F(z)\frac{z}{z-1}G(z)] = \frac{z-1}{z}[F(z)A(z)]$$

8. $c(t) = D_{\Delta x}\,{}_{\Delta\lambda}\int_0^x f(x-\lambda-\Delta\lambda)A(\lambda)\Delta\lambda$ $C(s) = s[F(s)\frac{G(s)}{s}] = s[F(s)A(s)]$

$$C(z) = \frac{z-1}{z}[F(z)\frac{z}{z-1}G(z)] = \frac{z-1}{z}[F(z)A(z)]$$

where

$F(s) = K_{\Delta x}$ Transform of the system input excitation function, f(x)
$f(x)$ = System input excitation function
$G(s)$ = System $K_{\Delta x}$ Transform transfer function
$A(s) = \frac{1}{s}G(s) = K_{\Delta x}$ Transform of the system response to a unit step function
$A(x)$ = System response to a unit step function
$A(0) = 0$, The system, by definition, is initally passive
$A'(x) = D_{\Delta x}A(x)$ = Discrete derivative of A(x)
$C(s) = K_{\Delta x}$ Transform of the system output function, c(x)
$c(x)$ = System output function
$\Delta\lambda = \Delta x$
$\lambda = n\Delta x$, $n = 0, 1, 2, 3, \ldots, \frac{x}{\Delta x}-1, x$
$\lambda = 0, \Delta\lambda, 2\Delta\lambda, 3\Delta\lambda, \ldots, t-\Delta\lambda, t$
$f'(\lambda) = D_{\Delta\lambda}f(\lambda) = \frac{f(\lambda+\Delta\lambda)-f(\lambda)}{\Delta\lambda}$ = discrete derivative of $f(\lambda)$
$\Delta f(\lambda) = f(\lambda+\Delta\lambda)-f(\lambda) = f'(\lambda)\Delta\lambda$

Comments - By definition, the discrete integral, ${}_{\Delta\lambda}\int_{0}^{x} h(\lambda)\Delta\lambda$ , is the sum, ${}_{\Delta\lambda}\sum_{\lambda=0}^{x-\Delta\lambda} h(\lambda)\Delta\lambda$ .

Though the integration is from $\lambda = 0$ thru x, the summation that defines this integration is from $\lambda = 0$ thru $x-\Delta\lambda$.

The system response to a $\Delta x$ width unit amplitude pulse input excitation is $A'(x-\Delta x)\Delta x$.

**52** $K_{\Delta x}$ Transform Transfer Function Response to a Unit Amplitude Pulse Input

g(x)

$r(x) = U(x-n\Delta x) - U(x-[n+1]\Delta x)$ — [ $K_{\Delta x}$ Transform Transfer Function ] — $c(x) = g(x-[n+1]\Delta x)U(x-[n+1]\Delta x)\Delta x$

$R(s) = (1+s\Delta x)^{-(n+1)}\Delta x$ G(s) $C(s) = (1+s\Delta x)^{-(n+1)}\Delta x\, G(s)$

Input Transfer Function Output

$c(x) = g(x-[n+1]\Delta x)U(x-[n+1]\Delta x)\Delta x$

where

$r(x) = U(x-n\Delta x) - U(x-[n+1]\Delta x)$ = Unit Amplitude Pulse Input
$R(s) = (1+s\Delta x)^{-(n+1)}\Delta x$ = Unit Amplitude Pulse Input $K_{\Delta x}$ Transform
$G(s)$ = $K_{\Delta x}$ Transform Transfer Function
$g(x) = K_{\Delta x}^{-1}[G(s)]$ = Inverse $K_{\Delta x}$Transform of G(s)
c(x) = Output response to the Input Unit Amplitude Pulse
$C(s) = K_{\Delta x}[c(t)]$ = Output $K_{\Delta x}$ Transform, $(1+s\Delta x)^{-(n+1)}\Delta x\, G(s)$
$\Delta x$ = Interval between successive discrete independent variable values
$x = 0, \Delta x, 2\Delta x, 3\Delta x, \ldots$ = discrete independent variable values
$n = 0, 1, 2, 3, \ldots$

**53** Z Transform Transfer Function Response to a Unit Area Impulse Input

g(x)

$r(x) = \delta(x-n\Delta x)$ — [ Z Transform Transfer Function ] — $c(x) = g(x-n\Delta x)U(x-n\Delta x)$

$R(z) = z^{-n}$ G(z) $C(z) = z^{-n}G(z)$

Input Transfer Function Output

$c(x) = g(x-n\Delta x)U(x-n\Delta x)$

where

$r(x) = \delta(x-n\Delta x)$ = Unit Impulse Input
$R(z) = z^{-n}$ = Unit Impulse Input Z Transform
G(z) = Z Transform Transfer Function

$g(x) = Z^{-1}[G(s)]$ = Inverse Z Transform of G(z)
c(x) = output response to the Input Unit Impulse
$C(z) = Z[c(t)]$ = Output Z Transform, $z^{-n}G(z)$
$\Delta x = T$ = Interval between successive discrete independent variable values
x = 0, T, 2T, 3T, … = discrete independent variable values
n = 0, 1, 2, 3, …

**54** $K_{\Delta x}$ Transform Transfer Function Response to a Unit Step Input

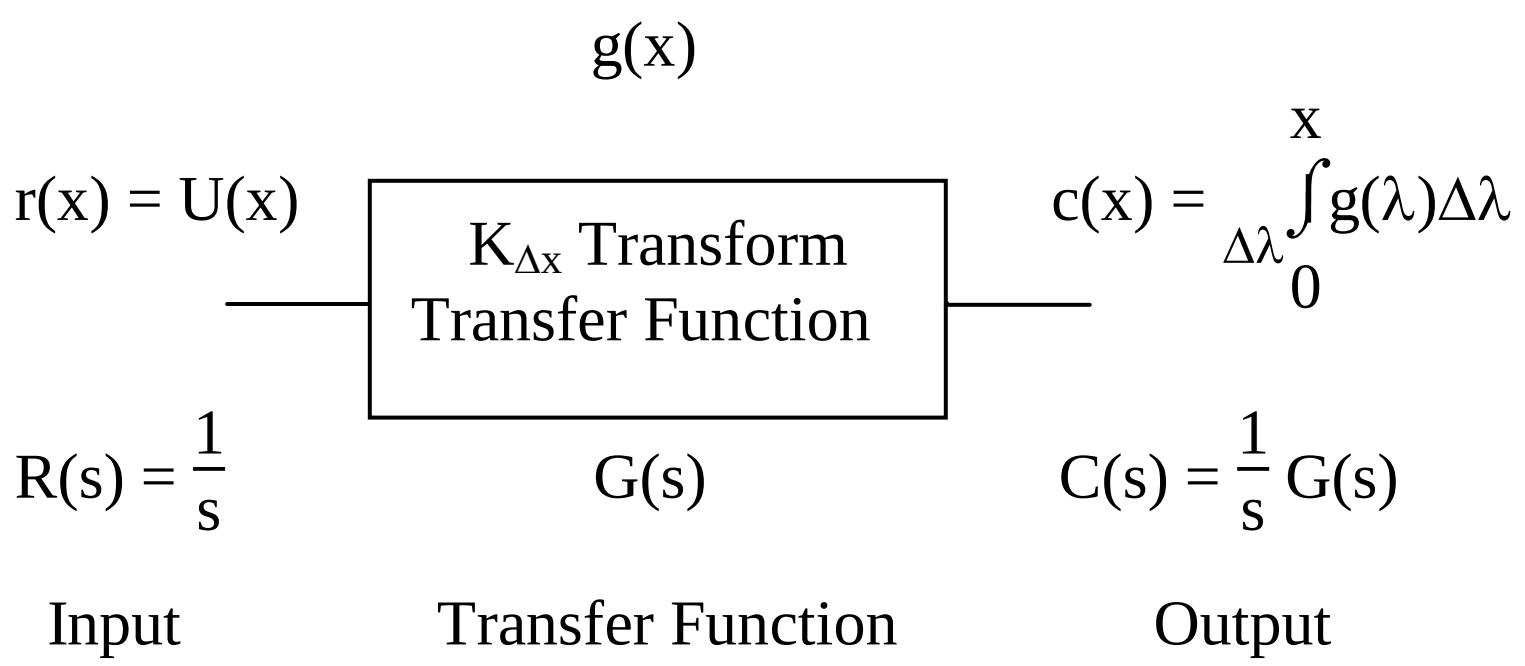

$$c(x) = {}_{\Delta\lambda}\int_0^x g(\lambda)\Delta\lambda \qquad 1)$$

where
r(x) = U(x) = Unit Step Input
$R(s) = \frac{1}{s}$ = Unit Step Input $K_{\Delta x}$ Transform
G(s) = $K_{\Delta x}$ Transform Transfer Function
$g(x) = K_{\Delta x}^{-1}[G(s)]$ = Inverse $K_{\Delta x}$Transform of G(s)
c(x) = output response to the input Unit Step
$C(s) = K_{\Delta x}[c(t)]$ = Output $K_{\Delta x}$ Transform, $\frac{1}{s}G(s)$
$\Delta x = \Delta\lambda$ = Interval between successive discrete independent variable values
x = 0, Δx, 2Δx, 3Δx, … = discrete independent variable values

**55** Z Transform Transfer Function Response to a Unit Step Input

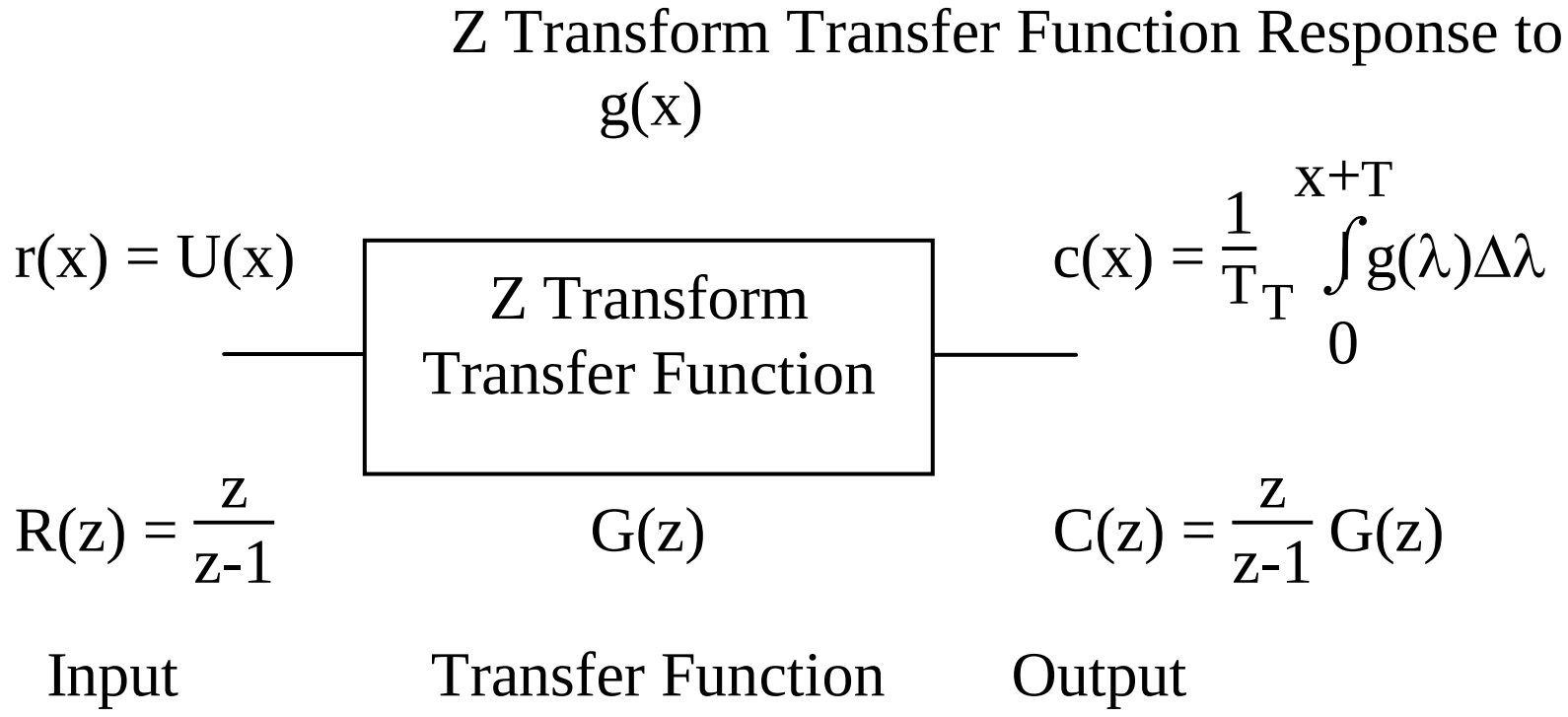

$$c(x) = \frac{1}{T}\,{}_{T}\int_0^{x+T} g(\lambda)\Delta\lambda \qquad 1)$$

where
r(x) = U(x) = Unit Step Input

$R(z) = \frac{z}{z-1}$ = Unit Step Input Z Transform

G(z) = Z Transform Transfer Function

$g(x) = Z^{-1}[G(s)]$ = Inverse Z Transform of G(z)

c(x) = output response to the input Unit Step

$C(z) = Z[c(t)]$ = Output Z Transform, $\frac{z}{z-1} G(z)$

$\Delta x = \Delta\lambda = T$ = Interval between successive discrete independent variable values

x = 0, T, 2T, 3T, … = discrete independent variable values

**56**

$K_{\Delta x}$ Transform Sample and Hold Switch

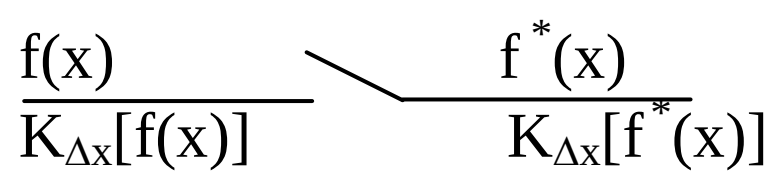

$K_{\Delta x}[f^*(x)] = K_{\Delta x}[f(x)]$

where

f(x) = Switch input function

$f^*(x)$ = Switch output function

$f^*(x)$ = Sample and hold shaped waveform

$0 \le x < \infty$

$\Delta x$ = interval between samples

$f^*(x) = f(n\Delta x)$ for $n\Delta x \le x < [n+1]\Delta x$ , n = 0, 1, 2, 3, …

Sample and Hold Shaped Waveform

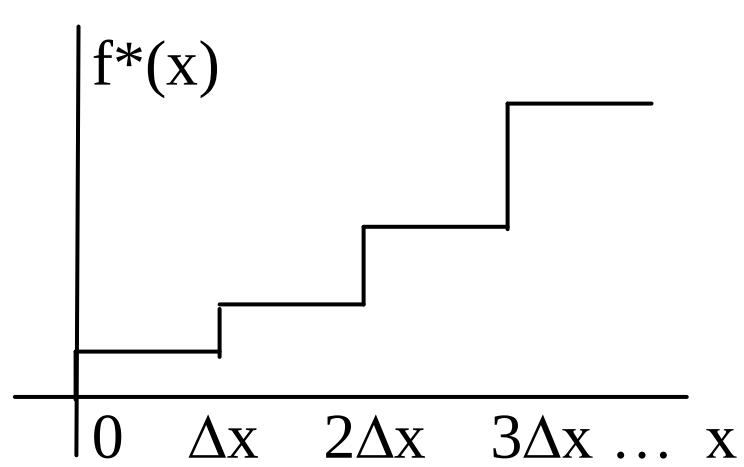

<u>Comments</u> – 1) $f^*(x) = f(x)$ at x = 0, Δx, 2Δx, 3Δx, …

2) $K_{\Delta x}[f^*(x)] = K_{\Delta x}[f(x)]$. If two functions have the same values at x = 0, Δx, 2Δx, 3Δx, …, their $K_{\Delta x}$ Transforms are the same.

3) If the input of a $K_{\Delta x}$ Transform sample and hold switch is designated as f(x) = F(x), the switch's output is designatd as $f^*(x)$ = F(x) also. However, it is understood that the two waveforms generally are not the same. The switch output is a sample and hold shaped waveform with values equal to the switch input values at the sampling instants x = 0, Δx, 2Δx, 3Δx, …

4) If two or more synchronously actived sampling switches are connected in series, the combination of switches acts as one switch.

**57** Finding the Laplace Transform of the output of a sample and hold switch using $K_{\Delta x}$ Transforms

$$L[f^*(x)] = \frac{1}{\mathsf{s}} sK_{\Delta x}[f(x)]\,\Big|_{s = \frac{e^{\mathsf{s}\Delta x} - 1}{\Delta x}} \quad , \quad K_{\Delta x} \text{ Transforms are a function of s}$$

or

$$F^*(\mathsf{s}) = \frac{1}{\mathsf{s}} sK_{\Delta x}[L^{-1}[F(\mathsf{s})]]\,\Big|_{s = \frac{e^{\mathsf{s}\Delta x} - 1}{\Delta x}} \quad , \quad K_{\Delta x} \text{ Transforms are a function of s}$$

or

$$L[f^*(x)] = \frac{f(0)}{\mathsf{s}} + \frac{1}{\mathsf{s}}[\,K_{\Delta x}[D_{\Delta x}f(x)]\,\Big|_{s = \frac{e^{\mathsf{s}\Delta x} - 1}{\Delta x}}\,] \quad , \quad K_{\Delta x} \text{ Transforms are a function of s}$$

where

- $f(x)$ = Switch input function, a continuous function of time, x
- $L[f(x)] = F(\mathsf{s})$ , the Laplace Transform of the switch input function, f(x), a function of $\mathsf{s}$
- $K_{\Delta x}[f(x)]$ , The $K_{\Delta x}$ Transform of the input function f(x) , a function of s
- $f^*(x)$ = Switch output discrete function of time where $x = 0, \Delta x, 2\Delta x, 3\Delta x, \ldots$
- $L[f^*(x)] = F^*(\mathsf{s})$ , the Laplace Transform of the switch output function, $f^*(x)$, a function of $\mathsf{s}$
- $\Delta x$ = Sampling interval
- $x = 0, \Delta x, 2\Delta x, 3\Delta x$
- $\mathsf{s}$ = The Laplace Transform function variable
- s = The $K_{\Delta x}$ Transform function variable

$$L[f^*(x)] = \frac{1}{\mathsf{s}} sK_{\Delta x}[f(x)]\,\Big|_{s = \frac{e^{\mathsf{s}\Delta x} - 1}{\Delta x}} = {}_{1}\sum_{n=0}^{\infty} f(n\Delta x)[\frac{e^{-n\mathsf{s}\Delta x}}{\mathsf{s}} - \frac{e^{-(n+1)\mathsf{s}\Delta x}}{\mathsf{s}}]$$

$$K_{\Delta x}[f^*(x)] = K_{\Delta x}[f(x)]$$

<u>Comment</u> - Note that the Laplace Transform of the output of a sample and hold switch, $L[f^*(t)]$, is a function of the $K_{\Delta t}$ Transform of the input, f(t), to the sample and hold switch. See a diagram of a sample and hold switch shown in the row above.

**58** Equation to convert the $K_{\Delta t}$ Transform of a sample and hold shaped waveform function into its Laplace Transform equivalent

$$L[f^*(t)] = \frac{1}{\mathsf{s}} sK_{\Delta t}[f^*(t)]\,\Big|_{s = \frac{e^{\mathsf{s}\Delta t} - 1}{\Delta t}} \quad , \quad K_{\Delta t} \text{ Transforms are a function of s}$$

where

- $\Delta t$ = sampling interval
- $t = 0, \Delta t, 2\Delta t, 3\Delta t, \ldots$
- $f^*(t)$ = Sample and hold shaped waveform function
- s = $K_{\Delta t}$ Transform variable
- $\mathsf{s}$ = Laplace Transform variable

| | |
|---|---|
| **59** | Equation to convert the Z Transform of a sample and hold shaped waveform function into its Laplace Transform equivalent<br><br>$L[f^*(t)] = \frac{1}{s}\frac{z-1}{z} Z[f^*(t)]\Big\|_{z=e^{sT}}$<br><br>where<br>$T = \Delta t$ = sampling interval<br>$t = 0, T, 2T, 3T, \ldots$<br>$f^*(t)$ = Sample and hold shaped waveform function<br>$z = e^{sT}$ = Z Transform variable<br>$s$ = Laplace Transform variable |
| **60** | Equation to convert the Laplace Transform of a sample and hold shaped waveform function into its $K_{\Delta t}$ Transform equivalent<br><br>$K_{\Delta t}[f^*(t)] = \frac{1}{s}[sL[f^*(t)]]\Big\|_{e^{s\Delta t} = 1+s\Delta t}$<br><br>where<br>$\Delta t$ = sampling interval<br>$t = 0, \Delta t, 2\Delta t, 3\Delta t, \ldots$<br>$f^*(t)$ = Sample and hold shaped waveform function<br>$s = K_{\Delta t}$ Transform variable<br>$s = \frac{1}{\Delta t}\ln(1+s\Delta t)$ , Laplace Transform variable |
| **61** | Equation to convert the Laplace Transform of a sample and hold shaped waveform function into its Z Transform equivalent<br><br>$Z[f^*(t)] = \frac{z}{z-1}[sL[f^*(t)]\,)]\Big\|_{e^{sT} = z}$<br><br>where<br>$T = \Delta t$ = sampling interval<br>$t = 0, T, 2T, 3T, \ldots$<br>$f^*(t)$ = Sample and hold shaped waveform function<br>$z = e^{sT}$ = Z Transform variable<br>$s$ = Laplace Transform variable |
| **62** | Finding the Laplace Transform of the output of a sample and hold switch using Z Transforms<br><br>$L[f^*(t)] = \frac{1}{s}\frac{z-1}{z} Z[f(t)]\Big\|_{z=e^{sT}}$ , Z Transforms are a function of z<br><br>or<br><br>$F^*(s) = \frac{1}{s}\frac{z-1}{z} Z[L^{-1}[F(s)]\Big\|_{z=e^{sT}}$ , Z Transforms are a function of z<br><br>or<br><br>$L[f^*(t)] = \frac{1}{s}[\,f(0) + Tz^{-1} Z[D_T[f(t)]\,]\Big\|_{z=e^{sT}}\,]$ , Z Transforms are a function of z |
| | where<br>$f(t)$ = Sample and hold switch input function, a continuous function of time, t<br>$Z[f(t)] = F(z)$ ,  The Z Transform of the switch input function, f (t), a function of z<br>$f^{*}(t)$ = Sample and hold switch output function, a sample and hold function of time, t<br>$L[f^{*}(t)] = F^{*}(s)$ ,  the Laplace Transform of the switch output function, $f^{*}(x)$, a function of $s$<br>$T = \Delta t$ = Sampling interval<br>$x = 0, T, 2T, 3T, \ldots$<br>$F(s)$ = Laplace Transform of the switch input function, f(t)<br>$F^{*}(s)$ = Laplace Transform of the switch sample and hold output function, $f^{*}(t)$<br>z = The Z Transform function variable<br>$z = e^{sT}$<br>$D_T[f(t)] = \dfrac{f(t+\Delta t) - f(t)}{T}$ ,  Discrete derivative of f(t) with a sampling interval of T<br>$s$ = Laplace Transform variable<br><br>$L[f^{*}(t)] = \dfrac{1}{s}\,\dfrac{z-1}{z}\,Z[f(t)]\Big\|_{z = e^{sT}} = {}_1\sum_{n=0}^{\infty} f(n\Delta t)\left[\dfrac{e^{-ns\Delta t}}{s} - \dfrac{e^{-(n+1)s\Delta t}}{s}\right]$<br>$Z[f^{*}(t)] = Tz^{-1}Z[f(t)]$<br><br><u>Comment</u> - Note that the Laplace Transform of the output of a sample and hold switch, $L[f^{*}(t)]$, is a function of the Z Transform of the input, f(t), to the sample and hold switch. See a diagram of a sample and hold switch shown two rows above. |
| **63** | $Tz^{-1}$ = Z Transform of a unit amplitude pulse of width T initiated at t = 0. |
| **64** | $K_{\Delta t}[f(t)]\Big\|_{s = \frac{e^{s\Delta t} - 1}{\Delta t}} = Tz^{-1}Z[f(t)]\Big\|_{z = e^{sT}}$<br>where<br>$T = \Delta t$ = sampling interval<br>s = $K_{\Delta t}$ Transform variable<br>z = Z Transform variable<br>$s$= Laplace Transform variable |
| **65** | $K_{\Delta t}$ Transform sample and hold sampling of the input and output to a continuous time system<br><br>$K_{\Delta t}[c^{*}(t)] = K_{\Delta t}[g(t)]K_{\Delta t}[f^{*}(t)]$<br>or<br>$C^{*}(s) = G^{*}(s)F^{*}(s)$<br>where<br>$g(t) = L^{-1}[G(s)]$<br><br>Fig. 1 represents the sample and hold sampled system of Fig. 2. The inputs of both Fig. 1 and Fig. 2 are the same for $t = 0, \Delta t, 2\Delta t, 3\Delta t, \ldots$ |

$K_{\Delta t}$ Transform system transfer function

Fig. 1

Sample and hold sampling

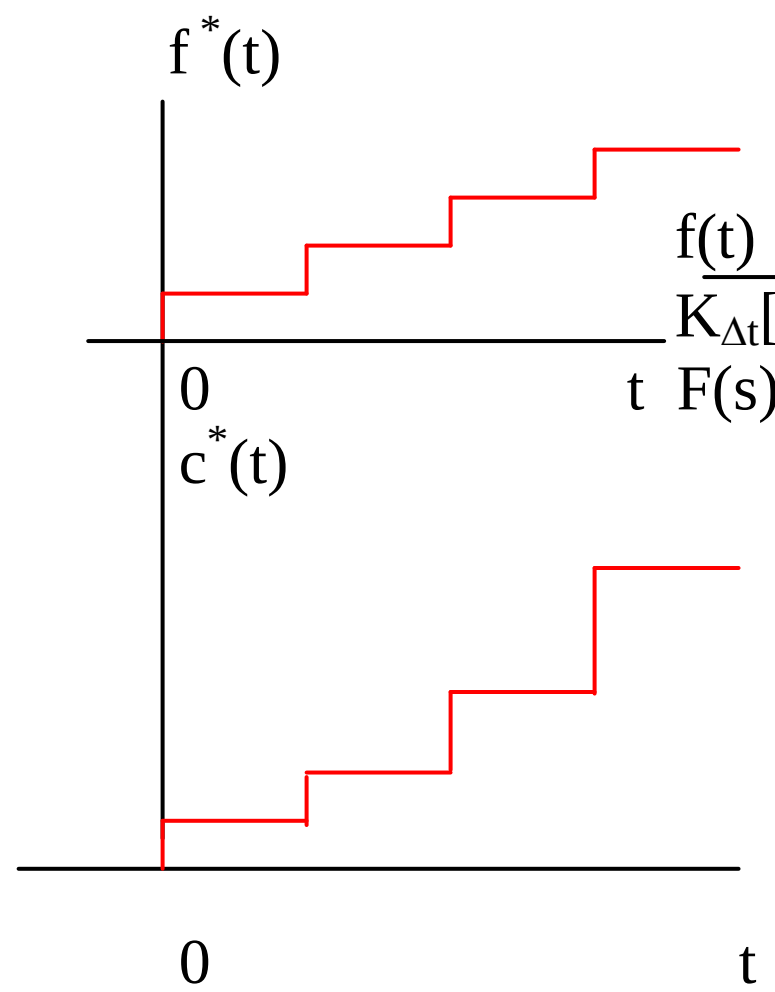

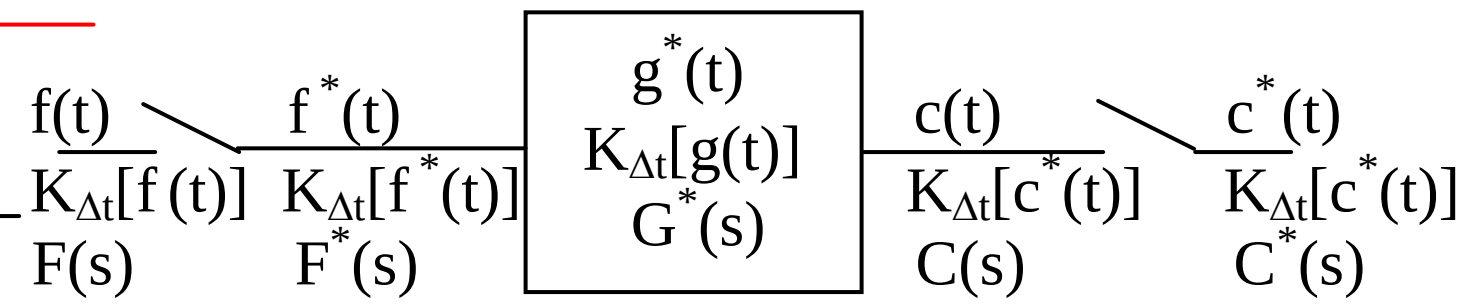

$K_{\Delta t}[c^*(t)] = K_{\Delta t}[g^*(t)]\,K_{\Delta t}[f^*(t)]$

$C^*(s) = G^*(s)F^*(s)$

$G^*(s) = K_{\Delta t}[L^{-1}[G(s)]]$

$g^*(t) = g(t) = L^{-1}[G(s)]$

$G^*(s) = K_{\Delta t}[g(t)] = Tz^{-1}Z[g(t)]\,\big|_{z\,=\,1+s\Delta t}$

$G^*(s) = Tz^{-1}G^*(z)\,\big|_{z\,=\,1+s\Delta t}$

$K_{\Delta t}[f^*(t)] = K_{\Delta t}[f(t)]$

The two switches shown in Fig. 1 above are synchronous sample and hold switches.

Laplace Transform system transfer function

Fig. 2

No sampling

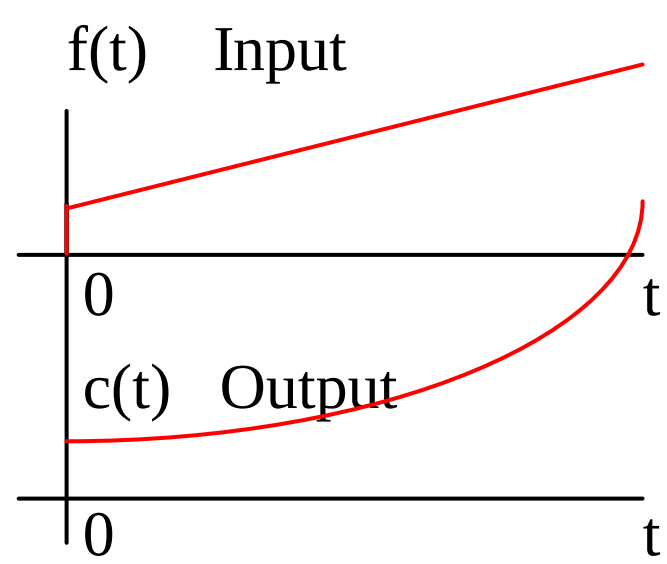

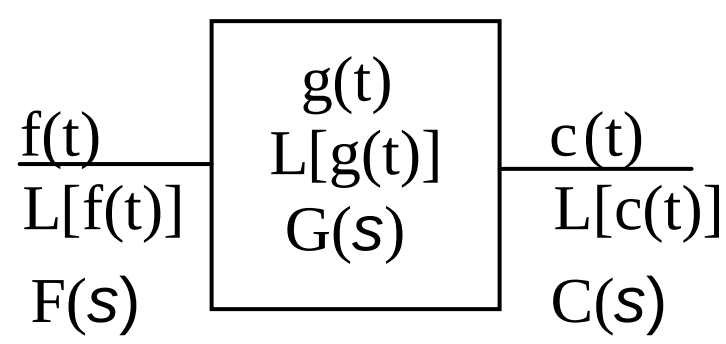

$L[c(t)] = L[g(t)]\,L[f(t)]$

$C(s) = G(s)F(s)$

$G(s) = L[g(t)]$

$g(t) = L^{-1}[G(s)]$

where

- $f(t)$ = system unsampled input function
- $f^*(t)$ = system sample and hold sampled input function
- $F^*(s) = K_{\Delta t}[f^*(t)]$ , The $K_{\Delta t}$ Transform of the sampled system input function, $f^*(t)$
- $g(t) = L^{-1}[G(s)]$ , Inverse Laplace Transform of the system Laplace transfer function, $G(s)$
- $c(t)$ = System unsampled output function
- $c^*(t)$ = System sample and hold sampled output function
- $G^*(s) = K_{\Delta t}[g(t)]$ = System $K_{\Delta t}$ Transform transfer function
- $g^*(t) = K_{\Delta t}^{-1}[G^*(s)]$ , Inverse Laplace Transform of the system $K_{\Delta t}$ Transform transfer function, $G^*(s)$
- $C^*(s) = K_{\Delta t}[c^*(t)]$ , $K_{\Delta t}$ Transform of the sampled output function, $c^*(t)$
- $g^*(t) = g(t) = L^{-1}[G(s)]$
- $G^*(s) = K_{\Delta t}[L^{-1}[G(s)]]$
- $F(s) = L[f(t)]$ , The Laplace Transform of the system input function, $f(t)$
- $C(s) = L[c(t)]$ , The Laplace Transform of the system output function, $c(t)$
- $G(s) = L[g(t)]$ = System Laplace Transform transfer function

$G^*(z) = Z[g(t)] =$ Z Transform of the function $g(t)$

$z = e^{sT} = 1+s\Delta t$

$s = K_{\Delta t}$ Transform variable

$s$ = Laplace Transform variable

To find $G^*(s)$ from $G(s)$ there is a table in the Appendix that may help. See TABLE 3a, The Conversion of Calculus Function Laplace Transforms to Equivalent Function $K_{\Delta t}$ Transforms.

**66** Laplace Transform sample and hold sampling of the input and output to a continuous time system

$L[c^*(t)] = K_{\Delta t}[g(t)]L[f^*(t)]$ , $K_{\Delta t}[g(t)]$ is a function of s

or

$C^*(s) = G^*(s)F^*(s)$

where

$g(t) = L^{-1}[G(s)]$

$$s = \frac{e^{s\Delta t} - 1}{\Delta t}$$

Fig. 1 represents the sample and hold sampled system of Fig. 2.

$K_{\Delta t}$ Transform system transfer function

Fig. 1 $f^*(t)$ Sample and hold sampling

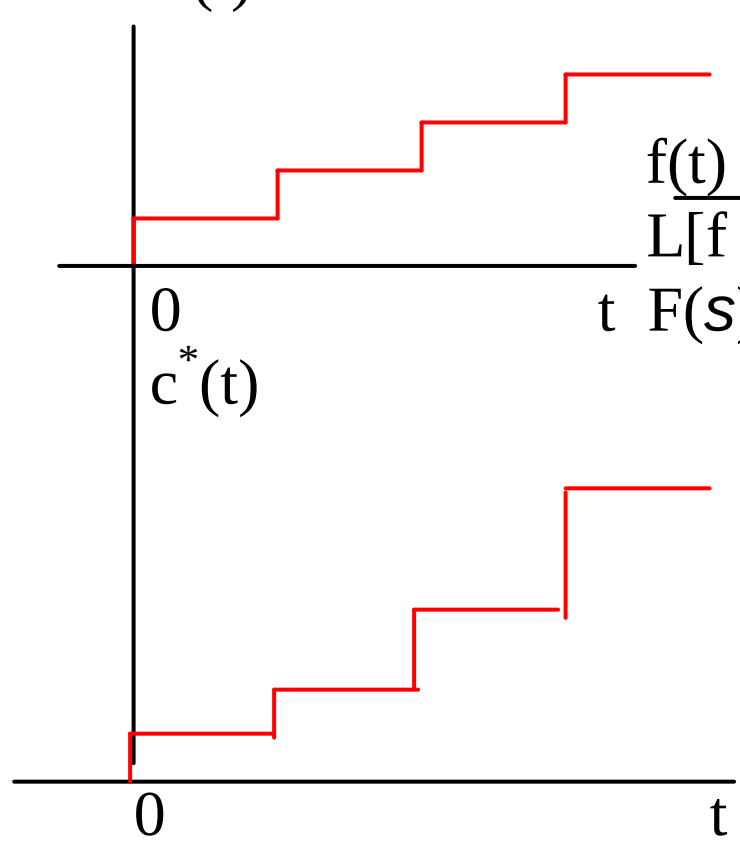

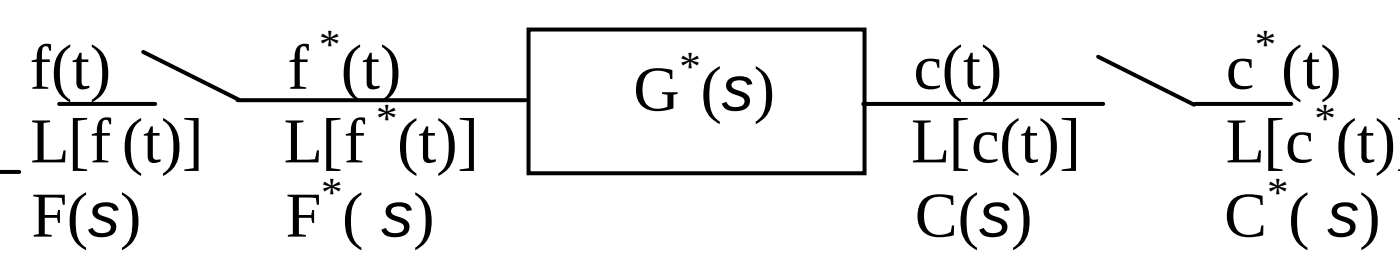

$L[c^*(t)] = K_{\Delta t}[g(t)]L[f^*(t)]$

$C^*(s) = G^*(s)F^*(s)$

$$F^*(s) = L[f^*(t)] = \frac{1}{s} sK_{\Delta t}[L^{-1}[F(s)]]\Big|_{s = \frac{e^{s\Delta t} - 1}{\Delta t}}$$

$$G^*(s) = K_{\Delta t}[L^{-1}[G(s)]]\Big|_{s = \frac{e^{s\Delta t} - 1}{\Delta t}}$$

or

$$G^*(s) = Tz^{-1}Z[L^{-1}[G(s)]]\Big|_{z = e^{s\Delta t}}$$

The two switches shown in Fig. 1 above are synchronous sample and hold switches.

Laplace Transform system transfer function

Fig. 2 No sampling

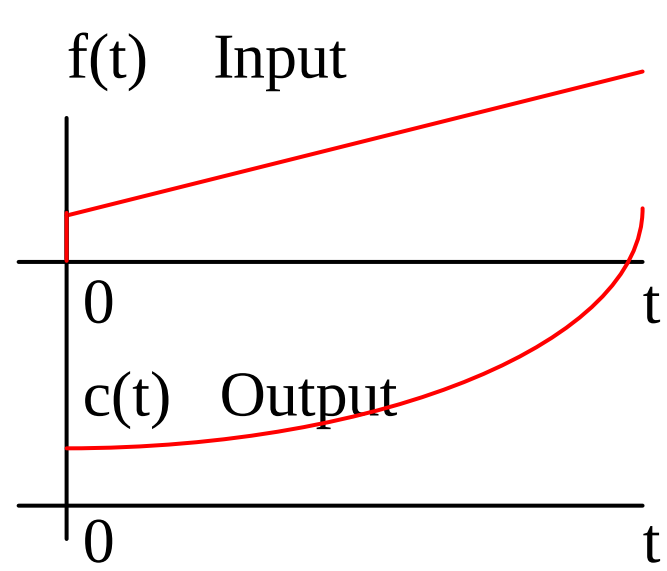

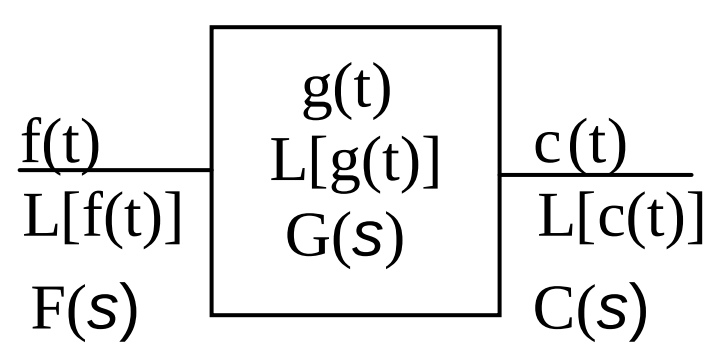

$L[c(t)] = L[g(t)]\, L[f(t)]$

| | |
|---|---|
| | $C(s) = G(s)F(s)$<br>$G(s) = L[g(t)]$<br>$g(t) = L^{-1}[G(s)]$<br>where<br>$f(t)$ = system unsampled input function<br>$F(s) = L[f(t)]$ , The Laplace Transform of the system input function, $f(t)$<br>$g(t) = L^{-1}[G(s)]$ , Inverse Laplace Transform of the system Laplace transfer function, $G(s)$<br>$G(s) = L[g(t)]$ = System Laplace Transform transfer function<br>$c(t)$ = System unsampled output function<br>$C(s) = L[c(t)]$ , Laplace Transform of the output function, $c(t)$<br>$f^*(t)$ = system sample and hold sampled input function<br>$F^*(s) = L[f^*(t)] = \frac{1}{s} sK_{\Delta t}[L^{-1}[F(s)]]$ , The Laplace Transform of the sampled input function, $f^*(t)$<br>$c^*(t)$ = System sample and hold sampled output function<br>$C^*(s) = L[c^*(t)]$ , Laplace Transform of the sampled output function, $c^*(t)$<br>$G^*(s) = K_{\Delta t}[L^{-1}[G(s)]] \Big\vert_{s = \frac{e^{s\Delta t} - 1}{\Delta t}}$<br>or<br>$G^*(s) = Tz^{-1}Z[L^{-1}[G(s)]] \Big\vert_{z = e^{s\Delta t}}$<br>$s = \frac{e^{s\Delta t} - 1}{\Delta t}$<br>$\Delta t$ = Sampling interval<br>$s = K_{\Delta t}$ Transform variable<br>$s$ = Laplace Transform variable<br><br>To find the $K_{\Delta t}$ transform of a function there are Tables 2, 3, 3a, and 3b in the Appendix that may help. |
| **67** | Z Transform sample and hold sampling of the input and output to a continuous time system<br><br>$Z[c^*(t)] = Tz^{-1}Z[g(t)]\,Z[f^*(t)]$<br>or<br>$C^*(z) = Tz^{-1}G^*(z)F^*(z)$<br>where<br>$g(t) = L^{-1}[G(s)]$<br>$G^*(z) = Z[g(t)]$ |

Fig. 1 represents the sample and hold sampled system of Fig. 2.

Z Transform system transfer function

Sample and hold sampling

Fig. 1

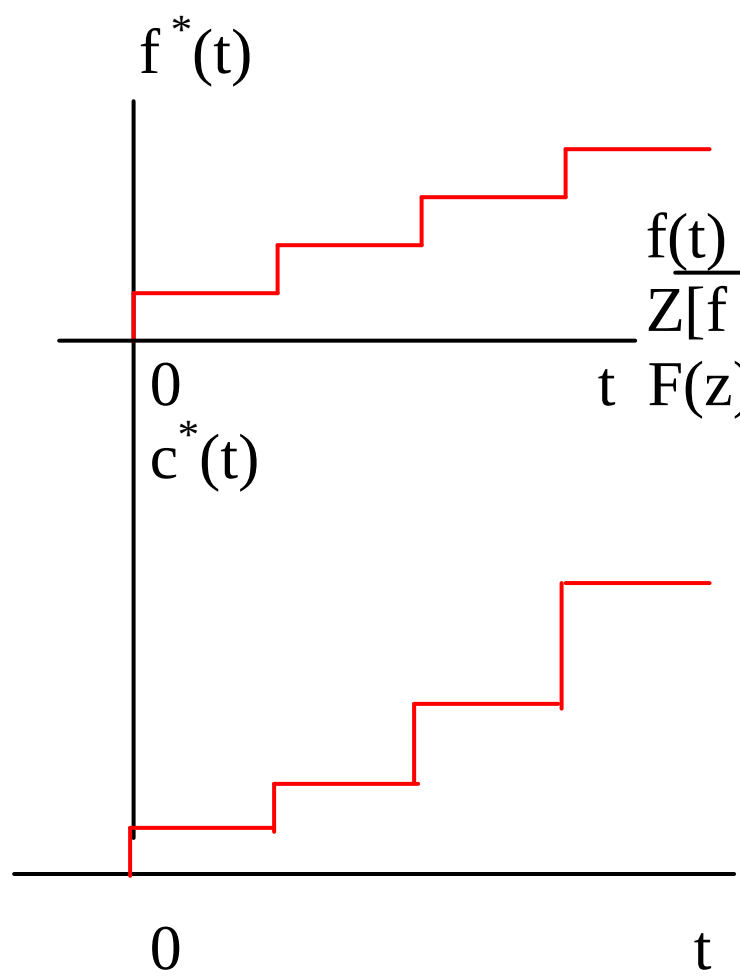

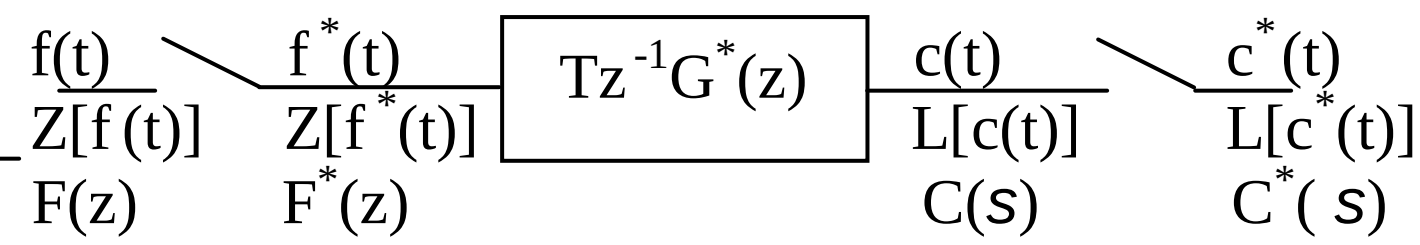

$Z[c^*(t)] = Tz^{-1}Z[g(t)]\,Z[f^*(t)]$

$C^*(z) = Tz^{-1}G^*(z)F^*(z)$

$F^*(z) = Z[f^*(t)]$

$G^*(z) = Z[L^{-1}[G(s)]]$

$g(t) = L^{-1}[G(s)]$

The two switches shown in Fig. 1 above are synchronous sample and hold switches.

Laplace Transform system transfer function

No sampling

Fig. 2

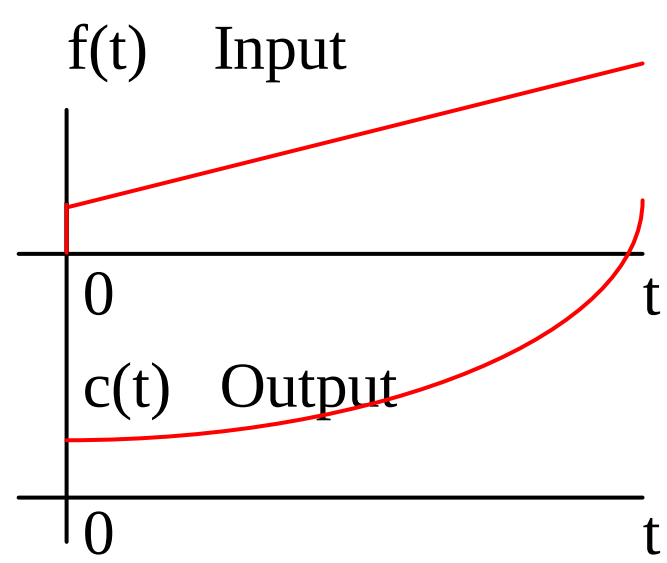

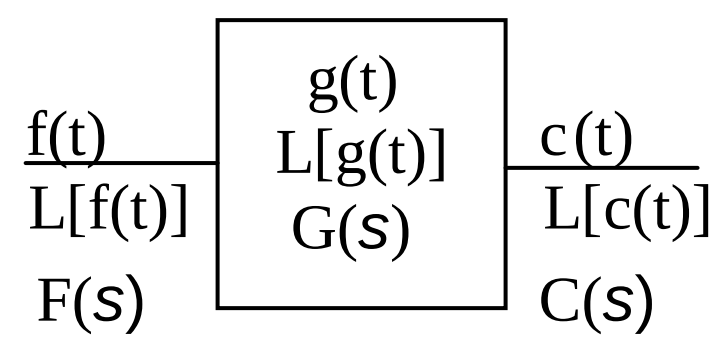

$L[c(t)] = L[g(t)]\,L[f(t)]$

$C(s) = G(s)F(s)$

$G(s) = L[g(t)]$

$g(t) = L^{-1}[G(s)]$

where

- $f(t)$ = system unsampled input function
- $F(z) = Z[f(t)]$, The Z Transform of the system input function, $f(t)$
- $g(t) = L^{-1}[G(s)]$, Inverse Laplace Transform of the system Laplace transfer function, $G(s)$
- $G(s) = L[g(t)]$ = System Laplace Transform transfer function
- $c(t)$ = System unsampled output function
- $C(z) = Z[c(t)]$, Z Transform of the output function, $c(t)$
- $f^*(t)$ = system sample and hold sampled input function
- $F^*(z) = Z[f^*(t)]$, The Z Transform of the sampled input function, $f^*(t)$
- $c^*(t)$ = System sample and hold sampled output function
- $C^*(z) = Z[c^*(t)]$, Z Transform of the sampled output function, $c^*(t)$
- $G^*(z) = Z[g(t)]$
- $\Delta t = T$ = Sampling interval
- $z = e^{s\Delta t}$, Z Transform variable
- $s$ = Laplace Transform variable

**68** Equivalent Sample and Hold $K_{\Delta t}$ Transform and Z Transform Transfer Functions #1

$K_{\Delta t}$ Transform Transfer Function Diagram     Z Transform Transfer Function Diagram

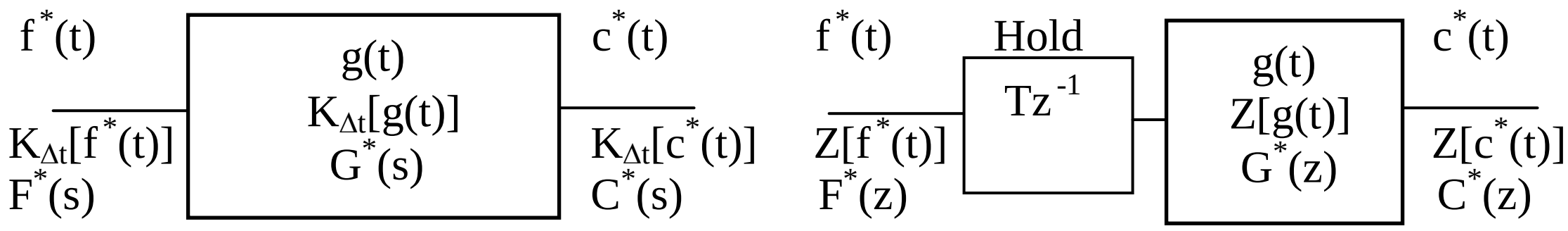

$G^*(s) = K_{\Delta t}[L^{-1}[G(s)]]$     $g(t) = L^{-1}[G(s)]$     $G^*(z) = Z[L^{-1}[G(s)]]$

For the same input to both transfer functions the same output will be obtained.

Equivalent $K_{\Delta t}$ Transform and Z Transform Equations

$K_{\Delta t}[c^*(t)] = K_{\Delta t}[f^*(t)]\,K_{\Delta t}[g(t)]$   or   $C^*(s) = F^*(s)\,G^*(s)$

$Z[c^*(t)] = Z[f^*(t)]Tz^{-1}Z[g(t)]$   or   $C^*(z) = F^*(z)Tz^{-1}G^*(z)$

**69** Equivalent Sample and Hold $K_{\Delta t}$ Transform and Z Transform Transfer Functions #2

$K_{\Delta t}$ Transform Transfer Function Diagram     Z Transform Transfer Function Diagram

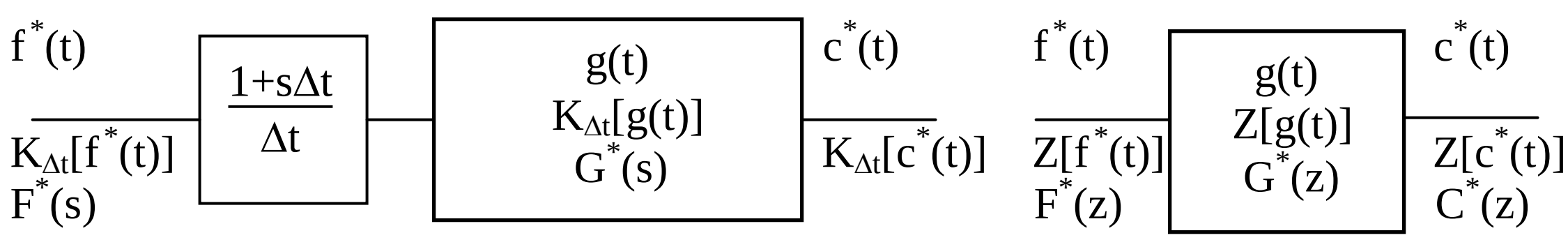

$G^*(s) = K_{\Delta t}[L^{-1}[G(s)]]$     $g(t) = L^{-1}[G(s)]$     $G^*(z) = Z[L^{-1}[G(s)]]$

$$\frac{1+s\Delta t}{\Delta t} = \left.\frac{z}{T}\right|_{\substack{z = 1+s\Delta t \\ T = \Delta t}}$$ , This term cancels the $K_{\Delta t}$ Transform sample and hold.

For the same input to both transfer functions the same output will be obtained.

Equivalent $K_{\Delta t}$ Transform and Z Transform Equations

$K_{\Delta t}[c^*(t)] = K_{\Delta t}[f^*(t)]\dfrac{1+s\Delta t}{\Delta t}K_{\Delta t}[g(t)]$   or   $C^*(s) = F^*(s)\dfrac{1+s\Delta t}{\Delta t}G^*(s)$

$Z[c^*(t)] = Z[f^*(t)]Z[g(t)]$   or   $C^*(z) = F^*(z)\,G^*(z)$

**70** Comparison of the $K_{\Delta t}$ Transform and the Z Transform

$K_{\Delta t}$ Transform Series Expansion

$$K_{\Delta t}[f(t)] = \Delta t[f(0)(1+s\Delta t)^{-1} + f(\Delta t)(1+s\Delta t)^{-2} + f(2\Delta t)(1+s\Delta t)^{-3} + f(3\Delta t)(1+s\Delta t)^{-4} + \ldots]$$

Z Transform Series Expansion

$$Z[f(t)] = f(0)z^{-0} + f(T)z^{-1} + f(2T)z^{-2} + f(3T)z^{-3} + \ldots$$

$K_{\Delta t}$ Transform to Z Transform Conversion

$$Z[f(t)] = F(z) = \frac{z}{T} F(s)\Big|_{s = \frac{z-1}{T}}$$

$T = \Delta t$ sampling period

$t = nT$

$n = 0,1,2,3,\ldots$

$Z[f(t)] = F(z)$ Z Transform of f(t)

$K_{\Delta t}[f(t)] = F(s)$ $K_{\Delta t}$ Transform of f(t)

Z Transform to $K_{\Delta t}$ Transform Conversion

$$K_{\Delta t}[f(t)] = F(s) = Tz^{-1} F(z)\Big|_{z = 1+s\Delta t}$$

$\Delta t = T$ sampling period

$t = n\Delta t$

$n = 0, 1, 2, 3, \ldots$

$K_{\Delta t}[f(t)] = F(s)$ $K_{\Delta t}$ Transform of f(t)

$Z[f(t)] = F(z)$ Z Transform of f(t)

where

f(t) = function of t

$\Delta t = T$ = Interval between successive values of t

**71** Modified $K_{\Delta t}$ Transform

$$K\Delta t_M[f(t)] = K_{\Delta t}[f(t+m\Delta t)]$$

or

$$K\Delta t_M[f(t)] = \Delta t[f(m\Delta t)(1+s\Delta t)^{-1} + f(\Delta t+m\Delta t)(1+s\Delta t)^{-2} + f(2\Delta t+m\Delta t)(1+s\Delta t)^{-3} + f(3\Delta t+m\Delta t)(1+s\Delta t)^{-4} + \ldots]$$

or

$$K\Delta t_M[f(t)] = TZ_M[f(t+mT)]\Big|_{\substack{z = 1+s\Delta t \\ \Delta t = T}}$$

or

$$K\Delta t_M[f(t)] = Tz^{-1}Z[f(t+mT)]\Big|_{\substack{z = 1+s\Delta t \\ \Delta t = T}}$$

where

$K\Delta t_M[f(t)]$ = Modified $K_{\Delta t}$ Transform of the function, f(t)

$Z_M[f(t)]$ = Modified Z Transform of f(t)

$K_{\Delta t}[f(t)] = K_{\Delta t}$ Transform of f(t)

$Z[f(t)]$ = Z Transform of f(t)

$0 \le m \le 1$

$T = \Delta t$, Interval between successive values of t

$t = t + m\Delta t$

$t = mT, T+mT, 2T+mT, 3T+mT, \ldots$

**72** Modified Z Transform

$$Z_M[f(t)] = z^{-1}{}_1 \sum_{k=0}^{\infty} f(kT+mT)z^{-k}$$

or

$$Z_M[f(t)] = [f(mT)z^{-1} + f(1T+mT)z^{-2} + f(2T+mT)z^{-3} + f(3T+mT)z^{-4} + \ldots]$$

or

$$Z_M[f(t)] = z^{-1}Z[f(t+mT)]$$

or

$$Z_M[f(t)] = z^{-1}\left[\frac{1}{T}\int_{0}^{\infty} f(t+mT)z^{-\frac{t}{\Delta t}}\Delta t\right]_T$$

or

$$Z_M[f(t)] = \frac{1}{T}K_{\Delta t}[f(t+m\Delta t)]\Big|_{\substack{s = \frac{z-1}{T} \\ T = \Delta t}}$$

or

$$Z_M[f(t)] = \frac{1}{T}K\Delta t_M[f(t)]\Big|_{\substack{s = \frac{z-1}{T} \\ T = \Delta t}}$$

where

$Z_M[f(t)]$ = Modified Z Transform of the function, f(t)
$K\Delta t_M[f(t)]$ = Modified $K_{\Delta t}$ Transform of f(t)
$Z[f(t)]$ = Z Transform of f(t)
$K_{\Delta t}[f(t)] = K_{\Delta t}$ Transform of f(t)
$0 \le m \le 1$
$T = \Delta t$, Interval between successive values of t
$t = t + m\Delta t$
$t = mT, T+mT, 2T+mT, 3T+mT, \ldots$

**73** Finding the Inverse $K_{\Delta t}$ Transform Using the $K_{\Delta t}$ Transform Asymptotic Series and Discrete Maclaurin Series

$$F(s) = K_{\Delta t}[f(t)] = \sum_{n=1}^{\infty} D_{\Delta t}^{\,n-1}f(0)s^{-n} = f(0)s^{-1} + D_{\Delta t}^{\,1}f(0)s^{-2} + D_{\Delta t}^{\,2}f(0)s^{-3} + D_{\Delta t}^{\,3}f(0)s^{-4} + \ldots$$

where

$D_{\Delta t}^{\,m}f(0) = D_{\Delta t}^{\,m}f(t)|_{t=0}$ , m = 0,1,2,3,… , Discrete derivatives of f(t) evaluated at t = 0
f(t) = function of t
$F(s) = K_{\Delta t}[f(t)] = K_{\Delta t}$ Transform of the function, f(t)
$\Delta t$ = Interval between successive values of t

The asymptotic expansion of the $K_{\Delta t}$ Transform, F(s), is obtained by the division of its denominator into its numerator or by the use of a computer program.

Note – The $K_{\Delta t}$ Transform Asymptotic Series is a series expansion of a $K_{\Delta t}$ Transform at $x = \infty$.

| | |
|---|---|
| | Finding $f(t) = K_{\Delta t}^{-1}[F(s)]$ , the Inverse $K_{\Delta t}$ Transform<br><br>The evaluated values of the initial condition discrete derivatives and the value of $\Delta t$ are entered into a discrete Maclaurin Series.<br><br>$$f(t) = \sum_{m=0}^{\infty} \frac{1}{m!} D_{\Delta t}^{m} f(t)\Big\vert_{t=0} [t]_{\Delta t}^{m}$$<br><br>Expanding<br><br>$$f(t) = f(0) + \frac{D_{\Delta t}^{1}f(0)}{1!}\, t + \frac{D_{\Delta t}^{2}f(0)}{2!}\, t(t-\Delta t) + \frac{D_{\Delta t}^{3}f(0)}{3!}\, t(t-\Delta t)(t-2\Delta t) + \frac{D_{\Delta t}^{4}f(0)}{4!}\, t(t-\Delta t)(t-2\Delta t)(t-3\Delta t) + \ldots$$<br><br>The above discrete Maclaurin Series is used to evaluate $f(t)$<br><br>If the discrete Maclaurin Series is truncated at the rth term, the series will calculate exact values for $f(t)$ at $t = 0, \Delta t, 2\Delta t, 3\Delta t, \ldots, (r-1)\Delta t$. In general, for larger values of t there will be some error and series convergence must be investigated. The discrete Maclaurin Series can be used to interpolate the value of $f(t)$ for $t = m\Delta t$ where m are the real numbers but not integers between 0 and r-1. |
| **74** | Finding the Inverse $K_{\Delta t}$ Transform Using the $K_{\Delta t}$ Transform Asymptotic Series and the Definition of Discrete Derivatives<br><br>$$F(s) = K_{\Delta t}[f(t)] = \sum_{n=1}^{\infty} D_{\Delta t}^{n-1}f(0)s^{-n} = f(0)s^{-1} + D_{\Delta t}^{1}f(0)s^{-2} + D_{\Delta t}^{2}f(0)s^{-3} + D_{\Delta t}^{3}f(0)s^{-4} + \ldots$$<br>where<br>$D_{\Delta t}^{m}f(0) = D_{\Delta t}^{m}f(t)\vert_{t=0}$ , $m = 0,1,2,3,\ldots$ , Discrete derivatives of $f(t)$ evaluated at $t = 0$<br>$f(t)$ = function of t<br>$F(s) = K_{\Delta t}[f(t)] = K_{\Delta t}$ Transform of the function, $f(t)$<br>$\Delta t$ = Interval between successive values of t<br><br>The asymptotic expansion of the $K_{\Delta t}$ Transform, $F(s)$, is obtained by the division of its denominator into its numerator or by the use of a computer program.<br><br><u>Note</u> – The $K_{\Delta t}$ Transform Asymptotic Series is a series expansion of a $K_{\Delta t}$ Transform at $x = \infty$.<br><br>Finding $f(t)$ at $t = 0, \Delta t, 2\Delta t, 3\Delta t, \ldots$<br><br>$f(0) = D_{\Delta t}^{0}f(0)$<br><br>$$\frac{f(\Delta t) - f(0)}{\Delta t} = D_{\Delta t}^{1}f(0)$$ |

$$\frac{f(2\Delta t) - 2f(\Delta t) + f(0)}{\Delta t^2} = D_{\Delta t}{}^2 f(0)$$

$$\frac{f(3\Delta t) - 3f(2\Delta t) + 3f(\Delta t) - f(0)}{\Delta t^3} = D_{\Delta t}{}^3 f(0)$$

$$\frac{f(4\Delta t) - 4f(3\Delta t) + 6f(2\Delta t) - 4f(\Delta t) + f(0)}{\Delta t^4} = D_{\Delta t}{}^4 f(0)$$

.
.
.

The coefficients of the above equations can be obtained by the Binomial Expansion of $(1+1)^n$ where n = 0, 1, 2, 3, …

**75** Finding the Inverse $K_{\Delta t}$ Transform Using the $K_{\Delta t}$ Transform Asymptotic Series and the Inverse $K_{\Delta t}$ Transform Formula

$$K_{\Delta t}{}^{-1}[K_{\Delta t}[f(t)]] = f(n\Delta t) = \sum_{p=0}^{n} [{}_nC_p \Delta t^p]\, a_{p+1} \quad , \quad \text{The Inverse } K_{\Delta t} \text{ Transform Formula}$$

where

$K_{\Delta t}[f(t)] = F(s) = a_1 s^{-1} + a_2 s^{-2} + a_3 s^{-3} + a_4 s^{-4} + \dots$ , Asymptotic Series form of F(s)

$F(s) = K_{\Delta t}[f(t)] = K_{\Delta t}$ Transform of the function, f(t)

$K_{\Delta t}{}^{-1}[F(s)] = f(t) =$ Inverse $K_{\Delta t}$ Transform of the function, F(s)

n = 0, 1, 2, 3, 4, …

p = 0, 1, 2, 3, … , n

$nC_p = \dfrac{n!}{r!(n-r)!}$

0! = 1

f(t) = function of t

$\Delta t$ = interval between successive values of t

$t = 0, \Delta t, 2\Delta t, 3\Delta t, \dots$

$D_{\Delta t}{}^m f(0) = a_{m+1}$ , $m = 0,1,2,3,\dots$ , Initial conditions of f(t)

**76** Finding the Inverse $K_{\Delta t}$ Transform Using the Z Transform

$$F(z) = \frac{z}{\Delta t} F(s)\Big|_{s = \frac{z-1}{\Delta t}} = \sum_{n=0}^{\infty} f(n\Delta t) z^{-n} = f(0) + f(1\Delta t)z^{-1} + f(2\Delta t)z^{-2} + f(3\Delta t)z^{-3} + f(4\Delta t)z^{-3} + \dots$$

where

$F(s) = K_{\Delta t}[f(t)]$ , $K_{\Delta t}$ Transform of the function, f(t)

$F(z) = Z[f(t)]$ , Z Transform of the function, f(t)

f(t) = function of t

$\Delta t$ = Interval between successive values of t

$t = n\Delta t$ , $n = 0,1,2,3,\dots$

$z = 1+s\Delta t$

| | |
|---|---|
| | For the $K_{\Delta t}$ Transform, F(s), the maximum order of the numerator polynomial is one less than the order of the denominator polynomial. For the Z Transform, F(z), the maximum order of the numerator polynomial is the same as the order of the denominator polynomial.<br><br>Finding f(t) at t = 0, Δt, 2Δt, 3Δt, …<br><br>The Z Transform series expansion , which is derived from a $K_{\Delta t}$ Transform, F(s), finds the values of f(t) where $f(t) = Z^{-1}[F(z)]$. The coefficients of the Z Transform series expansion represent the values of f(t) at t = nΔt , n = 0,1,2,3,… |
| **77** | The $K_{\Delta t}$ Transform Unit Amplitude Pulse Series<br><br>$F(s) = \sum_{n=0}^{\infty} f(n\Delta t)[\Delta t(1+s\Delta t)^{-n-1}] = f(0)[\Delta t(1+s\Delta t)^{-1}] + f(\Delta t)[\Delta t(1+s\Delta t)^{-2}] + f(2\Delta t)[\Delta t(1+s\Delta t)^{-3}] + f(3\Delta t)[\Delta t(1+s\Delta t)^{-4}] + \dots$<br>where<br>f(t) = function of t<br>$F(s) = K_{\Delta t}[f(t)]$ , $K_{\Delta t}$ Transform of the function f(t)<br>Δt = Interval between successive values of t<br>$\Delta t(1+s\Delta t)^{-n-1}$ = $K_{\Delta t}$ Transform of a unit amplitude pulse where the pulse rises at t = nΔt and falls at t = (n+1)Δt<br>n = 0,1,2,3,…<br><br>The coefficients of the $K_{\Delta t}$ Transform Unit Amplitude Pulse Series represent the values of f(t) for t = nΔt , n = 0,1,2,3,…<br><br>Finding f(t) at t = 0, Δt, 2Δt, 3Δt, …<br><br>The values of f(t) at t = 0, Δt, 2Δt, 3Δt, … are the coefficients of the $K_{\Delta t}$ Transform Unit Amplitude Pulse Series. |
| **78** | Finding the Inverse Laplace Transform Using the Laplace Transform Asymptotic Series and Maclaurin Series<br><br>$F(s) = L[f(t)] = \sum_{n=1}^{\infty} \frac{d^{n-1}}{dt^{n-1}} f(0)s^{-n} = f(0)s^{-1} + \frac{d}{dt} f(0)s^{-2} + \frac{d^2}{dt^2} f(0)s^{-3} + \frac{d^3}{dt^3} f(0)s^{-4} + \dots$<br>where<br>$f^{(m)}(0) = \frac{d^m}{dt^m} f(0) = \frac{d^m}{dt^m} f(t)\|_{t=0}$ , m = 0,1,2,3,… , Derivatives of f(t) evaluated at t = 0<br>f(t) = function of t<br>L[f(t)] = Laplace Transform of the function, f(t) |
| | Finding $f(t) = L^{-1}[F(s)]$, the Inverse Laplace Transform<br><br>The evaluated values of the initial condition derivatives are entered into a Maclaurin Series.<br><br>$$f(t) = \sum_{m=0}^{\infty} \frac{1}{m!} \frac{d^m}{dt^m} f(t)\Big\|_{t=0} t^m = \sum_{m=0}^{\infty} \frac{1}{m!} f^{(m)}(0)\, t^m$$<br><br>Expanding<br><br>$$f(t) = f(0) + \frac{1}{1!} f^{(1)}(0)\, t + \frac{1}{2!} f^{(2)}(0)\, t^2 + \frac{1}{3!} f^{(3)}(0)\, t^3 + \frac{1}{4!} f^{(4)}(0)\, t^4 + \dots$$<br><br>The above Maclaurin Series is used to evaluate f(t). |
| **79** | Finding the Inverse Z Transform Using the $K_{\Delta t}$ Transform Asymptotic Series and the Discrete Maclaurin Series<br><br>$K_{\Delta t}[f(t)] = Tz^{-1} Z[f(t)]\Big\|_{z = 1+s\Delta t,\ \Delta t = T}$, Z Transform to $K_{\Delta t}$ Transform Conversion<br><br>$$K_{\Delta t}[f(t)] = \sum_{n=1}^{\infty} D_{\Delta t}^{n-1} f(0) s^{-n} = f(0)s^{-1} + D_{\Delta t}^{1} f(0)s^{-2} + D_{\Delta t}^{2} f(0)s^{-3} + D_{\Delta t}^{3} f(0)s^{-4} + \dots$$<br>where<br>$D_{\Delta t}^{m} f(0) = D_{\Delta t}^{m} f(t)\|_{t=0}$, $m = 0,1,2,3,\dots$, Discrete derivatives of f(t) evaluated at t = 0<br>f(t) = function of t<br>Z[f(t)] = Z Transform of the function f(t)<br>$K_{\Delta t}[f(t)] = K_{\Delta t}$ Transform of the function, f(t)<br>$\Delta t$ = Interval between successive values of t<br><br>The asymptotic expansion of the $K_{\Delta t}$ Transform, F(s), is obtained by the division of its denominator into its numerator or by the use of a computer program.<br><br><u>Note</u> – The $K_{\Delta t}$ Transform Asymptotic Series is a series expansion of a $K_{\Delta t}$ Transform at $x = \infty$.<br><br>Finding $f(t) = Z^{-1}[f(t)] = K_{\Delta t}^{-1}[f(t)]$, the Inverse Z Transform<br><br>The evaluated values of the initial condition discrete derivatives and the value of $\Delta t$ are entered into a discrete Maclaurin Series. |

$$f(t) = \sum_{m=0}^{\infty} \frac{1}{m!} D_{\Delta t}^{m} f(t)\Big|_{t=0} [t]_{\Delta t}^{m}$$

Expanding

$$f(t) = f(0) + \frac{D_{\Delta t}^{1}f(0)}{1!}\, t + \frac{D_{\Delta t}^{2}f(0)}{2!}\, t(t-\Delta t) + \frac{D_{\Delta t}^{3}f(0)}{3!}\, t(t-\Delta t)(t-2\Delta t) + \frac{D_{\Delta t}^{4}f(0)}{4!}\, t(t-\Delta t)(t-2\Delta t)(t-3\Delta t) + \ldots$$

The above discrete Maclaurin Series is used to evaluate f(t)

If the discrete Maclaurin Series is truncated at the rth term, the series will calculate exact values for f(t) at $t = 0, \Delta t, 2\Delta t, 3\Delta t, \ldots, (r-1)\Delta t$. In general, for larger values of t there will be some error and series convergence must be investigated. The discrete Maclaurin Series can be used to interpolate the value of f(t) for $t = m\Delta t$ where m are the real numbers but not integers between 0 and r-1.

<u>Comment</u>

The Inverse Z Transform of f(t), $Z^{-1}[f(t)]$, is usually obtained from the Z Transform equation shown below.

$$F(z) = \sum_{n=0}^{\infty} f(n\Delta t)z^{-n} = f(0) + f(1\Delta t)z^{-1} + f(2\Delta t)z^{-2} + f(3\Delta t)(z^{-3} + f(4\Delta t)z^{-3} + \ldots$$

However, the discrete Maclaurin series representation of f(t) provides more information concerning f(t). Besides providing the values of f(t) at $t = 0, \Delta t, 2\Delta t, 3\Delta t, \ldots$, it provides interpolated values of f(t) between $t = 0, \Delta t, 2\Delta t, 3\Delta t, \ldots$ and the discrete derivatives at $t = 0$.

Also, the discrete Maclaurin Series can be obtained from the values of f(t) obtained from the above Z Transform series. The discrete derivatives at $t = 0$ can be calculated using the following equations.

$$D_{\Delta t}^{0}f(0) = f(0)$$

$$D_{\Delta t}^{1}f(0) = \frac{f(\Delta t) - f(0)}{\Delta t}$$

$$D_{\Delta t}^{2}f(0) = \frac{f(2\Delta t) - 2f(\Delta t) + f(0)}{\Delta t^{2}}$$

$$D_{\Delta t}^{3}f(0) = \frac{f(3\Delta t) - 3f(2\Delta t) + 3f(\Delta t) - f(0)}{\Delta t^{3}}$$

$$D_{\Delta t}^{4}f(0) = \frac{f(4\Delta t) - 4f(3\Delta t) + 6f(2\Delta t) - 4f(\Delta t) + f(0)}{\Delta t^{4}}$$

.
.
.

| | |
|---|---|
| | The coefficients of the previous equations can be obtained by the Binomial Expansion of $(1+1)^n$ where n = 0, 1, 2, 3, … |
| **80** | Heavyside Expansion Equations<br><br>1) $f(x) = K_{\Delta x}^{-1}[f(s)] = \sum_{n=1}^{m} \frac{p(a_n)}{q^{'}(a_n)} e_{\Delta x}(a_n,x) = \sum_{n=1}^{m} \frac{p(a_n)}{Q_n(a_n)} e_{\Delta t}(a_n,x)$<br><br>where<br>$K_{\Delta t}[f(x)] = f(s)$<br>$f(s) = \frac{p(s)}{q(s)}$<br>$q(s) = (s-a_1)(s-a_2)(s-a_3) \dots (s-a_m)$<br>p(s) = polynomial of degree < m<br>$e_{\Delta x}(a_n,x) = (1+a_n\Delta x)^{\frac{x}{\Delta x}}$<br>$\Delta x$ = x increment<br>$x = 0, \Delta x, 2\Delta x, 3\Delta x, \dots$<br>$Q_n(s)$ = product of all factors of q(s) except the factor, $s-a_n$<br>$q^{'}(s)$ = first derivative of q(s)<br><br>2) $f(x) = K_{\Delta x}^{-1}[f(s)] = \sum_{n=1}^{r} \frac{p^{(r-n)}(a)\,(1+a\Delta x)^{1-n}}{(r-n)!\,(n-1)!} [x]_{\Delta x}^{n-1} e_{\Delta x}(a,x)$<br><br>where<br>$K_{\Delta x}[f(x)] = f(s)$<br>$f(s) = \frac{p(s)}{(s-a)^r}$<br>p(s) = polynomial of degree < r<br>r = 1, 2, 3, … , the order of the pole at a<br>$e_{\Delta x}(a,x) = (1+a\Delta x)^{\frac{x}{\Delta x}}$<br>$\Delta x$ = x increment<br>$x = 0, \Delta x, 2\Delta x, 3\Delta x, \dots$ |

**Note:** $x = 0$ is $x = 0^+$

# TABLE 3

## $K_{\Delta x}$ Transforms
## (Compared to the Laplace Transform)

$$K_{\Delta x}[f(x)] = {}_{\Delta x}\int_0^{\infty} (1+s\Delta x)^{-\left(\frac{x+\Delta x}{\Delta x}\right)} f(x)\Delta x$$

$L[f(x)] = \lim_{\Delta x\to 0} K_{\Delta x}[f(x)] \equiv K_0[f(x)]$ , $L[f(x)]$ is the Laplace Transform

$x = m\Delta x$ , $m = 0,1,2,3,\dots$ , $\Delta x \neq 0$

<u>Comment</u> – There are additional $K_{\Delta x}$ Transforms in Table 3a and Table 3b

| # | $f_{\Delta x}(x)$ | $K_{\Delta x}[f_{\Delta x}(x)]$ | f(x) | L[f(x)] |
|---|---|---|---|---|
| **1** | 1<br>or<br>U(x) | $\frac{1}{s}$ | 1<br>or<br>U(x) | $\frac{1}{s}$ |
| **2** | $U(x - n\Delta x)$ | $\frac{(1+s\Delta x)^{-n}}{s}$ | $U(x - n\Delta x)$ | $\frac{e^{-sn\Delta x}}{s}$ |
| **3** | c | $\frac{c}{s}$ | c | $\frac{c}{s}$ |
| **4** | x | $\frac{1}{s^2}$ | x | $\frac{1}{s^2}$ |
| **5** | $e_{\Delta x}(a,x)$<br><br>or<br><br>$e_{\Delta x}(a,x)$ variable subscript identity<br>$e_{m\Delta x}(\frac{(1+a\Delta x)^m-1}{m\Delta x},t)$<br>$m = 1,2,3,\dots$ | $\frac{1}{s-a}$<br>root $s = a$<br><br>The related Z Transform is:<br>$\frac{z}{z-(1+aT)}$<br>$T = \Delta x$ | $e^{ax}$ | $\frac{1}{s-a}$<br>root $s = a$ |
| **6** | $\sin_{\Delta x}(b,x)$ | $\frac{b}{s^2+b^2}$<br>roots $s = jb, -jb$<br><br>The related Z Transform is:<br>$\frac{bTz}{z^2-2z+(1+b^2T^2)}$ | sinbx | $\frac{b}{s^2+b^2}$<br>roots $s = jb, -jb$ |

| # | $f_{\Delta x}(x)$ | $K_{\Delta x}[f_{\Delta x}(x)]$ | f(x) | L[f(x)] |
|---|---|---|---|---|
| **7** | $\cos_{\Delta x}(b,x)$ | $\frac{s}{s^2+b^2}$<br>roots $s = jb, -jb$<br><br>The related Z Transform is:<br>$\frac{z(z-1)}{z^2-2z+(1+b^2T^2)}$ | cosbx | $\frac{s}{s^2+b^2}$<br>roots $s = jb, -jb$ |
| **8** | $1+a\Delta x \neq 0$<br>$e_{\Delta x}(a,x)\sin_{\Delta x}(\frac{b}{1+a\Delta x}, x)$<br>or<br>$(1+a\Delta x)^{\frac{x}{\Delta x}}\sin_{\Delta x}(\frac{b}{1+a\Delta x}, x)$<br>$1+a\Delta x = 0$<br>$[b\Delta x]^{\frac{x}{\Delta x}} \sin\frac{\pi x}{2\Delta x}$ | $\frac{b}{(s-a)^2+b^2}$<br>roots $s = a+jb, a-jb$<br><br>The related Z Transform is:<br>$\frac{bTz}{z^2-2[1+aT]z+[(1+aT)^2+b^2T^2]}$ | $e^{ax}\sin bx$ | $\frac{b}{(s-a)^2+b^2}$<br>roots $s = a+jb, a-jb$ |
| **9** | $1+a\Delta x \neq 0$<br>$e_{\Delta x}(a,x)\cos_{\Delta x}(\frac{b}{1+a\Delta x}, x)$<br>or<br>$(1+a\Delta x)^{\frac{x}{\Delta x}}\cos_{\Delta x}(\frac{b}{1+a\Delta x}, x)$<br>$1+a\Delta x = 0$<br>$[b\Delta x]^{\frac{x}{\Delta x}} \cos\frac{\pi x}{2\Delta x}$ | $\frac{s-a}{(s-a)^2+b^2}$<br>roots $s = a+jb, a-jb$<br><br>The related Z Transform is:<br>$\frac{z^2- (1+aT)z}{z^2-2[1+aT]z+[(1+aT)^2+b^2T^2]}$ | $e^{ax}\cos bx$ | $\frac{s-a}{(s-a)^2+b^2}$<br>roots $s = a+jb, a-jb$ |
| **10** | $\sinh_{\Delta x}(b,x)$ | $\frac{b}{s^2-b^2}$<br>roots $s = b, -b$<br><br>The related Z Transform is:<br>$\frac{bTz}{z^2-2z+(1-b^2T^2)}$ | sinhbx | $\frac{b}{s^2-b^2}$<br>roots $s = b, -b$ |
| **11** | $\cosh_{\Delta x}(b,x)$ | $\frac{s}{s^2-b^2}$<br>roots $s = b, -b$<br><br>The related Z Transform is:<br>$\frac{z(z-1)}{z^2-2z+(1-b^2T^2)}$ | coshbx | $\frac{s}{s^2-b^2}$<br>roots $s = b, -b$ |

| # | $f_{\Delta x}(x)$ | $K_{\Delta x}[f_{\Delta x}(x)]$ | f(x) | L[f(x)] |
|---|---|---|---|---|
| **12** | $[x]^{n}_{\Delta x}$<br>or<br>$\prod_{m=1}^{n}(x-[m-1]\Delta x)$<br>n = 1,2,3,… | $\frac{n!}{s^{n+1}}$<br>The related Z Transform is:<br>$\frac{n!T^{n}z}{(z-1)^{n+1}}$ | $x^n$ | $\frac{n!}{s^{n+1}}$ |
| **13** | $x^2$ | $\frac{2}{s^3} + \frac{\Delta x}{s^2}$ | $x^2$ | $\frac{2}{s^3}$ |
| **14** | $x^3$ | $\frac{6}{s^4} + \frac{6\Delta x}{s^3} + \frac{\Delta x^2}{s^2}$ | $x^3$ | $\frac{6}{s^4}$ |
| **15** | $\sin\frac{ax}{\Delta x}$ | $\frac{\frac{\sin a}{\Delta x}}{(s + \frac{1-\cos a}{\Delta x})^2 + (\frac{\sin a}{\Delta x})^2}$ | sinax | $\frac{a}{s^2+a^2}$ |
| **16** | sinax | $\frac{\frac{\sin a\Delta x}{\Delta x}}{(s + \frac{1-\cos a\Delta x}{\Delta x})^2 + (\frac{\sin a\Delta x}{\Delta x})^2}$ | sinax | $\frac{a}{s^2+a^2}$ |
| **17** | $\cos\frac{ax}{\Delta x}$ | $\frac{s + \frac{1-\cos a}{\Delta x}}{(s + \frac{1-\cos a}{\Delta x})^2 + (\frac{\sin a}{\Delta x})^2}$ | cosax | $\frac{s}{s^2+a^2}$ |
| **18** | cosax | $\frac{s + \frac{1-\cos a\Delta x}{\Delta x}}{(s + \frac{1-\cos a\Delta x}{\Delta x})^2 + (\frac{\sin a\Delta x}{\Delta x})^2}$ | cosax | $\frac{s}{s^2+a^2}$ |
| **19** | sinhax | $\frac{\frac{\sinh a\Delta x}{\Delta x}}{(s + \frac{1-\cosh a\Delta x}{\Delta x})^2 - (\frac{\sinh a\Delta x}{\Delta x})^2}$ | sinhax | $\frac{a}{s^2-a^2}$ |
| **20** | coshax | $\frac{s + \frac{1-\cosh a\Delta x}{\Delta x}}{(s + \frac{1-\cosh a\Delta x}{\Delta x})^2 - (\frac{\sinh a\Delta x}{\Delta x})^2}$ | coshax | $\frac{s}{s^2-a^2}$ |

| # | $f_{\Delta x}(x)$ | $K_{\Delta x}[f_{\Delta x}(x)]$ | f(x) | L[f(x)] |
|---|---|---|---|---|
| **21** | xsinax | $\dfrac{(\frac{\sin a\Delta x}{\Delta x})(2s+\Delta t s^2)}{[(s + \frac{1-\cos a\Delta x}{\Delta x})^2 + (\frac{\sin a\Delta x}{\Delta x})^2]^2}$ | xsinax | $\dfrac{2as}{[s^2 + a^2]^2}$ |
| **22** | xcosax | $\dfrac{A - B}{[(s + \frac{1-\cos a\Delta x}{\Delta x})^2 + (\frac{\sin a\Delta x}{\Delta x})^2]^2}$ where $A = (s + \frac{1-\cos a\Delta x}{\Delta x})^2 - (\frac{\sin a\Delta x}{\Delta x})^2$ $B = \frac{1-\cos a\Delta x}{\Delta x}[\Delta x s^2 + 4s + 2(\frac{1-\cos a\Delta x}{\Delta x})]$ | xcosax | $\dfrac{s^2 - a^2}{[s^2 + a^2]^2}$ |
| **23** | xsinhax | $\dfrac{(\frac{\sinh a\Delta x}{\Delta x})(2s+\Delta t s^2)}{[(s + \frac{1-\cosh a\Delta x}{\Delta x})^2 - (\frac{\sinh a\Delta x}{\Delta x})^2]^2}$ | xsinhax | $\dfrac{2as}{[s^2 - a^2]^2}$ |
| **24** | xcoshax | $\dfrac{A - B}{[(s + \frac{1-\cosh a\Delta x}{\Delta x})^2 - (\frac{\sinh a\Delta x}{\Delta x})^2]^2}$ where $A = (s + \frac{1-\cosh a\Delta x}{\Delta x})^2 + (\frac{\sinh a\Delta x}{\Delta x})^2$ $B = \frac{1-\cosh a\Delta x}{\Delta x}[\Delta x s^2 + 4s + 2(\frac{1-\cosh a\Delta x}{\Delta x})]$ | xcoshax | $\dfrac{s^2 + a^2}{[s^2 - a^2]^2}$ |
| **25** | $e^{ax}$ sinbx | $\dfrac{\frac{e^{a\Delta x}\sin b\Delta x}{\Delta x}}{(s - \frac{e^{a\Delta x}\cos b\Delta x - 1}{\Delta x})^2 + (\frac{e^{a\Delta x}\sin b\Delta x}{\Delta x})^2}$ | $e^{ax}$ sinbx | $\dfrac{b}{(s-a)^2 + b^2}$ |
| **26** | $e^{ax}$ cosbx | $\dfrac{s - \frac{e^{a\Delta x}\cos b\Delta x - 1}{\Delta x}}{(s - \frac{e^{a\Delta x}\cos b\Delta x - 1}{\Delta x})^2 + (\frac{e^{a\Delta x}\sin b\Delta x}{\Delta x})^2}$ | $e^{ax}$ cosbx | $\dfrac{s-a}{(s-a)^2 + b^2}$ |
| **27** | $b^x$ | $\dfrac{1}{s - \frac{b^{\Delta x} - 1}{\Delta x}}$ | $b^x$ | $\dfrac{1}{s - \ln b}$ |
| **28** | $e^{ax}$ | $\dfrac{1}{s - \frac{e^{a\Delta x} - 1}{\Delta x}}$ | $e^{ax}$ | $\dfrac{1}{s - a}$ |
| **29** | $xe^{ax}$ | $\dfrac{e^{a\Delta x}}{(s - \frac{e^{a\Delta x} - 1}{\Delta x})^2}$ | $xe^{ax}$ | $\dfrac{1}{(s - a)^2}$ |
| **30** | $x^2e^{ax}$ | $\dfrac{e^{a\Delta x}(s\Delta x + e^{a\Delta x} + 1)}{(s - \frac{e^{a\Delta x} - 1}{\Delta x})^3}$ | $x^2e^{ax}$ | $\dfrac{2}{(s - a)^3}$ |

| # | $f_{\Delta x}(x)$ | $K_{\Delta x}[f_{\Delta x}(x)]$ | f(x) | L[f(x)] |
|---|---|---|---|---|
| 31 | $[x]^{n}_{\Delta x}\, e^{ax}$ <br> or <br> $[\prod_{m=1}^{n}(x-[m-1]\Delta x)]e^{ax}$ <br> n = 1,2,3, … | $\dfrac{e^{na\Delta x}\, n!}{(s - \frac{e^{a\Delta x}-1}{\Delta x})^{n+1}}$ | $x^n e^{ax}$ | $\dfrac{n!}{(s-a)^{n+1}}$ |
| 32 | $e^{\frac{ax}{\Delta x}}$ | $\dfrac{1}{(s - \frac{e^{a}-1}{\Delta x})}$ | $e^{ax}$ | $\dfrac{1}{s-a}$ |
| 33 | $xe^{\frac{ax}{\Delta x}}$ | $\dfrac{e^{a}}{(s - \frac{e^{a}-1}{\Delta x})^2}$ | $xe^{ax}$ | $\dfrac{1}{(s-a)^2}$ |
| 34 | $\frac{1}{\Delta x}[U(x-n\Delta x)-U(x-n\Delta x-\Delta x)]$ <br> n = 0,1,2,3,… <br> Unit Area Pulse <br> Comments – Interval = Δx <br> Amplitude = $\frac{1}{\Delta x}$ <br> Δx times a Unit Area Pulse = a Unit Amplitude Pulse | $(1+s\Delta x)^{-n-1}$ <br> or <br> $(1+s\Delta x)^{-(\frac{x+\Delta x}{\Delta x})}$ <br> or <br> $\frac{1}{s\Delta x}[(1+s\Delta x)^{-n}-(1+s\Delta x)^{-n-1}]$ <br> $x = n\Delta x$ | $\frac{1}{\Delta x}[U(x-n\Delta x)-U(x-n\Delta x-\Delta x)]$ | $\dfrac{e^{-sn\Delta x}}{\Delta x}\left(\dfrac{1-e^{-s\Delta x}}{s}\right)$ |
| 35 | $[U(x-n\Delta x)-U(x-n\Delta x-\Delta x)]$ <br> n = 0,1,2,3,… <br> Unit Amplitude Pulse <br> Comment – Interval = Δx <br> Ampitude = 1 | $(1+s\Delta x)^{-n-1}\Delta x$ <br> or <br> $(1+s\Delta x)^{-(\frac{x+\Delta x}{\Delta x})}\Delta x$ <br> or <br> $\frac{1}{s}[(1+s\Delta x)^{-n}-(1+s\Delta x)^{-n-1}]$ <br> $x = n\Delta x$ | $[U(x-n\Delta x)-U(x-n\Delta x-\Delta x)]$ | $e^{-sn\Delta x}\left(\dfrac{1-e^{-s\Delta x}}{s}\right)$ |

| # | $f_{\Delta x}(x)$ | $K_{\Delta x}[f_{\Delta x}(x)]$ | f(x) | L[f(x)] |
|---|---|---|---|---|
| **36** | $[x]_{\Delta x}^{n}\, e_{\Delta x}(a,x)$<br>or<br>$e_{\Delta x}(a,x)\ \prod_{m=1}^{n}(x-[m-1]\Delta x)$<br>n=1,2,3,…<br>or<br>$e^{\alpha x}\ \prod_{m=1}^{n}(x-[m-1]\Delta x)$<br>n=1,2,3,…<br>$\alpha = \frac{1}{\Delta x}\ln(1+a\Delta x)$ | $\frac{(1+a\Delta x)^n\, n!}{(s-a)^{n+1}}$<br><br>The related Z Transform is:<br>$\frac{n!(1+aT)^n T^n z}{[z-(1+aT)]^{n+1}}$ | $x^n e^{ax}$ | $\frac{n!}{(s-a)^{n+1}}$ |
| **37** | $xe_{\Delta x}(a,x)$ | $\frac{1+a\Delta x}{(s-a)^2}$<br><br>The related Z Transform is:<br>$\frac{(1+aT)Tz}{z^2-2(1+aT)z+(1+aT)^2}$ | $xe^{ax}$ | $\frac{1}{(s-a)^2}$ |
| **38** | $x^2 e_{\Delta x}(a,x)$ | $\frac{(1+a\Delta x)(s\Delta x+a\Delta x+2)}{(s-a)^3}$ | $x^2 e^{ax}$ | $\frac{2}{(s-a)^3}$ |
| **39** | $x \sin_{\Delta x}(b,x)$ | $\frac{2bs + b\Delta x(s^2-b^2)}{(s^2+b^2)^2}$<br><br>The related Z Transform is:<br>$\frac{T^2bz[z^2 - (1+b^2T^2)]}{[z^2-2z+(1+b^2T^2)]^2}$ | x sinbx | $\frac{2bs}{(s^2+b^2)^2}$ |
| **40** | $x \cos_{\Delta x}(b,x)$ | $\frac{(s^2-b^2)-2\Delta x b^2 s}{(s^2+b^2)^2}$<br><br>The related Z Transform is:<br>$\frac{Tz[z^2-2(1+b^2T^2)z+(1+b^2T^2)]}{[z^2-2z+(1+b^2T^2)]^2}$ | x cosbx | $\frac{s^2-b^2}{(s^2+b^2)^2}$ |

| # | $f_{\Delta x}(x)$ | $K_{\Delta x}[f_{\Delta x}(x)]$ | f(x) | L[f(x)] |
|---|---|---|---|---|
| **41** | $x\,\sinh_{\Delta x}(b,x)$ | $\frac{2bs + b\Delta x(s^2+b^2)}{(s^2-b^2)^2}$<br>The related Z Transform is:<br>$\frac{T^2bz[z^2 - (1-b^2T^2)]}{[z^2-2z+(1-b^2T^2)]^2}$ | x sinhbx | $\frac{2bs}{(s^2-b^2)^2}$ |
| **42** | $x\,\cosh_{\Delta x}(b,x)$ | $\frac{(s^2+b^2)+2\Delta x b^2 s}{(s^2-b^2)^2}$<br>The related Z Transform is:<br>$\frac{Tz[z^2-2(1-b^2T^2)z+(1-b^2T^2)]}{[z^2-2z+(1-b^2T^2)]^2}$ | x coshbx | $\frac{s^2+b^2}{(s^2-b^2)^2}$ |
| **43** | $[x]^n_{\Delta x}\, e_{\Delta x}(a, x-n\Delta x)$<br>or<br>$\prod_{m=1}^{n}(x-[m-1]\Delta x)e_{\Delta x}(a, x-n\Delta x)$<br>n = 1,2,3, … | $\frac{n!}{(s-a)^{n+1}}$ | $x^n e^{ax}$ | $\frac{n!}{(s-a)^{n+1}}$ |
| **44** | $[x]^n_{\Delta x}\, \sin_{\Delta x}(b, x-n\Delta x)$<br>or<br>$\prod_{m=1}^{n}(x-[m-1]\Delta x)\sin_{\Delta x}(b, x-n\Delta x)$<br>n = 1,2,3, … | $\frac{n!}{2j}[\frac{(s+jb)^{n+1} - (s-jb)^{n+1}}{(s^2+b^2)^{n+1}}]$ | $x^n$ sinbx | $\frac{n!}{2j}[\frac{(s+jb)^{n+1} - (s-jb)^{n+1}}{(s^2+b^2)^{n+1}}]$ |
| **45** | $[x]^n_{\Delta x}\, \cos_{\Delta x}(b,x-n\Delta x)$<br>or<br>$\prod_{m=1}^{n}(x-[m-1]\Delta x)\cos_{\Delta x}(b, x-n\Delta x)$<br>n = 1,2,3, … | $\frac{n!}{2}[\frac{(s+jb)^{n+1}+(s-jb)^{n+1}}{(s^2+b^2)^{n+1}}]$ | $x^n$ cosbx | $\frac{n!}{2}[\frac{(s+jb)^{n+1}+(s-jb)^{n+1}}{(s^2+b^2)^{n+1}}]$ |
| **46** | $U(x-n\Delta x)e_{\Delta x}(a,x-n\Delta x)$<br>n = 0,1, 2, 3, … | $\frac{(1+s\Delta x)^{-n}}{s-a}$ | $U(x-\Delta x)e^{a(x-\Delta x)}$ | $\frac{e^{-s\Delta x}}{s-a}$ |

| # | $f_{\Delta x}(x)$ | $K_{\Delta x}[f_{\Delta x}(x)]$ | $f(x)$ | $L[f(x)]$ |
|---|---|---|---|---|
| **47** | $\sum_{n=0}^{\infty} \frac{a_n}{n!}[x]^n_{\Delta x}$ <br> or <br> $a_o + \sum_{n=1}^{\infty} \frac{a_n}{n!} \prod_{m=0}^{n-1}(x-m\Delta x)$ <br> $a_n = D^n_{\Delta x} y(x) \|_{x=0}$ | $\sum_{n=0}^{\infty} \frac{a_n}{s^{n+1}}$ | $a_o + \sum_{n=1}^{\infty} \frac{a_n x^n}{n!}$ <br> $a_n = \frac{d^n}{dx^n} y(x) \|_{x=0}$ | $\sum_{n=0}^{\infty} \frac{a_n}{s^{n+1}}$ |
| **48** | $\frac{\sin_{\Delta x}(b,x) + bx\cos_{\Delta x}(b,x)}{b}$ | $\frac{2s^2-b^2s\Delta x}{(s^2+b^2)^2}$ | $\frac{\sin bx + bx\cos ax}{b}$ | $\frac{2s^2}{(s^2+b^2)^2}$ |
| **49** | $\frac{\sinh_{\Delta x}(b,x) + bx\cosh_{\Delta x}(b,x)}{b}$ | $\frac{2s^2+b^2s\Delta x}{(s^2-b^2)^2}$ | $\frac{\sinh bx + bx\cosh ax}{b}$ | $\frac{2s^2}{(s^2-b^2)^2}$ |
| **50** | $\frac{\sin_{\Delta x}(b,x) - bx\cos_{\Delta x}(b,x)}{b^3}$ | $\frac{2+s\Delta x}{(s^2+b^2)^2}$ | $\frac{\sin bx - bx\cos ax}{b^3}$ | $\frac{2}{(s^2+b^2)^2}$ |
| **51** | $\frac{-\sinh_{\Delta}x(b,x) + bx\cosh_{\Delta x}(b,x)}{b^3}$ | $\frac{2+s\Delta x}{(s^2-b^2)^2}$ | $\frac{-\sinh bx + bx\cosh bx}{b^3}$ | $\frac{2}{(s^2-b^2)^2}$ |
| **52** | $\cos_{\Delta x}(b,x) - \frac{bx}{2}\sin_{\Delta x}(b,x)$ | $\frac{s^3 + b^2\Delta x(s^2-b^2)}{(s^2+b^2)^2}$ | $\cos bx - \frac{bx}{2}\sin bx$ | $\frac{s^3}{(s^2+b^2)^2}$ |
| **53** | $\cosh_{\Delta x}(b,x) + \frac{bx}{2}\sinh_{\Delta x}(b,x)$ | $\frac{s^3 - b^2\Delta x(s^2+b^2)}{(s^2-b^2)^2}$ | $\cosh bx + \frac{bx}{2}\sinh bx$ | $\frac{s^3}{(s^2-b^2)^2}$ |
| **54** | $\frac{e_{\Delta x}(b,x) - e_{\Delta x}(a,x)}{b-a}$ | $\frac{1}{(s-b)(s-a)}$ | $\frac{e^{bx} - e^{ax}}{b-a}$ | $\frac{1}{(s-b)(s-a)}$ |
| **55** | $\frac{be_{\Delta x}(b,x) - ae_{\Delta x}(a,x)}{b-a}$ | $\frac{s}{(s-b)(s-a)}$ | $\frac{be^{bx} - ae^{ax}}{b-a}$ | $\frac{s}{(s-b)(s-a)}$ |
| **56** | Square Wave <br> 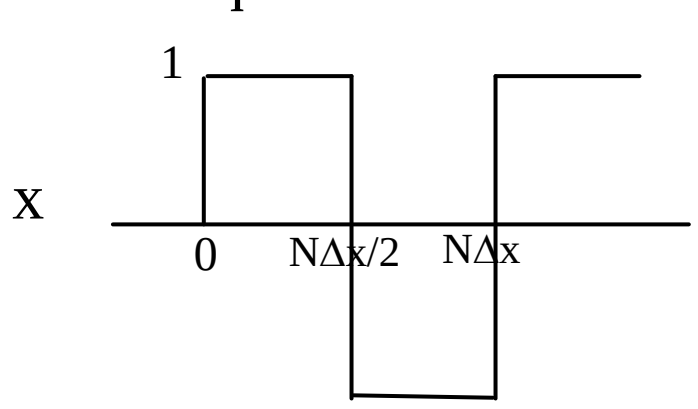  <br> SW(x) | Period = NΔx <br> SW(x) = SW(x+NΔx) <br> $\frac{1}{s}\frac{1-(1+s\Delta x)^{-\frac{N}{2}}}{1+(1+s\Delta x)^{-\frac{N}{2}}}$ <br> or <br> $\frac{1}{s}\tanh_{\Delta x}(s,\frac{N\Delta x}{4})$ | Period = a <br> SW(x) | $\frac{1}{s}\tanh(\frac{as}{4})$ |

| # | $f_{\Delta x}(x)$ | $K_{\Delta x}[f_{\Delta x}(x)]$ | f(x) | L[f(x)] |
|---|---|---|---|---|
| **57** | Triangular Waveform 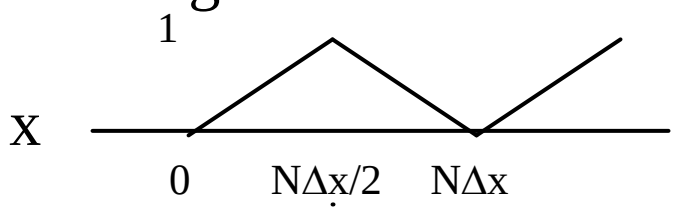  TW(x) | Period = NΔx<br>TW(x) = TW(x+NΔx)<br>$\frac{2}{N\Delta x s^2}\frac{1-(1+s\Delta x)^{-\frac{N}{2}}}{1+(1+s\Delta x)^{-\frac{N}{2}}}$<br>or<br>$\frac{2}{N\Delta x s^2}\tanh_{\Delta x}(s,\frac{N\Delta x}{4})$ | Period = a<br>TW(x) | $\frac{2}{as^2}\tanh(\frac{as}{4})$ |
| **58** | Saw Tooth Waveform 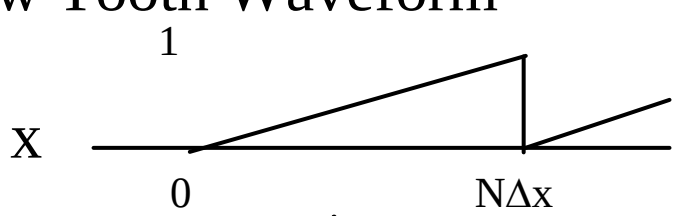  ST(x) | Period = NΔx<br>$\frac{1}{N\Delta x s^2} - \frac{(1+s\Delta x)^{-N}}{s[1-(1+s\Delta x)^{-N}]}$ | Period = a<br>ST(x) | $\frac{1}{as^2} - \frac{e^{-as}}{s(1-e^{-as})}$ |
| **59** | $\frac{1}{(\frac{x}{\Delta x})!}$<br>x = 0,Δx,2Δx,3Δx,… | $\Delta x(1+s\Delta x)^{-1}e^{(1+s\Delta x)^{-1}}$ | | |
| **60** | $\frac{1}{x}U(x-\Delta x)$<br>x = Δx,2Δx,3Δx,… | $(1+s\Delta x)^{-1}\ln(\frac{1+s\Delta x}{s\Delta x})$ | | |
| **61** | Impulse Function<br>(Impulse at x = 0)<br>δ(x) | $\Delta x(1+s\Delta x)^{-1}$<br>or<br>$\frac{1}{s+\frac{1}{\Delta x}}$ | δ(0) | 1 |
| **62** | Impulse Function<br>(Impulse at x = 0)<br>δ(x) | $\frac{K_{\Delta x}[f(x)]}{Z[f(x)]}$<br>where<br>z = 1+sΔx<br>Δx = T | δ(0) | 1 |
| **63** | Impulse Function<br>(Impulse at x = nΔx)<br>n = 0,1,2,3…<br>δ(x) | $\Delta x(1+s\Delta x)^{-n-1}$ | δ(nΔx) | $e^{-n\Delta xs}$ |

| # | $f_{\Delta x}(x)$ | $f(s) = K_{\Delta x}[f_{\Delta x}(x)]$ | Value Equations |
|---|---|---|---|
| **64** | $K_1 e_{\Delta x}(a,x) + K_2 e_{\Delta x}(b,x)$<br>or<br>$K_1(1+a\Delta x)^{\frac{x}{\Delta x}} + K_2(1+b\Delta x)^{\frac{x}{\Delta x}}]$<br><br>where<br>$x = 0, \Delta x, 2\Delta x, 3\Delta x, 4\Delta x, \ldots$<br>$\Delta x$ = x increment | The roots of the denominator are real<br><br>$\frac{K_1}{s-a} + \frac{K_2}{s-b} = \frac{Cs+D}{s^2+As+B}$<br><br>where<br>$s = a, b$, the roots of $s^2 + As + B$ | a,b,A,B,C,D = real constants<br>$K_1$, $K_2$ = real constants<br>$a = -\frac{A}{2} - \sqrt{(\frac{A}{2})^2 - B}$<br>$b = -\frac{A}{2} + \sqrt{(\frac{A}{2})^2 - B}$<br>$K_1 = \frac{D + aC}{a-b}$<br>$K_2 = -\frac{D + bC}{a-b}$ |
| **65** | $e_{\Delta x}(a,x)[K_1\cos_{\Delta x}(\frac{b}{1+a\Delta x},x) + K_2\sin_{\Delta x}(\frac{b}{1+a\Delta x},x)]$<br>or<br>$(1+a\Delta x)^{\frac{x}{\Delta x}}[K_1\cos_{\Delta x}(w,x) + K_2\sin_{\Delta x}(w,x)]$<br>or<br>$[(1+a\Delta x)^2+b^2\Delta x^2]^{\frac{x}{2\Delta x}}[K_1\cos\beta\frac{x}{\Delta x} + K_2\sin\beta\frac{x}{\Delta x}]$<br><br>where<br>$x = 0, \Delta x, 2\Delta x, 3\Delta x, 4\Delta x, \ldots$<br>$\Delta x$ = x increment | The roots of the denominator are complex<br><br>$\frac{K_1(s-a) + K_2 b}{(s-a)^2+b^2} =$<br>$\frac{K_1(s-a) + K_2 w(1+a\Delta x)}{(s-a)^2 + (w[1+a\Delta x])^2}$<br>$= \frac{Cs+D}{s^2+As+B}$<br><br>where<br>$s = a \pm jb$, the roots of $s^2 + As + B$ | a,b,A,B,C,D = real constants<br>$K_1$, $K_2$ = real constants<br>$a = -\frac{A}{2}$<br>$b = \sqrt{B - (\frac{A}{2})^2} = w(1+a\Delta x)$<br>$w = \frac{b}{1+a\Delta x}$<br>$K_1 = C$<br>$K_2 = \frac{D - \frac{CA}{2}}{\sqrt{B - (\frac{A}{2})^2}} = \frac{D + Ca}{b}$<br>$\beta = \begin{cases} \tan^{-1}\frac{b\Delta x}{1+a\Delta x} & \text{for } 1+a\Delta x \ge 0 \\ \pi + \tan^{-1}\frac{b\Delta x}{1+a\Delta x} & \text{for } 1+a\Delta x < 0 \end{cases}$ |

# TABLE 3a

## Conversion of Calculus Function Laplace Transforms to Equivalent Function

## $K_{\Delta x}$ Transforms

| # | $K_{\Delta x}$ Transform $\Delta x \to 0$<br>Laplace Transform<br>F(x) where $0 \le x < \infty$ | $K_{\Delta x}$ Transform<br>Generalized Laplace Transform<br>F(x) where $x = 0, \Delta x, 2\Delta x, 3\Delta x, \ldots$ | F(x)<br>Calculus<br>Function |
|---|---|---|---|
| **1** | $\frac{1}{s}$ | $\frac{1}{s}$ | 1 |
| **2** | $\frac{1}{s-a}$ | $\frac{1}{s - \frac{e^{a\Delta x}-1}{\Delta x}}$ | $e^{ax}$ |
| **3** | $\frac{b}{s^2+b^2}$ | $\frac{\frac{\sin b\Delta x}{\Delta x}}{(s + \frac{1-\cos b\Delta x}{\Delta x})^2 + (\frac{\sin b\Delta x}{\Delta x})^2}$ | sinbx |
| **4** | $\frac{b}{s^2-b^2}$ | $\frac{\frac{\sinh b\Delta x}{\Delta x}}{(s + \frac{1-\cosh b\Delta x}{\Delta x})^2 - (\frac{\sinh b\Delta x}{\Delta x})^2}$ | sinhbx |
| **5** | $\frac{s}{s^2+b^2}$ | $\frac{s + \frac{1-\cos b\Delta x}{\Delta x}}{(s + \frac{1-\cos b\Delta x}{\Delta x})^2 + (\frac{\sin b\Delta x}{\Delta x})^2}$ | cosbx |
| **6** | $\frac{s}{s^2-b^2}$ | $\frac{s + \frac{1-\cosh b\Delta x}{\Delta x}}{(s + \frac{1-\cosh b\Delta x}{\Delta x})^2 - (\frac{\sinh b\Delta x}{\Delta x})^2}$ | coshbx |

| # | $K_{\Delta x}$ Transform $\Delta x \to 0$<br>**Laplace Transform**<br>**F(x) where $0 \le x < \infty$** | $K_{\Delta x}$ Transform<br>**Generalized Laplace Transform**<br>**F(x) where x = 0, $\Delta x$, $2\Delta x$, $3\Delta x$, …** | **F(x)**<br>**Calculus**<br>**Function** |
|---|---|---|---|
| **7** | $\frac{b}{(s-a)^2 + b^2}$ | $\frac{\frac{e^{a\Delta x}\sin b\Delta x}{\Delta x}}{(s - \frac{e^{a\Delta x}\cos b\Delta x - 1}{\Delta x})^2 + (\frac{e^{a\Delta x}\sin b\Delta x}{\Delta x})^2}$ | $e^{ax}\sin bx$ |
| **8** | $\frac{b}{(s-a)^2 - b^2}$ | $\frac{\frac{e^{a\Delta x}\sinh b\Delta x}{\Delta x}}{(s - \frac{e^{a\Delta x}\cosh b\Delta x - 1}{\Delta x})^2 - (\frac{e^{a\Delta x}\sinh b\Delta x}{\Delta x})^2}$ | $e^{ax}\sinh bx$ |
| **9** | $\frac{s-a}{(s-a)^2 + b^2}$ | $\frac{s - \frac{e^{a\Delta x}\cos b\Delta x - 1}{\Delta x}}{(s - \frac{e^{a\Delta x}\cos b\Delta x - 1}{\Delta x})^2 + (\frac{e^{a\Delta x}\sin b\Delta x}{\Delta x})^2}$ | $e^{ax}\cos bx$ |
| **10** | $\frac{s-a}{(s-a)^2 - b^2}$ | $\frac{s - \frac{e^{a\Delta x}\cosh b\Delta x - 1}{\Delta x}}{(s - \frac{e^{ax}\cosh b\Delta x - 1}{\Delta x})^2 - (\frac{e^{a\Delta x}\sinh b\Delta x}{\Delta x})^2}$ | $e^{ax}\cosh bx$ |
| **11** | $e^{-sn\Delta x}f(s)$ | $(1+s\Delta x)^{-n}f(s)$ | $U(x-n\Delta x)F(x-n\Delta x)$<br>n = 0,1,2,3,… |
| **12** | $\frac{1}{s^2}$ | $\frac{1}{s^2}$ | $x$ |
| **13** | $\frac{2}{s^3}$ | $\frac{2}{s^3} + \frac{\Delta x}{s^2}$ | $x^2$ |
| **14** | $\frac{6}{s^4}$ | $\frac{6}{s^4} + \frac{6\Delta x}{s^3} + \frac{\Delta x^2}{s^2}$ | $x^3$ |
| **15** | $\frac{1}{(s-a)^2}$ | $\frac{e^{a\Delta x}}{(s - \frac{e^{a\Delta x} - 1}{\Delta x})^2}$ | $xe^{ax}$ |
| **16** | $\frac{2}{(s-a)^3}$ | $\frac{e^{a\Delta x}(s\Delta x + e^{a\Delta x} + 1)}{(s - \frac{e^{a\Delta x} - 1}{\Delta x})^3}$ | $x^2e^{ax}$ |

| # | $K_{\Delta x}$ Transform $\Delta x \to 0$ Laplace Transform F(x) where $0 \le x < \infty$ | $K_{\Delta x}$ Transform Generalized Laplace Transform F(x) where $x = 0, \Delta x, 2\Delta x, 3\Delta x, \ldots$ | F(x) Calculus Function |
|---|---|---|---|
| **17** | $\frac{2bs}{[\, s^2 + b^2 \,]^2}$ | $\frac{\frac{\sin b\Delta x}{\Delta x}(2s+\Delta x s^2)}{[(s + \frac{1-\cos b\Delta x}{\Delta x})^2 + (\frac{\sin b\Delta x}{\Delta x})^2]^2}$ | xsinbx |
| **18** | $\frac{2bs}{[\, s^2 - b^2 \,]^2}$ | $\frac{\frac{\sinh b\Delta x}{\Delta x}(2s+\Delta x s^2)}{[(s + \frac{1-\cosh b\Delta x}{\Delta x})^2 - (\frac{\sinh b\Delta x}{\Delta x})^2]^2}$ | xsinhbx |
| **19** | $\frac{s^2 - b^2}{[\, s^2 + b^2 \,]^2}$ | $\frac{[(s + \frac{1-\cos b\Delta x}{\Delta x})^2 - (\frac{\sin b\Delta x}{\Delta x})^2] - (\frac{1-\cos b\Delta x}{\Delta x})[\Delta x s^2 + 4s + 2(\frac{1-\cos b\Delta x}{\Delta x})]}{[(s + \frac{1-\cos b\Delta x}{\Delta x})^2 + (\frac{\sin b\Delta x}{\Delta x})^2]^2}$ | xcosbx |
| **20** | $\frac{s^2 + b^2}{[\, s^2 - b^2 \,]^2}$ | $\frac{[(s + \frac{1-\cosh b\Delta x}{\Delta x})^2 + (\frac{\sinh b\Delta x}{\Delta x})^2] - (\frac{1-\cosh b\Delta x}{\Delta x})[\Delta x s^2 + 4s + 2(\frac{1-\cosh b\Delta x}{\Delta x})]}{[(s + \frac{1-\cosh b\Delta x}{\Delta x})^2 - (\frac{\sinh b\Delta x}{\Delta x})^2]^2}$ | xcoshbx |

# TABLE 3b

## Conversions from $K_{\Delta x}$ Transforms to Calculus Functions

| # | $K_{\Delta x}$ Transform Generalized Laplace Transform of F(x) | F(x) Calculus Functions<br><br>F(t) where t = 0, Δx, 2Δx… | Coefficient Definitions |
|---|---|---|---|
| **1** | $\frac{1}{s}$ | 1 | |
| **2** | $\frac{1}{s-a}$ | $e^{\alpha x}$<br><br>The Equivalent non-Calculus Function is:<br>$e_{\Delta x}(\frac{e^{\alpha\Delta x}-1}{\Delta x}, x)$ | Denominator polynomial a real root<br>$\alpha = \frac{\ln(1+a\Delta x)}{\Delta x}$ |
| **3** | $\frac{b}{(s-a)^2+b^2}$ | $e^{\alpha x}\sin\beta x$<br><br>The Equivalent non-Calculus Function is:<br>$(1+a\Delta x)^{\frac{x}{\Delta x}}\sin_{\Delta x}(\frac{b\Delta x}{1+a\Delta x}, x)$ | Denominator polynomial complex conjugate roots = a ± jb<br>$\alpha = \frac{1}{2\Delta x}\ln[(1+a\Delta x)^2 + (b\Delta x)^2]$<br><br>$\beta = \begin{cases} \frac{1}{\Delta x}\tan^{-1}\left\lvert\frac{b\Delta x}{1+a\Delta x}\right\rvert & \text{for } 1+a\Delta x > 0 \quad b\Delta x \ge 0 \\ -\frac{1}{\Delta x}\tan^{-1}\left\lvert\frac{b\Delta x}{1+a\Delta x}\right\rvert & \text{for } 1+a\Delta x > 0 \quad b\Delta x < 0 \\ \frac{1}{\Delta x}[\pi - \tan^{-1}\left\lvert\frac{b\Delta x}{1+a\Delta x}\right\rvert] & \text{for } 1+a\Delta x < 0 \quad b\Delta x \ge 0 \\ \frac{1}{\Delta x}[-\pi + \tan^{-1}\left\lvert\frac{b\Delta x}{1+a\Delta x}\right\rvert] & \text{for } 1+a\Delta x < 0 \quad b\Delta x < 0 \end{cases}$<br><br>$1+a\Delta x \neq 0$<br><br>$0 \le \tan^{-1}\left\lvert\frac{b\Delta x}{1+a\Delta x}\right\rvert < \frac{\pi}{2}$<br><br><u>Note</u> – For sinβx, $(1+a\Delta x)^2 + (b\Delta x)^2 = 1$ |

| # | $K_{\Delta x}$ Transform Generalized Laplace Transform of F(x) | F(x) Calculus Functions<br><br>F(t) where t = 0, Δx, 2Δx… | Coefficient Definitions |
|---|---|---|---|
| **4** | $\frac{s-a}{(s-a)^2+b^2}$ | $e^{\alpha x}\cos\beta x$<br><br>The Equivalent non-Calculus Function is:<br><br>$(1+a\Delta x)^{\frac{x}{\Delta x}}\cos_{\Delta x}(\frac{b\Delta x}{1+a\Delta x}, x)$ | Denominator polynomial complex conjugate roots = $a \pm jb$<br><br>$\alpha = \frac{1}{2\Delta x}\ln[(1+a\Delta x)^2+(b\Delta x)^2]$<br><br>$\beta = \begin{cases} \frac{1}{\Delta x}\tan^{-1}\left\lvert\frac{b\Delta x}{1+a\Delta x}\right\rvert & \text{for } 1+a\Delta x>0 \quad b\Delta x\ge 0 \\ -\frac{1}{\Delta x}\tan^{-1}\left\lvert\frac{b\Delta x}{1+a\Delta x}\right\rvert & \text{for } 1+a\Delta x>0 \quad b\Delta x< 0 \\ \frac{1}{\Delta x}[\pi-\tan^{-1}\left\lvert\frac{b\Delta x}{1+a\Delta x}\right\rvert] & \text{for } 1+a\Delta x<0 \quad b\Delta x\ge 0 \\ \frac{1}{\Delta x}[-\pi+\tan^{-1}\left\lvert\frac{b\Delta x}{1+a\Delta x}\right\rvert] & \text{for } 1+a\Delta x<0 \quad b\Delta x< 0 \end{cases}$<br><br>$1+a\Delta x \ne 0$<br><br>$0 \le \tan^{-1}\left\lvert\frac{b\Delta x}{1+a\Delta x}\right\rvert < \frac{\pi}{2}$<br><br><u>Note</u> – For $\cos\beta x$, $\alpha = 0$, $(1+a\Delta x)^2+(b\Delta x)^2 = 1$ |
| **5** | $(1+s\Delta x)^{-n}f(s)$ | $U(x-n\Delta x)F(x-\Delta x)$ | The Unit Step Function, $U(x-n\Delta x)$<br>1 for $x \ge n\Delta x$<br>0 for $x < n\Delta x$<br><br>n = 0,1,2,3,…<br>f(s) = $K_{\Delta x}$ Transform of F(x) |
| **6** | $\frac{n!}{s^{n+1}}$ | $\prod_{m=1}^{n}(x-[m-1]\Delta x)$ | n = 1,2,3,… |
| **7** | $\frac{(1+a\Delta t)^n n!}{(s-a)^{n+1}}$ | $e^{\alpha x}\prod_{m=1}^{n}(x-[m-1]\Delta x)$ | Denominator polynomial a multiple real root<br><br>$\alpha = \frac{\ln(1+a\Delta x)}{\Delta x}$<br><br>n = 1,2,3,… |

# TABLE 3c

## Modified $K_{\Delta x}$ Transforms

| # | Calculus Function Laplace Transforms | f(x) Calculus Functions Interval Calculus Functions $x = 0, \Delta x, 2\Delta x \ldots$ | $K_{\Delta x}$ Transform of f(x) F(s) $x = 0, \Delta x, 2\Delta x, \ldots$ | Modified $K_{\Delta x}$ Transform of f(x) F(s,m) $0 \le m \le 1, \quad x = m\Delta x, \Delta x+m\Delta x, 2\Delta x$ |
|---|---|---|---|---|
| 1 | $\frac{1}{s}$ | 1 | $\frac{1}{s}$ | $\frac{1}{s}$ |
| 2 | $\frac{1}{s^2}$ | x | $\frac{1}{s^2}$ | $\frac{1}{s^2} + \frac{m\Delta x}{s}$ |
| 3 | | $x(x-\Delta x)$ | $\frac{2}{s^3}$ | $\frac{2}{s^3} + \frac{2m\Delta x}{s^2} + \frac{m(m-1)\Delta x^2}{s}$ |
| 4 | $\frac{2}{s^3}$ | $x^2$ | $\frac{2}{s^3} + \frac{\Delta x}{s^2}$ | $\frac{2}{s^3} + \frac{(2m+1)\Delta x}{s^2} + \frac{m^2\Delta x^2}{s}$ |
| 5 | | $e_{\Delta x}(a,x)$ | $\frac{1}{s-a}$ | $\frac{(1+a\Delta x)^m}{s-a}$ |
| 6 | $\frac{1}{s-a}$ | $e^{ax}$ or $e_{\Delta x}(\frac{e^{a\Delta x}-1}{\Delta x},x)$ | $\frac{1}{s - \frac{e^{a\Delta x}-1}{\Delta x}}$ | $\frac{e^{am\Delta x}}{s - \frac{e^{a\Delta x}-1}{\Delta x}}$ |
| 7 | | $\sin_{\Delta x}(b,x)$ | $\frac{b}{s^2+b^2}$ | $\frac{\sin_{\Delta x}(b,m\Delta x)s + b\cos_{\Delta x}(b,m\Delta x)}{s^2+b^2}$ |
| 8 | $\frac{b}{s^2+b^2}$ | sinbx or $e_{\Delta x}(\frac{\cos b\Delta x - 1}{\Delta x}, x)\sin_{\Delta x}(\frac{\tan b\Delta x}{\Delta x}, x)$ | $\frac{\frac{\sin b\Delta x}{\Delta x}}{(s + \frac{1-\cos b\Delta x}{\Delta x})^2 + (\frac{\sin b\Delta x}{\Delta x})^2}$ | $\frac{\cos bm\Delta x(\frac{\sin b\Delta x}{\Delta x}) + \sin bm\Delta x(s + \frac{1-\cos b\Delta x}{\Delta x})}{s^2 + 2(\frac{1-\cos b\Delta x}{\Delta x})s + 2(\frac{1-\cos b\Delta x}{\Delta x^2})}$ |

| # | Calculus Function Laplace Transforms | f(x) Calculus Functions **Interval Calculus Functions** $x = 0, \Delta x, 2\Delta x \ldots$ | $K_{\Delta x}$ Transform of f(x) F(s) $x = 0, \Delta x, 2\Delta x, \ldots$ | Modified $K_{\Delta x}$ Transform of f(x) F(s,m) $0 \le m \le 1, \quad x = m\Delta x, \Delta x + m\Delta x, 2\Delta x$ |
|---|---|---|---|---|
| **9** | | $\cos_{\Delta x}(b,x)$ | $\frac{s}{s^2+b^2}$ | $\frac{\cos_{\Delta x}(b,m\Delta x)s - b\sin_{\Delta x}(b,m\Delta x)}{s^2 + b^2}$ |
| **10** | $\frac{s}{s^2+b^2}$ | $\cos bx$<br>or<br>$e_{\Delta x}(\frac{\cos b\Delta x - 1}{\Delta x}, x) \cos_{\Delta x}(\frac{\tan b\Delta x}{\Delta x}, x)$ | $\frac{s + \frac{1-\cos b\Delta x}{\Delta x}}{(s + \frac{1-\cos b\Delta x}{\Delta x})^2 + (\frac{\sin b\Delta x}{\Delta x})^2}$ | $\frac{\cos bm\Delta x(s + \frac{1-\cos b\Delta x}{\Delta x}) - \sin bm\Delta x(\frac{\sin b\Delta x}{\Delta x})}{s^2 + 2(\frac{1-\cos b\Delta x}{\Delta x})s + 2(\frac{1-\cos b\Delta x}{\Delta x^2})}$ |
| **11** | | $e_{\Delta x}(a,x)\sin_{\Delta x}(\frac{b}{1+a\Delta x},x)$<br>$a \neq -\frac{1}{\Delta x}$ | $\frac{b}{(s-a)^2+b^2}$ | $\frac{(1+a\Delta x)^m [b\cos_{\Delta x}(\frac{b}{1+a\Delta x},m\Delta x) + \sin_{\Delta x}(\frac{b}{1+a\Delta x},m\Delta x)(s-a)]}{(s-a)^2 + b^2}$ |
| **12** | $\frac{b}{(s-a)^2 + b^2}$ | $e^{ax}\sin bx$<br>or<br>$e_{\Delta x}(\frac{e^{a\Delta x}\cos b\Delta x - 1}{\Delta x}, x)\sin_{\Delta x}(\frac{\tan b\Delta x}{\Delta x}, x)$ | $\frac{\frac{e^{a\Delta x}\sin b\Delta x}{\Delta x}}{(s - \frac{e^{a\Delta x}\cos b\Delta x - 1}{\Delta x})^2 + (\frac{e^{a\Delta x}\sin b\Delta x}{\Delta x})^2}$ | $\frac{e^{am\Delta x}[\cos m\Delta x(\frac{e^{a\Delta x}\sin b\Delta x}{\Delta x}) + \sin bm\Delta x(s + \frac{1-e^{a\Delta x}\cos b\Delta x}{\Delta x})]}{s^2 + 2(\frac{1-e^{a\Delta x}\cos b\Delta x}{\Delta x})s + (\frac{1-2e^{a\Delta x}\cos b\Delta x + e^{2a\Delta x}}{\Delta x^2})}$ |
| **13** | | $e_{\Delta x}(a,x)\cos_{\Delta x}(\frac{b}{1+a\Delta x},x)$<br>$a \neq -\frac{1}{\Delta x}$ | $\frac{s-a}{(s-a)^2+b^2}$ | $\frac{(1+a\Delta x)^m [\cos_{\Delta x}(\frac{b}{1+a\Delta x},m\Delta x)(s-a) - b\sin_{\Delta x}(\frac{b}{1+a\Delta x},m\Delta x)]}{(s-a)^2 + b^2}$ |
| **14** | $\frac{s-a}{(s-a)^2 + b^2}$ | $e^{ax}\cos bx$<br>or<br>$e_{\Delta x}(\frac{e^{a\Delta x}\cos b\Delta x - 1}{\Delta x}, x)\cos_{\Delta x}(\frac{\tan b\Delta x}{\Delta x}, x)$ | $\frac{s - \frac{e^{a\Delta x}\cos b\Delta x - 1}{\Delta x}}{(s - \frac{e^{a\Delta x}\cos b\Delta x - 1}{\Delta x})^2 + (\frac{e^{a\Delta x}\sin b\Delta x}{\Delta x})^2}$ | $\frac{e^{am\Delta x}[\cos bm\Delta x(s + \frac{1-e^{a\Delta x}\cos b\Delta x}{\Delta x}) - \sin bm\Delta x(\frac{e^{a\Delta x}\sin b\Delta x}{\Delta x})]}{s^2 + 2(\frac{1-e^{a\Delta x}\cos b\Delta x}{\Delta x})s + (\frac{1-2e^{a\Delta x}\cos b\Delta x + e^{2a\Delta x}}{\Delta x^2})}$ |

1. $K\Delta t_M[y(t)] = K_{\Delta t}[y(t+m\Delta t]$ Definition of the Modified $K_{\Delta t}$ Transform

2. $K\Delta t_M[y(t)] = Tz^{-1}Z[y(t+mT)] \,|_{z = 1+s\Delta t}$ Modified $K_{\Delta t}$ Transform from the Z Transform
   $\Delta t = T$ sampling period

   or

   $K\Delta t_M[y(t)] = TZ[y(t+[m-1]T)] \,|_{z = 1+s\Delta t}$ Modified $K_{\Delta t}$ Transform from the Z Transform
   $\Delta t = T$ sampling period

3. $K\Delta t_M[y(t)] = TZ_M[y(t)] \,|_{z = 1+s\Delta t}$ Modified $K_{\Delta t}$ Transform from the Modified Z Transform
   $\Delta t = T$ sampling period

4. $Z_M[y(t)] = z^{-1}Z[y(t+mT)]$ Definition of the Modified Z Transform

   or

   $Z_M[y(t)] = Z[y(t+[m-1]T)]$ Definition of the Modified Z Transform

5. $Z_M[y(t)] = \frac{1}{T} K\Delta t_M[y(t)] \,|_{s = \frac{z-1}{T}}$ Modified Z Transform from the Modified $K_{\Delta t}$ Transform
   $T = \Delta t$ sampling period

   where
   $t = 0, \Delta t, 2\Delta t, 3\Delta t, \ldots$
   $0 \le m < 1$
   $\Delta t = T$

# TABLE 3d

## Z Transforms of Some Interval Calculus Functions

| # | F(x)<br>Interval Calculus Functions<br><br>x = 0, Δx, 2Δx… | $K_{\Delta x}$ Transform<br>of F(x)<br><br>x = 0, Δx, 2Δx, … | Z Transform<br>of F(x)<br><br>T = Δx, x = 0, T, 2T, … |
|---|---|---|---|
| **1** | $[x]^n_{\Delta x}$<br>or<br>$\prod_{m=1}^{n}(x-[m-1]\Delta x)$<br>n = 1,2,3,… | $\frac{n!}{s^{n+1}}$ | $\frac{n!T^n z}{(z-1)^{n+1}}$ |
| **2** | $e_{\Delta x}(a,x)$ | $\frac{1}{s-a}$ | $\frac{z}{z-(1+aT)}$ |
| **3** | $xe_{\Delta x}(a,x)$ | $\frac{1+a\Delta x}{(s-a)^2}$ | $\frac{(1+aT)Tz}{z^2-2(1+aT)z+(1+aT)^2}$ |
| **4** | $[x]^n_{\Delta x}\, e_{\Delta x}(a,x)$<br>or<br>$e_{\Delta x}(a,x)\prod_{m=1}^{n}(x-[m-1]\Delta x)$<br>n=1,2,3,… | $\frac{(1+a\Delta x)^n\, n!}{(s-a)^{n+1}}$ | $\frac{n!(1+aT)^n T^n z}{[z-(1+aT)]^{n+1}}$ |
| **5** | $\sin_{\Delta x}(b,x)$ | $\frac{b}{s^2+b^2}$ | $\frac{bTz}{z^2-2z+(1+b^2T^2)}$ |

| # | F(x)<br>Interval Calculus Functions<br><br>**x = 0, Δx, 2Δx…** | $K_{\Delta x}$ Transform of F(x)<br><br>**x = 0, Δx, 2Δx, …** | Z Transform of F(x)<br><br>**T = Δx, x = 0, T, 2T, …** |
|---|---|---|---|
| **6** | $x\sin_{\Delta x}(b,x)$ | $\frac{2bs + b\Delta x(s^2-b^2)}{(s^2+b^2)^2}$ | $\frac{T^2bz[z^2 - (1+b^2T^2)]}{[z^2-2z+(1+b^2T^2)]^2}$ |
| **7** | $\cos_{\Delta x}(b,x)$ | $\frac{s}{s^2+b^2}$ | $\frac{z(z-1)}{z^2-2z+(1+b^2T^2)}$ |
| **8** | $x\cos_{\Delta x}(b,x)$ | $\frac{(s^2-b^2)-2\Delta xb^2s}{(s^2+b^2)^2}$ | $\frac{Tz[z^2-2(1+b^2T^2)z+(1+b^2T^2)]}{[z^2-2z+(1+b^2T^2)]^2}$ |
| **9** | $e_{\Delta x}(a,x)\sin_{\Delta x}(\frac{b}{1+a\Delta x},x)$<br>$a \neq -\frac{1}{\Delta x}$ | $\frac{b}{(s-a)^2+b^2}$ | $\frac{bTz}{z^2-2[1+aT]z+[(1+aT)^2+b^2T^2]}$ |
| **10** | $e_{\Delta x}(a,x)\cos_{\Delta x}(\frac{b}{1+a\Delta x},x)$<br>$a \neq -\frac{1}{\Delta x}$ | $\frac{s-a}{(s-a)^2+b^2}$ | $\frac{z^2- (1+aT)z}{z^2-2[1+aT]z+[(1+aT)^2+b^2T^2]}$ |
| **11** | $\sinh_{\Delta x}(b,x)$ | $\frac{b}{s^2-b^2}$ | $\frac{bTz}{z^2-2z+(1-b^2T^2)}$ |
| **12** | $x \sinh_{\Delta x}(b,x)$ | $\frac{2bs + b\Delta x(s^2+b^2)}{(s^2-b^2)^2}$ | $\frac{T^2bz[z^2 - (1-b^2T^2)]}{[z^2-2z+(1-b^2T^2)]^2}$ |
| **13** | $\cosh_{\Delta x}(b,x)$ | $\frac{s}{s^2-b^2}$ | $\frac{z(z-1)}{z^2-2z+(1-b^2T^2)}$ |
| **14** | $x \cosh_{\Delta x}(b,x)$ | $\frac{(s^2+b^2)+2\Delta xb^2s}{(s^2-b^2)^2}$ | $\frac{Tz[z^2-2(1-b^2T^2)z+(1-b^2T^2)]}{[z^2-2z+(1-b^2T^2)]^2}$ |

| # | F(x)<br>Interval Calculus Functions<br><br>$x = 0, \Delta x, 2\Delta x\ldots$ | $K_{\Delta x}$ Transform<br>of F(x)<br><br>$x = 0, \Delta x, 2\Delta x, \ldots$ | Z Transform<br>of F(x)<br><br>$T = \Delta x,\quad x = 0, T, 2T, \ldots$ |
|---|---|---|---|
| **15** | $K_{\Delta x}[f(x-n\Delta x)U(x-n\Delta x)]$<br><br>$U(x-n\Delta x) = \begin{Bmatrix} 1 & x \geq n\Delta x \\ 0 & x < n\Delta x \end{Bmatrix}$<br><br>$n = 0,1,2,3,\ldots$<br><br>Unit Step Function | $(1+s\Delta x)^{-n}K_{\Delta t}[f(t)]$ | $z^{-n}Z[f(t)]$<br><br>$n = 0,1,2,3,\ldots$ |
| **16** | $U(x-n\Delta x)-U(x-[n+1]\Delta x))$<br>$n = 0,1,2,3,\ldots$<br><br>Unit Amplitude Pulse<br><br>Pulse Interval = $\Delta x$<br>Pulse Amplitude = 1<br>Pulse initiation at $x = n\Delta x$<br>Pulse ends at $x = (n+1)\Delta x$ | $(1+s\Delta x)^{-n-1}\Delta x$<br>or<br>$(1+s\Delta x)^{-(\frac{x+\Delta x}{\Delta x})}\Delta x$<br>or<br>$\frac{1}{s}[(1+s\Delta x)^{-n}-(1+s\Delta x)^{-n-1}]$<br>$x = n\Delta x$ | $Tz^{-n-1}$<br><br>$n = 0,1,2,3,\ldots$ |
| **17** | $f(x)$ | $\int_{0}^{\infty}{}_{\Delta x}\,(1+s\Delta x)^{-(\frac{x+\Delta x}{\Delta x})}f(x)\Delta x$ | $\frac{1}{T}\int_{0}^{\infty}{}_{T}\,z^{-\frac{x}{\Delta x}}f(x)\Delta x$ |
| **18** | $xf(x)$ | $-(1+s\Delta x)\frac{d}{ds}K_{\Delta x}[f(x)] - \Delta x K_{\Delta x}[f(x)]$ | $-Tz\frac{d}{dz}Z[f(x)]$<br><br>or<br><br>$-Tz^2\frac{d}{dz}[\frac{Z[f(x)]}{z}] - Tz[f(x)]$ |

| # | F(x) Interval Calculus Functions<br><br>$x = 0, \Delta x, 2\Delta x\ldots$ | $K_{\Delta x}$ Transform of F(x)<br><br>$x = 0, \Delta x, 2\Delta x, \ldots$ | Z Transform of F(x)<br><br>$T = \Delta x, \quad x = 0, T, 2T, \ldots$ |
|---|---|---|---|
| **19** | $D_{\Delta x}f(x)$ | $sK_{\Delta x}[f(x)] - f(0)$ | $\frac{z-1}{T} Z[f(x)] - \frac{z}{T} f(0)$<br><br>$D_T f(x) = \frac{f(x+T) - f(x)}{T}$ |
| **20** | $D^n{}_{\Delta x}f(x)$ | $s^n K_{\Delta x}[f(x)] - s^{n-1}D^0{}_{\Delta x}f(0) - s^{n-2}D^1{}_{\Delta x}f(0) - s^{n-3}D^2{}_{\Delta x}f(0) - \ldots - s^0 D^{n-1}{}_{\Delta x}f(0)$<br><br>or<br><br>$s^n K_{\Delta x}[f(x)] - \sum_{m=1}^{n} s^{n-m} D_{\Delta x}{}^{m-1} f(0)$<br><br>n = 1,2,3,… | $\left(\frac{z-1}{T}\right)^n Z[f(x)] - \sum_{m=1}^{n} \frac{z}{T} \left(\frac{z-1}{T}\right)^{n-m} D_T{}^{m-1} f(0)$<br><br>n = 1,2,3,… |
| **21** | ${}_{\Delta x}\int_0^x f(x)\Delta x$ | $\frac{1}{s} K_{\Delta x}[f(x)]$ | $\frac{T}{z-1} Z[f(x)]$ |

$K_{\Delta x}$ Transform to Z Transform Conversion

$$Z[f(x)] = \frac{z}{T} K_{\Delta x}[f(x)]\Big|_{s = \frac{z-1}{T}}$$ Z[f(x)] = F(z) Z Transform

$T = \Delta x$

$x = nT$, $n = 0,1,2,3,\ldots$ $K_{\Delta x}[f(x)] = f(s)$ $K_{\Delta x}$ Transform

Comment - Often the variable t is used instead of the variable x

# TABLE 4
## Fundamental Interval Calculus Functions and Definitions

| | $f_{\Delta x}(x)$ | Definition of $f_{\Delta x}(x)$ | $f(x)$ | Definition of $f(x)$ |
|---|---|---|---|---|
| **1** | $[x]^n_{\Delta x}$<br>$n = 0,1,2,3, \ldots$ | $1 \quad n = 0$<br>$\prod_{m=1}^{n}(x-(m-1)\Delta x) , \quad n = 1,2,3\ldots$<br>$x = x + r\Delta x \ , \quad r = \text{integer}$ | $x^n$ | $1 \quad n = 0$<br>$\prod_{m=1}^{n} x \quad n = 1,2,3,\ldots$ |
| **1a** | $[x]^{-n}_{\Delta x}$<br>$n = 1,2,3, \ldots$ | $\frac{1}{[x]^n_{\Delta x}}$<br>$x = x + r\Delta x \ , \quad r = \text{integer}$ | $x^{-n}$ | $1 \quad n = 0$<br>$\prod_{m=1}^{n} \frac{1}{x} \quad n = 1,2,3, ..$ |
| **2** | $e_{\Delta x}(a,x)$ | $(1+a\Delta x)^{\frac{x}{\Delta x}}$ | $e^{ax}$ | $\lim_{\Delta x \to 0}(1+a\Delta x)^{\frac{x}{\Delta x}}$ |
| **3** | $e_{\Delta x}(a,x)$ | $\left[e^{\ln(1+a\Delta x)}\right]^{\frac{x}{\Delta x}}$ | $e^{ax}$ | $\lim_{\Delta x \to 0}(1+a\Delta x)^{\frac{x}{\Delta x}}$ |
| **4** | $\sin_{\Delta x}x$ | $\frac{(1+j\Delta x)^{\frac{x}{\Delta x}} - (1-j\Delta x)^{\frac{x}{\Delta x}}}{2j}$ | $\sin x$ | $\frac{e^{jx}-e^{-jx}}{2j}$ |
| **5** | $\sin_{\Delta x}x$ | $\sin_{\Delta x}(1,x)$ | $\sin x$ | $\frac{e^{jx}-e^{-jx}}{2j}$ |
| **6** | $\cos_{\Delta x}x$ | $\frac{(1+j\Delta x)^{\frac{x}{\Delta x}} + (1-j\Delta x)^{\frac{x}{\Delta x}}}{2}$ | $\cos x$ | $\frac{e^{jx}+e^{-jx}}{2}$ |
| **7** | $\cos_{\Delta x}x$ | $\cos_{\Delta x}(1,x)$ | $\cos x$ | $\frac{e^{jx}+e^{-jx}}{2}$ |
| **8** | $\sin_{\Delta x}(a,x)$ | $\frac{(1+ja\Delta x)^{\frac{x}{\Delta x}} - (1-ja\Delta x)^{\frac{x}{\Delta x}}}{2j}$ | $\sin ax$ | $\frac{e^{jax}-e^{-jax}}{2j}$ |
| **9** | $\sin_{\Delta x}(a,x)$ | $(1+[a\Delta x]^2)^{\frac{x}{2\Delta x}} \sin\beta\frac{x}{\Delta x}$<br>$\beta = \tan^{-1}a\Delta x$<br>or<br>$(1+[a\Delta x]^2)^{\frac{x}{2\Delta x}} \sin\beta x$<br>$\beta = \frac{1}{\Delta t} \tan^{-1}a\Delta x$ | $\sin ax$ | $\frac{e^{jax}-e^{-jax}}{2j}$ |

| | $f_{\Delta x}(x)$ | Definition of $f_{\Delta x}(x)$ | f(x) | Definition of f(x) |
|---|---|---|---|---|
| **10** | $\sin_{\Delta x}(a,x)$ | $(\sec\beta)^{\frac{x}{\Delta x}} \sin\beta\frac{x}{\Delta x}$ <br> $\beta = \tan^{-1}a\Delta x$ | sinax | $\frac{e^{jax}-e^{-jax}}{2j}$ |
| **11** | $\cos_{\Delta x}(a,x)$ | $\frac{(1+ja\Delta x)^{\frac{x}{\Delta x}} + (1-ja\Delta x)^{\frac{x}{\Delta x}}}{2}$ | cosax | $\frac{e^{jax}+e^{-jax}}{2}$ |
| **12** | $\cos_{\Delta x}(a,x)$ | $(1+[a\Delta x]^2)^{\frac{x}{2\Delta x}} \cos\beta\frac{x}{\Delta x}$ <br> $\beta = \tan^{-1}a\Delta x$ <br> or <br> $(1+[a\Delta x]^2)^{\frac{x}{2\Delta x}} \cos\beta x$ <br> $\beta = \frac{1}{\Delta t}\tan^{-1}a\Delta x$ | cosax | $\frac{e^{jax}+e^{-jax}}{2}$ |
| **13** | $\cos_{\Delta x}(a,x)$ | $(\sec\beta)^{\frac{x}{\Delta x}} \cos\beta\frac{x}{\Delta x}$ <br> $\beta = \tan^{-1}a\Delta x$ | cosax | $\frac{e^{jax}+e^{-jax}}{2}$ |
| **14** | $e_{\Delta x}(a,x)\sin_{\Delta x}(\frac{b}{1+a\Delta x},x)$ | $\frac{(1+a\Delta x+jb\Delta x)^{\frac{x}{\Delta x}} - (1+a\Delta x-jb\Delta x)^{\frac{x}{\Delta x}}}{2j}$ | $e^{ax}\sin bx$ | $\frac{e^{a+jbx}-e^{a-jbx}}{2j}$ |
| **15** | $e_{\Delta x}(a,x)\sin_{\Delta x}(\frac{b}{1+a\Delta x},x)$ | $[(1+a\Delta x)^2+(b\Delta x)^2]^{\frac{x}{2\Delta x}} \sin\beta\frac{x}{\Delta x}$ <br> $\beta = \begin{cases} \tan^{-1}\left\lvert\frac{b\Delta x}{1+a\Delta x}\right\rvert & \text{for } 1+a\Delta x > 0 \quad b\Delta x \geq 0 \\ -\tan^{-1}\left\lvert\frac{b\Delta x}{1+a\Delta x}\right\rvert & \text{for } 1+a\Delta x > 0 \quad b\Delta x < 0 \\ \pi - \tan^{-1}\left\lvert\frac{b\Delta x}{1+a\Delta x}\right\rvert & \text{for } 1+a\Delta x < 0 \quad b\Delta x \geq 0 \\ -\pi + \tan^{-1}\left\lvert\frac{b\Delta x}{1+a\Delta x}\right\rvert & \text{for } 1+a\Delta x < 0 \quad b\Delta x < 0 \end{cases}$ <br> $1+a\Delta x \neq 0$ <br> $0 \leq \tan^{-1}\left\lvert\frac{b\Delta x}{1+a\Delta x}\right\rvert < \frac{\pi}{2}$ | $e^{ax}\sin bx$ | $\frac{e^{a+jbx}-e^{a-jbx}}{2j}$ |

| | $f_{\Delta x}(x)$ | Definition of $f_{\Delta x}(x)$ | $f(x)$ | Definition of $f(x)$ |
|---|---|---|---|---|
| **16** | $e_{\Delta x}(a,x)\sin_{\Delta x}(\frac{b}{1+a\Delta x},x)$ | $e_{\Delta x}(a,x)(\sec\beta)^{\frac{x}{\Delta x}}\sin\beta\frac{x}{\Delta x}$<br>$\beta = \begin{cases} \tan^{-1}\left\lvert\frac{b\Delta x}{1+a\Delta x}\right\rvert & \text{for } 1+a\Delta x > 0 \quad b\Delta x \ge 0 \\ -\tan^{-1}\left\lvert\frac{b\Delta x}{1+a\Delta x}\right\rvert & \text{for } 1+a\Delta x > 0 \quad b\Delta x < 0 \\ \pi - \tan^{-1}\left\lvert\frac{b\Delta x}{1+a\Delta x}\right\rvert & \text{for } 1+a\Delta x < 0 \quad b\Delta x \ge 0 \\ -\pi + \tan^{-1}\left\lvert\frac{b\Delta x}{1+a\Delta x}\right\rvert & \text{for } 1+a\Delta x < 0 \quad b\Delta x < 0 \end{cases}$<br>$1+a\Delta x \neq 0$<br>$0 \le \tan^{-1}\left\lvert\frac{b\Delta x}{1+a\Delta x}\right\rvert < \frac{\pi}{2}$ | $e^{ax}\sin bx$ | $\frac{e^{a+jbx}-e^{a-jbx}}{2j}$ |
| **17** | $e_{\Delta x}(a,x)\cos_{\Delta x}(\frac{b}{1+a\Delta x},x)$ | $\frac{(1+a\Delta x+jb\Delta x)^{\frac{x}{\Delta x}}+(1+a\Delta x-jb\Delta x)^{\frac{x}{\Delta x}}}{2}$ | $e^{ax}\cos bx$ | $\frac{e^{a+jbx}-e^{a-jbx}}{2j}$ |
| **18** | $e_{\Delta x}(a,x)\cos_{\Delta x}(\frac{b}{1+a\Delta x},x)$ | $[(1+a\Delta x)^2+(b\Delta x)^2]^{\frac{x}{2\Delta x}}\cos\beta\frac{x}{\Delta x}$<br>$\beta = \begin{cases} \tan^{-1}\left\lvert\frac{b\Delta x}{1+a\Delta x}\right\rvert & \text{for } 1+a\Delta x > 0 \quad b\Delta x \ge 0 \\ -\tan^{-1}\left\lvert\frac{b\Delta x}{1+a\Delta x}\right\rvert & \text{for } 1+a\Delta x > 0 \quad b\Delta x < 0 \\ \pi - \tan^{-1}\left\lvert\frac{b\Delta x}{1+a\Delta x}\right\rvert & \text{for } 1+a\Delta x < 0 \quad b\Delta x \ge 0 \\ -\pi + \tan^{-1}\left\lvert\frac{b\Delta x}{1+a\Delta x}\right\rvert & \text{for } 1+a\Delta x < 0 \quad b\Delta x < 0 \end{cases}$<br>$1+a\Delta x \neq 0$<br>$0 \le \tan^{-1}\left\lvert\frac{b\Delta x}{1+a\Delta x}\right\rvert < \frac{\pi}{2}$ | $e^{ax}\cos bx$ | $\frac{e^{a+jbx}+e^{a-jbx}}{2}$ |
| **19** | $e_{\Delta x}(a,x)\cos_{\Delta x}(\frac{b}{1+a\Delta x},x)$ | $e_{\Delta x}(a,x)(\sec\beta)^{\frac{x}{\Delta x}}\cos\beta\frac{x}{\Delta x}$<br>$\beta = \begin{cases} \tan^{-1}\left\lvert\frac{b\Delta x}{1+a\Delta x}\right\rvert & \text{for } 1+a\Delta x > 0 \quad b\Delta x \ge 0 \\ -\tan^{-1}\left\lvert\frac{b\Delta x}{1+a\Delta x}\right\rvert & \text{for } 1+a\Delta x > 0 \quad b\Delta x < 0 \\ \pi - \tan^{-1}\left\lvert\frac{b\Delta x}{1+a\Delta x}\right\rvert & \text{for } 1+a\Delta x < 0 \quad b\Delta x \ge 0 \\ -\pi + \tan^{-1}\left\lvert\frac{b\Delta x}{1+a\Delta x}\right\rvert & \text{for } 1+a\Delta x < 0 \quad b\Delta x < 0 \end{cases}$ | $e^{ax}\cos bx$ | $\frac{e^{a+jbx}+e^{a-jbx}}{2}$ |

| | $f_{\Delta x}(x)$ | Definition of $f_{\Delta x}(x)$ | f(x) | Definition of f(x) |
|---|---|---|---|---|
| | | $1+a\Delta x \neq 0$ <br> $0 \leq \tan^{-1}\left\lvert\frac{b\Delta x}{1+a\Delta x}\right\rvert < \frac{\pi}{2}$ | | |
| **20** | $\tan_{\Delta x}(a,x)$ | $\frac{(1+ja\Delta x)^{\frac{x}{\Delta x}} - (1-ja\Delta x)^{\frac{x}{\Delta x}}}{(1+ja\Delta x)^{\frac{x}{\Delta x}} + (1-ja\Delta x)^{\frac{x}{\Delta x}}}$ | tanax | $\frac{e^{jax}-e^{-jax}}{e^{jax}+e^{-jax}}$ |
| **21** | $\tan_{\Delta x}(a,x)$ | $\tan\beta\frac{x}{\Delta x}$ <br> $\beta = \tan^{-1}a\Delta x$ | tanax | $\frac{e^{jax}-e^{-jax}}{e^{jax}+e^{-jax}}$ |
| **22** | $\sinh_{\Delta x}x$ | $\frac{(1+\Delta x)^{\frac{x}{\Delta x}} - (1-\Delta x)^{\frac{x}{\Delta x}}}{2}$ | sinhx | $\frac{e^{x}-e^{-x}}{2}$ |
| **23** | $\sinh_{\Delta x}x$ | $\sinh_{\Delta x}(1,x)$ | sinhx | $\frac{e^{x}-e^{-x}}{2}$ |
| **24** | $\cosh_{\Delta x}x$ | $\frac{(1+\Delta x)^{\frac{x}{\Delta x}} + (1-\Delta x)^{\frac{x}{\Delta x}}}{2}$ | coshx | $\frac{e^{x}+e^{-x}}{2}$ |
| **25** | $\cosh_{\Delta x}x$ | $\cosh_{\Delta x}(1,x)$ | coshx | $\frac{e^{x}+e^{-x}}{2}$ |
| **26** | $\sinh_{\Delta x}(a,x)$ | $\frac{(1+a\Delta x)^{\frac{x}{\Delta x}} - (1-a\Delta x)^{\frac{x}{\Delta x}}}{2}$ | sinhax | $\frac{e^{ax}-e^{-ax}}{2}$ |
| **27** | $\sinh_{\Delta x}(a,x)$ | $(1-[a\Delta x]^2)^{\frac{x}{2\Delta x}} \sinh\beta\frac{x}{\Delta x}$ <br> $\beta = \tanh^{-1}a\Delta x$ | sinhax | $\frac{e^{ax}-e^{-ax}}{2}$ |
| **28** | $\cosh_{\Delta x}(a,x)$ | $\frac{(1+a\Delta x)^{\frac{x}{\Delta x}} + (1-a\Delta x)^{\frac{x}{\Delta x}}}{2}$ | coshax | $\frac{e^{ax}+e^{-ax}}{2}$ |
| **29** | $\cosh_{\Delta x}(a,x)$ | $(1-[a\Delta x]^2)^{\frac{x}{\Delta x}} \cosh\beta\frac{x}{\Delta x}$ <br> $\beta = \tanh^{-1}a\Delta x$ | coshax | $\frac{e^{ax}+e^{-ax}}{2}$ |
| **30** | $\tanh_{\Delta x}(a,x)$ | $\frac{(1+a\Delta x)^{\frac{x}{\Delta x}} - (1-a\Delta x)^{\frac{x}{\Delta x}}}{(1+a\Delta x)^{\frac{x}{\Delta x}} + (1-a\Delta x)^{\frac{x}{\Delta x}}}$ | tanhax | $\frac{e^{ax}-e^{-ax}}{e^{ax}+e^{-ax}}$ |
| **31** | $\tanh_{\Delta x}(a,x)$ | $\tanh\beta\frac{x}{\Delta x}$ <br> $\beta = \tanh^{-1}a\Delta x$ | tanhax | $\frac{e^{ax}-e^{-ax}}{e^{ax}+e^{-ax}}$ |

**Note :** $\lim_{\Delta x\to 0}f_{\Delta x}(x) = f(x)$

# TABLE 5
## Interval Calculus Equations and Identities

| # | Equations and Identities | Comments |
|---|---|---|
| **1** | <u>Interval Calculus Function $\Delta x \to 0$ Limits</u><br><br>$\lim_{\Delta x \to 0}[\ln_{\Delta x}x - \ln_{\Delta x}1] = \ln x$ , n=1<br>$\lim_{\Delta x \to 0}\text{lnd}(n,\Delta x,x) = \frac{1}{(n-1)x^{n-1}}$ , $n \neq 1$<br>$\lim_{\Delta x \to 0}e_{\Delta x}(a,x) = e^{ax}$<br>$\lim_{\Delta x \to 0}\sin_{\Delta x}(a,x) = \sin ax$<br>$\lim_{\Delta x \to 0}\cos_{\Delta x}(a,x) = \cos ax$<br>$\lim_{\Delta x \to 0}e_{\Delta x}x = \lim_{\Delta x \to 0}e_{\Delta x}(1,x) = e^x$<br>$\lim_{\Delta x \to 0}\sin_{\Delta x}x = \lim_{\Delta x \to 0}\sin_{\Delta x}(1,x) = \sin x$<br>$\lim_{\Delta x \to 0}\cos_{\Delta x}x = \lim_{\Delta x \to 0}\cos_{\Delta x}(1,x) = \cos x$ | All Interval Calculus functions become their Calculus counterpart functions as $\Delta x$ goes to zero.<br><br>Note that only as $\Delta x$ goes to zero does (a,x) become ax.<br><br>Note that the counterpart function of lnx is not $\ln_{\Delta x}x$, it is $\ln_{\Delta x}x - \ln_{\Delta x}1$. |
| **2** | $e_{\Delta x}(a,x) = (1+a\Delta x)^{\frac{x}{\Delta x}}$<br>$x = 0,\Delta x,2\Delta x,3\Delta x,\ldots$ | |
| **3** | $e_{\Delta x}(a,x) = e_{m\Delta x}(\frac{(1+a\Delta x)^m-1}{m\Delta x},x)$<br>$x = 0,\Delta x,2\Delta x,3\Delta x,\ldots$<br>$m = 1,2,3,\ldots$ | $e_{\Delta x}(a,x)$ variable subscript identity |
| **4** | $e_{\Delta x}(a,x) = e_{2\Delta x}(\frac{a}{2},2x)$<br>$x = 0,\Delta x,2\Delta x,3\Delta x,\ldots$ | |
| **5** | $e_{\Delta x}(a,x) = e_{m\Delta x}(\frac{a}{m},mx)$<br>$x = 0,\Delta x,2\Delta x,3\Delta x,\ldots$<br>$m = 1,2,3,\ldots$ | |
| **6** | $e_{\Delta x}(a,mx) = e_{\Delta x}(\frac{(1+a\Delta x)^m-1}{\Delta x},x)$<br>$x = 0,\Delta x,2\Delta x,3\Delta x,\ldots$<br>$m = 1,2,3,\ldots$ | |
| **7** | $e_{\Delta x}(a,x) = \left[e^{\ln(1+a\Delta x)}\right]^{\frac{x}{\Delta x}}$ | |
| **8** | $e_{\Delta x}(a,-x) = (1+a\Delta x)^{-\frac{x}{\Delta x}}$ | |
| **9** | $\frac{1}{e_{\Delta x}(a,x)} = e_{\Delta x}(a,-x)$ | |
| **10** | $\frac{1}{e_{\Delta x}(a,-x)} = e_{\Delta x}(a,x)$ | |
| **11** | $e_{\Delta x}(a,x)e_{\Delta x}(b,x) = e_{\Delta x}(a+b+ab\Delta x,x)$ | |
| **12** | $e_{\Delta x}(a,x)e_{\Delta x}(b,-x) = e_{\Delta x}(\frac{a-b}{1+b\Delta x},x)$ | |
| **13** | $e_{\Delta x}(a+b,x) = e_{\Delta x}(a,x)\, e_{\Delta x}(\frac{b}{1+a\Delta x},x)$ | |
| **14** | $e_{\Delta x}(a,x+n\Delta x) = (1+a\Delta x)^{n}\, e_{\Delta x}(a,x)$ | |
| **15** | $e_{\Delta x}(a,x)e_{\Delta x}(a,y) = e_{\Delta x}(a,x+y)$ | |
| **16** | $e_{\Delta x}(a,n\Delta x) = (1+a\Delta x)^{n}$ | |
| **17** | $e_{\Delta x}(a,0) = 1$ | |
| **18** | $e_{\Delta x}(0,x) = 1$ | |
| **19** | $e_{\Delta x}(\frac{e^{a\Delta x}-1}{\Delta x},x) = e^{ax}$ | Identity |
| **20** | $e_{\Delta x}(a,x) = \frac{1+\tanh_{\Delta x}(a, \frac{x}{2})}{1-\tanh_{\Delta x}(a, \frac{x}{2})}$ | |
| **21** | $e_{\Delta x}(ja,x) = \frac{1+j\tan_{\Delta x}(a, \frac{x}{2})}{1-j\tan_{\Delta x}(a, \frac{x}{2})}$ | |
| **22** | $e_{\Delta x}(j,x) = \cos_{\Delta x}x + j\sin_{\Delta x}x$ | |
| **23** | $e_{\Delta x}(-j,x) = \cos_{\Delta x}x - j\sin_{\Delta x}x$ | |
| **24** | $e_{\Delta x}(ja,x) = \cos_{\Delta x(}a,x) + j\sin_{\Delta}x(a,x) = (1+ja\Delta x)^{\frac{x}{\Delta x}}$ | |
| **25** | $e_{\Delta x}(-ja,x) = \cos_{\Delta x}(a,x) - j\sin_{\Delta}x(a,x) = (1-ja\Delta x)^{\frac{x}{\Delta x}}$ | |
| **26** | $\frac{e_{\Delta x}(a+jb,x) - e_{\Delta x}(a-jb,x)}{2j} = e_{\Delta x}(a,x)\,\sin_{\Delta x}(\frac{b}{1+a\Delta x},x)$ | |
| **27** | $\frac{e_{\Delta x}(a+jb,x) + e_{\Delta x}(a-jb,x)}{2} = e_{\Delta x}(a,x)\,\cos_{\Delta x}(\frac{b}{1+a\Delta x},x)$ | |
| **28** | $(1+\Delta x^2)^{\frac{x}{\Delta x}} = \sin_{\Delta x}{}^2(1,x) + \cos_{\Delta x}{}^2(1,x) = \sin_{\Delta x}{}^2 x + \cos_{\Delta x}{}^2 x$ | |
| **29** | $(1+[a\Delta x]^2)^{\frac{x}{\Delta x}} = \sin_{\Delta x}{}^2(a,x) + \cos_{\Delta x}{}^2(a,x)$ | |
| **30** | $(1-[a\Delta x]^2)^{\frac{x}{\Delta x}} = \cosh_{\Delta x}{}^2(a,x) - \sinh_{\Delta}x^2(a,x)$ | |
| **31** | $\sin_{\Delta x}(a,x) = (\sec\beta)^{\frac{x}{\Delta x}} \sin\beta\frac{x}{\Delta x}$<br>$\beta = \tan^{-1}a\Delta x$ | Other forms of this equation are found in Table 4 |
| **32** | $e_{\Delta x}(0,x) = 1$ | |
| **33** | $e_{\Delta x}(a,0) = 1$ | |
| **34** | $\sin_{\Delta x}(0,x) = 0$ | |
| **35** | $\sin_{\Delta x}(b,0) = 0$ | |
| **36** | $\cos_{\Delta x}(0,x) = 1$ | |
| **37** | $\cos_{\Delta x}(b,0) = 1$ | |
| **38** | $\sin_{\Delta x}(b,\Delta x) = b\Delta x$ | |
| **39** | $\cos_{\Delta x}(b,\Delta x) = 1$ | |
| **40** | $\sin_{\Delta x}(b,2\Delta x) = 2b\Delta x$ | |
| **41** | $\cos_{\Delta x}(b,2\Delta x) = 1 - b^2\Delta x^2$ | |
| **42** | $\sin_{\Delta x}(b,m\Delta x) = [1+(b\Delta x)^2]^{\frac{m}{2}}\sin(m\beta)$<br>$\beta = \tan^{-1}b\Delta t$<br>$m = 0,1,2,3, \ldots$ | |
| **43** | $\cos_{\Delta x}(b,m\Delta x) = [1+(b\Delta x)^2]^{\frac{m}{2}}\cos(m\beta)$<br>$\beta = \tan^{-1}b\Delta t$<br>$m = 0,1,2,3, \ldots$ | |
| **44** | $\sin_{\Delta x}(b,x+m\Delta x) = \cos_{\Delta x}(b,m\Delta x)\sin_{\Delta x}(a,x) + \sin_{\Delta x}(b,m\Delta x)\cos_{\Delta x}(a,x)$<br>$m = 0,1,2,3, \ldots$ | |
| **45** | $\cos_{\Delta x}(a,x+m\Delta x) = \cos_{\Delta x}(a,m\Delta x)\cos_{\Delta x}(a,x) - \sin_{\Delta x}(a,m\Delta x)\sin_{\Delta x}(a,x)$<br>$m = 0,1,2,3, \ldots$ | |
| **46** | $\sin_{\Delta x}(a,x) = (\sec\beta)^{\frac{x}{\Delta x}}\left[\dfrac{2\tan_{\Delta x}(a,\frac{x}{2})}{1+\tan_{\Delta x}^2(a,\frac{x}{2})}\right]$<br>$\beta = \tan^{-1}a\Delta x$ | $\sec^2\beta = 1+[a\Delta x]^2$ |
| **47** | $\sin_{\Delta x}(-a,x) = -\sin_{\Delta x}(a,x)$ | |
| **48** | $\sin_{\Delta x}(a,x) = -j\sinh_{\Delta x}(ja,x)$ | |
| **49** | $\sin_{\Delta x}(ja,x) = j\sinh_{\Delta x}(a,x)$ | |
| **50** | $\sin_{\Delta x}2x = 2\sin_{\Delta x}x\cos_{\Delta x}x$ | |
| **51** | $\sin_{\Delta x}(a,2x) = 2\sin_{\Delta x}(a,x)\cos_{\Delta x}(a,x)$ | |
| **52** | $\sin_{\Delta x}(x+y) = \sin_{\Delta x}x\cos_{\Delta x}y + \cos_{\Delta x}x\sin_{\Delta x}y$ | |
| **53** | $\sin_{\Delta x}(a,x+y) = \sin_{\Delta x}(a,x)\cos_{\Delta x}(a,y) + \sin_{\Delta x}(a,y)\cos_{\Delta x}(a,x)$ | |
| **54** | $\sin_{\Delta x}(a,x+\Delta x) = \sin_{\Delta x}(a,x) + a\Delta x\cos_{\Delta x}(a,x)$ | |
| **55** | $\sin_{\Delta x}^2(a,x) = \dfrac{(1+[a\Delta x]^2)^{\frac{x}{\Delta x}} - \cos_{\Delta x}(a,2x)}{2}$ | |
| **56** | $\cos_{\Delta x}(a,x) = (\sec\beta)^{\frac{x}{\Delta x}}\cos\beta\dfrac{x}{\Delta x}$<br>$\beta = \tan^{-1}a\Delta x$ | Other forms of this equation are found in Table 4 |
| **57** | $\cos_{\Delta x}(0,x) = 1$ | |
| **58** | $\cos_{\Delta x}(a,x) = (\sec\beta)^{\frac{x}{\Delta x}}\left[\dfrac{1 - \tan_{\Delta x}{}^{2}(a, \frac{x}{2})}{1 + \tan_{\Delta x}{}^{2}(a, \frac{x}{2})}\right]$<br>$\beta = \tan^{-1}a\Delta x$ | $\sec^{2}\beta = 1+[a\Delta x]^{2}$ |
| **59** | $\cos_{\Delta x}(-a,x) = \cos_{\Delta x}(a,x)$ | |
| **60** | $\cos_{\Delta x}(a,x) = \cosh_{\Delta x}(ja,x)$ | |
| **61** | $\cos_{\Delta x}(ja,x) = \cosh_{\Delta x}(a,x)$ | |
| **62** | $\cos_{\Delta x}2x = \cos_{\Delta x}{}^{2}x - \sin_{\Delta x}{}^{2}x$ | |
| **63** | $\cos_{\Delta x}(a,2x) = \cos_{\Delta x}{}^{2}(a,x) - \sin_{\Delta x}{}^{2}(a,x)$ | |
| **64** | $\cos_{\Delta x}(x+y) = \cos_{\Delta x}x \cos_{\Delta x}y – \sin_{\Delta x}x \sin_{\Delta x}y$ | |
| **65** | $\cos_{\Delta x}(a,x+y) = \cos_{\Delta x}(a,x) \cos_{\Delta x}(a,y) – \sin_{\Delta x}(a,x) \sin_{\Delta x}(a,y)$ | |
| **66** | $\cos_{\Delta x}(a,x+\Delta x) = \cos_{\Delta x}(a,x) - a\Delta x\sin_{\Delta x}(a,x)$ | |
| **67** | $\cos_{\Delta x}{}^{2}(a,x) = \dfrac{(1+[a\Delta x]^{2})^{\frac{x}{\Delta x}} + \cos_{\Delta x}(a,2x)}{2}$ | |
| **68** | $[\sin_{\Delta x}{}^{2}x + \cos_{\Delta x}{}^{2}x]^{a} = \sin_{\frac{\Delta x}{a}}{}^{2}(a,x) + \cos_{\frac{\Delta x}{a}}{}^{2}(a,x)$ | |
| **69** | $\sin_{\Delta x}(a,x \pm\frac{\pi\Delta x}{2b}) = \pm (\sec b)^{\pm\frac{\pi}{2b}} \cos_{\Delta x}(a,x)$<br>$b = \tan^{-1}a\Delta x$ | |
| **70** | $\cos_{\Delta x}(a,x \pm\frac{\pi\Delta x}{2b}) = \mp (\sec b)^{\pm\frac{\pi}{2b}} \sin_{\Delta x}(a,x)$<br>$b = \tan^{-1}a\Delta x$ | |
| **71** | $\sin_{\Delta x}(a,x \pm \frac{\pi\Delta x}{b}) = -(\sec b)^{\pm\frac{\pi}{b}} \sin_{\Delta x}(a,x)$<br>$b = \tan^{-1} a\Delta x$ | |
| **72** | $\cos_{\Delta x}(a,x \pm \frac{\pi\Delta x}{b}) = -(\sec b)^{\pm\frac{\pi}{b}} \cos_{\Delta x}(a,x)$<br>$b = \tan^{-1} a\Delta x$ | |
| **73** | $[x]^{n}_{\Delta x} = (x)(x-\Delta x)(x-2\Delta x)\ldots(x-[n-1]\Delta x)$ , $n = 1,2,3,\ldots$<br>or<br>$[x]^{n}_{\Delta x} = \prod_{m=1}^{n}(x-[m-1]\Delta x)$ , $n = 1,2,3,\ldots$<br>$[x]^{0}_{\Delta x} = 1$<br>$[x]^{-m}_{\Delta x} = \frac{1}{[x]^{m}_{\Delta x}}$ , m = integers | |
| **74** | $[x]^{n}_{\Delta x} = [x]^{m}_{\Delta x}[x-m\Delta x]^{n-m}_{\Delta x}$ , m,n = positive integers , $n \geq m$<br>or<br>$[x]^{n+m}_{\Delta x} = [x]^{m}_{\Delta x}[x-m\Delta x]^{n}_{\Delta x}$ , m,n = positive integers | |
| **75** | $[x]^{n}_{\Delta x} \sin_{\Delta x}(b,x) = (1+[b\Delta x]^2)^n \cos_{\Delta x}(b,-n\Delta x)[x]^{n}_{\Delta x} \sin_{\Delta x}(b,x-n\Delta x) -$<br>$(1+[b\Delta x]^2)^n \sin_{\Delta x}(b,-n\Delta x)\,[x]^{n}_{\Delta x} \cos_{\Delta x}(b,x-n\Delta x)$ | |
| **76** | $[x]^{n}_{\Delta x} \cos_{\Delta x}(b,x) = (1+[b\Delta x]^2)^n \sin_{\Delta x}(b,-n\Delta x)\,[x]^{n}_{\Delta x} \sin_{\Delta x}(b,x-n\Delta x) +$<br>$(1+[b\Delta x]^2)^n \cos_{\Delta x}(b,-n\Delta x)\,[x]^{n}_{\Delta x} \cos_{\Delta x}(b,x-n\Delta x)$ | |
| **77** | $K_1\sin_{\Delta x}(a,x) + K_2\cos_{\Delta x}(a,x) =$<br>$\sqrt{K_1^2 + K_2^2}\,(\sec b)^{\frac{\beta}{b}} \cos_{\Delta x}(a,x - \frac{\beta}{b}\Delta x)$<br>$K_1, K_2, a$ = constants<br>$b = \tan^{-1} a\Delta x$<br>$\beta = \tan^{-1}\frac{K_1}{K_2}$ | Sum of $K_1\sin_{\Delta x}(a,x)$ and $K_2\cos_{\Delta x}(a,x)$ |
| **78** | $g_{\Delta x}(w,ax) = g_{\frac{\Delta x}{a}}(aw,x)$ | $g_{\Delta x}(p,q)$ is an Interval Calculus function |
| **79** | $\tan_{\Delta x}(a,x) = \tan\beta^{\frac{x}{\Delta x}}$<br>$\beta = \tan^{-1}(a\Delta x)$ | $\tan^{-1}a\Delta x = \frac{i}{2}\ln\frac{i+a\Delta x}{i-a\Delta x}$ |
| **80** | $\tan_{\Delta x}(a,x) = \dfrac{\pi}{\text{lnd}(1,1,1-\frac{bx}{\pi}) - \text{lnd}(1,1,\frac{bx}{\pi})}$<br>$b = \dfrac{\tan^{-1}(a\Delta x)}{\Delta x}$<br>for $bx = n\pi$, □ n=integer, $\tan_{\Delta x}(a,x)=0$ | |
| **81** | $\tan_{\Delta x}(a,x) = \dfrac{2\tan_{\Delta x}(a,\frac{x}{2})}{1-\tan_{\Delta x}^{2}(a,\frac{x}{2})}$ | |
| **82** | $\tan_{\Delta x}(-a,x) = -\tan_{\Delta x}(a,x)$ | |
| **83** | $\tan_{\Delta x}(a,x) = -j\tanh_{\Delta x}(ja,x)$ | |
| **84** | $\tan_{\Delta x}^{2}(a,x) = \dfrac{(\sec\beta)^{\frac{2x}{\Delta x}} - \cos_{\Delta x}(a,2x)}{(\sec\beta)^{\frac{2x}{\Delta x}} + \cos_{\Delta x}(a,2x)}$<br>$\beta = \tan^{-1}(a\Delta x)$ | $\sec^2\beta = 1+[a\Delta x]^2$ |
| **85** | $\csc_{\Delta x}(a,x) = \dfrac{1}{\sin_{\Delta x}(a,x)}$ | |
| **86** | $\sec_{\Delta x}(a,x) = \dfrac{1}{\cos_{\Delta x}(a,x)}$ | |
| **87** | $\text{cotan}_{\Delta x}(a,x) = \dfrac{1}{\tan_{\Delta x}(a,x)}$ | |
| **88** | $\sinh_{\Delta x}(a,x) = (1-[a\Delta x]^2)^{\frac{x}{2\Delta x}}\sinh\beta^{\frac{x}{\Delta x}}$<br>$\beta = \tanh^{-1}(a\Delta x)$ | |
| **89** | $\sinh_{\Delta x}(a,x) = (1-[a\Delta x]^2)^{\frac{x}{2\Delta x}}\left[\frac{2\tanh_{\Delta x}(a, \frac{x}{2})}{1 - \tanh_{\Delta x}^2(a, \frac{x}{2})}\right]$ | |
| **90** | $\sinh_{\Delta x}(-a,x) = -\sinh_{\Delta x}(a,x)$ | |
| **91** | $\sinh_{\Delta x}(a,x) = -j\sin_{\Delta x}(ja,x)$ | |
| **92** | $\cosh_{\Delta x}(a,x) = (1-[a\Delta x]^2)^{\frac{x}{2\Delta x}}\cosh\beta\frac{x}{\Delta x}$<br>$\beta = \tanh^{-1}(a\Delta x)$ | $1-a\Delta x \neq 0$ |
| **93** | $\cosh_{\Delta x}(a,x) = (1-[a\Delta x]^2)^{\frac{x}{2\Delta x}}\left[\frac{1 + \tanh^2_{\Delta x}(a, \frac{x}{2})}{1 - \tanh_{\Delta x}^2(a, \frac{x}{2})}\right]$ | |
| **94** | $\cosh_{\Delta x}(-a,x) = \cosh_{\Delta x}(a,x)$ | |
| **95** | $\cosh_{\Delta x}(a,x) = \cos_{\Delta x}(ja,x)$ | |
| **96** | $\tanh_{\Delta x}(a,x) = \tanh\beta\frac{x}{\Delta x}$<br>$\beta = \tanh^{-1}a\Delta x$ | |
| **97** | $\tanh_{\Delta x}(a,x) = \frac{-j\pi}{\mathrm{lnd}(1,1,1-\frac{bx}{\pi}) - \mathrm{lnd}(1,1,\frac{bx}{\pi})}$<br>$b = \frac{\tan^{-1}(ja\Delta x)}{\Delta x} = \frac{j\tanh^{-1}(a\Delta x)}{\Delta x}$<br>for $bx = 0$ □ $\tanh_{\Delta x}(a,x)=0$ | |
| **98** | $\tanh_{\Delta x}(a,x) = \frac{2\tanh_{\Delta x}( a,\frac{x}{2})}{1 + \tanh_{\Delta x}^2( a,\frac{x}{2})}$ | |
| **99** | $\tanh_{\Delta x}(-a,x) = - \tanh_{\Delta x}(a,x)$ | |
| **100** | $\tanh_{\Delta x}(a,x) = -j\tan_{\Delta x}(ja,x)$ | |
| **101** | $\mathrm{csch}_{\Delta x}(a,x) = \frac{1}{\sinh_{\Delta x}(a,x)}$ | |
| **102** | $\text{sech}_{\Delta x}(a,x) = \dfrac{1}{\cosh_{\Delta x}(a,x)}$ | |
| **103** | $\text{cotanh}_{\Delta x}(a,x) = \dfrac{1}{\tanh_{\Delta x}(a,x)}$ | |
| **104** | $1 = [x]_{\Delta x}^{0}$ | Identity |
| **105** | $x = [x]_{\Delta x}^{1}$ | Identity |
| **106** | $x^2 = [x]_{\Delta x}^{2} + \Delta x[x]_{\Delta x}^{1}$ | Identity |
| **107** | $x^3 = [x]_{\Delta x}^{3} + 3\Delta x[x]_{\Delta x}^{2} + \Delta x^2[x]_{\Delta x}^{1}$ | Identity |
| **108** | $A^x = e_{\Delta x}( \dfrac{A^{\Delta x} -1}{\Delta x}, x )$ | Identity |
| **109** | $e^{ax} = e_{\Delta x}( \dfrac{e^{a\Delta x} -1}{\Delta x}, x ) = [1 + ( \dfrac{e^{a\Delta x} -1}{\Delta x} )\Delta x]^{\frac{x}{\Delta x}}$<br>or<br>$e^{ax} = e_{\Delta x}( \dfrac{e^{-a\Delta x} -1}{\Delta x}, -x ) = [1 + ( \dfrac{e^{-a\Delta x} -1}{\Delta x} )\Delta x]^{-\frac{x}{\Delta x}}$<br>where<br>$\dfrac{x}{\Delta x}$ = integer | Identity |
| **110** | $e^{-ax} = e_{\Delta x}( \dfrac{e^{-a\Delta x} -1}{\Delta x}, x ) = [1 + ( \dfrac{e^{-a\Delta x} -1}{\Delta x} )\Delta x]^{\frac{x}{\Delta x}}$<br>or<br>$e^{-ax} = e_{\Delta x}( \dfrac{e^{a\Delta x} -1}{\Delta x}, -x ) = [1 + ( \dfrac{e^{a\Delta x} -1}{\Delta x} )\Delta x]^{-\frac{x}{\Delta x}}$<br>where<br>$\dfrac{x}{\Delta x}$ = integer | Identity |
| **111** | $e^{jax} = e_{\Delta x}( \dfrac{e^{ja\Delta x} -1}{\Delta x}, x) = [1 + ( \dfrac{e^{ja\Delta x} -1}{\Delta x} )\Delta x]^{\frac{x}{\Delta x}}$<br>or<br>$e^{jax} = e_{\Delta x}( \dfrac{e^{-ja\Delta x} -1}{\Delta x}, -x ) = [1 + ( \dfrac{e^{-ja\Delta x} -1}{\Delta x} )\Delta x]^{-\frac{x}{\Delta x}}$<br>where<br>$\dfrac{x}{\Delta x}$ = integer | Identity |
| **112** | $e^{-jax} = e_{\Delta x}(\frac{e^{-ja\Delta x} -1}{\Delta x}, x) = [1 + (\frac{e^{-ja\Delta x} -1}{\Delta x})\Delta x]^{\frac{x}{\Delta x}}$<br>or<br>$e^{-jax} = e_{\Delta x}(\frac{e^{ja\Delta x} -1}{\Delta x}, -x) = [1 + (\frac{e^{ja\Delta x} -1}{\Delta x})\Delta x]^{-\frac{x}{\Delta x}}$<br>where<br>$\frac{x}{\Delta x}$ = integer | Identity |
| **113** | $\sin bx = [\cos b\Delta x]^{\frac{x}{\Delta x}}\sin_{\Delta x}(\frac{\tan b\Delta x}{\Delta x}, x)$<br>or<br>$\sin bx = e_{\Delta x}(\frac{\cos b\Delta x -1}{\Delta x}, x)\sin_{\Delta x}(\frac{\tan b\Delta x}{\Delta x}, x)$<br>or<br>$\sin bx = [1+(\beta\Delta x)^2]^{-\frac{x}{2\Delta x}}\sin_{\Delta x}(\beta, x)$<br>where<br>$\beta = \frac{\tan b\Delta x}{\Delta x}$<br>$\cos b\Delta x \neq 0$<br>$x = m\Delta t$ ,  $m = 0,1,2,3, \ldots$<br>or<br>$\sin bx = \frac{1}{2j} e_{\Delta x}(A+jB,x) - \frac{1}{2j} e_{\Delta x}(A-jB,x)$<br>where<br>$A = \frac{\cos b\Delta x - 1}{\Delta x}$<br>$B = \frac{\sin b\Delta x}{\Delta x}$<br>$\frac{x}{\Delta x}$ = integer | Identity<br><br>If $\cos b\Delta x = 0$, function simplification with term cancellation yields sinbx<br><br>$\sin bx = \frac{1}{2j} e^{jbx} - \frac{1}{2j} e^{-jbx}$ |
| **114** | $e^{ax}\sin bx = [e^{a\Delta x}\cos b\Delta x]^{\frac{x}{\Delta x}}\sin_{\Delta x}(\frac{\tan b\Delta x}{\Delta x}, x)$<br>or<br>$e^{ax}\sin bx = e_{\Delta x}(\frac{e^{a\Delta x}\cos b\Delta x-1}{\Delta x}, x)\sin_{\Delta x}(\frac{\tan b\Delta x}{\Delta x}, x)$<br>where<br>$\frac{x}{\Delta x}$ = integer<br>$\cos b\Delta x \neq 0$ | Identity<br><br>If $\cos b\Delta x = 0$, function simplification with term cancellation yields $e^{ax}\sin bx$ |
| **115** | $\cos bx = [\cos b\Delta x]^{\frac{x}{\Delta x}} \cos_{\Delta x}(\frac{\tan b\Delta x}{\Delta x}, x)$<br>or<br>$\cos bx = e_{\Delta x}(\frac{\cos b\Delta x - 1}{\Delta x}, x) \cos_{\Delta x}(\frac{\tan b\Delta x}{\Delta x}, x)$<br>or<br>$\cos bx = [1+(\beta\Delta x)^2]^{-\frac{x}{2\Delta x}} \cos_{\Delta x}(\beta, x)$<br>where<br>$\beta = \frac{\tan b\Delta x}{\Delta x}$<br>$\cos b\Delta x \neq 0$<br>$x = m\Delta t$ , $m = 0,1,2,3, \ldots$<br>or<br>$\cos bx = \frac{1}{2} e_{\Delta x}(A+jB,x) + \frac{1}{2} e_{\Delta x}(A-jB,x)$<br>where<br>$A = \frac{\cos b\Delta x - 1}{\Delta x}$<br>$B = \frac{\sin b\Delta x}{\Delta x}$<br>$\frac{x}{\Delta x}$ = integer $\frac{x}{\Delta x}$ = integer | Identity<br><br>If $\cos b\Delta x = 0$, function simplification with term cancellation yields cosbx<br><br>$\cos bx = \frac{1}{2} e^{jbx} + \frac{1}{2} e^{-jbx}$ |
| **116** | $e^{ax}\cos bx = [e^{a\Delta x}\cos b\Delta x]^{\frac{x}{\Delta x}} \cos_{\Delta x}(\frac{\tan b\Delta x}{\Delta x}, x)$<br>or<br>$e^{ax}\cos bx = e_{\Delta x}(\frac{e^{a\Delta x}\cos b\Delta x - 1}{\Delta x}, x)\cos_{\Delta x}(\frac{\tan b\Delta x}{\Delta x}, x)$<br>where<br>$\frac{x}{\Delta x}$ = integer<br>$\cos b\Delta x \neq 0$ | Identity<br><br>If $\cos b\Delta x = 0$, function simplification with term cancellation yields $e^{ax}\cos bx$ |
| **117** | $\sinh bx = [\cosh b\Delta x]^{\frac{x}{\Delta x}} \sinh_{\Delta x}(\frac{\tanh b\Delta x}{\Delta x}, x)$<br>or<br>$\sinh bx = e_{\Delta x}(\frac{\cosh b\Delta x - 1}{\Delta x}, x)\sinh_{\Delta x}(\frac{\tanh b\Delta x}{\Delta x}, x)$<br>or<br>$\sinh bx = \frac{1}{2j} e_{\Delta x}(\frac{e^{b\Delta x} - 1}{\Delta x}, x) - \frac{1}{2j} e_{\Delta x}(\frac{e^{-b\Delta x} - 1}{\Delta x}, x)$ | Identity<br><br>$\sinh bx = \frac{1}{2j} e^{bx} - \frac{1}{2j} e^{-bx}$ |
| | | |
| **118** | $e^{ax}\sinh bx = [e^{a\Delta x}\cosh b\Delta x]^{\frac{x}{\Delta x}}\sinh_{\Delta x}(\frac{\tanh b\Delta x}{\Delta x}, x)$<br>or<br>$e^{ax}\sinh bx = e_{\Delta x}(\frac{e^{a\Delta x}\cosh b\Delta x - 1}{\Delta x}, x)\sinh_{\Delta x}(\frac{\tanh b\Delta x}{\Delta x}, x)$ | Identity |
| **119** | $\cosh bx = [\cosh b\Delta x]^{\frac{x}{\Delta x}}\cosh_{\Delta x}(\frac{\tanh b\Delta x}{\Delta x}, x)$<br>or<br>$\cosh bx = e_{\Delta x}(\frac{\cosh b\Delta x - 1}{\Delta x}, x)\cosh_{\Delta x}(\frac{\tanh b\Delta x}{\Delta x}, x)$<br>or<br>$\cosh bx = \frac{1}{2}e_{\Delta x}(\frac{e^{b\Delta x} - 1}{\Delta x}, x) + \frac{1}{2}e_{\Delta x}(\frac{e^{-b\Delta x} - 1}{\Delta x}, x)$ | Identity<br><br>$\cosh bx = \frac{1}{2}e^{bx} + \frac{1}{2}e^{-bx}$ |
| **120** | $e^{ax}\cosh bx = [e^{a\Delta x}\cosh b\Delta x]^{\frac{x}{\Delta x}}\cosh_{\Delta x}(\frac{\tanh b\Delta x}{\Delta x}, x)$<br>or<br>$e^{ax}\cosh bx = e_{\Delta x}(\frac{e^{a\Delta x}\cosh b\Delta x - 1}{\Delta x}, x)\cosh_{\Delta x}(\frac{\tanh b\Delta x}{\Delta x}, x)$ | Identity |
| **121** | $\sin_{\Delta x}(b,x-\Delta x) = \cos_{\Delta x}(b,-\Delta x)\sin_{\Delta x}(b,x) + \sin_{\Delta x}(b,-\Delta x)\cos_{\Delta x}(b,x)$<br>or<br>$\sin_{\Delta x}(b,x-\Delta x) = \frac{1}{1+[b\Delta x]^2}\sin_{\Delta x}(b,x) - \frac{b\Delta x}{1+[b\Delta x]^2}\cos_{\Delta x}(b,x)$ | Identity |
| **122** | $\cos_{\Delta x}(b,x-\Delta x) = -\sin_{\Delta x}(b,-\Delta x)\sin_{\Delta x}(b,x) + \cos_{\Delta x}(b,-\Delta x)\cos_{\Delta x}(b,x)$<br>or<br>$\cos_{\Delta x}(b,x-\Delta x) = \frac{b\Delta x}{1+[b\Delta x]^2}\sin_{\Delta x}(b,x) + \frac{1}{1+[b\Delta x]^2}\cos_{\Delta x}(b,x)$ | Identity |
| **123** | $\sin_{\Delta x}(b,x-n\Delta x) = \cos_{\Delta x}(b,-n\Delta x)\sin_{\Delta x}(b,x) + \sin_{\Delta x}(b,-n\Delta x)\cos_{\Delta x}(b,x)$<br>$n = 1,2,3,\ldots$ | Identity |
| **124** | $\cos_{\Delta x}(b,x-n\Delta x) = -\sin_{\Delta x}(b,-n\Delta x)\sin_{\Delta x}(b,x) + \cos_{\Delta x}(b,-n\Delta x)\cos_{\Delta x}(b,x)$<br>$n = 1,2,3,\ldots$ | Identity |
| **125** | $\sin_{\Delta x}(b,x) = \sin_{\Delta x}(b,x-\Delta x) + [b\Delta x]\cos_{\Delta x}(b,x-\Delta x)$ | Identity |
| **126** | $\cos_{\Delta x}(b,x) = -[b\Delta x]\sin_{\Delta x}(b,x-\Delta x) + \cos_{\Delta x}(b,x-\Delta x)$ | Identity |
| **127** | $\sin_{\Delta x}(b,x) = (1+[b\Delta x]^2)^n\cos_{\Delta x}(b,-n\Delta x)\sin_{\Delta x}(b,x-n\Delta x) - (1+[b\Delta x]^2)^n\sin_{\Delta x}(b,-n\Delta x)\cos_{\Delta x}(b,x-n\Delta x)$<br>$n = 1,2,3,\ldots$ | Identity |
| **128** | $\cos_{\Delta x}(b,x) = (1+[b\Delta x]^2)^n\sin_{\Delta x}(b,-n\Delta x)\sin_{\Delta x}(b,x-n\Delta x) + (1+[b\Delta x]^2)^n\cos_{\Delta x}(b,-n\Delta x)\cos_{\Delta x}(b,x-n\Delta x)$<br>$n = 1,2,3,\ldots$ | Identity |
| **129** | $\ln[e_{\Delta x}(a,x)] = [\frac{\ln(1+a\Delta x)}{\Delta x}]\,x$ | |
| **130** | $e^{jb\frac{x}{\Delta x}} = \cos b\frac{x}{\Delta x} + j\sin b\frac{x}{\Delta x}$<br>$b = \frac{2\pi}{m}$ , $m$ = integer , $m \neq 0$<br>$x = 0, \Delta x, 2\Delta x, 3\Delta x, \ldots$ | |
| **131** | $e^{j\frac{\pi x}{2\Delta x}} = [\,j\,]^{\frac{x}{\Delta x}}$<br>where<br>$x,\Delta x$ = real or complex values | $\ln(j) = \frac{\pi}{2}j$ |
| **132** | $e^{-j\frac{\pi x}{2\Delta x}} = [-j\,]^{\frac{x}{\Delta x}}$<br>where<br>$x,\Delta x$ = real or complex values | $\ln(-j) = -\frac{\pi}{2}j$ |
| **133** | $\sin\frac{\pi x}{2\Delta x} = \frac{[\,j\,]^{\frac{x}{\Delta x}} - [-j\,]^{\frac{x}{\Delta x}}}{2j}$<br>where<br>$x,\Delta x$ = real or complex values | $\ln(j) = \frac{\pi}{2}j$<br>$\ln(-j) = -\frac{\pi}{2}j$ |
| **134** | $\cos\frac{\pi x}{2\Delta x} = \frac{[\,j\,]^{\frac{x}{\Delta x}} + [-j\,]^{\frac{x}{\Delta x}}}{2}$<br>where<br>$x,\Delta x$ = real or complex values | $\ln(j) = \frac{\pi}{2}j$<br>$\ln(-j) = -\frac{\pi}{2}j$ |
| **135** | $x! = [x]_1^x = [1]_{-1}^x = x(x-1)(x-2)(x-3)\ldots(2)(1)$<br>$x = 0,1,2,3,\ldots$ | x factorial |
| **136** | ${}_nP_x = [n]_1^x = n(n-1)(n-2)(n-3)\ldots(n-x+1)$<br>$x = 0,1,2,3,\ldots n$ | Permutations<br>n things in groups of x |
| **137** | ${}_nC_x = \frac{[n]_1^x}{[x]_1^x} = \frac{[n-x+1]_{-1}^x}{[1]_{-1}^x} = \frac{n(n-1)(n-2)(n-3)\ldots(n-x+1)}{x!} = \frac{n!}{x!(n-x)!}$<br>$x = 0,1,2,3,\ldots n$ | Combinations<br>n things in groups of x |
| **138** | $\ln e_{\Delta x}(a,x) = [\frac{\ln(1+a\Delta x)}{\Delta x}]x$ | |

# TABLE 5a

## Some Commonly Used Interval Calculus/Calculus Function Identities

| # | Calculus Function<br>$x = 0, \Delta x, 2\Delta x, 3\Delta x, \ldots$ | Interval Calculus Function<br>$x = 0, \Delta x, 2\Delta x, 3\Delta x, \ldots$ |
|---|---|---|
| **1** | $e^{ax}$ | $e_{\Delta x}(\frac{e^{a\Delta x}-1}{\Delta x}, x)$<br>or<br>$[1 + (\frac{e^{a\Delta x}-1}{\Delta x})\Delta x]^{\frac{x}{\Delta x}}$ |
| **2** | $\sin bx$ | $[\cos b\Delta x]^{\frac{x}{\Delta x}}\sin_{\Delta x}(\frac{\tan b\Delta x}{\Delta x}, x)$<br>or<br>$e_{\Delta x}(\frac{\cos b\Delta x - 1}{\Delta x}, x)\sin_{\Delta x}(\frac{\tan b\Delta x}{\Delta x}, x)$<br>$\cos b\Delta x \neq 0$ |
| **3** | $\cos bx$ | $[\cos b\Delta x]^{\frac{x}{\Delta x}}\cos_{\Delta x}(\frac{\tan b\Delta x}{\Delta x}, x)$<br>or<br>$e_{\Delta x}(\frac{\cos b\Delta x - 1}{\Delta x}, x)\cos_{\Delta x}(\frac{\tan b\Delta x}{\Delta x}, x)$<br>$\cos b\Delta x \neq 0$ |
| **4** | $e^{ax}\sin bx$ | $[e^{a\Delta x}\cos b\Delta x]^{\frac{x}{\Delta x}}\sin_{\Delta x}(\frac{\tan b\Delta x}{\Delta x}, x)$<br>or<br>$e_{\Delta x}(\frac{e^{a\Delta x}\cos b\Delta x - 1}{\Delta x}, x)\sin_{\Delta x}(\frac{\tan b\Delta x}{\Delta x}, x)$<br>$\cos b\Delta x \neq 0$ |
| **5** | $e^{ax}\cos bx$ | $[e^{a\Delta x}\cos b\Delta x]^{\frac{x}{\Delta x}}\cos_{\Delta x}(\frac{\tan b\Delta x}{\Delta x}, x)$<br>or<br>$e_{\Delta x}(\frac{e^{a\Delta x}\cos b\Delta x - 1}{\Delta x}, x)\cos_{\Delta x}(\frac{\tan b\Delta x}{\Delta x}, x)$<br>$\cos b\Delta x \neq 0$ |

| # | Calculus Function<br>x = 0, Δx, 2Δx, 3Δx, … | Interval Calculus Function<br>x = 0, Δx, 2Δx, 3Δx, … |
|---|---|---|
| | | |
| **6** | $\sinh bx$ | $[\cosh b\Delta x]^{\frac{x}{\Delta x}}\sinh_{\Delta x}(\frac{\tanh b\Delta x}{\Delta x}, x)$<br>or<br>$e_{\Delta x}(\frac{\cosh b\Delta x - 1}{\Delta x}, x)\sinh_{\Delta x}(\frac{\tanh b\Delta x}{\Delta x}, x)$ |
| **7** | $\cosh bx$ | $[\cosh b\Delta x]^{\frac{x}{\Delta x}}\cosh_{\Delta x}(\frac{\tanh b\Delta x}{\Delta x}, x)$<br>or<br>$e_{\Delta x}(\frac{\cosh b\Delta x - 1}{\Delta x}, x)\cosh_{\Delta x}(\frac{\tanh b\Delta x}{\Delta x}, x)$ |
| **8** | $e^{ax}\sinh bx$ | $[e^{a\Delta x}\cosh b\Delta x]^{\frac{x}{\Delta x}}\sinh_{\Delta x}(\frac{\tanh b\Delta x}{\Delta x}, x)$<br>or<br>$e_{\Delta x}(\frac{e^{a\Delta x}\cosh b\Delta x - 1}{\Delta x}, x)\sinh_{\Delta x}(\frac{\tanh b\Delta x}{\Delta x}, x)$ |
| **9** | $e^{ax}\cosh bx$ | $[e^{a\Delta x}\cosh b\Delta x]^{\frac{x}{\Delta x}}\cosh_{\Delta x}(\frac{\tanh b\Delta x}{\Delta x}, x)$<br>or<br>$e_{\Delta x}(\frac{e^{a\Delta x}\cosh b\Delta x - 1}{\Delta x}, t)\cosh_{\Delta x}(\frac{\tanh b\Delta x}{\Delta x}, x)$ |
| **10** | $e^{ax}$<br>$a = \frac{\ln(1+\alpha\Delta x)}{\Delta x}$<br>$\alpha = \frac{e^{a\Delta x} - 1}{\Delta x}$ | $e_{\Delta x}(\alpha,x)$<br>or<br>$(1+\alpha\Delta x)^{\frac{x}{\Delta x}}$ |
| **11** | $(\sec\beta)^{\frac{x}{\Delta x}} \sin\beta\frac{x}{\Delta x}$<br>$\beta = \tan^{-1}b\Delta x$ | $\sin_{\Delta x}(b,x)$ |
| **12** | $(\sec\beta)^{\frac{x}{\Delta x}} \cos\beta\frac{x}{\Delta x}$<br>$\beta = \tan^{-1}b\Delta x$ | $\cos_{\Delta x}(b,x)$ |

| # | **Calculus Function**<br>**$x = 0, \Delta x, 2\Delta x, 3\Delta x, \ldots$** | **Interval Calculus Function**<br>**$x = 0, \Delta x, 2\Delta x, 3\Delta x, \ldots$** |
|---|---|---|
| **13** | $[(1+a\Delta x)^2+(b\Delta x)^2]^{\frac{x}{2\Delta x}} \sin\beta \frac{x}{\Delta x}$<br>or<br>$e_{\Delta x}(a,x)(\sec\beta)^{\frac{x}{\Delta x}} \sin\beta\frac{x}{\Delta x}$<br>$\beta = \begin{cases} \tan^{-1}\left\|\frac{b\Delta x}{1+a\Delta x}\right\| & \text{for } 1+a\Delta x > 0 \quad b\Delta x \ge 0 \\ -\tan^{-1}\left\|\frac{b\Delta x}{1+a\Delta x}\right\| & \text{for } 1+a\Delta x > 0 \quad b\Delta x < 0 \\ \pi - \tan^{-1}\left\|\frac{b\Delta x}{1+a\Delta x}\right\| & \text{for } 1+a\Delta x < 0 \quad b\Delta x \ge 0 \\ -\pi + \tan^{-1}\left\|\frac{b\Delta x}{1+a\Delta x}\right\| & \text{for } 1+a\Delta x < 0 \quad b\Delta x < 0 \end{cases}$<br>$1+a\Delta x \neq 0$<br>$0 \le \tan^{-1}\left\|\frac{b\Delta x}{1+a\Delta x}\right\| < \frac{\pi}{2}$ | $e_{\Delta t}(a,x)\sin_{\Delta x}(\frac{b}{1+a\Delta x}, x)$ |
| **14** | $[(1+a\Delta x)^2+(b\Delta x)^2]^{\frac{x}{2\Delta x}} \cos\beta\frac{x}{\Delta x}$<br>or<br>$e_{\Delta x}(a,x)(\sec\beta)^{\frac{x}{\Delta x}} \cos\beta\frac{x}{\Delta x}$<br>$\beta = \begin{cases} \tan^{-1}\left\|\frac{b\Delta x}{1+a\Delta x}\right\| & \text{for } 1+a\Delta x > 0 \quad b\Delta x \ge 0 \\ -\tan^{-1}\left\|\frac{b\Delta x}{1+a\Delta x}\right\| & \text{for } 1+a\Delta x > 0 \quad b\Delta x < 0 \\ \pi - \tan^{-1}\left\|\frac{b\Delta x}{1+a\Delta x}\right\| & \text{for } 1+a\Delta x < 0 \quad b\Delta x \ge 0 \\ -\pi + \tan^{-1}\left\|\frac{b\Delta x}{1+a\Delta x}\right\| & \text{for } 1+a\Delta x < 0 \quad b\Delta x < 0 \end{cases}$<br>$1+a\Delta x \neq 0$<br>$0 \le \tan^{-1}\left\|\frac{b\Delta x}{1+a\Delta x}\right\| < \frac{\pi}{2}$ | $e_{\Delta x}(a,x)\cos_{\Delta x}(\frac{b}{1+a\Delta x}, x)$ |

| # | Calculus Function<br><br>**x = 0, Δx, 2Δx, 3Δx, …** | Interval Calculus Function<br><br>**x = 0, Δx, 2Δx, 3Δx, …** |
|---|---|---|
| **15** | $(1-[b\Delta x]^2)^{\frac{x}{2\Delta x}} \sinh\beta\frac{x}{\Delta x}$<br>$\beta = \tanh^{-1}b\Delta x$ | $\sinh_{\Delta x}(b,x)$ |
| **16** | $(1-[b\Delta x]^2)^{\frac{x}{2\Delta x}} \cosh\beta\frac{x}{\Delta x}$<br>$\beta = \tanh^{-1}b\Delta x$ | $\cosh_{\Delta x}(b,x)$ |

# TABLE 6

## Some Interval Calculus Difference in Differential Form, and Integral Equations

| | Operation | Description |
|---|---|---|
| | **Difference Equations in Differential Form** | |
| **1** | $H_{\Delta x}u(x) = u(x+\Delta x)$ | $+\Delta x$ Increment |
| **2** | $H_{-\Delta x}u(x) = u(x-\Delta x)$ | $-\Delta x$ Decrement |
| **3** | $\Delta u(x) = u(x+\Delta x) – u(x)$ <br> or <br> $\Delta u(x) = \Delta x \left[\frac{u(x+\Delta x) - u(x)}{\Delta x}\right] = \Delta x\, D_{\Delta x}u(x)$ | Difference |
| **4** | $D_{\Delta x}u(x) = \left[ \frac{H_{\Delta x} -1}{\Delta x} \right]u(x) = \frac{u(x+\Delta x) - u(x)}{\Delta x}$ , Discrete differentiation | Discrete Calculus $+\Delta x$ Derivative |
| **5** | $D_{-\Delta x}u(x) = \left[ \frac{H_{-\Delta x} -1}{-\Delta x} \right]u(x) = \frac{u(x-\Delta x) - u(x)}{-\Delta x}$ | Discrete Calculus $-\Delta x$ Derivative |
| **6** | $D_{-\Delta x}f(x) = D_{\Delta x}H_{-\Delta x}f(x)$ | Equality |
| **7** | $D_{\Delta x}f(x) = D_{-\Delta x}H_{\Delta x}f(x)$ | Equality |
| **8** | $D_{\Delta x}H_{\Delta x}f(x) = H_{\Delta x}D_{\Delta x}f(x)$ | Operator Commutative Property |
| **9** | $D_{-\Delta x}H_{-\Delta x}f(x) = H_{-\Delta x}D_{-\Delta x}f(x)$ | Operator Commutative Property |
| **10** | $D_{\Delta x}[u(x)v(x)] = v(x)D_{\Delta x}u(x) + D_{\Delta x}v(x)u(x+\Delta x)$ <br> or <br> $D_{\Delta x}[u(x)v(x)] = v(x)D_{\Delta x}u(x) + u(x)D_{\Delta x}v(x) + D_{\Delta x}u(x)D_{\Delta x}v(x)\Delta x$ | Derivative of a Product |
| **11** | $D_{\Delta x}\left[\frac{u(x)}{v(x)}\right] = \frac{v(x)D_{\Delta x}u(x) - u(x)D_{\Delta x}v(x)}{v(x)H_{\Delta x}v(x)}$ | Derivative of a Division |

| | Operation | Description |
|---|---|---|
| **12** | Discrete Function Chain Rule<br>$D_{\Delta x}\left[F(v)\,\middle|_{v=g_{\Delta x}(x)}\right] = \left[D_{\Delta v}F(v)\right]\Big|_{\substack{v=g_{\Delta x}(x)\\ \Delta v = \Delta x D_{\Delta x} g_{\Delta x}(x)}} \quad \left[D_{\Delta x}g_{\Delta x}(x)\right]$<br><br>where<br>$F(v)$ = function of v<br>$g_{\Delta x}(x)$ = discrete Interval Calculus function of x<br>$v = g_{\Delta x}(x)$<br>$\Delta x,\Delta v$ = interval increments | Chain Rule |
| **13** | $D_{\Delta x}\, f(x+m\Delta x) = D_{\Delta v} f(v)\Big|_{\substack{v = x+m\Delta x\\ \Delta v = \Delta x\\ m = 0,1,2,3,\ \ldots}}$ | |
| **14** | $\Delta u(x) = D_{\Delta x}u(x)\Delta x$ | |
| **15** | $D_{\Delta x}ku(x) = kD_{\Delta x}u(x)$ | |
| **16** | $D_{\Delta x}k = 0$ | |
| **17** | $D_{\Delta x}x = 1$ | |
| **18** | $D_{\Delta x}[x(x-\Delta x)] = 2x$ | |
| **19** | $D_{\Delta x}[x(x-\Delta x)(x-2\Delta x)] = 3x(x-\Delta x)$ | |
| **20** | $D_{\Delta x}[x(x-\Delta x)(x-2\Delta x)(x-3\Delta x)] = 4x(x-\Delta x)(x-2\Delta x)$ | |
| **21** | $D_{\Delta x}[x(x+\Delta x)] = 2(x+\Delta x)$ | |
| **22** | $D_{\Delta x}[x(x+\Delta x)(x+2\Delta x)] = 3(x+\Delta x)(x+2\Delta x)$ | |
| **23** | $D_{\Delta x}[x(x+\Delta x)(x+2\Delta x)(x+3\Delta x)] = 4(x+\Delta x)(x+2\Delta x)(x+3\Delta x)$ | |
| **24** | $D_{\Delta x}\left[\frac{1}{x}\right] = -\frac{1}{(x+\Delta x)x}$ | |

| | Operation | Description |
|---|---|---|
| **25** | $D_{\Delta x}\left[\frac{1}{x(x-\Delta x)}\right] = -\frac{2}{(x+\Delta x)x(x-\Delta x)}$ | |
| **26** | $D_{\Delta x}\left[\frac{1}{x(x-\Delta x)(x-2\Delta x)}\right] = -\frac{3}{(x+\Delta x)x(x-\Delta x)(x-2\Delta x)}$ | |
| **27** | $D_{\Delta x}\left[\frac{1}{x(x+\Delta x)}\right] = -\frac{2}{x(x+\Delta x)(x+2\Delta x)}$ | |
| **28** | $D_{\Delta x}\left[\frac{1}{x(x+\Delta x)(x+2\Delta x)}\right] = -\frac{3}{x(x+\Delta x)(x+\Delta x)(x+3\Delta x)}$ | |
| **29** | $D_{\Delta x}[\mathit{x}]^{n}_{\Delta x} = \begin{cases} 0 & \text{for } n = 0 \\ n[\mathit{x}]^{n-1}_{\Delta x} & \text{for } n = 1 \quad 2 \quad 3 \ldots \\ n[\mathit{x}+\Delta x]^{n-1}_{\Delta x} & \text{for } n = -1 \quad -2 \quad -3 \ldots \end{cases}$<br>where<br>$[\mathit{x}]^{0}_{\Delta x} = 1$<br>$[\mathit{x}]^{n}_{\Delta x} = \prod_{m=1}^{n}(\mathit{x}-(m-1)\Delta x)$ , n = 1,2,3…<br>$[\mathit{x}]^{-m}_{\Delta x} = \frac{1}{[\mathit{x}]^{m}_{\Delta x}}$ , m = integer<br>$\mathit{x} = x + p\Delta x$ , p = integer<br>$x = x_0 + r\Delta x$ , r = integers | General Equation |
| **29a** | $D_{\Delta x}\left[\prod_{n=a}^{b}(x-n\Delta x)\right] = (b-a+1)\prod_{m=a}^{b-1}(x-m\Delta x)$<br>n = a, a+1, a+2,…, b-1, b<br>m = a, a+1, a+2,…, b-2, b-1<br>a,b = integers<br>$b \geq a$<br>$\prod_{n=a}^{a-1}(x-n\Delta x) = 1$<br>$x = x_0 + r\Delta x$, r = integers<br>b-a+1 = the order of the polynomial being differentiated | General Equation |

| | Operation | Description |
|---|---|---|
| **29b** | $D_{\Delta x}[\prod_{n=a}^{b}(x+n\Delta x)] = (b-a+1)\prod_{m=a+1}^{b}(x+m\Delta x)$<br>$n = a, a+1, a+2,\ldots, b-1, b$<br>$m = a+1, a+2, a+3,\ldots, b-1, b$<br>a,b = integers<br>$b \geq a$<br>$\prod_{n=a+1}^{a}(x-n\Delta x) = 1$<br>$x = x_0 + r\Delta x$, r = integers<br>b-a+1 = the order of the polynomial being differentiated | General Equation |
| **30** | $D_{\Delta x}\left[\dfrac{1}{\prod_{n=a}^{b}(x-n\Delta x)}\right] = -\left[\dfrac{(b-a+1)}{\prod_{m=a-1}^{b}(x-m\Delta x)}\right]$<br>$n = a, a+1, a+2,\ldots, b-1,b$<br>$m = a-1, a, a+1, \ldots, b-1, b$<br>a,b = integers<br>$b \geq a$<br>$x = x_0 + r\Delta x$, r = integers<br>b-a+1 = the order of the polynomial being differentiated | General Equation |
| **30a** | $D_{\Delta x}\left[\dfrac{1}{\prod_{n=a}^{b}(x+n\Delta x)}\right] = -\left[\dfrac{(b-a+1)}{\prod_{m=a}^{b+1}(x+m\Delta x)}\right]$<br>$n = a, a+1, a+2,\ldots, b-1,b$<br>$m = a, a+1, a+2,\ldots, b, b+1$<br>a,b = integers<br>$b \geq a$<br>$x = x_0 + r\Delta x$, r = integers<br>b-a+1 = the order of the polynomial being differentiated | General Equation |
| **31** | $D_{\Delta x}\,\mathrm{lnd}(n,\Delta x,x) = \pm\dfrac{1}{x^n}$ ,   + for n = 1,   – for $n \neq 1$<br>A computer program, LNDX, is available to calculate the function, lnd(n,Δx,x) | $\ln_{\Delta x}x \equiv \mathrm{lnd}(1,\Delta x,x)$ |
| **32** | $D_{\Delta x}\ln_{\Delta x}x = \dfrac{1}{x}$ | $\ln_{\Delta x}x \equiv \mathrm{lnd}(1,\Delta x,x)$ |
| **33** | $D_{\Delta x}f(g(x)) = D_{\Delta p}f(p)D_{\Delta x}g(x)$<br>$p = g(x)$ | Derivative of a function |

| | Operation | Description |
|---|---|---|
| **34** | $D_{2\Delta x} = D_{\Delta x} + \frac{\Delta x}{2} D_{\Delta x}^2$<br>or<br>$s_{2\Delta x} = s + \frac{\Delta x}{2} s^2$<br>where<br>$D_{\Delta x}u(x) = \frac{u(x+\Delta x) - u(x)}{\Delta x}$<br>$D_{2\Delta x}u(x) = \frac{u(x+2\Delta x) - u(x)}{2\Delta x}$<br>$s = s_{\Delta x} = D_{\Delta x}$ , discrete derivative operator | Discrete derivative operator equality |
| **35** | $D_{3\Delta x} = D_{\Delta x} + \Delta x D_{\Delta x}^2 + \frac{\Delta x^2}{3} D_{\Delta x}^3$<br>or<br>$s_{3\Delta x} = s + \Delta x s^2 + \frac{\Delta x^2}{3} s^3$<br>where<br>$D_{3\Delta x}u(x) = \frac{u(x+3\Delta x) - u(x)}{3\Delta x}$<br>$s = s_{\Delta x} = D_{\Delta x}$ , discrete derivative operator | Discrete derivative operator equality |
| **36** | $D_{4\Delta t} = D_{\Delta t} + \frac{3}{2}\Delta t D_{\Delta t}^2 + \Delta t^2 D_{\Delta t}^3 + \frac{1}{4}\Delta t^3 D_{\Delta t}^4$<br>or<br>$s_{4\Delta t} = D_{\Delta t} + \frac{3}{2}\Delta t s^2 + \Delta t^2 s^3 + \frac{1}{4}\Delta t^3 s^4$<br>where<br>$D_{4\Delta x}u(x) = \frac{u(x+4\Delta x) - u(x)}{4\Delta x}$<br>$s = s_{\Delta x} = D_{\Delta x}$ , discrete derivative operator | Discrete derivative operator equality |
| **37** | $D_{m\Delta x} = D_{\Delta x} + \sum_{n=2}^{\infty} \frac{\Delta x^{n-1}}{n!} \prod_{p=1}^{n-1}(m-p) D_{\Delta x}^n$ , $m = 2,3,4,\ldots$<br>or<br>$D_{m\Delta x} = D_{\Delta x} + \frac{(m-1)}{2!}\Delta x D_{\Delta x}^2 + \frac{(m-1)(m-2)}{3!}\Delta x^2 D_{\Delta x}^3 + \frac{(m-1)(m-2)(m-3)}{4!}\Delta x^3 D_{\Delta x}^4 + \ldots$ ,<br>$m = 2,3,4,\ldots$<br>or<br>$s_{m\Delta x} = s + \frac{(m-1)}{2!}\Delta x s^2 + \frac{(m-1)(m-2)}{3!}\Delta x^2 s^3 + \frac{(m-1)(m-2)(m-3)}{4!}\Delta x^3 s^4 + \ldots$ , $m = 2,3,4,\ldots$ | Discrete derivative operator equality |

| | Operation | Description |
|---|---|---|
| | where<br>$D_{m\Delta x}u(x) = \dfrac{u(x+m\Delta x) - u(x)}{m\Delta x}$<br>$s = s_{\Delta x} = D_{\Delta x}$ , discrete derivative operator | |
| **38** | $D_{2\Delta x}{}^2 = D_{\Delta x}{}^2 + \Delta x D_{\Delta x}{}^3 + \dfrac{\Delta x^2}{4} D_{\Delta x}{}^4$<br>or<br>$s_{2\Delta x}{}^2 = s^2 + \Delta x s^3 + \dfrac{\Delta x^2}{4} s^4$<br>where<br>$D_{2\Delta x}{}^2 u(x) = \dfrac{u(x+4\Delta x) - 2u(x+2\Delta x) + u(x)}{(2\Delta x)^2}$<br>$s = s_{\Delta x} = D_{\Delta x}$ , discrete derivative operator | Discrete derivative operator equality |
| **39** | $D_q D_r = D_r D_q$<br>where<br>$q,r = n\Delta x$ , $n = 0,1,2,3,\ldots$<br>$D = \dfrac{d}{dx} = D_0$<br>$D_a u(x) = \dfrac{u(x+a) - u(x)}{a}$ , $a = q,r$ , discrete derivative of u(x) | Discrete derivative operator equality |
| **40** | $D_q D_r D_s = D_q D_s D_r = D_r D_s D_q = D_r D_q D_s = D_s D_r D_q = D_s D_q D_r$<br>where<br>$q,r,s = n\Delta x$ , $n = 0,1,2,3,\ldots$<br>$D = \dfrac{d}{dx} = D_0$<br>$D_a u(x) = \dfrac{u(x+a) - u(x)}{a}$ , $a = q,r,s$ , discrete derivative of u(x) | Discrete derivative operator equality |
| **41** | $f(x+\Delta x) = f(x) + D_{\Delta x} f(x) \Delta x$ | |
| **42** | $D_{\Delta x} \sin\pi\dfrac{x}{\Delta x} = -\dfrac{2}{\Delta x}\sin\pi\dfrac{x}{\Delta x}$ | |
| **43** | $D_{\Delta x} \cos\pi\dfrac{x}{\Delta x} = -\dfrac{2}{\Delta x}\cos\pi\dfrac{x}{\Delta x}$ | |
| **44** | $D_{\Delta x} \sin\dfrac{\pi x}{2\Delta x} = \dfrac{1}{\Delta x}\left( \cos\dfrac{\pi x}{2\Delta x} - \sin\dfrac{\pi x}{2\Delta x} \right)$ | |
| **45** | $D_{\Delta x} \cos\dfrac{\pi x}{2\Delta x} = -\dfrac{1}{\Delta x}\left( \sin\dfrac{\pi x}{2\Delta x} + \cos\dfrac{\pi x}{2\Delta x} \right)$ | |

| | Operation | Description |
|---|---|---|
| **46** | $D_{\Delta x}e_{\Delta x}x = e_{\Delta x}x$ | $e_{\Delta x}x = (1+\Delta x)^{\frac{x}{\Delta x}}$ |
| **47** | $D_{\Delta x}\sin_{\Delta x}x = \cos_{\Delta x}x$ | |
| **48** | $D_{\Delta x}\cos_{\Delta x}x = -\sin_{\Delta x}x$ | |
| **49** | $D_{\Delta x}e_{\Delta x}(a,x) = ae_{\Delta x}(a,x)$<br>or<br>$D_{\Delta x}[1+a\Delta x]^{\frac{x}{\Delta x}} = a[1+a\Delta x]^{\frac{x}{\Delta x}}$ | $e_{\Delta x}(a,x) = (1+a\Delta x)^{\frac{x}{\Delta x}}$ |
| **50** | $D_{\Delta x}e_{\Delta x}(a,x+m\Delta x) = ae_{\Delta x}(a,x+m\Delta x)$<br>$m = 0, 1, 2, 3, \ldots$ | |
| **51** | $D_{2\Delta x}e_{\Delta x}(a,x) = (a + \frac{\Delta x}{2}a^2)\, e_{\Delta x}(a,x)$<br>$x = 0, \Delta x, 2\Delta x, 3\Delta x, \ldots$<br>$m = 1, 2, 3, \ldots$ | |
| **52** | $D_{m\Delta x}e_{m\Delta x}(\frac{a}{m},mx) = [\frac{(1+a\Delta x)^m-1}{m\Delta x}]e_{\Delta x}(a,x)$<br>or<br>$D_{m\Delta x}e_{m\Delta x}(\frac{a}{m},mx) = [\frac{(1+a\Delta x)^m-1}{m\Delta x}]e_{m\Delta x}(\frac{a}{m},mx)$<br>$x = 0, \Delta x, 2\Delta x, 3\Delta x, \ldots$<br>$m = 1, 2, 3, \ldots$ | |
| **53** | $D_{\Delta x}e_{\Delta x}(a,2x) = [2a+a^2\Delta x]e_{\Delta x}(a,2x)$<br>or<br>$D_{\Delta x}[1+a\Delta x]^{\frac{2x}{\Delta x}} = [2a+a^2\Delta x][1+a\Delta x]^{\frac{2x}{\Delta x}}$ | |
| **54** | $D_{\Delta x}e_{\Delta x}(a,mx) = [\frac{(1+a\Delta t)^m-1}{\Delta t}]e_{\Delta x}(a,mx)$ , $m = 0,1,2,3,\ldots$<br>or<br>$D_{\Delta x}[1+a\Delta x]^{\frac{mx}{\Delta x}} = [\frac{(1+a\Delta t)^m-1}{\Delta t}][1+a\Delta x]^{\frac{mx}{\Delta x}}$ , $m = 0,1,2,3,\ldots$ | |

| | Operation | Description |
|---|---|---|
| **55** | $D_{\Delta x}\sin_{\Delta x}(b,x) = b\cos_{\Delta x}(b,x)$ | |
| **56** | $D_{\Delta x}\sin_{\Delta x}(b,x+\Delta x) = b\cos_{\Delta x}(b,x+\Delta x)$ | |
| **57** | $D_{2\Delta x}\sin_{\Delta x}(b,x) = b\cos_{\Delta x}(b,x) - \frac{b^2\Delta t}{2}\sin_{\Delta x}(b,t)$ | |
| **58** | $D_{\Delta x}\sin_{\Delta x}{}^2(a,x) = a[\sin_{\Delta x}(a,x+\Delta x) + \sin_{\Delta x}(a,x)]\cos_{\Delta x}(a,x)$<br>or<br>$D_{\Delta x}\sin_{\Delta x}{}^2(a,x) = a[2\sin_{\Delta x}(a,x) + a\Delta x\cos_{\Delta x}(a,x)]\cos_{\Delta x}(a,x)$ | |
| **59** | $D_{\Delta x}\cos_{\Delta x}(b,x) = - b\sin_{\Delta x}(b,x)$ | |
| **60** | $D_{\Delta x}\cos_{\Delta x}(b,x+\Delta x) = - b\sin_{\Delta x}(b,x+\Delta x)$ | |
| **61** | $D_{2\Delta x}\cos_{\Delta x}(b,x) = -b\sin_{\Delta x}(b,x) - \frac{b^2\Delta t}{2}\cos_{\Delta t}(b,t)$ | |
| **62** | $D_{\Delta x}\cos_{\Delta x}{}^2(a,x) = -a[\cos_{\Delta x}(a,x+\Delta x) + \cos_{\Delta x}(a,x)]\sin_{\Delta x}(a,x)$<br>or<br>$D_{\Delta x}\cos_{\Delta x}{}^2(a,x) = -a[2\cos_{\Delta x}(a,x) - a\Delta x\sin_{\Delta x}(a,x)]\sin_{\Delta x}(a,x)$ | |
| **63** | $D_{\Delta x}\sinh_{\Delta x}x = \cosh_{\Delta x}x$ | |
| **64** | $D_{\Delta x}\cosh_{\Delta x}x = \sinh_{\Delta x}x$ | |
| **65** | $D_{\Delta x}\sinh_{\Delta x}(a,x) = a\cosh_{\Delta x}(a,x)$ | |
| **66** | $D_{\Delta x}\cosh_{\Delta x}(a,x) = a\sinh_{\Delta x}(a,x)$ | |
| **67** | $D_{\Delta x}e_{\Delta x}(a,-x) = -ae_{\Delta x}(a,-x-\Delta x)$ | |
| **68** | $D_{\Delta x}(1+a\Delta x)^{-\frac{x}{\Delta x}} = -a(1+a\Delta x)^{-\frac{x+\Delta x}{\Delta x}}$ | |
| **69** | $D_{\Delta x}k^{-\frac{x}{\Delta x}} = \left(\frac{k^{-1}-1}{\Delta x}\right)k^{-\frac{x}{\Delta x}}$ , k = constant | |
| **70** | $D_{\Delta x}e^{\frac{ax}{\Delta x}} = \left(\frac{e^a-1}{\Delta x}\right)e^{\frac{ax}{\Delta x}}$ | |

| | Operation | Description |
|---|---|---|
| **71** | $D_{\Delta x}\sin\frac{ax}{\Delta x} = [\frac{\cos a-1}{\Delta x}]\sin\frac{ax}{\Delta x} + [\frac{\sin a}{\Delta x}]\cos\frac{ax}{\Delta x}$ | |
| **72** | $D_{\Delta x}\cos\frac{ax}{\Delta x} = [\frac{\cos a-1}{\Delta x}]\cos\frac{ax}{\Delta x} - [\frac{\sin a}{\Delta x}]\sin\frac{ax}{\Delta x}$ | |
| **73** | $D_{\Delta x}\tan\frac{ax}{\Delta x} = [\frac{\sin a}{\Delta x}]\sec\frac{ax}{\Delta x}\sec a(\frac{x}{\Delta x} + 1)$ | |
| **74** | $D_{\Delta x}\csc\frac{ax}{\Delta x} = -\csc a(\frac{x}{\Delta x} + 1)[(\frac{\cos a-1}{\Delta x}) + (\frac{\sin a}{\Delta x})\cot an\frac{ax}{\Delta x}]$ | |
| **75** | $D_{\Delta x}\sec\frac{ax}{\Delta x} = \sec a(\frac{x}{\Delta x} + 1)[(\frac{\cos a-1}{\Delta x}) + (\frac{\sin a}{\Delta x})\tan\frac{ax}{\Delta x}]$ | |
| **76** | $D_{\Delta x}\cot an\frac{ax}{\Delta x} = -[\frac{\sin a}{\Delta x}]\csc\frac{ax}{\Delta x}\csc a(\frac{x}{\Delta x} + 1)$ | |
| | **<u>Integral Equations</u>** Indefinite integrals | |
| **77** | ${}_{\Delta x}\int(1)\Delta x = x + k$ , k = constant of integration | |
| **78** | ${}_{\Delta x}\int x\Delta x = \frac{x(x-\Delta x)}{2} + k$ | |
| **79** | ${}_{\Delta x}\int x(x-\Delta x)\Delta x = \frac{x(x-\Delta x)(x-2\Delta x)}{3} + k$ | |
| **80** | ${}_{\Delta x}\int x(x+\Delta x)\Delta x = \frac{(x-\Delta x)x\,(x+\Delta x)}{3} + k$ | |
| **81** | ${}_{\Delta x}\int\frac{1}{x(x-\Delta x)}\Delta x = -\frac{1}{x-\Delta x} + k$ | |
| **82** | ${}_{\Delta x}\int\frac{1}{x(x-\Delta x)(x-2\Delta x)}\Delta x = -\frac{1}{2(x-\Delta x)(x-2\Delta x)} + k$ | |
| **83** | ${}_{\Delta x}\int\frac{1}{x(x-\Delta x)(x-2\Delta x)(x-3\Delta x)}\Delta x = -\frac{1}{3(x-\Delta x)(x-2\Delta x)(x-3\Delta x)} + k$ | |
| **84** | ${}_{\Delta x}\int\frac{1}{x(x+\Delta x)}\Delta x = -\frac{1}{x} + k$ | |
| **85** | ${}_{\Delta x}\int\frac{1}{x(x+\Delta x)(x+2\Delta x)}\Delta x = -\frac{1}{2x(x+\Delta x)} + k$ | |
| **86** | ${}_{\Delta x}\int\frac{1}{x(x+\Delta x)(x+2\Delta x)(x+3\Delta x)}\Delta x = -\frac{1}{3x(x+\Delta x)(x+2\Delta x)} + k$ | |

| | Operation | Description |
|---|---|---|
| **87** | $\int_{\Delta x} [x]_{\Delta x}^{n} \Delta x = \begin{cases} \dfrac{[x]_{\Delta x}^{n+1}}{n+1} + k & \text{for } n = 0 \quad 1 \quad 2 \quad 3 \ldots \\ \dfrac{[x-\Delta x]_{\Delta x}^{n+1}}{n+1} + k & \text{for } n = -2 \quad -3 \ldots \end{cases}$ <br> where <br> $[x]_{\Delta x}^{0} = 1$ <br> $[x]_{\Delta x}^{n} = \prod_{m=1}^{n}(x-(m-1)\Delta x)$ , $n = 1,2,3\ldots$ <br> $[x]_{\Delta x}^{-m} = \dfrac{1}{[x]_{\Delta x}^{m}}$ , m = integer <br> $x = x + p\Delta x$ , p = integer <br> $x = x_0 + r\Delta x$ , r = integers | General Equation |
| **87a** | $\int_{\Delta x} \prod_{n=a}^{b}(x-n\Delta x)\Delta x = \dfrac{\prod_{m=a}^{b+1}(x-m\Delta x)}{b-a+2} + k$ <br> n = a, a+1, a+2,…, b-1, b <br> m = a, a+1, a+2,…, b, b+1 <br> a,b = integers <br> $b \geq a$ <br> $x = x_0 + r\Delta x$, r = integers <br> b-a+2 = one more than the order of the polynomial being integrated | General Equation |
| **87b** | $\int_{\Delta x} \prod_{n=a}^{b}(x+n\Delta x)\Delta x = \dfrac{\prod_{m=a-1}^{b}(x+m\Delta x)}{b-a+2} + k$ <br> n = a, a+1, a+2,…, b-1, b <br> m = a-1, a, a+1,…, b-1, b <br> a,b = integers <br> $b \geq a$ <br> $x = x_0 + r\Delta x$, r = integers <br> b-a+2 = one more than the order of the polynomial being integrated | General Equation |

| | Operation | Description |
|---|---|---|
| **88** | $$\int_{\Delta x} \frac{1}{\prod_{n=a}^{b}(x-n\Delta x)}\Delta x = -\frac{1}{(b-a)\prod_{m=a+1}^{b}(x-m\Delta x)} + k$$ <br> n = a, a+1, a+2,…, b-1, b <br> m = a+1, a+2, a+3,…, b-1, b <br> a,b = integers <br> b > a <br> $x = x_0 + r\Delta x$, r = integers <br> b-a = one less than the order of the denominator polynomial of the function being integrated | General Equation |
| **88a** | $$\int_{\Delta x} \frac{1}{\prod_{n=a}^{b}(x+n\Delta x)}\Delta x = -\frac{1}{(b-a)\prod_{m=a}^{b-1}(x+m\Delta x)} + k$$ <br> n = a, a+1, a+2,…, b-1, b <br><br> m = a, a+1, a+2,…, b-2, b-1 <br> a,b = integers <br> b > a <br> $x = x_0 + r\Delta x$, r = integers <br> b-a = one less than the order of the denominator polynomial of the function being integrated | General Equation |
| **89** | $$\int_{\Delta x} \frac{1}{x^n}\Delta x = \pm\,\mathrm{lnd}(n,\Delta x,x) + k\ , \quad +\text{ for } n = 1, \quad -\text{ for } n \neq 1$$ <br> A computer program, LNDX, is available to calculate the function, lnd(n,Δx,x) | |
| **90** | $$\int_{\Delta x} \frac{1}{x}\Delta x = \ln_{\Delta x}x + k \equiv \mathrm{lnd}(1,\Delta x,x) + k$$ | |
| **91** | $$\int_{\Delta x\, x_1}^{x_2} v(x)D_{\Delta x}u(x)\Delta x = u(x)v(x)\Big\vert_{x_1}^{x_2} - \int_{\Delta x\, x_1}^{x_2} D_{\Delta x}v(x)u(x+\Delta x)\Delta x$$ <br> or <br> $$\int_{\Delta x\, x_1}^{x_2} v(x)D_{\Delta x}u(x)\Delta x = u(x)v(x)\Big\vert_{x_1}^{x_2} - \int_{\Delta x\, x_1}^{x_2} u(x)D_{\Delta x}v(x)\Delta x - \Delta x \int_{\Delta x\, x_1}^{x_2} D_{\Delta x}u(x)D_{\Delta x}v(x)\Delta x$$ <br> where <br> $$D_{\Delta x}u(x) = \frac{u(x+\Delta x)-u(x)}{\Delta x}$$ | Interval Calculus Integration by Parts |

| | Operation | Description |
|---|---|---|
| **92** | ${}_{\Delta x}\int cf(x)\Delta x = c\ {}_{\Delta x}\int f(x)\Delta x$ , c = constant | |
| **93** | ${}_{\Delta x}\int e_{\Delta x}(x)\Delta x = e_{\Delta x}(x) + k$ | $e_{\Delta x}(x) = e_{\Delta x}(1,x)$ |
| **94** | ${}_{\Delta x}\int e_{\Delta x}(a,x)\Delta x = \frac{1}{a}\ e_{\Delta x}(a,x) + k$ | |
| **95** | ${}_{2\Delta x}\int e_{\Delta x}(a,x)\Delta x = \frac{1}{2a+a^2\Delta x}\ e_{\Delta x}(a,x) + k$<br>$x = 0, \Delta x, 2\Delta x, 3\Delta x, \ldots$ | |
| **96** | ${}_{m\Delta x}\int e_{\Delta x}(a,x)\Delta x = \frac{\Delta x}{[(1+a\Delta x)^m-1]}\ e_{\Delta x}(a,x) + k$<br>or<br>${}_{m\Delta x}\int e_{\Delta x}(a,x)\Delta x = \frac{\Delta x}{[(1+a\Delta x)^m-1]}\ e_{m\Delta x}(\frac{a}{m}, mx) + k$<br>$x = 0, \Delta x, 2\Delta x, 3\Delta x, \ldots$<br>$m = 1, 2, 3, \ldots$ | |
| **97** | ${}_{\Delta x}\int xe_{\Delta x}(a,x)\Delta x = \frac{e_{\Delta x}(a,x)}{a^2}[ax-1-a\Delta x] + k$ | |
| **98** | ${}_{\Delta x}\int e_{\Delta x}(a,-x)\Delta x = -\frac{1+a\Delta x}{a}e_{\Delta x}(a,-x) + k$ | |
| **99** | ${}_{\Delta x}\int (1+a\Delta x)^{-\frac{x}{\Delta x}}\ \Delta x = -\frac{1+a\Delta x}{a}(1+a\Delta x)^{-\frac{x}{\Delta x}} + k$ | |
| **100** | ${}_{\Delta x}\int \sin_{\Delta x}x\ \Delta x = -\cos_{\Delta x}x + k$ | |
| **101** | ${}_{\Delta x}\int \cos_{\Delta x}x\ \Delta x = \sin_{\Delta x}x + k$ | |
| **102** | ${}_{\Delta x}\int \sin_{\Delta x}(b,x)\Delta x = -\frac{1}{b}\cos_{\Delta x}(b,x) + k$ | |
| **103** | ${}_{2\Delta x}\int \sin_{\Delta x}(b,x)\Delta x = -\frac{2}{4b+b^3\Delta x^2}\cos_{\Delta x}(b,x) - \frac{b\Delta x}{4b+b^3\Delta x^2}\sin_{\Delta x}(b,x) + k$<br>$x = 0, \Delta x, 2\Delta x, 3\Delta x, \ldots$ | |
| **104** | ${}_{\Delta x}\int \cos_{\Delta x}(b,x)\Delta x = \frac{1}{b}\ \sin_{\Delta x}(b,x) + k$ | |
| **105** | ${}_{2\Delta x}\int \cos_{\Delta x}(b,x)\Delta x = \frac{2}{4b+b^3\Delta x^2}\ \sin_{\Delta x}(b,x) - \frac{b\Delta x}{4b+b^3\Delta x^2}\ \cos_{\Delta x}(b,x) + k$<br>$x = 0, \Delta x, 2\Delta x, 3\Delta x, \ldots$ | |

| | Operation | Description |
|---|---|---|
| **106** | ${}_{\Delta x}\int e_{\Delta x}(a,x)\sin_{\Delta x}(b,x)\Delta x = \frac{e_{\Delta x}(a,x)}{a^2+b^2(1+a\Delta x)^2}[a\sin_{\Delta x}(b,x) - b(1+a\Delta x)\cos_{\Delta x}(b,x)] + k$ | |
| **107** | ${}_{\Delta x}\int e_{\Delta x}(a,x)\cos_{\Delta x}(b,x)\Delta x = \frac{e_{\Delta x}(a,x)}{a^2+b^2(1+a\Delta x)^2}[a\cos_{\Delta x}(b,x) + b(1+a\Delta x)\sin_{\Delta x}(b,x)] + k$ | |
| **108** | ${}_{\Delta x}\int e_{\Delta x}(a,x)\sinh_{\Delta x}(b,x)\Delta x = \frac{e_{\Delta x}(a,x)}{a^2-b^2(1+a\Delta x)^2}[a\sinh_{\Delta x}(b,x) - b(1+a\Delta x)\cosh_{\Delta x}(b,x)] + k$ | |
| **109** | ${}_{\Delta x}\int e_{\Delta x}(a,x)\cosh_{\Delta x}(b,x)\Delta x = \frac{e_{\Delta x}(a,x)}{a^2-b^2(1+a\Delta x)^2}[a\cosh_{\Delta x}(b,x) - b(1+a\Delta x)\sinh_{\Delta x}(b,x)] + k$ | |
| **110** | ${}_{\Delta x}\int \sinh_{\Delta x}x\,\Delta x = \cosh_{\Delta x}x + k$ | $\sinh_{\Delta x}x = \sinh_{\Delta x}(1,x)$ |
| **111** | ${}_{\Delta x}\int \cosh_{\Delta x}x\,\Delta x = \sinh_{\Delta x}x + k$ | $\cosh_{\Delta x}x = \cosh_{\Delta x}(1,x)$ |
| **112** | ${}_{\Delta x}\int \sinh_{\Delta x}(a,x)\Delta x = \frac{1}{a}\cos_{\Delta x}(a,x) + k$ | |
| **113** | ${}_{\Delta x}\int \cosh_{\Delta x}(a,x)\Delta x = \frac{1}{a}\sinh_{\Delta x}(a,x) + k$ | |
| **114** | ${}_{\Delta x}\int (1+a\Delta x)^{-\frac{x}{\Delta x}}\Delta x = -\frac{1+a\Delta x}{a}(1+a\Delta x)^{-\frac{x}{\Delta x}} + k$ | |
| **115** | ${}_{\Delta x}\int z^{-\frac{x}{\Delta x}}\Delta x = -\left(\frac{z\Delta x}{z-1}\right)z^{-\frac{x}{\Delta x}} + k$<br>where<br>$z-1 \neq 0$ | |
| **116** | ${}_{\Delta x}\int c^{-\frac{x}{\Delta x}}\Delta x = \frac{\Delta x}{c^{-1}-1}c^{-\frac{x}{\Delta x}} + k$<br>where<br>c = constant<br>$c^{-1}-1 \neq 0$ | |
| **117** | ${}_{\Delta x}\int x^2\Delta x = \frac{x(x-\Delta x)(x-2\Delta x)}{3} + \frac{x(x-\Delta x)}{2}\Delta x + k$ | |
| **118** | ${}_{\Delta x}\int x^3\Delta x = \frac{x(x-\Delta x)(x-2\Delta x)(x-3\Delta x)}{4} + x(x-\Delta x)(x-2\Delta x)\Delta x + \frac{x(x-\Delta x)}{2}\Delta x^2 + k$ | |
| **119** | ${}_{\Delta x}\int x^4\Delta x = \frac{x(x-\Delta x)(x-2\Delta x)(x-3\Delta x)(x-4\Delta x)}{5}+\frac{3}{2}x(x-\Delta x)(x-2\Delta x)(x-3\Delta x)\Delta x+\frac{7}{3}x(x-\Delta x)(x-2\Delta x)\Delta x^2$<br>$+\frac{1}{2}x(x-\Delta x)\Delta x^3 + k$ | |

| | Operation | Description |
|---|---|---|
| **120** | $_{\Delta x}\int \sin\frac{\pi x}{\Delta x}\Delta x = -\frac{\Delta x}{2}\sin\frac{\pi x}{\Delta x} + k$ | |
| **121** | $_{\Delta x}\int \cos\frac{\pi x}{\Delta x}\Delta x = -\frac{\Delta x}{2}\cos\frac{\pi x}{\Delta x} + k$ | |
| **122** | $_{\Delta x}\int \sin\frac{\pi x}{2\Delta x}\Delta x = -\frac{\Delta x}{2}(\sin\frac{\pi x}{2\Delta x} + \cos\frac{\pi x}{2\Delta x}) + k$ | |
| **123** | $_{\Delta x}\int \cos\frac{\pi x}{2\Delta x}\Delta x = \frac{\Delta x}{2}(\sin\frac{\pi x}{2\Delta x} - \cos\frac{\pi x}{2\Delta x}) + k$ | |
| **124** | $_{\Delta x}\int A^x\Delta x = \begin{cases} (\frac{\Delta x}{A^{\Delta x}-1})A^x + k & A^{\Delta x}-1 \neq 0 \\ x + k & A^{\Delta x}-1 = 0 \end{cases}$ | |
| **125** | $_{\Delta x}\int b^{ax}\Delta x = \begin{cases} (\frac{\Delta x}{b^{a\Delta x}-1})b^{ax} + k & b^{a\Delta x}-1 \neq 0 \\ x + k & b^{a\Delta x}-1 = 0 \end{cases}$ | |
| **126** | $_{\Delta x}\int e^{ax}\Delta x = \begin{cases} (\frac{\Delta x}{e^{a\Delta x}-1})e^{ax} + k & e^{a\Delta x}-1 \neq 0 \\ x + k & e^{a\Delta x}-1 = 0 \end{cases}$ | |
| **127** | $_{\Delta x}\int \sin ax\Delta x = \begin{cases} \frac{\Delta x}{2}[\frac{\sin a(x-\Delta x) - \sin ax}{1-\cos a\Delta x}] + k & 1-\cos a\Delta x \neq 0 \\ k & 1-\cos a\Delta x = 0 \end{cases}$<br>or<br>$_{\Delta x}\int \sin ax\Delta x = \begin{cases} -\frac{\Delta x}{2}[\sin ax + (\frac{\sin a\Delta x}{1-\cos a\Delta x})\cos ax] + k & 1-\cos a\Delta x \neq 0 \\ k & 1-\cos a\Delta x = 0 \end{cases}$ | |
| **128** | $_{\Delta x}\int \cos ax\Delta x = \begin{cases} \frac{\Delta x}{2}[\frac{\cos a(x-\Delta x) - \cos ax}{1-\cos a\Delta x}] + k & 1-\cos a\Delta x \neq 0 \\ x + k & 1-\cos a\Delta x = 0 \end{cases}$<br>or<br>$_{\Delta x}\int \cos ax\Delta x = \begin{cases} -\frac{\Delta x}{2}[\cos ax - (\frac{\sin a\Delta x}{1-\cos a\Delta x})\sin ax] + k & 1-\cos a\Delta x \neq 0 \\ k & 1-\cos a\Delta x = 0 \end{cases}$ | |
| **129** | $_{\Delta x}\int \sinh ax\Delta x = \begin{cases} \frac{\Delta x}{2}[\frac{\sinh a(x-\Delta x) - \sinh ax}{1-\cosh a\Delta x}] + k & 1-\cosh a\Delta x \neq 0 \\ k & 1-\cosh a\Delta x = 0 \end{cases}$<br>or<br>$_{\Delta x}\int \sinh ax\Delta x = \begin{cases} -\frac{\Delta x}{2}[\sinh ax + (\frac{\sinh a\Delta x}{1-\cosh a\Delta x})\cosh ax] + k & 1-\cosh a\Delta x \neq 0 \\ k & 1-\cosh a\Delta x = 0 \end{cases}$ | |

| | Operation | Description |
|---|---|---|
| | | |
| **130** | $_{\Delta x}\int \cosh ax \Delta x = \begin{cases} \frac{\Delta x}{2}\left[\frac{\cosh a(x-\Delta x) - \cosh ax}{1-\cosh a\Delta x}\right] + k & 1-\cosh a\Delta x \neq 0 \\ x + k & 1-\cosh a\Delta x = 0 \end{cases}$<br>or<br>$_{\Delta x}\int \cosh ax \Delta x = \begin{cases} -\frac{\Delta x}{2}\left[\cosh ax + \left(\frac{\sinh a\Delta x}{1-\cosh a\Delta x}\right)\sinh ax\right] + k & 1-\cosh a\Delta x \neq 0 \\ k & 1-\cosh a\Delta x = 0 \end{cases}$ | |
| **131** | $_{\Delta x}\int \ln(1+\frac{\Delta x}{x})\Delta x = \Delta x \ln x + k$ | |
| **132** | $_{1}\int \frac{x}{(x+1)!}\Delta x = -\frac{1}{x!} + k \quad , \quad \Delta x = 1$ | |
| **133** | $_{1}\int \ln x \Delta x = \ln\Gamma(x) + k \quad , \quad \Delta x = 1$ | |
| **134** | $_{\Delta x}\int \ln_{\Delta x} x \Delta x = (x-\Delta x)\ln_{\Delta x} x - x + k$ | $\ln_{\Delta x} x \equiv \mathrm{lnd}(1,\Delta x,x)$ |
| **135** | $_{\Delta x}\int \mathrm{lnd}(n,\Delta x,x)\Delta x = (x-\Delta x)\mathrm{lnd}(n,\Delta x,x) \pm \mathrm{lnd}(n-1,\Delta x,x) + k$<br>where<br>$+$ for $n = 1, 2$<br>$-$ for $n \neq 1, 2$ | |
| **136** | $_{\Delta x}\int g(x)\Delta x = h(x) + k$<br>where<br>$g(x) = \sum_{n=0}^{N} a_n x^n$ , a polynomial of order N<br>$h(x) = \sum_{n=0}^{N+1} b_n x^n$ , a polynomial of order N+1<br>$[H_{\Delta x}-1]h(x) = g(x)\Delta x$<br>k = constant of integration<br>$\Delta x$ = x increment<br>$a_n, b_n$ = constants<br>$H_{\Delta x}h(x) = h(x+\Delta x)$<br><br>The coefficients, $b_n$ , $n = 0, 1, 2, 3, \ldots , N$ are calculated from the equation, $[H_{\Delta x}-1]h(x) = g(x)\Delta x$, by equating the coefficients of the terms of like powers of x. | Integral of a Polynomial |

| | Operation | Description |
|---|---|---|
| | | |
| **137** | $\sum_{\Delta x\;x=x_1}^{x_2} f(x) = \frac{1}{\Delta x}\;{}_{\Delta x}\!\int_{x_1}^{x_2+\Delta x} f(x)\Delta x = -\frac{1}{\Delta x}\;{}_{-\Delta x}\!\int_{x_2}^{x_1+(-\Delta x)} f(x)\Delta x$<br>$\Delta x$ = x increment<br>$\Delta x = \frac{x_2 - x_1}{n}$ , n =number of intervals between $x_1$ and $x_2$<br>f(x) = function of x<br>$x_1$ = intial x value of the summation<br>$x_2$ = final x value of the summation<br>$x = x_1, x_1+\Delta x, x_1+2\Delta x, x_1+3\Delta x, \ldots , x_2-\Delta x, x_2$ | Summation Evaluation |
| **138** | $D_{\Delta x}\;{}_{\Delta x}\!\int f(x)\Delta x = f(x)$ | Discrete derivative of a discrete indefinite integral |
| **139** | $\lim_{\Delta x \to 0}\;{}_{\Delta x}\!\int_{x_1}^{x_2} f(x)\Delta x = \int_{x_1}^{x_2} f(x)dx$ | |
| | <u>**Integral Equations**</u><br>**Definite integrals** | |
| **140** | ${}_{\Delta x}\!\int_{x_1}^{x_2} b^x \Delta x = \begin{cases} \left(\frac{\Delta x}{b^{\Delta x}-1}\right)(b^{x_2} - b^{x_1}) & b^{\Delta x} - 1 \neq 0 \\ x_2 - x_1 & b^{\Delta x} - 1 = 0 \end{cases}$ | |
| **141** | ${}_{\Delta t}\!\int_{0}^{T} \sin\left(\frac{2\pi n t}{T}\right)\Delta t = 0$<br>where<br>n = integers | |
| **142** | ${}_{\Delta t}\!\int_{0}^{T} \cos\left(\frac{2\pi n t}{T}\right)\Delta t = \begin{cases} 0 & 1-\cos\left(\frac{2\pi n \Delta t}{T}\right) \neq 0 \\ T & 1-\cos\left(\frac{2\pi n \Delta t}{T}\right) = 0 \end{cases}$<br>where<br>n = integers | |

| | Operation | Description |
|---|---|---|
| | | |
| **143** | $${}_{\Delta x}\int_{x_1}^{x_2} f(x)\Delta x = -{}_{-\Delta x}\int_{x_2}^{x_1} H_{-\Delta x}f(x)\Delta x = -{}_{-\Delta x}\int_{x_2}^{x_1} f(x-\Delta x)\Delta x = -{}_{-\Delta x}\int_{x_2+(-\Delta x)}^{x_1+(-\Delta x)} f(x)\Delta x$$ <br> where <br> $\Delta x = \frac{x_2-x_1}{n}$ , n = the number of x increments between a and b <br> $\Delta x$ = x increment <br> f(x) = function of x <br> $H_{-\Delta x}f(x) = f(x-\Delta x)$ <br> $x = x_1, x_1+\Delta x, x_1+2\Delta x, \ldots, x_2-\Delta x, x_2$ <br> $${}_{\Delta x}\int_{x_1}^{x_2} f(x)\Delta x = \Delta x \sum_{\Delta x\; x=x_1}^{x_2-\Delta x} f(x)$$ | Changing the direction of Integration |
| **144** | $$D_{\Delta x}\, {}_{\Delta x}\int_{x}^{x+c} f(x)\Delta x = {}_{\Delta x}\int_{0}^{c} \vartheta_{\Delta x} f(x+x)\Delta x$$ <br> where <br> $f(x)$ = a function of $x$ <br> $\Delta x$ = $x$ increment <br> $\Delta x$ = x increment <br> $\Delta x = \Delta x$ <br> $c = m\Delta x$ , m=0,1,2,3,… <br> $\vartheta_{\Delta x}$ = discrete partial derivative operator. The partial derivative is with respect to x with an x increment of Δx. | Discrete derivative of a discrete definite integral |
| **145** | $${}_{\Delta x}\int_{0}^{\infty} z^{-\frac{x}{\Delta x}}\Delta x = \frac{z\Delta x}{z-1}$$ | |
| **146** | $\zeta(n,\Delta x,x) = \frac{1}{\Delta x}\,\mathrm{lnd}(n,\Delta x,x)$ , for all n <br><br> $$\zeta(n,\Delta x,x_i) = \frac{1}{\Delta x}\,\mathrm{lnd}(n,\Delta x,x_i) = \sum_{\Delta x\; x=x_i}^{\pm\infty} \frac{1}{x^n} , \quad \mathrm{Re}(n)>1$$ <br> $+\infty$ for $\mathrm{Re}(\Delta x)>0$ or $\{\mathrm{Re}(\Delta x)=0$ and $\mathrm{Im}(\Delta x)>0\}$ <br> $-\infty$ for $\mathrm{Re}(\Delta x)<0$ or $\{\mathrm{Re}(\Delta x)=0$ and $\mathrm{Im}(\Delta x)<0\}$ <br><br> and | General Zeta Function |

| | Operation | Description |
|---|---|---|
| | $\zeta(n,\Delta x,x)\Big\vert_{x_1}^{x_2} = \frac{1}{\Delta x}\,\mathrm{lnd}(n,\Delta x,x)\Big\vert_{x_1}^{x_2} = \pm \sum_{\Delta x\; x=x_1}^{x_2-\Delta x} \frac{1}{x^n}$ ,  + for n = 1,  − for $n \neq 1$<br><br>where<br>$\zeta(n,\Delta x,x)$ = General Zeta Function<br>$\Delta x$ = x interval<br>x = real or complex variable<br>$x_i,x_1,x_2,\Delta x,n$ = real or complex constants<br>Any summation term where x = 0 is excluded | |
| **147** | $\int_{\Delta x\; x_i}^{\infty} \frac{1}{x^n}\,\Delta x = \Delta x \sum_{\Delta x\; x=x_i}^{\infty} \frac{1}{x^n} = \mathrm{lnd}(n,\Delta x,x_i) = \Delta x\,\zeta(n,\Delta x,x_i)$ ,  $\mathrm{Re}(n) > 1$<br>$\zeta(n,\Delta x,x_i)$ = General Zeta Function<br>Any summation term where x = 0 is excluded | |
| **148** | $\int_{1\; x_i}^{\infty} \frac{1}{x^n}\,\Delta x = \sum_{1\; x=x_i}^{\infty} \frac{1}{x^n} = \mathrm{lnd}(n,1,x_i) = \zeta(n,x_i)$ ,  $\mathrm{Re}(n)>1$<br>$\zeta(n,x)$ = Hurwitz Zeta Function<br>Any summation term where x = 0 is excluded | |
| **149** | $\int_{1\; 1}^{\infty} \frac{1}{x^n}\,\Delta x = \sum_{1\; x=1}^{\infty} \frac{1}{x^n} = \mathrm{lnd}(n,1,1) = \zeta(n)$<br>$\zeta(n)$ = Riemann Zeta Function<br>$\mathrm{Re}(n)>1$ | |
| **150** | n = 1<br>$\int_{\Delta x\; x_1}^{x_2} \frac{1}{x}\,\Delta x = \Delta x \sum_{\Delta x\; x=x_1}^{x_2-\Delta x} \frac{1}{x} = \mathrm{lnd}(1,\Delta x,x)\Big\vert_{x_1}^{x_2} \equiv \ln_{\Delta x} x\,\Big\vert_{x_1}^{x_2}$<br><br>Any summation term where x = 0 is excluded | |
| **151** | $n \neq 1$<br>$\int_{\Delta x\; x_1}^{x_2} \frac{1}{x^n}\,\Delta x = \Delta x \sum_{\Delta x\; x=x_1}^{x_2-\Delta x} \frac{1}{x^n} = -\,\mathrm{lnd}(n,\Delta x,x)\Big\vert_{x_1}^{x_2}$ | |

| | Operation | Description |
|---|---|---|
| | Any summation term where x = 0 is excluded | |
| **152** | n = 1<br>$_{\Delta x}\int_{x_1}^{x_2} \frac{1}{x-a}\,\Delta x = \Delta x \; _{\Delta x}\sum_{x=x_1}^{x_2-\Delta x} \frac{1}{x-a} = \text{lnd}(1,\Delta x,x-a)\Big\|_{x_1}^{x_2} \equiv \ln_{\Delta x}(x-a)\Big\|_{x_1}^{x_2}$<br>where<br>a = constant<br>Any summation term where x = 0 is excluded | |
| **153** | n ≠ 1<br>$_{\Delta x}\int_{x_1}^{x_2} \frac{1}{(x-a)^n}\,\Delta x = \Delta x \; _{\Delta x}\sum_{x=x_1}^{x_2-\Delta x} \frac{1}{(x-a)^n} = -\,\text{lnd}(n,\Delta x,x-a)\Big\|_{x_1}^{x_2}$<br>a = constant<br>Any summation term where x = 0 is excluded | |
| **154** | n = any value<br>$_{\Delta x}\int_{x_i}^{\infty} (-1)^{\frac{x-x_i}{\Delta x}} \frac{1}{x^n}\,\Delta x = \Delta x \; _{\Delta x}\sum_{x=x_i}^{\infty} (-1)^{\frac{x-x_i}{\Delta x}} \frac{1}{x^n} = \frac{\alpha(n)}{2^n}\left[\text{lnd}(n,\Delta x,\frac{x_i}{2}) - \text{lnd}(n,\Delta x,\frac{x_i+\Delta x}{2})\right]$<br>where $\alpha(n) = \begin{cases} 1 & n \neq 1 \\ -1 & n = 1 \end{cases}$<br>$x = x_i, x_i+\Delta x, x_i+2\Delta x, x_i+3\Delta x, \ldots$<br>$\Delta x$ = x increment<br>$x \neq 0$ where<br>n = integers | |
| **155** | n = any value<br>$_{\Delta x}\int_{x_1}^{x_2+\Delta x} (-1)^{\frac{x-x_1}{\Delta x}} \frac{1}{x^n}\,\Delta x = \Delta x \; _{\Delta x}\sum_{x=x_1}^{x_2} (-1)^{\frac{x-x_1}{\Delta x}} \frac{1}{x^n} = \frac{\alpha(n)}{2}\left[-\text{lnd}(n,2\Delta x,x)\Big\|_{x_1}^{x_2+\Delta x} + \text{lnd}(n,2\Delta x,x)\Big\|_{x_1+\Delta x}^{x_2+2\Delta x}\right]$<br>where $\alpha(n) = \begin{cases} 1 & n \neq 1 \\ -1 & n = 1 \end{cases}$<br>$x = x_1, x_1+\Delta x, x_1+2\Delta x, x_1+3\Delta x, \ldots, x_2-\Delta x, x_2$<br>$x_2 = x_1 + (2m+1)\Delta x$ ,  m=0,1,2,3,…<br>$\Delta x$ = x increment<br>$x \neq 0$ | |

| | Operation | Description |
|---|---|---|
| **156** | $\ln_1 x = \int_0^\infty \frac{e^{-t} - e^{-xt}}{1 - e^{-t}}\,dt$ <br> and <br> ${}_1\int_{x_1}^{x_2} \frac{1}{x}\Delta x = \Delta x \; {}_1\sum_{x=x_1}^{x_2-1} \frac{1}{x} = \ln_1 x \,\Big\|_{x_1}^{x_2} = \int_0^\infty \frac{e^{-x_1 t} - e^{-x_2 t}}{1 - e^{-t}}\,dt$ <br> where <br> $\Delta x, x_1, x_2$ = real or complex values <br> $x = x_1, x_1+1, x_1+2, x_1+3, \ldots, x_2-1, x_2$ <br> $\Delta x = 1$ <br> <u>Note</u> – $\ln_{\Delta x} x = \ln_1 \frac{x}{\Delta x} = \ln_1 \mathrm{x}$ where $\mathrm{x} = \frac{x}{\Delta x}$ | $\ln_1 x \equiv \mathrm{lnd}(1,1,x)$ |
| **157** | ${}_{\Delta\theta}\int_0^{2\pi} e^{jn\theta}\,\Delta\theta = \begin{cases} 0 & e^{jn\Delta\theta} - 1 \neq 0 \\ 2\pi & e^{jn\Delta\theta} - 1 = 0 \end{cases}$ <br> where <br> $\Delta\theta = \frac{2\pi}{m}$ <br> $\theta = 0, \Delta\theta, 2\Delta\theta, 3\Delta\theta, \ldots, 2\pi - \Delta\theta, 2\pi$ <br> n = integers <br> m = 1,2,3,… | |
| **158** | ${}_{\Delta x}\int_{x_1}^{x_2} f(x)\Delta x = \Delta x \; {}_{\Delta x}\sum_{x=x_1}^{x_2-\Delta x} f(x) = -\sum_{n=0}^{\infty} D^n f(x_1) \frac{\mathrm{lnd}(-n,\Delta x,x-x_1)}{n!} \Big\|_{\Delta x\; x_1}^{x_2}$ <br> where <br> f(x) = continuous function of x <br> $D^n f(x_1)$ = nth derivative of f(x) evaluated at $x = x_1$ <br> $x = x_1, x_1+\Delta x, x_1+2\Delta x, \ldots, x_2-\Delta x, x_2$ <br> The Taylor Series of f(x) must be convergent. | |

| | Operation | Description |
|---|---|---|
| **159** | $_{\Delta\lambda}\int_{0}^{x} g(x-\lambda-\Delta\lambda)\Delta\lambda = {}_{\Delta\lambda}\int_{0}^{x} g(\lambda)\Delta\lambda$ | |
| **160** | $_{T}\int_{0}^{x+T} g(x-\lambda)\Delta\lambda = {}_{T}\int_{0}^{x+T} g(\lambda)\Delta\lambda$ | |
| **161** | $D_{\Delta x}\,{}_{\Delta t}\int_{a}^{x} f(t)\Delta t = f(x)$ , $\Delta t = \Delta x$ , a = constant | |
| | **<u>Summation Equations</u>** | |
| | See the summation equations in Table 17 | |
| | **<u>Area Calculation Equations</u>** | |
| **162** | $A = \Delta x_{\Delta x}\sum_{x=x_1}^{x_2-\Delta x} f(x) = {}_{\Delta x}\int_{x_1}^{x_2} f(x)\Delta x - \Delta x[\sum MV[f(x)]$ , $x = x_1, x_2, x_3, \ldots, x_m$ ]<br>any $\frac{1}{0}$ term excluded<br>where<br>$MV[f(x)] = \lim_{\varepsilon\to 0}\left[\frac{f(x+\varepsilon) + f(x-\varepsilon)}{2}\right]$ , Mean Value of f(x)<br>$x = x_1, x_1+\Delta x, x_1+2\Delta x, x_1+3\Delta x, \ldots, x_2-\Delta x, x_2$<br>$x_r$ = those values of x, $x_1<x< x_2$, for which f(x) poles occur<br>r = 1,2,3,…,m<br>m = the number of poles between the limits $x_1$,$x_2$<br>$f(x_r)$ = summation terms divided by 0 which are excluded, r = 1,2,3,…,m<br>For all x = real values<br>A = area under the f(x) curve between the limits $x_1$ and $x_2$ as defined by the summation, $\Delta x_{\Delta x}\sum_{x=x_1}^{x_2-\Delta x} f(x)$, where if a summation term has a pole at $x_r$, it is excluded from the summation<br>$MV[f(x_r)]$ = Mean Value of the function, f(x), calculated at $x=x_r$<br>$\Delta x$,x = real or complex values<br>There can be no pole at the limits, $x_1$ and $x_2$<br>If there are no poles between $x_1$ and $x_2$, the area equation term, $\Delta x\sum DV(x)$, is left out.<br>The values of $MV[{}_{\Delta x}\int_{x_1}^{x_2} f(x)\Delta x]$, $x_1$, $x_2$, and $\sum MV[f(x)]$ , $x = x_1, x_2, x_3, \ldots, x_m$, must be finite. For the above conditions, A will be a finite value. | General Area Calculation Formla |

| | Operation | Description |
|---|---|---|
| **163** | $A = \Delta x \sum_{\Delta x}\limits_{x=x_1}^{x_2-\Delta x} \frac{1}{(x-a)^n} = \pm \,\mathrm{lnd}(n,\Delta x,x-a) \Big\vert_{x_1}^{x_2}$ , + for n = 1 , – for $n \neq 1$<br>any $\frac{1}{0}$ term excluded<br>where<br>$f(x) = \frac{1}{(x-a)^n}$<br>$x = x_1, x_1+\Delta x, x_1+2\Delta x, x_1+3\Delta x, \ldots , x_2-\Delta x, x_2$<br>There can be no pole at the limits, $x_1$ and $x_2$<br>$n,\Delta x,x,a$ = real or complex values<br>For all x = real values<br>A = area under the f(x) curve between the limits $x_1$ and $x_2$ as defined by the summation, $\Delta x \sum_{\Delta x}\limits_{x=x_1}^{x_2-\Delta x} f(x)$, where if a summation term has a pole at x = a, it is excluded from the summation | Area Calculation Formula for the function,<br>$f(x) = \frac{1}{(x-a)^n}$ |
| **164** | **Calculating the area under a continuous function using discrete integration**<br>$A = \lim_{\Delta x\to 0} \int_{\Delta x}\limits_{x_1}^{x_2} f(x)\,\Delta x = \int_{x_1}^{x_2} f(x)\,\Delta x$<br>where<br>$x_1 \leq x \leq x_2$<br>$x = x_1, x_1+\Delta x, x_1+2\Delta x, x_1+3\Delta x, \ldots , x_2-\Delta x, x_2$ , discrete values of x<br>f(x) = continuous function of x<br>A = area under the function, f(x), within the limits $x_1,x_2$<br>$\Delta x$ = x increment , $\Delta x \to 0$<br>There are no poles at any x | |
| **165** | **Estimating the area under a continuous function using discrete integration**<br>$A_e = \int_{\Delta x}\limits_{x_1}^{x_2} f(x)\,\Delta x$<br>where<br>$x_1 \leq x \leq x_2$<br>$x = x_1, x_1+\Delta x, x_1+2\Delta x, x_1+3\Delta x, \ldots , x_2-\Delta x, x_2$ , discrete values of x<br>f(x) = continuous function of x<br>$A_e$ = estimated area under the function, f(x), within the limits $x_1,x_2$<br>$\Delta x$ = x increment<br>Accuracy increases as $\Delta x \to 0$<br>There are no poles at any x<br>$A_e$ is exact for f(x) being a constant | |

| | Operation | Description |
|---|---|---|
| **166** | **Estimating the area under a continuous function using the discrete integration implementation of the Trapazoidal Rule**<br>$A_e = \underset{\Delta x}{\int_{x_1}^{x_2}} f(x)\,\Delta x + f(x)\frac{\Delta x}{2}\Big\vert_{\Delta x\; x_1}^{x_2}$<br>where<br>$x_1 \leq x \leq x_2$<br>$x = x_1, x_1+\Delta x, x_1+2\Delta x, x_1+3\Delta x, \ldots, x_2-\Delta x, x_2$ , discrete values of x<br>$f(x)$ = continuous function of x<br>$A_e$ = estimated area under the function, f(x), within the limits $x_1,x_2$<br>$\Delta x$ = x increment<br>Accuracy increases as $\Delta x \rightarrow 0$<br>There are no poles at any x<br>$A_e$ is exact for f(x) being a polynomial up to the first order | |
| **167** | **Estimating the area under a continuous function using the discrete integration implementation of the Modified Trapazoidal Rule**<br>$A_e = \underset{\Delta x}{\int_{x_1}^{x_2}} f(x)\,\Delta x + \left[\frac{f(x-\Delta x)+12f(x)-f(x+\Delta x)}{24}\right]\Delta x\,\Big\vert_{\Delta x\; x_1}^{x_2}$<br>where<br>$x = x_1, x_1+\Delta x, x_1+2\Delta x, x_1+3\Delta x, \ldots, x_2-\Delta x, x_2$<br>$f(x)$ = continuous function of x<br>$A_e$ = estimated area under the function, f(x), within the limits $x_1,x_2$<br>$\Delta x$ = x increment<br>Accuracy increases as $\Delta x \rightarrow 0$<br>There are no poles at any x<br>$A_e$ is exact for f(x) being a polynomial up to the third order<br>This rule is slightly less accurate than Simpson's Rule | |
| **168** | **Estimating the area under a continuous function using the discrete integration implementation of Simpson's Rule**<br>$A_e = \frac{4}{3}\underset{\frac{\Delta x}{2}}{\int_{x_1}^{x_2}} f(x)\,\Delta x - \frac{1}{3}\underset{\Delta x}{\int_{x_1}^{x_2}} f(x)\,\Delta x + f(x)\frac{\Delta x}{6}\Big\vert_{\Delta x\; x_1}^{x_2}$<br>where<br>$x_1 \leq x \leq x_2$<br>$x = x_1, x_1+\frac{\Delta x}{2}, x_1+\Delta x, x_1+\frac{3\Delta x}{2}, x_1+2\Delta x, \ldots, x_2-\frac{\Delta x}{2}, x_2$ , discrete values of x<br>$x_2 = x_1+n\Delta x$ , $n = 0,1,2,3, \ldots$<br>$f(x)$ = continuous function of x<br>$A_e$ = estimated area under the function, f(x), within the limits $x_1,x_2$<br>$\Delta x$ = x increment<br>Accuracy increases as $\Delta x \rightarrow 0$<br>There are no poles at any x<br>$A_e$ is exact for f(x) being a polynomial up to the third order | |

| | Operation | Description |
|---|---|---|
| **169** | **Estimating the area under a continuous function using the discrete integration implementation of Bode's Rule**<br>$A_e = \frac{64}{45}\ {}_{\Delta x}\!\int_{x_1}^{x_2} f(x)\,\Delta x - \frac{20}{45}\ {}_{2\Delta x}\!\int_{x_1}^{x_2} f(x)\,\Delta x + \frac{1}{45}\ {}_{4\Delta x}\!\int_{x_1}^{x_2} f(x)\,\Delta x + \frac{14}{45}\, f(x)\,\Delta x \,\Big\vert_{\Delta x}{}_{x_1}^{x_2}$<br>where<br>$x_1 \leq x \leq x_2$<br>$x = x_1, x_1+\Delta x, x_1+2\Delta x, x_1+3\Delta x, \ldots , x_2-\Delta x, x_2$ , discrete values of x<br>$x_2 = x_1+4n\Delta x$ , $n = 0,1,2,3, \ldots$<br>f(x) = continuous function of x<br>$A_e$ = estimated area under the function, f(x), within the limits $x_1,x_2$<br>$\Delta x$ = x increment<br>Accuracy increases as $\Delta x \rightarrow 0$<br><br>There are no poles at any x<br>$A_e$ is exact for f(x) being a polynomial up to the fifth order | |
| **170** | **Estimating the area under a continuous function using the discrete integration implementation of a non-uniform spacing rule**<br>$A_e = \frac{5}{120}\ {}_{\Delta x}\!\int_{x_1}^{x_2} f(x)\,\Delta x - \frac{128}{120}\ {}_{\frac{\Delta x}{2}}\!\int_{x_1}^{x_2} f(x)\,\Delta x + \frac{243}{120}\ {}_{\frac{\Delta x}{3}}\!\int_{x_1}^{x_2} f(x)\,\Delta x + \frac{11}{120}\, f(x)\,\Delta x \,\Big\vert_{\Delta x}{}_{x_1}^{x_2}$<br>where<br>$x_1 \leq x \leq x_2$<br>$x = x_1, x_1+\frac{\Delta x}{3}, x_1+\frac{\Delta x}{2}, x_1+\frac{2\Delta x}{3}, x_1+\Delta x, \ldots , x_2-\frac{\Delta x}{3}, x_2$ , discrete values of x<br>$x_2 = x_1+n\Delta x$ , $n = 0,1,2,3, \ldots$<br>f(x) = continuous function of x<br>$A_e$ = estimated area under the function, f(x), within the limits $x_1,x_2$<br>$\Delta x$ = x increment<br>Accuracy increases as $\Delta x \rightarrow 0$<br>There are no poles at any x<br>$A_e$ is exact for f(x) being a polynomial up to the fifth order<br>This rule is slightly less accurate than Bode's Rule | |
| | **<u>Complex Plane Discrete Contour Integration/Summation Equations</u>** | |
| | <u>Explanation of Discrete Contour Integration/Summation in the Complex Plane</u><br>1. The complex plane is divided into complex plane grid squares of equal area.<br>2. Allowable values occur only at the vertices of the complex plane grid squares.<br>3. A discrete closed contour in the complex plane is composed of only horizontal and vertical vectors positioned head to tail. These equal in length vectors are the sides of complex plane grid squares.<br>4. A closed discrete contour is comprised of vectors positioned head to tail | |

| | Operation | Description |
|---|---|---|
| | enclosing an integer number of grid squares.<br>5. Closed contour vectors may be integrated/summed in either a clockwise or counterclockwise direction. | |
| **171** | Summation Complex plane discrete closed contour area calculation equations<br><br>$$A_c = = -\frac{2j}{A_g}\,{}_1\!\sum_{k=0}^{v-1}\left(\frac{z_{k+1}+z_k}{2}\right)^3(z_{k+1}-z_k) = \frac{2j}{A_g}\,{}_1\!\sum_{k=0}^{v-1}\left(\frac{z_{k+1}+z_k}{2}\right)^3(z_{k+1}-z_k)$$<br>$\frac{z_{k+1}+z_k}{2} = z_k + \frac{z_{k+1}-z_k}{2}$ , an alternate form which may be used for calculation<br>where<br>$A_c$ = the area enclosed in a closed contour in the complex plane<br>$A_g$ = $\delta x\delta$ , the area of a complex plane grid square<br>$\delta$ = the length of the side of a complex plane grid square<br>v = the number of counterclockwise or clockwise pointing vectors comprising the complex plane closed contour being summed<br>$z_k$ = the coordinate value of the kth point on the closed contour<br>$z_{k+1} - z_k$ = the kth contour vector value, $z_{k+1}$ at the vector head, $z_k$ at the vector tail | |
| **172** | Integration Complex plane discrete closed contour area calculation equations (functions of k)<br><br>$$A_c = -\frac{2j}{A_g}\,{}_1\!\oint_0^v \left[\left(z_k + \frac{\Delta z}{2}\right)^3 \Delta z\right] \Delta k = \frac{2j}{A_g}\,{}_1\!\oint_0^v \left[\left(z_k + \frac{\Delta z}{2}\right)^3 \Delta z\right] \Delta k$$<br>where<br>$A_c$ = the area enclosed by the vectors forming a closed discrete contour in the complex plane<br>v = the number of horizontal and vertical vectors forming the closed contour<br>$z_k$ = f(k) , function(s) of k<br>k = 0, 1, 2, 3,…, v<br>$z_k$ = the coordinate value of the contour kth vector tail point<br>$A_g$ = $\delta x\delta$ , the area of a complex plane grid square<br>$\delta$ = the length of the side of a complex plane grid square<br>$\Delta z$ = real or imaginary vector values<br>$\Delta z = \begin{cases} +\delta \text{ for a horizontal left to right pointing vector} \\ -\delta \text{ for a horizontal right to left pointing vector} \\ +j\delta \text{ for a vertical down to up pointing vector} \\ -j\delta \text{ for a vertical up to down pointing vector} \end{cases}$<br>$z_k + \Delta z$ = the coordinates of the contour kth vector head point<br>$\Delta k = 1$ | |

| | Operation | Description |
|---|---|---|
| **173** | <u>Integration</u> Complex plane discrete closed contour area calculation equations (functions of z)<br><br>$A_c = \frac{2j}{A_g} \underset{\Delta z}{\oint_{z_i}^{z_f}} ( z + \frac{\Delta z}{2} )^3 \Delta z = - \frac{2j}{A_g} \underset{\Delta z}{\oint_{z_i}^{z_f}} ( z + \frac{\Delta z}{2} )^3 \Delta z$<br>where<br>$A_c$ = the area enclosed by the vectors forming a closed discrete contour in the complex plane<br>$z$ = the complex plane coordinates of the tail points of the contour vectors<br>$z_i$ = the initial value of z<br>$z_f$ = the final value of z<br>$\Delta z$ = real or imaginary vector values<br>$\Delta z = \begin{cases} +\delta \text{ for a horizontal left to right pointing vector} \\ -\delta \text{ for a horizontal right to left pointing vector} \\ +j\delta \text{ for a vertical down to up pointing vector} \\ -j\delta \text{ for a vertical up to down to pointing vector} \end{cases}$<br>$z + \Delta z$ = the coordinates of the head points of the contour vectors<br>$A_g = \delta x \delta$, the area of a complex plane grid square<br>$\delta$ = the length of the side of a complex plane grid square<br><br>Should a numerical value be substituted for the integral symbol $\Delta z$, it would be $\delta = \lvert \Delta z \rvert$. | |
| **174** | $N = \frac{A_c}{A_g}$<br>where<br>N = the number of grid squares enclosed in a closed contour in the complex plane<br>$A_c$ = the area enclosed in a discrete closed contour in the complex plane<br><br>$A_g$ = the area of a complex plane grid square | |
| **175** | Complex plane integratation along a contour straight line segment composed of head to tail vectors<br><br>Complex plane contour straight line segment<br><br>$z_i \xrightarrow{\Delta z} \xrightarrow{\Delta z} \xrightarrow{\Delta z} \; - \; - \; \xrightarrow{\Delta z} \xrightarrow{\Delta z} z_f$<br><br>$\underset{\Delta z}{\int_{z_i}^{z_f}} ( z + \frac{\Delta z}{2} )^3 \Delta z = \frac{1}{4} z^4 - \frac{1}{8} z^2 \Delta z^2 \Big\rvert_{z_i}^{z_f} = \frac{z^2}{8} ( 2z^2 - \Delta z^2 ) \Big\rvert_{z_i}^{z_f}$ | |

| | Operation | Description |
|---|---|---|
| | where<br>$z$ = the complex plane coordinates of the tail points of the contour vectors<br>$\Delta z$ = real or imaginary vector values<br>$\Delta z = \begin{cases} +\delta \text{ for horizontal left to right pointing vectors} \\ -\delta \text{ for horizontal right to left pointing vectors} \\ +j\delta \text{ for vertical down to up pointing vectors} \\ -j\delta \text{ for vertical up to down pointing vectors} \end{cases}$<br>$\delta$ = the length of the side of a complex plane grid square<br>$z + \Delta z$ = the coordinates of the head points of the contour vectors<br>$z_i$ = the initial value of z<br>$z_f$ = the final value of z | |
| **176** | Discrete complex plane closed contour area equation<br>$A_c = +\frac{2j}{A_g} \oint_{\Delta z} ( z + \frac{\Delta z}{2})^3 \Delta z = -\frac{2j}{A_g} \oint_{\Delta z} ( z + \frac{\Delta z}{2})^3 \Delta z = \left\| \frac{j}{2} \sum_{p=0}^{N} (-1)^p c_p^2 \right\|$<br>where<br>$A_c$ = The area enclosed within a discrete complex plane closed contour<br>N = 3,5,7,9,11,…<br>N = the number of the discrete complex plane contour corner points minus one<br>$c_p$ = the coordinates of the corner points of the discrete complex plane contour<br>p = 0,1, 2, 3, 4,…, N , corner point designations<br>\|r\| = absolute value of r<br>The contour integration can begin at any of the contour corner points. The integration initial point is designated as $c_0$. The following points progressing along the contour in the direction of integration are successively designated $c_1, c_2, c_3, \ldots$<br><br><u>Comments</u> - The value of the summation, $\frac{j}{2} \sum_{p=0}^{N} (-1)^p c_p^2$ , is real.<br>The summation is negative if the contour between $c_0$ and $c_1$ is horizontal and p increases counterclockwise around the contour or if the contour between $c_0$ and $c_1$ is vertical and p increases clockwise around the contour else the summation is positive.<br><br><u>Note</u> - The complex plane discrete closed contour area equation shows the contour enclosed area to be solely a function of the complex plane coordinates of the contour corner points. | |

| | Operation | Description |
|---|---|---|
| **177** | <u>Summation</u> Discrete closed contour summations which evaluate to zero<br><br>${}_{1}\overset{v-1}{\underset{k=0}{\sum\delta}}\left(\frac{z_{k+1}+z_k}{2}\right)^r(z_{k+1}-z_k) \;=\; {}_{1}\overset{v-1}{\underset{k=0}{\sum\delta}}\left(\frac{z_{k+1}+z_k}{2}\right)^r(z_{k+1}-z_k) = 0$<br><br>$\frac{z_{k+1}+z_k}{2} = z_k + \frac{z_{k+1}-z_k}{2}$ , an alternate form which may be used for calculation<br><br>where<br>$r = 0,1,2$<br>$\delta =$ the length of the side of a complex plane grid square<br>$v =$ the number of counterclockwise or clockwise pointing vectors comprising the complex plane closed contour being summed<br>$z_k =$ the coordinate value of the kth point on the closed contour<br>$z_{k+1} - z_k =$ the kth contour vector value, $z_{k+1}$ at the vector head, $z_k$ at the vector tail | |
| **178** | <u>Integration</u> Discrete closed contour integrals which evaluate to zero (functions of k)<br><br>${}_{1}\overset{v}{\underset{0}{\oint\delta}}\left[\left(z_k+\frac{\Delta z}{2}\right)^r\Delta z\right]\Delta k \;=\; {}_{1}\overset{v}{\underset{0}{\oint\delta}}\left[\left(z_k+\frac{\Delta z}{2}\right)^r\Delta z\right]\Delta k = 0$<br><br>where<br>$r = 0,1,2$<br>$v =$ the number of horizontal and vertical vectors forming the contour<br>$z_k = f(k)$ , function(s) of k<br>$k = 0, 1, 2, 3,\ldots, v$<br>$z_k =$ the coordinate value of the tail point of the kth contour vector<br>$\delta =$ the length of the side of a complex plane grid square<br>$\Delta z =$ real or imaginary vector values<br>$\Delta z = \begin{cases} +\delta & \text{for a horizontal left to right pointing vector} \\ -\delta & \text{for a horizontal right to left pointing vector} \\ +j\delta & \text{for a vertical down to up pointing vector} \\ -j\delta & \text{for a vertical up to down pointing vector} \end{cases}$<br>$z_k + \Delta z =$ the complex plane coordinate of the head point of the kth contour vector<br>$z_k + \frac{\Delta z}{2} =$ the complex plane coordinate of the center of the kth contour vector<br>$\Delta k = 1$ | |

| | Operation | Description |
|---|---|---|
| **179** | Integration Discrete closed contour integrals which evaluate to zero (functions of z)<br><br>$\underset{\Delta z}{\oint} ( z + \frac{\Delta z}{2})^r \, \Delta z \;=\; \underset{\Delta z}{\oint} ( z + \frac{\Delta z}{2})^r \, \Delta z \;=\; 0$<br>where<br>$r = 0,1,2$<br>$\delta =$ the length of the side of a complex plane grid square<br>$z =$ the complex plane coordinates of the tail points of the contour vectors<br>$z + \Delta z =$ the coordinates of the head points of the contour vectors<br>$z + \frac{\Delta z}{2} =$ the complex plane coordinates of the centers of the contour vectors<br>$\Delta z =$ real or imaginary value, vector values<br>$\Delta z = \begin{cases} +\delta \text{ for a horizontal left to right pointing vector} \\ -\delta \text{ for a horizontal right to left pointing vector} \\ +j\delta \text{ for a vertical down to up pointing vector} \\ -j\delta \text{ for a vertical up to down pointing vector} \end{cases}$<br><br>Should a numerical value be substituted for the integral symbol $\Delta z$, it would be $\delta = \lvert\Delta z\rvert$. | |
| | **$K_{\Delta x}$ and $J_{\Delta x}$ Transform Equations** | |
| **180** | The $K_{\Delta t}$ Transform<br><br>$F(s) = K_{\Delta t}[f(t)] = \underset{\Delta t}{\int_0^\infty} [1+s\Delta t]^{-\frac{t+\Delta t}{\Delta t}} f(t) \, \Delta t$<br>The $K_{\Delta x}$ Inverse Transform<br>$f(t) = K_{\Delta t}^{-1}[F(s)] = \frac{1}{2\pi j} \oint_c [1+s\Delta t]^{\frac{t}{\Delta t}} F(s) \, ds$<br>$s = \frac{e^{(\gamma+jw)\Delta t} - 1}{\Delta t} , \quad -\frac{\pi}{\Delta t} \le w < \frac{\pi}{\Delta t}$<br>where<br>$t = m\Delta t$ , $m = 0,1,2,3,\ldots$<br>$\Delta t =$ sampling interval, t increment<br>$f(t) =$ function of t<br>$f(t) = 0$ for $t < 0$<br>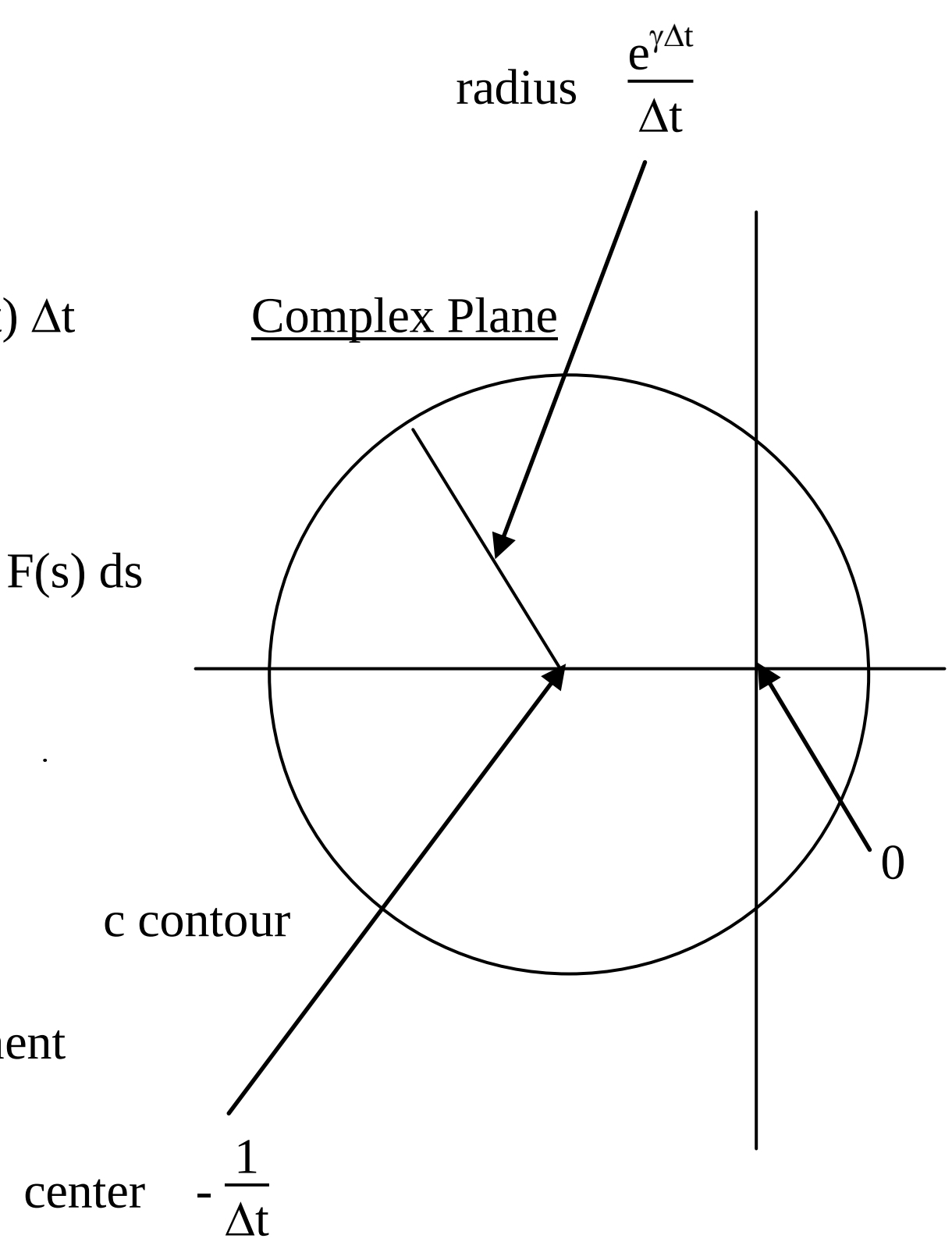 | The $K_{\Delta t}$ Transform is closely related to the $J_{\Delta x}$ Transform, $K_{\Delta x}[f(t)] = \frac{J_{\Delta x}[f(t)]}{(1+s\Delta t)}$ |

| | Operation | Description |
|---|---|---|
| | $F(s) = K_{\Delta t}[f(t)]$, a function of s<br>$\gamma$, w = real value constants<br>$\gamma > 0$<br><br>C, the complex plane contour of integration, is a circle of radius, $\frac{e^{\gamma\Delta t}}{\Delta t}$ , with center at $-\frac{1}{\Delta t}$<br><br>$\gamma$ is chosen so that the contour encloses all poles of F(s). If so, the f(t) complex plane closed contour integral can be evaluated using residue theorem methodology. The function of s, $[1+s\Delta t]^{\frac{t}{\Delta t}}$ F(s), for many F(s), can be expressed as a convergent Laurent Series for each pole as required by residue theory.<br><br><u>Notes</u> - The $K_{\Delta t}$ Transform becomes the Laplace Transform for $\Delta t \to 0$.<br>The c contour of the $K_{\Delta t}$ and $J_{\Delta t}$ Transforms is the same. | |
| **181** | <u>The $J_{\Delta t}$ Transform</u><br><br>$F(s) = J_{\Delta t}[f(t)] = \int_{0}^{\infty} [1+s\Delta t]^{-\frac{t}{\Delta t}} f(t)\, \Delta t$ <br><br><u>The $J_{\Delta t}$ Inverse Transform</u><br><br>$f(t) = J_{\Delta t}^{-1}[F(s)] = \frac{1}{2\pi j} \oint_c [1+s\Delta t]^{\frac{t}{\Delta t}-1} F(s)\, ds$<br><br>$s = \frac{e^{(\gamma+jw)\Delta t} - 1}{\Delta t}, \quad -\frac{\pi}{\Delta t} \le w < \frac{\pi}{\Delta t}$<br><br>where<br>$t = m\Delta t$ , m = 0,1,2,3,…<br>$\Delta t$ = sampling period, t increment<br>f(t) = function of t<br>$F(s) = J_{\Delta t}[f(t)]$ , a function of s<br>f(t) = 0 for t < 0<br>$\gamma$, w = real value constants<br>$\gamma > 0$<br><br>C, the complex plane contour of integration, is a circle of radius, $\frac{e^{\gamma\Delta t}}{\Delta t}$ , with center at $-\frac{1}{\Delta t}$ | The $J_{\Delta t}$ Transform is closely related to the $K_{\Delta x}$ Transform,<br>$J_{\Delta x}[f(t)] = (1+s\Delta x)K_{\Delta x}[f(t)]$ |

<table>
<tr><th></th><th>Operation</th><th>Description</th></tr>
<tr><td></td><td>

$\gamma$ is chosen so that the contour encloses all poles of F(s). If so, the f(t) complex plane closed contour integral can be evaluated using residue theorem methodology. The function of s, $[1+s\Delta t]^{\frac{t}{\Delta t}-1}$ F(s), for many F(s), can be expressed as a convergent Laurent Series for each pole as required by residue theory.

<u>Notes</u> - The $J_{\Delta t}$ Transform becomes the Laplace Transform for $\Delta t \to 0$.
The c contour of the $J_{\Delta t}$ and $K_{\Delta t}$ Transforms is the same. See the c contour diagram in row 132 above.

</td><td></td></tr>
<tr><td>182</td><td>

<u>Fourier Series for sample and hold shaped waveforms</u>

An example of a discrete Interval Calculus sample and hold shaped waveform

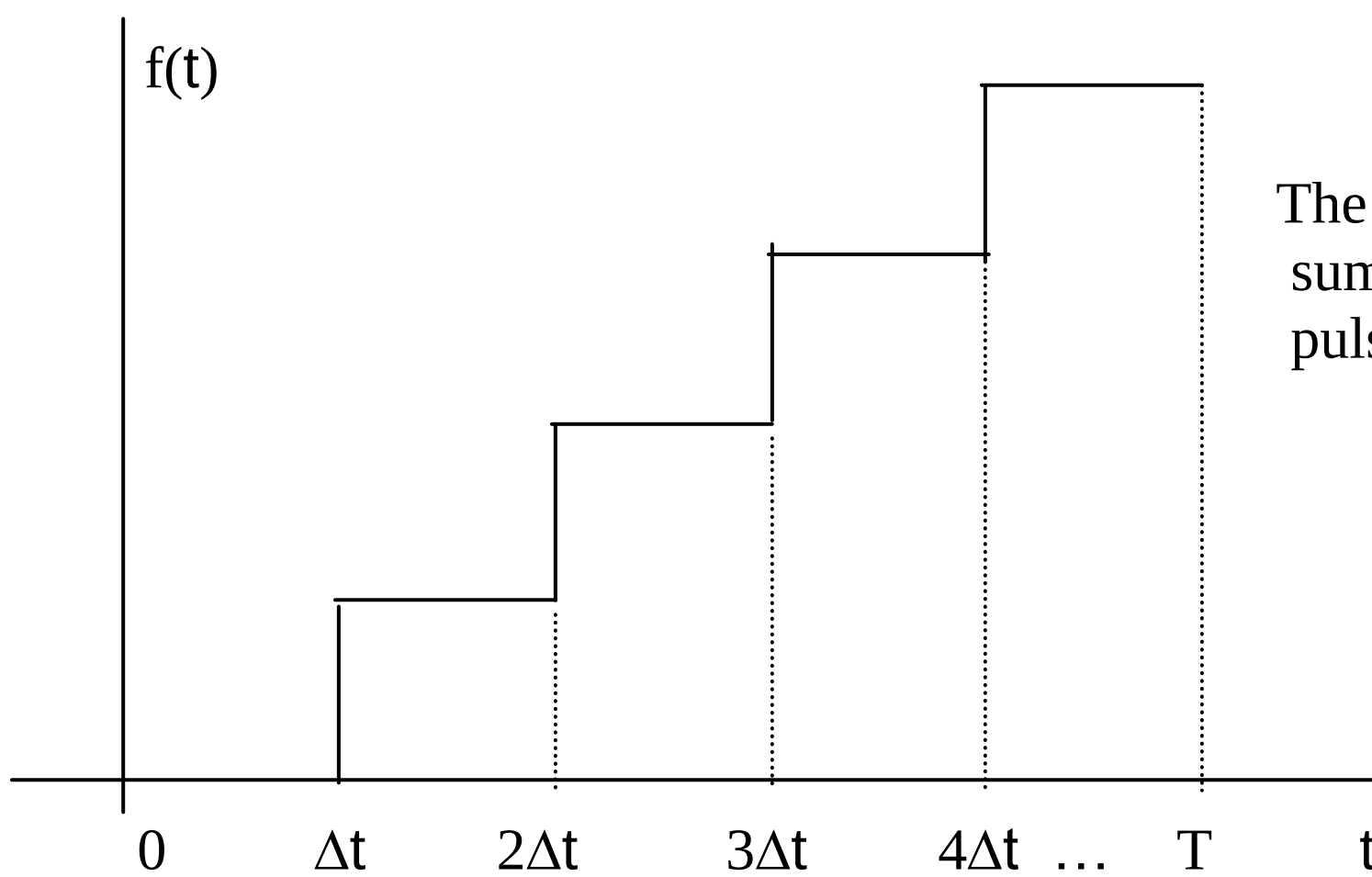

The waveform is the sum of consecutive pulses of equal width

Fourier Series of Sample and Hold Shaped Waveforms

1) $F(t) = \frac{a_0}{2} + \sum_{n=1}^{\infty} ( a_n \cos\frac{2\pi n}{T} t + b_n \sin\frac{2\pi n}{T} t )$

Fourier Series Coefficient Calculation Equations Generalized for Sample and Hold Shaped Waveforms

</td><td></td></tr>
</table>

| | Operation | Description |
|---|---|---|
| | 2) $a_0 = \frac{2}{T}\ {}_{\Delta t}\int_0^T f(t)\Delta t$ <br><br> 3) $a_n = \frac{1}{\pi n}\ {}_{\Delta t}\int_0^T f(t)\Delta \sin\frac{2\pi nt}{T} = \frac{1}{\pi n}\ {}_{\Delta t}\sum_{t=0}^{T-\Delta t} f(t)\,\Delta \sin\frac{2\pi n}{T}t$ <br><br> 4) $a_n = \frac{1}{2\pi nj}\,[\ {}_{\Delta t}\int_0^T f(t)\Delta e^{j\frac{2\pi nt}{T}} - {}_{\Delta t}\int_0^T f(t)\Delta e^{-j\frac{2\pi nt}{T}}\ ]$ <br><br> 5) $a_n = \frac{1}{\pi n}\ {}_{\Delta t}\int_0^T f(t)\Delta[e_{\Delta t}(\frac{\cos\frac{2\pi n}{T}\Delta t - 1}{\Delta t}, t)\sin_{\Delta t}(\frac{\tan\frac{2\pi n}{T}\Delta t}{\Delta t}, t)]$ <br><br> 6) $a_n = \frac{1}{2\pi nj}\,[\ {}_{\Delta t}\int_0^T f(t)\Delta e_{\Delta t}(\frac{e^{j\frac{2\pi n}{T}\Delta t} - 1}{\Delta t}, t) - {}_{\Delta t}\int_0^T f(t)\Delta e_{\Delta t}(\frac{e^{-j\frac{2\pi n}{T}\Delta t} - 1}{\Delta t}, t)\ ]$ <br><br> 7) $b_n = -\frac{1}{\pi n}\ {}_{\Delta t}\int_0^T f(t)\Delta \cos\frac{2\pi nt}{T} = -\frac{1}{\pi n}\ {}_{\Delta t}\sum_{t=0}^{T-\Delta t} f(t)\Delta \cos\frac{2\pi nt}{T}$ <br><br> 8) $b_n = -\frac{1}{2\pi n}\,[\ {}_{\Delta t}\int_0^T f(t)\Delta e^{j\frac{2\pi nt}{T}} + {}_{\Delta t}\int_0^T f(t)\Delta e^{-j\frac{2\pi nt}{T}}\ ]$ <br><br> 9) $b_n = -\frac{1}{\pi n}\ {}_{\Delta t}\int_0^T f(t)\Delta[e_{\Delta t}(\frac{\cos\frac{2\pi n}{T}\Delta t - 1}{\Delta t}, t)\cos_{\Delta t}(\frac{\tan\frac{2\pi n}{T}\Delta t}{\Delta t}, t)]$ <br><br> 10) $b_n = -\frac{1}{2\pi n}\,[\ {}_{\Delta t}\int_0^T f(t)\Delta e_{\Delta t}(\frac{e^{j\frac{2\pi n}{T}\Delta t} - 1}{\Delta t}, t) + {}_{\Delta t}\int_0^T f(t)\Delta e_{\Delta t}(\frac{e^{-j\frac{2\pi n}{T}\Delta t} - 1}{\Delta t}, t)\ ]$ <br> where <br> $t = 0, \Delta t, 2\Delta t, 3\Delta t, \ldots, T-\Delta t, T$ <br> $\Delta t = \frac{T}{m}$ = interval between successive values of t <br> $n = 1, 2, 3, \ldots$ | |

| | Operation | Description |
|---|---|---|
| | $T$ = positive constant<br>$m$ = number of intervals within the range, 0 to $T$<br>$m$ = positive integer<br>$f(t)$ = discrete sample and hold shaped waveform function of $t$<br>$\left(\frac{e^{j\frac{2\pi n}{T}\Delta t}-1}{\Delta t}\right)$ = constant<br>$\left(\frac{e^{-j\frac{2\pi n}{T}\Delta t}-1}{\Delta t}\right)$ = constant<br><br>$\Delta$ = Difference operator<br>$D_{\Delta t} = \frac{\Delta}{\Delta t}$ = discrete derivative<br>$\Delta g(t) = g(t+\Delta t) - g(t)$<br>$D_{\Delta t}\, g(t) = \frac{g(t+\Delta t) - g(t)}{\Delta t}$<br><br><u>Comment</u> – The previously specified equations containing the difference operator, $\Delta$, can be rewritten using the discrete derivative operator, $D_{\Delta t}$ .<br>For example, the equation,<br>$b_n = -\frac{1}{\pi n} \int_{\Delta t\ 0}^{T} f(t)\Delta\cos\frac{2\pi nt}{T}$, can be rewritten,<br>$b_n = -\frac{1}{\pi n} \int_{\Delta t\ 0}^{T} f(t)D_{\Delta t}\cos\frac{2\pi nt}{T}\Delta t$ . | |

# TABLE 7

## Equations for the Evaluation of Summations and Functions

| # | Definition or Relationship | Comments |
|---|---|---|
| | **<u>General Equations</u>** | |
| **1** | **The lnd(n,Δx,x) Series**<br><br><u>For all n</u><br><br>$\text{lnd}(n,\Delta x,x) \approx \left[\frac{1+\alpha(n)}{2}\right]\left[\ln\left(\frac{x}{\Delta x}-\frac{1}{2}\right)+\gamma\right]+\left[\frac{1-\alpha(n)}{2}\right]\left[\frac{1}{(n-1)(x-\frac{\Delta x}{2})^{n-1}}+k\right]$<br><br>$+\ \alpha(n)\sum_{m=1}^{\infty}\frac{\Gamma(n+2m-1)\left(\frac{\Delta x}{2}\right)^{2m}C_m}{\Gamma(n)(2m+1)!\left(x-\frac{\Delta x}{2}\right)^{n+2m-1}}$<br><br>accuracy improves rapidly for increasing $\lvert\frac{x}{\Delta x}\rvert$<br><br>where<br><br>$\alpha(n) = \begin{cases} 1 & n = 1 \\ -1 & n \neq 1 \end{cases}$<br><br>n,Δx,x = real or complex values<br>k = constant of integration for n≠1<br>γ = constant of integration for n=1, Euler's Constant .577215664…<br>Δx = x increment<br>$C_m$ = Series constants, m = 1,2,3,…<br><br>$C_1 = +1$ $\quad C_5 = +\frac{2555}{3}$<br><br>$C_2 = -\frac{7}{3}$ $\quad C_6 = -\frac{1414477}{105}$<br><br>$C_3 = +\frac{31}{3}$ $\quad C_7 = +286685$<br><br>$C_4 = -\frac{381}{5}$ $\quad$ … | $\text{lnd}(n,\Delta x,x_i) =$<br><br>$\int_{x_i}^{\infty}{}_{\Delta x}\frac{1}{x^n}\Delta x =$<br><br>$\Delta x \sum_{x=x_i}^{\infty}{}_{\Delta x}\frac{1}{x^n}$ |
| **2** | **The $\ln_{\Delta x} x \equiv \mathrm{lnd}(1,\Delta x,x)$ Series**<br><u>For n =1</u><br><br>$\ln_{\Delta x} x \approx \gamma + \ln\left(\frac{x}{\Delta x} - \frac{1}{2}\right) + \sum_{m=1}^{\infty} \frac{(2m-1)!\, C_m}{(2m+1)!\, 2^{2m} \left(\frac{x}{\Delta x} - \frac{1}{2}\right)^{2m}}$<br><br>accuracy improves rapidly for increasing $\lvert\frac{x}{\Delta x}\rvert$<br>$n,\Delta x,x,x_i$ are real or complex values<br>$\gamma$ = Euler's Constant, .577215664…<br>$x = x_i + n\Delta x$ ,    n = integers<br>$x_i$ = a value of x<br>$\Delta x$ = x increment<br>$C_m$ = series coefficients<br><br>$C_1 = +1$    $C_5 = +\frac{2555}{3}$<br>$C_2 = -\frac{7}{3}$    $C_6 = -\frac{1414477}{105}$<br>$C_3 = +\frac{31}{3}$    $C_7 = +286685$<br>$C_4 = -\frac{381}{5}$    … | This equation is used for n=1 only.<br><br>$\ln_{\Delta x} x \equiv \mathrm{lnd}(1,\Delta x,x)$<br><br>$\ln_{\Delta x} x_2 - \ln_{\Delta x} x_1 =$<br>$\int_{\Delta x\; x_1}^{x_2} \frac{1}{x}\,\Delta x =$<br>$\Delta x \sum_{\Delta x\; x=x_1}^{x_2} \frac{1}{x^n}$ |
| **2a** | $\ln_{\Delta x} x \approx \gamma + \ln\left(\frac{x}{\Delta x} - \frac{1}{2}\right) + \sum_{m=1}^{\infty} \frac{C_m}{(2m+1)(2m)\left(\frac{2x}{\Delta x} - 1\right)^{2m}}$<br><br>accuracy improves rapidly for increasing $\lvert\frac{x}{\Delta x}\rvert$<br>$n,\Delta x,x,x_i$ are real or complex values<br>$\gamma$ = Euler's Constant, .577215664…<br>$x = x_i + n\Delta x$ ,    n = integers<br>$x_i$ = a value of x<br>$\Delta x$ = x increment<br>$C_m$ = series coefficients | |
| **2b** | $$\ln_{\Delta x}x \approx \gamma - \sum_{m=0}^{\infty} \frac{C_m}{(2m+1)!} \left(\frac{\Delta x}{2}\right)^{2m} \frac{d^{2m}}{dx^{2m}} \ln\left(\frac{x}{\Delta x} - \frac{1}{2}\right)$$ <br> accuracy improves rapidly for increasing $\left\lvert\frac{x}{\Delta x}\right\rvert$ <br> $n,\Delta x,x,x_i$ are real or complex values <br> $\gamma$ = Euler's Constant, .577215664… <br> $x = x_i + n\Delta x$ ,  n = integers <br> $x_i$ = a value of x <br> $\Delta x$ = x increment <br> $C_m$ = series coefficients <br> $C_0 = -1$ | |
| **2c** | $$\ln_{\Delta x}x \approx \gamma + \ln\left(\frac{x}{\Delta x} - \frac{1}{2}\right) + \sum_{m=1}^{\infty} \frac{K_m}{\left(\frac{x}{\Delta x} - \frac{1}{2}\right)^{2m}}$$ <br> accuracy improves rapidly for increasing $\left\lvert\frac{x}{\Delta x}\right\rvert$ <br> $n,\Delta x,x,x_i$ are real or complex values <br> $\gamma$ = Euler's Constant, .577215664… <br> $x = x_i + n\Delta x$ ,  n = integers <br> $x_i$ = a value of x <br> $\Delta x$ = x increment <br> $K_m$ = series coefficients <br> $K_1 = +\frac{1}{24}$   $K_2 = -\frac{7}{960}$ <br> $K_3 = +\frac{31}{8064}$   $K_4 = -\frac{127}{30720}$ <br> $K_5 = +\frac{511}{67584}$   $K_6 = -\frac{1414477}{67092480}$ <br> $K_7 = +\frac{8191}{98304}$ <br> $K_m = \frac{(2m-1)!\,C_m}{2^{2m}(2m+1)!} = \frac{C_m}{2^{2m}(2m+1)(2m)}$ ,  m = 1,2,3,… | |
| **2d** | $\ln_{\Delta x}x \approx \gamma + \ln\left(\frac{x}{\Delta x} - \frac{1}{2}\right) + \frac{1}{24}\frac{1}{(\frac{x}{\Delta x} - \frac{1}{2})^2} - \frac{7}{960}\frac{1}{(\frac{x}{\Delta x} - \frac{1}{2})^4} + \frac{31}{8064}\frac{1}{(\frac{x}{\Delta x} - \frac{1}{2})^6}$ $- \frac{127}{30720}\frac{1}{(\frac{x}{\Delta x} - \frac{1}{2})^8} + \frac{511}{67584}\frac{1}{(\frac{x}{\Delta x} - \frac{1}{2})^{10}} - \frac{1414477}{67092480}\frac{1}{(\frac{x}{\Delta x} - \frac{1}{2})^{12}} + \frac{8191}{98304}\frac{1}{(\frac{x}{\Delta x} - \frac{1}{2})^{14}} + \ldots$<br><br>accuracy improves rapidly for increasing $\lvert\frac{x}{\Delta x}\rvert$<br><br>$n,\Delta x,x,x_i$ are real or complex values<br><br>$\gamma$ = Euler's Constant, .577215664…<br><br>$x = x_i + n\Delta x$ , n = integers<br>$x_i$ = a value of x<br>$\Delta x$ = x increment | This is 2c expanded |
| **3** | **The lnd(n,Δx,x) n≠1 Series**<br>For $n \neq 1$, $n \neq -1, -2, -3, -4, \ldots$<br><br>$\text{lnd}(n,\Delta x,x) \approx -\sum_{m=0}^{\infty}\frac{\Gamma(n+2m-1)\left(\frac{\Delta x}{2}\right)^{2m} C_m}{\Gamma(n)(2m+1)!\left(x - \frac{\Delta x}{2}\right)^{n+2m-1}} + k$<br><br>accuracy improves rapidly for increasing $\lvert\frac{x}{\Delta x}\rvert$<br>$n,\Delta x,x,x_i$ are real or complex values<br>k = the constant of integration<br>$x = x_i + r\Delta x$ , r = integers<br>$x_i$ = a value of x<br>$\Delta x$ = x increment<br><br>$C_0 = -1$     $C_4 = -\frac{381}{5}$<br>$C_1 = +1$     $C_5 = +\frac{2555}{3}$<br>$C_2 = -\frac{7}{3}$     $C_6 = -\frac{1414477}{105}$<br>$C_3 = +\frac{31}{3}$     $C_7 = +286685$ | This equation is used for $n \neq 1$ only.<br><br>$\text{lnd}(n,\Delta x,x_i) =$<br>$\int_{x_i}^{\infty}{}_{\Delta x}\frac{1}{x^n}\,\Delta x =$<br><br>$\Delta x \sum_{x=x_i}^{\infty}{}_{\Delta x}\frac{1}{x^n}$ |
| **3a** | For n = -1, -2, -3, -4, …<br>$\text{lnd}(n,\Delta x,x) = -\sum_{m=0}^{\infty} \frac{\Gamma(n+2m-1)\left(\frac{\Delta x}{2}\right)^{2m} C_m}{\Gamma(n)(2m+1)!\left(x-\frac{\Delta x}{2}\right)^{n+2m-1}}$<br>$\Delta x, x, x_i$ are real or complex values<br>$x = x_i + r\Delta x$ ,  r = integers<br>$x_i$ = a value of x<br>$\Delta x$ = x increment<br>Notes – The constant of integration is equal to 0<br>This equation is an equality | |
| **3b** | For $n \neq 1$<br>$\text{lnd}(n,\Delta x,x) \approx -\sum_{m=0}^{\infty} \frac{\left(\frac{\Delta x}{2}\right)^{2m} C_m}{(n-1)(2m+1)!} \frac{d^{2m}}{dx^{2m}}\left[\frac{1}{(x-\frac{\Delta x}{2})^{n-1}}\right] + k$<br>accuracy improves rapidly for increasing $\lvert\frac{x}{\Delta x}\rvert$<br>$n, \Delta x, x, x_i$ are real or complex values<br>k = the constant of integration<br>$x = x_i + r\Delta x$ ,  r = integers<br>$x_i$ = a value of x<br>$\Delta x$ = x increment<br>Notes – The constant of integration is equal to 0 for n = -1,-2,-3,-4,…<br>This equation is an equality for n = -1,-2,-3,-4,… | |
| **3c** | For $n \neq 1$<br>$\text{lnd}(n,\Delta x,x) \approx -\left[\frac{1}{(x-\frac{\Delta x}{2})^{n}}\right]\left[C_0\left(\frac{1}{n-1}\right)\left(x-\frac{\Delta x}{2}\right) + \sum_{m=1}^{\infty} C_m\left(\frac{(\frac{\Delta x}{2})^{2m}}{(2m+1)(2m)}\right)\prod_{p=1}^{2m-1}\left(\frac{n+p-1}{p(x-\frac{\Delta x}{2})}\right)\right] + k$<br>accuracy improves rapidly for increasing $\lvert\frac{x}{\Delta x}\rvert$<br>$n, \Delta x, x, x_i$ are real or complex values | This equation for $\text{lnd}(n,\Delta x,x)$ has the best form for computer calculations. Smaller numbers are used in the computations. |
| | k =constant of integration<br><br>$x = x_i + r\Delta x$ ,   r = integers<br><br>$x_i$ = a value of x<br><br>$\Delta x$ = x increment<br><br><u>Note</u> – The constant of integration is equal to 0 for n = -1,-2,-3,-4,… | |
| **4** | <u>The General Zeta Function</u><br><br>$\zeta(n,\Delta x,x) = \frac{1}{\Delta x}\,\text{lnd}(n,\Delta x,x)$<br><br>$\zeta(n,\Delta x,x_i) = \frac{1}{\Delta x}\,\text{lnd}(n,\Delta x,x_i) = \sum_{\Delta x\; x=x_i}^{\pm\infty} \frac{1}{x^n}$ ,   Re(n)>1<br><br>$+\infty$ for Re(Δx)>0 or {Re(Δx)=0 and Im(Δx)>0}<br>$-\infty$ for Re(Δx)<0 or {Re(Δx)=0 and Im(Δx)<0}<br><br>and<br><br>$\zeta(n,\Delta x,x)\Big\vert_{x_1}^{x_2} = \frac{1}{\Delta x}\,\text{lnd}(n,\Delta x,x)\Big\vert_{x_1}^{x_2} = \pm \sum_{\Delta x\; x=x_1}^{x_2-\Delta x} \frac{1}{x^n}$ ,<br><br>+ for n = 1,  – for $n \neq 1$,   Re(n)>1<br>where<br>$\Delta x$ = x interval<br>x = real or complex variable<br>$x_i, x_1, x_2, \Delta x, n$ = real or complex constants<br>Any summation term where x = 0 is excluded | The Zeta Function is calculated using the lnd(n,Δx,x) Function |
| **4a** | <u>The Hurwitz Zeta Function</u><br><br>$\zeta(n,x_i) = \text{lnd}(n,1,x_i)$ ,   $n \neq 1$<br><br>$\zeta(n,x_i) = \text{lnd}(n,1,x_i) = \sum_{1\; x=x_i}^{\infty} \frac{1}{x^n}$ ,   Re(n)>1<br><br>where<br>x = real or complex variable<br>$x_i$,n = real or complex constants<br>Any summation term where x = 0 is excluded<br><br>The Hurwitz Zeta Function is a special case of the General Zeta Function where $\Delta x = 1$, $x_1 = x_i$, and $x_2 \rightarrow \infty$ | A special case of the General Zeta Function |

<table>
<tr>
<td>4b</td>
<td><u>The Riemann Zeta Function</u><br><br>$\zeta(n) = \text{lnd}(n,1,1)\ ,\quad n\neq 1$<br><br>$\zeta(n) = \text{lnd}(n,1,1) = {}_{1}\sum_{x=1}^{\infty} \frac{1}{x^n}\ ,\quad \text{Re}(n)>1$<br><br>where<br>x = real or complex variable<br>n = real or complex constant<br><br>The Riemann Zeta Function is a special case of the General Zeta Function where $\Delta x = 1$, $x_1 = 1$, and $x_2 \to \infty$</td>
<td>A special case of the General Zeta Function</td>
</tr>
<tr>
<td>4c</td>
<td><u>The N=1 Zeta Function</u><br><br>$\zeta(1,\Delta x,x) = \frac{1}{\Delta x}\text{lnd}(1,\Delta x,x) \equiv \frac{1}{\Delta x}\ln_{\Delta x} x\ ,\quad n=1$<br><br>$\zeta(1,\Delta x,x_f) = \frac{1}{\Delta x}\text{lnd}(1,\Delta x,x_f) \equiv \frac{1}{\Delta x}\ln_{\Delta x} x_f = {}_{\Delta x}\sum_{x=\Delta x}^{x_f-\Delta x}\frac{1}{x}\ ,\quad n=1$<br><br>and<br><br>$\zeta(1,\Delta x,x)\Big|_{x_1}^{x_2} = \frac{1}{\Delta x}\text{lnd}(1,\Delta x,x)\Big|_{x_1}^{x_2} \equiv \frac{1}{\Delta x}\ln_{\Delta x}x\Big|_{x_1}^{x_2} = {}_{\Delta x}\sum_{x=x_1}^{x_2-\Delta x}\frac{1}{x}\ ,\quad n=1$<br><br>where<br>$\Delta x$ = x interval<br>x = real or complex variable<br>$x_1, x_2, x_f, \Delta x$ = real or complex constants<br>Any summation term where x = 0 is excluded<br><br>$\text{lnd}(1,\Delta x,\Delta x) = 0$<br>$\text{lnd}(1,\Delta x,0) = 0$<br><br>$\ln_{\Delta x}x$ is an optional form used to emphasize the similarity to the natural logarithm.<br><br>The N=1 Zeta Function is a special case of the General Zeta Function where $x_1 = \Delta x$, $x_2 = x_f$, n = 1</td>
<td>A very special case of the General Zeta Function</td>
</tr>
</table>

| | | |
|---|---|---|
| **5** | For n = 1,2,3, …<br><br>$${}_{\Delta x}\sum_{x=x_1}^{x_2} x^n = -\frac{1}{\Delta x}\sum_{m=0}^{N}\frac{\left(\frac{\Delta x}{2}\right)^{2m} C_m}{(n+1)(2m+1)!}\frac{d^{2m}}{dx^{2m}}\left[\left(x-\frac{\Delta x}{2}\right)^{n+1}\right]\Big|_{x_1}^{x_2+\Delta x}$$<br>where<br>$x = x_1, x_1+\Delta x, x_1+2\Delta x, \ldots, x_2-\Delta x, x_2$<br>$\frac{x_2 - x_1}{\Delta x}$ = integer<br>$x_1, x_2, x, \Delta x$ are real or complex values<br>$\Delta x$ = x increment<br>$N = \mathrm{Int}(\frac{n+1}{2})$ , Int(p) = integer value of p | $${}_{\Delta x}\sum_{x=x_1}^{x_2} x^n =$$ $$\frac{1}{\Delta x}\ {}_{\Delta x}\int_{x_1}^{x_2} x^n \Delta x$$ |
| **5a** | $${}_{\Delta x}\sum_{x=x_1}^{x_2} f(x) = -\sum_{m=0}^{\infty}\frac{\left(\frac{\Delta x}{2}\right)^{2m} C_m}{\Delta x(2m+1)!}\frac{d^{2m}}{dx^{2m}}\int f\left(x-\frac{\Delta x}{2}\right)dx\,\Big|_{x_1}^{x_2+\Delta x}$$<br>$x = x_1, x_1+\Delta x, x_1+2\Delta x, \ldots, x_2-\Delta x, x_2$<br>$\frac{x_2 - x_1}{\Delta x}$ = integer<br>$x_1, x_2, x, \Delta x$ = real or complex values<br>$\Delta x$ = x increment<br>f(x) = function of x<br><u>Necessary Conditions</u><br>f(x) has a Taylor/Maclaurin Series representation which is convergent for all x on the line in the complex plane delimited by $x_1$ and $x_2 + \Delta x$. Also, for a Sum Formula itself must converge. For f(x) = g(ax), a real value function of ax (a=constant), the smaller the quantity, $(a\Delta x)^2$, the more probable formula convergence will occur. (A sufficient condition for absolute convergence for g(ax) = $e^{ax}$, sinax, cosax, sinhax, and coshax is $(a\Delta x)^2 < 39$). | Sum Formula form #1<br><br>$${}_{\Delta x}\sum_{x=x_1}^{x_2} f(x) =$$ $$\frac{1}{\Delta x}\ {}_{\Delta x}\int_{x_1}^{x_2} f(x)\Delta x$$ |

| | | |
|---|---|---|
| | <u>Comment</u><br>Even if the function, f(x), satisfies the series convergence criteria, the Sum Formula may diverge or converge then diverge as more series terms are added. This may occur if the rate of change of f(x) at either x value limit increases rapidly in magnitude. f(x) = tanx is such a problematic function.<br>or<br>$$f(x) \equiv \sum_{p=1}^{P} a_p(x+c_p)^{n_p} , \quad p =1,2,3,\ldots,P$$<br>where $a_p,c_p,n_p$ = constants<br>$n_p$ = positive integers , $n_p \geq 0$<br>P = the number of terms of f(x) | |
| **5b** | $$\sum_{\Delta x\; x=x_1}^{x_2} f(x) = \frac{1}{\Delta x}\int_{x_1}^{x_2+\Delta x} f(x - \frac{\Delta x}{2})dx + \sum_{m=1}^{\infty} A_m\left(\frac{\Delta x}{2}\right)^{2m-1} \frac{d^{2m-1}}{dx^{2m-1}} f(x - \frac{\Delta x}{2}) \Big\vert_{x_1}^{x_2+\Delta x}$$<br>where $x = x_1, x_1+\Delta x, x_1+2\Delta x, \ldots, x_2-\Delta x, x_2$<br>$x_1,x_2,x,\Delta x$ = real or complex values<br>$\frac{x_2 - x_1}{\Delta x}$ = integer<br>$\Delta x$ = x increment<br>f(x) = function of x<br><u>Necessary Conditions</u><br>f(x) has a Taylor/Maclaurin Series representation which is convergent for all x on the line in the complex plane delimited by $x_1$ and $x_2 + \Delta x$. Also, the Sum Formula itself must converge. For f(x) = g(ax), a real value function of ax (a=constant), the smaller the quantity, $(a\Delta x)^2$, the more probable formula convergence will occur. (A sufficient condition for absolute convergence for $g(ax) = e^{ax}$, sinax, cosax, sinhax, and coshax is $(a\Delta x)^2 < 39$). | Sum Formula form #2<br>$\sum_{\Delta x\; x=x_1}^{x_2} f(x) =$<br>$\frac{1}{\Delta x} \int_{\Delta x\; x_1}^{x_2} f(x)\Delta x$ |

Comment

Even if the function, f(x), satisfies the series convergence criteria, the Sum Formula may diverge or converge then diverge as more series terms are added. This may occur if the rate of change of f(x) at either x value limit increases rapidly in magnitude. f(x) = tanx is such a problematic function.

or

$$f(x) \equiv \sum_{p=1}^{P} a_p(x+c_p)^{n_p}\ , \quad p = 1,2,3,\ldots,P$$

where $a_p, c_p$ = real or complex values

$n_p$ = positive integers , $n_p \geq 0$

P = the number of terms of f(x)

$A_m$ = constants

$$A_1 = -\frac{1}{12} \qquad A_2 = +\frac{7}{720}$$

$$A_3 = -\frac{31}{30240} \qquad A_4 = +\frac{127}{1209600}$$

$$A_5 = -\frac{73}{6842880} \qquad A_6 = +\frac{1414477}{1307674368000}$$

$$A_7 = -\frac{8191}{74724249600}$$

$$A_m = \frac{-C_m}{2(2m+1)!}\ , \quad m = 1,2,3,\ldots$$

$$\frac{A_{m+1}}{A_m} \approx -\frac{1}{10}$$

$$\lim_{m\to\infty}\frac{A_{m+1}}{A_m} = -.101321183642336\ldots$$

| | | |
|---|---|---|
| **5c** | $$\sum_{\Delta x\, x=x_1}^{x_2} f(x) = \frac{1}{\Delta x}\int_{x_1-\frac{\Delta x}{2}}^{x_2+\frac{\Delta x}{2}} f(x)dx + \sum_{m=1}^{\infty} B_m \Delta x^{2m-1} \frac{d^{2m-1}}{dx^{2m-1}} f(x) \Big\|_{x_1-\frac{\Delta x}{2}}^{x_2+\frac{\Delta x}{2}}$$ <br> $x = x_1, x_1+\Delta x, x_1+2\Delta x, \ldots, x_2-\Delta x, x_2$ <br> $x_1,x_2,x,\Delta x$ = real or complex values <br> $\frac{x_2 - x_1}{\Delta x}$ = integer <br> $\Delta x$ = x increment <br> $f(x)$ = function of x <br> <u>Necessary Conditions</u> <br> f(x) has a Taylor/Maclaurin Series representation which is convergent for all x on the line in the complex plane delimited by $x_1 - \frac{\Delta x}{2}$ and $x_2 + \frac{\Delta x}{2}$ . Also, the Sum Formula itself must converge. For f(x) = g(ax), a real value function of ax (a=constant), the smaller the quantity, $(a\Delta x)^2$, the more probable formula convergence will occur. (A sufficient condition for absolute convergence for $g(ax) = e^{ax}$, sinax, cosax, sinhax, and coshax is $(a\Delta x)^2 < 39$). <br> <u>Comment</u> <br> Even if the function, f(x), satisfies the series convergence criteria, the Sum Formula may diverge or converge then diverge as more series terms are added. This may occur if the rate of change of f(x) at either x value limit increases rapidly in magnitude. f(x) = tanx is such a problematic function. <br> or <br> $$f(x) \equiv \sum_{p=1}^{P} a_p(x+c_p)^{n_p}, \quad p = 1,2,3,\ldots,P$$ <br> where $a_p,c_p$ = real cor complex values <br> $n_p$ = positive integers , $n_p \geq 0$ | Sum Formula form #3 <br> $$\sum_{\Delta x\, x=x_1}^{x_2} f(x) = \frac{1}{\Delta x}\, {}_{\Delta x}\!\int_{x_1}^{x_2} f(x)\Delta x$$ |
| | P = the number of terms of f(x)<br>$B_m$ = constants<br>$B_1 = -\frac{1}{24}$ $B_2 = +\frac{7}{5760}$<br>$B_3 = -\frac{31}{967680}$ $B_4 = +\frac{127}{154828800}$<br>$B_5 = -\frac{73}{3503554560}$ $B_6 = +\frac{1414477}{2678117105664000}$<br>$B_7 = -\frac{8191}{612141052723200}$<br>$B_m = \frac{-C_m}{2^{2m}(2m+1)!}$ , m = 1,2,3,…<br>$\frac{B_{m+1}}{B_m} \approx -\frac{1}{40}$<br>$\lim_{m\to\infty}\frac{B_{m+1}}{B_m} = -.025330295910584\ldots$<br>$B_m = -\frac{1}{2}\frac{B_{2m}}{(2m)!}$ m =1 , $B_{2m}$ = even Bernoulli Constants<br>$B_m \approx -\frac{B_{2m}}{(2m)!}$ m = 2,3,4,… , ∞ | |
| **5d** | $\Delta x\sum_{x=x_1}^{x_2} f(x) \cong \frac{1}{\Delta x}\int_{x_1-\frac{\Delta x}{2}}^{x_2+\frac{\Delta x}{2}} f(x)dx + \sum_{m=1}^{\infty} B_m\,\Delta x^{2m-1}\frac{d^{2m-1}}{dx^{2m-1}}f(x)\Big\|_{x_1-\frac{\Delta x}{2}}^{x_2+\frac{\Delta x}{2}}$<br>accuracy improves rapidly for increasing $\left\|\frac{x}{\Delta x}\right\|$<br>This same formula as in 5c may be valid even if f(x) can not be represented by a Taylor/Maclaurin Series. However, equation accuracy will depend on the quantity, $\left\|\frac{x}{\Delta x}\right\|$ (a larger value improves accuracy). The approximate equality (≅) approaches equality (=) as $\Delta x\to 0$.<br>$x = x_1, x_1+\Delta x, x_1+2\Delta x, \ldots, x_2-\Delta x, x_2$<br>$x_1,x_2,x,\Delta x$ = real or complex values<br>$\frac{x_2 - x_1}{\Delta x}$ = integer<br>$\Delta x$ = x increment<br>f(x) = function of x | Sum Formula<br>form #3<br>$\Delta x\sum_{x=x_1}^{x_2} f(x) =$<br>$\frac{1}{\Delta x}\;\Delta x\int_{x_1}^{x_2} f(x)\Delta x$ |

Necessary Conditions

f(x) has a Laurent Series representation,

$$f(x) = \sum_{n=-\infty}^{\infty} a_n (x+c)^n$$

where $a_n$ , c = real or complex constants,

or

a Mittag - Leffler Series representation,

$$f(x) = f(0) + \sum_{n=1}^{\infty} a_n \left( \frac{1}{x-c_n} + \frac{1}{x+c_n} \right)$$

where $a_n$ , $c_n$ = real or complex constants,

which is convergent for all x on the line in the complex plane delimited by $x_1 - \frac{\Delta x}{2}$ and $x_2 + \frac{\Delta x}{2}$ . Also, the Sum Formula itself must converge.For f(x) = g(ax), a real value function of ax (a=constant), the smaller the quantity, $(a\Delta x)^2$, the more probable formula convergence will occur.

Comment

Even if the function, f(x), satisfies the series convergence criteria, the Sum Formula may diverge or converge then diverge as more series terms are added. This may occur if the rate of change of f(x) at either x value limit increases rapidly in magnitude. f(x) = tanx is such a problematic function.

or

$$f(x) \equiv \sum_{p=1}^{P} a_p (x+c_p)^{n_p} \quad , \quad p = 1,2,3,\ldots P$$

where $a_p, c_p, n_p$ = real or complex values

At least one $n_p$ value is not a positive integer or 0

P = the number of terms of f(x)

| | | |
|---|---|---|
| **6** | $$\sum_{\Delta x\, x=x_1}^{x_2} (-1)^{\frac{x-x_1}{\Delta x}} f(x) = -\frac{1}{2} f(x) \Big|_{x_1}^{x_2+\Delta x} + \sum_{m=1}^{\infty} H_m (2\Delta x)^{2m-1} \frac{d^{2m-1}}{dx^{2m-1}} f(x) \Big|_{x_1}^{x_2+\Delta x}$$ <br> $x = x_1, x_1+\Delta x, x_1+2\Delta x, \ldots, x_2-\Delta x, x_2$ <br> $x_2 = x_1 + (2p-1)\Delta x$ , $p=1,2,3,\ldots$ <br> $x_1,x_2,x,\Delta x$ = real or complex values <br> $\Delta x$ = x increment <br> $f(x)$ = function of x <br> <u>Necessary Conditions</u> <br> f(x) has a Taylor/Maclaurin Series representation which is convergent for all x on the line in the complex plane delimited by $x_1$ and $x_2+\Delta x$. Also, the Altenrating Sign Sum Formula itself must converge. For $f(x) = g(ax)$, a real value function of ax (a=constant), the smaller the quantity, $(a\Delta x)^2$, the more probable formula convergence will occur. (A sufficient condition for absolute convergence for $g(ax) = e^{ax}$, sinax, cosax, sinhax, and coshax is $(2a\Delta x)^2 < 39$). <br> <u>Comment</u> <br> Even if the function, f(x), satisfies the series convergence criteria, the Alternating Sign Sum Formula may diverge or converge then diverge as more series terms are added. This may occur if the rate of change of f(x) at either x value limit increases rapidly in magnitude. f(x) = tanx is such a problematic function. <br> or <br> $$f(x) \equiv \sum_{p=1}^{P} a_p(x+c_p)^{n_p} , \quad p = 1,2,3,\ldots,P$$ <br> where $a_p,c_p$ = real or complex values <br> $n_p$ = positive integers , $n_p \geq 0$ <br> P = the number of terms of f(x) | Alternating Sign Sum Formula <br> $$\sum_{\Delta x\, x=x_1}^{x_2} (-1)^{\frac{x-x_1}{\Delta x}} f(x) =$$ <br> $$\frac{1}{\Delta x} \int_{\Delta x\, x_1}^{x_2} (-1)^{\frac{x-x_1}{\Delta x}} f(x)\, \Delta x$$ |

$H_1 = +\frac{1}{8}$ $\qquad$ $H_2 = -\frac{1}{384}$

$H_3 = +\frac{1}{15360}$ $\qquad$ $H_4 = -\frac{17}{10321920}$

$H_5 = +\frac{31}{743178240}$ $\qquad$ $H_6 = -\frac{691}{653996851200}$

$H_7 = +\frac{5461}{204047017574400}$

$$H_m = \frac{B_{2m} + \frac{C_m}{2^{2m}(2m+1)}}{(2m)!} \quad , \quad m = 1,2,3, \ldots$$

$B_{2m}$ = even Bernoulli Constants

$C_m$ = Interval Calculus Summation Constants

$B_{2m}$, $C_m$, $H_m$ = real constants

$$\frac{H_{m+1}}{H_m} \approx -\frac{1}{40}$$

$$\lim_{m\to\infty}\frac{H_{m+1}}{H_m} = -\ .025330295910584\ldots$$

**6a**

$$\sum_{\Delta x\; x=x_1}^{x_2} (-1)^{\frac{x-x_1}{\Delta x}} f(x) \cong -\frac{1}{2} f(x)\Big|_{x_1}^{x_2+\Delta x} + \sum_{m=1}^{\infty} H_m (2\Delta x)^{2m-1} \frac{d^{2m-1}}{dx^{2m-1}} f(x)\Big|_{x_1}^{x_2+\Delta x}$$

accuracy improves rapidly for increasing $\left|\frac{x}{\Delta x}\right|$

This same formula as in 6 may be valid even if f(x) can not be represented by a Taylor/Maclaurin Series. However, equation accuracy will depend on the quantity, $\left|\frac{x}{\Delta x}\right|$ (a larger value improves accuracy). The approximate equality ($\cong$) approaches equality (=) as $\Delta x \to 0$.

$x = x_1, x_1+\Delta x, x_1+2\Delta x, \ldots, x_2-\Delta x, x_2$

$x_2 = x_1 + (2p-1)\Delta x \ , \qquad p=1,2,3,\ldots$

$x_1, x_2, x, \Delta x$ = real or complex values

$\Delta x$ = x increment

f(x) = function of x

Necessary Conditions

f(x) has a Laurent Series representation,

Alternating Sign Sum Formula

$$\sum_{\Delta x\; x=x_1}^{x_2} (-1)^{\frac{x-x_1}{\Delta x}} f(x) = \frac{1}{\Delta x}\int_{\Delta x\; x_1}^{x_2} (-1)^{\frac{x-x_1}{\Delta x}} f(x)\, \Delta x$$

$$f(x) = \sum_{n=-\infty}^{\infty} a_n (x+c)^n$$

where $a_n, c$ = real or complex constants,

or

a Mittag - Leffler Series representation,

$$f(x) = f(0) + \sum_{n=1}^{\infty} a_n \left( \frac{1}{x-c_n} + \frac{1}{x+c_n} \right)$$

where $a_n, c_n$ = real or complex constants,

which is convergent for all x on the line in the complex plane delimited by $x_1$ and $x_2+\Delta x$. Also, the Sum Formula itself must converge.For f(x) = g(ax), a real value function of ax (a=constant), the smaller the quantity, $(a\Delta x)^2$, the more probable formula convergence will occur.

Comment

Even if the function, f(x), satisfies the series convergence criteria, the Sum Formula may diverge or converge then diverge as more series terms are added. This may occur if the rate of change of f(x) at either x value limit increases rapidly in magnitude. f(x) = tanx is such a problematic function.

or

$$f(x) \equiv \sum_{p=1}^{P} a_p (x+c_p)^{n_p} \quad , \quad p = 1,2,3,\ldots P$$

where $a_p, c_p, n_p$ = real or complex values

At least one $n_p$ value is not a positive integer or 0

P = the number of terms of f(x)

| | | |
|---|---|---|
| **7** | <u>The Discrete Function Taylor Series is:</u><br>$$f(x_s) = \left[ \sum_{n=0}^{r} \frac{1}{n!} D^n_{\Delta x} f(x_0) \prod_{m=0}^{n-1} (x_s - x_0 - m\Delta x) \right] + $$<br>Polynomial series terms<br>$$\frac{1}{r!} \int_{x_0}^{x_s}{}_{\Delta x} D^{r+1}_{\Delta x} f(x) \prod_{m=0}^{r-1} (x_s - x - [m+1]\Delta x)\Delta x$$<br>or Remainder term, $R_r(x_s)$<br>For $x_s \geq x_0 + \Delta x$<br>$$\mathbf{f}(x_s) = \left[ \sum_{n=0}^{r} \frac{1}{n!} D^n_{\Delta x} f(x_0) \prod_{m=0}^{n-1} (x_s - x_0 - m\Delta x) \right] + $$<br>Polynomial series terms<br>$$\frac{1}{r!} \sum_{x=x_0}^{x_s - \Delta x}{}_{\Delta x} D^{r+1}_{\Delta x} f(x) \prod_{m=0}^{r-1} (x_s - x - [m+1]\Delta x)\Delta x$$<br>Remainder term, $R_r(x_s)$<br>where<br>$\Delta x$ = x increment<br>$\Delta x$ is a positive value<br>The discrete integral and summation subscript is the value of $\Delta x$<br>r = 0 , 1, 2, 3, …<br>r = one less than the number of polynomial series terms<br>$x_0$ = initial value of x<br>$x_s = x_0 + N\Delta x$<br>N = 0, 1, 2, 3, …<br>f(x) = function of x<br>$x = x_0, x_0+\Delta x, x_0+2\Delta x, x_0+3\Delta x, \ldots, x_s-\Delta x, x_s$<br>$D^n_{\Delta x} f(x_0)$ = nth discrete derivative of the function, f(x), evaluated at $x_0$<br>$\prod_{m=0}^{-1} (q) = 1$ where q = real number<br>For $x_s = x_0$ the value of the integral is 0<br>For $0 \leq N < r+1$ the remainder term is zero<br>For $N \geq r+1$, in general, the remainder term is non-zero | Discrete Function Taylor Series |
| | <u>Comment</u> – The Discrete Taylor Series is derived from the expansion of $f(x_0+h) - f(x_0) = -{}_{\Delta t}\int_0^h D_{\Delta t} f(x_0+h-t)\,\Delta t$ using the discrete integration by parts formula, then changing the variable from t to x using $x = x_0 = h - t$ and $D^n_{\Delta t} f(x_0+h-t) = (-1)^n D^n_{\Delta x} f(x)$. | |
| **7a** | $$f(x) = \sum_{n=0}^{\infty} a_n [x-x_0]^n_{\Delta x}$$ or $$f(x) = \sum_{n=0}^{\infty} a_n \prod_{m=0}^{n-1} (x-x_0-m\Delta x)$$ $$a_n = \frac{1}{n!} D_{\Delta x}{}^n f(x)\Big\vert_{x=x_0}$$ or $$f(x) = f(x_0) + \frac{D_{\Delta x} f(x_0)}{1!}(x-x_0) + \frac{D_{\Delta x}{}^2 f(x_0)}{2!}(x-x_0)(x-x_0-\Delta x) + \frac{D_{\Delta x}{}^3 f(x_0)}{3!}(x-x_0)(x-x_0-\Delta x)(x-x_0-2\Delta x) + \dots$$ where<br>$D_{\Delta x} f(x) = \dfrac{f(x+\Delta x) - f(x)}{\Delta x}$, discrete derivative of f(x)<br>$D_{\Delta x}{}^n f(x)$ = nth discrete derivative of f(x)<br>$\prod_{m=0}^{-1} (x-x_0-m\Delta x) = 1$<br>$[x-x_0]^0_{\Delta t} = 1$<br>$x = x_0+N\Delta x$ , $N = 0,1,2,3,\dots$<br>f(x) = function of x<br>$x_0$ = initial value of x<br>$\Delta x$ = x increment | Discrete Function Infinite Taylor Series |
| **7b** | $$f(x) = \sum_{n=0}^{M} a_n [x-x_0]^n_{\Delta x}$$ or | Discrete Function Truncated Taylor Series |
| | $f(x) = \sum_{n=0}^{M} a_n \prod_{m=0}^{n-1}(x-x_0-m\Delta x)$<br><br>$a_n = \frac{1}{n!} D_{\Delta x}{}^{n} f(x)\Big\vert_{x=x_0}$<br><br>or<br><br>$f(x) = f(x_0) + \frac{D_{\Delta x} f(x_0)}{1!}(x-x_0) + \frac{D_{\Delta x}{}^{2} f(x_0)}{2!}(x-x_0)(x-x_0-\Delta x)$<br>$+ \frac{D_{\Delta x}{}^{3} f(x_0)}{3!}(x-x_0)(x-x_0-\Delta x)(x-x_0-2\Delta x) + \ldots$<br>where<br>$D_{\Delta x} f(x) = \frac{f(x+\Delta x) - f(x)}{\Delta x}$ , discrete derivative of f(x)<br>$D_{\Delta x}{}^{n} f(x)$ = nth discrete derivative of f(x)<br>$\prod_{m=0}^{-1}(x-x_0-m\Delta x) = 1$<br>$[x-x_0]_{\Delta t}^{0} = 1$<br>$x = x_0+N\Delta x$ , $N = 0,1,2,3,\ldots$<br>f(x) = function of x<br>$x_0$ = initial value of x<br>$\Delta x$ = x increment<br>M+1 = the number of series terms , $M = 0,1,2,3\ldots$<br>For $0 \le N < M+1$ there is no series error<br>For $N \ge M+1$ there is, in general, series error | |
| | **<u>Specific Equations</u>** | |
| **8** | $e_{\Delta x}(a,x) = \sum_{n=0}^{\infty} \frac{a^n}{n!}[x]_{\Delta x}^{n}$<br><br>or | $e_{\Delta x}x = e_{\Delta x}(1,x)$ |
| | $e_{\Delta x}(a,x) = 1 + \sum_{n=1}^{\infty} \frac{a^n}{n!} \prod_{m=0}^{n-1} (x-m\Delta x)$ | |
| **9** | $e_{\Delta x}(ja,x) = (1+ja\Delta x)^{\frac{x}{\Delta x}} = 1 + \sum_{n=1}^{\infty} a^{2n} \frac{(-1)^n}{(2n)!} [x]_{\Delta x}^{2n} + j \sum_{n=1}^{\infty} a^{2n-1} \frac{(-1)^{n-1}}{(2n-1)!} [x]_{\Delta x}^{2n-1}$ <br> or <br> $e_{\Delta x}(ja,x) = (1+ja\Delta x)^{\frac{x}{\Delta x}} = 1 + \sum_{n=1}^{\infty} a^{2n} \frac{(-1)^n}{(2n)!} \prod_{m=0}^{2n-1} (x-m\Delta x) + j \sum_{n=1}^{\infty} a^{2n-1} \frac{(-1)^{n-1}}{(2n-1)!} \prod_{m=0}^{2(n-1)} (x-m\Delta x)$ | $e_{\Delta x} = \cos_{\Delta x}(a,x)$ <br> $+j\sin_{\Delta x}(a,x)$ |
| **10** | $\sin_{\Delta x}(a,x) = \sum_{n=1}^{\infty} a^{2n-1} \frac{(-1)^{n-1}}{(2n-1)!} [x]_{\Delta x}^{2n-1}$ <br> or <br> $\sin_{\Delta x}(a,x) = \sum_{n=1}^{\infty} a^{2n-1} \frac{(-1)^{n-1}}{(2n-1)!} \prod_{m=0}^{2(n-1)} (x-m\Delta x)$ | $\sin_{\Delta x}x = \sin_{\Delta x}(1,x)$ |
| **11** | $\cos_{\Delta x}(a,x) = 1 + \sum_{n=1}^{\infty} a^{2n} \frac{(-1)^n}{(2n)!} [x]_{\Delta x}^{2n}$ <br> or <br> $\cos_{\Delta x}(a,x) = 1 + \sum_{n=1}^{\infty} a^{2n} \frac{(-1)^n}{(2n)!} \prod_{m=0}^{2n-1} (x-m\Delta x)$ | $\cos_{\Delta x}x = \cos_{\Delta x}(1,x)$ |
| **12** | $\sinh_{\Delta x}(a,x) = \sum_{n=1}^{\infty} \frac{a^{2n-1}}{(2n-1)!} [x]_{\Delta x}^{2n-1}$ <br> or | $\sinh_{\Delta x}x = \sinh_{\Delta x}(1,x)$ |
| | $\sinh_{\Delta x}(a,x) = \sum_{n=1}^{\infty} \frac{a^{2n-1}}{(2n-1)!} \prod_{m=0}^{2(n-1)} (x-m\Delta x)$ | |
| **13** | $\cosh_{\Delta x}(a,x) = 1 + \sum_{n=1}^{\infty} \frac{a^{2n}}{2n!} [x]_{\Delta x}^{2n}$ <br> or <br> $\cosh_{\Delta x}(a,x) = 1 + \sum_{n=1}^{\infty} \frac{a^{2n}}{2n!} \prod_{m=0}^{2n-1} (x-m\Delta x)$ | $\cosh_{\Delta x}x = \cosh_{\Delta x}(1,x)$ |
| **14** | $\ln_1(1+x) = \sum_{m=1}^{\infty}(-1)^{m+1}\zeta(m+1)x^m \;,\quad \lvert x\rvert < 1$ <br> where $\zeta(r)$ = Riemann Zeta Function | |
| **15** | $\frac{e_{\Delta x}(s,x)}{(s-a)^m} = \frac{(1+s\Delta x)^{\frac{x}{\Delta x}}}{(s-a)^m} = e_{\Delta x}(a,x) \sum_{n=0}^{\infty} \frac{1}{n!} \frac{[x]_{\Delta x}^{n}}{(1+a\Delta x)^n} \frac{1}{(s-a)^{m-n}}$ <br> where <br> $m = 1,2,3,\ldots$ | Laurent Series of $\frac{e_{\Delta x}(s,x)}{(s-a)^m}$ |
| **16** | The General Exponential Orthogonal Series <br> $f(x) = \sum_{n=0}^{\infty} A_n y_n(x)$ <br> $y_n(x) = \sum_{m=0}^{n} \frac{(-1)^m nCm\Gamma(n+m+k)\Gamma(k)}{\Gamma(n+k)\Gamma(m+k)} e^{-mbx} = \sum_{m=0}^{n} g(n,m)e^{-mbx}$ <br> $g(n,m) = \frac{(-1)^m nCm\Gamma(n+m+k)\Gamma(k)}{\Gamma(n+k)\Gamma(m+k)}$ <br> or | |
| | $y_n(x) = \left[\frac{\Gamma(k)}{\Gamma(n+k)}\right]\frac{1}{p^{k-1}}\frac{d^n}{dp^n}p^{k-1}[p(1-p)]^n$ , $p = e^{-bx}$ <br><br> $A_n = \left[\frac{\Gamma(n+k)}{n!\Gamma(k)}\right]^2(2n+k)b\int_0^\infty e^{-kbx}y_n(x)f(x)dx$ <br><br> or <br><br> $A_n = \left[\frac{\Gamma(n+k)}{n!\Gamma(k)}\right]^2(2n+k)b\sum_{m=0}^{n}g(n,m)F(b[k+m])$ <br><br> where <br><br> $F(s) = L[f(x)] = \int_0^\infty e^{-sx}f(x)dx$ , $s = b[k+m]$ Laplace Transform <br><br> b,k = constants <br> $0 < b < \infty$ <br> $0 < k < \infty$ <br> $nCm = \frac{n!}{m!(n-m)!}$ the number of combinations of n things taken m at a time <br> f(x) = function of x <br> $\Gamma(k)$ = Gamma Function <br> $n = 0, 1, 2, 3, \dots \infty$ <br><br> $f(0) = \sum_{n=0}^{\infty}\frac{(-1)^n A_n}{\left[\frac{\Gamma(n+k)}{n!\Gamma(k)}\right]}$ <br><br> $f(\infty) = \sum_{n=0}^{\infty}A_n$ <br><br> The differential equation from which $y_n(x)$ is derived is as follows: <br><br> $(e^{bx}-1)\frac{d^2}{dx^2}y_n(x) - b[(k-1)e^{bx} - k]\frac{d}{dx}y_n(x) + b^2n[n+k]y_n(x) = 0$ <br><br> The solution is $c_n y_n(x)$ where $c_n$ is a constant | |

# TABLE 8

## Lnd(n,Δx,x) Function Definitions and Relationships

| # | Definition or Relationship | Comments |
|---|---|---|
| **1** | $\displaystyle {}_{\Delta x}\!\int_{x_1}^{x_2} \frac{1}{x^n}\,\Delta x = \Delta x \;{}_{\Delta x}\!\!\sum_{x=x_1}^{x_2-\Delta x} \frac{1}{x^n} = \pm\, \mathrm{lnd}(n,\Delta x,x)\,\Big\|_{x_2}^{x_1}$ , + for n=1, – for n≠1 | Integral to Summation Relationship |
| **2** | $\ln_{\Delta x} x \equiv \mathrm{lnd}(1,\Delta x,x)$ | The $\ln_{\Delta x} x$ form is sometimes used to emphasize the lnd(n,Δx,x) function's natural log similarities when n=1. |
| **3** | $\displaystyle {}_{\Delta x}\!\int \frac{1}{x^n}\,\Delta x = \pm\, \mathrm{lnd}(n,\Delta x,x) + k$ , + for n=1, – for n≠1<br>k = constant of integration | Discrete integration of the function, $\frac{1}{x^n}$ |
| **4** | $\displaystyle D_{\Delta x}\, \mathrm{lnd}(n,\Delta x,x) = \pm \frac{1}{x^n}$ , + for n = 1, – for n≠1 | $\ln_{\Delta x} x \equiv \mathrm{lnd}(1,\Delta x,x)$ |
| **5** | $\displaystyle D_{\Delta x}\, \mathrm{lnd}(n,\Delta x,x+a) = \pm \frac{1}{(x+a)^n}$ , + for n = 1, – for n≠1 | |
| **6** | $\displaystyle \frac{d^m}{dx^m} \ln_{\Delta x} x \Big\|_{x=x_i} = (-1)^{m+1} \Delta x\, m! \;{}_{\Delta x}\!\!\sum_{x=x_i}^{\infty} \frac{1}{x^{m+1}}$<br>where m = integer, m ≥1 | |
| **7** | $\displaystyle \frac{d^m}{dx^m} \ln_{\Delta x} x = \beta(m)(-1)^{m+1} m!\, \mathrm{lnd}(m+1,\Delta x,x)$<br>where m = 0,1,2,3,…<br>$\beta(m) = \begin{cases} -1 & \text{for } m = 0 \\ +1 & \text{for } m \neq 0 \end{cases}$ | |
| **8** | $\displaystyle \frac{d^m}{dx^m} \mathrm{lnd}(n,\Delta x,x) = (-1)^m \frac{\Gamma(n+m)}{\Gamma(n)} \mathrm{lnd}(n+m,\Delta x,x)$<br>where n ≠ 1, m = 0,1,2,3,… | |

| # | Definition or Relationship | Comments |
|---|---|---|
| **9** | $-\mathrm{lnd}(n,\Delta x,x+\Delta x)\Big\|_{x_1}^{x_2} = \mathrm{lnd}(n,-\Delta x,x)\Big\|_{x_2}^{x_1}$ | |
| **10** | <u>n = 1</u><br>$\ln_{\Delta x}x - \ln_{\Delta x}(\Delta x - x) = \begin{cases} 0\ , & \frac{x}{\Delta x} = \text{integer or integer} + .5 \\ -\pi\cot\mathrm{an}\frac{\pi x}{\Delta x}\ , & \frac{x}{\Delta x} \neq \text{integer} \end{cases}$<br>where<br>$x, \Delta x$ = real or complex values<br>$\Delta x$ = x increment | $\ln_{\Delta x}x \equiv \mathrm{lnd}(1,\Delta x,x)$ |
| **11** | <u>n = 2,3,4,5…</u><br>$(-1)^n\mathrm{lnd}(n,\Delta x,x) + \mathrm{lnd}(n,\Delta x,\Delta x - x) = \begin{cases} 0, \quad \frac{x}{\Delta x} = \text{integer or integer} + .5 \\ \quad n = 3,5,7,9\ldots \\ 2\mathrm{lnd}(n,\Delta x,\Delta x) = \frac{2}{(\Delta x)^{n-1}}\mathrm{lnd}(n,1,1) = \frac{2}{(\Delta x)^{n-1}}\zeta(n) \\ \quad \frac{x}{\Delta x} = \text{integer} \\ \quad n = 2,4,6,8\ldots \\ \frac{-\pi}{(n-1)!}\frac{d^{n-1}}{dx^{n-1}}\cot\mathrm{an}\frac{\pi x}{\Delta x}\ , \quad \frac{x}{\Delta x} \neq \text{integer} \\ \quad n = 2,3,4,5\ldots \end{cases}$<br>where<br>$x, \Delta x$ = real or complex values<br>$\Delta x$ = x increment<br>$\zeta(n)$ = Riemann Zeta Function | |
| **12** | $\ln_1 x = \ln_1(1-x)$<br>where x = integer or integer + .5 | |
| **13** | $\ln_1(x-1) = \ln_1 x - \frac{1}{x-1}$ | |

| # | Definition or Relationship | Comments |
|---|---|---|
| **14** | $\mathrm{lnd}(n,\Delta x,x+\Delta x) = \mathrm{lnd}(n,\Delta x,x) - \frac{\Delta x}{x^n}$ | |
| **15** | $\mathrm{lnd}(n,\Delta x,x) = \frac{1}{\Delta x^{n-1}}\,\mathrm{lnd}(n,1,\frac{x}{\Delta x}) = \frac{1}{\Delta x^{n-1}}\,\zeta(n,1,\frac{x}{\Delta x}) = \frac{1}{\Delta x^{n-1}}\,\zeta(n,\frac{x}{\Delta x})$<br>where<br>$x,\Delta x$ = real values with the condition that $x \neq$ positive real value when $\Delta x$ = negative real value<br>n = complex values<br>or<br>$x = \Delta x$ = complex value<br>n = complex value<br>or<br>$x,\Delta x$ = complex values<br>n = integer<br><br>$\zeta(n,x)$ = Hurwitz Zeta Function , $n \neq 1$<br>$\zeta(n,1,x)$ = General Zeta Function , n = 1 is allowed<br><br><u>Note</u> – A complex value, x+jy, can be a real value, x, an imaginary value, jy, or an integer, N+j0.<br><br><u>Comment</u> – For value combinations not specified above, the equality of the stated equation may not be valid. | |
| **16** | ${}_{\Delta x}\int_{\Delta x}^{x} \frac{1}{x}\,\Delta x = \ln_{\Delta x}x$ | |
| **17** | $\lim_{\Delta x\to 0}(\ln_{\Delta x}x_i - \ln_{\Delta x}1) = \ln x_i = \lim_{\Delta x\to 0}\,{}_{\Delta x}\int_{1}^{x_i} \frac{1}{x}\,\Delta x = \int_{1}^{x_i} \frac{1}{x}\,dx$ | |
| **18** | $\ln_{\Delta x}x = \ln_1\frac{x}{\Delta x}$ | |
| **19** | $\mathrm{lnd}(n,\Delta x,\Delta x) = \frac{1}{\Delta x^{n-1}}\,\mathrm{lnd}(n,1,1) = \frac{1}{\Delta x^{n-1}}\,\zeta(n)$<br>where<br>$n,\Delta x$ = complex values<br><br>$\zeta(n)$ = Riemann Zeta Function , $n \neq 1$<br><br><u>Note</u> – A complex value, x+jy, can be a real number, x, an imaginary value, jy, or an integer, N+j0. | |

| # | Definition or Relationship | Comments |
|---|---|---|
| **20** | $lnd(1,\Delta x,x) \equiv ln_{\Delta x}x \approx ln(\frac{x}{\Delta x} - \frac{1}{2}) + \gamma$<br>where<br>$\gamma$ = Euler's Number , .577215664…<br>Accuracy increases fairly rapidly of increasing $\mid\frac{x}{\Delta x}\mid$ | |
| **21** | $ln_{\Delta x}\Delta x = 0$ | |
| **22** | $lnd(n,\Delta x,x) = \frac{1}{(\alpha\Delta x)^{n-1}}\, lnd(n,\alpha,\frac{x}{\alpha\Delta x})$<br>where $\alpha = 1$ for {Re(Δx)>0} or {Re(Δx)=0 and Im(Δx)>0}<br>$\alpha = -1$ for {Re(Δx)<0} or {Re(Δx)=0 and Im(Δx)<0} | |
| **23** | $lnd(n,\Delta x,x) = \frac{1}{2}\,[lnd(n,2\Delta x,x+\Delta x) + lnd(n,2\Delta x,x)]$<br>where $n \neq 1$ | |
| **24** | $lnd(n,\Delta x,x) = \frac{1}{2^{1-n}}\, lnd(n,2\Delta x,2x)$<br>where n = all values | |
| **25** | $lnd(n,\Delta x,x_i) = lnd(n,\Delta x,x_p) - \alpha(n)\Delta x \sum_{\Delta x}{}_{x=x_i}^{x_p-\Delta x} \frac{1}{x^n}$<br>where<br>$x = x_i + m\Delta x$ , m = 0,1,2,3,…<br>$x_i,x_p$ = values of x<br>or in another form<br>$lnd(n,\Delta x,x) = lnd(n,\Delta x,x+M\Delta x) - \alpha(n)\Delta x \sum_{1}{}_{m=1}^{M} \frac{1}{(x+[M-m]\Delta x)^n}$<br>where<br>M = 1,2,3,…<br>Δx = x increment<br>n,Δx,x = real or complex values | |

| # | Definition or Relationship | Comments |
|---|---|---|
| | $\alpha(n) = \begin{cases} +1 & \text{for } n = 1 \\ -1 & \text{for } n \neq 1 \end{cases}$<br><br>Note - Either of these equations can be used with the lnd(n,Δx,x) Series to evaluate the function, lnd(n,Δx,x), when the quantity, $\left\vert \frac{x}{\Delta x} \right\vert$, is not sufficiently large for good series accuracy. | |
| **26** | lnd(n,Δx,x) – lnd(n,-Δx,x-Δx) = K<br><br>where<br>K = constant value for all x on a straight line x locus in the complex plane<br>K can change as a result of a change in n, Δx, or the set of x values<br>$x = x_i + m\Delta x$ , m = integers<br>$x_i$ = value of x<br>Δx = x increment<br>n = all values (including the special case n=1)<br>$n,\Delta x,x,x_i$ = real or complex values | |
| | **Equations and diagrams related to the lnd(n,Δx,x) n≠1 function and the lnd(n,Δx,x) n≠1 Series constants of integration, $K_1,K_2,K_3,k_1,k_2,k_3$** | |
| **27** | $\text{lnd}(n,\Delta x,x) \approx -\sum_{m=0}^{\infty} \frac{\Gamma(n+2m-1)\left(\frac{\Delta x}{2}\right)^{2m} C_m}{\Gamma(n)(2m+1)!\left(x - \frac{\Delta x}{2}\right)^{n+2m-1}} + K \quad , \quad n\neq 1$<br><br>$K = K_1$ or $K_2$ or $K_3$ or $k_1$ or $k_2$ or $k_3$<br><br>$\sum_{\Delta x\, x=x_1}^{x_2} \frac{1}{x^n} = -\frac{1}{\Delta x} \text{lnd}(n,\Delta x,x) \Big\vert_{x_1}^{x_2+\Delta x}$<br><br>The $K_2,K_3,k_1,k_2$ constants may change their value as a result of a change to either n, Δx, or the set of the summation x values (i.e. $x = x_1, x_1+\Delta x, x_1+2\Delta x, x_1+3\Delta x,\ldots, x_2$ ). The $K_1,k_3$ constants always have a 0 value. | |

<table>
<tr><th>#</th><th>Definition or Relationship</th><th>Comments</th></tr>
<tr><td>28</td><td>Re(Δx)≠0 and Im(Δx) ≠ 0<br>X Locus Segment Designations<br>Complex Plane im(x)<br>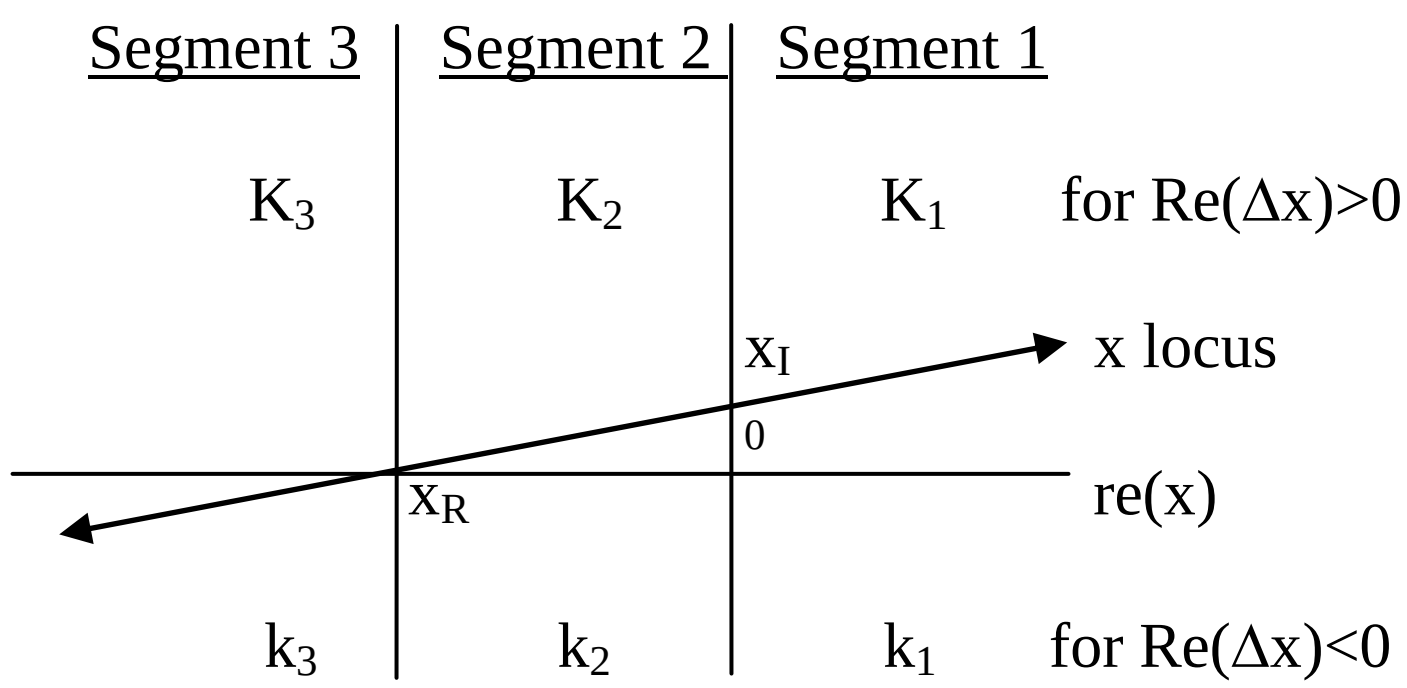<br>$x = x_i + m\Delta x$ , m = integer , This is the x locus<br>$x_i$ = value of x<br>$x_R$ = real axis x locus crossover point<br>$x_I$ = imaginary axis x locus crossover point</td><td>Use the constants, $K_1$, $K_2$, $K_3$ when Re(Δx)>0.<br><br>Use the constants, $k_1$, $k_2$, $k_3$ when Re(Δx)<0.<br><br>The constant of integration subscript indicates the x locus segment in which the constant is used.<br><br>The x locus axis crossover points, $x_R$ and $x_I$, determine the x locus segments.</td></tr>
<tr><td>29</td><td>Re(Δx)≠0 and Im(Δx) = 0<br>X Locus Segment Designations<br><br>Complex Plane im(x)<br>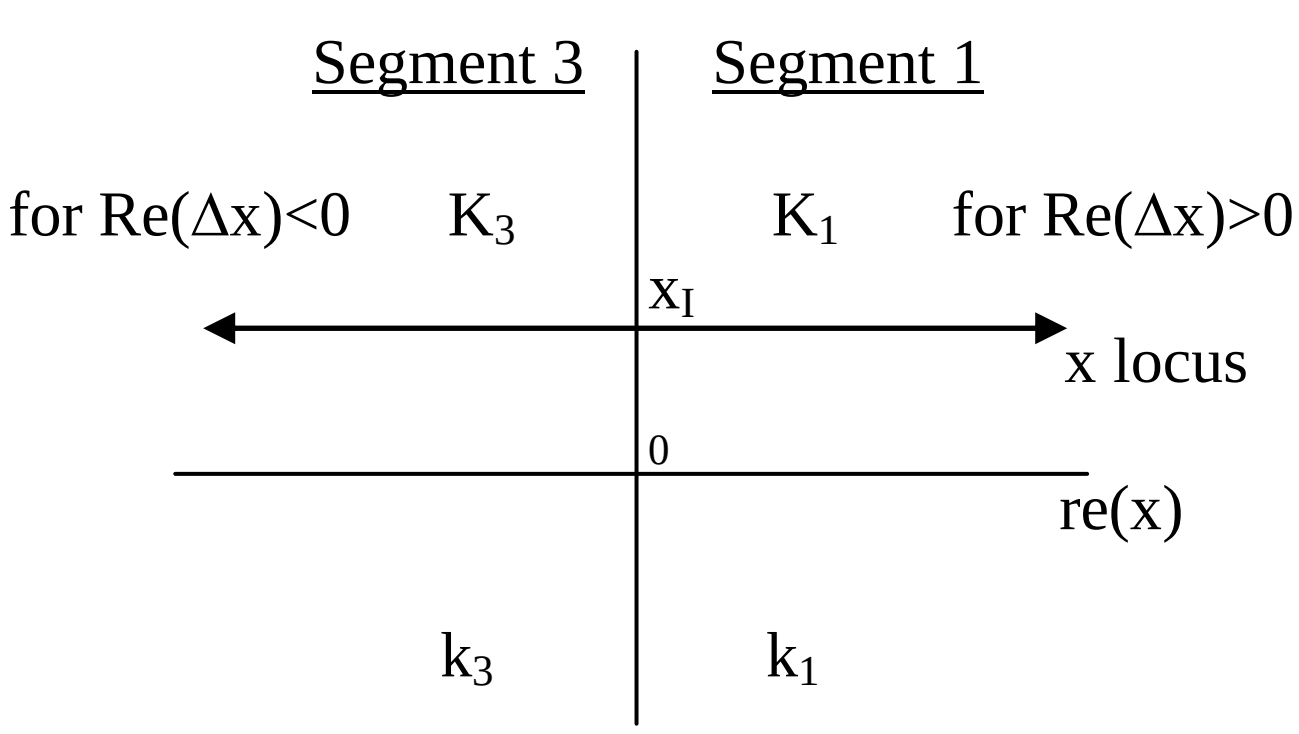<br>$x = x_i + m\Delta x$ , m = integer , This is the x locus<br>$x_i$ = value of x<br>$x_I$ = imaginary axis x locus crossover point</td><td>Use the constants, $K_1$, $K_3$, when Re(Δx)>0.<br><br>Use the constants, $k_1$, $k_3$, when Re(Δx)<0.<br><br>The constant of integration subscript indicates the x locus segment in which the constant is used.<br><br>The imaginary axis crossover point, $x_I$, determines the x locus segments.<br><br>Here there is no x locus segment 2.</td></tr>
</table>

| # | Definition or Relationship | Comments |
|---|---|---|
| **30** | <u>$Re(\Delta x) = 0$ and $Im(\Delta x) \neq 0$</u><br>X Locus Segment Designations<br><br><u>Complex Plane</u><br>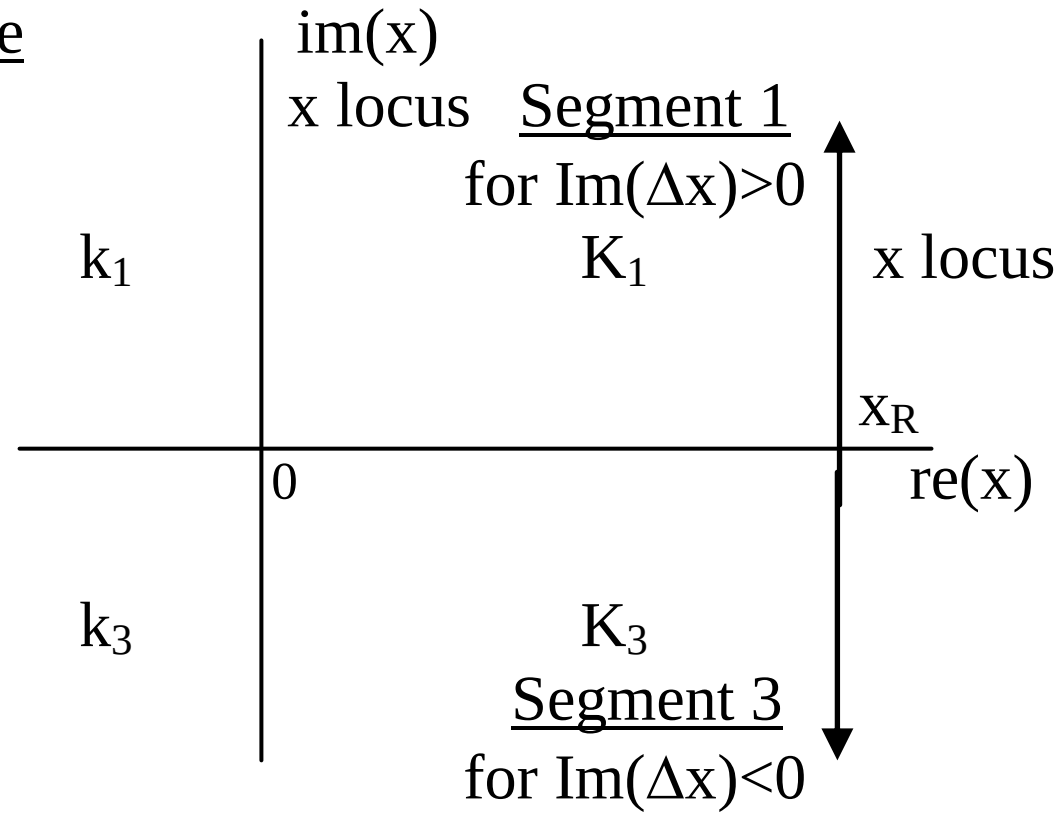<br><br>$x = x_i + m\Delta x$ , $m$ = integer , This is the x locus<br>$x_i$ = value of x<br>$x_R$ = real axis x locus crossover point | Use the constants, $K_1$, $K_3$, when $Re(\Delta x)=0$ and $Im(\Delta x)>0$.<br><br>Use the constants, $k_1$, $k_3$, when $Re(\Delta x)=0$ and $Im(\Delta x)<0$.<br><br>The constant of integration subscript indicates the x locus segment in which the constant is used.<br><br>The x locus real axis crossover point, $x_R$, determines the x locus segments.<br><br>Here there is no x locus segment 2. |
| **31** | <u>Comment</u> on the magnitude, $\left\lvert \frac{x}{\Delta x} \right\rvert$<br><br>The magnitude, $\left\lvert \frac{x}{\Delta x} \right\rvert$ , is often referred to as having to be large for the $lnd(n,\Delta x,x)$ $n \neq 1$ Series to be accurate. The larger the magnitude, the greater the accuracy. The quantity, $\frac{x}{\Delta x}$ , should have a significant real component.<br>Experience with the use of this magnitude has provided some indication as to what the term "large" actually means. It has been found that for good accuracy, $Re(x)$ should be equal to or larger than $100Re(\Delta x)$. | |
| **32** | <u>Relating the K,k constants to the lnd(n,Δx,x) function</u><br><br>$K_r - k_r = K_1 - k_1 = K_2 - k_2 = K_3 - k_3 = lnd(n,\Delta x,x) - lnd(n,-\Delta x,x-\Delta x)$<br>where<br>x is any x value on the x locus<br>$x = x_i + m\Delta x$ , $m$ = integer , This is the x locus<br>$x_i$ = value of x | |

| # | Definition or Relationship | Comments |
|---|---|---|
| | $K_r$ = lnd(n,Δx,x) n≠1 Series constants of integration for Re(Δx) > 0 or Re(Δx) = 0 and Im(Δx) > 0<br>$k_r$ = lnd(n,Δx,x) n≠1 Series constants of integration for Re(Δx) < 0 or Re(Δx) = 0 and Im(Δx) < 0<br>r = 1, 2, or 3<br>n ≠ 1<br>Note – For horizontal or vertical x locii there is no x locus segment 2 | |
| **33** | Calculation of the lnd(n,Δx,x) n≠1 Series constants of integration<br>$K_1 = k_3 = 0$<br>$K_3 = -k_1$ = lnd(n,Δx,x) – lnd(n,-Δx,x-Δx)<br>$K_2 - k_2$ = lnd(n,Δx,x) – lnd(n,-Δx,x-Δx) | There are not enough equations here to solve for $K_2$ and $k_2$. |
| **34** | Equations involving the lnd(n,Δx,x) n≠1 function and its constants of integration<br>1) $\Delta x \sum\limits_{\Delta x\; x=x_1}^{x_2} \frac{1}{x^n} = {}_{\Delta x}\int_{x_1}^{x_2+\Delta x} \frac{1}{x^n}\,\Delta x = -\mathrm{lnd}(n,\Delta x,x)\Big\|_{x_1}^{x_2+\Delta x}$<br>$= -\mathrm{lnd}(n,\Delta x,x_2+\Delta x) + \mathrm{lnd}(n,\Delta x,x_1)]$<br>$x = x_1, x_1+\Delta x, x_1+2\Delta x, x_1+3\Delta x, \ldots, x_2-\Delta x, x_2$<br>$x_1,x_2$ may be infinite values<br>n≠1<br>2) $\Delta x \sum\limits_{\Delta x\; x=x_i}^{\pm\infty} \frac{1}{x^n} = {}_{\Delta x}\int_{x_i}^{\pm\infty} \frac{1}{x^n}\,\Delta x = \mathrm{lnd}(n,\Delta x,x_i) \approx \mathrm{lnd}_f(n,\Delta x,x_i) + \begin{cases} K_r \\ k_r \end{cases}$<br>$x = x_i, x_i+\Delta x, x_i+2\Delta x, x_i+3\Delta x + x_i+4\Delta x, x_i+5\Delta x, \ldots$<br>Re(n)>1<br>Accuracy increases rapidly as $\left\|\frac{x}{\Delta x}\right\|$ increases in value.<br>3) $\Delta x \sum\limits_{\Delta x\; x=\mp\infty}^{\pm\infty} \frac{1}{x^n} = \mathrm{lnd}(n,\Delta x,x_i) - \mathrm{lnd}(n,-\Delta x,x_i-\Delta x) = \pm(K_r - k_r)$<br>Re(n)>1 | These equations apply for all values of n, Δx, and x unless otherwise specified |

| # | Definition or Relationship | Comments |
| --- | --- | --- |
| | 4) $\text{lnd}_f(n,\Delta x,x) = -\sum_{m=0}^{\infty} \frac{\Gamma(n+2m-1)\left(\frac{\Delta x}{2}\right)^{2m} C_m}{\Gamma(n)(2m+1)!\left(x - \frac{\Delta x}{2}\right)^{n+2m-1}}$<br>The absolute value, $\left\lvert\frac{x}{\Delta x}\right\rvert$, must be large for good accuracy | |
| | 5) $\text{lnd}(n,\Delta x,x) \approx \text{lnd}_f(n,\Delta x,x) + \begin{cases} K_r \\ k_r \end{cases}$, x is within the x locus segment, r<br>The absolute value, $\left\lvert\frac{x}{\Delta x}\right\rvert$, must be large for good accuracy | |
| | 6) $\left.\begin{matrix} K_r \\ k_r \end{matrix}\right\} \approx \text{lnd}(n,\Delta x,x) - \text{lnd}_f(n,\Delta x,x)$, x is within the x locus segment, r<br>The absolute value, $\left\lvert\frac{x}{\Delta x}\right\rvert$, must be large for good accuracy | |
| | 7) $\text{lnd}(n,\Delta x,x_i) = \Delta x \sum_{x=x_i}^{x_p-\Delta x}{}_{\Delta x} \frac{1}{x^n} + \text{lnd}_f(n,\Delta x,x_p) + \begin{cases} K_r \\ k_r \end{cases}$ | |
| | 8) $K_1 = 0$ | |
| | 9) $k_3 = 0$ | |
| | 10) $K_3 = -k_1$ | |
| | 11) $\text{lnd}(n,\Delta x,x) - \text{lnd}(n,-\Delta x,x-\Delta x) = \pm(K_r - k_r)$ | |
| | 12) $\text{lnd}(n,\Delta x,x) - \text{lnd}(n,-\Delta x,x-\Delta x) = -[\text{lnd}(n,-\Delta x,x) - \text{lnd}(n,\Delta x,x+\Delta x)]$ | |
| | 13) $k_1 = \pm\left[-\Delta x \sum_{x=\mp\infty}^{\pm\infty}{}_{\Delta x} \frac{1}{x^n}\right]$, $\text{Re}(n)>1$ | |

| # | Definition or Relationship | Comments |
|---|---|---|
| | 14) $K_3 = \pm\,[\Delta x \sum_{\Delta x\; x=\mp\infty}^{\pm\infty} \frac{1}{x^n}]\ ,\quad Re(n)>1$<br><br>15) $K_2 - k_2 = \pm\,[\Delta x \sum_{\Delta x\; x=\mp\infty}^{\pm\infty} \frac{1}{x^n}]\ ,\quad Re(n)>1$<br><br>E16) $K_r - k_r = \pm\,[\Delta x \sum_{\Delta x\; x=\mp\infty}^{\pm\infty} \frac{1}{x^n}]\ ,\quad Re(n)>1\ ,\quad r = 1,2,\text{ or }3$<br><br>E17) $\left.\begin{matrix} K_r \\ k_r \end{matrix}\right\} = \Delta x \sum_{\Delta x\; x=x_i}^{\pm\infty} \frac{1}{x^n} - \ln d_f(n,\Delta x,x_i)\ ,\quad Re(n)>1\ ,\quad r = 1,2,\text{ or }3$<br>$x_i$ is within the x locus segment, r<br>Accuracy increases rapidly as $\|\frac{x}{\Delta x}\|$ increases in value.<br>$x_i$ to $+\infty$ for $Re(\Delta x)>0$ or $\{Re(\Delta x)=0$ and $Im(\Delta x)>0\}$<br>$x_i$ to $-\infty$ for $Re(\Delta x)<0$ or $\{Re(\Delta x)=0$ and $Im(\Delta x)<0\}$<br><br>where<br>$x = x_i + m\Delta x$ , m = integers , This is an x locus in the complex plane<br>$x_i,x_p$ = values of x<br>$\Delta x$ = x increment<br>+ for $Re(\Delta x)>0$ or $\{Re(\Delta x)=0$ and $Im(\Delta x)>0\}$<br>– for $Re(\Delta x)<0$ or $\{Re(\Delta x)=0$ and $Im(\Delta x)<0\}$<br>$-\infty$ to $+\infty$ for $Re(\Delta x)>0$ or $\{Re(\Delta x)=0$ and $Im(\Delta x)>0\}$<br>$+\infty$ to $-\infty$ for $Re(\Delta x)<0$ or $\{Re(\Delta x)=0$ and $Im(\Delta x)<0\}$<br>x = real or complex values<br>$n,\Delta x,x_i,x_p Kr,k_r$ = real or complex constants<br>r = 1,2,3 , The x locus segment designations<br>$K_r$ = constant of integration for $Re(\Delta x)>0$ or $\{Re(\Delta x) = 0$ and $Im(\Delta x) > 0\}$<br>$k_r$ = constant of integration for $Re(\Delta x)<0$ or $\{Re(\Delta x) = 0$ and $Im(\Delta x) < 0\}$<br>$n \neq 1$<br><br>The x locus line is the straight line in the complex plane through all of the plotted x points. | |

| # | Definition or Relationship | Comments |
|---|---|---|
| **35** | $K_r - k_r = K_1 - k_1 = K_2 - k_2 = K_3 - k_3 = \Delta x \sum_{\Delta x \; x=\mp\infty}^{\pm\infty} \frac{1}{x^n} =$ <br><br> $lnd(n,\Delta x,x) - lnd(n,-\Delta x,x-\Delta x) = \int_{\Delta x \; x=\mp\infty}^{\pm\infty} \frac{1}{x^n}\,\Delta x$ <br><br> where <br> $Re(n) > 1$ <br> $x = x_i + m\Delta x$ , m = integer This is the x locus in the complex plane. <br> $x_i$ = value of x <br> $K_r = lnd(n,\Delta x,x)$ $n\neq1$ Series constants of integration for $Re(\Delta x) > 0$ or $Re(\Delta x) = 0$ and $Im(\Delta x) > 0$ <br> $k_r = lnd(n,\Delta x,x)$ $n\neq1$ Series constants of integration for $Re(\Delta x) < 0$ or $Re(\Delta x) = 0$ and $Im(\Delta x) < 0$ <br> r = 1,2,3 <br> $\Delta x$ = x increment <br><br> Comment - The summations, $\sum_{\Delta x \; x=\mp\infty}^{\pm\infty} \frac{1}{x^n}$ and $\sum_{\Delta x \; x=\mp\infty}^{\pm\infty} (-1)^{\frac{x-x_i}{\Delta x}} \frac{1}{x^n}$ are defined entirely in terms of the K,k constants. (The $lnd_f(n,\Delta x,x)$ series terms cancel out. See row 33 above.) | This equation applies only for Re(n)>1 |
| **36** | $lnd(n,\Delta x,x+\Delta x) - lnd(n,\Delta x,x) = \alpha(n)(\frac{1}{x^n})\Delta x$ <br><br> or <br><br> $\frac{1}{\Delta x}[lnd(n,\Delta x,x+\Delta x) - lnd(n,\Delta x,x)] = \alpha(n)(\frac{1}{x^n}) = D_{\Delta x}lnd(n,\Delta x,x)$ <br><br> where <br> $\Delta x$ = x increment <br> $n,\Delta x,x$ = real or complex values <br> $\alpha(n) = \begin{cases} +1 & \text{for } n = 1 \\ -1 & \text{for } n \neq 1 \end{cases}$ <br><br> Note that the value of the function, $\frac{1}{\Delta x}[lnd(n,\Delta x,x+\Delta x) - lnd(n,\Delta x,x)]$, remains | |

| # | Definition or Relationship | Comments |
|---|---|---|
| | constant at the value, $\alpha(n)(\frac{1}{x^n})$, for all values of Δx. In effect, the functions, $\frac{1}{\Delta x}$ [lnd(n,Δx,x+Δx) – lnd(n,Δx,x)] and $D_{\Delta x}$lnd(n,Δx,x), which are equal, are not functions of Δx. | |
| **37** | $D_{\Delta x}$lnd(n,Δx,x) = $D_{-\Delta x}$lnd(n,-Δx,x) = 0<br>or<br>$D_{\Delta x}$lnd(n,Δx,x) – $D_{\Delta x}$lnd(n,-Δx,x-Δx) = 0<br>or<br>[lnd(n,Δx,x+Δx)–lnd(n,Δx,x)] + [lnd(n,-Δx,x-Δx)–lnd(n,-Δx,x)] = 0<br><br>These equations are valid on all segments of an x locus in the complex plane (i.e. Segment 1, Segment 2, and Segment 3). | |
| **38** | <u>n = -1,-2,-3,…</u><br>lnd(n,Δx,x) – lnd(n,Δx,x-Δx) = 0. | This equation is believed to be valid but, presently, there is no formal proof . |
| | **Definitions, hypothesis, and conditions related to lnd(n,Δx,x) n≠1 Series snap characteristic** | |
| **39** | <u>Definition of lnd(n,Δx,x) n≠1 Series Snap</u><br>Lnd(n,Δx,x) n≠1 Series Snap is a characteristic of this series whereby an abrupt and rapid value transition sometimes takes place near the series' x locus complex plane axis crossover point(s). The resulting error in the calculated lnd(n,Δx,x) function value is correctable by a re-evaluation of the series constant of integration. | |
| **40** | <u>The lnd(n,Δx,x) n≠1 Series Snap Hypothesis</u><br><br>Lnd(n,Δx,x) n≠1 Series snap will occur, if it occurs at all, only at an x locus transition across a complex plane axis and, except where series snap occurs, the lnd(n,Δx,x) n≠1 Series constant of integration will not change. | This hypothesis is not yet proven. |
| **41** | <u>Necessary and sufficient condition</u><br>For no lnd(n,Δx,x) n≠1 Series snap to occur, $K_1=K_2=K_3=k_1=k_2=k_3=0$<br>where<br>$x = x_i + m\Delta x$ , m = integer This is the x locus in the complex plane.<br>$x_i$ = value of x<br>The x locus is neither horizontal nor vertical in the complex plane. | |

| # | Definition or Relationship | Comments |
|---|---|---|
| **42** | <u>Necessary and sufficient condition</u><br>For no lnd(n,Δx,x) n≠1 Series snap to occur, $K_1=K_3=k_1=k_3=0$<br>where<br>$x = x_i + m\Delta x$ , m = integer This is the x locus in the complex plane.<br>$x_i$ = value of x<br>The x locus is either horizontal or vertical in the complex plane. | |
| **43** | <u>Necessary condition</u><br>For no lnd(n,Δx,x) n≠1 Series snap to occur, lnd(n,Δx,x)–lnd(n,Δx,x-Δx)=0<br>where<br>$x = x_i + m\Delta x$ , m = integer This is the x locus in the complex plane.<br>$x_i$ = value of x<br>The x locus is neither horizontal nor vertical in the complex plane. | |
| **44** | <u>Necessary and sufficient condition</u><br>For no lnd(n,Δx,x) n≠1 Series snap to occur, lnd(n,Δx,x)–lnd(n,Δx,x-Δx)=0<br>where<br>$x = x_i + m\Delta x$ , m = integer This is the x locus in the complex plane.<br>$x_i$ = value of x<br>The x locus is either horizontal or vertical in the complex plane | |

# TABLE 9

## Miscellaneous Relationships, Equations, and Evaluations

| # | Relationships, Equations, or Evaluations | Comments |
|---|---|---|
| | Relationships | |
| **1** | $\int_{x_1}^{x_2} f(x)dx = \lim_{\Delta x \to 0} {}_{\Delta x}\int_{x_1}^{x_2} f(x)\Delta x$ | Integration |
| **2** | $\frac{d}{dx} = D = \lim_{\Delta x \to 0} D_{\Delta x}$ | Differentiation |
| **3** | ${}_{\Delta x}\sum_{x=x_1}^{x_2} f(x) = \frac{1}{\Delta x}\, {}_{\Delta x}\int_{x_1}^{x_2+\Delta x} f(x)\Delta x$ | Summation |
| **4** | $K_{\Delta x}$ Transform to Z Transform Conversion<br>$Z[f(x)] = F(z) = \frac{z}{T} f(s)\Big\|_{s = \frac{z-1}{T}}$ , $Z[f(x)] = F(z)$ Z Transform<br>$T = \Delta x$<br>$x = nT$ , $n = 0,1,2,3,\ldots$ $K_{\Delta x}[f(x)] = f(s)$ $K_{\Delta x}$Transform | $K_{\Delta x}$ Transform to<br>Z Transform Conversion |
| **5** | Z Transform to $K_{\Delta x}$ Transform Conversion<br>$K_{\Delta x}[f(x)] = f(s) = \frac{\Delta x}{1+s\Delta x} F(z)\Big\|_{z = 1+s\Delta x}$ , $K_{\Delta x}[f(x)] = f(s)$ $K_{\Delta x}$ Transform<br>$\Delta x = T$<br>$x = n\Delta x$ , $n = 0,1,2,3,\ldots$<br>$Z[f(x)] = F(z)$ Z Transform | Z Transform to<br>$K_{\Delta x}$ Transform Conversion |
| **6** | Continuous and discrete exponential function equivalence using the Laplace and $K_{\Delta x}$ Transforms<br><br>$L[e^{ax}] = \frac{1}{s-a}$ , Continuous mathematics $(\Delta x \to 0)$<br>where<br>$x$ = real values<br><br>$K_{\Delta x}[e^{ax}] = \frac{1}{s - \frac{e^{a\Delta x} - 1}{\Delta x}}$ , Discrete mathematics ($\Delta x$ does not go to 0)<br>where<br>$x = x_i + m\Delta x$ , m=integers<br>$x_i$ = a value of x | |

| # | Relationships, Equations, or Evaluations | Comments |
|---|---|---|
| | The continuous $e^{ax}$ function Laplace Transform can be converted to its equivalent discrete function $K_{\Delta x}$ Transform by changing the Laplace Transform denominator constant, a, to $\frac{e^{a\Delta x}-1}{\Delta x}$ . | |
| **7** | Homogeneous Differential Equation<br><br>The following homogeneous differential equation has the same solution irrespective of the value of $\Delta x$, including $\Delta x \to 0$.<br><br>$\left(D_{\Delta x} - \frac{e^{r_1\Delta x}-1}{\Delta x}\right)\left(D_{\Delta x} - \frac{e^{r_2\Delta x}-1}{\Delta x}\right)\left(D_{\Delta x} - \frac{e^{r_3\Delta x}-1}{\Delta x}\right)\ldots\left(D_{\Delta x} - \frac{e^{r_n\Delta x}-1}{\Delta x}\right) y(x) = 0$<br><br>The solution is:<br><br>$y(x) = K_1e^{r_1x} + K_2e^{r_2x} + K_3e^{r_3x} + \ldots + K_ne^{r_nx} =$<br><br>$K_1\left(1 + \left[\frac{e^{r_1\Delta x}-1}{\Delta x}\right]\Delta x\right)^{\frac{x}{\Delta x}} + K_2\left(1 + \left[\frac{e^{r_2\Delta x}-1}{\Delta x}\right]\Delta x\right)^{\frac{x}{\Delta x}} +$<br>$K_3\left(1 + \left[\frac{e^{r_3\Delta x}-1}{\Delta x}\right]\Delta x\right)^{\frac{x}{\Delta x}} + \ldots + K_n\left(1 + \left[\frac{e^{r_n\Delta x}-1}{\Delta x}\right]\Delta x\right)^{\frac{x}{\Delta x}}$<br><br>where<br>$x = x_i + m\Delta x$ m = integers<br>$x_i$ = a value of x<br>$\frac{e^{r_p\Delta x}-1}{\Delta x}$ = the differential equation roots , p = 1,2,3,…n<br>$r_p$ = the differential equation roots for $\Delta x \to 0$<br>$r_p$ = real or complex conjugate values<br>n = order of the differential equation<br>$D_{\Delta x}y(x) = \frac{y(x+\Delta x)-y(x)}{\Delta x}$<br>$\lim_{\Delta x\to 0}D_{\Delta x}y(x) = \lim_{\Delta x\to 0}\frac{y(x+\Delta x)-y(x)}{\Delta x} = \frac{d}{dx}y(x)$<br>The sets of initial conditions are consistent<br><br>Given a homogeneous differential equation with roots $r_1, r_2, r_3, \ldots, r_n$ , and a value of $\Delta x$, a related differential difference equation (with the same solution) can be obtained in accordance with the differential equation presented above. Using this equation with a consistent set of initial conditions, the values of y(x) can be obtained for $x = x_i + m\Delta x$. | |

| # | Relationships, Equations, or Evaluations | Comments |
|---|---|---|
| | Equations | |
| **8** | $\zeta(n,\Delta x,x) = \frac{1}{\Delta x}\,\mathrm{lnd}(n,\Delta x,x)$ , for all n <br><br> $\zeta(n,\Delta x,x_i) = \frac{1}{\Delta x}\,\mathrm{lnd}(n,\Delta x,x_i) = {}_{\Delta x}\sum_{x=x_i}^{\pm\infty}\frac{1}{x^n} = \frac{1}{\Delta x}\,{}_{\Delta x}\int_{x_i}^{\pm\infty}\frac{1}{x^n}\,\Delta x$ , Re(n)>1 <br><br> $+\infty$ for Re(Δx)>0 or {Re(Δx)=0 and Im(Δx)>0} <br> $-\infty$ for Re(Δx)<0 or {Re(Δx)=0 and Im(Δx)<0} <br><br> and <br><br> $\zeta(n,\Delta x,x)\Big\vert_{x_1}^{x_2} = \frac{1}{\Delta x}\,\mathrm{lnd}(n,\Delta x,x)\Big\vert_{x_1}^{x_2} = \pm\,{}_{\Delta x}\sum_{x=x_1}^{x_2-\Delta x}\frac{1}{x^n} = \pm\frac{1}{\Delta x}\,{}_{\Delta x}\int_{x_1}^{x_2}\frac{1}{x^n}\,\Delta x$ , <br><br> + for n = 1, − for n ≠ 1 <br><br> ${}_{\Delta x}\sum_{x=x_1}^{x_1-\Delta x}\frac{1}{x^n} = 0$ <br><br> or <br><br> $\zeta(n,\Delta x,x)\Big\vert_{x_1}^{x_2+\Delta x} = \frac{1}{\Delta x}\,\mathrm{lnd}(n,\Delta x,x)\Big\vert_{x_1}^{x_2+\Delta x} = \pm\,{}_{\Delta x}\sum_{x=x_1}^{x_2}\frac{1}{x^n} = \pm\frac{1}{\Delta x}\,{}_{\Delta x}\int_{x_1}^{x_2+\Delta x}\frac{1}{x^n}\,\Delta x$ <br><br> where <br> Δx = x interval <br> x = real or complex variable <br> $x_i,x_1,x_2,\Delta x,n$ = real or complex constants <br> Any summation term where x = 0 is excluded | General Zeta Function <br><br> General Zeta Function summation to infinity <br><br> General Zeta Function summation between finite limits |
| **8a** | $\zeta(n) = \mathrm{lnd}(n,1,1)$ , n≠1 <br><br> $\zeta(n) = {}_{1}\sum_{x=1}^{\infty}\frac{1}{x^n} = {}_{1}\int_{1}^{\infty}\frac{1}{x^n}\,\Delta x = \mathrm{lnd}(n,1,1)$ , Re(n)>1 <br><br> where <br> ζ(n) is the Riemann Zeta Function <br> n = real or complex value | Riemann Zeta Function <br><br> Riemann Zeta Function summation to infinity |

| # | Relationships, Equations, or Evaluations | Comments |
|---|---|---|
| **8b** | $\zeta(n,x) = \text{lnd}(n,1,x)\ ,\ \ n\neq 1$<br><br>$\zeta(n,x_i) = {}_{1}\sum_{x=x_i}^{\infty}\frac{1}{x^n} = {}_{1}\int_{x_i}^{\infty}\frac{1}{x^n}\,\Delta x = \text{lnd}(n,1,x_i)\ ,\ \ \text{Re}(n)>1$<br><br>where<br>$\zeta(n,x)$ is the Hurwitz Zeta Function<br>$n,x_i$ = real or complex values<br>Any summation term where x = 0 is excluded<br>and<br>$\zeta(n,x)\Big\|_{x_1}^{x_2} = \text{lnd}(n,1,x)\Big\|_{x_1}^{x_2} = -\,{}_{1}\sum_{x=x_1}^{x_2-1}\frac{1}{x^n} = -\,{}_{1}\int_{x_1}^{x_2}\frac{1}{x^n}\,\Delta x\ ,$<br><br>${}_{1}\sum_{x=x_1}^{x_1-1}\frac{1}{x^n} = 0$<br>or<br>$\zeta(n,x)\Big\|_{x_1}^{x_2+1} = \text{lnd}(n,1,x)\Big\|_{x_1}^{x_2+1} = -\,{}_{1}\sum_{x=x_1}^{x_2}\frac{1}{x^n} = -\,{}_{1}\int_{x_1}^{x_2+1}\frac{1}{x^n}\,\Delta x$<br><br>where<br>x = real or complex variable<br>$x_i,x_1,x_2,n$ = real or complex constants<br>Any summation term where x = 0 is excluded | Hurwitz Zeta Function<br><br>Hurwitz Zeta Function summation to infinity<br><br>Hurwitz Zeta Function summation between finite limits |
| **8c** | $\zeta(1,\Delta x,x) = \frac{1}{\Delta x}\text{lnd}(1,\Delta x,x) \equiv \frac{1}{\Delta x}\ln_{\Delta x}x\ ,\ \ n=1$<br><br>$\zeta(1,\Delta x,x_f) = \frac{1}{\Delta x}\text{lnd}(1,\Delta x,x_f) \equiv \frac{1}{\Delta x}\ln_{\Delta x}x_f = {}_{\Delta x}\sum_{x=\Delta x}^{x_f-\Delta x}\frac{1}{x} = \frac{1}{\Delta x}\ {}_{\Delta x}\int_{\Delta x}^{x_f}\frac{1}{x}\,\Delta x$<br><br>and<br>$\zeta(1,\Delta x,x)\Big\|_{x_1}^{x_2} = \frac{1}{\Delta x}\text{lnd}(1,\Delta x,x)\Big\|_{x_1}^{x_2} \equiv \frac{1}{\Delta x}\ln_{\Delta x}x\ \Big\|_{x_1}^{x_2} = {}_{\Delta x}\sum_{x=x_1}^{x_2-\Delta x}\frac{1}{x} = \frac{1}{\Delta x}\ {}_{\Delta x}\int_{x_1}^{x_2}\frac{1}{x}\,\Delta x$<br><br>where | N=1 Zeta Function |

| # | Relationships, Equations, or Evaluations | Comments |
|---|---|---|
| | $\Delta x$ = x interval<br>x = real or complex variable<br>$x_1, x_2, x_f, \Delta x$ = real or complex constants<br>Any summation term where x = 0 is excluded<br><br>$\text{lnd}(1,\Delta x,\Delta x) = 0$<br>$\text{lnd}(1,\Delta x,0) = 0$<br><br>$\ln_{\Delta x}x$ is an optional form used to emphasize the similarity to the natural logarithm. | |
| **9** | $\zeta(n,\Delta x,x) = \frac{1}{\Delta x}\,\text{lnd}(n,\Delta x,x) = \pm\frac{1}{\Delta x}\,[\int_{\Delta x} \frac{1}{x^n}\,\Delta x + k(n,\Delta x)]$ ,<br><br>– for $n\neq 1$, + for $n=1$ | General Zeta Function indefinite integral |
| **10** | $\zeta(0,\Delta x,x) = \frac{1}{\Delta x}\,\text{lnd}(0,\Delta x,x) = \frac{1}{\Delta x}\,[-x + \frac{\Delta x}{2}]$ | |
| **11** | $\zeta(-1,\Delta x,x) = \frac{1}{\Delta x}\,\text{lnd}(-1,\Delta x,x) = \frac{1}{\Delta x}\,[-\frac{x(x-\Delta x)}{2} - \frac{\Delta x^2}{12}]$ | |
| **12** | $\zeta(-2,\Delta x,x) = \frac{1}{\Delta x}\,\text{lnd}(-2,\Delta x,x) = \frac{1}{\Delta x}\,[-\frac{x(x-\Delta x)(x-2\Delta x)}{3} - \Delta x\left(\frac{x(x-\Delta x)}{2}\right)]$ | |
| **13** | $\frac{d}{dn}\zeta(0,x) \equiv \frac{d}{dn}\zeta(n,x)\Big\vert_{n=0} = \lim_{n\to 0}\frac{\text{lnd}(n,1,x) - \text{lnd}(0,1,x)}{n}$<br><br><u>Note</u> - $\ln\Gamma(x) - \frac{1}{2}\ln(2\pi) = \frac{d}{dn}\zeta(0,x)$ | |
| **14** | $\frac{d}{dx}\,\text{lnd}(n,\Delta x,x) = -n\,\text{lnd}(n+1,\Delta x,x)$ , $n\neq 0$ | Derivative of the $\text{lnd}(n,\Delta x,x)$ Function |
| **15** | $D_{\Delta x}\text{lnd}(n,\Delta x,x) = \pm\frac{1}{x^n}$ , – for $n\neq 1$, + for $n=1$ | Discrete derivative of the $\text{lnd}(n,\Delta x,x)$ Function |
| **16** | $\frac{d}{dx}\zeta(n,\Delta x,x) = -n\,\zeta(n+1,\Delta x,x)$ , $n\neq 0$ | Derivative of the General Zeta Function |
| **17** | $D_{\Delta x}\zeta(n,\Delta x,x) = \pm[\frac{1}{\Delta x}]\,\frac{1}{x^n}$ , – for $n\neq 1$, + for $n=1$<br><br><u>Note</u> - $\zeta(n,\Delta x,x) = \frac{1}{\Delta x}\,\text{lnd}(n,\Delta x,x)$ | Discrete derivative of the General Zeta Function |
| **18** | $\frac{d}{dx}\zeta(n,x) = -n\,\zeta(n+1,x)$ , $n\neq 0,1$ | Derivative of the Hurwitz Zeta Function |

| # | Relationships, Equations, or Evaluations | Comments |
|---|---|---|
| **19** | $\displaystyle \mathop{\int}_{\Delta x}\limits_{x_1}^{x_2} \frac{1}{(x+a)^n}\,\Delta x = \pm\,\Delta x\,\zeta(n,\Delta x,x+a)\Big\vert_{x_1}^{x_2} = \pm\,\text{lnd}(n,\Delta x,x+a)\Big\vert_{x_1}^{x_2}$ | Discrete definite Integral of the function, $\frac{1}{(x+a)^n}$<br>General Zeta Function |
| **20** | $\displaystyle {}_{\Delta x}\!\int \frac{1}{x^n}\,\Delta x = \pm\;\Delta x\,\zeta(n,\Delta x,x) + k(n,\Delta x)$ ,   $-$ for $n\neq 1$, $+$ for $n=1$ | Indefinite discrete integral of the function, $\frac{1}{x^n}$<br>General Zeta Function |
| **20a** | $D_{\Delta x}[\Delta x\,\zeta(n,\Delta x,x)] = \pm\frac{1}{x^n}$ ,   $-$ for $n\neq 1$, $+$ for $n=1$ | Discrete derivative of the General Zeta Function |
| **20b** | $\displaystyle {}_{1}\!\int \frac{1}{x^n}\,\Delta x = -\zeta(n,x) + k(n)\ ,\ \ n \neq 1$ | Indefinite discrete integral of the function, $\frac{1}{x^n}$<br>Hurwitz Zeta Function |
| **20c** | $D_1\zeta(n,x) = -\frac{1}{x^n}$ | Discrete derivative of the Hurwitz Zeta Function |
| **20d** | $\displaystyle {}_{1}\!\int \frac{1}{x}\,\Delta x = \psi(x) + k$ | Indefinite discrete integral of the function, $\frac{1}{x}$<br>Digamma Function |
| **20e** | $D_1\psi(x) = \frac{1}{x}$ | Discrete derivative of the Digamma Function |
| **21** | $\psi(x) = \text{lnd}(1,1,x) - \gamma\ ,\ \ x \neq 0,-1,-2,-3,\ldots$<br>or<br>$\psi(x) = \frac{d}{dx}\ln\Gamma(x) = \frac{\Gamma(x)'}{\Gamma(x)}$<br>where<br>$\gamma$ = Euler's Constant, .5772157…<br>$\Gamma(x)$ = Gamma Function<br>$\Gamma(x)' = \frac{d}{dx}\Gamma(x)$ | Digamma Function |
| **21a** | $\psi(\frac{x}{\Delta x}) = \text{lnd}(1,\Delta x,x) - \gamma = \text{lnd}(1,1,\frac{x}{\Delta x}) - \gamma,\ \ \frac{x}{\Delta x} \neq 0,-1,-2,-3,\ldots$<br>where<br>$\gamma$ = Euler's Constant, .5772157…<br>$\Delta x$ = real or complex constant | |

| # | Relationships, Equations, or Evaluations | Comments |
|---|---|---|
| **22** | $\psi(x) \approx \ln(x - \frac{1}{2})$<br>Accuracy increases fairly rapidly for increasing $\lvert x \rvert$<br>Let $\lvert x \rvert \geq 10$ | |
| **22a** | $\ln x \approx \psi(x + \frac{1}{2})$<br>Accuracy increases fairly rapidly for increasing $\lvert x \rvert$<br>Let $\lvert x \rvert \geq 10$ | |
| **23** | $\sum_{\Delta x\, x=x_1}^{x_2} \frac{1}{x} = \psi(x)\Big\rvert_{x_1}^{x_2+\Delta x} = \frac{1}{\Delta x}[\psi(\frac{x_2}{\Delta x}+1) - \psi(\frac{x_1}{\Delta x})]$<br>where<br>$x_1 \neq 0,-\Delta x,-2\Delta x,-3\Delta x,\ldots$<br>$x_2 \neq -\Delta x,-2\Delta x,-3\Delta x,\ldots$<br>x = real or complex variable<br>$x_1,x_2$ = real or complex constants<br><br>or<br><br>$\sum_{\Delta x\, x=x_1}^{x_2-\Delta x} \frac{1}{x} = \psi(x)\Big\rvert_{x_1}^{x_2} = \frac{1}{\Delta x}[\psi(\frac{x_2}{\Delta x}) - \psi(\frac{x_1}{\Delta x})]$<br>where<br>$x_1,x_2 \neq 0,-\Delta x,-2\Delta x,-3\Delta x,\ldots$<br>x = real or complex variable<br>$x_1,x_2$ = real or complex constants | The Digamma Function summation between finite limits |
| **23a** | $\psi(x) = \zeta(1,1,x) - \gamma = \sum_{1\, r=1}^{x-1} \frac{1}{r} - \gamma = \text{lnd}(1,1,x) - \gamma$<br>where<br>$\sum_{1\, r=1}^{0} \frac{1}{r} = 0$<br>$x = 1, 2, 3, \ldots$<br>$\gamma$ = Euler's Constant, .5772157…<br>$\psi(1) = -\gamma$ | The Digamma Function summation where $x = 1,2,3,\ldots$ |

| # | Relationships, Equations, or Evaluations | Comments |
|---|---|---|
| **23b** | $\displaystyle \Delta x \sum_{x=\Delta x}^{x} \frac{1}{x} \approx \frac{1}{\Delta x}\left[\ \ln\left(\frac{x}{\Delta x}+\frac{1}{2}\right)+\gamma\ \right]$ <br> where <br> $\gamma$ = Euler's Number , .577215664… <br> Accuracy increases fairly rapidly for increasing $\left\|\frac{x}{\Delta x}\right\|$ <br> Let $\left\|\frac{x}{\Delta x}\right\| \geq 10$ | |
| **24** | $\psi^{(m)}(x) = (-1)^{m+1}m!\,\zeta(1+m,x) = (-1)^{m+1}m!\,\text{lnd}(1+m,1,x)$ <br> or <br> $\psi^{(m)}(x) = \frac{d^m}{dx^m}\psi(x) = \frac{d^{m+1}}{dx^{m+1}}\ln\Gamma(x)$ <br> where <br> m = 1,2,3,… <br> x = real or complex values | The Polygamma Functions |
| **25** | $\text{lnd}(n,\Delta x,x+\Delta x) = \text{lnd}(n,\Delta x,x) \pm \frac{\Delta x}{x^n}$ , – for n≠1, + for n=1 | lnd(n,Δx,x) Function Recurrence Equation |
| **25a** | $\zeta(n,\Delta x,x+\Delta x) = \zeta(n,\Delta x,x) \pm \frac{1}{x^n}$ , – for n≠1, + for n=1 | General Zeta Function Recurrence Equation |
| **25b** | $\zeta(n,x+1) = \zeta(n,x) - \frac{1}{x^n}$ <br> or <br> $\zeta(n,1,x+1) = \zeta(n,1,x) - \frac{1}{x^n}$ | Hurwitz Zeta Function Recurrence Equation |
| **25c** | $\psi(x+1) = \psi(x) + \frac{1}{x}$ <br> where <br> $\psi(x) = \text{lnd}(1,1,x) - \gamma$ , Digamma Function <br> $\gamma$ = Euler's Constant, .5772157… | Digamma Function Recurrence Equation |
| **25d** | $\psi^{(n)}(x+1) = \psi^{(n)}(x) + \frac{(-1)^n n!}{x^{n+1}}$ , n=1,2,3,… <br> where <br> $\psi^{(n)}(x) = (-1)^{n+1}n!\,\text{lnd}(n+1,1,x)$ , Polygamma Functions | Polygamma Functions Recurrence Equation |

| # | Relationships, Equations, or Evaluations | Comments |
|---|---|---|
| **26** | $\text{lnd}(1,1,x) \equiv \ln_1 x = {}_1\!\sum_{x=1}^{x-1} \frac{1}{x}$<br>where<br>$x = 2,3,4,5,\ldots$<br>$\text{lnd}(1,1,0) = \text{lnd}(1,1,1) = 0$ | Calculation of $\ln_1 x$, $\Delta x = 1$ |
| **27** | $\text{lnd}(1,-1,x) \equiv \ln_{-1} x = -\, {}_{-1}\!\sum_{x=-1}^{x+1} \frac{1}{x}$<br>where<br>$x = -2,-3,-4,-5,\ldots$<br>$\text{lnd}(1,-1,0) = \text{lnd}(1,-1,-1) = 0$ | Calculation of $\ln_{-1} x$, $\Delta x = -1$ |
| **28** | $\ln_1 x \equiv \text{lnd}(1,1,x) = {}_1\!\sum_{x=1}^{\infty} \frac{x-1}{x(x+x-1)}$<br>where<br>$x \neq 0,-1,-2,-3,\ldots$ | Calculation of $\ln_1 x$<br><br>Convergence is slow |
| **29** | $\text{lnd}(m+1,\Delta x,x) = \frac{\beta(m)(-1)^{m+1}}{m!} \, {}_1\!\sum_{x=1}^{\infty} \frac{d^m}{dx^m}\left[\frac{\frac{x}{\Delta x}-1}{x\left(x+\frac{x}{\Delta x}-1\right)}\right]$<br>where<br>$m = 0,1,2,3,\ldots$<br>$\frac{x}{\Delta x} \neq 0,-1,-2,-3,\ldots$<br>$\beta(m) = \begin{cases} -1 & \text{for } m = 0 \\ +1 & \text{for } m \neq 0 \end{cases}$<br>$\Delta x$ = x increment | Calculation of $\text{lnd}(m+1,\Delta x,x)$<br>$m = 0,1,2,3,\ldots$<br><br>Convergence improves for increasing m |
| **30** | $\ln_{\Delta x} x \equiv \text{lnd}(1,\Delta x,x)$<br>where<br>$\Delta x$ = x increment | Alternate notation for the function, $\text{lnd}(1,\Delta x,x)$. The function, $\text{lnd}(1,\Delta x,x)$, has similarities to the natural logarithm. |

| # | Relationships, Equations, or Evaluations | Comments |
|---|---|---|
| **31** | $\ln_{\Delta x} x_f = \Delta x \,{}_{\Delta x}\!\sum_{x=\Delta x}^{x_f-\Delta x} \frac{1}{x} = \ln_1 \frac{x_f}{\Delta x} = {}_1\!\sum_{x=1}^{\frac{x_f}{\Delta x}-1} \frac{1}{x} = {}_{\Delta x}\!\int_{\Delta x}^{x_f} \frac{1}{x}\,\Delta x$<br>where<br>$x = \Delta x, 2\Delta x, 3\Delta x, \ldots, x_f$<br>$\frac{x}{\Delta x} = 1,2,3,\ldots$<br>$x, x_f, \Delta x$ = real or complex values | Calculation of $\ln_{\Delta x} x$ |
| **32** | $\mathrm{lnd}(n,\Delta x,x) = \frac{1}{\Delta x^{n-1}} \zeta(n,\frac{x}{\Delta x})$ , $n \neq 1$<br>where<br>$\zeta(n,x)$ is the Hurwitz Zeta Function<br>$\Delta x$ = x increment<br>$x,\Delta x$ = real values with the condition that $x \neq$ positive real value when $\Delta x$ = negative real value<br>n = complex value<br>or<br>$x = \Delta x$ = complex values<br>n = complex value<br>or<br>$x,\Delta x$ = complex values<br>n = integer | |
| **33** | $\lim_{\Delta x \to \pm 0} \mathrm{lnd}(n,\Delta x,x) = \frac{1}{(n-1)x^{n-1}}$<br>where<br>$n \neq 1$<br>$\Delta x \to \begin{cases} +0 & \text{for } \mathrm{Re}(x)>0 \\ -0 & \text{for } \mathrm{Re}(x)<0 \end{cases}$ | |
| **34** | $\mathrm{lnd}(1,\Delta x,x) - \mathrm{lnd}(1,\Delta x,x-1) = {}_1\!\sum_{n=2}^{\infty} \mathrm{lnd}(n,\Delta x,x)$<br>where<br>$\|x\|>1$<br>$x,\Delta x$ = real or complex values<br>$\Delta x$ = x increment<br>$n = 2,3,4,\ldots$ | |

| # | Relationships, Equations, or Evaluations | Comments |
|---|---|---|
| **35** | $\mathrm{lnd}(1,\Delta x,x+1) - \mathrm{lnd}(1,\Delta x,x) = {}_1\sum_{n=2}^{\infty} (-1)^n\, \mathrm{lnd}(n,\Delta x,x)$<br><br>where<br>$\|x\|>1$<br>$x,\Delta x$ = real or complex values<br>$\Delta x$ = x increment<br>$n = 2,3,4,\ldots$ | |
| **36** | $\mathrm{lnd}(n,-\Delta x,x-\Delta x) = (-1)^{n+1}\, \mathrm{lnd}(n,\Delta x,\Delta x-x)$<br><br>where<br>n = integers<br>$\Delta x$ = x increment<br>$x = x_i + m\Delta x$ , m = integers<br>$x_i$ = value of x<br>$x,x_i,\Delta x$ = real or complex values | |
| **37** | $\mathrm{lnd}(n,-\Delta x,-x) = (-1)^{n+1}\, \mathrm{lnd}(n,\Delta x,x)$<br><br>where<br>n = integers<br>$\Delta x$ = x increment<br>$x = x_i + m\ x$ , $\Delta$ m = integers<br>$x_i$ = value of x<br>$x,x_i,\Delta x$ = real or complex values | |
| **38** | $\mathrm{lnd}(n,-\Delta x,x-\Delta x) = \mathrm{lnd}(n,\Delta x,x)$<br><br>where<br>$n = -1,-2,-3,-4,\ldots$<br>$\Delta x$ = x increment<br>$x = x_i + m\Delta x$ , m = integers<br>$x_i$ = value of x<br>$x,x_0,\Delta x$ = real or complex values | |
| **39** | $\mathrm{lnd}(n,\Delta x,\Delta x-x) = (-1)^{n+1}\, \mathrm{lnd}(n,\Delta x,x)$<br><br>where<br>$n = -1,-2,-3,-4,\ldots$<br>$\Delta x$ = x increment<br>$x = x_i + m\Delta x$ , m = integers<br>$x_i$ = value of x<br>$x,x_i,\Delta x$ = real or complex values | |

| # | Relationships, Equations, or Evaluations | Comments |
|---|---|---|
| **40** | Definitions of the function, lnd(n,Δx,x)<br><br>$D_{\Delta x}\text{lnd}(1,\Delta x,x) \equiv D_{\Delta x}\ln_{\Delta x}x = +\frac{1}{x}$ ,    n =1<br><br>$D_{\Delta x}\text{lnd}(n,\Delta x,x) = -\frac{1}{x^n}$ ,    n ≠ 1<br><br>$\text{lnd}(1,\Delta x,x) \equiv \ln_{\Delta x}x = \Delta x \ {}_{\Delta x}\!\sum_{\mathit{x}=\Delta x}^{x-\Delta x} \frac{1}{\mathit{x}} = {}_{\Delta x}\!\int_{\Delta x}^{x} \frac{1}{\mathit{x}}\,\Delta \mathit{x}$ ,   n=1<br><br>$\text{lnd}(n,\Delta x,x) = \Delta x \ {}_{\Delta x}\!\sum_{\mathit{x}=x}^{\pm\infty} \frac{1}{\mathit{x}^n} = {}_{\Delta x}\!\int_{x}^{\pm\infty} \frac{1}{\mathit{x}^n}\,\Delta \mathit{x}$ ,   n≠1<br><br>where<br>$\Delta \mathit{x} = \mathit{x}$ increment<br>Δx = the value of $\Delta \mathit{x}$<br>x,Δx,$\mathit{x}$,$\Delta \mathit{x}$,n = real or complex values<br>x to +∞   for Re(Δx)>0 or {Re(Δx)=0 and Im(Δx)>0}<br>x to -∞   for Re(Δx)<0 or {Re(Δx)=0 and Im(Δx)<0} | Definitions of the lnd(n,Δx,x) function |
| **41** | $\text{lnd}(n,\Delta x,x) = \frac{1}{\Delta x^{n-1}}\,\text{lnd}(n,1,\frac{x}{\Delta x})$<br><br>where<br>x,Δx = real values with the condition that x ≠ positive real value<br>when Δx = negative real value<br>n = complex value<br>or<br>x = Δx = complex values<br>n = complex value<br>or<br>x,Δx = complex values<br>n = integer | |
| | Evaluations | |
| **42** | $\ln_{\Delta x}\Delta x \equiv \text{lnd}(1,\Delta x,\Delta x) = 0$ | |
| **43** | $\ln_{\Delta x}0 \equiv \text{lnd}(1,\Delta x,0) = 0$ | |

| # | Relationships, Equations, or Evaluations | Comments |
|---|---|---|
| **44** | $\sum_{n=2}^{\infty}\zeta(n)-1 = 1$<br>$\zeta(n)$ = Riemann Zeta Function | |
| **45** | $\sum_{n=2}^{\infty}(-1)^n[\zeta(n)-1] = \frac{1}{2}$<br>$\zeta(n)$ = Riemann Zeta Function | |
| **46** | $\sum_{n=1}^{\infty}\zeta(2n)-1 = \frac{3}{4}$<br>$\zeta(n)$ = Riemann Zeta Function | |
| **47** | $\sum_{n=1}^{\infty}\zeta(2n+1)-1 = \frac{1}{4}$<br>$\zeta(n)$ = Riemann Zeta Function | |
| **48** | $\sum_{n=1}^{\infty}(-1)^{n+1}[\zeta(2n)-1] = .5766740474685817413405\text{0}$<br>$\zeta(n)$ = Riemann Zeta Function | |
| **49** | $\ln x = \int_1^x \frac{1}{x}\,dx = \ln x\Big\|_1^x = \lim_{\Delta x\to 0} \int_{1\ \Delta x}^{x} \frac{1}{x}\,\Delta x = \lim_{\Delta x\to 0} \ln_{\Delta x} x\Big\|_1^x$<br>$\equiv \lim_{\Delta x\to 0} \text{lnd}(1,\Delta x,x)\Big\|_1^x = \lim_{\Delta x\to 0} \Delta x \sum_{x=1\ \Delta x}^{x-\Delta x} \frac{1}{x}$ , $x \neq$ negative real value<br>$\ln x = \lim_{n\to\infty} [\ \text{lnd}(1,\frac{x-a}{n},x) - \text{lnd}(1,\frac{x-a}{n}, a)\ ] + \ln(a)$ ,  x = any value<br>$a = \begin{cases} +1 & \text{for x = positive real value or} \\ & \text{for x = non-real complex value} \\ -1 & \text{for x = negative real value or} \\ & \text{for x = non-real complex value} \end{cases}$<br>$x \neq 0$ | lnx calculation<br><br>$\ln_{\Delta x}x \equiv \text{lnd}(1,\Delta x,x)$ |

<table>
<tr><th>#</th><th>Relationships, Equations, or Evaluations</th><th>Comments</th></tr>
<tr><td></td><td><u>Note</u> - Contrary to the general rule, $\lim_{\Delta x \to 0} \ln_{\Delta x} x$ is not equal to lnx.</td><td></td></tr>
<tr><td>50</td><td>$\pi = -\text{lnd}(1,4,1+4m) + \text{lnd}(1,-4,-3+4m)$<br>where<br>m = integers</td><td>π calculation</td></tr>
<tr><td>50a</td><td>$\pi = \lim_{m\to 0} [\text{lnd}(1+m,-1,N-1) - \text{lnd}(1+m,1,N)]i$<br>where<br>N = integer or integer + .5</td><td>π calculation</td></tr>
<tr><td>50b</td><td>$\pi = \lim_{m\to 0} [\text{lnd}(1+mi,-1,N-1) - \text{lnd}(1+mi,1,N)]i$<br>where<br>N = integer or integer + .5</td><td>π calculation</td></tr>
<tr><td>50c</td><td>$\pi = \lim_{n\to 0} \frac{1}{1+i} [\text{lnd}(1\pm n,4,1+4m) - \text{lnd}(1\pm n,-4,-3+4m)]$<br>where<br>m = integers</td><td>π calculation</td></tr>
<tr><td>50d</td><td>$\pi = \lim_{n\to 0} \frac{1}{1+i} [\text{lnd}(1\pm ni,4,1+4m) - \text{lnd}(1\pm ni,-4,-3+4m)]$<br>where<br>m = integers</td><td>π calculation</td></tr>
<tr><td>50e</td><td>$\pi = {}_1\sum_{x=m+\frac{5}{4}}^{\infty} \frac{1}{x} + {}_{-1}\sum_{x=m+\frac{1}{4}}^{-\infty} \frac{1}{x}$<br>where<br>m = integers</td><td>π calculation<br><br>slow convergence</td></tr>
<tr><td>51</td><td>$\gamma = \lim_{m\to 0}\left[\text{lnd}(1+m,1,1) - \frac{1}{m}\right]$<br>where<br>γ = .5772157… , Euler’s constant<br>m = real numbers</td><td>γ calculation</td></tr>
</table>

| # | Relationships, Equations, or Evaluations | Comments |
| --- | --- | --- |
| **51a** | $\gamma = \lim_{m\to 0}\left[ \text{lnd}(1+mi,1,1) + \frac{i}{m} \right]$<br>where<br>$\gamma = .5772157\ldots$ , Euler's constant<br>m = real numbers | γ calculation |
| **52** | $\sum_{x=1}^{\infty} {}_{1}(-1)^{x-1}\frac{1}{x} = \ln 2$ | |
| **53** | $\sum_{x=-\infty}^{\infty} {}_{\Delta x}\frac{1}{x^n} = \sum_{x=+\infty}^{-\infty} {}_{-\Delta x}\frac{1}{x^n} = 0$<br>where<br>$x = x_i + m\Delta x$ , m = integers<br>$x_i$ = value of x<br>n = 1,3,5,7,…<br>$\frac{x_i}{\Delta x}$ = integer or integer + .5<br>Any summation term where x = 0 is excluded | |
| **54** | $D_{\Delta x}\left[ \sum_{x=-\infty}^{\infty} {}_{\Delta x}\frac{1}{x^n} \right] = D_{\Delta x}\left[ \sum_{x=+\infty}^{-\infty} {}_{-\Delta x}\frac{1}{x^n} \right] = 0$<br>where<br>$x = x_i + m\Delta x$ , m = integers<br>$x_i$ = value of x<br>n = 2,4,6,8,…<br>$\frac{x_i}{\Delta x}$ = integer or integer + .5<br>Any summation term where x = 0 is excluded<br>$D_{\Delta x}$ is the discrete derivative operator<br>$D_{\Delta x}f(x) = \frac{f(x+\Delta x)-f(x)}{\Delta x}$ | |

| # | Relationships, Equations, or Evaluations | Comments |
|---|---|---|
| **55** | $\ln x = \lim_{m\to 0} \frac{2}{m}\left[\frac{x^m - 1}{x^m + 1}\right]$<br>where<br>m = real numbers<br>$m \to \pm 0$<br>x = real or complex numbers | |
| **55a** | $\ln x = \lim_{m\to 0} \left[\frac{x^m - 1}{m}\right]$<br>where<br>m = real numbers<br>$m \to \pm 0$<br>x = real or complex numbers | For the same value of m, this equation is approximately half as accurate as the equation in row 55 |
| **56** | $\log_b x = \lim_{m\to 0} \frac{2}{m \ln b}\left[\frac{x^m - 1}{x^m + 1}\right]$<br>where<br>m = real numbers<br>$m \to \pm 0$<br>x,b = real or complex numbers | |
| **57** | $e^x = \lim_{n\to\infty} \left[\frac{2 + \frac{x}{n}}{2 - \frac{x}{n}}\right]^n$<br>where<br>n = real numbers<br>$n \to \pm\infty$<br>x = real or complex number | |
| **58** | $b^x = \lim_{n\to\infty} \left[\frac{\frac{2}{\ln b} + \frac{x}{n}}{\frac{2}{\ln b} - \frac{x}{n}}\right]^n$<br>where<br>n = real numbers<br>$n \to \pm\infty$<br>x,b = real or complex numbers | |

| # | Relationships, Equations, or Evaluations | Comments |
|---|---|---|
| **59** | $$\sum_{\Delta x\; x=\mp\infty}^{\pm\infty} \frac{1}{x} = \frac{\pi}{\Delta x}\operatorname{cotan}\left(\frac{\pi x_i}{\Delta x}\right)$$ where<br>$x = x_i + m\Delta x$ , $m$ = integers<br>$x_i$ = value of x<br>$\frac{x_i}{\Delta x} \neq$ integer<br>$x, x_i, \Delta x$ = real or complex values | |
| **60** | $$\lim_{\Delta x\to 0} \sum_{\Delta x\; x=x_1}^{x_2} (-1)^{\frac{x-x_1}{\Delta x}} f(x) = \frac{f(x_1) - f(x_2)}{2}$$ where<br>$\frac{x_2 - x_1}{\Delta x}$ = integer<br>$x_2 = x_1 + (2r-1)\Delta x$ , $r = 1,2,3,\ldots$ | |

# TABLE 10
## Formula Constant Calculation Equations

| | |
|---|---|
| **1** | **$B$n n=0,1,2,3,… Bernoulli Constant Calculation** |

$D_0 = +1$

$D_1 = -\frac{1}{2}$

$D_2 = +\frac{1}{12}$

For n =2,3,4,5,…

$$D_{2n} = -\frac{1}{(2n+1)!} + \frac{1}{2(2n)!} - \sum_{m=1}^{n-1} \frac{D_{2m}}{(2n+1-2m)!}$$

$D_{2n-1} = 0$

$B_0 = +1$

$B_1 = -\frac{1}{2}$

$B_2 = +\frac{1}{6}$

For n =2,3,4,5,…

$B_{2n} = (2n)!D_{2n}$

$B_{2n-1} = 0$

Or calculated from the $C_n$ Constants

$C_0 = -1$

$C_1 = +1$

For n =1,2,3,4,…

$$B_{2n} = \frac{1 - 2n \displaystyle\sum_{m=1}^{n} \frac{C_m \prod_{s=0}^{2(m-1)}(2n-1-s)}{(2m+1)!}}{2^{2n}}$$

| | |
|---|---|
| **2** | **<u>Cn n=0,1,2,3,… Cn Constant Calculation</u>**<br>For n=1,2,3,4,…<br>$C_0 = -1$<br>$B_0 = 1$<br>$B_1 = -\frac{1}{2}$<br>$B_{2n}$ = even Bernoulli Constants<br>$B_{2n+1} = 0$ , odd Bernoulli Constants<br>$$\sum_{m=1}^{n} \frac{C_m \prod_{s=0}^{2(m-1)}(2n-1-s)}{(2m+1)!} = \frac{1-2^{2n}B_{2n}}{2n}$$<br>$$\sum_{m=1}^{n} \frac{C_m \prod_{s=0}^{2(m-1)}(2n-1-s)}{(2m+1)!} = \sum_{m=1}^{n} \frac{C_m(2n-1)!}{(2m+1)![2(n-m)]!}$$ |
| **3** | **<u>An n=0,1,2,3,… An Constant Calculation</u>**<br>$A_n = \frac{-C_n}{2(2n+1)!}$ , n = 1,2,3,…<br>where $C_n$ = $C_n$ Constants |
| **4** | **<u>Bn n=0,1,2,3,… Bn Constant Calculation</u>**<br>$B_n = \frac{-C_n}{2^{2n}(2n+1)!}$ , n = 1,2,3,…<br>where $C_n$ = $C_n$ Constants |
| **5** | **<u>Hn n=0,1,2,3,… Hn Constant Calculation</u>**<br>$H_n = \frac{B_{2n} + \frac{C_n}{2^{2n}(2n+1)}}{(2n)!}$ , n = 1,2,3, …<br>where $B_n$ = Bernoulli Constants<br>$C_n$ = $C_n$ Constants |

# TABLE 11
## Formula Constants

| | |
|---|---|
| **1** | **$C_n$ , n=1,2,3,… Constants**<br>$C_0 = -1$<br>$C_1 = +1$<br>$C_2 = -\frac{7}{3}$<br>$C\ 3 = +\frac{31}{3}$<br>$C\ 4 = -\frac{381}{5}$<br>$C\ 5 = +\frac{2555}{3}$<br>$C\ 6 = -\frac{1414477}{105}$<br>$C\ 7 = +286685$<br>$C\ 8 = -\frac{118518239}{15}$<br>$C\ 9 = +\frac{5749691557}{21}$<br>$C\ 10 = -\frac{640823941499}{55}$<br>$C\ 11 = +\frac{1792042792463}{3}$<br>$C\ 12 = -\frac{9913827341556185}{273}$<br>$C\ 13 = +2582950540044537$<br>$C\ 14 = -\frac{3187598676787461083}{15}$<br>$C\ 15 = +\frac{4625594554880206790555}{231}$<br>$C\ 16 = -\frac{182112049520351725262389}{85}$<br>$C\ 17 = +\frac{7749759553478559072551395}{3}$ |

$$C\ 18 = -\frac{9041858456194752424958344698 91}{25935}$$

$$C\ 19 = +\ 5235038863532181733472316757$$

$$C\ 20 = -\frac{14353174239884589601263410372223 7}{165}$$

$$C\ 21 = +\frac{33427300696841208116528825914877 41}{21}$$

$$C\ 22 = -\frac{73447208499508830514216203097846 7283}{23}$$

$$C\ 23 = +\frac{20985757843117067182095330601636 553591}{3}$$

$$C\ 24 = -\frac{55261733892726368667839334271075 79759250083}{3315}$$

$$C\ 25 = +\frac{47377713207329531930725194944660 08540099675}{11}$$

$$C\ 26 = -\frac{18040648148285033994746255593954 05480381292691}{15}$$

$$C\ 27 = +\frac{14440774190032405703502681013504 24763530320442 2955}{399}$$

$$C\ 28 = -\frac{16972613470291337144855739969414 34891001574251519417}{145}$$

$$C\ 29 = +\frac{12175371827101734898203263637136 40808095644117832 9291}{3}$$

$$C\ 30 = -\frac{70053421038731765784608637372075 77729051163404483391458991 17}{465465}$$

**2** **<u>$A_n$ , n=1,2,3,… Constants</u>**

$$A\ 1 = -\frac{1}{12}$$

$$A\ 2 = +\frac{7}{720}$$

$$A\ 3 = -\frac{31}{30240}$$

$$A\ 4 = +\frac{127}{1209600}$$

$$A\ 5 = -\frac{73}{6842880}$$

$$A\ 6 = +\frac{1414477}{1307674368000}$$

| | |
|---|---|
| | $A\ 7 = -\dfrac{8191}{74724249600}$<br>$A\ 8 = +\dfrac{16931177}{1524374691840000}$<br>$A\ 9 = -\dfrac{5749691557}{5109094217170944000}$<br>$A\ 10 = +\dfrac{91546277357}{802857662698291200000}$<br>$A\ 11 = -\dfrac{3324754717}{287777551824322560000}$<br>$A\ 12 = +\dfrac{1982765468311237}{1693824136731743669452800000}$<br>$A\ 13 = -\dfrac{22076500342261}{186134520519971831808000000}$<br>$A\ 14 = +\dfrac{65053034220152267}{54133236696365522170675200000000}$<br>$A\ 15 = -\dfrac{925118910976041358111}{7597902916460400683578420101120000000}$<br>$A\ 16 = +\dfrac{16555640865486520478399}{1341967268361837003852811862016000000000}$<br>$A\ 17 = -\dfrac{8089941578146657681}{6472037976754201114671946137600000000}$ |
| **3** | <u>**$B_n$ , n=1,2,3,… Constants**</u><br>$B\ 1 = -\dfrac{1}{24}$<br>$B\ 2 = +\dfrac{7}{5760}$<br>$B\ 3 = -\dfrac{31}{967680}$<br>$B\ 4 = +\dfrac{127}{154828800}$<br>$B\ 5 = -\dfrac{73}{3503554560}$<br>$B\ 6 = +\dfrac{1414477}{2678117105664000}$<br>$B\ 7 = -\dfrac{8191}{612141052723200}$ |
| | $B\ 8 = +\frac{16931177}{49950709902213120000}$<br>$B\ 9 = -\frac{5749691557}{669659197233029971968000}$<br>$B\ 10 = +\frac{91546277357}{420928638260761696665600000}$<br>$B\ 11 = -\frac{3324754717}{603513268363481705349120000}$<br>$B\ 12 = +\frac{1982765468311237}{14208826703980998799521113702400000}$<br>$B\ 13 = -\frac{22076500342261}{6245638111639999472316973056000000}$<br>$B\ 14 = +\frac{65053034220152267}{726564003867240624328165356994560000000}$<br>$B\ 15 = -\frac{925118910976041358111}{407909306804755512687816982320742662144000000}$<br>$B\ 16 = +\frac{16555640865486520478399}{2881852764958272713015226472499792314368000000000}$<br>$B\ 17 = -\frac{8089941578146657681}{5559438289725860403624550885933103185920000000}$ |
| **4** | **<u>$H_n$ , n=1,2,3,… Constants</u>**<br>$H\ 1 = +\frac{1}{8}$<br>$H\ 2 = -\frac{1}{384}$<br>$H\ 3 = +\frac{1}{15360}$<br>$H\ 4 = -\frac{17}{10321920}$<br>$H\ 5 = +\frac{31}{743178240}$<br>$H\ 6 = -\frac{691}{653996851200}$<br>$H\ 7 = +\frac{5461}{204047017574400}$<br>$H\ 8 = -\frac{929569}{1371195958099968000}$ |
| | $H\ 9 = +\frac{3202291}{186482650301595648000}$ |
| | $H\ 10 = -\frac{221930581}{5102165312251656929280000}$ |
| | $H\ 11 = +\frac{4722116521}{428581886229139182059520000}$ |
| | $H\ 12 = -\frac{56963745931}{204105820641830048114933760000}$ |
| | $H\ 13 = +\frac{14717667114151}{2081879370546666490772324352000000}$ |
| | $H\ 14 = -\frac{2093660879252671}{11691834544990079012177373560832000000}$ |
| | $H\ 15 = +\frac{86125672563201181}{18987539301063888315776054662791168000000}$ |
| | $H\ 16 = -\frac{129848163681107301953}{1130138339199322632554990773529330319360000000}$ |
| | $H\ 17 = +\frac{868320396104950823611}{2983565215486211749945175642117432043110400000000}$ |

# TABLE 12

## Differential Difference Equation Solutions

| | |
|---|---|
| 1 | **$D_{\Delta x}^1 y(x) + Ay(x) = 0$  First Order Differential Difference Equation with root, a**<br><br>$y(x) = Ke_{\Delta x}(a,x) = K(1+a\Delta x)^{\frac{x}{\Delta x}}$<br>where<br>$x = m\Delta x$ ,  m = integers<br>K,a = real constants<br>a = equation root<br>$\Delta x$ = x increment |
| 2 | **$D_{\Delta x}^2 y(x) + AD_{\Delta x}^1 y(x) + By(x) = 0$  Second Order Differential Difference Equation with real roots, a,b**<br><br>$y(x) = K_1 e_{\Delta x}(a,x) + K_2 e_{\Delta x}(b,x) = K_1(1+a\Delta x)^{\frac{x}{\Delta x}} + K_2(1+b\Delta x)^{\frac{x}{\Delta x}}$<br>where<br>$x = m\Delta x$ ,  m = integers<br>$K_1,K_2$,a,b = real constants<br>a,b = equation real roots<br>$\Delta x$ = x increment |
| 3 | **$D_{\Delta x}^2 y(x) + AD_{\Delta x}^1 y(x) + By(x) = 0$  Second Order Differential Difference Equation with complex roots, a+jb, a-jb and $1+a\Delta x = 0$**<br><br>$y(x) = (b\Delta x)^{\frac{x}{\Delta x}} [\ K_1 \cos \frac{\pi x}{2\Delta x} + K_2 \sin \frac{\pi x}{2\Delta x}\ ]$<br>where<br>$x = m\Delta x$ ,  m = integers<br>$K_1,K_2$,a,b = real constants<br>$1+a\Delta x = 0$<br>a+jb, a-jb = equation roots<br>$\Delta x$ = x increment |
| **4** | **$D_{\Delta x}^{2}y(x) + AD_{\Delta x}^{1}y(x) + By(x) = 0$ Second Order Differential Difference Equation with complex roots, a+jb, a-jb and $1+a\Delta x \neq 0$**<br><br>$y(x) = e_{\Delta x}(a,x)\,[\,K_1\cos_{\Delta x}(\frac{b}{1+a\Delta x}, x) + K_2\sin_{\Delta x}(\frac{b}{1+a\Delta x}, x)\,]$<br><br>or<br><br>$y(x) = (1+a\Delta x)^{\frac{x}{\Delta x}}\,[\,K_1\cos_{\Delta x}(\frac{b}{1+a\Delta x}, x) + K_2\sin_{\Delta x}(\frac{b}{1+a\Delta x}, x)\,]$<br><br>or<br><br>$y(x) = (1+a\Delta x)^{\frac{x}{\Delta x}}\,[\,K_1(1+[\frac{b\Delta x}{1+a\Delta x}]^2)^{\frac{x}{2\Delta x}}\cos\frac{\beta x}{\Delta x} + K_2(1+[\frac{b\Delta x}{1+a\Delta x}]^2)^{\frac{x}{2\Delta x}}\sin\frac{\beta x}{\Delta x}\,]$<br><br>where<br>$x = m\Delta x$ ,  m = integers<br>$K_1,K_2,a,b$ = real constants<br>$\beta = \begin{cases} \tan^{-1}\frac{b\Delta x}{1+a\Delta x} & \text{for } 1+a\Delta x \geq 0 \\ \pi+\tan^{-1}\frac{b\Delta x}{1+a\Delta x} & \text{for } 1+a\Delta x < 0 \end{cases}$<br>$1+a\Delta x \neq 0$<br>a+jb, a-jb = equation roots<br>$\Delta x$ = x increment |

# TABLE 13

## The Undetermined Coefficient Method for Solving Differential Difference Equations

**f(x) - The Solution of a Differential Difference Equation**

For the differential difference equation

$D_{\Delta x}{}^{n}f(x) + a_{n-1}D_{\Delta x}{}^{n-1}f(x) + a_{n-2}D_{\Delta x}{}^{n-2}f(x) + \ldots + a_1D_{\Delta x}f(x) + a_0f(x) = Q(x)$

$f(x) = f_C(x) + f_P(x)$ , the differential difference equation general solution

with the characteristic polynomial, h(r),

$h(r) = r^n + a_{n-1}r^{n-1} + a_{n-2}\,r^{n-2} + \ldots + a_2r^2 + a_1r + a_0 = (r-r_n)(r-r_{n-1})(r-r_{n-2})\ldots(r-r_3)(r-r_2)(r-r_1)$

**fc(x) - The Complementary Solution of the Differential Difference Equation**

For the related homogeneous differential difference equation

$D_{\Delta x}{}^{n}f_C(x) + a_{n-1}D_{\Delta x}{}^{n-1}f_C(x) + a_{n-2}D_{\Delta x}{}^{n-2}f_C(x) + \ldots + a_1D_{\Delta x}f_C(x) + a_0f_C(x) = 0$

where

n = order of the homogeneous equation
$a_{n-1}, a_{n-1}, \ldots, a_0$ = real constants
$f_C(x)$ = complementary solution, the general solution to the related homogeneous differential difference equation
$D_{\Delta x}{}^{n}f(x)$ = nth discrete derivative of the function, f(x)
$x = x_i + p\Delta x$ , p = integers
$x_i$ = initial value of x
$\Delta x$ = x increment

with the characteristic polynomial:

$h(r) = r^n + a_{n-1}r^{n-1} + a_{n-2}\,r^{n-2} + \ldots + a_2r^2 + a_1r + a_0 = (r-r_n)(r-r_{n-1})(r-r_{n-2})\ldots(r-r_3)(r-r_2)(r-r_1) = 0$

with n roots

$r = r_1, r_2, r_3, \ldots, r_{n-1}, r_n$

group the roots and catagorize them in the following way:

v = one of the different unique root or complex conjugate root pair values
m = multiplicity of each different unique root or complex comjugate root pair value, m = 1,2,3,…
g = number of groups of different unique root and root pair values ($g \le n$)

Note – A root of multiplicity 1 is a single root, a root of multiplicity 2 is a double root, etc.

W(v,m,x) = real value general solution function associated with the root(s), v, of multiplicity, m

The W(v,m,x) functions are shown in the table below

$$f_C(x) = \sum_{s=1}^{g} W_s(v_s,m_s,x)$$ , the general solution to the homogeneous differential difference equation

In the following table

$a,b,A_pB_p$ = real constants

$$[x]_{\Delta x}^{0} = 1$$

$$[x]_{\Delta x}^{q} = \sum_{u=1}^{q}(x-[u-1]\Delta x)\ ,\ \ q = 1,2,3,\ldots$$

$x = x_i + p\Delta x$ , p =integers

$x_i$ = initial value of x

$\Delta x$ = x increment

The $A_p$, $B_p$ constants in the differential difference equation general solution, f(x), are evaluated from the problem specified initial conditions:

$$f(0),\ f(\Delta x),\ f(2\Delta x),\ f(3\Delta x),\ \ldots,\ f(n\Delta x)$$

or

$$f(0),\ D_{\Delta x}f(0),\ D_{\Delta x}^{2}f(0),\ D_{\Delta x}^{3}f(0),\ \ldots,\ D_{\Delta x}^{n}f(0)$$

The Possible h(r) = 0 Root Outcomes

| v Root/Root Pair Value | W(v,m,x) Function | W(v,m,x) Calculation Function |
|---|---|---|
| a | $[\sum_{p=1}^{m} A_p[x]_{\Delta x}^{p-1}]\ e_{\Delta x}(a,x)$ | $e_{\Delta x}(a,x) = [1+a\Delta x]^{\frac{x}{\Delta x}}$ |
| $0 \pm jb$ | $[\sum_{p=1}^{m} A_p[x]_{\Delta x}^{p-1}]\ \sin_{\Delta x}(b,x) +$ $[\sum_{p=1}^{m} B_p[x]_{\Delta x}^{p-1}]\ \cos_{\Delta x}(b,x)$ | $\sin_{\Delta x}(b,x) = \frac{[1+jb\Delta x]^{\frac{x}{\Delta x}} - [1-jb\Delta x]^{\frac{x}{\Delta x}}}{2j}$ $\cos_{\Delta x}(b,x) = \frac{[1+jb\Delta x]^{\frac{x}{\Delta x}} + [1-jb\Delta x]^{\frac{x}{\Delta x}}}{2}$ |
| $-\frac{1}{\Delta x} \pm jb$ | $[\sum_{p=1}^{m} A_p[x]_{\Delta x}^{p-1}][b\Delta x]^{\frac{x}{\Delta x}} \sin\frac{\pi x}{2\Delta x}$ + $[\sum_{p=1}^{m} B_p[x]_{\Delta x}^{p-1}][b\Delta x]^{\frac{x}{\Delta x}} \cos\frac{\pi x}{2\Delta x}$ | $[b\Delta x]^{\frac{x}{\Delta x}} \sin\frac{\pi x}{2\Delta x} = \frac{[+jb\Delta x]^{\frac{x}{\Delta x}} - [-jb\Delta x]^{\frac{x}{\Delta x}}}{2j}$ $[b\Delta x]^{\frac{x}{\Delta x}} \cos\frac{\pi x}{2\Delta x} = \frac{[+jb\Delta x]^{\frac{x}{\Delta x}} + [-jb\Delta x]^{\frac{x}{\Delta x}}}{2}$ |
| $a \pm jb$ $a \neq -\frac{1}{\Delta x}$ | $[\sum_{p=1}^{m} A_p[x]_{\Delta x}^{p-1}]\ e_{\Delta x}(a,x) \sin_{\Delta x}(\frac{b}{1+a\Delta x}, x)$ + $[\sum_{p=1}^{m} B_p[x]_{\Delta x}^{p-1}]\ e_{\Delta x}(a,x) \cos_{\Delta x}(\frac{b}{1+a\Delta x}, x)$ | $e_{\Delta x}(a,x) \sin_{\Delta x}(\frac{b}{1+a\Delta x}, x) = \frac{[1+(a+jb)\Delta x]^{\frac{x}{\Delta x}} - [1+(a-jb)\Delta x]^{\frac{x}{\Delta x}}}{2j}$ $e_{\Delta x}(a,x) \cos_{\Delta x}(\frac{b}{1+a\Delta x}, x) = \frac{[1+(a+jb)\Delta x]^{\frac{x}{\Delta x}} + [1+(a-jb)\Delta x]^{\frac{x}{\Delta x}}}{2}$ |

**$f_P(x)$ - The Particular Solution of the Differential Difference Equation**

if $Q(x) = [\sum_{p=0}^{N} A_p[x]_{\Delta x}^{p}]\ e_{\Delta x}(u,x) \sin_{\Delta x}(w,x) + [\sum_{p=0}^{N} B_p[x]_{\Delta x}^{p}]\ e_{\Delta x}(u,x) \cos_{\Delta x}(w,x)$, the particular solution to this equation is:

$$f_P(x) = [\sum_{p=N_h}^{N+N_h} \mathcal{A}_p[x]_{\Delta x}^{p}]\ e_{\Delta x}(u,x) \sin_{\Delta x}(w,x) + [\sum_{p=N_h}^{N+N_h} \mathcal{B}_p[x]_{\Delta x}^{p}]\ e_{\Delta x}(u,x) \cos_{\Delta x}(w,x)$$

where

$u \neq -\frac{1}{\Delta x}$

N = the polynomial order of the summations in Q(x)

$N_h$ = the number of times the Q(x) complex conjugate related roots, $u+jw(1+u\Delta x)$ and $u-jw(1+u\Delta x)$, appear in the characteristic polynomial, h(r)

$A_p, B_p, \mathcal{A}_p, \mathcal{B}_p$, u, w = real value constants

n = order of the differential difference equation

$f_P(x)$ = particular solution to a differential difference equation

$\Delta x$ = x increment

or

if $Q(x) = [\sum_{p=0}^{N} A_p[x]_{\Delta x}^{p}][w\Delta x]^{\frac{x}{\Delta x}} \sin\frac{\pi x}{2\Delta x} + [\sum_{p=0}^{N} B_p[x]_{\Delta x}^{p}][w\Delta x]^{\frac{x}{\Delta x}} \cos\frac{\pi x}{2\Delta x}$, the particular solution to this equation is:

$$f_P(x) = [\sum_{p=N_h}^{N+N_h} \mathcal{A}_p[x]_{\Delta x}^{p}][w\Delta x]^{\frac{x}{\Delta x}} \sin\frac{\pi x}{2\Delta x} + [\sum_{p=N_h}^{N+N_h} \mathcal{B}_p[x]_{\Delta x}^{p}][w\Delta x]^{\frac{x}{\Delta x}} \cos\frac{\pi x}{2\Delta x}$$

where

N = the polynomial order of the summations in Q(x)

$N_h$ = the number of times the Q(x) complex conjugate related roots, $-\frac{1}{\Delta x} + jw$ and $-\frac{1}{\Delta x} - jw$ appear in the characteristic polynomial, h(r)

$A_p, B_p, \mathcal{A}_p, \mathcal{B}_p, w$ = real value constants

n = order of the differential difference equation

$f_P(x)$ = particular solution to a differential difference equation

$\Delta x$ = x increment

If $Q(x) = Q_1(x) + Q_2(x) + Q_3(x) + \ldots + Q_g(x)$, then $f_P(x)$ is the sum of the particular solutions of each $Q_r(x)$ where r = 1,2,3,…,g.

The undetermined coefficients of the $f_P(x)$ function are evaluated by introducing the obtained $f_P(x)$ function into the differential difference equation and solving for the coefficient values which will yield an identity.

# TABLE 13a

## Differential Difference Equation Q(x) Functions and their Corresponding Undetermined Coefficient Particular Solution Functions

$$D_{\Delta x}^{n}f(x) + a_{n-1}D_{\Delta x}^{n-1}f(x) + a_{n-2}D_{\Delta x}^{n-2}f(x) + \ldots + a_1 D_{\Delta x}f(x) + a_0 f(x) = Q(x)$$

| # | For these Interval Calculus Functions in Q(x)<br><br>The function specified should be replaced in Q(x) by its identity if its identity is listed.<br>$A_p, B_p, k_p$, A, a, b, k, u, w, Δx are constants | Put these Interval Calculus functions with undetermined coefficients in the particular solution, $f_P(x)$.<br><br>$N_h$ = the number of times the Q(x) related root(s) appear in the characteristic polynomial, h(r) | Related Root(s) |
|---|---|---|---|
| **1** | $\sum_{p=0}^{N} k_p[x]_{\Delta x}^{p}$ | $\sum_{p=N_h}^{N+N_h} A_p[x]_{\Delta x}^{p}$ | r = 0 |
| **1a** | k | $\sum_{p=N_h}^{N_h} A_p[x]_{\Delta x}^{p}$ | r = 0 |
| **1b** | kx<br>$= k[x]_{\Delta x}^{1}$ | $\sum_{p=N_h}^{1+N_h} A_p[x]_{\Delta x}^{p}$ | r = 0 |
| **1c** | $kx^2$<br>$= k\Delta x[x]_{\Delta x}^{1} + k[x]_{\Delta x}^{2}$ | $\sum_{p=N_h}^{2+N_h} A_p[x]_{\Delta x}^{p}$ | r = 0 |
| **1d** | $kx^3$<br>$= k(\Delta x^2)\,[x]_{\Delta x}^{1} + k(3\Delta x)\,[x]_{\Delta x}^{2} +$<br>$k[x]_{\Delta x}^{3}$ | $\sum_{p=N_h}^{3+N_h} A_p[x]_{\Delta x}^{p}$ | r = 0 |
| **1e** | $k[x]_{\Delta x}^{n}$<br>n = 1,2,3,… | $\sum_{p=N_h}^{n+N_h} A_p[x]_{\Delta x}^{p}$ | r = 0 |

| # | For these Interval Calculus Functions in Q(x)<br><br>The function specified should be replaced in Q(x) by its identity if its identity is listed.<br><br>$A_p, B_p, k_p$, A, a, b, k, u, w, Δx are constants | Put these Interval Calculus functions with undetermined coefficients in the particular solution, $f_P(x)$.<br><br>$N_h$ = the number of times the Q(x) related root(s) appear in the characteristic polynomial, h(r) | Related Root(s) |
|---|---|---|---|
| **2** | $[\sum_{p=0}^{N} k_p[x]_{\Delta x}^{p}]\ e_{\Delta x}(u,x)$ | $[\sum_{p=N_h}^{N+N_h} A_p[x]_{\Delta x}^{p}]\ e_{\Delta x}(u,x)$ | $r = u$<br><br>$u \neq -\frac{1}{\Delta x}$ |
| **2a** | $k e_{\Delta x}(u,x)$ | $[\sum_{p=N_h}^{N_h} A_p[x]_{\Delta x}^{p}]\ e_{\Delta x}(u,x)$ | $r = u$<br><br>$u \neq -\frac{1}{\Delta x}$ |
| **2b** | $kx e_{\Delta x}(u,x)$ | $[\sum_{p=N_h}^{1+N_h} A_p[x]_{\Delta x}^{p}]\ e_{\Delta x}(u,x)$ | $r = u$<br><br>$u \neq -\frac{1}{\Delta x}$ |
| **2c** | $k[x]_{\Delta x}^{n} e_{\Delta x}(u,x)$ | $[\sum_{p=N_h}^{n+N_h} A_p[x]_{\Delta x}^{p}]\ e_{\Delta x}(u,x)$ | $r = u$<br><br>$u \neq -\frac{1}{\Delta x}$ |
| **3** | $[\sum_{p=0}^{N} k_p[x]_{\Delta x}^{p}]\ \sin_{\Delta x}(w,x)$ | $[\sum_{p=N_h}^{N+N_h} A_p[x]_{\Delta x}^{p}]\ \sin_{\Delta x}(w,x)\ +$<br>$[\sum_{p=N_h}^{N+N_h} B_p[x]_{\Delta x}^{p}]\ \cos_{\Delta x}(w,x)$ | $r = \pm jw$ |
| **3a** | $k\sin_{\Delta x}(w,x)$ | $[\sum_{p=N_h}^{N_h} A_p[x]_{\Delta x}^{p}]\ \sin_{\Delta x}(w,x)\ +$<br>$[\sum_{p=N_h}^{N_h} B_p[x]_{\Delta x}^{p}]\ \cos_{\Delta x}(w,x)$ | $r = \pm jw$ |

| # | For these Interval Calculus Functions in Q(x)<br><br>The function specified should be replaced in Q(x) by its identity if its identity is listed.<br><br>$A_p, B_p, k_p$, A, a, b, k, u, w, Δx are constants | Put these Interval Calculus functions with undetermined coefficients in the particular solution, $f_P(x)$.<br><br>$N_h$ = the number of times the Q(x) related root(s) appear in the characteristic polynomial, h(r) | Related Root(s) |
|---|---|---|---|
| **3b** | $kx\sin_{\Delta x}(w,x)$ | $[\sum_{p=N_h}^{1+N_h} A_p[x]_{\Delta x}^{p}]\sin_{\Delta x}(w,x) +$<br>$[\sum_{p=N_h}^{1+N_h} B_p[x]_{\Delta x}^{p}]\cos_{\Delta x}(w,x)$ | $r = \pm jw$ |
| **3c** | $k[x]_{\Delta x}^{n}\sin_{\Delta x}(w,x)$ | $[\sum_{p=N_h}^{n+N_h} A_p[x]_{\Delta x}^{p}]\sin_{\Delta x}(w,x) +$<br>$[\sum_{p=N_h}^{n+N_h} B_p[x]_{\Delta x}^{p}]\cos_{\Delta x}(w,x)$ | $r = \pm jw$ |
| **4** | $[\sum_{p=0}^{N} k_p[x]_{\Delta x}^{p}]\cos_{\Delta x}(w,x)$ | $[\sum_{p=N_h}^{N+N_h} A_p[x]_{\Delta x}^{p}]\sin_{\Delta x}(w,x) +$<br>$[\sum_{p=N_h}^{N+N_h} B_p[x]_{\Delta x}^{p}]\cos_{\Delta x}(w,x)$ | $r = \pm jw$ |
| **4a** | $k\cos_{\Delta x}(w,x)$ | $[\sum_{p=N_h}^{N_h} A_p[x]_{\Delta x}^{p}]\sin_{\Delta x}(w,x) +$<br>$[\sum_{p=N_h}^{N_h} B_p[x]_{\Delta x}^{p}]\cos_{\Delta x}(w,x)$ | $r = \pm jw$ |

| # | For these Interval Calculus Functions in Q(x)<br><br>The function specified should be replaced in Q(x) by its identity if its identity is listed.<br><br>$A_p, B_p, k_p$, A, a, b, k, u, w, Δx are constants | Put these Interval Calculus functions with undetermined coefficients in the particular solution, $f_P(x)$.<br><br>$N_h$ = the number of times the Q(x) related root(s) appear in the characteristic polynomial, h(r) | Related Root(s) |
|---|---|---|---|
| **4b** | $kx\cos_{\Delta x}(w,x)$ | $[\sum_{p=N_h}^{1+N_h} A_p[x]_{\Delta x}^{p}]\sin_{\Delta x}(w,x) +$<br>$[\sum_{p=N_h}^{1+N_h} B_p[x]_{\Delta x}^{p}]\cos_{\Delta x}(w,x)$ | $r = \pm jw$ |
| **4c** | $k[x]_{\Delta x}^{n}\cos_{\Delta x}(w,x)$ | $[\sum_{p=N_h}^{n+N_h} A_p[x]_{\Delta x}^{p}]\sin_{\Delta x}(w,x) +$<br>$[\sum_{p=N_h}^{n+N_h} B_p[x]_{\Delta x}^{p}]\cos_{\Delta x}(w,x)$ | $r = \pm jw$ |
| **5** | $[\sum_{p=0}^{N} k_p[x]_{\Delta x}^{p}]e_{\Delta x}(u,x)\sin_{\Delta x}(w,x)$ | $[\sum_{p=N_h}^{N+N_h} A_p[x]_{\Delta x}^{p}]e_{\Delta x}(u,x)\sin_{\Delta x}(w,x) +$<br>$[\sum_{p=N_h}^{N+N_h} B_p[x]_{\Delta x}^{p}]e_{\Delta x}(u,x)\cos_{\Delta x}(w,x)$ | $r = u \pm jw(1+u\Delta x)$<br>$u \neq -\frac{1}{\Delta x}$ |
| **5a** | $ke_{\Delta x}(u,x)\sin_{\Delta x}(w,x)$ | $[\sum_{p=N_h}^{N_h} A_p[x]_{\Delta x}^{p}]e_{\Delta x}(u,x)\sin_{\Delta x}(w,x) +$<br>$[\sum_{p=N_h}^{N_h} B_p[x]_{\Delta x}^{p}]e_{\Delta x}(u,x)\cos_{\Delta x}(w,x)$ | $r = u \pm jw(1+u\Delta x)$<br>$u \neq -\frac{1}{\Delta x}$ |

| # | For these Interval Calculus Functions in Q(x)<br><br>The function specified should be replaced in Q(x) by its identity if its identity is listed.<br><br>$A_p, B_p, k_p$, A, a, b, k, u, w, Δx are constants | Put these Interval Calculus functions with undetermined coefficients in the particular solution, $f_P(x)$.<br><br>$N_h$ = the number of times the Q(x) related root(s) appear in the characteristic polynomial, h(r) | Related Root(s) |
|---|---|---|---|
| **5b** | $kxe_{\Delta x}(u,x)\ \sin_{\Delta x}(w,x)$ | $[\sum_{p=N_h}^{1+N_h} A_p[x]_{\Delta x}^{p}]\ e_{\Delta x}(u,x)\ \sin_{\Delta x}(w,x)\ +$<br>$[\sum_{p=N_h}^{1+N_h} B_p[x]_{\Delta x}^{p}]\ e_{\Delta x}(u,x)\ \cos_{\Delta x}(w,x)$ | $r = u \pm jw(1+u\Delta x)$<br>$u \neq -\frac{1}{\Delta x}$ |
| **5c** | $k[x]_{\Delta x}^{n}\ e_{\Delta x}(u,x)\ \sin_{\Delta x}(w,x)$ | $[\sum_{p=N_h}^{n+N_h} A_p[x]_{\Delta x}^{p}]\ e_{\Delta x}(u,x)\ \sin_{\Delta x}(w,x)\ +$<br>$[\sum_{p=N_h}^{n+N_h} B_p[x]_{\Delta x}^{p}]\ e_{\Delta x}(u,x)\ \cos_{\Delta x}(w,x)$ | $r = u \pm jw(1+u\Delta x)$<br>$u \neq -\frac{1}{\Delta x}$ |
| **6** | $[\sum_{p=0}^{N} k_p[x]_{\Delta x}^{p}]\ e_{\Delta x}(u,x)\ \cos_{\Delta x}(w,x)$ | $[\sum_{p=N_h}^{N+N_h} A_p[x]_{\Delta x}^{p}]\ e_{\Delta x}(u,x)\ \sin_{\Delta x}(w,x)\ +$<br>$[\sum_{p=N_h}^{N+N_h} B_p[x]_{\Delta x}^{p}]\ e_{\Delta x}(u,x)\ \cos_{\Delta x}(w,x)$ | $r = u \pm jw(1+u\Delta x)$<br>$u \neq -\frac{1}{\Delta x}$ |
| **6a** | $ke_{\Delta x}(u,x)\ \cos_{\Delta x}(w,x)$ | $[\sum_{p=N_h}^{N_h} A_p[x]_{\Delta x}^{p}]\ e_{\Delta x}(u,x)\ \sin_{\Delta x}(w,x)\ +$<br>$[\sum_{p=N_h}^{N_h} B_p[x]_{\Delta x}^{p}]\ e_{\Delta x}(u,x)\ \cos_{\Delta x}(w,x)$ | $r = u \pm jw(1+u\Delta x)$<br>$u \neq -\frac{1}{\Delta x}$ |

| # | For these Interval Calculus Functions in Q(x)<br><br>The function specified should be replaced in Q(x) by its identity if its identity is listed.<br><br>$A_p, B_p, k_p$, A, a, b, k, u, w, Δx are constants | Put these Interval Calculus functions with undetermined coefficients in the particular solution, $f_P(x)$.<br><br>$N_h$ = the number of times the Q(x) related root(s) appear in the characteristic polynomial, h(r) | Related Root(s) |
|---|---|---|---|
| **6b** | $kxe_{\Delta x}(u,x)\cos_{\Delta x}(w,x)$ | $[\sum_{p=N_h}^{1+N_h} A_p[x]_{\Delta x}^{p}]\, e_{\Delta x}(u,x)\sin_{\Delta x}(w,x)\ +$<br>$[\sum_{p=N_h}^{1+N_h} B_p[x]_{\Delta x}^{p}]\, e_{\Delta x}(u,x)\cos_{\Delta x}(w,x)$ | $r = u \pm jw(1+u\Delta x)$<br>$u \neq -\frac{1}{\Delta x}$ |
| **6c** | $k[x]_{\Delta x}^{n}\, e_{\Delta x}(u,x)\cos_{\Delta x}(w,x)$ | $[\sum_{p=N_h}^{n+N_h} A_p[x]_{\Delta x}^{p}]\, e_{\Delta x}(u,x)\sin_{\Delta x}(w,x)\ +$<br>$[\sum_{p=N_h}^{n+N_h} B_p[x]_{\Delta x}^{p}]\, e_{\Delta x}(u,x)\cos_{\Delta x}(w,x)$ | $r = u \pm jw(1+u\Delta x)$<br>$u \neq -\frac{1}{\Delta x}$ |
| **7** | $[\sum_{p=0}^{N} k_p[x]_{\Delta x}^{p}][w\Delta x]^{\frac{x}{\Delta x}}\sin\frac{\pi x}{2\Delta x}$ | $[\sum_{p=N_h}^{N+N_h} A_p[x]_{\Delta x}^{p}][w\Delta x]^{\frac{x}{\Delta x}}\sin\frac{\pi x}{2\Delta x}\ +$<br>$[\sum_{p=N_h}^{N+N_h} B_p[x]_{\Delta x}^{p}][w\Delta x]^{\frac{x}{\Delta x}}\cos\frac{\pi x}{2\Delta x}$ | $r = -\frac{1}{\Delta x} \pm jw$ |
| **7a** | $k[w\Delta x]^{\frac{x}{\Delta x}}\sin\frac{\pi x}{2\Delta x}$ | $[\sum_{p=N_h}^{N_h} A_p[x]_{\Delta x}^{p}][w\Delta x]^{\frac{x}{\Delta x}}\sin\frac{\pi x}{2\Delta x}\ +$<br>$[\sum_{p=N_h}^{N_h} B_p[x]_{\Delta x}^{p}][w\Delta x]^{\frac{x}{\Delta x}}\cos\frac{\pi x}{2\Delta x}$ | $r = -\frac{1}{\Delta x} \pm jw$ |

| # | **For these Interval Calculus Functions in Q(x)**<br><br>**The function specified should be replaced in Q(x) by its identity if its identity is listed.**<br>**$A_p, B_p, k_p$, A, a, b, k, u, w, $\Delta x$ are constants** | **Put these Interval Calculus functions with undetermined coefficients in the particular solution, $f_P(x)$.**<br><br>**$N_h$ = the number of times the Q(x) related root(s) appear in the characteristic polynomial, h(r)** | **Related Root(s)** |
|---|---|---|---|
| **7b** | $kx[w\Delta x]^{\frac{x}{\Delta x}}\sin\frac{\pi x}{2\Delta x}$ | $[\sum_{p=N_h}^{1+N_h} A_p[x]_{\Delta x}^{p}][w\Delta x]^{\frac{x}{\Delta x}}\sin\frac{\pi x}{2\Delta x}\ +$<br>$[\sum_{p=N_h}^{1+N_h} B_p[x]_{\Delta x}^{p}][w\Delta x]^{\frac{x}{\Delta x}}\cos\frac{\pi x}{2\Delta x}$ | $r = -\frac{1}{\Delta x} \pm jw$ |
| **7c** | $k[x]_{\Delta x}^{n}[w\Delta x]^{\frac{x}{\Delta x}}\sin\frac{\pi x}{2\Delta x}$ | $[\sum_{p=N_h}^{n+N_h} A_p[x]_{\Delta x}^{p}][w\Delta x]^{\frac{x}{\Delta x}}\sin\frac{\pi x}{2\Delta x}\ +$<br>$[\sum_{p=N_h}^{n+N_h} B_p[x]_{\Delta x}^{p}][w\Delta x]^{\frac{x}{\Delta x}}\cos\frac{\pi x}{2\Delta x}$ | $r = -\frac{1}{\Delta x} \pm jw$ |
| **8** | $[\sum_{p=0}^{N} k_p[x]_{\Delta x}^{p}][w\Delta x]^{\frac{x}{\Delta x}}\cos\frac{\pi x}{2\Delta x}$ | $[\sum_{p=N_h}^{N+N_h} A_p[x]_{\Delta x}^{p}][w\Delta x]^{\frac{x}{\Delta x}}\sin\frac{\pi x}{2\Delta x}\ +$<br>$[\sum_{p=N_h}^{N+N_h} B_p[x]_{\Delta x}^{p}][w\Delta x]^{\frac{x}{\Delta x}}\cos\frac{\pi x}{2\Delta x}$ | $r = -\frac{1}{\Delta x} \pm jw$ |
| **8a** | $k[w\Delta x]^{\frac{x}{\Delta x}}\cos\frac{\pi x}{2\Delta x}$ | $[\sum_{p=N_h}^{N_h} A_p[x]_{\Delta x}^{p}][w\Delta x]^{\frac{x}{\Delta x}}\sin\frac{\pi x}{2\Delta x}\ +$<br>$[\sum_{p=N_h}^{N_h} B_p[x]_{\Delta x}^{p}][w\Delta x]^{\frac{x}{\Delta x}}\cos\frac{\pi x}{2\Delta x}$ | $r = -\frac{1}{\Delta x} \pm jw$ |

| # | For these Interval Calculus Functions in Q(x)<br><br>The function specified should be replaced in Q(x) by its identity if its identity is listed.<br><br>$A_p$, $B_p$, $k_p$, A, a, b, k, u, w, Δx are constants | Put these Interval Calculus functions with undetermined coefficients in the particular solution, $f_P(x)$.<br><br>$N_h$ = the number of times the Q(x) related root(s) appear in the characteristic polynomial, h(r) | Related Root(s) |
|---|---|---|---|
| **8b** | $kx[w\Delta x]^{\frac{x}{\Delta x}}\cos\frac{\pi x}{2\Delta x}$ | $[\sum_{p=N_h}^{1+N_h} A_p[x]_{\Delta x}^{p}][w\Delta x]^{\frac{x}{\Delta x}}\sin\frac{\pi x}{2\Delta x}$ +<br>$[\sum_{p=N_h}^{1+N_h} B_p[x]_{\Delta x}^{p}][w\Delta x]^{\frac{x}{\Delta x}}\cos\frac{\pi x}{2\Delta x}$ | $r = -\frac{1}{\Delta x} \pm jw$ |
| **8c** | $k[x]_{\Delta x}^{n}[w\Delta x]^{\frac{x}{\Delta x}}\cos\frac{\pi x}{2\Delta x}$ | $[\sum_{p=N_h}^{n+N_h} A_p[x]_{\Delta x}^{p}][w\Delta x]^{\frac{x}{\Delta x}}\sin\frac{\pi x}{2\Delta x}$ +<br>$[\sum_{p=N_h}^{n+N_h} B_p[x]_{\Delta x}^{p}][w\Delta x]^{\frac{x}{\Delta x}}\cos\frac{\pi x}{2\Delta x}$ | $r = -\frac{1}{\Delta x} \pm jw$ |
| **9** | $[\sum_{p=0}^{N} k_p[x]_{\Delta x}^{p}]\sinh_{\Delta x}(w,x)$ | $[\sum_{p=N_h}^{N+N_h} A_p[x]_{\Delta x}^{p}]\sinh_{\Delta x}(w,x)$ +<br>$[\sum_{p=N_h}^{N+N_h} B_p[x]_{\Delta x}^{p}]\cosh_{\Delta x}(w,x)$ | $r = \pm w$ |
| **9a** | $k\sinh_{\Delta x}(w,x)$ | $[\sum_{p=N_h}^{N_h} A_p[x]_{\Delta x}^{p}]\sinh_{\Delta x}(w,x)$ +<br>$[\sum_{p=N_h}^{N_h} B_p[x]_{\Delta x}^{p}]\cosh_{\Delta x}(w,x)$ | $r = \pm w$ |

| # | For these Interval Calculus Functions in Q(x)<br><br>The function specified should be replaced in Q(x) by its identity if its identity is listed.<br><br>$A_p, B_p, k_p$, A, a, b, k, u, w, Δx are constants | Put these Interval Calculus functions with undetermined coefficients in the particular solution, $f_P(x)$.<br><br>$N_h$ = the number of times the Q(x) related root(s) appear in the characteristic polynomial, h(r) | Related Root(s) |
|---|---|---|---|
| **9b** | $kx\sinh_{\Delta x}(w,x)$ | $[\sum_{p=N_h}^{1+N_h} A_p[x]_{\Delta x}^{p}]\sinh_{\Delta x}(w,x) +$<br>$[\sum_{p=N_h}^{1+N_h} B_p[x]_{\Delta x}^{p}]\cosh_{\Delta x}(w,x)$ | $r = \pm w$ |
| **9c** | $k[x]_{\Delta x}^{n}\sinh_{\Delta x}(w,x)$ | $[\sum_{p=N_h}^{n+N_h} A_p[x]_{\Delta x}^{p}]\sinh_{\Delta x}(w,x) +$<br>$[\sum_{p=N_h}^{n+N_h} B_p[x]_{\Delta x}^{p}]\cosh_{\Delta x}(w,x)$ | $r = \pm w$ |
| **10** | $[\sum_{p=0}^{N} k_p[x]_{\Delta x}^{p}]\cosh_{\Delta x}(w,x)$ | $[\sum_{p=N_h}^{N+N_h} A_p[x]_{\Delta x}^{p}]\sinh_{\Delta x}(w,x) +$<br>$[\sum_{p=N_h}^{N+N_h} B_p[x]_{\Delta x}^{p}]\cosh_{\Delta x}(w,x)$ | $r = \pm w$ |
| **10a** | $k\cosh_{\Delta x}(w,x)$ | $[\sum_{p=N_h}^{N_h} A_p[x]_{\Delta x}^{p}]\sinh_{\Delta x}(w,x) +$<br>$[\sum_{p=N_h}^{N_h} B_p[x]_{\Delta x}^{p}]\cosh_{\Delta x}(w,x)$ | $r = \pm w$ |

| # | For these Interval Calculus Functions in Q(x)<br><br>The function specified should be replaced in Q(x) by its identity if its identity is listed.<br>$A_p, B_p, k_p$, A, a, b, k, u, w, Δx are constants | Put these Interval Calculus functions with undetermined coefficients in the particular solution, $f_P(x)$.<br><br>$N_h$ = the number of times the Q(x) related root(s) appear in the characteristic polynomial, h(r) | Related Root(s) |
|---|---|---|---|
| **10b** | $kx\cosh_{\Delta x}(w,x)$ | $[\sum_{p=N_h}^{1+N_h} A_p[x]_{\Delta x}^{p}]\sinh_{\Delta x}(w,x)\ +$<br>$[\sum_{p=N_h}^{n+N_h} B_p[x]_{\Delta x}^{p}]\cosh_{\Delta x}(w,x)$ | $r = \pm w$ |
| **10c** | $k[x]_{\Delta x}^{n}\cosh_{\Delta x}(w,x)$ | $[\sum_{p=N_h}^{n+N_h} A_p[x]_{\Delta x}^{p}]\sinh_{\Delta x}(w,x)\ +$<br>$[\sum_{p=N_h}^{n+N_h} B_p[x]_{\Delta x}^{p}]\cosh_{\Delta x}(w,x)$ | $r = \pm w$ |
| **11** | $[\sum_{p=0}^{N} k_p[x]_{\Delta x}^{p}]\,A^x$<br>$= [\sum_{p=0}^{N} k_p[x]_{\Delta x}^{p}]\,e_{\Delta x}(\frac{A^{\Delta x}-1}{\Delta x}, x)$ | $[\sum_{p=N_h}^{N+N_h} A_p[x]_{\Delta x}^{p}]\,e_{\Delta x}(\frac{A^{\Delta x}-1}{\Delta x}, x)$ | $r = \frac{A^{\Delta x}-1}{\Delta x}$ |
| **11a** | $kA^x$<br>$= ke_{\Delta x}(\frac{A^{\Delta x}-1}{\Delta x}, x)$ | $[\sum_{p=N_h}^{N_h} A_p[x]_{\Delta x}^{p}]\,e_{\Delta x}(\frac{A^{\Delta x}-1}{\Delta x}, x)$ | $r = \frac{A^{\Delta x}-1}{\Delta x}$ |
| **11b** | $kxA^x$<br>$= kxe_{\Delta x}(\frac{A^{\Delta x}-1}{\Delta x}, x)$ | $[\sum_{p=N_h}^{1+N_h} A_p[x]_{\Delta x}^{p}]\,e_{\Delta x}(\frac{A^{\Delta x}-1}{\Delta x}, x)$ | $r = \frac{A^{\Delta x}-1}{\Delta x}$ |

| # | For these Interval Calculus Functions in Q(x)<br><br>The function specified should be replaced in Q(x) by its identity if its identity is listed.<br><br>$A_p$, $B_p$, $k_p$, A, a, b, k, u, w, Δx are constants | Put these Interval Calculus functions with undetermined coefficients in the particular solution, $f_P(x)$.<br><br>$N_h$ = the number of times the Q(x) related root(s) appear in the characteristic polynomial, h(r) | Related Root(s) |
|---|---|---|---|
| **12** | $[\sum_{p=0}^{N} k_p[x]_{\Delta x}^{p}]\, e^{ax}$<br><br>$= [\sum_{p=0}^{N} k_p[x]_{\Delta x}^{p}]\, e_{\Delta x}(\frac{e^{a\Delta x}-1}{\Delta x}, x)$ | $[\sum_{p=N_h}^{N+N_h} A_p[x]_{\Delta x}^{p}]\, e_{\Delta x}(\frac{e^{a\Delta x}-1}{\Delta x}, x)$ | $r = \frac{e^{a\Delta x}-1}{\Delta x}$ |
| **12a** | $ke^{ax}$<br><br>$= ke_{\Delta x}(\frac{e^{a\Delta x}-1}{\Delta x}, x)$ | $[\sum_{p=N_h}^{N_h} A_p[x]_{\Delta x}^{p}]\, e_{\Delta x}(\frac{e^{a\Delta x}-1}{\Delta x}, x)$ | $r = \frac{e^{a\Delta x}-1}{\Delta x}$ |
| **12b** | $kxe^{ax}$<br><br>$= kxe_{\Delta x}(\frac{e^{a\Delta x}-1}{\Delta x}, x)$ | $[\sum_{p=N_h}^{1+N_h} A_p[x]_{\Delta x}^{p}]\, e_{\Delta x}(\frac{e^{a\Delta x}-1}{\Delta x}, x)$ | $r = \frac{e^{a\Delta x}-1}{\Delta x}$ |

| # | For these The function specified should be replaced in Q(x) by its identity if its identity is listed.<br>k = constant | Put these Interval Calculus functions with undetermined coefficients in the particular solution, $f_P(x)$.<br>$N_h$ = Interval Calculus Functions in Q(x)<br>the number of times the Q(x) related root(s) appear in the characteristic polynomial, h(r) | Related Root(s) |
|---|---|---|---|
| **13** | $[\sum_{p=0}^{N} k_p[x]_{\Delta x}^{p}]\ \text{sinbx}$<br>$= [\sum_{p=0}^{N} k_p[x]_{\Delta x}^{p}]\ e_{\Delta x}(\frac{\cos b\Delta x - 1}{\Delta x}, x)\ \sin_{\Delta x}(\frac{\tan b\Delta x}{\Delta x}, x)$<br>where<br>$\frac{x}{\Delta x}$ = integer<br>$\cos b\Delta x \neq 0$ | $[\sum_{p=N_h}^{N+N_h} A_p[x]_{\Delta x}^{p}]\ e_{\Delta x}(\frac{\cos b\Delta x-1}{\Delta x}, x)\ \sin_{\Delta x}(\frac{\tan b\Delta x}{\Delta x}, x) +$<br>$[\sum_{p=N_h}^{N+N_h} B_p[x]_{\Delta x}^{p}]\ e_{\Delta x}(\frac{\cos b\Delta x-1}{\Delta x}, x)\ \cos_{\Delta x}(\frac{\tan b\Delta x}{\Delta x}, x)$ | $r = \frac{\cos b\Delta x - 1}{\Delta x} \pm j\ \frac{\sin b\Delta x}{\Delta x}$<br>$\cos b\Delta x \neq 0$ |
| **13a** | ksinbx<br>$= k e_{\Delta x}(\frac{\cos b\Delta x - 1}{\Delta x}, x)\ \sin_{\Delta x}(\frac{\tan b\Delta x}{\Delta x}, x)$<br>where<br>$\frac{x}{\Delta x}$ = integer<br>$\cos b\Delta x \neq 0$ | $[\sum_{p=N_h}^{N_h} A_p[x]_{\Delta x}^{p}]\ e_{\Delta x}(\frac{\cos b\Delta x-1}{\Delta x}, x)\ \sin_{\Delta x}(\frac{\tan b\Delta x}{\Delta x}, x) +$<br>$[\sum_{p=N_h}^{N_h} B_p[x]_{\Delta x}^{p}]\ e_{\Delta x}(\frac{\cos b\Delta x-1}{\Delta x}, x)\ \cos_{\Delta x}(\frac{\tan b\Delta x}{\Delta x}, x)$ | $r = \frac{\cos b\Delta x - 1}{\Delta x} \pm j\ \frac{\sin b\Delta x}{\Delta x}$<br>$\cos b\Delta x \neq 0$ |

| # | For these The function specified should be replaced in Q(x) by its identity if its identity is listed.<br><br>k = constant | Put these Interval Calculus functions with undetermined coefficients in the particular solution, $f_P(x)$.<br><br>$N_h$ = Interval Calculus Functions in Q(x)<br><br>the number of times the Q(x) related root(s) appear in the characteristic polynomial, h(r) | Related Root(s) |
|---|---|---|---|
| **13b** | kxsinbx<br><br>$= kxe_{\Delta x}(\frac{\cos b\Delta x - 1}{\Delta x}, x)\sin_{\Delta x}(\frac{\tan b\Delta x}{\Delta x}, x)$<br><br>where<br><br>$\frac{x}{\Delta x}$ = integer<br><br>$\cos b\Delta x \neq 0$ | $[\sum_{p=N_h}^{1+N_h} A_p[x]_{\Delta x}^{p}]\, e_{\Delta x}(\frac{\cos b\Delta x-1}{\Delta x}, x)\sin_{\Delta x}(\frac{\tan b\Delta x}{\Delta x}, x) +$<br><br>$[\sum_{p=N_h}^{1+N_h} B_p[x]_{\Delta x}^{p}]\, e_{\Delta x}(\frac{\cos b\Delta x-1}{\Delta x}, x)\cos_{\Delta x}(\frac{\tan b\Delta x}{\Delta x}, x)$ | $r = \frac{\cos b\Delta x - 1}{\Delta x} \pm j\,\frac{\sin b\Delta x}{\Delta x}$<br><br>$\cos b\Delta x \neq 0$ |
| **14** | $[\sum_{p=0}^{N} k_p[x]_{\Delta x}^{p}]\cos bx$<br><br>$= [\sum_{p=0}^{N} k_p[x]_{\Delta x}^{p}]\, e_{\Delta x}(\frac{\cos b\Delta x - 1}{\Delta x}, x)\cos_{\Delta x}(\frac{\tan b\Delta x}{\Delta x}, x)$<br><br>where<br><br>$\frac{x}{\Delta x}$ = integer<br><br>$\cos b\Delta x \neq 0$ | $[\sum_{p=N_h}^{N+N_h} A_p[x]_{\Delta x}^{p}]\, e_{\Delta x}(\frac{\cos b\Delta x-1}{\Delta x}, x)\sin_{\Delta x}(\frac{\tan b\Delta x}{\Delta x}, x) +$<br><br>$[\sum_{p=N_h}^{N+N_h} B_p[x]_{\Delta x}^{p}]\, e_{\Delta x}(\frac{\cos b\Delta x-1}{\Delta x}, x)\cos_{\Delta x}(\frac{\tan b\Delta x}{\Delta x}, x)$ | $r = \frac{\cos b\Delta x - 1}{\Delta x} \pm j\,\frac{\sin b\Delta x}{\Delta x}$<br><br>$\cos b\Delta x \neq 0$ |

| # | For these The function specified should be replaced in Q(x) by its identity if its identity is listed.<br><br>k = constant | Put these Interval Calculus functions with undetermined coefficients in the particular solution, $f_P(x)$.<br><br>$N_h$ = Interval Calculus Functions in Q(x)<br><br>the number of times the Q(x) related root(s) appear in the characteristic polynomial, h(r) | Related Root(s) |
|---|---|---|---|
| **14a** | kcosbx<br>$= k e_{\Delta x}( \frac{\cos b\Delta x -1}{\Delta x}, x) \cos_{\Delta x}( \frac{\tan b\Delta x}{\Delta x}, x )$<br>where<br>$\frac{x}{\Delta x}$ = integer<br>$\cos b\Delta x \neq 0$ | $[ \sum_{p=N_h}^{N_h} A_p[x]_{\Delta x}^{p} ]\ e_{\Delta x}(\frac{\cos b\Delta x-1}{\Delta x}, x) \sin_{\Delta x}( \frac{\tan b\Delta x}{\Delta x}, x ) +$<br>$[ \sum_{p=N_h}^{N_h} B_p[x]_{\Delta x}^{p} ]\ e_{\Delta x}(\frac{\cos b\Delta x-1}{\Delta x}, x) \cos_{\Delta x}( \frac{\tan b\Delta x}{\Delta x}, x )$ | $r = \frac{\cos b\Delta x -1}{\Delta x} \pm j \frac{\sin b\Delta x}{\Delta x}$<br>$\cos b\Delta x \neq 0$ |
| **14b** | kxcosbx<br>$= kx e_{\Delta x}( \frac{\cos b\Delta x -1}{\Delta x}, x) \cos_{\Delta x}( \frac{\tan b\Delta x}{\Delta x}, x )$<br>where<br>$\frac{x}{\Delta x}$ = integer<br>$\cos b\Delta x \neq 0$ | $[ \sum_{p=N_h}^{1+N_h} A_p[x]_{\Delta x}^{p} ]\ e_{\Delta x}(\frac{\cos b\Delta x-1}{\Delta x}, x) \sin_{\Delta x}( \frac{\tan b\Delta x}{\Delta x}, x ) +$<br>$[ \sum_{p=N_h}^{1+N_h} B_p[x]_{\Delta x}^{p} ]\ e_{\Delta x}(\frac{\cos b\Delta x-1}{\Delta x}, x) \cos_{\Delta x}( \frac{\tan b\Delta x}{\Delta x}, x )$ | $r = \frac{\cos b\Delta x -1}{\Delta x} \pm j \frac{\sin b\Delta x}{\Delta x}$<br>$\cos b\Delta x \neq 0$ |

| # | For these The function specified should be replaced in Q(x) by its identity if its identity is listed.<br><br>k = constant | Put these Interval Calculus functions with undetermined coefficients in the particular solution, $f_P(x)$.<br><br>$N_h$ = Interval Calculus Functions in Q(x)<br><br>the number of times the Q(x) related root(s) appear in the characteristic polynomial, h(r) | Related Root(s) |
|---|---|---|---|
| **15** | $[\sum_{p=0}^{N} k_p[x]_{\Delta x}^{p}]\ e^{ax}\sin bx$<br>=<br>$[\sum_{p=0}^{N} k_p[x]_{\Delta x}^{p}]\ e_{\Delta x}(\frac{e^{a\Delta x}\cos b\Delta x - 1}{\Delta x}, x)\ \sin_{\Delta x}(\frac{\tan b\Delta x}{\Delta x}, x)$<br>where<br>$\frac{x}{\Delta x}$ = integer<br>$\cos b\Delta x \neq 0$ | $[\sum_{p=N_h}^{N+N_h} A_p[x]_{\Delta x}^{p}]\ e_{\Delta x}(\frac{e^{a\Delta x}\cos b\Delta x-1}{\Delta x}, x)\ \sin_{\Delta x}(\frac{\tan b\Delta x}{\Delta x}, x) +$<br><br>$[\sum_{p=N_h}^{N+N_h} B_p[x]_{\Delta x}^{p}]\ e_{\Delta x}(\frac{e^{a\Delta x}\cos b\Delta x-1}{\Delta x}, x)\ \cos_{\Delta x}(\frac{\tan b\Delta x}{\Delta x}, x)$ | $r = \frac{e^{a\Delta x}\cos b\Delta x - 1}{\Delta x} \pm j\,\frac{e^{a\Delta x}\sin b\Delta x}{\Delta x}$<br>$\cos b\Delta x \neq 0$ |
| **15a** | $ke^{ax}\sin bx$<br>$= ke_{\Delta x}(\frac{e^{a\Delta x}\cos b\Delta x-1}{\Delta x}, x)\ \sin_{\Delta x}(\frac{\tan b\Delta x}{\Delta x}, x)$<br>where<br>$\frac{x}{\Delta x}$ = integer<br>$\cos b\Delta x \neq 0$ | $[\sum_{p=N_h}^{N_h} A_p[x]_{\Delta x}^{p}]\ e_{\Delta x}(\frac{e^{a\Delta x}\cos b\Delta x-1}{\Delta x}, x)\ \sin_{\Delta x}(\frac{\tan b\Delta x}{\Delta x}, x) +$<br><br>$[\sum_{p=N_h}^{N_h} B_p[x]_{\Delta x}^{p}]\ e_{\Delta x}(\frac{e^{a\Delta x}\cos b\Delta x-1}{\Delta x}, x)\ \cos_{\Delta x}(\frac{\tan b\Delta x}{\Delta x}, x)$ | $r = \frac{e^{a\Delta x}\cos b\Delta x - 1}{\Delta x} \pm j\,\frac{e^{a\Delta x}\sin b\Delta x}{\Delta x}$<br>$\cos b\Delta x \neq 0$ |

| # | For these The function specified should be replaced in Q(x) by its identity if its identity is listed.<br><br>k = constant | Put these Interval Calculus functions with undetermined coefficients in the particular solution, $f_P(x)$.<br><br>$N_h$ = Interval Calculus Functions in Q(x)<br><br>the number of times the Q(x) related root(s) appear in the characteristic polynomial, h(r) | Related Root(s) |
|---|---|---|---|
| **15b** | $kxe^{ax}\sin bx$<br><br>$= kxe_{\Delta x}(\frac{e^{a\Delta x}\cos b\Delta x-1}{\Delta x}, x) \sin_{\Delta x}(\frac{\tan b\Delta x}{\Delta x}, x)$<br><br>where<br><br>$\frac{x}{\Delta x}$ = integer<br><br>$\cos b\Delta x \neq 0$ | $[\sum_{p=N_h}^{1+N_h} A_p[x]_{\Delta x}^{p}]\, e_{\Delta x}(\frac{e^{a\Delta x}\cos b\Delta x-1}{\Delta x}, x) \sin_{\Delta x}(\frac{\tan b\Delta x}{\Delta x}, x) +$<br><br>$[\sum_{p=N_h}^{1+N_h} B_p[x]_{\Delta x}^{p}]\, e_{\Delta x}(\frac{e^{a\Delta x}\cos b\Delta x-1}{\Delta x}, x) \cos_{\Delta x}(\frac{\tan b\Delta x}{\Delta x}, x)$ | $r = \frac{e^{a\Delta x}\cos b\Delta x -1}{\Delta x} \pm j\, \frac{e^{a\Delta x}\sin b\Delta x}{\Delta x}$<br><br>$\cos b\Delta x \neq 0$ |
| **16** | $[\sum_{p=0}^{N} k_p[x]_{\Delta x}^{p}]\, e^{ax}\cos bx$<br><br>=<br><br>$[\sum_{p=0}^{N} k_p[x]_{\Delta x}^{p}]\, e_{\Delta x}(\frac{e^{a\Delta x}\cos b\Delta x -1}{\Delta x}, x) \cos_{\Delta x}(\frac{\tan b\Delta x}{\Delta x}, x)$<br><br>where<br><br>$\frac{x}{\Delta x}$ = integer<br><br>$\cos b\Delta x \neq 0$ | $[\sum_{p=N_h}^{N+N_h} A_p[x]_{\Delta x}^{p}]\, e_{\Delta x}(\frac{e^{a\Delta x}\cos b\Delta x-1}{\Delta x}, x) \sin_{\Delta x}(\frac{\tan b\Delta x}{\Delta x}, x) +$<br><br>$[\sum_{p=N_h}^{N+N_h} B_p[x]_{\Delta x}^{p}]\, e_{\Delta x}(\frac{e^{a\Delta x}\cos b\Delta x-1}{\Delta x}, x) \cos_{\Delta x}(\frac{\tan b\Delta x}{\Delta x}, x)$ | $r = \frac{e^{a\Delta x}\cos b\Delta x -1}{\Delta x} \pm j\, \frac{e^{a\Delta x}\sin b\Delta x}{\Delta x}$<br><br>$\cos b\Delta x \neq 0$ |

| # | For these The function specified should be replaced in Q(x) by its identity if its identity is listed.<br><br>k = constant | Put these Interval Calculus functions with undetermined coefficients in the particular solution, $f_P(x)$.<br><br>$N_h$ = Interval Calculus Functions in Q(x)<br><br>the number of times the Q(x) related root(s) appear in the characteristic polynomial, h(r) | Related Root(s) |
|---|---|---|---|
| **16a** | $ke^{ax}\cos bx$<br><br>$= ke_{\Delta x}(\frac{e^{a\Delta x}\cos b\Delta x-1}{\Delta x}, x)\ \cos_{\Delta x}(\frac{\tan b\Delta x}{\Delta x}, x)$<br><br>where<br><br>$\frac{x}{\Delta x}$ = integer<br><br>$\cos b\Delta x \neq 0$ | $[\sum_{p=N_h}^{N_h} A_p[x]_{\Delta x}^{p}]\ e_{\Delta x}(\frac{e^{a\Delta x}\cos b\Delta x-1}{\Delta x}, x)\ \sin_{\Delta x}(\frac{\tan b\Delta x}{\Delta x}, x) +$<br><br>$[\sum_{p=N_h}^{N_h} B_p[x]_{\Delta x}^{p}]\ e_{\Delta x}(\frac{e^{a\Delta x}\cos b\Delta x-1}{\Delta x}, x)\ \cos_{\Delta x}(\frac{\tan b\Delta x}{\Delta x}, x)$ | $r = \frac{e^{a\Delta x}\cos b\Delta x - 1}{\Delta x} \pm j\,\frac{e^{a\Delta x}\sin b\Delta x}{\Delta x}$<br><br>$\cos b\Delta x \neq 0$ |
| **16b** | $kxe^{ax}\cos bx$<br><br>$= kxe_{\Delta x}(\frac{e^{a\Delta x}\cos b\Delta x-1}{\Delta x}, x)\ \cos_{\Delta x}(\frac{\tan b\Delta x}{\Delta x}, x)$<br><br>where<br><br>$\frac{x}{\Delta x}$ = integer<br><br>$\cos b\Delta x \neq 0$ | $[\sum_{p=N_h}^{1+N_h} A_p[x]_{\Delta x}^{p}]\ e_{\Delta x}(\frac{e^{a\Delta x}\cos b\Delta x-1}{\Delta x}, x)\ \sin_{\Delta x}(\frac{\tan b\Delta x}{\Delta x}, x) +$<br><br>$[\sum_{p=N_h}^{1+N_h} B_p[x]_{\Delta x}^{p}]\ e_{\Delta x}(\frac{e^{a\Delta x}\cos b\Delta x-1}{\Delta x}, x)\ \cos_{\Delta x}(\frac{\tan b\Delta x}{\Delta x}, x)$ | $r = \frac{e^{a\Delta x}\cos b\Delta x - 1}{\Delta x} \pm j\,\frac{e^{a\Delta x}\sin b\Delta x}{\Delta x}$<br><br>$\cos b\Delta x \neq 0$ |

| # | For these The function specified should be replaced in Q(x) by its identity if its identity is listed.<br><br>k = constant | Put these Interval Calculus functions with undetermined coefficients in the particular solution, $f_P(x)$.<br><br>$N_h$ = Interval Calculus Functions in Q(x)<br><br>the number of times the Q(x) related root(s) appear in the characteristic polynomial, h(r) | Related Root(s) |
|---|---|---|---|
| 17 | $[\sum_{p=0}^{N} k_p [x]_{\Delta x}^{p}]$ sinhbx<br>=<br>$[\sum_{p=0}^{N} k_p [x]_{\Delta x}^{p}]\, e_{\Delta x}(\frac{\cosh b\Delta x - 1}{\Delta x}, x)\, \sinh_{\Delta x}(\frac{\tanh b\Delta x}{\Delta x}, x)$ | $[\sum_{p=N_h}^{N+N_h} A_p [x]_{\Delta x}^{p}]\, e_{\Delta x}(\frac{\cosh b\Delta x - 1}{\Delta x}, x)\, \sinh_{\Delta x}(\frac{\tanh b\Delta x}{\Delta x}, x) +$<br><br>$[\sum_{p=N_h}^{N+N_h} B_p [x]_{\Delta x}^{p}]\, e_{\Delta x}(\frac{\cosh b\Delta x - 1}{\Delta x}, x)\, \cosh_{\Delta x}(\frac{\tanh b\Delta x}{\Delta x}, x)$ | $r = \frac{\cosh b\Delta x - 1}{\Delta x} \pm \frac{\sinh b\Delta x}{\Delta x}$ |
| 17a | ksinhbx<br>$= k e_{\Delta x}(\frac{\cosh b\Delta x - 1}{\Delta x}, x)\, \sinh_{\Delta x}(\frac{\tanh b\Delta x}{\Delta x}, x)$ | $[\sum_{p=N_h}^{N_h} A_p [x]_{\Delta x}^{p}]\, e_{\Delta x}(\frac{\cosh b\Delta x - 1}{\Delta x}, x)\, \sinh_{\Delta x}(\frac{\tanh b\Delta x}{\Delta x}, x) +$<br><br>$[\sum_{p=N_h}^{N_h} B_p [x]_{\Delta x}^{p}]\, e_{\Delta x}(\frac{\cosh b\Delta x - 1}{\Delta x}, x)\, \cosh_{\Delta x}(\frac{\tanh b\Delta x}{\Delta x}, x)$ | $r = \frac{\cosh b\Delta x - 1}{\Delta x} \pm \frac{\sinh b\Delta x}{\Delta x}$ |

| # | For these The function specified should be replaced in Q(x) by its identity if its identity is listed. <br><br> k = constant | Put these Interval Calculus functions with undetermined coefficients in the particular solution, $f_P(x)$. <br><br> $N_h$ = Interval Calculus Functions in Q(x) <br><br> the number of times the Q(x) related root(s) appear in the characteristic polynomial, h(r) | Related Root(s) |
|---|---|---|---|
| **17b** | kxsinhbx <br><br> $= kxe_{\Delta x}(\frac{\cosh b\Delta x - 1}{\Delta x}, x) \sinh_{\Delta x}(\frac{\tanh b\Delta x}{\Delta x}, x)$ | $[\sum_{p=N_h}^{1+N_h} A_p[x]_{\Delta x}^{p}] e_{\Delta x}(\frac{\cosh b\Delta x-1}{\Delta x}, x) \sinh_{\Delta x}(\frac{\tanh b\Delta x}{\Delta x}, x) +$ <br><br> $[\sum_{p=N_h}^{1+N_h} B_p[x]_{\Delta x}^{p}] e_{\Delta x}(\frac{\cosh b\Delta x-1}{\Delta x}, x) \cosh_{\Delta x}(\frac{\tanh b\Delta x}{\Delta x}, x)$ | $r = \frac{\cosh b\Delta x - 1}{\Delta x} \pm \frac{\sinh b\Delta x}{\Delta x}$ |
| **18** | $[\sum_{p=0}^{N} k_p[x]_{\Delta x}^{p}] \cosh bx$ <br><br> = <br><br> $[\sum_{p=0}^{N} k_p[x]_{\Delta x}^{p}] e_{\Delta x}(\frac{\cos bh\Delta x - 1}{\Delta x}, x) \cosh_{\Delta x}(\frac{\tanh b\Delta x}{\Delta x}, x)$ | $[\sum_{p=N_h}^{N+N_h} A_p[x]_{\Delta x}^{p}] e_{\Delta x}(\frac{\cosh b\Delta x-1}{\Delta x}, x) \sinh_{\Delta x}(\frac{\tanh b\Delta x}{\Delta x}, x) +$ <br><br> $[\sum_{p=N_h}^{N+N_h} B_p[x]_{\Delta x}^{p}] e_{\Delta x}(\frac{\cosh b\Delta x-1}{\Delta x}, x) \cosh_{\Delta x}(\frac{\tanh b\Delta x}{\Delta x}, x)$ | $r = \frac{\cosh b\Delta x - 1}{\Delta x} \pm \frac{\sinh b\Delta x}{\Delta x}$ |

| # | For these The function specified should be replaced in Q(x) by its identity if its identity is listed.<br><br>k = constant | Put these Interval Calculus functions with undetermined coefficients in the particular solution, $f_P(x)$.<br><br>$N_h$ = Interval Calculus Functions in Q(x)<br><br>the number of times the Q(x) related root(s) appear in the characteristic polynomial, h(r) | Related Root(s) |
|---|---|---|---|
| **18a** | kcoshbx<br><br>$= ke_{\Delta x}(\frac{\cosh b\Delta x - 1}{\Delta x}, x) \cosh_{\Delta x}(\frac{\tanh b\Delta x}{\Delta x}, x)$ | $[\sum_{p=N_h}^{N_h} A_p[x]_{\Delta x}^p] e_{\Delta x}(\frac{\cosh b\Delta x-1}{\Delta x}, x) \sinh_{\Delta x}(\frac{\tanh b\Delta x}{\Delta x}, x) +$<br><br>$[\sum_{p=N_h}^{N_h} B_p[x]_{\Delta x}^p] e_{\Delta x}(\frac{\cosh b\Delta x-1}{\Delta x}, x) \cosh_{\Delta x}(\frac{\tanh b\Delta x}{\Delta x}, x)$ | $r = \frac{\cosh b\Delta x - 1}{\Delta x} \pm \frac{\sinh b\Delta x}{\Delta x}$ |
| **18b** | kxcoshbx<br><br>$= kxe_{\Delta x}(\frac{\cosh b\Delta x - 1}{\Delta x}, x) \cosh_{\Delta x}(\frac{\tanh b\Delta x}{\Delta x}, x)$ | $[\sum_{p=N_h}^{1+N_h} A_p[x]_{\Delta x}^p] e_{\Delta x}(\frac{\cosh b\Delta x-1}{\Delta x}, x) \sinh_{\Delta x}(\frac{\tanh b\Delta x}{\Delta x}, x) +$<br><br>$[\sum_{p=N_h}^{1+N_h} B_p[x]_{\Delta x}^p] e_{\Delta x}(\frac{\cosh b\Delta x-1}{\Delta x}, x) \cosh_{\Delta x}(\frac{\tanh b\Delta x}{\Delta x}, x)$ | $r = \frac{\cosh b\Delta x - 1}{\Delta x} \pm \frac{\sinh b\Delta x}{\Delta x}$ |

# TABLE 14

## The Method of Related Functions Table and Block Diagram Description

**Functions with the same $K_{\Delta x}$ Transform**

$A_m, B_m, a, b, c$ are constants

| | Discrete Functions | Laplace/$K_{\Delta x}$ Transform |
|---|---|---|
| **1** | Discrete calculus<br>$y(x)$<br>Calculus<br>$y(x)$<br>See Note 9 | $Y(s)$ |
| **2** | Discrete calculus<br>$D_{\Delta x}^{n}y(x) + A_{n-1}D_{\Delta x}^{n-1}y(x) + A_{n-2}D_{\Delta x}^{n-2}y(x) + \ldots + A_1 D_{\Delta x}y(x) + A_0 y(x)$<br>Calculus<br>$\frac{d^n}{dx^n}y(x) + A_{n-1}\frac{d^{n-1}}{dx^{n-1}}y(x) + A_{n-2}\frac{d^{n-2}}{dx^{n-2}}y(x) + \ldots + A_1\frac{d}{dx}y(x) + A_0 y(x)$<br>See Note 10 | $s^n y(s) + A_{n-1}s^{n-1}y(s) + A_{n-2}s^{n-2}y(s) + \ldots + A_1 sy(s) + A_0 y(s) - B_{n-1}s^{n-1} - B_{n-2}s^{n-2} - B_{n-3}s^{n-3} - \ldots - B_1 s - B_0$<br>$B_m$, m=0,1,2,…,n-1 are intial condition constants |
| **3** | Discrete calculus<br>$c$<br>Calculus<br>$c$ | $\frac{c}{s}$ |
| **4** | Discrete calculus<br>$x = [x]_{\Delta x}^{1}$<br>Calculus<br>$x$ | $\frac{1}{s^2}$ |
| **5** | Discrete calculus<br>$[x]_{\Delta x}^{n}$<br>Calculus<br>$x^n$<br>$n = 0,1,2,3, \ldots$ | $\frac{n!}{s^{n+1}}$ |
| **6** | Discrete calculus<br>$e_{\Delta x}(a,x)$<br>Calculus<br>$e^{ax}$ | $\frac{1}{s - a}$<br>root $s = a$ |

| | Discrete Functions | Laplace/$K_{\Delta x}$ Transform |
|---|---|---|
| **7** | Discrete calculus<br>$\sin_{\Delta x}(b,x)$<br>Calculus<br>$\sin bx$ | $\frac{b}{s^2+b^2}$<br>roots s = jb, -jb |
| **8** | Discrete calculus<br>$\cos_{\Delta x}(b,x)$<br>Calculus<br>$\cos bx$ | $\frac{s}{s^2+b^2}$<br>roots s = jb, -jb |
| **9** | Discrete calculus<br>$a \neq -\frac{1}{\Delta x}$<br>$e_{\Delta x}(a,x)\ \sin_{\Delta x}(\frac{b}{1+a\Delta x}, x)$<br>$a = -\frac{1}{\Delta x}$<br>$[b\Delta x]^{\frac{x}{\Delta x}} \sin\frac{\pi x}{2\Delta x}$<br>Calculus<br>$e^{ax}\sin bx$ | $\frac{b}{(s-a)^2+b^2}$<br>roots s = a+jb, a-jb |
| **10** | Discrete calculus<br>$a \neq -\frac{1}{\Delta x}$<br>$e_{\Delta x}(a,x)\ \cos_{\Delta x}(\frac{b}{1+a\Delta x}, x)$<br>$a = -\frac{1}{\Delta x}$<br>$[b\Delta x]^{\frac{x}{\Delta x}} \cos\frac{\pi x}{2\Delta x}$<br>Calculus<br>$e^{ax}\cos bx$ | $\frac{s-a}{(s-a)^2+b^2}$<br>roots s = a+jb, a-jb |
| **11** | Discrete calculus<br>$\sinh_{\Delta x}(b,x)$<br>Calculus<br>$\sinh bx$ | $\frac{b}{s^2-b^2}$<br>roots s = b, -b |
| **12** | Discrete calculus<br>$\cosh_{\Delta x}(b,x)$<br>Calculus<br>$\cosh bx$ | $\frac{s}{s^2-b^2}$<br>roots s = b, -b |
| **13** | Discrete calculus<br>$x\sin_{\Delta x}(b,x)$<br>Calculus<br>$x\sin bx+[b\Delta x]x\cos bx$ | $\frac{2bs+[b\Delta x](s^2-b^2)}{(s^2+b^2)^2}$ |

| | Discrete Functions | Laplace/$K_{\Delta x}$ Transform |
|---|---|---|
| **14** | Discrete calculus<br>$x\cos_{\Delta x}(b,x)$<br>Calculus<br>$x\cos bx-[b\Delta x]x\sin bx$ | $\frac{(s^2-b^2)-[b\Delta x]2bs}{(s^2+b^2)^2}$ |
| **15** | Discrete calculus<br>$x\sin_{\Delta x}(b,x-\Delta x)$<br>or<br>$\frac{-[b\Delta x]x\cos_{\Delta x}(b,x)+x\sin_{\Delta x}(b,x)}{1+[b\Delta x]^2}$<br>Calculus<br>$x\sin bx$ | $\frac{2bs}{(s^2+b^2)^2}$ |
| **16** | Discrete calculus<br>$x\cos_{\Delta x}(b,x-\Delta x)$<br>or<br>$\frac{x\cos_{\Delta x}(b,x)+[b\Delta x]x\sin_{\Delta x}(b,x)}{1+[b\Delta x]^2}$<br>Calculus<br>$x\cos bx$ | $\frac{s^2-b^2}{(s^2+b^2)^2}$ |
| **17** | Discrete calculus<br>$[x]_{\Delta x}^{n}e_{\Delta x}(a,x-n\Delta x)$<br>or<br>$\frac{1}{(1+a\Delta x)^n}[x]_{\Delta x}^{n}e_{\Delta x}(a,x)$<br>Calculus<br>$x^n e^{ax}$<br>$n = 0,1,2,3,\ \ldots$ | $\frac{n!}{(s-a)^{n+1}}$ |
| **18** | Discrete calculus<br>$[x]_{\Delta x}^{n}\ \sin_{\Delta x}(b,x-n\Delta x)$<br>or | $\frac{n!}{2j}\left[\frac{(s+jb)^{n+1}-(s-jb)^{n+1}}{(s^2+b^2)^{n+1}}\right]$ |

| | Discrete Functions | Laplace/$K_{\Delta x}$ Transform |
|---|---|---|
| | $[x]^{n}_{\Delta x} \cos_{\Delta x}(b,-n\Delta x) \sin_{\Delta x}(b,x) +$<br>$[x]^{n}_{\Delta x} \sin_{\Delta x}(b,-n\Delta x) \cos_{\Delta x}(b,x)$<br>Calculus<br>$x^{n} \sin bx$<br>$n = 0,1,2,3, \ldots$ | |
| **19** | Discrete calculus<br>$[x]^{n}_{\Delta x} \cos_{\Delta x}(b, x-n\Delta x)$<br>or<br>$[x]^{n}_{\Delta x} \cos_{\Delta x}(b,-n\Delta x)\cos_{\Delta x}(b,x) -$<br>$[x]^{n}_{\Delta x} \sin_{\Delta x}(b,-n\Delta x)\sin_{\Delta x}(b,x)$<br>Calculus<br>$x^{n} \cos bx$<br>$n = 0,1,2,3, \ldots$ | $\frac{n!}{2}\left[\frac{(s+jb)^{n+1} + (s-jb)^{n+1}}{(s^2+b^2)^{n+1}}\right]$ |
| **20** | Discrete calculus<br>$x\sinh_{\Delta x}(b,x-\Delta x)$<br>or<br>$\frac{-[b\Delta x]x\cosh_{\Delta x}(b,x)+x\sinh_{\Delta x}(b,x)}{1-[b\Delta x]^2}$<br>Calculus<br>$x\sinh bx$ | $\frac{2bs}{(s^2-b^2)^2}$ |
| **21** | Discrete calculus<br>$x\cosh_{\Delta x}(b,x-\Delta x)$<br>or<br>$\frac{x\cosh_{\Delta x}(b,x)-[b\Delta x]x\sinh_{\Delta x}(b,x)}{1-[b\Delta x]^2}$<br>Calculus<br>$x\cosh bx$ | $\frac{s^2+b^2}{(s^2-b^2)^2}$ |

| | Discrete Functions | Laplace/$K_{\Delta x}$ Transform |
|---|---|---|
| **22** | Discrete calculus<br>$x\sinh_{\Delta x}(b,x)$<br>Calculus<br>$x\sinh bx+[b\Delta x]x\cosh bx$ | $\frac{2bs+[b\Delta x](s^2+b^2)}{(s^2-b^2)^2}$ |
| **23** | Discrete calculus<br>$x\cosh_{\Delta x}(b,x)$<br>Calculus<br>$x\cosh bx+[b\Delta x]x\sinh bx$ | $\frac{(s^2+b^2)+[b\Delta x]2bs}{(s^2-b^2)^2}$ |
| **24** | Discrete calculus<br>$e^{ax}$<br>Calculus<br>$e^{\frac{e^{a\Delta x}-1}{\Delta x}x}$<br>See Note 1 at the end of this table | $\frac{1}{s-\frac{e^{a\Delta x}-1}{\Delta x}}$ |
| **25** | Discrete calculus<br>$A^x$<br>Calculus<br>$e^{\frac{A^{\Delta x}-1}{\Delta x}x}$<br>See Note 2 at the end of this table | $\frac{1}{s-\frac{A^{\Delta x}-1}{\Delta x}}$ |
| **26** | Discrete calculus<br>$[x]_{\Delta x}^{n}\, e^{a(x-n\Delta x)}$<br>or<br>$e^{-an\Delta x}[x]_{\Delta x}^{n}\, e^{ax}$<br>Calculus<br>$x^n e^{\frac{e^{a\Delta x}-1}{\Delta x}x}$<br>$n = 0,1,2,3,\ \ldots$ | $\frac{n!}{(s-\frac{e^{a\Delta x}-1}{\Delta x})^{n+1}}$ |
| **27** | Discrete calculus<br>$\sin bx$<br>Calculus<br>$e^{\frac{\cos b\Delta x-1}{\Delta x}x}\sin(\frac{\sin b\Delta x}{\Delta x}x)$<br>See Note 3 at the end of this table | $\frac{\frac{\sin b\Delta x}{\Delta x}}{(s-\frac{\cos b\Delta x-1}{\Delta x})^2+(\frac{\sin b\Delta x}{\Delta x})^2}$ |

| | Discrete Functions | Laplace/$K_{\Delta x}$ Transform |
|---|---|---|
| **28** | Discrete calculus<br>$\cos bx$<br>Calculus<br>$e^{\frac{\cos b\Delta x-1}{\Delta x}x}\cos(\frac{\sin b\Delta x}{\Delta x}x)$<br>See Note 4 at the end of this table | $\frac{s-\frac{\cos b\Delta x-1}{\Delta x}}{(s-\frac{\cos b\Delta x-1}{\Delta x})^2+(\frac{\sin b\Delta x}{\Delta x})^2}$ |
| **29** | Discrete calculus<br>$e^{ax}\sin bx$<br>Calculus<br>$e^{\frac{e^{a\Delta x}\cos b\Delta x-1}{\Delta x}x}\sin(\frac{e^{a\Delta x}\sin b\Delta x}{\Delta x}x)$<br>See Note 5 at the end of this table | $\frac{\frac{e^{a\Delta x}\sin b\Delta x}{\Delta x}}{(s-\frac{e^{a\Delta x}\cos b\Delta x-1}{\Delta x})^2+(\frac{e^{a\Delta x}\sin b\Delta x}{\Delta x})^2}$ |
| **30** | Discrete calculus<br>$e^{ax}\cos bx$<br>Calculus<br>$e^{\frac{e^{a\Delta x}\cos(b\Delta x)-1}{\Delta x}x}\cos(\frac{e^{a\Delta x}\sin(b\Delta x)}{\Delta x}x)$<br>See Note 6 at the end of this table | $\frac{s-\frac{e^{a\Delta x}\cos b\Delta x-1}{\Delta x}}{(s-\frac{e^{a\Delta x}\cos b\Delta x-1}{\Delta x})^2+(\frac{e^{a\Delta x}\sin b\Delta x}{\Delta x})^2}$ |
| **31** | Discrete calculus<br>$\sinh bx$<br>Calculus<br>$e^{\frac{\cosh b\Delta x-1}{\Delta x}x}\sinh(\frac{\sinh b\Delta x}{\Delta x}x)$<br>See Note 7 at the end of this table | $\frac{\frac{\sinh b\Delta x}{\Delta x}}{(s-\frac{\cosh b\Delta x-1}{\Delta x})^2-(\frac{\sinh b\Delta x}{\Delta x})^2}$ |
| **32** | Discrete calculus<br>$\cosh bx$<br>Calculus<br>$e^{\frac{\cosh b\Delta x-1}{\Delta x}x}\cosh(\frac{\sinh b\Delta x}{\Delta x}x)$<br>See Note 8 at the end of this table | $\frac{s-\frac{\cosh b\Delta x-1}{\Delta x}}{(s-\frac{\cosh b\Delta x-1}{\Delta x})^2-(\frac{\sinh b\Delta x}{\Delta x})^2}$ |

Notes

Notes 1 thru 8 are identities

1. $e^{ax} = e_{\Delta x}(\frac{e^{a\Delta x}-1}{\Delta x}, x)$

2. $A^x = e_{\Delta x}(\frac{A^{\Delta x}-1}{\Delta x}, x)$

3. $\sin bx = e_{\Delta x}(\frac{\cos b\Delta x - 1}{\Delta x}, x)\sin_{\Delta x}(\frac{\tan b\Delta x}{\Delta x}, x)$, $\frac{x}{\Delta x}$ = integer , $\cos(b\Delta x) \neq 0$

4. $\cos bx = e_{\Delta x}(\frac{\cos b\Delta x - 1}{\Delta x}, x) \cos_{\Delta x}(\frac{\tan b\Delta x}{\Delta x}, x)$ , $\frac{x}{\Delta x}$ = integer , $\cos(b\Delta x) \neq 0$

5. $e^{ax} \sin bx = e_{\Delta x}(\frac{e^{a\Delta x}\cos b\Delta x - 1}{\Delta x}, x)\sin_{\Delta x}(\frac{\tan b\Delta x}{\Delta x}, x)$ , $\frac{x}{\Delta x}$ = integer , $\cos(b\Delta x) \neq 0$

6. $e^{ax} \cos bx = e_{\Delta x}(\frac{e^{a\Delta x}\cos b\Delta x - 1}{\Delta x}, x)\sin_{\Delta x}(\frac{\tan b\Delta x}{\Delta x}, x)$ , $\frac{x}{\Delta x}$ = integer , $\cos(b\Delta x) \neq 0$

7. $\sinh bx = e_{\Delta x}(\frac{\cosh b\Delta x - 1}{\Delta x}, x) \sinh_{\Delta x}(\frac{\tanh b\Delta x}{\Delta x}, x)$

8. $\cosh bx = e_{\Delta x}(\frac{\cosh b\Delta x - 1}{\Delta x}, x) \cosh_{\Delta x}(\frac{\tanh b\Delta x}{\Delta x}, x)$

9. The table upper discrete calculus function, y(x), is associated with discrete differential equations and the lower Calculus function, *y*(x), is associated with Calculus differential equations. The relationship between the two functions is $\lim_{\Delta x \to 0} f_1(x) = \lim_{\Delta x \to 0} f_2(x)$.

10. For the Related Functions discrete calculus and Calculus differential equations, the value of the respective initial conditions is the same (i.e. y(0) = *y*(0), $D_{\Delta x}y(0) = \frac{d}{dx}y(0)$, $D_{\Delta x}^2\, y(0) = \frac{d^2}{dx^2}y(0)$, etc.).

Description of the Method of Related Functions

Using the Method of Related Functions, a discrete calculus differential equation can be solved by first converting it into its related Calculus differential equation using the related functions in Table 14. This related Calculus differential equation is then solved by means of any of the usual Calculus methods. When the solution function to the Calculus differential equation is obtained, it is converted into its related discrete calculus function using the related functions in Table 14. The resulting function is the solution to the original discrete calculus differential equation to be solved.

The Method of Related Functions is shown in detail in the block diagram which follows. In this block diagram the function relationships to the Laplace/$K_{\Delta x}$ Transform are shown in dotted lines. However, in the application of the method, the transform relationships are not directly involved.

Block Diagram The Method of Related Functions used to solve differential difference equations

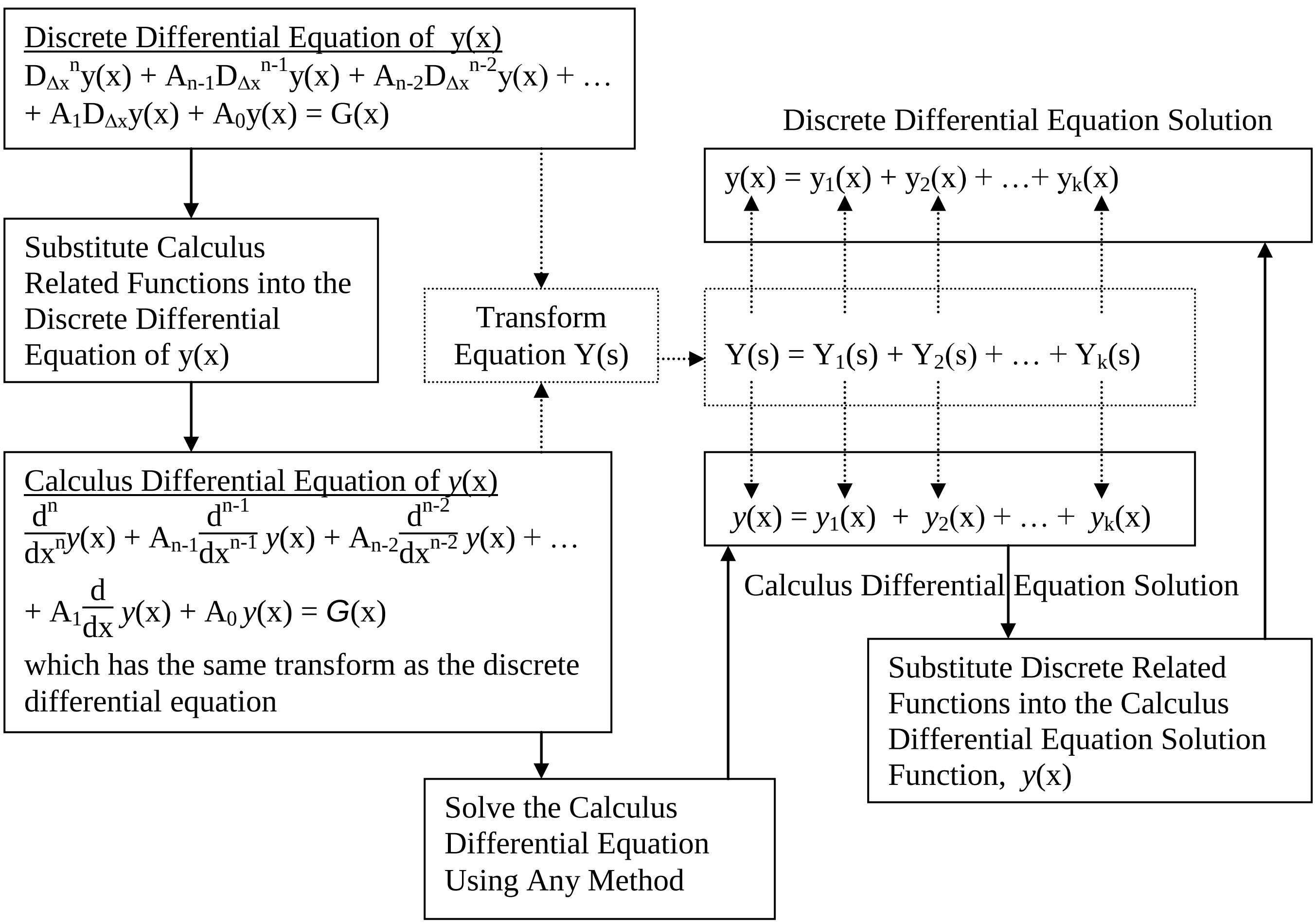

The block diagram solid lines represent the Method of Related Functions process. The block diagram dotted lines show the relationship of the method process to the Laplace/$K_{\Delta x}$ Transform. The process itself does not involve the use of transforms unless the Laplace Transform is used to find the solution to the related Calculus differential equation.

# TABLE 15
## Stability Criteria

| | |
|---|---|
| **1** | **$D_{\Delta x}^{1}y(x) + Ay(x) = 0$ First Order Differential Difference Equation with root, a**<br><br>$y(x) = Ke_{\Delta x}(a,x) = K(1+a\Delta x)^{\frac{x}{\Delta x}}$<br>where<br>$x = m\Delta x$ , $m = 0,1,2,3,\ldots$<br>K,a = real constants<br>a = equation root<br>$\Delta x$ = x increment , $\Delta x > 0$<br><br><u>Criteria for stability</u> Complex Plane<br>$-\frac{2}{\Delta x} < a < 0$ $-\frac{2}{\Delta x}$ ——— $-\frac{1}{\Delta x}$ ——— \| 0 x axis<br>or limit center point limit<br>$-1 < 1+a\Delta x < +1$ |
| **2** | **$D_{\Delta x}^{2}y(x) + AD_{\Delta x}^{1}y(x) + By(x) = 0$ Second Order Differential Difference Equation with real roots, a,b**<br><br>$y(x) = K_1e_{\Delta x}(a,x) + K_2e_{\Delta x}(b,x) = K_1(1+a\Delta x)^{\frac{x}{\Delta x}} + K_2(1+b\Delta x)^{\frac{x}{\Delta x}}$<br>where<br>$x = m\Delta x$ , $m = 0,1,2,3,\ldots$<br>$K_1,K_2$,a,b = real constants<br>a,b = equation real roots<br>$\Delta x$ = x increment , $\Delta x > 0$<br><br><u>Criteria for stability</u> Complex Plane<br>$-\frac{2}{\Delta x} < a,b < 0$ $-\frac{2}{\Delta x}$ ——— $-\frac{1}{\Delta x}$ ——— \| 0 x axis<br>or limit center point limit<br>$-1 < 1+a\Delta x < +1$ and $-1 < 1+b\Delta x < +1$ |

**3** **$D_{\Delta x}^{2}y(x) + AD_{\Delta x}^{1}y(x) + By(x) = 0$ Second Order Differential Difference Equation with complex roots, a+jb, a-jb and $1+a\Delta x = 0$**

$$y(x) = (b\Delta x)^{\frac{x}{\Delta x}} \left[ K_1 \cos \frac{\pi x}{2\Delta x} + K_2 \sin \frac{\pi x}{2\Delta x} \right]$$

where

$x = m\Delta x$ , $m = 0,1,2,3,\ldots$

$K_1, K_2, a, b =$ real constants

$1+a\Delta x = 0$ or $a = -\frac{1}{\Delta x}$

$a+jb$, $a-jb$ = equation roots

$\Delta x$ = x increment , $\Delta x > 0$

Complex Plane

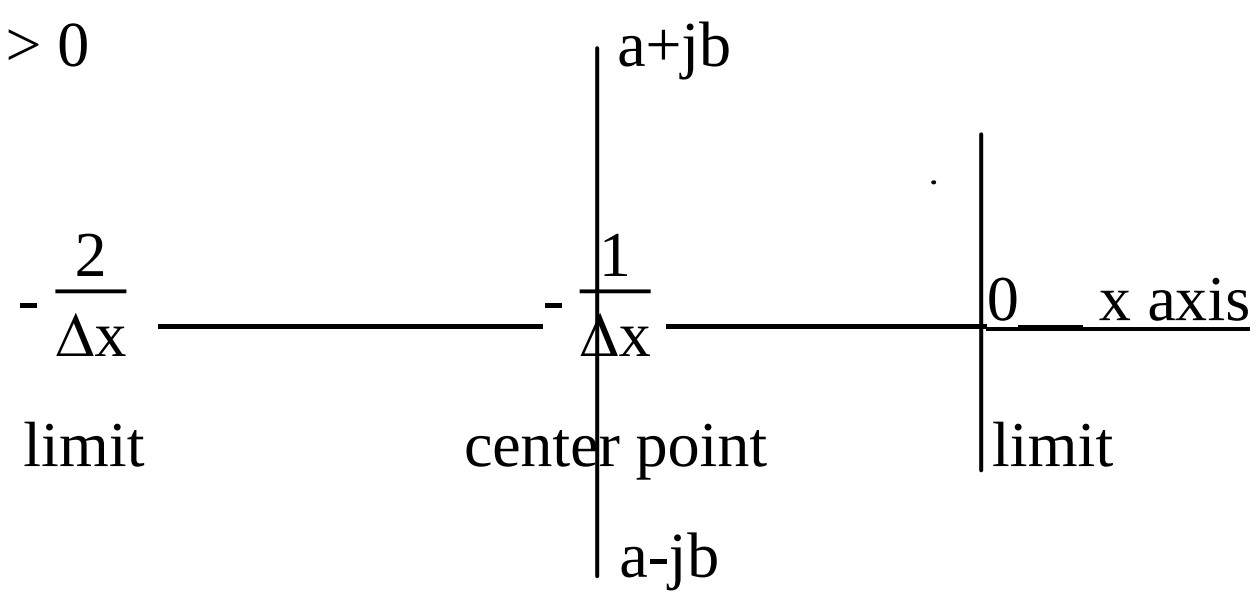

Criteria for stability

$$-\frac{1}{\Delta x} < b < +\frac{1}{\Delta x}$$

or

$$-1 < b\Delta x < +1$$

**4** **$D_{\Delta x}^{2}y(x) + AD_{\Delta x}^{1}y(x) + By(x) = 0$ Second Order Differential Difference Equation with complex roots, a+jb, a-jb and $1+a\Delta x \neq 0$**

$$y(x) = e_{\Delta x}(a,x) \left[ K_1 \cos_{\Delta x}\left(\frac{b}{1+a\Delta x}, x\right) + K_2 \sin_{\Delta x}\left(\frac{b}{1+a\Delta x}, x\right) \right] , \quad 1 + a\Delta x \neq 0$$

or

$$y(x) = (1+a\Delta x)^{\frac{x}{\Delta x}} \left[ K_1 \cos_{\Delta x}\left(\frac{b}{1+a\Delta x}, x\right) + K_2 \sin_{\Delta x}\left(\frac{b}{1+a\Delta x}, x\right) \right] , \quad 1 + a\Delta x \neq 0$$

or

$$y(x) = (1+a\Delta x)^{\frac{x}{\Delta x}} \left[ K_1 \left(1+\left[\frac{b\Delta x}{1+a\Delta x}\right]^2\right)^{\frac{x}{2\Delta x}} \cos\frac{\beta x}{\Delta x} + K_2 \left(1+\left[\frac{b\Delta x}{1+a\Delta x}\right]^2\right)^{\frac{x}{2\Delta x}} \sin\frac{\beta x}{\Delta x} \right] , \quad 1 + a\Delta x \neq 0$$

$$y(x) = \left[ \sqrt{(1+a\Delta x)^2 + (b\Delta x)^2} \right]^{\frac{x}{\Delta x}} \left[ K_1 \cos\frac{\beta x}{\Delta x} + K_2 \sin\frac{\beta x}{\Delta x} \right]$$

where

$x = m\Delta x$ , $m = 0,1,2,3,\ldots$

$K_1, K_2, a, b =$ real constants

$$\beta = \begin{cases} \tan^{-1}\frac{b\Delta x}{1+a\Delta x} & \text{for } 1+a\Delta x \geq 0 \\ \pi + \tan^{-1}\frac{b\Delta x}{1+a\Delta x} & \text{for } 1+a\Delta x < 0 \end{cases}$$

$a+jb$, $a-jb$ = equation roots

$\Delta x$ = x increment , $\Delta x > 0$

<u>Criteria for stability</u>

$0 \le (1+a\Delta x)^2 + (b\Delta x)^2 < 1$

or

$(a + \frac{1}{\Delta x})^2 + b^2 < (\frac{1}{\Delta x})^2$ (equation of a circle)

This above equation specifies that, for stability, the equation roots, a+jb and a-jb, must lie within a circle within the left half of the complex plane centered at $x = -\frac{1}{\Delta x}$ with a radius of $\frac{1}{\Delta x}$.

5

## General Criteria for Stability

**For stability, all $K_{\Delta x}$ Transform (transfer function) denominator roots must lie within the circle situated in the left half of the complex plane centered at $x = -\frac{1}{\Delta x}$ with a radius of $\frac{1}{\Delta x}$.**

$y(s) = K_{\Delta x}[y(x)]$ , The $K_{\Delta x}$ Transform of y(x)

$y(s) = \frac{K_1}{s+a}$ or $y(s) = \frac{K_1 s + K_2}{(s + a + jb)(s + a - jb)}$

**Root Location Stability Diagram**

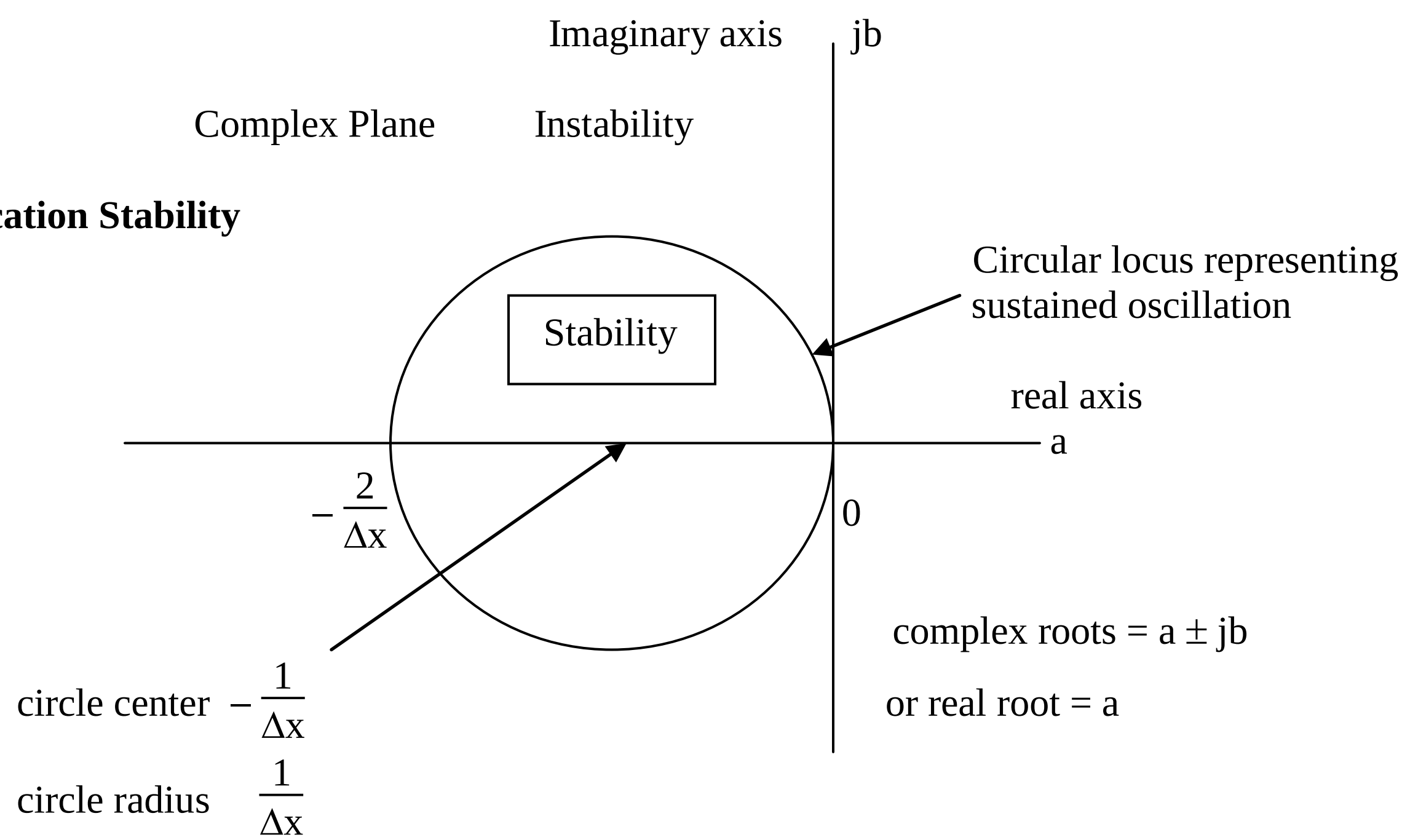

| | |
|---|---|
| 6 | **<u>The Modified Nyquist Criterion as it applies to discrete (sampled) variable control systems</u>** |

**($\Delta x \neq 0$)**

**Evaluate, in a counterclockwise direction, A(s) for s = all points on the perimeter of the Critical Circle of radius $\frac{1}{\Delta x}$ centered at $-\frac{1}{\Delta x}$. If a pole of A(s) lies on the Critical Circle perimeter, semicircle it very closely from outside the Critical Circle. A pole of A(s) on the Critical Circle perimeter is considered to be within the Critical Circle.**

**<u>Note</u> – $\Delta x$ , the x increment, may be designated $\Delta t$, a time increment, for time related variables.**

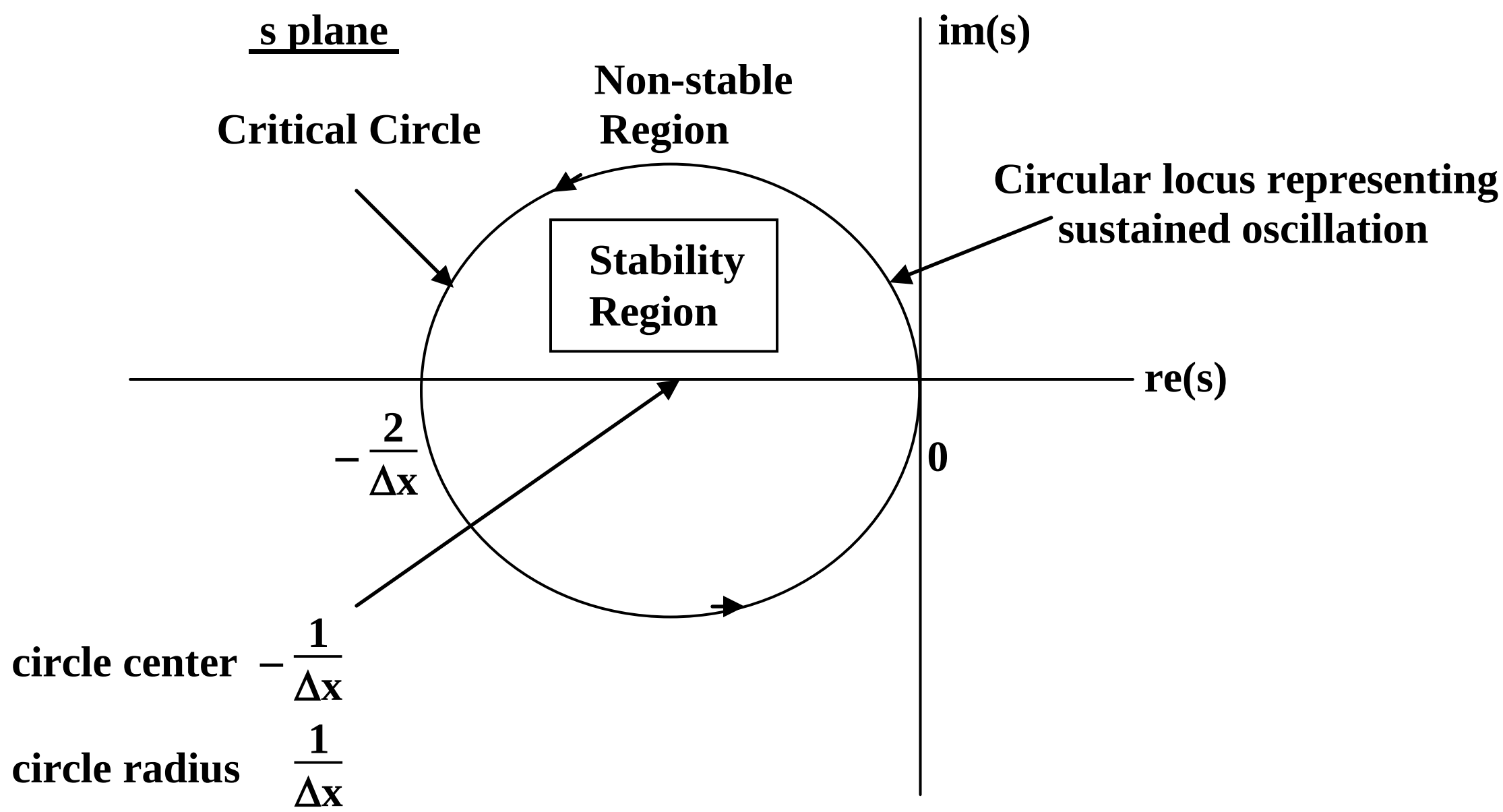

**<u>For control system stability</u>**

**where**

**N = the number of poles of 1+A(s) (or A(s)) and the number of zeros of 1+A(s).**

**P = the number of poles of 1+A(s) (or A(s)) located outside the Critical Circle which lies within the left half of the s plane.**

**If N poles of A(s) are located within the Critical Circle of the s plane, there must be no encirclement of the A(s) plane –1 point.**

**If N-P poles of A(s) are located within the Critical Circle of the s plane, there must be N-P counterclockwise encirclements of the A(s) plane –1 point.**

**For each zero of 1+A(s) located on the Critical Circle perimeter, there will be a transition through the –1 point of the A(s) plane.**

**<u>Notes</u> – The poles of A(s) are the same as those of 1+A(s).**

**For the functions, A(s) and 1+A(s), there can be no more zeros than poles.**

**7** **<u>The Modified Nyquist Criterion as it applies to continuous variable control systems</u> ($\Delta x \rightarrow 0$)**

**Evaluate, in a counterclockwise direction, A(s) for s = all points on the perimeter of the left half of the s plane. If a pole of A(s) lies on the perimeter, semicircle it very closely from outside the left half of the s plane. A pole of A(s) on the s plane imaginary axis is considered to be within the left half of the s plane.**

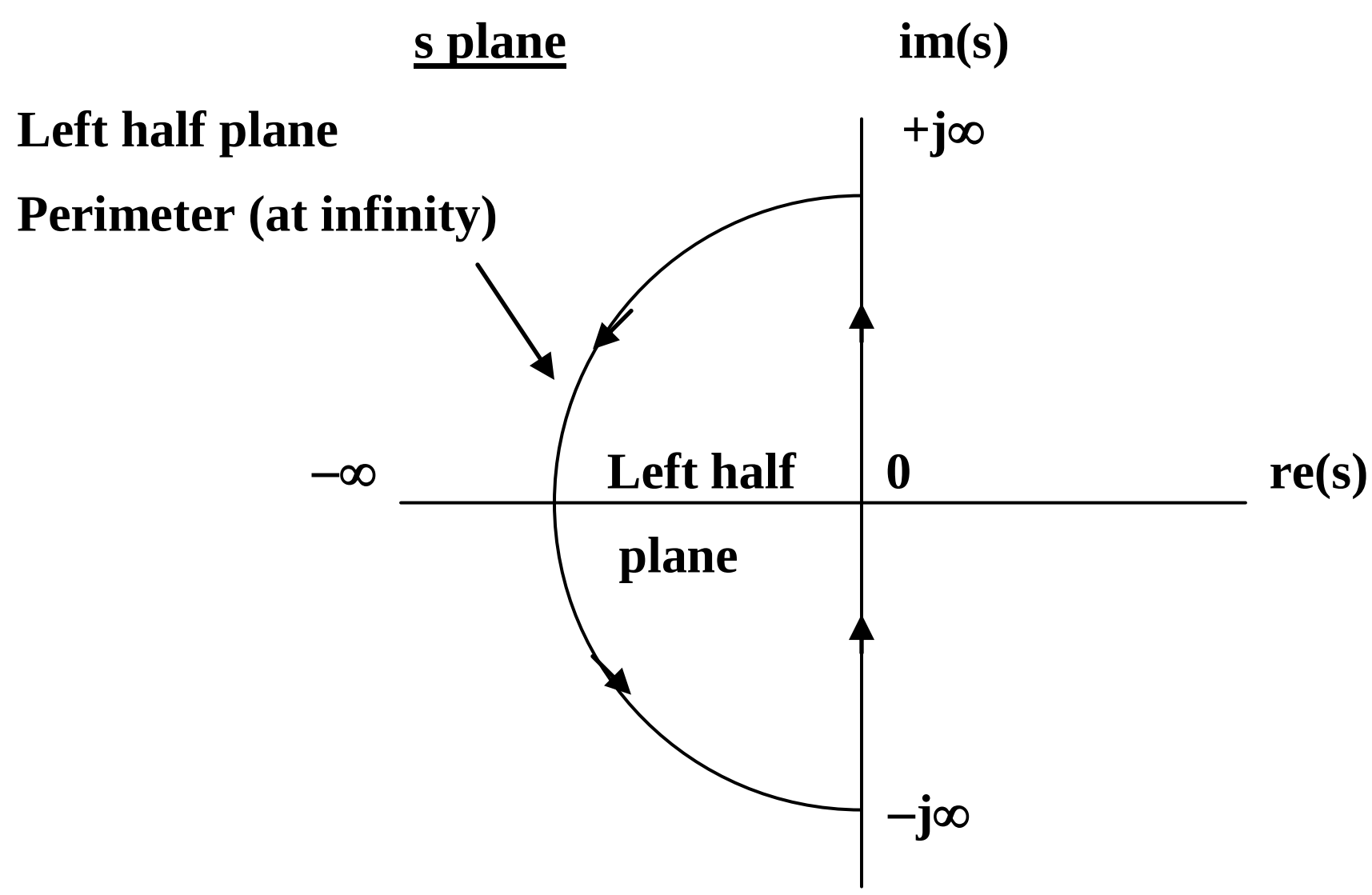

**<u>For control system stability</u>**

**where**

**N = the number of poles of 1+A(s) (or A(s)) and the number of zeros of 1+A(s).**

**P = the number of poles of 1+A(s) (or A(s)) outside the left half of the s plane.**

**If N poles of A(s) are located within the left half of the s plane, there must be no encirclement of the A(s) plane –1 point.**

**If N-P poles of A(s) are located within the left half of the s plane, there must be N-P counterclockwise encirclements of the A(s) plane –1 point.**

**For each zero of 1+A(s) located on the imaginary axis of the s plane, there will be a transition through the –1 point of the A(s) plane.**

|  | **<u>Notes</u> – The poles of A(s) are the same as those of 1+A(s).**<br>**For the functions, A(s) and 1+A(s), there can be no more zeros than poles.**<br><br>**<u>Comment</u> – The Modified Nyquist Criterion as applied to continuous variable control systems yields exactly the same A(s) plot and stability analysis data as the Nyquist Criterion.** |
|---|---|

# TABLE 16

## A listing of Zeta Function, Digamma Function, Polygamma Function, and lnd(n,Δx,x) Function Equations

### Zeta Function, Digamma Function, Polygamma Function, and lnd(n,Δx,x) Function Equations

**The General Zeta Function**

**1** The relationship of the General Zeta Function to the lnd(n,Δx,x) function

$$\zeta(n,\Delta x,x) = \frac{1}{\Delta x}\,\mathrm{lnd}(n,\Delta x,x)$$

where
- $x$ = real or complex variable
- $n,\Delta x$ = real or complex constants
- $\Delta x$ = x increment

**1a** The General Zeta Function Summation Between Finite Limits

$${}_{\Delta x}\sum_{x=x_1}^{x_2}\frac{1}{x^n} = \pm\,\zeta(n,\Delta x,x)\Big|_{x_1}^{x_2+\Delta x} = \pm\frac{1}{\Delta x}\,\mathrm{lnd}(n,\Delta x,x)\Big|_{x_1}^{x_2+\Delta x}, \quad -\text{ for } n\neq 1,\ +\text{ for } n=1$$

where
- $x$ = real or complex variable
- $x_1, x_2, \Delta x, n$ = real or complex constants
- $\Delta x$ = x increment
- Any summation term where $x = 0$ is excluded

or

$$\zeta(n,\Delta x,x)\Big|_{x_1}^{x_2} = \pm\,{}_{\Delta x}\sum_{x=x_1}^{x_2-\Delta x}\frac{1}{x^n} = \frac{1}{\Delta x}\,\mathrm{lnd}(n,\Delta x,x)\Big|_{x_1}^{x_2}, \quad -\text{ for } n\neq 1,\ +\text{ for } n=1$$

where

$${}_{\Delta x}\sum_{x=x_1}^{x_1-\Delta x}\frac{1}{x^n} = 0$$

- $x$ = real or complex variable
- $x_1, x_2, \Delta x, n$ = real or complex constants
- $\Delta x$ = x increment
- Any summation term where $x = 0$ is excluded

| | **Zeta Function, Digamma Function, Polygamma Function, and lnd(n,Δx,x) Function Equations** |
|---|---|
| | The following equation, which relates the lnd(n,Δx,x) Function to the Hurwitz Zeta Function, may be used.<br><br>$\mathrm{lnd}(n,\Delta x,x) = \frac{1}{\Delta x^{n-1}}\,\zeta(n,\frac{x}{\Delta x})$<br><br>where<br>x,Δx = real values with the condition that x ≠ positive real value when Δx = negative real value<br>n = complex values<br>or<br>x = Δx = complex value<br>n = complex value<br>or<br>x,Δx = complex values<br>n = integer<br><br>ζ(n,x) = Hurwitz Zeta Function , n ≠ 1<br><br><u>Comments</u> – A complex value, x+jy, can be a real value, x, an imaginary value, jy, or an integer, N+j0.<br><br>For value combinations not specified above, the equality of the stated equation may not be valid. |
| 1b | <u>The General Zeta Function Definite Discrete Integral Between Finite Limits</u><br><br>$_{\Delta x}\int_{x_1}^{x_2} \frac{1}{x^n}\,\Delta x = \Delta x\, {}_{\Delta x}\sum_{x=x_1}^{x_2-\Delta x} \frac{1}{x^n} = \pm\,\Delta x\,\zeta(n,\Delta x,x)\,\Big\|_{x_1}^{x_2} = \pm\,\mathrm{lnd}(n,\Delta x,x)\,\Big\|_{x_1}^{x_2}$ , – for n≠1, + for n=1<br><br>where<br>$_{\Delta x}\sum_{x=x_1}^{x_1-\Delta x} \frac{1}{x} = 0$<br>x = real or complex variable<br>$x_1,x_2,\Delta x,n$ = real or complex constants<br>Δx = x increment<br>Any summation term where x = 0 is excluded |

| | **Zeta Function, Digamma Function, Polygamma Function, and lnd(n,Δx,x) Function Equations** |
|---|---|
| 1c | The General Zeta Function Indefinite Discrete Integral<br><br>$\int_{\Delta x} \frac{1}{x^n}\Delta x = \pm\,\Delta x\,\zeta(n,\Delta x,x) + k(n,\Delta x) = \pm\,\mathrm{lnd}(n,\Delta x,x) + k(n,\Delta x)$ ,   $-$ for $n\neq1$, $+$ for $n=1$<br>where<br>x = real or complex variable<br>Δx,n = real or complex constants<br>Δx = x increment<br>k = constant of integration, a function of n,Δx |
| 1d | The General Zeta Function Summation to Infinity<br><br>$\sum_{\Delta x\; x=x_i}^{\pm\infty} \frac{1}{x^n} = \zeta(n,\Delta x,x_i) = \frac{1}{\Delta x}\,\mathrm{lnd}(n,\Delta x,x_i)$ ,   $Re(n)>1$<br><br>where<br>x = real or complex variable<br>$x_i$,n,Δx = real or complex constants<br>Δx = x increment<br>+∞ for Re(Δx)>0 or [ Re(Δx)=0 and Im(Δx)>0 ]<br>–∞ for Re(Δx)<0 or [ Re(Δx)=0 and Im(Δx)<0 ]<br>Any summation term where x = 0 is excluded |
| 1e | The General Zeta Function Definite Discrete Integral to Infinity<br><br>$\int_{\Delta x\; x=x_i}^{\pm\infty} \frac{1}{x^n}\Delta x = \Delta x \sum_{\Delta x\; x=x_i}^{\pm\infty} \frac{1}{x^n} = \Delta x\,\zeta(n,\Delta x,x_i) = \mathrm{lnd}(n,\Delta x,x_i)$ ,   $Re(n)>1$<br>where<br>x = real or complex variable<br>$x_i$,n,Δx = real or complex constants<br>Δx = x increment<br>+∞ for Re(Δx)>0 or [ Re(Δx)=0 and Im(Δx)>0 ]<br>–∞ for Re(Δx)<0 or [ Re(Δx)=0 and Im(Δx)<0 ]<br>Any summation term where x = 0 is excluded |

| | **Zeta Function, Digamma Function, Polygamma Function, and lnd(n,Δx,x) Function Equations** |
|---|---|
| 1f | The General Zeta Function Discrete Derivative<br><br>$D_{\Delta x}\,\zeta(n,\Delta x,x) = \pm\,\frac{1}{\Delta x}\,\frac{1}{x^n}\,,\quad -\text{ for } n\neq 1,\ +\text{ for } n=1$<br><br>where<br>x = real or complex variable<br>Δx,n = real or complex constants<br>Δx = x increment |
| 1g | The General Zeta Function term to term relationship<br><br>$\zeta(n,\Delta x,x+\Delta x) - \zeta(n,\Delta x,x) = \pm\,\frac{1}{x^n}\,,\quad -\text{ for } n\neq 1,\ +\text{ for } n=1$<br><br>where<br>x = real or complex variable<br>Δx,n = real or complex constants<br>Δx = x increment |
| 1h | The General Zeta Function Recursion Equation<br><br>$\zeta(n,\Delta x,x+\Delta x) = \zeta(n,\Delta x,x) \pm \frac{1}{x^n}\,,\ -\text{ for } n\neq 1,\ +\text{ for } n=1,$<br><br>where<br>x = real or complex variable<br>Δx,n = real or complex constants<br>Δx = x increment |
| | **The Hurwitz Zeta Function** |
| 2 | The relationship of the Hurwitz Zeta Function to the lnd(n,Δx,x) function<br><br>$\zeta(n,x) = \zeta(n,1,x) = \text{lnd}(n,1,x)\ ,\ \ n \neq 1$<br><br>where<br>x = real or complex variable<br>n = real or complex constant |
| 2a | The Hurwitz Zeta Function Summation Between Finite Limits<br><br>$\sum_{1\ x=x_1}^{x_2} \frac{1}{x^n} = -\,\zeta(n,x)\Big\vert_{x_1}^{x_2+1} = -\,\zeta(n,1,x)\Big\vert_{x_1}^{x_2+1} = -\,\text{lnd}(n,1,x)\Big\vert_{x_1}^{x_2+1}\ ,\ \ n \neq 1$<br><br>where<br>x = real or complex variable<br>$x_1, x_2$,n = real or complex constants<br>Any summation term where x = 0 is excluded<br><br>or |

| | Zeta Function, Digamma Function, Polygamma Function, and lnd(n,Δx,x) Function Equations |
|---|---|
| | $$\zeta(n,x)\Big|_{x_1}^{x_2} = -{}_1\sum_{x=x_1}^{x_2-1}\frac{1}{x^n} = \zeta(n,1,x)\Big|_{x_1}^{x_2} = \mathrm{lnd}(n,1,x)\Big|_{x_1}^{x_2}, \quad n\neq 1$$ where $${}_{\Delta x}\sum_{x=x_1}^{x_1-1}\frac{1}{x^n} = 0$$ x = real or complex variable<br>$x_1$,$x_2$,n = real or complex constants<br>Any summation term where x = 0 is excluded |
| 2b | The Hurwitz Zeta Function Definite Discrete Integral Between Finite Limits<br>$${}_1\int_{x_1}^{x_2}\frac{1}{x^n}\,\Delta x = {}_1\sum_{x=x_1}^{x_2-1}\frac{1}{x^n} = -\,\zeta(n,x)\Big|_{x_1}^{x_2} = -\,\zeta(n,1,x)\Big|_{x_1}^{x_2} = -\,\mathrm{lnd}(n,1,x)\Big|_{x_1}^{x_2}, \quad n\neq 1$$ where $${}_1\sum_{x=x_1}^{x_1-1}\frac{1}{x} = 0$$ x = real or complex variable<br>$x_1$,$x_2$,n = real or complex constants<br>Δx = 1, x increment<br>Any summation term where x = 0 is excluded |
| 2c | The Hurwitz Zeta Function Indefinite Discrete Integral<br>$${}_1\int\frac{1}{x^n}\,\Delta x = -\,\zeta(n,x) + k(n) = -\,\zeta(n,1,x) + k(n) = -\,\mathrm{lnd}(n,1,x) + k(n)\,, \quad n\neq 1$$ where<br>x = real or complex variable<br>n = real or complex constant<br>Δx = 1, x increment<br>k = constant of integration, a function of n |

| | Zeta Function, Digamma Function, Polygamma Function, and lnd(n,Δx,x) Function Equations |
|---|---|
| 2d | The Hurwitz Zeta Function Summation to Infinity<br><br>$$\sum_{x=x_i}^{\infty}{}_{1}\frac{1}{x^n} = \zeta(n,x_i) = \zeta(n,1,x_i) = \text{lnd}(n,1,x_i)\,, \quad \text{Re}(n)>1$$<br>where<br>x = real or complex variable<br>$x_i$,n = real or complex constants<br>Δx = 1, x increment<br>Any summation term where x = 0 is excluded |
| 2e | The Hurwitz Zeta Function Definite Discrete Integral to Infinity<br><br>$$\int_{x=x_i}^{\infty}{}_{1}\frac{1}{x^n}\Delta x = \sum_{x=x_i}^{\infty}{}_{1}\frac{1}{x^n} = \zeta(n,x_i) = \zeta(n,1,x_i) = \text{lnd}(n,1,x_i)\,, \quad \text{Re}(n)>1$$<br>where<br>x = real or complex variable<br>$x_i$,n= real or complex constants<br>Δx = 1, x increment<br>Any summation term where x = 0 is excluded |
| 2f | The Hurwitz Zeta Function Discrete Derivative<br><br>$$D_1\zeta(n,x) = -\frac{1}{x^n}\,, \quad n\neq1$$<br>where<br>x = real or complex variable<br>n = real or complex constant |
| 2g | The Hurwitz Zeta Function term to term relationship<br><br>$$\zeta(n,x+1) - \zeta(n,x) = -\frac{1}{x^n}\,, \quad n\neq1$$<br>where<br>x = real or complex variable<br>n = real or complex constant |

| | **Zeta Function, Digamma Function, Polygamma Function, and lnd(n,Δx,x) Function Equations** |
|---|---|
| 2h | The Hurwitz Zeta Function Recursion Equation<br><br>$\zeta(n,x+1) = \zeta(n,x) - \frac{1}{x^n}$ , $n\neq1$<br>where<br>x = real or complex variable<br>n = real or complex constant |
| | **The Riemann Zeta Function** |
| 3 | The relationship of the Riemann Zeta Function to the lnd(n,Δx,x) function<br><br>$\zeta(n) = \zeta(n,1,1) = \text{lnd}(n,1,1)$ , $n \neq 1$<br>where<br>n = real or complex constant |
| 3a | The Riemann Zeta Function Summation to Infinity<br><br>$\zeta(n) = \sum_{1 \atop x=1}^{\infty} \frac{1}{x^n} = \zeta(n,1,1) = \text{lnd}(n,1,1)$ , Re(n)>1<br>where<br>x = real or complex variable<br>n = real or complex constant |
| 3b | The Riemann Zeta Function Definite Discrete Integral to Infinity<br><br>$\int_{1 \atop x=1}^{\infty} \frac{1}{x^n}\Delta x = \sum_{1 \atop x=1}^{\infty} \frac{1}{x^n} = \zeta(n) = \zeta(n,1,1) = \text{lnd}(n,1,1)$ , Re(n)>1<br>where<br>x = real or complex variable<br>n = real or complex constant<br>Δx = 1, x increment |
| | **The Polygamma Functions** |
| 4 | The relationship of the Polygamma Functions to the lnd(n,Δx,x) function<br><br>$\psi^{(m)}(x) = (-1)^{m+1}m!\,\zeta(1+m,x) = (-1)^{m+1}m!\,\zeta(1+m,1,x) = (-1)^{m+1}m!\,\text{lnd}(1+m,1,x)$<br>where<br>m = 1,2,3,…<br>x = real or complex value |

| | Zeta Function, Digamma Function, Polygamma Function, and lnd(n,Δx,x) Function Equations |
|---|---|
| 4a | The Polygamma Functions Summation Between Finite Limits<br>$$\sum_{\Delta x\; x=x_1}^{x_2}\frac{1}{x^{1+m}} = \frac{(-1)^m}{m!\Delta x^{m+1}}\psi^{(m)}\left(\frac{x}{\Delta x}\right)\Big|_{x_1}^{x_2+\Delta x} = -\zeta(1+m,\Delta x,x)\Big|_{x_1}^{x_2+\Delta x} = -\frac{1}{\Delta x}\,\mathrm{lnd}(1+m,\Delta x,x)\Big|_{x_1}^{x_2+\Delta x}$$<br>where<br>x = real or complex variable<br>$x_1, x_2, \Delta x$ = real or complex constants<br>Δx = x increment<br>m = 1,2,3,…<br>Any summation term where x = 0 is excluded<br>or<br>$$\psi^{(m)}\left(\frac{x}{\Delta x}\right)\Big|_{x_1}^{x_2} = (-1)^{m+1}m!\Delta x^{m+1}\zeta(1+m,\Delta x,x)\Big|_{x_1}^{x_2} = (-1)^m m!\Delta x^{m+1}\sum_{\Delta x\; x=x_1}^{x_2-\Delta x}\frac{1}{x^{1+m}} = (-1)^{m+1}m!\Delta x^m\,\mathrm{lnd}(1+m,\Delta x,x)\Big|_{x_1}^{x_2}$$<br>where<br>$$\sum_{\Delta x\; x=x_1}^{x_1-\Delta x}\frac{1}{x^{1+m}} = 0$$<br>x = real or complex variable<br>$x_1, x_2, \Delta x$ = real or complex constants<br>Δx = x increment<br>m = 1,2,3,…<br>Any summation term where x = 0 is excluded |
| 4b | The Polygamma Functions Definite Discrete Integral Between Finite Limits<br>$$\int_{\Delta x\; x_1}^{x_2}\frac{1}{x^{1+m}}\,\Delta x = \Delta x\sum_{\Delta x\; x=x_1}^{x_2-\Delta x}\frac{1}{x^{1+m}} = \frac{(-1)^m}{m!\Delta x^m}\psi^{(m)}\left(\frac{x}{\Delta x}\right)\Big|_{x_1}^{x_2} - \Delta x\,\zeta(1+m,\Delta x,x)\Big|_{x_1}^{x_2} = -\,\mathrm{lnd}(1+m,\Delta x,x)\Big|_{x_1}^{x_2}$$<br>where<br>$$\sum_{\Delta x\; x=x_1}^{x_1-\Delta x}\frac{1}{x} = 0$$<br>x = real or complex variable<br>$x_1, x_2, \Delta x$ = real or complex constants<br>Δx = x increment<br>m = 1,2,3,…<br>Any summation term where x = 0 is excluded |

## Zeta Function, Digamma Function, Polygamma Function, and lnd(n,Δx,x) Function Equations

4c The Polygamma Functions Indefinte Discrete Integral

$$\int_{\Delta x} \frac{1}{x^{1+m}}\,\Delta x = \frac{(-1)^m}{m!\Delta x^m}\,\psi^{(m)}\left(\frac{x}{\Delta x}\right) + k(\Delta x) = -\,\Delta x\,\zeta(1+m,\Delta x,x) + k(\Delta x) = -\,\mathrm{lnd}(1+m,\Delta x,x) + k(\Delta x)$$

where

- $x$ = real or complex variable
- $\Delta x$ = real or complex constant
- $\Delta x$ = $x$ increment
- $m = 1,2,3,\ldots$
- $k$ = constant of integration, a function of $\Delta x$

4d The Polygamma Function Summation to Infinity

$$\sum_{\Delta x \atop x=x_i}^{\pm\infty} \frac{1}{x^{1+m}} = \frac{(-1)^{m+1}}{m!\Delta x^{m+1}}\,\psi^{(m)}\left(\frac{x_i}{\Delta x}\right) = \zeta(1+m,\Delta x,x_i) = \frac{1}{\Delta x}\,\mathrm{lnd}(1+m,\Delta x,x_i)$$

where

- $x$ = real or complex variable
- $x_i,m,\Delta x$ = real or complex constants
- $+\infty$ for $\mathrm{Re}(\Delta x)>0$ or [ $\mathrm{Re}(\Delta x)=0$ and $\mathrm{Im}(\Delta x)>0$ ]
- $-\infty$ for $\mathrm{Re}(\Delta x)<0$ or [ $\mathrm{Re}(\Delta x)=0$ and $\mathrm{Im}(\Delta x)<0$ ]
- $m = 1,2,3,\ldots$
- $\Delta x$ = $x$ increment
- Any summation term where $x = 0$ is excluded

4e The Polygamma Functions Definite Discrete Integral to Infinity

$$\int_{\Delta x \atop x_i}^{\pm\infty} \frac{1}{x^{1+m}}\,\Delta x = \frac{(-1)^{m+1}}{m!\Delta x^m}\,\psi^{(m)}\left(\frac{x_i}{\Delta x}\right) = \Delta x \sum_{\Delta x \atop x=x_i}^{\pm\infty} \frac{1}{x^{1+m}} = \Delta x\,\zeta(1+m,\Delta x,x_i) = \mathrm{lnd}(1+m,\Delta x,x_i)$$

where

- $x$ = real or complex variable
- $x_i,m,\Delta x$ = real or complex constants
- $+\infty$ for $\mathrm{Re}(\Delta x)>0$ or [ $\mathrm{Re}(\Delta x)=0$ and $\mathrm{Im}(\Delta x)>0$ ]
- $-\infty$ for $\mathrm{Re}(\Delta x)<0$ or [ $\mathrm{Re}(\Delta x)=0$ and $\mathrm{Im}(\Delta x)<0$ ]
- $m = 1,2,3,\ldots$
- $\Delta x$ = $x$ increment
- Any summation term where $x = 0$ is excluded

| | **Zeta Function, Digamma Function, Polygamma Function, and lnd(n,Δx,x) Function Equations** |
|---|---|
| 4f | The Polygamma Functions Discrete Derivative<br><br>$D_1\psi^{(m)}(x) = (-1)^m m!\frac{1}{x^{m+1}}$<br>where<br>m = 1,2,3,…<br>x = real or complex variable |
| 4g | The Polygamma Functions term to term relationship<br>$\psi^{(m)}(x+1) - \psi^{(m)}(x) = (-1)^m m!\frac{1}{x^{m+1}}$<br>where<br>m = 1,2,3,…<br>x = real or complex value |
| 4h | The Polygamma Functions Recursion Equation<br><br>$\psi^{(m)}(x+1) = \psi^{(m)}(x) + (-1)^m m!\frac{1}{x^{m+1}}$<br>where<br>m = 1,2,3,…<br>x = real or complex value |
| | **The Digamma Function** |
| 5 | The relationship of the Digamma Function to the lnd(n,Δx,x) function<br><br>$\psi(x) = \text{lnd}(1,1,x) - \gamma$ , $x \neq 0,-1,-2,-3,\ldots$<br>or<br>$\psi(\frac{x}{\Delta x}) = \text{lnd}(1,\Delta x,x) - \gamma = \text{lnd}(1,1,\frac{x}{\Delta x}) - \gamma = \Delta x\zeta(1,\Delta x,x) - \gamma$ , $x \neq 0,-\Delta x,-2\Delta x,-3\Delta x,\ldots$<br><br>where<br>x = real or complex values<br>Δx = x increment<br>γ = Euler's Constant, .5772157…<br>ψ(x) is infinite for x = 0,-1,-2,-3,…<br>$\psi(\frac{x}{\Delta x})$ is infinite for $x = 0,-\Delta x,-2\Delta x,-3\Delta x,\ldots$<br>ζ(1,Δx,x) is the N=1 Zeta Function<br><br>Comments – The Digamma Function, ψ(x), has first order poles at x ≠ 0,-1,-2,-3,….<br>At these values of x, where $\psi(x) = \frac{1}{0}$ , the value of ψ(x) is infinite.<br><br>The lnd(1,1,x) function minus Euler's number is equal to the Digamma Function for all values of x except for x = 0,-1,-2,-3,… For x = 0,-1,-2,-3,…, the lnd(1,1,x) function |

differs from the $\psi(x)$ function. The equation used to calculate the lnd(1,1,x) function when
$x = 0,-1,-2,-3,\ldots$ is $\mathrm{lnd}(1,1,x) = \mathrm{lnd}(1,1,1-x)$. In addition, any summation calculated using the lnd(1,1,x) function, excludes any summation term with a division by zero. These differences were introduced into the lnd(1,1,x) function to make possible the summation of integers along the real axis of the complex plane.

For $x = 0,-1,-2,-3,\ldots$, if the program, LNDX, is selected to find $\psi(x)$ using the equation, $\psi(x) = \mathrm{lnd}(1,1,x) - \gamma$, the resulting calculated value for $\psi(x)$ is $\psi(x) = \psi(1-x)$.

5a The Digamma Function Summation Between Finite Limits

$$\sum_{\Delta x\; x=x_1}^{x_2} \frac{1}{x} = \frac{1}{\Delta x}\psi\left(\frac{x}{\Delta x}\right)\Big|_{x_1}^{x_2+\Delta x} = \zeta\left(1,1,\frac{x}{\Delta x}\right)\Big|_{x_1}^{x_2+\Delta x} = \frac{1}{\Delta x}\mathrm{lnd}\left(1,1,\frac{x}{\Delta x}\right)\Big|_{x_1}^{x_2+\Delta x} = \frac{1}{\Delta x}\mathrm{lnd}(1,\Delta x,x)\Big|_{x_1}^{x_2+\Delta x}$$

where

$x$ = real or complex variable

$x_1,x_2,\Delta x$ = real or complex constants

$\Delta x$ = x increment

$x_1 \neq 0,-\Delta x,-2\Delta x,-3\Delta x,\ldots$ , These $x_1,x_2$ exclusions apply only if the Digamma Function is used.

$x_2 \neq -\Delta x,-2\Delta x,-3\Delta x,\ldots$ , If the Digamma Function is not used, any term of the summation with a division by zero will be excluded.

or

$$\psi\left(\frac{x}{\Delta x}\right)\Big|_{x_1}^{x_2} = \Delta x \sum_{\Delta x\; x=x_1}^{x_2-\Delta x} \frac{1}{x} = \Delta x\,\zeta\left(1,1,\frac{x}{\Delta x}\right)\Big|_{x_1}^{x_2} = \mathrm{lnd}\left(1,1,\frac{x}{\Delta x}\right)\Big|_{x_1}^{x_2} = \mathrm{lnd}(1,\Delta x,x)\Big|_{x_1}^{x_2}$$

where

$$\sum_{\Delta x\; x=x_1}^{x_1-\Delta x} \frac{1}{x} = 0$$

$x$ = real or complex variable

$x_1,x_2,\Delta x,n$ = real or complex constants

. $\Delta x$ = x increment

$x_1,x_2 \neq 0,-\Delta x,-2\Delta x,-3\Delta x,\ldots$ , These $x_1,x_2$ exclusions apply only if the Digamma Function is used. If the Digamma Function is not used, any term of the summation with a division by zero will be excluded.

| | **Zeta Function, Digamma Function, Polygamma Function, and lnd(n,Δx,x) Function Equations** |
|---|---|
| 5b | The Digamma Function Definite Discrete Integral Between Finite Limits<br><br>$$\int_{\Delta x\,x_1}^{x_2} \frac{1}{x}\,\Delta x = \Delta x \sum_{\Delta x\,x=x_1}^{x_2-\Delta x} \frac{1}{x} = \psi\left(\frac{x}{\Delta x}\right)\Big|_{x_1}^{x_2} = \Delta x\,\zeta\left(1,1,\frac{x}{\Delta x}\right)\Big|_{x_1}^{x_2} = \text{lnd}\left(1,1,\frac{x}{\Delta x}\right)\Big|_{x_1}^{x_2} = \text{lnd}(1,\Delta x,x)\Big|_{x_1}^{x_2}$$<br>where<br>$$\sum_{\Delta x\,x=x_1}^{x_1-\Delta x} \frac{1}{x} = 0$$<br>$x$ = real or complex variable<br>$x_1, x_2, \Delta x$ = real or complex constants<br>$\Delta x$ = x increment<br>$x_1, x_2 \neq 0, -\Delta x, -2\Delta x, -3\Delta x, \ldots$, These $x_1, x_2$ exclusions apply only if the Digamma Function is used. If the Digamma Function is not used, any term of the integral with a division by zero will be excluded. |
| 5c | The Digamma Function Indefinite Discrete Integral<br><br>$$\int_{\Delta x} \frac{1}{x}\,\Delta x = \psi\left(\frac{x}{\Delta x}\right) + k(\Delta x) = \Delta x\,\zeta\left(1,1,\frac{x}{\Delta x}\right) + k(\Delta x) = \text{lnd}\left(1,1,\frac{x}{\Delta x}\right) + k(\Delta x) = \text{lnd}(1,\Delta x,x) + k(\Delta x)$$<br>where<br>$x$ = real or complex variable<br>$\Delta x$ = real or complex constant<br>$\Delta x$ = x increment<br>$x \neq 0, -\Delta x, -2\Delta x, -3\Delta x, \ldots$ , This x value exclusion applies only if the Digamma Function is used. If the Digamma Function is not used, any term of the integral with a division by zero will be excluded.<br>$k$ = constant of integration, a function of $\Delta x$ |
| 5d | The Digamma Function Discrete Derivative<br><br>$$D_1 \psi(x) = \frac{1}{x}$$<br>where<br>$x$ = real or complex variable<br>$n$ = real or complex constant |

| | **Zeta Function, Digamma Function, Polygamma Function, and lnd(n,Δx,x) Function Equations** |
|---|---|
| 5e | The Digamma Function term to term relationship<br><br>$\psi(x+1) - \psi(x) = \frac{1}{x}$<br>where<br>x = real or complex variable |
| 5f | The Digamma Function Recursion Equation<br>$\psi(x+1) = \psi(x) + \frac{1}{x}$<br>where<br>x = real or complex variable |

# TABLE 17

## List of Functions Calculated by the lnd(n,Δx,x) Function
**(where n, Δx, x are real or complex values)**

# Primary Relationships

| # | Function Calculated | Means of Calculation | Comments |
|---|---|---|---|
| **1** | tan(ax) | $\dfrac{\pi}{\text{lnd}(1,1,1-\frac{ax}{\pi}) - \text{lnd}(1,1,\frac{ax}{\pi})}$<br>for ax = nπ□ n=integer, tan(ax)=0 | For tan(x),<br>a = 1 |
| **2** | sin(x) | $\dfrac{2\tan(\frac{x}{2})}{1+\tan^2(\frac{x}{2})}$ | Use #1, a = ½ |
| **2a** | csc(πx)<br>$= \dfrac{1}{\sin(\pi x)}$ | $2e^{\lim_{n\to 0}\frac{\text{lnd}(n,1,x)-\text{lnd}(0,1,x)+\text{lnd}(n,1,1-x)-\text{lnd}(0,1,1-x)}{n}}$<br><br><u>Comment</u> – Computer calculation of the above equation requires more and more computational accuracy as x approaches integer values, the x values at which the poles of csc(πx) are located. For x very nearly equal to an integer value, the following equation can be used with good accuracy:<br>$\csc\pi x \approx \pm\dfrac{1}{(x-n)\pi}$<br>where<br>x = value very nearly equal to an integer, n<br>+ for n even , – for n odd | |

| # | Function Calculated | Means of Calculation | Comments |
|---|---|---|---|
| **3** | $\cos(x)$ | $\dfrac{1-\tan^2(\frac{x}{2})}{1+\tan^2(\frac{x}{2})}$ | Use #1, a = ½ |
| **4** | $e^{jx}$ | 1) $\dfrac{1+j\tan(\frac{x}{2})}{1-j\tan(\frac{x}{2})}$<br>2) $-jx[\lim_{\Delta x\to 0}\text{lnd}(1-jx,\Delta x,e)]$ | Use #1, a = ½<br>For small Δx, may get computer overflow |
| **5** | $\tanh(ax)$ | $\dfrac{-j\pi}{\text{lnd}(1,1,1-\frac{bx}{\pi}) - \text{lnd}(1,1,\frac{bx}{\pi})}$<br>b = ja<br>for ax = 0, tanh(ax)=0 | For tanh(x), a = 1 |
| **6** | $\sinh(x)$ | $\dfrac{2\tanh(\frac{x}{2})}{1-\tanh^2(\frac{x}{2})}$ | Use #5, a = ½ |
| **7** | $\cosh(x)$ | $\dfrac{1+\tanh^2(\frac{x}{2})}{1-\tanh^2(\frac{x}{2})}$ | Use #5, a = ½ |
| **8** | $e^{x}$ | 1) $\dfrac{1+\tanh(\frac{x}{2})}{1-\tanh(\frac{x}{2})}$<br>2) $-x[\lim_{\Delta x\to 0}\text{lnd}(1-x,\Delta x,e)]$ | Use #5, a = ½<br>For small Δx, may get computer overflow |

| # | Function Calculated | Means of Calculation | Comments |
|---|---|---|---|
| 9 | $\tan_{\Delta x}(a,x)$ <br> $= \tan[\tan^{-1}(a\Delta x)]\frac{x}{\Delta x}$ | 1) $\dfrac{\pi}{\mathrm{lnd}(1,1,1-\frac{bx}{\pi}) - \mathrm{lnd}(1,1,\frac{bx}{\pi})}$ <br> $b = \dfrac{\tan^{-1}(a\Delta x)}{\Delta x}$ <br> for $bx = n\pi$, n=integer, $\tan_{\Delta x}(a,x)=0$ <br> 2) $\tan_{\Delta x}(a,x) = \dfrac{2\tan_{\Delta x}(a, \frac{x}{2})}{1-\tan_{\Delta x}^{2}(a, \frac{x}{2})}$ | |
| 10 | $\sin_{\Delta x}(a,x)$ | $[1+(a\Delta x)^2]^{\frac{x}{2\Delta x}} \dfrac{2\tan_{\Delta x}(a, \frac{x}{2})}{1+\tan_{\Delta x}^{2}(a, \frac{x}{2})}$ | Use # 9, $x \to \frac{x}{2}$ |
| 11 | $\cos_{\Delta x}(a,x)$ | $[1+(a\Delta x)^2]^{\frac{x}{2\Delta x}} \dfrac{1-\tan_{\Delta x}^{2}(a, \frac{x}{2})}{1+\tan_{\Delta x}^{2}(a, \frac{x}{2})}$ | Use # 9, $x \to \frac{x}{2}$ |
| 12 | $e_{\Delta x}(ja,x)$ | $\cos_{\Delta x}(a,x) + j\sin_{\Delta x}(a,x)$ | Use # 10, #11 |

| # | Function Calculated | Means of Calculation | Comments |
|---|---|---|---|
| **13** | $\tanh_{\Delta x}(a,x)$<br>$= \tanh[\tanh^{-1}(a\Delta x)]^{\frac{x}{\Delta x}}$ | 1) $\dfrac{-j\pi}{\mathrm{lnd}(1,1,1-\frac{bx}{\pi}) - \mathrm{lnd}(1,1,\frac{bx}{\pi})}$<br>$b = \dfrac{\tan^{-1}(ja\Delta x)}{\Delta x} = \dfrac{j\tanh^{-1}(a\Delta x)}{\Delta x}$<br>for $bx = 0$, $\tanh_{\Delta x}(a,x)=0$<br>2) $\tanh_{\Delta x}(a,x) = \dfrac{2\tanh_{\Delta x}(a, \frac{x}{2})}{1+ \tanh_{\Delta x}^{2}(a, \frac{x}{2})}$ | |
| **14** | $\sinh_{\Delta x}(a,x)$ | $[1-(a\Delta x)^2]^{\frac{x}{2\Delta x}} \dfrac{2\tanh_{\Delta x}(a, \frac{x}{2})}{1- \tanh_{\Delta x}{}^{2}(a, \frac{x}{2})}$ | Use #13, $x \to \frac{x}{2}$ |
| **15** | $\cosh_{\Delta x}(a,x)$ | $[1-(a\Delta x)^2]^{\frac{x}{2\Delta x}} \dfrac{1+\tanh_{\Delta x}{}^{2}(a, \frac{x}{2})}{1- \tanh_{\Delta x}{}^{2}(a, \frac{x}{2})}$ | Use #13, $x \to \frac{x}{2}$ |
| **16** | $e_{\Delta x}(a,x)$ | $\cosh_{\Delta x}(a,x) + \sinh_{\Delta x}(a,x)$ | Use #14, #15 |
| **17** | $\ln(x)$ | 1a) $\lim_{n\to\infty} [\,\mathrm{lnd}(1,\frac{x-a}{n},x) - \mathrm{lnd}(1,\frac{x-a}{n}, a)\,] + \ln(a)$<br>1b) $\lim_{n\to\infty} [\,\mathrm{lnd}(1,1, \frac{nx}{x-a}) - \mathrm{lnd}(1,1, \frac{na}{x-a})\,] + \ln(a)$<br>1c) $\lim_{n\to\infty} [\,\psi(\frac{nx}{x-a}) - \psi(\frac{na}{x-a})\,] + \ln(a)$ | $\ln_{\Delta x}x \equiv \mathrm{lnd}(1,\Delta x,x)$ |

| # | Function Calculated | Means of Calculation | Comments |
|---|---|---|---|
| | | where<br>$a = \begin{cases} +1 & \text{for } x = \text{positive real value or} \\ & \text{for } x = \text{non-real complex value} \\ -1 & \text{for } x = \text{negative real value or} \\ & \text{for } x = \text{non-real complex value} \end{cases}$<br>x = real or complex value<br>$x \neq -1,0,1$<br>n = positive integers, $n \rightarrow \infty$<br>$\psi(x)$ = Digamma Function, $x \neq 0,-1,-2,-3, \ldots$<br>$\ln(-1) = \pm\pi i$<br>ln(0) is infinite<br>$\ln(1) = 0$<br><br>or<br><br>2) $\lim_{n\rightarrow 0} \left[ \frac{-\text{lnd}(-n,1,x+1) + \text{lnd}(-n,1,x) - 1}{n} \right]$<br>n = real number, $n \neq 0$<br>$x \neq 0$<br>or<br><br>3) $\lim_{n\rightarrow 0} [\text{lnd}(n,1,x+1) - \text{lnd}(0,1,x+1) - \text{lnd}(n,1,x) + \text{lnd}(0,1,x)]$<br><br>n = real number, $n \neq 0$<br>$x \neq 0$<br>or<br><br>4) $\lim_{\Delta x\rightarrow 0}[\text{lnd}(1,\Delta x,x) - \text{lnd}(1,\Delta x,1)]$ | |

| # | Function Calculated | Means of Calculation | Comments |
|---|---|---|---|
| **18** | $x^N$ | 1) $\alpha(N)[\,\mathrm{lnd}(-N,1,x+1) - \mathrm{lnd}(-N,1,x)\,]$<br>$\alpha(N) = +1$ for $N = -1$, $\alpha(N) = -1$ for $N \neq -1$<br>2) $-N[\lim_{\Delta x \to 0} \mathrm{lnd}(1-N,\Delta x,x)]$ | For small $\Delta x$, may get computer overflow |
| **19** | $N^x$ | $\dfrac{1+\tanh[\frac{\ln(N)}{2}]x}{1-\tanh[\frac{\ln(N)}{2}]x}$ | Use #5, #17 $a = \dfrac{\ln(N)}{2}$ |
| **20** | 1 $\sum_{x=x_1}^{x_2} \ln(1+x)$ | 1) $\left[x(\ln x-1) + \dfrac{\mathrm{lnd}(1,1,x)}{2} + \sum_{n=2}^{\infty}(-1)^n \dfrac{\mathrm{lnd}(n,1,x)}{n+1}\,\right]\Big\|_{x_1}^{x_2+1}$<br>2) $\left[x(\ln x-1) + \dfrac{\mathrm{lnd}(1,1,x+\frac{1}{2})}{2} + \sum_{n=2}^{\infty}\left\{\dfrac{\mathrm{lnd}(2n-2,1,x+\frac{1}{2})}{(2n-1)(2)^{2n-2}} - \dfrac{\mathrm{lnd}(2n-1,1,x+\frac{1}{2})}{(2n-1)(2)^{2n-1}}\right\}\right]\Big\|_{x_1}^{x_2+1}$<br>$\mathrm{Re}(1/x) \geq -1$, $1/x \neq -1$ | Use #17<br>(1) has slower convergence than (2)<br>Re(x) is the real part of x |
| **21** | $\ln\Gamma(x)$ | 1) $\ln\Gamma(x)$<br>where<br>$\Gamma(x) = \dfrac{(2\pi)^x}{2\cos(\frac{\pi x}{2})}\left[\dfrac{\mathrm{lnd}(1-x,1,1)}{\mathrm{lnd}(x,1,1)}\right]$, $x \neq 0,1,3,5,7,9,\ldots$<br>$\Gamma(x)$ is infinite for $x = 0,-1,-2,-3,\ldots$<br>$\Gamma(1) = 1$<br>First use $\Gamma(x) = (x-1)\,\Gamma(x-1)$ if $x = 3,5,7,9,\ldots$ | Use #3, #17, #19<br>For small \|x\| use 1) or 2) |

| # | Function Calculated | Means of Calculation | Comments |
| --- | --- | --- | --- |
| | | or<br><br>2) $\lim_{n\to 0} \frac{\text{lnd}(n,1,x) - \text{lnd}(0,1,x)}{n} + \frac{\ln 2\pi}{2}$<br>n = real number, $n\to 0$<br><br>or<br><br>3) $\eta + (x\text{-}1)[\ln(x\text{-}1)\text{-}1] + \frac{\text{lnd}(1,1,x\text{-}1)}{2} + \sum_{n=2}^{\infty}(-1)^n \frac{\text{lnd}(n,1,x\text{-}1)}{n+1}$<br>$\|x\text{-}1\|>1$, $\eta = 0.6303307007539063114770736\,91$<br><br>or<br><br>4) $\eta + (x\text{-}1)[\ln(x\text{-}1)\text{-}1] + \frac{\text{lnd}(1,1,x\text{-}\frac{1}{2})}{2} + \sum_{n=2}^{\infty}\left[\frac{\text{lnd}(2n\text{-}2,1,x\text{-}\frac{1}{2})}{(2n\text{-}1)(2)^{2n-2}} - \frac{\text{lnd}(2n\text{-}1,1,x\text{-}\frac{1}{2})}{(2n\text{-}1)(2)^{2n-1}}\right]$<br>$\text{re}(1/x) \geq -1$, $1/x \neq -1$, $\eta = 0.630330700753906311477073691364\ldots$ | (4) has faster convergence than (3) |

| # | Function Calculated | Means of Calculation | Comments |
|---|---|---|---|
| **22** | $\Gamma(x)$ | 1) $\frac{(2\pi)^x}{2\cos(\frac{\pi x}{2})}\left[\frac{\text{lnd}(1-x,1,1)}{\text{lnd}(x,1,1)}\right]$, $x \neq 0,1,3,5,7,9,\ldots$<br>where<br>$\Gamma(x)$ is infinite for $x = 0,-1,-2,-3,\ldots$<br>$\Gamma(1) = 1$<br>First use $\Gamma(x) = (x-1)\,\Gamma(x-1)$ if $x = 3,5,7,9,\ldots$<br>or<br>2) $\sqrt{2\pi}\, e^{\lim_{n\to 0}\frac{\text{lnd}(n,1,x)-\text{lnd}(0,1,x)}{n}}$<br>n = real number<br>or<br>3) $e^{\eta+(x-1)[\ln(x-1)-1]+\frac{\text{lnd}(1,1,x-1)}{2}+\sum_{n=2}^{\infty}(-1)^n\frac{\text{lnd}(n,1,x-1)}{n+1}}$<br>$\lvert x-1\rvert>1$, $\eta = 0.6303307007539063114770736913664$<br>or<br>4) $e^{\eta+(x-1)[\ln(x-1)-1]+\frac{\text{lnd}(1,1,x-\frac{1}{2})}{2}+\sum_{n=2}^{\infty}\left[\frac{\text{lnd}(2n-2,1,x-\frac{1}{2})}{(2n-1)(2)^{2n-2}}-\frac{\text{lnd}(2n-1,1,x-\frac{1}{2})}{(2n-1)(2)^{2n-1}}\right]}$<br>$1/x \geq -1$, $1/x \neq -1$, $\eta = 0.6303307007539063114770736913664\ldots$ | Gamma Function<br><br>Use #3, #17, #21<br><br>For small \|x\| use 1) or 2)<br><br>(4) has faster convergence than (3) |

| # | Function Calculated | Means of Calculation | Comments |
|---|---|---|---|
| **23** | $\Gamma(x)\Gamma(1-x)$ | $\dfrac{\pi}{\sin(\pi x)}$ | Use #2 |
| **24** | $\dfrac{d\Gamma(x)}{dx}$ | $\psi(x)\Gamma(x)$ | Use #22, #25 |
| **25** | $\psi(x)$ | 1) $\text{lnd}(1,1,x) - \gamma\ ,\ \ x \neq 0,-1,-2,-3,\ldots$<br>2) $-\lim_{n\to 0} [\text{lnd}(1+n,1,x) - \frac{1}{n}]$<br>3) $\frac{d}{dx} \ln\Gamma(x)$<br>4) $\dfrac{\Gamma(x)'}{\Gamma(x)}$<br>where<br>$\Gamma(x)' = \dfrac{d\Gamma(x)}{dx}$<br>$\gamma$ = Euler's constant, .5772157…<br>$\psi(x)$ is infinite for $x = 0,-1,-2,-3,\ldots$<br>Some Digamma Function equations are:<br>$D_1\psi(x) = \frac{1}{x}$ , Discrete derivative<br>$\psi(x+1) = \psi(x) + \frac{1}{x}$ , Recursion equation<br>$\psi(1-x) = \psi(x) + \pi\,\text{cotan}(\pi x)$ , Reflection equation | Digamma Function |

| # | Function Calculated | Means of Calculation | Comments |
|---|---|---|---|
| **26** | $\psi^{(n)}(x)$ | 1) $(-1)^{n+1}\, n!\, \mathrm{lnd}(n+1,1,x)$ <br> 2) $(-1)^{n+1}\, n!\, \zeta(n+1,x)$ <br> 3) $\psi^{(n)}(x) = \frac{d^n}{dx^n}\psi(x) = \frac{d^{n+1}}{dx^{n+1}}\ln\Gamma(x)$ <br> where n = 1, 2, 3, 4, … | Polygamma Functions |
| **27** | $\eta(x)$ | $(1-2^{1-x})\, \mathrm{lnd}(x,1,1)$ <br> $n\neq 1$ | Dirichlet Eta Function |
| **28** | $E_{2n}$ | $2^{4n+2}\mathrm{lnd}(-2n,1,\frac{1}{4})$ <br> where n = 0, 1, 2, 3,… | Euler Numbers |
| **29** | $B_n$ | $(-1)^{n+1}\, n\, \mathrm{lnd}(1-n,1,1)$ <br> where n = 1, 2, 3, 4, … | Bernoulli Constants |
| **30** | $\zeta(n,\Delta x,x)$ | 1) $\frac{1}{\Delta x}\mathrm{lnd}(n,\Delta x,x)$ , n = any value <br> 2) $\frac{1}{(\Delta x)^n}\mathrm{lnd}(n,1,\frac{x}{\Delta x})$ , x,Δx = real values with the condition that <br> x ≠ positive real value when <br> Δx = negative real value <br> n = complex value <br> or <br> x = Δx = complex values <br> n = complex value <br> or <br> x,Δx = complex values <br> n = integer | General Zeta Function |

| # | Function Calculated | Means of Calculation | Comments |
|---|---|---|---|
| | | 3) $\frac{1}{(\Delta x)^{n-1}}\zeta(n,1,\frac{x}{\Delta x})$ , n,Δx,x conditions are the same as listed in 2)<br><br>4) $\frac{1}{(\Delta x)^{n-1}}\zeta(n,\frac{x}{\Delta x})$ , n,Δx,x conditions are the same as listed in 2)<br><br>$\zeta(n,\Delta x)$ = Hurwitz Zeta Function , $n \neq 1$<br>If evaluating a summation, any summation term with x = 0 is excluded<br><br><u>Note</u> – A complex value, x+jy, can be a real number, x, an imaginary value, jy, or an integer, N+j0.<br><br><u>Comments</u> – For value combinations not specified above, the equality of the stated equation may not be valid.<br><br>The General Zeta Function, unlike the Riemann and Hurwitz Zeta Functions, is defined for n=1. Note the N=1 Zeta Function in row 30c, the Digamma Function in row 25, and the $\frac{1}{x}$ function summation in row 31.<br><br>The General Zeta Function and the Hurwitz Zeta Function are closely related. With the exception of the case where n=1, $\zeta(n,x) = \zeta(n,1,x)$.<br><br>Some General Zeta Function equations are:<br>$\zeta(n,\Delta x,x) = \frac{1}{\Delta x}\text{lnd}(n,\Delta x,x)$ , General Zeta Function<br>$\zeta(n,x) = \zeta(n,1,x) = \text{lnd}(n,1,x)$ , Hurwitz Zeta Function<br>$\zeta(n) = \zeta(n,1,1) = \text{lnd}(n,1,1)$ , Riemann Zeta Function | |

| # | Function Calculated | Means of Calculation | Comments |
| --- | --- | --- | --- |
| | | $\psi(x) = \zeta(1,1,x) - \gamma$ , $x \neq 0,-1,-2,\ldots$ Digamma Function<br>$D_{\Delta x}\,\zeta(n,\Delta x,x) = \pm \frac{1}{\Delta x}\frac{1}{x^n}$ , Discrete derivative<br>$\frac{d}{dx}\,\zeta(n,\Delta x,x) = -n\,\zeta(n+1,\Delta x,x)$ , $n\neq 0$ Derivative<br>$\zeta(n,\Delta x,x+\Delta x) = \zeta(n,\Delta x,x) \pm \frac{1}{x^n}$ Recursion Equation<br>where + for $n = 1$, – for $n \neq 1$ | |
| **30a** | $\zeta(n)$ | 1) $\text{lnd}(n,1,1)$<br>2) $\frac{1}{2(1-2^{1-n})}\,[\text{lnd}(n,2,1) - \text{lnd}(n,2,2)]$<br>3) $(\Delta x)^{n-1}\text{lnd}(n,\Delta x,\Delta x)$ , n = integers<br>$n \neq 1$ | Riemann Zeta Function<br>$\zeta(n) = \zeta(n,1,1)$ |
| **30b** | $\zeta(n,x)$ | $\text{lnd}(n,1,x)$<br>$n \neq 1$<br>If evaluating a summation, any summation term with x = 0 is excluded<br>Some Hurwitz Zeta Function equations are:<br>$D_1\zeta(n,x) = -\frac{1}{x^n}$ , $n \neq 1$ Discrete derivative<br>$\frac{d}{dx}\,\zeta(n,x) = -n\,\zeta(n+1,x)$ , $n\neq 0,1$ Derivative<br>$\zeta(n,x+1) = \zeta(n,x) - \frac{1}{x^n}$ , $n \neq 1$ Recursion Equation | Hurwitz Zeta Function<br>$\zeta(n,x) = \zeta(n,1,x)$ |

| # | Function Calculated | Means of Calculation | Comments |
|---|---|---|---|
| **30c** | $\zeta(1,\Delta x,x)$ | 1) $\frac{1}{\Delta x}\,\text{lnd}(1,\Delta x,x)$<br><br>2) $\frac{1}{\Delta x}\,\text{lnd}(1,1,\frac{x}{\Delta x})$<br><br>n = 1<br>If evaluating a summation, any summation term with x = 0 is excluded<br><br><u>Comments</u> – The N=1 Zeta Function Function is closely related to the Digamma Function.<br>$\zeta(1,\Delta x,x) = \frac{1}{\Delta x}[\psi(\frac{x}{\Delta x}) + \gamma]\ ,\ \ \frac{x}{\Delta x} \neq 0,-1,-2,-3,\ldots$<br>$\zeta(1,1,x) = \psi(x) + \gamma\ ,\ \ x \neq 0,-1,-2,-3,\ldots$ | N=1 Zeta Function |
| **31** | $\sum_{\Delta x\ x=x_1}^{x_2} \frac{1}{x}$ | 1) $\frac{1}{\Delta x}\,\text{lnd}(1,\Delta x,x)\Big\vert_{x_1}^{x_2+\Delta x}$<br>or<br>2) $\frac{1}{\Delta x}\,\text{lnd}(1,1,\frac{x}{\Delta x})\Big\vert_{x_1}^{x_2+\Delta x}$<br><br>Any summation term with x = 0 is excluded<br><u>Comments</u> –<br>$\text{lnd}(1,1,x) = \psi(x) + \gamma\ ,\ \ x \neq 0,-1,-2,-3,\ldots$<br>$\psi(x)$ <u>is infinite</u> for $x = 0,-1,-2,-3,\ldots$<br>$\text{lnd}(1,1,x)$ <u>is not infinite</u> for $x = 0,-1,-2,-3,\ldots$<br>$\text{lnd}(1,1,x) = \text{lnd}(1,1,1-x)$ , $x$ = integer or integer + .5<br>$\text{lnd}(1,1,1-x) - \text{lnd}(1,1,x) = \begin{cases} 0\ ,\ x = \text{integer} \\ \pi\cot(\pi x)\ ,\ x \neq \text{integer} \end{cases}$<br><br>See row 31a | |

| # | Function Calculated | Means of Calculation | Comments |
|---|---|---|---|
| **31a** | $\sum_{\Delta x \; x=x_1}^{x_2} \frac{1}{x}$ | $\frac{1}{\Delta x}\psi(\frac{x}{\Delta x})\Big\vert_{x_1}^{x_2+\Delta x} = \zeta(1,1,\frac{x}{\Delta x})\Big\vert_{x_1}^{x_2+\Delta x} = \frac{1}{\Delta x}\text{lnd}(1,1,\frac{x}{\Delta x})\Big\vert_{x_1}^{x_2+\Delta x} = \frac{1}{\Delta x}\text{lnd}(1,\Delta x,x)\Big\vert_{x_1}^{x_2+\Delta x}$<br>where<br>x = real or complex variable<br>$x_1,x_2,\Delta x$ = real or complex constants<br>$\gamma$ = Euler's Constant, .5772157…<br>$\Delta x$ = x increment<br>$x_1 \neq 0,-\Delta x,-2\Delta x,-3\Delta x,\ldots$ , These exclusions apply only if the Digamma Function is used.<br>$x_2 \neq -\Delta x,-2\Delta x,-3\Delta x,\ldots$ , If the Digamma Function is not used, any term of the integral with a division by zero will be excluded.<br>Comments - $\psi(x)$ is infinite for $x = 0,-1,-2,-3,\ldots$<br>$\psi(x) = \text{lnd}(1,1,x) - \gamma$ , $x \neq 0,-1,-2,-3,\ldots$<br>$\psi(1-x) - \psi(x) = \pi\text{cotan}(\pi x)$<br>If $x_1$ or $x_2$ is equal to any of the above disallowed values, use the equations of row 31 which allow these values.<br>See row 31 | Digamma Function summation between finite limits |
| **31b** | $SHP_m = \sum_{\Delta x \; x=x_1}^{x_1+(m-1)\Delta x} \frac{1}{x}$ | $\frac{1}{\Delta x}\text{lnd}(1,\Delta x,x)\Big\vert_{x_1}^{x_1+m\Delta x}$<br>where<br>m = the number of terms in the harmonic progression<br>$\Delta x$ = common difference<br>$x_m = x_1 + (m-1)\Delta x$<br>$\frac{1}{x_1}$ = first term of the harmonic progression summation<br>$\frac{1}{x_m}$ = last term of the harmonic progression summation<br>$SHP_m$ = the harmonic progression sum | Sum of m terms of a Harmonic Progression |

| # | Function Calculated | Means of Calculation | Comments |
|---|---|---|---|
| **31c** | $_{\Delta x}\int_{x_1}^{x_2} \frac{1}{x+a}\,\Delta x$<br>or<br>$\Delta x \,_{\Delta x}\sum_{x=x_1}^{x_2-\Delta x} \frac{1}{x+a}$ | 1) $\text{lnd}(1,\Delta x,x+a)\Big\|_{x_1}^{x_2}$<br>or<br>2) $\text{lnd}(1,1,\frac{x+a}{\Delta x})\Big\|_{x_1}^{x_2}$<br><br>Any summation term with x+a = 0 is excluded<br><br>or<br>3) $\psi(\frac{x+a}{\Delta x})\Big\|_{x_1}^{x_2}$, $\frac{x+a}{\Delta x} \neq 0,-1,-2,-3,\ldots$<br><br>where $\psi(x)$ = Digamma Function | Discrete definite integral of the function, $\frac{1}{x+a}$<br><br>Eq 1) and Eq 2 can be used for all real and complex values of x,a,Δx |
| **31d** | $_{\Delta x}\int \frac{1}{x+a}\,\Delta x$ | 1) $\text{lnd}(1,\Delta x,x+a) + k(\Delta x)$<br>or<br>2) $\text{lnd}(1,1,\frac{x+a}{\Delta x}) + k(\Delta x)$<br>or<br>3) $\psi(\frac{x+a}{\Delta x}) + k(\Delta x)$, $\frac{x+a}{\Delta x} \neq 0,-1,-2,-3,\ldots$<br><br>where<br>$\psi(x)$ = Digamma Function<br>k = constant of integration, a function of Δx | Discrete indefinite integral of the function, $\frac{1}{x+a}$ |
| **32** | $_{\Delta x}\sum_{x=x_1}^{x_2} \frac{1}{x^n}$ | $\pm\,\zeta(n,\Delta x,x)\Big\|_{x_1}^{x_2+\Delta x} = \pm\,\frac{1}{\Delta x}\,\text{lnd}(n,\Delta x,x)\Big\|_{x_1}^{x_2+\Delta x}$, $-$ for $n \neq 1$, $+$ for $n = 1$<br><br>Any summation term where x = 0 is excluded | General Zeta Function summation between finite limits |

| # | Function Calculated | Means of Calculation | Comments |
|---|---|---|---|
| **32a** | $\sum_{\Delta x}{}_{x=x_1}^{x_2} \frac{1}{x^n}$ | $-\frac{1}{\Delta x}\,\mathrm{lnd}(n,\Delta x,x)\Big\vert_{x_1}^{x_2+\Delta x}\ ,\ n \neq 1$<br>Any summation term with x = 0 is excluded | For n = 1 use #31 |
| **32b** | $\sum_{1}{}_{x=x_1}^{x_2} \frac{1}{x^n}$ | $-\zeta(n,x)\Big\vert_{x_1}^{x_2+1} = -\zeta(n,1,x)\Big\vert_{x_1}^{x_2+1} = -\mathrm{lnd}(n,1,x)\Big\vert_{x_1}^{x_2+1}\ ,\ n \neq 1$<br>Any summation term where x = 0 is excluded | Hurwitz Zeta Function summation between finite limits |
| **32c** | $\sum_{\Delta x}{}_{x=x_1}^{x_2} \frac{1}{x^{1+m}}$ | $\frac{(-1)^m}{m!\Delta x^{m+1}}\,\psi^{(m)}(\frac{x}{\Delta x})\Big\vert_{x_1}^{x_2+\Delta x} = -\zeta(1+m,\Delta x,x)\Big\vert_{x_1}^{x_2+\Delta x} = -\frac{1}{\Delta x}\,\mathrm{lnd}(1+m,\Delta x,x)\Big\vert_{x_1}^{x_2+\Delta x}$<br>where<br>x = real or complex variable<br>$x_1, x_2, \Delta x$ = real or complex constants<br>$\Delta x$ = x increment<br>m = 1,2,3,…<br>Any summation term where x = 0 is excluded | Polygamma Functions summation between finite limits |
| **32d** | $SAP_m = \sum_{\Delta x}{}_{x=x_1}^{x_1+(m-1)\Delta x} x$ | $-\frac{1}{\Delta x}\,\mathrm{lnd}(-1,\Delta x,x\ )\Big\vert_{x_1}^{x_1+m\Delta x}$<br>where<br>m = the number of terms in the harmonic progression<br>$\Delta x$ = common difference<br>$x_m = x_1 + (m-1)\Delta x$<br>$x_1$ = first term of the arithmetic progression summation<br>$x_m$ = last term of the arithmetic progression summati<br>$SHP_m$ = the arithmetic progression sum | Sum of m terms of an Arithmetic Progression |

| # | Function Calculated | Means of Calculation | Comments |
|---|---|---|---|
| **32e** | $\int_{\Delta x\, x_1}^{x_2} \frac{1}{(x+a)^n} \Delta x$<br>or<br>$\Delta x \sum_{\Delta x\, x=x_1}^{x_2-\Delta x} \frac{1}{(x+a)^n}$ | 1) $-\,\mathrm{lnd}(n,\Delta x,x+a)\Big\vert_{x_1}^{x_2}$, $n \neq 1$<br>Any summation term with x+a = 0 is excluded<br>or<br>2) $-\frac{1}{(\Delta x)^{n-1}}\,\zeta(n,\frac{x+a}{\Delta x})\Big\vert_{x_1}^{x_2}$, $n \neq 1$<br>where<br>x+a,Δx = real values with the condition that<br>x+a ≠ positive real value when Δx = negative real value<br>n = complex value<br>or<br>x+a = Δx = complex values<br>n = complex value<br>or<br>x+a,Δx = complex values<br>n = integer<br><br>$\zeta(n,x)$ = Hurwitz Zeta Function , $n \neq 1$<br><br>Any summation term with x+a = 0 is excluded<br><br><u>Note</u> – A complex value, x+jy, can be a real number, x, an imaginary value, jy, or an integer, N+j0. | Discrete definite Integral of the function, $\frac{1}{(x+a)^n}$<br><br>For n = 1 use #31c |
| **32f** | $\int_{\Delta x} \frac{1}{(x+a)^n} \Delta x$ | $-\,\mathrm{lnd}(n,\Delta x,x+a) + k(n,\Delta x)$, $n \neq 1$<br>where<br>k = constant of integration, a function of n,Δx | Discrete indefinite Integral of the function, $\frac{1}{(x+a)^n}$<br><br>For n = 1 use #31d |

| # | Function Calculated | Means of Calculation | Comments |
|---|---|---|---|
| **33** | $\sum_{\Delta x\, x=x_1}^{x_2} (-1)^{\frac{x-x_1}{\Delta x}} \frac{1}{x^n}$ | $\frac{\alpha(n)}{2\Delta x}[-\mathrm{lnd}(n,2\Delta x,x) \Big\vert_{x_1}^{x_2+\Delta x} + \mathrm{lnd}(n,2\Delta x,x) \Big\vert_{x_1+\Delta x}^{x_2+2\Delta x}]$<br>$x = x_1, x_1+\Delta x, x_1+2\Delta x, x_1+3\Delta x, \ldots, x_2-\Delta x, x_2$<br>$x_2 = x_1 + (2m-1)\Delta x$ , $m = 1,2,3,\ldots$<br>$\Delta x$ = x increment<br>n = any value<br>$\alpha(n) = \begin{cases} 1 & n \neq 1 \\ -1 & n = 1 \end{cases}$<br>$x \neq 0$ | |
| **34** | $\sum_{\Delta x\, x=x_1}^{x_2} f(x)$<br>or<br>$\frac{1}{\Delta x} \int_{\Delta x\, x_1}^{x_2+\Delta x} f(x)\Delta x$ | 1) $\frac{1}{\Delta x} \sum_{r=1}^{R} \alpha_r(n_r) a_r \mathrm{lnd}(n_r,\Delta x,x+c_r) \Big\vert_{x_1}^{x_2+\Delta x}$<br>for $f(x) = \sum_{r=1}^{R} \frac{a_r}{(x+c_r)^{n_r}}$<br>R = the number of f(x) terms<br>$a_r$, $c_r$ = constants<br>$\alpha_r(n_r) = +1$ for $n_r = 1$<br>$\alpha_r(n_r) = -1$ for $n_r \neq 1$<br>Any summation term with $x+c_r = 0$ is excluded<br>or | Partial Fraction Expansion<br><br>Sum of a series of $\frac{a_r}{(x+c_r)^{n_r}}$ terms |

| # | Function Calculated | Means of Calculation | Comments |
|---|---|---|---|
| | | 2) $-\frac{1}{\Delta x}\sum_{x=0}^{\infty} f(x_1)^{(n)} \frac{\text{lnd}(-n,\Delta x,x-x_1)}{n!} \Big\vert_{x_1}^{x_2+\Delta x}$<br>where<br>$x = x_1, x_1+\Delta x, x_1+2\Delta x, x_1+3\Delta x), \ldots , x_2-\Delta x, x_2$<br>f(x) is a continuous function of x<br>f(x) can be expressed as a convergent Taylor Series about the point $x_1$<br>The series region of convergence contains the points $x_1$ thru $x_2+\Delta x$<br>$f(x_1)^{(n)}$ is the nth derivative of the function, f(x), evaluated at the point $x_1$<br><br>Comment – There are more ways to evaluate this summation using equations which do not involve the lnd(n,Δx,x) function. The following discrete calculus equation can be used,<br>$\sum_{\Delta x\, x=x_1}^{x_2} f(x) = \frac{1}{\Delta x}\sum_{n=0}^{\infty} D_{\Delta x}^{\ n} f(x_1) \frac{[x-x_1]_{\Delta x}^{n+1}}{(n+1)n!} \Big\vert_{x_1}^{x_2+\Delta x}$ ,<br>or several other equations listed in Table 7. | |
| **35** | $\sum_{\Delta x\, x=x_i}^{\pm\infty} \frac{1}{x^n}$ | $\zeta(n,\Delta x,x_i) = \frac{1}{\Delta x}\text{lnd}(n,\Delta x,x_i)$<br><br>+∞ for Re(Δx)>0 or [ Re(Δx)=0 and Im(Δx)>0 ]<br>-∞ for Re(Δx)<0 or [ Re(Δx)=0 and Im(Δx)<0 ]<br>$x = x_i + m\Delta x,\ m = 0,1,2,3,\ldots$<br>Re(n)>1<br><br>Any summation term with x = 0 is excluded | General Zeta Function summation to infinity |

| # | Function Calculated | Means of Calculation | Comments |
|---|---|---|---|
| **35a** | $\sum_{x=x_i}^{\infty} {}_1 \frac{1}{x^n}$ | $\zeta(n,x_i) = \zeta(n,1,x_i) = \text{lnd}(n,1,x_i)$<br><br>$x = x_i + m,\ \ m = 0,1,2,3,\ldots$<br>$\text{Re}(n)>1$<br>Any summation term with $x = 0$ is excluded | Hurwitz Zeta Function summation to infinity |
| **35b** | $\sum_{x=1}^{\infty} {}_1 \frac{1}{x^n}$ | $\zeta(x) = \zeta(x,1,1) = \ \ \text{lnd}(n,1,1)$<br><br>$x = 1,2,3,\ldots$<br>$\text{Re}(n)>1$ | Riemann Zeta Function summation to infinity |
| **36** | $\sum_{x=x_i}^{\pm\infty} {}_{\Delta x} (-1)^{\frac{x-x_i}{\Delta x}} \frac{1}{x^n}$ | 1) $\frac{1}{\Delta x}[\ -\text{lnd}(n,2\Delta x,x_i+\Delta x) + \text{lnd}(n,\Delta x,x_i)\ ]$<br>where $n\neq1$<br>or<br>2) $\frac{1}{\Delta x}[\ \text{lnd}(n,\Delta x,x_i) - 2^{1-n}\,\text{lnd}(n,\Delta x,\frac{x_i+\Delta x}{2})\ ]$<br>where $n\neq1$<br><br>or<br><br>3) $\frac{1}{\Delta x}[\ 2^{1-n}\,\text{lnd}(n,\Delta x,\frac{x_i}{2}) - \text{lnd}(n,\Delta x,x_i)\ ]$<br>where $n\neq1$<br>or<br>4) $\pm2^{-n}\,D_{\Delta x}\text{lnd}(n,\Delta x,\frac{x_i}{2})$, + for $n=1$, – for $n\neq1$<br>where $n=1$ is allowed<br>or | Alternating sum $x_i$ to $+\infty$ or $x_i$ to $-\infty$ |

| # | Function Calculated | Means of Calculation | Comments |
|---|---|---|---|
| | | 5) $\frac{\alpha(n)}{2^n\Delta x}$ [ $\text{lnd}(n,\Delta x,\frac{x_i}{2})$ – $\text{lnd}(n,\Delta x,\frac{x_i+\Delta x}{2})$ ]<br>where n =1 is allowed<br>or<br>6) $\frac{\alpha(n)}{2\Delta x}$ [ $-\text{lnd}(n,2\Delta x,x_i+\Delta x) + \text{lnd}(n,2\Delta x,x_i)$ ]<br>where n = 1 is allowed<br><br>Re(n)>0<br>$\alpha(n) = -1$ for n = 1<br>$\alpha(n) = +1$ for $n \neq 1$<br>$+\infty$ for Re(Δx)>0 or [ Re(Δx)=0 and Im(Δx)>0 ]<br>$-\infty$ for Re(Δx)<0 or [ Re(Δx)=0 and Im(Δx)<0 ]<br>$x = x_i + m\Delta x$, m = 0,1,2,3, …<br>$x\neq 0$<br><br><u>Comment</u> – From Eq 2) where Δx = 1 and $x_i$ = 1, the well known equation for the Riemann Zeta Function is derived.<br><br>$$\zeta(n) = \zeta(n,1,1) = \text{lnd}(n,1,1) = \frac{1}{1-2^{1-n}}\sum_{x=1}^{\infty}(-1)^{\frac{x-1}{\Delta x}}\frac{1}{x^n}$$<br>Re(n)>0 | |

| # | Function Calculated | Means of Calculation | Comments |
| --- | --- | --- | --- |
| **37** | $\sum_{\Delta x,\ x=\mp\infty}^{\pm\infty} \frac{1}{x^n}$ | $\frac{\alpha(n)}{\Delta x}$ [lnd(n,Δx,$x_i$) - lnd(n,-Δx,$x_i$-Δx)]<br><br>Re(n)≥1<br>α(n) = -1 for n = 1<br>α(n) = +1 for n ≠ 1<br>-∞ to +∞ for Re(Δx)>0 or [Re(Δx)=0 and Im(Δx)>0]<br>+∞ to -∞ for Re(Δx)<0 or [Re(Δx)=0 and Im(Δx)<0]<br>x = $x_i$+mΔx, m = 0,1,2,3,…<br>$x_i$ = value of x<br>n,x,$x_i$,Δx = real or complex values<br><br>Any summation term with x = 0 is excluded<br><br><u>Comment</u> – If the value of $x_i$ is changed to $x_i$ + rΔx, r = integers, the summation value remains the same. | Sum –∞ to +∞ or<br>Sum +∞ to –∞ |

| # | Function Calculated | Means of Calculation | Comments |
|---|---|---|---|
| **37a** | $\sum_{\Delta x\ x=\mp\infty}^{\pm\infty} \frac{1}{x}$ <br> $x = x_i + m\Delta x$ <br> m = integers <br> $x_i$ = value of x <br><br> or <br><br> $\sum_{\Delta x\ x=x_i}^{\pm\infty} \frac{1}{x} - \sum_{\Delta x\ x=\Delta x - x_i}^{\pm\infty} \frac{1}{x_i}$ | 0 <br> for $\frac{x_i}{\Delta x}$ = integer or integer + ½ <br><br> Any summation term with x = 0 is excluded <br> ----------------------------------- <br><br> $\frac{\pi}{\Delta x}\,\mathrm{cotan}(\frac{\pi x_i}{\Delta x})$ <br> for $\frac{x_i}{\Delta x} \neq$ integer <br><br> ----------------------------------- <br><br> $\frac{1}{\Delta x}$ [-lnd(1,Δx,$x_i$) + lnd(1,Δx,Δx-$x_i$)] <br> for $\frac{x_i}{\Delta x}$ = all values <br><br> -∞ to +∞ for Re(Δx)>0 or [Re(Δx)=0 and Im(Δx)>0] <br> +∞ to -∞ for Re(Δx)<0 or [Re(Δx)=0 and Im(Δx)<0] <br> x, $x_i$, Δx = real or complex values <br> Δx = x increment <br><br> Any summation term with x = 0 is excluded <br><br> <u>Comment</u> – If the value of $x_i$ is changed to $x_i$ + rΔx, r = integers, the summation value remains the same. | Sum –∞ to +∞ or <br> Sum +∞ to –∞ <br><br> or <br><br> Summations <br><br> Δx-$x_i$ and $x_i$ to +∞ or <br> Δx-$x_i$ and $x_i$ to -∞ |

| # | Function Calculated | Means of Calculation | Comments |
|---|---|---|---|
| 37b | $\sum_{x=\mp\infty}^{\pm\infty}{}_{\Delta x} \frac{1}{x^n}$<br>$x = x_i+m\Delta x$<br>m = integers<br>$x_i$ = value of x | 0<br>for $\frac{x_i}{\Delta x}$ = integer or integer + ½<br>n = 1, 3, 5, 7 …<br><br>Any summation term with x = 0 is excluded<br>------------------------------------------------------------<br><br>$\frac{2}{(\Delta x)^n}\,\text{lnd}(n,1,1) = \frac{2}{(\Delta x)^n}\,\zeta(n)$<br>for $\frac{x_i}{\Delta x}$ = integer<br>n = 2, 4, 6, 8, …<br><br>Any summation term with x = 0 is excluded<br>------------------------------------------------------------<br><br>$\frac{-\pi}{(n-1)!\Delta x}\frac{d^{n-1}}{dx^{n-1}}\,\text{cotan}(\frac{\pi x}{\Delta x})\,\Big\rvert_{x=-x_i}$<br>for $\frac{x_i}{\Delta x} \neq$ integer<br>n = 1, 2, 3, 4 …<br>------------------------------------------------------------<br>$\frac{1}{\Delta x}\,[(-1)^n\,\text{lnd}(n,\Delta x,-x_i) + \text{lnd}(n,\Delta x,x_i+\Delta x)]$<br>n = 2, 3, 4 …<br>for $\frac{x_i}{\Delta x}$ = all values<br><br>$-\infty$ to $+\infty$ for $Re(\Delta x)>0$ or $[Re(\Delta x)=0$ and $Im(\Delta x)>0]$ | Sum $-\infty$ to $+\infty$ or<br>Sum $+\infty$ to $-\infty$<br><br>Similar to #36 but n must be positive integers,<br>n = 1, 2, 3, 4, … |

| # | Function Calculated | Means of Calculation | Comments |
|---|---|---|---|
| | | $+\infty$ to $-\infty$ for $Re(\Delta x)<0$ or $[Re(\Delta x)=0$ and $Im(\Delta x)<0]$<br>$\Delta x$ = x increment<br>$\zeta(n)$ = Riemann Zeta Function<br>x, $x_i$, $\Delta x$ = real or complex values<br><br>Any summation term with x = 0 is excluded<br><br>Comment – If the value of $x_i$ is changed to $x_i + r\Delta x$, r = integers, the summation value remains the same. | |
| **38** | $\sum_{\Delta x,\ x=\mp\infty}^{\pm\infty} (-1)^{\frac{x-x_i}{\Delta x}} \frac{1}{x^n}$ | $\frac{\alpha(n)}{2\Delta x}$ [-lnd(n,2$\Delta x$,$x_i$+$\Delta x$)+lnd(n,2$\Delta x$,$x_i$)-lnd(n,-2$\Delta x$,$x_i$-2$\Delta x$)+lnd(n,-2$\Delta x$,$x_i$-$\Delta x$)]<br><br>$Re(n)>0$<br>$\alpha(n) = -1$ for $n = 1$<br>$\alpha(n) = +1$ for $n \neq 1$<br>$-\infty$ to $+\infty$ for $Re(\Delta x)>0$ or $[Re(\Delta x)=0$ and $Im(\Delta x)>0]$<br>$+\infty$ to $-\infty$ for $Re(\Delta x)<0$ or $[Re(\Delta x)=0$ and $Im(\Delta x)<0]$<br>$x = x_i + m\Delta x$, m = integers<br>$x_i$ = value of x<br>$x\neq 0$ | Alternating sum $x_i$ to $+\infty$ or $x_i$ to $-\infty$ |
| **39** | $\beta(w,z)$ | 1) $\sqrt{2\pi}\, e^{\lim_{n\to 0} \frac{lnd(n,1,w)-lnd(0,1,w)+lnd(n,1,z)-lnd(0,1,z)-lnd(n,1,w+z)+lnd(0,1,w+z)}{n}}$<br><br>2) $\frac{\Gamma(w)\Gamma(z)}{\Gamma(w+z)}$ | Beta Function<br><br>Use #19, #22 |

| # | Function Calculated | Means of Calculation | Comments |
|---|---|---|---|
| **40** | $\pi$ | -lnd(1,4,4m+5) + lnd(1,-4,4m+1)<br>m = integer | Pi |
| **41** | $\gamma$ | $\lim_{n\to 0}$ [ lnd(1+n,1,1) – $\frac{1}{n}$ ]<br>n = real numbers<br><br>$\lim_{x\to\infty}$ [ lnd(1,1,x) – ln(x – $\frac{1}{2}$) ] | Euler’s Constant |

In this table the equation constants and variables can be any real or complex value unless otherwise specified.

# Secondary Relationships

| # | Function Calculated | Means of Calculation | Comments |
|---|---|---|---|
| **1** | csc x | $\frac{1}{\sin x}$ | See P2 (P = Primary Relationships) |
| **2** | sec x | $\frac{1}{\cos x}$ | See P3 |
| **3** | cotan x | $\frac{1}{\tan x}$ | See P1 |
| **4** | csch x | $\frac{1}{\sinh x}$ | See P6 |
| **5** | sech x | $\frac{1}{\cosh x}$ | See P7 |
| **6** | cotanh x | $\frac{1}{\tanh x}$ | See P5 |
| **7** | $\log_{10} x$ | $(\log_{10} e)\ln x$ | See P17 |
| **8** | $\sin^{-1} x$ | $-j \ln[\ (1-x^2)^{1/2} + jx\ ],\ \ x^2 \leq 1$ | See P17, P18 |
| **9** | $\cos^{-1} x$ | $-j \ln[\ x + j(1-x^2)^{1/2}\ ],\ \ x^2 \leq 1$ | See P17, P18 |
| **10** | $\tan^{-1} x$ | $\frac{j}{2} \ln(\frac{1 - jx}{1 + jx}),\ \ x$ real | See P17 |

| # | Function Calculated | Means of Calculation | Comments |
|---|---|---|---|
| **11** | $\csc^{-1}x$ | $\sin^{-1}(\frac{1}{x})$, $x \geq 0$ | See S8 (S = Secondary Relationships) |
| **12** | $\sec^{-1}x$ | $\cos^{-1}(\frac{1}{x})$, $x \geq 0$ | See S9 |
| **13** | $\text{cotan}^{-1}x$ | $\tan^{-1}(\frac{1}{x})$, $x \geq 0$ | See S10 |
| **14** | $\sinh^{-1}x$ | $\ln[x + (x^2 + 1)^{1/2}]$ | See P17, P18 |
| **15** | $\cosh^{-1}x$ | $\ln[x + (x^2 - 1)^{1/2}]$, $x \geq 1$ | See P17, P18 |
| **16** | $\tanh^{-1}x$ | $\frac{1}{2}\ln(\frac{1 + x}{1 - x})$, $0 \leq x^2 < 1$ | See P17 |
| **17** | $\text{csch}^{-1}x$ | $\sinh^{-1}(\frac{1}{x})$ | See S14 |
| **18** | $\text{sech}^{-1}x$ | $\cosh^{-1}(\frac{1}{x})$ | See S15 |
| **19** | $\text{cotanh}^{-1}x$ | $\tanh^{-1}(\frac{1}{x})$ | See S16 |

| # | Function Calculated | Means of Calculation | Comments |
|---|---|---|---|
| **20** | $\csc_{\Delta x}(a,x)$ | $\frac{1}{\sin_{\Delta x}(a,x)}$ | See P10<br>(P = Primary Relationships) |
| **21** | $\sec_{\Delta x}(a,x)$ | $\frac{1}{\cos_{\Delta x}(a,x)}$ | See P11 |
| **22** | $\text{cotan}_{\Delta x}(a,x)$ | $\frac{1}{\tan_{\Delta x}(a,x)}$ | See P9 |
| **23** | $\text{csch}_{\Delta x}(a,x)$ | $\frac{1}{\sinh_{\Delta x}(a,x)}$ | See P14 |
| **24** | $\text{sech}_{\Delta x}(a,x)$ | $\frac{1}{\cosh_{\Delta x}(a,x)}$ | See P15 |
| **25** | $\text{cotanh}_{\Delta x}(a,x)$ | $\frac{1}{\tanh_{\Delta x}(a,x)}$ | See P13 |

# FUNCTIONS CALCULATED BY THE FUNCTION, LND(N,DX,X)

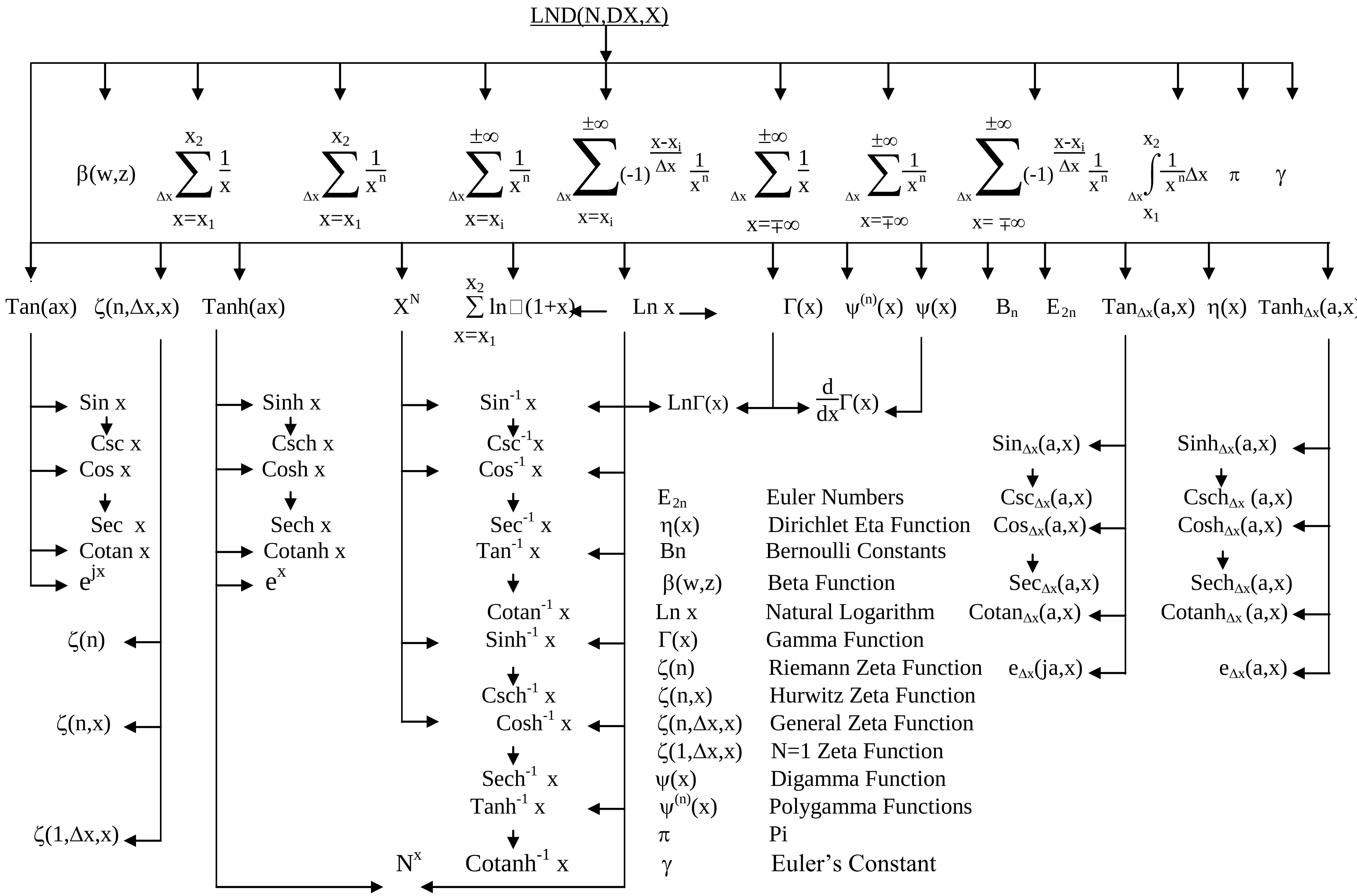

# TABLE 18

## Complex Plane Closed Contour Area Calculation Equations

| # | Closed Contour Area Calculation Equation | Initial Direction of Summation/ Integration | Comments | Equation Applicable to Complex Plane Closed Contours Composed of | | |
|---|---|---|---|---|---|---|
| | | | | **Vertical/Horizontal Straight Line Vectors Only S1** | **Straight Line Vectors with any Slope S2** | **Continuous Curves S3** |
| **1** | $A = \frac{ja}{2} \sum_{n=1}^{N} (-1)^{n-1} z_n^2$ | Horizontal or vertical then clockwise or counterclockwise | Use when all sides of a complex plane closed contour are alternating head to tail horizontal and vertical vectors<br><br>From this basic equation all of the following equations can be derived. | Yes | No | No |
| **2** | $A = -\frac{ja}{2} \sum_{n=1}^{\frac{N}{2}} (z_{2n}^2 - z_{2n-1}^2)$ | Horizontal or vertical then clockwise or counterclockwise | Use when all sides of a complex plane closed contour are alternating head to tail horizontal and vertical vectors | Yes | No | No |

| # | Closed Contour Area Calculation Equation | Initial Direction of Summation/ Integration | Comments | Equation Applicable to Complex Plane Closed Contours Composed of | | |
|---|---|---|---|---|---|---|
| | | | | Vertical/Horizontal Straight Line Vectors Only S1 | Straight Line Vectors with any Slope S2 | Continuous Curves S3 |
| **3** | $A = -a \sum_{n=1}^{\frac{N}{2}} (x_{2n-1}y_{2n-1} - x_{2n}y_{2n})$ or $A = -a \sum_{n=1}^{\frac{N}{2}} \begin{vmatrix} x_{2n-1} & y_{2n} \\ x_{2n} & y_{2n-1} \end{vmatrix}$ | Horizontal or vertical then clockwise or counterclockwise | Use when all sides of a complex plane closed contour are alternating head to tail horizontal and vertical vectors | Yes | No | No |
| **4** | $A = a \sum_{n=1}^{\frac{N}{2}} y_{2n-1}\Delta x_{2n-1}$ | Horizontal only then clockwise or counterclockwise | Use when all sides of a complex plane closed contour are alternating head to tail horizontal and vertical vectors | Yes | No | No |
| **5** | $A = a \sum_{n=1}^{\frac{N}{2}} x_{2n-1}\Delta y_{2n-1}$ | Vertical only then clockwise or counterclockwise | Use when all sides of a complex plane closed contour are alternating head to tail horizontal and vertical vectors | Yes | No | No |

| # | Closed Contour Area Calculation Equation | Initial Direction of Summation/ Integration | Comments | Equation Applicable to Complex Plane Closed Contours Composed of | | |
|---|---|---|---|---|---|---|
| | | | | Vertical/Horizontal Straight Line Vectors Only S1 | Straight Line Vectors with any Slope S2 | Continuous Curves S3 |
| **6** | $A = a \sum_{n=1}^{\frac{N}{2}} y_{2n}\Delta x_{2n-1}$ | Horizontal only then clockwise or counterclockwise | Use when all sides of a complex plane closed contour are alternating head to tail horizontal and vertical vectors | Yes | No | No |
| **7** | $A = a \sum_{n=1}^{\frac{N}{2}} x_{2n}\Delta y_{2n-1}$ | Vertical only then clockwise or counterclockwise | Use when all sides of a complex plane closed contour are alternating head to tail horizontal and vertical vectors | Yes | No | No |
| **8** | $A = -ja \sum_{n=1}^{\frac{N}{2}} \bar{z}_{2n-1}\Delta x_{2n-1}$ | Horizontal only then clockwise or counterclockwise | Use when all sides of a complex plane closed contour are alternating head to tail horizontal and vertical vectors | Yes | No | No |
| **9** | $A = a \sum_{n=1}^{\frac{N}{2}} \bar{z}_{2n-1}\Delta y_{2n-1}$ | Vertical only then clockwise or counterclockwise | Use when all sides of a complex plane closed contour are alternating head to tail horizontal and vertical vectors | Yes | No | No |

| # | Closed Contour Area Calculation Equation | Initial Direction of Summation/ Integration | Comments | Equation Applicable to Complex Plane Closed Contours Composed of | | |
|---|---|---|---|---|---|---|
| | | | | **Vertical/Horizontal Straight Line Vectors Only S1** | **Straight Line Vectors with any Slope S2** | **Continuous Curves S3** |
| **10** | $A = \frac{a}{2}\sum_{n=1}^{\frac{N}{2}}(y_{2n-1}\Delta x_{2n-1} - x_{2n}\Delta y_{2n})$ | Horizontal only then clockwise or counterclockwise | Use when all sides of a complex plane closed contour are alternating head to tail horizontal and vertical vectors | Yes | No | No |
| **11** | $A = \frac{a}{2}\sum_{n=1}^{\frac{N}{2}}(x_{2n-1}\Delta y_{2n-1} - y_{2n}\Delta x_{2n})$ | Vertical only then clockwise or counterclockwise | Use when all sides of a complex plane closed contour are alternating head to tail horizontal and vertical vectors | Yes | No | No |
| **12** | $A = -\frac{ja}{2}\sum_{n=1}^{\frac{N}{2}}[(x_{2n+1}+jy_{2n-1})^2 - (x_{2n-1}+jy_{2n-1})^2]$ | Horizontal only then clockwise or counterclockwise | Use when all sides of a complex plane closed contour are alternating head to tail horizontal and vertical vectors | Yes | No | No |
| **13** | $A = -\frac{ja}{2}\sum_{n=1}^{\frac{N}{2}}[(x_{2n-1}+jy_{2n+1})^2 - (x_{2n-1}+jy_{2n-1})^2]$ | Vertical only then clockwise or counterclockwise | Use when all sides of a complex plane closed contour are alternating head to tail horizontal and vertical vectors | Yes | No | No |

| # | Closed Contour Area Calculation Equation | Initial Direction of Summation/ Integration | Comments | Equation Applicable to Complex Plane Closed Contours Composed of | | |
|---|---|---|---|---|---|---|
| | | | | **Vertical/Horizontal Straight Line Vectors Only S1** | **Straight Line Vectors with any Slope S2** | **Continuous Curves S3** |
| **14** | $A = \frac{jb}{4}\sum_{n=1}^{N}[(x_{n+1}+jy_n)^2 - (x_n+jy_{n+1})^2]$ | Clockwise or counterclockwise | Use when all sides of a complex plane closed contour are head to tail vectors of any slope | Yes | Yes | Can use for approximation by applying to selected points on a continuous closed contour |
| **15** | $A = -b\sum_{n=1}^{N}\bar{y}_n\Delta x_n$ | Clockwise or counterclockwise | Use when all sides of a complex plane closed contour are head to tail vectors of any slope | Yes | Yes | Can use for approximation by applying to selected points on a continuous closed contour |
| **16** | $A = b\sum_{n=1}^{N}\bar{x}_n\Delta y_n$ | Clockwise or counterclockwise | Use when all sides of a complex plane closed contour are head to tail vectors of any slope | Yes | Yes | Can use for approximation by applying to selected points on a continuous closed contour |
| **17** | $A = jb\sum_{n=1}^{N}\bar{z}_n\Delta x_n$ | Clockwise or counterclockwise | Use when all sides of a complex plane closed contour are head to tail vectors of any slope | Yes | Yes | Can use for approximation by applying to selected points on a continuous closed contour |

| # | Closed Contour Area Calculation Equation | Initial Direction of Summation/ Integration | Comments | Equation Applicable to Complex Plane Closed Contours Composed of | | |
|---|---|---|---|---|---|---|
| | | | | **Vertical/Horizontal Straight Line Vectors Only S1** | **Straight Line Vectors with any Slope S2** | **Continuous Curves S3** |
| **18** | $A = b\sum_{n=1}^{N}\bar{z}_n\Delta y_n$ | Clockwise or counterclockwise | Use when all sides of a complex plane closed contour are head to tail vectors of any slope | Yes | Yes | Can use for approximation by applying to selected points on a continuous closed contour |
| **19** | $A = \frac{jb}{2}\sum_{n=1}^{N}\bar{z}_n\Delta z_n^*$ | Clockwise or counterclockwise | Use when all sides of a complex plane closed contour are head to tail vectors of any slope | Yes | Yes | Can use for approximation by applying to selected points on a continous closed contour |
| **20** | $A = -\frac{jb}{2}\sum_{n=1}^{N}\bar{z}_n^*\Delta z_n$ | Clockwise or counterclockwise | Use when all sides of a complex plane closed contour are head to tail vectors of any slope | Yes | Yes | Can use for approximation by applying to selected points on a continuous closed contour |
| **21** | $A = \frac{b}{2}\sum_{n=1}^{N}(x_n y_{n+1} - x_{n+1} y_n)$ or $A = \frac{b}{2}\sum_{n=1}^{N}\begin{vmatrix} x_n & y_n \\ x_{n+1} & y_{n+1} \end{vmatrix}$ | Clockwise or counterclockwise | Use when all sides of a complex plane closed contour are head to tail vectors of any slope | Yes | Yes | Can use for approximation by applying to selected points on a continuous closed contour |

| # | Closed Contour Area Calculation Equation | Initial Direction of Summation/ Integration | Comments | Equation Applicable to Complex Plane Closed Contours Composed of | | |
|---|---|---|---|---|---|---|
| | | | | **Vertical/Horizontal Straight Line Vectors Only S1** | **Straight Line Vectors with any Slope S2** | **Continuous Curves S3** |
| **22** | $A = \frac{b}{2}\sum_{n=1}^{N}(\bar{x}_n\Delta y_n - \bar{y}_n\Delta x_n)$ | Clockwise or counterclockwise | Use when all sides of a complex plane closed contour are head to tail vectors of any slope | Yes | Yes | Can use for approximation by applying to selected points on a continous closed contour |
| **23** | $A = \frac{b}{2}\sum_{n=1}^{N}(x_{n+1}\Delta y_n - y_{n+1}\Delta x_n)$ | Clockwise or counterclockwise | Use when all sides of a complex plane closed contour are head to tail vectors of any slope | Yes | Yes | Can use for approximation by applying to selected points on a continous closed contour |
| **24** | $A = \frac{b}{2}\sum_{n=1}^{N}(x_n\Delta y_n - y_n\Delta x_n)$ | Clockwise or counterclockwise | Use when all sides of a complex plane closed contour are head to tail vectors of any slope<br><br>This is the discrete form of Green's Theorem for the calculation of area | Yes | Yes | Can use for approximation by applying to selected points on a continous closed contour |

| # | Closed Contour Area Calculation Equation | Initial Direction of Summation/ Integration | Comments | Equation Applicable to Complex Plane Closed Contours Composed of | | |
|---|---|---|---|---|---|---|
| | | | | **Vertical/Horizontal Straight Line Vectors Only S1** | **Straight Line Vectors with any Slope S2** | **Continuous Curves S3** |
| **25** | $A = \frac{b}{2} \oint_c [xdy - ydx]$ | Clockwise or counterclockwise | Use when a complex plane closed contour is a continuous closed curve<br><br>This is Green's Theorem for the calculation of area | No | No | Yes |
| **26** | $A = -b \oint_c ydx$ | Clockwise or counterclockwise | Use when a complex plane closed contour is a continuous closed curve | No | No | Yes |
| **27** | $A = b \oint_c xdy$ | Clockwise or counterclockwise | Use when a complex plane closed contour is a continuous closed curve | No | No | Yes |
| **28** | $A = jb \oint_c zdx$ | Clockwise or counterclockwise | Use when a complex plane closed contour is a continuous closed curve | No | No | Yes |

| # | Closed Contour Area Calculation Equation | Initial Direction of Summation/ Integration | Comments | Equation Applicable to Complex Plane Closed Contours Composed of | | |
|---|---|---|---|---|---|---|
| | | | | **Vertical/Horizontal Straight Line Vectors Only S1** | **Straight Line Vectors with any Slope S2** | **Continuous Curves S3** |
| **29** | $A = b\oint_{c} z\,dy$ | Clockwise or counterclockwise | Use when a complex plane closed contour is a continuous closed curve | No | No | Yes |
| **30** | $A = \frac{jb}{2}\oint_{c} z\,dz^{*}$ | Clockwise or counterclockwise | Use when a complex plane closed contour is a continuous closed curve | No | No | Yes |
| **31** | $A = -\frac{jb}{2}\oint_{c} z^{*}dz$ | Clockwise or counterclockwise | Use when a complex plane closed contour is a continuous closed curve | No | No | Yes |

Variable and Constant Definitions

Summation and integration direction constants a and b

$$a = \begin{cases} +1 & \text{For initial direction horizontal then clockwise} \\ -1 & \text{For initial direction horizontal then counterclockwise} \\ -1 & \text{For initial direction vertical then clockwise} \\ +1 & \text{For initial direction vertical then counterclockwise} \end{cases}$$

$$b = \begin{cases} -1 & \text{For clockwise direction} \\ +1 & \text{For counterclockwise direction} \end{cases}$$

Categorization of Complex Plane Closed Contour Shapes

1. S1 – Discrete complex plane closed contours composed entirely of alternating horizontal and vertical vectors connected head to tail forming a closed loop
2. S2 – Discrete complex plane closed contours composed of vectors of any slope connected head to tail forming a closed loop
3. S3 – Continuous complex plane closed contours composed of infinitesimal vectors of any slope connected head to tail forming a closed loop

Examples of S1, S2, and S3 Complex Plane Closed Contours

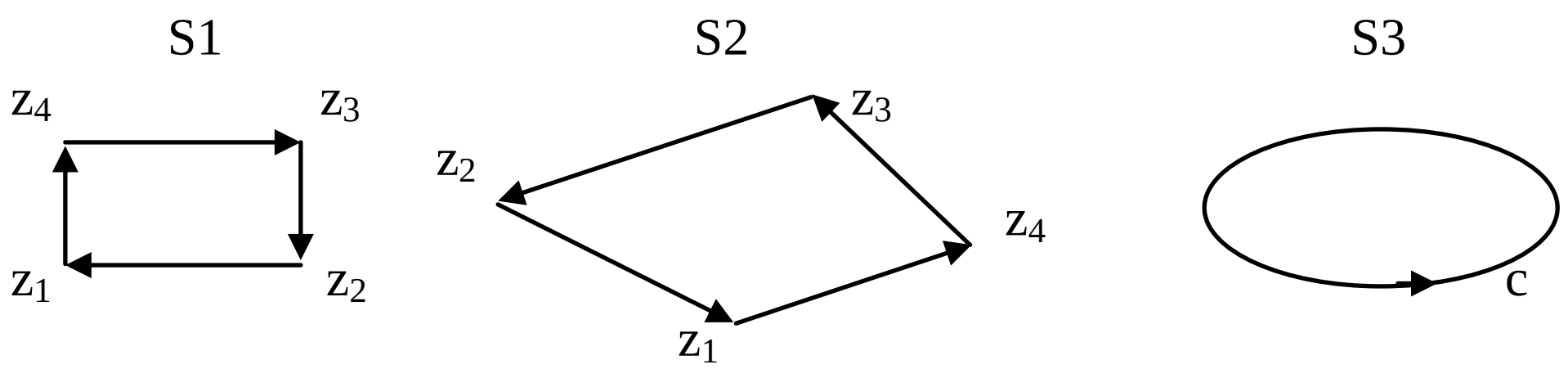

Clockwise Discrete Closed Contour | Counterclockwise Discrete Closed Contour | Counterclockwise Continuous Closed Contour

$z_1, z_2, z_3, z_4$ = Closed contour corner points defining the S1 and S2 discrete closed contours
c = The S3 continous closed contour

Comment 1 - N = the number of corner points that define a discrete closed contour in the complex plane
For S1 closed contours N = 4, 6, 8, 10, …
For S2 closed contours N = 3, 4, 5, 6, …
For S3 closed contours N is infinite

Comment 2 - An S1 and S3 closed contour is a special case of the S2 closed contour.

$$A = \frac{1}{2}(A_H + A_V)$$

where

$A$ = the area enclosed within a complex plane S2 closed contour

$A_H, A_V$ = the areas of two complex plane S1 closed contours constructed from the S2 closed contour corner points

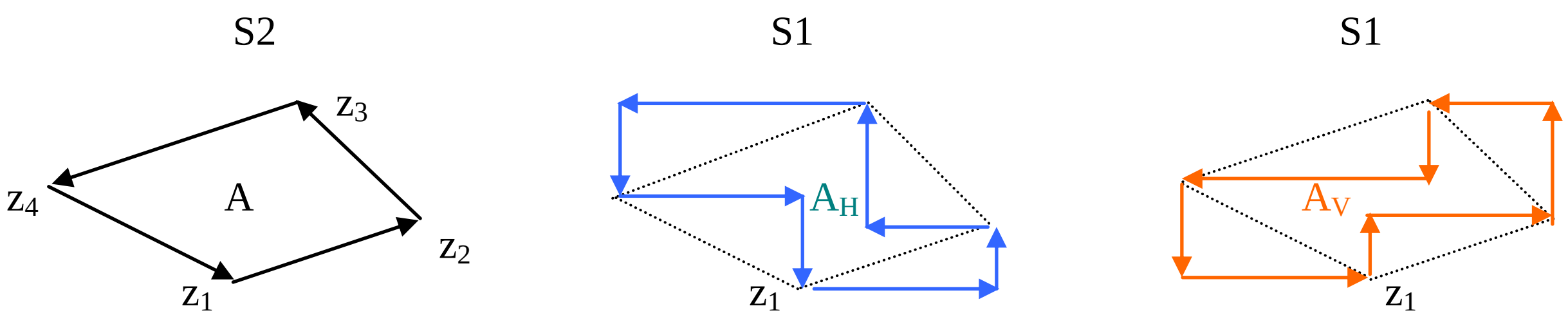

Counterclockwise | Counterclockwise, Initial Direction horizontal | Counterclockwise, Initial direction vertical

It has been shown that the area enclosed within an S2 closed contour can be calculated from half the sum of the areas enclosed within two S1 closed contours constructed from the same corner points of the S2 closed contour. Note the above complex plane S2 and two S1 closed contours.

Notation and Variable Definitions

$m$ = defined positive integers

$x_m$ = real values

$y_m$ = real values

$z_m = x_m + jy_m$ , the set of complex plane corner points defining a discrete closed contour

Note 1 – $z_m$ can be represented by a discrete function of m. For example, $z_m = e^{\frac{j2\pi(m-1)}{M}}$ where $m = 1, 2, 3, \ldots, M$

$$\bar{x}_m = \frac{x_{m+1} + x_m}{2} , \quad \bar{y}_m = \frac{y_{m+1} + y_m}{2} , \quad \bar{z}_m = \frac{z_{m+1} + z_m}{2}$$

$\Delta x_m = x_{m+1} - x_m$ , $\Delta y_m = y_{m+1} - y_m$ , $\Delta z_m = z_{m+1} - z_m$

$x_{N+1} = x_1$ , $y_{N+1} = y_1$ , $z_{N+1} = z_1$

$x(\theta)$ = A continuous real value function of $\theta$

$y(\theta)$ = A continuous real value function of $\theta$

$z(\theta) = x(\theta) + jy(\theta)$ , A continuous complex value function of $\theta$ defining a closed contour in the complex plane (ex: $z(\theta) = e^{j\theta}$)

A = The area enclosed within an S1, S2, or S3 closed contour in the complex plane. The area, A, is always a positive value

$A_H$ = The area enclosed within the S1 closed contour constructed from an S2 closed contour by starting at an initial S2 corner point, $z_1$, and then drawing alternately horizontal then vertical vectors connecting the S2 corner points. The two S1 closed contours contructed from the same S2 corner points must have the same sense, clockwise or counterclockwise.

$A_V$ = The area enclosed within the S1 closed contour constructed from an S2 closed contour by starting at an initial S2 corner point, $z_1$, and then drawing alternately vertical then horizontal vectors connecting the S2 corner points. The two S1 closed contours contructed from the same S2 corner points must have the same sense, clockwise or counterclockwise.

N = The number of corner points on the complex plane closed contour (i.e. The number of corner points defining the complex plane closed contour).

* = complex conjugate designation (ex: $[a+jb]^* = a-jb$ or $[e^{j\theta}]^* = e^{-j\theta}$ )

$\oint_c$ = The integration sign representing complex plane integration along the continuous closed contour, c

<u>Note 2</u> - Consecutive contour points, designated $z_1$, $z_2$, $z_3$, …, are connected by straight line vectors. These vectors are connected head to tail to form a closed loop.

Comment 3 - There are several interesting closed contour summations and integrals that equal zero. They are as follows:
For S1 complex plane closed contours

1. $\sum_{n=1}^{\frac{N}{2}} [\,(x_{2n}^2 - x_{2n-1}^2) = 0$ , $x_N = x_1$ , $x_{2n+1} = x_{2n}$ , $y_{2n} = y_{2n-1}$ , Initial direction of summation = horizontal

2. $\sum_{n=1}^{\frac{N}{2}} [\,(y_{2n}^2 - y_{2n-1}^2) = 0$ , $y_N = y_1$ , $y_{2n+1} = y_{2n}$ , $x_{2n} = x_{2n-1}$ , Initial direction of summation = vertical

For S2 complex plane closed contours

3. $\sum_{n=1}^{N} (x_{n+1}^2 - x_n^2) = 0$ , $x_{N+1} = x_1$

4. $\sum_{n=1}^{N} (y_{n+1}^2 - y_n^2) = 0$ , $y_{N+1} = y_1$

5. $\sum_{n=1}^{N} \overline{x}_n \Delta x_n = 0$

6. $\sum_{n=1}^{N} \overline{y}_n \Delta y_n = 0$

7. $\sum_{n=1}^{N} \bar{z}_n \Delta z_n = 0$

8. $\sum_{n=1}^{N} \bar{z}_n^{\,*} \Delta z_n^{\,*} = 0$

For S3 complex plane closed contours

9. $\oint_c x(\theta)dx(\theta) = 0$

10. $\oint_c y(\theta)dy(\theta) = 0$

11. $\oint_c z(\theta)dz(\theta) = 0$

12. $\oint_c z(\theta)^{*}dz(\theta)^{*} = 0$

<u>Comment 4</u> - The closed contour summations previously listed have another form, a discrete integral form. The discrete integral form is Particularly useful when a function is available to describe the closed contour for which the enclosed area is to be calculated

For example: $A = \frac{jb}{2}\sum_{n=1}^{N} \bar{z}_n \Delta z_n^{\,*} = \frac{jb}{2}\;{}_1\!\int_{1}^{N+1} \bar{z}_n \Delta z_n^{\,*}$

# TABLE 19

## Integral Calculus Conventions

| | |
|---|---|
| **1** | A summation term with a division by zero is generally excluded from the summation. |
| **2** | The symbol, $\Delta x$, when it appears in an equation, is used to represent the x increment. |
| **3** | The value of $\Delta x$, designated as, $_{\Delta x}$, is specified on the Interval Calculus mathematical operation symbols, for integration, $_{\Delta x}\int$, summation, $_{\Delta x}\sum$, product, $_{\Delta x}\prod$, and limit specification, $\Big\vert_{x_1}^{x_2}{}_{\Delta x}$. The $\Delta x$ value symbol, $_{\Delta x}$, may be replaced by its equivalent numerical value. |
| **4** | If no $\Delta x$ value appears on the $_{\Delta x}\sum$ or $_{\Delta x}\prod$ mathematical operation symbols, the value of $\Delta x$ is considered to be 1. If no $\Delta x$ value appears on the $_{\Delta x}\int$ mathematical operation symbol, the value of $\Delta x$ is infinitesimally small (i.e. $\Delta x \to 0$). If no $\Delta x$ value appears on the $\Big\vert_{x_1}^{x_2}{}_{\Delta x}$ mathematical operation symbol, the value of $\Delta x$ is determined from other operation symbols in the equation in which this symbol appears. In many cases, the value of $\Delta x$ is not needed. Only the upper and lower limit values, $x_1$ and $x_2$, are needed. |
| **5** | For the mathematical operation symbols $_{\Delta x}\int_{x_1}^{x_2}$, $_{\Delta x}\sum_{x=x_1}^{x_2}$, $_{\Delta x}\prod_{x_1}^{x_2}$, and limit specification, $\Big\vert_{x_1}^{x_2}{}_{\Delta x}$, the lower limit value (here $x_1$) is the initial value of x and the upper limit value (here $x_2$) is the final value of x. Also, the value of $\Delta x$, $_{\Delta x}$, must be consistent with the values of $x_1$ and $x_2$ (i.e. $x_1+n\Delta x = x_2$ where n = the number of x intervals between $x_1$ and $x_2$). |

# TABLE 20

## The lnd(n,Δt,t) Function Summation Relationship to the Real Value of n

The lnd(1,Δt,t) function is used to replace the natural log which has limited use in discrete mathematics. Much of Interval Calculus deals with the lnd(n,Δt,t) function which is defined by the equation, $D_{\Delta t}\text{lnd}(n,\Delta t,t) = \pm\frac{1}{t^n}$ , + for n = 1, − for n ≠ 1. The function, lnd(n,Δt,t) has three variables, n,Δt,t, which can be either real or complex. It is a very useful function in discrete mathematics. Note the discription of the lnd(n,Δt,t) function below. The lnd(n,Δx,x) function was derived to evaluate the commonly occurring discrete (variable) integration, $\int_{\Delta t\; t_1}^{t_2}\frac{1}{t^n}\,\Delta t = \Delta t\sum_{\Delta t\; t=t_1}^{t_2-\Delta t}\frac{1}{t^n}$ where t = 0,Δt,2Δt,3Δt, … and n,Δt,t are real or complex.

For Re(n) > 1

$$\int_{\Delta t\; t_1}^{t_2}\frac{1}{t^n}\,\Delta t = \Delta t\sum_{\Delta t\; t=t_1}^{t_2-\Delta t}\frac{1}{t^n} = -\,\text{lnd}(n,\Delta t,t)\Big|_{t_1}^{t_2}$$

$t_2$, can be infinite.

Comment

$\lim_{t\to\infty}\text{lnd}(n,\Delta t,t) = 0$ where Re(n) > 1, Re(Δx) > 0

$$\int_{\Delta t\; t}^{\infty}\frac{1}{t^n}\,\Delta t = \Delta t\sum_{\Delta t\; t=t}^{\infty}\frac{1}{t^n} = +\,\text{lnd}(n,\Delta t,t)\ ,\quad \text{Re}(n) > 1,\ \text{R}(\Delta x) > 0$$

It is observed that the function, lnd(n,Δt,t), provides an evaluation of the Zeta Functions.

For Re(n) = 1

$$\int_{\Delta t\; t_1}^{t_2}\frac{1}{t^n}\,\Delta x = \Delta t\sum_{\Delta t\; t=t_1}^{t_2-\Delta t}\frac{1}{t^n} = +\,\text{lnd}(n,\Delta t,t)\Big|_{t_1}^{t_2}$$

For 0 < Re(n) < 1

$$\int_{\Delta t\, t_1}^{t_2} \frac{1}{t^n}\,\Delta t = \Delta t \sum_{\Delta t\, t=t_1}^{t_2-\Delta t} \frac{1}{t^n} = -\,\mathrm{lnd}(n,\Delta t,t)\Big|_{t_1}^{t_2}$$

For Re(n) = 0

$$\int_{\Delta t\, t_1}^{t_2} \frac{1}{t^n}\,\Delta x = \Delta t \sum_{\Delta t\, t=t_1}^{t_2-\Delta t} \frac{1}{t^n} = -\,\mathrm{lnd}(n,\Delta t,t)\Big|_{t_1}^{t_2}$$

For Re(n) < 0

$$\int_{\Delta t\, t_1}^{t_2} t^{-n}\,\Delta x = \Delta t \sum_{\Delta t\, t=t_1}^{t_2-\Delta t} t^{-n} = -\,\mathrm{lnd}(n,\Delta t,t)\Big|_{t_1}^{t_2}$$

The computer program that evaluates the lnd(n,Δt,t) function is named LNDX and made available in this document (See the Abstract on page 1). From the lnd(n,Δt,t) function , one can evaluate sixty important mathematical functions including the Riemann, Hurwitz, Polygamma, Digamma, and Gamma functions, the Calculus trigonometric and hyperbolic engineering math functions, and discrete analogs of the Calculus engineering math functions. Moreover, in the course of writing the computer program to calculate the discrete natural log, much interesting and novel math has been found which involves neither the Euler-MacLaurin summation nor Bernoulli polynomials. The lnd(n,Δx,x) function can evaluate summations, calculate the area under sample and hold shaped curves, and solve differential difference equations and difference equations.

# CALCULATION PROGRAMS

## 1) UBASIC Program to Calculate the Function, LND(n,Δx,x)

Comment – UBASIC is a freeware BASIC interpreter written by Yuji Kida at Rikkyo University in Japan. UBASIC is specialized for mathematical computing. The following programs were written and tested for Version 8.21 and are downloadable from https://www.box.com/s/c59uvksuwomnw8rsyj92. Note that at time of writing (Jan 2013), UBASIC Version 8.74 is the latest stable release, and it is downloadable at: archives.math.utk.edu/software/msdos/number.theory/UBASIC

"LNDX.UB"

```
10  'Calculation of ln(DX,X) and LND(N,DX,X) -- N,DX,X can be complex
20  point 20 'Set decimal digits used in calculations
30  word 280 'Set word length of long variables
40  A=50 'Print A decimal digits in the printouts
50  B=25 'B+1 is the number of LNDf(N,DX,X) series terms used
60  C=100 'NN=C, the countup/countdown number
70  SUM=0 'SUM is a summation variable. Initialize SUM to 0
80  CT=0 'Initialize CT to 0. CT is a counter variable, CT=0 or 1
90  F2=0 'Flag that selects LND(N,DX,X) computation method, Tk=2, when
    F2=1
100 FG=0 'Initialize "formula only" flag to 0. FG=1 selects formula only
110 FF=0 'Flag used to avoid a division by 0
120 FB=0 'Flag to bypass a second input of N and DX, FB=1 disables bypass
130 SN=0.00000000000000001 'SN is a small number used for calculations
140 LN=100000000000000000 'LN is a large number used for calculations
150 X1=0:X2=0 'Initialize the X locus cross over points to 0
160 print
170 print "The following inputs can be real or complex, A+B#i"
180 print "DX can not equal 0"
190 print
200 input "What is N ";N
210 input "What is DX ";DX
```

```
220  input "What is X ";X
230  print
240  'Identify some characteristics of DX
250  if and{re(DX)=0,im(DX)=0} then
260  :print "DX can not equal 0"
270  :goto 3980
280  :endif
290  if re(DX)=0 then FD=1 'The real part of DX is 0
300  :else FD=0 'The real part of DX is not equal to 0
310  :endif
320  dim C(70),GR(70),G(70),J(70) 'Define arrays
330  'Calculate the LNDf(N,DX,X) series formula coefficients
340  C(0)=-1 'Cm coefficients
350  C(1)=1
360  C(2)=-7/3
370  C(3)=31/3
380  C(4)=-381/5
390  C(5)=2555/3
400  C(6)=-1414477/105
410  C(7)=286685
420  C(8)=-118518239/15
430  C(9)=5749691557/21
440  C(10)=-640823941499/55
450  C(11)=1792042792463/3
460  C(12)=-9913827341556185/273
470  C(13)=2582950540044537
480  C(14)=-3187598676787461083/15
490  C(15)=4625594554880206790555/231
500  C(16)=-182112049520351725262389/85
510  C(17)=774975953478559072551395/3
520  C(18)=-904185845619475242495834469891/25935
530  C(19)=5235038863532181733472316757
```

```
540   C(20)=-143531742398845896012634103722237/165
550   C(21)=33427300696841208116528825914877741/21
560   C(22)=-734472084995088305142162030978467283/23
570   C(23)=2098575784311706718209533060163655359l/3
580   C(24)=-55261733892726368667839334271075797592500083/3315
590   C(25)=4737771320732953193072519494466008540099675/11
600   if and{re(N)=1,im(N)=0} then goto 3000 'Calculate ln(DX,X)=ln1(X/DX)
610   'Calculate LND(N,DX,X) where N<>1
620   J(0)=1/(N-1) 'Calculate series remainder component terms
630   for Y=1 to B
640   J(Y)=((DX/2)^(2*Y))/(2*Y+1)/(2*Y)
650   next Y
660   if and{re(N)<0,(re(N)-fix(re(N)))=0,im(N)=0} then FG=1 'formula only
670   if F2=1 then goto 2020
680   gosub 760 'Calculate LND(N,DX,X)
690   E=P
700   SCn=Cn 'Save the value of Cn
710   if re(N)<=0 then goto 740 'Don't calculate K
720   G=Z1
730   gosub 2120 'Calculate LND(N,DX,-inf.) or LND(N,DX,+inf.), find K
740   Cn=SCn 'Restore the value of Cn
750   goto 3580 'Print calculated values
760   'LND(N,DX,X) calculation subroutine starts here
770   'Inputs are X,DX,N and the output is P
780   Tk=1 'LND(N,DX,X) Calculation Method #1
790   'Calculate the X locus X axis crossover point, XC
800   if im(DX)<>0 then XC=re(X)-im(X)*(re(DX)/im(DX))
810   :elseif or{im(DX)=0,re(DX)=0} then XC=0
820   :endif
830   'Define the crossover points
840   if or{and{re(DX)>0,XC>0},and{re(DX)<0,XC<0}} then
850   :X1=XC:X2=0 else X1=0:X2=XC
```

```
860   :elseif XC=0 then X1=XC,X2=XC
870   :endif
880   NN=C
890   if FD=0 then IN=NN*re(DX) elseif FD=1 then IN=NN*im(DX)
900   'Calculate Z1 thru Z4 which will define X ranges
910   Z1=X1+IN
920   Z2=X1-IN
930   Z3=X2+IN
940   Z4=X2-IN
950   Z5=X1-14*IN
960   G=Z1 'Initialize G
970   'Define the calculation conditions
980   if FD=0 then V=re(X) else if FD=1 then V=im(X)
990   'Condition #1
1000  if and{re(DX)>0,V>=Z5} then goto 1060
1010  :elseif and{re(DX)<0,V<=Z5} then goto 1060
1020  :elseif and{re(DX)=0,im(DX)>0,V>=Z5} then goto 1060
1030  :elseif and{re(DX)=0,im(DX)<0,V<Z5} then goto 1060
1040  :else goto 1100
1050  :endif
1060  Cn=1
1070  Z=X
1080  gosub 2360 'Calculate LND(N,DX,X)
1090  goto 2010
1100  'Condition #2
1110  if and{re(DX)>=0,abs(X1-X2)<=abs(10*IN),V<Z5} then goto 1150
1120  :elseif and{re(DX)<=0,abs(X1-X2)<=abs(10*IN),V>Z5} then goto 1150
1130  :else goto 1350
1140  :endif
1150  Cn=2
1160  'Calculate C2
1170  G=Z1
```

```
1180 ZI=Z5 'Calculate XZ5
1190 gosub 2940
1200 XZ5=XO
1210 Z=XZ5
1220 gosub 2360 'Calculate LNDc(N,DX,XZ5)
1230 A1=P
1240 FG=1 'Set the "formula only" flag
1250 Z=XZ5
1260 gosub 2360 'Calculate LNDf(N,DX,XZ5)
1270 A2=P
1280 C2=A1-A2 'Calculate C2
1290 FG=1 'Set the "formula only" flag
1300 Z=X
1310 gosub 2360 'Calculate LNDf(N,DX,X)
1320 A3=P
1330 P=A3+C2 'Calculate LND(N,DX,X)
1340 goto 2010
1350 'Conditions #3, #4, #5
1360 'Calculate C1
1370 G=Z1 'Define G
1380 ZI=Z2 'Calculate XZ2
1390 gosub 2940
1400 XZ2=XO
1410 Z=XZ2
1420 gosub 2360 'Calculate LNDc(N,DX,XZ2)
1430 A1=P
1440 FG=1 'Set the "formula only flag"
1450 Z=XZ2
1460 gosub 2360 'Calculate LNDf(N,DX,XZ2)
1470 A2=P
1480 C1=A1-A2 'Calculate C1
1490 'Calculate C2
```

```
1500 G=Z3 'Define G
1510 ZI=Z4 'Calculate XZ4
1520 gosub 2940
1530 XZ4=XO
1540 Z=XZ4
1550 gosub 2360 'Calculate LNDc(N,DX,XZ4)
1560 A3=P
1570 FG=1 'Set the "formula only" flag
1580 Z=XZ4
1590 gosub 2360 'Calculate LNDc(N,DX,XZ4)
1600 A4=P
1610 A5=A3-A4
1620 C2=A5+C1 'Calculate C2
1630 'Condition #3
1640 if and{re(DX)>=0,V<=Z2,V>Z3} then goto 1680
1650 :elseif and{re(DX)<=0,V>=Z2,V<Z3} then goto 1680
1660 :else goto 1750
1670 :endif
1680 Cn=3
1690 FG=1 'Set the "formula only" flag
1700 Z=X
1710 gosub 2360 'Calculate LNDf(N,DX,X)
1720 A6=P
1730 P=A6+C1
1740 goto 2010
1750 'Condition #4
1760 if and{re(DX)>=0,V<=Z3,V>Z4} then goto 1800
1770 :elseif and{re(DX)<=0,V>=Z3,V<Z4} then goto 1800
1780 :else goto 1870
1790 :endif
1800 Cn=4
1810 G=Z3 'Define G
```

```
1820 Z=X
1830 gosub 2360 'Calculate LNDf(N,DX,X)
1840 A7=P
1850 P=A7+C1
1860 goto 2010
1870 'Condition #5
1880 if and{re(DX)>=0,V<=Z4} then goto 1920
1890 :elseif and{re(DX)<=0,V>=Z4} then goto 1920
1900 :else goto 1990
1910 :endif
1920 Cn=5
1930 FG=1 'Set the "formula only" flag
1940 Z=X
1950 gosub 2360 'Calculate LNDf(N,DX,X)
1960 A8=P
1970 P=A8+C2
1980 goto 2010
1990 print "There is a program error"
2000 goto 3690
2010 return
2020 Tk=2 'LND(N,DX,X) Calculation Method #2
2030 if re(N)<=0 then goto 2060 'Don't calculate K
2040 'Calculate K
2050 gosub 2120 'Calculate LND(N,DX,-inf.) or LND(N,DX,+inf.), find K
2060 'Calculate LND(N,DX,X)
2070 Z=X
2080 gosub 2360 'Find LND(N,DX,X)
2090 E=P
2100 'Print LND(N,DX,X)
2110 goto 3580
2120 'Subroutine to calculate LND(N,DX,-inf.) or LND(N,DX,+inf.)
2130 'Inputs are X,DX,N,G and the output is K
```

```
2140 'A reference value of X, XA, must be calculated
2150 'Calculate XA for Re(DX)<>0 or Re(DX)=0
2160 if FD=0 then
2170 :XA=re(X)+re(DX)*fix(-re(X)/re(DX)) 'Real part
2180 :XA=XA+(im(X)+im(DX)*fix(-re(X)/re(DX)))*#i 'Imaginary part added
2190 :elseif FD=1 then
2200 :XA=re(X)+(im(X)+fix(-im(X)/im(DX))*im(DX))*#i
2210 :endif
2220 SX=X 'Save the value of X
2230 if F2=1 then Z=XA else X=XA
2240 if F2=1 then gosub 2360 else gosub 760 'Find LND(N,DX,XA)
2250 K1=P
2260 DX=-DX 'Negate the value of DX
2270 G=-G 'Negate the value of G
2280 if F2=1 then Z=XA+DX else X=XA+DX
2290 if F2=1 then gosub 2360 else gosub 760 'Find LND(N,-DX,XA-DX)
2300 K2=P
2310 DX=-DX 'Restore the value of DX
2320 G=-G 'Restore the value of G
2330 K=K1-K2
2340 X=SX 'Restore the value of X
2350 return
2360 'Subroutine to calculate LNDc(N,DX,X) and LNDf(N,DX,X)
2370 'Inputs are Z,DX,N,G and the output is P
2380 if F2=0 then goto 2510 'The following code is used only for Tk=2
2390 'Calculate the X locus X axis crossover point, XC
2400 if im(DX)<>0 then XC=re(X)-im(X)*(re(DX)/im(DX))
2410 :elseif im(DX)=0 then XC=im(X)
2420 :endif
2430 'Define the crossover points
2440 if or{and{re(DX)>0,XC>0},and{re(DX)<0,XC<0}} then
2450 :X1=XC else X1=0
```

```
2460 :endif
2470 NN=C
2480 if FD=0 then G=X1+NN*re(DX)
2490 :elseif FD=1 then G=X1+NN*im(DX) 'Number selected for higher accuracy
2500 :endif
2510 H=1000000 'The maximum number of calculation loops
2520 M=0 'Initialize loop counter
2530 TS=0 'Set the initial value of the truncated series sum to 0
2540 'Check if the accuracy will be OK or if the "formula only" flag is set
2550 if and{FD=0,re(DX)>0,re(Z)>=G} then goto 2690
2560 :elseif and{FD=0,re(DX)<0,re(Z)<=G} then goto 2690
2570 :elseif and{FD=1,im(DX)>0,im(Z)>=G} then goto 2690
2580 :elseif and{FD=1,im(DX)<0,im(Z)<=G} then goto 2690
2590 :elseif FG=1 then goto 2690
2600 :endif
2610 'The Count-up/Count-down routine below is needed for computation
2620 while or{and{FD=0,re(DX)>0,re(Z)<G},and{FD=0,re(DX)<0,re(Z)>G},and{FD=1,im(DX)>0,im(Z)<G},and{FD=1,im(DX)<0,im(Z)>G}}
2630 M=M+1 'Count series calculation loops
2640 if and{re(N)>0,abs(Z)<0.00000000001} then goto 2660 'Exclude 1/0 terms
2650 if re(N)>0 then TS=TS+DX/(Z^N) else TS=TS+DX*(Z^(-N)) 'Evaluate sum
2660 Z=Z+DX 'Increase the series Z values
2670 if M>H-1 goto 2920 'Test for excessive calculation looping
2680 wend
2690 if (Z-DX/2)=0 then 'Avoid a division by 0
2700 :SZ=Z 'Save the value of Z
2710 :Z=Z+DX 'Increment the value of Z
2720 :FF=1 'Set the division by 0 avoidance flag to 1
2730 :endif
2740 G(0)=(Z-DX/2) 'Calculate series remainder component terms
```

```
2750 G(1)=N/(Z-DX/2)
2760 for Y=2 to B
2770 G(Y)=G(Y-1)*(N+2*Y-3)*(N+2*Y-2)/(2*Y-2)/(2*Y-1)/((Z-DX/2)^2)
2780 next Y
2790 D=C(0)*J(0)*G(0)
2800 for Y=1 to B
2810 D=D+C(Y)*J(Y)*G(Y)
2820 next Y
2830 R=-D/((Z-DX/2)^N) 'The series remainder term is calculated
2840 if FF=1 then 'If FF=1, recalculate the series remainder term
2850 :R=R+DX/((Z-DX)^N) 'The series remainder term is recalculated
2860 :Z=SZ 'Restore the value of Z
2870 :FF=0 'Reset the division by O avoidance flag to 0"
2880 :endif
2890 P=TS+R 'Calculate LND(N,DX,Z) Z=inputted subroutine value
2900 FG=0 'Assure that the "formula only" flag is turned off
2910 return
2920 print "The truncated series loop count has exceeded ";M;" counts"
2930 goto 3690
2940 'Subroutine to find an X locus point with a nearly matching real
     value
2950 'Input is ZI and the output is XO
2960 if re(DX)=0 then XO=X+(fix((ZI-im(X))/im(DX)))*DX
2970 :else XO=X+(fix((ZI-re(X))/re(DX)))*DX
2980 :endif
2990 return
3000 'CALCULATION OF ln(DX,X)=LND(1,DX,X)
3010 'ln(DX,X)=ln1(R) where R=X/DX
3020 Tk=3 'ln(DX,X)=LND(1,DX,X) Calculation Method #3
3030 R=X/DX
3040 PI=#pi
3050 'Identify the sign of Re(R)
```

```
3060 if re(R)>=0 then Cn=6:goto 3180
3070 :else Cn=7
3080 :if and{im(R)=0,abs(R-fix(R))<10.0^(-24)} then T=0
3090 :elseif im(PI*R)>10^(2.5) then T= +PI*#i 'for very large positive PI*R
3100 :elseif im(PI*R)<-10^(2.5) then T= -PI*#i 'for very large negative
     PI*R
3110 :else T=-PI/tan(PI*R)
3120 :endif
3130 :endif
3140 Z=1-R 'Re(R)<0
3150 gosub 3220 'Find ln1(1-R)
3160 Y=P+T 'Calculate ln1(R) 'If R is an integer, T=0
3170 goto 3530
3180 Z=R 'Re(R)>0
3190 gosub 3220 'Calculate ln1(R)
3200 Y=P
3210 goto 3530
3220 'CALCULATE ln(DX,X)=LND(1,DX,X) -- Subroutine starts here
3230 G=20000 'Number selected for 44 place accuracy for X<G, G>0
3240 K=#euler
3250 if Z=1 then P=0:goto 3520 'Special case
3260 L=G-abs(Z) 'Is the magnitude of Z less than or equal to G
3270 if re(Z)<0 then goto 3300 'If re(Z) is negative
3280 if L<=0 then gosub 3440:goto 3520 'abs(Z) is greater than
     or equal to G
3290 'abs(Z) is less than G
3300 FX=re(Z)-fix(re(Z)) 'Get the fractional part of the real part of Z
3310 HX=G+FX+im(Z)*#i 'High Z selected for formula desired accuracy
3320 W=G-fix(re(Z)) 'Calculate the count down number
3330 Q=HX
3340 SZ=Z 'Save the value of Z
3350 Z=HX
```

```
3360 gosub 3440
3370 Z=SZ 'Restore the value of Z
3380 for I=1 to W 'Count down to obtain the desired LN1(Z)
3390 Q=Q-1
3400 if Q=0 then goto 3420:endif '1/Q=1/(+-0) is defined as 1/Q=0
3410 P=P-1/Q
3420 next I
3430 goto 3520
3440 'ln(DX,X)=LND(1,DX,X) Formula Subroutine
3450 P=log(Z-0.5)+K 'Formula to calculate ln1(Z), Re(Z)>0 and Z large
3460 S=0 'Calculate ln1(Z) series terms
3470 for M=1 to B
3480 S=S+C(M)/(2*M+1)/(2*M)/((2*Z-1)^(2*M))
3490 next M
3500 P=P+S
3510 return
3520 return
3530 if CT=0 then Cn1=Cn else Cn2=Cn 'For calculation path checking
3540 IY=alen(int(abs(Y)))+2 'Specify the number of integer digits to print
3550 print "ln(";DX;",";X;")= ";using(IY,A),Y 'Print to screen
3560 'lprint "ln(";DX;",";X;")= ";using(IY,A),Y 'Have printer print
3570 goto 3690
3580 if CT=0 then Cn1=Cn else Cn2=Cn 'For calculation path checking
3590 IE=alen(int(abs(E)))+2 'Specify the number of integer digits to print
3600 print "LND(";N;",";DX;",";X;") = ";using(IE,A),E 'Print to screen
3610 'lprint "LND(";N;",";DX;",";X;")= ";using(IE,A),E 'Have printer print
3620 if re(N)<=0 then goto 3690 'Don't print K
3630 if or{re(DX)>0,and{re(DX)=0,im(DX)>0}} then FI=" -inf. "
3640 :else FI=" +inf. "
3650 :endif
3660 IK=alen(int(abs(K)))+2 'Specify the number of integer digits to print
3670 print "LND(";N;",";DX;",";FI;") = ";using(IK,A),K 'Print to screen
```

```
3680 'lprint "LND(";N;",";DX;",";FI;") = ";using(IK,A),K 'Have printer print
3690 if and{re(N)=1,im(N)=0} then SUM=SUM+Y*(-1)^CT else SUM=SUM+E*(-1)^CT
3700 if CT=1 then goto 3850
3710 print
3720 strinput "Do you wish to enter another value of X (Y/N) ";AA
3730 CT=CT+1 'Count the number of X entries
3740 if or{asc(AA)=asc("y"),asc(AA)=asc("Y")} then goto 3750
     else goto 3980
3750 NF=N 'Save the first value of N
3760 DXF=DX 'Save the first value of DX
3770 XF=X 'Save the first value of X
3780 print
3790 if FB=0 then goto 3820
3800 input "What is N ";N 'Input the second value of N
3810 input "What is DX ";DX 'Input the second value of DX
3820 input "What is X ";X 'Input the second value of X
3830 print
3840 if and{re(N)=1,im(N)=0} then goto 3000 else goto 600
3850 print
3860 if and{re(N)=1,im(N)=0} then
3870 :SUM=-SUM
3880 :IS=alen(int(abs(SUM)))+2 'Specify the no. of integer digits to print
3890 :print "ln(";DX;",";X;")-ln(";DXF;",";XF;") = ";using(IS,A),SUM
3900 :else
3910 :IS=alen(int(abs(SUM)))+2 'Specify the no. of integer digits to print
3920 :print "-LND(";N;",";DX;",";X;")+LND(";NF;",";DXF;",";XF;") =
     ";using(IS,A),SUM
3930 :endif
3940 print
3950 if and{N=NF,DX=DXF} then goto 3960 else goto 3980 ' print SUM?
3960 ISD=alen(int(abs(SUM/DX)))+2 'Specify the no. integer digits to print
3970 print "SUM(";N;",";DX;",";XF;"to ";X-DX;") = ";using(ISD,A),SUM/DX
3980 end
```

## 2) UBASIC Program to Calculate the Cn, n=0,1,2,… Constants
### "CNCALC.UB"

```
10  'This program calculates the Cn constants from the Bernoulli constants
20  dim A(60),B(60),C(60),D(60)
30  I=30 'I is the number of Cn constants to be calculated
40  'Calculate the even Bernoulli constants, Bm, m=2,4,6,...,2I
50  D(0)=1
60  D(1)=-1//2
70  D(2)=1//12
80  B(2)=D(2)*!(2)
90  for N=2 to I
100  E1=-1//!(2*N+1)+1//2//!(2*N)
110  SUM=0
120  for M=1 to N-1
130  SUM=SUM+D(2*M)//!(2*N+1-2*M)
140  next M
150  E2=SUM
160  D(2*N)=(E1-E2)
170  B(2*N)=D(2*N)*!(2*N)
180  next N
190  'Calculate the Cn constants n=0,1,2,3,...,I
200  C(0)=-1
210  print
220  print "C 0 = ";C(0)
230  A(1)=0
240  for N=1 to I
250  S=0
260  if N=1 then goto 310
270  for M=1 to N-1
280  S=S+C(M)*!(2*N-1)//!(2*N-2*M)//!(2*M+1)
290  next M
```

```
300 A(N)=S
310 D=(1-(2^(2*N))*B(2*N))//(2*N)
320 C(N)=(2*N+1)*(2*N)*(D-A(N))
330 print "C";N;"= ";C(N)
340 next N
350 end
```

# 3) UBASIC Program to Calculate the An, n=1,2,3,… Constants

## "ANCALC.UB"

```
10  'This program calculates the An constants from the Cn constants
12  point 20
14  word 280
20  dim A(60),B(60),C(60),D(60)
30  I=17 'I is the number of An constants to be calculated
40  'Calculate the even Bernoulli constants, Bm, m=2,4,6,...,2I
50  D(0)=1
60  D(1)=-1//2
70  D(2)=1//12
80  B(2)=D(2)*!(2)
90  for N=2 to I
100  E1=-1//!(2*N+1)+1//2//!(2*N)
110  SUM=0
120  for M=1 to N-1
130  SUM=SUM+D(2*M)//!(2*N+1-2*M)
140  next M
150  E2=SUM
160  D(2*N)=(E1-E2)
170  B(2*N)=D(2*N)*!(2*N)
180  next N
190  'Calculate the Cn constants n=0,1,2,3,...,I
200  C(0)=-1
210  print
220  print "C 0 = ";C(0)
230  A(1)=0
240  for N=1 to I
250  S=0
260  if N=1 then goto 310
270  for M=1 to N-1
```

```
280 S=S+C(M)*!(2*N-1)//!(2*N-2*M)//!(2*M+1)
290 next M
300 A(N)=S
310 D=(1-(2^(2*N))*B(2*N))//(2*N)
320 C(N)=(2*N+1)*(2*N)*(D-A(N))
330 print "A";N;"= ";-C(N)//2//!(2*N+1)
340 next N
350 end
```

# 4) UBASIC Program to Calculate the Bn, n=1,2,3,… Constants

"BNCALC.UB"

```
10  'This program calculates the Bn constants from the Cn constants
12  point 20
14  word 280
20  dim A(60),B(60),C(60),D(60)
30  I=17 'I is the number of An constants to be calculated
40  'Calculate the even Bernoulli constants, Bm, m=2,4,6,...,2I
50  D(0)=1
60  D(1)=-1//2
70  D(2)=1//12
80  B(2)=D(2)*!(2)
90  for N=2 to I
100  E1=-1//!(2*N+1)+1//2//!(2*N)
110  SUM=0
120  for M=1 to N-1
130  SUM=SUM+D(2*M)//!(2*N+1-2*M)
140  next M
150  E2=SUM
160  D(2*N)=(E1-E2)
170  B(2*N)=D(2*N)*!(2*N)
180  next N
190  'Calculate the Cn constants n=0,1,2,3,...,I
200  C(0)=-1
210  print
230  A(1)=0
240  for N=1 to I
250  S=0
260  if N=1 then goto 310
270  for M=1 to N-1
280  S=S+C(M)*!(2*N-1)//!(2*N-2*M)//!(2*M+1)
```

```
290 next M
300 A(N)=S
310 D=(1-(2^(2*N))*B(2*N))//(2*N)
320 C(N)=(2*N+1)*(2*N)*(D-A(N))
330 print "B";N;"= ";-C(N)//(2^(2*N))//!(2*N+1)
340 next N
350 end
```

# 5) UBASIC Program to Calculate the Hn, n=1,2,3,… Constants

## "HNCALC.UB"

```
10  'This program calculates the Hn constants from the Bernoulli constants
20  dim A(60),B(60),C(60),H(60)
30  B(2)=1//6
40  B(4)=-1//30
50  B(6)=1//42
60  B(8)=-1//30
70  B(10)=5//66
80  B(12)=-691//2730
90  B(14)=7//6
100  B(16)=-3617//510
110  B(18)=43867//798
120  B(20)=-174611//330
130  B(22)=854513//138
140  B(24)=-236364091//2730
150  B(26)=8553103//6
160  B(28)=-23749461029//870
170  B(30)=8615841276005//14322
180  B(32)=-7709321041217//510
190  B(34)=2577687858367//6
200  B(36)=-26315271553053477373//1919190
210  B(38)=2929993913841559//6
220  B(40)=-261082718496449122051//13530
230  B(42)=1520097643918070802691//1806
240  B(44)=-27833269579301024235023//690
250  B(46)=596451111593912163277961//282
260  B(48)=-5609403368997817686249127547//46410
270  B(50)=495057205241079648212477525//66
280  C(0)=-1
```

```
290 print
300 print "C 0 = ";C(0)
310 A(1)=0
320 for N=1 to 25
330 S=0
340 if N=1 then goto 390
350 for M=0 to N-2
360 S=S+C(M+1)*!(2*N-1)//!(2*N-2-2*M)//!(3+2*M)
370 next M
380 A(N)=S
390 D=(1-(2^(2*N))*B(2*N))//(2*N)
400 C(N)=(2*N+1)*(2*N)*(D-A(N))
410 print "C";N;"= ";C(N)
420 next N
430 print:print
440 print "Hn Constants"
450 print
460 for I=1 to 17
470 H(I)=(B(2*I)+C(I)//(2*I+1)//2^(2*I))//!(2*I)
480 print "H";I;"= ";H(I)
490 next I
500 end
```

# 6) UBASIC Program to Plot a Nyquist Diagram and Calculate Phase Margin

"PLPM.UB"

```
10 'EXTENDED NYQUIST CRITERION COMPLEX PLANE POLAR PLOTTING PROGRAM
20 P=20 '% of left half plane circle calculated
30 D=5 '-D to +D is the x axis range, -.7D to +.7D is the y axis range
40 DX=0.1 'Delta X increment
50 N=20000 'N is one half of the number of points calculated
60 SN=0.00000000000001 'sn is a small number, 10^(-14)
70 C=0 'If C=1 no longer search for a curve negative x axis crossover
80 Q=0 'Number of points calculation counter initialized to 0
90 ST=0 'ST=1 if there is a curve negative x axis crossover x<-1
100 gosub *AXES(-D,D,-0.7*D,0.7001*D,1,30,"polar","white","Nyquist Plot")
110 for M=N to -N step -1
120 H="line":R=2
130   X=(#pi*M/N)*P/100 'Calculate S around the left half plane circle
140   S=(-1+exp(X*#i))/DX
150 '----------------FUNCTION OF S TO PLOT-----------------------
160 F=1.7*((1+S*DX)^3)/(S+1)/(S+0.7)/(S+SN)
170 '-----------------------------------------------------------
180 'Identify a curve crossover of the negative x axis where x<-1
190 V1=arg(F) 'Present value of the argument of the function, F
200 'V2 is the past value of V1, the argument of the function, F
210 Q=Q+1 'The number of points calculation counter
220 if Q=1 then V2=V1 'For the first point calculated, V2=V1
230 E=abs(V1-V2) 'The value of E spikes at an x axis crossover
240 V2=V1 'V2 is given the value of V1 to be the V1 previous value
250 'Calculate Phase Margin for a negative x axis crossover -1<=x<=0
```

```
260 if and{E>=1,re(F)<-1,im(F)<0.1,im(F)>-0.1} then ST=1
270 if C=1 goto 340
280 if and{abs(F)>=1,arg(F)<=0,ST<>1} then
290 :H="both":R=4:C=1 'Place a red dot where abs(F)=1
300 :G=fix(100*(#pi+arg(F))*180/#pi)/100 'Calculate Phase Margin
310 :W#=str(G)+" deg"
320 :glocate mapx(-D)+3,mapy(0.7*D)+35:gprint W#; 'Print Phase Margin
330 :endif
340   gosub *PLOT(abs(F),arg(F),"polar",H,R)
350 next M
360 end
370 '========== AXES SUBROUTINE ==========
380 'X1,X2: Lower and upper limit of x-axis
390 'Y1,Y2: Lower and upper linit of y-axis
400 'XTIC,YTIC: tick mark spacing.  No tics if <=0
410 '   if polar: XTIC refers to R-axis & YTIC is in degrees
420 S%:"r"=RectangularGrid,"p"=PolarGrid
430 'IBACK#: "b"=black background, "w"=white background
440 'TITLE#: Plot title
450 '
460 *AXES(X1,X2,Y1,Y2,XTIC,YTIC,IAXES#,IBACK#,TITLE#)
470 screen 23
480 gsize:if (left(IBACK#)="w" or left(IBACK#)="W") then gcolor -1 else
    gcolor -7
490 cls 3
500 ' Set range
510 Vxh=630:Vx1=320-Vxh\2:Vx2=Vx1+Vxh
520 Vyh=440:Vy1=240-Vyh\2:Vy2=Vy1+Vyh
530 XCSCALE=(X2-X1)*Vyh/Vxh:YCSCALE=Y2-Y1
540 ' Set view port and window.
550 if (left(IBACK#)="w" or left(IBACK#)="W") then view (Vx1,Vy1)-
    (Vx2,Vy2),15,7 else view (Vx1,Vy1)-(Vx2,Vy2),0,7
```

```
560  window (X1,Y2)-(X2,Y1)
570  if (left(IAXES#)="p" or left(IAXES#)="P") then goto 640
580  ' Rectangular grid
590   if (YTIC>0) then for I=1 to -int(Y1/YTIC):for J=0 to 100:pset
      (X1+(X2-X1)*J/100,-I*YTIC),8:next J:next I
600   if (YTIC>0) then for I=1 to int(Y2/YTIC):for J=0 to 100:pset
      (X1+(X2-X1)*J/100,I*YTIC),8:next J:next I
610   if (XTIC>0) then for I=1 to -int(X1/XTIC):for J=0 to 100:pset (-
      I*XTIC,Y1+(Y2-Y1)*J/100),8:next J:next I
620   if (XTIC>0) then for I=1 to int(X2/XTIC):for J=0 to 100:pset
      (I*XTIC,Y1+(Y2-Y1)*J/100),8:next J:next I
630  goto 680
640  ' Polar grid
650  RMAX=sqrt(max(X1^2+Y1^2,X2^2+Y1^2,X1^2+Y2^2,X2^2+Y2^2))
660  if (XTIC>0) then for I=1 to int(RMAX/XTIC):for J=0 to
     120:psetMax(X1,min(X2,I*XTIC*cos(J*2*#pi/120))),
     max(Y1,min(Y2,I*XTIC*sin(J*2*#pi/120)))),8:next J:next I
670  if (YTIC>0) then for I=1 to int(360/YTIC):for J=0 to 100:pset
     (max(X1,min(X2,J*RMAX/100*cos(I*#pi*YTIC/180))),
     max(Y1,min(Y2,J*RMAX/100*sin(I*#pi*YTIC/180)))),8:next J:next I
680  ' Draw X- and Y- axes.
690  if (Y1<=0 and Y2>=0) then line (X1,0)-(X2,0),1
700  if (X1<=0 and X2>=0) then line (0,Y1)-(0,Y2),1
710  'Write scales
720  glocate mapx(X1)+3,mapy((Y1+Y2)/2):gprint int(X1*10)/10;
730  glocate mapx(X2)-8*len(str(int(X2*10)/10)),mapy((Y1+Y2)/2):
     gprint int(X2*10)/10;
740  glocate mapx((X1+X2)/2),mapy(Y1)-15:gprint int(Y1*10)/10;
750  glocate mapx((X1+X2)/2),mapy(Y2)+2:gprint int(Y2*10)/10;
760  ' Display the title
770  glocate mapx(X1)+3,mapy(Y2)+2:gprint TITLE#;
780  return
```

```
790 '========= PLOT SUBROUTINE =========
800 'X=X coordinate (or r in polar)
810 'Y=Y coordinate (or theta in polar)
820 'IPOLAR: "r"=use rectangular coordinates, "p"=use polar coordinates
830 'ILINE:  "p"=use points, "l"=use lines, "b"=use both
840 'ICOLOR: COLOR (0-15)
850 '
860 *PLOT(X,Y,IPOLAR#,ILINE#,ICOLOR%)
870 local X1,Y1
880 X1=X:Y1=Y
890 if (left(IPOLAR#)="p" or left(IPOLAR#)="P") then
    X=X1*cos(Y1):Y=X1*sin(Y1)
900 if (left(ILINE#)="p" or left(ILINE#)="P") then ILINE#=0:goto 940
910 if (left(ILINE#)="l" or left(ILINE#)="L") then ILINE#=1:goto 940
920 if (left(ILINE#)="b" or left(ILINE#)="B") then ILINE#=2:goto 940
930 ILINE#=1
940 if (FIRSTPOINT%=1 and (ILINE#=1 or ILINE#=2)) then line –
    (X,Y),ICOLOR%
950 if (FIRSTPOINT%=0 and (ILINE#=1 or ILINE#=2)) then pset
    (X,Y),ICOLOR%:FIRSTPOINT%=1
960 if (ILINE#=0 or ILINE#>=2) then circle
    (X,Y),(XCSCALE/200,YCSCALE/200),ICOLOR%,"f",ICOLOR%
970 if (ILINE#=0) then FIRSTPOINT%=0
980 return
```

# 7) UBASIC Program to Plot a Difference Equation

"PLTX.UB"

```
 10  'DIFFERENCE EQUATION PLOTTING PROGRAM
 20  DX=0.1 ' The Delta X increment
 30  F0=0 'The initial value of y0, y(0)
 40  F1=0 'The initial value of y1, Y(DX)
 50  F2=0 'The initial value of y2, Y(2DX)
 60  XL=0 'Leftmost x value
 70  XR=50 'Rightmost x value
 80  YL=-10 'Lowest y value
 90  YU=10 'Uppermost y value
100  R=5 'Step input value
110  A=1 'A K Transform pole value, (S+A)
120  B=0.7 'A K Transform pole value, (S+B)
130  P=fix(XR/DX) 'The number of points to plot
140  K=1.7 'K Transform gain constant
150  gosub *AXES(XL,XR,YL,YU,10,10,"rectangle","white","Function Plot")
160  for M=0 to P
170    X=M*DX
180  FT=K*(DX^3)*R-A*B*(DX^2)*(F1-F0)-(A+B)*DX*(F2-2*F1+F0)
190  F3=(FT-(-3*F2+3*F1-F0))/(1+K*(DX^3))
200    gosub *PLOT(X,F0,"rectangle","line",2)
210  F0=F1
220  F1=F2
230  F2=F3
240  next M
250  end
260  '========== AXES SUBROUTINE ==========
270  'X1,X2: Lower and upper limit of x-axis
280  'Y1,Y2: Lower and upper linit of y-axis
290  'XTIC,YTIC: tick mark spacing.  No tics if <=0
300  '   if polar: XTIC refers to R-axis & YTIC is in degrees
310  'IAXES%: "r"=rectangular grid, "p"=polar grid
320  'IBACK#: "b"=black background, "w"=white background
330  'TITLE#: Plot title
340  '
350  *AXES(X1,X2,Y1,Y2,XTIC,YTIC,IAXES#,IBACK#,TITLE#)
360  screen 23
370  gsize:if (left(IBACK#)="w" or left(IBACK#)="W") then gcolor -1 else
     gcolor -7
380  cls 3
390  ' Set range
400  Vxh=630:Vx1=320-Vxh\2:Vx2=Vx1+Vxh
```

```
410  Vyh=440:Vy1=240-Vyh\2:Vy2=Vy1+Vyh
420  XCSCALE=(X2-X1)*Vyh/Vxh:YCSCALE=Y2-Y1
430  ' Set view port and window.
440  if (left(IBACK#)="w" or left(IBACK#)="W") then view (Vx1,Vy1)-
     (Vx2,Vy2),15,7 else view (Vx1,Vy1)-(Vx2,Vy2),0,7
450  window (X1,Y2)-(X2,Y1)
460  if (left(IAXES#)="p" or left(IAXES#)="P") then goto 530
470  ' Rectangular grid
480   if (YTIC>0) then for I=1 to -int(Y1/YTIC):for J=0 to 100:pset
     (X1+(X2-X1)*J/100,-I*YTIC),8:next J:next I
490   if (YTIC>0) then for I=1 to int(Y2/YTIC):for J=0 to 100:pset
      (X1+(X2-X1)*J/100,I*YTIC),8:next J:next I
500   if (XTIC>0) then for I=1 to -int(X1/XTIC):for J=0 to 100:pset
      (-I*XTIC,Y1+(Y2-Y1)*J/100),8:next J:next I
510   if (XTIC>0) then for I=1 to int(X2/XTIC):for J=0 to 100:pset
      (I*XTIC,Y1+(Y2-Y1)*J/100),8:next J:next I
520   goto 570
530  ' Polar grid
540  RMAX=sqrt(max(X1^2+Y1^2,X2^2+Y1^2,X1^2+Y2^2,X2^2+Y2^2))
550  if (XTIC>0) then for I=1 to int(RMAX/XTIC):for J=0 to 120:pset
     (max(X1,min(X2,I*XTIC*cos(J*2*#pi/120))),
     max(Y1,min(Y2,I*XTIC*sin(J*2*#pi/120)))),8:next J:next I
560  if (YTIC>0) then for I=1 to int(360/YTIC):for J=0 to 100:pset
     (max(X1,min(X2,J*RMAX/100*cos(I*#pi*YTIC/180))),
     max(Y1,min(Y2,J*RMAX/100*sin(I*#pi*YTIC/180)))),8:next J:next I
570  ' Draw X- and Y- axes.
580  if (Y1<=0 and Y2>=0) then line (X1,0)-(X2,0),1
590  if (X1<=0 and X2>=0) then line (0,Y1)-(0,Y2),1
600  'Write scales
610  glocate mapx(X1)+3,mapy((Y1+Y2)/2):gprint int(X1*10)/10;
620  glocate mapx(X2)-8*len(str(int(X2*10)/10)),mapy((Y1+Y2)/2):gprint
     int(X2*10)/10;
630  glocate mapx(X1),mapy(Y1)-15:gprint int(Y1*10)/10;
640  glocate mapx(X2),mapy(Y2)+2:gprint int(Y2*10)/10;
650  ' Display the title
660  glocate mapx(X1)+3,mapy(Y2)+2:gprint TITLE#;
670  return
680  '========== PLOT SUBROUTINE ==========
690  'X=X coordinate (or r in polar)
700  'Y=Y coordinate (or theta in polar)
710  'IPOLAR: "r"=use rectangular coordinates, "p"=use polar coordinates
720  'ILINE:  "p"=use points, "l"=use lines, "b"=use both
730  'ICOLOR: COLOR (0-15)
740  '
750  *PLOT(X,Y,IPOLAR#,ILINE#,ICOLOR%)
760  local X1,Y1
```

```
770  X1=X:Y1=Y
780  if (left(IPOLAR#)="p" or left(IPOLAR#)="P") then
     X=X1*cos(Y1):Y=X1*sin(Y1)
790  if (left(ILINE#)="p" or left(ILINE#)="P") then ILINE#=0:goto 830
800  if (left(ILINE#)="l" or left(ILINE#)="L") then ILINE#=1:goto 830
810  if (left(ILINE#)="b" or left(ILINE#)="B") then ILINE#=2:goto 830
820  ILINE#=1
830  if (FIRSTPOINT%=1 and (ILINE#=1 or ILINE#=2)) then line –
     (X,Y),ICOLOR%
840  if (FIRSTPOINT%=0 and (ILINE#=1 or ILINE#=2)) then pset
     (X,Y),ICOLOR%:FIRSTPOINT%=1
850  if (ILINE#=0 or ILINE#>=2) then circle
     (X,Y),(XCSCALE/200,YCSCALE/200),ICOLOR%,"f",ICOLOR%
860  if (ILINE#=0) then FIRSTPOINT%=0
870  return
```

# Explanation of the use of the lnd(n,Δx,x) function calculation program

# Lnd(n,Δx,x) Function Calculation Program

## (Program Disk is provided)

The envelope on the back binder cover contains a disk with the computer programming language, UBASIC, on it. UBASIC runs on the MS-DOS operating system. This disk also contains the lnd(n,Δx,x) function calculation program, LNDX. This program is written in the UBASIC programming language.

To run the lnd(n,Δx,x) function calculation program, do the following:

Call up the LNDX program in MS-DOS in the following manner:

| | |
|---|---|
| 1. C:\>E: push enter key | select the computer disk (E or other drive input) |
| 2. E:> appears | |
| 3. E:>UBASIC push enter key | select the UBASIC program |
| 4. OK | appears |
| 5. push the F1 key | to load a UBASIC program |
| 6. load " | appears |
| 7. LNDX push the enter key | select the LNDX UBASIC program |
| 8. run push the enter key | enter run to run the LNDX UBASIC program |
| 9. n push enter key<br>Δx push enter key<br>x push enter key | enter n,Δx,x    <u>Note</u> – To enter a complex number, 2+3i, enter 2+3#i |
| 10. lnd(n,Δx,x) | appears |
| 11. yes/no | enter yes for integrals with two limits, the first entry is the low limit |
| 12. x push the enter key | if yes, the high integration limit is entered |
| 13. integration calculation<br>sum calculation | if yes, appears |
| 14. system push the enter key | After a calculation is run, to exit to MS-DOS |
| 15. E:\> exit push the enter key | To exit MS-DOS |

# Simple examples to show how to use the lnd(n,Δx,x) function calculation program, LNDX

## Primary Equations:

$$\int_{\Delta x\; x_1}^{x_2} \frac{1}{(x-a)^n}\,\Delta x = \Delta x \sum_{\Delta x\; x=x_1}^{x_2-\Delta x} \frac{1}{(x-a)^n} = \pm\, \mathrm{lnd}(n,\Delta x,(x-a))\,\Big|_{x_1}^{x_2} \quad + \text{ for } n = 1,\ - \text{ for } n \neq 1 \qquad 1)$$

Often used to integrate or sum partial fraction expansions

$$\int_{\Delta x\; x_i}^{\infty} \frac{1}{x^n}\,\Delta x = \Delta x \sum_{\Delta x\; x=x_i}^{\infty} \frac{1}{x^n} = \mathrm{lnd}(n,\Delta x,x_i) \quad \text{for } n \neq 1 \qquad 2)$$

Often used to calculate the Riemann Zeta Function (Re(n)>1, Δx=1, x=1) and the Hurwitz Zeta function (Re(n)>1, Δx=1)

Note -- All variables can be real or complex values

**Example #1**
**Calculation using Eq 1**

**Simple Calculations**

**Left to right integration/summation**

$n = 1$, $\Delta x = .5$ $x_1 = -1$, $x_2 = 1.5$, $x_2-\Delta x = 1$

$$\int_{.5\; -1}^{1.5} \frac{1}{x}\,\Delta x = .5 \sum_{.5\; x=-1}^{1} \frac{1}{x} = \mathrm{lnd}(1,.5,1.5) - \mathrm{lnd}(1,.5,-1)$$

$$\sum_{.5\; x=-1}^{1} \frac{1}{x} = 2(1.5 - 1.5) = 0$$

Checking

$$\frac{1}{-1} + \frac{1}{-.5} + \frac{1}{0} + \frac{1}{.5} + \frac{1}{1} = 0$$

Good check

Note - As with the Hurwitz Zeta Funcion, any division by zero summation term is excluded.

**Right to left integration/summation**

n = 1, $\Delta x$ = -.5 $x_1$ = 1, $x_2$ = -1.5, $x_2$-$\Delta x$ = -1

$$ {}_{-.5}\int_{1}^{-1.5} \frac{1}{x}\,\Delta x = -.5\ {}_{-.5}\sum_{x=1}^{-1} \frac{1}{x} = \text{lnd}(1,-.5,-1.5) - \text{lnd}(1,-.5,1) $$

$$ {}_{-.5}\sum_{x=1}^{-1} \frac{1}{x} = -2(1.5 - 1.5) = 0 $$

Checking

$$ \frac{1}{1} + \frac{1}{.5} + \frac{1}{0} + \frac{1}{-.5} + \frac{1}{-1} = 0 $$

Good check

<u>Note</u> - As with the Hurwitz Zeta Funcion, any division by zero summation term is excluded.

**Example #2**
**Calculation using Eq 1**

n = 2.3+3.7i, $\Delta x$ = 1+i, $x_1$ = -1-3i, $x_2$ = 3+i, $x_2$-$\Delta x$ = 2

$$ {}_{1+i}\int_{-1-3i}^{3+i} \frac{1}{(x-2+i)^{2.3+3.7i}}\,\Delta x = (1+i)\ {}_{1+i}\sum_{x=-1-3i}^{2} \frac{1}{(x-2+i)^{2.3+3.7i}} = -\text{lnd}(2.3+3.7i,1+i,x-2+i)\Big|_{-1-3i}^{3+i} $$

$$ {}_{1+i}\sum_{x=-1-3i}^{2} \frac{1}{(x-2+i)^{2.3+3.7i}} = \frac{1}{1+i}\,[-\text{lnd}(2.3+3.7i,1+i,1+2i) + \text{lnd}(2.3+3.7i,1+i,-3-2i)] $$

$$ {}_{1+i}\sum_{x=-1-3i}^{2} \frac{1}{(x-2+i)^{2.3+3.7i}} = \frac{1}{1+i}\,[(-3.0442574254041-12.2176183323067i)+ $$

$$ (+155626.32956516287-24853.3785053376592i)] $$

$$ {}_{1+i}\sum_{x=-1-3i}^{2} \frac{1}{(x-2+i)^{2.3+3.7i}} = 65378.8445920337525-90244.4407157037186i $$

Checking

$$ {}_{1+i}\sum_{x=-1-3i}^{2} \frac{1}{(x-2+i)^{2.3+3.7i}} = \frac{1}{(-3-2i)^{2.3+3.7i}} + \frac{1}{(-2-i)^{2.3+3.7i}} + \frac{1}{(-1)^{2.3+3.7i}} + \frac{1}{(i)^{2.3+3.7i}} $$

$$ {}_{1+i}\sum_{x=-1-3i}^{2} \frac{1}{(x-2+i)^{2.3+3.7i}} = 65378.8445920337525-90244.4407157037186i $$

Good check

**Calculation using Eq 2**

$$ {}_{1}\int_{1}^{\infty} \frac{1}{x^2}\,\Delta x = {}_{1}\sum_{x=1}^{\infty} \frac{1}{x^{2}} = \operatorname{lnd}(2,1,1) = \zeta(2) = 1.644934066848226 $$

Checking

From The Handbook of Mathematical Tables

$$ {}_{1}\sum_{x=1}^{\infty} \frac{1}{x^{2}} = 1.644934066848226 \qquad \text{Good check} $$

# Acknowledgements

I have finally been able to complete my work due to some much needed advice from Professor Victor Kac of the Department of Mathematics at MIT and Professor Delfim Torres of the Department of Mathematics at the University of Aveiro in Portugal. I especially thank Professor Torres for endorsing this document for submission to the arXiv. Finally, I would like to thank my son, William Kaminsky, a doctoral student in the Department of Physics at MIT, for help in editing this document.

# Bibliographic Notes

The "Interval Calculus" put forth in this document was initially conceived and largely developed independently of the mathematics and engineering research communities. It was only toward the very end of the writing of this document that the author sought to correspond with mathematicians and academic engineers (please see the above Acknowledgements) as well as to review the contemporary literature to see what parallels exist to Interval Calculus.

There exist two strands of parallel research in the contemporary literature. These two strands are most often called "Quantum Calculus" (Kac and Cheung 2002) and "Time Scale Calculus" (Bohner and Peterson 2001), respectively. Whereas the Interval Calculus in this document deals exclusively with analysis on discrete domains of real or complex numbers separated by a fixed distance $\Delta x$, both Quantum Calculus and Time Scale Calculus deal with analysis on substantially more general domains. Quantum Calculus not only contains an "h-calculus" that considers domains that are additively defined like Interval Calculus --- that is, domains where if x[n] is the nth item, then x[n] = x[0] + nh --- but it also contains a "q-calculus" that considers domains that are multiplicatively defined --- that is, if x[n] denotes the nth item in the domain, then $x[n] = q^n x[0]$. Superficially, the h-calculus and the q-calculus appear quite different. However, the fact that one can map any bounded interval of the real line to the unit circle of the complex plane means that they are closely connected via q = exp(ih). Time Scale Calculus is even more general than Quantum Calculus, considering domains that are arbitrary closed subsets of the real numbers.

Despite the fact that both Quantum Calculus and Time Scale Calculus in principle subsume the Interval Calculus described in this document, the author believes that this Interval Calculus in practice has substantial value. Interval Calculus focuses on the most commonly encountered real world discrete variable problems --- that is, those with a fixed discretization. Moreover, the present document aims to be useful both as a pedagogical textbook for the advanced undergraduate coming from the standard engineering curriculum as well as a reference handbook for the working engineer. More specifically, the first 8 chapters of the present document are loosely based on Wylie and Barrett's standard survey textbook of advanced engineering mathematics (Wylie and Barrett 1995), and the appendix of the present document is loosely based on such standard reference tomes as Abramowitz and Stegun (1965) and the *CRC Standard Mathematical Tables and Formulae* (Zwillinger 2011).

An Annotated List of Key References for Quantum and Time Scales Calculus

All three calculi --- Interval, Quantum, and Time Scales --- have antecedents in the classical literature on the Calculus of Finite Differences (e.g., Jordan 1965 and references therein). Indeed, Baez (2002) asserts that Quantum Calculus in both its h-calculus and q-calculus varieties goes back to the infancy of rigorous calculus, at least as far back as Gauss. However, since such classical literature is overwhelmingly concerned with deriving series approximations and interpolation formulas for continuous calculus expressions as opposed to making a unified framework for calculus on both continuous and discrete structures, we shall not delve any deeper into it.

Time Scales Calculus is generally credited to have started with the Ph.D. thesis of Stefan Hilger (1988), where it was called "Analysis on Measure Chains." The seminal papers in the Time Scales literature on Laplace transformations are Hilger (1999) and Bohner and Peterson (2002). Except for notation, the $K_{\Delta t}$ transformation of Interval Calculus described in the present document is identical to the Laplace transformations of the latter reference (Bohner and Peterson 2002) when those are specialized to the case of real or complex numbers with a fixed discretization.

Both Quantum Calculus and Time Scales Calculus have garnered textbook-length introductions geared to undergraduates comfortable with rigorous real analysis at the level of a first semester course for mathematics majors. In the case of Quantum Calculus, this textbook is Victor Kac and Pokmon Cheung's *Quantum Calculus* (2002). In the case of Time Scales Calculus, this textbook is Martin Bohner and Allan Peterson's *Dynamic Equations on Time Scales: An Introduction with Applications* (2001). Beyond the previously discussed basic difference among Interval, Quantum, and Time Scales Calculi in the types of discrete domains they treat, some notable contrasts between these textbooks and the present document are:

1) The present document does not have a semester of rigorous real analysis as a prerequisite.

2) The present document treats both the analytical and especially the numerical evaluation of discrete integrals in much greater depth than either *Quantum Calculus* or *Dynamic Equations on Time Scales*.

3) The present document along with *Quantum Calculus* thoroughly treats Taylor approximation and the evaluation of the remainder in finite order Taylor approximations. (*Quantum Calculus* presents this via Bernoulli polynomials whereas the present document presents this via another family of polynomials that was found to be more natural.) *Dynamic Equations on Time Scales* does not treat the remainders to Taylor approximations.

4) The present document thoroughly treats discrete Laplace transforms through the lens of control systems engineering. *Quantum Calculus* does not present any Laplace transforms. In contrast, *Dynamic Equations on Time Scales* does discuss Laplace transforms as well as linear dynamical systems. However, it does not have the present document's emphasis on linear control theory.

In regard to point (4), when one moves to the research literature, one can find expositions both of Laplace transforms for Quantum Calculus (Bohner and Guseinov 2010), as well as control theory on Time Scales. In fact, the control theory literature on Time Scales is rapidly growing, encompassing not just the generalization of linear control theory to time scales (Davis et al. 2009; Pawluszewicz and Torres 2010), but also the generalization of optimal control theory (Girejko et al. 2011) and even topics at the present frontier of the conventional control theory literature like control based on fractional calculus (Pooseh et al. 2013).